%% file: mt-arxiv.tex
\documentclass{mt}
\usepackage{textcomp}
\usepackage[Symbol]{upgreek}
\usepackage[all,arc,curve,color,frame,line]{xy}
\usepackage{amssymb}
\usepackage{epic}
\usepackage{mathrsfs} 
\usepackage{accents}
\usepackage{slashbox}
\usepackage{listliketab}
\usepackage{stmaryrd}
\usepackage{etoolbox}
\usepackage[refpage,notocbasic]{nomencl}
\usepackage{mathtools}
\usepackage{tikz}
\usetikzlibrary{decorations.pathreplacing}
\usetikzlibrary{arrows,matrix,positioning}
\usepackage{calligra}
\usepackage[T1]{fontenc}

\usepackage{hyperref}

\DeclareRobustCommand{\gobblefive}[5]{} % if not using hyperref, use 4 instead of 5
\newcommand*{\SkipTocEntry}{\addtocontents{toc}{\gobblefive}}

\newcommand\hidevp{$v_P$}

\makeindex
\makenomenclature
\renewcommand{\marginpar}[2][]{}

\DeclareMathAlphabet{\mathbf}{OML}{cmm}{b}{it}

\numberwithin{section}{chapter}
\numberwithin{equation}{section}
\numberwithin{figure}{chapter}
\numberwithin{table}{chapter}

\renewcommand{\theequation}{\thesection.\arabic{equation}}
\renewcommand\labelenumi{(\theenumi)}

\newcommand{\bbold}{\mathbb}
\newcommand{\cal}{\mathcal}

\newcommand{\rom}{\textup}

\renewcommand\epsilon{\varepsilon}

\def\R { {\bbold R} }
\def\Q { {\bbold Q} }
\def\Z { {\bbold Z} }
\def\C { {\bbold C} }
\def\N { {\bbold N} }
\def\E {{\mathcal E}}
\let\cedille\c
\def\c {\operatorname{c}}
\def\I {\operatorname{I}}
\def \G{G}

\def \sol{\operatorname{sol}}
\def \order{\operatorname{order}}
\def \val{\operatorname{mul}}

\def \exc {{\mathscr E}}
\def \ex{\operatorname{e}}
\def \wr {\operatorname{wr}}
\def \Wr {\operatorname{Wr}}
\def \Frac {\operatorname{Frac}}
\def \alg {{\operatorname{a}}}

\def \sch{\operatorname{s}}
\def \Sch{\operatorname{S}}
\def \ca{\mathcal}
\def \ac{\operatorname{ac}}
\def \af{\mathbf{f}}
\def \ag{\mathbf{g}}
\def \f{\operatorname{f}}
\let \accentv\v
\def \v{\operatorname{v}}
\def \re{\operatorname{r}}
\def \tp{\operatorname{tp}}

\def \a{\operatorname{a}}
\def \d{\operatorname{d}}

\def \ev{\operatorname{e}}
\def \bar {\overline}
\def \<{\langle}
\def \>{\rangle}

\def \tilde {\widetilde}

\def \Gd{\Q\Gamma}

\def \coker{\operatorname{coker}}
\def \ind{\operatorname{index}}

\def \hat {\widehat}

\def \supp {\operatorname{supp}}
\def \Spec {\operatorname{Spec}}

\def \RxLE {\R [[  x^{\R}  ]] ^{\operatorname {LE}}  }
\def \RxE {\R [[  x^{\R}  ]] ^{\operatorname {E}}  }
\def \((  {(\!(}
\def \)) {)\!)}
\def \T{\mathbb{T}}

\newcommand\psile{\triangledown}

\def \dd{\operatorname{ddeg}}
\def \dv{\operatorname{dmul}}

\def \Hom{\operatorname{Hom}}

\def \res{\operatorname{res}}

\def \k {{{\boldsymbol{k}}}}
\def \flatter{\mathrel{\prec\!\!\!\prec}}
\DeclareMathSymbol{\precequ}{\mathrel}{symbols}{"16}
\DeclareMathSymbol{\succequ}{\mathrel}{symbols}{"17}
\def \flattereq{\mathrel{\precequ\!\!\!\precequ}}

\def \comp{\mathrel{-{\hskip0.06em\!\!\!\!\!\asymp}}}

\def \nasymp{\not\asymp}

\newcommand{\claim}[2][\!\!]{\medskip\noindent {\sc Claim #1:} {\it #2}\medskip}
\newcommand{\claimnoskip}[2][\!\!]{\medskip\noindent\hskip\normalparindent {\sc Claim #1:} \textit{#2}\/}
\newcommand{\case}[2][\!\!]{\medskip\noindent {\sc Case #1:} {\it #2}\/}
\newcommand{\subcase}[2][\!\!]{\medskip\noindent {\sc Subcase #1:} {\it #2}\/}

\newtheorem{theorem}{Theorem}[section]
\newtheorem{lemma}[theorem]{Lemma}
\newtheorem{prop}[theorem]{Proposition}
\newtheorem{cor}[theorem]{Corollary}

\newtheorem{theoremintro}{Theorem}
\newtheorem{corintro}{Corollary}

\theoremstyle{definition}

\newtheorem{definition}[theorem]{Definition}

\theoremstyle{remark}
\newtheorem*{example}{Example}
\newtheorem*{examples}{Examples}
\newtheorem{exampleNumbered}[theorem]{Example}
\newtheorem{examplesNumbered}[theorem]{Examples}
\newtheorem*{exampleUnnumbered}{Example}
\newtheorem{exampleintro}{Example}
\newtheorem*{notation}{Notation}
\newtheorem*{notations}{Notations}
\newtheorem*{remarks}{Remarks}
\newtheorem*{remark}{Remark}

\newtheorem{remarkNumbered}[theorem]{Remark}

\newcommand{\abs}[1]{\lvert#1\rvert}
\newcommand{\dabs}[1]{\lVert#1\rVert}

\def \sgn {\operatorname{sign}}

\def \RGLE {\R[[ G^{\text{LE}}]] }

\def \fG{{\mathfrak G}}
\def \fH{{\mathfrak H}}
\def \fM {{\mathfrak M}}
\def \fN{{\mathfrak N}}
\def \fR{{\mathfrak R}}
\def \fd {{\mathfrak d}}
\def \fm {{\mathfrak m}}
\def \fn {{\mathfrak n}}
\def \fg{{\mathfrak g}}
\def \fr{{\mathfrak r}}
\def \fv {{\mathfrak v}}
\def \fw {{\mathfrak w}}

\def \CM {{{C}}[[ \fM ]] }
\def \CN {{{C}}[[ \fN ]] }

\def \RlognE {\R [[  \ell_n^\R  ]] ^{\operatorname {E}}  }
\def \RlognplusE {\R [[  \ell_{n+1}^\R  ]] ^{\operatorname {E}}  }

\def \fM {{\mathfrak M}}
\def \fd {{\mathfrak d}}
\def \fm {{\mathfrak m}}
\def \fn {{\mathfrak n}}
\def \fv {{\mathfrak v}}
\def \fw {{\mathfrak w}}

\def \CM {{{C}}[[ \fM ]] }

\def \Ric{\operatorname{Ri}}
\def \Cl{\operatorname{Cl}}
\def \id{\operatorname{id}}
\def \Log{\operatorname{L}}

\def \GL{\operatorname{GL}}
\def \Zero{\operatorname{Z}}
\def \Ideal{\operatorname{I}}

\let\oldi\i
\let\oldj\j
\renewcommand\i{\relax\ifmmode{\boldsymbol{i}}\else\oldi\fi}
\renewcommand\j{\relax\ifmmode{\boldsymbol{j}}\else\oldj\fi}

\def \bomega{{\boldsymbol{\omega}}}
\def \btau{{\boldsymbol{\tau}}}

\renewcommand\leq{\leqslant}
\renewcommand\geq{\geqslant}
\renewcommand\preceq{\preccurlyeq}
\renewcommand\succeq{\succcurlyeq}
\renewcommand\le{\leq}
\renewcommand\ge{\geq}
\renewcommand\frak{\mathfrak}

\DeclareFontFamily{U}{fsy}{}
\DeclareFontShape{U}{fsy}{m}{n}{<->s*[.9]psyr}{}
\DeclareSymbolFont{der@m}{U}{fsy}{m}{n}
\DeclareMathSymbol{\der}{\mathord}{der@m}{182}

\DeclareSymbolFont{der@m}{U}{fsy}{m}{n}
\DeclareMathSymbol{\derdelta}{\mathord}{der@m}{100}

\newcommand\restrict{\upharpoonright}

\newcommand\dotprec{\mathrel{\dot\prec}}
\newcommand\dotsucc{\mathrel{\dot\succ}}
\newcommand\dotpreceq{\mathrel{\dot\preceq}}
\newcommand\dotasymp{\mathrel{\dot\asymp}}

\makeatletter
\newcommand{\raisemath}[1]{\mathpalette{\raisem@th{#1}}}
\newcommand{\raisem@th}[3]{\raisebox{#1}{$#2#3$}}
\makeatother
\newcommand\simflat{\sim^{\raisemath{0.1em}{\flat}}}

\newcommand\wt{\operatorname{wt}}
\newcommand\wv{\operatorname{wm}}

\newcommand\bsigma{{\boldsymbol{\sigma}}}
\newcommand\m{\mathfrak m}

\let\polishl\l
\renewcommand\l{{\boldsymbol l}}

\newcommand\Ricmu{\mu}
\newcommand\Ricnu{\nu}

\newcommand\Pmu{\operatorname{dwm}}
\newcommand\Pnu{\operatorname{dwt}}
\newcommand\dwv{\operatorname{dwm}}
\newcommand\dwt{\operatorname{dwt}}
\newcommand\dval{\operatorname{dmul}}

\newcommand\ddeg{\operatorname{ddeg}}

\newcommand\nwt{\operatorname{nwt}}
\newcommand\ndeg{\operatorname{ndeg}}
\newcommand\nval{\operatorname{nmul}}

\newcommand\degree{\operatorname{d}}

\DeclareSymbolFont{imag@m}{OT1}{cmr}{m}{ui}
\DeclareMathSymbol{\imag}{\mathord}{imag@m}{105}

\DeclareFontFamily{OMS}{smallo}{}
\DeclareFontShape{OMS}{smallo}{m}{n}{<->s*[.65]cmsy10}{}
\DeclareSymbolFont{smallo@m}{OMS}{smallo}{m}{n}
\DeclareMathSymbol{\smallo}{\mathord}{smallo@m}{79}

\DeclareFontFamily{OMS}{largerdot}{}
\DeclareFontShape{OMS}{largerdot}{m}{n}{<->s*[.8]cmsy10}{}
\DeclareSymbolFont{largerdot@m}{OMS}{largerdot}{m}{n}
\DeclareMathSymbol{\largerdot}{\mathord}{largerdot@m}{15}

\newcommand{\tr}{\operatorname{\mathfrak{tr}}}

\newcommand\Der{\operatorname{der}}
\newcommand\TrDer{\operatorname{trder}}
\newcommand\TrAut{\operatorname{TrAut}}
\newcommand\TrEnd{\operatorname{tr}}
\newcommand\Aut{\operatorname{Aut}}

\newcommand\End{\operatorname{End}}
\newcommand\Lie{\operatorname{Lie}}

\newcommand\diag{\operatorname{diag}}
\newcommand\companion{{\operatorname{c}}}

\DeclareMathSymbol{\llambda}{\mathord}{der@m}{108}
\DeclareMathSymbol{\rrho}{\mathord}{der@m}{114}

\newcommand\lgb{{[\hspace{-0.24em}[\hspace{-0.14em}[\hspace{-0.24em}[}}
\newcommand\rgb{{]\hspace{-0.24em}]\hspace{-0.14em}]\hspace{-0.24em}]}}
\def \upg{\upgamma}
\def \Upg{\Upgamma}
\def \upl{\uplambda}
\def \Upl{\Uplambda}
\def \upo{\upomega}
\def \Upo{\Upomega}

\def \Upd{\Updelta}
\def \upu{\upgamma}
\def \Upu{\Upgamma}

\def\HL{\Upl}

\def\HLO{\Upl\Upo}

\newcommand{\Iota}{\mathrm{I}}

\newcommand\llb{\llbracket}
\newcommand\rrb{\rrbracket}
\newcommand\lla{\langle\hspace{-0.25em}\langle}
\newcommand\rra{\rangle\hspace{-0.25em}\rangle}
\newcommand\blla{\big\langle\hspace{-0.35em}\big\langle}
\newcommand\brra{\big\rangle\hspace{-0.35em}\big\rangle}

\newcommand{\equationqed}[1]{\[\pushQED{\qed}#1 \qedhere\popQED\]\let\qed\relax}
\newcommand{\alignqed}[1]{\begin{align*}\pushQED{\qed} #1 \qedhere\popQED\end{align*}\let\qed\relax}

\def \ome{\omega}

\newcommand*{\longhookrightarrow}{\ensuremath{\lhook\joinrel\relbar\joinrel\rightarrow}}
 
\makeatletter
\everydisplay=\expandafter{\the\everydisplay
  \ifdim\@totalleftmargin>\z@
    \displayindent=\dimexpr+\leftmargin\relax
    \displaywidth=\columnwidth
  \fi
}
\makeatother

\DeclareRobustCommand{\gobblefour}[4]{}

\newbool{PUP}
\boolfalse{PUP}

\begin{document}

\frontmatter
\title{Asymptotic Differential Algebra\\ and Model Theory of Transseries}

\author[Aschenbrenner]{Matthias Aschenbrenner}
\address{Department of Mathematics\\
University of California, Los Angeles\\
Los Angeles, CA 90095\\
U.S.A.}
\email{matthias@math.ucla.edu}

\author[van den Dries]{Lou van den Dries}
\address{Department of Mathematics\\
University of Illinois at Urbana-Cham\-paign\\
Urbana, IL 61801\\
U.S.A.}
\email{vddries@math.uiuc.edu}

\author[van der Hoeven]{Joris van der Hoeven}
\address{\'Ecole Polytechnique\\
91128 Palaiseau Cedex\\
France}
\email{vdhoeven@lix.polytechnique.fr}

%\begin{abstract}
%Version of \today. \dots
%\end{abstract}

\date{\today.}

\maketitle

\include{mt-quotes}

\tableofcontents

\include{mt-in}

\include{mt-1}

\include{mt-2}

\include{mt-3}

\include{mt-4}

\include{mt-5}

\include{mt-6}

\include{mt-7}

\include{mt-8}

\include{mt-9n}

\include{mt-10}

\include{mt-11}

\include{mt-12}

\include{mt-13}

\include{mt-14}

\include{mt-15}

\include{mt-16}

\include{mt-trans}

\include{mt-modth}

\backmatter

\include{mt-bib}

\include{mt-index}

\end{document}

%% file: mt-quotes.tex
\ifbool{PUP}{

\begin{bookepigraph}
Had the apparatus [of transseries and analyzable functions] been introduced for the sole purpose of solving Dulac's ``conjecture,'' one might legitimately question the wisdom and cost-effectiveness of such massive investment in new machinery. However, [these notions] have many more applications, actual or potential, especially in the study of analytic singularities. But their chief attraction is perhaps that of giving concrete, if partial, shape to G.~H.~Hardy's dream of an \textit{all-inclusive, maximally stable algebra of ``totally formalizable functions.''} \\[-1.5em]
\epigraphsource{{\small ---~Jean \'Ecalle, \textit{Six Lectures on Transseries, Analysable Functions and the Constructive Proof of Dulac's Conjecture}.}}
\end{bookepigraph}

\vskip-3.5em
 
\begin{bookepigraph}
The virtue of model theory is its ability to organize succinctly the sort of tiresome algebraic details associated with elimination theory.\\[-1.5em]
\epigraphsource{\small ---~Gerald Sacks, \textit{The Differential Closure of a Differential Field}.}
\end{bookepigraph}

\vskip-3.5em

\begin{bookepigraph}
Les analystes $p$-adiques se fichent tout autant que les g\'eom\`etres alg\'e\-bristes~\dots, des gammes \`a plus soif sur les valuations compos\'ees, les groupes ordonn\'es baroques, sous-groupes pleins desdits et que sais-je. Ces gammes m\'eritent tout au plus d'enrichir les exercices de Bourbaki, tant que personne ne s'en sert. \\[-1.5em]
\epigraphsource{{\small ---~Alexander Grothendieck, letter to Serre dated October 31, 1961.}}
\end{bookepigraph}

\vskip-3.5em

\begin{bookepigraph}
I don't like either writing or reading two-hundred page papers. It's not my idea of fun.\\[-1.5em]
\epigraphsource{{\small ---~John H.~Conway, quoted in \textit{Genius at Play: The Curious Mind of John Hor\-ton Con\-way}\/ by Siobhan Roberts.}}
\end{bookepigraph}

}{

\thispagestyle{empty}
\vspace*{\fill} 

\begin{quote}
Had the apparatus [of transseries and analyzable functions] been introduced for the sole purpose of solving Dulac's ``conjecture,'' one might legitimately question the wisdom and cost-effectiveness of such massive investment in new machinery. However, [these notions] have many more applications, actual or potential, especially in the study of analytic singularities. But their chief attraction is perhaps that of giving concrete, if partial, shape to G.~H.~Hardy's dream of an {\it all-inclusive, maximally stable algebra of ``totally formalizable functions.''}\\

{\small ---~Jean \'Ecalle, {\it Six Lectures on Transseries, Analysable Functions and the Constructive Proof of Dulac's Conjecture}.}
\end{quote}

\vskip3em

\begin{quote}
The virtue of model theory is its ability to organize succinctly the sort of tiresome algebraic details associated with elimination theory.\\

{\small ---~Gerald Sacks, {\it The Differential Closure of a Differential Field}.}
\end{quote}

\vskip3em

\begin{quote}
Les analystes $p$-adiques se fichent tout autant que les g\'eom\`etres alg\'e\-bristes~\dots, des gammes \`a plus soif sur les valuations compos\'ees, les groupes ordonn\'es baroques, sous-groupes pleins desdits et que sais-je. Ces gammes m\'eritent tout au plus d'enrichir les exercices de Bourbaki, tant que personne ne s'en sert. \\

{\small ---~Alexander Grothendieck, letter to Serre dated October 31, 1961.}
\end{quote}

\vskip3em

\begin{quote}
I don't like either writing or reading two-hundred page papers. It's not my idea of fun.\\

{\small ---~John H.~Conway, quoted in {\it Genius at Play: The Curious Mind of John Hor\-ton Conway}\/ by Siobhan Roberts.}
\end{quote}

\vspace*{\fill} 

}

%% file: mt-in.tex
\ifbool{PUP}{\begin{thepreface}}{\chapter*{Preface}}

\noindent
We develop here the algebra
and model theory of the \textit{differential field of transseries,}\/
a fascinating mathematical structure obtained
by iterating a construction going back more than a century to Levi-Civita and Hahn. It was introduced about thirty years ago as an
exponential ordered field by Dahn and G\"oring in connection with Tarski's problem on the real field with exponentiation, and independently by \'Ecalle in his proof of the Dulac Conjecture on plane analytic vector fields. 
%It has since found further use in model theory \cite{DMM1} and analysis \cite{Costin}.
% with precursors in the generalized power series fields of Levi-Civita~\cite{Levi-Civita} and Hahn~\cite{Hahn}. It
%was introduced  by Dahn and G\"oring~\cite{DG} in connection with Tarski's problem on the real field with exponentiation, and independently by \'Ecalle~\cite{Ecalle1} in connection with the Dulac Conjecture on plane analytic vector fields. 
%It has since found further use in model theory \cite{DMM1} and analysis \cite{Costin}.
%Transseries can represent  
%real-valued functions, as Laurent series represent (germs of) meromorphic functions. 
%Transseries are also rooted
%in the asymptotic expansions of Poincar\'e and Stieltjes, and 

The analytic aspects of transseries have a precursor in Borel's summation of divergent series.
%series~\cite{Poincare86, Stieltjes86, Borel}.
Indeed, \'Ecalle's theory of accelero-summation vastly extends Borel summation, and associates to each {\em accelero-summable\/} transseries an
{\em analyzable\/} function. In this way many non-oscillating
real-valued functions that arise naturally (for example, as solutions of algebraic differential equations) can be represented faithfully by transseries.

For about twenty years we have studied the differential field of trans\-series within the broader program of developing 
\textit{asymptotic differential algebra.}\/ We have recently obtained  decisive positive results 
on its model theory, and we describe these results in an {\em Introduction and
Overview.}\/ That introduction assumes some rudimentary knowledge
of differential fields, valued fields, and model theory, but
no acquaintance with transseries.
It is intended to familiarize readers 
%Assuming only rudimentary knowledge
%of valued fields, differential fields, and model theory, but
%no acquaintance with transseries, we describe these results %below in an ``Introduction and
%Overview.'' An attentive reading of this will make the reader familiar 
with the main issues in this book and with the
terminology that we frequently~use.

Initially, Joris van der Hoeven in Paris and
Matthias Aschenbrenner and Lou van den Dries in Urbana on the other side of the Atlantic worked independently, but around 2000 we decided to join forces. In 2011 we arrived at a rough outline for
proving some precise conjectures: see our programmatic survey
{\em Toward a model theory for transseries.}\/
All the conjectures stated in that paper (with one minor change) did turn out to be true,
even though some seemed to us at the time rather optimistic. 

Why is this book so long? For one, several problems we faced had no short solutions. Also, we have chosen to work in a setting that is sufficiently flexible for further
developments, as we plan to show in a later volume. Finally, we have tried to be reasonably self-contained by assuming
only a working knowledge of basic algebra: groups, rings, modules, fields. 
Occasionally we refer to Lang's~{\em Algebra}.

%The length of this volume is partly explained by the 
%inherent difficulties that had to be overcome. It is 
%also due to our choice of working 
%in a setting that is sufficiently flexible for further
%developments, and to our wish to be reasonably self-%contained. We assume
%only a working knowledge of basic algebra: groups, rings, %modules, fields. 
%Occasionally we refer to Lang's textbook \cite{Lang}. 

\ifbool{PUP}{\enlargethispage*{2\baselineskip}}{}

After the {\em Introduction and Overview} this book consists of~16 chapters and~2 appendices. Each chapter has an introduction and is divided into
sections. 
%The sections within a chapter are numbered. 
Each section has subsections, the last one often consisting of (partly historical) notes and comments. Many chapters state in the beginning some assumptions---sometimes just notational in nature--- that are in force throughout
that chapter, and of course the reader should be aware of
those in studying a particular chapter, since we do not repeat these assumptions when stating theorems, etc. The same holds for many sections and subsections.
The end of the volume has a list of symbols and an index.

\ifbool{PUP}{\section*{Acknowledgments}
}{}
%}{\skiptocentry\section*{Acknowledgments}} in earlier version submitted to PUP

\noindent
Part of this work was carried out while some of the authors were in residence at various times at the
Fields Institute (Toronto), the
Institut des Hautes \'Etudes Scientifiques (Bures-sur-Yvette), the
Isaac Newton Institute for Mathematical Sciences (Cambridge), and the
Mathematical Sciences Research Institute (Berkeley). The support and hospitality of these institutions is gratefully acknowledged.

Aschenbrenner's work was partially supported by  the National Science Foundation under grants DMS-0303618, DMS-0556197, and DMS-0969642. Visits by van der Hoeven to Los Angeles were partially supported by the UCLA Logic~Center.

We thank the following copyright holders for permission to reproduce the text in the epigraphs in the front of this book: Springer Science and Business Media, New York, for the quote by Jean \'Ecalle from \cite{Ecalle2}; the American Mathematical Society  for the quote by Gerald Sacks from \cite{SacksDC},
\copyright~1972 American Mathematical Society;
Professor Jean-Pierre Serre for the quote by Alexander Grothendieck from~\cite{GrothendieckSerre}; and Siobhan Roberts for the quote by John~H. Conway that appears in her book \textit{Genius at Play: The Curious Mind of John Horton Conway}\/~\cite{Roberts} \textcopyright~Siobhan Roberts, published by Bloomsbury Publishing, Inc., 2016.

We thank David Marker and Angus Macintyre for their interest and steadfast moral support over the years. To Santiago Camacho, Andrei Gabrielov, Tigran Hakob\-yan, Elliot Kaplan, Nigel Pynn-Coates, Chieu Minh Tran, and especially to Allen Gehret, we are indebted for numerous comments on and corrections to the 
manuscript. We are also grateful to Philip Ehrlich for
setting us right on some historical points, and to
the anonymous reviewers for useful suggestions and for spotting some
errors. We are of course solely responsible for any remaining inadequacies.    

Finally, we thank our editor, Vickie Kearn, and the other staff at Princeton University Press, notably Nathan Carr and Glenda Krupa, for helping us to bring this  book into its final form.

\bigskip

\bigskip

\hskip15em \textsl{Matthias Aschenbrenner,} Los Angeles  

\medskip

\hskip15em \textsl{Lou van den Dries,} Urbana  

\medskip

\hskip15em \textsl{Joris van der Hoeven,} Paris 

\bigskip

\hskip15em September 2015  

\ifbool{PUP}{\end{thepreface}}{}

\chapter*{Conventions and Notations}

\ifbool{PUP}{\addcontentsline{toc}{chapter}{Conventions and Notations}}{}

\noindent
Throughout, $m$ and $n$ range over the set $\N=\{0,1,2,\dots\}$ of
natural numbers. For sets~$X$,~$Y$ we distinguish between
$X\subseteq Y$, meaning that $X$ is a subset of $Y$, and~$X\subset Y$, meaning that $X$ is a proper subset of $Y$.

\nomenclature[A]{$A^{\ne}$}{$A\setminus\{0\}$, for an additively written abelian group $A$}
\nomenclature[A]{$R^\times$}{group of units of a ring $R$}

For an (additively written) abelian group $A$ we set
$A^{\ne}:= A\setminus \{0\}$. By {\em ring\/} we mean an
associative but possibly non-commutative ring with identity $1$. Let~$R$ be a ring. A {\em unit\/} of $R$ is a $u\in R$ with a right-inverse (an $x\in R$ with $ux=1$) and a left-inverse (an $x\in R$ with $xu=1$). If $u$ is a unit of $R$, then $u$ has only 
one right-inverse and only one left-inverse, and these coincide. With respect to multiplication
the units of $R$ form a group 
$R^\times$ with identity $1$. Thus the multiplicative group of
a field~$K$ is
$K^\times= K\setminus\{0\}=K^{\ne}$.  
Subrings and ring morphisms preserve $1$. 

\index{ring}
\index{ring!unit}
\index{unit}
\index{domain}
\index{domain!integral}
\index{integral!domain}

A {\em domain\/} is a ring with $1\ne 0$ such that for all $x$,~$y$ in the ring, if $xy=0$, then~$x=0$ or $y=0$. Usually domains are commutative, but not always. 
However, an {\em integral\/} domain is always commutative, that is, a subring of a field. 

Let $R$ be a ring. An $R$-module is a left $R$-module unless
specified otherwise, and the scalar $1\in R$ acts as the identity on any $R$-module. Let $M$ be an $R$-module and~$(x_i)_{i\in I}$ a family in $M$. A family $(r_i)_{i\in I}$ in $R$ is {\em admissible\/} if $r_i=0$ for all but finitely 
many~${i\in I}$. An
{\em $R$-linear combination of 
$(x_i)$\/} is an $x\in M$ such that $x=\sum_i r_ix_i$ of $M$ for some admissible family $(r_i)$ in $R$. We say that
{\em $(x_i)$  generates~$M$}\/ if every element of $M$ 
is an $R$-linear combination of $(x_i)$. We say that
$(x_i)$ is {\em $R$-dependent} (or {\em linearly dependent over~$R$}) if $\sum_i r_i x_i=0$ for some admissible
family $(r_i)_{i\in I}$ in $R$ 
with~${r_i\ne 0}$ for some~${i\in I}$; for $I=\{1,\dots,n\}$ we also abuse language by expressing this as:  {\em $x_1,\dots, x_n$ are 
$R$-dependent}. We say that~$(x_i)$ is {\em $R$-independent}
(or {\em linearly independent over~$R$}) if~$(x_i)$ is not $R$-dependent. We call~{\em $M$ free on $(x_i)$} (or {\em $(x_i)$ a basis of $M$}) if $(x_i)$ generates~$M$ and $(x_i)$ is $R$-independent. Sometimes we use this terminology for sets $X\subseteq M$ to mean that for some
(equivalently, for every)
index set~$I$ and bijection $i\mapsto x_i\colon I \to X$ the family~$(x_i)$ has the corresponding property.

\index{linear combination}
\index{independent!linearly}
\index{dependent!linearly}
\index{admissible!family}
\index{family!admissible}
\index{linearly!independent}
\index{linearly!dependent}
\index{R-dependent@$R$-dependent}
\index{R-independent@$R$-independent}
\index{module!free}

Let $K$ be a commutative ring. A {\em $K$-algebra\/} is defined to be a ring $A$ together with a ring morphism $\phi\colon K \to A$ that takes its values
in the center of $A$; we then refer to~$\phi$ as the structural morphism of the $K$-algebra $A$, and construe
$A$ as a $K$-module by $\lambda a:= \phi(\lambda) a$ for 
$\lambda\in K$ and $a\in A$. 

\index{algebra}

 Given a field extension $F$ of a field $K$ and a family 
 $(x_i)$ in $F$ we use the expressions~{\em $(x_i)$ is algebraically \textup{(}in\textup{)}dependent over $K$} and {\em $(x_i)$ is a transcendence basis of~$F$ over $K$} in a way similar to the
above linear analogues; likewise, a set~${X\subseteq F}$ can be
referred to as being a transcendence basis of $F$ over $K$. 
   
\index{algebraically!independent}
\index{algebraically!dependent}
\index{independent!algebraically}
\index{dependent!algebraically}
\index{transcendence!basis}

When a vector space $V$ over a field $C$ is given, then \textit{subspace of $V$}\/ 
means \textit{vector subspace of $V$.}\/

\begingroup

\ifbool{PUP}{

\clearpage
\mbox{}
\thispagestyle{empty}
\clearpage

\chapter*{Leitfaden}
\addcontentsline{toc}{chapter}{Leitfaden}

\def\leitfaden#1#2{ {\parbox{#1}{\rightskip=0pt plus .2\hsize
\footnotesize #2}}}

\begin{center}
\vskip-1.5em
\begin{tikzpicture}[->,>=stealth',shorten >=0.2em,shorten <=0.2em]
  \node[draw] (2) at (0,0) {\leitfaden{11em}{2.~Valued Abelian Groups}};
      \node[draw, below left = 7.5em of 2,xshift=3.5em] (1)  {\leitfaden{7.25em}{1.~Some Commutative\\ \phantom{1.~}Algebra}};
  \node[draw, below = 2em of 2] (3)  {\leitfaden{11em}{3.~Valued Fields}} edge[semithick,<-] (2) edge[semithick,<-, bend left=-18] (1);
  \node[draw, below = 2em of 3] (4)  {\leitfaden{11em}{4.~Differential Polynomials}}   edge[semithick,<-] (3)  edge[semithick,<-] (1);
  \node[draw, below = 2em of 4] (5)  {\leitfaden{11em}{5.~Linear Differential Polynomials}}  edge[semithick,<-] (4)  edge[semithick,<-, bend right=-13] (1);
  \node[draw, below = 2em of 5] (6)  {\leitfaden{11em}{6.~Valued Differential Fields}}  edge[semithick,<-] (5);
  \node[draw, below = 2em of 6] (7)  {\leitfaden{11em}{7.~Differential-Henselian Fields}}  edge[semithick,<-] (6);
  \node[draw, right = 2em of 7] (8)  {\leitfaden{10.25em}{8.~Differential-Henselian Fields\\ \phantom{8.}~with Many Constants}}  edge[semithick,<-] (7);
  \node[draw, below = 2em of 7] (9)  {\leitfaden{11em}{9.~Asymptotic Fields and \\ \phantom{9.~}Asymptotic Couples}}  edge[semithick,<-] (7);
  \node[draw, below = 2em of 9] (10)  {\leitfaden{11em}{10.~$H$-Fields}}  edge[semithick,<-] (9);

  \node[draw, below = 2em of 10] (11) {\leitfaden{11em}{11.~Eventual Quantities,\\ \phantom{11.}~Immediate Extensions, \\ \phantom{11.}~and Special Cuts}}  edge[semithick,<-] (10);
  \node[draw, below left = 2.8em of 11] (12) {\leitfaden{6.5em}{12.~Triangular\\ \phantom{12.~}Automorphisms}};
  \node[draw, below = 2em of 11] (13) {\leitfaden{11em}{13.~The Newton Polynomial}}  edge[semithick,<-] (11)  edge[semithick,<-] (12);
  \node[draw, below = 2em of 13] (14) {\leitfaden{11em}{14.~Newtonian Differential Fields}}  edge[semithick,<-] (13);
  \node[draw, below right = 3em and -4em of 14] (16) {\leitfaden{8.5em}{16.~Quantifier Elimination}}  edge[semithick,<-, bend left=-13] (14);
  
  \node[draw, below left = 3em and -6em of 14] (15) {\leitfaden{12em}{15.~Newtonianity of Directed Unions}}  edge[semithick,<-, bend right=-14] (14);

\end{tikzpicture}
\end{center}

\vfill
}
{

\chapter*{Leitfaden}

\def\leitfaden#1#2{ {\parbox{#1}{\rightskip=0pt plus .2\hsize
\small #2}}}

\vfill

\begin{center}
\hspace*{-6em}
\begin{tikzpicture}[->,>=stealth',shorten >=0.2em,shorten <=0.2em]
  \node[draw] (2) at (0,0) {\leitfaden{14.5em}{2.~Valued Abelian Groups}};
      \node[draw, below left = 7.5em of 2,xshift=2.5em] (1)  {\leitfaden{13em}{1.~Some Commutative Algebra}};
  \node[draw, below = 2em of 2] (3)  {\leitfaden{14.5em}{3.~Valued Fields}} edge[semithick,<-] (2) edge[semithick,<-, bend left=-11] (1);
  \node[draw, below = 2em of 3] (4)  {\leitfaden{14.5em}{4.~Differential Polynomials}}   edge[semithick,<-] (3)  edge[semithick,<-] (1);
  \node[draw, below = 2em of 4] (5)  {\leitfaden{14.5em}{5.~Linear Differential Polynomials}}  edge[semithick,<-] (4)  edge[semithick,<-, bend right=-13] (1);
  \node[draw, below = 2em of 5] (6)  {\leitfaden{14.5em}{6.~Valued Differential Fields}}  edge[semithick,<-] (5);
  \node[draw, below = 2em of 6] (7)  {\leitfaden{14.5em}{7.~Differential-Henselian Fields}}  edge[semithick,<-] (6);
  \node[draw, right = 2em of 7] (8)  {\leitfaden{12.5em}{8.~Differential-Henselian Fields\\ \phantom{8.}~with Many Constants}}  edge[semithick,<-] (7);
  \node[draw, below = 2em of 7] (9)  {\leitfaden{14.5em}{9.~Asymptotic Fields and \\ \phantom{9.~}Asymptotic Couples}}  edge[semithick,<-] (7);
  \node[draw, below = 2em of 9] (10)  {\leitfaden{14.5em}{10.~$H$-Fields}}  edge[semithick,<-] (9);

  \node[draw, below = 2em of 10] (11) {\leitfaden{14.5em}{11.~Eventual Quantities, Immediate \phantom{11.}~Extensions,  and Special Cuts}}  edge[semithick,<-] (10);
  \node[draw, below left = 2.8em of 11] (12) {\leitfaden{13em}{12.~Triangular Automorphisms}};
  \node[draw, below = 2em of 11] (13) {\leitfaden{14.5em}{13.~The Newton Polynomial}}  edge[semithick,<-] (11)  edge[semithick,<-] (12);
  \node[draw, below = 2em of 13] (14) {\leitfaden{14.5em}{14.~Newtonian Differential Fields}}  edge[semithick,<-] (13);
  \node[draw, below right = 3em and -5em of 14] (16) {\leitfaden{11em}{16.~Quantifier Elimination}}  edge[semithick,<-, bend left=-13] (14);
  
  \node[draw, below left = 3em and -7em of 14] (15) {\leitfaden{15em}{15.~Newtonianity of Directed Unions}}  edge[semithick,<-, bend right=-14] (14);

\end{tikzpicture}
\end{center}

\vskip2em

\vfill

}

\endgroup

\ifbool{PUP}{

\clearpage
\mbox{}
\thispagestyle{empty}
\clearpage

\chapter*[Dramatis Person\ae]{\fontsize{1.2em}{1em}\calligra\pmb{Dramatis Person\ae}}
\addcontentsline{toc}{chapter}{Dramatis Person\ae} 
}{
\chapter*{{\for{toc}{Dramatis Person\ae}}{\except{toc}{\fontsize{1.2em}{1em}\calligra Dramatis Person\ae}}} 
\markboth{Dramatis Person\ae}{Dramatis Person\ae}

}

\noindent
We summarize here the definitions of some notions prominent in our work, together with a list of attributes that apply to them. We include the 
page number where each concept is first introduced. We let 
$m$,~$n$,~$r$ range over~$\N=\{0,1,2,\dots\}$.
Below~$K$ is a field, possibly equipped with further structure. We let $a$,~$b$,~$f$,~$g$,~$y$,~$z$ range over elements of $K$, and we let~$Y$ be an indeterminate over $K$.
If~$K$ comes equipped with a valuation, then we let~$\mathcal O$ be the valuation ring of $K$, and we
freely employ the dominance relations on~$K$
introduced in Section~\ref{sec:valued fields}. If $K$ comes equipped with a derivation~$\der$, then we also write~$f',f'',\dots,f^{(n)},\dots$ for $\der f,\der^2 f,\dots,\der^n f,\dots$, and~$f^\dagger=f'/f$ for the logarithmic derivative of any $f\ne 0$; in this case~${C=\{f:f'=0\}}$ denotes the
constant field of $K$, and $c$ ranges over $C$.
% and we assume  $Y$ is a {\it differential}\/ indeterminate over $K$. 
We let $\Gamma$ be an ordered abelian group, and let $\alpha$,~$\beta$,~$\gamma$ range over $\Gamma$.

%\ifbool{PUP}{}{}
\SkipTocEntry\section*{Valued Fields}
\chaptermark{Dramatis Person\ae}

\noindent
Let $K$ be a \textit{valued field,}\/ that is, a field equipped with a valuation on it; p.~\pageref{p:valued fields}.

\medskip
\noindent
\textit{Complete}\/: every cauchy sequence in $K$ has a limit in $K$; p.~\pageref{p:complete}.

\medskip
\noindent
\textit{Spherically complete}\/:
every pseudocauchy sequence in $K$ has a pseudolimit in $K$; p.~\pageref{p:sph complete}. ``Spherically complete'' is equivalent to ``maximal'' as defined below.

%\medskip
%\noindent
%{\it Step-complete}\/:
%every special pseudocauchy sequence in $K$ has a pseudolimit in $K$. 
%(p.~\pageref{p:step-complete}.)

%\medskip
%\noindent
%{\it Fluent}\/:
%every fluent pseudocauchy sequence in $K$ has a pseudolimit in $K$. 
%(p.~\pageref{p:fluent}.)

\medskip
\noindent
\textit{Maximal}\/: there is no proper immediate valued field extension of $K$; p.~\pageref{p:maximal}.

\medskip
\noindent
\textit{Algebraically maximal}\/: there is no proper immediate algebraic valued field extension of $K$; p.~\pageref{p:algebraically maximal}.

\medskip
\noindent
\textit{Henselian}\/: for every $P\in K[Y]$ with $P\preceq 1$, $P(0)\prec 1$, and $P'(0)\asymp 1$, there exists~$y\prec 1$ with $P(y)=0$; p.~\pageref{sec:henselian valued fields}.

%\ifbool{PUP}{}{}
\SkipTocEntry\section*{Differential Fields}
\chaptermark{Dramatis Person\ae}

\noindent
Let $K$ be a \textit{differential field}, that is, a field of characteristic zero equipped with a derivation on it; p.~\pageref{p:differential field}.

%\medskip
%\noindent
%{\it Weakly differentially closed}\/:
%for all $P\in K[Y,\dots, Y^{(r)}]\setminus K$ there exists
%$y$ such that $P(y,y',\dots,y^{(r)})=0$. (p.~\pageref{p:weakly %differentially closed}.)

\medskip
\noindent
\textit{Linearly surjective}\/: for all $a_0,\dots,a_r\in K$ with $a_r\neq 0$
there exists $y$ such that $$a_0y+a_1y'+\cdots+a_ry^{(r)}=1;\qquad\text{ 
p.~\pageref{p:linearly surjective}.}$$

\medskip
\noindent
\textit{Linearly closed}\/:
for all $r\ge 1$ and $a_0,\dots,a_r\in K$ there are $b_0,\dots,b_{r-1},b\in K$ 
with $a_0Y+a_1Y'+\cdots+a_rY^{(r)} = 
b_0(Y'+bY)+b_1(Y'+bY)'+\cdots+b_{r-1}(Y'+bY)^{(r-1)}$; p.~\pageref{p:linearly closed}.

\medskip
\noindent
\textit{Picard-Vessiot closed \textup{(or \textit{pv-closed}\/)}}:\/
for all $r\ge 1$ and $a_0,\dots,a_r\in K$ with $a_r\neq 0$
there are $C$-linearly independent $y_1,\dots,y_r$ such that $a_0y_i+a_1y_i'+\cdots+a_ry_i^{(r)}=0$ for $i=1,\dots,r$;
p.~\pageref{p:pv-closed}.

\medskip
\noindent
\textit{Differentially closed}\/:
for all $P\in K[Y,\dots, Y^{(r)}]^{\ne}$ and
$Q\in K[Y,\dots, Y^{(r-1)}]^{\ne}$ such that $\frac{\partial P}{\partial Y^{(r)}}\neq 0$
there is $y$ with $P(y,y',\dots,y^{(r)})=0$ and $Q(y,y',\dots,y^{(r-1)})\neq 0$; p.~\pageref{p:differentially closed}.

%\ifbool{PUP}{}{}
\SkipTocEntry\section*{Valued Differential Fields}
\chaptermark{Dramatis Person\ae}

\noindent
Let $K$ be a \textit{valued differential field}, that is, a differential field 
equipped with a valuation on it; p.~\pageref{p:valued diff field}.

\medskip
\noindent
\textit{Small derivation}\/:
$f\prec 1 \Rightarrow f'\prec~1$; 
p.~\pageref{p:small}.

\medskip
\noindent
\textit{Monotone}\/: $f\prec 1 \Rightarrow f'\preceq f$; p.~\pageref{p:monotone}.

\medskip
\noindent
\textit{Few constants}\/: $c\preceq 1$ for all $c$;
p.~\pageref{p:many/few constants}.

\medskip
\noindent
\textit{Many constants}\/: for every
$f$ there exists~$c$ with $f\asymp c$; p.~\pageref{p:many/few constants}.

\medskip
\noindent
\textit{Differential-henselian \textup{(or \textit{$\d$-henselian}\/)}}\/:
$K$  has small derivation and: 
\begin{list}{*}{\setlength\leftmargin{2.5em}}
\item[(DH1)]
for all $a_0,\dots,a_r \preceq 1$ in $K$ with $a_r\asymp 1$
there exists $y\asymp 1$ such that $$a_0y+a_1y'+\cdots+a_ry^{(r)} \sim 1;$$
\item[(DH2)] for every $P\in  K[Y,Y',\dots,Y^{(r)}]$ with $P\preceq 1$, $P(0)\prec 1$, and
$\frac{\partial P}{\partial Y^{(n)}}(0)\asymp 1$ for some~$n$,  there exists $y\prec 1$ such that $P(y,y',\dots,y^{(r)})=0$; p.~\pageref{sec:dh1}.
\end{list}

%\ifbool{PUP}{}{}
\SkipTocEntry\section*{Asymptotic Fields}
\chaptermark{Dramatis Person\ae}

\noindent
Let $K$ be an \textit{asymptotic field}\/, that is, a valued differential field such that for all nonzero  $f,g\prec 1$: $f\prec g\Longleftrightarrow f'\prec g'$;
p.~\pageref{As-Fields,As-Couples}.

\medskip
\noindent
\textit{$H$-asymptotic \textup{(or \textit{of $H$-type}\/)}}\/:     $0\neq f\prec g\prec 1\Rightarrow f^\dagger\succeq g^\dagger$;
p.~\pageref{p:H-asymptotic field}.

\medskip
\noindent
\textit{Differential-valued \textup{(or \textit{$\d$-valued}\/)}}\/: for
 all $f\asymp 1$  
there exists $c$ with $f\sim c$;
p.~\pageref{p:dv}.

%\medskip
%\noindent
%{\it Pre-differential-valued \textup{(or {\it pre-$\d$-valued}\/)}}\/:  $f\preceq 1\ \&\ 0\neq g\prec 1\Rightarrow f'\prec g^\dagger$. (p.~\pageref{p:pdv}.)

\medskip
\noindent
\textit{Grounded}\/: there exists nonzero $f\nasymp 1$ such that 
$g^\dagger \succeq f^\dagger$ for all nonzero $g\nasymp 1$; p.~\pageref{p:grounded asymptotic field}.

\medskip
\noindent
\textit{Asymptotic integration}\/:
for all $f\neq 0$ there exists $g\nasymp 1$ with $g'\asymp f$;
p.~\pageref{p:asymptotic integration}.

%\medskip
%\noindent
%{\it Rational asymptotic integration}\/:
%for all $f\neq 0$ and $n\geq 1$ there exists $g\nasymp 1$  
%such that 
%$(g')^n \asymp f$.
%(p.~\pageref{p:rational as int}.)

{\sloppy

\medskip
\noindent
\textit{Asymptotically maximal}\/: $K$ has no proper immediate asymptotic field extension;
p.~\pageref{p:asymptotically maximal}.

}

\medskip
\noindent
\textit{Asymptotically $\d$-algebraically maximal}\/: $K$ has no proper immediate differential-al\-ge\-braic asymptotic field extension;
p.~\pageref{p:asymptotically d-algebraically maximal}.

\medskip
\noindent
\textit{$\upl$-free}\/: $H$-asymptotic, ungrounded, and for all $f$ there exists $g\succ 1$ with $f-g^{\dagger\dagger}\succeq g^\dagger$; p.~\pageref{p:upl-free}.

\medskip
\noindent
\textit{$\upo$-free}\/:
$H$-asymptotic, ungrounded, and for all $f$ there exists a $g\succ 1$ such that ${f-\omega(g^{\dagger\dagger})\succeq (g^\dagger)^2}$, where $\omega(z):=-(2z'+z^2)$; p.~\pageref{p:upo-free}.

\medskip
\noindent
\textit{Newtonian}\/: $H$-asymptotic, ungrounded, and every $P\in K[Y,Y',\dots,Y^{(r)}]^{\neq}$
of Newton degree~$1$
has a zero in $\mathcal O$; p.~\pageref{ch:newtonian fields}. (See p.~\pageref{p:ndeg} for Newton degree.)

%\ifbool{PUP}{}{}
\SkipTocEntry\section*{Ordered Valued Differential Fields}
\chaptermark{Dramatis Person\ae}

\noindent
Let $K$ be an \textit{ordered valued differential field\/}, that is, a valued 
differential field equipped with an ordering in the usual sense of 
\textit{ordered field}\/;
p.~\pageref{p:ovdf}.

\medskip
\noindent
\textit{Pre-$H$-field}\/: 
$\mathcal O$ is convex in the ordered field $K$, and
for all $f$:   
$$f>\mathcal O\Longrightarrow f'>0;\qquad  \text{p.~\pageref{p:pre-H-field}.}$$

\medskip
\noindent
\textit{$H$-field}\/: $\mathcal O$ is the convex hull of $C$ 
in the ordered field~$K$, and for all $f$: 
  $$f>C\Longrightarrow f'>0, \quad f\asymp 1\Longrightarrow \text{
there exists~$c$ with $f\sim c$;}\qquad \text{p.~\pageref{p:H-field}.}$$

\medskip
\noindent
\textit{Liouville closed}\/:
$K$ is a real closed $H$-field and
for all $f$,~$g$ there exists $y\neq 0$ such that $y'+fy=g$;
p.~\pageref{p:Liouville closed}.

%\medskip
%\noindent
%{\it Schwarz closed}\/:  
%$K$ is a Liouville closed, $K=\omega(K)\cup\sigma(K^\times)$,
%and $\omega(K)$ is closed downward, where $\omega(z):=-(2z'+z^2)$ and $\sigma(y):=\omega(-y^\dagger)+y^2$ for $y\neq 0$.
%(p.~\pageref{p:Schwarz closed}.)

%\ifbool{PUP}{}{}
\SkipTocEntry\section*{Asymptotic Couples}
\chaptermark{Dramatis Person\ae}

\noindent
Let $(\Gamma, \psi)$ be an \textit{asymptotic couple\/}, that is, the 
ordered abelian group $\Gamma$ is equipped with a map 
$\psi\colon\Gamma^{\neq}\to\Gamma$ such that 
for all $\alpha,\beta\neq 0$:
\begin{list}{*}{\setlength\leftmargin{2.5em}}
\item[(AC1)] $\alpha+\beta\ne 0\ \Rightarrow\ \psi(\alpha+\beta)\ge \min\!\big(\psi(\alpha), \psi(\beta)\big)$;
\item[(AC2)] $\psi(k\alpha)=\psi(\alpha)$  for all $k\in \Z^{\ne}$;
\item[(AC3)] $\alpha>0 \ \Rightarrow\ \alpha':=\alpha + \psi(\alpha) > \psi(\beta)$;
\end{list}
p.~\pageref{sec:ascouples}. For $\gamma\ne 0$ we set $\gamma':= \gamma + \psi(\gamma)$. 

\medskip
\noindent
\textit{$H$-asymptotic \textup{(or \textit{of $H$-type}\/)}}\/:  
$0<\alpha<\beta\ \Rightarrow\ \psi(\alpha)\geq\psi(\beta)$;
p.~\pageref{p:H-ascouples}.

\medskip
\noindent
\textit{Grounded}\/: $\Psi:= \{\psi(\alpha):\ \alpha\neq 0\}$ 
has a largest element; p.~\pageref{p:grounded as couple/small der}.

\medskip
\noindent
\textit{Small derivation}\/: $\gamma>0\Rightarrow\gamma'>0$; p.~\pageref{p:grounded as couple/small der}.

\medskip
\noindent
\textit{Asymptotic integration}\/:
for all $\alpha$ there exists $\beta\neq 0$ with 
$\alpha=\beta'$;
p.~\pageref{p:asymptotic integration}.

%\medskip
%\noindent
%{\it Rational asymptotic integration}\/:
%for all $\alpha$ and $n\geq 1$ there exists $\beta\neq 0$  with 
%$\alpha=n\,\beta'$.
%(p.~\pageref{p:rational as int}.)

\ifbool{PUP}{

\clearpage
\mbox{}
\thispagestyle{empty}
\clearpage 
 
\makehalftitle

%\mbox{}
%\chaptermark{}
%\clearpage

}{}

\mainmatter

\ifbool{PUP}{
\clearpage

\chapter*[Introduction and Overview]{Introduction and Overview}
\addcontentsline{toc}{chapter}{Introduction and Overview}\markboth{\MakeUppercase{Introduction and Overview}}{\MakeUppercase{Introduction and Overview}}
}{\chapter*{Introduction and Overview}}

\ifbool{PUP}{
\section*{A Differential Field with No Escape}
\addcontentsline{toc}{section}{A Differential Field with No Escape}
}{\section*{A Differential Field with No Escape}}

\noindent
Our principal object of interest is the \textit{differential field}\/ $\T$ of \textit{transseries}. 
Transseries are formal
series in an indeterminate $x>\R$, such as
\begin{align} \tag{1}
  \varphi(x)\ & =\  - 3 \ex^{\ex^x} + \ex^{\textstyle\frac{\ex^x}{\log x} +
  \frac{\ex^x}{\log^2 x} + \frac{\ex^x}{\log^3 x} + \cdots} - x^{11} + 7
  \label{intro-ex}\\
  &   \hspace{1em} + \frac{\pi}{x} + \frac{1}{x \log x} + \frac{1}{x \log^2
  x} + \frac{1}{x \log^3 x} + \cdots \nonumber\\
  &   \hspace{1em} + \frac{2}{x^2} + \frac{6}{x^3} + \frac{24}{x^4} +
  \frac{120}{x^5} + \frac{720}{x^6} + \cdots \nonumber\\
  &   \hspace{1em} + \ex^{- x} + 2 \ex^{- x^2} + 3 \ex^{- x^3} + 4 \ex^{-
  x^4} + \cdots, \nonumber
\end{align}
where $\log^2 x := (\log x)^2$, etc.
% $x \rightarrow \infty$. 
As in this example, each transseries is a (possibly transfinite) 
sum, with
terms written from left to right, in asymptotically decreasing order. Each
term is the product of a real coefficient and a \textit{transmonomial}\/. Appendix~\ref{app:trans} contains the
inductive construction of~$\T$, including the definition of ``transmonomial'' and other
notions about transseries that occur in this introduction. For expositions
of $\T$ with proofs, see~\cite{DMM2,Edgar,JvdH}. 
In~\cite{DMM2},~$\T$ is denoted by~$\R\(( x^{- 1}\)) ^{\operatorname{LE}}$, and its
elements are called \textit{logarithmic-exponential series}.
At this point we just mention that transseries can be
added and multiplied in the natural way,
and that with these operations, $\mathbb T$ is a field containing~$\R$ as a subfield.
Transseries can also be differentiated term by term, subject to $r'=0$ for each $r\in\R$ and $x'=1$. In this way~$\mathbb T$ acquires the structure
of a \textit{differential field.}\/ 
%(We invite the reader to write out the derivative~$\varphi'(x)$ of the transseries $\varphi(x)$ from \eqref{intro-ex} as a transseries.)

\subsection*{Why transseries?}
Transseries naturally arise in solving differential equations at infinity and studying the asymptotic behavior of
their solutions, where 
ordinary power series, Laurent series, or even
Puiseux series in $x^{-1}$ are inadequate. Indeed, functions as simple as~$\ex^x$ or~$\log x$
cannot be expanded with respect to the asymptotic scale $x^{\R}$ of
real powers of $x$ at~$+\infty$. For merely solving algebraic equations, no
exponentials or logarithms are needed: it is classical that the
fields of Puiseux series over~$\R$ and~$\C$ are 
real closed and algebraically closed, respectively.

One approach to asymptotics 
with respect to more general scales was initiated by Hardy~\cite{Hardy11,Hardy2},
inspired by earlier work of du~Bois-Reymond~\cite{duBoisReymond} in the late 19th century. Hardy considered 
\textit{logarithmico-exponential functions}\/: real-valued functions
 built up from constants and the variable
$x$ using addition, multiplication, division, exponentiation and taking
logarithms. 
He showed that such a function, when defined on some interval $(a, + \infty)$,  has eventually constant sign (no oscillation!), 
and so the germs at 
$+\infty$ of these functions form an ordered field $H$ with derivation $\frac{d}{dx}$.
Thus $H$ is what Bourbaki~\cite{Bourbaki} calls a \textit{Hardy field}\/:
a
sub\textit{field}\/~$K$ of the ring of germs at $+\infty$ of differentiable functions
$f\colon (a, +\infty)\to \R$ with $a\in \R$, closed under taking derivatives; for more precision, see Section~\ref{As-Fields,As-Couples}. 
Each Hardy field is naturally an ordered differential field.
The Hardy field $H$ is rather special: every $f\in H$ satisfies an algebraic differential equation over $\R$. But $H$ lacks some closure 
properties that are
desirable for a comprehensive theory. For instance, 
$H$ has no antiderivative of $\ex^{x^2}$ 
(by Liouville; see~\cite{Rosenlicht72}), and the functional inverse of $(\log x)( \log \log x)$ doesn't lie in
$H$, and is not even asymptotic 
to any element of~$H$:~\cite{DMM1,vdH:PhD}; see also~\cite{Ressayre00}.  

With $\T$ and transseries we go beyond $H$ and logarithmico-exponential functions by admitting \textit{infinite}\/ sums.
It is important to be aware, however, that by virtue of its inductive construction,
$\T$ does not contain, for example, the series
\[ x + \log x + \log \log x + \log \log \log x + \cdots, \]
which does make sense in a suitable extension of $\T$. Thus 
$\T$ allows only certain kinds of infinite sums.
Nevertheless, it turns out that
the  differential field~$\T$ enjoys many remarkable 
closure
properties that $H$ lacks. For instance, 
$\T$ is closed under natural operations of
exponentiation, integration, composition, compositional inversion, and 
the resolution of \textit{feasible}\/ algebraic differential equations
(where the meaning of \textit{feasible}\/ can be made explicit). 
This makes $\T$ of interest for different areas of mathematics:

\ifbool{PUP}{}{\medskip}

\subsubsection*{Analysis} In connection with the Dulac Problem,
  $\T$ is sufficiently rich for modeling the asymptotic behavior of
  so-called Poincar{\'e} return maps. This analytically deep result is a crucial part of {\'E}calle's solution of the Dulac
  Problem~{\cite{Ecalle0,Ecalle1,Ecalle2}}. (At the end of this introduction we discuss this in more detail.)
  
\ifbool{PUP}{}{\medskip}

\subsubsection*{Computer algebra} Many transseries are concrete enough to compute with them, in the 
  sense of computer algebra~\cite{vdH:PhD,Sh04}. Moreover, many of
  the closure properties mentioned above can be made effective. This allows for the automation
  of an important part of asymptotic calculus for functions of one variable.

\ifbool{PUP}{}{\medskip}
  
\subsubsection*{Logic} Given an o-minimal expansion  of the real field, the germs at $+\infty$ of its  definable one-variable functions 
form a Hardy field, which in many cases %(for example, for the exponential field of reals) 
can be embedded into $\T$. This gives useful information about
the possible asymptotic behavior of these definable functions;  
see \cite{AvdD4,Miller} for more about this connection.

\medskip
\noindent
Soon after the introduction of $\T$ in the 1980s  it was  suspected that $\T$ might well be a
kind of \textit{universal domain}\/ for the differential algebra of Hardy fields and
similar ordered differential fields,
analogous to the role of the algebraically closed field~$\C$ as a universal domain for algebraic geometry of characteristic $0$
(Weil~\cite[Chapter~X, \S{}2]{Weil}),  and of~$\R$,~$\Q_p$, and~$\C \(( t\)) $  in related ordered and valued settings. This is corroborated by the strong closure properties enjoyed by $\T$.
See in particular 
p.~148 of
\'Ecalle's book~\cite{Ecalle1} for eloquent expressions of this idea. The present volume and the next
substantiate the \textit{universal domain}\/ nature of the differential field~$\T$,
using the language of \textit{model theory.}\/
%Structures enjoying strong closure properties, such as $\T$, often also have  a nice model theory.
The model-theoretic properties of
the classical fields~$\C$,~$\R$,~$\Q_p$ and~$\C \(( t\)) $ are well established thanks to
Tarski, Seidenberg, Robinson, Ax~\&~Kochen, Er{\accentv s}ov, Cohen, Macintyre, Denef, 
and others; see~\cite{Ta,Se,Robinson56,AxKochen,AxKochen2,Ersov, Cohen, Mac, Denef}. 
Our goal is to analyze likewise the differential field~$\T$, which
comes with a definable ordering and valuation, and in this book we achieve
 this goal.

\subsection*{The ordered and valued differential field $\T$}

\noindent
For what follows, it will be convenient to quickly survey some of the most
distinctive features of $\T$.
Appendix~\ref{app:trans} contains precise definitions and further details.
%and Appendix~\ref{app:...} explains model-theoretic concepts.

\medskip
\noindent
Each transseries $f = f (x)$ can be uniquely decomposed as a sum
$$ f\ =\ f_{\succ} + f_{\asymp} + f_{\prec},$$
where $f_{\succ}$ is the \textit{infinite part of~$f$}, $f_{\asymp}$ is its
constant term (a real number), and $f_{\prec}$ is its \textit{infinitesimal
part}. In the example~\eqref{intro-ex} above,
\begin{align*}
  \varphi_{\succ} &\ = \ - 3 \ex^{\ex^x} + \ex^{\textstyle\frac{\ex^x}{\log x} +
  \frac{\ex^x}{\log^2 x} + \frac{\ex^x}{\log^3 x} + \cdots} - x^{11},\\
  \varphi_{\asymp} &\ =\  7,\\
  \varphi_{\prec} &\ =\  \frac{\pi}{x} + \frac{1}{x \log x} + \cdots .
\end{align*}
In this example, $\varphi_{\succ}$ happens to be a finite sum, but 
this is not a necessary feature of transseries: take for example  
$f:= \frac{\ex^x}{\log x} +
\frac{\ex^x}{\log^2 x} + \frac{\ex^x}{\log^3 x} + \cdots$, with
$f_{\succ}=f$.
Declaring a
transseries to be positive iff its dominant (=~leftmost)
coefficient is positive turns~$\T$ into an \textit{ordered}\/
field extension of~$\R$ with $x >\R$. In our example~\eqref{intro-ex}, the
dominant transmonomial of~$\varphi(x)$ is $\ex^{\ex^x}$ and its dominant
coefficient is~$-3$, whence~$\varphi(x)$ is negative; in fact, $\varphi(x)<\R$.

\medskip
\noindent
The inductive definition of $\T$ involves constructing a certain exponential
operation $\exp \colon \T \to \T^\times$, with $\exp (f)$ also written as $\ex^f$, and
\[ \exp (f)\ =\ \exp (f_{\succ}) \cdot \exp (f_{\asymp}) \cdot \exp (f_{\prec})\
   =\ \exp (f_{\succ}) \cdot \ex^{f_{\asymp}} \cdot \sum_{n = 0}^{\infty}
   \frac{f_{\prec}^n}{n!} \]
where the first factor $\exp (f_{\succ})$ is a transmonomial,
the second factor $\ex^{f_{\asymp}}$ is the real number obtained by exponentiating the real
number $f_{\asymp}$ in the usual way, and the third factor $\exp (f_{\prec}) =
\sum_{n = 0}^{\infty} \frac{f_{\prec}^n}{n!}$ is expanded as a series in the
usual way. Conversely, each transmonomial is of the form $\exp (f_{\succ})$
for some transseries $f$. 
Viewed as an exponential field, $\T$ is an {\em elementary\/} extension of the
exponential field of real numbers; see~\cite{DMM1}. 
In particular, $\T$ is
real closed, and so its ordering is existentially definable (and universally
definable) from its ring operations:
\begin{equation}\tag{2}
  f \geqslant 0\ \Longleftrightarrow\ f = g^2  \text{ for some $g$.}
  \label{square-equiv}
\end{equation}
However, as emphasized above, our main interest is in $\T$ as a \textit{differential field,} with
derivation $f \mapsto f'$ on $\T$ defined termwise, with $r' = 0$ for $r \in
\R$, $x' = 1$, $(\ex^f)' = f' \ex^f$, and $(\log f)' = f' / f$ for $f > 0$.
Let us fix here some notation and terminology in force throughout this volume:
a \textit{differential field} is a field $K$ of characteristic $0$ together
with a single derivation $\der \colon K \to K$; if $\der$ is clear from the context
we often write $a'$ instead of~$\der (a)$, for $a \in K$. The
\textit{constant field} of a differential field $K$ is the subfield
\[ C_K\ :=\  \{a \in K\ : a' = 0\} \]
of $K$, also denoted by $C$ if $K$ is clear from the context. 
The constant field of $\T$ turns out to be $\R$, that is,
\[ \R\ =\ \{f \in \T :\  f' = 0\} . \]
%and  infinite in size, since $\ex^{\ex^x}$ is infinite.)
By an
\textit{ordered differential field}\/ we mean a differential field equipped
with a total ordering on its underlying set making it an \textit{ordered
field}\/ in the usual sense of that expression. So~$\mathbb T$ is an \textit{ordered differential field.}\/
More important than the ordering is the \textit{valuation}\/ on $\T$ with valuation ring
\[ \mathcal{O}_{\T}\ :=\ \big\{ f \in \T :\  \text{$|f| \leqslant r$ for some  $r
   \in \R$} \big\}\ =\ \{f \in \T :\  f_{\succ} = 0\}, \]
a convex subring of $\T$. The unique maximal ideal of $\mathcal O_{\T}$ is
$$\smallo_{\T}\ :=\ \big\{ f\in \T: \text{$|f| \leqslant r$ for all  $r>0$ 
   in $\R$}\big\}\ =\ \{f \in \T :\  f=f_{\prec}\}$$
and thus $\mathcal O_{\T}=\R+\smallo_{\T}$.
Its very definition shows that $\mathcal{O}_{\T}$ is
\textit{existentially}\/ definable in the differential field~$\T$. However,
$\mathcal{O}_{\T}$ is \textit{not}\/ universally definable in the differential
field~$\T$: Corollary~\ref{notmodelcomplete}. 
In light of the model completeness conjecture discussed below, it is therefore advisable to add the
valuation as an extra primitive, and so in the rest of this introduction \textit{we construe~$\T$
as an ordered and valued differential field, with valuation given by~$\mathcal{O}_{\T}$.}\/ By a \textit{valued differential field}\/ we mean throughout a
differential field $K$ equipped with a valuation ring of $K$ that contains the
prime subfield $\Q$ of $K$.

\subsection*{Grid-based transseries}

When referring to transseries we have in
mind the \textit{well-based transseries of finite logarithmic and
exponential depth}\/ of \cite{vdH:PhD}, also called \textit{logarithmic-exponential series}\/ in~\cite{DMM2}. The construction of the field $\T$ in
Appendix~\ref{app:trans} allows variants, and 
we briefly comment on one of
them.

Each transseries $f$ is an infinite sum $f=\sum_{\fm} f_{\fm}\fm$
where each $\fm$ is a transmonomial and $f_{\fm}\in\R$.  
The \textit{support}\/ of such a transseries $f$ is the set $\operatorname{supp}(f)$ of transmonomials~
$\mathfrak{m}$ for which the coefficient $f_{\mathfrak{m}}$ is
nonzero. For instance, the transmonomials in the support of the transseries $\varphi$ of 
example~\eqref{intro-ex} are
\ifbool{PUP}{
\begin{multline*}
\ex^{\ex^x},\ \ex^{\textstyle\frac{\ex^x}{\log x} +
   \frac{\ex^x}{(\log x)^2} + \frac{\ex^x}{(\log x)^3} + \cdots},\ x^{11}, \\
   1,\ \frac{1}{x},\ \frac{1}{x \log x},\ \ldots,\ \frac{1}{x^2},\ \frac{1}{x^3},
   \ldots,\ \ex^{-x},\ \ex^{-x^2},\ \ldots.
\end{multline*}
}{
\[ \ex^{\ex^x}, \ex^{\textstyle\frac{\ex^x}{\log x} +
   \frac{\ex^x}{(\log x)^2} + \frac{\ex^x}{(\log x)^3} + \cdots}, x^{11}, 1,
   \frac{1}{x}, \frac{1}{x \log x}, \ldots, \frac{1}{x^2}, \frac{1}{x^3},
   \ldots, \ex^{-x}, \ex^{-x^2},\ldots. \]}
%The set of transmonomials is totally ordered by the asymptotic dominance  relation
%\[ \mathfrak{m} \preccurlyeq \mathfrak{n} \quad \Longleftrightarrow \quad    \mathfrak{m} \in    \mathcal{O} \mathfrak{n} \quad \Longleftrightarrow \quad   \text{$\mathfrak{m} \leq c\mathfrak{n}$ for some $c\in \R^{>0}$.} \]
By imposing various restrictions on the kinds of permissible supports, the
construction from Appendix~\ref{app:trans} yields various interesting differential subfields of
$\T$.

To define multiplication on $\T$, supports should be {\em{well-based}}:
every nonempty subset of the support of a 
transseries $f$
should contain an asymptotically dominant element. 
%Transseries with this kind of support are said to be {\em{well-based}}. 
So well-basedness is a minimal requirement on supports. A much stronger  
condition on $\supp(f)$ is as follows: there are transmonomials
$\mathfrak{m}$ and  
$\mathfrak{n}_1, \ldots, \mathfrak{n}_k \in\smallo_{\T}$ ($k\in \N$) such that
$$  \operatorname{supp} f \ \subseteq \ \big\{ \mathfrak{m}\, \mathfrak{n}_1^{i_1} \cdots
  \mathfrak{n}_k^{i_k}:\  i_1, \ldots, i_k \in \N \big\} .
$$
Supports of this kind are called {\em{grid-based}}. Imposing this constraint
all along, the construction from Appendix~\ref{app:trans} builds the
differential subfield $\T_{\operatorname{g}}$ of {\em{grid-based}} transseries of~$\T$.
Other suitable restrictions on the support yield other interesting differential 
subfields of $\T$.

The differential field $\T_{\operatorname{g}}$ of grid-based transseries has been studied in detail
in~{\cite{JvdH}}.
% (see also {\cite[{\textsection}7.7]{DMM2}}). 
In particular, that book contains a kind of algorithm for
solving {\em{algebraic differential equations}} over $\T_{\operatorname{g}}$.
These equations are of the form
\begin{equation}\tag{3}  P \big(y, \ldots, y^{(r)}\big) \ = \ 0,  \label{alg-diff-eq}
\end{equation}
where $P \in \T_{\operatorname{g}} [Y, \ldots, Y^{(r)}]$ is a nonzero polynomial in $Y$ and a
finite number of its formal derivatives $Y', \ldots, Y^{(r)}$. We note here that by combining results
from~\cite{JvdH} and the present volume, any  solution~$y \in
\T$ to~\eqref{alg-diff-eq} is actually grid-based.  
Thus transseries outside $\T_{\operatorname{g}}$ such as~$\varphi(x)$ from \eqref{intro-ex} or
$\zeta (x) =
1 + 2^{- x} + 3^{- x} + \cdots$
%$$\ =\ \sum_{n=1}^{\infty}\ex^{-x\log n},$$ 
are differentially transcendental
over~$\T_{\operatorname{g}}$; see the
\textit{Notes and comments}\/ to Section~\ref{mctnl} for more details,
and Grigor$'$ev-Singer~\cite{GrigorievSinger} for an earlier result of this kind.

%Which supports we allow is related to which infinite sums  are allowed to exist. A family $(f_i)$ of well-based transseries is said to be {\it summable} if $\bigcup_{i \in I} \operatorname{supp} f_i$ is well-based and for every transmonomial $\mathfrak{m}$ the set $\{ i \in I : (f_i)_{\mathfrak{m}}\neq 0 \}$ is finite. By~\cite{KKS} there is no field of well-based transseries closed under under exponentiation, taking logarithms of positive elements, as well as infinite sums. For instance, the field~$\T$ constructed in Appendix~\ref{app:transseries} is closed under exponentiation and taking logarithms of positive elements, but it does not contain series such as
%\[ x + \log x + \log \log x + \log \log \log x + \cdots, \]
%so it is not closed under infinite sums. Nevertheless, $\T$ is the increasing union of a sequence of differential subfields of transseries  closed under infinite sums.
%One may also construct a strictly 
%increasing ordinal sequence
%$\T_{\alpha}$ of fields of transseries which are 
%all closed under
%infinite summation and logarithm, and such that 
%$\exp \T_{\alpha}
%\subseteq \T_{\alpha + 1}$ for all $\alpha$.

\subsection*{Model completeness}
One reason that ``geometric'' fields like $\C$, $\R$, $\Q_p$ are more manageable than ``arithmetic'' fields like $\Q$ is that the
former are {\em model complete}; see Appendix~\ref{app:modth} 
for this and other basic model-theoretic notions used in this volume. 
A consequence of the model completeness of
$\R$ is that any finite system of polynomial equations over $\R$ (in any number of unknowns) with a solution
in an ordered field extension of $\R$, has a solution in $\R$ itself.
By the $\R$-version of \eqref{square-equiv} we can also allow polynomial inequalities in such a system. (A related fact: if such a system has real algebraic coefficients, then it has a
real algebraic solution.)

For a more geometric view of model completeness we first specify an algebraic subset of~$\R^n$ to be the set of common zeros,
\[ \big\{ y = (y_1, \ldots, y_n) \in \R^n :\  P_1 (y) = \cdots = P_k (y) =
   0 \big\}, \]
of finitely many polynomials $P_1, \ldots, P_k \in \R [Y_1, \ldots, Y_n]$. 
Define a subset of $\R^m$ to be {\em{subalgebraic}} if it is the image of
an algebraic set in $\R^n$ for some~$n \geqslant m$ under the
projection map
\[ (y_1, \ldots, y_n) \mapsto (y_1, \ldots, y_m)\ \colon\ \R^n \to
   \R^m . \]
Then a consequence of the model completeness of $\R$ is that the complement in
$\R^m$ of any subalgebraic set is again
subalgebraic. Model completeness of $\R$ is a little stronger in that only polynomials with integer coefficients should be involved.

A nice analogy between $\R$ and $\T$ is the 
following intermediate value property,
announced in~\cite{vdH:ivt} and established for~$\T_{\operatorname{g}}$
in~\cite{JvdH}: Let $P (Y) = p (Y, \ldots, Y^{(r)})$ be a differential
polynomial over~$\T$, that is, with coefficients in $\T$, and let $f$, $h$ be
transseries with $f < h$; then $P (g)$ takes on all values strictly between $P
(f)$ and $P (h)$ for transseries~$g$ with $f < g < h$.
Underlying this opulence of $\T$ is a more robust 
property that we call \textit{newtonianity,}\/ which is
analogous to henselianity for valued fields.
The fact that $\T$ is newtonian implies, for instance, that any
differential equation $y'=Q(y, y', \dots, y^{(r)})$ with 
$Q\in x^{-2}\mathcal{O}_{\T}[Y, Y',\dots, Y^{(r)}]$
%\[ y' + 3 xy\ =\ \ex^{- \ex^x}  
%(y^2 + y' y'' - y'' y''') - \ex^{-
%   x^{15}} \]
has an infinitesimal solution $y\in\smallo_{\T}$. The
definition of ``newtonian'' is rather subtle, and is discussed later in
this introduction.

Another way that $\R$ and $\T$ are similar concerns
the factorization of linear differential operators: 
any linear differential operator $A = \der^r + a_1\der^{r-1}\cdots + a_r$ of order $r\ge 1$ with coefficients $a_1,\dots,a_r\in\T$, is a product of such operators of order one and order two, with coefficients in $\T$.  Moreover, any
linear differential equation $y^{(r)}+a_1y^{(r-1)}+\cdots+a_ry=b$ ($a_1,\dots,a_r,b\in\T$) has a solution $y\in\T$
(possibly $y=0$).
In particular, every transseries $f$ has a transseries integral $g$, that is,~$f = g'$.
(It is noteworthy that a \textit{convergent}\/ transseries can very well have a
\textit{divergent}\/ transseries as an integral; for example, the transmonomial~$\frac{\ex^x}{x}$ has as an integral the divergent transseries $\sum_{n =
0}^{\infty} n! \frac{\ex^x}{x^{n + 1}}$. The analytic aspects of transseries
are addressed by {\'E}calle's theory of {\em{analyzable
functions}}~{\cite{Ecalle1}}, where
genuine functions are associated to transseries such as $\sum_{n = 0}^{\infty} n!
\frac{\ex^x}{x^{n + 1}}$, using the process of accelero-summation,  a far
reaching generalization of Borel summation; these analytic issues are not addressed in the present volume.)

These strong closure properties make it plausible to conjecture that
$\T$ is model complete, as a valued differential field. 
This and some other conjectures to be mentioned in this introduction go back some~20 years, and are proved in the present volume. 
To state model completeness of $\T$ geometrically we use the terms
%let us use the terms  
\textit{$\d$-algebraic}\/ and \textit{$\d$-polynomial}\/ to abbreviate \textit{differential-algebraic}\/ and \textit{differential polynomial}\/ and we define a \textit{$\d$-algebraic set}\/ in $\T^n$ to be 
the set of common zeros,
$$\big\{f=(f_1,\dots, f_n)\in \T^n:\ P_1(f)=\cdots = P_k(f)=0\big\}$$
of some $\d$-polynomials $P_1,\dots, P_k$ in differential indeterminates $Y_1,\dots, Y_n$, 
$$P_i(Y_1,\dots, Y_n)\ =\ p_i\big(Y_1,\dots,Y_n,\, Y_1',\dots, Y_n',\, Y_1'',\dots, Y_n'',\, Y_1''',\dots, Y_n''',\,\dots\big)$$ 
over $\T$. 
We also define an \textit{H-algebraic set}\/ to be the
intersection of a $\d$-algebraic set with a set of the form
\[ \big\{ y = (y_1, \ldots, y_n) \in \T^n :\ \text{$y_i \in
   \smallo_{\T}$ for all $i\in I$}   \big\} \quad\text{where $I \subseteq \{ 1, \ldots, n \}$,} \]
and we finally define a
subset of $\T^m$ to be \textit{sub-$H$-algebraic}\/ if it is the
image of an $H$-algebraic set in $\T^n$ for some $n\ge m$ under the
projection map $$(f_1,\dots, f_n) \mapsto (f_1,\dots, f_m)\colon \T^n \to \T^m.$$ It follows from the model completeness of $\T$
that the complement in $\T^m$ of any sub-$H$-algebraic set is again sub-$H$-algebraic, in analogy with Gabrielov's     ``theorem of the complement'' for real subanalytic sets~\cite{Gabrielov}.
(The model completeness of~$\T$ is a little stronger: it
is equivalent to this ``complement'' formulation where
the defining $\d$-polynomials of the $\d$-algebraic sets
involved have integer
coefficients.) A consequence is that for subsets of 
$\T^m $,
$$ \text{sub-$H$-algebraic}\ =\ \text{definable in $\T$.}$$ 

\smallskip\noindent
The usual model-theoretic approach to establishing that a given structure is model complete consists of two steps. (There is also a preliminary choice to be made of {\em primitives}; our choice for $\T$:  its ring operations, its derivation, its ordering, and its valuation.) 
The first step is to record the basic compatibilities between primitives; ``basic'' here means in practice that they 
are also satisfied by the {\em substructures\/} 
of the structure of interest. For the more familiar structure of the ordered field $\R$ of real numbers, these basic compatibilities are the ordered field axioms. 
The second and harder step is to find some closure properties satisfied by our structure that together with these basic compatibilities can be shown to imply 
{\em all\/} its elementary properties. In the model-theoretic treatment of~$\R$, it turns out that this job is done by the closure properties defining {\em real closed fields}: every positive element has a square root, and every odd degree polynomial has a zero.

\subsection*{$H$-fields}
For $\T$ we try to capture the first step of the axiomatization
 by the notion
of an \textit{$H$-field.}\/ 
%(One minor difference to the case of $\R$ is that ordered valued  differential subfields
%of $H$-fields are not automatically $H$-fields themselves.)
We chose the prefix $H$ in honor of E.~Borel, H.~Hahn, G.~H.~Hardy, and F.~Hausdorff,
who pioneered our subject about a century ago {\cite{Borel,Hahn,Hardy,Hau}},
and who share the initial~H, except for Borel.
To define $H$-fields, let $K$ be an ordered differential field (with constant field $C$) and
set
\begin{align*}
  \mathcal{O} &\ :=\ \big\{ a \in K :\  |a| \leqslant c \text{ for some $c>0$
   in $C$} \big\} && \text{(a convex subring of $K$),} \\
  \smallo &\ :=\  \big\{ a \in K :\  |a| < c \text{ for all $c > 0$ in $C$}
  \big\}. && 
\end{align*}
These notations should remind the reader
of Landau's big~O and small~o. The elements of~$\smallo$ are thought of as 
\textit{infinitesimal,}\/  the elements of~$\mathcal{O}$ as 
\textit{bounded,}\/  and those of $K\setminus\mathcal{O}$ as \textit{infinite.}\/ Note that $\mathcal{O}$ is definable in
the ordered differential field~$K$, and is a valuation ring of~$K$ with (unique) maximal ideal~$\smallo$.
We define $K$ to be an \textit{$H$-field}\/ if it satisfies the two conditions below:
\begin{list}{*}{\setlength\leftmargin{2.5em}}
\item[(H1)] for all $a\in K$, if $a> C$, then $a'>0$, 
\item[(H2)] $\mathcal{O}= C+\smallo$.
\end{list}
By~(H2) the constant field $C$ can be identified canonically with the residue
field~$\mathcal O/\smallo$ of~$\mathcal{O}$. 
As we did with $\T$ we construe an $H$-field $K$ as an
\textit{ordered valued differential field.}\/
An $H$-field $K$ is said to have \textit{small derivation}\/ if~$\der \smallo
\subseteq \smallo$ (and thus~$\der\mathcal O\subseteq\smallo$). If~$K$ is an $H$-field and $a \in K$, $a > 0$, then
$K$ with its derivation~$\der$ replaced by~$a \der$ is also an $H$-field. Such changes of derivation play
a major role in our work.

Among $H$-fields with small derivation are $\T$ and its ordered differential subfields
containing $\R$, and any Hardy field
containing $\R$. Thus $\R (x)$, $\R
(x, \ex^x, \log x)$ as well as Hardy's larger field of 
logarithmico-exponential functions are $H$-fields. 
%An interesting ``hybrid'' class of $H$-fields mediating between these two classes of $H$-fields was introduced
%in~\cite{vdH:hfsol}: a {\it transserial Hardy field}\/ is
%an ordered differential subfield of $\T$ which also carries the structure of a Hardy field, 
%assuming some natural compatibility conditions. (In our later work, we plan
%to generalize the results from~\cite{vdH:hfsol}  by investigating how $H$-fields can be embedded in
%Hardy fields and fields of transseries.)

%Several variations on this definition are possible. For instance, one might consider Hardy fields that contain only germs of infinitely differentiable functions, or even more strict, only germs of real analytic functions (as is the case for Hardy's field $H$).

\subsection*{Closure properties} Let $\operatorname{Th}(\mathbf{M})$ be the first-order theory of an
$\mathcal L$-structure
$\mathbf M$, that is, $\operatorname{Th}(\mathbf M)$ is the set of $\mathcal{L}$-sentences that are true in $\mathbf M$; see Appendix~\ref{app:modth} for details.
In terms of $H$-fields, we can now make the model completeness conjecture more precise, as was done in~\cite{AvdD2}: 
$$\operatorname{Th}(\T)\ =\ \text{model companion of the theory of $H$-fields with small derivation},$$
where $\T$ is construed as an ordered and valued differential field.
This amounts to
adding to the earlier model completeness of $\T$ the claim that
any $H$-field with small derivation can be embedded as an ordered valued differential field into some ultrapower of $\T$.
Among the consequences of this conjecture is that any finite
system of algebraic differential equations over $\T$ (in several unknowns) has a solution in $\T$ whenever it has one in some $H$-field extension
of $\T$. 
It means that the concept of ``$H$-field'' is intrinsic to the differential field $\T$.
It also suggests studying systematically the extension theory
of $H$-fields: A.~Robinson taught us that for a theory to have a model companion at all---a rare phenomenon---is equivalent to certain embedding and extension properties of its class of
models. Here it helps to know that $H$-fields fall under the so-called {\em differential-valued fields} (abbreviated
as {\em $\d$-valued fields\/} below) of Rosenlicht, who began a study of these valued differential fields and their extensions in the early 1980s; see \cite{Rosenlicht2}.  (A {\em $\d$-valued\/} field is defined to be a valued differential field such
that $\mathcal{O}=C + \smallo$, and
$a'b\in b'\smallo$ for all $a,b\in \smallo$; here $\mathcal{O}$ is the valuation ring with maximal ideal $\smallo$, and $C$ is the constant field.) 
Most of our work is actually in the setting of valued differential fields where no field ordering is given, since even for $H$-fields the valuation is a more robust and useful feature than its field ordering. 

Besides developing the extension theory of $H$-fields we need to isolate the relevant \textit{closure properties}\/ of $\T$. 
First, $\T$ is real closed, but that property does not involve the derivation. Next, $\T$ is closed under integration and, by its very construction, also under exponentiation.  In terms of the derivation this gives two natural closure 
properties of~$\T$:
\[       \forall a \exists b\, (a = b'), \qquad     \forall a \exists b\, (b \neq 0 \mathbin{\&} ab = b').
\]
An $H$-field $K$ is said to be \textit{Liouville closed}\/ if it is real
closed and satisfies these two sentences; cf. Liouville~\cite{Liouville1,Liouville2}. 
So $\T$ is Liouville closed. It was shown in~\cite{AvdD2} that
any $H$-field has a \textit{Liouville closure,}\/ that is, a minimal
Liouville closed $H$-field extension. 
If $K$ is a Hardy field
containing~$\R$ as a subfield, then it has a unique Hardy field extension that
is also a Liouville closure of $K$, but it can happen that an $H$-field~$K$ 
has two Liouville closures that are not isomorphic over $K$; it cannot have
more than two. Understanding this ``fork in the road'' and dealing with it is
fundamental in our work.
 Useful notions in this connection are
{\em comparability classes}, {\em groundedness}, and
{\em asymptotic integration.\/} We discuss this briefly below for 
$H$-fields. (Parts of Chapters~\ref{ch:asymptotic differential fields} and~\ref{evq} treat these notions for a much larger class of valued differential fields.) 
Later in this introduction we encounter an important but rather hidden closure property,
called {\em $\upo$-freeness}, which rules over the fork in the road. Finally, there is the very powerful closure property of  newtonianity that we already mentioned earlier.

\subsection*{Valuations and asymptotic relations}

Let $K$ be an $H$-field, let $a$,~$b$ range over~$K$, and let $v \colon K \to \Gamma_{\infty}$ be the
(Krull) valuation on $K$ associated to $\mathcal{O}$, with value group
$\Gamma = v (K^{\times})$ and $\Gamma_{\infty} := \Gamma \cup
\{\infty\}$ with $\Gamma < \infty$.
Recall that $\Gamma$ is an ordered abelian group,  additively written
as is customary in valuation theory. Then
\[ va < vb \quad\Longleftrightarrow\quad |a| > c |b|  \text{ for all $c > 0$ in $C$.} \]
Thinking of elements of $K$ as
germs of functions at $+\infty$, we also adopt Hardy's notations from
asymptotic analysis:
\[ a \succ b, \hspace{2em} a \succcurlyeq b, \hspace{2em} a \prec b,
   \hspace{2em} a \preccurlyeq b, \hspace{2em} a \asymp b, \hspace{2em} a \sim
   b \hspace{2em} \]
are defined to mean, respectively,
\[ va < vb, \hspace{1em} va \leqslant vb, \hspace{1em} va > vb, \hspace{1em}
   va \geqslant vb, \hspace{1em} va = vb, \hspace{1em} v (a - b) > va. \]
(Some of these notations from~\cite{Hardy2} actually go back to du
Bois-Reymond~\cite{dBR71}.) 
%One should be aware of the ``reversal'' in
%the equivalence $a \succ b \Longleftrightarrow v a < v b$. 
Note that $a \succ
1$ means that~$a$ is infinite, that is, $|a| > C$, and $a \prec 1$ means
that~$a$ is infinitesimal, that is, $a \in \smallo$. It is crucial that the asymptotic relations above can be differentiated, provided we restrict to nonzero
$a,b$ with $a\not\asymp 1$, $b\not\asymp 1$: 
$$a\succ b\ \Longleftrightarrow\ a'\succ b',\qquad 
 a\asymp b\ \Longleftrightarrow\ a'\asymp b', \qquad a\sim b\ \Longleftrightarrow\ a'\sim b'.$$
%From time to time, we will also use Landau's notations $O (a)$ and $o (a)$ to mean, respectively, ``any element~$b$ with $b \preccurlyeq a$'' and ``any element~$b$ with $b \prec a$''. \marginpar{Really used?}
For $a \ne 0$ we let $a^{\dagger} := a' / a$ be its logarithmic
derivative, so $(ab)^{\dagger} = a^{\dagger} + b^{\dagger}$ for~${a, b \ne 0}$.
Elements $a, b \succ 1$ are said to be \textit{comparable}\/ if $a^{\dagger}
\asymp b^{\dagger}$; if $K$ is a Hardy field containing~$\R$ as subfield, or
$K = \T$, this is equivalent to the existence of an~$n \geqslant 1$ such that
$|a| \leqslant |b|^n$ and $|b| \leqslant |a|^n$. Comparability is an
equivalence relation on the set of infinite elements of $K$, and the comparability classes~Cl$(a)$
of~$K$ are totally ordered by Cl$(a) \leqslant \operatorname{Cl} (b) :
\Longleftrightarrow a^{\dagger} \preccurlyeq b^{\dagger}$.

\begin{exampleUnnumbered} 
For $K=\T$, set $e_0=x$ and $e_{n+1}= \exp(e_n)$.
Then the sequence $(\operatorname{Cl}(e_n))$ is strictly increasing and cofinal in the 
set of comparability classes. More important are the $\ell_n$
defined recursively by $\ell_0=x$, and $\ell_{n+1}=\log \ell_n$. Then the sequence $\operatorname{Cl}(\ell_0)>\operatorname{Cl}(\ell_1)>\operatorname{Cl}(\ell_2)>\cdots > \operatorname{Cl}(\ell_n)>\cdots$ is coinitial in the
set of comparability classes of $\T$. For later use it is worth
noting at this point that
$$\ell_n^\dagger\ =\ \frac{1}{\ell_0\cdots \ell_n},\qquad -\ell_n^{\dagger\dagger}\ =\ \frac{1}{\ell_0} + \frac{1}{\ell_0\ell_1} + \cdots + \frac{1}{\ell_0\ell_1\cdots \ell_n}.
$$
\end{exampleUnnumbered}

\noindent
We call $K$ {\em grounded\/} if $K$ has a smallest comparability class. Thus $\T$ is not grounded.
If~$\Gamma^{>}$
contains an element~$\alpha$ such that for every~$\gamma \in \Gamma^{>}$ we
have $n \gamma \geqslant \alpha$ for some~${n \geqslant 1}$, then~$K$ is grounded;
this condition on $\Gamma$ is in particular
satisfied if~${\Gamma \ne \{0\}}$ and $\Gamma$ has finite archimedean rank.  
%Also, if $K$ is of finite transcendence degree over $C$, then $K$ is grounded. 
If $K$ is grounded, then $K$ has only one Liouville closure (up to isomorphism over $K$).

The $H$-field $K$ is said to have {\em asymptotic integration\/} if $K$ satisfies
$\forall a \exists b (a\asymp b')$, equivalently,~${\{vb':\ b\in K\}=\Gamma_{\infty}}$. It is obvious that every Liouville closed $H$-field has asymptotic integration; in particular, $\T$ has asymptotic integration. In general,
at most one $\gamma\in \Gamma$ lies outside
${\{vb':\ b\in K\}}$; if~$K$ is grounded, then such a $\gamma$ exists, by results in Section~\ref{AbstractAsymptoticCouples}, and so $K$ cannot have asymptotic integration.

\ifbool{PUP}{
\section*{Strategy and Main Results}
\addcontentsline{toc}{section}{Strategy and Main Results}
}{\section*{Strategy and Main Results}}

\noindent
Model completeness of
$\T$ concerns finite systems of algebraic differential equations over $\T$ with asymptotic side conditions
in several differential indeterminates.

Robinson's strategy for establishing model completeness applied to $\T$ requires us to 
move beyond $\T$ to consider $H$-fields and their extensions. If we are lucky---as we are in this case---it will
suffice to consider extensions of $H$-fields by one element~$y$ at a time.
This leads to equations $P(y)=0$ with an asymptotic side condition $y\prec g$. Here $P\in K\{Y\}$ is a univariate differential
polynomial with coefficients in an $H$-field~$K$ with $g\in K^\times$, and
$K\{Y\}=K[Y, Y', Y'',\dots]$ is the differential domain
 of $\d$-polynomials in the differential indeterminate $Y$ over $K$. The key issue: when is there a solution in some $H$-field extension of $K$? 
A detailed study of such equations in the special case
$K = \T_{\operatorname{g}}$ and where we only look for solutions in $\T_{\operatorname{g}}$ itself was undertaken in~\cite{JvdH}, using an assortment of techniques (for instance, various fixpoint theorems)
heavily based on the particular structure of $\T_{\operatorname{g}}$.
Generalizing these results to suitable $H$-fields
is an important guideline in our work.

\subsection*{Differential Newton diagrams}

Let $K$ be an $H$-field, and consider a $\d$-algebraic equation with asymptotic side condition,
\begin{equation}\tag{4}
  \label{aade} P (y)\ =\ 0, \hspace{2em} y\ \prec\ g,
\end{equation}
where $P \in K \{ Y \}$, $P\ne 0$, and $g \in K^\times$; we look for nonzero solutions in $H$-field extensions of $K$. For the sake of concreteness we take $K=\T_{\operatorname{g}}$ and
look for nonzero solutions in $\T_{\operatorname{g}}$, focusing on 
the example below:
\begin{equation}\tag{5}
  \label{tr-aade-ex} 
  \ex^{- \ex^x} y^2 y'' + y^2 - 2 xyy' - 7\ex^{-x}y' - 4 + \textstyle\frac{1}{\log x}\ =\ 0, \hspace{2em} y\ \prec\ x.
\end{equation}
We sketch briefly how~\cite{JvdH} goes about 
solving \eqref{tr-aade-ex}.   
%Let $r$ be the order of the $\d$-polynomial $P$, let
%$\i=(i_0,\dots,i_r)$ and $\j$ range over $\N^{1+r}$, and set
%$Y^\i:=Y^{i_0}(Y')^{i_1}\cdots (Y^{(r)})^{i_r}$, and likewise %with $y\in K^\times$ in place of $Y$. Thus
%$$P = \sum_{\i} P_{\i}Y^{\i}\quad\text{
%where $P_\i\in K$ and $P_\i=0$ for all but finitely many $\i$.}$$
First of all, we need to find the possible \textit{dominant terms}\/ of solutions $y$. This is done by considering possible
cancellations. 
For example, $y^2$ and $-4$ might be the terms of least valuation in the left side of~\eqref{tr-aade-ex}, with all other terms having greater valuation, so negligible compared to~$y^2$ and~$-4$. This yields a cancellation $y^2\sim 4$, so $y\sim 2$ or $y\sim -2$, giving~$2$ and~$-2$ as potential dominant terms of a solution $y$.
  
Another case: $\ex^{-\ex^x}y^2y''$ and $y^2$ are the terms of least valuation. Then we get a cancellation 
$\ex^{- \ex^x} y^2 y'' \sim - y^2$, that is, $y'' \sim - \ex^{\ex^x}$, which leads to $y \sim - \ex^{\ex^x} /
\ex^{2 x}$. But this possibility is discarded, since \eqref{tr-aade-ex} also requires $y \prec x$. (On the other hand, if the asymptotic condition in \eqref{tr-aade-ex}
had been $y\prec \ex^{\ex^x}$, we would have kept 
$- \ex^{\ex^x} /
\ex^{2 x}$ as a potential dominant term of a solution $y$.)

What makes things work in these two cases is that the cancellations arise from terms of {\em different\/} degrees in $y, y', y'',\dots$.  
Such cancellations are reminiscent of the more familiar setting of algebraic equations 
where the dominant monomials of solutions 
can be read off from a
\textit{Newton diagram\/} and the corresponding dominant coefficients are zeros of the
corresponding \textit{Newton polynomials}\/;  see Section~\ref{Newton Diagrams}. This method still works in our $\d$-algebraic setting, for cancellations among terms of different degrees, but requires the construction of so-called {\em equalizers}.
%for example, the transmonomial $\ex^{\ex^x} / \ex^{2 x}$ equalizes the homogeneous parts of degree~$2$ and~$3$ in \eqref{tr-aade-ex}.  
 
\medskip\noindent
A different situation arises for cancellations between terms of the {\em same\/} degree. Consider for example the case that
$y^2$ and $-2xyy'$ have least valuation among the terms in the left side of~\eqref{tr-aade-ex}, with all other terms of higher valuation. Then $y^2 \sim 2 xyy'$, so $y^{\dag} \sim \frac{1}{2x}$. Now 
$y^{\dag}= \frac{1}{2x}$ gives $y=c x^{1/2}$ with $c \in \R^\times$, but the weaker condition $y^{\dag} \sim \frac{1}{2x}$
only gives $y=ux^{1/2}$ with $u\ne 0$, $u^\dagger\prec x^{-1}$,
that is, $|v(u)| < |v(x)|/n$ for all $n\ge 1$. Substituting $ux^{1/2}$ for $y$ in \eqref{tr-aade-ex} and considering~$u$ as the new unknown, the condition on $v(u)$ forces $u\asymp 1$, so after all we do get 
$y\sim c x^{1/2}$ with
$c\in \R^\times$, giving~$cx^{1/2}$ as a potential dominant term of a solution $y$.  It is important to note that here an {\em integration constant\/} $c$ gets introduced. 

Manipulations as we just did are similar to rewriting an equation $H(y)=0$ with {\em homogeneous\/} nonzero $H\in K\{Y\}$ of positive degree as a (Riccati) equation~${R(y^\dagger)=0}$ with $R$ of {\em lower} order than $H$.

This technique can be shown to work in general for cancellations among terms of the same degree, provided we are also allowed to transform the equation to an equivalent one by applying a suitable iteration of the {\em upward shift}\/ $f(x)\mapsto f(\ex^x)$. (For reasonable $H$-fields $K$ one can apply instead compositional conjugation by positive active elements; see below for \textit{compositional conjugation}\/ and \textit{active.}\/)

\begin{samepage}
 
\medskip\noindent
Having determined a possible dominant term $f=c\fm$ of a solution
of~\eqref{aade}, where ${c\in\R^\times}$ and $\fm$~is a transmonomial, we next perform a so-called \textit{refinement}\/
\begin{equation}\label{aade refinement}\tag{6}
 P(f+y)\ =\ 0,\hspace{2em} y\ \prec\ f
\end{equation}
of~\eqref{aade}.
%which consists of a change of variables together with the imposition of a~new asymptotic side condition. 
For instance, taking $f = 2$, the
equation~\eqref{tr-aade-ex} transforms into
\begin{eqnarray*}
  \ex^{- \ex^x}  y^2  y'' + y^2 - 2 x y 
  y' + 4 \ex^{- \ex^x}  y  y'' &  & \\
  {}+ 4 y - (4 x + 7\ex^{-x})  y' + 4 \ex^{- \ex^x} 
  y'' + \textstyle\frac{1}{\log x} & = & 0, \hspace{2em} y \prec 2.
\end{eqnarray*}
Now apply the same procedure to this refinement, to find the ``next'' term. \looseness=-1

\end{samepage}

\medskip\noindent
Roughly speaking, this yields an infinite process to obtain all possible
asymptotic expansions of solutions to any asymptotic equation.
How do we make this into a finite process?
For this, it is useful
to introduce the \textit{Newton degree}\/ of \eqref{aade}. This notion is
similar to the Weierstrass degree of a multivariate power series and
corresponds to the degree of the asymptotically significant part of the
equation. If the Newton degree is~$0$, then \eqref{aade}
has no solution. The Newton degree of~\eqref{tr-aade-ex} turns out to
be~$2$: this has to do with the fact that
$\ex^{- \ex^x} y^2
y'' \prec y^2$ whenever $y \prec x$. We shall return soon to the precise
definition of Newton degree for differential polynomials over rather general $H$-fields. As to the resolution of asymptotic  
equations over $K=\T_{\operatorname{g}}$, the following key facts were established in~\cite{JvdH}:
\begin{itemize}
  \item The Newton degree stays the same or decreases under refinement. 
%cannot strictly increase under refinement.
  \item If the Newton degree of the refinement \eqref{aade refinement} equals that of  \eqref{aade}, we employ so-called \textit{unravelings}\/; these resemble the \textit{Tschirnhaus transformations}\/ that overcome similar obstacles in the algebraic setting.
  Combining unravelings with refinements as described above, we arrive after finitely many steps
   at an asymptotic equation of Newton degree~$0$ or $1$.  
  \item The $H$-field $\T_{\operatorname{g}}$ is newtonian,
  that is, 
  any asymptotic equation over $\T_{\operatorname{g}}$ of Newton degree~$1$ has a solution
  in $\T_{\operatorname{g}}$.
\end{itemize}
All in all, we have for any given
asymptotic equation over $\T_{\operatorname{g}}$ a more or less finite procedure for gaining an overview of the entire space of solutions in $\T_{\operatorname{g}}$. 
    
\medskip
\noindent
To define the Newton degree of an
asymptotic equation \eqref{aade} over rather general $H$-fields, we first
need to introduce the \textit{dominant part}\/ of~$P$ and then, based on a process called \textit{compositional conjugation,}\/ the \textit{Newton polynomial}\/ of $P$.

\subsection*{The dominant part}
Let $K$ be an $H$-field. 
We extend the valuation $v$ of $K$ to the integral domain
$K\{Y\}$ %the differential domain $K \{Y\}$ of $\d$-poly\-no\-mials in the differential indeterminate $Y$ over $K$ 
by setting
\[ vP\ =\ \min \{va:\  \text{$a$ is a coefficient of $P$}\}, \]
and we extend the binary relations $\asymp$ and $\sim$ on $K$ to $K\{Y\}$ accordingly. 
It is also convenient to fix a monomial set $\mathfrak{M}$ in $K$, that is, a
subset $\mathfrak{M}$ of $K^{>}$ that is mapped bijectively by $v$ onto the
value group $\Gamma$ of $K$. This allows us to define the \textit{dominant
part~$D_P (Y)$} of a nonzero $\d$-polynomial $P (Y)$ over $K$ to be the
unique element of $C \{Y\} \subseteq K \{Y\}$ with $P\sim  \mathfrak{d}_PD_P$, where $\mathfrak{d}_P \in \mathfrak{M}$ is the {\em{dominant
monomial}} of $P$ determined by $P\asymp  \mathfrak{d}_P$. 
(Another choice of monomial
set would just multiply $D_P$ by some positive constant.) For $K=\T$ we
always take the set of transmonomials as our monomial set. 

\begin{exampleintro}
  \label{dom-ex} Let $K = \T$. For $P = x^5 + (2 + \ex^x) Y + (3 \ex^x + \log x) (Y')^2$, we have $\mathfrak{d}_P= \ex^x$ and 
  $D_P = Y + 3 (Y')^2$.
  For $Q = Y^2 - 2xYY'$  we have $D_Q  = - 2
YY'$.
\end{exampleintro}

\noindent
For $K$ with small derivation we can use $D_P$ to get near the zeros 
$a\asymp 1$ of
$P$: if $P (a) = 0$, $a \asymp 1$, then $D_P (c) = 0$
where $c$ is the unique constant with~${a \sim c}$.
We need to understand, however, the behavior of $P (a)$ not only for $a \asymp 1$,
that is,~${va=0}$, but also for ``sufficiently flat'' elements $a \in K$, that is, for $va$ approaching~$0\in \Gamma$. For instance, in $\T$,
the iterated logarithms $$\ell_0=x,\quad \ell_1=\log x,\quad \ell_2=\log \log x,\   \ldots$$  satisfy $v(\ell_n) \to 0$ in $\Gamma_{\T}$ and likewise $v(1/\ell_n) \to 0$. 
The \textit{dominant term~$\mathfrak{d}_PD_P$}\/ of~$P$ often
provides a good approximation for $P$ when evaluating at sufficiently flat
elements, but not always: for $K=\T$ and $Q$ as in Example~\ref{dom-ex} we note that  for $y = \ell_2$ we have:  $y^2 = \ell_2^2\succ  2 xyy' =  2 \ell_2 / \ell_1$, so $Q(y)\sim y^2\not\asymp (\mathfrak d_Q D_Q)(y)$.

In order to approximate $P(y)$ by $(\mathfrak{d}_P D_P)(y)$ for  
sufficiently flat~$y$, we need one more ingredient: \textit{compositional
conjugation.}\/ For $K=\T$ and $Q$ as in Example~\ref{dom-ex}, this amounts to a change of
variables $x = \ex^{\ex^{\tilde x}}$,
%\exp \exp \tilde x$, 
so that $Q (y) = y^2 - 2y (dy/d\tilde{x})\ex^{-\tilde x}$ for $y\in\T$. With respect to this new variable $\tilde x$, the
dominant term~$Y^2$ of the adjusted $\d$-polynomial $Y^2 - 2YY' \ex^{-\tilde x}$ is then an adequate 
approximation of $Q$ when
evaluating at sufficiently flat elements of $\T$. Such changes of variable do not make sense
for general $H$-fields, but as it turns out, compositional conjugation is a good substitute.

\subsection*{Compositional conjugation}
We define this for an arbitrary
differential field $K$. For $\phi \in K^{\times}$ we let $K^{\phi}$ be the
differential field obtained from $K$ by replacing its derivation $\der$ by the
multiple $\phi^{- 1}  \der$. Then a differential polynomial $P (Y) \in K
\{Y\}$ defines the same function on the common underlying set of $K$ and~$K^{\phi}$ as a certain differential polynomial $P^{\phi} (Y) \in K^{\phi}
\{Y\}$: for $P = Y'$, we have $P^{\phi} (Y) = \phi Y'$ (since over $K^{\phi}$
we evaluate $Y'$ according to the derivation $\phi^{- 1}  \der$), for $P =
Y''$ we have $P^{\phi} (Y) = \phi' Y' + \phi^2 Y''$ (with $\phi' = \der
\phi$), and so on. This yields a ring isomorphism
\[ P \mapsto P^{\phi}\colon K \{Y\} \to K^{\phi} \{Y\} \]
that is the identity on the common subring $K [Y]$. It is also an automorphism
of the common underlying $K$-algebra of $K \{Y\}$ and $K^{\phi} \{Y\}$, and
studied as such in Chapter~\ref{ch:triangular automorphisms}. We call
$K^{\phi}$ the \textit{compositional conjugate of $K$ by $\phi$,}\/ and
$P^{\phi}$ the \textit{compositional conjugate of $P$ by $\phi$.}\/ Note that
$K$ and $K^{\phi}$ have the same constant field $C$. If $K$ is an $H$-field
and $\phi \in K^{>}$, then so is $K^{\phi}$. It pays to note how things change
under compositional conjugation, and what remains invariant.

\subsection*{The Newton polynomial}
Suppose now that $K$ is an $H$-field with asymptotic integration. For $\phi
\in K^{>}$ we say that~$\phi$ is \textit{active}\/ (in $K$) if $\phi
\succcurlyeq a^{\dagger}$ for some nonzero~$a \nasymp 1$ in $K$; equivalently,
the derivation $\phi^{- 1}  \der$ of $K^{\phi}$ is small. Let $\phi \in K^{>}$
range over the active elements of~$K$ in what follows, fix a monomial set
$\mathfrak{M} \subseteq K^{>}$ of~$K$, and let $P \in K \{Y\}$, $P \ne 0$. 
The dominant part
$D_{P^{\phi}}$ of $P^{\phi}$ lies in $C \{Y\}$, and we show in Section~\ref{sec:The Dominant Part of a Differential Polynomial} that it  eventually stabilizes as
$\phi$ varies: there is a differential polynomial $N_P \in C \{Y\}$ and an
active~$\phi_0 \in K^{>}$ such that for all $\phi \preccurlyeq \phi_0$,
\[ D_{P^{\phi}}\ =\ c_{\phi} N_P, \hspace{2em} c_{\phi} \in C^{>} . \]
We call $N_P$ the \textit{Newton polynomial}\/ of $P$. It is of
course only determined up to a factor from $C^{>}$, but this ambiguity is
harmless. The (total) degree of $N_P$ is called the \textit{Newton degree}\/
of $P$.

\begin{exampleintro}
 Let $K = \T$. Then $f\in K^>$ is active iff $f\succeq \ell_n^\dagger=\frac{1}{\ell_0\ell_1\cdots\ell_n}$ for some~$n$.
 If~$P$ is as in Example~\ref{dom-ex},
 then for each $\phi$,
$$P^{\phi}\ =\
  x^5 + (2 + \ex^x) Y + \phi^2  (3 \ex^x + \log x) (Y')^2,$$
  so $D_{P^{\phi}} = Y$ if $\phi \prec 1$. This yields $N_P = Y$, so $P$ has
  Newton degree~$1$. It is an easy exercise to
  show that for $Q = Y^2 - 2xYY'$ we have $N_Q = Y^2$ .
\end{exampleintro}

\noindent
A crucial result in~\cite{JvdH} (Theorem~8.6) says that if $K =\T_{\operatorname{g}}$, then $N_P \in \R [Y]
  (Y')^{\N}$. A major step in our work was to isolate a robust class of $H$-fields $K$ with asymptotic integration for which likewise $N_P\in C[Y](Y')^{\N}$ for all nonzero $P\in K\{Y\}$. This required several completely new tools to
be discussed below.

\subsection*{The special cuts $\upu$, $\upl$ and $\upo$}
Recall that  $\ell_n$ denotes the $n$th iterated logarithm of~$x$ in $\T$, so $\ell_0 = x$ and $\ell_{n + 1} = \log \ell_n$. We  
introduce the elements %\marginpar{Why not $\upg$ etc. instead of $\upu$?}
\[ \begin{array}{rcccl}
     \upu_n & = & \ell_n^{\dag} & = & \displaystyle\frac{1}{\ell_0 \cdots \ell_n}\\[1em]
     \upl_n & = & - \upu_n^{\dag} & = & \displaystyle\frac{1}{\ell_0} +
     \frac{1}{\ell_0 \ell_1} + \cdots + \frac{1}{\ell_0 \ell_1 \cdots
     \ell_n}\\[1em]
     \upo_n & = & - 2 \upl_n' - \upl_n^2 & = & \displaystyle\frac{1}{\ell_0^2} +
     \frac{1}{\ell_0^2 \ell_1^2} + \cdots + \frac{1}{\ell_0^2 \ell_1^2 \cdots
     \ell_n^2} 
   \end{array} \]
of $\T$. As $n\to \infty$ these elements approach their formal limits
\begin{eqnarray*}
  \upu_{\T} & = & \frac{1}{\ell_0 \ell_1 \ell_2 \cdots}\\
  \upl_{\T} & = & \frac{1}{\ell_0} + \frac{1}{\ell_0 \ell_1} +
  \frac{1}{\ell_0 \ell_1 \ell_2} + \cdots\\
  \upo_{\T} & = & \frac{1}{\ell_0^2} + \frac{1}{\ell_0^2 \ell_1^2} +
  \frac{1}{\ell_0^2 \ell_1^2 \ell_2^2} + \cdots ,
\end{eqnarray*}
which for now are just suggestive expressions. 
Indeed, our field $\T$ of transseries {\em{of finite logarithmic
and exponential depth}} does not contain any pseudolimit of the pseudocauchy 
sequence $(\upl_n)$, nor of the pseudocauchy sequence $(\upo_n)$. There
are, however, immediate 
$H$-field extensions of $\T$ where such pseudolimits exist, and if we let
$\upl_{\T}$ be such a pseudolimit of $(\upl_n)$, then in some further 
$H$-field extension we have an element suggestively denoted by 
$\exp(\int {-\upl_{\T}})$ that can play the role of $\upu_{\T}$. 

\medskip\noindent
Even though $\upu_{\T}$, $\upl_{\T}$ and
$\upo_{\T}$ are not in $\T$, we can take them as elements of some $H$-field
extension of $\T$, as indicated above, and so we obtain sets
\begin{eqnarray*}
  \Upu (\T) & = & \{ a \in \T:\ a >
  \upu_{\T} \}\\
  \Upl (\T) & = & \{ a \in \T:\ a <
  \upl_{\T} \}\\
  \Upo(\T) & = & \{ a \in \T:\ a < \upo_{\T} \}
\end{eqnarray*}
that can be shown to be definable in $\T$. For instance,
\begin{eqnarray*}
  \Upu (\T) & = & \big\{ a \in \T:\ \forall b \in \T\ ( 
  b \succ 1 \Rightarrow a \neq b^{\dag}) \big\}\\
  & = & \{ - a' :\  a \in \T,\ a \geqslant 0 \} .
\end{eqnarray*}
In other words, $\upu_{\T}$, $\upl_{\T}$ and
$\upo_{\T}$ realize {\em{definable cuts}} in $\T$.

\medskip\noindent
For any ungrounded $H$-field $K\ne C$ we can build
a sequence $(\ell_{\rho})$ of elements
$\ell_{\rho} \succ 1$, indexed by the ordinals $\rho$ less than some infinite 
limit ordinal, such that
$$\sigma > \rho \ \Rightarrow\  \ell_{\sigma}^{\dag}
\prec \ell_{\rho}^{\dag}, \qquad v (\ell_{\rho}) \rightarrow 0 \text{ in $\Gamma$.}$$ 
These $\ell_{\rho}$ play in $K$ the role that the iterated logarithms 
$\ell_n$ play in $\T$. In analogy with~$\T$ they yield the elements
$$  \upu_\rho\ :=\ \ell_{\rho}^{\dag}, \qquad
     \upl_\rho\ :=\  - \upu_{\rho}^{\dag}, \qquad
     \upo_\rho\ :=\  - 2 \upl_n' - \upl_n^2,$$
of $K$, and $(\upl_{\rho})$ and $(\upo_{\rho})$ are pseudocauchy sequences. As with
$\T$ this gives rise to definable sets $\Upu(K)$, $\Upl(K)$ and
$\Upo(K)$ in $K$. The fact mentioned earlier that $\T$ does not contain 
$\upu_{\T}$,
$\upl_{\T}$ or $\upo_{\T}$ turns out to be very significant: in
general, we have $$\upu_K \in K\quad \Rightarrow\quad \upl_K \in K
\quad\Rightarrow\quad \upo_K \in K$$ and each of the four mutually exclusive cases
\[ \upu_K \in K, \qquad
     \upu_K \notin K\ \&\ \upl_K \in K, \qquad 
     \upl_K \notin K\ \&\  \upo_K \in K, \qquad  \upo_K \notin K \]
can occur; see Section~\ref{sec:specialH}. Here we temporarily abuse notations, since we should explain what we mean by 
$\upu_K\in K$ and the like; see the next subsections. 

\subsection*{On gaps and forks in the road}
Let $K$ be an $H$-field.
We say that an element~${\upu\in K}$ is a {\em{gap}} in $K$ if
for all $a \in K$ with $a \succ 1$ we have
\[ a^{\dag} \succ \upu \succ {(1 / a)'}.\] 
The existence of such a gap is the formal
counterpart to the informal statement that $\upu_K \in K$.
% We say that $K$ is {\em{$\upu$-free}} if no such gap exists. 
If $K$ has a gap $\upu$, then $\upu$ has no primitive in $K$, so $K$ is not closed under integration.
If $K$ has trivial derivation (that is, $K=C$), then $K$ has a gap~$\upu=1$. There are also~$K$ with $K\ne C$ (even Hardy fields) that have a gap. Not having a gap is equivalent to being grounded or having asymptotic integration. 

We already mentioned the result from~\cite{AvdD2} that $K$ may have two Liouville
closures that are not isomorphic over $K$ (but fortunately not more than two).
Indeed, if $K$ has a gap $\upu$, then in one Liouville closure
all primitives of $\upu$ are infinitely large,
whereas in the other $\upu$ has an infinitesimal primitive.
Even if $K$ has no gap, the above fork in the road can arise more indirectly: Assume that $K$ has asymptotic integration and $\upl\in K$ is such that for all $a \in K^\times$ with $a \prec 1$,
\[ a'^{\dag}\  <\  -\upl\  <\  a^{\dag\dag}.\] 
Then $K$ has no element $\upu\neq 0$
with $\upl=-\upu^{\dag}$, but $K$ has an $H$-field extension
$K \langle \upu \rangle$ generated by an element
$\upu$ with $\upl=-\upu^{\dag}$, and any such $\upu$ is a gap
in $K \langle \upu \rangle$. It follows again
that $K$ has two Liouville closures that are not $K$-isomorphic. 
 
 For real closed $K$ with asymptotic integration,
 the existence of such an element~$\upl$ corresponds to the informal statement that $\upu_K\notin K\ \&\ \upl_K \in K$. We define~$K$ to be
{\em{$\upl$-free}} if $K$ has asymptotic integration and satisfies the sentence
\[ \forall a \exists b \big[b \succ 1\ \&\ a - b^{\dagger \dagger} \succcurlyeq
   b^{\dagger}\big] . \]
It can be shown that for real closed $K$ with asymptotic integration, $\upl$-freeness is equivalent to the nonexistence of an element $\upl$ as
above. More generally, $K$ is $\upl$-free iff $K$ has asymptotic integration and $(\upl_{\rho})$ has no pseudolimit in $K$.

\subsection*{The property of $\upo$-freeness}

Even $\upl$-freeness might not prevent a fork in the road for some $\d$-algebraic extension. Let $K$ be an $H$-field, and define
\[ \omega = \omega_K \colon K \to K, \hspace{2em} \omega (z) := - 2 z' - z^2.
\]
Assume that $K$ is $\upl$-free and $\upo\in K$
is such that for all $b\succ 1$ in $K$,
\[\upo - \omega (b^{\dag\dag}) \prec (b^{\dag})^2.\] 
Then the first-order differential equation $\omega
(z) = \upo$ admits no solution in $K$, but~$K$ has an $H$-field extension $K \langle
\upl \rangle$ generated by a solution $z=\upl$ to $\omega
(z) = \upo$ such that~$K\langle \upl\rangle$ is no longer
$\upl$-free (and with a fork in its road towards Liouville closure).

For $\upl$-free $K$ the existence of an element $\upo$ as above corresponds to the informal statement that $\upl_K \notin K\ \&\  \upo_K \in K$. We say that $K$ is {\em{$\upo$-free}} if no such
$\upo$ exists, more precisely,
$K$ has asymptotic integration and satisfies the sentence
%{\tmabbr{i.e.}}
\[ \forall a \exists b \big[b \succ 1\ \&\ a - \omega (b^{\dagger \dagger})
   \succcurlyeq (b^{\dagger})^2\big].  \]
(It is easy to show that if $K$ is $\upo$-free, then it is 
$\upl$-free.) For $K$ with asymptotic integration, $\upo$-freeness is equivalent to the pseudocauchy sequence
$(\upo_{\rho})$ not having a pseudolimit in $K$.  Thus $\T$ is $\upo$-free. More generally, if $K$ has
asymptotic integration and is a union of grounded $H$-subfields, 
then $K$ is $\upo$-free by Corollary~\ref{recam}.

Much deeper and very useful is that if $K$ is an $\upo$-free $H$-field and $L$
is a~$\d$-algebraic $H$-field extension of $K$, then $L$ is
also $\upo$-free and has no comparability class smaller than all those of $K$; this is part of Theorem~\ref{upoalgebraic}. Thus
the property of $\upo$-freeness is very robust: if $K$ is
$\upo$-free, then forks in the road towards Liouville closure
no longer occur, even for $\d$-algebraic $H$-field extensions of $K$ (Corollary~\ref{upoliou}).
There are, however, Liouville closed $H$-fields that are not
$\upo$-free; see~{\cite{ADH}}.

Another important consequence of $\upo$-freeness is that  
Newton polynomials of differential polynomials 
then take the same simple shape as those over $\T_{\operatorname{g}}$:

\begin{theoremintro} If $K$ is $\upo$-free and $P\in K\{Y\}$, $P\ne 0$, then $N_P\in C[Y](Y')^{\N}$.
\end{theoremintro}

\noindent
%These deeper results about $\upo$-freeness are related to~\cite[Lemme~7.4]{Ecalle1} (stated there without proof). 
The proof in Chapter~\ref{ch:The Dominant Part and the Newton Polynomial} depends heavily on 
Chapter~\ref{ch:triangular automorphisms},
where we determine the invariants of certain automorphism groups of polynomial algebras in infinitely many variables $Y_0, Y_1, Y_2,\dots$ over a field of characteristic zero.

The function $\omega$ and the notion of $\upo$-freeness are 
closely related to second order linear differential equations
over $K$. More precisely~(Riccati), for~${y \in K^{\times}}$, ${4 y'' + fy = 0}$
is equivalent to $\omega (z) = f$ with $z := 2 y^{\dagger}$; so the
second-order linear differential equation $4 y'' + fy = 0$ reduces in a~way to
a first-order (but non-linear) differential equation 
$\omega (z) = f$. (The
factor~$4$ is just for convenience, to get simpler expressions below.)

\begin{exampleUnnumbered}
  The differential equation $y'' = - y$ has no solution $y \in \T^{\times}$,
  whereas the Airy equation $y'' = xy$ has two $\R$-linearly independent
  solutions in~$\T$~{\cite[Chapter~11, (1.07)]{Olver}}. Indeed, in
  Sections~\ref{sec:behupo} and~\ref{sec:special sets} we show that for
  $f \in \T$, the differential equation $4 y'' + fy = 0$ has a
  solution $y \in \T^{\times}$ if and only if $f < \upo_{\T}$, that is,
  $f < \upo_n=\frac{1}{\ell_0^2} + \frac{1}{\ell_0^2 \ell_1^2} + \cdots +
  \frac{1}{\ell_0^2 \ell_1^2 \cdots \ell_n^2}$ for some $n$. This fact
  reflects classical results {\cite{Hartman,Hille48}} on the question: for
  which logarithmico-exponential functions $f$ (in Hardy's sense) does the
  equation $4y'' + fy = 0$ have a non-oscillating real-valued solution (more
  precisely, a nonzero solution in a Hardy field)?
\end{exampleUnnumbered}

\subsection*{Newtonianity}

This is the most consequential elementary property of $\T$. An $\upo$-free
$H$-field~$K$ is said to be \textit{newtonian}\/ if every
$\d$-polynomial~$P(Y)$ over $K$ of Newton degree~$1$ has a zero in~$\mathcal{O}$. This turns out to be the correct analogue for valued
differential fields like $\T$ of the property of being henselian for a valued
field. We chose the adjective \textit{newtonian}\/ since it is this property that
allows us to develop in Chapter~\ref{ch:The Dominant Part and the Newton Polynomial} a Newton diagram method for differential polynomials. It is good to keep in mind that the role of newtonianity
in the results of Chapters~\ref{ch:newtonian fields}, \ref{ch:newtdirun}, and \ref{ch:QE} is more or less analogous to
that of henselianity in the theory of valued fields and as the key
condition in the Ax-Kochen-Er{\accentv s}ov results.

We already mentioned the result from~\cite{JvdH} that~$\T_{\operatorname{g}}$ is newtonian. 
That $\T$ is newtonian 
%In~\cite{vdH:PhD}, it has also been shown that $\T$ is newtonian. 
is a consequence of the following 
analogue in Chapter~\ref{ch:newtdirun} of the familiar valuation-theoretic fact that
spherically complete valued fields are henselian:

\begin{theoremintro}
  \label{duscnewt}If $K$ is an $H$-field, $\der K = K$, and $K$ is a directed
  union of spherically complete grounded $H$-subfields, 
 then~$K$ is \textup{(}$\upo$-free and\textup{)} newtonian.
\end{theoremintro}

\begin{exampleUnnumbered}
Let $K=\T$ 
and consider for $\alpha\in \R$ the differential polynomial 
$$P(Y)\ =\ Y'' - 2Y^3 - xY - \alpha \in \T\{Y\}.$$
For $\phi\in \T^{\times}$ we have
$(Y'')^\phi = \phi^2 Y''+\phi'Y'$ for $\phi\in \T^{\times}$, so
$$P^\phi\  =\ \phi^2 Y'' + \phi' Y' - 2Y^3 -xY-\alpha.$$
Now $\phi^2,\phi'\prec 1\prec x$ for active $\phi\prec 1$ in $\T^{>}$. Hence
$N_{P}\in \R^\times Y$, so $P$ has Newton degree~$1$. Thus the Painlev\'e II equation $y''=2y^3+xy+\alpha$ has a solution $y\in\mathcal{O}_{\T}$. (It is known that
$P$ has a zero $y\preceq 1$ in  the differential  subfield $\R(x)$ of $\T$ iff $\alpha\in\Z$; see for example \cite[Theorem~20.2]{GLS}.)
\end{exampleUnnumbered}

\noindent
The main results of Chapter~\ref{ch:newtonian fields} amount for
$H$-fields to the following:  

\begin{theoremintro}
  If $K$ is a newtonian $\upo$-free $H$-field with divisible value group,
  then~$K$ has no proper immediate $\d$-algebraic $H$-field extension.
\end{theoremintro}

%{corunnumbered}

\begin{corintro}
  Let $K$ be a real closed newtonian $\upo$-free $H$-field, and let $K^{\alg}
  = K [\imag]$ \textup{(}where $\imag^2 = - 1$\textup{)} be its algebraic
  closure. Then:
  \begin{enumerate}
    \item[\textup{(i)}] each $\d$-polynomial in $K^{\alg} \{Y\}$ of positive degree has a
    zero in $K^{\alg}$;
    
    \item[\textup{(ii)}] each linear differential operator in $K^{\alg} [\der]$ of positive
    order is a composition of such operators of order~$1$; 
    
    \item[\textup{(iii)}] each $\d$-polynomial in $K\{Y\}$ of odd degree has a zero in $K$; and
    
    \item[\textup{(iv)}] each linear differential operator in $K [\der]$ of positive order is
    a composition of such operators of order~$1$ and order~$2$.
  \end{enumerate}
\end{corintro}

\begin{theoremintro}
  \label{upodivnewtonization} If $K$ is an $\upo$-free $H$-field with divisible
  value group, then $K$ has an immediate $\d$-algebraic newtonian $H$-field
  extension, and any such extension embeds over $K$ into every $\upo$-free
  newtonian $H$-field extension of $K$.
\end{theoremintro}

\noindent
An extension of $K$ as in Theorem~\ref{upodivnewtonization} is
minimal over $K$ and thus unique up to isomorphism over $K$. We call such an extension a \textit{newtonization}\/ of $K$.

\begin{theoremintro}
  \label{newtonliouvilleclosure}If $K$ is an $\upo$-free $H$-field, then $K$
  has a $\d$-algebraic newtonian Liouville closed $H$-field extension that
  embeds over $K$ into every $\upo$-free newtonian Liouville closed $H$-field
  extension of $K$.
\end{theoremintro}

\noindent
An extension of $K$ as in
Theorem~\ref{newtonliouvilleclosure} is minimal over $K$ and thus unique up to isomorphism over $K$. We call such an extension a \textit{Newton-Liouville closure}\/ of $K$.

\subsection*{The main theorems}

We now come to the results in Chapter~\ref{ch:QE}, which in our view justify this volume. First, the various elementary conditions we have discussed
axiomatize a model complete theory. To be precise, construe $H$-fields in the
natural way as $\mathcal{L}$-structures where
$\mathcal{L}\ :=\ \{0, 1, {+}, {-}, {\, \cdot \,} , \der, {\leqslant},
   {\preccurlyeq}\}$, and 
let  $T^{\mathrm{nl}}$ be the
  $\mathcal{L}$-theory whose models are the newtonian $\upo$-free Liouville
  closed $H$-fields.

\begin{theoremintro}
  \label{nlmodelcompleteness}
  $T^{\mathrm{nl}}$ is model complete. 
  \end{theoremintro}

\noindent
The theory $T^{\mathrm{nl}}$ is not complete and has exactly two completions, namely $T^{\mathrm{nl}}_{\mathrm{small}}$ (small derivation) and $T^{\mathrm{nl}}_{\mathrm{large}}$ (large derivation). Thus newtonian $\upo$-free Liouville closed $H$-fields with small derivation have the same elementary properties as~$\mathbb T$. 

Every $H$-field with small derivation can be embedded into a  model of $T^{\mathrm{nl}}_{\mathrm{small}}$; thus Theorem~\ref{nlmodelcompleteness} yields the strong version of the model completeness conjecture from \cite{AvdD2} stated earlier in this introduction.
As $T^{\mathrm{nl}}_{\mathrm{small}}$ is complete and effectively axiomatized, it is
decidable. In particular,
there is an algorithm which, for any given $\d$-polynomials
$P_1, \ldots, P_m$ in indeterminates $Y_1, \ldots, Y_n$ with coefficients from~$\Z [x]$, decides whether there is a tuple $y \in \mathbb{T}^n$ such that ${P_1
(y) = \cdots = P_m (y) = 0}$. Such an algorithm with $\mathbb{T}$ replaced by
its differential subring $\R [[x^{- 1}]]$ is due to Denef and
Lipshitz~{\cite{DL84}}, but no such algorithm can exist with $\mathbb{T}$
replaced by $\R \(( x^{-1}\)) $ or by any of various other natural $H$-subfields of~$\mathbb{T}$~{\cite{AvdD3,GrigorievSinger}}. 

\medskip
\noindent
Theorem~\ref{nlmodelcompleteness} is the main step towards
an elimination of quantifiers, in a slightly
extended language: Let $\mathcal{L}^{\iota}_{\Upl, \Upo}$ be $\mathcal{L}$
augmented by the unary function symbol~$\iota$ and the unary predicates~$\Upl$,~$\Upo$, and extend $T^{\mathrm{nl}}$ to the $\mathcal{L}^{\iota}_{\Upl,
\Upo}$-theory $T^{\mathrm{nl}, \iota}_{\Upl, \Upo}$ by adding as defining
axioms for these new symbols the universal closures of
\begin{align*}
  {}\big[a \ne 0  & \longrightarrow a \cdot \iota (a) = 1\big]\  \&\ \big[a = 0
  \longrightarrow \iota (a) = 0\big],\\
  \Upl (a)\ & \longleftrightarrow\ \exists y \big[y \succ 1\ \&\ a = - y^{\dagger
  \dagger}\big],\\
  \Upo (a)\ & \longleftrightarrow\ \exists y \big[y \ne 0\ \&\ 4 y'' + ay = 0\big].
\end{align*}
For a model $K$ of $T^{\mathrm{nl}}$ this makes the sets $\Upl (K)$ and $\Upo
(K)$ downward closed with respect to the ordering of $K$. For example, for $f
\in \T$,
\begin{align*}
  f \in \Upl (\T)\ & \Longleftrightarrow\ f < \upl_n = \frac{1}{\ell_0} + \frac{1}{\ell_0
  \ell_1} + \cdots + \frac{1}{\ell_0 \ell_1 \cdots \ell_n}\  \text{ for some $n$,}\\
  f \in \Upo (\T)\ & \Longleftrightarrow\ f < \upo_n = \frac{1}{\ell_0^2} +
  \frac{1}{\ell_0^2 \ell_1^2} + \cdots + \frac{1}{\ell_0^2 \ell_1^2 \cdots
  \ell_n^2}\  \text{ for some $n$,}
\end{align*}
that is, $\Upl (\T)$ and $\Upo (\T)$ are the cuts in $\T$ determined by $\upl_\T$, $\upo_\T$ introduced earlier. 
We can now state what we view as the main result of this volume:

\begin{theoremintro}
  \label{intQE}The theory $T^{\mathrm{nl}, \iota}_{\Upl, \Upo}$ admits
  elimination of quantifiers.
\end{theoremintro}

\noindent
We cannot omit here either $\Upl$ or $\Upo$. In Chapter~\ref{ch:QE} we do
include for convenience one more unary predicate $\I$ in
$\mathcal{L}^{\iota}_{\Upl, \Upo}$: for a model $K$ of $T^{\mathrm{nl}}$ and
$a \in K$,
\[ \I (a)\ \longleftrightarrow\ 
\exists y \big[a \preccurlyeq y'\ \&\ y \preccurlyeq 1\big]\
   \longleftrightarrow\ a = 0 \vee \big[a \ne 0\ \&\ \neg \Upl (- a^{\dagger})\big], \]
where the first equivalence is the defining axiom for $\I$, and the second
shows that~$\I$ is superfluous in Theorem~\ref{intQE}. We note here that this
predicate $\I$ governs the solvability of first-order linear differential
equations \textit{with asymptotic side condition.}\/ More precisely, for $K$
as above and $f\in K$, $g, h \in K^\times$, the following are equivalent: \ifbool{PUP}{\enlargethispage{\baselineskip}}{}
\begin{enumerate}
  \item[(a)] there exists $y \in K$ with $y' = fy + g$ and $y \prec h$;
  \item[(b)] $\big[ (f - h^{\dagger}) \in \I (K) \text{ and } (g / h) \in \I
  (K) \big]  \text{ or } \big[ (f - h^{\dagger}) \notin \I (K) \text{ and } (g
  / h) \prec f - h^{\dagger} \big]$.
\end{enumerate}
This equivalence is part of Corollary~\ref{qelinI} and exemplifies
Theorem~\ref{intQE} (but is not derived from that theorem, nor used in its proof).

In the proof of Theorem~\ref{intQE}, and throughout the construction of
suitable $H$-field extensions, the predicates $\I$, $\Lambda$ and
$\Omega$ act as switchmen. Whenever a~fork in the road occurs due to the
presence of a gap $\upg$, then $\Iota (\upg)$ tells us to take the branch where $\int
\upg \preceq 1$, while $\neg \Iota(\upg)$ forces $\int \upg \succ 1$.
Likewise, the predicates $\Lambda$ and $\Omega$ control what happens when
adjoining elements $\upg$ and $\upl$ with $\upg^{\dag} = -
\upl$ and $\omega (\upl) = \upo$.

From the above defining axioms for $\Upl$ and $\Upo$ it is clear that these
predicates are (uniformly) existentially definable in models of
$T^{\mathrm{nl}}$. 
By model completeness of~$T^{\mathrm{nl}}$ they are also uniformly
\textit{universally}\/ definable in these models;
Section~\ref{sec:lHLO} deals with such algebraic-linguistic issues. 
%In view of this existential and universal definability, the model completeness of $T^{\mathrm{nl}}$ now follows from Theorem~\ref{intQE}. The completeness in Theorem~\ref{nlmodelcompleteness} also follows easily from this fundamental elimination result.

\medskip
\noindent
Next we list some more intrinsic consequences of our elimination theory.

\begin{corintro}\label{qe-top} Let $K$ be a newtonian $\upo$-free Liouville closed $H$-field, and suppose the set $X\subseteq K^n$ is definable. Then $X$ has empty interior in $K^n$ \textup{(}with respect to the order topology on $K$ and the product topology on
$K^n$\textup{)} if and only if for some nonzero $P\in K\{Y_1,\dots, Y_n\}$ we have $X\subseteq \big\{y\in K^n:\ P(y)=0\big\}$.
\end{corintro}

\noindent
In~(i) below the intervals are in the sense of the ordered field $K$.

\begin{corintro}
  \label{qe-cor}Let $K$ be a newtonian $\upo$-free Liouville closed $H$-field.
  Then:
  \begin{enumerate}
    \item[\textup{(i)}] $K$ is o-minimal at infinity: if $X \subseteq K$ is definable in
    $K$, then for some $a \in K$, either $(a, + \infty) \subseteq X$, or $(a,
    + \infty) \cap X = \emptyset$;
    
    \item[\textup{(ii)}] if $X \subseteq K^n$ is definable in $K$, then $X \cap C^n$ is
    semialgebraic in the sense of the real closed constant field $C$ of $K$;
    
    \item[\textup{(iii)}] $K$ has $\mathrm{NIP}$. \textup{(}See Appendix~\ref{app:modth} for
    this very robust property.\textup{)}
  \end{enumerate}
\end{corintro}

\noindent
It is hard to imagine obtaining these results for $K = \T$ without Theorem~\ref{intQE}. 
Item~(i) relates to classical bounds on solutions of algebraic differential equations over Hardy fields; see \cite[Section~3]{AvdD3}.
To illustrate item (ii) of Corollary~\ref{qe-cor}, we note that the set of real parameters $(\lambda_0,\dots,\lambda_n)\in\R^{n+1}$ for which the system
$$\lambda_0 y+\lambda_1 y'+\cdots+\lambda_n y^{(n)}=0,\qquad 0\neq y \prec 1$$
has a solution in $\T$ is a semialgebraic subset of $\R^{n+1}$; in fact, it agrees with the
set of all  $(\lambda_0,\dots,\lambda_n)\in\R^{n+1}$ such that the polynomial
$\lambda_0 +\lambda_1 Y+\cdots+\lambda_n Y^n\in \R[Y]$ has a negative zero in $\R$; see Corollary~\ref{cor:linconstcoeff, Schwarz closed}.
To illustrate item (iii), let $Y=(Y_1,\dots,Y_n)$ be a tuple of distinct differential indeterminates; 
for an $m$-tuple $\sigma=(\sigma_1,\dots,\sigma_m)$ of elements of $\{ {\prec}, {\asymp}, {\succ} \}$ we
say that $P_1,\dots,P_m\in \T\{Y\}$ \textit{realize~$\sigma$}\/ if there exists~$a\in\T^n$ such that $P_i(a)\,\sigma_i\, 1$  holds for $i=1,\dots,m$.
Then a special case of (iii) says that for fixed $d,n,r\in\N$,
the number of tuples $\sigma\in\{ {\prec}, {\asymp}, {\succ} \}^m$ 
realized by some $P_1,\dots,P_m\in \T\{Y\}$ of degree at most~$d$ and order at most~$r$
grows only polynomially
with~$m$, even though the total number of tuples is $3^m$. 
These manifestations of (ii) and (iii), though instructive, are perhaps a bit misleading, since they can be obtained without
appealing to (ii) and (iii).

\medskip
\noindent
In the course of proving Theorem~\ref{nlmodelcompleteness} we
also get:

\begin{theoremintro}
  \label{intnlcon}If $K$ is a newtonian $\upo$-free Liouville closed
  $H$-field, then $K$ has no proper $\d$-algebraic $H$-field extension with
  the same constant field.
\end{theoremintro}

\noindent
For $K=\T_{\operatorname{g}}$ this yields:
every $f\in\mathbb T\setminus\mathbb T_{\operatorname{g}}$ is $\d$-transcendental over $\mathbb T_{\operatorname{g}}$.

\medskip\noindent
We can also enlarge $\T$.
For example, the series $\sum_{n = 0}^{\infty}
e_n^{- 1}$, with $e_n$ the $n$th iterated exponential of $x$,
does not lie in $\T$ but does lie in a~certain completion~$\T^{\operatorname{c}}$ of
$\T$. This completion $\T^{\operatorname{c}}$ is naturally an ordered valued differential field extension of $\T$, and by Corollary~\ref{uponewlicomp} we have $\T \preccurlyeq \T^{\operatorname{c}}$.

\ifbool{PUP}{
\section*{Organization}
\addcontentsline{toc}{section}{Organization}
}{\section*{Organization}}

\noindent
Here we discuss the somewhat elaborate organization of this volume into chapters, some technical ingredients not mentioned so far, and some material that goes beyond the
setting of $H$-fields. Indeed, the supporting algebraic theory deserves to be developed in a broad way, and there are more notions to keep track of than one might expect.    

\subsection*{Background chapters}
To make our work more accessible and self-contained, we provide
in the first five chapters background on commutative algebra, valued
abelian groups, valued fields, differential fields, and linear differential
operators. This material has many
sources, and we thought it would be convenient to have it available all in one
place. In addition we have an appendix with the construction of
$\T$, and an appendix exposing the (small) part of model theory that we 
need.

The basic facts on Hahn products, pseudocauchy sequences and spherical
completeness in these early chapters are used throughout
the volume. Some readers might prefer to skip in a~first
reading cauchy sequences, completeness (for valued abelian groups and valued
fields) and step-completeness, which are not needed for the main results in
this volume (but see Corollary~\ref{uponewlicomp}). Some parts, like
Sections~\ref{sec:valued vector spaces} and~\ref{sec:systems}, fit naturally
where we put them, but are mainly intended for use in the next volume. On the other
hand, Section~\ref{Compositional Conjugation} on compositional conjugation is
elementary and frequently referred to in subsequent chapters, but this material seems virtually absent from the literature.

\subsection*{Valued differential fields}
We also profited from examining arbitrary valued differential fields $K$ with
small de\-ri\-va\-tion, that is, $\der \smallo \subseteq \smallo$ for the maximal
ideal $\smallo$ of the valuation ring~$\mathcal{O}$ of $K$. This yields the
continuity of the derivation $\der$ with respect to the valuation topology and
gives $\der \mathcal{O} \subseteq \mathcal{O}$, and so induces a derivation on
the residue field. To our surprise, we could establish in
Chapters~\ref{ch:valueddifferential} and \ref{sec:dh1} some useful facts in
this very general setting when the induced derivation on the residue field is nontrivial,
for example the Equalizer Theorem~\ref{theq}. We need this result in deriving
an ``eventual'' version of it for $\upo$-free $H$-fields in
Chapter~\ref{ch:The Dominant Part and the Newton Polynomial}, which in turn is
crucial in obtaining our main results, via its role in constructing an
appropriate Newton diagram for $\d$-polynomials.

\subsection*{Asymptotic couples}
A useful gadget is the \textit{asymptotic couple}\/ of an $H$-field $K$. This
is the value group~$\Gamma$ of $K$ equipped with the map $\gamma \mapsto
\gamma^{\dagger} \colon \Gamma^{\ne} \to \Gamma$ defined by: if $\gamma = vf$, $f
\in K^{\times}$, then $\gamma^{\dagger} = v (f^{\dagger})$. This map is a
valuation on $\Gamma$, and we extend it to a map $\Gamma \to \Gamma_{\infty}$
by setting $0^{\dagger} := \infty$. Two key facts are that
$\alpha^{\dagger} < \beta + \beta^{\dagger}$ for all $\alpha, \beta > 0$ in
$\Gamma$, and $\alpha^{\dagger} \geqslant \beta^{\dagger}$ whenever $0 <
\alpha \leqslant \beta$ in $\Gamma$. The condition on an $H$-field of having
small derivation can be expressed in terms of its asymptotic couple; the same
holds for having a gap, for being grounded, and for having asymptotic integration, but not for being $\upo$-free.

Asymptotic couples were introduced by Rosenlicht~{\cite{Rosenlicht2}} for $\d$-valued fields. In Chapter~\ref{ch:valueddifferential} we assign to {\em any\/} valued differential field with small derivation an asymptotic couple, with good effect. Asymptotic couples play also an important role in
Chapters~\ref{ch:asymptotic differential fields}, \ref{ch:H}, \ref{evq}, \ref{ch:The Dominant Part and the Newton Polynomial}, and \ref{ch:QE}.

\subsection*{Differential-henselian fields}
Valued differential fields with small derivation include the so-called
\textit{monotone}\/ differential fields defined by the condition $a' \preccurlyeq a$. In
analogy with the notion of a \textit{henselian}\/ valued field,
Scanlon~\cite{Scanlon} introduced \textit{differential-henselian}\/ monotone
differential fields. Using the Equalizer Theorem we extend this notion and
basic facts about it to arbitrary valued differential fields with small
derivation in Chapter~\ref{sec:dh1}. (We abbreviate
\textit{differential-henselian}\/ to \textit{$\d$-henselian.}\/) This
material plays a role in Chapter~\ref{ch:newtonian fields}, using the
following relation between \textit{$\d$-henselian}\/ and
\textit{newtonian}\/: an $\upo$-free $H$-field $K$ is newtonian iff for every
active $\phi \in K^{>}$ the compositional conjugate $K^{\phi}$ is
$\d$-henselian, with the valuation~$v$ on $K^{\phi}$ replaced by the coarser
valuation~$\pi \circ v$ where $\pi \colon v (K^{\times}) = \Gamma \to \Gamma /
\Delta$ is the canonical map to the quotient of $\Gamma$ by its convex
subgroup $$\Delta\ :=\ \{\gamma \in \Gamma :\ \gamma^{\dagger} > v
\phi\}.$$
We pay particular attention to two special cases: $v (C^{\times}) = \{0\}$
(few constants), and $v (C^{\times}) = \Gamma$ (many constants). The first
case is relevant for newtonianity, the second case is considered in a short
Chapter~\ref{ch:monotonedifferential}, where we present Scanlon's extension of
the Ax-Kochen-Er{\accentv s}ov theorems to $\d$-henselian valued fields with many
constants, and add some things on definability.

While $\d$-henselianity is defined in terms of solving differential
equations in one unknown, it implies the solvability of suitably non-singular
systems of~$n$ differential equations in $n$ unknowns: this is shown at the
end of Chapter~\ref{sec:dh1}, and has a nice consequence for newtonianity:
Proposition~\ref{newtalgpres}.

%Perhaps this is a good point to mention that besides~\cite{Scanlon}, there are only scant results on
%the model-theoretic properties of valued differential fields exhibiting some nontrivial interaction between valuation and
%derivation.
%(The basic model theory of ordered differential fields, where no interaction between ordering and derivation is assumed,
%was established by Singer~\cite{Singer78-1}; these results were generalized by
%Guzy and Point~\cite{Guzy05,GuzyPoint10} to other expansions of differential fields, including
%valued differential fields.)

\subsection*{Asymptotic differential fields}
To keep things simple we confined most of the exposition above to $H$-fields, but
this setting is a bit too narrow for various technical reasons. For example,
a~differential subfield of an $H$-field with the induced ordering is not
always an $H$-field, and passing to an algebraic closure like $\T
[\imag]$ destroys the ordering, though $\T [\imag]$ is still a $\d$-valued
field. On occasion we also wish to change the valuation of an $H$-field or
$\d$-valued field by coarsening. For all these reasons we introduce in
Chapter~\ref{ch:asymptotic differential fields} the class of
\textit{asymptotic differential fields}\/, which is larger and more flexible
than Rosenlicht's class of $\d$-valued fields. Many basic facts about
$H$-fields and $\d$-valued fields do have good analogues for asymptotic
differential fields. This is shown in Chapter~\ref{ch:asymptotic differential
fields}, which also contains a lot of basic material on asymptotic couples. 
Chapter~\ref{ch:H} deals more specifically with $H$-fields.

%allow ourselves to refer to {\cite{AvdD2,AvdD3}} for some
%facts whose proofs go through with only minor changes in this generalized
%setting. Except for these references and one appeal 
%in Chapter~\ref{ch:QE} to
%a result in~\cite{AvdD}, we have avoided explicit 
%dependence on our earlier published work in this area.

\subsection*{Immediate extensions}
Indispensable for attaining our main results is the fact that every $H$-field
with divisible value group and with asymptotic integration has a spherically
complete immediate $H$-field extension. This is part of
Theorem~\ref{thm:immediate}, and proving it about five years ago removed a
bottleneck. It provides the only way known to us of extending every $H$-field
to an $\upo$-free $H$-field. Possibly more important than
Theorem~\ref{thm:immediate} itself are the tools involved in its proof. In
view of the theorem's content, it is ironic that models of $T^{\operatorname{nl}}$ are
never spherically complete, in contrast to all prior positive results on
elementary theories of valued fields with or without extra structure,
cf.~\cite{AxKochen,AxKochen2,BMS,Ersov,Scanlon}.

\subsection*{The differential Newton diagram method}
Chapters~\ref{ch:The Dominant Part and the Newton Polynomial} and \ref{ch:newtonian fields} present the differential Newton diagram method in the
general context of asymptotic fields that satisfy suitable technical
conditions, such as $\upo$-freeness.
Before tackling these chapters,
the reader may profit from first studying our exposition of the Newton diagram method
for ordinary one-variable polynomials over henselian valued fields of equicharacteristic zero in Section~\ref{Newton Diagrams}.
Some of the issues encountered there (for example, the \textit{unraveling}\/ technique) appear again, albeit in more
intricate form, in the differential context of these chapters.
In the proofs of a few crucial facts about the special
cuts $\upl$ and $\upo$ in Chapter~\ref{ch:The Dominant Part and the Newton Polynomial} we use some results from
the preceding Chapter~\ref{ch:triangular automorphisms} on triangular automorphisms.
Chapter~\ref{ch:triangular automorphisms} is a bit special in being essentially independent of the earlier chapters.
 
\subsection*{Proving newtonianity}
Chapter~\ref{ch:newtdirun} contains the proof of Theorem~\ref{duscnewt}, and thus establishes that $\T$ is a model of
our theory $T^{\operatorname{nl}}_{\operatorname{small}}$. This theorem is also useful in other contexts: In~\cite{BM},
Berarducci and Mantova construct a derivation on Conway's field~{\bf No} of surreal numbers~\cite{Con76,Gon86}
turning it into a Liouville closed $H$-field with constant field~$\R$. From Theorem~\ref{duscnewt} and the completeness
of $T^{\operatorname{nl}}_{\small}$ it follows that
{\bf No} with this derivation and $\mathbb T$ are elementarily equivalent, 
as we show in \cite{ADHNo}.

\subsection*{Quantifier elimination}
In Chapter~\ref{ch:QE} we first prove Theorem~\ref{nlmodelcompleteness}
on model completeness, next we consider $H$-fields equipped with a $\Upl\Upo$-structure, and then deduce Theorem~\ref{intQE} about quantifier elimination with various interesting consequences,
such as Corollaries~\ref{qe-top} and~\ref{qe-cor}. The introduction to this chapter gives an overview of the proof and the role of various embedding and extension results in it.

\ifbool{PUP}{
\section*{The Next Volume}
\addcontentsline{toc}{section}{The Next Volume}
}{\section*{The Next Volume}}

\noindent
The present volume focuses on achieving quantifier elimination
(Theorem~\ref{intQE}), and so we left out various things we did since~1995 
that were not needed for that. In
a~second volume we intend to cover these things as required for developing our work further.
Let us briefly survey some highlights of what is to come.

\subsection*{Linear differential equations}
We plan to consider linear differential equations in much 
greater
detail, comprising the corresponding differential Galois theory, in connection
with constructing the linear surjective closure of a differential field,
factoring linear differential operators over suitable algebraically closed
$\d$-valued fields, and explicitly constructing the Picard-Vessiot extension
of such an operator.
Concerning the latter, the complexification
$\T [ \imag ]$ of $\T$ is no longer closed under
exponential integration, since oscillatory ``transmonomials'' such as
$\ex^{\imag x}$ are not in~$\T \left[ \imag \right]$. Adjoining these oscillatory transmonomials to $\T [ \imag ]$ yields a
$\d$-valued field that
contains a Picard-Vessiot extension of~$\T$ for each operator in $\T[\der]$. 

\subsection*{Hardy fields}
We also wish to pay more attention to Hardy fields, and this will bring
up analytic issues. For example, every Hardy field
containing~$\R$ can be shown to extend to an $\upo$-free Hardy field. Using methods
from~{\cite{vdH:hfsol}}, we also hope to prove that it always extends to a
newtonian $\upo$-free Hardy field. Indeed, that paper proves
among other things the following pertinent result (formulated here with our
present terminology): Let $\T_{\operatorname{g}}^{\operatorname{da}}$ consist of the grid-based
transseries that are $\d$-algebraic over~$\R$. Then $\T_{\operatorname{g}}^{\operatorname{da}}$ is a
newtonian $\upo$-free Liouville closed $H$-subfield of $\T_{\operatorname{g}}$ and is isomorphic
over~$\R$ to a Hardy field containing $\R$.

\subsection*{Embedding into fields of transseries}
Another natural question we expect to deal with is whether
every $H$-field can be given some kind of transserial structure.
This can be made more precise in terms of the axiomatic definition of a
\textit{field of transseries}\/ in terms of a transmonomial group 
$\mathfrak{M}$ in Schmeling's
thesis~{\cite{Schm01}}. For instance, one axiom there is that for all
$\mathfrak{m} \in \mathfrak{M}$ we have
$\operatorname{supp} \log \mathfrak{m} \subseteq \mathfrak{M}^{\succ}$. We hope that any $H$-field can be embedded
into such a field of transseries. This would be a natural
counterpart of Kaplansky's theorem~\cite{Kaplansky} embedding
certain valued fields
into Hahn fields, and would make it possible to
think of $H$-field elements as generalized transseries.

%\subsection*{Variants of $\T$}
%We already mentioned that variants of the
%construction of~$\T$ lead to various interesting differential %subfields of~$\T$, such as $\T_{\operatorname{g}}$.
%In the next volume we intend to elaborate on the model-%theoretic properties of transseries with various other restrictions 
%on the support.  
%We intend to elaborate on transseries with various restrictions  on the support and to consider also the complex transseries from~{\cite{vdH:osc}}.

\subsection*{More on the model theory of $\T$}
In the second volume we hope to deal with further issues around $\T$
of a model-theoretic nature: for example, identifying the induced structure on
its value group (conjectured to be given by its $H$-couple, as specified in
{\cite{AvdD}}); and determining the definable closure of a subset of a model
of~$T^{\operatorname{nl}}$, in order to get a handle on what functions are definable
in $\T$. 

A by-product of the present volume is a full description of several
important $1$-types over a given model of $T^{\operatorname{nl}}$, but the entire
space of such $1$-types remains to be surveyed.
Theorem~\ref{intnlcon} suggests that the model-theoretic notions of \textit{non-orthogonality to~$C$}\/ or \textit{$C$-internality}\/ may be significant for models of $T^{\operatorname{nl}}$; 
see also \cite{ADHdim}.

\ifbool{PUP}{
\section*{Future Challenges}
\addcontentsline{toc}{section}{Future Challenges}
}{\section*{Future Challenges}}

\noindent
We now discuss a few more open-ended avenues of inquiry.

%We already mentioned some issues that would be nice to resolve. Here are some more.

\subsection*{Differentiation and exponentiation}
 The restriction to $\smallo_{\T}$ of the exponential
function on $\T$ is easily seen to be definable in~$\T$, but by part (ii) of Corollary~\ref{qe-cor}, the restriction to $\R$ of this exponential
function is not definable in $\T$. This raises the question whether our
results can be extended to the differential field $\T$ 
{\em with exponentiation}, or
with some other extra o-minimal structure on it.

\subsection*{Logarithmic transseries}
A transseries is \textit{logarithmic}\/ if all transmonomials in it are of the
form $\ell_0^{r_0} \cdots \ell_n^{r_n}$ with $r_0, \ldots, r_n \in \R$. (See Appendix~\ref{app:trans}.) The
logarithmic transseries make up an $\upo$-free newtonian $H$-subfield
$\T_{\log}$ of $\T$ that is not Liouville closed. We conjecture that
$\T_{\log}$ as a valued differential field is model complete.
The asymptotic couple of $\T_{\log}$ has been successfully analyzed by Gehret~\cite{Geh},
and turns out to be model-theoretically tame, in particular, has NIP~\cite{Geh2}. (There is also the notion of a transseries being
\textit{exponential.}\/ The exponential transseries form a
real closed $H$-subfield $\T_{\exp}$ of $\T$ in which the set~$\Z$ is
existentially definable, see~\cite{AvdD3}. 
It follows that the differential field~$\T_{\exp}$ does not have a reasonable model theory: it is as complicated as
so-called \textit{second-order arithmetic.}\/)

\subsection*{Accelero-summable transseries}
The paper \cite{vdH:hfsol} on transserial Hardy fields yields on the one hand a method to associate a genuine function to a suitable formal transseries, and in the other direction also provides means to associate concrete
asymptotic expansions to elements of Hardy fields. We expect that more can be done in this direction.

{\'E}calle's theory of analyzable functions has a more canonical
procedure that associates a function to an {\em accelero-summable\/}
transseries. These transseries make up an $H$-subfield~$\T_{\operatorname{as}}$ of~$\T$. This procedure has the advantage that it does not only preserve
the ordered field structure, but also composition, functional inversion, and
several other operations. In its full generality, however, {\'E}calle's theory requires sophisticated analytic tools, and is beyond the scope of this volume.
It is clear that
$\T_{\operatorname{as}}$ is analytically more important than $\T$, but
the latter might help in understanding the former. The $H$-subfield  $\T_{\operatorname{as}}$ of $\T$ contains $\R$, is $\upo$-free and Liouville closed. Is it newtonian? In view of Theorem~\ref{intnlcon}, a positive answer would 
confirm \'Ecalle's belief~\cite[p.~148]{Ecalle1} that any
solution in $\T$ of an algebraic differential equation
$P(Y)=0$ over $\T_{\operatorname{as}}$ with $P\ne 0$ lies
in~$\T_{\operatorname{as}}$. 

\subsection*{Beyond $H$-fields} The derivation of a
differentially closed field $K$ cannot be continuous with respect to a nontrivial valuation on $K$; see Section~\ref{Miscellaneous Facts about Asymptotic Fields}.
This sets a limit for the study of valued differential fields with a reasonable interaction between valuation and derivation.
However, one may close off the
$\d$-valued field~$\T[\imag]$ under exponential integration,
by adding oscillatory
transmonomials  recursively. This results in valued differential fields of \textit{complex
transseries}\/ over which a version of the Newton diagram method for $\T_{\operatorname{g}}$ goes through; see~\cite{vdH:osc}. It would be interesting to find out more about the model theory of these rich valued differential fields.

\ifbool{PUP}{
\section*{A Historical Note on Transseries}
\addcontentsline{toc}{section}{A Historical Note on Transseries}
}{\section*{A Historical Note on Transseries}}

\noindent
The differential field of transseries was first defined and extensively used 
in \'Ecalle's solution of Dulac's problem, which is about
plane analytic vector fields. Its solution shows in particular
that a plane polynomial vector field admits only a finite number of limit 
cycles. 
%yields a partial answer to weak version of the second part of 
%Hilbert's 16th Problem \cite{Hilbert1900} asks whether any polynomial 
%vector field in the real plane admits only a finite number of limit cycles.
(A \textit{limit cycle}\/ of a planar vector field is a periodic
trajectory with an annular neighborhood not containing any other periodic trajectory.)
 It was long believed that in 1923 Dulac~\cite{Dulac23}
had given a proof of this finiteness statement, until Il$'$yashenko~\cite{Ilyashenko81} found a gap in~1981: Dulac was operating with asymptotic expansions of germs of functions as if they faithfully represented these germs. To justify this in Dulac's case is not easy: it was done
independently
by \'Ecalle~\cite{Ecalle1} and Il$'$yashenko~\cite{Ilyashenko91},
and required fundamental new ideas (and hundreds of pages). We briefly sketch here the role of transseries in \'Ecalle's approach.

\medskip
\noindent
Suppose towards a contradiction that some polynomial vector field 
on $\R^2$ has infinitely many limit cycles. Classical facts about planar vector fields such as the Poincar\'e-Bendixson Theorem 
allow us to reduce to the case where infinitely many of these limit cycles accumulate at a so-called
\textit{polycycle}\/ of the vector field; see~\cite[Theorem~24.22]{IlYa} for details. Such a polycycle~$\sigma$ consists of finitely
many trajectories    
$$S_1\to S_2,\ \dots\ ,S_{r-1}\to S_r,\ S_r\to S_1 \qquad\text{(the edges)}$$ between singularities 
$S_1,\dots,S_r$ (the vertices) of the vector field; see Figure~\ref{fig:limit cycle} where $r=3$. 
%Repetitions among the vertices $S_1,\dots,S_r$ are allowed whereas the edges are distinct.
Draw lines $\ell_1,\dots,\ell_r$ that cross these edges $S_r\to S_1,\,S_1\to S_2,\,\dots,\,$ 
$S_{r-1}\to S_r$ transversally at points $O_1,\dots, O_r$. For
any trajectory~$\varphi$ of the vector field that is
sufficiently close to $\sigma$ we consider the successive
points where~$\varphi$ meets $\ell_1,\ell_2,\dots,\ell_r,\ell_1$ and denote their
distances to $O_1,\dots, O_r, O_1$ by $t_1,t_2,\dots,t_r,t_{r+1}$. 
The behavior of the vector field near $S_i$ yields for some $\varepsilon_i>0$ a real analytic function
  $g_i\colon (0,\varepsilon_i)\to (0,\infty)$ 
such that $t_{i+1} = g_i(t_{i})$. %Each~$\Phi_i$ can be computed
%essentially in terms of the local behavior 
%of the vector field at $S_i$.
We have $g_i(t)\to 0$ as $t\to 0^+$ but $g_i$ does not necessarily extend 
analytically to~$0$.

\begin{figure} 
\begin{center}
\includegraphics{mt-dulac.eps}%[natwidth=68.5mm,natheight=57mm]{mt-dulac.eps}
\end{center}
\caption{A polycycle $\sigma$ and a close trajectory $\varphi$.}
\label{fig:limit cycle}
\end{figure}

\medskip
\noindent
The composition $f:=g_r\circ\cdots\circ g_1$ is defined on some
interval $(0,\varepsilon)$ and is called the \textit{Poincar\'e return map}\/ 
of the polycycle~$\sigma$ (relative to our choice of $\ell_1,\dots,\ell_r$). We have $t_{r+1}=f(t_1)$, for trajectories  close enough to $\sigma$.
Thus $f(t)=t$ corresponds to a periodic trajectory, and so it suffices to show that either $f$ is the identity,
or $f(t)\neq t$ for
all sufficiently small~$t>0$. One can even ask whether this non-oscillation 
property holds for Poincar\'e return maps of polycycles of plane 
{\em analytic\/} (not necessarily polynomial) vector fields; this is {\em Dulac's problem}. 

\medskip
\noindent
It is convenient to work at infinity by setting $x=t^{-1}$ and replace these functions $f,g_1,\dots, g_r$
of $t$ by functions $F, G_1,\dots, G_r$ of $x$ with
$F=G_r\circ \dots \circ G_1$. Dulac~\cite{Dulac23} provides formal series expansions $\tilde{G_i}$ of the $G_i$, which are rather simple transseries, usually divergent, and which by formal composition yields an often complicated 
transseries expansion 
$\tilde{F}=\tilde{G_r}\circ \cdots \circ \tilde{G_1}$ of $F$. 
 
\'Ecalle is able to reconstitute the germs $G_i$ and $F$ from their formal counterparts~$\tilde{G_i}$ and $\tilde{F}$
by developing a delicate analytic machinery of 
{\em accelero-summation}.
More precisely, he constructs an (accelero-summation) operator 
$\tilde{G}(x) \mapsto G$ whose domain of definition is a certain 
differential subfield $\T_{\operatorname{as}}$ of $\T$ and whose values are germs of real analytic functions at $+\infty$;
it assigns to each $\tilde{G_i}$ the germ $G_i$.
Moreover, $\T_{\operatorname{as}}$ is closed under composition, and accelero-summation preserves real constants, addition, multiplication, differentiation, composition, and the (total) field ordering: 
if  $\tilde{G}(x)\in \T_{\operatorname{as}}$, $\tilde{G}(x)>0$,
then $G(x)>0$ for all sufficiently large real $x$; here the $x$ in 
$\tilde{G}(x)$ is an indeterminate, while in $G(x)$
it ranges over real numbers. 
Applying this operator to $\tilde{F}(x) -x$ yields the desired result:  
either $F(x)=x$ for all large enough $x$, or
$F(x)<x$ for all large enough $x$, or $F(x)>x$ for all 
large enough $x$.

 Accelero-summation is very powerful---the solution of the Dulac problem is just one application---and much of \'Ecalle's book~\cite{Ecalle1} consists of
developing it in various directions. Unorthodox summations occur already in Euler's study of divergent series~(\cite[p.~220]{Euler1760}, see also \cite[\S{}5.3]{Varadarajan}), but even the remarkable generalization of this work by Emil Borel more than a century later~\cite{Borel} is not adequate to reconstitute the above $F$ from $\tilde{F}$.

\medskip\noindent
 \'Ecalle~\cite[\S{}1.9]{Ecalle1} also indicates another approach to 
Dulac's problem in which acce\-le\-ro-summation would play a smaller
role. It involves the 
group of
formal Laurent series $x(1+a_1 x^{-1}+a_2x^{-2}+\cdots)\in \T$  (with respect to composition) and the group
$\{ \ldots, \log \log x, \log x,
x, \ex^x, \ex^{\ex^x}, \ldots \}$ generated under composition by $\ex^x\in \T$. 
An intriguing open question is whether there exist nontrivial relations
between these two groups. In other words, is the group they generate under 
composition their free product? 
% the existence of 
%nontrivial relations between these groups.
% $f_1 \circ g_1 \circ \cdots \circ f_r \circ g_r = x$ with
%$f_i \in \mathcal{S} \setminus \{ x \}$, $g_i \in \mathcal{E} \setminus \{ x
%\}$ and $r\ge 1$.

\medskip
\noindent
Dulac's problem is often mentioned in connection with
Hilbert's 16th Problem, whose second part asks for a uniform bound (only depending on the degrees of the polynomials involved) on the number of limit cycles of a polynomial vector field in~$\R^2$. This remains open, and is part of Smale's list~\cite{Smale} of mathematical problems ``for the next century.''
 
\medskip
\noindent
The exponential field of transseries (without the derivation) was also introduced independently by Dahn and G\"o\-ring~\cite{DG}.
Motivated by Tarski's problem on the real exponential field, they saw $\T$ as a candidate for a non-standard model of the theory of this structure. 
This idea was vindicated
in~\cite{DMM1}, in the wake of Wilkie's solution~\cite{Wilkie} of the ``geometric'' part of Tarski's problem.

%% file: mt-1.tex
\chapter{Some Commutative Algebra}\label{ch:commalg}

\noindent
This chapter will enable us to give
a self-contained proof of 
Johnson's Theorem~\ref{thm:johnson} on regular solutions of systems of 
algebraic differential equations. (Johnson's result will be used in Chapter~\ref{sec:dh1}
on differential-henselian fields.) Sections~\ref{sec:regular local rings},~\ref{sec:differentials}, and~\ref{sec:derivations on field exts} below contain the facts on regular local rings and 
K\"ahler differentials needed for Theorem~\ref{thm:johnson}. 
We also include material needed
in Chapter~\ref{ch:valuedfields}, as well as material that would 
otherwise be ad hoc
parts of Chapters~\ref{ch:differential polynomials} and~\ref{ch:lindifpol}. Nothing in the present chapter is new; our aim is a clear and efficient exposition. 

\medskip\noindent 
We recall here a common notational convention, to be used throughout this chapter and later in these notes. 
Let $R$ be a commutative ring and $M$ an $R$-module. When~$U$
and~$V$ are given as additive subgroups of $R$ and $M$, respectively, then we set 
$$UV\ :=\ \left\{\sum_{i=1}^n r_i x_i\in M:\ r_1,\dots, r_n\in U,\ x_1,\dots, x_n\in V \right\},$$ 
the additive subgroup of $M$ generated by the products $rx$
with $r\in U$ and $x\in V$. Let~$I$ be an ideal of $R$.
Then by this convention $IM$ is the submodule of $M$ generated by the $rx$ with $r\in I$ and $x\in M$. If $A$ is a
commutative $R$-algebra, then~$IA$ turns out to be the ideal of $A$ generated by the image of $I$ in $A$. Construing~$R$ as an $R$-module as usual and given also an ideal $J$ of 
$R$, our convention yields the ideal $IJ\subseteq I\cap J$ of 
$R$ generated by the products $rs$
with $r\in I$ and $s\in J$. This allows us to define the ideal~$I^n$ of $R$ recursively by 
$I^{0}:= R$ and $I^{n+1}:=I^nI$. 

This notation applies also to subrings $K$ and $L$ of $R$, and then $KL$ is the subring~$K[L]=L[K]$ of $R$ generated by $K\cup L$. Often $R$ is a field and $K$,~$L$ are subfields of $R$
with a common subfield $C$ of $K$ and $L$ such that
$K$ or $L$ is algebraic over $C$; in that case $KL$ is a subfield of $R$, called the {\em compositum of $K$ and $L$}. 

In a few places $U$ and $V$ are only given as {\em subsets\/} of $R$ (not necessarily additive subgroups), and then $UV$ denotes just the set of products $rs$ with $r\in U$, $s\in V$. 

\index{field!compositum}
\nomenclature[A]{$KL$}{compositum of the fields $K$ and $L$}

\section{The Zariski Topology and Noetherianity}\label{sec:ztn}

\noindent
{\em Throughout this section $R$ is a commutative ring}. The set of prime ideals of 
$R$ is called the {\bf spectrum} of $R$, denoted by $\Spec(R)$. So for each
$\frak{p}\in \Spec(R)$ we have the integral domain $R/\frak{p}$.
Below we shall view the elements of $R$ as functions on~$\Spec(R)$, and make 
$\Spec(R)$ into a space, with a
{\em Nullstellensatz\/} for these functions.  This space is particularly 
nice when $R$ is noetherian.  

\index{ring!spectrum}
\index{spectrum}
\nomenclature[A]{$\Spec(R)$}{set of prime ideals of the commutative ring $R$}

\subsection*{The Zariski topology} Let $f\in R$ and $\frak{p}\in \Spec(R)$,
and think of the residue class $f+ \frak{p}\in R/\frak{p}$ as the 
{\em value of 
$f$ at the point $\frak{p}$}. Then ``$f$ vanishes at $\frak{p}$'' just 
means that $f\in \frak{p}$.   
 For a set $S\subseteq R$, the subset
$$\Zero(S)\ :=\ \big\{ \mathfrak p\in\Spec(R):\ \mathfrak p\supseteq S\big\}$$
of $\Spec(R)$ is then the set of points in $\Spec(R)$ at which
all $f\in S$ vanish: the set of common zeros of the 
``functions'' $f\in S$. Thus $\Zero(\emptyset)=\Spec(R)$, and 
$\Zero(R)=\emptyset$. 
For $f_1,\dots,f_n\in R$ we also write $\Zero(f_1,\dots,f_n)$ instead of 
$\Zero\!\big(\{f_1,\dots,f_n\}\big)$. 
If~$I$ is the ideal of $R$ generated by $S\subseteq R$, 
then $\Zero(S)=\Zero(I)$.
The sets $\Zero(S)$ with $S\subseteq R$
are the closed sets of a topology on $\Spec(R)$, called its 
{\bf Zariski topology}: \index{topology!Zariski} use that
$$\textstyle\bigcap_{\lambda\in\Lambda} \Zero(S_\lambda)\ =\ \Zero\left(\bigcup_{\lambda\in\Lambda} S_\lambda\right)$$
for any family $(S_\lambda)_{\lambda\in\Lambda}$ of subsets of $R$, and that
for $S_1,S_2\subseteq R$,
$$\Zero(S_1)\cup\Zero(S_2)\ =\ \Zero(S_1S_2), \qquad S_1S_2\ :=\ \{f_1f_2: f_1\in S_1,\ f_2\in S_2\}.$$
The sets $\operatorname{D}(f):=\Spec(R)\setminus \Zero(f)$ with $f\in R$
form a base of open sets for this Zariski topology.
The Zariski topology on $\Spec(R)$ is not in general hausdorff. 

\medskip\noindent
Let $\varphi\colon R\to S$ be a morphism of commutative rings. Then for each $\mathfrak q\in\Spec(S)$ we have
$\varphi^{-1}(\mathfrak q)\in\Spec(R)$, giving rise to an 
inclusion-preserving map
$$\varphi^*\colon \Spec(S)\to\Spec(R), \qquad \varphi^*(\mathfrak q):=\varphi^{-1}(\mathfrak q)$$
with $(\varphi^*)^{-1}\big(\!\Zero(I)\big)=\Zero\!\big(\varphi(I)\big)$ for ideals $I$ of $R$. 
So $\varphi^*$ is continuous, and if $\varphi$ is surjective with kernel $I$, 
then $\varphi^*$ is a ho\-meo\-morphism onto its image $\Zero(I)$.

\subsection*{Radical ideals} Let  $I$ be an ideal of~$R$.  The {\bf radical} of $I$ (more precisely, the nilradical of $I$) is the ideal
$$\sqrt{I}\ :=\ \{f\in R:\ \text{$f^n\in I$ for some $n$}\}$$
of $R$. Note that $\sqrt{I}\supseteq I$, and $\sqrt{I}=R$ iff $I=R$.
One says that $I$ is {\bf radical} if~$I=\sqrt{I}$. Every prime ideal of $R$ is radical, 
the radical of an ideal of $R$ is radical, and the intersection $\bigcap_{\lambda\in\Lambda} I_\lambda$ of a family $(I_\lambda)_{\lambda\in\Lambda}$ of radical ideals of $R$ is radical.
Moreover, if $I$ is radical and $S$ is a subset of $R$, then 
$$(I:S)\ :=\ \{f\in R:\  fS\subseteq I \}$$
is a radical ideal of $R$ containing $I$.
We call
$$\operatorname{nil}(R)\ :=\ \sqrt{ \{0\} }\ =\ \{a\in R:\ \text{$f^n=0$ for some $n$}\}$$
the {\bf nilradical} of $R$. If $\operatorname{nil}(R)$ is finitely generated, then, as is easily verified,  $\operatorname{nil}(R)^n=\{0\}$ for some $n\ge 1$. 
One says that $R$ is {\bf reduced} if $\operatorname{nil}(R)=\{0\}$, equivalently,~$\{0\}$ is a radical ideal of $R$.

\nomenclature[Az]{$\sqrt{I}$}{radical of an ideal $I$}
\nomenclature[Az]{$(I:S)$}{$\{r\in R: rS\subseteq I \}$}
\nomenclature[A]{$\operatorname{nil}(R)$}{nilradical of $R$}

\index{ideal!radical}
\index{ring!reduced}

\medskip
\noindent
For $X\subseteq \Spec(R)$, the radical ideal of all $f\in R$ that vanish at each point of $X$ is
$$\Ideal(X)\ :=\ \bigcap_{\mathfrak p\in X}\mathfrak p,$$
which by convention is $R$ for $X=\emptyset$.
Note that $\Zero\!\big(\!\Ideal(X)\big)$ is the closure of $X$ in the space $\Spec(R)$.
Here is an abstract version of Hilbert's Nullstellensatz:

\begin{prop}\label{prop:radical is intersection of primes}
For each ideal $I$ of $R$ we have $\Ideal\!\big(\!\Zero(I)\big)=\sqrt{I}$.
\end{prop}

\noindent
Towards the proof we first define a {\bf multiplicative subset of $R$\/} 
to be a set $S\subseteq R$ such that $1\in S$, and $fg\in S$ for all 
$f,g\in S$. 

\index{multiplicative subset}
\index{subset!multiplicative}

\begin{lemma}[Krull]\label{lem:Krull}
Let $S$ be a multiplicative subset of $R$. Let $I$ be an ideal of~$R$ disjoint from $S$ and maximal with these properties. Then $I$ is prime.
\end{lemma}
\begin{proof}
Let $f,g\in R$ be such that $fg\in I$. Assume towards a contradiction that 
$f,g\notin I$. Then $Rf+I$ and $Rg+I$ are ideals
of $R$ properly containing $I$, so we have $s,t\in S$ with $s\in Rf+I$, $t\in Rg+I$. Then we have a contradiction: 
\begin{equation*}
st \in S\cap (Rf+I)(Rg+I)\ \subseteq\ S\cap I\ =\ \emptyset.  \qedhere
\end{equation*} 
\end{proof}

\noindent
Using also Zorn, this lemma applied to $S=\{1\}$ gives the well-known 
fact that every ideal $I\ne R$ of $R$ is contained in a prime ideal of $R$.

\begin{proof}[Proof of Proposition~\ref{prop:radical is intersection of primes}]
As already noted, $\Ideal\!\big(\!\Zero(I)\big)$ is radical, which in view of
$\Ideal\!\big(\!\Zero(I)\big)\supseteq I$ gives $\Ideal\!\big(\!\Zero(I)\big)\supseteq \sqrt{I}$.
Conversely, suppose $f\in R$, $f\notin \sqrt{I}$. Then
$f\notin\Ideal\!\big(\!\Zero(I)\big)$, since Zorn and Lemma~\ref{lem:Krull}
yield a prime ideal $\mathfrak p\supseteq I$ of $R$ disjoint from the multiplicative 
subset $\{f^n:n=0,1,2,\dots\}$ of~$R$, so  $f\notin \mathfrak p\in \Zero(I)$.
\end{proof}

\begin{cor}\label{cor:radical is intersection of primes, 1}
The map $I \mapsto \Zero(I)$ is an inclusion-reversing bijection 
from the set of radical ideals of $R$ onto the
set of closed subsets of $\Spec(R)$, with inverse given by  
$X \mapsto \Ideal(X)$.
\end{cor}

\begin{cor} \label{cor:radical is intersection of primes, 2}
$\operatorname{nil}(R)=\bigcap_{\mathfrak p\in\Spec(R)} \mathfrak p$.
\end{cor}

\noindent
We now continue with some purely topological considerations that will help 
in coming to grips with non-traditional spaces like $\Spec(R)$.

%\begin{cor}
%For each $f\in R$, the subspace $X=\operatorname{D}(f)$ of $\Spec(R)$ is  compact: for every family $\{U_i\}_{i\in I}$ of open subsets of $X$  with $X=\bigcup_{i\in I} U_i$ there is a finite $I_0\subseteq I$ such that $X=\bigcup_{i\in I_0} U_i$. In particular, $\Spec(R)=\operatorname{D}(1)$ is compact.
%\end{cor}
%\begin{proof}
%Let $f\in R$. It suffices to show that if  $G\subseteq R$  with $\operatorname{D}(f) \subseteq \bigcup_{g\in G} \operatorname{D}(g)$, then there is a finite $G_0\subseteq G$ with $\operatorname{D}(f) \subseteq \bigcup_{g\in G_0} \operatorname{D}(g)$. To see this note that   Corollary~\ref{cor:radical is intersection of primes, 1} yields $f\in \sqrt{\sum_{g\in G} Rg}$. Hence $f\in \sqrt{\sum_{g\in G_0} Rg}$ for some finite $G_0\subseteq G$.
%\end{proof}

\subsection*{Irreducibility} {\em Let $X$ and $Y$ be topological spaces}.
Call $X$ {\bf irreducible} \index{irreducible!topological space}\index{topological space!irreducible} if $X\neq\emptyset$ and~$X$ is not the union of 
two proper closed subsets; note that then $X$ is connected and each 
nonempty open subset of $X$ is irreducible and dense in $X$.
A subset of $X$ is called irreducible if it is irreducible as a subspace of~$X$.
It is easy to see that if~$X\neq\emptyset$, then 
\begin{align*}
\text{$X$ is irreducible}&\quad\Longleftrightarrow\quad  
\text{\parbox{20em}{any two nonempty open subsets of $X$ have a nonempty intersection}}\\
&\quad\Longleftrightarrow\quad 
\text{every nonempty open subset of $X$ is dense in $X$.}
\end{align*}
Thus if $X$ is a subspace of $Y$, then 
$X$ is irreducible iff its closure $\operatorname{cl}(X)$ in $Y$ is irreducible.
In particular, for any $x\in X$ the subspace $\operatorname{cl}\!\big(\{x\}\big)$ of $X$ is irreducible.
If $f\colon X\to Y$ is continuous
and $X$ is irreducible, then $f(X)$ is irreducible.

\begin{definition}
An {\bf irreducible component} of $X$ is a maximal irreducible subset of $X$. \index{component!irreducible}\index{topological space!irreducible component}
\end{definition}

\noindent
Thus the irreducible components of $X$ are closed. One-point spaces are
irreducible, so every topological space is the union of its
irreducible components, in view of:

\begin{prop}\label{prop:irred comp}
Every irreducible subset of $X$ is contained in an irreducible
component of $X$.
\end{prop}
\begin{proof}
By Zorn, it suffices to verify that the union $Y$ of a family $\{Y_\lambda\}_{\lambda\in \Lambda}$ of irreducible subsets of
$X$, linearly ordered by inclusion, is irreducible. For this, let $U_i\subseteq X$  be open with $U_i\cap Y\neq\emptyset$,
$i=1,2$. For $i=1,2$, take $\lambda_i\in \Lambda$ with $U_i\cap Y_{\lambda_i}\neq\emptyset$.
Assuming $Y_{\lambda_1}\subseteq Y_{\lambda_2}$, we get $U_1\cap U_2\cap Y_{\lambda_2}\neq\emptyset$,
since $Y_{\lambda_2}$ is irreducible, and thus $U_1\cap U_2\cap Y\neq\emptyset$.
\end{proof}

\begin{cor}\label{cor:irred comp}
Suppose $X$ is a finite union of irreducible subsets. Then
$X$ has only finitely many irreducible components
and no irreducible component of $X$ is contained in the union of the others.
\end{cor}
\begin{proof}
By Proposition~\ref{prop:irred comp} we have $X=X_1\cup \cdots \cup X_m$
where $X_1,\dots,X_m$ are irreducible components of~$X$. It remains to note 
that if $Y$ is any irreducible subset of~$X$, then $Y\subseteq X_i$ for some 
$i$, using $Y=\bigcup_{i} (Y\cap X_i)$. 
\end{proof}

\noindent
We can now identify the irreducible closed subsets of $\Spec(R)$:

\begin{lemma}\label{lem:irred comp}
Let $X$ be a closed subset of $\Spec(R)$. Then $X$ is irreducible if and only 
if $\Ideal(X)$ is a prime ideal of $R$.
\end{lemma}
\begin{proof}
Suppose $X$ is irreducible, and $f,g\in R$, $fg\in \Ideal(X)$. Then 
for all $\mathfrak p\in X$ we have $f\in\mathfrak p$ or $g\in\mathfrak p$, 
so $X=\big(X\cap \Zero(f)\big)\cup \big(X\cap \Zero(g)\big)$, hence $X\subseteq \Zero(f)$ or $X\subseteq \Zero(g)$, and thus $f\in \Ideal(X)$ or $g\in \Ideal(X)$. Conversely,
assume $\Ideal(X)$ is prime, and $X=X_1\cup X_2$ with closed $X_1$, $X_2$.
Then $\Ideal(X)=\Ideal(X_1\cup X_2)=\Ideal(X_1)\cap \Ideal(X_2)$, hence
$\Ideal(X)=\Ideal(X_1)$ or $\Ideal(X)=\Ideal(X_2)$, so $X=X_1$ or $X=X_2$. 
\end{proof}

\noindent 
Thus $X\mapsto \Ideal(X)$ is an inclusion-reversing 
bijection from the set of
closed irreducible subsets of $\Spec(R)$ onto the set of prime ideals of $R$. 

\index{ideal!minimal prime divisor}
\index{minimal!prime divisor}

\medskip
\noindent
Let $I$ be an ideal of $R$. A {\bf minimal prime divisor of $I$} is 
a point $\mathfrak p\in \Zero(I)$ such that no $\mathfrak{q}\in \Zero(I)$
is strictly contained (as a set) in $\mathfrak{p}$. 
For $I=\{0\}$, the minimal prime
divisors of $I$ are exactly the minimal prime ideals of $R$ (with respect to
inclusion).
The above bijection $X\mapsto \Ideal(X)$ maps the set of irreducible components of the closed subset $\Zero(I)$ of $\Spec(R)$ onto the set of 
minimal prime divisors of $I$. 
Thus by Proposition~\ref{prop:irred comp}, every $\mathfrak{p}\in \Zero(I)$ 
contains a minimal prime divisor of $I$. 

\subsection*{Prime avoidance}
For use in Sections~\ref{sec:PIT} and~\ref{sec:regular local rings} we show: 

\begin{lemma}\label{lem:prime avoidance}
Let $I,J_1,\dots,J_n$ be ideals of $R$ with
$n\geq 2$ and $I\subseteq J_1\cup\cdots\cup J_n$. Assume 
$J_1,\dots,J_{n-2}$ are prime. Then $I\subseteq J_j$ for some 
$j\in\{1,\dots,n\}$.
\end{lemma}
\begin{proof}
We argue by induction on $n=2,3,\dots$.
Let $n=2$, and assume that neither $I\subseteq J_1$ nor $I\subseteq J_2$.
Take $a_1\in I\setminus J_2$, $a_2\in I\setminus J_1$. Then $a_1\in J_1$, $a_2\in J_2$, hence
$a_1+a_2\in I\setminus (J_1\cup J_2)$. Next, let $n\ge 3$ and let
$j$, $k$ range over $\{1,\dots,n\}$.
It suffices to find $k$ with $I\subseteq \bigcup_{j\neq k} J_j$, since then we are done by inductive hypothesis.
Towards a contradiction, assume that $I\not\subseteq \bigcup_{j\neq k} J_j$ for all $k$. 
For each $k$, take 
$a_k\in I\setminus \left(\bigcup_{j\neq k} J_j\right)$. Then $a_k\in J_k$ for each $k$, hence 
$a:=\prod_{k\neq 1} a_k\in I\cap \bigcap_{k\neq 1} J_k$; but $a\notin J_1$, since $J_1$ is prime.
Thus $a_1+a\in I$ and $a_1+a\notin J_k$ for all $k$, a contradiction.
\end{proof}

\begin{cor}\label{cor:prime avoidance}
Let $I, J$ be ideals of $R$ such that $I\not\subseteq J$, and let
 $\frak{p}_1,\dots,\frak{p}_n\in \Spec(R)$ be such that
$I\setminus J\subseteq \frak{p}_1\cup\cdots\cup \frak{p}_n$. Then 
$I\subseteq \frak{p}_j$ for some $j\in\{1,\dots,n\}$.
\end{cor}
\begin{proof}
Use $I\subseteq \frak{p}_1\cup\cdots\cup \frak{p}_n\cup J$ and 
Lemma~\ref{lem:prime avoidance}.
\end{proof}

\subsection*{Chain conditions}
In this subsection the set $S$ is \textit{partially}\/ ordered by~$\leq$, with
$$s < s'\ :\Longleftrightarrow\ \text{$s\le s'$ and $s\ne s'$,}$$ for $s, s'\in S$. 
A {\bf maximal} element of a set $X\subseteq S$ is an $x\in X$ 
such that there is no $x'\in X$ with $x< x'$.
The following are equivalent:
\begin{enumerate}
\item Every nonempty subset of $S$ has a maximal element;
\item there is no sequence $s_1< s_2 < \cdots < s_n < s_{n+1} <\cdots$ in $S$.
\end{enumerate}
We say that $S$ satisfies the {\bf ascending chain condition} 
(or has {\bf acc}) if 
$S$ fulfills one of these equivalent conditions. We say that $S$ satisfies the {\bf descending chain condition} (or has {\bf dcc}) if $S$ with the reversed partial
ordering $\geq$ has acc. Thus if the ordering of~$S$ is total,
then $S$ has dcc~iff $S$ is well-ordered.

\index{maximal!elements of a partially ordered set}
\index{ascending chain condition}
\index{descending chain condition}
\index{chain conditions}
\index{acc}
\index{dcc}

\subsection*{Noetherian rings}
Call $R$ {\bf noetherian} if the set of its ideals, partially ordered 
by inclusion, has acc. If $R$ is noetherian, then so is its image under any 
ring morphism.

\index{ring!noetherian}
\index{noetherian!ring}

\begin{lemma} \label{lem:noeth} $R$ is noetherian $\Longleftrightarrow$ every ideal of $R$ is finitely generated.
\end{lemma}
\begin{proof}
If the ideal $I$ of $R$ is not finitely generated, then we obtain 
inductively a sequence $(r_n)$ in $I$ with $r_{n+1}\notin (r_0,\dots,r_n)$ for all $n$, giving a strictly increasing infinite sequence
$(r_0) \subset (r_0, r_1) \subset (r_0, r_1, r_2) \subset \cdots$
of ideals of $R$. 
Conversely, if every ideal of $R$ is finitely generated, 
it follows easily that there cannot exist a strictly increasing infinite 
sequence of ideals of $R$.  
\end{proof}

\noindent
Fields and, more generally, principal ideal domains (like the ring of integers)
are noetherian. In Chapter~\ref{ch:valuedfields} we study valuation rings, which are only noetherian in a very special (but important) case, as we
now explain. We define here
a {\bf valuation ring\/} to be an integral domain $R$ such that
for all $a,b\in R$: $\ a\in bR$ or $b\in aR$. 

Suppose that $R$ is a valuation ring.  For any $a_1,\dots, a_n\in R$, $n\ge 1$, we have $(a_1,\dots, a_n)R=a_jR$ for some $j$; in
particular, every finitely generated ideal of $R$ is 
principal. The set $R\setminus R^{\times}$ of nonunits of
$R$ is clearly the largest proper ideal of $R$, and is thus the unique maximal ideal $\m_R$ of $R$.

Fields are valuation rings, but are viewed as trivial within the class of valuation rings. Next in complication are {\em discrete valuation rings}:
a {\bf discrete valuation ring\/} (or {\bf DVR}) is an integral domain $R$ with an element $t$ such that $t\ne 0$, $t\notin R^\times$, and $R^{\ne}=R^\times t^{\N}$. 
Note that every DVR is indeed a valuation ring. 
The power series ring~$\k[[t]]$ in an indeterminate $t$ over a field $\k$ is clearly a DVR. 

Suppose $R$ is a DVR. Then $R$ is a PID (principal ideal domain):  with $t$ as in the definition of 
\textit{discrete valuation ring}\/,
any ideal $I\ne \{0\}$ of $R$ is clearly generated by $t^m$ where~$m$ is minimal such that $t^m\in I$. In particular, a DVR
is noetherian:
 
\index{valuation ring!discrete}
\index{DVR}
\index{discrete valuation ring}
\index{valuation ring}
\index{ring!valuation}

\begin{cor}
Let $R$ be a valuation ring that is not a field. Then
$$\text{$R$ is noetherian}\ \Longleftrightarrow\ 
\text{$R$ is a \rm{PID}}\ \Longleftrightarrow\ 
\text{$R$ is a \rm{DVR}.}$$
\end{cor}
\begin{proof}
Suppose $R$ is a PID. Take $t\in R$ such that $\m_R=Rt$. Then $t$ is prime in~$R$.
If~$t'$ is also prime in $R$, then, $R$ being a valuation ring, $t=at'$ or 
$t'=at$ with $a\in R$, and in either case, $a\in R^{\times}$. 
Thus,  $R$ being factorial, every $r\in R^{\neq}$ has the form $r=ut^n$ with $u\in R^\times$. So $R$ is a DVR. The rest is clear.
\end{proof}

\noindent
Many rings of natural origin are noetherian, by a famous result:

\begin{prop}[Hilbert's Basis Theorem]\label{prop:HBT}
If $R$ is noetherian, then so is the ring~$R[X]$ of polynomials in the indeterminate $X$ over $R$. 
\end{prop}
\begin{proof}
Suppose $I$ is an ideal of $R[X]$ that is not finitely generated. Set~${f_0:=0}$, and with $f_0,f_1,\dots,f_n\in I$, take $f_{n+1}\in I\setminus (f_1,\dots,f_n)$ of minimal degree. For ${n\geq 1}$, let
$r_n\in R^{\neq}$ and $d_n\in\N$ be such that
$f_n=r_nX^{d_n}+\text{lower degree terms}$. Then $d_1\leq d_2\leq\cdots$.
Also, $(r_1,\dots, r_n)\ne (r_1,\dots, r_n, r_{n+1})$: otherwise
$r_{n+1} = \sum_{i=1}^n a_ir_i$ with all $a_i\in R$, so $f_{n+1}-\sum_{i=1}^n a_i X^{d_{n+1}-d_i} f_i$ has smaller degree than $f_{n+1}$, contradicting the choice of $f_{n+1}$.
Hence $R$ is not noetherian.
\end{proof}

\noindent
Thus if $R$ is noetherian, then so is every finitely generated commutative $R$-algebra.

\subsection*{Noetherian spaces}
{\em Let $X$, $Y$ be topological spaces}. 
Call $X$ {\bf noetherian} if its collection of closed sets satisfies the descending chain
condition: there is no strictly descending infinite sequence
$X_0\supset X_1\supset\cdots$ of closed subsets of $X$; equivalently, each nonempty collection
of closed subsets of $X$ has a minimal element with respect to inclusion.
If $R$ is noetherian, then the space $\Spec(R)$ is noetherian.

\index{topological space!noetherian}
\index{noetherian!topological space}

\begin{remarks}
A noetherian space is quasicompact (every covering by open subsets has a finite subcovering), but in the absence of being
hausdorff this is less useful than some other facts:
\begin{enumerate}
\item each subspace of a noetherian space is noetherian;
\item if $X$ is noetherian and $f\colon X\to Y$ is continuous, then $f(X)\subseteq Y$ is noetherian;
\item if $X$ is covered by finitely many noetherian subspaces, then $X$ is
noetherian.
\end{enumerate} 
\end{remarks}

\noindent
Suppose $X$ is noetherian. Then $X$ is a finite union of 
irreducible closed
subsets: if not, take a minimal closed subset $Y$ of $X$ that is not a finite
union of irreducible closed subsets. Then $Y\neq\emptyset$ and $Y$ is not irreducible, so
$Y = Y_1 \cup Y_2$ with $Y_1$, $Y_2$ proper closed subsets of $Y$. Each $Y_i$ is a finite
union of irreducible closed subsets, and so is~$Y$, a contradiction.
Thus by Corollary~\ref{cor:irred comp}:

\begin{cor}\label{cor:irred comp 1} 
If $\Spec(R)$ is noetherian \textup{(}in particular, if $R$ is noetherian\textup{)} 
and~$I$ is an ideal of $R$, then $I$ has only finitely 
many minimal prime divisors.
\end{cor}

\noindent
So if $\Spec(R)$ is noetherian, then $R$ has only finitely many minimal
prime ideals. By a {\bf zero divisor\/} of $R$ we mean an $a\in R$ such that $ab=0$ for
some $b\in R^{\ne}$.

\index{ring!zero divisor}
\index{zero divisor}

\begin{cor}\label{cor:irred comp 2}
Suppose $\Spec(R)$ is noetherian and $\mathfrak p_1,\dots,\mathfrak p_n$ are the minimal prime ideals of $R$.
Then $\operatorname{nil}(R) = \mathfrak p_1 \cap\cdots\cap\mathfrak p_n$, 
the elements of $\mathfrak p_1 \cup\cdots\cup\mathfrak p_n$ are zero divisors of $R$, and if $R$ is reduced, $\mathfrak p_1 \cup\cdots\cup\mathfrak p_n$ is the set of zero divisors of $R$.
\end{cor}
\begin{proof}
The first statement holds by Corollary~\ref{cor:radical is intersection of primes, 2}.
The second statement holds clearly for $n=0,1$, so assume $n\geq 2$ and
$\mathfrak p_1,\dots,\mathfrak p_n$ are distinct. Let $i$, $j$ range over $\{1,\dots,n\}$, and let 
$r\in\mathfrak p_j$. For $i\ne j$, take $s_i\in \mathfrak{p}_i\setminus \mathfrak{p}_j$. Then $s:= \prod_{i\ne j} s_i\in \bigcap_{i\neq j}\mathfrak p_i$, 
$s\notin\mathfrak p_j$. Then $rs\in\bigcap_i\mathfrak p_i$, so
we have $m\geq 1$ with $(rs)^m=0$. Since $s\notin\mathfrak p_j$, we have $s^m\neq 0$, so we get
$i\in \N$ with $i< m$ and $r^i s^m\neq 0$, $r^{i+1}s^m=0$. Thus~$r$ is a zero divisor of $R$. 
For the third statement, assume that $R$ is reduced, that is, $\mathfrak p_1 \cap\cdots\cap\mathfrak p_n=\{0\}$.
Let $r\in R$ be a zero divisor. Take $s\in R^{\neq}$ with $rs=0$, and
then take $j\in\{1,\dots,n\}$ with $s\notin\mathfrak p_j$. From $rs=0\in\mathfrak p_j$ we obtain $r\in\mathfrak p_j$.
\end{proof}

\subsection*{Krull dimension} \index{topological space!Krull dimension}
\index{dimension!Krull}
\index{Krull dimension!topological space}
This is a notion of dimension suitable for noetherian spaces. Let $X$,~$Y$ be topological spaces. In this subsection we take suprema and infima in~$\N \cup \{-\infty, +\infty\}$. We define the (Krull) {\bf dimension} of $X$ to be the supremum of the set of $n$ for which there is a 
strictly increasing sequence
$X_0 \subset X_1 \subset \cdots \subset X_n$
of irreducible closed subsets of $X$. In particular, the  dimension of $X$ is $-\infty$ iff~$X = \emptyset$. We denote the  dimension of $X$ by~$\dim(X)$. (If $X$ is a nonempty hausdorff space, then $\dim(X) = 0$, so  Krull 
dimension is of no interest for hausdorff spaces.) 
If~$X$ is a subspace of $Y$, then $\dim(X) \leq \dim(Y)$. 
(Use that the closure in $Y$ of an irreducible subset is irreducible.)
Moreover, 
$$\dim(X)\ =\ \sup\big\{\!\dim(Y) : \text{$Y$ is an irreducible component of $X$}\big\}.$$

\subsection*{Some special prime ideals} 
Let $K$ be a field, and
let a family 
$X=(X_\lambda)_{\lambda\in \Lambda}$ of distinct indeterminates be given. If 
$\Lambda$ is finite, then the ring $K[X]$ is noetherian by Proposition~\ref{prop:HBT}. If
$\Lambda$ is infinite, then $K[X]$ is not noetherian. However, such rings appear later as rings of differential polynomials, and in Section~\ref{sec:diff rings} we use:

\begin{lemma}\label{spprid} Let $I$ be an ideal of $K[X]$ generated by homogeneous polynomials of degree $1$. Then $I$ is a prime ideal of $K[X]$.
\end{lemma}
\begin{proof} Take a $K$-independent set $H\subseteq K[X]$ of homogeneous polynomials of degree~$1$ that generates the ideal $I$. Then $H$ is part of a basis $B$ of the $K$-linear subspace of~$K[X]$ generated by the $X_\lambda$. Take a bijection $\lambda\mapsto b_{\lambda}\colon \Lambda \to B$. Then $P(X) \mapsto
P\big((b_{\lambda})\big)$ is an automorphism of the $K$-algebra $K[X]$. Replacing $I$ by the inverse image of $I$ under this automorphism, we arrange that $H=\{X_\lambda:\ \lambda\in \Lambda_0\}$, with $\Lambda_0\subseteq \Lambda$.
Hence $K[X]/I$ is isomorphic as a $K$-algebra to $K[X_\lambda: \lambda\in \Lambda\setminus \Lambda_0]$.  
\end{proof}

\subsection*{Notes and comments} 
 The notion of {\em noetherian ring}\/ and the basic facts about it are due to E.~Noe\-ther~\cite{Noether}. 
(The short proof of Hilbert's Basis Theorem~\cite{Hilbert} 
given above is from~\cite{Sarges}.)
Krull's Lemma~\ref{lem:Krull} is from \cite[p.~732]{Krull29a}.
%Artinian rings were introduced in \cite{Artin27}.
Antecedents for the Zariski topology on $\Spec(R)$ include Stone~\cite{Stone36,Stone37},
Jacobson~\cite{Jacobson45}, and 
Zariski~\cite{Zariski52}. Irreducible and noetherian spaces come
from Serre~\cite{Serre55}.

\section{Rings and Modules of Finite Length}\label{sec:finite length}

\noindent
The algebra in this section is not just commutative: {\em $R$ is a ring, possibly not commutative, and $M$,~$N$ range over $R$-modules}.  

\subsection*{Composition series} An {\bf $M$-series}\index{series!$M$-series} (of length $m$) 
is a strictly increasing sequence
\begin{equation}\label{eq:comp series}
\{0\}=M_0\subset M_1\subset\cdots\subset M_m=M
\end{equation}
of submodules of $M$. A {\bf refinement\/} \index{refinement!$M$-series} of an $M$-series 
\eqref{eq:comp series} is an $M$-series
$$\{0\}=N_0\subset N_1\subset\cdots\subset N_{n}=M$$ such that
$\{M_0,\dots,M_m\}\subseteq \{N_0,\dots, N_n\}$ (so $m\le n$).
Two $M$-series
$$\{0\}=M_0\subset M_1\subset\cdots\subset M_m=M, \quad \{0\}=N_0\subset N_1\subset\cdots\subset N_{n}=M$$
of $M$ are said to be {\bf equivalent} \index{equivalence!$M$-series} if $m=n$ and there is a permutation 
$i\mapsto i'$ of $\{1,\dots,m\}$ such that $M_{i}/M_{i-1}\cong N_{i'}/N_{i'-1}$, as $R$-modules, 
for $i=1,\dots,m$. The next result is known as the
Schreier Refinement Theorem.

\begin{theorem}\label{thm:Schreier} Any two $M$-series have equivalent refinements.
\end{theorem}

\noindent
See \cite[Chapter~I, \S{}3]{Lang} for a proof of an analogue of this theorem for groups which adapts in a straightforward
way to modules.

\medskip
\noindent
We call $M$ {\bf simple} (or {\bf irreducible}) \index{module!irreducible}\index{irreducible!module}\index{module!simple}\index{simple!module} if $M\neq\{0\}$ and $M$ has no submodules other than~$\{0\}$ and $M$. If $M$ is simple, then $M=Re$
for each $e\in M^{\neq}$. 
Hence $M$ is simple iff $M\cong R/\m$ for some maximal left ideal $\m$ of~$R$.
It is clear that if $M$ and~$N$ are simple $R$-modules, then every $R$-linear map $\varphi\colon M\to N$ with $\varphi\neq 0$ is an isomorphism of $R$-modules (Schur's Lemma).

\medskip
\noindent
A {\bf composition series} of $M$ is an $M$-series \eqref{eq:comp series} such that $M_{i+1}/M_{i}$ is simple for $i=0,\dots,m-1$.
Every $M$-series equivalent to a composition series of~$M$ is itself a composition series of $M$, and a composition series of~$M$ cannot be 
refined to a strictly longer $M$-series.
Thus Theorem~\ref{thm:Schreier} yields the Jordan-H\"older Theorem:

\index{series!composition}
\index{composition!series}

\begin{cor}\label{cor:JH}
Any two composition series of $M$ are equivalent. 
\end{cor}

\subsection*{Euler-Poincar\'e maps}
Let $A$ be an abelian group. Suppose that for certain $R$-modules 
$M$ there is defined a quantity $\chi(M)\in A$, such that: \begin{enumerate}
\item $\chi(\{0\})$ is defined and equal to $0\in A$,
\item if $0 \to K\to M\to N\to 0$ is an exact sequence of $R$-modules, then $\chi(M)$ is defined if and only if both $\chi(K)$ and $\chi(N)$ are defined, and in this case, $\chi(M) = \chi(K) + \chi(N)$.
\end{enumerate}
(Such an assignment $M\mapsto \chi(M)$ is called an {\bf Euler-Poincar\'e map} 
on $R$-modules.)
Clearly if $M\cong N$, then $\chi(M)$ is defined iff $\chi(N)$ is defined, and in this case $\chi(M)=\chi(N)$. For example, if $R=\Z$, then setting
$\chi(M):=\log |M|$ for a finite abelian group~$M$ defines
an Euler-Poincar\'e map on abelian groups with values in the additive group $A=\R$.

\index{Euler-Poincar\'e map}
\index{map!Euler-Poincar\'e}

\begin{lemma}\label{alternatingexact, length}
Let an exact sequence
$$0\longrightarrow M_0 \xrightarrow{\ \varphi_0\ } M_1 \xrightarrow{\ \varphi_1\ } \cdots \longrightarrow M_n \longrightarrow 0$$
of $R$-modules be given such that $\chi(M_i)$ is defined for $i=0,\dots,n$.
Then 
$$\sum_{i=0}^n (-1)^i \chi(M_i)=0.$$
\end{lemma}
\begin{proof}
By induction on $n$. The cases $n=0,1,2$ are obvious. Suppose $n\geq 3$,
and put $N:=\operatorname{im} \varphi_1 = \ker\varphi_2$. Then we have exact sequences
$$0\longrightarrow M_0 \xrightarrow{\ \varphi_0\ } M_1\xrightarrow{\ \varphi_1\ } N\longrightarrow 0, \qquad
0\longrightarrow N \xrightarrow{\ \subseteq\ } M_2 \longrightarrow \cdots \longrightarrow M_n \longrightarrow 0.$$
By the first sequence $\chi(N)$ is defined and 
$\chi(M_0)-\chi(M_1)+\chi(N)=0$. Now use the inductive hypothesis on the second sequence.
\end{proof}

\subsection*{Modules of finite length}
If $M$ has a composition series, then all composition series of $M$ have a common length, called the {\bf length} of~$M$, denoted by~$\ell(M)$. If~$M$ does not have a composition series, we put $\ell(M):=\infty$.
Note: ${\ell(M)=0}$ iff~${M=\{0\}}$.
If
$N$ is a submodule of $M$, then $M$ has finite length if and only if~$N$ and~$M/N$ have finite length, in which case $\ell(M)=\ell(N)+\ell(M/N)$. Hence ${M\mapsto \ell(M)}$, defined
for $M$ of finite length, is a $\Z$-valued Euler-Poincar\'e map on $R$-modules.

\index{module!length}
\index{length}

\begin{example}
Suppose $R=K$ is a field. Then $\ell(M)=\dim_K M$. 
In this way Lem\-ma~\ref{alternatingexact, length} 
contains Lemma~\ref{alternatingexact} as a special case. 
\end{example}

\index{module!noetherian}
\index{module!artinian}
\index{noetherian!module}
\index{artinian module}

\noindent
Call $M$ {\bf noetherian} if the set of submodules of $M$, partially ordered by inclusion, has acc. Call
$M$ {\bf artinian} if this partially ordered set has dcc. 
As in the proof of Lemma~\ref{lem:noeth}, $M$ is noetherian iff all its submodules are finitely generated. 

\begin{lemma}\label{lem:cc and finite length}
$M$ is noetherian and artinian if and only if $M$ has finite length.
Thus if $M$ has finite length, then $M$ is finitely generated.
\end{lemma}
\begin{proof} Previous remarks make it clear that
if $\ell(M)<\infty$, then $M$ is noetherian and artinian.
Conversely, suppose $M$ is artinian. 
If~$M\neq \{0\}$,
let $M_1$ be a minimal nonzero submodule of $M$; if $M\neq M_1$, let $M_2$ be a submodule of $M$ which is minimal
among the submodules of $M$ properly containing~$M_1$, and so on. This yields
a strictly increasing sequence
$\{0\}=M_0\subset M_1\subset M_2\subset\cdots$ of submodules of $M$. If $M$ is also noetherian, this construction stops with $M_m=M$ for some $m$, and then
we have a composition series of~$M$.
\end{proof}

\subsection*{Rings of finite length} {\em In this subsection $R$ is commutative}. In the phrases ``$R$ has  finite length'' and ``$R$ is artinian'' we view $R$ as an $R$-module in the 
usual way.

\begin{lemma}\label{lem:artinian=>dim 0}
Suppose $R$ is artinian. Then every prime ideal of $R$ is maximal.
\end{lemma}
\begin{proof}
Let $\mathfrak p$ be a prime ideal of $R$; then $R/\mathfrak p$ is an
artinian integral domain, and replacing $R$ by $R/\mathfrak p$ we may
assume that $R$ is an integral domain, and need to show that~$R$ is a field.  Let $r\in R^{\neq}$. 
Considering the chain $Rr\supseteq Rr^2\supseteq\cdots$, we obtain~$n\geq 1$ and $s\in R$ such that
$r^n=r^{n+1}s$; hence $1=rs$. 
\end{proof}

\begin{prop}\label{prop:dim 0 rings}
The following conditions on $R$ are equivalent:
\begin{enumerate}
\item[\textup{(i)}] $R$ is noetherian and every prime ideal of $R$ is maximal;
\item[\textup{(ii)}] every finitely generated $R$-module has finite length;
\item[\textup{(iii)}] $R$ has finite length;
\item[\textup{(iv)}] $R$ is noetherian and artinian.
\end{enumerate}
\end{prop}
\begin{proof} If $R=\{0\}$, then all four conditions are trivially satisfied, so
let $R\ne \{0\}$. 
Assume (i) holds, and let $\mathfrak m_1,\dots,\mathfrak m_n$ be the minimal prime ideals of~$R$; so each $\mathfrak m_i$ is maximal.  Put $\mathfrak n:=\mathfrak m_1\cdots\mathfrak m_n$, so $\mathfrak n\subseteq \mathfrak m_1\cap\cdots\cap\mathfrak m_n=\operatorname{nil}(R)$, giving $m\geq 1$ with $\mathfrak n^m=\{0\}$.
Let $M$ be a finitely generated $R$-module. 
Consider the chain
$$M\ \supseteq\ \mathfrak m_1 M\ \supseteq\ \mathfrak m_1\mathfrak m_2 M\ \supseteq\ \cdots\supseteq\ \mathfrak m_1\cdots\mathfrak m_n M\ =\ \mathfrak n M$$
of submodules of $M$.
For $i=0,\dots,n-1$, the quotient $\mathfrak m_1\cdots\mathfrak m_i M/\mathfrak m_1\cdots\mathfrak m_{i+1}M$ is a finitely generated vector space over the field
$R/\mathfrak m_{i+1}$, hence has finite length as an $R$-module.
So $M/\mathfrak nM$ has finite length. Likewise, $\mathfrak nM/\mathfrak n^2M$
has finite length, and so $M/\mathfrak n^2M$ has finite length. 
Proceeding this way we see that $M/\mathfrak n^mM$ has finite length. 
As $\mathfrak n^m=\{0\}$, so does $M$, showing (ii).
The implication (ii)~$\Rightarrow$~(iii) is obvious, (iii)~$\Longleftrightarrow$~(iv) is 
Lemma~\ref{lem:cc and finite length}, and  
Lemma~\ref{lem:artinian=>dim 0} gives
(iv)~$\Rightarrow$~(i). 
\end{proof}

\subsection*{Notes and comments}  
The Schreier Refinement Theorem is from \cite{Schreier}, and the Jordan-H\"older Theorem
from~\cite{Jordan, Hoelder89}. 
%Artinian rings were introduced in~\cite{Artin27}. 
In connection with Proposition~\ref{prop:dim 0 rings} we mention that every artinian commutative ring is automatically noetherian: theorem of Akizuki~\cite{Akizuki}; see 
also~\cite[Theorem~3.2]{Matsumura}.

\section{Integral Extensions and Integrally Closed Domains}\label{sec:intext}

\noindent
In this section we establish some facts for use in Chapter~\ref{ch:valuedfields}. But first a reminder on matrices. Let $R$ be a  
commutative ring, $n\geq 1$, and let $R^{n\times n}$ be the $R$-algebra of $n\times n$ matrices over $R$, with multiplicative identity $\I_n$, the $n\times n$ identity matrix. 

\index{matrix!adjoint}
\index{adjoint!matrix}

Let a matrix $T =(T_{ij})\in R^{n\times n}$ be given. It has determinant
$\det T\in R$, and transpose $T^{\operatorname{t}}\in R^{n\times n}$. For $i,j=1,\dots,n$, let $T^{ij}$ be the determinant of the
$(n-1)\times (n-1)$ matrix obtained from $T$ by removing the $i$th row and the $j$th column from $T$; by convention this determinant equals $1$ for $n=1$. 
The matrix $$T^*\ :=\ \big((-1)^{i+j} T^{ij}\big)^{\operatorname{t}} \in R^{n\times n} $$ is called the {\bf adjoint of $T$} and satisfies
$$TT^*\  =\  T^*T\ =\ (\det T)\I_n \qquad(\text{Cramer's Rule, Laplace expansion}).$$

\subsection*{Integral extensions}
{\em In this subsection $A$ is a subring of the commutative ring~$B$.}\/
An element $b\in B$ is said to be {\bf integral over $A$} if
there are $n\ge 1$ and $a_1,\dots, a_n\in A$ such that
$b^n + a_1b^{n-1} + \cdots + a_{n-1}b + a_n=0$.  \index{integral!element}\index{element!integral}

\begin{samepage}

\begin{lemma}\label{lem:integral element}
Let $b\in B$. The following are equivalent:
\begin{enumerate}
\item[\textup{(i)}] $b$ is integral over $A$;
\item[\textup{(ii)}] the submodule $A[b]$ of the $A$-module $B$ is finitely generated;
\item[\textup{(iii)}] some subring of $B$ contains $A[b]$ and is finitely generated as an $A$-module.
\end{enumerate}
\end{lemma}

\end{samepage}

\begin{proof} If $b^n+ a_{1}b^{n-1}+ \cdots + a_{n-1}b + a_n=0$ with $n\ge 1$ and $a_1,\dots, a_n\in A$, then clearly 
$A[b]=A + Ab + \cdots + Ab^{n-1}$. This gives (i)~$\Rightarrow$~(ii), and (ii)~$\Rightarrow$~(iii) is trivial. 
To show (iii)~$\Rightarrow$~(i), let $C\supseteq A[b]$ be a subring of $B$ that is finitely generated as an $A$-module. 
Take $c_1,\dots, c_n\in C$ such that $C=Ac_1+ \cdots + Ac_n$, $n\ge 1$.
Then
\begin{align*} bc_1\ &=\ a_{11}c_1+ \cdots + a_{1n}c_n\\
  \vdots \ & \hskip3em           \vdots \hskip3em  \vdots\hskip2em \vdots\  \\
            bc_n\ &=\ a_{n1}c_1+ \cdots + a_{nn}c_n
\end{align*}
for certain $a_{ij}\in A$. Then for the $n\times n$ matrix
$T=b\I_n - (a_{ij})$ we have $Tc=0$, where $c$ is 
the column vector $(c_1,\dots, c_n)^{\operatorname{t}}\in C^n$ and likewise,
$0=(0,\dots,0)^{\operatorname{t}}$ in $C^n$.
Multiplying $Tc=0$ on the left by $T^*$ we obtain
$\det(T)=0$. This gives an equality $b^n+ a_{1}b^{n-1}+ \cdots + a_{n-1}b + a_n=0$ with $a_1,\dots, a_n\in A$.       
\end{proof}

\noindent
We say that $B$ is {\bf integral over $A$}  if every element of $B$ is integral over $A$.\index{integral!ring extension}\index{extension!integral}\index{ring!integral extension}

\begin{cor}
Let $b_1,\dots,b_m\in B$. Then the following are equivalent:
\begin{enumerate}
\item[\textup{(i)}] each $b_i$ is integral over $A$;
\item[\textup{(ii)}] the submodule $A[b_1,\dots,b_m]$ of the $A$-module $B$ is finitely generated;
\item[\textup{(iii)}] $A[b_1,\dots,b_m]$ is integral over $A$.
\end{enumerate}
\end{cor}
\begin{proof} Suppose $P_i(b_i)=0$ for monic $P_i\in A[X]$
of degree $d_i\ge 1$, $i=1,\dots,m$. 
Then $A[b_1,\dots, b_m]$ is generated as an $A$-module by
the products $b_1^{j_1}\cdots b_m^{j_m}$ with $0\le j_1< d_1,\dots, 0\le j_m< d_m$. This gives (i)~$\Rightarrow$~(ii). The implication (ii)~$\Rightarrow$~(iii) follows from (iii)~$\Rightarrow$~(i) in Lemma~\ref{lem:integral element},
and  (iii)~$\Rightarrow$~(i) is obvious.
\end{proof}

\begin{cor}
Let  $c$ be an element in a commutative ring extension of $B$.
Suppose~$c$ is integral over $B$ and $B$ is integral over~$A$. Then $c$ is integral over $A$.
\end{cor}
\begin{proof}
Take $b_1,\dots,b_n\in B$ such that $c$ is integral over the ring $A[b_1,\dots,b_n]$; the latter
is a finitely generated $A$-module, and
$A[b_1,\dots,b_n,c]$ is finitely generated as an $A[b_1,\dots,b_n]$-module, hence also as an $A$-module.
So $c$ is integral over $A$.
\end{proof}

\index{closure!integral}
\index{integral!closure}
\index{integrally closed}
\index{closed!integrally}

\noindent
We say that $A$ is {\bf integrally closed in $B$\/} if every $b\in B$ that is integral over $A$ already lies in $A$. 
The set of elements of $B$ that are integral over $A$ is called the {\bf integral closure of $A$ in $B$}.\/
By the previous two corollaries,
the integral closure of $A$ in $B$ is a subring of~$B$ that contains~$A$, is integral over $A$, and 
integrally closed in $B$.

\subsection*{Prime ideals under integral extensions}
{\em In this subsection $A$ is a subring of the commutative 
ring $B$, and $B$ is integral over $A$}.
  
\begin{lemma}\label{lem:intersect nzd}
Assume $1\ne 0$ in $B$, and let $J$ be an ideal of $B$ containing an
element that is not a zero divisor of $B$. Then $J\cap A\neq\{0\}$.
\end{lemma}
\begin{proof}
Suppose $b\in J$ is not a zero divisor of $B$. Take $n\ge 1$ minimal such that there exist $a_1,\dots,a_n\in A$ with $b^n+a_{1}b^{n-1}+\cdots+a_n=0$. For such $a_1,\dots,a_n$ we have 
$0\ne a_n\in bB\cap A\subseteq J\cap A$.
\end{proof}

\begin{cor}\label{cor:primes incomp}
Let $\mathfrak q,\mathfrak q'\in\Spec(B)$, and suppose 
$\mathfrak q\subseteq\mathfrak q'$ and
$\mathfrak q\cap A=\mathfrak q'\cap A$. Then~$\mathfrak q=\mathfrak q'$.
\end{cor}
\begin{proof}
Let $\mathfrak p:=\mathfrak q\cap A=\mathfrak q'\cap A\in\Spec(A)$, and let $\overline{A}$ be the image of $A/\mathfrak p$ in $\overline{B}:=B/\mathfrak{q}$ under the natural
embedding $A/\mathfrak p\to B/\mathfrak q$. Then the domain $\overline{B}$ is integral over its subring~$\overline{A}$, and $\mathfrak q'/\mathfrak q$ is a prime ideal of $\overline{B}$ that intersects $\overline{A}$ trivially. By Lemma~\ref{lem:intersect nzd} this yields $\mathfrak q=\mathfrak q'$.
\end{proof}

\begin{lemma}\label{lem:CS}
Let $I$ be an ideal of $A$. Then $IB\cap A\subseteq\sqrt{I}$.
\end{lemma}
\begin{proof}
Let $b\in IB\cap A$. 
Take a finitely generated subalgebra $C$ of the $A$-algebra~$B$
with $b\in IC$; then $C$ is also finitely generated as
$A$-module. Proceeding as in (iii)~$\Rightarrow$~(i) in the proof
of Lemma~\ref{lem:integral element}, and using the notations there, we get all $a_{ij}\in I$, and so obtain
$b^n + a_1 b^{n-1} + \cdots + a_n = 0$ where $n\geq 1$ and $a_1,\dots,a_{n}\in I$,
so $b^n\in I$.
\end{proof}

\begin{cor}\label{cor:CS}
For each $\mathfrak p\in\Spec(A)$ there is a $\mathfrak q\in\Spec(B)$ with $\mathfrak p=\mathfrak q\cap A$.
\end{cor}
\begin{proof}
Let $\mathfrak p\in\Spec(A)$, and set $S:=A\setminus\mathfrak p$.
Then $\mathfrak pB\cap S=\emptyset$ by Lemma~\ref{lem:CS},
so Lemma~\ref{lem:Krull} gives $\mathfrak q\in\Spec(B)$
with $\mathfrak pB\subseteq \mathfrak q$ and $\mathfrak q\cap S=\emptyset$, 
hence $\mathfrak q\cap A=\mathfrak p$.
\end{proof}

%\begin{cor}
%Let $\iota\colon A\to B$ be the natural inclusion. Then the %continuous map
%$\mathfrak q\mapsto \iota^*(\mathfrak q)=
%\mathfrak q\cap A\colon \Spec(B)\to\Spec(A)$ is surjective
%and closed \textup{(}i.e., maps closed subsets of 
%$\Spec(B)$ to closed subsets of $\Spec(A)$\textup{)}.
%\end{cor}
%\begin{proof}
%Let $\Zero(J)$, where $J$ is an ideal $B$, 
%be a closed subset of
%$\Spec(B)$. Set $I:=J\cap A$ and let $\pi_J\colon B\to B/J$ 
%and $\pi_I\colon A\to A/I$ be the natural surjections. 
%Let also $\overline{\iota}\colon A/I\to B/J$
%be the ring embedding with $\pi_J\circ \iota = 
%\overline{\iota}\circ\pi_I$.
%Then $\overline{\iota}^*$ is
%surjective by Corollary~\ref{cor:CS}, and
%$\iota^*\circ\pi_J^* = \pi_I^*\circ\overline{\iota}^*$, thus
%$$\iota^*\big(\Zero(J)\big)=
%\iota^*\big(\pi_J^*(\Spec B/J)\big)=
%\pi_I^*\big(\overline{\iota}^*(\Spec B/J)\big)=
%\pi_I^*(\Spec A/I)=\Zero(I)$$
%is closed.
%\end{proof}

\begin{cor}\label{cor:prime ideal chains} Let $\mathfrak{p}, \mathfrak{p}'\in \Spec(A)$ and $\mathfrak{q}\in \Spec(B)$ be such that $\mathfrak{p}\subseteq \mathfrak{p}'$ and 
$\mathfrak{p}=\mathfrak{q}\cap A$. Then there exists
$\mathfrak{q}'\in \Spec(B)$ such that 
$\mathfrak{q}\subseteq \mathfrak{q}'$ and $\mathfrak{p}'=\mathfrak{q}' \cap A$.
\end{cor}
\begin{proof} Apply Corollary~\ref{cor:CS} to the prime ideal
$\mathfrak{p}'/\mathfrak{p}$ of $A/\mathfrak{p}$ and the natural embedding $A/\mathfrak{p} \to B/\mathfrak{q}$.
\end{proof}

\begin{cor}
Let $\mathfrak q\in\Spec(B)$ and $\mathfrak p=\mathfrak q\cap A$.
Then: $$ \text{$\mathfrak p$ is a maximal ideal of~$A$}\ \Longleftrightarrow\    \text{$\mathfrak q$ is a maximal ideal of $B$.}$$
\end{cor}
\begin{proof} Corollary~\ref{cor:primes incomp} gives 
``$\Rightarrow$'' and Corollary~\ref{cor:prime ideal chains} yields ``$\Leftarrow$.''
\end{proof} 

\noindent
Thus any maximal ideal $\mathfrak p$ of $A$
is contained in some prime ideal $\mathfrak q$ of $B$; 
for any such~$\mathfrak p, \mathfrak q$ we have: 
$\mathfrak p=\mathfrak q\cap A$,
$\mathfrak q$ is a maximal ideal of $B$, and the field
$B/\mathfrak q$ is algebraic over the field $A/\mathfrak p$
 (after identifying the latter with its natural image in $B/\mathfrak q$).   

\subsection*{An application} In Section~\ref{sec:diff rings} we shall need the following basic fact about integral domains that are finitely generated over an infinite field:

\begin{lemma}\label{lem:unit alg}
Let the infinite field $K$ be a subring of the integral domain
$B$. Assume that $B$ is finitely generated as a $K$-algebra and 
$x\in B$
is transcendental over~$K$. Then $x-c\in B^\times$ for only finitely many 
$c\in K$.  
\end{lemma}
\begin{proof} 
Let $F$ be the fraction field of $B$.
Take $x_1,\dots,x_n\in B$ such that $B=K[x_1,\dots,x_n]$, $x=x_1$,
and $x_1,\dots,x_r$ is a transcendence basis of $F$  over~$K$, where ${1\le r\le n}$. Take $g\in K[x_1,\dots,x_r]^{\neq}$ such that
$x_{r+1},\dots, x_n\in F$ are integral over the subring
$A:=K[x_1,\dots,x_r,g^{-1}]$ of $F$. Then $B[g^{-1}]\subseteq F$ is 
integral over $A$, so every maximal ideal of $A$ extends to a
maximal ideal of $B[g^{-1}]$. Thus every $K$-algebra morphism $A\to K$ 
extends to a $K$-algebra morphism
$B[g^{-1}]\to  K^\alg$, where $K^\alg$ is an algebraic closure of $K$.
Treating $x_1,x_2,\dots, x_r$ as indeterminates over $K$ and using that $K$ is infinite, we
have $c_2,\dots, c_r\in K$ such that
$g(x_1, c_2,\dots, c_r)\ne 0$ in $K[x_1]$. Let $c\in K$ be such that $g(c,c_2,\dots, c_r)\ne 0$. All but finitely many elements of $K$ satisfy this condition, so it only remains to 
show that $x-c\notin A^\times$. Taking $c_1:= c$ we get a 
$K$-algebra morphism $B\to K$ sending $x_i$ to $c_i$ for $i=1,\dots,r$. 
It extends to a $K$-algebra morphism $A[g^{-1}]\to K^\alg$ sending
$x-c$ to $0$, so $x-c\notin A^{\times}$. 
\end{proof}

\subsection*{Integrally closed domains}
We define an {\bf integrally closed domain\/}  
to be an integral domain
that is integrally closed in its fraction field. 
The ring of integers, and polynomial rings (in any set of variables) over fields are integrally closed domains,
since they are factorial domains (also called {\em unique factorization domains}):

\index{integrally closed!domain}

\begin{lemma}
Factorial domains are integrally closed domains. 
\end{lemma}
\begin{proof} Let $A$ be a factorial domain with fraction field $K$, and suppose $f\in K^{\times}$ is integral over $A$. Then $f^n + a_1f^{n-1} + \cdots + a_{n-1}f+a_n=0$ with coefficients
$a_1,\dots, a_n\in A$. We arrange $f=a/b$ where $a,b\in A^{\ne}$ and no irreducible element of~$A$ divides both
$a$ and $b$ in $A$. Then 
$a^n + a_1a^{n-1}b + \cdots +a_{n-1}ab^{n-1}+ a_nb^n=0$, so~$b$ divides $a^{n}$, which forces $b\in A^\times$, and thus 
$f\in A$.
\end{proof}

\noindent
Next we characterize the integral closure of an integrally closed domain in a field extension of its fraction field:

\begin{lemma}\label{lem:int closure}
Suppose $A$ is an integrally closed domain with fraction field $K$ and ${L\supseteq K}$ is a field extension. An element of $L$ is integral over $A$ if and only if it is algebraic over $K$ and
its minimum polynomial over $K$ has its coefficients in $A$.
\end{lemma}

\begin{proof}
Suppose $b$ is integral over $A$; then clearly $b$ is algebraic over $K$.
Let~$P\in K[X]$ be the minimum polynomial of $b$ over $K$, say of degree $n$. Let $K^{\alg}$ be an algebraic closure of 
$K$. Then in $K^{\alg}[X]$ we have
$P=(X-b_1)\cdots (X-b_n)$ where each $b_i\in K^{\alg}$ is a conjugate of $b$, that is, of the form $\sigma(b)$ for
some embedding $\sigma\colon K(b) \to K^{\alg}$ over~$K$.  
Hence every $b_i$ is integral over~$A$, so the coefficients of~$P$ 
are integral over $A$, and thus lie in $A$, since $A$ is integrally closed. 
\end{proof}

\begin{lemma}\label{lem:fg int closure}
Suppose that $A$ is an integrally closed domain with fraction field $K$,
${L\supseteq K}$ is a separable field extension of finite degree, and $B$ is the integral closure of~$A$ in $L$.
Then $L=K(x)$ for some $x\in B$.
For any such $x$ with minimum polynomial~$P$ over~$K$
we have $B \subseteq P'(x)^{-1}A[x]$.
\end{lemma}

\begin{proof}
The Primitive Element Theorem~\cite[Chapter~V, Theorem~4.6]{Lang} gives $x\in L$ with $L=K(x)$.
Multiplying $x$ by a suitable element of $A^{\ne}$ we get 
$x\in B$.
Let~$P$ be the minimum polynomial of $x$ over $K$; then
$P\in A[X]$ by Lemma~\ref{lem:int closure}.
Take a field extension $M$ of $L$ such that $M|K$ is a Galois extension of finite degree.
Set $G:=\Aut(M|K)$, $H:=\Aut(M|L)$, and take a coset decomposition $$G\ =\ \sigma_1 H\cup\cdots \cup \sigma_n H, \qquad n=[L:K]=\deg P.$$ So $\sigma_1|L,\dots, \sigma_n|L$ are the distinct $K$-embeddings of $L$ into $M$. We take $\sigma_1=\id_M$. 
Note that $P=\prod_{i=1}^n \big(X-\sigma_i(x)\big)$.
For $i=1,\dots,n$ put
$Q_i:=P\big/\big(X-\sigma_i(x)\big)\in M[X]$. Since $x\in B$, we have $b_0,\dots,b_{n-2}\in B$ such that \begin{align*}
Q_1\ &=\ P/(X-x)\ =\ X^{n-1}+b_{n-2}X^{n-2}+\cdots+b_0\quad\text{ hence}\\ 
Q_i\ &=\ \sigma_i(Q_1)\ =\ X^{n-1}+\sigma_i(b_{n-2})X^{n-2}+\cdots+\sigma_i(b_0).
\end{align*}
Also note that $Q_1(x)=P'(x)$ and $Q_i(x)=0$ for $i\geq 2$. Let $b\in B$. Then
\begin{align*}
	     Q_1(x)b\  
		&=\ \sum_{i=1}^n Q_i(x)\sigma_i(b)\
		=\ \sum_{i=1}^n \big(x^{n-1}+\sigma_i(b_{n-2})x^{n-2}+\cdots+\sigma_i(b_0)\big)\sigma_i(b) \\
		&=\ \left(\textstyle\sum\limits_i \sigma_i(b)\right)x^{n-1} + \left(\textstyle\sum\limits_i \sigma_i(b_{n-2}b)\right)x^{n-2}+\cdots+
\textstyle\sum\limits_i \sigma_i(b_0b).
\end{align*}
The coefficients $\sum_i \sigma_i(b),\sum_i \sigma_i(b_{n-2}b),\dots,\sum_i \sigma_i(b_0b)$ are in $K$
since they are tra\-ces of elements of $L$, and they are integral over $A$, hence they are in $A$.
Therefore $P'(x)b=Q_1(x)b\in A[x]$, so  $b\in P'(x)^{-1}A[x]$. 
\end{proof}

\subsection*{Notes and comments}
The notions of
\textit{integral element}\/ and \textit{integral closure}\/ arose from that of an \textit{algebraic integer}\/ 
and the \textit{ring of integers}\/ of an algebraic number field (19th century: Gauss, Dirichlet, Kummer, Dedekind).
The general theory was developed by E.~Noe\-ther~\cite{Noether27}, Krull~\cite{Krull}, and
Cohen-Seidenberg~\cite{CohenSeidenberg}.

\section{Local Rings}\label{sec:local rings}

\noindent
{\em Throughout this section $A$ is a commutative ring}.
Call $A$ {\bf local} \index{ring!local} \index{local ring} if $A$ has exactly one maximal ideal.
Thus valuation rings are local rings. Often we denote 
a local ring~$A$ by $(A,\m)$, thus indicating its maximal ideal~$\m$. A local ring $(A,\m)$ has {\bf residue field} $\k:=A/\m$ \index{local ring!residue field}\index{residue field!local ring} with a surjective ring morphism $a\mapsto\overline{a}:=a+\m\colon A\to\k$.
Local rings have rather good properties compared to arbitrary 
commutative rings. 
This is exemplified by the generation of modules over local rings, the first topic of this section. Next we describe the maximal ideals of certain integral extensions of a local ring.
We then discuss the process of localization, which can often simplify
a problem by reduction to a local ring issue.
Finally, we consider regular sequences, as needed in the study
of regular local rings in Section~\ref{sec:regular local rings}.

\subsection*{Nakayama's Lemma}
{\em In this subsection $M$ is a finitely generated $A$-module}. Thus any set of 
generators of~$M$ has a finite subset generating $M$. 
We denote by~$\mu(M)$ the least $m$ such that $M$ is generated by a subset
of size $m$.
We say that $G\subseteq M$ is a {\bf minimal set of generators} of $M$ if 
$G$ generates $M$, but no proper subset of $G$ does. Any set of generators of~$M$ contains a minimal set of generators of $M$,
and every set of generators of $M$ of size~$\mu(M)$ is minimal.

\index{minimal!set of generators of a module}
\index{module!minimal set of generators}
\nomenclature{$\mu(M)$}{minimal number of generators of the module $M$}

\begin{lemma}[Nakayama]\label{lem:NAK}
Suppose the ideal $I$ of $A$ is contained in every maximal ideal of $A$,
and $IM=M$. Then $M=\{0\}$.
\end{lemma}
\begin{proof}
Towards a contradiction, assume $n=\mu(M)\geq 1$. Take  $x_1,\dots,x_n\in M$ with  $M=Ax_1+\cdots+Ax_n$. Since $IM=M$, we can take $a_1,\dots,a_n\in I$ such that 
$x_n=a_1x_1+\cdots+a_nx_n$. Then 
$(1-a_n)x_n\in Ax_1+\cdots+Ax_{n-1}$ where 
$1-a_n\in A^\times$. Hence $M=Ax_1+\cdots+Ax_{n-1}$, a contradiction.
\end{proof}

\noindent
In the rest of this subsection $(A,\mathfrak m)$ is a local ring with residue field $\k=A/\m$.
We view $\overline{M}:=M/\m M$ as a $\k$-linear space in the natural way, and for $x\in M$ we set $\overline{x}:=x+\m M\in\overline{M}$. 

\begin{cor}\label{cor:NAK, 1}
Let $N$ be a submodule of $M$.
If $M=N+\mathfrak m M$, then $M=N$.
\end{cor}
\begin{proof}
Apply Lemma~\ref{lem:NAK} to $\m$ and $M/N$ in place of $I$ and $M$.
\end{proof}

\begin{cor}\label{cor:NAK, 2}
Let $x_1,\dots,x_n\in M$. Then
$$M\ =\ A\,x_1+\cdots+A\,x_n\qquad\Longleftrightarrow\qquad\overline{M} = \k \,\overline{x_1}+\cdots+\k\,\overline{x_n}.$$
\end{cor}
\begin{proof}
Apply the previous corollary to $N=A\,x_1+\cdots+A\,x_n$.
\end{proof}

\noindent
Familiar properties of bases of $\k$-vector spaces and the preceding corollary yield:

\begin{cor}\label{cor:NAK, 3}
Let $x_1,\dots,x_m\in M$ be distinct. Then:
\begin{enumerate}
\item[\textup{(i)}] $\{x_1,\dots,x_m\}$ is a minimal set of generators of $M$ if and only if
$\overline{x_1},\dots,\overline{x_m}$ is a basis of the $\k$-linear space $\overline{M}$;
thus $\mu(M)=\dim_{\k} \overline{M}$;
\item[\textup{(ii)}] $\{x_1,\dots,x_m\}$ is contained in a minimal set of generators of $M$ if and only if $\overline{x_1},\dots,\overline{x_m}$ are $\k$-linearly independent.
\end{enumerate}
\end{cor}

\noindent
%The following corollary is used in the proof of 
%Proposition~\ref{prop:reg => domain} below.
Let $x\in\fm\setminus\fm^2$; then $A^*:=A/Ax$ is a local
ring with maximal ideal $\m^*:=\m/Ax$. For $a\in A$, put $a^*:=a+Ax\in A^*$.

\begin{cor}\label{cor:NAK, 4}
Let $x_1,\dots,x_n\in\m$, and suppose $x_1^*,\dots,x_n^*$ are distinct, and
$\{x_1^*,\dots,x_n^*\}$ is a minimal
set of generators of $\m^*$. Then $x\notin \{x_1,\dots,x_n\}$, and
$\{x_1,\dots, x_n,x\}$ is a minimal set of
generators of $\m$. 
%Thus $\mu(\m)=\mu(\m^*)+1$.
\end{cor}
\begin{proof}
Clearly we have $x\notin \{x_1,\dots, x_n\}$, and $x_1,\dots,x_n,x$ generate  $\fm$.
For the  set $\{x_1,\dots, x_n,x\}$ of
generators of~$\m$ to be minimal, it is enough by Corollary~\ref{cor:NAK, 3} that the 
residue classes 
$x_1+\m^2,\dots, x_n+\m^2, x+\m^2\in \m/\m^2$
are $\k$-linearly independent. Thus, let $a_1,\dots,a_n,a\in A$ be such that
$$a_1x_1+\cdots+a_nx_n + ax\in\m^2.$$
Then in $A^*$  we have
$$a_1^*x_1^*+\cdots+a_n^*x_n^*\in(\fm^*)^2$$
and hence $a_1^*,\dots,a_n^*\in\fm^*$ by hypothesis and 
Corollary~\ref{cor:NAK, 3},
so $a_1,\dots,a_n\in\fm$. This yields $ax\in\fm^2$ and thus $a\in\fm$, since $x\notin\fm^2$.
\end{proof}

\begin{cor}\label{cor:NAK, 5}
Suppose $\fm$ is finitely generated. Then either:
\begin{enumerate}
\item[\textup{(i)}] $\m^n \neq \m^{n+1}$ for all $n$, or
\item[\textup{(ii)}] $\m^n = \{0\}$ for some $n$, in which case $\Spec(A)=\{\m\}$.
\end{enumerate}
\end{cor} 
\begin{proof}
Suppose $\m^n = \m^{n+1}$. Then $\m^n=\{0\}$ by Nakayama's Lemma, and so for $\mathfrak p\in\Spec(A)$ we have
$\m^n\subseteq\mathfrak p$, hence 
$\m\subseteq \mathfrak p$.
\end{proof}

\noindent
Here is a characterization of DVRs among local rings:

\begin{lemma}\label{NAKDVR} Suppose $A$ is noetherian,
$t\in A$ is not a zero divisor of $A$, and $\m=tA$. Then $A$ is a $\operatorname{DVR}$. 
\end{lemma}
\begin{proof} Using that $t$ is not a zero divisor we get
$\m I=I$ for $I:= \bigcap_n \m^n=\bigcap_n t^n A$, so $I=\{0\}$ by Nakayama's Lemma. Therefore, given
$a\in A^{\ne}$, we have $n$ with $a\in t^nA\setminus t^{n+1}A$, and then $a=t^nu$, $u\in A^{\times}$. Thus $A$ is a DVR.
\end{proof}

\subsection*{The maximal ideals of an integral extension of a local ring} Let $A$ be a local ring with maximal ideal $\mathfrak{m}$ and $X$ an indeterminate. We extend the surjective ring morphism $a\mapsto\overline{a}:=a+\mathfrak m\colon A\to \k:= A/\mathfrak{m}$ to a surjective ring morphism $P\mapsto\overline{P}\colon A[X]\to \mathbf k[X]$ that sends $X$ to $X$. With these notations we have:

\begin{lemma} \label{lem:Dedekind}
Let $P\in A[X]$ be monic of positive degree, and $A[x]:=A[X]/(P)$ with $x:=X+(P)$.
Suppose $\overline{P}\in \mathbf k[X]$ factors as
$$\overline{P}\ =\ \overline{P_1}^{e_1} \cdots \overline{P_n}^{e_n}$$
where each $e_i\ge 1$, each $P_i\in A[X]$ is monic, and $\overline{P_1},\dots,\overline{P_n}\in \mathbf k[X]$ are irreducible and distinct. Then 
$\mathfrak m_1,\dots,\mathfrak m_n$ with
$$\mathfrak m_i\ :=\ \mathfrak m A[x] + P_i(x)A[x] $$ 
are the distinct maximal ideals of $A[x]$. 
% and $A[x]/\mathfrak m_i\cong 
%\mathbf k_A[X]/(\overline{P}_i)$ as rings.
\end{lemma}
\begin{proof}
Consider for each $i$ the composite of the natural surjections
$$A[x]=A[X]/(P) \to \k[X]/(\overline{P})\to \k[X]/(\overline{P_i}).$$
It is easy to see that the composite map has kernel $\mathfrak m_i$; since $\k[X]/(\overline{P_i})$ is a field,~$\mathfrak m_i$ is a maximal ideal of $A[x]$. If $1\le i < j\le n$, then $\overline{P_i}\not\equiv 0 \mod \overline{P_j}$ in $\k[X]$, so $\mathfrak m_i\neq\mathfrak m_j$. It remains to show that the $\mathfrak m_i$ are the only maximal ideals of $A[x]$. First note that~$A[x]$ is integral over $A$, so each maximal ideal of $A[x]$ contains $\mathfrak{m}$. The polynomial $P-\prod_iP_i^{e_i}$ is in $\mathfrak{m}A[X]$, and $P(x)=0$, so 
$\prod_i P_i(x)^{e_i}\in \mathfrak{m}A[x]$. Thus each maximal ideal of $A[x]$ contains some $P_i(x)$, and thus equals $\mathfrak m_i$ for some $i$.
\end{proof}

\subsection*{Localization} 
%See \cite[Chapter~II, \S{}4]{Lang} for (easy) proofs of some 
%facts not proved in this subsection.
Let $S$ be a multiplicative subset of $A$. 
Recall the construction of the localization $S^{-1}A$ of $A$ at $S$:  
this is the ring whose elements are the equivalence classes
of the equivalence relation $\sim$ on $A\times S$ defined by
$$ (a_1,s_1) \sim (a_2,s_2) \quad :\Longleftrightarrow \quad \text{$a_1 s_2 s= 
a_2 s_1 s$ for some $s\in S$},$$
with the equivalence class of $(a,s)\in A\times S$ denoted by 
$\frac{a}{s}$ or $a/s$, and with addition and multiplication given by
$$\frac{a_1}{s_1} + \frac{a_2}{s_2}\ =\ \frac{a_1s_2+a_2s_1}{s_1s_2}, \qquad
\frac{a_1}{s_1}\cdot \frac{a_2}{s_2}\ =\ \frac{a_1a_2}{s_1s_2}.$$
Then $S^{-1}A$ has zero element $0/1$ and identity $1/1$, and
$S^{-1}A=\{0\}\Longleftrightarrow 0\in S$.
The map 
$\iota=\iota_A^S\colon A\to S^{-1}A$ sending $a\in A$ to  $a/1$ 
is a ring morphism
with $$\iota(S)\ \subseteq\ (S^{-1}A)^\times, \qquad \ker\iota\ =\ \{a\in A:\ \text{$as=0$ for some $s\in S$}\}.$$ 
We recall the key universal property of $\iota\colon A\to S^{-1}A$: 
for every ring morphism $\phi\colon A\to B$ into a commutative ring $B$ with 
$\phi(S)\subseteq B^\times$ there is a unique ring morphism
$\phi'\colon S^{-1}A \to B$ such that $\phi=\phi'\circ \iota$.  
%If $T\subseteq A$ is also multiplicative with 
%$S\subseteq T$, then
%we have a unique ring morphism $S^{-1}A\to T^{-1}A$ which %sends $\iota^S(a)$ to $\iota^T(a)$, for $a\in A$.

\nomenclature[A]{$S^{-1}A$}{localization of $A$ at its multiplicative subset $S$}
\index{localization}

\medskip
\noindent
Let $I$ be an ideal of $A$. Then
$ S^{-1}I := \left\{\textstyle\frac{a}{s}: a\in I,\ s\in S\right\}$
is the ideal of $S^{-1}A$ generated by 
$\iota(I)$. The ideal
$$\iota^{-1}(S^{-1}I)\ =\ \{ a\in A:\ \text{$as\in I$ for some $s\in S$} \}$$
of $A$ contains $I$, with
$\iota^{-1}\big(S^{-1}I)=A$ iff $I\cap S\neq\emptyset$.
If $J$ is an ideal of $S^{-1}A$ and $I=\iota^{-1}(J)$, then 
$J=S^{-1}I$. Hence if $A$ is noetherian, then so is $S^{-1}A$. 
An ideal $J$ of $S^{-1}A$ is prime iff 
the ideal $\iota^{-1}(J)$ of $A$ is prime and disjoint from $S$. Thus: 

\begin{cor}
The map $\iota^*\colon\Spec(S^{-1}A)\to\Spec(A)$ is a homeomorphism onto its image $\big\{\mathfrak p\in\Spec(A):\ \mathfrak p\cap S=\emptyset\big\}$, with inverse
given by
$\mathfrak{p}\mapsto S^{-1}\mathfrak{p}$. 
\end{cor}

\begin{example}  
For $S=\{1,s,s^2,\dots\}$ with $s\in A$ we denote $S^{-1}A$ also by $A_s$, and
%and call $R_s$  the localization of $R$ at $s$. 
the image of $\iota^*$ in this case is $\operatorname{D}(s)=\big\{\mathfrak p\in\Spec(A):s\notin\mathfrak p\big\}$.
\end{example}

\noindent
Localization produces local rings: Let $\mathfrak p\in\Spec(A)$; then  $S=A\setminus\mathfrak p$ is a multiplicative subset of $A$ with $0\notin S$, and the localization $S^{-1}A$ of $A$ 
at $S$ is then denoted by~$A_{\mathfrak p}$ and called the localization 
of $A$ at $\mathfrak p$. 
The ring $A_{\mathfrak p}$ is indeed local, with maximal ideal~$\mathfrak p{A_{\mathfrak p}}$: by the preceding corollary $\iota^*$ maps 
$\Spec A_{\mathfrak p}$ bijectively onto
$$\big\{\mathfrak q\in\Spec(A):\ \mathfrak q\subseteq\mathfrak p\big\}.$$ 
The  morphism $\iota\colon A\to A_{\mathfrak p}$ induces an embedding 
$A/\mathfrak p\to A_{\mathfrak p}/\mathfrak p{A_{\mathfrak p}}$ of the domain~$A/\mathfrak p$ into the field $F:=A_{\mathfrak p}/\mathfrak p{A_{\mathfrak p}}$,
and $F$ is the fraction field of the image of this embedding
in the sense that every $f\in F$ equals $s^{-1}a$ for some $a$,~$s$ in the image of this embedding,~$s\ne 0$.   
%If $A$ is local with maximal ideal $\mathfrak p$, 
%then $\iota$  is an isomorphism.  
%If $\mathfrak q\subseteq \mathfrak p$ are prime ideals 
%of $R$, then we have
%a ring morphism $R_{\mathfrak p}\to R_{\mathfrak q}$ sending $r/1\in R_{\mathfrak p}$ to $r/1\in R_{\mathfrak q}$, 
%for $r\in R$.

%\begin{example}
%If $R$ is an integral domain, then $\mathfrak p=\{0\}$ 
%is a prime ideal of $R$,   and $R_{\mathfrak p}$ is 
%the field of fractions of $R$.
%\end{example}

\subsection*{Localization in integral domains} Suppose $A$ is an integral domain
and $S$ is a multiplicative subset of $A$ with $0\notin S$.
Then we have the ring embedding 
$$ a/s\mapsto s^{-1}a\ :\ S^{-1}A \to \Frac{A} \qquad(a\in A,\ s\in S),$$
via which, throughout this volume, we identify $S^{-1}A$ with a subring of the fraction
field $\Frac{A}$ of $A$. (Thus $S^{-1}A=\Frac{A}$ for 
$S= A^{\ne}$.) Note that if $A$ is a PID, then so is 
$S^{-1}A$. For $\mathfrak{p}\in \Spec(A)$ and $S=A\setminus \mathfrak{p}$ we obtain the local domain $A_{\mathfrak{p}}$. 
Thus if~$A$ is a PID and
$\{0\}\ne \mathfrak{p}\in \Spec(A)$, then 
$A_{\mathfrak p}$ is a DVR by Lemma~\ref{NAKDVR}.

\begin{lemma}\label{lem:local ring incl} Let  $\mathfrak m$ be a maximal ideal of the integral domain $A$, and let  $B$ be a local ring
such that $A\subseteq B\subseteq A_{\mathfrak m}$. Then $B=A_{\mathfrak m}$.  
\end{lemma}
\begin{proof}
Let $\mathfrak n$ be the maximal ideal of $B$ and $s\in A\setminus\m$; it is enough to show that then $s\in B^\times$,
that is, $s\notin \mathfrak n$.
Take $t\in A$ with $st\equiv 1\bmod\m$, so $st-1\in\m$. Since $\m\subseteq \m A_{\m}\cap B\subseteq\mathfrak n$,
this yields $st-1\in\mathfrak n$, so $s\notin\mathfrak n$.
\end{proof}

\noindent
We can now complete Lemma~\ref{lem:fg int closure} as follows:

\begin{cor}\label{cor:fg int closure}
Let $A$ be an integrally closed domain with fraction field 
$K$, let $L\supseteq K$ be a separable field extension of finite degree. Let $B$ be the integral closure of~$A$ in $L$,
let $x\in B$ with $L=K(x)$ have minimum polynomial $P$ over~$K$, and let $\m$ and $\mathfrak n$ be maximal ideals of $A[x]$ and
$B$, respectively, such that $P'(x)\notin\m=\mathfrak n \cap A[x]$. 
Then we have $A[x]_{\mathfrak m}=B_{\mathfrak{n}}$. 
\end{cor}
\begin{proof}
From $P'(x)\notin\m$ and Lemma~\ref{lem:fg int closure} we get $B\subseteq  A[x]_{\mathfrak m} \subseteq B_{\mathfrak n}$.
Now Lemma~\ref{lem:local ring incl} applied to $B$, $A[x]_{\mathfrak m}$, $\mathfrak n$ in place
of $A$, $B$, $\m$, yields $A[x]_{\mathfrak m}=B_{\mathfrak n}$.
\end{proof}

\noindent
{\em In the next four lemmas $A$ is an integral domain and $S$ is a multiplicative subset of $A$
with $0\notin S$}. So for $\mathfrak p\in\Spec(A)$ with $\mathfrak p\cap S=\emptyset$ we have $S^{-1}\mathfrak{p}\in \Spec(S^{-1}A)$. We set $K:=\Frac(A)$. %\marginpar{clarify $S$}

\begin{lemma}\label{locloc} 
If $\mathfrak p\in\Spec(A)$ and $\mathfrak p\cap S=\emptyset$,
then $(S^{-1}A)_{S^{-1}\mathfrak p}=A_{\mathfrak p} $ in $K$.
\end{lemma} 

\noindent
Let $X$ be an indeterminate over $A$ in the next two lemmas.

\begin{lemma} 
We have $(S^{-1}A)[X] = S^{-1}\big(A[X]\big)$ inside  $K(X)$. 
\end{lemma}

\begin{lemma}\label{cor:loc of polynomial ring}
Let $\mathfrak P\in\Spec\!\big(A[X]\big)$. Then $\mathfrak p:=\mathfrak P\cap A\in\Spec(A)$, $\mathfrak{P}$ generates a prime ideal 
$\mathfrak{Q}$ in $A_{\mathfrak p}[X]\subseteq K[X]$, and
$A[X]_{\mathfrak P} = A_{\mathfrak p}[X]_{\mathfrak Q}$ inside 
$K(X)$.
\end{lemma}
\begin{proof} The set $S:= A\setminus \frak{p}$ is a multiplicative subset
of $A$ as well as of $A[X]$, and 
$\mathfrak P\cap S=\emptyset$. So $\mathfrak{P}$ does generate a prime ideal
$\mathfrak{Q}=S^{-1}\mathfrak{P}$ in $S^{-1}\big(A[X]\big)= A_{\mathfrak p}[X]$.
Now apply Lemma~\ref{locloc} to $A[X]$ and  $\mathfrak{P}$ in the role of
$A$ and  $\frak{p}$.
\end{proof}

\noindent
In Section~\ref{sec:johnson} we use:
 
\begin{lemma}\label{lem:expansion-contraction}
Let $I$ be an ideal of $A$ and $\mathfrak p\in\Zero(I)$. Then the ideal
$$\iota^{-1}(IA_{\mathfrak p})\ =\ \{a\in A:\ \text{$as\in I$ for some $s\in A\setminus\mathfrak p$}\}$$
of $A$ satisfies $I\subseteq \iota^{-1}(IA_{\mathfrak p})\subseteq\mathfrak p$ and is contained in 
every prime ideal~$\mathfrak q$ of $A$ with $I\subseteq\mathfrak q\subseteq\mathfrak p$.
If $\mathfrak q:=\iota^{-1}(IA_{\mathfrak p})$ itself is prime, then $\mathfrak q A_{\mathfrak q}=I A_{\mathfrak q}$.
\end{lemma}
\begin{proof}
The first statement is clear. Suppose $\mathfrak q:=\iota^{-1}(IA_{\mathfrak p})$ is prime.
We have $\mathfrak q A_{\mathfrak p}=I A_{\mathfrak p}$, and applying the ring morphism $A_{\mathfrak p}\to A_{\mathfrak q}$ yields $\mathfrak q A_{\mathfrak q}=I A_{\mathfrak q}$.
\end{proof}

\subsection*{Localization of modules} Let $S$ be a multiplicative subset of $A$. Every $A$-module~$M$
gives rise to the $S^{-1}A$-module $S^{-1}M$, whose 
elements are the formal fractions $x/s$ ($x\in M$, $s\in S$), with
$$x/s\ =\ 0\ \Longleftrightarrow\ \text{$tx=0$ for some $t\in S$}.$$ 
Addition of these fractions is given by
$$(x_1/s_1) + (x_2/s_2)\ =\ (s_2x_1+s_1x_2)/s_1s_2  \qquad(x_1, x_2\in M,\ s_1,s_2\in S)$$
and their multiplication by scalars from $S^{-1}A$ is given by 
$$(a/s)(x/t)\ =\ ax/st \qquad(a\in A,\ s,t\in S,\ x\in M).$$ 
Any $S^{-1}A$-module is construed as an $A$-module via $\iota\colon A\to S^{-1}A$. The map $\iota_M=\iota_M^S\colon M \to S^{-1}M$ defined by $\iota_M(x)=x/1$ is
$A$-linear, and for any $A$-linear map $f\colon M \to N$ into an $S^{-1}A$-module $N$
there is a unique $S^{-1}A$-linear map $f'\colon S^{-1}M\to N$ such that 
$f=f'\circ\iota_M$.

\subsection*{Regular sequences}  
Let $M$ be an $A$-module. Call $a\in A$ a {\bf zero divisor on $M$} if
$ax=0$ for some $x\in M^{\ne}$. Thus $a\in A$ is not a zero divisor on $M$ iff 
the $A$-linear map $x\mapsto ax\colon M\to M$ is injective.
Let $\vec r=(r_1,\dots,r_n)\in A^n$, $n\geq 1$. 

\index{module!zero divisor}
\index{zero divisor}

\begin{definition}
The sequence $\vec{r}$ is called {\bf regular on $M$}, 
if 
\begin{enumerate}
\item $M\neq r_1M+\cdots+r_nM$; and
\item $r_{j}$ is not a zero divisor on $M/(r_1M+ \cdots+r_{j-1}M)$, 
for $j=1,\dots,n$.
\end{enumerate}
\end{definition}

\index{sequence!regular}
\index{regular!sequence}

\noindent
This notion, applied to the $A$-module $A$,
is motivated by:

\begin{example} Let $K$ be a commutative ring with $0\ne 1$, let $X_1,\dots, X_n$ be 
distinct indeterminates, $n\ge 1$, $A=K[X_1,\dots,X_n]$. Then 
$(X_1,\dots,X_n)$ is regular on $A$.  
\end{example}

\noindent
It is easy to see that for $1\le m < n$,
\begin{equation}\label{eq:split regular sequ}
\text{$\vec r$ is regular on $M$}\ \Longleftrightarrow\  
\begin{cases} &
\parbox{13em}{$(r_1,\dots,r_m)$ is regular on $M$ and $(r_{m+1},\dots,r_n)$ is
regular on $M/(r_1M+ \cdots+r_{m}M)$.}
\end{cases}
\end{equation}
The following are also easy to verify:

\begin{lemma}\label{lem:reg sequence}
If $\vec r$ is regular on $A$ and $X$ is an indeterminate, then
$\vec r$ is regular on the $A[X]$-module $A[X]$.
\end{lemma}

\begin{lemma}\label{lem:reglocal}
If $\vec r=(r_1,\dots, r_n)$ is regular on $A$ and $r_1,\dots, r_n\in \mathfrak p$ where $\mathfrak p\in\Spec(A)$,
then $\big(\frac{r_1}{1},\dots,\frac{r_n}{1}\big)$ is regular on the 
$A_{\mathfrak p}$-module 
$A_{\mathfrak p}$.
\end{lemma}

\noindent
{\em In the next lemma and its corollary $X$ is an indeterminate and $I$ is
an ideal of $A[X]$ 
such that $\m:=I\cap A$ is a maximal ideal of $A$}.

\begin{lemma}\label{Xregzero}
If $I\neq \m A[X]$, then
$I=\m A[X] + PA[X]$ where $P\in A[X]$ is not a zero divisor on
the $A[X]$-module 
$A[X]/\m A[X]$.
\end{lemma}

{\sloppy
\begin{proof} Extend the canonical map $A\to A/\m$ to the ring morphism
$A[X]\to (A/\mathfrak m)[X]$ by sending $X$ to $X$. Since
the kernel of this extension is $\mathfrak m A[X]$, we obtain a
ring isomorphism $A[X]/\m A[X]\xrightarrow{\cong}(A/\m)[X]$.
Now use that the polynomial ring $(A/\m)[X]$  over the 
field $A/\m$ is a domain and that its ideals are principal.
\end{proof}
}
\noindent
From Lemmas~\ref{lem:reg sequence} and~\ref{Xregzero} we obtain: 

\begin{cor}\label{cor:reg poly}
If $\vec r=(r_1,\dots,r_n)$ is regular on $A$ and $\m=r_1A+\cdots+r_nA$, and 
$I\ne \m A[X]$,
then there exists $P\in A[X]$ such that $I=r_1A[X]+\cdots+r_nA[X]+PA[X]$ and
$(r_1,\dots,r_n,P)$ is regular on the $A[X]$-module $A[X]$.
\end{cor}

\subsection*{Notes and comments} 
Nakayama's Lemma in its various forms is due to Krull, Azumaya, and Nakayama;
% Krull and Azumaya~\cite{Azumaya}; 
see the discussion in \cite[pp.~212--213]{Nagata62}. Early studies of
local rings are 
Krull~\cite{Krull38} and Chevalley~\cite{Chevalley43}. Lemma~\ref{lem:Dedekind} is essentially 
due to Kummer~\cite{Kummer} and Dedekind~\cite{Dedekind}.
Localization goes back to Grell~\cite{Grell}, and in the
generality above, to
Chevalley~\cite{Chevalley44}
and Uzkov~\cite{Uzkov}. 
%Localization of non-commutative rings was first 
%investigated by Ore~\cite{Ore31}.

\section{Krull's Principal Ideal Theorem}\label{sec:PIT}

\noindent
{\em Throughout this section $R$ is a commutative ring.}\/   The {\bf height} of a prime ideal~$\mathfrak p$ of $R$, denoted by $\operatorname{ht}(\mathfrak p)$, is the supremum, in $\N\cup\{\infty\}$, of the lengths $n$ of strictly 
increasing sequences 
$\mathfrak p_0 \subset \mathfrak p_1 \subset \cdots \subset \mathfrak p_n=\mathfrak{p}$ of prime ideals of $R$.
The (Krull) {\bf dimension} of~$R$~is 
$$\dim R\ :=\ \sup\big\{\!\operatorname{ht}(\mathfrak p):\  \mathfrak p\in \Spec(R)\big\} $$
with supremum in $\N\cup\{-\infty, +\infty\}$. So
$\dim R = \dim \Spec(R)$
by Lemma~\ref{lem:irred comp}, with $\dim \Spec(R)$ defined as at the end of 
Section~\ref{sec:ztn}. More generally, for each ideal $I$ of~$R$, we have $\dim R/I = \dim \Zero(I)$. Note:
 $$\dim R = -\infty\ \Longleftrightarrow\  \text{$1=0$ in $R$}\ \Longleftrightarrow\ \Spec(R)=\emptyset.$$ 
%more generally, for each ideal $I$ of $R$, we have $\dim R/I = \dim \Zero(I)$.

\index{height}
\index{ring!Krull dimension}
\index{dimension!Krull}
\index{Krull dimension!commutative ring}

\begin{examples}
A prime ideal of $R$ has height $0$ 
if and only if it is  minimal. The ring~$R$ has dimension $0$ if and only if $R\neq\{0\}$ and every prime ideal of $R$ is maximal; thus an integral domain has dimension $0$ if and only if it is a field.
See Proposition~\ref{prop:dim 0 rings}  for a characterization of $0$-dimensional noetherian commutative rings.
\end{examples}

\noindent
Clearly
$\operatorname{ht}(\mathfrak p) + \dim(R/\mathfrak p) \leq \dim R$.
If $R$ is local with maximal ideal $\m$, then $\operatorname{ht}(\m)=\dim R$.
Therefore, if $\mathfrak p$ is a prime ideal of $R$, then $\operatorname{ht}(\mathfrak p)=\operatorname{ht}(\mathfrak p R_{\mathfrak p})=\dim R_{\mathfrak p}$.

\medskip
\noindent
The following is clear from the definitions:
 
\begin{lemma}\label{lem:ht under morphisms}
Let $\varphi \colon R\to S$ be a surjective ring morphism
and $\mathfrak p$ be a prime ideal of $R$ with $\mathfrak p\supseteq\ker\varphi$. Then
$\varphi(\mathfrak p)$ is a prime ideal of $S$, and $\operatorname{ht}(\mathfrak p)\ge \operatorname{ht}(\varphi(\mathfrak{p}))$.
\end{lemma}

\noindent
From Corollaries~\ref{cor:primes incomp},~\ref{cor:CS}, and~\ref{cor:prime ideal chains} we obtain: 

\begin{cor} If $R$ is a subring of the commutative ring $S$, and $S$ is integral over~$R$, then $\dim(R)=\dim(S)$.
\end{cor}

\noindent
\textit{In the rest of this section we assume that $R$ is noetherian.}\/
The following theorem is a key result about Krull dimension in this 
noetherian setting:

\begin{theorem}[Principal Ideal Theorem]\label{thm:PIT}
If $x\in R$ and $\mathfrak p$ is a minimal prime divisor of $Rx$,
then $\operatorname{ht}(\mathfrak p)\leq 1$.
\end{theorem}

\noindent
Since $\operatorname{ht}(\mathfrak p)=\dim R_{\mathfrak p}$ for $\mathfrak{p}\in \Spec(R)$, this immediately follows from the case where~$R$ is a local ring and $\mathfrak p$ is its maximal ideal, treated in the next lemma:

\begin{lemma} 
Suppose $(R,\mathfrak m)$ is a local ring, $x\in R$, and $\mathfrak m$ is a minimal prime divisor of $Rx$. Then $\dim R\leq 1$.
\end{lemma}
\begin{proof}
Let $\mathfrak q\neq\mathfrak m$ be a prime ideal of $R$;
we need to show that $\operatorname{ht}(\mathfrak q)=0$. Let 
$\iota\colon R\to R_{\mathfrak q}$ be the morphism given by $\iota(r)=\frac{r}{1}$. Set
$\mathfrak q^{(n)} := \iota^{-1}(\mathfrak q^n R_{\mathfrak q})$ 
(the \textit{$n$th symbolic power}\/ of $\mathfrak q$)
and consider the
chain 
$$\mathfrak q^{(1)}+Rx\ \supseteq\ \mathfrak q^{(2)}+Rx\ \supseteq \cdots $$
of ideals of $R$. The ring $\overline{R}:=R/Rx$ has exactly one prime ideal (namely $\mathfrak m/Rx$), so $\dim\overline{R}=0$. Then Proposition~\ref{prop:dim 0 rings} gives $n$ such that $\mathfrak q^{(n)} + Rx = \mathfrak q^{(n+1)} + Rx$.

\claim{$\mathfrak q^{(n)} =  \mathfrak q^{(n+1)}+\mathfrak m \mathfrak q^{(n)}$.}

\noindent
To prove this claim, let $a\in \mathfrak q^{(n)}$. Then $a=b+rx$ where $b\in \mathfrak q^{(n+1)}$, $r\in R$.
Since $rx\in \mathfrak q^{(n)}$ and $x\notin\mathfrak q$, we obtain $r\in\mathfrak q^{(n)}$ and hence
$a=b+rx\in \mathfrak q^{(n+1)}+\mathfrak m \mathfrak q^{(n)}$.

\medskip
\noindent
The claim and Corollary~\ref{cor:NAK, 1} yields $\mathfrak q^{(n)} = \mathfrak q^{(n+1)}$ and hence
$\mathfrak q^n R_{\mathfrak q} = \mathfrak q^{n+1} R_{\mathfrak q}$.
Thus 
$\operatorname{ht}(\mathfrak q)=\dim R_{\mathfrak q}=0$ by Corollary~\ref{cor:NAK, 5}, as required.
\end{proof}

\noindent
Combining Theorem~\ref{thm:PIT} with Corollary~\ref{cor:irred comp 2} gives:
 
\begin{cor}\label{cor:PIT}
If $x\in R$ is not a zero divisor and  $\mathfrak p$ is a minimal prime 
divisor of~$Rx$, then $\operatorname{ht}(\mathfrak p)=1$.
\end{cor}

\noindent
Below we extend Theorem~\ref{thm:PIT} to non-principal ideals. First a lemma:

\begin{lemma}\label{lem:GPIT}
Let $\mathfrak p,\mathfrak q_1,\dots,\mathfrak q_m$ be prime ideals of $R$ 
such that $\mathfrak p\not\subseteq\bigcup_{i=1}^m \mathfrak q_i$, and
$\operatorname{ht}(\mathfrak p)\ge n\ge 1$. 
Then there is a strictly increasing sequence
$\mathfrak p_0\subset \mathfrak p_1\subset\cdots\subset \mathfrak p_n$ of prime ideals of $R$ with $\mathfrak p_n=\mathfrak p$ and $\mathfrak p_j\not\subseteq\bigcup_{i=1}^m \mathfrak q_i$ for $j=1,\dots,n$.
\end{lemma}
\begin{proof}
By induction on $n$. The case $n=1$ being trivial, suppose~$n\geq 2$.
Since $\operatorname{ht}(\mathfrak p)\geq n\geq 2$, we can take a prime ideal $\mathfrak q\subseteq \mathfrak{p}$ of $R$  with   
$\operatorname{ht}(\mathfrak p/\mathfrak q)\geq 2$ in $R/\mathfrak{q}$ and $\operatorname{ht}(\mathfrak q)\geq n-2$.
Lemma~\ref{lem:prime avoidance} provides $x\in\mathfrak p\setminus\left(\mathfrak q\cup\bigcup_{i=1}^m \mathfrak q_i\right)$.
Take a minimal prime divisor $\mathfrak p'\subseteq \mathfrak{p}$ of $xR+\mathfrak q$. 
Since
$\operatorname{ht}(\mathfrak p'/\mathfrak q)\leq 1$ by Theorem~\ref{thm:PIT}, we have
$\mathfrak p'\neq\mathfrak p$. From $x\in\mathfrak p'$ we get
$\mathfrak p'\not\subseteq \bigcup_{i=1}^m \mathfrak q_i$. Now apply the inductive hypothesis
to $\mathfrak p'$.
\end{proof}

\begin{theorem}[Generalized Principal Ideal Theorem]\label{thm:GPIT}
If $\mathfrak p$ is a minimal prime divisor of an ideal of $R$  generated by 
$n$ elements, then $\operatorname{ht}(\mathfrak p)\leq n$.
\end{theorem}
\begin{proof}
By induction on $n$. The case $n=0$ is trivial. Let $n\geq 1$, and
let $\mathfrak p$ be a minimal prime divisor of 
$I=Rx_1+\cdots+Rx_n$, $x_1,\dots, x_n\in R$.
Let $\mathfrak q_1,\dots,\mathfrak q_m$ be the minimal prime divisors of $J:=Rx_1+\cdots+Rx_{n-1}$; thus $\operatorname{ht}(\mathfrak q_i)\leq n-1$ for $i=1,\dots,m$ by inductive assumption.
If $\mathfrak p\subseteq\bigcup_{i=1}^m \mathfrak q_i$, then for some $i\in\{1,\dots,m\}$ we have $\mathfrak p\subseteq\mathfrak q_i$
and hence $\operatorname{ht}(\mathfrak p)\leq\operatorname{ht}(\mathfrak q_i)\leq n-1$, and we are done. So assume $\mathfrak p\not\subseteq\bigcup_{i=1}^m \mathfrak q_i$.
Towards a contradiction, suppose that $\operatorname{ht}(\mathfrak p)\geq n+1$.
Lemma~\ref{lem:GPIT} yields a strictly increasing sequence 
$\mathfrak p_0\subset \mathfrak p_1\subset\cdots\subset\mathfrak p_{n+1}$
of prime ideals of $R$ with $\mathfrak p_{n+1}=\mathfrak p$ and $\mathfrak p_j\not\subseteq\bigcup_{i=1}^m \mathfrak q_i$ for $j=1,\dots,n$.
Now consider the ring $\overline{R}:=R/J$, whose minimal prime ideals are the  
$\overline{\mathfrak q_i}:=\mathfrak q_i/J$ ($i=1,\dots,m$). Setting $\overline{\mathfrak p}:=\mathfrak p/J$,
we have $\operatorname{ht}(\overline{\mathfrak p})\leq 1$ by Theorem~\ref{thm:PIT}.
Therefore $\overline{\mathfrak p}$ is a minimal prime divisor of the ideal $(\mathfrak p_{1}+J)/J$ of $\overline{R}$, since $\mathfrak p_{1}\not\subseteq\bigcup_{i=1}^m \mathfrak q_i$. Hence $\mathfrak p$ is a minimal prime divisor
of~$\mathfrak p_{1}+J$.
Thus the prime ideal $\mathfrak p/\mathfrak p_{1}$ of $R/\mathfrak p_{1}$ has height $\geq n$ and is a
minimal prime divisor of the
ideal generated by $x_1+\mathfrak p_{1},\dots,x_{n-1}+\mathfrak p_{1}$ in $R/\mathfrak p_{1}$,  
in contradiction to the inductive hypothesis.
\end{proof}

\noindent
Important consequences of Theorem~\ref{thm:GPIT} are the following:

\begin{cor} Every prime ideal of $R$ has finite height. 
\end{cor}

\begin{cor}\label{cor:GPIT, 1}
If $R$ is local with maximal ideal $\m$, then $\dim(R)\leq \mu(\m)$.
\end{cor}

\noindent
Here is a generalization of Corollary~\ref{cor:PIT}: 

\begin{cor}\label{cor:GPIT, 2} If $m\ge 1$, the sequence $(r_1,\dots,r_m)\in R^m$ is 
regular on $R$, and~$\mathfrak p$ is a minimal prime divisor of
$Rr_1+\cdots+Rr_m$, then $\operatorname{ht}(\mathfrak p)=m$.
\end{cor}

\noindent
This follows from Theorem~\ref{thm:GPIT} and an easy induction on $m$, using:

\begin{lemma}
Let $\mathfrak p\in\Spec(R)$, suppose $x\in \mathfrak p$ is not a zero divisor, and set
$\overline{R}:=R/Rx$, $\overline{\mathfrak p}:=\mathfrak p/Rx\in\Spec(\overline{R})$.
Then $\operatorname{ht}(\mathfrak p) \geq \operatorname{ht}(\overline{\mathfrak p})+1$.
\end{lemma}
\begin{proof}
For $n=\operatorname{ht}(\overline{\mathfrak p})$ we get a strictly increasing sequence
 $\mathfrak p_0\subset \mathfrak p_1 \subset\cdots \subset \mathfrak p_n=\mathfrak{p}$ of prime ideals of $R$ with $Rx\subseteq \mathfrak p_0$.
Corollary~\ref{cor:irred comp 2} and $x\in\mathfrak p_0$ give that
 $\mathfrak p_0$ is not a minimal prime ideal of $R$. Thus
$\operatorname{ht}(\mathfrak p)\ge n+1$. 
\end{proof}

\noindent
The Generalized Principal Ideal Theorem has a converse:

\begin{prop}\label{prop:GPIT converse}
Let $\mathfrak p$ be a prime ideal of $R$ and $n=\operatorname{ht}(\mathfrak p)$. Then there are $x_1,\dots,x_n\in\mathfrak p$
such that $\mathfrak p$ is a minimal prime divisor of
$Rx_1+\cdots+Rx_n$.
\end{prop}
\begin{proof}
By induction on $n$. The case $n=0$ being trivial, suppose $n\geq 1$.
Let $\mathfrak q_1,\dots,\mathfrak q_m$ be the minimal prime ideals of $R$;
then $\mathfrak p\not\subseteq\mathfrak q_j$ for $j=1,\dots,m$.  
By Lemma~\ref{lem:prime avoidance} we get $x_1\in\mathfrak p$ with
$x_1\notin\bigcup_{j=1}^m \mathfrak q_j$.
Put $\overline{R}:=R/Rx_1$ and $\overline{\mathfrak p}:=\mathfrak p/Rx_1$.
Then $\operatorname{ht}(\overline{\mathfrak p})\leq n-1$ by the choice of $x_1$. So inductively we have $x_2,\dots,x_n\in R$ with
$\overline{\mathfrak p}$ a minimal prime divisor of the ideal of $\overline{R}$ generated by
the cosets $x_j+Rx_1$, $j=2,\dots,n$. Thus $\mathfrak p$ is a minimal prime 
divisor of $Rx_1+\cdots+Rx_n$.
\end{proof}

\begin{cor}\label{cor:GPIT, converse}
Let $\mathfrak p$ be a prime ideal of $R$ of height $n$. Let $x\in\mathfrak p$ and
$\overline{R}:=R/Rx$, $\overline{\mathfrak p}:=\mathfrak p/Rx$. Then $\operatorname{ht}(\overline{\mathfrak p})=n$ or $\operatorname{ht}(\overline{\mathfrak p})=n-1$.
\end{cor}
\begin{proof}
By Lemma~\ref{lem:ht under morphisms} it remains to show that
$\operatorname{ht}(\overline{\mathfrak p})\geq n-1$.
Suppose that $\operatorname{ht}(\overline{\mathfrak p})\leq n-2$.
Then $n\ge 2$, and Proposition~\ref{prop:GPIT converse} gives $x_1,\dots,x_{n-2}\in\mathfrak p$
such that $\overline{\mathfrak p}$ is a minimal prime divisor of 
the ideal of $\overline{R}$ generated by the cosets $x_j+Rx$, $j=1,\dots,n-2$.
Then $\mathfrak p$ is a minimal prime divisor of $Rx+Rx_1+\cdots +Rx_{n-2}$, 
hence 
$\operatorname{ht}(\mathfrak p)\leq n-1$ by Theorem~\ref{thm:GPIT}, 
a contradiction.
\end{proof}

\subsection*{Notes and comments}
The notion of Krull dimension and the Principal Ideal Theorem are due to 
Krull~\cite{Krull28, Krull29b}. It seems that Kronecker~\cite[p.~80]{Kronecker} already knew the Generalized Principal Ideal Theorem for polynomial rings over a field.
Not every noetherian ring has finite Krull dimension; see~\cite[Appendix~A1]{Nagata62}.

\section{Regular Local Rings}\label{sec:regular local rings}

\noindent
{\em In this section $(R,\m)$ is a noetherian local ring with residue field $\k=R/\m$}.

\subsection*{Definition and basic properties}
The {\bf embedding dimension} of $R$ is defined as $\operatorname{edim}(R):=\mu(\m)$, the minimal number of generators of the maximal ideal $\m$ of~$R$. 
Corollary~\ref{cor:NAK, 3} gives  $\operatorname{edim}(R)=\dim_{\k} \m/\m^2$,
and Corollary~\ref{cor:GPIT, 1} says that $\dim(R)\leq\operatorname{edim}(R)$. We
call $(R,\m)$
a {\bf regular} local ring if $\dim(R)=\operatorname{edim}(R)$.
If $\dim R=0$, then $R$ is regular iff $R$ is a field. 

\index{dimension!embedding}
\index{embedding!dimension}
\index{local ring!regular}
\index{regular!local ring}
\index{ring!regular local}

\begin{prop}\label{prop:reg => domain}
Every regular local ring is an integral domain.
\end{prop}

\noindent
In the proof of this proposition we use:

\begin{lemma}\label{lem:reg => domain}
Suppose $(R,\m)$ is a regular local ring, and $x\in\m\setminus\m^2$. Then $\overline{R}:=R/Rx$
is a regular local ring with $\dim R = 1+\dim\overline{R}$.
\end{lemma}
\begin{proof}
%We let $\overline{\m}:=\m/Rx$ denote the maximal ideal of $\overline{R}$.
We have 
$\dim R =  \operatorname{edim} R = 1+ \operatorname{edim} \overline{R}
\ge 1+\dim  \overline{R}$, using Corollary~\ref{cor:NAK, 4} for the second equality. 
Corollary~\ref{cor:GPIT, converse} gives
$\dim R \le 1+ \dim\overline{R}$.
\end{proof}

\begin{proof}[Proof of Proposition~\ref{prop:reg => domain}] 
Let $(R,\m)$ be a regular local ring;
we show by induction on the dimension $n$ of $R$ that $R$ is a domain.
This is clear if $\dim R=0$, so assume~$n\geq 1$. 
Let $\mathfrak p_1,\dots,\mathfrak p_m$ be the minimal prime ideals of $R$.
If $\mathfrak m\setminus\mathfrak m^2\subseteq \bigcup_{i=1}^m \mathfrak p_i$, then Corollary~\ref{cor:prime avoidance}
gives $\mathfrak m\subseteq\mathfrak p_i$ for some $i\in\{1,\dots,m\}$, which contradicts $\operatorname{ht}(\mathfrak m)=n\geq 1$.
Take $x\in\m\setminus\m^2$ with $x\notin \bigcup_{i=1}^m\mathfrak p_i$.
By Lemma~\ref{lem:reg => domain}, 
$\overline{R}:=R/Rx$ is a regular local ring of dimension $n-1$, hence an integral domain by inductive hypothesis. So $Rx$ is a prime ideal of $R$, and we can take $i\in\{1,\dots,m\}$ with~$\mathfrak p_i\subseteq Rx$. We claim that $\mathfrak p_i=\{0\}$ (and hence $R$ is an integral domain).
To see this, let $a\in\mathfrak p_i$, and take $r\in R$ such that $a=rx$.
Since $x\notin\mathfrak p_i$, we have $r\in\mathfrak p_i$, so $a=rx\in \mathfrak p_i\m$.
This yields $\mathfrak p_i = \m\mathfrak p_i$, hence $\mathfrak p_i=\{0\}$ by Nakayama's Lemma.
\end{proof}

\begin{cor}\label{cor:reg => domain}
Suppose $(R,\m)$ is a regular local ring. Let $I\neq R$ be an ideal of~$R$. Then the following are equivalent:
\begin{enumerate}
\item[\textup{(i)}] there is a minimal set $G$ of generators of $\m$ with distinct
$x_1,\dots,x_m\in G$ such that $I=Rx_1+\cdots +Rx_m$;
\item[\textup{(ii)}] $R/I$ is a regular local ring $($and thus $I$ is a prime ideal of $R)$. 
\end{enumerate}
In that case and with $m$ as in \rm{(i)}, we have $\dim R = m+ \dim R/I$.  
\end{cor}
\begin{proof} Let $m$ be as in (i). Then
induction on $m$ and Lemma~\ref{lem:reg => domain} gives 
(i)~$\Rightarrow$~(ii) and $\dim R = m + \dim R/I$.
To show (ii)~$\Rightarrow$~(i),
suppose $R/I$ is a regular local ring. Let $x_1,\dots,x_d\in\m$, where $d:=\dim R/I=\operatorname{edim} R/I$, be such that
$x_1+I,\dots,x_d+I$ are the distinct elements of a minimal set of 
generators of the maximal ideal $\m/I$ of~$R/I$. The natural exact
sequence of $\k$-linear maps
$$0\longrightarrow I/I\cap\m^2 \longrightarrow \m/\m^2 \longrightarrow (\m/I)/(\m/I)^2 \longrightarrow 0$$
and Corollary~\ref{cor:NAK, 3} yield $x_{d+1},\dots,x_n\in I$  such that
$x_1,\dots,x_n$ are the distinct elements of a minimal set of generators of $\m$.
Put $I':=Rx_{d+1}+\cdots+Rx_n$. By what we showed before, $R':=R/I'$ is a regular local ring
of dimension $d$ and $I'$ is a prime ideal of $R$. Since $I'\subseteq I$, this yields $I'=I$.
\end{proof}

\noindent
Next, a useful characterization of regularity:

\begin{prop}\label{prop:char of regularity} Assume $R$ is not a field. Then
the local ring $(R,\m)$ is regular if and only if there is an $n\ge 1$ and a
sequence $\vec r= (r_1,\dots, r_n)\in R^n$ such that~$\vec r$ is regular on $R$ and
$\m=Rr_1+ \cdots +Rr_n$. Moreover, if $r_1,\dots, r_n$ are the distinct elements
of a minimal set of generators of $\m$ and $\dim R=n$, then the sequence
$(r_1,\dots,r_n)$ is regular on $R$.
\end{prop}
\begin{proof}
Suppose $n\ge 1$, the sequence 
$(r_1,\dots,r_n)\in R^n$ is regular on~$R$, and $\m=Rr_1+\cdots+Rr_n$. Then by Corollaries~\ref{cor:GPIT, 1} and \ref{cor:GPIT, 2}
$$\dim(R)\ =\ \operatorname{ht}(\m)\ =\ n\ \geq\ \operatorname{edim}(R)\ \geq\ \dim(R),$$
so $\dim(R)=\operatorname{edim}(R)$.
Conversely, suppose that $(R,\m)$ is a regular local ring and $r_1,\dots, r_n$ are the distinct elements
of a minimal set of generators of $\m$. We show by induction on~$n$ that 
then the sequence $(r_1,\dots,r_n)$ is regular on $R$. Since $R$ is an 
integral domain, $r_1$ is not a zero divisor on $R$. Also,
$r_1\in\m\setminus\m^2$, so $\overline{R}:=R/Rr_1$ is regular of dimension $n-1$, by Lemma~\ref{lem:reg => domain}. We are done if $n=1$, 
so let~${n\ge 2}$. 
Put $\overline{r}:=r+Rr_1$ for $r\in R$. Then $\overline{r_2},\dots,\overline{r_n}$ are the distinct elements of a minimal set of generators of 
the maximal ideal $\m/Rr_1$ of $\overline{R}$. Hence the sequence $(\overline{r_2},\dots,\overline{r_n})$ is regular on $\overline{R}$, by 
inductive hypothesis. Thus $(r_1,\dots,r_n)$ is regular on~$R$, by~\eqref{eq:split regular sequ}.
\end{proof}

\subsection*{Examples of regular local rings}
First the $1$-dimensional case: 
 
\begin{lemma}\label{regdim1dvr} 
$R$ is regular and $\dim(R)=1\ \Longleftrightarrow\ R$ is a $\rm{DVR}$.
\end{lemma}
\begin{proof} The direction $\Leftarrow$ is clear. 
Assume $R$ is regular (hence an integral domain) and $\dim(R)=1$. Then $\m=tR$ with $0\ne t\in R$, so $R$ is a DVR
by Lemma~\ref{NAKDVR}.
%Then $\m^n\neq\m^{n+1}$ for all $n$ by 
%Corollary~\ref{cor:NAK, 5}.
%Moreover, if $I$ is an ideal of $R$, $I\neq\{0\}$ 
%and $I\neq R$, then $\sqrt{I}=\fm$ by 
%Proposition~\ref{prop:radical is intersection of primes} 
%(since $\m$ is the only non-zero prime ideal of $R$), 
%hence $\m^n\subseteq I$ for some $n$.
%This yields $\bigcap_n \fm^n=\{0\}$. Therefore,
%taking $t\in R$ with $\m=Rt$, every $r\in R^{\neq}$ equals 
%$ut^n$ for a unique pair $(u,n)$, $u\in R^\times$. 
%Thus $R$ is a DVR.
\end{proof}

\begin{prop}\label{prop:loc of poly rings are regular local}
Let $R$ be a regular local ring, $X$ an indeterminate over $R$ and~$\mathfrak P$ a prime ideal of $R[X]$ with
$\mathfrak P\cap R = \m$. Then $R[X]_{\mathfrak P}$ is a regular local ring.
\end{prop}
\begin{proof} Set $n:=\dim R$. If $n=0$, then $R$ is a field, and so 
$R[X]_{\mathfrak P}$ is either a field (iff
$\mathfrak P = \{0\}$) or a DVR. Assume $n\ge 1$.
Then Proposition~\ref{prop:char of regularity} gives a sequence 
$(r_1,\dots,r_n)\in R^n$ that is regular on $R$ with
$\m=Rr_1+\cdots + Rr_n$.
Then the desired result follows from Corollary~\ref{cor:reg poly}, 
Lemma~\ref{lem:reglocal}, and Proposition~\ref{prop:char of regularity}. 
\end{proof}

\noindent
The following consequence will be
used in Section~\ref{sec:johnson}:

\begin{cor}\label{cor:loc of poly rings are regular local}
Let $K$ be a $\rm{PID}$, 
$X_1,\dots,X_n$ distinct indeterminates,
and $\mathfrak P$ a prime ideal of $K[X_1,\dots,X_n]$. Then $K[X_1,\dots,X_n]_{\mathfrak P}$ is a regular local ring.
\end{cor}
\begin{proof}
By induction on $n$. The case $n=0$ reduces to Lemma~\ref{regdim1dvr}. Assume $n\geq 1$, put $A:=K[X_1,\dots,X_{n-1}]$
and $\mathfrak p:=A\cap\mathfrak P$. Inductively, $A_{\mathfrak p}$ is a regular local ring. Also, $K[X_1,\dots,X_n]_{\mathfrak P}= A_{\mathfrak p}[X_n]_{\mathfrak Q}$ for some $\mathfrak Q\in\Spec\!\big(A_{\mathfrak p}[X_n]\big)$, by Corollary~\ref{cor:loc of polynomial ring}. Hence $K[X_1,\dots,X_n]_{\mathfrak P}$ is a regular local ring by Proposition~\ref{prop:loc of poly rings are regular local}.
\end{proof}

\subsection*{Notes and comments} 
Regular local rings were studied by Krull~\cite{Krull38}, which contains Proposition~\ref{prop:reg => domain}. The terminology comes from~\cite{Chevalley43}. The 
algebraic-geometric significance of regular local rings is largely due to Zariski~\cite{Zariski47}.

\section{Modules and Derivations}\label{sec:modules}

\noindent
The rest of this chapter is dominated by the general idea of linearization, and by the useful dualities that come with it. Thus we consider tensor products of modules, derivations, and (dually) differentials. 
% to prepare for  Sections~\ref{sec:differentials} and 
%\ref{sec:johnson}. 
%below, in this section
%we recall some basic facts about modules.
{\em Throughout this section we fix a commutative ring $R$ and let $M$, $N$ range over $R$-modules}.

\subsection*{Modules of morphisms}
Let $A$, $B$, $C$ be $R$-modules. 
We denote by $\Hom_R(A,M)$ (or $\Hom(A,M)$ if $R$ is clear from the context) the set of $R$-linear maps $A\to M$, made into an $R$-module in the natural way.
Let $\phi\colon A\to B$ be $R$-linear. Then 
$$\beta\mapsto \beta\circ \phi\colon\Hom(B,M)\to \Hom(A,M)$$
is an $R$-linear map, denoted by $\phi^*$ or $\Hom(\phi,M)$. If $\psi\colon B\to C$ is also $R$-linear,
then $(\psi\circ \phi)^*=\phi^*\circ\psi^*$.

\nomenclature[A]{$\Hom_R(A,M)$}{module of $R$-linear maps $A\to M$}

\begin{lemma}\label{lem:Hom and exact sequ}
The following are equivalent for $\phi$ and $R$-linear $\psi\colon B \to C$:
\begin{enumerate}
\item[\textup{(i)}] the sequence 
$ A  \xrightarrow{\ \phi\ }  B \xrightarrow{\ \psi\ }  C \longrightarrow  0$
is  exact;
\item[\textup{(ii)}] for each $M$ the sequence
$$\Hom(A,M)  \xleftarrow{\ \phi^*\ }  \Hom(B,M) \xleftarrow{\ \psi^*\ }  \Hom(C,M) \longleftarrow 0$$
of induced morphisms is exact.
\end{enumerate}
\end{lemma}
\begin{proof}
The direction (i)~$\Rightarrow$~(ii) is entirely routine.
Assume (ii). Exactness at~$C$ in~(i) means surjectivity of $\psi$.
Let $\alpha\colon C\to C/\operatorname{im}\psi$ be the canonical map. Since
$\alpha\circ \psi=0$, exactness in (ii) at $\Hom(C,M)$
for $M=C/\operatorname{im}\psi$ gives $\alpha=0$, so~$\operatorname{im}\psi=C$. Next, applying
$(\psi\circ \phi)^*=0$ to 
$\operatorname{id}_{C}$ gives $\psi\circ \phi=0$. It only remains to
show that $\operatorname{im}\phi\supseteq \ker \psi$. Let $\beta\colon B \to M:=B/\operatorname{im}\phi$ be the canonical map.
Then $\beta\circ \phi=0$, so $\beta=\gamma\circ \psi$ for some
$\gamma\in \Hom(C, M)$, so $\ker \psi\subseteq \ker \beta= \operatorname{im}\phi$.
\end{proof}

\subsection*{Tensor products} The tensor product $M\otimes_R N$ of the $R$-modules $M$ and $N$ is an $R$-module
with an $R$-bilinear map 
$$(x,y) \mapsto x\otimes_R y\ \colon\ M\times N \to M\otimes_R N$$
that is universal: for any $R$-module $B$ and $R$-bilinear map $\beta\colon M\times N \to B$ there is a unique $R$-linear $b\colon M\otimes_R N \to B$ such that
$b(x\otimes y)=\beta(x,y)$ for all $x\in M$, $y\in N$. 
 For use below we recall the construction of $M\otimes N$: Let $x$,~$x'$ range over $M$, $y$,~$y'$ over $N$, and $r$ over $R$. Let $F$ be the free $R$-module with basis $M\times N$, and let~$G$ be the submodule of $F$ generated by all elements of the form
\begin{align*} (x+x',y)-(x,y)-(x',y), \qquad &(x,y+y')-(x,y)-(x,y'),\\
     (rx,y)-r(x,y),\qquad &(x,ry)-r(x,y).
\end{align*}
Put $M\otimes_R N:=F/G$; then 
$$(x,y) \mapsto x\otimes y:=(x,y)+G\ \colon\ M\times N\to M\otimes N$$ is an $R$-bilinear map having the desired universal property.
%See for example \cite[Chapter~XVI]{Lang} or 
%\cite[Appendix~A]{Matsumura}. 
We drop the subscript~$R$
in expressions like $M\otimes_R N$ and 
$x\otimes_R y$ when~$R$ is clear from the context.
We shall need the following variant of the universal property of $M\otimes N=M\otimes_R N$:

\begin{lemma}\label{univtensadd}
Let $\phi\colon M\times N\to Z$ be a biadditive map into an abelian group $Z$ such that
$\phi(rx,y)=\phi(x,ry)$ for all $r\in R$, $x\in M$, $y\in N$. Then there is a unique
additive map $f\colon M\otimes N\to Z$ such that $f(x\otimes y)=\phi(x,y)$ for all $x\in M$, $y\in N$.
\end{lemma}
\begin{proof} With the notations in the above construction
of $M\otimes N$,  
%Let $x$,~$x'$ range over $M$, $y$,~$y'$ 
%over $N$, and $r$ over $R$.
%Recall the construction of $M\otimes N$: Let $F$ be the free %$R$-module with basis $M\times N$, and let $G$ be the %submodule of $F$ generated by all elements of the form
%\begin{align*} (x,y+y')-(x,y)-(x,y'), 
%\qquad &(x,y+y')-(x,y)-(x,y'),\\
%     (rx,y)-r(x,y),\qquad &(x,ry)-r(x,y).
%\end{align*}
%Put $M\otimes N:=F/G$; then 
%$$(x,y) \mapsto x\otimes y:=(x,y)+G\ \colon\ 
%M\times N\to M\otimes N$$ is an $R$-bilinear map having 
%the desired universal property. Now let 
%$\phi\colon M\times N\to Z$ be 
% a biadditive map into an abelian group $Z$ such that
%$\phi(rx,y)=\phi(x,ry)$ for all $r$,~$x$,~$y$. 
we have the internal direct sum decomposition
$F=\bigoplus_{(x,y)} R\cdot(x,y)$, and 
for given $x$, $y$, the map
$r\mapsto \phi(rx,y)\colon R\to Z$ is
additive. Thus we have an additive map $\hat\phi\colon F\to Z$ 
with $\hat\phi\big(r(x,y)\big)=\phi(rx,y)$ for all $r$, $x$, $y$.
Clearly $G\subseteq\ker\hat\phi$, so $\hat\phi$ induces an additive map $f\colon M\otimes N=F/G\to Z$
such that $f(x\otimes y)=\phi(x,y)$ for all $x$, $y$. 
There can be at most one such map $f$, since every element 
of $M\otimes N$ is a finite sum of elements of the form 
$x\otimes y$.
\end{proof}

\noindent
We use this variant to construct derivations. 
A {\bf derivation\/} on $R$ is a map $\der \colon R\to R$
such that $\der(a+b)=\der(a)+\der(b)$ and 
$\der(ab)=\der(a)b+a\der(b)$ for all $a,b\in R$. The map $R \to R$ sending 
every element of $R$ to $0$ is a derivation on $R$, called the 
{\bf trivial\/} derivation (on $R$). If $\der$ is a derivation on $R$ and 
$a\in R$, then $a\der$ is also a derivation on $R$.
Below $\der$ is a derivation on $R$. It gives a subring
$\{a\in R:\ \der(a)=0\}$ of $R$. If a set~${X\subseteq R}$
generates a subring $R_0$ of $R$ with $\der(X)\subseteq R_0$,
then $\der(R_0)\subseteq R_0$. A {\bf $\der$-compatible derivation} on   $M$
is an additive map $d\colon M\to M$ such that $d(rx)=\der(r)x+rd(x)$ for all $r\in R$ and $x\in M$.

\index{derivation}
\index{derivation!trivial}
\nomenclature[M]{$\der$}{derivation on a ring}
\index{module!$\der$-compatible derivation}
\index{derivation!$\der$-compatible}

\begin{cor}\label{cor:der on tensor prod, 1}
Let $d_M$ and $d_N$ be $\der$-compatible derivations on~$M$ and $N$.
Then there is a unique $\der$-compatible derivation $\derdelta$ on $M\otimes N$
such that
\begin{equation}\label{eq:der on tensor prod}
\derdelta(x\otimes y)\ =\ d_M(x) \otimes y + x\otimes d_N(y)\quad\text{for all $x\in M$, $y\in N$.}
\end{equation}
\end{cor}
\begin{proof}
The map
$$(x,y)\mapsto \phi(x,y):=d_M(x) \otimes y + x\otimes d_N(y)\ \colon\ M\times N\to M\otimes N$$
is biadditive and satisfies $\phi(rx,y)=\phi(x,ry)$ for all $r\in R$, $x\in M$, $y\in N$.
Hence Lemma~\ref{univtensadd} yields an additive map $\derdelta\colon M\otimes N\to M\otimes N$
satisfying \eqref{eq:der on tensor prod}, and this map is easily checked to be a $\der$-compatible derivation on
$M\otimes N$. 
\end{proof}

\noindent
Now let 
$A$ be a commutative $R$-algebra. Then $A\otimes_R M$ is not only an $R$-module, but even an  
$A$-module inducing the given $R$-module structure, with 
$$a\cdot (b\otimes x)\ =\ ab\otimes x \qquad (a,b\in A,\ x\in M).$$
Let $(x_i)$ be a family of elements in $M$. If the $R$-module $M$ is free on 
$(x_i)$, then the $A$-module $A\otimes_R M$ is free on $(1\otimes x_i)$; the converse holds if the structural morphism $R \to A$ is injective.

\begin{cor}\label{cor:der on tensor prod, 2}
Let $\der_A$ be a derivation on $A$ such that $\der_A(r1_A)=\der(r)1_A$ for all~$r\in R$. Then $\der_A$ is a $\der$-compatible derivation on the $R$-module $A$, and if $d$ is a $\der$-compatible derivation on $M$, then the $\der$-compatible derivation $\derdelta$ on $A\otimes_R M$
with
$$\derdelta(a\otimes x)\ =\ \der_A(a) \otimes x + a\otimes d(x)\quad\text{for all $a\in A$, $x\in M$}$$
is $\der_A$-compatible.
\end{cor}

\noindent
Let $B$ be a second commutative $R$-algebra. Then
the $A$-module $A\otimes_R B$ is even a (commutative) $A$-algebra with
multiplicative identity $1_A\otimes 1_B$ and multiplication given by
$$(a_1\otimes b_1)\cdot (a_2\otimes b_2)\ =\ a_1a_2\otimes b_1b_2 \qquad (a_1,a_2\in A,\ b_1,b_2\in B).$$
We have ring morphisms $a\mapsto a\otimes 1_B\colon A\to A\otimes B$ and
$b\mapsto 1_A\otimes b\colon B\to A\otimes B$, which are injective if $R$ is a field and $1_A\ne 0$, $1_B\ne 0$.

\begin{cor}\label{cor:der on tensor prod, 3}
Let $\der_A$ and $\der_B$ be derivations on the rings $A$
and $B$, respectively, such that $\der_A(r1_A)=\der(r)1_A$ and $\der_B(r1_B)=\der(r)1_B$ for all $r\in R$.
Then there is a unique derivation $\derdelta$ on the ring $A\otimes_R B$ such that
$$\derdelta(a\otimes b)\ =\ \der_A(a)\otimes b + a\otimes\der_B(b)\quad\text{for all $a\in A$, $b\in B$.}$$
\end{cor}

\subsection*{Some useful isomorphisms}
Let $A$ be a commutative $R$-algebra. Let $E$ be an $A$-module, viewed as an
$R$-module via the structural morphism $R \to A$. For an $R$-linear map
$\phi\colon M\to E$ and $a\in A$ we define the $R$-linear map $a\phi\colon M \to E$
by $a\phi(x):= a\cdot \phi(x)$~($x\in M$). This makes  
$\Hom_R(M,E)$ an $A$-module inducing its given $R$-module structure. 
For any $R$-linear $\phi\colon M \to E$ we have the $A$-linear map
$$\phi_A\ \colon\ A\otimes_R M\to E, \qquad \phi_A(a\otimes x)=a\phi(x)\quad(a\in A,\ x\in M),$$
and it is routine to check that the resulting map 
\begin{equation}\label{eq:iso tensor, 1}
\phi\mapsto \phi_A\ \colon\ \Hom_R(M,E)\to \Hom_A(A\otimes_R M,E)
\end{equation}
is an isomorphism of $A$-modules.

\medskip
\noindent
Let $I$ be an ideal of $R$ and consider $A:=R/I$ as an $R$-algebra in 
the obvious way. Then $M/IM$ is naturally an $A$-module, and we have an $A$-linear map
\begin{equation}\label{eq:iso tensor, 3}
x+IM\mapsto 1_A\otimes x\ \colon\ M/IM\to  A\otimes_R M  \qquad(x\in M).
\end{equation}
This map \eqref{eq:iso tensor, 3} is an isomorphism of $A$-modules: its inverse
is the $A$-linear map $A\otimes_R M\to M/IM$ induced by the $R$-bilinear map
$A\times M\to M/IM$ that sends $(r+I,x)$ to $rx+IM$ for $r\in R$, $x\in M$.

\medskip\noindent
Let $S$ be a multiplicative subset of $R$, and consider $A:=S^{-1} R$ as an $R$-algebra via the canonical map
$\iota\colon R \to S^{-1}R=A$. Then we have an isomorphism
\begin{equation}\label{eq:iso tensor, 4}
A \otimes_R M\ \to\ S^{-1}M
\end{equation} of $A$-modules that
sends $(r/s)\otimes x$ to $(rx)/s$ for $r\in R$, $s\in S$, $x\in M$.

\medskip
\noindent
Let $h\colon A \to B$ be a morphism between commutative $R$-algebras $A$ and $B$.
We have the $A$-module $A\otimes_R M$, and so the $B$-module $B\otimes_A (A\otimes_R M)$, by construing $B$ as an $A$-algebra via $h$. We also have a $B$-linear
map
\begin{equation}\label{eq:iso tensor, 2}
B\otimes_R M\ \to\  B\otimes_A (A\otimes_R M), \qquad b\otimes_R x\mapsto b\otimes_A(1\otimes_R x) \ \text{ for $x\in M$.} 
\end{equation}
It is easy to check (by explicit construction of an inverse) that this
map is in fact an isomorphism of $B$-modules: {\em transitivity of base change}.

\subsection*{Rational rank} Let $M$ be an abelian group, in other words, a $\Z$-module. Then we define the {\bf rational rank} $\operatorname{rank}_{\Q} M$ of $M$ by 
$$\operatorname{rank}_{\Q} M\ :=\ \dim_{\Q} \Q\otimes_{\Z} M,$$
if the $\Q$-linear space $\Q\otimes_{\Z} M$ has finite dimension, and set $\operatorname{rank}_{\Q} M:= \infty$ otherwise. Using for example \eqref{eq:iso tensor, 4} it follows that this rational rank is the largest $m$ for which there are $\Z$-independent
$x_1,\dots, x_m\in M$, if such a largest $m$ exists, and is $\infty$ otherwise. In particular, the free abelian group
$\Z^m$ has rational rank $m$. It is easy to check that the
rational rank, when restricted to abelian groups of finite rational rank,
yields an Euler-Poincar\'e map with values in $\Z$.  

\nomenclature[A]{$\operatorname{rank}_\Q(M)$}{rational rank of the abelian group $M$}
\index{abelian group!rational rank}
\index{rank!rational}
\index{rational!rank}

\subsection*{Independence at a prime}
Let $f_i\in M$ for $i\in I$. Given $\frak{p}\in \Spec(R)$, the family~$(f_i)$ is said to be {\bf independent at~$\frak{p}$} 
if the family $(f_i+\frak{p}M)$ in the $R/\frak{p}$-module~${M/\frak{p}M}$ is linearly independent over~$R/\frak{p}$. The fact below is used in Section~\ref{sec:johnson}:

\index{independent!at a prime}

\begin{lemma}\label{lem:independence under generalization} Suppose the 
$R$-module $M$ is free and $\frak{p}, \frak{q}\in \Spec(R)$, 
$\frak{p}\supseteq \frak{q}$.
If~$(f_i)$ is independent
at $\frak{p}$, then
$(f_i)$ is independent at $\frak{q}$.
\end{lemma}
\begin{proof} Let $f_1,\dots, f_m\in M$ with $m\ge 1$ be such that
$f_1+\frak{p}M,\dots, f_m+\frak{p}M$ are linearly independent over
$R/\frak{p}$. It suffices to show that then $f_1+\frak{q}M,\dots, f_m + \frak{q}M$ are linearly independent over $R/\frak{q}$. Take a basis of $M$ and
take distinct $e_1,\dots, e_n$ from that basis such that
for $i=1,\dots,m$ we have $f_i= f_{i1}e_1+\cdots +f_{in}e_n$ with all
$f_{ij}\in R$. The vectors $e_1+\frak{p}M,\dots, e_n+\frak{p}M$ of $M/\frak{p}M$ are linearly independent over~$R/\frak{p}$, so $m\le n$, and some $m\times m$ submatrix
of $(f_{ij})$ has determinant $D\notin \frak{p}$. Then 
$D\notin\frak{q}$, so $f_1+\frak{q}M,\dots, f_m + \frak{q}M$ 
are linearly independent over $R/\frak{q}$.   
\end{proof} 

\subsection*{Notes and comments} 
%For Corollaries~\ref{cor:der on tensor prod, 1}--\ref{cor:der on tensor prod, 3} see, e.g., \cite[p.~49]{CF}.
Lem\-ma~\ref{lem:independence under generalization} is taken from \cite{JJ}.

\section{Differentials}\label{sec:differentials}

\noindent
Below we define K\"ahler differentials and prove some basic facts
about them, as needed in Section~\ref{sec:johnson}.
{\em Throughout, $K$ is a commutative ring, $A$ is a commutative $K$-algebra, and
$M$ is an $A$-module}. (Later we only use the case where $K$ is a field of characteristic zero, but that doesn't simplify things.)
We let $a$, $b$ range over $A$.

\subsection*{$K$-derivations}
A map $\Delta\colon A\to M$ is said to be a {\bf $K$-derivation} if 
$\Delta$ is $K$-linear and $\Delta(ab)=a\Delta(b)+b\Delta(a)$ for all $a$, $b$. 
A $K$-derivation $A\to A$ is called a {\bf $K$-derivation on $A$.}
So every $K$-derivation on $A$ is a derivation on $A$ in the sense of Section~\ref{sec:modules},
and a derivation on a commutative ring $R$ is the same thing as a $\Z$-derivation on $R$, where $R$ is viewed as a $\Z$-algebra.

\index{K-derivation@$K$-derivation}
\index{derivation!$K$-derivation}
\index{module!$K$-derivation}

In the rest of this subsection $\Delta\colon A\to M$ is a $K$-derivation, so
$\Delta(\kappa\cdot 1)=0$ for all~$\kappa\in K$, and $A^\Delta:=\ker\Delta$ is a subalgebra of $A$. The rules below for computing with~$\Delta$ are easy
to  verify:

\begin{lemma}
For all $n\ge 1$ and $a_1,\dots,a_n\in A$,
$$\Delta(a_1\cdots a_n)\ =\ \sum_{j=1}^n a_1\cdots a_{j-1}a_{j+1}\cdots a_n \, \Delta(a_j).$$
In particular, $\Delta(a^n)=na^{n-1}\Delta(a)$ for $n\geq 1$, and if
$I$ is an ideal of $A$ and $n\geq 1$, then $\Delta(I^n)\subseteq I^{n-1}M$.
\end{lemma}

\noindent
Let $X=(X_1,\dots, X_n)$ be a tuple of distinct indeterminates, $n\ge 1$. 
Consider a 
polynomial $P=\sum_{\i}P_{\i}X^{\i}\in A[X]$ with
$\i=(i_1,\dots, i_n)$ ranging over a finite subset of $\N^n$ and with
$P_{\i}\in A$ and $X^{\i}:= X_1^{i_1}\cdots X_n^{i_n}$ for all $\i$. 

\begin{lemma}\label{derpolA} For ${\mathbf a}=(a_1,\dots, a_n)\in A^n$, and setting 
${\mathbf a}^{\i}:= a_1^{i_1}\cdots a_n^{i_n}$,
$$\Delta\big(P({\mathbf a})\big) \ =\  \sum_{\i}{\mathbf a}^{\i}\Delta(P_{\i}) + \sum_{j=1}^n \frac{\partial P}{\partial X_j}({\mathbf a}) \Delta(a_j).$$
\end{lemma}

\noindent
In particular, for polynomials over $K$:

\begin{cor}\label{cor:deriv of poly expr}
Given $F\in K[X]$ and ${\mathbf a}=(a_1,\dots,a_n)\in A^n$, we have
$$\Delta\big(F({\mathbf a})\big)\ =\ \sum_{j=1}^n
  \frac{\partial F}{\partial X_j}\!({\mathbf a})\,\Delta(a_j).$$
\end{cor}

%\begin{cor}\label{cor:deriv of poly expr, 1}
%Suppose $A$ is a field, and let $P\in K[X]^{\neq}$ and $a$  satisfy $P(a)=0$, $P'(a)\neq 0$. Then $\Delta(a)=0$.
%\end{cor}

\noindent
If $A$ is a field, then $A^\Delta$ is a subfield of~$A$:

\begin{lemma}[Quotient Rule]\label{lem:quotient rule}
Let $s\in A^\times$. Then  
$$\Delta(s^{-1}a)\ =\ s^{-2}\big(s\Delta(a)-a\Delta(s)\big),$$
in particular, $\Delta(s^{-1})=-s^{-2}\Delta(s)$.
\end{lemma}

\begin{lemma}\label{lem:quotient rule, 1} Let $S$ be a multiplicative subset of 
$A$; so $S^{-1}M$ is an $S^{-1}A$-module. 
There is exactly one $K$-derivation $S^{-1}\Delta\colon S^{-1}A\to S^{-1}M$ making the diagram
$$\xymatrix{ 
S^{-1}A \ar[r]^{S^{-1}\Delta} & S^{-1}M \\
A \ar[r]^\Delta \ar[u]^\iota & M \ar[u]^{\iota_M}  }$$
commutative; it is given by
\begin{equation}\label{eq:DeltaS}
S^{-1}\Delta(a/s)\ =\ \big(s\Delta(a)-a\Delta(s)\big)\big/s^2\quad (s\in S).
\end{equation}
\end{lemma}
\begin{proof} The main thing to check is that the value of the right-hand side
in \eqref{eq:DeltaS} does not change when replacing $a$,~$s$ by $at$,~$st$ for $t\in S$.
\end{proof}

\begin{cor}\label{cor:derivation into fraction field} If $A$ is an integral domain with fraction field $F$,
then any $K$-derivation
 $\der\colon A\to F$  extends uniquely to a $K$-derivation on $F$. % \marginpar{Need to replace all references to {\tt lem:derivation into fraction field} by references to this corollary.}
\end{cor}

\subsection*{The module of $K$-derivations}
The set $\Der_K(A,M)$ of $K$-derivations $A\to M$ is a submodule of the $K$-module $\Hom_K(A,M)$ of $K$-linear maps $A\to M$. Given $\Delta\in\Der_K(A,M)$, 
we also have $a\Delta\in\Der_K(A,M)$ with $a\Delta\colon A\to M$ defined by $(a\Delta)(b):=a\Delta(b)$. This scalar multiplication
$$(a,\Delta)\mapsto a\Delta\ \colon\ A\times\Der_K(A,M)\to \Der_K(A,M)$$
makes $\Der_K(A,M)$ an $A$-module. The $K$-module structure on 
$\Der_K(A,M)$ induced by this $A$-module structure is the already given one.  Set 
$\Der_K(A):=\Der_K(A,A)$. 

%\begin{example}
%Suppose $A$ and $K$ are fields with $A$ separable algebraic over $K$; then $\Der_K(A,M)=\{0\}$ for each $M$, by Corollary~\ref{cor:deriv of poly expr, 1}.
%\end{example}

\medskip
\noindent
Let $N$ also be an $A$-module, and let $\phi\colon M\to N$ be $A$-linear.
Then we have an $A$-linear map
$$\Der_K(A,\phi)\ \colon\ \Der_K(A,M)\to \Der_K(A,N), \qquad \Delta\mapsto \phi\circ \Delta.$$
Let $P$ be a third $A$-module, and let $\psi\colon N\to P$ be $A$-linear. Then 
$$\Der_K(A,\psi\circ\phi)\ =\ \Der_K(A,\psi)\circ\Der_K(A,\phi).$$

\begin{lemma}\label{lem:exact sequ 1}
Suppose $0\longrightarrow M\xrightarrow{\ \phi\ } N\xrightarrow{\ \psi\ } P$ is exact. Then
$$0\longrightarrow \Der_K(A,M)\xrightarrow{\ \Der_K(A,\phi)\ } \Der_K(A,N)\xrightarrow{\ \Der_K(A,\psi)\ } \Der_K(A,P)$$
is an exact sequence of $A$-modules and $A$-linear maps.
\end{lemma}
\begin{proof}
That $\Der_K(A,\phi)$ is injective and
maps $\Der_K(A,M)$ into $\ker \Der_K(A,\psi)$ is clear. If $\Delta\in\Der_K(A,N)$ and 
$\psi\circ\Delta=0$, then
$\Delta(A)\subseteq\ker\psi=\operatorname{im}\phi$ and hence $\Delta\in\Der_K(A,\phi)\big(\!\Der_K(A,M)\big)$.
\end{proof}

\noindent
Next, let $B$ be a $K$-algebra, $h\colon A\to B$ a $K$-algebra morphism, 
and $N$ a $B$-module.
If $\Delta\colon B\to N$ is a $K$-derivation, then $\Delta\circ h\colon A\to N$ is also a $K$-derivation,
where we view $N$ as an $A$-module via~$h$.
So we obtain an $A$-linear~map
$$ \Der_K(h,N)\ \colon\
\Der_K(B,N)\to \Der_K(A,N), \qquad \Delta\mapsto \Delta\circ h,$$
where the $B$-module $\Der_K(B,N)$ is viewed as an $A$-module via $h$. 
We note that $\Der_A(B,N)\subseteq \Der_K(B,N)$, not just as sets, but also as $A$-modules.

\begin{lemma}\label{lem:exact sequ 2}
The sequence
$$0\longrightarrow \Der_A(B,N)\longhookrightarrow \Der_K(B,N)\xrightarrow{\ \Der_K(h,N)\ } 
\Der_K(A,N)$$
of $A$-modules and $A$-linear maps is exact.
\end{lemma}
\begin{proof}
Clearly $\Der_A(B,N)\subseteq\ker \Der_K(h,N)$. 
If $\Delta\in\Der_K(B,N)$ and $\Delta\circ h=0$, then $\Delta$ is an 
$A$-derivation and
so $\Delta \in \Der_A(B,N)$.
\end{proof}

\noindent
Let $I$ be an ideal of $A$ with $IM=\{0\}$, set $B:= A/I$, and let
$h\colon A \to B$ be the canonical map. We construe the $A$-modules $M$ and $I/I^2$
as $B$-modules in the usual way. Note that if $\Delta\in \Der_K(A,M)$ and $\Delta(I)=\{0\}$, then $\Delta$ induces a $K$-derivation $h(a) \mapsto \Delta(a)\colon B \to M$.  
For 
$\Delta\in\Der_K(A,M)$ we have $\Delta(I^2)=\{0\}$, 
so the restriction of $\Delta$ to $I$ induces a $B$-linear map $\overline{\Delta}\colon I/I^2\to M$. This yields an $A$-linear~map
$$\Delta\mapsto \overline{\Delta}\ \colon\ \Der_K(A,M)\to \Hom_{B}(I/I^2,M)$$
where the $B$-module $\Hom_{B}(I/I^2,M)$ is viewed as an $A$-module via $h$. 

\begin{lemma}\label{lem:exact sequ 3}
The sequence
$$0\longrightarrow \Der_K(B,M)\xrightarrow{\ \Der_K(h,M)\ } \Der_K(A,M)\xrightarrow{\ \Delta\mapsto \overline{\Delta}\ } \Hom_{B}(I/I^2,M)$$
of $A$-modules and $A$-linear maps is exact.
\end{lemma}
\begin{proof}
From $\Der_A(B,M)=\{0\}$ and Lemma~\ref{lem:exact sequ 2} we get that
$\Der_K(h,M)$ is injective. If $\Delta\in \operatorname{im}\Der_K(h,M)$, then 
$\overline{\Delta}=0$.
Conversely, if $\Delta\in\Der_K(A,M)$ and $\overline{\Delta}=0$,
then $\Delta(I)=\{0\}$, so $\Delta$ induces a $K$-derivation $B=A/I\to M$, hence
$\Delta\in\operatorname{im}\Der_K(h,M)$. 
\end{proof}

\subsection*{The universal $K$-derivation}
 A {\bf universal} $K$-derivation of $A$ is an $A$-module~$\Omega$ together with a $K$-de\-ri\-va\-tion $\d\colon A \to \Omega$ such that for any $A$-module $M$ and $K$-derivation $\Delta\colon A \to M$ there
is a unique $A$-linear map $\phi\colon \Omega\to M$ with
$\Delta=\phi\circ \d$:
$$\xymatrix{& \Omega \ar@{-->}[dr]^{\phi} & \\
A \ar[ur]^{\d} \ar[rr]^{\Delta} &&  M}$$
If $\d_1\colon A\to\Omega_1$ and $\d_2\colon A\to\Omega_2$ are universal $K$-derivations of $A$, then the unique $A$-linear map $\phi\colon\Omega_1\to\Omega_2$
with $\d_2=\phi\circ \d_1$ is an isomorphism of $A$-modules, 
by the usual argument. So there is at most one universal $K$-derivation of~$A$ up to canonical isomorphism.
To get a universal $K$-derivation of $A$, take some bijection $d\colon A \to d(A)$ and let
$F$ be the free $A$-module on the set $d(A)$. Let $G$ be the $A$-submodule of $F$ generated by the elements 
$$d(ab)-ad(b)-bd(a),\quad d(a+b)-d(a)-d(b),\quad d(\kappa a) - \kappa d(a)  \ (\kappa\in K).$$
Let $\Omega$ be the quotient module
$F/G$, and let $\d\colon A \to \Omega$ be given by $\d(a):= d(a)+G$.
It is clear that then $\Omega$ with $\d\colon A \to \Omega$
is a universal $K$-derivation of~$A$. For definiteness, we set
$\Omega_{A|K}:= \Omega$ and let $\d\colon A \to \Omega_{A|K}$ be the universal $K$-derivation just defined. Note that the $A$-module
$\Omega_{A|K}$ is generated by its elements~$\d(a)$.  Instead of $\d(a)$, we also write $\d a$, or $\d_{A|K}(a)$ if we need to indicate the dependence on $A$ and $K$. 
The $A$-module $\Omega_{A|K}$ is called  the module of (K\"ahler) {\bf differentials} of the $K$-algebra~$A$,
and $\d a$ is called the {\bf differential} of $a$. 

\index{K-derivation@$K$-derivation!universal}
\index{universal $K$-derivation}
\nomenclature[A]{$\Omega_{A\vert{}K}$}{$A$-module of K\"ahler differentials of the $K$-algebra $A$}
\nomenclature[A]{$\d a$}{differential of $a$ (in $\Omega_{A\vert{}K}$)}
\index{differential}
\index{K\"ahler differentials}

It may help to think of the elements 
of $A$ as $K$-valued functions on a space, with $\Omega_{A|K}$ 
as a kind of cotangent bundle for this space. In any case, the 
$K$-derivation $\d\colon A \to \Omega_{A|K}$ can be useful in linearizing problems involving $A$.

\medskip\noindent
Note that if $A$ is generated as a $K$-algebra by its elements 
$a_i$ with $i\in I$, then by Corollary~\ref{cor:deriv of poly expr}
the $A$-module $\Omega_{A|K}$ is generated by the differentials $\d a_i$ of
these elements.
Using both Corollary~\ref{cor:deriv of poly expr} and 
the Quotient Rule~\ref{lem:quotient rule} we get:

\begin{cor}\label{cor: OKL} If $K$ is a field and $L=K(x_1,\dots,x_n)$
a field extension of $K$,
then the vector space $\Omega_{L|K}$ over $L$ is generated by $\d x_1,\dots,\d x_n$, so $\dim_L \Omega_{L|K}\leq n$.
\end{cor}

\medskip
\noindent
Each $A$-module $M$ gives an isomorphism of $A$-modules
\begin{equation}\label{eq:Hom Der}
\Hom_A(\Omega_{A|K},M)\ \xrightarrow{\ \cong\ }\ \Der_K(A,M),\qquad \phi\mapsto \phi\circ \d.
\end{equation}
In particular we have an isomorphism of $A$-modules
\begin{equation}\label{eq:Hom Der 2}
\Omega_{A|K}^*\ :=\ \Hom_A(\Omega_{A|K},A)\ \xrightarrow{\ \cong\ }\ \Der_K(A),\qquad \phi\mapsto \phi\circ \d,
\end{equation} 
so we have an $A$-bilinear map 
$$\<\,\ ,\ \>\ \colon\ \Omega_{A|K}\times \Der_K(A) \to A$$
given by $\<\omega,\Delta\>:=\phi(\omega)$ for $\omega\in\Omega_{A|K}$, $\Delta\in\Der_K(A)$, where
$\phi\colon\Omega_{A|K}\to A$ is the $A$-linear map such that $\phi\circ\d=\Delta$.
Thus for $\Delta\in\Der_K(A)$,
\begin{equation}\label{eq:Hom Der 3}
\<\d a,\Delta\>\ =\ \Delta(a).
\end{equation}

\begin{lemma}\label{lem:diff module of polynomial algebra}
Suppose $A=K[X]$ where $X=(X_i)_{i\in I}$ is a family of distinct
indeterminates over $K$. Then the $A$-module $\Omega_{A|K}$ is free with basis
 $(\d X_i)_{i\in I}$, and for $P\in A$
we have  $\frac{\partial P}{\partial X_i}=0$ for all but finitely many $i\in I$ and
\begin{equation}\label{eq:universal der, poly}
\d P\  =\ \sum_{i\in I} \frac{\partial P}{\partial X_i}\,\d X_i.
\end{equation}
\end{lemma}
\begin{proof} Let $i$,~$j$ range over $I$. 
Corollary~\ref{cor:deriv of poly expr} yields
 \eqref{eq:universal der, poly} for each $P\in A$. Thus the $A$-module $\Omega_{A|K}$ is 
generated by the $\d X_i$. Next, $\partial/\partial X_j\in \Der_K(A)$, and by~\eqref{eq:Hom Der 3} we have 
$\<\d X_i,\partial/\partial X_j\>=\delta_{ij}$ (Kronecker delta). This 
``orthogonality'' shows that the $A$-module $\Omega_{A|K}$ is free on  $(\d X_i)_{i\in I}$.   
\end{proof}

\noindent
Let $S$ be a multiplicative subset of $A$. The $K$-derivation
$\d\colon A \to \Omega:= \Omega_{A|K}$ yields by Lemma~\ref{lem:quotient rule, 1}
a $K$-derivation 
$$S^{-1}\d\ \colon\ S^{-1}A \to S^{-1}\Omega$$ 
into the $S^{-1}A$-module  $S^{-1}\Omega$. 

\begin{lemma}\label{lem:quotient rule, 2} The $K$-derivation 
$S^{-1}\d\colon S^{-1}A \to S^{-1}\Omega$ is a 
universal $K$-deri\-va\-tion of $S^{-1}A$. In particular,
$\Omega_{S^{-1}A|K}\ \cong\ S^{-1}\Omega_{A|K}$ as $S^{-1}A$-modules.
\end{lemma}
\begin{proof}
Let $N$ be an $S^{-1}A$-module and $D\in\Der_K\!\big(S^{-1}A, N\big)$. The
$S^{-1}A$-mo\-dule~$S^{-1}\Omega$ is generated by the $S^{-1}\d(a/1)=(\d a)/1$, so there is at most one
$S^{-1}A$-linear map $\psi\colon S^{-1}\Omega \to N$ with $\psi\circ S^{-1}\d=D$; our job 
is to find such $\psi$. The $K$-derivation $\Delta:=D\circ \iota\colon A \to N$
yields an $A$-linear $\phi\colon \Omega\to N$ with
 $\Delta= \phi\circ \d$. From $\phi$ we get an $S^{-1}A$-linear 
$\psi\colon S^{-1}\Omega \to N$ with $\psi(\omega/1)=\phi(\omega)$ for all 
$\omega\in \Omega$. This map $\psi$ has the desired property.
\end{proof}

\noindent
Next we give a useful alternative construction of the module of K\"ahler differentials of~$A$. Consider the $K$-al\-ge\-bra~${A\otimes_K A}$, with identity
$1\otimes 1$ and multiplication
given by $$(a_1\otimes a_2)\cdot (b_1\otimes b_2)\ =\ a_1b_1\otimes a_2b_2.$$
Then $A\otimes_KA$ is an $A$-algebra with structural morphism
$$a\mapsto a\otimes 1\colon A \to A\otimes_K A \quad\text{ (so $a(b_1\otimes b_2)=ab_1\otimes b_2$)}$$
which induces its given $K$-algebra structure.
We have an $A$-algebra morphism 
$$\mu\colon A\otimes_K A \to A,\qquad \mu(a\otimes b)\ =\ ab.$$ 
Set $J:= \ker\mu$. Then $1\otimes a-a\otimes 1\in J$ for all $a\in A$, and we 
get a $K$-linear map
$$d_A\colon A \to J/J^2, \qquad d_A(a)\ :=\ (1\otimes a - a\otimes 1)+J^2,$$
which is easily verified to be a $K$-derivation, using
$$1\otimes ab-b\otimes a - a\otimes b + ab\otimes 1\ =\ 
(1\otimes a - a\otimes 1)\cdot(1\otimes b - b\otimes 1)\ \in\ J^2.$$  
The universal property of $\d$ thus yields an $A$-linear map 
$$\Omega_{A|K} \to J/J^2\ \colon\ \d a \mapsto  d_A(a)\ =\ (1\otimes a - a\otimes 1)+J^2.$$ 

\begin{lemma} \label{isoIsq} This map $\Omega_{A|K} \to J/J^2$ is an isomorphism of $A$-modules.
\textup{(}Thus the $K$-derivation $d_A$ is universal.\textup{)}
\end{lemma}
\begin{proof} First, $J$ is generated as an $A$-module
by the $1\otimes a-a\otimes 1$: let $s\in J$ and take $a_1,\dots, a_n, b_1,\dots, b_n\in A$ with $s=\sum_i a_i \otimes b_i$; then $\mu(s)=\sum_i a_ib_i=0$, so
$$ s\ =\ \sum_i a_i(1\otimes b_i - b_i\otimes 1),$$
which proves the generation claim. Therefore, as an $A$-module, $J/J^2$ is generated  by the~$d_A(a)$. 
It remains to find an $A$-linear map
$\phi\colon J/J^2\to \Omega_{A|K}$ with $\phi\circ d_A=\d$. The $K$-bilinear map $(a,b) \mapsto a\d b\colon A\times A \to \Omega_{A|K}$ yields a $K$-linear map 
$\psi\colon A\otimes_K A \to \Omega_{A|K}$ with $\psi(a\otimes b)=a \d b$ for all
$a$,~$b$. Then $\psi$ is actually $A$-linear, 
$\psi(1\otimes a - a\otimes 1)=\d a$ for all $a$, and
$\psi$ sends each product $(1\otimes a-a\otimes 1)\cdot(1\otimes b-b\otimes 1)$ to $0$, and so $\psi$ induces an $A$-linear
map $\phi\colon J/J^2 \to \Omega_{A|K}$ with $\phi\circ d_A=\d$.
\end{proof}

\begin{cor}\label{cor:derivation on Omega} Let $\der$ be a derivation on $K$ and $\der_A$ a derivation on $A$ such that $\der_A(\kappa 1)=\der(\kappa)1$ for all $\kappa\in K$.
Then there is a unique $\der_A$-compatible derivation $\der^*$ on
the $A$-module $\Omega_{A|K}$ such that $\der^*(\d a)=\d\der_A(a)$ for all $a$. 
\end{cor}
\begin{proof} Since $\Omega_{A|K}$ is generated as an $A$-module by the
$\d a$, there can be at most one such $\der^*$.
 Corollary~\ref{cor:der on tensor prod, 3} yields a derivation  $\derdelta$ on the ring 
$A\otimes_K A$ such that 
$$\derdelta(a\otimes b)\ =\ \der_A(a)\otimes b + a\otimes\der_A(b)\qquad\text{for all $a$,~$b$.}$$
Then $\derdelta(J)\subseteq J$, so
$\derdelta(J^2)\subseteq J^2$, hence~$\derdelta$ induces a  
$\der_A$-compatible derivation $\der^*$ on~$J/J^2$ with  
$\der^*(d_A(a))=d_A(\der_A(a))$ for all $a$. Now use  Lemma~\ref{isoIsq}.
\end{proof}

\subsection*{The fundamental exact sequences}
Let $L$ be a commutative ring, $B$ be a commutative $L$-algebra, and suppose we are given
a commutative diagram
$$\xymatrix{ A\ar[r]^{h} & B \\
K \ar[r]\ar[u] & L\ar[u]}$$ 
of ring morphisms, where the vertical arrows are the structural morphisms.

\begin{lemma}\label{lem:functorial map}
With the $B$-module $\Omega_{B|L}$ construed as an $A$-module via $h$, there is exactly one $A$-linear map $\Omega_h\colon \Omega_{A|K} \to \Omega_{B|L}$ such that
the diagram below commutes:
$$\xymatrix{\Omega_{A|K} \ar[r]^{\Omega_h} & \Omega_{B|L} \\
A \ar[r]^h \ar[u]^{\d_{A|K}} & B \ar[u]_{\d_{B|L}}}$$
\end{lemma}
\begin{proof}
The map $\d_{B|L}\circ h$ is a $K$-derivation of $A$; so the existence and uniqueness of~$\Omega_h$
follow from the universal property of $\d_{A|K}$.
\end{proof}

\noindent
Now let $B$ be a commutative $K$-algebra and $h\colon A\to B$ a $K$-algebra morphism. Then we are in the situation above with $L=K$ and the identity
map $K \to L$,
so Lemma~\ref{lem:functorial map} yields an $A$-linear map 
$$\Omega_h\ \colon\ \Omega_{A|K} \to \Omega_{B|K}, \quad \Omega_h(\d a)\ =\ \d h(a),$$
which in turn (see previous section) yields a $B$-linear map
$$\alpha:=(\Omega_h)_B \colon B\otimes_A \Omega_{A|K} \to \Omega_{B|K},
\quad \alpha(1\otimes\d a)=\d h(a)\ \text{ for all $a$.}$$
Applying Lemma~\ref{lem:functorial map}
to another diagram yields a 
$B$-linear map
$$\beta\ \colon\ \Omega_{B|K}\to\Omega_{B|A}, \qquad 
\beta(\d_{B|K} f)\ =\ \d_{B|A} f\ \text{ for $f\in B$.}$$

\begin{prop}[First fundamental exact sequence]\label{prop:first exact sequ}
The following sequence of $B$-modules and $B$-linear maps is exact:
$$%\begin{equation}\label{eq:1st exact}
B\otimes_A \Omega_{A|K} \xrightarrow{\ \alpha\ } \Omega_{B|K} \xrightarrow{\ \beta\ } \Omega_{B|A}\longrightarrow 0
$$%\end{equation}
\end{prop}
\begin{proof}
Let $E$ be a $B$-module. Then \eqref{eq:iso tensor, 1} yields
an isomorphism
$$ \Hom_A\!\big(\Omega_{A|K},E\big) \xrightarrow{\ \cong\ }
\Hom_B\!\big(B\otimes_A \Omega_{A|K},E\big)$$
of $B$-modules. Identifying these two $B$-modules via this isomorphism, 
we have a commutative diagram of $A$-linear maps:
$$\xymatrix@C+=3.4em{
\Hom_A(\Omega_{A|K},E)\ar[d]^\cong  & \ar[l]_{\alpha^*}\ar[d]^\cong \Hom_B(\Omega_{B|K},E) & \ar[l]_{\beta^*}\ar[d]^\cong  \Hom_B(\Omega_{B|A},E) & \ar[l] 0 \\ 
\Der_K(A,E)  &  \ar[l]_{\Der_K(h,E)} \Der_K(B,E) & \ar[l]_\supseteq \Der_A(B,E)  & \ar[l] 0
}$$  
where the vertical arrows are isomorphisms given by \eqref{eq:Hom Der}
and the bottom row is exact by Lemma~\ref{lem:exact sequ 2}. Hence the top row is also exact.
Now use Lemma~\ref{lem:Hom and exact sequ}.
\end{proof}

\noindent
Next, let $I$ be an ideal of $A$, set $B:= A/I$, and let $h\colon A \to B$ be
the canonical map. To determine $\ker\alpha$ in this case, first note that 
$\d_{A|K}(I^2)\subseteq I\Omega_{A|K}$.
Composing the natural surjection
$\Omega_{A|K}\to \Omega_{A|K}/I\Omega_{A|K}$ with  the restriction of
$\d_{A|K}$  to $I$    yields an $A$-linear map
$I\to   \Omega_{A|K}/I\Omega_{A|K}$ whose kernel contains $I^2$,  hence in turn induces a $B$-linear map
$\gamma\colon I/I^2\to \Omega_{A|K}/I\Omega_{A|K}$. 
Identifying the $B$-modules $\Omega_{A|K}/I\Omega_{A|K}$ and $B\otimes_A \Omega_{A|K}$ as in~\eqref{eq:iso tensor, 3}, we get 
$$ \gamma\ :\ I/I^2\to B\otimes_A \Omega_{A|K}, \quad 
\gamma(a+I^2)\ =\ 1_B\otimes\d a\ \text{ for $a\in I$.}$$

\begin{prop}[Second fundamental exact sequence]\label{prop:second exact sequ}
The following sequence of $B$-modules and $B$-linear maps is exact:
$$I/I^2 \xrightarrow{\ \gamma\ } B\otimes_A \Omega_{A|K} \xrightarrow{\ \alpha\ } \Omega_{B|K}\longrightarrow 0$$
\end{prop}
\begin{proof}
Let $E$ be a $B$-module. As in the proof of 
Proposition~\ref{prop:first exact sequ} we get a commutative diagram of 
$A$-linear maps:
$$\xymatrix@C+=3.5em{
\Hom_B(I/I^2,E) \ar[d]^= & \ar[l]_{\gamma^*}  \ar[d]^\cong \Hom_A(\Omega_{A|K},E) & \ar[l]_{\alpha^*}  \ar[d]^\cong \Hom_B(\Omega_{B|K},E) & \ar[l] 0 \\
\Hom_B(I/I^2,E)  &  \ar[l]_{\quad \overline{\Delta}\,\mapsfrom\, \Delta} \Der_K(A,E) & \ar[l]_{\Der_K(h,E)} \Der_K(B,E)  & \ar[l] 0
}$$ 
with exact bottom row by Lemma~\ref{lem:exact sequ 3}. Now use Lemma~\ref{lem:Hom and exact sequ}. 
\end{proof}

\noindent
Let $R$ be a commutative $K$-algebra, $\mathfrak p\in\Spec(R)$, and set $A=R_{\mathfrak p}$.
Then $A$ is a $K$-algebra and a local ring with maximal ideal $\m=\mathfrak p A$. Let $F:=A/\m$ be its residue field.
%(so $F\cong\Frac(R/\mathfrak p)$). 
Then Lemma~\ref{lem:quotient rule, 2} and the isomorphisms
\eqref{eq:iso tensor, 4} and \eqref{eq:iso tensor, 2} give
$$F\otimes_A \Omega_{A|K}\ \cong\ F\otimes_A (A\otimes_R \Omega_{R|K})\ \cong\ F\otimes_{R} \Omega_{R|K}\qquad\text{($F$-linear isomorphisms)}$$
with $1\otimes \d(r/1)\in F\otimes_A \Omega_{A|K}$ corresponding to $1\otimes \d r\in F\otimes_{R} \Omega_{R|K}$
for $r\in R$. In view of Proposition~\ref{prop:second exact sequ} this leads to an exact sequence 
\begin{equation}\label{eq:second exact sequ}
\m/\m^2 \xrightarrow{\ \gamma_0\ } F\otimes_{R} \Omega_{R|K} \xrightarrow{\ \alpha_0\ } \Omega_{F|K}\longrightarrow 0
\end{equation}
of vector spaces over $F$ and $F$-linear maps, with $\gamma_0\big((r/1)+\m^2\big) = 1\otimes\d r$ for $r\in\mathfrak{p}$
and $\alpha_0(1\otimes \d r)=\d h(r/1)$
for $r\in R$, where $h\colon A \to F$ is the canonical map.

\subsection*{Notes and comments}
The module $\Omega_{A|K}$ of differentials was introduced by K\"ahler~\cite{Kaehler}. The description of $\Omega_{A|K}$ via Lemma~\ref{isoIsq} is due to Cartier~\cite{Cartier}; see~\cite{Kunz03} for the history.
For more information on $\Omega_{A|K}$, see~\cite{Kunz86}.
Corollary~\ref{cor:derivation on Omega} is from~\cite{JJ69-1}.
If $K$ is a field of characteristic zero, then the $F$-linear map~$\gamma_0$ 
from~\eqref{eq:second exact sequ}
is injective; see \cite[Corollary~6.5~(a)]{Kunz86} or \cite[Theorems~25.2 and 26.9]{Matsumura}.
See Chapter~\ref{ch:triangular automorphisms} below for $K$-derivations on non-commutative $K$-algebras.

\section{Derivations on Field Extensions}\label{sec:derivations on field exts}

\noindent
{\em Let $K$ be a field and $L$ a field extension of $K$}. Below we indicate ways to
extend a derivation on $K$ to a derivation on $L$.
%We use this to establish a fact appealed to in the proof of Proposition~\ref{prop:hensel, multivar}. 
In the case that $\operatorname{char}(K)=0$
we also relate bases for the vector space $\Omega_{L|K}$ over $L$ to 
transcendence bases of $L$ over $K$.

\subsection*{Extending derivations}
In this subsection we fix a derivation $\der$ on $K$.
Let $X=(X_1,\dots, X_n)$ be a tuple of distinct indeterminates
and let $P=\sum_{\i}P_{\i}X^{\i}\in K[X]$ be a polynomial: the 
$\i=(i_1,\dots, i_n)$ range over a finite subset of $\N^n$ and
$P_{\i}\in K$, $X^{\i}:= X_1^{i_1}\cdots X_n^{i_n}$ for all~$\i$. Set
$P^{\der}:= \sum_{\i}\der(P_{\i})X^{\i}\in K[X]$.  \nomenclature[O]{$P^\der$}{result of applying the derivation $\der$ to the coefficients of~$P$} It is easy to verify that~${P\mapsto P^\der}$ is a derivation on~$K[X]$ extending $\der$.
The following identity is a special case of Lemma~\ref{derpolA} (but we allow $n=0$ here) and is 
used frequently:  

\begin{lemma}\label{derpol} For all $a=(a_1,\dots, a_n)\in K^n$,
$$\der\big(P(a)\big) \ =\  P^\der(a) + \sum_{j=1}^n \frac{\partial P}{\partial X_j}(a)
\cdot \der(a_j).$$
\end{lemma}

\noindent
Recall: $L$ is separably algebraic over $K$ iff for every
$a\in L$ we have $P(a)=0$ and $P'(a)\ne 0$ for some $P\in K[X]$. Thus if 
$\operatorname{char}(K)=0$ and $L$ is algebraic over $K$, then $L$ is separably algebraic over $K$. 

\begin{lemma}\label{lem:derivative on alg ext}
Suppose that $L$  is separably algebraic over $K$. Then $\der$ extends uniquely to a
derivation on~$L$.
\end{lemma}
\begin{proof}
Let $\der$ be extended to a derivation on $L$. Let us denote this extended derivation also by $\der$. For any $a\in L$ with minimum polynomial $P\in K[X]$
over $K$ we have $P(a)=0$ and $P'(a)\neq 0$, and hence by \ref{derpol},
\begin{equation}\label{eq:derivative on alg ext}
\der(a)\ =\ \frac{-P^{\der}(a)}{P'(a)}.
\end{equation}
So there is at most one such derivation on $L$. To construct such a derivation, let $a\in L$; it suffices to get a derivation on $K[a]=K(a)$ extending $\der$. Let $P\in K[X]$ be the minimum polynomial of $a$ over $K$, so the ring morphism
$$ Q\mapsto Q(a)\ \colon\ K[X]\to K[a]$$
has kernel $PK[X]$. Consider the additive map
$$d\colon K[X]\to K[a]=K(a),\qquad d(Q)\ =\ Q^{\der}(a)+Q'(a)\frac{-P^{\der}(a)}{P'(a)}.$$
An easy computation shows that $d(Q)=0$ for $Q\in PK[X]$, so $d$ induces an additive map $K[a]\to K[a]$
sending $Q(a)$ to $Q^{\der}(a)+Q'(a)\frac{-P^{\der}(a)}{P'(a)}$, for all $Q\in K[X]$. This last map is easily checked to be a derivation on $K[a]$ extending~$\der$.
\end{proof}

\noindent
Let $L$ be as in Lemma~\ref{lem:derivative on alg ext}, and let 
$\der$ be extended to a derivation 
on $L$. Then any subfield of $L$ containing $K$ is closed under 
this derivation of $L$, by \eqref{eq:derivative on alg ext}.

\begin{lemma}\label{lem:Xy} Let $x_1,\dots, x_n, y_1,\dots, y_n\in L$ be such
that $x_1,\dots, x_n$ are algebraically independent over $K$, set
$x=(x_1,\dots, x_n)$, and consider
$L$ as a $K[x]$-module via the inclusion $K[x]\to L$. 
Then there is a unique $\Z$-derivation
$K[x] \to L$ that extends $\der$ and sends $x_j$ to $y_j$ for $j=1,\dots,n$.
\end{lemma} 
\begin{proof} By Lemma~\ref{derpolA} there is at most one such extension. 
The map $$ P(x)\ \mapsto\  P^\der(x) + 
\sum_{j=1}^n \frac{\partial P}{\partial X_j}(x)y_j\ \colon\ K[x] \to L$$
is easily verified to be an extension as required.
\end{proof}

\noindent
Here is an easy consequence of  Lemma~\ref{lem:Xy} and Corollary~\ref{cor:derivation into fraction field}:

\begin{cor}\label{lem:derivative on trans ext}
Suppose that $L=K(x)$ where $x=(x_i)_{i\in I}$ is a family in~$L$ that is 
algebraically independent over $K$.
Then there is for each family $(y_i)_{i\in I}$ in $L$ a unique 
extension of $\der$ to a derivation on $L$ with $\der(x_i)=y_i$ for all $i\in I$. 
\end{cor}

\begin{lemma}\label{lem:derivative on purely insep ext}
Suppose that $K$ has characteristic $p > 0$ and $L = K(a)$ where
$a^p = c \in K\setminus K^p$ and $\der(c) = 0$. Then there is for each 
$b \in L$ a unique extension of~$\der$ to a 
derivation
on $L$ with $\der(a) = b$.
\end{lemma}
\begin{proof} Let $b \in L$. The minimum polynomial of $a$ over $K$ is 
$X^p - c$, so the ring
morphism $K[X] \to K[a]\colon Q \mapsto  Q(a)$ is surjective with kernel $(X^p-c)K[X]$. 
View $K[a] = K(a) = L$ as a $K[X]$-algebra via this morphism. Then the map
$$ d\ :\ K[X] \to K[a],\qquad d(Q)\ : =\  Q^{\der}(a) + Q'(a)b,$$
is a $\Z$-derivation, and $d(X^p-c) = -\der(c) = 0$, hence $d$ induces a derivation on 
$K[a]$ sending $a$
to $b$ and extending $\der$. Uniqueness is clear by Lemma~\ref{derpol}.
\end{proof}

\noindent
Call $L$ {\em  separably generated\/} over $K$ if $L$ is separably algebraic over $K(B)$ for some transcendence base $B$ of $L|K$. 
In particular, if $\operatorname{char} K=0$, then $L$ is separably generated over $K$. 

\begin{cor}\label{cor:characterization of sep alg} Suppose $L$ is separably generated over $K$. Then there exists an extension of
$\der$ to a derivation on $L$; there is a unique such extension iff
$L$ is algebraic over $K$.
\end{cor}
\begin{proof} Let $B$ be a transcendence basis for $L|K$ such that
$L$ is separably algebraic over $K(B)$.
%, and let $F$ be the 
%subfield of~$L$ consisting of the elements of $L$ that are separably algebraic over $K(B)$. 
We can extend~$\der$ 
to a derivation on $K(B)$ by Corollary~\ref{lem:derivative on trans ext},
in more than one way if $B\ne \emptyset$. Any such extension then 
extends further to a derivation on~$L$ by Lemma~\ref{lem:derivative on alg ext}.
%If $F\ne L$, then any derivation on $F$ extends in more than one way to
%a derivation on the purely inseparable extension $L$ of $F$, by 
%Lemma~\ref{lem:derivative on purely insep ext}. 
\end{proof}

\noindent
For $\der=0$ we have the following variant:

\begin{cor}\label{bettersepalg} Suppose $L$ is finitely generated as a field over $K$.
Then: 
\begin{align*}
 \text{$L$ is separably algebraic over $K$}\ \Longleftrightarrow\ &
\begin{cases}
\ \ \parbox{16em}{any derivation on $L$ that extends\\
the trivial derivation on $K$ is trivial.}
\end{cases}
\end{align*} 
\end{cor}
\begin{proof} Lemma~\ref{lem:derivative on alg ext} gives $\Rightarrow$. If $\operatorname{char} K=0$, then Corollary~\ref{cor:characterization of sep alg} gives $\Leftarrow$.  
Suppose $\operatorname{char} K=p>0$ and $L$ is not separably algebraic over $K$. We have a tower $K=K_0\subseteq K_1\subseteq \cdots \subseteq K_n=L$ of intermediate fields with $K_{i+1}=K_{i}(x_i)$ for $i < n$,
such that either
(1) $x_i$ is transcendental over $K_{i}$, or (2) $x_i^p\in K_i,\ x_i\notin K_{i}$, or (3)
$x_i$ is separably algebraic over $K_{i}$. Take $i<n$ maximal such that 
(1) or (2) holds. Then by Corollary~\ref{lem:derivative on trans ext} and
Lemma~\ref{lem:derivative on purely insep ext} the trivial derivation on $K_i$ extends to more than one derivation on $K_{i+1}$, each of which extends to a derivation on $L$.
     \end{proof}

\noindent
The next result will be used in the proof of
Proposition~\ref{prop:hensel, multivar}. For an 
$n$-tuple $P=(P_1,\dots,P_n)\in K[X]^n$ we define the
$n\times n$ matrix
$$P'\ :=\ \left(\frac{\partial P_i}{\partial X_j}\right)_{i,j=1,\dots,n} \qquad\text{(entries in $K[X]$)}.$$
Evaluating its entries at a point $a\in L^n$ gives the 
$n\times n$ matrix $P'(a)$ over $L$.

\begin{cor}\label{cor:regsepalg}
Let $P=(P_1,\dots,P_n)\in K[X]^n$ and $a=(a_1,\dots,a_n)\in L^n$ be such that $P(a)=0$  and $\det P'(a)\neq 0$. Then $a_1,\dots,a_n$ are separably algebraic over~$K$.
\end{cor}
\begin{proof}
Let $\Delta\in\Der_K K(a)$. Then 
Corollary~\ref{cor:deriv of poly expr} and the assumptions $P(a)=0$  and $\det P'(a)\neq 0$ give $\Delta(a_1)=\cdots=\Delta(a_n)=0$, so $\Delta=0$. Hence $K(a)|K$ is separably algebraic by
Corollary~\ref{bettersepalg}.
\end{proof}

\subsection*{K\"ahler differentials for field extensions}
If $L$ is separably algebraic over $K$, then $\Omega_{L|K}=\{0\}$ by 
Lemma~\ref{lem:derivative on alg ext} and the isomorphism~\eqref{eq:Hom Der 2}. We now generalize this fact for $K$ of characteristic zero:

\begin{lemma}\label{lem:lindep in Omega} Suppose $\operatorname{char} K=0$. Let $(x_i)_{i\in I}$ be a family in $L$. Then the family $(x_i)$ is algebraically independent 
over $K$ if and only if the family  $(\d x_i)$ in $\Omega_{L|K}$ is linearly independent over $L$. 
\end{lemma}
\begin{proof}
Suppose $(x_i)$ is algebraically independent over $K$.  
Corollary~\ref{cor:characterization of sep alg} and its proof yield for each $i\in I$ a 
$K$-derivation $\der_i$ on $L$ such that $\der_i(x_j)=\delta_{ij}$ (Kronecker delta) for all $j\in I$. Thus  $(\d x_i)$ is linearly
independent over $L$ by~\eqref{eq:Hom Der 3}. 

Next, assume
$x_1,\dots, x_n\in L$ are algebraically dependent over $K$. 
Set $x:=(x_1,\dots,x_n)$ and
let $P\in K[X_1,\dots,X_n]$ be a nonzero polynomial such that $P(x)=0$, and 
of minimal degree with these properties. This minimality gives $i$ with
$\frac{\partial P}{\partial X_i}\neq 0$ and
so $\frac{\partial P}{\partial X_i}(x)\neq 0$.
Then by Corollary~\ref{cor:deriv of poly expr},  
$$\sum_{i=1}^n \frac{\partial P}{\partial X_i}(x)\,\d x_i\ =\  \d P(x)\ =\  0\ \text{ in $\Omega_{L|K}$,}$$
so $\d x_1,\dots,\d x_n$ are linearly dependent over $L$.
\end{proof}

\noindent 
In particular, if  $\operatorname{char} K=0$, then a family $(x_i)_{i\in I}$ in
$L$ is a transcendence basis for~$L|K$ iff $(\d x_i)$ is a basis
of the vector space $\Omega_{L|K}$ over $L$, and so
$$\dim_L  \Omega_{L|K}\ =\ \operatorname{trdeg}(L|K).$$
In Section~\ref{sec:johnson} we prove an analogue of this for extensions of
differential fields.

\subsection*{Notes and comments}
This section stems from A.~Weil~\cite[Chapter~I, \S{}5]{Weil}.

%% file: mt-2.tex
\chapter{Valued Abelian Groups}

\noindent
Our main objects of interest are fields like $\mathbb{T}$
that are equipped with a compatible valuation and
derivation. To analyze these objects we
need valuation theory, 
%some of it not well-known, 
and so
we include two chapters containing the purely valuation-theoretic tools.
Since these tools come from a variety of sources scattered over the 
literature and some are new, we include 
all but the most routine proofs.
%, for the benefit of readers without 
%prior exposure to valuation theory. 
%We hope that this way, the present chapter can serve as 
%an introduction to the part
%of valuation theory relevant to our study of~$\mathbb T$, 
%even for a reader without prior exposure to this subject. 
%In accordance with our general set-up we shall be self-contained as far as definitions and statements of results is concerned, but in such an established subject we can refer to the literature for many proofs. 

After introducing some terminology concerning ordered sets in Section~\ref{sec:ordered sets}, this chapter treats valued abelian groups (Sections~\ref{sec:valued abelian gps},~\ref{sec:valued vector spaces}) and
ordered abelian groups (Section~\ref{sec:oag}); the latter will occur as 
value groups of valued fields.
Valued abelian groups arise in our work because
the logarithmic derivative map on a 
valued differential field like $\mathbb T$ induces a 
valuation on the value group that turns out to be very useful. Moreover, the notion of 
a \textit{pseudocauchy sequence}\/
makes perfect sense in the general setting of valued abelian groups, and
the basic facts about these sequences yield a natural proof of a 
generalized Hahn Embedding Theorem which
can serve as a model for later proofs of several much deeper 
embedding theorems. 

%Some basic texts on valuation theory are \cite{Prestel-Engler, %Ribenboim}, and the comprehensive (as yet to be published) book 
%\cite{FVK}. For more on ordered sets see \cite{Harzheim}, and on %ordered algebraic structures~\cite{Fuchs, PC}. For ordered fields %in particular see \cite{BCR,KS}. 

\input{mt-2-1}

\input{mt-2-2}

\input{mt-2-3}

\input{mt-2-4}

%% file: mt-2-1.tex
\section{Ordered Sets}\label{sec:ordered sets}

\noindent
%Here we define the basic notions concerning ordered sets that are 
%extensively used later. 
This section serves mainly to fix notations and terminology.
By convention {\em ordered set\/} means 
{\em totally ordered set\/} unless specified otherwise. This agrees with 
the usual meaning of {\em ordered abelian group\/} and
{\em ordered field}. \index{ordering}\index{ordered set}\index{set!ordered}

\medskip
\noindent
Let $S$ be an ordered set. We denote the ordering on $S$ by 
$\leq$, the corresponding strict ordering by $<$, and the reversals of $\le$ and $<$ by $\ge$ and $>$, respectively. 
We extend these notations to sets $A,B\subseteq S$ by
$A\le B:\Longleftrightarrow \text{$a \le b$ for all $(a,b)\in A\times B$}$, and
$A < B:\Longleftrightarrow \text{$a < b$ for all $(a,b)\in A\times B$}$. 
Also for $b\in S$, $A\le b:\Longleftrightarrow A\le \{b\}$, and so on. 
We view a subset of $S$ as ordered by the induced ordering.
We put 
$S_\infty := S \cup \{\infty\}$,
$\infty\notin S$, with the ordering on $S$ extended to a (total)
ordering on~$S_\infty$ by $S<\infty$. Occasionally, we even take two distinct
elements $-\infty, \infty \notin S$, and extend the ordering on $S$ to an ordering
on $S\cup \{-\infty, \infty\}$ by $-\infty < S < \infty$. For $A\subseteq S$ and $b\in S$ we  
set $A^{\le b}:=\{a\in A:\ a\le b\}$; similarly for $<$, $\ge$, and $>$ instead of $\le$.
For $a,b\in S$ we put \nomenclature[E]{$[a,b]=[a,b]_S$}{convex hull of $\{a,b\}$ in $S$}
$$[a,b]\ =\ [a,b]_S\ :=\ \{x\in S:\ a\leq x\leq b\}.$$
A subset $C$ of $S$ is said to be {\bf convex} in $S$ if for all 
$a,b\in C$ we have $[a,b]\subseteq C$.\index{convex} For $A\subseteq S$ we let \nomenclature[E]{$\operatorname{conv}(A)$}{convex hull of $A$} 
$$\operatorname{conv}(A)\ :=\ \{ x\in S:\ \text{$a\leq x\leq b$ for some 
$a,b\in A$} \}$$
be the {\bf convex hull of $A$} \index{convex!hull}\index{hull!convex} in $S$, that is, the smallest convex 
subset of $S$ containing~$A$.  For $-\infty\leq a<b\leq\infty$ we also set
$$(a,b)\ =\ (a,b)_S\ 	:=\ \{x\in S:\ a < x < b\}.$$
The sets of the form $(a,b)$ are called {\bf intervals} in $S$. 
The intervals in $S$ form a basis for a hausdorff topology on $S$, 
the {\bf interval topology} or {\bf order topology} on $S$. 

\nomenclature[E]{$(a,b)=(a,b)_S$}{interval in $S$}
\index{order!topology}
\index{interval}
\index{interval!topology}
\index{topology!interval}
\index{topology!order}

\medskip
\noindent
The ordered set $S$ is said to be {\bf dense} if for all $a<b$ in $S$ we have $(a,b)\neq\emptyset$, and
{\bf without endpoints} if for all $a,b\in S$ we have $(a,\infty)\neq\emptyset$ and $(-\infty,b)\neq\emptyset$.

\index{ordered set!dense}
\index{ordered set!without endpoints}

\medskip
\noindent
Let $S'$ be also an ordered set and $f\colon S\to S'$ a map.
Then $f$ is said to be 

\storestyleof{enumerate} \begin{listliketab}
\newcounter{tabenum}\setcounter{tabenum}{0} \newcommand{\nextnum}{(\addtocounter{tabenum}{1}\thetabenum)} 
\begin{tabular}{Lll}
\nextnum &{\bf increasing} if for all $a,b\in S$:  & $a\leq b\Rightarrow f(a)\leq f(b)$; \\
\nextnum &{\bf strictly increasing} if for all $a,b\in S$: & $a< b\Rightarrow f(a)< f(b)$; \\
\nextnum &{\bf decreasing} if for all $a,b\in S$:  & $a\leq b\Rightarrow f(a)\geq f(b)$; \\
\nextnum &{\bf strictly decreasing} if for all $a,b\in S$: & $a< b\Rightarrow f(a)>f(b)$.\\
\end{tabular}
\end{listliketab}

\index{increasing map}
\index{strictly!increasing}
\index{decreasing map}
\index{strictly!decreasing}
\index{map!increasing}
\index{map!strictly increasing}
\index{map!decreasing}
\index{map!strictly decreasing}
\index{isomorphism!ordered sets}
\index{ordered set!isomorphism}
\index{intermediate value property}
\index{ordered set!well-ordered}
\index{well-ordered}
\index{ordered set!order type}
\index{order!type}
\index{type!order}
\nomenclature[E]{$\operatorname{ot}(S)$}{order type of the ordered set $S$}

\smallskip
\noindent
An increasing bijection $S\to S'$ is called an {\bf isomorphism} of 
ordered sets. 
We say that~$f$ has the {\bf intermediate value property} if for all $a<b$ in 
$S$, 
\begin{align*} f(a) < f(b)\ &\Longrightarrow\ \big(f(a),f(b)\big)_{S'}\subseteq f\big((a,b)_S\big),\\
 f(a) > f(b)\ &\Longrightarrow\ \big(f(b),f(a)\big)_{S'}\subseteq 
f\big((a,b)_S\big).
\end{align*}
The ordered set $S$ is said to be {\bf well-ordered} if every nonempty 
subset of $S$ has a smallest element; equivalently, there is no infinite 
sequence $$a_0>a_1>\cdots>a_n>\cdots$$ in $S$.
%$(a_n)$ in $S$ with $a_0>a_1>\cdots>a_n>\cdots$. 
We say that two ordered sets {\bf have the same order type} if they are 
isomorphic.
If $S$ is well-ordered then there is a unique ordinal number,  
denoted by~$\operatorname{ot}(S)$, with the same order type as $S$.
(As usual, an ordinal is construed here in von~Neumann's sense as the set
of all smaller ordinals.) 
If $S$ is well-ordered and $S'$ is an ordered subset of $S$, then 
$S'$ is well-ordered with $\operatorname{ot}(S')\leq\operatorname{ot}(S)$.

\begin{lemma}\label{lem:wo under increasing maps}
Let $S'$ be an ordered set and $f\colon S\to S'$ increasing and surjective. If~$S$ is well-ordered, then so is $S'$, with
$\operatorname{ot}(S)\geq\operatorname{ot}(S')$.
\end{lemma}
\begin{proof}
For each $s'\in S'$ pick $g(s')\in f^{-1}(s')$. 
Then the map $g\colon S'\to S$ is strictly increasing, and the claims follow.
\end{proof}

\noindent
A subset $A$ of $S$ is said to be a {\bf cut} \index{cut} in $S$, or {\bf downward closed} \index{downward!closed} in $S$, if for all $a\in A$ and $s\in S$ we have $s< a \Rightarrow s\in A$. The empty subset of $S$ and $S$ itself are cuts in~$S$; these are called the {\bf trivial} cuts in $S$.  The collection of cuts in an ordered set is totally ordered by inclusion (with largest element $S$ and smallest element $\emptyset$), and 
the union and intersection of any set of cuts are also cuts.
We say that an element~$x$ of an ordered set extending $S$ {\bf realizes} 
the cut $A$ in $S$ if $A<x<S\setminus A$.

%A non-trivial cut in $S$ without a largest element is 
%called a {\bf dedekind cut} in $S$.

Dually, $A\subseteq S$ is {\bf upward closed} \index{upward!closed} in $S$ if for all $a\in A$ and 
$s\in S$ we have $a<s \Rightarrow s\in A$. (So $A\subseteq S$ is upward closed in $S$ iff $S\setminus A$ is downward closed in~$S$.)
For $A\subseteq S$ we put 
\begin{align*}
A^{\downarrow} &:= \{ s\in S: \text{$s\leq a$ for some $a\in A$}\}, \\
A^{\uparrow} &:= \{ s\in S: \text{$a\leq s$ for some $a\in A$}\};
\end{align*}
clearly $A^{\downarrow}$ is the smallest downward closed subset of $S$ containing $A$, and 
 $A^{\uparrow}$ is the smallest upward closed subset of $S$ containing $A$.
A subset $A$ of $S$ is said to be {\bf cofinal} \index{cofinal} in $S$ if $A^\downarrow=S$, i.e.,
if for each $s\in S$ there is $a\in A$ with $s\le a$, and $A$ is said to be {\bf coinitial} \index{coinitial} in $S$ if $A^\uparrow=S$, that is,
if for each $s\in S$ there is $a\in A$ with~$s\ge a$. 
The relation of cofinality is transitive: if $A\subseteq S$ is cofinal in $S$ and $B\subseteq A$ is cofinal in $A$, then $B$ is cofinal in $S$. A map 
$f\colon S \to S'$ into an ordered set~$S'$ is said to be {\bf cofinal} if~$f(S)$ is cofinal in $S'$.
%; similarly with ``coinitial'' in place of ``cofinal.''

\index{cofinality}
\index{coinitiality}
\index{ordered set!cofinality}
\index{ordered set!coinitiality}
\nomenclature[E]{$\operatorname{cf}(S)$}{the cofinality of an ordered set $S$}
\nomenclature[E]{$\operatorname{ci}(S)$}{the coinitiality of an ordered set $S$}

\nomenclature[E]{$A^{\uparrow}$}{smallest upward closed subset containing $A$}
\nomenclature[E]{$A^{\downarrow}$}{smallest downward closed subset containing $A$}
\index{regular!ordinal}

\medskip
\noindent
As a consequence of the Axiom of Choice, there always exists a well-ordered cofinal subset of $S$. 
The {\bf cofinality} of $S$, denoted by $\operatorname{cf}(S)$, is defined to be  
the smallest ordinal~$\lambda$ such that there exists a cofinal subset of $S$ of order type $\lambda$. This cofinality~$\lambda$ is actually a cardinal; here we identify as usual an ordinal $\alpha$ with the set of ordinals~${<\alpha}$, and a cardinal~$\kappa$ with the least ordinal of cardinality $\kappa$.
Dually, the {\bf coinitiality} of $S$, denoted by $\operatorname{ci}(S)$, is the cofinality of $S$ equipped with the reversed ordering, that is, the smallest ordinal $\lambda$ such that there exists a coinitial subset of~$S$ of reversed
order type~$\lambda$. %The cofinality of an ordinal $\alpha$ can then be characterized as the smallest ordinal $\lambda\leq\alpha$ such that $\lambda$ is cofinal in $\alpha$. 
Thus $\operatorname{cf}(\alpha)\leq\alpha$ for every ordinal $\alpha$. An ordinal $\alpha$ such that $\operatorname{cf}(\alpha)=\alpha$ is called {\bf regular.} Using transitivity of cofinality one easily sees that $\operatorname{cf}(S)$ is regular.

\begin{lemma}\label{lem:cofinality}
Let $S'$ be a cofinal subset of $S$. Then $\operatorname{cf}(S')= \operatorname{cf}(S)$.
\end{lemma}
\begin{proof} Transitivity of cofinality gives $\operatorname{cf}(S')\geq \operatorname{cf}(S)$. 
To get $\operatorname{cf}(S)\geq \operatorname{cf}(S')$ we can replace $S'$ by a well-ordered cofinal subset of order type $\operatorname{cf}(S')$
and assume that $S'$ is well-ordered.
Take a strictly increasing and cofinal sequence $\big(s(\gamma)\big)$ in $S$ 
indexed by the ordinals $\gamma<\lambda:=\operatorname{cf}(S)$. Define the sequence $\big(s'(\gamma)\big)$ in $S'$ indexed by the same ordinals by
$$s'(\gamma) := \min\!\big\{ s'\in S': s'\geq s(\gamma)\big\}.$$
The sequence  $\big(s'(\gamma)\big)$ is increasing, hence by Lemma~\ref{lem:wo under increasing maps} the order type of its image~$s'(\lambda)$ satisfies
$\lambda\geq \operatorname{ot}\!\big(s'(\lambda)\big)$. Since $s'(\lambda)$ is well-ordered and cofinal in $S'$, we also have $\operatorname{ot}\!\big(s'(\lambda)\big)\geq\operatorname{cf}(S')$, and thus 
$\operatorname{cf}(S)=\lambda \geq \operatorname{cf}(S')$ .
\end{proof}

\begin{cor}\label{cor:cofinality, 2}
Let $\alpha$ be a regular ordinal. Then every subset of $\alpha$ of order type~$\alpha$ is cofinal in $\alpha$.
%there is no ordinal $\beta<\alpha$ having a 
%cofinal subset of order type $\alpha$.
\end{cor}
\begin{proof}
Let $\beta$ be an ordinal and $S'$ a cofinal subset of $\beta$ with $\operatorname{ot}(S')=\alpha$. Then $\alpha=\operatorname{cf}(\alpha)=\operatorname{cf}(S')=\operatorname{cf}(\beta)\leq\beta$, where the third equality holds by 
Lemma~\ref{lem:cofinality}.
\end{proof}

\begin{cor} \label{cor:cofinality}
Let $S'$ be an ordered set without a largest element, and let $f\colon S\to S'$ be an increasing and cofinal map. Then $\operatorname{cf}(S)=\operatorname{cf}(S')$.
\end{cor}
\begin{proof}
Define an equivalence relation $\sim$ on $S$ by $x\sim y:\Longleftrightarrow f(x)=f(y)$. Pick a set $R$ of representatives for this equivalence relation, 
and observe that $f$ restricts to an isomorphism $R\to f(R)=f(S)$ of ordered 
sets and $R$ is cofinal in $S$. Hence 
$\operatorname{cf}(S)=\operatorname{cf}(R)=\operatorname{cf}\!\big(f(S)\big)=\operatorname{cf}(S')$, by Lemma~\ref{lem:cofinality}. 
\end{proof}

\noindent
A {\bf well-indexed sequence} is a sequence $(a_\rho)$
whose terms $a_\rho$ are indexed by the elements $\rho$
of an infinite well-ordered set without a last element. 
Restricting a well-indexed sequence $(a_\rho)$ to a cofinal subset of
its index set yields a well-indexed sequence; it is what we call a 
{\bf cofinal subsequence\/} of $(a_\rho)$. Given a set $A$,
a {\bf well-indexed sequence in $A$\/} is a well-indexed sequence  
whose terms are all in $A$.

\index{sequence!well-indexed}
\index{well-indexed}
\index{sequence!cofinal subsequence}

\subsection*{Notes and comments}
The material in this section is standard fare from the theory of ordered sets. For more on ordered sets see \cite{Harzheim}. The proof of Corollary~\ref{cor:cofinality} is from the proof of the ``Fact'' on p.~22 of \cite{vdDries-Tame-II}.

%% file: mt-2-2.tex
\section{Valued Abelian Groups}\label{sec:valued abelian gps}

\noindent
We introduce the basic notions concerning valued abelian groups, treat  
pseudocauchy and cauchy sequences in this context, and prove a fixpoint theorem.

\subsection*{Valued abelian groups} 
Let $G$ be an abelian (additively written) group. A {\bf valuation\/} on $G$ 
is a function
$v\colon G \to S_{\infty}$ where $S$ is an ordered set, such that for all $x,y\in G$ the following 
conditions are satisfied: 
\begin{list}{*}{\setlength\leftmargin{3em}}
\item[(VA1)] $v(x)=\infty \Longleftrightarrow x=0$;
\item[(VA2)] $v(-x)=v(x)$;
\item[(VA3)] $v(x+y) \ge \min\big(v(x),v(y)\big)$.
\end{list}
Let $v\colon G \to S_{\infty}$ be a valuation on $G$. Note that if 
$x,y \in G$ and $vx < vy$, then $v(x+y) = vx$.  
(Use that $vx = v(x+y-y) \geq \min\big(v(x+y),vy\big)$.)
More generally, if $x_1,\ldots,x_n \in G$, $n \geq 1$, and 
$v(x_1)<v(x_2),\ldots,v(x_n)$, then $v(x_1+\cdots+x_n) = v(x_1)$.
Another elementary fact that is often used is that the condition $v(x-y) > vx$
for $x,y\in G$ (which implies $x,y\ne 0$) defines an equivalence relation on 
$G^{\ne}$.

\index{valuation!on an abelian group}
\nomenclature[G]{$v$}{valuation on an abelian group}

\medskip\noindent
Assume now that the valuation
$v\colon G \to S_{\infty}$ is surjective. For 
$s\in S$ and $a\in G$ we define the 
{\bf open ball\/} $B_a(s)$ \nomenclature[G]{$B_a(s)$}{open ball centered at $a$ with 
 radius $s$} and the {\bf closed ball\/} $\bar{B}_a(s)$ \nomenclature[G]{$\bar{B}_a(s)$}{closed ball centered at $a$ with radius $s$} by  \marginpar{Do you really want to have the center $a$ in the subscript 
and the radius $s$ in parentheses, and not the other way around 
(which is more common, I think)?}
\begin{align*} 
B_a(s)\ &:=\ \big\{x\in G:\ v(x-a)> s\big\},\\
\bar{B}_a(s)\ &:=\ \big\{x\in G:\ v(x-a) \ge s\big\}.
\end{align*}
We refer to $B_a(s)$ as the {\bf open ball centered at $a$ with 
valuation radius $s$}, and to $\bar{B}_a(s)$ as the {\bf closed ball
centered at $a$ with valuation radius $s$}.  \index{ball} Note that $B(s):= B_0(s)$ \nomenclature[G]{$B(s)$}{$B_0(s)$} is a subgroup of $G$, and the 
$B_a(s)=a+B(s)$ are its cosets. Likewise, $\bar{B}(s):= \bar{B}_0(s)$ \nomenclature[G]{$\bar{B}(s)$}{$\bar{B}_0(s)$} is
a subgroup of $G$, and the $\bar{B}_a(s)=a+\bar{B}(s)$ are its cosets. 
Viewing balls in this way as cosets, it is clear that if $D$, $E$ are balls 
(of any kind) with nonempty intersection, then 
$D \subseteq E$ or $E \subseteq D$. Likewise, any point in a ball can serve as a center of that ball.

\medskip\noindent
The open balls form a basis for a topology on $G$, the {\bf $v$-topology}.  (If the valuation~$v$ is understood from the context we also speak of the {\bf valuation topology} on $G$.)
It is easy to see that balls (of any kind) are both open and closed in the 
$v$-topology, and that $G$ with the $v$-topology is a hausdorff 
topological group. If $S$ has a largest element, then the $v$-topology on 
$G$ is discrete. This is in particular the case if $v$ is 
{\bf trivial}, that is, $S$ is a singleton. If $S=\emptyset$, then
of course $G=\{0\}$. 

\index{trivial!valuation}
\index{valuation!trivial}
\index{valuation!topology}
\index{topology!valuation}

\nomenclature[G]{$G(s)$}{$B(s)/\bar{B}(s)$}

We have $B(s)\subseteq \bar{B}(s)$, and
we let $G(s):= \bar{B}(s)/B(s)$ be the corresponding
quotient group, which is nontrivial. 
Here is a useful and suggestive cardinality bound
on $G$ in terms of the cardinalities of these quotients:

\begin{lemma}\label{krullgravett} $|G|\ \le \prod_{s\in S}|G(s)|$.
\end{lemma}
\begin{proof} Let $s\in S$, and let $\cal{B}_s$ be the set of all
open balls with valuation radius $s$. Each $B\in \cal{B}_s$ is contained 
in a unique closed ball $\bar{B}$ of valuation radius $s$, and
$\bar{B}$ is the disjoint union of exactly $|G(s)|$-many $E\in \cal{B}_s$.
Thus we have a map $f_s\colon \cal{B}_s\to G(s)$ such that for all 
$D,E\in  \cal{B}_s$, if
$\bar{D}=\bar{E}$ and $D\ne E$, then  $f_s(D)\ne f_s(E)$.
We now associate to each $a\in G$ the element $\tilde{a}$ of the cartesian 
product set $\prod_{s\in S}G(s)$ such that
$\tilde{a}(s)=f_s\big(B_a(s)\big)$ for all $s\in S$. Then the map 
$$a\mapsto \tilde{a}\ \colon\ G \to \prod_{s\in S}G(s)$$
is injective: if $a,b\in G$ and $a\ne b$, then for $s:= v(a-b)$
we have $\tilde{a}(s) \ne \tilde{b}(s)$ since~$a$ and~$b$ are 
in the same closed 
ball of valuation radius $s$ but in different open balls of valuation 
radius $s$. 
\end{proof}   

\begin{example}[Hahn products]
Let $(G_s)_{s\in S}$ be a family of nontrivial abelian groups 
indexed by an ordered set $S$. For any element $g=(g_s)$ of the 
product group $\prod_sG_s$ we define its support by \index{support!Hahn product}\nomenclature[G]{$\supp g$}{support of an element $g\in\prod_s G_s$}
$$ \supp g\ :=\ \{s\in S:\ g_s\ne 0\}.$$
We define the {\bf Hahn product} \index{Hahn product} of the family $(G_s)$ to be the 
subgroup $H\big[(G_s)\big]$ of~$\prod_s G_s$ consisting of the $g\in \prod_s G_s$ with
well-ordered support. 
The (surjective) valuation \index{Hahn product!valuation} on this Hahn product $G:=H\big[(G_s)\big]$ \nomenclature[G]{$H[(G_s)]$}{Hahn product of the family $(G_s)$} given by
$$v\ \colon\ G\to S_{\infty}, \qquad v(g)\ :=\ \min (\supp g)\text{ for $g\ne 0$,}$$ 
is called the {\bf Hahn valuation\/} \index{valuation!Hahn} of $G$. Let $G$ be 
equipped with its Hahn valuation. 
Given $\sigma\in S$ we have 
the obvious projection map 
$$ (g_s) \mapsto g_{\sigma}\ \colon\ G \to G_{\sigma},$$
which restricts to a surjective group morphism 
$\bar{B}(\sigma) \to  G_{\sigma}$ with kernel
$B(\sigma)$, thus inducing a group isomorphism $G(\sigma) \to  G_{\sigma}$.
If all $G_s$ are equal, say $G_s=A$ for each $s\in S$, where $A$ is an abelian group, then we set 
$H[S,A]\ :=\ H\big[(G_s)\big]$. \nomenclature[G]{$H[S,A]$}{Hahn product $H[(G_s)]$ where $G_s=A$ for each $s\in S$}
\end{example}

\noindent
We define a {\bf valued abelian group\/} to be a triple $(G,S,v)$ with
$G$ an abelian group, $S$ an ordered set, and $v\colon G \to S_{\infty}$ a 
surjective valuation. 
We call the ordered set $S$ the {\bf value set} of the valued abelian group $(G,S,v)$.
A Hahn product is considered as a valued abelian group by equipping 
it with its Hahn valuation. When an ambient valued abelian group 
$(G, S, v)$
is given, then we set for $x,y\in G$, 
\begin{align*} x\preceq y\ &:\Longleftrightarrow\ v(x)\ge v(y), \qquad x\prec y\ :\Longleftrightarrow\ v(x)> v(y),\\
x\succeq y\ &:\Longleftrightarrow\ v(x)\le v(y), \qquad x\succ y\ :\Longleftrightarrow\ v(x)< v(y),\\
x\asymp y\ &:\Longleftrightarrow\ v(x)= v(y), \qquad x\sim y\ :\Longleftrightarrow\ v(x-y)> v(x).
\end{align*}
These relational notations are shorter and often more suggestive than notations using~$v$. In particular, $\sim$ is an equivalence relation on $G^{\ne}$.  

\index{dominance relation}

\medskip
\noindent
Let  $(G, S, v)$ and $(G', S', v')$ be
valued abelian groups. Then we say that $(G, S, v)$ is a 
{\bf valued subgroup} of $(G', S', v')$, or $(G', S', v')$ {\bf extends} 
$(G, S, v)$, or $(G', S', v')$ is an {\bf extension} of 
$(G, S, v)$, if
$G$ is a subgroup of $G'$, $S$ is an ordered subset of $S'$, and $v(x)=v'(x)$
for all $x\in G$; notation: $(G, S, v) \subseteq (G', S', v')$.

\index{valued abelian group}
\index{group!valued abelian}
\index{abelian group!valued}
\index{extension!valued abelian groups}
\nomenclature[G]{$(G,S,v)$}{valued abelian group}

Assume $(G, S, v) \subseteq (G', S', v')$.  Then we have for each 
$s\in S$ a natural
group embedding $G(s) \to G'(s)$ sending, 
for $a\in \bar{B}(s)$, the
coset $a+ B(s)$ of the open ball~$B(s)$ in $G$ to the coset
$a+ B'(s)$ of the open
ball $B'(s)=\{x\in G': v'(x)>s\}$ in $G'$. 
This extension of valued abelian groups is said to be {\bf immediate\/} if
$S= S'$ and for each $s\in S$ the group embedding 
$G(s) \to G'(s)$ is bijective; equivalently, for each $0\neq x'\in G'$ there is
$x\in G$ with $x\sim x'$.
For example, if $G$ is dense in $G'$ in the $v'$-topology, then $(G',S',v')$ is an immediate extension of $(G,S,v)$.
Such immediate extensions, in the setting of
asymptotic differential fields, will play an 
important role in what follows. A tool for coming to grips with
immediate extensions is the notion of pseudoconvergence, to which we turn
in the next subsection. At this stage we can make one useful definition:
Call a valued abelian group $(G,S,v)$ {\bf maximal\/} if it has no 
proper immediate valued abelian group extension. Hahn products are maximal: 
this is Corollary~\ref{hahnmax} below.

\index{valued abelian group!maximal}
\index{maximal!valued abelian group}

\begin{cor} 
Every valued abelian group has an immediate
valued abelian group extension that is maximal.
\end{cor}
\begin{proof} Use Zorn and Lemma~\ref{krullgravett}.
\end{proof} 

\noindent
By imposing extra structure on our valued abelian groups, we
can add to this existence result a corresponding uniqueness property;
see Corollaries~\ref{eqmaxsph} and~\ref{isomaximm}.

\subsection*{Pseudoconvergence}  \index{pseudoconvergence}
Fix a valued abelian group $(G,S,v)$. Let $(a_\rho)$ 
be a  well-indexed sequence in $G$, and $a\in G$.
Then $(a_\rho)$ is said to {\bf pseudoconverge to $a$} 
(notation: $a_\rho \leadsto a$), \nomenclature[G]{$a_\rho\leadsto a$}{$(a_\rho)$ pseudoconverges to $a$} if $\big(v(a-a_\rho)\big)$ is eventually strictly increasing, that is,
for some index $\rho_0$ we have $a-a_\sigma \prec a-a_\rho$ whenever
$\sigma > \rho > \rho_0$. We also say in that case that {\bf $a$ is a 
pseudolimit of $(a_\rho)$}. \index{pseudolimit} Note that if  $a_\rho \leadsto a$, then $a_\rho+ b \leadsto a+b$ 
for each~$b\in G$, and that
$$a_\rho \leadsto 0\ \Longleftrightarrow\   
\text{$(va_\rho)$ is eventually strictly increasing.}$$

\begin{lemma}\label{pc1} 
Let $(a_\rho)$ be a well-indexed sequence in $G$ such that 
$a_\rho \leadsto a$ where ${a\in G}$. With $s_\rho:= v(a-a_\rho)$, 
we have: \begin{enumerate}
\item[\textup{(i)}] either $a \prec a_\rho$ eventually, or $a\sim a_\rho$ eventually;
\item[\textup{(ii)}] $(va_\rho)$ is either eventually strictly increasing, or
eventually constant;
\item[\textup{(iii)}] for each $b\in G$:\quad 
$a_\rho \leadsto b\ \Longleftrightarrow\ v(a-b) > s_\rho$ 
eventually.
\end{enumerate}
\end{lemma}
\begin{proof} Let $\rho_0$ be as in the definition of  
``$a_\rho \leadsto a$.'' Suppose $a_\rho\preceq a$, where 
$\rho > \rho_0$. Then for $\sigma > \rho$ we have
$a-a_{\sigma}\prec a-a_\rho \preceq a$, so 
$a\sim a_\sigma$. This proves~(i). 
 Now~(ii) follows from (i) by noting that if $a\prec a_\rho$,  then $va_\rho=v(a-a_\rho)$, and if $a\sim a_{\rho}$, then
 $va_\rho=va$. We leave (iii) to the reader.
\end{proof}

\noindent 
If $(a_\rho)$ is a well-indexed sequence in $G$
and $a\in G'$ where $(G',S',v')$ is a valued abelian group extending
$(G,S,v)$, then
``$a_\rho \leadsto a$'' is to be interpreted in this valued
extension, that is, by considering the sequence $(a_\rho)$ in
$G$ as a sequence in $G'$.

\begin{lemma}\label{pc2} Suppose $(G',S',v')$ 
is an immediate valued abelian group 
extension of~$(G,S,v)$, and let $a\in G'\setminus G$. 
Then there is a well-indexed
sequence $(a_\rho)$ in $G$ such that $a_\rho \leadsto a$ and $(a_\rho)$
has no pseudolimit in $G$. 
\end{lemma}
\begin{proof} We claim that the subset 
$\big\{v'(a-x): x\in G\big\}$ of $S$ 
has no largest element. To see this, let $x\in G$; we shall find $y\in G$ 
such that $v'(a-y) > v'(a-x)$. We have $s:= v'(a-x)\in S$ and $G(s)=G'(s)$,
so we can take $b\in G$ with $a-x\in b+B'(s)$, hence $v'(a-y)>s$ for $y:=x+b$,
as claimed. It follows that we can take a well-indexed sequence  
$(a_\rho)$ in $G$
such that the sequence $\big(v'(a-a_\rho)\big)$ is strictly increasing and cofinal
in $\big\{v'(a-x): x\in G\big\}$. Thus $a_\rho \leadsto a$. If
$a_\rho\leadsto g\in G$, then $v'(a-g) > v'(a-a_\rho)$ 
for all $\rho$ by Lemma~\ref{pc1}(iii), hence
$v'(a-g) > v'(a-x)$ for all $x\in G$, a contradiction.   
\end{proof}

\subsection*{Pseudocauchy sequences} As before, fix a valued abelian 
group $(G,S,v)$.
To capture within $(G,S,v)$ that a well-indexed sequence
$(a_\rho)$ in $G$ has a pseudolimit in some valued abelian group 
extension we make
the following definition. 

A {\bf pseudocauchy sequence in $G$\/}, 
more precisely, in $(G,S,v)$, is a well-indexed sequence
$(a_\rho)$ in $G$ such that
for some index $\rho_0$ we have 
$$\tau > \sigma > \rho > \rho_0\ \Longrightarrow\ 
a_\tau-a_\sigma \prec a_\sigma-a_\rho.$$ We also write \textit{pc-sequence}\/ for
\textit{pseudocauchy sequence.}\/ A cofinal subsequence of a pc-sequence
in $G$ is a pc-sequence in $G$.  

\index{pseudocauchy sequence}
\index{sequence!pseudocauchy}
\index{pc-sequence}

\begin{lemma}\label{pc3} 
Let $(a_\rho)$ be a well-indexed sequence in $G$. Then
 $(a_\rho)$ is a pc-sequence in $G$ if and only if 
 $(a_\rho)$ has a pseudolimit in some valued abelian group 
extension of~$(G,S,v)$. In that case, $(a_\rho)$ has even a pseudolimit 
in some elementary extension of the two-sorted structure $(G,S,v)$. \rm{(See \ref{sec:eleqelsub} for \textit{elementary extension.}\/)}
\end{lemma}
\begin{proof} Suppose $(a_\rho)$ is a pc-sequence in $G$, and
let $\rho_0$ be as in the definition of  \textit{pseudocauchy sequence.}\/
We refer to \ref{sec:sat} for the notion of type. Consider the type in the variable $x$
consisting of the formulas 
$$ x-a_\sigma \prec x-a_\rho \qquad(\sigma > \rho > \rho_0).$$ 
Every finite subset of this type
is realized by $a_\tau$ for any sufficiently large $\tau$. Thus
we can realize this  type 
by a suitable $a\in G'$ for some elementary extension  
$(G', S',v')$
of~$(G,S,v)$, and then~$a_\rho\leadsto a$.

For the converse, suppose $a_\rho\leadsto a$ where $a\in G'$ and  
$(G',S',v')$ is a valued abelian group extension of $(G, S,v)$.
Let $\rho_0$ be as in the definition of pseudolimit, and let
$\sigma > \rho > \rho_0$. Then $a_\sigma-a_\rho = (a_\sigma-a) -(a_\rho -a)$,
so $a_\sigma-a_\rho\asymp a-a_\rho$. So if in addition $\tau > \sigma$, then
$a_\tau - a_\sigma\asymp a-a_\sigma \prec a-a_\rho \asymp a_\sigma-a_\rho$.   
\end{proof} 

\noindent
It will be relevant to us that the last part of Lemma~\ref{pc3} goes through 
with the same proof and the corresponding notion of elementary extension
when $(G,S,v)$ has extra (first-order) structure.
For example, we shall apply it to valued differential fields: 
fields equipped with both a valuation and a derivation, as defined later.

\begin{cor}\label{impca} If every pc-sequence in $G$ has a pseudolimit in $G$, 
then the valued abelian group $(G,S,v)$ is maximal. 
\end{cor}
\begin{proof} Immediate from Lemma~\ref{pc2}.
\end{proof} 

\noindent
The converse of Corollary~\ref{impca} holds when our valued abelian group 
has suitable extra structure; see Corollary~\ref{eqmaxsph}.    
It is easy to check that every pc-sequence in a Hahn product pseudoconverges
in it. Thus by Corollary~\ref{impca}:

\begin{cor}\label{hahnmax} Any Hahn product is maximal as a 
valued abelian group.
\end{cor}

\noindent
Lemma~\ref{pc3} and part~(ii) of Lemma~\ref{pc1} yield:

\begin{cor} \label{pc4} 
If $(a_\rho)$ is a pc-sequence in $G$, 
then $(va_\rho)$ is either eventually strictly increasing 
\textup{(}so $a_{\rho} \leadsto 0$\textup{)}, or
eventually constant \textup{(}so $a_{\rho} \not\leadsto 0$\textup{)}.
\end{cor}

\noindent
Let $(a_\rho)$ be a pc-sequence in $G$, pick $\rho_0$ 
as above, and put 
$$s_\rho :=v(a_{\rho'}-a_\rho)\in S\qquad\text{for
$\rho'>\rho > \rho_0$;}$$ 
this depends only on $\rho$ as the notation
suggests. Then $(s_\rho)_{\rho > \rho_0}$ is strictly increasing.
% and $s_\rho=v(a_{\rho+1}-a_\rho)$, with $\rho+1$
%the immediate successor of $\rho$. We call the sequence $(v(a_{\rho+1}-a_\rho))%$ in $S_\infty$ the {\bf gauge} of $(a_\rho)$.

\begin{lemma}\label{eqpslim} Let $a\in G$. Then the following are equivalent:
\begin{enumerate} 
\item[\textup{(i)}] $a_\rho\leadsto a$;
\item[\textup{(ii)}] $v(a-a_\rho)=s_\rho$ for all $\rho>\rho_0$;
\item[\textup{(iii)}] $v(a-a_\rho)\geq s_\rho$ for all $\rho>\rho_0$;
\item[\textup{(iv)}] $v(a-a_\rho)\geq s_\rho$ eventually.
\end{enumerate}
\end{lemma}
\begin{proof}
Suppose $a_\rho\leadsto a$, and let $\rho>\rho_0$. Since $v(a-a_\sigma)$ is 
eventually strictly increasing, $v(a-a_\sigma)\notin\{ v(a-a_\rho),s_\rho \}$ for sufficiently large $\sigma>\rho$; for such $\sigma$,
$$s_\rho=v(a_\sigma-a_\rho)=\min\big(v(a-a_\rho),v(a-a_\sigma)\big)=v(a-a_\rho).$$ 
This shows (i)~$\Rightarrow$~(ii), and (ii)~$\Rightarrow$~(iii)~$\Rightarrow$~(iv) are trivial. For (iv)~$\Rightarrow$~(i), suppose $\rho_1\geq\rho_0$ is such that $v(a-a_\rho)\geq s_\rho$ for all $\rho>\rho_1$.
It suffices to show that then $v(a-a_\rho)=s_\rho$ for $\rho>\rho_1$. Suppose towards a contradiction that $\rho>\rho_1$ is such that $v(a-a_\rho)>s_\rho$. Take any $\sigma>\rho$; then $v(a-a_{\sigma})> s_{\rho}$, so
$$  s_\rho\ =\ v(a_{\sigma}-a_{\rho})\ =\ v\big((a-a_\rho)-(a-a_\sigma)\big)\ >\ s_\rho,$$
a contradiction.
\end{proof}

\noindent
The {\bf width} of $(a_\rho)$ is an upward closed subset of $S_\infty$, namely 
$$  \{s \in S_{\infty}:\ s>s_\rho\ \mbox{for all}\ \rho > \rho_0\}\  \qquad(\text{independent of the choice of $\rho_0$}).$$    
Its significance is that if $a,b\in G$ and $a_\rho \leadsto a$, then by 
Lemma~\ref{pc1}, 
$$a_\rho \leadsto b\ \Longleftrightarrow\ \text{$v(a-b)$ is in the width of $(a_\rho)$.}$$

\index{width of a pc-sequence}
\index{pc-sequence!width}

\subsection*{Fixpoint theorem} 
As before, $(G,S,v)$ is a valued abelian group.
To visualize the property that every pc-sequence has a pseudolimit, 
define a {\bf nest of balls\/} in~$G$ to be a
collection of balls in~$G$ any two of which meet. So a nest of balls in $G$ is 
(totally) ordered by inclusion.
We call~$(G,S,v)$
{\bf spherically complete\/} if every nonempty nest of closed balls in $G$ 
has a point in its intersection. (For example, if~$S$ is well-ordered under the reversed ordering, then trivially $(G, S,v)$ is spherically complete.)

\label{p:sph complete}
\index{nest!balls}
\index{valued abelian group!spherically complete}
\index{spherically complete}

\begin{lemma}\label{fp1} $(G, S,v)$ is spherically complete if and only if 
each pc-sequence in~$G$ has a pseudolimit in $G$.
\end{lemma}
\begin{proof}
Let $(a_\rho)$ be a pc-sequence in $G$, and let $\rho_0$ be as in the 
definition of pc-sequence.
Then $\mathcal B=\big\{\overline{B}_{a_\rho}(s_\rho):\ \rho >\rho_0 \big\}$ 
is a nonempty nest of closed balls in~$G$, and 
$\bigcap\mathcal B=\{a\in G:a_\rho\leadsto a\}$ by 
Lemma~\ref{eqpslim}.
Thus if $(G, S,v)$ is spherically complete, then
each pc-sequence in $G$ has a pseudolimit in $G$. Conversely, 
suppose the latter condition holds,
and let $\mathcal B$ be a nonempty nest of closed balls in $G$; we need to 
show that $\bigcap\mathcal B\ne \emptyset$. 
Eliminating a trivial case,
assume that $\mathcal{B}$ has no smallest ball in it. 
Replacing~$\mathcal B$  by a coinitial subset (under inclusion), we
arrange that $\mathcal{B}=\{B_{\rho}:\ \rho< \lambda\}$ for some infinite 
limit ordinal $\lambda$, with $B_{\rho}$ strictly containing~$B_{\sigma}$ whenever $\rho< \sigma < \lambda$. For each~$\rho<\lambda$ 
we take $a_{\rho}\in B_{\rho}\setminus B_{\rho +1}$. It follows that
$(a_{\rho})$ is a pc-sequence. Take~$a\in G$ such that 
$a_{\rho} \leadsto a$. Then $a\in \bigcap\mathcal B$.
\end{proof}

\noindent
A routine argument shows:

\begin{lemma}\label{fp2} If $X\subseteq G$ and 
$B_x(s)\subseteq X$ whenever $x,y\in X,\ x\ne y,\ s=v(x-y)$, 
then any pseudolimit in $G$ of any pc-sequence in $X$ lies in $X$.
\end{lemma}

\noindent
For example, any ball (open or closed) in $G$ satisfies the hypothesis of
Lemma~\ref{fp2}, and so does, trivially, $X=G$.  
A map $f\colon X\to X$ with $X\subseteq G$ is said to be 
{\bf contractive} if for all 
distinct $x,y\in X$ we have $f(x)-f(y) \prec x-y$.

\index{map!contractive}
\index{contractive}

\begin{theorem} \label{thm:fixpoint}
Let $f\colon X\to X$ be a contractive map, where
$\emptyset \ne X\subseteq G$ and each pc-sequence in $X$ has a pseudolimit 
in $X$. Then $f$ has a unique fixpoint.
\end{theorem} 
\begin{proof} It is clear that $f$ has at most one fixpoint. Take any point 
$x_0\in X$ and make it the initial term of a sequence
$(x_{\lambda})_{\lambda< \nu}$ in $X$ indexed by the ordinals less than an ordinal~$\nu>0$, such that \begin{enumerate}
\item [(1)] $x_{\lambda} \ne x_{\mu}$ whenever $\lambda < \mu < \nu$;
\item [(2)]$x_{\lambda+1}=f(x_{\lambda})$ whenever $\lambda < \lambda +1< \nu$;
\item[(3)] $x_{\lambda''}-x_{\lambda'} \prec x_{\lambda'}-x_{\lambda}$ whenever $\lambda < \lambda' < \lambda'' < \nu$.
\end{enumerate}
If $\nu=\mu+1$ is a successor ordinal, and $x_{\mu}$ is not yet a fixpoint of $f$, then we set $x_{\nu}:= f(x_{\mu})$, and then the extended sequence
$(x_{\lambda})_{\lambda < \nu +1}$ satisfies (1)--(3) with $\nu+1$
instead of $\nu$. If $\nu$ is a limit ordinal, then we let $x_{\nu}$ be a pseudolimit in $X$ of the pc-sequence $(x_{\lambda})_{\lambda< \nu}$, and then the extended sequence $(x_{\lambda})_{\lambda < \nu +1}$ satisfies again (1)--(3) with $\nu+1$ instead of $\nu$. This building process must come to a halt 
by producing a fixpoint of $f$.
\end{proof}

\begin{cor}\label{fixcor} Suppose $(G, S,v)$ is spherically complete and
$e\colon G \to G$ is a contractive map. Then $\operatorname{id}_G + e\colon G \to G$
is bijective. 
\end{cor}

\begin{proof} Injectivity is immediate from $e$ being contractive.
To get sur\-jec\-ti\-vi\-ty,
let ${a\in G}$, and define $f\colon G \to G$ by $f(x)=a-e(x)$. Then
$f$ is contractive. Let~$x$ be the fixpoint of $f$. Then $a=x+e(x)$.
\end{proof}

\noindent
The variant below is also useful. Let $\emptyset\ne P\subseteq S$, 
$\emptyset\ne X\subseteq G$, and
let a map $f\colon X\to X$ be given. A {\bf $P$-fixpoint} of $f$ is a point
$x\in X$ such that $v(f(x)-x) \ge s$ for some $s\in P$. We call $f$ 
{\bf contractive up to $P$} if for all 
$x,y\in X$ with $v(x-y) < P$ we have $f(x)-f(y)\prec x-y$.

\index{map!contractive up to $P$}
\index{contractive!up to $P$}
\index{map!$P$-fixpoint}

\begin{lemma} Suppose every pc-sequence in $X$ has a pseudolimit in $X$, and
$f$ is contractive up to $P$. Then $f$ has a $P$-fixpoint. 
\end{lemma}
\begin{proof} Take any point $x_0\in X$ and make it the initial term 
of a sequence
$(x_{\lambda})_{\lambda< \nu}$ in~$X$ indexed by the ordinals less than an ordinal 
$\nu>0$, such that \begin{enumerate}
\item[(1)] $v(x_{\mu}-x_{\lambda})< P$ whenever $\lambda < \mu < \nu$;
\item[(2)] $x_{\lambda+1}=f(x_{\lambda})$ whenever $\lambda < \lambda +1< \nu$;
\item[(3)] $x_{\lambda''}-x_{\lambda'}\prec x_{\lambda'}-x_{\lambda}$ whenever $\lambda < \lambda' < \lambda'' < \nu$.
\end{enumerate}
We now continue as in the proof of Theorem~\ref{thm:fixpoint}, 
replacing ``fixpoint'' by ``$P$-fixpoint'' throughout.
\end{proof}

\subsection*{Generalizing closed balls} Let $(G, S, v)$ 
be a valued abelian group. We define a {\bf union-closed ball} in $G$ to be the union of a nonempty nest of closed balls in~$G$. We 
can assume the closed balls in such a nest to have a common
center: Let $\{B_i:\ i\in I\}$ be a nest of closed balls in $G$
where~$I$ is a nonempty index set and $B_i=\overline{B}_{a_i}(s_i)$, with $a_i\in G$ and 
$s_i\in S$ for $i\in I$; then 
$\bigcup_i B_i = \bigcup_i \overline{B}_{a}(s_i)$ for any $a\in \bigcup_i B_i$. Using this fact, one easily shows:

\index{ball!union-closed}
\index{union-closed}

\begin{lemma}\label{genclosedballs} Let $B_1$ and $B_2$ be 
union-closed balls in $G$
such that $B_1\cap B_2\ne \emptyset$. Then $B_1\subseteq B_2$ or $B_2\subseteq B_1$. If in addition $B_1$ is properly contained in $B_2$, then there is a closed
ball $B$ in $G$ such that $B_1\subseteq B \subseteq B_2$.
\end{lemma}

\noindent
Note that if $S$ has no largest element, then every open ball
in $G$ is a union-closed ball.
A {\bf nest of union-closed balls} in $G$ is a collection of
union-closed balls in $G$ any two of which meet. Such a
nest is (totally) ordered by inclusion.

\index{nest!union-closed balls}

\begin{cor}\label{corgenclosedballs} If $(G,S, v)$ is spherically complete, then any
nonempty nest of union-closed balls in $G$ has a point in its
intersection.
\end{cor}
\begin{proof} Let $(B_i)$ be a family of union-closed balls in $G$ indexed by a nonempty well-ordered set
$I$ such that $I$ has no largest element and
$B_i$ properly contains $B_j$ whenever $i < j$ in $I$.
For each $i$, let $s(i)$ be the immediate successor of $i$ in $I$, and
use Lemma~\ref{genclosedballs} to get a closed ball $D_i$ in $G$ such that $B_{s(i)} \subseteq D_i \subseteq B_i$. Then $\{D_i:\ i\in I\}$ is a nest of closed balls,
and $\bigcap_i D_i = \bigcap_i B_i$. 

The corollary is an easy consequence of this construction. 
\end {proof}

\subsection*{Equivalence of pc-sequences} Let $(G,S,v)$ be a valued abelian group.

\begin{lemma}\label{pc5} Given pc-sequences $(a_{\rho})$ and $(b_{\sigma})$ in $G$
\textup{(}with possibly different index sets\textup{)}, the
following conditions are equivalent: \begin{enumerate}
\item[\textup{(i)}] $(a_{\rho})$ and $(b_{\sigma})$ have the same pseudolimits in every
valued abelian group extension of $(G,S,v)$;
\item[\textup{(ii)}] $(a_{\rho})$ and $(b_{\sigma})$ have the same width, and have a 
common pseudolimit in some valued abelian group extension of $(G,S,v)$;
\item[\textup{(iii)}] there are arbitrarily large $\rho$ and $\sigma$ such that 
for all $\rho'> \rho$ and $\sigma'> \sigma$ we have 
$a_{\rho'}-b_{\sigma'}\prec a_{\rho'}-a_{\rho}$, and there are arbitrarily large 
$\rho$ and $\sigma$ such that for all $\rho'> \rho$ and $\sigma'> \sigma$ we 
have $a_{\rho'}-b_{\sigma'}\prec b_{\sigma'}-b_{\sigma}$.
\end{enumerate}
\end{lemma}
\begin{proof} Suppose $s\in S$ is in the width of $(a_{\rho})$ but not in
the width of $(b_{\sigma})$. Take $x\in G$ with $vx=s$ and take $a$ in some
valued abelian group extension of $(G,S,v)$ such that $a_{\rho} \leadsto a$.
Then also $a_\rho \leadsto a+x$, but $a$ and $a+x$ cannot both be pseudolimits 
of $(b_{\sigma})$. This argument proves (i)~$\Rightarrow$~(ii).
To prove (ii)~$\Rightarrow$~(iii), assume~(ii), and take~$a$ in some
valued abelian group extension of $(G,S,v)$ such that $a_{\rho} \leadsto a$
and $b_{\sigma} \leadsto a$. Let~$\rho_0$ be so large that
$v(a-a_\rho)$ is strictly increasing for $\rho > \rho_0$ and 
$a-a_\rho\asymp a_{\rho'}-a_{\rho}$ for all $\rho'>\rho > \rho_0$. 
Likewise, let $\sigma_0$ be so
large that 
$v(a-b_\sigma)$ is strictly increasing for $\sigma > \sigma_0$, and
$a-b_{\sigma}\asymp b_{\sigma'}-b_\sigma$ for all $\sigma'>\sigma > \sigma_0$.
Take any $\rho > \rho_0$ and then take $\sigma> \sigma_0$ such that 
$a-b_\sigma \prec a-a_{\rho}$. Then for $\rho'>\rho$ and 
$\sigma'>\sigma$ we have
$a_{\rho'}-b_{\sigma'}=(a_{\rho'}-a)+(a-b_{\sigma'})$, with
$$a_{\rho'}-a\prec a-a_{\rho}\asymp a_{\rho'}-a_{\rho}, \qquad a-b_{\sigma'}\prec a-b_{\sigma}\prec a_{\rho'}-a_{\rho},$$
so $a_{\rho'}-b_{\sigma'}\prec a_{\rho'}-a_{\rho}$. Likewise, we obtain the 
second part of (iii).
 
For (iii)~$\Rightarrow$~(i), assume (iii), and let $a$ be a pseudolimit of 
$(a_\rho)$ in a 
valued abelian group extension of $(G,S,v)$. Take $\rho_0$ and $\sigma_0$ 
such that
$v(a-a_\rho)$ is strictly increasing for $\rho > \rho_0$, and 
$s_{\sigma}:=v(b_{\sigma'}-b_{\sigma})$ for $\sigma' >\sigma> \sigma_0$ depends only
on $\sigma$ (not on $\sigma'$) and is strictly increasing as a function of $\sigma>\sigma_0$.
Let any $\sigma > \sigma_0$ be given. Then
$a-b_{\sigma}=(a-a_{\rho}) + (a_{\rho}-b_{\sigma})$ and 
$$a_{\rho}-b_{\sigma}=(a_{\rho}-b_{\sigma'})+(b_{\sigma'}-b_{\sigma})\asymp b_{\sigma'}-b_{\sigma}, \text{ and }v(b_{\sigma'}-b_{\sigma})=s_{\sigma}$$
for all sufficiently large $\rho> \rho_0$ and $\sigma'>\sigma$.
As $v(a-a_{\rho})$ is strictly increasing for~$\rho > \rho_0$, and
$v(a-b_{\sigma})$ and $v(a_{\rho}-b_{\sigma})=s_{\sigma}$ depend only on $\sigma$
(for big enough $\rho$)
it follows that
$v(a-b_{\sigma})= s_{\sigma}$. Thus $b_{\sigma} \leadsto a$.
\end{proof}

\begin{remark}
Lemma~\ref{pc5} goes through with the same proof when the valued abelian group $(G,S,v)$ 
is equipped with extra first-order structure and in (i) and (ii) 
the extensions are required to be elementary extensions of the 
expanded structure.
\end{remark}

\noindent
Two pc-sequences $(a_{\rho})$ and $(b_{\sigma})$ in $G$ are said to be 
{\bf equivalent\/} if the conditions of Lemma~\ref{pc5} are satisfied.
This equivalence relation on the class of pc-sequences in~$G$
is relative to the ambient valued abelian group
$(G,S,v)$. However, condition~(iii) of Lemma~\ref{pc5} shows that for pc-sequences
$(a_{\rho})$ and $(b_{\sigma})$ in $G$ and any valued abelian group extension
$(G',S',v')$ of $(G,S,v)$ we have: $(a_{\rho})$ and $(b_{\sigma})$ are 
equivalent with respect to $(G,S,v)$ iff  $(a_{\rho})$ and $(b_{\sigma})$ are
equivalent with respect to $(G',S',v')$. Note that any cofinal subsequence
of a pc-sequence $(a_{\rho})$ in $G$ is equivalent to $(a_{\rho})$.

\index{equivalence!pc-sequences}
\index{pc-sequence!equivalence}

\medskip\noindent
Let $a$ be an element of a valued abelian group extension of $(G,S,v)$
with $a\notin G$. For convenience, denote the valuation of that extension also
by $v$, and set $$v(a-G):= \big\{v(a-g):\ g\in G\big\},$$
a nonempty subset of the value set of that extension. The next lemmas collect some basic facts about this situation. The easy proofs are left to the reader.

\begin{lemma}\label{pcmax} The following are equivalent:
\begin{enumerate}
\item[\textup{(i)}] $v(a-G)$ has no largest element;
\item[\textup{(ii)}] $a_{\rho}\leadsto a$ for some pc-sequence $(a_{\rho})$ in $G$ without
pseudolimit in $G$.
\end{enumerate}
If $a$ lies in an immediate valued abelian group extension of $(G,S,v)$, then \textup{(i)} holds. If~\textup{(i)} holds, then $v(a-G)$ is a downward closed subset of $S$.
\end{lemma}

\begin{lemma}\label{pcnonmax} Let $(a_{\rho})$ be a well-indexed sequence in $G$ such that $v(a-a_{\rho})$ is strictly increasing as a function of $\rho$. Then $a_{\rho} \leadsto a$, and we have:
$$\big(v(a-a_{\rho})\big) \text{ is cofinal in $v(a-G)$}\ \Longleftrightarrow\
\text{$(a_{\rho})$  has no pseudolimit in $G$.}$$
\end{lemma}

\noindent
A {\bf divergent\/} pc-sequence in $G$ is a pc-sequence in $G$ without a pseudolimit in $G$. 

\index{pc-sequence!divergent}

\begin{cor}\label{eqpccor}
Any two divergent pc-sequences in $G$ with a common 
pseudolimit in an immediate valued abelian group extension of $(G,S,v)$ 
are equivalent.
\end{cor}

\subsection*{Coarsenings} Let surjective valuations $v\colon G \to S_{\infty}$ and 
$v'\colon G \to S'_\infty$ on the abelian group $G$ be given.
We say that $v'$ is a {\bf coarsening} of $v$ 
(or $v'$ is {\bf coarser} than~$v$, or $v$ is {\bf finer} than~$v'$) if for 
all $a,b\in G$ we have $v(a)\leq v(b)\Rightarrow v'(a)\leq v'(b)$; 
equivalently, there exists an increasing surjection 
$i\colon S_{\infty}\to S'_{\infty}$ such that $v'=i\circ v$. Note that for such a map~$i$ we have $i(\infty)=\infty$ and $i(S)=S'$, since
$v'g\ne \infty$ for $g\in G^{\ne}$. 
We call $v$ and $v'$ {\bf equivalent} if $v$ is coarser than $v'$ and
$v'$ is coarser than~$v$; equivalently, there 
exists an isomorphism $i\colon S_{\infty}\to S_{\infty}'$ of ordered sets such that $v'=i\circ v$ (and such $i$ is then uniquely determined).
If $v$ and $v'$ are equivalent, then the $v$-topology agrees with 
the $v'$-topology. If $v'$ is coarser than $v$ and 
$S'$ has no largest element, then the $v$-topology also agrees with the
$v'$-topology. 

\index{coarsening!valuation on an abelian group}
\index{coarser}
\index{finer}
\index{equivalence!valuations on an abelian group}

\begin{lemma}\label{coarsepc} Suppose $v'$ is coarser than $v$, and $(a_{\rho})$ is a pc-sequence
in $G$ with respect to $v'$. Then  $(a_{\rho})$ is a pc-sequence with respect 
to $v$, and for all $a\in G$, 
$$\text{$a_{\rho} \leadsto a$ with respect to v}\ \Longleftrightarrow\ \text{$a_{\rho} \leadsto a$ with respect to $v'$.}$$
\end{lemma}
 
\noindent
To verify $\Rightarrow$, use for example Lemma~\ref{eqpslim}.

\begin{cor}\label{coco} If $(G, S, v)$ is spherically complete and 
$v'$ is coarser than $v$, then~$(G, S', v')$ is spherically complete.
\end{cor}  

\subsection*{Cauchy sequences} 
Let $(G,S,v)$ be a valued abelian group and  $(a_\rho)$ a well-indexed sequence in $G$. We say that $(a_\rho)$ is a {\bf cauchy sequence} (or a {\bf c-sequence}) in $G$ if for every $s\in S$ there is~$\rho_0$ such that $v(a_\rho-a_{\rho'})>s$ for all $\rho,\rho'>\rho_0$. 
Thus if $S$ has a largest element, then~$(a_{\rho})$ is a c-sequence iff $a_{\rho}$ is eventually constant.
A c-sequence $(a_\rho)$ in $G$ remains a c-sequence in every extension $(G',S',v')$ of $(G,S,v)$ with $S$ cofinal in $S'$. 
For~$a\in G$ we say that $(a_\rho)$ {\bf converges to $a$} if for each $s\in S$ there is some $\rho_0$ such that $v(a-a_\rho)>s$ for all $\rho>\rho_0$; in symbols: $a_\rho\to a$.
We say that $(a_\rho)$ {\bf converges in $G$} if $a_\rho\to a$ for some $a\in G$.
Note that if $(a_\rho)$ converges in some extension of $G$, then $(a_\rho)$ is a c-sequence.
There is at most one $a\in G$ with $a_\rho\to a$, and if there is such an $a$ we call it the {\bf limit of~$(a_\rho)$}.  
Clearly each pc-sequence of width $\{\infty\}$ is a c-sequence. (For a partial converse of this statement see Lemma~\ref{lem:c-sequence width} below.) 
If  $(a_\rho)$ is a pc-sequence of width~$\{\infty\}$ and $a\in G$, then we 
have the equivalence
$a_\rho\leadsto a\ \Longleftrightarrow\ a_\rho\to a$. The following is obvious.

\index{cauchy sequence!in a valued abelian group}
\index{c-sequence!in a valued abelian group}
\index{sequence!cauchy}
\index{sequence!convergence}
\index{sequence!limit}
\index{limit of a sequence}
\nomenclature[G]{$a_\rho\to a$}{$(a_\rho)$ converges to $a$}

\begin{lemma}\label{lem:addition of c-sequences}
Let $(a_\rho)$, $(b_\rho)$ be c-sequences in $G$ with the same index set, and~$a,b\in G$. Then

\begin{enumerate}
\item[\textup{(i)}] $(-a_\rho)$ is a c-sequence, and if $a_\rho\to a$, then $-a_\rho\to -a$;
\item[\textup{(ii)}] $(a_\rho+b_\rho)$ is a c-sequence, and if $a_\rho\to a$ and $b_\rho\to b$, then $a_\rho+b_\rho\to a+b$.
\end{enumerate}
\end{lemma}

\noindent
We also have an analogue of Corollary~\ref{pc4} for c-sequences:

\begin{lemma} \label{lem:pc4 analogue}
Let $(a_\rho)$ be a c-sequence in $G$. Then either
for every $s\in S$ there is a $\rho_0$ such that $va_\rho>s$ for all $\rho>\rho_0$ \textup{(}that is, $a_\rho\to 0$\textup{)}, or $va_\rho$ takes an eventually constant value in $S$ \textup{(}so $a_\rho\not\to 0$\textup{)}. 
\end{lemma}
\begin{proof}
Suppose $a_\rho\not\to 0$. Take $s\in S$ such that for every $\rho$ there is some $\rho'>\rho$ with $va_{\rho'}\leq s$. Also, $(a_\rho)$ being a c-sequence, take $\rho_0$ such that $v(a_\rho-a_{\rho'})>s$ for all $\rho,\rho'>\rho_0$. Let $\rho'>\rho>\rho_0$; we claim that $va_\rho=va_{\rho'}$. To see this take $\rho''>\rho'$ such that $va_{\rho''}\leq s$. Then
$v(a_\rho-a_{\rho''})>s$ and hence $va_\rho=va_{\rho''}$, and similarly $va_{\rho'}=va_{\rho''}$, whence the claim.
\end{proof}

\noindent
For every c-sequence $(a_\rho)$ in $G$ with $a_\rho\not\to 0$ there is an $s\in S$ with $va_\rho=s$ eventually, by Lemma~\ref{lem:pc4 analogue}; we call $s$ the eventual valuation of $(a_\rho)$, and if $a_\rho\to 0$ we define the eventual valuation of $(a_\rho)$ to be $\infty$. Note that if $a$ is an element in an extension of~$G$ with $a_\rho\to a$, then the eventual valuation of $(a_\rho)$ is $va$.

\medskip
\noindent
A cofinal subsequence $(b_\sigma)$ of a c-sequence $(a_\rho)$ is a c-sequence, with $a_\rho\to a$ iff $b_\sigma\to a$, for all $a\in G$. The next lemma allows us to restrict our attention to c-sequences indexed by the ordinals less than $\operatorname{cf}(S)$:

\begin{lemma} \label{lem:cofinality c-sequences}
Suppose $(a_\rho)$ is a c-sequence that is not eventually constant. Then the index set of $(a_\rho)$ has cofinality $\operatorname{cf}(S)$.
\end{lemma}
\begin{proof}
The assumption implies that $S\ne \emptyset$ and $S$ has no largest element.
Choose a sequence $(s_\gamma)$ in $S$, indexed by the ordinals 
$\gamma<\operatorname{cf}(S)$, which is strictly increasing and cofinal 
in $S$. For each $\gamma$, define
$$I_\gamma := \big\{ \rho:\ \text{$v(a_{\rho_1}-a_{\rho_2})>s_\gamma$ for all 
$\rho_1, \rho_2\ge\rho$} \big\},$$
so $I_\gamma\neq\emptyset$ and $I_{\gamma}$ is upward closed in the set of indices $\rho$.  Set $\rho(\gamma):=\min I_\gamma$. The well-indexed sequence 
$\big(\rho(\gamma)\big)$ is increasing, since $I_{\gamma'}\subseteq I_\gamma$ for $\gamma'\geq \gamma$. Also $\bigcap_\gamma I_\gamma=\emptyset$, since $(a_\rho)$ is 
not ultimately constant. Therefore, given any index $\rho$ we can choose  
$\gamma$ with $\rho\notin I_\gamma$, and then $\rho(\gamma) > \rho$.
Hence the sequence $\big(\rho(\gamma)\big)$ is cofinal in the set of indices $\rho$, and
thus by Corollary~\ref{cor:cofinality} the set of indices $\rho$ has  
cofinality $\operatorname{cf}\!\big(\!\operatorname{cf}(S)\big)=\operatorname{cf}(S)$.
\end{proof}

\noindent
We call an extension $(G',S',v')\supseteq (G,S,v)$ of valued abelian groups {\bf dense}\index{valued abelian group!dense extension}\index{dense extension!valued abelian groups}\index{extension!valued abelian groups!dense} if $G$ is dense in $G'$ (in the $v'$-topology on $G'$). 
Every dense extension of valued abelian groups is immediate.
Hence by Lemma~\ref{pc2} and its proof:

\begin{lemma}\label{lem:pc2 analogue}
Given a dense extension $(G',S',v')\supseteq (G,S,v)$ and $a\in G'\setminus G$,  there is a divergent pc-sequence $(a_\rho)$ of width~$\{\infty\}$ in $G$, 
indexed by the ordinals~$\rho<\operatorname{cf}(S)$, such that 
$a_\rho\leadsto a$.
\end{lemma}

\index{valued abelian group!complete}
\index{complete!valued abelian group}
\label{p:complete}

\noindent
One says that $G$ is {\bf complete} if every c-sequence in $G$ converges in $G$. Thus if~$S$ has a largest element, then
$G$ is automatically complete. If $v'$ is a coarsening of the valuation~$v$ of $G$ and the value set $S'$ of $v'$ has no largest element, then
Lemma~\ref{coarsepc} goes through with c-sequences instead of pc-sequences, and 
$a_{\rho}\to a$ instead of $a_{\rho} \leadsto a$. Therefore:

\begin{cor}\label{coarsecomplete} If $G$ is complete with respect to $v$, then $G$ is complete with respect to any coarsening of $v$.
\end{cor} 

\noindent
By Lemma~\ref{lem:pc2 analogue}, if $G$ is complete, then $G$ has no proper dense 
extension.
On the other hand, every c-sequence in $G$ converges in some dense extension 
of $G$:

\begin{theorem}\label{thm:completion}
Every valued abelian group has a dense complete extension. 
\end{theorem}

\begin{proof}
If $S=\emptyset$ (so $G=\{0\}$), or $S$ has a largest element, 
then $G$ is complete and~$G^{\operatorname{c}}:=G$ 
has the required properties. For the rest of the proof we 
assume $S$ is nonempty and has no largest element
(so $\operatorname{cf}(S)$ is infinite).
Let $G^{\operatorname{cs}}$ be the set of all c-sequences~$(a_\rho)$ in $G$ indexed by the ordinals $\rho<\operatorname{cf}(S)$. By  Lemma~\ref{lem:addition of c-sequences}, $G^{\operatorname{cs}}$ is an abelian group under componentwise addition of sequences, and
$$N := \big\{ (a_\rho)\in G^{\operatorname{cs}} : a_\rho\to 0 \big\}$$
is a subgroup of $G^{\operatorname{cs}}$; we let $G^{\operatorname{c}}:=G^{\operatorname{cs}}/N$ be the quotient group.
For $(a_\rho),(b_\rho)\in G^{\operatorname{cs}}$ with $b_\rho\to 0$ the eventual valuations of $(a_\rho)$ and $(a_\rho+b_\rho)$ are the same, hence we can define a function
$$v^{\operatorname{c}}\colon G^{\operatorname{c}}\to S_\infty,\qquad  v^{\operatorname{c}}\big((a_\rho)+N\big)=\text{eventual valuation of $(a_\rho)$.}$$
It is easy to check that $v^{\operatorname{c}}$ is a valuation of the abelian group $G^{\operatorname{c}}$.
The map which associates to each element $a\in G$ the coset of $N$ containing the constant sequence~$(a)\in G^{\operatorname{cs}}$ is an embedding $G\to G^{\operatorname{c}}$ of groups, and identifying $G$ with its image under this embedding, the valued group $(G^{\operatorname{c}}, S, v^{\operatorname{c}})$ is an extension of~$(G,S, v)$.

Let $(a_\rho)\in G^{\operatorname{cs}}$ and $s\in S$. Since $(a_\rho)$ is a c-sequence, we can take $\rho_0$ such that  $v(a_\rho-a_{\rho_0})>s$ for all $\rho>\rho_0$; thus the eventual value of the c-sequence $(a_\rho-a_{\rho_0})$  is larger than $s$, hence $a_{\rho_0}\in G\cap B_{(a_\rho)+N}(s)$. Thus $G$ is dense in~$G^{\operatorname{c}}$. By a similar argument, $a_\rho\to (a_\rho)+N$ in $G^{\operatorname{c}}$. Therefore, by Lemma~\ref{lem:cofinality c-sequences}, \textit{every}\/ c-sequence in~$G$ converges in~$G^{\operatorname{c}}$.
It now follows that $G^{\operatorname{c}}$ is complete: Let $(b_\sigma)$ be a c-sequence in $G^{\operatorname{c}}$; we need to show that $(b_\sigma)$ has a limit in $G^{\operatorname{c}}$.
By what we have shown, applied to $G^{\operatorname{c}}$ in place of $G$, $(b_\sigma)$ has a limit  $b'$ in some dense valued abelian group extension~$G'$ of~$G^{\operatorname{c}}$. Since $G'\supseteq G$ is also dense,
by Lemma~\ref{lem:pc2 analogue} there is a c-sequence $(a_\rho)$ in~$G$ with $a_\rho\to b'$ in $G'$. This c-sequence converges in $G^{\operatorname{c}}$, thus $b'\in G^{\operatorname{c}}$.
\end{proof}

\begin{cor}
$G$ is complete iff $G$ has no proper dense extension.
\end{cor}

\noindent
So ``spherically complete $\Rightarrow$ maximal $\Rightarrow$ complete'' 
by Corollaries~\ref{fp1} and~\ref{impca}. Thus any Hahn product is complete, by Corollary~\ref{hahnmax}.

\medskip
\noindent
Corollary~\ref{cor:completion} below says that any two dense complete extensions of $G$ are isomorphic over $G$ by a unique isomorphism. We obtain this from a more general observation (Lemma~\ref{lem:extend to closure}) about the extension of continuous maps.

\medskip
\noindent
Let $(a_\rho)$, $(b_\sigma)$ be c-sequences in $G$ (with possibly different index sets). We call~$(a_\rho)$ and~$(b_\sigma)$   equivalent (in symbols: $(a_\rho)\sim (b_\sigma)$) if for each $s\in S$ there are $\rho_0$,~$\sigma_0$ such that $v(a_\rho-b_\sigma)>s$ for all $\rho>\rho_0$, $\sigma>\sigma_0$.
So if $(a_\rho)$ and $(b_\sigma)$ have the same limit in some valued abelian group extension of $G$, then
$(a_\rho)\sim (b_\sigma)$; conversely, if $(a_\rho)\sim (b_\sigma)$, then 
$(a_\rho)$ and $(b_\sigma)$ have the same limit in some dense valued abelian group extension of $G$ (by Theorem~\ref{thm:completion}).
The relation $\sim$ is an equivalence relation on the class of c-sequences in $G$. If $(b_\sigma)$ is a  cofinal subsequence of $(a_\rho)$, then $(a_\rho)\sim (b_\sigma)$. Moreover, given $a\in G$ we have
$a_\rho\to a$ iff $(a_\rho) \sim (a)$, where $(a)$ denotes an arbitrary constant c-sequence all of whose terms equal $a$.

\medskip
\noindent
Let $X\subseteq G$; by a c-sequence in $X$ we mean a c-sequence $(a_\rho)$ in $G$ such that $a_\rho\in X$ for each $\rho$. In the next two lemmas we consider a map $f\colon X\to G_1$ where $(G_1,S_1,v_1)$ is a valued abelian group. 
We say that $f$ is {\bf uniformly continuous} if for each $s_1\in S_1$ there is
an $s\in S$ such that for all $a,b\in X$, if 
$v(a-b)>s$, then~$v\big(f(a)-f(b)\big)>s_1$.
We say that $f$ is {\bf cauchy-continuous} (for short: {\bf c-continuous}) if for every c-sequence~$(a_\rho)$ in $X$, the image sequence $\big(f(a_\rho)\big)$ is a c-sequence in~$G_1$.
Clearly we have the implications 
$$\text{uniformly continuous $\ \Rightarrow\ $ c-continuous 
$\ \Rightarrow\ $ continuous.}$$ Conversely, if $G$ is complete, $X$ is closed,
and $f$ is continuous, then $f$ is c-continuous. 

\index{map!uniformly continuous}
\index{uniformly continuous}
\index{map!cauchy-continuous}
\index{cauchy-continuous}
\index{map!c-continuous}
\index{c-continuous}

\begin{lemma}
Suppose $f$ is c-continuous, $(a_\rho)$ and $(b_\sigma)$ are c-sequences in $X$,
and $(a_\rho)\sim (b_\sigma)$. Then $\big(f(a_\rho)\big) \sim \big(f(b_\sigma)\big)$.
\end{lemma}
\begin{proof}
By passing to cofinal subsequences of $(a_\rho)$, $(b_\sigma)$, we arrange that
the two sequences have the same ordered index set, by Lemma~\ref{lem:cofinality c-sequences}. Order the set of pairs~$(\rho,i)$, where $i=0,1$, by $(\rho,i)<(\rho',i')$ iff $\rho<\rho'$ or $\rho=\rho'$ and $i=0$, $i'=1$. Then the sequence $(c_{(\rho,i)})$ with $c_{(\rho,0)}=a_\rho$ and $c_{(\rho,1)}=b_\rho$ is a c-sequence in $X$ equivalent to both $(a_\rho)$ and $(b_\rho)$. Since $f$ is c-continuous, $\big(f(c_{(\rho,i)})\big)$ is a c-sequence in $G_1$, and  $\big(f(c_{(\rho,i)})\big)$ is equivalent to both $\big(f(a_\rho)\big)$ and $\big(f(b_\sigma)\big)$.
\end{proof}

\begin{lemma}\label{lem:extend to closure}
Suppose that for each c-sequence $(a_\rho)$ in $X$, the sequence~$\big(f(a_\rho)\big)$ converges in $G_1$. Then $f$ has a unique extension to a continuous map $X'\to G_1$,  where~$X'$ is the closure of $X$ in $G$ in the $v$-topology. This extension is uniformly continuous if $f$ is uniformly continuous.
$$\xymatrix@R-0.5em@C+1em{
X' \ar@{-->}[r]      				& G_1 \\
X \ar[u]^{\subseteq}\ar[ur]_{f}	& 
}$$
\end{lemma}
\begin{proof}
We only show existence of such an extension, since uniqueness then follows 
easily.
Let $a'\in X'$ be given; choose a c-sequence $(a_\rho)$ in $X$  such that $a_\rho\to a'$, and  let $f'(a')$ be the limit of $\big(f(a_\rho)\big)$ in $G_1$.
Note that $f'(a')$ does not depend on our choice of the sequence $(a_\rho)$: if $(b_\sigma)$ is another c-sequence in $X$ such that $b_\sigma\to a'$, then 
$(a_\rho)\sim (b_\sigma)$ and hence
$\big(f(a_\rho)\big)\sim \big(f(b_\sigma)\big)$ by the previous lemma, so
$\big(f(a_\rho)\big)$ and $\big(f(b_\sigma)\big)$ have the same limit in $G_1$. The map $f'\colon X'\to G_1$ so defined clearly extends $f$ and is continuous. It is easy to verify that if $f$ is uniformly continuous, then so is~$f'$.
\end{proof}

%\marginpar{Lemma~\ref{lem:extend to closure} is to extend the derivation of an asymptotic field to its completion.}

\noindent
Let $(G',S',v')$ be a valued abelian group. A {\bf morphism\/} 
$(G,S,v)\to (G',S',v')$ of valued abelian groups is a
pair $(i,j)$, where $i$ is a group morphism $G\to G'$ and~$j$ is an 
increasing map $S_\infty\to S'_\infty$ such that 
$j\big(v(x)\big)=v'\big(i(x)\big)$ for all $x\in G$
(so $j(\infty)=\infty$).  \index{valued abelian group!morphism} \index{morphism!valued abelian groups}
Such a morphism $(i,j)$ is
an {\bf embedding\/} if $j$ is injective (in which case~$i$ is also injective). \index{valued abelian group!embedding} \index{embedding!valued abelian groups}
%We call $i$ the embedding of abelian groups corresponding 
%to the embedding $(i,j)$ and $j$ the embedding of ordered sets 
%corresponding to $(i,j)$. 
If $(G,S,v)\subseteq (G',S',v')$, then we have an embedding
$$(G,S,v)\to (G',S',v')$$ of valued abelian groups
whose components are the natural inclusions 
$G\to G'$ and $S_\infty\to S'_\infty$. A morphism 
$(i,j)\colon (G,S,v)\to (G',S',v')$ of valued abelian groups
is said to be an {\bf isomorphism\/}
if $i\colon G\to G'$ and $j\colon S_\infty\to S'_\infty$ are bijections; note that then
$(i^{-1},j^{-1})\colon (G',S',v')\to (G,S,v)$ is also 
an isomorphism (of valued abelian groups). An embedding
$(i,j)\colon (G,S,v)\to (G',S',v')$ determines an isomorphism
$(i,j)\colon (G,S,v)\to (iG, jS, v'|_{iG})$ onto a valued subgroup of 
$(G',S',v')$. \index{valued abelian group!isomorphism} \index{isomorphism!valued abelian groups}

\begin{cor}\label{cor:extend to closure}
Let $(G,S,v)\xrightarrow{(i,j)} (G',S',v')$ be a morphism of valued abelian groups where $i$ is continuous. 
Let $(G_1,S_1,v_1)$ be a dense valued abelian group extension  of $(G,S,v)$ and 
$(G_1',S_1',v_1')$ a valued abelian group extension of $(G',S',v')$ such that 
$S'$ is cofinal in $S_1'$ and every c-sequence in $G'$ converges in $G_1'$.
$$\xymatrix@R-0.5em@C+1em{
(G_1,S_1,v_1) \ar@{-->}[r]			& (G_1',S_1',v_1') \\
(G,S,v) \ar[r]^{(i,j)}\ar[u] 	& (G',S',v') \ar[u]
}
$$
Then there is a unique extension of $(i,j)$ to a morphism $(G_1,S_1,v_1)\to (G_1',S_1',v_1')$ whose valued group component is continuous.
\end{cor}
\begin{proof} Since $i$ is actually uniformly continuous and 
$S'$ is cofinal in $S_1'$,
the group morphism $G\xrightarrow{\,i\,} G'\xrightarrow{\,\subseteq\,} G_1'$ is 
uniformly continuous. Hence it extends uniquely to a
continuous map $i_1\colon G_1\to G_1'$. Then $i_1$ is a group morphism 
since it restricts to a group morphism $G\to G'$ on the dense subgroup 
$G$ of $G_1$.
Moreover, for every c-sequence~$(a_\rho)$ in $G$ with limit $a\in G_1$ we have $i(a_\rho)\to i_1(a)$ in $G_1'$, hence 
\begin{align*}
j\big(v_1(a)\big)\ &=\ j\big(\text{eventual valuation of $(a_\rho)$}\big)\\
&=\ \text{eventual valuation of $i(a_\rho)$}\ =\ v_1'\big(i_1(a)\big).
\end{align*}
Thus $(i_1,j)$ is a morphism of valued abelian groups as desired.
\end{proof}

\noindent
This gives the uniqueness of dense complete extensions: 

\begin{cor}\label{cor:completion}
Let $G'$ be a dense complete extension of $G$ and $G_1$ a 
dense extension of $G$. Then there is a unique embedding $G_1\to G'$ over $G$, 
and if $G_1$ is in addition complete, this embedding is an isomorphism.
\end{cor}

\noindent
Thanks to this corollary we may refer to a dense complete extension of $G$
as the {\bf completion} of $G$, to be denoted by $G^{\operatorname{c}}$. 
By Corollary~\ref{cor:extend to closure}, taking completions is
functorial in the following way: \index{completion!valued abelian group} \index{valued abelian group!completion} \nomenclature[G]{$G^{\operatorname{c}}$}{completion of $G$}

\begin{cor}\label{cor:completion functor}
If $(G',S',v')$ is a valued abelian group extension of $(G,S,v)$ and
$S$ is cofinal in $S'$, then the inclusion $G\to G'$ extends uniquely to a valued abelian group embedding $G^{\operatorname{c}}\to G'^{\operatorname{c}}$.
\end{cor}

%\noindent
%The identity map $G\to G$ is c-continuous if $G$ is equipped with 
%the valuation $v$ on the left-hand side and  with any nontrivial coarsening of %$v$ on the right-hand side, and vice versa. Also, $G$ is complete with respect %to $v$ iff $G$ is complete with respect to any of its nontrivial coarsenings.
%More precisely, by Corollary~\ref{cor:extend to closure}, the completion 
%of $G$ %with respect to $v$ is essentially the same as the completion 
%of $G$ with respect to any of its nontrivial coarsenings:

%\begin{cor}\label{cor:completion coarsening}
%Let $v_1\colon G\to (S_1)_\infty$ be a nontrivial valuation on $G$ 
%and let $j\colon S_\infty\to (S_1)_\infty$ be an increasing surjection 
%such that $v_1=j\circ v$. Let 
%$(G_1^{\operatorname{c}},v_1^{\operatorname{c}},S_1)$ be 
%the completion of $(G,v_1,S_1)$.
%Then there is a unique group isomorphism 
%$i\colon G^{\operatorname{c}}\to G_1^{\operatorname{c}}$
%which is the identity on $G$, such that
%$(i,j)$ is a morphism $(G^{\operatorname{c}},v^{\operatorname{c}},S)\to (G_1^{\%operatorname{c}},v_1^{\operatorname{c}},S_1)$.
%\end{cor}

\noindent
Next we clarify the relationship between pc-sequences and c-sequences: %\marginpar{This is used, e.g., in the material on step-completeness.}

\begin{lemma}\label{lem:c-sequence width}
Each c-sequence in $G$ that is not eventually constant has a cofinal subsequence which is a pc-sequence of width~$\{\infty\}$. \textup{(}Hence $G$ is complete iff every pc-sequence of width $\{\infty\}$ has a pseudolimit in $G$.\textup{)}
\end{lemma}
\begin{proof}
Let $(a_\rho)$ be a c-sequence in $G$ that is not eventually constant.
Thus $\operatorname{cf}(S)$ is infinite. Take $a\in G^{\operatorname{c}}$ 
such that $a_{\rho} \to a$.
By passing to a suitable cofinal subsequence we arrange that
$a_\rho\ne a$ for all $\rho$. By Lemma~\ref{lem:cofinality c-sequences} we 
also arrange that the ordered index set of $(a_\rho)$ is the set of all ordinals 
$\rho< \operatorname{cf}(S)$. Then the subset
$$\big\{ v(a-a_\rho): \rho<\operatorname{cf}(S) \big\}$$ of $S$ is
cofinal in $S$.
An obvious transfinite recursion
yields an ordinal~$\lambda$ and a strictly increasing sequence 
$\big(\rho(\gamma)\big)_{\gamma < \lambda}$ in the set of indices 
$\rho< \operatorname{cf}(S)$, 
such that $v\big(a-a_{\rho(\gamma)}\big)$ as a function of 
$\gamma< \lambda$ is strictly increasing and cofinal in $S$.
Then $\operatorname{cf}(\lambda)=\operatorname{cf}(S)$ by Corollary~\ref{cor:cofinality}.
By passing to a cofinal subsequence and 
reindexing we replace~$\lambda$ by~$\operatorname{cf}(\lambda)$, and then 
$\lambda=\operatorname{cf}(S)$.
It remains to note that then $\big\{\rho(\gamma):\ \gamma < \lambda\big\}$ is 
cofinal in $\big\{\rho:\ \rho< \operatorname{cf}(S)\big\}$ 
by Corollary~\ref{cor:cofinality, 2}.
\end{proof}

%\noindent
%We say that a nest $\mathcal B$ of $\cup$-closed balls in $G$ 
%is a %{\bf cauchy nest} if for each $s_0\in S$ there are $a\in G$, 
%$s\in S^{> s_0}$ with $\overline{B}_a(s)\supseteq B$ for some 
%$B\in\mathcal B$. A cauchy nest of $\cup$-closed balls in $G$ 
%has at most one point in its intersection. 
%\index{cauchy nest} \index{nest!cauchy}

%\begin{lemma}\label{cor:cauchy nests} The following conditions 
%on $(G, S, v)$ are equivalent: \begin{enumerate}
%\item every nonempty cauchy nest of closed balls in $G$ 
%has nonempty intersection;
%\item $G$ is complete;
%\item every nonempty cauchy nest of $\cup$-closed balls in $G$ has %nonempty intersection.
%\end{enumerate}
%\end{lemma}
%\begin{proof}
%Let $(a_\rho)$ be a pc-sequence of width $\{\infty\}$ in $G$, 
%and let $\rho_0$ be as in the definition of pc-sequence.
%Then $\mathcal B=\big\{\overline{B}_{a_\rho}(s_\rho):\ 
%\rho>\rho_0\big\}$ is a nonempty cauchy nest of closed balls in $G$, 
%with $\bigcap\mathcal B=\{a\in G:a_\rho\leadsto a\}$.
%Thus $(1) \Rightarrow (2)$ by Lemma~\ref{lem:c-sequence width}.
%It remains to show $(2) \Rightarrow (3)$, so  assume (2),
%and let $\mathcal B$ be a nonempty cauchy nest of $\cup$-closed balls %in $G$. To get a point in its intersection we can arrange by the
%proof of Corollary~\ref{corgenclosedballs} that $\mathcal{B}$ %consists of closed balls.
%Then we construct a point in its intersection just as in the 
%proof of Lemma~\ref{fp1}, noting that the pc-sequence 
%constructed in that proof will in our case have width $\{\infty\}$. 
%\end{proof}

\noindent
For use in Section~\ref{sec:pc in valued fields} we also note:

\begin{lemma}\label{lem:complete approx} Suppose $(G,S,v)$ is complete and
$(G',S',v')$ is a valued abelian group extension of $(G,S,v)$.  
Let $a'\in G'\setminus G$ be such that $v'(a'-G)$ $($a subset of $S')$ has no 
largest element. Then $v'(a'-G)\subseteq S$, and $v'(a'-G)$ 
is not cofinal in $S$.
\end{lemma}
\begin{proof} First, given $g\in G$, take $a\in G$ such that 
$v'(a'-g)< v'(a'-a)$. Then $v'(a'-g)=v(a-g)$. Thus $v'(a'-G)\subseteq S$. 
Towards a contradiction, assume that $v'(a'-G)$ 
is cofinal in $S$.  Take a well-indexed sequence $(a_\rho)$ in $G$ such that 
$\big(v'(a'-a_\rho)\big)$ is strictly increasing and cofinal in $v'(a'-G)$, 
and hence in $S$. Then $(a_\rho)$ is a pc-sequence of width $\{\infty\}$ in $G$,
and thus $a_\rho\leadsto a\in G$, so $v'(a'-a)> v'(a'-G)$, a contradiction.
\end{proof}

\subsection*{Direct product of valued abelian groups}
Let $(G_i,S_i,v_i)$ be a valued abelian group for $i=1,\dots,n$, where $n\ge 1$.
Suppose $S$ is an ordered set containing each~$S_i$ as an ordered subset, with $S=S_1\cup\cdots\cup S_n$.
Set $G:=G_1\times\cdots\times G_n$ (direct product of abelian groups), and
for each $g=(g_1,\dots,g_n)\in G$, set
$$vg\  :=\ \min\big( v_1(g_1),\dots, v_n(g_n) \big).$$
Then $v\colon G\to S_\infty$ is a valuation on $G$. We call the valued abelian group $(G,S,v)$ the {\bf direct product} of the valued abelian groups $(G_i,S_i,v_i)$.
Given a closed ball $B=\overline{B}_{g}(s)$ in $(G,S,v)$, with $g=(g_1,\dots,g_n)\in G$, $s\in S$, we have
\begin{align*} B\ &=\ B_1\times\cdots\times B_n\quad\text{where $B_i=\big\{x\in G_i:v_i(x-g_i)\geq s\big\}$, so}\\
     B_i\ &=\ \bigcup\big\{\overline{B}_{g_i}(s_i):\  s_i\in S_i,\ s_i\ge s\big\} \quad
\text{if there is an $s_i\ge s$ in $S_i$, and}\\ 
    B_i\ &=\  \{g_i\} \quad \text{if there is no $s_i\ge s$ in $S_i$}.
\end{align*}

\index{valued abelian group!direct product}

\begin{lemma}\label{lem:direct product of vag}
If each $(G_i,S_i,v_i)$ is spherically complete, then so is $(G,S,v)$.
\end{lemma} 
\begin{proof}
Let $\mathcal B\neq\emptyset$ be a nest of closed balls in $(G,S,v)$. For each 
$B\in\mathcal B$, let $B_1,\dots, B_n$ be as above. Let $i\in \{1,\dots,n\}$
be given. If $B_i$ is a singleton for some $B\in \mathcal B$, then
$\bigcap\{ B_i:B\in\mathcal B\}$ is a singleton. 
Otherwise, by the remarks preceding the lemma, 
$\{B_i:B\in\mathcal B\}$ is a
nonempty nest of union-closed balls in $(G_i,S_i,v_i)$. 
 If each $(G_i,S_i,v_i)$ is spherically complete, then
 Corollary~\ref{corgenclosedballs} yields for each $i$ a point 
 $g_i\in \bigcap\{ B_i:B\in\mathcal B\}$, and thus a point
 $g=(g_1,\dots, g_n)\in \bigcap \mathcal B$. 
\end{proof}

\begin{lemma}\label{productcomplete} Suppose each $S_i$ is cofinal in $S$. Then:
$$ \text{$(G, S, v)$  is complete }\ \Longleftrightarrow\  \text{ each
 $(G_i,S_i,v_i)$ is complete.} $$
 \end{lemma} 
\begin{proof} Let $(g_{\rho})$ be a well-indexed sequence in $G$, with $g_\rho=(g_{\rho 1},\dots,g_{\rho n})$ and $g_{\rho i}\in G_i$ for each $i$ and $\rho$. By the cofinality assumption, $(g_{\rho})$ is a c-sequence in $G$ iff~$(g_{\rho i})$ is for every $i$ a c-sequence in $G_i$. This
 yields the desired equivalence.
\end{proof}

\noindent
We now consider the case where all $(G_i,S_i,v_i)$ are equal, say $(G_i,S_i,v_i)=(H,S,v)$ for $i=1,\dots,n$.
Then the valuation topology on $G$ with respect to the valuation $v\colon G\to S_\infty$ on  $G$ defined above coincides with the product topology on $G=H^n$ with respect to  the valuation topology of $H$.  
By the previous lemma, the dense extension~$(H^{\operatorname{c}})^n$ of $G=H^n$ is complete, so we may take $(H^{\operatorname{c}})^n$ as the completion~$G^{\operatorname{c}}$ of~$G$. Note also that by
Lemma~\ref{lem:extend to closure}:

\begin{cor}\label{lem:continuous extension to completion}
Every c-continuous map $X\to H$ \textup{(}$X\subseteq H^n$\textup{)} has a 
unique continuous extension $X'\to H^{\operatorname{c}}$ to the closure $X'$ of $X$ in $(H^{\operatorname{c}})^n$.
\end{cor}

\subsection*{Notes and comments}
The material in this section is classical valuation theory, some of it generalized from valued
fields to valued abelian groups. An early explicit mention of valued abelian groups is in \cite[Chapter~IV, \S{}4]{Fuchs}.
Lemma~\ref{krullgravett} is due to Krull~\cite{Krull}; the short proof given here was found by Gravett~\cite{Gravett-1}. The concept of pseudocauchy sequence and its basic properties, for rank~$1$ valued fields, are due to Ostrowski~\cite{Ostrowski}. (See also Section~\ref{sec:pc in valued fields} below.)
The Fixpoint Theorem~\ref{thm:fixpoint} is due to Prie\ss{}-Crampe~\cite{PC2} (in a more general setting of ultrametric spaces). 

Equivalence of pc-sequences will show up
in several later chapters. This notion of equivalence seems to have been used first in \cite{Scanlon}, and the flexibility it provides
was exploited heavily in \cite{BMS} in connection with valued difference fields.

The completion of a valued abelian group $G$
can also be introduced as the completion of the hausdorff topological group
$G$ (with respect to the valuation topology) equipped with the natural extension of the valuation of $G$ to this completion by continuity (see, e.g., \cite[\S{}II.4]{PC}); we chose the approach via c-sequences in order to stress the 
analogy with
pc-sequences and spherical completeness. 
Lemma~\ref{lem:c-sequence width} appears in~\cite[II, \S{}4, Lemma~11]{PC}, with a corrected proof (different from the one here) given in~\cite[Theorem~5.7]{SS}.

%% file: mt-2-3.tex
\section{Valued Vector Spaces}\label{sec:valued vector spaces} 

\noindent
In this section we define valued vector spaces, characterize spherically 
complete valued vector spaces, and then prove a version of the Hahn 
Embedding Theorem for valued vector spaces. 
We pay special attention to Hahn spaces, which are
particularly well-behaved valued vector spaces.

\medskip\noindent
Let $C$ be a field. A {\bf valued vector space over $C$}
is a vector space $G$ over $C$ equipped with a surjective valuation 
$v\colon G \to S_{\infty}$ on its underlying additive group such that
$v(\lambda a)=v(a)$ for all $\lambda\in C^\times$ and $a\in G$. 
Note that then for each $s\in S$ the subgroups $B(s)=B_0(s)$ and $\bar{B}(s)=\bar{B}_0(s)$ of $G$ are subspaces of $G$, and hence $G(s)=\bar{B}(s)/B(s)$ is in a natural way a vector space over $C$. 

\index{valued vector space}
\index{vector space!valued}

\begin{example}[Hahn products]
Let $S$ be an ordered set and $(G_s)_{s\in S}$ a family of vector spaces
$G_s\ne \{0\}$ over $C$. Then the Hahn product $H=H\big[(G_s)\big]$ of $(G_s)$ is a subspace of the vector space $\prod_s G_s$ over $C$, and $H$ with its 
Hahn valuation is a valued vector space over $C$. For each $s\in S$ the natural group isomorphism $H(s)\to G_s$ is $C$-linear.
\end{example}

\noindent
{\em For the remainder of this section $(G,S,v)$ and $(G',S',v')$ are valued 
vector spaces over $C$}. We call $(G,S,v)$ a {\bf valued subspace} 
of~$(G',S',v')$, or~$(G',S',v')$ an {\bf extension} of~$(G,S,v)$, if, as valued abelian groups,
$(G,S,v)$ is a valued subgroup of $(G',S',v')$, and 
$G$ is a subspace (not merely a subgroup) of $G'$. 

\index{valued vector space!subspace}

\medskip\noindent
An {\bf embedding} \index{embedding!valued vector spaces} \index{valued vector space!embedding} of valued vector spaces over $C$
is an embedding 
$$(i,j)\ :\ (G,S,v)\to (G',S',v')$$ of
valued abelian groups 
as defined in Section~\ref{sec:valued abelian gps} such that 
$i\colon G\to G'$ is $C$-linear. Such an embedding 
$(i,j)\colon (G,S,v)\to (G',S',v')$ is said to be an {\bf isomorphism} \index{isomorphism!valued vector spaces}  \index{valued vector space!isomorphism} (of valued vector spaces over $C$) if $i\colon G\to G'$ and 
$j\colon S_\infty\to S'_\infty$ are bijections; note that then
$(i^{-1},j^{-1})\colon (G',S',v')\to (G,S,v)$ is also 
an isomorphism. Such an embedding
$(i,j)\colon (G,S,v)\to (G',S',v')$ determines an isomorphism
$$(i,j)\ :\  (G,S,v)\to (iG, jS, v'|{iG})$$ onto a valued subspace of 
$(G',S',v')$.  

\medskip\noindent
We consider $C$ as the $1$-dimensional valued vector space over $C$
whose valuation is trivial. Note that any $1$-dimensional valued vector 
space over $C$ is isomorphic to the valued vector space $C$. 

\medskip\noindent
Below $(G,S,v)$ is usually referred to as \textit{the valued vector space $G$}\/, 
in particular when
$S$ and $v$ are clear from the context. Also, when $S=S'$, as ordered sets,
then by an embedding $i\colon G \to G'$ of valued vector spaces over $C$ we mean
an embedding $(i,\id_{S_{\infty}})\colon (G,S,v) \to (G',S,v')$ of valued vector spaces over 
$C$.

\subsection*{Maximal valued vector spaces} 
Note that if $(a_{\rho})$ is a pc-sequence in $G$ and $a_{\rho} \leadsto a\in G$,
then for $g\in G$ and $\lambda\in C^\times$ we have a pc-sequence
$(g+\lambda a_{\rho})$ in $G$ with $g+\lambda a_{\rho}\leadsto g+\lambda a$.

\begin{lemma}\label{divimm} Let $(a_{\rho})$ be a divergent pc-sequence in 
$G$. Then $G$ has an immediate extension $G\oplus Ca$  with 
$a_{\rho} \leadsto a$, such that any extension $G\oplus Cb$ of 
$G$ with $a_{\rho} \leadsto b$ gives an isomorphism 
$G\oplus Ca \to G\oplus Cb$ of valued vector spaces over $C$ 
that is the identity
on $G$ and sends $a$ to $b$.
\end{lemma}
\begin{proof} Take an extension 
$G\oplus Ca$ of $G$ such that $a_{\rho} \leadsto a$. (This is possible
by Lemma~\ref{pc3} and the remark following its proof.) 
We claim that this extension is immediate. Indeed, we have $s\in S$ such that
$$v(a- a_{\rho})\ >\ va\ =\ va_{\rho}\ =\ s,\ \text{ eventually,}$$
in particular, $va\in S$, and, eventually, 
$a+B(s)=a_{\rho} + B(s)$ in $G\oplus Ca$. Instead of~$a$, this 
also works for each of its affine transforms 
$g+ \lambda a$ with $g\in G$ and $\lambda\in C^\times$, since 
$(g+\lambda a_{\rho})$ is a divergent pc-sequence in $G$ with
$g+\lambda a_{\rho}\leadsto g+\lambda a$, for such~$g$,~$\lambda$.
This proves the immediacy claim. If $G\oplus Cb$ is also an extension
with $a_{\rho} \leadsto b$, then by the above, for any $g\in G$ and
$\lambda\in C^\times$,
$$g+\lambda a\ \asymp\ g+\lambda a_{\rho}\ \asymp\ g+\lambda b,\ 
\text{ eventually,}$$
and so we have an isomorphism $G\oplus Ca \to G\oplus Cb$ of valued vector spaces over $C$ that is the identity
on $G$ and sends $a$ to $b$. 
\end{proof}

\noindent
Call $G$ {\bf maximal\/} if it has no proper immediate (valued vector space) 
extension. From Corollary~\ref{impca} and Lemmas~\ref{fp1},~\ref{divimm}
we obtain: \index{valued vector space!maximal} \index{maximal!valued vector space}

\begin{cor}\label{eqmaxsph} The following are equivalent: 
\begin{enumerate}
\item[\textup{(i)}] $G$ is maximal as a valued abelian group;
\item[\textup{(ii)}] $G$ is maximal;
\item[\textup{(iii)}] each pc-sequence in $G$ has a pseudolimit in $G$;
\item[\textup{(iv)}] $G$ is spherically complete.
\end{enumerate}
\end{cor}

\noindent
Using Zorn, we obtain from Lemmas~\ref{krullgravett},~\ref{pc2}, 
and~\ref{divimm}:

\begin{cor}\label{isomaximm} The valued vector space $G$ has a maximal 
immediate extension, 
and all maximal immediate extensions of $G$ are isomorphic over $G$.
\end{cor}

%Let $(i,j)$ be such an embedding of valued vector spaces over $C$. 
%Then each $s\in S$
%gives an embedding $G(s)\to G'(j(s))$ of vector spaces 
%over $C$ sending, for $a\in\bar{B}(s)$, the coset $a+B(s)$ of 
%$B(s)$ to the coset $i(a)+B'(j(s))$ of $B'(j(s))$. We say 
%that $(i,j)$ is {\bf immediate} if $j: S \to S'$ and all 
%induced group embeddings $G(s)\to G'(j(s))$ are bijective. 
%\marginpar{commented out material checked but not needed?}

\noindent
The next result is a form of the Hahn Embedding Theorem and  
identifies the maximal immediate 
extension of $G$ from Corollary~\ref{isomaximm}.

\begin{prop}\label{prop:Hahn for valued vector spaces} 
There exists an embedding $i\colon G\to H\big[(G(s))\big]$ of valued vector spaces 
over $C$ such that for all $g\in G^{\ne}$,
$$ i(g)_s\ =\ g+B(s)\in G(s)\quad  \text{ where $s=v(g)$.}$$
Note that then $H\big[(G(s))\big]$ is an immediate extension of its valued 
subspace $iG$.
\end{prop}
\begin{proof} Set $H:=H\big[(G(s))\big]$, and let $H'$ be the subspace of $H$ 
consisting of all $g\in\prod_s G(s)$ with \textit{finite}\/ support. 
The valuation of $H$ restricts to a surjective valuation $H'\to S_\infty$, 
making $H$ an immediate extension of its valued subspace $H'$. 

Choose for each $s\in S$ a $C$-linear right-inverse 
$f_s\colon G(s)\to \bar{B}(s)$ of the natural projection 
$\bar{B}(s) \to \bar{B}(s)/B(s)=G(s)$, and 
let $f$ be the $C$-linear map $H'\to G$ such that $f((g_s))=\sum_s f_s(g_s)$.
Then $G$ is an immediate extension of its valued subspace 
$G':= f(H')$, and we have an isomorphism $i'\colon G'\to H'$ of valued vector
spaces over~$C$ given by $i'(g)=h$ iff $f(h)=g$, for all $g\in G', h\in H'$.

Take a maximal immediate extension $\tilde{G}$ of $G$. This is also a maximal
immediate extension of $G'$. Since $H$ is a maximal immediate extension
of $H'$, it follows from Corollary~\ref{isomaximm} that $i'$ extends to an isomorphism
$\tilde{G}\to H$. Then the restriction of this isomorphism to
$G$ gives an embedding $i\colon G \to H$ as in the proposition.
\end{proof}

\noindent
In Section~\ref{sec:valued abelian gps} we defined the completion $G^{\operatorname{c}}$ of $G$ as a valued abelian group. We make~$G^{\operatorname{c}}$ into a vector
space over $C$ as follows.
Given $\lambda\in C$, 
the map $$x\mapsto \lambda x\ \colon\ G \to G$$ is uniformly continuous, so extends 
uniquely to a continuous map  $G^{\operatorname{c}} \to G^{\operatorname{c}}$
by Lemma~\ref{cor:extend to closure}. Denoting the image of
$x\in G^{\operatorname{c}}$
under this map also by $\lambda x$ and varying~$\lambda$ we have an
operation  
$$ (\lambda,x)\mapsto \lambda x\ \colon\  C\times G^{\operatorname{c}} \to  
G^{\operatorname{c}}$$ that makes $G^{\operatorname{c}}$ a 
valued vector space over $C$ extending $G$.

%Inspection of the proof of Theorem~\ref{thm:completion}) shows that 
%$G^{\operatorname{c}}$ is in a natural way a vector space over $C$ 
%making $G^{\operatorname{c}}$ a valued vector space extension on $G$.
%\marginpar{above commented out sentence has been checked
%and gives an alternative way of making $G^{\operatorname{c}}$ a vector space}

\subsection*{Valuation bases}
Let $G_0$ be a subspace of $G$, and let $\mathcal B\subseteq G$.
We say that $\mathcal B$ is {\bf valuation-independent over~$G_0$} if $0\notin\mathcal B$, and
for every family $(\lambda_b)_{b\in \mathcal B}$ of elements of $C$, with 
$\lambda_b\ne 0$ for only finitely many $b\in \mathcal B$, and every $g_0\in G_0$, 
$$v\left(\sum_{b\in \mathcal B} \lambda_b b + g_0\right)\  =\  
\min \big( \{vb:\ b\in\mathcal B,\ \lambda_b\ne 0\}\cup \big\{v(g_0)\big\} \big).$$
This gives a linearly independent family $(b+G_0)_{b\in \mathcal B}$ in the $C$-vector space~$G/G_0$. For a set $X\subseteq G$ we let~$\langle X\rangle$ be the $C$-linear subspace of $G$ generated by $X$.
If $\mathcal B$ is valuation-independent over $G_0$ and $\langle\mathcal B \cup G_0\rangle = G$,
then we call~$\mathcal B$ a {\bf valuation basis} of~$G$ 
over $G_0$.
If $\mathcal B$ is valuation-independent over $\{0\}$, then we just say that~$\mathcal B$ is valuation-independent, and if $\mathcal B$ is a valuation basis of $G$ 
over~$\{0\}$, then we just say that $\mathcal B$ is a valuation basis of $G$.
Thus:  \index{valuation!independence} \index{valuation!basis} \index{independent!valuation}

\begin{lemma}\label{lem:val basis}
Let $\mathcal B, \mathcal B'\subseteq G$ be disjoint. Then $\mathcal B\cup\mathcal B'$ is valuation-independent over~$G_0$ if and only if $\mathcal B$ is valuation-independent over $G_0$ and $\mathcal B'$ is valuation-in\-de\-pen\-dent over $\langle\mathcal B\cup G_0\rangle$.
\end{lemma}

\noindent
Every valuation basis of $G$ over $G_0$ is clearly maximal among the 
subsets of $G$ that are valuation-independent over $G_0$.
It is also easy to see that an increasing union of subsets of $G$, each valuation-independent over $G_0$, is valuation-independent over~$G_0$; hence by Zorn, there exists a maximal (possibly empty) subset of $G$ which is valuation-independent over $G_0$.

\begin{lemma}\label{lem:val basis, 1} The extension
$G_0\subseteq G$ is immediate if and only if  $G$ has no nonempty subset that is valuation-independent over $G_0$.
\end{lemma}
\begin{proof}
Suppose $G_0\subseteq G$ is immediate, and let $g\in G^{\neq}$. 
Take $g_0\in G_0$ such that $g\sim g_0$. Then $v(g-g_0)\neq \min( vg,vg_0)$, so $\{g\}$ is not valuation-independent over $G_0$. On the other hand, suppose $G_0\subseteq G$ is not immediate.
Then we can take $g\in G^{\neq}$ such that for each $g_0\in G_0$ we have $v(g-g_0)\le vg$  and hence $v(g-g_0)=\min(vg,vg_0)$; then $\{g\}$ is valuation-independent over $G_0$.
\end{proof} 

\begin{cor}\label{cor:val basis}
There exists a subspace $G_1\supseteq G_0$ of $G$ such that $G_1$ admits a valuation basis over $G_0$, and $G_1\subseteq G$ is immediate. 
\end{cor}
\begin{proof}
Let $\mathcal B\subseteq G$ be maximal valuation-independent over $G_0$ and $G_1:={\langle \mathcal B\cup G_0\rangle}$. Then $G_1$ admits a valuation basis over $G_0$. Moreover, if $\emptyset\neq\mathcal B'\subseteq G$ is valuation-independent over $G_1$, then $\mathcal B\cup\mathcal B'$ is valuation-independent over $G_0$, contradicting maximality of $\mathcal B$. Hence by Lemma~\ref{lem:val basis, 1}, $G_1\subseteq G$ is immediate.
\end{proof}

\begin{cor}\label{cor139}
Suppose the set $v(G\setminus G_0)$ is reverse well-ordered. Then $G$ has a valuation basis over $G_0$.
\end{cor}
\begin{proof}
Take $G_1$ as in the previous lemma. For a contradiction assume $G\neq G_1$.
Then the nonempty subset $v(G\setminus G_1)$ of $v(G\setminus G_0)$ is also reverse well-ordered, so we can take $g\in G\setminus G_1$ such that $vg=\max v(G\setminus G_1)$. Since $G_1\subseteq G$ is immediate there is some $g_1\in G_1$ with $v(g-g_1)>vg$, and $g-g_1\in G\setminus G_1$, a contradiction.
\end{proof}

\noindent
We denote the dimension of a vector space~$V$ over $C$ by $\dim_C V$. \nomenclature[A]{$\dim_C V$}{dimension of the vector space $V$ over the field $C$}

\begin{lemma}\label{lem:val basis, 2}
$\dim_C G/G_0 \geq \abs{v(G)\setminus v(G_0)}$.
\end{lemma}
\begin{proof}
If $\mathcal B$ is a subset of $G$ such that $vb\notin v(G_0)$ for each $b\in\mathcal B$ and
$b\not\asymp b'$ for all $b\neq b'$ in $\mathcal B$, then $\mathcal B$ is $C$-linearly independent over $G_0$.
\end{proof}

\noindent
Note that by Corollary~\ref{cor139} and Lemma~\ref{lem:val basis, 2}, 
every finite-dimensional valued vector space over $C$ has a valuation basis. 
More generally:

\begin{cor}\label{cor:brown}
Every valued vector space over $C$ of countable dimension has a valuation basis.
\end{cor}
\begin{proof}
Suppose $G$ has dimension $\aleph_0$. Then $G=\bigcup_n G_n$ where 
$$G_0=\{0\}, \quad G_n\subseteq G_{n+1}, \quad \dim_C G_{n+1}/G_n=1  \qquad(\text{for every $n$}).$$ Then each $G_n$ is finite-dimensional, so $v(G_n)$
is finite by Lemma~\ref{lem:val basis, 2}. Thus $G_n$ as a valued subspace of~$G$ is spherically complete and hence maximal (Corollary~\ref{eqmaxsph}), so $G_n\subseteq G_{n+1}$ is not immediate. By Corollary~\ref{cor:val basis} and since
$\dim_C G_{n+1}/G_n=1$, we can take $b_{n}\in G_{n+1}$ such that $\{b_{n}\}$ is a valuation basis 
for $G_{n+1}$ over $G_n$. Then $\{b_n:n\geq 0\}$ is a valuation basis of $G$, by Lemma~\ref{lem:val basis}.
\end{proof}

\noindent
In Section~\ref{sec:max} we shall need ``good'' right-inverses to certain linear maps:

\begin{lemma}\label{valvecrightinv} Suppose $G$ is maximal
and $A\colon G \to G$ is a $C$-linear map such that $g\succeq A(g)$ for all $g\in G$, and for all $h\in G$ there is
$g\in G$ with $A(g)=h$ and $g\asymp h$. Then there exists a
$C$-linear map $B\colon G \to G$ such that $A\circ B=\operatorname{id}_G$ and $B(h)\asymp h$ for all $h\in G$.
\end{lemma}
\begin{proof} Assume $H_0\ne G$ is a $C$-linear subspace of $G$ and
$B_0\colon H_0\to G$ is a $C$-linear map such that $A\circ B_0$ is the
inclusion $H_0\to G$
%=\operatorname{id}_{H_0}$ 
and $B_0(h)\asymp h$ for all $h\in H_0$. 

\claim{There exists $h_1\in G\setminus H_0$ such that $B_0$ extends to a $C$-linear map 
$$B_1\ :\  H_1 \to G, \qquad H_1:= H_0+ C h_1,$$ for which 
$A\circ B_1$ is the inclusion $H_1\to G$
%=\operatorname{id}_{H_1}$ 
and $B_1(h)\asymp h$ for all $h\in H_1$.}

\noindent
It is clear that Lemma~\ref{valvecrightinv} follows from this claim. 
To prove the claim we first pick an element $b\in G\setminus H_0$ and distinguish two cases:

\case[1]{$H_0+Cb$ is not an immediate extension of $H_0$ \textup{(}where both are viewed as valued subspaces of $G$\textup{).}} Then we set 
$H_1:= H_0+ Cb$.
In view of Lemma~\ref{lem:val basis, 1} we can take $h_1\in H_1\setminus H_0$ 
such that $\{h_1\}$ is a valuation basis of $H_1$ over $H_0$. Take any $g\in G$ with $A(g)=h_1$ and
$g\asymp h_1$, and let $B_1\colon H_1\to G$ be the $C$-linear extension
of~$B_0$ with $B_1(h_1)=g$. Then $B_1$ has the claimed property.

\case[2]{$H_0+C b$ is an immediate extension of $H_0$.} Take a divergent pc-se\-quence~$(b_{\rho})$ in $H_0$ such that $b_{\rho} \leadsto b$. Then $(a_{\rho}):=\big(B_0(b_{\rho})\big)$ is a divergent pc-sequence in~$B_0(H_0)$, and has a pseudolimit $a$ in $G$. Then $A(a_{\rho})=b_{\rho} \leadsto h_1:= A(a)$ by
Lemma~\ref{eqpslim}, so $h_1\notin H_0$.   
Set $H_1:= H_0+ Ch_1$ and let $B_1\colon H_1\to G$ be the $C$-linear extension
of~$B_0$ with $B_1(h_1)=a$. Then $B_1$ has the claimed property.
To see this, consider any $h= h_0+\lambda h_1\in H_1$
with $h_0\in H_0$ and $\lambda\in C^\times$. Then $\big(B_0(h_0)+\lambda a_{\rho}\big)$ is a divergent pc-sequence in $B_0(H_0)$ and
$B_0(h_0) + \lambda a_{\rho}\leadsto B_0(h_0) + \lambda a$,
so $$B_1(h)\ =\ B_0(h_0)+ \lambda a\ \asymp\ B_0(h_0) + \lambda a_{\rho},\  \text{ eventually}$$
by Lemma~\ref{pc1}. Applying $A$ to $B_0(h_0)+\lambda a_{\rho}$ we
get $h_0 + \lambda b_{\rho}\leadsto h_0+\lambda h_1=h$, and
so $h\asymp h_0+\lambda b_{\rho}$, eventually. Since 
$B_0(h_0) + \lambda a_{\rho}\asymp h_0+\lambda b_{\rho}$
for all $\rho$, we get $B_1(h)\asymp h$,
as required. 
\end{proof}

\noindent
The rest of this section is mainly intended for use in a later volume.

\subsection*{Hahn spaces} 
If for all $a,b\in G^{\ne}$ with $a\asymp b$ there
exists $\lambda\in C$ such that $a\sim \lambda b$, then we call~$G$ a
{\bf Hahn space over~$C$.} 
Equivalently, $G$ is a Hahn space over $C$ iff all vector spaces~$G(s)$ have dimension~$1$. In particular, an immediate extension of a Hahn space is a Hahn space, and a valued subspace of a Hahn space is again a Hahn space. 
Given a family~$(G_s)$ of $1$-dimensional vector spaces over $C$, indexed by the elements $s$ of an ordered set $S$, the Hahn product $H\big[(G_s)\big]$ of $(G_s)$ is a Hahn space over $C$. 
By Proposition~\ref{prop:Hahn for valued vector spaces}, every Hahn space embeds into a Hahn product of this kind: \index{Hahn space}

\begin{cor}\label{Hahn space}
If $G$ is a Hahn space over $C$, then there is an embedding 
$i\colon G\to H[S,C]$ of valued vector spaces over $C$ such that $H[S,C]$
is an immediate extension of its valued subspace $iG$.
\end{cor}

\noindent
Suppose $G$ is a Hahn space over $C$ and $\mathcal B\subseteq G$ with $0\notin \mathcal B$.  Then $\mathcal B$ is valuation-independent iff $b\not\asymp b'$ for all $b\neq b'$ in $\mathcal B$, and $\mathcal B$ is maximal valuation-independent iff it is valuation-independent and $v(\mathcal B)=S$. Thus if $G$ 
admits a valuation basis, then $\dim_C G=\abs{S}$, and $G$ is
determined (up to isomorphism of valued vector spaces) by the ordered set $S$. 
By Corollary~\ref{cor:brown}:

\begin{lemma} \label{Hahn valuations, lemma}
Suppose $G$ is a Hahn space of countable dimension as a vector space over $C$.
Then $G$ has a basis $\mathcal B$ with $v(\mathcal B)=S$ and $b\not\asymp b'$ for all $b\neq b'$ in~$\mathcal B$.
\end{lemma}
%\begin{proof}
%This follows by a well-known argument.  Note first that if $y_1,\dots, y_m\in G^{\ne}$  and $vy_1 < \dots < vy_m$, then $y_1,\dots,y_m$ are linearly independent over $C$.  Hence  the size of $v(G^{\ne})$ is at most $d$. We claim that $v(G^{\ne})$ has size $d$. Assume $d >0$. Take $y\in G^{\ne}$ such that $vy=\min v(G^{\ne})$. Then there is for any $a\in G$ a unique $\lambda(a)\in C$ such that $vy < v\bigl(a-\lambda(a)y\bigr)$. The set $\bigl\{a-\lambda(a)y: a\in G\bigr\}$ is a $C$-linear subspace of $G$ of dimension $d -1$, so we can proceed inductively.
%\end{proof}

\subsection*{Scalar extension} 
\textit{In this subsection $G$ is a Hahn space over $C$, and $D$ is an extension field of $C$.}\/ We construe $G_{D}:=D\otimes_{C}G$ as a vector space over $D$ in the usual way. Identifying $g\in G$ with
$1\otimes g\in G_{D}$ makes
$G$ a $C$-linear subspace of $G_{D}$. 
%The following is very useful:

\begin{lemma}\label{LemmaStandardForm} Let $h\in G_{D}$. There are $d_i\in D^{\ne}$, $g_i\in G^{\neq}$, $i=1,\dots,m$, such that 
\begin{equation}\label{StandardForm}
h\ =\ \sum_{i=1}^m d_i g_i,  \qquad g_1\succ \cdots \succ g_m. 
\end{equation}
\end{lemma}
\begin{proof}
We have $h=\sum_{j=1}^n e_j h_j$ ($e_j\in D$, $h_j\in G$).
We use induction on~$n$ to get~$h$ in the form \eqref{StandardForm}.
If $n=1$ we can do this with $m=0$ or $m=1$, so suppose $n>1$.
We can assume $h_1\succeq \cdots \succeq h_n\ne 0$.
Take $c_j\in C$ with $h_j-c_jh_1\prec h_1$ for $j=2,\dots,n$.
Then
$$h\ =\ \left(e_1+\sum_{j=2}^n c_je_j\right) h_1+\sum_{j=2}^n e_j(h_j-c_jh_1).$$
It remains to apply the inductive hypothesis to $\sum_{j=2}^n e_j(h_j-c_jh_1)$. 
\end{proof}

\begin{cor}\label{ScalarExtension}
There is a unique valuation~$w$ on the abelian group $G_{D}$ that
extends the valuation $v$ of $G$ and makes $G_{D}$ a valued vector space over~$D$. It makes~$G_D$ a Hahn space over $D$ with $w(G_D)=v(G)$.
\end{cor}
\begin{proof} For such a valuation $w$ and $h\in G_D^{\ne}$
as in \eqref{StandardForm} we have $w(h)=v(g_1)$. As to
existence:
Corollary~\ref{Hahn space} gives an
embedding $G \to H[S,C]$ of valued vector spaces over $C$; it extends to 
a $D$-linear injective map $G_{D}\to H[S,D]$.  The Hahn valuation of~$H[S,D]$ yields in this way a valuation $w$ on $G_D$ as claimed. \end{proof}

\subsection*{Fredholm operators} In order to discuss some properties of linear operators on Hahn spaces, we begin with some
generalities about Fredholm operators.
Let $V$ and $V'$ be  vector spaces over $C$ and $A\colon V\to V'$
be a $C$-linear map, with kernel
$\ker A := \big\{v\in V:A(v)=0\big\}$ and cokernel $\coker A := V'/A(V)$.
One says that~$A$ is a {\bf Fredholm operator}
if both $C$-vector spaces $\ker A$ and
$\coker A$ are finite-dimensional, and in that case we define the 
{\bf index} of $A$ by
$$\ind A := \dim_C\ker A - \dim_C\coker A\qquad \text{(an integer).}$$
If both $V$ and $V'$ are finite-dimensional, then every $C$-linear map $V\to V'$
is a Fredholm operator of index $\dim_C V-\dim_C V'$. \index{Fredholm operator} \index{operator!Fredholm} \index{index!Fredholm operator}
\nomenclature[A]{$\ker A$}{kernel of a linear map $A$}
\nomenclature[A]{$\coker A$}{cokernel of a linear map $A$}
\nomenclature[A]{$\ind A$}{index of a Fredholm operator $A$}

\begin{prop}\label{Fredholm}
If $A\colon V\to V'$ and $B\colon V'\to V''$ are Fredholm operators, then so is
$B\circ A\colon V\to V''$, and
$\ind B\circ A = \ind B+\ind A$. 
\end{prop}

\noindent
Towards the proof of this proposition, we first note an easy lemma:

\begin{lemma}\label{lem:Fredholm direct sum}
Let $A\colon V\to V'$ and $B\colon W\to W'$ be $C$-linear.  
Then $A$ and $B$ are Fredholm operators iff 
$A\oplus B\colon V\oplus W\to V'\oplus W'$ is one; in that case 
$$\ind(A\oplus B)\ =\ \ind A + \ind B.$$
\end{lemma}

\noindent
We also need two small items from linear algebra; for the first one, see for example~\cite[Lem\-ma~III.9.1]{Lang}; the second one is a special case of Lemma~\ref{alternatingexact, length}.

\begin{lemma}[Snake Lemma] \label{lem:snake} \marginpar{for the moment assume
without proof}
Assume that the following diagram of vector spaces over $C$ and $C$-linear maps is commutative with exact rows: 
\[ \xymatrix{ 
			& V'\ar[d]^{F'} \ar[r]^{i}	& V \ar[d]^{F} \ar[r]^{p}	& V'' \ar[d]^{F''} \ar[r]	& 0 \\
0 \ar[r]		& W'\ar[r]^{j}			& W \ar[r]^{q}					& W''						& }
\]
Then we have an exact sequence of $C$-linear maps
\[ 
\ker F' \longrightarrow	 \ker F \longrightarrow	 \ker F'' \xrightarrow{\ \delta\ }	 \coker F' \longrightarrow	 \coker F \longrightarrow	 \coker F''
\]
where the maps besides the ``connecting'' map $\delta$ are the natural ones. 
\end{lemma}

\begin{lemma}\label{alternatingexact}
For each exact sequence
$$0\longrightarrow V_1 \longrightarrow V_2 \longrightarrow \cdots \longrightarrow V_n \longrightarrow 0$$
of finite-dimensional vector spaces over $C$ we have
$$\sum_{i=1}^n (-1)^i \dim_C V_i=0.$$
\end{lemma}

\noindent
Suppose in Lemma~\ref{lem:snake} the map $i$ is injective, and 
the map $q$ is surjective. Then by Lemmas~\ref{lem:snake} and~\ref{alternatingexact}, if
two of the three maps $F$, $F'$ and $F''$ are Fredholm operators, then so is the third, and $\ind F=\ind F'+\ind F''$.

\begin{proof}[Proof of Proposition~\ref{Fredholm}]
Let $A\colon V\to V'$ and $B\colon V'\to V''$ be Fredholm operators. Let
$$i\colon V \to V\oplus V',\quad p\colon V\oplus V'\to V',\quad j\colon V' \to V''\oplus V',\quad 
q\colon V''\oplus V' \to V'$$ be the
$C$-linear maps given by 
\begin{align*}
i(x)		&=\big(x,A(x)\big), 		& p(x,x')	&=A(x)-x', &&  \\
j(x')	&=\big(B(x'),x'\big),	& q(x'',x')	&=x''-B(x') & \quad &(x\in V, x'\in V', x''\in V''). 
\end{align*} 
Then the diagram
\[ \xymatrix{ 
0 \ar[r]		& V\ar[d]^{A} \ar[r]^i	& V\oplus V' \ar[d]^{(B\circ A)\oplus\id_{V'}} \ar[r]^p	& V' \ar[d]^{B} \ar[r]	& 0 \\
0 \ar[r]		& V'\ar[r]^{j}		& V''\oplus V' \ar[r]^q								& V'' \ar[r]				& 0}
\]
commutes and has exact rows.
%\begin{align*}
%i(x)		&=(x,A(x)), 		& p(x,x')	&=A(x)-x', &&  \\
%j(x')	&=(B(x'),x'),	& q(x'',x')	&=x''-B(x') & \quad &(x\in V, x'\in V', %x''\in V''). 
%\end{align*}
It now remains to use the remark preceding this proof and Lemma~\ref{lem:Fredholm direct sum}.
\end{proof}

\subsection*{Linear operators on Hahn spaces}
{\em In this subsection $G$ is a Hahn space over~$C$, and $A\colon G\to G$ is a $C$-linear operator. We set
$\ker^{\ne} A := (\ker A)^{\ne}$}.

\begin{lemma} \label{Hahn valuations, cor} 
Suppose that $\ker A$ is finite-dimensional. 
Then there exists a $C$-linear
subspace $M$ of $G$ such that
$$G=M\oplus \ker A, \qquad
v(M)\cap v(\ker^{\ne} A)=\emptyset.$$
\end{lemma}
\begin{proof}
By Corollary~\ref{Hahn space} 
we may assume that $G$ is a valued subspace of the Hahn product $H[S,C]$. 
Let $y_1,\dots,y_d\in\ker^{\ne} A$ be such that
$vy_1>\cdots>vy_d$ are the distinct elements of $v(\ker^{\ne} A)$. 
Then
$$M := \big\{ a\in G : \supp(a)\cap 
\{vy_1,\dots,vy_d\}=\emptyset\big\}$$
is a $C$-linear subspace of $G$ such that
that $v(M)\cap v(\ker^{\ne} A)=\emptyset$ (in particular
$M\cap\ker A=\{0\}$). Given $a\in G$ we let $i(a)$ denote the largest
$i\in\{1,\dots,d\}$  such that $vy_i\in\supp(a)$, if such an
$i$ exists, and $i(a):=0$ otherwise. We show that for every $a\in G$ we have
$a\in M+\ker A$ by induction on $i(a)$. If $i(a)=0$, then $a\in M$.
If $i(a)=i>0$, then there exists $\lambda\in C^\times$ such that 
$vy_i\notin\supp(a-\lambda y_i)$,
hence $i(a-\lambda y_i)<i$ and therefore $a=(a- \lambda y_i)+ \lambda y_i\in M+\ker A$ by
inductive hypothesis. Thus $M$ has the required properties.
\end{proof}

\begin{cor} \label{v(A(N))}
If $\ker A$ is finite-dimensional, then
$$v\big(A(G)\big) = \big\{ v\big(A(y)\big) : y\in G,\ vy\notin v(\ker^{\ne} A)\big\}.$$
\end{cor}

\noindent
The set $\cal F_G$ of Fredholm operators $G\to G$ is a monoid under 
composition (by Proposition~\ref{Fredholm}), with identity element $\id_G$. In the lemma below we let 
$\cal S_G$ be the submonoid of 
$\cal F_G$ consisting of all surjective Fredholm operators $G\to G$. \nomenclature[A]{$\cal F_G$}{monoid of Fredholm operators on~$G$}
\nomenclature[A]{$\cal S_G$}{monoid of surjective Fredholm operators on~$G$} \nomenclature[G]{$A^{-1}$}{distinguished right-inverse of $A$}

\begin{lemma}\label{inverse}
There exists a map $A\mapsto A^{-1}\colon \cal S_G\to\cal F_G$ such that 
for all $A\in\cal S_G$:
$$A\circ A^{-1}=\id_G, \quad v\big(A^{-1}(G)\big)\cap v(\ker^{\ne} A)=\emptyset,$$
and if $B\colon G\to G$ is a bijective $C$-linear operator, then 
$$(B\circ A)^{-1}=A^{-1}\circ B^{-1}.$$
\end{lemma}
\begin{proof}
As in the proof of Lemma~\ref{Hahn valuations, cor}
we may assume that $G$ is a valued subspace of $H[S,C]$. For each $A\in\cal F_G$ we define
the subspace $M_A=M$ as in that proof.
Let $A\in\cal S_G$. For $g\in G$ we let $A^{-1}(g)$ be the unique $f\in M_A$
with $A(f)=g$. This yields a $C$-linear operator $A^{-1}\colon G\to G$
with $A\big(A^{-1}(g)\big)=g$ and $v\big(A^{-1}(g)\big)\notin v(\ker^{\ne} A)$ for all $g\in G$. Note that $A^{-1}$ is injective and $A^{-1}(G)=M_A$; in 
particular, $\ker A^{-1}=\{0\}$ and $\coker A^{-1}=G/M_A\cong \ker A$ are both 
finite-dimensional, hence $A^{-1}\in\cal F_G$. Let in addition $B\colon G\to G$
be a bijective $C$-linear operator; then $\ker B\circ A =\ker A$, hence 
$M_{B\circ A}=M_A$.
So if $g\in G$, then $f:=A^{-1}\big(B^{-1}(g)\big)$ satisfies $f\in M_{B\circ A}$ and
$(B\circ A)(f)=g$, hence $f=(B\circ A)^{-1}(g)$.
Thus $(B\circ A)^{-1}=A^{-1}\circ B^{-1}$.
\end{proof}

\subsection*{Notes and comments}
Corollaries~\ref{eqmaxsph} and~\ref{isomaximm} are in Gra\-vett~\cite{Gravett-2}.
(For valued abelian \textit{groups,}\/ the equivalence of maximality and spherical completeness remains true, whereas the uniqueness in Corollary~\ref{isomaximm} does
not go through; see~\cite{Fleischer}.)
Proposition~\ref{prop:Hahn for valued vector spaces} is from Conrad~\cite{Conrad};
it generalizes the classical embedding theorem for ordered abelian groups (Corollary~\ref{cor:Hahn classical} below) of Hahn~\cite{Hahn}. The proof given here follows~\cite{Gravett-2}. (See \cite{Ehrlich-2} for a detailed discussion of Hahn's work and its various spinoffs.)
Corollary~\ref{cor:brown} is due to Brown~\cite{Brown}. Our presentation follows \cite{Kuhlmann96}.
%, where it is shown more generally that if $G_0$ a valued subspace of 
%$G$ such that $\dim_C (G/G_0)$ is countable and $G_0$ admits a 
%valuation basis, then $G$ admits a valuation basis.
Our notion of Hahn space generalizes that of \cite{AvdD} where Hahn spaces
are a special kind of ordered vector spaces; see Section~\ref{sec:oag} below. Lemma~\ref{LemmaStandardForm} is implicit in the proof of \cite[Theorem 3]{Rosenlicht2} and is a variant of \cite[Lemma~2.1]{AvdD}.
The concept of Fredholm operator and Proposition~\ref{Fredholm} are well-known in operator theory; see~\cite{Sarason}.

%% file: mt-2-4.tex
\section{Ordered Abelian Groups}\label{sec:oag}

\noindent
An {\bf ordered abelian group} is an abelian group $G$ (usually written additively) equipped with a (total) ordering $\leq$ on its underlying set such
that for all $x,y,z\in G$,
$$x\leq y \quad\Rightarrow\quad x+z \leq y+z.$$
Each ordered abelian group is a hausdorff topological group with respect to
its order topology.
\index{ordered abelian group}
\index{abelian group!ordered}
\index{group!ordered abelian}
{\em In the rest of this section $G$ is an ordered abelian group}, and we let~$a$,~$b$,~$c$ range over $G$. We put $\abs{a}:=\max(a,-a)$. We set $G^{<}:= G^{<0}$ and $G^{>}:= G^{>0}$, and similarly with $\leq$ and $\geq$ in place of $<$ and $>$, respectively. 
We extend the
addition on $G$ to $G_\infty$ by
$a+\infty=\infty+a=\infty+\infty=\infty$. 

\begin{example}[Hahn products]
Let $(G_s)_{s\in S}$ be a family of ordered abelian groups $G_s\ne \{0\}$ indexed by an ordered set $S$. Then 
there is a unique ordering making
the Hahn product $H\big[(G_s)\big]$ of $(G_s)$ an ordered abelian group such that
$$0<g \quad \Longleftrightarrow\quad 0<g_{vg},$$
for all nonzero $g=(g_s)\in H\big[(G_s)\big]$; call it the {\bf Hahn ordering\/} of
$H\big[(G_s)\big]$.
\end{example}

\index{Hahn product!ordering}
\index{ordering!Hahn}

\noindent
Every ordered abelian group is torsion-free. 
We shall consider each torsion-free abelian group $A$ as a subgroup of the
divisible abelian group $\Q A=\Q\otimes_\Z A$ (the divisible hull of $A$) via the embedding
$x\mapsto 1\otimes x$. We also equip $\Q G$ with the
unique ordering that makes it an ordered abelian group containing $G$ as an ordered subgroup.
Any torsion-free abelian group $A$ embeds as a group into a Hahn product~$H[S,\Q]$ for some ordered set $S$: picking a basis
$S$ of the vector space $\Q A$ over~$\Q$ and an order
on the set $S$ yields group inclusions and a group isomorphism
$$A\ \subseteq\ \Q A\ =\ \bigoplus_{s\in S} \Q s\ \subseteq\ H[(\Q s)]\ \cong\ H[S,\Q].$$
Thus by the example on Hahn products above, every torsion-free abelian group can be expanded to an ordered abelian group.

\index{hull!divisible}
\nomenclature[A]{$\Q G$}{divisible hull $\Q\otimes_\Z G$ of an abelian group $G$}

\subsection*{Convex valuations and archimedean classes}
Let $v\colon G\to S_\infty$ be a surjective valuation on $G$. We say that $v$ is {\bf compatible} with the ordering of $G$ or a {\bf convex valuation} of $G$ if one of the following  equivalent conditions is satisfied:
\begin{list}{*}{\setlength\leftmargin{3em}}
\item[(C1)] for all $a$, $b$, if $0<a<b$, then $va\geq vb$;
\item[(C2)] for every $s\in S$ the subgroup $B(s)=\{a:va>s\}$  of $G$ is convex;
\item[(C3)] for every $s\in S$ the subgroup $\overline{B}(s)=\{a:va\geq s\}$  of $G$ is convex.
\end{list}
Note that then $va\leq vb$ whenever $n\abs{a}\ge \abs{b}$, and
$v(ka)=va$ when $k\in\mathbb Z^{\ne}$. Hence if~$v$ is convex, then 
$v$ has a unique extension to a convex valuation $\Q G\to S_\infty$. 

\index{valuation!convex}
\index{valuation!standard}
\index{ordered abelian group!standard valuation}

\medskip
\noindent
The {\bf archi\-mede\-an class\/} of $a$ is defined by
$$[a]:=\big\{g\in G:
\text{$|a| \leq n|g|$ and 
$|g| \leq n|a|$ for some $n\ge 1$}\big\}.$$
The archimedean classes partition $G$. Each archi\-medean class $[a]$ with $a\neq 0$ is the disjoint union of the two convex sets $[a]\cap G^{<}$ and $[a]\cap G^{>}$.
We order the set 
$$[G]:=\big\{[a]:a\in G\big\}$$ of archimedean classes 
by $$[a]<[b] \quad:\Longleftrightarrow\quad 
\text{$n\abs{a}<\abs{b}$ for all $n\ge 1$}.$$ 
We also write $a=o(b)$ instead of $[a]<[b]$. 
We have $[0]<[a]$ for all $a\in G^{\ne}$, and 
$$[a] \leq [b] \quad\Longleftrightarrow\quad 
\text{$\abs{a}\leq n\abs{b}$ for some $n\ge 1$}.$$ 
We also write $a=O(b)$ instead of $[a]\le [b]$. 
Equipping $[G]$ with the reversed ordering on $[G]$ so as to make $[0]$ its largest element, the map $x\mapsto [x]$ 
becomes a convex valuation of $G$, called the {\bf standard valuation} of $G$.  The surjective convex valuations on $G$ are precisely the coarsenings of 
the standard valuation of $G$. 

\index{archimedean!class}
\nomenclature[I]{$[a]$}{archimedean class of an $a$}
\nomenclature[I]{$[G]$}{set of archimedean classes of  $G$}
\nomenclature[I]{$a=o(b)$}{$[a]<[b]$}
\nomenclature[I]{$a=O(b)$}{$[a]\leq [b]$}

\begin{lemma}\label{lem:Massaza}
Let $v$ be a convex valuation on $G$ such that $S=v(G^{\ne})$ has no largest 
element. Then the $v$-topology on $G$ agrees with the order topology on $G$.
\end{lemma}
\begin{proof}
If $a>0$, then $B(va)\subseteq (-a,a)$. If $s\in S$, then $(-b,b)\subseteq B(s)$ for any $b$ with $vb>s$.
\end{proof}

\noindent
Let $v$ be a convex valuation on $G$ and
$(a_{\rho})$ a pc-sequence in $G$ with respect to $v$. We wish to view the pseudolimits of $(a_{\rho})$ in terms of the ordering of $G$. For simplicity, assume $v(a_{\tau}-a_{\sigma}) > v(a_{\sigma}-a_{\rho})$ for all indices $\tau > \sigma > \rho$, in particular, $a_{\rho}\ne a_{\sigma}$ for all $\rho\ne \sigma$. For each $\rho$, let $s(\rho)$ be the immediate successor of $\rho$ in the set of indices, and set $d_{\rho}:= a_{s(\rho)}-a_{\rho}$. We divide the
set of indices into two disjoint subsets: 
$$L\ :=\ \{\lambda:\ a_{\lambda} < a_{s(\lambda)}\}, \qquad R\ :=\ \{\rho:\ a_{s(\rho)} < a_{\rho}\}.$$
Then $a_{\lambda} < a_{\sigma}< a_{\lambda} + 2d_{\lambda}$ whenever $\lambda\in L$ and $\lambda<\sigma$,
and $a_{\rho} + 2d_{\rho} < a_{\sigma}< a_{\rho}$ whenever $\rho \in R$ and $\rho < \sigma$.
Thus $a_{\lambda} < a_{\rho}$ for all $\lambda\in L$ and $\rho\in R$.
Set 
\begin{align*} P\ &:=\ \{a_{\lambda}:\ \lambda\in L\}\cup \{a_{\rho}+2d_{\rho}:\ \rho\in R\},\\
Q\ &:=\ \{a_{\rho}:\ \rho\in R\}\cup \{a_{\lambda}+2d_{\lambda}:\ \lambda\in L\}.
\end{align*}  

\begin{lemma}\label{pcordgp} We have $P < Q$, and for every $a\in G$,
$$ a_{\rho} \leadsto a\ \Longleftrightarrow\ P < a < Q.$$
In particular, if $G$ has no least positive element, and $\kappa$ is a cardinal such that the set of indices $\rho$ has cardinality $< \kappa$ and  $G$
is $\kappa$-saturated as an ordered set, then~$(a_{\rho})$ has a pseudolimit in $G$. \rm{(See \ref{sec:sat} for the notion of $\kappa$-saturation.)}
\end{lemma}
\begin{proof} We get $P<Q$ from the inequalities already stated.
Suppose $a_{\rho}\leadsto a\in G$. Then for each $\rho$ we have
$v(a-a_{s(\rho)})> vd_{\rho}$, so $a_{\rho} < a < a_{\rho}+2d_{\rho}$
when $\rho\in L$, and $a_{\rho}+2d_{\rho} < a < a_{\rho}$ when $\rho\in R$.
This gives $P < a < Q$. Next, assume $P<a<Q$. Then 
$v(a-a_{\rho})\ge vd_{\rho}$ for all $\rho$, and so 
$a_{\rho}\leadsto a$ by Lemma~\ref{eqpslim}. 
\end{proof}

\noindent
We say that $G$ is {\bf archimedean} \index{ordered abelian group!archimedean} if $[G^{\neq}]:=[G]\setminus \{[0]\}$ \nomenclature[I]{$[G^{\neq}]$}{set of nonzero archimedean classes of $G$} is a singleton, that is, $G\ne \{0\}$, and for all nonzero $a$ and $b$ there is some $n$ such that $\abs{a}\leq n\abs{b}$. Thus if $G$ is archimedean, then there is no nontrivial convex valuation on $G$. Moreover:

\begin{lemma}[H\"older] \label{lem:hoelder}
If $G$ is archimedean and $e\in G^{>}$, then there is a unique embedding of $G$ into the ordered additive group of $\R$ sending $e$ to $1$.
\end{lemma}

\noindent
We leave the easy proof of this fact to the reader.

\begin{example}
Let $(G_s)_{s\in S}$ be a family of ordered abelian groups $G_s\ne \{0\}$ indexed by an ordered set $S$, and $H=H\big[(G_s)\big]$ the  Hahn product of $(G_s)$. Then the valuation~$v$ of $H$ is convex. If every ordered abelian group $G_s$ is archimedean, then~$v$ is equivalent to the standard valuation of $H$. 
\end{example}

\noindent
Let $G_0$ be a subgroup of $G$. We view $G_0$ as an ordered subgroup, and
$[G_0]$ is accordingly identified with an ordered subset of $[G]$ in the obvious way. Then the standard valuation of $G$ restricts to the standard valuation of $G_0$, and if $G_0$ is dense in $G$ (in the order topology on $G$), then $[G_0]=[G]$. Also
$[\Q G]=[G]$.

\begin{lemma}\label{archextclass} Let $G_0$ be a subgroup of $G$ and $a\in G$. Then: \begin{enumerate}
\item[\textup{(i)}] $[a]\notin [G_0]  \Rightarrow  
[G_0 + \Z a]=[G_0] \cup \big\{[a]\big\}, 
[a]\le [g] \text{ for all $g\in (G_0 + \Z a)  \setminus G_0$;}$
\item[\textup{(ii)}] $[G_0+\Z a]\ne [G_0]\, \Rightarrow\, [G_0+\Z a]=[G_0] \cup \big\{[g]\big\} \text{ for some $g\in G_0+\Z a$.}$
\end{enumerate}
\end{lemma}

\noindent
We leave the proof of this lemma to the reader; (ii) follows easily from (i).

\begin{lemma}\label{exarchextclass} Let $C$ be a cut in the ordered set $[G^{\ne}]$. Then there is an ordered group extension $G+\Z x$ of $G$ such that: \begin{enumerate}
\item[\textup{(i)}] $x>0$ and $[x]$ realizes the cut $C$;
\item[\textup{(ii)}] for any ordered abelian group extension $H$ of $G$ and $y\in H^{>}$ such that $y>0$ and $[y]$ realizes $C$ there is a unique ordered group embedding $G+\Z x\to H$ that is the identity on $G$ and sends $x$ to $y$.
\end{enumerate} 
\end{lemma}
\begin{proof} Consider $G$ in the usual way as a subgroup of an abelian group 
$G\oplus \Z x$ with $nx\ne 0$ for all $n\ge 1$. The ordering of $G$ extends uniquely to an abelian group ordering  
on $G\oplus \Z x$ such that for all $g\in G$ and $n\ge 1$, 
$$g+nx>0 \Leftrightarrow \text{$g\ge 0$  or $[g]\in C$,}\qquad g-nx>0 \Leftrightarrow g>0 \text{ and $[g]>C$.}$$
Then $G\oplus \Z x=\G+ \Z x$ with its element $x$ has the desired properties.
\end{proof}

\subsection*{Convex subgroups and rank}
It is easy to see that for a subgroup $H$ of $G$, the following are equivalent:
\begin{enumerate}
\item $H$ is convex;
\item for all $g,h\in G$, if $0\leq \abs{g}\leq \abs{h}$ and $h\in H$, then $g\in H$;
\item for all $g,h\in G$, if $[g]\leq [h]$ and $h\in H$, then $g\in H$.
\end{enumerate}
As a consequence, the map $H\mapsto [H]$ is an inclusion-preserving bijection from the set of convex subgroups of $G$ onto the set of nonempty cuts in~$[G]$.
(In particular, the set of convex subgroups of $G$ is totally ordered by inclusion.)
The {\bf rank} of~$G$, denoted by $\operatorname{rank}(G)$, is defined to be $n$ if there are exactly $n$ nontrivial convex subgroups of~$G$,
and defined to be $\infty$ if there are infinitely many convex subgroups of $G$.
Thus~$G$ has rank $1$ iff $G$ is archimedean.
If $G_0$ is an ordered subgroup of~$G$, then clearly $\operatorname{rank}(G_0)\leq \operatorname{rank}(G)$. 

\index{rank!ordered abelian group}
\index{ordered abelian group!rank}
\nomenclature[I]{$\operatorname{rank}(G)$}{rank of  $G$}

\medskip
\noindent
Given a subgroup $H$ of $G$, the convex hull $\operatorname{conv}(H)$ of $H$ in $G$ is a convex subgroup of~$G$ containing $H$; in fact $\operatorname{conv}(H)$ is the smallest convex subgroup of $G$ containing~$H$.
In particular, for $g\in G$,
$$\operatorname{conv}(\Z g) = \big\{x\in G:[x]\leq [g]\big\}$$ 
is the smallest convex subgroup of $G$ containing $g$. Convex subgroups of $G$ of the form $\operatorname{conv}(\Z g)$, where $g\in G$, are said to be {\bf principal.} Note that if
$\operatorname{rank}(G)< \infty$, then every convex subgroup of $G$ is principal. Thus
$$ \operatorname{rank}(G)< \infty\ \Longleftrightarrow\ [G] \text{ is finite,}$$
and if $[G]$ is finite, then $ \operatorname{rank}(G)= |[G^{\ne}]|$.  

\index{convex!subgroup}
\index{principal!convex subgroup}

\medskip
\noindent
Recall from Section~\ref{sec:modules} that the {\bf rational rank} of $G$, denoted by $\operatorname{rank}_\Q(G)$, is
the dimension of~$\Q G$ as a vector space over $\Q$ (which by convention is $\infty$ 
if this vector space is not finitely generated). Clearly, if $(g_i)_{i\in I}$ is a family in $G^{\ne}$ 
and $[g_i]\ne [g_j]$ for all~$i\ne j$, then $(g_i)$ is $\Q$-linearly independent in $\Q G$.
Thus $\operatorname{rank}(G)\leq \operatorname{rank}_\Q(G)$. 

\index{rank!rational}
\index{rational!rank}
\index{abelian group!rational rank}

\medskip
\noindent
Sometimes we use the following easy observation: 

\begin{lemma}\label{lem:convex hull of sequ}
Let $(s_\rho)$ be a strictly increasing well-indexed sequence in $G^{>}$.
Then the following are equivalent:
\begin{enumerate}
\item[\textup{(i)}] $(s_{\rho})$ is cofinal in some convex subgroup of $G$;
\item[\textup{(ii)}] for each $\rho$ there is $\sigma>\rho$ with $s_{\sigma}\ge 2s_\rho$.
\end{enumerate}
\end{lemma}

\subsection*{Ordered quotient groups}
Let $H$ be a convex subgroup of $G$ with quotient group $\dot{G}:=G/H$ and canonical surjective group morphism $$x\mapsto\dot x:=x+H\ \colon\  G\to\dot G.$$
We equip $\dot G$ with the unique ordering making $\dot G$ an ordered abelian group such that for all $a,b$, if $a\le b$, then $\dot a \leq \dot b$; this is done by declaring $\dot a > \dot 0$ iff $a>H$.
It is easy to verify that for $a,b\notin H$ we have $[a] \leq [b] \Longleftrightarrow [\dot a]\leq [\dot b]$, hence $x\mapsto\dot x$ induces an isomorphism $[G]\setminus [H]\to [\dot G^{\neq}]$ of ordered sets.
Hence 
$$\operatorname{rank}(G) = \operatorname{rank}(H) + \operatorname{rank}(\dot G) .$$

\begin{example}
Suppose $H$, $H'$ are subgroups of $G$ with $H$ convex in $G$ and $G=H\oplus H'$ (internal direct sum of subgroups of $G$).
Then $G$ is ordered by the reverse lexicographic ordering: for $h\in H$, $h'\in H'$, $h+h'\geq 0$ iff $h'>0$ or $h'=0$, $h\geq 0$. The restriction of the morphism $G\to\dot G=G/H$ to $H'$ is an isomorphism $H'\xrightarrow{\cong}\dot G$ of ordered abelian groups, and so
$\operatorname{rank}(G) = \operatorname{rank}(H) + \operatorname{rank}(H')$.
\end{example}

\begin{example}
Suppose $G$ is divisible. Then every convex subgroup of $G$ is also divisible. Hence if $H$ is a convex subgroup of~$G$, then we can take a (divisible) subgroup~$H'$ of~$G$ such that $G=H\oplus H'$ (internal direct sum of subgroups of $G$), and $G$ is ordered by the reverse lexicographic ordering as described in the previous example.
\end{example}

\begin{lemma}\label{lem:uniqueness of ordering} 
Let $v\colon G\to S_\infty$ be a convex valuation on $G$ 
and let $(G',S,v')$ be  an immediate extension of the valued abelian 
group $(G,S,v)$. Then just one ordering of~$G'$ makes $G'$ an ordered 
abelian group extension of $G$ and $v'$ a convex valuation.
\end{lemma}
\begin{proof} Let $a'\in (G')^{\ne}$, and take $a\in G$ such that $a'\sim a$. For an ordering
of~$G'$ as in the lemma, we would have $a'>0\Longleftrightarrow a>0$. 
Conversely, this equivalence
determines an ordering as claimed.
\end{proof}

\begin{lemma}\label{archimedeanquotients} Let $v\colon G\to S_{\infty}$ 
be the standard valuation, so
$S_{\infty}=[G]$ with the reversed ordering. For $s\in S$ we have convex
subgroups $B(s) \subseteq \bar{B}(s)$ of $G$, and $G(s)= \bar{B}(s)/B(s)$
with the quotient ordering is archimedean.
\end{lemma}

\subsection*{Steady functions and slow functions} 
For use in Section~\ref{AbstractAsymptoticCouples} we establish here a useful condition for functions on ordered abelian groups to have the intermediate value property.
Let $v\colon G \to S_{\infty}$ be a convex valuation on $G$. Thus for all 
$x,y\in G$, if $vx>vy$, then $x=o(y)$. For $y\in G$ we let $o_v(y)$ stand for any element $x\in G$ with $vx> vy$, and accordingly, for $x\in G$ we use
$x=o_{v}(y)$ as a suggestive notation for $vx > vy$, sometimes preferred over $x\prec y$, which has the same meaning.  

Let $U$ be a nonempty convex subset of $G$. A function $i\colon U \to G$ is said to be {\bf $v$-steady\/} if $i$ has the intermediate value property and
$i(x)-i(y)=x-y + o_{v}(x-y)$ for all distinct $x,y\in U$. (Note that then $i$ is strictly increasing.) In particular, a~$v$-steady function $G \to G$ is bijective.
If $a\in U$ and the restrictions of $i\colon U \to G$ to~$U^{\le a}$ and~$U^{\ge a}$ are $v$-steady, then $i$ is $v$-steady.
A function $\eta\colon U \to G$
is said to be {\bf $v$-slow on the right\/} if for all $x,y,z\in U$,
\begin{enumerate}
\item[(s1)]$\eta(x)-\eta(y)=o_v(x-y)$ if $x\ne y$;
\item[(s2)] $\eta(y)=\eta(z)$ if $x<y<z$ and $z-y=o_v(z-x)$.
\end{enumerate} 
We define $\eta\colon U \to G$ to be {\bf $v$-slow on the left\/} in the same way except that in clause~(s2) we replace
``$x<y<z$'' by ``$x>y>z$.'' 
%We can now extend Lemma 2.2 in \cite{AvdD}:

\index{map!$v$-steady} \index{v-steady@$v$-steady}
\index{map!$v$-slow} \index{v-slow@$v$-slow}

\begin{lemma}\label{sloste} Suppose $i\colon U \to G$ is $v$-steady and $\eta\colon U \to G$
is $v$-slow on the left or on the right. Then $i+\eta\colon U \to G$ is $v$-steady.
\end{lemma}
\begin{proof} It is clear that $(i+\eta)(x)-(i+\eta)(y)=x-y + o_{v}(x-y)$ for all distinct $x,y\in U$, and thus $i+\eta$ is strictly increasing. To prove that $i+\eta$ has the intermediate value property, assume that $\eta$ is
$v$-slow on the right,
let $a,b\in U$ with $a<b$, put $c:= b-a$ and define $i_1, \eta_1\colon [0,c] \to G$
by $$i_1(x)=i(a+x)-i(a), \qquad \eta_1(x)=\eta(a+x)-\eta(a).$$ 
Then $i_1$ is $v$-steady and $\eta_1$ is $v$-slow on the right,
and it suffices to prove the intermediate value property for $x\mapsto i_1(x)+\eta_1(x)$. So we can assume $U=[0,c]$ and
$i(0)=\eta(0)=0$. Let $0<y<i(c)+\eta(c)$; it suffices to get $x\in (0,c)$ with $i(x)+\eta(x)=y$. Note: $i(c)=c+o_{v}(c),\ \eta(c)=o_{v}(c)$. We distinguish two cases: 

\case[1]{$c-y=o_{v}(c)$.} Then $0<y-\eta(c) < i(c)$, so $y-\eta(c)=i(x)$
with $0<x<c$. It follows easily that $c-x=o_v(c)$, so $\eta(x)=\eta(c)$ by (s2), and thus
$i(x)+\eta(x)=y$.

\case[2]{$c-y\ne o_v(c)$.} Then $0<y<c$ and $0<y-\eta(y) <i(c)$, so $y-\eta(y)=i(x)$
with $0<x<c$. It follows easily that $v(x)=v(y)$, $x-y=o_v(x)$, and
$x-y=o_v(y)$ so $\eta(x)=\eta(y)$ by (s2), and thus
$i(x)+\eta(x)=y$. 
\end{proof}

\begin{lemma}\label{LatticeOrdered}
Suppose the functions $\eta_1,\eta_2\colon U \to G$ are $v$-slow on the right. Then
$-\eta_1, \eta_1+\eta_2, \min(\eta_1,\eta_2)$ are also $v$-slow on the right. 
\end{lemma}
\begin{proof} It is clear that $\eta$ is $v$-slow on the right for $\eta=-\eta_1$ and 
for $\eta=\eta_1 + \eta_2$. Put $\eta:=\min(\eta_1,\eta_2)$. Then (s2) holds. 
Let $a,b\in U$, $a<b$.
If for some $i\in\{1,2\}$, we have $\eta(a)=\eta_i(a)$ 
and $\eta(b)=\eta_i(b)$,
then
$\eta(a)-\eta(b)=\eta_i(a)-\eta_i(b)=o_{v}(a-b)$. 
Suppose 
$\eta_1(a)\leq\eta_2(a)$ and $\eta_2(b)\leq
\eta_1(b)$. If $\eta(a)\geq\eta(b)$, this yields
$$0\leq \eta(a)-\eta(b)=\eta_1(a)-\eta_2(b)\leq
\eta_2(a)-\eta_2(b)=o_{v}(a-b),$$
and if $\eta(a)\leq\eta(b)$, then
$$0\leq \eta(b)-\eta(a)=\eta_2(b)-\eta_1(a)\leq
\eta_1(b)-\eta_1(a)=o_{v}(a-b).$$
The remaining case $\eta_2(a)\leq\eta_1(a)$ and $\eta_1(b)\leq
\eta_2(b)$ follows by symmetry.
\end{proof}

\subsection*{Completeness} {\em In this subsection $G^{>}$ has no smallest element.}\/ 
Note that this condition holds for divisible ordered abelian groups.
An easy consequence of this assumption on $G$ is that
for each $\varepsilon\in G^{>}$ there exists a $\delta\in G^{>}$ such that $2\delta < \varepsilon$. Below we often use this fact without comment.

\medskip\noindent
Let $(a_\rho)$ be a well-indexed sequence in $G$. 
%Write $a_\rho=(a_{1\rho},\dots,a_{n\rho})$ with $a_{i\rho}\in G$ for all $i$, $\rho$.
We say that $(a_\rho)$ is a {\bf cauchy sequence} (or a {\bf c-sequence}) in $G$ if for every $\varepsilon\in G^{>}$ there is $\rho_0$ such that $\abs{a_\rho-a_{\rho'}}<\varepsilon$ for all~$\rho,\rho'>\rho_0$. 
%Here and below
%$$\abs{a}=\max\big(\abs{a_1},\dots,\abs{a_n}\big)\qquad\text{for $a=(a_1,\dots,a_n)\in G^n$.}$$
For $a\in G$ we say that $(a_\rho)$ {\bf converges to $a$} if for each $\varepsilon\in G^{>}$ there is an index $\rho_0$ such that $\abs{a-a_\rho}<\varepsilon$ for all $\rho>\rho_0$; in symbols: $a_\rho\to a$.
%Thus $(a_\rho)$ is a c-sequence iff $(a_{i\rho})$ is a c-sequence for each $i$, 
%and for $a=(a_1,\dots,a_n)%\in G^n$ we have
%$a_\rho\to a$ iff $a_{i\rho}\to a_i$ for each $i$.
We say that~$(a_\rho)$ {\bf converges in $G$} if $a_\rho\to a$ for some $a\in G$. There is at most 
one $a\in G$ with $a_\rho\to a$. %and such $a$ is called the {\bf limit of~$(a_\rho)$.}
If $(a_\rho)$ converges in some ordered abelian group extension~$G'$ of $G$ such that~$G^{>}$ is coinitial in $(G')^{>}$, then $(a_\rho)$ is a c-sequence in $G$. 
If $(a_\rho)$ is a c-sequence in~$G$ with its standard valuation, then $(a_\rho)$ is a c-sequence in the ordered abelian group~$G$; similarly, for $a\in G$, if $a_\rho\to a$ in $G$ with its standard valuation, then  $a_\rho\to a$ in the ordered abelian group $G$. The converses of these implications hold if~$[G^{\neq}]$ with its natural ordering has no smallest element.

\index{cauchy sequence!in an ordered abelian group}
\index{c-sequence!in an ordered abelian group}
\index{sequence!cauchy}
\index{sequence!convergence}
\index{sequence!limit}
\index{limit of a sequence}
\nomenclature[I]{$a_\rho\to a$}{$(a_\rho)$ converges to $a$}

It is easy to check that Lemma~\ref{lem:addition of c-sequences} goes through for this notion of c-sequence in an ordered abelian group. We also have an analogue of Lemma~\ref{lem:cofinality c-sequences} with the same proof except for trivial rewordings: 

\begin{lemma}\label{lem:cofinality c-sequences, 2}
Let $(a_\rho)$ be a c-sequence in $G$ which is not eventually constant. Then the index set of $(a_\rho)$ has cofinality $\operatorname{ci}(G^{>})$.
\end{lemma}

\index{ordered abelian group!dense extension}
\index{dense extension!ordered abelian groups}
\index{extension!ordered abelian groups!dense}

\noindent
A {\bf dense} extension of $G$ is an extension $G'\supseteq G$ of ordered abelian groups such that~$G$ is dense in $G'$ in the order topology; then $G^{>}$ is coinitial in $(G')^{>}$. Thus:

\begin{lemma}\label{densecseq}
Given a dense extension $G'\supseteq G$ and $a\in G'\setminus G$, there is a c-se\-quence~$(a_\rho)$ in $G$, indexed by the ordinals $\rho<\operatorname{ci}(G^{>})$, such that $a_\rho\to a'$.
\end{lemma}  

\noindent
We say that $G$ is {\bf complete} if every c-sequence in $G$ converges in $G$.
A function $f\colon X \to G$ with  $X\subseteq G$ is said to be {\bf uniformly
continuous\/} if for each $\varepsilon\in G^{>}$ there exists a $\delta\in 
G^{>}$ such that for all $x,y\in X$, if $|x-y|< \delta$, then 
$|f(x)-f(y)| < \varepsilon$. Completeness and uniform continuity are
related in the usual way: 

\index{ordered abelian group!complete}
\index{complete!ordered abelian group}
\index{map!uniformly continuous}

\begin{lemma} \label{ucextord} Suppose $G$ is complete and  
$f\colon X \to G$ with  $X\subseteq G$ is uniformly
continuous. Then $f$ extends uniquely to a continuous map from the topological closure of $X$ in $G$ into $G$, and this extension is also uniformly continuous.
\end{lemma}

\noindent
Here is the analogue of Theorem~\ref{thm:completion}:

\begin{theorem}\label{thm:completion oag}
$G$ has a dense complete extension. \textup{(}Hence $G$ is complete iff $G$ has no proper dense extension.\textup{)}
\end{theorem} 

\begin{proof}
We mimic the construction in the proof of Theorem~\ref{thm:completion}. Let $G^{\operatorname{cs}}$ be the set of c-sequences $(a_\rho)$ in $G$ indexed by the ordinals $\rho<\operatorname{ci}(G^{>})$, made into an abelian group via the componentwise addition of such sequences. Then 
$$N\ :=\ \big\{ (a_\rho)\in G^{\operatorname{cs}}:\  a_\rho\to 0\big\}$$
is a subgroup of $G^{\operatorname{cs}}$, and we let $G^{\operatorname{d}}:=G^{\operatorname{cs}}/N$.  
If $(a_\rho)\in G^{\operatorname{cs}}\setminus N$, then for some~$\varepsilon\in G^{>}$ we have
$a_\rho>\varepsilon$ eventually, or $a_\rho<-\varepsilon$ eventually.  Using this fact we make~$G^{\operatorname{d}}$ into an ordered abelian group such that for all $(a_\rho)\in G^{\operatorname{cs}}$, 
$$(a_\rho)+N > 0 \qquad\Longleftrightarrow\qquad \text{for some $\varepsilon\in G^{>}$, $a_\rho > \varepsilon$ eventually.}$$
The map which associates to each $a\in G$ the coset $(a)+N$, where $(a)\in G^{\operatorname{cs}}$ is the constant sequence with value $a$, is an embedding  $G\to G^{\operatorname{d}}$ of ordered abelian groups; we identify $G$ with an ordered subgroup of $G^{\operatorname{d}}$ via this embedding. 

It is routine to show that $G$ is dense in $G^{\operatorname{d}}$, and that for
$(a_\rho)\in G^{\operatorname{cs}}$ we have $a_\rho \to a:=(a_\rho)+N$ in $G^{\operatorname{d}}$.
Thus by Lemma~\ref{lem:cofinality c-sequences, 2}, every c-sequence in $G$ converges in $G^{\operatorname{d}}$. 
As in the proof of Theorem~\ref{thm:completion} one now argues that $G^{\operatorname{d}}$ is complete, using
the previous lemma in place of Lemma~\ref{lem:pc2 analogue}.
\end{proof} 

\noindent
Below $G^{\operatorname{d}}$ is a dense complete extension of $G$ (not necessarily the one constructed in the proof of Theorem~\ref{thm:completion oag}). Using
Lemma~\ref{densecseq}, we have:

\begin{cor} For any dense extension $G'$ of $G$ there exists a 
unique ordered abelian group embedding $G'\to G^{\operatorname{d}}$ over $G$. 
\end{cor}

\noindent
It follows that the ordered abelian group extension $G^{\operatorname{d}}$ of $G$
is determined uniquely (up-to-unique-isomorphism-over $G$) by the property of 
being a dense complete extension of $G$. So there is no harm in referring to 
it as the {\bf completion of $G$.} 

\nomenclature[I]{$G^{\operatorname{d}}$}{completion of the ordered abelian group $G$}
\index{completion!ordered abelian group}

\subsection*{Ordered vector spaces}
An {\bf ordered field} \index{ordered field} \index{field!ordered} is a field $C$ with a (total) ordering $\leq$ of $C$ such that for all $x,y,z\in C$,
$$x\leq y\ \Rightarrow\ x+z\leq y+z,\qquad x\leq y\ \&\ z\geq 0 \ \Rightarrow\ xz\leq yz.$$
(Thus the ordering $\leq$ makes the additive group of $C$ an ordered abelian group.) Just one
ordering makes
the field $\Q$ of rational numbers into an ordered field, and it is given by 
$q\ge 0 \Longleftrightarrow \text{$q=m/n$ for some $m$,~$n$ with $n\ne 0$}$. This is how we view $\Q$ as
an ordered field. If $C$ is an ordered field, then $x^2\ge 0$ for all $x\in C$,
$\operatorname{char}(C)=0$, and we always take $\Q$ as an ordered subfield 
of $C$. Also just one ordering makes
the field $\R$ of real numbers into an ordered field, and it is given by $r\ge 0 \Longleftrightarrow r=x^2 \text{ for some $x\in \R$}$; this is how $\R$ is considered as an ordered field. If $C$ is an ordered field, then $C$ is a 
topological field with respect to the order topology, that is, the 
addition and multiplication maps $C\times C \to C$ 
(with the product topology on $C\times C$), and
the map $x\mapsto x^{-1}\colon C^\times \to C$ are continuous.

\index{ordered vector space} \index{vector space!ordered}

\medskip\noindent
Let $C$ be an ordered field, and let our ordered abelian group $G$ come equipped with an operation $(\lambda,x)\mapsto \lambda x\colon C\times G\to G$ making $G$ a vector space over $C$. Then $G$ is said
to be an {\bf ordered vector space over $C$} if for all $\lambda\in C$ and $x\in G$: 
$$\lambda> 0 \ \&\ x > 0 \quad \Rightarrow \quad \lambda x> 0.$$
Note that then $|\lambda x|=|\lambda|\cdot |x|$ for all $\lambda\in C$ and $x\in G$. 
Every divisible ordered abelian group, viewed as a vector space over $\Q$, is an ordered vector space over the ordered field $\Q$. We view $C$ as an
ordered vector space over $C$ in the obvious way. If~$G$ is a $1$-dimensional
ordered vector space over $C$ and $e\in G^{>}$, then $\lambda \mapsto \lambda e\colon C \to G$ is an isomorphism of ordered vector spaces over $C$.

\medskip\noindent
Let $G$ be an ordered vector space over $C$, and let $P$ be a cut in $G$, 
that is, $P$ is
a downward closed subset of $G$. We set $Q:= G\setminus P$, so $P < Q$. 
Take a vector space extension $G+ Cb\supseteq G$ over 
$C$ with $b\notin G$. There is a unique ordering
on $G+ Cb$ that extends the ordering of $G$ and makes $G+ Cb$
into an ordered vector space over~$C$, such that
$P< b < Q$: in this ordering we have for $\lambda\in C^{>}$ and $g,h\in G$,
$$g+\lambda b< h\Leftrightarrow\ b<\lambda^{-1}(h-g)\ \Leftrightarrow\ \lambda^{-1}(h-g)\in Q.$$
In the next two lemmas we consider $G+ Cb=G\oplus Cb$ as an ordered vector space
over~$C$ in this way. The first lemma is obvious:

\begin{lemma}\label{uniorvec} For any vector $b'$ in any ordered vector space
extension 
$G'\supseteq G$ over~$C$ with $P< b' < Q$, there is a unique embedding 
$G+Cb \to G'$ of ordered vector spaces over~$C$ that is the identity on $G$ 
and sends $b$ to $b'$.
\end{lemma}

\begin{lemma}\label{gendense} $G$ is dense in $G+C b$ if and only if
\begin{enumerate}
\item[\textup{(i)}] for each $\varepsilon \in G^{>}$ there are $p\in P,\ q\in Q$ with
$q-p< \varepsilon$, and 
\item[\textup{(ii)}] $P$ has no largest element and $Q$ has no least element.
\end{enumerate}
\end{lemma}
\begin{proof} It is clear that if $G$ is dense in $G+C b$, then (i) and (ii)
hold. Assume~(i) and~(ii). For each $a\in G$ there is by (ii) an $\varepsilon\in G^{>}$
with $|a-b|>\varepsilon$. For each 
$\varepsilon \in G^{>}$ there is by (i) an $a\in G$ with $|a-b|<\varepsilon$.   
These properties of $b$ are 
inherited by every affine transform $a+\lambda b$ with $a\in G$
and $\lambda\in C^\times$, so
there is for every $\delta\in (G+Cb)^{>}$ an $\varepsilon\in G^{>}$
with $\delta>\varepsilon$. Thus every neighborhood of $b$ in $G+Cb$
contains a vector from $G$; this property of $b$ is inherited by every
vector in $G+C^{\times}b$ and trivially holds for vectors in $G$. Thus
$G$ is dense in $G+C b$.
\end{proof}

\noindent
Suppose $G$ is an ordered vector space over $C$ with $G\ne \{0\}$; so $G^{>}$
has no least element. Given any $\lambda\in C$, 
the map $x\mapsto \lambda x\colon G \to G$ is uniformly continuous, so extends 
uniquely to a continuous map  $G^{\operatorname{d}} \to G^{\operatorname{d}}$
by Lemma~\ref{ucextord}. Denoting the image of any $x\in G^{\operatorname{d}}$
under this map also by $\lambda x$ and varying $\lambda$ we have an
operation  
$$ (\lambda,x)\mapsto \lambda x\ \colon\  C\times G^{\operatorname{d}} \to  
G^{\operatorname{d}}$$ that makes $G^{\operatorname{d}}$ into an ordered vector 
space over $C$ that extends $G$ as an ordered vector space over $C$. 
In particular, the completion of any nontrivial divisible ordered abelian group is divisible.

\subsection*{Valued ordered vector spaces}
Let $C$ be an ordered field. A {\bf valued ordered vector space over $C$} \index{vector space!valued ordered}
is an ordered vector space $G$ over $C$ with a convex valuation 
$v\colon G\to S_\infty$ such that $(G,S,v)$ is a valued vector space over $C$.
Note that then for $s\in S$ we have the convex subspaces 
$B(s)\subseteq\overline{B}(s)$ of $G$, and so $G(s)=\overline{B}(s)/B(s)$ 
with the quotient ordering is an ordered vector space over $C$.

\index{valued ordered vector space}
\index{ordered vector space!valued}

Let $(G_s)$ be a family of ordered vector spaces $G_s\ne \{0\}$ over $C$, 
indexed by the elements $s$ of an 
ordered set $S$. Then the Hahn product $H\big[(G_s)\big]$ 
with its Hahn ordering and its Hahn valuation is a valued ordered 
vector space over~$C$, and this is how such a Hahn product will be construed
as a valued ordered 
vector space over~$C$. In this context we can add something to
Proposition~\ref{prop:Hahn for valued vector spaces}:

\begin{prop}\label{prop:Hahn for valued ordered vector spaces}
Let $G$ be a valued ordered vector space over $C$, and let 
$$i\ \colon\  G \to H:=H\big[(G(s))\big]$$ 
be as in Proposition~\ref{prop:Hahn for valued vector spaces}. Then $i$ 
preserves order: if $g\in G^{>}$, then $i(g)\in H^{>}$.
\end{prop}

\noindent
%Proposition~\ref{prop:Hahn for valued ordered vector spaces}, 
With this, the embedding of an ordered abelian group into its divisible hull 
and Lemmas~\ref{lem:hoelder} and~\ref{archimedeanquotients} yield the classical Hahn Embedding Theorem:

\begin{cor}\label{cor:Hahn classical} Any ordered abelian group $G$ embeds 
into the Hahn ordered group $H[S,\R]$ where $S=[G^{\neq}]$ with
reversed ordering.
\end{cor}

\noindent
Here is a complement to Lemma~\ref{lem:uniqueness of ordering}:

\begin{cor} \label{cor:ordering in immediate extensions}
Let $G$ be a valued ordered vector space over $C$ and~${G'\supseteq G}$ 
an immediate extension of valued vector spaces, with valuation $v'$. Then 
just one ordering on $G'$ makes $G'$ an ordered abelian group extension 
of~$G$ with convex $v'$. With~$v'$ and this ordering, $G'$ is a valued 
ordered vector space over $C$.
%if $G$ is dense in $G'$ in the valuation topology, then $G$ is dense in 
%$G'$ in the order topology.
\end{cor} 

\noindent
This follows easily by inspecting the proof of Lemma~\ref{lem:uniqueness of ordering}.

\subsection*{The $C$-valuation of an ordered vector space over $C$}
Let $G$ be an ordered vector space over the ordered field $C$. The {\bf $C$-archi\-mede\-an class\/} of $a\in G$ \index{archimedean!class}\index{C-archimedean class@$C$-archimedean class}\nomenclature[I]{$[a]_C$}{$C$-archimedean class of $a$}  is
$$[a]_C\ :=\ \big\{g\in G:\ 
\text{$\abs{a} \leq \lambda\abs{g}$ and 
$\abs{g} \leq \lambda\abs{a}$ for some $\lambda\in C^{>}$}\big\}.$$
For $G$ as an ordered vector space over the ordered subfield $\Q$ of $C$ we get $[a]_\Q=[a]$ for $a\in G$, where $[a]$ is the archimedean class of $a$ as defined earlier in this section. 
Let $[G]_C$ be the set of $C$-archimedean classes. \nomenclature[I]{$[G]_C$}{set of all $C$-archimedean classes of $G$}\nomenclature[I]{$[G^{\neq}]_C$}{set of all nonzero $C$-archimedean classes of $G$} Then $[G]_C$ is a partition
of $G$, and we linearly order $[G]_C$ by
\begin{align*}
[a]_C<[b]_C &\quad:\Longleftrightarrow\quad 
\text{$\lambda\abs{a}<\abs{b}$ for all $\lambda\in C^{>}$} \\
&\quad\hskip0.3em \Longleftrightarrow\quad \text{$[a]_C\neq [b]_C$ and $\abs{a}<\abs{b}$.}
\end{align*} 
Thus $[0]_C=\{0\}$ is the smallest $C$-archimedean class. The map 
$$x\mapsto [x]_C\ :\  G \to [G]_C$$ with the reversed ordering
on $[G]_C$ (so with $[0]_C$ as the largest element)
is a convex valuation of $G$ making $G$ a valued ordered vector space over $C$.
This valuation is
called the {\bf $C$-valuation} of the ordered vector space $G$, and is a coarsening of
the standard valuation of the ordered additive group of $G$ defined earlier in this section.
%~\ref{sec:oag}.
%The valuations on $G$ making $G$ a valued ordered vector space 
%over $C$ are precisely the coarsenings of the $C$-valuation of $G$.
A subspace $G_0$ of~$G$ is considered as an ordered vector space over $C$ by restricting the ordering to $G_0$, and $[G_0]_C$ then becomes an ordered subset of $[G]_C$ via the obvious identification; the $C$-valuation of $G$ then restricts to the $C$-valuation of $G_0$.

\index{C-valuation@$C$-valuation}
\index{valuation!$C$-valuation}
\index{ordered vector space!$C$-valuation}
\index{ordered vector space!$C$-archimedean}

We say that $G$ is {\bf $C$-archimedean} if $[G^{\neq}]_C:=[G]_C\setminus \big\{[0]\big\}$ is a singleton.
Viewing~$G$ as a valued ordered vector space over $C$ by the $C$-valuation, and given  $s\in [G^{\neq}]_C$, the 
ordered vector space $G(s)=\overline{B}(s)/B(s)$ over $C$ is $C$-archimedean.

\begin{example}
Let $(G_s)$ be a family of $C$-archimedean ordered vector spaces over $C$
indexed by the elements of an ordered set $S$. 
Then the Hahn valuation of the Hahn product $H:=H\big[(G_s)\big]$ is 
equivalent to the $C$-valuation of $H$.
\end{example} 

\noindent
For use in the proof of Lemma~\ref{newprop, lemma 1} we record the following:

\begin{lemma}\label{divimm, standard val}
Let $(a_\rho)$ be a divergent pc-sequence in $G$ with respect to the $C$-valuation of $G$. Then $G$ has an extension $G\oplus Ca$ as an ordered vector space over~$C$ such that $[G\oplus Ca]_C=[G]_C$, and $a_\rho\leadsto a$ with respect to the 
$C$-valuation of~${G\oplus Ca}$. 
\end{lemma}
\begin{proof}
Lemma~\ref{divimm} and Corollary~\ref{cor:ordering in immediate extensions} 
provide an immediate
valued ordered vector space extension $G\oplus Ca$ of $G$ with 
$a_\rho\leadsto a$. It follows easily that the valuation on~${G\oplus Ca}$ 
is the $C$-valuation.
\end{proof}

\noindent
Let $G\oplus Ca$ be as in Lemma~\ref{divimm, standard val}, and let
$G\oplus Cb$ also be as in that lemma, with~$b$ instead of $a$.
Then there is an isomorphism $G\oplus Ca\to G\oplus Cb$
of ordered vector spaces over $C$ which is the identity on $G$ and 
sends $a$ to $b$.

\subsection*{Completion of valued ordered vector spaces}
Let $G$ be a valued ordered vector space over an ordered field $C$ 
such that the value set $S=v(G^{\ne})$ is not empty and has no largest element.
Then the valuation topology and the order topology on $G$ coincide, by
Lemma~\ref{lem:Massaza}. Thus  
$G$ is complete as a valued abelian group if and only if 
$G$ is complete as an ordered abelian group. 

Recall that $G^{\operatorname{c}}$ is the completion of $G$ as a valued abelian 
group. We construe it as a valued vector space 
over $C$ as in Section~\ref{sec:valued vector spaces}. 
Corollary~\ref{cor:ordering in immediate extensions} then gives a
unique ordering on $G^{\operatorname{c}}$ making it a valued ordered 
vector space over $C$ extending the valued order vector space $G$. 
With this ordering $G^{\operatorname{c}}$ is also complete 
as an ordered abelian group, and
$G$ is dense in $G^{\operatorname{c}}$ in the order topology. 

On the other hand, $G^{\operatorname{d}}$ is the completion of $G$ as an ordered abelian group, and we 
construe it here as an ordered vector space over the ordered field $C$ as
indicated at the end of the subsection on ordered vector spaces. 
Thus by the uniqueness property of $G^{\operatorname{d}}$ we have:

\begin{cor}\label{lem:Gc vs Gd}
There is a unique isomorphism $G^{\operatorname{c}}\to G^{\operatorname{d}}$ of ordered vector spaces over $C$ which is the identity on $G$.
\end{cor}

\subsection*{Hahn spaces over ordered fields}
Let $C$ be an ordered field and $G$ an ordered vector space over $C$. Call $G$ a {\bf Hahn space} over $C$ if $G$ with its 
$C$-valuation is a Hahn space over $C$ as defined in Section~\ref{sec:valued vector spaces}, that is, for all $a,b\in G^{\neq}$ with $[a]_C=[b]_C$ there exists 
$\lambda\in C^\times$ such that $[a-\lambda b]_C<[a]_C$. 

\index{Hahn space}

\begin{examples}
\mbox{}

\begin{enumerate}
\item Any $1$-dimensional ordered vector space over $C$ is a Hahn space over $C$.
\item The ordered vector space $\Q+\Q\sqrt 2 \subseteq \R$ over $\Q$ is not a
Hahn space over $\Q$.
\item Any ordered vector space over the ordered field $\R$ is a Hahn space over
$\R$.
\item If $(G_s)$ is a family of Hahn spaces $G_s\ne \{0\}$ over $C$ indexed by the elements~$s$ of an ordered set, then the Hahn product $H\big[(G_s)\big]$ as a vector
space over $C$ with its Hahn ordering is a Hahn space over $C$.
\end{enumerate}
\end{examples}

\noindent
Using Proposition~\ref{prop:Hahn for valued ordered vector spaces} we obtain 
easily: 

\begin{cor}
If $G$ is a Hahn space over $C$ equipped with its $C$-valuation, then there is an embedding $i\colon G\to H[S,C]$ of valued ordered vector spaces over $C$  such that $H[S,C]$ is an immediate extension of its valued subspace $iG$.
\end{cor}

\subsection*{Notes and comments}
The fact that torsion-free abelian groups are orderable is in
\cite{Levi}.
Lemma~\ref{lem:Massaza} is in~\cite{Massaza},
Lemma~\ref{lem:hoelder} is from~\cite{Hoelder}. 
The archimedean property was isolated by Stolz~\cite{Stolz}. (See \cite{Ehrlich-1,Fisher}.)
The material on steady and slow functions extends a result in \cite{AvdD}.
The completion of an ordered abelian group can also be obtained by realizing 
certain kinds of Dedekind cuts; see~\cite{Cohen-Goffman}. (The particular cuts required were 
already considered by Veronese~\cite{Veronese}.) Completing 
valued abelian groups and ordered abelian group 
are special cases of completing a uniform space; see~\cite[Section~II.3.7]{Bourbaki-Top} and~\cite[Lemmas~II.4.8 and III.4.5]{PC}.
%Ehrlich~\cite{Ehrlich} argues that proper dedekind cuts 
%should be called ``veronese cuts'' in honor of Veronese's 
%work~\cite{Veronese} on non-archimedean geometry.
The notion of ``Hahn space over an ordered field'' is from~\cite{AvdD}.
For more on ordered algebraic structures, see~\cite{Fuchs, PC}.

%% file: mt-3.tex
\chapter{Valued Fields}\label{ch:valuedfields}

\noindent
In this chapter we assume familiarity with basic field theory, including
the rudiments of the theory of ordered fields (which are summarized, without proofs, in Section~\ref{sec:valued ordered fields} below).
In Section~\ref{sec:valued fields} 
we take up valued fields and establish their basic properties, focusing in particular on extensions of valued fields.  Next, in Section~\ref{sec:pc in valued fields} we study  pseudocauchy sequences in valued fields;
these sequences will be needed later in constructing solutions of
algebraic differential equations in immediate extensions of suitable 
valued differential fields.  
%elementarily equivalent to such as $\mathbb T$ and in suitable extensions. 
In Section~\ref{sec:henselian valued fields} we consider 
{\em henselian\/} valued fields. 
(In Chapter~\ref{sec:dh1} we study for valued differential fields the analogous notion of {\em differential-henselian}.) 
Whenever it simplifies matters we restrict our attention 
in this section to valued fields whose residue field is of 
characteristic zero, since this is the only case that arises later.
In Section~\ref{sec:decomposition} we show how to decompose a 
valuation on a field into simpler ones. This leads to a 
study of various special types of pseudocauchy sequences. The valuation 
of $\mathbb T$ is compatible with its natural ordering.
Consequently, Section~\ref{sec:valued ordered fields} of this chapter contains basic facts about fields with compatible ordering and valuation. 
In Section~\ref{sec:modth val fields} we review some basic model theory of valued fields.
In the final Section~\ref{Newton Diagrams} we consider the Newton diagram and Newton tree of a polynomial over a valued field.

\input{mt-3-1}

\input{mt-3-2}

\input{mt-3-3}

\input{mt-3-4}

\input{mt-3-5}

\input{mt-3-6}

\input{mt-3-7}

%% file: mt-3-1.tex
\section{Valuations on Fields}\label{sec:valued fields}

\noindent
Let $A$ be an integral domain.
A {\bf valuation}  on $A$ is a map $v\colon A^{\ne} \to \Gamma$, where $\Gamma$ is an ordered abelian group, 
such that for all $x,y\in A^{\ne}$:
\begin{list}{*}{\setlength\leftmargin{3em}}
\item[(V1)] $v(xy)\ =\ v(x)+v(y)$;
\item[(V2)] $v(x+y)\ \geq\ \min(v(x),v(y))$ when $x+y\ne 0$.
\end{list}
Let $v\colon A^{\ne}\to\Gamma$ be a valuation on $A$. Then $v1=v(-1)=0$, since $v$ restricted to the group $A^\times$ of units of $A$ is a group morphism. Hence $vx=v(-x)$ for all $x\in A^{\ne}$. By convention we extend $v$ to 
$v\colon A \to \Gamma_\infty$ by $v(0):= \infty$, which makes $v$ a 
valuation on the additive group of $A$ as defined in Section~\ref{sec:valued abelian gps}.
We can extend $v$ uniquely to a valuation $v\colon K^\times\to\Gamma$ on the fraction field $K$ of $A$, by
$$v(x/y)\ =\ vx-vy\qquad (x,y\in A^{\neq}).$$
Thus $v(K^\times)$ is a subgroup of $\Gamma$.
When we refer to a valuation $v\colon K^\times \to\Gamma$ on a field~$K$, we 
assume from now on that $v(K^\times)=\Gamma$, unless specified otherwise. We call $\Gamma=v(K^\times)$ the {\bf value group} of $v$. The valuation $v$ on $K$ is trivial iff $\Gamma=\{0\}$.

\index{valuation!on an integral domain}
\index{value group!valuation}
\index{valuation!value group}
\nomenclature[K]{$\smallo=\smallo_K$}{maximal ideal of $\mathcal O$}

\medskip
\noindent
Let $v\colon K^\times\to\Gamma$ be a valuation on the field $K$.
We set
$$\mathcal O\	:=\ \{ x\in K:\ vx\geq 0\},\ \qquad
\smallo\ 		:=\  \{ x\in K:\  vx>0\}.$$
If we need to indicate the dependence on $v$, we write $\mathcal{O}_v$ and $\smallo_v$
instead of $\mathcal{O}$ and $\smallo$. Thus $\mathcal O$ is a subring of $K$ and $\smallo$ is an ideal of $\mathcal O$.  In fact, $vx=0\Longleftrightarrow v(x^{-1})=0$, for $x\in K^\times$, so $\mathcal O^\times=\mathcal O\setminus\smallo$. 
Note that for $x,y\in \mathcal{O}$ we have $x\in \mathcal{O}y$ iff $vx\ge vy$.
Therefore $\mathcal O$ is a valuation ring as defined in Section~\ref{sec:ztn}, and in particular a local ring with maximal ideal~$\smallo$.
We call~$\mathcal O$ the {\bf valuation ring} \index{valuation!valuation ring}\index{valuation ring!valuation on a field} of $v$.
The {\bf residue field} \index{residue field!valuation on a field}\index{valuation!residue field} of $v$ is the field
$\k_v:=\mathcal{O}/\smallo$. 
With the notation of Section~\ref{sec:valued abelian gps} and working in
the valued abelian group $(K,\Gamma,v)$ we have 
$\mathcal{O}=\overline{B}(0)$, $\smallo=B(0)$,
$K(0)=\mathcal{O}/\smallo$, and for each $a\in K^\times$ the group isomorphism
$$x\mapsto xa\ \colon\ \overline{B}(0)\to \overline{B}(va)$$ 
induces a group isomorphism $K(0)\to K(va)$.
The $v$-topology on $K$ makes $K$ into a hausdorff topological field.

\subsection*{Local rings}
Recall that a commutative ring is said to be {\bf local} if it has exactly one maximal ideal. Thus a commutative ring $A$ is local iff $A\setminus A^\times$ is an ideal of $A$, and in this case, $\mathfrak m_A:=A\setminus A^\times$ is the (unique) maximal ideal of $A$. We let $\mathbf k_A:=A/\mathfrak m_A$ be the {\bf residue field} of the local ring~$A$. 
Given local rings $A$, $B$, we say that {\bf $B$ lies over $A$} if $A$ is a subring of $B$ and $\mathfrak m_A\subseteq\mathfrak m_B$ (and hence $\mathfrak m_B\cap A=\mathfrak m_A$). In this case we have an induced embedding $\mathbf k_A=A/\mathfrak m_A\to B/\mathfrak m_B=\mathbf k_B$ of residue fields, by means of which we identify $\mathbf k_A$ with a subfield of $\mathbf k_B$. We then have a commutative diagram
$$\xymatrix@R-0.5em@C+1em{
B \ar[r]      				& \k_B \\
A \ar[r] \ar[u]^\subseteq	& \k_A \ar[u]^\subseteq
}$$
where the horizontal arrows are the residue morphisms.

\index{local ring}
\index{local ring!residue field}
\index{local ring!lying over}
\index{ring!local}
\index{residue field!local ring}
\index{valuation ring!field}

\medskip\noindent
A {\bf valuation ring of a field $K$\/} is a subring $A$ of $K$ such that for each $x\in K^\times$ we have $x\in A$ or $x^{-1}\in A$. Each valuation ring of a field is clearly a valuation ring as defined in Section~\ref{sec:ztn}. 
Let $K\subseteq L$ be a field extension, and let $A, B$ be valuation rings of  $K$ and
$L$, respectively. Then $B\cap K$ is a valuation ring of $K$, and 
$$\text{$B$ lies over $A$}\ \Longleftrightarrow\ A=B\cap K.$$

\subsection*{Correspondence between valuation rings and valuations}
Let $K$ be a field. The valuation ring of a valuation on $K$ is a valuation ring of $K$. Conversely, to a valuation ring $A$ of $K$ we associate a valuation on $K$ as follows: Consider the (abelian) quotient group $\Gamma_A=K^\times/A^\times$, written additively. The binary relation $\geq$ on $\Gamma_A$ defined by
$$xA^\times \geq yA^\times\ :\Longleftrightarrow\ x/y\in A\qquad (x,y\in K^\times)$$
makes $\Gamma_A$ into an ordered abelian group, and the natural map
$$v_A\colon K^{\times}\to \Gamma_A,\qquad v_A x := 
xA^\times,$$
is a valuation on $K$. The valuation ring of $v_A$ is $A$, and
every valuation on $K$ with valuation ring $A$ is equivalent to $v_A$, 
as defined in Section~\ref{sec:valued abelian gps}. In fact, for valuations
$v\colon K^\times \to \Gamma=v(K^\times)$ and $v'\colon K^\times\to \Gamma'=v'(K^\times)$ on 
the field $K$ we have: $v$ and~$v'$ are equivalent as defined in Section~\ref{sec:valued abelian gps}, if and only if $v$ and $v'$ have the same valuation ring, and also if and only if there is an isomorphism $i\colon \Gamma \to \Gamma'$ of
ordered abelian groups such that $v'=i\circ v$.

\medskip
\noindent
Let $K\subseteq L$ be a field extension, $A$ a valuation ring of $K$, and $B$ be a valuation ring of $L$ lying over $A$.
Then $\mathfrak m_A=\mathfrak m_B\cap K$, so we have an induced embedding $\mathbf k_A=A/\mathfrak m_A\to B/\mathfrak m_B=\mathbf k_B$ of residue fields, by means of which we identify 
$\mathbf k_A$ with a subfield of $\mathbf k_B$. Also, $B^\times\cap K=A^\times$, so we have an induced embedding $v_A(x)\mapsto v_B(x)\colon\Gamma_A\to\Gamma_B$ ($x\in K^\times$) of ordered abelian groups, by means of which~$\Gamma_A$ is identified with an ordered subgroup of $\Gamma_B$. Then $v_A(x)=v_B(x)$ for all $x\in K$.

\label{p:valued fields}

\medskip
\noindent
A {\bf valued field} is just a field equipped with one of its valuation rings.
Let $K$ be a valued field. Unless we specify otherwise 
we let $\mathcal O$ be the distinguished
valuation ring of $K$, with $v$ the corresponding valuation on $K$, and 
$\smallo$ the maximal ideal of $\mathcal O$.
 The {\bf value group} of $K$ is $\Gamma=v(K^\times)$. If we need to indicate
the dependence on~$K$ we attach a subscript $K$, so 
$\mathcal{O}=\mathcal{O}_K$, $v=v_K$, and so on. \nomenclature[K]{$\mathcal O=\mathcal O_K$}{valuation ring of a valued field $K$}\nomenclature[K]{$v=v_K$}{valuation of a valued field $K$} Sometimes (especially when $K$ is an asymptotic differential field, 
see Section~\ref{As-Fields,As-Couples})
we prefer more relational notation: 
define the
binary relations $\preceq$ and $\prec$ on $K$ by
$$ a\preceq b\ :\Longleftrightarrow\ va \ge vb, \qquad a\prec b\ :\Longleftrightarrow\ va > vb,$$
just as we already did for valued abelian groups. Accordingly, 
$$K^{\preceq 1}\ :=\ \{a\in K:\ a\preceq 1\}$$ is then the valuation ring of $K$, and $K^{\prec 1} :=\{a\in K:\ a\prec 1\}$ is the maximal ideal 
of~$K^{\preceq 1}$. Likewise, $K^{\succ 1}\ :=\ \{a\in K:\ a\succ 1\}$.
The {\bf residue field} of $K$ is 
$$\res(K)\ :=\ \mathbf k_{\mathcal O}\ =\ \mathcal{O}/\smallo\ =\ K^{\preceq 1}/K^{\prec 1},$$
and for $a\in \mathcal{O}$ we let $\res(a)$, or sometimes $\overline{a}$, be the
residue class $a+\smallo\in \res(K)$. 
The residue fields of the valued  fields of interest in our work (such as $\mathbb T$) are of characteristic zero. 
One says that $K$ is of {\bf equicharacteristic zero} if $\res(K)$ has characteristic zero.
Note that then $K$ also has characteristic zero, and the
valuation ring $\mathcal O$ of $K$ contains a subfield (i.e., a subring that happens to be a field):
the unique ring morphism $\Z\to\mathcal O$ is injective, and identifying $\Z$ with its image in $\mathcal O$ under this morphism, each nonzero integer is a unit in $\mathcal O$, so the fraction field $\Q\subseteq K$ of $\Z$ is
contained in $\mathcal O$.
Note that if $\mathcal O$ contains a subfield $C$, then $K$ as a vector space over~$C$ with the valuation $v$ is a valued vector space over $C$ (as defined in Sec\-tion~\ref{sec:valued vector spaces}).

\index{valued field}
\index{field!valued}
\index{value group!valued field}
\index{valued field!value group}
\index{valued field!valuation ring}
\index{valued field!residue field}
\index{valuation ring!valued field}
\index{residue field!valued field}
\index{valued field!equicharacteristic zero}
\nomenclature[K]{$K^{\preceq 1}$}{valuation ring of a valued field $K$}
\nomenclature[K]{$K^{\succ 1}$}{the complement of $K^{\preceq 1}$ in a valued field $K$}
\nomenclature[K]{$K^{\prec 1}$}{maximal ideal of the valuation ring of a valued field $K$}
\nomenclature[K]{$\res(K)$}{residue field of a valued field $K$}
\nomenclature[K]{$\res(a)=\overline{a}$}{residue class of $a\in K^{\preceq 1}$}
\nomenclature[K]{$\Gamma=\Gamma_K$}{value group of a valued field $K$}

\medskip
\noindent
Let $K$ be a valued field.
A {\bf valued field extension} of $K$ is a field extension~$L$ of~$K$ equipped with a valuation ring of $L$ that lies over $\mathcal O$;
in this situation we also call $K$ a {\bf valued subfield} of $L$.
Thus a 
valued field extension $K\subseteq L$ gives rise to a field extension
$\res(K) \subseteq \res(L)$ and to an ordered group extension 
$\Gamma \subseteq \Gamma_L$.

\index{extension!valued fields}
\index{valued field!extension}
\index{valued field!valued subfield}

\subsection*{Correspondence between dominance relations and valuations}
The binary relation $\preceq$ introduced above is an example of a dominance relation:

\begin{definition}
A {\bf dominance relation} \index{dominance relation}\nomenclature[R]{$f\preceq g$}{$f$ is dominated by $g$} on a field $K$ is a binary relation 
$\preceq$ on~$K$ such that for all $f,g,h\in K$:
\begin{list}{*}{\setlength\leftmargin{1em}}
\item[(D1)] $1\not\preceq 0$;
\item[(D2)] $f\preceq f$;
\item[(D3)] $f\preceq g \text{ and } g\preceq h \Rightarrow f\preceq h$;
\item[(D4)] $f\preceq g$ or $g\preceq f$;
\item[(D5)] $f\preceq g \Leftrightarrow fh\preceq gh$, provided $h\ne 0$;
\item[(D6)] $f\preceq h \text{ and } g\preceq h \Rightarrow f+g \preceq h$.
\end{list}
If $f\preceq g$, we say that $f$ is {\bf dominated} by $g$, or $g$
{\bf dominates} $f$. 
\end{definition}

\noindent
Thus, if $v$ is a valuation on $K$ with valuation ring $K^{\preceq 1}$, we obtain a
dominance relation on $K$ by setting, for $f,g\in K$:
\begin{equation}\label{domval}
f\preceq g \ :\Longleftrightarrow\  vf\geq vg \ \Longleftrightarrow\ 
\text{$f=gh$ for some $h\in K^{\preceq 1}$.}
\end{equation}
Conversely, if $\preceq$ is a dominance relation on $K$, then clearly
$K^{\preceq 1}:=\{f\in K:f\preceq 1\}$ is a valuation ring of $K$,
and if $v$
denotes the corresponding valuation on $K$, 
then the equivalence \eqref{domval}
holds, for all $f,g\in K$. We call $v$ the valuation {\bf associated} to
the dominance relation $\preceq$. This yields a one-to-one
correspondence between dominance relations on $K$ and valuation rings of $K$.
(That is why the reverse of a dominance relation is sometimes called a
{\em valuation divisibility}.) We shall use the valuation and 
dominance terminologies interchangeably and switch between them without further comment.

\index{valuation!associated to a dominance relation}
\index{dominance relation!associated valuation}

\begin{notations}
Let $K$ be a field with 
dominance relation $\preceq$ on $K$, let
$v\colon K^{\times} \to\Gamma$ be the associated valuation, and let 
$f,g$ denote elements of $K$. We define, just as we did for valued abelian groups in Section~\ref{sec:valued abelian gps}: \nomenclature[R]{$f\prec g$}{$f$ is strictly dominated by $g$}
$$f\prec g \quad :\Longleftrightarrow\quad vf>vg 
\quad \Longleftrightarrow\quad \text{$f\preceq g$ and 
$g\not\preceq f$.}$$
If $f\prec g$, we say that  $f$ is {\bf strictly dominated} by $g$.
If $f\preceq g$ and $g\preceq f$, then we say that $f$ and $g$ are
{\bf asymptotic}, written as $f\asymp g$, \nomenclature[R]{$f\asymp g$}{$f$ and $g$ are asymptotic} and if $f-g\prec f$, then~$f$ and $g$ are said to be {\bf equivalent}, written as $f\sim g$. \nomenclature[R]{$f\sim g$}{$f$ and $g$ are equivalent} The relations $\asymp$ and $\sim$ are equivalence relations on~$K$ and~$K^\times$,
respectively. Note that
\begin{align*}
f\asymp g\quad &\Longleftrightarrow\quad vf=vg, \\ 
f\sim g  \quad &\Longleftrightarrow\quad v(f-g)>vf.
\end{align*}
In particular, $f\sim g \Rightarrow f\asymp g$. Define
\begin{align*}
f\preceq 1 &\quad:\Longleftrightarrow\quad \text{$f$ is {\bf bounded},} \\
f\prec   1 &\quad:\Longleftrightarrow\quad \text{$f$ is 
{\bf infinitesimal},}\\
f\succ   1 &\quad:\Longleftrightarrow\quad \text{$f$ is {\bf infinite}.}
\end{align*}
The elements of the maximal ideal $K^{\prec 1}$ of the valuation ring 
$K^{\preceq 1}$ are exactly the infinitesimals of $K$.
\end{notations}

\nomenclature[R]{$f\preceq 1$}{$f$ is bounded}
\nomenclature[R]{$f\prec 1$}{$f$ is infinitesimal}
\nomenclature[R]{$f\succ 1$}{$f$ is infinite}
\index{bounded}
\index{infinitesimal}
\index{infinite}

%\noindent
%If $K\subseteq L$ is a field extension, $\preceq$ a
%dominance relation on $K$ and $\preceq_L$ a dominance relation on $L$
%such that $f\preceq g \Longleftrightarrow f\preceq_L g$ for all
%$f,g\in K$, then $(K,\preceq)\subseteq (L,\preceq_L)$ is called an {\bf
%extension} of fields with dominance relation. (Usually the dominance relation 
%on $L$ is also denoted just by the symbol $\preceq$.)

\subsection*{Well-based series} Let $\frak{M}$ be an ordered set whose ordering is $\preceq$. We think of the elements of $\mathfrak{M}$ as {\em monomials\/} and accordingly denote its elements by $\fm$, $\fn$, and so on. A set $\mathfrak{S}\subseteq \frak{M}$ 
is said to be {\bf well-based}\/
if it is well-ordered for the reverse ordering~$\succeq$, that is, 
there is no strictly increasing infinite sequence 
$\fm_0 \prec \fm_1 \prec \frak{m}_2 \prec \cdots$ in $\mathfrak{S}$.
Clearly the union of two well-based subsets of $S$ is well-based. 

Next, let $C$ be an (additive) abelian group whose elements are to be thought of as coefficients. In this spirit, a function $f\colon \frak{M} \to C$ will be denoted as a series
$\sum_{\mathfrak{m} \in \mathfrak{M}} f_{\mathfrak{m}} \mathfrak{m}$,
with $f_{\mathfrak{m}}=f(\mathfrak{m})$, with
{\bf support\/} 
$\supp f := \{ \mathfrak{m} \in \mathfrak{M} : f_{\mathfrak{m}}
\neq 0 \}$. Then 
\[ C [[ \mathfrak{M} ]]\ :=\ \left\{ f\colon  \frak{M} \to C :\  
   \text{$\supp f \subseteq \mathfrak{M}$ is well-based} \right\} \]
is a subgroup of the additive group of all functions $\mathfrak M\to C$ with pointwise addition. 
Indeed, $C[[\mathfrak M]]$ is just the Hahn product $H[\mathfrak{M},C]$ of Section~\ref{sec:valued abelian gps}, with respect to the reverse
ordering of~$\mathfrak{M}$. For nonzero $f \in C [[ \mathfrak{M} ]]$ we define 
\[ \mathfrak{d} ( f )\ :=\  \max_{\preccurlyeq} \supp f, \]
the {\bf dominant monomial} of $f$. 
\index{series!well-based}
\index{well-based series}
\index{well-based series!support}
\index{well-based series!dominant monomial}
\index{support!well-based}
\index{dominant!monomial}
\index{monomial!dominant}
\nomenclature[K]{$\frak d(f)$}{dominant monomial of the well-based series $f$}
\nomenclature[K]{$\supp(f)$}{support of the series $f$}

\medskip\noindent
{\em In the rest of this subsection $\frak{M}$
is an ordered abelian group and $C$ is a \textup{(}coefficient\textup{)} field.}
Thus  $C [[ \frak{M} ]] $
is a subspace of the $C$-vector space of all functions $\mathfrak M\to C$. We take $\frak{M}$ as a {\em multiplicative\/} group, in view of the role of its elements as monomials.
We leave the proof of the next result to the reader.

\begin{lemma}\label{lem:neumann}
Let $\mathfrak S_1,\mathfrak S_2\subseteq\mathfrak M$ be well-based. Then 
for each $\mathfrak m\in \mathfrak M$ there are only finitely many $(\mathfrak m_1,\mathfrak m_2)\in \mathfrak S_1\times \mathfrak S_2$ such that $\mathfrak m=\mathfrak m_1\cdot\mathfrak m_2$, and the set $$\mathfrak S_1\cdot \mathfrak S_2\ :=\ \{\mathfrak m_1\cdot\mathfrak m_2:\ \mathfrak m_1\in \mathfrak S_1,\ \mathfrak m_2\in \mathfrak S_2\}$$ 
is well-based.
\end{lemma}

\noindent
Lemma~\ref{lem:neumann} gives a binary operation on $C[[\mathfrak M]]$ by 
$$f\cdot g\ :=\ \sum_{\mathfrak m\in\mathfrak M} \left(\sum_{\mathfrak n_1\cdot\mathfrak n_2=\mathfrak m} f_{\mathfrak n_1}g_{\mathfrak n_2}\right)\mathfrak m.$$
With this operation as multiplication, $C[[\mathfrak M]]$ is a domain, with subfield $C$ via the identification $c\mapsto f$ with
$f_{1}=c$ and $f_{\mathfrak{m}}=0$ for all $\mathfrak{m}\ne 1$.
We identify the
group~$\mathfrak{M}$ with a subgroup of~$C[[\mathfrak M]]^{\times}$ via
$\mathfrak m \mapsto f$, with $f_{\mathfrak m}=1$ and $f_{\mathfrak{n}}=0$ for all
$\mathfrak{n}\ne \mathfrak{m}$.

\medskip\noindent
Let $\Gamma$ be an additive copy of the group $\mathfrak M$, with group isomorphism 
$$\mathfrak m\mapsto v\mathfrak m\ \colon\ \mathfrak M\to\Gamma,$$
and equip $\Gamma$ with the ordering $\leq$ such
that for all $\mathfrak m,\mathfrak n\in\mathfrak M$: $\mathfrak m\succeq \mathfrak n \Longleftrightarrow v\mathfrak m\leq v\mathfrak n$.
Then the map
$$v \colon\ C[[\mathfrak M]]\to\Gamma_\infty, \qquad  vf\ =\ \begin{cases} v\big(\mathfrak d(f)\big) & \text{if $f\neq 0$} \\
\infty & \text{if $f=0$}\end{cases}$$
is a valuation on the $C$-vector space $C[[\mathfrak M]]$, making $C[[\mathfrak M]]$ a spherically complete Hahn space over $C$. Moreover, 
$v$ is a valuation on the domain $C[[\mathfrak M]]$. The binary
relation $\preceq$ on $C[[\mathfrak M]]$ associated to $v$ satisfies, for nonzero $f,g\in C[[\mathfrak M]]$:
$$f\preceq g \qquad\Longleftrightarrow\qquad \mathfrak d(f)\preceq\mathfrak d(g) \text{ (in the ordered set $\mathfrak M$).}$$

\begin{lemma}
$C[[\mathfrak M]]$ is a field.
\end{lemma}
\begin{proof}
Let $f\in C[[\mathfrak M]]$, $f\neq 0$; to get $g\in C[[\mathfrak M]]$ with $fg=1$, divide $f$ by $f_{\mathfrak d(f)}\mathfrak d(f)$ to arrange $f=1+\epsilon$ with $\epsilon\prec 1$.
The map $\Phi\colon C[[\mathfrak M]]\to C[[\mathfrak M]]$ given by $\Phi(x)=1-\epsilon x$ is contractive: if $x,y\in C[[\mathfrak M]]$, $x\ne y$, then $\Phi(x)-\Phi(y)=\epsilon(y-x) \prec x-y$. Since the valued additive group $C[[\mathfrak M]]$ is spherically complete, Theorem~\ref{thm:fixpoint} gives $g\in C[[\mathfrak M]]$ with $1-\epsilon g=\Phi(g)=g$, and so $fg=1$. 
\end{proof}

\noindent
The valuation ring of the valued field $K=C[[\mathfrak M]]$ is
$$\mathcal O\ =\ \big\{ f\in K :\  \supp(f) \subseteq \mathfrak M^{\preceq 1} \big\}.$$
The map sending $f\in \mathcal O$ to its constant term $f_1$ is a surjective ring morphism $\mathcal O\to C$ with kernel
$$\smallo\ =\ \big\{ f\in K :\  \supp(f) \subseteq \mathfrak M^{\prec 1} \big\},$$
hence induces an isomorphism $\res(K)=\mathcal O/\smallo\to C$. 
We call $K$ the valued field of {\bf well-based series with coefficients in $C$ and monomials from~$\mathfrak{M}$.} A valued field of the form $C[[\mathfrak M]]$ is also referred to as a {\bf Hahn field over $C$}. 

\nomenclature[K]{$C[[\frak M]]$}{valued field of  well-based series with coefficients in $C$ and monomials from~$\mathfrak{M}$}
\nomenclature[K]{$C\(( t^\Gamma\)) $}{valued field $C[[x^\Gamma]]$ where $x^\gamma=t^{-\gamma}$}
\index{field!Hahn}

\bigskip\noindent
Often we prefer an alternative notation for Hahn fields, as follows. Let
$C$ be a field and~$\Gamma$ an additive ordered abelian group.
Then $C\(( t^\Gamma\)) $ is the field consisting of the formal series $f=\sum_{\gamma} c_{\gamma}t^\gamma$  (summation over all $\gamma\in \Gamma$) with all coefficients~$c_\gamma\in C$, such that $\{\gamma:\ c_\gamma\ne 0\}$
is a well-ordered subset of $\Gamma$, and with the obvious addition, and multiplication according to $t^{\alpha}t^\beta=t^{\alpha+\beta}$ for $\alpha,\beta\in \Gamma$. Note that~$C\(( t^\Gamma\)) $ is $C[[x^\Gamma]]$ in our original notation, with $x^\gamma:=t^{-\gamma}$
for $\gamma\in \Gamma$, and with dominance relation $\preceq$
on the multiplicative group $x^\Gamma=t^\Gamma$ given by:
$x^{\alpha} \preceq x^\beta\Longleftrightarrow \alpha\le\beta$.
This yields the
Hahn field $C\(( t^{\Gamma}\)) $ over $C$ with valuation $v\colon C\(( t^{\Gamma} \)) {}^\times \to \Gamma$ given by $$v(f)\ =\ \min\{\gamma\in \Gamma:\ c_{\gamma}\ne 0\}$$
for nonzero $f= \sum_{\gamma} c_{\gamma}t^\gamma$ as above.
When $\Gamma$ is the ordered group $\Z$ of integers, this gives the valued field $C\(( t^\Z\)) =C\(( t\)) =C[[x^{\Z}]]$ \nomenclature[K]{$C\(( t\)) $}{valued field of Laurent series with coefficients in $C$}\index{series!Laurent} of Laurent series in $t:=t^1=x^{-1}$ over $C$, with
valuation ring $C[[t]]$ in conventional notation.

\subsection*{The overrings of a valuation ring}
Let $A$ be a domain with fraction field $K$. Recall from Section~\ref{sec:local rings} that, given a prime ideal $\mathfrak p$ of $A$,   the localization of  $A$ with respect to $\mathfrak p$
is the subring
$$A_{\mathfrak p}= \left\{a/s:\ a,s\in A,\ s\notin\mathfrak p\right\}$$
of $K$, that $A_{\mathfrak p}$ is a local ring with maximal 
ideal $\mathfrak p A_{\mathfrak p}$ generated by $\mathfrak p$, and
that $\mathfrak p A_{\mathfrak p}\cap A=\mathfrak p$.
Note also that for prime ideals $\mathfrak p$, $\mathfrak p'$ of $A$ we have:
$$A_{\mathfrak p}\subseteq A_{\mathfrak p'}\ \Longleftrightarrow\ \mathfrak p\supseteq\mathfrak p'.$$

\noindent
{\em In the rest of this subsection $A$ is a valuation ring of $K$ with associated valuation $v\colon K^{\times}\to\Gamma$ where $\Gamma=\Gamma_A$}. The following is easy to verify:

\index{localization}
\nomenclature[A]{$A_{\mathfrak p}$}{localization of  $A$ with respect to its prime ideal~$\mathfrak p$}

\begin{lemma}\label{lem:overrings vs prime ideals}
If $B$ is a subring of $K$ containing $A$, then $B$ is a valuation ring of~$K$, $\mathfrak{m}_B\subseteq \mathfrak{m}_A$, and $B=A_{\mathfrak p}$ for a unique prime ideal~$\mathfrak p$ of~$A$, namely $\mathfrak p = \mathfrak{m}_B$. 
\end{lemma}

\noindent
Combining the previous lemma with the next lemma shows that the collection of overrings of $A$ in its fraction field is totally ordered by inclusion. The proof of this lemma is also a routine verification left to the reader.

\begin{lemma}\label{lem:convex subgroups vs prime ideals} 
We have an inclusion-reversing bijection from the set of convex subgroups of $\Gamma$ onto the set of prime ideals of $A$ given by
$$\Delta\ \mapsto\ \mathfrak p_\Delta\ :=\ \{ x\in K:\ vx>\Delta\},$$
with inverse
$$\mathfrak p\ \mapsto\ \Delta_{\mathfrak p}\ :=\ \{ \gamma\in\Gamma:\ \text{$\abs{\gamma}<vx$ for all $x\in \mathfrak p$} \}.$$
\end{lemma}

\noindent
The rank of the valuation ring $A$ is defined to be the rank of the ordered abelian group~$\Gamma$, and similarly, by the rank of a valued field we mean the rank of its value group. (We use analogous terminology for rational rank in place of rank.)
By the previous lemma, if $A$ has finite rank $r$, then $K$ has exactly $r+1$ 
subrings containing~$A$ and they form a tower: $A=B_0\subseteq B_1\subseteq\cdots\subseteq B_{r}=K$.

\begin{cor}\label{cor:archmax} The following conditions on 
$\Gamma$ and $A$ are equivalent: \begin{enumerate}
\item[\textup{(i)}] $\Gamma$ is archimedean;
\item[\textup{(ii)}] $A$ is a proper subring of $K$ and is maximal with respect to inclusion among the proper subrings of $K$.
\end{enumerate}
\end{cor}

\noindent
Let $\mathfrak p$ be a prime ideal of $A$ and $B=A_{\mathfrak p}$, with associated valuation $v_B\colon K^{\times}\to \Gamma_B$, and let $\Delta=\Delta_{\mathfrak p}$ with ordered quotient group $\dot\Gamma=\Gamma/\Delta$. The inclusion $A\subseteq B$ gives rise to an ordered group morphism $\Gamma=K^\times/A^\times\to K^\times/B^\times=\Gamma_B$ with kernel~$\Delta$, so we have an isomorphism $\Gamma_B\to\dot\Gamma$ (of ordered groups) which fits into the 
commutative diagram
$$\xymatrix@R-0.5em@C+1em{
K^\times \ar[r]^v \ar[d]_{v_B}      & \Gamma\ar[d] \\
\Gamma_B \ar[r]^{\cong}       & \dot\Gamma
}$$
where the arrow on the right is the natural surjection $\Gamma\to\dot\Gamma$.

\subsection*{Degree, residue degree, and ramification index}
Let $K\subseteq L$ be a field extension. Then $[L:K]$ is its {\bf degree}, 
that is, the dimension of~$L$ as a vector space over $K$, with the convention that $[L:K]=\infty$ if this dimension is infinite. Likewise, given an extension $\Gamma\subseteq\Gamma'$ of abelian groups we let $[\Gamma':\Gamma]$ be its index, which by convention is $\infty$ if $\Gamma'/\Gamma$ is infinite. We also have the 
transcendence degree $\operatorname{trdeg}(L|K)$ of $L$ over $K$, set equal to
$\infty$ if $\operatorname{trdeg}(L|K)$ is not finite. 
{\em Below in this subsection~$K\subseteq L$ is a valued field extension.}\/ We call $[\res(L):\res(K)]$ the {\bf residue degree} of~$L$ over~$K$ and $[\Gamma_L:\Gamma]$ the {\bf ramification index} of $L$ over $K$. 

\index{index!abelian group extension}
\index{extension!abelian groups!index}
\index{extension!fields!transcendence degree}
\index{extension!fields!degree}
\index{extension!valued fields!residue degree}
\index{extension!valued fields!ramification index}
\index{residue degree}
\index{ramification index}
\index{index!ramification}
\index{transcendence!degree}
\index{degree!field extension}
\index{degree!transcendence}
\index{degree!residue}

\nomenclature[A]{$[L:K]$}{degree of the field extension $L\supseteq K$}
\nomenclature[A]{$\operatorname{trdeg}(L\lvert K)$}{transcendence degree of the field extension $L\supseteq K$}
\nomenclature[A]{$[\Gamma':\Gamma]$}{index of the abelian group extension $\Gamma'\supseteq\Gamma$}

\begin{prop}\label{prop:ef inequ}
Let $b_1,\dots,b_m\in\mathcal O_L$ be such that $\overline{b_1},\dots,
\overline{b_m}\in\res(L)$ are
linearly independent over $\res(K)$, $m\geq 1$.  
Likewise, let $c_1,\dots,c_n\in L^\times$ be such
that $vc_1,\dots,vc_n\in\Gamma_L$ lie in distinct cosets of $\Gamma$, $n\geq 1$.
Then 
\begin{equation}\label{eq:ef inequ}
v\left(\sum_{i,j} a_{ij} b_ic_j \right)\ =\ \min_{i,j} v(a_{ij}c_j) \qquad\text{\textup{(}all $a_{ij}\in K$\textup{)}.}
\end{equation}
In particular, the family $(b_ic_j)$ is $K$-linearly independent.
\end{prop}

\begin{proof}
First we show that $v(a_1b_1+\cdots+a_mb_m)=\min_i v(a_i)$ for $a_1,\dots,a_m\in K$. We can
assume that some $a_i\neq 0$, and then dividing by an $a_i$ of minimum valuation, we can reduce
to the case that $v(a_i)\geq 0$ for all $i$ and $v(a_i)=0$ for some $i$. We must show that
then $v(a_1b_1+\cdots+a_mb_m)=0$. This follows by reduction mod $\smallo_L$ using the linear
independence of $\overline{b_1},\dots,\overline{b_m}$ over $\res(K)$. 
To show \eqref{eq:ef inequ},
we can assume that for each $j$ there is an $i$ with $a_{ij}\neq 0$. So for any $j$,
\begin{multline*}
\gamma_j\ :=\ v\left(\sum_{i} a_{ij}b_ic_j\right)\ =\ v\left( \left(\sum_i a_{ij}b_i\right)c_j\right)\ =\ 
v\left(\sum_i a_{ij}b_i\right) + vc_j\ = \\
\min_i v(a_{ij})+vc_j \in \Gamma+vc_j.
\end{multline*}
Now for $j\neq j'$ we have $\Gamma+vc_{j}\neq\Gamma+vc_{j'}$, so $\gamma_{j}\neq\gamma_{j'}$.
Hence
\equationqed{v\!\left(\sum_{i,j} a_{ij}b_ic_j\right)\ =\ \min_j\gamma_j\ =\ \min_j\!\big(\!\min_i v(a_{ij})+v(c_j)\big)\ =\ 
\min_{i,j} v(a_{ij}c_j).}
\end{proof}

\begin{cor}\label{cor:ef inequ}
$[L:K] \geq  \big[\!\res(L):\res(K)\big]\cdot [\Gamma_L:\Gamma]$.
\end{cor}

\noindent
Under suitable extra assumptions on $K$, the inequality in Corollary~\ref{cor:ef inequ} is an equality; see Corollary~\ref{cor:ef equ}.
For algebraic extensions the previous corollary yields:

\begin{cor}\label{cor:ef inequ, 2}
If $[L:K]=n$, then $\big[\!\res(L):\res(K)\big]\leq n$ \textup{(}so $\res(L)$ is algebraic over $\res(K)$\textup{)}, and $[\Gamma_L:\Gamma]\leq n$ \textup{(}so $m\Gamma_L\subseteq\Gamma$ for some $m\in\{1,\dots,n\}$\textup{)}.
\end{cor}

\begin{remark}
Degree, residue degree and ramification index are multiplicative:
for valued field extensions $K\subseteq L\subseteq M$ with $[M:K]<\infty$,
\begin{align*}
[M:K]\				&=\ [M:L]\cdot[L:K],\\ 
\big[\!\res(M):\res(K)\big]\ 	&=\ \big[\!\res(M):\res(L)\big]\cdot\big[\!\res(L):\res(K)\big],\\
\big[\Gamma_M:\Gamma\big]\ 	&=\ [\Gamma_M:\Gamma_L]\cdot[\Gamma_L:\Gamma].
\end{align*}
\end{remark}

\noindent
Corollary~\ref{cor:ef inequ} also has an analogue for the transcendence degree 
of valued field extensions. To see this, note:

{\sloppy

\begin{lemma}\label{lem:trdeg inequ}
Let $x_1,\dots,x_m\in\mathcal O_L$ be such that the residue classes 
$\overline{x_1},\dots,\overline{x_m}\in\res(L)$ are algebraically independent over $\k=\res(K)$, and let $y_1,\dots,y_n\in L^\times$ be such that the cosets 
$vy_1+\Gamma,\dots,vy_n+\Gamma\in\Gamma_L/\Gamma$
are $\Z$-linearly independent. Then $x_1,\dots,x_m,y_1,\dots,y_n$ are algebraically independent over $K$. 
\end{lemma}

}
\noindent
This lemma is immediate from Proposition~\ref{prop:ef inequ}, and in turn now entails what is sometimes called the
{\em Zariski-Abhyankar Inequality}:

\begin{cor}\label{cor:ZA inequality}
$\operatorname{trdeg}(L|K)\ \geq\ \operatorname{rank}_\Q (\Gamma_L/\Gamma) + \operatorname{trdeg}\!\big(\!\res(L)|\res(K)\big)$. 
%where the inequality is to be understood as an inequality of cardinal numbers 
%if one of these quantities is infinite. 
\end{cor} 

\noindent
The following consequence of this inequality is used in Section~\ref{sec:modth val fields}:

\begin{cor}\label{cor:ZA app} Let $E\subseteq F$ be a field extension with $\operatorname{trdeg}(F|E)=1$ and let~$A$ be a valuation ring of $F$ such that $E\subseteq A\ne F$. Then
$A$ is maximal among the proper subrings of $F$.
\end{cor}
\begin{proof} By Corollary~\ref{cor:ZA inequality} applied to
$K:= E$ with the trivial valuation, and $L:= F$ with the valuation given by $A$, the value group $\Gamma_L$ has 
$\operatorname{rank}_{\Q}(\Gamma_L)=1$, and so is archimedean. It remains to use Corollary~\ref{cor:archmax}. 
\end{proof}

%\noindent
%Let $\sigma\in\Aut(L|K)$. Then
%$\mathcal O_L':=\sigma^{-1}(\mathcal O_L)$ is  a valuation ring of 
%$L$ lying over $\mathcal O$;
%let $v'\colon L^\times\to\Gamma':=L^\times/(\mathcal O_L')^\times$ be 
%the valuation on $L$ associated to $\mathcal O_L'$.
%The automorphism $\sigma$ induces an isomorphism 
%$$\overline{\sigma}\colon\Gamma'=L^\times/(\mathcal O_L')^\times\to 
%L^\times/\mathcal O_L^\times=\Gamma_L$$ of ordered groups, and the diagram
%$$\xymatrix@R-1em@C+2em{
%& \Gamma' \ar[dd]^{\overline{\sigma}} \\
%L^\times \ar[ru]^{v'} \ar[rd]_{v_L\circ\sigma} \\
%& \Gamma_L }$$
%commutes. In particular, if $\mathcal O_L'=\mathcal O_L$, then 
%$\Gamma'=\Gamma_L$, $\overline{\sigma}=\id_{\Gamma_L}$, and 
%$v_L=v_L\circ\sigma$.

%\subsection*{Integral extensions}

%\noindent
%If $\mathfrak p$ and $\mathfrak q$ are prime ideals of $A$ 
%and $B$, respectively, we say that {\bf $\mathfrak q$ lies %over $\mathfrak p$} if $\mathfrak q\cap A=\mathfrak p$.  
%\marginpar{Still need to deal with this.}

\subsection*{Integral closure and valuations}
Next we relate integrality to valuations. 
{\em In this subsection we fix a field $K$ and a local subring $A$ of $K$ with maximal ideal $\mathfrak m=\mathfrak m_A$}.
The first proposition and Zorn imply that there is always a valuation ring of $K$ lying over $A$. (See the beginning of
this section for the meaning of {\em lying over} for local rings and in particular for valuation rings.)

\begin{prop}\label{prop:Krull}
Consider the class of all local subrings of $K$ lying over $A$, partially ordered by $B\leq B' :\Longleftrightarrow \text{$B'$ lies over $B$.}$ Any maximal element of this class is a valuation ring of $K$.
\end{prop}

\noindent
In the proof of this proposition we use the following lemma. 

\begin{lemma}
Suppose $x\in K^\times$ is such that
$1\in\mathfrak m A[x^{-1}]+x^{-1}A[x^{-1}]$. Then~$x$ is integral over $A$.
\end{lemma}

\begin{proof}
We have $1=a_nx^{-n}+\cdots+a_1x^{-1}+a_0$ where $a_1,\dots,a_n\in A$ and $a_0\in\mathfrak m$. Multiplying both sides by $x^n$ yields:
$$x^n(1-a_0) + \text{(terms of lower degree in $x$)}\ =\ 0.$$
Since $1-a_0$ is a unit in $A$, it follows that $x$ is integral over $A$.
\end{proof}

\begin{proof}[Proof of Proposition~\ref{prop:Krull}]
Replacing $A$ by a maximal element of the class of local subrings of $K$ lying over $A$, we arrange that $A$ is the only local subring of~$K$ lying over $A$; we need to show that then $A$ is a valuation ring of~$K$.
For this, let $x\in K^\times$. 
Suppose first that $x$ is integral over $A$. Then $A[x]$ is integral over~$A$, so Corollaries~\ref{cor:primes incomp} and~\ref{cor:CS} give a maximal ideal $\mathfrak n$ of $A[x]$ with $\mathfrak{n} \cap A=\mathfrak m$, and then~$A[x]_{\mathfrak n}$ is a local subring of $K$ which lies over $A$.
By maximality of $A$ this gives~${x\in A}$.
Next, suppose $x$ is not integral over $A$. Then by the lemma above we have $1\notin\mathfrak m A[x^{-1}]+x^{-1}A[x^{-1}]$, so we have a maximal ideal $\mathfrak n$ of $A[x^{-1}]$ such that $\mathfrak n\supseteq\mathfrak m$, 
and thus $\mathfrak n\cap A=\mathfrak m$. The local subring $A[x^{-1}]_{\mathfrak n}$ of $K$ lies over $A$, and maximality of $A$ yields $x^{-1}\in A$.
\end{proof}

\begin{cor}\label{cor:Krull}
Let $L\supseteq K$ be a field extension. 
For each valuation ring of~$K$ there exists a valuation ring of $L$ which lies over it.
\end{cor}

\begin{lemma} \label{lem:val rings int closed}
Each valuation ring is integrally closed. 
\end{lemma}
\begin{proof}
Suppose that $A$ is a valuation ring of $K$, let $v=v_A$ be the associated valuation, and let $x\in K$ satisfy  $x^n+a_{n-1}x^{n-1}+\cdots+a_0=0$ (all $a_i\in A$). 
If $vx<0$, then $v(x^n)=nvx < ivx+v(a_i)=v(a_ix^i)$ for $i=0,\dots,n-1$ and hence 
$v(x^n + \cdots + a_0)=nvx\ne \infty$, a contradiction.
\end{proof}

\noindent
From Lemma~\ref{lem:val rings int closed} we obtain: 

\begin{cor}\label{cor:acvf}
If $K$ is algebraically closed and $A$ is a valuation ring of $K$, 
then the field $\k_A=A/\mathfrak m$ is algebraically closed and 
the abelian group $\Gamma_A$ is divisible.
\end{cor}

\noindent
In view of Corollary~\ref{cor:ef inequ, 2} this gives: 

\begin{cor} Let $L$ be an algebraic closure of $K$. Let
$A$, $B$ be valuation rings of $K, L$, respectively, such that $B$ 
lies over $A$. \textup{(}So $(K,A)$ is a valued subfield of~$(L,B)$.\textup{)}
Then $\k_B$ is an algebraic closure of $\k_A$, and 
$\Gamma_B$ is a divisible hull of $\Gamma_A$. 
\end{cor}

\begin{prop} \label{prop:Krull, 2}
The integral closure of $A$ in $K$ equals the intersection of the valuation 
rings of $K$ lying over $A$.
\end{prop}
\begin{proof}
Any $x\in K$ integral over $A$ lies in every valuation ring of $K$ containing~$A$ as a subring. Next, suppose $x\in K^\times$ is not integral over $A$. Then, by the lemma above, $\mathfrak m A[x^{-1}]+x^{-1}A[x^{-1}]$ is a proper ideal of $A[x^{-1}]$. So we can take a maximal ideal~$\mathfrak n\supseteq\mathfrak m$ of $A[x^{-1}]$ that contains $x^{-1}$. This yields a local subring $A[x^{-1}]_{\mathfrak n}$ of $K$ that lies over $A$. Let $V$ be a maximal element of the class of local subrings of $K$ lying over~$A[x^{-1}]_{\mathfrak n}$. Then $V$ is a valuation ring of $K$, by Proposition~\ref{prop:Krull}, $V$ lies over $A$, and $x^{-1}\in\mathfrak m_V$, so $x\notin V$.
\end{proof}

\subsection*{Valuations and algebraic field extensions}
{\em In this subsection $K$ is a field and~$A$ is a valuation ring of $K$ with 
maximal ideal $\mathfrak m=\mathfrak m_A$. Also, $L$ is an algebraic field 
extension of $K$, and $B$ is the integral closure of~$A$ in~$L$.}

\begin{prop} \label{prop:Krull, 3} The valuation rings of $L$ lying over $A$
are exactly the $B_{\mathfrak q}$ with~$\mathfrak{q}$ a maximal
ideal of $B$.
\end{prop} 
\begin{proof} Let 
$V$ be a valuation ring of $L$ lying over $A$. Valuation rings are integrally closed, so $B\subseteq V$. Set $\mathfrak q:=\mathfrak m_V\cap B$, so
$\mathfrak{q}\cap A=\mathfrak{m}$, hence $\mathfrak q$ is a maximal ideal of $B$ by Corollary~\ref{cor:primes incomp}.

\claim{$V=B_{\mathfrak q}$.} 

\noindent
To prove this claim, note first that
$B\setminus\mathfrak q\subseteq V\setminus\mathfrak m_V$, hence $B_{\mathfrak q}\subseteq V$. For the other inclusion
let $x\in V^{\neq}$. We have a relation
$a_nx^n+\cdots+a_0=0$ where $a_0,\dots,a_n\in A$, $a_n\neq 0$. Take $s\in\{0,\dots,n\}$ maximal such that
$v_A(a_s)=\min_i v_A(a_i)$, and put $b_i:=a_i/a_s$. Dividing by $a_sx^s$ yields
$$\underbrace{(b_nx^{n-s}+\cdots+b_{s+1}x+1)}_{y}+x^{-1}\underbrace{(b_{s-1}+\cdots+b_0/x^{s-1})}_{z} = 0.$$
So $x=-y^{-1}z$. Since $b_i\in\mathfrak m$ for $i=s+1,\dots,n$ we have $y\in V^\times$,  thus~$y\notin\mathfrak q$.
Hence to get $x\in B_{\mathfrak q}$ it suffices to show that $y,z\in B$. This will follow from Proposition~\ref{prop:Krull, 2} if we show that $y$, $z$ lie in every valuation ring of $L$ lying over~$A$. If such a ring contains~$x$, it also contains $y$, hence it contains $z=-yx$. If such a ring does not contain $x$, then it contains $x^{-1}$, and thus $z=b_0/x^{s-1} + \cdots$. This finishes the proof of
our claim. 

Conversely, let $\mathfrak{q}$ be a maximal ideal of $B$. It is clear that
$B_{\mathfrak{q}}$ lies over $A$, and it remains to show that $B_{\mathfrak q}$ is a valuation ring of $L$. Proposition~\ref{prop:Krull} gives a valuation ring $V$ of $L$ that lies over
$B_{\mathfrak q}$, and so $V=B_{\mathfrak{p}}$ where
$\mathfrak{p}$ is a maximal ideal of $B$. As $B_{\mathfrak q}\subseteq V$,
this gives $\mathfrak p\subseteq \mathfrak q$, so 
$\mathfrak p= \mathfrak q$, and thus  $B_{\mathfrak q}=V$.  
\end{proof}

\noindent
Proposition~\ref{prop:Krull, 3} yields a bijection $\mathfrak q\mapsto B_{\mathfrak q}$ from the set of maximal ideals of $B$ onto the set of valuation rings of $L$ lying over $A$. In particular, distinct valuation rings of $L$ lying
over $A$ are incomparable with respect to inclusion.
The next proposition concerns the behavior of valuations under \textit{normal}\/ field extensions:

\begin{prop}\label{prop:transitivity}
Suppose $L\supseteq K$ is normal. Then for any valuation rings $V$, $V'$ of $L$ lying over $A$ there is some $\sigma\in\Aut(L|K)$ such that $\sigma(V)=V'$. 
\end{prop}

\noindent
In the proof we use the Chinese Remainder Theorem~\cite[II,~\S 2]{Lang}: \textit{Let $I_1,\dots,I_n$ \textup{(}$n\ge 1$\textup{)} be ideals in
the commutative ring $R$ such that $I_i+I_j=R$ for all $i$, $j$ with $1\leq i<j\leq n$, and let $a_1,\dots,a_n\in R$. Then there exists an $x\in R$ such that $x\equiv a_i\mod I_i$ for $i=1,\dots,n$.}\/ \index{theorem!Chinese Remainder Theorem} \index{Chinese Remainder Theorem} \label{p:CRT}

\begin{proof}[Proof of Proposition~\ref{prop:transitivity}]
We assume $[L:K]<\infty$ below, since the proposition follows from its validity in that special case.

Let $\mathfrak q$ and $\mathfrak q'$ be maximal ideals of $B$.
%with $\mathfrak{q}\cap A = \mathfrak{q}' \cap A =\mathfrak m$;
By Proposition~\ref{prop:Krull, 3} it suffices to show that 
$\sigma(\mathfrak q)=\mathfrak q'$ for some 
$\sigma\in G:=\Aut(L|K)$.
Assume towards a contradiction that $\sigma(\mathfrak q)\neq\mathfrak q'$ for all $\sigma\in G$. Thus the two sets of maximal ideals of~$B$, $\big\{\sigma(\mathfrak q):\sigma\in G\big\}$ and $\big\{\sigma(\mathfrak q'):\sigma\in G\big\}$,
are disjoint; since $[L:K]<\infty$, these two sets are also finite. By the Chinese Remainder Theorem we get $x\in B$ such that
$$x\equiv 0\bmod \sigma(\mathfrak q), \quad x\equiv 1\bmod \sigma(\mathfrak q')\qquad\text{for all $\sigma\in G$.}$$
Hence $\sigma x\in\mathfrak q\setminus\mathfrak q'$ for all $\sigma\in G$. Recall that
$$\operatorname{N}_{L|K}(x)\ :=\ \left(\prod_{\sigma\in G} \sigma x\right)^\ell$$
lies in $K$, where $\ell=1$ if $\operatorname{char} K=0$, and $\ell=p^e$ for some $e\in\mathbb N$ if  $\operatorname{char} K=p>0$. Each $\sigma x$ is integral over $A$, so $\operatorname{N}_{L|K}(x)\in A$ since $A$ is integrally closed. Also 
$\operatorname{N}_{L|K}(x)\in\mathfrak q\setminus\mathfrak q'$ because $\mathfrak q$ and $\mathfrak q'$ are prime ideals, hence $\operatorname{N}_{L|K}(x)\in A\cap\mathfrak q=\mathfrak m\subseteq\mathfrak q'$, so 
$\operatorname{N}_{L|K}(x)\in\mathfrak q'$, 
 a contradiction.
\end{proof}

\begin{cor}
Suppose $[L:K]<\infty$. There are only finitely many valuation rings of $L$ lying over~$A$.
\end{cor}
\begin{proof}
Take a field extension $L'\supseteq L$ such that $[L':K]<\infty$ and
$L'\supseteq K$ is normal, apply Proposition~\ref{prop:transitivity} to $L'$, 
and use that $\Aut(L'|K)$ is finite.
\end{proof}

\noindent
In the next result $K$ is the valued field with valuation ring $A$. 

\begin{cor}\label{cor:Krull, 2}
Let $K^\alg$ be an algebraic closure of $K$ equipped with a valuation ring of $K^\alg$ lying over $A$. Then any valued field embedding $K\to F$, where $F$ is an algebraically closed valued field, extends to a valued field embedding $K^\alg\to F$.
\end{cor}
\begin{proof}
Let $K\to F$ be a valued field embedding, and
let $j\colon K^\alg\to F$ be a field embedding that extends the field embedding $K\to F$. Let $V$ be the valuation ring of $F$. Then $j^{-1}(V)$ is a valuation ring of $K^\alg$ lying over $A$, and so we have a field automorphism $\sigma$
of $K^\alg$ over $K$ such that $\sigma^{-1}j^{-1}(V)$ is the valuation ring of $K^\alg$. Then~$j\sigma$ is a valued field embedding $K^\alg\to F$ as desired.
\end{proof}

\noindent
Although for our purpose we need to consider only valued fields of
characteristic zero, we note
for the sake of completeness:

\begin{lemma} \label{lem:purely insep}
Suppose $\operatorname{char}(K)=p>0$ and $L\supseteq K$ is purely inseparable. 
Then there is a unique valuation ring of $L$ that lies over $A$.
\end{lemma}
\begin{proof} This unique valuation ring is 
$\big\{x\in L:\ \text{$x^{p^n}\in A$ for some $n$}\big\}$.
\end{proof}  

\begin{lemma} \label{lem:unique val, 1}
Let $V$ be a valuation ring of $L$ lying over $A$, with $v$ as its associated 
valuation. Suppose $\sigma\in\Aut(L|K)$ is such that $\sigma(V)=V$.
Then $v\circ\sigma=v$.
\end{lemma}
\begin{proof} Otherwise we have $x\in K^\times$ with $v(\sigma (x))< v(x)$, 
that is, $\sigma(x)/x\notin V$. 
By induction we get $\sigma^n(x)/\sigma^{n-1}(x)\notin V$ for all $n\ge 1$,
and thus $v(\sigma^n(x)) < v(x)$ for all~$n\ge 1$. Taking $n\ge 1$ 
such that $\sigma^n(x)=x$ we get a contradiction. 
\end{proof} 

\begin{cor} \label{cor:unique val, 1}
Let $L$ be a normal extension of $K$ and $V$ a valuation ring of~$L$ lying over $A$, with associated valuation $v$. Suppose $V$ is the 
only valuation ring of~$L$ lying over $A$. Then
for $x\in L^\times$ with
minimum polynomial 
$$a_0+a_1X+\cdots+a_{n-1}X^{n-1}+X^n\qquad (a_0,\dots,a_{n-1}\in K,\ n\geq 1)$$
over $K$ we have $v(x)=\frac{1}{n}v(a_0)$.
\end{cor}
\begin{proof} Use Lemma~\ref{lem:unique val, 1} and the fact
that $a_0$ or $-a_0$ is a product $x_1\cdots x_n$ of 
conjugates $x_i$ of $x$.
\end{proof}

\noindent
Recall that a Galois extension of a field $E$ is an algebraic field extension
of $E$ that is both normal and separable over $E$.
 Let $L$ be a Galois extension of $K$
and $G:=\Aut(L|K)$. Let $V$ be a valuation ring of $L$ lying over $A$. Then 
the  subgroup
$$G^{\operatorname{d}}\ :=\ \big\{\sigma\in G:\  \sigma(V)=V \big\}$$
of $G$ is called the {\bf decomposition group} of $V$ over $K$, and
the fixed field $L^{\operatorname{d}}$ of~$G^{\operatorname{d}}$ is called the {\bf decomposition field} of $V$ over $K$.

\nomenclature[K]{$G^{\operatorname{d}}$}{decomposition group of a valuation ring}
\nomenclature[K]{$L^{\operatorname{d}}$}{decomposition field of $L\supseteq K$}
\index{decomposition!group}
\index{decomposition!field}
\index{group!decomposition}
\index{field!decomposition}

\begin{cor} \label{cor:decomposition field}
Suppose $L$ is a Galois extension of $K$.  If $V'\neq V$ is another valuation ring of $L$ lying over $A$, then $V\cap L^{\operatorname{d}}\neq V'\cap L^{\operatorname{d}}$. Moreover, $L^{\operatorname{d}}$ is the smallest subfield of $L$ containing $K$ and having this property.
\end{cor}
\begin{proof}
If $V'$ is a valuation ring of $L$ lying over $A$ with $V\cap L^{\operatorname{d}}=V'\cap L^{\operatorname{d}}$, then by Proposition~\ref{prop:transitivity} there is some $\sigma\in\Aut(L|L^{\operatorname{d}})$ such that $\sigma(V)=V'$; but $\Aut(L|L^{\operatorname{d}})=G^{\operatorname{d}}$, hence for such $\sigma$ we have $\sigma(V)=V$ and thus $V=V'$. Suppose~$L'$ is any subfield of $L$ with $K\subseteq L'$ which also has the indicated property. Then $\Aut(L|L')\subseteq G^{\operatorname{d}}$: if $\sigma\in\Aut(L|L')$ then both $\sigma(V)$ and $V$ are valuation rings of $L$ lying over $V\cap L'$, hence $\sigma(V)=V$, and thus $\sigma\in G^{\operatorname{d}}$. Therefore $L^{\operatorname{d}}\subseteq L'$. 
\end{proof}

\subsection*{Adjoining roots}
Let $K$ be a valued field with residue field $\k=\res(K)$. For use later in this section and in Section~\ref{sec:henselian valued fields} we show:

\begin{lemma}\label{lem:pth root}
Let $p$ be a prime number, and $x$ an element in a field extension of~$K$ such that $x^p=a\in K^\times$ where $va\notin p\Gamma$. Then $X^p-a$ is the minimum polynomial of~$x$ over $K$, and $v$ extends uniquely to a valuation $w\colon K(x)^{\times}\to\Delta$ with $\Delta\subseteq\Q\Gamma$ \textup{(}as ordered groups\textup{)}. The residue field of $w$ remains $\k$, and $[\Delta:\Gamma]=p$, with
$$\Delta\ =\ \bigcup_{i=0}^{p-1} \Gamma+iw(x)\qquad \text{\textup{(}disjoint union\textup{)}.}$$
\end{lemma}
\begin{proof}
Let $w\colon K(x)^\times\to\Delta$ with $\Delta\subseteq\Q\Gamma$ be a valuation extending $v$. (By Corollaries~\ref{cor:Krull} and \ref{cor:ef inequ, 2} there is such an extension.) Since $va\notin p\Gamma$, the elements
$$w(x^0)\ =\ 0,\quad w(x^1)\ =\ \frac{va}{p},\ \dots,\ 
w(x^{p-1})\ =\ \frac{(p-1)va}{p}$$
of $\Delta$ are in distinct cosets of $\Gamma$, so $1,x,\dots,x^{p-1}$ are $K$-linearly independent, and thus $X^p-a$ is the minimum polynomial of $x$ over $K$. Also, for an arbitrary nonzero element $b=b_0+b_1x+\cdots+b_{p-1}x^{p-1}$ of $K(x)$ (all $b_i\in K$), 
$$w(b)\ =\ \min\left\{v(b_i)+\frac{iva}{p}:\ i=0,\dots,p-1\right\},$$
showing the uniqueness of $w$. This also proves the claims made by the lemma about~$\Delta$. The residue field of $w$ remains $\k$ by
Corollary~\ref{cor:ef inequ}.
\end{proof}

\subsection*{Simple transcendental extensions}
Let $K$ be a valued field with value group $\Gamma$ and residue field $\res(K)$. So far we mainly considered algebraic extensions, but in the rest of this
section $L=K(x)$ is a field extension with $x$ transcendental over $K$. 
We shall indicate some valuations of $L$ that extend the given valuation
$v$ of $K$.

\begin{lemma}\label{lem:rella}
Let $\Delta$ be an ordered abelian group extension of 
$\Gamma$ and $\delta\in \Delta$. Then~$v$ extends uniquely to a valuation $v\colon L^\times \to \Delta$ of $L$ such that
\begin{equation}\label{eq:rella}
va\ =\ \min_i \big(v(a_i)+i\delta\big) \qquad\text{for $a\ =\ \sum_i a_i x^i\in K[x]^{\neq}$ \textup{(}all $a_i\in K$\textup{)}.}
\end{equation}
Here we do not require $v\colon L^\times \to \Delta$ to be surjective.  \end{lemma}
\begin{proof}
There is clearly at most one such extension as in the lemma. To prove existence, define $v\colon K[x]^{\neq}\to\Delta$ by \eqref{eq:rella}.
It suffices to show that $v$ is a valuation on $K[x]$ (hence extends to a valuation on $L$). It is clear that (V2) is satisfied. To show~(V1), let $a=\sum_i a_ix^i\ ,b=\sum_j b_jx^j\in K[x]^{\neq}$ ($a_i,b_j\in K$). Then
$$ab\  =\ \sum_n c_nx^n\qquad\text{where $\ c _n=\sum_{i+j=n} a_ib_j$,}$$
and for all $n$,
$$v\big(c_n\big)+n\delta\ \geq\ \min_{i+j=n} \big( (v(a_i)+i\delta) + (v(b_j)+j\delta)\big)\ \geq\  va+vb.$$
Take $i_0$ minimal such that $v(a_{i_0})+i_0\delta=\min_i \big(v(a_i)+i\delta\big)$ and take $j_0$ minimal such that $v(b_{j_0})+j_0\delta=\min_j \big(v(b_j)+j\delta\big)$. 
Now set $n_0=i_0+j_0$ and consider
$$c_{n_0} x^{n_0}\  =\ a_{i_0}b_{j_0}x^{n_0} + \underbrace{\text{terms $a_ib_jx^{n_0}$ with $i+j=n_0$, and $i<i_0$ or $j<j_0$.}}_{\text{\normalsize each has valuation $>v(a_{i_0})+v(b_{j_0})+n_0\delta$}}$$
Thus
$v( c_{n_0} ) + n_0\delta = \big(v(a_{i_0})+i_0\delta\big)+\big(v(b_{j_0})+j_0\delta\big)=va+vb$.
\end{proof}

\noindent
Next we consider in more detail two special cases. 

\begin{lemma}\label{lem:lift value group ext}
Let $\Delta$ be an ordered abelian group extension of $\Gamma$, and
$\delta\in \Delta$ such that $n\delta\notin\Gamma$ for all $n\ge 1$. Then $v$ extends uniquely to a $($not necessarily surjective$)$ valuation $v\colon L^\times \to \Delta$ of $L$ with $vx=\delta$. The value group of this valuation is
the internal direct sum  $\Gamma\oplus\Z\delta$ in $\Delta$, and its
residue field is $\k:=\res(K)$. 
\end{lemma}
\begin{proof}
If $v$ is an extension as in the lemma, then 
$$v(1)\ =\ 0,\quad v(x)\ =\ \delta,\quad  v(x^2)\ =\ 2\delta,\  \dots~$$
lie in different cosets of $\Gamma$, so by Proposition~\ref{prop:ef inequ}, $v$ satisfies \eqref{eq:rella}, and thus~$v$ must be the extension
of Lemma~\ref{lem:rella}. Let $v$ be the extension described in
Lemma~\ref{lem:rella}. Then $vx=\delta$, and $v(L^\times)=\Gamma\oplus\Z\delta$.
To see that the residue field is still $\k$, let $a\in L$ with $va=0$,
so $a=b/c$ where $b,c\in K[x]^{\neq}$. 
Then $b=dx^m(1+r)$, $c=ex^n(1+s)$ where $d,e\in K^\times$ and 
$r,s\prec 1$ in $L$. Since $vd+m\delta=vb=vc=ve+n\delta$ we have $vd=ve$ and $m=n$. Thus $a=(d/e)(1+t)$ where  $t\prec 1$ in $L$, so 
$\overline{a}=\overline{d/e}\in\k$.
\end{proof}

\noindent
The valuation on $L$ described in the next lemma is called the 
{\bf gaussian extension} of $v$ to $L$ (with respect to $x$). \index{valuation!gaussian extension} \index{gaussian extension} \index{extension!gaussian}

\begin{lemma}\label{lem:gauss}
There is a unique valuation ring of $L$ that makes 
$L$ a valued field extension of $K$ with $x\preceq 1$ and 
$\overline{x}$ transcendental
over $\k:=\res(K)$. Equipping $L$ with this valuation ring we have $\Gamma_L=\Gamma$ and $\res(L)=\k(\overline{x})$, and
\begin{equation}\label{eq:gauss}
v_L(a)\ =\ \min_i v(a_i)\qquad\text{for $a\ =\ \sum_i a_i x^i\in K[x]$ \textup{(}$a_i\in K$\textup{)}.}
\end{equation}
\end{lemma}
\begin{proof}
If $L$ is equipped with a valuation ring as in the lemma with associated valuation $v_L$, then $1,\overline{x},\overline{x}^2,\dots$ are linearly independent over $\k$, so \eqref{eq:gauss} holds by Proposition~\ref{prop:ef inequ}, which determines $v_L$ uniquely, with $\Gamma_L=\Gamma$. As to existence, apply Lemma~\ref{lem:rella} with $\Delta=\Gamma$ and $\delta=0$ to 
equip $L$ with the valuation ring whose associated valuation $v_L$ satisfies \eqref{eq:gauss}. Then clearly $v_L(x)=0$, and $\overline{x}$ is transcendental over $\k$. To get $\res(L)=\k(\overline{x})$, let $a\in L$ with $v_L(a)=0$; we claim that $\overline{a}\in\k(\overline{x})$.
Now $a=b/c$ where $b,c\in K[x]^{\neq}$ with $v_L(b)=v_L(c)=0$. Then 
$\overline{b}\neq 0$, $\overline{c}\neq 0$, and $\overline{b}, \overline{c}\in \k[\overline{x}]$,
and so $\overline{a}=\overline{b}/\overline{c}\in \k(\overline{x})$ as desired.
\end{proof}

\noindent
Let $X$ be an indeterminate, and equip the field $K(X)$ 
with the gaussian extension of $v$ with respect to $X$.
A polynomial $P\in K[X]$ is said to be {\bf primitive} if $vP=0$.
The proof of the next lemma is obvious. \index{polynomial!primitive} \index{primitive polynomial}

\begin{lemma}
If $P\in K[X]\setminus K$ is primitive and $P=P_1\cdots P_n$ with $P_i\in K[X]$ for $i=1,\dots,n$, then there are $a_1,\dots,a_n\in K^\times$ such that each polynomial $Q_i=a_iP_i$ is primitive and $P=Q_1\cdots Q_n$.
\end{lemma}

\noindent
Note that $\mathcal{O}[X]^{\times}=\mathcal{O}^\times$, since $\mathcal{O}$ is a domain.

\begin{lemma}\label{lem:gauss, 2}
Suppose $P\in \mathcal{O}[X]$, $P\notin \mathcal{O}$. Then 
$$\text{$P$ is irreducible in $\mathcal O[X]$}\ \Longleftrightarrow\ 
\text{$P$ is primitive and irreducible in $K[X]$.}$$
\end{lemma}
\begin{proof}
Take $a\in\mathcal O$ with $va=vP$. Then $P=a(a^{-1}P)$ and $a^{-1}P$ is primitive, so if $P$ is irreducible in $\mathcal O[X]$, then $P$ is primitive, hence irreducible in $K[X]$ by the previous lemma. The direction $\Leftarrow$ is obvious.
\end{proof}

\subsection*{Prescribing value group and residue field extensions}
Let $K$ be a valued field with value group $\Gamma$ and residue field $\k=\res(K)$.
We finish this section with the following existence result:

\begin{prop}\label{prop:prescribed}
Let $\Gamma'$ be an ordered abelian group extension of $\Gamma$ and $\k'$ a field extension of $\k=\res(K)$. Then there is a valued field extension $K'$ of $K$ with value group $\Gamma'$ and with residue field isomorphic to $\k'$ over $\k$.
\end{prop}

\noindent
This follows from Lemmas~\ref{lem:pth root} and \ref{lem:lift value group ext} (for the value group extension) and \ref{lem:gauss} and the following fact (for the residue field extension).

\begin{lemma}\label{extresa}
Let $\xi$ lie in an algebraic closure of $\k$, let $P\in\mathcal O[X]$ be monic such that $\overline{P}\in\k[X]$ is the minimum polynomial of $\xi$ over $\k$, and let $a$ be a zero of $P$ in an algebraic closure of $K$. Then $v$ extends uniquely to a \textup{(}not necessarily surjective\textup{)} valuation 
$v'\colon K(a)^\times \to \Q\Gamma$ on $K(a)$. The residue field $\k'$ of $v'$ is $\k$-isomorphic to~$\k(\xi)$, $v'(K(a)^\times)=\Gamma$, and
$[K(a):K]=[\k':\k]$.
\end{lemma}
\begin{proof}
By Lemma~\ref{lem:gauss, 2}, $P$ is the minimum polynomial of $a$ over $K$, hence 
$$[K(a):K]\ =\ \deg P\ =\ \deg\overline{P}\ =\ [\k(\xi):\k].$$
Let $v'\colon K(a)^\times \to \Q\Gamma$ be a valuation on $K(a)$ extending $v$,
with residue field $\k'\supseteq \k$. Then $v'a\ge 0$ since $a$ is integral over $\mathcal O$.
Now $\overline{P}(\overline{a})=0$ in $\k'$, hence there exists a $\k$-isomorphism $\k(\overline{a})\to\k(\xi)$; in particular, $[K(a):K]=[\k(\overline{a}):\k]$ and thus $\k'=\k(\overline{a})$ and $v'(K(a)^\times)=\Gamma$, by Corollary~\ref{cor:ef inequ}. The uniqueness of $v'$ follows from Proposition~\ref{prop:ef inequ}.
\end{proof}

\subsection*{Notes and comments}
Everything in this section is classical valuation theory.
The notations $\prec$ and $\sim$  were introduced by  du Bois-Reymond~\cite{dBR71,dBR72} in asymptotic analysis.
(See \cite{Ehrlich-1, Ehrlich, Fisher} for discussions of his work.) 

The valued fields $C[[\mathfrak M]]$ occur in Hahn~\cite{Hahn}; 
a variant with $\mathfrak M=x^{\R}$ is in Levi-Civita~\cite{Levi-Civita};
generalizations of Hahn fields were considered by Mal$'$cev~\cite{Malcev}, Neumann~\cite{Neumann}, and Higman~\cite{Higman}; see \cite{Ehrlich-2} for some history. 

As to our notation $C[[\mathfrak M]]$ for Hahn fields, another popular notation is $C\(( \mathfrak{M} \)) $, used for example in \cite{DMM2}. But here we only use the latter for
$\mathfrak{M}=t^\Gamma$ as specified at the end of the subsection
on well-based series.  
% w specified at the end of the su 
%with~``$C\(( \mathfrak{M} \)) $'' 
% the more traditionally based use of $C\(( t^{\Gamma} \)) $ 
%where~$\Gamma$ is an additive ordered abelian group and
%$\gamma \mapsto t^\gamma: \Gamma\to \mathfrak{M}, 
%with the Laurent series field $C\(( t\)) =C\(( t^{\Z}\)) $ 
%as special case;
% in our present system we can use this traditional 
%notation without ambiguity,
%keeping in mind that $C\(( t^{\Gamma} \)) =C[[x^\Gamma]]$, 
%with $t^\gamma=x^{-\gamma}$ for $\gamma\in \Gamma$.  

 %Lemma~\ref{lem:neumann} is from~\cite{Neumann}. 
Zariski~\cite[Chapter~VI, \S{}10]{ZS} contains a special case of Corollary~\ref{cor:ZA inequality}, and Abhyankar~\cite{Abhyankar} a more general version for noetherian local rings.
Propositions~\ref{prop:Krull} and \ref{prop:Krull, 2}--\ref{prop:transitivity} and their corollaries are due to Krull~\cite{Krull}, building on work by Deu\-ring~\cite{Deuring} and
Ostrowski~\cite{Ostrowski} for rank~$1$~valuations.  
In that setting Ostrowski defined the ramification index and residue degree, and proved Corollary~\ref{cor:ef inequ} and 
Lemma~\ref{lem:purely insep}.
Lemma~\ref{lem:rella} stems from Ostrowski~\cite{Ostrowski} and Rella~\cite{Rella}, and 
Proposition~\ref{prop:prescribed} from Mac~Lane~\cite{MacLane}.
For the history of valuation theory before Krull's~\cite{Krull}, see
Roquette~ \cite{Roquette02}.

%% file: mt-3-2.tex
\section{Pseudoconvergence in Valued Fields}\label{sec:pc in valued fields}

\noindent
Let $K$ be a valued field. Viewing the additive group of $K$ as a valued abelian group, the material on pseudoconvergent and pseudocauchy sequences in Section~\ref{sec:valued abelian gps} applies.
Let $(a_\rho)$ be a well-indexed sequence in $K$, and $a\in K$.
Recall that  $(a_\rho)$ is said to pseudoconverge to $a$
(notation: $a_\rho \leadsto a$)
if $v(a-a_\rho)$ is eventually strictly increasing. We also say in that case that $a$ is a 
pseudolimit of $(a_\rho)$.
Note that if  $a_\rho \leadsto a$ and $b\in K$, then $a_\rho+ b \leadsto a+b$, and $a_\rho b \leadsto ab$ if $b\neq 0$. More generally:

\begin{prop}\label{prop:pseudocontinuity}
If $a_\rho\leadsto a$, and $P\in K[X],\ P\notin K$,
% be nonconstant  \textup{(}that is, $P\notin K$\textup{)}. 
then $P(a_\rho)\leadsto P(a)$.
\end{prop}

\noindent
The proof of this proposition is based on Taylor expansion of polynomials: for a polynomial $P\in K[X]$ of degree at most $d$ there are unique polynomials $P_{(i)}\in K[X]$ such that in the ring $K[X,Y]$ of polynomials over $K$ in the distinct indeterminates~$X$,~$Y$, the identity
\begin{equation}\label{eq:Taylor}
P(X+Y)\ =\ \sum_{i=0}^d P_{(i)}(X)\cdot Y^i
\end{equation}
holds. For convenience we also set $P_{(i)}=0$ for $i>d$, so $P_{(0)}=P$ and $P_{(1)}=P'$ (the formal derivative of $P$). If $\operatorname{char}(K)=0$ then $P_{(i)}=\frac{1}{i!}P^{(i)}$ where $P^{(i)}$ is the usual $i$th formal derivative of $P$. Although this will not be used until the next section, we already note here a useful identity:

\begin{lemma} \label{lem:iterated Taylor}
$P_{(i)(j)}={i+j \choose i}P_{(i+j)}$ \textup{(}$i,j\in\N$\textup{)}.
\end{lemma}

\medskip
\noindent
We also use a fact on ordered abelian groups, with the easy proof left to the reader:

\begin{lemma}\label{lem:linear functions}
For each $i$ in a finite nonempty set $I$ let $\beta_i\in\Gamma$ and $n_i\in\N^{\ge 1}$, and let $\lambda_i\colon\Gamma\to\Gamma$ be the linear function given by $\lambda_i(\gamma)=\beta_i+n_i\gamma$. Assume that $n_i\neq n_j$ for all distinct $i,j\in I$. Let $\rho\mapsto\gamma_\rho$ be a strictly increasing function from an infinite linearly ordered set without largest element into $\Gamma$. Then there is an $i_0\in I$ such that if $i\in I$ and $i\neq i_0$, then $\lambda_{i_0}(\gamma_\rho)<\lambda_i(\gamma_\rho)$ eventually.
\end{lemma}

\begin{proof}[Proof of Proposition~\ref{prop:pseudocontinuity}]
Assume $a_\rho\leadsto a$, and $P\in K[X]$ with $P\notin K$ has degree at most $d$. Substituting $a$ for $X$ and $a_\rho-a$ for $Y$ in \eqref{eq:Taylor} yields
\begin{align*} P(a_\rho)-P(a)\ &=\ \sum_{i=1}^d P_{(i)}(a)(a_\rho-a)^i, \quad \text{and for $i=1,\dots,d$,}\\
v\big(P_{(i)}(a)(a_\rho-a)^i\big)\ &=\ \beta_i+i\gamma_\rho\quad\text{where $\beta_i:=v\big(P_{(i)}(a)\big)$ and
$\gamma_\rho:=v(a_\rho-a)$.}
\end{align*}
Since $P=P(a)+\sum_{i=1}^d P_{(i)}(a)(X-a)^i$ and $P$ is not constant, there is an $i\in\{1,\dots,d\}$
with $P_{(i)}(a)\neq 0$. 
%Now $$v\big(P_{(i)}(a)(a_\rho-a)^i\big)&=\beta_i+i\gamma_\rho\quad
%\text{where $\beta_i:=v\big(P_{(i)}(a)\big)$ and $\gamma_\rho:=v(a_\rho-a)$.}$$
Since $(\gamma_\rho)$ is eventually strictly increasing, Lemma~\ref{lem:linear functions} yields $i_0\in \{1,\dots,d\}$ with $P_{(i_0)}(a)\neq 0$  such that for every $i\in\{1,\dots,d\}$ with $i\neq i_0$ we have
$\beta_{i_0}+i_0\gamma_\rho < \beta_i+i\gamma_\rho$ eventually.
Then $v\big(P(a_\rho)-P(a)\big)=\beta_{i_0}+i_0\gamma_\rho$ eventually; in particular, the sequence
$\big(v\big(P(a_\rho)-P(a)\big)\big)$ is eventually strictly increasing, that is, $P(a_\rho)\leadsto P(a)$.
\end{proof}

\noindent
%Recall also that $(a_\rho)$ is called a pseudo-cauchy sequence 
%in $K$ (or pc-sequence in~$K$, for short) if
%for some index $\rho_0$ we have 
%$$\tau > \sigma > \rho > \rho_0\ \Longrightarrow\ 
%a_\tau-a_\sigma \prec a_\sigma-a_\rho.$$
%and in that case we also say that {\bf $(a_\rho)$ is pc.}
The sequence $(a_\rho)$ is pc in $K$ if and only if 
$(a_\rho)$ has a pseudolimit in some valued field extension of 
$K$, and in that case, $(a_\rho)$ has even a pseudolimit 
in some elementary extension of the valued field $K$. (See the remark following Lemma~\ref{pc3}.) 
Together with the previous proposition this immediately yields:

\begin{cor}\label{cor:pc and polynomials}
Suppose $(a_\rho)$ is a pc-sequence in $K$, and  $P\in K[X]$ is nonconstant. Then $\big(P(a_\rho)\big)$ is a pc-sequence in $K$.
\end{cor}

\noindent
If $(a_\rho)$ is a pc-sequence, then  $(v(a_\rho))$ 
is either eventually strictly increasing, or eventually constant. Hence:

\begin{cor}\label{cor:pseudolimit 0}
Let $(a_\rho)$ and $(b_\rho)$ be pc-sequences in $K$.  Then:
$$a_\rho b_\rho \leadsto 0 \qquad\Longleftrightarrow\qquad \text{$a_\rho
\leadsto 0$ or $b_\rho\leadsto 0$.}$$
\end{cor}

\noindent 
A valued field extension $L\supseteq K$ is said to be {\bf immediate} if $\res(L)=\res(K)$ and $\Gamma_L=\Gamma$ (equivalently, the extension $L\supseteq K$ of valued additive groups is immediate as defined in Section~\ref{sec:valued abelian gps}). \index{extension!valued fields!immediate} \index{immediate extension}

\medskip\noindent
In the rest of this subsection we assume that $(a_\rho)$ is a pc-sequence in $K$,
and we shall prove that then $(a_\rho)$ has a pseudolimit in an 
immediate valued field extension of $K$.
To prepare for that, let $P\in K[X]$ be nonconstant. Then by Corollaries~\ref{pc4} and \ref{cor:pc and polynomials} there are two possibilities:  
\begin{align*}\text{either }&\big(v(P(a_\rho))\big)\ \text{ is eventually strictly increasing (equivalently,
$P(a_\rho)\leadsto 0$),}\\
\text{or }&\big(v(P(a_\rho))\big)\ \text{ is eventually constant (equivalently,
$P(a_\rho)\not\leadsto 0$).}
\end{align*} 
We say that $(a_\rho)$ is of {\bf algebraic type over $K$} if the first possibility is realized for some nonconstant 
$P\in K[X]$, and then such a $P$ of least degree is called a {\bf minimal polynomial of $(a_\rho)$ over $K$.}
By Corollary~\ref{cor:pseudolimit 0} such a minimal polynomial of $(a_\rho)$ over $K$ is irreducible.
We say that $(a_\rho)$ is of {\bf transcendental type over $K$} if the second possibility is realized for all
nonconstant $P\in K[X]$. Note that for $a\in K$ we have $a_\rho\leadsto a$ iff $P(a_\rho)\leadsto 0$ for $P=X-a$; in particular, if $(a_\rho)$ is of algebraic type over~$K$ and diverges in $K$, then a minimal polynomial of $(a_\rho)$ over~$K$ has degree at least $2$, and
if $(a_\rho)$ is of transcendental type over $K$, then $(a_\rho)$ diverges in $K$.
A pc-sequence in $K$ of transcendental type over $K$ determines an essentially unique immediate extension:

\index{pc-sequence!algebraic type} \index{pc-sequence!transcendental type} \index{pc-sequence!minimal polynomial}
\index{minimal!polynomial} 

\begin{lemma}\label{lem:Kaplansky, 1}
Suppose $(a_\rho)$ is of transcendental type over $K$. The
valuation $v$ on~$K$ extends uniquely to a valuation $v\colon K(X)^{\times} \to \Gamma$ such that
\begin{equation}\label{eq:Kaplansky, 1}
vP = \text{eventual value of $v\big(P(a_\rho)\big)$}\qquad \text{for each $P\in K[X]$.}
\end{equation}
With this valuation $K(X)$ is an immediate valued field extension of $K$ in which $a_\rho\leadsto X$.
Moreover, if $a_\rho\leadsto a$ in a valued field extension of $K$, then there is a valued field
isomorphism $K(X)\to K(a)$ over $K$ that sends $X$ to $a$.
\end{lemma}
\begin{proof}
It is clear that defining $vP$ for $P\in K[X]$ as in \eqref{eq:Kaplansky, 1}, we have a valuation on~$K[X]$ and thus on $K(X)$. The value group of this valuation is still $\Gamma$, and one checks easily that $a_\rho\leadsto X$. To verify that the residue field of $K(X)$ is the residue field of $K$ we first note that because the value groups are equal, each $R\in K(X)$ with $vR=0$ has the form $R=P/Q$ where $P,Q\in K[X]$ with $vP=vQ=0$. So it is enough to consider a nonconstant $P\in K[X]$ with $vP=0$, and find $b\in K$ with $v(P-b)>0$. We have $0=vP=v\big(P(a_\rho)\big)$ eventually, and $\big(v(P-P(a_\rho))\big)$ is eventually strictly increasing, so $v\big(P-P(a_\rho)\big)>0$, eventually. Thus $b=P(a_\rho)$ for big enough $\rho$ will do the job.

Finally, suppose $a_\rho\leadsto a$ with $a$ in a valued field extension of $K$ (whose valuation we continue to denote by $v$ as usual). For nonconstant $P\in K[X]$ we have $P(a_\rho)\leadsto P(a)$ and thus $v\big(P(a)\big)=v\big(P(a_\rho)\big)$, eventually; in particular, $P(a)\neq 0$ and $v\big(P(a)\big)=vP\in\Gamma$. Thus $a$ is transcendental over $K$ and the field isomorphism $K(X)\to K(a)$ over $K$ that sends $X$ to $a$ is even a valued field isomorphism.
\end{proof}

\noindent
Here is an analogue for pc-sequences of algebraic type:

\begin{lemma}\label{lem:Kaplansky, 2}
Suppose $(a_\rho)$ is divergent of algebraic type over $K$. Let $\mu(X)$ 
be a minimal polynomial
of $(a_\rho)$ over $K$, and $a$ a zero of $\mu$ in an extension field of~$K$.  
Then~$v$ extends uniquely to a valuation $v\colon K(a)^{\times} \to \Gamma$ such that
$$v\big(P(a)\big) = \text{eventual value of $v\big(P(a_\rho)\big)$}\qquad \text{for each $P\in K[X]$ of degree $<\deg\mu$.}$$
With this valuation $K(a)$ is an immediate valued field extension of $K$, and $a_\rho\leadsto a$.
Moreover, if $\mu(b)=0$ and $a_\rho\leadsto b$ in a valued field extension of 
$K$, then there is a valued field
isomorphism $K(a)\to K(b)$ over $K$ that sends $a$ to $b$.
\end{lemma}
\begin{proof}
Much of the proof duplicates the proof for the case of transcendental type. A difference is in how we obtain the multiplicative law for $v\colon K(a)^\times\to\Gamma$ as defined above. Let $s,t\in K(a)^\times$. Then $s=S(a)$, $t=T(a)$ with nonzero $S,T\in K[X]$ of degree less than $\deg\mu$, and $ST=Q\mu+R$ with $Q,R\in K[X]$ and $\deg R<\deg \mu$, so $st=R(a)$, and thus eventually
$$vs = v\big(S(a_\rho)\big),\quad vt = v\big(T(a_\rho)\big), \quad v(st) = v\big(R(a_\rho)\big).$$
Also
$$vs+vt = v\big(S(a_\rho)T(a_\rho)\big) = v\big(Q(a_\rho)\mu(a_\rho)+R(a_\rho)\big),\qquad\text{eventually.}$$
Since $\big(v\big(Q(a_\rho)\mu(a_\rho)\big)\big)$ is either eventually strictly increasing or eventually $\infty$, this forces $vs+vt=v(R(a_\rho))$, eventually, so $vs+vt=v(st)$.
\end{proof}

\noindent
The previous two lemmas together now immediately imply:

\begin{cor} \label{cor:kaplansky}
Every pc-sequence in $K$ has a pseudolimit in an immediate valued field extension of $K$.
\end{cor}

\noindent
A minimal polynomial of a pc-sequence in $K$ of algebraic type over $K$ is in 
general non-unique, even if we require it to be monic. For example, adding 
to the constant term of a minimal polynomial $\mu\in K[X]$ of $(a_\rho)$ 
over $K$ any element of $K$ whose valuation lies in the width of the 
pc-sequence $\big(\mu(a_\rho)\big)$ does not change its status as a minimal polynomial 
of $(a_\rho)$ over $K$.

\subsection*{Maximal and algebraically maximal valued fields}
A valued field is said to be {\bf maximal} if it has no proper immediate 
valued field extension.
By Lemma~\ref{krullgravett} and Zorn, every valued field has an immediate 
extension which is maximal. \index{valued field!maximal} \index{maximal!valued field}

\label{p:maximal}

\begin{cor}\label{cor:maximal valued fields}
The following are equivalent:
\begin{enumerate}
\item[\textup{(i)}] $K$ is maximal;
\item[\textup{(ii)}] the valued additive group of $K$ is maximal;
\item[\textup{(iii)}] $K$ is spherically complete;
\item[\textup{(iv)}] every pc-sequence in $K$ has a pseudolimit in $K$.
\end{enumerate}
\end{cor}
\begin{proof} 
The implication (i)~$\Rightarrow$~(iv) follows from 
Corollary~\ref{cor:kaplansky}, (iii)~$\Leftrightarrow$~(iv) from Lemma~\ref{fp1},  
(iv)~$\Rightarrow$~(ii) from Corollary~\ref{impca}, and (ii)~$\Rightarrow$~(i) is trivial.
\end{proof}

\noindent
For any field $C$ and ordered abelian group $\mathfrak M$, the valued field $C[[\mathfrak M]]$ is maximal, by Corollaries~\ref{hahnmax} and \ref{cor:maximal valued fields}.
Thus by Corollaries~\ref{cor:ef inequ, 2} and~\ref{cor:acvf}:

\begin{cor}\label{cor:McL} 
Let $C$ be a field and $\mathfrak M$ be an ordered abelian group.
The field~$C[[\mathfrak M]]$ is algebraically closed iff $C$ is algebraically closed and $\mathfrak M$ is divisible.
\end{cor}

\noindent
By Lemma~\ref{pc2}, if $a\notin K$ is an element in an immediate valued field extension of~$K$, then there is a divergent pc-sequence $(a_\rho)$ in $K$ 
such that $a_\rho\leadsto a$. The following lemma complements this result:

\begin{lemma}\label{lem:pc2 for algebraic type}
Let $a$ in an immediate valued field extension of $K$ be algebraic over~$K$, and $a\notin K$. Then there is a divergent pc-sequence $(a_\rho)$ in $K$ of algebraic type over $K$ such that $a_\rho\leadsto a$. 
\end{lemma}
\begin{proof}
Take a divergent pc-sequence $(a_\rho)$ in $K$ such that $a_\rho\leadsto a$.
Let $P\in K[X]$ be the minimum polynomial of $a$ over $K$.  By the Taylor identity \eqref{eq:Taylor} we have
$$P(a_\rho)=P(a_\rho)-P(a)=(a_\rho-a)\cdot Q(a_\rho),\qquad Q\in K(a)[X]$$
and thus
$$v\big(P(a_\rho)\big) = v(a_\rho-a) + v\big(Q(a_\rho)\big).$$
Since $\big(v(a_\rho-a)\big)$ is eventually strictly increasing and $\big(v(Q(a_\rho))\big)$ is eventually strictly increasing or eventually constant, $\big(v(P(a_\rho))\big)$ is eventually strictly increasing, so~$(a_\rho)$ is of algebraic type over $K$.
\end{proof}

\noindent
A valued field is {\bf algebraically maximal} if it has no immediate proper algebraic valued field extension. Thus every maximal valued field is algebraically maximal, and every algebraically closed valued field is algebraically maximal. Lemmas~\ref{lem:Kaplansky, 2} and \ref{lem:pc2 for algebraic type} yield:  \index{valued field!algebraically maximal} \index{algebraically!maximal}

\label{p:algebraically maximal}

\begin{cor}\label{cor:algebraically maximal}
$K$ is algebraically maximal if and only if each pc-sequence in~$K$ of algebraic type over $K$ has a pseudolimit in $K$.
\end{cor}

\noindent
By Zorn, $K$ has an immediate valued field extension that is algebraically maximal and algebraic over $K$.  By Corollary~\ref{cor:alg max implies henselian} and
Theorem~\ref{thm:henselization} such an extension is unique up to unique isomorphism over $K$, if $K$ has equicharacteristic zero.
%, then such an extension is unique up-to-unique-isomorphism over $K$; 
%see Corollary~\ref{cor:alg max implies henselian} and
%Theorem~\ref{thm:henselization}.
 
\subsection*{Completion of valued fields}
{\em Let $K$ be a valued field}. We say that a valued field extension~$L\supseteq K$ is {\bf dense} if $K$ is dense in $L$ in the valuation topology on~$L$. Every dense extension of valued fields is immediate. The following complements Theorem~\ref{thm:completion}:\index{valued field!dense extension}\index{dense extension!valued fields}\index{extension!valued fields!dense}

\begin{theorem}\label{thm:completion valued fields}
There is a dense valued field extension $K^{\operatorname{c}}\supseteq K$ such that any dense valued field extension $L\supseteq K$ embeds uniquely over $K$ into 
$K^{\operatorname{c}}$.
%unique embedding $L\to K^{\operatorname{c}}$ which is the identity on $K$.
\end{theorem}

\noindent
The proof of this theorem uses the following routine lemma:

\begin{lemma}\label{lem:multiplication of c-sequences}
Let $(a_\rho)$, $(b_\rho)$ be c-sequences in $K$ with the same index set. Then: 
%and $a,b\in K$. Then
\begin{enumerate}
\item[\textup{(i)}]  $(a_\rho\cdot b_\rho)$ is a c-sequence;
\item[\textup{(ii)}] if $a,b\in K$ and $a_\rho\to a$, $b_\rho\to b$, then $a_\rho\cdot b_\rho\to a\cdot b$;
\item[\textup{(iii)}] if  $a_\rho\neq 0$ for all $\rho$ and $a_\rho\not\to 0$, then $(1/a_\rho)$ is a c-sequence.
\end{enumerate}
\end{lemma}
%\begin{proof}
%To show (1), let $\varepsilon\in\Gamma$ be given. By 
%Lemma~\ref{lem:pc4 analogue} either for each $\gamma\in\Gamma$ 
%we have $v(a_\rho)>\gamma$ eventually, or $(v(a_\rho))$ is 
%eventually constant. Hence there is some $\rho_0$ with 
%$v(b_\rho-b_\sigma)+v(a_\sigma)>\varepsilon$ for all 
%$\rho,\sigma>\rho_0$. Similarly, there is some $\rho_1$ 
%such that $v(a_\rho-a_\sigma)+v(b_\rho)>\varepsilon$ for all 
%$\rho,\sigma>\rho_1$. Hence for $\rho,\sigma>\max(\rho_0,\rho_1)$ we have
%\begin{multline*}
%v(a_\rho b_\rho-a_\sigma b_\sigma)=v\big( (a_\rho-a_\sigma)b_\rho + 
%a_\sigma(b_\rho-b_\sigma) \big) \geq \\
%\min\big( v(a_\rho-a_\sigma)+v(b_\rho), v(b_\rho-b_\sigma)+
%v(a_\sigma)\big) > \varepsilon.
%\end{multline*}
%This shows that $(a_\rho b_\rho)$ is a c-sequence; similarly one 
%shows that if $a_\rho\to a$ and $b_\rho\to b$ then 
%$a_\rho\cdot b_\rho\to a\cdot b$. For (2) suppose  
%$a_\rho\neq 0$ for all $\rho$ and $a_\rho\not\to 0$. 
%Then there is some $\alpha\in\Gamma$ such that 
%$v(a_\rho)=\alpha$ eventually. We have
%$v(1/a_\rho-1/a_\sigma) = v(a_\sigma-a_\rho)-2\alpha$ for all 
%sufficiently large $\rho$, $\sigma$; hence $(1/a_\rho)$ is a c-sequence.
%\end{proof}

\begin{cor} \label{cor:completion valued fields, 1}
Suppose $L=K(a_i: i\in I)$ is a valued field extension of~$K$ such that $\Gamma$ is cofinal
in $\Gamma_L$ and each generator $a_i$ is the limit in $L$ of a c-sequence in~$K$. Then $L\supseteq K$ is dense.
\end{cor}

\begin{proof}[Proof of Theorem~\ref{thm:completion valued fields}]
Let $K^{\operatorname{c}}$ be the completion of the valued additive group of $K$. (See Section~\ref{sec:valued abelian gps}.)
By Lemmas~\ref{lem:continuous extension to completion} and \ref{lem:multiplication of c-sequences},
the multiplication $(x,y)\mapsto x\cdot y$ and inversion $x\mapsto 1/x$ ($x\neq 0$) on $K$ have unique extensions to continuous maps 
$K^{\operatorname{c}}\times K^{\operatorname{c}}\to  K^{\operatorname{c}}$ and 
$(K^{\operatorname{c}})^{\neq}\to  K^{\operatorname{c}}$,
with the product topology on $K^{\operatorname{c}}\times K^{\operatorname{c}}$, and~$K^{\operatorname{c}}$ is a field extension of $K$ with the first map as
multiplication and the second map as inversion. By Lemmas~\ref{lem:pc4 analogue} and~\ref{lem:multiplication of c-sequences}, the valuation of 
$K^{\operatorname{c}}$ is a valuation on the field~$K^{\operatorname{c}}$.
Given  a dense valued field extension $L$ of $K$ there is a unique embedding $L\to K^{\operatorname{c}}$ of valued additive groups over $K$, by Corollary~\ref{cor:completion}, and by continuity, this embedding is a valued field embedding.  
\end{proof}

\noindent 
The properties of the
valued field extension $K^{\operatorname{c}}$ of $K$ postulated in Theorem~\ref{thm:completion valued fields} determine $K^{\operatorname{c}}$ up to a  unique (valued field) isomorphism over $K$. We call $K^{\operatorname{c}}$  the {\bf completion of $K$.} Note that by construction,  $K^{\operatorname{c}}$ is indeed complete: every c-sequence in  $K^{\operatorname{c}}$ converges in  $K^{\operatorname{c}}$.
Hence: \index{valued field!completion} \index{completion!valued field} \nomenclature[K]{$K^{\operatorname{c}}$}{completion of the valued field $K$}

\begin{cor} \label{cor:completion valued fields, 2}
A valued field  is complete iff it has no proper dense valued field extension.
\textup{(}So every maximal valued field is complete.\textup{)}
\end{cor}

{\sloppy 
\noindent
For Hahn fields we make the above more
concrete. Consider a (valued) Hahn field~$C[[\mathfrak{M}]]$
over the field $C$ with $\mathfrak{M}\ne \{1\}$. Since
$C[[\mathfrak{M}]]$ is spherically complete, it is complete. For $S\subseteq C[[\mathfrak{M}]]$, set 
$\supp S:= \bigcup_{f\in S}\supp f$, and let $\operatorname{cl}(S)$ be the closure of $S$ in~$C[[\mathfrak{M}]]$ with respect to its valuation topology. For $f=\sum_{\fm}f_{\fm} \fm$ in~$C[[\mathfrak{M}]]$ and $\fn\in \mathfrak{M}$ we define the {\bf truncation $f_{|\fn}$ of $f$ at $\fn$} by
$$f_{|\fn}\ :=\ \sum_{\fm \succ \fn} f_{\fm}\fm,$$ 
so
$f\in \operatorname{cl}\big(\{f_{|\fn}:\ \fn\in \mathfrak{M}\}\big)$.
A set $S\subseteq C[[\mathfrak{M}]]$ is said to be {\bf truncation closed\/} if for all $f\in S$ and $\fn\in \mathfrak{M}$ we have $f_{|\fn}\in S$. Note that if $S$ is a
truncation closed $C$-linear subspace of $C[[\mathfrak{M}]]$,
then $\supp S\subseteq S$. We do not use this here, but it is worth mentioning that many subsets of $C[[\mathfrak{M}]]$ of a natural origin are truncation closed; see~\cite{vdDtrunc}. From the observations above we get:
}

\index{truncation}
\index{closed!truncation}
\nomenclature[K]{$f_{\lvert\fn}$}{truncation of $f\in C[[\mathfrak{M}]]$ at $\fn\in\fM$}

\begin{lemma} \label{clhahn} If $S$ is a truncation closed subset of $C[[\mathfrak{M}]]$, then
$$\operatorname{cl}(S)\ =\ \big\{f\in C[[\mathfrak{M}]]:\ f_{|\fn}\in S \text{ for every }\fn\in \mathfrak{M}\big\},$$
and so $\operatorname{cl}(S)$ is also truncation closed. 
\end{lemma}

\noindent
Let $K$ be a (valued) subfield of $C[[\mathfrak{M}]]$. Then
$\operatorname{cl}(K)$ is also a subfield of $C[[\mathfrak{M}]]$.  Assume further that
$\frak{M}\subseteq K$. Then the valuation topology of $C[[\mathfrak{M}]]$ induces on~$K$ the valuation topology of the valued field $K$, and likewise with $\operatorname{cl}(K)$ instead of~$K$. Thus the valued field extension
$\operatorname{cl}(K)\supseteq K$ is dense. As a valued abelian group 
$\operatorname{cl}(K)$ is complete by Lemma~\ref{lem:c-sequence width}, so the unique valued field embedding $\operatorname{cl}(K) \to K^{\operatorname{c}}$ over~$K$ is an isomorphism.  
Therefore it is reasonable to call 
$\operatorname{cl}(K)$ the {\bf completion of~$K$ in~$C[[\mathfrak{M}]]$}.   
\index{completion!valued subfield of a Hahn field} \nomenclature[K]{$\operatorname{cl}(K)$}{completion in $C[[\mathfrak{M}]]$ of its valued subfield field $K$}

\begin{cor}\label{complhahn} Suppose $K$ is a truncation closed subfield of $C[[\mathfrak{M}]]$ such that $\mathfrak{M}\subseteq K$. Then the completion of $K$ in $C[[\mathfrak{M}]]$ is 
$$\big\{f\in C[[\mathfrak{M}]]:\ \text{$f_{|\fn}\in K$ for every $\fn\in \mathfrak{M}$}\big\}.$$
\end{cor}

\begin{exampleNumbered}\label{ex:LC}
Let $L := \big\{ f\in C[[\fM]]: \text{$\supp f_{|\fn}$ is finite for all $\fn\in\fM$} \big\}$.
Then $C[\fM]\subseteq L$, $C[\fM]$ is dense in $L$ (for the valuation topology on $C[[\fM]]$), and $L$ is a truncation closed subalgebra of the $C$-algebra $C[[\fM]]$.
Hence $\operatorname{cl}(C[\fM])=\operatorname{cl}(L)=L$ by Lemma~\ref{clhahn}; in particular, $L$ is complete as a valued abelian group.

Suppose now that $\fM$ is archimedean. Then $(1-\varepsilon)^{-1}=1+\varepsilon + \varepsilon^2+ \cdots\in L$ for
$\varepsilon\in L$ with $\varepsilon \prec 1$, and so $L$ is a
subfield
of $C[[\fM]]$. Hence $C(\fM)\subseteq L$, and
thus~$L$ is the completion of $C(\fM)$ in $C[[\fM]]$.
\end{exampleNumbered}

\noindent
Let $(\mathfrak{M}_i)_{i\in I}$ be a family of ordered subgroups of $\mathfrak{M}$ such that for all $i,j\in I$ there is a~$k\in I$ with $\mathfrak{M}_i, \mathfrak{M}_j\subseteq \mathfrak{M}_k$
(automatic if the $\mathfrak{M}_i$ are convex in
$\mathfrak{M}$) and $\mathfrak{M}=\bigcup_i \mathfrak{M}_i$. Then the hypothesis of Corollary~\ref{complhahn} holds for $K:=\bigcup_i C[[\mathfrak{M}_i]]$. If in addition
the~$\mathfrak{M}_i$ are convex in $\mathfrak{M}$, then
$\operatorname{cl}(K)=C[[\mathfrak{M}]]$.

\medskip\noindent
The valued field $K$ is said to be {\bf discrete} if $\Gamma\cong\Z$ (as ordered abelian groups), equivalently, the valuation ring $\mathcal{O}$ is a DVR as defined in Section~\ref{sec:ztn}. In that case, every pc-sequence in $K$ has width $\{\infty\}$. Hence for discrete valued fields, the properties {\em complete}, {\em spherically complete}, and {\em maximal\/} are all equivalent. \index{valued field!discrete}

\medskip
\noindent
From Corollary~\ref{cor:completion functor} we obtain:

\begin{cor}\label{cor:completion functor, fields}
If $L\supseteq K$ is a valued field extension and $\Gamma$ is cofinal in $\Gamma_L$, then the inclusion $K\to L$ extends uniquely to a valued field embedding $K^{\operatorname{c}}\to L^{\operatorname{c}}$.
\end{cor}

\noindent
Next we show the continuity of roots of a polynomial in the valuation topology. First an elementary fact: given $n\ge 1$ and $a_1,\dots, a_n, b_1,\dots,b_n\in K$ (not necessarily distinct), there are permutations
$\sigma$ and $\tau$ of $\{1,\dots,n\}$ such that for $i=1,\dots,n$,
$$v\big(a_{\sigma(i)}-b_{\tau(i)}\big)\ =\ \max_{i\leq j\leq n} v\big(a_{\sigma(i)}-b_{\tau(j)}\big).$$ 
To see this, pick $i_1,j_1\in \{1,\dots,n\}$ such that 
$$v(a_{i_1}-b_{j_1})\ =\ \max\big\{v(a_i-b_j):\ 1\le i,j\le n\big\},$$ and set $\sigma(1)=i_1$, $\tau(1)=j_1$. Now continue inductively with the $a_i$ with $i\ne i_1$ and
the $b_j$ with $j\ne j_1$. In the next lemma $K[X]$ is equipped with the gaussian extension of the valuation of $K$.

\begin{lemma}\label{continuityofroots}
Let
$P=\prod_{i=1}^n\, (X-a_i)$ with $n\ge 1$ and $a_1,\dots,a_n\in K$.
Then for each $\gamma\in\Gamma$ there is a $\beta\in\Gamma$ with the property
that if 
$Q=\prod_{i=1}^n\, (X-b_i)$ with $b_1,\dots,b_n\in K$ and
$v(P-Q)>\beta$, and $v(a_i-b_i)=\max_{i\leq j\leq n} v(a_i-b_j)$ for each~$i$,
then $v(a_i-b_i)>\gamma$ for each~$i$.
\end{lemma}
\begin{proof}
We proceed by induction on $n$, the case $n=1$ being trivial. Suppose $n>1$ and let 
$\gamma\in\Gamma$ be given. 
Set $\tilde{P}:=P/(X-a_1)$, and
by inductive hypothesis choose $\tilde{\beta}\in\Gamma$ such that 
if $\tilde{Q}=\prod_{i=2}^n (X-b_i)$ (all $b_i\in K$)
with $v(\tilde{P}-\tilde{Q})>\tilde\beta$, and  $v(a_i-b_i) \geq v(a_i-b_j)$ for~$2\leq i\leq j\leq n$, then $v(a_i-b_i)>\gamma$ for each $i\geq 2$.

Let now $b_1,\dots,b_n\in K$ with  
$v(a_i-b_i)\geq v(a_i-b_j)$ for~$i\leq j\leq n$, and 
$Q=\prod_{i=1}^n\, (X-b_i)$.
Euclidean Division in $K[X]$ shows that the coefficients of $\tilde{P}$ are polynomials (with integer coefficients) in $a_1$ and the coefficients $P_i$ of $P$, and the same polynomials, with $a_1$ and $P_i$ replaced by $b_1$ and $Q_i$, respectively, yield the coefficients of $\tilde{Q}:=\prod_{i=2}^n (X-b_i)$.
Since polynomial functions are continuous in the valuation topology, we can therefore
choose $\beta_0\ge \gamma$ in $\Gamma$ such that if $v(a_1-b_1)>\beta_0$ and 
$v(P-Q)>\beta_0$, then
$v(\tilde{P}-\tilde{Q})>\tilde\beta$.
We have
$$v\big(Q(a_1)\big)\ =\ \sum_{i=1}^n v(a_1-b_i)\ \leq\ nv(a_1-b_1).$$
Using $P(a_1)=0$ we get
$$v\big(Q(a_1)\big)\ =\ v\big((P-Q)(a_1)\big)\ \geq\ v(P-Q)+\min\big(0,nv(a_1)\big).$$
Hence if $v(P-Q) > \beta_1:= n\beta_0-\min\big(0,nv(a_1)\big)$, then $v(a_1-b_1)>\beta_0$.
Thus $\beta:=\max(\beta_0,\beta_1)$ has the required property.
\end{proof}

\begin{cor}\label{cor:completion of acf} If $K$ is algebraically closed, then so is 
its completion $K^{\operatorname{c}}$.
\end{cor}
\begin{proof} Assume $K$ is algebraically closed. Let $L$ be
an algebraic closure of the completion $K^{\operatorname{c}}$ of $K$,
equipped with a valuation ring of $L$ lying over the valuation ring of~$K^{\operatorname{c}}$. To get $K^{\operatorname{c}}$ algebraically
closed, it is enough to show that $K$ is dense in~$L$, 
since then $K^{\operatorname{c}}=L$. So let $a\in L$, and let
$\gamma\in\Gamma$ be given. Let $P\in K^{\operatorname{c}}[X]$ be the
minimum polynomial of $a$ over $K^{\operatorname{c}}$. Choose
$\beta\in\Gamma$ as in Lemma~\ref{continuityofroots} applied to~$L$ in place of~$K$. As $K$ is
dense in $K^{\operatorname{c}}$ we can take a
monic $Q\in K[X]$ of the same degree as $P$ with $v(P-Q)>\beta$. Since $K$ is
algebraically closed, it now follows from Lemma~\ref{continuityofroots}
that there is $b\in K$ with $Q(b)=0$ and
$v(a-b)>\gamma$.
\end{proof}

\subsection*{Valued vector spaces over valued fields}
To study extensions of complete valued fields of finite degree, we temporarily move to the  setting 
of valued vector spaces over valued fields. {\em Below $K$ is a valued field}. \index{valued vector space!over a valued field}

\begin{definition}
A {\bf valued vector space over $K$} is a vector space $G$ over $K$, together
with a surjective valuation $v\colon G\to S_\infty$ on the additive group of $G$, and an action of the value group $\Gamma$ of $K$ on the value set $S$,
$$(\gamma,s)\mapsto \gamma+s\ \colon\ \Gamma\times S\to S,$$ 
such that for all
$\alpha, \beta\in \Gamma$, $s, s'\in S$, and $a\in K^\times,\ g\in G^{\ne}$,
\begin{list}{*}{\setlength\leftmargin{1em}}
\item[(VS1)] $\alpha\leq \beta\ \Longleftrightarrow\  \alpha+s\leq \beta+s$;
\item[(VS2)] $s\le s'\ \Rightarrow\ \alpha+s \le \alpha + s'$;
\item[(VS3)] $v(ag) = va+vg$.
\end{list}
% $a,b\in K^\times$ and $g,h\in G^{\neq}$:
%\begin{list}{*}{\setlength\leftmargin{1em}}
%\item[(VS1)] $va\leq vb\ \&\ vg\leq vh \Rightarrow va+vg\leq vb+vh$;
%\item[(VS2)] $v(ag)=va+vg$.
%\end{list}
Note that here the same symbol $v$ denotes the valuation of the valued abelian group~$G$ and of the valued field $K$. 
It is convenient to extend the action of $\Gamma$ on $S$ to a map $\Gamma_\infty\times S_\infty\to S_{\infty}$ by $\infty+s=\gamma+\infty=\infty$ for all $\gamma\in\Gamma_\infty$, $s\in S_\infty$. Then (VS2) and (VS3) hold for all $\alpha,\beta\in \Gamma_{\infty}$,
$s, s'\in S_{\infty}$, and $a\in K$, 
$g\in G$, and so does the forward direction $\Rightarrow$ of (VS1).
\end{definition}

\noindent
Below we denote such a valued vector space by $(G,S,v)$, or simply by~$G$. For a (vector) subspace
$H$ of $G$ we have $\Gamma+ v(H^{\ne})\subseteq v(H^{\ne})$, and so $H$ is a valued vector
space by restricting the valuation of $G$ to $H$ and the action of $\Gamma$ to $v(H^{\ne})$. 

\begin{examples}
\mbox{}
\begin{enumerate}
\item If the valuation of $K$ is trivial, so $\Gamma=\{0\}$, then a valued vector space over~$K$ is 
the same as
a valued vector space $(G,S,v)$ over the \textit{field $K$}\/ as defined in Section~\ref{sec:valued vector spaces}, together with the  trivial action of $\Gamma$ on $S$.
\item The valued field $K$ is in a natural way a valued vector space $(K,\Gamma,v)$ over itself, with the action being the addition on $\Gamma$. More generally, if $L\supseteq K$ is a valued field extension,  with the valuation of $L$ denoted by $v$,  then $(L,\Gamma_L,v)$, with the action of $\Gamma$ on $\Gamma_L$ by addition, is a valued vector space over $K$.
\item Given a vector space $G\ne \{0\}$ over $K$ and a basis $\mathcal B$ of $G$, we turn $G$ into a valued vector space $(G,\Gamma,v_{\mathcal B})$ over $K$ by setting $v_{\mathcal{B}}g=\min\{va_b : b\in\mathcal B\}$ for $g=\sum_{b\in\mathcal B} a_b b$, where $a_b\in K$ for all $b\in\mathcal B$, $a_b=0$ for all but finitely many $b\in\mathcal B$, with the action of $\Gamma$ on $\Gamma$ given by addition.
%\item Let $n\ge 1$ and let
%$(G_i,S,v_i)$ be a valued vector space over $K$ for $i=1,\dots,n$ where the ordered set $S$ and the action
%of $\Gamma$ on $S$ associated to $(G_i, S, v_i)$
%does not depend on $i$. 
%Let $(G,S,v)$ be the direct product of the valued abelian groups $(G_i,S_i,v_i)$, so 
%$G=G_1\times\cdots\times G_n$ and
%$vg=\min(v_1g_1,\dots,v_ng_n)$ for $g=(g_1,\dots,g_n)\in G$. Then $(G,S,v)$ with the given action of 
%$\Gamma$ on $S$
%is a valued vector space over $K$. {\bf commented out example checked, but seems not needed} 
\end{enumerate}
\end{examples}

\noindent
Let $(G,S,v)$ be a valued vector space over $K$. The notions ``valuation-independent'' and ``valuation basis'' from Section~\ref{sec:valued vector spaces} extend
to this more general setting:
$\mathcal B\subseteq G$ is called {\bf valuation-independent} if $0\notin\mathcal B$, and for every family $(a_b)_{b\in\mathcal B}$ in~$K$, with $a_b=0$ for all but finitely many $b\in\mathcal B$, we have
$$v\left(\sum_{b\in\mathcal B} a_b b \right)\ =\ \min\big(\{v(a_b b):b\in\mathcal B\}\cup\{\infty\}\big).$$
Every valuation-independent subset of $G$ is $K$-linearly independent. (For example,
the set $\{b_ic_j:\ 1\le i\le m,\ 1\le j\le n\}$ in Proposition~\ref{prop:ef inequ} is valuation-independent in $L$ as valued vector space over $K$.)  \index{valuation!independence} \index{valuation!basis} \index{independent!valuation}

\begin{lemma}\label{lem:val independence}
Let $\mathcal B\subseteq G$ be valuation-independent with span $H$, and $g\in G\setminus H$. Then 
$\mathcal B\cup\{g\}$ is valuation-independent iff $vg\geq v(h+g)$ for all $h\in H$.
\end{lemma}
\begin{proof}
The forward direction being obvious, suppose that $vg\geq v(h+g)$ for all $h\in H$. It suffices to show that for all  $h\in H$  and $a\in K^\times$ we have $v(h+ag) = \min(vh,va+vg)$. For this, after dividing by $a$, we may assume $a=1$. If $vh\neq vg$, then $v(h+g)=\min(vh,vg)$, and if $vh=vg$, then $\min(vh,vg)\leq v(h+g)\leq vg$ by assumption, hence $v(h+g)=vg=\min(vh,vg)$.
\end{proof}

\noindent
A {\bf valuation basis} of $G$ is a valuation-independent vector space basis of $G$.  So in Example~(3)
above, $\mathcal B$ is a valuation basis of $(G,\Gamma,v_{\mathcal B})$.

Let $\mathcal B$ be a valuation basis of $G$ and suppose $\abs{\mathcal B}=n\ge 1$ is finite, say $\mathcal B=\{b_1,\dots,b_n\}$. Then $(G,S,v)$ is isomorphic as a valued abelian group to the direct product of the valued abelian groups $(K,S_i,v_i)$, $i=1,\dots,n$, where ${S_i=\Gamma+vb_i}$ as an ordered subset of $S$, $v_ia=va+vb_i$ for $a\in K$. 
Each valued abelian group~$(K,S_i,v_i)$ in turn is isomorphic to the valued abelian group 
$(K,\Gamma,v)$.
Thus, by Lemma~\ref{lem:direct product of vag}:

\begin{lemma}\label{lem:sc vvs}
If $G$ has a finite valuation basis and $K$ is spherically complete, then~$G$ is spherically complete.
\end{lemma}

\begin{cor}\label{cor:sc vvs}
Suppose $G$ has finite dimension as a $K$-vector space, and $K$ is spherically complete. Then $G$ has a valuation basis, hence is spherically complete.
\end{cor}

\begin{proof}
Let $\mathcal B\subseteq G$ be valuation-independent, with span $H$.  Then $H$ as a valued vector space
over $K$ is spherically complete by the lemma above, so we are done if $G=H$. Suppose $G\neq H$ and take $g\in G\setminus H$. Then
$v(g-H)\subseteq S$ has a largest element $v(g-h_0)$, $h_0\in H$, by 
Lemma~\ref{pcmax}. Replacing $g$ by $g-h_0$ we
arrange that $vg\geq v(g-h)$ for all $h\in H$. Then by Lemma~\ref{lem:val independence}, the set $\mathcal B\cup\{g\}$ is valuation-independent. Continuing this way we build a valuation basis of $G$.
\end{proof}

\noindent
With completeness instead of spherical completeness, we have: 

\begin{prop}\label{prop:vvs}
Suppose $K$ is complete, $G$ is finite-dimensional as a vector space over $K$, and  
$\Gamma+s$ is cofinal in $S$ for all $s\in S$. Then $(G,S,v)$ is complete.
\end{prop}

\begin{proof}
Let $\mathcal B$ be a basis for $G$, where $n=\abs{\mathcal B}\ge 1$, say $\mathcal B=\{b_1,\dots,b_n\}$. Then~$\mathcal B$ is a valuation basis of $(G,\Gamma,v_{\mathcal B})$.
As $K$ is complete, $(G,\Gamma,v_{\mathcal B})$ is complete by Lem\-ma~\ref{productcomplete} and the considerations preceding 
Lem\-ma~\ref{lem:sc vvs}. 

\claim{$v_{\mathcal B}$ is {\em equivalent\/} to $v$, in the sense that
there exist $s,t\in S$ with
$$v_{\mathcal B}x + s\ \leq\  vx\  \leq\  v_{\mathcal B}x + t\qquad\text{for all $x\in G$.}$$}

\vskip-1.5em
\noindent
(It follows that $(G,S,v)$ is complete: if $(x_{\rho})$ is a c-sequence in $(G,S,v)$, then
$(x_{\rho})$ is a c-sequence in $(G,\Gamma,v_{\mathcal B})$, which gives $x\in G$
with $x_{\rho}\to x$ in $(G,\Gamma,v_{\mathcal B})$, and then $x_{\rho}\to x$ in $(G,S,v)$.)
To prove the claim, we set $s:=\min(vb_1,\dots,vb_n)$, and then 
$vx\geq v_{\mathcal B}x+s$ for all $x\in G$. 
We show by induction on $n$ that there exists $t\in S$ with $vx\leq v_{\mathcal B}x+t$ 
for all $x\in G$. For $n=1$ this holds for $t:=s$, so assume $n>1$.
For $i=1,\dots,n$, let $G_i$ be the span of $\mathcal B\setminus\{b_i\}$.
Inductively, we can assume that $G_i$ is complete as a valued subspace of $G$, so by Lemma~\ref{lem:complete approx} we can take $t_i\in S$ such that $v(b_i+g)\leq t_i$ for all $g\in G_i$.
Then for $x=\sum_i a_i b_i\in G$ (all $a_i\in K$), we get
$vx \leq va_i+t_i$ for all $i$, so $vx\leq v_{\mathcal B}x+t$ where $t=\max_i t_i$.
\end{proof}

\noindent
{\em In the next two corollaries $L\supseteq K$ is a valued field extension
with $[L:K]<\infty$}. Note that then $\Gamma$ is cofinal in $\Gamma_L$, by Corollary~\ref{cor:ef inequ, 2}. Thus
by Corollary~\ref{cor:sc vvs} and the previous proposition:

\begin{cor}\label{cor:finite ext of complete}
If $K$ is spherically complete, then $L$ is spherically complete. If~$K$ is complete, then $L$ is complete.
\end{cor}

\noindent
By Corollary~\ref{cor:completion functor, fields} we have a unique valued field embedding
$K^{\operatorname{c}}\to L^{\operatorname{c}}$ extending the inclusion $K\to L$.
We view $K^{\operatorname{c}}$ as a valued subfield of $L^{\operatorname{c}}$ via this embedding, and this gives us the subfield $K^{\operatorname{c}}L$ of  
$L^{\operatorname{c}}$.

\begin{cor}\label{cor:completion finite ext} We have 
$L^{\operatorname{c}}=K^{\operatorname{c}}L$ and 
$\mathcal O_{L^{\operatorname{c}}} = \mathcal O_{K^{\operatorname{c}}}\mathcal O_{L}$. 
%The valuation ring of $L^{\operatorname{c}}$ is the only valuation ring 
%of $L^{\operatorname{c}}$ lying over both $\mathcal O_L$ and 
%$\mathcal O_{K^{\operatorname{c}}}$. 
\end{cor}
\begin{proof} Since $[K^{\operatorname{c}}L:K^{\operatorname{c}}]\leq [L:K]<\infty$,
the valued subfield  $LK^{\operatorname{c}}$ of  
$L^{\operatorname{c}}$ is complete by Corollary~\ref{cor:finite ext of complete}.
As $L$ is dense in $L^{\operatorname{c}}$, so is $K^{\operatorname{c}}L$, and thus
$K^{\operatorname{c}}L=L^{\operatorname{c}}$.   

Next, take a basis of the vector space $L$ over $K$ that is contained in
$\mathcal{O}_L$  and extract from it
a basis $\mathcal{B}=\{b_1,\dots,b_n\}\subseteq \mathcal{O}_L$ of 
$L^{\operatorname{c}}$ as a vector space
over $K^{\operatorname{c}}$, where $[L^{\operatorname{c}}:K^{\operatorname{c}}]=n$.
This gives us the valuation $v_{\mathcal{B}}$ on the vector space~$L^{\operatorname{c}}$ over~$K^{\operatorname{c}}$.
The proof of Proposition~\ref{prop:vvs} gives $\alpha\in \Gamma_L$ such that
$v_{\mathcal{B}}x\ge vx+\alpha$ for all $x\in L^{\operatorname{c}}$. 
Let any $a\in \mathcal O_{L^{\operatorname{c}}}$ be given. By density we get
$b\in\mathcal O_{L}$ with $v(a-b)\ge -\alpha$, so  
$v_{\mathcal{B}}(a-b)\ge 0$, that is, $a=b+\sum_{i=1}^n a_ib_i$ with
all $a_i\in \mathcal{O}_{K^{\operatorname{c}}}$, and so $a\in \mathcal O_{K^{\operatorname{c}}} \mathcal O_{L}$, as promised.
\end{proof}

\subsection*{Notes and comments}
Proposition~\ref{prop:pseudocontinuity} for rank $1$ is in Ostrowski~\cite{Ostrowski}.
The classification of pc-sequences into those of algebraic and transcendental type is due to Kaplansky~\cite{Kaplansky}, who also proved Lemmas~\ref{lem:Kaplansky, 1} and \ref{lem:Kaplansky, 2} as well as (ii)~$\Leftrightarrow$~(iii) in Corollary~\ref{cor:maximal valued fields}.
Corollary~\ref{cor:McL} is in Mac~Lane~\cite{MacLane39}. The completion of~$\R(t^\R)$ in~$\R \(( t^\R \)) $ (see Example~\ref{ex:LC}) was first considered 
by Levi-Civita~\cite{Levi-Civita,Levi-Civita2}, and later used for restricted versions of ``non-standard analysis'' by Laugwitz~\cite{Laugwitz} and Lightstone and Robinson~\cite{LightstoneRobinson,Robinson73a}. 
Co\-rol\-lary~\ref{cor:completion of acf} for rank~$1$ is
from the early days of valuation theory~\cite{Kuerschak, Rychlik}.
An early source for results like Corollary~\ref{cor:sc vvs} is Theorem 11 in Chapter 2 of~\cite{Schilling}. 
%due to Krull~\cite{Krull}.
The proof of Corollary~\ref{cor:completion finite ext} via Proposition~\ref{prop:vvs} stems from~\cite{Roquette58} and is credited there to E.~Artin. The notion of valuation-independence was introduced by Baur~\cite{Baur} (with different terminology).

%% file: mt-3-3.tex
\section{Henselian Valued Fields}\label{sec:henselian valued fields}

\noindent
{\em In this section $K$ is a valued field with residue field $\k=\res(K)$}.
Recall that for a polynomial $P\in\mathcal O[X]$ we let $\overline{P}$ be
the polynomial
with coefficients in $\k$ obtained by replacing each coefficient of $P$ by its residue class.
The definition of henselianity isolates a key algebraic property of maximal valued fields (see Corollaries~\ref{cor:alg max implies henselian} and~\ref{cor:alg max equals henselian}): given $P\in \mathcal O[X]$, every non-singular zero of $\overline{P}$ in~$\k$ can be lifted to a zero of the original polynomial $P$ in $\mathcal O$. More precisely:

\begin{definition}\label{def:henselian}
We call $K$ {\bf henselian} if for every polynomial $P\in \mathcal O[X]$ and $\alpha\in\k$ with
$\overline{P}(\alpha)=0$ and $\overline{P}{'}(\alpha)\neq 0$ there is 
$a\in\mathcal O$
with $P(a)=0$ and $\overline{a}=\alpha$. 
\end{definition}

\index{henselian valued field} \index{valued field!henselian}

\noindent
By the next lemma, the $a$ in this definition is unique.
In Chapter~\ref{sec:dh1} we introduce a version of henselianity for 
valued differential fields lacking such uniqueness.

\begin{lemma}\label{lem:henselian, uniqueness}
Let $P\in\mathcal O[X]$ and $\alpha\in\k$ be such that $\overline{P}(\alpha)=0$ and $\overline{P}{'}(\alpha)\neq 0$. Then there is at most one $a\in\mathcal O$ with $P(a)=0$ and $\overline{a}=\alpha$.
\end{lemma}
\begin{proof}
Suppose $a\in\mathcal O$ satisfies $P(a)=0$ and $\overline{a}=\alpha$. Then $\overline{P}{}'(\overline{a})=\overline{P}{}'(\alpha)\neq 0$,
hence $P'(a)\in\mathcal O^\times$. Taylor expansion in $x\in\smallo$ around $a$ gives
\begin{align*}
P(a+x)\ =\ P(a)+P'(a)x+bx^2\ &=\ P'(a)x+bx^2 \qquad (b\in\mathcal O) \\
							&=\  P'(a)x\left(1+P'(a)^{-1}bx\right).
\end{align*}
Since $P'(a)\left(1+P'(a)^{-1}bx\right)\in \mathcal O^\times$ this gives $P(a+x)=0$ iff $x=0$.
\end{proof}

\noindent
Next we show that algebraically maximal valued fields are henselian:

\begin{lemma}\label{lem:hensel}
Let $P\in\mathcal O[X]$ and $a\in\mathcal O$ be such that $P(a)\prec 1$, $P'(a)\asymp 1$, and~$P$ has no zero in
$a+\smallo$. Then there is a pc-sequence $(a_\rho)$ in $K$ such that $P(a_\rho)\leadsto 0$, $a_\rho\equiv a\bmod\smallo$ for each $\rho$, and $(a_\rho)$ has no pseudolimit in $K$.
\end{lemma}

\begin{proof}
Starting with $a$, one step of the classical Newton process for approximating the zeros of polynomials yields an element $b$ of $\mathcal O$ such that 
$v(b-a)=v\big(P(a)\big)>0$ and $v\big(P(b)\big)\geq 2v\big(P(a)\big)>0$: use Taylor expansion to write
\begin{align*}
P(a+x)\ 	&=\  P(a)+P'(a)x+\text{terms of higher degree in $x$} \\
		&=\  P'(a)\left( P'(a)^{-1}P(a)+x+\text{terms of higher degree in $x$} \right);
\end{align*}
setting $x=-P'(a)^{-1}P(a)$ then yields
$$P(a+x)\ =\  P'(a)\left(\text{multiple of $P(a)^2$}\right),$$
hence $b=a+x$ has the required property. Note that $b\equiv a\bmod\smallo$ and thus $P'(b)\equiv P'(a)\bmod\smallo$, in particular $P'(b)\asymp 1$. Hence the hypothesis on $a$ in the statement of the lemma also applies to $b$. We now iterate this process. More precisely, let $\lambda$ be a nonzero ordinal and $(a_\rho)$ a sequence in $a+\smallo$ indexed by the ordinals $\rho<\lambda$ such that $a_0=a$, and $v(a_\sigma-a_\rho)=v\big(P(a_\rho)\big)$ and $v\big(P(a_\sigma)\big)\geq 2v\big(P(a_\rho)\big)$ whenever $\lambda > \sigma > \rho$. (For $\lambda=2$ we have such a sequence with $a_0=a$, $a_1=b$.)
If $\lambda=\mu+1$ is a successor ordinal, then we construct the next term $a_\lambda\in a+\smallo$ by Newton approximation so that $v(a_\lambda-a_\mu)=v\big(P(a_\mu)\big)$ and $v\big(P(a_\lambda)\big)\geq 2v\big(P(a_\mu)\big)$.

Suppose now that $\lambda$ is a limit ordinal. Then $(a_\rho)$ is clearly a pc-sequence and $P(a_\rho)\leadsto 0$.
If $(a_\rho)$ has no pseudolimit in $K$, then we are done. If $(a_\rho)$ has a pseudolimit in $K$, let $a_\lambda$
be such a pseudolimit. Then $P(a_\rho)\leadsto P(a_\lambda)$. Since $\big(v(P(a_\rho))\big)$ is strictly increasing, this yields
$v\big(P(a_\lambda)\big)\geq v\big(P(a_{\rho+1})\big)\geq 2v\big(P(a_\rho)\big)$ for each index $\rho<\lambda$. It is also clear that
$v(a_\lambda-a_\rho)=v(a_{\rho+1}-a_\rho)=v\big(P(a_\rho)\big)$ for each $\rho<\lambda$, in particular, $a_\lambda\in a+\smallo$.
Thus we have extended our sequence by one more term. This building process must come to an end.
\end{proof}

\begin{cor}\label{cor:alg max implies henselian}
Each algebraically maximal valued field is henselian. \textup{(}Thus each maximal valued field and each algebraically closed valued field is henselian.\textup{)}
\end{cor}

\noindent
This follows from Corollary~\ref{cor:algebraically maximal} and the previous lemma.
For equicharacteristic zero valued fields, see Corollary~\ref{cor:alg max equals henselian} for a converse.

\medskip
\noindent
Let $P\in\mathcal O[X]$ and $a\in\mathcal O$ satisfy the hypothesis of Lemma~\ref{lem:hensel}, let $(a_\rho)$ be the pc-sequence constructed in the proof of that lemma, and consider the strictly increasing sequence $(\gamma_\rho)$, where $\gamma_\rho=v(a_\rho-a_\sigma)$ with $\sigma>\rho$. Then by construction of $(a_\rho)$ we have
$\gamma_\sigma\geq 2\gamma_\rho$ for $\sigma>\rho$; thus $(\gamma_\rho)$ is cofinal in a convex subgroup $\neq\{0\}$ of $\Gamma$ (Lemma~\ref{lem:convex hull of sequ}). So if $\operatorname{rank}(\Gamma)=1$, then $(a_\rho)$ has width $\{\infty\}$. Thus:

\begin{cor}\label{cor:henselian complete rank 1}
Each complete  valued field of rank $1$ is hen\-selian.
\end{cor}

\noindent
Corollary~\ref{cor:henselian complete rank 1} is commonly known as {\em Hensel's Lemma}. \index{Hensel's Lemma} 
%For a generalization see Section~\ref{sec:decomposition} below; 
%there we will need the following complement to Lemma~\ref{lem:hensel}:

\begin{lemma}\label{lem:hensel, complement}
Let $P\in\mathcal O[X]$ and $a\in\mathcal O$ be such that 
$P(a)\prec 1$, $P'(a)\asymp 1$, and let $x$ in a valued field extension $L$ of 
$K$ be a zero of $P$ with $v(x-a)>0$. Then $v\big(P(b)\big)=v(x-b)$ for all 
$b\in L$ with $v(a-b)>0$.
%Let also $(a_\rho)$ be a well-indexed sequence in $a+\smallo$. Then
%$v(P(a_\rho))=v(x-a_\rho)$ for all $\rho$,
%hence $a_\rho\leadsto x$ iff $P(a_\rho)\leadsto 0$.
%If $a_\rho\leadsto x$ and $P$ is irreducible in $K[X]$, 
%$then $P$ is a minimal polynomial of $(a_\rho)$ over $K$.
\end{lemma}
\begin{proof}
Taylor expansion yields
\begin{align*}
P(x+Y)\	&=\  P(x) + P'(x)\cdot Y+\text{terms of higher degree in $Y$} \\
		&=\ P'(x)\cdot Y \cdot (1+Q)\qquad\text{where $Q\in\mathcal O_L[Y]$, $Q(0)=0$.}
\end{align*}
Let $b\in L,\ v(a-b)>0$. Substituting $b-x$ for $Y$ yields 
$$P(b)\ =\ P'(x)\cdot (b-x)\cdot \big(1+Q(b-x)\big)\qquad
\text{with $Q(b-x)\prec 1$,}$$
which gives the desired conclusion. 
%This shows the first statement.
%Suppose now $a_\rho\leadsto x$ and $P$ is irreducible in $K[X]$. Then for 
%each nonconstant $Q\in K[X]$ with $\deg Q<\deg P$ we have $Q(x)\neq 0$, 
%and  $Q(a_\rho)\leadsto Q(x)$ by Proposition~\ref{prop:pseudocontinuity}, 
%hence $Q$ cannot be a minimal polynomial of $(a_\rho)$ over $K$. Thus $P$ 
%is a minimal polynomial of $(a_\rho)$ over $K$.
\end{proof}

\begin{cor}\label{cor:Kc henselian}
The following conditions on $K$ are equivalent:
\begin{enumerate}
\item[\textup{(i)}] the completion $K^{\operatorname{c}}$ of $K$ is henselian;
\item[\textup{(ii)}] for every polynomial $P\in\mathcal O[X]$ and $a\in\mathcal O$ with
$P(a)\prec 1$ and $P'(a)\asymp 1$
and every $\gamma\in\Gamma^{>}$ there exists $b\in \mathcal O$ with $v\big(P(b)\big)>\gamma$
and $\overline{a}=\overline{b}$.
\end{enumerate}
In particular, if $K$ is henselian, then so is $K^{\operatorname{c}}$.
\end{cor}
\begin{proof}
Suppose $K^{\operatorname{c}}$ is henselian, $P\in\mathcal O[X]$, $a\in\mathcal O$,
$P(a)\prec 1$, $P'(a)\asymp 1$, and $\gamma\in \Gamma^{>}$.
Take $x\in K^{\operatorname{c}}$ with $x\preceq 1$, $P(x)=0$, and $\overline{x}=\overline{a}$. Next, 
take $b\in K$ such that $v(b-x)>\gamma$. Then 
by Lemma~\ref{lem:hensel, complement} applied to $L=K^{\operatorname{c}}$ we obtain
$v\big(P(b)\big)=v(x-b)>\gamma$.
Thus (i)~$\Rightarrow$~(ii).
Assume (ii),
and let $Q\in \mathcal O_{K^{\operatorname{c}}}[X]$  and $\alpha\in \k$ be such that $\overline{Q}(\alpha)=0$ and $\overline{Q}{'}(\alpha)\neq 0$. We extend the valuation of $K^{\operatorname{c}}$ to its algebraic closure and use Corollary~\ref{cor:alg max implies henselian} to get a zero $x$ of $Q$ in this algebraic
closure such that $\overline{x}=\alpha$. Let 
$\gamma\in\Gamma^{>}$.
Since~$K$ is dense in~$K^{\operatorname{c}}$, we can 
take $P\in K[X]$ with $v(P-Q)>\gamma$. Then $P\in\mathcal O[X]$ and $\overline{P}=\overline{Q}$.
By~(ii) we have $b\in\mathcal O$ with 
$v\big(P(b)\big)>\gamma$
and $\overline{b}=\alpha$. Hence by Lemma~\ref{lem:hensel, complement} again
we have $v(x-b)>\gamma$.
By Corollary~\ref{cor:completion valued fields, 1} and \ref{cor:completion valued fields, 2} this yields $x\in K^{\operatorname{c}}$, so (i) holds.
\end{proof}

\subsection*{Lifting the residue field}
Suppose $C$ is a subfield of the valuation ring $\mathcal O$ of~$K$. Then $C$ is mapped onto a sub\-field~$\overline{C}$ of $\k=\res(K)$ under the residue map 
$$x\mapsto\overline{x}\ \colon\ \mathcal O\to 
\mathcal O/\smallo=\k.$$ In case $\overline{C}=\k$ we call $C$ a {\bf lift} of $\k$ (in $\mathcal O$). For example, the subfield $C$ of a Hahn field~$C[[\mathfrak M]]$ is a lift of the residue field of this Hahn field. 
We say that the residue field of~$K$ can be lifted if there is a lift of $\k$ in $\mathcal O$. 

\index{lifting!residue field}
\index{lift}
\index{residue field!lift}
\index{residue field!lifting}

\begin{prop}\label{prop:lift}
Suppose $K$ is henselian of equicharacteristic zero. Then the residue field of $K$ can be lifted; in fact, every maximal subfield of $\mathcal O$ is a lift of $\k$.
\end{prop}

\noindent
In the proof we use:

\begin{lemma}\label{lem:max subfield}
Let $C$ be a maximal subfield of $\mathcal O$. Then $C$ is algebraically closed in~$K$, and the field extension $\k\supseteq \overline{C}$ is algebraic.
\end{lemma}
\begin{proof}
Since $\mathcal O$ is integrally closed in $K$, it contains the algebraic closure of $C$ in~$K$; hence $C$ is algebraically closed in $K$. Let $\xi\in\k$ be transcendental over $\overline{C}$, and take $x\in\mathcal O$ such that $\overline{x}=\xi$. For $P\in C[X]$, $P\neq 0$, we have $\overline{P}\in\overline{C}[X]$, $\overline{P}\neq 0$, hence $\overline{P(x)}=\overline{P}(\xi)\neq 0$, so $P(x)\in\mathcal O^\times$. In particular, the subring $C[x]$ of $\mathcal O$ is mapped isomorphically onto the subring $\overline{C}[\xi]$ of $\k$ by the residue map. Thus $C[x]$ is a domain with fraction field $C(x)$ inside $\mathcal O$, and $C(x)$ is mapped isomorphically onto the subfield~$\overline{C}(\xi)$ of $\k$. Thus $C(x)$ is a subfield of $\mathcal O$ that strictly contains $C$, a contradiction.
\end{proof}

\begin{proof}[Proof of Proposition~\ref{prop:lift}]
Since $\operatorname{char}\k=0$, the valuation ring $\mathcal O$ of $K$ contains a subfield. By Zorn there is a maximal subfield $C$ of $\mathcal O$.  Suppose $x\in\mathcal O$ and $\overline{x}\in\k\setminus\overline{C}$.  By the previous lemma, $\xi:=\overline{x}$ is algebraic over $\overline{C}$. Let $P\in C[X]$ be a monic polynomial such that its image $\overline{P}\in\overline{C}[X]$ is the minimum polynomial of $\xi$ over~$\overline{C}$. Then $P$ is irreducible since $\overline{P}$ is. 
%For each $\varepsilon\in\smallo$ the residue map 
%$\mathcal O\to\k$ restricts to a surjective ring morphism
%$C[x+\varepsilon] \to \overline{C}[\xi]$,
%where $C[x+\varepsilon]$ is a subring of $\mathcal O$ 
%and $\overline{C}[\xi]$ is a subfield of $\k$. 
Since $\operatorname{char}(\k)=0$, the irreducible polynomial $\overline{P}\in\overline{C}[X]$ is separable, so $\overline{P}{}'(\xi)\neq 0$. As $K$ is henselian, this gives $\varepsilon\in\smallo$ such that $P(x+\varepsilon)=0$. 
%The map $C[x+\varepsilon] \to \overline{C}[\xi]$ is injective and 
Then $C[x+\varepsilon]$ is a subfield of $\mathcal O$, contradicting the maximality of~$C$. Thus $\k=\overline{C}$, so $C$ is a lift of $\k$.
\end{proof}

\subsection*{Characterizations of henselianity}
The following elementary lemma contains useful reformulations of the henselianity condition:
 
\begin{lemma} \label{lem:char henselian}
The following are equivalent:
\begin{enumerate}
\item[\textup{(i)}] $K$ is henselian;
\item[\textup{(ii)}] each polynomial 
$$1+X+P_2X^2+\cdots+P_nX^n\qquad\text{where $n\geq 2$, $P_2,\dots,P_n\in\smallo$}$$  
has a zero in $\mathcal O$ $($of course, such a zero must lie in 
$-1+\smallo)$;
\item[\textup{(iii)}] each polynomial 
$$Y^n+Y^{n-1}+Q_{n-2}Y^{n-2}+\cdots+Q_0\quad \text{where $n\geq 2$, $Q_0,\dots,Q_{n-2}\in\smallo$}$$
has a zero in $\mathcal O^\times$;
\item[(iv)] given a polynomial $P\in\mathcal O[X]$ and $a\in\mathcal O$ such that $P(a)\prec P'(a)^2$, there is $b\in\mathcal O$ such that $P(b)=0$ and $b-a\prec P'(a)$. \textup{(}Newton version\textup{)}
\end{enumerate}
\end{lemma}
\begin{proof}
Assume $K$ is henselian and let $P=1+X+P_2X^2+\cdots+P_nX^n$ with $n\geq 2$ and $P_2,\dots,P_n\in\smallo$. Then for $a=-1$ we have $P(a)\prec 1$ and 
$P'(a)\asymp 1$. Thus~$P$ has a zero in $\mathcal O$. This shows (i)~$\Rightarrow$~(ii). For (ii)~$\Leftrightarrow~$(iii), use the substitution $X=1/Y$. Suppose now that (ii) holds,  let $P$, $a$ be as in the hypothesis of (iv). Note that 
$P'(a)\neq 0$ and set $c:=P(a)/P'(a)^2$ (so $c\prec 1$). 
Let $x\in\mathcal O$ and consider the expansion:
\begin{align*}
P(a+x)\ 	&=\ P(a)+P'(a)x+\sum_{i\geq 2} P_{(i)}(a)x^i \\
			&=\ cP'(a)^2 + P'(a)x+\sum_{i\geq 2} P_{(i)}(a)x^i.
\end{align*} 
Set $x=cP'(a)y$ where $y\in\mathcal O$. Then
$$P(a+x)\ =\ cP'(a)^2 \left(1+y+\sum_{i\geq 2} a_i y^i\right)$$
where the $a_i\in\smallo$ do not depend on $y$. By (ii) choose $y\in\mathcal O$ such that
$$1+y+\sum_{i\geq 2} a_iy^i\ =\ 0.$$
This yields an element $b=a+x=a+cP'(a)y$ as required. This shows (ii)~$\Rightarrow$~(iv), and (iv)~$\Rightarrow$~(i) is clear.
\end{proof}

\noindent
Given an algebraic closure $K^\alg$ of $K$, the valuation $v\colon K^\times \to \Gamma $ extends to a valuation $(K^\alg)^\times \to \Q\Gamma$. The following proposition shows among other things that uniqueness
of such an extension is equivalent to $K$ being henselian:

\begin{prop}\label{prop:char henselian}
The following are equivalent:
\begin{enumerate}
\item[\textup{(i)}] $K$ is henselian;
\item[\textup{(ii)}] for each algebraic field extension $L\supseteq K$ there is a unique valuation ring of $L$ lying over $\mathcal{O}$;
\item[\textup{(iii)}] for each monic polynomial $P\in\mathcal O[X]$ which is irreducible in $K[X]$ there is some~$m\geq 1$ and a monic $Q\in\mathcal O[X]$ with $\overline{Q}$ is irreducible in $\k[X]$ and~$\overline{P}=\overline{Q}^m$;
\item[\textup{(iv)}] given monic $P,Q,R\in\mathcal O[X]$ with $\overline{P}=\overline{Q}\cdot\overline{R}$ and $\overline{Q}$, $\overline{R}$  relatively prime in~$\k[X]$, there are monic $Q^*,R^*\in\mathcal O[X]$ with $P=Q^*R^*$ and $\overline{Q^*}=\overline{Q}$, $\overline{R^*}=\overline{R}$.
\end{enumerate}
\end{prop}

\begin{proof}
Suppose (ii) fails. Then we have a field extension $L\supseteq K$ with $[L:K]<\infty$ and more than one valuation ring of $L$ lying over $\mathcal O$. By Lemma~\ref{lem:purely insep} we can replace $L$ by the separable closure of $K$ in $L$ and arrange that $L$ is separable over~$K$. Replacing $L$ by its normal closure over $K$ (in an algebraic closure of $K$), we can also assume that $L$
is a Galois extension of $K$. Let $\mathcal O_L$ be a valuation ring of~$L$ lying over~$\mathcal O$, and let $G^{\operatorname{d}}$ be the decomposition group and $L^{\operatorname{d}}$ the decomposition field of $\mathcal O_L$ over~$K$; see Corollary~\ref{cor:decomposition field}.
Let $\mathcal O_1,\dots,\mathcal O_m$ be the distinct valuation rings of~$L$ lying over $\mathcal O$, with $\mathcal O_1=\mathcal O_L$, and let $\smallo_i$ be the maximal ideal of $\mathcal O_i$. Then $m>1$ and $\mathcal O_1\cap L^{\operatorname{d}}\neq\mathcal O_i\cap L^{\operatorname{d}}$ for $i=2,\dots, m$.
Let $B$ be the integral closure of $\mathcal O$ in~$L^{\operatorname{d}}$. Then $\smallo_1\cap B,\dots, \smallo_m \cap B$ are maximal ideals of $B$, by Proposition~\ref{prop:Krull, 3}, and $\smallo_1\cap B\ne \smallo_i\cap B$ for $i=2,\dots,m$. Next, the Chinese Remainder Theorem provides an element $x\in B$ such that $x\in\smallo_1$ and $x\notin\smallo_i$ for $i=2,\dots,m$. Let 
$$P\ :=\ X^n+a_{n-1}X^{n-1}+\cdots+a_0\qquad (a_i\in K)$$ 
be the minimum polynomial of $x$ over $K$. From $m>1$ we get $x\notin K$, so $n\geq 2$ and~$P$ does not have a zero in $K$. By Lemma~\ref{lem:int closure} we have $P\in\mathcal O[X]$; we claim that $a_0\in\smallo$ and $a_{1}\notin\smallo$.
To see this, let $x=x_1,x_2,\dots,x_n$ be the distinct conjugates of~$x$ under the action of $\Aut(L|K)$. Now $a_0=(-1)^n x(x_2\cdots x_n)$ and $x_2,\dots, x_n$ are integral over $\mathcal{O}$, so $x_2\cdots x_n\in B$, and thus $a_0\in xB\cap \mathcal{O} \subseteq \smallo$.
Let $j\in\{2,\dots,n\}$ and take $\sigma\in\Aut(L|K)$ with $x_j=\sigma(x)$; since $x_j\neq x$ and $x\in L^{\operatorname{d}}$ we have $\sigma\notin G^{\operatorname{d}}$, hence $\sigma^{-1}(\mathcal O_1)=\mathcal O_i$ where $i\in\{2,\dots,m\}$. For such $i$ we have $x\notin\smallo_i=\sigma^{-1}(\smallo_1)$, so $x_j=\sigma(x)\notin\smallo_1$. 
Therefore the sum $\sum_{j=1}^n x_1\cdots x_{j-1}x_{j+1}\cdots x_n$ in
$$a_1\ =\ (-1)^{n-1}\sum_{j=1}^n x_1\cdots x_{j-1}x_{j+1}\cdots x_n,$$
has precisely one term, namely $x_2\cdots x_n$, that misses the factor $x_1=x$, and so this is the only term in the sum not in $\smallo_1$; thus $a_1\notin \smallo$. It follows that $K$ is not henselian. 
This proves (i)~$\Rightarrow$~(ii), since we have established its contrapositive. 

Suppose now that (ii) holds. Let $K^\alg$ be an algebraic closure of $K$. Then by
Proposition~\ref{prop:Krull, 3}, the integral closure $B$ of $\mathcal O$ in  $K^\alg$ is local. 
Let $P\in\mathcal O[X]$ be monic and irreducible in $K[X]$, and let $x\in K^\alg$ be a zero of $P$. Then the subring~$\mathcal O[x]$ of $B$ is local, since distinct maximal ideals of $\mathcal{O}[x]$
would extend to distinct maximal ideals of~$B$. Now $P$ is monic and irreducible, so
$\mathcal O[x]\cong\mathcal O[X]/(P)$ as ring extensions of~$\mathcal{O}$. By Lemma~\ref{lem:Dedekind}, $\overline{P}=\overline{Q}^m$ for
some $m\ge 1$ and monic irreducible polynomial $\overline{Q}\in\k[X]$.
This shows (ii)~$\Rightarrow$~(iii).

Suppose (iii) holds; let $P,Q,R\in\mathcal O[X]$ be monic such that $\overline{P}=\overline{Q}\cdot\overline{R}$ and $\overline{Q}$, $\overline{R}$ are relatively prime in $\k[X]$. We have $P=P_1\cdots P_m$ where each $P_i\in K[X]$ is monic and irreducible. The coefficients of $P_i$ are elementary symmetric functions in the zeros of $P_i$, and these zeros are among the zeros of $P$ and hence integral over~$\mathcal O$; thus $P_i\in\mathcal O[X]$ for each $i$.
By (iii), $\overline{P_i}$ is a power of an irreducible polynomial in~$\k[X]$. Since  $\overline{Q}$, $\overline{R}$ are relatively prime, either $\overline{P_i}$ divides $\overline{Q}$ or $\overline{P_i}$ divides $\overline{R}$. Let~$Q^*$ be the product of those $P_i$ such that $\overline{P_i}$ divides $\overline{Q}$, and let $R^*$ be the product of the remaining~$P_i$. Then $Q^*$, $R^*$ have the required properties, showing (iv). 

Assume (iv); to derive (i), it suffices, by the equivalence of (i) and (iii) in Lem\-ma~\ref{lem:char henselian}, that each polynomial $P=X^n+X^{n-1}+a_{n-2}X^{n-2}+\cdots+a_0$ with $n\geq 2$ and $a_0,\dots,a_{n-2}\in\smallo$ has a zero in $\mathcal O^\times$. Now $\overline{P}$ factors as
$(X+1)X^{n-1}$. By~(iv) take monic $Q^*,R^*\in\mathcal O[X]$ with $P=Q^*R^*$, $\overline{Q^*}=\overline{Q}=X+1$, $\overline{R^*}=X^{n-1}$.
Then $Q^*=X-a$ with $a\in\mathcal O$, so
$\overline{a}=-1$, hence $a\in\mathcal O^\times$ and $P(a)=0$.
\end{proof}

\noindent
The equivalence of (i) and (ii) in this proposition has an important consequence:

\begin{cor}\label{cor:alg ext of henselian}
If $K$ is henselian and $L$ is an algebraic valued field extension of~$K$, then $L$ is henselian and $\mathcal{O}_L$ is the integral closure of 
$\mathcal{O}$ in $L$.
\end{cor}

\noindent
To get a variant of property (iii) in Proposition~\ref{prop:char henselian} 
for not necessarily monic polynomials we use the following lemma:

\begin{lemma}
Suppose $K$ is henselian, and let $P\in\mathcal O[X]$ be irreducible in $K[X]$ such that 
$\deg \overline{P}\ge 1$. Then  $\deg(\overline{P})=\deg(P)$.
\end{lemma}
\begin{proof}
Let $K^\alg$ be an algebraic closure of $K$ and $\mathcal O^\alg$ the unique valuation ring of~$K^\alg$ lying over $\mathcal O$; we continue to denote the associated valuation $(K^\alg)^\times \to \Q\Gamma$ on~$K^\alg$  by~$v$. Then $v\circ\sigma=v$ for all $\sigma\in\Aut(K^\alg|K)$, by Lemma~\ref{lem:unique val, 1}.  In $K^\alg[X]$,
\begin{equation}\label{eq:char henselian}
P\ =\ a\prod_{j=1}^n (X-x_j)\qquad\text{where $a\in \mathcal{O},\ a\ne 0$, $x_1,\dots,x_n\in K^\alg$.}
\end{equation}
Since for all $i$, $j$ there is $\sigma\in\Aut(K^\alg|K)$ with $\sigma(x_i)=x_j$, we have $\gamma\in\Q\Gamma$ with $v(x_i)=\gamma$ for all $i$. We claim that $\gamma\ge 0$. Suppose $\gamma<0$, and let
$$\prod_{j=1}^n(X-x_j)\ =\ X^n + b_1X^{n-1} + \cdots + b_{n-1}X + b_n, \qquad b_1,\dots, b_{n}\in K.$$ Then $v(b_n)=n\gamma$, and $v(b_i)\ge i\gamma>n\gamma= v(b_n)$ for $i=1,\dots, n-1$, and this holds also
for $b_0:= 1$. But this contradicts $\deg \overline{P}\ge 1$, and thus $\gamma\ge 0$, as claimed. 
Applying the residue morphism to both sides in \eqref{eq:char henselian} and using $\overline{P}\ne 0$ now yields $\overline{a}\neq 0$, that is, $\deg(\overline{P})=n=\deg(P)$.
\end{proof}

\begin{cor}
Suppose $K$ is henselian, $P\in\mathcal O[X]$ is irreducible in $K[X]$, and 
$\deg\overline{P}\ge 1$. Then 
there exist $m\geq 1$, $a\in\mathcal O^\times$, and monic $Q\in\mathcal O[X]$ 
such that $\overline{Q}$ is irreducible in $\k[X]$ and $\overline{P}=\overline{a}\cdot\overline{Q}^m$.
%\item if $P,Q,R\in\mathcal O[X]$, $\overline{P}=
%\overline{Q}\cdot\overline{R}$ and 
%$\overline{Q}, \overline{R}\in \k[X]^{\ne}$ are relatively prime, 
%then there exist $Q^*,R^*\in\mathcal O[X]$ with 
%$P=Q^*R^*$, $\overline{Q^*}=\overline{Q}$, $\overline{R^*}=\overline{R}$, 
%and $\deg Q^*=\deg \overline{Q}$.
%\end{enumerate}
\end{cor}
\begin{proof} The leading coefficient $a$ of $P$ lies in $\mathcal{O}^\times$ 
by the previous lemma, and so the desired result
follows from Proposition~\ref{prop:char henselian}(iii) applied to 
the monic polynomial $P/a\in \mathcal O[X]$. 
%Let $P$, $Q$, $R$ be as in the hypothesis of (2); to get the conclusion 
%of (2), we may assume that $\deg \overline{Q}=\deg Q$. By 
%Corollary~\ref{cor:gauss} we have $P=P_1\cdots P_m$ where each 
%$P_i\in\mathcal O[X]$ is irreducible in $K[X]$. 
%By (1), if $\overline{P_i}\notin\k$ then $\deg \overline{P_i}=\deg P_i$ 
%and $\overline{P_i}$ is a $\k$-multiple of a power of an irreducible 
%polynomial in $\k[X]$. Therefore, since 
%$\overline{P}=\overline{P_1}\cdots\overline{P_m}=\overline{Q}\overline{R}$, 
%if $\overline{P_i}\notin\k$ then $\overline{P_i}$ divides exactly 
%one of $\overline{Q}$ or $\overline{R}$. Let $Q_0$ be the product 
%of the polynomials $P_i$ such that $\overline{P_i}\notin\k$ and 
%$\overline{P_i}$ divides $\overline{Q}$, with $Q_0=1$ if there is 
%no such $P_i$. Then $\overline{Q_0}=\overline{c}\cdot\overline{Q}$ 
%for some $c\in\mathcal O^\times$,
%and $\deg Q_0=\deg\overline{Q}$. Set $Q^*:=c^{-1}Q_0\in\mathcal O[X]$ and
%take $R^*\in K[X]$ such that $P=Q^*R^*$; then  
%$R^*\in\mathcal O[X]$ and $\overline{R^*}=\overline{R}$.
\end{proof}

\noindent
We also note a consequence of Corollary~\ref{cor:unique val, 1} and 
Proposition~\ref{prop:char henselian}:

\begin{cor}\label{cor:unique val, 2}
If $K$ is henselian and $x$ an element of an algebraic valued 
field extension of $K$, with
minimum polynomial 
$$a_0+a_1X+\cdots+a_{n-1}X^{n-1}+X^n\qquad (a_0,\dots,a_{n-1}\in K,\ n\geq 1)$$
over $K$, then $v(x)=\frac{1}{n}v(a_0)$.
\end{cor}

\noindent
The henselian axiom concerns polynomials over $\mathcal O$ in a single variable, but implies an analogue for multivariate polynomials. To discuss this, let $n\ge 1$, and let
$P=(P_1,\dots,P_n)$ be an $n$-tuple of polynomials $P_i \in K[Y_1,\dots,Y_n]$. 
For $a\in K^n$ we have $P(a):=\big(P_1(a),\dots,P_n(a)\big)\in K^n$, and we recall from Section~\ref{sec:derivations on field exts} that 
$$P'(a)\ :=\ \left(\frac{\partial P_i}{\partial Y_j}(a)\right)_{i,j=1,\dots,n} \qquad\text{(an $n\times n$-matrix 
with entries in $K$)}.$$
We equip the additive abelian group $K^n$ with the valuation 
$$v\colon K^n\to\Gamma_\infty,\qquad
v(a)\ :=\ \min_i v(a_i)\quad\text{for $a=(a_1,\dots,a_n)\in K^n$.}$$
We consider here the elements of $K^n$ as column vectors, that is, as
$n\times 1$-matrices with entries in $K$. In particular, this holds for
the above $a\in K^n$ and $P(a)\in K^n$.

\begin{prop}\label{prop:hensel, multivar} Let $K$ be henselian, let $P_1,\dots, P_n\in \mathcal{O}[Y_1,\dots,Y_n]$, $n\ge 1$, and let $a\in \mathcal{O}^n$. Then:
\begin{enumerate}
\item[\textup{(i)}] if $v\big(P(a)\big)>0$ and $v\big(\det P'(a)\big)=0$, then there is a unique $b\in\mathcal O^n$ with $P(b)=0$ and $v(a-b)>0$;
\item[\textup{(ii)}] Newton version: if $v\big(P(a)\big)> 2v\big(\det P'(a)\big)$, then there is a unique $b\in\mathcal O^n$ with $P(b)=0$ and $v(a-b)>v\big(\det P'(a)\big)$. 
%\textup{(}Newton version\textup{)}
\end{enumerate}
\end{prop}
\begin{proof} To prove (i), 
let $v\big(P(a)\big)>0$ and $v\big(\det P'(a)\big)=0$. 
The square matrix $J:=P'(a)$ has entries in $\mathcal O$ with 
$\det J\in\mathcal O^\times$, so has inverse $J^{-1}$ with entries in 
$\mathcal O$. For $y\in \smallo^n$ we have by Taylor expansion
$$P(a+y)\ =\ P(a) + J\cdot y + Q(y), \qquad Q=(Q_1,\dots, Q_n),$$
where the $Q_i\in \mathcal{O}[Y_1,\dots, Y_n]$ are independent of $y$
and contain only monomials of degree $\ge 2$. Thus for $y\in \smallo^n$, 
$$ J^{-1}\cdot P(a+y)\ =\ J^{-1}\cdot P(a) + y + R(y)\qquad R=(R_1,\dots, R_n),$$
where the $R_i\in \mathcal{O}[Y_1,\dots, Y_n]$ are independent of $y$
and contain only monomials of degree $\ge 2$. Thus with $c:= J^{-1}\cdot P(a)$
we have $c\in \smallo^n$, and it remains to show that there is exactly one
$y\in \smallo^n$ with $c+y+R(y)=0$, that is, $-c-R(y)=y$. Consider the map
$f\colon \smallo^n\to \smallo^n$ given by $f(y)=-c-R(y)$.
For $x,y\in \smallo^n$ we have $f(x)-f(y)=R(y)-R(x)$, and so
$v\big(f(x)-f(y)\big) > v(x-y)$ if $x\ne y$. If~$K$ is maximal as a valued field, then the valued abelian group $K^n$ is spherically complete, by Lemma~\ref{lem:direct product of vag} and Corollary~\ref{cor:maximal valued fields}, and so $f$ has a unique fixpoint 
by Theorem~\ref{thm:fixpoint}, and (i) holds. Suppose $K$ is not maximal.
Then we take some algebraically closed maximal valued field extension 
$L$ of $K$ and obtain a unique~$y\in \smallo_L^n$ with $P(a+y)=0$.
Then the entries $y_i$ of the vector $y$ are separably algebraic over $K$ by Corollary~\ref{cor:regsepalg}. 
%(which is also \cite[Proposition~VIII.5.3]{Lang}). 

Let~$K^\alg$ be the algebraic closure of~$K$ in $L$, so 
$\mathcal{O}^{\alg}:= \mathcal{O}_L\cap K^\alg$ is the unique valuation
ring of $K^\alg$ lying over $\mathcal{O}$. Let~$\smallo^{\alg}$ be the maximal ideal
of $\mathcal{O}^{\alg}$, so $y\in (\smallo^{\alg})^n$.
Then for each $\sigma\in\Aut(K^\alg|K)$,
$$P\big(a+\sigma(y)\big)\ =\ \sigma\big(P(a+y)\big)\ =\ 0, \quad \sigma(y)\in (\sigma \smallo^{\alg})^n= (\smallo^{\alg})^n,$$
so $\sigma(y)=y$ by uniqueness of $y$. Therefore $y\in\smallo^n$.
This proves (i). 

For (ii), set $J:= P'(a)$, and put $d:= \det J$, and assume 
$v\big(P(a)\big)> 2v(d)$. 
As in the proof of (i) we have $Q_1,\dots, Q_n\in \mathcal{O}[Y_1,\dots, Y_n]$, each having only monomials of degree $\ge 2$, such that for all $y\in K^n$,
$$ P(a+y)\ =\ P(a) + J\cdot y + Q(y), \qquad Q=(Q_1,\dots, Q_n).$$
Note that $Q_i(dY_1,\dots, dY_n)=d^2R_i(Y_1,\dots, Y_n)$ with $R_i\in \mathcal{O}[Y_1,\dots, Y_n]$.
Then
$$P(a+dy)\ =\ P(a)+ J\cdot dy + d^2R(y), \qquad R=(R_1,\dots,R_n),$$
for all $y\in K^n$. Now $J^{-1}$ has all its entries in $d^{-1}\mathcal{O}$, so for all $y\in K^n$,
$$J^{-1}\cdot P(a+dy)\ =\ J^{-1}\cdot P(a) + dy + dS(y), \qquad S=(S_1,\dots, S_n)$$
with $S_i\in \mathcal{O}[Y_1,\dots, Y_n]$ independent of $y$ and all its monomials of degree $\ge 2$.
So 
$$d^{-1}J^{-1}\cdot P(a+dy)\ =\ d^{-1}J^{-1}\cdot P(a) + y + S(y),$$
for all $y\in K^n$, and $v\big(d^{-1}J^{-1}\cdot P(a)\big)> 0$. Thus (i) gives a unique $y\in\smallo^n$ with $P(a+dy)=0$, which is the conclusion of (ii).
\end{proof}

\noindent
In the corollary below we let $X_1,\dots,X_m, Y_1,\dots,Y_n$ be distinct indeterminates and set $X=(X_1,\dots,X_m)$, $Y=(Y_1,\dots,Y_n)$.
Given an $n$-tuple  $P=(P_1,\dots,P_n)$ of polynomials $P_i\in K[X,Y]$ and
$a\in K^m$, $b\in K^n$ we define the $n\times n$-matrix
$$\frac{\partial P}{\partial Y}(a,b)\ :=\ \left(\frac{\partial P_i}{\partial Y_j}(a,b)\right)_{i,j=1,\dots,n}$$
with entries in $K$.

\begin{cor}[Implicit Function Theorem]\label{cor:IFT}
Suppose $K$ is henselian. Let $P=(P_1,\dots,P_n)$ where $P_i\in \mathcal O[X,Y]$,
$m,n\geq 1$.  Let $a\in \mathcal O^m$, $b\in \mathcal O^n$ be such that 
$$P(a,b)\ =\ 0, \qquad 
\det\frac{\partial P}{\partial Y}(a,b)\ \neq\  0,$$
and set $\delta:=v\left(\det\frac{\partial P}{\partial Y}(a,b)\right)\in\Gamma$.
Then  for all $x\in\mathcal O^m$ with $v(x-a)>2\delta$
there is a unique $y\in\mathcal O^n$ such that 
$P(x,y)=0$ and $v(y-b)>\delta$.
\end{cor}
\begin{proof}
Let $x\in\mathcal O^m$, $v(x-a)>2\delta$. Taylor expansion yields
$$P(x,b)\ =\  P(a,b) + Q(x-a)\ =\  Q(x-a),\qquad Q=(Q_1,\dots,Q_n)$$
where the $Q_i\in \mathcal{O}[X_1,\dots, X_m]$ contain only monomials of degree~$\geq 1$, and so
$v\big(P(x,b)\big)\geq v(x-a)>2\delta$. Similarly,
$$\det\frac{\partial P}{\partial Y}(x,b)\ =\ \det\frac{\partial P}{\partial Y}(a,b) + R(x-a)\quad\text{
where  $R\in \mathcal O[X]$, $R(0)=0$,}$$   so
$v\left(\det\frac{\partial P}{\partial Y}(x,b)\right) = v\left(\det\frac{\partial P}{\partial Y}(a,b)\right)=\delta$. Now apply Proposition~\ref{prop:hensel, multivar}(ii) to   $\big(P_1(x,Y),\dots,P_n(x,Y)\big)$ and~$b$ in place of~$P$ and~$a$.
\end{proof}

\subsection*{Hensel configuration and algebraic maximality} 
For $K$ to be henselian is a condition on polynomials over its valuation ring $\mathcal O$. It is convenient to have an equivalent condition in terms of polynomials over $K$. Let $P\in K[X]$ be of degree at most~$d$, and let $a\in K$, so
$$P_{+a}(X)\ :=\ P(a+X)\ =\ \sum_{i=0}^d P_{(i)}(a)X^i\ =\  P(a) + P'(a)X + \sum_{i=2}^d P_{(i)}(a)X^i.$$
We say that $P$ is in {\bf hensel configuration} at $a$ if $P'(a)\neq 0$, and either $P(a)=0$ or 
$P(a)\ne 0$ and $\gamma:= v\big(P(a)\big)-v\big(P'(a)\big)$ satisfies
$$v\big(P(a)\big)\ <\  v\big(P_{(i)}(a)\big)+i\gamma\quad\text{for $i=2,\dots,d$.}$$

\index{hensel configuration} \index{configuration!hensel}

\begin{lemma}\label{lem:hensel config}
Suppose $K$ is henselian and $P$ as above is in hensel configuration at $a\in K$. Set $\gamma:= v\big(P(a)\big)-v\big(P'(a)\big)\in \Gamma_{\infty}$. Then there is a unique $b\in K$ such that $P(a+b)=0$ and $v(b)\geq \gamma$; this $b$ satisfies $v(b)=\gamma$.
\end{lemma}
\begin{proof} 
If $P(a)=0$, then $b:=0$ works, so assume $P(a)\neq 0$. Take $g\in K$ such that 
$vg=\gamma$, and set $Q:=P_{+a}$, $h:=P(a)$, and  
$$\tilde{Q}(X)\ :=\ Q(gX)/h\  =\ 1 + \underbrace{\big(P'(a)g/h\big)}_{\asymp 1}X +  \sum_{i=2}^d \underbrace{\big(P_{(i)}(a)g^i/h\big)}_{\prec 1}X^i.$$
Since $K$ is henselian and the image of $\tilde{Q}$ under the natural surjection $\mathcal O[X]\to\k[X]$ has degree $1$,
the polynomial $\tilde{Q}$ has a unique zero $u\in\mathcal O$. We have $u\asymp 1$ and $P(a+ug)=0$, hence $b:=ug$ has the required properties. 
\end{proof}

\begin{prop}\label{prop:hensel config}
Suppose $K$ is of equicharacteristic zero. Let $(a_\rho)$ be a pc-sequence in $K$ with $a_\rho\leadsto a\in K$, and set $\gamma_\rho:= v(a-a_{\rho})$. 
Let $P\in K[X]$ of degree~$\le d$ be such that $P(a_\rho)\leadsto 0$ and $P_{(i)}(a_\rho)\not\leadsto 0$ for $i= 1,\dots,d$. Then for $i=2,\dots,d$, we have, eventually,
$$v\big(P(a_{\rho})\big)\ =\ v\big(P(a_\rho)-P(a)\big)\ =\  v\big(P'(a)\big) + \gamma_\rho\ <\  v\big(P_{(i)}(a)\big)+i\gamma_\rho.$$
\end{prop}
\begin{proof}
The proof of Proposition~\ref{prop:pseudocontinuity} yields a unique $i_0\geq 1$
such that for each~$i\geq 1$ with $i\neq i_0$,
$$v\big(P(a_\rho)-P(a)\big)\ =\  v\big(P_{(i_0)}(a)\big) + i_0\gamma_\rho\ <\  v\big(P_{(i)}(a)\big)+i\gamma_\rho,\quad\text{eventually.}$$
Now $P(a_\rho)\leadsto 0$, so $v\big(P(a_{\rho})\big)=v\big(P(a_\rho)-P(a)\big)$, eventually, and for $i\geq 1$, $i\neq i_0$:
$$v\big(P(a_\rho)\big)\ =\  v\big(P_{(i_0)}(a)\big)+i_0\gamma_\rho\ <\   v\big(P_{(i)}(a)\big)+i\gamma_\rho,\quad\text{eventually.}$$
We claim that $i_0=1$. Let $i>1$; our claim will then follow by deriving
$$v\big(P'(a)\big)+\gamma_\rho\ \le\  v\big(P_{(i)}(a)\big)+i\gamma_\rho,\quad\text{eventually.}$$
Now the proof of Proposition~\ref{prop:pseudocontinuity} applied to $P'$ instead of $P$ also yields
$$v\big(P'(a_\rho)-P'(a)\big)\ \leq\ v\big(P'_{(j)}(a)\big)+j\gamma_\rho,\quad\text{eventually}$$
for all $j\geq 1$. Since $v\big(P'(a_\rho)\big)=v\big(P'(a)\big)$ eventually, this yields
$$v\big(P'(a)\big)\ \leq\ v\big(P'_{(j)}(a)\big)+j\gamma_\rho,\quad\text{eventually}$$
for all $j\geq 1$.
Using $P'_{(j)}=(1+j)P_{(1+j)}$ (Lemma~\ref{lem:iterated Taylor}), this gives for $j\geq 1$:
$$v\big(P'(a)\big)\ \leq\  v\big(P_{(1+j)}(a)\big)+j\gamma_\rho,\quad\text{eventually.}$$
For $j=i-1$, this yields
$$v\big(P'(a)\big)+\gamma_\rho\ \le\  v\big(P_{(i)}(a)\big)+i\gamma_\rho,\quad\text{eventually.}$$
Thus $i_0=1$ as claimed.  
\end{proof}

\begin{cor}\label{cor:hensel config}
Suppose $K$ is of equicharacteristic zero.
Let $(a_\rho)$ be a pc-sequence in $K$, and let $P\in K[X]$ be such that $P(a_\rho)\leadsto 0$ and $P_{(i)}(a_\rho)\not\leadsto 0$ for all~$i\geq 1$.
Then $P$ is in hensel configuration at $a_{\rho}$, eventually, and
in any henselian valued field extension of $K$ there is a unique $b$ such that $a_\rho\leadsto b$ and $P(b)=0$.
\end{cor}
\begin{proof}
Let $a$ be a pseudolimit of $(a_\rho)$ in some valued field extension of $K$ (whose valuation we continue to denote by $v$), and set $\gamma_\rho := v(a-a_\rho)$. 
Since for each~${i\geq 1}$ we have $v\big(P_{(i)}(a_\rho)\big)=v\big(P_{(i)}(a)\big)$, eventually, Proposition~\ref{prop:hensel config} shows that $P$ is in hensel configuration at $a_{\rho}$, eventually.
Let $K'$ be a henselian valued field extension of $K$. After deleting an initial segment of the sequence $(a_\rho)$ we can assume that $(\gamma_\rho)$ is strictly increasing, and that $v(a_\sigma-a_\rho)=\gamma_\rho$ whenever $\sigma>\rho$, and that $P'(a_\rho)\neq 0$ for all $\rho$. Likewise, by~Lemma~\ref{lem:hensel config} and Proposition~\ref{prop:hensel config} we can assume that for every $\rho$ there is a unique $b_\rho\in K'$ such that $P(b_\rho)=0$ and $v(a_\rho-b_\rho)=v\big(P(a_\rho)\big)-v\big(P'(a_\rho)\big)$. Proposition~\ref{prop:hensel config} also shows that we can assume that $v\big(P(a_\rho)\big)-v\big(P'(a_\rho)\big)=\gamma_\rho$ for all $\rho$, so $v(a_\rho-b_\rho)=\gamma_\rho$ for all $\rho$.
The uniqueness of $b_\rho$ yields $b_\sigma=b_\rho$ whenever $\sigma>\rho$. Thus all $b_\sigma$ are equal to a single~$b$, which has the desired properties.
\end{proof}

\begin{cor}\label{cor:alg max equals henselian}
Suppose $K$ is of equicharacteristic zero. Then $K$ is henselian if and only if $K$ is algebraically maximal.
\end{cor}
\begin{proof}
We already know  that algebraically maximal valued fields are henselian (Co\-rol\-lary~\ref{cor:alg max implies henselian}).
Assume $K$ is henselian, and let $(a_\rho)$ be a pc-sequence in $K$ of algebraic type over $K$. Take a monic minimal
polynomial $P$ of $(a_\rho)$ over $K$. Then $P(a_\rho)\leadsto 0$ and $P_{(i)}(a_\rho)\not\leadsto 0$ for each $i\geq 1$, so $P$ is in hensel configuration at $a_{\rho}$, eventually, and thus Lemma~\ref{lem:hensel config} gives $a\in K$ such that $P(a)=0$. Since $P$ is irreducible, we obtain $P=X-a$, so $a_\rho\leadsto a$.
\end{proof}

\noindent
This leads to a partial converse of some earlier results (see~\ref{cor:acvf} and \ref{cor:alg max implies henselian}):

\begin{cor}\label{achensacdiv}
Suppose $K$ is of  equicharacteristic zero. Then $K$ is algebraically closed if and only if $K$ is henselian, $\res(K)$ is algebraically closed, and~$\Gamma$ is divisible.
\end{cor}

\begin{exampleNumbered}[from Laurent series to Puiseux series]\label{ex:Puiseux}\index{series!Puiseux}
Let $C$ be a field. The elements of the valued subfield
$$\operatorname{P}(C)\ :=\ \bigcup_{n\geq 1} C \(( t^{\frac{1}{n}\Z} \)) $$
of the Hahn field $C \(( t^\Q \)) $ are called {\bf Puiseux series} over~$C$ (in $t$). Note that $\operatorname{P}(C)$ has the same residue field (isomorphic to $C$) as its extension
$C \(( t^\Q \)) $ and the same value group $\Q$, so $C \(( t^\Q \)) $ is an
immediate extension of $\operatorname{P}(C)$.
For $n\geq 1$ we have
$$C \(( t^{\frac{1}{n}\Z} \)) \ =\ C \(( t \)) + C \(( t \)) a
+\cdots+ C \(( t \)) a^{n-1}\ =\ C \(( t \)) [a], \qquad a\ :=\ t^{\frac{1}{n}},$$
so $C \(( t^{\frac{1}{n}\Z} \)) $ is algebraic of degree $n$
over $C \(( t \)) $. It follows that
$\operatorname{P}(C)$ is algebraic over the Laurent series field $C \(( t \)) = C\(( t^{\Z} \)) $. Now $C \(( t^{\frac{1}{n}\Z} \)) $ is henselian for every~${n\ge 1}$, 
by Corollary~\ref{cor:alg max implies henselian}, so $\operatorname{P}(C)$ is henselian.
Hence by Corollary~\ref{achensacdiv}, if~$C$ is algebraically closed of characteristic zero, then $\operatorname{P}(C)$
is algebraically closed, and thus~$\operatorname{P}(C)$ is the algebraic closure of~$C \(( t \)) $ in~$C \(( t^\Q \)) $.
\end{exampleNumbered}

%\begin{example}
%Let $C$ be a field of characteristic zero and let $\mathfrak M$ be an ordered abelian group (written multiplicatively). Then the valued field $C[[\mathfrak M]]$ is algebraically closed iff $C$ is algebraically closed and $\mathfrak M$ is divisible.
%\end{example}

\begin{cor}\label{cor:pc-limits in simple trans exts}
Suppose $K$ is henselian of equicharacteristic zero. Let~$L=K(y)$ be a valued field extension of $K$ with
$nv(y)\notin \Gamma$ for all $n\ge 1$. Then no divergent pc-sequence in $K$ has a pseudolimit in $L$.
\end{cor}
\begin{proof} Let $(a_\rho)$ be a divergent pc-sequence in $K$ with
$a_{\rho}\leadsto a\in L$. Then~$(a_\rho)$ is of transcendental type over $K$, 
by Corollaries~\ref{cor:algebraically maximal} and \ref{cor:alg max equals henselian},  so~$K(a)$ is an immediate extension of $K$ by 
Lemma~\ref{lem:Kaplansky, 1}. But $a$ is transcendental over $K$ and
so $y$ is algebraic over $K(a)$, and thus $nv(y)\in\Gamma_{K(a)}=\Gamma$ for some $n\ge 1$, a contradiction.
\end{proof}

\subsection*{Henselization}
Let $K$ be a valued field. A {\bf henselization of $K$} is a henselian valued field extension~$K^{\operatorname{h}}$ of $K$ such that any valued field embedding $K\to L$ into a henselian valued field $L$ extends uniquely to an embedding $K^{\operatorname{h}}\to L$. \index{henselization of a valued field} \index{valued field!henselization} \nomenclature[K]{$K^{\operatorname{h}}$}{henselization of $K$}

\begin{prop} Every valued field has a henselization.
\end{prop}
\begin{proof}
Fix an algebraic closure $K^{\operatorname{a}}$ of the underlying field of $K$ and pick a valuation ring $\mathcal{O}^{\operatorname{a}}$ of $K^{\operatorname{a}}$ that lies
over $\mathcal{O}$. Then $(K^{\operatorname{a}}, \mathcal{O}^{\operatorname{a}})$ is a henselian valued field extension of $K$, and the intersection $E$ of the subfields
$F\supseteq K$ of $K^a$ such that~$(F, \mathcal{O}^a \cap F)$ is henselian is itself a henselian valued field extension of $K$ with valuation ring $\mathcal{O}_E:=
\mathcal{O}^a\cap E$. We claim that $E$ is a henselization
of $K$. To see this, let a valued field embedding $i\colon K \to L$ into a henselian valued field $L$ be given. We have to show that~$i$ extends uniquely to a valued field embedding $E \to L$.
We can arrange that $K$ is a
valued subfield of $L$ with $i$ the inclusion map. Replacing~$L$ by
the relative algebraic closure of $K$ in $L$ we also reduce to the case that
$L$ is algebraic over $K$. Take an algebraic closure $L^a$ of $L$ and take a valuation ring $V$ of $L^a$ that lies over $\mathcal{O}_L$. Then
Proposition~\ref{prop:transitivity} gives a valued field isomorphism $\sigma\colon (K^{\operatorname{a}}, \mathcal{O}^{\operatorname{a}}) \to (L^{\operatorname{a}},V)$ over~$K$, and so $\sigma$ maps $E$ onto the smallest
henselian valued subfield of $(L^{\operatorname{a}}, V)$ containing~$K$, and
thus $\sigma(E)\subseteq L$. So we have a valued field embedding $\sigma|E\colon E\to L$ over~$K$. Let $j\colon E \to L$ also be a valued field embedding over $K$; it is enough to
show that then $j=\sigma|E$. Set $F:=\big\{a\in E: j(a)=\sigma(a)\big\}$. Then
$F\supseteq K$ is a subfield of $E$, and we take it as a valued subfield of
$E$. Then $F$ is henselian: let $P\in \mathcal{O}_F[X]$ and $a\in \mathcal{O}_F$ be such that $P(a)\prec 1$ and $P'(a)\asymp 1$. Then~$P$ has a unique zero
$b\in E$ with $a- b\prec 1$, and likewise, the image of $P$ under~$j$ in $L[X]$
has a unique zero~$c$ in~$L$ with $j(a)- c\prec 1$, and so $j(b)=c$. For the same
reason, $\sigma(b)=c$, and so $b\in F$. We have now shown that $F$ is henselian. Hence $E=F$ by the minimality of $E$,
and thus~${j=\sigma|E}$. 
\end{proof}

\noindent
If $K_1$ and $K_2$ are henselizations of $K$, then the unique embedding $K_1\to K_2$
over $K$ is an isomorphism; thus there is no harm in referring to \textit{the}\/ henselization of $K$. 
Of course, if $K$ is henselian, then it is its own henselization. 
We already know that $K$ has an immediate algebraically maximal (and thus henselian) valued field extension that is algebraic over $K$. Hence:

\begin{cor}
The henselization of $K$ is an immediate valued field extension of~$K$ and is algebraic over $K$.
\end{cor}

\noindent
If $K'$ is any henselian valued field extension of $K$, then by the {\bf henselization of~$K$ in~$K'$} we mean the henselization $K^{\operatorname{h}}$ of $K$ with $K\subseteq K^{\operatorname{h}}\subseteq K'$ as valued fields. 

%\begin{lemma}
%Let $K^\alg$ be an algebraic closure of $K$ equipped with a valuation ring lying over that of $K$, and suppose $K$ has henselization $K^{\operatorname{h}}$  in~$K^\alg$. Let $L$ be a valued subfield of $K^\alg$ containing $K$, having  henselization $L^{\operatorname{h}}$ in~$K^\alg$. Then $L^{\operatorname{h}}=LK^{\operatorname{h}}$. 
%\end{lemma}
%\begin{proof}
%By Corollary~\ref{cor:alg ext of henselian}, $LK^{\operatorname{h}}$ is a henselian valued field extension of $L$; hence $L^{\operatorname{h}}\subseteq LK^{\operatorname{h}}$. Conversely, $L^{\operatorname{h}}$ also contains the valued subfields $L$ and $K^{\operatorname{h}}$, hence we have $L^{\operatorname{h}}=LK^{\operatorname{h}}$. 
%\end{proof}

%\noindent
%Note that if in the situation of the previous lemma the field extension $L\supseteq K$ is finite, then so is $L^{\operatorname{h}}\supseteq K^{\operatorname{h}}$.

\medskip
\noindent
We can now describe the absolute Galois group of the henselization of $K$: 

\begin{lemma}\label{lem:abs Galois gp of henselization}
Let $K^\alg$ be an algebraic closure of $K$ and $\mathcal{O}^{\alg}$ a valuation ring of~$K^\alg$ lying over $\mathcal O$. Let $K^{\operatorname{h}}$ be the henselization of $K$ in~$(K^\alg,\mathcal{O}^{\alg})$. Then
$$\Aut(K^\alg|K^{\operatorname{h}})\ =\ \big\{\sigma\in\Aut(K^\alg|K):\ \sigma(\mathcal O^\alg)=\mathcal O^\alg\big\}.$$
\end{lemma}
\begin{proof}
Let $\sigma\in \Aut(K^\alg|K^{\operatorname{h}})$. Then $\sigma(\mathcal O^\alg)$ and $\mathcal O^\alg$ are valuation rings of $K^\alg$ lying over the valuation ring of $K^{\operatorname{h}}$, hence they are equal by 
Proposition~\ref{prop:char henselian}. Conversely, let $\sigma\in\Aut(K^\alg|K)$ with $\sigma(\mathcal O^\alg)=\mathcal O^\alg$. Then $\sigma$ is an automorphism of the valued field $(K^\alg,\mathcal{O}^{\alg})$, so its valued subfields $K^{\operatorname{h}}$ and $\sigma(K^{\operatorname{h}})$ are both henselizations of~$K$ in~$(K^\alg,\mathcal{O}^{\alg})$, with isomorphism $\sigma|K^{\operatorname{h}}$ between them.  Therefore $K^{\operatorname{h}}=\sigma(K^{\operatorname{h}})$ and~${\sigma|K^{\operatorname{h}}=\id}$.
\end{proof}

\noindent
This leads to the following top-down construction of $K^{\operatorname{h}}$ when $\operatorname{char}(K)=0$:

\begin{cor}\label{cor:abs Galois gp of henselization}
Assume $\operatorname{char}(K)=0$, and let $K^\alg$, $\mathcal{O}^{\alg}$,  and $K^{\operatorname{h}}$ be as in the lemma above.
Then $K^{\operatorname{h}}$ is the decomposition field of $\mathcal O^\alg$ over $K$.
\end{cor}

%From Corollary~\ref{cor:Krull, 2} and the previous lemma we obtain a description of the part of the %henselization of $K$ contained in a given normal extension of $K$:

%\begin{cor}\label{cor:abs Galois gp of henselization}
%Let $K^\alg$ and $K^{\operatorname{h}}$ be as in the lemma above.
%Let $N$ be a normal valued field extension of $K$ inside $K^\alg$, and $G:=\Aut(N|K)$. 
%Then $K^{\operatorname{h}}\cap N$ is the decomposition field 
%$N^{\operatorname{d}}$ of $\mathcal O_N$ over $K$, i.e., the fixed field of the decomposition group 
%$$G^{\operatorname{d}}=\big\{\sigma\in G:\sigma(\mathcal O_N)=\mathcal O_N\big\}$$
%of $\mathcal O_N$ over $K$. 
%\end{cor}

%\noindent
%In the equicharacteristic zero case, we also have:

\begin{theorem}\label{thm:henselization}
Suppose $K$ is of equicharacteristic zero and $L$ is an immediate henselian valued 
field extension of $K$ and $L$ is algebraic over $K$. 
Then $L$ is a henselization of $K$.
\end{theorem} 

\begin{proof} By Corollary~\ref{cor:alg max equals henselian} the henselization $K^{\operatorname{h}}$ of $K$ in $L$ is algebraically maximal, and $L$ is an immediate
algebraic extension of $K^{\operatorname{h}}$, so $K^{\operatorname{h}}=L$. \end{proof}

%\begin{cor}\label{cor:henselization}
%Suppose the valued field $K$ has equicharacteristic zero. Then any algebraically maximal immediate valued %field extension of $K$ that is algebraic over $K$ is a henselization of $K$.
%\end{cor}

\noindent
We can now prove the following useful fact:
  
\begin{lemma}\label{lem:min poly immediate ext}
Suppose $K$ has equicharacteristic zero, and $L\supseteq K$ is an immediate valued field extension
with $[L:K]<\infty$. Then $L=K(x)$ for some $x\prec 1$ in $L$ 
whose minimum polynomial over $K$ has the form 
$$X^n + a_{n-1}X^{n-1} + \cdots + a_1X+a_0\qquad\text{with all $a_i\preceq 1$, $a_1\asymp 1$, $a_0 \prec 1$.}$$
\end{lemma}
\begin{proof}
Let $K^\alg$ be an algebraic closure of $K$ containing $L$ as a subfield and take a valuation ring $\mathcal O^{\alg}$ of $K^\alg$ lying above $\mathcal O_L$.
Let $L^{\operatorname{h}}$ be the henselization of $L$ in~$(K^\alg,\mathcal O^{\alg})$; then $L^{\operatorname{h}}$ is also the henselization of $K$ in $(K^\alg,\mathcal O^{\alg})$, by Theorem~\ref{thm:henselization}.
Let $N$ be a Galois extension of $K$ inside $K^\alg$ which contains $L$ with $[N:K]<\infty$. Then $\mathcal O_N=\mathcal O^{\alg}\cap N$ is a valuation ring of 
$N$. Let $G:=\Aut(N|K)$ and let 
$$G^{\operatorname{d}}=\big\{\sigma\in G:\sigma(\mathcal O_N)=\mathcal O_N\big\}$$
be the decomposition group of $\mathcal O_N$ over $K$, with fixed field~$N^{\operatorname{d}}$.
By Corollary~\ref{cor:abs Galois gp of henselization} we have $N^{\operatorname{d}}=L^{\operatorname{h}}\cap N\supseteq L$; in particular, $G^{\operatorname{d}}\subseteq\Aut(N|L)$.
Let~$A$ be the integral closure of $\mathcal O$ in $L$ and $B$ the integral closure of $\mathcal O$ in $N$, and~$\mathfrak q$ the maximal ideal of~$B$ such that $B_{\mathfrak q}=\mathcal O_N$.
$$\xymatrix@R-1em@C+0.5em{
						&													&	N \ar@{-}[d] & \ar@{<->}[dd]^{G}	\\
						& \mathcal O_N\ar@{-}[ur]	\ar@{-}[d]	& L				&	\\
B\ar@{-}[ur]	\ar@{-}[d]	& \mathcal O_L\ar@{-}[ur]							& K\ar@{-}[u] &					\\
A\ar@{-}[ur]	\ar@{-}[d]	&													&  &\\
\mathcal O\ar@{-}[uurr]	&													& &
}$$

\claim{Suppose $\sigma\in G$ and $\sigma(\mathfrak q)\cap A=\mathfrak q\cap A$.
Then $\sigma\in\Aut(N|L)$.}\\[-2em]

\begin{proof}[Proof of Claim] From $\mathcal{O}_N=B_{\mathfrak{q}}$ we get 
$\sigma(\mathcal{O}_N)=B_{\sigma(\mathfrak{q})}$, so $\mathcal O_N$ and $\sigma(\mathcal O_N)$ lie  over $A_{\mathfrak{q}\cap A}=\mathcal O_L$. Since $N\supseteq L$ is normal, Proposition~\ref{prop:transitivity} gives $\tau\in\Aut(N|L)$ with $\tau(\mathcal O_N)=\sigma(\mathcal O_N)$, so $\tau^{-1}\sigma\in G^{\operatorname{d}}\subseteq\Aut(N|L)$,  thus $\sigma\in\Aut(N|L)$.
\end{proof}

\noindent
By the claim, the (finite) sets 
$$\big\{\sigma(\mathfrak q)\cap A:\sigma\in\Aut(N|L)\big\}, \qquad \big\{\sigma(\mathfrak q)\cap A:\sigma\in G\setminus\Aut(N|L)\big\}$$ 
of maximal ideals of $A$ are disjoint. Hence the Chinese Remainder Theorem (below Proposition~\ref{prop:transitivity})
%(p.~\pageref{p:CRT})
gives~${x\in A}$ such that 
$$\text{$x\in\sigma(\mathfrak q)$ for all 
$\sigma\in\Aut(N|L)$,} \qquad x\notin\sigma(\mathfrak q)\ \text{ for all $\sigma\in G\setminus \Aut(N|L)$.}$$ We show that $x$ has the required properties.
Let $H$ be the subgroup of $G$ whose fixed field in $N$ is $K(x)$. From $K(x)\subseteq L$ we get $\Aut(N|L)\subseteq \Aut(N|K(x))=H$. On the other hand,
$x\in\mathfrak q$ gives
$x=\sigma(x)\in\sigma(\mathfrak q)$ for all $\sigma\in H$, so
$H\subseteq\Aut(N|L)$; thus $H=\Aut(N|L)$ and $L=K(x)$. 

Now fix a coset decomposition $G=\sigma_1 H\cup\cdots\cup\sigma_n H$ of $G$ with respect to its subgroup $H$, where $\sigma_1=\id$ and $n=[G:H]=[L:K]$. Then $\sigma_1|L,\dots,\sigma_n|L$ are the distinct $K$-embeddings of $L$ into $N$. The minimum polynomial of $x$ over $K$ is
$$P(X)\ =\ \prod_{i=1}^n \big(X-\sigma_i(x)\big)\ =\ X^n+a_{n-1}X^{n-1}+\cdots+a_1X+a_0 \qquad \text{(all $a_j\in K$).}$$
An argument similar to the one at the end of the proof of the implication (i)~$\Rightarrow$~(ii) in Proposition~\ref{prop:char henselian}  now shows that 
$a_j\preceq 1$ for all $j$ and $a_0\prec 1$, $a_1\asymp 1$.
%Since each $\sigma_i(x)$ is integral over $\mathcal O$, so are the $P_j$, hence $P_j\in A\cap K=\mathcal O=K^{\preceq 1}$. Now
%$$P_0 = (-1)^n\prod_{i=1}^n \sigma_i(x)=x\cdot \text{(element of $A$}) \in \mathfrak q\cap K=\smallo=K^{\prec 1}.$$
%Moreover
%$$P_1 = (-1)^{n-1}\sum_{i=1}^n \sigma_1(x)\cdots\sigma_{i-1}(x)\sigma_{i+1}(x)\cdots\sigma_n(x),$$
%and precisely one term in this sum, namely $\sigma_2(x)\cdots\sigma_n(x)$, misses the factor $\sigma_1(x)=x$, hence this is the only term in this sum which is not in $\mathfrak q$; thus $P_1\notin K^{\prec 1}$.
\end{proof}

\noindent
The next lemma describes the valuation ring $\mathcal{O}_L$ of $L$ in Lemma~\ref{lem:min poly immediate ext}:

\begin{lemma}\label{lem:min poly immediate ext, addition} 
Let $L=K(x)$ be a separable algebraic field extension of $K$ where $x\prec 1$ and the minimum polynomial of $x$ over $K$ has the form 
$$X^n + a_{n-1}X^{n-1} + \cdots + a_1X+a_0\qquad\text{with all $a_i\preceq 1$, $a_1\asymp 1$, $a_0 \prec 1$.}$$
Then $\smallo(x):= \smallo\mathcal O[x] + x\mathcal O[x]$
is a maximal ideal of $\mathcal O[x]$
and $\mathcal O_L = \mathcal O[x]_{\smallo(x)}$. 
\end{lemma}

\begin{proof}
By Lemma~\ref{lem:Dedekind}, $\smallo(x)$ is a maximal ideal
of $\mathcal O[x]$. Let $B$ be the integral closure of $\mathcal O$ in $L$. Then $\mathcal O_L=B_{\mathfrak{n}}$ for 
a maximal ideal $\mathfrak{n}$ of $B$, by Proposition~\ref{prop:Krull, 3}. Hence $x\in \mathfrak{n}$, so
$\smallo(x)=\mathfrak{n}\cap \mathcal O[x]$, and thus 
$\mathcal O[x]_{\smallo(x)}=B_{\mathfrak{n}}$ by
Corollary~\ref{cor:fg int closure}. 
\end{proof}

\subsection*{Monomial groups} A {\bf monomial group\/} of the valued field~$K$ is a subgroup $\frak{M}$
of $K^\times$ that is mapped bijectively onto $\Gamma$ by $v$. This notion is equivalent to that of a {\bf cross-section~$s$} of the valued field $K$, that is, a group morphism $s\colon\Gamma\to K^\times$ such that $v(s\gamma)=\gamma$
for every~$\gamma\in\Gamma$. Indeed, for such a cross-section $s$, the set~$s(\Gamma)$ is a monomial group, and every monomial group of $K$ is of this form for a unique
cross-section $s$ of $K$. For example, if $C$ is a field and $\mathfrak M$ an ordered abelian multiplicative group and $K=C[[\mathfrak M]]$ the associated Hahn field, then the group $\mathfrak M$, identified with a subgroup of $K^\times$ in the natural way, is a monomial group of $K$. Other examples for monomial groups are provided by the following lemma: \index{valued field!monomial group} \index{monomial! group}\index{group!monomial}  \index{valued field!cross-section} \index{cross-section}

\begin{lemma}\label{divmon}
Let $G$ be a divisible subgroup of $K^\times$ with $v(G)=\Gamma$. Then there exists a cross-section $s$ of $K$ with $s(\Gamma)\subseteq G$.
\end{lemma}
\begin{proof}
The inclusion $\mathcal O^\times\cap G\to G$ and the restriction of the valuation 
$v\colon K^\times \to \Gamma$ to $G$ yield an 
exact sequence $$1 \to \mathcal O^\times\cap G\to G \to \Gamma \to 0$$ of abelian 
groups. Since $G$ is divisible, so is $\mathcal O^\times\cap G$,  hence this exact sequence splits, which gives a section $s$ as claimed. 
\end{proof}

\noindent
Thus algebraically closed valued fields have monomial groups.
Likewise, if $K$ is real closed (see Section~\ref{sec:valued ordered fields} below), then $K$ has a monomial group contained in $K^{>}$. It is also easy to see that if $K$ is henselian with residue field $\k=\res(K)$ of characteristic zero, and the abelian groups $\k^\times$ and $\Gamma$ are divisible, then $K^\times$ is divisible, hence $K$ has a monomial group.

\begin{lemma}\label{exmongr} Let $K$ be an algebraically closed valued field and $\frak{M}_0$
a monomial group of a valued subfield $K_0$ of $K$ over which $K$ is algebraic. 
Then $\frak{M}_0$ extends to a monomial group $\frak{M}$ of $K$.
\end{lemma}
\begin{proof} Let $s_0\colon \Gamma_0=v(K_0^\times)\to K_0^\times$ 
be the cross-section of $K_0$ with $s_0(\Gamma_0)=\frak{M}_0$. Since the abelian
group $K^\times$ is divisible, $s_0$ extends to a group morphism 
$$s\ \colon\ \Gamma=v(K^\times)\to K^\times.$$ Using $\Gamma=\Q\Gamma_0$ it follows easily
that $\frak{M}:=s(\Gamma)$ is a monomial group of $K$.
\end{proof} 

\noindent
Next we want to show that every valued field (possibly equipped with additional structure) has an elementary extension with a monomial group. This requires a
digression on abelian groups. Let $A$,~$B$ be (additively written) abelian
groups. Call $A$ a {\bf pure subgroup} of $B$ if $A$ is a subgroup of $B$
such that $A\cap nB = nA$ for all~$n\ge 1$. If $A$ is a subgroup
of $B$ and
$B/A$ is torsion-free, then $A$ is a pure subgroup of $B$. In case $B$ is itself torsion-free, then
$$\text{$A$ is a pure subgroup of $B$}\ \Longleftrightarrow\ \text{$B/A$ is torsion-free.}$$
Also, if $A$ is a direct summand of $B$ (that is, $A$ is a subgroup of $B$ and $B=A\oplus B'$, internally, for some subgroup $B'$ of $B$), then
$A$ is a pure subgroup of $B$. \index{subgroup!pure}  \index{pure subgroup}

\begin{lemma}\label{pu1} Suppose $A$ is a pure subgroup of $B$, and the group
$B/A$ is finitely generated. Then $B=A\oplus B'$, internally, for some subgroup $B'$ of
$B$.
\end{lemma}
\begin{proof} By the fundamental theorem on finitely generated abelian groups we have 
$$B/A\ =\ \Z(b_1+A) \oplus \dots \oplus \Z(b_m+A)$$ for
suitable $b_1,\dots,b_m\in B$. If $b_i+A$ has finite order $n_i\ge 1$ in $B/A$,
then $n_ib_i\in A$, so $n_ib_i=n_ia_i$ with $a_i\in A$, and replacing
$b_i$ by $b_i-a_i$ we arrange $n_ib_i=0$. With this adjustment
of the $b_i$s one checks easily that $B=A\oplus B'$ where
\equationqed{ B'\ =\ \Z b_1 + \dots + \Z b_m\ =\ \Z b_1 \oplus \dots \oplus \Z b_m.}
\end{proof}

\begin{cor}\label{pu2} Suppose $A$ is a pure subgroup of $B$,
$e_{ij}\in \Z$ for $1\le i\le m$ and $1\le j \le n$, and $a_1,\dots,a_m\in A$.
Suppose the system of equations
\begin{align*} e_{11}x_1 + \dots +e_{1n}x_n\ &=\ a_1\\
  \vdots \hskip3em           \vdots \hskip3em  \vdots    \ &\hskip2em \vdots\  \\
                  e_{m1}x_1 + \dots + e_{mn}x_n\ &=\ a_m
\end{align*}
has a solution in $B$, that is, there are $x_1,\dots, x_n\in B$ for which the above equations hold. Then this system has a solution in $A$.
\end{cor}

\begin{lemma}\label{pu3}  Suppose $A$ is a pure subgroup of $B$,
and $b\in B$. Then there is a pure subgroup $A'$ of $B$ that contains
$A$ and $b$ such that $A'/A$ is countable.
\end{lemma}
\begin{proof} By the Downward L\"{o}wenheim-Skolem Theorem (see \ref{prop:LS downwards}) we can take a
subgroup~$A'$ of $B$ that contains $A$ and $b$ such that $A'/A$ is countable
and $A'/A \preccurlyeq B/A$. It follows easily that $A'$ is a pure subgroup of $B$.
\end{proof}

\begin{prop}\label{pu4} Let $h\colon A \to U$ be a group morphism into an $\aleph_1$-saturated 
\textup{(}additive\textup{)} abelian group $U$ and suppose $A$ is a pure subgroup of $B$.
Then $h$ extends to a group morphism $B \to U$.
\end{prop}
\begin{proof} By the previous lemma we reduce to the case that 
$B/A$ is
countable. (This reduction does not use that $U$ is $\aleph_1$-saturated.) Let $b_0,b_1, b_2,\ldots$ generate~$B$ over~$A$. If we can find
elements $u_0,u_1, u_2,\ldots \in U$ such that
$\sum_{i=0}^n e_iu_i=h(a)$ whenever $\sum_{i=0}^n e_ib_i=a$ (all $e_i\in \Z$, $a\in A$),
then we can extend $h$ as desired
by sending $b_n$ to~$u_n$ for each $n$. Note that each
finite subset of this countable set of constraints on $(u_0,u_1,u_2,\dots)$
can be satisfied, by Corollary~\ref{pu2}. The desired result follows.
\end{proof}

\begin{cor}\label{pu5} Let $U$ be an $\aleph_1$-saturated \textup{(}additive\textup{)} abelian group and a pure subgroup of $B$. Then $U$ is a direct summand of $B$.
\end{cor}
\begin{proof} Apply Proposition~\ref{pu4} to the identity map $U\to U$.
\end{proof}

\noindent
For our valued field $K$, construed as a field with valuation ring, we get:

\begin{lemma}\label{xs1} If $K$ is $\aleph_1$-saturated, then it has
a cross-section.
\end{lemma}
\begin{proof} With $U=\mathcal O^\times$,
the natural inclusion $U \to K^\times$ and the map
$v\colon K^\times \to \Gamma$ yield the exact sequence of abelian groups
$$ 1 \to U\to  K^\times \to  \Gamma \to 0.$$
Since $\Gamma$ is torsion-free, $U$ is a pure subgroup
of $K^\times$. If $K$ is $\aleph_1$-saturated, then so is the group
$U$, and thus the above exact sequence splits by Corollary~\ref{pu5}.
\end{proof}

\noindent
The following variant of this lemma is also sometimes useful: 

\begin{lemma}\label{xs2}
Let $K$ be $\aleph_1$-saturated, and let
$E$ be a valued subfield of $K$ such that~$\Gamma_E$ is pure in $\Gamma$ and $\aleph_1$-saturated. Let $s_E$ be a cross-section of $E$.
Then $s_E$ extends to a cross-section of $K$.
\end{lemma}
\begin{proof} By Lemma~\ref{xs1} we have a cross-section $s$ of $K$.
By Corollary~\ref{pu5} we have an internal direct sum decomposition $\Gamma=\Gamma_E \oplus \Delta$ with $\Delta$ a subgroup of $\Gamma$. This gives a cross-section of $K$ that coincides with
$s_E$ on $\Gamma_E$ and with $s$ on $\Delta$.
\end{proof}

\subsection*{Uniqueness of maximal immediate extensions}
\textit{In this subsection we assume that $K$ is a valued field of equicharacteristic zero.}\/ By Corollaries~\ref{isomaximm} and \ref{cor:maximal valued fields} any two maximal immediate valued field extensions of $K$ are isomorphic over $K$ \textit{as valued abelian additive groups}\/; we can now improve on this:

\begin{cor}\label{cor:unique max imm ext, 1}
Any two maximal immediate valued field extensions of $K$ are isomorphic over $K$ as valued fields.
\end{cor}
\begin{proof}
Let $K_1$, $K_2$ be maximal immediate valued field extensions of $K$. Below, each subfield of $K_i$ ($i=1,2$) is viewed as a valued field by taking as valuation ring the intersection of the valuation ring of $K_i$ with the subfield. Consider fields~$L_1$,~$L_2$ with $K\subseteq L_i\subseteq K_i$ for $i=1,2$, and suppose that we have a valued field isomorphism $L_1\cong L_2$. Note that $L_1=K_1$ iff $L_2=K_2$.

Suppose $L_1\neq K_1$ (and hence $L_2\neq K_2$). It suffices to show that then we can extend the isomorphism $L_1\cong L_2$ to a valued field isomorphism $L_1'\cong L_2'$ where $L_i'$ is a field with $L_i\subseteq L_i'\subseteq K_i$ and $L_i'\neq L_i$, for $i=1,2$. If $L_1$ is not henselian, then we can take for $L_i'$ the henselization of $L_i$ in $K_i$. 
Thus, suppose $L_1$ is henselian. Take $b\in K_1\setminus L_1$ and take a divergent pc-sequence $(a_\rho)$ in $L_1$ such that $a_\rho\leadsto b$. Since $L_1$ is algebraically maximal by Corollary~\ref{cor:alg max equals henselian}, $(a_\rho)$ is of transcendental type over $L_1$. The image of $(a_\rho)$ under our isomorphism $L_1\cong L_2$ has a pseudolimit $c\in K_2$. By Lemma~\ref{lem:Kaplansky, 1} we can extend this isomorphism to an isomorphism $L_1(b)\cong L_2(c)$.
\end{proof}

\noindent
Assuming also that $K$ has a monomial group we can identify the maximal immediate extension of $K$ with a Hahn field:

\begin{cor}\label{cor:unique max imm ext, 2}
Suppose $K$ is maximal and $\mathfrak M$ is a monomial group of~$K$. Then there is a valued field isomorphism $K\cong \k[[\mathfrak M]]$ 
that induces the identity on $\k=\res(K)$ and on $\mathfrak M$.
\end{cor}
\begin{proof}
By Corollary~\ref{cor:alg max implies henselian}, $K$ is henselian. Thus by Proposition~\ref{prop:lift} there is a lift~$C$ of the residue field $\k$ in $\mathcal O$. Let $i$ be the field isomorphism~${c\mapsto\overline{c}\colon C\to\k}$. Then $i$ extends uniquely to a valued field embedding of the valued subfield $C(\mathfrak M)$ of $K$ in\-to~$\k[[\mathfrak M]]$ that is the identity on $\mathfrak M$. Clearly $K$ is an immediate extension of~$C(\mathfrak M)$. 
By the previous corollary, this embedding extends to a valued field iso\-mor\-phism~${K\cong \k[[\mathfrak M]]}$. 
\end{proof}

\noindent
The proof of Corollary~\ref{cor:unique max imm ext, 1} also yield a useful 
embedding property of the maximal immediate extension of $K$. To formulate this, set $\Gamma:= v(K^\times)$ and
let $L$ be a valued field extension of~$K$. We say that 
$L$ is {\bf $\Gamma$-maximal} if $L$ is henselian and every
pc-sequence in $L$ of length~$\le \operatorname{cf}(S)$ for some
$S\subseteq \Gamma$ pseudoconverges in $L$. In particular, if $L$ is maximal,
then $L$ is $\Gamma$-maximal. Also, if $L$ is henselian and
$|\Gamma|^{+}$-saturated, then~$L$ is $\Gamma$-maximal. \index{valued field!$\Gamma$-maximal} \index{maximal!$\Gamma$-maximal valued field} 

\begin{cor} If $M$ is a maximal immediate valued field extension of $K$,
then~$M$ embeds as a valued field over $K$ into any
$\Gamma$-maximal valued field extension of $K$. 
\end{cor}

\begin{figure}[h]
\begin{tikzpicture}
\begin{scope}[xshift=2em,on grid]
%\draw[help lines,step=2em,gray!20] (0,0) grid (20,10);
\node (Kmax) at (4,8) {$K^{\operatorname{max}}$};
\node[below left = 7em of Kmax] (Kc)  {};
\node[below right = 14em of Kc] (K)  {$K$};
\node[above right = 7em of K] (Kh) {$K^{\operatorname{h}}$};
\node[above right = 7em of Kh] (Khunr) {$(K^{\operatorname{h}})^{\operatorname{unr}}$};
\node[above right = 7em of Khunr] (Kalg) {$K^{\alg}$};
\draw (Kmax) -- (Kh) -- (K);
\draw (Kh) -- (Khunr) -- (Kalg);
\coordinate[above = of Kmax] (a);
\coordinate[left = of Kc] (b);
\coordinate[below = of K] (c) ;
\coordinate[right = of Kh] (d) ;
\draw[dotted] (a)--(b)--(c)--(d)--(a);
\coordinate[below = 1em of c] (cprime);
\coordinate[above left = 6em of cprime] (e);
\coordinate[above right = 26em of e] (f);
\coordinate[below right = 6em of f] (g);
\draw[dashed] (e)--(f)--(g)--(cprime)--(e);
\node[below = 9em of Kc] {immediate extensions};
\node[below = 9em of Khunr] {algebraic extensions};
\end{scope}
\end{tikzpicture}
\caption{Inclusion diagram of some valued field extensions of a valued field $K$ of equicharacteristic zero.} \label{fig:valued field exts}
\end{figure}
 
\noindent
Figure~\ref{fig:valued field exts} shows various valued field extensions of $K$ discussed in this section.
Here, $K^{\max}$ denotes a maximal immediate extension of $K$,
and $K^{\operatorname{h}}$ a henselization of $K$ in~$K^{\max}$. We introduce $(K^{\operatorname{h}})^{\operatorname{unr}}$ in the following subsection.

\subsection*{Unramified extensions}
\textit{In this subsection $K$ is henselian.}\/  We fix an algebraic closure 
$K^\alg$ of $K$, and
equip it with the unique valuation ring $\mathcal{O}^{\alg}$ making 
$K^{\alg}$ a valued field extension of $K$. Any algebraic
field extension of $K$ we mention is assumed to be a subfield of 
$K^\alg$, and regarded as a
{\em valued\/} subfield of $K^{\alg}$. 

Let $L$ be an algebraic field extension of $K$. If $[L:K]<\infty$, then we 
say that $L$ is {\bf unramified over $K$}
if the associated extension $\res(L)\supseteq\res(K)$ of residue fields is 
separable and $[L:K] = \big[\!\res(L):\res(K)\big]$
(and thus $\Gamma_L=\Gamma$, by Corollary~\ref{cor:ef inequ}). Suppose
$[L:K] < \infty$ and $E$ is an intermediate field: $K$ is a subfield of $E$ and $E$ is a subfield of $L$. Then by the 
multiplicativity of degrees and properties of separability, 
$$\text{$L$ is unramified over $K$}\ \Longleftrightarrow\ \text{$L$ is unramified over $E$ and $E$ is unramified over $K$.}$$ 
Without assuming $[L:K]< \infty$ we call $L$ {\bf unramified over $K$} 
if every intermediate field $E$ with $[E:K]<\infty$ is unramified over $K$
(and thus $\Gamma_L=\Gamma$). 

\index{extension!valued fields!unramified} 
\index{unramified extension of valued fields}

\begin{lemma} \label{lem:unram}
Let $L'\supseteq K$ be an algebraic field extension, and $L$, $K'$ subfields of~$L'$ containing $K$, such
that $L'=LK'$ and $L$ is unramified over $K$.
$$\xymatrix@R+1em@C+1em{
L \ar[r]      & L'  \\
K \ar[r]     \ar[u]  & K' \ar[u]
}$$
Then $L'$ is unramified over $K'$. 
\end{lemma}
\begin{proof}
We may assume $[L: K]<\infty$. Take $x\in\mathcal O_L$ such that $\res(L)=\res(K)(\overline{x})$, and let $P\in K[X]$ be the minimum polynomial of $x$ over $K$. Then $P\in\mathcal O[X]$ and
\begin{align*}
\big[\!\res(L):\res(K)\big]\, &\leq\, \deg\overline{P}\, =\, \deg P\\ &=\, \big[K(x):K\big]\, \leq\, [L:K]\ =\ \big[\!\res(L):\res(K)\big],
\end{align*}
hence $L=K(x)$ (and thus $L'=K'(x)$), and $\overline{P}$ is the minimum polynomial of~$\overline{x}$ over~$\res(K)$. Let $\mathcal{O}'$ be the valuation ring
of $K'$; so $\mathcal{O}'$ is the integral closure of $\mathcal{O}$ in $K'$.
Let $Q\in K'[X]$ be the minimum polynomial of $x$ over $K'$. Then
$Q\in \mathcal{O}'[X]$ and $P=QR$ with monic $R\in \mathcal{O}'[X]$, so 
$\overline{Q}$ divides $\overline{P}$ in $\res(K')[X]$, and hence 
$\overline{Q}\in \res(K')[X]$ is separable. Thus 
$\overline{Q}\in \res(K')[X]$ is irreducible by the irreducibility of 
$Q\in K'[X]$ and the lifting property 
of Proposition~\ref{prop:char henselian}(iv). Therefore,
\begin{align*}
\big[\!\res(L'):\res(K')\big]\ &\leq\ [L':K']\ =\ \deg Q\ =\ \deg\overline{Q}\\
   &=\ \big[\!\res(K')(\overline{x}):\res(K')\big]\ \leq\ \big[\res(L'):\res(K')\big].
\end{align*}
It follows that $\res(L')=\res(K')(\overline{x})$ and $[L':K']=\big[\res(L'):\res(K')\big]$, and so $L'$ is unramified over $K'$.
\end{proof}

\begin{cor} Let $L\supseteq K$ be an algebraic field extension which is unramified over~$K$. Let $E$ be an intermediate field. Then $L$ is unramified over $E$.
\end{cor}
\begin{proof} Lemma~\ref{lem:unram} applied to $K'=E$ shows that $L$ is unramified
over $E$.  
\end{proof}

\begin{cor} Suppose $L,L'\supseteq K$ are algebraic field extensions unramified over $K$.
Then the compositum $LL'$ is also unramified over $K$.
\end{cor}
\begin{proof} 
We can reduce to the case that $[L:K]<\infty$ and $[L':K]<\infty$.  It follows from Lemma~\ref{lem:unram}
that $LL'$ is unramified over $L'$. Then the equivalence just before Lemma~\ref{lem:unram} gives that
$LL'$ is unramified over $K$.
\end{proof}

\noindent
Let $L$ be an algebraic field extension of $K$. By the previous corollary,
the composite of the intermediate fields that are unramified over $K$ is itself unramified over
$K$, and is called the {\bf maximal unramified extension of $K$ in $L$.}  \index{maximal!unramified extension of a valued field} \index{valued field!maximal unramified extension} \nomenclature[K]{$K^{\operatorname{unr}}$}{maximal unramified extension of  $K$}

\begin{lemma}\label{lem:max unram}
Let $M$ be the maximal unramified extension of $K$ in $L$. Then $\res(M)$
is the separable algebraic closure of $\res(K)$ in $\res(L)$.
\end{lemma}
\begin{proof}
Let $\xi\in\res(L)$ be separable over $\res(K)$; we need to show $\xi\in\res(M)$. Let $P\in\mathcal O[X]$ be monic such that $\overline{P}$ is the minimum polynomial of $\xi$ over $\res(K)$. Then $P\in K[X]$ is irreducible. As $L$ is henselian, $P$ has a zero $x\in\mathcal O_L$ with $\overline{x}=\xi$. 
It follows that 
$$\big[K(x):K\big]\ =\ \big[\!\res(K)(\xi):\res(K)\big],$$ hence $K(x)$ is unramified over $K$, so $K(x)\subseteq M$, 
and thus $\xi=\overline{x}\in\res(M)$. 
\end{proof}

\noindent
We let $K^{\operatorname{unr}}$ be the maximal unramified extension of $K$ in $K^{\alg}$. By the previous lemma, the residue field of $K^{\operatorname{unr}}$ is the separable algebraic closure of $\res(K)$
in its algebraic closure $\res(K^{\alg})$.
Note that if $L\supseteq K$ is any algebraic field extension, then $K^{\operatorname{unr}}\cap L$ is the maximal unramified extension of $K$ in $L$.

\begin{prop}\label{prop:ramified}
Suppose $K$ has equicharacteristic zero and $L\supseteq K$ is a field extension with $[L:K]<\infty$
and $\res(K)=\res(L)$. Then there is a chain of subfields
$$K=K_0\ \subseteq\  K_1\ \subseteq\cdots\subseteq\ K_{n}= L$$
of $L$ such that for $i=1,\dots,n$ we have $K_{i}=K_{i-1}\big(b^{1/p}\big)$ for some
prime number $p$ and some $b\in K_{i-1}^\times$ with
$vb\notin p\Gamma_{K_{i-1}}$. Moreover, $[L:K]\ =\ [\Gamma_L:\Gamma]$.
\end{prop}
\begin{proof}
The quotient group $\Gamma_L/\Gamma$ is finite, so we have a chain of subgroups
$$\Gamma=\Gamma_0\ \subseteq\ \Gamma_{1}\ \subseteq \cdots \subseteq\ \Gamma_{n}=\Gamma_L$$
of $\Gamma_L$ such that $[\Gamma_{i}:\Gamma_{i-1}]$ is prime for $i=1,\dots,n$.
Set $K_0=K$, let $i\in\{1,\dots,n\}$ and suppose $K_{i-1}$ is an intermediate field
with $\Gamma_{K_{i-1}}=\Gamma_{i-1}$. With $[\Gamma_{i}:\Gamma_{i-1}]=p$, take
$\alpha\in\Gamma_i$ such that $\alpha+\Gamma_{i-1}$ has order $p$ in $\Gamma_i/\Gamma_{i-1}$.
Then $\alpha=va$ where $a\in L^\times$, and $p\alpha=vb$ with $b\in K_{i-1}^\times$.
We arrange $a^{p}=b$ as follows:
since $\res(L)=\res(K)$ we have $a^{p}/b=cu$ where $c\in K$ and $u\in L$ with $c\asymp 1$
and $u\sim 1$. As $L$ is henselian, the polynomial $X^{p}-u\in\mathcal O_L[X]$ has a zero $x$ in $\mathcal O_L$ with $x\sim 1$; replacing~$a$ by~$a/x$ and~$b$ by $bc$ we have $a^p=b$. With $K_i:=K_{i-1}(a)$,
this yields $\Gamma_{K_i}=\Gamma_i$ by Lemma~\ref{lem:pth root}. In this way we construct a chain $K_0\subseteq K_1\subseteq \cdots \subseteq K_n$ of subfields of $L$.  Then the extension $L\supseteq K_n$ is immediate, so $L=K_n$ by Corollary~\ref{cor:alg max equals henselian}. Lemma~\ref{lem:pth root}
also gives $[L:K]\ =\ [\Gamma_L:\Gamma]$ by multiplicativity.
\end{proof}

\noindent
Proposition~\ref{prop:ramified}  together with Lemma~\ref{lem:max unram} 
yield a key fact:

\begin{cor}\label{cor:ef equ} Suppose $K$ has equicharacteristic zero and $L\supseteq K$ is a field extension with $[L:K]<\infty$. Then
$$%\begin{equation}\label{eq:ef equ}
[L:K]\  =\  [\Gamma_L:\Gamma]\cdot \big[\!\res(L):\res(K)\big].
$$%\end{equation}
\end{cor}

\begin{proof} Let $M$ be the maximal unramified extension of $K$ in $L$. Then by definition 
$\big[\!\res(M):\res(K)\big]=[M:K]$. Next, $\res(M)=\res(L)$ by Lemma~\ref{lem:max unram}, so $[L:M]=[\Gamma_L:\Gamma_M]$ by Proposition~\ref{prop:ramified}. It remains to note that $\Gamma_M=\Gamma$. 
\end{proof}

\index{extension!valued fields!purely ramified} 
\index{purely ramified}

\noindent
Let $L\supseteq K$ be a  field extension with $[L:K]<\infty$. 
One says that $L$ is  {\bf purely ramified over $K$} if $[L:K]=[\Gamma_L:\Gamma]$.
Thus by Corollaries~\ref{cor:ef inequ} and~\ref{cor:ef equ}, 
if  $K$ has equicharacteristic zero, then 
\begin{align*}
\text{$L$ is unramified over $K$} &\ \Longleftrightarrow\  \Gamma_L=\Gamma, \\
  \text{$L$ is purely ramified over $K$} &\ \Longleftrightarrow\  \res(L)=\res(K).
\end{align*}

\subsection*{Notes and comments} Hensel~\cite{Hensel} introduced the valued field of $p$-adic numbers and showed it to be hen\-se\-lian. 
Its generalization stated in
Corollary~\ref{cor:henselian complete rank 1}  is due to Rychl{\'\i}k~\cite{Rychlik}.
The characterizations of henselianity collected in Lemma~\ref{lem:char henselian} and  
Proposition~\ref{prop:char henselian} are due to Nagata~\cite{Nagata},  Rayner~\cite{Rayner}, and Rim~\cite{Rim}. 
The notion of henselization and the existence of the henselization of a 
valued field 
is also due to Nagata. (The henselization of rank~$1$ valued fields 
had been constructed already by Ostrowski~\cite{Ostrowski}.) 
%The proof of Proposition~\ref{prop:hensel, multivar} given above is from \cite{FVK-hensel} and credited there to F.~Pop. 

Puiseux series (see Example~\ref{ex:Puiseux}) stem from Newton~\cite{Newton} and Pui\-seux~\cite{Puiseux50}. 
It is well-known that for a field $C$ of characteristic $p>0$, the Puiseux series field~$\operatorname{P}(C)$ is not algebraically closed in its extension $C \(( t^{\Q} \)) $: 
the polynomial $X^p-X-x$ with $x=t^{-1}$ has a zero ${x^{1/p}+x^{1/p^2}+\cdots}\in C \(( t^{\Q} \)) \setminus \operatorname{P}(C)$. (This example is from Chevalley~\cite[p.~64]{Chevalley51}. See~\cite{Kedlaya,Kedlaya15} for more on this issue.) 
%For a description of the algebraic closure of 
%$C[[x^\Z]]$ in~$C[[x^\Q]]$ where
%$C$ is an algebraically closed field of positive %characteristic see \cite{Kedlaya}.

For Lemma~\ref{lem:min poly immediate ext}, see \cite[Proposition~1.1]{Scanlon}. 
The facts on abelian groups in \ref{pu1}--\ref{pu5} are taken from \cite[Chapter~V, \S{}5]{Cherlin} where this material is treated
for modules over any ring. Lemma~\ref{xs1} is also in \cite{Cherlin}. 
Cross-sections are from~\cite{AxKochen}. The notion of pure subgroup and basic facts about it
are due to Pr\"ufer~\cite{Pruefer}.
Proposition~\ref{prop:lift} goes back to~\cite{Hasse-Schmidt} for discrete valuations, and for general valuations to~\cite{Kaplansky,MacLane}.
Corollaries~\ref{cor:unique max imm ext, 1} and \ref{cor:unique max imm ext, 2} are due to Kaplansky~\cite{Kaplansky}. The maximal unramified extension plays a key role in the Galois theory of valuations, of which we only needed a few facts here; see~\cite{Prestel-Engler, FVK} for more. 

Our use of ``henselian'' for valued fields is common in valuation theory, but clashes a bit with its use in the larger area of commutative algebra: there it would be the {\em valuation ring\/} of the valued field that is henselian, according to
Azumaya~\cite{Azumaya} who defined a local ring $\mathcal O$ with residue field $\k$ to be henselian if it satisfies
the condition in Definition~\ref{def:henselian}. 
Some results in this section generalize to this local ring setting, for example Proposition~\ref{prop:hensel, multivar}, but in that case the proof is more difficult: \cite[(18.5.11), (b)]{EGAIV}; see also \cite[(1.9)]{ArtinApprox}, \cite[Theorem~4.2,~(d$'$)]{Milne}.

%% file: mt-3-4.tex
\section{Decomposing Valuations}\label{sec:decomposition}

\noindent
A major theme in our work is constructing zeros of differential polynomials over
valued differential fields. We
approximate/construct these zeros by pc-sequences that have special properties. In this section we introduce these 
properties in a purely valuation-theoretic setting. A basic tool for analyzing such  pc-sequences is the two-pronged process of coarsening and specializing a valuation, which we treat first. {\em Throughout this section $K$ is a valued field with valuation ring $\mathcal{O}$, 
valuation $v\colon K^\times\to\Gamma$  and  dominance relation $\preceq$}. 

\index{coarsening!valuation on a field}
\index{Delta-coarsening@$\Delta$-coarsening}
\nomenclature[K]{$\dot v=v_\Delta$}{coarsening of  $v$ by $\Delta$}
\nomenclature[Rx]{$\preceq_\Delta$ or $\dot{\preceq}$}{dominance relation associated to $\dot v=v_\Delta$}

\subsection*{Coarsening and specialization} 
Let $\Delta$ be a convex subgroup of
$\Gamma=v(K^\times)$. Then we have the ordered quotient group~$\dot{\Gamma}=\Gamma/\Delta$
and the valuation 
$$\dot{v}=v_\Delta \colon
K^\times\to\dot{\Gamma}$$ on the field $K$, defined by 
$\dot{v}(f):=v(f) + \Delta$. We call $\dot{v}$ 
the {\bf coarsening of $v$ by~$\Delta$}, or the 
{\bf $\Delta$-coarsening of $v$}; the valued field 
$(K, \dot{v})$ whose valuation is
$\dot{v}$ instead of $v$ is called the {\bf $\Delta$-coarsening of $K$}. 
Note that $\dot v$ is indeed a coarsening of $v$ as defined in 
Section~\ref{sec:valued abelian gps}. Thus if 
$\Delta\neq\Gamma$, then the $v$-topology agrees with the 
$\dot v$-topology.
The dominance relation on $K$ corresponding to the coarsened
valuation~$\dot{v}$ is denoted by~$\preceq_\Delta$, or by $\dot{\preceq}$ if $\Delta$ is understood from the context, so 
$f\preceq g \Rightarrow f\ \dot{\preceq}\ g$ for $f,g\in K$.
The valuation ring of $\dot{v}$ is
$$\dot{\mathcal{O}}\ :=\ \{a\in K: \text{$va\ge \delta$  for some  $\delta\in \Delta$} \},$$
which has $\mathcal{O}$ as a subring and has maximal ideal \nomenclature[K]{$\dot{\mathcal O}$}{valuation ring of the coarsening $\dot v=v_\Delta$ of a valuation~$v$} \nomenclature[K]{$\dot{\smallo}$}{maximal ideal of $\dot{\mathcal O}$}
$$\dot{\smallo}\ :=\ \{a\in K: va>\Delta\}\subseteq \smallo.$$ 
Put $\dot K:=\dot{\mathcal{O}}/\dot{\smallo}$, 
the residue field of $\dot{\mathcal{O}}$, and put 
$\dot{a}:= a + \dot{\smallo} \in \dot{K}$ for 
$a\in \dot{\mathcal{O}}$. 
Then for $a\in \dot{\mathcal{O}}$ with $a\notin \dot{\smallo}$ 
the value 
$va$ depends only on the residue class $\dot a\in \dot K$, which gives
the valuation 
\begin{equation}\label{eq:valuation on dot K}
v\colon \dot K^{\times} \to \Delta, \qquad v\dot a:= va
\end{equation}
with valuation ring $\mathcal{O}_{\dot K}\ =\ \{\dot{a}: a\in \mathcal{O}\}$. The maximal ideal of
$\mathcal{O}_{\dot K}$ is $$\smallo_{\dot K}\ =\ \{\dot{a}:a\in \smallo\}.$$ Throughout 
$\dot{K}$ stands for the valued field $(\dot{K},\mathcal{O}_{\dot K})$.
The composed map 
$$\mathcal{O} \to \mathcal{O}_{\dot K} \to \res(\dot K)=
\mathcal{O}_{\dot K}/\smallo_{\dot K}$$
has kernel $\smallo$, and thus 
induces a field isomorphism 
$\res(K) \overset{\cong}{\longrightarrow} \res(\dot K)$, and
we identify~$\res(K)$ and $\res(\dot K)$ via this map.
We call $\dot K$ a {\bf specialization} of the
valued field $K$. The following is now straightforward.

\nomenclature[K]{$\dot{K}$}{valued residue field of the coarsening $\dot v=v_\Delta$ of $v$}
\index{valued field!specialization}
\index{specialization}

\begin{lemma}\label{copc} Let $(a_{\rho})$ be a well-indexed sequence in $\dot{\mathcal{O}}$
such that $(\dot{a}_{\rho})$ is a pc-sequence in $\dot{K}$. 
Then $(a_{\rho})$ is a pc-sequence in $K$, and for all $a\in \dot{\mathcal{O}}$,
$$ a_{\rho} \leadsto a\quad \text{in $K$}\ \Longleftrightarrow\ \dot{a}_{\rho}\leadsto \dot{a}\quad \text{in $\dot{K}$.}$$
\end{lemma}

%\begin{remark}
%The process of passing from the valued field $K$ to the valued fields 
%$(K,\dot{\mathcal O})$ and $\dot K$ can be reversed, but this will 
%not be needed.
%\end{remark}

\noindent
We picture the valuation $v$ on $K$ as follows:

\begin{center}
\begin{tikzpicture}
  \matrix (m) [matrix of math nodes, row sep=3em,
    column sep=3em]{
    		& K 			& 						\\
\Gamma	& 			& \operatorname{res}(K)	\\ };
  \draw[-stealth]
    (m-1-2) edge[thick] (m-2-1) edge[thick] (m-2-3);
\end{tikzpicture}
\end{center}

\noindent
The greater the angle between the arrows, the larger the valuation ring. Given our convex subgroup $\Delta$ of $\Gamma$, this gives rise to the following diagram which displays the original valuation $v$, its coarsening by $\Delta$, and the corresponding specialization: 

\begin{center}
\begin{tikzpicture}
  \matrix (m) [matrix of math nodes, row sep=3em,
    column sep=3em]{
    			& K 			& 						&  \\
\dot\Gamma	& 			& \dot K					&  \\
\Gamma		& 			& \operatorname{res}(K)	& \operatorname{res}(\dot K) \\
\Delta		&			&						& \\ };
  \draw[-stealth,>=]
    (m-1-2) edge[densely dotted,very thick] (m-2-1) edge[densely dotted, very thick] (m-2-3) 
    (m-3-3) edge[stealth-stealth] node[above] {$\cong$} (m-3-4);
  \draw[->>,=stealth]  	(m-3-1) edge (m-2-1);  
  \draw[-stealth] (m-4-1) edge node[left]{$\subseteq$} (m-3-1);
   \draw[->]
    (m-2-3) edge[dotted,thick] (m-3-4) edge[dotted,thick] (m-4-1);
  \draw[-stealth]
    (m-1-2) edge[thick] (m-3-1) edge [-,line width=6pt,draw=white] (m-3-3) edge[thick] (m-3-3);
\end{tikzpicture}
\end{center}

\noindent
The original valuation $v\colon K^\times\to\Gamma$ can often be understood 
in terms of the two (usually simpler)
valuations $\dot v\colon K^\times\to\dot{\Gamma}$ and 
$v\colon \dot K^\times \to \Delta$.
Here is an example:

\begin{lemma}\label{lem:henselian decomposition}
$K$ is henselian $\ \Longleftrightarrow\ (K,\dot{\mathcal O})$ and $\dot K$ are henselian.
\end{lemma}
\begin{proof}
We extend the natural surjection $x\mapsto\dot x\colon \dot{\mathcal O}\to \dot{\mathcal O}/\dot\smallo=\dot K$ to
the ring morphism $Q\mapsto \dot{Q}\colon \dot{\mathcal O}[X]\to \dot K[X]$ which maps the indeterminate $X$ to itself.

Suppose first that $K$ is henselian. Then $(K,\dot{\mathcal O})$ clearly is  henselian, using $\dot{\smallo}\subseteq\smallo$, ${\mathcal O}\subseteq{\dot{\mathcal O}}$, and the equivalence of (i) and (ii) in Lemma~\ref{lem:char henselian}. 
To show that $\dot{K}$ is also henselian, we again use condition (ii) for henselianity from Lemma~\ref{lem:char henselian}. Thus, let $P=1+X+P_2X^2+\cdots+P_nX^n\in {\mathcal O}_{\dot K}[X]$ where $n\geq 2$ and $P_i\in\smallo_{\dot K}$ for $i=2,\dots,n$; we have to show that $P$ has a zero in ${\mathcal O}_{\dot K}$. Take $Q=1+X+Q_2X^2+\cdots+Q_nX^n\in\dot{\mathcal O}[X]$ with $\dot Q=P$. Then $Q_i\in\smallo$ for $i=2,\dots,n$, hence
$Q$ has a zero $x\in\mathcal O$, and then $\dot x\in {\mathcal O}_{\dot K}$ is a zero of~$P$.

Conversely, suppose $(K,\dot{\mathcal O})$ and $\dot K$ are 
henselian.
Let $P\in\mathcal O[X]$ and $a\in\mathcal O$ be such that $P(a)\in \smallo$ and 
$P'(a)\in \mathcal O^{\times}$; we need to show that $P$ has a zero $b\in\mathcal O$ with $b\equiv a\bmod\smallo$.
Since $\dot K$ is henselian, we can take $c\in\mathcal O$ (so $\dot c\in\mathcal O_{\dot K}$) such that 
$\dot P(\dot c)=0$ and $\dot c\equiv\dot a\bmod \smallo_{\dot{K}}$. Then
$P(c)\in \dot{\smallo}$ and
$c\equiv a\bmod \smallo$, so 
$P'(c)\in \mathcal O^{\times}\subseteq \dot{\mathcal O}^{\times}$. 
Since~$(K,\dot{\mathcal O})$ is henselian, this gives $b\in\dot{\mathcal O}$ with $P(b)=0$ and $b\equiv c\bmod\dot\smallo$. 
Then $b$ is a zero of $P$ as required.
\end{proof}

\subsection*{Coarsening and valued field extensions}
Let $\Delta$ be a convex subgroup of
$\Gamma$ with ordered quotient group $\dot{\Gamma}=\Gamma/\Delta$. As described above, this gives rise to the coarsened valuation
$\dot{v}=v_\Delta\colon K^{\times} \to \dot\Gamma$ with residue field
$\dot{K}=\dot{\cal{O}}/\dot{\smallo}$, the latter equipped with  the valuation  $v\colon\dot K^{\times}\to \Delta$ such that $v\dot a=va$ for all 
$a\in\dot{\mathcal O}\setminus \dot{\smallo}$. 

\medskip
\noindent
Next, let $L$ be a valued field extension of $K$, and let $\Delta_L$ be the convex hull of $\Delta$ in~$\Gamma_L$. 
Then the natural inclusion $\Gamma\to\Gamma_L$ induces an embedding 
$$\dot\Gamma=\Gamma/\Delta\to\dot\Gamma_L:=\Gamma_L/\Delta_L$$ of ordered abelian groups, and identifying $\dot\Gamma$ with its image under this embedding,
the coarsening $\dot v=v_{\Delta_L}$ of the valuation $v$ of $L$ by $\Delta_L$ extends the coarsening $\dot v=v_\Delta$ of the valuation of $K$ by $\Delta$, so we have a valued field extension 
$(L, \dot{\mathcal{O}}_L)$ of $(K, \dot{\mathcal{O}})$.  
As usual, we identify $\res(K)$ with a subfield of $\res(L)$ and
the residue field $\dot K$ of $v_\Delta$ with a subfield of the residue field $\dot L$ of $v_{\Delta_L}$. Then the valuation $v\colon \dot L^{\times}\to \Delta_L$ with $v\dot a=va$ for all $a\in\dot{\mathcal O}_L\setminus \dot{\smallo}_L$ restricts to the valuation  $v\colon\dot K^{\times}\to \Delta$ on $\dot K$. Thus the diagrams
$$\xymatrix@R-0.5em@C-0.3em{
L^\times \ar[r]^{v_{\Delta_L}}					& \dot\Gamma_L			 & \dot{L}^\times \ar[r]^v					& \Delta_L	 & \res(L) \ar[r]^\cong      				& \res(\dot L) \\
K^\times \ar[r]^{v_\Delta} \ar[u] 	& \dot\Gamma \ar[u]	 &  \dot{K}^\times \ar[r]^v\ar[u]	& \Delta \ar[u]	 & \res(K) \ar[r]^\cong \ar[u]		& \res(\dot K) \ar[u]
}$$
commute. (Here the vertical arrows are the natural inclusions.)

\begin{lemma}\label{lem:coarsening, finite ext}
Suppose $[L:K]<\infty$. Then $[\dot L:\dot K]\leq [L:K]<\infty$. Thus if 
$\dot K$ is spherically complete \textup{(}respectively, complete\textup{)}, then so is $\dot L$.
\end{lemma}
\begin{proof}
The first statement follows from Corollary~\ref{cor:ef inequ, 2}. 
The second statement follows from the first and 
Corollary~\ref{cor:finite ext of complete}.
\end{proof}

\noindent
If $\Gamma_L=\Gamma$, then $\Delta_L=\Delta$, and if $\res(L)=\res(K)$, then 
$\res(\dot L)=\res(\dot K)$.
In particular, if $L$ (with valuation ring $\mathcal O_L$) is an 
immediate extension of $K$, then 
$\dot L$ is an immediate extension of 
$\dot K$. 

\begin{lemma}\label{coarseningcompletion} Suppose $L$ is a completion
of the valued field $K$ and $\Delta\ne \Gamma$. 
Then the coarsening $(L, \dot{\mathcal{O}}_L)$ is a completion of
$(K, \dot{\mathcal{O}})$.
\end{lemma}
\begin{proof} Note: $L$ is an immediate extension of $K$. Lemma~\ref{coarsecomplete} gives completeness of 
$(L, \dot{\mathcal{O}}_L)$. Use $\Delta\ne \Gamma$ to get $K$ dense in $L$ with respect to $\dot{v}$. 
\end{proof}

\noindent
Conversely, suppose $(L,\dot{\mathcal O}_L)$ is an {\em immediate\/} 
valued field extension of $(K, \dot{\mathcal O})$, and let~$\dot v$ denote
both the valuation of $(K, \dot{\mathcal O})$ and of $(L,\dot{\mathcal O}_L)$.
Thus  $\dot{K}=\dot{L}$ after the usual identification, where $\dot L$ is
the residue field of $(L,\dot{\mathcal O}_L)$. 

We define a map 
$v\colon L^\times \to \Gamma$ extending the valuation 
$v\colon K^\times \to \Gamma$ 
as follows. 
For $f\in L^\times$, take $g\in K^\times$ and $u\in L^\times$ such that
$f=gu$ and $\dot{v}(u)=0$; then $\dot{u}\in
\dot{L}^\times=\dot{K}^\times$, so $v(\dot{u})\in \Delta$; it is 
easy to check that $v(g)+v(\dot{u})\in \Gamma$ depends only on~$f$ and not
on the choice of $g$, $u$; now put $v(f):= v(g)+v(\dot{u})$. 

\begin{lemma}\label{immcoarse} The map $v\colon L^\times \to \Gamma$ 
is a valuation on $L$, 
its coarsening by $\Delta$ 
is exactly the given valuation $\dot{v}$ on $L$, and $L$ with 
this valuation $v$ is an immediate extension of the valued field $K$.
\end{lemma}
\begin{proof} It is easy to check that 
$v\colon L^\times \to \Gamma$ is a group morphism extending $v$ on~$K^\times$, 
and
that $\dot{v}(f)=v(f)+ \Delta\in \dot{\Gamma}$ for $f\in L^\times$. 
Also, if $f\in L^\times$ and $\dot{v}f>0$, then $vf>0$ and $v(1+f)=0$. 
Next, for $f_1,f_2\in L^\times$ with $f_1+f_2\ne 0$ one shows that 
$v(f_1+f_2) \ge \min\big\{vf_1,vf_2\big\}$ by distinguishing the cases
$\dot{v}f_1=\dot{v}f_2$ and $\dot{v}f_1 < \dot{v}f_2$. The rest 
of the lemma follows easily.
\end{proof}

\begin{cor}\label{cor:max decomp}
$K$ is maximal $\ \Longleftrightarrow\  (K,\dot{\mathcal O})$ and $\dot K$ are maximal.
\end{cor}
\begin{proof} Assume $K$ is maximal. Then $(K,\dot{\mathcal O})$ is maximal by 
 Corollary~\ref{coco}, and~$\dot K$ is maximal by 
Lemma~\ref{copc}. Conversely, assume $(K,\dot{\mathcal O})$ and $\dot K$ are maximal, and let~$L$ be an immediate valued field
extension of $K$.
Then with the notations above,~$\dot L$ is an
immediate extension of $\dot K$, so $\dot L = \dot K$. 
It follows that $(L,\dot{\mathcal O}_L)$ is an immediate extension
of $(K,\dot{\mathcal O})$, and so $K=L$. 
\end{proof}

\index{pc-sequence!special}
\index{pc-sequence!$\Delta$-special}
\index{special!pc-sequence}
\index{Delta-special@$\Delta$-special pc-sequence}

\subsection*{Step-completeness}
Let $(a_\rho)$ be a pc-sequence in $K$, and set 
$$s(\rho)\ :=\ \text{immediate successor of $\rho$.}$$ Then 
$\gamma_\rho := v(a_{s(\rho)}-a_\rho)\in \Gamma$ for all sufficiently large $\rho$, and
$(\gamma_\rho)$ is eventually strictly increasing. Also, if 
$a_\rho\leadsto a\in L$, where $L$ is a valued field
extension of $K$, then $\gamma_\rho=v(a-a_\rho)$ for all sufficiently 
large $\rho$. Given a nontrivial convex subgroup~$\Delta$ of~$\Gamma$ we say that $(a_\rho)$ is {\bf $\Delta$-special} if  
$(\gamma_\rho)$ is eventually cofinal in $\Delta$, that is, $\gamma_\rho\in \Delta$ for
all sufficiently large $\rho$, and for each $\delta\in \Delta$ we have
$\gamma_\rho > \delta$ for all sufficiently large $\rho$. 
Thus $(a_\rho)$ is $\Delta$-special iff the width of $(a_\rho)$ is $\{\gamma\in\Gamma_\infty:\gamma>\Delta\}$; in particular,
if~$(a_\rho)$ is $\Delta$-special, then so is every pc-sequence in $K$ equivalent to $(a_\rho)$.
We call~$(a_\rho)$ {\bf special} if it is $\Delta$-special for some
nontrivial convex subgroup $\Delta$ of $\Gamma$. 
We say that~$K$ is {\bf step-complete} if each special 
pc-sequence in $K$ has a pseudolimit in~$K$. By 
Lemma~\ref{lem:c-sequence width}, this is the same as
requiring that for each convex subgroup~$\Delta\ne \{0\}$ of $\Gamma$
the residue field~$\dot K$ of the coarsened valuation $\dot v=v_\Delta$ 
is complete with respect to the valuation $v$ defined in 
\eqref{eq:valuation on dot K}. Therefore: 

\index{valued field!step-complete}
\index{step-complete valued field}
\label{p:step-complete}

\begin{cor}\label{lem:step-complete decomp}
Let $\Delta$ be a convex subgroup of $\Gamma$, with ordered quotient group 
$\dot\Gamma=\Gamma/\Delta$,
corresponding coarsening $\dot v=v_\Delta\colon K^{\times}\to\dot\Gamma$ of $v$, 
and residue field $\dot K$ of~$\dot v$ with valuation 
$v\colon\dot K^{\times}\to\Delta$ as above. Then 
$$\text{$K$  is step-complete }\ \Longleftrightarrow\ 
\text{$(K,\dot{\mathcal O})$ and $\dot K$  are step-complete.}$$
\end{cor}

\noindent
We have the following implications for our valued field $K$:
$$ \text{spherically complete} \quad\Longrightarrow\quad 
\text{step-complete} \quad\Longrightarrow\quad  \text{complete.}
$$
If $K$ has rank~$1$, then the second arrow may be reversed.
If $K$ is discrete, then both arrows may be reversed. Thus by Corollaries~\ref{cor:max decomp} and~\ref{lem:step-complete decomp} and induction on $n$:

\begin{lemma}\label{lem:McL}
If $K$ is step-complete and $\Gamma$ is isomorphic to the lexicographically 
ordered group $\Z^n$, then $K$ is spherically complete. 
\end{lemma}

\noindent
We note a few more basic facts about step-complete valued fields:

\begin{cor}\label{lem:step-complete is henselian}
If $K$ is step-complete, then $K$ is henselian.
\end{cor} 
\begin{proof}
Let $P\in\mathcal O[X]$ and $a\in\mathcal O$ be such that $P(a)\prec 1$, $P'(a)\asymp 1$ and $P$ has no zero in $a+\smallo$. Then by the discussion preceding Corollary~\ref{cor:henselian complete rank 1}, there is a special divergent
pc-sequence in $K$, so $K$ is not step-complete.
\end{proof}

\begin{cor}\label{lem:finite ext of step-complete} Suppose $K$ is
step-complete and $L$ is a valued field extension of~$K$ with $[L:K]<\infty$.
Then $L$ is step-complete.
\end{cor}
\begin{proof} By Corollary~\ref{cor:ef inequ, 2},
every convex subgroup of $\Gamma_L$ is the convex hull in 
$\Gamma_L$ 
of a convex subgroup of $\Gamma$. 
Now apply Lemma~\ref{lem:coarsening, finite ext}.
\end{proof}

\noindent
Every valued field has a step-complete valued field extension; in the case of equi-characteristic zero there is such an extension with a semiuniversal property:

\begin{prop}\label{prop:step-completion}
Suppose $K$ has equicharacteristic zero. Then
there exists a step-complete valued field extension $K^{\operatorname{sc}}$ of 
$K$ which embeds over $K$ into any step-complete valued field extension of $K$.
\end{prop}
\begin{proof}
Fix a maximal immediate extension $M$ of $K$; for a valued subfield $L$ of~$M$ we denote by $L^{\operatorname{h}}$ the henselization of $L$ inside $M$. Starting with $K_0=K^{\operatorname{h}}$, we construct an increasing sequence $(K_\rho)$ of henselian valued subfields of $M$, indexed by
ordinals~$\rho$,  as follows. Let an ordinal $\mu>0$ and
an increasing sequence $(K_\rho)_{\rho< \mu}$ of 
henselian valued subfields of $M$ be given. If $\mu$ is a limit ordinal, we set $K_\mu:=\bigcup_{\rho<\mu} K_\rho$. Suppose $\mu$ is a successor ordinal, say $\mu=\lambda+1$. If $K_\lambda$ is step-complete, then we are done and set 
$K^{\operatorname{sc}}:=K_\lambda$. So assume $K_\lambda$ is not step-complete. Then we 
take a convex subgroup $\Delta\neq\{0\}$ of $\Gamma$ and a divergent $\Delta$-special pc-sequence $(a_\rho)$ in~$K_\lambda$;
note that~$(a_\rho)$ is of transcendental type, by 
%Lemma~\ref{lem:Kaplansky, 1} and 
Corollary~\ref{cor:alg max equals henselian}.
Take $a\in M$ with $a_\rho\leadsto a$ and set $K_\mu:=K_\lambda(a)^{\operatorname{h}}$. This construction has to terminate eventually. By Lemmas~\ref{lem:Kaplansky, 1} and \ref{lem:step-complete is henselian} and the universal property of henselization, $K^{\operatorname{sc}}$ has the required property.
\end{proof}

\noindent
Suppose $K$ has equicharacteristic zero. Call a valued field extension $K^{\operatorname{sc}}$ of~$K$ as in 
Proposition~\ref{prop:step-completion} a {\bf step-completion} of $K$. Each step-completion of $K$ 
is an immediate extension of $K$, since every maximal immediate extension of $K$ is step-complete. 
Does $K$ have up to isomorphism of valued fields
over $K$ a unique step-completion? (The answer is ``yes'' if $K$ has rank~$1$, and also if 
$\Gamma=\Z^n$, ordered lexicographically, by 
Corollary~\ref{cor:unique max imm ext, 1} and Lemma~\ref{lem:McL}.)

\index{valued field!step-completion}
\index{step-completion of a valued field}
\nomenclature[K]{$K^{\operatorname{sc}}$}{a step-completion of the valued field $K$}

\subsection*{Fluent and jammed pc-sequences}
Let $(a_\rho)$ be a pc-sequence in $K$. Set
$$\gamma_\rho\ :=\ v(a_{s(\rho)}-a_\rho) \quad \text{ where }s(\rho):= \text{immediate successor of $\rho$.}$$ 
So $\gamma_\rho\in \Gamma$ for all sufficiently large $\rho$, and $(\gamma_\rho)$ is eventually strictly increasing. 

\medskip\noindent
Let $\Delta$ be a convex subgroup of $\Gamma$. We say that $(a_\rho)$ is {\bf $\Delta$-fluent} \index{Delta-fluent@$\Delta$-fluent pc-sequence}\index{pc-sequence!$\Delta$-fluent} if $(a_\rho)$ remains a pc-sequence for the coarsened valuation $v_{\Delta}$; equivalently, for some index $\rho_0$,
$$\gamma_{\rho'} - \gamma_\rho > \Delta\quad\text{for all $\rho'> \rho > \rho_0$.}$$
The following is almost immediate but good to keep in mind.

\begin{lemma}\label{pcDeltainv} Suppose $(a_{\rho})$ is $\Delta$-fluent and $a\in K$. Then
$$a_{\rho} \leadsto a \text{ in $K$}\ \Longleftrightarrow\ a_{\rho} \leadsto a \text{ in the $\Delta$-coarsening of $K$}.$$
\end{lemma}

\noindent
We say that $(a_\rho)$ is {\bf $\Delta$-jammed} \index{pc-sequence!$\Delta$-jammed}\index{Delta-jammed@$\Delta$-jammed pc-sequence} if for some index $\rho_0$,
$$\gamma_{\rho'}-\gamma_\rho \in \Delta\quad \text{ for all $\rho'>\rho > \rho_0$.}$$
It is easily checked that if $(a_\rho)$ is not $\Delta$-jammed, 
then it has a cofinal $\Delta$-fluent subsequence. Let 
$W\subseteq \Gamma_{\infty}$ be the width of $(a_{\rho})$. Then $(a_{\rho})$ is
$\Delta$-jammed iff for some $\gamma\in \Gamma$ with $\gamma <  W$
we have $\gamma'-\gamma\in \Delta$ for all $\gamma'$ with
$\gamma < \gamma' < W$. Thus being $\Delta$-jammed depends only on the width. 
Therefore, if $(a_{\rho})$ is $\Delta$-jammed, then so is
every equivalent pc-sequence in $K$.

\medskip\noindent
We say that $(a_\rho)$ is {\bf fluent} \index{pc-sequence!fluent}\index{fluent!pc-sequence} if it is $\Delta$-fluent for some 
nontrivial convex subgroup $\Delta$ of $\Gamma$, and we say that $(a_\rho)$ is 
{\bf jammed} \index{pc-sequence!jammed}\index{jammed pc-sequence} if it is $\Delta$-jammed for every nontrivial convex 
subgroup~$\Delta$ of $\Gamma$. If $(a_\rho)$ is not jammed, then it has a 
fluent cofinal subsequence. 

\medskip\noindent
Let us say that $K$ is {\bf fluent} \index{fluent!valued field}\index{valued field!fluent}\label{p:fluent} if every fluent pc-sequence in 
$K$ has a pseudolimit in~$K$; equivalently, every divergent pc-sequence 
in $K$ is jammed; equivalently, every coarsening of $K$ by a nontrivial convex
subgroup of $\Gamma$ is maximal as a valued field. Every maximal valued field is fluent.
If $K$ is fluent, then so is the residue 
field (with its induced valuation) of any coarsening of $K$. 
The notions introduced above are only relevant when
$[\Gamma^{\neq}]$ has no smallest element.

\begin{lemma}\label{specfluent} Suppose $[\Gamma^{\ne}]$ has no smallest element. Then each special pc-se\-quence in $K$ has a cofinal fluent subsequence.
\end{lemma}
\begin{proof} Let $(a_{\rho})$ be a $\Delta'$-special pc-sequence in $K$ where
$\Delta'$ is a nontrivial convex subgroup of $\Gamma$. Take a nontrivial 
proper convex subgroup $\Delta$ of $\Delta'$. It is enough to observe that 
$(a_{\rho})$ is not $\Delta$-jammed, and so has a cofinal $\Delta$-fluent
subsequence. 
\end{proof}

\begin{cor}\label{maxflstep} Suppose $[\Gamma^{\ne}]$ has no smallest 
element and $K$ is fluent. Then $K$ is step-complete, and thus 
henselian.
\end{cor}
\begin{proof} By Lemma~\ref{specfluent} every special pc-sequence in $K$ has a pseudolimit in $K$. It remains to appeal to Lemma~\ref{lem:step-complete is henselian}. 
\end{proof}

\begin{cor}
Suppose $K$ is fluent and $L$ is a valued field extension of $K$ with $[L:K]<\infty$. Then $L$ is fluent. 
\end{cor}
\begin{proof} Let $\Delta'$ be a nontrivial convex subgroup of
$\Gamma_L$; we show that the $\Delta'$-coarsening
of $L$ is maximal. 
It follows from Corollary~\ref{cor:ef inequ, 2} that
$\Delta:= \Delta'\cap \Gamma$ is a nontrivial convex subgroup of $\Gamma$ and
$\Delta'=\Delta_L$. As the $\Delta$-coarsening of $K$ is maximal, it remains to appeal to Corollary~\ref{cor:finite ext of complete}.
\end{proof}

\subsection*{$\Delta$-immediate extensions} \index{extension!valued fields!$\Delta$-immediate}
\index{valued field!$\Delta$-immediate extension}
\index{Delta-immediate@$\Delta$-immediate extension}
 In this subsection 
{\em extension\/} means
{\em valued field extension}.
Let $\Delta$ be a convex subgroup of $\Gamma$. A {\bf $\Delta$-immediate 
extension of $K$} is an immediate extension $L$ of $K$ such that each $a\in L\setminus K$ 
is a pseudolimit of a divergent $\Delta$-fluent pc-sequence in $K$; 
equivalently,
it is an immediate extension $L$ of $K$ such that
no $a\in L\setminus K$ is a pseudolimit of a divergent $\Delta$-jammed
pc-sequence in $K$. If $L$ is a $\Delta$-immediate extension of $K$ and $M$ is a $\Delta$-immediate extension of $L$, then $M$ is a $\Delta$-immediate 
extension of $K$:
this uses the characterization of ``$\Delta$-immediate extension'' in terms of
$\Delta$-jammed pc-sequences. We now consider still another useful 
characterization in terms of coarsening.

\medskip
\noindent
Let $(K,\dot{\cal{O}})$ be the field $K$ with the coarsened valuation 
$\dot{v}=v_{\Delta}$ and its valuation ring~$\dot{\cal{O}}$.
Given an immediate extension $L$ of $K$ and using again $\Delta$ to coarsen,
we obtain a valued field extension
$(L,\dot{\cal{O}}_L)$ of $(K,\dot{\cal{O}})$. 

\begin{lemma} Let $L$ be an immediate extension of $K$. Then $L$ is a 
$\Delta$-immediate extension of $K$ iff the extension 
$(L,\dot{\cal{O}}_L)$ of $(K,\dot{\cal{O}})$ is immediate.
\end{lemma}
\begin{proof} Assume $L|K$ is $\Delta$-immediate. Let $a\in \dot{\cal{O}}_L$;
it is enough to find $b\in\dot{\cal{O}}$ such that $v(a-b)>\Delta$.
For this we can assume $a\notin K$. Take a divergent 
$\Delta$-fluent pc-sequence $(a_{\rho})$ in $K$ such that $a_{\rho} \leadsto a$
and $\gamma_\rho:= v(a_{s(\rho)}-a_{\rho})=v(a-a_{\rho})$ for all~$\rho$.
Take $\rho$ such that $a\sim a_{\rho}$ (so $\gamma_{\rho}\ge \delta$
for some $\delta\in \Delta$), and $\gamma_{s(\rho)}-\gamma_{\rho}>\Delta$.
Then for $b:= a_{s(\rho)}$ we have 
$v(a-b)=\gamma_{s(\rho)}> \Delta + \gamma_{\rho}$, so  $v(a-b)>\Delta$.

For the converse, assume the extension $(L,\dot{\cal{O}}_L)$ of $(K,\dot{\cal{O}})$ is 
immediate. Let~$a\in L\setminus K$, and take a divergent
pc-sequence $(a_{\rho})$ in $(K,\dot{\cal{O}})$ such that $a_{\rho} \leadsto a$, in~$(L,\dot{\cal{O}}_L)$. Then $(a_{\rho})$ is a $\Delta$-fluent pc-sequence
with $a_{\rho} \leadsto a$ in $L$. If $b\in K$ were such that 
$a_{\rho} \leadsto b$ in $K$, then also $a_{\rho} \leadsto b$ in
$(K,\dot{\cal{O}})$, and so no such $b$ exists.
\end{proof}

\begin{cor} Let $L$ be a $\Delta$-immediate extension of $K$. Then 
$L$ is also a $\Delta$-immediate extension of any valued field $F$ with
$K\subseteq F\subseteq L$. 
\end{cor}

\begin{cor}\label{maxdeltafluent} 
The following conditions on $K$ are equivalent: \begin{enumerate}
\item[\textup{(i)}] $K$ has no proper $\Delta$-immediate extension;
\item[\textup{(ii)}] every $\Delta$-fluent pc-sequence in $K$ pseudoconverges in $K$;
\item[\textup{(iii)}] $(K, \dot{\mathcal{O}})$ is a maximal valued field.
\end{enumerate}
If $K$ has any of these properties, then $(K, \dot{\mathcal{O}})$ is henselian.
\end{cor}
\begin{proof}
We have (ii)~$\Longleftrightarrow$~(iii) by Corollary~\ref{cor:maximal valued fields} applied to 
the $\Delta$-coarsening of~$K$ in place of $K$, and (iii)~$\Rightarrow$~(i) is clear. Let $(L,\dot{\mathcal O}_L)$ be an immediate extension of~$(K,\dot{\mathcal O})$. Then by Lemma~\ref{immcoarse} the valuation of $K$ extends to a valuation $v\colon L\to\Gamma_\infty$ making~$L$ an immediate extension of $K$ and $\dot{\mathcal O}_L$ the valuation ring of the coarsening of this valuation by $\Delta$. Equipped with this valuation, $L$ is a $\Delta$-immediate extension of $K$.
Thus (i)~$\Rightarrow (K,\dot{\mathcal O})=(L,\dot{\mathcal O}_L)$. 
This shows (i)~$\Rightarrow$~(iii). If
any of the conditions~(i),~(ii),~(iii) holds, then $(K, \dot{\mathcal{O}})$ is henselian by Corollary~\ref{cor:alg max implies henselian}.
\end{proof}

\noindent
By Zorn there exists a $\Delta$-immediate extension of $K$ that
has no proper $\Delta$-immediate extension. Such an extension of $K$ is called a
 {\bf maximal $\Delta$-immediate extension} of~$K$; under 
certain conditions it is unique up to isomorphism over $K$.
To discuss this, let us fix some
maximal $\Delta$-immediate extension $K^{\Delta}$ of $K$. 
Then $K^{\Delta}$ is also a maximal $\Delta$-immediate extension of 
any intermediate
valued field $L$ with $K\subseteq L\subseteq K^{\Delta}$, and~$K^{\Delta}$ is henselian with respect to
its $\Delta$-coarsened valuation.

\index{extension!valued fields!maximal $\Delta$-immediate}
\index{valued field!maximal $\Delta$-immediate extension}
\index{maximal!$\Delta$-immediate extension}
\nomenclature[K]{$K^\Delta$}{maximal $\Delta$-immediate extension of the valued field $K$}

\begin{prop}\label{undelta} Suppose $K$ has equicharacteristic zero. 
Then every maximal $\Delta$-immediate extension of $K$ is isomorphic over $K$ 
to the valued field $K^{\Delta}$.
\end{prop}

\begin{proof}
Let $L$ be a maximal $\Delta$-immediate extension of $K$. Then both the $\Delta$-coar\-se\-ning of $L$ and the $\Delta$-coarsening of $K^\Delta$ are maximal immediate extensions of the $\Delta$-coarsening of $K$.
Then Corollary~\ref{cor:unique max imm ext, 1} gives an isomorphism $L\to K^\Delta$ of these $\Delta$-coar\-se\-nings that is the identity on $K$. Since for each 
$b\in L^\times$ we have $b=a(1+\varepsilon)$ with $a\in K^\times$ and 
$\varepsilon\in \dot{\smallo}_L$, this isomorphism is also an isomorphism with respect to the valuations of $L$ and $K^{\Delta}$, respectively.
\end{proof}

\noindent
The following lemma is obvious, and is used in the next subsection.

\begin{lemma}\label{deldel} If $L$ is a $\Delta$-immediate extension of $K$,
then any maximal $\Delta$-immediate extension of $L$ is also a maximal
$\Delta$-immediate extension of $K$.
\end{lemma}

\subsection*{Fluent completion}
{\em In this subsection we assume that $K$ has equicharacteristic zero, 
$\Gamma\ne \{0\}$ and $[\Gamma^{\ne}]$ has no smallest element}. By 
{\em extension\/} we mean {\em valued field extension}, and {\em embedding\/} 
means {\em valued field embedding}.

Let $L$ be an extension of $K$. 
A pc-sequence $(a_{\rho})$ in $L$ is said to be {\bf $K$-fluent} if there is 
a nontrivial convex subgroup 
$\Delta$ of $\Gamma$ such 
that $(a_{\rho})$ has length $\le\operatorname{cf}(S)$ for some 
$S\subseteq \Gamma/\Delta$, and $(a_{\rho})$ is $\Delta_L$-fluent. 
We say that $L$ is $K$-fluent if $L$ is henselian
and every $K$-fluent pc-sequence
in $L$ pseudoconverges in $L$. Note that if $L$ is fluent, then $L$ is
$K$-fluent. 
Also, if $L$ is $\Gamma$-maximal, then $L$ is $K$-fluent. 
The following is a variant of Proposition~\ref{undelta}.

\index{pc-sequence!$K$-fluent}
\index{K-fluent@$K$-fluent pc-sequence}

\begin{lemma}\label{fldelta} Let $\Delta$ be a nontrivial convex subgroup of $\Gamma$ and let
$K^{\Delta}$ be as above.
% a maximal $\Delta$-immediate extension of $K$. 
Then $K^{\Delta}$ embeds over $K$ into any $K$-fluent extension of $K$.
\end{lemma}  
\begin{proof} Let $L$ be an intermediate valued field: $K\subseteq L \subseteq K^{\Delta}$, assume $L\ne K^{\Delta}$, and let $i\colon L \to F$ be an embedding over 
$K$ into a
$K$-fluent extension $F$ of $K$. It suffices to show that $i$ 
can be extended to an embedding into $F$ of some intermediate valued field 
that strictly contains $L$. 
Note that $i$ remains an embedding when replacing~$L$ by its 
$\Delta$-coarsening, and $F$ by its $\Delta_F$-coarsening. Since the 
$\Delta$-coarsening of $K^{\Delta}$ is henselian and the $\Delta_F$-coarsening
of $F$ is henselian, we can extend $i$ and arrange that the
$\Delta$-coarsening of $L$ is henselian. 
Take $a\in K^{\Delta}\setminus L$. Take a divergent $\Delta$-fluent
pc-sequence~$(a_{\rho})$ in $L$ such that $a_{\rho}\leadsto a$.
By passing to a cofinal subsequence we arrange that $(a_{\rho})$
has length  $\le\operatorname{cf}(S)$ for some 
$S\subseteq \Gamma/\Delta$. Then~$(a_{\rho})$ is of transcendental type over the $\Delta$-coarsening of
$L$. In view of Corollary~\ref{cor:pc and polynomials} and Lemma~\ref{pcDeltainv}, $(a_{\rho})$ is of transcendental type over $L$.  Also $\big(i(a_{\rho})\big)$ is a $K$-fluent pc-sequence in $F$, so 
$i(a_{\rho}) \leadsto b$ with $b\in F$, hence we can extend $i$ to an
embedding $L(a) \to F$ sending~$a$ to~$b$. 
\end{proof}

\begin{prop}\label{fluentcompletion} $K$ has an immediate fluent extension that embeds 
over $K$ into any $K$-fluent extension of $K$.
\end{prop}

\noindent
We call an extension of $K$ as in this proposition a {\bf fluent completion} of $K$.

\index{valued field!fluent completion}
\index{fluent!completion}

\begin{proof} Fix a decreasing coinitial sequence $(\Delta_{\alpha})$
of nontrivial convex subgroups of~$\Gamma$ indexed by the ordinals
$\alpha< \lambda$ for some infinite limit ordinal $\lambda$, where ``coinitial''
means that every nontrivial convex subgroup $\Delta$ of $\Gamma$ includes 
$\Delta_{\alpha}$ for some $\alpha$ (and thus for all sufficiently large 
$\alpha$). For each $\alpha$ we pick a maximal $\Delta_{\alpha}$-immediate
extension~$K^{\Delta_{\alpha}}$ of~$K$. We arrange this so that
$K^{\Delta_{\alpha}}$ is a valued subfield of $K^{\Delta_{\beta}}$ whenever
$\alpha < \beta < \lambda$: by transfinite recursion and using 
Lemma~\ref{deldel},
take $K^{\Delta_{\alpha+1}}$ to be a maximal $\Delta_{\alpha+1}$-immediate extension 
of $K^{\Delta_{\alpha}}$ for $\alpha< \lambda$, and if $\beta< \lambda$ is an 
infinite limit ordinal, take $K^{\Delta_{\beta}}$ to be a maximal
$\Delta_{\beta}$-immediate extension of the $\Delta_{\beta}$-immediate extension
$\bigcup_{\alpha< \beta}K^{\Delta_{\alpha}}$ of $K$. Put 
$$  K^{\operatorname{f}}\ :=\  \bigcup_{\alpha< \lambda} K^{\Delta_{\alpha}}.$$
Then  $K^{\operatorname{f}}$ is an immediate extension of $K$ with the following properties: 
\begin{enumerate}
\item[(1)] $K^{\operatorname{f}}$ is henselian;
\item[(2)] every fluent pc-sequence in $K$ pseudoconverges in $K^{\operatorname{f}}$;
\item[(3)] every $a\in K^{\operatorname{f}}\setminus K$ is a pseudolimit of some divergent fluent pc-sequence in $K$;
\item[(4)] $K^{\operatorname{f}}$ embeds over $K$ into any $K$-fluent 
extension of $K$.
\end{enumerate}
To get (1), let $n\ge 2$, $a_2,\dots, a_n\in  K^{\operatorname{f}}$, $a_2,\dots, a_n \prec 1$; our job is to show that then the polynomial 
$1+X + a_2X^2 + \dots + a_nX^n$ has a zero in the valuation ring of 
$K^{\operatorname{f}}$.  Take $\alpha< \lambda$ so large that
$a_2,\dots,a_n \in  K^{\Delta_{\alpha}}$ and $v(a_2),\dots, v(a_n)> \Delta_{\alpha}$. 
(This is possible since $[\Gamma^{\ne}]$ has no smallest element.)
Since $K^{\Delta_{\alpha}}$ is henselian
with respect to its $\Delta_{\alpha}$-coarsened valuation, 
our polynomial does have a zero $b\in K^{\Delta_{\alpha}} $ with 
$v(1+b)> \Delta_{\alpha}$, and so $v(b)\ge 0$. This proves (1).

To get (4), let $L$ be any $K$-fluent 
extension of $K$. Now use Lemma~\ref{fldelta} and transfinite 
recursion to obtain for each $\alpha< \lambda$ an embedding $i_{\alpha}\colon K^{\Delta_{\alpha}} \to L$ such that~$i_{\beta}$ extends $i_{\alpha}$ whenever 
$\alpha < \beta < \lambda$. 

\index{valued field!semifluent completion}
\index{semifluent completion}

We define a {\em semifluent completion\/} of $K$ to be 
an immediate extension of $K$ with the properties (1)--(4) of 
$K^{\operatorname{f}}$. So we have shown that $K$ has a semifluent completion. 
Let any ordinal $\nu>0$ be given. We now build an 
increasing sequence~$(K_\mu)$ of immediate extensions of $K$, indexed by
the ordinals $\mu < \nu$, such that
$K_0=K$, $K_{\mu +1}$ is a semifluent completion of $K_{\mu}$
whenever $\mu < \mu+1 < \nu$, and $K_{\mu}=\bigcup_{\alpha< \mu} K_{\alpha}$
whenever~$\mu< \nu$ is an infinite limit ordinal. With $\nu$ large enough, 
we obtain $\mu < \mu+1< \nu$ such that $K_{\mu}=K_{\mu+1}$, and it is easy to
check that then $K_{\mu}$ is a fluent completion of $K$.      
\end{proof}

\noindent
Are all fluent completions of $K$ isomorphic over $K$, for every $K$?
We do not know the answer, but we can proceed without. Here is an important property of fluent completions, used in Sections~\ref{sec:uplfree} and~\ref{sec:behupo}.  

\begin{lemma}\label{lem:no jammed pseudolimits} Suppose $(a_\rho)$ is a jammed pc-sequence in $K$ with a pseudolimit in a fluent completion of $K$.
Then $(a_\rho)$ has a pseudolimit in $K$.
\end{lemma}
\begin{proof} Every fluent completion of $K$ embeds over $K$ into the
fluent completion~$K_{\mu}$ from the proof of 
Proposition~\ref{fluentcompletion}. By 
transfinite induction this gives a reduction to the case that $(a_\rho)$ has a 
pseudolimit in a semifluent completion of $K$ as defined in that proof. 
Since every pc-sequence in $K$ equivalent to $(a_\rho)$ is also 
jammed, it remains to use property~(3) of semifluent completions and 
Corollary~\ref{eqpccor}. 
\end{proof}

\subsection*{Approximation by special pc-sequences}
The main goal of this subsection is the result below about approximating elements in the henselization $K^{\operatorname{h}}$ of the valued field $K$. It generalizes the fact (immediate from Corollary~\ref{cor:henselian complete rank 1}) that for $K$ of rank $1$ every element in $K^{\operatorname{h}}$ is the limit of a c-sequence in $K$. 

\begin{prop}[F.-V.~Kuhlmann]\label{prop:approx henselization}
For each $x\in K^{\operatorname{h}}\setminus K$ there is a divergent special pc-sequence $(x_\rho)$ in $K$ and some $a\in K^\times$ such that $x_\rho\leadsto x/a$.
\end{prop}

\noindent
We give the proof after a number of auxiliary results. 
\textit{In the rest of this subsection,~$\Delta$ is a nontrivial convex subgroup of $\Gamma$.}\/
We begin with a reformulation of the conclusion of this proposition:

\begin{lemma}
Let $x$ be an element of a valued field extension of $K$ with $x\notin K$, and let $a\in K^\times$.  Then the following are equivalent:
\begin{enumerate}
\item[\textup{(i)}] there is a divergent $\Delta$-special pc-sequence $(x_\rho)$ in $K$ such that $x_\rho\leadsto x/a$;
\item[\textup{(ii)}] the coset $\alpha+\Delta$, where $\alpha=va$, is a cofinal subset of $v(x-K)$.
\end{enumerate}
\end{lemma}
\begin{proof}
Let $(x_\rho)$ be a divergent $\Delta$-special pc-sequence in $K$ such that $x_\rho\leadsto x/a$.
We may assume that $v(x/a-x_\rho)\in\Delta$ for all $\rho$, so $v(x-ax_\rho)\in \alpha+\Delta$ for each~$\rho$, and $\big(v(x-ax_\rho)\big)$ is cofinal in $\alpha+\Delta$. Thus $\alpha+\Delta\subseteq v(x-K)$ since $v(x-K)$ is downward closed.
Since the pc-sequence $(ax_\rho)$ in $K$ diverges and $ax_\rho\leadsto x$, the sequence $\big(v(x-ax_\rho)\big)$ is cofinal in $v(x-K)$. Hence $\alpha+\Delta$ is cofinal in $v(x-K)$.
Conversely, suppose $\alpha+\Delta\subseteq v(x-K)$ and $\alpha+\Delta$ is cofinal in $v(x-K)$.
Since $\Delta\neq\{0\}$, the set $\alpha+\Delta$, and hence also the set $v(x-K)$, does not have a largest element.
Choose a well-indexed sequence $(y_\rho)$ in $K$ such that $v(x-y_\rho)\in\alpha+\Delta$ for each $\rho$ and
$\big(v(x-y_\rho)\big)$ is strictly increasing and cofinal in $\alpha+\Delta$, and hence in $v(x-K)$. Then $(x_\rho)$ where $x_\rho:=y_\rho/a$ for each $\rho$ is a divergent $\Delta$-special pc-sequence in $K$ such that $x_\rho\leadsto x/a$.
\end{proof}

\noindent
Let $x$ be an element of a valued field extension of $K$. We say that $x$ is {\bf almost $\Delta$-special over $K$} if $x\notin K$ and for some  $a\in K^\times$
the equivalent conditions (i) and (ii) in the previous lemma hold. We say that $x$ is  {\bf  $\Delta$-special over $K$} if $x\notin K$ and conditions (i) and (ii) hold for $a=1$.
So if $x$ is almost $\Delta$-special over $K$, then  some $K^\times$-multiple of~$x$ is  $\Delta$-special over $K$. 
Moreover, for $a,b\in K$, $a\neq 0$, 
$$v(ax+b-K) = va+v(x-K),$$
hence $ax+b$ is almost $\Delta$-special over $K$ iff $x$ is.
We call $x$ {\bf almost special over $K$} if~$x$ is almost $\Delta$-special over $K$ for some $\Delta$, and similarly we say that $x$ is {\bf  special over~$K$} if 
$x$ is  $\Delta$-special over $K$ for some $\Delta$.
%From Lemma~\ref{specfluent} we obtain:

%\begin{lemma}
%Suppose $[\Gamma^{\neq}]$ does not have a smallest element. If the element $x$ in %a valued field extension of $K$ is  special over $K$, then there is no divergent %jammed pc-sequence in $K$ with pseudolimit $x$.
%\end{lemma}

\index{element!almost $\Delta$-special over $K$}
\index{element!$\Delta$-special over $K$}
\index{element!almost special over $K$}
\index{element!special over $K$}

\medskip\noindent
In the next lemmas we consider valued field extensions $L$ of $K$. For such $L$ we let~$\Delta_L$ be the convex hull of $\Delta$ in $\Gamma_L$, and $\dot K$ is the valued residue field of the coarsening $v_\Delta$ of the valuation $v$ in $K$ by $\Delta$, viewed as a valued subfield of the valued residue field $\dot L$ of the coarsening $v_{\Delta_L}$ of the valuation on $L$ by $\Delta_L$ as usual.

\begin{lemma}\label{lem:Delta-special}
Let $x\in L$ and $vx\in\Delta_L$. Then $x$ is  $\Delta$-special over $K$ iff
$\dot x\notin\dot K$ and~$\dot x$ is the limit in $\dot{L}$ of a c-sequence in $\dot K$.
\end{lemma}
\begin{proof}
Let $(x_\rho)$ be a divergent $\Delta$-special pc-sequence in $K$ with $x_\rho\leadsto x$. 
Then $vx=v(x_\rho)$ for sufficiently large $\rho$, and after passing to a cofinal subsequence, we may assume that this holds for all $\rho$. Then $(\dot x_\rho)$ is a c-sequence in $\dot K$ and $\dot{x}_{\rho} \to \dot{x}$. Also, $\dot{x}\notin \dot{K}$: otherwise, $\dot{x}=\dot{a},\ a\in \dot{\mathcal O}$, and then $x_\rho\leadsto a$ by Lemma~\ref{eqpslim}, a contradiction.  Conversely, suppose $(x_\rho)$ is a well-indexed sequence in $\dot{\mathcal O}$ such that~$(\dot x_\rho)$ is a c-sequence in $\dot K$  and $\dot x_\rho\to\dot x\notin\dot K$.
After passing to a cofinal subsequence, we may assume that $(\dot x_\rho)$ is a pc-sequence of width $\{\infty\}$ in $\dot K$, by Lemma~\ref{lem:c-sequence width}; thus~$(x_\rho)$ is a $\Delta$-special pc-sequence in $K$ with $x_\rho\leadsto x$. If $a\in K$ is a pseudolimit of $(x_\rho)$, then $va=vx_\rho$ eventually, in particular $va\in\Delta$, and
$\dot x_\rho\to\dot a= \dot x$, a contradiction. Hence~$(x_\rho)$ does not have a pseudolimit in $K$, and so~$x$ is $\Delta$-special
over $K$.
\end{proof}

\noindent
Note that if $x\in L$ is  $\Delta$-special over $K$ with $vx\in\Delta_L$, then by the previous lemma the valued field extension $\dot K(\dot x)\supseteq \dot K$
is dense, and hence every $y\in L$ with $vy\in\Delta_L$ and $\dot y\in \dot K(\dot x)\setminus \dot K$ is  $\Delta$-special over $K$.

\medskip
\noindent
The following lemma indicates a source of  special elements over $K$. 

\begin{lemma}\label{lem:hens, special}
Let $P\in\mathcal O[X]$ and $a\in\mathcal O$ be such that $P(a)\prec 1$, $P'(a)\asymp 1$, and~$P$ has no zero in $a+\smallo$. Let $x$ in a valued field extension $L$ of $K$ be a zero of $P$ with $x-a\prec 1$. Then $vx\in\Delta$ and
$x$ is $\Delta$-special, for some $\Delta$. 
\end{lemma}
\begin{proof}
By the discussion preceding Corollary~\ref{cor:henselian complete rank 1}, we have
$\Delta$ and a $\Delta$-special divergent pc-sequence $(a_\rho)$ in $K$ such that $P(a_\rho)\leadsto 0$ and $a_\rho\equiv a\bmod\smallo$ for all $\rho$. 
Taylor expansion around $x\in L$ yields
\begin{align*}
P(x+Y)	&= P(x) + P'(x)\cdot Y+\text{terms of higher degree in $Y$} \\
		&= P'(x)\cdot Y \cdot (1+Q)\qquad\text{where $Q\in\mathcal O_L[Y]$, $Q(0)=0$.}
\end{align*}
Substituting $a_\rho-x$ for $Y$ yields 
$$P(a_\rho)=P'(x)\cdot (a_\rho-x)\cdot \big(1+Q(a_\rho-x)\big)\qquad\text{with $Q(a_\rho-x)\prec 1$.}$$
Thus $v\big(P(a_\rho)\big)=v(a_\rho-x)$ for all $\rho$, hence $a_{\rho} \leadsto x$, and so $x$ is $\Delta$-special over $K$. Also, $vx=va_{\rho}< v(x-a_{\rho})\in \Delta$, eventually, so $vx\in\Delta$.
\end{proof}

\noindent
We say that a valued field extension $L$ of $K$ is {\bf almost $\Delta$-special} if every element of~$L\setminus K$ is almost $\Delta$-special over $K$. A valued field extension $L$ of $K$ is {\bf almost special} if every element of $L\setminus K$ is almost special over $K$. (So Proposition~\ref{prop:approx henselization} says that the henselization of $K$ is an almost special extension of $K$.) Thus every almost special valued field extension is immediate. The next lemma shows that if $M\supseteq L$ and $L\supseteq K$ are almost special valued field extensions, then so is $M\supseteq K$:

\index{extension!valued fields!almost special}
\index{extension!valued fields!almost $\Delta$-special}

\begin{lemma}\label{lem:almost special transitive}
Let $L$ be an almost special valued field extension of $K$ and $x$ an element of a valued field extension of $L$ such that $x$ is almost special over $L$. Then~$x$ is almost special over $K$.
\end{lemma}
\begin{proof}
Clearly $v(x-K)\subseteq v(x-L)$, and if equality holds, then $x$ is almost special over $K$. So suppose $v(x-K) \ne v(x-L)$. 
We can take $y\in L$ with $v(x-K)<v(x-y)$; then $y\notin K$, so $y$ is almost special over~$K$. For each $a\in K$ we have $v(y-a)=v\big((y-x)+(x-a)\big)=v(x-a)$, hence $v(y-K)=v(x-K)$. Thus $x$ is almost special over $K$.
\end{proof}

\noindent
We now define a condition on an element $x$ of a valued field extension of~$K$ which ensures that $K(x)\supseteq K$ is almost special: call an element $x$ in a valued field extension~$L$ of $K$ {\bf very $\Delta$-special over $K$} if 
\begin{list}{*}{\setlength\leftmargin{3em}}
\item[(VS1)] $vx\in\Delta_L$ and $x$ is $\Delta$-special over $K$;
\item[(VS2)] for all $n\ge 1$, if $1,x,\dots,x^n$ are $K$-linearly independent, then $1,\dot x,\dots,\dot x^n$ in~$\dot{L}$ are $\dot K$-linearly independent.
\end{list}
Here  $\dot K$ is the valued residue field of the coarsening~$v_\Delta$ of the valuation $v$ in $K$ by~$\Delta$, viewed as a valued subfield of the valued residue field $\dot L$ of the coarsening~$v_{\Delta_L}$ of the valuation on $L=K(x)$ by the convex hull $\Delta_L$ of $\Delta$ in $\Gamma_L$. Note that given~(VS1),  
condition (VS2)  expresses that $\big[K(x):K\big]=\big[\dot K(\dot x):\dot K\big]\in \N\cup\{\infty\}$. Thus a very $\Delta$-special $x$ over $K$ is transcendental over $K$ iff $\dot x$ is transcendental over 
$\dot K$.

\index{element!very $\Delta$-special over $K$}

\begin{lemma}\label{lem:very special}
Let $x$ be very $\Delta$-special over $K$. Then $K(x)\supseteq K$ is almost $\Delta$-special.
\end{lemma}
\begin{proof}
Let $y\in K(x)\setminus K$; we need to show that $y$ is almost $\Delta$-special over~$K$.
We first assume that $x$ is algebraic over $K$.  Then $y=P(x)$ with $P\in K[X]$ of degree~$m$ with $1\le m<\big[K(x):K\big]$. Replacing $y$ by $ay+b$ and $P$ by $aP+b$, for suitable $a,b\in K$, $a\neq 0$, we may assume that $vP=0$ and $P(0)=0$.
Thus $P\in\dot{\mathcal O}[X]$ with $\deg \dot P\ge 1$. By (VS2), 
$1,\dot x,\dots,\dot x^m$ are $\dot K$-linearly independent, hence $\dot y=\dot P(\dot x)\notin\dot K$.
Thus by (VS1) and the remark following Lemma~\ref{lem:Delta-special}, $y=P(x)$ is  $\Delta$-special over~$K$.

Now suppose that $x$ is transcendental over $K$; then $\dot x$ is transcendental over~$\dot K$ by~(VS2). So $y=P(x)/Q(x)$ with $P,Q\in K[X]$, $Q\neq 0$. 
After multiplying~$P$,~$Q$ and~$y$ by suitable elements of $K^\times$ we may assume that $vP=vQ=0$.
Now $P=\sum_i P_i X^i$, $Q=\sum_j Q_jX^j$ with $P_i,Q_j\in\mathcal O$ for all $i$, $j$; set $n=\deg \dot Q$. Then $Q_n\neq 0$, so with $b:=P_n/Q_n$ we have $P_n-bQ_n=0$. Take $a\in K^\times$ such that $v(P-bQ)=va$. Then with $R:=a^{-1}(P-bQ)\in\mathcal O[X]$,
we have $v(R)=0$ and $R_n=0$, so
$\dot R(\dot x)/\dot Q(\dot x)\notin \dot K$. By (VS1) and the remark after Lemma~\ref{lem:Delta-special},
$R(x)/Q(x)$, and hence also $y=b+ aR(x)/Q(x)$,
are $\Delta$-special over~$K$.
\end{proof}

\noindent
An element of a valued field extension of $K$ is said to {\bf very special over $K$} if it is very $\Delta$-special over $K$ for some
$\Delta$.
Together with Lemma~\ref{lem:almost special transitive}, the previous lemma implies:

\begin{cor}\label{cor:almost special transitive}
Let $L$ be an almost special valued field extension of $K$ and~$x$ an element of a valued field extension of $L$ such that $x$ is very special over $L$. Then $L(x)\supseteq K$ is almost special.
\end{cor}

\noindent
Here is another consequence of Lemma~\ref{lem:very special}:

\begin{cor}\label{cor:transcendental type, almost special}
Suppose  $K$ is henselian of equicharacteristic zero. Let~$(a_\rho)$ be a divergent $\Delta$-special pc-sequence in $K$, and let $x$ be a pseudolimit of $(a_\rho)$ in a valued field extension of~$K$. Then $x-a$ is very $\Delta$-special over $K$, for some $a\in K$; in particular, $K(x)\supseteq K$ is almost $\Delta$-special.
\end{cor}
\begin{proof}
Take $\rho_0$ such that $v(x-a_\rho)=v(a_\rho-a_\sigma)\in\Delta$ for $\rho_0\leq \rho<\sigma$, and set $a:=a_{\rho_0}$;
then $(b_\rho)_{\rho\geq\rho_0}:=(a_\rho-a)_{\rho\geq\rho_0}$ is a divergent $\Delta$-special pc-sequence in~$K$ with pseudolimit $y:=x-a$, and $vy\in\Delta$ and
$v(b_\rho)\in\Delta$ for $\rho_0\leq\rho$.
The sequence~$(\dot b_\rho)$ in $\dot K$ is a divergent pc-sequence in $\dot K$ (of width $\{\infty\}$) with pseudolimit~$\dot y$. Since~$K$ is henselian of equicharacteristic zero, so is $\dot K$, by Lemma~\ref{lem:henselian decomposition}. Thus~$(b_\rho)$ is of transcendental type over $K$ and 
$(\dot b_\rho)$ is of transcendental type over $\dot K$, by Lemma~\ref{lem:Kaplansky, 2} and Corollary~\ref{cor:alg max equals henselian}. Hence $y$ is transcendental over $K$ and $\dot y$ is transcendental over~$\dot K$ by Lemma~\ref{lem:Kaplansky, 1}. Thus $y$ is very $\Delta$-special over $K$, and so $K(x)=K(y)\supseteq K$ is almost $\Delta$-special, by Lemma~\ref{lem:very special}.
\end{proof}

%\no indent
%An obvious transfinite induction in combination with Lemmas~\ref{lem:almost special transitive} and \ref{lem:very special} yields:

%\begin{cor}
%Let $(a_\lambda)_{\lambda<\mu}$ be a sequence of elements of some valued field extension of $K$, indexed by all ordinals $\lambda$ less than some ordinal $\mu$. Suppose that for each $\nu<\mu$, the element $a_{\nu}$ is very special over $K(a_{\lambda}:\lambda<\nu)$. Then the valued field extension $K(a_\lambda:\lambda<\mu)\supseteq K$ is almost special.
%\end{cor}

\noindent
In analogy to Lemma~\ref{lem:hens, special}, we have:

\begin{lemma}\label{lem:hens, very special}
Let $P\in\mathcal O[X]$ be monic and let $a\in\mathcal O$ be such that $P(a)\prec 1$, $P'(a)\asymp 1$, and $P$ has no zero in $a+\smallo$. 
Suppose that for each $\Delta$, either $\dot P\in\mathcal O_{\dot K}[X]$ is irreducible over $\dot K$ or has a zero in $\dot a+\smallo_{\dot K}$.
Then each zero $x$ of $P$ in any valued field extension of $K$ with $x-a\prec 1$ is very special over $K$.
\end{lemma}
\begin{proof}
Let $x$ be an element in a valued field extension of $K$ with $P(x)=0$ and $x-a\prec 1$.
By Lemma~\ref{lem:hens, special} we can take $\Delta$ such that $x$ is $\Delta$-special over $K$ and $vx\in\Delta$. Then $\dot P(\dot a)\prec 1$, $\dot P'(\dot a)\asymp 1$, as well as
$\dot P(\dot x)=0$, $\dot x - \dot a \prec 1$, and $\dot x\notin\dot K$ (by Lemma~\ref{lem:Delta-special}). 
Thus, by Lemma~\ref{lem:henselian, uniqueness}, $\dot P$ has no zero in $\dot a+\smallo_{\dot K}$, so
by hypothesis,  $\dot P$ is irreducible. It follows that $\dot P$ is the minimum polynomial of $\dot x$ over $\dot K$. With $L:=K(x)$, we now have
$$[L:K]\ \leq\ \deg P\ =\  \deg \dot P\  =\  
\big[\dot K(\dot x):\dot K\big]\  \leq\ [\dot{L}:\dot K]\ \leq\  [L:K]$$
and thus $[L:K]=\big[\dot K(\dot x):\dot K\big]$. Hence $x$ is very $\Delta$-special over $K$.
\end{proof}

\noindent
We can now give the proof of Proposition~\ref{prop:approx henselization}:

\begin{proof}[Proof of Proposition~\ref{prop:approx henselization}]
Let $L$ be a valued subfield of $K^{\operatorname{h}}$ containing $K$ such that $L\supseteq K$ is almost special. If $L=K^{\operatorname{h}}$, then we are done, so suppose otherwise; then $L$ is not henselian, by the minimality of henselizations. Take a monic polynomial $P\in\mathcal O_L[X]$ of minimal degree such that for some $a\in\mathcal O_L$ we have $P(a)\prec 1$ and $P'(a)\asymp 1$, and $P$ has no zero in $a+\smallo_L$. Fix such an $a$. We claim that
for every $\Delta$ either $\dot{P}\in \dot{L}[X]$ is irreducible, or $\dot{P}$ has a zero in $\dot a+\smallo_{\dot L}$. To prove this claim, suppose towards a contradiction that $\Delta$ is such that $\dot P\in \dot L[X]$ is reducible and has no zero in $\dot a+\smallo_{\dot L}$. Then we have monic $Q,R\in \mathcal{O}_L[X]$ of degree~$\ge 1$ such that $\dot{P}=\dot{Q}\dot{R}$. Then
$\dot{Q}(\dot{a})\prec 1$ or $\dot{R}(\dot{a})\prec 1$, say $\dot{Q}(\dot{a})\prec 1$. Then also $\dot{Q}'(\dot{a}) \asymp 1$ and $\dot{Q}$ has no zero in~$\dot a+\smallo_{\dot L}$. Hence $Q(a)\prec 1$, $Q'(a)\asymp 1$, and $Q$ has no zero in~$a+\smallo_L$, and $\deg Q< \deg P$, contradicting the minimality of $\deg P$.
This proves the claim. Take a zero $x\in K^{\operatorname{h}}$ of $P$ with $x-a\prec 1$. By Lemma~\ref{lem:hens, very special}, $x$ is very special over $L$, and so by Lemma~\ref{lem:very special}, the valued field extension $L(x)\supseteq L$ is almost special. Hence $L(x)\supseteq K$ is almost special, by Lemma~\ref{lem:almost special transitive}. It now remains to appeal to Zorn. 
\end{proof}

\noindent
We use Proposition~\ref{prop:approx henselization} to refine Proposition~\ref{prop:step-completion}:

\begin{cor}
Let $K$ have equicharacteristic zero. Then any step-com\-ple\-tion of $K$ is almost special over $K$.
\end{cor}
\begin{proof} For a cardinal $\kappa$ we let $\kappa^+$ be the next bigger 
cardinal.
Fix a maximal immediate valued field extension $M$ of $K$.
In the proof of Proposition~\ref{prop:step-completion} we constructed a step-completion $K^{\operatorname{sc}}$ of $K$  as the union of an  increasing sequence $(K_\lambda)_{\lambda<\nu}$ of henselian valued subfields of $M$, indexed by all ordinals~$\lambda$ less than some ordinal $\nu<\abs{K}^+$, as follows: $K_0=K^{\operatorname{h}}$, and for $0<\mu<\nu$, if~$\mu$ is a limit ordinal, then $K_\mu=\bigcup_{\lambda<\mu} K_\lambda$, and if $\mu$ is a successor ordinal, $\mu=\lambda+1$,
then $K_\mu=K_{\lambda}(a)^{\operatorname{h}}$ where $a\in M$ is a pseudolimit of a special divergent pc-sequence in $K_{\lambda}$. (Here the superscript $\operatorname{h}$ refers to henselization in $M$.)

Now transfinite induction using Proposition~\ref{prop:approx henselization} and Co\-rol\-lary~\ref{cor:transcendental type, almost special}  implies
that $K_\mu\supseteq K$ is almost special for each $\mu < \nu$. So $K^{\operatorname{sc}}$ is almost special over $K$.
Any step-completion of $K$ embeds into $K^{\operatorname{sc}}$ over $K$, and is 
therefore
almost special over $K$ as well.
\end{proof}

\subsection*{Notes and comments}
Lemma~\ref{lem:henselian decomposition} is in Nagata~\cite{Nagata}.
Step-completeness was defined by Krull~\cite[p.~177]{Krull} in an ideal-theoretic way; Ribenboim~\cite{Ribenboim58} has the connection to special pc-sequences, see also~\cite{Ribenboim}. Lemma~\ref{lem:step-complete decomp} is from \cite{Ribenboim58},
Lemma~\ref{lem:McL} is in Mac~Lane~\cite{MacLane}, Lemma~\ref{lem:step-complete is henselian} in Krull~\cite{Krull}, and Lemma~\ref{lem:finite ext of step-complete} in Ribenboim~\cite{Ribenboim60}. The proofs of these facts given here follow \cite{Warner}.
According to \cite{MacLane} the question posed after Proposition~\ref{prop:step-completion} goes back to Krull. It is discussed, in the context of valued ordered fields, in \cite{Tressl}.
Proposition~\ref{prop:approx henselization} is proved in \cite{FVK-approx}.

%% file: mt-3-5.tex
\section{Valued Ordered Fields}\label{sec:valued ordered fields}

\noindent
The basic facts on ordered and real closed fields, due to Artin and Schreier, are summarized without proof in Theorems~\ref{thm:AS} and 
~\ref{thm:AS, 2} below. Next, we focus on
the relevant valuations on ordered fields. Indeed,
the more robust features of a non-archimedean ordered field are 
better described in terms of a valuation than in terms of the ordering.
Throughout this section $K$ is (at least) a field. 

\subsection*{Ordered fields} Recall that in Section~\ref{sec:oag} we defined and briefly discussed ordered fields. 
{\em In this subsection $K$ is an ordered field, so $\Q\subseteq K$}.

Note that $K^{>}$ with the induced ordering is an ordered multiplicative group. Obviously, the ordered additive group of $K$ is archimedean iff 
$\operatorname{conv}(\Q)=K$; in this case we also call the ordered field $K$ {\bf archimedean}. Lemma~\ref{lem:hoelder} has an analogue for ordered fields:

\index{ordered field!archimedean}
\index{archimedean!ordered field}

\begin{lemma}[H\"older] \label{lem:hoelder, 2}
If $K$ is archimedean, then the unique embedding $K\to\R$ of ordered additive groups sending $1\in K$ to $1\in \R$  is a ring morphism.
\end{lemma}

\noindent
We note the following easy bound on zeros of polynomials over ordered fields:

\begin{lemma}\label{boundrealalg}
Let $P=X^d + a_{d-1}X^{d-1} + \cdots + a_0$ with all $a_i\in K$. Set
$M := 1 + \abs{a_{d-1}} + \cdots +\abs{a_0}$.
Then for all $x\in K$ with $\abs{x}\geq M$ we have $P(x) = x^d(1+\epsilon)$
with $\epsilon\in K$, $|\epsilon|< 1$, so $P(x)\ne 0$. 
\end{lemma}

\begin{cor}\label{cor:boundrealalg} Let
$K(x)$ be a field extension of $K$ with $x$ transcendental over~$K$. Then there is a unique ordering of $K(x)$ that makes $K(x)$ an ordered field extension of~$K$ with $x>K$.
\end{cor}
\begin{proof} For any such ordering on $K(x)$ and monic polynomial 
$P\in K[X]$ we have $P(x)>0$ by Lemma~\ref{boundrealalg}, and so there can be at most one such ordering. This also shows how to define such an ordering; alternatively, the existence of such an ordering follows by a routine model-theoretic compactness argument.
\end{proof}

\subsection*{Real closed fields}
Call a field $K$ {\bf orderable} if some ordering of $K$ makes $K$ an ordered field; note that then $\operatorname{char}(K)=0$. No algebraically closed field is orderable.

\medskip\noindent
Call $K$ {\bf euclidean} if $x^2+y^2\ne -1$ for all $x,y\in K$, and 
$$K\ =\ \{x^2:\ x\in K\} \cup \{-x^2:\ x\in K\}.$$  
%$$a\leq b \quad :\Longleftrightarrow\quad \text{there is some $y\in K$ with $b-a=y^2$}$$
%defines an ordering of $K$ making $K$ an ordered field. 
If $K$ is euclidean, then $K$ is an ordered field for a unique ordering, namely  
$$a\geq 0\quad \Longleftrightarrow\quad a=x^2\text{ for some $x\in K$.}$$

\index{field!orderable}
\index{field!euclidean}
\index{orderable field}
\index{euclidean!field}
\nomenclature[A]{$\imag$}{an element of a field with $\imag^2=-1$}

\begin{theorem}[Artin \& Schreier] \label{thm:AS}
The following conditions on $K$ are equivalent:
\begin{enumerate}
\item[\textup{(i)}] $K$ is orderable, and $K$ has no orderable proper algebraic field extension;
\item[\textup{(ii)}] $K$ is euclidean, and every $P\in K[X]$ of odd degree has a zero in $K$;
\item[\textup{(iii)}] $K$ is not algebraically closed, and $K(\imag)$, where $\imag^2=-1$, is algebraically closed;
\item[\textup{(iv)}] $K$ is not algebraically closed and has an algebraically closed field extension $L\supseteq K$ with 
$[L:K]<\infty$.
\end{enumerate}
\end{theorem}

\noindent
Call $K$ {\bf real closed} if it satisfies the (equivalent) conditions of Theorem~\ref{thm:AS}. 

\index{closed!real}
\index{field!real closed}
\index{real closed field}

\begin{cor}\label{cor:AS}
Let $K'$ be a subfield of the real closed field $K$. Then: 
$$ \text{$K'$ is real closed } \Longleftrightarrow\ 
 \text{$K'$ is algebraically closed in $K$.}$$
\end{cor}
\begin{proof}
If $K'$ is real closed, then $K'$ is algebraically closed in every orderable 
field extension, in particular, in $K$. Conversely, suppose $K'$ is algebraically closed in~$K$. Then $K'$ is euclidean and every $P\in K'[X]$ of odd degree has a zero in $K'$, since~$K'$ inherits these properties from $K$.
\end{proof}

\noindent
Here is an obvious consequence of Corollary~\ref{cor:AS}:

\begin{cor}\label{cor:cor:AS} Let 
$(K_i)_{i\in I}$ with $I\ne \emptyset$ be a family of real
closed subfields of a real closed field $K$. Then
$\bigcap_i K_i$ is a real closed subfield of $K$.
\end{cor}

\noindent
The archetypical example of a real closed field is the field $\R$ of real numbers. By the previous corollary, the algebraic closure of $\Q$ in $\R$ (known as the field of real algebraic numbers) is also real closed.
Below we always consider a real closed field as equipped with the unique ordering making it an ordered field. Here are some further properties of real closed fields.
%\medskip
%\noindent
%To condition (2) in Theorem~\ref{thm:AS} we can add some further equivalent conditions:

\begin{prop}\label{prop:AS}
Suppose $K$ is a real closed field and $P\in K[X]$. Then:
\begin{enumerate}
\item[\textup{(i)}] $P$ is monic and irreducible in $K[X]$ iff $P=X-a$ for some $a\in K$, or $P=(X-a)^2+b^2$ for some $a,b\in K$ with $b\ne 0$;
\item[\textup{(ii)}] the map $x\mapsto P(x)\colon K\to K$ has the intermediate value property.
\end{enumerate}
\end{prop}
\begin{proof} The quadratic polynomials 
in (i) take only values in $K^{>}$, and are thus irreducible. Conversely,
suppose $P$ is monic, irreducible, and of degree $>1$. Then $P$ is the minimum polynomial over $K$ of $a+b\imag$ for some $a,b\in K$ with $b\ne 0$, so 
$$P\ =\ \big(X-(a+b\imag)\big)\cdot\big(X-(a-b\imag)\big)\ =\ (X-a)^2 + b^2.$$ This proves (i). 
As to (ii), we can reduce to the case that $P(a)< 0 < P(b)$ with $a< b$ in $K$; it is enough to
show that then $P(x)=0$ for some $x\in K$ with $a<x<b$. The existence of such an $x$ 
follows easily from (i) by
factoring $P$ into linear factors, and quadratic factors taking only values $>0$.  
\end{proof}

\begin{cor}\label{cor:propAS} Let $K$ be a real closed field and $K(x)$ a field extension of $K$ with~$x$ transcendental over $K$. Let $A$ be a cut in $K$. Then there is a unique ordering of~$K(x)$ that makes it an ordered field extension of $K$ for which $x$ realizes the cut $A$.
\end{cor}
\begin{proof} Let $P\in K[X]$ be monic. Then $P(X)=Q(X)\prod_{i=1}^n(X-a_i)$
where $a_1,\dots, a_n\in K$ and $Q(X)$ is a product of monic irreducible quadratic polynomials in~$K[X]$ as described in Proposition~\ref{prop:AS}(i). Hence for any ordering of $K(x)$ with the indicated properties we have: $P(x)>0$ iff the number of 
$i\in \{1,\dots,n\}$ such that $a_i\notin A$ is even; so there can be at most one such ordering. This also shows how to define such an ordering; alternatively, the existence of such an ordering follows by a routine model-theoretic compactness argument.
\end{proof}

\noindent
Now let $K$ be an ordered field. Then a {\bf real closure} of $K$ is a real closed algebraic field extension of $K$ whose ordering extends the ordering of $K$. Here is the key result on this notion:

\index{ordered field!real closure}
\index{real closure of an ordered field}

\begin{theorem}[Artin \& Schreier] \label{thm:AS, 2} $K$ has a real closure. 
If $K'$ is a real closure of~$K$,
then every ordered field embedding $K\to L$ into a real closed field $L$ has a unique extension to an ordered field embedding $K'\to L$.
\end{theorem}

\noindent
Therefore, if $K_1,\ K_2$ are real closures of $K$, then there is a unique isomorphism $K_{1}\xrightarrow{\cong} K_{2}$ over $K$. Thus we can speak of 
\textit{the}\/ real closure of $K$, denoted by $K^{\operatorname{rc}}$.

Note that by Lemma~\ref{boundrealalg} there is for each
$a\in K^{\operatorname{rc}}$ an element $b\in K^{>0}$ such that $|a|\le b$.
It is tempting to jump to the conclusion that $K$ is dense in  
$K^{\operatorname{rc}}$ (with respect to the order topology), and this jump
is indeed a notorious source of error in
the subject. A counterexample is the ordered field $K=\R(x)$ with $x>\R$:
in its real closure the interval $(\sqrt{x}-1, \sqrt{x}+1)$ contains
no element of $K$.  

Suppose $K$ is given as an ordered subfield of the real closed 
(ordered) field~$F$.
Then the algebraic closure $\{a\in F: \text{$a$ is algebraic over $K$}\}$ of $K$ in $F$ is a real closure of $K$, by Corollary~\ref{cor:AS}, and is called 
the {\bf real closure of $K$ in $F$}; it is clearly the only field extension of $K$ inside $F$ that is a real closure of $K$.

\subsection*{Convex valuations}
{\em In this subsection $K$ is an ordered field}.

\begin{lemma}
Let $A$ be a subring of $K$.
\begin{enumerate}
\item[\textup{(i)}] The convex hull of $A$ in $K$ is a subring of $K$;
\item[\textup{(ii)}] $A$ is convex in $K$ iff $[0,1]\subseteq A$; and
\item[\textup{(iii)}] if $A$ is convex, then $A$ is a valuation ring of $K$.
\end{enumerate}
\end{lemma}
\begin{proof}
Part (i) is clear. In (ii), $[0,1]\subseteq A$ is clearly necessary for convexity of~$A$; conversely, if $[0,1]\subseteq A$ and $a\in A$, $x\in K$ with $0<x<a$, then $xa^{-1}\in A$, hence $x=xa^{-1}\cdot a\in A$.
For (iii), suppose $A$ is convex.  Let $x\in K^\times$; if $\abs{x}\leq 1$, then $x\in A$, and if $\abs{x}>1$ then $\abs{x^{-1}}<1$, so
$x^{-1}\in A$. Thus $A$ is a valuation ring of~$K$.
\end{proof}

\begin{lemma}\label{lem:convex valuation ring}
Let $\mathcal O$ be a valuation ring of $K$. 
The following are equivalent:
\begin{enumerate}
\item[\textup{(i)}] $\mathcal O$ is convex;
\item[\textup{(ii)}] $\smallo$ is convex;
\item[\textup{(iii)}] $\smallo\subseteq (-1,1)$;
\item[\textup{(iv)}] $|a|< 1/n$ for all $a\in \smallo$ and $n\ge 1$;
\item[\textup{(v)}] $\operatorname{conv}(\Q)\subseteq\mathcal O$.
\end{enumerate}
\end{lemma}
\begin{proof}
For (i)~$\Rightarrow$~(ii), assume (i). Let $a,x\in K$, $0<x<a$ and $a\in\smallo$. From $a^{-1}\notin\mathcal O$ and $0<a^{-1}<x^{-1}$ we get $x^{-1}\notin\mathcal O$, so $x\in\smallo$.
The implications (ii)~$\Rightarrow$~(iii) and (iii)~$\Rightarrow$~(iv) are obvious. For (iv)~$\Rightarrow$~(v), assume (iv), let $a\in\operatorname{conv}(\Q)$, $a>0$, and take $n\ge 1$ with $a\le n$. Then $a^{-1} \ge 1/n$, so $a^{-1}\notin\smallo$, and thus $a\in \mathcal O$. Finally, (v)~$\Rightarrow$~(i) follows from part (ii) of the previous lemma. 
\end{proof} 

\noindent
We say that a valuation $v$ on $K$ is {\bf convex} if its valuation ring 
$\mathcal O$ satisfies the equivalent conditions in the previous lemma. Thus a valuation on the ordered field~$K$ is convex iff it is convex as a valuation on the ordered additive group of $K$ as defined in Section~\ref{sec:oag}.
In terms of the dominance relation $\preceq$ associated to $v$:
$$  \text{$v$ is convex}\ \Longleftrightarrow\ \text{ for all $x,y\in K$  with 
$\abs{x}\leq \abs{y}$  we have   $x\preceq y$.}$$
Suppose $v$ is a convex valuation on $K$. Then $\smallo$ is a convex subgroup
of the ordered additive group $\mathcal{O}$, and the resulting ordering on the quotient group 
$\mathcal{O}/\smallo$ makes the residue field an ordered field. (Whenever a convex valuation on
an ordered field is given we regard the residue field as an ordered field in this way.)
If $\dot v$ is a coarsening of $v$ by a convex subgroup of the value group of $v$, then $\dot v$ is also convex, and so is the valuation $v$ on the ordered residue field $\dot K$ of $\dot v$.
By Lemma~\ref{lem:Massaza}, if $v$ is nontrivial, then the $v$-topology on $K$ coincides with the order topology on $K$.

\index{valuation!convex}

\medskip
\noindent
Every ordered field carries a canonical convex valuation:

\begin{example} The standard valuation $v\colon K \to [K]$ with the reverse 
ordering on $[K]$
is given by
$vx=[x]$. It is not just a valuation of the ordered additive group of $K$, 
but even a convex
valuation $K^\times \to [K^\times]$ on the ordered field $K$, where 
$[K^\times]$ is made an 
ordered abelian (additive) group
by $[x]+[y]=[x\cdot y]$ for  $x,y\in K^{\times}$ and  $0=[1]$.
The valuation ring of $v$ is $\mathcal O=\operatorname{conv}(\Q)$ of $K$.
\end{example}

\noindent
The ordered field analogue of Lemma~\ref{lem:uniqueness of ordering} follows from that lemma and its proof:

\begin{cor}\label{lem:ordering immediate ext}
Let $\mathcal O$ be a convex subring of $K$, and $(L,\mathcal O_L)$ 
an immediate valued field extension of $(K,\mathcal O)$. Then just one ordering on $L$ makes $L$ an ordered field extension of $K$ such that
 $\mathcal O_L$ is convex. Moreover, $\mathcal O_L$ is the convex hull of $\mathcal O$ with respect to this ordering on $L$. 
\end{cor}

\noindent
%This follows from Lemma~\ref{lem:uniqueness of ordering} and its proof.
Here is another useful extension result:

\begin{lemma}\label{extror} Let $\mathcal{O}$ be a convex subring of $K$ and $(K(y), \mathcal{O}_y)$ a valued field extension of $(K, \mathcal{O})$ such that 
$nvy\notin \Gamma$ for all $n\ge 1$. Then just one ordering on~$K(y)$ makes
it an ordered field extension of $K$ such that $y>0$ and $\mathcal{O}_y$ is convex. 
\end{lemma}
\begin{proof} By Corollary~\ref{cor:ef inequ, 2} and Lemma~\ref{lem:lift value group ext}, $y$ is transcendental over $K$, the value group of
$(K(y), \mathcal{O}_y)$ is $\Gamma\oplus \Z vy$ (internal direct sum), 
and $\res(K)=\res(K(y))$. Let a field ordering on $K(y)$ be given with $y>0$ that extends the ordering of $K$ and with respect to which
the valuation of $K(y)$ is convex. Let $f\in K(y)^\times$. Then
$f=y^kgu$ with $k\in \Z,\ g\in K^{>}$, and $u\asymp 1$ in $K(y)$, and
so $\res(u)\in \res(K)^{\times}$, hence 
$$\text{$f>0$ in $K(y)$}\ \Longleftrightarrow\ \text{$\res(u)>0$ in $\res(K)$.}$$  This equivalence shows that there can only be one such
ordering. It is routine to check that this equivalence also yields a definition of such an ordering.
\end{proof}

\noindent
Of course, Lemma~\ref{extror} goes through with $y<0$ instead of $y>0$.

\begin{lemma}\label{lem:KW}
Every henselian valuation ring of $K$ is convex.
\end{lemma}
\begin{proof}
Suppose $\mathcal O$ is a henselian valuation ring of $K$.
Let $a\in\smallo$. Then for $x\in\mathcal O$ we have $x^2+x+a\equiv x(x+1)\bmod\smallo$,
so by Proposition~\ref{prop:char henselian} the polynomial $P:=X^2+X+a\in\mathcal O[X]$ has two distinct zeros in $\mathcal O$. Hence the discriminant $1-4a$ of $P$ is positive, that is, $a<\frac{1}{4}$.
This also holds for $-a$ in place of $a$, so $|a|<1$.
Thus~$\mathcal O$ is convex, by (i)~$\Leftrightarrow$~(iii) in Lemma~\ref{lem:convex valuation ring}.
\end{proof} 

\begin{lemma}\label{imagconvex} Let $\mathcal O$ be a convex subring of $K$. Then $\mathcal O + \mathcal O\imag$ is the unique valuation ring of $K[\imag]$ that lies over $\mathcal{O}$. The maximal ideal of $\mathcal O + \mathcal O\imag$ is $\smallo + \smallo \imag$.
\end{lemma}
\begin{proof} Let $a,b\in K$ and $a+b\imag \notin\mathcal O + \mathcal O\imag$.
Then $|a|> \mathcal{O}$ or $|b|> \mathcal O$, so 
$$\frac{1}{a+b\imag}\ =\ \frac{a}{a^2+b^2} -\frac{b}{a^2+b^2}\imag\  \in\ \smallo + \smallo\imag\ \subseteq\ \mathcal O + \mathcal O\imag.$$
Thus $\mathcal O + \mathcal O\imag$ is a valuation ring of $K[\imag]$ lying over $\mathcal{O}$. Since $\imag$ is integral over $\mathcal O$, any valuation ring of $K[\imag]$ lying over $\mathcal{O}$ includes $\mathcal O + \mathcal O\imag$
and thus equals  $\mathcal O + \mathcal O\imag$.
\end{proof} 

\noindent
We have already seen that if a valued field is algebraically closed, then its value group is divisible and its residue field is also algebraically closed (Corollary~\ref{cor:acvf}). Here is an analogue for real closed valued fields:

\begin{theorem} \label{thm:KW}
Suppose $K$ is real closed, and $K$ is equipped with a valuation ring~$\mathcal O$ of $K$. Then 
the value group $\Gamma$ of $K$ is divisible, and $\res(K)$ is either real closed or algebraically closed. Moreover,
the following are equivalent:
\begin{enumerate}
\item[\textup{(i)}] $\res(K)$ is real closed;
\item[\textup{(ii)}] $K$ is henselian;
\item[\textup{(iii)}] $\mathcal O$ is convex.
\end{enumerate}
\end{theorem}
\begin{proof}
For every $a\in K^{>}$ and $n\ge 1$ the polynomial $X^n-a$ has a zero in $K$, so~$\Gamma$ is divisible.
Equip  the algebraic closure $K^{\operatorname{a}}=K(\imag)$ of $K$,  $\imag^2=-1$,  with a valuation ring of $K^{\operatorname{a}}$ lying over $\mathcal O$. Then $\res(K^{\operatorname{a}})$ is an algebraic closure of $\res(K)$ and 
$$\big[\!\res(K^{\operatorname{a}}):\res(K)\big]\  
\leq\ [K^{\operatorname{a}}:K]\ =\ 2.$$
Hence either $\res(K^{\operatorname{a}})=\res(K)$, in which case $\res(K)$ is algebraically closed, or  
$\big[\!\res(K^{\operatorname{a}}):\res(K)\big]=2$, in which case $\res(K)$ is real closed (by Theorem~\ref{thm:AS}).

If $\res(K)$ is real closed, then the valued field $K$ is algebraically maximal, and hence henselian by
Corollary~\ref{cor:alg max implies henselian}.
Thus (i)~$\Rightarrow$~(ii) holds, and (ii)~$\Rightarrow$~(iii) follows from Lemma~\ref{lem:KW}. For (iii)~$\Rightarrow$~(i), note that if $\mathcal O$ is convex, then $\res(K)$ is orderable, so not algebraically closed, and hence
real closed.
\end{proof}

\subsection*{Ordered fields with a convex valuation}
{\em In this subsection $K$ is an ordered field with a convex valuation ring
$\mathcal{O}$ of $K$}, so $K$ is an ordered and a valued field.  
This includes the case $\mathcal{O}=K$, where the corresponding valuation is trivial; we can then make it nontrivial by considering an ordered field extension $K(x)$ with $x>K$ as in Corollary~\ref{cor:boundrealalg}, and taking
the convex hull~$\mathcal{O}_x$ of~$K$ in $K(x)$: then 
$(K(x), \mathcal{O}_x)$ is an ordered and valued field extension of $(K, \mathcal{O})$ for
$\mathcal{O}=K$, with $\mathcal{O}_x\ne K(x)$. 

\begin{cor}\label{convexlift}
Suppose $K$ is real closed and $C$ is a maximal subfield of $\mathcal O$. Then~$C$ is a lift of $\res(K)$, and $\mathcal O=\operatorname{conv}(C)$.
\end{cor}
\begin{proof} Proposition~\ref{prop:lift} and Theorem~\ref{thm:KW}
imply that $C$ is a lift of~$\res(K)$. For $x\in\mathcal O$ we have $c\in C$ with $x-c\in \smallo$, so $c-1<x<c+1$,  thus $x\in \operatorname{conv}(C)$. 
%  By Lemma~\ref{lem:max subfield} there is some monic  $P=\sum_i P_iX^i\in C[X]$ of degree $d>0$ such that %$\overline{P}(\overline{x})=0$. Thus $P(x)\in\smallo\subseteq\mathfrak m_A$, and if we had $x\notin A$ %then $x^{-1}\in\mathfrak m_A$,  hence $1+(P_{n-1}x^{-1}+\cdots+P_0x^{-d})=x^{-d}P(x)\in\mathfrak m_A$, a %contradiction.
\end{proof}

\begin{cor}\label{cvrc}
There exists a unique convex valuation ring of the real closure~$K^{\operatorname{rc}}$ of $K$ lying over the valuation ring $\mathcal O$ of $K$. Equipping $K^{\operatorname{rc}}$ with this valuation ring we have 
$\Gamma_{K^{\operatorname{rc}}}=\Q\Gamma$ and
$\res(K^{\operatorname{rc}})=\res(K)^{\operatorname{rc}}$.
\end{cor}
\begin{proof} Let  $\mathcal{O}^{\operatorname{rc}}$ be 
the convex hull of $\mathcal O$ in $K^{\operatorname{rc}}$. Then $\mathcal{O}^{\operatorname{rc}}$ is a convex valuation ring of~$K^{\operatorname{rc}}$ lying over $\mathcal O$. By the remarks following the
proof of Proposition~\ref{prop:Krull, 3},  $\mathcal{O}^{\operatorname{rc}}$ is the only convex valuation ring of $K^{\operatorname{rc}}$ lying over $\mathcal O$.
Turn $K^{\operatorname{rc}}$ into a valued field with valuation ring
$\mathcal{O}^{\operatorname{rc}}$. Then  $\Gamma_{K^{\operatorname{rc}}}$ is divisible, so $\Gamma_{K^{\operatorname{rc}}}=\Q\Gamma$.
Since $\res(K^{\operatorname{rc}})$ is an algebraic ordered field extension of $\res(K)$ and $\res(K^{\operatorname{rc}})$ is real closed (by Theorem~\ref{thm:KW}), we have $\res(K^{\operatorname{rc}})=\res(K)^{\operatorname{rc}}$.
\end{proof} 

\begin{cor}\label{cor:Prestel} For our valued field $K$ we have:
$$\text{$K$ is real closed }\ \Longleftrightarrow\  
\text{$K$ is henselian,
$\res(K)$ is real closed, and
$\Gamma$ is divisible}.$$
\end{cor}
\begin{proof}
The direction $\Rightarrow$ follows from Theorem~\ref{thm:KW}.
Conversely, suppose the conditions on the right side are satisfied. Corollary~\ref{cvrc} gives a 
valuation ring of $K^{\operatorname{rc}}$ making $K^{\operatorname{rc}}\supseteq K$  an immediate algebraic extension of valued fields. It now follows from Corollary~\ref{cor:alg max equals henselian} that $K=K^{\operatorname{rc}}$ is real closed.
\end{proof}

\begin{example}
Let $C$ be an ordered field and $\mathfrak M$ a (multiplicative) ordered abelian group. Viewing $C[[\mathfrak M]]$ as a Hahn product, its Hahn ordering 
is given by 
$$f>0\ \Longleftrightarrow\ f_{\frak{m}}>0\qquad   (f\in C[[\mathfrak M]]^{\ne},\ \frak{m}=\frak{d}(f)).$$
Then $C[[\mathfrak M]]$ with its Hahn ordering is called an {\bf ordered Hahn field\/}; it is indeed an ordered field containing 
$C$ as an ordered subfield and with $\mathfrak{M}$ as an ordered subgroup of
$C[[\mathfrak M]]^{>}$. The 
 Hahn valuation on  $C[[\mathfrak M]]$ has valuation ring $\operatorname{conv}(C)$ with respect to the
Hahn ordering, and is thus a convex valuation. Therefore:
$$\text{$C[[\mathfrak M]]$ is real closed } \Longleftrightarrow\  \text{$C$ is real closed and $\mathfrak M$ is divisible.}$$
Thus the Hahn field $\R\(( t^{\Q}\)) $ is real closed. The field $\operatorname{P}(C)$ of Puiseux series over $C$ (Example~\ref{ex:Puiseux}) is an ordered subfield of the ordered
Hahn field $C \(( t^\Q \)) $.
If $C$ is real closed, then $\operatorname{P}(C)$ is the real
closure of its ordered subfield $C\(( t \)) $ in $C \(( t^\Q \)) $.
\end{example}

\begin{cor}\label{cor:completion of rcf}
Suppose $K$ is real closed. Then the completion $K^{\operatorname{c}}$ of the valued field $K$ is real closed, and the valuation of  $K^{\operatorname{c}}$ is convex.
\end{cor}
\begin{proof} 
Since $K\subseteq K^{\operatorname{c}}$ is an immediate extension, the residue field of $K^{\operatorname{c}}$ is real closed, by Theorem~\ref{thm:KW}, hence 
$K^{\operatorname{c}}$ is not algebraically closed, by Corollary~\ref{cor:acvf}.
Extend the valuation of $K$ to the algebraic closure $K(\imag)$, $\imag^2=-1$, of $K$.
By Corollary~\ref{cor:completion of acf}, 
$K(\imag)^{\operatorname{c}}$ is algebraically closed, and by Corollary~\ref{cor:completion finite ext}, $K(\imag)^{\operatorname{c}}=K^{\operatorname{c}}(\imag)$. 
Hence $K^{\operatorname{c}}$ is real closed, and by Theorem~\ref{thm:KW}, its valuation is convex.
\end{proof}

\subsection*{Completion of ordered fields}
{\em In this subsection $K$ is an ordered field}.
Call an ordered field extension $L\supseteq K$ {\bf dense} if $K$ is dense in $L$, in the order topology on~$L$.
We have an analogue of Theorem~\ref{thm:completion valued fields} for 
\textit{ordered}\/ fields:

\index{extension!ordered fields!dense}
\index{dense extension!ordered fields}

\begin{theorem}\label{thm:completion ordered fields}
There is a dense ordered field extension $K^{\operatorname{d}}\supseteq K$ such that any dense ordered field extension $L\supseteq K$ embeds uniquely over $K$ into 
$ K^{\operatorname{d}}$.
\end{theorem} 
\begin{proof} Let $K^{\operatorname{d}}$ be the completion of the ordered additive group of $K$; see Section~\ref{sec:oag}. It is easy to check
that Lemma~\ref{lem:multiplication of c-sequences}
goes through for c-sequences in the ordered field~$K$, and
from this we obtain that the multiplication $(x,y)\mapsto x\cdot y$ and inversion $x\mapsto 1/x$ ($x\neq 0$) on $K$ have unique extensions to continuous maps
$$K^{\operatorname{d}}\times K^{\operatorname{d}}\to  K^{\operatorname{d}}, 
\qquad (K^{\operatorname{d}})^{\neq}\to  K^{\operatorname{d}},$$ with the product topology on $K^{\operatorname{d}}\times K^{\operatorname{d}}$,
and that $K^{\operatorname{d}}$ is an ordered field extension of~$K$ 
with the first map as multiplication and the second map as inversion. 
It is now routine to check that this
ordered field extension has the desired properties.
\end{proof}

\noindent
The properties of the ordered field extension $K^{\operatorname{d}}$ of $K$ 
postulated in Theorem~\ref{thm:completion ordered fields} determine 
$K^{\operatorname{d}}$ up to a  unique (ordered field) isomorphism over $K$. 
We call $K^{\operatorname{d}}$  the {\bf completion of $K$.} Note that by 
construction,  $K^{\operatorname{d}}$ is indeed complete: every c-sequence 
in  $K^{\operatorname{d}}$ converges in  $K^{\operatorname{d}}$. When is 
$K^{\operatorname{d}}$ real closed?

\index{completion!ordered field}
\index{ordered field!completion}
\nomenclature[K]{$K^{\operatorname{d}}$}{completion of the ordered field $K$}

\begin{prop}\label{drc} $K^{\operatorname{d}}$ is real closed if and only if $K$ is dense in its real closure.
\end{prop}

\noindent
Consider first the case that  $K$ is 
archimedean. Then $K$ is isomorphic to a unique ordered subfield of $\R$,
and identifying $K$ with this subfield we can take $K^{\operatorname{d}}=\R$,
and~$K^{\operatorname{rc}}$ as a subfield of $\R$.
Thus $K^{\operatorname{d}}$ is real closed and  $K$ is dense in 
$K^{\operatorname{rc}}$. In general, if~$K^{\operatorname{d}}$ is real closed, 
then clearly $K$ is dense in $K^{\operatorname{rc}}$. For the converse, we can
assume~$K$ is not archimedean, and so $K$ has a convex subring  
$\mathcal O\ne K$.

\medskip\noindent
Accordingly, consider $K$ below as equipped with a convex valuation 
ring $\mathcal O\ne K$ of~$K$. 
Let $K^{\operatorname{c}}$ be the completion of the valued field $K$. Then 
$K^{\operatorname{c}}\supseteq K$ is an immediate valued field extension, so by Lemma~\ref{lem:ordering immediate ext} we can take a unique ordering on $K^{\operatorname{c}}$ making~$K^{\operatorname{c}}$ an ordered field extension of $K$ such that
the valuation ring of~$K^{\operatorname{c}}$ is convex with respect to
this ordering. The valuation topology
and the order topology on $K^{\operatorname{c}}$ coincide, 
by Lemma~\ref{lem:Massaza}, so with this ordering 
$K^{\operatorname{c}}$ is also complete as an ordered field, and $K$ is dense in
$K^{\operatorname{c}}$ in the order topology. Thus: 

\begin{lemma}
There is a unique isomorphism $K^{\operatorname{c}}\to K^{\operatorname{d}}$ of ordered fields which is the identity on $K$.
\end{lemma}
 
\noindent
In view of Corollary~\ref{cor:completion of rcf}, this yields:

\begin{cor} \label{cor:completion of rcf, 2}
If $K$ is real closed, then so is $K^{\operatorname{d}}$.
\end{cor}

\noindent
To conclude the proof of Proposition~\ref{drc} it remains to use this 
corollary and to note that if $K$ is dense in an ordered field extension $L$, 
then we can take $K^{\operatorname{d}}=L^{\operatorname{d}}$.   

\begin{example}
Let $C$ be an ordered field. By Example~\ref{ex:LC} 
the completion in the Hahn field $C\(( t^\Q \)) =C[[x^\Q]]$ ($t=x^{-1}$) of its subfield $C(t^\Q)$ is
$$L\ :=\ \big\{ f\in C[[x^\Q]]:\  \text{$\supp f|_{x^q}$ is finite for all $q\in\Q$} \big\}.$$
Now consider $L$ as an ordered subfield of the ordered Hahn field 
$C\(( t^\Q \)) $. Then $L$ is an ordered field extension of
the ordered subfield $C(t^\Q)$ of $C\(( t^\Q \)) $, and is as such also a completion of this ordered field $C(t^\Q)$. 
If $C$ is real closed, then 
so is $L$, by
Corollaries~\ref{cor:henselian complete rank 1} and~\ref{cor:Prestel}.
\end{example}

\subsection*{Notes and comments}
%Lemma~\ref{boundrealalg} is often attributed 
%to Cauchy~\cite{Cauchy}.
The notion of a real closed field and Theorems~\ref{thm:AS} and \ref{thm:AS, 2} are due to Artin and Schreier~\cite{Artin-Schreier,Artin-Schreier-2}. 
For proofs of these theorems, see for example~\cite{KS} or~\cite{Prestel}; for the equivalence of (iii) and (iv) in Theorem~\ref{thm:AS}, see also Leicht~\cite{Leicht}. 
Lemma~\ref{lem:hoelder, 2} is from~\cite{Hoelder}.
Theorem~\ref{thm:KW} is in Knebusch and Wright~\cite{KW}, and Corollary~\ref{cor:Prestel} in Prestel~\cite{Prestel}.
Corollary~\ref{cor:completion of rcf, 2} has been noticed by various authors; see \cite{Baer,Hauschild,Scott}.

%% file: mt-3-6.tex
\section{Some Model Theory of Valued Fields}  
\label{sec:modth val fields}

\noindent
In this section we assume familiarity with Appendix~\ref{app:modth}. We establish here quantifier elimination for
\textit{algebraically closed fields with a nontrivial valuation}\/ and \textit{real closed fields
with a nontrivial convex valuation.}\/
These well-known results foreshadow
the deeper elimination theorem about the valued differential field 
$\mathbb{T}$
in Chapter~\ref{ch:QE}. We also need Proposition~\ref{indsemi} below in Section~\ref{sec:embth}. 
 
\subsection*{Algebraically closed valued fields} We augment the language
$\{ 0, 1, {-}, {+}, {\,\cdot\,}\}$ of rings by a binary relation symbol $\preceq$ to obtain the language $\mathcal L_{\preceq}$.\nomenclature[Ba]{$\mathcal L_{\preceq}$}{language of integral domains with a dominance relation} Let $\operatorname{ACVF}$ be the $\mathcal L_{\preceq}$-theory whose models are
the structures $(K,{\preceq})$ where~$K$ is an algebraically closed field and $\preceq$
is a nontrivial dominance relation on $K$.\nomenclature[Bf]{$\operatorname{ACVF}$}{theory of algebraically closed fields with nontrivial dominance relation in the language $\mathcal L_{\preceq}$} Here a dominance relation on a field $K$ is said to be {\bf trivial} if its corresponding valuation
ring is $K$. %, equivalently if the value group of its associated valuation is $\{0\}$.

\begin{theorem}\label{thm:ACVF QE}
$\operatorname{ACVF}$ has \textup{QE}.
\end{theorem}

\noindent
First some remarks about dominance relations.
Let $R$ be an integral domain. We define a {\bf dominance relation} on $R$ \index{dominance relation} to be a binary relation $\preceq$ on $R$ such that conditions (D1)--(D6) in Section~\ref{sec:valued fields} hold for all $f,g,h\in R$.
The {\bf trivial} \index{dominance relation!trivial}\index{trivial!dominance relation} dominance relation 
$\preceq_{\operatorname{t}}$  on $R$ is the one with
$r \preceq_{\operatorname{t}} s$ for all $r,s\in R$ with~$s\neq 0$.
For any dominance relation $\preceq$ on $R$ there is a unique dominance relation $\preceq_F$
on $F=\Frac(R)$ such that $(R,{\preceq})\subseteq (F,{\preceq_F})$; it is given by
$$\frac{r_1}{s} \preceq_F \frac{r_2}{s}\quad\Longleftrightarrow\quad r_1\preceq r_2
\qquad (r_1,r_2,s\in R,\ s\neq 0).$$

\begin{exampleNumbered}\label{ex:dr on Z}
The nontrivial dominance relations on $\Z$ are exactly the dominance relations
$\preceq_p$, where $p$ is a prime number:
$$a \preceq_p b  \quad\Longleftrightarrow\quad  \text{for all $n$: $b\in p^n\Z\Rightarrow a\in p^n\Z$.}$$ 
\end{exampleNumbered}

\noindent
Every substructure of a model of $\operatorname{ACVF}$ is a pair
$(R,{\preceq})$ with $R$ an integral domain and~$\preceq$ a dominance relation on $R$.
Conversely, for any integral domain $R$ with a dominance relation $\preceq$ on it, $(R,{\preceq})$ is a substructure of a model
of $\operatorname{ACVF}$: first extend~$(R,{\preceq})$ to  $(F,{\preceq_F})$ as above; then extend $\preceq_F$ to a dominance relation on the algebraic closure of $F$; in case a valuation (in the form of a dominance relation) is trivial, adjoin a transcendental 
to make it nontrivial.

\medskip
\noindent
Let $(R,{\preceq})$ be an integral domain with a dominance relation on it. Let
$$i\ \colon\ (R,{\preceq})\to (K,{\preceq_K})$$ be an embedding into an algebraically closed field $K$
with dominance relation~$\preceq_K$ on~$K$. Let $\preceq^{\alg}$ be a dominance relation on
the algebraic closure $F^{\alg}$ of $F=\Frac(R)$ such that $(R,{\preceq})\subseteq (F^{\alg},{\preceq^{\alg}})$.
Then by Corollary~\ref{cor:Krull, 2} we can extend~$i$ to an embedding $(F^{\alg},{\preceq_{F^{\alg}}})\to (K,{\preceq_K})$. 

\begin{proof}[Proof of Theorem~\ref{thm:ACVF QE}]
By the remarks preceding this proof and~\ref{prop:QE test, 2} it suffices to show the following:
Let $E$, $F$ be nontrivially valued algebraically closed fields such that $F$ is $\abs{K}^+$-saturated, where $K$ is a proper algebraically closed subfield of $E$; view subfields of $E$ as valued subfields of $E$, and 
let $i\colon K\to F$ be a valued field
embedding. Then there exists $x\in E\setminus K$ and an extension of~$i$ to a valued field embedding
$j\colon K(x)\to F$. 

To find such $x$ and $j$ we distinguish three cases.
To simplify notation we identify the valued field $K$ with the valued field $iK$ via $i$.

\case[1]{$\res(K)\neq\res(E)$.}
Take $x\in \mathcal O_E$ such that $\overline{x}\notin\res(K)$.
Since $\res(K)$ is algebraically closed, $\overline{x}$ is transcendental over $\res(K)$.
Also $x\notin K$, so $x$ is transcendental over $K$. By the saturation assumption
on $F$ we can find $y\in\mathcal O_F$ with $\overline{y}\notin\res(K)$.
So $\overline{y}$ is transcendental over $\res(K)$ and $y$ is transcendental
over~$K$. Then the field embedding $j\colon K(x)\to F$ over $K$ with $j(x)=y$ is a valued
field embedding by the uniqueness part of Lemma~\ref{lem:gauss}.

\case[2]{$\Gamma\neq\Gamma_E$.} 
Note that $\Gamma$ is divisible.
Take any $\alpha\in\Gamma_E\setminus\Gamma$.
Since $F$ is $|K|^+$-saturated, so is $\Gamma_F$ as an ordered set.
Also, $\Gamma_F\neq \{0\}$. Hence we can take $\beta\in\Gamma_F$ realizing the same
cut in $\Gamma$ as $\alpha$. 
So we have an isomorphism $\Gamma+\Z\alpha\to\Gamma+\Z\beta$ of ordered abelian groups
over $\Gamma$ which sends $\alpha$ to $\beta$.
Take $x\in E^\times$, $y\in F^\times$ with $vx=\alpha$, $vy=\beta$.
Then $x\notin K$, so $x$ is transcendental over~$K$; likewise,~$y$ is transcendental
over~$K$. By the uniqueness part of Lemma~\ref{lem:lift value group ext}, the field embedding 
$j\colon K(x)\to F$ over $K$ with $j(x)=y$ is even a valued field embedding.

\case[3]{$\res(K)=\res(E)$ and $\Gamma=\Gamma_E$.}
Take any $x\in E\setminus K$. Then the valuation on~$K(x)$ is uniquely determined
by the valuation on $K$ and by the map 
$$a\mapsto v(x-a)\ \colon\ K\to\Gamma,$$
since each monic $f\in K[x]$ factors as $f=\prod_{i=1}^n (x-a_i)$ with all $a_i\in K$,  
so $vf=\sum_i v(x-a_i)$. It follows that for $y\in F\setminus K$ with
$v(x-a)=v(y-a)$ for all $a\in K$, the field embedding
$K(x)\to F$ over $K$ that sends $x$ to $y$ is a valued field embedding.
Such an element~$y$ exists by saturation and the next general lemma.
\end{proof}

\begin{lemma}
Let $K\subseteq L$ be a valued field extension such that $\res K=\res L$.
Let $a_1,\dots,a_n\in K$, $n\geq 1$, and let $x\in L\setminus K$ be such that
$v(x-a_i)\in\Gamma$ for $i=1,\dots,n$. Then there exists $a\in K$ such that 
$v(x-a_i)=v(a-a_i)$ for $i=1,\dots,n$.
\end{lemma}
\begin{proof}
Any $a\in K$ such that $v(a-x)>v(a-a_i)$ for $i=1,\dots,n$ has the desired property.
We may assume $v(x-a_1)\geq v(x-a_i)$ for $i=2,\dots,n$. Since $v(x-a_1)\in\Gamma$
we can take $b\in K$ such that $v(x-a_1)=vb$. So $v\big(\frac{x-a_1}{b}\big)=0$,
and since $\res K=\res L$, we have $\frac{x-a_1}{b}=c+\varepsilon$ with $c\in K$,
$c\asymp 1$, $\varepsilon\prec 1$. Then $a=a_1+bc$ works because $x-a=b\varepsilon$ and
$v(b\varepsilon)>v(x-a_i)$.
\end{proof}

\begin{cor}
$\operatorname{ACVF}$ is the model completion of the ${\mathcal L}_{\preceq}$-theory of pairs
$(R, \preceq)$ where $R$ is an integral domain and $\preceq$ is a dominance
relation on $R$. 
\end{cor}

\noindent
For a valued field $K$ with residue field $\k$,
the pair $(\operatorname{char} K,\operatorname{char}\k)$ 
is among the following, where $p$ is a prime number:
\begin{list}{*}{\setlength\leftmargin{2.5em}}
\item[$(0,0)$:] equicharacteristic $0$;
\item[$(0,p)$:] mixed characteristic $p$;
\item[$(p,p)$:] equicharacteristic $p$.
\end{list}
Each of these actually occurs:
if $\k$ is a field of characteristic $p$, with $p=0$ or $p$ a prime number, then 
the Hahn fields $\k\(( t^\Gamma\)) $ have equicharacteristic $p$; if $p$ is a prime number, then the unique dominance relation $\preceq$ on $\Q$ such that
$(\Z, {\preceq_p})\subseteq (\Q, {\preceq})$ yields a valued field of
mixed characteristic $p$. 
%For each prime number $p$, the valued field $\Q_p$ of $p$-adic numbers is of %mixed characteristic $p$.

\index{valued field!characteristic}
\index{characteristic of a valued field}

The
{\bf characteristic} of a valued field $K$ is $(\operatorname{char} K,\operatorname{char}\k)$; it equals the characteristic of any valued field
extension of $K$.
Let $(m,n)$ be the characteristic of some valued field; define
$\operatorname{ACVF}_{(m,n)}$ as the $\mathcal L_{\preceq}$-theory whose models are the structures
$(K, \preceq)\models\operatorname{ACVF}$ that are of characteristic $(m,n)$
as a valued field. 

\begin{cor}
$\operatorname{ACVF}_{(m,n)}$ is complete. 
\end{cor}
\begin{proof} Let $p$ range over prime numbers. 
Construing valued fields as fields with a dominance relation, the classification of dominance relations on $\Z$ 
from~\ref{ex:dr on Z} gives: $(\Z,{\preceq_{\operatorname{t}}})$ embeds into every
valued field of characteristic $(0,0)$, and $(\Z,{\preceq_p})$ into every
valued field of characteristic $(0,p)$. Clearly, $(\mathbb F_p,{\preceq_{\operatorname{t}}})$ embeds into every valued field of
characteristic $(p,p)$. Now use Corollary~\ref{cor:QE and sentences}.
\end{proof}

\subsection*{Real closed valued fields} An {\bf ordered integral domain\/} is an integral domain~$R$ with a (total) ordering $\le$ of $R$ such that for all $x,y,z\in R$,
$$x\le y\ \Rightarrow\ x+z\le y+z, \quad x\le y\ \&\ z\ge 0\ \Rightarrow\ xz\le yz.$$
Given an ordered integral domain $(R, {\le})$ there is a unique field ordering~$\le_F$ of its fraction field $F$ such that $(R, {\le}) \subseteq (F, {\le_F})$; we call $(F, {\le_F})$ the {\bf ordered fraction field\/} of $(R, {\le})$. 
We augment the language
$\{ 0, 1, {-}, {+}, {\,\cdot\,}\}$ of rings by a binary relation symbol $\le$ to obtain the language $\mathcal{L}_{\operatorname{OR}}$ of ordered rings, and we construe ordered integral domains as $\mathcal{L}_{\operatorname{OR}}$-structures in the natural way.

We augment $\mathcal L_{\operatorname{OR}}$ by a binary relation symbol~$\preceq$
to obtain the language $\mathcal L_{\operatorname{OR},\preceq}$,
and construe each valued ordered field as an $\mathcal L_{\operatorname{OR},\preceq}$-structure in the obvious way.
A dominance relation $\preceq$ on an ordered integral domain~$R$ is said to be  
{\bf convex} if for all~$r,s\in R$ we have: $0\leq r\leq s\Rightarrow r\preceq s$.
So a dominance relation on an ordered field is convex iff 
its corresponding valuation ring is convex. If $\preceq$ is a convex dominance relation on the ordered integral domain $R$, then $\preceq_F$ as in the  previous subsection is a convex dominance relation on its ordered fraction field $F$. 

\index{dominance relation!convex}
\index{convex!dominance relation}
\nomenclature[Bb]{$\mathcal L_{\operatorname{OR}}$}{language of ordered rings}
%\nomenclature[Bex52]{$\mathcal L_{\operatorname{OR}}$}{language of ordered rings}
\nomenclature[Bc]{$\mathcal L_{\operatorname{OR},\preceq}$}{language of ordered integral domains with a dominance relation}
\nomenclature[Bg]{$\operatorname{RCVF}$}{theory of real closed ordered fields with nontrivial convex dominance relation in the language $\mathcal L_{\operatorname{OR},\preceq}$}

Let $\operatorname{RCVF}$ be the $\mathcal L_{\operatorname{OR},\preceq}$-theory whose models are
the $\mathcal L_{\operatorname{OR},\preceq}$-structures~$(K,{\preceq})$ where~$K$ is a real closed ordered field and $\preceq$
is a nontrivial convex dominance relation on~$K$. The substructures of the models of
$\operatorname{RCVF}$ are exactly the pairs~$(R,{\preceq})$ where $R$ is an ordered integral domain and $\preceq$ is
a convex dominance relation on $R$; this follows easily from the remarks above and Corollary~\ref{cvrc}.  

Let $R$ be an ordered integral domain, $\preceq$ a convex
dominance relation on $R$, and $i\colon (R,{\preceq})\to (K,{\preceq_K})$ an embedding
into a real closed ordered field $K$ with a convex dominance relation $\preceq_K$ on $K$. By Corollary~\ref{cvrc} there is a unique convex dominance relation $\preceq_{F^{\operatorname{rc}}}$
on the real closure $F^{\operatorname{rc}}$ of
the ordered fraction field $F$ of~$R$ with $(R,{\preceq})\subseteq (F^{\operatorname{rc}},{\preceq_{F^{\operatorname{rc}}}})$; it follows that
the unique extension of $i$ to an ordered field embedding $F^{\operatorname{rc}}\to K$ is also an
embedding $(F^{\operatorname{rc}}, {\preceq_{F^{\operatorname{rc}}}}) \to (K,{\preceq_K})$.

\begin{theorem}\label{thm:RCVF QE}
$\operatorname{RCVF}$ has \textup{QE}.
\end{theorem}
\begin{proof}
By the above remarks and ~\ref{prop:QE test, 2} it suffices to show the following:
Let $E, F\models \operatorname{RCVF}$ be such that~$F$ is $\abs{K}^{+}$-saturated where $K$ is a proper real closed subfield of $E$;
view subfields of $E$ as valued ordered subfields of~$E$, and let
$i\colon K\to F$ be an
embedding of ordered valued fields. Then there is an $x\in E\setminus K$ and an extension of~$i$ to an embedding
$j\colon K(x)\to F$ of ordered valued fields. 

To find such $x$ and $j$ we distinguish the same three cases as in the proof of Theorem~\ref{thm:ACVF QE}.
To simplify notation we identify the valued ordered field $K$ with the 
valued ordered field $iK$ via $i$.

\case[1]{$\res(K)\neq\res(E)$.}
Take $x\in \mathcal O_E$ such that $\overline{x}\notin\res(K)$.
Since $\res(K)$ is real closed, $\overline{x}$ is transcendental over $\res(K)$.
Also $x\notin K$, so $x$ is transcendental over~$K$. By the saturation assumption
on $F$ we can find $y\in\mathcal O_F$ such that $\overline{y}$ realizes the same cut
in $\res(K)$ as $\overline{x}$.  
In particular $\overline{y}\notin\res(K)$, 
so $\overline{y}$ is transcendental over~$\res(K)$ and~$y$ is transcendental
over~$K$. 
The field embedding $j\colon K(x)\to F$ over $K$ with $j(x)=y$ is a valued
field embedding by the uniqueness part of Lemma~\ref{lem:gauss}.
It is easy to check that $x\in \mathcal{O}_E$ and $y\in \mathcal{O}_F$ realize the same cut in $K$, so $j$ preserves order by Corollary~\ref{cor:propAS}.

\case[2]{$\Gamma\neq\Gamma_E$.} 
Take any $x\in E^>$ with $\alpha:=vx\notin\Gamma$.
As in the proof of Case~2 of Theorem~\ref{thm:ACVF QE} we get $y\in F^>$ such that $\beta:=vy$ realizes the same cut in $\Gamma$ as~$\alpha$. As in that proof, $x$ 
and $y$ are transcendental over~$K$, and the field embedding 
$j\colon K(x)\to F$
over $K$ with $j(x)=y$ is a valued field embedding. It is easy to check that 
$x$ and $y$ realize the same cut in $K$, so $j$ preserves order as in Case 1.

\case[3]{$\res(K)=\res(E)$ and $\Gamma=\Gamma_E$.}
Take any $x\in E\setminus K$. With $v$ extended to a valuation on $E[\imag]$
we have $v\big(x-(a+b\imag)\big)=\min\big(v(x-a), vb\big)$ for $a,b\in K$ by Lemma~\ref{imagconvex}. Thus we can proceed as in the proof of  Case~3 of Theorem~\ref{thm:ACVF QE} to obtain
a valued field embedding $j\colon K(x)\to F$ over $K$. By Corollary~\ref{lem:ordering immediate ext}, $j$ preserves order.
\end{proof}

\noindent
The separation in three cases in the proofs of Theorem~\ref{thm:ACVF QE} and \ref{thm:RCVF QE}, according to whether the residue field extends, the value group extends, or neither,
is a common feature of many proofs for QE or model completeness of (expansions of) valued fields;
we see this again in proving Theorem~\ref{thm:tamepairs} below, and in establishing QE for~$\mathbb T$ in Chapter~\ref{ch:QE}. The case of immediate extensions
is usually the hardest.

In view of the remarks preceding Theorem~\ref{thm:RCVF QE} we have:

\begin{cor}\label{cor:RCVF QE} $\operatorname{RCVF}$ is the model completion of the 
$\mathcal{L}_{\operatorname{OR}, \preceq}$-theory of ordered
integral domains $R$ equipped with a convex dominance
relation on $R$.
\end{cor} 

\begin{cor}
$\operatorname{RCVF}$ is complete. 
\end{cor}
\begin{proof} Use that the ordered ring of integers with its trivial dominance relation embeds into
every model of $\operatorname{RCVF}$. 
\end{proof}

\begin{cor}\label{cor:RCVF NIP}
$\operatorname{RCVF}$ has \textup{NIP}.
\end{cor}
\begin{proof}
Let $(K,{\preceq})\models \operatorname{RCVF}$ and suppose that the 
relation $R\subseteq K^m\times K^n$ is  $0$-definable in~$(K,{\preceq})$; we need to show that $R$ is dependent. 
Let $K^*$ be a $\abs{K}^+$-saturated elementary extension of the ordered field $K$,
and take $a^*>0$ in $K^*$ such that
$\mathcal O=(-a^*,a^*)_{K^*}\cap K$.
Let $\varphi(x,y)$ be a quantifier-free 
$\mathcal L_{\operatorname{OR}, \preceq}$-formula that defines~$R$ in~$(K,{\preceq})$, with $x=(x_1,\dots, x_m)$,  $y=(y_1,\dots, y_n)$. Boolean combinations
of dependent relations are dependent (Lemma~\ref{lem:VC for boolean combinations}). So we can assume  that~$\varphi(x,y)$ is  of the form ``$P(x,y)\leq Q(x,y)$'' or 
``$P(x,y)\preceq Q(x,y)$'' where $P,Q\in \Z[x,y]$. 
We associate to $\varphi$  a quantifier-free formula $\varphi^*$ in the language of ordered rings:
in the first case we take $\varphi^*:=\varphi$, and in the second case $\varphi^*$ expresses
$$ P(x,y)=Q(x,y)=0 \vee \big(Q(x,y)\neq 0 \ \&\ \abs{P(x,y)}<a^*\abs{Q(x,y)}\big).$$
Let $R^*$ be the subset of $(K^*)^{m+n}$ defined by 
$\varphi^*$. Then  $R^*\cap K^{m+n}=R$. 
Since $R^*$ is dependent by Corollary~\ref{cor:RCF NIP}, so is $R$.
\end{proof}

\subsection*{Tame pairs} A {\bf tame pair} (tacitly: of real closed fields) is a pair $(K,C)$ where $K\models \operatorname{RCVF}$ and $C$ is a real closed subfield of $K$ such that $\mathcal{O}=C+\smallo$, where~$\mathcal{O}$ is the valuation ring of $K$ (corresponding to the distinguished dominance relation~$\preceq$ on $K$); note that then $\mathcal{O}$ is the convex hull of $C$ in $K$ and 
$C$ is a lift of $\res(K)$.
Conversely, if $K\models \operatorname{RCVF}$, then by Proposition~\ref{prop:lift} and Zorn, $\mathcal O$ contains a lift of~$\res(K)$, and for any such lift $C$ we have a tame pair $(K,C)$. 

It is worth noting that a tame pair $(K,C)$ has a {\em standard part map\/} 
$\operatorname{st}\colon\mathcal O\to C$: it assigns to $a\in\mathcal O$ the unique $c\in C$ such that $a-c\in\smallo$; thus $\operatorname{st}\colon\mathcal O\to C$ is a ring morphism, and if $a,b\in \mathcal{O}$, $a\le b$, then $\operatorname{st}(a)\le \operatorname{st}(b)$. 

Let $\mathcal L_{\operatorname{tame}}$ be the language 
$\mathcal{L}_{\operatorname{OR}, \preceq}$ augmented by a unary relation symbol $U$.
We view each tame pair $(K,C)$ as an $\mathcal L_{\operatorname{tame}}$-structure in the natural way, interpreting~$U$ as
the underlying set of $C$. Note that for tame pairs 
$(K,C)$ and $(L,C_L)$ with
$(K,C)\subseteq(L,C_L)$, the standard part map of $(L, C_L)$
extends the standard part map of $(K, C)$. 
We let~$\operatorname{RCF}_{\operatorname{tame}}$  be the $\mathcal L_{\operatorname{tame}}$-theory
whose models are exactly the tame pairs~$(K,C)$.

\index{real closed field!tame pair}
\index{tame pair}
%\index{standard part}
\nomenclature[Bd]{$\mathcal L_{\operatorname{tame}}$}{language of tame pairs}
\nomenclature[Bh]{$\operatorname{RCF}_{\operatorname{tame}}$}{theory of nontrivial tame pairs in the language $\mathcal L_{\operatorname{tame}}$}

\begin{lemma}\label{lem:tamepairs} Assume $\mathbf K=(K,C)$ and $\mathbf E=(E,C_E)$ are models
of~$\operatorname{RCF}_{\operatorname{tame}}$ such that $\mathbf K\subseteq\mathbf E$. Let $e\in C_E$, and let $K(e)^{\operatorname{rc}}$ and $C(e)^{\operatorname{rc}}$ be the real closures of~$K(e)$ and~$C(e)$ in $E$ and $C_E$, respectively. Then $C_E\cap K(e)^{\operatorname{rc}} = C(e)^{\operatorname{rc}}$, and
$$\mathbf K_1\ :=\  \big(K(e)^{\operatorname{rc}}, C(e)^{\operatorname{rc}}\big) \models \operatorname{RCF}_{\operatorname{tame}}, \qquad
\mathbf K\ \subseteq\ \mathbf K_1\ \subseteq\ \mathbf E,$$
where $K(e)^{\operatorname{rc}}$ is construed as a valued
ordered subfield of $E$. 
\end{lemma}
\begin{proof} This is trivial if $e\in C$. Assume 
$e\notin C$, so $e\notin K$. Also,  
$$K_1:= K(e)^{\operatorname{rc}}\models \operatorname{RCVF}, \qquad
C_1:= C(e)^{\operatorname{rc}}\subseteq C_E\cap K_1 \subseteq \mathcal{O}_1 
:= \mathcal{O}_E\cap K_1.$$ 
By
Lemma~\ref{lem:gauss} for $x:=e$, and Corollary~\ref{cvrc}, 
$C_1$ is a lift of the residue field of~$K_1$, so $C_1=C_E\cap K_1$ and $\mathbf K_1\models \operatorname{RCF}_{\operatorname{tame}}$. Thus
$\mathbf K \subseteq \mathbf K_1\subseteq\mathbf E$.
\end{proof}  

\begin{theorem}\label{thm:tamepairs} $\operatorname{RCF}_{\operatorname{tame}}$ is model complete.
\end{theorem}
\begin{proof}
Let $\mathbf K=(K,C)$, $\mathbf E=(E,C_E)$,  $\mathbf F=(F,C_F)$ be models
of~$\operatorname{RCF}_{\operatorname{tame}}$, where $\mathbf K\subseteq\mathbf E$ and $\mathbf K\subseteq\mathbf F$,
and $\mathbf F$ is $\abs{E}^+$-saturated. By Corollary~\ref{cor:modelcomplete test} and Zorn it is enough to show that
there is a substructure $\mathbf K_1$ of $\mathbf E$ that properly contains~$\mathbf K$, is a model of $\operatorname{RCF}_{\operatorname{tame}}$, 
and embeds over $K$ into $\mathbf F$. To simplify notation
we view subfields of $E$, respectively $F$, as valued ordered subfields of~$E$, respectively~$F$.
As before we distinguish three cases:

\case[1]{$C\neq C_E$.}
Take $e\in C_E\setminus C$; then $e\notin K$.  The saturation assumption on~$\mathbf F$ gives $f\in C_F$ realizing the same cut in $K$ as $e$. Arguing as in Case 1 in the proof of Theorem~\ref{thm:RCVF QE} this yields an ordered and valued field embedding $j\colon K(e)\to F$ over~$K$ sending
$e$ to $f$. With notations as in Lemma~\ref{lem:tamepairs} and its proof we can extend $j$ to an ordered field embedding $j_1\colon K_1\to F$, which
by Corollary~\ref{cvrc} is also a valued field embedding. 
On the $\mathbf F$-side, $j_1(K_1)$ is the real closure
of~$K(f)$ in~$F$ and $j_1(C_1)$ is the real closure of $C(f)$ in $C_F$, so 
$j_1$ embeds $\mathbf K_1$ into $\mathbf F$. 

\case[2]{$\Gamma\neq\Gamma_E$.} We first argue as in Case 2 in the proof of
Theorem~\ref{thm:RCVF QE}, then take real closures as in Case 1 above, and follow
the argument there, 
using instead of Lemma~\ref{lem:tamepairs} the statement about residue fields in Lemma~\ref{lem:lift value group ext}.

\case[3]{$C=C_E$ and $\Gamma=\Gamma_E$.}
Argue as in Case~3 in the proof of Theorem~\ref{thm:RCVF QE}, and then extend 
$j\colon K(x)\to F$ to the real closure of $K(x)$ inside $E$. 
\end{proof}

\begin{cor}\label{cor:RCFtame complete}
$\operatorname{RCF}_{\operatorname{tame}}$ is complete.
\end{cor}
\begin{proof} 
Take $x$ in an ordered field extension of $\Q$ with $x>\Q$, and let $\Q^{\operatorname{rc}}$ be the real closure of $\Q$ in $\Q(x)^{\operatorname{rc}}$. Consider $\Q(x)^{\operatorname{rc}}$
as a valued field whose valuation ring is the convex hull of
$\Q$ in $\Q(x)^{\operatorname{rc}}$.
Then $\Q(x)^{\operatorname{rc}}\models \operatorname{RCVF}$, and
$\Q^{\operatorname{rc}}$ is a maximal subfield of the 
valuation ring of $\Q(x)^{\operatorname{rc}}$, so
$\big(\Q(x)^{\operatorname{rc}},\Q^{\operatorname{rc}}\big)\models \operatorname{RCF}_{\operatorname{tame}}$. We claim
that $\big(\Q(x)^{\operatorname{rc}},\Q^{\operatorname{rc}}\big)$ embeds into every model of $\operatorname{RCF}_{\operatorname{tame}}$. To prove this claim, let $(K, C)\models\operatorname{RCF}_{\operatorname{tame}}$ and take any $y\in K^{>}$ with $y\succ 1$. Then
$y>\Q\subseteq K$, so we have an ordered field embedding $j\colon \Q(x)^{\operatorname{rc}}\to K$ with  $j(x)=y$. Consider any subfield of~$K$ as an ordered subfield of $K$, and let 
$\Q(y)^{\operatorname{rc}}$ be the real closure of $\Q(y)$ in~$K$. Then $j\big(\Q(x)^{\operatorname{rc}}\big)= \Q(y)^{\operatorname{rc}}$. The $j$-image of the valuation ring
of $\Q(x)^{\operatorname{rc}}$ is the convex hull of $\Q$ in
$\Q(y)^{\operatorname{rc}}$. This convex hull is a convex valuation ring of
$\Q(y)^{\operatorname{rc}}$; it is contained in $\mathcal{O}\cap \Q(y)^{\operatorname{rc}}$, but does not contain $y$, and thus it equals $\mathcal{O}\cap \Q(y)^{\operatorname{rc}}$ by Corollary~\ref{cor:ZA app}, so $j$ is a valued field embedding.

Finally, $j(\Q^{\operatorname{rc}})$ is
the real closure of $\Q$ in $K$, so $j(\Q^{\operatorname{rc}})\subseteq C\cap \Q(y)^{\operatorname{rc}}$. Now $C\cap \Q(y)^{\operatorname{rc}}$ is a real closed subfield of $\Q(y)^{\operatorname{rc}}$ by Corollary~\ref{cor:cor:AS}, and
$y\notin C\cap \Q(y)^{\operatorname{rc}}$, so $j(\Q^{\operatorname{rc}})= C\cap \Q(y)^{\operatorname{rc}}$. Thus $j$ embeds $\big(\Q(x)^{\operatorname{rc}},\Q^{\operatorname{rc}}\big)$ into $(K, C)$. 
\end{proof}

\begin{prop}\label{indsemi} 
Let $\mathbf K=(K,C)\models\operatorname{RCF}_{\operatorname{tame}}$.
If $X\subseteq K^n$ is definable \textup{(}with parameters\textup{)} in  $\mathbf K$, then $X\cap C^n$ is semialgebraic in the sense of~$C$. 
Thus if $X\subseteq K^n$ is semialgebraic in the sense of
$K$, then $X\cap C^n$ is semialgebraic in the sense of~$C$.
\end{prop}
\begin{proof} \marginpar{proof requires basics about RCF}
%Let $\mathcal L$ be the language of rings and $\mathcal L'=\mathcal %L_{\operatorname{tame}}$.
Let $x=(x_1,\dots,x_n)$ be a tuple of distinct variables.
Note that $C\preceq K$ as $\mathcal{L}_{\operatorname{OR}}$-structures, by \ref{thm:RCF QE}. 
Thus we have to show for every
$\mathcal L_{\operatorname{tame},K}$-formula $\varphi'(x)$ that $U(x)\wedge\varphi'(x)$ is equivalent in $\mathbf K$
to $U(x)\wedge\varphi(x)$ for some $\mathcal L_{\operatorname{OR},C}$-formula $\varphi(x)$, where
$U(x):=U(x_1)\wedge\cdots\wedge U(x_n)$. By Lemma~\ref{lem:separation}, this reduces to showing: \marginpar{\bf accept this reduction for the moment}
\begin{itemize}
\item[($\ast$)]
\textit{Let $\mathbf E=(E,C_E)$ and $\mathbf F=(F,C_F)$ be elementary
extensions of~$\mathbf K$ and suppose $e\in C_E^n$ and $f\in C_F^n$
realize the same type over $C$ in the real closed fields
$C_E$ and~$C_F$, respectively. Then $e$ and $f$ realize the same type over~$K$ in~$\mathbf E$ and~$\mathbf F$, respectively.}\/
\end{itemize}
Taking real closures in $C_E$ and $C_F$,  
the hypothesis of ($\ast$) yields a unique ordered field isomorphism $C(e)^{\operatorname{rc}}\to C(f)^{\operatorname{rc}}$ over~$C$ sending $e$ to $f$. We claim that this isomorphism extends to an ordered field isomorphism $K(e)^{\operatorname{rc}}\to K(f)^{\operatorname{rc}}$ over~$K$ (taking real closures in $E$ and~$F$),
and that \begin{align*}
C_E\cap K(e)^{\operatorname{rc}}\ &=\ C(e)^{\operatorname{rc}}, & \big(K(e)^{\operatorname{rc}},C(e)^{\operatorname{rc}}\big)&\ \models\ \operatorname{RCF}_{\operatorname{tame}},\\
 C_F\cap K(f)^{\operatorname{rc}}\ &=\ C(f)^{\operatorname{rc}}, & \big(K(f)^{\operatorname{rc}},C(f)^{\operatorname{rc}}\big)&\ \models\ \operatorname{RCF}_{\operatorname{tame}}.
 \end{align*}
In view of Theorem~\ref{thm:tamepairs}, the
conclusion of ($\ast$) follows from this claim.  

By induction on $n$ we reduce the claim to the case $n=1$. If 
$e\in C$, the claim holds trivially, so assume $e\notin C$. Then $f\notin C$, and  $e$ and $f$ realize the same cut in~$C$,
and therefore the same cut in $K$. Thus we have an ordered field isomorphism 
$K(e)^{\operatorname{rc}}\to K(f)^{\operatorname{rc}}$ over~$K$ sending $e$ to $f$, and the rest holds by Lemma~\ref{lem:tamepairs}. 
%this is also an isomorphism $\big(K(e)^{\operatorname{rc}}, %C(e)^{\operatorname{rc}}\big)\to \big(K(f)^{\operatorname{rc}}, %C(f)^{\operatorname{rc}}\big)$. (Also need ~\ref{thm:tamepairs} %again?)
 %Case~1 of the proof of Theorem~\ref{thm:tamepairs} 
 %we already showed
%that $\operatorname{st}_i K(c_i)^{\operatorname{rc}} = %C(c_i)^{\operatorname{rc}}$ for $i=1,2$.
\end{proof}
\marginpar{section checked except for and
boldfaced comments}

\subsection*{Notes and comments} The first explicitly model-theoretic result
on valued fields is due to A.~Robinson~\cite{Robinson56}: model completeness of ACVF (close to Theorem~\ref{thm:ACVF QE}).
Theorems~\ref{thm:RCVF QE} and~\ref{thm:tamepairs} are from Cherlin-Dickmann~\cite{CherlinDickmann} and Macintyre~\cite{Macintyre68}. 
Model completeness of $\operatorname{ACVF}$ and $\operatorname{RCVF}$ also follow
from the more general results of Ax~\&~Kochen~\cite{AxKochen,AxKochen2} and Er{\accentv s}ov~\cite{Ersov},
in view of \ref{ex:DOAb},~\ref{thm:ACF QE},~\ref{thm:RCF QE}. (See \cite[Section~1.3]{CherlinDickmann}.) 
Corollary~\ref{cor:RCVF QE} was noticed by Becker~\cite{Becker83}.
Proposition~\ref{indsemi} is a special case of  \cite[Proposition~8.1]{vdD}; see also~\cite{Marker-Steinhorn}.
The proof that $\mathbb T$ has NIP in Section~\ref{sec:embth} is modeled
on that of Corollary~\ref{cor:RCVF NIP} above.
%By Delon~\cite{Delon}, a henselian valued field of equicharacteristic zero has NIP iff its residue field has NIP.
For proofs that $\operatorname{ACVF}$ and $\operatorname{RCF}_{\operatorname{tame}}$ have NIP, see~\cite{Delon} and \cite{GH}, respectively. 
%(For $\operatorname{RCF}_{\operatorname{tame}}$ this also follows
%from Corollary~\ref{cor:RCFtame complete} and Proposition~\ref{}.) 
%\marginpar{Reference is to 15.6.4.}

In this section we construed valued fields as one-sorted structures by
encoding the valuation as a dominance relation.
It is often more informative to view valued fields as three-sorted structures, with sorts for the underlying field, for the value group, and for the residue field; for example, see~\cite{Pas89}. In
Chapter~\ref{ch:monotonedifferential} we use this setting in dealing with certain valued \textit{differential}\/ fields.

%% file: mt-3-7.tex
\section{The Newton Tree of a Polynomial over a Valued Field}
\label{Newton Diagrams}

\noindent
In this section we use Newton diagrams to construct a
{\em Newton tree}\/ for any given nonconstant polynomial over a
henselian valued field of equicharacteristic zero. Such a Newton
tree is a finite tree of {\em approximate zeros}\/ of the polynomial, and
induces a partition of the field into finitely many simple pieces
on each of which the polynomial behaves in a simple way. 
In Chapters~\ref{ch:The Dominant Part and the Newton Polynomial} and~\ref{ch:newtonian fields} we
use more delicate Newton diagrams for differential polynomials
over suitable valued differential fields, and in preparation for this, the reader may find the exposition of the Newton diagram method
for ordinary polynomials below helpful. 
The present section also has an application in the proof of Proposition~\ref{prop:vP(y), order 1} below, and will be useful in the next volume. This section does not depend on Sections~\ref{sec:decomposition}--\ref{sec:modth val fields}.  

\begin{notation}  Let $F$ be a field and $P\in F[Y]^{\ne}$. Then $\deg P\in \N$ denotes the {\bf degree}
of $P$, and $\val P\in \N$ denotes the {\bf multiplicity} of
$P$ at $0$: the largest $m$ such that $P\in Y^mF[Y]$.
For any $f\in F$ we let $P_{+ f}:= P(f+Y)$ be the {\bf additive conjugate\/} of~$P$ by~$f$, and 
$P_{\times f}:= P(fY)$ the 
{\bf multiplicative conjugate\/} of $P$ by~$f$.   
\end{notation}

\index{degree!polynomial}
\index{polynomial!degree}
\index{multiplicity!polynomial}
\index{conjugation!additive}
\index{conjugation!multiplicative}
\nomenclature[A]{$\deg P$}{degree of $P$}
\nomenclature[A]{$\val P$}{multiplicity of $P$ at $0$}
\nomenclature[A]{$P_{+h}$}{additive conjugate $P(Y+h)$ of $P$ by $h$}
\nomenclature[A]{$P_{\times h}$}{multiplicative conjugate $P(hY)$ of $P$ by $h$}

\noindent
{\em Throughout this section we fix a valued field $K$
with valuation ring $\mathcal O\ne K$, 
residue field $\k$ of characteristic $0$,
residue map $a\mapsto\bar{a}\colon \mathcal O=K^{\preceq 1}\to \k$, 
value group~$\Gamma$ and valuation 
$v\colon 
K^{\times} \to \Gamma$. We also choose a subset $\fM$ of $K^\times$
which is mapped bijectively onto $\Gamma$ by $v$. \textup{(}We do not assume that $\fM$ is a subgroup of $K^\times$.\textup{)} 
For $\gamma\in\Gamma$ we let $s\gamma$ be the unique element of $\fM$ with $v(s\gamma)=\gamma$.
``Equivalence'' in this section refers to the equivalence relation $\sim$ on
$K^{\times}$ induced by the valuation. We let
$f$,~$g$,~$y$,~$z$ range over $K$, and $\fm$, $\fn$ over $\fM$. Finally, we fix throughout a polynomial $P(Y)\in K [Y]^{\neq}$: $\ P\ =\ a_0 + a_1Y + \dots + a_nY^n,\ a_0,a_1,\dots,a_n\in K,\ a_n\ne 0.$}

\index{dominant!monomial}
\index{monomial!dominant}
\index{dominant!part}
\index{dominant!multiplicity}
\index{dominant!degree}
\index{multiplicity!dominant}
\index{degree!dominant}
\nomenclature[K]{$\mathfrak d_P$}{dominant monomial of $P$}
\nomenclature[K]{$D_P$}{dominant part of $P$}
\nomenclature[K]{$\dval P$}{dominant multiplicity of $P$}
\nomenclature[K]{$\ddeg P$}{dominant degree of $P$}

\subsection*{Dominant part} The {\bf dominant monomial} of $P$
is the unique element $\mathfrak d_P$ of $\fM$ with $\mathfrak d_P\asymp P$. Then 
$\frak d_P^{-1}P\in \mathcal O[Y]$, and we call
the  polynomial
$$D_P\ :=\ \sum_{i} \bar{(a_i/\frak d_P)}\, Y^i \in \k[Y]$$
the {\bf dominant part} of $P$.   
Clearly $D_P$ is nonzero with 
$$\val P\ \leq\ \val D_P\ \leq\ \deg D_P\ \le\ \deg P.$$ 
We call $\dval P:=\val D_P$ the {\bf dominant multiplicity} of $P$ at $0$ and
$\ddeg P:=\deg D_P$ the {\bf dominant degree} of $P$. 
Note that
\begin{align*}
P_{\times\fm}\	&=\ \sum_i a_i \fm^i\,Y^i, \quad\text{so} \\
D_{P_{\times\fm}}\	&=\ \sum_i \bar{(a_i \fm^i/\mathfrak d)}\, Y^i\qquad\text{where
$\mathfrak d=\mathfrak d_{P_{\times \fm}}$.}
\end{align*}

\begin{lemma}\label{lem:ddeg monotone} 
$\fm\prec \fn\ \Rightarrow\  \dval P_{\times \fm}\leq \ddeg P_{\times \fm} \leq \dval P_{\times \fn}\leq \ddeg P_{\times \fn}.$
\end{lemma}
\begin{proof}
Suppose $\fm\prec \fn$; it suffices to show that then $\ddeg P_{\times \fm}\leq \dval P_{\times \fn}$.
Let $d=\dval P_{\times \fn}$. Then for $i>d$ we have $a_i\fn^i \preceq a_d\fn^d$ and so
$$a_i\fm^i\ =\ (a_i\fn^i) \cdot (\fm/\fn)^i\ \prec\ (a_d\fn^d) \cdot  (\fm/\fn)^d\ =\ a_d\fm^d,$$
and thus $\ddeg P_{\times \fm} \leq d$ as required.
\end{proof}

\subsection*{Approximate zeros}
An {\bf approximate zero} of $P$ is an element
$y$ such that
$$v\bigl(P(y)\bigr)\ >\ \min_i v(a_iy^i)\quad 
\text{ (in particular $n\ge 1$ and $y\ne 0$)};$$ 
equivalently, $y\neq 0$ and $D_{P_{\times\fm}}(c)=0$, where $\fm\asymp y$ and $c=\bar{(y/\fm)}$. In this case the polynomial $D_{P_{\times\fm}}\in\k[Y]$ is not homogeneous, since $c\ne 0$, and so
with $\delta := \min_i v(a_iy^i)$, 
there are at least two elements $i\in \{0,\dots,n\}$ 
such that $v(a_iy^i)=\delta$. We say that $\fm$ is a {\bf starting monomial} for $P$ if $D_{P_{\times\fm}}$ is not
homogeneous.
Note that if $y$ is an approximate zero of $P$ and $y \sim z$, then $z$ is an
approximate zero of~$P$.
If $P(y)=0$ and $y\ne 0$, then $y$ is an approximate zero of $P$. 

\index{approximate!zero}
\index{monomial!starting}
\index{starting monomial}

\begin{remark} 
We have $P_{+f}(Y)= P(f+Y)=b_0 + b_1Y + \dots + b_nY^n$ with
$$b_i\ =\ \sum_{j=0}^{n-i}\binom{i+j}{i}a_{i+j}f^j,\ \qquad
\text{in particular, $b_0=P(f)$.}$$
Suppose that $f\preceq 1$; then $P_{+f}\asymp P$: the identities above give $P_{+f}\preceq P$, and hence
$P\preceq P_{+f}$ by
$P=(P_{+f})_{+g}$ for $g=-f$. 
Next, suppose that $vf=\beta \le 0$ and~$f$ is \textit{not}\/ 
an approximate zero of $P$. Then 
$v(P_{+f})$ depends only on $\beta$, not on $f$:
 $$v(P_{+f})\ =\ v\bigl(P(f)\bigr)\ =\ \min_i v(a_i) +i\beta.$$
To see this, put $\gamma:= v(b_0)= \min_i v(a_i) +i\beta$, so for $0\le i \le n$ 
and $0\le j \le n-i$,
 $$v(a_{i+j}f^j)\ =\ v(a_{i+j})+j\beta\ 
=\  v(a_{i+j}) + (i+j)\beta - i\beta\ \ge\ \gamma - i\beta.$$ 
Hence $v(b_i)\ge \gamma - i\beta \ge \gamma$ for $i>0$. The assertion follows.
\end{remark}

\subsection*{Geometric interpretation} Recall that 
$\Gd$ is the divisible
hull of $\Gamma$. We refer to $\Z \times \Gd$ as
the \textit{plane\/}, and to its elements as \textit{points\/}. The {\bf abscissa\/} \index{abscissa}
of a point~$(i,\alpha)$ is the integer $\operatorname{abscis}(i,\alpha):=i$. We view
$\Z\times \{0\}$ as the horizontal axis of the plane, and 
$\{0\}\times \Gd$ as its vertical axis. For $\beta\in \Gd$, 
define the additive function
$$L_\beta\colon \Z \times \Gd \to \Gd, \qquad 
L_\beta(i,\alpha):=\alpha +i\beta.$$
Given any $\beta,\delta\in \Gd$ we refer to the set
$$\bigl\{(i,\alpha)\in \Z\times \Gd:\ L_\beta(i,\alpha)=\delta\bigr\}$$  
as \textit{the line $L_\beta=\delta$.}\/ A point 
$(i,\alpha)$ is said to lie \textit{above\/} (respectively
\textit{on\/}, respectively \textit{below\/}) this line if
$L_\beta(i,\alpha) > \delta$ (respectively $L_\beta(i,\alpha) =\delta$, 
respectively
$L_\beta(i,\alpha) < \delta$). Since these are the only kind of lines we need
to consider, by ``line'' we shall always mean a line of the form 
$L_\beta=\delta$ as above. Note that if $p$,~$q$ are points on the line~$\ell$ and $p \ne q$, then 
$\operatorname{abscis}(p)\ne \operatorname{abscis}(q)$.
Each line contains infinitely many 
points, and is of the form $L_\beta=\delta$ for exactly one pair 
$(\beta,\delta)\in \Gd \times \Gd$.
For any two points $(i_1,\alpha_1)$ and $(i_2,\alpha_2)$ with $i_1 \ne i_2$
there is exactly one line 
containing both points, namely $L_\beta=\delta$ with 
$\beta=-\bigl( \frac{\alpha_2 - \alpha_1}{i_2-i_1}\bigr)$, 
$\delta=\alpha_1 + i_1\beta$.
We call $\beta$ the {\bf antislope} of the line 
$L_\beta=\delta$. (Its slope is $-\beta$.) 

\index{antislope}
\index{line!antislope}
\index{line!slope}
\index{slope!line}
\index{Newton!diagram}
\index{Newton!diagram!edges}
\index{edge}
\index{vertex}
\nomenclature[K]{$\mathcal N(P)$}{Newton diagram of $P$}

%\begin{definition}
\medskip\noindent
We define the {\bf Newton diagram} of $P$ to be the finite 
nonempty set of points
$$\mathcal N (P)\ :=\ 
\bigl\{\bigl(i,v(a_i)\bigr):\ i=0,\dots,n,\ a_i\ne 0\bigr\}\
\subseteq\ \Z \times \Gamma\ \subseteq\ \Z \times \Gd.$$
 An {\bf edge} of $\mathcal N(P)$ is a line $\ell$ 
that contains at least two points of $\mathcal N(P)$ and such that all points of
$\mathcal N(P)$ lie on or above $\ell$. In this case, the point 
of $\ell \cap \mathcal N(P)$ with least abscissa is called the 
{\bf left vertex} \index{line!left vertex} of $\ell$ in $\mathcal N(P)$, and the point of
$\ell \cap \mathcal N(P)$ with largest abscissa is called the 
{\bf right vertex} \index{line!right vertex} of $\ell$ in $\mathcal N(P)$. 
Figure~\ref{fig:newtondiag1} shows a Newton diagram with $6$ points, $4$ edges, and $n=6$; this diagram is missing a point with abscissa $4$, since~$a_4=0$. 
%\end{definition}

\begin{figure}[h]
\input{mt-newtondiag1.tex}
\caption{Picture of a Newton diagram.}\label{fig:newtondiag1}
\end{figure}

\medskip\noindent
Note that if~$\ell$ and~$\ell'$ are edges of $\mathcal N(P)$ and 
$\ell \ne \ell'$, then $\operatorname{antislope}(\ell) \ne
\operatorname{antislope}(\ell')$.  The {\bf antislopes} \index{Newton!diagram!antislopes} of  $\mathcal N(P)$ are by definition the antislopes of the
edges of $\mathcal N(P)$. If~$\ell$ and~$\ell'$ are edges 
of $\mathcal N(P)$ with antislopes $\beta>\beta'$, then
$$\text{abscis(right vertex of }\ell)\ \le\ 
\text{abscis(left vertex of }\ell').$$
In more graphic terms: the antislope increases when
moving from right to left in the diagram along the edges. 
Now let $\ell$ be a line given by $L_\beta=\delta$, where $\beta,\delta\in\Q\Gamma$. Then for $i=0,\dots,n$,
\begin{align*}
\text{$\big(i,v(a_i)\big)$ lies above $\ell$}&\quad\Longleftrightarrow\quad v(a_i)+i\beta>\delta, \\
\text{$\big(i,v(a_i)\big)$ lies on $\ell$}&\quad\Longleftrightarrow\quad v(a_i)+i\beta=\delta, \\
\text{$\big(i,v(a_i)\big)$ lies below $\ell$}&\quad\Longleftrightarrow\quad v(a_i)+i\beta<\delta.
\end{align*}
Hence 
$$\left.\text{\parbox{20em}{the line
$\ell$ contains at least one point of $\mathcal N(P)$ and all points of $\mathcal N(P)$ lie on or above $\ell$}}\ \right\}
\quad\Longleftrightarrow\quad \delta=\min_i v(a_i)+i\beta.$$
The antislope $\beta$ of $\ell$ does not necessarily lie in $\Gamma$, but suppose it does, and suppose also that
$\delta=\min_i v(a_i)+i\beta$. Then, with $\fm=s\beta$, 
$$D_{P_{\times\fm}}\ =\ \sum_i \bar{(a_i \fm^i/\mathfrak d)}\, Y^i\qquad\text{where
$\mathfrak d=\mathfrak d_{P_{\times\fm}}=s\delta$,}$$
and so $\ell$ is an edge of $\mathcal N(P)$ if and only if 
$\fm$ is a starting monomial for $P$. 
%In fact, map $\fm\mapsto v\fm$ is a bijection between the set of starting monomials for $P$ and the set of antislopes of $\mathcal N(P)$ which lie in~$\Gamma$.
If $\ell$ is an edge of $\mathcal N(P)$, then we set
$P_\beta:=D_{P_{\times\fm}}$, and we have 
$$\text{abscis(left vertex of $\ell$)}\ =\ \val P_\beta, \qquad 
\text{abscis(right vertex of $\ell$)}\ =\ \deg P_\beta.$$
The relation to the notion of approximate zero is as follows:
\begin{enumerate}
\item If $y\asymp\fm$ is an approximate zero of $P$, then $\mathcal N(P)$ has an antislope $\beta:=v\fm$, and $\bar{(y/\fm)}$ 
is a zero of $P_\beta$. 
\item If $\beta=v\fm\in\Gamma$ is an antislope of $\mathcal N(P)$ and $c\in \k^{\times}$ is a zero of $P_\beta$, then $P$ has an
approximate zero $y\asymp\fm$ with $\bar{(y/\fm)}=c$.
\end{enumerate}
Note that $y$ in (2) is determined up to equivalence % in the sense of $\sim$ 
by $\beta$ and~$c$. Clearly there are at most~$n$ antislopes of $\mathcal N(P)$.
In this way the geometric interpretation of approximate zeros shows that
$P$ has, up to equivalence, at most $n$ approximate zeros. 

\index{closed!$\preceq$-closed}
\index{subset!$\preceq$-closed}
\index{asymptotic equation}
\index{equation!asymptotic}
\index{asymptotic equation!solution}
\index{asymptotic equation!approximate solution}
\index{approximate!solution}
\index{asymptotic equation!dominant degree}
\index{asymptotic equation!primary dominant part}
\index{dominant!part!primary}
\index{primary dominant part}
\index{dominant!degree!on $\mathcal E$}
\index{dominant!degree!asymptotic equation}
\nomenclature[K]{$\ddeg\eqref{ast}$}{dominant degree of \eqref{ast}}
\nomenclature[K]{$P_{\beta\eqref{ast}}$}{primary dominant part of \eqref{ast}}
\nomenclature[K]{$\ddeg_{\E}P$}{dominant degree of $P$ on $\E$}
\nomenclature[K]{$\beta\eqref{ast}$}{smallest antislope of $\mathcal N(P)$ which is
$\geq v\fm$ for some $\fm\in\E$, if there is one}

\subsection*{Asymptotic equations} 
A set $\E\subseteq K^\times$ is called {\bf $\preceq$-closed} if $\E\ne \emptyset$, and $f\in\E$ whenever $0\neq f\preceq g\in\E$. 
Let us consider an 
{\bf asymptotic equation}
\begin{equation}\label{ast}\tag{E}
P(Y)=0, \quad  \quad Y\in\E
\end{equation}
where $\E\subseteq K^\times$ is $\preceq$-closed. 
A {\bf solution} of \eqref{ast} is an element $y$ of $\E$
such that $P(y)=0$. 
An {\bf approximate solution} of \eqref{ast} is an approximate
zero~$y\in\E$ of~$P$. 
If $\mathcal N(P)$ has an antislope $\geq v\fm$ for some $\fm\in\E$,
let $\beta\eqref{ast}$ be the least among these antislopes, and define
the {\bf dominant degree}
of \eqref{ast} to be the abscissa of the right vertex of the edge of $\mathcal N(P)$ with antislope $\beta\eqref{ast}$, denoted by $\ddeg\eqref{ast}$; thus $\val P<\ddeg\eqref{ast}\leq\deg P$.
If in addition $\beta\eqref{ast}\in\Gamma$, then $\ddeg\eqref{ast}=\deg P_{\beta\eqref{ast}}$, 
and we call $P_{\beta\eqref{ast}}$ the {\bf primary dominant part} 
of \eqref{ast}. If all antislopes of $\mathcal N(P)$ are $<v(\E)$, then 
we define the dominant degree~$\ddeg\eqref{ast}$ of~\eqref{ast} to be $0$. (In that case~\eqref{ast} has no approximate solution.)
Related to $\ddeg \eqref{ast}$ is the {\bf dominant degree} of $P$ on $\E$, defined to be the natural number
$$\ddeg_{\E} P := \max \big\{\ddeg P_{\times\fm}:\fm\in\E\big\}.$$
Note that for all $\fm$,
$$\val P\ =\ \val P_{\times\fm}\ \leq\ 
\val D_{P_{\times\fm}}\ \leq\ \deg D_{P_{\times\fm}}\ \leq\
\deg P_{\times\fm}\ =\ \deg P,$$
and so $\val P \leq \ddeg_{\E} P\leq \deg P$. Moreover:

\begin{lemma} If $\ddeg \eqref{ast}=0$, then $\val P = \ddeg_{\E} P$. If on the other hand $\ddeg \eqref{ast}>0$, then $\val P < \ddeg\eqref{ast} = \ddeg_{\E} P$.
%$\ddeg\eqref{ast} \leq \ddeg_{\E} P$, and
%$$ \ddeg\eqref{ast} = \ddeg_{\E} P
%\ \Longleftrightarrow\ \ddeg_\E P>\val P\ 
%\Longleftrightarrow\ \ddeg\eqref{ast}>0.$$ 
\end{lemma}

\begin{proof}
Let $K^\alg$ be an algebraic closure of $K$, equipped with a valuation extending that of $K$,
let $\fM^\alg\supseteq\fM$ be a subset of $(K^\alg)^\times$ which maps bijectively onto the value group $\Q\Gamma$ of $K^\alg$ under
this valuation, and let 
$$\E^\alg\ :=\ \{a\in (K^\alg)^\times:\ \text{$a\preceq\fm$ for some $\fm\in\E$}\}$$
be the smallest $\preceq$-closed subset of $(K^\alg)^\times$ containing $\E$.
Then $\ddeg\eqref{ast}$ does not change if  we replace $K$, $\fM$, $\E$ by
$K^\alg$, $\fM^\alg$, $\E^\alg$, and $\ddeg_{\E} P=\ddeg_{\E^\alg} P$ by Lemma~\ref{lem:ddeg monotone}.
Thus we may assume that $\Gamma$ is divisible, 
and do so below. 

We establish the first implication by proving the contrapositive. So assume that $m:=\val P< d:= \ddeg_\E P$, and take 
$\fm\in\E$ with $\ddeg P_{\times\fm}=d$.
Then $\big(m,v(a_m)\big)$ lies on or above the line 
through $\big(d,v(a_d)\big)$ with antislope $v\fm$.
Hence $\mathcal N(P)$ has an edge with left vertex 
$\big(m,v(a_m)\big)$ and antislope $\geq v\fm$, and thus
$\ddeg\eqref{ast} > 0$.

For the second implication, assume $\ddeg\eqref{ast} > 0$.
It is clear that then
$\val P < \ddeg\eqref{ast} \leq \ddeg_{\E} P$. 
Assume towards a contradiction that $\fm\in\E$ is such that 
$d:=\ddeg\eqref{ast}<i:=\ddeg P_{\times\fm}$.
Then   $v(a_d)+dv\fm \geq v(a_i)+iv\fm$.
Also, the edge of~$\mathcal N(P)$ with antislope $\beta=\beta\eqref{ast}$ 
passes through the point $\big(d,v(a_d)\big)$ but
not through $\big(i,v(a_i)\big)$,
so $v(a_i)+i\beta>v(a_d)+d\beta$. 
Let $\ell_1$ be the line passing through $\big(d,v(a_d)\big)$ and $\big(i,v(a_i)\big)$; then the antislope $\beta_1=-\left(\frac{v(a_i)-v(a_d)}{i-d}\right)$
of $\ell_1$ satisfies $v\fm\leq\beta_1<\beta$.
Thus $\mathcal N(P)$ has an edge with left vertex $\big(d,v(a_d)\big)$ and antislope $\beta^*$ such that
$v\fm\le \beta_1 \le \beta^* < \beta$,
 contradicting the minimality of~$\beta$.
%This proves ``$\ddeg\eqref{ast}>0 \Rightarrow \ddeg\eqref{ast} = \ddeg_{\E} P$.''
\end{proof}

%\noindent
%We say that \eqref{ast} is {\bf quasilinear} if 
%$\ddeg\eqref{ast} = 1$.

\begin{remarkNumbered}\label{rem:henselian, unique sol}
Suppose $K$ is henselian and $\ddeg\eqref{ast} = 1$. 
Then~\eqref{ast} has a unique solution. To see this, 
let $\beta=\beta\eqref{ast}$, and note that $\beta\in\Gamma$. Let $\fm=s\beta$, $\mathfrak d=\mathfrak d_{P_{\times\fm}}$.
The image of $\mathfrak d^{-1} P_{\times\fm}\in \mathcal{O}[Y]$ in $\k[Y]$ is 
$P_\beta=D_{P_{\times\fm}}$, which has degree~$1$ and multiplicity~$0$, and thus has a (unique) zero in~$\k^{\times}$. 
Now $K$ being henselian,
it follows that $\mathfrak d^{-1} P_{\times\fm}$ has a unique zero $y$ in~$\mathcal{O}$. For this $y$ we have $y\asymp 1$, and
so~$\fm y$ is a solution of \eqref{ast}. 
It is the only solution of  \eqref{ast} because $\mathcal N(P)$
has only one antislope that is $\geq v\fn$ for some 
$\fn\in\E$.
\end{remarkNumbered}

\subsection*{Refinements}
A {\bf refinement} of \eqref{ast} is an asymptotic equation of the form
\begin{equation}\label{astast}\tag{E$'$}
P_{+f}(Y)=0, \quad \quad Y\in\E'
\end{equation}
where $f\in \E\cup\{0\}$ and $\E'\subseteq\E$ is $\preceq$-closed. 
If $y$ is a solution of~\eqref{astast} and $f+y\neq 0$, then $f+y$ is a solution of \eqref{ast}.\index{asymptotic equation!refinement}\index{refinement} 
Moreover, if 
$y\not\sim-f$ is an approximate solution
of~\eqref{astast}, then $f + y$ is an approximate solution 
of \eqref{ast}. To see this, 
let $y\not\sim-f$ be an approximate solution
of \eqref{astast}. We have
$$ P_{+f}(Y)\ =\ b_0 + b_1Y + \dots + b_nY^n, \qquad
b_i\ =\ \frac{P^{(i)}(f)}{i!}\ =\ \sum_{j=i}^n \binom{j}{i}a_j f^{j-i}.$$
From $y\not\sim-f$ we get $vf,vy\geq v(f+y)$, and so
\begin{align*}
v\bigl(P(f + y)\bigr)\ &=\ v\bigl( P_{+f}(y)\bigr)\ >\ \min_i v(b_iy^i)\ 
      =\ \min_i v\left( \sum_{j=i}^n \binom{j}{i}a_j f^{j-i}y^i\right)\\
    &\ge\ \min_i\min_{j\ge i} v\left(a_j f^{j-i}y^i\right)\ \ge\ \min _j v\bigl(a_j(f + y)^j\bigr).  
%      &\ge\ \min_i v\left( \sum_{j=i}^n \binom{j}{i}a_j 
% (f+y)^{j}\right)\ \ge\ \min _j v\bigl(a_j(f + y)^j\bigr).
\end{align*}
For $g\neq 0$ and $\E=\{y: 0\neq y\prec g\}$ we also indicate \eqref{ast} by
$$P(Y)=0,\quad \quad Y\prec g,$$
and we set $\ddeg_{\prec g} P:=\ddeg_{\E} P$. This notation is used in the next lemmas.

\begin{lemma} \label{Refinement-Lemma}
Let $f$ be an approximate solution of \eqref{ast}. Put
$\beta:=vf$, $\fm:=s\beta$, $\mu:=\text{multiplicity of the zero 
$\bar{f/\fm}$ of $P_\beta$}$, and
$P_{+f}(Y)=b_0 + b_1Y + \dots + b_nY^n$ with $b_0,\dots,b_n \in K$. 
Consider the refinement
\begin{equation}\label{astphi}\tag{E$_{+f}$}
P_{+f}(Y)=0, \quad \quad Y\prec f
\end{equation}
of \eqref{ast}. Then
\begin{enumerate}
\item[\textup{(i)}] $b_{\mu}\ne 0$, and all points of $\mathcal N(P_{+f})$ with abscissa 
$<\mu$ lie above the line
through $\bigl(\mu,v(b_\mu)\bigr)$ with antislope $\beta$, and all points of 
$\mathcal N(P_{+f})$ with abscissa $>\mu$ lie on or above that line;
\item[\textup{(ii)}] if $b_i= 0$ for all $i<\mu$, then \eqref{astphi} has 
dominant degree $0$;
\item[\textup{(iii)}] if $b_i\ne 0$ for some $i<\mu$, then 
\eqref{astphi} has dominant degree $\mu$.
\end{enumerate}
\end{lemma}
\begin{proof} Note that
$b_i=\frac{P^{(i)}(f)}{i!}$. Set
$\delta:= \min_i v(a_i f^i)$, $\mathfrak d:= s\delta$, and
$$Q(Y):= \mathfrak d^{-1} P_{\times\fm}(Y) = \sum_i (a_i \fm^i/\mathfrak d)\, Y^i\in \mathcal{O}[Y].$$ 
Hence $Q^{(i)}(f/\fm)=(\fm^i/\mathfrak d)P^{(i)}(f)$, so 
$b_i=\frac{\mathfrak d}{i!\fm^i}Q^{(i)}(f/\fm)$, $i=0,\dots,n$. Note that
\begin{align*}
v\bigl(Q^{(i)}(f/\fm)\bigr)\ &>\ 0\quad \text{for $0\le i<\mu$,}\\ 
v\bigl(Q^{(\mu)}(f/\fm)\bigr)\ &=\ 0, \\ 
v\bigl(Q^{(j)}(f/\fm)\bigr)\ &\ge\ 0\quad \text{for $\mu < j \le n$.}
\end{align*} 
Thus 
\begin{align*}
 v(b_i)\ &>\  \delta - i\beta\quad \text{for $0 \le i < \mu$,}\\
 v(b_\mu)\ &=\  \delta - \mu\beta, \\
 v(b_j)\ &\ge\  \delta - j\beta\quad \text{for $\mu < j \le n.$}
\end{align*}
These inequalities give (i). Also, if 
$0\ne y \prec f$ and $\mu < j \le n$, then
$$v(b_\mu y^\mu)\ =\ (\delta - \mu\beta)+\mu vy\ =\ 
\delta + \mu\bigl(vy - \beta\bigr)\ <\  
\delta + j\bigl(vy -\beta\bigr)\ \le\ v(b_j y^j).$$
Hence the dominant degree of \eqref{astphi} is at most
$\mu$. If $b_i=0$ for all
$i<\mu$ (equivalently, $f$ is a zero of multiplicity $\mu$ of $P$), 
then \eqref{astphi} has dominant degree $0$.
Suppose $b_i\ne 0$ for some $i<\mu$. The inequalities above show that 
each line through the point $\bigl(\mu,v(b_\mu)\bigr)$ and some point 
$\bigl(i,v(b_i)\bigr)$ 
with $i<\mu$ and $b_i\ne 0$ has anti\-slope~$> \beta$. Among these lines, let
$\ell'$ be the one with minimal antislope $\beta' > \beta$. 
One checks easily that
all points $\bigl(j,v(b_j)\bigr)$ with $\mu < j\le n$ and $b_j\ne 0$ 
lie above
$\ell'$, so $\ell'$ is an edge of $\mathcal N(P_{+f})$. 
Thus \eqref{astphi} has dominant degree $\mu$ and 
$\beta\eqref{astphi} = \beta'$.
\end{proof} 

\begin{figure}[h]
\input{mt-newtondiag2.tex}
\caption{Behavior of Newton diagrams under refinement.}\label{fig:newtondiag2}
\end{figure} 

\noindent
Figure~\ref{fig:newtondiag2} illustrates Lemma~\ref{Refinement-Lemma} above in the case where our value group $\Gamma$ is the ordered abelian group~$\Q$, $\E=\{y\in K^\times: vy>1/2\}$, and
$$P\ =\ t^2 + 2t^{-1}Y + t^{-2}Y^2 + t^{-3/2}Y^3 + t^{-2} Y^5 + t^{3/2}Y^6
\quad\text{where $t\in\fM$, $vt=1$.}$$
In the Newton diagram on the left, the dotted line has antislope $1/2$, and the edges of antislope $>1/2$ have been highlighted. On the right is the Newton diagram   of
the additive conjugate $P_{+f}$ of $P$ by the approximate zero $f=-2t$ of $P$. 

\medskip
\noindent
In the previous lemma the dominant degree 
of \eqref{astphi} is less than the dominant degree of \eqref{ast}, 
except when $\beta=\beta\eqref{ast}$ and $\mu=\deg (P_\beta)$, in which case
\eqref{ast} and~\eqref{astphi} have the same dominant degree. The next lemma
handles this exceptional case. 
In the proof we tacitly use that if $n\geq 1$
(so $P'\neq 0$), then the Newton diagram~$\mathcal N(P')$ is obtained by shifting each point of $\mathcal N(P)$ 
not on the vertical axis
 $\{0\}\times \Gd$ by one unit to the left:
$$\mathcal N(P')\ =\ \bigl\{\bigl(i-1,v(a_i)\bigr):\ i\ge 1,\ a_i\ne 0\bigr\}.$$ 
%We say that \eqref{ast} is {\bf unraveled} 
%if \eqref{ast} has dominant degree $d\geq 1$ and 
%either $\beta\eqref{ast}\notin\Gamma$ or 
%$\beta\eqref{ast}\in\Gamma$ and $P_{\beta\eqref{ast}}$
%has no zero in $\k$ of multiplicity $d$.

\begin{lemma}\label{Special-Lemma} 
Let $K$ be henselian. Suppose that $\eqref{ast}$ has dominant degree $d\geq 1$, $\beta\eqref{ast}\in\Gamma$, and $P_{\beta\eqref{ast}}$
has a zero in $\k$ of multiplicity $d$. Then:
\begin{enumerate}
\item[\textup{(i)}] $P^{(d-1)}$ has a unique zero $f\in\E$;
\item[\textup{(ii)}] $vf = \beta\eqref{ast}$ and $f$ is an approximate solution 
of \eqref{ast}; moreover, if $y$ is any approximate solution of \eqref{ast}, then $y\sim f$;
\item[\textup{(iii)}] if the refinement \eqref{astphi} 
of \eqref{ast} still has dominant degree $d$ with 
$\beta\eqref{astphi}\in \Gamma$, 
then $P_{\beta\eqref{astphi}}$ has no zero  in $\k$ of multiplicity $d$. 
\end{enumerate} 
\end{lemma}
\begin{proof} To keep notations simple, put $\beta=\beta\eqref{ast}$. Then 
$P_\beta =c_1(Y-c_2)^d$
with $c_1,c_2\in \k$, $c_1\ne 0$. Also $c_2\ne 0$, since 
$\val P_\beta< \deg P_\beta$.  Set 
$$\delta\ :=\ \min_i v(a_i) + i\beta, \quad \fm\ :=\ s\beta,\quad \mathfrak d\ :=\ s\delta,
\quad Q(Y)\ :=\ \mathfrak d^{-1}P_{\times\fm}(Y) \in \mathcal{O}[Y],$$
so $Q$ has image $P_\beta$ in $\k[Y]$. 
Note that $$Q^{(d-1)}(Y)\ =\ \frac{\fm^{d-1}}{\mathfrak d}P^{(d-1)}(\fm Y),$$
with image $d!c_1(Y-c_2)$ in $\k[Y]$. It follows that the asymptotic equation 
$$ P^{(d-1)}(Y)\ =\ 0, \quad \quad Y\in\E $$ has dominant degree $1$ with 
primary dominant part $c(Y-c_2)$ for some~$c\in \k^\times$. 
Hence~(i) holds by Remark~\ref{rem:henselian, unique sol}. It also follows that $f=\fm g$ where $g$ is 
the unique zero of $Q^{(d-1)}$ in $\mathcal{O}$.
From $\bar{g}=c_2$ we obtain that $vf= \beta$ and that 
$f$ is an approximate zero of~\eqref{ast}. Since  
$\beta=\beta\eqref{ast}$ and $P_\beta(0)\ne 0$,
$\beta$ is the only antislope of~$\mathcal N(P)$ that is~$\geq v\fn$ for some $\fn\in\E$. Hence, if
$y$ is an approximate solution of \eqref{ast}, then
$vy=\beta$, so
$\bar{y/\fm}=c_2$ and $y\sim f$. 
As to (iii),  this follows from $P^{(d-1)}(f)=0$, that is,
the coefficient of $Y^{d-1}$ in $P_{+f}(Y)$ is $0$. 
\end{proof}

\noindent
The exceptional case treated in this lemma does actually occur: 

\begin{example}
Let $K=\R\(( t^\Q\)) $ with its usual valuation 
$v\colon K^\times \to \Gamma=\Q$, so
$vt=1$; see Section~\ref{sec:valued fields}. We  take $\fM=t^\Q$,
and identify the residue field $\k$ of $K$ with $\R$
in the usual way. We set 
%$x:=t^{-1}$ and
\begin{align*}
P(Y)\ &:=\ (g^2-t^2) -2gY + Y^2\ =\ (Y-g)^2-t^2,
\text{ where}\\ 
g\ &:=\ \sum_{k=1}^\infty t^{1-(1/k)}\ =\ 1+t^{\frac{1}{2}}+
t^{\frac{2}{3}}+t^{\frac{3}{4}}+\cdots \in K.
\end{align*}
Then ${\mathcal N}(P)=\bigl\{(0,0),\ (1,0),\ (2,0)\bigr\}$.
Let $\E:=\{y\in K^\times:\ y\prec t^{-1}\}$; then
the asymptotic equation \eqref{ast} 
has $\beta\eqref{ast}=0$, $\ddeg\eqref{ast}=2$,
$P_{\beta\eqref{ast}}(Y)=(Y-1)^2$ in $\R[Y]$, and $1$ is an
approximate solution of \eqref{ast}. The polynomial
$$P_{+1}(Y)\ =\ \big((g-1)^2-t^{2}\big)-2(g-1)Y+ Y^2\ =\ \big(Y-(g-1)\big)^2-t^{2}$$
has Newton diagram ${\mathcal N}(P_{+1})=
\bigl\{(0,1),\ (1,\frac{1}{2}),\ (2,0)\bigr\}$,
and the refinement 
$$P_{+1}(Y)\ =\ 0,\qquad Y\prec 1$$
of \eqref{ast} has dominant degree $2$ with primary dominant
part $(Y-1)^2$. 
In fact, for every $m\geq 1$ and $g_m:=\sum_{k=1}^m t^{1-(1/k)}$, 
the refinement
$$P_{+g_m}(Y)\ =\ 0,\qquad Y\prec t^{1-(1/m)}$$
of \eqref{ast} still has primary dominant part $(Y-1)^2$. On the other 
hand, 
the polynomial $P'(Y)=2(Y-g)$ has the unique zero $f:=g$ in $K$. We obtain 
$P_{+f}(Y)=Y^2-t^{2}$, $\beta\eqref{astphi}=1$, so the primary dominant part of \eqref{astphi}
is $Y^2-1$, with zeros $1$ and $-1$ in~$\R$, each of 
multiplicity $1$. See Figure~\ref{fig:newtondiag3}.
\end{example}

\begin{figure}[h]
\input{mt-newtondiag3.tex}
\caption{Unraveling an asymptotic equation.}\label{fig:newtondiag3}
\end{figure}

\index{complexity!asymptotic equation}
\index{asymptotic equation!complexity}
\nomenclature[K]{$c\eqref{ast}$}{complexity of \eqref{ast}}

\noindent
The {\bf complexity} of \eqref{ast} is the pair
$c\eqref{ast}=(d,\ell)\in \N\times\{0,1\}$ defined as follows: 
\begin{align*} d\ &=\ \text{dominant degree
of \eqref{ast}},\\
 \ell\ &=\ \begin{cases} 1 & 
\text{\parbox{25em}{if $d\ge 1$, $\beta\eqref{ast}\in \Gamma$, and $P_{\beta\eqref{ast}}$ has a zero in $\k$ of multiplicity $d$,}}\\ 
  0 & \text{otherwise.}
\end{cases}
\end{align*} 
In particular, if $d=0$, then $c\eqref{ast}=(0,0)$. 
We order 
the cartesian product $\N\times\{0,1\}$ lexicographically. In
Lemma \ref{Refinement-Lemma}, if $\beta \ne \beta\eqref{ast}$ or 
$c\eqref{ast}=(d,0)$ with $d\geq 1$,  
then $c\eqref{ast} > c\eqref{astphi}$. 
In Lemma~\ref{Special-Lemma} we have
$c\eqref{ast}=(d,1)$, and either $c\eqref{astphi}=(d,0)$ or $c\eqref{astphi}=(0,0)$; thus $c\eqref{ast}> c\eqref{astphi}$ in Lemma~\ref{Special-Lemma} as well. 

\medskip
\noindent
Passing to refinements \eqref{astphi} of \eqref{ast} for suitably chosen~$f$ in order to lower the complexity of \eqref{ast} in some sense \emph{unravels} the asymptotic equation. In Sections~\ref{sec:aseq}, \ref{sqe}, and \ref{sec:unravelers} we likewise try to unravel asymptotic \textit{differential}\/ equations.

\subsection*{Newton trees} In this subsection $K$ is henselian. The
lemmas above and the decrease
in complexity upon refinement
yield a family $(f_{\sigma})_{\sigma \in \Sigma}$
of elements $f_{\sigma}\in K^{\times}$ indexed by a finite (possibly empty) 
set $\Sigma$ of finite sequences $\sigma = (\sigma_1,\dots,\sigma_{k})$, where~$k\ge 1$ and each 
$\sigma_j$ is a positive integer, such that: 
\begin{enumerate}
\item $\Sigma$ is closed under taking initial segments of positive length, that is,
whenever $\sigma = (\sigma_1,\dots,\sigma_{k})\in \Sigma$ and $1\le j \le k$, then
$\sigma|j :=(\sigma_1,\dots,\sigma_{j})\in \Sigma$; 
\item  $\bigl\{i\in \N^{\geq 1}: (i)\in \Sigma\bigr\}= \{1,\dots,q\}$ 
with $q\in \N$; the elements $f_{(1)},\dots,f_{(q)}$ are 
approximate zeros of $P$, and any approximate zero of $P$ is equivalent to
  $f_{(i)}$ for exactly one $i\in \{1,\dots,q\}$;
\item for each $\sigma = (\sigma_1,\dots,\sigma_{k})\in \Sigma$ we have 
$$\bigl\{i\in \N^{\geq 1}:\ \sigma\ast i=(\sigma_1,\dots,\sigma_{k}, i)\in \Sigma\bigr\}\ =\ 
\{1,\dots,q_{\sigma}\}, \quad q_{\sigma}\in \N;$$ 
the elements $f_{\sigma\ast 1},\dots,f_{\sigma\ast q_{\sigma}}$ are 
approximate solutions of 
\begin{equation}\tag{E$_{\sigma}$}\label{phie}
P_{+f(\sigma)}(Y)\ =\ 0, \quad \quad Y\prec f_{\sigma}
\end{equation}
where $f(\sigma):=\sum_{j=1}^k f_{\sigma|j}$;
any approximate solution of \eqref{phie} is equivalent to
$f_{\sigma\ast i}$
for exactly one $i\in \{1,\dots,q_{\sigma}\}$;  and 
if $c \eqref{phie}= (d,1)$,
then $q_{\sigma}=1$ and $P_{+f(\sigma)}^{(d-1)}(f_{\sigma\ast 1})=0$.
\end{enumerate}
We call such a family $(f_\sigma)_{\sigma\in \Sigma}$
a {\bf Newton tree for $P$ in $K$}. (Strictly speaking, it is a forest
rather than a tree, but below we shall give it a root.) Note that 
items~(2) and (3) contain the instructions for growing a Newton tree 
for $P$ in $K$, 
which is in general not unique. Besides $K$
and the polynomial $P$ this notion of Newton tree also involves 
the ``monomial'' set $\fM$, but changing $\fM$ does not affect whether
the primary dominant part of \eqref{ast}  has a zero in 
$\k$ of multiplicity equal to its degree. Thus a Newton tree for $P$ in~$K$ 
with respect to $\fM$ remains a Newton tree for $P$ in~$K$ with another 
choice of $\fM$.

\index{Newton!tree}

\medskip\noindent
Let  $(f_\sigma)_{\sigma\in \Sigma}$ be
a Newton tree for $P$ in $K$. If
$\sigma=(\sigma_1,\dots,\sigma_{k})\in \Sigma$, then 
$$f_{\sigma|1}\succ f_{\sigma|2} \succ \dots \succ f_{\sigma|k}=f_\sigma,$$
so $f(\sigma) \sim f_{\sigma|1}$, $f(\sigma)$ is an approximate zero of $P$,
and whenever $y-f(\sigma) \prec f_\sigma$ and $1\le j < k$, then 
$y-f(\sigma|j)$ is an approximate solution of $(\operatorname{E}_{\sigma|j})$.
It follows that if $y$ is an approximate 
zero of $P$, then there is a unique $\sigma\in \Sigma$ such that 
$y-f(\sigma) \prec f_\sigma$ and 
$y-f(\sigma)$ is not an approximate zero of $P_{+f(\sigma)}$: take $\sigma\in \Sigma$
such that $y-f(\sigma) \prec f_\sigma$ and
$\sigma$ is not a proper initial segment
of any $\sigma'\in \Sigma$ with $y-f(\sigma') \prec f_{\sigma'}$.  

The case that $y$ is {\em not\/} an approximate zero of $P$ is put under this
roof as follows: put $\Sigma_0 := \Sigma\cup\{\emptyset\}$, set 
$f_\emptyset:=\infty$,  $f(\emptyset):=0$, and let
$y \prec \infty$ for each $y$, by convention. 
(Thus the forest $\Sigma$ becomes a 
single tree $\Sigma_0$ with root $\emptyset$.) 
Then there is for each $y$ a unique $\sigma\in \Sigma_0$
such that $y-f(\sigma) \prec f_\sigma$ and $y-f(\sigma)$ is not an 
approximate zero of~$P_{+f(\sigma)}$.  
See Figure~\ref{fig:newtontree}.

\begin{figure}[h]
\input{mt-newtontree.tex}
\caption{Newton tree.}\label{fig:newtontree}
\end{figure}

This leads to a piecewise
uniform description of $v\bigl(P(y)\bigr)$ in terms of functions of the
form $v\bigl(y-f(\sigma)\bigr)$, namely:

\begin{lemma} \label{piecewise uniform description}
Let $\sigma\in \Sigma_0$ and $P_{+f(\sigma)}= \sum_i b_iY^i$. If
$y - f(\sigma) \prec f_\sigma$ and $y-f(\sigma)$ is not an approximate zero of 
$P_{+f(\sigma)}$, then
$$v\bigl(P(y)\bigr)\ =\ 
\min \bigl\{v(b_i)+i\cdot v\bigl(y-f(\sigma)\bigr)\colon\ i=0,\dots,n\bigr\}
\qquad \text{\textup{(}$0\cdot \infty:= 0$ in $\Gamma_{\infty}$\textup{)}.}$$
\end{lemma}

\index{disk!special}
\index{disk!special!with holes}
\index{special!disk}
\index{special!disk with holes}

\noindent
For $\sigma\in \Sigma_0$, let $G_\sigma$ be the set of all $y$ 
satisfying the condition in this lemma:
$$G_\sigma\ :=\ \big\{y:\ \text{$y -f(\sigma) \prec f_\sigma$ and $y-f(\sigma)$ 
is not an approximate zero of $P_{+f(\sigma)}$}\big\}.$$  
For a more geometric description
of $G_\sigma$, define a {\bf special disk in $K$\/} to be 
a subset of $K$ that is either $K$ itself
or an open ball of the form $\big\{y: v(y-f) > vg\big\}$ where 
$f,g\in K^\times$ and $vf\le vg$. A {\bf special disk in $K$ with holes\/} is a
subset of $K$ of the form $D\setminus (D_1 \cup \dots\cup D_k)$
where $D$ is a special disk in $K$ and $D_1,\dots,D_k$ are disjoint 
special disks in $K$ properly contained in $D$. We can now summarize some
of the above as follows:

\begin{cor}\label{holes} For each $\sigma\in \Sigma_0$ the set
$G_\sigma$ is a special disk in $K$ with holes, and each $y$ 
belongs to $G_\sigma$ for exactly one $\sigma\in \Sigma_0$.
The zeros of $P$ in $K^{\times}$ are among the $f(\sigma)$ with $\sigma\in\Sigma$.
\end{cor}

\subsection*{Behavior under extension} Suppose $K$ is henselian, and 
let $L$ be a henselian valued field extension 
of $K$ such that the residue field of $K$ is algebraically closed 
in the residue field of $L$, and $\Gamma_L \cap \Gd =\Gamma$
(inside $\Q\Gamma_L$). Then any approximate zero of $P$ in
$L$ is equivalent to some element in~$K$.
It follows that a Newton tree
for $P$ in $K$ remains a Newton tree for $P$ in $L$. 

More generally, let $(f_\sigma)_{\sigma\in \Sigma}$ be a Newton tree for $P$ in $K$
and let $L$ be any henselian valued field extension of $K$. 
Then we can extend $\Sigma$ to a finite set $\Sigma_L\supseteq \Sigma$ of 
finite sequences of positive integers and extend $(f_\sigma)_{\sigma\in \Sigma}$ to
a Newton tree $(f_\sigma)_{\sigma\in \Sigma_L}$ of $P$ in $L$. 
This is because if a polynomial in $\k[Y]$ has degree $d\geq 1$ and has a zero of
multiplicity $d$ in an extension field of $\k$, then this zero
lies in~$\k$.

\subsection*{Notes and comments}
Newton diagrams were introduced by Newton~\cite{Newton} (1676) in constructing fractional power series solutions $y=y(x)$ to polynomial equations $P(x,y)=0$. Puiseux \cite{Puiseux50} rediscovered this in a complex-analytic setting:
the field of convergent Puiseux series with complex coefficients is algebraically closed; see \cite[\S{}17]{Abhyankar76},   \cite[Section~8.3]{BrieskornKnoerrer}, \cite[p.~396]{Chrystal}, and \cite[\S{}38]{HenselLandsberg}.
The use of derivatives as in the proof of Lemma~\ref{Special-Lemma} goes back to
Smith~\cite{Smith75}, but
Tschirnhaus~\cite{Tschirnhaus} already observed that for a monic
polynomial $P\in F[Y]$ of degree $n\geq 2$ over a field~$F$ and $f\in F$ with $P^{(n-1)}(f)=0$ one has
$P_{+f}=Y^n+\text{(terms of degree $\leq n-2)$}$. 
Dumas~\cite{Dumas} applied Newton diagrams to algebraic equations over 
the field of $p$-adic numbers, and 
Ostrowski~\cite{Ostrowski33, Ostrowski} and Rella~\cite{Rella} extended this to more general henselian valued fields.

Early uses of Newton diagrams in the study of algebraic differential equations
over the differential field $\operatorname{P}(\C)$  of Puiseux series over $\C$ are by Briot-Bouquet~\cite{BriotBouquet}
and by Fine~\cite{Fine89,Fine90}. 
For recent uses of Newton diagrams in connection with
differential equations over $\C[[x^{\R}]]$ and various differential subfields of it we refer to 
\cite{Cano93, DellaDoraDicrescenzoTournier, Farto, GrigorievSinger, Malgrange,Ramis}.

The Newton diagrams for differential polynomials in our Chapters~\ref{ch:The Dominant Part and the Newton Polynomial} and~\ref{ch:newtonian fields} do not seem to be related to the
Newton polyhedrons
for polynomials in several indeterminates over a discrete valued field in~\cite{Smirnov,Thaler}.

%% file: mt-newtondiag1.tex
\vskip1em

\begin{tikzpicture}
\coordinate (Origin)   at (0,0);
    \coordinate (XAxisMin) at (-0.5,0);
    \coordinate (XAxisMax) at (7.5,0);
    \coordinate (YAxisMin) at (0,-3.5);
    \coordinate (YAxisMax) at (0,4.5);
	\clip (-1.5,-3.5) rectangle (7.5,4.5);
    \draw [thick,-latex] (XAxisMin) -- (XAxisMax);% Draw x axis
    \draw [thick,-latex] (YAxisMin) -- (YAxisMax);% Draw y axis
    \foreach \i in {1,...,7} {
\draw[style=help lines, dashed] (\i,-4)--(\i,4.5);
}
  %  \draw[style=help lines,dashed] (-0.5,-4) grid[step=1] (8,4.5);
 \draw (0,-0.05) node [below, rectangle, fill=white] {\tiny $0$};  
 \draw (1,-0.05) node [below, rectangle, fill=white] {\tiny $1$};
 \draw (2,-0.05) node [below, rectangle, fill=white] {\tiny $2$};
 \draw (3,-0.05) node [below, rectangle, fill=white] {\tiny $3$};
 \draw (4,-0.05) node [below, rectangle, fill=white] {\tiny $4$};
 \draw (5,-0.05) node [below, rectangle, fill=white] {\tiny $5$};
 \draw (6,-0.05) node [below, rectangle, fill=white] {\tiny $6$};
 \draw (7.5,-0.25) node [below left] {$\Z$};
  \draw (0.25,4.25) node [right] {$\Q\Gamma$};
\node[draw,circle,inner sep=0.2em,fill] at (0,2) {};
\node[draw,circle,inner sep=0.2em,fill] at (1,-1) {};
\node[draw,circle,inner sep=0.2em,fill] at (2,-2) {};
\node[draw,circle,inner sep=0.2em,fill] at (3,-1.5) {};
\node[draw,circle,inner sep=0.2em,fill] at (5,-2) {};
\node[draw,circle,inner sep=0.2em,fill] at (6,1.5) {};
\draw[thick, dashed] (1.5,-2) -- (5.5,-2);
\draw[thick, dashed, domain=4.9:6.10] plot (\x, {-2+3.5*(\x-5)});
\draw[thick, dashed, domain=0.65:2.25] plot (\x, {-2-1*(\x-2)});
\draw[thick, dashed, domain=-0.25:1.25] plot (\x, {-1-3*(\x-1)});

\draw (0.1,2.1) node [right, rectangle, fill=white] {\small $(0,v(a_0))$};
\draw (1,-0.6) node [right, rectangle, fill=white] {\small $(1,v(a_1))$};
\draw (2,-2.3) node [below, rectangle, fill=white] {\small $(2,v(a_2))$};
\draw (3,-1.4) node [above, rectangle, fill=white] {\small $(3,v(a_3))$};
\draw (5.25,-2) node [right, rectangle, fill=white] {\small $(5,v(a_5))$};
\draw (6.1,1.5) node [right, rectangle, fill=white] {\small $(6,v(a_6))$};
\draw (5.5,-2.7) node [right, rectangle, fill=white] {$\mathcal N(P)$};
\end{tikzpicture}

%% file: mt-newtondiag2.tex
\vskip1em
\begin{tikzpicture}[scale=0.65]

\coordinate (Origin)   at (0,0);
    \coordinate (XAxisMin) at (-0.5,0);
    \coordinate (XAxisMax) at (7.5,0);
    \coordinate (YAxisMin) at (0,-3.5);
    \coordinate (YAxisMax) at (0,4.5);
	\clip (-1.5,-3.5) rectangle (7.5,4.5);
    \draw [thick,-latex] (XAxisMin) -- (XAxisMax);% Draw x axis
    \draw [thick,-latex] (YAxisMin) -- (YAxisMax);% Draw y axis
%        \foreach \i in {1,...,7} {
%\draw[style=help lines, dashed] (\i,-4)--(\i,4.5);
%}
 \draw[style=help lines,dashed] (-0.5,-4) grid[step=1] (8,4.5);
 \draw (0,-0.1) node [below, rectangle, fill=white] {\tiny $0$};  
 \draw (1,-0.1) node [below, rectangle, fill=white] {\tiny $1$};
 \draw (2,-0.1) node [below, rectangle, fill=white] {\tiny $2$};
 \draw (3,-0.1) node [below, rectangle, fill=white] {\tiny $3$};
 \draw (4,-0.1) node [below, rectangle, fill=white] {\tiny $4$};
 \draw (5,-0.1) node [below, rectangle, fill=white] {\tiny $5$};
 \draw (6,-0.1) node [below, rectangle, fill=white] {\tiny $6$};
 \draw (7.7,-0.25) node [below left] {$\Z$};
  \draw (-0.05,3.5) node [right] {$\Q$};
\node[draw,circle,inner sep=0.2em,fill] at (0,2) {};
\node[draw,circle,inner sep=0.2em,fill] at (1,-1) {};
\node[draw,circle,inner sep=0.2em,fill] at (2,-2) {};
\node[draw,circle,inner sep=0.2em,fill] at (3,-1.5) {};
\node[draw,circle,inner sep=0.2em,fill] at (5,-2) {};
\node[draw,circle,inner sep=0.2em,fill] at (6,1.5) {};
\draw[thick, dashed] (2,-2) -- (5,-2);
\draw[thick, dashed, domain=5:6] plot (\x, {-2+3.5*(\x-5)});
\draw[ultra thick, domain=1:2] plot (\x, {-2-1*(\x-2)});
\draw[ultra thick, domain=0:1] plot (\x, {-1-3*(\x-1)});
\draw[thick, dotted, domain=-0.5:5.5] plot (\x, {-2-0.5*(\x-2)});
\draw[thick, ->] (2.25,1) .. controls (1.7,-0.2) .. (1.5, -1.3);
\draw (0.9,1) node [right, rectangle, fill=white] {antislope $\beta\eqref{ast}=1$};
\draw (5.5,-2.5) node [right, rectangle, fill=white] {$\mathcal N(P)$};
\end{tikzpicture}
\begin{tikzpicture}[scale=0.65]

\coordinate (Origin)   at (0,0);
    \coordinate (XAxisMin) at (-0.5,0);
    \coordinate (XAxisMax) at (7.5,0);
    \coordinate (YAxisMin) at (0,-3.5);
    \coordinate (YAxisMax) at (0,4.5);
	\clip (-1.5,-3.5) rectangle (7.5,4.5);
    \draw [thick,-latex] (XAxisMin) -- (XAxisMax);% Draw x axis
    \draw [thick,-latex] (YAxisMin) -- (YAxisMax);% Draw y axis
    \draw[style=help lines,dashed] (-0.5,-4) grid[step=1] (8,4.5);
%    \foreach \i in {1,...,7} {
%\draw[style=help lines, dashed] (\i,-4)--(\i,4.5);
%}
 \draw (0,-0.1) node [below, rectangle, fill=white] {\tiny $0$};  
 \draw (1,-0.1) node [below, rectangle, fill=white] {\tiny $1$};
 \draw (2,-0.1) node [below, rectangle, fill=white] {\tiny $2$};
 \draw (3,-0.1) node [below, rectangle, fill=white] {\tiny $3$};
 \draw (4,-0.1) node [below, rectangle, fill=white] {\tiny $4$};
 \draw (5,-0.1) node [below, rectangle, fill=white] {\tiny $5$};
 \draw (6,-0.1) node [below, rectangle, fill=white] {\tiny $6$};
 \draw (7.7,-0.25) node [below left] {$\Z$};
  \draw (-0.05,3.5) node [right] {$\Q$};
\node[draw,circle,inner sep=0.2em,fill] at (0,3/2) {};
\node[draw,circle,inner sep=0.2em,fill] at (1,-1) {};
\node[draw,circle,inner sep=0.2em,fill] at (2,-2) {};
\node[draw,circle,inner sep=0.2em,fill] at (3,-1.5) {};
\node[draw,circle,inner sep=0.2em,fill] at (4,-1) {};
\node[draw,circle,inner sep=0.2em,fill] at (5,-2) {};
\node[draw,circle,inner sep=0.2em,fill] at (6,3/2) {};
\draw[thick, dashed] (1,-1) -- (2,-2) -- (5,-2) -- (6,1.5);

%\draw[thick, dashed, domain=1:2] plot (\x, {-1*(\x-1)});
\draw[ultra thick, domain=0:1] plot (\x, {1.5-2.5*(\x)});
\draw[ thick, dotted, domain=-0.5:5.5] plot (\x, {-1-1*(\x-1)});
 \draw (1.4,2.3) node [right, rectangle, fill=white] {antislope $\beta\eqref{astphi}>1$};
\draw[thick, ->] (1.8,1.8) -- (0.5, 1);
  \draw (2.95,-2.6) node [right, rectangle, fill=white] {antislope $\beta\eqref{ast}=1$};
\draw[thick, ->] (5.25,-3) .. controls (4.5,-3.5) .. (3.5, -3.2);
\draw (5.4,-1.75) node [right, rectangle, fill=white] {$\mathcal N(P_{+f})$};
\end{tikzpicture}

%% file: mt-newtondiag3.tex
\vskip1em

\begin{tikzpicture}

\coordinate (Origin)   at (0,0);
    \coordinate (XAxisMin) at (-0.5,0);
    \coordinate (XAxisMax) at (3.5,0);
    \coordinate (YAxisMin) at (0,-0.5);
    \coordinate (YAxisMax) at (0,3.5);
	\clip (-1.5,-0.75) rectangle (3.5,3.5);
    \draw [thick,-latex] (XAxisMin) -- (XAxisMax);% Draw x axis
    \draw [thick,-latex] (YAxisMin) -- (YAxisMax);% Draw y axis
    \draw[style=help lines,dashed] (-0.5,-4) grid[step=1] (8,4.5);
 \draw (0,-0.15) node [below, rectangle, fill=white] {\tiny $0$};  
 \draw (1,-0.15) node [below, rectangle, fill=white] {\tiny $1$};
 \draw (2,-0.15) node [below, rectangle, fill=white] {\tiny $2$};
 \draw (3,-0.15) node [below, rectangle, fill=white] {\tiny $3$};
 %\draw (4.05,0) node [below right] {\tiny $4$};
 %\draw (5.05,0) node [below right] {\tiny $5$};
 %\draw (6.05,0) node [below right] {\tiny $6$};
 \draw (3.6,-0.25) node [below left] {$\Z$};
  \draw (0.25,3.3) node [right] {$\Q$};
\node[draw,circle,inner sep=0.2em,fill] at (0,0) {};
\node[draw,circle,inner sep=0.2em,fill] at (1,0) {};
\node[draw,circle,inner sep=0.2em,fill] at (2,0) {};

%\draw[thick, dashed, domain=1:2] plot (\x, {-1*(\x-1)});
\draw[ultra thick, domain=0:1] (0,0) -- (2,0);

\draw (2.2,1) node [right, rectangle, fill=white] {$\mathcal N(P)$};
\end{tikzpicture}
\begin{tikzpicture}

\coordinate (Origin)   at (0,0);
    \coordinate (XAxisMin) at (-0.5,0);
    \coordinate (XAxisMax) at (3.5,0);
    \coordinate (YAxisMin) at (0,-0.5);
    \coordinate (YAxisMax) at (0,3.5);
	\clip (-1.5,-0.75) rectangle (3.5,3.5);
    \draw [thick,-latex] (XAxisMin) -- (XAxisMax);% Draw x axis
    \draw [thick,-latex] (YAxisMin) -- (YAxisMax);% Draw y axis
    \draw[style=help lines,dashed] (-0.5,-4) grid[step=1] (8,4.5);
 \draw (0,-0.15) node [below, rectangle, fill=white] {\tiny $0$};  
 \draw (1,-0.15) node [below, rectangle, fill=white] {\tiny $1$};
 \draw (2,-0.15) node [below, rectangle, fill=white] {\tiny $2$};
 \draw (3,-0.15) node [below, rectangle, fill=white] {\tiny $3$};
 %\draw (4.05,0) node [below right] {\tiny $4$};
 %\draw (5.05,0) node [below right] {\tiny $5$};
 %\draw (6.05,0) node [below right] {\tiny $6$};
 \draw (3.6,-0.25) node [below left] {$\Z$};
  \draw (0.25,3.3) node [right] {$\Q$};
  \draw[ultra thick, domain=0:1] (0,2) -- (2,0);

\node[draw,circle,inner sep=0.2em,fill] at (0,2) {};
%\node[draw,circle,inner sep=0.2em,fill=white] at (1,1) {};
\node[draw,circle,inner sep=0.2em,fill] at (2,0) {};

%\draw[thick, dashed, domain=1:2] plot (\x, {-1*(\x-1)});

\draw (2.1,1) node [right, rectangle, fill=white] {$\mathcal N(P_{+f})$};
\end{tikzpicture}

%% file: mt-newtontree.tex
\vskip1em
\begin{tikzpicture}
\tikzstyle{level 1}=[sibling distance=8em] \tikzstyle{level 2}=[sibling distance=2.2em] \tikzstyle{level 3}=[sibling distance=3em]
\node {$f_{\emptyset}$} 
child{ node {$f_{(1)}$} 
child{edge from parent[dotted]}  child{  edge from parent[draw=none]} child{ edge from parent[dotted]}}
child{ node {$f_{(2)}$}  child{ node {\small $f_{(2,1)}$} child{edge from parent[dotted]} child{node {$f_{(2,1,\dots,\sigma_k)}$} edge from parent[dashed]} child{ edge from parent[dotted]}} child{node {\small $f_{(2,2)}$}}  child{ node {$\dots\ $} edge from parent[draw=none]} child{node {\small $f_{(2,q_{(2)})}$}}  }  
child{ node {$\dots$} edge from parent[draw=none] } 
child{ node {$\small f_{(q)}$}   child{ edge from parent[dotted]}  child{   edge from parent[draw=none]} child{ edge from parent[dotted]  } };
\end{tikzpicture}

%% file: mt-4.tex
\chapter{Differential Polynomials}\label{ch:differential polynomials}

\noindent
Our differential fields are of characteristic zero with one 
distinguished derivation. We prove here some basic facts about
these
differential fields and their differential field extensions. In our work we often decompose differential polynomials into parts
of a special form and operate formally on differential
polynomials in various ways: additive and multiplicative conjugation, Ritt division, composition. Here we study these decompositions
and operations in their natural setting.
% except for the Riccati transform and the key operation of
%compositional conjugation, which are considered in the next chapter.

We also consider valued differential fields, the
property of continuity of the derivation with respect to the valuation topology, and the gaussian extension of the valuation to the 
ring of differential polynomials. Valued differential fields will be studied further
in Chapter~\ref{ch:valueddifferential}. We finish this chapter with some basic results on
simple differential rings, and on differentially closed fields.

\input{mt-4-1}

\input{mt-4-2}

\input{mt-4-3}
\input{mt-4-4}

\input{mt-4-5}

\input{mt-4-6}

\input{mt-4-7}

%% file: mt-4-1.tex
\section{Differential Fields and Differential Polynomials}
\label{Differential Fields and Differential Polynomials}

\noindent
Recall that derivations on commutative rings were introduced in 
Section~\ref{sec:modules}. When we say below that a commutative ring $K$ {\em contains $\Q$\/} we are abusing language: the meaning of this phrase is that there exists a (necessarily unique) ring embedding of the field $\Q$ into $K$; as usual we identify $\Q$ in that case with its image in 
$K$ under this embedding. 

\subsection*{Differential rings} 
A {\bf differential ring} is by definition a commutative ring $K$ 
containing $\Q$,
equipped with a derivation $\der$ on $K$; when $\der$ is clear from 
the context, we set $a':=\der(a)$, and similarly, $a^{(n)}:=\der^n(a)$, with
$\der^n$ the $n$th iterate of $\der$.
The Leibniz identity $(ab)'=a'b+ab'$ then extends as follows:
for elements $a_1,\dots,a_m$ ($m\geq 1$) of a differential ring,
$$%\begin{equation}\label{eq:Leibniz}
(a_1\cdots a_m)^{(n)}\ =\ \sum \frac{n!}{i_1!\cdots i_m!}\, a_1^{(i_1)}\cdots a_m^{(i_m)},
$$%\end{equation}
where the sum on the right is over all 
$(i_1,\dots,i_m)\in\N^m$ with $i_1+\cdots+i_m=n$. Given a differential ring $K$
and an element $a$ of a differential ring extension $L$ of $K$ we let 
$K\{a\}:=K[a, a', a'',\dots]$ denote the smallest subring of $L$ containing $K\cup\{a\}$ and closed under $\der$: {\em the differential ring generated
by $a$ over $K$}. 

\index{differential ring}
\index{ring!differential}
\index{ring!constants}
\index{differential ring!ring of constants}
%\index{derivative}
\nomenclature[M]{$a^{(n)}$}{$n$th derivative $a^{(n)}=\der^n(a)$ of $a$}
\nomenclature[M]{$K\{a\}$}{differential ring generated by $a$ over $K$}

\medskip\noindent
Let $K$ be a differential ring.
The subring $\{a\in K: a'=0\}$ of $K$ contains $\Q$. It is called the 
{\bf ring of constants} of $K$, and denoted by $C_K$ (or just $C$ if $K$ is clear from the context). If $c\in C$,
then $(ca)'=ca'$ for all $a\in K$. For $a\in K^\times$ we put $a^\dagger:= a'/a$, the {\bf logarithmic derivative\/} \index{derivative!logarithmic}\nomenclature[M]{$a^\dagger$}{logarithmic derivative $a^\dagger= a'/a$ of a unit $a$} of $a$. Let $a,b\in K^\times$; then 
$(ab)^\dagger = a^\dagger + b^\dagger$, in particular, $(1/a)^\dagger=-a^\dagger$,
and $a^\dagger=b^\dagger$ iff $a=bc$ for some $c\in C$.

\subsection*{Localization} 
Let $R$ be a differential ring and $S$ a multiplicative subset of $R$ with~$0\notin S$.
Then there is a 
unique derivation on $S^{-1}R$ making $S^{-1}R$ into a differential ring
such that the natural map $r\mapsto r/1: R\to S^{-1}R$ is a morphism of differential rings; 
it is given by
$$(r/s)'\  =\  (r's-rs')/s^2\quad\text{for $r\in R$, $s\in S$, }$$
and we always consider $S^{-1}R$ as a differential ring in this way.
In particular, if $R$ is a  differential integral domain (that is, 
a differential ring whose underlying
ring is an integral domain), then the derivation~$\der$ of~$R$ extends uniquely to a derivation on the fraction field of $R$.

\index{localization}

%\begin{lemma}\label{lem:derivation into fraction field} 

\subsection*{Differential fields}
A {\bf differential field} is a
differential ring whose underlying ring is a field (of characteristic $0$
since it contains $\Q$).
If $K$ is a differential field, then $C=C_K=\{c\in K:c'=0\}$
is a subfield of $K$, called the {\bf field of constants} of 
$K$. \index{field!differential}\index{differential field}\index{differential field!field of constants}\index{field!constants}\nomenclature[M]{$C_K$}{field of constants of the differential field $K$} If $K$ is a differential field and $L$ is a field extension of $K$ that is algebraic over~$K$, then we always consider $L$ as the differential field extension of $K$ obtained as in Lemma~\ref{lem:derivative on alg ext} by extending the derivation of $K$ to a derivation on $L$.

\label{p:differential field}

\begin{lemma}\label{lem:C alg closed in K}
Let $K$ be a differential field. Any element of $K$ 
that is algebraic over~$C$ lies in $C$. Thus if $K$ is algebraically closed as a field, then so is $C$, and if~$K$ is real closed, then so is $C$.
\end{lemma}
\begin{proof}
Suppose $a\in K$ is algebraic over~$C$, with minimum polynomial $P\in C[X]$ over~$C$. Then $P'(a)\neq 0$ since $\operatorname{char}(C)=0$, and by
the case $n=1$ of Lemma~\ref{derpol} we obtain $0=P'(a)a'$, so $a'=0$, hence $a\in C$.
\end{proof}

%\begin{lemma}\label{lem:derivative on alg ext}

\begin{lemma}\label{lem:const field of algext} Let $K$ be a differential field, $L$ a differential field extension of $K$, and suppose $a\in C_L$ is algebraic over~$K$.
Then $a$ is algebraic over~$C$.
\end{lemma}
\begin{proof} Let $P(X)$ be the minimum polynomial of $a$ over~$K$.
Then $$P(a)'\ =\ P^{\der}(a)+ P'(a)\cdot a'\ =\ P^{\der}(a)\ =\ 0,$$ with $\deg P^{\der}<\deg P$, so $P^{\der}=0$, and thus $P\in C[X]$.
\end{proof}

\noindent
Let $K$ be a field of characteristic zero. Then $P\mapsto P'=\partial P/\partial X$ is a
derivation on the
ring $K[X]$ of polynomials in the indeterminate $X$ over~$K$. This
derivation extends to a derivation $R\mapsto \partial R/\partial X$ on the fraction field $K(X)$, with $K$ as its field of constants.
Taking logarithmic derivatives gives the following:

\begin{lemma}\label{lem:log der of rational fn}
Let $R=a\cdot(X-b_1)^{k_1}\cdots (X-b_n)^{k_n}$ where $a\in K^\times$,
$b_1,\dots,b_n\in K$, $k_1,\dots,k_n\in\Z$. Then
$$\frac{\partial R/\partial X}{R}\ =\ \sum_{i=1}^n \frac{k_i}{X-b_i}.$$
\end{lemma}

\noindent
Assume next that $K$ is a differential field. Then its derivation
extends to the derivation $P\mapsto P^\der$ on $K[X]$ with $X^{\der}=0$, and we extend this further to a
derivation $R\mapsto R^\der$ on~$K(X)$. With these notations
we get by a routine computation: 

\begin{cor}\label{cor:derpol}
Let $P,Q\in K[X]$ with $Q\neq 0$, and set $R:=P/Q\in K(X)$. Then
$$R(a)'\ =\ R^\der(a) + (\partial R/\partial X)(a)\cdot a'\  \text{ for $a\in K$  with $Q(a)\ne 0$.}$$
\end{cor}

%(This follows easily from the case $n=1$ of Lemma~\ref{derpol}.) 

%\noindent
%The following lemma is used in the proof of Lemma~\ref{compconjval} below:

%\begin{lemma}\label{closed under der}
%Let $K$ be a valued differential field, $\phi\in K$, and $v(\phi')\geq v\phi$. 
%Then for every $n$, the $\mathcal{O}$-submodules
%$\phi^n\mathcal{O}$ and $\phi^n \mathfrak{m}(\mathcal{O})$ of $K$ are 
%closed under $\der$.
%\end{lemma}
%\begin{proof} We can assume $n>0$ and $\phi\ne 0$. Let
%$a = b \phi^n$ where $b\in \mathcal{O}$. Then
%$a'=b'\phi^n + nb\phi^{n-1}\phi'$ with $v(b'\phi^n), 
%v(b\phi^{n-1}\phi') \geq nv\phi$, with strict inequality if 
%$b\in \mathfrak{m}(\mathcal{O})$. 
%Hence $a'\in\phi^n \mathcal{O}$, and 
%$a'\in\phi^n \mathfrak{m}(\mathcal{O})$ if $b\in \mathfrak{m}(\mathcal{O})$. 
%\end{proof}

\subsection*{Differential automorphisms}
Let $K$ be a differential ring. An {\bf automorphism} of $K$ is 
an automorphism $\sigma$ of the ring~$K$ such that $\sigma\circ\der=\der\circ\sigma$, with $\der$ the derivation of $K$.
Let $L$ be a differential ring extension of  $K$. 
The set of ring automorphisms~$\sigma$ of~$L$ with 
$\sigma(a)=a$ for all $a\in K$ 
is a group under composition, denoted by $\Aut(L|K)$. 
The set of automorphisms $\sigma$ of the differential ring $L$ 
with $\sigma(a)=a$ for all $a\in K$ 
is a subgroup of $\Aut(L|K)$, denoted by $\Aut_\der(L|K)$.
If both $K$ and~$L$ are differential fields and the field extension~$L|K$ is algebraic, then $\Aut_\der(L|K)=\Aut(L|K)$. If~$L$ is an integral domain and $F$ is the differential fraction field of $L$, then every
$\sigma\in\Aut_\der(L|K)$ extends uniquely to an automorphism $\sigma_F\in\Aut_\der(F|K)$, and the map
$\sigma\mapsto\sigma_F\colon \Aut_\der(L|K)\to \Aut_\der(F|K)$ is an embedding of groups.

\index{automorphism!differential ring}
\index{differential ring!automorphism}
\nomenclature[A]{$\Aut(L\lvert K)$}{group of automorphisms of $L$ over~$K$}
\nomenclature[M]{$\Aut_\der(L\lvert K)$}{group of automorphisms of a differential ring $L$ over~$K$}

\subsection*{Differential polynomials}
Let $K$ be a differential ring with derivation $\der$ and $Y$ a
differential indeterminate over~$K$. Then $K\{Y\}$ denotes the ring of
differential polynomials in $Y$ over~$K$. As a ring, $K\{Y\}$
is just the polynomial ring $K[Y,Y',Y'',\dots]$ in the distinct
 indeterminates~$Y^{(n)}$ (${n\in\N}$) over~$K$, where as usual we write $Y$, $Y'$, $Y''$ 
instead
of $Y^{(0)}$, $Y^{(1)}$, $Y^{(2)}$. We consider
$K\{Y\}$ as the differential ring whose derivation,
extending the derivation of~$K$ and also denoted by
$\der$, is given by $\der(Y^{(n)})=Y^{(n+1)}$
for every $n$. 
Given $P\in K\{Y\}$, the smallest $r\in\N$ such that 
$P\in K\big[Y,Y',\dots,Y^{(r)}\big]$ is called the
{\bf order} of the differential polynomial $P$, and the
{\bf degree} $\deg P$ of $P$ is its
(total) degree as an element of the polynomial ring $K[Y,Y',\dots]$
(with $\deg 0= -\infty$). It is easy to check that if $P\in K\{Y\}$ with
$P\notin K$ has order $r$, then $P'$ has order $r+1$ and degree $1$
in $Y^{(r+1)}$. 
For $P\in K\{Y\}$ and~$y$ an element of a differential ring extension of
$K$, we let $P(y)$ be the element of that extension obtained
by substituting $y,y',\dots$ for $Y,Y',\dots$ in $P$,
respectively. In particular, we have $P=P(Y)$ for $P\in K\{Y\}$. 

\index{polynomial!differential}
\index{differential polynomial}
\index{differential polynomial!order}
\index{differential polynomial!degree}
\index{degree!differential polynomial}
\index{order!differential polynomial}

\nomenclature[M]{$K\{Y\}$}{ring of differential polynomials in  $Y$ over   $K$}
\nomenclature[O]{$\order(P)$}{order of  $P$}
\nomenclature[O]{$\deg(P)$}{degree of  $P$}

Let $L$ be a differential ring extension of $K$ and $y\in L$. Then the map
$$P \mapsto P(y)\colon\  K\{Y\} \to L$$ is a differential ring morphism, and is in 
fact the unique differential ring morphism $K\{Y\} \to L$
that is the identity on $K$ and sends $Y$ to $y$. With $L=K\{Y\}$ this gives
a composition operation $(P, Q) \mapsto P\circ Q:= P(Q)\colon K\{Y\}^2 \to K\{Y\}$. \nomenclature[O]{$P\circ Q=P(Q)$}{composition of  $P$ with $Q$} \index{composition!differential polynomials} From the uniqueness we easily get $(P\circ Q)(y)=P\big(Q(y)\big)$ for all 
$P,Q\in K\{Y\}$ and all $y$ in all 
differential ring extensions of $K$, and thus the associativity of 
composition:
$$(P\circ Q)\circ R\ =\ P\circ (Q\circ R) \quad \text{ for $P,Q,R\in K\{Y\}$.}$$
Let  $Y_1,\dots,Y_m$ be distinct differential indeterminates.
We define the differential ring $K\{Y_1,\dots,Y_m\}$ inductively by letting it be 
$K$ for $m=0$ and setting 
$$K\{Y_1,\dots,Y_m\}\ :=\
K\{Y_1,\dots,Y_{m-1}\}\{Y_m\} \qquad (m>0).$$
 This is a domain if $K$ is a domain.
{\em In the rest of this subsection $K$ is a differential field}. The fraction field of 
$K\{Y_1,\dots,Y_m\}$ is denoted
by $K\langle Y_1,\dots,Y_m\rangle$, is given the unique derivation that 
extends the derivation of $K\{Y_1,\dots,Y_m\}$, and
is called the field of {\bf differential rational
functions} in $Y_1,\dots,Y_m$ with coefficients in $K$.

\index{differential rational function}
\index{rational!differential function}
\nomenclature[M]{$K\{Y_1,\dots,Y_m\}$}{ring of differential polynomials in  $Y_1,\dots,Y_m$ over~$K$}
\nomenclature[M]{$K\langle Y_1,\dots,Y_m\rangle$}{field of differential rational functions in  $Y_1,\dots,Y_m$ with coefficients in the differential field $K$}

\medskip\noindent
Let $L$ be a differential field extension of $K$ and $y\in L$.
Then 
$$K\langle y\rangle\ :=\ K(y,y',y'',\dots)$$ denotes the differential
subfield of $L$ generated over~$K$ by $y$. Likewise,
$K\langle y_1,\dots,y_{m}\rangle$ is the differential subfield
of $L$ generated over~$K$ by elements $y_1,\dots,y_{m}\in L$. We say that
$y$ is {\bf differentially algebraic} over~$K$ if $P(y)=0$ for some
$P\in K\{Y\}^{\ne}$, equivalently, $K\langle y\rangle$ has finite 
transcendence degree over~$K$. We also say that $L$ is 
{\bf differentially algebraic} over~$K$ if each element of $L$ is differentially algebraic over~$K$. If $L=K\langle y_1,\dots,y_{m}\rangle$
and each $y_i\in L$ is differentially algebraic over~$K$, then $L$ is
differentially algebraic over~$K$. If $L$ is
differentially algebraic over~$K$ and~$M$ is a differential field extension of
$L$ that is differentially algebraic over~$L$, then~$M$ is
differentially algebraic over~$K$. If $y$ is not differentially algebraic over~$K$,
then $y$ is said to be {\bf differentially transcendental} over~$K$. 
Thus the differential indeterminate~$Y$, as an element of 
$K\langle Y \rangle$, is differentially transcendental over~$K$. 
\index{differentially!algebraic}
\nomenclature[M]{$K\langle y_1,\dots,y_{m}\rangle$}{differential field generated by $y_1,\dots,y_m$ over~$K$}
\index{differentially!transcendental} 
The equality $Y'=Y\cdot Y^\dagger$ shows that $Y$ is differentially 
algebraic over~$K\langle Y^\dagger \rangle$, so~$Y^\dagger$ must be 
differentially transcendental over~$K$ as well. More generally:

\begin{lemma}\label{dtrdalg} Suppose $y\in L$ is $\d$-transcendental over~$K$, and $z\in K\<y\>,\ z\notin K$.
Then $y$ is $\d$-algebraic over~$K\<z\>$, and so
$z$ is $\d$-transcendental over~$K$. 
\end{lemma}

\noindent
Here and below we use the prefix $\d$ to abbreviate {\em differential\/} or {\em differentially.}\/ Thus \textit{$\d$-algebraic}\/ stands for \textit{differentially algebraic.}

\begin{proof} We have $z=P(y)/Q(y)$ where $P, Q\in K\{Y\}^{\ne}$ and $P\ne aQ$ for all $a\in K$. Set 
$F(Y):= P(Y)-zQ(Y)\in K\<z\>\{Y\}$. From $z\notin K$ we get 
$F\ne 0$. Then $F(y)=0$ yields that
$y$ is $\d$-algebraic over~$K\<z\>$.
\end{proof}

%We shall often use the prefix $\d$ to abbreviate {\em %differential\/} or {\em differentially}. Thus ``$\d$-%algebraic" stands for ``differentially algebraic.'' 

\subsection*{Minimal annihilators}
{\em In this subsection, $K$ is a differential field and $y$ is an element of a differential field extension of~$K$}.
A {\bf minimal annihilator\/} of~$y$ over~$K$ is an irreducible 
$P\in K\{Y\}^{\ne}$, say of order $r$, such that $P(y)=0$ and $Q(y)\ne 0$ 
for every $Q\in K\{Y\}^{\ne}$ of order at most~$r$
with $\deg_{Y^{(r)}}Q < \deg_{Y^{(r)}}P$.
  
\index{minimal!annihilator}
\index{annihilator!minimal}

\begin{lemma}\label{lem:min annihilator}
Let $P\in K\{Y\}^{\ne}$ of order $r$ be irreducible with $P(y)=0$ and
$Q(y)\neq 0$ for all $Q\in K\{Y\}^{\ne}$ of order~$<r$. Then $P$ is a minimal annihilator of $y$ over~$K$. 
\end{lemma}
\begin{proof}
Towards a contradiction,
assume that $P$ is not a minimal annihilator of $y$ over~$K$. Take $Q\in K\{Y\}^{\ne}$ of order $r$ with $Q(y)=0$ and $\deg_{Y^{(r)}}Q<\deg_{Y^{(r)}}P$.
So $d:=\deg_{Y^{(r)}}Q\geq 1$, and we
choose $Q$ with these properties such that $d$ is minimal. By polynomial division
in   $K\big[Y,Y',\dots,Y^{(r)}\big]$ we obtain $A\in K\big[Y,Y',\dots,Y^{(r-1)}\big]{}^{\ne}$ and $B\in K\big[Y,Y',\dots,Y^{(r)}\big]$  with
$AP=BQ+R$ where $\deg_{Y^{(r)}} R<d$. Then $R(y)=0$, so $R=0$ by minimality of $d$, hence
$P$ divides $BQ$. So $P$ divides $B$, since $P$ is irreducible and $\deg_{Y^{(r)}}Q<\deg_{Y^{(r)}}P$.
Hence $B\in P\,K\{Y\}$, so $AB\in AP\,K\{Y\}=BQ\,K\{Y\}$ and thus
$0\neq A\in Q\,K\{Y\}$, contradicting   $\order(A)<r=\order(Q)$. 
\end{proof}

\begin{cor}\label{cor:min annihilator}
Suppose that $y$ is $\d$-algebraic over~$K$. Then $y$ has a minimal annihilator over~$K$. 
Such a minimal annihilator of $y$ is unique up to multiplication by a factor
from~$K^\times$. 
\end{cor}
\begin{proof}
Take $P\in K\{Y\}^{\neq}$ of minimal order $r$ such that $P(y)=0$. Replacing $P$
by some factor we arrange that $P$ is irreducible.
Then by the preceding lemma, $P$ is a  minimal annihilator of~$y$ over~$K$.
Let $P^*\in K\{Y\}^{\ne}$ also be a minimal annihilator of~$y$ over~$K$. Then
$r=\order(P)=\order(P^*)$ and $d:=\deg_{Y^{(r)}} P = \deg_{Y^{(r)}} P^*$, so we 
have
$A,A^*\in K\{Y\}^{\ne}$ of order $<r$ such that $Q:=AP+A^*P^*$ has degree~$<d$ in $Y^{(r)}$. Then $Q(y)=0$, so $Q=0$, and thus by irreducibility of $P$ and $P^*$ we get 
$P^*=aP$ with $a\in K^\times$.
\end{proof}

\begin{cor}\label{cor:order lower bound}
Let $y$ be $\d$-algebraic over~$K$ with minimal annihilator 
$P\in K\{Y\}$ over~$K$. Suppose $Q\in K\{Y\}^{\neq}$ is such
that $Q(y)=0$ and $\order(Q)\leq \order(P)$. Then 
$Q\in P\,K\{Y\}$, and hence $\order(Q) = \order(P)$.
\end{cor}
\begin{proof} Replacing $Q$ by a suitable factor we arrange that
$Q$ is irreducible, and thus a minimal annihilator of $y$ over~$K$ by Lemma~\ref{lem:min annihilator}. Now apply
Corollary~\ref{cor:min annihilator}. 
\end{proof}

\noindent
Let $P\in K\{Y\}^{\ne}$ be irreducible. Then there is an element $a$ in a
differential field extension of $K$ with $P$ as a minimal annihilator over~$K$.
For any such $a$, if $b$ in a
differential field extension of $K$ also has $P$ as a minimal annihilator over 
$K$, then there is a differential field isomorphism $K\<a\> \to K\<b\>$ over~$K$
that sends $a$ to $b$.

The following two lemmas contain the proofs of these facts, which we include here since they will serve as templates for similar more involved proofs later, for differential fields equipped with additional structure. 

\begin{lemma}\label{lem:construct diff alg extension}
Let $P\in K\{Y\}^{\neq}$ be irreducible. Then there is an element $a$ in a differential  field extension of $K$ such that $a$ is $\d$-algebraic over~$K$ with minimal annihilator~$P$ over~$K$.
\end{lemma}
\begin{proof}
Let $r=\order(P)$, and
take a polynomial $p\in K[Y_0,\dots,Y_r]$ with $P=p\big(Y,Y',\dots,Y^{(r)}\big)$, so $p$ is irreducible. Consider the integral domain
$$K[y_0,\dots,y_r]\ :=\ K[Y_0,\dots,Y_r]/(p),\qquad \text{$y_i:=Y_i+(p)$ for $i=0,\dots,r$,}$$
and let $F=K(y)$, where $y=(y_0,\dots,y_r)$, be its fraction field. Note that~${p(y)=0}$ and 
$q(y)\neq 0$ for $q\in K[Y_0,\dots,Y_r]$ with $\deg_{Y_r} q < \deg_{Y_r} p$; in particular, we have~${(\partial p/\partial Y_r)(y)\neq 0}$. We are going to extend $\der$ to a derivation on $F$ such that $y_i'=y_{i+1}$ for $i=0,\dots,r-1$. We first set
$$y_{r+1}\ :=\ -\frac{p^\der(y)+\sum_{i=0}^{r-1} (\partial p/\partial Y_i)(y)\cdot y_{i+1}}{(\partial p/\partial Y_r)(y)} \quad\text{in $K(y)$,}$$
which by Lemma~\ref{derpol} is the value that $y_r'$ will necessarily have, for any derivation on~$F$ extending $\der$ with $y_i'=y_{i+1}$ for $i=0,\dots,r-1$. Next we define the additive map $d\colon K[Y_0,\dots,Y_r]\to F$ by
$$d(f)\  =\  f^\der(y)+\sum_{i=0}^r \frac{\partial f}{\partial Y_i}(y)\cdot y_{i+1}.$$
As in the proof of Lemma~\ref{lem:derivative on alg ext} we check that the kernel of $d$ contains $(p)$ and that the induced additive map
$K[y]\to F$ is a derivation into $F$, which by Corollary~\ref{cor:derivation into fraction field} extends uniquely to a derivation on $F$.
This derivation extends $\der$, and setting $a:=y_0$ we have $a^{(i)}=y_i$ for $i=0,\dots,r$. This $a$ has the desired property.
\end{proof}

\begin{lemma}\label{lem:uniqueness of diff alg extension}
Let $a$ and $b$ in differential field extensions of $K$ be $\d$-algebraic over~$K$ with common minimal annihilator $P\in K\{Y\}$ over~$K$. Then there is a differential field isomorphism $K\<a\>\to K\<b\>$ over~$K$ sending $a$ to $b$.
\end{lemma}
\begin{proof}
Take $p\in K[Y_0,\dots,Y_r]$ such that $P=p(Y,Y',\dots,Y^{(r)})$, where $r:=\order(P)$.
Put $\vec{a}:=\big(a,a',\dots,a^{(r)}\big)$ and $\vec{b}:=\big(b,b',\dots,b^{(r)}\big)$. The ring morphism
$$f\mapsto f(\vec{a}\,)\ \colon\ K[Y_0,\dots,Y_r]\to K[\vec{a}\,]$$ has kernel $(p)$, and likewise with $b$ instead of $a$, so we have a ring isomorphism $K[\vec{a}\,]\to K[\vec{b}\,]$ that sends $a^{(i)}$ to $b^{(i)}$ for $i=0,\dots,r$. 
Since $P$ is a minimal annihilator of~$a$ over~$K$, we have $(\partial p/\partial Y_r)(\vec{a}\,)\neq 0$, so
Lemma~\ref{derpol} gives
$$a^{(r+1)}\ =\ -\frac{p^\der(\vec{a}\,)+\sum_{i=0}^{r-1} (\partial p/\partial Y_i)(\vec{a}\,)\cdot a^{(i+1)}}{(\partial p/\partial Y_r)(\vec{a}\,)}\in K(\vec{a}\,),$$
so $K[\vec{a}\,]'\subseteq K(\vec{a}\,)$ (as sets) and hence $K[\vec{a}]\subseteq K\{a\}\subseteq K(\vec{a}\,)$ (as rings). Likewise with $b$ instead of $a$, so the ring isomorphism $K[\vec{a}\,]\to K[\vec{b}\,]$ from above extends to a differential field isomorphism $K\<a\>\to K\<b\>$ as desired.
 \end{proof}

\index{differential polynomial!separant}

\noindent
Let $P\in K\{Y\}^{\neq}$ be irreducible of order $r\in\N$, and set 
$$ S\ :=\ \partial P/\partial Y^{(r)}\in K\big[Y, \dots, Y^{(r)}\big].$$ 
(In Section~\ref{Operations on Differential Polynomials} we call $S$ the {\em separant\/} of $P$.)  
Let $a$ be an element of a differential field extension of $K$ with minimal annihilator $P$ over~$K$.
The proofs above show that then $S(a)\neq 0$ and
that the subring $K\left[a,a',\dots,a^{(r)},\textstyle\frac{1}{S(a)}\right]$
of  $K\<a\>$
is closed under the derivation of $K\<a\>$. Therefore,
$$
K\{a\}\ \subseteq\ K\left[a,a',\dots,a^{(r)},\textstyle\frac{1}{S(a)}\right], \qquad K\<a\>\ =\ K\big(a,a',\dots,a^{(r)}\big)$$
with $a,a',\dots,a^{(r-1)}$ as a transcendence basis for the field extension $K\<a\>\supseteq K$.
Thus if $b$ is an element of a differential field extension of $K$ with $P(b)=0$ and $S(b)\neq 0$, then we have a differential ring morphism 
$K\left[a,a',\dots,a^{(r)},\textstyle\frac{1}{S(a)}\right]\to K\<b\>$ over~$K$ sending $a$ to $b$, which restricts to a differential ring morphism $K\{a\}\to K\{b\}$.

\begin{cor}\label{cor:order of min ann}
Let $y$ be $\d$-algebraic over~$K$. Then
$$\operatorname{trdeg}\!\big(K\<y\>|K\big)\ =\ \text{order of a minimal annihilator of $y$ over~$K$.}$$
\end{cor}

\subsection*{Differential transcendence bases}
Let $L\supseteq K$ be an extension of differential fields. A set 
$E\subseteq L$ generates over~$K$ the differential
field $K\<E\>\subseteq L$, and we put
$$\operatorname{cl} E\ :=\ \big\{f\in L:\  \text{$f$ is  $\d$-algebraic over 
$K\<E\>$}\big\}.$$ 
Then $\operatorname{cl} E$ is a differential subfield of $L$.

\begin{lemma}\label{lem:pregeom}
Let $E\subseteq L$ and $a,b\in L$. Then:
\begin{enumerate}
\item[\textup{(i)}] $E\subseteq\operatorname{cl} E$;
\item[\textup{(ii)}] $\operatorname{cl} E = \bigcup\{\operatorname{cl} E_0: \text{$E_0\subseteq E$ is finite} \}$;
\item[\textup{(iii)}] $\operatorname{cl}(\operatorname{cl} E)=\operatorname{cl} E$;
\item[\textup{(iv)}] $a\notin \operatorname{cl} E,\ b\notin \operatorname{cl}\big(E\cup\{a\}\big)\ 
\Rightarrow\ a\notin \operatorname{cl}\big(E\cup\{b\}\big)$.
\end{enumerate}
\end{lemma}
\begin{proof}
Parts (i)--(iii) are clear. For (iv), put 
$F:=  \operatorname{cl} E$, and
assume $a\notin F$ and~$b$ is $\d$-transcendental over~$F\<a\>$. Then
the family
$\big(a^{(m)}\big)$ is algebraically independent over~$F$ and the family $\big(b^{(n)}\big)$ is 
algebraically independent over~$F\<a\>$. It follows that their ``union''  
$\big(a^{(m)}, b^{(n)}\big)$ is
algebraically independent over~$F$. Hence $\big(a^{(m)}\big)$ is algebraically 
independent over~$F\<b\>$, that is, $a$ is $\d$-transcendental
over~$F\<b\>$.   
\end{proof}

\noindent
Lemma~\ref{lem:pregeom} says that the operation
$E\mapsto \operatorname{cl} E$ (for $E\subseteq L$) is a \textit{pregeometry}\/
on~$L$. This yields a notion of independence, analogous to linear
independence in vector spaces and algebraic independence for field extensions.
Below we just formulate this independence notion in our situation and state 
the relevant facts about it. (These facts are special cases of
generalities about pregeometries that can be found in many places,
for example in \cite[V, \S{}5, Exercice~14]{Bourbaki-Alg},
in \cite[VII, \S{}2]{Cohn65}, and in \cite[Section~8.1]{Marker}.) % \cite[Appendix~C]{TentZiegler}.  

Let $E$ range over subsets of $L$ in the rest of this subsection.
We say that~$E$ is {\bf $\d$-algebraically independent} over~$K$ if $x\notin \operatorname{cl}(E\setminus\{x\})$
for all $x\in E$.\index{differentially!algebraically independent}
\index{independent!d-algebraically@$\d$-algebraically} So~$E$ is $\d$-algebraically independent over~$K$ iff 
the family
$\big(x^{(n)}\big)_{x\in E,\ n=0,1,2,\dots}$ is algebraically independent over~$K$.
Call $E$ a {\bf $\d$-transcendence basis} of $L$ over~$K$ if~$E$ is 
$\d$-algebraically independent over~$K$ and $L$ is 
$\d$-algebraic over~$K\<E\>$.\index{differential transcendence!basis} If $E$ is $\d$-algebraically independent over~$K$,
then $E$ is contained in a $\d$-transcendence basis of $L$ over~$K$. If 
$L=\operatorname{cl} E$,
then $E$ contains a $\d$-transcendence basis of $L$ over~$K$.  
All $\d$-transcendence bases of $L$ over~$K$ have the same cardinality, which
is called the {\bf differential transcendence degree} of $L$
over~$K$, denoted by $\operatorname{trdeg}_\der L|K$.\index{differential transcendence!degree}\index{degree!differential transcendence}\nomenclature[M]{$\operatorname{trdeg}_\der(L\lvert{}K)$}{differential transcendence degree of the differential field extension $L\supseteq K$} So $\operatorname{trdeg}_\der L|K = 0$ iff $L$ is $\d$-algebraic over~$K$.
Hence $\operatorname{trdeg}_\der L|K = 0$ if the derivation of $L$ is trivial. 
If $M$ is a differential field extension of~$L$, then
$$\operatorname{trdeg}_\der M|K = \operatorname{trdeg}_\der M|L + \operatorname{trdeg}_\der L|K.$$

\subsection*{The zeros of a linear differential polynomial}
Let $K$ be a differential field and suppose the  differential polynomial
$A(Y)\in K\{Y\}$ is homogeneous of degree~$1$ and order $r\in \N$.
Then the map $y\mapsto A(y)\colon K\to K$ is $C$-linear, so the set of zeros
$$Z(A)\ :=\ \big\{ y\in K:\  A(y) = 0 \big\}$$
of $A$ in $K$ is a $C$-linear subspace of $K$.  Towards proving 
the well-known fact
that the dimension of this $C$-linear space is at most~$r$, we define the
{\bf Wronskian matrix} of a $(1+n)$-tuple
$(y_0,\dots,y_n)$ of elements of a differential ring by
$$\Wr(y_0,\dots,y_n)\  :=\ \begin{pmatrix} y_0 & y_1 & \cdots & y_n \\
y_0' & y_1' & \cdots & y_n' \\
\vdots & \vdots & \ddots & \vdots \\
y_0^{(n)} & y_1^{(n)} & \cdots & y_n^{(n)}\end{pmatrix},$$ 
with determinant
$$\wr(y_0,\dots,y_n)\ :=\ \det\Wr(y_0,\dots,y_n),$$
the {\bf Wronskian} of $(y_0,\dots,y_n)$.
We first show:

\index{Wronskian}
\nomenclature[M]{$\Wr(y_0,\dots,y_n)$}{Wronskian matrix of $y_0,\dots,y_n$}
\nomenclature[M]{$\wr(y_0,\dots,y_n)$}{Wronskian determinant $\det\Wr(y_0,\dots,y_n)$}

\begin{lemma}\label{lem:wronskian}
Let $y_0,\dots,y_n\in K$. Then
$$\wr(y_0,\dots,y_n)=0\quad\Longleftrightarrow\quad \text{$y_0,\dots,y_n$ are
 $C$-linearly dependent.}$$
\end{lemma}
\begin{proof}
Suppose $c_0,\dots,c_n\in C$ are not all zero and $\sum_{i=0}^n c_iy_i=0$. Taking successive derivatives yields $\sum_{i=0}^n c_i y_i^{(j)} = 0$ for each $j$, showing that the columns of $\Wr(y_0,\dots,y_n)$ are linearly dependent (over~$C$), so $\wr(y_0,\dots,y_n)=0$. This yields~``$\Leftarrow$,'' and we prove ``$\Rightarrow$'' by induction on $n$. The case $n=0$ being trivial, suppose $\wr(y_0,\dots,y_n)=0$ and $n>0$. Thus there are $a_0,\dots,a_n\in K$, not all zero, such that $\sum_{i=0}^n a_i y_i^{(j)}=0$ for $j=0,\dots,n$. After reindexing and normalization, we may assume that $a_0=1$, so
\begin{equation}\label{eq:wronskian}
y_0^{(j)} + \sum_{i=1}^n a_i y_i^{(j)}\ =\ 0\qquad\text{for $j=0,\dots,n$.}
\end{equation}
Taking derivatives in \eqref{eq:wronskian} for $j=0,\dots,n-1$ yields
$$ y_0^{(j+1)}+\sum_{i=1}^n a_i y_i^{(j+1)}+\sum_{i=1}^n a_i' y_i^{(j)}\ =\ 0\ 
=\ \sum_{i=1}^n a_i' y_i^{(j)} \qquad\text{for $j=0,\dots,n-1$.}$$
Hence, if $a_i'\neq 0$ for some  $i\in\{1,\dots,n\}$, then $\wr(y_1,\dots, y_n)=0$, so $y_1,\dots,y_n$ are $C$-linearly dependent by inductive hypothesis. 
If $a_i'=0$ for all $i\in\{1,\dots,n\}$, then $y_0,\dots,y_n$ are $C$-linearly dependent by \eqref{eq:wronskian} for $j=0$.
\end{proof}

\noindent
Thus if $y_0,\dots,y_n\in K$ are linearly independent over~$C$ and $L$ is a differential field extension of $K$, then $y_0,\dots, y_n$ remain linearly independent over~$C_L$. 

\begin{cor}\label{cor:wronskian}
$\dim_C Z(A) \leq r$.
\end{cor}
\begin{proof}
We may assume that $A(Y)=Y^{(r)}+\sum_{i=0}^{r-1} a_i Y^{(i)}$ where $a_0,\dots,a_{r-1}\in K$.
Let $y_0,\dots,y_r\in Z(A)$. Then the last row of
$$\Wr(y_0,\dots,y_r) = \begin{pmatrix} y_0 & y_1 & \cdots & y_r \\
y_0' & y_1' & \cdots & y_r' \\
\vdots & \vdots & \ddots & \vdots \\
-\sum a_i y_0^{(i)} & -\sum a_i y_1^{(i)} & \cdots & -\sum a_i y_r^{(i)}\end{pmatrix}$$
is a $K$-linear combination of the preceding rows, so $\wr(y_0,\dots, y_r)=0$. 
Now Lem\-ma~\ref{lem:wronskian} tells us that $y_0,\dots,y_r$ are $C$-linearly 
dependent.
\end{proof}

\begin{cor}\label{cor:wronskian, 1} 
Let $y_0,\dots,y_r\in K$ be linearly independent over~$C$ and let $z_0,\dots,z_r\in K$. Then there is
a unique homogeneous $A\in K\{Y\}$ of degree~$1$ and order at most~$r$ such that $A(y_i)=z_i$ for $i=0,\dots,r$.
\end{cor}
\begin{proof}
By Lemma~\ref{lem:wronskian}, the matrix $W:=\Wr(y_0,\dots,y_r)$ is invertible.
Take $a_0,\dots, a_r\in K$ such that $(a_0,\dots,a_r)W=(z_0,\dots,z_r)$, and set
$$A\ :=\ a_rY^{(r)}+a_{r-1}Y^{(r-1)}+\cdots+a_0Y.$$
Then $A(y_i)=z_i$ for $i=0,\dots,r$. If $B\in K\{Y\}$ is also homogeneous of degree $1$ and order at most $r$
with $B(y_i)=z_i$ for $i=0,\dots,r$, then $\dim_C Z(A-B)> r$ and hence $A=B$ by Corollary~\ref{cor:wronskian}.
\end{proof}

\noindent
We record a few more useful and well-known facts about Wronskians. Let 
$$A=Y^{(r)}+a_{r-1}Y^{(r-1)}+\cdots+a_0 Y\qquad\text{  where $a_0,\dots,a_{r-1}\in K$, $r\geq 1$.}$$
%First, an immediate consequence of Lemma~\ref{lem:wronskian} 
%and its Corollary~\ref{cor:wronskian}:

\begin{cor}\label{cor:wronskian, 2} Let $y_1,\dots,y_r$ be a basis of $Z(A)$ $($so $\dim_C Z(A)=r)$. Then
$$A(Y)\ =\ \frac{\wr(y_1,\dots,y_r,Y)}{\wr(y_1,\dots,y_r)}\quad \text{   in $K\<Y\>$.}$$
\end{cor}
\begin{proof} The right-hand quotient equals $B(Y)=Y^{(r)}+b_{r-1}Y^{(r-1)}+\cdots+b_0Y$ with $b_0,\dots, b_{r-1}\in K$, so $A-B$ is homogeneous of degree $1$
and of order $< r$. Since $y_i\in Z(A-B)$ for $i=1,\dots,r$, this gives
$A=B$ by Corollary~\ref{cor:wronskian}.
\end{proof}

\noindent
The following fact is known as Abel's identity:

\begin{lemma}\label{lem:abel} 
Let $y_1,\dots,y_r\in Z(A)$ and $w:=\wr(y_1,\dots,y_r)$. Then
$$w'\  =\ -a_{r-1}w.$$
\end{lemma}
\begin{proof}
Expressing the determinant $w$ in the usual way as a sum of $r!$ products, and
differentiating gives
$w'=w_1+\cdots+w_r$ where $w_i$ is the determinant of the matrix~$W_i$ obtained from $\Wr(y_1,\dots,y_r)$ by differentiating its $i$th row. For  $i=1,\dots,r-1$, each matrix~$W_i$  contains two identical rows, so $w_i=0$. Thus 
$$w'\ =\ w_r\ =\ \det W_r\ =\ \det \begin{pmatrix} 
y_1			& \cdots		& y_r \\
y_1'			& \cdots		& y_r'\\
\vdots		& \ddots		& \vdots \\
y_1^{(r-2)}	& \cdots		& y_r^{(r-2)} \\
y_1^{(r)}	& \cdots		& y_r^{(r)}\end{pmatrix}.$$
Using row operations and $A(y_i)=0$ for $i=1,\dots,r$, this equals
$$\det \begin{pmatrix} 
y_1					& \cdots		& y_r \\
y_1'					& \cdots		& y_r'\\
\vdots				& \ddots		& \vdots \\
y_1^{(r-2)}			& \cdots		& y_r^{(r-2)} \\
-a_{r-1}y_1^{(r-1)}	& \cdots		& -a_{r-1}y_r^{(r-1)}\end{pmatrix}\ =\ -a_{r-1}w,$$
as required.
\end{proof}

\begin{lemma}\label{lem:Wr under transforms}
Let $M$ be an $n\times n$ matrix over~$C$, $n\ge 1$,  and let 
$y=(y_1,\dots,y_n)$ and $z=(z_1,\dots,z_n)\in K^n$ be row vectors with $z=yM$. Then
$$\Wr(z_1,\dots,z_n)=\Wr(y_1,\dots,y_n)\cdot M, \quad \wr(z_1,\dots,z_n)=\wr(y_1,\dots,y_n)\cdot\det M.$$
\end{lemma}
\begin{proof}
We have $(z_1^{(i)},\dots,z_n^{(i)})=(y_1^{(i)},\dots,y_n^{(i)})\cdot M$ for each $i\in\N$.
\end{proof}

%\subsection*{Algebraic independence of iterated logarithms}
%In this subsection $K$ is a differential field, and 
%$Y$ a differential indeterminate over~$K$.
%An {\bf iterated logarithm sequence} in $K$ is a sequence 
%$(\ell_n)_{n\in\N}$ of nonzero elements of $K$ such that
%$\ell_{n+1}'=\ell_n^\dagger$ for all $n$. Such an iterated logarithm 
%sequence is {\bf standard} if $\ell_0'=1$. If $K$ is closed under 
%integration (for every $f\in K$ there is $y\in K$ with $y'=f$), 
%then for every $f\in K^\times$ there is an iterated logarithm 
%sequence $(\ell_n)$ with $\ell_0=f$.
%We have the following useful fact (see Theorems~6.1, 6.2, and 
%Lemma~1.3 in \cite{Babakhanian}):

%\begin{lemma}\label{criterion for P=0}
%Let $(\ell_n)_{n\in\N}$ be a standard iterated logarithm sequence in $K$.
%If $P\in C\{Y\}$ has order at most $r$ and $P(\ell_r)=0$, then $P=0$.
%\end{lemma}

\ifbool{PUP}{\enlargethispage*{0.9\baselineskip}}{}

\subsection*{Notes and comments}
Abstract differential fields were introduced by Baer~\cite{Baer27}. 
In the older literature one also finds \textit{hypertranscendental}~\cite{Maillet} and \textit{transcendentally transcendental}~\cite{Moore} instead of the terminology 
\textit{differentially transcendental} introduced by Kolchin~\cite{Kolchin44}.
The notion of $\d$-algebraic independence and the basic facts about it are due to Raudenbush~\cite{Raudenbush33}.
Many of the foundational results in differential algebra are due to Ritt and Kolchin. A comprehensive reference on this subject is Kolchin's book~\cite{Kolchin-DA}.  In some places, like~\cite{Kaplansky76}, a {\em differential ring\/} is a commutative ring~$K$ equipped with a derivation on~$K$, without requiring that $K$ contains $\Q$ as a subring.
(What we call a differential ring is called a {\em Ritt algebra\/} in~\cite{Kaplansky76}.)
A history of differential algebra can be found in \cite{Buium-Cassidy}.
Lemma~\ref{lem:abel} dates back to~\cite{Abel}.

%% file: mt-4-2.tex
\section{Decompositions of Differential Polynomials}
\label{Decompositions of Differential Polynomials}

\noindent
{\em In this section 
$K$ is a differential ring and $Y$ 
a differential indeterminate}. We let~$\der$
denote the derivation of $K$ as well as of $K\{Y\}$. Let $P=P(Y)\in K\{Y\}$ 
be a differential polynomial of order at most $r\in \N$.
We indicate some useful ways of expressing $P$ as a sum of 
differential polynomials of a special form.

\nomenclature[O]{$P^{\i}$}{$P^{i_0}(P')^{i_1}\cdots (P^{(r)})^{i_r}$, for  $\i =(i_0,\dots,i_r)\in \N^{1+r}$}
\index{differential polynomial!decomposition!natural}
\index{decomposition of a differential polynomial!natural}
\nomenclature[O]{$P_{\i}$}{coefficient of $Y^{\i}$ in the natural decomposition of  $P$}

\subsection*{Natural decomposition} 
For $\i =(i_0,\dots,i_r)\in \N^{1+r}$ and $Q\in K\{Y\}$, put 
$$Q^{\i} := Q^{i_0}(Q')^{i_1}\cdots (Q^{(r)})^{i_r}.$$
In particular, $Y^{\i} = Y^{i_0}(Y')^{i_1}\cdots (Y^{(r)})^{i_r}$, and
$y^{\i} = y^{i_0}(y')^{i_1}\cdots (y^{(r)})^{i_r}$ for $y\in K$.
Let $P_{\i}\in K$ be the coefficient of $Y^{\i}$ in $P$; then
$$P(Y)=\sum_{\i} P_{\i}Y^{\i} \quad\text{\bf (natural decomposition)}.$$

\subsection*{Decomposition into homogeneous parts} 
For each $\i$ as above we put $|\i| :=i_0 + \cdots + i_r$, and for 
$d\in \N$ we let 
$$P_d(Y):= \sum_{|\i| =d}P_{\i}Y^{\i},$$
the homogeneous part of degree $d$ of $P$. Note that then
$$ P(Y) = \sum_d P_d(Y) \quad 
\text{\bf (decomposition into homogeneous parts)}.$$
For $d=0$ and $d=1$ this gives $P_0=P(0)\in K$, and
$$ P_1\ =\ \sum_{n=0}^r \frac{\partial P}{\partial Y^{(n)}}(0)Y^{(n)}, \quad P_1\in \sum_{n=0}^r K\cdot Y^{(n)},$$
where $\frac{\partial P}{\partial Y^{(n)}}$ is the formal partial derivative of $P$ with respect 
to the variable $Y^{(n)}$, which has nothing to do with the derivation
$\der$ of $K\{Y\}$. Given also $Q\in K\{Y\}$ we have
$(PQ)_1=P(0)Q_1 + Q(0)P_1$.
For nonzero $P$ we have $$\deg(P)=\max\{d:P_d\neq 0\},$$ and we define 
the {\bf multiplicity of $P$ at $0$} as
$$\val(P):=\min\{d:P_d\neq 0\}.$$ We also set $\val(0):=\infty>\N$. 
Then for $Q\in K\{Y\}$, with the usual addition in $\N\cup\{\infty\}$:
$$\val(PQ)\ \ge\ \val(P)+\val(Q),\qquad \val(P+Q)\ \geq\ \min\!\big(\!\val(P),\val(Q)\big),$$
and if $K$ is an integral domain, then $\val(PQ)= \val(P)+\val(Q)$. We say that $P$ is {\bf homogeneous of degree $d$} if $P=P_d$. 

\nomenclature[Cb]{$\lvert\i\rvert$}{degree $i_0+\cdots+i_r$ of $\i=(i_0,\dots,i_r)\in\N^{1+r}$}
\index{differential polynomial!decomposition!into homogeneous parts}
\index{decomposition of a differential polynomial!into homogeneous parts}
\nomenclature[O]{$P_d$}{homogeneous part of degree $d$ of $P$}
\index{multiplicity!differential polynomial}
\index{differential polynomial!homogeneous}
\index{differential polynomial!multiplicity}
\index{homogeneous!differential polynomial}

\nomenclature[O]{$\val(P)$}{multiplicity of $P$ at $0$}
\nomenclature[O]{$\lvert\i\rvert$}{degree  $\lvert\i\rvert =i_0 + \cdots + i_r$ of $\i =(i_0,\dots,i_r)\in \N^{1+r}$}

\subsection*{Decomposition along orders} Let
$S^*$ be the set of words on a set $S$. \nomenclature[Ca]{$S^*$}{set of (finite) words on a set $S$} For any word ${\bsigma} =\sigma_1\cdots\sigma_n\in \{0,\dots,r\}^*$ 
of length $|\bsigma|= n$ \nomenclature[Cc]{$\lvert\bsigma\rvert$}{length $n$ of  $\bsigma=\sigma_1\cdots\sigma_n\in\N^*$} and $Q\in K\{Y\}$ we put 
$$Q^{[\bsigma]}\ :=\ Q^{(\sigma_1)}\cdots Q^{(\sigma_n)},$$
so for each permutation $\pi$ of $\{1,\dots,n\}$ we have 
$Q^{[\bsigma]} =Q^{[\pi(\bsigma)]}$, where
$\pi(\bsigma) := \sigma_{\pi(1)}\cdots\sigma_{\pi(n)}$. 
Thus 
$$P(Y)\ =\ \sum_{\bsigma} P_{[\bsigma]}Y^{[\bsigma]}
\quad\text{\bf (decomposition along orders)}$$
where the coefficients $P_{[\bsigma]} \in K$ are uniquely determined
by the requirements that $P_{[\bsigma]}=P_{[\pi(\bsigma)]}$
for each permutation $\pi$ of $\{1,\dots,\abs{\bsigma}\}$ (and of course
$P_{[\bsigma]}=0$ for all but finitely many $\bsigma$). For example, with $P=YY'$ and $r=1$, the words $\bsigma$ on $\{0,1\}^*$
with $|\bsigma|=\deg P=2$ are $00$, $01$, $10$, and $11$, with
$Y^{[00]}=Y^2$, $Y^{[01]}=Y^{[10]}=YY'$, and $Y^{[11]}=(Y')^2$, so
$P_{[01]}=P_{[10]}=1/2$, and $P_{[00]}=P_{[11]}=0$.

To find expressions relating the $P_{\i}$'s and
$P_{[\bsigma]}$'s, we first note that $Y^{[\bsigma]}=Y^{\i(\bsigma)}$, where
$\i(\bsigma)=(i_0,\dots,i_r)$ with 
$$i_k\ =\ \bigl|\bigl\{j\in\{1,\dots,|\bsigma|\}:\  \sigma_j=k\bigr\}\bigr|.$$ 
Thus $|\i(\bsigma)|=|\bsigma|$.
Using this notation we have
$$ P_{\i}\ =\ \sum_{\i(\bsigma)=\i} P_{[\bsigma]} .$$
Putting $\i!:= i_0!\cdots i_r!$ and $\binom{n}{\i} := \frac{n!}{\i!}$ we have
$$P_{[\bsigma]}\ =\ \frac{P_{\i(\bsigma)}}{\bigl|\bigl\{\btau:\ 
\abs{\btau}= \abs{\bsigma},\ \i(\btau)=\i(\bsigma)\bigr\}\bigr|}\
=\  {\binom{\abs{\bsigma}}{\i(\bsigma)}}^{-1}P_{\i(\bsigma)}.$$
%We also put $\abs{\bomega}:=\omega_1+\cdots+\omega_n$ for
%$\bomega=\omega_1\cdots\omega_n\in\{0,\dots,r\}^*$.

\index{differential polynomial!decomposition!along orders}
\index{decomposition of a differential polynomial!along orders}
\nomenclature[O]{$P^{[\bsigma]}$}{$P^{(\sigma_1)}\cdots P^{(\sigma_n)}$, for ${\bsigma} =\sigma_1\cdots\sigma_n\in \N^*$}
\nomenclature[O]{$P_{[\bsigma]}$}{coefficient of $Y^{[\bsigma]}$ in the decomposition along orders of~$P$}
\nomenclature[Cb]{$\i"!$}{$i_0"!\cdots i_r"!$ for $\i=(i_0,\dots,i_r)\in\N^{1+r}$}

\subsection*{Taylor expansion}
 For $\i=(i_0,\dots,i_r)$ and 
$\j=(j_0,\dots,j_r)$ in $\N^{1+r}$ we define \nomenclature[Cb]{$\i\leq\j$}{partial ordering on $\N^{1+r}$}
$$\i\leq\j \quad :\Longleftrightarrow\quad
i_0 \le j_0,\dots,i_r\le j_r,$$
and if $\i\leq\j$ we put $\binom{\j}{\i}:=
\binom{j_0}{i_0}\cdots\binom{j_r}{i_r}$. \nomenclature[Cb]{$\binom{\j}{\i}$}{$\binom{j_0}{i_0}\cdots\binom{j_r}{i_r}$, for $\i=(i_0,\dots,i_r)\leq\j=(j_0,\dots,j_r)$ in $\N^{1+r}$} It is also 
convenient to define, for $\i\in\N^{1+r}$: \nomenclature[O]{$P_{(\i)}$}{$\frac{P^{(\i)}}{\i"!}$, for  $\i\in\N^{1+r}$}\nomenclature[O]{$P^{(\i)}$}{$\frac{\partial^{\lvert\i\rvert}P}{\partial^{i_0}Y\cdots \partial^{i_r}Y^{(r)}}$, for $\i=(i_0,\dots,i_r)\in\N^{1+r}$}
$$P_{(\i)}\ :=\ \frac{P^{(\i)}}{\i !} \qquad\text{where
$P^{(\i)}\ :=\ 
\frac{\partial^{|\i|}P}{\partial^{i_0}Y\cdots \partial^{i_r}Y^{(r)}}$.}$$
(Here the right-hand side is an iterated partial derivative of $P$ with respect 
to the variables $Y, \dots, Y^{(r)}$ and has nothing to do with the derivation
$\der$ of $K\{Y\}$.) Thus, if $\abs{\i}=0$, then $P_{(\i)}=P$, 
and if $\abs{\i}=1$,
then $P_{(\i)}$ is one of the $\frac{\partial P}{\partial Y^{(k)}}$. We have
$$ \deg P_{(\i)}\ \le\ \deg P - |\i|, \qquad   
\text{so $\ P_{(\i)}=0\ $ if $\ |\i|>\deg P$.}$$ 
Repeated application of Lemma~\ref{lem:iterated Taylor} gives $(P_{(\i)})_{(\j)} = \binom{\i+\j}{\j}P_{(\i+\j)}$ for
$\i,\j\in\N^{1+r}$, where $\N^{1+r}$ is construed as a monoid under 
pointwise addition. 
%(construing $\N^{1+r}$ as a monoid under pointwise addition): 
%$$(P_{(\i)})_{(\j)} = \binom{\i+\j}{\j}P_{(\i+\j)}\qquad\text{for all
%$\i,\j\in\N^{1+r}$.}$$
Let $Z$ be a new differential indeterminate. Then 
$$P(Y+Z)\ =\ \sum_{\i} P_{(\i)}(Y)Z^{\i}
\quad\text{\bf (Taylor expansion)}.$$
Here is an application:

\index{differential polynomial!Taylor expansion}

\begin{lemma}\label{notthin} Suppose $K$ is a differential field and 
$y'\ne 0$ for some $y\in K$. If $P\ne 0$, then $P(y)\ne 0$ for
some $y\in K$.
\end{lemma}
\begin{proof} If $y\in K$ is algebraic over $C$, then $y\in C$. Thus $K$ has 
infinite dimension as a vector space over $C$. It follows that
if $a_0,\dots, a_r\in K$ are not all zero, then 
$a_0b+ \cdots + a_r b^{(r)}\ne 0$ for some $b\in K$, by
Corollary~\ref{cor:wronskian}.
Assume $P\notin K$ 
has order~$r$ and (total) degree $d\geq 1$. 
%$P\ne 0$, and set $d(P):= \sum_i \deg_{Y^{(i)}} P$. 
%We can assume $P\notin K$ and 
Then we have for all $b\in K$,
$$ P(Y+b)-P(Y)\ =\ \sum_{k=0}^r \frac{\partial P}{\partial Y^{(k)}}\cdot b^{(k)} + \sum_{|\i|\ge 2} P_{(\i)}(Y)\cdot b^{\i}\ =:\  Q_b(Y),$$
with $\deg(Q_b)< d$. Take $\i=(i_0,\dots, i_r)\in \N^{1+r}$ such that
$|\i|=d$ and $Y^{\i}$ occurs in~$P$. Take $k\in \{0,\dots,r\}$ 
with $i_k\ne 0$. Then $\frac{\partial P}{\partial Y^{(k)}}$ contains the 
monomial $Y^{\j}$ with $\j=(i_0,\dots, i_{k-1}, i_k-1, i_{k+1},\dots, i_r)$
%$$(Y^{(k)})^{(i_k-1)}(Y^{(k+1)})^{i_{k+1}}\cdots       
%  (Y^{(r)})^{i_r}$$ 
of degree $d-1$. This monomial does not occur in 
$\sum_{|\i|\ge 2} P_{(\i)}(Y)\cdot b^{\i}$, so the coefficient of this monomial in $Q_b(Y)$ is
$a_0b+ \cdots + a_r b^{(r)}$
where $a_0,\dots, a_r\in K$ are
independent of $b$, with $a_k\ne 0$. Now take $b\in K$ such that 
$a_0b+ \cdots + a_r b^{(r)}\ne 0$, so $Q_b\ne 0$. 
We can assume inductively that $Q_b(y) \ne 0$ for some $y\in K$. Then $P(y+b)\ne 0$ or $P(y)\ne 0$ for such $y$, by the identity above. 
\end{proof}

\index{thin}
\index{differential field!thin subset}
\index{subset!thin} 

\noindent
Suppose $K$ is a differential field with nontrivial derivation, that is, $y'\ne 0$ for some ${y\in K}$. Call a set $S\subseteq K$ {\bf thin} (in $K$) if there is an $F\in K\{Y\}^{\ne}$ such that $F(y)=0$ for all $y\in S$. Note that a finite union of thin subsets of
$K$ is thin, and that $K$ is not thin by Lemma~\ref{notthin}. 
Given any $F\in K\{Y\}^{\ne}$ the set
$$ \big\{y\in K^\times:\ F(y^\dagger)=0\big\}$$
is thin, 
since in $K\langle Y\rangle $ we have $0\ne Y^nF(Y^\dagger)\in K\{Y\}$ 
for some $n$. 

\begin{lemma}\label{thindesc} Let $L|K$ be a differential field extension
of the differential field $K$. Assume $y'\ne 0$ for some $y\in K$ and 
$S\subseteq L$ is thin in $L$. Then $S\cap K$ is thin in $K$.
\end{lemma}
\begin{proof} Take $P\in L\{Y\}^{\ne}$ such that 
$S\subseteq \big\{y\in L:\ P(y)=0\big\}$. We have
$$P\ =\ b_1P_1+ \cdots + b_nP_n$$ where
$b_1,\dots, b_n\in L$ are linearly independent over $K$, 
$P_1,\dots, P_n\in K\{Y\}$ are nonzero, and~$n\ge 1$. Then 
$S\cap K\subseteq \big\{y\in K:\ P_1(y)=0\big\}$.
\end{proof}

\subsection*{Decomposition into isobaric parts}
In this subsection
$P,Q\in K\{Y\}$ are of order at most $r\in \N$. Also,
$\i=(i_0,\dots, i_r)\in \N^{1+r}$, $\|\i\|:= i_1+2i_2 + \cdots + ri_r$, we let~$\bsigma$ 
range over words in $\{0,\dots,r\}^*$, and  
for $\bsigma=\sigma_1\cdots\sigma_d$ we put
$\|\bsigma\|:=\sigma_1+\cdots+\sigma_d$,  
so $\|\bsigma\|=\|\i(\bsigma)\|$.
Define the {\bf weight} $\wt(P)\in \N\cup \{-\infty\}$ of $P$ by
$$\wt(P)\ :=\ \max \big\{ \|\bsigma\|\ : P_{[\bsigma]}\neq 0 \big\} \text{ if $P\neq 0$,} \qquad \wt(0):= -\infty < \N. $$
In particular, 
$$\wt(Y^{\i})\ =\ \|\i\|, \quad \wt(P_{(\i)})\ \le\ \wt(P) - \|\i\|,$$ 
and
$$\wt(P+Q)\ \leq\ \max\!\big(\!\wt(P), \wt(Q)\big).$$
Define the {\bf weighted multiplicity} $\wv(P)\in \N\cup\{\infty\}$ of $P$ at $0$ by
$$\wv(P)\ :=\ \min \big\{ \|\bsigma\| :\  P_{[\bsigma]}\neq 0 \big\} \text{ if $P\neq 0$,} \qquad \wv(0):= \infty > \N.$$
In particular, $\wv(Y^{\i})=\wt(Y^{\i})=\|\i\|$, and 
$  \wv(P+Q)\geq \min\!\big(\!\wv(P), \wv(Q)\big)$.
We say that 
$P$ is {\bf isobaric} if $P_{[\bsigma]}=0$ for all $\bsigma$ with
$\|\bsigma\|\ne \wt(P)$. Note that $0$ is isobaric, and that if $P\ne 0$, then
$P$ is isobaric iff $\wv(P)=\wt(P)$. (Although $\deg 0=\wt 0 =-\infty$, we do consider $0$ to be homogeneous
of any degree $d\in \N$ and to be isobaric of any weight $w\in \N$.) Degree and weight behave as follows under 
differentiation: $(Y^\i)'$ is homogeneous of the same degree $|\i|$ as 
$Y^\i$, and is isobaric of weight $\wt(Y^\i) +1=\|\i\|+1$.  Thus if  $P\in C\{Y\}$ is isobaric of weight $w\in \N$, then~$P'$ is isobaric of weight $w+1$. 
For $w\in\N$, let
$$P_{[w]} := \sum_{\|\bsigma\|=w} P_{[\bsigma]} Y^{[\bsigma]}$$
be the {\bf isobaric part of $P$ of weight $w$.} 
Then $P_{[w]}=0$ for $w>\wt(P)$, and
$$P=\sum_w P_{[w]} \qquad\text{\bf (decomposition into isobaric
parts)}.$$
Note that $(PQ)_{[w]}=\sum_{i+j=w} P_{[i]}Q_{[j]}$.
In particular, if $P$ and $Q$ are isobaric, then so is $PQ$. Moreover:

\index{differential polynomial!decomposition!into isobaric parts}
\index{decomposition of a differential polynomial!into isobaric parts}
\nomenclature[Cc]{$\lvert\lvert\bsigma\rvert\rvert$}{weight $\lvert\lvert\bsigma\rvert\rvert=\sigma_1+\cdots+\sigma_d$ of  $\bsigma=\sigma_1\cdots\sigma_d\in\N^*$}
\nomenclature[Cb]{$\lvert\lvert\i\rvert\rvert$}{weight $\lvert\lvert\i\rvert\rvert=i_1+2i_2 + \cdots + ri_r$ of $\i=(i_0,\dots, i_r)\in \N^{1+r}$}
\nomenclature[O]{$\wt(P)$}{weight of  $P$}
\nomenclature[O]{$\wv(P)$}{weighted multiplicity of  $P$}
\index{weight}
\index{weighted multiplicity}
\index{isobaric!differential polynomial}
\index{differential polynomial!isobaric}
\nomenclature[O]{$P_{[w]}$}{isobaric part of~$P$ of weight~$w$}

\begin{lemma}\label{w and W} If $K$ is an integral domain, then 
$$\wv(PQ)=\wv(P)+\wv(Q), \qquad \wt(PQ)=\wt(P)+\wt(Q).$$
\end{lemma}

\subsection*{Decomposition into subhomogeneous parts}
In Chapter~\ref{ch:The Dominant Part and the Newton Polynomial} we use yet another decomposition of $P$. For $\i=(i_0,\dots,i_r)\in\N^{1+r}$ define the 
{\bf subdegree of $\i$} \index{subdegree}\nomenclature[Cb]{$\lvert\i\rvert'$}{subdegree $\lvert\i\rvert'=i_1+\cdots+i_r$ of $\i=(i_0,\dots, i_r)\in \N^{1+r}$} as $|\i|':=i_1+\cdots+i_r$, and for $d\in\N$, set
$P_{|d|'} := \sum_{|\i|'=d} P_{\i}\,Y^{\i}$,  so\nomenclature[O]{$P_{\lvert d\rvert'}$}{subhomogeneous part of $P$ of subdegree~$d$}\index{differential polynomial!decomposition!into subhomogeneous parts}\index{decomposition of a differential polynomial!into subhomogeneous parts}
$$  P = \sum_d P_{|d|'} \quad\text{\bf (decomposition into subhomogeneous parts)}.
$$
The {\bf subdegree of $P$}, $\operatorname{sdeg}(P)$,
is the largest $d$ with $P_{|d|'}\neq 0$ if $P\ne 0$, and $\operatorname{sdeg}(P):=-\infty$ if $P=0$. \index{differential polynomial!subdegree}\nomenclature[O]{$\operatorname{sdeg}(P)$}{subdegree of  $P$}
 We say that $P$ is {\bf subhomogeneous} of subdegree $d$ if $P=P_{|d|'}$. \index{differential polynomial!subhomogeneous}
\index{subhomogeneous} Thus $P$ is subhomogeneous of subdegree $d$ iff there are homogeneous  
$P_j\in K\{Y'\}\subseteq K\{Y\}$
of degree $d$ such that $P=\sum_j P_j(Y)\, Y^j$. Suppose $P\ne 0$ is subhomogeneous of subdegree $d$; then
$\wt(P)\geq d$, with equality iff $P\in K[Y]\cdot (Y')^d$.

\subsection*{Logarithmic decomposition} In this subsection we consider 
iterated logarithmic derivatives. Let $K$ be a 
differential field and $y\in K$. We set $y^{\<0\>}:= y$, and
inductively, if $y^{\<n\>}\in K$ is defined and nonzero, 
$y^{\<n+1\>}:=  (y^{\<n\>})^{\dagger}$ (and otherwise $y^{\<n+1\>}$ is not defined). The results we state below follow by easy inductions on 
$n$ using Lemma~\ref{derpol}. First, if $y^{\<n\>}$ 
is defined, then
$$  y^{(n)}\ =\ y\cdot L_n(y^{\<1\>},\dots,y^{\<n\>} )\ =\ 
y^{\<0\>}\cdot L_n(y^{\<1\>},\dots,y^{\<n\>} )    $$
where the polynomial $L_n\in \Z[X_1,\dots, X_n]$ depends only on $n$ (not on $y$
or $K$) and is homogeneous of degree $n$ with nonnegative coefficients:
\begin{align*}  L_0\ &=\ 1\\
               L_1\ &=\ X_1\\
              L_2\ &=\ X_1^2 + X_1X_2\\
             L_3\ &=\ X_1^3 + 3X_1^2X_2 + X_1X_2^2 + X_1X_2X_3\\
                 &\ \ \vdots
\end{align*}
The $L_n$ are given by the recursion
$$ L_0= 1, \qquad L_{n+1}\ =\ X_1L_n + \sum_{j=1}^n X_jX_{j+1}\frac{\partial L_n}{\partial X_j}.$$
So if  $y^{\<n\>}$ is defined, then 
$\Q\big[y,y',\dots, y^{(n)}\big]\subseteq \Q\big[y^{\<0\>}, y^{\<1\>},\dots,y^{\<n\>}\big]$. Similarly: \begin{enumerate}
\item[(1)] if $y^{\<n\>}$ is defined, then $\Q\big(y,y',\dots, y^{(n)}\big)\ =\ \Q\big(y^{\<0\>}, y^{\<1\>},\dots,y^{\<n\>}\big)$;
\item[(2)] if $y,y', \dots,y^{(n)}$ are algebraically independent over $\Q$, 
then $y^{\<n\>}$ is defined, and so $y^{\<0\>}, y^{\<1\>},\dots,y^{\<n\>}$ are algebraically independent over $\Q$ by (1).
\end{enumerate} 
Thus in $K\langle Y \rangle$ each $Y^{\<n\>}$ is
defined,  
$Y^{\<0\>}, Y^{\<1\>},\dots,Y^{\<n\>}$
are algebraically independent over $K$, and 
$K\langle Y \rangle=K\big(Y^{\<n\>}: n=0,1,2,\dots\big)$. 
If $y^{\<n\>}$ is defined and $\i=(i_0,\dots, i_n)\in \mathbb{N}^{1+n}$, 
we set $$  y^{\<\i\>}\ :=\ 
(y^{\<0\>})^{i_0}(y^{\<1\>})^{i_1}\cdots (y^{\<n\>})^{i_n}. $$
Using the $L_n$ it follows that $P$ has a unique decomposition
\[\begin{array}{lll} P\ =\ &\sum_{\i}P_{\<\i\>}Y^{\<\i\>}&  
\text{{\bf (logarithmic decomposition)},}
\end{array}\]
with $\i$ ranging over $\N^{1+r}$, all $P_{\<\i\>}\in K$, and
$P_{\<\i\>}=0$ for all but finitely many $\i$.

\index{differential polynomial!decomposition!logarithmic}
\index{decomposition of a differential polynomial!logarithmic}

\nomenclature[O]{$P_{\<\i\>}$}{coefficient of $Y^{\<\i\>}$ in the logarithmic decomposition of~$P$}
\nomenclature[M]{$a^{\<n\>}$}{$n$th iterated logarithmic derivative of $a$}

\subsection*{Notes and comments}
Lemma~\ref{notthin} appears in \cite[p.~35]{Ritt}.

%% file: mt-4-3.tex
\section{Operations on Differential Polynomials}\label{Operations on Differential Polynomials}

\noindent
{\em In this section $K$ is a differential ring with derivation $\der$ and 
$Y$ is a differential indeterminate.}\/ We discuss here additive and 
multiplicative conjugation of differential polynomials, Ritt reduction, 
composition, and substituting powers. The operation of
compositional conjugation is more involved and is treated in
Section~\ref{Compositional Conjugation}.

\subsection*{Additive and multiplicative conjugation}
Let $P=P(Y)\in K\{Y\}$ be of order at most $r$, and let $h, h_1, h_2\in K$.
We define the differential polynomials
\[
\begin{array}{lll}
P_{+h}\     =\ P_{+h}(Y)\     :=\ & P(h+Y)\in K\{Y\}   
&\text{{\bf (additive conjugate of $P$)},}\\
P_{\times h}\ =\ P_{\times h}(Y)\ :=\ & P(hY)\in K\{Y\}   
&\text{{\bf (multiplicative conjugate of $P$)}.} 
\end{array}\]
Thus
$$P_{+(h_1 + h_2)} =(P_{+h_1})_{+h_2},\quad
P_{\times h_1h_2} =(P_{\times h_1})_{\times h_2},\quad 
(P_{\times h_1})_{+h_2} = (P_{+h_1h_2})_{\times h_1}.$$
These conjugations commute with $\der$, that is,
$$(P_{+h})'\ =\ (P')_{+h}, \qquad 
(P_{\times h})'\ =\ (P')_{\times h}.$$
{\sloppy Additive conjugation with $h$ is an automorphism $Q\mapsto Q_{+h}$ of the 
differential ring~$K\{Y\}$, with inverse $Q\mapsto Q_{-h}:=Q_{+(-h)}$,
and if $h$ is a unit of $K$, then also multiplicative conjugation with $h$ is
an automorphism $Q\mapsto Q_{\times h}$ of this differential ring, with inverse 
$Q\mapsto Q_{/h}:=Q_{\times h^{-1}}$. These automorphisms are the identity on~$K$. 
If $P$ is homogeneous of degree $d$, then so is $P_{\times h}$. Thus (without restricting $P$): 
$$(P_{\times h})_d\ =\  (P_d)_{\times h}, \qquad \text{($d\in \N$).}$$
Note also that $(P_{+h})_1=\sum_{n=0}^r \frac{\partial P}{\partial Y^{(n)}}(h)Y^{(n)}$.
}
\index{conjugation!additive}
\index{conjugation!multiplicative}
\index{differential polynomial!conjugation!additive}
\index{differential polynomial!conjugation!multiplicative}

\nomenclature[O]{$P_{+h}$}{additive conjugate $P(Y+h)$ of  $P$ by $h$}
\nomenclature[O]{$P_{\times h}$}{multiplicative conjugate $P(hY)$ of  $P$ by $h$}

\begin{example}
Consider a homogeneous linear $A\in K\{Y\}$:
$$A(Y)\ =\ a_0Y + a_1Y' + \cdots + a_nY^{(n)} \qquad
(a_0,\dots,a_n\in K).$$
Then
$$A_{\times h}(Y)\ =\ b_0Y+ b_1Y'+\cdots+b_nY^{(n)}$$
where
$$b_i\ =\ \sum_{j=i}^n \binom{j}{i} a_j h^{(j-i)},$$
in particular $b_0=a_0h+a_1h'+\cdots+a_nh^{(n)}=A(h)$ and $b_n=a_nh$.
\end{example}

\nomenclature[Cc]{$\bsigma \le \btau$}{partial ordering on $\N^*$}
\nomenclature[Cc]{$\binom{\btau}{\bsigma}$}{$\binom{\tau_1}{\sigma_1}\cdots\binom{\tau_d}{\sigma_d}$,
for $\bsigma=\sigma_1\cdots\sigma_d\leq \btau=\tau_1\cdots\tau_d$ in $\N^*$}

\noindent
We now derive expressions for the coefficients of the additive and 
multiplicative conjugates of $P$. 
%For $\i=(i_0,\dots,i_r)$ and 
%$\j=(j_0,\dots,j_r)$ in $\N^{r+1}$ we define
%$$\i\leq\j \quad :\Longleftrightarrow\quad
%i_0 \le j_0,\dots,i_r\le j_r,$$
%and if $\i\leq\j$ we put $\binom{\j}{\i}:=
%\binom{j_0}{i_0}\cdots\binom{j_r}{i_r}$. 
Given words
$\bsigma = \sigma_1\cdots \sigma_m$ and
$\btau = \tau_1\cdots \tau_n$ in $\{0,\dots,r\}^*$, define
$$\bsigma \le \btau \quad :\Longleftrightarrow\quad  m=n \text{ and } 
\sigma_1 \le \tau_1,\dots,\sigma_n\le \tau_n,$$
and if $\bsigma \le \btau$, put 
$$\binom{\btau}{\bsigma} 
:=\binom{\tau_1}{\sigma_1}\cdots\binom{\tau_n}{\sigma_n}, \qquad
\btau-\bsigma:=(\tau_1-\sigma_1)\cdots (\tau_n-\sigma_n)\in
\{0,\dots,r\}^*.$$ 
With these notations we have:

\begin{lemma}\label{Conjugation-Lemma}
Let $\i=(i_0,\dots,i_r)\in\N^{1+r}$ and 
$\bsigma\in\{0,\dots,r\}^*$. Then
\begin{align}
(P_{+h})_{\i}\ & =\ \sum_{\j \ge \i}\binom{\j}{\i}P_{\j}h^{\j - \i},
\label{Additive-Conj}\\
(P_{\times h})_{[\bsigma]}\ &=\ \sum_{\btau \ge \bsigma}
\binom{\btau}{\bsigma}P_{[\btau]} h^{[\btau - \bsigma]}.
\label{Multiplicative-Conj}
\end{align}
\end{lemma}
\begin{proof}
By  Taylor expansion,
$$
P_{+h}(Y)\ =\ \sum_\j P_{(\j)}(h)Y^\j,$$
hence $(P_{+h})_{\i} = P_{(\i)}(h)$.
Now $P(Y)=\sum_{\j} P_{\j}Y^{\j}$ gives 
$$P_{(\i)}\ =\ \sum_{\j} P_{\j}\frac{(Y^{\j})^{(\i)}}{\i!}\ 
 =\  \sum_{\j \ge \i} \binom{\j}{\i} P_{\j}Y^{\j - \i} ,$$
since $(Y^{\j})^{(\i)}= \frac{\j!}{(\j - \i)!}Y^{\j-\i}$ if $\j \ge \i$
and $(Y^{\j})^{(\i)}=0$ otherwise,
as an easy induction shows. This yields \eqref{Additive-Conj}.
A simple induction on $n=\abs{\btau}$ shows that for
$\btau = \tau_1\cdots \tau_n\in \{0,\dots,r\}^*$ we have
$$(ab)^{[\btau]}\ =\ \sum_{\bsigma\le \btau} \binom{\btau}{\bsigma} 
a^{[\btau - \bsigma]}b^{[\bsigma]}\qquad\text{for $a,b\in K\{Y\}$.}$$
Hence
\begin{align*} P(hY)\ =\ \sum_{\btau} P_{[\btau]}(hY)^{[\btau]}\ &=\
\sum_{\btau} \left( P_{[\btau]} \sum_{\bsigma\le \btau} 
\binom{\btau}{\bsigma} h^{[\btau - \bsigma]}Y^{[\bsigma]}\right)\\
&=\ \sum_{\bsigma}\left(\sum_{\btau \ge \bsigma}
\binom{\btau}{\bsigma}P_{[\btau]} h^{[\btau - \bsigma]}\right)Y^{[\bsigma]},
\end{align*}
and this yields \eqref{Multiplicative-Conj}.
\end{proof}

\begin{cor} We have $\deg P_{+h} = \deg P$ and $\wt P_{+h} = \wt P.$ If
$h$ is a unit of~$K$, then $\deg P_{\times h} = \deg P$ and 
$\wt P_{\times h} = \wt P$. 
\end{cor} 

\noindent
The following identities have routine proofs, using the Chain Rule:
\begin{align}
\frac{\partial P_{+h}}{\partial Y^{(i)}}\ &=\ \left(\frac{\partial P}{\partial Y^{(i)}}\right)_{+h},\label{eq:partial ders and operations, 1} \\
\frac{\partial P_{\times h}}{\partial Y^{(i)}}\ &=\ \sum_{j\geq i} {j\choose i} h^{(j-i)} \left(\frac{\partial P}{\partial Y^{(j)}}\right)_{\times h}. \label{eq:partial ders and operations, 2}
\end{align}

\noindent
Suppose now that $P\neq 0$. We have  
$$0\leq \val P_{+h} \leq \deg P, \qquad
\val P_{+h}>0 \Longleftrightarrow P(h)=0.$$
We call $\val P_{+h}$ the
{\bf multiplicity} of $P$ at $h$. 
Note that if $P\in K[Y]$, then 
$$P(Y)\ =\ Q(Y)\cdot(Y-h)^\mu,\ \text{ where $Q(Y)\in K[Y]$, $Q(h)\ne 0$,  $\mu=\val P_{+h}$.}$$

\index{differential polynomial!multiplicity}
\index{multiplicity}

\begin{lemma}\label{lem:multiplicity}
Suppose $0\neq P\in K[Y](Y')^\N$ and $h'=0$, and let $\mu\in\N$. Then 
\begin{align*}\val P_{+h}=\mu \quad\Longleftrightarrow\quad &\text{there are $i,j\in \N$
and $Q\in K[Y]$ such that $i+j=\mu$,}\\
&\text{$P(Y) =  Q(Y)\cdot(Y-h)^i\cdot (Y')^j$, and $Q(h)\ne 0$}. 
\end{align*}
\end{lemma}
\begin{proof}
The implication $\Leftarrow$ is obvious. For
the direction $\Rightarrow$, note that we have
$P(Y)=R(Y)\cdot(Y')^j$ with $j\in\N$ and $R\in K[Y]$. Set $i:=\val R_{+h}$; then $\val P_{+h} = i + j$ and
$R(Y)=Q(Y)\cdot(Y-h)^i$ where $Q\in K[Y],\ Q(h)\neq 0$.
\end{proof}

\begin{cor}\label{cor:multiplicity}
For $0\neq P\in K[Y](Y')^\N$ and $h'=0$ we have:
$$\val P_{+h}=\deg P \ \Longleftrightarrow\ \text{$P =  a \cdot (Y-h)^i\cdot (Y')^j$ where  $a\in K^{\neq}$ and $i+j=\deg P$.}$$
\end{cor}

%\noindent
%We also note
%(construing $\N^{r+1}$ as a monoid under pointwise addition): 
%$$(P_{(\i)})_{(\j)} = \binom{\i+\j}{\j}P_{(\i+\j)}\qquad\text{for all
%$\i,\j\in\N^{r+1}$.}$$

%For $\btau \in \{0,\dots,r\}^*$ we have $(P_d)_{[\btau]}= P_{[\btau]}$ 
%if $l(\btau)=d$, and 
%$(P_d)_{[\btau]}=0$ otherwise. Using this fact it follows easily
%from \eqref{Multiplicative-Conj} that 
%$$(P_{\times h})_d = (P_d)_{\times h}.$$
%using the fact that
%$(P_d)_{[\btau]}= P_{[\btau]}$ if $l(\btau)=d$, and 
%$(P_d)_{[\btau]}=0$ otherwise.

\subsection*{Complexity and Ritt division} {\em Let $K$ be a differential field,
$P,Q\in K\{Y\}$, and~$a\in K$}. For $P\notin K$, denote the order of $P$ by $r_P$, the degree of $P$ in $Y^{(r_P)}$
by~$s_P$, and the total degree of $P$ by $t_P$ (so $s_P, t_P\ge 1$), and define the {\bf complexity} of~$P$ to be the triple $\c(P)=(r_P, s_P, t_P)\in \N^3$.\index{complexity!differential polynomial}
\index{differential polynomial!complexity}\nomenclature[O]{$\c(P)$}{complexity of $P$} For $P\in K$ we set $\c(P)=(0,0,0)$. We order $\N^3$ lexicographically.
Note that if $P,Q\notin K$, then $\c(P), \c(Q) < \c(PQ)$. 
For $P\notin K$ and $r=r_P$, $s=s_P$ we have 
$$P\ =\ F_0 + F_1\cdot Y^{(r)} + \cdots + F_s\cdot (Y^{(r)})^s, \qquad F_0,\dots, F_s\in K\big[Y,\dots, Y^{(r-1)}\big],$$ and we define the {\bf initial}
of $P$ to be $I_P:= F_s$, \index{initial}\index{differential polynomial!initial} and the {\bf separant} of $P$ to be\index{separant}\index{differential polynomial!separant}\nomenclature[O]{$I_P$}{initial of $P$}
\nomenclature[O]{$S_P$}{separant of $P$}
$$S_P\ :=\ \frac{\partial P}{\partial Y^{(r)}}\ =\ 
\sum_{i=1}^s  iF_i\cdot (Y^{(r)})^{i-1},$$
so $I_P\ne 0,\ S_P\ne 0$, and $\c(I_P) < \c(P),\ \c(S_P) < \c(P)$.
The quantities above transform as follows under various conjugations:
\begin{align*} &\c(P_{+a})=\c(P),\quad I_{P_{+a}}=I_{P,+a},\  \qquad \quad 
S_{P_{+a}}=S_{P,+a},\\
\text{ and for $a\ne 0$,}\quad &\c(P_{\times a})=\c(P),\quad I_{P_{\times a}}=a^{s_P}\cdot I_{P,\times a}, \quad S_{P_{\times a}}=a\cdot S_{P,\times a}.
\end{align*}

\begin{lemma}\label{lem:comp1}
Suppose $P\notin K$, and set $r=r_P$. Then for $n\geq 1$, 
$$P^{(n)}\ =\ G_n + S_PY^{(r+n)}\qquad\text{with $G_n\in K\big[Y,Y',\dots,Y^{(r+n-1)}\big]$}.$$
In particular,
$P^{(n)}\notin K$ and $P^{(n)}$ has order $r+n$, for every $n$. 
\end{lemma}
\begin{proof}
Let
$s=s_P$ be the degree of $P$ in $Y^{(r)}$, so $s\geq 1$, and
$$P\ =\ F_0+F_1\cdot Y^{(r)}+\cdots + F_s\cdot (Y^{(r)})^s,\quad  F_0,\dots,F_s\in K\big[Y,Y',\dots,Y^{(r-1)}\big].$$
Then $P'=G_1+ S_PY^{(r+1)}$ with
$$G_1:=F_0'+F_1'\cdot Y^{(r)}+\cdots+ F_s'\cdot (Y^{(r)})^s\in K\big[Y,Y',\dots,Y^{(r)}\big].$$
This gives the case $n=1$, and the rest goes by induction
on $n$.
\end{proof}

\begin{cor}\label{PdividesP'} Assume $P\notin K$. Then $P'\notin PK\{Y\}$.
\end{cor}
\begin{proof} With the notations above and $r=r_P$ we have $P'=S_PY^{(r+1)} + G_1$. Assume $P'=P Q$. Then 
$Q=AY^{(r+1)}+ B$ with $A,B\in K[Y,\dots, Y^{(r)}]$, so
$S_P=PA$, contradicting $\operatorname{c}(S_P) < \operatorname{c}(P)$.
\end{proof}

\begin{cor}\label{gaussdagger} Assume $F\in K\<Y\>$, $F\notin K$. Then $F^\dagger\notin K\{Y\}$.
\end{cor}
\begin{proof} We have $F=P/Q$ with relatively prime
$P,Q\in K\{Y\}^{\ne}$. Towards a contradiction, assume
$F^\dagger=R\in K\{Y\}$. Then $P'Q-PQ'=PQR$, so $P$ divides 
$P'Q$ in~$K\{Y\}$, hence $P$ divides $P'$ in~$K\{Y\}$, and so $P\in K$ by Corollary~\ref{PdividesP'}. Likewise we
get $Q\in K$, so
$F\in K$, a contradiction.
\end{proof}

\noindent
Next an analogue for differential polynomials of division with remainder: \index{differential polynomial!Ritt division}\index{Ritt division}\index{theorem!Ritt division}
%The proof is constructive.

\begin{theorem}[Ritt division]\label{thm:Ritt division}
Let $P\notin K$. Then
  $$ I_P^k S_P^l Q\ =\ A_0P+A_1P' + \dots + A_nP^{(n)} + R$$
for some $k,l\in \N$, $A_0,\dots, A_n\in K\{Y\}$, and 
$R\in K\{Y\}$ with $\c(R) < \c(P)$. 
\end{theorem}

\noindent
We break up the proof into two steps, and first show:

\begin{lemma}\label{lem:Ritt division}
Let $P\notin K$. There are $l\in \N$ and $A_0,\dots, A_n\in K\{Y\}$ such
that
  $$ S_P^l Q\ =\ A_0P+A_1P' + \dots + A_nP^{(n)} + R, \qquad \order(R) \leq \order(P).$$
\end{lemma}
\begin{proof}
Let $r:=r_P$ and $Q\in K[Y,\dots, Y^{(r+i)}]$. By induction on~$i$ we show: 
$$S_P^l Q\ \equiv\ R\mod (P,P',\dots,P^{(i)})\quad\text{in the ring $K\big[Y,Y',\dots,Y^{(r+i)}\big]$}$$
for suitable $l\in\N$ and
$R\in K\{Y\}$ with $\order(R) \leq r$. If $i=0$, then we can take 
$l=0$, $R=Q$. So suppose $i>0$. From Lemma~\ref{lem:comp1} we obtain
$$P^{(i)}\ =\ S_P Y^{(r+i)}+G_i,\quad G_i\in K\big[Y,Y',\dots,Y^{(r+i-1)}\big].$$
As $Q=\sum_{k=0}^s F_k\cdot (Y^{(r+i)})^k$ with all $F_k\in K\big[Y,Y',\dots,Y^{(r+i-1)}\big]$, we get \\
\begin{align*} S_P^s Q\ &=\ \sum_{k=0}^s F_k\cdot (P^{(i)}-G_i)^k S_P^{s-k}\ \equiv\ Q_i\mod(P^{(i)}) \quad \text{in }K\big[Y,\dots,Y^{(r+i)}\big]\\
 &\text{ with }Q_i\ =\ \sum_{k=0}^s F_k\cdot (-G_i)^k S_P^{s-k}\in K\big[Y,Y',\dots,Y^{(r+i-1)}\big].
\end{align*}
The inductive assumption gives $l\in\N$ and $R\in K\{Y\}$ with $\order(R)\leq r$ and 
$$S_P^l Q_i\ \equiv\ R\mod \big(P,P',\dots,P^{(i-1)}\big)\quad\text{in the ring $K\big[Y,Y',\dots,Y^{(r+i-1)}\big]$.}$$
In combination with the previous congruence, this gives
$$S_P^{l+s} Q\ \equiv\ R\mod \big(P,P',\dots,P^{(i)}\big)\quad\text{in the ring $K\big[Y,Y',\dots,Y^{(r+i)}\big]$,} $$
which finishes the induction.
\end{proof}

\begin{proof}[Proof of Theorem~\ref{thm:Ritt division}]
Take $l, A_0,\dots,A_n,R$ as in Lem\-ma~\ref{lem:Ritt division}. We are done if $\c(R)<\c(P)$, so assume $\c(R)\ge \c(P)$. Then $r:=r_P=r_R$
and $s_R\ge s_P$.
Ordinary division by $P$ for polynomials in $Y^{(r)}$ 
over $K\big[Y,Y',\dots,Y^{(r-1)}\big]$
yields
$$I_P^k R\ =\ A^*P+R^*,\quad k=1+s_R-s_P,\ A^*,R^*\in K\big[Y,Y',\dots,Y^{(r)}\big],\ s_{R^*}<s_P.$$
Hence
$$I_P^k S_P^l Q \ =\ A_0^*P+A_1^*P'+\cdots+A_n^*P^{(n)}+R^*$$
where $A_0^*:=I^k_P A_0+A^*$, $A_i^*:=I^k_P A_i$ for $i=1,\dots,n$, and $\c(R^*)<\c(P)$.
\end{proof}

\noindent
Additive conjugation in the Ritt division of the theorem above gives:
$$ I_{P_{+a}}^k S_{P_{+a}}^l Q_{+a}\ =\ (A_0)_{+a}P_{+a}+(A_1)_{+a}P_{+a}' + \dots + (A_n)_{+a}P_{+a}^{(n)} + R_{+a}.$$

\noindent
Corollary~\ref{cor:order lower bound} and Lemma~\ref{lem:Ritt division} yield:

\begin{cor}\label{cor:rittvariant}
Let $y$ be an element of a differential field extension of $K$,  $\d$-algebraic over~$K$, and let $P\in K\{Y\}^{\neq}$ be a minimal annihilator of $y$ over $K$. Then for $Q\in K\{Y\}$ we have: 
$Q(y)=0$ if and only if there exist  $l\in\N$ and $A_0,\dots,A_n\in K\{Y\}$ such that $S^l_P Q = A_0P+A_1P'+\cdots+A_nP^{(n)}$.
\end{cor}

\index{composition!differential polynomials}

\subsection*{Composition of differential polynomials}
{\em Let $K$ be a differential field, and $P, Q\in K\{Y\}$}. We prove here some
elementary facts about $P(Q)$.

\begin{lemma}\label{lem:comp2}
For each $d\in\N$ we have
$(P^{(n)})_d=(P_d)^{(n)}$. If $P$ is homogeneous, then so is~$P^{(n)}$.
If $P\notin K$, then
$\deg P^{(n)}=\deg P$ and $\wt P^{(n)}=(\wt P)+n$, and
$$ \val P^{(n)}\geq \val P, \qquad \wv P^{(n)}\geq\wv P.$$
\end{lemma}
\begin{proof}
For any word $\bsigma=\sigma_1\cdots\sigma_d\in\{0,\dots,r\}^*$ of length $\abs{\bsigma}=d$ we have
$$(Y^{[\bsigma]})'\ =\ \sum_{i=1}^d Y^{(\sigma_i+1)}\big(Y^{[\bsigma]}/Y^{(\sigma_i)}\big),$$
so $(P')_d=(P_d)'$, which gives the claim about homogeneous $P$. 
Let $P\notin K$. Then $P'\ne 0$ by Lemma~\ref{lem:comp1}, so $\deg P' =\deg P$. With $w:= \wt P$
we have $\wt P' \le w+1$. It remains to show that $(P')_{[w+1]}\ne 0$. 
Let $K_c$ be the field $K$ with the trivial derivation. 
Then $K_{c}\{Y\}$ and $K\{Y\}$ have the same underlying ring,
and $(P')_{[w+1]}$ equals in this underlying ring
the derivative of
$P_{[w]}$ in $K_c\{Y\}$. As $P_{[w]}\notin K$, this gives 
$(P')_{[w+1]}\ne 0$.
Verifying the claims about $\val P^{(n)}$ and $\wv P^{(n)}$ is routine.
\end{proof}

\noindent
Thus we let $P^{(n)}_d$ denote both $(P^{(n)})_d$ and $(P_d)^{(n)}$.

\begin{lemma}\label{lem:comp3}
If $P\ne 0$ and $Q\notin K$, then $P(Q)\neq 0$.
\end{lemma}
\begin{proof} If $Q\notin K$, then
$Q$ is $\d$-transcendental over $K$ by Lemma~\ref{dtrdalg}.\end{proof}

\noindent
Thus if $Q\notin K$, then $P\mapsto P\circ Q= P(Q)$ is an injective 
differential ring endomorphism of $K\{Y\}$, and so extends uniquely to a 
differential field embedding  
$$F \mapsto F\circ Q= F(Q)\colon\  K\<Y\>\to K\<Y\>.$$

\begin{cor}\label{degundercomp}
Suppose $P\ne 0$, $Q\notin K$, $\deg(P)=d$, $\deg(Q)=e$. Then
$$\deg(P\circ Q)\ =\ d \, e, \qquad (P\circ Q)_{de}\ =\ P_d \circ Q_e.$$
If $P$ and $Q$ are homogeneous, then so is $P\circ Q$.
\end{cor}
\begin{proof}
Let $\i=(i_0,\dots,i_r)\in\N^{1+r}$.
By Lemma~\ref{lem:comp2} we have $\deg Q^{(n)}=e$ for all~$n$, so
\begin{align*}
Q^{\i}\ 	&=\ Q^{i_0}(Q')^{i_1}\cdots(Q^{(r)})^{i_r}\\
			&=\ Q_e^{i_0}(Q_e')^{i_1}\cdots(Q_e^{(r)})^{i_r}+R\quad\text{where $\deg(R)<\abs{\i}e$.}
\end{align*}
Thus $\deg (P\circ Q)\leq d\,e$, and $(P\circ Q)_{de}=P_d\circ Q_e$, and so 
by Lemma~\ref{lem:comp3} we get $\deg(P\circ Q) = d\,e$.
\end{proof}

\begin{cor}\label{wtundercomp}
Suppose $P\ne 0$ has order $\le r$ and $Q\notin K$. Then
$$\wv(P\circ Q)\ \geq\ (\val P)(\wv Q).$$
Moreover, with $\i$ ranging over $\N^{1+r}$,
$$\wt(P\circ Q)\ =\ \max\big\{\abs{\i}(\wt Q)+\dabs{\i}:\ P_{\i}\neq 0\big\},$$
and so if $P$ is homogeneous or isobaric, then
$$\wt(P\circ Q)\ =\ (\deg P)(\wt Q)+(\wt P).$$
\end{cor}
\begin{proof} Let $\i$ range over $\N^{1+r}$. Lemma~\ref{lem:comp2} gives $\wv Q^{(\i)} \geq \abs{\i}(\wv Q)$ for all $\i$, so $\wv(P\circ Q) \geq (\val P)(\wv Q)$. Set 
$$w:=\wt Q,\qquad \mu:=\max\big\{\abs{\i}w+\dabs{\i}:P_{\i}\neq 0\big\}.$$
By Lemma~\ref{lem:comp2}, $\wt Q^{(n)}=w+n$, so with 
$Q_{[w+n]}^{(n)}:= (Q^{(n)})_{[w+n]}$ we have 
\begin{align*}
Q^{(\i)}\ 	&=\ Q^{i_0}(Q')^{i_1}\cdots(Q^{(r)})^{i_r}\\
			&=\ Q_{[w]}^{i_0}(Q_{[w+1]}')^{i_1}\cdots(Q_{[w+r]}^{(r)})^{i_r}+S\quad\text{where $\wt(S)<\abs{\i}w+\dabs{\i}$}.
\end{align*}
Thus $\wt (P\circ Q)\leq \mu$.  Let $K_c$ be the field $K$ with the trivial
derivation. Then~$K\{Y\}$ and~$K_c\{Y\}$ have the same underlying ring, 
and $Q_{[w+n]}^{(n)}\in K\{Y\}$ equals in this underlying ring the $n$th derivative of $Q_{[w]}$ in $K_c\{Y\}$. So
$(P\circ Q)_{[\mu]}\in K\{Y\}$ equals in this underlying ring the composition
$\tilde{P}\circ (Q_{[w]})$ computed in $K_c\{Y\}$, where
$$ \tilde{P}\ :=\ \sum_{\abs{\i}w+\dabs{\i}=\mu} P_{\i}\, Y^{\i}\ \neq\ 0.$$
Thus $\wt(P\circ Q) = \mu$ by Lemma~\ref{lem:comp3}. \end{proof}

\begin{cor}\label{cor:comp is associative} Let $F\in K\<Y\>$, and
suppose $P, Q\notin K$. Then 
$$P\circ Q\notin K, \qquad
F\circ (P\circ Q)\ =\ (F\circ P)\circ Q.$$
\end{cor}
\begin{proof} We have $P\circ Q\notin K$ by Corollary~\ref{degundercomp}. The 
associativity of the composition operator on $K\{Y\}$ then leads to
$F\circ (P\circ Q)\ =\ (F\circ P)\circ Q$. 
\end{proof}

\subsection*{Substituting powers}
{\em Let $K$ be a differential field and $c\in C^\times$}.  We consider here the
effect of substituting a power $y^c$ in a differential polynomial, where $y$ and $y^c$ are nonzero elements of a differential field
extension $L$ of $K$ such that $(y^c)^\dagger=cy^\dagger$.

\begin{lemma}\label{lem:derivatives of powers}
For each $n$ there is a homogeneous and isobaric $E_n\in C\{Y\}$ of order~$n$, degree $n$, and weight $n$, such that for any
differential field extension $L$ of~$K$ and any elements $y,z\in L^\times$,
$$z^\dagger\ =\  cy^\dagger\ \Longrightarrow\ z^{(n)}\ =\ z \cdot\frac{E_n(y)}{y^n}.$$
\end{lemma}
\begin{proof}
By induction on $n$. We can take $E_0=1$, and $E_1=cY'$.
Suppose $E_n\in C\{Y\}$ has the required property for a certain 
$n\geq 1$, and let $L$ be a differential field extension of $K$
and let $y,z\in L^\times$ satisfy $z^\dagger= cy^\dagger$.
Then 
\begin{align*}
z^{(n+1)}\ 	&=\ z'\cdot \frac{E_n(y)}{y^n} + z\cdot \left(\frac{E_n(y)}{y^n}\right)' \\
				&=\ z  \cdot \frac{cy'E_n(y)}{y^{n+1}}  + z \cdot \frac{y^n E_n(y)'-ny^{n-1}y'E_n(y)}{y^{2n}}\\
				&=\ z\cdot \frac{(c-n)y'E_n(y) + yE_n(y)'}{y^{n+1}}.
\end{align*}
Thus $E_{n+1}(Y):=(c-n)Y'E_n(Y) + YE_n(Y)'\in C\{Y\}$ has the desired property for $n+1$ in place of $n$, in view of  Lemmas~\ref{lem:comp1} and \ref{lem:comp2}. 
\end{proof}

\begin{cor}
Let $P\in K\{Y\}^{\neq}$ be homogeneous and isobaric, and $d=\deg(P)$, $w=\wt(P)$. Then there is 
a homogeneous and isobaric $E\in K\{Y\}^{\neq}$ of degree~$w$ and weight~$w$ such that for any
differential field extension $L$ of $K$ and ${y,z\in L^\times}$,  
$$z^\dagger\ =\ cy^\dagger\ \Longrightarrow\ P(z)\ =\ z^d\cdot \frac{E(y)}{y^w}
.$$
\end{cor}

\noindent
This follows easily from Lemma~\ref{lem:derivatives of powers}.

\begin{cor}\label{cor:derivatives of powers}
Let $P\in K\{Y\}^{\neq}$ be homogeneous, $d=\deg(P)$, $w=\wt(P)$. Then
there is a homogeneous $E\in K\{Y\}^{\ne}$ of degree $w$ and weight $w$ such that for any
differential field extension $L$ of $K$ and $y,z\in L^\times$,
$$z^\dagger\ =\ cy^\dagger\ \Longrightarrow\ 
P(z)\ =\ z^d\cdot \frac{E(y)}{y^{w}}.$$
\end{cor}
\begin{proof} Apply the previous corollary to the isobaric parts of $P$.
\end{proof}

\noindent 
Corollary~\ref{cor:derivatives of powers} in the case $c=-1$ yields:

\begin{cor}\label{cor:P(1/Y)}
If $P\in K\{Y\}^{\ne}$ is homogeneous, then
$$P(1/Y)\ =\ \frac{E(Y)}{Y^{d+w}}\qquad \text{in $K\<Y\>$, $d:= \deg P$,  $w:= \wt P$,}$$
where $E\in K\{Y\}^{\ne}$ is homogeneous of degree $w$, and $\wt(E)=w$.
\end{cor}

\subsection*{Notes and comments}
Ritt's division theorem (Theorem~\ref{thm:Ritt division} above) appears in \cite[p.~5]{Ritt};
see also \cite[Lemma~7.3]{Kaplansky76}.

%% file: mt-4-4.tex
\section{Valued Differential Fields and Continuity}\label{Valdifcon}

\label{p:valued diff field}

\noindent
We define a {\bf valued differential field\/}
to be a differential field $K$ equipped with a valuation ring 
$\mathcal{O}\supseteq \Q$ of $K$. \index{differential field!valued}\index{valued field!differential}\index{field!valued differential}\index{valued differential field} Here are some basic examples:

\begin{enumerate}
\item Let $\k$ be a differential field and $\Gamma$ an ordered abelian group. Then we make the Hahn field $\k\(( t^{\Gamma}\)) $ into a
differential field extension of $\k$ by 
$$\der\!\left(\sum a_{\gamma}t^{\gamma}\right)\ =\ \sum a_{\gamma}'\ t^\gamma\qquad\text{(so $t'=0$)}.$$
We refer to this valued differential field as the 
{\em Hahn differential field~$\k\(( t^{\Gamma}\)) $}. It 
has constant field  $C_{\k}\(( t^\Gamma\)) $. 
\item Let $\k$ be a differential field, and let $g\in \k\(( t^{\Q}\)) $. We
make the Hahn field~$\k\(( t^{\Q}\)) $ into a valued differential field and a differential field extension of $\k$ by
$$\der\!\left(\sum a_{q}t^{q}\right)\ =\ \sum a_{q}'t^q + \left(\sum qa_q t^{q-1}\right)g\qquad\text{(so $t'=g$)}.$$
\item Let $C$ be a field of characteristic zero (so $C\supseteq \Q$, with the usual identification). We make the
field $C(x)$ of rational functions in $x$ over $C$
into a differential field
with constant field $C$ and $x'=1$. Let $v\colon C(x)^\times \to \Z$ be the discrete valuation on the field $C(x)$ given by $v(C^\times)=\{0\}$ and $v(x)=-1$. Then $C(x)$ equipped with the valuation ring $\mathcal{O}_v$ of 
$v$ is a valued differential field. 
\item $\T$ and $\T_{\log}$: see Appendix~\ref{app:trans}.
\end{enumerate}

\noindent
In the rest of this section we fix a 
valued differential field $K$. We begin with some easy observations:

\begin{lemma}\label{easyobs} 
Let $a,b\in K$, $b\ne 0$. Then
$$ a'\preceq a,\ b'\preceq b\ \Longrightarrow\ (ab)'\preceq ab, \quad (1/b)'\preceq 1/b,\quad  (a/b)'\preceq a/b.$$
In particular, if $b'\preceq b$, then $(b^k)'\preceq b^k$ for every $k\in \Z$. 
\end{lemma} 

\noindent
Our valued differential fields often have the property that the 
maximal ideal of its valuation ring is closed under the 
derivation: then the valuation ring is also closed under the derivation 
and various other useful inclusions hold:

\begin{lemma}\label{closed under der}
Suppose that $\der \smallo\subseteq \smallo$.
Then also $\der \mathcal{O}\subseteq \mathcal{O}$. Moreover, 
let $a\in K$ be such that $a'\preceq a$. Then $\der(a^m\smallo)\subseteq a^m\smallo$ and $\der(a^m\mathcal{O})\subseteq a^m\mathcal{O}$ for all $m$,
hence
$$\der^n(a^m\smallo)\subseteq a^m\smallo,\quad \der^n(a^m\mathcal{O})\subseteq a^m\mathcal{O}\qquad\text{for all $m$,~$n$,}$$
and if $a'\prec a$, then 
$\der^n(a)\prec a$ for all $n\geq 1$.
\end{lemma}
\begin{proof}
Let $x\in \mathcal{O}$, and suppose $x'\notin \mathcal{O}$. Set 
$y:=1/x'$, $z:=x/x'$. Then $y,z\in \smallo$, so 
$y',z'\in \smallo$. But $z=xy$, so
$z'=xy'+1$, a contradiction. This proves the first assertion. For  the second,
let $a\in K$ with $a'\preceq a$. Then
$\der(a\smallo) \subseteq a'\smallo + a\der\smallo\subseteq a\smallo$ and
$\der(a\mathcal{O}) \subseteq a'\mathcal{O} + a\der \mathcal{O} \subseteq
a\mathcal{O}$, so we are done for $m=1$. Now use that 
$a' \preceq a$ gives $(a^m)'\preceq a^m$ for all $m$ by the previous lemma.
If $a'\prec a$, then
$\der(a)\in a\smallo$ and hence $\der^n(a) \in \der^{n-1}(a\smallo) \subseteq a\smallo$ for $n\geq 1$.
\end{proof}

\noindent
If $\der\smallo\subseteq \smallo$, then by Lemma~\ref{closed under der} the
residue field $\k=\mathcal{O}/\smallo$ is a differential field with
derivation $a+\smallo \mapsto a' + \smallo$ ($a\in \mathcal{O}$). (From 
$\Q\subseteq \mathcal{O}$ we get that $\k$ is of characteristic zero, and
so with the obvious identification $\Q$ is a subfield of $\k$.)
The derivation
$\der$ of~$K$ is said to be {\bf small} \label{p:small} \index{valued differential field!with small derivation}
\index{derivation!small} if
$\der\smallo\subseteq \smallo$; in that case we refer to $\k$ with the induced derivation as the {\bf differential residue field\/} of $K$. 
Note that if $K$ has small derivation, then so does every valued differential subfield of $K$. 

The derivation of the Hahn differential
field $\k\(( t^{\Gamma}\)) $ of Example (1) above is small, and
its differential residue field is isomorphic to $\k$. In Example (2)
above, the derivation of $\k\(( t^\Q\)) $ is small if $vg\ge 1$, in which case its differential residue field is again isomorphic to $\k$.
The derivation of $C(x)$ in Example (3) is small, and its differential residue field is isomorphic to $C$ with the trivial derivation. The derivations of~$\T$ and~$\T_{\log}$ of Example (4) are also small, and in both cases the differential residue field is isomorphic to $\R$ with the trivial derivation.

\begin{lemma}\label{Cohn-Lemma}
Suppose $\der \smallo\subseteq \smallo$, and $y\in K$. Then
$$ y\preceq 1 \Rightarrow (y')^2 \preceq y, \qquad y \succeq 1 \Rightarrow (y')^2 \preceq y^3.$$ 
If moreover $\der \cal{O} \subseteq \smallo$, then  
$$0\ne y\preceq 1 \Rightarrow (y')^2 \prec y, \qquad y \succeq 1 \Rightarrow (y')^2 \prec y^3.$$ 
\end{lemma}
\begin{proof}
Let $y\preceq 1$, and
suppose for a contradiction that $(y')^2 \succ y$. Then 
$y=\varepsilon(y')^2$ with $\varepsilon \prec 1$; differentiating yields 
$y'=\varepsilon'(y')^2+2\varepsilon y'y''$, contradicting 
$y'\succ \varepsilon'(y')^2$ and
$y' \succ \varepsilon y'y''$.
Next, suppose $y \succeq 1$; then $y^{-1}\preceq 1$, hence
$(y'/y^2)^2=\big((y^{-1})'\big)^2\preceq y^{-1}$, so $(y')^2\preceq y^3$.
The proof of the part assuming $\der \cal{O} \subseteq \smallo$ is similar.
\end{proof}

\begin{cor}\label{cor:Cohn-Lemma}
Suppose $\der\smallo\subseteq \smallo$.
Let $\Delta$ be a convex subgroup of $\Gamma$, and $\dot{\smallo}$ the maximal ideal of the valuation ring $\dot{\mathcal{O}}$ of the $\Delta$-coarsening of $K$. Then $\der \dot{\smallo}\subseteq \dot{\smallo}$.
\end{cor}
\begin{proof} If $y\in \dot{\smallo}$, then $(y')^2\in \dot{\smallo}$
by Lemma~\ref{Cohn-Lemma}, so $y'\in \dot{\smallo}$.
\end{proof}

\subsection*{Continuity} We now give the valued differential field
$K$ its valuation topology. While not part of our definition of 
{\em valued differential field}, we are only interested in the case that
 $\der\colon K \to K$ is continuous. Note that 
if $\der$ is continuous, then for every differential polynomial 
$P\in K\{Y\}$ the function $y\mapsto P(y)\colon K \to K$ is continuous.
Here is a slightly stronger version of this fact:

\begin{lemma}\label{c-condifpol} Suppose $\der\colon K \to K$ is continuous. Then for each $P(Y)\in K\{Y\}$ the function $y\mapsto P(y)\colon K \to K$ is c-continuous.
\end{lemma}
\begin{proof} Since $\der$ is continuous at $0$ and additive, $\der$  is even uniformly continuous. Hence for every $n$ the map $y\mapsto y^{(n)}\colon K \to K$
is uniformly continuous and thus c-continuous. The c-continuity of differential polynomial functions on $K$ now follows easily from
Lemmas~\ref{lem:addition of c-sequences} and \ref{lem:multiplication of c-sequences}.
\end{proof}

\begin{lemma}\label{cd1} If $\der\smallo\subseteq \smallo$, then
$\der\colon K \to K$ is continuous.
\end{lemma}
\begin{proof} Assume $\der\smallo\subseteq \smallo$. If $\gamma\in \Gamma^{>}$,
and $y\in K$, $vy> 2\gamma$, then $vy' > \gamma$, 
by Lemma~\ref{Cohn-Lemma}. So $\der$ is continuous at $0$, and
since $\der$ is additive, it is continuous.
\end{proof}

\noindent
Up to a factor from $K^\times$ the condition $\der\smallo\subseteq \smallo$
captures exactly continuity of $\der$: 

\begin{lemma}\label{cd2} The following conditions on $K$ are equivalent: \begin{enumerate}
\item[\textup{(i)}] $\der\colon K \to K$ is continuous;
\item[\textup{(ii)}] for some $a\in K^\times$ we have $\der \smallo \subseteq a\smallo$;
\item[\textup{(iii)}] for some $a\in K^\times$ we have
$\derdelta\smallo \subseteq \smallo$, for $\derdelta:= a^{-1}\der$. 
\end{enumerate}
\end{lemma}
\begin{proof} As to (i)~$\Rightarrow$~(ii), assume $\der$ is continuous.
Take $\beta\in \Gamma$ such that $f'\in \smallo$ for all $f\in K$ with
$vf> \beta$, and take $b\in K$ with $vb=\beta$. Then for $g\in \smallo$ we have
$v(bg)> \beta$, so $(bg)'=b'g+bg'\in \smallo$, hence 
$g'\in b^{-1}\smallo+b^\dagger\smallo$. Taking $a\in K^\times$ with $va= \min(-vb, vb^\dagger)$ we obtain $\der\smallo\subseteq a\smallo$. The implication
(ii)~$\Rightarrow$~(iii) is obvious, and (iii)~$\Rightarrow$~(i) follows from 
Lemma~\ref{cd1} applied to $\derdelta$.
\end{proof}

\noindent
A favorable situation is when $\der$ is small {\em and\/} the
derivation of the differential residue field $\k$ is nontrivial: this
often allows properties of the differential residue field $\k$ to be
lifted to useful information about the valued differential field $K$. 

In this connection, consider the set $A$ of all $a\in K^\times$ such that
$\der \smallo\subseteq a\smallo$ (that is, $a^{-1}\der$ is small), and the derivation of the differential residue field of $(K, a^{-1}\der)$ is nontrivial. Here $(K, a^{-1}\der)$ denotes the valued differential
field $K$ with derivation~$a^{-1}\der$ instead of~$\der$.
If $a\in A$ and $a\asymp b$ in $K$, then clearly $b\in A$. We wish to record the following observation:

\begin{lemma}\label{obun} 
If $A\ne \emptyset$, then $v(A)$ consists of just one element.
\end{lemma}
\begin{proof} Suppose $\der$ is small and the
derivation of the differential residue field $\k$ is nontrivial, and let
$a\in K^\times$. It suffices to note that if $a\prec 1$, then 
$a^{-1}\der$ is no longer small, and if $a\succ 1$, then $a^{-1}\der$ is small, but $a^{-1}\der\mathcal{O}\subseteq a^{-1}\mathcal{O}\subseteq \smallo$, so the derivation of the differential residue field of $(K, a^{-1}\der)$ is trivial.
\end{proof}

\begin{lemma}\label{c2} Assume that $C\subseteq \mathcal{O}$.
Let $P(Y)=
a_0Y+\cdots+a_nY^{(n)}$, with all $a_i\in K$, $a_n\neq 0$. Then each level set
$P^{-1}(s)$ \rom{(}$s\in K$\rom{)} is a discrete subset of $K$.
\end{lemma}
\begin{proof} Such a level set is empty or a translate of the $C$-linear subspace $V:=P^{-1}(0)$ of $K$. Any $y_0,\dots, y_n\in K$ with $y_0 \succ y_1 \succ \cdots \succ
y_n \succ 0$ are linearly independent over $C$. Since $\dim_C V \le n$, there are no such $y_0,\dots, y_n\in V$.
\end{proof}

\begin{lemma}\label{c3} Suppose $\der$ is nontrivial, and
the valuation of $K$ is nontrivial. Then
no differential polynomial $P(Y)\in K\{Y\}^{\ne}$
vanishes identically on any nonempty open subset of $K$.
\end{lemma}
\begin{proof} Assume towards a contradiction that
$P(Y)\in K\{Y\}^{\ne}$ vanishes identically on
$\big\{y\in K:\ v(y-a)\ge \gamma\big\}$, $a\in K$, $\gamma\in \Gamma$. 
Take $g\in K$ with $vg=\gamma$ and set $Q:= P(a+gY)\in K\{Y\}^{\ne}$.
Then $Q$ vanishes identically on $\mathcal{O}$, and so 
$Q(Y)\cdot Q(Y^{-1})\cdot Y^n$ is for sufficiently large $n$ a nonzero
differential polynomial that vanishes identically on~$K$. As the derivation
of $K$ is nontrivial, this is impossible by Lemma~\ref{notthin}.
\end{proof}

\subsection*{Completion}
Recall from Section~\ref{sec:pc in valued fields} that any valued field $E$ can be {\em completed\/}: it is
dense in a valued field extension
$E^{\operatorname{c}}$ such that for each valued field extension~$E\subseteq F$ with $E$ dense
in~$F$ there is a unique valued field embedding $F\to E^{\operatorname{c}}$ that is the
identity on~$E$. (Here ``dense'' is with respect to the relevant valuation topology.) These
properties determine~$E^{\operatorname{c}}$ up to unique valued field isomorphism over $E$,
and~$E^{\operatorname{c}}$ is called the {\bf completion of~$E$}.
Recall that $E^{\operatorname{c}}|E$ is an immediate extension. 

\index{completion!valued differential field}

\begin{lemma}\label{ExtendingDerivation} Suppose $\der\colon K \to K$ is continuous.
Then there is a unique continu\-ous derivation on $K^{\operatorname{c}}$ 
extending the deri\-vation
of $K$. 
\end{lemma}
\begin{proof} The derivation of $K$ being additive, it is even uniformly 
continuous, that is, for each $\gamma\in \Gamma$ there is 
$\delta\in \Gamma$
such that whenever $x,y\in K$ and $v(x-y) > \delta$, then
$v(x'-y') >\gamma$. It follows that $\der$
extends uniquely to a continuous map 
$K^{\operatorname{c}}\to K^{\operatorname{c}}$ (Lemma~\ref{lem:extend to closure}); this
map is a derivation.  
\end{proof}

\noindent
Suppose $K$ is as in the lemma. Then we consider $K^{\operatorname{c}}$ as the valued differential field
whose derivation is the unique continuous derivation on $K^{\operatorname{c}}$ that
extends the one of~$K$. If $K\subseteq L$ is a valued differential field extension
such that $K$ is dense in $L$ and the derivation of $L$ is continuous
with respect to the valuation topology, then the unique valued field embedding
$L \to K^{\operatorname{c}}$ that is the identity on $K$ is a differential 
field embedding. 

\begin{cor}\label{completionsmall} If $K$ has small derivation, then so does $K^{\c}$.
\end{cor}

\subsection*{Traces and norms} Let $L|K$ be an extension of 
valued differential fields with $[L:K]<\infty$ and put $n:=[L:K]$. We let $K^\alg$ be an 
algebraic closure of $K$ equipped with the unique derivation extending the 
derivation of $K$. 
For $f\in L$, let
$$
\operatorname{tr}_{L|K}(f)\ :=\ \sum_{i=1}^n \sigma_i(f), \qquad
\operatorname{N}_{L|K}(f)\ :=\ \prod_{i=1}^n \sigma_i(f)
$$
denote the trace of $f$ in $L|K$ and the norm of $f$ in $L|K$, respectively. 
Here $\sigma_1,\dots,\sigma_n$ are the distinct field embeddings 
$L\to K^\alg$ which are the identity on $K$. Note that each~$\sigma_i$ 
is an embedding $L\to K^\alg$ of differential fields. Hence for all $f\in L$:
$$\operatorname{tr}_{L|K}(f')\ =\ \operatorname{tr}_{L|K}(f)', \qquad
  \operatorname{tr}_{L|K}(f^\dagger)\ =\
\operatorname{N}_{L|K}(f)^\dagger\ \text{ if $f\neq 0$.}$$
We pick a valuation ring of $K^\alg$ lying over the
valuation ring of $K$ to make $K^\alg$ a valued field extension of $K$.
If $K$ is henselian, then each $\sigma_i$ is a valued field embedding,
and after identifying $L$ with a valued subfield of $K^\alg$ via $\sigma_1$, 
say, we have $\sigma_i f \asymp f$ for $f\in L$ and $i=1,\dots,n$, hence 
(for $f\in L$):
$$ \operatorname{tr}_{L|K}(f)\ \preceq\ f, \qquad 
\operatorname{N}_{L|K}(f)\ \asymp\ f^n.$$
It follows that for henselian $K$ and $f\in L^\times$ we have
$$ \operatorname{N}_{L|K}(f)^\dagger\ =\  \operatorname{tr}_{L|K}(f^\dagger)\ \preceq f^\dagger.$$

\subsection*{Monotonicity} We say that the valued differential field
$K$ is {\bf monotone\/} if $a'\preceq a$ for all $a\prec 1$ in $K$; in that case the derivation of $K$ is small, and $a'\preceq a$
for all~$a\in K$, by Lemma~\ref{easyobs}. The Hahn differential field 
$\k\(( t^{\Gamma}\)) $ of Example~(1) above is monotone. The valued differential field $\k\(( t^{\Q}\)) $ of Example~(2) above is monotone if~$vg\ge 1$. Also $C(x)$ of Example~(3) is monotone. In Example~(4),~$\T_{\log}$ is monotone, but $\T$ is not. 
Note that if~$K$ is monotone, then so is every valued differential subfield of $K$, and every coarsening of $K$ by a convex subgroup of~$\Gamma=v(K^\times)$. 

\label{p:monotone}

\index{valued differential field!monotone}
\index{monotone!valued differential field}

\subsection*{Many constants and few constants} For later use we introduce the following conditions, where as usual $C=C_K$ and $\Gamma=v(K^\times)$: \begin{enumerate}
\item $K$ has {\bf many constants\/} if $v(C^\times)=\Gamma$;
\item $K$ has {\bf few constants\/} if $v(C^\times)=\{0\}$, equivalently, $C\subseteq \mathcal{O}$.
\end{enumerate}
The Hahn differential field $\k\(( t^{\Gamma}\)) $ of Example (1) above has many constants. The valued differential fields of Examples (3) and (4)
have few constants. If $K$ has small derivation and many constants,
then $K$ is monotone. 

\label{p:many/few constants}
\index{valued differential field!many constants}

If $K$ has many constants, and $L$ is a valued differential field extension of $K$
with $\Gamma_L=\Gamma$, then $L$ also has many constants.
%Conversely, every valued differential field  $K$
%has a valued differential field extension  $L$ with many constants and $\Gamma_L=\Gamma$.
%(By Corollary~\ref{lem:derivative on trans ext} 
%and Lemma~\ref{lem:gauss}.)
Although we don't really need this, here are two easy results on 
extensions with possibly bigger value group.

\begin{lemma} Suppose $K$ is henselian, $K$ has many constants, 
$\k$ is real closed or algebraically closed, and $L$ is a valued differential field extension of $K$ and algebraic over $K$. Then $L$ has many constants.
\end{lemma}
\begin{proof} Let $a\in L^{\times}$, and take $n\ge 1$ such that 
$nv(a)\in \Gamma$. Then $a^n\asymp c$ with~$c\in C^\times$, so $a^n=bc$
with $b\in L^\times$, $b\asymp 1$. Now $L$ is also henselian with real closed or 
algebraically closed residue field, so $b=d^n$ or $-b=d^n$ with $d\in L^\times$.
Then $c=(a/d)^n$ or $-c=(a/d)^n$, and thus $a/d\in C_L$ and $v(a)=v(a/d)$.
\end{proof}

\begin{lemma} Suppose $K$ is monotone with many constants, and let  
$\Gamma+ \Z\beta$ be an 
ordered group extension of $\Gamma=v(K^\times)$ such that $n\beta\notin \Gamma$ for all $n\ge 1$. Let $K(b)$ be a field extension of $K$ with $b$ 
transcendental over $K$ and make $K(b)$ into a valued differential 
field extension of $K$ by requiring $b'=0$ and by extending the valuation
of $K$ to a valuation $v\colon K(b)^\times \to \Gamma+ \Z\beta$ with
$vb=\beta$. Then $K(b)$ is monotone, $C_{K(b)}=C(b)$, and $K(b)$ 
has many constants.
\end{lemma}
\begin{proof} For $a_0,\dots, a_n\in K$, not all zero, and $P=\sum_i a_iY^i$ 
we have
\begin{align*} vP(b)\ &=\ v(a_0+ a_1b + \cdots + a_n b^n)\ =\
\min_i \big(v(a_i)+i\beta\big),\\
  P(b)'\ &=\  P^{\der}(b),\ \text{ so }\ vP(b)'\ =\ \min_i\big(v(a_i')+i\beta\big),
\end{align*}
from which we get $P(b)'\preceq P(b)$. Thus $K(b)$ is monotone. 

Let $P,Q\in K[Y]^{\ne}$ be coprime, with $Q$ monic, such that 
$P(b)/Q(b)\in C_{K(b)}$.
Then $P^{\der}(b)Q(b)=P(b)Q^{\der}(b)$, so $Q|Q^{\der}$ in $K[Y]$. Since
$\deg Q^{\der}< \deg Q$, this gives~$Q^{\der}=0$, and so $P^{\der}=0$ as well.
Then $P,Q\in C[Y]$, so $P(b)/Q(b)\in C(b)$. This gives  $C_{K(b)}=C(b)$, and also shows that $K(b)$ has many constants. 
\end{proof}

\subsection*{Classifying pc-sequences} Let $(a_{\rho})$ be a pc-sequence in $K$. 

\index{pc-sequence!differential-algebraic type}
\index{pc-sequence!differential-transcendental type}
\index{pc-sequence!minimal differential polynomial}
\index{minimal!differential polynomial}

\medskip\noindent
We say that $(a_\rho)$ is of {\bf differential-algebraic type over $K$} (for short: {\bf $\d$-algebraic type over $K$}) if
$G(b_{\lambda})\leadsto 0$ for some $G(Y)\in K\{Y\}$ and some pc-sequence
$(b_{\lambda})$ in~$K$ equivalent to $(a_{\rho})$. A {\bf minimal differential
polynomial of} $(a_{\rho})$ over $K$ is a differential polynomial $G(Y)\in K\{Y\}$ with the following properties: \begin{enumerate}
\item $G(b_{\lambda})\leadsto 0$ for some pc-sequence $(b_{\lambda})$ in $K$ equivalent to $(a_{\rho})$ (so $G\notin K$);
\item $H(b_{\lambda})\not\leadsto 0$ whenever $H\in K\{Y\}$ has lower complexity than $G$ and the pc-sequence $(b_{\lambda})$ in $K$ is equivalent to $(a_{\rho})$.
\end{enumerate}
Thus $(a_{\rho})$ is of $\d$-algebraic type over $K$ if and only if 
$(a_{\rho})$ has a minimal differential polynomial over $K$.
If $G$ is a minimal differential polynomial of $(a_{\rho})$ over $K$ and $a\in K$, then $G_{+a}$ is a minimal differential polynomial of
$(a_{\rho}-a)$.

\medskip\noindent
We say that $(a_{\rho})$ is of {\bf differential-transcendental type over $K$} (for short: {\bf $\d$-trans\-cen\-den\-tal type over $K$}) if it is not of $\d$-algebraic type over $K$, that is, $G(b_{\lambda})\not\leadsto 0$ for each $G\in K\{Y\}$ and each pc-sequence 
$(b_{\lambda})$ in $K$ equivalent to $(a_{\rho})$.

\subsection*{Notes and comments}
Lemmas~\ref{easyobs},~\ref{closed under der}, and~\ref{Cohn-Lemma} are 
partly borrowed from Section~2 of \cite{Cohn}, but the terminology there
is different from ours. The property of monotonicity appears in \cite{Cohn}; see also \cite[Proposition~2.2]{Morrison}. Monotonicity 
together with having many constants is a key assumption in Scanlon's~\cite{Scanlon}.

%% file: mt-4-5.tex
\section{The Gaussian Valuation}\label{The Gaussian Valuation}

\noindent
Let $K$ be
a valued differential field. Then we extend its valuation $v$
to the valuation 
$v\colon K\<Y\> \to \Gamma_{\infty}$ on $K\<Y\>$  by 
  $$v(P):=\min_{\i} v(P_{\i}) = 
\min_{\bomega}v(P_{[\bomega]}) \in \Gamma_{\infty},$$
for $P(Y)=\sum_{\i} P_{\i} Y^{\i}\in K\{Y\}$. Because of its familiar
connection to Gauss's lemma on unique factorization, this extended valuation 
is called the {\bf gaussian extension} of~$v$. (See also Section~\ref{sec:valued fields}.) 
Note that for $P,Q\in K\{Y\}$ we have
$v(P+Q) = \min\bigl(v(P), v(Q)\bigr)$ whenever $v(P)\ne v(Q)$, and also
whenever $P$, $Q$ have no common monomials, that is, for all $\i$ either 
$P_{\i}=0$ or $Q_{\i}=0$. Hence $v(P)=\min_d v(P_d)$.

\index{valuation!gaussian extension}
\index{gaussian extension}
\index{extension!gaussian}

\medskip\noindent
\textit{Recall that $\mathcal{O}$ is the valuation ring of $K$ with
maximal ideal $\smallo$. In the rest of this section we assume 
$\der \smallo\subseteq \smallo$ \textup{(}and hence $\der \mathcal{O}\subseteq \mathcal{O}$\textup{)}. Also, $\phi$, $f$, $g$ range over $K$.}\/

\subsection*{The function $v_P$} 

\begin{lemma}\label{v-under-conjugation} Let $P\in K\{Y\}$. Then we have:
\begin{enumerate}
\item[\textup{(i)}] if $f\preceq 1$, then $P_{+f}\asymp P$; if $f\prec 1$ and $P\ne 0$, then $P_{+f}\sim P$;
\item[\textup{(ii)}] the element $v(P_{\times f})$ of $\Gamma_\infty$ depends only on $vf$;
\item[\textup{(iii)}] if $P\neq 0$ and $P(0)=0$, then the map
$$vf \mapsto v(P_{\times f})\ \colon\ \Gamma_\infty\to\Gamma_\infty$$
is strictly increasing.
\end{enumerate}
\end{lemma}
\begin{proof} 
Assume $f\preceq 1$. From \eqref{Additive-Conj} we obtain
$P_{+f}\preceq P$. Since $P=(P_{+f})_{+g}$ with $g:=-f$, this also gives
$P \preceq P_{+f}$. Next, assume $f\prec 1$ and $P\ne 0$. 
Then \eqref{Additive-Conj} yields $P_{+f}=P+Q$ with $Q\prec P$, and so $P_{+f}\sim P$. This proves (i). 

Suppose that $v\phi=0$.
By \eqref{Multiplicative-Conj}, this gives
$v(P_{\times \phi}) \ge v(P)$.
Using $P=(P_{\times \phi})_{\times \phi^{-1}}$ we can reverse this 
inequality to
get $v(P) \ge v(P_{\times \phi})$. Hence $v(P)=v(P_{\times \phi})$.
Suppose now that $vf=vg$, and write $f=g\phi$ with $v\phi=0$.
By what we just proved and using $P_{\times f} =(P_{\times g})_{\times \phi}$ 
we obtain
$v(P_{\times f})=v(P_{\times g})$. This shows (ii).
Next, suppose $P\neq 0$, $P(0)=0$ and $vf > vg$. Write
$f= g\phi$ where $v\phi>0$. Then the identity~\eqref{Multiplicative-Conj} 
restricted to $\bsigma$ of length~$>0$, together with
the assumption that $\der\smallo \subseteq \smallo$, yields
$v(P_{\times \phi}) > v(P)$. Using
$P_{\times f} =(P_{\times g})_{\times \phi}$, we conclude 
$v(P_{\times f}) > v(P_{\times g})$.
\end{proof}

\noindent
For $P\in K\{Y\}$ we define $v_P\colon \Gamma_\infty \to \Gamma_\infty$
by $$v_P(\gamma) := v(P_{\times f}) \quad \text{ whenever $vf=\gamma$.}$$
(Thus for $P=0$ the function $v_P$ takes the constant value $\infty$.)
For $P,Q\in K\{Y\}$ we have $$(PQ)_{\times f}= P_{\times f}Q_{\times f}
\quad\text{and}\quad (P+Q)_{\times f}=P_{\times f}+Q_{\times f},$$ so
$$v_{PQ}(\gamma)=v_P(\gamma) + v_Q(\gamma)\quad\text{and}\quad
v_{P+Q}(\gamma)\geq\min\bigl(v_P(\gamma),v_Q(\gamma)\bigr).$$ 
Note also that 
$$v_P(\gamma)=\min_d v_{P_d}(\gamma).$$ 
Thus in some sense the properties
of the functions $v_P$ reduce to the case that $P$ is homogeneous.
In Chapter~\ref{ch:valueddifferential} we consider these functions $v_P$ in more detail,  
and then the following results will be very useful. The first
one says that if the derivation of the differential residue field $\k$ is nontrivial
and $v_P(\alpha)=\beta$, then $v(P(f))=\beta$ for ``almost all'' $f$ with $vf=\alpha$. More precisely:

\nomenclature[V]{$v_P(\gamma)$}{gaussian valuation of $P_{\times f}$, where  $vf=\gamma$}

\begin{lemma}\label{vpnice} Assume the derivation of $\k$ is nontrivial. Let $P\in K\{Y\}^{\ne}$ and $\alpha, \beta\in \Gamma$ be such that
$v_P(\alpha)=\beta$, and let $a\in K$ be such that $va=\alpha$. Then
there is a thin set $S\subseteq \k$ such that $vP(ay)=\beta$ for all $y\asymp 1$
in $K$ with $\overline{y}\notin S$. 
\end{lemma}
\begin{proof} Take $b\in K$ with $vb=\beta$, so $vF=0$ for 
$F:= b^{-1}P_{\times a}$. Then the thin set $S:=\big\{\bar{y}\in \k:\ \bar{F}(\bar{y})=0\big\}\subseteq \k$ has the property that
$F(y)\asymp 1$ for all $y\asymp 1$ in $K$ with $\bar{y}\notin S$, and
so $vP(ay)=\beta$ for all such $y$. 
\end{proof}

\begin{lemma}\label{vplemma} Suppose $g\in K^\times$ and $g'\preceq g$.
Then $g^{(n)}\preceq g$ for all $n$, and with $\gamma:= vg$ we have
$v_P(\gamma)=v(P)+d\gamma$ for homogeneous $P\in K\{Y\}^{\ne}$ 
of degree $d$. 
\end{lemma}
\begin{proof} 
From $g' \preceq g$ and Lemma~\ref{closed under der} we obtain 
$g^{(n)}\preceq g$ for all $n$.
Next, consider the case $d=1$ and $P=Y^{(n)}$. Then 
$$P_{\times g}=(gY)^{(n)}= gY^{(n)} + ng'Y^{(n-1)} + \cdots + g^{(n)}Y,$$
so $v_P(\gamma)=\gamma$. Thus for $\i=(i_0,\dots, i_r)\in \N^{1+r}$, $d=|\i|$, and $P=Y^{\i}$,
$$ P_{\times g}\ =\ g^dY^{\i} + R, \quad R\in K\{Y\},$$
where $vR\ge d\gamma$ and all monomials in $R$ are lower than $Y^{\i}$
in the antilexicographic ordering on the set of monomials. 
For general homogeneous $P\in K\{Y\}^{\ne}$ of degree~$d$ we take among the terms $aY^{\i}$ in $P$ with $va=vP$ the one
for which $Y^{\i}$ is largest in the antilexicographic ordering 
on the set of monomials. For this term $aY^{\i}$ in $P$ we have $|\i|=d$ and
$$ P_{\times g}\ =\ ag^{d}Y^{\i} +  R + S, \quad R,S\in K\{Y\},$$
where $vR\ge vP + d\gamma$, all monomials in $R$ are antilexicographically 
lower than $Y^{\i}$, and $vS > vP + d\gamma$.
\end{proof}

\begin{cor}\label{monvp} If $K$ is monotone and $P\in K\{Y\}^{\ne}$ is 
homogeneous of degree~$d$, then
$v_P(\gamma)\ =\ v(P) + d\gamma$ for all $\gamma\in \Gamma$.
\end{cor}

\noindent
Let $\Delta$ be a convex subgroup of $\Gamma$, and let 
$\dot{v}\colon K^\times \to \dot{\Gamma}=\Gamma/\Delta$ be the $\Delta$-coarsening 
of~$K$ with valuation ring $\mathcal{O}$ and maximal ideal $\dot{\smallo}$ of
$\mathcal{O}$.
Then $\der \dot{\smallo}\subseteq \dot{\smallo}$ by Corollary~\ref{cor:Cohn-Lemma}.
The following is obvious. 

\begin{lemma}\label{coarsevp} For $P\in K\{Y\}^{\ne}$ and $\gamma\in \Gamma$ we 
have $$\dot{v}_P(\dot{\gamma})\ =\ v_P(\gamma) + \Delta.$$
\end{lemma}

\subsection*{Dominant weight}
In this subsection we assume $P,Q\in K\{Y\}^{\ne}$, and set
\begin{align*}
\Pnu(P)\ &:=\ \max \big\{ \|\bsigma\| :\  v(P_{[\bsigma]})=v(P) \big\},\\
\Pmu(P)\ &:=\ \min \big\{ \|\bsigma\| :\  v(P_{[\bsigma]})=v(P) \big\}, 
\end{align*}
so $$\wv(P)\ \leq\ \Pmu(P)\ \leq\ \Pnu(P)\ \leq\ \wt(P).$$ 
We call $\Pnu(P)$ the 
{\bf dominant weight} of $P$ and $\Pmu(P)$ the {\bf dominant weighted multiplicity} of $P$. These quantities will be needed in Section~\ref{evtbeh}.
% connection with Newton diagrams of differential polynomials. 

\index{differential polynomial!dominant!weight}
\index{differential polynomial!dominant!weighted multiplicity}
\index{weight!dominant}
\index{multiplicity!dominant weighted}
\index{weighted multiplicity!dominant}

\nomenclature[V]{$\dwt(P)$}{dominant weight of $P$}
\nomenclature[V]{$\dwv(P)$}{dominant weighted multiplicity of $P$}

\begin{lemma} \label{mu and nu and polynomial multiplication}
$$\Pnu(PQ)=\Pnu(P)+\Pnu(Q), \quad \Pmu(PQ)=\Pmu(P)+\Pmu(Q).$$
\end{lemma}
\begin{proof}
Take $f$ with $vf=v(P)$, so $P_0:=f^{-1}P\in \mathcal{O}\{Y\}$;
then
$\Pnu(P)=\wt(\bar{P_0})$ and $\Pmu(P)=\wv(\bar{P_0})$ where $\bar{P_0}$ is the
image of $P_0$ under the ring morphism 
$$\mathcal{O}\{Y\}\to\k\{Y\} $$ that extends 
the residue
map $a\mapsto\bar{a}\colon \mathcal{O}\to \k:= \res(K)$ and sends $Y^{(n)}$ to $Y^{(n)}$
for all $n$. The lemma now
follows from Lemma~\ref{w and W}.
\end{proof}

\begin{lemma}\label{dwtP}
The quantity $\Pnu(P_{\times g})$ depends only on $vg$, for $g\ne 0$.
\end{lemma}
\begin{proof} Take 
$\bsigma$ such that $\Pnu(P)=\|\bsigma\|$ and $v(P_{[\bsigma]})=v(P)$, and 
let $vg=0$. Then by 
\eqref{Multiplicative-Conj} we have
$v\bigl((P_{\times g})_{[\bsigma]}\bigr)=v(P_{[\bsigma]})=v(P)=v(P_{\times g})$, so 
$\Pnu(P_{\times g}) \geq \Pnu(P)$. To reverse this inequality, use 
$P=(P_{\times g})_{\times g^{-1}}$. This proves the lemma for $vg=0$, and the 
general case follows easily from this special case. 
\end{proof}

\noindent
Thus for $\gamma\in\Gamma$ we can define
$\Pnu_P(\gamma):=\Pnu(P_{\times g})$ where $\gamma=vg$.

\begin{lemma} Suppose $P$ is homogeneous and 
$g\in K^\times$ satisfies $g'\preceq g$. Then we have 
$\dwt(P_{\times g})=\dwt(P)$.
\end{lemma}
\begin{proof} By Lemma~\ref{vplemma} we have $g^{(n)}\preceq g$
for all $n$ and $v(P_{\times g})=v(P)+d\gamma$, where 
$\gamma:= vg$ and $d:= \deg P$. Therefore, if $|\bsigma|=d$ and 
$\|\bsigma\|>\dwt(P)$, then by formula \eqref{Multiplicative-Conj} we have
$v\big((P_{\times g})_{[\bsigma]}\big)> v(P_{\times g})$.
Now take $\bsigma$ such that $|\bsigma|=d$, $\|\bsigma\|=\dwt(P)$, and $v(P_{[\bsigma]})=v(P)$.
Then by \eqref{Multiplicative-Conj},
$$(P_{\times g})_{[\bsigma]}\ =\ P_{[\bsigma]}g^d + f, \qquad
vf\ >\  v(P_{[\bsigma]}g^d)\ =\ v(P) + d\gamma\  =\ v(P_{\times g}),$$ 
so $v\big((P_{\times g})_{[\bsigma]}\big)=v(P_{\times g})$ and thus 
$\dwt(P_{\times g})=\dwt(P)$.
\end{proof}

\noindent
The following is proved just like Lemma~\ref{dwtP}:

\begin{lemma} \label{mu_P}
If $\der \mathcal{O}\subseteq \smallo$, then $\Pmu(P_{\times g})$ 
depends only on $vg$, for $g\neq 0$.
\end{lemma}
%\begin{proof}
%Assume $\der \mathcal{O}\subseteq \smallo$. Take 
%$\bomega$ such that $\Pmu(P)=\|\bomega\|$ and $v(P_{[\bomega]})=v(P)$, and 
%let $vg=0$. Then by 
%\eqref{Multiplicative-Conj} we have
%$v\bigl((P_{\times g})_{[\bomega]}\bigr)=v(P_{[\bomega]})$, so 
%$\Pmu(P_{\times g}) \leq \Pmu(P)$. To reverse this inequality, use 
%$P=(P_{\times g})_{\times g^{-1}}$. This proves the lemma for $vg=0$, and the 
%general case follows easily from this special case. 
%\end{proof}

\noindent
If $\der \mathcal{O}\subseteq \smallo$ and $\gamma\in\Gamma$, then we define
$\Pmu_P(\gamma):=\Pmu(P_{\times g})$ where $\gamma=vg$.

\nomenclature[V]{$\dwt_P(\gamma)$}{dominant weight of $P_{\times g}$, where   $vg=\gamma$}
\nomenclature[V]{$\dwv_P(\gamma)$}{dominant weighted multiplicity of $P_{\times g}$, where $vg=\gamma$}

\begin{lemma}\label{mu_P, more} Suppose 
$\der\mathcal O\subseteq\smallo$ and
$P$ is homogeneous. Let
$g\in K^\times$, $g'\prec g$. Then $g^{(n)}\prec g$ for all $n\geq 1$, and $\Pmu(P_{\times g})=\Pmu(P)$.
\end{lemma}
\begin{proof} By Lemma~\ref{closed under der} we have $g^{(n)}\prec g$ for all $n\geq 1$. By Lemma~\ref{vplemma} we have $v(P_{\times g})=v(P)+d\gamma$, where 
$\gamma:= vg$ and $d:=\deg P$.  
 Take $r\in \N$ such that $P$ has order $\le r$. Below, 
 $\bsigma, \btau\in \{0,\dots,r\}^*$ are subject to $|\bsigma|=|\btau|=d$. Then  
\begin{equation}\label{eq:btau > bsigma}
\btau\geq\bsigma\text{ and }\btau\neq\bsigma\quad\Rightarrow\quad v\big(P_{[\btau]}g^{[\btau-\bsigma]}\big)\ >\ v(P)+ d\gamma.
\end{equation}
So by  \eqref{Multiplicative-Conj} and \eqref{eq:btau > bsigma}, if $\dabs{\bsigma}<\Pmu(P)$, then 
$v\big((P_{\times g})_{[\bsigma]}\big) > v(P_{\times g})$.
Take $\bsigma$ such that $\|\bsigma\|=\Pmu(P)$, and $v(P_{[\bsigma]})=v(P)$.
Then by \eqref{Multiplicative-Conj} and  \eqref{eq:btau > bsigma},
$$(P_{\times g})_{[\bsigma]}\ =\ P_{[\bsigma]}g^d + f, \qquad
vf\ >\  v\big(P_{[\bsigma]}g^d\big)\ =\ v(P) + d\gamma\  =\ v(P_{\times g}),$$ 
so $v\big((P_{\times g})_{[\bsigma]}\big)=v(P_{\times g})$ and thus 
$\Pmu(P_{\times g})=\Pmu(P)$.
\end{proof}

%% file: mt-4-6.tex
\section{Differential Rings}\label{sec:diff rings}

\noindent
We gather here some basic facts about differential rings to
be used at various places. This concerns differential ideals,
simple differential rings, and linear disjointness over constant fields. {\em Throughout this section $R$ is a differential ring}. Recall from Section~\ref{Differential Fields and Differential Polynomials} that this includes $\Q$ being a subring of $R$.

\subsection*{Radical differential  ideals}
A {\bf differential ideal} of $R$ is an ideal $I$ of $R$ 
with $f'\in I$ for all $f\in I$. If $I$ and $J$ are differential ideals of $R$, then so are~$I\cap J$ and~$IJ$.
Note that if $I$ is a differential ideal of~$R$ with $I\ne R$, then $a+I\mapsto a'+I$ is
a derivation on $R/I$, making $R/I$ a differential ring and 
the natural surjection $R\to R/I$ a morphism of differential rings. Conversely, if~$S$ is a differential ring and $h\colon R\to S$ is a morphism of differential rings, then $\ker h$ is a differential   
ideal of $R$.  %If $I$ is a differential ideal of $R$, then $I$ is also a %differential ideal of $R^\phi$  for each $\phi\in R^\times$. 
Given a subset $S$ of $R$, we denote by $[S]$ the smallest differential ideal of~$R$ containing~$S$; so~$[S]$ is the ideal generated by the derivatives~$s^{(n)}$ of the elements~$s\in S$. A {\bf radical differential ideal of $R$\/} is a differential ideal of $R$ that is radical as an ideal of $R$.

\nomenclature[Mx]{$[S]$}{differential ideal generated by $S$}

\index{ideal!differential}
\index{differential ideal}

\medskip
\noindent
{\em In the rest of this subsection $I$ is a differential ideal of $R$}.
Suppose that $S$ is a multiplicative subset of $R$ with
$0\notin S$. Recall from Section~\ref{Differential Fields and Differential Polynomials} the differential ring
morphism $r\mapsto \iota(r):=\frac{r}{1}\colon R\to S^{-1}R$. The ideal $S^{-1}I$ of $S^{-1}R$
generated by $\iota(I)$ is a differential ideal, hence 
$$\iota^{-1}(S^{-1}I)\ =\ \{r\in R:\ \text{$rs\in I$ for some $s\in S$}\}$$
is a differential ideal of $R$ containing $I$.

\begin{lemma}\label{lem:radical of diff ideal}
Let $a\in R$ with $a^n\in I$, where $n\geq 1$. Then $(a')^{2n-1}\in I$.
\end{lemma}
\begin{proof}
By induction on $k=1,\dots,n$ we show that $a^{n-k}(a')^{2k-1}\in I$. The case $k=1$ holds since
$a^{n-1}a'=\frac{1}{n}(na^{n-1}a')=\frac{1}{n}(a^n)'\in I$. Suppose $1\le k < n$ and
$a^{n-k}(a')^{2k-1}\in I$. Differentiating yields
$$(n-k)a^{n-k-1}(a')^{2k}+(2k-1)a^{n-k}(a')^{2k-2}a''\in I,$$
and then multiplying by $\frac{1}{n-k}a'$ gives $a^{n-k-1}(a')^{2k+1}\in I$, as required.
\end{proof}

\noindent
Thus the radical of a differential ideal of $R$  
is a differential ideal of~$R$. 

\begin{lemma}\label{lem:radical prod}
Suppose $I$ is radical and $a,b\in R$, $ab\in I$. Then $a'b\in I$. 
\end{lemma}
\begin{proof} Use $a'b+ab'=(ab)'\in I$ and multiply by $a'b$.
\end{proof}

\begin{cor}
If $I$ is radical and $S\subseteq R$, then $(I:S)$
is a radical differential ideal of $R$.
\end{cor}

\noindent
Given $S\subseteq R$, the smallest radical differential ideal of $R$ containing $S$ is $\sqrt{[S]}$.

\begin{cor}\label{cor:radical mult}
Let $S,T\subseteq R$. Then $\sqrt{[S]}\cdot\sqrt{[T]} \subseteq \sqrt{[ST]}$.
\end{cor}
\begin{proof}
Let $r\in R$. By the previous corollary, $\big(\sqrt{[rT]}:r\big)$ is a radical differential ideal of $R$; it contains $T$
and thus also $\sqrt{[T]}$. Thus $r\sqrt{[T]}\subseteq\sqrt{[rT]}$ for all~$r\in R$.
Hence the radical differential ideal $\big(\sqrt{[ST]} : \sqrt{[T]}\,\big)$ of $R$ contains $S$ and therefore~$\sqrt{[S]}$; it follows that $\sqrt{[S]}\cdot\sqrt{[T]} \subseteq \sqrt{[ST]}$.
\end{proof}

\noindent
A differential ideal $I$ of $R$ is said to be {\bf prime} if it is prime as an
ideal of $R$, that is,~$R/I$ is an integral domain. A differential ideal of $R$ is said to be
{\bf maximal} if it is proper (not equal to $R$) and maximal among proper
differential ideals of $R$ (which does not imply that it is a maximal ideal
of $R$; see the next subsection).
Here is a differential analogue of a well-known fact from commutative algebra:

\index{differential ideal!prime}
\index{differential ideal!maximal}

\begin{prop}\label{prop:radical is intersection of primes, differential}
Every radical differential ideal of $R$ is the intersection in $R$ of a collection of prime differential ideals of $R$.
\end{prop}

\noindent
In the proof we use:

\begin{lemma}\label{lem:disjoint}
Let $S$ be a multiplicative subset of $R$. Let $I$ be a radical differential ideal of $R$ disjoint from $S$ and
maximal with these properties. Then $I$ is prime.
\end{lemma}
\begin{proof}
Let $a,b\in R$ with $ab\in I$, and suppose for a contradiction that $a,b\notin I$. 
Then $\sqrt{[I,a]}$ and $\sqrt{[I,b]}$ are radical differential ideals which properly contain $I$;
thus we have $s,t\in S$ with $s\in \sqrt{[I,a]}$  and $t\in \sqrt{[I,b]}$. 
Then by Corollary~\ref{cor:radical mult},
$$st \in S\cap \sqrt{[I,a]}\sqrt{[I,b]}\ \subseteq\ S\cap I\ =\ \emptyset,$$
a contradiction.
\end{proof}

\begin{proof}[Proof of Proposition~\ref{prop:radical is intersection of primes, differential}]
Let $I$ be a radical differential ideal of $R$ and $r\in R\setminus I$. Take a radical differential ideal $J$ of $R$ with $J\supseteq I$, disjoint from $S:=\{1,r,r^2,\dots\}$, and maximal with these properties; then $J$ is prime.
\end{proof}

\begin{cor}\label{cor:max diff => prime}
Every maximal differential ideal of $R$ is prime.
\end{cor}

\begin{cor}\label{amalgdifferfld} Let $K$ be a differential field, and let $E$ and 
$F$ be differential field extensions of $K$. Then there is a differential
field $L$ with differential field embeddings $E\to L$ and $F\to L$ that agree on $K$.
\end{cor}
\begin{proof} Corollary~\ref{cor:der on tensor prod, 3} makes the commutative
$K$-algebra $E\otimes_K F$ into a differential ring $R$ such that
the maps $a\mapsto a\otimes 1\colon E\to R$ and
$b\mapsto 1\otimes b\colon F \to R$ are (injective) differential ring morphisms. Take a prime differential ideal $\mathfrak{p}$
of $R$. Composing these maps $E\to R$ and $F\to R$
with the natural map $R \to R/\mathfrak{p}$ and extending~$R/\mathfrak{p}$ to its differential fraction field $L$
yields the desired result.  
\end{proof}

\noindent
Let $K$ be a differential field and $R=K\{Y_1,\dots,Y_n\}$, where $Y_1,\dots,Y_n$ are distinct differential indeterminates over~$K$. We have a 
Differential Nullstellensatz:

\index{theorem!Differential Nullstellensatz}

\begin{cor}\label{cor:diff Nullstellensatz}
Let $P_1,\dots,P_m,Q\in R$. Then the following are equivalent:
\begin{enumerate}
\item[\textup{(i)}] $Q\notin \sqrt{[P_1,\dots,P_m]}$;
\item[\textup{(ii)}] there exists a differential field extension $L$ of $K$ and $\vec{y}\in L^n$
such that $P_i(\vec{y})=0$ for $i=1,\dots,m$ and $Q(\vec{y})\neq 0$.
\end{enumerate}
\end{cor}
\begin{proof}
Suppose (i) holds. Then we have a prime differential ideal $I$ of $R$ with 
$I\supseteq \sqrt{[P_1,\dots,P_m]}$ and $Q\notin I$. Identifying $K$  in the 
obvious way with a differential subfield of the differential integral domain 
$R/I$, we get (ii) with $L$
the differential fraction field of $R/I$ and $\vec{y}=(y_1,\dots,y_n)\in L^n$ 
where $y_j=Y_j+I\in R/I$ for $j=1,\dots,n$. The direction 
(ii)~$\Rightarrow$~(i) is obvious.
\end{proof}

\begin{remark} 
If $P_1,\dots,P_m$ in the previous corollary are homogeneous of degree~$1$, 
then $\sqrt{[P_1,\dots,P_m]}=[P_1,\dots,P_m]$,
since as an ideal of $R$, the latter 
is generated by homogeneous differential polynomials of degree $1$, and 
thus prime by Lemma~\ref{spprid}.
\end{remark} 

\subsection*{Simple differential rings}
We say that $R$ is {\bf simple} if 
the only differential ideals of $R$ are $\{0\}$ and $R$, that is, $\{0\}$ is a maximal differential ideal of $R$.
Simple differential rings are integral domains, by Corollary~\ref{cor:max diff => prime}. If $c$ is a constant of $R$, then $cR$ is a
differential ideal of $R$. It follows that the ring of constants $C_R$ of a simple differential ring $R$ is a field. 

\index{differential ring!simple}
\index{simple!differential ring}

\begin{lemma}\label{lem:integration simple}
Suppose $R$ is simple and $r\in R\setminus \der(R)$. Let $x$ be an 
element of a differential ring extension of $R$ with $x'=r$. Then the 
differential ring $R[x]$ is simple and $x$ is transcendental over $R$.
\end{lemma}
\begin{proof}
We equip the polynomial ring $R[X]$ over $R$ with the  derivation extending the derivation of $R$ such that $X'=r$.
We claim that $R[X]$ is simple. Let $I$ be a nonzero differential ideal of $R[X]$, and let $n$ be the smallest degree of a nonzero element of $I$.
Let $J$ be the 
set of all $a\in R$ such that there is some $P\in I$ with $P=aX^n+\text{terms of degree $<n$}$.
It is easy to see that $J$ is a nonzero differential ideal of $R$, hence~$1\in J$.
Let  $P=X^n+a_{n-1}X^{n-1}+\cdots+a_0\in I$ with $a_0,\dots,a_{n-1}\in R$.
Suppose~$n\geq 1$. We have $P'\in I$ and 
$$P'=(n r+a_{n-1}')X^{n-1}+\text{terms of degree $<n-1$,}$$ 
hence $P'=0$ and
thus $n r+a_{n-1}'=0$, so $y:=-\frac{1}{n}a_{n-1}\in R$ satisfies $y'=r$, a contradiction. Hence $n=0$, so $P=1\in I$.
Thus $R[X]$ is simple. The ring morphism $\phi\colon R[X]\to R[x]$ over $R$  
sending $X$ to $x$ is a surjective morphism of differential rings, so $\phi$ is an isomorphism.
\end{proof}

\begin{lemma}\label{lem:exp integration simple}
Suppose $R$ is simple, $r\in R$, and $y'\ne mry$ for all $y\in R^{\neq}$ and~${m\ge 1}$. Let $x$ be a unit of a differential ring extension $A$ of $R$ with $x'=rx$. Then the differential subring $R[x,x^{-1}]$ of $A$ is simple, and $x$ is transcendental over~$R$.
\end{lemma}
\begin{proof}
Let $X$ be an indeterminate over $R$ and equip $R[X,X^{-1}]$ with the derivation extending the derivation of $R$ such that $X'=rX$. 
We claim that $R[X, X^{-1}]$ is simple. (As in the proof of Lemma~\ref{lem:integration simple}, it follows that then $x$ is transcendental over~$R$ and 
$R[x,x^{-1}]$ is simple.) For this, let $I\neq \{0\}$ be a differential ideal of~$R[X,X^{-1}]$; then $I\cap R[X]\neq\{0\}$. Let $n$ be the smallest degree of a nonzero element of $I\cap R[X]$. Let $J$ be the  set of all $a\in R$ such that there is some $P\in I\cap R[X]$ with $P=aX^n+\text{terms of degree $< n$}$. Then $J$ is a nonzero differential ideal of $R$, hence $1\in J$. Let  $P=X^n+a_{n-1}X^{n-1}+\cdots+a_0\in I$ with $a_0,\dots,a_{n-1}\in R$. If~$n=0$,
then $P=1\in I$. Suppose $n\geq 1$. Then  
$$P'\ =\ nrX^n + \big(a_{n-1}'+(n-1)a_{n-1}r\big)X^{n-1}+\cdots+(a_1'+a_1r)X+a_0' \in I,$$ 
hence $nrP-P'\in I$ has degree less than $n$, so $nrP=P'$, hence $a_{i}'=(n-i) r a_{i}$ for $i=0,\dots,n-1$. By hypothesis we obtain $a_0=\cdots=a_{n-1}=0$, that is, $P=X^n$. Hence $1=X^{-n}P\in I$.
\end{proof}

\noindent
For $R$,~$r$,~$x$ as in Lemma~\ref{lem:exp integration simple} the differential subring $R[x]$ of $R[x,x^{-1}]$ is not simple, since it has $xR[x]$ as a differential ideal. 

\medskip\noindent
{\em In the rest of this subsection $K$ is a differential field with constant field $C$}.

\begin{cor}\label{cor:unit alg}
Suppose $R$ is a simple differential ring extension of $K$.  
Then every constant element of the differential fraction field of $R$ lies in $R$. If $R$ is finitely generated as a $K$-algebra, then every such element is algebraic over $C$. 
\end{cor}
\begin{proof}
Let $a$ be a constant element of the differential fraction field of~$R$. Then $I:=\{b\in R: ab\in R\}$ is a nonzero differential ideal of $R$, so $I=R$ and hence $a\in R$.
Thus for every $c\in C$ with $a\neq c$, $(a-c)R$ is a nonzero differential ideal of~$R$, so
$a-c\in R^\times$. Assume also that $R$ is finitely generated as a $K$-algebra. Then by Lemma~\ref{lem:unit alg}, $a$ is algebraic over $K$, so $a$ is algebraic over $C$ by 
Lemma~\ref{lem:const field of algext}. 
%Let
%$P(X)\in K[X]$ be the minimum polynomial of $a$ over $K$. Differentiating the
%equality $P(a)=0$ shows $P(X)\in C[X]$. 
\end{proof}

\begin{cor}\label{cor:no new exps under integration}
Let $a\in K\setminus \der(K)$. 
Let $y$ be an element of a differential field extension of $K$ with $y'=a$. 
Then every $r\in K(y)^\times$ with $r^\dagger\in K$ lies in $K$.
\end{cor}
\begin{proof}
By Lemma~\ref{lem:integration simple}, $y$ is transcendental over $K$
and $K[y]$ is simple.
Let $r\in K(y)^\times$ be such that $s:=r^\dagger\in K$. Take
coprime $f,g\in K[y]^{\neq}$ and $k\in\Z$ such that $r=(f/g)y^k$ and $y$ does not divide $fg$ in $K[y]$. To show $r\in K$ we can multiply~$r$ with an element of $K^\times$ and arrange that
$f$,~$g$ are monic in $y$. Then $f'\in K[y]$ is of lower degree in $y$
than $f$, and $g'\in K[y]$ is of lower degree in $y$ than~$g$, and
$$r^\dagger\ =\ (f'/f) - (g'/g)+(ka/y)\ =\ s.$$ 
Multiplying both sides by $fgy$ yields
$$ (f'g-fg')y+fg\cdot ka = fg \cdot y s  .$$
Hence $fg|(f'g-fg')y$, so
$f|f'gy$ and $g|fg'y$, and thus $f|f'$ and $g|g'$ in $K[y]$.
Hence $f'=g'=0$, which by Corollary~\ref{cor:unit alg} yields
$f,g\in C\subseteq K$, so $f=g=1$.   
Then $k(a/y)=s\in K$, which in view of $a\ne 0$ gives $k=0$. 
Therefore $r=1\in K$.
\end{proof}

\begin{cor}\label{cor:solution with alg constant field ext}
Let $P\in K\{Y\}^{\neq}$ be irreducible of order $r$, and let
$Q\in K\{Y\}^{\neq}$ have order~$<r$. Then there is an element $y$ of a differential field extension of~$K$ such that $P(y)=0$, $S_P(y)\neq 0$, $Q(y)\neq 0$,
and $C_{K\<y\>}$ is algebraic over $C$.   
\end{cor}
\begin{proof}
Lemma~\ref{lem:construct diff alg extension} gives an element $a$ in a differential field extension of $K$ with minimal annihilator $P$ over $K$. 
The differential subring $R:=K\big\{a, S(a)^{-1}, Q(a)^{-1}\big\}$ of $K\<a\>$ with $S:= S_P$ satisfies 
$$ R\ =\ K\big[a,a',\dots,a^{(r)},S(a)^{-1},Q(a)^{-1}\big].$$
Let $M$ be a maximal differential ideal of $R$. Then $\overline{R}:=R/M$ 
is a simple differential ring, finitely generated as a $K$-algebra. 
For $y:=a+M\in\overline{R}$ we get $P(y)=0$, $S(y)\neq 0$, $Q(y)\neq 0$
in $K\<y\>$, and  $C_{K\<y\>}$ is algebraic over $C$ by
Corollary~\ref{cor:unit alg}.
\end{proof}

\noindent
We use this to derive a technical fact needed in the next volume: %Section~\ref{sec:ls closure}:

\begin{lemma}\label{prop:condition to be min ann}
Let $P\in K\{Y\}^{\neq}$ be irreducible of order $r$, and let
$Q\in K\{Y\}^{\neq}$ have order $<r$. Suppose $P$ is a minimal annihilator over $K$ of every element $y$ of every differential field extension $L$ of $K$
such that $P(y)=0$, $Q(y)\neq 0$, and $C_{L}$ is algebraic over $C$.
Let $a$ in a differential field extension of $K$ satisfy $P(a)=0$, $Q(a)\neq 0$, and $S_P(a)\neq 0$. Then $P$ is a minimal annihilator of $a$ over $K$.  
\end{lemma}
\begin{proof} Put $S:=S_P$.
Corollary~\ref{cor:solution with alg constant field ext}
and its proof give an element $y$ in a differential field extension
of $K$ such that $P(y)=0$, $S(y)\ne 0$, $Q(y)\ne 0$, $C_{K\<y\>}$ is
algebraic over $C$, and the differential subring 
$R:=K\{y, S(y)^{-1}, Q(y)^{-1}\}$ of~$K\<y\>$ is simple. Then $P$ is a minimal annihilator of $y$ by the hypothesis of the lemma. 
The remarks after Lemma~\ref{lem:uniqueness of diff alg extension} give
a surjective differential ring morphism 
$R\to K\big\{a, S(a)^{-1},Q(a)^{-1}\big\}$
over $K$ sending $y$ to $a$. Since $R$ is simple, this morphism is an 
isomorphism.
\end{proof}

\subsection*{Linear disjointness of constant fields}
Let $K$,~$L$ be subfields of a field $F$, and let $C$ be a subfield of both $K$ and $L$. We say that {\bf $K$ is 
linearly disjoint from~$L$ over~$C$} if every finite set of elements of $K$ which is $C$-linearly independent
is also $L$-linearly independent. Equivalently, $K$ is linearly disjoint from $L$ over $C$ iff
the ring morphism $K\otimes_C L\to F$ with $a\otimes b\mapsto ab$ is injective.
In particular, if $K$ is linearly disjoint from $L$ over $C$, then $L$ is linearly disjoint from $K$ over $C$. For proofs of this and other facts on linear disjointness, see \cite[VIII,~\S{}3]{Lang}.

{\em In the rest of this subsection, $K$ is a differential field, and $L$ is a differential field extension of~$K$}.
The following lemma shows that $K$ is linearly disjoint from the constant field $C_L$ of $L$ over the constant field $C$ of $K$.

\index{linearly!disjoint}

\begin{lemma}\label{lem:lin disjoint over constants}
Let $R$ be a differential ring extension of $K$. The ring morphism 
$$K\otimes_C C_R\to R\qquad\text{ with $\ a\otimes c\mapsto ac\ $ for $a\in K$, $c\in C_R$,}$$ 
is injective. $($Its image is the differential subring $K[C_R]$ of $R$.$)$
\end{lemma}
\begin{proof}
Suppose not. Then we have $c_1,\dots,c_n\in C_R$ ($n\geq 1$), linearly independent over $C$, and $a_1,\dots,a_n\in K$,
not all zero, such that $\sum_{i=1}^n a_ic_i=0$. Let $n$ be minimal. Then $a_1\in K^\times$, so $a_1c_1\ne 0$, hence $n\geq 2$. Arranging
$a_1=1$ and differentiating gives $\sum_{i=2}^n a_i'c_i=0$, so
$a_i\in C$ for $i=1,\dots,n$ by the minimality of~$n$, but this 
contradicts the $C$-linear independence of $c_1,\dots,c_n$.
\end{proof}

\begin{cor}\label{cor:lin disjoint over constants, 1}
Let $A\in K\{Y\}$ be homogeneous of degree $1$ and let $g\in K$. If $A(f)=g$ for some $f\in K[C_L]$, then 
$A(f)=g$ for some $f\in K$.
\end{cor}
\begin{proof}
Let $B$ be a $C$-vector space basis of $C_L$ with $1\in B$, let $b$ range over~$B$, and assume $f\in K[C_L]$. By Lemma~\ref{lem:lin disjoint over constants} we have a
unique admissible family $(f_b)$ in $K$ such that 
$f=\sum_b f_b b$. If $f\in K[C_L]$ and $A(f)=g$, then $\sum_b A(f_b)b=g$ and hence $A(f_1)=g$.
\end{proof}

\noindent
In the same way we obtain:

\begin{cor}\label{cor:lin disjoint over constants, 2}
Let $D$ be a subring of $C_L$ containing $C$. Then $K[D]$ is a 
differential subring of $L$ with $D$ as its ring of constants. 
\end{cor}

\begin{cor}\label{cor:lin disjoint over constants, 3}
Let $R$ be a differential ring extension of $K$ and $D$ a subring of~$C_R$ containing $C$ such that $R=K[D]$.
Let $J$ be a $K$-linear subspace of $R$ such that $\der J\subseteq J$, and
set $I:=J\cap D$. Then $J=IK$. If $D$ is a field, then $R$ is simple.
\end{cor}
\begin{proof}
Let $B$ be a $C$-vector space basis of $D$, and let $b$ range over $B$. For any $f\in R$ there is  
by Lemma~\ref{lem:lin disjoint over constants} a unique admissible family $(f_b)$ in $K$ such that $f=\sum_b f_b b$, and define the length $\ell(f)$ of $f$ by $\ell(f):=|\{b:\ f_b\neq 0\}|$.
By induction on $\ell(f)$ we show that if $f\in J$, then $f\in IK$. This is clear if $\ell(f)=0$, so
suppose $f\in J$ satisfies $\ell(f)\geq 1$. Take $b_1\in B$ with $f_{b_1}\neq 0$;
we may assume $f_{b_1}=1$. We may also assume $f\notin D$ and thus take $b_2\in B$ with $f_{b_2}\notin C$. Now $f'=\sum_b f_b' b$ and so $\ell(f')<\ell(f)$, hence $f'\in IK$.
Setting $a:=f_{b_2}$, we similarly have $(a^{-1}f)'\in IK$ and thus $(a^{-1})'f=(a^{-1}f)'-a^{-1}f'\in IK$.
Since $a^{-1}\notin C$, we obtain $f\in IK$. This proves $J=IK$. If $D$ is a field, it follows that $R$ is simple as a differential ring. 
\end{proof}

\noindent
Thus with $R$ and $D$ as in
 Corollary~\ref{cor:lin disjoint over constants, 3} the map $I\mapsto IR$ is a bijection from the set of ideals of $D$ onto the set of differential ideals
of $R$, with inverse $J\mapsto J\cap D$.

From Corollaries~\ref{cor:lin disjoint over constants, 2}, 
~\ref{cor:lin disjoint over constants, 3}, and
~\ref{cor:unit alg} we obtain:

\begin{cor}\label{cor:lin disjoint over constants, 4}
Let $D$ be a subfield of $C_L$ containing $C$. Then $K(D)$ is a 
differential subfield of $L$ with $D$ as its field of constants. 
\end{cor}

\subsection*{Adjoining constants}
Let $K$ be a differential field and $L$ an extension field of $K$ with
subfield $D\supseteq C$ such that $L=K(D)$. If the field $L$ has a derivation extending that of $K$ with constant field~$D$, then by the previous subsection $K$ and $D$ are linearly disjoint over 
$C$. In Section~\ref{sec:H-fields} we need a converse of this fact:

\begin{lemma}\label{lem:constfieldext}
Suppose $K$ and $D$ are linearly disjoint over $C$. Then
there is a unique derivation on the field $L$ that extends that of $K$ and is trivial on
$D$; this derivation has~$D$ as its constant field.
\end{lemma}
\begin{proof} There is clearly at most one such derivation. For existence, let $B$ be a transcendence
basis of $D$ over $C$. Linear disjointness gives algebraic independence of~$B$ over $K$; see \cite[VIII, Proposition~3.2]{Lang}. So the derivation of $K$ extends to a derivation of the field $K(B)$ that annuls each
element of $B$ (by Corollary~\ref{lem:derivative on trans ext}), hence to a derivation of the algebraic extension~$L$ of $K(B)$,
and this derivation will annul each element of $D$, since $D$ is algebraic over~$C(B)$ (by Corollary~\ref{cor:characterization of sep alg}).
By Corollary~\ref{cor:lin disjoint over constants, 4} the constant field of this derivation on $L$ equals $D$.
\end{proof}

\subsection*{Notes and comments} 
Proposition~\ref{prop:radical is intersection of primes, differential} is 
due to Ritt; the proof as given is from~\cite{Kaplansky76}. The equivalence of 
(i) and (ii) in Corollary~\ref{cor:diff Nullstellensatz} is in
Raudenbush~\cite{Raudenbush}, with a constructive proof in 
Seidenberg~\cite{Seidenberg}. 
Corollary~\ref{cor:unit alg} is in Levelt~\cite[Lem\-ma~A3]{Levelt}.
Lemma~\ref{lem:lin disjoint over constants} is from \cite[p.~292]{Rosenlicht1}.
For Corollary~\ref{cor:lin disjoint over constants, 2} see \cite[Corollary~5 on p.~768]{Kolchin-GTDF},
and for Corollary~\ref{cor:lin disjoint over constants, 3}, see \cite[Proposition~1]{Levelt}.

%% file: mt-4-7.tex
\section{Differentially Closed Fields}
\label{sec:DCF}

\noindent
This section presents L.~Blum's proof of A.~Robinson's theorem 
that the theory of differential fields has a model completion, namely the theory of \textit{differentially closed fields}\/; see \ref{sec:mc} and \ref{sec:qe} for the concept of model completion.  We shall meet differentially closed fields briefly in Chapter~\ref{ch:lindifpol} and Section~\ref{Miscellaneous Facts about Asymptotic Fields}. We also 
construct \textit{differential closures\/} 
of differential fields, but do not use these later.

\medskip\noindent
The language of differential rings is the language  $\mathcal L_\der=\{0,1,{-},{+},{\,\cdot\,},\der\}$
obtained by augmenting the language $\{0,1,{-},{+},{\,\cdot\,}\}$ of rings by an extra unary function symbol~$\der$.
We view differential rings as structures for this language in the obvious way. A differential ring morphism is a morphism with respect to this language. Likewise, a differential subring of a differential ring $R$ is a substructure of $R$ with respect to the language $\mathcal L_\der$. 

Let $K$ be a differential field. Note that $\Z$ with its trivial derivation is a substructure of $K$ that is not a differential subring of $K$, in view of our convention that differential rings contain $\Q$. In fact,
the substructures of $K$ are exactly the
subrings of $K$ closed under the derivation
of $K$ and equipped with the restriction of this
derivation to the subring. 

\begin{definition} \label{p:differentially closed}
A differential field $K$ is said to be {\bf differentially closed}
if for all $P\in K\big[Y,\dots, Y^{(r)}\big]{}^{\ne}$ and
$Q\in K\big[Y,\dots, Y^{(r-1)}\big]{}^{\ne}$ such that $Y^{(r)}$ occurs in $P$
there is $y\in K$ with $P(y)=0$ and $Q(y)\neq 0$. (In this definition, we may restrict to irreducible $P$.)
Taking $P\in K[Y]$ we see that if $K$ is differentially closed, 
then $K$ is algebraically closed.
\end{definition}

\index{differential field!differentially closed}
\index{closed!differentially}
\index{differentially!closed}

\noindent
Applying Corollary~\ref{cor:solution with alg constant field ext} 
iteratively, we see that every 
differential field $K$ has a differential field extension $L$ such that
$L$ is differentially closed and $\d$-algebraic over~$K$, and
$C_L$ is algebraic over $C$.

\medskip
\noindent
Let  $\operatorname{DCF}$ be the theory of differentially closed fields, formulated in the  language $\mathcal L_\der$.

\nomenclature[Bi]{$\operatorname{DCF}$}{theory of differentially closed fields in the language $\mathcal L_\der$}
\nomenclature[Be]{$\mathcal L_\der$}{language of differential rings}

\begin{theorem}\label{thm:DCF-qe}
$\operatorname{DCF}$ admits quantifier elimination.
\end{theorem}
\begin{proof}
Let $E$ and $F$ be differentially closed fields such that $F$ is $\abs{E}^+$-saturated, and let $R$ be a proper substructure of $E$ and $\phi \colon R \to F$ be
an embedding; we shall ex\-tend~$\phi$ to an embedding of a differential subring of $E$, properly containing~$R$, into~$F$. If $R$ is not a field, then we can extend $\phi$ to the fraction field of $R$ inside~$E$ (which is also a differential subfield of $E$). So assume $R$ is a differential
field~$K$. Take any $a \in E \setminus K$. 
Consider first the case that $a$ is $\d$-transcendental over $K$. Using that $F$ is  $\abs{E}^+$-saturated with nontrivial
derivation
we get a $b \in F$ that is $\d$-transcendental over the differential
subfield $\phi(K)$ of $F$. Then $\phi$ extends to an embedding
$K\{a\} \to F$ sending $a$ to $b$.
Next assume that $a$ is $\d$-algebraic
over $K$, and let $P\in K\{Y\}$, of order $r\in\N$, be a minimal annihilator of $a$ over $K$, so
$$P(a) = 0,\qquad \text{$Q(a) \neq 0$ for all $Q \in K\{Y\}^{\neq}$ of order $< r$.}$$
Using that $F$ is  $\abs{E}^+$-saturated and differentially closed 
we can take $b \in F$ such that this equation and these inequations hold with $b$ instead of $a$ and with $P$ and the $Q$'s replaced by their
$\phi$-images. Then by Lemma~\ref{lem:uniqueness of diff alg extension} we can extend $\phi$ to an embedding of $K\{a\}$
into $F$ that sends $a$ to $b$.
\end{proof}

\begin{cor}\label{cor:DCF-qe}
$\operatorname{DCF}$ is complete, and $\operatorname{DCF}$ is the model completion of the theory of
differential fields.
\end{cor}
\begin{proof}
The field $\Q$ with the trivial derivation embeds into every differentially closed field, hence the first part follows
from Theorem~\ref{thm:DCF-qe}. The second statement follows from that theorem in view of the fact that every differential field extends to a differentially closed field.
\end{proof}

\noindent
As a consequence, we can now add to the Differential 
Nullstellensatz~\ref{cor:diff Nullstellensatz} two conditions, 
each equivalent to condition (i) in that theorem:

\index{theorem!Differential Nullstellensatz}

\begin{enumerate} 
\item[(iii)] for some differentially closed differential field extension $L$ of $K$ and some 
$\vec{y}\in L^n$, we have
$P_i(\vec{y})=0$ for $i=1,\dots,m$ and $Q(\vec{y})\neq 0$;
\item[(iv)] for every differentially closed differential field extension $L$ of $K$, there is  
$\vec{y}\in L^n$
such that $P_i(\vec{y})=0$ for $i=1,\dots,m$ and $Q(\vec{y})\neq 0$.
\end{enumerate}

\begin{cor}\label{cor:induced structure on C}
Suppose $K$ is a differentially closed field. If $X\subseteq K^n$ is definable in $K$, then 
$X\cap C^n$ is definable in the algebraically closed field $C$.
\end{cor}
\begin{proof}
Let $P=P(Y_1,\dots, Y_n)\in K\{Y_1,\dots, Y_n\}$ be a differential polynomial in the distinct indeterminates $Y_1,\dots, Y_n$ over $K$. Upon removing from~$P$ the monomials involving any $Y_i^{(r)}$ with $r\ge 1$ we obtain an ordinary polynomial $p\in K[Y_1,\dots, Y_n]$ such that for all
$y_1,\dots, y_n\in C$,
$$P(y_1,\dots,y_n)\ =\ p(y_1,\dots,y_n).$$
By QE for $\operatorname{DCF}$ (Theorem~\ref{thm:DCF-qe}), this fact reduces our job to showing that
for any polynomial $p\in K[Y_1,\dots, Y_n]$ the subset
$\big\{ c\in C^n: p(c)=0\big\}$ of $C^n$ is definable in 
the algebraically closed field $C$;
for this, use the following Lemma~\ref{lem:induced Zariski top}.
\end{proof}

\begin{lemma}\label{lem:induced Zariski top}
Let $E\subseteq F$ be a field extension, $Y_1,\dots,Y_n$ distinct indeterminates, and $p\in F[Y_1,\dots,Y_n]$. Then there are polynomials
$p_1,\dots,p_m\in E[Y_1,\dots,Y_n]$ such that for all $a\in E^n$,
$$p(a)=0 \quad \Longleftrightarrow\quad p_1(a)=\cdots=p_m(a)=0.$$
\end{lemma}
\begin{proof}
Take a basis $b_1,\dots,b_m$ of the $E$-linear subspace of $F$ generated by
the coefficients of $P$. Then $p=b_1p_1+\cdots+b_mp_m$ with $p_1,\dots,p_m\in E[Y_1,\dots,Y_n]$,
and then $p_1,\dots,p_m$ have the desired property.
\end{proof}

\subsection*{Differential closures} This subsection is not used later. 

\begin{lemma}\label{lem:construct diff closure}
Let $K$ be a differential field, let $P\in K\{Y\}^{\neq}$ be irreducible of order~$r$, and let $Q\in K\{Y\}^{\neq}$ have
order~$<r$. Then there is an element $y$ of a differential field extension of $K$ with
$P(y)=0$, $Q(y)\neq 0$, such that $K\<y\>$ embeds over $K$
into any differentially closed differential field extension of $K$.
\end{lemma}
\begin{proof}
Let $S=S_P$, and let $a$ and $R=K\big\{a,S(a)^{-1},Q(a)^{-1}\big\}$ be as in the proof of
Corollary~\ref{cor:solution with alg constant field ext}.
Let $M$ be a maximal differential ideal of $R$. Then
$\overline{R}:=R/M$ is a simple differential ring, hence an integral domain; we identify
$K$ with its image in $\overline{R}$ via the natural embedding $K\to\overline{R}$.
With $y:=a+M\in\overline{R}$ we have $P(y)=0$, $Q(y)\neq 0$.
Since the ring $R$ is noetherian,
we can take $P_1,\dots,P_m\in K\{Y\}$ of order~$\leq r$ such that $P_1(a),\dots,P_m(a)$ generate the ideal~$M$ in $R$.
Let $L$ be a differentially closed differential field extension of~$K$.
By the equivalence of~(iii) and~(iv) above, we can take $z\in L$ with $P(z)=0$, $P_i(z)=0$ for $i=1,\dots,m$, $S(z)\neq 0$, $Q(z)\neq 0$. The remarks after Lemma~\ref{lem:uniqueness of diff alg extension} give
a differential ring morphism $R\to L$ over $K$ sending $a$ to $z$.
Now $M$ is part of the kernel of this differential ring morphism,
so we get a differential ring
morphism $\overline{R}\to L$ over $K$, which is an embedding since $\overline{R}$ is simple.
\end{proof}

\noindent
Let $K$ be a differential field. A {\bf differential closure} of $K$ is a differential 
field extension~$K^{\operatorname{dc}}$ of $K$
which is differentially closed and
such that every embedding of~$K$ into a differentially closed field $L$ extends to an embedding of $K^{\operatorname{dc}}$ into $L$.
Every differential field has a differential closure: this follows by  a straightforward
transfinite construction using Lemma~\ref{lem:construct diff closure}.
For the sake of completeness we also mention the following results
without proof.

\index{differential field!differential closure}
\index{closure!differential}
\nomenclature[M]{$K^{\operatorname{dc}}$}{differential closure of the differential field $K$}

\begin{theorem}\label{thm:diff closure}
Any two differential closures of $K$ are isomorphic over $K$.
\end{theorem}

\noindent
Thanks to this theorem, we may refer to \textit{the}\/ differential closure
$K^{\operatorname{dc}}$ of $K$. Note that~$K^{\operatorname{dc}}$ is 
$\d$-algebraic over $K$, and that the constant field of $K^{\operatorname{dc}}$ is the 
algebraic closure of the constant field $C$ of $K$ inside~$K^{\operatorname{dc}}$. 

In contrast to the corresponding notions for fields,
differential fields always have proper $\d$-algebraic 
extensions, and the differential closure of a differential 
field~$K$ is not always minimal over $K$:

\begin{prop}\label{prop:diff closure not minimal}
If the derivation of $K$ is trivial, then $K^{\operatorname{dc}}$ properly contains
a differentially closed differential subfield containing $K$.
\end{prop}

\subsection*{Notes and comments} 
Most of the material above stems from Robinson~\cite{Robinson59} 
and Blum~\cite{Blum}. Theorem~\ref{thm:diff closure} is 
in Shelah~\cite{Shelah72a}. Proposition~\ref{prop:diff closure not minimal} is due, independently, to Kolchin~\cite{Kolchin74}, Rosenlicht~\cite{Rosenlicht74}, and Shelah~\cite{Shelah73}.
See \cite{MTF} for an exposition of \ref{thm:diff closure} and \ref{prop:diff closure not minimal}. 

 Singer~\cite{Singer78-1} shows that the theory of
{\em ordered\/} differential fields has a model completion.   
(No relation between ordering and derivation is imposed, and the relevant language is that of differential rings with a binary relation symbol~$\leq$ for the ordering.) Given a model $K$ of this model completion,
Singer~\cite{Singer78-2} also showed that the differential field~$K[\imag]$, where $\imag^2=-1$, is
differentially closed (and hence~$K[\imag]$ is the differential closure of $K$). This fact together with  Corollary~\ref{cor:induced structure on C} shows that such $K$, as a subset of~$K[\imag]$, is not definable in the differential field $K[\imag]$. By Proposition~\ref{prop:defining T} below this
 is in contrast to what happens for $K=\T$.

By Michaux~\cite{Michaux}, the theory
of valued differential fields (no relation between derivation and valuation being imposed)
has a model completion. Similar results for certain theories of topological differential fields are in
Guzy and Point~\cite{GuzyPoint10}.

%% file: mt-5.tex
\chapter{Linear Differential Polynomials}\label{ch:lindifpol} 

\noindent
Linear differential polynomials and their zero sets play a special role in our work. 
Their theory is of course much better understood than for
differential polynomials in general. For us this fact is particularly
relevant since the property of a valued differential field being differential-henselian (studied in Chapter~\ref{sec:dh1}) involves the linear part of an arbitrary differential polynomial in an essential way. Also, the key operation of compositional conjugation
on arbitrary differential polynomials is defined by transformations
of linear differential polynomials; see Section~\ref{Compositional Conjugation}.

Accordingly we consider in this chapter homogeneous linear differential polynomials and the corresponding linear operators, and prove various basic results on them as needed later. In particular,
we study the property of a linear differential operator over a differential field $K$ of defining a surjective map $K \to K$, and the transformation of a system of linear differential equations in several unknowns to an equivalent system of several linear differential equations in a single unknown. 

In the final section of this chapter we apply this material on linear differential polynomials (plus some commutative algebra) to prove a result on zeros of systems of
non-linear algebraic differential equations.

\input{mt-5-1}

\input{mt-5-2}

\input{mt-5-3}

\input{mt-5-4}

\input{mt-5-5}

\input{mt-5-6}

\input{mt-5-7}

\input{mt-5-8}

\input{mt-5-9}

%% file: mt-5-1.tex
\section{Linear Differential Operators}
\label{Linear Differential Operators}

\noindent
This section contains definitions and basic results 
about linear differential operators. \textit{Throughout
$K$ is a differential ring with derivation $\der$ and ring of constants $C$.
We let~$a$,~$b$, sometimes with subscripts, range over $K$.}

\index{linear differential operator}
\index{operator!linear differential}
\nomenclature[Q]{$K[\der]$}{ring of linear differential operators over $K$}

\subsection*{The ring $K[\der]$}
A homogeneous linear differential polynomial 
$$A(Y)\ =\ a_0Y + a_1Y' + \dots + a_nY^{(n)}\in K\{Y\} $$ defines a $C$-linear 
operator
$y \mapsto A(y)=a_0y + a_1y' + \dots + a_ny^{(n)}$ on $K$ and
the composition of two such operators is again an operator of this form.
In its role as an operator the above $A$ is more conveniently written as
$a_0 + a_1\der + \dots + a_n\der^n$. Formally,~$K[\der]$ is 
a ring that contains $K$ as a
subring and with a distinguished element, 
denoted~$\der$, such that 
$K[\der]$, as a left-module over $K$, is free with basis 
$$ \der^0, \der^1, \der^2, \der^3,\dots, \qquad 
\text{ with $\der^0=1$, $\der^1=\der$,
$\der^m \ne \der^n$ whenever $m\ne n$,}$$
and such that $\der a=a\der + a'$. (So $K[\der]$ is commutative 
iff $a'=0$ for all $a$.) This description determines $K[\der]$ up to 
isomorphism over $K$ as a ring extension of $K$ with a distinguished element $\der$. Note the following identity:
\begin{equation}\label{D^na} \der^n a\ =\
\sum_{i=0}^n \binom{n}{i}a^{(n-i)}\der^i\ =\ 
a^{(n)} + na^{(n-1)}\der + \dots + a\der^n.
\end{equation}
The element $\der$ of $K[\der]$ should not be confused with the 
derivation $\der$ of $K$: $\der\ne 0$ in~$K[\der]$ even if the 
derivation $\der$ is trivial on $K$. Every $A\in K[\der]$ has the form
$$%\begin{equation}\label{DiffOp}
A\ =\ a_0 + a_1\der + \dots + a_n\der^n,
$$%\end{equation}
and for such $A$ we put 
$$A(y)\ :=\ a_0y + a_1y' + \dots + a_ny^{(n)}$$ for $y$ in a
differential ring extension of $K$. Multiplication of elements of 
$K[\der]$ corresponds to composition:
$ (AB)(y)=A\bigl(B(y)\bigr)$ for $A,B\in K[\der]$ and $y$ in a 
differential ring extension of $K$. The map $K[\der] \to \Hom_C(K,K)$ that associates to each 
$A\in K[\der]$ the $C$-linear operator $y\mapsto A(y)$ on $K$ is itself
$C$-linear. If $K$ is a differential field and $C\ne K$, then this
map $K[\der] \to \Hom_C(K,K)$ is injective, by Lemma~\ref{notthin}. 
Multiplication on the right by an 
element of $K$ 
corresponds to multiplicative conjugation of the corresponding linear 
differential polynomial: 
\begin{align*}\text{if }\quad(a_0 + a_1\der + \dots + a_n\der^n)b\ &=\
b_0 + b_1\der + \dots + b_n \der^n,\\
\text{then }\quad \bigl(a_0Y + a_1Y' + \dots + a_nY^{(n)}\bigr)_{\times b}\ &=\
b_0Y + b_1Y' + \dots + b_nY^{(n)}.
\end{align*}
If $L$ is a differential ring, then any differential ring morphism $K\to L$ extends uniquely to a ring morphism
$K[\der] \to L[\der]$ sending $\der\in K[\der]$ to $\der\in L[\der]$. Thus 
every automorphism $\sigma$ of the differential ring $K$ extends to an automorphism, also denoted by $\sigma$, of the ring $K[\der]$ with $\sigma \der=\der$
(so $\sigma\big(A(y)\big)=(\sigma A)(\sigma y)$ 
for $A\in K[\der]$ and $y\in K$). 

\medskip\noindent
Call $A\in K[\der]$ {\bf monic of order $n$} if it has the form 
$\der^n + a_{n-1}\der^{n-1} + \dots + a_0$. 
If $A,B\in K[\der]$ are monic of order $m$, $n$, then $AB$
is monic of order $m+n$. To $P\in K\{Y\}$ we associate 
its {\bf linear part} $L_P \in K[\der]$, 
$$L_P\ :=\ \sum_n \frac{\partial P}{\partial Y^{(n)}}(0)\der^n \qquad(\text{so }\ L_{P_{+a}}\ =\ \sum_n \frac{\partial P}{\partial Y^{(n)}}(a)\der^n). $$ 
Hence $L_P(y)=P_1(y)$
for all $y$ in all differential ring extensions of $K$, and then for
$P,Q\in K\{Y\}$ we have $L_{PQ}=P(0)L_Q + Q(0)L_P$, 
$ L_{P_{\times a}}=L_P\cdot a$ and $\der L_P= L_{P'}$. 
%Also $L_{P_{+a}}=\sum_n \frac{\partial P}
%{\partial Y^{(n)}}(a)\der^n$.

\index{linear differential operator!monic}
\index{differential polynomial!linear part}
\index{linear part}
\index{linear differential operator!kernel}
\nomenclature[O]{$L_P$}{linear part of $P$}
\nomenclature[Q]{$\ker A$}{kernel of  $A$}

\medskip\noindent
Let $A\in K[\der]$. Then its {\bf kernel}
$$\ker A\ :=\ \{y\in K:\ A(y)=0\}$$ is a $C$-submodule of $K$. If we need to stress the dependence of $\ker A$ on $K$, we write~$\ker_K A$. If~$L$ is a differential ring extension of $K$, then $\ker_K A$ is a $C$-submodule of $\ker_L A$, and $\ker_K A=K\cap\ker_L A$. 

\medskip
\noindent
Let $a$ be a unit of $K$. We define the {\bf twist} of $A$ 
by $a$ to be 
$A_{\ltimes a}:=a^{-1}Aa\in K[\der]$. \index{linear differential operator!twist}\nomenclature[Q]{$A_{\ltimes a}$}{twist of  $A$ by $a$} In particular,
$ \der_{\ltimes a}= \der + a^\dagger$.
We have $\ker A=a\ker A_{\ltimes a}$. Note that if $A$ is monic of order~$n$, then so is $A_{\ltimes a}$. The map $B \mapsto B_{\ltimes a}$ is an 
automorphism of the ring $K[\der]$; it is the identity on $K$, with inverse $B \mapsto B_{\ltimes a^{-1}}$. 

\subsection*{Right-inverses of linear differential operators} 
Suppose now that $\der^{-1}$ is a {\em right-inverse\/} to the derivation $\der$ of $K$, that is, $\der^{-1}\colon K \to K$ satisfies $\der \circ \der^{-1}=\operatorname{id}_K$. (Think of $\der^{-1}$ as
integration.) Then, if $b\in K$ and $b=a^\dagger$ with~$a$ a unit of $K$, we obtain a right-inverse $(\der + b)^{-1}:= a^{-1} \der^{-1} a$ to $\der+b$, since
$$(\der + b)\circ a^{-1}\der^{-1} a\ =\ (a^{-1}\der a)\circ a^{-1}\der^{-1} a\ =\ \operatorname{id}_K.$$ Here $a^{-1} \der^{-1} a\colon K \to K$ sends each $x\in K$ to 
$a^{-1}\der^{-1}(ax)$. Next, if $A\in K[\der]$ is monic of order $n\ge 1$ and factors as
$A=(\der+b_1)\cdots (\der + b_n)$ with each $b_i=a_i^\dagger$ and~$a_i$ a unit of $K$, then the operator
$A$ on $K$ has right-inverse 
$$A^{-1}\ :=\ (a_n^{-1}\der^{-1}a_n) \circ \dots \circ (a_1^{-1}\der^{-1}a_1)$$
in the sense that $A\circ A^{-1}=\operatorname{id}_K$.

\subsection*{The derivative of a linear differential operator} 
For an operator
$$A\ =\ a_0+a_1\der+\cdots+a_n\der^n\in K[\der]$$ we define its {\bf derivative} $A' \in K[\der]$ by
$$A'\ :=\ a_1 + 2a_2\der + \cdots + na_n\der^{n-1} \qquad   
(\text{so $A'=0$ if $n=0$}).$$
The map $A\mapsto A'$ is a derivation on the ring $K[\der]$:

\index{derivative!linear differential operator}
\nomenclature[Q]{$A'$}{derivative of  $A$}

\begin{lemma} Let $A,B\in K[\der]$. Then 
$$(A+B)'\ =\ A'+B', \qquad (AB)'\ =\ A'B+AB'.$$
\end{lemma}
\begin{proof}
The first identity is clear. The second follows from its 
special cases
$$(A\der^j)'\ =\ A'\der^j + jA\der^{j-1}, \qquad (aB)'\ =\ aB', \qquad 
(\der^jb)'\ =\ j\der^{j-1}b,$$
where $j\in\N^{\ge 1}$.
(For the last equality, use \eqref{D^na}.) 
\end{proof}

\noindent
For $A\in K[\der]$ and $i\in \N$ we define $A^{(i)}\in K[\der]$ by 
$$A^{(0)}\ :=\ A, \qquad A^{(i+1)}\ :=\ (A^{(i)})'.$$ 
Induction on $i$ yields that if
$A=\sum_{j=0}^n a_j\der^j$, then
\begin{equation}\label{A^(i)}
\frac{A^{(i)}}{i!}\ =\ \sum_{j=i}^n \binom{j}{i}a_j \der^{j-i},
\end{equation}
and thus for all $f\in K$:
\begin{equation}\label{Af}
Af\ =\ \sum_{i=0}^{n} \frac{A^{(i)}(f)}{i!}\der^i.
\end{equation}
In the next lemma, $k$ is an integer, and
$$\binom{k}{n}\ :=\ \frac{1}{n!} \cdot k(k-1)\cdots (k-n+1)
\qquad \text{(a rational number).}$$
So $\binom{k}{n}=\frac{k!}{(k-n)!\cdot n!}\in\N^{\ge 1}$ if $k\geq n$ and 
$\binom{k}{n}=0$ if $0\leq k<n$.

\begin{lemma}\label{A(fx^k)}
Let $A\in K[\der]$, $f\in K$, and suppose $x\in K^\times$ satisfies $x'=1$. Then
$$A(fx^k)\ =\ \sum_{i\geq 0} \binom{k}{i} A^{(i)}(f)\cdot x^{k-i}.$$
\end{lemma}
\begin{proof} By the identity~\eqref{Af}  we have 
$$A(fx^k)\ =\ (Af)(x^k)\ =\ \sum_{i\ge 0}\frac{A^{(i)}(f)}{i!}(x^k)^{(i)}.$$
It remains to note that for every $i\in\N$, 
\equationqed{(x^k)^{(i)}\ =\ 
k(k-1)\cdots (k-i+1)\cdot x^{k-i}\ =\ i!\cdot\binom{k}{i}x^{k-i}.}
\end{proof}

%\begin{cor}\label{A(fx^k), cor}
%Let $x$ be a unit of a differential ring extension of $K$ with $x'=1$, and 
%let $$A=a_j\der^j+a_{j+1}\der^{j+1}+\cdots+a_n\der^n\quad\text{ with
%$0\leq j\leq n$, $a_j,\dots,a_n\in K$, $a_j\neq 0$.}$$
%If $0\leq k<j$ then $A(x^k)=0$, and if $k<0$ or $k\geq j$ then
%$$A(x^k)=(a+g)\cdot x^{k-j}$$
%for some nonzero $a$ and some $g\in x^{-1}K[x^{-1}]$.
%\end{cor}
%\begin{proof}
%By Lemma~\ref{A(fx^k)}
%$$A(x^k)=\sum_{m=j}^n 
%k(k-1)\cdots (k-m+1)a_m\cdot x^{k-m},$$
%from which the corollary follows immediately.
%\end{proof}

\subsection*{The order of a linear differential operator} \textit{In this subsection we assume that~$K$ is a differential field.}\/ Let $A,B\in K[\der]$. Define the {\bf order} \index{order!linear differential operator} \index{linear differential operator!order} of $A$ by: 
$$\order 0 = -\infty, \qquad
\order A = n \text{ for $A=a_0 + a_1\der + \dots + a_n\der^n$,  
$a_n\ne 0$.}$$
It is easy to check that $\order AB = \order A + \order B$; 
in particular, if $A\ne 0$ and $B\ne 0$, then 
$AB\ne 0$. With $[A,B]:=AB-BA$ it follows easily from  
\eqref{D^na} that $\order\,[A,B]<\order A+\order B$ for
nonzero $A,B$.

\medskip
\noindent
\textit{In the rest of this subsection we assume $K\neq C$.}

\begin{lemma}\label{lem:commuting with scalar}
Let $A\in K[\der]$ and $A\notin K$. Then $[A,b]\ne 0$ for some $b\in K$.
% iff $A\in K$.
\end{lemma}
\begin{proof} We have $A=a_0+a_1\der+\cdots+a_n\der^n$, $n\ge 1$, $a_n\neq 0$.
Take $b\in K$ with $b'\neq 0$. Then $bA = a_nb\der^n + a_{n-1}b\der^{n-1}+\text{terms of order $<n-1$}$, and by \eqref{D^na},
$Ab= a_nb\der^n + (a_n n b' + a_{n-1}b)\der^{n-1}+
\text{terms of order $<n-1$,}$
so $bA \ne Ab$.
\end{proof}

\noindent
In view of $\der a=a'+a\der$ we obtain:

\begin{cor}\label{cor:center of K[der]}
Let $A\in K[\der]$. Then $[A,B]=0$ for all $B\in K[\der]$ iff $A\in C$.
\end{cor}

\noindent
An {\bf ideal} of $K[\der]$ is an additive subgroup $I$ of $K[\der]$ such that
$FA\in I$ and $AF\in I$ whenever $F\in K[\der]$ and $A\in I$.

\index{ideal!ring of differential operators}

\begin{cor}\label{cor:K[der] is simple}
The ring $K[\der]$ has no ideals except $\{0\}$ and $K[\der]$.
\end{cor}
\begin{proof}
Let $I$ be an ideal of $K[\der]$, and suppose $I\neq \{0\}$. Take a nonzero $A\in I$ of minimal order.
Then for all $b\in K$ we have $[A,b]\in I$ and $\order\,[A,b]<\order A$, hence $[A,b]=0$.
Thus $A\in K$ by Lemma~\ref{lem:commuting with scalar} and hence $I=K[\der]$.
\end{proof}

\subsection*{Euclidean division in $K[\der]$} \textit{In this subsection 
$K$ is a differential field.}\/
The ring~$K[\der]$ admits {\em division with remainder on the left}: for 
$A,B\in K[\der]$, $B\neq 0$,
there exist unique $Q,R\in K[\der]$ with $A=QB+R$ and 
$\order R < \order B$.
It follows that for every left ideal $I$ of $K[\der]$ there exists an
$A\in K[\der]$ such that $I=K[\der] A$. (A \textit{left ideal\/} of $K[\der]$
is an additive subgroup $I$ of  $K[\der]$ such that $FA\in I$ whenever
$F\in K[\der]$ and $A\in I$.) Using division with remainder on the left and induction on $\min(\order A, \order B)$ it follows easily that $K[\der]A \cap K[\der]B\ne \{0\}$ for $A,B\in K[\der]^{\ne}$.
Thus for $A_1,\dots,A_r\in K[\der]^{\ne}$, $r\ge 1$, we can define the 
\textit{least common left multiple\/} of $A_1,\dots,A_r$ to be the
unique monic $A\in K[\der]$ such that 
$$K[\der] A\ =\  K[\der] A_1 \cap\cdots\cap K[\der] A_r.$$
Note that if $A,B\in K[\der]$, then $1\notin K[\der]A+ K[\der]B$ iff 
$A$ and $B$ have a nontrivial common right divisor, that is, a
$D\in K[\der]$ of positive order such that there are $A_1, B_1\in K[\der]$
with $A=A_1D$ and $B=B_1D$.  

\index{ideal!ring of differential operators}
\index{linear differential operator!Euclidean division}

\begin{lemma}\label{fa1} Let $A,B\in K[\der]$ and $a$ satisfy
$(\der -a)A\in K[\der]B$ and ${\order(B)\ge 2}$. Then 
$1\notin K[\der]A+ K[\der]B$.
\end{lemma}
\begin{proof} Suppose $1=\alpha A + \beta B$ with $\alpha, \beta\in K[\der]$. 
Then
$\alpha=q\cdot(\der-a)+b$ with $q\in K[\der]$. Also $(\der -a)A=\gamma B$ with
$\gamma\in K[\der]$, hence 
$$1\ =\ bA + (q\gamma + \beta)B.$$ 
Then $b\ne 0$ since 
$\order(B)>0$, and so $b^{-1}(\der-c)b=\der -a$ for
$c=a+b^{\dagger}$. Hence $b^{-1}(\der -c)\in K[\der]B$, contradicting 
$\order(B)\ge 2$.
\end{proof}

\begin{lemma}\label{fa2} Let $A,B,D\in K[\der]$ and $a$ be such that
$$K[\der]A + K[\der]B\ =\ K[\der]D\ \ne\ K[\der](\der -a)A + K[\der]B.$$
Then  $K[\der](\der -a)A + K[\der]B = K[\der](\der-c)D$ for some $c\in K$.
\end{lemma}
\begin{proof} Since $A,B\in K[\der]D$ we can divide by $D$ on the right and
reduce to the case that $D=1$, and so $K[\der](\der-a)A + K[\der]B=K[\der]E$ 
with $\order(E)\ge 1$. But also 
$1\in K[\der]A+ K[\der]B\subseteq K[\der]A+K[\der]E$, so $\order(E)\ge 2$ 
would contradict the previous lemma (with $E$ instead of $B$). Thus 
$\order(E)=1$.
\end{proof}

\noindent
For $A=\sum_i a_i\der^i\in K[\der]$ we define its {\bf adjoint} $A^*\in K[\der]$ by
$$ A^*\ :=\ \sum_i(-1)^i\der^ia_i,$$
so $a^*=a$ and $\der^*=-\der$.

\begin{lemma}
The map $A\mapsto A^*$ is an involution: for all $A,B\in K[\der]$,
$$ (A+B)^*\ =\ A^*+B^*, \qquad (AB)^*\ =\ B^*A^*, \qquad A{}^*{}^*\ =\ A.$$
\end{lemma}
\begin{proof}
The first identity is obvious. 
Let $K[\der]^{\operatorname{opp}}$ be the opposite ring of $K[\der]$; it has the same underlying additive abelian group as $K[\der]$, and its multiplication $\ast$ is given by
$A\ast B:=BA$.
Then $K[\der]^{\operatorname{opp}}$ contains $K$ as a subring and as a left $K$-module is free
with basis $(-\der)^0,(-\der)^1,(-\der)^2,\dots$, by \eqref{D^na}. Moreover
$(-\der)\ast a = a\ast (-\der) + a'$ for all $a\in K$. So  
the $K$-linear map $K[\der] \to K[\der]^{\operatorname{opp}}$ 
of left $K$-modules sending $\der^n$ to $(-\der)^n$ for each $n$, is a ring isomorphism. This map is nothing but $A \mapsto A^*$,  so $(AB)^* = (A^*)\ast (B^*) = B^*A^*$
for all $A,B\in K[\der]$.
%Now $A{}^*{}^*=A$ for all $A\in K[\der]$ follows by applying $%(\ )^*$ to both sides of the defining equation for~$A^*$.
\end{proof}

\noindent
Using this involution ``left'' results imply similar ``right'' results. For
example, we obtain in this way {\em division with remainder on the right}: 
for $A,B\in K[\der],\ B\ne 0$, 
there exist unique $Q,R\in K[\der]$ such that $A=BQ+R$ and 
$\order(R) < \order(B)$.
Hence for every right ideal $I$ of $K[\der]$ there is $B\in K[\der]$ such that
$I=BK[\der]$. 
%(A {\it right ideal\/} of $K[\der]$ is an additive 
%subgroup $I$ of  $K[\der]$ such that $AF\in I$ whenever
%$A\in I$ and $F\in K[\der]$.) 
The right analogues of the lemmas 
above now follow easily:

\begin{lemma}\label{fa3} Let $A,B\in K[\der]$ and $a$ satisfy 
$A(\der -a)\in BK[\der]$ and $\order(B)\ge 2$. Then 
$1\notin AK[\der]+ BK[\der]$.
\end{lemma}

\begin{lemma}\label{fa4} Let $A,B,D\in K[\der]$ and $a$ be such that
$$AK[\der]+ BK[\der]\ =\  D K[\der]\ \ne\ A(\der -a) K[\der]+ BK[\der].$$
Then  $A(\der -a)K[\der] + B K[\der] = D(\der-c)K[\der]$ for some $c\in K$.
\end{lemma}

\noindent
We also note:

\begin{lemma}\label{lem:ideals under extension}
Let $L$ be a differential field extension of $K$, and $A\in K[\der]$. Then 
$L[\der]A\cap K[\der]=K[\der]A$ and $AL[\der]\cap K[\der]=AK[\der]$. 
\end{lemma}
\begin{proof}
We may assume that $A\neq 0$.
Let $B\in L[\der]A\cap K[\der]$. Then $B=QA$ where $Q\in L[\der]$. Uniqueness of division with remainder on the left gives $Q\in K[\der]$. Thus $L[\der]A\cap K[\der]=K[\der]A$, and 
$AL[\der]\cap K[\der]=AK[\der]$ follows likewise.
\end{proof}

\index{adjoint!linear differential operator}
\index{linear differential operator!adjoint}
\nomenclature[Q]{$A^*$}{adjoint of  $A$}

\subsection*{The kernel of a linear differential operator} \textit{Let $K$ be 
a differential field in this subsection.}\/
Each $A\in K[\der]$ acts as a $C$-linear operator on $K$, and its kernel
$$\ker A\ =\ \bigl\{y\in K: A(y)=0\bigr\}$$
is a $C$-linear subspace of $K$, which 
is of dimension $\le n$ 
if $\order A \le n$, $A\ne 0$ (Co\-rol\-lary~\ref{cor:wronskian}).
If $A\in K[\der]^{\neq}$, $\dim_C \ker A=\order A$, and $L$ is a
differential field extension of $K$ with $C_L=C$, then $\ker_L A=\ker A$.
For 
$A,B\in K[\der]^{\ne}$ we have an exact sequence
$$0 \to \ker B \to \ker AB \to \ker A$$
of $C$-linear spaces, where the map on the right is given by $y\mapsto B(y)$, so
$$\dim_{C} \ker AB\ \le\ \dim_C \ker A + \dim_C \ker B, $$ 
with equality if and only if $B(K)\supseteq \ker A$. 
In particular: 

\begin{lemma}\label{lem:factoring and fund systems}
Let $A,B\in K[\der]^{\neq}$, and $m=\order(A)$, $n=\order(B)$. 
Then
$$\dim_C \ker AB\ =\ m+n \quad\Rightarrow\quad \text{$\dim_C \ker A = m$ and $\dim_C \ker B = n$.}$$
\end{lemma}

\noindent
We leave the proof of the next lemma to the reader:

\begin{lemma} \label{special form}
Let $A=a_0+a_1\der+\cdots+a_n\der^n\in K[\der]$ be of order 
$n\geq 1$. Suppose 
$g\in K^\times$ satisfies $g^\dagger=-(na_n)^{-1}a_{n-1}$, and put
$\tilde{A}=a_n^{-1}g^{-1}Ag\in K[\der]$. Then 
$$\tilde{A}\ =\ \tilde{a}_0+\tilde{a}_1\der+\cdots+
\tilde{a}_{n-2}\der^{n-2}+\der^n $$
for suitable $\tilde{a}_0,\dots,\tilde{a}_{n-2}\in K$. Moreover, we have 
an isomorphism 
$$y\mapsto gy\ \colon\  \ker \tilde{A}\to\ker A$$ of $C$-linear spaces. 
\end{lemma}

\begin{example}
If $A=a_0+a_1\der+\der^2$, then 
$\tilde{A}=\big(a_0 - \frac{1}{2}a_1' -\frac{1}{4}a_1^2\big)+\der^2$.
\end{example}

\begin{lemma}\label{lem:homog vs inhomog}
Let $A\in K[\der]$, $b\in K^\times$, and 
$B:=(\der-b^\dagger) A$. Let also $y\in K$. 
\begin{enumerate}
\item[\textup{(i)}] If $A(y)=b$, then  $B(y)=0$ 
and $\ker B = C y\oplus \ker A$.
\item[\textup{(ii)}] If $A(y)\neq 0$ and $B(y)=0$, then $A(cy)=b$ for some $c\in C^\times$.
\end{enumerate}
\end{lemma}
\begin{proof}
Part (i) is easy to see, and for (ii) note that if $A(y)\neq 0$ and $B(y)=0$, then $A(y)^\dagger=b^\dagger$
and hence $A(cy)=cA(y)=b$ for some $c\in C^\times$.
\end{proof}

\noindent 
The following lemma uses kernels to formulate a criterion for a differential operator to be contained in a given left ideal of $K[\der]$:

\begin{lemma}\label{lem:divisibility and kernels}
Let $A,B\in K[\der]^{\neq}$, $m=\order A$, $n=\order B$. 
\begin{enumerate}
\item[\textup{(i)}] If $\dim_C \ker A = m$ and $\ker A\subseteq \ker B$, then $B\in K[\der]A$;
\item[\textup{(ii)}] if $\dim_C \ker B = n$ and $B\in K[\der]A$, then $\ker A\subseteq\ker B$ and $\dim_C \ker A=m$.
\end{enumerate}
Thus if $\dim_C \ker A = m$ and $\dim_C \ker B = n$, then 
$$B\in K[\der]A\quad \Longleftrightarrow\quad \ker A\subseteq\ker B.$$
\end{lemma}
\begin{proof}
For (i), suppose $\dim_C \ker A = m$ and $\ker A\subseteq \ker B$. Take $Q,R\in K[\der]$ with $B=QA+R$, $\order R<m$;
then $\ker A\subseteq\ker R$, so $\order R<m\leq\dim_C\ker R$ and hence $R=0$. Part (ii) is immediate from
Lemma~\ref{lem:factoring and fund systems}.
\end{proof}

\subsection*{Linear differential equations with constant coefficients}
\textit{In this subsection~$K$ is a differential field.}\/ The ring
$C[\der]$ is commutative, and so we have a $C$-algebra isomorphism $P(Y)\mapsto P(\der)\colon C[Y]\to C[\der]$.  

\textit{Let $A\in C[\der]^{\ne}$, and take $P=P(Y)\in C[Y]^{\ne}$ with $A=P(\der)$, so $\order A=\deg P$. We also assume there is given an element $x\in K$ with $x'=1$, and that for each $c\in C$ there exists 
$a\in K^\times$ with $a'=ca$. For each $c\in C$ we pick such an $a$ and denote it suggestively by $\ex^{cx}$.}\/
Note that $x$ is transcendental over $C$, by 
Lemma~\ref{lem:C alg closed in K}. For $f\in C[x]$ we let $\deg_x f\in\N\cup\{-\infty\}$ be the degree of $f$ viewed as a polynomial in $x$ over $C$. Note also that for $c\in C$ we have $(\der-c)\ex^{cx}=\ex^{cx}\der$ in~$K[\der]$, and thus  $(\der-c)^i\ex^{cx}=\ex^{cx}\der^i$ in $K[\der]$ for all $i\in \N$.

\begin{lemma}\label{cor:linconstcoeff, 1}
Let $f\in C[x]$ and let $c\in C$ be such that $P(c)\neq 0$. Then
$$A(\ex^{cx}f)\ =\ \ex^{cx}g\qquad\text{for some $g\in C[x]$ with $\deg_x f=\deg_x g$.}$$
\end{lemma}
\begin{proof}
Take $a_0,\dots,a_n\in C$ such that $P(Y) = \sum_{i=0}^n a_i (Y-c)^i$. Then 
$A=\sum_{i=0}^n a_i (\der-c)^i$. Using the above identity
for $(\der-c)^i\ex^{cx}$ we get
$$A(\ex^{cx}f)\ =\  \ex^{cx}\left(\sum_{i=0}^n a_i f^{(i)}\right)\ =\ \ex^{cx}g\ \text{ for }\ g\ =\ \sum_{i=0}^n a_if^{(i)}\in C[x],$$
and $\deg_x g=\deg_x f$, since $a_0=P(c)\neq 0$.
\end{proof}

\begin{cor}\label{cor:linconstcoeff, 2}
Let $c_1,\dots,c_n\in C$ be distinct, where $n\geq 1$. Then the elements $\ex^{c_1 x}, \dots,\ex^{c_n x}$ of $K$ are linearly independent over $C(x)$.
\end{cor}
\begin{proof}
By induction on $n$. The case $n=1$ being obvious, let $n\geq 2$, and let $f_1,\dots,f_n\in C[x]$ satisfy $\sum_{i=1}^n  \ex^{c_i x} f_i=0$. Take $d\in\N$ with $d>\deg_x f_n$. Then by Lemma~\ref{cor:linconstcoeff, 1} applied to $(Y-c_n)^d$
instead of $P$ we get
$$0\ =\ (\der-c_n)^d\left(\sum_{i=1}^n  
\ex^{c_i x} f_i\right)\ =\ 
\sum_{i=1}^{n-1} \ex^{c_i x}g_i + \ex^{c_n x} f_n^{(d)}\ =\ \sum_{i=1}^{n-1} \ex^{c_i x}g_i$$
where $g_i\in C[x]$ with $\deg_x f_i= \deg_x g_i$ for $i=1,\dots,n-1$.
By inductive hypothesis we have $g_i=0$ and hence $f_i=0$, for $i=1,\dots,n-1$, and thus $f_n=0$.
\end{proof}

\noindent
Zeros of $P(Y)$ in $C$ yield elements of the kernel of $A$:

\begin{prop}\label{prop:linconstcoeff}
Let $c_1,\dots,c_n\in C$ be distinct zeros of $P$, of respective multiplicities $m_1,\dots,m_n\in\N^{\ge 1}$, $n\geq 1$.
Then $A(\ex^{c_i x} x^j)=0$ for $1\le i\le n$, $0\le j <m_i$, and the family $\big(\ex^{c_i x} x^j \big)_{1\le i\le n,\ 0\le j < m_i}$
 is linearly independent over $C$.
\end{prop}
\begin{proof}
Let $1\le i \le n$, $0 \le j< m_i$.
Then $P=Q\cdot(Y-c_i)^{m_i}$ where $Q\in C[Y]$. Set $B:= Q(\der)\in C[\der]$.
Then $A=B\cdot (\der-c_i)^{m_i}$, so
$$A(\ex^{c_i x} x^j)\ =\ B\big( \ex^{c_i x} (x^j)^{(m_i)} \big)\ =\  B(0)\ =\ 0.$$
The linear independence statement is immediate from Corollary~\ref{cor:linconstcoeff, 2}.
\end{proof}

\begin{notation} 
 Let $\operatorname{m}(P)$ denote the number of zeros of~$P$ in~$C$, counted with 
 multiplicity; thus $\operatorname{m}(P)=\sum_{c\in C} \operatorname{mul}(P_{+c})$, where only
 finitely many terms in this sum are nonzero.  
 Note that 
 $\operatorname{m}(P)=\operatorname{m}(P_{+a})$ for $a\in C$, and  
 $\operatorname{m}(P)=\operatorname{m}(P_{\times b})$ for $b\in C^\times$. Also, $\operatorname{m}(P)=\operatorname{m}(P_1)+\cdots+\operatorname{m}(P_n)$ if $P=P_1\cdots P_n$, $P_1,\dots,P_n\in C[Y]^{\neq}$ ($n\geq 1$).
 We also set $\operatorname{m}(A):=\operatorname{m}(P)$; so 
 $\operatorname{m}(A)=\operatorname{m}(A_1)+\cdots+\operatorname{m}(A_n)$ if
 $A=A_1\cdots A_n$, $A_1,\dots,A_n\in C[\der]^{\neq}$ ($n\geq 1$). 
\end{notation}

\nomenclature[O]{$\operatorname{m}(P)$}{$\sum_{c\in C} \operatorname{mul}(P_{+c})$, for $P\in C[Y]$}
\nomenclature[Q]{$\operatorname{m}(A)$}{$\operatorname{m}(P)$ for $P\in C[Y]$ with $A=P(\der)\in C[\der]$}

\noindent
From Proposition~\ref{prop:linconstcoeff} we obtain:
\begin{equation}\label{eq:mledim}\operatorname{m}(A)\ \leq\ \dim_C \ker A.
\end{equation}
If $C$ is algebraically closed, then $\order(A)=\operatorname{m}(A) = \dim_C \ker A$. 
Lemma~\ref{lem:linconstcoeff, omega} below gives another situation where  
$\operatorname{m}(A) = \dim_C \ker A$.

\subsection*{The type of a linear differential operator}
\textit{In this subsection $K$ is a differential field.}\/
We say that $A,B\in K[\der]$ have {\bf the same type} if the $K[\der]$-modules
$K[\der]/K[\der]A$ and $K[\der]/K[\der]B$ are isomorphic. 
For $a\in K^\times$ and $A\in K[\der]$ we have $K[\der]A=
K[\der]aA$, and the automorphism $B\mapsto Ba$ of the (left) $K[\der]$-module $K[\der]$ maps $K[\der]A$ onto~$K[\der]Aa$, so $A$, $aA$, and $Aa$ have the same type.  
For $A\in K[\der]^{\neq}$, the order of $A$ equals the $K$-vector space dimension of $K[\der]/K[\der]A$, so if $A,B\in K[\der]^{\neq}$ have the same type, then $\order(A)=\order(B)$.

\index{type!linear differential operator}
\index{linear differential operator!type}

\begin{lemma}\label{lem:same type}
Let $A,B\in K[\der]^{\neq}$. 
Then $A$ and $B$ have the same type if and only if $\order(A)=\order(B)$ and there
is $R\in K[\der]$ of order less than $\order(A)$ with $1\in K[\der]R+K[\der]A$ and $BR\in K[\der]A$.
\end{lemma}
\begin{proof}
Let $\phi\colon K[\der]/K[\der]B\to K[\der]/K[\der]A$ be an isomorphism of $K[\der]$-mod\-ules. 
Set $n:=\order(A)=\order(B)$. Then 
$\phi(1+K[\der]B)=R+K[\der]A$ with $R\in K[\der]$ of order 
$<n$.  One checks easily that then $1\in K[\der]R+K[\der]A$ and 
$BR\in K[\der]A$.
Conversely, assume $\order(A)=\order(B)=n$, and $1\in K[\der]R+K[\der]A$ and $BR\in K[\der]A$ with $R\in K[\der]$ of order $<n$. Then the $K[\der]$-linear map $K[\der]\to K[\der]/K[\der]A$ sending~$1$ to $R+K[\der]A$ is surjective and its kernel contains~$K[\der]B$, hence induces a $K[\der]$-linear map
$K[\der]/K[\der]B\to K[\der]/K[\der]A$, which is bijective since $\dim_K K[\der]/K[\der]A=n=\dim_K K[\der]/K[\der]B$.
\end{proof}

\begin{example}
Let $A,B\in K[\der]$ be monic of order $1$. Then $A$ and $B$ have the same type iff
$B=A_{\ltimes a}$ for some $a\in K^\times$.
\end{example}

\begin{cor}
Let $A,B\in K[\der]^{\neq}$ have the same type, and let $R\in K[\der]$ be as in Lemma~\ref{lem:same type}. Then
$R(\ker A)\subseteq\ker B$, and the map $x\mapsto R(x)\colon\ker A\to\ker B$ is an isomorphism of $C$-linear spaces.
\end{cor}
\begin{proof}
Take $L_1,L_2,L\in K[\der]$ with $1=L_1R+L_2A$ and $BR=LA$. Then for each $x\in \ker A$ we have $B\big(R(x)\big)=L\big(A(x)\big)=L(0)=0$, so $R(x)\in\ker B$, and $x=L_1\big(R(x)\big)$; hence we have an injective $C$-linear
map $x\mapsto R(x)\colon\ker A\to\ker B$. By symmetry, we also have an injective $C$-linear map $\ker B\to\ker A$,
showing that $\dim_C \ker A=\dim_C \ker B$, so  $x\mapsto R(x)\colon\ker A\to\ker B$ is surjective.
\end{proof}

\subsection*{Factorization in $K[\der]$} \textit{In this subsection 
$K$ is a differential field.}\/ We call $A\in K[\der]$ of positive order 
{\bf irreducible} \index{linear differential operator!irreducible} if there are
no $A_1,A_2\in K[\der]$ of positive order with $A=A_1A_2$. If 
$A\in K[\der]$ has positive order, then $A=A_1\cdots A_r$ 
with $r\in \N^{\ge 1}$ and each $A_i\in K[\der]$ irreducible.
We say that $A\in K[\der]$ 
{\bf splits} over $K$ \index{linear differential operator!splits} if
$A\ne 0$ and $A=c(\der-a_1)\cdots (\der  -a_n)$ for some~$c\in K^\times$ and 
some $a_1,\dots,a_n$. If $A,B\in K[\der]$
split over $K$, so does $AB$ 
(use twisting to eliminate factors $c\in K^\times$). 
In Section~\ref{sec:newtonization}  
we shall prove under natural assumptions on $K$ that each 
nonzero $A\in K[\der]$ splits over~$K$. 
In this connection the following facts are relevant.

\begin{lemma} \label{linear factor}
Let $A\in K[\der]$ be of order $n\ge 1$, and $u^\dag=a\in K$ with nonzero $u$ from some differential field extension of $K$. Then 
$$A\in K[\der] (\der-a)\ \Longleftrightarrow\ 
A(u)=0.$$ 
\end{lemma} 
\begin{proof} If $A=B\cdot(\der-a) + b$ with $B\in K[\der]$, $b\in K$, then
$A(u)=bu$.
\end{proof} 

\noindent
In particular, if $A\in K[\der]$ and $A(u)=0$, $u\in K^{\times}$, then 
$A\in K[\der](\der-u^\dagger)$.
%The main result of this subsection is the following:

\begin{prop}\label{splitfactor} Suppose $A,B,D\in K[\der]$, $A=BD$, and $A$ splits over $K$. 
Then $B$ and $D$ split over $K$.
\end{prop}
\begin{proof} We shall use Lemma~\ref{fa4} to show that $B$ splits over $K$. (A similar use of Lemma~\ref{fa2} shows that $D$ splits over $K$.)
Let $A=(\der-a_1)\cdots (\der  -a_n)$, $n\ge 1$, and for $i=0,\dots,n$, take the unique monic
$B_i\in K[\der]$ such that
$$B_iK[\der]\ =\ (\der-a_1)\cdots (\der-a_i)K[\der] + BK[\der],$$
so $B_0=1$ and $B_n=bB$, $b\ne 0$. By Lemma~\ref{fa4} we have, for $0\le i<n$, either
$B_{i+1}=B_i$, or $B_{i+1}=B_i(\der - c)$ for some $c\in K$.
Thus each $B_i$ splits over $K$, and so does $B$.
\end{proof}

\noindent
The next lemma shows how this proposition applies to linear parts of 
differential polynomials.
If $y$ in a differential field extension of $K$ is 
$\d$-algebraic over $K$ and~$P$ is a minimal annihilator of $y$ 
over $K$ (so $P\in K\{Y\}^{\ne}$ has minimal complexity
subject to $P(y)=0$), then $S_P(y)\ne 0$ and so $L_{P_{+y}}$ has order $r_P$.

\begin{lemma}\label{fa7} Let $E$ and $F$ be differential subfields of $K$ 
with $E\subseteq F \subseteq K$, and let $f\in K$ be $\d$-algebraic over $E$ with minimal annihilator $P$ over $E$ and $Q$ over $F$. 
Then $L_{P_{+f}}\in K[\der]L_{Q_{+f}}$, so if $L_{P_{+f}}$ splits over $K$, then so does $L_{Q_{+f}}$. 
\end{lemma}
\begin{proof} From $P(f)=0$ and Corollary~\ref{cor:rittvariant} we obtain
$$  H P\ =\ A_0Q+A_1Q' + \dots + A_nQ^{(n)}, \qquad H\ =\ S_Q^l,$$
where $l\in \N$, $A_0,\dots, A_n\in F\{Y\}$. By additive conjugation,
$$H_{+ f} P_{+  f}\ =\ (A_0)_{+ f} Q_{+ f} + \dots + (A_n)_{+f} Q^{(n)}_{+ f}.$$
Taking linear parts and using $ P(f)=Q( f)=0$ gives
\begin{align*} H(f)L_{ P_{+ f}}\ &=\ A_0(f)L_{Q_{+f}} + \dots + A_n(f) L_{Q^{(n)}_{+f}}\\
      &=\ \big(A_0(f)+ \cdots +
        A_n( f)\der^n\big)L_{Q_{+ f}}.
        \end{align*} 
It only remains to note that $H(f)\ne 0$.         
\end{proof}

\subsection*{Factorization into irreducibles and composition series} \textit{In this subsection~$K$ is a differential field, $R=K[\der]$, and $A\in R^{\neq}$ is monic.}\/
We discuss in what sense a  factorization of $A$ into irreducible elements of~$R$ is unique.
Given $A_1,A_2\in R$ with $A=A_1A_2$, we have $RA\subseteq RA_2$, which gives an  exact sequence 
$$0\longrightarrow R/RA_1 \longrightarrow R/RA \longrightarrow R/RA_2 \longrightarrow 0$$
of $R$-linear maps, 
with $R/RA_1 \to R/RA$ induced by the
$R$-linear endomorphism $B\mapsto BA_2$ of $R$. 
Conversely, let $M_1$ be a submodule of the $R$-module $M=R/RA$. Take monic $A_2\in R$ such that $RA_2$ is the kernel
of the composition 
$$R\to R/RA=M\to M/M_1$$ 
of the natural surjections. 
Then $A\in RA_2$, we have monic $A_1\in R$
with $A=A_1A_2$, and so we obtain a commutative
diagram of $R$-linear maps
$$\xymatrix@C+2em{	0 \ar[r]	& M_1 \ar@{^{(}->}[r] 							& M \ar[r] 				& M/M_1 \ar[r] 			& 0 \\
			0 \ar[r]	& R/RA_1 \ar[u] \ar[r] 	& R/RA \ar[u]^{=} \ar[r]	& R/RA_2 \ar[u]\ar[r]	& 0}$$
where the map $R/RA_1\to R/RA$ is induced by the endomorphism $B \mapsto BA_2$ of the $R$-module $R$, and the vertical arrows are isomorphisms, the rightmost one being $B+RA_2\mapsto B+M_1$ for $B\in R$.   
Thus $A$ is irreducible (in $R$) iff the $R$-module~$R/RA$ is irreducible. Using this recursively, a factorization $A=A_1\cdots A_m$ of~$A$ into irreducible $A_1,\dots,A_m\in R$ gives rise to a composition series
$$\{0\}\ =\ M_0\ \subset\ M_1\ \subset\ \cdots\ \subset\
M_m\ =\ M$$
of length $m$ of the $R$-module $M=R/RA$ and for $i=1,\dots,m$
isomorphisms $M_i \cong R/RA_1\cdots A_{i}$ 
and $M_i/M_{i-1}\cong R/RA_{i}$ of $R$-modules. 

\medskip
\noindent
Usually, $A$ has many factorizations into irreducible monic operators; e.g., if $K=C(x)$ where $x\notin C$ and $x'=1$, then $\der^2=(\der+f^\dagger)(\der-f^\dagger)$ for all $f=ax+b$ with~$a,b\in C$, not both zero.
However, the \textit{number}\/ and the \textit{types}\/ of the irreducible factors in such a factorization do not depend on the particular factorization:

\begin{cor}\label{cor:uniqueness of factorization}
If $A_1,\dots, A_m,B_1,\dots, B_n$ are irreducible elements of $R$ such that $A_1\cdots A_m=B_1\cdots B_n$, then $m=n$ and there is a permutation $i\mapsto i'$ of $\{1,\dots,m\}$ such that
$A_i$ and $B_{i'}$ are of the same type for $i=1,\dots,m$.
\end{cor}

\noindent
This corollary  follows from the remarks preceding it in combination with the Jordan-H\"older Theorem (Corollary~\ref{cor:JH}).
Note that Corollary~\ref{cor:uniqueness of factorization} gives rise to another proof of Proposition~\ref{splitfactor}. The remarks above also yield:

\begin{cor}\label{cor:splitting and comp series}
The linear differential operator $A$ splits over $K$ if and only if the $R$-module $M=R/RA$ has a composition
series 
$$\{0\}\ =\ M_0\ \subset\ M_1\ \subset\ \cdots\ \subset\ M_m\ =\ M$$
with $\dim_K M_i/M_{i-1}=1$ for $i=1,\dots,m$.
\end{cor}

\subsection*{Linear closedness and linear surjectivity} \textit{In this subsection
 we let $K$ be a differential field, and we let $r$ range over $\N^{\ge 1}$.}\/
We define $K$ to be {\bf $r$-linearly closed} 
if each nonzero $A\in K[\der]$ of order $\le r$ splits over~$K$.
We define $K$ to be {\bf linearly closed} \label{p:linearly closed} if it is $r$-linearly closed for 
each $r$. 
If the derivation of $K$ is trivial, that is, $K=C$, then 
$K[\der]=C[\der]$ is the usual (commutative) polynomial ring with $\der$
as an indeterminate over $C$, so in that case $K$ is linearly closed iff~$C$ is algebraically closed. 
%Does it follow more generally that if $K$ is linearly closed, 
%then $C$ is algebraically closed?
The property of being linearly closed has a nice first-order axiomatization: 

\index{differential field!linearly closed}
\index{linearly!closed}

\begin{lemma} Let $n\ge 1$ be given. Then there exists a
differential polynomial 
$$d_n(Y_1,\dots, Y_n,Y)\ \in\ \Q\{Y_1,\dots, Y_n,Y\} \qquad  
(\text{with all its coefficients in $\Z$})$$ 
such that for each differential ring
$R$ and all $a_1,\dots, a_n,a\in R$,
$$d_n(a_1,\dots, a_n,a)=0\ \Longleftrightarrow\ \der^n + a_1\der^{n-1} + \cdots + a_n = (\der-a)B \text{ for some $B\in R[\der]$.}$$ 
\end{lemma} 
\begin{proof} Just note that in the ring $\Q\{Y_1,\dots, Y_n,Y\}[\der]$ we have
$$\der^n + Y_1\der^{n-1} + \cdots + Y_n\ =\ (\der-Y)B + d_n$$
with $B\in \Q\{Y_1,\dots, Y_n,Y\}[\der]$ and $d_n\in \Q\{Y_1,\dots, Y_n,Y\}$. 
\end{proof}

\noindent
Thus $K$ is linearly closed iff for each $n\ge 2$,
$$ K \models\ \forall y_1\cdots \forall y_n \exists y\  d_n(y_1,\dots, y_n,y)=0.$$
We say that $K$ is {\bf $r$-linearly surjective} if $A(K)=K$
for each nonzero 
$A\in K[\der]$ of order
at most $r$, equivalently,
for all $a_0,\dots,a_r\in K$ such that $a_i\ne 0$ for some~$i$, the 
inhomogeneous linear differential equation 
$1+a_0y + \dots + a_ry^{(r)}=0$ has a solution in $K$. We say that $K$ is {\bf linearly surjective} if it is $r$-linearly 
surjective for each $r$. \label{p:linearly surjective}\index{differential field!linearly surjective}\index{linearly!surjective} One shows easily:

\begin{lemma}\label{lem:linsurproperties} Linear surjectivity has the following properties:
\begin{enumerate}
\item[\textup{(i)}] if $A(K)=K$ for every monic irreducible
$A\in K[\der]$ of order $\leq r$, then $K$ is $r$-linearly surjective;
\item[\textup{(ii)}] if $K$ is linearly closed and $1$-linearly surjective, then 
$K$ is linearly surjective;
\item[\textup{(iii)}] if $K$ is a directed union of $r$-linearly surjective differential
subfields, then $K$ is $r$-linearly surjective;
\item[\textup{(iv)}] if $K$ is $r$-linearly surjective and $A_1,\dots,A_n\in K[\der]^{\neq}$ \textup{(}$n\geq 1$\textup{)} are of order at most~$r$, then
$\dim_C \ker A_1\cdots A_n\ =\  \dim_C \ker A_1 + \cdots + \dim_C \ker A_n$.
\end{enumerate}
\end{lemma}

\noindent
Note that (ii) follows from (i).
Next we establish a descent property of linear surjectivity. Let $L$ be a differential
field extension of $K$ with $d:=[L:K]<\infty$. Let~$K^\alg$ be an algebraic closure of $K$ 
with the unique derivation extending the
derivation of~$K$. For $y\in L$, let
$$\operatorname{tr}_{L|K}(y)\ :=\  \sum_{i=1}^d \sigma_i(y)$$
be the trace of $y$ in $L|K$. Here $\sigma_1,\dots,\sigma_d$
are the distinct field embeddings $L\to K^\alg$ which are the identity on $K$.
The map 
$$y\mapsto \tau_{L|K}(y) := \frac{1}{d}\operatorname{tr}_{L|K}(y)\ \colon\ 
L\to K$$ 
is $K$-linear and the identity on $K$. Moreover, $\tau_{L|K}(f')=\tau_{L|K}(f)'$
for each $f\in L$. It follows that if
$A\in K[\der]$, $g\in K$, and the equation $A(y)=g$ has a solution $f\in L$, 
then it has the solution $ \tau_{L|K}(f)$ in $K$. Since every 
algebraic extension field of $K$ is a directed union of algebraic extension 
fields of finite degree, this gives the following: 

\begin{cor}\label{linsuralgdesc} If $K$ has an $r$-linearly surjective 
algebraic 
differential field extension, then $K$ is $r$-linearly surjective. 
\end{cor}   

\noindent
In Section~\ref{sec:systems} below we show that algebraic differential field extensions 
of linearly surjective differential fields are linearly surjective. This 
requires some generalities on systems of linear differential equations. In Section~\ref{sec:differential modules} we 
prove the same with {\em linearly closed\/} instead of {\em linearly surjective}. 

\medskip
\noindent 
By convention we consider any $K$ as being $0$-linearly surjective.

\subsection*{Picard-Vessiot closed differential fields} 
\textit{In this subsection $K$ is a differential field, and $r$ ranges over $\N^{\geq 1}$.}\/
We assume here familiarity with the facts stated about zeros of linear
differential polynomials in 
Section~\ref{Differential Fields and Differential Polynomials}.
We say that $K$ is {\bf $r$-Picard-Vessiot closed} (or {\bf $r$-pv-closed}) if
$\dim_C \ker A = \order A$ for all $A\in K[\der]^{\neq}$ of order at most $r$.
% every nonzero $A\in K[\der]$
%of order at most $r$ has a fundamental system of zeros in $K$. 
If 
$K$ is a directed union of $r$-pv-closed differential
subfields, then $K$ is $r$-pv-closed. Also, $K$ is $1$-pv-closed iff $(K^\times)^\dagger=K$.
We  say that $K$ is {\bf Picard-Vessiot closed} (or {\bf pv-closed}) if
$K$ is $r$-pv-closed for all $r$. This subsection, with some additional material in Section~\ref{sec:differential modules}, will mainly get used in the next volume.

\label{p:pv-closed}
\index{differential field!Picard-Vessiot closed}
\index{differential field!pv-closed}
\index{Picard-Vessiot closed}

%\marginpar{In Lemma~3.9.1 add the sentence: If $K$ is $r$-pv-closed, then so is $K^\phi$.}

\begin{lemma}\label{le:r-pv closed} The following two conditions on $K$ are equivalent: 
\begin{enumerate}
\item[\textup{(i)}] $K$ is $r$-linearly closed and $1$-pv-closed;
\item[\textup{(ii)}] $\ker A\neq \{0\}$ for all $A\in K[\der]\setminus K$ of order at most $r$.
\end{enumerate}
\end{lemma}
\begin{proof}
Use Lemma~\ref{linear factor}.
\end{proof}

\begin{lemma}\label{lem:r-pv closed}
Suppose $r\geq 2$. Then the following are equivalent:
\begin{enumerate}
\item[\textup{(i)}] $K$ is $r$-pv-closed;
\item[\textup{(ii)}] $K$ is $r$-linearly closed, $r$-linearly surjective, and $1$-pv-closed;
\item[\textup{(iii)}] $K$ is $r$-linearly closed, $1$-linearly surjective, and $1$-pv-closed.
\end{enumerate}
In particular, 
$$\text{$K$ is pv-closed}\quad\Longleftrightarrow\quad \text{$K$ is linearly closed, linearly surjective, and $(K^\times)^\dagger=K$.}$$
\end{lemma}
\begin{proof}
To show (i)~$\Rightarrow$~(ii),
suppose $K$ is $r$-pv-closed. Then $K$ is $r$-linearly closed and $1$-pv-closed, by the previous lemma.
To show that $K$ is $r$-linearly surjective, it suffices to show 
that $K$ is $1$-linearly surjective. Let $A\in K[\der]^{\neq}$ with $\order A=1$, and $b\in K^\times$; we need to show $b\in A(K)$.
Put $B:=(\der-b^\dagger)\cdot A$, so $\order B=2\leq r$;
then $\dim_C \ker A = 1 < 2 = \dim_C\ker B$,
so there exists $y\in K$ with $A(y)\neq 0$, $B(y)=0$, and then $A(cy)=b$ for some $c\in C$, by 
Lemma~\ref{lem:homog vs inhomog}(ii). 

The implication (ii)~$\Rightarrow$~(iii) is trivial, and  (iii)~$\Rightarrow$~(i) follows from the fact that if $A,B\in K[\der]^{\neq}$
with $B(K)=K$, then $\dim_C \ker AB = \dim_C\ker A+\dim_C\ker B$.
\end{proof}

\noindent
By Lemmas~\ref{le:r-pv closed} and~\ref{lem:r-pv closed}, every 
differentially closed field is linearly closed, linearly surjective, and pv-closed.

\begin{lemma}\label{lem:pv closed => alg closed}
Suppose $K$ is pv-closed and $C$ is algebraically closed. Then $K$ is algebraically closed.
\end{lemma}
\begin{proof}
Let $P\in K[Y]$ have degree $\ge 2$, $P\notin YK[Y]$. Take an algebraic closure $K^\alg$ of $K$ and give it the unique derivation extending the derivation of~$K$. Then $C_{K^{\alg}}=C$ 
%constant field of  $K^\alg$ is $C$. 
by Lemma~\ref{lem:const field of algext}.
Let $V$ be the $C$-linear subspace of  $K^\alg$ generated by
the zeros of $P$ in~$K^\alg$. Let
$y_1,\dots,y_n$ be a basis for $V$ (so $n\geq 1$), 
and let 
$$A(Y)\ :=\ \frac{\wr(y_1,\dots,y_n,Y)}{\wr(y_1,\dots,y_n)}\ =\ Y^{(n)} + a_1Y^{(n-1)} + \cdots + a_nY, \quad a_1,\dots, a_n\in K^\alg.$$
Then  $V=Z(A)$. If also $V=Z(B)$, $B=Y^{(n)}+b_1Y^{(n-1)} + \cdots + b_nY$ where $b_1,\dots, b_n\in K^{\alg}$, then $A=B$ by 
Corollary~\ref{cor:wronskian, 2}. Each $\sigma\in\Aut(K^\alg|K)=\Aut_\der(K^\alg|K)$ satisfies $\sigma(V)=V$ and so
fixes the coefficients of $A$.  
Thus $A\in K\{Y\}$.
Since $K$ is pv-closed, this yields $V\subseteq K$; in particular
$P(y)=0$ for some $y\in K$.
\end{proof}

\subsection*{Relating  $K[\imag]$  and  $K$} \textit{In this subsection $K$ is a 
differential field in which $-1$ is not a square in $K$.}\/
It follows that $K$ has a (differential) field extension 
$K[\imag],\ \imag^2=-1$. The results below indicate how
factorization in $K[\der]$ relates to factorization in~$K[\imag][\der]$.
This will be used in later chapters where $K$ is real closed,
and thus $K[\imag]$ is algebraically closed.

\begin{lemma} If every $A\in K[\der]^{\ne}$ splits over $K[\imag ]$, then $K[\imag]$ is linearly closed.
\end{lemma}
\begin{proof} The ``complex conjugation'' automorphism $a+b\imag\mapsto \bar{a+b\imag}:= a-b\imag$ ($a,b\in K$) of the differential field $K[\imag]$ extends uniquely to an automorphism
$A\mapsto \bar{A}$ of the ring $K[\imag][\der]$ with $\bar{\der}=\der$. Then
$\bar{A}=A\Longleftrightarrow A\in K[\der]$, for all $A\in K[\imag][\der]$. Let
$A\in K[\imag][\der]^{\ne}$ be monic. Let $B$ be the least common left multiple 
of $A$
and $\bar{A}$. Then $\bar{B}=B$, so $B\in K[\der]$. If
$B$ splits over $K[\imag]$, then so does $A$, by Proposition~\ref{splitfactor}.
\end{proof}

\begin{lemma} \label{first or second order}
Suppose $u$ is a nonzero element in a 
differential field extension of~$K[\imag]$ such that 
$u^\dagger\in K[\imag]$. Then there is $B\in K[\der]$ of order $2$ such that 
$B(u)=0$.
\end{lemma}
\begin{proof}
Write $u^\dagger=a+\imag b$, so $u'=au+\imag bu$.
Differentiating this relation yields
$$u''\ =\ (a'+a^2 -b^2)u + \imag (2ab+b')u.$$
If $b=0$, then we can take $B=\der^2-(a'+a^2)$. Suppose that 
$b\neq 0$. Then we eliminate~$\imag$ by forming a suitable linear 
combination of $u''$ and $u'$:
$$bu''-(2ab+b')u'\ =\ (-ab' + a'b-a^2b-b^3)u. $$ 
It follows that we can take
$B:= b\der^2 -(2ab+b')\der + (ab'-a'b+a^2b + b^3)$.  
\end{proof}

\begin{lemma}\label{fa5} Suppose $A\in K[\der]$ is monic and irreducible over $K$ of order~$2$
and splits over $K[\imag]$. Then 
$$A=\big(\der-(a-b\imag+b^\dagger)\big)\cdot\big(\der-(a+b\imag)\big)\quad\text{for some $a\in K$, $b\in K^\times$.}$$
\end{lemma}
\begin{proof} We have $A=(\der - c)\big(\der - (a+b\imag)\big)$ with
$c\in K[\imag]$ and $a,b\in K$. Then $b\ne 0$ and $c=a-b\imag+b^\dagger$, and
so $A$ has the desired form. Note also that
\equationqed{A\ =\ \der^2 - (2a+b^\dagger)\der + (-a'+a^2 +ab^\dagger + b^2).}
\end{proof}

{\sloppy
\begin{lemma}\label{fa6} Suppose $A\in K[\der]$ is monic and splits over 
$K[\imag]$.
Then $A=A_1\cdots A_m$ with all $A_i\in K[\der]$ monic and irreducible of order $1$ or order $2$.
\end{lemma} 
}
\begin{proof} We proceed by induction on $\order(A)$ and can assume $A$ has order $n>1$. By Lemma~\ref{linear factor} we have a nonzero $u$ in some differential field extension of $K[\imag]$
such that $u^\dagger\in K[\imag]$ and $A(u)=0$. By 
Lemma~\ref{first or second order} this
gives a monic $D\in K[\der]$ of order $1$ or $2$ such that $D(u)=0$. Take such 
$D$ of
least order. Then $A=BD+R$ with $B,R\in K[\der]$ and $\order(R)<\order(D)$, 
hence
$R(u)=0$, so $R=0$ and thus $A=BD$. Now $B$ and $D$ split over $K[\imag]$ by 
Lemma~\ref{fa4}, and so we can apply the inductive assumption to $B$. 
\end{proof}

\noindent
For $f=a+b\imag\in K[\imag]$ ($a,b\in K$) we set $\bar{f}:= a-b\imag$, $\Re(f):= a$, 
$\Im(f):= b$, so $f\mapsto \bar{f}$ is a differential field
automorphism of $K[\imag]$.
Now let $F$ be a differential
subfield of~$K$, and suppose $f\in K[\imag]$ is $\d$-algebraic over 
$F[\imag]$ with
minimal annihilator~$P(Y)$ over~$F[\imag]$. Then
$f$, $\bar{f}$,  and $\Re(f)$ are $\d$-algebraic over $F$. Set 
$r:=r_P$; then
$L_{P_{+f}}$ has order $r$, with coefficient $S_P(f)$ of $\der^r$.
  
\nomenclature[A]{$\Re(f)$}{real part of $f\in K[\imag]$}
\nomenclature[A]{$\Im(f)$}{imaginary part of $f\in K[\imag]$}

\begin{lemma}\label{fa8} Suppose $L_{P_{+f}}$ splits over $K[\imag]$, and 
let $R$ be a minimal annihilator of~$\Re f$ over $F[\imag]\langle f \rangle$. Then 
$L_{R_{+\Re f}}$ splits over $K[\imag ]$.
\end{lemma}
\begin{proof} Set $Q:= R_{\times 1/2, +f}$, so $R=Q_{-f,\times 2}$. From $\Re(f)=\frac{1}{2}(f+\bar{f})$ we get $Q(\bar f)=0$. Thus $Q$ is a minimal
annihilator of $\bar f$ over $F[\imag]\langle f \rangle$. Also, $\bar P$ is a minimal annihilator of $\bar f$ over $F[\imag]$, hence $L_{Q_{+\bar f}}$ splits over $K[\imag]$ by Lemma~\ref{fa7}. Now 
$$ Q_{+\bar f}\ =\  R_{+\Re f, \times 1/2}, \qquad L_{Q_{+\bar f}}\ =\ L_{R_{+\Re f, \times 1/2}},$$
so $L_{R_{+\Re f}}=L_{R_{+\Re f,\times 1/2, \times 2}}=2L_{R_{+\Re f,\times 1/2}}$ splits 
over $K[\imag]$.
\end{proof}

\begin{lemma}\label{fa9} Suppose $K[\imag]$ is $r$-linearly closed, and 
$S\in K\{Y\}$ is a minimal annihilator of $\Re f$ over $F$. 
Then $L_{S_{+\Re f}}$ splits over $K[\imag]$.
\end{lemma}
\begin{proof} Since $K[\imag]$ is $r$-linearly closed, $L_{P_{+f}}$ splits over
$K[\imag]$. Let $R$ be a minimal annihilator of $\Re f$ over $F[\imag]\langle f \rangle$. Then
$L_{R_{+\Re f}}$ splits over $K[\imag]$ by Lemma~\ref{fa8}. We now apply 
Lemma~\ref{fa7} with 
$\Re f$, $S$, $R$ in the role of $f$, $P$, $Q$, and obtain 
$$L_{S_{+\Re f}} = A\cdot L_{R_{+\Re f}}\ \text{ with }A\in K[\imag][\der],\quad \order A= r_S - r_R.$$
So it is enough to show that $A$ splits over $K[\imag]$.
Let $Q: =R_{\times 1/2, +f}$ be as in the proof of Lemma~\ref{fa8}, so $r_Q=r_R$.
By considering the inclusions among the fields
$$F[\imag],\quad E:=F[\imag]\langle f \rangle,\quad E\langle  \bar f \rangle= 
E\langle  \Re f\rangle,\quad F[\imag]\langle \Re f \rangle$$
and the corresponding transcendence degrees, we obtain
$$r + r_R\ =\  r+ r_Q\ =\  r_S + \operatorname{trdeg}\big(E\langle  \Re f \rangle | F[\imag]\langle \Re f \rangle\big),$$
so $r_S - r_R = r - \operatorname{trdeg}\big(E\langle \Re f \rangle | F[\imag]\langle \Re f \rangle\big) \le r$. Thus $A$ splits over $K[\imag]$.
\end{proof}

\subsection*{Linear elements} \textit{In this subsection $K$ is a differential field, $y$ lies in a differential field extension of $K$, and $r\in \N$.}\/ We say that $y$ is {\bf $r$-linear} over $K$ if there is an  $A\in K[\der]^{\neq}$ of order~$\leq r$ such that $A(y)=0$.
Note that $0\in K$ is $0$-linear over $K$, and each $a\in K^\times$ is $1$-linear over~$K$
(take $A=\der-a^\dagger$). If $y$ is $r$-linear over $K$, then
$$y^{(n)}\ \in\ Ky + \cdots+ Ky^{(r-1)}\ \text{ for all $n$,}$$
and thus $K\<y\>=K(y,y',\dots,y^{(r-1)})$ as fields.
Lemma~\ref{lem:homog vs inhomog}(i) shows that $y$ is $(r+1)$-linear over $K$ if  $A(y)\in K$ for some  $A\in K[\der]^{\neq}$ of order~$\leq r$.
We say that~$y$ is {\bf linear} over $K$ if it is $r$-linear over $K$
for some $r$. Thus $y$ is linear over~$K$
iff $A(y)\in K$ for some  $A\in K[\der]^{\neq}$. ``Algebraic over $K$'' implies  ``linear over~$K$'':
 
%\index{$r$-linear}
\index{element!linear over $K$}

\begin{lemma}\label{lem:algebraic implies linear}
Suppose $r=[K(y):K]<\infty$.
Then $y$ is $r$-linear over $K$.
\end{lemma}
\begin{proof}
Let $P\in K[Y]$ be the minimum polynomial of $y$ over $K$. 
From $P'(y)\neq 0$ we get $y'=-P^\der(y)/P'(y) \in K(y)$, by Lemma~\ref{derpol}. Hence
$y^{(n)}\in K(y)$ for all $n$; so a nontrivial $K$-linear combination of $y,y',\dots,y^{(r)}$ is zero, that is, $A(y)=0$ for some  $A\in K[\der]^{\neq}$ of order~$\leq r$.
\end{proof}

\noindent
Sums, products, and derivatives of linear elements remain linear:

\begin{lemma}\label{lem:algebraic combinations of linear elements} Let
$L$ be a differential field extension of $K$ and let $y,z\in L$ be
$m$-linear and $n$-linear over $K$, respectively. Then $y'$ is $m$-linear over $K$, $y+z$ is $(m+n)$-linear over $K$, and $yz$ is $mn$-linear over $K$.
\end{lemma}
\begin{proof} All $y^{(i)}$ ($i\ge 0$) lie in the $K$-linear space
$Ky + \cdots + Ky^{(m-1)}$ of dimension at most~$m$, so  
$y',\dots, y^{(m+1)}$ are $K$-linearly dependent, and thus $y'$ is 
$m$-linear over~$K$. Also
$z^{(i)}\in Kz + \cdots + Kz^{(n-1)}$ for all $i$, and
so all $(y+z)^{(i)}$ lie in the $K$-linear space  
$Ky + \cdots + Ky^{(m-1)}+ Kz + \cdots + Kz^{(n-1)}$ of dimension at most~$m+n$.
Hence $y+z, (y+z)', \dots, (y+z)^{(m+n)}$ are $K$-linearly dependent, and thus 
$y+z$ is $(m+n)$-linear over $K$.
Finally, $(yz)^{(i)}\in\sum_{j< m, k< n} Ky^{(j)}z^{(k)}$ for all $i$, and so by the same reasoning  $yz$ is $mn$-linear over $K$. 
\end{proof}

\begin{cor}\label{cor:algebraic combinations of linear elements} If $L$ is a differential field extension of~$K$ and $y_1,\dots,y_n\in L$ are linear over $K$, 
then all elements of $K\{y_1,\dots,y_n\}$ are linear over $K$.
\end{cor}

\noindent
There is no such result for reciprocals:

\begin{lemma}\label{lem:reciprocal of antider}
Let $a\in K\setminus\der(K)$, and $y'=a$.
Then $y$ is transcendental over $K$ and $2$-linear over $K$, but $z=1/y$ is not linear over $K$.
\end{lemma}
\begin{proof}
Suppose towards a contradiction that $y$ is algebraic over $K$, and let $P=\sum_{i=0}^r a_i Y^i$ be the minimum polynomial of $y$ over $K$ ($r\ge 1$, $a_i\in K$ for $i=0,\dots,r$, $a_r=1$).
By Lemma~\ref{derpol}, we then have $P^\der(y)+P'(y)a=0$ and
hence $$P^\der(Y)+P'(Y)a\ =\ 0\  \text{ in $K[Y]$.}$$ In particular, $a_{r-1}'+ra=0$, so $a=(-a_{r-1}/r)'\in \der(K)$, a contradiction. Thus~$y$, and hence also $z=1/y$, are transcendental over $K$.
To see that $z$ is not linear over~$K$, one first shows by an easy induction on $n\ge 1$ that
$$z^{(n)}\ =\  n!(-a)^n\, z^{n+1}+f_n\qquad\text{where $f_n\in Kz^2+\cdots+Kz^n$.}$$
Thus from $A(z)=0$ with monic $A\in K[\der]$ of order~$r$ we obtain 
$$r!(-a)^r\,z^{r+1}+g(z)\ =\ 0\qquad\text{where $g(z)\in Kz+\cdots+Kz^r$,}$$
contradicting the transcendence of $z$ over $K$.
\end{proof}

\subsection*{Notes and comments} Early algebraic studies of the ring $K[\der]$ are Ore~\cite{OreI, OreII}, but the {\em type\/} of a linear differential operator
occurs already in Poincar\'e~\cite{Poincare}.
Lemma~\ref{lem:divisibility and kernels} is \cite[Lemma~2.1]{Singer96}.  
Proposition~\ref{prop:linconstcoeff} goes back to Euler~\cite{Euler43}.
Corollary~\ref{cor:uniqueness of factorization} is due to Landau~\cite{Landau} and Loewy~\cite{Loewy03}.
Lemma~\ref{lem:algebraic combinations of linear elements} and its Corollary~\ref{cor:algebraic combinations of linear elements} are from~\cite{Singer79}, but 
have a classical origin \cite[Sec.~167]{Schlesinger}.
In connection with Lemma~\ref{lem:reciprocal of antider} we note that by~\cite{HarrisSibuya}, if
$y$ is a nonzero element of a differential field extension of $K$ such that both~$y$ and~$1/y$ are linear over $K$, then $y^\dagger$ is algebraic over $K$. See \cite[Section~2.4]{Singer96} for some history of the algebraic study of linear differential operators.

%% file: mt-5-2.tex
\section{Second-Order Linear Differential Operators}\label{sec:secondorder}

\noindent
Lemma~\ref{fa6} above is a first indication that linear differential operators of order~$1$ and~$2$ are going to play a special role. In preparation for this, we now make a few basic observations concerning linear differential operators
of order $2$. 

\medskip
\noindent
Let $K$ be a differential field, and let $A\in K[\der]$ be of order $2$. Lemma~\ref{special form} and the example following it show that if $(K^\times)^\dagger=K$, then in order to describe the kernel $\ker A$ of $A$, one may 
reduce to the case that  
$$A\ =\ 4\der^2+f\qquad (f\in K)$$ 
which we assume in the rest of this section. Note that $\dim_C\ker A\leq 2$, and $A$ is irreducible in $K[\der]$ iff $A$ does not split over $K$.
We introduce the function \nomenclature[M]{$\omega(z)$}{$\omega(z)=-(2z'+z^2)$}
$$z\mapsto\omega(z):=-(2z'+z^2)\ \colon\ K\to K.$$
Then for $y\in K^\times$ we have $\omega(2y^\dagger)=-4y''/y$, hence
$A(y)=y(f-\omega(2y^\dagger))$, and thus
$$A(y)=0 \quad\Longleftrightarrow\quad 4y''+fy=0 \quad\Longleftrightarrow\quad \omega(2y^\dagger)=f.$$
Hence
$$\dim_C\ker A\geq 1 \ \Longleftrightarrow\  f\in\omega\big(2(K^\times)^\dagger\big).$$
Another way to formulate this involves the set 
$$\Upo(K)\ :=\ \{f\in K:\  \text{$4y''+fy=0$ for some $y\in K^\times$}\}.$$
Thus $\Upo(K)$ is existentially definable in $K$, and if  
$(K^\times)^\dagger= K$, then $\Upo(K)=\omega(K)$.

\nomenclature[M]{$\Upo(K)$}{set of $f\in K$ such that $4y''+fy=0$ for some $y\in K^\times$}

If $A$ splits over $K$, then $A=4(\der+a)(\der-a)$ and $f=\omega(2a)$ for some
$a\in K$. Conversely, if $z\in K$ and $\omega(z)=f$, then  
$A=4 \left(\der + \frac{z}{2}\right)\left(\der-\frac{z}{2}\right)$.
Together with Lemma~\ref{linear factor}, the above yields
$$\dim_C\ker A\geq 1 \ \Rightarrow\  \text{$A$ splits over $K$}  \ \Longleftrightarrow\  f\in\omega(K).$$
%If $z\in K$ and $\omega(z)=f$, then  $A$ factors as
%$$A\ =\ 4 \left(\der + \frac{z}{2}\right)\left(\der-\frac{z}%{2}\right).$$
If $(K^\times)^\dagger=K$
(equivalently, $\dim_C \ker B=1$ for each $B\in K[\der]$ of or\-der~$1$), then 
$$\dim_C\ker A\geq 1 \ \Longleftrightarrow\ \text{$A$ splits over $K$}  \ \Longleftrightarrow\  f\in\omega(K).$$
The following computations involving $\omega$ will be used several times later on:

\begin{lemma}\label{varrholemma, claim, 1}
Let $w,y,z\in K$ be such that $y=z-w\neq 0$. Then
$$\omega(w)-\omega(z)\ =\ y\cdot\bigl(2(y^\dagger + w)+y\bigr).$$
\end{lemma}

\vskip-1.5em

\begin{flalign*}
\hskip\normalparindent&\text{\sc Proof.} & \omega(w)-\omega(z)\ &=\ 2(z-w)'+2(z-w)w+(z-w)^2 &\\
&& &=\ (z-w)\cdot\bigl(2(z-w)^\dagger + 2w + (z-w)\bigr)\\
&& &=\ y\cdot\bigl(2y^\dagger+2w+y\bigr)\\
&& &=\ y\cdot\bigl(2(y^\dagger + w)+y\bigr). & \qed
\end{flalign*}

\begin{cor}\label{cor:formula for omega}
Let $y\in K^\times$, $z=-y^\dagger$. Then
$\omega(z+y)=\omega(z)-y^2$.
\end{cor}

\noindent
Membership in $\omega(K)$ governs solvability of certain differential equations:

\begin{cor}
Let $f,g,h\in K$, $f\neq 0$. Then
$$\exists z\in K\, \big[ z'=fz^2+gz+h\big] \quad\Longleftrightarrow\quad
4fh-(g+f^\dagger)^2+2(g+f^\dagger)' \in \omega(K).$$
\end{cor}
\begin{proof}
Let $z$ range over $K$, and set $f_*:=-2f$, $\ g_*:=-2g$, $\ h_*:=-2h$. Then
\begin{align*}
z'\	=\ fz^2+gz+h\  &\Longleftrightarrow\
-2z'\ =\ f_*z^2+g_* z+h_*\\  &\Longleftrightarrow\ 
-2f_*z'\ =\ (f_*z)^2+g_*(f_*z)+f_*h_*\\ 
&\Longleftrightarrow\
-2\big((f_*z)'-f_*'z\big)\ =\ (f_*z)^2+g_*(f_*z)+f_*h_*\\ 
&\Longleftrightarrow 
-2(f_*z)'\	=\ (f_*z)^2 + (g_*-2f_*^\dagger)(f_*z)+f_*h_*.
\end{align*}
So with $g_{**}:=g_*-2f_*^\dagger=-2(g+f^\dagger)$ and $h_{**}:=f_*h_*=4fh$,
$$\exists z \, \big[ z'=fz^2+gz+h\big] \quad\Longleftrightarrow\quad
\exists z \, \big[ -2z'=z^2+g_{**}z+h_{**}\big].$$
Now completing the square yields
$$-2z'=z^2+g_{**}z+h_{**} \quad\Longleftrightarrow\quad
-2\left(z+\textstyle\frac{g_{**}}{2}\right)' = 
\left(z+\textstyle\frac{g_{**}}{2}\right)^2+\left(h_{**}-\textstyle\frac{g_{**}^2}{4}-g_{**}'\right),$$
and 
$h_{**}-\textstyle\frac{g_{**}^2}{4}-g_{**}'=4fh-(g+f^\dagger)^2+2(g+f^\dagger)'$.
\end{proof}

\subsection*{Relation to the Schwarzian} 
We now consider the function
$$y\mapsto \sch(y):=\omega(z+y)=\omega(z)-y^2\ \colon \ K^\times\to K\qquad\text{where $z=-y^\dagger$.}$$
For a nonconstant element $u$ of $K$,
$$\Sch(u)\  := \ (u^{\prime\dagger})' - \frac{1}{2} (u^{\prime\dagger})^2\ =\
\frac{u'''}{u'} - \frac{3}{2}\left(\frac{u''}{u'}\right)^2$$
is known as the {\bf Schwarzian derivative} of $u$. It is related to $\sch(u^\dagger)$ as follows:

\nomenclature[M]{$\sch(y)$}{$\omega(-y^\dagger)-y^2$, for  $y\neq 0$}
\nomenclature[M]{$\Sch(u)$}{Schwarzian derivative of $u$}
\index{derivative!Schwarzian}
\index{Schwarzian derivative}

\begin{lemma}
Let $u\in K\setminus C$. Then
$2\Sch(u) = \sch(u^\dagger)$.
\end{lemma}
\begin{proof}
To see this note that
$u^{\prime\dagger} = (u u^\dagger)^\dagger = u^\dagger+u^{\dagger\dagger}$
and hence
\begin{align*}
(u^{\prime\dagger})'\ &=\ (u^{\dagger})'+(u^{\dagger\dagger})',\\
(u^{\prime\dagger})^2\ &=\  (u^{\dagger})^2+(u^{\dagger\dagger})^2+
2u^\dagger u^{\dagger\dagger}\ =\ (u^{\dagger})^2+(u^{\dagger\dagger})^2+2(u^\dagger)',
\end{align*}
so
\alignqed{
	2\Sch(u)\ 	&=\  2(u^{\prime\dagger})' - (u^{\prime\dagger})^2 \\
				&=\ 2(u^{\dagger})'+2(u^{\dagger\dagger})' - (u^{\dagger})^2-(u^{\dagger\dagger})^2-2(u^\dagger)' \\
				&=\ 2(u^\dagger{}^\dagger)'-(u^\dagger{}^\dagger)^2-(u^\dagger)^2 = \sch(u^\dagger).}
\end{proof}

\noindent
Note that for each $y\in K^\times$ we have  $\sch(-y)=\sch(y)$; in particular,
for $u\in K\setminus C$ we have $\Sch(u)=\Sch(1/u)$. More generally, 
consider the (left) action of the group $\GL_2(C)$ on the set $K\setminus C$ 
given by
$$(T, u)\ \mapsto Tu\ :=\ \frac{ t_{11}u+t_{12} }{ t_{21}u+t_{22}}  
\quad \text{ for }  
T=\left(\begin{smallmatrix} t_{11} & t_{12} \\ t_{21} & t_{22}\end{smallmatrix}\right)\in \GL_2(C).$$ 
Then $\Sch$ is invariant under this action:

\begin{lemma}\label{lem:S(Tu)} Let $T\in \GL_2(C)$ and $u\in K\setminus C$. Then
$\Sch(Tu)\ =\ \Sch(u)$.
\end{lemma}
\begin{proof} Let $T=\left(\begin{smallmatrix} t_{11} & t_{12} \\ t_{21} & t_{22}\end{smallmatrix}\right)$.
Set $u_*:=Tu$, $y:=t_{21}u+t_{22}$, $c:=t_{21}$ (so $y'=cu'$), and $d:=t_{11}t_{22}-t_{12}t_{21}$ (so $d\in C^\times$). Then
$$u_*'\ =\ \frac{d}{y^2} u',\qquad u_*^{\prime\dagger}\ =\ u^{\prime\dagger}-2c \frac{u'}{y},$$
hence
\begin{align*}
(u_*^{\prime\dagger})'\  &=\ (u^{\prime\dagger})' + 2c^2 \left(\frac{u'}{y}\right)^2 - 2c\frac{u''}{y}, \\
(u_*^{\prime\dagger})^2\  &=\ (u^{\prime\dagger})^2 + 4c^2  \left(\frac{u'}{y}\right)^2  - 4c \frac{u''}{y},
\end{align*}
and thus $\Sch(u_*)=\Sch(u)$.
\end{proof}

\begin{lemma}\label{lem:wrscu}
Let $y_1,y_2\in K^\times$ be such that $\wr(y_1,y_2)\in C^\times$. Then 
$$u:=y_1/y_2\notin C, \qquad 2\sch(u^\dagger)\ =\  
 \omega(2y_1^\dagger) + \omega(2y_2^\dagger) .$$
\end{lemma}
\begin{proof} 
Set $c:=\wr(y_1,y_2)$, so $c=y_1y_2'-y_1'y_2=-y_1y_2u^\dagger\in C^\times$ and hence $u\notin C$
and $u^\dagger=-c/(y_1y_2)$. Therefore  $u^\dagger=y_1^\dagger-y_2^\dagger$ and  
$-u^{\dagger\dagger}=y_1^\dagger+y_2^\dagger$, so
\begin{align*}
2\sch(u^\dagger)\			&=\ 
2\omega(-u^{\dagger\dagger})-2(u^\dagger)^2 \\
			&=\  -4\big( (y_1^\dagger)'+(y_2^\dagger)' \big) - 2(y_1^\dagger+y_2^\dagger)^2 - 2(y_1^\dagger-y_2^\dagger)^2 \\
			&=\ -2\big( 2(y_1^\dagger)'+2(y_2^\dagger)' \big) - (2y_1^\dagger)^2 - (2y_2^\dagger)^2, 		
\end{align*}
and the lemma follows.
\end{proof}

\noindent
We can use the function $\sch$ to detect whether $\dim_C \ker A=2$:

\begin{cor}\label{cor:s, 1}
If $y_1,y_2\in\ker A$ are $C$-linearly independent, then $\sch(u^\dagger)=f$ for $u=y_1/y_2\in K\setminus C$.
Conversely, if  $u\in K\setminus C$, $\sch(u^\dagger)=f$, $r\in K^\times$, 
$r^\dagger=\frac{1}{2}u^{\prime\dagger}$, then
$y_1=\frac{u}{r}$ and $y_2=\frac{1}{r}$ form a basis of the $C$-linear space $\ker A$ with $y_1/y_2=u$.
\end{cor}

\noindent
Note that for $u\in K\setminus C$ and $r\in K^\times$, we have: 
$r^\dagger=\frac{1}{2}u^{\prime\dagger}\Longleftrightarrow r^2\in C^\times u'$.

\begin{proof} For the first claim, let $y_1,y_2\in\ker A$ be $C$-linearly 
independent, and set $w:= \wr(y_1, y_2)$, $u:= y_1/y_2$. Then $w'=0$ by 
Lemma~\ref{lem:abel}, 
so $w\in C^\times$, and thus $\sch(u^\dagger)=f$ using Lemma~\ref{lem:wrscu} and
$\omega(2y_1^\dagger)=\omega(2y_2^\dagger)=f$.
For the second claim, let $u\in K\setminus C$, $\sch(u^\dagger)=f$, $r\in K^\times$, $r^\dagger=\frac{1}{2}u^{\prime\dagger}$, and set
$y_1=\frac{u}{r}$ and $y_2=\frac{1}{r}$.
Then 
$$2y_1^\dagger\ =\ 2u^\dagger-u^{\prime\dagger}\ =\ -u^{\dagger\dagger}+u^\dagger$$ and hence
$\omega(2y_1^\dagger)=\sch(u^\dagger)=f$; also
$$2y_2^\dagger\ =\ -u^{\prime\dagger}\ =\ -u^{\dagger\dagger}-u^\dagger\ =\ -(1/u)^{\dagger\dagger}+(1/u)^\dagger,$$ 
so
$\omega(2y_2^\dagger)=\sch\big((1/u)^\dagger\big)=\sch(-u^\dagger)=\sch(u^\dagger)=f$. 
Thus~$A(y_1)=A(y_2)=0$.
\end{proof}

\begin{cor}\label{cor:omega(K)=s(K)}
Suppose $(K^\times)^\dagger=K$. Then 
$$\dim_C\ker A = 2 \ \Longleftrightarrow\  f\in \sch(K^\times).$$
If in addition $K$ is $1$-linearly surjective, then  $\omega(K)=\sch(K^\times)$.
\end{cor}
\begin{proof}
The first part follows from Corollary~\ref{cor:s, 1}. 
Suppose now that~$K$ is also $1$-linearly surjective.
By the remark preceding Lemma~\ref{lem:factoring and fund systems}, if $A$ splits over~$K$, then
$\dim_C \ker A=2$; this gives $\omega(K)=\sch(K^\times)$ by the first part. 
\end{proof}

\begin{cor}\label{cor:s, 2}
Suppose $(K^\times)^\dagger=K$, and
let $u,u_*\in K\setminus C$. Then 
$$\sch(u^\dagger)=\sch(u_*^\dagger)\ \Longleftrightarrow\ \text{$u_*=Tu$ for some $T\in\GL_2(C)$.}$$
\end{cor}
\begin{proof}
The forward direction is a consequence of Corollary~\ref{cor:s, 1}, and the backward direction follows from Lemma~\ref{lem:S(Tu)}.
\end{proof}

\subsection*{The function $\sigma$}
We now assume that $-1$ is not a square in $K$, and
work in the differential field extension~$K[\imag]$ ($\imag^2=-1$) of~$K$.
Note that for $a,b\in K$ we have
$$ A(a+b\imag)\ =\ A(a) + A(b)\imag,$$
so if the equation $A(y)=0$ has a nonzero solution in $K[\imag]$, then
it has a nonzero solution in $K$. Thus
$$\omega\big( 2(K[\imag]^\times)^\dagger \big)\cap K\ =\ \omega\big( 2(K^\times)^\dagger \big),$$
and by Corollary~\ref{cor:omega(K)=s(K)}: 

\begin{cor}\label{d=0or2}
Suppose $K$ is $1$-linearly surjective and $(K^\times)^\dagger=K$, and let $d:=\dim_C\ker A$. Then either 
$d=0$ or $d=2$, and $d$ also equals the  $C[\imag]$-vector space dimension of
the kernel of $A$ viewed as an element of $K[\imag][\der]$.
\end{cor}

\nomenclature[M]{$\sigma(y)$}{$\omega(-y^\dagger)+y^2$, for $y\neq 0$}

\noindent
Consider now also the function $\sigma\colon K^\times \to K$ defined by
$$ \sigma(y)\ :=\ \sch(y\imag)\ =\ \omega(z+y\imag)\ =\ \omega(z)+y^2 \qquad\text{where $z=-y^\dagger$}.$$
If $f=\sigma(y)$ with $y\in K^\times$, and $z=-y^\dagger$, then $A$ factors over $K[\imag]$ as
$$A\ =\ 4\left(\der+\frac{z+y\imag}{2}\right)\left(\der-\frac{z+y\imag}{2}\right).$$
Note also that $\sigma(y)=\sigma(-y)$ for $y\in K^\times$. For $a,b\in K$ we have
$$\omega(a+b\imag)\ =\ \omega(a)+b^2- 2(ab +b')\imag,$$
so if $b\neq 0$, then:  $\omega(a+b\imag)\in K\Leftrightarrow a=-b^\dagger$. Hence
$$\omega\big(K[\imag]\big)\cap K\ =\ \omega(K) \cup \sigma(K^\times),$$
and thus
\begin{equation}\label{eq:A splits over K[i]}
\text{$A$ splits over $K[\imag]$}\ \Longleftrightarrow \ f\in\omega(K)\cup\sigma(K^\times).
\end{equation}

\medskip\noindent
Suppose that 
$A$ splits over~$K[\imag]$. Then $f\in \omega(K)$ or 
$f\in\sigma(K^\times)$.
We now indicate a differential field extension $L$ of $K$
such that $\ker_L A$ has dimension $2$ as a vector space over $C_L$, and
we specify a basis for this vector space.

\case[1]{$f\in\omega(K)$.} Let $z\in K$ be such that
$f=\omega(z)$.  
Take a nonzero element~$y_1$ in a differential field extension of 
$K$ with $2y_1^\dagger=z$, and next
an element~$y_2$ of a differential field extension of $K\<y_1\>=K(y_1)$ with $y_2'-(z/2)y_2=1/y_1$. Then 
$L := K\<y_1, y_2\> = K(y_1, y_2)$
is a differential field extension of $K$ such that $y_1, y_2\in \ker_L A$, and $\wr(y_1, y_2)=y_1y_2'-y_1'y_2=1$. So $L$ is as promised, with
$y_1, y_2$ as a basis of the vector space $\ker_L A$ over $C_L$. 
Note that if $(K^\times)^\dagger=K$ and $K$ is $1$-linearly surjective, then we can take $y_1, y_2\in K$, so that $L=K$. 

\case[2]{$f\in\sigma(K^\times)$.}
Let $y\in K^\times$ satisfy $\sigma(y)=f$. First take an element 
$r$ in a
field extension of $K[i]$ with $r^2=y$, and next an element $e(y)\ne 0$
in a differential field extension of $K[\imag,r]$ such that
$e(y)^\dagger=\frac{1}{2}y\imag $. Thinking of $e(y)$ as $\exp(\frac{1}{2}\imag \int y)$, set  
$$e(-y)\ :=\ e(y)^{-1}, \qquad y_1\ :=\ \frac{e(y)}{r}, \qquad
  y_2\ :=\ \frac{e(-y)}{r}, \qquad L\ :=\ K[\imag]\<y_1, y_2\>.$$
Then $2y_1^\dagger= -y^\dagger +y\imag $ and $2y_2^\dagger= 
-y^\dagger -y\imag $, so $L=K[\imag](y_1, y_2)$ and 
$$\omega(2y_1^\dagger)\ =\ \sigma(y)\ =\ \sigma(-y)\ =\ \omega(2y_2^\dagger)\ =\ f.$$
Also $(y_1/y_2)^\dagger= y_1^\dagger - y_2^\dagger= y\imag \ne 0$, so
$y_1, y_2$ are linearly independent over $C_L$. 
Thus~$L$ is as we promised, with
$y_1, y_2$ as a basis of the vector space $\ker_L A$ over $C_L$.

%We can also take
%$\frac{y_1+y_2}{2}$ (to be thought of as $\frac{\cos\phi}{r}$)
%and $\frac{y_1-y_2}{2i}$
%(to be thought of as $\frac{\sin\phi}{r}$) as $C[\imag]$-linearly %independent zeros of $A$. 
 
\subsection*{Application to linear differential operators with constant coefficients}
\textit{In this subsection we assume that $C$ is real closed, and that $\imag$ is an element of a differential field
extension of $K$ with $\imag^2=-1$.}\/
We have 
$$\omega(C)\ =\ \big\{{-c^2}:\ c\in C\big\}\ =\ C^{\leq},\qquad \sigma(C^\times)\ =\ \big\{c^2:\ c\in C^\times\big\}\ =\ C^>,$$ 
so
$C  = \omega(C)\cup\sigma(C^\times)$.
Hence if $f\in C$, then $A$ splits over $K[\imag]$, by~\eqref{eq:A splits over K[i]}.
For the next lemma, recall the definition of $\operatorname{m}(B)$ for $B\in C[\der]^{\ne}$ from Section~\ref{Linear Differential Operators}.

 \begin{lemma}\label{lem:linconstcoeff, omega} 
 Suppose $K$ is $1$-linearly surjective, $(K^\times)^\dagger=K$, and $\omega(K)\cap C \subseteq C^{\leq}$.
 Let $B\in C[\der]^{\neq}$. Then  
 $\dim_C\ker B = \operatorname{m}(B)$.
 \end{lemma}
 \begin{proof}
We have $\operatorname{m}(B)\leq\dim_C\ker B$ by~\eqref{eq:mledim}. To show
 $\dim_C\ker B \leq \operatorname{m}(B)$,
we can assume $B$ is monic of order $\ge 1$.  
Then Proposition~\ref{prop:AS}(i) yields $B=B_1\cdots B_n$ ($n\geq 1$), 
 with
  $B_i\in C[\der]$ monic of order~$1$ or $2$ for $i=1,\dots,n$. Hence
  $$\operatorname{m}(B)\ =\ \operatorname{m}(B_1)+\cdots+\operatorname{m}(B_n), \quad  
  \dim_C \ker B\ \leq\ \dim_C \ker B_1 + \cdots + \dim_C \ker B_n.$$
This gives a reduction to the case that $B$ has order~$1$ or~$2$. 
  If $\order(B)=1$ then $\ker B\neq\{0\}$ (since $(K^\times)^\dagger=K$), so $\dim_C\ker B=1=\operatorname{m}(B)$.
  Suppose $\order(B)=2$, say
   $B=\der^2+b\der+c$ with $b,c\in C$. Put $\tilde B:=\der^2+(c-\frac{b^2}{4})$;
  then $\dim_C \ker B = \dim_C \ker \tilde B$; see Lemma~\ref{special form} and the example following it.
 Also  
  $\operatorname{m}(B)=\operatorname{m}({\tilde B})$.  So replacing~$B$ by~$\tilde B$ we arrange $b=0$. Then $\operatorname{m}(B) = 0$ if $c > 0$ and $\operatorname{m}(B) = 2$ if $c\leq 0$. By Corollary~\ref{d=0or2} applied to $A=4B$ we either have
  $\dim_C \ker B=0$ or $\dim_C \ker B=2$.
From $(K^\times)^\dagger=K$ and $\omega(K)\cap C = C^{\leq}$, we get
     $$\dim_C\ker B\geq 1 \quad \Longleftrightarrow \quad 4c\in\omega(K)\quad \Longleftrightarrow\quad c\leq 0,$$
and this yields $\dim_C\ker B=\operatorname{m}(B)$.
  \end{proof}

\subsection*{Notes and comments}
The Schwarzian derivative plays a role in the analytic theory of linear differential equations, where versions of
Corollaries~\ref{cor:s, 1} and \ref{cor:s, 2} are well-known; see \cite[Chapter~10]{Hille}.

%% file: mt-5-3.tex
\section{Diagonalization of Matrices}\label{sec:diagonalization}

\noindent
Given a differential field $K$, the ring $K[\der]$ is euclidean
as defined below. Here we establish a basic result about 
matrices over any euclidean ring. In the next section we use this
to reduce
a system of linear differential equations in several unknowns 
to several linear differential equations in one unknown.
Throughout this section $R$ is a (possibly non-commutative) ring with $1\neq 0$
and we let $a$,~$b$,~$c$ range over $R$.
We say that $R$ is a {\bf domain} 
if for all $a$,~$b$, if $ab=0$, then $a=0$ or $b=0$. 
 
\index{domain}
\index{ring!domain}

\subsection*{Total divisibility}
We say that  
{\bf $a$ totally divides~$b$} (notation: $a\|b$) if
$Rb\subseteq aR$ and $bR\subseteq Ra$ (and thus
$RbR\subseteq aR\cap Ra$). Note:
\begin{enumerate}
\item $a\|0$,
\item $0\|b\Longleftrightarrow b=0$,
\item $a\|a\Longleftrightarrow aR=Ra$, 
\item $a\|b\ \&\ b\|c \Rightarrow a\|c$, 
\item if $a$ is a unit of $R$, then $a\|b$ for all $b$, and
\item given $a$ the set $\{b: a\|b\}$ is a (two-sided) ideal of $R$.
\end{enumerate}
The ring $R$ is said to be {\bf simple} if the only (two-sided) ideals of $R$ are $\{0\}$ and $R$.
If $R$ is simple and $a\|b$, then $b=0$ or $a\in R^{\times}$. 

\nomenclature[A]{${a\lvert\rvert b}$}{$a$ totally divides $b$}
\index{ring!simple}
\index{simple!ring}

\subsection*{Degree functions} Addition on $\N$ and the usual ordering
on $\N$ are extended as usual to $\N\cup\{-\infty\}$.
A {\bf degree function} on $R$ is a map 
$\operatorname{d}\colon R\to\N\cup\{-\infty\}$ such that for all $a$,~$b$,
\begin{list}{*}{\setlength\leftmargin{2.5em}}
\item[(D1)] $\operatorname{d}(a)=-\infty \Longleftrightarrow a=0$;
\item[(D2)] $\operatorname{d}(-a)=\operatorname{d}(a)$;
\item[(D3)] $\operatorname{d}(a+b)\leq \max\{ \operatorname{d}(a),\operatorname{d}(b) \}$; and
\item[(D4)] $\operatorname{d}(ab)=\operatorname{d}(a)+\operatorname{d}(b)$.
\end{list}
%Here we extend the addition on $\N$ to a binary operation on 
%$\N\cup\{-\infty\}$ by $n+(-\infty)=(-\infty)+n
%=(-\infty)+(-\infty):=-\infty$ for all $n$, and we 
%extend the usual ordering on $\N$ to an ordering on 
%$\N\cup\{-\infty\}$ such that $-\infty<n$ for all $n$.
Let $\operatorname{d}$ be a degree function on $R$. Then $\operatorname{d}(1)=0$ by (D1) and (D4). 
It follows from (D1), (D2), (D3) that
$\operatorname{d}(a+b) = \max\!\big\{\!\operatorname{d}(a),\operatorname{d}(b) \big\}$  if $\operatorname{d}(a)\neq \operatorname{d}(b)$.
By (D1) and (D4) the set $R^{\neq}$ is closed under multiplication, so
$R$ is a domain.

\index{degree function!on a ring}
\index{ring!degree function}

\subsection*{Euclidean rings}
Let $\operatorname{d}$ be a degree function on $R$. Note that
then for all $a,b$ with~$a\neq 0$, there is at most one pair 
$(q,r)\in R^2$ with $b=qa+r$ and $\operatorname{d}(r)<\operatorname{d}(a)$, and also
at most one pair~$(q^*,r^*)\in R^2$ with $b=aq^*+r^*$ and $\operatorname{d}(r^*)<\operatorname{d}(a)$.
We say that $R$ is {\bf left euclidean with respect to $\operatorname{d}$} if
for all $a,b$, $a\neq 0$, there is $q\in R$ with $\operatorname{d}(b-qa)<\operatorname{d}(a)$; similarly,~$R$ is {\bf right euclidean with respect to $\operatorname{d}$} if
for all $a,b$, $a\neq 0$, there is $q^*\in R$ with $\operatorname{d}(b-aq^*)<\operatorname{d}(a)$. We  say that $R$ is {\bf euclidean with respect to $\operatorname{d}$} if it is both left and right euclidean with respect to~$\operatorname{d}$.
Note that if $R$ is euclidean with respect to $\operatorname{d}$, then $R^\times=\big\{a:\operatorname{d}(a)=0\big\}$. 

\index{ring!euclidean}
\index{euclidean!ring}

If $R$ is a domain, we call $a$ {\bf irreducible} \index{irreducible!element of a domain} (in $R$) if 
$a\notin R^\times$,
and there are no $a_1,a_2\in R\setminus R^\times$ with $a=a_1a_2$. Note that if
$a$ is irreducible in the domain $R$, then $a\ne 0$, and $au$,~$ua$ are also 
irreducible for $u\in R^{\times}$.  
Call $R$
{\bf euclidean} if it is euclidean with respect to some degree function on $R$.
If $R$ is euclidean, then $R$ is a domain, every nonzero element of $R\setminus R^\times$ equals
$a_1\cdots a_n$ for some $n\ge 1$ and irreducible $a_1,\dots, a_n$ in $R$,
%and $R$ is a PID (principal ideal domain), that is, $R$ is a
%domain, 
every left ideal of $R$ is principal (of the form $Ra$), and every right ideal of $R$ is principal (of the form $aR$). 

\index{principal!ideal}

\begin{lemma}\label{lem:submodule of free}
Suppose every left ideal of $R$ is principal. Then every submodule~$M$ of the left $R$-module $R^n$ is generated by $n$ elements.
\end{lemma} 
\begin{proof}
The case $n=0$ holds trivially, so suppose $n\geq 1$, and identify $R^{n-1}$ with $R^{n-1}\times\{0\}\subseteq R^n$ in
the natural way. Inductively, the module 
$M\cap R^{n-1}$ is generated by elements $b_1,\dots,b_{n-1}$. Let $e:=(0,\dots,0,1)\in R^n$, and
consider the left ideal $I:=\{a\in R:ae\in M+R^{n-1}\}$ of $R$. Let $a\in R$ with $I=Ra$, and pick
$b_n\in M$ with $ae\in b_n+R^{n-1}$.
Then $b_1,\dots,b_n$ generate $M$.
\end{proof}

\subsection*{Ore domains}
A {\bf right Ore domain} is a domain $R$ such that
for all $a,b\in R^{\neq}$ we have $aR\cap bR\neq \{0\}$. A {\bf left Ore
domain} is defined similarly. If $R$ is both a left and right Ore domain,
then~$R$ is called an {\bf Ore domain}. \index{domain!Ore}\index{Ore domain}\index{ring!Ore domain}
 (For example, every integral domain is an
Ore domain.) Call~$R$ {\bf right noetherian} if every right ideal of~$R$ is finitely generated (as a right $R$-module);
likewise with \textit{left}\/ instead of \textit{right.}\/\index{ring!noetherian}
%For use in Section~\ref{sec:diff modules, model theory} we note:  

\begin{lemma}\label{lem:Ore} Let $R$ be a right noetherian domain.
Then $R$ is a right Ore domain.
\end{lemma}
\begin{proof}
Let $a,b\in R^{\ne}$. Since the right ideal of $R$ generated by the $a^mb$
is finitely generated, we have $n$ such that $a^{n+1}b= \sum_{i=0}^n a^ibc_i$
with $c_0,\dots, c_n\in R$. By canceling some power $a^m$ with $m\le n$
on both sides
we arrange $c_0\ne 0$. Then 
$$0\ \ne\ bc_0\ =\ a^{n+1}b-\sum_{i=1}^na^ibc_i\in aR \cap bR,$$
as desired.\end{proof}

\noindent
Similarly, every left noetherian domain is left Ore. In particular, 
every euclidean ring is an Ore domain.

\medskip\noindent
Now let $R$ be a left Ore domain and 
$M$ an $R$-module.
Then  $$M_{\operatorname{tor}}\ :=\ \big\{x\in M: \text{$ax=0$ for some $a\in R^{\neq}$}\big\}$$
is a submodule of $M$: for $x,y\in M_{\operatorname{tor}}$, take $a,b\in R^{\neq}$ with $ax=by=0$; take $r,s,t\in R$ with $ra=sb=t\neq 0$; then
$t(x+y)=rax+sby=0$, so $x+y\in M_{\operatorname{tor}}$; it follows likewise that if $x\in M_{\operatorname{tor}}$ and $a\in R$, then $ax\in M_{\operatorname{tor}}$.
We call $M_{\operatorname{tor}}$ the {\bf torsion submodule} of $M$.
Call $M$ a {\bf torsion module} if $M_{\operatorname{tor}}=M$ and
{\bf torsion-free} if $M_{\operatorname{tor}}=\{0\}$. So  $M_{\operatorname{tor}}$ is 
a torsion module and
$M/M_{\operatorname{tor}}$ is torsion-free. 

\index{module!torsion submodule}
\index{module!torsion module}
\index{module!torsion-free}
\index{torsion!submodule}
\index{torsion!module}
\nomenclature[A]{$M_{\operatorname{tor}}$}{torsion submodule of $M$}

\subsection*{Diagonalization} 
\textit{In this subsection $R$ is euclidean and $m,n\geq 1$.}\/ When $m$,~$n$ are clear from the context, we let $0$ denote the
zero element of the additive group~$R^{m\times n}$ of $m\times n$ matrices 
over $R$. For $m=n$ we denote the
multiplicative group of units of the ring $R^{n\times n}$ by 
$\operatorname{GL}_n(R)$. 
An $m\times n$ matrix $A=(a_{ij})$ over $R$ is said to be {\bf diagonal} if $a_{ij}=0$ for all $i\in\{1,\dots,m\}$ and
$j\in\{1,\dots,n\}$ with $i\neq j$. Given $m\times n$ matrices $A$ and $B$ over $R$, we say that $A$ and $B$ are
{\bf equivalent} (over $R$) if there are $P\in \operatorname{GL}_m(R)$ and $Q\in \operatorname{GL}_n(R)$ 
with $A=PBQ$; in symbols: $A\sim B$. Clearly $\sim$ is an equivalence relation on the set of $m\times n$ matrices  over $R$, with the equivalence 
class of $0$ being $\{0\}$.
We are going to show that every matrix  over a euclidean ring is equivalent to a diagonal matrix; more precisely:

\nomenclature[A]{$\operatorname{GL}_n(R)$}{group of invertible $n\times n$ matrices over a ring $R$}
\index{matrix!diagonal}
\index{diagonal!matrix}
\index{equivalence!matrices}

\begin{theorem}\label{thm:smith nf}
Every $m\times n$ matrix over $R$ is equivalent to a diagonal $m\times n$ matrix $D=(d_{ij})$ over $R$  such that 
$d_{ii}\|d_{jj}$ for $1\leq i< j\leq \min\{m,n\}$.
\end{theorem}

\noindent
The proof involves row operations on an $m\times n$ matrix $A$ over~$R$:
\begin{list}{*}{\setlength\leftmargin{2.5em}}
\item[(R1)] interchange two rows;
\item[(R2)] add a left multiple of the $i$th row to the $j$th row ($i\ne j$);
\item[(R3)] multiply a row on the left by a unit of $R$.
\end{list} 
If $B$ arises from $A$ by applying one of the operations (R1)--(R3), 
then  $B=PA$ for some  $P\in \operatorname{GL}_m(R)$, so $A\sim B$. 
We also have column operations on $A$:
\begin{list}{*}{\setlength\leftmargin{2.5em}}
\item[(C1)] interchange two columns;
\item[(C2)] add a right multiple of the $i$th column to the $j$th column ($i\ne j$);
\item[(C3)] multiply a column on the right by a unit of $R$.
\end{list} 
If $B$ arises from $A$ by applying one of the operations (C1)--(C3), then  
$B=AQ$
for some  $Q\in \operatorname{GL}_n(R)$, hence $A\sim B$.

\medskip
\noindent
We fix a degree function $\operatorname{d}$ on $R$ with respect to which $R$ is euclidean. Given an $m\times n$ matrix $A=(a_{ij})$ over $R$, let
$\operatorname{d}(A)$ be the minimum of the $\operatorname{d}(a_{ij})$ with 
$a_{ij}\ne 0$ if $A\ne 0$, and $d(A):= -\infty$ if $A=0$.

\begin{lemma}\label{lem:min 2x2}
Let $A=\left(\begin{smallmatrix} a & 0 \\ 0 & b \end{smallmatrix}\right)$ where $a,b\in R$, $a\neq 0$,
and suppose
$\operatorname{d}(a)\le \operatorname{d}(B)$ for all $2\times 2$ matrices
$B$ over $R$ with $A\sim B$. Then $a\|b$.
\end{lemma}
\begin{proof}
Let $c\in R$, and take $q,r\in R$ with $bc=qa+r$ and $\operatorname{d}(r)<\operatorname{d}(a)$.
First applying the operation (R2) and then (C2) we see that $A\sim \left(\begin{smallmatrix} a & 0 \\ r & b \end{smallmatrix}\right)$, hence $r=0$. This shows $bR\subseteq Ra$. Similarly one obtains $Rb\subseteq aR$.
\end{proof}

\noindent
Towards the proof of Theorem~\ref{thm:smith nf}, consider an $m\times n$ 
matrix $A$ over $R$.

\smallskip

\claim{$A$ is equivalent to a diagonal matrix $D$ over $R$ with 
$\operatorname{d}(D)\le \operatorname{d}(A)$.} 

\begin{proof}[Proof of Claim] Set $k:=\min\{m,n\}\in\N^{\geq 1}$ and $d:=\operatorname{d}(A)\in\N\cup\{-\infty\}$. We order the set~${\N^{\geq 1}\times(\N\cup\{-\infty\})}$ lexicographically and prove
the claim by induction on $(k,d)$.
If $A=0$, then  the claim holds trivially, so assume $A\neq 0$.
Applying the operations (R1) and (C1) we first replace $A$ by an equivalent $m\times n$ matrix, without changing~$(k,d)$,
to reduce to the case that $\operatorname{d}(a_{11})=d$.
Now using (R2) and (C2) and euclidean division by $a_{11}$, we see that $A\sim B$ where $B=(b_{ij})$ is an $m\times n$ matrix over $R$
with $b_{11}=a_{11}$ (hence $\operatorname{d}(B)\leq d$) and $\operatorname{d}(b_{i1})<\operatorname{d}(a_{11})$ for $i=2,\dots,m$
and $\operatorname{d}(b_{1j})<\operatorname{d}(a_{11})$ for $j=2,\dots,n$.
If $m>1$ and $b_{i1}\neq 0$ for some $i\in\{2,\dots,m\}$, then $\operatorname{d}(B)<\operatorname{d}(a_{11})=d$,
and we can apply the inductive hypothesis to $B$.
Similarly, if $n>1$ and $b_{1j}\neq 0$ for some $j\in\{2,\dots,n\}$, then $\operatorname{d}(B)<\operatorname{d}(A)$,
and the inductive hypothesis applies to $B$.
Thus we may assume $b_{i1} = 0$ for all $i\in\{2,\dots,m\}$ and $b_{1j} = 0$ for all $j\in\{2,\dots,n\}$.
If $k=1$, then $B$ is already diagonal, so suppose~$k>1$. Then we have an  $(m-1)\times (n-1)$ matrix $B'$ over $R$ such that $$B\ =\
\left(
\begin{array}{cc}
a_{11}	&   \\[-0.5em] 
		& \boxed{ \begin{matrix} & & \\ & B' & \\ & & \end{matrix}}\end{array}\right).
$$
The claim now follows by applying the inductive hypothesis to $B'$. \end{proof} 

\begin{proof}[Proof of Theorem~\ref{thm:smith nf}]
Let $A$ be an $m\times n$ matrix over $R$.
To show that $A$ is equivalent over~$R$ to a diagonal matrix as in Theorem~\ref{thm:smith nf}, we can assume $A\neq 0$.
The claim gives a diagonal matrix $D=(d_{ij})$ over $R$ with $A\sim D$ and
$\operatorname{d}(D)\le \operatorname{d}(B)$ for all
$m\times n$ matrices $B$ over $R$ with $A\sim B$.
By applying the operations~(R1) and~(C1) to~$D$ we arrange
$\operatorname{d}(D)=\operatorname{d}(d_{11})$. Set $k:=\min\{m,n\}$.
We are done if $k=1$, so assume~$k>1$.
By Lemma~\ref{lem:min 2x2} and the minimality of $\operatorname{d}(D)$ we have
$d_{11}\|d_{ii}$ for $i=2,\dots,k$. We 
can assume inductively that the $(m-1)\times (n-1)$-matrix
$(d_{ij})_{i,j\ge 2}$ is equivalent to a diagonal matrix
as in  Theorem~\ref{thm:smith nf}, the entries
of which will then be in the (two-sided) ideal of $R$ generated by the 
$d_{ii}$ with
$i=2,\dots,k$, and so $d_{11}$ totally divides each of those entries. \end{proof}

\begin{cor}\label{cor:smith nf, 1}
Suppose $R$ is simple. Let $A\neq 0$ be an $m\times n$ matrix over~$R$, and
$k:=\min\{m,n\}$. Then $A$ is equivalent to a
diagonal $m\times n$ matrix $D=(d_{ij})$ over $R$ such that for some
$r\in\{1,\dots,k\}$ we have
$$ d_{11}\ =\ \cdots\ =\ d_{r-1,r-1}\ =\ 1,\qquad d_{rr}\neq 0, \quad d_{r+1,r+1}\ =\ \cdots\ =\ d_{kk}\ =\ 0.$$
\end{cor}
\begin{proof}
Theorem~\ref{thm:smith nf} gives a
diagonal $m\times n$ matrix $B=(b_{ij})$  over $R$ with $A\sim B$ and 
$b_{ii}\|b_{jj}$ for $1\leq i< j\leq k$. Since $A\neq 0$, we have $b_{rr}\neq 0$ for some $r\in\{1,\dots,k\}$;
take $r$ maximal with this property. For $i\in\{1,\dots,r-1\}$ we have $b_{ii}\|b_{rr}$, and so $b_{ii}$ is a unit of $R$ (as $R$ is simple). Now use (R3).
\end{proof}

\noindent
The next result is an application of Theorem~\ref{thm:smith nf} to finitely generated modules over euclidean rings. Here and below, 
\textit{module}\/ means \textit{left module.}\/

\begin{cor}\label{cor:smith nf, 2}
Let $M$ be a finitely generated  $R$-module. Then
$$M\cong (R/Rd_1)\oplus\cdots\oplus (R/Rd_r)\oplus R^s$$
where $d_1,\dots,d_r\in R^{\neq}$, $r,s\in\N$, and $d_1\|d_2\|\cdots\|d_r$.
If $R$ is simple, then
$$M\cong (R/Rd)\oplus R^s\  \qquad (d\in R^{\neq},\ s\in\N).$$
\end{cor}
\begin{proof}
Take $n\geq 1$ and a surjective $R$-linear map $R^n\to M$.
If the kernel is trivial, then the above holds with $r=0$,~$s=n$ (and $d=1$
for simple~$R$).
Assume the kernel is not trivial. Then Lemma~\ref{lem:submodule of free}
yields an $m\geq 1$ and an $m\times n$ matrix~$A$ over~$R$ such that 
$M\cong R^n/N$
where $N=\{yA:y\in R^m\}$ and the elements of $R^m$ and $R^n$
are viewed as row vectors.
Applying Theorem~\ref{thm:smith nf} and  Corollary~\ref{cor:smith nf, 1} to~$A$
yields the desired result with $r+s=n$.
\end{proof}

\noindent
Hence a finitely generated $R$-module $M$ is torsion-free iff $M\cong R^n$ for some $n$.

\subsection*{Independence and rank} Let $R$ be a euclidean domain, and let $M$ range over $R$-modules.
We refer to the Conventions and Notations in the beginning of this volume for the terminology in dealing with linear (in)dependence over $R$.

\medskip\noindent
Let $m,n\geq 1$ and let $A=(a_{ij})$ be an $m\times n$ matrix over $R$. \index{independent} \index{dependent} \index{R-independent@$R$-independent} \index{R-dependent@$R$-dependent} We say that the rows of~$A$ are $R$-independent if $A_1,\dots,A_m$ are $R$-independent, where
$A_i=(a_{i1},\dots,a_{in})$ is the $i$th row, considered as a vector of the (left) 
$R$-module $R^n$. For $r_1,\dots, r_m\in R$ we may view $(r_1,\dots, r_m)$ as a $1\times m$ matrix over $R$, and so
$$r_1A_1+\cdots + r_mA_m\ =\ (r_1,\dots, r_m)A.  $$
Therefore, if the rows of $A$ are 
$R$-independent and $A\sim B$ for the $m\times n$ matrix $B$ over $R$,
then the rows of $B$ are $R$-independent. Suppose $A\sim D$,
where $D=(d_{ij})$ is a  diagonal $m\times n$ matrix over $R$. Then
the rows of $A$ are $R$-independent iff
$m\le n$ and $d_{ii}\ne 0$ for $i=1,\dots,m$. 
For $R$-linear $f\colon R^m \to R^n$ this yields:

\begin{cor}\label{cor:rank 1}
If $f$ is injective, then~$m\leq n$.
\end{cor}
\begin{proof} Let $e_1,\dots, e_m$ be the usual basis vectors of $R^m$. Assume $f$ is injective. Let $A$ be the
$m\times n$ matrix with $i$th row $A_i:= f(e_i)$. Then the rows
of $A$ are $R$-independent, so $m\le n$ by 
Theorem~\ref{thm:smith nf} and the remarks above. 
\end{proof}

\begin{cor}\label{cor:rank 2}
If $f$ is surjective, then~$m\geq n$.
\end{cor}
\begin{proof}  Let $e_1,\dots, e_n$ be the usual basis vectors of $R^n$. Assume $f$ is surjective. Take $b_j\in R^m$ with
$f(b_j)=e_j$ for $j=1,\dots,n$. Then $R^n\cong Rb_1+\cdots + Rb_n\subseteq R^m$, so $n\le m$ by Corollary~\ref{cor:rank 1}. 
\end{proof}

\index{rank!module}
\index{module!rank}

\noindent
Of course these two corollaries also hold for $m=0$ or $n=0$,
which we allow below. By these two corollaries, if $R^m\cong R^n$ as $R$-modules, then $m=n$.
So for each finitely generated free $R$-module $M$ we may define the {\bf rank} of $M$ to be the unique~$n$ such that $M\cong R^n$; in this case,
\begin{align*}
\operatorname{rank} M\	&=\ \text{largest $m$ such that there are $R$-independent $x_1,\dots, x_m\in M$} \\
						&=\ \text{least $m$ such that $M$ is generated by some $x_1,\dots, x_m\in M$}.
\end{align*}
Now let $M$ be any finitely generated $R$-module. Then the finitely generated torsion-free $R$-module $M/M_{\operatorname{tor}}$
is free, and we set $\operatorname{rank}(M):=\operatorname{rank}(M/M_{\operatorname{tor}})$;
thus~$M$ is a torsion module iff $\operatorname{rank}(M)=0$. Also,
$$\operatorname{rank} M\ 	=\  \text{largest $m$ such that there are $R$-independent $x_1,\dots, x_m\in M$}.$$
Clearly if $M$, $N$ are finitely generated $R$-modules, then
$$\operatorname{rank}(M\oplus N)\ =\ \operatorname{rank}(M)+\operatorname{rank}(N).$$
In fact, the next lemma shows that $M \mapsto \operatorname{rank}(M)$, for finitely generated $M$, is a $\Z$-valued Euler-Poincar\'e map on $R$-modules in the sense of Section~\ref{sec:finite length}.

\begin{lemma}\label{lem:additivity of rank}
Let an exact sequence of $R$-modules and $R$-linear maps be given:  
$$0\longrightarrow K\xrightarrow{\ \iota\ } M\xrightarrow{\ \pi\ } N\longrightarrow 0.$$
If $M$ is finitely generated, then so are $K$ and~$N$, and
$\operatorname{rank}(M) = \operatorname{rank}(K) + \operatorname{rank}(N)$.
\end{lemma}
\begin{proof}
Suppose $M$ is finitely generated. Then $N$ is finitely generated, and by Lem\-ma~\ref{lem:submodule of free}, so is $K$.
To prove the rank formula, assume until further notice 
that~$N$ is a torsion module. Suppose $x_1,\dots, x_m\in M$ are
$R$-independent. We have $r_1,\dots, r_m\in R^{\ne}$ with 
$r_1x_1,\dots, r_mx_m\in \ker \pi= \iota(K)$, and as
$r_1x_1,\dots, r_mx_m$ are also $R$-independent, we get
$m\le \operatorname{rank}(K)$. This shows $\operatorname{rank}(M)\le \operatorname{rank}(K)$. It is obvious that 
$\operatorname{rank}(K)\le \operatorname{rank}(M)$, and so
the two ranks are equal. 
%The $R$-modules
%$\overline{K}:=K/K_{\operatorname{tor}}$ and  
%$\overline{M}:=M/M_{\operatorname{tor}}$ 
%are finitely generated and free, with 
%$m:=\operatorname{rank}(\overline{K})=\operatorname{rank}(K)$ %and
%$n:=\operatorname{rank}(\overline{M})=\operatorname{rank}(M)$, %and 
%$\overline{N}:=N/\pi(M_{\operatorname{tor}})$ is a torsion %module.
%We have $\iota(K_{\operatorname{tor}}) = %M_{\operatorname{tor}}\cap \iota(K)$, hence~$\iota$ induces an
%injective $R$-linear map $\overline{\iota}\colon \overline{K}  %\to \overline{M}$. Hence $m\le n$. Let 
%$\overline{\pi}\colon \overline{M} \to \overline{N}$ be the %surjective
%$R$-linear map induced by $\pi$. Then we get the exact %sequence
%$$0\longrightarrow \overline{K}
%\xrightarrow{\ \overline{\iota}\ } 
%\overline{M}\xrightarrow{\ \overline{\pi}\ } 
%\overline{N}\longrightarrow 0$$
%{\bf proof checked sofar} Theorem~\ref{thm:smith nf} 
%provides a basis $b_1,\dots, b_m$ of $\overline{K}$,
%a basis $e_1,\dots, e_n$ of $\overline{M}$, and 
%$r_1,\dots, r_m\in R$ such that 
%$\overline{\iota}(b_i)=r_i e_i$ for $i=1,\dots,m$.
%It follows that $m=n$, that is,
%$\operatorname{rank}(\overline{K})=
%\operatorname{rank}(\overline{M})$.

We now drop the assumption that $N$ is a torsion module. The canonical map
$\nu\colon N \to N/N_{\operatorname{tor}}$ yields the $R$-linear surjection  
$\nu\circ \pi\colon M\to N/N_{\operatorname{tor}}$ with kernel
$M_1:=\pi^{-1}(N_{\operatorname{tor}})$. Since 
the $R$-module $N/N_{\operatorname{tor}}$ is free, we have
$M\cong M_1\oplus (N/N_{\operatorname{tor}})$, so
$\operatorname{rank}(M)=\operatorname{rank}(M_1)+\operatorname{rank}(N)$.
We have an exact sequence
$$0\longrightarrow K\xrightarrow{\ \iota\ } M_1\xrightarrow{\ \pi|M_1\ } N_{\operatorname{tor}}\longrightarrow 0$$
and so $\operatorname{rank}(K)=\operatorname{rank}(M_1)$ by what we showed before.
\end{proof}

\begin{cor}\label{cor:additivity of rank} Suppose $M$ is generated by elements $x_1,\dots, x_n$, and $m\le n$ is such that $x_1,\dots, x_m$ are $R$-independent
and $x_1,\dots,x_m, x_j$ are $R$-dependent for all~$j$ with
$m < j \le n$. Then $\operatorname{rank}(M)=m$.
\end{cor}
\begin{proof} Let $K:= Rx_1+\cdots + Rx_m$. Then $M/K$ is a torsion module and $K\cong R^m$,
so $\operatorname{rank}(K)=m$. Now use
Lemma~\ref{lem:additivity of rank}.
\end{proof}

\subsection*{Notes and comments}
Lemma~\ref{lem:Ore} is due to Goldie~\cite{Goldie}.
Theorem~\ref{thm:smith nf} for $R=\Z$ was proved by Smith~\cite{Smith}. The version here, for euclidean $R$, is due
to Wedderburn~\cite{Wedderburn} and Jacobson~\cite{Jacobson};
our presentation follows~\cite[Section~1.4]{CohnIdeal}. 
%The theorem goes through
%for $R$ a PID; see Teichm\"uller~\cite{Teichmuller}. %Corollaries~\ref{cor:rank 1} and ~\ref{cor:rank 2} hold even %for left noetherian $R$; see ~\cite[\S{}1]{Lam}. 

%% file: mt-5-4.tex
\section{Systems of Linear Differential Equations}\label{sec:systems}

\noindent
In this section $K$ is a differential field.
We apply the previous section to $R=K[\der]$.
The ring $K[\der]$ is euclidean with respect to the degree function $\operatorname{d}$ on 
$K[\der]$ given by $\operatorname{d}(A):=\order(A)$.
% if $A\neq 0$ and $\operatorname{d}(0):=-\infty$.
By Corollary~\ref{cor:K[der] is simple}, if $C\neq K$, then $K[\der]$ is simple.

\subsection*{Inhomogeneous equations}
In this subsection, $R=K[\der]$.
Let $m,n\ge 1$, and let~$A$ be an $m\times n$ matrix over $R$. Given a column vector $f=(f_1,\dots,f_n)^{\operatorname{t}}\in K^n$ 
we let $A(f)$ be the column vector in~$K^m$ with $i$th entry 
$$\sum_{j=1}^n A_{ij}(f_j)\quad  (i=1,\dots,m),$$ not to be confused with the
matrix product $Af\in R^m$ that has $i$th entry 
$\sum_{j=1}^n A_{ij}f_j$. The map $f\mapsto A(f)\colon K^n \to K^m$ is $C$-linear.
It is easy to check that if $B$ is an $n\times p$ matrix over $R$ with $p\in \N^{\ge 1}$, and
$g=(g_1,\dots, g_p)^{\operatorname{t}}\in K^p$, then $(AB)(g)=A\big(B(g)\big)$.
For the $n\times n$ identity matrix $I$ over $R$ we have $I(f)=f$
for all $f\in K^n$. 

\medskip\noindent 
Let $a=(a_1,\dots,a_m)^{\operatorname{t}}\in K^m$ be a column vector.
The pair $(A,a)$ gives rise to a system $A(y)=a$ of linear differential equations over $K$.
A {\bf solution} in $K$ to this system 
is a column vector $f=(f_1,\dots,f_n)^{\operatorname{t}}\in K^n$ such that
$A(f)=a$. We say that the system $A(y)=a$ is {\bf $R$-consistent} if
for every row vector $r=(r_1,\dots, r_m)\in R^m$ such that $rA=0$ (a matrix product)
we have $r(a)=0$, where
$$r(a)\ :=\ r_1(a_1) + \cdots + r_m(a_m).$$ 
If $A(y)=a$ has a solution in $K$, then $A(y)=a$ is clearly $R$-consistent. 

\index{system of linear differential equations!$R$-consistent}
\index{system of linear differential equations!solution}

\medskip\noindent
We can increase $m$ while keeping $n$ fixed by adding extra zero rows to $A$ and extra zero entries to $a$; in this way we can arrange
that $m\ge n$. Such a change in $A$, $a$ does not change the solutions to the system
in $K$ nor its $R$-consistency status. So we assume below that $m\ge n$. 

\medskip\noindent
Using Theorem~\ref{thm:smith nf}, take $P\in \operatorname{GL}_m(R)$ and $Q\in \operatorname{GL}_n(R)$ such that $PAQ$ is diagonal, and set $b=P(a)\in K^m$. Then each column $z\in K^n$ with
$(PAQ)(z)=b$ yields a column $y=Q(z)\in K^n$ with $A(y)=a$; this gives a
bijective correspondence $z\mapsto Q(z)$ between the solutions of $(PAQ)(z)=b$
in $K$ and the solutions of $A(y)=a$ in $K$. 
 It is also easy to see that
$A(y)=a$ is $R$-consistent iff $(PAQ)(z)=b$ is $R$-consistent.   
Put $B_i=(PAQ)_{ii}\in R$ for $i=1,\dots, n$.
Then the system $(PAQ)(z)=b$ with $z=(z_1,\dots, z_n)^{\operatorname{t}}\in K^n$
takes the form
$$B_1(z_1)\ =\ b_1, \dots, B_n(z_n)\ =\ b_n,\qquad 0\ =\ b_i\ \text{ for $n<i\le m$,}$$
that is, $(PAQ)(z)=b$ has no solution in $K$ if $b_i\ne 0$ for some $i$ with
$n<i\le m$, and otherwise its solutions are the columns $z\in K^n$ such that
$$B_1(z_1)\ =\ b_1, \dots, B_n(z_n)\ =\ b_n.$$ 
The system
$(PAQ)(z)=b$ is $R$-consistent if and only if for all $r_1,\dots, r_m\in R$ with
$r_1B_1=\dots = r_nB_n=0$ we have $r_1(b_1) + \cdots + r_m(b_m)=0$, which in turn is equivalent to $b_i=0$ for all $i\in \{1,\dots,n\}$ with $B_i=0$, and $b_i=0$ for all $i\in\{n+1,\dots,m\}$.  Thus:

\begin{lemma}\label{consistent}
If $K$ is linearly surjective and a system $A(y)=a$ as above is $R$-con\-sis\-tent,
then it has a solution in $K$. 
\end{lemma}

\noindent
Next we  consider a more common way of presenting a system of linear
differential equations, namely as a matrix differential equation $y'=Ay+b$ where
$A$ is an $n\times n$ matrix over $K$ (not over $K[\der]$ as above), with $n\ge 1$, and
$b=(b_1,\dots, b_n)^{\operatorname{t}}\in K^n$ is a given column vector. A solution of this
equation in $K$ is a column $f=(f_1,\dots, f_n)^{\operatorname{t}}\in K^n$ such that
$f'=Af+b$ where $f':= (f'_1,\dots, f'_n)^{\operatorname{t}}\in K^n$.
Let such an equation $y'=Ay+b$ be given. Define the $n\times n$-matrix $B$ over
$R=K[\der]$ by $B_{ij}:=-A_{ij}$ for~$i\ne j$ and $B_{ii}:=\der - A_{ii}$.
Then the equation $y'=Ay+b$ has clearly the same solutions in $K$ as
the system $B(y)=b$. We claim that the system $B(y)=b$ is automatically $R$-consistent. To see why, let
$r=(r_1,\dots, r_n)\in R^n$ be a row vector and $rB=0$ in $R^n$. Then
$r_j\der=rA_j$ for $j=1,\dots,n$ where $A_j$ is the $j$th column of~$A$. 
Hence $r_1=\dots = r_n=0$, since otherwise an equality $r_j\der=rA_j$
with nonzero $r_j$ of highest order gives a contradiction. This proves our claim.

\begin{cor}\label{matrixdifcon}
If $K$ is linearly surjective, then each matrix differential equation $y'=Ay+b$ over $K$ as above has a solution in $K$. 
\end{cor}

\noindent
Let $n\ge 1$, $a_1,\dots, a_n\in K$, $L=\der^n + a_1\der^{n-1} + \cdots + a_n\in K[\der]$, and $b\in K$. The solutions of the equation
$L(z)=b$ in $K$ are the $f\in K$ such that $L(f)=b$. Setting $y_0:= z$, this equation is equivalent to the system 
\begin{equation}\label{eq:system}
\tag{$\ast$} y_0'\ =\ y_1, \dots, y_{n-2}'\ =\ y_{n-1}, \quad y_{n-1}'\ =\ -(a_1y_{n-1} + \cdots + a_ny_0)+b  
\end{equation}
in the unknowns $y_0,\dots,y_{n-1}$, in the sense that $z\mapsto (z, z',\dots, z^{(n-1)})$ maps the set of solutions in
$K$ of $L(z)=b$ bijectively onto the set of solutions $(y_0,\dots, y_{n-1})\in K^n$ of \eqref{eq:system}. The system \eqref{eq:system} can be written as a matrix equation
$y'= A_Ly+ (0,\dots,0,b)^{\operatorname{t}}$ where $A_L$ is the  $n\times n$ matrix
$$A_L:=\begin{pmatrix}
0		& 1		& 0		& \cdots		& 0 \\
0		& 0		& 1		& \cdots		& 0 \\
\vdots	& \vdots	& \vdots	& \ddots		& \vdots \\
0		& 0		& 0		& \cdots		& 1\\
-a_n		& -a_{n-1}	& -a_{n-2}	& \cdots		& -a_1
\end{pmatrix}$$
over $K$, called the  {\bf companion matrix} of $L$.
Thus the converse of Corollary~\ref{matrixdifcon} is also valid. We can now derive:

\index{matrix!companion}
\index{linear differential operator!companion matrix}
\index{companion!matrix}
\nomenclature[Q]{$A_L$}{companion matrix of $L$}

\begin{cor} \label{cor:lin surj under alg extensions}
Suppose $K$ is linearly surjective and $E$ is a differential
field extension of $K$ and algebraic over $K$. Then $E$ is linearly surjective. 
\end{cor}
\begin{proof} Let $n\ge 1$, let $A$ be an $n\times n$ matrix over $E$, and $b\in E^n$ a column vector. We have to show that the matrix equation $y'=Ay+b$ has a
solution in $E$. The entries of $A$ and $b$ lie in a finite degree extension
of $K$ inside $E$, so we can arrange that $E$ is of finite degree over $K$, say with basis $e_1,\dots, e_m$ over $K$. Writing the~$e_i'$, the $e_ie_j$, and 
the entries of $A$ and $b$ as $K$-linear combinations of $e_1,\dots, e_m$ and making the substitution 
$$ y_j\ =\ z_{j1}e_1 +\cdots + z_{jm}e_m\quad (j=1,\dots,n),$$
one obtains an $mn\times mn$ matrix $A^{\diamond}$ over $K$ and a column vector
$b^{\diamond}\in K^{mn}$ such that any solution of the matrix equation $z'=A^{\diamond}z+b^{\diamond}$
in $K$ yields a solution of $y'=Ay+b$ in $E$. 
\end{proof}

\subsection*{Independence and finite-dimensionality}
Let $n\ge 1$ and let $Y=(Y_1,\dots, Y_n)$ be a tuple of distinct differential indeterminates.
Throughout, $r$ ranges over $\N$.  

\medskip\noindent
Let $K[\der]^n$ be the
(left) $K[\der]$-module of row vectors $(L_1,\dots, L_n)$ with components $L_j\in K[\der]$, and 
$L(L_1,\dots, L_n)=(LL_1,\dots, LL_n)$ for
$L\in K[\der]$. In order to relate homogeneous differential polynomials in $K\{Y\}$ of degree $1$ to vectors in $K[\der]^n$ we consider the $K$-linear space  
$$K\{Y\}_1\ :=\ \sum_{j=1}^n \sum_{r} KY_j^{(r)}$$
of homogeneous differential polynomials in $K\{Y\}$ of degree $1$, and
the $K$-linear bijection
$$A \mapsto A^{\der}:=(A^{\der}_1,\dots, A^{\der}_n)\ \colon\  K\{Y\}_1\to K[\der]^n$$
of (left) $K$-linear spaces, such that for $A=\sum_{j=1}^n\sum_{r} a_{jr}Y_j^{(r)}$, all $a_{jr}\in K$, we have $A^{\der}_j:= \sum_{r} a_{jr}\der^r\in K[\der]$ for $j=1,\dots,n$. Let $A\in K\{Y\}_1$.
Then clearly $$A(y)\ =\ \sum_{j=1}^n A_j^{\der}(y_j)\quad\text{ for $y=(y_1,\dots, y_n)\in K^n$.}$$ 
Also $A'\in K\{Y\}_1$ and $(A')^{\der}= \der A^{\der}$ in the (left) $K[\der]$-module $K[\der]^n$, and so $(A^{(r)})^{\der}=\der^rA^\der$ by induction on $r$. Thus
for $A_1,\dots, A_m\in K\{Y\}_1$ the following are equivalent: \begin{enumerate}
\item the family $\big(A_i^{(r)}\big)$ ($i=1,\dots,m$, $r=0,1,\dots$) is linearly dependent over $K$;
\item there exist $L_1,\dots, L_m\in K[\der]$, not all equal to $0$, such
that  $$L_1A_1^{\der} + \cdots + L_mA_m^{\der}\ =\ 0,$$
that is, $A_1^{\der},\dots, A_m^{\der}$ are $K[\der]$-dependent as defined in the previous section. 
\end{enumerate}

\noindent
If these conditions are satisfied we say that $A_1,\dots, A_m$ are
{\bf $\d$-dependent}; if not, we say that $A_1,\dots, A_m$ are {\bf $\d$-independent}.
(We might add {\em over $K$}, but if $A_1,\dots, A_m$ are $\d$-independent
over $K$, then $A_1,\dots, A_m$ are $\d$-independent over any differential field extension.) If $A_1,\dots, A_m$ are $\d$-independent, then
$m\le n$ by Section~\ref{sec:diagonalization}.  

\index{dependent}
\index{independent}

\medskip\noindent
Let $A_1,\dots, A_m\in K\{Y\}_1$ be given, where $m\ge 1$.
The above $K$-linear bijection $K\{Y\}_1\to K[\der]^n$
maps $A_i$ for $i=1,\dots,m$ to the vector
$$(A_{i,1}^\der,\dots, A_{i,n}^\der)\in K[\der]^n,$$ which we take as the
$i$th row of the $m\times n$ matrix $A^\der=(A_{i,j}^\der)$ over
$K[\der]$. For this matrix and $y=(y_1,\dots, y_n)^{\operatorname{t}}\in K^n$ we have
$$\big(A_1(y),\dots, A_m(y)\big)^{\operatorname{t}}\ =\ A^{\der}(y).$$
Consider the case $m=n$, so we are given $A_1,\dots, A_n\in K\{Y\}_1$. Let $A^{\der}$ be the corresponding $n\times n$ matrix
over $K[\der]$. The zero set 
\begin{align*} Z(A_1,\dots, A_n)\ &:=\ \big\{y\in K^n:\
A_1(y)=\cdots =A_n(y)=0\big\}\\ &\ =\ \big\{y\in K^n:\ A^{\der}(y)=0\big\}
\end{align*} is a $C$-linear subspace of $K^n$. When is it finite-dimensional?

\begin{lemma}\label{fdzind} Assume $K\ne C$. Then the vector space $Z(A_1,\dots, A_n)$ over $C$ is finite-dimensional iff $A_1,\dots, A_n$ are $\d$-in\-de\-pen\-dent. If $A_1,\dots, A_n$ are $\d$-in\-de\-pen\-dent and $K$ is linearly surjective, then for every $a\in K^n$ there is
a $y\in K^n$ with~$A(y)=a$.
\end{lemma}
\begin{proof} Take $P, Q\in \operatorname{GL}_n\!\big(K[\der]\big)$ such that
 $PA^{\der}Q=D$ is diagonal. The $C$-linear bi\-jec\-tion $z\mapsto Q(z)\colon K^n \to K^n$
maps the $C$-linear space $\big\{z\in K^n: D(z)=0\big\}$ onto
the $C$-linear space $\big\{y\in K^n: A^{\der}(y)=0\big\}$.
Thus the latter vector space over~$C$ is finite-dimensional iff
$D_{ii}\ne 0$ for $i=1,\dots,n$, which in turn is equivalent to 
$A_1,\dots, A_n$ being $\d$-independent, using the equivalence of (1) and (2) above. The last claim of the lemma now follows because if
$K$ is linearly surjective and $D_{ii}\ne 0$ for
$i=1,\dots,n$, then there is for every $b\in K^n$ a $z\in K^n$ with
$D(z)=b$.  
\end{proof}

\begin{cor}\label{indlininhom}
Let $A_1,\dots,A_n\in K\{Y\}_1$ be $\d$-independent and  $b_1,\dots,b_n\in K$. Let $L$ be a differential field extension of $K$ and let $y=(y_1,\dots,y_n)^{\operatorname{t}}\in L^n$ be such that $A_i(y)=b_i$ for $i=1,\dots,n$.
Then $y_1,\dots, y_n$ are linear over $K$.
\end{cor}
\begin{proof} Take $P, Q\in \operatorname{GL}_n(K[\der])$ such that
$PA^{\der}Q=D$ is diagonal. Then $D_{ii}\ne 0$ for 
$i=1,\dots,n$ and for 
$z=(z_1,\dots,z_n)^{\operatorname{t}}:= Q^{-1}(y)$ we have $D(z)=P(b)$.
Hence each $z_i$ is linear over $K$, and thus the components
of $y=Q(z)$ are linear over $K$.
\end{proof}

\noindent
Note that we have an isomorphism
$$A=(A_1,\dots, A_m)\mapsto A^\der\ :\  K\{Y\}_1^m\to K[\der]^{m\times n}$$ of left $K$-modules. Let $A=(A_1,\dots, A_m)\in K\{Y\}_1^m$, and recall that then 
$$A^\der(y)\ =\ \big(A_1(y),\dots, A_m(y)\big)^{\operatorname{t}}$$ for all $y=(y_1,\dots, y_n)^{\operatorname{t}}\in K^n$. Also let a tuple $B=(B_1,\dots, B_n)\in K\{Z\}_1^n$ be given,
with $Z=(Z_1,\dots, Z_p)$ a tuple of $p$ distinct indeterminates, $p\ge 1$, and let $B^\der$ be the corresponding $n\times p$ matrix over $K[\der]$. Substituting $B_j$ for $Y_j$ in $A$ yields 
$$ \bigg(A_1\big(B_1(Z),\dots,B_n(Z)\big),\dots, A_m\big(B_1(Z),\dots, B_n(Z)\big)\bigg)\in K\{Z\}_1^m$$
whose corresponding $m\times p$ matrix over $K[\der]$ is $A^\der B^\der$, as is easily verified.

\medskip\noindent
Now let $L$ be a differential field extension of $K$ such that
$[L:K]=n$. Let $b_1,\dots, b_n$ be a basis of $L$ as a vector space over
$K$. Let $X$ be a new differential indeterminate, and
$Y=(Y_1,\dots, Y_n)$ be as before. 
%Thinking of $Z$ as $b_1Y_1+\cdots + b_nY_n$ we consider
%the differential ring morphism $L\{Z\}\to  L\{Y\}$ over $L$ sending
%$Z$ to $b_1Y_1 + \cdots + b_nY_n$. 
Let $A=A(X)\in L\{X\}_1$. Then
$$A(b_1Y_1+\cdots + b_nY_n)\ =\ A_1(Y)b_1+\cdots + A_n(Y)b_n\in L\{Y\}_1$$
with uniquely determined $A_1,\dots, A_n\in K\{Y\}_1$. We now have:

\begin{lemma}\label{fextind} If $K\ne C$ and $A\ne 0$, then $A_1,\dots, A_n$ are $\d$-independent.
\end{lemma}
\begin{proof} The $K$-linear map 
$$y=(y_1,\dots, y_n)\mapsto b_1y_1+ \cdots + b_ny_n\colon K^n \to L$$ 
maps $Z(A_1,\dots, A_n)\subseteq K^n$ bijectively onto 
$Z(A)\subseteq L$. If $A\ne 0$, then $Z(A)$ has finite dimension as a vector space over $C_L$, and so $Z(A_1,\dots, A_n)$ has finite dimension as a vector space over $C$ in view of
$[C_L:C]\le n$. Thus the desired result follows, in view of
Lemma~\ref{fdzind}. 
\end{proof}

\noindent
Let $P\in K\{Y\}$. The homogeneous part of $P$ of degree $1$ is $$P_1\ :=\ \sum_{j=1}^n\left( \sum_r \frac{\partial P}{\partial Y_j^{(r)}}(0) Y_j^{(r)}\right) \in K\{Y\}_1.$$
Note that $(P_1)'\ =\ (P')_1$. Let $y=(y_1,\dots,y_n)\in K^n$. We set
$$P_{+y}\ :=\ P(y+Y)\ =\ P(y_1+Y_1,\dots, y_n+Y_n)\in K\{Y\}$$
and we note that $(P_{+y})' = (P')_{+y}$, and
$$(P_{+y})_1\ =\ \sum_{j=1}^n \left(\sum_r \frac{\partial P}{\partial Y_j^{(r)}}(y)Y_j^{(r)}\right).$$
Let $P_1,\dots,P_m\in K\{Y\}$ and $y\in K^n$. Then we say that $P_1,\dots,P_m$  are {\bf $\d$-dependent at~$y$} if
$(P_{1,+y})_1, \dots,(P_{m,+y})_1$ are $\d$-dependent, and {\bf $\d$-independent at~$y$} otherwise. \index{dependent!at $y$}\index{independent!at $y$} Note that if $L$ is a differential field extension of $K$, then $P_1,\dots,P_m$  are $\d$-dependent at~$y$ with respect to $K$ iff they are
$\d$-dependent at~$y$ with respect to $L$. 
If $P_1,\dots,P_m$  are $\d$-independent at some point of
$K^n$, then $m\leq n$. 

\medskip\noindent
The following lemma gives a simple sufficient condition for $\d$-independence. First, for $P\in K\{Y\}$ and $\vec r=(r_1,\dots,r_n)\in\N^n$ we say that 
$P$ has order at most $\vec r$
if $P\in K\big[Y_j^{(r)}:\ 1\le j\le n,\ 0\le r\le r_j\big]$. 
Thus if
$A\in K\{Y\}_1$ has order at most~$\vec r\in \N^n$, then 
$A=\sum_{j=1}^n \big(\sum_{r=0}^{r_j} a_{jr}Y_j^{(r)}\big)$
with all $a_{jr}=\partial A\big/\partial Y_j^{(r)}\in K$.

\begin{lemma}\label{lem:rubel} 
Let $\vec r=(r_1,\dots,r_n)\in\N^n$, and
assume $A_1,\dots,A_m\in K\{Y\}_1$ have order at most 
$\vec{r}$  and the $m\times n$ matrix $\big(\partial A_i\big/\partial Y_j^{(r_j)}\big)$ over~$K$ has rank $m$. Then $A_1,\dots,A_m$ 
are $\d$-independent.
\end{lemma}
\begin{proof} Suppose not. Take $s\in \N$ minimal such that
$\big(A_i^{(r)}\big)_{1\le i\le m,\ 0\le r\le s}$ is linearly dependent over $K$. 
Take $f_{ir}\in K$ for $1\le i \le m$ and $0\le r \le s$ 
such that $f_{is}\neq 0$ for some $i$ and $\sum_{i,r} f_{ir} A_i^{(r)} = 0$.
Now $A_i^{(r)}$ has order at most $(r_1+r,\dots,r_n+r)$
and $\partial A_i^{(r)}\big/\partial Y_j^{(r_j+r)}=\partial A_i\big/\partial Y_j^{(r_j)}$. 
Comparing coefficients of $Y_j^{(r_j+s)}$ yields
$$f_{1s} \big(\partial A_1\big/\partial Y_j^{(r_j)}\big)+\cdots+ f_{ms} \big(\partial A_m\big/\partial Y_j^{(r_j)}\big)=0\qquad
(j=1,\dots,n),$$ 
so $\big(\partial A_i\big/\partial Y_j^{(r_j)}\big)$ has rank $<m$.
\end{proof}

\noindent
We cannot reverse Lemma~\ref{lem:rubel}: take $m=n=2$, $\vec r=(1,1)$ and 
$A_1:=Y_1'+Y_2'$, $A_2:=Y_1$; then $A_1, A_2$ are $\d$-independent, but $\big(\partial A_i\big/\partial Y_j'\big)$~has~ 
only rank~$1$.

\subsection*{Homogeneous equations}
In this subsection we let $n\geq 1$.

\begin{notation}
Let $R$ be a differential ring. For an $n\times n$ matrix $A=(a_{ij})$ over $R$ we set
$A':=(a_{ij}')$. Then $A\mapsto A'$ is a derivation on the ring of $n\times n$ matrices over $R$:
$$(A+B)'\ =\ A'+B',\quad (AB)'\ =\ A'B+AB'\qquad\text{for $n\times n$ matrices $A$, $B$ over $R$,}$$
and also $(A^{\operatorname{t}})'=(A')^{\operatorname{t}}$ for such $A$.
For $A\in \operatorname{GL}_n(R)$ we have $
  (A^{-1})'=-A^{-1}A'A^{-1}.$
For a column vector $f=(f_1,\dots,f_n)^{\operatorname{t}}\in R^n$ we set 
$f'=(f_1',\dots,f_n')^{\operatorname{t}}$. Then
$$(Af)'\ =\ A'f+Af'\quad\text{for column vectors $f\in R^n$ and $n\times n$ matrices $A$ over $R$.}$$
\end{notation}

\noindent
Let $A$ be an $n\times n$ matrix over $K$. Given a differential ring extension $R$ of $K$, the solutions to the  matrix differential equation $y'=Ay$ form a $C_R$-submodule of $R^n$, which we denote by $\sol_R(A)$; we also set $\sol(A):=\sol_K(A)$. 

\nomenclature[M]{$\sol_R(A)$}{set of solutions of $y'=Ay$ over $R$}

\begin{lemma}\label{lem:lin indep}
If $f_1,\dots,f_r\in \sol(A)$ are $K$-linearly dependent, then $f_1,\dots,f_r$ are $C$-linearly dependent.
\textup{(}In particular, $\dim_C \sol(A)\leq n$.\textup{)}
\end{lemma}
\begin{proof}
Let $f_1,\dots,f_r\in\sol(A)$ with $r\ge 1$ be such that
$f_r=\sum_{i=1}^{r-1} a_if_i$ where $a_1,\dots,a_{r-1}\in K$, and  
$f_1,\dots,f_{r-1}$ are $K$-linearly independent. It is enough to show that
then $f_1,\dots, f_r$ are  $C$-linearly dependent. Now
$$0\ =\ f_r'-Af_r\ =\ \sum_{i=1}^{r-1} a_i'f_i + \sum_{i=1}^{r-1} a_i(f_i'-Af_i)\
=\
 \sum_{i=1}^{r-1} a_i'f_i.$$
Then $a_i'=0$ for $i=1,\dots,r-1$, hence $a_i\in C$ for $i=1,\dots,r-1$.
\end{proof}

\noindent
Some other aspects of homogeneous linear equations are better understood in the setting of differential modules; see the next section.

%% file: mt-5-5.tex
\section{Differential Modules}\label{sec:differential modules}

\noindent
{\em Throughout this section $K$ is a differential field.}\/
We define here differential modules over $K$, and use these to show, among other things, that
if $K$ is linearly closed (respectively, pv-closed), then so is any algebraic differential field extension of $K$.

\medskip\noindent
Let $R$ be a differential ring. If $M$ is a (left) $R[\der]$-module, then the additive map $\der_M\colon M \to M$ given by
$\der_M(x)=\der x$ for $x\in M$ (with $\der\in R[\der]$ 
in the module product $\der x$) satisfies $\der_M(ax)=a'x+a\der_M(x)$ for $a\in R\subseteq R[\der]$, $x\in M$, so 
$\der_M$ is a $\der$-compatible derivation on
the $R$-module $M$ as defined in Section~\ref{sec:modules}.
Conversely, let $M$ be an $R$-module and 
$\der_M$ a $\der$-compatible derivation on the $R$-module $M$.
Then there is a unique (left) $R[\der]$-module with $M$ as its underlying
$R$-module such that $\der_M\colon M \to M$ equals the multiplication
$x\mapsto \der x$ by the scalar $\der\in R[\der]$.

In particular, the derivation $\der$ of $R$ is a 
$\der$-compatible derivation on the $R$-module $R$,
so this makes $R$ an $R[\der]$-module with
$\der a=a'$ for $a\in R$. The differential ideals
of $R$ are exactly the submodules of this 
$R[\der]$-module $R$.  

Note also that for $R[\der]$-modules $M$ and $N$, a map $f: M \to N$
is $R[\der]$-linear iff it is $R$-linear and $f(\der x)=\der f(x)$ for all $x\in M$. 
For $R=K$ this suggests the following notion which turns out to be very useful: 

{\sloppy

\begin{definition}
A {\bf differential module} over $K$ is a finite-dimensional vector space~$M$ over $K$
together with a $\der$-compatible derivation on $M$; we construe such~$M$ as a (left) 
$K[\der]$-module as indicated in the remarks preceding this definition. \index{module!differential} \index{differential module} The {\bf dimension} of a differential module $M$ over $K$ is the dimension $\dim_K M$ of $M$ as a vector space over~$K$.
\end{definition}
}

\noindent
Let $A=(a_{ij})$ be an $n\times n$ matrix over $K$. We make the
$K$-linear space $K^n$ with standard basis $e_1,\dots,e_n$ into a differential 
module $M_A$ over $K$ by requiring
$$\der(e_j)\  =\  -\sum_{i=1}^n a_{ij} e_i \qquad(j=1,\dots,n).$$
This determines $M_A$, and we call $M_A$ the differential module 
associated to $A$.
Note that for $e=\sum_j f_j e_j$ ($f_1,\dots, f_n\in K$) we have
$$\der (e)\  =\ \sum_j f_j'e_j - \sum_i \left( \sum_j a_{ij} f_j \right) e_i\ =\ e'-Ae,$$
from which it follows that $\sol(A)=\{e\in M_A:\ \der(e)=0\}$.
Conversely, if $M$ is a differential module over $K$ of dimension $n$
with basis $b_1,\dots,b_n$ as a $K$-linear space, then there is a unique $n\times n$ matrix $A=(a_{ij})$ over $K$ such that the $K$-linear map $M\to M_A$ with $b_i\mapsto e_i$ for $i=1,\dots,n$ is an isomorphism of $K[\der]$-modules,
and we call this $A$ the  matrix associated to $M$ with respect to the basis $b_1,\dots,b_n$.

\begin{example}
For $L\in K[\der]^{\ne}$ the submodule $K[\der]L$ of the (left)
$K[\der]$-mod\-ule~$K[\der]$ yields the quotient module $K[\der]/K[\der]L$,
which has dimension $\operatorname{order}(L)$ as
a $K$-linear space, so  $K[\der]/K[\der]L$ is a differential module over $K$.
\end{example}

\noindent
Let $A$, $B$ be  $n\times n$ matrices over $K$.
A matrix $P\in\operatorname{GL}_n(K)$ defines an isomorphism $e\mapsto Pe\colon M_A\to M_B$ of $K[\der]$-modules if and only if $\der(Pe_j)=P\der(e_j)$ for
$j=1,\dots,n$, if and only if
$BP-PA=P'$; in this case, any differential ring extension~$R$ of $K$ yields an isomorphism 
$f\mapsto Pf\colon \sol_R(A)\to\sol_R(B)$ of $C_R$-modules.
We call the matrix differential equations $y'=Ay$ and $y'=By$ {\bf equivalent} if $M_A\cong M_B$.
Figure~\ref{fig:diff operators} illustrates the various incarnations of homogeneous linear differential equations.

\index{equivalence!matrix differential equations}
\nomenclature[Q]{$M_A$}{differential module 
associated to an $n\times n$ matrix $A$}

\begin{figure}[h!]
$$\xymatrix{&
\txt<9.25em>{monic operators $L\in K[\der]$ of order $n$} \ar[dl]_{\hskip-2.5em L\mapsto \big[y'=A_Ly\big]} \ar[dr]^{\hskip1.5em L\mapsto K[\der]/K[\der]L}& \\ 
\txt<9.25em>{matrix equations $y'=Ay$, with $A$ an $n\times n$ matrix over $K$}\ar[rr]^{\big[y'=Ay\big] \mapsto M_A} & & \txt<9.25em>{$n$-dimensional differential modules over $K$}
}
$$
\caption{The correspondence between linear differential operators, matrix differential equations, and differential modules.}
\label{fig:diff operators}
\end{figure}

\medskip
\noindent
Given a differential module $M$ over $K$, the elements $f\in M$ with $\der(f)=0$ are said to be
{\bf horizontal}. \index{horizontal} \index{differential module!horizontal elements} \index{differential module!horizontal} The set of horizontal elements of $M$ is a finite-dimensional $C$-linear subspace of $M$. Indeed, if $A$ is
an $n\times n$ matrix over $K$, then $\sol(A)$ is the $C$-linear subspace of the underlying $K$-vector space $K^n$ of $M_A$ consisting of the horizontal elements of $M_A$, 
and $\dim_C \sol(A)\leq n$, with equality iff the $K$-linear space $M_A$ has a basis consisting of horizontal elements
(by Lemma~\ref{lem:lin indep}).

\begin{exampleNumbered}\label{ex:K as K[der]-module}
Turn $K$ into a $K[\der]$-module with scalar multiplication
$$(L,f)\mapsto L(f)\  \colon\  K[\der]\times K\to K.$$ 
Then $K$ becomes a differential module over $K$, and $1\in K$ is horizontal.
\end{exampleNumbered}

\noindent
Let $M$ be a differential module over $K$. If $f\in M^{\ne}$ is horizontal, then
$Kf$ is a submodule of the $K[\der]$-module $M$, and
$a\mapsto af\colon K \to Kf$ is an isomorphism with the above differential module 
$K$. We call $M$ {\bf horizontal} if $M$ is isomorphic as a $K[\der]$-module 
to a direct sum of copies of the above differential module $K$; equivalently, 
$M$ has a basis consisting of horizontal elements.

\medskip
\noindent
Corollary~\ref{cor:smith nf, 2} yields the following for any (left)
$K[\der]$-module $M$:

\begin{cor} \label{cor:cyclic vector, 1}
Suppose $C\neq K$. Then $M$ is a differential module over $K$ if and only if 
$M\cong K[\der]/K[\der]L$ for some monic $L\in K[\der]^{\neq}$.
\end{cor}

\noindent
A {\bf cyclic vector\/} of a differential module $M$ over $K$ of dimension
$n$ is a vec\-tor~${e\in M}$ such that $e, \der e,\dots, \der^{n-1}e$ is a basis
of $M$ as a $K$-linear space. For example, if $L=\der^n+a_1\der^{n-1}+\cdots + a_n\in K[\der]$ ($a_1,\dots,a_n\in K$), then the differential
module $M:= K[\der]/K[\der]L$ of dimension $n$ has cyclic vector
$e:= 1+K[\der]L$ with $Le=0$, and the matrix of $M$ with respect to the basis
$e, \der e,\dots, \der^{n-1}e$ is $-A_L^{\operatorname{t}}$.
Conversely, a cyclic vector $e$ of a differential module $M$ over $K$ of dimension $n$ with $Le=0$ and
$L\in K[\der]$ of order $n$ yields an isomorphism
$K[\der]/K[\der]L\ \to M$ of differential modules over $K$ sending $1+K[\der]L$ to $e$.
If $C\neq K$, then by Corollary~\ref{cor:cyclic vector, 1} every differential module $M\ne \{0\}$ over $K$ has
a cyclic vector. In the next subsection we explain 
the role of cyclic vectors in connection with the 
correspondences in Figure~\ref{fig:diff operators} above.

\index{cyclic vector}

\medskip
\noindent
If $C\ne K$, then by Corollary~\ref{cor:splitting and comp series} and the existence of cyclic vectors 
the differential field
$K$ is linearly closed if and only if every differential module $M\neq\{0\}$ over $K$ 
has a $1$-dimensional differential submodule; it is
easy to see that this equivalence also holds if $C=K$. 
  
\begin{cor}\label{cor:lc going up}
Suppose $K$ is linearly closed. Let $L$ be a differential field extension of $K$ that is algebraic over $K$. Then $L$ is
also linearly closed. 
\end{cor}
\begin{proof}
It suffices to consider the case $[L:K]<\infty$. For this case,
let $M\neq\{0\}$ be a differential module over $L$ and let $M_K$ denote $M$ viewed as a $K[\der]$-module. Then~$M_K$ is a differential module over $K$. Hence we can take a $1$-dimensional differential submodule $N$ of $M_K$.
Then the $L$-linear subspace $LN$ of $M$ generated by $N$ is an $L[\der]$-submodule of $M$ with $\dim_L LN=1$.
%We have a surjective $L$-linear map $L\otimes_K N\to LN$, so $1\leq\dim_L LN\leq\dim_L L\otimes_K N=1$.
\end{proof}

\subsection*{Duality}
Let $M$ and $N$ be $K[\der]$-modules.
The $K$-linear space $\Hom_K(M,N)$ of all $K$-linear maps $M\to N$ is
made into a $K[\der]$-module by defining
$$(\der\phi)(f)\ =\  \der(\phi f)-\phi(\der f)\qquad\text{for $\phi\in\Hom_K(M,N)$ and $f\in M$.}$$
If $M$ and $N$ are differential modules over $K$, then so is $\Hom_K(M,N)$, with
$$\dim_K \Hom_K(M,N)\ =\ \dim_K M\cdot\dim_K N,$$
and if in addition $M$ and $N$ are horizontal, then so is $\Hom_K(M,N)$.
A special case of this construction is the {\bf dual} $M^* := \Hom_K(M,K)$ of $M$. 
Here we view $K$ as a horizontal differential module over $K$ as explained in Example~\ref{ex:K as K[der]-module}.
Writing $\< \phi, f\> := \phi(f)\in K$ for $\phi\in M^*$ and $f\in M$, we have
the identity
\begin{equation}\label{eq:der on dual}
\der \<\phi, f\>\ =\ \<\der\phi,f\>+\<\phi,\der f\>.
\end{equation}
The natural $K$-linear map $f\mapsto \<-,f\>: M\to M^{**}$ is a morphism of $K[\der]$-modules,
and if $M$ is a differential module over $K$, then this morphism is an isomorphism.
In the rest of this subsection, we let  $M\neq \{0\}$ be a  differential module over $K$ and $e_0,\dots,e_{n-1}$ be a basis of $M$, and we let $A$ be the matrix associated to $M$ with respect to $e_0,\dots,e_{n-1}$, so $A=(a_{ij})_{0\le i\le n-1,\ 0\le j\le n-1}$ with $\der(e_j)=-\sum_i a_{ij}e_i$ for all $j$.  
We also let $A^*$ be the matrix associated to $M^*$ with respect to
the basis $e_0^*,\dots,e_{n-1}^*$ of~$M^*$ dual to $e_0,\dots,e_{n-1}$.

\index{differential module!dual}
\nomenclature[Q]{$M^*$}{dual of the $K[\der]$-module $M$}

\begin{lemma}
$e_0,\dots,e_{n-1}$ are horizontal iff $e_0^*,\dots,e_{n-1}^*$ are horizontal.
\end{lemma}
\begin{proof}
Let $i$ and $j$ range over $\{0,\dots,n-1\}$.
From \eqref{eq:der on dual} we obtain
\begin{equation}\label{eq:der on dual basis}
\< \der e_i^*,e_j \> + \<e_i^*,\der e_j\>\ =\ \der\<e_i^*,e_j\>\ =\ \der\delta_{ij}=0,
\end{equation}
where $\delta_{ij}$ is the Kronecker delta.
If each $e_j$ is horizontal, then 
for each $i$,
$$\<\der e_i^*,e_j\>\ =\ -\<e_i^*,\der e_j\>\ =\ -\<e_i^*,0\>\ =\ 0\quad\text{for each $j$,}$$
so $\der e_i^*=0$. This shows the forward direction, and the
converse goes likewise.
\end{proof}

\noindent
By the previous lemma, if $M$ is horizontal, then so is $M^*$.

\begin{lemma}
$A^*=-A^{\operatorname{t}}$.
\end{lemma}
\begin{proof}
By \eqref{eq:der on dual basis} we obtain for $j,k=0,\dots,n-1$, 
$$\<\der e_j^*, e_k\>\	=\ - \<e_j^*, \der e_k\>\ 
						=\ \sum_i a_{ik} \<e_j^*,e_i\>\ =\ a_{jk}$$
and thus $\der e_j^* = \sum_i a_{ji} e_i^*$ as required by the lemma.
\end{proof}

\noindent
The matrix differential equation $y'=A^*y$ is called
the {\bf adjoint equation} of the matrix differential equation $y'=Ay$.
In the next lemma, given $L\in K[\der]$, we let $L^*\in K[\der]$ denote the adjoint of $L$ as defined in Section~\ref{Linear Differential Operators}.

\index{system of linear differential equations!adjoint}
\index{adjoint!system of linear differential equations}
\nomenclature[Q]{$y'=A^*y$}{adjoint equation of  $y'=Ay$}

\begin{lemma}\label{lem:Gabber}
Suppose $e_0$ is a cyclic vector of $M$ and $e_i=\der^i e$ for $i=0,\dots,n-1$.
Let $L\in K[\der]$ be monic of order $n$ 
with $Le=0$, $e:= e_0$.  Then $e^*:=e_{n-1}^*$ is a cyclic vector of $M^*$
and $L^*e^*=0$.
\end{lemma}
\begin{proof} We have $L=\der^n+\sum_{i=1}^{n}a_{n-i}\der^{n-i}$ 
with $a_0,\dots,a_{n-1}\in K$. Recall that $\der^*=-\der$,
so $L^*=(\der^*)^n+\sum_{i=1}^{n} (\der^*)^{n-i} a_{n-i}$.
By the previous lemma, $A^*=A_L$, so
$$\der e_i^* + e_{i-1}^*\ 	=\  a_i e_{n-1}^* \qquad\text{for $i=0,\dots,n-1$,}$$
where $e_{-1}^*:=0$.
These identities give $M^*=Ke^*+K\der e^*+\cdots+K\der^{n-1}e^*$, so
$e^*$ is a cyclic vector of $M^*$. By induction on $i$ they also yield 
$$\big( (\der^*)^i + \sum_{j=1}^i(\der^*)^{i-j}a_{n-j}\big) e^*\ =\ e^*_{n-i-1}\qquad\text{for $i=0,\dots,n$.}$$
Taking $i=n$ we obtain $L^*e^*=0$ as claimed.
\end{proof}

\begin{cor} Let $L\in K[\der]$ be monic of order $n\ge 1$. Then
$$(K[\der]/K[\der]L)^*\ \cong\ K[\der]/K[\der]L^*\ \cong\ M_{A_L},$$
as differential modules over $K$.
\end{cor}
\begin{proof} Set $M:= K[\der]/K[\der]L$ and $e:= 1+K[\der]L\in M$. Then 
$M$ has matrix~$-A_L^{\operatorname{t}}$
with respect to the basis $e,\der e,\dots,\der^{n-1}e$, so
$M^*$ has matrix $A_L$ with respect to the dual basis, and thus
$M^*\cong M_{A_L}$. By the lemma above, 
$M^*\cong K[\der]/K[\der]L^*$.
\end{proof}

\begin{cor}\label{cor:cyclic vector, 2}
If $C\neq K$ and $n\ge 1$, then any matrix differential equation $y'=Ay$, where $A$ is an $n\times n$ matrix over $K$ is equivalent to a matrix differential equation $y'=A_Ly$ where $L\in K[\der]$ is monic of order $n$.
\end{cor}

\begin{proof} Assume $C\ne K$ and $A$ is an $n\times n$ matrix over $K$,
$n\ge 1$. 
By Corollary~\ref{cor:cyclic vector, 1} we have $M_A\cong K[\der]/K[\der]L^*$
where $L\in K[\der]$ is monic of order $n$; then $M_A\cong M_{A_{L}}$ by the preceding corollary.
\end{proof}

\noindent
In the next subsection we use the following {\em lemma on integrating factors}. This concerns the situation where $e_0,\dots, e_{n-1}$ are horizontal,
$b\in M$ is given, and we wish to integrate $b$, that is, determine $f\in M$ such that $\der f = b$. This can be done if we can integrate the ``Fourier coefficients'' $\<e_i^*,b\>$ of $b$: 

\begin{lemma} \label{lem:integrating factors}
Let $e_0,\dots,e_{n-1}$ be horizontal, and $b\in M$.
Suppose $a_i\in K$ are such that $a_i'=\<e_i^*,b\>$, for $i=0,\dots,n-1$. 
Then $f:=\sum_i a_ie_i$ satisfies $\der f=b$.
\end{lemma}
\begin{proof} From $\der e_i=0$ we get
\equationqed{\der f\	=\ \sum_i \der (a_i e_i)\ 
		=\ \sum_i (a_i'e_i+a_i\der e_i)\ 
		=\ \sum_i \<e_i^*,b\>e_i\  =\  b.}
\end{proof}

\subsection*{Fundamental matrices}
In this subsection $n$ ranges over $\N^{\geq 1}$.
Let $A$ be an $n\times n$ matrix over~$K$.
Note that if $F$ is an $n\times n$ matrix over $K$, then $F'=AF$ iff the columns of $F$ lie in $\sol(A)$,
and in this case,
$F\in \operatorname{GL}_n(K)$  iff the columns of $F$ form a basis of the $C$-linear space~$\sol(A)$.
Let $R$ be a differential ring extension of~$K$. A matrix $F\in \operatorname{GL}_n(R)$ is called a {\bf fundamental matrix} for the differential equation
$y'=Ay$ if $F'=AF$. We say that~$R$ {\bf contains a fundamental matrix for $y'=Ay$} if  $\operatorname{GL}_n(R)$ contains one.

\index{system of linear differential equations!fundamental matrix}
\index{matrix!fundamental}
\index{fundamental matrix}

\begin{lemma}\label{lem:all fund matrices}
Let $F,G\in \operatorname{GL}_n(R)$ be fundamental matrices for $y'=Ay$. Then $F^{-1}G$ lies in the subgroup $\operatorname{GL}_n(C_R)$ of  
$\operatorname{GL}_n(R)$. Thus  
%has coefficients in $C_R$. Hence if $C_R=C$ or $R$ is an integral domain, then 
the set of all fundamental matrices for $y'=Ay$ in $\operatorname{GL}_n(R)$ is 
$F\operatorname{GL}_n(C_R)$.
\end{lemma}
\begin{proof}
Set $P:=F^{-1}G$; then $FP=G$, so
$$AG\ =\ G'\ =\ F'P+FP'\ =\ AFP+FP'\ =\ AG+FP',$$
hence $P'=0$, and thus $P$ has all its entries in $C_R$. Likewise, the inverse 
$G^{-1}F$ of~$P$ has all its entries in $C_R$, so 
$P\in \operatorname{GL}_n(C_R)$.
%Under the extra stated hypothesis on $R$, we have  
%$R^\times\cap C_R = C_R^\times$, and so $P\in \operatorname{GL}_n(C_R)$.
\end{proof}

\begin{lemma}\label{lem:adjoint, 1}
Let   $F\in\operatorname{GL}_n(R)$. Then $F$ is a fundamental matrix for $y'=Ay$ iff $G:=(F^{\operatorname{t}})^{-1}$ is a fundamental matrix for the adjoint equation $y'=A^*y$.
\end{lemma}
\begin{proof}
First note that $G\in \operatorname{GL}_n(R)$.
We  have $G'=-G(F')^{\operatorname{t}}G$, and a straightforward computation now shows that $G'=A^*G\Longleftrightarrow F'=AF$. %\marginpar{Uncomment to see the computation.}
%\begin{align*}
%G'=A^*G	&\quad\Longleftrightarrow\quad -G(F')^{\operatorname{tr}}G=-A^{\operatorname{tr}}G \\
%		&\quad\Longleftrightarrow\quad G(F')^{\operatorname{tr}}=A^{\operatorname{tr}} \\
%		&\quad\Longleftrightarrow\quad F'G^{\operatorname{tr}}=A \\
%		&\quad\Longleftrightarrow\quad F'=A(G^{\operatorname{tr}})^{-1} = AF. 
%\end{align*}
\end{proof}

\begin{lemma}\label{lem:abel, generalized}
Let $F\in\operatorname{GL}_n(R)$ be a fundamental matrix for $y'=Ay$. Then 
$$(\det F)'\ =\ \operatorname{tr} A \cdot \det F.$$
\end{lemma}
\begin{proof} Let $i$,~$j$ range over $\{1,\dots,n\}$, and
let $f_j\in R^n$ be the $j$th row of $F$. By expanding the determinant and the product formula:
$$(\det F)'\ 	=\  \sum_i \det F_i\qquad\text{where $F_i\ =\ \left(\begin{smallmatrix} f_1 \\ \vdots \\  f_i' \\   \vdots \\ f_n\end{smallmatrix}\right)$.}$$
We have $f_i'=\sum_j a_{ij} f_j$ since $F'=AF$. Subtracting from the $i$th row of $F_i$ the linear combination $\sum_{j \neq i} a_{ij} f_j$ replaces therefore
the $i$th row of $F_i$ by $a_{ii}f_i$, and so $\det F_i=a_{ii}\det F$.
Hence  $(\det F)'=\sum_i a_{ii}\det F=\operatorname{tr} A\cdot \det F$.	\end{proof}

\noindent
Consider a monic operator 
$$L\ =\ \der^n+a_{n-1}\der^{n-1}+\cdots+a_0\in K[\der]\qquad (a_0,\dots,a_{n-1}\in K),$$
with companion matrix $A_L$.
The map
\begin{equation}\label{eq:ker L iso sol(A_L)}
f\mapsto \big(f,f',\dots,f^{(n-1)}\big)^{\operatorname{t}}\ \colon\ \ker_R L 
\to \sol_R(A_L)
\end{equation}
is an isomorphism of $C_R$-modules.
Thus, given an $n\times n$ matrix $F$ over $R$, we have:
$F$ is a fundamental matrix for $y'=A_Ly$ iff there are $f_1,\dots,f_n\in \ker_R L$ such that $F=\operatorname{Wr}(f_1,\dots,f_n)$
and $\operatorname{wr}(f_1,\dots,f_n)\in R^\times$.
%, and in this case $\operatorname{tr} A=-a_{n-1}$. Thus 
%Lemma~%\ref{lem:abel, generalized} generalizes Lemma~\ref{lem:abel}. 

\begin{lemma}\label{prop:fund matrix, fund system}
Suppose $C\neq K$ and
$R$ contains a fundamental matrix for each equation $y'=A_Ly$ where $L\in K[\der]$ is monic of order $n$. Then
$R$ contains a fundamental matrix for each equation $y'=Ay$, where $A$ is an $n\times n$ matrix over~$K$.
\end{lemma}
\begin{proof}
Let $A$ be an $n\times n$ matrix over $K$. Corollary~\ref{cor:cyclic vector, 2} 
provides a monic $L\in K[\der]$ of order $n$ such that
$y'=Ay$ is equivalent to $y'=A_Ly$. Take $P\in\operatorname{GL}_n(K)$ such that 
$P\sol_R(A_L)=\sol_R(A)$, and let $F\in \operatorname{GL}_n(K)$ be a fundamental matrix for $y'=A_Ly$. 
Then $PF\in\operatorname{GL}_n(R)$ is a fundamental matrix for~$y'=Ay$.
\end{proof}

\noindent
In the rest of this subsection we assume that our differential ring extension $R$ of~$K$ is an integral domain whose
ring of constants $C_R$ is a field that equals the field of constants of the differential fraction field of $R$. Let $L\in K[\der]$ be monic of
order $n$ as displayed above.
Then $\ker_R L$ is a $C_R$-linear subspace of $R$ and
$\dim_{C_R} \ker_R L\leq n$. 
Moreover, by Lemma~\ref{lem:wronskian},  
$y_1,\dots,y_n\in R$ are $C_R$-linearly independent iff 
$\wr(y_1,\dots,y_n)\neq 0$.
In particular, if $\operatorname{Wr}(y_1,\dots,y_n)$ is a fundamental matrix for $y'=A_Ly$,
then $y_1,\dots,y_n$ is a basis of the $C_R$-linear space $\ker_R L$. 
The following lemma contains a partial converse.  

\begin{lemma}\label{lem:fundamental system}
Suppose  that $\ker_R (\der-a_{n-1})\neq \{0\}$. Let $y_1,\dots,y_n\in\ker_R L$ be such that
$w:=\operatorname{wr}(y_1,\dots,y_n)\neq 0$. Then $w\in R^\times$.
\end{lemma}
\begin{proof}
By  Lemma~\ref{lem:abel} we have $w'+a_{n-1}w=0$; so  $w^{-1}\in \Frac(R)$ satisfies $(w^{-1})'-a_{n-1}w^{-1}=0$.
By assumption, we also have $v\in R^{\neq}$ with $v'-a_{n-1}v=0$. Hence $w^{-1}=cv$
for some $c\in C_R^\times$; in particular, $w\in R^\times$.
\end{proof}

\begin{cor}\label{cor:fund matrix, fund system}
Suppose $C\neq K$, and let $r\in\N^{\geq 1}$.
The following are equivalent:
\begin{enumerate}
\item[\textup{(i)}] $R$ contains a fundamental matrix for $y'=Ay$, for every $n\le r$ and $n\times n$ matrix~$A$ over $K$;
\item[\textup{(ii)}] $\dim_{C_R} \sol_R A = n$ for every $n\le r$ and $n\times n$ matrix $A$ over $K$;
\item[\textup{(iii)}] $\dim_{C_R} \ker_R L = \order L$ for every $L\in K[\der]^{\neq}$ of order at most~$r$;
\item[\textup{(iv)}] $R$ contains a fundamental matrix for $y'=A_Ly$, for every monic
$L\in K[\der]$ of positive order at most~$r$.
\end{enumerate}
\end{cor}
\begin{proof}
We clearly have (i)~$\Rightarrow$~(ii). 
The implication (ii)~$\Rightarrow$~(iii) is immediate from the isomorphism \eqref{eq:ker L iso sol(A_L)}, 
(iii)~$\Rightarrow$~(iv) holds by Lem\-ma~\ref{lem:fundamental system}, and
(iv)~$\Rightarrow$~(i) by Lem\-ma~\ref{prop:fund matrix, fund system}.
(Only the last implication used the hypothesis $C\neq K$.)
\end{proof}

\noindent
Next an analogue of Corollary~\ref{matrixdifcon} for homogeneous linear differential equations: 

\begin{cor}\label{matrixdifcon, homog}
Let $r\in\N^{\geq 1}$. The following are equivalent:
\begin{enumerate}
\item[\textup{(i)}] $K$ contains a
fundamental matrix for $y'=Ay$, for every $n\le r$ and
$n\times n$ matrix $A$ over $K$;
\item[\textup{(ii)}] $\dim_C \sol A=n$ for every $n\le r$ and $n\times n$ matrix $A$ over $K$;
\item[\textup{(iii)}] $K$ is $r$-pv-closed;
\item[\textup{(iv)}] $K$ contains a fundamental matrix for $y'=A_Ly$, for every
monic $L\in K[\der]$ of positive order at most~$r$;
\item[\textup{(v)}] every differential module $M$ over $K$ with $\dim_K M\leq r$ is horizontal.
\end{enumerate}
\end{cor}
\begin{proof}
The equivalence of (i)--(iv) follows from Corollary~\ref{cor:fund matrix, fund system}.
The equivalence of (iv) and (v) involves Lemma~\ref{lem:lin indep}.
(Each of (i)--(v) implies  $C\ne K$.)
\end{proof}

\noindent
Corollary~\ref{cor:pv-closed and alg extensions} below is an analogue of Corollary~\ref{cor:lin surj under alg extensions}, and follows easily from Corollary~\ref{matrixdifcon, homog} and the next lemma:

\begin{lemma}\label{lem:fund matrix over alg ext}
Suppose $C\neq K$.
Let $E$ be a differential subring of $R$ containing~$K$ such that
$E$ is an algebraic field extension of $K$. Suppose 
$R$ contains a fundamental matrix for $y'=Ay$, for every $n$ and 
$n\times n$ matrix $A$ over $K$.
Then~$R$ also contains a fundamental matrix for $y'=Ay$, for
every $n$ and $n\times n$ matrix $A$ over $E$.
\end{lemma}
\begin{proof} Let $n\ge 1$ and let 
$A$ be an $n\times n$ matrix over $E$; by  
Corollary~\ref{cor:fund matrix, fund system} it suffices to show that $\dim_{C_R} \sol_R A = n$.
We may assume that $E$ is of finite degree over $K$, say with basis $e_1,\dots, e_m$ over $K$. 
Take the $mn\times mn$ matrix $A^{\diamond}$ over $K$ from the
proof of Corollary~\ref{cor:lin surj under alg extensions}. As
in that proof we see that for all $z\in\sol_R A^{\diamond}$ we have
$Bz\in \sol_R A$, where $B$ is the $n\times mn$ matrix
$$B\ =\  \begin{pmatrix} 
 e_1 & e_2 & \cdots & e_m &        &     &     &        &    \\
     &     &        &     & \ddots &     &     &        &    \\
     &     &        &     &        & e_1 & e_2 & \cdots & e_m
\end{pmatrix}$$
over $E$. 
Let $Z\in \operatorname{GL}_{mn}(R)$ be a fundamental matrix for $z'=A^{\diamond} z$. 
Then the $n\times mn$ matrix $BZ$ over $R$ has rank
$n$ over $\Frac(R)$, and its $mn$ columns lie in $\sol_R A$, so 
$\sol_R A$ has dimension at least $n$ as a vector space over $C_R$.
By Lemma~\ref{lem:lin indep} this dimension is also at most $n$, so equal to $n$.
\end{proof}

\noindent
The situations considered in the next two corollaries, with $E\ne K$ in the first one, actually occur, by a fact mentioned in the \textit{Notes and comments}\/ to this section. 

\begin{cor}\label{cor:pv-closed and alg extensions}
If $K$ is pv-closed and $E$ is a differential field extension of $K$ and algebraic over $K$, then $E$ is pv-closed. 
\end{cor}

\noindent
From this corollary we now obtain a variant of Lemma~\ref{lem:pv closed => alg closed}:

\begin{cor}\label{corpvrcf}
Suppose $K$ is pv-closed and $C$ is real closed. Then $K$ is real closed.
\end{cor}
\begin{proof}
The differential field extension $E=K[\imag]$, where $\imag^2=-1$, is algebraic over~$K$ and has constant field $C_E=C[\imag]$.
By the previous corollary, $E$ is pv-closed, and 
since $C$ is real closed, $C_E$ is algebraically closed. Hence $E$ is algebraically closed by
Lemma~\ref{lem:pv closed => alg closed} and so  $K$ is real closed by Theorem~\ref{thm:AS}.
\end{proof}

\noindent
Let $A$ be an $n\times n$ matrix over $K$, and identify the $K$-linear space $M=M_A$ with its dual $M^*$
via the isomorphism $f\mapsto \<f,{-}\>\colon M\to M^*$, 
where $\<\hskip0.5em ,\hskip0.5em \>$ denotes the
usual inner product $M\times M \to K$ given by
$$\<f,g\>\ :=\ \sum_i f_ig_i \qquad (f=(f_1,\dots,f_n),\ g=(g_1,\dots,g_n)\in K^n)$$
on $M=K^n$. Let $F\in\operatorname{GL}_n(K)$ be a fundamental matrix of $y'=Ay$,
with columns $f_1,\dots,f_n$, and let  $G=(F^{\operatorname{t}})^{-1}$, with
columns $g_1,\dots,g_n$. Then $\<g_i,f_j\>=\delta_{ij}$ for $i,j=1,\dots,n$, so $g_1,\dots,g_n$ is the
basis of $M^*$ dual to the horizontal basis $f_1,\dots,f_n$ of $M$.
The next result goes under the name of \textit{variation of constants.}\/

\begin{lemma}\label{lem:adjoint, 2}
Let $b\in K^n$ be a column vector, and
suppose $a_i\in K$ satisfies $a_i'=\<g_i,b\>$, for $i=1,\dots,n$. Then 
$f:=\sum_i a_i f_i\in K^n$ satisfies $f'=Af+b$.
\end{lemma}

\noindent
This is immediate from Lemma~\ref{lem:integrating factors}. 
By Lemma~\ref{lem:r-pv closed}, if $K$ is pv-closed, then~$K$ is linearly surjective. 
The following corollary complements this fact:  if one is already given a basis of the kernel of a monic differential operator $L\in K[\der]^{\neq}$, then solving inhomogeneous linear differential equations $L(y)=b$ reduces to integration.

\begin{cor}\label{cor:adjoint}
Let $L\in K[\der]^{\neq}$ be monic of order $n$ and 
suppose $y_1,\dots,y_n$ is a basis of $\ker L$.
Let $(z_1,\dots,z_n)$ be the bottom row of the matrix $(W^{\operatorname{t}})^{-1}$, where
$W=\Wr(y_1,\dots,y_n)$. 
Let $b\in K$, and
suppose $a_i\in K$ satisfies $a_i'=bz_i$, for $i=1,\dots,n$. Then
$y:=\sum_i a_iy_i\in K$ satisfies $L(y)=b$.
\end{cor}
\begin{proof}
Apply the previous lemma with the companion matrix $A_L$ of $L$ and 
the column vector $(0,\dots,0,b)^{\operatorname{t}}$ in place of  $A$ and $b$.
\end{proof}

\begin{example}
Let $L=\der-f$ where $f\in K$, let $y_1\in K^\times$ be such that $y_1'=fy_1$,
and let $b\in K$. Then any $a\in K$ with  $a'=b/y_1$ yields $L(y)=b$ for $y:=ay_1$.
\end{example}

\subsection*{Notes and comments}
Lemma~\ref{consistent} is from \cite[pp.~91--94]{CF}.
Corollary~\ref{cor:cyclic vector, 1} was apparently first shown by Loewy \cite{Loewy} for the case where $K$ is a field of
meromorphic functions, and the argument used here (via Corollary~\ref{cor:smith nf, 1}) was  suggested in~\cite{Adjamagbo}; see \cite{ChKo} for a historical discussion.
Lemma~\ref{lem:Gabber} is credited to O.~Gabber in \cite[Lem\-ma~(1.5.3)]{Katz}.
The notion of adjoint equation and Lemmas~\ref{lem:adjoint, 1},~\ref{lem:adjoint, 2}, and Corollary~\ref{cor:adjoint} stem from the classical literature
on linear differential equations, see
\cite[Chapter~III, \S{}10]{Poole}.

\marginpar{fact from \cite{CHP, PillayKamensky} taken on faith} By \cite{CHP,PillayKamensky}, any real closed $K$ has a pv-closed differential field extension~$L$ with~$C_L=C$; note that any such $L$ is real closed by Corollary~\ref{corpvrcf} applied to~$L$.

%% file: mt-5-6.tex
\section{Linear Differential Operators in the Presence of a Valuation}\label{sec:ldopv}

\noindent
{\em In this section $K$ is a valued differential field with valuation $v$ and derivation~$\der$, such that $\der \smallo \subseteq\smallo$}. We let 
$a$, $b$ (possibly with subscripts) range over $K$. 

Note that $\der$ induces a derivation on the residue field
$\k=\mathcal{O}/\smallo$ which we also denote by~$\der$, so 
$\der(a+\smallo)=\der(a) + \smallo$ for $a\in \mathcal{O}$. We have 
$\mathcal{O}[\der]$ as a 
subring of $K[\der]$, with a ring homomorphism
$$\sum_i a_i\der^i \mapsto 
\sum_i (a_i+\smallo)\der^i\ \colon\ \mathcal{O}[\der] \to \k[\der].$$
We extend $v$ to a valuation on the additive group of
$K[\der]$ by setting 
$$v(A)\ :=\ \min_i v(a_i)\in\Gamma_\infty\ \text{ for }\ A=\sum_i a_i\der^i \in K[\der].$$ Note that then $v(aA)= va + v(A)$. Also
$v(Aa)$ depends, for fixed $A\in K[\der]$, only on $\gamma=va$, 
by Lemma~\ref{v-under-conjugation}(ii), but generally 
in a more complicated way; thus we can define 
$v_A(\gamma):=v(Aa)$ for $\gamma=va$.

\begin{lemma}\label{Functoriality-Easy} Let $A,B\in K[\der]^{\ne}$
and $\gamma\in\Gamma$. Then:
\begin{enumerate}
\item[\textup{(i)}] if $v(B)=0$, then $v(AB)=v(A)$;
\item[\textup{(ii)}] $v_{AB}(\gamma)= v_A\bigl(v_B(\gamma)\bigr)$.
\end{enumerate}
\end{lemma}
\begin{proof} For (i) we use $v(aA)=va+v(A)$ to reduce to the case
$v(A)=v(B)=0$. Then $v(AB)=0$ follows by applying the homomorphism
$\mathcal{O}[\der] \to \k[\der]$. For (ii),
take $y\in K^\times$ with $vy=\gamma$, so $v_B(\gamma)=v(By)$ and
$v_{AB}(\gamma)=v(ABy)$. Now $By=b\hat{B}$ with $vb=v(By)$ and
$v(\hat{B})=0$. Hence $ABy=Ab\hat{B}$, so by (i):
$$v_{AB}(\gamma)\ =\ v(ABy)\ =\ v(Ab\hat{B})\ =\ v(Ab)\ 
=\ v_A(vb)\ =\ v_A\bigl(v(By)\bigl)\ =\ v_A\bigl(v_B(\gamma)\bigr)$$
as claimed.
\end{proof}

\begin{cor}\label{speclinclosed} If $K$ is $r$-linearly closed, then so is $\k$.
\end{cor} 
\begin{proof} Let $A\in \mathcal{O}[\der]$ be monic, and $A=BD$ with $B,D\in K[\der]$ and $B, D$ of order~$\ge 1$. 
Take $a\in K^\times$ with $va=-vD$. Then $A=(Ba^{-1})\cdot (aD)$ with $v(aD)=0$, and thus $v(Ba^{-1})=vA=0$. Therefore, if $A\in \mathcal{O}[\der]$ is monic of order $\ge 1$ and its image in $\k[\der]$ is irreducible, then $A$ is irreducible in $K[\der]$. The lemma is an easy consequence of this fact. 
\end{proof}

\begin{cor}\label{Gauss}
Let $A,B\in K[\der]$ be monic and $AB\in \mathcal{O}[\der]$. Then $A,B\in \mathcal{O}[\der]$.
\end{cor}
\begin{proof}
By Lemma~\ref{Functoriality-Easy}(ii) and 
Lemma~\ref{v-under-conjugation}(iii) we have 
$$0\ =\ v(AB)\ =\ v_{AB}(0)\ =\ v_A\big(v_B(0)\big)\ 
\le\ v_A(0)\ \le\ 0,$$
so $v(A)=v_A(0)=0$, and hence $v(B)=v_B(0)=0$ by 
Lemma~\ref{v-under-conjugation}(iii). 
\end{proof}

\noindent
For $A=\sum_i a_i\der^i \in K[\der]$ we have $a_0=A(1)$ and 
we put $$\Pmu(A)\ :=\ \min\bigl\{i:v(a_i)=v(A)\bigr\}\in \N,$$
so $\Pmu(A)=0 \Longleftrightarrow v(a_0) = v(A)$.
Here is a special case of Lemma~\ref{mu_P}:

\begin{lemma}\label{mu-lemma} Suppose that
$\der\mathcal{O}\subseteq \smallo$. Then, given any
$A\in K[\der]^{\ne}$, the quantity
$\Pmu(Ag)$ for $g\in K^{\times}$ depends only on $vg$.
\end{lemma}

\noindent
If $\der\mathcal{O}\subseteq \smallo$, then we define 
for $A\in K[\der]$, $A\ne 0$:
$$\Pmu_A(\gamma)\ :=\ \Pmu(Ag),\ \text{ where }
g\in K^{\times},\ vg=\gamma.$$

\begin{lemma}\label{nLM-formula} Suppose that 
$\der\mathcal{O}\subseteq \smallo$. 
Let $A,B\in K[\der]^{\ne}$. Then
$$ \Pmu_{AB}(\gamma)\ =\ \Pmu_A\bigl(v_B(\gamma)\bigr)+\Pmu_B(\gamma)\qquad 
(\gamma\in \Gamma).$$
\end{lemma}
\begin{proof}
The induced derivation
on $\k$ is trivial, so $\k[\der]$ is commutative. Therefore, if $v(A)=v(B)=0$, then 
$\Pmu(AB)=\Pmu(A) + \Pmu(B)$.  
With $b$, $y$, $\gamma=vy$, and~$\hat B$ as in the proof of part (ii) of 
Lemma~\ref{Functoriality-Easy},
write $Ab=a \hat{A}$ where
$va=v_A(vb)=v_A\bigl(v_B(\gamma)\bigr)$ and $v(\hat A)=0$. Then
$ABy=a \hat{A}\hat{B}$, so 
$$\Pmu_{AB}(\gamma)\ =\ \Pmu(\hat{A} \hat{B})\ =\ 
\Pmu(\hat{A}) + \Pmu(\hat{B}).$$
The claim now follows from $\Pmu(\hat B)=\Pmu(b\hat{B})=\Pmu(By)=\Pmu_B(\gamma)$ 
and 
\equationqed{\Pmu(\hat{A})\ =\ \Pmu(a\hat{A})\ =\ \Pmu(Ab)\ 
=\ \Pmu_A(vb)\ =\ \Pmu_A\bigl(v_B(\gamma)\bigr).}
\end{proof}

\noindent
Suppose that $C\subseteq\mathcal O$ and $A\in K[\der]^{\ne}$.
If $y_1,\dots, y_r\in\ker^{\neq} A$ and
 $y_1\succ \dots \succ y_r$, then $y_1,\dots,y_r$ are $C$-linearly independent. So $\abs{v(\ker^{\neq} A)}\leq\dim_C \ker A$. Moreover:

\begin{lemma} \label{valuation-of-basis}  Suppose that 
$\mathcal{O}=C+\smallo$, and let $A\in K[\der]^{\ne}$. 
Then 
\begin{enumerate}
\item[\textup{(i)}] $v(\ker^{\ne} A)$ is finite, of size $\dim_C\ker A$;
\item[\textup{(ii)}] if $v(\ker^{\ne} A)= \{vy_1, \dots , vy_m\}$ where
$y_1,\dots, y_m\in\ker^{\ne} A $, $vy_1 < \dots < vy_m$, 
then $y_1,\dots,y_m$ is a basis of $\ker A$.
\end{enumerate}
\end{lemma} 
\begin{proof} This follows from Lemma~\ref{Hahn valuations, lemma}, in view of the fact that $K$ with the valuation~$v$ is a Hahn space over $C$.
\end{proof}

\noindent
The assumptions $\der\mathcal{O}\subseteq \smallo$ 
and $\mathcal{O}=C+\smallo$ of the lemmas above are satisfied 
if $K$ is {\em differential-valued\/} as defined in 
Section~\ref{As-Fields,As-Couples}.

%\begin{lemma} \label{valuation-of-basis}  Suppose that 
%$\mathcal{O}=C+\smallo$, and $A\in K[\der]$ has order $n$. 
%Then 
%\begin{enumerate}
%\item $v(\ker^{\ne} A)$ is finite, of size $\dim_C\ker A$;
%\item if $v(\ker^{\ne} A)= \{vy_1, \dots , vy_m\}$ where
%$y_1,\dots, y_m\in\ker^{\ne} A $, $vy_1 < \dots < vy_m$, 
%then $y_1,\dots,y_m$ is a basis of $\ker A$;
%\item $\gamma\in v(\ker^{\ne} A) 
%\Rightarrow \Pmu_A(\gamma)>0$.
%\end{enumerate}
%\end{lemma} 
%\begin{proof}  
%Items (1) and (2) follow from 
%Lemma~\ref{Hahn valuations, lemma}, 
%since $K$ with the valuation $v$ is a Hahn space over $C$. 
%As to 
%item (3), note that if $g\in \ker^{\ne} A$ 
%and $\gamma=vg$, then 
%$v_A(\gamma)= v(Ag) < \infty = v\bigl(A(g)\bigr)$, so  
%$\Pmu_A(\gamma)>0$.
%\end{proof}

%\noindent
%The assumptions $\der\mathcal{O}\subseteq \smallo$ 
%and $\mathcal{O}=C+\smallo$ of the lemmas above are 
%satisfied if $K$ is differential-valued as defined in 
%Section~\ref{As-Fields,As-Couples}. 
%Lemma~\ref{valuation-of-basis} is useful in solving an %equation $A(y)=0$: if we know where $\Pmu_A>0$, then 
%(3) helps to pin down the valuation of the potential %solutions.

\medskip
\noindent
For $A\in K[\der]^{\ne}$, $A=\sum_i a_i\der^i$, we define
$$\Pnu(A)\ :=\ \max\bigl\{i:v(a_i)=v(A)\bigr\}.$$
Then the analogues of Lemmas~\ref{mu-lemma} and \ref{nLM-formula} hold 
without the assumption~${\der\mathcal{O}\subseteq \smallo}$: 
$\Pnu(Ag)$ for $g\in K^{\times}$ depends only on $vg$, and setting
$\Pnu_A(\gamma) := \Pnu(Ag)$ when $\gamma=vg$, $g\in K^{\times}$, we have
$\Pnu_{AB}(\gamma)=\Pnu_A\bigl(v_B(\gamma)\bigr)+\Pnu_B(\gamma)$ for
$A, B \in K[\der]^{\ne}$ and $\gamma\in \Gamma$.
(Same proofs.)

\subsection*{Exceptional values}
\textit{In this subsection we let $\gamma$ range over $\Gamma$ and $y$ range over~$K^\times$.}\/
Let $A\in K[\der]^{\ne}$,
$A = a_0 + a_1\der + \cdots + a_n\der^n$, 
$a_0,\dots, a_n\in K$,
so $A(1)=a_0$.
Recall that $\Pmu(A)=\min\big\{i:\ v(a_i)=v(A)\big\}$, so $\Pmu(A) >0$ is equivalent to $A(1)\prec A$.
We call $\gamma$ an {\bf exceptional value} for $A$ if there is a $y$
such that $vy=\gamma$ and $A(y)\prec Ay$, and we let~$\exc(A)$ be the set of exceptional values for~$A$,
in other words
 $$\exc(A)\	:= \ \big\{ vy:\ A(y)\prec Ay  \big\} \ =\ \big\{vy:\ \Pmu(Ay)>0\big\}.$$ 
For $a\in K^\times$ we have
$$%\begin{equation}\label{E-mult, 2}
\exc(aA)\ =\ \exc(A), \qquad \exc(Aa)\ =\ \exc(A) - va.
$$%\end{equation} 

\begin{lemma}\label{lem:exc(A) and ker}
$v(\ker^{\ne} A) \subseteq \exc(A)$.
\end{lemma}
\begin{proof}
If $g\in \ker^{\ne} A$, then 
$A(g)=0\prec Ag$, so $vg\in\exc(A)$.
\end{proof}

\index{exceptional value}
\index{linear differential operator!exceptional value}
\nomenclature[Q]{$\exc(A)$}{set of exceptional values of $A$}

\noindent
Lemma~\ref{lem:exc(A) and ker} is useful in solving a differential equation
$A(y)=0$: knowing $\exc(A)$ helps to pin down the valuation of the potential solutions.

\begin{lemma}\label{asymptotic relations under A, lemma}
Let $f_1,f_2\in K^\times$ with $v(f_2)\notin\exc(A)$. Then:
$$f_1\asymp f_2\ \Rightarrow\ A(f_1)\asymp A(f_2), \qquad f_1\prec f_2\ \Rightarrow\ A(f_1)\prec A(f_2),$$  
and thus $f_1\sim f_2\ \Rightarrow\ A(f_1)\sim A(f_2)$.
\end{lemma}
\begin{proof}
If $f_1\asymp f_2$, then  $A(f_1)\asymp Af_1 \asymp Af_2 \asymp A(f_2)$,
and if $f_1\prec f_2$, then $A(f_1)\preceq Af_1\prec Af_2\asymp A(f_2)$. 
\end{proof}

\noindent
If $\der\mathcal O\subseteq\smallo$, then
$\exc(A) =\big\{ \gamma: \Pmu_A(\gamma) >0 \big\}$, so for $B\in K[\der]^{\ne}$,
by Lemma~\ref{nLM-formula}:
$$%\begin{equation}\label{E-mult, 1}
\exc(AB)\ =\ v_B^{-1}\big(\exc(A)\big)\cup\exc(B).
$$%\end{equation}

\noindent
Note that $\Gamma^\flat:=\{va:\ a'\prec a\}$ is a subgroup of $\Gamma$. By Lemma~\ref{mu_P, more} we have:

\begin{lemma} If $\der\mathcal{O}\subseteq \smallo$, then 
$\exc(A)$ is a union of cosets of $\Gamma^\flat$.
\end{lemma} 

\begin{example}
Suppose $\der\mathcal{O}\subseteq \smallo$ and  $\order(A)=1$. Then
$$vy\in\exc(A)\quad\Longleftrightarrow\quad (a_0/a_1)+y^\dagger\prec 1.$$
Hence either $\exc(A)=\emptyset$ or $\exc(A)$ is a coset of $\Gamma^\flat$. 

As a special case, the valued differential field $K=C(x)$ from
   Example~(3) at the beginning of 
   Section~\ref{Valdifcon} satisfies $y^\dagger\prec 1$
   for all $y$, so
   $\der\mathcal{O}\subseteq \smallo$. For this $K$ and $a\in K$, if  $a\prec 1$, then
   $\exc(\der+a)=\Gamma^\flat$, and if $a\succeq 1$, then $\exc(\der+a)=\emptyset$.
\end{example}

\noindent
To indicate the dependence of~$\exc(A)$ on $K$ we write it as $\exc_K(A)$.
If $L$ is a valued differential field extension of $K$ with small derivation, then $\exc_K(A)\subseteq\exc_L(A)$, and if in addition $\der\mathcal{O}_L\subseteq \smallo_L$, then 
$\exc_K(A)=\exc_L(A)\cap\Gamma$ by Lemma~\ref{mu-lemma}.

\subsection*{Linear surjectivity in the presence of a valuation}
%In this subsection $K$ is a valued differential field  with $\der \smallo \subseteq\smallo$. Let $\k=\mathcal{O}/\smallo$ be its differential residue field.  
Let $A\in K[\der]^{\ne}$. Recall that if $a\in K^{\times}$ and $\Pmu(Aa)=0$,
then $$v\big(A(a)\big)\ =\ v(Aa)\ =\ v_A(va).$$ 
We call $A$ {\bf neatly surjective} if for all $b\in K^\times$ there is 
$a\in K^\times$ such that $A(a)=b$ and $v_A(va)=vb$; note that then 
$A\colon K \to K$ and $v_A\colon \Gamma\to \Gamma$ are surjective.
In Chapter~\ref{ch:valueddifferential} we prove that $v_A\colon \Gamma\to \Gamma$ {\em is\/}
indeed surjective, but this is a rather difficult result, so we prefer to
add it explicitly as an assumption in some results below.  

\index{linear differential operator!neatly surjective}
\index{neatly surjective}

\begin{lemma}\label{neatlysur} Neat surjectivity has the following basic properties: \begin{enumerate}
\item[\textup{(i)}] each $\phi\in K^\times \subseteq K[\der]^{\ne}$ is neatly surjective;
\item[\textup{(ii)}] if $A,B\in K[\der]^{\ne}$ are neatly surjective, then so is $AB$;
\item[\textup{(iii)}] if $K$ is linearly closed and each $A\in K[\der]^{\ne}$ of order $1$
is neatly surjective, then every $A\in K[\der]^{\ne}$ is neatly surjective;
\item[\textup{(iv)}] if $A\in K[\der]^{\ne}$, $A(K)=K$, and $\Pmu(Aa)=0$ for all 
$a\in K^\times$, then $A$ is neatly surjective.
\end{enumerate}
\end{lemma} 

\begin{lemma}\label{nslinsur} If $K$ is $1$-linearly surjective, 
$\der$ is neatly surjective, and $\smallo\subseteq (K^\times)^\dagger$, then
every $A\in K[\der]$ of order $1$ is neatly surjective $($so $\k$ is
$1$-linearly surjective$)$.
\end{lemma}
\begin{proof} Let $A=\der + \phi\in K[\der]$, and suppose $A(K)=K$. 
If $\Pmu(Aa)=0$ for all~$a\in K^\times$, then $A$ is neatly surjective by Lemma~\ref{neatlysur}(iv).
So let $\Pmu(Aa)>0$, with~$a\in K^\times$.
Then $a^{-1}Aa=\der+ (a^\dagger + \phi)$ with $a^\dagger + \phi\in \smallo$, 
and so, assuming $\smallo\subseteq (K^\times)^\dagger$,  we have
$a^\dagger + \phi=b^\dagger$ with $b\in K^\times$, hence 
$a^{-1}Aa=b^{-1}\der b$. Assuming next that $\der$ is neatly
surjective, $A$ is neatly surjective by
Lemma~\ref{neatlysur}(i)(ii). 
\end{proof}

\noindent
If every $A\in K[\der]^{\ne}$ of order~$\le n$ is neatly surjective,
then $\k$ is $n$-linearly surjective. We have the following approximate 
converse to this fact: 

\begin{lemma}\label{apprsol} 
Suppose $\k$ is $n$-linearly surjective, $A\in K[\der]^{\ne}$ has
order $\le n$ and $v_A\colon \Gamma\to \Gamma$ is surjective. 
Then there is for each $b\in K^\times$ an $a\in K^\times$
such that $A(a) \sim b$ and $v_A(va)=vb$.
\end{lemma}
\begin{proof} Let $b\in K^\times$ and take 
$\alpha\in \Gamma$ with $vb = v_A(\alpha)$.
Take $\phi\in K^\times$ with $v\phi = \alpha$, so $B:= b^{-1}A\phi$ 
satisfies $v(B)=0$. Take $u\in K$ with $u\asymp 1$ such that
$B(u)\sim 1$ (possible since $\k$
is $n$-linearly surjective). Then $A(u\phi) \sim b$, so $a=u\phi$ works.
\end{proof}

\noindent
This suggests a way to ``neatly'' solve equations $A(y)=g$:

\begin{cor}\label{cor1628} Assume $K$ is spherically complete, $\k$ is $n$-linearly surjective, $A\in K[\der]^{\ne}$ has order $\le n$, and $v_A\colon \Gamma\to \Gamma$ is surjective. Then $A$ is neatly surjective.
\end{cor}
\begin{proof} Let $g\in K^\times$; we wish to find $f\in K^\times$ such that
$A(f)=g$ and $v_A(vf)=vg$. Take $f_0\in K^\times$ with $A(f_0)\sim g$ and $v_A(vf_0)=vg$. If $A(f_0)=g$, we are done. Suppose $A(f_0)\ne g$. Then we take
$y\in K^\times$ such that $A(y) \sim g-A(f_0)\prec g$ and $v_A(vy)=v\big(g-A(f_0)\big)$
and set $f_1:=f_0 + y$. Then $g-A(f_1) \prec g-A(f_0)$, and $f_0\sim f_1$. If
$A(f_1)=g$, we are done, and otherwise we continue as before, with $f_1$ instead of
$f_0$. In general, we have a sequence
$(f_{\lambda})_{\lambda < \rho}$ in $K^\times$, indexed by an ordinal~$\rho>0$, with $f_0$ as chosen initially, such that the following conditions hold: \begin{enumerate}
\item $v\big(g-A(f_{\lambda})\big)$ is strictly increasing as a function of
$\lambda$, 
\item $v_A\big(v(f_{\mu}-f_{\lambda})\big)\ =\ v\big(g- A(f_{\lambda})\big)$ for $\lambda < \mu < \rho$.
\end{enumerate}
In particular, $f_{\lambda} \sim f_0$ for all $\lambda$.
Consider first the case that $\rho=\nu+1$ is a successor ordinal. If $A(f_{\nu})=g$,
we are done, so assume $A(f_{\nu})\ne g$. The same way we got~$f_1$ from $f_0$ we take $f_{\nu +1}\in K^\times$ with 
$g-A(f_{\nu +1}) \prec g-A(f_{\nu})$ and ${v_A\big(v(f_{\nu+1}-f_{\nu})\big)=v\big(g-A(f_{\nu})\big)}$.
Then the extended sequence $(f_{\lambda})_{\lambda < \rho+1}$ has the above properties with~$\rho+1$ instead of $\rho$.
 
 Suppose $\rho$ is a limit ordinal. Then $(f_{\lambda})$ is a pc-sequence, so $f_\lambda\leadsto f_{\rho}\in K$.  Then the extended sequence $(f_{\lambda})_{\lambda < \rho+1}$ has the above properties with $\rho+1$ instead of $\rho$.
 Eventually, this building process must result in an $f$ as desired. 
\end{proof} 

\noindent
Neat surjectivity is preserved under adjoining $\imag$, as we explain now. 
Assume $-1$ is not a square in the residue field $\k$.
Then $K[\imag]$ is a valued differential field with valuation given by 
$v(a+b\imag )=\min(va, vb)$ for $a,b\in K$, so its valuation ring is 
$\mathcal{O} + \mathcal{O}\imag$, with maximal ideal $\smallo + \smallo \imag$ satisfying
$\der(\smallo + \smallo \imag)\subseteq \smallo + \smallo \imag$.
The value group of $K[\imag]$ is again $v(K^\times)=\Gamma$, and: 

\begin{lemma} Let $A\in K[\der]^{\ne}$. Then $A\colon K \to K$ is neatly surjective
if and only if $A\colon K[\imag] \to K[\imag]$ is neatly surjective.
\end{lemma}
\begin{proof} Assume first that $A\colon K \to K$ is neatly surjective. Let
$f\in K[\imag]^\times$ and take~${\gamma\in \Gamma}$ such that $v_A(\gamma)=vf$.
Then $f=g+h\imag$, with $g,h\in K$, and $vf=\min(vg, vh)$. Take $a,b\in K$ such that
$A(a)=g$, $v_A(va)=vg$, $A(b)=h$, $v_A(vb)=vh$. Then $\gamma=\min (va, vb)= v(a+b\imag )$, and $A(a+b\imag )=f$. Thus $A\colon K[\imag]\to K[\imag]$ is neatly surjective.

Next, assume $A\colon K[\imag] \to K[\imag]$ is neatly surjective. 
Let $f\in K^\times$, and
take $a,b\in K$ such that $A(a+b\imag)=f$ and $v_A(\gamma)=vf$ with $\gamma=\min(va, vb)$. Then $A(a)=f$ and $A(b)=0$. If $va \le vb$, then $va=\gamma$, and 
if $va > vb$, then $\gamma=v(a+b)$ and $A(a+b)=f$. 
\end{proof}

\subsection*{Independence and diagonalization in the presence of a valuation} In the next two lemmas $Y=(Y_1,\dots,Y_n)$ is a tuple of
distinct differential indeterminates, $n\ge 1$. These lemmas involve
independence notions defined earlier in this chapter.

\begin{lemma}\label{indepreduc} Suppose $A_1,\dots, A_m$ are in $\mathcal{O}\{Y\}_1$ and their images $\overline{A}_1,\dots, \overline{A}_m$ in $\k\{Y\}_1$ are $\d$-independent \textup{(}over $\k$\textup{)}. Then $A_1,\dots, A_m$ are $\d$-independent.
\end{lemma}
\begin{proof} This follows easily from the definition of
$\d$-dependence in Section~\ref{sec:systems}.
\end{proof}

\begin{lemma}\label{diagval}
Let $n\ge 1$ and let $A$ be an $n\times n$ matrix over $\mathcal{O}[\der]$ such that its image $\overline{A}$ as an $n\times n$ matrix over $\k[\der]$ has $\k[\der]$-independent rows.
Then there are
$S, T\in \operatorname{GL}_n\big(\mathcal{O}[\der]\big)$ such that $SAT=(D_{ij})$ 
with $D_{ij}\prec 1$ for all $i\ne j$ in $\{1,\dots,n\}$, and
$D_{ii}\asymp 1$ for $i=1,\dots,n$. 
\end{lemma}
\begin{proof} Diagonalize
$\overline{A}$ by row and column operations, and
lift each step to a row or column operation over $\mathcal{O}[\der]$. 
\end{proof}

\subsection*{Notes and comments} For archimedean $\Gamma$ one can find Lemma~\ref{Functoriality-Easy}(ii) and Lemma~\ref{nLM-formula} in Section~1.6 of~\cite{Robba}.

%% file: mt-5-7.tex
\section{Compositional Conjugation}\label{Compositional Conjugation}

\noindent
{\em In this section $K$ is a differential ring with derivation $\der$
and $\phi\in K^\times$}. We define~$K^\phi$ to be the 
differential ring with the same underlying ring 
as~$K$ but with the derivation~$\derdelta$ given by 
$\derdelta(a)= \phi^{-1}\cdot \der(a)$
for $a\in K$, so $\derdelta=\phi^{-1}\der$. This gives rise to the ring~$K^\phi\{Y\}$ of differential polynomials over $K^\phi$. 
Thus $K^\phi\{Y\}$ has the same underlying ring as~$K\{Y\}$ but its derivation extends $\derdelta=\phi^{-1} \der$. We denote this extended
derivation also by $\derdelta$, so
$\derdelta(Y^{(n)})=Y^{(n+1)}$ for all $n$.
For a differential polynomial $P\in  K^\phi\{Y\}$ written 
as an ordinary polynomial $P=p(Y,Y',\dots,Y^{(n)})\in K\big[Y,Y',\dots,Y^{(n)}\big]$,
we have 
$P(y)=p\bigl(y,\derdelta(y),\dots,\derdelta^n(y)\bigr)$ 
for $y\in K$.

\nomenclature[M]{$\derdelta$}{derivation $\derdelta=\phi^{-1}\der$ of  $K^\phi$}

\subsection*{Transformation formulas}
In order to relate the differential rings $K\{Y\}$ and 
$K^\phi\{Y\}$ we take $\derdelta$ as the element $\phi^{-1} \der$
of $K[\der]$ and express the
powers of $\der$ 
as linear combinations of powers of $\derdelta$ in the ring $K[\der]$, with $\phi^{(k)}=\der^k(\phi)$:
\begin{align*}
\der^1\ &=\ \phi \derdelta\\
\der^2\ &=\ \phi^2 \derdelta^2 + \phi' \derdelta\\
\der^3\  &=\ \phi^3 \derdelta^3 + 
3\phi\phi' \derdelta^2 + \phi'' \derdelta\\
\der^4\ &=\ \phi^4\derdelta^4 + 6\phi^2\phi'\derdelta^3+\big(4\phi\phi''+3 (\phi')^2\big)\derdelta^2 + \phi^{(3)}\derdelta \\
                      &\vdots \\
\der^n\  &=\  F^n_n(\phi) \derdelta^n + F^n_{n-1}(\phi) \derdelta^{n-1} + 
\dots + F^n_{1}(\phi) \derdelta
\end{align*}
where the differential polynomials $F_k^n(X)\in \Q\{X\}\subseteq K\{X\}$ \nomenclature[O]{$F^n_k$}{used in expressing $\der^n$ in terms of $\derdelta$} have nonnegative integer coefficients and are independent of $K$ and $\phi$:
$F_n^n(X)=X^n$ and  $F_1^n(X)=X^{(n-1)}$ for $n\ge 1$, and
we have the recursion formula 
\begin{equation}\label{recursion for F} F_k^{n+1}(X)\ =\ F_k^n(X)' + X F_{k-1}^n(X) \qquad (k=1,\dots,n)
\end{equation}
where $F_0^n:= 0$ for $n\ge 1$, to make the last identity true for $k=1$.
For later use we also set $F_0^0:= 1$.
So $F_k^n$ is homogeneous
of degree $k$ and isobaric of weight $n-k$ for $0\le k \le n$, and of order 
$n-k$ for $1\le k\le n$.
   The identities relating the powers of $\der$ and $\derdelta$ suggest that we consider the ring morphism 
$$P(Y) \mapsto  {P^\phi}(Y) \colon K\{Y\} \to  K^\phi\{Y\}$$ that is 
the identity on $K[Y]$ and sends $Y^{(n)}$, for each $n\ge 1$, to
$$(Y^{(n)})^{\phi}\ :=\ F^n_n(\phi) Y^{(n)} + F^n_{n-1}(\phi) Y^{(n-1)} + 
\dots + F^n_{1}(\phi) Y'.$$ 
In particular, $Y^\phi=Y,\quad (Y')^\phi=\phi Y', \quad (Y'')^\phi=\phi^2 Y''+\phi' Y'$.
Note that $P\mapsto P^\phi$ is a ring isomorphism: 
it maps each subring $K\big[Y,\dots,Y^{(n)}\big]$ of $K\{Y\}$ isomorphically
onto the same subring of $K^\phi\{Y\}$.
We call $K^\phi$ the {\bf compositional conjugate of~$K$ by~$\phi$} and $P^\phi$ the {\bf compositional conjugate of $P$ by~$\phi$.}
 
\nomenclature[M]{$K^\phi$}{compositional conjugate of $K$ by $\phi\in K^\times$}
\nomenclature[O]{$P^\phi$}{compositional conjugate of  $P$  by $\phi$}
\index{differential field!compositional conjugate}
\index{differential polynomial!conjugation!compositional}
\index{conjugation!compositional}

\begin{lemma} \label{comp-conjugation lemma}
Let $P\in K\{Y\}$ and $y\in K$. Then 
\begin{align*}
P(y)\ =\  {P^\phi}(y), 
&\qquad (\der P)^\phi\ =\ \phi\cdot\derdelta({P^\phi}),\\
({P^\phi})_{+y}\ =\ {(P_{+y})^\phi}, &\qquad 
({P^\phi})_{\times y}\ =\ {(P_{\times y})^\phi}.
\end{align*}
\end{lemma}
\begin{proof}
For $P=Y^{(n)}$ these identities follow from the definitions,
using induction on~$n$ if necessary. The general case reduces easily to this special case.
\end{proof}

\noindent
In view of the last two identities we let $P^\phi_{+y}$ denote both 
$({P^\phi})_{+y}$ and ${(P_{+y})^\phi}$, and let~$P^\phi_{\times y}$ denote both 
$({P^\phi})_{\times y}$ and ${(P_{\times y})^\phi}$.

\medskip\noindent
If $\theta$ is also a unit of $K$, then 
$K^{\phi\theta} = (K^\phi)^\theta$, and we have ring isomorphisms
$$\begin{matrix}
P \mapsto { P^\phi} & \colon & K\{Y\} &\longrightarrow 
& K^\phi\{Y\}\\
         Q \mapsto { Q^\theta} &\colon & 
       K^\phi\{Y\}& \longrightarrow & K^{\phi\theta}\{Y\}\\   
        P \mapsto {P^{\phi\theta}} &\colon & K\{Y\} 
&\longrightarrow &  K^{\phi\theta}\{Y\}
\end{matrix}
$$
and the last map is the composition of the preceding two: 
$$%\begin{equation}\label{eq:repeated compconj}
{P^{\phi\theta}}={({P^\phi})^\theta}, \quad P\in K\{Y\}.
$$%\end{equation}
For $P=Y^{(n)}$ this is easily checked by an induction on $n$, using the second 
identity of Lemma~\ref{comp-conjugation lemma}; the identity for all $P$ 
is then an easy consequence. Note also that for~$\phi=1$ we have $K^\phi=K$ and
$P^\phi=P$.

Another induction on $n$ using the second identity of 
Lemma~\ref{comp-conjugation lemma} 
yields 
$$\bigl(\der^n Q\bigr)^\phi\  =\  P^\phi\big(Q^\phi\big)\ \text{ for } 
P=Y^{(n)}\in K\{Y\} \text{ and all }Q\in K\{Y\}.$$ 
It follows that $P(Q)^\phi=P^\phi(Q^\phi)$
for all  $P,Q\in K\{Y\}$. 
Thus we have a ring isomorphism
$$A\mapsto A^\phi\ \colon\ K[\der]\to  K^\phi[\derdelta], \quad a^\phi= a \text{ for }a\in K, \quad \der^\phi=\phi\derdelta.$$
It sends $\der^n$ for $n\ge 1$ to 
$F^n_n(\phi) \derdelta^n + F^n_{n-1}(\phi) \derdelta^{n-1} + 
\dots + F^n_{1}(\phi) \derdelta$, and $A^{\phi}(y)=A(y)$
for all $A\in K[\der],\ y\in K$. We identify the rings
 $K[\der]$ and $K^\phi[\derdelta]$ via this isomorphism. Thus for
$A=a_0+a_1\der+\cdots+a_r\der^r\in K[\der]$
with $a_0,\dots,a_r\in K$ we have
\begin{align*} A=A^\phi\  &=\ b_0+b_1\derdelta+\cdots+
b_r\derdelta^r\in  K^\phi[\derdelta], \text{  where}\\
b_i\  &=\ \sum_{j=i}^r F^j_{i}(\phi) a_j \quad \text{for $i=0,\dots,r$.}
\end{align*}
Here are some easy consequences:

\begin{lemma}\label{linclsurcom} Let $K$ be a differential field
 and $r\in \N^{\ge 1}$. If $K$ is $r$-linearly closed \textup{(}respectively $r$-linearly surjective, respectively $r$-pv-closed\textup{)}, 
then so is $K^{\phi}$. 
\end{lemma}

\noindent
The following identities for the powers $\der^n$ in the ring 
$K^\phi[\derdelta]$ will also be useful:
\begin{align*}
\der^1\  &=\ \phi\cdot \derdelta \\
\der^2\  &=\ \phi^2\cdot\derdelta^2 + \phi\cdot\derdelta(\phi)\cdot\derdelta \\
\der^3\  &=\  \phi^3\cdot\derdelta^3 + 3\phi^2\cdot\derdelta(\phi)\cdot\derdelta^2 
+ \big(\phi^2\cdot\derdelta^2(\phi)+\phi\cdot \derdelta(\phi)^2\big)\cdot\derdelta \\
&\vdots \\
\der^n\  &=\  G^n_n(\phi)\cdot \derdelta^n + 
G^n_{n-1}(\phi)\cdot \derdelta^{n-1} + \dots + G^n_{1}(\phi)\cdot \derdelta
\end{align*}
where the differential polynomials $G_k^n(X)\in \Q\{X\}\subseteq K^\phi\{X\}$ 
have nonnegative integer coefficients and are independent of~$K$ and $\phi$:
$G^n_n(X) = X^n$ for $n\ge 1$, and
\begin{equation}\label{recursion for G}
G^{n+1}_k\ =\ X\cdot\big(\derdelta(G^n_{k})+ G^n_{k-1} \big) \qquad(k=1,\dots,n)
\end{equation}
where $G_0^n:= 0$ for $n\ge 1$, to make the last identity true for $k=1$.
For later use we also set $G_0^0:= 1$. The differential 
polynomial $G^n_k\in K^\phi\{X\}$ is 
homogeneous of degree $n$ and isobaric of weight $n-k$, and $G^n_k(\phi)=F^n_k(\phi)$, for $0\le k \le n$.

\nomenclature[O]{$G^n_k$}{used in expressing $\der^n$ in terms of $\derdelta$}

\begin{lemma}\label{compderderdelta} $\Q[\phi,\dots, \der^n(\phi), \phi^{-1}]=\Q[\phi,\dots, \derdelta^n(\phi), \phi^{-1}]$, as rings. Thus for a differential field $K$, the differential subfields 
$\Q\<\phi\>=\Q(\der^n(\phi): n=0,1,2,\dots)$ of~$K$ 
and $\Q\<\phi\>=\Q(\derdelta^n(\phi): n=0,1,2,\dots)$ of $K^{\phi}$ have the same underlying field. 
\end{lemma}

\noindent
In the rest of this subsection $P\in K\{Y\}$ has order at most $r\in \N$, and 
$\bsigma,\btau\in \{0,\dots,r\}^*$ have equal length. For 
$\bsigma \le \btau$, 
$\bsigma=\sigma_1\cdots\sigma_d,\ \btau=\tau_1\cdots\tau_d$ we define 
$$F^{\btau}_{\bsigma}\ :=\ F^{\tau_1}_{\sigma_1} \cdots F^{\tau_d}_{\sigma_d} \in \Q\{X\},$$
so $F^{\btau}_{\bsigma}$ is homogeneous of degree $\|\bsigma\|$, and isobaric of weight $\|\btau\|-\|\bsigma\|$; in particular, 
$F^\varepsilon_\varepsilon:=1$ for the the empty word $\varepsilon$ in 
$\{0,\dots,r\}^*$, and $F^{\btau}_{\bsigma}=0$ if $\tau_i> \sigma_i=0$ for some
$i\in \{1,\dots,d\}$.

\begin{lemma} \label{conj-lemma}
$\quad (P^\phi)_{[\bsigma]}\ =\  \sum_{\btau \geq \bsigma} F^{\btau}_{\bsigma}(\phi) P_{[\btau]}.$
\end{lemma}
\begin{proof}
For $\btau=\tau_1\cdots\tau_d$ we have
$$(Y^{[\btau]})^\phi\ =\ 
(Y^{(\tau_1)})^\phi \cdots (Y^{(\tau_d)})^\phi\ =\
\sum_{\bsigma\leq \btau} F^{\btau}_{\bsigma} (\phi) Y^{[\bsigma]}.$$
With $P=\sum_{\btau} P_{[\btau]}Y^{[\btau]}$, this gives
\begin{align*} P^\phi\ &=\ \sum_{\btau} P_{[\btau]}(Y^{[\btau]})^\phi\ =\ \sum_{\btau} P_{[\btau]}\sum_{\bsigma\le \btau} F^{\btau}_{\bsigma} (\phi) Y^{[\bsigma]}\\
&=\ \sum_{\bsigma}\left(\sum_{\btau \ge \bsigma}F^{\btau}_{\bsigma} (\phi) P_{[\btau]}\right)Y^{[\bsigma]},
\end{align*}
from which the lemma follows.
\end{proof}

%For example, if $\phi'=0$ then $(P^\phi)_{[\bomega]}=\phi^{\abs{\bomega}} P_{[\bomega]}$ for every $\bomega$.

\begin{cor}\label{corconj-lemma} $P$ and $P^\phi$ have the same order, degree, weight, and complexity. 
If~$P$ is homogeneous, then so is $P^\phi$. If $P$ is 
subhomogeneous, then so is $P^\phi$.
\end{cor}
\begin{proof}
The differential polynomials $(Y^{(n)})^\phi$ are homogeneous of degree $1$, so
if~$P$ is homogeneous of degree $d$, then $P^\phi$ is too. Considering
the homogeneous parts of~$P$, this
yields $\deg(P^\phi)=\deg(P)$. Similarly one shows that if $P$ is subhomogeneous
of subdegree $d$, then $P^\phi$ is too.
Lemma~\ref{conj-lemma} shows that if $\|\bsigma\| > \wt(P)$, 
then $(P^\phi)_{[\bsigma]}=0$.
Hence  $\wt(P^\phi)\leq \wt(P)$, and since $P=(P^\phi)^{\phi^{-1}}$,
we obtain  $\wt(P^\phi)=\wt(P)$. In the same way we get
$\order(P^\phi) =\order(P)$ and $\c(P^\phi)=\c(P)$.
\end{proof}

\subsection*{Compositional conjugation and upward shift} In this subsection we assume familiarity with Appendix~\ref{app:trans}. We
study here $K^\phi$ for $K=\mathbb{T}$ and special $\phi$.

\index{differential polynomial!upward shift}
\index{upward!shift}
\nomenclature[O]{$P{\uparrow}$}{upward shift of the differential polynomial $P$ over $\mathbb T$}

\begin{example}
Consider $\mathbb{T}$ with its usual derivation $\der=\frac{d}{dx}$, and set $t=\frac{1}{x}$. The {\em upward shift}\/ of $f=f(x)\in \mathbb{T}$ is defined by
$ f{\uparrow}\ :=\  f(\ex^x)$.
Then $f\mapsto f{\uparrow}$ 
is an isomorphism $\mathbb{T}^t = (\mathbb{T} ,x\der)\to
(\mathbb{T},\der)$ of differential fields, with inverse $f\mapsto f{\downarrow}=
f(\log x)$. Let $P(Y)\in \mathbb{T}\{Y\}$, and let
$P{\uparrow}$ be the differential polynomial in $\mathbb{T}\{Y\}$ obtained by
applying $f\mapsto f{\uparrow}$ to the coefficients of $P^{t}\in \mathbb{T}^t\{Y\}$. Thus $(Y'){\uparrow}=\ex^{-x}Y'$. We have
$P(y){\uparrow} = P{\uparrow}(y{\uparrow})$
for $y\in \mathbb{T}$. (It suffices to check this for $P=Y^{(n)}$.)
Note:
$$ (P_{+a}){\uparrow}\ =\ (P{\uparrow})_{+a{\uparrow}}, \quad  
(P_{\times a}){\uparrow}\ =\ (P{\uparrow})_{\times a{\uparrow}} \qquad(P\in\mathbb{T}\{Y\},\ a\in \mathbb{T}),$$
and $P\mapsto P{\uparrow}$ is an $\R[Y]$-algebra automorphism of
$\mathbb{T}\{Y\}$ agreeing with $f\mapsto f{\uparrow}$ on~$\mathbb{T}$. 
It is shown in \cite{JvdH} that upward shifting is a very useful tool 
for analyzing the zero sets of differential polynomials over $\mathbb{T}$. 
Although in an arbitrary differential field the operation 
$f\mapsto f{\uparrow}$ is not available, in cases of interest
we have a nonzero element~$x$ with $x'=1$, and then compositional conjugation 
by $t=\frac{1}{x}$ makes sense and is a substitute for upward shifting, 
especially for differential polynomials with constant 
coefficients. This is worked out in Section~\ref{appdifpol}, 
%Chapter~\ref{ch:triangular automorphisms}, 
where we relate compositional conjugation by such~$t$ to more
general compositional conjugations.
\end{example}

\noindent
In the rest of this subsection $x$ is a unit of $K$ with $x'=1$, we set 
$t=\frac{1}{x}$, and work in~$K^t$ whose derivation is $\derdelta:=x\der$, so
 $\derdelta(x)=x$ and $\derdelta(t)=-t$. 

\medskip\noindent
Let $n\brack k$ ($k=0,\dots,n$) be the (unsigned) Stirling
numbers of the first kind; see \cite[\S\S{}5.5, 6.3]{Comtet-Book}, \cite[\S{}6.1]{GKP}.
They are defined recursively by ${n\brack 0}:=0$ for $n\ge 1$, ${n\brack n}:=1$, 
$${n+1\brack k}={n\brack k-1}+n\cdot {n\brack k} \quad\text{for $k=1,\dots,n$.}$$
Thus ${n\brack k}\in \N^{\ge 1}$ for $k=1,\dots,n$. We also set ${n\brack k}:=0$ for $k>n$. See Table~\ref{tab:stirling first kind}. As the table suggests, ${n+1\brack n}=\binom{n+1}{2}$ and ${n+1\brack 1}=n!$, which is easily verified.

\index{Stirling!numbers of the first kind!unsigned}
\nomenclature[Cd]{$n\brack k$}{unsigned Stirling numbers of the first kind}

\begin{table}[h]
\begin{center}
\begin{tabular}{c|llllll{c}}
\backslashbox[1em]{$k$}{$n$}
& 0 & 1 & 2 & 3 & 4 & 5 &\\
\hline
$0$		& $1$	& $0$	& $0$	& $0$	&  $0$ 	& $0$  	& \dots \\  
$1$		&     	& $1$  	& $1$  	& $2$  	&  $6$ 	& $24$	& \dots \\
$2$		&     	&      	& $1$  	& $3$  	&  $11$	& $50$	& \dots \\
$3$		&     	&      	&      	& $1$  	&  $6$	& $35$	& \dots \\
$4$		&     	&      	&      	&      	&  $1$	& $10$	& \dots \\
$5$		&     	&      	&      	&      	&		& $1$	& \dots \\
$\vdots$	& 		&      	&      	&      	&      	&    	& $\ddots$ 
\end{tabular}

\bigskip
\end{center}
\caption{Unsigned Stirling numbers of the first kind.}\label{tab:stirling first kind}
\end{table}

\noindent
For $k,n\geq 0$ we let $s(n,k):=(-1)^{n-k}{n\brack k}$, so
$s(n,0)=0$ for $n\ge 1$, $s(n,n)=1$, and 
\begin{equation}\label{eq:recurrence for s(n,k)}
s(n+1,k)\ =\ s(n,k-1)-n\cdot s(n,k)\qquad\text{for $k=1,\dots,n$.}
\end{equation}
The $s(n,k)$ are 
known as the signed Stirling numbers of the first kind.

\index{Stirling!numbers of the first kind!signed}
\nomenclature[Cd]{$s(n,k)$}{signed Stirling numbers of the first kind}

\begin{lemma}\label{lem:Stirling}
For $k=0,\dots,n$ we have
\begin{equation}\label{Stirling}
F^n_k(t)\ =\ s(n,k)\cdot t^{n},
\end{equation}
and the sum of the coefficients of $G^n_k$ equals ${n\brack k}$.
\end{lemma}

\begin{proof}
An easy induction on $n$ using \eqref{recursion for F} shows \eqref{Stirling}.
Since $\derdelta(t)=-t$ and each differential polynomial~$G^n_k$ is isobaric of weight $n-k$ and homogeneous of degree $n$, we have $G^n_k(t)=(-1)^{n-k}c_{n,k}t^n$ where $c_{n,k}=\text{sum of the coefficients of $G^n_k$}$. Since $G^n_k(t)=F^n_k(t)$ we get 
$c_{n,k}={n\brack k}$.
\end{proof}

\begin{exampleNumbered}\label{ex:example for G}
For $n\ge 1$ we have $G^n_{n-1} = {n\brack n-1}\, X^{n-1}\, X'$, and thus
%$$G^n_{n-1} = c(n,n-1)\, X^{n-1}\, X'.$$
$$  F^n_{n-1}(\phi)\  = {n\brack n-1}\, \phi^{n-1}\, \phi^\dagger.$$
To see this, use \eqref{recursion for G} to show that $G^n_{n-1}$ is an integer multiple of $X^{n-1}\, X'$, and then apply the previous lemma.
\end{exampleNumbered}

\noindent
For $\btau=\tau_1\cdots\tau_d\geq \bsigma=\sigma_1\cdots\sigma_d$ we have
$$F^{\btau}_{\bsigma}(t) =  s(\btau,\bsigma) \cdot t^{\|\btau\|}
\qquad\text{where $s(\btau,\bsigma) :=
s(\tau_1,\sigma_1) \cdots s(\tau_d,\sigma_d)$,}$$
hence
$$%\begin{equation}\label{eq:uparrow formula}
(P^t)_{[\bsigma]}\ =\ \sum_{\btau\geq\bsigma} s(\btau,\bsigma) t^{\|\btau\|} P_{[\btau]}
$$%\end{equation}
by Lemma~\ref{conj-lemma}. 

\subsection*{Compositional conjugation in the algebraic closure of $K\<Y\>$} The material in this subsection and the next one will
only be used in Section~\ref{eveqth}.  

Let $K$ be a differential field.
Since  $P\mapsto P^\phi\colon K\{Y\}\to K^\phi\{Y\}$ is a ring isomorphism, it extends uniquely to a
field isomorphism  $P\mapsto P^\phi\colon K\<Y\>\to K^\phi\<Y\>$. Fix an algebraic closure $F=K\<Y\>^\alg$ of the field~$K\<Y\>$.
Since $K\<Y\>=K^\phi\<Y\>$ as \textit{fields,}\/ this isomorphism $K\<Y\>\to K^\phi\<Y\>$ may be viewed as an automorphism of the field $K\<Y\>$, and we extend it to
an automorphism $P\mapsto P^\phi$ of $F$. The derivation $\der$ of~$K\<Y\>$ and the derivation $\derdelta$ of $K^\phi\<Y\>$ 
both extend uniquely to derivations of $F$, also denoted by~$\der$ and~$\derdelta$, respectively,
with the same constant field $C_F=C^\alg$. Note that $P\mapsto P^\phi\colon F \to F$ is the identity on $K$, and maps $C_F$ onto itself. 

\medskip
\noindent
The second identity in Lemma~\ref{comp-conjugation lemma} continues to hold if $P$ is an element of $F$:

\begin{lemma}
For each $P\in F$, we have $(\der P)^\phi=\phi\cdot\derdelta(P^\phi)$.
\end{lemma}
\begin{proof} Let $\sigma$ be the automorphism $P\mapsto P^\phi$ of $F$. 
The derivation $\phi\derdelta$ of $F$ yields a derivation 
$\sigma^{-1}\circ (\phi\derdelta)\circ \sigma$ on $F$ that agrees with 
$\der$ on $K\{Y\}$
by the second identity in Lemma~\ref{comp-conjugation lemma}. Hence
$\der=\sigma^{-1}\circ (\phi\derdelta)\circ \sigma$, and thus
$\sigma\circ\der=(\phi\derdelta)\circ \sigma$. 
\end{proof} 

\begin{cor}\label{cor:compconj and composition}
Let $P\in K\{Y\}$ and $Q\in F$. Then $P(Q)^\phi=P^\phi(Q^\phi)$.
\end{cor}
\begin{proof}
An  induction on $n$ using the previous lemma yields the identity
$(\der^n Q)^\phi=(Y^{(n)})^\phi(Q^\phi)$ and this gives
what we want. 
\end{proof}

\begin{lemma}\label{lem:composition in K<Y>a}
Suppose $Q\in F$ is transcendental over $K$, that is, $Q\notin K^{\alg}\subseteq F$. Then $Q$ is even $\d$-transcendental over $K$.
\end{lemma}
\begin{proof} Since $Y$ is $\d$-transcendental over $K$, it
is enough to show that $Y\in F$ is $\d$-algebraic over $K\<Q\>$.
As $Q$ is algebraic over $K\<Y\>$, we can take $R_0,\dots,R_n\in K\{Y\}$
with $R_n\neq 0$, such that
$$R_n Q^n+R_{n-1}Q^{n-1}+\cdots +R_0\ =\ 0.$$
Put 
$$R\ :=\ R_n(Z)Q^n+R_{n-1}(Z)Q^{n-1}+\cdots+R_0(Z)\in K\<Q\>\{Z\},$$ 
where $Z$ is a differential indeterminate 
different from $Y$. Take $r$ with $R_0,\dots, R_n$ of order $\le r$ and take $\i\in\N^{1+r}$ with $(R_n)_{\i}\neq 0$.
Then in $K\<Q\>$ we have the equality
$$R_{\i}\ =\ (R_n)_{\i}Q^n+(R_{n-1})_{\i}Q^{n-1}+\cdots+(R_0)_{\i},$$
so $R_{\i}\ne 0$ since $Q$ is transcendental over $K$. Hence $R\neq 0$ in  $K\<Q\>\{Z\}$, and ${R(Y)=0}$. Thus
$Y$ is $\d$-algebraic over $K\<Q\>$. 
\end{proof}

\subsection*{Compositional conjugation and rational powers} Let $K$ be a differential field. As before, we fix an algebraic 
closure $F$ of $K\<Y\>$, taking it as a differential field extension of
$K\<Y\>$, and extend the ring isomorphism $P\mapsto P^\phi\colon K\{Y\}\to K^\phi\{Y\}$ to an automorphism $P\mapsto P^\phi$ of the field $F$. 

\begin{notation}
For elements $f$ and $g$ of a differential field $E$ we define 
$$f =_{\operatorname{c}} g\ :\Longleftrightarrow\ f=c\cdot g \text{ for some $c\in C_E^\times$,}$$
so if $f,g\ne 0$, then $f =_{\operatorname{c}} g\Longleftrightarrow f^\dagger=g^\dagger$.
\end{notation}

\nomenclature[M]{$f =_{\operatorname{c}} g$}{$f=c\cdot g$ for some constant $c\neq 0$}

\noindent
We use this for $E:=F$.
For $P,Q\in F$ we have $P =_{\operatorname{c}} Q \Rightarrow P^\phi =_{\operatorname{c}} Q^\phi$.

\medskip
\noindent
Next, we extend the usual power map
$$(P,k)\mapsto P^k\colon F^\times\times\Z  \to F^\times$$
of the multiplicative group $F^\times$ to a map
$$(P,q)\mapsto P^q\colon F^\times\times\Q \to F^\times$$
such that 
$$(P^{q})^\dagger=qP^\dagger\qquad\text{for all $q\in\Q$ and $P\in F^\times$.}$$
Note that  
$$(P^{k/l})^l =_{\operatorname{c}} P^k\qquad\text{for $k,l\in\Z$, $l\neq 0$.}$$ 
The following rules are also easy to verify, for $P,Q\in F^\times$ and $q,q_1,q_2\in\Q$:
\begin{enumerate}
\item $P =_{\operatorname{c}} Q\ \Rightarrow\ P^q =_{\operatorname{c}} Q^q$;
\item $P^{q_1}P^{q_2}\ =_{\operatorname{c}}\ P^{q_1+q_2}$;
\item $(P^{q_1})^{q_2}\ =_{\operatorname{c}}\ P^{q_1  q_2}$;
\item $(PQ)^q\ =_{\operatorname{c}}\ P^q Q^q$.
\end{enumerate}
We allow of course $\phi\in \Q^\times\subseteq K^\times$, but a ``power'' 
$P^q$ should not be confused with the ``compositional conjugate'' 
$P^\phi$ of $P\in F^\times$ by $\phi$. 
Powers behave as follows under compositional conjugation:

\begin{lemma}\label{lem:powers and compconj}
Let $P\in F^\times$, $q\in\Q$. Then $(P^q)^\phi =_{\operatorname{c}} (P^\phi)^q$.
\end{lemma}
\begin{proof} For $k\in \Z$ we have $(P^k)^\phi = (P^\phi)^k$
since $P^k$ has its usual meaning in the multiplicative group $F^\times$.
Take $n\ge 1$ such that $nq=k\in \Z$, and set $Q:= (P^q)^\phi,\ R:=  (P^\phi)^q$.
Then $(P^q)^n =_{\operatorname{c}} P^k$ gives 
$$ Q^n\ =\  \big((P^q)^n\big)^{\phi}\  =_{\operatorname{c}}\ (P^k)^{\phi}\ =\ 
(P^\phi)^k\  =_{\operatorname{c}}\ R^n,$$
so $Q =_{\operatorname{c}}\ R$. 
\end{proof}

\noindent
Recall that by convention $G^0_0=1$ and $G^n_0=0$ for $n\ge 1$. We have
\begin{equation}\label{ynphi} (Y^{(n)})^\phi\ =\ \sum_{m=0}^n G^n_m(\phi)Y^{(m)} \in K^\phi\{Y\}.\end{equation}
Let $q\in \Q$ and suppose $\phi^q\in K$. Then in $K^\phi\{Y\}$,
$$(Y^{(m)})_{\times \phi^q}\  =\  \sum_{i=0}^m {m\choose i}\derdelta^{m-i}(\phi^q)Y^{(i)}, $$
and thus by identity  \eqref{ynphi} we have in $K^{\phi}\{Y\}$,
\begin{equation}\label{yncphi} (Y^{(n)})^\phi_{\times \phi^q}\ =\ \sum_{i=0}^n \left(\sum_{m=i}^n {m\choose i} G^n_m(\phi)\derdelta^{m-i}(\phi^q)\right)Y^{(i)}.\end{equation}
%Next we consider in more detail the case $q=1$, with $\phi^1=\phi\in K^\times$.

\medskip\noindent
We consider in more detail the case $q=1$, with $\phi^1=\phi\in K^\times$, since
Section~\ref{cutsvalgrp} requires us to deal with
$P^\phi_{\times\phi}$ where compositional and multiplicative conjugation are 
combined. 
The key to this is the fact that 
$P^\phi_{\times\phi}(Y')=P(Y')^{\phi}$. More formally, for $P\in K\{Y\}^{\ne}$, set \nomenclature[O]{$P^{\times \phi}$}{$P^{\times 1,\phi}$} $P^{\times \phi}:=P^\phi_{\times\phi}\in K^\phi\{Y\}$, and $P^{\times}:=P(Y')\in K\{Y'\}$. \nomenclature[O]{$P^{\times}$}{$P(Y')$} Then
\begin{equation}\label{eq:Ptimesl}
(P^{\times})^\phi\ =\ P^\phi\big( (Y')^\phi \big)\ =\ 
P^\phi( \phi Y')\ =\  P^{\times \phi}(Y')
\end{equation}
in $K^\phi\<Y\>$. Setting \nomenclature[O]{$R^n_k$}{$G^{n+1}_{k+1}$}
$$R^n_i(Y)\ 	:=\  \sum_{m=i}^n {m\choose i} G^n_m(Y) \, Y^{(m-i)}\in \Q\{Y\}\subseteq K^\phi\{Y\} \qquad (i=0,\dots,n),$$
the identity \eqref{yncphi} gives the following identity in $K^{\phi}\{Y\}$:
$$(Y^{(n)})^{\times \phi}\ =\  R^n_0(\phi) Y + R^n_1(\phi) Y' + \cdots + R^n_n(\phi) Y^{(n)}.$$
We have $Y^{(n+1)}=(Y^{(n)})^{\times}$ and hence by \eqref{eq:Ptimesl} applied to $P=Y^{(n)}$,
$$(Y^{(n+1)})^\phi\ =\ (Y^{(n)})^{\times\phi}(Y')\ =\  R^n_0(\phi) Y' + R^n_1(\phi) Y'' + \cdots + 
R^n_n(\phi) Y^{(n+1)}.$$
In view of \eqref{ynphi} this gives $R^n_i(\phi)=G^{n+1}_{i+1}(\phi)$ for $i=0,\dots,n$. Lemma~\ref{compderderdelta} shows that $\phi\in K^{\phi}$ is $\d$-transcendental over $\Q$ for suitable $K, \phi$, and thus
$$%\begin{equation}\label{eq:Rni=Gn+1,i+1}
R^n_i\ =\ G^{n+1}_{i+1}\ \text{ in }\Q\{Y\} \qquad\text{for $i=0,\dots,n$.}
$$%\end{equation}
Next, let $q\in \Q$ be such that $\phi^q\in K$. For $P\in K\{Y\}^{\ne}$ we set
$$ P^{\times q,\phi}\ :=\ P^\phi_{\times\phi^q}\in K^\phi\{Y\}, \qquad P^{\times q}\ :=\  P\big( (Y')^q \big)\in F.$$ 
Now suppose that $P\in K\{Y\}^{\ne}$ is homogeneous; let
$d=\deg(P)$ and $w=\wt(P)$. 
Then Corollary~\ref{cor:derivatives of powers} gives
$$P^{\times q}\ =_{\operatorname{c}}\ (Y')^{dq-w}\cdot E(Y')$$
with homogeneous $E\in K\{Y\}^{\neq}$, $\deg(E)=\wt(E)=w$. 
Computing in $F$ gives 
\begin{align*}
(P^{\times q})^\phi\ 	&=\ \, P^\phi\big( ((Y')^q)^\phi \big) \\
			&=_{\operatorname{c}}\  P^\phi \big( ((Y')^\phi)^q  \big) \\
					&=\ \, P^\phi \big( (\phi Y')^q \big) \\
			&=_{\operatorname{c}}\  P^\phi \big( \phi^q (Y')^q \big)\ 
			=\  P^{\times q,\phi}\big( (Y')^q \big),
\end{align*}
using Corollary~\ref{cor:compconj and composition} for the first equality; for the second, use the homogeneity of~$P^{\phi}$ and Lemma~\ref{lem:powers and compconj}. For future reference we summarize:

\nomenclature[O]{$P^{\times q,\phi}$}{$P^\phi_{\times\phi^q}$}
\nomenclature[O]{$P^{\times q}$}{$P\big((Y')^q\big)$}

\begin{lemma}\label{lem:Ptimesqphi}
Let $P\in K\{Y\}^{\ne}$ be homogeneous, 
$d=\deg(P)$, $w=\wt(P)$, $q\in \Q$. Then there is a homogeneous  $E\in K\{Y\}^{\neq}$
with $\deg(E) = \wt(E) = w$ such that for all $\phi\in K^\times$ with
$\phi^q\in K$ we have \textup{(}in $F$\textup{)}:
$$P^{\times q,\phi}\big( (Y')^q \big)\  =_{\operatorname{c}}\ (\phi Y')^{dq-w}\cdot E^{\times \phi}(Y'). $$
\end{lemma}

\noindent
Suppose $q\in\N^{\geq 1}$. If $P\in K\{Y\}^{\neq}$ is homogeneous, then $P^{\times q}=P\big((Y')^q\big)\in K\{Y\}$ is homogeneous with $\deg(P^{\times q})=q\,\deg(P)$, by Corollary~\ref{degundercomp}.  Hence for arbitrary $P\in K\{Y\}^{\neq}$ and $d\in\N$ we have $(P^{\times q})_{qd} = (P_d)^{\times q}$.

\subsection*{Notes and comments}
An explicit formula for the transformation coefficients $F^n_k$ appears in \cite[Theorem~1]{Todorov}
and a formula for the $G^n_k$ in \cite[(8)]{Comtet}. (Combining Propositions~\ref{prop:explicit fm for Bell polys} 
and \ref{prop:itmatrix F} below also gives a closed form expression for the~$F^n_k$.)
The formula for $F^n_k$ in the special case of Lemma~\ref{lem:Stirling}  dates back to at least
the 1823 dissertation of Scherk~\cite{Scherk}; see~\cite[Appendix]{BF}.
The Stirling number $n\brack k$ counts the number of permutations of $n$ objects with $k$ disjoint cycles; these numbers were introduced by Stirling~\cite[p.~11]{Stirling} in 1730.

The $K$-algebra underlying both $K\{Y\}$ and $K^\phi\{Y\}$ is $A=K[Y,Y',\dots]$, and $P\mapsto P^\phi$ is an automorphism of $A$; in Chapter~\ref{ch:triangular automorphisms} we study this automorphism. 
%operation $P\mapsto P^\phi$ as an automorphism of the $K$-algebra $A$.
%If $C\neq K$, then the map sending $\phi\in K^\times$ to $P\mapsto P^\phi$ is not a group morphism $K^\times\to\Aut(A)$; in fact, the identity \eqref{eq:repeated compconj} expresses that the map that associates to $(\phi,\theta)\in K^\times\times K^\times$ the  automorphism
%$P\mapsto P^\phi\colon K^\theta\{Y\}\to K^{\phi\theta}\{Y\}$ of %$A$ is a {\it cocycle}\/ for the action of $K^\times$ on itself %by multiplication, taking values in $\Aut(A)$.

%% file: mt-5-8.tex
\section{The Riccati Transform}\label{The Riccati Transform}

\noindent
In this section $K$ is a differential 
{\em field}, $y$ ranges over $K^{\times}$, and $z:=y^\dag$.
Then $$\frac{y^{(n)}}{y}\ =\  R_n(z)$$ where the
differential polynomial $R_n(Z)\in K\{Z\}$ 
has nonnegative integer coefficients and is independent of $K$ and $y$:
\begin{align*}
         R_0(Z)\ &=\ 1\\
         R_1(Z)\ &=\ Z\\
         R_2(Z)\ &=\ Z^2 + Z'\\
         R_3(Z)\ &=\ Z^3 + 3ZZ' + Z''\\
                 &\ \ \vdots
\end{align*}
These $R_n$ are defined by the recursion 
\begin{equation}\label{Definition of R_n}
R_0(Z)\ :=\ 1, \qquad R_{n+1}(Z)\ :=\ ZR_n(Z) + R_n(Z)'.
\end{equation}
An easy induction yields $R_n(Z) = Z^n + A_n(Z)$ with $\deg(A_n) \leq n-1$ and where every monomial in $A_n$ has the form $Z^{i_0}(Z')^{i_1}\cdots (Z^{(k)})^{i_k}$ with $1\le k \leq n-1,\ i_k\geq 1$. 
Note that $R_2(Z) = -\frac{1}{4}\omega(2Z)$ where $\omega(Z)=-(Z^2+2Z')$ is as in Section~\ref{sec:secondorder}.
The~$R_n$ are related to the $F^n_k$ from Section~\ref{Compositional Conjugation} as follows:

\index{Riccati polynomial}
\nomenclature[O]{$R_n$}{$n$th Riccati polynomial}

\begin{lemma} $R_n(Z) = F^n_n(Z) + F^n_{n-1}(Z) + \cdots + F^n_0(Z)$, so the homogeneous and isobaric parts of $R_n$ are given by $R_n(Z)_k=R_n(Z)_{[n-k]}=F^n_k(Z)$ for $k=0,\dots,n$.
\end{lemma}
\begin{proof} Let $Y$ be a differential indeterminate over $\Q$, let $\der$ be the usual derivation of~$K\<Y\>$, set $\phi:= Y^\dagger\in \Q\<Y\>$, and 
take $\derdelta=\phi^{-1}\der$.  Then $\derdelta^k(Y)=Y$ for all $k\in\N$, and thus
$$R_n(\phi)\ =\ \frac{Y^{(n)}}{Y}\ =\ \frac{\sum_{k=0}^n F^n_k(\phi)\,\derdelta^k(Y)}{Y}\ =\ \sum_{k=0}^n F^n_k(\phi).$$
It remains to note that $\phi$ is differentially transcendental over $\Q$. 
\end{proof}

\noindent
Hence $R_n(Z)_{[0]}=F^n_n(Z)=Z^n$, so the $R_n(Z)$ are linearly independent over $K$. If~$n\ge 1$, then $R_n(Z)_1=R_n(Z)_{[n-1]}=F^n_1(Z)=Z^{(n-1)}$, so
$R_n$ has order $n-1$.

\medskip
\noindent
For each $n$ we have the following binomial identity in the 
differential polynomial ring~$\Q\{ Z_1,Z_2 \}$:
\begin{equation}\label{R_n-property 1}
R_n(Z_1 + Z_2)\ =\ \sum_{i+j=n}\binom{n}{i}R_i(Z_1)R_j(Z_2).
\end{equation}
To see this, note that in the differential field $\Q\<Y_1,Y_2\>$,
\begin{multline*}
R_n\bigl(Y_1^\dag + Y_2^\dag\bigr)\
=\  R_n\bigl((Y_1Y_2)^\dag \bigr)\
=\ \frac{(Y_1Y_2)^{(n)}}{Y_1Y_2}\
=\ \frac{\sum_{i=0}^n\binom{n}{i}Y_1^{(i)}Y_2^{(n-i)}}{Y_1Y_2}\\ 
=\  \sum_{i+j=n}\binom{n}{i}\frac{Y_1^{(i)}}{Y_1}\frac{Y_2^{(j)}}{Y_2}\
=\ \sum_{i+j=n}\binom{n}{i}R_i(Y_1^\dag)R_j(Y_2^\dag).
\end{multline*}

\noindent
By induction on $m$ this gives:

\begin{cor}\label{binomric} For all $z_1,\dots, z_m \in K$ we have
$$ R_n(z_1+ \cdots +z_m)\ =\ \sum\frac{n!}{i_1!\cdots i_m!}R_{i_1}(z_1)\cdots R_{i_m}(z_m)$$
where the sum on the right is over all tuples $(i_1,\dots, i_m)\in \N^m$ with $i_1+ \cdots +i_m=n$.
\end{cor}

\begin{definition}\label{Riccati-Def}
Let $\Ric\colon K\{Y\}\to K\{Z\}$ be the $K$-algebra homomorphism
such that $\Ric\big(Y^{(n)}\big)= R_n(Z)$ for all $n$. 
We call $\Ric(P)$ the {\bf Riccati transform} of $P$.
If $P\in K\{Y\}$ is homogeneous of
degree $d$, then
$$\Ric(P)(z)\ =\ \frac{P(y)}{y^d}.$$
\end{definition}

\index{Riccati transform!differential polynomial}
\index{differential polynomial!Riccati transform}
\nomenclature[O]{$\Ric(P)$}{Riccati transform of  $P$}

\begin{remark}
If $P\in K\{Y\}^{\ne}$ is homogeneous of order $r\ge 1$, then
$$\Ric(P)\ne 0, \quad \operatorname{order} \Ric(P)\ =\ r-1, \quad \deg \Ric(P)\ \le\ \wt P\ \le\  r\cdot \operatorname{sdeg} P.$$ 
\end{remark}

\noindent
A straightforward computation yields:

\begin{lemma}
If $P\in K\{Y\}$ is homogeneous of degree $d$ and $g\in K$,  then
$$\Ric(P)_{+g}(Z)\ =\ \sum_{\bsigma} \left(\sum_{\btau\geq\bsigma}
P_{[\btau]}\binom{\btau}{\bsigma}R_{\btau-\bsigma}(g)\right)R_{\bsigma}(Z).$$
Here $\bsigma=\sigma_1\cdots\sigma_d$ and 
$R_{\bsigma}:=R_{\sigma_1}\cdots R_{\sigma_d}=\Ric(Y^{[\bsigma]})$.
\end{lemma}

\begin{cor}\label{Riccati-Mult}
Let $P\in K\{Y\}$ be homogeneous of degree $d$ and let $h\in K^\times$. Then
$$\Ric(P_{\times h})\ =\ h^d\cdot\Ric(P)_{+h^\dagger}.$$
\end{cor}
\begin{proof}
With $g:=h^\dagger$ and $\bsigma$ ranging over words 
$\sigma_1\cdots\sigma_d$ of length $d$ we have
\begin{align*}
h^d\cdot \Ric(P)_{+g}(Z)\ &=\ 
\sum_{\bsigma} \left(\sum_{\btau\geq\bsigma}
P_{[\btau]}\binom{\btau}{\bsigma}h^d 
R_{\btau-\bsigma}(g)\right)R_{\bsigma}(Z)\\ 
&=\
\sum_{\bsigma} \left(\sum_{\btau\geq\bsigma}
P_{[\btau]}\binom{\btau}{\bsigma}h^{[\btau-\bsigma]}\right)R_{\bsigma}(Z)\ =\ 
\Ric(P_{\times h})(Z)
\end{align*}
by the previous lemma and \eqref{Multiplicative-Conj}.
\end{proof}

\noindent
Here is how the Riccati transform interacts with compositional conjugation:

\begin{lemma}\label{rica1}\label{Ri and compositional conjugation} 
Let $\phi\in K^\times$ and $P\in K\{Y\}$. Then
$$  \Ric(P^\phi)\ =\ \Ric(P)^\phi_{\times\phi}\ \text{ in $K^\phi\{Z\}$.}$$
\end{lemma}
\begin{proof} We can reduce to the case that $P$ is homogeneous of degree $d$.
Then for $y\in K^\times$ and $z=y^\dagger$, and setting $\derdelta=\phi^{-1}\der$,
we have
\begin{align*} \Ric(P)^\phi_{\times\phi}\big(\derdelta(y)/y\big)\ &=\ \Ric(P)^\phi(z)\ =\ \Ric(P)(z)\\
&=\ \frac{P(y)}{y^d}\ =\ \frac{P^\phi(y)}{y^d}\ =\ \Ric(P^\phi)\big(\derdelta(y)/y\big).
\end{align*}
This remains true for $y$ in differential field extensions of $K^\phi$, and so
for all $a$ in all such extensions we have 
$\Ric(P^\phi)(a)=\Ric(P)^\phi_{\times\phi}(a)$.
\end{proof}

\subsection*{The Riccati transform of a linear differential operator}
Let
$$A\ =\ a_0 + a_1\der + \dots + a_n\der^n\in K[\der] 
\qquad (a_0,\dots,a_n\in K).$$ 
The {\bf Riccati transform} $\Ric(A)$ of $A$ is defined to be
the Riccati transform of the corresponding differential polynomial: 
$$\Ric(A)(Z)\ :=\ a_0R_0(Z) + a_1R_1(Z) + \dots + a_nR_n(Z)\in K\{Z\}.$$ 
Using \eqref{R_n-property 1} and \eqref{A^(i)} we have for $g\in K$:
$$\Ric(A)_{+g}(Z)\ =\ b_0R_0(Z) + b_1R_1(Z) + \cdots + b_nR_n(Z),$$
where
$$b_i\ =\ \sum_{j=i}^{n} \binom{j}{i} a_{j} R_{j-i}(g)\ =\  
\frac{1}{i!}\Ric(A^{(i)})(g).$$
In particular, $b_0=\Ric(A)(g)$ and $b_n=a_n$.

\begin{lemma} Let $A\in K[\der]$, set $R:= \Ric(A)$, and let $a\in K$. Then 
$$ A\in K[\der](\der-a)\ \Longleftrightarrow\ R(a)=0.$$
\end{lemma}
\begin{proof} Immediate from Lemma~\ref{linear factor}.
\end{proof}  

\index{Riccati transform!linear differential operator}
\index{linear differential operator!Riccati transform}
\nomenclature[Q]{$\Ric(A)$}{Riccati transform of $A$}

\noindent
We define a {\bf Riccati polynomial} \index{Riccati polynomial} over $K$ to be a differential polynomial
$$R(Z)=a_0R_0(Z) + a_1R_1(Z) + \dots + a_nR_n(Z) 
\qquad (a_0,\dots,a_n\in K),$$
that is, a Riccati polynomial over $K$ is the Riccati transform of some
$A\in K[\der]$. It follows that if $R$ is a Riccati polynomial over 
$K$ and $g\in K$, then $R_{+g}$ is a Riccati polynomial over 
$K$ of the same order.

\begin{lemma}\label{lem:twist coeffs} Let $A\in K[\der]$, and let
$h$ be an element of a differential field extension of $K$ with $h\neq 0$ and 
$h^\dagger=g\in K$. Then
$A_{\ltimes h}\in K[\der]$.
\end{lemma}
\begin{proof}
By Corollary~\ref{Riccati-Mult} we have 
$\Ric(A_{\ltimes h})=\Ric(A)_{+g}\in K\{Z\}$. Next one shows by an easy
induction on $n$ and using $R_n=Z^n+\text{lower degree terms}$: 
if $L$ is a differential field extension of $K$, and
$R=a_0R_0+\cdots+a_nR_n$ ($a_0,\dots,a_n\in L$) is a Riccati polynomial 
over $L$ with $R\in K\{Z\}$, then $a_0,\dots,a_n\in K$.
\end{proof}

\index{differential field!weakly $r$-differentially closed}
\index{differential field!weakly differentially closed}
\index{weakly!differentially closed}
\index{differentially!closed!weakly}
\label{p:weakly differentially closed}

\noindent
The decrease in order under Riccati transformation is the basis of
several inductive proofs. The next lemma is an example. 
For $r\in \N$ we say that 
$K$ is {\bf weakly $r$-differentially closed} if every 
$P\in K\{Y\}\setminus K$ of order $\le r$ has a zero
in $K$.   (Thus~$K$ is weakly $0$-differentially closed iff $K$ is algebraically closed.)
We also say that~$K$ is {\bf weakly differentially closed} if $K$ is 
weakly $r$-differentially closed for each $r\in\N$. Clearly differentially closed $\Rightarrow$ weakly differentially closed.

\begin{lemma}\label{fa0} Suppose $K$ is weakly $r$-differentially closed, $r\in \N$. Then $K$ is $(r+1)$-linearly closed.
\end{lemma}
\begin{proof} By induction on $r$. The case $r=0$ is obvious, so let $r>0$, and let
$A\in K[\der]$ have order $r+1$. Then its Riccati transform $\Ric(A)\in K\{Z\}$ has order $r$, so has a zero $a\in K$, which gives $A=B\cdot (\der - a)$ with $B\in K[\der]$ of
order $r$. Now apply the inductive assumption to $B$. \end{proof}

\subsection*{Riccati transforms in the presence of a valuation}
{\em In this subsection $K$ is a valued differential field such that
$\der \smallo\subseteq  \smallo$, and thus $\der\mathcal{O} \subseteq \mathcal{O}$.}

\begin{lemma}\label{valric} Let $P\in K\{Y\}$ be homogeneous. Then 
$v(P)=v(\Ric(P))$.
\end{lemma}
\begin{proof} By $K$-linearity of the Riccati transform we can reduce 
to the case~$vP=0$. Then $\Ric(P)\in \mathcal{O}\{Z\}$. Let $\bar{P}$ be the
image of $P\in \mathcal{O}\{Y\}$ in $\k\{Y\}$. Then $\bar{P}$ is 
homogeneous and 
nonzero, so $\Ric(\bar{P})\ne 0$. Since $\Ric (\bar{P})$ is the image of
$\Ric(P)$ under the natural map $\mathcal{O}\{Z\}\to  \k\{Z\}$, we
get $v\big(\Ric(P)\big)=0$.       
\end{proof}

\noindent
In view of Corollary~\ref{Riccati-Mult}, this lemma yields:

\begin{cor}\label{corvalric} If $P\in K\{Y\}^{\ne}$ is homogeneous of degree $d$, and 
$a\in K^\times$, then
$$ v_P(\alpha)\ =\ d\alpha + v(R_{+a^\dagger}), \qquad \text{ with $\alpha:= va$, $R:= \Ric(P)$.}$$
\end{cor}

\noindent
For a nonzero Riccati polynomial 
$R=a_0R_0 + \dots + a_nR_n$ ($a_0,\dots,a_n\in K$), put
$$\Ricmu(R)\ :=\ \min\bigl\{i:v(a_i)=v(R)\bigr\}, \qquad 
\Ricnu(R)\ :=\ \max\bigl\{i:v(a_i)=v(R)\bigr\}.$$ 
Note that for nonzero $A\in K[\der]$ and $R=\Ric(A)$, we have
$$%\begin{equation}\label{mu-1}
v(A)\ =\ v(R), \quad \Pmu(A)\ =\ \Ricmu(R), \quad \Pnu(A)\ =\ \Ricnu(R),
$$%\end{equation}
and by Corollary~\ref{Riccati-Mult}:
$$%\begin{equation}\label{mu-2}
v(Ay)=vy+v(R_{+z}), \quad
\Pmu(Ay)=\Ricmu(R_{+z}), \quad \Pnu(Ay)=\Ricnu(R_{+z}) .
$$%\end{equation}
The quotients of $K$ by its $\Q$-linear subspaces $\smallo$ and $\mathcal{O}$ are
$K/\smallo$ and $K/\mathcal{O}$, with canonical maps 
$$g \mapsto g+\smallo
\colon K\to K/\smallo, \qquad g \mapsto g+\mathcal{O}
\colon K\to K/\mathcal{O}.$$ If $R$ is a nonzero Riccati polynomial over $K$,
then by Lemma~\ref{v-under-conjugation}(i), the
value $v(R_{+g})$ for $g\in K$ depends only on the coset $g+\mathcal{O}$.

\nomenclature[V]{$\mu(R)$}{dominant multiplicity of $R$}
\nomenclature[V]{$\nu(R)$}{dominant weight of $R$}

\begin{lemma} \label{mu and nu}
Let $R$ be a nonzero
Riccati polynomial over $K$. Then $\Ricmu(R_{+g})$ and $\Ricnu(R_{+g})$ depend
only on $g+\smallo$ and  $g+\mathcal{O}$, respectively,
for $g\in K$.
\end{lemma}
\begin{proof} Let $R= a_0R_0 + \dots + a_nR_n$ with $a_0,\dots,a_n\in K$, and
let $g\in K^{\prec 1}$. Then $R_j(g)\in \smallo$ for $j>0$,
since $\smallo$ is closed under $\der$.
Now
$$R_{+g}(Z)\ =\  b_0R_0(Z) + b_1R_1(Z) + \cdots + b_nR_n(Z)$$
where
$$b_i\ =\ \sum_{j=i}^n \binom{j}{i} a_j R_{j-i}(g)\ 
=\   a_i + (i+1)a_{i+1}R_1(g) +\cdots + \binom{n}{i}a_nR_{n-i}(g).
$$
Hence $v(b_i)\geq v(R)$ for all $i$, and
$v(b_i)=v(a_i)=v(R)$ for $i=\Ricmu(R)$; 
so $\Ricmu(R_{+g})\leq\Ricmu(R)$.
Using $R=(R_{+g})_{-g}$ we get $\Ricmu(R_{+g})=\Ricmu(R)$.
The general case follows easily from this special case $vg>0$. The argument 
for $\Ricnu(R_{+g})$ is similar.
\end{proof}

\noindent
For a nonzero Riccati polynomial $R$ over $K$ we define 
\begin{align*} &\Ricmu_R\colon K/\smallo \to \N, \quad \quad
 \Ricmu_R(g+\smallo) := \Ricmu(R_{+g}) \quad \quad (g\in K),\\
&\Ricnu_R\colon K/\mathcal{O} \to \N, \quad \quad
\Ricnu_R(g+\mathcal{O}) := \Ricnu(R_{+g}) \quad \quad (g\in K).
\end{align*}

\nomenclature[V]{$\mu_R(g+\smallo)$}{$\mu(R_{+g})$}
\nomenclature[V]{$\nu_R(g+\mathcal{O})$}{$\nu(R_{+g})$}

\subsection*{Notes and comments}
Riccati polynomials are named after J.~Riccati (1676--1754), who
studied the differential equation $Z'+aZ^2 = b x^r$ for $a,b,r\in \R$ in~\cite{Riccati}.

The results in Section~\ref{sec:newtonization} and Chapter~\ref{ch:newtdirun} below yield that 
$\T[\imag]$ and $\T_{\log}[\imag]$
($\imag^2=-1$) are weakly differentially closed.

%% file: mt-5-9.tex
\section{Johnson's Theorem}\label{sec:johnson}

\noindent
Below $Y=(Y_1,\dots, Y_n)$ is a tuple of distinct $\d$-indeterminates, $n\ge 1$. The goal of this section is to prove 
the following result due to 
Joseph Johnson~\cite{JJ}.  

\begin{theorem}\label{thm:johnson} 
Let $K\subseteq L$ be a differential field extension, let 
$P_1,\dots, P_n\in K\{Y\}$, and let $y=(y_1,\dots, y_n)\in L^n$ be such that
$P_1(y)=\cdots = P_n(y)=0$ and $P_1,\dots,P_n$ are 
$\d$-independent at $y$. Then $y_1,\dots, y_n$ are
$\d$-algebraic over $K$.
\end{theorem} 

\noindent
Theorem~\ref{thm:johnson} is an analogue of
Corollary~\ref{cor:regsepalg}, with differential field extensions replacing
field extensions. Accordingly, we begin by establishing systematically analogues of results in Sections~\ref{sec:modules},~\ref{sec:differentials}, and~\ref{sec:derivations on field exts}
for differential $K$-algebras $A$ instead of $K$-algebras,
with $\Omega_{A|K}$ construed as an $A[\der]$-module instead of
just an $A$-module, and differential transcendence degree
instead of transcendence degree. 

Note also that Theorem~\ref{thm:johnson} partially extends
Corollary~\ref{indlininhom} to non-linear systems of algebraic differential equations.

Johnson considers more generally fields of arbitrary characteristic with finitely many commuting derivations, as is usual in the Kolchin tradition. In our setting where differential fields are of characteristic $0$ with a single
derivation, things are a bit simpler, but nevertheless we need some nontrivial facts on regular local rings from Chapter~\ref{ch:commalg}.
We give the proof of Theorem~\ref{thm:johnson} at the end of this section; we precede it with some generalities on the tensor product of  $K[\der]$-modules and the module of differentials of a given differential field extension.

We shall apply Theorem~\ref{thm:johnson} in 
Section~\ref{sec:max} to differential-henselian fields. 

\subsection*{Tensor products of $K[\der]$-modules} Let $K$ be
a differential ring, let $M$, $N$ be $K[\der]$-modules, and consider the $K$-module $T:=M\otimes_K N$. Corollary~\ref{cor:der on tensor prod, 1} gives a $\der$-compatible
derivation on $T$ making $T$ a $K[\der]$-module such that
$$\der(x\otimes y)\ =\ (\der x)\otimes y+x\otimes (\der y) \qquad\text{for all
$x\in M$, $y\in N$.}$$
If $K$ is a differential field and $M$ and $N$ are differential modules over $K$, then $T$ is a differential module over $K$ with
$\dim_K T = \dim_K M \cdot \dim_K N$,
and if in addition~$M$ and~$N$ are horizontal, then so is $T$.
%The $K$-linear map 
%$$M^*\otimes_K M\to K\quad\text{ with $\phi\otimes f\mapsto \<\phi,f\>=\phi(f)$ for  $\phi\in M^*$, $f\in M$}$$ 
%considered in Section~\ref{sec:systems} is even $K[\der]$-linear.

\medskip
\noindent
Let $A$ be a {\bf differential $K$-algebra,} that is, $A$ is
a differential ring equipped with a differential ring morphism
$K \to A$; the latter also makes $A$ a $K$-algebra.  
We view $A$ as a $K[\der]$-module with $\der a=a'$ for 
$a\in A$. By Corollary~\ref{cor:der on tensor prod, 2} we obtain the $A[\der]$-module $A\otimes_K M$, called the {\bf base change} of $M$ to $A$. 
Note that   
$$\der(a\otimes x)\ =\ a'\otimes x + a\otimes (\der x) \qquad\text{for $a\in A$, $x\in M$.}$$
Let $(x_i)$ be a family in $M$. Then the $K[\der]$-module $M$ is free on $(x_i)$ iff the $K$-module~$M$ is free on
$(\der^n x_i)_{i,n}$. If the $K[\der]$-module $M$
is free on $(x_i)$, then the $A[\der]$-module $A\otimes_K M$
is free on $(1\otimes x_i)$: use that $\der^n(1\otimes x)=1\otimes \der^n x$ for $x\in M$. The differential ring morphism $K \to A$ extends to a ring morphism $K[\der] \to A[\der]$ sending $\der$ to $\der$, and any $A[\der]$-module is construed as a $K[\der]$-module via this ring morphism~$K[\der]\to A[\der]$.  

\index{algebra!differential}
\index{differential algebra}
\index{differential module!base change}
\index{base change}

Let $I$ be a differential ideal of $K$, and $A=K/I$, viewed as a differential $K$-algebra via the canonical map $K \to A$.  Then $IM$ is a $K[\der]$-submodule of $M$, which makes~$M/IM$ a $K[\der]$-module as well as an $A$-module. We give 
$M/IM$ its unique structure of $A[\der]$-module that induces its given $A$-module structure as well as its $K[\der]$-module structure coming from the natural ring morphism $K[\der]\to A[\der]$. 
Then the isomorphism 
$A\otimes_K M\xrightarrow{\cong} M/IM$ of \eqref{eq:iso tensor, 3} is $A[\der]$-linear. 

The next lemma is for use in the second volume. 

\begin{lemma}\label{lem:iso base change}
Let $K$ be a differential field, $M$ a differential module over~$K$,
and~$F$ a differential field extension of~$K$. Then the base change
$F\otimes_K M$ of $M$ to $F$ is a differential module over $F$ with $\dim_{F} F\otimes_K M = \dim_K M$.
Moreover, if $M=K[\der]/K[\der]L$ where  $L\in K[\der]^{\neq}$, then 
$F\otimes_K M\cong F[\der]/F[\der]L$
as $F[\der]$-modules. 
\end{lemma}
\begin{proof}
The first part holds by the remarks preceding the lemma. Suppose $M=K[\der]/K[\der]L$ where $L\in K[\der]^{\neq}$ has
order $r\ge 1$. Then $e=1+K[\der]L$ is a cyclic vector of~$M$, since  $M=Ke\oplus K\der e\oplus\cdots\oplus K\der^{r-1} e$ (internal direct sum of $K$-linear subspaces) and $Le=0$. Setting $e_F:=1\otimes e$ we  have $ \der^n e_F= 1\otimes \der^n e$ for each $n$ and so $F\otimes_K M = Fe_F\oplus F\der e_F\oplus\cdots\oplus F\der^{r-1} e_F$  with $Le_F=1\otimes Le=0$.
\end{proof}

\noindent
If $K$ is a differential field, and $A$, $B$ are differential $K$-algebras, then the derivation on the ring $A\otimes_K B$
provided by Corollary~\ref{cor:der on tensor prod, 3} is the unique derivation on this ring for which the (injective)
maps $a\mapsto a\otimes 1\colon A\to A\otimes_K B$ and $b\mapsto 1\otimes b\colon B\to A\otimes_K B$ are
morphisms of differential rings.

\subsection*{The module of differentials as a $K[\der]$-module}
Let $K$ be a differential ring and~$A$ a differential $K$-algebra.
The map $\der^*$ of Corollary~\ref{cor:derivation on Omega} makes the $A$-module~$\Omega_{A|K}$ 
of K\"ahler differentials
an $A[\der]$-module (and thus a $K[\der]$-module), and
the universal $K$-derivation $\d_{A|K}\colon A\to \Omega_{A|K}$ is then $K[\der]$-linear. Below we view~$\Omega_{A|K}$ as an $A[\der]$-module in this way; so 
$\der\omega=\der^*(\omega)$ for $\omega\in \Omega_{A|K}$.  

Let $B$ be a second differential $K$-algebra and let $h\colon A\to B$ be a morphism of differential $K$-algebras. Then $B$ is via $h$ also a differential $A$-algebra, and so we have the $B[\der]$-module $\Omega_{B|A}$, and the base change 
$B\otimes_A \Omega_{A|K}$ of the $A[\der]$-module~$\Omega_{A|K}$ to~$B$ is a $B[\der]$-module. By Proposition~\ref{prop:first exact sequ} we have an exact sequence
$$B\otimes_A \Omega_{A|K} \xrightarrow{\ \alpha\ } \Omega_{B|K} \xrightarrow{\ \beta\ } \Omega_{B|A}\longrightarrow 0$$
of $B$-modules and $B$-linear maps, where $\alpha(1\otimes \d a)=\d h(a)$ for all $a\in A$ and $\beta(\d_{B|K}b)=\d_{B|A}b$ for all $b\in B$. One verifies easily that $\alpha$ and $\beta$ are $B[\der]$-linear.

Next, let $I\ne A$ be a differential ideal of $A$, let 
$B:= A/I$
be the quotient differential $K$-algebra with derivation
$\der_B$, with the canonical morphism
$h\colon A \to B$ of differential
$K$-algebras. Then the $\der_B$-compatible derivation
$a+ I^2\mapsto a' + I^2$ (for~$a\in I$) on the $B$-module $I/I^2$ makes the latter a $B[\der]$-module. Proposition~\ref{prop:second exact sequ} yields the exact sequence of $B$-modules and $B$-linear maps
$$I/I^2 \xrightarrow{\ \gamma\ } B\otimes_A \Omega_{A|K} \xrightarrow{\ \alpha\ } \Omega_{B|K}\longrightarrow 0,$$
with $\gamma(a+I^2)=1\otimes\d a$ for $a\in I$. It is easy to check that $\gamma$ is even $B[\der]$-linear.

Finally, let $R$ be a differential $K$-algebra, let 
$\mathfrak p$ be a 
differential prime ideal of~$R$, and set $A:=R_{\mathfrak p}$.
Then $A$ is a differential $K$-algebra and a local ring
whose maximal ideal $\m=\mathfrak p A$ is a differential ideal of $A$, with differential residue field $F:=A/\m$. It is easy to check that the maps $\alpha_0$ and
$\gamma_0$ in the exact sequence 
$$\m/\m^2 \xrightarrow{\ \gamma_0\ } F\otimes_{R} \Omega_{R|K} \xrightarrow{\ \alpha_0\ } \Omega_{F|K}\longrightarrow 0$$
of $F$-linear spaces and $F$-linear maps from \eqref{eq:second exact sequ} are $F[\der]$-linear.

\subsection*{Modules of differentials for differential field extensions}
In the rest of this section $K$ is a differential field and $L$ is a differential field extension of $K$. For any family $(a_i)_{i\in I}$ in $L$, the following conditions are 
equivalent: 
\begin{enumerate}
\item the family $(a_i)$ in $L$ is $\d$-algebraically independent over $K$;
\item the family $\big(a_i^{(n)}\big)_{i,n}$ in $L$ is algebraically independent over $K$;
\item the family $\big(\!\d a_i^{(n)}\big)_{i,n}$ in the vector space $\Omega_{L|K}$ is linearly independent over $L$;
\item the family $(\d a_i)$ in the $L[\der]$-module $\Omega_{L|K}$ is $L[\der]$-independent. 
\end{enumerate}
For (1)~$\Leftrightarrow$~(2) we refer to the subsection on differential transcendence bases in Section~\ref{Differential Fields and Differential Polynomials}, (2)~$\Leftrightarrow$~(3)
follows from Lemma~\ref{lem:lindep in Omega}, and (3)~$\Leftrightarrow$~(4) holds because
$\d a^{(n)}=\der^n \d a$ for $a\in L$. 
In particular, (1)~$\Leftrightarrow$~(4) yields for $a\in L$:
\begin{align*} a \text{ is $\d$-algebraic over $K$}\quad 
&\Longleftrightarrow\quad 
\d a\in (\Omega_{L|K})_{\operatorname{tor}}.\\
\hskip-1.5em\text{Thus:}\ \qquad \qquad
\text{$L$ is $\d$-algebraic over $K$} \quad  
&\Longleftrightarrow\quad \text{$\Omega_{L|K}$ is a torsion $L[\der]$-module.}%\qquad \qquad
\end{align*}
Suppose now that $L=K\<x_1,\dots,x_n\>$. Then 
$L=K\big(x_j^{(m)}: m\in \N,\ j=1,\dots,n\big)$, so the vector space $\Omega_{L|K}$ over $L$ is generated by the $\d\!\big(x_j^{(m)}\big)=\der^m \d x_j$ with $m\in \N$ and $j=1,\dots, n$. Hence as an $L[\der]$-module, $\Omega_{L|K}$ is generated by $\d x_1,\dots,\d x_n$. Moreover:

{\sloppy
\begin{cor}\label{rank=trdegder}
$\operatorname{rank} \Omega_{L|K} = \operatorname{trdeg}_\der L|K$.
\end{cor}
\begin{proof} By permuting the indices we arrange
that $m\le n$ is such that $x_1,\dots, x_m$ is a $\d$-transcendence base of $L$ over $K$; in particular, $\operatorname{trdeg}_\der L|K = m$.
Also, the elements $\d x_1,\dots, \d x_m$ of $\Omega_{L|K}$ are $L[\der]$-independent by (1)~$\Leftrightarrow$~(4), and likewise,
if~${m < j \le n}$, then $x_1,\dots, x_m, x_j$ are $\d$-algebraically dependent, so $\d x_1,\dots, \d x_m, \d x_j$ are
$L[\der]$-dependent. Thus $\operatorname{rank} \Omega_{L|K}=m$
by Corollary~\ref{cor:additivity of rank}.
\end{proof}
}

\subsection*{Independence at a prime}
Let $R$ be the differential $K$-algebra $K\{Y_1,\dots, Y_n\}$.
By Lemma~\ref{lem:diff module of polynomial algebra} the $R$-mod\-ule~$\Omega_{R|K}$ is free on its generating family $\big(\!\d(Y_j^{(r)})\big)_{j,r}$. As an $R[\der]$-module with
$\der \d P=\d(P')$ for~$P\in R$, it is free on $\d Y_1,\dots,\d Y_n$. 

\index{independent!at a prime}
\index{d-independent@$\d$-independent}

Let $P_1,\dots, P_m\in R$ and
$\frak{p}\in \Spec(R)$. Then we say that $P_1,\dots, P_m$ are
{\bf $\d$-independent at $\frak{p}$} if the family
$\big(\!\d(P_i^{(r)})\big)_{i,r}$ in the $R$-module $\Omega_{R|K}$ is independent at $\frak{p}$ as defined in Section~\ref{sec:modules}. 
If $P_1,\dots, P_m$
are $\d$-independent at $\frak{p}$ and $\frak{p}\supseteq \frak{q}\in \Spec(R)$, then $P_1,\dots, P_m$ are 
$\d$-independent at~$\frak{q}$, by Lemma~\ref{lem:independence under generalization}. 
In view of the isomorphism~\eqref{eq:iso tensor, 3}, the following are equivalent:
\begin{enumerate}
\item $P_1,\dots, P_m$ are $\d$-independent at $\frak{p}$; 
\item the family $\big(1\otimes \d(P_i^{(r)})\big)_{i,r}$ in $(R/\mathfrak p)\otimes_R \Omega_{R|K}$ is $(R/\mathfrak p)$-linearly independent.
\end{enumerate}
The $R$-module $\Omega_{R|K}$ is free on $\big(\!\d(Y_j^{(r)})\big)_{j,r}$, so the $(R/\frak{p})$-module $(R/\mathfrak p)\otimes_R \Omega_{R|K}$ is free on 
$\big(1\otimes\d(Y_j^{(r)})\big)_{j,r}$, and likewise
the vector space 
$F\otimes_R \Omega_{R|K}$ over  $F:=\Frac(R/\frak{p})$ is free on $\big(1\otimes\d(Y_j^{(r)})\big)_{j,r}$. Thus (1) and (2)
are also equivalent to \begin{enumerate}
\item[(3)] the family $\big(1\otimes \d(P_i^{(r)})\big)_{i,r}$ in $F\otimes_R \Omega_{R|K}$ is $F$-linearly independent.
\end{enumerate}
If $\frak{p}$ is also a differential ideal of $R$,
then (2) is equivalent to (4), and (3) to (5): 
\begin{enumerate}
\item[(4)] $1\otimes \d P_1,\dots, 1\otimes \d P_m\in (R/\mathfrak p)\otimes_R \Omega_{R|K}$ are $(R/\mathfrak p)[\der]$-independent;
\item[(5)] $1\otimes \d P_1,\dots, 1\otimes \d P_m\in F\otimes_R \Omega_{R|K}$ are $F[\der]$-independent.
\end{enumerate}
Suppose $\frak{p}\in \Spec(R)$ is a differential ideal of $R$. The inclusion $K \to R$ composed with the canonical map 
$R\to R/\frak{p}$ yields an injective differential ring
morphism $K \to R/\frak{p}$ via which we identify $K$ below with a differential subring of $R/\frak{p}$. 

\begin{lemma}\label{dindeplemma} Let $y\in L^n$ and let $\mathfrak p$ be the kernel of the differential $K$-algebra morphism $R\to L$ that
sends $Y_j$ to $y_j$ for $j=1,\dots,n$. Then $\mathfrak p$ is a differential prime ideal of $R$, and this morphism $R\to L$ induces a differential ring isomorphism 
$R/\mathfrak{p}\to K\{y\}$, which extends to a differential
field isomorphism $F\to K\<y\>$, where 
$F:= \Frac(R/\mathfrak{p})$. Moreover, we have the equivalence
$$\text{$P_1,\dots,P_m$ are $\d$-independent at $y$}\quad\Longleftrightarrow\quad\text{$P_1,\dots, P_m$
 are $\d$-independent at $\frak{p}$.}$$
 \end{lemma}
\begin{proof} Obvious, except for the final equivalence. For that we use the relevant definitions, 
the equality~\eqref{eq:universal der, poly}, the isomorphism
of free $F[\der]$-modules
$$(A_1,\dots,A_n)\mapsto A_1\cdot(1\otimes\d Y_1)+\cdots+A_n\cdot(1\otimes\d Y_n)\ \colon\ F[\der]^n\to F\otimes_R \Omega_{R|K},$$
and the above equivalence (1)~$\Leftrightarrow$~(5).  
\end{proof}

\subsection*{An abstract version of Johnson's theorem}
We keep the notations introduced in the previous subsection.
Theorem~\ref{thm:johnson} will be derived from the next result.

\begin{prop}\label{prop:johnson}
Let $P_1,\dots, P_n\in R$, $I:=[P_1,\dots, P_n]$, and let $\frak{p}\in \Spec(R)$ be such that $I\subseteq \frak{p}$ and $P_1,\dots, P_n$ are $\d$-independent at $\frak{p}$. Set
$$\frak{q}\ :=\ \{a\in R:\ \text{$ab\in I$ for some $b\in R\setminus \frak{p}$}\}.$$
Then $\frak{q}$ is a differential ideal of $R$ with the following properties: \begin{enumerate}
\item[\textup{(i)}] $\frak{q}$ is the smallest prime
ideal of $R$ containing $I$ and contained in $\frak{p}$;
\item[\textup{(ii)}] the differential fraction field of $R/\frak{q}$ is $\d$-algebraic over $K$.
\end{enumerate}
\end{prop}
\begin{proof} We consider $R$ and its localization 
$A:= R_{\frak{p}}$ as differential subrings of the differential field $\Frac(R)= K\<Y\>$. 
Clearly $\mathfrak q$ is a differential ideal of $R$ with $I\subseteq\mathfrak q\subseteq\mathfrak p$,
and any prime ideal of~$R$ containing $I$ and contained in $\mathfrak p$ contains $\mathfrak q$;
so to finish the proof of (i) it remains to show that $\mathfrak q$ is prime.
Note that $A$ is a local domain with maximal ideal $\frak{m}:=\frak{p}A$, and $(IA)\cap R=\frak{q}$. We shall prove that $IA$ is a prime ideal of $A$ (so
$\frak{q}$ is prime as a consequence). To
reduce to a more finitary situation, set for $r\in \N$:
$$R_r:=K\big[Y_j^{(k)}:\ j=1,\dots, n,\ k=0,\dots,r\big],\quad \frak{p}_r:= \frak{p}\cap R_r, \quad A_r:= (R_r)_{\frak{p}_r},$$
so $A_r$ is a local subring of $A$ with maximal ideal 
$\frak{m}_r:=\frak{p}_rA_r=\frak{m}\cap A_r$.
For $IA$ to be a prime ideal of $A$, it is clearly enough to show:

\claim{Let $r, N\in \N$ and let distinct
$$G_1,\dots, G_N\in R_r\cap \big\{P_i^{(k)}:\ i=1,\dots,n,\ k\in \N\big\}$$
be given. Then $G_1,\dots, G_N$ generate a prime ideal of $A_r$.} 

\noindent
To prove this claim, set $F := A/\frak{m}$ and $F_r := A_r/\frak{m}_r$,
so $F_r$ is a subfield of the field~$F$ after the usual identification.
We have an exact sequence
\begin{equation}\label{eq:johnson}
\m/\m^2 \xrightarrow{\ \gamma_0\ } F\otimes_R \Omega_{R|K} \longrightarrow \Omega_{F|K}\longrightarrow 0
\end{equation}
of $F$-linear spaces, where
$\gamma_0(a+\m^2) = 1\otimes\d a$ for $a\in \frak{p}$. Composing $\gamma_0$ with
the natural $F_r$-linear map $\m_r/\m_r^2\to\m/\m^2$ yields an $F_r$-linear map
$$ \frak{m}_r/\frak{m}_r^2 \to F\otimes_R\Omega_{R|K}, \qquad 
a+\frak{m}_r^2 \mapsto 1\otimes \d a 
\text{ for }a\in\frak{p}_r.$$ 
Since $P_1,\dots,P_n$ are $\d$-independent at $\frak{p}$, the natural isomorphism
$\Frac(R/\mathfrak p)\cong F$ and the equivalence (1)~$\Leftrightarrow$~(3) preceding Lemma~\ref{dindeplemma} 
yield that in $F\otimes_R \Omega_{R|K}$ the elements
$1\otimes \d G_1,\dots, 1\otimes \d G_N$ 
are linearly independent over~$F$. Therefore the elements
$G_1+ \frak{m}_r^2,\dots, G_N+\frak{m}_r^2$ of $\frak{m}_r/\frak{m}_r^2$ are linearly independent over $F_r$. Now~$A_r$ is a regular local ring by Corollary~\ref{cor:loc of poly rings are regular local},
so $G_1,\dots, G_N$ generate a prime ideal in~$A_r$ by Corollaries~\ref{cor:NAK, 3}(ii) and \ref{cor:reg => domain}.
This concludes the proof of~(i). 

\medskip
\noindent
Towards proving (ii), first note that by Lemma~\ref{lem:expansion-contraction} we have
$\mathfrak q R_{\mathfrak q} = I R_{\mathfrak q}$ and
$$\mathfrak q\ =\ \{a\in R:\ \text{$ab\in I$ for some $b\in R\setminus \frak{q}$}\}.$$
To show (ii) we replace $\mathfrak p$ by $\mathfrak q$ to arrange that
$\mathfrak p$ is a differential ideal of $R$ and $\mathfrak{m}=\mathfrak p A = I A$ is the differential ideal of $A$ generated by $P_1,\dots, P_n$.
Then \eqref{eq:johnson} is an exact sequence of $F[\der]$-modules and $F[\der]$-linear maps. As $P_1,\dots,P_n$ are $\d$-independent at~$\frak{p}$, the elements
$1\otimes\d P_1,\dots, 1\otimes\d P_n$ of $F\otimes_R\Omega_{R|K}$ are $F[\der]$-independent. Since the $F[\der]$-module $\m/\m^2$ is generated by $P_1+\m^2,\dots,P_n+\m^2$, the $F[\der]$-linear map~$\gamma_0$ yields that the $F[\der]$-module
$\m/\m^2$ is free on $P_1+ \m^2,\dots, P_n+ \m^2$, and 
$\gamma_0$ is injective. The $F[\der]$-module $F\otimes_R \Omega_{R|K}$ being free on 
$1\otimes \d Y_1,\dots, 1\otimes\d Y_n$, taking ranks
of $F[\der]$-modules, we get 
$\operatorname{rank}(\Omega_{F|K}) = 0$
from \eqref{eq:johnson} and Lemma~\ref{lem:additivity of rank}.
As the canonical field isomorphism $\Frac(R/\mathfrak p)\cong F$ respects the natural derivations on these fields,
this ``rank $0$'' fact gives (ii) in view of Corollary~\ref{rank=trdegder}.
\end{proof}

\begin{proof}[Proof of Theorem~\ref{thm:johnson}] We are given
$P_1,\dots, P_n\in K\{Y\}$ and $y\in L^n$ such that 
$P_1(y)=\cdots=P_n(y)=0$ and $P_1,\dots, P_n$ are $\d$-independent at $y$; we have to show that then $K\<y\>$ is
$\d$-algebraic over $K$. Let $\mathfrak p$ be the kernel of the differential $K$-algebra morphism $R\to L$ that
sends $Y_j$ to $y_j$ for $j=1,\dots,n$. Then 
$P_1,\dots, P_n$ are $\d$-independent at $\frak{p}$ by 
Lemma~\ref{dindeplemma}. Let $\mathfrak{q}$ be as in
Proposition~\ref{prop:johnson}. Then~(ii) of that proposition
yields an element $Q_j\in \mathfrak{q} \cap K\{Y_j\}^{\ne}$
for $j=1,\dots,n$. Then $Q_j\in \mathfrak{p}$ for 
$j=1,\dots,n$, so $K\<y\>$ is $\d$-algebraic over $K$.
\end{proof}

\noindent
We finish with an application of Theorem~\ref{thm:johnson} and Lemma~\ref{lem:rubel}:

\begin{cor}\label{cor:rubel} Suppose $K\subseteq L$ is a differential field extension, 
$P_1,\dots, P_n$ in~$K\{Y\}$ have order at most $(r_1,\dots,r_n)\in\N^n$, and $y=(y_1,\dots, y_n)\in L^n$ satisfies
$$P_1(y)\ =\ \cdots\ =\ P_n(y)\ =\ 0, \qquad 
\det\left(\big(\partial P_i\big/\partial Y_j^{(r_j)}\big)
(y)\right)\ \neq\ 0.$$ 
Then $y_1,\dots, y_n$ are $\d$-algebraic over $K$.
\end{cor}
\begin{proof} By the hypothesis $P_1,\dots, P_n$ are $\d$-independent at $y$.
\end{proof}

\begin{example}
Let  $Q_1,\dots,Q_n\in K\{Y\}$ have order at most $(r_1,\dots,r_n)\in\N^n$. If $L$ is a differential field extension of $K$ and $y=(y_1,\dots,y_n)\in L^n$ 
satisfies
$$y_i^{(r_i+1)}\ =\ Q_i(y)\qquad (i=1,\dots,n),$$ 
then $y_1,\dots,y_n$ are $\d$-algebraic over $K$. (Take $P_i:=Y_i^{(r_i+1)}-Q_i$.)
\end{example}

\subsection*{Notes and comments}
Proposition~\ref{prop:johnson} is from \cite{JJ}. 
See also \cite{JJ69-1} for
more information on K\"ahler differentials in the context of (partial) differential fields. An analytic version of
Corollary~\ref{cor:rubel} is due to Ru\-bel~\cite{Rubel82}, whose proof is elementary and independent of Johnson's theorem.

%% file: mt-6.tex
\chapter{Valued Differential Fields}\label{ch:valueddifferential}

\setcounter{theorem}{0}

\noindent
{\em Throughout this chapter, $K$ is a valued differential field whose derivation is
small in the sense that $\der\smallo\subseteq \smallo$}. In Chapter~\ref{ch:differential polynomials} we 
already derived some
consequences of this assumption, and here we go much further with it, 
sometimes under mild additional conditions. As noted before, 
one thing we get from $\der\smallo\subseteq \smallo$ is that 
$\der\mathcal{O}\subseteq \mathcal{O}$, and so
the residue field $\k=\mathcal{O}/\smallo$ is naturally a differential field.
Sometimes we impose as extra condition that the derivation of $\k$ is 
nontrivial. While this extra condition is not satisfied for $K=\T$,
it does hold in suitable coarsenings of compositional conjugates
of $\T$, and this is how the results of this chapter apply
to $\T$, and other valued differential fields of interest.
 
\medskip\noindent
Section~\ref{sec:asbevp} considers the asymptotic behavior of 
the function $v_P\colon \Gamma \to \Gamma$ for homogeneous $P\in K\{Y\}^{\ne}$. In 
Section~\ref{sec:algebraicext} we show among other things
that the derivation of any valued differential 
field extension of $K$ that is algebraic over $K$ is also small. In 
Section~\ref{sec:resext} we show how differential field extensions 
of the residue differential field $\k$ give rise to 
valued differential field extensions of $K$ 
with small derivation and the same value group. We also study there monotone extensions.  

\medskip\noindent
In Section~\ref{sec:dvvg} we show how the derivation and the valuation on
$K$ jointly give rise to a valuation on $\Gamma$. This helps 
in improving the estimates on $v_P$ from Section~\ref{sec:asbevp}. 
The ordered abelian group $\Gamma$ with this valuation is an 
{\em asymptotic couple\/} in the sense of Rosenlicht, and for frequent 
later use we prove in Section~\ref{sec:ascouples} some basic facts about 
asymptotic couples in general. 
In Section~\ref{sec: dominant} we define the {\em dominant part\/} of a differential
polynomial, and establish its basic properties. The key facts about
the functions $v_P$, the valuation on $\Gamma$, and dominant parts are then
used to prove the important Equalizer Theorem in Section~\ref{equalize}:

\index{theorem!Equalizer Theorem}

\begin{theorem}\label{theq} Let $P,Q\in K\{Y\}^{\ne}$ be
homogeneous of degree $d$ and $e$, respectively, with~${d>e}$, and suppose  
$(d-e)\Gamma=\Gamma$. Then there is a unique
$\alpha\in \Gamma$ such that $v_P(\alpha)=v_Q(\alpha)$.
\end{theorem}

\noindent 
We shall use this mainly for $d=1$, $e=0$, or when $\Gamma$ is divisible, and
in either case the condition $(d-e)\Gamma=\Gamma$ is trivially satisfied.
(For $d=1$, $e=0$ the theorem says that $v_A\colon \Gamma \to \Gamma$ is bijective
for each $A\in K[\der]^{\ne}$.)

\medskip\noindent
In Section~\ref{sec:pceq} we consider a pc-sequence
$(a_{\rho})$ in $K$ and a differential polynomial $G\in K\{Y\}\setminus K$.
If the derivation of $\k$ is 
nontrivial, then we show how to replace~$(a_{\rho})$ by an equivalent pc-sequence to arrange that  
$G(a_\rho)$ is also a pc-sequence. 
We also prove some variants of this 
important fact, and use this in Section~\ref{sec:cimex} to construct immediate extensions of
such $K$
in analogy with Kaplansky's lemmas~\ref{lem:Kaplansky, 1} and~\ref{lem:Kaplansky, 2}. This has the effect that
any such $K$ has a spherically complete immediate valued differential field extension with small
derivation. 

\subsection*{Notes and comments} The key result in this chapter
that is needed later in connection with the model theory of $\T$ is the Equalizer Theorem of Section~\ref{equalize}. 
%The construction of spherically complete
%immediate extensions in Section~\ref{sec:cimex}  
Other sections are clearly part of any reasonably broad theory of valued differential fields, and play a role in the next chapter on differential-henselian fields. Sections~\ref{sec:pceq} and~\ref{sec:cimex} are inspired by analogous results on certain kinds of valued difference fields due to B\'elair, Macintyre, and Scanlon \cite{BMS}.

\section{Asymptotic Behavior of \protect\hidevp}\label{sec:asbevp}

\noindent
Here is a key consequence of our assumption $\der\smallo\subseteq \smallo$:

\begin{lemma}\label{cohnest} If $y\in \smallo$, then $(y')^{n+1}\preceq y^n$ for all $n$.
\end{lemma}
\begin{proof} Suppose there is a counterexample. Take $n$ least with 
$(y')^{n+1} \succ y^n$ for some $y\in \smallo$, and fix such
$y$. Then $n\ge 1$ and $y^{n}=\varepsilon (y')^{n+1}$, $\varepsilon\prec 1$, so
$$ny^{n-1}y'\ =\ \varepsilon'(y')^{n+1} + (n+1)\varepsilon (y')^{n}y''.$$
After dividing by $y'$ this gives
$$ny^{n-1}\ =\ \varepsilon'(y')^n + (n+1)\varepsilon (y')^{n-1}y''.$$ 
By the minimality of $n$ we have $\varepsilon' (y')^n \prec y^{n-1}$, so
$y^{n-1}\asymp \varepsilon(y')^{n-1}y''$, and hence 
$y^{n-1}y'\asymp \varepsilon(y')^{n}y''$, so by taking $n$th powers, 
$$y^{n(n-1)}(y')^n\ \asymp\ \varepsilon^n (y')^{n^2}(y'')^n.$$ 
Thus, using $y^n=\varepsilon(y')^{n+1}$,
$$ \varepsilon^{n-1}(y')^{(n+1)(n-1)}(y')^n\ \asymp\ \varepsilon^n(y')^{n^2}(y'')^n,$$
so $(y')^{n-1}\ \asymp\ \varepsilon (y'')^n$, contradicting the minimality 
of $n$. 
\end{proof} 

%\noindent
%In view of 
%Lemma~\ref{easyobs}, this gives: 

\begin{cor}\label{rk1cor} If $\Gamma$ has rank $1$, then $K$ is monotone. 
\end{cor}

\noindent
{\em In the remainder of this section $\alpha$,~$\beta$,~$\gamma$ range over $\Gamma$}.

\begin{cor}\label{vP-Lemma} Let $P\in K\{Y\}^{\ne}$ be homogeneous of degree $d$
and $\gamma\ne 0$. 
Then $v_P(\gamma)\ =\ v(P) + d\gamma+o(\gamma)$.  More generally, if $\alpha\ne  \beta$, then
$$v_P(\alpha)-v_P(\beta)\ =\ d\cdot (\alpha-\beta) + o(\alpha-\beta).$$
\end{cor}
\begin{proof} Let $\gamma>0$, take $g\in K^\times$
with $vg=\gamma$, and set $\Delta:=\big\{\delta\in\Gamma:\ \delta=o(\gamma)\big\}$. 
In the $\Delta$-coarsening of $K$ we have $\dot{v}g'\ge \dot{v}g$, so $\dot{v}_P(\dot{\gamma})=\dot{v}(P)+d\dot{\gamma}$ by Lemma~\ref{vplemma}, and thus $v_P(\gamma) =v(P) + d\gamma+o(\gamma)$.
Next, for $Q:= P_{\times g^{-1}}$,
$$v(P)\ =\ v_Q(\gamma)\ =\ v(Q) + d\gamma + o(\gamma), \qquad v(Q)\ =\ v_P(-\gamma),$$
so $v_P(-\gamma)=v(P) -d\gamma+ o(-\gamma)$. Thus the first part of
the lemma also holds if $\gamma< 0$. For the second part, let $\alpha\ne \beta$ and take
$a,b\in K^\times$ with $va=\alpha$, $vb=\beta$. Then 
$$v_P(\alpha) = v(P_{\times a}) = v(P_{\times b, \times (a/b)}) = v_{P_{\times b}}(\alpha-\beta),$$ 
so by the first part,
$$ v_P(\alpha)\ =\ v(P_{\times b}) + d(\alpha-\beta) + o(\alpha-\beta)\ =\ v_P(\beta)+d(\alpha-\beta) + o(\alpha-\beta),$$
which gives the desired result.
\end{proof}

\noindent
Here is some convenient notation. Let $\gamma>0$. Then
$\alpha\ge \beta+o(\gamma)$ is defined to mean that 
$\alpha\ge \beta+ \delta$ for some $\delta=o(\gamma)$
in $\Gamma$; equivalently, it means that $\alpha\ge \beta -\frac{1}{n}\gamma$ for all $n\ge 1$. We also declare that $\infty\ge \beta+ o(\gamma)$ for
$\infty\in \Gamma_{\infty}$. 
%In particular, for $y\in \smallo$, $y\ne 0$ we have
%$vy'\ge vy + o(vy)$.
%Here and below we use the convention that for 
%$\beta,\gamma\in \Gamma$ with
%$\gamma>0$, and $\alpha\in \Gamma_{\infty}$, the meaning 
%of $\alpha\ge \beta+o(\gamma)$ is that $\alpha \ge 
%\beta - \frac{1}{n}\gamma$ for all $n\ge 1$, that is,
%with $\Delta:=\big\{\delta\in \Gamma:\ 
%\delta=o(\gamma)\big\}$ we have
%$\alpha-\beta\in \Delta$ or $\alpha-\beta> \Delta$. 

\begin{cor} Let $P\in K\{Y\}^{\ne}$, $P(0)=0$, and $\gamma>0$. Then
$$v_P(\gamma)\ \ge \ v(P)+\gamma + o(\gamma).$$
\end{cor}

\begin{cor}\label{equalizer} Let $P, Q\in K\{Y\}^{\ne}$ be homogeneous 
of degrees $d$,~$e$.
Then
$$\gamma\ne 0\ \Longrightarrow\ v_P(\gamma)-v_Q(\gamma)\ =\ v(P)-v(Q) +(d-e)\gamma + o(\gamma).$$
Also, if $d>e$, then $v_P-v_Q\colon \Gamma \to \Gamma$ is strictly increasing.
\end{cor}
\begin{proof} The second part here follows from the second part of 
~\ref{vP-Lemma}.
\end{proof}

\begin{lemma}\label{alphaOgamma}
Let $P,Q\in K\{Y\}^{\neq}$ be homogeneous of the same degree and let $\gamma>0$ be such that
$v(P)\geq v(Q)+\gamma$. Then for $\alpha= O(\gamma)$ we have $v_P(\alpha) > v_Q(\alpha)$.
\end{lemma}
\begin{proof}
This is clear if $\alpha=0$, so let $0\neq\alpha=O(\gamma)$.  Then 
by Corollary~\ref{equalizer},
\equationqed{v_P(\alpha)-v_Q(\alpha)\ =\ 
v(P)-v(Q)+o(\gamma)\ \geq\ \gamma+o(\gamma)\ >\ 0.}
\end{proof}

\noindent
The next consequence of Lemma~\ref{alphaOgamma} is needed in
Section~\ref{sec:maxdh}.

\begin{cor}\label{cor:vA=vB}
Let $A,B\in K[\der]^{\neq}$, $B=A+E$, $\gamma>0$, $v(E)\geq v(A)+\gamma$, and $\alpha=O(\gamma)$. Then
$v_A(\alpha)=v_B(\alpha)$. Moreover, for any $a\in K^\times$ with $va=\alpha$ we have $A(a) \prec Aa\ \Longleftrightarrow\ B(a) \prec Ba$. Therefore, $\alpha\in\exc(A)\ \Longleftrightarrow\ \alpha\in\exc(B)$.
\end{cor}
\begin{proof}
By Corollary~\ref{alphaOgamma} we have $v_E(\alpha)>v_A(\alpha)$ and thus $v_A(\alpha)=v_B(\alpha)$.
Let $a\in K^\times$ be such that $va=\alpha$ and 
$A(a)\prec Aa$. Then $E(a) \preceq Ea \prec Aa\asymp Ba$, so
$B(a)=A(a)+E(a) \prec Ba$. We can interchange $A$ and $B$
in this argument.
%$v\big(A(a)\big)>v_A(\alpha)=v_B(\alpha)$. Then
%$v\big(E(a)\big)\geq v_E(\alpha)>v_A(\alpha)=v_B(\alpha)$, 
%so $v\big(B(a)\big)>v_B(\alpha)$, that is,
%$\alpha\in\exc(B)$.
\end{proof}

\noindent
To get useful information from Corollary~\ref{equalizer} concerning 
$v_F\colon \Gamma \to \Gamma$ when the differential polynomial $F\in K\{Y\}^{\ne}$ is not homogeneous, we
need two order-theoretic facts. To formulate those facts, let $X$ and $Y$ be nonempty ordered sets.

\begin{lemma} Let $f_0,\dots, f_n\colon X \to Y$ be strictly increasing bijections. Then the map $f\colon X\to Y$ given by $f(x):=\min_i f_i(x)$ is a
strictly increasing bijection.
\end{lemma}
\begin{proof} Define $g\colon Y\to X$ by $g(y):= \max_i f_i^{-1}(y)$. One verifies easily that then $g\circ f= \operatorname{id}_X$ and $f\circ g=\operatorname{id}_Y$.
\end{proof}

\noindent
Given maps $f,g\colon X \to Y$ we define $f<_c g$ to mean that
there are convex subsets $A< B < C$ in $X$ such that $X=A\cup B \cup C$,  $f>g$ on $A$, $f= g$ on $B$, and $f< g$ on~$C$. (Some of $A$, $B$, $C$ can be empty.)

\begin{lemma}\label{glueconvex} Let $f_0,\dots, f_n\colon X \to Y$ be such that $f_i <_c f_j$ for all $i<j$. Define $f\colon X \to Y$ by $f(x):=\min_i f_i(x)$. Then there are $i_0< \dots < i_m$ in $\{0,\dots,n\}$ and convex subsets 
$D_m<  \cdots < D_0$ of $X$ such that
\begin{enumerate}
\item[\textup{(i)}] $X=D_m\cup \dots \cup D_0$, and $D_m,\dots, D_0$ are nonempty;
\item[\textup{(ii)}] $f=f_{i_k}$ on $D_k$, for $k=0,\dots,m$.
\end{enumerate}
\end{lemma}
\begin{proof} For $n=0$ the lemma holds 
with $m=0$ and $D_0=X$. Let $n\ge 1$, and for 
$i=1,\dots,n$, take convex subsets $A_i< B_i < C_i$ in $X$ such that $X=A_i\cup B_i \cup C_i$, $f_0>f_i$ on $A_i$,
$f_0=f_i$ on $B_i$, and $f_0< f_i$ on $C_i$. 
Take $i^*\in \{1,\dots,n\}$ such that 
$C_{i^*}\subseteq C_{i}$ for $i=1,\dots,n$. Then $f_0< f_i$ on $D_0:=C_{i^*}$ for $i=1,\dots,n$, so if $D_0=X$, then the lemma holds with $m=0$, $i_0=0$. Suppose
that $D_0\ne X$. Then for $X':= X\setminus D_0$ we have $f_0\ge f_{i^*}$ on $X'$, and
then we can use as an inductive assumption that the lemma holds for the restrictions of $f_1,\dots, f_{n}$ to $X'$. 
\end{proof}

\begin{cor}\label{vPglueconvex} Let $F\in K\{Y\}^{\ne}$ have 
$F_{d_0}, \dots, F_{d_n}$ 
with $d_0 < \dots < d_n$ as its nonzero homogeneous parts. Set 
$f_i:= v_{F_{d_i}}\colon \Gamma \to \Gamma$ for $i=0,\dots,n$. Then there are $i_0< \dots < i_m$ in $\{0,\dots,n\}$ and convex
subsets $D_m < \cdots < D_0$ of $\Gamma$ such that
\begin{enumerate}
\item[\textup{(i)}] $\Gamma=D_m\cup \dots \cup D_0$, and $D_m,\dots, D_0$ are nonempty;
\item[\textup{(ii)}] $v_F(\gamma)=f_{i_k}(\gamma)$ for all $\gamma\in D_k$ and 
$k=0,\dots,m$.
\end{enumerate}
If $F(0)=0$ and each $f_i$ is bijective, then so is $v_F$.
\end{cor}
\begin{proof} By Corollary~\ref{equalizer} we have 
$f_i <_c f_j$ for $0\le i < j \le n$, in the 
sense defined earlier. We can now apply Lemma~\ref{glueconvex}, since
$v_F(\gamma)=\min_i f_i(\gamma)$
for all $\gamma$, and if $F(0)=0$, then each function 
$f_i$ is strictly increasing.
\end{proof}  
 
\subsection*{Notes and comments}
Lemma~\ref{cohnest} and Corollary~\ref{rk1cor} are in \cite{Cohn}, Section~2.
 
\section{Algebraic Extensions }\label{sec:algebraicext}

\noindent
Recall the standing assumption of this chapter that the derivation of 
$K$ is small, that is, $\der\smallo\subseteq \smallo$. At the end of this
section we prove:

\begin{prop}\label{Algebraic-Extensions} 
If $L$ is a valued differential field extension of $K$ and
$L$ is algebraic over $K$, then the derivation of $L$ is also small. 
\end{prop}

\noindent
We actually work in the more general setting of
an ambient valued differential field~$L$ with subfields 
$E\subseteq F$ such that $F$ is algebraic over $E$. In the lemmas 
below we consider $E$ and $F$ as {\em valued\/} subfields 
of $L$, and first deal with the case
that~$F|E$ is immediate, next the case that $E$ is henselian and
$F|E$ is unramified,
and finally the purely ramified case. It is important for 
later use in Section~\ref{sec:cimex} not to assume in the first two cases that $E$ or $F$ 
is closed under the derivation $\der$ of $L$. 

\subsection*{Immediate algebraic extensions}

\noindent
\begin{lemma}\label{key} Suppose $F|E$ is immediate, and 
$\der\smallo_E\subseteq \smallo_L$ and
$\der\mathcal{O}_E\subseteq \mathcal{O}_L$. 
Then $\der\smallo_F\subseteq \smallo_L$ and $\der\mathcal{O}_F\subseteq \mathcal{O}_L$. 
\end{lemma}
\begin{proof} Since we have $\mathcal{O}_F= \mathcal{O}_E +\smallo_F$,  it suffices to 
show that $\der\smallo_F\subseteq \smallo_L$. If $\mathcal{O}_E=E$, then $\mathcal{O}_F=F$ and we are done.
Assume $\mathcal{O}_E\ne E$. We can also arrange that~${1< {[F:E]}<\infty}$.
Then by Lemma~\ref{lem:min poly immediate ext} we have $F=E(y)$ where $y\in L^\times$, $y\prec 1$ and  
$y$ has minimum
 polynomial $$ P(Y)\ =\ Y^n+a_{n-1}Y^{n-1}+\cdots+a_1Y+a_0$$ 
over $E$ with
coefficients $a_i\in \mathcal{O}_E$, $a_1\asymp 1$,  $a_0\prec 1$.
(Note that $a_0\ne 0$, since $n={[F:E]}>1$.)

\claim[1]{We have $y\asymp a_0$, and $y' \preceq (aa_0)'$ for some $a\in \mathcal{O}_E$.}

\noindent
We get $y\asymp a_0$ from $a_1y\sim -a_0$. Next,
$P(y)=0$ gives $$0\ =\ P(y)'\ =\ P^{\der}(y) + P'(y)y',$$ and since
$P'(y)\sim a_1\asymp 1$, this gives
$$y'\asymp P^{\der}(y)=\sum_{i=0}^{n-1} a_i'y^i,$$ so we get $i\in \{0,\dots,n-1\}$
with $y'\preceq a_i'y^i$. If $i=0$, this gives $y'\preceq a_0'$.
If $i\ge 1$, then we get $y'\preceq a_i'a_0=(a_ia_0)'-a_ia_0'$, so
$y'\preceq (a_ia_0)'$, or $y'\preceq a_ia_0'\preceq a_0'$.

\medskip\noindent
By Lemma~\ref{lem:min poly immediate ext, addition}, 
$$\smallo(y)\ :=\ \smallo_E\mathcal{O}_E[y] +y\mathcal{O}_E[y]$$
is a maximal ideal of $\mathcal{O}_E[y]$, and 
\begin{align*} \mathcal{O}_F\ &=\  \mathcal{O}_E[y]_{\smallo(y)}\ =\
S^{-1}\mathcal{O}_E[y]\ \text{ with $S:=1+\smallo(y)$,}\\
\smallo_F\ &=\ \smallo(y) \mathcal{O}_F\ =\ S^{-1}\smallo(y).
\end{align*}

\claim[2]{For each $b\in  \mathcal{O}_E$ there exists $u\in  \mathcal{O}_E$ such that $(by)' \preceq (ua_0)'$.}

\noindent
To see why, let $b\in\mathcal{O}_E$, so $(by)'=b'y + by'$. 
Now $b'y\asymp b'a_0=(ba_0)'-ba_0'$, so $b'y\preceq (ba_0)'$ or
$b'y\preceq ba_0'\preceq a_0'$. Also, with $a$ as in Claim 1 we have
$by'\preceq (aa_0)'$. Thus $(by)'\preceq (ua_0)'$ for $u=b$ or $u=1$ or
$u=a$. 

\claim[3]{For each $\phi\in \smallo(y)$ there exists $\varepsilon\in \smallo_E$ such that $ \phi' \preceq \varepsilon'$.}

\noindent
This property holds for $\phi\in \smallo_E$, and for $\phi=y$ by Claim~1,
and it is inherited under taking sums and products of elements in
$\smallo(y)$. This yields Claim~3 using also Claim~2. 

\medskip\noindent
Let $f\in \smallo_F$. Then
$f=\phi/(1+e)$
with $\phi,e\in  \smallo(y)$; differentiating this 
quotient, we
obtain from Claim 3 that there is 
$\varepsilon\in \smallo_E$ such
that $f' \preceq \varepsilon'$.
\end{proof}

\noindent
Here is a further reduction of the problem of showing 
$\der\smallo_F\subseteq \smallo_L$ and $\der\mathcal{O}_F\subseteq \mathcal{O}_L$:

\begin{lemma}\label{keycor} Suppose $F|E$ is immediate and $S\subseteq \mathcal{O}_E$ satisfies \
$$\mathcal{O}_E\ =\ \{f/g:\ f,g\in S,\ g\asymp 1\}, \quad
\der S\subseteq \mathcal{O}_L, \quad \der (S\cap \smallo_E)\subseteq \smallo_L.$$
Then $\der\smallo_F\subseteq \smallo_L$ and $\der\mathcal{O}_F\subseteq \mathcal{O}_L$.
\end{lemma}
\begin{proof} Let $a\in \mathcal{O}_E$, so $a=f/g$ with $f,g\in S$, $g\asymp 1$.
Then $f', g'\preceq 1$, and thus
$a'=(f'g-fg')/g^2\preceq 1$. 
Also, if $a\prec 1$, then $f\prec 1$, so $f'\prec 1$, and thus $a'\prec 1$.
Now apply Lemma~\ref{key}. 
\end{proof}

\subsection*{Unramified and purely ramified algebraic extensions}

\begin{lemma}\label{Residue-Field-Extension}
Suppose $E$ is henselian, $\Gamma_F=\Gamma_E$,
and  $\der \smallo_E\subseteq \smallo_L$ and $\der\mathcal{O}_E\subseteq \mathcal{O}_L$.  Then $\der \smallo_F\subseteq \smallo_L$ and $\der\mathcal{O}_F\subseteq \mathcal{O}_L$.
\end{lemma}
\begin{proof} We can arrange $[F:E]=\bigl[\res(F):\res(E)\bigr]=n>1$.
Take~${y\asymp 1}$ in~$F$ such that $\res(F)=\res(E)[\bar{y}]$, where
$\bar{y}$ is the residue class of~$y$ in $\res(F)$. Corollary~\ref{cor:alg ext of henselian} gives 
$a_0,\dots,a_{n-1}\preceq 1$ in $E$ such that
$$P(Y)=Y^n+a_{n-1}Y^{n-1}+\cdots +a_1Y+a_0\in E[Y]$$ is the
minimum polynomial of $y$ over $E$, so its reduction $\bar{P}(Y)$ in
$\res(E)[Y]$ is the minimum polynomial of $\bar{y}$ over $\res(E)$; 
in particular $P'(y)\asymp 1$. Moreover, 
$$P'(y)y'\ =\ -\sum_{i=0}^{n-1}a_i'y^i,$$ 
hence $y' \preceq 1$. Let $f\in F$, and take 
$f_0,\dots,f_{n-1}\in E$ such that
\begin{align*}f\ &=\ f_0+f_1y+\cdots+f_{n-1}y^{n-1},\ \text{ so}\\ 
 f'\ &=\  \sum_{i=0}^{n-1} f_i'y^i + 
\left(\sum_{j=1}^{n-1} jf_jy^{j-1}\right)y'.
\end{align*}
If $f\prec 1$, then $f_0,\dots, f_{n-1}\prec 1$, and thus $f'\prec 1$. Likewise, if $f\preceq 1$, then $f_0,\dots, f_{n-1}\preceq 1$, and thus $f'\preceq 1$.
\end{proof}

\begin{lemma}\label{Purely ramified} Let $F=E\bigl(u^{1/p}\bigr)$
where $p$ is a prime number, $u\in E^{\times}$ and $vu\notin p\Gamma_E$, and
assume that $\der E\subseteq E$ and $\der\smallo_E\subseteq \smallo_E$.
Then $\der\smallo_F\subseteq \smallo_F$.
\end{lemma}
\begin{proof}
Let $u^{i/p} := (u^{1/p})^i$ for $i\in\Z$.
By Lemma~\ref{lem:pth root},  
$\res(E)=\res(F)$. Let $a\in F^\times$, $a\prec 1$; our job is to show that
$a'\prec 1$. We have  
$$a\ =\ a_0+ a_1u^{1/p} + \cdots +a_{p-1}u^{(p-1)/p}$$ 
with $a_0,\dots, a_{p-1}\in E$. Then $a_iu^{i/p}\prec 1$ for $i=0,\dots,p-1$,
so we may reduce to the case $a=a_iu^{i/p}$, and so $a^p=b\in E^\times$.
Then $pa^{p-1}a'=b'$, and so by Lemma~\ref{cohnest}, 
$$va'\ =\ vb' -(p-1)va\ =\ vb'-\big((p-1)/p\big)vb\  >\ 0,$$
as desired.
\end{proof}

\begin{proof}[Proof of Proposition~\ref{Algebraic-Extensions}] Let the
valued differential field extension $L$ of $K$ be algebraic over $K$. Note that any field between $K$ and $L$ is closed under the derivation of $L$.
By extending $L$ we arrange that $L$ is an algebraic closure of~$K$. Next, by  
Lemma~\ref{key}, we arrange that $K$ is henselian. We then reach $L$ in two steps. In the first step we pass from~$K$ to its maximal unramified extension
$K^{\operatorname{unr}}$ in $L$. Then the derivation on $K^{\operatorname{unr}}$ is small by 
Lemma~\ref{Residue-Field-Extension}. In the second step we apply Proposition~\ref{prop:ramified} and Lemma~\ref{Purely ramified} to $L$ as an
extension of~$K^{\operatorname{unr}}$. 
\end{proof}

\section{Residue Extensions}\label{sec:resext}

\index{gaussian extension}
\index{extension!gaussian}

\noindent
The differential field $K\<Y\>$ with the gaussian
extension of the valuation of $K$ is a valued differential field extension of $K$, with
differential subring $\mathcal{O}\{Y\}$ of $K\{Y\}$, and we claim that
\begin{enumerate} 
\item $\mathcal{O}\{Y\}\subseteq \mathcal{O}_{K\<Y\>}$ and 
$\der\smallo_{K\<Y\>} \subseteq \smallo_{K\<Y\>}$;
\item $Y\preceq 1$, and the image $y$ of $Y$ in the differential 
residue field of $K\<Y\>$
is $\d$-transcendental over the differential residue field $\k$ of $K$;
\item $\k\<y\>$ is the differential residue field of $K\<Y\>$.
\end{enumerate}
For (1), consider first the case of $P\in K\{Y\}$ with $vP>0$. Then $vP'>0$, 
because $P=gF$ with $g\in \smallo$ and $F\in \mathcal{O}\{Y\}$, so $P'=g'F+ gF'$.
Next, let $f(Y)\in \smallo_{K\<Y\>}$. Then $f=P/Q$ with $P,Q\in K\{Y\}$ and
$vQ=0$ and $vP>0$, so $f'=(P'Q-PQ')/Q^2$, and $v(P'Q)=vP' >0$ and 
$v(PQ')\ge vP>0$, so $vf'>0$. For (2), let $g(Y)\in \k\{Y\}^{\ne}$
and take $P\in \mathcal{O}\{Y\}$ with $\overline{P}=g$. Then $vP=0$ and so~$g(y)$, the image of $P$ under the residue map, is~$\ne 0$. For (3), let 
$f(Y)\in K\<Y\>$ and $vf=0$.
Then $f=P/Q$ with $P,Q\in K\{Y\}$ and $vP=vQ=0$, so the image of~$f$ in the residue field of $K\<Y\>$ is $\overline{P}(y)/\overline{Q}(y)\in \k\<y\>$.

\medskip\noindent
This valued differential field extension $K\<Y\>$ has the following universal
property:

\begin{lemma}\label{univdiftr} Let $L$ be any valued differential field extension of $K$
with~${\der\smallo_{L} \subseteq \smallo_{L}}$ and with an element $a\preceq 1$
whose image in the differential residue field of $L$
is $\d$-transcendental over $\k$. Then there is a unique
valued differential field embedding $K\<Y\>\to L$ over $K$ sending $Y$ to $a$.
\end{lemma}

\noindent
Here is an analogue of the above for \textit{$\d$-algebraic}\/ residue field extensions:

\begin{theorem}\label{resext} Let $K\<a\>$ be a
differential field extension of $K$ such that $a$ is $\d$-al\-ge\-bra\-ic over $K$ with
minimal annihilator $F$ over $K$ of order $r$. Assume that~${vF=0}$ and~${vI=0}$ where $I$ is the
initial of $F$, and that the image $\overline{F}$ of $F$ in $\k\{Y\}$ is irreducible. Then $v$ extends uniquely to a valuation 
$v\colon K\<a\>^\times \to \Gamma$ such that 
$\der\smallo_{K\<a\>} \subseteq \smallo_{K\<a\>}$, $a\preceq 1$, and $\res{a}$ is
$\d$-algebraic over $\k$ with minimal annihilator 
$\overline{F}$ over $\k$. For this extended valuation the differential residue 
field of $K\<a\>$ is $\k\<\res{a}\>$.   
\end{theorem}
\begin{proof} Let $F$ have order $r$, and set $d:= \deg_{Y^{(r)}}F$, so $d\ge 1$. Then $\overline{F}$ has order~$r$ and $d= \deg_{Y^{(r)}}\overline{F}$. Note that if $b\in K\<a\>^\times$, then $b=P(a)/Q(a)$ where $Q\in K\big[Y,\dots, Y^{(r-1)}\big]{}^{\ne}$ and $P\in K\big[Y,\dots, Y^{(r)}\big]{}^{\ne}$ with 
$\deg_{Y^{(r)}}P < d$.

Suppose an extension of the valuation of $K$ to a valuation $v$ of $K\<a\>^\times$ as in the theorem is given. Let 
$P\in K\big[Y,\dots, Y^{(r)}\big]{}^{\ne}$ and $\deg_{Y^{(r)}}P < d$;
we claim that then $vP(a)=vP$. To prove this claim we can multiply $P$ by an element of $K^\times$ and arrange that $vP=0$. Then $\overline{P}\ne 0$ and
$\deg_{Y^{(r)}}\overline{P} < d$, so $\res{P(a)}=\overline{P}(\res{a})\ne 0$, and
thus $vP(a)=0$, as claimed. Thus there can be at most one such extension.

To construct such an extension, let $b=P(a)/Q(a)$ be as above; we claim that then
$vP-vQ$ depends only on $b$ and not on the choice of $P$ and $Q$. To prove this claim,
let also $b=G(a)/H(a)$ where $H\in K\big[Y,\dots, Y^{(r-1)}\big]{}^{\ne}$ and $G\in K\big[Y,\dots, Y^{(r)}\big]{}^{\ne}$ with $\deg_{Y^{(r)}}G < d$. Then 
$$\deg_{Y^{(r)}} HP-QG\ <\  d, \qquad 
(HP-QG)(a)\ =\ 0,$$ so $HP=QG,$ and thus $vP-vQ=vG-vH$ as claimed.
This allows us to extend the valuation
$v$ on $K$ to a function $v\colon K\<a\>^\times \to \Gamma$ by setting
$vb:= vP-vQ$ when $b=P(a)/Q(a)$ for $b$, $P$, $Q$ as above. Next we claim that
$v\colon K\<a\>^\times \to \Gamma$ is a valuation on the field $K\<a\>$.
To prove this, let $b,c\in K\<a\>^{\ne}$, and $b=P(a)/Q(a)$ as above, and also $c=G(a)/H(a)$ with $H\in K\big[Y,\dots, Y^{(r-1)}\big]{}^{\ne}$ and $G\in K\big[Y,\dots, Y^{(r)}\big]{}^{\ne}$ with $\deg_{Y^{(r)}}G < d$. We can arrange $Q=H$. Then $b+c= (P+G)(a)/Q(a)$, 
and so, if $b+c\ne 0$, 
$$v(b+c)\ =\ v(P+G)-vQ\ \ge\ \min (vP-vQ, vG-vQ)\ =\ \min(vb, vc).$$
As to $v(bc)=vb+vc$, this holds if $G\in K\big[Y,\dots, Y^{(r-1)}\big]$, because
then $bc=(PG)(a)/Q^2(a)$. In particular, it holds for 
$c\in K\big(a,\dots, a^{(r-1)}\big)$. This allows us to reduce the general case to
the case that $b=P(a)$ and $c=G(a)$ with $P$,~$G$ as above satisfying the additional condition $vP=vG=0$ (so $vb=vc=0$).  
Since $vI=0$, division with remainder in $\mathcal{O}\big[Y,\dots, Y^{(r)}\big]$ gives 
$I^mPG=BF+R$ with $B, R\in \mathcal{O}\big[Y,\dots, Y^{(r)}\big]$ and $\deg_{Y^{(r)}} R < d$. Then $vR=0$, since $vR>0$ would give $\overline{I}{}^m\cdot \overline{P}\cdot \overline{G}= 
\overline{B} \cdot \overline{F}$ in $\k\{Y\}$, contradicting the irreducibility of
$\overline{F}$. Now $bc=R(a)/I^m(a)$, so $v(bc)=vR- mvI=0$, as required. We have 
now shown that~$v$ is a valuation on $K\<a\>$.

Towards proving that 
$\der\smallo_{K\<a\>} \subseteq \smallo_{K\<a\>}$ we first note that 
$a^{(i)}\preceq 1$ for $i< r$. We have 
$F_0,\dots, F_d\in \mathcal{O}\big[Y,\dots, Y^{(r-1)}\big]$ with $F_d=I$ such that
$$ F\ =\ F_d\cdot \big(Y^{(r)}\big)^d +F_{d-1}\cdot \big(Y^{(r)}\big)^{d-1} + \cdots + F_0.$$
Thus $a^{(r)}$ is integral over the subring 
$\mathcal{O}\big[a,\dots, a^{(r-1)}\big]_{I(a)}$ of the valuation ring~$\mathcal{O}_{K\<a\>}$ of $K\<a\>$, so
$a^{(r)}\preceq 1$. In view of $f_i:= F_i(a)$ and $v(F_i)\ge 0$, this gives
$f_i'\preceq 1$ for $i=0,\dots, d$. Now with $g:= a^{(r)}$ we have
$$0\ =\ F(a)'\ =\ f_d'g^d + f_{d-1}'g^{d-1} + \cdots + f_0' + \frac{\partial F}{\partial Y^{(r)}}(a) \cdot g',$$
with $\deg_{Y^{(r)}}\frac{\partial F}{\partial Y^{(r)}} < d$ and 
$v\frac{\partial F}{\partial Y^{(r)}} = 0$ (using $vI=0$), so 
$g'=a^{(r+1)}\preceq 1$. Now let  
$b=P(a)/Q(a)$ be as above with $b\in \smallo_{K\<a\>}$. We can arrange
$vP>0$ and $vQ=0$, and then 
$$b'\ =\ \frac{P(a)'Q(a)-P(a)Q(a)'}{Q(a)^2},$$
which in view of $a^{(i)}\preceq 1$ for $i=0,\dots, r+1$ yields
$P(a)' \prec 1$, $Q(a)'\preceq 1$, in addition to $P(a)\prec 1$ and 
$Q(a)\asymp 1$. Thus $b'\in \smallo_{K\<a\>}$, as promised. From $F(a)=0$ we 
get $\overline{F}(\res{a})=0$. Suppose $A\in \k\{Y\}^{\ne}$ 
has order at most $r$ and $\deg_{Y^{(r)}} A < d$. Take 
$P\in \mathcal{O}\{Y\}^{\ne}$ of order at most $r$ with $\deg_{Y^{(r)}}P < d$ 
such that $\overline{P}=A$. Then $vP=0$, so $vP(a)=0$, and thus 
$A(\res{a})=\res{P(a)}\ne 0$. Thus $\overline{F}$ is
a minimal 
annihilator of~$\res{a}$ over $\k$.

Now let $b=P(a)/Q(a)$ as before be such that $vb=0$. Then we can assume
$vP=vQ=0$, so $\res{b} = \overline{P}(\res{a})/\overline{Q}(\res{a})\in \k\<\res{a}\>$.
Thus the differential residue field of $K\<a\>$ is $\k\<\res{a}\>$.   
\end{proof}

\noindent
This gives a partial analogue of Proposition~\ref{prop:prescribed} for valued differential fields:

\begin{cor}\label{cor:resext}
Let $\k_L$ be a differential field extension of $\k$. Then $K$ has a valued differential field extension $L$ with $\der\smallo_L\subseteq\smallo_L$, with the same value group as $K$ and with differential residue field isomorphic to $\k_L$ over $\k$.
\end{cor}
\begin{proof} We can reduce to the case $\k_L=\k\<y\>$. 
If $y$ is $\d$-transcendental over $\k$, the corollary holds with 
$L=K\<Y\>$ equipped with the gaussian extension of the valuation of $K$. 
Next, suppose $y$ is $\d$-algebraic over $\k$, with minimal
annihilator $\overline{F}(Y)\in \k\{Y\}$ over~$\k$. Take $F\in \mathcal{O}\{Y\}$
to have image $\overline{F}$ in $\k\{Y\}$ and to have the same complexity as~$\overline{F}$.
Then  $vF=vI=0$, where $I$ is the initial of $F$, and~$F$ is irreducible in $K\{Y\}$. Take a differential field extension 
$L=K\<a\>$ of $K$ such that $a$
has minimal annihilator~$F$ over~$K$. Then $L$ with the valuation defined
in Theorem~\ref{resext} has the desired properties.
\end{proof}

\subsection*{Preserving monotonicity} Let $L$ be an ambient valued differential field with derivation
$\der$, and consider subfields of $L$ as valued subfields. Let $E$ be a subfield of~$L$ such that $a'\preceq a$ for all $a\in E$.
In particular $\der\smallo_E\subseteq \smallo_L$ and 
$\der\mathcal{O}_E\subseteq \mathcal{O}_L$. We do not assume $\der E \subseteq E$. After proving some lemmas we shall derive:

\begin{prop}\label{monalgext} If $b\in L$ is algebraic over $E$, 
then $b'\preceq b$.
\end{prop}

\begin{lemma}\label{sc2}\label{stas11} Let $F\supseteq E$ be a subfield of $L$
with $\der\mathcal{O}_F \subseteq \mathcal{O}_L$ and $\Gamma_E=\Gamma_F$. 
Then $b'\preceq b$ for all $b\in F$.
\end{lemma}
\begin{proof} Let $b\in F^\times$, and take 
$a\in E^\times$, $u\in F$, with
$b=au$ and $u\asymp 1$. Now use $a'\preceq a\asymp b$ and $u'\preceq 1$ to get
$b'=a'u+au'\preceq b$.
\end{proof}

\begin{cor}\label{cor:sc2}
If $L$ has small derivation and a monotone valued differential subfield $K$ with $\Gamma_L=\Gamma$, then $L$ is monotone.
\end{cor}
\begin{proof}
Apply Lemma~\ref{sc2} with $E=K$, $F=L$.
\end{proof}

\begin{cor}\label{cor:monotonicity under extensions} If $K$ is monotone, then so is the differential field $K\<Y\>$ with the gaussian valuation, as well as any extension $K\<a\>$ as in  
Theorem~\ref{resext},
and any extension $L$ as in Corollary~\ref{cor:resext}.
\end{cor}

{\sloppy
\begin{lemma}\label{lem:S1 on subspace}  
Let $\mathcal B$ be a valuation-independent subset of the valued vector space~$L$ over the 
valued field $E$, such that
$b'\preceq b$ for all $b\in\mathcal B$. Then $f'\preceq f$ for all $f$ in 
the $E$-linear span of $\mathcal B$ in $L$.
\end{lemma}
\begin{proof}
Let $f=\sum_{i=1}^n a_i b_i\neq 0$ with $a_1,\dots, a_n\in E$ and distinct 
$b_1,\dots, b_n\in\mathcal B$. Then $v(f)= \min_i v(a_ib_i)$ and
$f'= \sum a_i' b_i + \sum a_i b_i'$, and 
\begin{align*} v\left( \sum_i a_i' b_i\right)\ &\geq\ 
\min_i v(a_i'b_i)\ \geq\ \min_i v(a_ib_i)\ =\ v(f),\\ 
v\left(\sum_i a_i b_i' \right)\ &\geq\ \min_i v(a_ib_i')\  \geq\ \min_i v(a_ib_i)\ =\ v(f),
\end{align*}
so $v(f') \ge v(f)$.  
\end{proof}

}

\begin{lemma}\label{stas12} Let $y\in L$ be such that
 $y^p=a\in E^\times$ with $p$ a prime number
and $vy\notin \Gamma_E$. Then $b'\preceq b$ for all $b\in E(y)$. 
\end{lemma}
\begin{proof} We have $pvy=va$ and $py^{p-1}y'=a'$, so 
$$y'\ \asymp\ a'/y^{p-1}\ \preceq\ a/y^{p-1}\ =\ y.$$
The desired result now follows from Lemmas~\ref{lem:pth root} and~\ref{lem:S1 on subspace}, the latter applied to
$\mathcal{B}=\{1,y,\dots, y^{p-1}\}$.
\end{proof}

\begin{proof}[Proof of Proposition~\ref{monalgext}] By passing to an
algebraic closure we arrange that~$L$ is algebraically closed. Let $E^\alg$ be the algebraic closure of
$E$ in $L$, let $E^{\operatorname{h}}$ be the henselization of $E$ in $E^\alg$, and let
$F$ be the maximal unramified extension of $E^{\operatorname{h}}$ in~$E^\alg$.
Then $\der\mathcal{O}_F\subseteq \mathcal{O}_L$ by Lemmas~\ref{key} and~\ref{Residue-Field-Extension}, and so $b'\preceq b$ for all $b\in F$ by Lemma~\ref{stas11}. 
Now use Lemma~\ref{stas12}
to conclude that $b'\preceq b$ for all $b\in E^\alg$. 
\end{proof}  

\begin{cor}\label{cor:monalgmon} If $K$ is monotone, then every valued differential field
extension of $K$ that is algebraic over $K$ is also monotone.
\end{cor}

\begin{cor}\label{cor:S1 and trdeg 1} 
Suppose $E\subseteq \mathcal{O}_L$ and $L$ has transcendence degree $1$ 
over $E$.  Then $L^\phi$ is monotone for some $\phi\in L^\times$.
\end{cor}
\begin{proof} The case $\Gamma_L=\{0\}$ is trivial, so
assume $\Gamma_L\neq\{0\}$. Take $y\in L^{\times}$ with $y\prec 1$. 
Then $y$ is transcendental over $E$. Let $\mathcal{B}=\{1, y, y^2, \dots\}$
and note that~$\mathcal{B}$ is valuation-independent over $E$ and
$L$ is algebraic over $E(y)$. If $y'\preceq y$, then $b'\preceq b$ for
all $b\in E(y)$ by Lemmas~\ref{lem:S1 on subspace} and~\ref{easyobs}, and so 
$L$ is monotone by Proposition~\ref{monalgext}. Suppose that $y'\succ y$. Then we set $\phi=y^\dagger$, and 
$\derdelta=\phi^{-1}\der$, so $\derdelta(a) \preceq a$ for all $a\in E$.
Moreover $\derdelta(y)=y$. Thus the previous argument applies to
 $L^\phi$ instead of $L$, so $L^{\phi}$ is monotone.
\end{proof}

%\begin{exampleNumbered}\label{ex:S1 and trdeg 1}
%Suppose $E\subseteq\mathcal O_L$ and $L=E(y)$ where $y\nasymp %1$ and $y'\preceq y$.
%Then~$L$ is monotone.
%\end{exampleNumbered}

%\begin{cor}\label{cor:extend to many const}
%Let $K$ be a monotone valued differential field. Then $K$ has a monotone valued differential field extension with many constants and the same value group as $K$.
%\end{cor}
%\begin{proof}
%Let $g\in K^\times$; it suffices that then $K$ has a monotone valued differential field extension with the same value group as $K$ and constant $c\asymp g$.
%Put $F:=Y'+g^\dagger Y\in K\{Y\}$, and use Lemma~\ref{lem:construct diff alg extension} to get an element $a$ in a differential field
%extension of $K$ with minimal annihilator $a$ over $K$. Equip $K\<a\>$ with the valuation $v\colon K\<a\>^\times\to\Gamma$ from Theorem~\ref{resext}. Then $a\asymp 1$ and $K\<a\>$ is monotone by Corollary~\ref{cor:sc2} and $c:=ag\asymp g$ is a constant.
%\end{proof}

\subsection*{Notes and comments} Proposition~\ref{monalgext}  and Corollary~\ref{cor:S1 and trdeg 1} are slight
extensions of \cite[Theorem~4.1]{Morrison}, 
\cite[Theorem~3]{Duval} and \cite[Proposition~5.1]{Morrison}.

\section{The Valuation Induced on the Value Group}\label{sec:dvvg}

\noindent
Let $a$,~$b$ range over $K$ and 
$\alpha$,~$\beta$,~$\gamma$ over $\Gamma$. Note: if $a'\preceq a$ and 
$a\asymp b$, then 
$b'\preceq b$. Also, if $a'\succ a$ and $a\asymp b$, then
$b'\succ b$ and
$a^{\dagger} \sim b^{\dagger}$, and thus $a^{\dagger} \asymp b^{\dagger}$. This allows us to define a function
$\psile\colon \Gamma^{\ne} \to \Gamma$ as follows: for $\alpha\ne 0$, $\alpha=va$,
$$
 a'\preceq a\ \Rightarrow\ \psile(\alpha):=0, \qquad a'\succ a\ \Rightarrow\  \psile(\alpha):= v(a^\dagger)=v(a')-v(a)\ < 0.$$
Thus $\psile(\Gamma^{\ne})\subseteq \Gamma^{\le}$. We extend $\psile$ to all of $\Gamma$ by
$\psile(0):= \infty$. 

\nomenclature[T]{$\psile$}{valuation on the value group}

\begin{lemma}\label{valv} The function $\psile\colon \Gamma\to \Gamma_{\infty}$ has the following properties: \begin{enumerate}
\item[\textup{(i)}] $\psile$ is a \textup{(}non-surjective\textup{)} valuation on the abelian group 
$\Gamma$; 
\item[\textup{(ii)}] $\psile(k\alpha)=\psile(\alpha)$ for nonzero $k\in \Z$;
\item[\textup{(iii)}] $\psile(\alpha)=o(\alpha)$ for $\alpha\ne 0$;
\item[\textup{(iv)}] if $a'\succ a$ and $va=\alpha$, then $v(a^{(n)})=\alpha+n\psile(\alpha)$ for
all $n$. 
\end{enumerate}
\end{lemma}
\begin{proof} It is clear that $\psile(-\alpha)=\psile(\alpha)$.
Let $\alpha,\beta\ne 0$, and take $a$,~$b$ such that ${va=\alpha}$, $vb=\beta$. Then
$v(ab)=\alpha+\beta$, and by considering separately the cases 
$(ab)'\preceq ab$ and $(ab)'\succ ab$ we obtain 
$\psile(\alpha+\beta) \ge \min\!\big(\psile(\alpha),\psile(\beta)\big)$.
%$$\psile(\alpha+\beta)\ \ge\  v\big((ab)^\dagger\big)\ 
%=\ v(a^\dagger + b^\dagger)\ \ge\ \min (va^\dagger, vb^\dagger)\ 
%\ge\ \min\big(\psile(\alpha),\psile(\beta)\big).$$
This proves~(i), and~(ii) follows from $(a^k)^\dagger=ka^\dagger$ for $a\ne 0$ and
$k\in \Z$. As to (iii), replacing~$\alpha$ by $-\alpha$ if necessary, we arrange $\alpha>0$, and then taking $a$ with $va=\alpha$, we have 
$v(a^\dagger)\ge -\alpha/n$ for all $n\ge 1$, by Lemma~\ref{cohnest}, which gives the desired estimate. 

Let $a'\succ a$ and set $\alpha:=va$, so $\psile(\alpha) < 0$.
It is obvious that for $n=0,1$ we have $v\big(a^{(n)}\big)=\alpha+n\psile(\alpha)$. Assume 
this equality holds 
for a certain $n\ge 1$. Applying~(iii) to~$\psile(\alpha)$ in the role of 
$\alpha$ gives 
$\psile\big(\psile(\alpha)\big)=o\big(\psile(\alpha)\big)$, so $\psile(\alpha) < \psile\big(\psile(\alpha)\big)=\psile\big(n\psile(\alpha)\big)$, 
and thus $\psile\big(\alpha + n\psile(\alpha)\big)=\psile(\alpha) < 0$. Therefore,
$$v\big(a^{(n+1)}\big)\ =\ \alpha+n\psile(\alpha) + \psile\big(\alpha+n\psile(\alpha)\big)\ =\  \alpha+(n+1)\psile(\alpha),$$
which proves (iv).
\end{proof}  

\noindent
In the next section we show that $(\Gamma, \psile)$ is a so-called asymptotic couple,
and there we establish some useful general facts about asymptotic couples.
The above lemma is good enough for the application that follows now.

\subsection*{A more precise estimate on $v_P$} We first consider
the case $P=Y^{(n)}$.

\begin{lemma}\label{asivp2} 
Suppose $\psile(\gamma)< 0$. Then
\begin{align*} 
 v_{Y^{(n)}}(\gamma)\ &=\ \gamma + n\psile(\gamma),\\
  v_{Y^{[\bsigma]}}(\gamma)\ &=\ d\gamma + \|\bsigma\|\psile(\gamma)\ 
\text{ for $\bsigma=\sigma_1\dots \sigma_d$.} 
\end{align*}
\end{lemma}
\begin{proof} To obtain 
the first identity, take $g\in K^\times$ with $vg=\gamma$. Then  \begin{align*}\big(Y^{(n)}\big)_{\times g}\ &=\ (gY)^{(n)}=\sum_{i=0}^n \binom{n}{i}g^{(i)}Y^{(n-i)},\ \text{ so}\\ 
v_{Y^{(n)}}(\gamma)\ &=\
\min_{0\le i\le n} vg^{(i)}\ =\ \min_{0\le i \le n}\gamma + i\psile(\gamma)\ =\ \gamma  + n\psile(\gamma).
\end{align*}
The second identity follows from the first using 
$Y^{[\bsigma]}=Y^{(\sigma_1)}\cdots Y^{(\sigma_d)}$.  
\end{proof}

\noindent
Let $P\in K\{Y\}^{\ne}$ be homogeneous of degree $d$. Lemma~\ref{vplemma} says that if $\psile(\gamma)\ge 0$, then $v_P(\gamma)=v(P)+d\gamma$. Here we deal 
with the case $\psile(\gamma)<0$:

\begin{prop}\label{ps1} Suppose $\psile(\gamma) <0$. Then 
$$vP+ d\gamma + \wt(P)\psile(\gamma)\ \le\  v_P(\gamma)\ \le\  vP +d\gamma +\big(\!\wt(P)-\Pnu(P)\big)\cdot \abs{\psile(\gamma)}.$$ 
In particular, 
$v_P(\gamma) = vP + d\gamma + O\big(\psile(\gamma)\big)$.
\end{prop}
\begin{proof} If $P=aY^{\bsigma}$, $a\in K^\times$, then by Lemma~\ref{asivp2},
$$ v_P(\gamma)\ =\ va+ d\gamma + \|\bsigma\|\psile(\gamma).$$
This gives the lower bound for all $P$. For the upper bound we first 
prove: 

\claim{Let $\bsigma_0$ be such that $P_{[\bsigma_0]}\ne 0$. Let $\bsigma\ge \bsigma_0$ be maximal with the property
$vP_{[\bsigma]}\le vP_{[\bsigma_0]}+ \|\bsigma -\bsigma_0\|\cdot \abs{\psile(\gamma)}$, and let $vg=\gamma$. Then 
$v\left((P_{\times g})_{[\bsigma]}\right)\ =\ vP_{[\bsigma]}+d\gamma$.}

%\vskip-1.5em
\noindent
Towards proving the claim, recall that
$$ (P_{\times g})_{[\bsigma]}\ =\ \sum_{\btau\ge \bsigma} \binom{\btau}{\bsigma}
P_{[\btau]}g^{[\btau-\bsigma]}\ =\ P_{[\bsigma]}g^d + \sum _{\btau> \bsigma} 
\binom{\btau}{\bsigma}P_{[\btau]}g^{[\btau-\bsigma]}.$$  
Let $\btau > \bsigma$, and note that
maximality of $\bsigma$ gives 
$$ vP_{[\btau]}\ >\  vP_{[\bsigma_0]} + \|\btau-\bsigma_0\|\cdot \abs{\psile(\gamma)},$$
which in view of $\|\btau-\bsigma_0\|=\|\btau-\bsigma\| + \|\bsigma-\bsigma_0\|$ gives $vP_{[\btau]} > vP_{[\bsigma]}+\|\btau-\bsigma\|\cdot \abs{\psile(\gamma)}$. Thus, using $|\psile(\gamma)|=-\psile(\gamma)$, we get
\begin{align*} v\big(P_{[\btau]}g^{[\btau-\bsigma]}\big)\ &=\ vP_{[\btau]}+d\gamma + \|\btau-\bsigma\|\cdot \psile(\gamma)\\
    &>\ vP_{[\bsigma]} +  \|\btau-\bsigma\|\cdot\abs{\psile(\gamma)} + d\gamma + 
\|\btau-\bsigma\|\cdot \psile(\gamma)\\
    &=\ vP_{[\bsigma]} +d\gamma\ =\ v\big(P_{[\bsigma]}g^d\big), 
    \end{align*}
so $v(P_{\times g})_{[\bsigma]}\ =\  v(P_{[\bsigma]}g^d)\ =\ vP_{[\bsigma]}+d\gamma$
as claimed. Now, take $\bsigma_0$ such that $vP_{[\bsigma_0]}=vP$ and $\|\bsigma_0\|=\Pnu(P)$, and take $\bsigma\ge \bsigma_0$ as in the claim.  Then
$$ v_P(\gamma)\ \le\  vP_{[\bsigma]} + d\gamma 
\ \le\  v(P)+d\gamma + \|\bsigma-\bsigma_0\|\cdot \abs{\psile(\gamma)},$$
by the claim, which gives the desired upper bound. 
\end{proof}

\subsection*{Notes and comments} To our knowledge the valuation 
$\psile$ induced on the value group has not been used before in this generality.

\section{Asymptotic Couples}\label{sec:ascouples}

\noindent
{\em In this section $\Gamma$ is an arbitrary ordered abelian group, unless
we specify it as the value group of our valued differential field $K$. 
We let $\alpha$,~$\beta$,~$\gamma$ range over $\Gamma$}. 

\index{asymptotic couple}
\nomenclature[T]{$(\Gamma,\psi)$}{asymptotic couple}

\medskip\noindent
An {\bf asymptotic couple} is a pair $(\Gamma, \psi)$ with $\psi\colon \Gamma^{\ne} \to \Gamma$, such that for all $\alpha, \beta\ne 0$, 
\begin{list}{*}{\setlength\leftmargin{3em}}
\item[(AC1)] $\alpha+\beta\ne 0\ \Longrightarrow\ \psi(\alpha+\beta)\ge \min\big(\psi(\alpha), \psi(\beta)\big)$;
\item[(AC2)] $\psi(k\alpha)=\psi(\alpha)$  for all $k\in \Z^{\ne}$, in particular, $\psi(-\alpha)=\psi(\alpha)$;
\item[(AC3)] $\alpha>0 \ \Longrightarrow\ \alpha + \psi(\alpha) > \psi(\beta)$. 
\end{list}
If in addition for all $\alpha,\beta$, 
\begin{list}{*}{\setlength\leftmargin{3em}}
\item[(HC)] $0<\alpha\leq \beta \Rightarrow \psi(\alpha)\geq\psi(\beta)$,
\end{list}
then $(\Gamma,\psi)$ is said to be of {\bf $H$-type}, or to be
an {\bf $H$-asymptotic couple}. \label{p:H-ascouples}\index{asymptotic couple!H-type@$H$-type}\index{H-asymptotic@$H$-asymptotic!couple}\index{H-type@$H$-type!asymptotic couple}  (We chose the prefix $H$ because $H$-asymptotic 
couples originate from Hardy fields; for more on this, see Chapter~\ref{ch:asymptotic differential fields}.)
 Trivial examples of $H$-asymptotic couples are
obtained by taking any~$\Gamma$ and any constant function
$\psi\colon \Gamma^{\ne} \to \Gamma$. In the next lemma
$\Gamma$ is the value group of our valued differential field $K$.

\begin{lemma} Let the function $\psile\colon \Gamma^{\ne} \to \Gamma$ be as defined in 
the previous section. Then~$(\Gamma, \psile)$ is an asymptotic couple.
\end{lemma}
\begin{proof} Conditions (AC1) and (AC2) correspond to (i) and (ii) of
Lemma~\ref{valv}. As to (AC3), if $\alpha>0$, $\beta\ne 0$, then 
$\alpha + \psile(\alpha) > 0\ge \psile(\beta)$
by (iii) of Lemma~\ref{valv}.
\end{proof}

\noindent
Let $(\Gamma,\psi)$ be an asymptotic couple. 
When convenient we consider $\psi$ as extended 
to a function $\Gamma_{\infty}\to \Gamma_{\infty}$ by 
$\psi(0)=\psi(\infty):=\infty$. Then $\psi(\alpha+ \beta)\ge \min\big(\psi(\alpha), \psi(\beta)\big)$ holds for all 
$\alpha$,~$\beta$, and 
$\psi\colon \Gamma \to \Gamma_{\infty}$ is a (non-surjective) valuation on the abelian group~$\Gamma$, as defined in Section~\ref{sec:valued abelian gps}. In particular,
$$ \psi(\alpha) < \psi(\beta)\ \Longrightarrow\
\psi(\alpha + \beta)=\psi(\alpha).$$ 
Each $\alpha$ yields subgroups
$$\big\{\gamma:\  \psi(\gamma)\ge \alpha\big\}, \quad \big\{\gamma:\ \psi(\gamma)>\alpha\big\}$$
of $\Gamma$. 
If $(\Gamma,\psi)$ is of $H$-type, then these two subgroups are
convex in $\Gamma$, and for all~$\alpha$,~$\beta$, if $\psi(\alpha) > \psi(\beta)$, then $[\alpha]< [\beta]$, by (AC2) and (HC). In particular,
if~$(\Gamma,\psi)$ is of $H$-type, then $\psi\colon \Gamma \to \psi(\Gamma)$ 
is a convex valuation on the ordered abelian group~$\Gamma$, as defined in Section~\ref{sec:oag}.  

\medskip\noindent
We define the {\bf shift} $\psi+\alpha\colon \Gamma^{\ne} \to \Gamma$ by
$$(\psi+\alpha)(\gamma)\ :=\ \psi(\gamma) + \alpha.$$
Then $(\Gamma, \psi+ \alpha)$ is also an asymptotic couple (a {\bf shift} of $(\Gamma, \psi)$), and if $(\Gamma,\psi)$ is
of $H$-type, so is $(\Gamma,\psi + \alpha)$. The next lemma gives two other useful
constructions of asymptotic couples. 

\index{asymptotic couple!shift}

\begin{lemma}\label{lem:minpsipsi2} 
If $(\Gamma, \psi_1)$ and $(\Gamma, \psi_2)$ are asymptotic couples, then
so is $(\Gamma, \psi)$ with $\psi\colon \Gamma^{\ne}\to \Gamma$ given by
$\psi(\gamma) := \min\!\big(\psi_1(\gamma), \psi_2(\gamma)\big)$.
If $(\Gamma, \psi)$ is an asymptotic couple, $\gamma_0\in \Gamma$, 
$0\le \gamma_0 < \big\{\gamma+\psi(\gamma):\ \gamma>0\big\}$, then so is $(\Gamma, \psi_0)$ with  
$\psi_0\colon \Gamma^{\ne}\to \Gamma$ given by
$$\psi_0(\gamma)= \psi(\gamma)\ \text{ if $\psi(\gamma)\le 0$},\qquad \psi_0(\gamma)=\gamma_0\
\text{ if  $\psi(\gamma)>0$.}$$
\end{lemma}

\noindent
The lemma goes through with \textit{$H$-asymptotic}\/ in place of \textit{asymptotic.}\/ The proofs are routine verifications. Less routine is the following. 

\begin{lemma}\label{ro} Let $(\Gamma, \psi)$ be an asymptotic couple and 
$\alpha, \beta\ne 0$. Then 
$$n\big(\psi(\beta)-\psi(\alpha)\big)\  <\  |\alpha| \quad \text{ for all $n$.}$$
\end{lemma}
\begin{proof} Replacing $\alpha$ if necessary by $-\alpha$ we can assume $\alpha>0$ and likewise we arrange $\beta>0$. We can also assume 
$\psi(\beta) > \psi(\alpha)$. Next, replacing
$\psi$ by its shift $\psi - \psi(\alpha)$ we arrange $\psi(\alpha)= 0$. Then
$\psi(\beta) >0$, and our job is to show that $n\psi(\beta) < \alpha$ for all $n$.
Note that $\psi(\beta) < \alpha +\psi(\alpha)=\alpha$. Also 
$$\psi(\beta) <  \psi(\beta)+ \psi\big(\psi(\beta)\big),$$ hence
$\psi\big(\psi(\beta)\big)>0$, so $\psi\big(n\psi(\beta)\big) >0$ for all~${n\ge 1}$, and thus
$n\psi(\beta) \ne \alpha$ for all~${n\ge 1}~$. Assume towards a contradiction that
$n\psi(\beta) > \alpha$ for some $n\ge 1$. Then for the least such $n$ we have
$n\psi(\beta) = \alpha + \gamma$ with $0 < \gamma < \psi(\beta)$, so
$\psi(\gamma)=\psi\big(n\psi(\beta) -\alpha\big)=0$ since $\psi(\alpha)=0< \psi\big(n\psi(\beta)\big)$.
Therefore, $\psi(\beta) < \gamma+ \psi(\gamma)=\gamma$, a contradiction. 
\end{proof}

\noindent
Let $(\Gamma,\psi)$ be an asymptotic couple.
Using Lemma~\ref{ro}, $\psi$ extends 
 uniquely to a map $(\Q\Gamma)^{\ne} \to \Q\Gamma$, also denoted by $\psi$,
such that $(\Q\Gamma,\psi)$ is an asymptotic couple. Then
$\psi\big((\Q\Gamma)^{\ne}\big)=\psi(\Gamma^{\ne})$ and if $(\Gamma,\psi)$ 
is of $H$-type, so is 
$(\Q\Gamma,\psi)$. 

Suppose $\Gamma$ is divisible, and thus a vector space over $\Q$. 
Now $\psi$ is a valuation on~$\Gamma$, and so if $(\alpha_i)$ is a 
family of elements of $\Gamma^{\ne}$
with $\psi(\alpha_i) \ne \psi(\alpha_j)$ for all distinct indices~$i$,~$j$, then 
$(\alpha_i)$ is a $\Q$-linearly independent family. In particular,  
if $\dim_{\Q} \Gamma$ is finite, then $\psi(\Gamma^{\ne})$ is finite
of size at most $\dim_{\Q} \Gamma$.

\medskip\noindent  
The next two lemmas state other basic facts about asymptotic couples.

\begin{lemma}\label{BasicProperties} Let $(\Gamma,\psi)$ be an 
asymptotic couple. 
\begin{enumerate}
\item[\textup{(i)}] If $\alpha,\beta< \gamma+\psi(\gamma)$ for all $\gamma>0$, then
$\psi(\alpha-\beta)>\min(\alpha,\beta)$.
In particular, if $\alpha,\beta\neq 0$, then
$\psi\bigl(\psi(\alpha)-\psi(\beta)\bigr)>
\min\big(\psi(\alpha),\psi(\beta)\big)$.
\item[\textup{(ii)}] If $\alpha,\beta\neq 0$ and $\alpha\neq \beta$, then 
$\psi(\alpha)-\psi(\beta)=o(\alpha-\beta)$.
\item[\textup{(iii)}] The map $\gamma\mapsto \gamma+ \psi(\gamma)\colon \Gamma^{\ne} \to\Gamma$ is strictly
increasing.
\end{enumerate}
\end{lemma}
\begin{proof}
For (i), let $\alpha <\beta <\gamma+\psi(\gamma)$ for all 
$\gamma>0$. Then $\beta< (\beta-\alpha)+ \psi(\beta-\alpha)$, so
$\psi(\beta-\alpha)>\alpha$, as required.
For (ii),
let $\alpha,\beta\neq 0$ with $\gamma:=\alpha-\beta\neq 0$. We have to show that then
$n\bigl|\psi(\alpha)-\psi(\beta)\bigr|<|\gamma|$ for all $n$.
If $\psi(\gamma)>\psi(\beta)$, then $\psi(\alpha)=\psi(\beta)$ 
by axioms~(AC1) and (AC2). 
Suppose $\psi(\gamma)\leq
\psi(\beta)$. Then by axiom (AC1) again we have $\psi(\gamma) \le \psi(\alpha)$, hence by Lemma~\ref{ro}:
$$n\psi(\gamma)\leq n\psi(\beta)<n\psi(\gamma)+|\gamma|, \qquad
n\psi(\gamma)\leq n\psi(\alpha)<n\psi(\gamma)+|\gamma|.$$
Thus $n\bigl|\psi(\alpha)-\psi(\beta)\bigr|<|\gamma|$ in all cases.
Property (iii) follows easily from (ii). \end{proof}

\noindent
Item (ii) justifies thinking of $\psi(\gamma)$ as a {\em slowly varying\/} 
function of $\gamma\in \Gamma^{\ne}$. From (i) we obtain something needed in Section~\ref{The Shape of the Newton Polynomial}:

\begin{cor}\label{1+o(1)} Suppose $(\Gamma, \psi)$ is of $H$-type and $\alpha>0$ and $\beta$ are
such that $\alpha=\psi(\alpha)\le \beta < \gamma+\psi(\gamma)$ for all
$\gamma>0$. Then $\beta=\alpha+ o(\alpha)$.
\end{cor} 

\begin{lemma}\label{CharAsymptoticCouples}
Let $\Gamma$ be an ordered abelian group and $\psi\colon\Gamma^{\ne}\to\Gamma$.
The pair~$(\Gamma,\psi)$ is an asymptotic couple if and only if
\begin{enumerate}
\item[\textup{(i)}] $\psi(k\gamma)=\psi(\gamma)$ for all $k\in\Z^{\ne}$, 
$\gamma\in\Gamma^{\ne}$,
\item[\textup{(ii)}] $\psi(\alpha-\beta)\geq\min\big(\psi(\alpha),\psi(\beta)\big)$
for all $\alpha,\beta\in\Gamma^{\ne}$, $\alpha\neq \beta$, and
\item[\textup{(iii)}] the map $\gamma\mapsto \gamma+ \psi(\gamma)\colon\Gamma^{>}\to\Gamma$ is strictly increasing.
\end{enumerate}
\end{lemma}
\begin{proof}
If  $(\Gamma,\psi)$ is an
asymptotic couple, then clearly (i)--(iii) hold. Conversely, 
suppose $(\Gamma,\psi)$ satisfies (i)--(iii). To verify (AC3), let $\alpha,\beta\in \Gamma^{>}$; we have to derive that $\psi(\beta)<\alpha+\psi(\alpha)$. 
From $\alpha+\beta>\beta>0$ and (iii) we get $\alpha+\beta+\psi(\alpha+\beta)>
\beta+\psi(\beta)$. Also $\alpha+\beta+\psi(\beta)>\beta+\psi(\beta)$,
hence by (ii)
$$\alpha+\beta+\psi\bigl((\alpha+\beta)-\beta\bigr)>
\beta+\psi(\beta),$$
that is, $\alpha + \beta + \psi(\alpha) > \beta + \psi(\beta)$, so
$\alpha+\psi(\alpha)>\psi(\beta)$ as
required.
\end{proof}

\subsection*{Notes and comments} Rosenlicht 
introduced and studied
asymptotic couples in  \cite{Rosenlicht1, Rosenlicht2, Rosenlicht3}.
Lemma~\ref{ro} is  Theorem~5 in \cite{Rosenlicht2}; the proof here follows~\cite{psiremarks}.

\section{Dominant Part}\label{sec: dominant}

\noindent
{\em In this section we choose for every $P\in K\{Y\}^{\ne}$ an element 
$\frak d_P\in K^\times$ with $\frak d_P \asymp P$, such that
$\frak d_P=\frak d_Q$ whenever $P\sim Q$, $P,Q\in K\{Y\}^{\ne}$.
Let $P\in K\{Y\}^{\ne}$}.

\index{dominant!part}
\index{differential polynomial!dominant!part}
\nomenclature[V]{$D_P$}{dominant part of  $P$}
\nomenclature[V]{$\mathfrak d_P$}{dominant monomial of $P$}

\begin{definition}
We have $\frak d_P^{-1}P\in \mathcal{O}\{Y\}$; we call
the differential polynomial
$$D_P\ :=\ \bar{\frak d_P^{-1}P}\ =\  \sum_{\i} (\bar{P_{\i}/\frak d_P})\, Y^{\i}\ =\ 
\sum_{\bomega} (\bar{P_{[\bomega]}/\frak d_P})\, Y^{[\bomega]}\in \k\{Y\}$$
the {\bf dominant part} of $P$.   
Clearly $D_P$ is nonzero with 
\begin{align*} \deg D_P\ &\le\ \deg P,   \hskip-6em  & \order D_P\ &\le\ \order P,\\ 
\wv(D_P)\ &=\ \Pmu(P), \hskip-6em &\wt(D_P)\ &=\ \Pnu(P).
\end{align*}
If $P$ is homogeneous of degree $d$, respectively isobaric of weight $w$,
so is $D_P$. It will also be convenient to set $\frak d_Q:= 0\in K$ and
$D_Q:= 0\in \k\{Y\}$ for $Q=0\in K\{Y\}$. 
\end{definition}

\noindent
Another choice of $\frak d_P$ multiplies $D_P$
by a factor from $\k^\times$. In fact,
only quantities like $\deg D_P$ and $\val D_P$ that are independent 
of the choice of $\frak d_P$ will matter in this chapter.
Here are a few simple rules:

\begin{samepage}

\begin{lemma} \label{simple rules for D_P, new} 
Let $Q\in K\{Y\}$. Then
\begin{enumerate}
\item[\textup{(i)}] if $P \succ Q$, then $\frak d_{P+Q}=\frak d_P$ and $D_{P+Q} = D_P$;
\item[\textup{(ii)}] given any $a\in K^\times$ we have $D_{aP}=e\,D_P$ for some
$e\in \k^\times$;
\item[\textup{(iii)}] $\frak d_{PQ}=u^{-1}\,\frak d_P\frak d_Q$ and $D_{P Q} = \bar{u}\,D_P D_Q$, 
for some $u\asymp 1$ in $K$.
\end{enumerate}
\end{lemma}

\end{samepage}

\begin{proof}
Item (i) is clear; (ii) follows from (iii). For (iii) we may 
assume $Q\ne 0$. We have
$v(P Q)=v(P)+v(Q)$, so we have $u\asymp 1$ in $K$ with
$\frak d_{PQ}=u^{-1}\frak d_P\frak d_Q$,
hence with $\i$, $\j$, $\l$ ranging over $\N^{1+r}$ with 
$r=\max\{\order P, \order Q\}$:
\begin{align*}
D_{PQ}\ &=\ \sum_{\i} \bar{(P\cdot Q)_{\i}/\frak d_{PQ}}\ Y^{\i}\ =\ 
\bar{u}\sum_{\i} \bar{(P\cdot Q)_{\i}/(\frak d_P \frak d_Q)}\ Y^{\i}\\
             &=\  \bar{u}\sum_{\i} \left(\sum_{\j + \l=\i} \bar{P_{\j}/\frak d_P}\cdot 
           \bar{Q_{\l}/\frak d_Q}\right) Y^{\i}\ 
             =\  \bar{u}D_P D_Q. \qedhere
\end{align*} 
\end{proof}

\begin{lemma} \label{lem:composition and dominant parts}
Let $Q\in \mathcal O\{Y\}$ be such that $\overline{Q}\notin\k$. Then 
$$P(Q)\asymp P,\qquad D_{P(Q)}\in \k^\times\cdot D_P(\overline{Q}).$$
\end{lemma}
\begin{proof}
We have 
$\mathfrak d_P^{-1}P(Q) = \sum_{\i} (P_{\i}/\mathfrak d_P) Q^{\i}\in\mathcal O\{Y\}$, and
$$\overline{\mathfrak d_P^{-1}P(Q)}\ =\ \sum_{\i} (\overline{P_{\i}/\mathfrak d_P})\, \overline{Q}^{\i}\
=\ D_P(\overline{Q})\ \neq\  0$$
by Lemma~\ref{lem:comp3}, from which the claims follow.
\end{proof}

\noindent
We have $P=\frak d_P D + R$ with $D\in \mathcal{O}\{Y\}$, $\bar{D}=D_P$, and $R\in K\{Y\}$,
$R\prec P$. We also observe the following:

\begin{lemma}\label{lem:clean, new}
If $P=Q+R$, $Q,R\in K\{Y\}$, $R\prec Q$, 
then $P\sim Q$, $D_P=D_Q$. If $P= \frak d_P D + R$, $D, R\in K\{Y\}$, 
$R\prec P$, then $D\in \mathcal{O}\{Y\}$ and $\bar{D}=D_P$.
\end{lemma}

\noindent
Recall that $\val(P)$ is the least $d\in \N$ such that $P_d\ne 0$. Thus
$$\val D_P >0\ \Longleftrightarrow\ D_P(0)=0\ \Longleftrightarrow\ P(0) \prec P.$$
Note that the quantity $\val D_P$ does not depend on the choice 
of $\frak{d}_P$.
The main part of Lemma~\ref{dn1,new} is the implication $a\prec 1 \Rightarrow \deg D_{P_{\times a}}\le \val D_{P}$.

\begin{lemma}\label{dn1,new}  Let $a\in K$, $a\preceq 1$. Then: 
\begin{enumerate}
\item[\textup{(i)}] $D_{P_{+a}}\in \k^\times \cdot(D_P)_{+\bar{a}}$, and thus $\deg D_{P_{+a}}= \deg D_P$;
\item[\textup{(ii)}] if $a\asymp 1$, then $D_{P_{\times a}}\in \k^\times\cdot (D_P)_{\times\bar{a}}$, $\val D_{P_{\times a}}=\val D_{P}$,  
$\deg D_{P_{\times a}}=\deg D_P$;
\item[\textup{(iii)}] if $a\prec 1$, then $ D_{P_{+a}}=D_P$, 
$\deg D_{P_{\times a}}\le \val D_{P}$.
\end{enumerate}   
\end{lemma}
\begin{proof} Taking $Q=Y+a$ in 
Lemma~\ref{lem:composition and dominant parts} gives (i).
Taking $Q=aY$ in that lemma gives (ii). If $a\prec 1$, then
$P_{+a}\sim P$, which gives the
equality in (iii). For the inequality in (iii) we can 
assume $P\asymp 1$. Then
$d:=\val(D_P)$ gives  
\begin{align*} P\ &=\ \sum_{i<d} P_i + P_d + \sum_{i> d}P_i, \quad P_i\prec 1 \text{ for $i<d$,} 
\quad P_d \asymp 1, \quad P_i \preceq 1 \text{ for $i > d$,}\\   
P_{\times a}\ &=\ \sum_{i<d}\, (P_i)_{\times a} + (P_d)_{\times a} + 
\sum_{i> d}\,(P_i)_{\times a}.
\end{align*}
Assume $0\ne a\prec 1$. Then we obtain from Corollary~\ref{vP-Lemma} that for $i>d$,
$$v\big((P_d)_{\times a}\big)\ =\ dva+o(va)\ <\  v\big((P_i)_{\times a}\big)\ =\ v(P_i) + iva+ o(va).$$ 
The inequality in (iii) now follows, since $(P_{\times a})_i=(P_i)_{\times a}$ for all $i$. \end{proof}

\noindent
Set $\dv(P):=\val D_P$ and $\dd(P):=\deg D_P$. We call $\dv(P)$ 
the {\bf dominant multiplicity}
of $P$ at $0$, and $\dd(P)$ the {\bf dominant degree} of $P$.

\index{differential polynomial!dominant!multiplicity}
\index{differential polynomial!dominant!degree}
\index{dominant!multiplicity}
\index{dominant!degree}
\index{degree!dominant}
\nomenclature[V]{$\dv(P)$}{dominant multiplicity of $P$ at $0$}
\nomenclature[V]{$\dd(P)$}{dominant degree of $P$}

\begin{cor}\label{aftra} Let $a,b\in K$, $g\in K^\times$ 
be such that $a-b\preceq g$. Then 
$$\ddeg P_{+a,\times g}\ =\ \ddeg P_{+b,\times g}.$$
\end{cor}
\begin{proof} Just note that $d:= (b-a)/g\in \mathcal{O}$ and
$P_{+b,\times g}=P_{+a, \times g, +d}$.
\end{proof}

\begin{cor}\label{dn1cor,new} Let $\fm,\fn\in K^\times$. Then
$\val{P} =\val(P_{\times \fm}) \le \dd P_{\times\fm}$ and 
$$\fm \prec \fn\ \Longrightarrow\ \dv P_{\times\mathfrak m}\ \leq\
\dd P_{\times\fm}\ \leq\ \dv P_{\times\mathfrak n}\  
\le\ \dd P_{\times\fn}.$$ 
\end{cor}

\begin{lemma}\label{dndn,new} Let $a\in K^\times$, and suppose 
$\dd P_{\times a}=k$.
Then $$0\ \le\ v(P_k)-v(P)\ \le\ (1+\deg P)\,|va|. $$
\end{lemma}
\begin{proof} This is clear when $vP=vP_k$, so assume $vP < vP_k$, and take
$l$ with $vP=vP_l$. Then $l\ne k$. Recall that
$P_{\times a,k}=P_{k,\times a}$, and $P_{\times a,l}=P_{l,\times a}$. Set $\alpha:= va$, and consider first the case that $\psile(\alpha) \ge 0$. Then by Lemma~\ref{vplemma},  
$$v(P_{\times a,k})=vP_k + k\alpha, \quad  v(P_{\times a,l})=vP_l + l\alpha= vP+l\alpha,$$
so $v(P_k)+ k\alpha \le vP+l\alpha$, and thus $0 < vP_k-vP \le |(l-k)\alpha|$, which gives what we want. Next, assume $\psile(\alpha) < 0$. Then by Proposition~\ref{ps1},     
$$v(P_{\times a,k})=vP_k + k\alpha + o(\alpha), \quad  v(P_{\times a,l})=vP_l + l\alpha +o(\alpha)= vP+l\alpha + o(\alpha),$$
and we can argue as before. 
\end{proof}

\subsection*{The dominant degree} {\em As before, $P\in K\{Y\}^{\ne}$. 
Let 
$\fm$ and $\fn$ range over $K^\times$, and let 
$\E\subseteq K^\times$ be nonempty such that $\fm\in \E$ 
whenever $\fm \preceq \fn \in \E$}. 
The {\bf dominant degree} of $P$ on $\E$, $\dd_{\E}P$, is the natural
number given by 
$$ \dd_{\E}P\ :=\ \max\big\{\!\dd(P_{\times\fm}):\ \fm\in \E\big\}.$$
By Corollary~\ref{dn1cor,new} we have $\val P\le \dd_{\E}P$.
If $Q\in K\{Y\}$, $Q\ne 0$, then clearly 
$$\dd_{\E}PQ\ =\ \dd_{\E}P + \dd_{\E}Q.$$ 
For $\phi\in K^\times$ we have
$\dd_{\E} P_{\times \phi} = \dd_{\phi\E}P$. 

\index{differential polynomial!dominant!degree!on $\mathcal E$}
\index{dominant!degree!on $\mathcal E$}
\index{degree!dominant!on $\mathcal E$}
\nomenclature[V]{$\dd_{\mathcal E}(P)$}{dominant degree of $P$ on $\mathcal E$}

\begin{lemma} 
Suppose $v(\mathcal E)$ does not have a smallest element. Then 
$$\dd_{\mathcal E} P = \max\big\{\!\dv(P_{\times\mathfrak m}):\ \mathfrak m\in\mathcal E\big\}.$$
In particular, there exists $\frak m\in\mathcal E$ such that for all 
$\frak n\in \mathcal E$ with $\fm\prec\frak n$, the differential polynomial  $D_{P_{\times\mathfrak n}}$ is homogeneous of degree $\dd_{\mathcal E} P$.
\end{lemma}
\begin{proof} We have
$\dv P_{\times\mathfrak m}\leq \dd P_{\times\mathfrak m}\leq \dd_{\mathcal E}P$ for $\mathfrak m\in\mathcal E$. Take $\frak m\in\mathcal E$ with $\dd P_{\times\mathfrak m}=\dd_{\mathcal E} P$, and let $\mathfrak n\in\mathcal E$ be such that $\mathfrak m\prec\mathfrak n$. Then by Corollary~\ref{dn1cor,new} we have $\dd_{\mathcal E} P=\dd P_{\times\mathfrak m} =\dv P_{\times\mathfrak n}=\dd P_{\times \fn}$. This yields the claim.
\end{proof}

\begin{lemma}\label{dn4,new} Suppose that $f\in \E$. Then 
$\dd_{\E}P_{+f}=\dd_{\E} P$.
\end{lemma}
\begin{proof} It is enough that for $\fm\in \E$ with
$\fm\succeq f$ we have $\dd P_{+f,\times\fm}=\dd  P_{\times\fm}$, and this is a special case of Corollary~\ref{aftra}. 
%We have $P_{+f,\times\fm}=P_{\times \fm,+a}$ 
%with $a=f/\fm \preceq 1$.
%By Lemma~\ref{dn1,new}, (1), this gives the 
%desired equality of degrees. 
\end{proof}

\begin{cor}\label{cor:NP homog deg 1}
Suppose $\ddeg_{\E} P=1$. Let $y\in \E\cup\{0\}$ be a zero of $P$ and $\fm\in\E$.
Then 
$\val P_{+y,\times\fm}\ =\ \dv P_{+y,\times \fm}\ =\ \ddeg P_{+y,\times \fm}\ =\ 1.$
\end{cor}
\begin{proof}
Since $(P_{+y})_0=P(y)=0$, we have
$$1 \leq \val P_{+y} = \val P_{+y,\times\fm} \leq \dv P_{+y,\times\fm} \leq \ddeg P_{+y,\times\fm} \leq \ddeg_{\E} P_{+y} = 1,$$
using Lemma~\ref{dn4,new} for the last step.
\end{proof}

\noindent
Let $\E'\subseteq \E$ be nonempty such that $\fm\in \E'$ whenever
$\fm\preceq \fn\in \E'$. Then for $f\in \E$ we have 
$\dd_{\E'}P_{+f}\le \dd_{\E} P$ by Lemma~\ref{dn4,new}.
For $\gamma\in \Gamma$ and $\E=\{\fn: v\fn\ge \gamma\}$, we set 
$\dd_{\geq\gamma} P:= \dd_{\E} P$, so if $v\fm=\gamma$, then
$\dd_{\geq\gamma} P = \dd P_{\times \fm}$. 

\begin{cor}\label{dn7a,new} Let $a, b\in K$ and $\alpha, \beta\in \Gamma$ 
be such that $v(b - a)\ge \alpha$ as well as $\beta \ge \alpha$. 
Then $\dd_{\geq\beta} P_{+b}\ \le\ \dd_{\geq\alpha}P_{+a}$. 
\end{cor}
\begin{proof} Since $P_{+b}=P_{+a, +(b-a)}$, we have 
$\dd_{\geq\alpha}P_{+b}=\dd_{\geq\alpha}P_{+a}$ by Lem\-ma~\ref{dn4,new}. 
It remains to note that $\dd_{\geq\beta} P_{+b}\le \dd_{\geq\alpha}P_{+b}$. 
\end{proof}

%\noindent
%Here is how $\dd_{\E}$ behaves under extensions of valued differential fields.
%Let $L$ be a valued differential field extension of $K$ with small 
%derivation, and let $\mathfrak u$, $\mathfrak v$ range over~$L^\times$. 
%\marginpar{New lemma.}

%\begin{lemma}
%Let
%$\E_L := \{ \mathfrak u : \text{$\mathfrak u\preceq \mathfrak m$ for 
%some $\mathfrak m\in\mathcal E$} \}$.  
%Then $\E_L\cap K=\E$, and
%$\mathfrak u\in\mathcal E_L$ whenever 
%$\mathfrak u\preceq \mathfrak v\in\mathcal E_L$. Moreover,
%$\dd_{\mathcal E} P = \dd_{\mathcal E_L} P$.
%\end{lemma}
%\begin{proof}
%The first statement is obvious.
%Since $\mathcal E \subseteq\mathcal E_L$, we clearly have 
%$\dd_{\mathcal E} P \leq \dd_{\mathcal E_L} P$.
%To show the reverse inequality, let $\mathfrak u\in\mathcal E_L$, 
%and take $\mathfrak m\in\mathcal E$ with $\mathfrak u\preceq\mathfrak m$. 
%Now if $\mathfrak u\asymp\mathfrak m$, then $$\dd P_{\times \mathfrak u}= 
%\dd\, (P_{\times \mathfrak u})_{\times \mathfrak m/\mathfrak u}=
%\dd P_{\times \mathfrak m}$$ by Lemma~\ref{dn1,new},~(3), and 
%if $\mathfrak u\prec\mathfrak m$ then
%$\dd P_{\times \mathfrak u} < \dd P_{\times \mathfrak m}$ by 
%Corollary~\ref{dn1cor,new}.
%In both cases $\dd P_{\times \mathfrak u} \leq  \dd_{\mathcal E} P$.
%Since this holds for all $\mathfrak u\in\mathcal E_L$, we have 
%$\dd_{\mathcal E_L} P \leq \dd_{\mathcal E} P$ as required.
%\end{proof}

\noindent
Suppose $\Gamma\ne \{0\}$ and $\Gamma^{>}$ has no least element.
Then $\E=\{\fn: \fn \prec \fm\}$ is nonempty, and we set 
$\dd_{\prec \fm}P:= \dd_{\E}P$. Also 
$v\big(\{\fn: \fn \prec \fm\}\big)$ has no least element, hence 
$$ \dd_{\prec \fm}P\ =\ \max\big\{\!\dd P_{\times\fn}:\ \fn\prec \fm\big\}\ =\ \max\big\{\!\dv P_{\times\fn}:\ \fn\prec \fm\big\}.$$

\section{The Equalizer Theorem}\label{equalize}

\noindent
We stated this theorem in the introduction 
to this chapter. {\em In this section 
$a$, $b$, $f$, $g$, $y$ range over $K$, and
$\alpha$, $\beta$, $\gamma$ over $\Gamma$. We fix $P\in K\{Y\}^{\ne}$}.

\subsection*{Lemmas on equalizing} These lemmas are technical facts to be
used in the proof of the Equalizer Theorem. 

\begin{lemma}\label{dn5} Suppose $va=\alpha <0$, $\dd P_{\times a}=m$, $\wt(P)=w$, $m<n$. Then
\begin{align*}  \psile(\alpha)= 0\ &\Longrightarrow\ v(P_n)-v(P_m)\ >\  
(n-m)|\alpha|,\\
\psile(\alpha) < 0\ &\Longrightarrow\  v(P_n)-v(P_m)\ >\ (n-m)|\alpha|+ 2w\cdot \psile(\alpha).
\end{align*}
\end{lemma}
\begin{proof} Since $0\ne P_{\times a, m}=P_{m,\times a}$, we have $P_m\ne 0$, so
$v(P_n)- v(P_m)\in \Gamma_{\infty}$ is certainly defined. From $m<n$ we get 
$v(P_{m, \times a})< v(P_{n, \times a})$, so
\begin{align*} \psile(\alpha)= 0\ &\Longrightarrow\  v(P_m) + m\alpha\ <\  v(P_n)+ n\alpha,\\
\psile(\alpha) < 0\ &\Longrightarrow\  v(P_m) + m\alpha + w\cdot \psile(\alpha)\ <\  v(P_n)+ n\alpha -w\cdot \psile(\alpha),
\end{align*}
by Lemma~\ref{vplemma} and Proposition~\ref{ps1}, respectively.  
\end{proof}

\begin{lemma}\label{dn6} Suppose $va=\alpha < 0$, $vb=\beta \ge \alpha-(1/d)\alpha$ for some $d\ge 1$, and $\dd_{\geq\alpha}P = m$. Then $P_{+a,m}\ \sim\ P_{+(a+b),m}$.
\end{lemma}
\begin{proof} Set $w:= \wt P$ and $Q:= P_{+a}$, so $P_{+(a+b),m} = Q_{+b,m}$. We have to show $Q_m \sim Q_{+b,m}$ (which includes $Q_m\ne 0$, $Q_{+b,m}\ne 0$). By Lemma~\ref{dn4,new} the assumption $\dd_{\geq\alpha}P = m$ gives $\dd_{\geq\alpha}Q =m$, so $\dd Q_{\times a}=m$.
Define $\delta\in \Gamma$ by $\delta:= 0$ if $\psile(\alpha)\ge 0$, and
$\delta:= 2w \cdot \psile(\alpha)$ if $\psile(\alpha) < 0$; in  any case, $\delta=o(\alpha)$ by Lemma~\ref{valv}. By
Lemma~\ref{dn5} applied to $Q$ we have $Q_m\ne 0$ and  
\begin{equation}\label{eq:Qn vs Qm}
v(Q_n)-v(Q_m)\ > \ (n-m)|\alpha| +\delta\ \text{ for all $n>m$.} 
\end{equation}
We can assume $b\ne 0$. We have $Q_{+b}=\sum_n Q_{n,+b}$ and by Taylor expansion, 
$$ Q_{n,+b}\ =\ \sum_{|\i|\le n} Q_{n, (\i)}(Y)\cdot b^{\i},$$ 
where each $Q_{n,(\i)}$ is homogeneous of degree $n-|\i|$. It follows that
\begin{align*} Q_{n,+b,m}\ &=\ 0\ \text{ for $n<m$,} \quad \text{ and }\  Q_{m, +b,m}=Q_m,\ \text{ which gives}\\
Q_{+b,m}\ &=\ Q_m + \sum_{n>m} Q_{n,+b,m},\quad
 Q_{n,+b,m}\ =\ \sum_{|\i|=n-m} Q_{n,(\i)}\cdot b^{\i}\  \text{ for $n>m$.}
\end{align*}
Hence it is enough to get $v\big(Q_{n,(\i)}b^{\i}\big) > v(Q_m)$ for all $n > m$ and $|\i|=n-m$. For such $n$ and $\i$ with  $Q_{n,(\i)}\ne 0$ we have
$\|\i\|\le w$, so  $v(b^{\i})\ge (n-m)\beta + \epsilon$,
where $\epsilon:= 0\in \Gamma$ if $\psile(\beta)\ge 0$, and $\epsilon:= w\cdot \psile(\beta)$ if $\psile(\beta) < 0$; in any case, $\epsilon= o(\beta)$ if $\beta\ne 0$. Thus it suffices to get $v(Q_n) + (n-m)\beta +\epsilon> v(Q_m)$ for all $n>m$, that is, 
$$v(Q_n)-v(Q_m)\ >\  (n-m)(-\beta) - \epsilon\ \text{ for all $n>m$.}$$ 
It remains to use \eqref{eq:Qn vs Qm} and $(n-m)|\alpha| +\delta > (n-m)(-\beta) - \epsilon$ for $n>m$. 
\end{proof}

\begin{lemma}\label{eqsuffa} Let
$P, Q\in K\{Y\}^{\ne}$ be homogeneous of degrees~$d$,~$e$, respectively, with $d > e$, and assume $(d-e)\Gamma=\Gamma$. Set $\alpha:= (vQ - vP)/(d-e)$ and 
$\beta :=  \big(v_Q(\alpha)-v_P(\alpha)\big)\big/(d-e)$. Then
\begin{enumerate}
\item[\textup{(i)}] $\psile(\alpha) \ge 0\ \Longrightarrow\ v_P(\alpha)=v_Q(\alpha)$;
\item[\textup{(ii)}] $\alpha\ne 0\ \Longrightarrow\ \psile(\beta) > \psile(\alpha)$;
\item[\textup{(iii)}]  $\psile(\alpha) < 0,\ \psile(\beta) <0\ \Longrightarrow\
\psile(\beta) = o\big(\psile(\alpha)\big)$. 
\end{enumerate}
\end{lemma}
\begin{proof} If $\psile(\alpha) \ge 0$, then $v_P(\alpha)=vP+ d\alpha$ and
$v_Q(\alpha)=vQ + e\alpha$, by Lemma~\ref{vplemma},  
so $v_P(\alpha)=v_Q(\alpha)$. This proves (i). For (ii) and (iii), assume $\alpha\ne 0$. 
Then  $\psile(\alpha) = o(\alpha)$ by 
Lemma~\ref{valv}(iii).
If $\beta=0$, we are done, so assume $\beta\ne 0$. Then $\psile(\alpha) < 0$ by~(i). If $\psile(\beta) \ge 0$, then $\psile(\beta) > \psile(\alpha)$, and we are done. So assume $\psile(\beta) < 0$. Set $w:= \max(\wt P, \wt Q)$. By Proposition~\ref{ps1}, 
\begin{align*} v_Q(\alpha)\ &=\ vQ + e\alpha + \epsilon, \qquad |\epsilon| \le w|\psile(\alpha)|,\\
v_P(\alpha)\ &=\ vP + d\alpha + \delta, \qquad |\delta| \le w|\psile(\alpha)|,
\end{align*}
so $|\beta| = |\epsilon - \delta|/(d-e) \le  2w|\psile(\alpha)|/(d-e)$. 
From
$\psile(\beta) < 0$ we get $\psile(\beta) = o(\beta)$, so $\psile(\beta)=o\big(\psile(\alpha)\big)$. In view of $\psile(\alpha) < 0$, this also gives $\psile(\beta) > \psile(\alpha)$.
\end{proof}

\noindent
In the next lemma we consider $\Gamma$ as equipped with its valuation 
$\psile$.

\begin{lemma}\label{eqpseudo} Let $P, Q\in K\{Y\}^{\ne}$ be homogeneous of degrees~$d, e$, respectively, with $d > e$. Let
$(\alpha_{\rho})_{\rho< \nu}$ be a pc-sequence in $\Gamma$ such that 
$\alpha_{\rho}\leadsto \alpha$, with $\gamma_{\rho}:=\psile(\alpha_{\rho +1}-\alpha_{\rho}) <0$  for all $\rho$, and $\gamma_{\sigma} = o(\gamma_{\rho})$ whenever
$\rho < \sigma < \nu$. Then we have $v_P(\alpha_{\rho})-v_Q(\alpha_{\rho}) \leadsto
v_P(\alpha)-v_Q(\alpha)$. 
\end{lemma} 
\begin{proof} Set $i:= v_P-v_Q$. We have to show that 
$i(\alpha_{\rho}) \leadsto i(\alpha)$. We can arrange that
$\psile(\alpha-\alpha_{\rho}) =\gamma_{\rho}$ for all $\rho$. 
Take $a_{\rho}\in K$ with
$va_{\rho} = \alpha_{\rho}$. By Proposition~\ref{ps1},
\begin{align*} v_P(\alpha)\ &=\ v_P\big(\alpha_{\rho} +(\alpha-\alpha_{\rho})\big)\ =\ 
v_{P_{\times a_{\rho}}}(\alpha-\alpha_{\rho})\\
&=\ v(P_{\times a_{\rho}}) + d(\alpha-\alpha_{\rho}) + O(\gamma_{\rho})\ =\ v_P(\alpha_{\rho}) + d(\alpha-\alpha_{\rho}) + O(\gamma_{\rho}),
\end{align*}
and likewise with $v_Q(\alpha)$. Hence $i(\alpha)=i(\alpha_{\rho}) + (d-e)(\alpha-\alpha_{\rho}) + O(\gamma_{\rho})$. Then by Lemma~\ref{BasicProperties}(ii) and (AC2) we have
$$ \psile \big(i(\alpha)-i(\alpha_{\rho})\big) -\psile(\alpha-\alpha_{\rho})\ =\ o(\gamma_{\rho}),$$ 
that is, $\psile\big(i(\alpha)-i(\alpha_{\rho})\big)= \gamma_{\rho} + o(\gamma_{\rho})$.
It follows that $\psile\big(i(\alpha)-i(\alpha_{\rho})\big)$ as a function of $\rho$ is
strictly increasing, so $i(\alpha_{\rho}) \leadsto i(\alpha)$. 
\end{proof}    

\subsection*{Proof of the Equalizer Theorem} Let $P,Q\in K\{Y\}^{\ne}$ be homogeneous
of degrees $d$ and $e$ with $d>e$, and assume 
$(d-e)\Gamma=\Gamma$. We wish to find an equalizer for $P$, $Q$, \index{equalizer} that
is, an $\alpha$ such that $v_P(\alpha)=v_Q(\alpha)$. (Since $v_P-v_Q$ is strictly increasing, there is at most one equalizer for~$P$,~$Q$.) If $vP=vQ$, then 
$\alpha=0$ works. Assume
$vP\ne vQ$. Since for $\alpha\ne 0$ we have 
$v_P(\alpha)=vP + d\alpha + o(\alpha)$ and 
$v_Q(\alpha)=vQ + e\alpha + o(\alpha)$, we expect the $\alpha$ such that 
$vP + d\alpha=vQ+e\alpha$, that is, $\alpha=(vQ-vP)/(d-e)$, to
be a good approximation
to an equalizer. This leads to the following approximation scheme:
Set $\alpha_0=0$ and take this as the initial term of
a sequence $(\alpha_{\rho})_{\rho< \nu}$
in $\Gamma$ indexed by the ordinals $\rho$ less than a certain ordinal $\nu>0$,
such that the following conditions hold:
\begin{enumerate}
\item $v_P(\alpha_{\rho})\ne v_Q(\alpha_{\rho})$ for $\rho < \nu$;
\item $\alpha_{\rho+1}\ =\ \alpha_{\rho} + \big(v_Q(\alpha_{\rho}) - v_P(\alpha_{\rho})\big)/(d-e)$ when $\rho+1< \nu$;
\item $\psile\big(v_Q(\alpha_{\sigma})- v_P(\alpha_{\sigma})\big)\ >\ \psile\big(v_Q(\alpha_{\rho}) - v_P(\alpha_{\rho})\big)$ for
$\rho < \sigma < \nu$;
\item $\psile(\alpha_{\sigma}-\alpha_{\rho})\  =\  \psile\big(v_Q(\alpha_{\rho})-v_P(\alpha_{\rho})\big)$ for $\rho < \sigma < \nu$.
\end{enumerate}
Note that by (3) we have $\alpha_\rho\ne \alpha_{\sigma}$ whenever $\rho < \sigma < \nu$.
We also pick for each~$\rho< \nu $ an element $a_{\rho}\in K$ such that $va_{\rho}=\alpha_{\rho}$. 
Consider first the case that $\nu$ is a successor ordinal, $\nu=\mu +1$.
Then we set $$\alpha_{\nu}:= \alpha_{\mu} + \big(v_Q(\alpha_{\mu}) - v_P(\alpha_{\mu})\big)/(d-e).$$ If $v_P(\alpha_{\nu})=v_Q(\alpha_{\nu})$,
we are done. Assume $v_P(\alpha_{\nu})\ne v_Q(\alpha_{\nu})$. Then (1)--(4) continue to hold with $\nu$ replaced by $\nu+1$. This is clear for (1) and (2); for (3), apply Lemma~\ref{eqsuffa} to $P_{\times a_{\mu}}$ and
$Q_{\times a_{\mu}}$ in the role of $P$ and $Q$, 
so that $\alpha=\alpha_{\nu}-\alpha_{\mu}$. Now~(4) with $\nu +1$ instead of $\nu$ follows easily. 

\medskip\noindent
Next, consider the case that $\nu>0$ is a limit ordinal. Then by (3) and (4) 
we have a pc-sequence $(\alpha_{\rho})$ in $\Gamma$ (with respect to the valuation $\psile$ on $\Gamma$). Set 
$$\gamma_{\rho}\ :=\  \psile(\alpha_{\rho+1}-\alpha_{\rho})\ =\ \psile\big(v_Q(\alpha_{\rho})-v_P(\alpha_{\rho})\big).$$
Then $\gamma_{\rho} < 0$ for all $\rho$: if this would fail for a certain $\rho < \nu$, then Lemma~\ref{eqsuffa} applied to 
$P_{\times a_{\rho}}$ and
$Q_{\times a_{\rho}}$ in the role of $P$ and $Q$ (which corresponds to $\alpha=\alpha_{\rho+1}-\alpha_{\rho}$) gives $v_P(\alpha_{\rho+1})= v_Q(\alpha_{\rho+1})$, contradicting (1).
Likewise, (3) and Lemma~\ref{eqsuffa} yield: 
$$\rho < \sigma < \nu\ \Longrightarrow\ \gamma_{\rho}\ <\ \gamma_{\sigma}\ <\  0, \quad \gamma_{\sigma}\ =\ o(\gamma_{\rho}).$$
We shall construct a pseudolimit of $(\alpha_{\rho})$ in $\Gamma$, but before we
do this, let us assume that $\alpha_{\nu}\in \Gamma$ is such a pseudolimit. If 
$v_P(\alpha_{\nu})=v_Q(\alpha_{\nu})$, we are done, so assume 
$v_P(\alpha_{\nu})\ne v_Q(\alpha_{\nu})$. We claim that then (1)--(4) holds with $\nu+1$ in place of~$\nu$. This claim is clearly valid for (1) and (2). As to (4), let
$\rho < \nu$ be given, and take $\sigma$ with $\rho < \sigma < \nu$ so 
large that $\psile(\alpha_{\nu}-\alpha_{\sigma})= \gamma_{\sigma}$. Then 
$\gamma_{\sigma}>\gamma_{\rho}$ gives
$$ \psile(\alpha_{\nu}-\alpha_{\rho})\ =\
 \psile\big( (\alpha_{\nu}-\alpha_{\sigma})+(\alpha_{\sigma}-\alpha_{\rho})\big)\ =\ \gamma_{\rho}.$$
A similar argument using Lemma~\ref{eqpseudo} gives (3) with $\nu +1$ 
in place of $\nu$.

\medskip\noindent
We now turn to the construction of a pseudolimit of $(\alpha_{\rho})$. Note:
$$\rho < \sigma < \nu\ \Longrightarrow\ v(a_{\sigma}^\dagger -a_{\rho}^\dagger)\ =\ \psile(\alpha_{\sigma}-\alpha_{\rho})\ =\ \gamma_{\rho}.$$
Let $R:= \Ric(P)$ and $S:= \Ric(Q)$. By Corollary~\ref{dn7a,new},
\[\rho < \sigma < \nu\ \Longrightarrow\ 
\begin{cases}
\displaystyle \dd_{\geq\gamma_{\rho}}R_{+a_{\rho}^\dagger}\ \ge\ \dd_{\geq\gamma_{\sigma}}R_{+a_{\sigma}^\dagger},\\[1em]
\displaystyle\dd_{\geq\gamma_{\rho}}S_{+a_{\rho}^\dagger}\ \ge\ \dd_{\geq\gamma_{\sigma}}S_{+a_{\sigma}^\dagger},
\end{cases}
\]
%$$  \rho < \sigma < \nu\ \Longrightarrow\ \dd_{\geq\gamma_{\rho}}R_{+a_{\rho}^\dagger}\ \ge\ \dd_{\geq\gamma_{\sigma}}R_{+a_{\sigma}^\dagger},\ \ \dd_{\geq\gamma_{\rho}}S_{+a_{\rho}^\dagger}\ \ge\ \dd_{\geq\gamma_{\sigma}}S_{+a_{\sigma}^\dagger},$$
which gives $m$, $n$ and an index $\rho_0$ such that for all $\rho\ge \rho_0$,
$$\dd_{\geq\gamma_{\rho}}R_{+a_{\rho}^\dagger}=m, \qquad \dd_{\geq\gamma_{\rho}}S_{+a_{\rho}^\dagger}=n.$$   
Applying Lemma~\ref{dn6} to $R_{+a_{\rho_0}^\dagger}$ in the role of $P$ and 
$a:=a_{\rho}^\dagger-a_{\rho_0}^\dagger$ and $\alpha:= \gamma_{\rho_0}$, we obtain 
\begin{align*}
\rho > \rho_0,\ vb\ge \gamma_{\rho}\quad &\Longrightarrow & R_{+a_{\rho}^\dagger, m}\ &\sim\ R_{+a_{\rho}^\dagger + b, m},\ \text{ and likewise,}\\
 \rho > \rho_0,\ vb\ge \gamma_{\rho}\quad &\Longrightarrow & S_{+a_{\rho}^\dagger, n}\ &\sim\ S_{+a_{\rho}^\dagger + b, n},\ \text{ and therefore}\\
\rho_0 < \rho < \sigma < \nu\quad &\Longrightarrow &
R_{+a_{\rho}^\dagger, m}\ &\sim\ R_{+a_{\sigma}^\dagger, m},\quad S_{+a_{\rho}^\dagger, n}\ \sim\ 
S_{+a_{\sigma}^\dagger, n}.
\end{align*}
Thus the following element of $\Gamma$ does not depend on the choice of 
$\rho > \rho_0$:
\begin{equation}\label{eq:alpha}
\alpha\ :=\ \frac{1}{d-e}\big(v\big(S_{+a_{\rho}^\dagger, n}\big)- v\big(R_{+a_{\rho}^\dagger, m}\big)\big)\qquad (\rho > \rho_0).
\end{equation} 
We claim that $\alpha_{\rho} \leadsto \alpha$.
First, for $\rho > \rho_0$, 
\begin{align*} \alpha - \alpha_{\rho}\ &=\ \frac{1}{d-e}\bigg[\big(e\alpha_{\rho} + v\big(S_{+a_{\rho}^\dagger, n}\big)\big)- \big(d\alpha_{\rho} + v\big(R_{+a_{\rho}^\dagger, m}\big)\big)\bigg], \\
v_P(\alpha_{\rho})\ &=\ vP_{\times a_{\rho}}\ =\ d\alpha_{\rho} + v\big(R_{+a_{\rho}^\dagger}\big).
\end{align*}
Next, by Lemma~\ref{dndn,new} applied to $R_{+a_{\rho}^\dagger}$ instead of $P$  we have for $\rho > \rho_0$,
\begin{align*} 0\ &\le\  v\big(R_{+a_{\rho}^\dagger,m}\big)-v\big(R_{+a_{\rho}^\dagger}\big)\ \le\ 
\big(1+\wt(P)\big)|\gamma_{\rho}|,\ \text{ so}\\
0\ &\le\ d\alpha_{\rho} +v\big(R_{+a_{\rho}^\dagger,m}\big) - v_P(\alpha_{\rho})\ \le\ \big(1+\wt(P)\big)|\gamma_{\rho}|,\ \text{ and likewise,}\\
 0\ &\le\ e\alpha_{\rho} +v\big(S_{+a_{\rho}^\dagger,n}\big) - v_Q(\alpha_{\rho})\ \le\ 
\big(1+\wt(Q)\big)|\gamma_{\rho}|,\ \text{ and thus}\\
\alpha - \alpha_{\rho}\ &=\  \frac{1}{d-e}\big(v_Q(\alpha_{\rho}) - v_P(\alpha_{\rho})\big) + O(\gamma_{\rho})\ =\ \alpha_{\rho+1}-\alpha_{\rho} + O(\gamma_{\rho}).
\end{align*}
The above together with Lemma~\ref{BasicProperties}(ii) gives
$$\psile(\alpha - \alpha_{\rho}) -\psile(\alpha_{\rho+1}-\alpha_{\rho})\ =\ o(\gamma_{\rho}),$$
that is, $\psile(\alpha - \alpha_{\rho}) =\gamma_{\rho}+o(\gamma_{\rho})$, for
$\rho > \rho_0$. Now, for $\rho_0 < \rho < \sigma < \nu$ we have
$\gamma_{\sigma}=o(\gamma_{\rho})$, so 
$\psile(\alpha - \alpha_{\rho}) < \psile(\alpha - \alpha_{\sigma})$, which
proves our claim that $\alpha_{\rho} \leadsto \alpha$. 

\medskip\noindent
The arguments above show that if there were no equalizer for $P$, $Q$, we could
extend our sequence $(\alpha_{\rho})_{\rho< \nu}$ indefinitely, that is, 
make the ordinal $\nu$ as large as we want, which contradicts 
$|\nu|\le |\Gamma|$.
This concludes the proof of the Equalizer Theorem. \qed

\medskip\noindent
{\em In the rest of this section we keep assuming that $P, Q\in K\{Y\}^{\ne}$ are homogeneous of 
degrees $d$,~$e$ with $d>e$ {\rm{(}}so $d\ge 1${\rm{)}}, and that $(d-e)\Gamma=\Gamma$}. 

\begin{cor} The function
$v_P-v_Q\colon \Gamma \to \Gamma$ is bijective. If $d\Gamma=\Gamma$, then the function 
$v_P\colon\Gamma \to \Gamma$ is bijective.
\end{cor}
\begin{proof} $(v_P-v_Q)(\alpha)=\beta$ is equivalent to
$v_P(\alpha)=v_{bQ}(\alpha)$ where $vb=\beta$. Also, $v_P = v_P-v_Q$ for $Q=1$.
\end{proof}

\noindent
The transfinite part of the proof of the Equalizer Theorem is rather bizarre and can 
probably be eliminated. (Allen Gehret noticed that this transfiniteness is very mild: any sequence 
$(\alpha_{\rho})_{\rho < \nu}$ in the proof above must have
length $\nu \le \omega p$ for some~${p\in \omega}$; this is because 
$\nu=\omega^2$ leads to the equalities \eqref{eq:alpha} for $\rho_0 < \rho < \omega^2$, with $\rho_0=\omega q$, $q\in \omega$, and thus $\alpha_{\omega(q+1)}=\alpha_{\omega(q+2)}$, a contradiction.)   In any case, since the theorem is true, logical considerations suggest that an equalizer is obtainable
with a finite bound on the number of steps used, where the bound
depends only on the orders and degrees of $P$ and $Q$.

{\sloppy

Acting on this suggestion we indicate an algorithm
to compute the equalizer of~$P$,~$Q$. As before we set
$R:= \Ric(P)$ and $S:= \Ric(Q)$, so $\deg R \le \wt P$ and
$\deg S \le \wt Q$. We also assume that for each $\alpha\in \Gamma$ there is given a ``monomial''~$\fm_{\alpha}\in K^\times$ with $v\fm_{\alpha} =\alpha$. 
We define the function
$$ v_{P,Q}\ :\ \Gamma\to \Gamma, \qquad v_{P,Q}(\alpha)\ 
:=\ \alpha+ \frac{v_Q(\alpha)-v_P(\alpha)}{d-e}.$$
Define an {\bf equalizer sequence\/} for $P$,~$Q$ to be a sequence $\alpha_0,\dots, \alpha_N$ in $\Gamma$ with~${N\in \N}$ such that $\alpha_0=0$, and for each $i< N$, either $\alpha_{i+1}=v_{P,Q}(\alpha_i)$, or 
$$\alpha_{i+1}\ =\ \frac{v\big(S_{+\fm_{\alpha}^\dagger, n}\big)-v\big(R_{+\fm_{\alpha}^\dagger, m}\big)}{d-e}$$
for $\alpha:= \alpha_i$ and some $m\le \wt P$ and $n\le \wt Q$ with $R_{+\fm_{\alpha}^\dagger, m}\ne 0$, $S_{+\fm_{\alpha}^\dagger, n}\ne 0$.
%For $m\le \wt P$ and $n\le \wt Q$ we define the function
%$$ v_{P,Q}^{m,n}\ :\ \Gamma\to \Gamma, \qquad 
%v_{P,Q}^{m,n}(\alpha)\ :=\ 
%\frac{v(S_{+\fm_{\alpha}^\dagger, n})-
%v(R_{+\fm_{\alpha}^\dagger, m})}{d-e}.$$
%$\alpha_{i+1}=v_{P,Q}^{m,n}(\alpha_i)$ for some $m\le \wt P$ %and $n\le \wt Q$. 
}

\index{equalizer!sequence}
\index{sequence!equalizer}

\begin{cor}\label{eqalgo} 
Assume $P$ and $Q$ have order $\le r$. 
Then there exists an equalizer sequence 
$\alpha_0,\dots, \alpha_N$ for $P$,~$Q$ such that $\alpha_N$ is an equalizer for $P$,~$Q$ and such that $N\le N(r,d)$ with
$N(r,d)\in \N$ depending only on $r$ and $d$.
\end{cor}
\begin{proof} Recall that there is no strictly descending sequence $(\beta_n)$ in a well-ordered set. Thus, starting with the equalizer
obtained from the proof of the Equalizer Theorem, and going back appropriately in the (possibly transfinite) sequence of that proof, we obtain the reversal of
an equalizer sequence for $P$,~$Q$. A uniform bound as claimed
must exist by model-theoretic compactness.  
\end{proof}

\noindent 
Note that the validity of the algorithm implicit in Corollary~\ref{eqalgo} depends on the Equalizer Theorem, whose proof is
hardly constructive.

\subsection*{Notes and comments} Suppose $P\in K\{Y\}^{\ne}$ is homogeneous of degree~$1$. Then the Equalizer Theorem says that the function $v_P\colon \Gamma\to \Gamma$ is a bijection. 
We do have a very different and more constructive (but longer) proof of this fact if $K$ is $H$-asymptotic as defined in Section~\ref{As-Fields,As-Couples} below. This unpublished proof predates the material above, and also yields the definability of the function $v_P$
in the asymptotic couple $(\Gamma, \psile)$ when $K$ is a
Liouville closed $H$-field as defined in Section~\ref{sec:Liouville closed}.  

We also have a more constructive proof of the Equalizer Theorem for $P$,~$Q$ of order~$\le 1$ and 
$H$-asymptotic $K$.

\section{Evaluation at Pseudocauchy Sequences}\label{sec:pceq}

\noindent
{\em In this section we assume that the induced derivation on the residue field $\k$ is nontrivial.}\/ In addition $L$ denotes a valued differential field extension 
of $K$  with $\der \smallo_L \subseteq \smallo_L$,
and~$(a_\rho)$ is a pc-sequence in $K$
with $a_{\rho} \leadsto a\in L$, and $G(Y)\in L\{Y\}\setminus L$.
We set $\gamma_{\rho}:= v\big(a_{s(\rho)}-a_\rho\big)\in \Gamma_{\infty}$, where 
$s(\rho):= \text{immediate successor of $\rho$}$.

\begin{lemma}\label{pcc1} There is a pc-sequence $(b_\rho)$ in $K$ 
equivalent to 
$(a_\rho)$ such that $\big(G(b_\rho)\big)$ is a
pc-sequence and $G(b_\rho) \leadsto G(a)$. 
\end{lemma}

\begin{proof}
After removing some initial $\rho$'s we can assume
$\gamma_{\rho}=v(a-a_\rho)\in \Gamma$ for all~$\rho$ and 
$\gamma_{\rho'} > \gamma_{\rho}$ whenever $\rho'>\rho$. 
Take $g_\rho\in K$ with
$v(g_\rho)=\gamma_\rho$ and define $u_\rho\in L$ by 
$a_\rho - a = g_\rho u_\rho$, so $u_\rho\asymp 1$. Let $x_{\rho}\in K$ 
be such that $x_{\rho}\asymp 1$ and $u_{\rho} + x_{\rho}\asymp 1$ in~$L$. 
Put $y_{\rho}: = u_{\rho} + x_{\rho}$ and
$b_{\rho}:= a_{\rho} + g_{\rho} x_{\rho}\in K$, so
$b_\rho - a =g_{\rho}y_{\rho}$. It follows that~$(b_{\rho})$ pseudoconverges to $a$ and has the same width as $(a_\rho)$, so
by Lemma~\ref{pc5} it is a pc-sequence in $K$ equivalent to $(a_{\rho})$. We have 
$$G(b_\rho) - G(a)= 
\sum_{|\i|\ge 1} G_{(\i)}(a)(g_\rho y_\rho)^{\i}$$
where $G_{(\i)}=\frac{G^{(\i)}}{\i !}$. 
Put $g_{\i}:=G_{(\i)}(a)\in L$ for 
$|\i|\ge 1$. Then
\begin{align*} G(b_\rho)-G(a)\ &=\ 
\sum_{|\i|\ge 1} g_{\i}(g_\rho y_\rho)^{\i}\
=\ P(g_{\rho}y_{\rho})\ =\ P_{\times g_{\rho}}(y_\rho), \text{ where }\\
 P(Y)\ :=\ G(a+Y)-G(a)\ &=\  
\sum_{|\i|\ge 1} g_{\i}Y^{\i}\in L\{Y\}, \text{ so $\deg P\ge 1$, $P(0)=0$.}
\end{align*}
By Lemma~\ref{vpnice} with $L$ instead of $K$, we have for
each $\rho$ a thin set $T_{\rho}\subseteq \k_L$, independent
of the choice of $x_{\rho}$,  
such that $0\in T_{\rho}$ and for all $y\in L$,
$$ y\asymp 1,\ \bar{y}\notin T_{\rho}\ \Longrightarrow\  v\big(P(g_\rho y)\big)\ =\ v_P(\gamma_{\rho}).  $$
Note that  $v_P(\gamma_{\rho})$ is strictly increasing as a function of
$\rho$. By Lemma~\ref{thindesc} we can take for each $\rho$ a thin set $S_{\rho}$ in $\k$
such that for all $e\in \k$, if $\bar{u}_{\rho}+ e\in  T_{\rho}$, then 
$e\in S_{\rho}$. 
Therefore, by choosing the $x_{\rho}$ such that $\bar{x}_{\rho}\notin S_{\rho}$
we get $\bar{y}_{\rho}\notin T_{\rho}$, and then
$\big(G(b_{\rho})\big)$ is a pc-sequence and 
$G(b_{\rho}) \leadsto G(a)$.
\end{proof}

\noindent
Note that $v(b_{\rho}-a)=\gamma_{\rho}$ and $v\big(P(b_{\rho}-a)\big)=v_P(\gamma_{\rho})$, eventually, for $(b_{\rho})$ as in the above proof. In addition we
can arrange that for
$e=1,\dots, \deg P$, if $P_e\ne 0$, then $v\big(P_e(b_{\rho}-a)\big)=v_{P_e}(\gamma_{\rho})$, eventually. 
Indeed, all this works just as well for finitely many differential polynomials:

\begin{lemma}\label{pcc2} Let $\mathcal{H}$ be a finite subset of $L\{Y\}$. 
Then there is a pc-sequence $(b_\rho)$ in $K$ equivalent to $(a_\rho)$ 
such that for each $H\in \mathcal{H}$, if $H\notin L$, then
$\big( H(b_\rho)\big)$ is a pc-sequence with 
$H(b_\rho) \leadsto H(a)$.
\end{lemma}

\noindent 
When $G(a_\rho) \leadsto 0$ we can improve these lemmas as follows:

\begin{lemma}\label{pcc3} Suppose that $G(a_\rho) \leadsto 0$, and let 
$\mathcal{H}$ be a finite subset of $L\{Y\}$.  Then there
is a pc-sequence $(b_\rho)$ in $K$ that is equivalent to $(a_\rho)$, 
such that $G(b_\rho) \leadsto 0$, and for each $H\in\mathcal{H}$, if 
$H\notin L$, then $\big(H(b_\rho)\big)$ is a pc-sequence with 
$H(b_\rho) \leadsto H(a)$.
\end{lemma}
\begin{proof} We can assume $G\in \mathcal{H}$ and
$v\big(G(a_\rho)\big)$ strictly increases with $\rho$. 
We now make the same reductions as in the proof of
Lemma~\ref{pcc1}, introducing 
$g_\rho$, $u_{\rho}$, $x_{\rho}$, $y_{\rho}$, $b_{\rho}$, 
$P=  
\sum_{|\i|\ge 1} g_{\i}Y^{\i}\in L\{Y\}$, accordingly. 
The proof of that lemma shows how to
arrange that $(b_\rho)$ is a pc-sequence in $K$ such that  
$H(b_\rho) \leadsto H(a)$, for each $H\in \mathcal{H}$ with $H\notin L$, and 
$$ v\big(P(g_\rho y_{\rho})\big)\ =\ v_P(\gamma_{\rho}),\ \text{ eventually.}  $$
We claim that then $G(b_\rho) \leadsto 0$.  To prove this claim, note that for all
$\rho$,
$$G(a_\rho) - G(a)\ =\ \sum_{|\i|\ge 1} g_\i(a_\rho -a)^{\i}\ =\ P_{\times g_{\rho}}(u_{\rho})\ \preceq\ P_{\times g_{\rho}}(y_{\rho})\ =\ G(b_{\rho})-G(a).$$
There can only be one $\rho$ with $G(a_{\rho})-G(a)\prec G(a)$, because for such $\rho$ we have $v\big(G(a_{\rho})\big) =v\big(G(a)\big)$.
So $G(a_{\rho})-G(a)\succeq G(a)$, eventually, hence
$G(b_{\rho})-G(a)\succeq G(a)$, eventually, and thus $G(b_{\rho})-G(a)\succ G(a)$, eventually. It now follows that $v\big(G(b_{\rho})\big)=v\big(G(b_{\rho})-G(a)\big)$, eventually, so $G(b_{\rho})\leadsto 0$.  
\end{proof}

\section{Constructing Canonical Immediate Extensions}\label{sec:cimex}

\noindent
{\em In this section 
we assume the derivation on $\k$ is nontrivial}.
Thus Lemmas~\ref{pcc1},  \ref{pcc2}, and~\ref{pcc3} are available for our $K$. 

In Lemma~\ref{diftr}, Corollary~\ref{cor:diftr}, and Lemma~\ref{zda} we fix a pc-sequence $(a_\rho)$ in~$K$. Recall
from the end of Section~\ref{Valdifcon} the notions of $(a_{\rho})$
being of $\d$-algebraic type over~$K$, that of $G\in K\{Y\}$ being a minimal differential polynomial of $(a_{\rho})$ over~$K$, and
of~$(a_{\rho})$ being of $\d$-transcendental type over $K$. 
We are going to associate to~$K$ and the
pc-sequence $(a_{\rho})$ an immediate valued differential field extension~$K\<a\>$. If~$(a_{\rho})$ is of $\d$-transcendental type, this extension
is canonical, and if it is of $\d$-algebraic type, it is canonical
modulo a choice of minimal differential polynomial.

\begin{lemma}\label{diftr} 
Suppose $(a_\rho)$ is of $\d$-transcendental type
over $K$.  Then $K$ has an
immediate valued differential field extension $K\langle a \rangle$
such that: \begin{enumerate} 
\item[\textup{(i)}] $\der \smallo_{K\<a\>}\subseteq \smallo_{K\<a\>}$,
$a_\rho \leadsto a$, and $a$ is $\d$-transcendental over $K$;
\item[\textup{(ii)}] for any valued differential field extension $L$ of 
$K$ with $\der\smallo_L\subseteq \smallo_L$ and any $b\in L$ with
$a_\rho \leadsto b$ there is a unique valued differential field
embedding $K\langle a \rangle\ \longrightarrow\ 
L$ over $K$ that sends $a$ to $b$. 
\end{enumerate}
\end{lemma}
\begin{proof} Let $F$ be an elementary extension
of $K$ containing a pseudolimit $a$ of~$(a_\rho)$.  
Let  $K\langle a \rangle$
be the valued differential subfield of $F$ generated by $a$ over $K$.
Let $G(Y)\in K\{Y\}$, $G\notin K$.  
Lemma~\ref{pcc1} gives a pc-sequence
$(b_\rho)$ in $K$ equivalent to~$(a_{\rho})$ such that 
$G(b_\rho) \leadsto  G(a).$
Now, $G(b_\rho) \not\leadsto 0$, since $(a_\rho)$
is of $\d$-transcendental type.  So $G(a)\ne 0$ and 
$G(b_\rho)\sim G(a)$, eventually.
Since $G$ was arbitrary, we see that $a$ is $\d$-transcendental over $K$ and $K\langle a \rangle$ is an immediate extension of~$K$. From $\der\smallo_F\subseteq \smallo_F$ we get $\der \smallo_{K\<a\>}\subseteq \smallo_{K\<a\>}$.  

Let $L$ and $b$ be as in (ii). By the proof of Lemma~\ref{pcc1} 
we can arrange in the argument above, in addition to 
$G(b_\rho) \leadsto  G(a)$, that
$G(b_\rho) \leadsto  G(b)$; hence $v_{K\<a\>}\big(G(a)\big)= v_L\big(G(b)\big)\in \Gamma$. 
%with $v_a$ the valuation of $K\<a\>$.  
\end{proof} 

\noindent
The following is immediate from Lemma~\ref{diftr}:

{\sloppy
\begin{cor}\label{cor:diftr} Let $L$ be a valued differential field extension
of $K$ satisfying~${\der \smallo_L\subseteq \smallo_L}$, and assume 
$a_{\rho}\leadsto b$, where $b\in L$ is
$\d$-algebraic over $K$. Then~$(a_{\rho})$ is of 
$\d$-algebraic type over $K$.
\end{cor} 
}

\begin{lemma}\label{zda} Suppose $P$ is a minimal differential polynomial of 
$(a_{\rho})$ over $K$. Then~$K$ has an immediate valued differential 
field extension $K\<a\>$ such that:
\begin{enumerate} 
\item[\textup{(i)}] $\der \smallo_{K\<a\>}\subseteq \smallo_{K\<a\>}$,
$a_\rho \leadsto a$ and $P(a)=0$;
\item[\textup{(ii)}] for any valued differential field extension $L$ of 
$K$ with $\der\smallo_L\subseteq \smallo_L$ and any $b\in L$ with
$a_\rho \leadsto b$ and $P(b)=0$ there is a unique valued differential field
embedding $K\langle a \rangle \to
L$ over $K$ that sends $a$ to $b$. 
\end{enumerate}
\end{lemma}
\begin{proof} Let $P$ have order $r$ and take $p\in K[Y_0,\dots,Y_r]$ such that
$$P\ =\ p\big(Y, Y',\dots, Y^{(r)}\big).$$ Then $p$ is irreducible by 
Corollary~\ref{cor:pseudolimit 0}. Consider the domain 
$$K[y_0,\dots, y_r]\ :=\ K[Y_0,\dots, Y_r]/(p), \qquad y_i:=Y_i + (p) \text{ for $i=0,\dots,r$,}$$
and let $F=K(y_0,\dots,y_r)$ be its fraction field. 

{\sloppy
We extend the valuation $v$ on $K$ to a valuation $v\colon F^\times \to \Gamma$ as
follows.  Pick a pseu\-do\-li\-mit 
$e$ of $(a_\rho)$ in some valued differential field extension $E$ of $K$
with~${\der \smallo_E \subseteq \smallo_E}$. We let $v$ also denote the 
valuation of $E$. Let
$\phi\in F^\times$, so 
$$\phi\ =\ f(y_0,\dots,y_r)/g(y_0,\dots,y_{r-1})$$
with $f\in K[Y_0,\ldots,Y_r]$ of lower degree in $Y_r$
than $p$ and  $g\in K[Y_0,\ldots,Y_{r-1}]^{\ne}$. 
Set $\vec e:=\big(e, e',\dots, e^{(r)}\big)$, and also 
$\vec b=\big(b, b',\dots,b^{(r)}\big)$ for $b\in K$. 
}

\medskip\noindent
By Lemma~\ref{pcc2} we can take a pc-sequence $(b_{\rho})$ in $K$ equivalent to
$(a_{\rho})$ such that if $f\notin K$, then
$f(\vec{b}_{\rho})\leadsto f(\vec e)$, and if $g\notin K$, then
$g(\vec{b}_{\rho})\leadsto g(\vec e)$. Also $f(\vec{b}_{\rho})\not\leadsto 0$ 
and $g(\vec{b}_{\rho})\not\leadsto 0$ by the minimality 
of $P$, so  eventually, $f(\vec{b}_{\rho})\sim f(\vec e)$ and 
$g(\vec{b}_{\rho})\sim g(\vec e)$, in particular, $f(\vec e)\ne 0$ and $g(\vec e)\ne 0$, and $v\big(f(\vec{e})\big)\in \Gamma$ and 
$v\big(g(\vec{e})\big)\in \Gamma$. 

\claim{$f(\vec e)/g(\vec e)$ depends only on $\phi$ and not on the choice of $(f,g)$.}  

\noindent
To see why this claim is true, suppose that
also $$\phi=f_1(y_0,\dots,y_r)/g_1(y_0,\dots,y_{r-1})$$
with $f_1\in K[Y_0,\ldots,Y_r]$ of lower degree in $Y_r$
than $p$ and  $g_1\in K[Y_0,\ldots,Y_{r-1}]^{\neq}$. Then 
$fg_1\equiv f_1g \bmod p$ in $K[Y_0,\dots,Y_r]$, and thus $fg_1=f_1g$ since
$fg_1$ and $f_1g$ have lower degree in $Y_r$ than $p$. Thus
$f(\vec e)/g(\vec e)=f_1(\vec e)/g_1(\vec e)$, as promised. 
This proves the claim and allows us to define $v\colon F^\times \to \Gamma$
by 
$$v(\phi)\ :=\ v\big(f(\vec e)/g(\vec e)\big)\ =\ v\big(f(\vec{e})\big)-
v\big(g(\vec{e})\big).$$
Clearly, this $v$ extends the valuation of $K$. Let $\phi_1, \phi_2\in F^{\times}$. It is easy to check that if $\phi_1+\phi_2\ne 0$, then $v(\phi_1+\phi_2)\ge \min(v\phi_1, v\phi_2)$. 
We have
$$\phi_1=\frac{f_1(y_0,\dots,y_r)}{g_1(y_0,\dots,y_{r-1})}, \quad \phi_2=\frac{f_2(y_0,\dots,y_r)}{g_2(y_0,\dots,y_{r-1})}, \quad
 \phi_1\phi_2=\frac{f_3(y_0,\dots,y_r)}{g_3(y_0,\dots,y_{r-1})}
$$
where $f_1, f_2, f_3\in K[Y_0,\ldots,Y_r]$ have lower degree in $Y_r$
than $p$ and where $g_1,g_2, g_3$ are nonzero polynomials
in $K[Y_0,\ldots,Y_{r-1}]$. Then
$$\frac{f_1}{g_1}\frac{f_2}{g_2}\ =\ \frac{pq}{g_1g_2g_3} + \frac{f_3}{g_3}$$ 
with $q \in K[Y_0, \dots ,Y_r]$. Lemma~\ref{pcc3} gives a pc-sequence $(b_{\lambda})$ in $K$ equivalent to $(a_{\rho})$, such that $p(\vec{b}_{\lambda}) \leadsto 0$, and if $q\notin K$, then
$q(\vec{b}_{\lambda})\leadsto q(\vec e)$, and such that for $i=1,2,3$
we have $f_i(\vec{b}_{\lambda})\sim f_i(\vec e)$ and
$g_i(\vec{b}_{\lambda})\sim g_i(\vec e)$, eventually.
This gives $f_1(\vec{b}_\lambda)f_2(\vec{b}_{\lambda})/g_1(\vec{b}_{\lambda})g_2(\vec{b}_{\lambda})\sim f_3(\vec{b}_{\lambda})/g_3(\vec{b}_{\lambda})$, eventually, and thus 
$$v(\phi_1\phi_2)\ =\ v(\phi_1) + v(\phi_2).$$ 
Thus $v\colon F^\times \to \Gamma$ is indeed a valuation on $F$, and below we consider $F$ as a valued field accordingly. Now let $f\in K[Y_0,\ldots,Y_r]$ be of lower degree in $Y_r$ than $p$, $f\notin K$. Take a pc-sequence
$(b_{\rho})$ in $K$ equivalent to $(a_{\rho})$ such that $f(\vec{b}_{\rho})\leadsto f(\vec e)$. Then 
$$v\big(f(y_0,\dots, y_r)-f(\vec{b}_{\rho})\big)\ =\ v\big(f(\vec e)-f(\vec{b}_{\rho})\big),\ \text{ for all $\rho$.}$$
As $f(\vec e)\sim f(\vec{b}_{\rho})$, eventually, this gives
$f(y_0,\dots, y_r)\sim f(\vec{b}_{\rho})$, eventually. Thus $F$ is an immediate extension of $K$. 

\medskip\noindent
We now equip $F$ with the derivation extending the derivation of $K$
such that $y_i'=y_{i+1}$ for $0\le i < r$. Setting $a:= y_0$ we have $a^{(i)}=y_i$
for $i=0,\dots,r$, $K\<a\>=F$, and $P(a)=0$.  We claim that $a_{\rho}\leadsto a$. Consider first the case that $r=0$ and $p$ has degree $1$ in $Y_0$. Then $p=fY_0+g$
with $f,g\in K$, $f\ne 0$, so $a=y_0=-g/f \in K$, and we have a pc-sequence
$(b_{\lambda})$ in $K$ equivalent to $(a_{\rho})$ such that
$p(b_\lambda)\leadsto 0$, that is, $b_{\lambda}\leadsto a$, so
 $a_{\rho}\leadsto a$. If $r=0$ and $p$ has degree $>1$ in $Y_0$, or 
 $r>0$, then $a\notin K$ and $v(a-a_{\rho})=v(y_0-a_{\rho})=v(e-a_{\rho})$ for all $\rho$, which again gives $a_{\rho}\leadsto a$.

\medskip\noindent
To get $\der\smallo_{F}\subseteq \smallo_{F}$, we set
$$S\ :=\ \big\{g(a):\ g\in K\big[Y,\dots, Y^{(r-1)}\big],\ g(a)\preceq 1\big\}.$$
(If $r=0$, then we have $K\big[Y,\dots, Y^{(r-1)}\big]=K$, so $S=\mathcal{O}$.)
By Lemma~\ref{keycor} applied to~$K\big(a,\dots, a^{(r-1)}\big)$ in the role of $E$, it is enough to show
that $\der S\subseteq \mathcal{O}_{F}$ and~${\der(S\cap \smallo_F)\subseteq \smallo_F}$. We prove the first of these inclusions. The second follows in 
the same way. 

\medskip\noindent
Let $g\in K\big[Y,\dots,Y^{(r-1)}\big]\setminus K$ with $g(a)\preceq 1$;
we have to show $g(a)'\preceq 1$. We can assume $g(a)'\ne 0$. 
Take  $g_1(Y), g_2(Y)\in K\big[Y,\dots,Y^{(r-1)}\big]$ such that
$$g(Y)'\  =\  g_1(Y) + g_2(Y)Y^{(r)}  \quad \text{in $K\{Y\}$.}$$
Then $$g(a)'\ =\  
g_1(a) + g_2(a)a^{(r)},$$ and for all $y\in K$,
$$g(y)'\ =\   g_1(y) + g_2(y)y^{(r)}.$$
Take a pc-sequence $(b_{\lambda})$ in $K$ equivalent to $(a_{\rho})$ such
that $g(b_{\lambda})\leadsto g(e)$. Hence $g(b_{\lambda})\sim g(e)$,
eventually, so $g(b_{\lambda})\sim g(a)$, eventually, and thus 
$g(b_{\lambda})'\preceq 1$, eventually.
We now distinguish two cases:

\case[1]{$P$ has degree $>1$ in $Y^{(r)}$, or $g_2=0$.}
Then we can assume that $(b_{\lambda})$ has been chosen such that in addition we have, eventually, 
$$ g(b_{\lambda})'\ =\ g_1(b_{\lambda})+ g_2(b_{\lambda})b_{\lambda}^{(r)}\  
 \sim\ g_1(a) + g_2(a)a^{(r)}\ =\ g(a)'.$$
Therefore $g(a)'\preceq 1$, as desired.   

\case[2]{$P$ has degree $1$ in  $Y^{(r)}$, and $g_2 \ne 0$.}
Then 
$$g_1 + g_2Y^{(r)}\ =\ \frac{h_1P + h_2}{h}, \qquad h,h_1,h_2\in K[Y,\dots,Y^{(r-1)}],\ h, h_1\ne 0,$$ so
$0\ne g(a)'=h_2(a)/h(a)$, so $h_2\ne 0$. As in the previous case we can
assume
$(b_{\lambda})$ to have been chosen such that in addition we have, eventually, 
$$P(b_{\lambda})\leadsto 0, \quad h(b_{\lambda})\sim h(a), \quad
h_1(b_{\lambda})\sim h_1(a), \quad h_2(b_{\lambda})\sim h_2(a).$$
Now, eventually $g(b_{\lambda})'\preceq 1$, so eventually
$h_1(b_{\lambda})P(b_{\lambda}) + h_2(b_{\lambda})\preceq h(b_{\lambda})$.
From this it follows easily that $h_2(b_{\lambda})\preceq h(b_{\lambda})$, eventually, so $h_2(a)\preceq h(a)$, that is, $g(a)'\preceq 1$.
This finishes the proof of (i). We now turn to (ii). 

\medskip\noindent
Suppose $L$ is a valued differential field extension of $K$
with $\der \smallo_L\subseteq \smallo_L$, and $b\in L$ satisfies
$P(b)=0$ and $a_{\rho}\leadsto b$. Let $Q\in K[Y,\dots, Y^{(r)}]$, 
$Q\notin K$, and suppose $Q$ has lower degree in $Y^{(r)}$ than $P$. Then $Q(b)\ne 0$ by the argument in the beginning of the proof with 
$b$ in the role of $e$.
Thus $P$ is a minimal annihilator of $b$ over $K$. By
Lemma~\ref{pcc1} (or rather its proof) we have
 a pc-sequence $(b_{\rho})$ in $K$ equivalent to~$(a_{\rho})$ such that
$Q(b_{\rho})\leadsto Q(e)$ as well as
$Q(b_{\rho})\leadsto Q(b)$. As in the beginning of the
proof this gives 
$$v\big(Q(b_{\rho})\big)\ =\ v\big(Q(e)\big)\ =\ v\big(Q(a)\big), \quad \text{eventually,}$$ 
and also $v\big(Q(b_{\rho})\big)=v_L\big(Q(b)\big)$, eventually.   
In particular, $v\big(Q(a)\big)=v_L\big(Q(b)\big)$. 
Thus the
differential field embedding $K\langle a \rangle \to L$ over $K$ sending $a$ to $b$ is also a valued field embedding.
\end{proof}

\noindent
Here are two immediate consequences of Lemmas~\ref{diftr} and~\ref{zda}:

\begin{cor}\label{cor1diftrzda} If $K$ has no proper immediate valued differential field extension 
with small derivation, then $K$ is spherically complete.
\end{cor}

\begin{cor}\label{cor2diftrzda} $K$ has an immediate valued differential field extension with small
derivation that is spherically complete. 
\end{cor}

\noindent
By Corollary~\ref{cor:unique max imm ext, 1},  any two extensions of $K$ as in the last corollary are
isomorphic over $K$ as valued fields; we would like to improve this to
being isomorphic over $K$ as valued differential fields. In Section~\ref{umie}
we establish this improvement under some extra assumptions on $K$.
This involves an extension of the notion of ``henselian''
to the setting of valued differential fields, to be developed
in the next chapter.

%% file: mt-7.tex
\chapter{Differential-Henselian Fields}\label{sec:dh1}

\setcounter{theorem}{0}

\noindent
{\em In this chapter $K$ is a valued differential
field with small derivation}. As usual, $\Gamma:=v(K^\times)$ and
$\k:=\res(K)$, the latter a {\em differential\/} field. By an {\em extension of\/} $K$ we mean
a valued differential field extension of $K$ whose derivation is small.
 
In this setting we study {\em differential-henselianity\/} and establish useful results about this notion in analogy with various facts about henselian
valued fields in Section~\ref{sec:henselian valued fields}.
We say that $K$ is {\bf differential-henselian} (for short: {\bf $\d$-henselian}) if the two conditions
below are satisfied:
\begin{list}{*}{\setlength\leftmargin{3em}}
\item[(DH1)] $\k$ is linearly surjective;
\item[(DH2)] for every $P\in \mathcal O\{Y\}$ with $P_0\prec 1$ and
$P_1\asymp 1$, there is $y\prec 1$ in $K$ such that $P(y)=0$.
\end{list}
If $\Gamma=\{0\}$ this is the same as $K$ being linearly surjective. Note that 
for $P\in  \mathcal O\{Y\}$ the condition $P_0\prec 1$ and
$P_1\asymp 1$ in (DH2) implies $\dval P =1$.       

\index{valued differential field!differential henselian}
\index{differential-henselian}
\index{d-henselian@$\d$-henselian}

\index{maximal!valued differential field}
\index{differential-algebraically maximal}
\index{valued differential field!maximal}
\index{valued differential field!differential-algebraically maximal}

Define $K$ to be {\bf maximal\/} if $K$ has
no proper immediate extension. If the derivation of $\k$ is nontrivial, then by Corollary~\ref{cor1diftrzda},  
$$   \text{$K$ is maximal}\ \Longleftrightarrow\    \text{$K$ is spherically complete.}$$ By Zorn's Lemma, $K$ does have a maximal 
immediate extension.  
Define $K$ to be {\bf differential-algebraically maximal\/} (for short: {\bf $\d$-algebraically maximal\/}) if $K$ has no proper immediate $\d$-algebraic extension. If the derivation of $\k$ is nontrivial, then by Lemma~\ref{zda}, $K$ is $\d$-algebraically maximal iff there is no divergent
pc-sequence in $K$ of $\d$-algebraic type over $K$. By Zorn,
$K$ does have an immediate $\d$-algebraic extension that is
$\d$-algebraically maximal. 

It is obvious that
{\em $\d$-henselianity\/} can be formulated
in the language of valued differential fields as a first-order
axiom scheme. For {\em $\d$-algebraic maximality\/} this is not obvious,
and perhaps false. We study below how these two notions are related.

After preliminaries about $\d$-henselianity in Section~\ref{preldh} we prove the following in Section~\ref{sec:maxdh}, in analogy with one direction of Corollary~\ref{cor:alg max equals henselian}: 

\begin{theorem}\label{thm:damdh} If $\k$ is linearly surjective and 
$K$ is $\d$-algebraically maximal, then~$K$ is $\d$-henselian.
\end{theorem}

\noindent
This depends critically on the $d=1$ case of the Equalizer Theorem.
An immediate consequence is a differential analogue of Hensel's Lemma:

{\sloppy
\begin{cor}\label{cor:mdh} If $\k$ is linearly surjective and $K$ is spherically complete, then~$K$ is $\d$-henselian. 
\end{cor}
}
\noindent
In particular, if $\k$ is a linearly surjective differential field, then the Hahn differential field 
$\k\(( t^{\Gamma}\)) $ 
defined in Section~\ref{Valdifcon} is $\d$-henselian.

\medskip\noindent
For monotone $K$ with linearly surjective
$\k$ we prove in Section~\ref{umie} the 
uniqueness-up-to-isomorphism-over-$K$ of maximal immediate extensions. This will be crucial in the next chapter. In Section~\ref{sec:fewcon} we assume $C\subseteq \mathcal{O}$ and show that in the presence
of monotonicity (perhaps unnecessary) we have a converse to Theorem~\ref{thm:damdh}:

\begin{theorem}\label{thm:fcdifhensalgmax} If $C\subseteq \mathcal{O}$ and $K$ is monotone and $\d$-henselian, then $K$ is $\d$-al\-ge\-brai\-cal\-ly maximal.
\end{theorem}

\noindent
Finally, we consider differential-henselianity in several variables and obtain a partial analogue of Proposition~\ref{prop:hensel, multivar}. In detail, let $Y=(Y_1,\dots, Y_n)$ be a tuple of distinct differential indeterminates, $n\ge 1$.
Let $P_1(Y),\dots, P_n(Y)\in \mathcal{O}\{Y\}$.
We consider the system of equations
$$ P_1(Y)\ =\ \cdots\ =\ P_n(Y)\ =\ 0.$$ 
Let $A_i\in \mathcal{O}\{Y\}$ be the homogeneous part 
of $P_i$ of
degree $1$, and let $\bar{A}_i$ be its image in~$\k\{Y\}_1$. Recall
the notion of $\d$-independence defined in Section~\ref{sec:systems}.

\begin{theorem}\label{thm:maxndh} Suppose $\k$ is linearly surjective,
$K$ is $\d$-algebraically maximal, $\bar{A}_1,\dots, \bar{A}_n\in \k\{Y\}_1$ are $\d$-independent, and $P_1(0)\prec 1,\dots, P_n(0)\prec 1$, with $0:=(0,\dots,0)\in K^n$. Then there exists a tuple $y=(y_1,\dots, y_n)\in K^n$ such that
$$ P_1(y)\ =\ \cdots\ =\ P_n(y)\ =\ 0, \qquad y_1\prec 1,\ \dots,\ y_n\prec 1.$$
\end{theorem}

\noindent
In Section~\ref{sec:max} we prove this first for spherically complete $K$ by approximation
arguments, using also heavily the material in Section~\ref{sec:systems}. The theorem then follows by appealing to 
the result from Section~\ref{sec:johnson} due to Johnson~\cite{JJ}.

\section{Preliminaries on Differential-Henselianity}\label{preldh}

\noindent
Throughout this chapter we assume $r\in \N$. To allow certain kinds of inductive arguments, we define $K$ to be {\bf $r$-differential-henselian} (for short: {\bf $r$-$\d$-henselian}) if the two conditions
below are satisfied:
\begin{list}{*}{\setlength\leftmargin{3em}}
\item[(DH$r$1)] $\k$ is $r$-linearly surjective;
\item[(DH$r$2)] for every $P\in \mathcal O\{Y\}$ of order $\le r$ with $P_0\prec 1$ and
$P_1\asymp 1$, there is $y\prec 1$ in $K$ such that $P(y)=0$.
\end{list}
For $r=0$ this is the same as $K$ being henselian as a valued field. If $\Gamma=\{0\}$ it is the same as $K$ being $r$-linearly surjective. Note that $K$ is $\d$-henselian iff $K$ is $r$-$\d$-henselian for each $r$.

\begin{lemma}\label{difhen1} Suppose $K$ is $r$-$\d$-henselian, 
$P\in \mathcal O\{Y\}$ has order $\le r$, 
$P_1\asymp 1$, and $P_i\prec 1$ for all $i\ge 2$. Then 
$P(y)=0$ for some $y\in \mathcal O$. 
\end{lemma}
\begin{proof} The assumption on $P$ gives $\ddeg P =1$. 
Use (DH$r$1) to get $u\in \mathcal O$ with
$D_P(\bar{u})=0$, and thus $\val (D_P)_{+\bar{u}}=\deg (D_P)_{+\bar{u}}=1$.
Therefore $\dval P_{+u}=\ddeg P_{+u}=1$ by Lemma~\ref{dn1,new}(i). 
Also $P_{+ u}\asymp 1$ by Lemma~\ref{v-under-conjugation}(i). Now  
apply (DH$r$2) to $P_{+u}$
in the role of $P$ to get $y\in \smallo$ with $P(u+y)=0$.
\end{proof}

% $P(u)\prec 1$. By Lemma~\ref{v-under-conjugation},~(1) we have 
%$P_{+ u}\asymp 1$ and $P_{i,+u}\prec 1$ for $i=2,\dots,r$. Hence 
%$P_{1,+u} \asymp 1$, so we can apply (DH$r$2) to $P_{+u}$
%in the role of $P$ to get $y\in \smallo$ with $P(u+y)=0$.
%\end{proof}
  
\begin{cor}
Suppose $K$ is $r$-$\d$-henselian, and $A\in K[\der]$ with $v(A)=0$ has order at most $r$.
Then $A(\smallo)=\smallo$ and $A(\mathcal O)=\mathcal O$. 
\end{cor}
\begin{proof}
The first statement follows from (DH$r$2), and the second statement from Lemma~\ref{difhen1}.
\end{proof}

\noindent
For $P\in \mathcal O\{Y\}$ and $a\in \mathcal O$ we say that $P$ is in 
{\bf differential-hensel position at~$a$} 
(abbreviated as {\bf dh-position at~$a$})
if $P(a) \prec 1$ and $P_{+a,1}\asymp 1$. Note that then~$P$ is 
in dh-position at $b$ for each $b\in \mathcal O$ with $a-b\in \smallo$. 
Note also that if $K$ is 
$r$-$\d$-henselian, $P\in \mathcal O\{Y\}$ is of order $\le r$ and 
$P$ is in dh-position at $a\in \mathcal O$, then there is $b\in \mathcal O$ with 
$a-b\in \smallo$ such that $P(b)=0$. This gives the following analogue of
an important result (Proposition~\ref{prop:lift}) about henselian valued fields:

\index{differential-hensel!position}
\index{dh-position}

\begin{prop}\label{lift.res.field} Suppose $K$ is $\d$-henselian.
Then $\k$ can be lifted to a differential subfield of $K$, that is,
there is a differential subfield $F$ of $K$ such that $F\subseteq \mathcal O$
and~$F$ maps \textup{(}isomorphically\textup{)} onto $\k$ under the
residue map $\mathcal{O}\to \k$.
\end{prop} 
\begin{proof} This is the case $E=\Q$ of the following more general result:
Let $E\subseteq \mathcal{O}$ be a differential subfield of $K$; then there 
is a differential subfield $F$ of $K$ such that
$E\subseteq F \subseteq \mathcal{O}$ and $F$ maps onto $\k$ under the
residue map $\mathcal{O}\to \k$. 
Suppose $\res(E)\ne \k$. Take $a\in \mathcal{O}$ such that 
$\overline{a}\notin \res(E)$. If $P(a)\asymp 1$ for all nonzero
$P(Y)\in E\{Y\}$, then~$E\langle a\rangle$ is a proper
differential field extension of $E$ contained in~$\mathcal{O}$. Next,
consider the case that $P(a)\prec 1$  for some nonzero
$P(Y)\in E\{Y\}$. Pick such $P$ of minimal complexity, say $P$ has order $r$.  
Then $Q(a)\asymp 1$
for all nonzero $Q(Y)\in E\{Y\}$ of lower complexity, hence 
$\frac{\partial P}{\partial Y^{(r)}}(a)\asymp 1$, so
$P$ is in dh-position at $a$. 
This gives $b\in \mathcal{O}$ with $P(b)=0$ and 
$\overline{a}=\overline{b}$. Since $\overline{a}\notin\res(E)$, we have $b\notin E$. 
It follows that $E\langle b \rangle$ is a proper
differential field extension of $E$ contained in $\mathcal{O}$. 
We finish the proof by in\-vo\-king~Zorn.
\end{proof}

\subsection*{An embedding result}
We now relate $\d$-henselianity to the material on
residue extensions in Section~\ref{sec:resext}. Let~$L$ and~$F$ be 
extensions of $K$. Then 
$\k$ is a {\em differential\/}
subfield of $\k_L$ and of $\k_F$. Assume also that
$L$ and $F$ are $r$-$\d$-henselian, and let 
$i\colon \k_L \to \k_F$ be a differential field embedding over $\k$.

\begin{lemma}\label{resembdh} Suppose $\overline{a}\in \k_L$ is $\d$-algebraic 
over $\k$, with minimal
annihilator $\overline{P}(Y)\in \k\{Y\}$ of order $r$ over $\k$. 
Then there exist $b\in \mathcal{O}_L$ with $\res(b)=\overline{a}$,
$\k_{K\<b\>}=\k\<\overline{a}\>$,
and $v(K\<b\>^\times)=\Gamma$, and a valued differential field embedding
$j\colon K\<b\> \to F$, such that $\res(jy)=i\big(\!\res(y)\big)$ for all 
$y\in \mathcal{O}_{K\<b\>}$.
\end{lemma}
\begin{proof}
Take $P\in \mathcal{O}\{Y\}$ such that
$P$ has image $\overline{P}$ in $\k\{Y\}$ and $P$ has the same complexity 
as $\overline{P}$. Then $v(I)=0$ where $I$ is the initial of $P$. 
Take also $a\in \mathcal{O}_L$ with $\res(a)=\overline{a}$.
Then $\frac{\partial P}{\partial Y^{(r)}}(a) \asymp 1$, 
so $P$ is in dh-position at $a$. This gives $b\in \mathcal{O}_L$
with $P(b)=0$ and $\res(b)=\overline{a}$. If $Q\in K\{Y\}^{\ne}$ has lower 
complexity than $P$, then $Q(b)\ne 0$ by an easy reduction to the case 
$Q\asymp 1$. Thus $P$ is a minimal
annihilator of $b$ over $K$. Likewise, we get $f\in \mathcal{O}_F$
with $P(f)=0$ and $\res(f)=i(\overline{a})$. Then~$P$ is also a minimal 
annihilator of $f$ over $K$, and so Theorem~\ref{resext} gives
a valued differential field embedding $j\colon K\<b\> \to F$
with $jb=f$. This $j$ has the desired property.
\end{proof}

\noindent
In the above set-up, assume that $L$ and $F$ are even $\d$-henselian
(rather than $r$-$\d$-hen\-sel\-ian). Then by Lemma~\ref{resembdh} and the
result on $\d$-transcendental residue extensions preceding 
Theorem~\ref{resext}: 

\begin{cor} There exist a valued differential subfield 
$E\supseteq K$ of $L$ such that $\k_{E}=\k_L$ and $v(E^\times)=\Gamma$,
and a valued differential field embedding $j\colon E\to F$ such that
$\res(jb)=i(\res(b))$ for all $b\in \mathcal{O}_{E}$.
\end{cor}

\subsection*{Differential-henselianity and specialization}
We now consider the behavior of $\d$-henselianity with respect to
coarsening and specialization.
Let $\Delta$ be a convex subgroup of $\Gamma$, and $\dot{v}=v_{\Delta}$ the
$\Delta$-coarsening of $v$, with valuation ring $\dot{\mathcal{O}}$ and
maximal ideal $\dot{\smallo}$ of $\dot{\mathcal{O}}$. Then
$\der\dot{\smallo}\subseteq \dot{\smallo}$ by Corollary~\ref{cor:Cohn-Lemma}, and thus $\der \dot{\mathcal{O}}\subseteq \dot{\mathcal{O}}$, which gives a differential residue
field $\dot{K}$. Moreover $\der \smallo_{\dot{K}} \subseteq \smallo_{\dot{K}}$,
where for convenience $\der$ denotes also the derivation of $\dot{K}$.
Thus the residue field $\res(\dot{K})$ of $\dot{K}$ is naturally a differential field, and the canonical ring isomorphism $\res(K) \cong \res(\dot{K})$ is a differential ring isomorphism.

\begin{lemma}\label{inva4} If $K$ is $r$-$\d$-henselian, 
then so is $\dot{K}$.
\end{lemma}
\begin{proof} 
Let $P\in \mathcal O\{Y\}$ of order $\le r$ have image $Q$ in $\dot{K}\{Y\}$. Then
$Q$ has coefficients in the valuation ring of $\dot{K}$. It now remains to note that the differential residue field of~$K$ is isomorphic to the differential residue field of $\dot{K}$, and that 
if ${Q_0\prec 1}$, $Q_1\asymp 1$, then $P_0\prec 1$, $P_1\asymp 1$. 
\end{proof}

\begin{lemma}\label{inva5} Suppose $(K, \dot{\mathcal{O}})$ and $\dot{K}$ are $r$-$\d$-henselian. Then so is $K$.% is $r$-$\d$-hen\-sel\-ian.
\end{lemma}
\begin{proof} Since $K$ and $\dot{K}$ have isomorphic differential residue fields,
$K$ satisfies condition (DHr1). Next, let 
$P\in \mathcal{O}\{Y\}$ of order $\le r$ satisfy $P_0\prec 1$ and $P_1\asymp 1$; it is enough to find $b\prec 1$ in $K$
such that $P(b)=0$. Let $Q$ be the image of $P$ in 
$\dot{K}\{Y\}$; then $Q\in \mathcal{O}_{\dot{K}}\{Y\}$. 
%of $P$ under the differential ring morphism 
%$\mathcal{O}\{Y\} \to \mathcal{O}_{\dot{K}}\{Y\}$ which
%sends each $a\in \mathcal{O}$ to $\dot{a}\in \mathcal{O}%_{\dot{K}}$ and each $Y^{(i)}$ to $Y^{(i)}$. 
Then $Q_0\prec 1$ and $Q_1\asymp 1$, which gives 
$a\in \smallo$
with $Q(\dot{a})=0$. Then $P(a) \dotprec 1$. Now $P_{+a}\sim P$ by Lemma~\ref{v-under-conjugation},
and so $P_{+a,0}=P(a) \dotprec 1$ and $vP_{+a,1}=0$, from which
we get $\dot{v}P_{+a,1}=0$. As $(K, \dot{\mathcal{O}})$ is 
$r$-$\d$-henselian, this gives $y\in \dot{\smallo}$ such that
$P_{+a}(y)=0$, and then $b:= a+y$ satisfies $b\prec 1$ and $P(b)=0$.
\end{proof}

\noindent
So far we have not used
the Equalizer Theorem in studying $\d$-henselianity, but to
progress further we need the $d=1$ case of this theorem, which says that
for $A\in K[\der]^{\ne}$ the function $v_A\colon \Gamma\to \Gamma$ is
bijective. Recall also that this function 
is strictly increasing  and satisfies $v_A(\gamma)=v(A)+\gamma+ o(\gamma)$
for $\gamma\in \Gamma^{\ne}$.

\subsection*{Relation to neat surjectivity} If $K$ is $\d$-henselian,
then not only $\k$ but also $K$ itself is linearly surjective. More precisely:

\begin{lemma}\label{dhns} Suppose $K$ is $r$-$\d$-henselian. Then every
$A\in K[\der]^{\ne}$ of order $\le r$ is neatly surjective.
\end{lemma}
\begin{proof} Let $A\in K[\der]^{\ne}$ have order $\le r$, let $b\in K^\times$,
and take $\alpha\in \Gamma$ such that $v_A(\alpha)=\beta:= vb$. We have to find
$a\in K^\times$ such that $va=\alpha$ and $A(a)=b$. First, take any 
$\phi\in K^\times$ with $v\phi=\alpha$. Then $v(A\phi)=vb$, so $v(B)=0$ with 
$B:= b^{-1}A\phi$. Then Lemma~\ref{difhen1} gives $y\in \mathcal O$ with 
$B(y)=1$, so $y\asymp 1$ and $A(\phi y)=b$. Thus $a:=\phi y$
has the desired property.
\end{proof} 

\begin{cor}\label{dhns2} If $K$ is $1$-$\d$-henselian, then $\smallo=(1+\smallo)^{\dagger}$.
\end{cor}
\begin{proof} Let $a\in \smallo$. We look for $y\in \smallo$ such that
$y'/(1+y)=a$, that is, $ay-y'=-a$. Now use that if $a- \der$ is neatly 
surjective, the equation $ay-y'=-a$ does indeed have a solution $y$ 
in $\smallo$.
\end{proof}

\subsection*{Application to having many constants} For use in the next chapter we show that under certain conditions $\d$-henselianity yields many constants.

\begin{lemma}\label{1dconst} Suppose $K$ is $1$-$\d$-henselian and monotone, and let $b\in K$, $b'\prec b$. Then $b\asymp c$ for some $c\in C$.
\end{lemma}
\begin{proof} Assume that $a\in K^\times$. 

\claim{$y\asymp a$ and $y'=a$ for some $y\in K$.}

\noindent
To see this, note that $(aZ)'-a=a(a^\dagger Z + Z'-1)$ and $a^\dagger\preceq 1$, so by Lemma~\ref{difhen1} 
we get $z\in \mathcal{O}$
with $a^\dagger z+z'=1$. Then $z\asymp 1$, so $y:=az\asymp a$ satisfies $y'=a$.

\medskip\noindent
Returning to $b$, note first that if $b'=0$, then
$c:= b$ gives $c\in C$ with $c\asymp b$.
Assume $b'\ne 0$. Then with $a:= b'$
the claim above gives $y\asymp b'$ such that $y'=b'$, so
$c:=y-b\in C$ and $c\sim -b$. 
\end{proof}

{\sloppy
\begin{cor}\label{cor:1dconst} Suppose $K$ is $1$-$\d$-henselian, monotone, and 
$(\k^{\times})^{\dagger}=\k$. Then~$K$ has many constants
and $(K^{\times})^{\dagger}\ =\ (\mathcal{O}^\times)^\dagger\ =\ 
\mathcal{O}$.
\end{cor}
}
\begin{proof} Let $b\in K^\times$, and suppose $b'\asymp b$; to show
$K$ has many constants, it is enough by Lemma~\ref{1dconst} to get 
$y\asymp 1$ in $K$ such that $(by)'\prec b$, that is,~${(by)^{\dagger}\prec 1}$, that is, ${b^\dagger + y^\dagger \prec 1}$.
Since $b^\dagger \asymp 1$, the assumption of the corollary 
provides such a~${y\asymp 1}$. As to $(\mathcal{O}^\times)^\dagger\ =\ \mathcal{O}$, 
let $a\in \mathcal{O}$. If $a\prec 1$, then
use $\smallo=(1+\smallo)^\dagger$. Assume $a\asymp 1$. Take 
$u\asymp 1$ with $a\sim u^\dagger$. Then 
$a-u^\dagger =(1+\varepsilon)^\dagger,$ $\varepsilon\in \smallo$, so
$a= (u(1+\varepsilon))^\dagger\in  (\mathcal{O}^\times)^\dagger$. This proves the claimed equalities.
\end{proof}

\section{Maximality and Differential-Henselianity}\label{sec:maxdh}

\noindent
We begin with proving a converse to Lemma~\ref{difhen1}.

\begin{lemma}\label{Hdhencor} Suppose 
for all $P\in \mathcal O\{Y\}$ of order $\le r$ with
$P_1\asymp 1$ and $P_i\prec 1$ for all $i\ge 2$, there is $z\in \mathcal O$ with $P(z)=0$.
Then $K$ is $r$-$\d$-henselian.
\end{lemma}
\begin{proof} It is clear that (DH$r$1) holds. For (DH$r$2), 
let $P\in \mathcal O\{Y\}$ be of order~${\le r}$ with $P_0\prec 1$ and
$P_1\asymp 1$; we have to get $y\prec 1$ in $K$ such that $P(y)=0$. If $P_0=0$, then $y=0$ works. Assume $P_0\ne 0$. Take $\gamma\in \Gamma$ such that
$v_{P_1}(\gamma) = v(P_0)$. Then~$\gamma>0$, so 
$$v_{P_{1}}(\gamma) = \gamma + o(\gamma)\ <\  v_{P_{i}}(\gamma)\quad \text{ for all $i\ge 2$.}$$ Take $g\prec 1$ in $K$
with $vg=\gamma$. Then $P_0\asymp P_{1,\times g}$ and 
$P_{i,\times g}\prec P_{1,\times g}$ for all $i\ge 2$. Now take $h\in K^\times$ 
with $h\asymp  P_{1,\times g}$ and apply our hypothesis to $h^{-1}P_{\times g}$ 
in the role of~$P$ to get $z\in \mathcal O$
with $P(gz)=0$. Then $y:= gz\prec 1$ and $P(y)=0$.
\end{proof}

\noindent
In the next four lemmas 
we consider more closely a differential
polynomial $P\in \mathcal O\{Y\}$ of order $\le r$ and $a\in \mathcal O$ such that 
$P$ is in dh-position at $a$, that is, with $Q:=P_{+a}$ we have
$Q(0)=P(a)\prec 1$ and $Q_1\asymp 1$. Let $Q_1=a_0Y+a_1Y' + \cdots + a_rY^{(r)}$
($a_0,\dots,a_r\in \mathcal O$) and set 
$$A\ :=\ L_{P_{+a}}\ =\  a_0 + a_1\der + \dots + a_r\der^r\in K[\der],\ \text{ so $vA=0$.}$$ 
We define $v(P,a)$ to be the unique $\gamma\in \Gamma_{\infty}$ such that 
$v_A(\gamma)=v\big(P(a)\big)$, with $v(P,a)=\infty$ by convention if $P(a)=0$. 
Thus $v(P,a)>0$, and
by Corollary~\ref{vP-Lemma},
$$  P(a)\ \ne\ 0\ \Longrightarrow\  
   v\big(P(a)\big)\  =\  v(P,a) + o\big(v(P,a)\big).$$ 
Note: if $K$ is monotone, then $v(P,a)=v\big(P(a)\big)$ by Corollary~\ref{monvp}.

\begin{lemma}\label{newtona}
Suppose $K$ is $r$-$\d$-henselian. 
Then there is $b\in \mathcal{O}$ with $P(b)=0$ 
and $v(a-b)\ge v(P,a)$; any such $b$ satisfies
$v(a-b)=v(P,a)$.
\end{lemma} 
\begin{proof} If $P(a)=0$, then we must take $b=a$. Assume $P(a)\ne 0$, set
$\gamma:= v(P,a)\in \Gamma^{>}$, take $g\in K$ with 
$vg=\gamma$.
Then $P(a+gY)\ =\ Q_{\times g}$, and
\begin{align*} Q_{\times g}\ &=\ P(a) + Q_{1,\times g} + \sum_{j\ge 2} Q_{j,\times g}, \qquad v\big(Q_{1,\times g}\big)\ =\ v\big(P(a)\big),  \\
  v\big(Q_{j,\times g}\big)\ &\ge\  j\gamma + o(\gamma)\ \text{ for $j\ge 2$,  so}\\
P(a+gY)\ &=\ P(a)\cdot\left( 1+\sum_{i=0}^r b_iY^{(i)} + R(Y)\right), \quad \val R\ge 2,\\
b_0,\dots, b_r&\in \mathcal{O},\  b_i \asymp 1\  \text{ for some $i$,} \qquad
v(R)\ge \gamma + o(\gamma) >0. 
\end{align*}
Now $K$ is
$r$-$\d$-henselian, so Lemma~\ref{difhen1} provides
$y\in \mathcal O$ such that $$1+\sum_{i=0}^r b_iy^{(i)} + R(y)=0.$$
Any such $y$ satisfies $y\asymp 1$, and so for $b:=a+gy$ we have $P(b)=0$ and
$v(a-b)=v(P,a)$.
\end{proof}

\noindent
Under a weaker assumption on $K$ we can still draw a useful conclusion:

\begin{lemma}\label{newton}
Suppose $\k$ is $r$-linearly surjective and $P(a)\ne 0$.
Then there is~${b\in \mathcal{O}}$ such that $P$ is in
dh-position at $b$, $v(a-b)\ge v(P,a)$ and 
$P(b)\prec P(a)$; any such $b$ satisfies
$v(a-b)=v(P,a)$ and $v(P,b)> v(P,a)$.
\end{lemma}
\begin{proof} At the point in the proof of Lemma~\ref{newtona} where we invoke
that $K$ is $r$-$\d$-henselian, we choose instead
$y\in \mathcal{O}$ such that 
$1+\sum_{i=0}^r b_iy^{(i)} \prec 1$.
Then $y \asymp 1$, so $v(b - a)=v(P,a)$ and $P(b) \prec P(a)$ with
$b:= a+gy$. 
To show that $P$ is in dh-position at $b$, we use Taylor expansion. Let
$\i=(i_0,\dots,i_r)$ and $\j$ range over~$\mathbb{N}^{1+r}$, recall that
$P_{(\i)}= \frac{P^{(\i)}}{\i!}\in K\{Y\}$, and
$$P(a+Y)\ =\ P(a) + \sum_{|\i|\ge 1} P_{(\i)}(a)Y^{\i}, \qquad 
P_{+a,1}=\sum_{|\i|=1} P_{(\i)}(a)Y^{(\i)} .$$ 
Taylor expanding $P_{(\i)}$ at $a$ gives
$$P_{(\i)}(b)\ =\ P_{(\i)}(a)+
\sum\limits_{|\j|\geq 1}P_{(\i)(\j)}(a)\cdot (gy)^{\j}.$$
Since $ P_{(\i)(\j)}(a)=
{{\i+\j}\choose {\j}}P_{(\i+\j)}(a)$, this gives for $|\i|\ge 1$:
\begin{align*} P_{(\i)}(b)\ &\sim\ P_{(\i)}(a) \hskip-5em && 
\text{ if  $P_{(\i)}(a)\asymp 1$,}\\ 
   P_{(\i)}(b)\ &\prec\ 1 \hskip-5em &&
\text{ if  $P_{(\i)}(a)\prec 1$.}
\end{align*}
Thus $P$ is in dh-position at $b$.
It only remains to show $v(P,b)> v(P,a)$. The same way~$(P,a)$ gives
rise to $A=L_{P_{+a}}\in K[\der]$, the pair $(P,b)$ yields 
$B=L_{P_{+b}}\in K[\der]$, 
and the arguments above show that $B=A+E$
with $v(E)\ge \gamma + o(\gamma)$ where $\gamma:=v(P,a)$.
In combination with $v_A(\gamma) = \gamma + o(\gamma)$, this gives
$$v_E(\gamma)\ =\ v(E)+ \gamma+ o(\gamma)\ \ge\ 
2\gamma + o(\gamma)\ >\ v_A(\gamma),$$
Hence $v_B(\gamma) = v_A(\gamma)$, and so $P(b) \prec P(a)$ forces
$v(P,b) > v(P,a)$.  
\end{proof}

\noindent 
Without even assuming that $\k$ is $r$-linearly surjective, the arguments 
above and Corollary~\ref{cor:vA=vB} give something that will be useful later:

\begin{samepage}

\begin{lemma}\label{henseld1} Suppose $P(a)\neq 0$,
$b\in \mathcal{O}$, $v(a-b)\ge v(P,a)$ and $P(b) \prec P(a)$. Then
$v(a-b)=v(P,a)$, and for all $b^*\in  \mathcal{O}$ with
$v(b-b^*)>v(P,a)$ and $B^*:=L_{P_{+b^*}}$,
\begin{enumerate}
\item[\textup{(i)}] $P$ is in dh-position at $b^*$; 
\item[\textup{(ii)}] $P(b^*) \prec P(a)$ and $v(P,b^*) > v(P,a)$;
\item[\textup{(iii)}] for all $y\in K^\times$, if $vy=O\big(v(P,a)\big)$ and
$A(y)\prec Ay$, then $B^*(y)\prec B^*y$;
\item[\textup{(iv)}] $\big\{ \alpha\in\exc(A):\ 
\alpha=O\big(v(P,a)\big)\big\}\ \subseteq\  
\exc(B^*)$.
\end{enumerate}
\end{lemma}
\end{samepage}

\begin{lemma}\label{hensel.imm} Let $\k$ be $r$-linearly surjective, 
and suppose there is no $b\in K$ with $P(b)=0$ and 
$v(a-b) = v(P,a)$.
Then there exists a divergent pc-sequence $(a_\rho)$ in~$K$ such that
$P(a_\rho) \leadsto 0$.  
\end{lemma}  
\begin{proof} Let $(a_\rho)_{\rho<\lambda}$ be a sequence 
in $\mathcal{O}$ with $\lambda$ an ordinal $>0$, $a_0=a$, and
\begin{enumerate}
  \item $P$ is in dh-position at $a_\rho$, for all $\rho < \lambda$;
  \item $v(a_{\rho'} -a_\rho)=v(P,a_\rho)$ 
whenever $\rho<\rho'<\lambda$; and
  \item $P(a_{\rho'})\prec P(a_\rho)$ and 
$v(P,a_{\rho'})>v(P,a_\rho)$ whenever $\rho<\rho'<\lambda$.
\end{enumerate}
Note that there is such a sequence if $\lambda=1$. 
Suppose $\lambda= \mu +1$ is a successor ordinal. Then Lemma~\ref{newton} 
yields $a_\lambda\in K$ such that $v(a_\lambda - a_\mu)=v(P,a_\mu)$,
$P(a_\lambda)\prec P(a_\mu)$ and 
$v(P,a_\lambda)>v(P,a_\mu)$. Then the extended
sequence $(a_\rho)_{\rho<\lambda +1}$ has the above properties 
with $\lambda+1$ instead of $\lambda$. 

Suppose $\lambda$ is a limit ordinal. Then $(a_\rho)$ is
a pc-sequence and 
$P(a_\rho) \leadsto 0$. If~$(a_\rho)$ has no pseudolimit in $K$ 
we are done. Assume otherwise, and take
a pseudolimit $a_\lambda\in K$ of $(a_\rho)$. The extended 
sequence $(a_\rho)_{\rho<\lambda +1}$ clearly satisfies condition~(2) with~$\lambda+1$ instead of $\lambda$. Applying Lemma~\ref{henseld1} to
$a_\rho$, $a_{\rho+1}$ and $a_\lambda$ in the place of $a$, $b$ and~$b^*$, where
$\rho < \lambda$, we see that conditions (1) and (3) are also satisfied 
with $\lambda+1$ instead of~$\lambda$.
This building process must come to an end. 
%by producing a divergent pc-sequence $(a_\rho)$. 
\end{proof}

\begin{cor}\label{henseld3} If $\k$ is $r$-linearly surjective
and $K$ is spherically complete, then~$K$ is $r$-$\d$-henselian.
\end{cor}

\begin{cor}\label{henimm} If $\k$ is $r$-linearly surjective, then $K$ has an immediate $r$-$\d$-hen\-se\-lian extension. If 
$\k$ is linearly surjective, then $K$ has an immediate $\d$-hen\-sel\-ian extension.
\end{cor}
\begin{proof} For $r=0$, take the henselization. For $r\ge 1$, 
use Corollary~\ref{cor2diftrzda}.  
\end{proof}

\begin{proof}[Proof of Theorem~\ref{thm:damdh}]
Assume $\k$ is linearly surjective, $K$ is $\d$-algebraically maximal,
and suppose towards a contradiction that
 $P\in \mathcal{O}\{Y\}$, $P_0\prec 1$, $P_1\asymp 1$, and there is no
$b\in \smallo$ with $P(b)=0$. Then Lemma~\ref{hensel.imm} provides
a divergent pc-sequence $(a_{\rho})$ in $K$ with $P(a_{\rho})\leadsto 0$.
Thus $(a_{\rho})$ is of $\d$-algebraic type over $K$, and so 
Lemma~\ref{zda} yields a proper immediate $\d$-algebraic 
extension of $K$.
\end{proof}

\subsection*{Step-completeness and differential-henselianity}
The spherical completeness in Corollary~\ref{henseld3} can be replaced by
step-completeness, at the cost of a stronger linear hypothesis:
Call $\mathcal{O}$ {\bf $r$-linearly surjective} if
for all $a_0,\dots,a_r\in \mathcal{O}$ such that $a_i\asymp 1$ for some $i$, 
the inhomogeneous linear differential equation 
$$1+a_0y + \dots + a_ry^{(r)}=0$$ has a solution in $\mathcal{O}$. 
It is easy to check that $\mathcal{O}$ is $r$-linearly surjective
iff each $A\in K[\der]^{\ne}$ of order $\le r$ is neatly surjective.
Thus by Lemma~\ref{dhns}, if $K$ is
$r$-$\d$-henselian, then~$\mathcal{O}$ is  
$r$-linearly surjective. 
Call $\mathcal{O}$ {\bf linearly surjective} if it is 
$r$-linearly surjective for each~$r$. 

\index{linearly!surjective!valuation ring}

\begin{lemma}\label{newtona, linear}
Suppose $\mathcal O$ is $r$-linearly surjective. Let $P\in K\{Y\}$ have  order~$\leq r$
with $P\asymp 1$ and $\deg P=\ddeg P=1$, and let $a\in\mathcal O$ be such that $P(a)\prec 1$. Then $P$ is in dh-position at $a$, there is $b\in\mathcal O$ with $P(b)=0$ and 
$v(a-b)\geq v(P,a)$, and any such $b$ satisfies $v(a-b)=v(P,a)$.
\end{lemma}
\begin{proof}
By Lemma~\ref{v-under-conjugation}(i) and Lemma~\ref{dn1,new}(i) we see that $P$ is in dh-position at~$a$. Now argue as in the proof of 
Lemma~\ref{newtona}, using the fact that $R=0$ and the hypothesis that $\mathcal O$ is $r$-linearly surjective
in place of Lemma~\ref{difhen1}.
\end{proof}

\begin{lemma}\label{newtonsp}
Suppose $\mathcal{O}$ is $r$-linearly surjective,
$P\in \mathcal O\{Y\}$ of order $\le r$ is in dh-position at $a\in \mathcal O$, and $P(a)\ne 0$.
Then there is $b\in \mathcal{O}$ such that $P$ is in
dh-position at~$b$, $v(a-b)= v(P,a)$, and 
$v\big(P(b)\big)\ \ge\ 2v\big(P(a)\big) + o\big(v(P(a))\big)$. 
\end{lemma}
\begin{proof} 
Put $\gamma:= v(P,a)$. We follow the proof
of Lemma~\ref{newtona}. There we took $y\in \mathcal{O}$
such that $1+\sum_{i=0}^r b_iy^{(i)} + R(y) =0$, and here we take  
$y\in \mathcal{O}$ such that $1+\sum_{i=0}^r b_iy^{(i)} =0$. Then
$b:= a+gy$ gives $P(b)=P(a)R(y)$, so
\begin{align*}
v\big(P(b)\big)\ =\  v\big(P(a)\big) + v\big(R(y)\big)\  &\ge\ v\big(P(a)\big)+\gamma+o(\gamma)\\ 
\ &=\ 2v\big(P(a)\big) + o\big(v(P(a))\big). \qedhere
\end{align*}
\end{proof}

\noindent
Suppose $\mathcal O$, $P$, $a$, $b$ are as in Lemma~\ref{newtonsp}. Then
$$v(P,b)\ \ge\  2 v(P,a) + o(v(P,a)).$$ 
This is because $v(P,b)= v\big(P(b)\big) + o\big(v(P(b))\big)$ if $P(b)\ne 0$.

\begin{lemma}\label{hensel.immK} Suppose $\mathcal{O}$ is 
$r$-linearly surjective, 
and $P(Y)\in \mathcal{O}\{Y\}$ of order~${\le r}$ is in dh-position at 
$a\in \mathcal{O}$.
Suppose that there is no $b\in K$ with $P(b)=0$ and 
$v(a-b) = v(P,a)$.
Then there exists a divergent special pc-sequence $(a_\rho)$ in $K$ such that
 $P(a_\rho) \leadsto 0$.
 \end{lemma} 
\begin{proof}{\sloppy We follow the proof of Lemma~\ref{hensel.imm}
except that in condition~(3) on~$(a_\rho)_{\rho<\lambda}$
the second inequality is replaced by the stronger
$$ v(P,a_{\rho'})\ >\ (3/2)v(P,a_\rho)\ \text{  whenever $\rho<\rho'<\lambda$,}$$
the appeal to Lemma~\ref{newton} is replaced by an appeal to
 Lemma~\ref{newtonsp}, which yields $v(P,a_\lambda)>(3/2)v(P,a_\mu)$
instead of $v(P,a_\lambda)>v(P,a_\mu)$. Also, when  
$\lambda$ is a limit ordinal, $(a_\rho)$ is a {\em special}\/
pc-sequence.}
\end{proof} 

\begin{cor}\label{henseld4special} If $\mathcal{O}$ is 
$r$-linearly surjective and each special pc-sequence in~$K$ with a minimal differential polynomial over $K$ of order $\le r$ pseudoconverges in~$K$, then~$K$ is $r$-$\d$-henselian.
\end{cor} 

\begin{cor}\label{henseld4} If $\mathcal{O}$ is 
linearly surjective and $K$ is step-complete, then $K$ is
$\d$-henselian.
\end{cor}

\subsection*{Lifting zeros of linear differential operators} 
Let $A\in K[\der]$ have or\-der~$\le r$, with $vA=0$. Let $\bar{A}$ be the image of $A$ under the natural
map $\mathcal{O}[\der] \to \k[\der]$.

\begin{lemma}\label{dhres} Suppose $\mathcal{O}$ is $r$-linearly surjective. % $n\ge 1$.
Then the 
additive map
$$ a \mapsto \overline{a} \colon\ \mathcal{O}\cap \ker A\ \to\ \ker \bar{A}$$
is surjective. 
\end{lemma}
\begin{proof} 
Suppose $a\in \mathcal{O}$ and 
$\overline{A}(\overline{a})=0$. Then $A(a)\prec 1$, so 
we have $y\prec 1$ in $K$ with $A(y)=A(a)$. Then
$A(a-y)=0$ and $\overline{a-y}=\overline{a}$. 
\end{proof}  

\noindent
If  $\mathcal O$ is $r$-linearly surjective, then
$\dim_{C_{\k}} \ker \overline{A} \leq  \dim_C \ker A$ by 
Lemma~\ref{dhres}. 
Under an additional condition on $K$ we have equality:

\begin{prop}\label{prop:equality of dimensions}
Suppose $\mathcal O$ is $r$-linearly surjective and $K$ has many constants $($that is, $v(C^\times)=\Gamma)$. Then $\dim_{C_{\k}} \ker\overline{A}=\dim_C \ker A$.
\end{prop}

\begin{proof} {\sloppy
By Lemma~\ref{dhres} we can
take $f_1,\dots,f_m\in\mathcal O\cap\ker A$  whose re\-si\-dues~$\overline{f_1},\dots,\overline{f_m}$ form a basis of the $C_{\k}$-linear space $\ker \overline{A}$. Let $V$ be the $C$-linear subspace of $\ker A$ generated by $f_1,\dots,f_m$ and consider the  homogeneous differential polynomial
$$P(Y)\ :=\ \wr(Y,f_1,\dots,f_m)\in\mathcal O\{Y\}$$
of degree $1$ and order~$\leq m$.
Then by Lemma~\ref{lem:wronskian} we have 
$$V\ =\ \big\{y\in K:\ P(y)=0\big\}.$$
Moreover
$$P(Y)\ =\ \sum_{i=0}^m (-1)^i (\det W_i)\, Y^{(i)},$$
where $W_i$ is the $m\times m$ matrix whose $j$th column is 
$$\big(f_j,\dots,f_j^{(i-1)},f_j^{(i+1)},\dots,f_j^{(m)}\big)^{\operatorname{t}};$$ 
in particular, $W_m=\Wr(f_1,\dots,f_m)$. Since $\wr(\overline{f_1},\dots,\overline{f_m})\neq 0$, this gives~${\deg D_P=1}$
and $\order D_P=\order P=m\le r$. 
For all $y\in\mathcal O$,
$$\text{$P$ is in dh-position at $y$}\Longleftrightarrow P(y)\prec 1  \Longleftrightarrow  \wr(\overline{y},\overline{f_1},\dots,\overline{f_m})=0
\Longleftrightarrow \overline{A}(\overline{y})=0.$$
It is enough to show $V=\ker A$. Suppose
towards a contradiction that $V\neq\ker A$. Take $a\in \ker A$, $a\notin V$. Since $K$ has many constants we can multiply $a$ by a nonzero constant
to arrange $a\asymp 1$. From $\overline{A}(\overline{a})=0$
we get $0\ne P(a)\prec 1$, so $\gamma:= v(P,a)\in\Gamma^>$.
By Lemma~\ref{newtona, linear} we can take $b\in V$
with $v(a-b)=\gamma$. Next, take $c\in C$ such that $vc=-\gamma$ and put $g:=c(a-b)$; then $g\asymp 1$, 
$\overline{A}(\overline{g})=0$, so
Lemma~\ref{newtona, linear} gives $h\in V$ with $v(g-h)=v(P,g)>0$.
Put $b^*:=b+c^{-1}h\in V$. Then $a-b^* = (a-b)-c^{-1}h=c^{-1}(g-h)$, so
$$v(a-b^*)\ =\ \gamma+v(g-h)\ >\ \gamma\ =\ v(P,a),$$ 
in contradiction to Lemma~\ref{newtona, linear}.}
\end{proof}

\subsection*{Differential-henselianity and completion} By Corollary~\ref{completionsmall} the completion~$K^{\c}$ of $K$ is an immediate extension of $K$. Our aim in this subsection is:

\begin{prop}\label{dhenscompletion} If $K$ is $\d$-henselian and $\operatorname{cf}(\Gamma)=\omega$, then $K^{\c}$ is $\d$-henselian.
\end{prop}

\noindent
Let $\mathcal{O}^{\c}$ be the valuation ring of $K^{\c}$. The proof uses a variant of Lemma~\ref{newton}:

\begin{lemma}\label{variantnewton} Assume $K$ is $\d$-henselian and 
$P\in \mathcal{O}^{\c}\{Y\}$
is in dh-position at the point $a\in \mathcal{O}^{\c}$, with $P(a)\ne 0$.
Let $\alpha\in \Gamma^{>}$ be given. Then there is
$b\in \mathcal{O}^{\c}$ such that $P$ is in dh-position at $b$,
$v(a-b)=v(P,a)$ and $P(b)\prec P(a)$, $vP(b) \ge \alpha$, 
$v(P,b) > v(P,a)$.
\end{lemma}
\begin{proof} Take $g\in K^{\c}$ with $vg=\gamma:= v(P,a)$ as in the proof of Lemma~\ref{newtona}. As in that proof we get
\begin{align*}P(a+gY)\ &=\ P(a)\cdot\left( 1+\sum_{i=0}^r b_iY^{(i)} + R(Y)\right), \quad R\in \mathcal{O}^{\c}\{Y\},\ \val R\ge 2,\\
b_0,\dots, b_r&\in \mathcal{O}^{\c},\  b_i \asymp 1\  \text{ for some $i$,} \qquad
v(R)\ge \gamma + o(\gamma) >0. 
\end{align*}
As $K$ is $\d$-henselian and dense in $K^{\c}$ we have
$y\in \mathcal O$ such that $$v\left(1+\sum_{i=0}^r b_iy^{(i)} + R(y)\right)\ \ge\ \alpha.$$
Any such $y$ satisfies $y\asymp 1$ and so for $b:=a+gy$ we have
$P(b)\prec P(a)$ and $vP(b)\ge \alpha$ and $v(a-b)=v(P,a)>0$. Thus
$P$ is in dh-position at $b$, and so $v(P,b) > v(P,a)$ by Lemma~\ref{newton}. 
\end{proof} 

\begin{proof}[Proof of Proposition~\ref{dhenscompletion}]
Assume $K$ is $\d$-henselian and $\operatorname{cf}(\Gamma)=\omega$.
 Let $P\in \mathcal{O}^{\c}\{Y\}$ be given such that $0\ne P(0)\prec 1$ and
$P_1\asymp 1$; it suffices to show that then~$P$ has a zero in the maximal ideal  $\smallo^{\c}$ of $\mathcal{O}^{\c}$. 
Fix a strictly increasing cofinal sequence~$(\alpha_{n})$ in $\Gamma^{>}$ with $\alpha_0=vP(0)$. We proceed as in the proof of Lemma~\ref{hensel.immK}, with some differences: Let $(a_m)_{m<n}$ be a finite sequence 
in $\smallo^{\c}$ with $n\ge 1$ and
\begin{enumerate}
  \item $a_0=0$, and $P$ is in dh-position at $a_m$ for all $m < n$,
  \item $v(a_{m'} -a_m)=v(P,a_m)$ 
whenever $m<m'<n$,
  \item $P(a_{m'})\prec P(a_m)$, 
$v(P,a_{m'})>v(P,a_m)$ whenever $m<m'<n$,
\item $vP(a_{m})\ge  \alpha_{m}$ whenever $m<n$.
\end{enumerate}
Note that we have such a sequence for $n=1$. 
Let $n= \mu +1$. If $P(a_{\mu})=0$, we are done, so assume $P(a_{\mu})\ne 0$. Then Lemma~\ref{variantnewton} 
yields $a_n\in K$ such that 
$vP(a_{n})\ge \alpha_{n}$, 
$v(a_n - a_\mu)=v(P,a_\mu)$,
$P(a_n)\prec P(a_\mu)$ and 
$v(P,a_n)>v(P,a_\mu)$. Then the extended
sequence $(a_m)_{m<n +1}$ has the above properties 
with $n+1$ instead of~$n$. Iterating this extension procedure yields an
infinite sequence $(a_n)$ in $\smallo^c$. This is a c-sequence, and
so $a_n \to a\in \smallo^{\c}$. Then $P(a)=0$.
\end{proof}

\subsection*{Notes and comments}
Proposition~\ref{prop:equality of dimensions} is a variant of \cite[Theorem~4.1.6]{MAnderson}. 

\section{Differential-Hensel Configurations} 

\noindent
{\em In this section
we assume that $\Gamma \ne \{0\}$.}\/ The notion of {\em
dh-position}\/ is too closely tied to differential polynomials 
of a special form, so we relax it as follows.
Let $P\in K\{Y\}$ have order $\le r$, and $a\in K$. Then
$$ P(a+Y)=P(a) + A(Y) +R(Y), \qquad A,R\in K \{Y\},$$ 
where $A= \sum_{i=0}^r P_{(i)}(a)Y^{(i)}$ 
is homogeneous of degree $1$, and all terms in $R$ have degree $\ge 2$.
Let us say that $P$ is in {\bf differential-hensel configuration} 
(abbreviated as {\bf dh-configuration}) at $a$ if $A\ne 0$, and
there is $\gamma\in \Gamma$ such that
$v(P(a))\ge v_A(\gamma) < v_R(\gamma)$; equivalently, $A\ne 0$, 
and either $P(a)=0$ or there is $\gamma\in \Gamma$ such that
$v(P(a))=v_A(\gamma) < v_R(\gamma)$. 
To prove this equivalence, assume $v(P(a)) > v_A(\gamma) < v_R(\gamma)$, $\gamma\in \Gamma$. Increasing $\gamma$ to
$\gamma+\delta$ ($\delta\in \Gamma^{>}$) such that $v(P(a))=v_A(\gamma+ \delta)$
we note that $v_A(\gamma+\delta)=v_A(\gamma) + \delta +o(\delta)$
and $v_R(\gamma+\delta) \ge v_R(\gamma) + 2\delta + o(\delta)$, by Corollary~\ref{vP-Lemma},
so $v_A(\gamma+\delta) < v_R(\gamma+\delta)$.

For any extension $L$ of $K$, the Equalizer Theorem yields: $P$ is
in dh-con\-fi\-gu\-ra\-tion at~$a$ with respect to~$K$ iff $P$ is
in dh-con\-fi\-gu\-ra\-tion at~$a$ with respect to~$L$.

 For differential polynomials of order $0$, ``differential-hensel configuration'' agrees with ``hensel configuration'' as defined in Section~\ref{sec:henselian valued fields}.

\index{configuration!differential-hensel}
\index{dh-configuration}
\index{differential-hensel!configuration}

\medskip
\noindent
For $P$ in dh-configuration at $a$ we define $v(P,a)\in \Gamma_{\infty}$ as follows: if $P(a)\ne 0$, then $v(P,a)$
is the unique $\gamma\in \Gamma$ with 
$v\big(P(a)\big) = v_A(\gamma)$, and if $P(a)=0$, then
$v(P,a):=\infty$.
Note: if $P\in \mathcal O\{Y\}$ and $a\in \mathcal O$, then 
$$\text{$P$ is in dh-position at $a$}\ \Longrightarrow\  \text{$P$ is in dh-configuration at $a$.}$$
Suppose $P$ is in dh-configuration at $a$ with $P(a)\ne 0$. 
Take $g\in K^\times$ with $vg=v(P,a)$, and put
$G(Y):= P(a+gY)/P(a)$, so $G(Y)=1+B(Y) + S(Y)$ where $B\in K\{Y\}$ is homogeneous of
degree $1$ with $v(B)=0$ and all terms in $S\in K\{Y\}$ have degree~${\ge 2}$ with
$v(S)>0$. Assuming now that $\k$ is $r$-linearly surjective, we can take
$y\in K$ with $y\asymp 1$ and $1+B(y)\prec 1$. Then
$G$ is in dh-position at $y$. If $K$ is
$r$-$\d$-henselian we can take $y$ as above such that $G(y)=0$, 
and then
$b:=a+gy$ satisfies $P(b)=0$ and $v(a-b)=v(P,a)$. So we have shown:

\begin{lemma}\label{dh} If $K$ is  
$r$-$\d$-henselian and $P\in K\{Y\}$ of order~$\le r$ is in dh-con\-fi\-gu\-ra\-tion at $a\in K$, then $P(b)=0$ and $v(a-b)=v(P,a)$ for some $b\in K$.
\end{lemma} 

\begin{lemma}\label{dhddeg} Let $a\in K$ and $P\in K\{Y\}$. Then
$$ P \text{ is in dh-configuration at $a$}\ \Longleftrightarrow\  \operatorname{ddeg} P_{+a,\times g}=1\ \text{ for some $g\in K^{\times}$.}$$
\end{lemma}
%\begin{proof} We have $P_{+a}=P(a) + A + R$ 
%with $A, R\in K\{Y\}$,
%$A$ homogeneous of degree $1$, and all terms in $R$ 
%of degree $\ge 2$,
%so for $g\in K^\times$ we have
%$$P_{+a,\times g}\ =\  P(a) + A_{\times g} + R_{\times g}.$$ 
%Suppose $P$ is in dh-configuration at $a$. Then
%$A\ne 0$, and we can take $g\in K^{\times}$ such that
%$P(a) \preceq A_{\times g} \succ R_{\times g}$; 
%for such $g$ we have $\operatorname{ddeg} P_{+a,\times g}=1$.
%In particular, if $P(a)\ne 0$ and 
%$g\in K^{\times},\ vg=v(P,a)$, then
%$\operatorname{ddeg} P_{+a,\times g}=1$ and $P(a)\asymp %A_{\times g}$.

%Conversely, if $g\in K^{\times}$ and
%$\operatorname{ddeg} P_{+a,\times g}=1$, then $P(a) \preceq %A_{\times g} \succ R_{\times g}$, so $P$ is in dh-%configuration at $a$.
%\end{proof}

\noindent
The proof is straightforward and left to the reader.

\begin{lemma}\label{dhli4-} 
Let $P\in K\{Y\}$ be in dh-configuration at $a\in K$ and $P(a)\ne 0$.
Then~$P$ is in dh-configuration at $b$ for every $b\in K$ with
$v(a-b)\ge v(P,a)$.
\end{lemma}
\begin{proof} By passing to the algebraic closure of $K$ we arrange that
$\Gamma$ has no least positive element. Take $g\in K^\times$ with $vg=v(P,a)$. Then 
$$Q\ :=\ P_{+a,\times g}\ =\ P(a) + A + R, \quad
P(a)\asymp A\succ R$$ where $A\in K\{Y\}$ is homogeneous of degree $1$
and $\val R \ge 2$. Since $\Gamma$ has no least positive element, we can take $h\in K$ such that $h\succ 1$ and 
$A_{\times h}\succ R_{\times h}$. 
Let $b\in K$ be such
that $v(a-b) > v(gh)$. Then $b=a+ghy$ with $y\prec 1$, so
$$P_{+b,\times gh}\ =\ P_{+a, \times gh, +y}\ \sim\ P_{+a, \times gh}\ =\ Q_{\times h},$$
and thus $\ddeg P_{+b,\times gh}=1$.
\end{proof}

\noindent
If $P$ is in dh-configuration at $a$ (with $P\in K\{Y\}$, $a\in K$),
and $\phi\in K^\times$, then $P_{\times \phi}$ is in dh-configuration at $a/\phi$,
with $v(P_{\times \phi}, a/\phi)=v(P,a)-v\phi$.

\subsection*{Compositional conjugation}
Let $\phi\in K$ and $\phi\asymp 1$. Then the derivation of~$K^{\phi}$ is also 
small, and the residue differential
field of $K^{\phi}$ is $\k^{\overline{\phi}}$, where $\overline{\phi}$ is the 
residue class
$\phi + \smallo$ in $\k=\mathcal O/\smallo$. In particular,
$\k$ is $r$-linearly surjective iff $\k^{\overline{\phi}}$ is $r$-linearly
surjective. For $P\in K\{Y\}$ we have
$(P^{\phi})_d=(P_d)^{\phi}$ for all $d\in \N$, and $v(P^\phi)=v(P)$. Thus 
$\mathcal O$ is $r$-linearly surjective iff 
$\mathcal O{}^\phi$ is $r$-linearly surjective, and if
$K$ is $r$-$\d$-henselian, then so is $K^{\phi}$.

\subsection*{Coarsening and specialization} In this subsection $\Delta$ is a convex 
subgroup of $\Gamma$. Let $\dot{\mathcal O}$ be the valuation ring of the 
corresponding coarsening $\dot{v}\colon K^\times \to  \dot{\Gamma}=\Gamma/\Delta$.

\begin{lemma}\label{inva3} If $K$ is $r$-$\d$-henselian, then so is $(K, \dot{\mathcal{O}})$.
\end{lemma}
\begin{proof} 
%The case $\Gamma=\{0\}$ is trivial, so let 
Recall that $\Gamma\ne \{0\}$. Let $P\in \dot{\mathcal O}\{Y\}$ and $a\in \dot{\mathcal O}$, and suppose
$P$ is in dh-position at $a$ with respect to the valuation $\dot{v}$.

\claim{$P$ is in dh-configuration at $a$ with
respect to the original valuation $v$.}

\noindent
To see why, let $P(a+Y) = P(a) + L(Y) + R(Y)$ where  $\dot{v}(P(a))>0$,
$L\in \dot{\mathcal O}\{Y\}$ is homogeneous of degree $1$ with $\dot{v}(L)=0$,  
and all terms of $R\in \dot{\mathcal O}\{Y\}$ have degree~${\ge 2}$.
We can assume $P(a)\ne 0$. Take $\gamma\in \Gamma$ such that
$v\big(P(a)\big)= v_L(\gamma)$. Then  $\dot{v}\big(P(a)\big)= \dot{v}_L(\dot{\gamma})$
by Lemma~\ref{coarsevp}, so 
$\gamma > \Delta$. Then $v_L(\gamma)=v(L)+ \gamma + o(\gamma)$,
and $v_R(\gamma) \ge  v(R) + 2\gamma + o(\gamma)$ by Corollary~\ref{vP-Lemma}.  
Since $v(L)\in \Delta$ and $v(R) \ge \delta$ for some $\delta\in \Delta$,
it follows that $v_L(\gamma) <  v_R(\gamma)$. This proves the claim, and also gives $v(P,a)> \Delta$. Assume now that $K$ is 
$r$-$\d$-henselian. Then the claim and Lemma~\ref{dh} yield 
$b\in K$
with $P(b)=0$ and $v(a-b)>\Delta$, so $\dot{v}(a-b)>0$.
It remains to note that $\dot{K}$ is $r$-linearly surjective by Lemma~\ref{inva4}. 
\end{proof}

\noindent
In combination with Lemmas~\ref{inva4} and~\ref{inva5}, we obtain:

\begin{cor}
$K$ is $r$-$\d$-henselian if and only if both
$(K,\dot{\mathcal O})$ and $\dot{K}$ are $r$-$\d$-henselian. 
\end{cor}

\section{Maximal Immediate Extensions in the Monotone Case}\label{umie}

\noindent
{\em In this section $K$ is a monotone valued differential
field, the induced derivation 
on~$\k$ is nontrivial, and $\Gamma\ne \{0\}$.}\/ Note that any
extension~$L$ of~$K$ with $\Gamma_L=\Gamma$ is monotone, by Corollary~\ref{cor:sc2}.
%by Lemma~\ref{stas11}. 

We begin with proving the key result that will give uniqueness of maximal
immediate extensions when the differential residue field is linearly surjective.

\subsection*{The differential-henselian configuration theorem}
Proposition~\ref{dhconf} below is analogous to Proposition~\ref{prop:hensel config} and plays a similar critical~role.
% as Corollary~\ref{cor:hensel config} 
%in certain uniqueness issues.

Let $L$ be a monotone extension of $K$, and $(a_\rho)$ a pc-sequence in $K$ such that $a_{\rho} \leadsto a\in L$. We set $\gamma_{\rho}:= v(a_{s(\rho)}-a_\rho)$, and 
let $\alpha$,~$\beta$,~$\gamma$ range over $\Gamma$.

\begin{prop} \label{dhconf}
Let $G(Y)\in L\{Y\}\setminus L$ have order $\le r$ and suppose \begin{enumerate}
\item[\textup{(i)}]  $G(a_\rho) \leadsto 0$;
\item[\textup{(ii)}] for every $\i\in \N^{1+r}$ with $|\i|=1$ and every pc-sequence 
$(e_{\lambda})$ in $K$ equivalent to~$(a_\rho)$, we have
$G_{(\i)}(e_\lambda)\not \leadsto 0$.
\end{enumerate}
Then $G$ is in dh-configuration at $a$, and $v(G,a)> v(a-a_{\rho})$, eventually.
There is also an index $\rho_0$ such that for all $\rho> \rho_0$ and all 
$g\in L$,
$$ g\asymp a-a_{\rho}\ \Longrightarrow\ 
\dval G_{+a,\times g}\ =\ 
\operatorname{ddeg} G_{+a,\times g}\ =\ 1.$$
\end{prop}
\begin{proof} We let 
$\i$ range over the elements of $\N^{1+r}$ with
$|\i|=1$, and $\j$ over all elements in $\N^{1+r}$. 
By removing some initial $\rho$'s we can assume
$\gamma_{\rho}=v(a-a_\rho)\in \Gamma^{\ne}$ for all $\rho$, and 
$\gamma_{\rho'} > \gamma_{\rho}$ whenever $\rho'>\rho$. 
Next, set  
$$ g_{\j}:= G_{(\j)}(a), \quad P(Y)\ :=\  G(a+Y)-G(a)\ =\  \sum_{|\j|\ge 1} g_{\j}Y^{\j}\ =\ \sum_{e=1}^N P_e(Y)\in L\{Y\},$$
where $N:= \deg P \ge 1$.
Lemma~\ref{pcc3} and the remarks after Lemma~\ref{pcc1} provide a pc-sequence $(b_{\rho})$ in $K$ that is equivalent to $(a_{\rho})$ such that: \begin{enumerate}
\item[(1)] $G(b_{\rho}) \leadsto G(a), \quad G(b_{\rho}) \leadsto 0$,
\item[(2)] $G_{(\i)}(b_{\rho}) \leadsto G_{(\i)}(a)$ whenever 
$G_{(\i)}\notin L$,
\item[(3)] $v(b_{\rho}-a)=\gamma_{\rho}$, eventually, 
\item[(4)] $v\big(P_e(b_{\rho}-a)\big)=v_{P_e}(\gamma_{\rho})$, eventually, whenever $1\le e\le N$ and $P_e\ne 0$.
\end{enumerate}
If 
$e, e'\in  \{1,\dots,N\}$, $e\ne e'$, $P_e\ne 0$, $P_{e'}\ne 0$, then either
$v_{P_e}(\gamma_{\rho}) <v_{P_{e'}}(\gamma_{\rho})$, eventually, or $v_{P_{e'}}(\gamma_{\rho}) <v_{P_e}(\gamma_{\rho})$, eventually, by Corollary~\ref{equalizer}. 
This yields a $d\in \{1,\dots,N\}$
such that, after removing some initial~$\rho$'s, we have 
$$ e\in \{1,\dots,N\}, d\ne e\ \Longrightarrow\ v_{P_d}(\gamma_{\rho}) <v_{P_e}(\gamma_{\rho}) \text{ for all $\rho$.}$$   
Now $G(b_{\rho})-G(a)=P(b_{\rho}-a)$ for all $\rho$, and $G(b_{\rho})\succ G(a)$, eventually, and
$$G(a+Y)= G(a) + A(Y) + R(Y), \quad A:= P_1, \quad R:=P_2+ \cdots + P_N.$$
Suppose that $d=1$. Then $G(b_{\rho}) \sim A(b_{\rho}-a)$, eventually, 
and so, eventually, 
$$v\big(G(a)\big) > v\big(A(b_{\rho}-a)\big)=
v_A(\gamma_{\rho}) < 
v_R(\gamma_{\rho}),$$
so $G$ is in dh-configuration at $a$, and $v(G,a)> \gamma_{\rho}=v(a-a_{\rho})$, eventually. Also the second part of the desired conclusion follows easily.

\medskip\noindent
It only remains to show that $d=1$. For each $\i$ we have
$$
G_{(\i)}(b_\rho)-G_{(\i)}(a)\ =\ 
\sum_{\j>\i}\binom{\j}{\i}G_{(\j)}(a)(b_{\rho}-a)^{\j-\i}\ =\ 
\sum_{\j>\i}\binom{\j}{\i}g_{\j}(b_{\rho}-a)^{\j-\i}.$$
If $G_{(\i)}\ne 0$, then $G_{(\i)}(b_\rho)\sim G_{(\i)}(a)=g_{\i}$, eventually, so
\begin{align*}  g_{\i}\ &\succ\
\sum_{\j>\i}\binom{\j}{\i}g_{\j}(b_{\rho}-a)^{\j-\i} ,
\text{ eventually, and thus}\\
  g_{\i}(b_{\rho}-a)^{\i}\ &\succ\ 
\sum_{\j>\i}\binom{\j}{\i}g_{\j}(b_{\rho}-a)^{\j},\ \text{ eventually.}
\end{align*}
From $\sum_{\i}\binom{\j}{\i}\ =\ j_0 + j_1+ \cdots + j_r$ we get
$$ \sum_{\i}\left(\sum_{\j>\i}\binom{\j}{\i}g_{\j}(b_{\rho}-a)^{\j}\right)\ 
=\ \sum_{e=2}^N e\left(\sum_{|\j|=e} g_{\j}(b_{\rho}-a)^{\j}\ \right)\ =\ 
\sum_{e=2}^N eP_e(b_{\rho}-a).$$
Now $G_{(\i)}\ne 0$ for at least one $\i$, and if $G_{(\i)}= 0$ and $\j>\i$, then
$G_{(\j)} = 0$, so
$$ \min_{\i} v\big(g_{\i}(b_{\rho}-a)^{\i}\big)\  <\ 
v\left(\sum_{e=2}^N eP_e(b_{\rho}-a)\right)\ =\ \min\big\{v_{P_e}(\gamma_{\rho}):\ e=2,\dots,N\big\}, $$ 
eventually. Since $A=P_1$ we have $v\big(P_1(b_{\rho}-a)\big)=v_A(\gamma_{\rho})$ for all $\rho$, so
to obtain $d=1$ it is enough to get, cofinally in $\rho$, 
$$v_A(\gamma_{\rho})\ \le\  \min_{\i} v\big(g_{\i}(b_{\rho}-a)^{\i}\big) .$$ 
To derive $d=1$ from the above, we now use that $L$ is monotone, which gives 
\begin{align*} v\big(g_{\i}(b_{\rho}-&a)^{\i}\big)\ \ge\ v(g_{\i}) + \gamma_{\rho}\ \text{ for all $\i$ and
$\rho$, and so}\\  
v_A(\gamma_{\rho})\ &=\ v(A)+\gamma_{\rho}\ =\ \big(\min_{\i} v(g_{\i})\,\big) + \gamma_{\rho}\ \le\ \min_{\i} v\big(g_{\i}(b_{\rho}-a)^{\i}\big)
\end{align*}
for all $\rho$, so $d=1$.
\end{proof}

\noindent
In this proof the assumption that $K$ and $L$ are monotone 
is used only at the end when it is explicitly invoked. The earlier part
of the proof goes through if this assumption is dropped, and so 
with the weaker assumption the proposition holds for $\deg G = 1$.

\subsection*{Uniqueness of maximal immediate extensions}

\begin{lemma}\label{damls} Let $L$ be a monotone $\d$-algebraically maximal
extension of~$K$ such that $\k_L$ is linearly surjective. Let
$(a_\rho)$ be a divergent pc-sequence in $K$ with minimal 
differential polynomial $G(Y)$ over $K$. Then $a_\rho \leadsto b$ and $G(b)=0$ for some $b\in L$.
\end{lemma}
\begin{proof} Take a pseudolimit $a\in L$ of
$(a_\rho)$, and take a pc-sequence $(b_{\lambda})$ in $K$ equivalent to
$(a_{\rho})$ such that $G(b_{\lambda})\leadsto 0$. Then $G$ is in dh-configuration at $a$, and $v(G,a) > v(a-b_\lambda)$, eventually, by Proposition~\ref{dhconf}. Now $L$ is $\d$-henselian by Theorem~\ref{thm:damdh}, so by Lemma~\ref{dh} we have $b\in L$ such that
$$ v_L(a-b)\ =\ v(G,a)\ \text{  and }\ G(b)\ =\ 0.$$ 
This gives $v_L(a-b)> v(a-b_{\lambda})$, eventually, so $b_{\lambda}\leadsto b$, 
and thus
$a_\rho \leadsto b$. 
\end{proof} 

\noindent
Together with Lemmas~\ref{sc2},~\ref{diftr}, and~\ref{zda}
this yields:

\begin{theorem}\label{unique.max.imm.ext} Suppose $\k$ is linearly surjective.
Then any two maximal immediate extensions of $K$ are isomorphic over
$K$ and $\d$-henselian. Also, any two $\d$-algebraic immediate extensions of $K$ that are
$\d$-algebraically maximal are isomorphic over $K$ and $\d$-henselian. 
\end{theorem}

\noindent
We now state minor variants of these results using the notion of saturation 
from model theory (see \ref{sec:sat}), as needed in the proof of the 
Equivalence Theorem~\ref{embed11}.
Let $|X|$ denote the cardinality of a set $X$, and let 
$\kappa$ be a cardinal.

\begin{lemma} Let $K$ be $\d$-henselian. 
Let $E$ be a valued differential subfield of~$K$ such that the derivation of
$\k_E$ is nontrivial. Assume
$K$ is $\kappa$-saturated with $\kappa>|v(E^\times)|$. If~$(a_{\rho})$ is a divergent pc-sequence 
in $E$ with minimal
differential polynomial~$G(Y)$ over~$E$, then $a_\rho \leadsto b$ and $G(b)=0$ for some $b\in K$.
\end{lemma}
\begin{proof} Let $(a_{\rho})$ be a divergent pc-sequence in $E$ with minimal
differential polynomial $G(Y)$ over $E$. By saturation we have a pseudolimit $a\in K$ of
$(a_\rho)$. Take a pc-sequence $(b_{\lambda})$ in $E$ equivalent to
$(a_{\rho})$ such that $G(b_{\lambda})\leadsto 0$. By Proposition~\ref{dhconf}, $G$
is in dh-configuration at $a$ with $v(G,a)> v(a-b_{\lambda})$, eventually.
Take $b\in K$ with
$v(a-b)=v(G,a)$ and $G(b)=0$. Then $b_{\lambda} \leadsto b$, so
$a_{\rho}\leadsto b$.
\end{proof}

\noindent
In combination with Lemmas~\ref{diftr},~\ref{zda}, and~\ref{damls}, this yields:

\begin{cor}\label{immsat} Let $K$ and $E$ be as in the previous lemma, 
and assume also that~$\k_E$ is linearly surjective and $v(E^\times)\ne \{0\}$. 
Then any maximal immediate 
extension of $E$ can be embedded in $K$ over $E$.
\end{cor}

\noindent
Recall that a valued field of equicharacteristic $0$ is henselian if and only
if it is algebraically maximal. We already established an
analogue for valued differential fields in one direction:
$\d$-algebraically maximal valued differential fields with small derivation and linearly surjective 
differential residue field are $\d$-henselian. 
The converse fails, even in the monotone case: 

\begin{example} 
Let $\k$ be a countable linearly surjective differential field
and let $\Gamma$ be a countable ordered 
abelian group, $\Gamma\ne \{0\}$. Then $\k(t^\Gamma)$ is a countable 
valued differential 
subfield of the
Hahn differential field $\k\(( t^\Gamma\)) $. The latter is $\d$-henselian,
so we can take a countable $\d$-henselian $K$ such that 
$$\k(t^\Gamma)\ \subseteq\ K\ \subseteq\  \k\(( t^\Gamma\)) \ \quad\text{(as valued differential fields).}$$
The constant field of $\k\(( t^\Gamma\)) $ is $C_{\k}\(( t^\Gamma\)) $, which is 
uncountable. Take $c\in C_{\k}\(( t^\Gamma\)) $ with $c\notin K$. Then 
$K\<c\>=K(c)$ is a proper immediate $\d$-algebraic
extension of $K$, so $K$ is not $\d$-algebraically maximal.  
\end{example} 

\noindent
In this example $K$ has many constants. With few constants, we do have a converse in the monotone case, as we shall see in the next section.

\section{The Case of Few Constants}\label{sec:fewcon}

\noindent
In this section we consider in more detail the situation where $C\subseteq \mathcal{O}$, that is, the valuation is trivial on $C$.  
Key facts are Lemma~\ref{dhli3} and Proposition~\ref{noextrazeros} below.  
%As a consequence we derive the existence of a ``differential-%henselization'' of $K$, provided $\k$ is linearly surjective and 
%$C\subseteq \mathcal{O}$. Hmm, need condition that $C$ doesn't
%increase in going to nice immediate extensions

\subsection*{Valuation properties of linear differential operators} {\em In
this subsection we fix $r\ge 1$, and assume: $K$ is $r$-$\d$-henselian,
$A\in K[\der]^{\ne}$ has or\-der~${\le r}$}. 

\begin{lemma}\label{dhli1} Suppose 
$A(1)\prec A$. Then $A(y)=0$ for some $y\sim 1$ in $K$.
\end{lemma}
\begin{proof} We have $A=a_0+ a_1\der + \cdots + a_r\der^r$
with $a_0, a_1,\dots, a_r\in K$, so $A(1)=a_0$. We can arrange that
$a_0\in \smallo$ and $a_1,\dots,a_r\in \mathcal{O}$, $a_i\asymp 1$ for some
$i\in\{1,\dots, r\}$. With $R(Z):=\Ric(A)$ we have $R(Z)\in \mathcal{O}\{Z\}$,
$R(0)=a_0\prec 1$ and 
$$R(Z)_1\ =\ a_1Z+ \cdots + a_rZ^{(r-1)}\ \asymp\ 1,$$ so
we get $z\in \smallo$ with $R(z)=0$. By Corollary~\ref{dhns2} we can
take $y\in 1+\smallo$ such that $y^\dagger=z$, and then $A(y)=0$. 
\end{proof}

\noindent
The proof shows that in Lemma~\ref{dhli1} we can relax the assumption that $K$ is $r$-$\d$-henselian to $K$ being
$(r-1)$-$\d$-henselian, with
$(1+\smallo)^\dagger=\smallo$ in case $r=1$. The same holds
for Lemma~\ref{dhli2}, Corollaries~\ref{cor:dhli2} and~\ref{pcdhli}, and Lemma~\ref{dhli3}. 

\begin{lemma}\label{dhli2} Suppose $C\subseteq \mathcal{O}$. 
Then there are no $b_0\succ b_1\succ \dots \succ b_{r}$ in $K^\times$ 
with $A(b_i)\prec Ab_i$ for $i=0,\dots,r$.
\end{lemma}
\begin{proof} Let $b\in K^\times$ and $A(b) \prec Ab$. Then
$B(1)\prec B$ for $B=Ab$, so we have $y\sim 1$ in $K$ with
$B(y)=0$, and thus $A(by)=0$ with $by\sim b$.
It remains to note that if $b_0\succ b_1\succ \dots \succ b_r$ in $K^{\times}$, then $b_0,\dots, b_r$ are $C$-linearly independent. 
\end{proof}

\noindent
Here is a reformulation of Lemma~\ref{dhli2} and its proof in terms of the set   
$$\exc(A)\ =\ \big\{ vb:\ b\in K^\times,\ A(b)\prec Ab \big\}$$  
of exceptional values for $A$ that was introduced in Section~\ref{sec:ldopv}:

\begin{cor}\label{cor:dhli2} If $C\subseteq \mathcal{O}$, then 
$\exc(A)=v(\ker^{\neq} A)$, so $\abs{\exc(A)}\leq \dim_C \ker A$.
\end{cor}

\begin{cor}\label{pcdhli} Suppose $C\subseteq \mathcal{O}$. Let $(a_{\rho})$ be a well-indexed sequence in $K$ such that $a_{\rho} \leadsto a$, with $a\in K$. Then $A(a_{\rho}) \leadsto A(a)$.
\end{cor}
\begin{proof} Replacing $a_{\rho}$ by $a_{\rho} -a$ we reduce to the case
$a=0$. By omitting some initial terms from $(a_{\rho})$ we can further assume that $va_{\rho}$ is strictly increasing as a function of~$\rho$. 
So $v(Aa_{\rho})=v_A(va_{\rho})$ is strictly increasing as a function of $\rho$. It remains to note that $A(a_{\rho}) \prec Aa_{\rho}$ for only finitely many $\rho$, by Lemma~\ref{dhli2}, and  $A(a_{\rho})\asymp Aa_{\rho}$ for all other $\rho$. 
\end{proof}

\subsection*{Extension to dominant degree $1$} Lemma~\ref{dhli2} leads to:

\begin{lemma}\label{dhli3} Suppose $r\ge 1$, $K$ is $r$-$\d$-henselian, $C\subseteq \mathcal{O}$ and
$G\in K\{Y\}\setminus K$ has order $\le r$. Then there do not exist 
$y_0,\dots ,y_{r+1}\in K$ such that: \begin{enumerate}
\item[\textup{(i)}] $y_{i-1}-y_{i}\succ y_{i}-y_{i+1}$ for all $i\in\{1,\dots, r\}$, and
$y_r\ne y_{r+1}$;
\item[\textup{(ii)}] $G(y_0)=\cdots = G(y_{r+1})=0$;
\item[\textup{(iii)}] $\operatorname{ddeg} G_{+y_{r+1},\times g}=1$ and 
$y_0-y_{r+1}\preceq g$ for some $g\in K^{\times}$.
\end{enumerate}
\end{lemma}
\begin{proof} Towards a contradiction, suppose $y_0,\dots, y_{r+1}\in K$ satisfy (i), (ii), (iii). Set $a_i:= y_i-y_{r+1}$ for $i=0,\dots,r$ and
$P:= G_{+y_{r+1}}$. Then $a_i\sim y_i-y_{i+1}$ for $i=0,\dots,r$ and
so by (i) and (ii),
$$a_0\succ a_1\succ \dots \succ a_{r}\ne 0, \qquad P(a_0)=\cdots=P(a_{r})= P(0)=0.$$
Now $P=A+R$ with $A=P_1$ and $\val R\ge 2$. Taking $g$ as in (iii) we have
$g\succeq a_0$, so $a_i=gb_i$ with $b_i\preceq 1$ for $i=0,\dots,r$, and $b_0 \succ b_1 \succ\dots \succ b_r$.
Also $P_{\times g}=A_{\times g} + R_{\times g}$ with
$R_{\times g} \prec A_{\times g}$ by (iii), and for $i=0,\dots,r$,
\begin{align*} P(a_i)\ &=\ A_{\times g}(b_i)+ R_{\times g}(b_i)\ =\ 0,\ \text{ so}\\
A_{\times g}(b_i)\ &=\ -R_{\times g}(b_i)\ \preceq\ R_{\times gb_i}\ \prec\ A_{\times gb_i}
\end{align*}
which contradicts Lemma~\ref{dhli2}.
\end{proof}

\begin{prop}\label{noextrazeros} Let $r\ge 1$, and assume $K$ is $r$-$\d$-henselian and $C\subseteq \mathcal{O}$. Let 
$G\in K\{Y\}$ with $\order G\leq r$ and
$\ddeg G =1$. 
Let $E$ be an immediate extension of~$K$.
Then $G$ has the same zeros in $\mathcal{O}$ as in $\mathcal{O}_E$.
\end{prop}
\begin{proof} Note first that $\ddeg G_{+y}=1$ for all 
$y\in \mathcal{O}_E$. Towards a contradiction, suppose $G(\ell)=0$
with $\ell\in \mathcal{O}_E\setminus\mathcal{O}$. Now 
$\ell\preceq 1$ gives $\ddeg G_{+\ell}=1$, and from $G(\ell)=0$
it follows easily that $\ddeg G_{+\ell, \times g}=1$ for all $g\preceq 1$ in $E^\times$.  

\claim{Let $\gamma\in v(\ell-K)^{\geq 0}$. Then there is a $y\in \mathcal{O}$ with
$G(y)=0$ and $v(\ell-y)\ge \gamma$.}

\noindent
To prove this claim, take $a\in K$ and $g\in K^\times$ such that 
$v(\ell -a)=vg=\gamma$. Then by Corollary~\ref{aftra} and the 
observation preceding the claim,
$\ddeg G_{+a,\times g}\ =\ \ddeg G_{+\ell,\times g}\ =\ 1$, so
we get $b\in \mathcal{O}$ such that $G(a+gb)=0$, so $y:=a+gb$ satisfies the claim. Having proved the claim, we get
$y_0,\dots, y_r, y_{r+1}\in \mathcal{O}$ such that 
$$\ell -y_0 \succ \ell - y_1 \succ \cdots \succ  \ell-y_{r+1}, \qquad 
G(y_0)\ =\ G(y_1)\ =\ \cdots\ =\ G(y_{r+1})\ =\ 0,$$ 
contradicting Lemma~\ref{dhli3}: take $g=1$ in (iii).    
\end{proof}

\noindent
We also have the following variant:

\begin{lemma} Assume that $C\subseteq \mathcal{O}$. Let
 $A\in K[\der]^{\ne}$ be neatly surjective. Let $E$ be an immediate extension of $K$. Then $\ker_E A = \ker A$. 
\end{lemma} 
\begin{proof} Take $r\ge 1$ and $a_0,\dots, a_r\in K$ such that 
$A= a_0 + a_1\der + \cdots + a_r\der^r$. Set $G:= a_0Y + \cdots + a_rY^{(r)}\in K\{Y\}$. We claim that $G$ has the same zeros in~$\mathcal{O}$ as in
$\mathcal{O}_E$, and for this we follow the proof by contradiction of Proposition~\ref{noextrazeros}, deriving $\ddeg G_{+a,\times g}= 1$ as in
that proof. Now $G_{+a, \times g}(Y)=G(a) + G_{\times g}(Y)$, and
$G(a)=A(a)$, while $G_{\times g}$ corresponds to $Ag$. Thus
$v(Ag) \le v\big(A(a)\big)$, and as $A$ is neatly surjective, we get
$b\in \mathcal{O}$ with $A(gb)=-A(a)$, so $y=a+gb$ gives $A(y)=0$
and $v(\ell - y)\ge v(\ell - a)$.
  As in the proof of Proposition~\ref{noextrazeros}, this argument yields 
$y_0,\dots, y_r, y_{r+1}\in K$ such that 
$$\ell -y_0\succ \ell - y_1 \succ\cdots \succ  \ell-y_{r+1}, \qquad 
A(y_0)\ =\ A(y_1)\ =\ \cdots\ =\ A(y_{r+1})\ =\ 0,$$
so $y_0-y_{r+1}\succ y_1-y_{r+1} \succ \cdots \succ y_r-y_{r+1}$ are
$C$-linearly independent elements of~$\ker A$, which is impossible.
This proves our claim.
For $y\in \ker _E A$ with $y\succ 1$,  
take $f\in K^\times$ with $f\asymp y$, and note that then 
$f^{-1}y\in (\ker_E Af) \cap \mathcal{O}_E$. Applying the above to
$Af$ instead of $A$ we conclude that $f^{-1}y\in \mathcal{O}$, so $y\in K$. 
\end{proof}

\noindent
The following special case is worth recording:

\begin{cor}\label{sameconstants} If $C\subseteq \mathcal{O}$ and  
$\der \in K[\der]^{\ne}$ is neatly surjective, then $C_E=C$ for any immediate extension $E$ of $K$. 
\end{cor}

\noindent
It turns out that any $K$ satisfying the assumptions of Corollary~\ref{sameconstants} (and thus any $\d$-henselian $K$ with few constants) is {\em asymptotic\/} in the sense of Chapter~\ref{ch:asymptotic differential fields}, by Lem\-ma~\ref{neatderas}. For more on such $K$, see Section~\ref{Asymptotic-Fields-Basic-Facts}, in particular, Corollary~\ref{asexistdhenselization}. 
When does $K$ extend to a
$\d$-henselian valued differential field with few constants? (The conclusion of Lemma~\ref{dhli2} holds for such $K$, any~${r\ge 1}$, and any $A\in K[\der]^{\ne}$ of order~${\le r}$.)
Corollary~\ref{aspsidh} and subsequent remarks address this question.
%When does $K$ extend to a
%$\d$-henselian valued differential field with few constants?
%Note that the conclusion of Lemma~\ref{dhli2} holds for such $K$, 
%any~${r\ge 1}$, and any $A\in K[\der]^{\ne}$ of order 
%$\le r$. 

%\begin{cor}\label{existdhenselization} Assume $\der$ is %neatly surjective, $\k$ is linearly surjective, 
%and~$C\subseteq \mathcal{O}$. Then $K$ has an immediate %extension $L$ such that: \begin{enumerate}
%\item $L$ is $\d$-algebraic over $K$;
%\item $L$ is $\d$-henselian;
%\item no proper differential subfield of $L$ containing $K$ %is $\d$-henselian.
%\end{enumerate}
%\end{cor}
%\begin{proof} Corollary~\ref{henimm} yields an 
%immediate $\d$-henselian extension $F$ of $K$. 
%By keeping only the elements of 
%$F$ that are $\d$-algebraic over $K$ we arrange that 
%$F$ is also $\d$-algebraic over $K$. 
%Note that $C_F=C$ by Corollary~\ref{sameconstants}. 
%Let $L$ be the intersection inside
%$F$ of the collection of all differential subfields 
%$E$ of $F$ that contain~$K$ and are $\d$-henselian. 
%Applying Lemma~\ref{noextrazeros} to the extensions 
%$E\subseteq F$ shows that $L$ has the desired property.
%\end{proof} 

\subsection*{Few constants and monotonicity} Having few constants
leads to the promised converse of Theorem~\ref{thm:damdh} in the monotone case:

\begin{proof}[Proof of Theorem~\ref{thm:fcdifhensalgmax}]
Assume $C\subseteq \mathcal{O}$ and $K$ is monotone and
$\d$-hen\-sel\-ian; our job is to show that $K$ is $\d$-algebraically maximal. Towards a contradiction, let~$K\<a\>$ be an immediate
$\d$-algebraic extension of $K$ with $a\notin K$. 
Take a divergent pc-sequence~$(a_{\rho})$ in $K$ such that
$a_{\rho} \leadsto a$. Then 
$(a_{\rho})$ is of $\d$-algebraic type over $K$ by Corollary~\ref{cor:diftr}, so we have 
a minimal differential polynomial~$G(Y)$ of $(a_{\rho})$ 
over~$K$. Replacing $(a_{\rho})$ by an equivalent
pc-sequence in $K$ we arrange that $G(a_{\rho}) \leadsto 0$. Then the
assumptions of Proposition~\ref{dhconf} are satisfied, so taking $g_{\rho}\in K$
with $g_{\rho}\asymp a-a_{\rho}$ we have 
$\ddeg G_{+a, \times g_{\rho}}=1$ eventually, and we can arrange this
holds for all $\rho$. Then 
$\ddeg G_{+a_{\rho}, \times g_{\rho}}=1$ for all $\rho$
by Corollary~\ref{aftra}. This gives for each~$\rho$ an element
$z_{\rho}\in K$ with $G(z_{\rho})=0$ and 
$a_{\rho}-z_{\rho}\preceq g_{\rho}$, so $a-z_{\rho}\preceq g_{\rho}$. 

Take $r\ge 1$ such that $G$ has order $\le r$. Pick some index $\rho_0$ and set $g:= g_{\rho_0}$.
In view of Lemma~\ref{pcnonmax} the above yields indices $\rho_0 < \rho_1< \cdots < 
\rho_{r+1}$ such that $a-z_{\rho_i}\succ a-z_{\rho_j}$ whenever $0\le i < j \le r+1$. Set $y_i:= z_{\rho_i}$ for $i=0,\dots, r+1$. Then 
conditions (i) and (ii) of Lemma~\ref{dhli3} are satisfied. Also $a-y_{r+1}\prec a-y_0\preceq g$, hence
$y_0-y_{r+1}\preceq g$ and $a_{\rho_0}-y_{r+1}\preceq g$. In view of 
Corollary~\ref{aftra} the latter gives 
$\ddeg G_{+y_{r+1},\times g}=\ddeg G_{+a_{\rho_0}, \times g}=1$, so
condition (iii) in Lemma~\ref{dhli3} also holds, which contradicts
that lemma.
\end{proof}

\section{Differential-Henselianity in Several Variables}\label{sec:max}

\noindent
In this section we prove Theorem~\ref{thm:maxndh} in several stages.
After some preliminaries we first handle the case of
spherically complete $K$ with archimedean $\Gamma$; this goes by
diagonalization and successive approximation. Next we treat arbitrary
spherically complete~$K$ using a reduction to the previous case by iterated coarsening. The last stage is an appeal to 
a result in commutative differential algebra due to J.~Johnson.

In this section $n$ ranges over $\N^{\ge 1}$, and for $y=(y_1,\dots, y_n)\in K^n$ we set
$vy:= \min(vy_1,\dots, vy_n)$. This makes $K^n$ into a valued abelian group, and in particular, $vy\ge 0$ means $y\in \mathcal{O}^n$, and
$vy>0$ means $y\in \smallo^n$. Recall that if $K$ is spherically complete, then $K^n$ is spherically complete as a valued abelian group.

\subsection*{Notations and some easy equivalences}  
Let $Y=(Y_1,\dots, Y_n)$ be a tuple of distinct differential indeterminates, and equip $K\{Y\}$ with the
gaussian extension of the valuation of $K$. Let $P=(P_1,\dots, P_n)\in K\{Y\}^n$. Then we set
\begin{align*} vP\ &:=\ \min(vP_1,\dots, vP_n),\\  
P(y)\ &:=\ \big(P_1(y),\dots, P_n(y)\big)\in K^n\ \text{ for $y\in K^n$.}
\end{align*} 
A {\bf solution of
the system $P(Y)=0$ in $K$} is a point $y\in K^n$ such that $P(y)=0$. 
Also, a {\bf solution of $P(Y)=0$ in $\mathcal{O}$} is a point $y\in \mathcal{O}^n$ with $P(y)=0$, and likewise, with $\smallo$ instead of
$\mathcal{O}$. 
Recall that for $F\in K\{Y\}$ and $y=(y_1,\dots, y_n)\in K^n$,
$$F_{+y}\ =\ F(y+Y)\ =\ F(y_1+Y_1,\dots, y_n+Y_n)\in K\{Y\}.$$
We also set $P_{+y}:= (P_{1,+y},\dots, P_{n,+y})\in K\{Y\}^n$
for $y\in K^n$.  

\medskip\noindent
In the rest of this section, let $P\in \mathcal{O}\{Y\}^n$ and consider the associated system
\begin{equation}\label{eq:dh system}
P(Y)\ =\ 0. \tag{$\ast$}
\end{equation}
Let $A_i=P_{i,1}$ be the homogeneous part of $P_i$ of degree $1$, 
with image $\bar{A}_i$ in $\k\{Y\}_1$. 

\medskip\noindent
We say that
the system \eqref{eq:dh system} is in {\bf differential-hensel position} (dh-position)
if the differential polynomials $\bar{A}_1,\dots, \bar{A}_n$ are $\d$-independent and 
$P_1(0)\prec 1,\dots, P_n(0)\prec 1$. 
(For $n=1$ this means ``dh-position at $0$'' as defined in 
Section~\ref{preldh}.) Note that if \eqref{eq:dh system} is in dh-position and
$y\in \smallo^n$, then $P_{+y}=0$ is in dh-position. 

\index{differential-hensel!position}
\index{dh-position}

\index{diagonal!differential-hensel position}
\index{differential-hensel!position, diagonal}
\index{ddh-position}

We say that the system
\eqref{eq:dh system} 
is in {\bf diagonal differential-hensel position} (ddh-po\-si\-tion) if for 
$i=1,\dots, n$, 
$$ 0\ne \bar{A}_i\in \k\{Y_i\},\qquad  P_i(0)\prec 1.$$
Note that if \eqref{eq:dh system} is in ddh-position, then it is in dh-position.
Suppose now that \eqref{eq:dh system} is in ddh-position, so 
for $i=1,\dots,n$,
$$P_i(Y)\ =\ a_i +  A_i(Y) + R_i(Y),\qquad a_i\in \smallo,$$
where $R_i\in \mathcal{O}\{Y\}$ has only terms of degree $\ge 2$.
Set $a:=(a_1,\dots, a_n)$. Then a solution of \eqref{eq:dh system} cannot be much smaller than $a$: 

\begin{lemma} Suppose $a\ne 0$ and $y$ is a solution
of \eqref{eq:dh system} in $K$. Then $$vy\ \le\ va+ o(va).$$
\end{lemma}
\begin{proof} Suppose $vy\ge (1+\epsilon)va$ where $\epsilon\in \Q^{>}$. Then for $i=1,\dots,n$ we have
$A_i(y), R_i(y) \prec a$,
 which for $i$ with $va_i=va$ is impossible.
\end{proof}

\begin{lemma}\label{eqdhdhh} Assume $\Gamma\ne \{0\}$. The following are equivalent: 
\begin{enumerate}
\item[\textup{(i)}] every system \eqref{eq:dh system} such that
$\bar{A}_1,\dots, \bar{A}_n$ are $\d$-independent and the coefficients of the
monomials in $P$ of degree $\ge 2$ are in $\smallo$ has a solution
in $\mathcal{O}$;
\item[\textup{(ii)}] $\k$ is linearly surjective, and every \eqref{eq:dh system}
in ddh-position has a solution in $\smallo$;
\item[\textup{(iii)}] $\k$ is linearly surjective, and every \eqref{eq:dh system} 
in dh-position has a solution in $\smallo$.
\end{enumerate}
\end{lemma}
\begin{proof} Suppose (i) holds. For $n=1$ this yields by
Lemma~\ref{Hdhencor} that $K$ is $\d$-henselian, so
$\k$ is linearly surjective. Let \eqref{eq:dh system} be in ddh-position.
Using the notations above,  
$P_i=a_i + A_i +R_i$, $a_i\prec 1$ and
 $A_i(Y)=D_i(Y_i) + E_i(Y)$ with $D_i\in K\{Y_i\}_1$, $D_i\asymp 1$,
 and $E_i\in K\{Y\}_1$, $E_i\prec 1$, for $i=1,\dots,n$. 
To get~(ii) we need to find a solution of \eqref{eq:dh system} in $\smallo$. We can take nonzero $g\in \smallo$ with $vg>0$ so small that
$$vD_i(gY_i)\ \le\  va_i, \qquad vD_i(gY_i)\ <\  vE_i(gY), vR_i(gY) \qquad(i=1,\dots,n).$$ Taking nonzero 
$h_i\in \smallo$ with $vh_i=vD_i(gY_i)$, we obtain a system
$$h_1^{-1}P_1(gY)\ =\ \cdots\ =\ h_n^{-1}P_n(gY)\ =\ 0.$$ Considering
the homogeneous parts of degree $1$ of the $h_i^{-1}P_i(gY)$, 
this system has a solution $y$ in $\mathcal{O}$ by (i), and so 
$gy$ is a solution of \eqref{eq:dh system} in $\smallo$.

Next, assume (ii). Let \eqref{eq:dh system} be in dh-position.
By Lemma~\ref{diagval} we can transform this into a system
$Q(Y)=0$ in ddh-position such that any solution of $Q(Y)=0$ in
$\smallo$ gives rise to a solution of \eqref{eq:dh system} in $\smallo$.
Now $Q(Y)=0$ has a solution in $\smallo$ by (ii), and so~\eqref{eq:dh system} has a solution in $\smallo$. This gives (iii). 

Finally, assume (iii). Let \eqref{eq:dh system} be such that
$\bar{A}_1,\dots, \bar{A}_n$ are $\d$-independent and 
the coefficients of the
monomials in $P$ of degree~$\ge 2$ are in $\smallo$. To get (i) we need to find a solution of \eqref{eq:dh system} in $\mathcal{O}$. Since $\k$ is
linearly surjective we can use Lemma~\ref{fdzind} to get~${y\in \mathcal{O}^n}$ such that
$vP(y)>0$. Then $P_{+y}(Y)=0$ is a system in dh-position, and so
has a solution $b$ in $\smallo$ by (iii), and then $y+b$ is a solution of
\eqref{eq:dh system} in $\mathcal{O}$.
\end{proof}

\subsection*{The case of an archimedean value group} We define
$\mathcal{O}$ to be {\bf strongly linearly surjective} if for all
$n\ge 1$ and 
$D_1(Y_1)\in \mathcal{O}\{Y_1\}_1,\dots, D_n(Y_n)\in \mathcal{O}\{Y_n\}_1$ 
with $D_1,\dots, D_n\asymp 1$, all $E_1(Y),\dots, E_n(Y)\prec 1$ in
$\mathcal{O}\{Y\}_1$, and all $a=(a_1,\dots, a_n)\in \mathcal{O}^n$ there exists a $y\in \mathcal{O}^n$ such that $y\asymp a$ and
$$D_1(y_1) + E_1(y)\ =\ a_1,\ \dots,\  D_n(y_n)+E_n(y)\ =\ a_n.$$

\index{strongly!linearly surjective}
\index{linearly!surjective!strongly}

\begin{lemma}\label{maxmonstrls} Suppose $\Gamma:= v(K^\times)$ is archimedean, $\k$ is linearly surjective, and~$K$ is spherically complete. Then $\mathcal{O}$ is strongly linearly surjective.
\end{lemma}
\begin{proof} Note that $K$ is monotone by Corollary~\ref{rk1cor}.
In addition, $K$ is $\d$-hen\-sel\-ian. Therefore, if $A\asymp 1$ in $K[\der]$,
then by Corollary~\ref{monvp} and Lemma~\ref{dhns}
there is for each $a\in K$ an
element $y\in K$ such that $A(y)=a$ and $y\asymp a$.

 Let $D_1,\dots, D_n, E_1,\dots, E_n$ be as above. Define $\Q$-linear maps 
$$D, E\colon \mathcal{O}^n \to \mathcal{O}^n,\quad
D(y):=\big(D_1(y_1),\dots, D_n(y_n)\big),  \ E(y):=\big(E_1(y),\dots, E_n(y)\big).$$
From the fact stated in the beginning of the proof and by Lemma~\ref{valvecrightinv} we obtain a valuation preserving right-inverse
$D^*$ to $D$, that is, $D^*\colon \mathcal{O}^n\to \mathcal{O}^n$
is $\Q$-linear, $D\circ D^*=\operatorname{id}_{\mathcal{O}^n}$, and
$D^*(a)\asymp a$ for all $a\in \mathcal{O}^n$. Consider the
$\Q$-linear map $F:= (D+E)\circ D^*=\operatorname{id}_{\mathcal{O}^n}+ G\colon \mathcal{O}^n \to \mathcal{O}^n$, with $G:= E\circ D^*$.
Then $G(a)\prec a$ for all nonzero $a\in \mathcal{O}^n$, and so
by Corollary~\ref{fixcor}, 
$F\colon\mathcal{O}^n \to \mathcal{O}^n$ is a valuation-preserving $\Q$-linear bijection. Given $a\in \mathcal{O}^n$, we get $b\in \mathcal{O}^n$ with $F(b)=a$ and $b\asymp a$, so $y:= D^*(b)$ yields 
$(D+E)(y)=a$ and $y\asymp a$. 
\end{proof}

\noindent
Assume $\Gamma$ is archimedean. Let \eqref{eq:dh system} be in ddh-position, so for $i=1,\dots, n$ we have $P_i=a_i + D_i + E_i + R_i$ 
where $a_i\in \smallo$,
 $D_i\in \mathcal{O}\{Y_i\}_1, D_i\asymp 1$, $E_i\in \mathcal{O}\{Y\}_1,\ E_i\prec 1$,
and all terms in $R_i\in \mathcal{O}\{Y\}$ have degree $\ge 2$.   
We try to construct a solution~$y\prec 1$ in~$K$ to~\eqref{eq:dh system}. Set
$a:=(a_1,\dots, a_n)$.
If $a=0$, then $y=0$ is such a solution, so assume~$a\ne 0$. We associate to \eqref{eq:dh system} the {\em linear\/} system 
$$ a_1+D_1(Y_1)+E_1(Y)\ =\ \cdots\  =\  a_n+ D_n(Y_n)+ E_n(Y)\ =\ 0.$$
Suppose it
has a solution $y$ in $K$ with $vy=va$. Substituting $y + Y$ for $Y$ gives
\begin{align*} &P_{i,+y}(Y)\ =\ P_i(y+Y)\ =\ 
a_i^* + D_i(Y_i) + E_i^*(Y) + R_i^*(Y),\ \text{ where}\\
a_i^*\ &=\ R_i(y),\ E_i^*\ =\ E_i+ dR_i(y,Y),\ dR_i(y,Y)\ :=\ \sum_{j=1}^n\sum_{k\in \N} \frac{\partial R_i}{\partial Y_j^{(k)}}(y)Y_j^{(k)},
\end{align*}
and $R_i^*\in \mathcal{O}\{Y\}$ has only terms of degree $\ge 2$,  and thus
$$va_i^*\ \ge\ 2vy\ =\ 2va, \qquad 
v\big(dR_i(y,Y)\big)\ \ge\ vy\ =\ va.$$
In particular, $va^* \ge 2va$ for $a^*:=(a^*_1,\dots, a^*_n)$, and $E_i^*\prec 1$. Thus
the modified system $P_{+y}(Y)=0$ with $P_{+y}:=(P_{1,+y},\dots, P_{n,+y})$ is
still in ddh-position.

\begin{lemma}\label{monstrls} Suppose $\Gamma$ is archimedean, $\k$
is linearly surjective, $K$ is spherically complete,
and \eqref{eq:dh system} is in ddh-position. Then there is $f\in \smallo^n$ such that $P(f)=0$. 
\end{lemma} 
\begin{proof}
Let $(f_{i})_{i< m}$ be a sequence in $\smallo^n$
with $m\ge 1$, such that \begin{enumerate}
\item[(1)] $P_{+f_{i}}(Y)=0$ is in ddh-position, for all $i<m$,
\item[(2)] $v(f_{j}-f_{i})=vP(f_{i})$ whenever $i < j<m$,
\item[(3)] $vP(f_{j}) \ge 2vP(f_{i})$ whenever $i< j< m$.
\end{enumerate}
Taking $f_0=0$ shows that such a sequence exists for $m=1$.
Let $m=\mu+1$. By Lemma~\ref{maxmonstrls} we can take $y\in \smallo^n$ such that $vy=vP(f_{\mu})$ and
$y$ is a solution of the linear system associated to
$P_{+f_{\mu}}(Y)=0$. Since $P(f_{\mu})\ne 0$, the arguments preceding the lemma
applied to $P_{+f_{\mu}}(Y)$ instead of $P(Y)$ show that 
for $f_{m}=f_{\mu}+y$ we have $vP(f_{m})\ge 2vP(f_{\mu})$
and $P_{+f_{m}}(Y)=0$ is in ddh-position.
Then the extended sequence $(f_{i})_{i< m+1}$ has the above properties with $m+1$ instead of $m$. 

This construction yields a c-sequence $(f_i)_{i\in \N}$ in
$\smallo^n$ with a limit $f\in \smallo^n$, and then~$P(f)=0$.
\end{proof} 

\noindent
Combining Lemmas~\ref{eqdhdhh} and~\ref{monstrls} yields: 

\begin{cor}\label{monlinmaxdh} If $\Gamma$ is archimedean,  $\k$ is linearly surjective, and $K$ is spherically complete, then conditions \textup{(i)}, \textup{(ii)}, \textup{(iii)} of Lemma~\ref{eqdhdhh} are satisfied. 
\end{cor}

\subsection*{Reduction to the case of an archimedean value group} This reduction
rests on iterated coarsening in combination with the following lemma.

\begin{lemma}\label{corminmo} Assume $\k$ is linearly surjective, $K$ is spherically complete, and $\Gamma$ has a smallest nontrivial convex subgroup $\Delta$. Let
\eqref{eq:dh system} be in dh-position, with $vP(0)\in \Delta$. 
Then $vy\in \Delta$ and $vP(y)>\Delta$ for some
$y\in \smallo^n$. 
\end{lemma}
\begin{proof} The valued differential residue field $\dot{K}$
of the $\Delta$-coarsening of $K$ has linearly surjective residue field, has archimedean value group $\Delta$, is maximal as a valued field, and \eqref{eq:dh system} yields a system $\dot{P}(Y)=0$ in dh-position
over $\dot{K}$. Then Corollary~\ref{monlinmaxdh} yields $y=(y_1,\dots, y_n)\in \smallo^n$
such that for $\dot{y}:= (\dot{y}_1,\dots, \dot{y}_n)$ we have
$\dot{P}(\dot{y})=0$. Since $\dot{y}\ne 0$ we have $vy\in \Delta$, and
so $y$ has the desired property. 
\end{proof}

\begin{prop}\label{thm:maxnmdh} Suppose $\k$ is linearly surjective, 
and $K$ is spherically complete. Then every system \eqref{eq:dh system} 
in dh-position has a solution in $\smallo$. 
\end{prop}
\begin{proof} Since $K$ is $\d$-henselian by Corollary~\ref{cor:mdh}, the assumptions on $K$ are inherited by any coarsening of $K$ by a convex subgroup of $\Gamma$. This will be tacitly used in what follows.
Let \eqref{eq:dh system} be in dh-position; our job is to show that \eqref{eq:dh system} has a solution in $\smallo$. If $P(0)=0$ we are done, so assume $P(0)\ne 0$. Then $\Gamma\ne \{0\}$. 
Let $\lambda>0$ be an ordinal, $(f_{\rho})_{\rho< \lambda}$ a sequence in $\smallo^n$, and $(\Delta_{\rho})_{\rho< \lambda}$ an increasing sequence of convex subgroups of $\Gamma$, such that for all $\rho<\lambda$,  \begin{enumerate}
\item[(1)] $vP(f_{\rho}) > \Delta_{\rho}$;
\item[(2)] $\Delta_{\rho}$ is the largest convex subgroup of $\Gamma$ not containing $vP(f_{\rho})$;
\item[(3)] $vP(f_{\rho})\in \Delta_{\rho+1}$ whenever $\rho+1< \lambda$;
\item[(4)] $v(f_{\rho'}-f_{\rho})> \Delta_{\rho}$ whenever $\rho < \rho'<\lambda$;
\item[(5)] $v(f_{\rho+1}-f_{\rho})\in  \Delta_{\rho+1}$ whenever $\rho+1<\lambda$.
\end{enumerate}
Such sequences exist for $\lambda=1$: take $f_0=0$ and let $\Delta_{0}$ be the largest convex subgroup
of $\Gamma$ not containing $vP(0)$. 

Suppose $\lambda=\mu+1$ is a successor ordinal. If $P(f_{\mu})=0$ we are done, so assume $P(f_{\mu})\ne 0$. Using (1) for $\rho:=\mu$ and Lemma~\ref{indepreduc}, the system $P_{+f_{\mu}}(Y)=0$ is in dh-position with respect to the $\Delta_{\mu}$-coarsening of
$K$.  By (2) for $\rho=\mu$ we have a smallest convex subgroup $\Delta$ of $\Gamma$ that properly contains $\Delta_{\mu}$. Then 
$vP(f_{\mu})\in \Delta$, so Lemma~\ref{corminmo} (applied to the
$\Delta_{\mu}$-coarsening of $K$ instead of $K$) yields a 
$y\in \smallo^n$ such that $vy> \Delta_{\mu}$, $vy\in \Delta$, and $vP_{+f_{\mu}}(y)> \Delta$. If $P(f_{\mu}+y)=0$, we are done. Suppose $P(f_{\mu}+y)\ne 0$. Then we set $f_{\lambda}:= f_{\mu}+ y$ and let $\Delta_{\lambda}$ be the largest convex subgroup of $\Gamma$ not containing $vP(f_{\lambda})$, in particular, $\Delta\subseteq \Delta_{\lambda}$. It is clear that then the conditions~(1)--(5) are satisfied with $\lambda+1$ instead of $\lambda$.

Next assume that $\lambda$ is a limit ordinal. Then $(f_{\rho})_{\rho<\lambda}$ is a pc-sequence by (4) and~(5), and thus has a pseudolimit $f_{\lambda}$ in
$\smallo^n$. Then $v(f_{\lambda}-f_{\rho}) > \Delta_{\rho}$ for all
$\rho< \lambda$, so $vP(f_{\lambda}) > \Delta_{\rho}$
for all $\rho< \lambda$, by (1). If $P(f_{\lambda})=0$ we are done, so assume
$P(f_{\lambda})\ne 0$.  Let $\Delta_{\lambda}$ be the largest convex subgroup of $\Gamma$ not containing $vP(f_{\lambda})$. Then
$\Delta_{\lambda}\supseteq \Delta_{\rho}$ for all $\rho< \lambda$. Thus conditions (1)--(5) hold for $\lambda+1$ instead of $\lambda$. 

The construction above cannot continue indefinitely, so must end in producing a solution of \eqref{eq:dh system} in $\smallo$. 
\end{proof}

\subsection*{The final step and an application}
Let $Y=(Y_1,\dots, Y_n)$ be as before.

\begin{proof}[Proof of Theorem~\ref{thm:maxndh}] We have $\d$-algebraically maximal $K$ with linearly surjective $\k$, and a system
\eqref{eq:dh system} in dh-position; our job is to show that \eqref{eq:dh system}
has a solution in~$\smallo$. Take a maximal immediate extension
$L$ of $K$. Then Proposition~\ref{thm:maxnmdh} yields 
$y\in \smallo^n_L$ such that $P(y)=0$. By Lemma~\ref{indepreduc}
and Theorem~\ref{thm:johnson} the immediate extension~$K\<y\>$ of $K$ is 
$\d$-algebraic, so $y\in \smallo^n$.
\end{proof}

\noindent
Using also Lemma~\ref{fextind}, this has the following consequence:

\begin{cor}\label{dhenspresalg} Suppose $\k$ is linearly surjective, $K$ is $\d$-algebraically maximal, and
$L$ is an algebraic extension of $K$ with
$\Gamma_L=\Gamma$. Then $L$ is $\d$-henselian.
\end{cor}
\begin{proof} We can reduce to the case that $[L:K]=n<\infty$, and
as $K$ is henselian, this yields a basis $b_1,\dots, b_n$ of
the vector space $L$ over $K$ such that $b_i\asymp 1$ for all~$i$ and
$\overline{b}_1,\dots, \overline{b}_n$ is a basis of $\k_L$ over $\k$.
Note that $\k_L$ is linearly surjective, by 
Corollary~\ref{cor:lin surj under alg extensions}. 
Let $X$ be a single differential indeterminate and 
$F\in \mathcal{O}_L\{X\}$ such that $A:= F_1\asymp 1$ and $F(0)\prec 1$;
our job is to show that $F$ has a zero in $\smallo_L$. Making the
substitution $X=b_1Y_1+ \cdots + b_nY_n$ we have
\begin{align*} F(b_1Y_1+\cdots + b_nY_n)\ &=\ P_1(Y)b_1+\cdots + P_n(Y)b_n, \qquad P_1,\dots, P_n\in \mathcal{O}\{Y\},\\
A(b_1Y_1+\cdots + b_nY_n)\ &=\ A_1(Y)b_1+\cdots + A_n(Y)b_n, \qquad A_1,\dots, A_n\in \mathcal{O}\{Y\},\\
\overline{A}(\overline{b}_1Y_1+ \cdots + \overline{b}_nY_n)\ &=\ \overline{A}_1(Y)\overline{b}_1+ \cdots + \overline{A}_n(Y)\overline{b}_n, \qquad \overline{A}_1,\dots, \overline{A}_n\in \k\{Y\}.
\end{align*}
Also $A_i$ is the homogeneous part of degree $1$ of $P_i$, and
$\overline{A}_1,\dots, \overline{A}_n\in \k\{Y\}_1$ are 
$\d$-independent
by Lemma~\ref{fextind}. From $F(0)\prec 1$ we get 
$P_1(0)\prec 1,\dots, P_n(0)\prec 1$, and so by Theorem~\ref{thm:maxndh} the system $P_1(Y)=\cdots = P_n(Y)=0$ has
a solution $y=(y_1,\dots, y_n)\in \smallo^n$, which gives a zero
$y_1b_1+\cdots + y_nb_n\in \smallo_L$ of $F$. 
\end{proof}

\noindent
This corollary might not be optimal: we would like to drop
the hypothesis $\Gamma_L=\Gamma$, and it would also be nice to
weaken the assumption on $K$ to $\d$-henselianity, or strengthen the conclusion to $L$ being $\d$-algebraically maximal.

%% file: mt-8.tex
\chapter{Differential-Henselian Fields with Many Constants}\label{ch:monotonedifferential}

\setcounter{theorem}{0}

\noindent
The results in this brief chapter will not be used later, but 
complement the earlier generalities. The valued 
differential fields considered here are easier to analyze as
to their model-theoretic properties than $\T$, and so this may
provide an illuminating contrast with our later
focus on objects like $\T$. Note
that {\em $\d$-henselian\/} includes having small derivation, so $\d$-henselian valued differential fields with many constants are {\em monotone}: this is a strong restriction, opposite in spirit to the {\em asymptotic\/} condition imposed in later chapters.

\medskip\noindent
Our goal here is to derive
Scanlon's extension in \cite{Scanlon, Scanlon03} of the Ax-Kochen-Er\accentv{s}ov theorems
to $\d$-henselian valued differential fields with many constants.
%(We do not assume as in \cite{Scanlon} that the multiplicative 
%group of the residue field is divisible.)
This is largely an application of Chapter~\ref{sec:dh1}, in particular Section~\ref{umie}. 

Given structures $\mathbf{M}$ and $\mathbf{N}$ for the same 
(first-order) language,
$\mathbf{M} \equiv \mathbf{N}$ means that $\mathbf{M}$ and 
$\mathbf{N}$ are elementarily equivalent: they satisfy the same sentences in that language; see Appendix~\ref{app:modth}.
Among the
results to be established is the following:

\begin{theorem}\label{akesc} Suppose $K$ and $L$ are $\d$-henselian
valued differential fields with many constants. Then $K\equiv L$ as valued 
differential fields if and only if $\res{K}\equiv \res{L}$ as differential fields and 
  $\Gamma_K\equiv \Gamma_L$ as ordered abelian groups.
\end{theorem}

\noindent
In particular, if $K$ is a $\d$-henselian valued differential field with many constants, and with differential residue field $\k$
and value group $\Gamma$, then $K\equiv \k\(( t^{\Gamma}\)) $, where~$\k\(( t^{\Gamma}\)) $ is the Hahn differential field of Example~(1) in Section~\ref{Valdifcon}.

Theorem~\ref{akesc} refers to the logical notion of elementary equivalence. We
derive from it a purely algebraic result for which we have no other proof:

\begin{cor}\label{mondhalg} Let $K$ be a $\d$-henselian valued differential field with many 
constants. Then every valued differential field extension of $K$ that is 
algebraic over $K$ is also $\d$-henselian.
\end{cor}

\noindent
By Corollary~\ref{cor:monalgmon}, if $L$ is
a valued differential field extension of a monotone valued
differential field $K$
and $L$ is algebraic over $K$, then $L$ is monotone, so has small derivation. Here is an
easy proof
of Corollary~\ref{mondhalg}, using Theorem~\ref{akesc}:

\begin{proof} First some remarks on an arbitrary valued field extension
$E\subseteq F=E(a)$ where $E$ is henselian, and $[F:E]=n$. We can express
properties of the valued field~$F$ in terms of the valued field $E$, as follows.
Let  $$P(X)\ =\ X^n+ a_{n-1}X^{n-1}+ \cdots + a_1X + a_0$$
be the minimum polynomial of $a$ over $E$. Let $f\in F=E(a)$,
and let $(b_0,\dots, b_{n-1})$ be the unique tuple in $E^n$
such that $f=b_0+ b_1a + \cdots + b_{n-1}a^{n-1}$. For $1\le m\le n$ and
$f_0,\dots, f_{m-1}\in E$ we can then express ``$X^m + f_{m-1}X^{m-1} + \cdots + f_1X + f_0$ is the minimum polynomial of $f$ over 
$E$'' as a first-order condition
$$E \models C_{m,n}(a_0,\dots, a_{n-1},b_0,\dots, b_{n-1}, f_0,\dots, f_{m-1})$$
on the quantities indicated. By Corollary~\ref{cor:unique val, 2} we can then
express the valuation on~$F$ in terms of the valuation on $E$ by
$v(f)=v(f_0)/m$ with $f, f_0, \dots, f_{m-1}$ as above. 
If in addition $\der$ is a derivation on $E\supseteq \Q$ and 
$\der$ also stands for its unique extension to a derivation on
$F$, then, with $P':= \partial P/\partial X\in E[X]$, we have
$$\der(a)\ =\ -P^{\der}(a)\big/P'(a), \qquad \der(f)\ =\ \sum_{i=0}^{n-1} \der(b_i)a^i + \sum_{i=1}^{n-1}ib_i\der(a)^{i-1}.$$
Turning now to our $K$, let $n\ge 1$ be given. 
By the remarks above there is a set~$\Sigma_n$ of sentences in the language of
valued differential fields, independent of~$K$, such that
$K\models \Sigma_n$ if and only if every valued differential field extension 
$L$ of $K$ with $[L:K]=n$ is $\d$-henselian. 
Now with $\k$ the differential residue field of $K$ and
$\Gamma=v(K^\times)$ we have $K\equiv \k\(( t^\Gamma\)) $. But for the Hahn differential field $\k\(( t^{\Gamma}\)) $ we have $\k\(( t^{\Gamma}\)) \models \Sigma_n$,
since every valued differential field extension $L$ of $\k\(( t^{\Gamma}\)) $ with 
$\big[L:\k\(( t^{\Gamma}\)) \big]=n$  
is maximal as a valued field by Corollary~\ref{cor:finite ext of complete}, 
and hence $\d$-henselian by Corollaries~\ref{cor:lin surj under alg extensions} and~\ref{henseld3}.  Thus
$K\models \Sigma_n$. As this holds for all $n\ge 1$, we obtain the 
desired result.  
\end{proof}

\noindent
In model-theoretic work on valued
fields and their expansions we often prefer a many-sorted set-up. 
This can be useful even in stating results properly. Thus we
shall consider here
three-sorted structures
$$ \mathbf K\ =\ \big(K, \k, \Gamma; \pi, v  \big)$$
where $K$ and $\k$ are fields, $\Gamma$ is an ordered abelian 
group, $v\colon K^\times \to \Gamma$ is a (surjective) valuation making
$K$ into a valued field, and
$\pi\colon \mathcal{O} \to \k$ with $\mathcal{O}:= \mathcal{O}_v$ is a surjective
ring morphism.
Note that then $\pi$ has kernel $\smallo_v$, and thus induces a
field isomorphism
$\k_v \cong \k$ between the residue field $\k_v$ and $\k$ 
such that the diagram
\[
\xymatrix@C=10pt{
&\mathcal{O}_v \ar[dl] \ar[dr]^{\pi}&\\
\mathbf{k}_v \ar[rr]^{\sim}&&\mathbf{k}\\
}
\]
commutes. Let us call
$\Gamma$ the value group of $\mathbf{K}$, and $\k$ its residue field 
(even though the latter is only isomorphic to the residue 
field $\k_v$ of $\ca{O}_v$).
We shall refer to~$\mathbf{K}$ as a valued field, since it
represents a way to construe the (one-sorted) valued field~$(K, \ca{O}_v)$ as a three-sorted model-theoretic structure where
the residue field and the value group are more explicitly present.

To adapt this setting to valued {\em differential\/} fields we 
consider three-sorted structures~$\mathbf{K}$ as above, but with some additional features: 
$K$ and $\k$ are differential fields (rather than just fields), and
$v$ makes $K$ into a valued differential field with small derivation such that
$\pi$ is a differential ring morphism (rather than just a ring morphism).
Note that then the above field isomorphism $\k_v \cong \k$ is actually a
differential field isomorphism. When referring to $\mathbf{K}$ as a valued 
differential field, we assume these additional features are present.

In this three-sorted set-up we have field variables ranging over $K$,
residue variables ranging over $\k$, and value group variables ranging
over $\Gamma$. Given a possibly many-sorted language $\mathcal{L}$ and 
$\mathcal{L}$-structures $\mathbf{M}$ and $\mathbf{N}$, we have the usual
 notions of~$\mathbf{M}$ being a substructure
of $\mathbf{N}$ (notation: $\mathbf{M}\subseteq \mathbf{N}$), of $\mathbf{M}$ and $\mathbf{N}$ being elementarily equivalent (notation: $\mathbf{M}\equiv\mathbf{N}$), and of $\mathbf{M}$ being an elementary substructure
of $\mathbf{N}$ (notation: $\mathbf{M}\preccurlyeq \mathbf{N}$); see Appendix~\ref{app:modth}. 
The main result of this chapter is the Equivalence Theorem~\ref{embed11},
among whose consequences are the following:

\begin{theorem}\label{pelsub} Suppose $\mathbf{K}$ and $\mathbf{K}^*$ are
$\d$-henselian valued differential fields with many 
constants such that $\mathbf{K}\subseteq \mathbf{K}^*$. Then
$\mathbf{K}\preccurlyeq \mathbf{K}^*$ if and only if $\k \preccurlyeq \k^*$  as differential fields, and $\Gamma\preccurlyeq \Gamma^*$ as ordered abelian groups.
\end{theorem} 

\noindent
 By \textit{definable}\/ we mean
\textit{definable with parameters from the ambient structure.}\/

\begin{theorem}\label{defwitt} Suppose $\mathbf{K}$ is a $\d$-henselian valued differential field with many constants. Then each subset of 
$\k^m\times \Gamma^n$ definable in $\mathbf{K}$ is a finite union of rectangles $Y\times Z$ with $Y\subseteq \k^m$ definable in the differential
field $\k$ and $Z\subseteq \Gamma^n$ definable in the ordered abelian group
$\Gamma$.  
\end{theorem} 

\noindent
The next sections lead up to the statement and proof of the Equivalence 
Theorem in Section~\ref{et}, and to conclude this chapter we derive the above
consequences.

\subsection*{Notes and comments} We refer to \cite{Scanlon}
for Scanlon's original treatment and to \cite{Scanlon03} for an update. 
The latter also deals with the mixed characteristic case.

\section{Angular Components} \label{ang}

\noindent
Let
$\mathbf{K}=(K, \k, \Gamma; \pi,v)$ be a valued field as
explained in the introduction to this chapter, with valuation ring
$\mathcal{O}:=\{a\in K:\ va\ge 0\}$.
For a subfield
$E$ of $K$ we set $\mathcal{O}_E:= \mathcal{O}\cap E$, a valuation ring of
$E$, and
$\k_E:=\pi(\mathcal{O}_E)$, a subfield of $\k$. For such~$E$ we also let $E$ stand for the valued subfield
$(E,\mathcal{O}_E)$ of $(K, \mathcal{O})$, as well as for the substructure
$(E, \k_E, v(E^\times); \dots)$ of $\mathbf{K}$.

\medskip\noindent  
An {\bf angular component map} on $\mathbf{K}$ is a 
multiplicative group morphism 
$$\ac \colon K^{\times} \to \k^{\times}$$
such that $\ac(a)=\pi(a)$ whenever $a\asymp 1$ in $K$; we extend it to $\ac \colon K \to \k$ by setting $\ac(0)=0$ (so $\ac(ab)=\ac(a)\ac(b)$ for all $a,b\in K$),  and also refer to this extension as an angular 
component map on $\mathbf{K}$. 

\nomenclature[K]{$\ac$}{angular component map on a valued field}
\index{angular component map}
\index{valued field!angular component map}
\index{map!angular component}

\begin{example}
Construing a Hahn field
$\k\(( t^{\Gamma}\)) $ as a valued field
$$ \mathbf{K}\ =\ \big(\k\(( t^{\Gamma}\)) , \k, \Gamma; \pi, v\big)$$
in the natural way, we have the angular component map
$\ac\colon\k\(( t^{\Gamma}\)) \to \k$ on $\mathbf{K}$ given by
$\ac(ct^\gamma + g)=c$ for $c\in \k^\times$, $\gamma\in \Gamma$,
and $g\in \k\(( t^{\Gamma}\)) $ with $v(g)> \gamma$.
\end{example}

\noindent
A cross-section $s$ on the valued field $\mathbf{K}$ 
yields an angular component map $\ac$ on $\mathbf{K}$ by setting 
$\ac(x)= \pi\big(x/s(vx)\big)$ for $x\in K^\times$. 
Thus by Lemma~\ref{xs1}:

\begin{cor}\label{acm2} If the valued field $\mathbf{K}$
is $\aleph_1$-saturated, then
there exists an angular component map on $\mathbf{K}$. 
\end{cor}

\subsection*{Monotone valued differential fields with angular component}
The presence of an angular component map simplifies the proof of the 
Equivalence Theorem~\ref{embed11}, but in the aftermath 
we can often discard these maps again, by Corollary~\ref{acm3}.

\medskip\noindent    
Let $\mathbf{K}=(K, \k, \Gamma;\pi, v)$ be a 
monotone valued differential field. 
By an {\bf angular component map} on $\mathbf{K}$ we mean
an angular component map $\ac$ on $\mathbf{K}$ as a valued field such that the following conditions are satisfied:
\begin{enumerate}
\item[(1)] $\ac(c)\in C_{\k}\subseteq \k$ for all $c\in C$;
\item[(2)] $\ac(a)'=\ac(a')$ for all $a\in K^\times$ with $v(a)=v(a')$. 
\end{enumerate}
Examples are the Hahn differential fields $\k\(( t^\Gamma\)) $ with angular 
component map given by $\ac(a)=a_{\gamma_0}$ for
nonzero $a=\sum a_\gamma t^\gamma\in \k\(( t^\Gamma\)) $ and $\gamma_0:=va$.

\index{angular component map}
\index{valued differential field!monotone!angular component map}
\index{map!angular component}

\begin{lemma}\label{acm1} Suppose $\mathbf{K}$ has many constants. Then
each angular component map on its valued subfield $C$
extends uniquely to an angular component map on $\mathbf{K}$. 
\end{lemma}
\begin{proof} Given an angular component map $\ac$ on $C$
the claimed extension to $\mathbf{K}$, also denoted by $\ac$, is obtained as 
follows: for $x\in K^\times$
we have $x=uc$ with $u\in K^\times$, $u\asymp 1$, $c\in C^\times$; then
$\ac(x)=\pi(u)\ac(c)$. 
\end{proof}

\noindent
Here is an immediate consequence of Corollary~\ref{acm2} and Lemma~\ref{acm1}: 

\begin{cor}\label{acm3} 
If $\mathbf{K}$ has many constants and is $\aleph_1$-saturated, 
then there is an 
angular component map on $\mathbf{K}$.
\end{cor}

{\sloppy
\subsection*{Notes and comments}
Angular component maps were introduced in \cite{Denef}. In~\cite{Pas}
there is an example of a valued field of residue characteristic~$0$ that has no angular component map (and thus also no cross-section).
Angular components often facilitate access to the definable sets,
and this is also their role in this chapter. 

}

\section{Equivalence over Substructures}\label{et}

\noindent
In this section we consider three-sorted structures
$$ \mathbf{K}\ =\ (K, \k, \Gamma;  \pi, v, \ac )$$
where $(K, \k, \Gamma; \pi, v )$ is a monotone
valued differential field with an angular component map
$\ac\colon K \to \k$ on it.
Such a structure will be called a
{\bf monotone $\ac$-valued differential field}.
The main result of this chapter is Theorem~\ref{embed11}. It
tells us when two $\d$-henselian monotone $\ac$-valued differential 
fields with many
constants are elementarily 
equivalent over a common substructure.  
In Section~\ref{cqe1} we derive from it in the 
usual way some attractive consequences on the elementary theories of
such valued differential fields and on the induced structure on
the value group and differential residue field.  

\index{valued differential field!monotone!$\ac$-valued}
\index{monotone!valued differential field!$\ac$-valued}

\medskip\noindent
If $\mathbf{K}$ is $\d$-henselian, then by Theorem~\ref{lift.res.field} 
there is a differential
ring mor\-phism ${i\colon\k \to \ca{O}}$ such that $\pi\big(i(a)\big)=a$ for all $a \in K$;
we call such $i$ a {\bf lifting} of~$\k$ to~$\mathbf{K}$. This will play a minor
role in the proof of the Equivalence Theorem.

\index{lifting!differential residue field}
\index{residue field!lifting}
\index{residue field!lift}

\medskip\noindent
A {\bf good substructure} \index{substructure!good}\index{good!substructure} of  $\mathbf{K}=(K, \k, \Gamma; \pi, v,\ac)$ is a triple
$\mathbf{E}=(E, \k_{\mathbf{E}}, \Gamma_{\mathbf{E}})$ such that 
\begin{list}{*}{\setlength\leftmargin{2.5em}}
\item[(GS1)] $E$ is a differential subfield of $K$, 
\item[(GS2)] $\k_{\mathbf{E}}$ is a differential subfield of $\k$ with 
$\ac(E)\subseteq \k_{\mathbf{E}}$ (hence $\pi(\ca{O}_E)\subseteq \k_{\mathbf{E}}$),
\item[(GS3)] $\Gamma_{\mathbf{E}}$ is an ordered abelian subgroup of $\Gamma$ 
with $v(E^\times)\subseteq \Gamma_{\mathbf{E}}$.
\end{list}
For good substructures $\mathbf{E}=(E, \k_{\mathbf{E}}, \Gamma_{\mathbf{E}})$ and
 $\mathbf{F}=(F, \k_{\mathbf{F}}, \Gamma_{\mathbf{F}})$, we define 
$\mathbf{E}\subseteq \mathbf{F}$ to mean that 
$E \subseteq F$,  $\k_{\mathbf{E}} \subseteq  \k_{\mathbf{F}}$, and $\Gamma_{\mathbf{E}} \subseteq  \Gamma_{\mathbf{F}}$.
If $E$ is a differential subfield of~$K$ with 
$\ac(E)=\pi(\ca{O}_E)$, then 
$\big(E, \pi(\ca{O}_E), v(E^\times)\big)$ is a good substructure of $\mathbf{K}$, 
and if in addition $F\supseteq E$ is a differential subfield
of $K$ such that $v(F^\times)=v(E^\times)$, then $\ac(F)=\pi(\ca{O}_F)$.
In the remainder of this section 
$$\mathbf{K}\ =\ (K, \k, \Gamma; \pi, v, \ac), \qquad \mathbf{K}^*\ =\ (K^*, \k^*, 
\Gamma^*; \pi^*, v^*, \ac^*)$$ are monotone $\ac$-valued differential fields, with 
valuation rings $\ca{O}$ and $\ca{O}^*$, and
$$\mathbf{E}\ =\ (E,\k_{\mathbf{E}}, \Gamma_{\mathbf{E}}), \qquad \mathbf{E}^*\ =\ (E^*, \k_{\mathbf{E}^*}, \Gamma_{\mathbf{E}^*})$$ are good substructures of
$\mathbf{K}$, $\mathbf{K}^*$, respectively. We put 
$\ca{O}_{E^*}:= \ca{O}^*\cap E^*$.

\medskip
\noindent
A {\bf good map} $\mathbf{f}\colon \mathbf{E} \to \mathbf{E}^*$ \index{map!good}\index{good!map} is a triple $\mathbf{f}=(f, f_{\re}, f_{\v})$ consisting
of a differential field isomorphism
$f\colon E \to E^*$,  a differential field isomorphism $f_{\re}\colon \k_{\mathbf{E}} \to
\k_{\mathbf{E}^*}$, and an ordered group isomorphism $f_{\v}\colon \Gamma_{\mathbf{E}} \to
\Gamma_{\mathbf{E}^*}$, such that
\begin{list}{*}{\setlength\leftmargin{2.5em}}
\item[(GM1)]
$f_{\re}\big(\!\ac(a)\big)=\ac^*\!\big(f(a)\big)$ for all $a \in E$, and
$f_{\re}$ is elementary as a partial map between the
differential fields $\k$ and $\k^*$;
\item[(GM2)] $f_{\v}\big(v(a)\big)=v^*\big(f(a)\big)$ for all $a \in E^\times$, and
$f_{\v}$ is elementary as a partial map between the
ordered abelian groups $\Gamma$ and $\Gamma^*$. 
\end{list}
Let $\mathbf{f}\colon \mathbf{E} \to \mathbf{E}^*$ be a good map as above. Then the 
field part $f\colon E \to E^*$ of $\af$ is a valued differential field isomorphism, and $f_{\re}$ and
$f_{\v}$ agree on $\pi(\ca{O}_E)$ and $v(E^\times)$ with
the maps 
$\pi(\ca{O}_E) \to \pi^*(\ca{O}_{E^*})$ and $v(E^\times) \to v^*({E^*}^\times)$ induced by $f$. 
We say that a good map $\ag= (g, g_{\re}, g_{\v}) \colon \mathbf{F} \to \mathbf{F^*}$ 
{\bf extends} $\af$ if $\mathbf{E}\subseteq \mathbf{F}$, $\mathbf{E^*}\subseteq \mathbf{F^*}$, 
and $g$, $g_{\re}$, $g_{\v}$ extend $f$, $f_{\re}$, $f_{\v}$, respectively.  
The {\bf domain} of $\mathbf{f}$ is $\mathbf{E}$.
Note that if a good map $\mathbf{E}\to \mathbf{E}^*$ exists, then $\mathbf{k}\equiv\mathbf{k}^*$ as differential fields and $\Gamma\equiv \Gamma^*$ as ordered abelian groups.

\medskip
\noindent
The next two lemmas show that various parts of the
conditions (GM1) and (GM2) are automatically
satisfied by certain extensions of good maps.

\begin{lemma}\label{ac1} Let $\af\colon \mathbf{E} \to \mathbf{E}^*$ be a good map, and
suppose
$F\supseteq E$ and~${F^*\supseteq E^*}$ are differential subfields of $K$ and
$K^*$, respectively, such that  
$\pi(\ca{O}_F)\subseteq \k_{\mathbf{E}}$ and $v(F^\times)=v(E^\times)$. Let $g\colon F \to F^*$ be a valued differential field 
isomorphism such that $g$ extends~$f$ and $f_{\re}\big(\pi(u)\big)=\pi^*\big(g(u)\big)$
for all $u\in \ca{O}_F$. Then $\ac(F)\subseteq  \k_{\mathbf{E}}$
and $f_{\re}\big(\!\ac(a)\big)= \ac^*\!\big(g(a)\big)$ for all $a \in F$.  
\end{lemma}
\begin{proof} Let $a \in F$. Then $a=a_1u$ where $a_1 \in E$ and  
$u\in \ca{O}_F$, $v(u)=0$, so $\ac(a)=\ac(a_1)\pi(u)\in \k_{\mathbf{E}}$. 
It follows easily that $f_{\re}\big(\!\ac(a)\big)= \ac^*\!\big(g(a)\big)$.
\end{proof}

\noindent
In the same way we obtain: 

\begin{lemma}\label{ac2} Suppose $\pi(\ca{O}_E)=\k_{\mathbf{E}}$, let 
$\af\colon \mathbf{E} \to \mathbf{E}^*$ be a good map, and let
$F\supseteq E$ and $F^*\supseteq E^*$ be differential subfields of $K$ and
$K^*$, respectively, such that $v(F^\times)=v(E^\times)$. 
Let $g\colon F \to F^*$ be a valued differential field isomorphism extending~$f$.   
Then $\ac(F)= \pi(\ca{O}_F)$ and
$g_{\re}\big(\!\ac(a)\big)= \ac^*\!\big(g(a)\big)$ for all $a \in F$, where the map
$g_{\re}\colon\pi(\ca{O}_F)\to \pi^*(\ca{O}_{F^*})$ is induced by $g$
\textup{(}and thus extends~$f_{\re}$\textup{)}. 
\end{lemma}

\noindent
Lemma~\ref{ac2} is also useful in getting the following:

\begin{lemma}\label{ac2 variant}
Suppose $\pi(\mathcal{O}_E) = \mathbf{k}_{\mathbf{E}}$, let $\mathbf{f}\colon\mathbf{E}\to\mathbf{E}^*$ be a good map. Suppose that $F\supseteq E$ and $F^*\supseteq E^*$ are differential subfields of $K$ and $K^*$, respectively, and are \emph{immediate} extensions of $E$ and $E^*$, respectively. Suppose that $g\colon F\to F^*$ is a valued differential field isomorphism that extends~$f$. Then $\mathbf{g} = (g,f_{\re},f_{\v})$ is a good map that extends~$\mathbf{f}$.
\end{lemma}

\noindent
The following is useful in connection with having many constants:

\begin{lemma}\label{fixr} Let $b\in K^{\times}$. Then the following are
equivalent:
\begin{enumerate}
 \item[\textup{(i)}] There is $c \in C^{\times}$ such that $b\asymp c$. 
 \item[\textup{(ii)}] There is $a \in K^{\times}$
such that $a \asymp 1$ and $a^\dagger = b^\dagger$.
\end{enumerate}
\end{lemma}
\begin{proof} 
If $c \in C^{\times}$ and $b\asymp c$, then $a:=bc^{-1}$ satisfies $a \asymp 1$ and $a^\dagger = b^\dagger$. Conversely, if $a\in K^\times$ and  $a \asymp 1,\ a^\dagger = b^\dagger$, then $b\asymp c:=a^{-1}b\in C^\times$.
\end{proof}

\begin{theorem}\label{embed11} Suppose $\mathbf{K}$, $\mathbf{K}^*$ are $\d$-henselian with many constants. Then any good map
$\mathbf{E} \to \mathbf{E}^*$ is a partial elementary map between $\mathbf{K}$ and 
$\mathbf{K}^*$.
\end{theorem}    
\begin{proof} The theorem holds trivially for $\Gamma=\{0\}$, so assume 
that $\Gamma\ne \{0\}$. 
Let $\mathbf{f}=(f, f_{\re}, f_{\v})\colon \mathbf{E} \to \mathbf{E}^*$
be a good map. By passing to suitable elementary extensions of~$\mathbf{K}$ and $\mathbf{K}^*$ we arrange that
$\mathbf{K}$ and $\mathbf{K}^*$ are $\kappa$-saturated, where 
$\kappa$ is an uncountable cardinal such that 
$|\k_{\mathbf{E}}|,\ |\Gamma_{\mathbf{E}}| < \kappa$.  
Call a good substructure $\mathbf{E}_1=(E_1,\k_1, \Gamma_1)$ of~$\mathbf{K}$ 
{\em small\/} if $|\k_1|,\ |\Gamma_1|<\kappa$.
We shall prove that the good maps with small domain form a back-and-forth
system from $\mathbf{K}$ to $\mathbf{K}^*$. (This clearly suffices to obtain 
the theorem: Proposition~\ref{prop:bf, 2}.) In other words, we shall prove that under the 
present assumptions on~$\mathbf{E}$, $\mathbf{E}^*$ and $\mathbf{f}$, there is   
for each $a \in K$ a good map $\ag$ extending~$\af$ such that
$\ag$ has small domain $\mathbf{F}=(F,\dots)$ with $a\in F$. The most
delicate of the extension procedures we need comes from Corollary~\ref{immsat}
(to extend domains) and Theorem~\ref{unique.max.imm.ext} (to extend good maps).
In addition we have 
several other basic extension procedures:

\ifbool{PUP}{}{\bigskip}
\subsubsection*{\textup{(1)} 
Given $d\in \k$, arranging that $d\in \k_{\mathbf{E}}$.}
By saturation and the definition of ``good map'' this can be achieved without changing $f$, $f_{\v}$, $E$, $\Gamma_{\mathbf{E}}$ by extending $f_{\re}$ to a partial elementary map between $\k$ and $\k^*$ with
$d$ in its domain.

\ifbool{PUP}{}{\bigskip}
\subsubsection*{\textup{(2)} 
Given $\gamma\in \Gamma$, arranging that 
$\gamma\in \Gamma_{\mathbf{E}}$.} This follows in the same way.

\ifbool{PUP}{}{\bigskip}
\subsubsection*{\textup{(3)} 
Arranging $\k_{\mathbf{E}}=\pi(\ca{O}_E)$.} Suppose 
$d\in \k_{\mathbf{E}}$, $d\notin \pi(\ca{O}_E)$; set 
$e:= f_{\re}(d)$. 

Assume first that $d$ is $\d$-transcendental over $\pi(\ca{O}_E)$.
Pick $a\in \ca{O}$ and $b\in \ca{O}^*$ such
that $\overline{a}=d$ and $\overline{b}=e$. Then 
$v(E\<a\>^\times)=v(E^\times)$, and
Lemmas~\ref{univdiftr} and~\ref{ac1} yield a 
good map~$\ag=(g, f_{\re}, f_{\v})$ with small domain
$(E \langle a \rangle, \k_{\mathbf{E}}, \Gamma_{\mathbf{E}})$  
such that $\ag$ extends~$\af$ and $g(a)=b$.

Next, assume that $d$ is $\d$-algebraic
over $\pi(\ca{O}_E)$. By introducing a
% Let $\overline{P}(Y)\in \pi(\ca{O}_E)\{Y\}$ be a
minimal annihilator of $d$ over $\pi(\ca{O}_E)$, Lemmas~\ref{resembdh} and~\ref{ac1}
provide an element $a\in \ca{O}$ with $\res(a)=d$,
$\pi(\ca{O}_{E\<a\>})=\pi(\ca{O}_E)\<d\>$, and 
$v(E\<a\>^\times)=v(E^\times)$, and a good map $(g, f_{\re}, f_{\v})$ extending $\af$ 
with small domain $(E \langle a \rangle, \k_{\mathbf{E}},  \Gamma_{\mathbf{E}})$.

\bigskip\noindent
By iterating (3) we can arrange $\k_{\mathbf{E}}=\pi(\ca{O}_E)$; this condition 
is actually preserved in each of the extension procedures (4)--(8) below, 
as the reader may easily verify. We do assume in the rest of the proof that  
$\k_{\mathbf{E}}=\pi(\ca{O}_E)$. Let us say that $\mathbf{E}$ 
{\em has many constants\/} 
if $v(C_E^{\times})= v(E^\times)$.

\ifbool{PUP}{}{\bigskip}

\subsubsection*{\textup{(4)} Extending $\af$ to a good map whose domain has many constants.}
Let $\beta \in v(E^\times)$. Pick $b \in E^{\times}$ such
that $v(b)=\beta$. Since $\mathbf{K}$ has many constants, we can use 
Lemma~\ref{fixr}
to get $a \in K$ such that $a\asymp 1$ and $P(a)=0$ where
$$P(Y)\ :=\  Y'- b^{\dagger}Y\in \mathcal{O}_E\{Y\}.$$ Note that
$v(qa)=0$ and $P(qa)=0$ for all $q \in \mathbb{Q}^{\times} \subseteq
E^{\times}$. Hence by saturation we can arrange that
$\overline{a}$ is transcendental
over $\k_{\mathbf{E}}$. Then $\overline{P}(Y)$ is a minimal annihilator
of $\overline{a}$ over $\k_{\mathbf{E}}$. 
By Lemma~\ref{lem:gauss},
$$E\langle a \rangle=E(a), \qquad v\big(E(a)^\times\big)=v(E^\times), \qquad \pi\big(\ca{O}_{E(a)}\big)=\k_{\mathbf{E}}(\overline{a}).$$
We shall find a good map extending $\af$ with domain 
$\big(E(a), \k_{\mathbf{E}}(\overline{a}), \Gamma_{\mathbf{E}}\big)$. 
Consider the
differential polynomial $Q:= f(P)$, that is,
$$Q(Y)\ =\  Y' - f(b)^{\dagger}Y.$$
By saturation we
can find $e \in \k^*$ with $\overline{Q}(e)=0$ and a differential field
isomorphism $g_{\re}\colon \k_{\mathbf{E}}(\overline{a}) \to
\k_{\mathbf{E}^*}(e)$ that extends~$f_{\re}$, sends $\overline{a}$
to $e$ and is elementary as a partial map between the
differential fields $\k$ and $\k^*$. 
Using again Lemma~\ref{fixr} we find  $a^* \in K^*$
such that $a^*\asymp 1$ and $Q(a^*)=0$. Since 
$\overline{Q}(\res(a^*))=\overline{Q}(e)=0$, we can
multiply $a^*$ by an element in $C^*$ of valuation zero and
arrange $\res(a^*)=e$.
Then Theorem~\ref{resext} and Lemma~\ref{ac2} yield a 
good map $\mathbf{g}=(g, g_{\re}, f_{\v})$ where 
$g\colon E(a) \to E^*(a^*)$ extends~$f$ and sends $a$ to
$a^*$. The domain 
$\big(E(a), \k_{\mathbf{E}}(\overline{a}), \Gamma_{\mathbf{E}}\big)$  of~$\mathbf{g}$ is  small, and $vc=\beta$ for $c:=a^{-1}b\in C_{E(a)}$.

\bigskip\noindent
In the extension procedures (1)--(4) the
value group $v(E^\times)$ does not change, so if the domain $\mathbf{E}$ of $\af$
has many constants, then so does the domain of the extension of $\af$
constructed in each of (1)--(4). Also $\Gamma_{\mathbf{E}}$ does not change in
(1),~(3), and~(4), but at this stage we allow $\Gamma_{\mathbf{E}}\ne v(E^\times)$. 

\ifbool{PUP}{}{\bigskip}

\subsubsection*{\textup{(5)}
Arranging that $\mathbf{k}_{\mathbf{E}} = \pi(\mathcal{O}_E)$ is linearly surjective and $\mathbf{E}$ has many constants.} This can be done by repeated applications of (1), (3), and (4).

\ifbool{PUP}{}{\bigskip}

\subsubsection*{\textup{(6)} 
Arranging that $\mathbf{k}_{\mathbf{E}} = \pi(\mathcal{O}_E)$ and~$\mathbf{E}$ is $\d$-henselian \textup{(}by which we
mean that~$E$ as a valued differential
subfield of~$\mathbf{K}$ is $\d$-henselian\textup{)}.} 
After arranging (5) above, we can additionally arrange that $\mathbf E$ is $\d$-henselian by Corollary~\ref{immsat},  Theorem~\ref{unique.max.imm.ext}, and Lemma~\ref{ac2 variant}.

\ifbool{PUP}{}{\bigskip}

\subsubsection*{\textup{(7)} 
Towards arranging $\Gamma_{\mathbf{E}}=v(E^\times)$; the case of no torsion
modulo $v(E^\times)$.} 

Suppose $\gamma\in \Gamma_{\mathbf{E}}$ has no torsion
modulo $v(E^\times)$, that is, $n\gamma\notin v(E^\times)$ for all $n\geq 1$.
Take $a \in C$ such that
$v(a)=\gamma$. Let $i$ be a lifting 
of the differential residue field $\k$ to~$\mathbf{K}$. Since $\ac(a)\in C_{\k}^\times$,
we have $a/i\big(\!\ac(a)\big)\in C$ and $v\big(a/i\big(\!\ac(a)\big)\big)=\gamma$. 
So replacing~$a$ by $a/i\big(\!\ac(a)\big)$ we arrange that 
$v(a)=\gamma$ and $\ac(a)=1$. In the same way we obtain $a^* \in C_{K^*}$ such
that $v^*(a^*)=\gamma^*:= f_{\v}(\gamma)$ and $\ac^*(a^*)=1$. Then
Lemma~\ref{lem:lift value group ext} gives an
isomorphism $g\colon E(a) \to E^*(a^*)$ of valued fields extending~$f$ with $g(a)=a^*$. Then $(g, f_{\re},
f_{\v})$ is a good map extending $\af$ with small domain $\big(E(a),\k_{\mathbf{E}}, \Gamma_{\mathbf{E}}\big)$;
this domain has many constants if $\mathbf{E}$ does.

\ifbool{PUP}{}{\bigskip}
\subsubsection*{\textup{(8)} 
Towards arranging $\Gamma_{\mathbf{E}}=v(E^\times)$; the case of prime torsion
modulo $v(E^\times)$.} Here we assume that $\mathbf{E}$ has many constants and is $\d$-henselian. 

Let $\gamma\in \Gamma_{\mathbf{E}}\setminus v(E^\times)$ with
$\ell\gamma \in v(E^\times)$, where $\ell$ is a prime number. 
As $E$ has many constants
we can pick $b \in C_E$ such that $v(b)=\ell\gamma$. 
Since $E$ is
$\d$-henselian we have a lifting of its differential residue field 
$\k_{\mathbf{E}}$ to $\mathbf{E}$ and we can use this as in~(7) 
to arrange that
$\ac(b)=1$. We shall find $c \in C$ such that $c^\ell=b$ and $\ac(c)=1$. 
As in (7) we have $a \in C$ such that $v(a)=\gamma$ and 
$\ac(a)=1$. Then the polynomial
$P(Y):=Y^\ell-b/a^\ell\in \mathcal{O}[Y]$ satisfies $P(1)\prec 1$ and
$P'(1)\asymp 1$. This gives
$u \in K$ such that $P(u)=0$ and $u\sim 1$. Now let $c=au$.
Clearly $c^\ell=b$ (so $c\in C$) and $\ac(c)=1$.  Likewise we find $c^* \in C_{K^*}$ such that ${c^*}^\ell=f(b)$ 
and $\ac^*(c^*)=1$. Then Lemma~\ref{lem:pth root} gives an isomorphism $g\colon E(c) \to E^*(c^*)$ of valued fields extending~$f$ with $g(c)=c^*$, and so $(g, f_{\re},
f_{\v})$ is a good map extending $\af$ with small domain $\big(E(c),\k_{\mathbf{E}}, \Gamma_{\mathbf{E}}\big)$; this domain has many constants.

\bigskip\noindent
By iterating (7) and (8) we can assume in the rest of the proof that  
$\Gamma_{\mathbf{E}}=v(E^\times)$, and we shall do so. This condition is actually
preserved in the earlier extension procedures~(3) and~(4), as the reader 
may easily verify. Anyway, we can refer from now on to
$\Gamma_{\mathbf{E}}$ as the {\em value group\/} of $E$. 
Note that in (7) and (8) the
differential residue field does not change.

\bigskip\noindent
Let $a \in K$ be given. We need to extend $\af$ to a good map
whose domain is small and contains $a$. At this stage we can assume 
$\k_{\mathbf{E}}=\pi(\ca{O}_E)$, $\Gamma_{\mathbf{E}}=v(E^\times)$, and~$\mathbf{E}$ has many constants. By Lemma~\ref{lem:trdeg inequ},
$\big|\pi\big(\ca{O}_{E\<a\>}\big)\big|< \kappa$ and $\big|v\big(E\<a\>^\times\big)\big| < \kappa$.
Then~(1)--(8) and Corollary~\ref{immsat} plus
Theorem~\ref{unique.max.imm.ext} allow us to extend $\af$ to a good map
$\af_1=(f_1, f_{1,\re}, f_{1,\v})$ with small domain $\mathbf{E}_1\supseteq \mathbf{E}$ 
such that 
$\mathbf{E}_1=(E_1, \k_1, \Gamma_1)$ has many constants, $\k_1$ is linearly 
surjective, and
$$
\pi\big(\ca{O}_{E\<a\>}\big)\ \subseteq\ \k_1\ =\ \pi\big(\ca{O}_{E_1}\big), \qquad  v\big(E\<a\>^\times\big)\  \subseteq\ \Gamma_1\ =\ v(E_1^\times).$$
As we extended $\af$ to $\af_1$, we extend $\af_1$
to a good map $\af_2$ with small domain $\mathbf{E}_2\supseteq \mathbf{E}_1$ such that
$\mathbf{E}_2=(E_2, \k_2, \Gamma_2)$ has many constants, $\k_2$ linearly 
surjective, and
$$\pi\big(\ca{O}_{E_1\<a\>}\big)\ \subseteq\ \k_2\ =\ \pi\big(\ca{O}_{E_2}\big), \qquad
 v\big(E_1\<a\>^\times\big)\ \subseteq\ \Gamma_2\ =\ v(E_2^\times).$$ 
Continuing this way
and taking a union of the resulting domains and good maps gives a small
good substructure $\mathbf{E}_{\infty}\supseteq \mathbf{E}$, such that
$\mathbf{E}_{\infty}=(E_{\infty}, \k_{\infty}, \Gamma_{\infty})$ has many constants,
$\k_{\infty}$ is linearly surjective, and
$$
\k_{\infty}\ =\ \pi\big(\ca{O}_{E_{\infty}\<a\>}\big)\ =\ \pi(\ca{O}_{E_{\infty}}),\qquad \Gamma_{\infty}\ =\ v\big(E_{\infty}\<a\>^\times\big)\ =\ v(E_{\infty}^\times),$$
together with an extension of $\af$ to a good 
map $\af_{\infty}=(f_{\infty},\dots)$ with domain $\mathbf{E}_{\infty}$. Thus the valued
differential subfield $E_{\infty}\<a\>$ of $\mathbf{K}$ is an 
immediate extension of its 
valued differential subfield $E_{\infty}$. 
By Corollary~\ref{immsat} the valued differential subfield~$E_{\infty}\<a\>$ of~$\mathbf{K}$ has a maximal immediate
valued differential field extension $F$ inside~$\mathbf{K}$.    Then~$F$ is a maximal immediate extension of $E_{\infty}$ as well.
 This gives
a good substructure $\mathbf{F}=(F,\Gamma_{\infty}, \k_{\infty})$ of~$\mathbf{K}$. Likewise, the valued differential subfield $f_{\infty}(E_{\infty})$ of
$\mathbf{K}^*$ has a maximal immediate valued differential field extension~$F^*$ in~$\mathbf{K}^*$.
Use
Theorem~\ref{unique.max.imm.ext} and Lemma~\ref{ac2} to extend
$\af_{\infty}$ to a good map $\mathbf{F}\to \mathbf{F}^*=(F^*,\dots)$, and 
note that~${a\in F}$. 
\end{proof}

\section{Relative Quantifier Elimination}\label{cqe1}

\noindent
Here we derive various consequences of the Equivalence Theorem of 
Section~\ref{et}. 
Let~$\mathcal{L}$ be the 
three-sorted language of valued fields, with 
sorts $\f$ (the field sort), $\re$ (the residue sort), and $\v$ (the value group sort). 
This language consists of two copies of the one-sorted language $\mathcal L_{\operatorname{R}}=\{ 0, 1, {-}, {+}, {\,\cdot\,}\}$ of rings, one  
in the sort~$\f$ and one in the sort~$\re$, a copy of the language $\mathcal L_{\operatorname{OA}} = 
\{ {\leq}, 0, {-}, {+} \}$ of ordered abelian groups in the sort $\v$,  and function symbols $v$ of sort $\f\v$ (for the valuation) and~$\pi$ (for the residue morphism) of sort $\f\re$.
 We view a valued field $(K, \k, \Gamma;\dots)$ as an 
$\mathcal{L}$-structure in the natural way, with $\f$-variables ranging over $K$, 
$\re$-variables
over~$\k$, and $\v$-variables over $\Gamma$. We augment~$\mathcal{L}$ 
with a function symbol $\der$ of sort~$\f\f$, a function symbol $\overline{\der}$ of sort $\re \re$, and a function symbol $\ac$ of sort $\f\re$ to get the language
$\mathcal{L}(\der, \overline{\der}, \ac)$ of $\ac$-valued differential fields.
If we do not indicate otherwise, then in this section
$$\mathbf{K}\ =\ (K, \k, \Gamma; \dots), \qquad \mathbf{K}^*\ =\ (K^*, \k^*, \Gamma^*; \dots)$$
are $\d$-henselian monotone $\ac$-valued differential fields with
many constants;
they are considered as $\ca{L}(\der, \overline{\der}, \ac)$-structures in the obvious way.

\begin{cor}\label{comp00}  
$\mathbf{K} \equiv \mathbf{K}^*$ if and only if $\k \equiv \k^*$ as differential fields  and  $\Gamma \equiv \Gamma^*$ as ordered abelian groups.
\end{cor}
\begin{proof} The ``only if'' direction is obvious. Suppose $\k \equiv \k^*$ as differential fields, and $\Gamma \equiv \Gamma^*$ as ordered groups. This gives 
good substructures 
$\mathbf{E}:=(\mathbb{Q}, \mathbb{Q}, \{0\})$ of~$\mathbf{K}$, and $\mathbf{E}^*:=\big(\mathbb{Q},\mathbb{Q}, \{0\}\big)$ of $\mathbf{K}^*$, and an obvious good map $\mathbf{E} \to \mathbf{E}^*$. Now apply Theorem~\ref{embed11}.
\end{proof}

\noindent
Thus $\mathbf{K}$ is elementarily equivalent to the Hahn differential field
$\k\(( t^\Gamma\)) $ with angular component map as in 
Section~\ref{ang}. 

\begin{cor}\label{comp01} Suppose 
$\mathbf{E}=(E, \k_E, \Gamma_E;\dots)\subseteq \mathbf{K}$ is a $\d$-henselian $\ac$-valued differential subfield of $\mathbf{K}$ with many constants, such that $\k_E\preccurlyeq \k$ as differential fields, and $\Gamma_E \preccurlyeq \Gamma$ as ordered abelian groups. Then $\mathbf{E} \preccurlyeq \mathbf{K}$.
\end{cor}
\begin{proof} Take an elementary extension $\mathbf{K}^*$ of 
$\mathbf{E}$. Then $\mathbf{K}^*$ has many constants, $(E, \k_E, \Gamma_E)$ is a good substructure of both $\mathbf{K}$ and
$\mathbf{K}^*$, and the identity on $(E, \k_E, \Gamma)$ is a good map. Hence by Theorem~\ref{embed11} we have
$\mathbf{K} \equiv_{\mathbf{E}} \mathbf{K}^*$. Since $\mathbf{E} \preccurlyeq \mathbf{K}^*$, this gives $\mathbf{E}\preccurlyeq \mathbf{K}$.
\end{proof}

\noindent
The proofs of these corollaries use only a weak form of the 
Equivalence Theorem, but now we
turn to a result that uses its full strength: a relative elimination of 
quantifiers 
for the $\ca{L}(\der, \overline{\der}, \ac)$-theory $T$ of 
$\d$-henselian monotone $\ac$-valued differential
fields with many constants.
We specify that the function symbols $\pi$ and $v$ of $\ca{L}(\der,\overline{\der}, \ac)$ 
are to be 
interpreted as {\em total\/} functions in any $\mathbf{K}$ as follows: extend 
$\pi\colon \ca{O} \to \k$ to $\pi\colon K \to \k$ by $\pi(a)=0$ for $a\notin \ca{O}$, and extend 
$v\colon K^\times \to \Gamma$  to $v\colon K \to \Gamma$ by $v(0)=0$.
Let $\mathcal{L}_{\re}$ be the sublanguage of 
$\mathcal{L}(\der,\overline{\der}, \ac)$ involving only the sort $\re$, that is,
$\mathcal{L}_{\re}$ is a copy of the language of differential fields, 
with $\overline{\der}$ as its
symbol for the derivation operator. Let $\mathcal{L}_{\v}$ be the 
sublanguage of 
$\mathcal{L}(\der,\overline{\der}, \ac)$ involving only the sort $\v$, that is,
$\mathcal{L}_{\v}=\mathcal L_{\operatorname{OA}}$ is
the language of ordered abelian groups.
 
Let $x=(x_1, \dots, x_l)$ be 
a tuple of distinct $\f$-variables, 
$y=(y_1, \dots, y_m)$ a tuple of distinct $\re$-variables, 
and $z=(z_1, \dots, z_n)$ a tuple of distinct $\v$-variables. 
Set 
$$\Z\{x\}\ :=\ \big\{P\in \Q\{x_1,\dots, x_l\}:\ \text{all coefficients of $P$ are in $\Z$}\big\}.$$
Define a 
{\em special $\re$-formula in $(x,y)$\/} to be an $\ca{L}(\der,\overline{\der}, \ac)$-formula 
$$\psi(x,y)\ :=\ \psi_{\re}\big(\!\ac(P_1(x)),\dots, \ac(P_k(x)),y\big)$$
where $k\in \mathbb{N}$, $\psi_{\re}(u_1,\dots, u_k,y)$ is an $\mathcal{L}_{\re}$-formula, 
and $P_1(x),\dots, P_k(x)\in \mathbb{Z}\{x\}$. Also, a {\em special $\v$-formula in $(x,z)$\/} is an $\ca{L}(\der,\overline{\der}, \ac)$-formula 
$$\theta(x,z)\  :=\ \theta_{\v}\big(v(P_1(x)),\dots, v(P_k(x)),z\big)$$
where $k\in \mathbb{N}$, $\theta_{\v}(v_1,\dots, v_k,z)$ is an $\mathcal{L}_{\v}$-formula, 
and $P_1(x),\dots, P_k(x)\in \mathbb{Z}\{x\}$. Note that these special formulas do not
have quantified $\f$-variables. We can now state our relative quantifier
elimination:

\begin{cor}\label{qe} Every $\ca{L}(\der,\overline{\der}, \ac)$-formula
$\ \phi(x,y,z)\ $ is $T$-equivalent to
$$ \big(\psi_1(x,y)\wedge \theta_1(x,z)\big)\vee \cdots \vee \big(\psi_N(x,y)\wedge \theta_N(x,z)\big)$$
for some $N\in \N$ and some special
$\re$-formulas $\psi_1(x,y),\dots, \psi_N(x,y)$ in $(x,y)$, and some special $\v$-formulas $\theta_1(x,z),\dots, \theta_N(x,z)$ in $(x,z)$. 
\end{cor}
\begin{proof} Let $\Theta(x,y,z)$ be the set of $\ca{L}(\der,\overline{\der}, \ac)$-formulas 
$$ \big(\psi_1(x,y)\wedge \theta_1(x,z)\big)\vee \cdots \vee \big(\psi_N(x,y)\wedge \theta_N(x,z)\big)$$ displayed
in the statement of the corollary. It is clear that $\Theta$ is closed
under taking disjunctions, and easy to check that $\Theta$ is closed
under taking negations, modulo logical equivalence. Thus by 
Corollary~\ref{cor:separation} it is enough to show that every
$T$-realizable $(x,y,z)$-type is completely 
determined by its intersection with $\Theta$. This guides the argument that follows.
  
Let $\psi(x,y)$ and $\theta(x,z)$ range over special formulas as described
above. For $\mathbf{K}=(K, \Gamma, \k;\dots)\models T$ and 
$a\in K^l$, $r \in \k^m$, $\gamma \in \Gamma^n$, let
\begin{align*} \tp_{\re}^{\mathbf{K}}(a,r)\ &:=\ 
\big\{\psi(x,y):\  \mathbf{K} \models \psi(a,r)\big\}, \\
  \tp_{\v}^{\mathbf{K}}(a,\gamma)\ &:=\ \big\{\theta(x,z):\ 
\mathbf{K} \models \theta(a,\gamma)\big\}.
\end{align*}
Let $\mathbf{K}$ and $\mathbf{K}^*$ be any models of $T$, and let
$$(a,r,\gamma)\in K^l\times \k^m\times \Gamma^n, \qquad (a^*,r^*,\gamma^*)\in 
(K^*)^l\times (\k^*)^m\times (\Gamma^*)^n$$ be such that
$\tp^{\mathbf{K}}_{\re}(a,r)=\tp^{\mathbf{K}^*}_{\re}(a^*,r^*)$ and $\tp^{\mathbf{K}}_{\v}(a,\gamma)=\tp^{\mathbf{K}^*}_{\v}(a^*,\gamma^*)$. 
By the above, it suffices to show that under these assumptions we have 
$$\tp^{\mathbf{K}}(a, r, \gamma)\ =\ \tp^{\mathbf{K}^*}(a^*, r^*, \gamma^*).$$
%(See Corollary~\ref{cor:separation}.)
Let $\mathbf{E}:= (E,\Gamma_{\mathbf{E}}, \k_{\mathbf{E}})$ where 
$E:= \mathbb{Q}\langle a \rangle$, $\k_{\mathbf{E}}$ is the differential subfield of
$\k$ generated by $\ac(E)$ and $r$, and $\Gamma_{\mathbf{E}}$ is the ordered subgroup of $\Gamma$ generated by
$\gamma$ over $v(E^\times)$, so~$\mathbf{E}$ is a good substructure of
$\mathbf{K}$. Likewise we define the good substructure $\mathbf{E}^*$ of~$\mathbf{K}^*$.
For each $P(x)\in \mathbb{Z}\{x\}$ we have $P(a)=0$ iff $\ac\big(P(a)\big)=0$, and also 
$P(a^*)=0$ iff $\ac^*\!\big(P(a^*)\big)=0$. 
In view of this fact, the assumptions give us a good map $\mathbf{E} \to \mathbf{E}^*$ 
sending $a$ to
$a^*$, $\gamma$ to $\gamma^*$ and $r$ to $r^*$. It remains to apply 
Theorem~\ref{embed11}.
\end{proof}

\noindent
In the proof of the corollary above it is important that 
our notion of a good substructure~$\mathbf{E}=(E, \Gamma_{\mathbf{E}}, \k_{\mathbf{E}})$ 
did not require $\Gamma_{\mathbf{E}}=v(E^\times)$ or $ \k_{\mathbf{E}}=\pi(\ca{O}_E)$.
Related to it is that in Corollary~\ref{qe} we have a separation of 
$\re$- and $\v$-variables;  this
makes the next result almost obvious. 

\begin{cor}\label{comp02} Each subset of $\k^m\times \Gamma^n$ definable in $\mathbf{K}$ is a finite union of rectangles $Y\times Z$ with $Y\subseteq \k^m$
definable in the differential field $\k$ and $Z\subseteq \Gamma^n$ definable in 
the
ordered abelian group $\Gamma$.
\end{cor}
\begin{proof} By Corollary~\ref{qe} and using its notations it is enough to
observe that for $a\in K^l$, a
special $\re$-formula $\psi(x,y)$ in $(x,y)$, and a
special $\v$-formula $\theta(x,z)$ in $(x,z)$, 
the set $\big\{r\in \k^m: \mathbf{K}\models \psi(a,r)\big\}$ is definable
in the differential field $\k$, and the set $\big\{\gamma\in \Gamma^n: \mathbf{K}\models \theta(a,\gamma)\big\}$ 
is definable in the ordered abelian group $\Gamma$.   
\end{proof}

\noindent
Corollary~\ref{comp02} says in particular 
that the relations on $\k$ definable in 
$\mathbf{K}$ are
definable in the differential field $\k$, and likewise, the relations on 
$\Gamma$ definable in $\mathbf{K}$ are definable in the ordered abelian group 
$\Gamma$. 

\medskip\noindent
Theorems~\ref{akesc} and~\ref{defwitt} from the introduction to this chapter 
do not mention
angular component maps. To get these results from Corollaries~\ref{comp00} 
and~\ref{comp02} we first pass to suitable $\aleph_1$-saturated 
elementary extensions and then use
Corollary~\ref{acm3} to get the necessary angular component
maps. 

To get Theorem~\ref{pelsub} from Corollary~\ref{comp01}, we 
arrange likewise that $\mathbf{K}$ and $\mathbf{K^*}$ from that theorem
are $\aleph_1$-saturated (but not yet equipped with angular component maps).
Then we have a cross-section $s\colon \Gamma \to C^{\times}$ of the valued subfield
$C$ of $\mathbf{K}$. Use Lemma~\ref{xs2} 
to extend $s$ to a cross-section $s^*\colon \Gamma^* \to C_{K^*}^{\times}$ 
of the valued subfield~$C_{K^*}$ of $\mathbf{K}^*$. These cross-sections yield
angular component maps on the valued fields~$C$ and~$C^*$, which by
Lemma~\ref{acm1} extend uniquely to angular component maps on $\mathbf{K}$ and
$\mathbf{K}^*$. This allows us to use Corollary~\ref{comp01} to get
 $\mathbf{K} \preccurlyeq \mathbf{K}^*$.

\subsection*{Notes and comments}
Corollary~\ref{qe} is analogous to a result by Pas~\cite{Pas89}
for henselian valued fields of equicharacteristic zero.
Readers familiar with the model-theoretic properties of \textit{stable embeddedness}\/ and
\textit{orthogonality}\/ will observe that by
Corollary~\ref{comp02}, $\k$ and $\Gamma$ are stably embedded in $\mathbf{K}$, and 
$\k$ and $\Gamma$ are orthogonal in $\mathbf{K}$. 
We refer to Appendix~\ref{app:modth} for the definition of the model-theoretic property NIP, the
{\em Non-Independence Property}.
Using Corollary~\ref{qe} one can show (along the lines of \cite[proof of Theorem~A.15]{SimonNIP}) that 
if the differential residue field $\k$ has NIP 
(that is, its theory has NIP), 
then so does~$\mathbf{K}$; this also uses the fact that by~\cite{GurevichSchmitt} every ordered abelian group has NIP.

\section{A Model Companion}\label{sec:modelcompanionmonotone}

\noindent
Let $\mathcal L_{\der, \preceq}:= \{0,1,{-},{+}, {\,\cdot\,}, {\der}, {\preceq}\}$ be the one-sorted language of valued differential fields, and construe valued differential fields as $\mathcal L_{\der, \preceq}$-structures in the obvious way. 

\begin{prop}
The theory of $\d$-henselian valued differential fields with many constants, differentially closed differential residue field,
and nontrivial divisible value group is complete and model complete. It
is the model companion of the theory of monotone valued differential fields.
\end{prop}
\begin{proof}
For completeness and model completeness, use Theorems~\ref{akesc} and~\ref{pelsub}
in combination with
Corollary~\ref{cor:DCF-qe} and Example~\ref{ex:DOAb}.
Let $K$ be an arbitrary monotone valued differential field; to prove the 
model companion claim, it is enough to embed~$K$ into a $\d$-henselian valued differential field with many constants, differentially closed differential residue field,
and nontrivial divisible value group. By Example~(1) at the beginning of Section~\ref{p:valued diff field},  $K$ has a
monotone valued differential field extension~$K_1$ with a
nontrivial value group. By Corollary~\ref{cor:monalgmon}, $K_1$ has a monotone valued differential field extension $K_2$ with divisible value group. 
By Corollary~\ref{cor:monotonicity under extensions}, $K_2$
has a monotone valued differential field extension $K_3$ with differentially closed differential residue field and $\Gamma_{K_2}=\Gamma_{K_3}$.
Theorem~\ref{thm:damdh} then yields
an immediate $\d$-henselian valued differential field
extension $K_4$ of $K_3$. Then~$K_4$ is still monotone by Corollary~\ref{cor:sc2}, and thus has many constants by Corollary~\ref{cor:1dconst}. 
\end{proof}

%% file: mt-9n.tex
\chapter{Asymptotic Fields and Asymptotic Couples}\label{ch:asymptotic differential fields}

\noindent
The key restriction on valued differential fields in Chapter~\ref{ch:valueddifferential}
was the {\em continuity\/} of the derivation. (Strictly speaking, we assumed
the derivation $\der$ to be {\em small}, but continuity of $\der$
reduces by compositional conjugation to smallness of $\der$.)

In this chapter we introduce asymptotic differential fields:
valued differential
fields with a much stronger interaction of the valuation 
and derivation.
For brevity we just call them {\em asymptotic fields}.
They include Rosenlicht's differential-valued fields~\cite{Rosenlicht2} 
and share 
many of their basic properties. The advantage of
the class of asymptotic fields over its subclass of differential-valued 
fields is that the former is closed under taking 
valued differential subfields,
under coarsening, and even under specialization 
(subject to a mild restriction).

A key feature of an asymptotic field is its asymptotic couple, which is 
just its value group
with some extra structure induced by the derivation.  
In Section~\ref{As-Fields,As-Couples}
we define asymptotic fields, their asymptotic couples, and discuss
Hardy fields. In 
Section~\ref{AbstractAsymptoticCouples} we consider asymptotic couples  
independent of their connection to asymptotic fields. This is used
in Section~\ref{Applicationtodifferentialpolynomials} to describe 
the behavior of differential polynomials as functions on asymptotic fields.
In Section~\ref{Asymptotic-Fields-Basic-Facts} we
consider asymptotic fields with small derivation
and the operations of coarsening and specialization. 
In Section~\ref{pafae} we show that algebraic extensions of 
asymptotic fields are asymptotic. In Section~\ref{ImmExtAs} we 
adapt the results on immediate extensions from Section~\ref{sec:cimex}
to asymptotic fields. Section~\ref{caseorderone} treats differential 
polynomials of order one over $H$-asymptotic fields.
In Section~\ref{sec:extH} we return to asymptotic couples
and prove some useful extension results about them, and in Section~\ref{sec:cac}
we establish a property of {\em closed\/} $H$-asymptotic couples as needed
in Chapter~\ref{ch:QE}. 
The present chapter and the next include some new material but are 
mainly based on ~\cite{Rosenlicht2,AvdD,AvdD2,AvdD3}.

\medskip\noindent
Some terminology: When a valued differential field $K$ is given, then an {\em extension of~$K$\/} is a valued differential field extension of $K$. Likewise, an {\em extension\/} of an ordered valued differential field $K$ is an ordered valued differential field extension of $K$. The term ``embedding'' is used in a similar way:
when $K$ and $L$ are given as valued differential fields, then an
{\em embedding\/} of $K$ into $L$ is an embedding of valued differential fields, and when $K$ and $L$ are given as ordered valued 
differential fields, then an
{\em embedding\/} of $K$ into $L$ is an embedding of ordered valued differential fields. 

\label{p:ovdf}
\index{extension!valued differential fields}
\index{extension!ordered valued differential fields}

\input{mt-9n-1}
\input{mt-9n-2}

\input{mt-9n-3}

\input{mt-9n-4}

\input{mt-9n-5}

\input{mt-9n-6}

\input{mt-9n-7}

\input{mt-9n-8}
\input{mt-9n-9}

%% file: mt-9n-1.tex
\section{Asymptotic Fields and Their Asymptotic Couples}\label{As-Fields,As-Couples}

\noindent
In the first subsection we define asymptotic fields and differential-valued
fields, in the second subsection we
show how to visualize an asymptotic couple, in the third subsection
we introduce the asymptotic couple of an asymptotic field,
and in the fourth subsection we define comparability classes and the property 
of being grounded. In the last subsection we discuss Hardy fields as
examples.

\subsection*{Asymptotic fields} 
An {\bf asymptotic differential field}, or just {\bf asymptotic field\/}, 
is a valued differential field~$K$  such that for all 
$f,g\in K^\times$ with $f,g\prec 1$,
\begin{itemize}
\item[(A)] $f\prec g \Longleftrightarrow f'\prec g'$.
\end{itemize}
If in addition we have for all $f,g\in K^\times$ with $f,g\prec 1$,
\begin{itemize}
\item[(H)] $f\prec g \Longrightarrow f^\dagger\succeq g^\dagger$,
\end{itemize}
then we say that $K$ is an {\bf $H$-asymptotic field\/} \label{p:H-asymptotic field} or an asymptotic field
of {\bf $H$-type}. Our main interest is in $H$-asymptotic fields, 
but many things go through without the $H$-type assumption.
%An {\bf asymptotic extension} of an asymptotic field $K$ is by definition 
%a valued differential field extension of $K$ that is asymptotic.  

\index{extension!asymptotic}
\index{asymptotic extension}
\index{valued differential field!asymptotic}
\index{field!asymptotic}
\index{asymptotic field}
\index{asymptotic field!H-type@$H$-type}
\index{H-type@$H$-type!asymptotic field}
\index{H-asymptotic@$H$-asymptotic!field}

\begin{lemma}\label{neatderas} Let $K$ be a valued differential field such that 
$C\subseteq \mathcal{O}$, $\der \smallo\subseteq \smallo$, and~$\der$ is 
neatly surjective. Then $K$ is asymptotic.
\end{lemma}
\begin{proof} From $C\subseteq \mathcal{O}$ we obtain $C\cap \smallo=\{0\}$, so
the restriction of $\der$ to $\smallo$ is injective. From $\der \smallo\subseteq \smallo$ we obtain a strictly increasing map 
$v_{Y'}=v_{\der}: \Gamma \to \Gamma$ with $v_{\der}(0)=0$. 
Let $f,g\in \smallo^{\ne}$.
Then $f', g'\in \smallo^{\ne}$, so $v_{\der}(vf)=v(f')$ and $v_{\der}(vg)=v(g')$ by
the neat surjectivity of $\der$.
Thus $f\prec g\Longleftrightarrow f'\prec g'$. 
\end{proof}

\noindent
Suppose $K$ is an asymptotic field. Then
$C\cap \smallo=\{0\}$: if $0\ne c\in C\cap \smallo$, 
then $c^2\prec c\prec 1$, so $0=(c^2)'\prec c'=0$ by
(A), which is impossible. Thus the valuation $v$ is trivial on $C$, and
the map $\varepsilon\mapsto\varepsilon'\colon \smallo\to K$ is injective.
In particular, $C\subseteq \mathcal{O}$ and 
the residue map $a \mapsto \bar{a}\colon \mathcal{O} \to \k$ is injective on 
$C$.
The following three conditions on $K$ are clearly equivalent: \begin{enumerate}
\item $\mathcal{O}=C+\smallo$;
\item $\{\bar{a}:\ a\in C\}=\k$;
\item for all $f\asymp 1$ in $K$ there exists $c\in C$ with $f\sim c$.
\end{enumerate}
\noindent
We say that $K$ is {\bf differential-valued\/} (or {\bf $\d$-valued}, for short) if it satisfies these 
three
(equivalent) conditions. So the constant field of a $\d$-valued field is also
a lift of its residue field. 
The next lemma is clear from (3) above.

\label{p:dv}
\index{field!differential-valued}
%\index{field!$H$-differential-valued}
\index{differential-valued!field}
\index{valued differential field!differential-valued}
\index{d-valued@$\d$-valued!field}
%\index{$H$-differential-valued field}
%\index{valued differential field!$H$-differential-valued}

\begin{lemma}\label{dv}
If $L$ is an asymptotic extension of a $\d$-valued field $K$ with
$\res(K)=\res(L)$, then $L$ is $\d$-valued, with $C_L=C$. 
\end{lemma}

\noindent
%Differential-valued fields of $H$-type are also called 
%{\em $H$-differential-valued fields}. 
Our final results in Chapter~\ref{ch:QE} concern 
just $\d$-valued fields of $H$-type such as $\mathbb{T}$, but towards this goal 
we need the wider setting of $H$-asymptotic fields 
where we can coarsen and pass to suitable
differential subfields at our convenience.

\begin{figure}
\begin{center}
\begin{picture}(300,200)

\put(0,0){\line(0,1){200}}
\put(0,0){\line(1,0){300}}

\put(20,20){\line(0,1){160}}
\put(20,20){\line(1,0){260}}

\put(40,40){\line(0,1){120}}
\put(40,40){\line(1,0){220}}

\put(60,60){\line(0,1){80}}
\put(60,60){\line(1,0){180}}

\put(0,200){\line(1,0){300}}
\put(20,180){\line(1,0){260}}
\put(40,160){\line(1,0){220}}
\put(60,140){\line(1,0){180}}

\put(300,0){\line(0,1){200}}
\put(280,20){\line(0,1){160}}
\put(260,40){\line(0,1){120}}
\put(240,60){\line(0,1){80}}

\multiput(80,80)(0,4){10}{\line(0,1){2.5}}
\multiput(218.5,80)(0,4){10}{\line(0,1){2.5}}
\multiput(80,80)(4,0){35}{\line(1,0){2.5}}

\multiput(218.5,119)(-4,0){35}{\line(-1,0){2.5}}

\put(10,5){valued differential fields}
\put(30,25){asymptotic fields}
\put(50,45){$\d$-valued fields}
\put(70,65){$\d$-valued fields of $H$-type}
\put(90,85){$H$-fields}

\end{picture}
\end{center}
\caption{Classes of valued differential fields.}
\label{fig:classes of valued differential fields}
\end{figure}

\ifbool{PUP}{}{\medskip\noindent}
Figure~\ref{fig:classes of valued differential fields} indicates the inclusion relations among some
classes of valued differential fields, with the class of $H$-fields from the 
next chapter as
the smallest class. Strictly speaking, $H$-fields are more than just valued differential fields, since they also carry an ordering; that's why we use a dashed rectangle for this class. 

Any differential subfield of an asymptotic field $K$
with the restricted dominance relation is itself
an asymptotic field. If $K$ is an asymptotic
field and $\phi\in K^{\times}$, then its compositional conjugate
$K^\phi$ with the same dominance relation remains an 
asymptotic field with the same constant field. (These two statements 
remain true with {\em $H$-asymptotic\/} in place of {\em asymptotic}.) 
If $K$ is a 
$\d$-valued field and $\phi\in K^{\times}$, then 
$K^\phi$ (with same dominance relation) remains $\d$-valued.

An asymptotic field is said to be {\bf asymptotically maximal} if it has no proper immediate asymptotic extension. Likewise, 
an asymptotic field is {\bf asymptotically $\d$-algebraically maximal} 
if it has no proper immediate $\d$-algebraic asymptotic extension.
Any asymptotic field has, by Zorn, an 
immediate asymptotic extension that is asymptotically maximal, and also an immediate $\d$-algebraic asymptotic extension that is asymptotically $\d$-algebraically maximal. 
These notions will become important in Chapter~\ref{ch:newtonian fields}.

\label{p:asymptotically maximal}
\label{p:asymptotically d-algebraically maximal}
\index{asymptotically!maximal}
\index{asymptotically!d-algebraically maximal@$\d$-algebraically maximal}
\index{asymptotic field!asymptotically maximal}
\index{asymptotic field!asymptotically d-algebraically maximal@asymptotically $\d$-algebraically maximal}

\subsection*{Asymptotic couples} We defined asymptotic couples
in Section~\ref{sec:ascouples}. Consider an asymptotic couple $(\Gamma, \psi)$.
For reasons that will become clear in the next subsection we also 
use the notation
$$ \alpha^\dagger:= \psi(\alpha), \qquad \alpha':= \alpha + \psi(\alpha) \qquad (\alpha \in \Gamma^{\ne}).$$
Of course this is only used when $\psi$ is understood from
the context. Thus $\alpha'>\beta^\dagger$ for $\alpha\in \Gamma^{>},\ \beta\in \Gamma^{\ne}$. The following subsets of $\Gamma$ play special roles:
\begin{align*}  (\Gamma^{\ne})'\ &:=\ \{\gamma':\ \gamma\in \Gamma^{\ne}\}, 
\qquad (\Gamma^{>})'\ :=\ \{\gamma':\ \gamma\in \Gamma^{>}\},\\
   \Psi\ &:=\ \psi(\Gamma^{\ne})\ =\ \{\gamma^\dagger:\ \gamma\in \Gamma^{\ne}\}\ =\ \{\gamma^\dagger:\ \gamma\in \Gamma^{>}\}.
\end{align*}
We call $\Psi$ the $\Psi$-set of $(\Gamma,\psi)$, and set $\alpha^\dagger:=\infty\in \Gamma_{\infty}$ for
$\alpha=0\in \Gamma$.  
\input{mt-psipic.tex}
%(This set is denoted by $\Psi$ in \cite{ }, { }, { }, { }.) 

\medskip\noindent
Figure~\ref{psipic} gives a rough idea of the graphs of $\gamma'$ and $\gamma^\dagger$
as functions of $\gamma\in \Gamma^{\ne}$. We have found this picture very 
helpful in getting a ``feel'' for asymptotic couples. 
It suggests for example that $\gamma'$ is 
strictly increasing on $\Gamma^{\ne}$, and this is part (iii) of Lemma~\ref{BasicProperties}.
It also suggests that $\gamma'$ has the intermediate value property on~$\Gamma^{>}$ as well as on $\Gamma^{<}$. This is the case for 
$H$-asymptotic couples by Lemma~\ref{ivp} below, but not for all asymptotic 
couples; see Example 2.8 in 
\cite{AvdD2}.   
The picture is really meant just for
$H$-asymptotic couples, where $\gamma^\dagger$ is increasing for $\gamma<0$ and
decreasing for $\gamma>0$. Of course, we are unable to make the picture show,
for $H$-asymptotic $(\Gamma, \psi)$, that
$\psi$ is constant on each archimedean class $[\alpha]$, 
$\alpha\in \Gamma^{\ne}$.  

\nomenclature[T]{$\Psi$}{$\psi(\Gamma^{\neq})$}
\nomenclature[T]{$\gamma'$}{$\gamma'=\gamma+\psi(\gamma)$, for $\gamma\neq 0$}
\nomenclature[T]{$\gamma^\dagger$}{$\gamma^\dagger=\psi(\gamma)$, for $\gamma\neq 0$}

\subsection*{The asymptotic couple of an asymptotic field}
Let $K$ be an asymptotic 
field. It is clear that then $vg'$ is uniquely determined by $vg$
for $g\in K^\times$ and $vg\ne 0$, that is, the derivation of $K$ induces a function
$$\gamma \mapsto \gamma'\ \colon\ \Gamma^{\ne} \to \Gamma   \qquad \text{($\gamma=vg$, $\gamma'=vg'$, 
 $g$ as above)}$$ 
on its value group $\Gamma$. We also consider the logarithmic-derivative analogue:
$$ \gamma \mapsto \gamma^{\dagger}:= \gamma'- \gamma\ \colon\ \Gamma^{\ne} \to \Gamma,$$
that is, $\gamma^\dagger=v(g^\dagger)$ for $g$ as above.
The asymptotic couple of $K$ is just
the value group of $K$ equipped with this induced operation 
$\gamma \mapsto \gamma^{\dagger}$. To justify this terminology, we characterize asymptotic fields as follows:

\begin{prop}\label{CharacterizationAsymptoticFields}
Let $K$ be a valued differential field.
% with dominance relation $\preceq$. 
Then the conditions below are equivalent:
\begin{enumerate}
\item[\textup{(i)}] $K$ is an asymptotic field;
\item[\textup{(ii)}] there is an asymptotic couple $(\Gamma, \psi)$ with
$\Gamma:=v(K^\times)$ such that for 
all $g\in K^\times$ with $g\nasymp 1$ we have
$\psi(vg)=v(g^\dagger)$;
\item[\textup{(iii)}] for all $f,g\in K^\times$ with
$f,g\nasymp 1$ we have: $f\preceq g \Longleftrightarrow f'\preceq g'$;
\item[\textup{(iv)}] for all $f,g\in K^\times$ we have:
$$\begin{cases}
&f\prec 1,\ g\nasymp 1 \quad\Rightarrow\quad f'\prec g^\dagger,\\
&f\asymp 1,\ g\nasymp 1 \quad\Rightarrow\quad f'\preceq g^\dagger.
\end{cases}$$
\end{enumerate}
\end{prop}
\begin{proof}
We first show the equivalence of (ii) and (iv). Suppose that (ii) holds.
The first implication of (iv) is clear. For the second implication, 
let $f,g\in K^\times$, $f\asymp 1$, $g\nasymp 1$.
We have $g^\dagger\asymp (fg)^\dagger=f^\dagger+g^\dagger$.
Hence $f'\asymp f^\dagger\preceq g^\dagger$.
Conversely, suppose that (iv) holds. Then  $g^\dagger\ne 0$ for $g\nasymp 1$ in $K^\times$ by taking  $f=g$ or $f=g^{-1}$ in the first implication of (iv). 
Let $f,g\in K^\times$ with $f\asymp g
\prec 1$. Then
$f^\dagger - g^\dagger = (f/g)^\dagger \asymp (f/g)'$ with $f/g\asymp 1$,
so $f^\dagger-g^\dagger \preceq f^\dagger$ and 
$f^\dagger-g^\dagger\preceq g^\dagger$
by the second implication of (iv). It follows that $f^\dagger\asymp g^\dagger$.
Thus $v(f^\dagger)$ only depends on $vf$, for $f\in K^{\times}$ with 
$f\nasymp 1$. Now~(ii) follows, with (AC3) a consequence of
the first implication in (iv).

Note that 
(ii)~$\Rightarrow$~(iii) is a consequence of Lemma~\ref{BasicProperties}(iii),
and (iii)~$\Rightarrow$~(i) is trivial. Lemma~\ref{CharAsymptoticCouples} gives (i)~$\Rightarrow$~(ii).
\end{proof}

\noindent
If $K$ is an asymptotic field, we call $(\Gamma,\psi)$ as
defined in (ii) of Proposition~\ref{CharacterizationAsymptoticFields} 
the {\bf asymptotic couple of $K$}. An asymptotic field is of $H$-type iff
its asymptotic couple is of $H$-type. If $(\Gamma,\psi)$ is the 
asymptotic couple of the asymptotic field $K$ and $a\in K^{\times}$, 
then $(\Gamma, \psi^a)$ with $\psi^a:= \psi-va$ is the asymptotic couple 
of $K^a$.

\index{asymptotic couple!asymptotic field}
\index{asymptotic field!asymptotic couple}
\index{asymptotic couple!shift}

\begin{cor}\label{CharacterizationAsymptoticFields-Corollary}
Let $K$ be an asymptotic field and $f,g\in K$. Then
\begin{enumerate}
\item[\textup{(i)}] if $f\prec g\nasymp 1$, then $f'\prec g'$;
\item[\textup{(ii)}] if $f\preceq g\nasymp 1$, then $f\sim g \Longleftrightarrow f'\sim g'$;
\item[\textup{(iii)}] if $f\preceq g\nasymp 1$ and $g'\preceq g$, then $f^{(n)}\preceq g'$
for all $n\ge 1$;
\item[\textup{(iv)}] if $f\asymp 1$ and $0\neq g\nasymp 1$, then
$f^\dagger\not\sim g^\dagger$.
\end{enumerate}
\end{cor}
\begin{proof}
For (i), if $f\prec g$ and $f,g\nasymp 1$, then $f'\prec g'$
by Proposition~\ref{CharacterizationAsymptoticFields}(iii), 
and if $1\asymp f\prec g$, then $f'\preceq g^\dagger=g'/g
\prec g'$ by Proposition~\ref{CharacterizationAsymptoticFields}(iv). 
For (ii), suppose that~$g\nasymp 1$.
If $f\sim g$, then $1\nasymp g\succ f-g$, so $g'\succ f'-g'$ by (i),
i.e.,~$f'\sim g'$. Conversely, suppose that $f\preceq g\nasymp 1$ and 
$f'\sim g'$. If $f-g\nasymp 1$,
then  Proposition~\ref{CharacterizationAsymptoticFields}(iii) yields
$f \sim g$. 
If $f-g\asymp 1$, then  $g\succ 1$, and thus $f-g\asymp 1\prec g$.
Part (iii) follows by induction on $n$ from (i) and
Proposition~\ref{CharacterizationAsymptoticFields}(iii).
As to (iv), if  $f\asymp 1$ and $0\neq g\nasymp 1$, then
$g/f\asymp g\nasymp 1$, so 
$(g/f)^\dagger=g^\dagger-f^\dagger \asymp
g^\dagger$, hence~$f^\dagger\not\sim g^\dagger$.
%Part (5) follows from
%(4) in Proposition~\ref{CharacterizationAsymptoticFields}.
\end{proof}

\noindent
{\em In the next two corollaries we assume that $K$ is an asymptotic field}.

\begin{cor}\label{c1} The derivation of $K$ is continuous.
\end{cor}
\begin{proof} If $\Gamma=\{0\}$, then the
valuation topology on $K$ is discrete and so $\der$
is continuous. Assume $\Gamma \ne \{0\}$,
and take $g\in K^\times$ with $g\succ 1$. Then $a:= g'\ne 0$, and
so $\der \smallo \subseteq a\smallo$ by Corollary~\ref{CharacterizationAsymptoticFields-Corollary}(i). Thus $\der$ is continuous by Lemma~\ref{cd2}.
\end{proof}

\noindent
As a consequence of Corollary~\ref{c1} and Lemma~\ref{ExtendingDerivation}, 
the completion  $K^{\operatorname{c}}$ of our asymptotic field $K$ 
is naturally a valued differential field, but more is true:

\begin{cor}\label{completingasymp} The valued differential field $K^{\operatorname{c}}$ is an asymptotic field. If $K$ is $\d$-valued, then so is $K^{\operatorname{c}}$.
\end{cor}
\begin{proof} 
We claim that~$K^{\operatorname{c}}$ is asymptotic with the same asymptotic couple
$(\Gamma, \psi)$ as $K$, and to show this
we use (ii) of Proposition~\ref{CharacterizationAsymptoticFields}. Let 
$f\in K^{\operatorname{c}}$ and
$0\ne f\nasymp 1$. Density of $K$ in $K^{\operatorname{c}}$ gives $g\in K$ with
$f\sim g$ and $f'\sim g'$. Then $\psi(vf)=\psi(vg)=v(g^\dagger)=v(f^\dagger)$. This proves our claim. If $K$ is even $\d$-valued, then so is~$K^{\operatorname{c}}$ by Lemma~\ref{dv},
since $\res{K}=\res{K^{\operatorname{c}}}$. 
\end{proof}

\index{completion!asymptotic field}
\index{asymptotic field!completion}

\noindent
We say that an asymptotic couple $(\Gamma, \psi)$ has 
{\bf  asymptotic integration }
if $$\Gamma=(\Gamma^{\ne})'.$$ An asymptotic field is said to 
have {\bf asymptotic integration\/} if its asymptotic couple has asymptotic integration.
The reason for this terminology is that if $K$ is a $\d$-valued field,
then $K$ has asymptotic integration iff for all $a\in K^\times$ there is
$b\in K$ such that $a\sim b'$. 

\label{p:asymptotic integration}
\index{asymptotic integration!}
\index{asymptotic field!with asymptotic integration}
\index{asymptotic couple!with asymptotic integration}

\subsection*{Comparability and groundedness} Let
$K$ be an asymptotic field.
On the set~$K^\times$ we
define the binary relations $\comp$,~$\flatter$,~$\flattereq$ as follows: \nomenclature[R]{$f\comp g$}{$f$ and $g$ are comparable}
\nomenclature[R]{$f\flatter g$}{$f$ is flatter than $g$}
\nomenclature[R]{$f\flattereq g$}{$f$ is flatter than or comparable with $g$}
\begin{align*}
f\comp g \quad &:\Longleftrightarrow\quad f^\dagger \asymp g^\dagger, \\
f\flatter g \quad &:\Longleftrightarrow\quad f^\dagger \prec g^\dagger, \\
f\flattereq g \quad &:\Longleftrightarrow\quad f^\dagger \preceq g^\dagger.
\end{align*}
For the meaning of $\flattereq$ in
Hardy fields, see Corollary~\ref{hardyflateq}. Note: these relations on~$K^\times$ do not change when passing from $K$ to a compositional conjugate $K^{\phi}$ with ${\phi\in K^\times}$. 
For $f,g\in K^\times$ with $f,g\nasymp 1$ we say that
$f$ and~$g$ are {\bf comparable} if $f\comp g$,  and we say that
$f$ is {\bf flatter} than $g$ if 
$f\flatter g$. 
 Comparability
is an equivalence relation on $\{f\in K^\times: f\nasymp 1\}$. The equivalence
class of an element~$f$ from this set is called its {\bf comparability class},
written as $\Cl(f)$. The relation $\flatter$ induces a linear ordering
on the set of comparability classes by 
$$\Cl(f)<\Cl(g)\ :\Longleftrightarrow\
f\flatter g.$$ Let $(\Gamma,\psi)$ be the asymptotic couple of $K$. Then for $f\in K^\times$ with $f\nasymp 1$ we have:
$\Cl(f)$ is the smallest comparability class of $K$ iff $\psi(vf)$ is the largest element of~$\Psi$. Thus $K$ has a smallest comparability class iff $\Psi$ has a largest element. Note also that if~$K$ is of $H$-type and
$f,g\in K$, then $1\prec f\preceq g \Rightarrow f\flattereq g$. We define the asymptotic field $K$ to be {\bf grounded}
if $K$ has a smallest comparability class, equivalently,
$\Psi$ has a largest element; an asymptotic field that is not grounded is called {\bf ungrounded}. 

\label{p:grounded asymptotic field}
\nomenclature[Rx]{$\Cl(f)$}{comparability class of $f$}
\index{comparability class}
\index{comparable}
\index{flatter}
\index{asymptotic field!grounded}
\index{grounded!asymptotic field}

\subsection*{Hardy fields} Hardy fields, defined below, 
are $H$-asymptotic fields of a classical origin, and
will often serve as examples and counterexamples. They are not just 
asymptotic fields
but also carry a natural ordering.
%Basic references on this subject are \cite[Appendice]{Bourbaki} and \cite{Rosenlicht4}, which also include some proofs omitted here.

Let $\mathcal{G}$ be the ring
of germs at $+\infty$ of real-valued functions whose domain is
a subset of $\R$
containing an interval $(a, +\infty)$, $a\in \R$; the domain may vary 
and the ring operations are defined as usual. We call a germ $g\in \mathcal{G}$ {\em continuous}, respectively {\em differentiable}, if it is the germ of
a continuous, respectively differentiable, function $(a,+\infty)\to \R$ for some $a\in \R$; for differentiable $g\in \mathcal{G}$ we let
$g'\in \mathcal{G}$ denote the germ of the derivative of that function.
If $g\in \mathcal{G}$ is the germ of a real-valued function 
on some interval $(a, +\infty)$, $a\in \R$, then we simplify notation
by letting $g$ also denote this function if the resulting ambiguity is harmless.
With this convention, given a property~$P$ of real numbers
and $g\in \mathcal{G}$ we say that~$P\big(g(t)\big)$ holds eventually if~$P\big(g(t)\big)$ holds for all sufficiently large real $t$. 
We identify each real number $r$ with the germ at~$+\infty$ of the function 
$\R\to \R$ that takes the constant value $r$. This  makes the field $\R$ into a subring of $\mathcal{G}$.

\begin{definition}[N.~Bourbaki]
A {\bf Hardy field\/} is a subring $K$ of $\mathcal{G}$ such that $K$ is a field, all $g\in K$ are differentiable, and $g'\in K$ for all $g\in K$.
\end{definition}

\begin{example}
Every subfield of
$\R$ is a Hardy field. Given polynomials 
$p(x), q(x)\in \R[x]$ with $q\ne 0$, identify the rational function 
$p(x)/q(x)\in \R(x)$
with the germ at $+\infty$ of $t \mapsto p(t)/q(t)$, for $t\in \R$ with $q(t)\ne 0$. This makes the rational function field~$\R(x)$ into a Hardy field, with~$x$ the germ of the identity function on $\R$. 
% and so is any subfield of $\R$.  
\end{example}

\index{germs}
\index{Hardy field}
\index{field!Hardy}

\noindent
\textit{In the rest of this subsection $K$ is a Hardy field, and $f$, $g$ range over~$K$.}\/ We consider~$K$ as a differential field with 
derivation~$f\mapsto f'$. It has constant field 
$C=K\cap \R$, with intersection taken inside $\mathcal{G}$. 
Every nonzero $f$ has a multiplicative inverse in $K$, so eventually $f(t)\neq 0$, hence
either eventually $f(t)<0$ or eventually $f(t)>0$ (by eventual continuity of $f$).
We make~$K$ an ordered field by declaring 
$$f>0\ :\Longleftrightarrow\ \text{$f(t) > 0$, eventually.}$$
Thus $C$ is a common ordered subfield of $K$ and of $\R$. 
Since $f'\in K$, either $f'<0$, or $f'=0$, or $f'>0$,
and accordingly, $f$ is either eventually strictly decreasing, or eventually constant, or eventually strictly increasing, hence the limit $\lim\limits_{t\to\infty} f(t)$ always exists, as an element of the extended real line $\R\cup\{\pm\infty\}$.
Our Hardy field~$K$ is a valued field 
with convex valuation ring
$$ \mathcal{O}\ =\ \big\{f :\ \text{$|f|\le n$ for some $n$}\big\}\ =\ \big\{f :\ \text{$|f| \le c$ for some $c\in C$}\big\}.$$
The natural dominance relation $\preccurlyeq$ of a Hardy field and the derived
asymptotic re\-la\-tions~$\prec$,~$\asymp$,~$\sim$ from Section~\ref{sec:valued fields} have the following meaning in terms of limits:
% can be formulated in terms of limits:
%allow useful reinterpretations in terms of limits
%as $t\to+\infty$:

\begin{lemma}\label{Limits}
For  $g\ne 0$, we have: 
\begin{align*}
f\preccurlyeq g &\Longleftrightarrow
\lim\limits_{t\to+\infty} \frac{f(t)}{g(t)} \in\R,\quad \qquad f\prec g \Longleftrightarrow
\lim\limits_{t\to+\infty} \frac{f(t)}{g(t)} =0,\\
f\asymp g &\Longleftrightarrow
\lim\limits_{t\to+\infty} \frac{f(t)}{g(t)} \in\R^\times,\ \qquad f\sim g \Longleftrightarrow
\lim\limits_{t\to+\infty} \frac{f(t)}{g(t)} = 1.
\end{align*}
%\begin{enumerate}
%\item[\textup{(i)}] $f\preccurlyeq g \Longleftrightarrow
%\lim\limits_{t\to+\infty} \frac{f(t)}{g(t)} \in\R$,
%\item[\textup{(ii)}] $f\prec g \Longleftrightarrow
%\lim\limits_{t\to+\infty} \frac{f(t)}{g(t)} =0$,
%\item[\textup{(iii)}] $f\asymp g \Longleftrightarrow
%\lim\limits_{t\to+\infty} \frac{f(t)}{g(t)} \in\R^\times$, 
%\item[\textup{(iv)}] $f\sim g \Longleftrightarrow
%\lim\limits_{t\to+\infty} \frac{f(t)}{g(t)} = 1$.
%\end{enumerate}
\end{lemma}

\noindent
Thus $f\succ 1$ iff $\lim\limits_{t\to+\infty}
|f(t)|=+\infty$.

\begin{example}
Let $K=\R(x)$. Then for
$a_0,a_1,\dots,a_n\in\R$, $a_n\neq 0$ we have in $K$: 
$a_0+a_1x+\cdots+a_nx^n \sim a_n x^n$. Thus 
$\Gamma=\Z v(x)$ with $v(x)<0=v(1)$, and the valuation is
given by $v(p/q)=(\deg p - \deg q)v(x)$, for $p,q\in\R[x]^{\neq}$.
\end{example}

\noindent
%Here are some properties of Hardy fields pertaining to the interaction
%between the asymptotic relations and the derivation. 
Recall that $f^\dagger:=f'/f=\big(\log |f|\big)'$ for $f\neq 0$.

\begin{prop}\label{LHospital}
Let $f,g\neq 0$.
\begin{enumerate}
\item[\textup{(i)}] If $f\prec g$, then $f^\dagger < g^\dagger$. 
\item[\textup{(ii)}] If $f\prec g\prec 1$, then $f^\dagger\succeq g^\dagger$.
\item[\textup{(iii)}] If $f\asymp 1$ and $\R\subseteq K$, then $f\sim c$ for some $c\in\R^\times$.
\item[\textup{(iv)}] If $f\preccurlyeq 1$, $g\nasymp 1$, then $f'\prec g^\dagger$.
\end{enumerate}
\end{prop}
\begin{proof}
For $f=1$,
item (i) is clear: if $g>\mathcal{O}$, say, then $g$ is
ultimately strictly increasing, hence $g^\dagger=g'/g>0$.
In general, if $f\prec g$, then $1\prec g/f$ and hence $0<(g/f)^\dagger=g^\dagger-f^\dagger$.
As to (ii): if $f\prec g\prec 1$, then $f^\dagger < g^\dagger <1^\dagger=0$ by (i), so $f^\dagger\succeq g^\dagger$.
If $f\asymp 1$, then $c:=\lim\limits_{t\to+\infty}f(t)$ works in (iii).
Next, let $f\preccurlyeq 1$, $g\nasymp 1$. To get $f'\prec g^\dagger$, 
replace $f$ by $f+1$ and $g$ by $1/g$, if necessary, to arrange $f\asymp 1\prec g$. 
Then $f^k\prec g$ for all $k\in \Z$, so $kf^\dagger < g^\dagger$ for all $k\in \Z$,
by (i), hence $f^\dagger \prec g^\dagger$, and thus 
$f'\prec fg^\dagger \asymp g^\dagger$. This proves (iv).  
%  $\lim\limits_{t\to+\infty} f(t)\in\R^\times$ and 
%$\lim\limits_{t\to+\infty}
%g(t)\in\{0,\pm\infty\}$ by Lemma~\ref{Limits}. 
%Hence
%by l'Hospital's rule,
%$$f\ =\ fg/g\ \sim\ (fg)'/g'\ =\ (fg'+f'g)/g'\ =\ f+f'g/g'.$$
%So $f'g/g'\prec f\asymp 1$ and therefore $f'\prec g'/g=g^\dagger$.
\end{proof}

\noindent
By the equivalence of (i) and (iv) in Proposition~\ref{CharacterizationAsymptoticFields} and items (ii) and (iv) in Proposition~\ref{LHospital}, every Hardy field is $H$-asymptotic. Moreover, every Hardy field containing~$\R$ is $\d$-valued, by Proposition~\ref{LHospital}(iii).

\medskip
\noindent  
Next we indicate without proof three ways of extending our Hardy field $K$:\begin{enumerate}
\item let $K^{\operatorname{rc}}$ consist of the continuous germs
$y\in \mathcal{G}$ such that $P(y)=0$ for some $P(Y)\in K[Y]^{\ne}$; then $K^{\operatorname{rc}}$ is the
unique real closed Hardy field extension of $K$
that is algebraic over~$K$, and is thus a real closure of the ordered field $K$;  
\item any differentiable $h\in \mathcal{G}$ with $h'\in K$ yields a 
Hardy field $K(h)\supseteq K$;
\item $\ex^f$ generates a Hardy field extension $K(\ex^f)$ of $K$.
\end{enumerate} 
All three are proved in \cite{Rosenlicht4}, but (1) earlier in \cite{Robinson72}, and (2) and (3) in \cite[Appendice]{Bourbaki}. An easy consequence of (2) is that there is a 
Hardy field $K(\R)\supseteq K$ generated as a field 
over $K$ by $\R$. It has $\R$ as its constant field. A special case
of~(2) is that for $f\in K^{>}$ we have a Hardy field 
$K(\log f)$, since $(\log f)'=f^\dagger\in K$.

\begin{cor}
The derivation of $K$ is small.
\end{cor}
\begin{proof} We have $x'=1$ for the germ $x\in \mathcal{G}$
of the identity function on $\R$, so by extending $K$ we can arrange $x\in K$. Then $f\preccurlyeq 1$ gives $f\prec x$, and hence $f'\prec x'=1$, since $K$ is asymptotic.
\end{proof}

\begin{cor}\label{hardyflateq} Suppose that $f,g\succ 1$. Then
$$f\flattereq g\ \Longleftrightarrow\ \text{$|f|\le |g|^n$  for some $n$.}$$
\end{cor}
\begin{proof} Assume $f\flattereq g$. Then  $f^\dagger\preceq g^\dagger$, that is, $(\log{|f|})'\preceq (\log{|g|})'$.
Working in a Hardy field extension of $K$ containing $\log{|f|}$ and $\log{|g|}$ we have $\log{|f|}\succ 1$ and $\log{|g|}\succ 1$, so $\log{|f|}\preceq \log{|g|}$,
that is, $\log{|f|}\le n\log{|g|}$ for some $n$, and thus 
$|f|\le |g|^n$ for some $n$. The converse follows by reversing
this reasoning.
\end{proof}

\begin{example}[of a Hardy field that is not $\d$-valued]
The Hardy subfield $\Q(x)$ of $\R(x)$ is $\d$-valued with
constant field $\Q$, but its algebraic extension 
$$\Q\left(\sqrt{2+x^{-1}}\,\right)    \qquad\text{(a subfield of the Hardy field $\R(x)^{\operatorname{rc}}$)}$$ 
is a Hardy field with
the same constant field $\Q$ and with $\sqrt{2}$ in the residue field.
Thus $\Q\left(\sqrt{2+x^{-1}}\,\right)$ is not $\d$-valued, even 
though it is a valued differential subfield of the $\d$-valued Hardy field $\R(x)^{\operatorname{rc}}$.
\end{example}

\subsection*{Notes and comments} 
Rosenlicht introduced differential-valued fields and their asymptotic couples in~\cite{Rosenlicht2}; besides Hardy fields that paper has nice examples of a 
complex-analytic nature. To our knowledge the larger class of asymptotic 
fields has not been singled out previously for special attention.  

Basic references on Hardy fields are \cite[Appendice]{Bourbaki} and \cite{Rosenlicht4}.
Corollary~\ref{hardyflateq} is from \cite{Rosenlicht5}.
The next chapter defines $H$-fields and pre-$H$-fields
as ordered valued differential fields that share certain key 
elementary properties with Hardy fields. 

%% file: mt-psipic.tex
\begin{figure}[h!]
\begin{center}
\setlength{\unitlength}{0.240900pt}
\ifx\plotpoint\undefined\newsavebox{\plotpoint}\fi
\sbox{\plotpoint}{\rule[-0.175pt]{0.350pt}{0.350pt}}%
\begin{picture}(1400,800)(150,100)
%\begin{picture}(1500,900)(100,100)
%\tenrm
\put(264,472){\rule[-0.175pt]{282.335pt}{0.350pt}}
\put(850,158){\rule[-0.175pt]{0.350pt}{151.526pt}}
\put(749,787){\makebox(0,0)[l]{\shortstack{$\Gamma\,\uparrow$}}}
\put(1400,518){\makebox(0,0){$\rightarrow\ \Gamma$}}
%\put(944,556){\makebox(0,0)[l]{\scriptsize $(1,1)$}}
%\put(772,438){\makebox(0,0){\scriptsize -1}}
%\put(928,438){\makebox(0,0){\scriptsize 1}}
%\put(873,528){\makebox(0,0){\scriptsize 1}}
\put(850,556){\makebox(0,0){$\circ$}}
%\put(928,542){\makebox(0,0){\tiny $\bullet$}}
%\put(772,472){\makebox(0,0){\tiny $|$}}
%\put(928,472){\makebox(0,0){\tiny $|$}}
%\put(850,541){\makebox(0,0){---}}
\put(1163,822){\makebox(0,0)[l]{$\gamma'$}}
\put(1241,403){\makebox(0,0)[l]{$\gamma^\dagger$}}
\put(264,445){\usebox{\plotpoint}}
\put(264,445){\rule[-0.175pt]{31.317pt}{0.350pt}}
\put(394,446){\rule[-0.175pt]{8.672pt}{0.350pt}}
\put(430,447){\rule[-0.175pt]{5.541pt}{0.350pt}}
\put(453,448){\rule[-0.175pt]{2.891pt}{0.350pt}}
\put(465,449){\rule[-0.175pt]{2.891pt}{0.350pt}}
\put(477,450){\rule[-0.175pt]{2.891pt}{0.350pt}}
\put(489,451){\rule[-0.175pt]{2.891pt}{0.350pt}}
\put(501,452){\rule[-0.175pt]{1.445pt}{0.350pt}}
\put(507,453){\rule[-0.175pt]{1.445pt}{0.350pt}}
\put(513,454){\rule[-0.175pt]{1.325pt}{0.350pt}}
\put(518,455){\rule[-0.175pt]{1.325pt}{0.350pt}}
\put(524,456){\rule[-0.175pt]{1.445pt}{0.350pt}}
\put(530,457){\rule[-0.175pt]{1.445pt}{0.350pt}}
\put(536,458){\rule[-0.175pt]{1.445pt}{0.350pt}}
\put(542,459){\rule[-0.175pt]{1.445pt}{0.350pt}}
\put(548,460){\rule[-0.175pt]{1.445pt}{0.350pt}}
\put(554,461){\rule[-0.175pt]{1.445pt}{0.350pt}}
\put(560,462){\rule[-0.175pt]{0.964pt}{0.350pt}}
\put(564,463){\rule[-0.175pt]{0.964pt}{0.350pt}}
\put(568,464){\rule[-0.175pt]{0.964pt}{0.350pt}}
\put(572,465){\rule[-0.175pt]{0.964pt}{0.350pt}}
\put(576,466){\rule[-0.175pt]{0.964pt}{0.350pt}}
\put(580,467){\rule[-0.175pt]{0.964pt}{0.350pt}}
\put(584,468){\rule[-0.175pt]{0.662pt}{0.350pt}}
\put(586,469){\rule[-0.175pt]{0.662pt}{0.350pt}}
\put(589,470){\rule[-0.175pt]{0.662pt}{0.350pt}}
\put(592,471){\rule[-0.175pt]{0.662pt}{0.350pt}}
\put(595,472){\rule[-0.175pt]{0.964pt}{0.350pt}}
\put(599,473){\rule[-0.175pt]{0.964pt}{0.350pt}}
\put(603,474){\rule[-0.175pt]{0.964pt}{0.350pt}}
\put(607,475){\rule[-0.175pt]{0.723pt}{0.350pt}}
\put(610,476){\rule[-0.175pt]{0.723pt}{0.350pt}}
\put(613,477){\rule[-0.175pt]{0.723pt}{0.350pt}}
\put(616,478){\rule[-0.175pt]{0.723pt}{0.350pt}}
\put(619,479){\rule[-0.175pt]{0.578pt}{0.350pt}}
\put(621,480){\rule[-0.175pt]{0.578pt}{0.350pt}}
\put(623,481){\rule[-0.175pt]{0.578pt}{0.350pt}}
\put(626,482){\rule[-0.175pt]{0.578pt}{0.350pt}}
\put(628,483){\rule[-0.175pt]{0.578pt}{0.350pt}}
\put(631,484){\rule[-0.175pt]{0.723pt}{0.350pt}}
\put(634,485){\rule[-0.175pt]{0.723pt}{0.350pt}}
\put(637,486){\rule[-0.175pt]{0.723pt}{0.350pt}}
\put(640,487){\rule[-0.175pt]{0.723pt}{0.350pt}}
\put(643,488){\rule[-0.175pt]{0.578pt}{0.350pt}}
\put(645,489){\rule[-0.175pt]{0.578pt}{0.350pt}}
\put(647,490){\rule[-0.175pt]{0.578pt}{0.350pt}}
\put(650,491){\rule[-0.175pt]{0.578pt}{0.350pt}}
\put(652,492){\rule[-0.175pt]{0.578pt}{0.350pt}}
\put(655,493){\rule[-0.175pt]{0.578pt}{0.350pt}}
\put(657,494){\rule[-0.175pt]{0.578pt}{0.350pt}}
\put(659,495){\rule[-0.175pt]{0.578pt}{0.350pt}}
\put(662,496){\rule[-0.175pt]{0.578pt}{0.350pt}}
\put(664,497){\rule[-0.175pt]{0.578pt}{0.350pt}}
\put(667,498){\rule[-0.175pt]{0.530pt}{0.350pt}}
\put(669,499){\rule[-0.175pt]{0.530pt}{0.350pt}}
\put(671,500){\rule[-0.175pt]{0.530pt}{0.350pt}}
\put(673,501){\rule[-0.175pt]{0.530pt}{0.350pt}}
\put(675,502){\rule[-0.175pt]{0.530pt}{0.350pt}}
\put(678,503){\rule[-0.175pt]{0.482pt}{0.350pt}}
\put(680,504){\rule[-0.175pt]{0.482pt}{0.350pt}}
\put(682,505){\rule[-0.175pt]{0.482pt}{0.350pt}}
\put(684,506){\rule[-0.175pt]{0.482pt}{0.350pt}}
\put(686,507){\rule[-0.175pt]{0.482pt}{0.350pt}}
\put(688,508){\rule[-0.175pt]{0.482pt}{0.350pt}}
\put(690,509){\rule[-0.175pt]{0.578pt}{0.350pt}}
\put(692,510){\rule[-0.175pt]{0.578pt}{0.350pt}}
\put(694,511){\rule[-0.175pt]{0.578pt}{0.350pt}}
\put(697,512){\rule[-0.175pt]{0.578pt}{0.350pt}}
\put(699,513){\rule[-0.175pt]{0.578pt}{0.350pt}}
\put(702,514){\rule[-0.175pt]{0.578pt}{0.350pt}}
\put(704,515){\rule[-0.175pt]{0.578pt}{0.350pt}}
\put(706,516){\rule[-0.175pt]{0.578pt}{0.350pt}}
\put(709,517){\rule[-0.175pt]{0.578pt}{0.350pt}}
\put(711,518){\rule[-0.175pt]{0.578pt}{0.350pt}}
\put(714,519){\rule[-0.175pt]{0.578pt}{0.350pt}}
\put(716,520){\rule[-0.175pt]{0.578pt}{0.350pt}}
\put(718,521){\rule[-0.175pt]{0.578pt}{0.350pt}}
\put(721,522){\rule[-0.175pt]{0.578pt}{0.350pt}}
\put(723,523){\rule[-0.175pt]{0.578pt}{0.350pt}}
\put(726,524){\rule[-0.175pt]{0.578pt}{0.350pt}}
\put(728,525){\rule[-0.175pt]{0.578pt}{0.350pt}}
\put(730,526){\rule[-0.175pt]{0.578pt}{0.350pt}}
\put(733,527){\rule[-0.175pt]{0.578pt}{0.350pt}}
\put(735,528){\rule[-0.175pt]{0.578pt}{0.350pt}}
\put(738,529){\rule[-0.175pt]{0.530pt}{0.350pt}}
\put(740,530){\rule[-0.175pt]{0.530pt}{0.350pt}}
\put(742,531){\rule[-0.175pt]{0.530pt}{0.350pt}}
\put(744,532){\rule[-0.175pt]{0.530pt}{0.350pt}}
\put(746,533){\rule[-0.175pt]{0.530pt}{0.350pt}}
\put(749,534){\rule[-0.175pt]{0.578pt}{0.350pt}}
\put(751,535){\rule[-0.175pt]{0.578pt}{0.350pt}}
\put(753,536){\rule[-0.175pt]{0.578pt}{0.350pt}}
\put(756,537){\rule[-0.175pt]{0.578pt}{0.350pt}}
\put(758,538){\rule[-0.175pt]{0.578pt}{0.350pt}}
\put(761,539){\rule[-0.175pt]{0.723pt}{0.350pt}}
\put(764,540){\rule[-0.175pt]{0.723pt}{0.350pt}}
\put(767,541){\rule[-0.175pt]{0.723pt}{0.350pt}}
\put(770,542){\rule[-0.175pt]{0.723pt}{0.350pt}}
\put(773,543){\rule[-0.175pt]{0.964pt}{0.350pt}}
\put(777,544){\rule[-0.175pt]{0.964pt}{0.350pt}}
\put(781,545){\rule[-0.175pt]{0.964pt}{0.350pt}}
\put(785,546){\rule[-0.175pt]{0.723pt}{0.350pt}}
\put(788,547){\rule[-0.175pt]{0.723pt}{0.350pt}}
\put(791,548){\rule[-0.175pt]{0.723pt}{0.350pt}}
\put(794,549){\rule[-0.175pt]{0.723pt}{0.350pt}}
\put(797,550){\rule[-0.175pt]{1.445pt}{0.350pt}}
\put(803,551){\rule[-0.175pt]{1.445pt}{0.350pt}}
\put(809,552){\rule[-0.175pt]{1.325pt}{0.350pt}}
\put(814,553){\rule[-0.175pt]{1.325pt}{0.350pt}}
\put(820,554){\rule[-0.175pt]{1.445pt}{0.350pt}}
\put(826,555){\rule[-0.175pt]{1.445pt}{0.350pt}}
%\put(832,556){\rule[-0.175pt]{10.118pt}{0.350pt}}
\put(832,556){\rule[-0.175pt]{2.500pt}{0.350pt}}
\put(856,556){\rule[-0.175pt]{4.300pt}{0.350pt}}
\put(874,555){\rule[-0.175pt]{1.445pt}{0.350pt}}
\put(880,554){\rule[-0.175pt]{1.325pt}{0.350pt}}
\put(885,553){\rule[-0.175pt]{1.325pt}{0.350pt}}
\put(891,552){\rule[-0.175pt]{1.445pt}{0.350pt}}
\put(897,551){\rule[-0.175pt]{1.445pt}{0.350pt}}
\put(903,550){\rule[-0.175pt]{0.723pt}{0.350pt}}
\put(906,549){\rule[-0.175pt]{0.723pt}{0.350pt}}
\put(909,548){\rule[-0.175pt]{0.723pt}{0.350pt}}
\put(912,547){\rule[-0.175pt]{0.723pt}{0.350pt}}
\put(915,546){\rule[-0.175pt]{0.964pt}{0.350pt}}
\put(919,545){\rule[-0.175pt]{0.964pt}{0.350pt}}
\put(923,544){\rule[-0.175pt]{0.964pt}{0.350pt}}
\put(927,543){\rule[-0.175pt]{0.723pt}{0.350pt}}
\put(930,542){\rule[-0.175pt]{0.723pt}{0.350pt}}
\put(933,541){\rule[-0.175pt]{0.723pt}{0.350pt}}
\put(936,540){\rule[-0.175pt]{0.723pt}{0.350pt}}
\put(939,539){\rule[-0.175pt]{0.578pt}{0.350pt}}
\put(941,538){\rule[-0.175pt]{0.578pt}{0.350pt}}
\put(943,537){\rule[-0.175pt]{0.578pt}{0.350pt}}
\put(946,536){\rule[-0.175pt]{0.578pt}{0.350pt}}
\put(948,535){\rule[-0.175pt]{0.578pt}{0.350pt}}
\put(951,534){\rule[-0.175pt]{0.530pt}{0.350pt}}
\put(953,533){\rule[-0.175pt]{0.530pt}{0.350pt}}
\put(955,532){\rule[-0.175pt]{0.530pt}{0.350pt}}
\put(957,531){\rule[-0.175pt]{0.530pt}{0.350pt}}
\put(959,530){\rule[-0.175pt]{0.530pt}{0.350pt}}
\put(962,529){\rule[-0.175pt]{0.578pt}{0.350pt}}
\put(964,528){\rule[-0.175pt]{0.578pt}{0.350pt}}
\put(966,527){\rule[-0.175pt]{0.578pt}{0.350pt}}
\put(969,526){\rule[-0.175pt]{0.578pt}{0.350pt}}
\put(971,525){\rule[-0.175pt]{0.578pt}{0.350pt}}
\put(974,524){\rule[-0.175pt]{0.578pt}{0.350pt}}
\put(976,523){\rule[-0.175pt]{0.578pt}{0.350pt}}
\put(978,522){\rule[-0.175pt]{0.578pt}{0.350pt}}
\put(981,521){\rule[-0.175pt]{0.578pt}{0.350pt}}
\put(983,520){\rule[-0.175pt]{0.578pt}{0.350pt}}
\put(986,519){\rule[-0.175pt]{0.578pt}{0.350pt}}
\put(988,518){\rule[-0.175pt]{0.578pt}{0.350pt}}
\put(990,517){\rule[-0.175pt]{0.578pt}{0.350pt}}
\put(993,516){\rule[-0.175pt]{0.578pt}{0.350pt}}
\put(995,515){\rule[-0.175pt]{0.578pt}{0.350pt}}
\put(998,514){\rule[-0.175pt]{0.578pt}{0.350pt}}
\put(1000,513){\rule[-0.175pt]{0.578pt}{0.350pt}}
\put(1002,512){\rule[-0.175pt]{0.578pt}{0.350pt}}
\put(1005,511){\rule[-0.175pt]{0.578pt}{0.350pt}}
\put(1007,510){\rule[-0.175pt]{0.578pt}{0.350pt}}
\put(1010,509){\rule[-0.175pt]{0.482pt}{0.350pt}}
\put(1012,508){\rule[-0.175pt]{0.482pt}{0.350pt}}
\put(1014,507){\rule[-0.175pt]{0.482pt}{0.350pt}}
\put(1016,506){\rule[-0.175pt]{0.482pt}{0.350pt}}
\put(1018,505){\rule[-0.175pt]{0.482pt}{0.350pt}}
\put(1020,504){\rule[-0.175pt]{0.482pt}{0.350pt}}
\put(1022,503){\rule[-0.175pt]{0.530pt}{0.350pt}}
\put(1024,502){\rule[-0.175pt]{0.530pt}{0.350pt}}
\put(1026,501){\rule[-0.175pt]{0.530pt}{0.350pt}}
\put(1028,500){\rule[-0.175pt]{0.530pt}{0.350pt}}
\put(1030,499){\rule[-0.175pt]{0.530pt}{0.350pt}}
\put(1032,498){\rule[-0.175pt]{0.578pt}{0.350pt}}
\put(1035,497){\rule[-0.175pt]{0.578pt}{0.350pt}}
\put(1037,496){\rule[-0.175pt]{0.578pt}{0.350pt}}
\put(1040,495){\rule[-0.175pt]{0.578pt}{0.350pt}}
\put(1042,494){\rule[-0.175pt]{0.578pt}{0.350pt}}
\put(1045,493){\rule[-0.175pt]{0.578pt}{0.350pt}}
\put(1047,492){\rule[-0.175pt]{0.578pt}{0.350pt}}
\put(1049,491){\rule[-0.175pt]{0.578pt}{0.350pt}}
\put(1052,490){\rule[-0.175pt]{0.578pt}{0.350pt}}
\put(1054,489){\rule[-0.175pt]{0.578pt}{0.350pt}}
\put(1057,488){\rule[-0.175pt]{0.723pt}{0.350pt}}
\put(1060,487){\rule[-0.175pt]{0.723pt}{0.350pt}}
\put(1063,486){\rule[-0.175pt]{0.723pt}{0.350pt}}
\put(1066,485){\rule[-0.175pt]{0.723pt}{0.350pt}}
\put(1069,484){\rule[-0.175pt]{0.578pt}{0.350pt}}
\put(1071,483){\rule[-0.175pt]{0.578pt}{0.350pt}}
\put(1073,482){\rule[-0.175pt]{0.578pt}{0.350pt}}
\put(1076,481){\rule[-0.175pt]{0.578pt}{0.350pt}}
\put(1078,480){\rule[-0.175pt]{0.578pt}{0.350pt}}
\put(1081,479){\rule[-0.175pt]{0.723pt}{0.350pt}}
\put(1084,478){\rule[-0.175pt]{0.723pt}{0.350pt}}
\put(1087,477){\rule[-0.175pt]{0.723pt}{0.350pt}}
\put(1090,476){\rule[-0.175pt]{0.723pt}{0.350pt}}
\put(1093,475){\rule[-0.175pt]{0.964pt}{0.350pt}}
\put(1097,474){\rule[-0.175pt]{0.964pt}{0.350pt}}
\put(1101,473){\rule[-0.175pt]{0.964pt}{0.350pt}}
\put(1105,472){\rule[-0.175pt]{0.662pt}{0.350pt}}
\put(1107,471){\rule[-0.175pt]{0.662pt}{0.350pt}}
\put(1110,470){\rule[-0.175pt]{0.662pt}{0.350pt}}
\put(1113,469){\rule[-0.175pt]{0.662pt}{0.350pt}}
\put(1116,468){\rule[-0.175pt]{0.964pt}{0.350pt}}
\put(1120,467){\rule[-0.175pt]{0.964pt}{0.350pt}}
\put(1124,466){\rule[-0.175pt]{0.964pt}{0.350pt}}
\put(1128,465){\rule[-0.175pt]{0.964pt}{0.350pt}}
\put(1132,464){\rule[-0.175pt]{0.964pt}{0.350pt}}
\put(1136,463){\rule[-0.175pt]{0.964pt}{0.350pt}}
\put(1140,462){\rule[-0.175pt]{1.445pt}{0.350pt}}
\put(1146,461){\rule[-0.175pt]{1.445pt}{0.350pt}}
\put(1152,460){\rule[-0.175pt]{1.445pt}{0.350pt}}
\put(1158,459){\rule[-0.175pt]{1.445pt}{0.350pt}}
\put(1164,458){\rule[-0.175pt]{1.445pt}{0.350pt}}
\put(1170,457){\rule[-0.175pt]{1.445pt}{0.350pt}}
\put(1176,456){\rule[-0.175pt]{1.325pt}{0.350pt}}
\put(1181,455){\rule[-0.175pt]{1.325pt}{0.350pt}}
\put(1187,454){\rule[-0.175pt]{1.445pt}{0.350pt}}
\put(1193,453){\rule[-0.175pt]{1.445pt}{0.350pt}}
\put(1199,452){\rule[-0.175pt]{2.891pt}{0.350pt}}
\put(1211,451){\rule[-0.175pt]{2.891pt}{0.350pt}}
\put(1223,450){\rule[-0.175pt]{2.891pt}{0.350pt}}
\put(1235,449){\rule[-0.175pt]{2.891pt}{0.350pt}}
\put(1247,448){\rule[-0.175pt]{2.650pt}{0.350pt}}
\put(1258,447){\rule[-0.175pt]{5.782pt}{0.350pt}}
\put(1282,446){\rule[-0.175pt]{8.672pt}{0.350pt}}
\put(1318,445){\rule[-0.175pt]{28.426pt}{0.350pt}}
\put(518,158){\usebox{\plotpoint}}
\put(519,159){\usebox{\plotpoint}}
\put(520,160){\usebox{\plotpoint}}
\put(521,161){\usebox{\plotpoint}}
\put(522,162){\usebox{\plotpoint}}
\put(523,163){\usebox{\plotpoint}}
\put(524,164){\usebox{\plotpoint}}
\put(525,165){\usebox{\plotpoint}}
\put(526,166){\usebox{\plotpoint}}
\put(527,167){\usebox{\plotpoint}}
\put(528,168){\usebox{\plotpoint}}
\put(529,169){\usebox{\plotpoint}}
\put(530,170){\usebox{\plotpoint}}
\put(531,171){\usebox{\plotpoint}}
\put(532,172){\usebox{\plotpoint}}
\put(533,173){\usebox{\plotpoint}}
\put(534,174){\usebox{\plotpoint}}
\put(535,175){\usebox{\plotpoint}}
\put(536,176){\usebox{\plotpoint}}
\put(537,178){\usebox{\plotpoint}}
\put(538,179){\usebox{\plotpoint}}
\put(539,180){\usebox{\plotpoint}}
\put(540,181){\usebox{\plotpoint}}
\put(541,182){\usebox{\plotpoint}}
\put(542,183){\usebox{\plotpoint}}
\put(543,184){\usebox{\plotpoint}}
\put(544,185){\usebox{\plotpoint}}
\put(545,186){\usebox{\plotpoint}}
\put(546,187){\usebox{\plotpoint}}
\put(547,188){\usebox{\plotpoint}}
\put(548,189){\usebox{\plotpoint}}
\put(549,191){\usebox{\plotpoint}}
\put(550,192){\usebox{\plotpoint}}
\put(551,193){\usebox{\plotpoint}}
\put(552,194){\usebox{\plotpoint}}
\put(553,195){\usebox{\plotpoint}}
\put(554,196){\usebox{\plotpoint}}
\put(555,197){\usebox{\plotpoint}}
\put(556,198){\usebox{\plotpoint}}
\put(557,199){\usebox{\plotpoint}}
\put(558,200){\usebox{\plotpoint}}
\put(559,201){\usebox{\plotpoint}}
\put(560,202){\usebox{\plotpoint}}
\put(561,204){\usebox{\plotpoint}}
\put(562,205){\usebox{\plotpoint}}
\put(563,206){\usebox{\plotpoint}}
\put(564,207){\usebox{\plotpoint}}
\put(565,208){\usebox{\plotpoint}}
\put(566,209){\usebox{\plotpoint}}
\put(567,210){\usebox{\plotpoint}}
\put(568,211){\usebox{\plotpoint}}
\put(569,212){\usebox{\plotpoint}}
\put(570,213){\usebox{\plotpoint}}
\put(571,214){\usebox{\plotpoint}}
\put(572,215){\usebox{\plotpoint}}
\put(573,217){\usebox{\plotpoint}}
\put(574,218){\usebox{\plotpoint}}
\put(575,219){\usebox{\plotpoint}}
\put(576,220){\usebox{\plotpoint}}
\put(577,221){\usebox{\plotpoint}}
\put(578,223){\usebox{\plotpoint}}
\put(579,224){\usebox{\plotpoint}}
\put(580,225){\usebox{\plotpoint}}
\put(581,226){\usebox{\plotpoint}}
\put(582,227){\usebox{\plotpoint}}
\put(583,228){\usebox{\plotpoint}}
\put(584,230){\usebox{\plotpoint}}
\put(585,231){\usebox{\plotpoint}}
\put(586,232){\usebox{\plotpoint}}
\put(587,233){\usebox{\plotpoint}}
\put(588,235){\usebox{\plotpoint}}
\put(589,236){\usebox{\plotpoint}}
\put(590,237){\usebox{\plotpoint}}
\put(591,238){\usebox{\plotpoint}}
\put(592,240){\usebox{\plotpoint}}
\put(593,241){\usebox{\plotpoint}}
\put(594,242){\usebox{\plotpoint}}
\put(595,243){\usebox{\plotpoint}}
\put(596,245){\usebox{\plotpoint}}
\put(597,246){\usebox{\plotpoint}}
\put(598,247){\usebox{\plotpoint}}
\put(599,248){\usebox{\plotpoint}}
\put(600,249){\usebox{\plotpoint}}
\put(601,251){\usebox{\plotpoint}}
\put(602,252){\usebox{\plotpoint}}
\put(603,253){\usebox{\plotpoint}}
\put(604,254){\usebox{\plotpoint}}
\put(605,255){\usebox{\plotpoint}}
\put(606,256){\usebox{\plotpoint}}
\put(607,258){\usebox{\plotpoint}}
\put(608,259){\usebox{\plotpoint}}
\put(609,260){\usebox{\plotpoint}}
\put(610,261){\usebox{\plotpoint}}
\put(611,263){\usebox{\plotpoint}}
\put(612,264){\usebox{\plotpoint}}
\put(613,265){\usebox{\plotpoint}}
\put(614,266){\usebox{\plotpoint}}
\put(615,268){\usebox{\plotpoint}}
\put(616,269){\usebox{\plotpoint}}
\put(617,270){\usebox{\plotpoint}}
\put(618,271){\usebox{\plotpoint}}
\put(619,273){\usebox{\plotpoint}}
\put(620,274){\usebox{\plotpoint}}
\put(621,275){\usebox{\plotpoint}}
\put(622,276){\usebox{\plotpoint}}
\put(623,278){\usebox{\plotpoint}}
\put(624,279){\usebox{\plotpoint}}
\put(625,280){\usebox{\plotpoint}}
\put(626,281){\usebox{\plotpoint}}
\put(627,283){\usebox{\plotpoint}}
\put(628,284){\usebox{\plotpoint}}
\put(629,285){\usebox{\plotpoint}}
\put(630,286){\usebox{\plotpoint}}
\put(631,288){\usebox{\plotpoint}}
\put(632,289){\usebox{\plotpoint}}
\put(633,290){\usebox{\plotpoint}}
\put(634,291){\usebox{\plotpoint}}
\put(635,293){\usebox{\plotpoint}}
\put(636,294){\usebox{\plotpoint}}
\put(637,295){\usebox{\plotpoint}}
\put(638,296){\usebox{\plotpoint}}
\put(639,298){\usebox{\plotpoint}}
\put(640,299){\usebox{\plotpoint}}
\put(641,300){\usebox{\plotpoint}}
\put(642,301){\usebox{\plotpoint}}
\put(643,303){\usebox{\plotpoint}}
\put(644,304){\usebox{\plotpoint}}
\put(645,305){\usebox{\plotpoint}}
\put(646,306){\usebox{\plotpoint}}
\put(647,308){\usebox{\plotpoint}}
\put(648,309){\usebox{\plotpoint}}
\put(649,310){\usebox{\plotpoint}}
\put(650,311){\usebox{\plotpoint}}
\put(651,313){\usebox{\plotpoint}}
\put(652,314){\usebox{\plotpoint}}
\put(653,315){\usebox{\plotpoint}}
\put(654,316){\usebox{\plotpoint}}
\put(655,318){\usebox{\plotpoint}}
\put(656,319){\usebox{\plotpoint}}
\put(657,320){\usebox{\plotpoint}}
\put(658,322){\usebox{\plotpoint}}
\put(659,323){\usebox{\plotpoint}}
\put(660,324){\usebox{\plotpoint}}
\put(661,326){\usebox{\plotpoint}}
\put(662,327){\usebox{\plotpoint}}
\put(663,328){\usebox{\plotpoint}}
\put(664,330){\usebox{\plotpoint}}
\put(665,331){\usebox{\plotpoint}}
\put(666,332){\usebox{\plotpoint}}
\put(667,334){\rule[-0.175pt]{0.350pt}{0.350pt}}
\put(668,335){\rule[-0.175pt]{0.350pt}{0.350pt}}
\put(669,336){\rule[-0.175pt]{0.350pt}{0.350pt}}
\put(670,338){\rule[-0.175pt]{0.350pt}{0.350pt}}
\put(671,339){\rule[-0.175pt]{0.350pt}{0.350pt}}
\put(672,341){\rule[-0.175pt]{0.350pt}{0.350pt}}
\put(673,342){\rule[-0.175pt]{0.350pt}{0.350pt}}
\put(674,344){\rule[-0.175pt]{0.350pt}{0.350pt}}
\put(675,345){\rule[-0.175pt]{0.350pt}{0.350pt}}
\put(676,347){\rule[-0.175pt]{0.350pt}{0.350pt}}
\put(677,348){\rule[-0.175pt]{0.350pt}{0.350pt}}
\put(678,350){\usebox{\plotpoint}}
\put(679,351){\usebox{\plotpoint}}
\put(680,352){\usebox{\plotpoint}}
\put(681,354){\usebox{\plotpoint}}
\put(682,355){\usebox{\plotpoint}}
\put(683,356){\usebox{\plotpoint}}
\put(684,358){\usebox{\plotpoint}}
\put(685,359){\usebox{\plotpoint}}
\put(686,360){\usebox{\plotpoint}}
\put(687,362){\usebox{\plotpoint}}
\put(688,363){\usebox{\plotpoint}}
\put(689,364){\usebox{\plotpoint}}
\put(690,366){\usebox{\plotpoint}}
\put(691,367){\usebox{\plotpoint}}
\put(692,368){\usebox{\plotpoint}}
\put(693,369){\usebox{\plotpoint}}
\put(694,371){\usebox{\plotpoint}}
\put(695,372){\usebox{\plotpoint}}
\put(696,373){\usebox{\plotpoint}}
\put(697,374){\usebox{\plotpoint}}
\put(698,376){\usebox{\plotpoint}}
\put(699,377){\usebox{\plotpoint}}
\put(700,378){\usebox{\plotpoint}}
\put(701,379){\usebox{\plotpoint}}
\put(702,381){\usebox{\plotpoint}}
\put(703,382){\usebox{\plotpoint}}
\put(704,383){\usebox{\plotpoint}}
\put(705,385){\usebox{\plotpoint}}
\put(706,386){\usebox{\plotpoint}}
\put(707,387){\usebox{\plotpoint}}
\put(708,389){\usebox{\plotpoint}}
\put(709,390){\usebox{\plotpoint}}
\put(710,391){\usebox{\plotpoint}}
\put(711,393){\usebox{\plotpoint}}
\put(712,394){\usebox{\plotpoint}}
\put(713,395){\usebox{\plotpoint}}
\put(714,397){\usebox{\plotpoint}}
\put(715,398){\usebox{\plotpoint}}
\put(716,399){\usebox{\plotpoint}}
\put(717,401){\usebox{\plotpoint}}
\put(718,402){\usebox{\plotpoint}}
\put(719,403){\usebox{\plotpoint}}
\put(720,405){\usebox{\plotpoint}}
\put(721,406){\usebox{\plotpoint}}
\put(722,407){\usebox{\plotpoint}}
\put(723,409){\usebox{\plotpoint}}
\put(724,410){\usebox{\plotpoint}}
\put(725,411){\usebox{\plotpoint}}
\put(726,413){\usebox{\plotpoint}}
\put(727,414){\usebox{\plotpoint}}
\put(728,415){\usebox{\plotpoint}}
\put(729,417){\usebox{\plotpoint}}
\put(730,418){\usebox{\plotpoint}}
\put(731,419){\usebox{\plotpoint}}
\put(732,421){\usebox{\plotpoint}}
\put(733,422){\usebox{\plotpoint}}
\put(734,423){\usebox{\plotpoint}}
\put(735,425){\usebox{\plotpoint}}
\put(736,426){\usebox{\plotpoint}}
\put(737,427){\usebox{\plotpoint}}
\put(738,429){\usebox{\plotpoint}}
\put(739,430){\usebox{\plotpoint}}
\put(740,431){\usebox{\plotpoint}}
\put(741,433){\usebox{\plotpoint}}
\put(742,434){\usebox{\plotpoint}}
\put(743,435){\usebox{\plotpoint}}
\put(744,437){\usebox{\plotpoint}}
\put(745,438){\usebox{\plotpoint}}
\put(746,439){\usebox{\plotpoint}}
\put(747,441){\usebox{\plotpoint}}
\put(748,442){\usebox{\plotpoint}}
\put(749,444){\usebox{\plotpoint}}
\put(750,445){\usebox{\plotpoint}}
\put(751,446){\usebox{\plotpoint}}
\put(752,447){\usebox{\plotpoint}}
\put(753,449){\usebox{\plotpoint}}
\put(754,450){\usebox{\plotpoint}}
\put(755,451){\usebox{\plotpoint}}
\put(756,452){\usebox{\plotpoint}}
\put(757,454){\usebox{\plotpoint}}
\put(758,455){\usebox{\plotpoint}}
\put(759,456){\usebox{\plotpoint}}
\put(760,457){\usebox{\plotpoint}}
\put(761,459){\usebox{\plotpoint}}
\put(762,460){\usebox{\plotpoint}}
\put(763,461){\usebox{\plotpoint}}
\put(764,462){\usebox{\plotpoint}}
\put(765,464){\usebox{\plotpoint}}
\put(766,465){\usebox{\plotpoint}}
\put(767,466){\usebox{\plotpoint}}
\put(768,467){\usebox{\plotpoint}}
\put(769,469){\usebox{\plotpoint}}
\put(770,470){\usebox{\plotpoint}}
\put(771,471){\usebox{\plotpoint}}
\put(772,472){\usebox{\plotpoint}}
\put(773,474){\usebox{\plotpoint}}
\put(774,475){\usebox{\plotpoint}}
\put(775,476){\usebox{\plotpoint}}
\put(776,477){\usebox{\plotpoint}}
\put(777,478){\usebox{\plotpoint}}
\put(778,479){\usebox{\plotpoint}}
\put(779,480){\usebox{\plotpoint}}
\put(780,482){\usebox{\plotpoint}}
\put(781,483){\usebox{\plotpoint}}
\put(782,484){\usebox{\plotpoint}}
\put(783,485){\usebox{\plotpoint}}
\put(784,486){\usebox{\plotpoint}}
\put(785,487){\usebox{\plotpoint}}
\put(786,489){\usebox{\plotpoint}}
\put(787,490){\usebox{\plotpoint}}
\put(788,491){\usebox{\plotpoint}}
\put(789,492){\usebox{\plotpoint}}
\put(790,493){\usebox{\plotpoint}}
\put(791,494){\usebox{\plotpoint}}
\put(792,496){\usebox{\plotpoint}}
\put(793,497){\usebox{\plotpoint}}
\put(794,498){\usebox{\plotpoint}}
\put(795,499){\usebox{\plotpoint}}
\put(796,500){\usebox{\plotpoint}}
\put(797,501){\usebox{\plotpoint}}
\put(798,503){\usebox{\plotpoint}}
\put(799,504){\usebox{\plotpoint}}
\put(800,505){\usebox{\plotpoint}}
\put(801,506){\usebox{\plotpoint}}
\put(802,507){\usebox{\plotpoint}}
\put(803,508){\usebox{\plotpoint}}
\put(804,509){\usebox{\plotpoint}}
\put(805,510){\usebox{\plotpoint}}
\put(806,511){\usebox{\plotpoint}}
\put(807,512){\usebox{\plotpoint}}
\put(808,513){\usebox{\plotpoint}}
\put(809,515){\usebox{\plotpoint}}
\put(810,516){\usebox{\plotpoint}}
\put(811,517){\usebox{\plotpoint}}
\put(812,518){\usebox{\plotpoint}}
\put(813,519){\usebox{\plotpoint}}
\put(814,520){\usebox{\plotpoint}}
\put(815,522){\usebox{\plotpoint}}
\put(816,523){\usebox{\plotpoint}}
\put(817,524){\usebox{\plotpoint}}
\put(818,525){\usebox{\plotpoint}}
\put(819,526){\usebox{\plotpoint}}
\put(820,528){\usebox{\plotpoint}}
\put(821,529){\usebox{\plotpoint}}
\put(822,530){\usebox{\plotpoint}}
\put(823,531){\usebox{\plotpoint}}
\put(824,532){\usebox{\plotpoint}}
\put(825,533){\usebox{\plotpoint}}
\put(826,534){\usebox{\plotpoint}}
\put(827,535){\usebox{\plotpoint}}
\put(828,536){\usebox{\plotpoint}}
\put(829,537){\usebox{\plotpoint}}
\put(830,538){\usebox{\plotpoint}}
\put(831,539){\usebox{\plotpoint}}
\put(832,540){\usebox{\plotpoint}}
\put(833,541){\usebox{\plotpoint}}
\put(834,542){\usebox{\plotpoint}}
\put(835,543){\usebox{\plotpoint}}
\put(836,544){\usebox{\plotpoint}}
\put(837,545){\usebox{\plotpoint}}
\put(838,546){\usebox{\plotpoint}}
\put(839,547){\usebox{\plotpoint}}
\put(840,548){\usebox{\plotpoint}}
\put(841,549){\usebox{\plotpoint}}
\put(842,550){\usebox{\plotpoint}}
\put(843,551){\usebox{\plotpoint}}
\put(845,552){\usebox{\plotpoint}}
%\put(846,553){\usebox{\plotpoint}}
%\put(847,554){\usebox{\plotpoint}}
%\put(848,555){\usebox{\plotpoint}}
%\put(849,556){\usebox{\plotpoint}}
%\put(850,557){\usebox{\plotpoint}}
%\put(851,558){\usebox{\plotpoint}}
%\put(852,559){\usebox{\plotpoint}}
%\put(853,560){\usebox{\plotpoint}}
\put(854,561){\usebox{\plotpoint}}
\put(855,562){\usebox{\plotpoint}}
\put(857,563){\usebox{\plotpoint}}
\put(858,564){\usebox{\plotpoint}}
\put(859,565){\usebox{\plotpoint}}
\put(861,566){\usebox{\plotpoint}}
\put(862,567){\usebox{\plotpoint}}
\put(863,568){\usebox{\plotpoint}}
\put(865,569){\usebox{\plotpoint}}
\put(866,570){\usebox{\plotpoint}}
\put(867,571){\usebox{\plotpoint}}
\put(869,572){\usebox{\plotpoint}}
\put(870,573){\usebox{\plotpoint}}
\put(871,574){\usebox{\plotpoint}}
\put(872,575){\usebox{\plotpoint}}
\put(874,576){\usebox{\plotpoint}}
\put(875,577){\usebox{\plotpoint}}
\put(876,578){\usebox{\plotpoint}}
\put(877,579){\usebox{\plotpoint}}
\put(878,580){\usebox{\plotpoint}}
\put(880,581){\usebox{\plotpoint}}
\put(881,582){\usebox{\plotpoint}}
\put(882,583){\usebox{\plotpoint}}
\put(884,584){\usebox{\plotpoint}}
\put(885,585){\usebox{\plotpoint}}
\put(886,586){\usebox{\plotpoint}}
\put(888,587){\usebox{\plotpoint}}
\put(889,588){\usebox{\plotpoint}}
\put(891,589){\rule[-0.175pt]{0.361pt}{0.350pt}}
\put(892,590){\rule[-0.175pt]{0.361pt}{0.350pt}}
\put(894,591){\rule[-0.175pt]{0.361pt}{0.350pt}}
\put(895,592){\rule[-0.175pt]{0.361pt}{0.350pt}}
\put(897,593){\rule[-0.175pt]{0.361pt}{0.350pt}}
\put(898,594){\rule[-0.175pt]{0.361pt}{0.350pt}}
\put(900,595){\rule[-0.175pt]{0.361pt}{0.350pt}}
\put(901,596){\rule[-0.175pt]{0.361pt}{0.350pt}}
\put(903,597){\rule[-0.175pt]{0.361pt}{0.350pt}}
\put(904,598){\rule[-0.175pt]{0.361pt}{0.350pt}}
\put(906,599){\rule[-0.175pt]{0.361pt}{0.350pt}}
\put(907,600){\rule[-0.175pt]{0.361pt}{0.350pt}}
\put(909,601){\rule[-0.175pt]{0.361pt}{0.350pt}}
\put(910,602){\rule[-0.175pt]{0.361pt}{0.350pt}}
\put(912,603){\rule[-0.175pt]{0.361pt}{0.350pt}}
\put(913,604){\rule[-0.175pt]{0.361pt}{0.350pt}}
\put(915,605){\rule[-0.175pt]{0.413pt}{0.350pt}}
\put(916,606){\rule[-0.175pt]{0.413pt}{0.350pt}}
\put(918,607){\rule[-0.175pt]{0.413pt}{0.350pt}}
\put(920,608){\rule[-0.175pt]{0.413pt}{0.350pt}}
\put(921,609){\rule[-0.175pt]{0.413pt}{0.350pt}}
\put(923,610){\rule[-0.175pt]{0.413pt}{0.350pt}}
\put(925,611){\rule[-0.175pt]{0.413pt}{0.350pt}}
\put(927,612){\rule[-0.175pt]{0.482pt}{0.350pt}}
\put(929,613){\rule[-0.175pt]{0.482pt}{0.350pt}}
\put(931,614){\rule[-0.175pt]{0.482pt}{0.350pt}}
\put(933,615){\rule[-0.175pt]{0.482pt}{0.350pt}}
\put(935,616){\rule[-0.175pt]{0.482pt}{0.350pt}}
\put(937,617){\rule[-0.175pt]{0.482pt}{0.350pt}}
\put(939,618){\rule[-0.175pt]{0.482pt}{0.350pt}}
\put(941,619){\rule[-0.175pt]{0.482pt}{0.350pt}}
\put(943,620){\rule[-0.175pt]{0.482pt}{0.350pt}}
\put(945,621){\rule[-0.175pt]{0.482pt}{0.350pt}}
\put(947,622){\rule[-0.175pt]{0.482pt}{0.350pt}}
\put(949,623){\rule[-0.175pt]{0.482pt}{0.350pt}}
\put(951,624){\rule[-0.175pt]{0.442pt}{0.350pt}}
\put(952,625){\rule[-0.175pt]{0.442pt}{0.350pt}}
\put(954,626){\rule[-0.175pt]{0.442pt}{0.350pt}}
\put(956,627){\rule[-0.175pt]{0.442pt}{0.350pt}}
\put(958,628){\rule[-0.175pt]{0.442pt}{0.350pt}}
\put(960,629){\rule[-0.175pt]{0.442pt}{0.350pt}}
\put(961,630){\rule[-0.175pt]{0.578pt}{0.350pt}}
\put(964,631){\rule[-0.175pt]{0.578pt}{0.350pt}}
\put(966,632){\rule[-0.175pt]{0.578pt}{0.350pt}}
\put(969,633){\rule[-0.175pt]{0.578pt}{0.350pt}}
\put(971,634){\rule[-0.175pt]{0.578pt}{0.350pt}}
\put(974,635){\rule[-0.175pt]{0.482pt}{0.350pt}}
\put(976,636){\rule[-0.175pt]{0.482pt}{0.350pt}}
\put(978,637){\rule[-0.175pt]{0.482pt}{0.350pt}}
\put(980,638){\rule[-0.175pt]{0.482pt}{0.350pt}}
\put(982,639){\rule[-0.175pt]{0.482pt}{0.350pt}}
\put(984,640){\rule[-0.175pt]{0.482pt}{0.350pt}}
\put(986,641){\rule[-0.175pt]{0.578pt}{0.350pt}}
\put(988,642){\rule[-0.175pt]{0.578pt}{0.350pt}}
\put(990,643){\rule[-0.175pt]{0.578pt}{0.350pt}}
\put(993,644){\rule[-0.175pt]{0.578pt}{0.350pt}}
\put(995,645){\rule[-0.175pt]{0.578pt}{0.350pt}}
\put(998,646){\rule[-0.175pt]{0.578pt}{0.350pt}}
\put(1000,647){\rule[-0.175pt]{0.578pt}{0.350pt}}
\put(1002,648){\rule[-0.175pt]{0.578pt}{0.350pt}}
\put(1005,649){\rule[-0.175pt]{0.578pt}{0.350pt}}
\put(1007,650){\rule[-0.175pt]{0.578pt}{0.350pt}}
\put(1010,651){\rule[-0.175pt]{0.482pt}{0.350pt}}
\put(1012,652){\rule[-0.175pt]{0.482pt}{0.350pt}}
\put(1014,653){\rule[-0.175pt]{0.482pt}{0.350pt}}
\put(1016,654){\rule[-0.175pt]{0.482pt}{0.350pt}}
\put(1018,655){\rule[-0.175pt]{0.482pt}{0.350pt}}
\put(1020,656){\rule[-0.175pt]{0.482pt}{0.350pt}}
\put(1022,657){\rule[-0.175pt]{0.530pt}{0.350pt}}
\put(1024,658){\rule[-0.175pt]{0.530pt}{0.350pt}}
\put(1026,659){\rule[-0.175pt]{0.530pt}{0.350pt}}
\put(1028,660){\rule[-0.175pt]{0.530pt}{0.350pt}}
\put(1030,661){\rule[-0.175pt]{0.530pt}{0.350pt}}
\put(1032,662){\rule[-0.175pt]{0.482pt}{0.350pt}}
\put(1035,663){\rule[-0.175pt]{0.482pt}{0.350pt}}
\put(1037,664){\rule[-0.175pt]{0.482pt}{0.350pt}}
\put(1039,665){\rule[-0.175pt]{0.482pt}{0.350pt}}
\put(1041,666){\rule[-0.175pt]{0.482pt}{0.350pt}}
\put(1043,667){\rule[-0.175pt]{0.482pt}{0.350pt}}
\put(1045,668){\rule[-0.175pt]{0.482pt}{0.350pt}}
\put(1047,669){\rule[-0.175pt]{0.482pt}{0.350pt}}
\put(1049,670){\rule[-0.175pt]{0.482pt}{0.350pt}}
\put(1051,671){\rule[-0.175pt]{0.482pt}{0.350pt}}
\put(1053,672){\rule[-0.175pt]{0.482pt}{0.350pt}}
\put(1055,673){\rule[-0.175pt]{0.482pt}{0.350pt}}
\put(1057,674){\rule[-0.175pt]{0.482pt}{0.350pt}}
\put(1059,675){\rule[-0.175pt]{0.482pt}{0.350pt}}
\put(1061,676){\rule[-0.175pt]{0.482pt}{0.350pt}}
\put(1063,677){\rule[-0.175pt]{0.482pt}{0.350pt}}
\put(1065,678){\rule[-0.175pt]{0.482pt}{0.350pt}}
\put(1067,679){\rule[-0.175pt]{0.482pt}{0.350pt}}
\put(1069,680){\rule[-0.175pt]{0.482pt}{0.350pt}}
\put(1071,681){\rule[-0.175pt]{0.482pt}{0.350pt}}
\put(1073,682){\rule[-0.175pt]{0.482pt}{0.350pt}}
\put(1075,683){\rule[-0.175pt]{0.482pt}{0.350pt}}
\put(1077,684){\rule[-0.175pt]{0.482pt}{0.350pt}}
\put(1079,685){\rule[-0.175pt]{0.482pt}{0.350pt}}
\put(1081,686){\rule[-0.175pt]{0.482pt}{0.350pt}}
\put(1083,687){\rule[-0.175pt]{0.482pt}{0.350pt}}
\put(1085,688){\rule[-0.175pt]{0.482pt}{0.350pt}}
\put(1087,689){\rule[-0.175pt]{0.482pt}{0.350pt}}
\put(1089,690){\rule[-0.175pt]{0.482pt}{0.350pt}}
\put(1091,691){\rule[-0.175pt]{0.482pt}{0.350pt}}
\put(1093,692){\rule[-0.175pt]{0.413pt}{0.350pt}}
\put(1094,693){\rule[-0.175pt]{0.413pt}{0.350pt}}
\put(1096,694){\rule[-0.175pt]{0.413pt}{0.350pt}}
\put(1098,695){\rule[-0.175pt]{0.413pt}{0.350pt}}
\put(1099,696){\rule[-0.175pt]{0.413pt}{0.350pt}}
\put(1101,697){\rule[-0.175pt]{0.413pt}{0.350pt}}
\put(1103,698){\rule[-0.175pt]{0.413pt}{0.350pt}}
\put(1104,699){\rule[-0.175pt]{0.379pt}{0.350pt}}
\put(1106,700){\rule[-0.175pt]{0.379pt}{0.350pt}}
\put(1108,701){\rule[-0.175pt]{0.379pt}{0.350pt}}
\put(1109,702){\rule[-0.175pt]{0.379pt}{0.350pt}}
\put(1111,703){\rule[-0.175pt]{0.379pt}{0.350pt}}
\put(1112,704){\rule[-0.175pt]{0.379pt}{0.350pt}}
\put(1114,705){\rule[-0.175pt]{0.379pt}{0.350pt}}
\put(1115,706){\rule[-0.175pt]{0.361pt}{0.350pt}}
\put(1117,707){\rule[-0.175pt]{0.361pt}{0.350pt}}
\put(1119,708){\rule[-0.175pt]{0.361pt}{0.350pt}}
\put(1120,709){\rule[-0.175pt]{0.361pt}{0.350pt}}
\put(1122,710){\rule[-0.175pt]{0.361pt}{0.350pt}}
\put(1123,711){\rule[-0.175pt]{0.361pt}{0.350pt}}
\put(1125,712){\rule[-0.175pt]{0.361pt}{0.350pt}}
\put(1126,713){\rule[-0.175pt]{0.361pt}{0.350pt}}
\put(1128,714){\rule[-0.175pt]{0.361pt}{0.350pt}}
\put(1129,715){\rule[-0.175pt]{0.361pt}{0.350pt}}
\put(1131,716){\rule[-0.175pt]{0.361pt}{0.350pt}}
\put(1132,717){\rule[-0.175pt]{0.361pt}{0.350pt}}
\put(1134,718){\rule[-0.175pt]{0.361pt}{0.350pt}}
\put(1135,719){\rule[-0.175pt]{0.361pt}{0.350pt}}
\put(1137,720){\rule[-0.175pt]{0.361pt}{0.350pt}}
\put(1138,721){\rule[-0.175pt]{0.361pt}{0.350pt}}
\put(1140,722){\rule[-0.175pt]{0.361pt}{0.350pt}}
\put(1141,723){\rule[-0.175pt]{0.361pt}{0.350pt}}
\put(1143,724){\rule[-0.175pt]{0.361pt}{0.350pt}}
\put(1144,725){\rule[-0.175pt]{0.361pt}{0.350pt}}
\put(1146,726){\rule[-0.175pt]{0.361pt}{0.350pt}}
\put(1147,727){\rule[-0.175pt]{0.361pt}{0.350pt}}
\put(1149,728){\rule[-0.175pt]{0.361pt}{0.350pt}}
\put(1150,729){\rule[-0.175pt]{0.361pt}{0.350pt}}
\put(1152,730){\rule[-0.175pt]{0.361pt}{0.350pt}}
\put(1153,731){\rule[-0.175pt]{0.361pt}{0.350pt}}
\put(1155,732){\rule[-0.175pt]{0.361pt}{0.350pt}}
\put(1156,733){\rule[-0.175pt]{0.361pt}{0.350pt}}
\put(1158,734){\rule[-0.175pt]{0.361pt}{0.350pt}}
\put(1159,735){\rule[-0.175pt]{0.361pt}{0.350pt}}
\put(1161,736){\rule[-0.175pt]{0.361pt}{0.350pt}}
\put(1162,737){\rule[-0.175pt]{0.361pt}{0.350pt}}
\put(1164,738){\usebox{\plotpoint}}
\put(1165,739){\usebox{\plotpoint}}
\put(1166,740){\usebox{\plotpoint}}
\put(1168,741){\usebox{\plotpoint}}
\put(1169,742){\usebox{\plotpoint}}
\put(1170,743){\usebox{\plotpoint}}
\put(1172,744){\usebox{\plotpoint}}
\put(1173,745){\usebox{\plotpoint}}
\put(1174,746){\usebox{\plotpoint}}
\put(1176,747){\usebox{\plotpoint}}
\put(1177,748){\usebox{\plotpoint}}
\put(1178,749){\usebox{\plotpoint}}
\put(1179,750){\usebox{\plotpoint}}
\put(1180,751){\usebox{\plotpoint}}
\put(1182,752){\usebox{\plotpoint}}
\put(1183,753){\usebox{\plotpoint}}
\put(1184,754){\usebox{\plotpoint}}
\put(1185,755){\usebox{\plotpoint}}
\put(1186,756){\usebox{\plotpoint}}
\put(1188,757){\usebox{\plotpoint}}
\put(1189,758){\usebox{\plotpoint}}
\put(1191,759){\usebox{\plotpoint}}
\put(1192,760){\usebox{\plotpoint}}
\put(1193,761){\usebox{\plotpoint}}
\put(1195,762){\usebox{\plotpoint}}
\put(1196,763){\usebox{\plotpoint}}
\put(1197,764){\usebox{\plotpoint}}
\put(1199,765){\usebox{\plotpoint}}
\put(1200,766){\usebox{\plotpoint}}
\put(1201,767){\usebox{\plotpoint}}
\put(1203,768){\usebox{\plotpoint}}
\put(1204,769){\usebox{\plotpoint}}
\put(1205,770){\usebox{\plotpoint}}
\put(1207,771){\usebox{\plotpoint}}
\put(1208,772){\usebox{\plotpoint}}
\put(1209,773){\usebox{\plotpoint}}
\put(1211,774){\usebox{\plotpoint}}
\put(1212,775){\usebox{\plotpoint}}
\put(1213,776){\usebox{\plotpoint}}
\put(1215,777){\usebox{\plotpoint}}
\put(1216,778){\usebox{\plotpoint}}
\put(1217,779){\usebox{\plotpoint}}
\put(1219,780){\usebox{\plotpoint}}
\put(1220,781){\usebox{\plotpoint}}
\put(1221,782){\usebox{\plotpoint}}
\put(1223,783){\usebox{\plotpoint}}
\put(1224,784){\usebox{\plotpoint}}
\put(1225,785){\usebox{\plotpoint}}
\put(1226,786){\usebox{\plotpoint}}
\end{picture}

\end{center}
\caption{The graphs of $\gamma\mapsto\gamma'$ and $\gamma\mapsto\gamma^\dagger$.}\label{psipic}
\end{figure}

%% file: mt-9n-2.tex
\section{$H$-Asymptotic Couples}\label{AbstractAsymptoticCouples}

\noindent 
We begin this section with proving some basic facts about 
arbitrary asymptotic couples, but after that we focus on $H$-asymptotic couples.
Using an intermediate value property for suitable functions on 
ordered abelian groups, we derive a key trichotomy 
for $H$-asymptotic couples: Corollary~\ref{trich}. We also 
introduce a contraction map and study $\psi$-maps on suitable 
$H$-asymptotic couples.
We finish with a look at the device of coarsening asymptotic couples and 
asymptotic fields. 

{\em Throughout this section  $(\Gamma,\psi)$ is an asymptotic couple.}\/
 
\index{H-asymptotic@$H$-asymptotic!couple}

\subsection*{Further basic facts on asymptotic couples}

\ifbool{PUP}{In this subsection $\alpha$,~$\beta$,~$\gamma$ range over $\Gamma$. 
By axiom~(AC3) for asymptotic couples (see Section~\ref{sec:ascouples}) we have $\Psi<(\Gamma^>)'$.
Recall from Section~\ref{sec:ordered sets} that 
$$\Psi^{\downarrow}=\{\alpha:\text{$\alpha\leq\gamma$ for some $\gamma\in\Psi$}\}.$$
If $\alpha<0$, then $\alpha'=\alpha^\dagger+\alpha<\alpha^\dagger$.
In particular, $(\Gamma^{<})'\ \subseteq\ \Psi^\downarrow$.  
}{
In this subsection $\alpha$,~$\beta$,~$\gamma$ range over $\Gamma$. 
By axiom~(AC3) for asymptotic couples (see Section~\ref{sec:ascouples}) we have $\Psi<(\Gamma^>)'$.
Recall from Section~\ref{sec:ordered sets} that 
$\Psi^{\downarrow}=\{\alpha:\text{$\alpha\leq\gamma$ for some $\gamma\in\Psi$}\}$. 
If $\alpha<0$, then $\alpha'=\alpha^\dagger+\alpha<\alpha^\dagger$.
In particular, $(\Gamma^{<})'\ \subseteq\ \Psi^\downarrow$.  
}
 
\begin{theorem}\label{26} 
The set $\Gamma\setminus (\Gamma^{\neq})'$ has at most one
element. If $\Psi$ has a largest element $\max\Psi$, then $\Gamma\setminus (\Gamma^{\neq})'=\{\max\Psi\}$.
\end{theorem}

\noindent
For the proof we need the following lemmas:

\begin{lemma}\label{lem:phipsi} Suppose $\beta\ne 0$, $\alpha=\beta'$, $\gamma\neq\alpha$, and 
$\gamma\geq \beta^\dagger$. Then $(\alpha-\gamma)^\dagger\leq\gamma$.
\end{lemma}
\begin{proof}  We have $\alpha-\gamma=\beta+\beta^\dagger-\gamma\leq\beta$ with 
$\alpha-\gamma\ne 0$ and $\beta\neq 0$ and hence 
$(\alpha-\gamma)'\le \beta'$, that is,
$\alpha-\gamma+ (\alpha-\gamma)^\dagger\ \leq\ \alpha$,
 and thus $ (\alpha-\gamma)^\dagger \leq\gamma$, as claimed.
\end{proof}

\begin{lemma}\label{UsefulLemma} 
The following conditions on $\alpha$ are equivalent:
\begin{enumerate}
\item[\textup{(i)}] $\alpha\in (\Gamma^{\neq})'$;
\item[\textup{(ii)}] $ (\alpha-\gamma )^\dagger=\gamma$ for some $\gamma\in\Psi$;
\item[\textup{(iii)}] $ (\alpha-\gamma )^\dagger\leq \gamma$ for some $\gamma\in\Psi$.
\end{enumerate}
\end{lemma}
\begin{proof}
Assume
$\alpha=\beta'$, $\beta\neq 0$, and put $\gamma:=\beta^\dagger\in\Psi$. Then $\gamma\neq \alpha$,
so $(\alpha-\gamma )^\dagger\leq \gamma$ by Lemma~\ref{lem:phipsi}. This shows (i)~$\Rightarrow$~(iii).
For (iii)~$\Rightarrow$~(ii), we first reduce to the case $\alpha=0$ 
by passing from $\psi$ and $\gamma$ to $\psi-\alpha$ and $\gamma-\alpha$.
So assume towards a contradiction that $\gamma\in\Psi$,
$\gamma^\dagger\leq \gamma$, but $\delta^\dagger\neq\delta$ for all 
$\delta\in\Psi$.
Note that then $\gamma\neq 0$ and $\gamma^\dagger< \gamma<(\Gamma^>)'$ and 
$\gamma^{\dagger\dagger} \neq \gamma^\dagger$ (take $\delta:=\gamma^\dagger$).
Hence by Lemma~\ref{BasicProperties}(i),
$$\min\bigl(\gamma^\dagger,\gamma^{\dagger\dagger}\bigr)\ =\ (\gamma-\gamma^\dagger)^\dagger\ >\ \min \big(\gamma,\gamma^\dagger \big)\ =\ \gamma^\dagger,$$
a contradiction. For (ii)~$\Rightarrow$~(i),
assume $\gamma\in\Psi$ and
$(\alpha-\gamma )^\dagger=\gamma$. Then $\alpha-\gamma\neq 0$ and
$\alpha= (\alpha-\gamma )+\gamma= (\alpha-\gamma )+  (\alpha-\gamma )^\dagger=(\alpha-\gamma)'$, so $\alpha\in (\Gamma^{\neq})'$. 
\end{proof}

\begin{proof}[Proof of Theorem~\ref{26}]
Suppose $\alpha\neq \beta$  and $\alpha,\beta\notin (\Gamma^{\neq})'$.
Then Lem\-ma~\ref{UsefulLemma} gives for $\gamma\in\Psi$ that $(\alpha-\gamma)^\dagger>\gamma$ and
$(\beta-\gamma)^\dagger>\gamma$. Thus for
$\gamma:= (\alpha-\beta)^\dagger$, 
%=\bigl(\alpha-\gamma^\dagger\bigr)-\bigl(\beta-\gamma^\dagger\bigr)$, 
$$(\alpha-\beta)^\dagger\ =\ \big( (\alpha-\gamma)-(\beta-\gamma)\big)^\dagger\
\geq\ \min\bigl((\alpha-\gamma)^\dagger,(\beta-\gamma)^\dagger\bigr)\
>\ \gamma=(\alpha-\beta)^\dagger,$$
 a contradiction. Suppose that $\Psi$ has a largest element $\max\Psi$.    
If $\max\Psi=\alpha'$, 
$\alpha\neq 0$, then $\alpha<0$, hence
$\alpha^\dagger=\max\Psi-\alpha>\max\Psi\geq \alpha^\dagger$,
a contradiction. 
\end{proof}

\begin{cor}\label{asy1} There is at most one $\beta$ such that
$$ \Psi\  < \beta\  < (\Gamma^{>})'.$$
If $\Psi$ has a largest element, there is no such $\beta$.
\end{cor}
%\begin{proof} If $\Psi\  \le \alpha < \beta\  < (\Gamma^{>})'$,
%then $\gamma:= \beta-\alpha>0$ gives 
%$$ \gamma^\dagger \le \alpha =\beta-\gamma < \gamma'-\gamma=\gamma^\dagger,$$ 
%a contradiction. 
%\end{proof}

\noindent
An element $\beta$ as in Corollary~\ref{asy1} is called a {\bf gap} in 
$(\Gamma, \psi)$.
If $K$ is an asymptotic field with asymptotic couple $(\Gamma,\psi)$, 
then such an element is also called a gap in $K$. 
Call $(\Gamma, \psi)$ {\bf grounded} if $\Psi$ has a largest element, and {\bf ungrounded} otherwise. So an asymptotic field is grounded iff its asymptotic couple is grounded. 

\index{gap!in an asymptotic couple}
\index{gap!in an asymptotic field}
\index{asymptotic couple!gap}
\index{asymptotic field!gap}
\index{asymptotic couple!grounded}
\index{grounded!asymptotic couple}

\begin{lemma}\label{GapLemma} 
Suppose $K$ is an asymptotic field. Then:
\begin{enumerate}
\item[\textup{(i)}] $K$ has at most one gap;
\item[\textup{(ii)}] if $K$ is grounded, then $K$ has no gap;
\item[\textup{(iii)}] if $L$ is an $H$-asymptotic field extension of $K$ such that $\Gamma^{<}$ is
cofinal in $\Gamma_L^{<}$, then a gap in $K$ remains a gap in $L$;
\item[\textup{(iv)}] if $K$ is the union of a directed family $(K_i)$ of asymptotic subfields
such that no $K_i$ has a gap, then $K$ has no gap. 
\end{enumerate} 
\end{lemma}

\noindent
Theorem~\ref{26} and Lemma~\ref{BasicProperties}(iii)  yield:

\begin{lemma}\label{cofin} If $\Gamma\ne \{0\}$, then $(\Gamma^{>})'$ is cofinal in $\Gamma$, and $(\Gamma^{<})'$ is coinitial
in~$\Gamma$. 
\end{lemma}  
%\begin{proof} Assume $\Gamma\ne \{0\}$. Shifting $\psi$ does not affect the conclusion of the lemma, so we can arrange that $\beta^\dagger \ge 0$ for some $\beta>0$, and thus $(\Gamma^{>})'>0$. Let $\gamma>0$; we have to find $\alpha>0$ such that $\alpha'\ge \gamma$. It is easy to check that $\alpha:= 2\gamma$ works. This proves the desired cofinality.

%Next, let $\gamma<0$; we have to find $\alpha<0$ such that $\alpha'\le \gamma$. Again, $\alpha:= 2\gamma$ works if $\gamma' \le 0$, and this can be arranged by decreasing $\gamma$ so that $-\gamma > \Psi$. 
%\end{proof} 

\begin{lemma}\label{id+qpsi} Suppose $\Gamma$ is divisible and 
$(\Gamma, \psi)$ has asymptotic integration. Then, given $q\in \Q^\times$ and $\alpha$, the map $\gamma\mapsto \gamma + q\psi(\gamma-\alpha)\colon \Gamma^{\ne \alpha} \to \Gamma$
is bijective.
\end{lemma}
\begin{proof} First, reduce to the case $\alpha=0$. Let $\theta$
be the map in the lemma for $\alpha=0$. Now use that 
$\theta(q\gamma)=q\big(\gamma+\psi(q\gamma)\big)=q\big(\gamma+\psi(\gamma)\big)=q\gamma'$ for $\gamma\ne 0$.
\end{proof}

\noindent
Recall from Section~\ref{sec:ascouples} that $\psi$ extends uniquely to a map $(\Q\Gamma)^{\neq}\to\Gamma$, also denoted by~$\psi$, such that $(\Q\Gamma,\psi)$ is an asymptotic couple.
Here $\Q\Gamma$ denotes the ordered divisible hull of $\Gamma$.
Note that $\Psi=\psi(\Gamma^{\neq})=\psi \bigl( (\Q\Gamma)^{\neq} \bigr)$. 
Thus $(\Gamma, \psi)$ is grounded iff  $(\Q\Gamma,\psi)$ is 
grounded. Using Lemma~\ref{UsefulLemma}, we also get:

\begin{cor}\label{QGlemma} 
$(\Gamma^{\neq})'\ =\ ((\Q\Gamma)^{\neq})'\cap\Gamma$.
\end{cor}

\noindent
We say that~$(\Gamma,\psi)$ {\bf has small derivation} if $\alpha'>0$
for all $\alpha>0$. So an asymptotic field has small derivation iff its asymptotic couple has small derivation.

\label{p:grounded as couple/small der}
\index{asymptotic couple!with small derivation}
%\index{asymptotic field!with small derivation}

\begin{lemma}\label{PresInf-Lemma} 
$(\Gamma,\psi)$  has small derivation iff there is no
$\gamma<0$ with $\Psi \le \gamma$.
\end{lemma}
\begin{proof}
Suppose $\Psi\leq\gamma<0$. For  $\alpha<0$ we have $\alpha'<\alpha^\dagger\leq\gamma$.
Thus $\gamma,0\notin (\Gamma^<)'$, and hence $\gamma$ or $0$ is in $(\Gamma^>)'$ by Theorem~\ref{26}.
In each case $(\Gamma,\psi)$  does not have small derivation.
Conversely, suppose $(\Gamma,\psi)$ does not have small derivation. If~$\Gamma^{>}$ has a least element $\varepsilon$, then $\min (\Gamma^>)'=\varepsilon'$ and so $\Psi< \varepsilon'\leq 0$, hence $\Psi\leq -\varepsilon$.
If~$\Gamma^{>}$ does not have a least element, we take any $\alpha>0$ with $\alpha'\leq 0$, and then $\beta$ with $0<\beta<\alpha$ satisfies  
$\Psi \le \gamma:=\beta'<0$.
\end{proof}

\begin{lemma}\label{PresInf-Lemma2} 
Suppose that $(\Gamma,\psi)$ has small
derivation. Then:
\begin{enumerate}
\item[\textup{(i)}] $(\Q\Gamma, \psi)$  has small derivation;
\item[\textup{(ii)}] $\gamma^\dagger \le \gamma\ \Rightarrow\ \gamma^\dagger > -\gamma/n$ for all $n\ge 1$; 
\item[\textup{(iii)}] $\gamma^\dagger > \gamma\ \Rightarrow\ 
\gamma^\dagger > \gamma/n$ for all $n\ge 1$;
\item[\textup{(iv)}] $\gamma^\dagger \leq 0\ \Rightarrow\ 
\gamma^\dagger=o(\gamma)$.
\end{enumerate}
\end{lemma}
\begin{proof}
If $\Gamma^{>}$ has no least element, then $\Gamma^{>}$ is coinitial
in $(\Q\Gamma)^{>}$, therefore $(\Gamma^{>})'$ is
coinitial in $((\Q\Gamma)^{>})'$, which yields (i).
If $\Gamma^{>}$ has a least
element $\varepsilon$, then $\varepsilon'\geq \varepsilon > 0$, and
thus $((\Q\Gamma)^{>})'>\varepsilon^\dagger\geq 0$.
 
Let $n\ge 1$ in the rest of the proof. Suppose that $\gamma^\dagger\leq \gamma$. Then 
$$(-\gamma)'\ =\  -\gamma + \gamma^\dagger\ \le\ 0,$$
so $\gamma > 0$. Applying (i) yields
$(\gamma/n)+
\gamma^\dagger=(\gamma/n)+(\gamma/n)^\dagger>0$,
which yields (ii).  For (iii),
assume $\gamma^\dagger>\gamma$. 
If $\gamma<0$, then
$(-\gamma/n)+ \gamma^\dagger
=(-\gamma/n)+ (-\gamma/n)^\dagger>0$ 
by~(i), so
$\gamma^\dagger > \gamma/n$.
If $\gamma\ge 0$, then $\gamma^\dagger > \gamma\ge \gamma/n$.
This proves (iii). Part (iv) follows from (ii) by taking $\gamma>0$.
\end{proof}

\begin{lemma}\label{hasc} Let $(\Gamma,\psi)$ be an ungrounded $H$-asymptotic couple, and let $\alpha\in \Psi^{\downarrow}$. Then there are $\gamma_0\in \Psi^{>\alpha}$ and $\delta_0\in (\Gamma^>)'$ such that
the map
$$ \gamma \mapsto \psi(\gamma-\alpha)\ \colon\ \Gamma \to \Gamma_{\infty}$$
is constant on the set $[\gamma_0,\delta_0]:=\{\gamma:\gamma_0\leq\gamma\leq\delta_0\}$, with value $>\alpha$.
\end{lemma}
\begin{proof} Take $\beta_0\in \Gamma^{>}$ so small that 
$\psi(\beta_0)> \alpha$ and
$[\beta_0] < \big[\psi(\beta_0)-\alpha\big]$,
and set $\gamma_0:=\psi(\beta_0)$ and $\delta_0:=\beta_0'$.
Then
$$[\gamma_0-\alpha]\ =\  \big[\psi(\beta_0)-\alpha\big]\ 
=\ \big[\beta_0+\psi(\beta_0)-\alpha\big]\ 
=\ [\delta_0-\alpha],$$
so $\psi(\gamma_0-\alpha)=\psi(\delta_0-\alpha)$.
Since the map $\gamma\mapsto\psi(\gamma-\alpha)\colon \Gamma^{>\alpha}\to\Gamma$ is decreasing,
$\gamma_0$ and $\delta_0$ have the desired property. The map takes value
$>\alpha$ by Lemma~\ref{BasicProperties}(i).  
\end{proof} 

\begin{remark}
With $(\Gamma,\psi)$ and $\alpha$, $\gamma_0$, $\delta_0$ as in Lemma~\ref{hasc} and any $H$-asymptotic couple~$(\Gamma_1,\psi_1)$ extending $(\Gamma,\psi)$, the map $\gamma_1 \mapsto \psi_1(\gamma_1-\alpha)$
is constant on $[\gamma_0,\delta_0]_{\Gamma_1}$.
\end{remark}

\noindent
We let  $\psi(*-\alpha)$ denote the constant value of $\psi(\gamma-\alpha)$ for 
$\gamma\in [\gamma_0, \delta_0]$  in Lemma~\ref{hasc}; it does not depend on
the choice of $\gamma_0, \delta_0$ in that lemma.   
\nomenclature[T]{$\psi(*-\alpha)$}{value of $\psi(\gamma-\alpha)$ for all large enough $\gamma\in \Psi$}

\begin{cor}\label{corhasc} Let $(\Gamma,\psi)$ be as in 
Lemma~\ref{hasc}. If $\alpha\in \Psi^{\downarrow}$, then
$\alpha < \psi(*-\alpha)$. If
$\alpha, \beta\in \Psi^{\downarrow}$ and 
$\alpha\le \beta$, then  $\psi(*-\alpha)\le \psi(*-\beta)$.
\end{cor}

\noindent
At one point in Section~\ref{sec:cac} we also need a variant of the above for
$\alpha\in (\Gamma^{>})'$:

\begin{lemma}\label{has+} Let $(\Gamma, \psi)$ be an ungrounded
$H$-asymptotic couple and $\alpha\in (\Gamma^{>})'$. Then there are 
$\gamma_0\in \Psi$ and $\delta_0\in (\Gamma^>)'$ with $\delta_0<\alpha$
such that the map
$$ \gamma \mapsto \psi(\gamma-\alpha)\ \colon\ \Gamma \to \Gamma_{\infty}$$
is constant on the set $[\gamma_0,\delta_0]$.
\end{lemma}
\begin{proof} Take the unique $\beta>0$ with $\beta'=\alpha$ and set 
$\gamma_0:= \psi(\beta)$. Then $\gamma_0-\alpha=-\beta$, so 
$\psi(\gamma_0-\alpha)=\psi(\beta)$. Next, in $(\Q\Gamma, \psi)$ we have
$(\beta/2)' < \beta'=\alpha$ and 
$(\beta/2)'-\alpha=(\beta/2) + \psi(\beta)-\alpha=-\beta/2$, so
$\psi\big((\beta/2)'-\alpha\big)=\psi(\beta)$ as well. So $\gamma_0$ and any 
$\delta_0=\beta_0'$ with $0<\beta_0\le \beta/2$ and $\beta_0\in \Gamma$ have the required property.  
\end{proof}

\noindent
With $(\Gamma,\psi)$ and $\alpha$ as in Lemma~\ref{has+}, we
set $\psi(*-\alpha):=\psi(\beta)$ for the unique 
$\beta>0$ with $\beta'=\alpha$, so  $\psi(*-\alpha)$ is the constant value of
the map in that lemma.

\subsection*{Steady functions and slow functions on $H$-asymptotic couples} 
{\em In this subsection $(\Gamma,\psi)$ is of $H$-type,}\/ and 
$\alpha$,~$\beta$,~$\gamma$ 
range over~$\Gamma$. Our first goal is to show:

\begin{lemma}\label{ivp} The functions 
$$\gamma \mapsto \gamma'\colon \Gamma^{>} \to \Gamma, \qquad \gamma \mapsto \gamma'\colon
\Gamma^{<} \to \Gamma$$
have the intermediate value property. \textup{(}Recall from Lemma~\ref{BasicProperties} that these functions are strictly
increasing.\textup{)}
\end{lemma}

\noindent
The proof is based on the results about steady and slow functions from Section~\ref{sec:oag}. From this section we recall some terminology.
Let $v\colon \Gamma \to S_{\infty}$ be a convex valuation on the ordered abelian group $\Gamma$.
Then  
$v\alpha > v\beta\ \Rightarrow\ \alpha=o(\beta)$. We can replace here 
``$\Rightarrow$'' by ``$\Longleftrightarrow$'' by taking for $v$
the standard valuation
$\gamma \mapsto [\gamma]\colon \Gamma \to [\Gamma]$, which assigns to each $\gamma$ its archimedean class
$[\gamma]$, with the reversed natural ordering on $[\Gamma]$ so as to make $[0]$ its largest element. 
Also important for us is the convex valuation $\psi\colon \Gamma \to \Gamma_{\infty}$. 
We let $o_v(\beta)$ stand for any element $\alpha$ with $v\alpha> v\beta$.
 
Let $U$ be a nonempty convex subset of $\Gamma$. In Section~\ref{sec:oag} we defined a function $i\colon U \to \Gamma$ to be $v$-steady if $i$ has the intermediate value property and
$i(x)-i(y)=x-y + o_{v}(x-y)$ for all distinct $x,y\in U$.  
We also defined a function $\eta\colon U \to \Gamma$ to be  $v$-slow on the right if for all $x,y,z\in U$,
\begin{enumerate}
\item[(s1)]$\eta(x)-\eta(y)=o_v(x-y)$ if $x\ne y$;
\item[(s2)] $\eta(y)=\eta(z)$ if $x<y<z$ and $z-y=o_v(z-x)$.
\end{enumerate} 
In the same way we defined $\eta\colon U \to G$ to be $v$-slow on the left, except that in clause~(s2) we replace
``$x<y<z$'' by ``$x>y>z$.''
By Lemma~\ref{sloste}, the sum of a $v$-steady function $U\to\Gamma$ and a $v$-slow (on the right or on the left) function $U\to\Gamma$ is $v$-steady.

\begin{proof}[Proof of Lemma~\ref{ivp}]
With $v$ the standard valuation, the  identity function on~$\Gamma^{>}$ is clearly $v$-steady. Also, the restriction of $\psi$ to
$\Gamma^{>}$ is $v$-slow on the right: part~(ii) of Lemma~\ref{BasicProperties} shows
that (s1) is satisfied, and the $H$-type assumption implies easily that~(s2) is satisfied. Thus by Lemma~\ref{sloste} the map 
$$\gamma \mapsto \gamma'=\gamma+\psi(\gamma)\ \colon\  \Gamma^{>} \to \Gamma$$ is $v$-steady; in particular, it has the intermediate value property. The other part of Lemma~\ref{ivp} follows in the same way,
with $\Gamma^{<}$ and ``$v$-slow on the left'' in place of~$\Gamma^{>}$ and
``$v$-slow on the right.''
\end{proof}

\begin{lemma}\label{have1} For all $\gamma$ we have:
$\ \gamma \in (\Gamma^{<})'\ \Longleftrightarrow\ \gamma < \psi(\alpha) \text{ for some $\alpha > 0$.}$ \newline 
If $\Psi^{>0}\ne \emptyset$, then there is a unique element
$1\in \Gamma^{>}$ with $\psi(1)=1$.
\end{lemma}
\begin{proof} If $\gamma=(-\alpha)'$, $\alpha>0$, then $\gamma =\psi(\alpha)-\alpha< \psi(\alpha)$. Thus the forward direction of the equivalence holds.
The other direction holds trivially if $\Gamma=(\Gamma^{\neq})'$, so
assume $\Gamma\neq(\Gamma^{\neq})'$. Then Theorem~\ref{26} gives 
$\Gamma\setminus (\Gamma^{\neq})'=\{\beta\}$.
By Corollary~\ref{cofin} and Lemma~\ref{ivp}, $(\Gamma^<)'$ is downward closed in $\Gamma$ and $(\Gamma^>)'$ upward closed in $\Gamma$, so $(\Gamma^<)'=\Gamma^{<\beta}$ and $(\Gamma^>)'=\Gamma^{>\beta}$. Since $(\Gamma^<)'\subseteq\Psi^\downarrow<(\Gamma^>)'$, this yields $\Psi\leq\beta$. So if $\gamma<\psi(\alpha)$ for some $\alpha>0$, then $\gamma<\beta$ and so
$\gamma\in (\Gamma^{<})'$.

Suppose now that $\Psi^{>0}\ne \emptyset$. Then $0\in (\Gamma^<)'$, so there is some
$1\in\Gamma^>$ such that $0=(-1)'$, that is, $\psi(1)=1$; uniqueness of $1$ follows from Lemma~\ref{BasicProperties}(iii).
\end{proof}

\noindent
If $\Psi^{>0}\ne \emptyset$, then the element $1$ as in Lemma~\ref{have1} serves as a unit of length, and we
identify $\Z$ with the subgroup $\Z\cdot 1$ of $\Gamma$ via 
$k\mapsto k\cdot 1$. 

\begin{cor}\label{trich} $(\Gamma, \psi)$ has exactly one of the following three properties: \begin{enumerate}
\item[\textup{(i)}] there is $\beta$ such that $\Psi< \beta < (\Gamma^{>})'$, that is, $(\Gamma, \psi)$ has a gap;
\item[\textup{(ii)}] $\Psi$ has a largest element, that is, $(\Gamma, \psi)$ is grounded;
\item[\textup{(iii)}] $\Gamma=(\Gamma^{\ne})'$, that is, $(\Gamma, \psi)$ has asymptotic integration.
\end{enumerate}
\end{cor}
\begin{proof} We know from Lemma~\ref{asy1} that (i) and (ii) exclude each 
other. Also,~$\beta$ as in~(i) cannot be in $(\Gamma^{<})'$, since $(\Gamma^{<})'\subseteq\Psi^\downarrow$. If $\Psi$ has a largest element, this element is not in $(\Gamma^{\neq})'$, by Theorem~\ref{26}. Thus (i) as well as (ii) excludes (iii). If we are not in case (i) or (ii), then we are in case~(iii) 
by Lemma~\ref{have1}.
\end{proof}

\noindent
The order in which the three possibilities are listed is natural: The trivial 
$H$-asymptotic couple with $\Gamma= \{0\}$ falls under (i), with $\beta=0$.
When $\Gamma$ is not trivial, but divisible and of finite dimension as a vector space
over $\Q$, then $\Psi$ is nonempty and finite, and so we are in case (ii).
Case (iii) typically requires a more infinitary construction: for example, the asymptotic couple of $\mathbb{T}$ falls under (iii).
In all cases the set $\Psi$ has an upper bound in~$\Gamma$, with 
$\sup \Psi = \beta$ in case (i), but the set $\Psi$ has no
supremum in $\Gamma$ in case~(iii).

\medskip\noindent
In view of Corollary~\ref{QGlemma} it follows that a gap in $(\Gamma, \psi)$ 
remains a gap in $(\Q\Gamma,\psi)$, and that if $(\Q\Gamma, \psi)$ 
has asymptotic integration, then so does $(\Gamma,\psi)$. 
 But there are $(\Gamma,\psi)$ with asymptotic integration
such that $(\Q\Gamma, \psi)$ does not have asymptotic integration; see \cite[Example~12.9]{AvdD3}. %\marginpar{but needs work to see this gives an example}
This situation needs to be dealt with in parts of Section~\ref{sec:LO-cuts} where the next lemma will be used. This lemma concerns the special role of the set $2\Psi=\big\{2\psi(\gamma):\gamma\neq 0\big\}$ in that section. 

\begin{lemma}\label{lem:split gamma}
Suppose $(\Gamma,\psi)$ has asymptotic integration, and let $\gamma$ be given. Then the following conditions are equivalent:
\begin{enumerate}
\item[\textup{(i)}] $\gamma$ is a supremum of $2\Psi$ in the ordered set $\Gamma$;
\item[\textup{(ii)}] $\gamma>2\Psi$, and there are no $\alpha,\beta>\Psi$ with $\alpha+\beta=\gamma$;
\item[\textup{(iii)}] $\frac{1}{2}\gamma$ is a gap in $(\Q\Gamma,\psi)$;
\item[\textup{(iv)}] $\frac{1}{2}\gamma$ is a supremum of $\Psi$ in the ordered set $\Q\Gamma$.
\end{enumerate}
\end{lemma}

{\sloppy
\begin{proof} Assume $\gamma=\sup 2\Psi$, and
$\alpha,\beta>\Psi$ satisfy $\alpha+\beta=\gamma$.
Take $\alpha_1\in\Gamma$ with $\Psi < \alpha_1 < \alpha$.
Then $2\Psi < \alpha_1+\beta < \gamma$, contradicting (i).
This gives (i)~$\Rightarrow$~(ii).

To get (ii)~$\Rightarrow$~(iii), we prove the contrapositive. So assume
$\frac{1}{2}\gamma$ is not a gap in~$(\Q\Gamma,\psi)$.
Then either $\Psi\not<\frac{1}{2}\gamma$ or $\frac{1}{2}\gamma\not<  \big( (\Q\Gamma)^>\big){}'$.
If $\Psi\not<\frac{1}{2}\gamma$, then $\gamma\not>2\Psi$. Suppose  
$\frac{1}{2}\gamma\not<  \big( (\Q\Gamma)^>\big){}'$. Since $\Gamma^>$ is coinitial in $(\Q\Gamma)^>$, we have 
$\delta\in\Gamma^>$ with $\frac{1}{2}\gamma \geq \delta'$.
Then $\gamma-\delta'\geq \frac{1}{2}\gamma>\Psi$, so $\alpha:=\delta'$ and $\beta:=\gamma-\alpha$ 
satisfy $\alpha,\beta>\Psi$ and $\alpha+\beta=\gamma$.

Corollary~\ref{asy1} gives  (iii)~$\Rightarrow$~(iv), and
(iv)~$\Rightarrow$~(i) is clear.
 \end{proof}

\noindent
The $H$-asymptotic couple $(\Gamma,\psi)$ is said to have {\bf rational asymptotic integration} if $(\Q\Gamma,\psi)$ has asymptotic integration, and an $H$-asymptotic field is said to have rational asymptotic integration if its asymptotic couple has rational asymptotic integration.
}

\label{p:rational as int}
\index{asymptotic couple!with rational asymptotic integration}
\index{asymptotic field!with rational asymptotic integration}
\index{asymptotic integration!rational}
\index{rational!asymptotic integration}

\subsection*{Contraction} {\em In this subsection $(\Gamma,\psi)$ is 
$H$-asymptotic and ungrounded}, and we let $\alpha$,~$\beta$,~$\gamma$
range over $\Gamma$. So we are in 
case (i) or case (iii) of Corollary~\ref{trich}. By Lemma~\ref{have1} we 
can define the so-called contraction map $\chi\colon \Gamma^{<} \to \Gamma^{<}$ 
by $\chi(\alpha)'=\alpha^\dagger$. It has the following basic properties: \index{contraction}\index{asymptotic couple!contraction}\nomenclature[T]{$\chi$}{contraction map of an $H$-asymptotic couple}

\begin{lemma}\label{contracting} Let $\alpha, \beta<0$. Then: \begin{enumerate}
\item[\textup{(i)}] the map $\chi$ is increasing: $\alpha < \beta\ \Longrightarrow\ \chi(\alpha)\le \chi(\beta)$;
\item[\textup{(ii)}] $\alpha^\dagger=\beta^\dagger\ \Longrightarrow\ \chi(\alpha)=\chi(\beta)$;
\item[\textup{(iii)}] $\chi(\alpha)=o_{\psi}(\alpha)$, and thus $\chi(\alpha)=o(\alpha)$; 
\item[\textup{(iv)}] $\alpha \ne \beta\ \Longrightarrow\ \chi(\alpha)-\chi(\beta)=o(\alpha-\beta)$.
\end{enumerate}
Moreover, $\chi$, while defined in terms of $\psi$, does not 
change upon replacing $\psi$ by a shift $\psi + \gamma$. If $(\Gamma, \psi)$ has small derivation, then 
$$ \alpha^\dagger< 0,\ \psi(\alpha^\dagger)<0\ \Rightarrow\ \alpha^\dagger=o(\alpha),\ \psi(\alpha^\dagger) = o(\alpha^\dagger),\   \chi(\alpha)=\alpha^\dagger-\psi(\alpha^\dagger).$$
\end{lemma}
\begin{proof} Properties (i) and (ii) hold because the map 
$\alpha \mapsto \alpha'$ is strictly increasing and the map $\alpha \mapsto \alpha^\dagger$ is increasing (where $\alpha<0$). To get (iii), we note that
$$\chi(\alpha)<0, \qquad \chi(\alpha)'\ =\ \chi(\alpha)+ \chi(\alpha)^\dagger\ =\ \alpha^\dagger,$$
so $\chi(\alpha)^\dagger > \alpha^\dagger$.  As to (iv), this follows from (ii) if $\alpha^\dagger=\beta^\dagger$. If $\alpha^\dagger < \beta^\dagger$,
then $\beta=o(\alpha)$, so $\chi(\alpha)=o(\alpha-\beta)$ and
$\chi(\beta)=o(\alpha-\beta)$ by (iii), so 
$\chi(\alpha)-\chi(\beta)=o(\alpha-\beta)$. 

The invariance of $\chi$ under shifts follows easily from the definition of $\chi$. 
The last statement is easily deduced from Lemma~\ref{PresInf-Lemma2}(iv).
\end{proof}

%It is also contracting in the sense that
%$\chi(\alpha)=o_{\psi}(\alpha)$, and thus $\chi(\alpha)=o(\alpha)$, for
%$\alpha<0$. 
%A useful feature of $\chi$ is that while it depends on $\psi$, it does not 
%change upon replacing $\psi$ by a shift $\psi + \gamma$. Nevertheless, 
%if $(\Gamma, \psi)$ has small derivation, then for all $\alpha$,
%$$ \alpha< 0,\ \alpha^\dagger< 0,\ \psi(\alpha^\dagger)<0\ \Rightarrow\ 
%\alpha^\dagger=o(\alpha),\ \psi(\alpha^\dagger) = o(\alpha^\dagger),\   %
%\chi(\alpha)=\alpha^\dagger-\psi(\alpha^\dagger),$$
%as is easily deduced from Lemma~\ref{PresInf-Lemma2},~(4). 
%The dependence of $\chi$ on $\psi$ holds also in the sense that 
%if $\psi(\alpha)=\psi(\beta)$
%with $\alpha,\beta<0$, then $\chi(\alpha)=\chi(\beta)$.  

\begin{lemma}\label{Deltacontract}  Let $\Delta$ be a subgroup of $\Gamma$ with $\psi(\Delta^{\ne}) \subseteq \Delta$, and suppose $\alpha\in \Delta^{<}$ is such that $\alpha^\dagger$ is not maximal in $\psi(\Delta^{\ne})$. Then $\chi(\alpha)\in \Delta$.
\end{lemma}
\begin{proof} By the assumptions on $\Delta$ and $\alpha$ we can take $\beta\in \Delta^{<}$ such that $\beta' =\alpha^\dagger$.
Then $\beta=\chi(\alpha)$.
\end{proof}

%\begin{cor}\label{cor:no gaps in T}
%Suppose that for every $\gamma\in\Gamma^<$, the sequence
%$\big(\chi^n(\gamma)\big)$ is cofinal in $\Gamma^<$, and let $\Delta$ be a subgroup of $\Gamma$ with $\Delta\neq\{0\}$ and $\psi(\Delta^{\neq})\subseteq\Delta$. Then either~$\Delta^<$ is cofinal in $\Gamma^<$, or $\psi(\Delta^{\neq})$ has a largest element.
%\end{cor}
%\begin{proof}
%Take some $\gamma\in\Delta^<$.
% If $\chi^n(\gamma)\in \Delta$
%for all $n$, then~$\Delta^<$
%is cofinal in $\Gamma^<$. Suppose that $\chi^n(\gamma)\notin \Delta$ for some $n$, and 
%take $n$ minimal with this property. Then $n>0$ since $\gamma\in\Delta$, and so 
%$\psi\bigl(\chi^{n-1}(\gamma)\bigr)=\max \psi(\Delta^{\neq})$ by Lemma~\ref{Deltacontract}.
%\end{proof}

\noindent
The next result will soon be needed in Section~\ref{Applicationtodifferentialpolynomials}. 

\begin{lemma}\label{evtcon} 
%Assume $\Gamma\ne\{0\}$ is divisible. 
Assume $\Gamma\ne \{0\}$ and $d, e_0, e_1,\dots, e_n\in \Z$ are not all $0$, and let $\alpha\in \Gamma$. Then there exists a $\beta<0$ such that:
\begin{align*} \text{either for all }\gamma\in (\beta,0),\ \alpha + d\psi(\gamma) + e_0\gamma+ e_1\chi(\gamma) + \cdots + e_n\chi^n(\gamma)\ &<\  0,\\
\text{or for all }\gamma\in (\beta,0),\
\alpha + d\psi(\gamma) + e_0\gamma+ e_1\chi(\gamma) + \cdots + e_n\chi^n(\gamma)\ &>\  0.
\end{align*}
\end{lemma}
\begin{proof} Since $\Gamma^{<}$ is cofinal in $(\Q\Gamma)^{<}$ we can pass to the divisible hull and arrange that
$\Gamma$ is divisible. Consider first the case that $(\Gamma, \psi)$ has  asymptotic 
integration.

\subcase[1]{$d \ge 1$.} Then $\alpha + d\psi(\gamma)$ is 
increasing as a function of $\gamma<0$, and 
$\{\alpha+d\psi(\gamma): \gamma<0\}$ has no supremum in $\Gamma$ (using
divisibility of $\Gamma$). Thus we have $\beta\in \Gamma^{<}$ and $\epsilon\in \Gamma^{>}$ such that either $\alpha+d\psi(\gamma) < -\epsilon$ for all $\gamma\in (\beta,0)$, or  $\alpha+d\psi(\gamma) > \epsilon$ for all $\gamma\in (\beta,0)$.
By increasing $\beta$ if necessary, we have also 
$|e_0\gamma+ e_1\chi(\gamma) + \cdots + e_n\chi^n(\gamma)| < \epsilon/2$
for all $\gamma\in (\beta,0)$, and so this $\beta$ does the job.

\subcase[2]{$d \le -1$.} This follows by symmetry from Subcase~1.

\subcase[3]{$d = 0$.} If $\alpha\ne 0$, then the desired result
follows since for all $\gamma<0$,
$$|e_0\gamma+ e_1\chi(\gamma) + \cdots + e_n\chi^n(\gamma)|\ <\ \big(|e_0|+1\big)|\gamma|.$$
If $\alpha=0$, then we take $i$ least with
$e_i\ne 0$ and use that $e_i\chi^i(\gamma)$ is the dominant term
in the sum $e_0\gamma+ e_1\chi(\gamma) + \cdots + e_n\chi^n(\gamma)$ for $\gamma<0$.

\medskip\noindent
This concludes the asymptotic integration case. Before we continue, we 
make a general observation: for $\delta\in \Gamma$ the shift 
$(\Gamma, \psi-\delta)$
has the same $\chi$-map as $(\Gamma,\psi)$, and
$\alpha + d\psi(\gamma)=(\alpha+d\delta) + d(\psi-\delta)(\gamma)$
for $\gamma<0$, so we may replace $(\Gamma, \psi)$ by its shift~${(\Gamma, \psi-\delta)}$, keeping the same $d,e_0,\dots, e_n$ and replacing
$\alpha$ by $\alpha+d\delta$. The remaining case is that $\Psi$ has
a supremum in $\Gamma$. Then we arrange by shifting that this supremum is~$0$, so $\sup \Psi=0\notin \Psi$. If $\alpha\ne 0$, the desired result follows
from
$$|d\psi(\gamma) + e_0\gamma+ e_1\chi(\gamma) + \cdots + e_n\chi^n(\gamma)|\ <\ 
   \big(|e_0|+1\big)|\gamma| \qquad (\gamma< 0),$$
which in turn depends on Lemma~\ref{PresInf-Lemma2}(iv).
Suppose $\alpha=0$. If also $d=0$, then we reach the desired
result by arguing as in Subcase 3 above, so assume $d\ne 0$. If $e_0\ne 0$,
then we use that $e_0\gamma$ is the dominant term in the sum 
$$d\psi(\gamma) + e_0\gamma+ e_1\chi(\gamma) + \cdots + e_n\chi^n(\gamma).$$
So we can assume $\alpha=0$, $d\ne 0$, $e_0=0$. In this case we proceed by
induction on $n$. If $n=0$, the desired result holds, so let $n\ge 1$. 
Using $\psi(\gamma)=\chi(\gamma)+ \psi\big(\chi(\gamma)\big)$ and setting
$\gamma^*:=\chi(\gamma)$ for $\gamma<0$, we have
$$ d\psi(\gamma) + e_1\chi(\gamma) + \cdots + e_n\chi^n(\gamma)\ =\
d\psi(\gamma^*) + (d+e_1)\gamma^*  +\sum_{i=1}^{n-1}e_{i+1}\chi^i(\gamma^*),$$
and so a suitable inductive hypothesis takes care of this.
\end{proof}

\noindent
By Lemma~\ref{have1} the image of the strictly increasing map 
$\gamma\mapsto \gamma+\psi(\gamma)\colon \Gamma^{<}\to\Gamma$ is a
cofinal subset of $\Psi^\downarrow$. In contrast to this, we have:

\begin{lemma}\label{lem:not cofinal in Psi}
Suppose $(\Gamma, \psi)$ has asymptotic integration and $\Gamma$ is divisible. Let $\alpha\in \Gamma$, $d_0\in \N^{\ge 1}$, $d\in\N^{\ge 2}$ and the map 
$i\colon \Gamma^<\to \Gamma$ be such that for all $\gamma<0$,
$$i(\gamma)\ =\ \alpha+d_0\,\gamma+d\,\psi(\gamma)+\epsilon(\gamma)\quad\text{with $\epsilon(\gamma)=o(\gamma)$,}$$
and $i$ is increasing. 
Then there are $\beta<0$ and 
$\gamma_0\in \Psi$ and $\delta_0\in (\Gamma^{>})'$ such that the $i$-image
of the interval $(\beta,0)$ is disjoint from the interval $(\gamma_0,\delta_0)$.
\end{lemma}
\begin{proof} 
If $i(\beta)>\Psi$ for some $\beta<0$, we take such $\beta$ such that $i(\beta)\in (\Gamma^{>})'$, and then the lemma holds with $\delta_0:= i(\beta)$ and any $\gamma_0\in \Psi$.
So we can assume that $i(\Gamma^{<})\subseteq \Psi^{\downarrow}$, and
it remains to show that $i(\Gamma^<)$ is  not cofinal in $\Psi^\downarrow$.
Towards a contradiction, suppose~$i(\Gamma^<)$ is  cofinal in 
$\Psi^\downarrow$. Pick $\alpha_0\in \Gamma$ such that for
all $\gamma<0$,
$$i(\gamma)-\alpha_0\ =\ d_0\gamma + d(\psi(\gamma)-\alpha_0)+ \epsilon(\gamma),$$
in other words, $(1-d)\alpha_0 = \alpha$. 
Replacing $(\Gamma,\psi)$ by its shift $(\Gamma,\psi-\alpha_0)$ and $i$ by the map $i-\alpha_0$, we arrange that $\alpha=0$. We now distinguish two cases.

\case[1]{$(\Gamma,\psi)$ has small derivation.} Then $\Psi\cap \Gamma^{>}\ne \emptyset$, so we have $1\in \Gamma^>$ with $\psi(1)=1$. Identifying $\Q$ with its image under
$r\mapsto r\cdot 1\colon\Q\to\Gamma$, we get
$\frac{3}{2}=\frac{1}{2}+\psi\left(\frac{1}{2}\right)>\Psi$.
Let $-1\leq \gamma<0$. Then $\psi(\gamma) \ge \psi(-1) =1$ and 
$\epsilon(\gamma)> \gamma$, hence
$$i(\gamma)\ =\ d_0\,\gamma+d\,\psi(\gamma)+\epsilon(\gamma)\ >\ (d_0+1)\gamma+d.$$
Taking in addition $\gamma\geq\frac{1}{d_0+1}\left(\frac{3}{2}-d\right)$, we get
$i(\gamma) > \left(\textstyle\frac{3}{2}-d\right)+d=\textstyle\frac{3}{2}$, 
so $i(\gamma)\notin\Psi^\downarrow$, and we have a contradiction.

\case[2] {$(\Gamma, \psi)$ does not have small derivation.} Then by Lemma~\ref{PresInf-Lemma} we can take
$\delta>0$ such that $\Psi \le  -\delta$. For $\gamma<0$ we have 
\begin{align*} i(\gamma)\ &=\ \psi(\gamma) + \big(d_0\gamma + (d-1) \psi(\gamma) + \epsilon(\gamma)\big), \text{ and}\\ 
 d_0\gamma + (d-1) \psi(\gamma) + \epsilon(\gamma)\ &\le\ \psi(\gamma)\ \le \ -\delta,
 \end{align*}
so the $i$-image of $\Gamma^{<}$ is contained in $(\Psi-\delta)^{\downarrow}$, which is contained $\Psi^{\downarrow}$ but not cofinal in it, since $\psi(\delta)> \Psi-\delta$. 
\end{proof} 

%\noindent
%We cannot drop in the lemma above the assumption that $(\Gamma, \psi)$
%has asymptotic integration. 

{\sloppy
\begin{cor}\label{cor:not cofinal in Psi} Suppose $(\Gamma, \psi)$ has rational asymptotic integration. Let~${\alpha\in \Gamma}$ and 
$d_0\in \N^{\ge 1}$, $d\in \N^{\ge 2}$, $e_1,\dots, e_n\in \Z$. 
Let $i\colon \Gamma^<\to\Gamma$ be given by
\begin{equation}\label{eq:def of i}
i(\gamma)\ =\ \alpha+ d_0\,\gamma + d\psi(\gamma) +  \sum_{i=1}^{n} e_i\,\chi^{i}(\gamma).
\end{equation}
Then $i$ is strictly increasing, and there are $\beta<0$ and 
$\gamma_0\in \Psi$ and $\delta_0\in (\Gamma^{>})'$ such that the $i$-image
of the interval $(\beta,0)$ is disjoint from the interval $(\gamma_0,\delta_0)$.
\end{cor}

}
\begin{proof} By passing to the divisible hull of $\Gamma$ and extending $i$ according to \eqref{eq:def of i}, we arrange that $\Gamma$ is divisible.
Clearly $i$ is strictly increasing by Lemma~\ref{contracting}(iv).
Thus Lemma~\ref{lem:not cofinal in Psi} applies.
\end{proof}

%\begin{cor}\label{cor:not cofinal in Psi} Suppose 
%$(\Gamma, \psi)$ has rational asymptotic integration. 
%Let $\alpha\in\Gamma$ and  $d_0,d_1,\dots, d_{r}\in \N$, 
%$d_0\geq 1$, $d_1+\cdots+d_r\geq 2$. 
%Let $i\colon \Gamma^<\to\Gamma$ be given by
%\begin{equation}\label{eq:def of i}
%i(\gamma)\ =\ \alpha+ d_0\,\gamma +  \sum_{i=1}^{r} d_i\,\psi%% 
%   \big(\chi^{i-1}(\gamma)\big)
%\end{equation}
%Then $i$ is strictly increasing, and there are $\beta<0$ and 
%$\gamma_0\in \Psi$ and $\delta_0\in (\Gamma^{>})'$ such that 
%the $i$-image of the interval $(\beta,0)$ is disjoint from the interval 
%$(\gamma_0,\delta_0)$.
%\end{cor}
%\begin{proof} By passing to the divisible hull of $\Gamma$ and 
%extending $i$ according to \eqref{eq:def of i}, we arrange that 
%$\Gamma$ is divisible.
%Clearly $i$ is strictly increasing. Induction on $m$ gives 
%$\psi(\chi^m(\gamma))=\psi(\gamma)-\sum_{i=1}^m\chi^i(\gamma)$ for 
%$\gamma<0$. With $d:=d_1+\cdots+d_r$ and $e_i:=d_{i+1}+\cdots+d_r$ 
%for $i=1,\dots,r-1$, we obtain: 
%$$i(\gamma) \ =\ \alpha+ d_0\,\gamma + d\,\psi(\gamma) - 
%\sum_{i=1}^{r-1} e_i\,\chi^{i}(\gamma) \qquad (\gamma<0).$$
%Thus Lemma~\ref{lem:not cofinal in Psi} applies.
%\end{proof}

\subsection*{Some facts about $\psi$-functions} {\em In this subsection $(\Gamma, \psi)$ is of $H$-type with small derivation, $\Gamma\ne \{0\}$ and $\Gamma$ is divisible}; $\alpha$,~$\beta$,~$\gamma$ range over $\Gamma$. We set
$$\Gamma_{\psi}\ :=\  \big\{\gamma:\ \text{$\psi^n(\gamma)<0$  for all $n\ge 1$}\big\},$$
so $\psi(\Gamma_{\psi})\subseteq \Gamma_{\psi}^{<0}=
-\Gamma_{\psi}^{>0}$, and $\Gamma_{\psi}^{<0}$ is downward closed: if
$\alpha \le \beta\in \Gamma_{\psi}^{<0}$, then $\psi(\alpha) \le \psi(\beta)$, and so by induction $\alpha\in \Gamma_{\psi}^{<0}$. In particular, $\Gamma_{\psi}^{<0}$ is a convex (possibly empty)
subset of $\Gamma^{<}$. Note also that for $\gamma\in \Gamma_{\psi}$ we have
$$\psi^{n+1}(\gamma)\ =\ o\big(\psi^n(\gamma)\big).$$
Let $\alpha$ be given as well as natural numbers $r$
and $d,d_1\dots, d_r$. Then we have a function $i=i_{\alpha, d,d_1,\dots, d_r}\colon \Gamma_{\psi} \to \Gamma$,
$$
i(\gamma)\ :=\ \alpha + d\gamma + d_1\psi(\gamma) + \dots + d_r\psi^r(\gamma).
$$
We call $i$ a {\bf $\psi$-function of slope $d$ and order $r$}, \index{function@$\psi$-function}\index{H-asymptotic@$H$-asymptotic!couple!$\psi$-function}\index{slope!$\psi$-function}\index{order!$\psi$-function}\index{type!$\psi$-function} and also a {\bf $\psi$-function of type~$(\alpha,d,r)$\/} if we want to specify $\alpha$.
By part (i) of Lemma~\ref{BasicProperties}, if
$\alpha$,~$\beta$ are distinct and nonzero, then $\psi\big(\psi(\alpha)-\psi(\beta)\big)> \psi(\alpha -\beta)$, that is, 
$$\psi(\alpha)-\psi(\beta)\ =\ o_{\psi}(\alpha-\beta).$$ It follows that 
the function 
$$\gamma \mapsto  d_1\psi(\gamma) + \dots + d_r\psi^r(\gamma)\ \colon\ \Gamma_{\psi} \to \Gamma$$
is $\psi$-slow on the left when restricted to~$\Gamma_{\psi}^{<0}$ and $\psi$-slow
on the right when restricted to~$\Gamma_{\psi}^{>0}$. Also, if $d\ge 1$, then
for all distinct
$x,y\in \Gamma_{\psi}$,
$$i(x)-i(y)\ =\ d(x-y)+o_{\psi}(x-y)$$ and so $i$ is 
strictly increasing, and $i$ has the intermediate value property on
$\Gamma_{\psi}^{<0}$ as well as on $\Gamma_{\psi}^{>0}$. The next result
is needed in Section~\ref{Applicationtodifferentialpolynomials}.

\begin{lemma}\label{ps2}
Let $i$, $j$ be $\psi$-functions of slopes $d> e$. Then either
$i< j$ on $\Gamma_{\psi}^{<0}$, or there is a unique 
$\beta\in \Gamma_{\psi}^{<0}$ such that $i(\beta)=j(\beta)$; in the latter
case, $i(\gamma)>j(\gamma)$ if $\beta < \gamma \in \Gamma_{\psi}^{<0}$. Likewise, either $i>j$ on 
$\Gamma_{\psi}^{>0}$, or there is a unique 
$\beta\in \Gamma_{\psi}^{>0}$ such that $i(\beta)=j(\beta)$; in the latter
case, $i(\gamma) <j(\gamma)$ if $\beta > \gamma\in  \Gamma_{\psi}^{>0}$. 
\end{lemma}
\begin{proof} Use that $i-j$ is strictly increasing, and that it has the intermediate value property on $\Gamma_{\psi}^{<0}$ as well as on $\Gamma_{\psi}^{>0}$.
\end{proof}

\subsection*{Coarsening} 
Given also an asymptotic couple $(\Gamma_1,\psi_1)$, a {\bf morphism} 
$$h\colon (\Gamma,\psi)\to (\Gamma_1,\psi_1)$$ is a $\le$-preserving morphism
$h\colon\Gamma\to\Gamma_1$ of abelian groups such that
$$ h\bigl(\psi(\gamma)\bigr)= \psi_1\bigl(h(\gamma)\bigr) \text{ for all $\gamma\in
\Gamma\setminus\ker h$.}$$
Let $\Delta$ be a convex subgroup of $\Gamma$. Then we have the ordered quotient
group $\dot{\Gamma}:=\Gamma/\Delta$: if $\gamma \ge 0$ in $\Gamma$,
then $\dot{\gamma}:=\gamma + \Delta\ge 0$ in $\dot{\Gamma}$. 
Lemma~\ref{Coarsening-Prop} below yields an asymptotic couple 
$(\dot{\Gamma}, \dot{\psi})$.
  
\index{morphism!asymptotic couples}
\index{asymptotic couple!morphism}
\nomenclature[T]{$(\dot{\Gamma}, \dot{\psi})$}{coarsening of $(\Gamma,\psi)$ by a convex subgroup}

\begin{lemma}\label{Coarsening-Prop}
There is a unique map 
$\dot{\psi}=\psi_\Delta \colon\ \dot{\Gamma}^{\ne}\to \dot{\Gamma}$
such that $(\dot{\Gamma}, \dot{\psi})$ is an asymptotic couple and
$\gamma\mapsto \dot{\gamma} \colon \Gamma \to \dot{\Gamma}$ is a 
morphism $(\Gamma,\psi) \to
(\dot{\Gamma}, \dot{\psi})$. It is given by
$$\dot{\psi}(\dot{\gamma})\ =\ \psi(\gamma) + \Delta \qquad\text{for 
$\gamma\in
\Gamma\setminus\Delta$.}$$ 
If $\Delta \neq \{0\}$, then 
$\Psi_\Delta:=\dot{\psi}\bigl(\dot{\Gamma}^{\ne})$ has 
supremum $\psi(\delta)+\Delta$, where
$\delta\in\Delta^{\ne}$ is arbitrary.
If $(\Gamma,\psi)$ has small derivation, then 
$(\dot{\Gamma}, \dot{\psi})$ has small derivation.
If $(\Gamma,\psi)$ is $H$-asymptotic, then 
$(\dot{\Gamma}, \dot{\psi})$ is $H$-asymptotic.
\end{lemma}
\begin{proof}
It suffices to prove existence. In what follows,
let $\alpha,\beta\in\Gamma^{\ne}$ and $\alpha \neq \beta$; then by part (ii) of
Lemma~\ref{BasicProperties}, we have
$\bigl[\psi(\alpha)-\psi(\beta)\bigr]<[\alpha-\beta]$.
Hence, if $\alpha-\beta\in\Delta$, then
$\psi(\alpha)-\psi(\beta)\in\Delta$. Therefore we have a map
$$\dot{\psi}\colon \dot{\Gamma}^{\ne}\to \dot{\Gamma} , \quad
\dot{\gamma} \mapsto\psi(\gamma) + \Delta \quad 
(\gamma\in \Gamma \setminus \Delta).$$
If $\beta-\alpha>\Delta$, then
$\beta'-\alpha'=\bigl(\beta+\psi(\beta)\bigr)-\bigl(\alpha+\psi(\alpha)\bigr)>\Delta$, and
thus by Lemma~\ref{CharAsymptoticCouples},  $(\dot{\Gamma},\dot{\psi})$
is an asymptotic couple. 
To prove the claim about $\Psi_{\Delta}$, 
assume in addition that $\alpha\in \Delta$, $\alpha>0$, and $\beta>\Delta$. 
Then $\bigl(\beta+\psi(\beta)\bigr)-\psi(\alpha)>\Delta$, so in $\dot{\Gamma}$,
$$\psi(\alpha) + \Delta\ <\  \beta' + \Delta\ =\ (\beta+\Delta)'.$$
Moreover,  $\psi(\beta) < \alpha + \psi(\alpha)\in \psi(\alpha) + \Delta$, so  
$\Psi_{\Delta} \leq \psi(\alpha)+\Delta$. 

If $(\Gamma,\psi)$ has small derivation, then $\alpha>\Delta$ implies
$\alpha'=\alpha+ \alpha^\dagger > \Delta$ by Lemma~\ref{PresInf-Lemma2}(iv), 
and so $(\dot{\Gamma}, \dot{\psi})$ has small derivation.  
\end{proof}

\noindent
In connection with coarsening we sometimes use the following lemma.

\begin{lemma}\label{asymp-lemma1}
Assume $\Delta \ne \{0\}$. The following conditions are equivalent:
\begin{enumerate}
\item[\textup{(i)}] $\psi(\Delta^{ \ne})\cap\Delta\neq\emptyset$;
\item[\textup{(ii)}] $\psi(\Delta^{\ne})\subseteq\Delta$;
\item[\textup{(iii)}] $(\Delta^{\ne})'\cap\Delta \neq\emptyset$;
\item[\textup{(iv)}] $(\Delta^{\ne})'\subseteq\Delta$.
\end{enumerate}
If $(\Gamma,\psi)$ is of $H$-type and \textup{(i)} holds, then $\Psi^{> 0}\subseteq \Delta$.
\end{lemma}
\begin{proof}
Let $\delta\in\Delta^{\ne}$ be such that $\delta^\dagger\in\Delta$. Then we have for
$\delta_1\in\Delta^{\ne}$: $$\abs{\delta^\dagger-\delta_1^\dagger}\leq\abs{\delta-
\delta_1}\in\Delta,$$ so $\delta_1^\dagger\in\Delta$. The equivalences easily
follow. Now assume $(\Gamma,\psi)$ is of $H$-type, (i) holds, and
$\Psi^{>0}\neq \emptyset$.
Then $(\Delta, \psi|\Delta^{\ne})$ is an $H$-asymptotic couple with
$\psi(\delta)>0$ for some $\delta\in \Delta^{\ne}$, so by Lemma~\ref{have1} the unique element $1\in \Gamma^{>}$ with $\psi(1)=1$ lies in $\Delta$.
Then $\Psi^{>0}< 1+1\in \Delta$, and thus $\Psi^{>0}\subseteq \Delta$. 
\end{proof}

\noindent
Now let $K$ be an asymptotic field with valuation $v$ and
asymptotic couple $(\Gamma, \psi)$. The convex subgroup $\Delta$ of
$\Gamma$ then leads to the coarsening of $v$ by $\Delta$:
$$\dot{v}=v_\Delta \colon\
K^\times\to\dot{\Gamma}=\dot{v}(K^\times),\qquad \dot{v}(f):=v(f) + \Delta.$$
The dominance relation on $K$ corresponding to the coarsened
valuation $\dot{v}$ is denoted by $\dotpreceq$ or $\preceq_\Delta$, so 
$f\preceq g \Rightarrow f \dotpreceq g$ for $f,g\in K$. (See Section~\ref{sec:decomposition}.)  Let $(K, {\dotpreceq})$ be the valued differential
field $K$ with $v$ replaced by $\dot{v}$.  

\begin{cor}\label{coarsen}
The valued differential field $(K,{\dotpreceq})$ has the following properties: 
\begin{enumerate}
\item[\textup{(i)}] $(K,{\dotpreceq})$ is
an asymptotic field with asymptotic couple $(\dot{\Gamma},\dot{\psi})$;
\item[\textup{(ii)}] if $a,b\in K$, $a, b\not\asymp 1$, then $a \dotprec b\ \Longleftrightarrow\ a' \dotprec b'$.
\end{enumerate}
\end{cor}
\begin{proof} For (i), use Lemma~\ref{Coarsening-Prop} and 
Proposition~\ref{CharacterizationAsymptoticFields}. For (ii),
let $a,b\in K$ and $a, b\nasymp 1$. Then $a\prec b \Longleftrightarrow 
a' \prec b'$. Set $\alpha=va$, $\beta=vb$, so
$va'=\alpha + \alpha^\dagger$ and $vb'=\beta+ \beta^\dagger$. Now
Lemma~\ref{BasicProperties}(ii) gives
$a \dotprec b\ \Longleftrightarrow\ a' \dotprec b'$.  
\end{proof}

\subsection*{Notes and comments} Lemmas~\ref{evtcon},~\ref{lem:not cofinal in Psi}, Corollary~\ref{cor:not cofinal in Psi}, and the results on 
$\psi$-functions are new, but the rest of this section is essentially taken from~\cite{AvdD, AvdD2, AvdD3, Rosenlicht1}. (For example, Theorem~\ref{26} is 
from \cite{Rosenlicht1}, but its first part only for well-ordered $\Psi$. 
In the above generality, it is~\cite[Theorem 2.6]{AvdD2}.) 

A more explicit version of Lemma~\ref{hasc} is in~\cite{Geh}. 
The $H$-type assumption in
Lemma~\ref{ivp} cannot be dropped; see \cite[Examples~2.8,~2.9]{AvdD2}. 
Ordered abelian groups with a contraction map are studied in~\cite{FVK-chi1, FVK-chi2}
and \cite[Appendix~A]{Kuhlmann00}; see also~\cite[Section~5]{psiremarks}. 
Constructions of $\d$-valued fields with prescribed asymptotic couple can be 
found in \cite{KuMa2, Rosenlicht1, Rosenlicht3} and \cite[Section~11]{AvdD3}.

%% file: mt-9n-3.tex
\section{Application to Differential Polynomials}\label{Applicationtodifferentialpolynomials}

\noindent 
{\em In this section $K$ is an $H$-asymptotic field with divisible value group
$\Gamma\ne \{0\}$.}\/ We let~$\beta$ and $\gamma$ range over 
$\Gamma$. We also fix a differential polynomial $P\in K\{Y\}^{\ne}$ of 
order~$r$. We shall use here the $\psi$-functions from 
Section~\ref{AbstractAsymptoticCouples}
to describe for some $\beta<0$ the behavior of $vP(y)$ for $vy$ in the interval
 $(\beta,0)\subseteq \Gamma$.

\subsection*{Special case} Here is a key step in our work:

\begin{lemma}\label{ps6a} Assume $\der\smallo\subseteq \smallo$. Let $B$ be a nonempty convex subset of $\Gamma^{<}$ 
such that $\psi(B)\subseteq B$ \textup{(}so $B\subseteq \Gamma_{\psi}^{<0}$\textup{)}.  Then there are 
$\beta\in B$ and a $\psi$-function $i=i_{\alpha,d,d_1,\dots, d_r}$ 
such that $P_d\ne 0$, and for all $y\in K$ with $vy \in B^{>\beta}$:  \begin{enumerate}
\item[\textup{(i)}] $P(y)\sim P_d(y)$; and $P_d(y)\succ P_e(y)$ whenever $e\in \N$, $e\ne d$;
\item[\textup{(ii)}] $vP(y)\ =\ i(vy)$;
\item[\textup{(iii)}] $d_j\le d\cdot r!$ for $j=1,\dots,r$. 
\end{enumerate}
\end{lemma}

\noindent
The inductive proof
goes as follows: if $P$ is inhomogeneous, consider the homogeneous parts of 
$P$; reduce the homogeneous case of positive order to that of lower 
order by taking the Riccati transform.

\begin{proof}[Proof of Lemma~\ref{ps6a}] We proceed by induction on $r$.
For $r=0$ and homogeneous $P=aY^d$, $a\in K^\times$, the desired result 
holds with any $\beta\in B$ and $i=i_{\alpha,d}$, $\alpha=va$. Assume it holds 
whenever $P$ (of order $r$) is homogeneous. Suppose now that~$P$ of order~$r$ is not
necessarily homogeneous, and let $P_{d_0},\dots, P_{d_k}$ with $d_0 < \dots < d_k$
be the nonzero homogeneous parts of $P$. 
Take $\beta\in B$
and $\psi$-functions $i_0,\dots, i_k$ with
$i_p=i_{\alpha_p,d_{p},d_{p,1},\dots, d_{p,r}}$ and
$d_{p,1},\dots, d_{p,r}\le d_pr!$ for $p=0,\dots,k$, and
$$vP_{d_0}(y)=i_0(vy), \dots,  vP_{d_k}(y)=i_k(vy)\ \text{ whenever $y\in K$, $vy \in B^{>\beta}$,} $$ 
and such that for all distinct $p,q\in \{0,\dots,k\}$, 
either $i_p> i_q$ on $B^{>\beta}$, or $i_p < i_q $ on~$B^{>\beta}$.
This gives $p\in  \{0,\dots,k\}$ such that for all $y\in K$  
with $vy \in B^{>\beta}$,
$$P(y)\sim P_{d_p}(y), \qquad P_{d_p}(y)\succ P_e(y) \text{ whenever
 $e\in \N$, $e\ne d_p$.}$$
Suppose $H\in K\{Y\}$ has order $r+1$, and is homogeneous of
degree $d$. Put $G:= \Ric(H)$, so $G$ is of order $r$, $v(G)=v(H)$,
and $\deg G\le \wt(H) \le (r+1)d$.
Take $\beta$ and $i_p$ as above for $G$ in the role of $P$, and note that 
$vH(y)=dvy + i_p\big(\psi(vy)\big)$ for all $y\in K$ with $vy \in B^{>\beta}$.
Thus the lemma holds for $H$ instead of $P$. 
\end{proof}

\noindent
The same induction on $r$ shows that in Lemma~\ref{ps6a}, if $P\in C\{Y\}^{\ne}$
and $P$ is homogeneous of degree $d$, then we can take $i$ to be of type 
$(0,d,r)$. For later use it is important to find out more about the values of
$\alpha$ and $d$ in Lemma~\ref{ps6a} for general $P$. 
The next lemma provides such
information for special $K$ and $B$. 

\begin{lemma}\label{vpy} Suppose $\sup \Psi = 0$ and $0\notin \Psi$. Set 
$B:=\Gamma^{<}$; so $\psi(B)\subseteq B$. Then there are $\beta<0$ and a $\psi$-function $i=i_{\alpha,d,d_1,\dots, d_r}$ with the properties
described in Lemma~\ref{ps6a} such that in addition $\alpha=v(P)$ and 
$d=\max\!\big\{e\in \N:\ v(P)=v(P_e)\big\}$. 
\end{lemma}
\begin{proof} Follow the proof of Lemma~\ref{ps6a}, taking into account 
that $\Gamma_{\psi}=\Gamma^{\ne}$, and that for any $\psi$-functions $i_1=i_{\alpha_1,d,d_1,\dots, d_m}$ and $i_2=i_{\alpha_2,e,e_1,\dots, e_n}$, if
either $\alpha_1 < \alpha_2$, or
$\alpha_1=\alpha_2$ and $d>e$, then $i_1 < i_2$ on some interval $(\beta,0)$ with $\beta < 0$.
\end{proof} 

\begin{cor}\label{vpyc}  Suppose $\sup \Psi = 0$ and $0\notin \Psi$.
Assume also that $P\succ P_0$. 
Then there is $\beta<0$
such that $P(y_1) \succ P(y_2) \succ P(z)$ whenever $y_1, y_2,z\in K$ and 
$\beta < vy_1 < vy_2 < vz=0$. 
\end{cor}
\begin{proof} Take $\beta$ and $i=i_{\alpha,d,d_1,\dots, d_r}$ as in 
Lemma~\ref{vpy} above. Since $P\ne 0$ and $P\succ P_0$ we have $d\ge 1$, 
and so the desired result holds for this $\beta$. 
\end{proof}

\subsection*{A useful coarsening} In order to better pin down 
the value of $\alpha$ in Lemma~\ref{ps6a} we 
introduce a coarsening that gives a reduction to
the case $\sup \Psi=0\notin \Psi$. This will be used in the proof of 
Proposition~\ref{alpha=vev}. {\em In this subsection we assume 
$\der\smallo\subseteq \smallo$}.

\medskip\noindent
Let $B$ be a nonempty convex subset of $\Gamma^{<}$ such that 
$\psi(B)\subseteq B$, just as in Lemma~\ref{ps6a} \textup{(}so $B\subseteq \Gamma_{\psi}^{<}$\textup{)}. 
Then we have the convex subgroup 
$$\Delta=\Delta(B)\ :=\ \big\{\gamma:\  \psi(\gamma)> \psi(B)\big\}$$
of $\Gamma$. Note that also $\Delta= \big\{\gamma:\  
\psi(\gamma)> B\big\}$. We are now going to
use the coarsening $\dot{v}=v_{\Delta}$ of $v$ by $\Delta$. Recall from the
end of Section~\ref{AbstractAsymptoticCouples} that
$K$ with $\dot{v}$ is an $H$-asymptotic field with asymptotic couple 
$(\dot{\Gamma}, \dot{\psi})$ and $\Psi_{\Delta}:=\dot{\psi}(\dot{\Gamma}^{\ne})$. 

\begin{lemma}\label{alphadelta}  Let $\alpha$ be as in Lemma~\ref{ps6a}. Then 
$\dot{\alpha}=\dot{v}(P)$. Moreover, if $P\dotsucc P_0$, then there is $\beta\in B$ such that 
$P(y_1) \dotsucc P(y_2) \dotsucc P(z)$
whenever $y_1, y_2 ,z\in K$ and 
$vy_1, vy_2\in B^{>\beta}$, $v(z)\in \Delta$, $y_1\dotsucc y_2$.
\end{lemma}
\begin{proof} We first show that $\Delta=\big\{\gamma: B<-|\gamma|\big\}$.  
If $\beta\in B$ and $\gamma\in \Delta^{<}$, 
then $\psi(\beta) < \psi(\gamma)$, so $\beta < \gamma$. This gives
$B<\Delta$. Next, assume $\gamma\in \Gamma^{<}$ and $B<\gamma$. Then
$\beta < \psi(\beta)  \le \psi(\gamma)$ for all $\beta\in B$, so
$\gamma\in \Delta$. This proves our claim about $\Delta$.

Note that $\dot{\Gamma}=\Gamma/\Delta$ is
divisible, and $\dot{B}$ is a convex subset of $\dot{\Gamma}^{<}$
with $[\dot{\beta}, \dot{0})\subseteq \dot{B}$ for all $\beta\in B$.
So for $\beta\in B$ we have $\beta < \psi(\beta)\in B$, so 
$\dot{\beta} \le \dot{\psi}(\dot{\beta})=\psi(\beta)+\Delta < \dot{0}$.
Thus $\dot{0}\notin\Psi_{\Delta}$ and $\sup \Psi_{\Delta}=\dot{0}$.
We can now apply 
Lemma~\ref{vpy} and Corollary~\ref{vpyc} to $K$ equipped with $\dot{v}$,
and this gives the desired result.
\end{proof}

\subsection*{Reduction to the special case}
{\em In this subsection $K$ is ungrounded.}\/
 So we have the contraction map $\chi\colon \Gamma^{<} \to \Gamma^{<}$. We now consider a convex nonempty set~$B\subseteq \Gamma^{<}$ 
that satisfies a {\em weaker\/} condition than in Lemma~\ref{ps6a}, 
namely, $\chi(B)\subseteq B$. For example, this holds for $B=\Gamma^{<}$. 
(In this subsection we do not assume $\der\smallo\subseteq \smallo$.)

\begin{prop}\label{ps6b} There are
$\alpha\in \Gamma$ and $\beta\in B$, and $d_0, d_1,\dots, d_r\in \N$ such that
$d_0\le \deg P$, and for all $y\in K$,
$$vy\in B^{>\beta}\ \Longrightarrow\ vP(y)\ =\ \alpha+ d_0vy + \sum_{i=1}^r d_i\psi\big(\chi^{i-1}(vy)\big).$$
If $P(0)=0$ we can take $d_0\ge 1$. If $P$ is homogeneous, we can take 
$d_0=\deg P$.
\end{prop}

\begin{proof} Consider first the special case that  $\der\smallo\subseteq \smallo$ and  $\psi(B)<0$. Then we have ${\psi(B)\subseteq B}$: given $\gamma\in B$ we have $\psi(\gamma)<0$ and $\chi(\gamma)\in B$ and so $\psi\big(\chi(\gamma)\big)<0$, hence
$$\gamma < \psi(\gamma)=\chi(\gamma)+\psi\big(\chi(\gamma)\big)< \chi(\gamma)\in B,$$
and thus $\psi(\gamma)\in B$ since $B$ is convex. Also $\chi(\gamma)=\psi(\gamma)-\psi^2(\gamma)$ for $\gamma\in B$, by Lemma~\ref{contracting}.
%earlier remarks. 
This gives $\psi(\gamma)< \chi(\gamma) < \psi(\gamma)/2$ for
$\gamma\in B$, so $\psi^2(\gamma)=\psi\big(\chi(\gamma)\big)$ for
such $\gamma$, and thus by induction on $m$,
$$\psi^m(\gamma)=\psi\big(\chi^{m-1}(\gamma)\big) \qquad(\gamma\in B,\ m\ge 1).$$
This yields the desired result in view of Lemma~\ref{ps6a}. 

\medskip\noindent
The general case is reduced to this special case by
compositional conjugation. In this connection recall that, unlike
$\psi$, the contraction map $\chi$ is invariant under compositional conjugation.
Take an elementary extension $L$ of $K$ with an element~$\theta\in L^\times$
such that $\psi(B)< v\theta < (\Gamma_L^{>})'$ (so $L^{\theta}$ has
small derivation). 
Let $B_L$ be the convex hull of $B$ in
$\Gamma_L$. Then $\chi(B_L)\subseteq B_L$ and $\psi^\theta(B_L) < 0$, so we are
in the special case considered before, with $P^\theta$, $B_L$, and $L^\theta$ in place of $P$,
$B$, and $K$. This gives $\alpha_0\in \Gamma_L$, $\beta\in B$, and $d_0,\dots, d_r\in \N$, $d_0\le \deg P^{\theta}$, such that for all $y\in L$, 
$$ vy\in B_L^{> \beta}\ \Longrightarrow\  vP^{\theta}(y)\ =\ 
\alpha_0+ d_0vy + \sum_{i=1}^{r} d_i\psi^\theta\big(\chi^{i-1}(vy)\big), $$
and thus with $d:= d_1 +\dots + d_r$ we have for all $y\in L$ with 
$vy\in B_L^{> \beta}$, 
$$  vP(y)\ =\ vP^{\theta}(y)\ =\ (\alpha_0-dv\theta) + d_0vy + \sum_{i=1}^{r} d_i\psi\big(\chi^{i-1}(vy)\big).$$
This holds in particular for $y\in K$ with $vy\in B^{> \beta}$, and so 
$\alpha:=\alpha_0-dv\theta\in \Gamma$. 
\end{proof}

\noindent
By taking $B=\Gamma^{<}$ in Proposition~\ref{ps6b}, we obtain the existence of
some $\beta<0$ such that $P(y)\ne 0$ for all $y\in K$ with $\beta<vy <0$.

\begin{cor}\label{ps6e} Suppose $P(0)=0$. Then there are $\beta<0$ and a strictly increasing function $i\colon (\beta,0) \to \Gamma$ with the intermediate value property such that $vP(y)=i(vy)$ for all $y\in K$ with 
$\beta < vy <0$. 
\end{cor}
\begin{proof} With $B=\Gamma^{<}$ and notations from the proof of Proposition~\ref{ps6b}, the function
$$i\colon B_L \to \Gamma_L, \quad i(\beta):=\alpha + d_0\beta + \sum_{i=1}^r d_i\psi\big(\chi^{i-1}(\beta)\big)$$
is strictly increasing and has the intermediate value property.
Now use that $(\Gamma_L, \psi_L)$ is an elementary extension of $(\Gamma,\psi)$.
\end{proof}

\noindent
We have $\psi\big(\chi^m(\gamma)\big)=\psi(\gamma)-\sum_{i=1}^m \chi^i(\gamma)$ for 
$\gamma<0$.
This yields a formulation of Proposition~\ref{ps6b} for
$B=\Gamma^{<}$ that is more illuminating for $vy$ tending to $0$:

\begin{cor}\label{ps6c} There are 
$\alpha\in \Gamma$ and $d,d_0,e_1,\dots, e_{r-1}\in \N$ with $d_0\le \deg P$, $d\ge e_1\ge \dots \ge e_{r-1}$, such that for some $\beta<0$ we have, for all $y\in K$,
$$\beta < vy < 0\ \Longrightarrow\  vP(y)\ =\ \alpha+ d\psi(vy) + d_0vy  - \sum_{i=1}^{r-1} e_i\chi^{i}(vy).$$ 
If $P(0)=0$ we can take here $d_0\ge 1$. If $P$ is homogeneous, we can take $d_0=\deg P$.
\end{cor}

\noindent 
Here we adopt the convention that $\sum_{i=1}^{r-1} e_i\chi^{i}(\gamma) = 0$ for 
$r\le 1$ and $\gamma <0$. With~$\alpha$ and $d_0,\dots,d_r$ as in 
Proposition~\ref{ps6b} for $B=\Gamma^{<}$ we can take in  
Corollary~\ref{ps6c} the same values for $\alpha$ and $d_0$, and $d=d_1+\dots + d_r,\  e_i=d_{i+1}+ \dots + d_r$. Note:
$$y\in K,\ vy<0 \ \Longrightarrow\ \sum_{i=1}^{r-1} e_i\chi^{i}(vy) = o_{\psi}(vy).$$
By Lemma~\ref{evtcon} there is a unique tuple 
$(\alpha, d, d_0, e_1,\dots, e_{r-1})$ as in Corollary~\ref{ps6c}. 
Here is a slight extension of Corollary~\ref{ps6e}:

\begin{cor}\label{ps6f} Suppose $P(a)=0$, $a\in K$. Then for any $\gamma$ there are $\beta<\gamma$ and a strictly increasing $i\colon (\beta,\gamma) \to \Gamma$ with the intermediate value property such that $vP(a+y)=i(vy)$ for all $y\in K$ with 
$\beta < vy <\gamma$. 
\end{cor}
\begin{proof} Take $g\in K$ with $vg=\gamma$, put $Q:= P_{+a,\times g}=P(a+gY)$,
and apply Corollary~\ref{ps6e} to $Q$ in the role of $P$.
\end{proof}

\subsection*{Notes and comments} The material in this section is new and 
plays a key role in Chapter~\ref{evq}.

%% file: mt-9n-4.tex
\section{Basic Facts about Asymptotic Fields}\label{Asymptotic-Fields-Basic-Facts}

\noindent
Next we consider asymptotic fields in connection with having a
small derivation, $\d$-henselianity, coarsening, and
specialization.
{\em Throughout this section $K$ is an asymptotic field with
asymptotic couple $(\Gamma,\psi)$, valuation ring $\mathcal O$, maximal ideal $\smallo$ of $\mathcal O$,
and residue field~$\k=\mathcal O/\smallo$}. 

\subsection*{The case of a small derivation}
Recall that the derivation of a valued differential field is said to be small
if the derivative of every infinitesimal is infinitesimal.
Any asymptotic field has a compositional conjugate with small derivation, so this enables us to reduce 
problems to the case where the derivation is small. 

\medskip\noindent
Suppose $K$ has small derivation~$\der$.
By Lemma~\ref{closed under der}, not only do we have 
$\der \smallo\subseteq \smallo$, but even 
$\der\mathcal{O}\subseteq \mathcal{O}$, and so $\der$ 
induces a derivation on %the residue field 
$\k$.
As we saw in Sections~\ref{The Gaussian Valuation} and~\ref{sec:ldopv}
this allows us to associate to $P\in K\{Y\}^{\ne}$ and 
$A\in K[\der]^{\ne}$ functions $v_P\colon \Gamma \to \Gamma$ and
$v_A\colon \Gamma \to \Gamma$ such that  
$$v_P(\gamma)=v(P_{\times g}),  \qquad v_A(\gamma)=v(Ag) \qquad(g\in K^\times,\ vg=\gamma).$$ 
For example, for $A=\der$ we have
$$v_\der(\gamma)\ =\ \min(\gamma, \gamma') \text{ if  
$\gamma\in \Gamma^{\ne}$,} \qquad v_\der(0)\ =\ 0.$$
%Section~\ref{pafae} shows that the algebraic 
%closure of any asymptotic field
%is asymptotic; here the algebraic closure is equipped with 
%the unique derivation extending the given one, and any %valuation extending the given one. In view of
%part~(1) of Lemma~\ref{PresInf-Lemma2} it follows that the
%algebraic closure of $K$ has small derivation. 
An asymptotic field has small derivation iff its asymptotic couple has. Thus:

\begin{lemma} If $K$ has small derivation and $L$ is
an asymptotic field extension of~$K$ with $\Gamma_L=\Gamma$, then $L$ has small
derivation. 
\end{lemma}

\noindent
When referring to $\sup \Psi$ we mean the supremum of $\Psi$
in the ordered set $\Gamma_{\infty}$, and the condition
$\sup \Psi=\alpha$ (for $\alpha\in \Gamma_{\infty}$) includes the requirement that $\sup \Psi$ exists. 
Note that if $\sup \Psi=0$, then $K$ has small
derivation. While $\sup \Psi=0$ is a strong restriction, it is often
possible and useful to reduce to this case. (If
$\Gamma=\{0\}$, then $\Psi=\emptyset$, so $\sup \Psi =0$.)
In the next two lemmas we describe situations where $\sup \Psi=0$ holds. The first lemma is in the spirit of
Lemma~\ref{obun}. 

\begin{lemma}\label{ndsupPsi=0} Suppose $K$ has small derivation and the differential residue field $\k$ has nontrivial derivation. Then $\sup \Psi=0$. 
\end{lemma} 
\begin{proof} Take $f\asymp 1$ in $K$ such that $f'\asymp 1$.
Then $g^\dagger\succeq 1$ for all $g\in K^\times$ with $g\not\asymp 1$,
by Proposition~\ref{CharacterizationAsymptoticFields}, and thus
$\Psi\le 0$. It remains to use Lemma~\ref{PresInf-Lemma}. 
\end{proof}

%The next two lemmas relate the
%condition $\sup \Psi=0$ to $\der$ being neatly surjective. 

\begin{lemma}\label{supPsi=0} Suppose $\der\smallo =\smallo$. Then 
$(\Gamma^{>})'= \Gamma^{>}$, and thus $\sup \Psi=0$. Also,
$C$ maps onto the constant field $C_{\k}$ of 
$\k$ under the residue map $\mathcal{O} \to \k$.  
\end{lemma}

\noindent 
This is an easy consequence of Corollary~\ref{asy1}.
The condition $\der \smallo =\smallo$ is satisfied when~$K$ has small derivation and $\der$ is neatly surjective. The following partial 
converse to Lemma~\ref{supPsi=0} 
is useful in verifying the hypothesis of Lemma~\ref{nslinsur}.

\begin{lemma}\label{sup=0neatsur} Suppose $\sup \Psi=0$, $C$ maps onto the constant field $C_{\k}$ of 
$\k$ under the residue map $\mathcal{O} \to \k$, and $\der K = K$. Then 
$\der$ is neatly surjective.
% there is for each
%$b\in K^\times$ an $a\in K^\times$ such that $a'=b$ and $v_{\der}(va)=vb$.
\end{lemma} 
\begin{proof} The condition $\sup \Psi=0$ gives $v_{\der}(\gamma) = \gamma'$ 
for all $\gamma\in \Gamma$, where 
$0'=0\in \Gamma$ by convention. Hence the only  
nonobvious case is when $a, b\in K^\times$,
$a'=b$ and $a\asymp 1$, $b\prec 1$. Then the image $\bar{a}$ of $a$ in $\k$
satisfies $\bar{a}'=0$, so $\bar{a}=\bar{c}$ with $c\in C$, by our 
assumption on the constant field of $\k$. Hence $(a-c)'=b$ and $a-c\prec 1$,
so  $v_{\der}\big(v(a-c)\big)=vb$, as desired.
\end{proof}

\noindent
In cases like $K=\mathbb{T}$ where $\sup \Psi=0$ fails, it is often useful to arrange it by compositional conjugation and coarsening. 
%with $\Gamma\neq\{0\}$.

\begin{lemma}\label{impreservesas} 
Suppose $\sup \Psi =0$ and $L$ is an immediate
valued differential field extension
of $K$ with small derivation. Then $L$ is asymptotic.
\end{lemma}
\begin{proof} We use the equivalence of
(i) and (ii) in Proposition~\ref{CharacterizationAsymptoticFields}. Let $b\in L^\times$ be such that $b\nasymp 1$. Take $a\in K^\times$ with
$b=a(1+\varepsilon)$, $\varepsilon\prec 1$. Then 
$b^\dagger =a^\dagger + \frac{\varepsilon'}{1+\varepsilon}$ 
with $a^\dagger \succeq 1$ and 
$\frac{\varepsilon'}{1+\varepsilon}\prec  1$, so 
$b^\dagger \sim a^\dagger$, and thus $\psi(vb)=\psi(va)=v(a^\dagger)=v(b^\dagger)$. 
\end{proof}

\noindent
The condition $\sup \Psi=0$ arises naturally in some situations:

\begin{lemma}\label{lem:K<Y> asymptotic} Suppose that $\der\smallo\subseteq \smallo$ and that $K\<Y\>$ equipped with its gaussian valuation is asymptotic. Then $\sup\Psi=0$.
\end{lemma}
\begin{proof} We have $Y\asymp 1$ in $K\<Y\>$, so
$f'\prec 1\asymp Y'\preceq g^\dagger$ for all $f,g\in K^\times$ with $f\prec 1$ and $g\not\asymp 1$, by (i)~$\Leftrightarrow$~(iv) in 
Proposition~\ref{CharacterizationAsymptoticFields}. Thus $\sup\Psi=0$.
\end{proof}

\begin{lemma}\label{lem:K<a>asymp} Assume $\der\smallo\subseteq \smallo$.
Let $K\<a\>$ be as in Theorem~\ref{resext} 
with minimal annihilator $F$ of $a$ over $K$ of order $r\geq 1$, and
$\overline{F}\notin \k^\times  Y'$. 
If $K\<a\>$ is asymptotic, then $\sup\Psi=0$.
\end{lemma}
\begin{proof} We have $a\asymp 1$. We claim that $a'\asymp 1$.
This claim clearly holds for $r\ge 2$, and also when $r=1$ and 
$F$ has degree $\ge 2$ in $Y'$. Suppose $r=1$ and $F$ has degree~$1$ in $Y'$. Then $F=IY'+G$ with
$I, G\in \mathcal{O}[Y],\ vI=0$.  The assumption on $\overline{F}$ and its irreducibility in $\k\{Y\}$ yields $vG=0$, and so again $a'\asymp 1$. With the claim established, argue
as in the proof of Lemma~\ref{lem:K<Y> asymptotic} with $a$ replacing $Y$.
\end{proof}

\noindent
We now turn to the case $\max \Psi=0$:

\begin{lemma}\label{KgroundedLasymp} Assume $\max \Psi=0$.
Let $L$ be a valued differential field extension of $K$ with small derivation, $\Gamma_L=\Gamma$, and $(\k_L^\times)^\dagger\cap \k \subseteq (\k^\times)^\dagger$. Then $L$ is asymptotic.
\end{lemma}
\begin{proof} Let $g\in L^\times$, $g\prec 1$. By Proposition~\ref{CharacterizationAsymptoticFields} it is enough to show that $v(g^\dagger)=\psi(vg)$.
We have $u\in L$, $f\in K$
such that $g=uf$, $u\asymp 1$, and $f\asymp g$. 
Then $g^\dagger = u^\dagger +f^\dagger$,
and  $u^\dagger\preceq 1$, $f^\dagger\succeq 1$.
If $u^\dagger \prec 1$ or $f^\dagger \succ 1$, then $v(g^\dagger)=v(f^\dagger)=\psi(vf)=\psi(vg)$. 
It remains to consider the case that $u^\dagger\asymp 1$ and
$f^\dagger\asymp 1$. If $u^\dagger + f^\dagger\asymp 1$, we
have again $v(g^\dagger)=v(f^\dagger)=\psi(vf)=\psi(vg)$,
so assume $u^\dagger + f^\dagger\prec 1$. Then
$\bar{u}^\dagger\in (\k_L^\times)^\dagger\cap \k$, and so
we have $\phi\in K$ such that 
$\phi\asymp 1$ and $u^\dagger \sim \phi^\dagger$. Then
$0\ne \phi f\prec 1$ and $(\phi f)^\dagger=\phi^\dagger + f^\dagger \prec 1$, contradicting $\max \Psi=0$. 
\end{proof}

\begin{cor} Suppose that $\max \Psi=0$. Then $K\<Y\>$ equipped with its gaussian valuation is asymptotic.
\end{cor}
\begin{proof} Apply Lemma~\ref{KgroundedLasymp} to $L=K\<Y\>$,
and Corollary~\ref{gaussdagger} with $K$ as $\k$.
\end{proof}

\subsection*{Differential-henselian asymptotic fields}
{\em In this subsection we assume that the asymptotic field $K$ has small derivation. Let~$\k$ be its differential residue field.}\/ Recall that if $K$ is $1$-differential-henselian, then~$\mathcal{O}$ is $1$-linearly surjective.
Thus by the next lemma, if~$K$ is $1$-differential-henselian, then~${\sup \Psi=0}$.

\begin{lemma}\label{dhns3} Suppose $\mathcal{O}$ is $1$-linearly surjective. Then $\der\smallo=\smallo$, 
$\sup \Psi = 0$, and the field embedding 
$a\mapsto \bar{a}\colon C \to C_{\k}$ is bijective.
\end{lemma}
\begin{proof} From $v(\der)=0$ we get $v_{\der}(\Gamma^{>})=\Gamma^{>}$ by the Equalizer Theorem. Since $\mathcal{O}$ is $1$-linearly surjective, $\der$ is neatly surjective, so $\der\smallo=\smallo$. Now use Lemma~\ref{supPsi=0}.
\end{proof}

\noindent
Since $\sup \Psi = 0$ fails for $\mathbb{T}$, this asymptotic field 
is not $\d$-henselian. Chapter~\ref{ch:newtdirun} and Lemma~\ref{diffnewt}, however, show that $\mathbb{T}$ becomes $\d$-henselian if we coarsen its valuation by any convex subgroup of $\Gamma_{\T}$ that contains $v(x)$.

%\begin{lemma}\label{dhres} Suppose $K$ is 
%$n$-differential-henselian, $n\ge 1$, and
%identify $C$ and $C_{\k}$ via the field isomorphism 
%$a\mapsto \bar{a}%\colon C \to C_{\k}$.
%Let $A\in K[\der]$ have order $\le n$, with $vA=0$. Then the map
%$$ a \mapsto \overline{a} \colon\ \mathcal{O}\cap \ker A\ \to\ 
%\ker \bar{A}$$
%is $C$-linear and surjective.
%\end{lemma}
%\begin{proof} Let $A=a_0 + a_1\der+ \dots + a_n\der^n$ with 
%$a_0,\dots, a_n\in K$. Set $$P(Y)\ :=\ a_0Y+ \dots + a_nY^{(n)}.$$ 
%Suppose $a\in \mathcal{O}$, $v(a)=0$ and 
%$\overline{A}(\overline{a})=0$. Then $P$ is in dh-position 
%at $a$, so we have $y\in \mathcal{O}$ such that
%$P(y)=0$ and $\overline{y}=\overline{a}$. Then $A(y)=0$.
%\end{proof}  

\begin{cor}\label{asexistdhenselization} Assume $\k$ is linearly surjective. Then $K$ has an immediate asymptotic field extension $L$ such that: \begin{enumerate}
\item[\textup{(i)}] $L$ is $\d$-algebraic over $K$;
\item[\textup{(ii)}] $L$ has small derivation and is $\d$-henselian;
\item[\textup{(iii)}] no proper differential subfield of $L$ containing $K$ is $\d$-henselian.
\end{enumerate}
\end{cor}
\begin{proof} Corollary~\ref{henimm} yields an immediate valued differential field extension~$F$ of~$K$ with small derivation such that $F$ is $\d$-henselian. Then $F$ is asymptotic by Lem\-ma~\ref{impreservesas}. 
By keeping only the elements of 
$F$ that are $\d$-algebraic over $K$ we arrange that $F$ is  $\d$-algebraic over $K$. 

 Let $L$ be the intersection inside
$F$ of the collection of all differential subfields~$E$ of $F$ that
contain $K$ and are $\d$-henselian. 
Proposition~\ref{noextrazeros} applied to the
ex\-ten\-sions~${E\subseteq F}$ shows that $L$ has the desired property.
\end{proof} 

\noindent
We conjecture that there is only one $L$ as in 
Corollary~\ref{asexistdhenselization}, up to isomorphism over~$K$. For any $r\in \N$, Corollary~\ref{asexistdhenselization}
goes through with ``linearly surjective'' and ``$\d$-henselian'' replaced by
``$r$-linearly surjective'' and ``$r$-$\d$-henselian.''

%This raises two questions: 
%\begin{enumerate}
%\item[\textup{(i)}] Are there $\d$-henselian $K$ with 
%$\Gamma\ne \{0\}$ and $0\notin \Psi$? 
%\item[\textup{(ii)}] Are there $\d$-henselian $K$ with $0\in \Psi$?
%\end{enumerate}
%Both questions have a positive answer, as we now indicate.

%\marginpar{commented out material below has been checked}
%For (ii) we start with a $\d$-henselian $K$ such that 
%$\Gamma\ne \{0\}$ and $0\notin \Psi$. 
%Take $\gamma\in \Psi$, and set
%$\delta:= \gamma^\dagger$, so 
%$\gamma < \delta < 0$ and $\delta=o(\gamma)$. Set 
%$$\Delta\ :=\ \{\alpha\in \Gamma:\ |\alpha|\le n|\delta| 
%\text{ for some }n\}.$$ Then $\Delta$ is a convex 
%subgroup of $\Gamma$ and $\delta, \delta^\dagger\in \Delta$. 
%Also $\gamma\notin \Delta$ and $\gamma^\dagger\in \Delta$, 
%so the $\Psi$-set $\Psi_{\Delta}$ of the $\Delta$-coarsening %of $K$ has largest element $0\in \Gamma/\Delta$, and 
%this $%\Delta$-coarsening of
%$K$ is asymptotic and $\d$-henselian, by 
%Corollary~\ref{coarsen} and Lemma~\ref{inva3}.   

\subsection*{Specialization of asymptotic fields}
{\em Let $K$ have small derivation, and let~$\Delta$ be a convex subgroup
of $\Gamma=v(K^\times)$}. 
It follows from Corollary~\ref{cor:Cohn-Lemma} 
that
the valuation ring $\dot{\mathcal{O}}$ of the coarsening 
$\dot{v}=v_{\Delta}\colon K^\times \to \dot{\Gamma}$ of~$v$  
and the maximal ideal $\dot{\smallo}$
of~$\dot{\mathcal{O}}$
are closed under $\der$. Note that if $\sup \Psi=0$, then 
$\psi(\Delta^{\ne}) \subseteq \Delta$ and $\sup \psi(\Delta^{\ne})=0$ in $\Delta$, by 
Lemma~\ref{PresInf-Lemma2}(iv).

\medskip\noindent
{\em In the rest of this subsection we assume that 
$\psi(\Delta^{\ne}) \subseteq \Delta$}.
The residue 
field $\dot K= \dot{\mathcal{O}}/\dot{\smallo}$
of $\dot{\mathcal{O}}$ with its induced valuation 
$v\colon \dot K^\times \to \Delta$
and the induced derivation is then an asymptotic field
with asymptotic couple~$\big(\Delta,\psi|\Delta^{\ne}\big)$. (This follows from 
the equivalence of (i) and (ii) in 
Proposition~\ref{CharacterizationAsymptoticFields}.) Moreover, $\dot K$ 
has small derivation, the field isomorphism 
$\res(K) \to \res(\dot K)$ from 
Section~\ref{sec:decomposition} is a differential field isomorphism, and if $K$ is
of $H$-type, then so is $\dot{K}$.
%We identify $C$ with a subfield of $C_{\dot K}$ via the 
%field embedding $c\mapsto\dot c\colon C\to C_{\dot K}$.

\begin{lemma}\label{specasymptint} Suppose $K$ has asymptotic integration and $\Delta\ne \{0\}$.
Then $\dot{K}$ has asymptotic integration.
\end{lemma}
\begin{proof} Take $\alpha\in \Delta^{\ne}$, and note that $\alpha+\psi(\alpha)\in \Delta$. Let $\delta\in \Delta$, $\delta\ne \alpha+\psi(\alpha)$, and take $\beta\in \Gamma^{\ne}$ with $\beta + \psi(\beta)=\delta$. Then by 
Lemma~\ref{BasicProperties}(ii),  
$$\big(\alpha + \psi(\alpha)\big)-\big(\beta+\psi(\beta)\big)=\alpha-\beta + o(\alpha-\beta)\in \Delta,$$
so $\alpha-\beta\in \Delta$, and thus $\beta\in \Delta$. 
\end{proof}

\index{closed!under integration}
\index{closed!under logarithms}
\index{differential field!closed under integration}
\index{differential field!closed under logarithms}

\noindent
We say that the differential field $F$ is {\bf closed under integration}
if for every $f\in F$ there is a $y\in F$ with $y'=f$, and we say that
$F$ is {\bf closed under logarithms} if for every $f\in F^\times$ 
there is a $y\in F$ with $y'=f^\dagger$. So if $F$ is closed under
integration, then~$F$ is closed under logarithms. 
%Every Liouville closed
%$H$-field is closed under integration.
If the differential field $F$ is closed under integration 
(closed under logarithms)
and $\phi\in F^\times$, then $F^\phi$ is closed under integration
(closed under logarithms, respectively). 

\begin{lemma} Suppose $\Delta\ne \{0\}$.
If $K$ is closed under integration, then so is $\dot{K}$. If~$K$ is
closed under logarithms, then so is $\dot{K}$.
\end{lemma}
\begin{proof}
Take $\delta\in\Delta^{\ne}$; then $\psi(\delta)\in\Delta$.
Suppose $K$ is closed under integration, and
let $f\in \dot{\cal{O}}$. We need to find $y\in\dot{\cal{O}}$ with $y'-f\in
\dot{\smallo}$. If $vf>\Delta$, then we may take $y:=0$, so we may assume
$vf\in\Delta$. Take $y\in K$ with $y'=f$; we claim that $y\in\dot{\cal{O}}$.
Suppose $y\notin \dot{\cal{O}}$. Then $\gamma:=vy<\Delta$, 
in particular 
$[\gamma]>[\delta]$
and hence $\big[\psi(\gamma)-\psi(\delta)\big]<[\gamma-\delta]=[\gamma]$ by
Lemma~\ref{BasicProperties}(ii). 
Since $\gamma+\psi(\gamma)=vf\in
\Delta$, we therefore have $[\gamma]=\big[\gamma+\psi(\gamma)-\psi(\delta)\big]\in [\Delta]$ and hence $\gamma\in\Delta$, a contradiction. 
Thus $\dot{K}$ is closed under integration.
Similarly one shows that if $K$ is
closed under logarithms, then so is $\dot{K}$.
\end{proof} 

\noindent
Let $L$ be an asymptotic field extension of $K$  with small derivation. 
Let $\Delta_L$ be the convex hull of $\Delta$ in $\Gamma_L$, and let
$$ \dot{\mathcal{O}}_L:=
  \{y\in L: \text{$vy\ge \delta$  for some $\delta\in \Delta_L$}\}
$$
be the valuation ring of the corresponding coarsening 
$\dot{v}\colon L^\times \to \dot{\Gamma}_L= \Gamma_L/\Delta_L$.  
Then~$\dot{\mathcal{O}}_L$ and its maximal ideal 
$\dot{\smallo}_L$
are closed under the derivation $\der$ of $L$. Note that~$\dot{\mathcal{O}}_L$ lies over $\dot{\mathcal{O}}$. Let 
$\dot L:=  \dot{\mathcal{O}}_L/\dot{\smallo}_L$ be the residue 
field of $\dot{\mathcal{O}}_L$, equipped with the induced valuation
$v\colon {\dot L}^\times \to \Delta_L$. We have
$\psi_L(\Delta_L^{\ne}) \subseteq \Delta_L$, so the field 
$\dot L$ with its valuation $v$
and the induced derivation is an asymptotic field
with asymptotic couple~$\big(\Delta_L,\psi_L|\Delta_L^{\ne}\big)$. With the usual
identifications we have $\dot K\subseteq \dot L$, not only 
as (residue) fields, 
but also as asymptotic fields.

\subsection*{Flattening} For certain {\em definable} ways of coarsening we 
use special notation and terminology, borrowed from \cite{JvdH}. These
particular coarsenings are called flattenings and 
will be very useful in Chapters~\ref{ch:The Dominant Part and the Newton Polynomial} and~\ref{ch:newtonian fields}. (Some readers might prefer to skip
this subsection until it gets used in those chapters.)  
{\em Assume $K$ is of $H$-type, and let 
$\gamma$ range over $\Gamma$.}\/
We have the convex subgroup 
$$\Gamma^\flat\ :=\ 
\bigl\{\gamma :\  \psi(\gamma)>0\bigr\}$$
of $\Gamma$. 
For example, if $K=\mathbb T$, then
$$\Gamma^\flat\ =\ \{vh:\  h\in K^\times,\ h\flatter \ex^x\}.$$
Let
$$v^\flat\ \colon\ K^\times\to \Gamma^\sharp
\qquad\text{where $\ \Gamma^\sharp\ :=\ \Gamma/\Gamma^\flat$}$$ 
be the coarsening of the valuation $v$ by $\Gamma^\flat$, with associated dominance relation $\preceq^\flat$ on~$K$, and valuation ring $\mathcal{O}^\flat$. We call the valuation $v^\flat$ the {\bf flattening} of $v$. 
Note that the definition of $\Gamma^\flat$ depends on the derivation of $K$.

\index{flattening}
\index{valuation!flattening}
\nomenclature[Z]{$\Gamma^\flat$}{$\{\gamma:\psi(\gamma)>0\}$}
\nomenclature[Z]{$v^\flat$}{flattening of the valuation $v$ of an $H$-asymptotic field}
\nomenclature[Rx]{$\preceq^\flat$}{dominance relation associated to $v^\flat$}
\nomenclature[Z]{$\Gamma^\sharp$}{$\Gamma^\sharp=\Gamma/\Gamma^\flat$}

Let $\phi\in K^\times$, and use a subscript $\phi$ to indicate the flattened objects $\Gamma^\flat_\phi$, $v^\flat_\phi$, $\preceq^\flat_\phi$, $\prec^\flat_\phi$, $\asymp^\flat_\phi$, and $\simflat_\phi$ associated to the asymptotic field $K^\phi$. In particular
$$\Gamma^\flat_\phi\ =\  \big\{ \gamma :\  
\psi(\gamma) > v\phi \big\}.$$
Clearly $\Gamma^\flat_\phi=\{0\}$ iff
$v\phi\notin \big(\Gamma^{<}\big)'$, and $\Gamma^\flat_\phi=\Gamma$
iff $\Psi>v\phi$. The flattening $v^\flat_\phi$ of the valuation $v$ on 
$K^{\phi}$ will be useful for ungrounded
$K$ and $v\phi\in \Psi^{\downarrow}$. In that case, 
with~$v\phi$ increasing cofinally in $\Psi^{\downarrow}$, the convex subgroup
$\Gamma^\flat_\phi$ of $\Gamma$ becomes arbitrarily small, and so
$v^\flat_\phi$ (with value group $\Gamma/\Gamma^\flat_{\phi}$) approximates $v$ in some sense. 

We denote the gaussian extension of
$v^\flat_{\phi}$ to $K^\phi\<Y\>$ also by $v^\flat_{\phi}$. The binary relations $\preceq^\flat_\phi$, $\prec^\flat_\phi$, $\asymp^\flat_\phi$, and $\simflat_\phi$,
on $K$ extend likewise, with the same notations, to binary relations on $K^{\phi}\<Y\>$. As $K\<Y\>$ and 
$K^{\phi}\<Y\>$ have the same underlying field, we also use
the same notations for elements of $K\<Y\>$; for example, 
given $f\in K$ and $P\in K\<Y\>$, $f\preceq^\flat_\phi P$ means that $v^\flat_{\phi}(f) \ge v^\flat_{\phi}(P)$ in $\Gamma/\Gamma^\flat_{\phi}$, where in the latter inequality
$f$ and $P$ are taken as elements of $K^\phi\subseteq K^\phi\<Y\>$ and 
$K^{\phi}\<Y\>$, respectively.

\nomenclature[Z]{$\Gamma^\flat_\phi$}{$\{\gamma:\psi(\gamma)>v\phi\}$}
\nomenclature[Z]{$v^\flat_\phi$}{flattening of the valuation $v$ of $K^\phi$}
\nomenclature[Rx]{$\preceq^\flat_\phi$}{dominance relation associated to $v^\flat_\phi$}

\begin{lemma}
Suppose $\Phi\in K^\times$ with $\Phi\nasymp 1$ and 
$\Phi^\dagger\asymp\phi$, and $f,g\in K$. Then
\begin{align*}
f\preceq^\flat_\phi g &\quad\Longleftrightarrow\quad v^\flat_\phi(f)\geq v^\flat_\phi(g) 
\quad\Longleftrightarrow\quad \text{$f\preceq hg$ for some $h\in K^\times$ with $h\flatter\Phi$,}\\
f\prec^\flat_\phi g &\quad\Longleftrightarrow\quad v^\flat_\phi(f)> v^\flat_\phi(g) 
\quad\Longleftrightarrow\quad \text{$f\prec hg$ for all $h\in K^\times$ with $h\flatter\Phi$,}\\
f\asymp^\flat_\phi g &\quad\Longleftrightarrow\quad v^\flat_\phi(f)=v^\flat_\phi(g) 
\quad\Longleftrightarrow\quad \text{$f\asymp hg$  for some $h\in K^\times$ with $h\flatter \Phi$.}
\end{align*}
\end{lemma}
\begin{proof}
These equivalences follow from their validity for $g=1$, and in that case
we argue as follows:
\begin{align*}
f\preceq^\flat_\phi 1 &\quad\Longleftrightarrow\quad \text{$vf\geq vh$
for some $h\in K^\times$ with $vh\in\Gamma^\flat_\phi$} \\
&\quad\Longleftrightarrow\quad \text{$f\preceq h$ for some $h\in K^\times$ with $h\flatter \Phi$,}
\end{align*}
where we use that $1\flatter \Phi$ to get the last equivalence when $vh=0$.
We proceed likewise with $f\prec^\flat_\phi 1$ and $f\asymp^\flat_\phi 1$, treating the case $vh=0$ separately. 
\end{proof}

\noindent
The
valuation ring of $v^\flat_\phi$ is
$$\cal{O}^\flat_\phi\ =\ \big\{ f\in K:\  \text{$vf\geq 0$ or $v(f')>vf + v\phi$} \big\},$$
with maximal ideal 
$$\smallo^\flat_\phi\ =\ \big\{ f\in K:\  \text{$vf>0$ and $v(f')\leq vf + v\phi$} \big\}.$$
Note that if $\Gamma^\flat=\{0\}$, then the canonical map $\Gamma\to\Gamma^\sharp$ is an isomorphism
of ordered groups via which $\Gamma$ is identified with $\Gamma^\sharp$.
If $\Gamma^\flat\neq\{0\}$, then $K$ has small derivation and 
there exists a unique positive element $1$ in
$\Gamma^\flat$ with $\psi(1)=1$, so $\psi\big((\Gamma^\flat)^{\ne}\big)\subseteq\Gamma^\flat$
by Lemma~\ref{asymp-lemma1}(i). Note also that if $\Gamma^\flat\neq\{0\}$, 
%and $\Gamma^\flat\neq \Gamma$, 
then by Lemma~\ref{Coarsening-Prop} we have
$\sup \Psi^\flat=0$ in $\Gamma^\sharp$, where
$$\Psi^\flat:= \big\{v^\flat(f^\dagger):\ f\in K^\times,\ f\not\asymp^\flat 1\big\}.$$
If $\Gamma^\flat\neq\{0\}$, then 
we denote the asymptotic residue field of $v^\flat$ by
$K^\flat := \cal{O}^\flat/\smallo^\flat$, with asymptotic couple  
$\big(\Gamma^\flat,\psi|\Gamma^\flat\big)$.

\nomenclature[Z]{$\cal{O}^\flat_\phi$}{valuation ring of $v^\flat_\phi$}
\nomenclature[Z]{$\smallo^\flat_\phi$}{maximal ideal of the valuation ring of $v^\flat_\phi$}

\medskip\noindent
Let $\Phi\in K^\times$ with $\Phi\nasymp 1$. The binary relations 
$\preceq^\flat_{\Phi^\dagger}$, $\prec^\flat_{\Phi^\dagger}$, $\asymp^\flat_{\Phi^\dagger}$ on $K$, and on $K\<Y\>$,  are often denoted by
$\preceq_{\Phi}$, $\prec_{\Phi}$, $\asymp_{\Phi}$, respectively. Thus
for $f,g\in K$: 
\begin{align*}
f\preceq_\Phi g &\quad\Longleftrightarrow\quad  \text{$f\preceq hg$ for some $h\in K^\times$ with $h\flatter\Phi$,}\\
f\prec_\Phi g &\quad\Longleftrightarrow\quad    \text{$f\prec hg$ for all $h\in K^\times$ with $h\flatter\Phi$,}\\
f\asymp_\Phi g &\quad\Longleftrightarrow\quad \text{$f\asymp hg$  for some $h\in K^\times$ with $h\flatter \Phi$.}
\end{align*}
An advantage of this notation is that the relations $\preceq_{\Phi}$, $\prec_{\Phi}$, $\asymp_{\Phi}$ do not change in passing from $K$ to
a compositional conjugate $K^\phi$ ($\phi\in K^\times$). Note that $\Phi\nasymp_\Phi 1$ by Corollary~\ref{CharacterizationAsymptoticFields-Corollary}(iv), and if $\Phi\prec 1$, then $\Phi\prec_{\Phi} 1$. Also, for $f\nasymp 1$ in $K^{\times}$,
$$f\ \asymp_\Phi\ 1\ \Longleftrightarrow\ f\ \flatter\ \Phi.$$ 
Here are some further rules:

\nomenclature[Rx]{$\preceq_\Phi$}{dominance relation $\preceq^\flat_{\Phi^\dagger}$, for $\Phi\in K^\times$, $\Phi\nasymp 1$}

\begin{lemma}\label{symflat} Let $f,g\in K^\times,\ f, g\nasymp 1,\  f\prec_{g} g$. Then $f\prec_{f} g$.
\end{lemma}
\begin{proof} We have $f/g\prec^{\flat}_{g^\dagger} 1$, so $\psi(vf - vg)\le \psi(vg)$, hence $\psi(vf-vg)\le \psi(vf)$, which in view of $vf> vg$ yields $f/g\prec^{\flat}_{f^\dagger} 1$, so $f\prec_{f} g$.
\end{proof}

\begin{lemma}\label{lem:flat equivalence}
Let $\Phi_1,\Phi_2\in K^\times$ with $\Phi_1,\Phi_2\nasymp 1$ and $\Phi_1 \flattereq \Phi_2$. Then for $f,g\in K$:
$$f\preceq_{\Phi_1} g  \ \Rightarrow\ f\preceq_{\Phi_2} g\qquad\text{and}\qquad
f\prec_{\Phi_2} g \ \Rightarrow\  f\prec_{\Phi_1} g.$$
\end{lemma}

\begin{lemma}\label{cor:flat equivalence}
Suppose $\Phi\in K^\times$, $\Phi \nasymp^\flat 1$. 
Then for all $f,g\in K$:
$$f\preceq^\flat g  \ \Rightarrow\ f\preceq_{\Phi} g\qquad\text{and}\qquad
f\prec_{\Phi} g \ \Rightarrow\  f\prec^\flat g.$$
\end{lemma}
\begin{proof} The second part follows from the first. To prove the first part, arrange $f\neq 0$ and $g=1$.
Assume $f\preceq^{\flat} 1$. Then $f\preceq 1$ or 
$f^\dagger \prec 1$. From
$\Phi \nasymp^\flat 1$ we get $\Phi^\dagger \succeq 1$, so
$f\preceq 1$ or $f^{\dagger}\prec \Phi^\dagger$, and thus $f\preceq_{\Phi}1$.
\end{proof}

%\noindent
%Let $\Phi\in K^\times$, $\Phi\not\asymp 1$.
%Given $f,g\in K$, the truth value of each of the conditions
%$$f\asymp^{\flat}g,\ f\preceq^{\flat} g,\ f\prec^{\flat} g,\ $f\asymp_{\Phi}g,\ f\preceq_{\Phi} g,\ f\prec_{\Phi} g$$
%depends only on the pair $(vf, vg)$. This allows us 
%to extend each of
%the binary relations $\asymp^{\flat}$, $\preceq^{\flat}$, 
%$\prec^{\flat}$, 
%$\asymp_{\Phi}$, $\preceq_{\Phi}$,
%$\prec_{\Phi}$ on $K$ to a binary relation on $K\{Y\}$, 
%to be denoted by the same symbol: for $P,Q\in K\{Y\}$, we %declare $P\asymp^{\flat} Q$ to mean $f\asymp^{\flat} g$
%with $f,g\in K$ such that $vP=vf$, $vQ=vg$, and 
%likewise with the other relations.

\noindent
{\em In the next lemma and its corollaries below
we assume that $K$ is ungrounded and has small derivation. Also $P\in K\{Y\}^{\neq}$ and $g\in K^\times$, $g\nasymp 1$}.

\begin{lemma}\label{lem:Pg flat equivalence}
Suppose $P$ is homogeneous of degree $d$.  Then 
$P_{\times g} \asymp_g P\,g^d$.
\end{lemma}
\begin{proof}
If $g^\dagger\preceq 1$, then $P_{\times g}\asymp Pg^d$ by Lemma~\ref{vplemma}. 
If $g^\dagger\succ 1$, then $g^{\dagger\dagger}\prec g^\dagger$ by 
Lemma~\ref{PresInf-Lemma2}(iv), and thus
$P_{\times g} \asymp_g Pg^d$ by Proposition~\ref{ps1}.
\end{proof}

\begin{cor}\label{cor:Pg flat equivalence, 1}
Suppose $g\prec 1$ and $d:=\dv(P)=\val(P)$. 
Then
$$P_{\times g}\ \asymp_g  P\,g^d.$$ 
\end{cor}
\begin{proof}
We have $P_i\preceq P_d$ for $i\ge d$, and $g\prec_g 1$, hence
$P_i\,g^i\prec_g P_d\,g^d$ for $i>d$, and the claim follows from
Lemma~\ref{lem:Pg flat equivalence}.
\end{proof}

\noindent
For $F\in K\{Y\}$ and $d\in\N$ we put 
$$F_{\leq d}\ :=\ F_0+F_1+\cdots+F_d, \qquad F_{>d}\ :=\ F-F_{\leq d}.$$
With this notation, we have:

\begin{cor}\label{cor:Pg flat equivalence, 2}
Suppose $g\succ 1$ and $d:=\ddeg(P)=\ddeg(P_{\times g})$. Then $$gP_{>d}\ \preceq_g\  P.$$
\end{cor}
\begin{proof} The case $P_{>d}=0$ being trivial, assume $P_{>d}\ne 0$.
Take $i>d$ such that $P_i \asymp P_{>d}$. Then by Lemma~\ref{lem:Pg flat equivalence} and $g\succ 1$ we have
$$(P_{\times g})_i\ =\ (P_i)_{\times g}\ \asymp_g\ P_i\,g^i\ \succeq\ P_i\,g^{d+1}\ \asymp\ g^{d+1}P_{>d}$$ and so
$(P_{\times g})_{>d} \succeq_g  g^{d+1}P_{>d}$. Since $\ddeg(P_{\times g})=d$, we have  
$(P_{\times g})_d\succeq  (P_{\times g})_{>d}$, and by Lemma~\ref{lem:Pg flat equivalence} again,
$(P_{\times g})_d = (P_d)_{\times g} \asymp_g g^dP_d$. Since $\ddeg(P)=d$, we have $P\asymp P_d$, so
$$ g^dP\ \asymp\ g^dP_d\ \asymp_g\ (P_{\times g})_d\ \succeq\ (P_{\times g})_{>d}\ \succeq_g\ g^{d+1}P_{>d},$$ and thus
$P\succeq_g gP_{>d}$ as claimed.
\end{proof}

\begin{cor}\label{cor:Pg flat equivalence, 3}
Let $f\in K^\times$ with $f\nasymp 1$ and $f\flatter g$. Then $P_{\times f} \asymp_g P$.
\end{cor}
\begin{proof}
Take $d\in\N$ with $P_{\times f}\asymp (P_d)_{\times f}$.
Lemma~\ref{lem:Pg flat equivalence} gives 
$(P_d)_{\times f}\asymp_f P_df^d$, so 
$P_{\times f} \asymp_f P_d\, f^d$, which by Lemma~\ref{lem:flat equivalence} yields $P_{\times f}\asymp_g P_df^d$. From
$f\flatter g$ we also get $f\asymp_g 1$, so
$P_d\, f^d \asymp_g P_d \preceq P$, hence $P_{\times f}\preceq_g P$. Applying the same
argument to $P_{\times f}$, $f^{-1}$ in place of $P$, $f$ yields $P=(P_{\times f})_{\times f^{-1}}\preceq_g P_{\times f}$.
\end{proof}

%% file: mt-9n-5.tex
\section{Algebraic Extensions of Asymptotic Fields}
\label{pafae}
 
\noindent
In this section we prove analogues for asymptotic fields of
results in Section~\ref{sec:algebraicext} on algebraic extensions of 
valued differential fields with small derivation.  

\subsection*{Immediate algebraic extensions} {\em In this subsection we fix a valued differential field $L$ with a subfield $K$.
We consider~$K$ as a valued subfield of $L$ with valuation ring
$\mathcal{O}$, maximal ideal $\smallo$ of $\mathcal{O}$,
and value group~$\Gamma$.}\/ 
In order for $L$ to be asymptotic it is of course necessary that 
for all $f,g\in K^\times$ with $f,g\prec 1$ we have
$$f\prec g\quad \Longleftrightarrow\quad f'\prec g'.$$
Call $K$ {\bf asymptotic in  $L$} if this condition is satisfied. 
Below we prove:

\index{valued differential field!asymptotic in}

\begin{prop}\label{pa5} Suppose $K$ is asymptotic in $L$ and $L$ 
is an immediate 
algebraic extension of $K$. Then $L$ is asymptotic.
\end{prop}

\noindent
For later use it is important that in this proposition 
we do not require $K$ to be closed under the derivation of $L$. 
First some remarks and a lemma.
Suppose $K$ is asymptotic in $L$, and set $C:=\{a\in K:\ a'=0\}$, a subfield
of $K$. Then
$C\cap \smallo =\{0\}$ by the same argument as for
asymptotic fields. Thus the valuation $v$ is trivial on~$C$, and
the map $\varepsilon\mapsto\varepsilon'\colon \smallo\to L$ is injective.
A variant of Proposition~\ref{CharacterizationAsymptoticFields} goes through with the same proof:

\begin{lemma}\label{pa1} Suppose $v(L^\times)=v(K^\times)~(=\Gamma)$.
The following are equivalent:
\begin{enumerate}
\item[\textup{(i)}] $K$ is asymptotic in $L$;
\item[\textup{(ii)}] there is an asymptotic couple $(\Gamma,\psi)$ 
such that for all $f\in K^\times$ with $f\nasymp 1$ we have
$\psi(vf)=v(f^\dagger)$;
\item[\textup{(iii)}] for all $f,g\in K^\times$ with
$f,g\nasymp 1$ we have: $f\preceq g \Longleftrightarrow f'\preceq g'$;
\item[\textup{(iv)}] for all $f,g\in K^\times$ we have:
$$\begin{cases}
&f\prec 1,\ g\nasymp 1 \quad\Rightarrow\quad f'\prec g^\dagger,\\
&f\asymp 1,\ g\nasymp 1 \quad\Rightarrow\quad f'\preceq g^\dagger.
\end{cases}$$
\end{enumerate}
\end{lemma}

\noindent
If $K$ is asymptotic in $L$ and $v(L^\times)=v(K^\times)$, then we call $(\Gamma, \psi)$ as defined in 
(ii) of the lemma above the {\bf asymptotic couple of $K$ in $L$},  so
$$ \Psi\ :=\ \psi(\Gamma^{\ne})\ =\ \bigl\{v(g^{\dagger}): 1 \nasymp g\in K^{\times}\bigr\}.$$

\begin{proof}[Proof of Proposition~\ref{pa5}] 
Let $(\Gamma,\psi)$ be the
asymptotic couple of $K$ in $L$. We shall prove that $L$ is asymptotic
with asymptotic couple $(\Gamma,\psi)$. Let $g\in L^{\times}$, $g\not\asymp 1$; 
by Proposition~\ref{CharacterizationAsymptoticFields}, 
(i)~$\Leftrightarrow$~(ii), 
it suffices
to show that then 
$v(g^\dagger)=\psi(vg)$.

Before doing so,  note that  we may replace $L$ by a compositional conjugate: if ${\phi\in L^\times}$, then $K$ is
also asymptotic in $L^{\phi}$
with~${(\Gamma, \psi -v\phi)}$ as the asymptotic couple of~$K$ in $L^{\phi}$, and
with $\phi^{-1}g^\dagger$ as the
logarithmic derivative of~$g$ in $L^\phi$.

Take $f\in K^\times$ with $f\sim g$, so 
$g=f(1+\varepsilon)$ with $\varepsilon \in \smallo_L$.
Then $g^\dagger=f^\dagger + \frac{\varepsilon'}{1+\varepsilon}$.
Now $v(f^\dagger)=\psi(vf)=\psi(vg)\in \Psi$. Compositional conjugation by 
$\phi:= f^\dagger$  (and renaming) arranges that $v(f^\dagger)=0\in \Psi$, so
$(\Gamma, \psi)$ has small derivation. Then 
$\der \smallo\subseteq \smallo_L$
and~${\der\mathcal{O}\subseteq \mathcal{O}_L}$ by 
Lemma~\ref{pa1}(iv), so 
$\der\smallo_L\subseteq \smallo_L$ by Lemma~\ref{key}. This gives 
$f^\dagger \asymp 1\succ  \frac{\varepsilon'}{1+\varepsilon}$, so
$g^\dagger\sim f^\dagger$, and thus $v(g^\dagger)=v(f^\dagger)=\psi(vg)$. 
\end{proof}

\subsection*{Algebraic extensions}
{\em In this subsection we assume that $K$ is an 
asymptotic field
and $L|K$ is an extension of valued differential fields}.  We shall prove:

\begin{prop}\label{Algebraic-Extensions-Asymptotic}
If $L|K$ is algebraic, then $L$ is an asymptotic field.
\end{prop}

\begin{lemma}\label{trace and norm}
Suppose that $K$ is henselian and $L|K$ is of finite degree. 
Assume also that for every $f\preceq 1$ in $L$ there exists 
$a\preceq 1$ in $K$ such that
\begin{enumerate}
\item[\textup{(i)}] $f'\preceq a'$, and
\item[\textup{(ii)}] if $f'\asymp a'$ and $f\prec 1$, then $a\prec 1$.
\end{enumerate}
Then $L$ is an asymptotic field.
\end{lemma}
\begin{proof}
We verify that $L$ satisfies condition (iv) of 
Proposition~\ref{CharacterizationAsymptoticFields}. For this, 
let $f,g\in L^\times$ with $f\preceq 1$ and $g\nasymp 1$. Let $a\preceq 1$ in $K$
satisfy (i) and (ii). Put $b:=\operatorname{N}_{L|K}(g)\in K^\times$; 
then $b\nasymp 1$ and $b^\dagger\preceq g^\dagger$ by a result at the end of
Section~\ref{Valdifcon}. Hence if $f\asymp 1$, 
then by (i): 
$f'\preceq a'\preceq b^\dagger\preceq g^\dagger$. 
If $f \prec 1$, then we similarly get, using~(ii), that 
$f' \prec a'\preceq b^\dagger\preceq g^\dagger$, or $a\prec 1$ and 
$f'\asymp a'\prec b^\dagger\preceq g^\dagger$.
\end{proof}

\begin{lemma}\label{Residue-Field-Extension-Asymptotic}
Suppose that $K$ is henselian and $L|K$ is of finite degree, such that
$[L:K]=\bigl[\res(L):\res(K)\bigr]=n>1$. Then $L$ is an asymptotic field.
\end{lemma}
\begin{proof}
Take $y\asymp 1$ in $L$ such that $\res(L)=\res(K)[\bar{y}]$, where
$\bar{y}$ is the residue class of~$y$ in $\res(L)$. Let 
$a_0,\dots,a_{n-1}\preceq 1$ in $K$ be such that
$$P(Y)=Y^n+a_{n-1}Y^{n-1}+\cdots +a_1Y+a_0\in K[Y]$$ is the
minimum polynomial of $y$ over $K$, so its reduction $\bar{P}(Y)$ in
$\res(K)[Y]$ is the minimum polynomial of $\bar{y}$ over $\res(K)$; 
in particular $P'(y)\asymp 1$. Moreover, $P'(y)y'=-\sum_{i=0}^{n-1}a_i'y^i$, 
hence $y' \asymp a_{i_0}'$ with $i_0\in\{0,\dots,n-1\}$. Let 
$f\preceq 1$ in $L$. Then we have $f_0,\dots,f_{n-1}\preceq 1$ in $K$ such that
$$f=f_0+f_1y+\cdots+f_{n-1}y^{n-1}.$$
Then
$$f' = 
\sum_{i=0}^{n-1} f_i'y^i + \left(\sum_{j=1}^{n-1} jf_jy^{j-1}\right)y'$$
and hence 
$$v(f') \geq \gamma := \min\left\{\min_{0\leq i<n} v(f_i'), 
\min_{0<j<n} v(f_j)+v(a_{i_0}')\right\}.$$
We now define an element $a\in  \mathcal{O}$ as follows: If 
$\gamma=\min_{0\leq i<n} v(f_i')$, then $a:=f_{i_1}$ where 
$0\leq i_1<n$ and $v(f_{i_1}')=\gamma$;
if $\gamma<\min_{0\leq i<n} v(f_i')$, then $a:=a_{i_0}$. 
Then $a$ satisfies conditions (i) and (ii) of Lemma~\ref{trace and norm}, 
so $L$ is an asymptotic field.
\end{proof}

\noindent
The next lemma is a variant of Lemma~3.3 in \cite{AvdD2}:

\begin{lemma}\label{Lemma34}
Suppose $\res(K) = \res(L)$,
$T\supseteq K^\times$ is a subgroup of
$L^{\times}$ such that $L=K(T)$ \rom{(}as fields\rom{)}, each element of 
$K[T]\setminus \{0\}$ has the
form $t_1+\dots+t_k$ with $k\ge 1$, $t_1,\dots,t_k\in T$ and
$t_1 \succ t_i$ for $2\le i \le k$, and for all 
$a,b\in T$,
\begin{align}
  a,b\prec 1 
&\quad\Longrightarrow\quad a' \prec b^\dagger \label{Eq1} \\
  a\asymp 1,\ b\prec 1 
&\quad\Longrightarrow\quad a' \preceq b^\dagger. \label{Eq2}
\end{align}
Then $L$ is an asymptotic field.
\end{lemma}
\begin{proof}
Let $g\in L^\times$, $g\prec 1$. Then by the assumptions on $T$ we have
$$g=b\cdot\frac{1+\sum_{i=1}^m a_i}{1+\sum_{j=1}^n b_j}$$
where $b\in T$, and $a_i,b_j\in T$ with $a_i,b_j\prec 1$, 
for $1\leq i\leq m$, $1\leq j\leq n$. So $g\asymp b\prec 1$; we claim that 
$g'\asymp b'$. For this, note that $b^\dagger \ne 0$ by \eqref{Eq1}, so
$$g^\dagger = b^\dagger 
\left(1+\frac{\sum_{i=1}^m a_i'/b^\dagger}{1+\sum_{i=1}^m a_i}-
\frac{\sum_{j=1}^n b_j'/b^\dagger}{1+\sum_{j=1}^n b_j}\right).$$
By \eqref{Eq1} we have $a_i'/b^\dagger, b_j'/b^\dagger\prec 1$ for 
$1\leq i\leq m$, $1\leq j\leq n$. Therefore $g^\dagger \asymp b^\dagger$ 
and thus $g'\asymp b'$ as claimed.

Now if $f,g\in L^\times$ satisfy $f,g\prec 1$, then this claim and 
\eqref{Eq1} yield $f' \prec g^\dagger$. Suppose $f,g\in L^\times$ with
$f\asymp 1$ and $g\prec 1$. Take $b\in T$ with $b\asymp g$ and 
$b'\asymp g'$ as above. Since $\res(K)=\res(L)$ we have
$f=a+h$ with $a\in K$, $a\asymp 1$ and
$h\in L$, $h\prec 1$. By
the claim above we find $t\in T\cup\{0\}$ with $h\asymp t$ and $h'\asymp t'$. 
Now $a' \preceq b^\dagger$
by \eqref{Eq2} and $t'\prec b^\dagger$ by \eqref{Eq1}, so
$f'=a'+h'\preceq b^\dagger\asymp g^\dagger$. The implication 
(iv)~$\Rightarrow$~(i) of Proposition~\ref{CharacterizationAsymptoticFields} 
now yields that $L$ is an asymptotic field.
\end{proof}

\begin{lemma}\label{Purely ramified-Asymptotic}
Let $p$ be a prime number, and suppose that $L=K\bigl(u^{1/p}\bigr)$
where $u\in K^{\times}$ with $vu\notin p\Gamma$.
Then $L$ is an asymptotic field.
\end{lemma}
\begin{proof}
Let $u^{i/p} := (u^{1/p})^i$ for $i\in\Z$, and put 
$T:=\bigcup_{i=0}^{p-1} K^\times u^{i/p}$. Then $T$ is a multiplicative 
subgroup of $L^\times$. By Lemma~\ref{lem:pth root},  
$\res(K)=\res(L)$, and each element of $L$ has the form $t_1+\cdots+t_k$ with 
$t_1,\dots,t_k\in T$, $t_1\succ \cdots \succ t_k$. Now let $a,b\in T$; 
then $a^p\in K^\times$ and $(a^p)^\dagger=pa^\dagger\asymp a^\dagger$, 
and similarly $b^p\in K^\times$, $(b^p)^\dagger\asymp b^\dagger$. Hence 
if $a,b\prec 1$, then 
$v(a')=va+v\big((a^p)^\dagger\big)=\frac{1}{p}\alpha+\psi(\alpha)$ and 
$v(b^\dagger)=v\big((b^p)^\dagger)=\psi(\beta)$ where 
$\alpha:=pva,\beta:=pvb\in\Gamma$, and in the asymptotic couple 
$(\Q\Gamma,\psi)$ we have $\frac{1}{p}\alpha+\psi(\alpha)=
(\id+\psi)\big(\frac{1}{p}\alpha\big)>\psi(\beta)$, so $a' \prec b^\dagger$. 
If $a\asymp 1 \succ b$, then 
$a'\asymp (a^p)' \preceq (b^p)^\dagger\asymp b^\dagger$. 
Hence $L$ is an asymptotic field by Lemma~\ref{Lemma34}.
\end{proof}

\begin{proof}[Proof of Proposition~\ref{Algebraic-Extensions-Asymptotic}] 
We assume that $L|K$ is algebraic. We need to show that then $L$ is an
asymptotic field. The property of
being an asymptotic field is inherited by valued differential subfields, so
we can assume that $L$ is an algebraic closure of $K$. Next, by  
Lemma~\ref{pa5}, we can arrange that $K$ is henselian. We then reach $L$ in two steps. In the first step we pass from~$K$ to its maximal unramified extension
$K^{\operatorname{unr}}$ inside $L$. Then $K^{\operatorname{unr}}$ is asymptotic by 
Lemma~\ref{Residue-Field-Extension-Asymptotic}. In the second step we obtain~$L$ as a
purely ramified extension of~$K^{\operatorname{unr}}$ (Proposition~\ref{prop:ramified}), and now Lemma~\ref{Purely ramified-Asymptotic}  applies.  
\end{proof}

\noindent
%Lemma~\ref{PresInf-Lemma2},~(1), 
If $K$ has small derivation and $L|K$ is algebraic, then 
$L$ has small derivation, by Proposition~\ref{Algebraic-Extensions}. If $K$ is of $H$-type and 
$L|K$ is 
algebraic, then $L$ is of $H$-type. 

\subsection*{An application} {\em In this subsection $K$ is an ungrounded $H$-asymptotic field, and $\Gamma=v(K^\times)\ne \{0\}$}. For use in Chapter~\ref{evq} we extend slightly some results from 
Section~\ref{Applicationtodifferentialpolynomials} by dropping the assumption that $\Gamma$ is divisible:
 
%\noindent
%For use in Chapter~ \ref{evq} we also observe that some %results from 
%Section~\ref{Applicationtodifferentialpolynomials} go 
%through under a 
%weaker assumption than divisibility of $\Gamma=v(K^\times)$:

\begin{cor}\label{betterdifpol} 
%Assume $K$ is of $H$-type and 
%$\Gamma\ne \{0\}$ has no smallest positive element. 
Let $P\in K\{Y\}^{\ne}$
have order $r$. Then Corollary~\ref{ps6e}
goes through if we drop ``with the intermediate value property.''
Also Corollary~\ref{ps6c} goes through, for a unique tuple
$(\alpha, d, d_0, e_1,\dots, e_{r-1})$. In particular, there exists $\beta\in \Gamma^{<}$ such that $P(y)\ne 0$ for all $y\in K$ with $\beta < vy < 0$.
\end{cor}
\begin{proof} This follows by applying Corollaries~\ref{ps6e} and~\ref{ps6c} to the algebraic closure~$K^{\a}$ of $K$, taking into account
that the value group of $K^{\a}$ is $\Q\Gamma$, that 
$\Gamma^{<}$ is cofinal in $(\Q\Gamma)^{<}$, that $\Psi_{K^{\a}}=\Psi$, and that the $\chi$-map
of the $H$-asymptotic couple~$(\Q\Gamma, \psi)$ extends the
$\chi$-map of $(\Gamma, \psi)$. The ``uniqueness'' uses Lemma~\ref{evtcon}.
\end{proof} 

\noindent
Here is a striking consequence that will not
be needed, since in Chapter~\ref{ch:The Dominant Part and the Newton Polynomial} we prove more precise results for so-called $\upo$-free $K$ (using heavier machinery):

\begin{cor}\label{betterdifpolcor} Suppose $K$ is existentially closed in some grounded $H$-asymp\-to\-tic field extension of $K$. Then
Corollary~\ref{ps6c} holds with all $e_i=0$, for every~$P$ in~$K\{Y\}^{\ne}$ of order $r$. \rm{(See end of \ref{sec:ultra products} for \textit{existentially closed.}\/)}
\end{cor}
\begin{proof} The assumption on $K$ 
yields an elementary extension of $K$ with an element~$y\succ 1$ such that $K\<y\>$
has $\Cl(y)$ as its smallest comparability class. Let $P\in K\{Y\}^{\ne}$ have order $r$, let $\beta$ and $(\alpha, d, d_0, e_1,\dots, e_{r-1})$ be as in Corollary~\ref{ps6c}, and
suppose towards a contradiction that some $e_i\ne 0$. 
We have $\beta < vy <0$, and taking $a\in K^\times$ with
$va=\alpha$ we get an element $f=P(y)/ay^{d_0}(y^\dagger)^d$
of $K\<y\>^\times$ such that
$vf=-\sum_{i=1}^{r-1}e_i\chi^i(vy)$,
so $f\not\asymp 1$ and  $\Cl(f) < \Cl(y)$, a contradiction. 
\end{proof} 

\noindent
If $K$ is a directed union of grounded asymptotic subfields, then $K$ is existentially closed in some grounded $H$-asymptotic field extension of $K$; see~\ref{cor:ec in ext}. Note that
$K=\T$ and $K=\T_{\log}$ are such directed unions.

\subsection*{Notes and comments} The above material on algebraic 
extensions is a mild generalization of Section~3 in~\cite{AvdD2}.

%% file: mt-9n-6.tex
\section{Immediate Extensions of Asymptotic Fields}
\label{ImmExtAs}

\noindent  
We begin by recording analogues for asymptotic fields of the
results on maximal immediate extensions 
from Chapter~\ref{ch:valueddifferential}.
Next we prove
some easy technical facts that are often needed. We also show how to construct immediate extensions by (un)coarsening, and we finish
by constructing fluent completions.

\subsection*{Asymptotically maximal immediate extensions} 
{\em In this subsection $K$ is an asympotic field with small derivation and differential residue field $\k$ such that the derivation on $\k$ is nontrivial}.
This is the same condition as in Section~\ref{sec:cimex}. 
Recall that the derivation of any immediate asymptotic extension of
$K$ is small. 
Here are two easy consequences of Section~\ref{sec:cimex} and Lemmas~\ref{ndsupPsi=0} and~\ref{impreservesas}:

\begin{cor}\label{cor1asdiftrzda} If $K$ is asymptotically maximal, then $K$ is spherically complete.
\end{cor}

\begin{cor}\label{cor2asdiftrzda} $K$ has an immediate asymptotic 
extension that is spherically complete. 
\end{cor}

\subsection*{Easy technical lemmas}
{\em In this 
subsection $F$ is a valued differential field with a subfield $K$
that is asymptotic in $F$, and with $L$ as an intermediate field, so that
$K\subseteq L \subseteq F$.}\/ Note that then $L$ is asymptotic in $F$
if for each $f\in L^\times$ with $f\prec 1$ there is $a\in K^\times$ such that~$f\asymp a$ and~$f'\asymp a'$.

\begin{lemma}\label{pa3} 
Suppose $L|K$ is immediate and for all $a\in K^\times$ and
$f\in L$,
\begin{equation}\label{ImmediateEq}
\bigl(a\prec 1, f\prec
1\bigr)\quad\Longrightarrow\quad f' \prec a^\dagger.
\end{equation}
Then $L$ is asymptotic in $F$.
\end{lemma}
\begin{proof}
Let $a\in K^\times$, $h\in L^\times$ with $h\prec a\prec 1$; we claim
that $h'\prec a'$. To see this, set $g:= h/a\in L^\times$, so
$g \prec 1$.  Now apply
\eqref{ImmediateEq} with $f=g$ to get
$\frac{h}{a}-\frac{h'}{a'}=-\frac{g'a}{a'}\prec 1$; hence $h'/a'\prec 1$ and thus $h'\prec a'$.

Let now $f\in L^\times$, $f\prec 1$. Take $a\in K$
with $f\sim a$, so
$h=f-a \prec a\prec 1$. 
Then $h'\prec a'$, so $f'=a'+h'\asymp a'$ as required.
\end{proof}

\begin{lemma}\label{pa4}
Suppose $L|K$ is immediate, $U\supseteq K$ is a
$K$-linear subspace of $L$ with 
$L=\bigl\{u/w : u,w\in U, w\neq 0\bigr\}$, and for all
$a\in K^\times$ and  $u\in U$,
\begin{align}
\bigl(a\prec 1, u\prec 1\bigr)
\quad\Longrightarrow\quad u'\prec a^\dagger, \label{ImmediateEq1}\\
\bigl(a\prec 1, u\asymp 1\bigr)
\quad\Longrightarrow\quad u' \preceq a^\dagger. \label{ImmediateEq2}
\end{align}
Then $L$ is asymptotic in $F$.
\end{lemma}
\begin{proof}
We shall verify \eqref{ImmediateEq}. Let $f\in L$, $f\prec 1$. Then
$f=u/w$ with $u,w\in U$ and $w\neq 0$. After dividing $u$ and
$w$ by an element of $K$ asymptotic to $w$, we may assume $u\asymp f$
and $w\asymp 1$. We have $f'=u'/w -fw'/w$ with $u'/w \asymp u' \prec
a^\dagger$ and $fw'/w \asymp fw' \prec  a^\dagger$ for all
$a\in K^\times$ with $a\prec 1$, by \eqref{ImmediateEq1} and
\eqref{ImmediateEq2}, respectively; hence $f'\prec a^\dagger$ for all
such $a$. 
\end{proof}

\subsection*{Immediate extensions by coarsening} Let $K$ be an asymptotic field and $\Delta$ a convex subgroup
of $\Gamma=v(K^\times)$, giving rise to the coarsened valuation
$$\dot{v}\colon K^\times \to \dot{\Gamma}=\Gamma/\Delta$$ 
with residue field $\dot{K}=\dot{\cal{O}}/\dot{\smallo}$. 
Let $(L, \dot{v})$ be an immediate valued field
extension of~$(K, \dot{v})$, and denote its residue field by 
$\dot{L}$, so $\dot{L}=\dot{K}$ after the usual identification.
In Section~\ref{sec:decomposition} we extended the 
valuation $v$ on $K$ to a valuation
$v\colon L^\times \to \Gamma$
such that if $f\in L^\times$ and
$f=gu$ with $g\in K^\times$ and $u\in L^\times$, $\dot{v}(u)=0$, then
$v(f)=v(g)+v(\dot{u})$. By Lemma~\ref{immcoarse}, $L$ with $v$ is an immediate
valued field extension of the valued field~$K$, and its coarsening by
$\Delta$ is the $(L,\dot{v})$ we started with. 
Let $L$ also be equipped with a derivation $\der$ making 
$(L, \dot{v})$ an asymptotic extension of
$(K, \dot{v})$. Then:

\begin{lemma} 
\label{Un-coarsen} With the above valuation $v$ on $L$, the 
valued differential field 
$L$ is a $\Delta$-immediate asymptotic extension of $K$. 
\end{lemma}   
\begin{proof} Let $f\in L$ and $\dot{v}f>0$; we claim that 
then $v(f') > \Psi$.
To see this, note that $f=g(1+h)$ with $g\in K$, $h\in L$, and $\dot{v}h>0$,
since $(L, \dot{v})$ is an immediate extension of $(K, \dot{v})$. Then
$vh>0$, so $vg=vf> \Delta$. Now 
$$f'=g'(1+h)+ gh',\ \text{ with }\ v(g'(1+h))=vg' > \Psi.$$ 
Also, $0< \dot{v}f = \dot{v}g$ gives
$\dot{v}(f')=\dot{v}(g')= v(g')+\Delta$,
so $v(f')=v(g') + \delta$ with $\delta\in \Delta$, hence $vf' > \Psi + \delta$.
With $h$ instead of $f$ this gives $vh' > \Psi + \delta$ for some 
$\delta\in \Delta$, so $v(gh')> \Psi$ in view of $vg> \Delta$. Thus
$v(f') > \Psi$, as claimed.

Now $(L,v)$ is an immediate 
valued field extension of $K$ by Lemma~\ref{immcoarse}. Using this fact, 
the claim above, and Lemma~\ref{pa3} we show that $(L,v)$ is an
asymptotic field. 
Let $f\in L$ with $vf>0$; it is enough to show
that then $v(f') > \Psi$. As before, $f=g(1+h)$ with $g\in K$, $h\in L$, and
$\dot{v}h>0$. Then $vg>0$ and $f'=g'(1+h)+ gh'$ with $v(g'(1+h))=vg' > \Psi$.
By the claim above applied to $h$ instead of $f$ we have $v(h')>\Psi$,
so $v(gh') \ge v(h')> \Psi$. Thus $v(f') > \Psi$.
%, as desired.
\end{proof}

\subsection*{Fluency and $\Delta$-immediate extensions} 
Let $K$ be an asymptotic field, $\Delta$ a convex subgroup of $\Gamma$, and 
$\dot{v}=v_{\Delta}$ the coarsening of $v$ by $\Delta$. By 
Lemma~\ref{Un-coarsen} any proper immediate asymptotic extension of the asymptotic field $(K, \dot{v})$ yields by ``uncoarsening'' a proper
$\Delta$-immediate asymptotic extension of~$K$.

\begin{lemma}\label{flasy} Assume $\Delta\ne \{0\}$. Then the following are equivalent: \begin{enumerate}
\item[\textup{(i)}] $K$ has no proper $\Delta$-immediate asymptotic extension;
\item[\textup{(ii)}] every $\Delta$-fluent pc-sequence in $K$ pseudoconverges in $K$.
\end{enumerate}
\end{lemma}
\begin{proof} The direction (ii)~$\Rightarrow$~(i) follows from 
Corollary~\ref{maxdeltafluent}. As to (i)~$\Rightarrow$~(ii), we prove its 
contrapositive. Assume there exists a divergent $\Delta$-fluent pc-sequence in~$K$.
Take $\delta\in \Delta^{\ne}$ and $\phi\in K^\times$ with $v\phi=\psi(\delta)$. Then $\psi^{\phi}(\delta)=0$, so $0\in \psi^{\phi}(\Delta^{\ne}) \cap \Delta$. 
Replacing $K$ by its compositional conjugate $K^{\phi}$ 
we arrange $0\in \psi(\Delta^{\ne})\cap \Delta$, so $K$ has small derivation
and $\psi(\Delta^{\ne})\subseteq \Delta$. Then
$(K, \dot{v})$ has small derivation, and the derivation of the differential
residue field of $(K, \dot{v})$ is nontrivial. Hence $(K, \dot{v})$ has a proper immediate 
asymptotic extension by Corollary~\ref{cor1asdiftrzda}, and this yields a proper $\Delta$-immediate asymptotic extension of $K$ by Lemma~\ref{Un-coarsen}.
\end{proof}

\noindent
By Zorn there exists a $\Delta$-immediate asymptotic
extension of $K$ that has no proper $\Delta$-immediate asymptotic
extension. Such an extension of $K$ is called
a {\em maximal $\Delta$-immediate asymptotic extension\/} of $K$.
Let $K(\Delta)$ be a maximal $\Delta$-immediate asymptotic extension
of $K$. By Lemma~\ref{flasy}, if $\Delta\ne \{0\}$, then $K(\Delta)$ is also a maximal $\Delta$-immediate extension of $K$ as a valued field, and 
so its $\Delta$-coarsening~$(K, \dot{v})$ is henselian. Here is a weak ``asymptotic'' version of Proposition~\ref{fluentcompletion}:

\index{maximal!$\Delta$-immediate asymptotic extension}

\begin{prop}\label{asfl} Let $K$ be an asymptotic field such that $\Gamma\ne \{0\}$ and $[\Gamma^{\ne}]$ has no least element.
Then $K$ has an immediate asymptotic extension which, as a valued field extension of $K$, is a fluent completion of $K$.
\end{prop}
\begin{proof} We follow the construction in the proof of Proposition~\ref{fluentcompletion}. Fix a decreasing coinitial sequence $(\Delta_{\alpha})$
of nontrivial convex subgroups of~$\Gamma$ indexed by the ordinals~$\alpha< \lambda$ for some infinite limit ordinal $\lambda$. For each $\alpha$ we pick a maximal $\Delta_{\alpha}$-immediate
asymptotic extension $K(\Delta_{\alpha})$ of $K$. We arrange this so that
$K(\Delta_{\alpha})$ is a valued differential subfield of 
$K(\Delta_{\beta})$ whenever
$\alpha < \beta < \lambda$. Put 
$$  K^{\operatorname{f}}\ :=\  \bigcup_{\alpha< \lambda} K(\Delta_{\alpha}).$$
Then  $K^{\operatorname{f}}$ is an asymptotic extension of $K$, and $K^{\operatorname{f}}$ as a valued field extension of~$K$ is a semifluent completion of $K$ as defined in the proof of Proposition~\ref{fluentcompletion}. Iterating as in that proof the construction
$K  \leadsto K^{\operatorname{f}}$ we eventually arrive at an asymptotic extension of $K$ that is also a fluent completion of $K$.      
\end{proof}

%% file: mt-9n-7.tex
\section{Differential Polynomials of Order One}\label{caseorderone}

\noindent
{\em Throughout this section
$K$ is an $H$-asymptotic field with asymptotic couple
$(\Gamma, \psi)$, $\Gamma\ne \{0\}$. We let $\gamma$ range over $\Gamma$ and $y$ over $K^\times$, and assume 
$P\in K\{Y\}^{\ne}$ has order $\le 1$.}\/ For such $P$ we improve on what 
Sections~\ref{sec:cimex} and \ref{Applicationtodifferentialpolynomials} 
yield for differential polynomials of 
arbitrary order. 
%The results here are not used later in this volume. 

\subsection*{Behavior of $v\big(P(y)\big)$} The goal of this
subsection is to show the following:

\begin{prop}\label{prop:vP(y), order 1}
Suppose $P(0)=0$. Then there is a finite set $\Delta(P)\subseteq\Gamma$ and a finite partition of $\Gamma\setminus\Delta(P)$ into  
convex subsets of $\Gamma$ such that for each set $U$ in 
this partition we have a strictly increasing function $i_U\colon U\to\Gamma$ for which
$$v\big(P(y)\big)\ =\ i_U(vy)\  \text{ whenever $vy \in U$.}$$
\end{prop}

\index{map!slowly varying}
\index{slowly varying}

{\sloppy
\noindent
We derive this from a more precise result for homogeneous $P$.
We start with some generalities about functions on ordered abelian groups. Let $G$ be an ordered abelian group and $U\subseteq G$. 
A function $\eta\colon U \to G$ is said to be {\bf slowly varying} if
${\eta(\alpha)-\eta(\beta)}=o(\alpha - \beta)$
for all distinct $\alpha, \beta \in U$. Note that then the fun\-ction~${\alpha \mapsto  \alpha + \eta(\alpha)\colon U \to G}$ is strictly increasing.
A key example of a slowly varying function is of course~${\psi\colon \Gamma^{\ne} \to \Gamma}$. Note that each constant function $U\to G$ is
slowly varying, and that if $\eta_1, \eta_2\colon U\to G$ are slowly 
varying, so
are $\eta_1 + \eta_2, \eta_1 - \eta_2,$ and
$$\alpha \mapsto \min\!\big(\eta_1(\alpha), \eta_2(\alpha)\big)\colon\ U\to G.$$

\begin{lemma} \label{slowlyvarying}
Let $s\in K$ be given. Then there is a $\gamma_0\in \Gamma$ and
a slowly varying function 
$\eta\colon \Gamma \setminus \{\gamma_0\} \to \Gamma$
with the following properties: \begin{enumerate}
\item[\textup{(i)}] $v(y^\dagger-s)=\eta(vy)$ for all $y$ with $vy\ne \gamma_0$;
\item[\textup{(ii)}] for each $\alpha\in \Gamma$ the set
$\big\{\gamma:\ \gamma\ne \gamma_0,\ \eta(\gamma)\le \alpha\big\}$
is a union of finitely many dis\-joint convex subsets of $\Gamma$.
\end{enumerate}
\end{lemma}

}

\begin{proof} Suppose first that $s=a^\dagger$ with $a\in K^\times$.
Then $y^\dagger - s = (y/a)^\dagger$, so $v(y^\dagger - s) = \psi(vy - va)$
if $vy \ne va$. In this case we can take $\gamma_0=va$ and
$\eta(\gamma)= \psi(\gamma - va)$. 

Next assume that $s\ne y^\dagger$ for all $y$. Then we take a nonzero 
$\phi$ in an elementary extension $L$ of $K$ such that
$v(y^\dagger -s)\le v(\phi^\dagger -s)$ for all $y$. (This $\phi$ could be
an element of $K$.) Let $\gamma_1:=v(\phi^\dagger -s)\in \Gamma_L$.  
Then we claim that for all $y$ with $vy\ne v\phi$:
$$ v(y^\dagger - s)\ =\ \begin{cases} \psi_L(vy - v\phi) 
& \text{ if } \psi_L(vy - v\phi)\le \gamma_1\\
\gamma_1 & \text{ if }\psi_L(vy - v\phi)\ge \gamma_1. \end{cases}
$$
To see this, let $vy \ne v\phi$. From 
$v(y^\dagger - s) \le v(\phi^\dagger -s)$ we get
$y^\dagger - \phi^\dagger \not\sim s-\phi^\dagger$. Since
$$y^\dagger - s\ =\ (y^\dagger - \phi^\dagger) - (s-\phi^\dagger)\ \text{ and }\
y^\dagger - \phi^\dagger\  =\ (y/\phi)^\dagger,$$ 
this gives the claim. Note that the claim can also be expressed as:
$$v(y^\dagger-s) = \min\big\{\psi_L(vy-v\phi),\gamma_1\big\} 
\text{ whenever }vy\ne v\phi.$$
If $v\phi\in \Gamma$, then the lemma clearly
holds with $\gamma_0:= v\phi$, and if $v\phi\notin \Gamma$, then
the lemma holds for any $\gamma_0\in \Gamma$. 
\end{proof}

\begin{lemma}\label{lem:vP(y), order 1, homog}
Suppose $P$ is homogeneous of degree $d$. There is a finite set $\Delta(P)\subseteq\Gamma$ and a finite partition of $\Gamma\setminus\Delta(P)$ into  
convex subsets of $\Gamma$ such that for each set $U$ in this partition
we have a slowly varying function $\eta_U\colon U\to\Gamma$ for which
$$v\big(P(y)\big)\ =\ d\,vy+\eta_U(vy)\  \text{ whenever $vy \in U$.}$$
\end{lemma}

\begin{proof}
We have $P=a_0Y^d + a_1Y^{d-1}Y' + \dots +a_d(Y')^d$ with  $a_0,\dots,a_d\in K$,
so $P(y)=y^d(a_0 + a_1z + \dots + a_dz^d)$, where 
$z=y^\dagger$. The henselization of $K$ is still $H$-asymptotic with the same asymptotic couple, so we can
assume $K$ is henselian. Now apply Lemma~\ref{piecewise uniform description} and the surrounding 
remarks to the polynomial 
$a_0+ a_1Z + \dots + a_dZ^d$ to partition $K$ into 
sets $G_1,\dots,G_k$ such that for $i=1,\dots,k$: 
$G_i$ is a special disk in $K$ with holes, and we have  
$b_{i0},\dots, b_{id}, s_i\in K$ with
$$v\big(P(y)\big)\ =\ 
d\,vy + \min\big\{v(b_{ij}) + j\,v(z - s_i):\  0\le j\le d\big\}\qquad\text{if $z\in G_i$.}$$ 
In addition, for $i=1,\dots,k$,  either $s_i=0$ and the condition
$z\in G_i$ is of the form 
$$v(z- s_{i,1})\ \le\ v(s_{i,1}),\ \dots,\ 
v(z - s_{i,n(i)})\ \le\  v(s_{i,n(i)}),$$
with $s_{i,1},\dots,s_{i,n(i)}\in K^\times$,
or $s_i\ne 0$ and the condition $z\in G_i$ is of the form
$$v(z - s_i)\ >\  v(t_i),\  v(z- s_{i,1})\ \le\ v(t_{i,1}),\
\dots,\ v(z - s_{i,n(i)})\ \le\  v(t_{i,n(i)})$$
with $t_i,s_{i,1}, t_{i,1},\dots,s_{i,n(i)}, t_{i,n(i)}\in K^\times$ and $s_i\succeq t_i,\ s_{i,1}\succeq t_{i,1},\dots, s_{i,n(i)}\succeq t_{i,n(i)}$.
Now use Lemma~\ref{slowlyvarying}.
\end{proof}

\begin{proof}[Proof of Proposition~\ref{prop:vP(y), order 1}]
We have
$$P(Y)\ =\ P_1(Y) + \dots + P_n(Y)\ 
\text{ (decomposition into homogeneous parts)}.$$ 
Let
$d$ below range over the numbers $1,\dots,n$ for which $P_d\ne 0$; 
likewise with $e$. Applying Lemma~\ref{lem:vP(y), order 1, homog} to all $P_d$ simultaneously we
obtain a finite $\Delta(P)\subseteq\Gamma$ and a finite partition of~$\Gamma\setminus\Delta(P)$
into convex subsets of $\Gamma$
such that for each set $U$ in the partition and each $d$ we have a slowly varying function $\eta_{U,d}\colon U\to\Gamma$ 
for which
$$v\big(P_d(y)\big)\ =\ d\,vy+\eta_{U,d}(vy)\ \text{ whenever $vy\in U$.}$$
If $d<e$, then the function 
$$\gamma \mapsto 
\big(e\,\gamma+ \eta_{U,e}(\gamma)\big) -
\big(d\,\gamma+ \eta_{U,d}(\gamma)\big)\colon\ U \to \Gamma$$ is strictly increasing. Hence, after increasing $\Delta(P)$ and refining our 
partition if necessary, we can arrange that for each $U$ in the partition 
there is a $d=d_U$ such that for all~${e\ne d}$:
$$v\big(P_d(y)\big)< v\big(P_e(y)\big)  \text{ whenever $vy\in U$.}$$
Thus $v\big(P(y)\big) =  v\big(P_d(y)\big)$ for $d=d_U$ and $vy\in U$, and the proposition follows.
\end{proof}

\subsection*{Evaluation at pc-sequences} {\em In this subsection we fix a pc-sequence $(a_{\rho})$ in $K$}. 
Proposition~\ref{pca-extension}
below is an analogue of Kaplansky's theorem (Lemma~\ref{lem:Kaplansky, 2}) about pc-sequences of algebraic type.
It is stronger than Lemma~\ref{zda} in not requiring the derivation of 
$K$ to be small with nontrivial differential
residue field, but weaker in assuming that $K$ is $H$-asymptotic and $P$ has order $\le 1$. Moreover, we never need to replace in our situation $(a_{\rho})$ by an equivalent pc-sequence. 
%In addition to $P\in K\{Y\}^{\ne}$ being of order 
%$\le 1$ we also assume below that $P\notin K$. 

\begin{lemma}\label{psa} Assume $P\notin K$ and 
$a_\rho\leadsto a\in K$. Then $P(a_\rho) \leadsto P(a)$.
\end{lemma}

\begin{proof} 
Replacing $a_\rho$ by $a_\rho -a$ and $P(Y)$ by $P_{+a}(Y)-P(a)$ we can  assume
that $a_\rho \leadsto 0$, $P(0)=0$, and have to show that then
$P(a_\rho) \leadsto 0$. Proposition~\ref{prop:vP(y), order 1} gives a finite subset $\Delta(P)$ of $\Gamma$ and a finite partition of $\Gamma\setminus\Delta(P)$ into finitely
many convex subsets of $\Gamma$ such that for each set~$U$ in this partition,
$v\big(P(y)\big)$ is strictly increasing as a function of $vy\in U$. 
Taking $U$ in the partition such that $v(a_\rho)\in U$  eventually, 
we conclude that $P(a_\rho) \leadsto 0$.
\end{proof}

\begin{cor}\label{pca} If $P\notin K$, then $\big(P(a_\rho)\big)$ 
is a pc-sequence.
\end{cor}
\begin{proof} Use that  $(a_\rho)$ has a pseudolimit in
some elementary extension of $K$. 
\end{proof}

\begin{prop} \label{pca-extension} Assume $P\notin K[Y]$, 
$P(a_{\rho})\leadsto 0$, and $Q(a_{\rho})\not\leadsto 0$ for
every $Q\in K\{Y\}$ with $\operatorname{c}(Q) < \operatorname{c}(P)$.
Then $K$ has an immediate asymptotic extension
field~$K\langle a \rangle$ with the following properties: \begin{enumerate}
\item[\textup{(i)}] $P(a)=0$ and $a_\rho \leadsto a$;
\item[\textup{(ii)}] for any $H$-asymptotic extension field $L$ of $K$ 
and any $b\in L$ with $P(b)=0$ and $a_\rho \leadsto b$ there is a 
unique embedding
$K\langle a\rangle \to L$ over $K$ that sends $a$ to $b$. 
\end{enumerate}
\end{prop}
\begin{proof} Note first that $P\in K[Y,Y']$ is irreducible.  
Consider the domain 
$$K[y_0, y_1]\ =\ K[Y,Y']/(P), \quad y_0:= Y+(P), \quad y_1:= Y'+(P),$$ and let $K(y_0,y_1)$ be its field of fractions. 
We extend the valuation~$v$ on $K$ to a valuation 
$v\colon K(y_0,y_1)^\times \to \Gamma$ as
follows. Pick a pseudolimit $e$ of $(a_\rho)$ in some $H$-asymptotic field
extension of $K$. Let $\phi\in K(y_0,y_1)$, $\phi\ne 0$, so 
$\phi=f(y_0, y_1)/g(y_0)$ with $f\in K[Y,Y']$ of lower degree in $Y'$ 
than $P$ and $g\in K[Y]^{\neq}$. 

\claim{$v\big(f(e,e')\big), v\big(g(e)\big)\in \Gamma$, and
$v\big(f(e,e')\big)- v\big(g(e)\big)$ depends only on $\phi$ and not on 
the choice of $(f,g)$.}

\noindent
To see why this claim is true, note that 
$f(a_\rho, a_\rho') \not\leadsto 0$ by the minimality of $P$, 
and that $f(a_\rho, a_\rho') \leadsto f(e,e')$ if $f\notin K$, 
by Lemma~\ref{psa}. Hence
$$f(a_{\rho}, a_{\rho}') \sim f(e,e')\ \text{ and }\
v\big(f(a_{\rho}, a_{\rho}')\big)= v\big(f(e,e')\big), 
\quad \text{eventually}.$$
In particular, $v\big(f(e,e')\big) \in \Gamma$, and likewise, 
$v\big(g(e)\big)\in \Gamma$. 
Suppose that also 
$\phi=f_1(y_0, y_1)/g_1(y_0)$ with $f_1\in K[Y,Y']$ of lower degree in $Y'$ 
than $P$ and $g_1\in K[Y]^{\neq}$. Then $fg_1 \equiv f_1g \bmod P$ in 
$K[Y,Y']$, and thus $fg_1=f_1g$ since $fg_1$ and $f_1g$ have lower degree in 
$Y'$ than $P$. Hence 
$v\big(f(e,e')\big)- v\big(g(e)\big)=v\big(f_1(e,e')\big)- v\big(g_1(e)\big)$,
thus establishing the claim.

\medskip\noindent
This allows us to define $v\colon K(y_0,y_1)^\times \to \Gamma$ by
$$v\phi\ :=\  v\big(f(e,e')\big)- v\big(g(e)\big).$$
It is routine to check that this map $v$ is a valuation on the field 
$K(y_0,y_1)$, except maybe for the multiplicative law. 
For this, for $i=1,2$ let $\phi_i\in
K(y_0,y_1)^\times$, so $\phi_i=f_i(y_0,y_1)/g_i(y_0)$ with $f_i\in K[Y,Y']^{\neq}$ 
of lower degree in~$Y'$ than $P$ and $g_i\in K[Y]^{\neq}$. Then 
$$v(\phi_i)\ =\ v\big(f_i(a_\rho,a_\rho')\big)-v\big(g_i(a_\rho)\big)
\quad\text{eventually ($i=1,2$)}$$
and hence
$$v(\phi_1)+v(\phi_2)\ =\ v\big((f_1f_2)(a_\rho,a_\rho')\big)-
 v\big((g_1g_2)(a_\rho,a_\rho')\big)
 \quad\text{eventually.}$$
We have $f_1f_2=(qP)/h\ +\ r/h$
where $q,r\in K[Y,Y']$, $h\in K[Y]^{\neq}$ and $\deg_{Y'} r<\deg_{Y'} P$.
Then 
$\phi_1\phi_2=r(y_0,y_1)/s(y_0)$ where $s:=g_1g_2h\in K[Y]^{\neq}$, 
and hence $$v(\phi_1\phi_2)\ =\
v\big(r(a_\rho,a_\rho')\big)-v\big(s(a_\rho)\big)
\qquad\text{eventually.}$$
Now $(f_1f_2)(a_\rho,a_\rho')\not\leadsto 0$ since
$f_1(a_\rho,a_\rho')\not\leadsto 0$ and $f_2(a_\rho,a_\rho')\not\leadsto 0$. 
Moreover, $P(a_\rho,a_\rho')\leadsto 0$ 
and $r(a_\rho,a_\rho')\not\leadsto 0$; hence
$$v\big((f_1f_2)(a_\rho,a_\rho')\big)\ =\ 
v\big(r(a_\rho,a_\rho')\big)-v\big(h(a_\rho)\big)
\quad\text{eventually.}$$
Therefore $v(\phi_1\phi_2)=v(\phi_1)+v(\phi_2)$ as desired.
Obviously, the value group of the valuation $v$ on $K(y_0,y_1)$ is $\Gamma$.
Its residue field is that of $K$. To see this, let $\phi\in K(y_0,y_1)^\times$
and $v\phi=0$; we shall find $s\in K$ with $v(\phi -s)>0$. We have
$\phi=f(y_0,y_1)/g(y_0)$ with $f\in K[Y,Y']$ of lower
degree in $Y'$ than $P$ and $g\in K[Y]^{\neq}$. Multiplying
$f$ and $g$ by a suitable element of $K^\times$ we may
assume that $v\big(f(e,e')\big)=v\big(g(e)\big)=0$. By the above we have
$f(e,e')\sim f(a_\rho, a_\rho')$ and $g(e)\sim g(a_\rho)$, eventually, hence
$v\left(\frac{f(e,e')}{g(e)} - s\right) >0$ where 
$s:=\frac{f(a_\rho, a_\rho')}{g(a_\rho)}$ with large enough $\rho$,
so $v(\phi -s)>0$ for such $s$. 

\medskip\noindent
We now equip $K(y_0,y_1)$ with the derivation 
extending the derivation of $K$ such that $y_0'=y_1$. We also set $a:= y_0$, 
so $a'=y_1$ and $K(y_0, y_1)=K\langle a \rangle$. Then $P(a)=0$, trivially, 
and $a_\rho \leadsto a$, as is easily checked.

\medskip\noindent
To show that $K\langle a \rangle$ is an asymptotic field, we use 
Proposition~\ref{pa5} with $K(a)$ in the role of $K$ and 
$L= K\langle a \rangle$. By that proposition it is enough to show
that  $K(a)$ is asymptotic in $K\langle a \rangle$. To do that we apply Lemma~\ref{pa4} with $L=K(a)$, $F= K\langle a \rangle$ 
and $U=K[a]$.
Consider elements $u=g(a)$  with $g\in K[Y]\setminus K$, and 
$b\in K^\times$ such that $b\prec 1$; it suffices to show that then
$u' \prec b^\dagger$ if $u\prec 1$, and
$u' \preceq b^\dagger$ if $u\asymp 1$. If $u\asymp 1$, then we take
$s\in K$ with $u \sim s$ and use $s'\preceq b^\dagger$
to reduce to the case $u \prec 1$. So we assume $u\prec 1$.
Note that 
$v(u)= \text{eventual value of } v\big(g(a_\rho)\big)$,
so $v\big(g(a_\rho)\big)>0$ eventually, hence
$v\big( g(a_\rho)' \big)$ is eventually constant as a function of $\rho$, 
and $v\big( g(a_\rho)' \big)>v(b^\dagger)$ eventually.
Therefore it is enough to show:

\claim{$v(u')= \text{eventual value of $v\big(g(a_\rho)'\big)$}$.}

\noindent
Let $g=c_0 + c_1Y + \dots +c_nY^n$ ($c_i\in K$), and put
$g^{\der}:=c_0' + c_1'Y + \dots +c_n'Y^n$, so
$$ u'\ =\  g(a)'\ =\  g^{\der}(a) + \frac{\partial g}{\partial Y}(a)a', \qquad 
g(a_\rho)'\ =\  g^{\der}(a_\rho) +  \frac{\partial g}{\partial Y}(a_\rho)a_\rho'.$$
Therefore the claim holds if $P(Y,Y')$ is of degree $>1$ in $Y'$, so we can
assume 
$$P(Y,Y')\ =\ P_0(Y)+P_1(Y)Y'\qquad\text{where $P_0,P_1\in K[Y]$, $P_1\neq 0$.}$$
Put $R:=-P_0/P_1$. Then $P_1(a)\neq 0$ and $a'=R(a)$. 
Also $P_1(a_\rho)\neq 0$ eventually. We may assume that 
$P_1(a_\rho)\neq 0$ for all $\rho$. Then for each $\rho$ we have
\begin{equation}\label{formula for g(a_rho)'}
g(a_\rho)'\ =\  g^{\der}(a_\rho) + \frac{\partial g}{\partial Y}(a_\rho)R(a_\rho) + 
\frac{\partial g}{\partial Y}(a_\rho)\big(a_\rho'-R(a_\rho)\big).
\end{equation}
Now $$v\big( g^{\der}(a_\rho) + \textstyle\frac{\partial g}{\partial Y}(a_\rho)R(a_\rho) \big)\  =\  
v\big( g^{\der}(a) + \textstyle\frac{\partial g}{\partial Y}(a)R(a) \big)$$ eventually. 
Also
$a_\rho'-R(a_\rho)=P(a_\rho,a_\rho')/P_1(a_\rho)$, hence 
$v\big(\textstyle\frac{\partial g}{\partial Y}(a_\rho)\big(a_\rho'-R(a_\rho)\big)\big)$ is eventually 
strictly increasing. We have
$$v\big(\textstyle\frac{\partial g}{\partial Y}(a_\rho)\big(a_\rho'-R(a_\rho)\big)\big)\ >\ 
v\big( g^{\der}(a) + \textstyle\frac{\partial g}{\partial Y}(a)R(a) \big)$$
eventually: otherwise,
$$v\big(\textstyle\frac{\partial g}{\partial Y}(a_\rho)\big(a_\rho'-R(a_\rho)\big)\big)\ <\ 
v\big( g^{\der}(a) + \textstyle\frac{\partial g}{\partial Y}(a)R(a) \big)$$
eventually, and then, by \eqref{formula for g(a_rho)'},
$v\big(g(a_\rho)'\big)$ would be both eventually constant and eventually
strictly increasing. Hence
$$v\big(g(a_\rho)'\big)\ =\   v\big( g^{\der}(a) + 
\textstyle\frac{\partial g}{\partial Y}(a)R(a) \big)\ =\  v(u')$$
eventually, proving the claim. Thus $K\langle a \rangle$ is an 
asymptotic field.

\medskip\noindent
It remains to prove item (ii). But this follows easily from the above,
since any $b$ as in~(ii) is transcendental over $K$, and can serve
as $e$ in the arguments above. 
\end{proof}

\subsection*{Notes and comments} Suppose $K$ has small derivation and~$P$ (of order $\le 1$) is homogeneous
of degree $d\geq 1$, with $d\Gamma=\Gamma$. Then the function
$v_P\colon\Gamma\to\Gamma$ is a strictly
increasing bijection, by the Equalizer Theorem. For this case, however, we have a more constructive (unpublished) argument leading to the stronger result that 
the function 
$\gamma\mapsto \textstyle\frac{1}{d}v_P(\gamma)\colon\Gamma\to\Gamma$ is
$\psi$-steady. \marginpar{need to check this is correct}

%% file: mt-9n-8.tex
\section{Extending $H$-Asymptotic Couples}\label{sec:extH} 

\noindent
Let $(\Gamma, \psi)$ and $(\Gamma_1, \psi_1)$ be asymptotic couples. An
{\bf embedding} $$h \colon\ (\Gamma,\psi)
\to (\Gamma_1,\psi_1)$$ is an embedding $h \colon\Gamma\to\Gamma_1$
of ordered abelian groups
such that $$ h\bigl(\psi(\gamma)\bigr)\ =\ \psi_1\bigl(h(\gamma)\bigr)\ 
\text{ for $\gamma\in\Gamma^{\ne}$.}$$
If $\Gamma\subseteq \Gamma_1$ and the inclusion $\Gamma\hookrightarrow\Gamma_1$ is 
an embedding $(\Gamma,\psi)
\to (\Gamma_1,\psi_1)$, then we call~$(\Gamma_1,\psi_1)$ an {\bf extension} of $(\Gamma,\psi)$.

\medskip\noindent
{\em In the rest of this section we assume
$(\Gamma,\psi)$ to be an $H$-asymptotic couple, and $(\Gamma_1, \psi_1)$ to be 
an asymptotic couple, not necessarily of $H$-type}. Thus
$\psi$ is constant on every archimedean class of $\Gamma$: for 
$\alpha,\beta\in\Gamma^{\ne}$ with $[\alpha]=[\beta]$ we have $\psi(\alpha)=
\psi(\beta)$. 

 Most of the extension results in this section come from 
\cite[Section 2]{AvdD2}, but the restriction
to $H$-asymptotic $(\Gamma,\psi)$ leads to fewer case distinctions.

\index{asymptotic couple!embedding}
\index{asymptotic couple!extension}
\index{embedding!asymptotic couples}
\index{extension!asymptotic couples}

\begin{lemma}\label{extension4}
Let $i\colon\Gamma\to G$ be an embedding of ordered abelian groups inducing
a bijection $[\Gamma]\to [G]$. Then there is a
unique function $\psi_G\colon G^{\ne} \to G$ such that~$(G,\psi_G)$
is an $H$-asymptotic couple and $i\colon (\Gamma, \psi) \to (G,\psi_G)$ is an
embedding. 
%and $i\bigl(\psi(\gamma)\bigr)=\psi_1\bigl(i(\gamma)\bigr)$ for all 
%$\gamma\in\Gamma^{\ne}$. 
\end{lemma}
\begin{proof} Define $\psi_G(g):= i\bigl(\psi(\gamma)\bigr)$ for $g\in G^{\ne}$ and
$\gamma\in \Gamma^{\ne}$ with $[g]=\big[i(\gamma)\big]$. Then $\psi_G\colon G^{\ne} \to G$ has the required properties: to check (AC3), pass to the divisible hulls of $\Gamma$ and $G$. 
\end{proof}

\noindent
The next
lemma and its proof show how to remove a gap. 

\begin{lemma}\label{extension1}
Let $\beta$ be a gap in $(\Gamma, \psi)$. Then there is an 
$H$-asymptotic couple ${(\Gamma+\Z \alpha,
\psi^{\alpha})}$ extending $(\Gamma,\psi)$ such that \begin{enumerate}
\item[\textup{(i)}] $\alpha>0$ and $\alpha'=\beta$;
\item[\textup{(ii)}] if $i\colon (\Gamma,\psi)\to (\Gamma_1,\psi_1)$ is an embedding and $\alpha_1\in\Gamma_1$, $\alpha_1>0$, $\alpha_1'=i(\beta)$, then
$i$ extends uniquely to an embedding $j\colon\bigl(\Gamma+\Z \alpha,
\psi^{\alpha}\bigr)\to (\Gamma_1,\psi_1)$ with $j(\alpha)=\alpha_1$.
\end{enumerate}
\end{lemma}
\begin{proof} Suppose $(\Gamma+\Z \alpha,
\psi^{\alpha})$ is an asymptotic couple that extends $(\Gamma,\psi)$ and
satisfies~(i). At this point we do not assume $(\Gamma+\Z \alpha, \psi^{\alpha})$ to be of $H$-type. Then $\alpha'< (\Gamma^{>})'$ gives 
$0 <\alpha< \Gamma^{>}$. Since $\Psi$ has no largest element, 
$[\Gamma^{\ne}]$ has no least element, so $0 < n\alpha < \Gamma^{>}$ 
for all~$n\ge 1$, in particular,
$\Gamma + \Z\alpha=\Gamma\oplus \Z\alpha$. Hence
$$\psi^{\alpha}(\alpha)\ =\ \alpha'-\alpha\ =\ \beta-\alpha\ >\ \Psi,$$ 
and thus for all $\gamma\in \Gamma$ and $k\in \Z$ with $\gamma+k\alpha\ne 0$, 
\begin{equation}\label{eq:psialpha}
\psi^{\alpha}(\gamma+k \alpha)\ =\ \begin{cases}
\psi(\gamma), &\text{if $\gamma\neq 0$,}\\
\beta-\alpha, &\text{otherwise.}
\end{cases} 
\end{equation}
It easily follows that $(\Gamma+\Z\alpha, \psi^{\alpha})$ is of $H$-type,
and has the universal property~(ii). All this assumes the existence
of an asymptotic couple $(\Gamma+\Z\alpha, \psi^{\alpha})$ that extends~$(\Gamma, \psi)$
and satisfies (i). 

To get such a couple, take an ordered abelian group extension 
$\Gamma^{\alpha}=\Gamma+\Z \alpha$ of~$\Gamma$ such that
$0 < n\alpha < \Gamma^{>}$ for all~$n\ge 1$ and extend $\psi$ to 
$\psi^\alpha\colon (\Gamma^{\alpha})^{\ne} \to \Gamma^{\alpha}$ according to \eqref{eq:psialpha}; in particular $\alpha+ \psi^{\alpha}(\alpha)=\beta$.  
It remains to show that $(\Gamma+\Z \alpha,
\psi^{\alpha})$ is an asymptotic couple. It is tedious but routine to check that (AC1) is satisfied, and (AC2) holds trivially. As to (AC3), note first that 
$$\max\{\psi^\alpha(\gamma+k\alpha):\ \gamma\in \Gamma,\ k\in \Z,\ \gamma+k\alpha\ne 0\}\ =\ \beta-\alpha.$$ 
Thus, given $\gamma\in \Gamma$ and $k\in \Z$ with 
$\gamma+k\alpha>0$, we need only verify that 
$$\beta-\alpha\ <\ \psi^\alpha(\gamma+k \alpha)+(\gamma+k \alpha).$$ 
This is easy to do by distinguishing the cases $\gamma\ne 0$ and $\gamma=0$. 
% let $\gamma+r \alpha, \delta+s \alpha\in\Gamma^{\ne}$ ($\gamma,\delta\in\Gamma$,
%$r,s\in\Z$) with $\delta+s \alpha>0$ (hence $\delta\ge 0$); we have to
%show $\psi^\alpha(\gamma+r \alpha)<\psi^\alpha(\delta+s \alpha)+(\delta+s \alpha)$. The case  
%$\psi^\alpha(\gamma+r \alpha)\leq\psi^\alpha(\delta+s \alpha)$
%is obvious, so assume $\psi^\alpha(\gamma+r \alpha)>
%\psi^\alpha(\delta+s \alpha)$. Then $\delta\ne 0$ by \eqref{eq:psialpha}, so $\delta>0$. 
%If moreover $\gamma=0$, then
%$[0]<[\alpha]<[\Gamma^{\ne}]$ and $\beta<(\Gamma^{>})'$
%give $\beta-(s+1)\alpha<\delta+\psi(\delta)$, that is,
%$$\psi^\alpha(\gamma+r \alpha)-\psi^\alpha(\delta+s \alpha)\ =\ \beta-\alpha-\psi(\delta)\ <\ \delta+s \alpha,$$
%as required. Similarly, if $\gamma\neq 0$, we get
%$$\psi^\alpha(\gamma+r \alpha)\ =\ \psi(\gamma)\ <\ \psi(\delta)+(\delta+s \alpha),$$ by (AC3) for
%$(\Gamma,\psi)$, and $[0]<[\alpha]<[\Gamma^{\ne}]$.
\end{proof}
  
%For details we refer to the proof of Lemma~2.10 in~\cite{AvdD2}, which 
%is just our lemma for arbitrary asymptotic couples instead 
%of $H$-asymptotic couples. 

\noindent
The universal property (ii) determines $(\Gamma+\Z \alpha, \psi^{\alpha})$ 
up to isomorphism over $(\Gamma, \psi)$. Note also that 
$[\Gamma+\Z\alpha]=[\Gamma] \cup \{[\alpha]\}$, so for $\Psi^{\alpha}:=\psi^{\alpha}\big((\Gamma+\Z\alpha)^{\ne}\big)$ we have:
\begin{equation}\label{eq:Psialpha} \Psi^{\alpha}\ =\ \Psi\cup\{\beta-\alpha\}, \qquad \beta-\alpha\ =\ \max\Psi^{\alpha}.
\end{equation}
Lemma~\ref{extension1} goes through with $\alpha < 0$ and $\alpha_1 <0$ in place of $\alpha >0$ and $\alpha_1 > 0$, respectively. In the setting of this modified lemma
we have $\Gamma^{<} < n\alpha < 0$
for all~$n\ge 1$, and \eqref{eq:psialpha} goes through for $\gamma\in \Gamma$ and $k\in \Z$ with $\gamma + k\alpha \ne 0$, and \eqref{eq:Psialpha} goes through. So we have really two ways to
remove a gap, and this is a pervasive {\em fork in the road.}\/ In any case,
removal of a gap as above leads by \eqref{eq:Psialpha} to a grounded $H$-asymptotic couple, and this is the situation we consider next. 

\begin{lemma}\label{extas2} Assume $\Psi$ has a largest element $\beta$. 
Then there is an 
$H$-asymptotic couple $(\Gamma+\Z \alpha,
\psi^{\alpha})$ extending $(\Gamma,\psi)$ with $\alpha\ne 0$, $\alpha'=\beta$, such that for any embedding $i\colon (\Gamma,\psi)\to (\Gamma_1,\psi_1)$ and any $\alpha_1\in\Gamma_1^{\ne}$ with $\alpha_1'=i(\beta)$ there is a unique extension of~$i$ to an embedding $j\colon(\Gamma+\Z \alpha,
\psi^{\alpha})\to (\Gamma_1,\psi_1)$ with $j(\alpha)=\alpha_1$.
\end{lemma}
\begin{proof} Suppose $(\Gamma+\Z \alpha,
\psi^{\alpha})$ is an asymptotic couple extending $(\Gamma,\psi)$ with
$\alpha\ne 0$, $\alpha'=\beta$. In the divisible hull
of  $(\Gamma+\Z \alpha,
\psi^{\alpha})$ we have for $\gamma\in (\Q\Gamma)^{<}$, 
$$\gamma'\ =\ \gamma+\psi(\gamma)\ <\ \psi(\gamma)\ \le\ \beta\ =\ \alpha',$$ 
so $\gamma< \alpha< 0$ by Lemma~\ref{BasicProperties}(iii).
Hence $\Gamma^{<} < n\alpha < 0$ for all~$n\ge 1$, in particular,
$\Gamma + \Z\alpha=\Gamma\oplus \Z\alpha$. Also $\psi^{\alpha}(\alpha)= \beta-\alpha > \Psi$, and thus \eqref{eq:psialpha} holds
for all $\gamma\in \Gamma$ and $k\in \Z$ with $\gamma+k\alpha\ne 0$. 
It easily follows that $(\Gamma+\Z\alpha, \psi^{\alpha})$ is of $H$-type, and 
has the required universal property. All this assumes the existence
of an asymptotic couple~${(\Gamma+\Z\alpha, \psi^{\alpha})}$ extending 
$(\Gamma, \psi)$ with $\alpha\ne 0$, $\alpha'=\beta$, but the above also 
suggests how to construct it. The details are similar to those in the proof of Lemma~\ref{extension1} and are left to the reader. 
%That the resulting construction leads 
%to an asymptotic couple $(\Gamma+\Z\alpha, \psi^{\alpha})$ of $H$-type
%follows from \eqref{eq:psialpha}.  
\end{proof}

\medskip\noindent
Let $(\Gamma+\Z\alpha, \psi^{\alpha})$ be as in Lemma~\ref{extas2}. Then
$[\Gamma+\Z\alpha]=[\Gamma] \cup \{[\alpha]\}$, so \eqref{eq:Psialpha}
holds for $\Psi^{\alpha}:=\psi^{\alpha}\big((\Gamma+\Z\alpha)^{\ne}\big)$.
Thus our new $\Psi$-set $\Psi^{\alpha}$ still has a maximum, but
this maximum is larger than the maximum $\beta$ of the original 
$\Psi$-set $\Psi$. By iterating this construction indefinitely,
taking a union, and passing to the divisible hull, we obtain a divisible
$H$-asymptotic couple with asymptotic integration. Once we have a divisible $H$-asymptotic couple with asymptotic integration, we can create an extension
with a gap as follows:

\begin{lemma}\label{addgap} Suppose $(\Gamma, \psi)$ has asymptotic 
integration and $\Gamma$ is divisible. 
Then there is an $H$-asymptotic couple 
$(\Gamma+ \Q\beta, \psi_{\beta})$ extending $(\Gamma, \psi)$ such that: \begin{enumerate}
\item[\textup{(i)}] $\Psi< \beta < (\Gamma^{>})'$;
\item[\textup{(ii)}] for any divisible $H$-asymptotic $(\Gamma_1, \psi_1)$ extending 
$(\Gamma, \psi)$ and $\beta_1\in \Gamma_1$ with  $\Psi< \beta_1 < (\Gamma^{>})'$
there is a unique embedding 
$(\Gamma+ \Q\beta, \psi_{\beta})\to (\Gamma_1, \psi_1)$ of asymptotic 
couples that is the 
identity on $\Gamma$ and sends $\beta$ to $\beta_1$. 
\end{enumerate}
\end{lemma}
\begin{proof} Since $\Psi$ has no largest element, we can take an 
elementary ex\-ten\-sion $(\Gamma^*, \psi^*)$ of
$(\Gamma,\psi)$ with an element $\beta\in \Gamma^*$ such that
$\Psi< \beta < (\Gamma^{>})'$. Moreover, for each $\alpha\in \Gamma^{>}$ we have
$\alpha^\dagger < \beta < \alpha'$ and $\alpha'-\alpha^\dagger= \alpha$, so
$\Gamma$ is dense in $\Gamma+ \Q \beta$ by Lemma~\ref{gendense}. 
Hence $[\Gamma+\Q \beta]=[\Gamma]$, and thus 
$\psi^*$ maps $(\Gamma+\Q \beta)^{\ne}$ into $\Psi$, and with~$\psi_{\beta}$
the restriction of $\psi^*$ to $(\Gamma+\Q \beta)^{\ne}$ we have an 
$H$-asymptotic 
couple $(\Gamma+ \Q\beta, \psi_{\beta})$ extending $(\Gamma, \psi)$
satisfying (i). Let $(\Gamma_1, \psi_1)$ be a divisible $H$-asymptotic couple extending 
$(\Gamma, \psi)$ and $\beta_1\in \Gamma_1$, $\Psi< \beta_1 < (\Gamma^{>})'$.
By the universal property of Lemma~\ref{uniorvec} we have a unique
embedding $\Gamma+ \Q\beta \to \Gamma_1$
of ordered vector spaces over $\Q$ that is the identity on $\Gamma$ and sends
$\beta$ to $\beta_1$. It is routine to check that this embedding is also
an embedding $(\Gamma+ \Q\beta, \psi_{\beta})\to (\Gamma_1, \psi_1)$ 
of asymptotic couples.  
\end{proof}

\noindent
Let $(\Gamma, \psi)$ be divisible with asymptotic
integration and let $(\Gamma+ \Q\beta, \psi_{\beta})$ be an 
$H$-asymptotic couple as in Lemma~\ref{addgap}. If 
$(\Gamma + \Q\alpha, \psi_{\alpha})$ is also a divisible $H$-asymptotic couple 
extending $(\Gamma, \psi)$ with $\Psi< \alpha < (\Gamma^{>})'$, then
by (ii) we have an isomorphism 
$(\Gamma+ \Q\beta, \psi_{\beta})\to (\Gamma+ \Q\alpha, \psi_{\alpha})$ 
of asymptotic couples
that is the identity on $\Gamma$ and sends $\beta$ to $\alpha$.
In this sense, $(\Gamma + \Q\beta, \psi_{\beta})$ is unique up to isomorphism
over~$(\Gamma, \psi)$. Thus the construction of $(\Gamma+ \Q\beta, \psi_{\beta})$
in the proof of Lemma~\ref{addgap} gives the following extra information,
with $\Psi_{\beta}$ the set of values of $\psi_{\beta}$ on 
$(\Gamma+ \Q\beta)^{\ne}$: 

\begin{cor}\label{addgapcor} The set $\Gamma$ is dense 
in the ordered abelian group 
$\Gamma+ \Q\beta$, so 
$[\Gamma] = [\Gamma+ \Q\beta]$, 
$\Psi_{\beta}=\Psi$ and $\beta$ is a gap in 
$(\Gamma+ \Q\beta, \psi_{\beta})$.  
\end{cor}

\noindent
Here is a version of Lemma~\ref{addgap} for possibly non-divisible $\Gamma$:

\begin{cor}\label{addgapmore}
Suppose $(\Gamma,\psi)$ has rational asymptotic integration. Then there is an $H$-asymptotic couple $(\Gamma +\Z\beta,\psi_\beta)$ extending $(\Gamma,\psi)$ such that:
\begin{enumerate}
\item[\textup{(i)}] $\Psi<\beta<(\Gamma^>)'$, $k\beta\notin \Gamma$ for all $k\in \Z^{\ne}$, and $[\Gamma+\Z\beta]=[\Gamma]$;
\item[\textup{(ii)}] for any $H$-asymptotic couple $(\Gamma_1,\psi_1)$ extending $(\Gamma,\psi)$ and any 
 $\beta_1\in\Gamma_1$ with $\Psi<\beta_1<(\Gamma^>)'$ there is a unique  embedding $(\Gamma\oplus\Z\beta,\psi_\beta)\to (\Gamma_1,\psi_1)$ of asymptotic couples that is the identity on $\Gamma$ and sends $\beta$ to
 $\beta_1$. 
\end{enumerate}
Moreover, $\beta$ is a gap in $(\Gamma +\Z\beta,\psi_\beta)$ for  any such 
extension of $(\Gamma, \psi)$. 
 \end{cor}
\begin{proof}
By assumption, $(\Q\Gamma,\psi)$ has asymptotic integration, so we can take an $H$-asymptotic couple
$(\Q\Gamma+\Q\beta,\psi_\beta)$ extending $(\Q\Gamma,\psi)$ and satisfying (i) and (ii) in Lemma~\ref{addgap},
with $(\Q\Gamma,\psi)$ replacing $(\Gamma,\psi)$. Note that $k\beta\notin \Gamma$ for all $k\in \Z^{\ne}$.
By Corollary~\ref{addgapcor} we have $\psi_{\beta}\big((\Gamma+\Z\beta)^{\ne}\big)=\Psi$.
Denoting the restriction of $\psi_{\beta}$ to $(\Gamma+\Z\beta)^{\ne}$
also by $\psi_{\beta}$ we obtain therefore
an $H$-asymptotic couple $(\Gamma\oplus\Z\beta,\psi_{\beta})$ extending
$(\Gamma,\psi)$ satisfying (i) and (ii) in the present corollary.
\end{proof}

\noindent 
Recall from Section~\ref{sec:ordered sets} that a cut in an 
ordered set $S$ is just a downward closed subset of $S$, and that 
an element $a$ of an ordered set extending $S$ is said to
realize a cut $C$ in $S$ if $C<a<S\setminus C$ (so $a\notin S$).

\begin{lemma}\label{extension5}
Let $C$ be a cut in $[\Gamma^{\ne}]$ and let $\beta\in\Gamma$ be such that
$\beta<(\Gamma^{>})'$,
$\gamma^\dagger \leq\beta$ for all $\gamma\in\Gamma^{\ne}$ with $[\gamma]> C$,
and $\beta\leq \delta^\dagger$ for all $\delta\in\Gamma^{\ne}$ with $[\delta]\in C$. 
Then there 
exists an $H$-asymptotic couple $(\Gamma\oplus\Z\alpha,\psi^\alpha)$ 
extending
$(\Gamma,\psi)$, with $\alpha>0$, such that:
\begin{enumerate}
\item[\textup{(i)}] $[\alpha]$ realizes the cut $C$ in $[\Gamma^{\ne}]$, and 
$\psi^\alpha(\alpha)=\beta$;
\item[\textup{(ii)}] given any embedding $i$ of $(\Gamma,\psi)$ into an $H$-asymptotic couple
$(\Gamma_1,\psi_1)$ 
and any element $\alpha_1\in \Gamma_1^{>}$ such that $[\alpha_1]$ realizes 
the cut $\bigl\{\bigl[i(\delta)\bigr]:
[\delta]\in C\bigr\}$
in $\bigl[i(\Gamma^{\ne})\bigr]$ and $\psi_1(\alpha_1)=i(\beta)$, there is
a unique extension of $i$ to an embedding $j\colon 
{(\Gamma\oplus\Z\alpha,\psi^\alpha)}\to (\Gamma_1,\psi_1)$ with $j(\alpha)=
\alpha_1$.
\end{enumerate}
If $(\Gamma, \psi)$ has asymptotic integration, then so does
$(\Gamma\oplus \Z\alpha, \psi^\alpha)$. If $(\Gamma, \psi)$  has rational asymptotic integration, then so does $(\Gamma\oplus \Z\alpha, \psi^\alpha)$.
\end{lemma}
\begin{proof} By Lemma~\ref{exarchextclass} we can extend $\Gamma$ to an ordered abelian group
$\Gamma^\alpha:=\Gamma\oplus\Z\alpha$ with $\alpha>0$ such that $[\alpha]$ realizes the cut $C$ in $[\Gamma^{\neq}]$. Then $[\Gamma^\alpha]=[\Gamma]\cup\bigl\{[\alpha]
\bigr\}$. We extend $\psi\colon\Gamma^{\ne}\to\Gamma$ 
to $\psi^\alpha\colon (\Gamma^\alpha)^{\ne}\to\Gamma$ by 
$$\psi^\alpha(\gamma+k\alpha)\ :=\ \min\bigl\{\psi(\gamma),\beta\bigr\} \quad
\text{for $\gamma\in\Gamma$, $k\in\Z^{\ne}$.}$$ 
(So $\psi^\alpha\bigl((\Gamma^\alpha)^{\neq}\bigr)=\Psi\cup\{\beta\}$.)
A tedious but routine checking of cases shows that~$\psi^\alpha$ decreases on $(\Gamma^\alpha)^{>}$, and that
axioms~(AC1) and (AC2) for asymptotic couples hold for $(\Gamma^\alpha,
\psi^\alpha)$. 
To verify (AC3), let $\delta=\gamma+k\alpha$ and $\delta^*=\gamma^*+k^*\alpha$ be  elements of~$(\Gamma^\alpha)^>$ ($\gamma,\gamma^*\in\Gamma$, $k,k^*\in\Z$); we have to show that $\psi^\alpha(\delta^*)<\delta+\psi^\alpha(\delta)$. We can assume $[\delta^*]<[\delta]$, since otherwise $\psi^\alpha(\delta^*)\leq\psi^\alpha(\delta)<\delta+\psi^\alpha(\delta)$. We distinguish the following cases, using Lemma~\ref{BasicProperties}(i) in Cases 2 and 3:

\case[1]{$[\delta^*]=[\gamma^*]$, $[\delta]=[\gamma]$.}
Then
$\bigl[\psi(\gamma^*)-\psi(\gamma)\bigr]<[\gamma^*-\gamma]=[\gamma]=[\delta]$, hence
$\psi^\alpha(\delta^*)=\psi(\gamma^*)<\psi(\gamma)+\delta=\psi^\alpha(\delta)+\delta$.

\case[2]{$[\delta^*]=[\alpha]$, $[\delta]=[\gamma]$.} Then $[\gamma]>C$, so $\psi(\gamma)\leq\beta$.
Hence
$\psi\big(\beta-\psi(\gamma)\big)>\psi(\gamma)$,
so
 $\bigl[\beta-\psi(\gamma)\bigr]<[\gamma]=[\delta]$, and thus
$\psi^\alpha(\delta^*)=\beta < \psi(\gamma)+\delta=\psi^\alpha(\delta)+\delta$.

\case[3]{$[\delta^*]=[\gamma^*]$, $[\delta]=[\alpha]$.}
Then $[\gamma^*]\in C$, so 
$$\psi\big(\psi(\gamma^*)-\beta\big)\ >\ \min\big\{\psi(\gamma^*),\beta\big\}\ =\ \beta,$$ so $\big[\psi(\gamma^*)-\beta\big]\in C$ or $\psi(\gamma^*)=\beta$.
Hence
$\bigl[\psi(\gamma^*)-\beta\bigr]<[\alpha]=[\delta]$, and thus 
$\psi^\alpha(\delta^*)=\psi(\gamma^*)<\beta+\delta=\psi^\alpha(\delta)+\delta$.

\medskip
\noindent
So $(\Gamma^{\alpha}, \psi^{\alpha})$ is indeed an 
$H$-asymptotic couple
satisfying (i). That it satisfies (ii) follows easily from 
the universal property of  Lemma~\ref{exarchextclass}.
 
%Assume that $(\Gamma, \psi)$ has asymptotic integration. Then $\Psi$ has no
%maximum and 
%$\beta\in \Psi^{\downarrow}$, so 
%$\Psi^{\alpha}=\Psi\cup\{\beta\}$ has no maximum. Suppose towards a 
%contradiction that $(\Gamma^{\alpha}, \psi^{\alpha})$ has a gap 
%$\gamma\in \Gamma^{\alpha}$, so $\Psi^{\alpha} < \gamma < 
%\big((\Gamma^{\alpha})^{>}\big)'$. 
%Working in the divisible hull of $(\Gamma^{\alpha}, \psi^{\alpha})$.
%Inside $\Q\Gamma^{\alpha}$ we have 
%$\Q\Gamma \cap (\Gamma\oplus\Z\alpha)=\Gamma$, so
%$\gamma\notin \Q\Gamma$, hence $\Q\Gamma^{\alpha}=\Q\Gamma + \Q\gamma$.
%Note also that $[\Gamma^{\ne}]$ has no least element and equals 
%$[(\Q\Gamma)^{\ne}]$. Thus we can
%apply Lemma~\ref{gendense} to the extension 
%$\Q\Gamma^{\alpha}\supseteq \Q\Gamma$ of ordered vector spaces 
%over $\Q$, with $P$ as the downward closure of
%$\Psi$ in $\Q\Gamma$ and with $\gamma$ in the role of~$b$. 
%It follows that $\Q\Gamma$ is dense in $\Q\Gamma^{\alpha}$, 
%hence $[\Q\Gamma]=[\Q\Gamma^{\alpha}]$,
% contradicting $[\alpha]\notin [\Gamma]=[\Q\Gamma]$.

Assume $(\Gamma, \psi)$ has asymptotic integration; given $\gamma\in\Gamma$ and $k\in \Z^{\ne}$ we shall find
an antiderivative of $\gamma+k\alpha$ in $(\Gamma^{\alpha}, \psi^{\alpha})$. If
$\gamma = \beta$, then $k\alpha$ is such an antiderivative, so
assume $\gamma\neq \beta$. Lemma~\ref{lem:minpsipsi2} gives the
asymptotic couple $\big(\Gamma, \min(\psi,\beta)\big)$ whose $\Psi$-set
has maximum element $\beta$. Since $\gamma\neq \beta$, Theorem~\ref{26} gives
$\gamma^*\in\Gamma$ such that
$\gamma^*+\min\big(\psi(\gamma^*),\beta\big) = \gamma$ and so
$\gamma^*+k\alpha + \psi^{\alpha}(\gamma^*+k\alpha) = \gamma+k\alpha$, that is,
$\gamma^*+k\alpha$ is an antiderivative as required. 

 Preserving rational asymptotic integration is done likewise.
%but now $\gamma\in \Q\Gamma^{\alpha}$, and we get 
%$\gamma\notin \Q\Gamma$ from 
%$\Psi < \gamma < \big((\Q\Gamma)^{>}\big)'$. 
\end{proof}

\noindent
For $C=\emptyset$ and $\beta$ a gap in $(\Gamma, \psi)$, this gives:

\begin{cor} Let $\beta\in \Gamma$ be a gap in $(\Gamma, \psi)$. Then there exists an $H$-asymptotic couple $(\Gamma+\Z\alpha,\psi^\alpha)$ 
extending $(\Gamma,\psi)$, such that:
\begin{enumerate}
\item[\textup{(i)}] $0 < n\alpha < \Gamma^{>}$ for all $n\ge 1$, and 
$\psi^\alpha(\alpha)=\beta$;
\item[\textup{(ii)}] for any embedding $i$ of $(\Gamma,\psi)$ into an $H$-asymptotic couple
$(\Gamma_1,\psi_1)$ 
and any $\alpha_1\in \Gamma_1^{>}$ with $\psi_1(\alpha_1)=i(\beta)$, there is
a unique extension of $i$ to an embedding $j\colon 
(\Gamma+\Z\alpha,\psi^\alpha)\to (\Gamma_1,\psi_1)$ with $j(\alpha)=
\alpha_1$.
\end{enumerate}
\end{cor}

\subsection*{Notes and comments} As already mentioned, much of this 
section comes from \cite[Section 2]{AvdD2}. For example,
Lemma~\ref{extas2} is a special case of~\cite[Lemma~2.12]{AvdD2}, and
Lemma~\ref{extension5}  combines Lemma~2.15 of \cite{AvdD2} and 
a remark that follows the proof of that lemma.

%% file: mt-9n-9.tex
\section{Closed $H$-Asymptotic Couples}\label{sec:cac}

\noindent
An $H$-asymptotic couple 
$(\Gamma,\psi)$ is said to be
{\bf closed\/} if it is divisible with asymptotic integration and  $\Psi:=\psi(\Gamma^{\ne})$ is downward
closed. At the beginning of Section~\ref{sec:Liouville closed}
we indicate why the $H$-asymptotic couple
of $\T$ is closed.

By \cite{AvdD}, closed $H$-asymptotic couples admit quantifier elimination.
A step in that direction is 
Proposition~\ref{cac} below, which is needed in 
Section~\ref{sec:valgrp}. 
In this book we do not need quantifier elimination 
for closed $H$-asymptotic couples, but to give 
some orientation for what follows we mention a consequence of it:
$$\text{ closed = existentially closed \qquad(for $H$-asymptotic couples)}.$$
In the first subsection we prove an easy part of 
this fact, namely that every $H$-asymp\-to\-tic couple extends to a 
closed $H$-asymptotic couple. 

\index{H-asymptotic@$H$-asymptotic!couple!closed}
\index{asymptotic couple!closed}
\index{closed!H-asymptotic couple@$H$-asymptotic couple}

\subsection*{Embedding $H$-asymptotic couples into closed $H$-asymptotic couples} Let us construe an asymptotic couple $(\Gamma, \psi)$ 
as an ordered group $\Gamma$ equipped with a binary relation on $\Gamma$, 
namely 
the graph of $\psi\colon \Gamma^{\ne} \to \Gamma$. In this way the
$H$-asymptotic couples are exactly the models of a set of 
$\forall\exists$-sentences in the
language of ordered abelian groups augmented by a binary relation symbol. 
By Section~\ref{sec:mc} and Lemma~\ref{lem:ec models} this gives 
the notion of an $H$-asymptotic couple
being {\em existentially closed\/},
and the fact that every $H$-asymptotic couple extends to an existentially closed one.

\begin{lemma}\label{eccac} Existentially closed $H$-asymptotic couples are
closed. 
\end{lemma}
\begin{proof} Let $(\Gamma, \psi)$ be an existentially closed $H$-asymptotic 
couple. Remarks after the proof of 
Lemma~\ref{ro} show that $\Gamma$ is divisible. Next, it follows from
Lemma~\ref{extension1} that~${(\Gamma, \psi)}$ has no gap,
from Lemma~\ref{extas2} that  $(\Gamma, \psi)$ is not grounded (so it has asymptotic integration) and from
Lemma~\ref{extension5} that $\Psi$ is downward closed.  
\end{proof}

\subsection*{A closure property of closed $H$-asymptotic couples}
Let $(\Gamma, \psi)$ be an asymptotic couple. Recall from~\ref{sec:ascouples} that we ex\-ten\-ded ${\psi\colon \Gamma^{\ne}\to \Gamma}$ to a function $\psi\colon \Gamma_{\infty} \to \Gamma_{\infty}$ by $\psi(0)=\psi(\infty):= \infty$. For $\alpha_1,\dots,\alpha_n\in \Gamma$, $n\ge 1$, we define the function  
$\psi_{\alpha_1,\dots,\alpha_n}\colon \Gamma_{\infty}\to \Gamma_{\infty}$ by
recursion on $n$: 
$$\psi_{\alpha_1}(\gamma)\ :=\ \psi(\gamma-\alpha_1), \qquad
\psi_{\alpha_1,\dots, \alpha_{n}}(\gamma)\ :=\ \psi\big(\psi_{\alpha_1,\dots, \alpha_{n-1}}(\gamma)-\alpha_{n}\big) 
\text{ for $n\ge 2$.}$$     

\begin{prop}\label{cac} Let $(\Gamma, \psi)$ be a closed $H$-asymptotic couple 
and let $(\Gamma^*,\psi^*)$ be an $H$-asymptotic couple
extending $(\Gamma,\psi)$. Suppose $n\ge 1$, 
$\alpha_1,\dots, \alpha_n\in \Gamma$, $q_1,\dots, q_n\in \Q$ and 
$\gamma\in \Gamma^*$ are such that
\begin{align*} \psi^*_{\alpha_1,\dots,\alpha_n}(\gamma)&\ne \infty\quad (\text{so
$\psi^*_{\alpha_1,\dots,\alpha_i}(\gamma)\ne \infty$ for $i=1,\dots,n$}),\ \text{ and}\\
\gamma + q_1\psi^*_{\alpha_1}(\gamma) &+ \cdots + q_n\psi^*_{\alpha_1,\dots,\alpha_n}(\gamma)\in \Gamma\qquad   \text{\textup{(}in $\Q\Gamma^*$\textup{)}}. \end{align*}
Then $\gamma\in \Gamma$.
\end{prop} 

\noindent
In the rest of this section we establish Proposition~\ref{cac}.

\subsection*{Some lemmas} Let $D$ be a subset of an ordered abelian group $\Gamma$. We say that $D$ is {\bf bounded\/} if $D\subseteq [p,q]$ for
some $p\le q$ in $\Gamma$, and otherwise we call $D$ {\bf unbounded}. (These notions and the next one are with respect to the ambient $\Gamma$.)
A {\bf (convex) component of $D$} is by 
definition a nonempty convex subset $C$ of $\Gamma$ such that 
$C\subseteq D$ and $C$ is maximal with these properties.  \index{component!convex}
The components of $D$ partition the set $D$: for $d\in D$ the unique component of 
$D$ containing $d$ is 
$$\big\{\gamma\in D^{\le d}:\ [\gamma,d]\subseteq D\big\} \cup \big\{\gamma\in D^{\ge d}:\ 
[d,\gamma]\subseteq D\big\}.$$ 
Now let $(\Gamma, \psi)$ be an $H$-asymptotic couple,
$n\ge 1$, and let $\alpha$ be a sequence $\alpha_1,\dots, \alpha_n$ from $\Gamma$. 
We set $$D_{\alpha}\ :=\ \big\{\gamma\in \Gamma:\ \psi_{\alpha}(\gamma)\ne \infty\big\}.$$
Thus 
\begin{align*}
D_{\alpha}\	&=\ \Gamma\setminus \{\alpha_1\} &&\text{for $n=1$, and} \\
D_{\alpha}\	&=\ \big\{\gamma\in D_{\alpha'}:\ \psi_{\alpha'}(\gamma)\ne \alpha_n\big\} &&\text{for $n>1$ and $\alpha'=\alpha_1,\dots, \alpha_{n-1}$.} 
\end{align*}
One checks easily by induction on $n$ that for distinct 
$\gamma, \gamma'\in  D_{\alpha}$,
$$\psi_{\alpha}(\gamma)- \psi_{\alpha}(\gamma')\ =\ o(\gamma-\gamma').$$ 

\begin{lemma}\label{siivp} Assume $\Gamma$ is divisible, and let $q_1,\dots, q_n\in \Q$. Then the function
$$\gamma\ \mapsto\ \gamma+ q_1\psi_{\alpha_1}(\gamma) + q_2\psi_{\alpha_1,\alpha_2}(\gamma)+ \cdots + q_n \psi_{\alpha}(\gamma)\ :\ D_{\alpha} \to \Gamma$$
is strictly increasing. Moreover, this function has the intermediate value property on every
component of  $D_{\alpha}$.
\end{lemma}
\begin{proof} Let $\eta\colon D_{\alpha} \to \Gamma$ be the function given by
$$\eta(\gamma)\ :=\ q_1\psi_{\alpha_1}(\gamma) +  q_2\psi_{\alpha_1,\alpha_2}(\gamma)+ \cdots + q_n \psi_{\alpha}(\gamma).$$
The function $\gamma \mapsto \gamma+ \eta(\gamma)\colon D_{\alpha} \to \Gamma$ is strictly increasing since $\eta(\gamma)-\eta(\gamma')=o(\gamma-\gamma')$ for 
all distinct
$\gamma, \gamma'\in D_{\alpha}$. Let $C$ be a component of $D_{\alpha}$ with
$\alpha_1< C$, and let $\beta< \gamma_1 < \gamma_2$ in $C$, with 
$\gamma_2-\gamma_1\le \gamma_1-\beta$. Then
$$0 \ <\ \gamma_1-\alpha_1\ <\ \gamma_2-\alpha_1\ =\ 
(\gamma_2-\gamma_1)+(\gamma_1-\alpha_1)\ \le\ 2(\gamma_1-\alpha_1),$$
so $\psi(\gamma_1-\alpha_1)=\psi(\gamma_2-\alpha_1)$. Hence 
$\eta(\gamma_1)=\eta(\gamma_2)$, since $\eta(\gamma)$ depends only on 
$\psi(\gamma-\alpha_1)$. Using terminology of Section~\ref{sec:oag} we have shown that $\eta|C$ is $v$-slow on the right, where $v$ is the standard valuation of $\Gamma$. Thus by Lemma~\ref{sloste} the function $\gamma\mapsto \gamma+\eta(\gamma)$ on
$C$ has the intermediate value property. 
For the components $<\alpha_1$ of $D_{\alpha}$
we can argue likewise. 
\end{proof}

\noindent
{\em In the rest of this subsection the $H$-asymptotic couple $(\Gamma, \psi)$ is closed}.

\begin{lemma} The set $D_{\alpha}$ has at most $2^n$ components, and on each of these the function $\psi_{\alpha}$ is monotone and has
the intermediate value property.
\end{lemma}
\begin{proof} For $n=1$ we have $\alpha=\alpha_1$, and the two
components of $D_{\alpha}$ are $\Gamma^{<\alpha}$, on which $\psi_{\alpha}$ is increasing,
and $\Gamma^{>\alpha}$, on which $\psi_{\alpha}$ is decreasing. On each of these
$\psi_{\alpha}$ has the intermediate value property, since $\Psi$ is 
downward closed. Suppose the lemma holds for a certain 
$\alpha=\alpha_1,\dots, \alpha_n$, let $\alpha_{n+1}\in \Gamma$, and set
${\alpha+}:=\alpha_1,\dots, \alpha_n,\alpha_{n+1}$.
Consider a component $C$ of $D_{\alpha}$. Then  $\psi_\alpha$ is
monotone on $C$, say increasing on~$C$, and has the intermediate
value property on $C$. Put
$$C_1 \ :=\ \bigl\{\gamma\in C:\ \psi_\alpha(\gamma)<\alpha_{n+1}\bigr\},$$
and similarly define $C_2$ and $C_3$, with $=$ and $>$, respectively, replacing $<$.
%\begin{align*}
%  C_1\ &:=\ \bigl\{\gamma\in C:\ \psi_\alpha(\gamma)<\alpha_{n+1}\bigr\}, \\
%  C_2\ &:=\ \bigl\{\gamma\in C:\ \psi_\alpha(\gamma)=\alpha_{n+1}\bigr\}, \\
%  C_3\ &:=\ \bigl\{\gamma\in C:\ \psi_\alpha(\gamma)>\alpha_{n+1}\bigr\}.
%\end{align*}
Thus $C$ is the disjoint union of its convex subsets $C_1$, $C_2$ and
  $C_3$,
and $C_1<C_2<C_3$. Also $C\cap D_{\alpha+}= C_1\cup C_3$, 
$\psi_{\alpha+}$ is clearly increasing on $C_1$ and
decreasing on $C_3$, and has the intermediate value property on $C_1$ and on $C_3$. 
If both $C_1$ and $C_3$ are nonempty, then $C_2$ is
nonempty (because of the intermediate value property of $\psi_\alpha$ on
$C$), and thus~$C_1$ and~$C_3$ are the components of $D_{\alpha+}$ that
are contained in $C$. Otherwise $C$ only contributes one
component to $D_{\alpha+}$, or none at all, depending on whether just one or both
of $C_1$ and $C_3$ are empty.
\end{proof}

\noindent
At this point we need more 
terminology. Let $f\colon \Gamma_\infty\to \Gamma_\infty$ be a function, and 
let~$C$ be a
nonempty convex subset of $\Gamma$ on which $f$ does not take the value 
$\infty$. Let
$p,q\in \Gamma$, and let $S\subseteq \Gamma$ be downward closed. (We only use this for $f=\psi_\alpha$ with $C$ a component of $D_\alpha$, and $S=\Psi$.)
\begin{enumerate}
\item {\bf $f$ increases on $C$ from $p$ to $q$} if $f|C$ is increasing, $p\leq
q$, and $f(C)=[p,q]$;
\item {\bf $f$ decreases on $C$ from $p$ to $q$} if $f|C$ is decreasing, $p\geq
q$, and $f(C)=[q,p]$;
\item {\bf $f$ increases on $C$ from $p$ to $S$} if $f|C$ is increasing, and
$f(C)=S^{\ge p}$;
\item {\bf $f$ decreases on $C$ from $S$ to $q$} if $f|C$ is decreasing, and
$f(C)=S^{\ge q}$;
\item {\bf $f$ decreases on $C$ from $S$ to $-\infty$} 
if $f|C$ is decreasing, and $f(C)=S$;
\item {\bf $f$ decreases on $C$ from $p$ to $-\infty$} 
if $f|C$ is decreasing, and $f(C)=(-\infty,p]$.
\end{enumerate}
Next, let $(\Gamma^*,\psi^*)$ be a closed $H$-asymp\-to\-tic couple that extends 
our closed $H$-asymp\-to\-tic couple $(\Gamma,\psi)$. Besides 
$D_{\alpha}\subseteq \Gamma$ we now also have the set 
$D^*_{\alpha}\subseteq \Gamma^*$
with its components, taken
relative to the ambient $(\Gamma^*,\psi^*)$; note that   
$D^*_{\alpha}\cap \Gamma=D_{\alpha}$.

\begin{lemma}\label{LongLemma} The components $C$ of $D_{\alpha}$ have the 
following properties: 
\begin{enumerate}
\item[\textup{(i)}] $C$ is contained in a \textup{(}necessarily unique\textup{)} component
$C^*$ of $D^*_\alpha$, and the map ${C\mapsto C^*}$ is a bijection from the set of
components of $D_\alpha$ onto the set of components of $D^*_\alpha$, with 
$C^*\cap \Gamma=C$ for
each $C$;
\item[\textup{(ii)}] $D_\alpha$ has a \rom{(}necessarily unique\rom{)} unbounded component $C_\infty>\alpha_1$; the corresponding component $C^*_\infty$ of $D^*_\alpha$ is
$>\alpha_1$ and unbounded in $\Gamma^*$;
\item[\textup{(iii)}] for bounded $C$ there are $p,q\in \Gamma$ such that
one of the following holds:
\begin{itemize}
\item[\textup{(a)}] $\psi_\alpha$ increases on $C$ from $p$ to $q$ and
      $\psi^*_\alpha$ increases on $C^*$ from $p$ to $q$,
\item[\textup{(b)}] $\psi_\alpha$ decreases on $C$ from $p$ to $q$ and
      $\psi^*_\alpha$ decreases on $C^*$ from $p$ to $q$,
\item[\textup{(c)}] $\psi_\alpha$ increases on $C$ from $p$ to $\Psi$ and
      $\psi^*_\alpha$ increases on $C^*$ from $p$ to $\Psi^*$,
\item[\textup{(d)}] $\psi_\alpha$ decreases on $C$ from $\Psi$ to $q$ and
      $\psi^*_\alpha$ decreases on $C^*$ from $\Psi^*$ to $q$;
\end{itemize}
\item[\textup{(iv)}] for the unbounded component $C_{\infty}>\alpha_1$ of $D_\alpha$
and the
corresponding component $C_\infty^*$ of $D_\alpha^*$, one of the
following holds:
\begin{itemize}
\item[\textup{(a)}] $\psi_\alpha$ decreases on $C_\infty$ from $\Psi$ to~$-\infty$
      and $\psi_\alpha^*$ decreases on $C^*_\infty$ from $\Psi^*$ to~$-\infty$, 
\item[\textup{(b)}] there is $p\in \Gamma$ such that $\psi_\alpha$ decreases on $C_\infty$ from $p$ to $-\infty$ 
and $\psi^*_\alpha$ decreases on $C^*_\infty$ from $p$ to $-\infty$.
\end{itemize}
\end{enumerate}
\end{lemma}

\noindent
Before we start the proof we note that $\alpha_1\notin D_{\alpha}$ and 
$\psi_{\alpha}(\alpha_1+\gamma)= \psi_{\alpha}(\alpha_1-\gamma)$ for all $\gamma\in \Gamma$. Thus $\alpha_1+\gamma\mapsto \alpha_1-\gamma\colon \Gamma \to \Gamma$ 
maps each component $>\alpha_1$ of $D_{\alpha}$ onto a 
component $<\alpha_1$ of $D_{\alpha}$. The lemma therefore also gives a unique 
unbounded component $<\alpha_1$ of $D_{\alpha}$, with properties symmetric to 
those for $C_{\infty}$.
 
\begin{proof}
We proceed by induction on $n$. The case $n=1$ is easy to verify. Suppose the
lemma holds for a certain $\alpha=\alpha_1,\dots, \alpha_n$, let $\alpha_{n+1}\in \Gamma$, and set
${\alpha+}:=\alpha_1,\dots, \alpha_n,\alpha_{n+1}$. By the remark preceding this 
proof we only need to consider components $>\alpha_1$.
Let $C>\alpha_1$ be a component of $D_\alpha$, and $C^*$ the corresponding component of $D^*_\alpha$.
Set
$$C_1 := \bigl\{\gamma\in C:\ \psi_\alpha(\gamma)<\alpha_{n+1}\bigr\}.$$
Similarly  define $C_2$ and $C_3$, with $=$ and $>$, respectively, replacing $<$,
%\begin{align*}
%C_1\  &:=\  \bigl\{\gamma\in C:\ \psi_\alpha(\gamma)<\alpha_{n+1}\bigr\}, \\
%C_2\  &:=\ \bigl\{\gamma\in C:\ \psi_\alpha(\gamma)=\alpha_{n+1}\bigr\}, \\
%C_3\  &:=\ \bigl\{\gamma\in C:\ \psi_\alpha(\gamma)>\alpha_{n+1}\bigr\},
%\end{align*}
and define the sets $C^*_i$ for $i=1,2,3$ likewise, replacing $C$ 
by $C^*$
and $\psi_\alpha$ by $\psi^*_\alpha$. Hence $C^*_i\cap \Gamma=C_i$, for $i=1,2,3$. 
The components of
$D_{\alpha+}$ contained in $C$ are the nonempty sets among~$C_1$ and~$C_3$,
and the components of $D^*_{\alpha+}$ contained in $C^*$ are the
nonempty sets among $C^*_1$ and $C^*_3$. The inductive assumption
easily gives that for $i=1,2,3$ we have: $C_i\ne \emptyset\Leftrightarrow C^*_i\ne \emptyset$. This proves (i) and (ii).

\medskip\noindent 
Assume that $C$ is bounded in $\Gamma$ (and hence $C^*$ is bounded in $\Gamma^*$).
We also 
assume~$\psi_\alpha$ increases on $C$ and $\psi^*_\alpha$ increases on 
$C^*$.
(The case that $\psi_\alpha$ decreases on $C$ and~$\psi^*_\alpha$ decreases on
$C^*$ is similar and left to the reader.)
We distinguish cases:

\case[1]{There exist $p,q\in \Gamma$ such that $\psi_\alpha$ increases on $C$ from $p$ to $q$ and
$\psi^*_\alpha$ increases on $C^*$ from $p$ to $q$.} Fix such $p,q$. Then we have several subcases:
\begin{enumerate}
\item[(a)] \textit{$q\leq \alpha_{n+1}$.}\/ Then $C_3, C_3^*=\emptyset$. If $q<\alpha_{n+1}$, then
$C_1,C_1^*\neq\emptyset$, $C_2,C_2^*=\emptyset$, $\psi_{\alpha+}$ increases on $C_1$
from $\psi(p-\alpha_{n+1})$ to $\psi(q-\alpha_{n+1})$, and $\psi^*_{\alpha+}$ increases on~$C_1^*$
from $\psi(p-\alpha_{n+1})$ to $\psi(q-\alpha_{n+1})$. If $\alpha_{n+1}=q>p$, then
$C_1,C_1^*\neq\emptyset$, $C_2,C_2^*\neq\emptyset$, and $\psi_{\alpha+}$ increases on $C_1$
from $\psi(p-\alpha_{n+1})$ to $\Psi$, and $\psi^*_{\alpha+}$ increases on $C_1^*$ from
$\psi(p-\alpha_{n+1})$ to $\Psi^*$. If $\alpha_{n+1}=p=q$, then $C_1=\emptyset$
and $C_1^*=\emptyset$.
\item[(b)] \textit{$\alpha_{n+1}\leq p$ and $q\neq \alpha_{n+1}$.}\/ Then $C_1,C_1^*=\emptyset$ and
$C_3,C_3^*\neq\emptyset$. If $\alpha_{n+1}<p$, then $C_2,C_2^*=\emptyset$, and $\psi_{\alpha+}$
decreases on $C_3$ from ${\psi(p-\alpha_{n+1})}$ to ${\psi(q-\alpha_{n+1})}$, and~$\psi_{\alpha+}^*$
decreases on $C_3^*$ from $\psi(p-\alpha_{n+1})$ to $\psi(q-\alpha_{n+1})$. If 
$\alpha_{n+1}=p$,
then $C_2,C_2^*\neq\emptyset$, $\psi_{\alpha+}$ decreases on $C_3$ from $\Psi$ to
$\psi(q-\alpha_{n+1})$, and $\psi_{\alpha+}^*$ decreases on $C_3^*$ from~$\Psi^*$ to
$\psi(q-\alpha_{n+1})$.
\item[(c)] \textit{$p<\alpha_{n+1}<q$.}\/ Then $C_i,C_i^*\neq\emptyset$ ($i=1,2,3$).
%$C_1,C_1^*, C_2,C_2^*, C_3,C_3^*\neq\emptyset$. 
Here, $\psi_{\alpha+}$ increases on $C_1$ from $\psi(p-\alpha_{n+1})$
to $\Psi$, and $\psi_{\alpha+}^*$ increases on $C_1^*$ from $\psi(p-\alpha_{n+1})$ to~$\Psi^*$. Similarly, $\psi_{\alpha+}$ decreases on $C_3$ from $\Psi$ to
$\psi(q-\alpha_{n+1})$, and $\psi_{\alpha+}^*$ decreases on~$C_3^*$ from~$\Psi^*$ to
$\psi(q-\alpha_{n+1})$.
\end{enumerate}
\case[2]{There exists $p\in \Gamma$ such that $\psi_\alpha$ increases on $C$ from $p$ to $\Psi$,
and $\psi^*_\alpha$ increases on $C^*$ from $p$ to $\Psi^*$.} Fix such $p$, and note that then $p\in\Psi$. This case is
essentially treated as Case~1, using Lemmas~\ref{hasc} and~\ref{has+} and some notation introduced after their proofs. If, for example,
$\alpha_{n+1}<p$, so that $C_1,C_1^*=\emptyset$, $C_2,C_2^*=\emptyset$ and
$C_3,C_3^*\neq\emptyset$,
then
$\psi_{\alpha+}$
decreases on $C_3$ from $\psi(p-\alpha_{n+1})$ to
$\psi(*-\alpha_{n+1})$, and
$\psi^*_{\alpha+}$
decreases on $C^*_3$ from $\psi(p-\alpha_{n+1})$ to 
$\psi^*(*-\alpha_{n+1})$, which equals $\psi(*-\alpha_{n+1})$.
We leave the details to the reader.

\medskip
\noindent
Now suppose $C=C_\infty$ is the unbounded component $>\alpha_1$ of $D_\alpha$. 
Then
$C^*=C^*_\infty$ is the unbounded component $>\alpha_1$ of $D_\alpha^*$. We have two cases again:

\case[3]{$\psi_\alpha$ decreases on $C$ from $\Psi$ to $-\infty$, and $\psi_\alpha^*$ decreases
on $C^*$ from $\Psi^*$ to~$-\infty$.} If $\alpha_{n+1}>\Psi$, then
$C_1,C_1^*\neq\emptyset$ and
$C_2,C_2^*,C_3,C_3^*=\emptyset$, so $\psi_{\alpha+}$ decreases on~$C_1$ from
$\psi(*-\alpha_{n+1})$ to $-\infty$, $\psi^*_{\alpha+}$ decreases on $C^*_1$ from
$\psi(*-\alpha_{n+1})$ to $-\infty$, $C_1$ is the unbounded
component $>\alpha_1$ of $D_{\alpha+}$, and $C_1^*$ is the unbounded component~$>\alpha_1$ of~$
D^*_{\alpha+}$.
If, on the other hand, $\alpha_{n+1}\in\Psi$, then $C_1,C_1^*\neq\emptyset$, $C_3,C_3^*\neq
\emptyset$, $C_1 > C_2 > C_3$, $C_1^*> C_2^*> C_3^*$,
$\psi_{\alpha+}$ decreases on $C_1$ from $\Psi$ to $-\infty$, $\psi^*_{\alpha+}$
decreases
on~$C_1^*$ from~$\Psi^*$ to $-\infty$, $\psi_{\alpha+}$ increases on~$C_3$ from
$\psi(*-\alpha_{n+1})$ to $\Psi$, $\psi^*_{\alpha+}$ increases on~$C_3^*$ from
$\psi(*-\alpha_{n+1})$ to $\Psi^*$, the unbounded component~$>\alpha_1$ of $D_{\alpha+}$
is $C_1$, and the unbounded component $>\alpha_1$ of $D^*_{\alpha+}$ is $C_1^*$.

\case[4]{There is $p\in \Gamma$ such that $\psi_\alpha$ decreases on $C$ from $p$ to $-\infty$ and
$\psi_\alpha^*$ decreases on~$C^*$ from $p$ to $-\infty$.}
This case is treated like Case~3, except that we now have three subcases,
according to whether $\alpha_{n+1}>p$, $\alpha_{n+1}=p$, or $\alpha_{n+1}<p$.

\medskip
\noindent
This finishes the inductive step, hence the proof of the lemma.
\end{proof}
 
%\begin{cor}\label{CorPsia} For each component $C$ of $D_\alpha$ we have
%$\psi_\alpha^*(C^*)\cap \Gamma=\psi_\alpha(C)$.
%\end{cor} 

\noindent
Let $q_1,\dots,q_n\in\Q$ and let $\theta\colon D_\alpha\to \Gamma$ be given by
$$ \theta(\gamma)\ :=\ \gamma+q_1\psi_{\alpha_1}(\gamma)+q_2\psi_{\alpha_1,\alpha_2}(\gamma)+\cdots +q_n\psi_\alpha(\gamma).$$
 
\begin{lemma}\label{caca}
Let $C_\infty$ be the unbounded component $>\alpha_1$ of $D_\alpha$. Then $\theta$ is
not bounded from above on $C_\infty$: for any $\beta\in \Gamma$ there exists $\gamma\in C_\infty$ with $\theta(\gamma)>\beta$.
Similarly, $\theta$ is not bounded from below on the unbounded
component $<\alpha_1$ of $D_\alpha$.
\end{lemma}
\begin{proof}
Note that $\big[\theta(\gamma)-\theta(\delta)\big]=[\gamma-\delta]$
for all $\gamma, \delta\in D_\alpha$, and that $[\Gamma]$ has no maximum, by closedness of
$(\Gamma,\psi)$. Given $\beta\in \Gamma$ we pick $\delta\in C_\infty$, and $\gamma>\delta$ such that
$[\gamma]>[\delta],\bigl[\beta-\theta(\delta)\bigr]$. Then
$\bigl[\theta(\gamma)-\theta(\delta)\bigr]>\bigl[\beta-\theta(\delta)\bigr]$, in particular
$\theta(\gamma)>\beta$. 
\end{proof}

\subsection*{Proof of Proposition~\ref{cac}} 

{\sloppy
Recall the setting: 
$(\Gamma, \psi)$ is a closed $H$-asymptotic couple, 
$(\Gamma^*,\psi^*)$ is an $H$-asymptotic couple that
extends $(\Gamma,\psi)$; also, $\alpha$ is a sequence
$\alpha_1,\dots, \alpha_n$ with $n\ge 1$
and $\alpha_1,\dots, \alpha_n\in \Gamma$, and  $q_1,\dots,q_n\in\Q$. 
This yields the function $\theta\colon D_\alpha \to \Gamma$ as defined at the end
of the previous subsection, and likewise we have the function 
$\theta^*\colon D^*_\alpha\to \Gamma^*$ given by
$$\theta^*(\gamma)\ :=\ \gamma+q_1\psi^*_{\alpha_1}(\gamma)+q_2\psi^*_{\alpha_1,\alpha_2}(\gamma)+\cdots +q_n\psi^*_\alpha(\gamma) \qquad(\gamma\in  D^*_\alpha).$$ 
It is clear that $D^*_{\alpha}\cap \Gamma=D_{\alpha}$ and that $\theta^*$ extends 
$\theta$. The claim we have to establish is that every $\gamma\in  D^*_\alpha$
with $\theta^*(\gamma)\in \Gamma$ lies in $\Gamma$. 
We first note that by extending~$(\Gamma^*,\psi^*)$ further, if necessary, we
can arrange that the $H$-asymptotic couple~$(\Gamma^*,\psi^*)$ is also closed; this uses 
Lemmas~\ref{eccac} and~\ref{lem:ec models}. In view of
the results in the previous subsection it suffices to prove under these
conditions:
}

\begin{lemma}
Let $C$ be a component of $D_\alpha$, with corresponding 
component $C^*$ of~$D^*_\alpha$. If $\delta\in C^*\setminus C$,
then $\theta^*(\delta)\in \Gamma^*\setminus \Gamma$.
\end{lemma}
\begin{proof}
By induction on $n$. The case $n=1$ is easily checked using Lemmas~\ref{siivp} and~\ref{id+qpsi}.
Assume the lemma holds for a certain sequence $\alpha=\alpha_1,\dots,\alpha_n$,
and certain
$q_1,\dots,q_n\in\Q$. Let $\alpha_{n+1}\in \Gamma$, put
$\alpha+=\alpha_1,\dots,\alpha_n,\alpha_{n+1}$ and let
$q_{n+1}\in\Q$. Then we have corresponding functions 
$\theta_{+}\colon D_{\alpha+}\to \Gamma$
and $\theta^*_{+}\colon D^*_{\alpha+}\to \Gamma^*$ given by
\begin{align*}
\theta_{+}(\gamma)\ &:=\ \theta(\gamma)+q_{n+1}\psi_{\alpha+}(\gamma),\\
\theta^*_{+}(\gamma)\ &:=\ \theta^*(\gamma)+q_{n+1}\psi^*_{\alpha+}(\gamma).
\end{align*}
Let $C$ be a component of $D_\alpha$ with
corresponding component $C^*$ of $D^*_\alpha$. Define~$C_i$,~$C_i^*$ 
(for $i=1,2,3$) as
in the proof of Lemma~\ref{LongLemma}. Then the components of~$D_{\alpha+}$ contained in $C$ are the nonempty sets among
$C_1$, $C_3$,
and the components of~$D^*_{\alpha+}$ contained in $C^*$ are the
nonempty sets among $C_1^*$, $C_3^*$. We assume $\delta\in C_i^*\setminus C_i$ for $i=1$
or $i=3$, and have to
show that $\theta^*_{+}(\delta)\notin \Gamma$. If $\delta$ lies in the 
convex hull
of~$C_i$ in~$C_i^*$, that
is, if there are $p,q\in C_i$ such that $p<\delta<q$, then the injectivity
of $\theta^*_{+}$
and intermediate value
property of $\theta_{+}|[p,q]$ already guarantee that $\theta^*_{+}(\delta)\in
\Gamma^*\setminus \Gamma$, without use of the induction hypothesis. So from now on, we assume that $\delta$ does
not lie in the convex hull of $C_i$ in $C_i^*$. 

Suppose there exists an element
$\beta\in \Gamma$ lying strictly between $\delta$ and $\alpha_1$, and set
$\varepsilon:=\frac{1}{2}|\beta-\alpha_1|>0$. Then
$\psi^*_{\alpha_1}$ is constant on the segment
$$I\ =\ I_\beta\ :=\ \bigl\{\gamma\in \Gamma^*:\  \delta-\varepsilon\leq \gamma \leq \delta+\varepsilon\bigr\},$$
since $[\gamma-\alpha_1]=[\delta-\alpha_1]$ for all $\gamma\in I$. 
An easy induction on $k$ gives $I\subseteq D^*_{\alpha_1,\dots,\alpha_k}$ and 
$\psi^*_{\alpha_1,\dots,\alpha_k}$ is
constant on $I$,
for $k=1,\dots,n+1$. Hence $I\subseteq C^*_i$,
$\psi^*_{\alpha+}$ is constant on $I$, and 
$\theta^*_{+}(\gamma)=\theta^*_{+}(\delta)+\gamma-\delta$ for all $\gamma\in I$.
If $I\cap C_i\neq\emptyset$, say $\xi\in I\cap C_i$, then
$$\theta^*_{+}(\delta)\ =\ \theta^*(\delta)+q_{n+1}\psi^*_{\alpha+}(\delta)\ =\
\theta^*(\delta)+q_{n+1}\psi_{\alpha+}(\xi)\notin \Gamma,$$
since $\theta^*(\delta)\notin \Gamma$, by the induction hypothesis. 
Thus for the rest
of the proof we assume $I_{\beta}\cap C_i=\emptyset$ for all $\beta\in \Gamma$ strictly
between $\delta$ and $\alpha_1$.
Next, by the remark preceding the proof of Lemma~\ref{LongLemma} we arrange
that $\delta$, $C$ and $C^*$ are all
$> \alpha_1$.

We now first consider the case that $C$ is bounded, 
$\psi_\alpha$ increases on $C$,  $\psi^*_\alpha$ increases on $C^*$, 
and $i=1$. Then $C_1<C_2<C_3$ and
$C_1^*<C_2^*<C_3^*$. The following possibilities arise (see
proof of Lemma~\ref{LongLemma}):

\case[1]{We have $p\in \Gamma$ such that $\psi_{\alpha+}$ increases on $C_1$ from $p$ to $\Psi$ and 
$\psi^*_{\alpha+}$
increases on $C_1^*$ from $p$ to $\Psi^*$.}
By the proof of
Lemma~\ref{LongLemma} this gives $C_2\neq\emptyset$. Since $\delta$ is not in the
convex hull of $C_1$ in $C_1^*$, either $\delta>C_1$ or $\delta<C_1$.
\begin{enumerate}
\item[(a)] $\delta>C_1$. Taking any $\beta\in C_1$ we have 
$\alpha_1<\beta<\delta$, and thus $C_1<I<C_2$, with $I=I_\beta$ as defined above.
Let $\varepsilon=\frac{1}{2}|\beta-\alpha_1|\in \Gamma^{>}$ be as before, and
choose $\xi\in C_1$ so large that
$|\alpha_{n+1}-\psi_\alpha(\xi)|\leq\varepsilon$.
Hence, in $[\Gamma^*]$,
$$
\bigl[\psi^*_{\alpha+}(\delta)-\psi_{\alpha+}(\xi)\bigr]\ <\
 \bigl[\psi^*_\alpha(\delta)-\psi_\alpha(\xi)\bigr]\
\leq\ \bigl[\alpha_{n+1}-\psi_\alpha(\xi)\bigr]\ \leq\ [\varepsilon].
$$
Let $f(\gamma):=\theta^*(\gamma)+q_{n+1}\psi_{\alpha+}(\xi)$, for $\gamma\in I$. 
Then  $\theta^*_{+}(\delta)>f(\delta-\varepsilon)$ since
\begin{align*}
\theta^*_{+}(\delta)-f(\delta-\varepsilon)\ &=\ 
  \varepsilon+ \theta^*_{+}(\delta-\varepsilon)-f(\delta-\varepsilon) \\
&=\  \varepsilon+ q_{n+1}\bigl(\psi^*_{\alpha+}(\delta)-\psi_{\alpha+}(\xi)\bigr).
\end{align*}
Likewise $\theta^*_{+}(\delta)<f(\delta+\varepsilon)$. Hence, by the intermediate
value property for $f$ on $I$ (Lemma~\ref{siivp}) we get $\gamma\in I$ with
$f(\gamma)=\theta^*_{+}(\delta)$. Since $I\cap C_1=\emptyset$, we have 
$\gamma\notin \Gamma$, so
$f(\gamma)\notin \Gamma$ by inductive hypothesis. Thus 
$\theta^*_{+}(\delta)\notin \Gamma$.
\item[(b)] $\delta<C_1$. Then $\psi^*_{\alpha+}(\gamma)=p$ for all $\gamma$ such that 
$\delta\leq \gamma <C_1$. 
In particular
$\theta^*_{+}(\delta)=\theta^*(\delta)+q_{n+1}p\notin \Gamma$, 
by the induction hypothesis.
\end{enumerate}

\case[2]{We have $p,q\in \Gamma$ such that  $\psi_{\alpha+}$ increases on $C_1$ from $p$ to $q$ and $\psi^*_{\alpha+}$
increases on $C_1^*$ from $p$ to $q$.} Again,
either $\delta<C_1$
or $\delta>C_1$. Both subcases are treated
as in Case~1(b).

\medskip
\noindent
Next we consider the case that $C$ is bounded, 
$\psi_\alpha$ increases on $C$ and $\psi^*_\alpha$ increases on~$C^*$, 
and $i=3$. Either $\psi_{\alpha+}$ decreases on $C_3$ from $\Psi$ to
 $q$ and
$\psi^*_{\alpha+}$ decreases on~$C_3^*$ from $\Psi^*$ to $q$, for some 
$q\in \Gamma$, or
$\psi_{\alpha+}$ decreases on $C_3$ from $p$ to $q$ and~$\psi^*_{\alpha+}$
decreases on $C_3^*$ from $p$ to $q$, for some $p,q\in \Gamma$. The latter 
subcase is
treated as in Case~2 above. In the first
subcase, suppose that $\delta<C_3$. Then $C_2\neq\emptyset$, hence there exists~$\beta\in \Gamma$ with $\alpha_1<\beta<\delta$, and thus $C_2<I_{\beta}<C_3$, with $I_{\beta}$ as defined previously.
Now for any
$\varepsilon\in \Gamma^{>}$, in particular for 
$\varepsilon=\frac{1}{2}(\beta-\alpha_1)$, we can choose $\xi\in C_3$ 
such that
$|\alpha_{n+1}-\psi_\alpha(\xi)|\leq\varepsilon$. Now continue as in Case~1(a) above. If~${\delta>C_3}$, argue as in Case~1(b).
The case that $C$ is bounded and $\psi_\alpha$ is decreasing on $C$
can be handled in a similar way, and
is
left to the reader.

\medskip
\noindent
Now assume $C$ is unbounded and $i=1$. Then $\psi_{\alpha}$ decreases on $C$ by Lem\-ma~\ref{LongLemma}(iv), so $C_3 < C_2 < C_1$, and
$C_1$ is necessarily the unbounded component~${>\alpha_1}$ of $D_{\alpha+}$. We have the following cases:

\case[3]{$\psi_{\alpha+}$ decreases on $C_1$ from $p$ to $-\infty$ and
$\psi^*_{\alpha+}$ decreases on $C_1^*$ from $p$ to~$-\infty$, for some 
$p\in \Gamma$.}
Again, either $\delta<C_1$ or $\delta>C_1$. If $\delta<C_1$, proceed as in 
Case~1(b) above; if  $\delta>C_1$, then $\theta^*_{+}(\delta)\notin \Gamma$ follows from 
Lemmas~\ref{caca} and~\ref{siivp}.

\case[4]{$\psi_{\alpha+}$ decreases on $C_1$ from $\Psi$ to $-\infty$ and
$\psi^*_{\alpha+}$ decreases on $C_1^*$ from $\Psi^*$ to~$-\infty$.} 
If $\delta<C_1$, then inspection of the proof of Lemma~\ref{LongLemma} gives
$C_2\neq\emptyset$.
Hence there exists $\beta\in \Gamma$ with $\alpha_1<\beta<\delta$. Now adopt the argument in Case~(1)(a) above.
If $\delta>C_1$, we again apply Lemma~\ref{caca}.

\medskip
\noindent
Finally, consider the case that $C$ is unbounded and $i=3$. Then 
$\psi_{\alpha+}$ increases on~$C_3$ from $p$ to $\Psi$ and 
$\psi^*_{\alpha+}$
increases on $C_3$ from $p$ to $\Psi^*$, for some $p\in \Gamma$.
If $\delta>C_3$, note that any $\beta\in C_3$ will satisfy $\alpha_1<\beta<\delta$, and continue as in
Case~1(a). If $\delta<C_3$, argue as in Case~1(b).
This finishes the induction.
\end{proof}

\subsection*{Notes and comments} Proposition~\ref{cac}
is essentially Property~(B) on p.~333 of~\cite{AvdD}, proved there on
pp.~336--342.  
The $H$-asymptotic couples considered there are equipped 
with extra structure, but this
can be dropped, as pointed out in Section~6 of that paper.

%% file: mt-10.tex
\chapter{$H$-Fields}\label{ch:H}

\noindent
Valued differential subfields of differential-valued fields are not always differential-valued, as shown by an example at the end of Section~\ref{As-Fields,As-Couples}. They do satisfy an axiom
that defines the notion of a \textit{pre-differential-valued \textup{(}pre-$\d$-valued\textup{)} field}.\/ In Section~\ref{sec:pdv} we upgrade some basic facts on asymptotic fields to pre-$\d$-valued fields; for example, algebraic extensions of pre-$\d$-valued fields are 
pre-$\d$-valued, not just asymptotic. In Section~\ref{sec:integrals} we adjoin integrals (solutions $y=\int f$ of equations $y'=f$) to pre-$\d$-valued fields of $H$-type;
the expression~$\int f$ here is purely suggestive and 
we attach no formal meaning to it.
This is used in Section~\ref{sec:dv(K)} to show that every pre-$\d$-valued field of $H$-type has a canonical $\d$-valued extension. In Section~\ref{sec:exp integrals} we adjoin exponential integrals (solutions $y=\exp(\int g)$ of equations 
$y^\dagger =g$) to pre-$\d$-valued fields of $H$-type; again, the use of $\exp$ here is only suggestive.

Hardy fields are pre-$\d$-valued
fields, but also have a field ordering that interacts with the valuation and derivation. Axiomatizing this interaction yields the notion of a
\textit{pre-$H$-field}\/; \textit{$H$-fields\/} are $\d$-valued pre-$H$-fields. Our main goal (only reached in Chapter~\ref{ch:QE}) is the model theory of the particular $H$-field~$\T$, but this requires considering $H$-fields in general and their ordered valued differential subfields, the pre-$H$-fields. We begin their
study in Section~\ref{sec:H-fields}, showing among other things that each pre-$H$-field has a canonical $H$-field extension. Applying and adapting earlier adjunction results we show in Section~\ref{sec:Liouville closed} that every $H$-field can be extended  in a minimal way to one that is \textit{Liouville closed,}\/ that is, real closed and closed under integration and exponential integration; there are at most two such
minimal Liouville closed extensions, up-to-isomorphism. We finish this chapter
with some miscellaneous facts about asymptotic fields in
Section~\ref{Miscellaneous Facts about Asymptotic Fields}.

\section{Pre-Differential-Valued Fields}\label{sec:pdv}

\noindent 
{\em Throughout this section $K$ is a valued differential field}. We say that $K$ is {\bf pre-dif\-fe\-ren\-tial-valued\/} \label{p:pdv}\index{valued differential field!pre-differential-valued}\index{pre-differential-valued field}\index{field!pre-differential-valued} (for short: {\bf pre-$\d$-valued}) if the following holds:

\begin{list}{*}{\setlength\leftmargin{3em}}
\item[(PDV)] for all $f, g\in K^\times$, if $f\preceq 1$, $g\prec 1$, then $f'\prec g^\dagger$.
\end{list}
Any differential field with the trivial valuation is pre-$\d$-valued. If $K$ is pre-$\d$-valued, then so is
$K^{\phi}$ for $\phi\in K^\times$.
If $K$ is $\d$-valued, then $K$ is pre-$\d$-valued. Any valued differential subfield of a pre-$\d$-valued field is pre-$\d$-valued.

\begin{samepage}
\begin{lemma}\label{lem:asymp-predif}
The following conditions on $K$ are equivalent: \begin{enumerate}
\item[\textup{(i)}] $K$ is pre-$\d$-valued;
\item[\textup{(ii)}] for all $f,h\in K^\times$, if $f\preceq 1\nasymp h$, then  $f'\prec h^\dagger$;
\item [\textup{(iii)}] $K$ is asymptotic  and for all $f,g\in K^\times$, if
$f\asymp 1\nasymp g$, then $f' \prec g^\dagger$.
\end{enumerate}
\end{lemma}
\end{samepage}

\begin{proof}
The equivalence (i)~$\Leftrightarrow$~(ii) is clear. 
%(If $h\succ 1$ take $g=1/h$ and use $g^\dagger=-h^\dagger$.) 
The equivalence (ii)~$\Leftrightarrow$~(iii) follows from the equivalence of (i) and (iv) in Proposition~\ref{CharacterizationAsymptoticFields}.
%using in addition that if $K$ is asymptotic then $C\subseteq\mathcal O$.
\end{proof}

\noindent
From this lemma and $u^\dagger\asymp u'$ for $u\asymp 1$ in $K$ we
obtain:
 
\begin{cor}\label{cor:pdv dagger}
If $K$ is pre-$\d$-valued and $f,g\in K^\times$, $f\asymp g\not\asymp 1$, then $f^\dagger\sim g^\dagger$.
\end{cor}

\noindent
From Proposition~\ref{CharacterizationAsymptoticFields} and Lemma~\ref{lem:asymp-predif} we conclude:

\begin{cor}\label{asymp-predif} Ungrounded asymptotic fields are pre-$\d$-valued.
\end{cor}

\noindent
See Section~\ref{Miscellaneous Facts about Asymptotic Fields} for examples of asymptotic fields that are not pre-$\d$-valued.
The next lemma characterizes
pre-$\d$-valued fields as exactly the valued differential fields whose
valuation is trivial on the constant field and that obey a
valuation-theoretic version of l'Hospital's Rule:

\begin{lemma}\label{lem:char pdv}
The following conditions on $K$ are equivalent:
\begin{enumerate}
\item[\textup{(i)}] $K$ is pre-$\d$-valued;
\item[\textup{(ii)}] $C\subseteq\mathcal O$, and for all $f,g\in K^\times$, if $f\preceq g\prec 1$, then $\frac{f}{g}-\frac{f'}{g'}\prec 1$;
\item[\textup{(iii)}] $C\subseteq\mathcal O$, and for all $f,g\in K^\times$,
if $1\prec f\preceq g$, then $\frac{f}{g}-\frac{f'}{g'}\prec 1$.
\end{enumerate}
\end{lemma}
\begin{proof}
Suppose $K$ is pre-$\d$-valued. Then $K$ is asymptotic by Lemma~\ref{lem:asymp-predif}. Hence $C\subseteq\mathcal O$. Let $f,g\in K^\times$ and $f\preceq g\prec 1$.
Then $g'\ne 0$, and
$$\frac{f}{g}-\frac{f'}{g'}\ =\ \frac{fg' - f'g}{gg'}\ =\ -\frac{(f/g)'}{g^\dagger}.$$
As $f/g\preceq 1$ and $g\prec 1$, this gives $(f/g)'\prec g^\dagger$, and
we have shown (i)~$\Rightarrow$~(ii).
To prove (ii)~$\Rightarrow$~(iii), suppose (ii) holds, and let $f,g\in K$ with $1\prec f\preceq g$
and set $a:=1/g$, $b:=1/f$. Then $a\preceq b\prec 1$, hence $a^2/b \preceq a\prec 1$. Also,
$$\frac{f}{g}-\frac{f'}{g'}\ =\ \frac{1/b}{1/a}-\frac{(1/b)'}{(1/a)'}\ =\ \frac{a}{b}-\frac{a^2b'}{b^2a'}\ =\
\frac{(a^2/b)'}{a'}-\frac{a^2/b}{a},$$
which is infinitesimal by (ii).
Finally, we show  (iii)~$\Rightarrow$~(i).
Suppose (iii) holds, and let $f,h\in K^\times$ with $f\preceq 1\nasymp h$.
To get $f'\prec h^\dagger$, we may replace $f$ by $f+1$ and $h$ by $1/h$ to arrange
$f\asymp 1\prec h$. Then $1\prec fh\asymp h$, so
$$f'\frac{h}{h'}\ =\ \frac{(fh)'}{h'}- \frac{fh}{h}\ \prec\ 1$$
by (iii), and therefore $f'\prec h^\dagger$. Hence $K$ is pre-$\d$-valued by Lemma~\ref{lem:asymp-predif}.
\end{proof}

\begin{example}
Let $\k$ be a differential field, and consider the Hahn field $\k\(( t^\Q\)) $.
Let $g\in \k\(( t^{\Q}\)) $ and make
$\k\(( t^{\Q}\)) $ into a differential field extension of $\k$ by
$$\left(\sum a_{q}t^{q}\right)'\ =\ \sum a_{q}'t^q + \left(\sum qa_q t^{q-1}\right)g\qquad\text{(so $t'=g$.)}$$
Suppose $vg< 1$. Then $\k\(( t^\Q\)) $ is pre-$\d$-valued: just note
that for $f\in \k\(( t^\Q\)) ^\times$ we have $v(f')=vf+vg-1$ if $vf\neq 0$, and $v(f')>vg-1$ if $vf = 0$.
\end{example}

\subsection*{Pre-$\d$-valued fields by coarsening} 
\textit{In this subsection $K$ is an asymptotic field,  and $\Delta$ is a convex subgroup of
$\Gamma=v(K^{\times})$.}\/ With notation as at the end of Section~\ref{AbstractAsymptoticCouples},
we have the $\Delta$-coarsening~$(K,{\dotpreceq})$ of $K$, which is again an asymptotic field.
The asymptotic fields that arise naturally and
are not pre-$\d$-valued are
obtained by coarsening the valuation of a
pre-$\d$-valued field. This raises the question whether
every asymptotic field arises by coarsening the valuation of a
pre-$\d$-valued field. In Section~\ref{Miscellaneous Facts about Asymptotic Fields} we show that this is
not the case. The next result answers another question:
when is a coarsening of an asymptotic field pre-$\d$-valued?

\begin{prop}\label{Coarsenings-Pdv}
Suppose that $\{0\} \ne \Delta\neq\Gamma$. The following are equivalent:
\begin{enumerate}
\item[\textup{(i)}] $(K,\dotpreceq)$ is a pre-$\d$-valued field;
\item[\textup{(ii)}] the subset $\Psi_{\Delta}$ of $\dot{\Gamma}$ does not have a largest element;
\item[\textup{(iii)}] $\psi(\gamma)-\psi(\delta)\notin\Delta$ for all $\gamma\in\Gamma
\setminus\Delta$, $\delta\in\Delta^{\ne}$;
\item[\textup{(iv)}] $\psi(\gamma)-\psi(\delta) < \Delta$ for all $\gamma\in\Gamma
\setminus\Delta$, $\delta\in\Delta^{\ne}$.
\end{enumerate}
\end{prop}
\begin{proof}
The implication (iv)~$\Rightarrow$~(iii) is trivial,
(iii)~$\Rightarrow$~(ii) follows by Lemma~\ref{Coarsening-Prop}, and
(ii)~$\Rightarrow$~(i) by
Corollary~\ref{asymp-predif}.
To show (i)~$\Rightarrow$~(iv), let
$(K,\dotpreceq)$ be pre-$\d$-valued, and let
$\gamma\in\Gamma$, $\gamma >\Delta$, $\delta\in\Delta^{\ne}$.
Then $\gamma=vg$, $\delta=
vf$ with $f,g\in K^\times$, so
$\dot{v}(f') > \dot{v}( g^\dag)$, hence $\Delta>v(g^\dagger)-v(f')=
\psi(\gamma)-\bigl(\delta+\psi(\delta)\bigr)$. Thus
$\psi(\gamma)-\psi(\delta) < \Delta$.
\end{proof}

\noindent
In view of Lemma~\ref{asymp-lemma1} this yields:

\begin{cor}\label{Coarsenings-Pdv, application}
Suppose $\psi(\Delta^{\neq})\cap\Delta\neq\emptyset$. Then
$$\text{$(K,{\dotpreceq})$ is pre-$\d$-valued} \quad\Longleftrightarrow\quad
\psi(\Gamma\setminus\Delta)\ \subseteq\ \Gamma\setminus\Delta.$$
\end{cor}

\begin{exampleNumbered}\label{ex:pdv with gap, 2} Assume $K$ is of
$H$-type and $1\in\Gamma^>$ satisfies $\psi(1)=1$. Set
$$\Delta\ :=\ \big\{\gamma\in\Gamma :\ \text{$\psi^n(\gamma) \geq 0$
for some $n\geq 1$} \big\}.$$
Then $\Delta$ is a convex subgroup of $\Gamma$ with
$1\in\psi(\Delta^{\neq})\cap\Delta$, and~${\psi(\Gamma\setminus\Delta)\subseteq\Gamma\setminus\Delta}$.
There\-fore~${(K,{\dotpreceq})}$ is pre-$\d$-valued by
Corollary~\ref{Coarsenings-Pdv, application}. Moreover, any
$b\in K^\times$ with ${1\nasymp b\dotasymp 1}$ yields a gap
$\dot v(b')=\dot v(b^\dagger)=\psi\big(v(b)\big)+\Delta$ of
$(K,{\dotpreceq})$, by Lemma~\ref{Coarsening-Prop} and the equivalence
of (i) and (ii) in Proposition~\ref{Coarsenings-Pdv}.
Since $\Delta\ne \{0\}$, such $b$ exist. 

If $K$ is an  $\aleph_0$-saturated elementary extension  of $\T$,  then
$\Delta\neq\Gamma$ for this $\Delta$, since for each $n\geq 1$ and
$\gamma:= v\big(\!\exp^{n+1}(x)\big)$ we have $\psi^m(\gamma)<0$
for $m=1,\dots,n$.
\end{exampleNumbered}

\subsection*{Specialization} {\em In this subsection $K$ is an asymptotic field with small derivation and $\Delta$ is a convex subgroup of $\Gamma=v(K^\times)$ with $\psi(\Delta^{\neq})\subseteq\Delta$.}\/ We have the residue 
field $\dot K= \dot{\mathcal{O}}/\dot{\smallo}$
of the valuation ring $\dot{\mathcal{O}}$ of $(K,{\dotpreceq})$.
Recall from Section~\ref{Asymptotic-Fields-Basic-Facts} that 
$\dot K$
  with its induced valuation and derivation is an asymptotic field with small derivation having
asymptotic couple~${\big(\Delta,\psi|\Delta^{\ne}\big)}$. 
The residue map $a\mapsto \dot a\colon \dot{\mathcal{O}}\to\dot K$ restricts to a surjective differential ring morphism $\mathcal O\to\mathcal O_{\dot K}$, which in turn restricts to 
a field embedding $C\to C_{\dot K}$, and we identify $C$ with a subfield of $C_{\dot K}$ via this embedding. In addition we have a differential field
isomorphism $\res(K)\to\res(\dot K)$,  making the diagram
$$\xymatrix@R-0.5em@C-0.3em{ C\ \ar@{^{(}->}[r] \ar[d] & \mathcal O \ar[r] \ar[d]  & \res(K) \ar[d]^{\cong} \\
 C_{\dot K}\,  \ar@{^{(}->}[r] & \mathcal O_{\dot K} \ar[r] & \res(\dot K)
}$$
%$$\xymatrix@R-0.5em@C-0.3em{ C_{\dot K}\,  \ar@{^{(}->}[r] & \mathcal O_{\dot K} \ar[r] & \res(\dot K) \\
%C\ \ar@{^{(}->}[r] \ar[u] & \mathcal O \ar[r] \ar[u]  & \res(K) \ar[u]_{\cong}
%}$$
of differential ring morphisms commutative. 

\begin{lemma}\label{Specialization of asymptotic fields}
Suppose $K$ is pre-$\d$-valued. Then $\dot K$ is
pre-$\d$-valued.
If $K$ is $\d$-valued, then $\dot K$ is $\d$-valued
with $C_{\dot K}=C$.
\end{lemma}
\begin{proof}
Let $f\in \mathcal{O}$ and
$g\in \smallo\setminus \dot{\smallo}$, so
$\dot f\in {\mathcal O}_{\dot K}$ and
$\dot g \in \smallo_{\dot K}\setminus\{0\}$. Then $vg\in \Delta^{\ne}$, so
$vg'\in \Delta$, hence
$v(\dot {f}') \geq v(f')>v(g^\dagger)=
v({\dot g}^\dagger)$, so
$\dot K$ is pre-$\d$-valued. Now assume that
$K$ is $\d$-valued.
Then the composition $C\hookrightarrow\mathcal O\to\res(K)$ is an isomorphism, hence the composition   $C_{\dot K}\hookrightarrow {\mathcal O}_{\dot K}\to \res(\dot K)$
is also an isomorphism. Therefore $\dot K$ is
$\d$-valued with $C_{\dot K}=C$.
\end{proof}

\subsection*{The case $\sup \Psi=0$}
\textit{In this subsection $K$ is an asymptotic field and 
$(\Gamma,\psi)$, $\mathcal O$, $\smallo$, and 
$\k=\mathcal O/\smallo$ have the usual meaning.}\/ We complement here the remarks on small derivation in Section~\ref{Asymptotic-Fields-Basic-Facts}, focusing
on the case $\sup \Psi=0\notin \Psi$, where necessarily $K$ has small derivation, is ungrounded, and thus pre-$\d$-valued.

\begin{lemma}\label{smallunramasymp}
Suppose $\sup\Psi=0\notin\Psi$, and
$L$ is a valued differential field extension of $K$ with small derivation and $\Gamma_L=\Gamma$.
Then $L$ is pre-$\d$-valued.
\end{lemma}
\begin{proof}
Since $L$ has small derivation, we have $f'\preceq 1$ for all $f\preceq 1$
in $L$. So it is enough to show that $g^\dagger\succ 1$ whenever $g\in L^\times$ and $g\not\asymp 1$.
For such~$g$, take $a\in K^\times$ and $u\in L^\times$ with $g=au$ and $u\asymp 1$. Then $g^\dagger=a^\dagger+u^\dagger$ with $u^\dagger\asymp u'\preceq 1 \prec a^\dagger$ and hence $g^\dagger\sim a^\dagger\succ 1$.
\end{proof}

\noindent
Note:
$\sup \Psi=0\notin \Psi \Longleftrightarrow 0\in \Gamma \text{ is a gap in $K$.}$
Thus to get a $K$ with $\Gamma\ne \{0\}$ and $\sup \Psi=0\notin \Psi$ we can
start with any asymptotic field with a gap and nontrivial valuation---Example~\ref{ex:pdv with gap, 2}  yields pre-$\d$-valued fields of $H$-type with this
property---and then arrange by compositional conjugation that this gap is $0$. 

\begin{exampleNumbered}\label{ex:pdv with gap, 1}
Suppose $\sup\Psi=0\notin\Psi$. Let $a$ in a differential field extension of $K$ be transcendental over $K$ with $a'=1$; equip~$K(a)$ with the valuation in
Theorem~\ref{resext} for $F:=Y'-1$. Then by Lemma~\ref{smallunramasymp} the valued differential field~$K(a)$ is pre-$\d$-valued with value group $\Gamma$
and $a\asymp 1$. 
\end{exampleNumbered}

\noindent
In view of Corollary~\ref{cor:resext} we obtain from
Lemma~\ref{smallunramasymp}:

\begin{cor}\label{resextpred}
If $\sup\Psi=0\notin\Psi$ and
$\k_L$ is a differential field extension of~$\k$, then $K$ has a
pre-$\d$-valued field extension $L$ with small derivation, the same
value group as $K$, and differential residue field isomorphic to $\k_L$ over $\k$.
\end{cor}

\noindent
The assumption $\sup \Psi=0\notin \Psi$ in Lemma~\ref{smallunramasymp} may seem overly restrictive, but is in fact appropriate in view of the next two lemmas.

\begin{lemma}\label{lem:K<Y> pdv} Suppose $\der\smallo\subseteq \smallo$ and $K\<Y\>$ with its gaussian valuation is pre-$\d$-valued. Then $\sup\Psi=0\notin\Psi$.
\end{lemma}
\begin{proof} 
By Lemma~\ref{lem:K<Y> asymptotic} we have $\sup\Psi=0$.
Also, $f'\prec 1\asymp Y'\prec g^\dagger$ for all $f,g\in K^\times$ with $f\prec 1$ and $g\not\asymp 1$, and thus $0\notin \Psi$.
\end{proof}

\begin{lemma} Assume $\der\smallo\subseteq \smallo$.
Let $K\<a\>$ be as in Theorem~\ref{resext} 
with minimal annihilator $F$ of $a$ over $K$ of order $r\geq 1$, and
$\overline{F}\notin \k^\times  Y'$. 
If $K\<a\>$ is pre-$\d$-valued, then $\sup\Psi=0\notin\Psi$.
\end{lemma}
\begin{proof}
By Lemma~\ref{lem:K<a>asymp} and its proof we have $\sup\Psi=0$ and $a'\asymp 1$. Now argue as in the proof of Lemma~\ref{lem:K<Y> pdv} with $a$ replacing $Y$.
\end{proof} 

\noindent
In Section~\ref{Asymptotic-Fields-Basic-Facts} we already mentioned the significant constraint $\sup \Psi=0$
on $\d$-hen\-sel\-ian~$K$. This condition is necessary for
$K$ to be merely {\em embeddable\/} into a
$\d$-henselian asymptotic field; is it also sufficient?
Here is a partial answer.

\begin{cor}\label{aspsidh} Suppose $\sup \Psi=0\notin \Psi$.
Then $K$ has a $\d$-henselian pre-$\d$-valued field extension $L$ with $\Gamma_L=\Gamma$.
\end{cor}

{\sloppy

\begin{proof} Since $\k$ can be extended to a
linearly surjective differential field, Corol\-lary~\ref{resextpred} yields a pre-$\d$-valued field extension $E$ of $K$ with small derivation, the same
value group as $K$, and linearly surjective differential residue field $\k_E\supseteq \k$. Then Corollary~\ref{cor2diftrzda} and Lemma~\ref{smallunramasymp} yield
an immediate spherically complete pre-$\d$-valued field extension $L$ of $E$; any such $L$ is $\d$-henselian.
\end{proof}
}

\noindent
If $\max \Psi=0$, then $K$ also has a $\d$-henselian asymptotic field extension with the same value group, by Lemma~\ref{KgroundedLasymp} together with facts about the {\em linear surjective closure of a differential field\/}, a notion to be developed in the next volume.

\subsection*{Extensions of pre-differential-valued fields I}
\textit{In this subsection we assume that $L$ is an asymptotic field extension of the pre-$\d$-valued field $K$.}\/

\begin{lemma}\label{lem:pdv extension}
Suppose $\Psi$ is cofinal in $\Psi_L$, and for all $f\in L^{\preceq 1}$, either $f'\preceq a'$ for some $a\in K^{\preceq 1}$, or $f'\preceq h'$ for some $h\in L^{\prec 1}$.
Then $L$ is pre-$\d$-valued. 
\end{lemma}
\begin{proof}
Let $f,g\in L^\times$, $f\preceq 1$, $g\nasymp 1$.
Take $b\in K^\times$ with $b\nasymp 1$ and $b^\dagger \preceq g^\dagger$. If $a\in K^{\preceq 1}$ and $f'\preceq a'$, then $a' \prec b^\dagger$, so $f'\prec g^\dagger$. If 
$h\in L^{\prec 1}$ and $f'\preceq h'$, then $h'\prec g^\dagger$, so
again $f'\prec g^\dagger$. Thus~$L$ satisfies~(PDV).
\end{proof}

\begin{cor}  \label{cor:pdv extension, 1}
Suppose $\Psi$ is cofinal in $\Psi_L$,  and $\res(L)=\res(K)$. Then $L$ is pre-$\d$-valued.
\end{cor}
\begin{proof}
Let $f\in L^{\preceq 1}$. 
Take $a\in K$ such that $h:=f-a\prec 1$. 
Then $a\preceq 1$ and $f'=a'+h'$, so $f'\preceq a'$ or $f'\preceq h'$. Thus $L$ is pre-$\d$-valued by Lemma~\ref{lem:pdv extension}.
\end{proof}

\begin{cor} \label{cor:pdv extension, 2}
If $L|K$ is immediate, then $L$ is pre-$\d$-valued.
\end{cor}

\noindent 
Corollary~\ref{completingasymp} and the previous corollary yield:

\begin{cor}
The completion $K^{\operatorname{c}}$ of $K$ is pre-$\d$-valued.
\end{cor}

\subsection*{Extensions of pre-differential-valued fields II}
\textit{In this subsection $K$ is a pre-$\d$-valued field and $L|K$ is a valued differential field extension.}\/ We first record
a pre-$\d$-valued version of Lemma~\ref{Lemma34} with the same proof:

\begin{lemma}\label{Lemma3.3} 
Suppose $\res(K) = \res(L)$,
$T\supseteq K^\times$ is a subgroup of
$L^{\times}$ with $L=K(T)$, each element of 
$K[T]^{\ne}$ has the
form $t_1+\dots+t_k$ with $k\ge 1$, $t_1,\dots,t_k\in T$ and
$t_1 \succ t_i$ for $2\le i \le k$, and for all 
$a,b\in T$,
$$ a\preceq 1,\ b\prec 1
\quad\Longrightarrow\quad a' \prec b^\dagger. \label{Eq3}
$$
Then $L$ is pre-$\d$-valued. 
\end{lemma}

\noindent
If $L|K$ is algebraic, then $L$ is asymptotic by Proposition~\ref{Algebraic-Extensions-Asymptotic}, and so $\Psi_L$ is defined and equals~$\Psi$. Thus by Corollary~\ref{cor:pdv extension, 1}: 

\begin{lemma}  \label{Purely ramified-pdv}
If $L|K$ is algebraic and $\res(K)=\res(L)$, then $L$ is pre-$\d$-valued.
\end{lemma}

\begin{lemma}\label{Residue-Field-Extension-pdv}
Suppose that $K$ is henselian and $L|K$ is of finite degree such that
$[L:K]=\bigl[\res(L):\res(K)\bigr]$. Then $L$ is pre-$\d$-valued.
\end{lemma}
\begin{proof}
The proof of Lemma~\ref{Residue-Field-Extension-Asymptotic} shows that for all $f\in L^{\preceq 1}$ we have $f'\preceq a'$ for some
$a\in K^{\preceq 1}$. Hence $L$ is pre-$\d$-valued by
Lemma~\ref{lem:pdv extension}.
\end{proof}

\begin{prop}\label{Algebraic-Extensions-pdv}
Suppose that $L|K$ is algebraic. Then $L$ is pre-$\d$-valued.
\end{prop}
\begin{proof} We first arrange that $L$ is an algebraic closure of $K$. By Lemma~\ref{Purely ramified-pdv} we next arrange that $K$ is henselian. 
Using~\ref{Residue-Field-Extension-pdv} and 
\ref{Purely ramified-pdv} we now
obtain that $L$ is pre-$\d$-valued as in the proof of Proposition~\ref{Algebraic-Extensions-Asymptotic}.
\end{proof}

\begin{cor}\label{cor:Algebraic-Extensions-pdv}
If $K$ is $\d$-valued, then its algebraic closure $K^\alg$, equipped with the unique extension of the 
derivation of $K$ to a derivation on $K^\alg$ and
any valuation on~$K^\alg$ extending that of $K$, is also $\d$-valued.
\end{cor}
\begin{proof}
The constant field of $K^\alg$ is an algebraic closure of $C$ (Lemmas~\ref{lem:C alg closed in K},~\ref{lem:const field of algext}), and the residue field of $K^\alg$ is an algebraic closure of $\res(K)$ (Corollary~\ref{cor:ef inequ, 2}).
\end{proof}

\begin{remark}
The end of Section~\ref{As-Fields,As-Couples} has an example of
a $\d$-valued Hardy field $\Q(x)$ with an algebraic Hardy field extension that is not $\d$-valued.
\end{remark}

\subsection*{Notes and comments}
Pre-$\d$-valued fields are implicit in~\cite[Theorem 1]{Rosenlicht2} and explicit in \cite{AvdD2}.
The statement of~\cite[Theorem 6]{Rosenlicht2} is not quite correct; its proof does give Corollary~\ref{cor:Algebraic-Extensions-pdv}. 
Proposition~\ref{Algebraic-Extensions-pdv} is from \cite{AvdD2}.

\section{Adjoining Integrals}\label{sec:integrals}

\noindent
We perform here mild variants of adjunctions done in Sections~4 and~5 of \cite{AvdD2}. In contrast to that paper we consider only the case of $H$-type; this 
avoids some case distinctions.
\textit{In this section we assume that $K$ is an $H$-asymptotic field.}\/

\begin{lemma}\label{extensionpdvfields1}
Suppose $K$ is pre-$\d$-valued and $vs$ is
a gap in $K$, with $s\in K$.
Then~$K$ has a pre-$\d$-valued extension $K(y)$ of $H$-type with $y'=s$, $y\prec 1$, and $\res K(y)=\res K$, such that 
for any asymptotic extension $L$ of $K$ and $z\in L$ with $z'=s$ and $z\prec 1$ there is a unique $K$-embedding $K(y) \to L$ sending
$y$ to $z$.
\end{lemma}
\begin{proof} Suppose $K(y)$ is an asymptotic extension of $K$
with $y'=s$ and $y\prec 1$.  
Then  $\alpha:= vy>0$ and $\alpha'=vs$, so $0 < n\alpha < \Gamma^{>}$ for all $n\ge 1$, by Lemma~\ref{extension1} and its proof. Hence $y$ is transcendental over $K$ and the asymptotic couple
of~$K(y)$ is just the $H$-asymptotic couple
$(\Gamma + \Z\alpha, \psi^{\alpha})$ with $\alpha'=\beta$ from 
Lemma~\ref{extension1}, with $\beta:=vs$. Thus $K(y)$ is $H$-asymptotic and has the desired universal property. %stated in the present lemma. 

To construct such an extension, take a valued differential field extension $K(y)$ of~$K$ with
$y$ transcendental over $K$, $y'=s$, and $0<n\alpha<\Gamma^{>}$ for all $n\geq 1$, where $\alpha:=vy$. Note that then $\res K  = \res K(y)$ by Lemma~\ref{lem:lift value group ext}.
It remains to show that $K(y)$ is pre-$\d$-valued. 
We have $\Gamma_{K(y)}:=v\big(K(y)^\times\big)=\Gamma + \Z \alpha=\Gamma \oplus \Z\alpha$.
We verify the conditions of Lemma~\ref{Lemma3.3} for $L=K(y)$ and
$T:=  K^{\times}y^{\Z}$. It is clear that  
every element of $K[T]^{\neq}$ has the form
$t_1 + \dots + t_n$ with $n \ge 1$, $t_1,\dots,t_n\in T$, and $t_1\succ \cdots \succ t_n$. 
Let $a\in T$; 
it clearly suffices to show for $\beta:= vs$: $a\preceq 1\Rightarrow v(a')\ge \beta$, and $a\prec 1\Rightarrow v(a^\dagger) < \beta$. If $a\asymp 1$, then $a\in K^\times$, so $v(a')\ge \beta$; accordingly we assume $a\prec 1$ below. With $a=fy^k$, $f\in K^\times$, we have
$$a'\ =\ f'y^k +kfy^{k-1}y'\ =\ f'y^k + kas/y, \qquad a^\dagger\ =\ f^\dagger + ky^\dagger.$$ 
Either $f\prec 1$, or $f\asymp 1$ and $k>0$. Consider first the
case $f\prec 1$. Then $v(f'y^k)=v(f')+k\alpha > \beta$ and
$v(kas/y)\ge v(a)+\beta-\alpha\ge \beta$, so $v(a')\ge \beta$;
also $v(f^\dagger) < \beta-\alpha$ and $v(y^\dagger)=\beta-\alpha$, so $v(a^\dagger)<\beta$. Next assume $f\asymp 1$ and $k>0$. 
Then we get $v(a')\ge \beta$ as before, and $v(f^\dagger)=v(f')\ge \beta$ and so $v(a^\dagger)=\beta-\alpha< \beta$.
\end{proof}

\begin{remarks}
Let  $K$ and $K(y)$ be as in Lemma~\ref{extensionpdvfields1}. Then
$K(y)$ is grounded by \eqref{eq:Psialpha}. 
If $K$ is $\d$-valued, then so is $K(y)$ with $C_{K(y)}=C$, by Lemma~\ref{dv}.
\end{remarks}

\noindent
What if we replace $y\prec 1$ in 
Lemma~\ref{extensionpdvfields1} by $y\succ 1$? To get a concise
answer we restrict ourselves to $\d$-valued $K$:

\begin{lemma}\label{variant41}
Let $K$ be $\d$-valued, and let $vs$ be a gap in $K$, with $s\in K$.
Then $K$ has a $\d$-valued extension $K(y)$ of $H$-type with
$y'=s$ and 
$y\succ 1$ such that 
for any asymptotic extension $L$ of $K$ and $z\in L$ with
$z'=s$ and $z\succ 1$ there is a unique 
$K$-embedding $K(y) \to L$ sending $y$ to $z$.
\end{lemma}
\begin{proof}
Similar to that of Lemma~\ref{extensionpdvfields1}, using the remarks following~\eqref{eq:Psialpha} on the variant of Lemma~\ref{extension1} when $\alpha<0$.
With $\alpha$,~$\beta$,~$T$ as in the proof of Lemma~\ref{extensionpdvfields1}
and $a\in T$, one shows: $a\preceq 1\Rightarrow v(a')> \beta-\alpha$, and $a\prec 1\Rightarrow v(a^\dagger) \le \beta-\alpha$;
this uses the $\d$-valued assumption on $K$.
\end{proof}

\begin{remark}
In Lemma~\ref{variant41} we have $C_{K(y)}=C$ and $K(y)$ is grounded, by~\eqref{eq:Psialpha}. 
\end{remark}

%\noindent
%Next we adjoin an integral for 
%an element $s$ with $vs=\max \Psi$:

\begin{lemma}\label{pre-extas2}
Let $K$ be pre-$\d$-valued, and let $s\in K$ be such that $vs = \max \Psi$. Then $K$ has a pre-$\d$-valued extension $K(y)$ of $H$-type with $y'=s$ such that 
for any pre-$\d$-valued extension $L$ of $K$ and $z\in L$ with
$z'=s$ there is a unique $K$-embedding $K(y) \to L$ sending
$y$ to $z$.
\end{lemma}
\begin{proof} Suppose $K(y)$ is a pre-$\d$-valued extension of $K$
such that~$y'=s$. 
The assumption that $K(y)$ is 
pre-$\d$-valued gives $\alpha:=vy< 0$, and as
at the beginning of the proof of Lemma~\ref{extas2} (working in the divisible hull of the asymptotic couple of~$K(y)$) we see that $\Gamma^{<} < n\alpha < 0$ for all 
$n\ge 1$.
Thus~$y$ is transcendental over~$K$, and the asymptotic couple
of~$K(y)$ is the $H$-asymptotic couple $(\Gamma + \Z\alpha, \psi^{\alpha})$ with
$\alpha'=\beta:=vs$ from Lemma~\ref{extas2}. This makes it clear why~$K(y)$ has the universal property
stated in the lemma, and it also indicates a way to construct
such an extension~$K(y)$ along the lines of the proof of Lemma~\ref{extensionpdvfields1}.   
\end{proof}

\noindent
Let $K$, $s$, and $K(y)$ be as in Lemma~\ref{pre-extas2} and set $\alpha:= vy$. In view of the remarks following the proof of Lemma~\ref{extas2} we have 
$\Psi_{K(y)}=\Psi\cup\{\alpha^\dagger\}$ with 
$\Psi < \alpha^\dagger$. By Lemma~\ref{lem:lift value group ext} we have 
$\res  K(y)=\res K$, and by Corollary~\ref{cor:no new exps under integration}
and $s\notin\der K$ we have $C_{K(y)}=C$. 
If $K$ is $\d$-valued, then so is $K(y)$, by Lemma~\ref{dv}.

\begin{lemma}\label{var43} Suppose $K$ is henselian, 
$s\in K$, $vs\in (\Gamma^{>})'$, 
$s\notin \der\smallo$, 
and $$S\ :=\ \big\{v(s-a'):\ a\in \smallo\big\}$$ has no largest element. 
Let $L=K(y)$ be a field extension of
$K$ with $y$ transcendental over $K$, and let $L$ be equipped with the 
unique derivation extending the derivation of $K$ such that $y'=s$. Then there is a unique valuation of $L$ that makes it an 
$H$-asymptotic
extension of~$K$ with $y\nasymp 1$. With this valuation $L$ is an
immediate extension of $K$ with $y\prec 1$, and so 
$L$ is pre-$\d$-valued if $K$ is. 
\end{lemma}
\begin{proof} 
Let $\kappa = \operatorname{cf}(S)$, so $\kappa$ is an infinite cardinal, and let $\rho$, $\sigma$, $\tau$ range over or\-di\-nals~${<\kappa}$.
Let $(a_\rho)$
be a sequence in $\smallo$ such that $\bigl(v(s-a_\rho')\bigr)$ is
a strictly increasing sequence in $S$, and cofinal in $S$.
Then $s-a_\rho' \asymp (a_\sigma-a_\rho)'$
for $\sigma>\rho$, hence
${(a_\tau-a_\sigma)'\prec (a_\sigma-a_\rho)'}$ and so 
$a_\tau-a_\sigma \prec a_\sigma-a_\rho$
for $\tau>\sigma>\rho$; thus $(a_\rho)$ is a 
pc-sequence in $K$.
Also, $(a_\rho)$ has no 
pseudolimit in $K$: suppose $a_\rho\leadsto a\in K$; then
$a-a_\rho \asymp a_\sigma-a_\rho \prec 1$
for $\sigma>\rho$, so $a \in \smallo$;
moreover, $a' - a_\rho' \asymp a_\sigma'-a_\rho'$
for $\sigma>\rho$, hence $v(s-a') > v(s-a_\rho')$
for all $\rho$, contradicting the cofinality property of the $v(s-a_{\rho}')$. 
With $P(Y):=Y'-s$ we have $P(a_\rho) \leadsto 0$, and since
 $K$ is henselian, we have
$Q(a_\rho)\not\leadsto 0$  for all $Q(Y)\in K[Y]^{\neq}$.
Hence the hypotheses of Proposition~\ref{pca-extension} are satisfied.  The first part of that proposition then yields a valuation of $L$ that makes it an immediate $H$-asymptotic extension of $K$ with $y\prec 1$. 
Let any valuation of $L$ be given that makes it an $H$-asymptotic extension of $K$ with $y\nasymp 1$. Then $y\prec 1$, since $y\succ 1$ implies $vs=v(y')< (\Gamma^>)'$. Also $y'=s$ gives that $\big( v(y'-a_\rho') \big)$ is strictly increasing, and so is $\big( v(y-a_\rho) \big)$, and thus $a_\rho\leadsto y$. It remains to use 
Proposition~\ref{pca-extension}, in particular part (ii).
\end{proof}

\begin{remarks} Suppose $K$ and $s$ are as in Lemma~\ref{var43} 
and $E$ is an asymptotic  field extension of $K$ with an element $z\prec 1$ such that $z'=s$. Then $z$ is transcendental over~$K$. This is because
$a_{\rho}' \leadsto z'$ for the pc-sequence $(a_\rho)$ in the proof of that lemma, so $a_\rho \leadsto z$; it remains to recall that $(a_\rho)$ is of transcendental type over $K$. 
\end{remarks}

\noindent
The assumption in Lemma~\ref{var43} that $\big\{v(s-a'):a\in \smallo\big\}$ has no maximum is always satisfied for pre-$\d$-valued $K$, by part (iii) of the next lemma.  

\begin{lemma}\label{lem:remark5.1(2)}
Suppose $K$ is pre-$\d$-valued and $s\in K$. Then: \begin{enumerate}
\item[\textup{(i)}] if $\big\{ v(s-a'):\  a\in K\big\}$ has a maximum
$\beta\in \Gamma_{\infty}$, then $\beta\notin (\Gamma^{\neq})'$;
\item[\textup{(ii)}] if $K$ has asymptotic integration and $s\notin \der K$, then $\big\{ v(s-a'):\  a\in K\big\}$ has no maximum;
\item[\textup{(iii)}] if $vs\in (\Gamma^{>})'$, and $s\notin \der\smallo$, then $\big\{v(s-a'):a\in \smallo\big\}$ has no maximum.
\end{enumerate}
\end{lemma}
\begin{proof} Suppose $\max\{v(s-a'): a\in K\}=\beta\in (\Gamma^{\neq})'$.
Take $a\in K$ and $b\in K^\times$, $b\nasymp 1$, such that 
$\beta=v(s-a')=v(b')$. 
Next, take $u\in K$ with $u\asymp 1$ and $s-a'=ub'$. Then 
$u'\prec b^\dagger$ as $K$ is pre-$\d$-valued, hence
$s-(a+ub)'=-u'b\prec b'$, contradicting the maximality of 
$\beta$. This proves (i), and (ii) is an immediate consequence.
The proof of (iii) is like that of (i), with $(\Gamma^{>})'$ and $a,b\in \smallo$ instead of $(\Gamma^{\ne})'$ and $a,b\in K$.
\end{proof}

%\begin{remarks} Lemma~\ref{lem:remark5.1(2)} shows that if $K$ is 
%$\d$-valued, $s\in K$, $vs\in (\Gamma^{>})'$, and $s\notin \der\smallo$, then $\{v(s-a'):a\in \smallo\}$ has no maximum; the latter is an assumption in Lemma~\ref{var43}. Likewise, if
%$K$ is pre-$\d$-valued with asymptotic integration, $s\in K$, and $s\notin \der K$,
%then $\{v(s-a'): a\in K\}$ has no maximum.

%Suppose $K$ and $s$ are as in Lemma~\ref{var43} 
%and $E$ is an asymptotic  field extension of $K$ with an element $z\prec 1$ such that $z'=s$. Then $z$ is transcendental over $K$. This is because
%$a_{\rho}' \leadsto z'$ for the pc-sequence $(a_\rho)$ in the proof of that lemma, so $a_\rho \leadsto z$; it remains to recall that $(a_\rho)$ is of transcendental type over $K$. 
%\end{remarks}

%\begin{remark} If $K$ is pre-$\d$-valued and $s\in K$ is such that $vs\in (\Gamma^{>})'$, $s\ne a'$ for all $a\in \smallo$, then the assumption in  Lemma~\ref{var43} that $S$ has no largest element is automatically satisfied.  To see this, note first that $vs\in S$. Let $\gamma\in S$ with $\gamma \ge vs$, say $\gamma = v(s-a')$ with $a \in \m$. Since $(\Gamma^>)'$ is upward closed, there exists $b\in \smallo$ with $v(b')=\gamma$. Thus for  some $u\in K$ with $u\asymp 1$ we have $v(s-a'-ub') > \gamma$, and then $v(u'b) > v(b')=\gamma$, so  $v\bigl(s-a' - (ub)'\bigr)>\gamma$.
%\end{remark}

%\noindent
%We shall also use a variant of the lemma above:

\begin{lemma}\label{var51} Suppose $K$ is henselian, $s\in K$,
$v(s-a')< (\Gamma^{>})'$ for all $a\in K$, 
and $S:= \big\{v(s-a'):\ a\in K\big\}$ 
has no largest element. Let $L=K(y)$ be a field extension of
$K$ with $y$ transcendental over $K$, and let $L$ be equipped with the 
unique derivation extending the derivation of~$K$ such that $y'=s$. Then there is a unique valuation of~$L$ making it an 
$H$-asymptotic extension of $K$. With this valuation $L$ is an
immediate extension of $K$ with~${y\succ 1}$, and so $L$ is pre-$\d$-valued if $K$ is. 
\end{lemma}
\begin{proof} Note that $vs\in S$. Since $S$ has no largest element, each 
$\alpha\in S$ satisfies
$\alpha < \gamma$ for some $\gamma\in \Psi$. Set 
$\kappa= \operatorname{cf}(S)$,
and let $\rho$, $\sigma$, $\tau$ range over ordinals~$<\kappa$.
Take a sequence $(a_{\rho})$ in $K$ such that $v(s-a_{\rho}')$ is strictly increasing and 
cofinal in~$S$ as a function of $\rho$, and $s-a_{\rho}'\prec s$ for all
$\rho$. 
Then $s-a_\rho'\sim a_\sigma'-a_\rho'$
for $\sigma > \rho$, and hence
$(a_\tau-a_\sigma)'\prec (a_\sigma-a_\rho)'$
for $\tau>\sigma>\rho$. 
Note that $a_\sigma-a_\rho \succ 1$ for $\sigma>\rho$, since otherwise
$s-a_\sigma'\sim (a_{\sigma +1}-a_{\sigma})'\prec 
(a_{\sigma}-a_\rho)'$, so $v(s-a_\sigma')\in (\Gamma^>)'$, contradicting the hypothesis.  Hence
$a_\tau-a_\sigma \prec a_\sigma-a_\rho$
for $\tau> \sigma> \rho$, so $(a_\rho)$ is a 
pc-sequence in $K$. From $a_\sigma-a_\rho\succ 1$ for $\sigma<\rho$ we get $a_{\rho} \preceq 1$ for at 
most one $\rho$, and so we can arrange that $a_{\rho} \succ 1$
for all $\rho$. If $a_\rho\leadsto a\in K$, then for some index $\rho_0$ we have
$a-a_\rho \asymp a_\sigma-a_\rho \succ 1$ for $\sigma > \rho > \rho_0$, and so $a' - a_\rho' \asymp a_\sigma'-a_\rho'$
for $\sigma>\rho > \rho_0$, hence $v(s-a') > v(s-a_\rho')$
for all $\rho$, contradicting the cofinality 
property of the $v(s-a_{\rho}')$. Thus the pc-sequence $(a_\rho)$ in $K$ is divergent. 

Now apply Proposition~\ref{pca-extension} to $P(Y):= Y'-s$ as in the proof of Lemma~\ref{var43} to get a valuation of $L$ making it an $H$-asymptotic extension of $K$. Assume any such
valuation is given. Then $y\succ 1$, since $y\preceq 1$ implies $vs=v(y')\geq\Psi$, so  
$v(s-a')\in (\Gamma^>)'$ for some $a\in K$, a contradiction. Also, $y'=s$ gives that $\big( v(y'-a_\rho') \big)$ is strictly increasing, and so
is $\big( v(y-a_\rho) \big)$, and thus $a_\rho\leadsto y$. 
It remains to use 
Proposition~\ref{pca-extension}, in particular part (ii). 
\end{proof} 

\noindent 
By alternating the passage to a henselization with the extension procedures 
of Lemmas~\ref{var43} and~\ref{var51} we obtain the following:

\begin{prop}\label{asintint} Let $K$ be $\d$-valued with 
asymptotic integration. Then $K$ has an immediate asymptotic extension
$K(\int)$ such that:  
\begin{enumerate}
\item[\textup{(i)}] $K(\int)$ is henselian and closed under integration;
\item[\textup{(ii)}] $K(\int)$ embeds over $K$ into any henselian $\d$-valued $H$-asymptotic extension of $K$ that is closed under integration.
\end{enumerate}
\end{prop}

\nomenclature[Z]{$K(\int)$}{closure of $K$ under integration}

\begin{proof} Define an {\em integration tower\/} on $K$ to be a strictly 
increasing chain $(K_{\lambda})_{\lambda\le \mu}$ of immediate asymptotic extensions
of $K$, indexed by the ordinals less than or equal to some ordinal $\mu$, such 
that:  $K_0 = K$; if $\lambda$ is a limit ordinal, $0<\lambda \le \mu$, then 
$K_\lambda = \bigcup_{\iota<\lambda} K_\iota$; for $\lambda < \lambda +1 \le \mu$, {\em either\/} $K_{\lambda + 1}$ is a henselization of $K_\lambda$,
{\em or\/} $K_\lambda$ is already henselian,
$K_{\lambda + 1} = K_\lambda(y_{\lambda})$ with $y_\lambda\notin K_\lambda$
(hence $y_\lambda$ is transcendental over~$K_\lambda$), $y_\lambda'= s_\lambda \in K_\lambda$, and (1) or (2) below
holds, where $(\Gamma_\lambda, \psi_\lambda)$
denotes the asymptotic
couple of $K_\lambda$: 
\begin{enumerate}
\item
$v(s_\lambda)\in (\Gamma_\lambda^{>})'$, and 
$s_\lambda\neq\varepsilon'$ for all 
$\varepsilon\in K_\lambda^{\prec 1}$;
\item $S_\lambda :=\bigl\{v(s_\lambda - a'): a\in K_\lambda\bigr\} <
(\Gamma_\lambda^{>})'$.
\end{enumerate}
Note that in (1) we can arrange $y_{\lambda} \prec 1$ by subtracting a 
constant from $y_{\lambda}$, and that in~(2) the set
$S_\lambda$ has no largest element and $y_{\lambda} \succ 1$. This is relevant
in connection with the uniqueness of the valuations in 
Lemmas~\ref{var43} and~\ref{var51}.

Take a maximal integration tower $(K_{\lambda})_{\lambda\le \mu}$ on $K$, where ``maximal'' means that it cannot be extended to an integration tower 
$(K_{\lambda})_{\lambda\le \mu+1}$ on $K$. Then the top~$K_{\mu}$ of the tower 
has the properties stated about $K(\int)$.  
\end{proof}

\begin{cor}\label{minint} With $K$ and $K(\int)$ as in 
Proposition~\ref{asintint}, the only henselian asymptotic subfield 
of $K(\int)$ containing $K$ and closed under integration
is $K(\int)$.
\end{cor} 

{\sloppy
\begin{proof} Let $F$ be the intersection of all
henselian asymptotic subfields $L\supseteq K$ of~$K(\int)$ that are closed
under integration. Then $F$ itself is among these $L$, and so~$F$ 
is the smallest
such $L$. Since $F$ is a
$\d$-valued field extension of $K$ closed under integration,
 $K(\int)$ embeds into $F$ over $K$, and the image of any such embedding is
also among these $L$, and thus equals $F$. In particular, $K(\int)$ is 
$K$-isomorphic to~$F$, and so inherits the minimality property of $F$.    
\end{proof}
}

\noindent
Suppose $K$ is $\d$-valued and has
asymptotic integration. The minimality property of 
Corollary~\ref{minint} yields that an asymptotic extension 
$K(\int)$ as in Proposition~\ref{asintint} is unique
up to isomorphism over $K$: if $L$ is also an asymptotic extension of $K$ with the properties of $K(\int)$
as in that proposition, then there is a $K$-embedding of $L$ into $K(\int)$, 
and the image of this embedding is then $K(\int)$.   

\medskip
\noindent
Here is an obvious consequence of Proposition~\ref{asintint}:

\begin{cor}\label{asintint, cor}
Spherically complete $\d$-valued fields of $H$-type with 
asymptotic integration are closed under integration.
\end{cor}

\subsection*{Notes and comments} Lemmas~\ref{extensionpdvfields1}--\ref{lem:remark5.1(2)}
are variants of~\cite[4.1, 4.2, 4.3, 5.1]{AvdD2}. That paper also has \ref{asintint}--\ref{asintint, cor} without $H$-type assumption. Another proof of~\ref{asintint, cor}, also without assuming $H$-type, is in~\cite{FVK-hensel}.
If $C$ is algebraically closed, then the differential  field~$K(\int)$ from Proposition~\ref{asintint} contains the Picard-Vessiot antiderivative closure of~$K$  constructed in~\cite{Magid}.

\section{The Differential-Valued Hull}\label{sec:dv(K)}

\noindent
Valued differential subfields of $\d$-valued fields are pre-$\d$-valued. We prove here a strong converse in the $H$-type case. 
\textit{Throughout this section $K$ is  a pre-$\d$-valued field of $H$-type.}\/

{\sloppy

\begin{theorem}\label{thm:dv(K)}
The pre-$\d$-valued field $K$ of $H$-type has a $\d$-valued ex\-ten\-sion~$\operatorname{dv}(K)$ of $H$-type such that any embedding of $K$ into any $\d$-valued field $L$ of $H$-type extends uniquely to an embedding of $\operatorname{dv}(K)$ 
into~$L$.
\end{theorem}
}

\index{hull!differential-valued}
\index{differential-valued!hull}
\index{d-valued@$\d$-valued!hull}
\nomenclature[Z]{$\operatorname{dv}(K)$}{differential-valued hull of  $K$}

\noindent
The universal property in the theorem determines $\operatorname{dv}(K)$ up to unique isomorphism over~$K$ of valued differential field extensions of $K$. We call $\operatorname{dv}(K)$ the {\bf differential-valued hull of $K$.}

\begin{proof}[Proof of Theorem~\ref{thm:dv(K)}]
If $K$ is not $\d$-valued, then for some $b\asymp 1$ in $K$ we have $b'\notin \der\smallo$,
%\neq\varepsilon'$ for all $\varepsilon\in\smallo$, 
and as~$K$ is of $H$-type, either $v(b')\in (\Gamma^>)'$ or
$v(b')<(\Gamma^>)'$ for such $b$, and in the latter case $K$ has a gap~$v(b')$.
For the purpose of this proof, call $K$ \textit{nice}\/ if  there is no $b\asymp 1$ in~$K$ such that $v(b')$ is a gap in $K$. Thus if $K$ has no gap (in particular, if~$K$ is grounded), then $K$ is nice.
Also, if $K$ is $\d$-valued, then $K$ is nice. If $K$ is nice, so is any immediate pre-$\d$-valued extension.
If $K$ is nice, put $K_0:=K$;
otherwise, take $b\asymp 1$ in $K$ such that $v(b')$ is a gap in $K$, put $s:=b'$, and take $K_0:=K(y)$ with $y\prec 1$ as in Lemma~\ref{extensionpdvfields1}.
Thus $\res(K_0)=\res(K)$ and $K_0$ is nice. 
Starting with $K_0$ and iterating and 
alternating applications of Proposition~\ref{Algebraic-Extensions-pdv} and
Lemma~\ref{var43} we obtain an immediate 
henselian $\d$-valued extension $K_1$ of $K_0$ such that 
any embedding of $K$ into any henselian $\d$-valued field $L$ 
of $H$-type extends
to an embedding of $K_1$ into~$L$; this also uses the remarks
following Lemma~\ref{var43}.

Let $D$ be the constant field of $K_1$; so $D$ maps isomorphically
onto $\res(K_1) =
\res(K)$ under the residue map ${\mathcal O}_{K_1}\to\res(K_1)$. Put
$\operatorname{dv}(K):=K(D)$, a $\d$-valued subfield of 
$K_1$. To show that $\operatorname{dv}(K)$ has the desired universal property,
let $i\colon K\to L$ be
any embedding of $K$ into a $\d$-valued field $L$ of $H$-type.  Extend $i$  
to an embedding $i_1\colon K_1\to L^{\operatorname{h}}$. Then $i_1(D)\subseteq C_{L^{\operatorname{h}}}=C_L\subseteq L$, so  
$i_1\big(\!\operatorname{dv}(K)\big)\subseteq L$. This gives an embedding
$i_1|\operatorname{dv}(K)\colon\operatorname{dv}(K)\to L$ that extends $i$.
Given $d\in D$, any such embedding must send $d$ to the element
of $C_L$ whose residue class in $\res(L)$ equals the natural $i$-image of 
$\res(d)\in \res(K)$ in $\res(L)$. This gives the required uniqueness.  
\end{proof}

\noindent
The proof of the theorem above provides extra information: 
\begin{enumerate}
\item[(a)] $\res K=\res \operatorname{dv}(K)$;
\item[(b)] $\operatorname{dv}(K)=K(D)$ where $D$ is the constant field of $\operatorname{dv}(K)$;
\item[(c)] $\operatorname{dv}(K)$ is
$\d$-algebraic over $K$, as a consequence of (b).
\end{enumerate}

\begin{cor}\label{cor:value group of dv(K)}
The value group of $\operatorname{dv}(K)$ is as follows:
\begin{enumerate}
\item[\textup{(i)}] Assume $K$ has no gap; then $\Gamma_{\operatorname{dv}(K)}=\Gamma$.
\item[\textup{(ii)}] Assume $K$ has a gap $\beta$ and no $b\asymp 1$ in $K$ satisfies
$v(b')=\beta$; then $\Gamma_{\operatorname{dv}(K)}=\Gamma$.
\item[\textup{(iii)}] Assume $K$ has a gap $\beta$ and $b\asymp 1$ in $K$ satisfies
$v(b')=\beta$. Let $a$ be the unique element of
$\operatorname{dv}(K)$ with $a'=b'$ and $a\prec 1$. Then 
the asymptotic couple of $\operatorname{dv}(K)$ is 
$(\Gamma+\Z\alpha,\psi^\alpha)$ as
in Lemma~\ref{extension1}, with $\alpha:=va$.
\end{enumerate}
\end{cor}
\begin{proof}
We use the notation and terminology from the proof of Theorem~\ref{thm:dv(K)}.
If $K$ has no gap, or $K$ has a gap $\beta$ and there is no $b\asymp 1$ in $K$ with
$v(b')=\beta$, then $K$ is nice,
hence  $K_1$ is an immediate extension of $K_0=K$, and so 
$\operatorname{dv}(K)$ is an immediate extension of $K$.
This shows (i) and (ii).
Suppose $K$ has a gap $\beta$ and $b\asymp 1$ in $K$ satisfies
$v(b')=\beta$.  Then $K_0=K(y)$ where $y'=b'$ and $y\prec 1$, and
 $K_1$ is an immediate extension of $K_0$. 
 Thus $\operatorname{dv}(K)$ is an immediate extension of~$K_0$, and (iii) follows from Lemma~\ref{extensionpdvfields1}, taking $a:=y\in K(D)=\operatorname{dv}(K)$.
\end{proof}

\noindent
If $K$ is $\d$-valued, then of course $\operatorname{dv}(K)=K$,
and the assumption in (i) or (ii) of Corollary~\ref{cor:value group of dv(K)} holds. The assumption in Corollary~\ref{cor:value group of dv(K)}(iii) holds for 
$(K,{\dotpreceq})$ as in Example~\ref{ex:pdv with gap, 2}.
If $K$,~$a$ are as in Example~\ref{ex:pdv with gap, 1}, then
$K(a)$,~$0$,~$a$ also satisfy the assumption on $K$,~$\beta$,~$b$
in Corollary~\ref{cor:value group of dv(K)}(iii).

In Section~\ref{sec:LO-cuts} we shall need the following:

\begin{lemma}\label{alggapcase} Suppose $K$ has
asymptotic integration. Let $L$ be a pre-$\d$-valued field extension of $K$ of $H$-type and suppose $L$ is
algebraic over~$K$ and has a gap $\beta$.
Then there is no
$b\asymp 1$ in $L$ with $v(b')=\beta$. 
\end{lemma}
\begin{proof} Identify $\operatorname{dv}(K)$ with a valued differential subfield
of $\operatorname{dv}(L)$ via the embedding $\operatorname{dv}(K)\to \operatorname{dv}(L)$ that extends the
inclusion $K \to L$. Now $\operatorname{dv}(K)$ is an immediate extension of
$K$ by  Corollary~\ref{cor:value group of dv(K)} and item (\rm{a}) preceding that corollary, and $\operatorname{dv}(L)$ is algebraic over $\operatorname{dv}(K)$
by items (\rm{a}) and (\rm{b}) preceding the corollary. It follows that
$\Gamma_{\operatorname{dv}(L)}/\Gamma$ is a torsion group, and so
$\Gamma_{\operatorname{dv}(L)}/\Gamma_L$ is a torsion group.
Now apply part (iii) of Corollary~\ref{cor:value group of dv(K)} to $L$ in the role of $K$.
\end{proof} 

\subsection*{Notes and comments} 
Theorem~\ref{thm:dv(K)} without $H$-type restriction
is in~\cite{AvdD2}.

\section{Adjoining Exponential Integrals}\label{sec:exp integrals}

\noindent 
{\em In this section $K$ is an
 $H$-asymptotic field 
with asymptotic couple $(\Gamma, \psi)$, $\Gamma\ne \{0\}$, and $a$, $b$ range over $K^\times$ and $j$, $k$ over $\Z$.}\/ Given $s\in K$ we wish to adjoin an exponential integral $\exp(\int s)$ to
$K$.  
\textit{In Proposition~\ref{newprop} and
Lemmas~\ref{newprop-1}--\ref{newprop, lemma 2} below,  
we assume
$s\in K^\times$ is such that $s\ne a^\dagger$ for all~$a$, and
$K(f)$ is a field extension of~$K$ with $f$ transcendental over~$K$, 
equipped with the unique derivation extending the derivation of $K$ such 
that $f^\dagger=s$.}\/ If $K$ is algebraically closed,   
then the constant field of $K(f)$ is~$C$, by 
Lemmas~\ref{lem:C alg closed in K}, \ref{lem:exp integration simple}, and Corollary~\ref{cor:unit alg}.
%seems to be true also for real closed $K$.

\begin{prop} \label{newprop}
Suppose $K$ is algebraically closed and $\d$-valued. Then there is a valuation of $K(f)$ that 
makes it a $\d$-valued extension of $H$-type of $K$,  
and for any such valuation with asymptotic couple $(\Gamma_f,\psi_f)$ of $K(f)$,
$$ \Gamma_f=\Gamma\ \Longleftrightarrow\ 
\text{$v(s-a^\dagger)\in (\Gamma^{>})'$  for some $a$.}$$
\end{prop}

\noindent
The proof consists of several lemmas, which give more precise information.

\begin{lemma}\label{newprop-1} Let
$K(f)$ carry a valuation making it a $\d$-valued
extension of~$K$ with value group $\Gamma_f=\Gamma$. Then 
$v(s-a^\dagger) \in (\Gamma^{>})'$ for some $a$.
\end{lemma}
\begin{proof} Since $vf\in \Gamma$, we have $a$ and $g$ such that 
$f=ag$, $g\in K(f)$, and $g\asymp 1$. Then 
$s-a^\dagger=g^\dagger\asymp g'$. Thus 
$v(s-a^\dagger)=v(g') \in (\Gamma^{>})'$.
\end{proof}

\begin{lemma}\label{step4} 
Assume $K$ is henselian and pre-$\d$-valued,
and $vs \in (\Gamma^{>})'$. 
Then there is a unique valuation
of $K(f)$ that makes it an $H$-asymptotic extension of~$K$ 
with $f-1\not\asymp 1$. With this valuation $K(f)$ is pre-$\d$-valued, and an immediate extension of~$K$ with $f \sim 1$.
\end{lemma}

\begin{proof}
Put $S := \big\{v\big(s-(1+\epsilon)^\dagger\big): \epsilon\in \smallo\big\}\subseteq (\Gamma^{>})'$. 

\claim[1]{The set $S$ has no largest element.}

\noindent
To see this, note first that $vs\in S$. Let $\gamma\in S$ with
$\gamma \ge vs$, and take $\epsilon \in \smallo$ with $\gamma = v\big(s-(1+\epsilon)^\dagger\big)$. Since 
$(\Gamma^>)'$ is upward closed, we have $b\in
\smallo$ with $v(b')=\gamma$. Take $u\in K$ with $u\asymp 1$ and $s-(1+\epsilon)^\dagger=ub'$. Now
$v(u'b) > v(b')=\gamma$, so with $\delta\in \smallo$ such that
$(1+\epsilon)(1+ub)=1+\delta$ we get 
$$s-(1+\delta)^\dagger\ =\ s-(1+\epsilon)^\dagger - (1+ub)^\dagger\ =\ 
ub'-\frac{(ub)'}{1+ub}\ =\ \frac{u^2bb'-u'b}{1+ub},$$ 
hence 
$v\big(s-(1+\delta)^\dagger\big)>\gamma$. 
This proves Claim 1.

\medskip
\noindent
Let $\kappa = \operatorname{cf}(S)$, 
so $\kappa$ is an infinite cardinal; let $\rho$, $\sigma$, $\tau$ range over the ordinals~$<\kappa$.
Take a sequence $(\epsilon_\rho)$ in $\smallo$ such that 
$v\big(s-(1+\epsilon_\rho)^\dagger\big)$ 
is strictly increasing as a function of~$\rho$, and 
cofinal in $S$. In
particular, $(1+\epsilon_\rho)^\dagger \leadsto s$. 

\claim[2]{$(\epsilon_\rho)$ is a pc-sequence in $K$.}

\noindent
For $\rho<\sigma$ we have
$s-(1+\epsilon_\rho)^\dagger\sim  (1+\epsilon_\sigma)^\dagger-(1+\epsilon_\rho)^\dagger$. Also, for $\rho\ne \sigma$, 
$(\epsilon_\sigma-\epsilon_\rho)'\succ
\epsilon_\rho'(\epsilon_\sigma-\epsilon_\rho)$, hence
$$  (1+\epsilon_\sigma)^\dagger-(1+\epsilon_\rho)^\dagger\ 
=\ \frac{(1+\epsilon_\rho)(\epsilon_\sigma-\epsilon_\rho)'-
\epsilon_\rho'(\epsilon_\sigma-\epsilon_\rho)}{(1+ \epsilon_\rho)(1+\epsilon_\sigma)}\ \sim\ 
(\epsilon_\sigma-\epsilon_\rho)'.$$
It follows that
$(\epsilon_\tau-\epsilon_\sigma)'\prec (\epsilon_\sigma-\epsilon_\rho)'$
for $\tau > \sigma > \rho$. Hence
$\epsilon_\tau-\epsilon_\sigma \prec \epsilon_\sigma-\epsilon_\rho$
for $\tau > \sigma > \rho$. Thus $(\epsilon_\rho)$ is a  pc-sequence. 

\claim[3]{The pc-sequence $(\epsilon_\rho)$ has no 
pseudolimit in $K$.} 

\noindent
To see this, suppose $\epsilon_\rho\leadsto\epsilon\in K$. Then $\epsilon \in \smallo$ and $\epsilon-\epsilon_{\rho}\asymp \epsilon_{\sigma}-\epsilon_{\rho}$ for $\sigma>\rho$, so
$(\epsilon-\epsilon_\rho)'\asymp (\epsilon_\sigma-\epsilon_\rho)'$
for $\sigma > \rho$.
Computations as in the proof of Claim~2 then give
 $(1+\epsilon)^\dagger - (1+\epsilon_\rho)^\dagger \asymp 
(1+\epsilon_\sigma)^\dagger-(1+\epsilon_\rho)^\dagger$
for $\sigma>\rho$. Hence $(1+\epsilon_{\rho})^\dagger\leadsto (1+\epsilon)^\dagger$, which in view of $(1+\epsilon_{\rho})^\dagger\leadsto s$ gives $v\big(s-(1+\epsilon)^\dagger\big)>
v\big(s-(1+\epsilon_\rho)^\dagger\big)$
for all~$\rho$, contradicting the choice of $(\epsilon_{\rho})$ . This proves Claim~3. 

\medskip
\noindent
With $P(Y):= Y'-(1+Y)s$ we have $P(\epsilon_{\rho})\leadsto 0$
and $P(y)=0$ for $y:= f-1$. By the claims we can apply Proposition~\ref{pca-extension} to $P(Y)$ to get a valuation of $K(f)$ that makes it an immediate 
$H$-asymptotic extension of $K$ with 
$\epsilon_{\rho} \leadsto y$. 

As to uniqueness, assume $K(f)$ is given a 
valuation making it an $H$-asymptotic extension of $K$ with
$y\nasymp 1$. Then $y\prec 1$, since $y\succ 1$ gives
$y\sim y+1=f$, and so
$v(y^\dagger)=v(f^\dagger)=vs\in (\Gamma^>)'$, which is impossible. From $(1+\epsilon_{\rho})^\dagger\leadsto s=(1+y)^\dagger$ we obtain as in the proof of Claim 2 that $\epsilon_{\rho}\leadsto y$, and so part (ii) of Proposition~\ref{pca-extension} yields the desired uniqueness. 
\end{proof}

\begin{lemma}\label{newprop0} 
Let $K$ be algebraically closed and $\d$-valued, and suppose $a$ satisfies
$v(s-a^\dagger)\in (\Gamma^{>})'$. Then some valuation of $K(f)$ 
makes it a $\d$-valued extension of $H$-type of $K$; 
any such valuation makes $K(f)$ an immediate extension of $K$.
\end{lemma}
\begin{proof} By Lemma~\ref{step4} with $s-a^\dagger$
and $f/a$ instead of $s$ and $f$ we get a valuation of $K(f)$  
making it a $\d$-valued extension of $H$-type of $K$. Let~$K(f)$
be equipped with any such valuation.
Then $(f/a)^\dagger = s-a^\dagger$, hence $f\asymp a$, so  $f/a = c(1+z)$ with 
$c\in C^\times$ and $z\prec 1$, and then $(1+z)^\dagger =s-a^\dagger$. 
Thus by Lemma~\ref{step4}, the valued
field $K(f)=K(1+z)$ is an immediate extension of $K$.
\end{proof}

\noindent
This takes care of Proposition~\ref{newprop} in the case that
$v(s-a^\dagger)\in (\Gamma^{>})'$ for some $a$.
It remains to consider the case 
that $v(s-a^\dagger)< (\Gamma^{>})'$ for all $a$. 
This condition is for $\d$-valued $K$ equivalent to: 
\begin{equation}\label{eq:newprop}
s-a^\dagger\ \succ\  u'\quad 
\text{ for all $a$  and for all $u\in K^{\preceq 1}$.}
\end{equation}
{\em In the next two lemmas we presuppose that \eqref{eq:newprop} holds,}\/
but do {\em not\/} assume that~$K$ is algebraically closed or 
$\d$-valued, since the extra
generality will be of some use.
(Of course we keep the assumption that $K$ is an 
$H$-asymptotic field with asymptotic couple~$(\Gamma, \psi)$, $\Gamma\ne \{0\}$.) Put 
$$ S\ :=\ \big\{ v(s-a^\dagger) :\  a\in K^\times \big\}\ \subseteq\ \Gamma,$$ so
$S<(\Gamma^{>})'$. Taking $a=1$ shows that $vs\in S$. 
The next lemma is the most substantial result of this section.

\begin{lemma} \label{newprop, lemma 1}
Suppose $S$ has no maximum and $\Gamma$ is divisible. Then there is 
a unique valuation
on $K(f)$ that makes $K(f)$ an $H$-asymptotic extension of $K$ with 
asymptotic couple $(\Gamma_f, \psi_f)$. Moreover, for this valuation we have:
\begin{enumerate}
  \item[\textup{(i)}] $vf\notin \Gamma$, $\Gamma_f=\Gamma\oplus\Z {vf}$, $[\Gamma_f]=[\Gamma]$,  
and $\res K(f) = \res K$; 
\item[\textup{(ii)}] if $K$ has small derivation, then so does $K(f)$;
\item[\textup{(iii)}] if $K$ is pre-$\d$-valued, then so is $K(f)$; 
\item[\textup{(iv)}] if $K$ is $\d$-valued, then so is $K(f)$, with the same constant field.
\end{enumerate}
\end{lemma}
\begin{proof}
Since $S$ has no largest element, Lemma~\ref{have1} and 
Corollary~\ref{trich} yield $S\subseteq (\Gamma^{<})'$.
Let the infinite cardinal $\kappa$ be the cofinality of $S$, let $\rho$, $\sigma$, $\tau$ range over the ordinals~$<\kappa$, and let 
$(a_\rho)$ be a sequence in $K^\times$ 
such that the sequence $\bigl(v(s-a_\rho^\dagger)\bigr)$ is strictly
increasing and cofinal in $S$ with $v(s-a_\rho^\dagger)>vs$ for all
$\rho$. Then $s\sim a_\rho^\dagger$ for all $\rho$. Since $vs\in (\Gamma^{<})'$,
and $b^\dagger \asymp b'$ if $b\asymp 1$, this gives $a_\rho\not\asymp 1$ for all $\rho$ by (i) of Corollary~\ref{CharacterizationAsymptoticFields-Corollary}. Also,
$$s-a_\rho^\dagger\ \asymp\ 
a_\sigma^\dagger - a_\rho^\dagger\ \asymp\ a_\rho^\dagger - a_\sigma^\dagger\qquad\text{for $\rho<\sigma$,}$$
since $s$ is a pseudolimit of the pc-sequence $(a_\rho^\dagger)$. Hence
$$(a_\rho/a_\sigma)^\dagger\ =\ a_\rho^\dagger-a_\sigma^\dagger\ \succ\
a_\sigma^\dagger-a_\tau^\dagger\ =\ (a_{\sigma}/a_\tau)^\dagger \qquad\text{for $\rho<\sigma<\tau$.}$$
We have $v(a_\rho/a_\sigma) \neq 0$ for $\rho<\sigma$:
otherwise $s-a_\rho^\dagger \asymp (a_\rho/a_\sigma)^\dagger \asymp (a_\rho/a_\sigma)'$
and $v\big((a_\rho/a_\sigma)'\big)> (\Gamma^{<})'$
by (i) of Corollary~\ref{CharacterizationAsymptoticFields-Corollary},  
a contradiction. Put $\alpha_\rho:=v(a_\rho)$, so 
$\alpha_\rho-\alpha_\sigma=v(a_\rho/a_\sigma)$ and 
$\psi(\alpha_\rho-\alpha_\sigma)=v(s-a_\rho^\dagger)$ for $\rho<\sigma$. 
Thus $(\alpha_\rho)$ is a pc-sequence with respect to the valuation $\psi$ on 
the ordered abelian group~$\Gamma$.

\claim[1]{For $\alpha=va$ we have $\psi(\alpha-\alpha_{\rho})=
v(a^\dagger-a_{\rho}^\dagger)=v(a^\dagger-s)$, eventually.}

\noindent
This is because $a^\dagger-a_\rho^\dagger=(a^\dagger-s)+(s-a_\rho^\dagger)$, and
$a^\dagger-s \succ s-a_\rho^\dagger$, eventually. Since $\alpha\in \Gamma$ in 
Claim 1 is arbitrary, it follows that  $(\alpha_\rho)$ 
has no pseudolimit in $(\Gamma, \psi)$.

\medskip\noindent
Suppose now that $K(f)$ is equipped with a valuation that makes it an 
$H$-asymptotic field extension of $K$ with asymptotic couple
$(\Gamma_f,\psi_f)$. Then $\eta:=vf\notin\Gamma$:
otherwise $f=ua$ for some $a$ and $u\in K(f)$ with $u\asymp 1$; 
for such $a$ and $u$
we have $u^\dagger\asymp u' \prec b'$ for all $b\succ 1$
by (i) of Corollary~\ref{CharacterizationAsymptoticFields-Corollary}, 
in particular 
$u^\dagger \prec s-a_\rho^\dagger$ and therefore
$$a^\dagger-a_\rho^\dagger\ =\ 
f^\dagger -u^\dagger-a_{\rho}^\dagger\ =\ 
(s-a_\rho^\dagger) - u^\dagger\ \asymp\  s-a_\rho^\dagger$$
for all $\rho$, contradicting Claim~1. This
yields $\Gamma_f=\Gamma\oplus\Z\eta$. Now 
$\psi_f(\eta-\alpha_\rho)=v(s-a_\rho^\dagger)$ for all $\rho$, so 
$\alpha_\rho\leadsto \eta$ in $(\Gamma_f, \psi_f)$, and for 
$\alpha=va$, Claim 1 gives
$$\psi_f(\alpha-\eta)\ =\ v(a^\dagger-s)\ =\ v(a^\dagger-a_\rho^\dagger)\ =\ 
\psi(\alpha-\alpha_\rho),\quad \text{eventually}.$$
Hence the sequence
$\big([\eta-\alpha_\rho]\big)$ is strictly decreasing in $[\Gamma_f]$,
and 
$[\alpha-\eta] = [\alpha-\alpha_\rho]$ eventually. Hence 
$[\Gamma_f]=[\Gamma]$. It also follows that 
for all $\alpha\in\Gamma$,
$$\alpha<\eta \Longleftrightarrow \alpha<\alpha_\rho \text{ eventually,} \qquad 
\alpha>\eta \Longleftrightarrow \alpha>\alpha_\rho \text{ eventually.}$$
This determines the ordering on $\Gamma_f$. Hence there is at most
one valuation on $K(f)$ making $K(f)$ an $H$-asymptotic field extension of 
$K$.

We now construct such a valuation. The valuation $\psi:\Gamma^{\ne} \to \Psi$ 
is coarser than
the standard valuation of $\Gamma$, so $(\alpha_\rho)$ has no
pseudolimit in $\Gamma$ with respect to the standard valuation by
Lemma~\ref{coarsepc}. Hence Lemma~\ref{divimm, standard val} gives an element~$\eta$ in an ordered abelian group extending $\Gamma$
such that $\eta\notin \Gamma$,
the sequence~${\bigl([\eta-\alpha_\rho]\bigr)}$ is eventually strictly
decreasing, and $[\Gamma_f]=[\Gamma]$, where
$\Gamma_f:=\Gamma\oplus\Z\eta$. By Lemma~\ref{extension4} we have
a unique extension of $\psi\colon\Gamma^{\ne}\to\Gamma$ to a map
$\psi_f\colon\Gamma_f^{\ne}\to\Gamma$ such that~$(\Gamma_f,\psi_f)$ is an
asymptotic couple of $H$-type. Now $\Gamma$ is divisible, so $\Gamma^{>}$ is coinitial
in $\Gamma_f^{>}$. Therefore, 
if $(\Gamma, \psi)$ has small derivation, then so does $(\Gamma_f, \psi_f)$.

\claim[2]{$\psi_f(va +j\eta)=v(a^\dagger+js)\in S$ for all $a$ and 
all $j\neq 0$.}

\noindent
To prove Claim~2, let $a$ and $j\ne 0$ be given, and take $u,d\in K^{\times}$ such that $a=ud^{-j}$ and $u\asymp 1$ (which is possible since $\Gamma$ is divisible). Then $a^\dagger = -jd^\dagger + u^\dagger$ and $u^\dagger\asymp u'\prec s-d^\dagger$,
so
\begin{align*} v(a^\dagger+js)\ &=\ v(js-jd^\dagger+ u^\dagger)\ =\ v(s-d^\dagger)\in S, \\
\psi_f(va+j\eta)\ &=\ \psi_f\big(j(-vd+\eta)\big)\ =\ \psi_f(\eta-vd),
\end{align*} 
so it only remains
to show that $\psi_f(\eta-vd)=v(s-d^\dagger)$. Since
$\bigl([\eta-\alpha_\rho]\bigr)$ is strictly decreasing, 
there are arbitrarily large $\rho$ with $[\eta-vd]=
[\alpha_\rho-vd]$.  For those $\rho$ we have
$\alpha_\rho\neq vd$ (since $\eta\neq vd$), and thus
$$\psi_f(\eta-vd)\ =\  \psi(\alpha_\rho-vd)\ =\ 
v(a_\rho^\dagger-d^\dagger).$$
By cofinality of 
$\bigl(v(s-a_\rho^\dagger)\bigr)$ in $S$ we get
$v(a_\rho^\dagger-d^\dagger)=
v(s-d^\dagger)$
for all sufficiently large $\rho$. Hence
 $\psi_f(\eta-vd)=v(s-d^\dagger)$ 
as required.

\medskip\noindent
We now equip $K(f)$ with the valuation $v\colon K(f)^\times \to \Gamma_f$
that extends the valuation of~$K$ and such that $vf=\eta$. 

\claim[3]{$K(f)$ is an asymptotic field 
with asymptotic couple $(\Gamma_f,\psi_f)$.}

\noindent
By Claim~2 we have:
\begin{equation}\label{psi-f}
v\big((af^j)^\dagger\big)\ =\ \psi_f\big(va+j\eta\big) 
\qquad\text{if $va+j\eta\neq 0$.}
\end{equation} 
Therefore, if $va+j\eta,vb+k\eta>0$, then $$v\big((af^j)'\big)\ =\
(va+j\eta)'\ >\ \psi_f(vb+k\eta)\ =\ v\big((bf^k)^\dagger\big).$$ 
Now suppose $va+j\eta=0$; then $j=0$, so $v\big((af^j)'\big)=v(a')>S$. 
In particular, 
if $k\neq 0$, then $v\big((af^j)'\big)>v\big((bf^k)^\dagger\big)$ by 
\eqref{psi-f} and Claim~2. If $vb\neq 0$, then 
$v\big((af^j)'\big)=v(a')\geq v(b^\dagger)$.
Since $\Gamma_f\neq\Gamma$ we have $\res K(f)=\res K$.
Hence $K$, $L:=K(f)$ and $T:= K^\times f^{\Z}$ satisfy the hypotheses 
of Lemma~\ref{Lemma34}, so $L$ is an asymptotic field. By 
\eqref{psi-f}, its asymptotic couple is $(\Gamma_f,\psi_f)$. 
If $K$ is pre-$\d$-valued, then so is $K(f)$ by Corollary~\ref{cor:pdv extension, 1}. If~$K$ is $\d$-valued, then so is $K(f)$ with
$C_{K(f)}=C$, by Lemma~\ref{dv}.
\end{proof}

\noindent
The uniqueness part of Lemma~\ref{newprop, lemma 1}
strengthens Proposition~\ref{newprop} significantly 
in the case considered.

\begin{lemma}\label{newprop, lemma 2}
Suppose $S$ has a maximum and $\Gamma$ is divisible. Then there is a valuation on $K(f)$ 
making it an $H$-asymptotic extension of $K$ with
$[v(af)] \notin [\Gamma]$ for some $a$.
For any such valuation, with asymptotic couple $(\Gamma_f, \psi_f)$
of $K(f)$,  \begin{enumerate}
\item[\textup{(i)}] $vf\notin \Gamma$, $\Gamma_f=\Gamma\oplus \Z vf$,  
$\res K(f)=\res K$, and 
$\Psi_f=\Psi\cup \{\max S\} \subseteq \Gamma$;
\item[\textup{(ii)}] if $K$ has small derivation, then so does $K(f)$; and
\item[\textup{(iii)}] if $K$ is $\d$-valued, then so is $K(f)$, with the same 
constant field.
\end{enumerate}
\end{lemma}
\begin{proof}
Let $vt$ with $t=s-b^\dagger$ be the largest element of $S$. Then 
$v(t-a^\dagger)=v\big(s-(ab)^\dagger\big)\leq vt$ for all $a$. By renaming 
$f/b$ as $f$ and $t$ as $s$ the hypotheses of  
Lemma~\ref{newprop, lemma 2} remain valid; 
then $v(s-a^\dagger)\leq vs$ for all $a$. Also 
\begin{equation}\label{min-formula asymptotic}
v(a^\dagger+js)\ =\ 
\min\big(vs,v(a^\dagger)\big)\in S\qquad\text{for all $a$ and all $j\neq 0$.}
\end{equation}
To see this, let $a$ and $j\ne 0$ be given. Since $\Gamma$ is divisible we can
take $u,d\in K^{\times}$ such that $u\asymp 1$ and $a=ud^{-j}$. 
Then $a^\dagger = -jd^\dagger + u^\dagger$ and $u^\dagger\asymp u'\prec s-d^\dagger$,
so
$$v(a^\dagger+js)\ =\ v(js-jd^\dagger+ u^\dagger)\ =\ v(s-d^\dagger)\
=\ \min\big(vs,v(d^\dagger)\big)\ =\ 
\min\big(vs,v(a^\dagger)\big).$$
(To get the last two equalities, consider the cases $vs = v(d^\dagger)$ and $vs=v(a^\dagger)$ separately.) 
Since $(\Gamma,\psi)$ is of $H$-type, 
$$M\ :=\ \bigl\{\gamma\in\Gamma^{<}:\ \psi(\gamma) \leq vs\bigr\}$$ 
is downward closed  in $\Gamma$.
Let $\eta$ be an element of an ordered abelian group extending~$\Gamma$ such that $M<\eta<\Gamma\setminus M$, so $\eta < 0$. Put  
$\Gamma_f:=\Gamma\oplus \Z\eta$. Then $[\eta]\notin [\Gamma]$:
otherwise $n\gamma < \eta < \gamma$ where $\gamma\in \Gamma^{<}$ and $n\geq 2$,
so $\psi(\gamma)= \psi(n\gamma)\le vs$, a contradiction.
Hence $[\Gamma_{f}]=[\Gamma]\cup \bigl\{[\eta]\bigr\}$.
We can now apply Lemma~\ref{extension5} with 
$$C\ :=\ \big\{[\gamma]:\ \gamma\in \Gamma,\ M < \gamma < 0\big\}, \quad \beta: =vs.$$
This gives an $H$-asymptotic couple $(\Gamma_f,\psi_f)$ that extends $(\Gamma, \psi)$
with 
$$\psi_f\colon \Gamma_f^{\ne} \to \Gamma, \qquad
\psi_f(\gamma + j\eta)\ :=\ 
\min\big(\psi(\gamma),vs\big)  \quad
\text{ for $\gamma\in\Gamma$ and  $j\neq 0$.}$$
Equip $K(f)$ with the valuation
$v \colon K(f)^\times \to \Gamma_f$ extending the valuation of $K$ 
such that $vf=\eta$.  Since 
$\Gamma_f\neq\Gamma$ we have $\res K(f)=\res K$. If $va+j\eta\neq 0$, then 
$$v\big((af^j)^\dagger\big)\ =\ v(a^\dagger+js)\ =\ \psi_f(va+j\eta)$$ 
by \eqref{min-formula asymptotic}. Therefore, if $va+j\eta,vb+k\eta>0$, 
then 
$$v\big((af^j)'\big)\ =\ (va+j\eta)'\ >\ \psi_f(vb+k\eta)\ =\ 
v\big((bf^k)^\dagger\big).$$ 
Next suppose $va+j\eta=0$; then $va=0$, $j=0$, 
so $v\big((af^j)'\big)=v(a')$. Now $v(a')>vs$, 
and if $vb\neq 0$, then $v(a')\geq v(b^\dagger)$; hence 
$v\big((af^j)'\big)\geq v\big((bf^k)^\dagger\big)$ 
if $vb+k\eta\neq 0$. We have now shown that $K$, $L=K(f)$ and 
$T=K^\times f^{\Z}$ 
satisfy the hypotheses of Lemma~\ref{Lemma34}, so $K(f)$ is an 
asymptotic field.

Now let $K(f)$ be equipped with any valuation that makes
it an $H$-asymptotic extension of $K$ with asymptotic couple 
$(\Gamma_f, \psi_f)$
such that $[v(af)]\notin [\Gamma]$ for some~$a$. Fix such~$a$. 
Then $vf\notin \Gamma$, so $\Gamma_f=\Gamma\oplus \Z vf$ and 
$\res K(f)=\res K$. Also
$[\Gamma_f]=[\Gamma]
\cup \big\{[v(af)]\big\}$ by Lemma~\ref{archextclass}, and 
$\psi_f\big(v(af)\big)=v(a^\dagger+s)$, hence
$$\Psi_f\ =\ \Psi \cup \big\{v(a^\dagger+s)\big\}.$$
To finish (i) we show that $v(a^\dagger+s)=\max S$: clearly, 
$v(a^\dagger+s)\in S$, and
\begin{align*} v(s-b^\dagger)\ &=\ v\big((af/ab)^\dagger\big)\ =\ 
\psi_f\big(v(af)-v(ab)\big)\\ 
            &=\ \min\{\psi_f\big(v(af)\big),\ \psi\big(v(ab)\big)\}\ \le\ 
\psi_f\big(v(af)\big)\ =\ v(a^\dagger+s),
\end{align*}
where the third equality holds because $[v(af)]\ne [v(ab)]$ and
$(\Gamma_f, \psi_f)$ is of $H$-type. 
To prove (ii), assume that $K$ has small derivation. We 
wish to show that~$K(f)$ has small derivation. With $a$ as above we have 
$\psi_f\big(v(af)\big)=v(a^\dagger + s)\in \Psi_f$, so we are done if
$v(a^\dagger + s)\ge 0$. Suppose $v(a^\dagger + s)< 0$. Then
Lemma~\ref{PresInf-Lemma} gives $\gamma\in \Gamma^{<}$ such that 
$v(a^\dagger +s) < \psi(\gamma)$. Then
$v(af) < \gamma <0$ or $-v(af) < \gamma < 0$. Because $[\Gamma_f]=[\Gamma]\cup
\big\{[v(af)]\big\}$, it follows that $\Gamma^{<}$ is cofinal in 
$\Gamma_f^{<}$, and we are done by Lemma~\ref{PresInf-Lemma}. 
Lemma~\ref{dv} gives (iii). 
\end{proof} 

\noindent
Lemmas~\ref{newprop-1},~\ref{newprop0},~\ref{newprop, lemma 1},~\ref{newprop, lemma 2} above cover all claims of
Proposition~\ref{newprop}.

\begin{example}
Let $K=\C[[x^{\Q}]]$ with its usual derivation 
($x'=1$ and $C=\C$) and valuation $v\colon K^\times \to \Gamma=\Q$, $vx=-1$, so
$K$ is an algebraically closed $\d$-valued field of $H$-type 
with $\Psi= \{1\}$. Let $s=\imag\in\C$, $\imag^2=-1$, so $s\ne a^\dagger$
for all $a$. Let
$f$ be as in the beginning of this subsection: transcendental over $K$ 
with the derivation on~$K(f)$ extending that of $K$ with $f^\dagger=s$.
(One may think of $f$ as $\ex^{\imag x}$.) 
We have $S=\{0\}=\{vs\}$, and $K(f)$ with the valuation constructed 
in the proof of Lemma~\ref{newprop, lemma 2} is $\d$-valued and
satisfies $vf<\Gamma$ and $\Psi_f=\{0,1\}$.
\end{example}

\begin{cor} \label{newcor}
Suppose that $K$ is algebraically closed, $\d$-valued, and 
has small derivation.
Then $K$ has a $\d$-valued extension $L$ of $H$-type
with small derivation and the following properties:
\begin{enumerate}
\item[\textup{(i)}] $L$ is algebraically closed with constant field $C$;
\item[\textup{(ii)}] for each $s\in K$ there exists $f\in L^\times$
with $f^\dagger=s$;
\item[\textup{(iii)}] $L$ is algebraic over its subfield 
generated over $K$ by the $f\in L^\times$ with $f^\dagger\in K$;
\item[\textup{(iv)}] there exists a family $(f_i)_{i\in I}$
of elements of $L^\times$ such that $f_i^\dagger\in K$ for all $i\in I$, and
$\Gamma_L=\Gamma\oplus \bigoplus_{i\in I} \Q v(f_i)$ 
\textup{(}internal direct sum of $\Q$-linear subspaces\textup{)}.
\end{enumerate}
Moreover, any such $L$ has the following additional properties:
\begin{enumerate}
\item[\textup{(v)}] $\Psi_L\subseteq \Gamma$, and if $\Psi$ has a maximum $\alpha$, 
then $\Psi_L=\Gamma^{\leq\alpha}$;
\item[\textup{(vi)}] $L$ has no proper differential field extension $M$ satisfying \textup{(i),~(ii),~(iii)} with $M$ instead of $L$.
\end{enumerate}
\end{cor}
\begin{proof} 
Iterating the extension steps of Lemmas~\ref{newprop0},
\ref{newprop, lemma 1}, and \ref{newprop, lemma 2}, and alternating these 
steps by taking algebraic closures, we obtain a $\d$-valued 
field extension~$L$ of $K$, with small derivation, of $H$-type, and
satisfying (i)--(iv).

Let any such $L$ be given. To prove (v), consider any element of $\Psi_L$.
It equals~$v(f^\dagger)$ with $f\in L^\times$, $f\nasymp 1$. By (iv) we have a
nonzero $k\in \Z$ such that
$$f^k\asymp af_1^{k_1}\cdots f_m^{k_m}, \qquad f_1,\dots,f_m\in L^\times,\ k_1,\dots, k_m\in \Z.$$
Then $kf^\dagger\asymp a^\dagger + k_1f_1^{\dagger} + \cdots + k_mf_m^\dagger \in K$,
so $v(f^\dagger)\in \Gamma$, as claimed. 
Now assume also that $\Psi$ has largest element $\alpha$. Then by 
$\Psi_L\subseteq \Gamma$ and 
Lemma~\ref{asy1} we have 
$\Psi_L\subseteq \Gamma^{\leq\alpha}$. The reverse inclusion is clear from (ii).

If $M$ is a 
differential field extension of $L$ satisfying (i), (ii), (iii), then for 
any ele\-ment~$g\neq 0$ of~$M$ with $g^\dagger\in K$ we have $f\in L^\times$ such that
$f^\dagger=g^\dagger$, hence $g=cf$ for some $c\in C$, and thus $g\in L$.
This proves (vi). 
\end{proof}

\subsection*{Notes and comments}
Lemma~\ref{step4} is a variant of \cite[Lemma~5.2]{AvdD2}. The other results in
this section are new. 

\section{$H$-Fields and Pre-$H$-Fields}\label{sec:H-fields}

\noindent
An {\bf $H$-field\/} \label{p:H-field}\index{H-field@$H$-field}
\index{ordered differential field!H-field@$H$-field} is by definition an ordered differential field $K$ 
such that:
\begin{list}{*}{\setlength\leftmargin{3em}}
\item[(H1)] for all $f\in K$, if $f>C$, then $f'>0$;
\item[(H2)] $\mathcal{O}=C+\smallo$ where 
$\mathcal{O}= \big\{g\in K:\text{$\abs{g} \le c$ for some $c\in C$}\big\}$ and $\smallo$ is the maximal ideal of the convex subring $\mathcal{O}$ of $K$.
\end{list} 
Hardy fields containing $\R$ as a subfield are $H$-fields, by Proposition~\ref{LHospital}(i),~(iii). Also~$\mathbb{T}$ with its usual ordering and derivation is an $H$-field.
Any compositional conjugate~$K^{\phi}$ of an $H$-field $K$ with 
$\phi\in K^{>}$ is an $H$-field (with the same ordering).

We regard any $H$-field
$K$ as an ordered valued differential field by taking the valuation given by the valuation ring $\mathcal{O}$ defined in (H2). Note that in an 
$H$-field the ordering determines the valuation, but in our experience
the valuation is the more robust and useful feature. Nevertheless, 
the ordering deserves attention.

\begin{lemma}
Let $K$ be an $H$-field, $f,g\in K$, $f\preceq 1$, $0\neq g\nasymp 1$. Then $f' \prec g^\dagger$.
\end{lemma}
\begin{proof}
Subtracting a constant from~$f$ and using (H2), we may assume $f\prec 1$.
Replacing $g$ by $-g$ if necessary we may assume that $g>0$; and replacing~$g$ by $1/g$ if necessary we may assume that $g>C$.  Let $c\in C^>$. Then $c+f$, $c-f$ are~$>0$ and~$\asymp 1$, so $g(c+f),g(c-f)>C$ and hence $g'(c+f)+gf'>0$ and $g'(c-f)-gf'>0$ by~(H1). Dividing by $g'>0$ gives
$-c-f<f'/g^\dagger<c-f$. This holds for all $c\in C^>$, so $f'/g^\dagger\prec 1$.
\end{proof}

\noindent
Note that $H$-fields are $\d$-valued by the previous lemma and axiom (H2). 
% in particular, $K$ is a Hahn space when
%viewed as valued vector space over $C$.
%We also note that for $a,b\in K^{>}$ we have: 
%$a\asymp b$ if and only if
%$a$ and $b$ have the same $C$-archimedean class, 
%that is, there exist
%$c,d \in C^{>}$ such that $c a \le b \le d a$. 
%Thus $K$ as ordered vector space over $C$
%is a Hahn space in the sense of Section~\ref{sec:oag}.
The example at the end of Section~\ref{As-Fields,As-Couples} shows that a differential 
subfield of an $H$-field with the induced ordering and valuation is not
always an $H$-field; it is, however, always a pre-$H$-field in the following 
sense: \index{pre-$H$-field}\index{ordered differential field!pre-$H$-field}
A {\bf pre-$H$-field\/} is an ordered pre-$\d$-valued field $K$ whose
ordering, valuation, and derivation interact as follows: 
\begin{list}{*}{\setlength\leftmargin{3em}}
\item[(PH1)] the valuation ring $\mathcal{O}$ is convex with respect to 
the ordering;
\item[(PH2)] for all $f\in K$, if $f>\mathcal{O}$, then $f'>0$. \label{p:pre-H-field}
\end{list}
If $K$ is a pre-$H$-field, then so is any ordered valued differential subfield of $K$, and any compositional conjugate $K^{\phi}$ with $\phi\in K^{>}$. 
Hardy fields are pre-$H$-fields. 
Any ordered differential field with the trivial valuation is a pre-$H$-field.
Since we construe $H$-fields as ordered valued differential fields, they are
in particular pre-$H$-fields. 
By part (ii) of the next lemma, pre-$H$-fields are $H$-asymptotic fields:

\begin{lemma}\label{psidecreasing1}
Let $K$ be a pre-$H$-field and $f,g\in K^\times$. 
\begin{enumerate}
\item[\textup{(i)}] If $f\prec g$, then $f^\dagger < g^\dagger$. 
\item[\textup{(ii)}] If $f\prec g\prec 1$, then $f^\dagger\succeq g^\dagger$.
\end{enumerate}
\end{lemma}
\begin{proof}
Suppose $f\prec g$. Then $g/f\succ 1$ and hence  $g^\dagger=f^\dagger+(g/f)^\dagger>f^\dagger$ by (PH2).
This shows (i), and (ii) follows from (i) by taking inverses.
\end{proof}

\begin{lemma}\label{intersection}
Let $K$ be a $\d$-valued field and let $(K_i)_{i\in I}$   be a family of $\d$-valued
subfields of $K$ with $I\ne \emptyset$. Then 
$\bigcap_i K_i$
is a $\d$-valued subfield of $K$. 
\end{lemma}
\begin{proof}
With $\mathcal O_i$ the valuation ring of $K_i$, the valuation
ring of $\bigcap_i K_i$ is $\bigcap_i \mathcal O_i$. Now, given any 
$a\in \big(\bigcap_i \mathcal O_i\big)^{\ne}$, there are unique $c\in C$ and 
$c_i\in C_{K_i}$ such that $a\sim c$ and $a\sim c_i$. Hence all $c_i$ are equal to $c$, and thus $c\in \bigcap_i K_i$.
It follows that   $\bigcap_i K_i$  is $\d$-valued.
\end{proof}

\noindent
It is easy to check that for any pre-$H$-field $K$,
$$ \text{$K$ is an $H$-field}\ \Longleftrightarrow\  \text{$K$ is $\d$-valued.}$$
Thus Lemma~\ref{intersection} goes through if \textit{$\d$-valued field}\/ and \textit{$\d$-valued subfield}\/ are replaced by 
\textit{$H$-field}\/ and \textit{$H$-subfield}\/, respectively.

\subsection*{Algebraic extensions of pre-$H$-fields}
\textit{In this subsection $K$ is a pre-$H$-field.}\/

\begin{prop}\label{Algebraic-Extensions-preh}
Let $L|K$ be an algebraic extension of ordered valued differential fields such that $\mathcal O_L$ is the convex hull of $\mathcal O$ in $L$.  Then $L$ is a pre-$H$-field.
\end{prop}

\noindent
This follows from Proposition~\ref{Algebraic-Extensions-pdv} and the case $T=K^\times$ of the next lemma:

\begin{lemma}\label{prehunr}
Let $L$ be an ordered pre-$\d$-valued extension of $K$ with convex valuation ring $\mathcal O_L$. Let $T$ be a subgroup of $L^\times$ such that $\Gamma_L \subseteq \Q v(T)$ \textup{(}in $\Q\Gamma_L$\textup{)} and for each $t\in T$ with $t\succ 1$ we have $t^\dagger>0$. 
Then $L$ is a pre-$H$-field.
\end{lemma}
\begin{proof}
Since $L$ satisfies (PH1),  it remains to verify (PH2). For this, let $f\in L$ and $f>\mathcal O_L$.
Take $n\geq 1$ and $t\in T$ with $f^n\asymp t$. Then $nf^\dagger=(f^n)^\dagger \sim t^\dagger$ by Corollary~\ref{cor:pdv dagger}, and so $t^\dagger>0$ gives $f'>0$.
\end{proof}

\noindent
By Corollary~\ref{cvrc} and its proof, the convex
hull $\mathcal{O}^{\operatorname{rc}}$ of $\mathcal{O}$ in the
real closure $K^{\operatorname{rc}}$ of~$K$ is the only convex valuation ring of
$K^{\operatorname{rc}}$ that lies over $\mathcal{O}$. Consider $K^{\operatorname{rc}}$ as equipped with the valuation corresponding to $\mathcal{O}^{\operatorname{rc}}$. Then $K^{\operatorname{rc}}$ is a pre-$H$-field extension of $K$, by Proposition~\ref{Algebraic-Extensions-preh}. The constant field of 
$K^{\operatorname{rc}}$ is 
the real closure of the ordered subfield~$C$ of $K$ in $K^{\operatorname{rc}}$, by Lemmas~\ref{lem:C alg closed in K} 
and \ref{lem:const field of algext}. Consequently:

\begin{cor}\label{cor:H-field rc} If $K$ is an $H$-field, then so is $K^{\operatorname{rc}}$.
\end{cor}

\noindent
We equip the differential 
field extension $K[\imag]$ ($\imag^2=-1$) of $K$ with the valuation ring 
$\mathcal {O} + \mathcal{O}\imag$ of~$K[\imag]$; see Lemma~\ref{imagconvex}.
This makes $K[\imag]$ a valued differential
field extension of~$K$ with the same value group as $K$, and we have
$v(a+b\imag) = \min(va, vb)$ for all $a,b\in K$. Moreover:

\begin{lemma} The valued differential field $K[\imag]$ is pre-$\d$-valued and of $H$-type. 
If~$K$ is an $H$-field, then $K[\imag]$ is $\d$-valued. 
\end{lemma}
\begin{proof} $K[\imag]$ is pre-$\d$-valued by Proposition~\ref{Algebraic-Extensions-pdv}.
%Let $g\in K[\imag]^\times$, $g=a+b\imag$, $a,b\in K$. Then  $g\bar{g}= a^2+ b^2$, so $g^\dagger + {\bar{g}}^\dagger= (a^2+b^2)^\dagger$, which in view of $g^\dagger \asymp  {\bar{g}}^\dagger$ gives $g^\dagger \succeq (a^2+b^2)^\dagger$. As $K$ is pre-$\d$-valued, it follows easily that $K[\imag]$ is pre-$\d$-valued. 
Since $K$ and $K[\imag]$ have the 
same asymptotic couple, $K[\imag]$ is of $H$-type. For the last claim of the lemma, use that the constant field of $K[\imag]$ is $C+C\imag$. 
\end{proof}

\noindent
Our main interest is in the model theory of the $H$-field 
$\mathbb{T}$, but it can be useful to work in algebraically closed 
extensions like
$\mathbb{T}[\imag]$ that are not $H$-fields 
but still $\d$-valued of $H$-type. (In Section~\ref{Miscellaneous Facts about Asymptotic Fields} we show that the subset $\T$ of $\T[\imag]$ is definable
in the differential field $\T[\imag]$.)

\subsection*{Immediate extensions of pre-$H$-fields}
\textit{In this subsection $K$ is a pre-$H$-field.}\/
The next lemma is often used
and elaborates on 
Corollary~\ref{lem:ordering immediate ext}.

\begin{lemma}\label{prehimm} Let $L$ be an immediate asymptotic
extension of $K$. Then $L$ has a unique field ordering extending that of $K$ 
in which $\mathcal{O}_L$ is convex. With this ordering $L$ is a pre-$H$-field 
and $\mathcal{O}_L$ is the convex hull of $\mathcal{O}$ in $L$.
If $K$ is an $H$-field, then so is $L$ with this ordering.
\end{lemma} 

\begin{proof} $L$ is pre-$\d$-valued
by Corollary~\ref{cor:pdv extension, 2}.  If $K$ is an $H$-field,
then $L$ is $\d$-valued by Lemma~\ref{dv}.
Corollary~\ref{lem:ordering immediate ext} yields a unique ordering on $L$ as claimed.  With this ordering $L$ is a pre-$H$-field, by Lemma~\ref{prehunr} with $T=K^\times$.
\end{proof}

\noindent
We always consider any immediate asymptotic 
extension of $K$ as a pre-$H$-field according to Lemma~\ref{prehimm}.
In view of Corollary~\ref{completingasymp} we get:

\begin{cor}\label{completinghfields} The completion $K^{\operatorname{c}}$ of $K$ is a pre-$H$-field. If $K$ is an $H$-field, then so is $K^{\operatorname{c}}$.
\end{cor}

\subsection*{Adjoining integrals}
\textit{In this subsection $K$ is a pre-$H$-field.}\/
The next three results are pre-$H$-field versions of Lemmas~\ref{extensionpdvfields1}, \ref{variant41} and \ref{pre-extas2}. 

\begin{cor}\label{extensionpdvfields1, cor}
Let  $s$ and $K(y)$ be as in Lemma~\ref{extensionpdvfields1}.
Then there is a unique ordering on $K(y)$ making it a pre-$H$-field extension of $K$.
If $K$ is an $H$-field, then so is $K(y)$ with that unique ordering, and $C_{K(y)}=C$.
\end{cor}
\begin{proof} 
By passing from $s$, $y$ to $-s$, $-y$ if necessary, we
arrange that $s<0$.
Then $\smallo < y < K^{>\smallo}=\{f\in K:f>\smallo\}$ for any ordering of $K(y)$ as in the corollary.

Equip $K(y)$ with the unique ordering that according to Lemma~\ref{extror} makes it an ordered field extension of $K$
with $y>0$ and convex $\mathcal{O}_{K(y)}$. To show that then~$K(y)$ is a pre-$H$-field, we
verify the hypothesis of Lemma~\ref{prehunr} for $L=K(y)$ and $T=K^\times y^{\mathbb Z}$. So let $t=gy^j\succ 1$ with $g\in K^\times$ and $j\in \Z$. Then $t^\dagger=g^\dagger + js/y$, and
$\gamma+j\alpha<0$, where
$\gamma:= vg$ and $\alpha:= vy$. So either $\gamma<0$, or
$\gamma=0$ and $j<0$. If $\gamma<0$, then $v(g^\dagger)=\psi(\gamma) < v(s/y)=vs-\alpha$, that is, $g^\dagger\succ s/y$,
which in view of $g^\dagger>0$ gives $t^\dagger>0$. If $\gamma=0$ and $j<0$, then $v(g^\dagger)\ge vs>v(s/y)$, so $s/y\succ g^\dagger$, which by $j<0$ and $s/y<0$ yields again $t^\dagger>0$.  

For the $H$-field case, see the remark following Lemma~\ref{extensionpdvfields1}.
\end{proof}

\begin{cor}\label{variant41, cor}
Let the $H$-field $K$ and $s$ and $K(y)$ be as in Lem\-ma~\ref{variant41}.
Then there is a unique ordering on $K(y)$ making it a pre-$H$-field extension of $K$.
In fact,~$K(y)$ with that unique ordering is an $H$-field and $C_{K(y)}=C$.
\end{cor}
\begin{proof}
Similar to that of Corollary~\ref{extensionpdvfields1, cor}. (If $s>0$ and $K(y)$ is equipped with an ordering making $K(y)$ a pre-$H$-field extension of $K$, then 
$\mathcal O < y < K^{>\mathcal O}$.)
\end{proof}

\begin{cor}\label{pre-extas2, pre-H}
With $K$, $s$, and $K(y)$ as in Lemma~\ref{pre-extas2}, there is a unique ordering on $K(y)$ making it a pre-$H$-field extension of $K$. If $K$ is an $H$-field, then so is $K(y)$ with that ordering, and $C_{K(y)}=C$.
\end{cor}
\begin{proof} Like that of Corollary~\ref{variant41, cor}: here also
$s>0$ forces $\mathcal O < y < K^{>\mathcal O}$ for any ordering
on $K(y)$ making it a pre-$H$-field extension of $K$. For the $H$-field case, use the remarks following the proof of Lemma~\ref{pre-extas2}. 
\end{proof}

\noindent
As a pre-$H$-field, $K$ is pre-$\d$-valued and so has a $\d$-valued hull $\operatorname{dv}(K)$. In view of~\ref{extensionpdvfields1, cor} and~\ref{prehimm}, the construction of $\operatorname{dv}(K)$ in the proof of
Theorem~\ref{thm:dv(K)} yields an analogue for pre-$H$-fields:

\begin{cor}\label{cor:H(K)}
A unique field ordering on $\operatorname{dv}(K)$ makes $\operatorname{dv}(K)$ a pre-$H$-field extension of $K$. 
Let~$H(K)$ be $\operatorname{dv}(K)$ equipped with this ordering.
Then $H(K)$ is an $H$-field and embeds uniquely over $K$ into any $H$-field extension of $K$.
\end{cor}

\noindent 
Of course this universal property determines $H(K)$ up to unique isomorphism
(of ordered valued differential fields) over $K$.  We call $H(K)$ the {\bf $H$-field hull} of $K$. 

\index{hull!H-field@$H$-field}
\index{H-field@$H$-field!hull}
\nomenclature[Z]{$H(K)$}{$H$-field hull of the pre-$H$-field $K$}

\subsection*{Extending the constant field} Here we show that $\d$-valued fields have $\d$-valued extensions by constants; likewise for $H$-fields. First a variant of Lemma~\ref{pa4}:

\begin{lemma}\label{pU} Let $K$ be a pre-$\d$-valued field and
$L$ a valued differential field extension of $K$ such that $\Gamma_L=\Gamma$. Let $U\supseteq K$ be a
$K$-linear subspace of $L$ such that 
$L=\{u/w : u,w\in U,\, w\neq 0\}$, and for all $u\in U^{\ne}$ and
$a\in K^\times$,
$$ u\preceq 1,\ a\prec 1
\quad\Longrightarrow\quad u'\prec a^\dagger.$$
Then $L$ is pre-$\d$-valued.
\end{lemma}
\begin{proof}
Let $f,g\in L^\times$, $f\preceq 1$, $g\prec 1$. Then
$f=u/w$ with $u,w\in U^{\ne}$. Using $U\supseteq K$ and $\Gamma_L=\Gamma$ we arrange $u\asymp f$
and $w\asymp 1$. Also $g=ah$ with $a\in K^\times$, $a\prec 1$, and $h\asymp 1$. Then $f'=(u'-fw')/w\asymp u'-fw' \prec 
a^\dagger$. Thus $f'\prec a^\dagger$, and likewise 
$h^\dagger\asymp h'\prec a^\dagger$.
Hence $g^\dagger=a^\dagger + h^\dagger\sim a^\dagger$, so $f'\prec g^\dagger$. 
\end{proof}

\noindent
Let $K$ be a $\d$-valued field. Note that then $K$, 
as a valued vector space over $C$, is a Hahn space as defined in Section~\ref{sec:valued vector spaces}. Next, let $L$ be an extension field of $K$ with a 
subfield $D \supseteq C$ such that $K$ and $D$ are linearly disjoint over $C$ and $L = K(D)$. Then 
Lemma~\ref{lem:constfieldext} provides a unique derivation on $L$ that extends the one of $K$ and is trivial on $D$; this derivation has $D$ as its constant field; below we consider $L$ as a differential field in this way. 

\begin{prop}\label{prop:constant field ext}
There exists a unique valuation on the field $L$ extending that of $K$ and trivial on $D$.
This valuation has the same value group as $K$.
Equipped with this valuation, $L$ is $\d$-valued.
\end{prop}
\begin{proof} The linear disjointness gives a $D$-linear isomorphism $D\otimes_C K \to K[D]$ given by 
$d\otimes a \mapsto da$. By Corollary~\ref{ScalarExtension} this yields a valuation $v\colon K[D]^{\ne} \to \Gamma$ on the abelian group $K[D]$ that extends the valuation of $K$ and makes $K[D]$ a Hahn space over $D$. 
An easy consequence of Lemma~\ref{LemmaStandardForm} is that $v(fg)=v(f)+v(g)$ for $f,g\in K[D]^{\ne}$, so $v$ extends to a (field) valuation $v\colon L^\times \to \Gamma$. Corollary~\ref{ScalarExtension} also shows that $v$ is the only field valuation on $L$ that extends the valuation of $K$ and is trivial on $D$. Lemma~\ref{pU} for $U=K[D]$ and Lemma~\ref{LemmaStandardForm} show that $L$ with $v$ is pre-$\d$-valued.
Lemma~\ref{LemmaStandardForm} also gives for $f\in K[D]^{\ne}$ that
$f\sim da$ (with respect to $v$) for suitable $d\in D^\times$ and
$a\in K^\times$. Hence $L$ with $v$ is $\d$-valued.    
\end{proof}

\noindent
In the next proposition we consider the differential field $L$ to be equipped with the
valuation of Proposition~\ref{prop:constant field ext}. Thus $L$ extends the $\d$-valued field $K$. 

\begin{prop}\label{prop:Hconstants}
Let $K$ and $D$ be given orderings that
make $K$ an $H$-field and~$D$ an ordered field extension of $C$. Then there is a unique field ordering
of~$L$ extending the orderings of~$K$ and $D$ in which the valuation ring
of $L$ is convex. With this ordering~$L$ is an $H$-field.
\end{prop}
\begin{proof} Let $f\in K[D]^{\ne}$. Lemma~\ref{LemmaStandardForm} gives $f\sim da$ with $d\in D^{\ne}$ and $a\in K^{>}$ (and thus $f\asymp a$). We claim that the sign of $d$ in the ordered field $D$ depends only on~$f$, not on the choice of $d$ and $a$. To see this, suppose also $f\sim eb$ with $e\in D^{\ne}$ and $b\in K^{>}$. Then $d/e\sim b/a\asymp 1$
in the $\d$-valued field $L$, and $b/a>0$ in $K$. Hence $d/e\in C^{>}$, and so $d$ and $e$ have the same sign.  Thus we can make $K[D]$ uniquely into an ordered vector space over $D$ such that for any $f$ as above we have $f > 0$ iff $d > 0$. This ordering is clearly compatible with multiplication: if 
$0 < f,g \in K[D]$, then $0 < fg$. We extend this ordering to the fraction field $L$ of $K[D]$ to make $L$ an ordered field. Clearly this ordering on $L$ is the only candidate for meeting the requirements. It does extend the orderings of $K$ and $D$, 
and it is an easy exercise to check, first, that $\mathcal{O}_L$ is convex
in $L$ for this ordering, and next, that it is the convex hull of 
$D$ in $L$ for this ordering. With this ordering $L$ is a pre-$H$-field by Lemma~\ref{prehunr} applied to $T=K^\times$, and
thus an $H$-field, since $L$ is $\d$-valued. 
\end{proof}

\begin{cor}\label{cor:dconstantsval} Let $F$ be a pre-$H$-field, $E$ a differential subfield of $F$, and $F=E(C_F)$. Then  
$v(F^\times)=v(E^\times)+\Z\alpha$ for some 
$\alpha\in v(F^\times)$.
\end{cor}
\begin{proof} Extending $F$ to its $H$-field hull
we arrange that $F$ is an $H$-field. Next, consider $E$
as a pre-$H$-subfield of $F$ and identify its $H$-field hull $H(E)$ with an $H$-subfield of $F$ via its embedding over $E$ into $F$. Then $H(E)$ and $C_F$ are linearly disjoint over $C_{H(E)}$ by a remark preceding Lemma~\ref{lem:lin disjoint over constants}.
It remains to note that $v(F^\times)=v\big(H(E)^\times\big)$ by Proposition~\ref{prop:constant field ext}, and that
$v\big(H(E)^\times\big)=v(E^\times)+\Z\alpha$ for some 
$\alpha\in v\big(H(E)^\times\big)$, by Corollary~\ref{cor:value group of dv(K)}. 
\end{proof}

\subsection*{Adjoining exponential integrals}
\textit{In this subsection $K$ is a pre-$H$-field
with asymptotic couple $(\Gamma, \psi)$, $\Gamma\ne \{0\}$, and~$a$,~$b$ range over~$K^\times$ and $j$,~$k$ over~$\Z$. 
We assume
$s\in K^\times$ is such that $s\ne a^\dagger$ for each~$a$, and
we take a field extension $K(f)$ of~$K$ with $f$ transcendental over~$K$, 
equipped with the unique derivation extending  the derivation of $K$ such 
that $f^\dagger=s$.}\/ We prove here some pre-$H$-field versions
of extension lemmas in Section~\ref{sec:exp integrals}. 
%The ordering of $K$ always extends to a field ordering 
%of~$K(f)$; see, e.g., Corollary~\ref{cor:boundrealalg}. 
%In the next section we will see that  there is in fact an ordering of~$K(f)$ that makes~$K(f)$ a pre-$H$-field extension of~$K$. In preparation for this, we deal here with certain special cases, imposing extra assumptions on $K$ and $s$. We begin with versions of Lemmas~\ref{newprop0} and~\ref{newprop, lemma 1} for pre-$H$-fields:

%\begin{lemma}\label{newprop0, preH} 
%Suppose $K$ is henselian and $v(s-a^\dagger)\in (\Gamma^{>})'$. Then there is a unique valuation on $K(f)$ making it an $H$-asymptotic extension of~$K$ with $f\sim a$. Equipped with this valuation,
%$K(f)$ is an immediate extension of~$K$, and thus, by Lemma~\ref{prehimm}, 
%a pre-$H$-field extension of~$K$ for a unique ordering of~$K(f)$. 
%\end{lemma}
%\begin{proof}
%Lemma~\ref{step4} with $s-a^\dagger$, 
%$f/a$ instead of $s$ and $f$
%provides a unique valuation of $K(f)$ making it an $H$-asymptotic extension of $K$ with $f\sim a$.
%With this valuation $K(f)$ is  an immediate extension of~$K$ by that same lemma. 
%\end{proof}

\begin{lemma}\label{newprop0, preH} 
Suppose $K$ is henselian and $v(s)\in (\Gamma^{>})'$. Then there is a unique valuation on $K(f)$ making it an $H$-asymptotic extension of~$K$ with $f\sim 1$. Equipped with this valuation,
$K(f)$ is an immediate extension of~$K$, and thus, by Lemma~\ref{prehimm}, 
a pre-$H$-field extension of~$K$ for a unique ordering of~$K(f)$. 
\end{lemma}
\begin{proof} Immediate from Lemma~\ref{step4}.
\end{proof}

\noindent
Next we establish a pre-$H$-version of 
Lemma~\ref{newprop, lemma 1}, and accordingly we now assume that~${s-a^\dagger \succ u'}$ for all $a$ and all $u\in K^{\preceq 1}$, and put
$$ S\ :=\ \big\{ v(s-a^\dagger) :\  a\in K^\times \big\}\ \subseteq\ \Gamma,$$ 
so $S<(\Gamma^>)'$.

\begin{lemma}\label{prehexpint} 
Suppose $S$ has no maximum and $\Gamma$ is divisible. With the valuation on $K(f)$ of Lemma~\ref{newprop, lemma 1}, there is a unique field ordering on $K(f)$ making it a pre-$H$-field extension of $K$ with $f>0$.
\end{lemma}
\begin{proof} We can assume $s < 0$. (If $s>0$, replace $s$ and
$f$ by $-s$ and $f^{-1}$.) 

\claim{For all $a\prec 1$ we have:
$f\succ a\Longleftrightarrow s>a^\dagger$.}

\noindent
To prove this claim, let $a\prec 1$ and set $\alpha:=va$. Then
\begin{align*}
f\ \succ\ a\ &\Longleftrightarrow\ \eta\ <\ \alpha\ \Longleftrightarrow\ \alpha_{\lambda}\ <\ \alpha \text{ eventually}\\ 
&\Longleftrightarrow\ a_{\lambda}\ \succ\ a \text{ eventually}\\
&\Longleftrightarrow\ a_{\lambda}^\dagger\ >\ a^\dagger \text{ eventually}\ 
\Longleftrightarrow\ s\ >\  a^\dagger.
\end{align*}
We order $K(f)$ as indicated in Lemma~\ref{extror} with $f$ in the role of
$y$, and show that this makes $K(f)$ a pre-$H$-field extension of $K$.  
Now $K(f)$ is pre-$\d$-valued by item~(iii) of Lemma~\ref{newprop, lemma 1}, and so with $T$ as in the proof of that lemma, it remains to show by Lemma~\ref{prehunr} that
for all $t\in T$ with $t\succ 1$ we have $t^\dagger > 0$. 

Let $t\in T$, $t\succ 1$, $t=af^j$. Consider first the case $j>0$.
Take $d\in K^\times$ with $ad^j\asymp 1$.  Then $t\asymp (f/d)^j\succ 1$ gives $f \succ d$, so
$s > d^\dagger$ by the Claim, hence $js> jd^\dagger= -a^\dagger + g$ with $g\in K$, $vg> \Psi$, and thus
$a^\dagger + js > g$. Also $v(a^\dagger + js)\in S < vg$
by Claim~2 in the proof of Lemma~\ref{newprop, lemma 1}, and thus
$t^\dagger= a^\dagger + js >0$. In the same way we treat the case that $j < 0$. The case $j=0$ is trivial.   
\end{proof}

\noindent
Since $(-f)^\dagger=f^\dagger=s$, Lemma~\ref{prehexpint} goes through with $f<0$ instead of $f>0$.

\medskip\noindent
The next result, Lemma~\ref{lemma62}, is essential for dealing with 
``Liouville closure'' in Section~\ref{sec:Liouville closed}. The case where the set 
$\big\{v(s-a^\dagger): a\in K^\times\big\}$ in Lemma~\ref{lemma62} has no maximum 
is basically contained in Lemma~\ref{prehexpint}, but the
rather lengthy proof we give here (we have no other) 
makes no use of that lemma.

\begin{lemma}\label{lemma62}
Suppose $K$ is a real closed $H$-field, $s< 0$, and
 $v(s-a^\dagger)\in \Psi^{\downarrow}$ for all $a$.
 Then there is a unique pair 
consisting of a valuation of $L=K(f)$ and a field ordering on $L$ making it a 
pre-$H$-field extension of $K$ with $f>0$. With this valuation and
ordering $L$ is an $H$-field, and we have:
\begin{enumerate}
\item[\textup{(i)}]  $vf\notin \Gamma$, $\Gamma_L=\Gamma\oplus\Z vf$, $f\prec 1$;
\item[\textup{(ii)}] $C=C_L$ and $\Psi$ is cofinal in $\Psi_{L}$;
\item[\textup{(iii)}] a gap in $K$ remains a gap in $L$;
\item[\textup{(iv)}] if $L$ has a gap which is not in $\Gamma$, then $[\Gamma]=[\Gamma_{L}]$.
\end{enumerate} 
\end{lemma}
\begin{proof}
Suppose $L=K(f)$ is equipped with 
a valuation and an ordering making~$L$ a 
pre-$H$-field extension of $K$ with $f>0$.

\claimnoskip[1]{$vf\notin \Gamma$.} 
{Otherwise $f=au$
for some $a$ and some $u\in L$ with $u\asymp 1$. For such $a$, $u$ we 
have
$s-a^\dagger = u^\dagger$, so $v(s-a^\dagger)>\Psi$, contradicting an assumption.}

\claimnoskip[2]{$f\prec 1$.} 
This is because $f\succ 1$ would give $s=f^\dagger>0$.
 
\claimnoskip[3]{$f\succ b \Longleftrightarrow s>b^\dagger$.} 
This holds by Lemma~\ref{psidecreasing1}(i) and Claim 1.

\medskip
\noindent
Claim~3 shows how $vf$ determines a cut
in $\Gamma$. 
Thus in constructing a valuation of~$L$ and ordering of $L$
with the desired properties, the three claims above leave no choice: 
we equip $L$ with the unique valuation extending the valuation of $K$ 
such that $0<vf\notin \Gamma$ realizes the 
cut in $\Gamma$ described in Claim~3 above (Lemmas~\ref{uniorvec} and \ref{lem:lift value group ext}), and with the unique ordering 
extending the ordering of $K$ in which the valuation ring of $L$ is convex, 
and with $f>0$ (Lemma~\ref{extror}). We are going to show that
with this valuation and ordering~$L$ is an $H$-field. 

Put $\eta := vf$, so $\Gamma_L = \Gamma \oplus \Z\eta$. 
We claim that  $va+j\eta>0$
if and only if 
$$ \text{\em either}\/\quad \text{(1)\hskip0.5em $va=0$ and $j>0$},\quad \text{\em or}\quad
\text{(2)\hskip0.5em $va\neq 0$ and $a^\dagger+js<0$}.$$
Since $\eta>0$, this is clear if $va=0$, or $va\ne 0$ and $j=0$. 
Assume $va\neq 0$ and $j<0$. Let $d\in K$ be a solution to
the equation $d^j=1/|a|$. 
We have $va+j\eta>0$ if and only
if $\eta < vd$, which is equivalent to $s>(-1/j)\cdot a^\dagger$, by definition of
the cut in $\Gamma$ realized by $\eta=vf$, that is, to $a^\dagger+js<0$. If $va\ne 0$,
$j>0$, one argues similarly.

We observe that for $va+j\eta\ne 0$
we have $v(a^\dagger+js)\in\Psi^\downarrow$. To see this, we may assume $j\neq 0$; take $d\in K^\times$ such that 
$d^j=1/|a|$; then 
$v(a^\dagger+js)=v(s-d^\dagger)\in\Psi^\downarrow$.
Suppose $va+j\eta\ne 0$ and $va=vb$; then $v(a^\dagger+js)=
v(b^\dagger+js)$: $a=ub$ with $u\in K$, $u\asymp 1$ gives $v(u^\dagger)=v(u')>
\Psi$, hence $v(u^\dagger)>v(b^\dagger+js)$, and $a^\dagger+js=(b^\dagger+js)+u^\dagger$, and thus
$v(a^\dagger+js) = 
v(b^\dagger+js)$, as promised. 
We can therefore extend $\psi\colon \Gamma^{\neq} \to \Gamma$
to a map $\psi_L\colon \Gamma_L^{\neq} \to \Gamma$ by
$$\psi_L(va + j\eta)\ :=\ v(a^\dagger+js) \quad\text{for $va+j\eta\neq 0$.}$$

\claimnoskip[4]{The function $\psi_L$ is decreasing on
$\Gamma_L^{>}$.}  Let $va+j\eta>vb+k\eta>0$; our job is to show that
$v(a^\dagger+js) \leq v(b^\dagger+ks)$. The above
criterion ``either (1) or (2)'' for an element of
$\Gamma_L$ to be positive yields:
\begin{enumerate}
\item[(3)] either $va=vb$ and $j>k$, or $va\neq vb$ and $a^\dagger +js < b^\dagger + ks$;
\item[(4)] either $vb=0$ and $k>0$, or $vb\neq 0$ and $b^\dagger + ks <0$.
\end{enumerate}
This leads to  four cases:

\case[1]{$va=vb=0$ and $j>k>0$.} Then $v(a^\dagger+js)=vs=v(b^\dagger+ks)$.

\case[2]{$va=vb\neq 0$, $j>k$, and $b^\dagger + ks <0$.} Then $v(b^\dagger+js)\in\Psi^\downarrow$ gives $(a/b)^\dagger
\prec b^\dagger+js$, hence
$a^\dagger + js \sim b^\dagger + js < b^\dagger+ks<0$, so
$v(a^\dagger+js) \le v(b^\dagger+ks)$.

\case[3]{$va\neq vb$,  $a^\dagger +js < b^\dagger + ks$, $vb=0$, and $k>0$.}
Then $b^\dagger+ks \sim ks <0$ and thus $v(a^\dagger+js)
\leq v(b^\dagger+ks)$.

\case[4]{$va\neq vb$,  $a^\dagger +js < b^\dagger + ks$, $vb\neq 0$, and
$b^\dagger+ks<0$.} Then clearly  $a^\dagger +js < b^\dagger + ks < 0$, and thus $v(a^\dagger+js) \leq v(b^\dagger+ks)$.

\medskip\noindent
This proves Claim 4. Since $\psi_L(k\gamma)=\psi_L(\gamma)$ for 
$\gamma\in\Gamma_L^{\neq}$, $k\neq 0$, the function
$\psi_L$ is constant on every archimedean class of
$\Gamma_L$. 

\claimnoskip[5]{$(\Gamma_L,\psi_L)$ is an $H$-asymptotic couple and $\Psi$ is cofinal in $\Psi_L:=\psi_L(\Gamma_L^{\neq})$.}  
First assume $[\Gamma]=[\Gamma_L]$. For $va+j\eta\neq 0$, let $b$ be such that $[vb]=
 [va+j\eta]$; then $\psi_L(va+j\eta)=
\psi_L(vb)=v(b^\dagger)=\psi(vb)$. Thus $(\Gamma_L,\psi_L)$
is an $H$-asymptotic couple by Lemma~\ref{extension4} and its proof. 
Next suppose $[\Gamma]\neq [\Gamma_L]$, and take any 
$\gamma_0\in\Gamma_L^{>}$ such that $[\gamma_0]
\notin [\Gamma]$.
Then we have for $j\ne 0$,
\begin{align*}  [va]<[\gamma_0]\ \Rightarrow\  \psi_L(va+j\gamma_0)&=\psi_L(\gamma_0), \quad [\gamma_0] < [va]\ \Rightarrow\ \psi_L(va+j\gamma_0)=v(a^\dagger),\\ 
\text{ so }\ \psi_L(va+j\gamma_0)\ &=\ \min\bigl\{v(a^\dagger),\psi_L(\gamma_0)\bigr\}. 
\end{align*}
Set $\Gamma_0:=\Gamma\oplus\Z\gamma_0\subseteq\Gamma_L$. Then $\bigl(\Gamma_0,\psi_L|\Gamma_0^{\neq}\bigr)$
is an $H$-asymptotic couple: apply Lemmas~\ref{exarchextclass}
and~\ref{extension5} to the cut
$C:=[\Gamma^{\neq}]^{<[\gamma_0]}$ and $\beta:=\psi_L(\gamma_0)$, using also a fact from the proof of the latter.
Since $\Gamma_L\subseteq \Q\Gamma_0$, $(\Gamma_L,\psi_L)$ is an $H$-asymptotic couple as well. The cofinality is because $v(a^\dagger+js)\in \Psi^\downarrow$ whenever $va + j\eta\ne 0$. 

\claimnoskip[6]{A gap in $(\Gamma, \psi)$ remains a gap in
$(\Gamma_L,\psi_L)$. If~$(\Gamma_L,\psi_L)$ has a gap that is not in $\Gamma$, then $[\Gamma]=[\Gamma_L]$.} 
Suppose $\beta\in \Gamma$ is a gap in $(\Gamma, \psi)$, but not in
$(\Gamma_L,\psi_L)$. This gives $\alpha\in \Gamma_L^{>}$ with
$\beta\le\psi_L(\alpha)$ or $\alpha'\le \beta$. In either case,
$0 < \alpha < \Gamma^{>}$, so $\psi_L(\alpha)\ge \Psi$,
contradicting that $\Psi_L\subseteq \Psi^{\downarrow}$ and that
in the presence of a gap in $(\Gamma,\psi)$ the set $\Psi$ cannot have a 
maximum.

For the second part of Claim~6, suppose towards a contradiction that
$(\Gamma_L,\psi_L)$ has a gap that is not in $\Gamma$, but 
$[\Gamma]\ne [\Gamma_L]$. Then $\Psi_L$ has no maximum, so 
$\Psi$ has no maximum.
By the first part of Claim~6, $(\Gamma, \psi)$ has no gap. Thus
$(\Gamma,\psi)$ has rational asymptotic integration. With $\Gamma_0$ 
as in the proof
of Claim~5, we obtain from the last statement in Lemma~\ref{extension5} that
$(\Q\Gamma_L, \psi_L)=(\Q\Gamma_0,\psi_L)$ has asymptotic integration,
contradicting that $(\Q\Gamma_0,\psi_L)$ has a gap.

\medskip\noindent
By Claim~1 and Lemma~\ref{lem:lift value group ext} we have $\res K=\res L$.
Using this fact and Claim~5 above we can prove very quickly:

\claimnoskip[7]{$L$ is $\d$-valued.}
To see this, put $T:=K^{\times}f^{\Z}$, and let $t=af^j\in T$.
Then $vt=va+ j\eta$, and if $vt\neq 0$, then $\psi_L(vt)=v(a^\dagger+js)=
v(t^\dagger)$ and so
$v(t') =vt+\psi_L(vt)$.
Now use Lemma~\ref{Lemma3.3}  and Claim~5.

\medskip\noindent
To complete the proof of the lemma, it remains to show:

\claimnoskip[8]{$L$ is an $H$-field.}
Let $t\in T$, $t\succ 1$; 
by Lemma~\ref{prehunr} it is enough to derive $t^\dagger>0$. We can assume $t\notin K$ and $t=af^j$, so $j\ne 0$. Take $d\in K^\times$ with 
$d^j=1/|a|$. First assume $j>0$. Then $va+j\eta=vt<0$ gives $\eta<vd$,
so $s>d^\dagger=-a^\dagger/j$, by the definition of the cut in $\Gamma$
realized by $\eta=vf$. Thus $t^\dagger=a^\dagger+js>0$, as required. 
Suppose $j<0$. Then $vd<vf$, and we distinguish the cases $vd>0$ (similar to the case $j>0$), $vd=0$ (where we use $s<0$ and $vs<v(d^\dagger)=v(a^\dagger)$), and $vd<0$ (where we use $v(a^\dagger)=v(d^\dagger)<vs$ and $a^\dagger>0$).
\end{proof}

%\noindent
%Lemma~\ref{lemma62} also goes through with the condition $s<0$ replaced by $s>0$, if in~(i) we change ``$f\prec 1$'' 
%to ``$f\succ 1$''.
%Both versions of Lemma~\ref{lemma62} also hold with~$f>0$ replaced by $f<0$.

\subsection*{Notes and comments} A.~Robinson~\cite{Robinson73b}
derived some Hardy field asymptotics from first-order axioms about ordered differential fields, and he showed where his axioms
fall short. The concepts of
$H$-field and pre-$H$-field, introduced in~\cite{AvdD2}, do not suffer from that defect. A variant of the notion of $H$-field with an analogue of Corollary~\ref{completinghfields} is in~\cite{Kulchinovskii98,Kulchinovskii99}.
Hahn fields $\R\(( t^\Gamma\)) $ with $H$-field derivations respecting infinite sums are studied in \cite{AvdD3, KuMa2}; 
see also the survey~\cite{Ma}.

{\sloppy
%Lemma~\ref{prehimm} slightly extends Corollary~3.2 in~\cite{AvdD2}.
Corollary~\ref{cor:H(K)} is~\cite[Corollary 4.6]{AvdD2}, but its proof
there incorrectly claims a unique field ordering
on $\operatorname{dv}(K)$ extending that of $K$ in which the valuation ring of~$\operatorname{dv}(K)$ is convex. Propositions~\ref{prop:constant field ext} and \ref{prop:Hconstants} are \cite[Theorem~3]{Rosenlicht2} and 
\cite[Proposition~9.1]{AvdD3}, respectively.
Lemma~\ref{lemma62} is a combination of Lemma~5.3 and the remarks following it in \cite{AvdD2}.
%Lemma~\ref{pU} is \cite[Lemma~3.1]{AvdD2}.

}

\section{Liouville Closed $H$-Fields}\label{sec:Liouville closed}

\noindent
A pre-$H$-field $K$ is said to be
{\bf Liouville closed\/} \label{p:Liouville closed}\index{H-field@$H$-field!Liouville closed}\index{Liouville closed $H$-field} if it is a real closed $H$-field and for all~$a\in K$ there exist
$y,z\in K$ such that $y'=a$ and $z\ne 0$, $z^\dagger =a$; note that then
any equation $y'+ay=b$ with $a,b\in K$ has a solution in $K^\times$. 

\begin{examples}
The $H$-field~$\mathbb{T}$ is Liouville closed: see Appendix~\ref{app:trans}. 
A Hardy field containing $\R$ as a subfield is Liouville closed if and only if it is real closed, closed under integration and closed under exponentiation. 
For a Hardy field $K\supseteq \R$ we define~$\operatorname{Li}(K)$ as the smallest
Hardy field extension of $K$ that is real closed and 
closed under integration and exponentiation; thus 
$\operatorname{Li}(K)$ is the smallest
Liouville closed Hardy field containing $K$. 
If~$K$ is a Liouville closed $H$-field and $\phi\in K^{>}$, then $K^{\phi}$ is a Liouville closed $H$-field. 
\end{examples}

\noindent
Let $K$ be a Liouville closed $H$-field. Then $K$ is closed under integration (and thus under logarithms)
as defined in Section~\ref{Asymptotic-Fields-Basic-Facts}.
Also, $K$ is \textit{closed under powers}\/ in the following sense:
for all $c\in C$ and $f\in K^\times$ there exists $y\in K^\times$ such that $y^\dagger=cf^\dagger$; such $y$ behaves like $f^c$.
The $H$-asymptotic couple $(\Gamma, \psi)$ of $K$ is closed in the sense of Section~\ref{sec:cac}: it has asymptotic integration,
$\Gamma$ is divisible, and~$\Psi$ is downward closed. Thus $(\Gamma,\psi)$ has a contraction map $\chi\colon \Gamma^{<} \to \Gamma^{<}$ as defined in Section~\ref{AbstractAsymptoticCouples}. This 
contraction map is induced by a logarithm map on $K^{>}$ as follows:
choosing for each
$a\in K^{>}$ a ``logarithm'' \index{logarithm!on an asymptotic field} $\Log(a)\in K$, that is, 
$\Log(a)'=a^{\dagger}$, we have 
$$\chi(\alpha)\ =\ v\big(\Log(a)\,\big)\ \text{ 
whenever $\alpha\in \Gamma^{<}$, $\alpha=va$, $a\in K^{>}$.}$$ 

\begin{example}
For $K=\T$ and $a\in \T^>$ we may take $\Log(a)=\log(a)$; if $a>\R$, then the sequence $a,\log a, \log\log a,\dots$ of iterated logarithms of $a$ is coinitial in $\T^{>\R}$, by results in Appendix~\ref{app:trans}. Thus for any $\alpha\in \Gamma_{\T}^{<}$ the sequence $\big(\chi^n(\alpha)\big)$ is cofinal in $\Gamma_{\T}^{<}$.  
\end{example}

\noindent
For use in Section~\ref{sec:special sets} we mention two easy results:

\begin{lemma}\label{intersectliou} Let $K$ and $L$ be Liouville closed $H$-subfields
of an $H$-field $M$ such that $C_L=C_M$. Then $K\cap L$ is a Liouville
closed $H$-subfield of $M$. 
\end{lemma}
\begin{proof} It is easy to check that $K\cap L$ is a real closed 
$H$-subfield
of $M$ with constant field $C=C_K$. Let $a\in K\cap L$. First, take
$f\in K$ and $g\in L$ such that $f'=a=g'$. Then $f-g\in C_M$, so
$f-g\in L$, and thus $f\in K\cap L$. Next, take $f\in K^\times$ and
$g\in L^\times$ such that $f^\dagger=a=g^\dagger$. Then
$f/g\in C_M$, so $f/g\in L$, and thus $f\in K\cap L$.
\end{proof}

\noindent
Using Lemma~\ref{intersection} and Corollary~\ref{cor:cor:AS} we obtain likewise:

\begin{lemma}\label{intersectlioumany} If $K$ is an
$H$-field and $(K_i)_{i\in I}$ $(I\ne \emptyset)$ is a family of Liouville closed $H$-subfields of $K$, all with the same constants as $K$, then $\bigcap_i K_i$ is a Liouville closed $H$-subfield of $K$.
\end{lemma}

\subsection*{Completion} Recall that by Corollary~\ref{completingasymp} the completion of a $\d$-valued field is 
$\d$-valued.
In this subsection we show that taking the completion preserves some
further properties, like being a Liouville closed $H$-field. 
%Note first that if $K\subseteq L$ is an extension of asymptotic fields and $B$ is a subset of
%$L$ such that $\Gamma$ is cofinal in $\Gamma_L$ and for each $b\in B$ and~$\gamma\in\Gamma$ there is $a\in K$ with
%$v(a-b)>\gamma$, then $K$ is dense in the valued differential subfield~$K\< B \>$ of $L$; this follows from the continuity of the derivation of $L$ and Corollary~\ref{cor:completion valued fields, 1}.  

\begin{cor}\label{cor:Kc closed under integration}
Suppose the $\d$-valued field $K$ is closed under integration. Then
the completion $K^{\operatorname{c}}$ of $K$ is also closed under integration.
\end{cor}

\begin{proof}
Note that $(\Gamma, \psi)$ has asymptotic integration, in particular $\Gamma\neq\{0\}$. 
Let ${b\in K^{\operatorname{c}}}$. To get
$a\in K^{\operatorname{c}}$ with $a'=b$ we may subtract from $b$ an element of $K$ and arrange in this way that $vb>\Psi$. Take
a c-sequence $(b_{\rho})$ in $K$ with $b_{\rho}\to b$ and
$vb_\rho> \Psi$ for all~$\rho$. Take for each $\rho$ the unique $a_{\rho}\in K$ with $a_{\rho}'=b_{\rho}$ and $a_{\rho}\prec 1$.
Then~$(a_{\rho})$ is clearly a c-sequence in $K$, which gives 
$a\in K^{\operatorname{c}}$ with $a_{\rho} \to a$, and thus, taking de\-ri\-va\-tives,~$a'=b$.
\end{proof}

\begin{cor}\label{cor:Kc closed under exp integration}
Suppose the $\d$-valued field $K$ satisfies 
$(K^\times)^\dagger=K$. Then the completion $K^{\operatorname{c}}$ of $K$ also satisfies $((K^{\operatorname{c}})^\times)^\dagger=K^{\operatorname{c}}$.
\end{cor}
\begin{proof} Let $b\in K^{\operatorname{c}}$. To get nonzero
$a\in K^{\operatorname{c}}$ with $a^\dagger=b$ we may subtract from $b$ an element of $K$ and arrange in this way that $vb>\Psi$. Take
a c-sequence $(b_{\rho})$ in $K$ with $b_{\rho}\to b$ and
$vb_\rho> \Psi$ for all $\rho$. Take for each $\rho$ the unique $a_{\rho}\in \smallo$ with $(1+a_{\rho})^\dagger=b_{\rho}$.
Now $(1+\delta)^\dagger-(1+\epsilon)^\dagger\sim (\delta-\epsilon)'$ for distinct $\delta, \epsilon\in \smallo$, so
$(a_{\rho})$ is a c-sequence in $K$. This gives 
$a\in K^{\operatorname{c}}$ with $a_{\rho} \to a$, and thus $a\prec 1$ and $(1+a)^\dagger=b$.
\end{proof}

\noindent
Corollaries~\ref{cor:completion of rcf}, \ref{completinghfields},  \ref{cor:Kc closed under integration}, and \ref{cor:Kc closed under exp integration}  yield  
a result used in Section~\ref{sec:reldifhens}:

\begin{lemma}\label{lioucompletion} If $K$ is a 
Liouville closed $H$-field, then so is its completion $K^{\operatorname{c}}$. 
\end{lemma}

\subsection*{Liouville extensions} Let $K$ be a differential field. A {\bf Liouville extension} of $K$
is a differential field extension~$L$ of $K$ such that $C_L$ is algebraic over $C$ and for each $a\in L$ there are
$t_1,\dots,t_n\in L$ with $a\in K(t_1,\dots,t_n)$ and for $i=1,\dots,n$, 
\begin{enumerate}
\item $t_i$ is algebraic over $K(t_1,\dots,t_{i-1})$, or
\item $t_i'\in K(t_1,\dots,t_{i-1})$, or 
\item $t_i\ne 0$ and $t_i^\dagger \in K(t_1,\dots,t_{i-1})$.
\end{enumerate}

\noindent
We leave the routine proofs of the next two lemmas to the reader.

\begin{lemma}\label{transliou} Let $M|L$ and $L|K$ be differential field extensions. If $M|K$ is a Liou\-ville extension, then so is $M|L$.
If $M|L$ and $L|K$ are Liouville extensions, then so is~$M|K$.
\end{lemma}

\begin{lemma}\label{lliou} Let $M|K$ be a differential field extension such that $C_M$ is algebraic over $C$. Then the subfield of $M$ generated by any nonempty set of intermediate Liouville extensions of $K$ is also a Liouville extension of $K$. Thus there exists a largest differential subfield of $M$ that contains $K$ and is a Liouville extension of $K$. 
\end{lemma}

\noindent
The assumption in Lemma~\ref{lliou} that $C_M$ is algebraic over $K$ cannot be dropped: for a transcendental real number $r$, the Hardy subfields $\Q(x)$ and $\Q(x+r)$ of $\R(x)$ are both Liouville extensions of $\Q$ with $\Q$ as common field
of constants, but the Hardy subfield
$\Q(x, x+r) = \Q(x,r)$ of $\R(x)$ is not a Liouville extension of $\Q$.

\medskip\noindent
Let $K$ be a differential field, and 
$|K|:=\text{cardinality of $K$.}$ We observe:

\begin{lemma}\label{cardinality} 
Suppose $L$ is a Liouville extension of $K$.
Then $|L| = |K|$.
\end{lemma}
\begin{proof} Define a chain of differential subfields
$K=K_0 \subseteq K_1 \subseteq K_2 \subseteq \cdots$   
of $L$:
\[ K_{n+1} = \left\{ \begin{array}{ll}
                     \text{algebraic closure of } K_n  \text{ in } L
  & \mbox{for $ n\equiv 0 \pmod 3$} \\
                  K_n\bigl(\{a\in L : a'\in K_n\}\bigr) 
   & \mbox{for $ n\equiv 1 \pmod 3$}\\
                     K_n\bigl(\{a\in L^{\times} : a^\dagger \in K_n\}\bigr) 
     & \mbox{for $ n\equiv 2 \pmod 3$}
\end{array} \right. \] 
Clearly $\abs{K_n}=\abs{K}$ 
for all $n$ (by induction), and
$L=\bigcup_n K_n$, so $\abs{L}=\abs{K}$. 
\end{proof}

\begin{lemma}\label{Minimality, 1}
Let $K$ be a Liouville closed $H$-field. Then $K$ has no proper Liouville extension with the same 
constants as $K$.
\end{lemma}
\begin{proof}
Suppose $L$ is a proper Liouville extension of the 
differential field $K$ with the same constants as $K$. Up to
$K$-isomorphism the only
proper algebraic extension field of $K$ is $K(\imag)$ with $\imag^2=-1$, 
and as a differential field extension of $K$ it contains the constant 
$\imag\notin C$.
Hence $L$ must contain a solution $y\notin K$ to an equation $y'=a$ 
with $a\in K$, or a solution $z\notin K$ with $z\ne 0$ to an equation 
$z^\dagger=b$ with $b\in K$. But given $y$ as above, take $y_0\in K$ with
$y_0'=a$, and then $y-y_0\in C_L\setminus C$, 
a contradiction. Similarly, given $z$ as above, take $z_0\in K^{\times}$
with $z_0^\dagger = b$, and note that then $z/z_0\in C_L\setminus C$, 
a contradiction. 
\end{proof}

\subsection*{Liouville closure}
Let $K$ be an $H$-field. A {\bf Liouville closure} of $K$ is a Liouville closed $H$-field extension~$L$ of $K$ such that $L$ is also a Liouville extension of $K$. Note that if $L$ is a Liouville closure of $K$, then by Lemma~\ref{Minimality, 1} there is no proper $H$-subfield of $L$ that contains~$K$ and is Liouville closed.

\begin{lemma}\label{lem:max Liouville}
Let $K\subseteq M$ be an extension of $H$-fields such that $M$ is Liouville closed and $C_M$ is algebraic over $C$. Then $M$ has
a unique $H$-subfield $L\supseteq K$ that is a Liouville closure of $K$. 
\end{lemma}
\begin{proof} Lemma~\ref{lliou} gives a largest Liouville extension $L$ of $K$ in $M$. Then $C_L=C_M$, so $L$ is an $H$-subfield of $M$. It is also clear that $L$ is real closed, and that for any $a\in L$ and $y,z\in M$ with $y'=a$ and $z\ne 0$, $z^\dagger=a$ we have $y,z\in L$.
Thus $L$ is a Liouville closure of $K$; uniqueness follows from
Lemma~\ref{Minimality, 1}. 
\end{proof}

\begin{cor}\label{cor:HardyLiouvilleClosure} Let $K$ be a Hardy field containing $\R$ as a subfield. Then $\operatorname{Li}(K)$
is a Liouville closure of $K$. 
\end{cor}
\begin{proof} Use $C_{\operatorname{Li}(K)}=\R$ and Lemma~\ref{lem:max Liouville}. 
\end{proof}

%\noindent
%In general, the $H$-field $\hat K$ in Proposition~\ref{CharLiouvilleClosures} is not unique:

%\begin{example}
%The Hardy
%fields $K(x)$ and $K(x+\pi)$ are both 
%Liouville extensions of $K=\Q$ (with $\Q$ as field
%of constants), but $K(x, x+\pi) = K(x,\pi)$ is not, since the constant $\pi$ is not
%algebraic over $\Q$; hence there is more than one Liouville closure~$\hat K$ of $K$ with $K\subseteq\hat K\subseteq\operatorname{Li}(\R)$.
%\end{example}

%However, if $K$, $L$ are as in the previous proposition, and  $C_L$ is algebraic over $C$,
%then by Lemma~\ref{intersectliou} and Proposition~\ref{CharLiouvilleClosures} there is a {\it unique}\/
%$H$-field $\hat K$ such that $K\subseteq \hat K\subseteq L$
%and $\hat K$ is a Liouville closure of $K$, which we call the {\bf Liouville closure of $K$ inside~$L$.}
%The constant field of $\hat K$ is the real closure of the constant field of~$K$.

\subsection*{The main result about Liouville closures} 
%This is the following:

\begin{theorem}\label{thm:Liouville closures}
Let $K$ be an $H$-field. Then one of the following occurs: 
\begin{enumerate}
\item[\rom{(I)}] $K$ has exactly one Liouville closure up to isomorphism over $K$,
\item[\rom{(II)}] $K$ has exactly two Liouville closures up to isomorphism over $K$.
\end{enumerate}
\end{theorem}

\noindent
Moreover, for any $H$-field $K$ we have:
\begin{enumerate}
\item[(1)] If no Liouville 
$H$-field extension of $K$ has a gap, then $K$ falls under 
Case~(I).
Special case: $H$-subfields of $\T$ properly containing 
$\R$ fall under Case~(I).
\item[(2)]  If $K$ is grounded, then $K$ falls under Case~(I).
\item[(3)] Suppose $K$ has a gap $\gamma$. Then $K$ falls under Case~(II): in one Liouville closure $L_1$ of $K$, all 
$s\in K$ with $vs=\gamma$ have the form $b'$ with  
$b\succ 1$, while in another Liouville closure $L_2$ of $K$ 
all $s\in K$ with $vs=\gamma$ have the form $b'$
with $b\prec 1$.
\end{enumerate}
We prove Theorem~\ref{thm:Liouville closures} and statements (1), (2), (3) later in this section.

\subsection*{Liouville towers}
{\em In this subsection $K$ is an $H$-field}. A
{\bf Liouville tower on $K$} is a strictly increasing chain~$(K_\lambda)_{\lambda \le \mu}$ of $H$-fields, indexed by the 
ordinals less than or equal to some ordinal~$\mu$,
such that 
\begin{enumerate}
\item[(1)] $K_0 = K$;
\item[(2)] if $\lambda$ is a limit ordinal, $0<\lambda \le \mu$, then 
$K_\lambda = \bigcup_{\iota<\lambda} K_\iota$;
\item[(3)] for $\lambda < \lambda +1 \le \mu$, {\em either\/}
\begin{enumerate}
\item[(a)] $K_{\lambda}$ is not real closed and $K_{\lambda + 1}$ is a real closure of $K_\lambda$,
\end{enumerate}
{\em or\/} $K_\lambda$ is real closed,
$K_{\lambda + 1} = K_\lambda(y_{\lambda})$ with $y_\lambda\notin K_\lambda$
(so $y_\lambda$ is transcendental over $K_\lambda$), and one of
the following holds, with $(\Gamma_\lambda, \psi_\lambda)$
the asymptotic
couple of $K_\lambda$ and $\Psi_\lambda := \psi_\lambda(\Gamma_\lambda^{\neq})$: 
\begin{enumerate}
\item[(b)] $y_\lambda'= s_\lambda \in K_\lambda$ with $y_\lambda\prec 1$ and 
$v(s_\lambda)$ is a gap in $K_\lambda$,
\item[(c)] $y_\lambda'= s_\lambda \in K_\lambda$ with $y_\lambda\succ 1$ and 
$v(s_\lambda)$ is a gap in $K_\lambda$,
\item[(d)] $y_\lambda'= s_\lambda \in K_\lambda$ with
$v(s_\lambda)= \max \Psi_\lambda$,
\item[(e)] $y_\lambda'= s_\lambda \in K_\lambda$ with $y_{\lambda}\prec 1$,
$v(s_\lambda)\in (\Gamma_\lambda^{>})'$, and 
$s_\lambda\neq\varepsilon'$ for all 
$\varepsilon\in K_\lambda^{\prec 1}$,
\item[(f)] $y_\lambda'= s_\lambda \in K_\lambda$ such that
$S_\lambda :=\bigl\{v(s_\lambda - a'): a\in K_\lambda\bigr\} <
(\Gamma_\lambda^{>})'$, and 
$S_\lambda$ has no largest element,
\item[(g)] $y_\lambda^\dagger = s_\lambda \in K_\lambda$ with 
$y_\lambda\sim 1$,
$v(s_\lambda)\in (\Gamma_\lambda^{>})'$, and 
$s_\lambda\neq a^\dagger$ for all $a\in K_\lambda^\times$, 
\item[(h)]  
$y_\lambda^\dagger= s_\lambda \in K_\lambda^{<}$ with 
$y_\lambda > 0$, and $v(s_\lambda-a^\dagger)\in\Psi_\lambda^\downarrow$
for all $a\in K_\lambda^\times$.
\end{enumerate}
\end{enumerate}
The $H$-field $K_\mu$ is called the {\bf top} of the tower 
$(K_\lambda)_{\lambda\leq\mu}$.
Note that (a), (b), (c), (d) correspond to Corollaries~\ref{cor:H-field rc}, \ref{extensionpdvfields1, cor}, \ref{variant41, cor}, \ref{pre-extas2, pre-H}, respectively,
and (e), (f), (g), (h) to 
Lemmas~\ref{var43}, \ref{var51}, \ref{newprop0, preH}, \ref{lemma62}, respectively.

\begin{lemma}\label{towerfacts}
Let a Liouville tower on $K$ as above be given. Then: 
\begin{enumerate}
\item[\textup{(i)}] $K_\mu$ is a Liouville extension of $K$;
\item[\textup{(ii)}] the constant field $C_\mu$ 
of $K_\mu$ is a real closure of $C$ if $\mu >0$; 
\item[\textup{(iii)}] $\abs{K_\mu}=\abs{K}$, 
hence $\mu < \abs{K}^+$, where $|K|^+$ is the least cardinal 
$>|K|$.
\end{enumerate}
\end{lemma}
\begin{proof}
For (i) and (ii), use results of
Sections~\ref{sec:integrals} and \ref{sec:exp integrals} to show by induction on~$\lambda\le \mu$ that
$K_\lambda$ is a Liouville extension of $K$, and that the constant field of
$K_\lambda$ is a real closure of $C$ for $\lambda >0$. Item (iii) follows
from (i) by Lemma~\ref{cardinality}.
\end{proof}

\noindent
By Lemma~\ref{towerfacts}(iii) there is 
a {\em maximal\/}
Liouville tower $(K_\lambda)_{\lambda \le \mu}$ on $K$, ``maximal'' meaning
that it cannot be extended to a Liouville tower
$(K_\lambda)_{\lambda \le \mu +1}$ on $K$.

\begin{lemma}\label{CharLiouvilleClosures, prep}
Let $L$ be the top of a maximal Liouville tower on $K$. Then $L$ is Liouville closed, and hence a Liouville closure
of $K$. 
\end{lemma}
\begin{proof} Using (a) and~\ref{cor:H-field rc} we see
that $L$ is real closed. Likewise, $L$ has no gap by (b) and~\ref{extensionpdvfields1, cor}, and $L$ is not grounded by (d) and~\ref{pre-extas2, pre-H}. So $L$ has asymptotic integration. Then (e) with~\ref{var43}, and (f) with~\ref{var51} show
that $L$ is closed under integration. In the same way, (g) with
~\ref{newprop0, preH} and  (h) with \ref{lemma62} show that
$(L^\times)^\dagger=L$. Thus $L$ is Liouville closed. 
\end{proof}

\noindent
By the last two lemmas, each $H$-field has a Liouville closure. The following result (where $C_M$ is not necessarily algebraic over $C$) has a straightforward proof:

\begin{lemma}\label{internalliouvilletower} Let $M$ be a Liouville closed $H$-field extension of $K$ and $(K_{\lambda})_{\lambda\le \mu}$ a Liouville tower on $K$. Suppose this tower is in $M$ 
$($consists of $H$-subfields of $M)$, and maximal in $M$, that is, 
it cannot be cannot be extended to a Liouville tower
$(K_\lambda)_{\lambda \le \mu +1}$ on $K$ in $M$.
Then $(K_{\lambda})_{\lambda\le \mu}$ is a maximal Liouville tower on $K$.
\end{lemma} 

\noindent
From Lemmas~\ref{Minimality, 1},~\ref{CharLiouvilleClosures, prep}, and~\ref{internalliouvilletower} we obtain:

\begin{cor} A Liouville closed $H$-field extension $L$ of $K$
is a Liouville closure of $K$ iff no proper $H$-subfield of $L$ containing $K$ is Liouville closed.
\end{cor}

\subsection*{Uniqueness of Liouville closure}  
The uniqueness properties in Sections~\ref{sec:integrals}, \ref{sec:exp integrals}, and \ref{sec:H-fields}, such as Lemma~\ref{prehimm}, together with Lemma~\ref{lem:max Liouville}, yield:

\begin{lemma}\label{gap}
Let $K$ be an $H$-field and $(K_\lambda)_{\lambda \le \mu}$ a 
Liouville tower
on $K$ such that no $K_\lambda$ with $\lambda < \mu$ has a gap.
Then every embedding of $K$ into a Liouville closed $H$-field $L$
extends to an embedding of $K_\mu$ into $L$. If
$K_\mu$ is also Liouville closed, then~$K_\mu$ is up to isomorphism over $K$ the unique Liouville closure of $K$.
\end{lemma} 

\noindent
Combining Lemmas~\ref{CharLiouvilleClosures, prep} and~\ref{gap}
yields item (1) after Theorem~\ref{thm:Liouville closures}:

\begin{cor}\label{noligaps} If no Liouville $H$-field extension of the $H$-field
$K$ has a gap, then $K$ has up to isomorphism over $K$ a unique Liouville closure.
\end{cor}

\noindent
To apply this to $H$-subfields of $\T$ we first note:

\begin{lemma}\label{nogapsinT} No $H$-subfield of $\T$ properly containing~$\R$ has a gap. 
\end{lemma} 
\begin{proof} More generally, let $\Delta\ne \{0\}$ be a subgroup of $\Gamma_{\T}$ with $\psi(\Delta^{\ne})\subseteq \Delta$,
where $\psi=\psi_{\T}$;
we claim that then the asymptotic couple 
$\big(\Delta, \psi|\Delta^{\ne}\big)$ does not have a gap. This is clear if this asymptotic couple is grounded. Suppose it is ungrounded,
and take $\alpha\in \Delta^{<}$.
Then $\chi^n(\alpha)\in \Delta^{<}$
for all $n$, by Lemma~\ref{Deltacontract}, so $\Delta^{<}$ is cofinal in $\Gamma_{\T}^{<}$, and thus $\big(\Delta, \psi|\Delta^{\ne}\big)$ has asymptotic integration. 
\end{proof}

\noindent
Combining~\ref{CharLiouvilleClosures, prep},
~\ref{internalliouvilletower}, ~\ref{gap}, and~\ref{nogapsinT} gives:

\begin{cor}\label{cor:no two lcs in T}
If $K$ is an $H$-subfield of $\T$ properly containing~$\R$,
then any two Liouville closures of $K$ are isomorphic 
over $K$.
\end{cor}

\begin{cor}\label{cor:Liouville closures}
Let $K\supseteq\R$ be a Hardy field and $e\colon K\to\T$ an $H$-field embedding with $e|\R=\id_\R$. Then $e$ extends to
an $H$-field embedding $\operatorname{Li}(K)\to\T$.
\end{cor}
\begin{proof} By Corollary~\ref{cor:HardyLiouvilleClosure}, $\operatorname{Li}(K)$ is a Liouville closure of $K$. If $K=\R$, then $\R(x)\subseteq \operatorname{Li}(K)$ and we can extend $e$ to $\R(x)$ to reduce to the case that $K$ properly contains $\R$. In that case Corollary~\ref{cor:no two lcs in T} applies.
\end{proof}

\noindent
Let $K$ be an $H$-field.
To construct useful Liouville towers on $K$ with the property stated in 
Lemma~\ref{gap}, let  $\Lambda \subseteq \bigl\{\rom{(a)}, \rom{(b)}, \dots, \rom{(h)}
\bigr\}$ with $\rom{(a)}\in
\Lambda$. Then the definition of 
{\em $\Lambda$-tower on $K$\/} is identical to that of
{\em Liouville tower on $K$}, except that in clause (3) of
that definition only the items from $\Lambda$ occur.

\begin{lemma}\label{hnogapno}
Let $K$ be a real closed $H$-field without a gap, and let 
$(K_\lambda)_{\lambda \le \mu}$ be a $\Lambda$-tower on $K$ with $\rom{(h)}\notin\Lambda$.
Then no $K_\lambda$ with $\lambda \leq \mu$ has a gap.
\end{lemma}
\begin{proof}
By induction on $\lambda\leq\mu$ one shows easily that the asymptotic couple of $K_{\lambda}$ is divisible without a gap, or grounded. 
\end{proof}

\noindent
Lemma~\ref{hnogapno} indicates that the Liouville extension of type (h) considered in Lem\-ma~\ref{lemma62} is 
special:
while none of 
the extensions of type (b)--(g) can produce a gap that wasn't 
already there, (h)  can {\em create\/} a gap; Section~\ref{sec:specialH} below contains a concrete example.
In Section~\ref{sec:special cuts} we show that \textit{every}\/ real closed $H$-field with asymptotic integration has an immediate $H$-field extension $K$ with an element~$s$ satisfying the hypotheses of Lemma~\ref{lemma62}, and such that for any nonzero $y$ in any 
$H$-field extension of $K$
with $y^\dagger=s$, the pre-$H$-field $K(y)$ has a gap. (Thus no spherically complete $H$-field can be Liouville closed.) This explains the perhaps curious
arrangement of the results leading to the proof of the main theorem below.  

\medskip
\noindent
Next we turn to statement (2) after Theorem~\ref{thm:Liouville closures}:
 
\begin{prop}\label{MainResult1}
If $K$ is a grounded $H$-field, then all Liouville closures of $K$ are isomorphic over $K$.  
\end{prop}

\noindent
Towards its proof we first show:

\begin{lemma}\label{Approx1}
Let $K$ be a grounded $H$-field.
There exists a Liouville tower on $K$ with top $L$ such that:
\begin{enumerate}
\item[\textup{(i)}] every $H$-field in the tower, in particular 
$L$, is grounded;
%\textup{(}and hence has no gap\textup{)}; 
%\item[\textup{(ii)}] $L$ is grounded; and
\item[\textup{(ii)}] for every $a\in K$ there exist $y,z\in L$ with $y'=a$ and $z\neq 0$, $z^\dagger=a$.
\end{enumerate}
\end{lemma}
\begin{proof} Let $\Lambda :=\bigl\{\rom{(a), (e), (f), (g), (h)}\bigr\}$. 
Take a maximal 
$\Lambda$-tower $(K_\lambda)_{\lambda \le \mu}$ on~$K$. 
%(Here ``maximal''
%means that we cannot extend this tower to an $\alpha$-tower  
%$(K_\lambda)_{\lambda \le \mu +1}$ on $K$.)
Induction on $\lambda$ using Lemmas~\ref{var43}, \ref{var51}, \ref{newprop0, preH} and \ref{lemma62} shows that each~$\Psi_\lambda$ has maximum
$\max\Psi$. %In particular, no $K_\lambda$ has a gap. 
Maximality 
with respect to (a), (g), (h) and Lemmas~\ref{newprop0, preH} and \ref{lemma62} yield: $K_{\mu}$ is real closed and
for all $a\in K$ there is
$z\in K_{\mu}^\times$ with $z^\dagger=a$. 
 Take $s \in K_{\mu}$ with $vs=\max\Psi$. Lemma~\ref{pre-extas2} and the remarks following it give an $H$-field extension
$L :=K_{\mu}(y)$ of $K_{\mu}$ such that $y$ is
transcendental over~$K_\mu$ and $y'=s$. 
Then~$\Psi_L$ again has a maximum, namely 
$\psi_L(vy)>\max\Psi$. %In particular $L$ has no gap.
It only remains to show that $K\subseteq \der L$. Suppose $t\in K$ 
and $t\notin \der K_{\mu}$. Then the maximality property of the
tower with respect to (a), (e), (f) together with Lemma~\ref{lem:remark5.1(2)}(i) gives
 $\max\Psi=v(t-a')$ for some $a\in
K_{\mu}$. For such $a$ we have
$t-a'=cs+d$ with $c\in C_\mu$ and $d\in K_{\mu}$, 
$vd>\max\Psi$. Then $d=e'$ with $e\in K_{\mu}$, and so $t=(a+cy+e)'$ with $a+cy+e\in L$. 
\end{proof}

\noindent
Let $K$ be a grounded $H$-field, and let $\ell(K)$ be the real closure of an $H$-field extension~$L$ of 
$K$ as in Lemma~\ref{Approx1}. Then $\Psi_{\ell(K)}= \Psi_L$, so 
$\ell(K)$ 
is grounded as well. Thus we can iterate this operation, and 
form $\ell^2(K) := \ell\big(\ell(K)\big)$, and so on. Taking the union of the 
increasing sequence of $H$-fields $\ell^n(K)$ built in this way, and applying
Lemma~\ref{gap}, we obtain Proposition~\ref{MainResult1}.

\medskip
\noindent
We now prove statement (3) following Theorem~\ref{thm:Liouville closures}:

\begin{prop}\label{MainResult2}
Let $K$ be an $H$-field with a gap $\gamma\in\Gamma$. Then $K$ has 
Liouville closures $L_1$ and $L_2$, such 
that any embedding of $K$ into a Liouville closed $H$-field~$M$ extends to an embedding of $L_1$ or of $L_2$ into $M$,
depending on whether the image of~$\gamma$ in $\Gamma_M$ lies in
$(\Gamma_M^{<})'$ or in 
$(\Gamma_M^{>})'$. Each Liouville closure of $K$
is $K$-isomorphic to $L_1$ or to $L_2$, but $L_1$ and $L_2$ are not
$K$-isomorphic. 
\end{prop}

\begin{proof}  Take $s\in K$ such that
$vs=\gamma$.  Let $K_1:=K(y_1)$ and $K_2:=K(y_2)$ be $H$-field extensions
of $K$ with $y_i$
transcendental over $K$ and $y_i'=s$, for $i=1,2$, such that
$y_1\succ 1$ and $y_2\prec 1$.  (Such $K_i$ exist by 
Corollaries~\ref{extensionpdvfields1, cor} and \ref{variant41, cor}.) Then both~$K_1$ and $K_2$
 are grounded. For $i=1,2$ let $L_i$  be a Liouville closure
of $K_i$.  Let an embedding of $K$ into a
Liouville closed $H$-field $M$ be given. If the image of $\gamma$ in~$\Gamma_M$
lies in $(\Gamma_M^{<})'$, then we can extend
that embedding to an embedding of $K_1$ into $M$, and hence by Proposition~\ref{MainResult1}, to an embedding of $L_1$ into $M$. If the image of $\gamma$ in~$\Gamma_M$
lies in $(\Gamma_M^{>})'$, then we can similarly
extend
that embedding to an embedding of $L_2$ into~$M$. It is now routine to
show that $L_1$ and $L_2$ as defined here have all the properties claimed 
in the proposition.   
\end{proof}

\begin{example}
Let $K$ be an ordered field, and equip $K$ with the trivial derivation
and trivial valuation. Then $K$ is an $H$-field with $\Gamma=\{0\}$, and
has gap $0=v(1)$.  The two Liouville closures $L_1$ and $L_2$ of $K$
in Proposition~\ref{MainResult2} satisfy $0\in \Psi_{L_1}$ and $\Psi_{L_2}<0$. Replacing 
the derivation $\der$ of $L_1$ by a suitable multiple
$a\der$,  $a\in L_1^>$, we obtain a $K$-isomorphic copy of $L_2$.
(For a more interesting $H$-field with a gap, see Section~\ref{sec:specialH}.)
%A more interesting $H$-field with a gap is the $H$-field $L$ constructed in the ``Example'' after Lemma~\ref{gap}, with gap $v(y)$. 
\end{example}

\noindent
Propositions~\ref{MainResult1} and~\ref{MainResult2} above concern two special cases, and we now
turn to the general situation in the course of proving
 Theorem~\ref{thm:Liouville closures}. Let $K$ be an $H$-field. Take a maximal
Liouville tower $(K_\lambda)_{\lambda \le \mu}$ on $K$. Then $K_\mu$
is a Liouville closure of~$K$.  We have two cases:

\begin{enumerate} 
\item[(A)] No $K_\lambda$ has a gap. Then $K$ falls under Case~(I) by 
Lemma~\ref{gap}.
\item[(B)] Some $K_\lambda$ in the tower has a gap. Take $\lambda$
minimal with this property. Let $L_1$ and~$L_2$ be the two Liouville
closures of $K_\lambda$ as in Proposition~\ref{MainResult2}. Then~$L_1$ and~$L_2$ are also Liouville closures of $K$. Given any embedding
of $K$ into a Liouville closed $H$-field $M$, we can first extend
it to an embedding of $K_\lambda$ into $M$, and then by Proposition~\ref{MainResult2} to an embedding of $L_1$ or $L_2$ into $M$. Applying this to the inclusion of $K$ into any Liouville closure $M$ of $K$,
it follows that $L_1$ or $L_2$ is $K$-isomorphic to $M$.
Thus if $L_1$ and $L_2$ are $K$-isomorphic, then~$K$ falls under Case~(I), 
and otherwise $K$ falls under Case~(II).
\end{enumerate}

\noindent
This yields the following more precise version of  Theo\-rem~\ref{thm:Liouville closures}:

\begin{theorem}\label{MainResult3}
Let $K$ be an $H$-field. Then $K$ has at least one and at most two Liouville
closures, up to isomorphism over $K$. Any embedding of $K$ into a
Liouville closed $H$-field $M$ extends to an embedding of some Liouville
closure of $K$ into $M$. Moreover, if $K$ has two Liouville
closures, not isomorphic over $K$, then $K$ has a Liouville $H$-field 
extension $L$ with a gap such that $L$ embeds over $K$ into 
any Liouville closed $H$-field extension of $K$.
\end{theorem}

%\begin{cor}
%Let $K$ be an ungrounded $H$-field, and let $L$ be a 
%Liouville closure of $K$. Then
%$L$ is up to $K$-isomorphism
%the only Liouville closure of $K$ if and only if
% $\Gamma^{>}$ is coinitial in $\Gamma_L^{>}$.
%\end{cor}

\subsection*{Notes and comments}
The differential-algebraic notion of {\em Liouville extension\/} was motivated by Liouville's work~\cite{Liouville3, Liouville4} on explicit solutions of second-order linear differential equations; see  Kolchin~\cite[p.~5]{Kolchin48}  
%Ritt~\cite{Ritt27}, \cite[Chapter~VI]{Ritt48}, and
and Rosenlicht-Singer~\cite{Rosenlicht-Singer}. 

Rosenlicht~\cite[Theorem~3]{Rosenlicht7} considers
the Liouville closure $\operatorname{Li}(K)$ of the Hardy field $K=\R$; it contains Hardy's field
of logarithmico-exponential functions~\cite{Hardy11}.
Liouville closed $H$-fields and the results of this section are  from~\cite{AvdD2} (with errata at the end of \cite{AvdD3}), 
except for
Lem\-ma~\ref{lioucompletion}, which is \cite[Lemma~10.2]{AvdD3}. 
We point out that the material of the present section demands the $H$-field setting:
Hardy fields and $H$-subfields of $\T$ tend to obscure the fork in the road
caused by gaps. 

If an $H$-field $K$ has a Liouville $H$-field extension with a gap, 
then $K$ has two Liouville closures that are not isomorphic over $K$,
according to \cite{AvdD2}. In reviewing that paper we realized that it
doesn't contain a proof of that claim, but   
Allen Gehret has since provided us with one; we do not use the claim in this book.   

In the next chapter we introduce a first-order condition ($\upo$-freeness) on an $H$-field that makes it fall under Case~(I) of Theorem~\ref{thm:Liouville closures}; see Corollary~\ref{upoliou}.
% (valid if $K$ is an ungrounded $H$-subfield of $\T$ properly 

A miniature version of ``Liouville closure'' is {\em closure under powers\/},
studied in the setting of Hardy fields in~\cite{Rosenlicht6}, and for $H$-fields in
\cite[Sections~7,~8]{AvdD3}: every 
$H$-field has a closure under powers, and up to isomorphism it has at most two.

%then its value group~$\Gamma$ can be turned into an ordered vector space over $C$ in a natural way, and
%if $K$ is Liouville closed with small derivation, then the pair
%$(\Gamma,\psi)$ is a closed $H$-couple in the sense of \cite{AvdD}; see \cite[Section~7]{AvdD3}.
%In \cite[Section~8]{AvdD3} it is shown that each $H$-field has a minimal $H$-field extension which is closed under powers,
%and an analogue of Theorem~\ref{thm:Liouville closures} holds: every $H$-field $K$ has (up to isomorphism over
%$K$) at least one and at most two such closures under powers.

\section{Miscellaneous Facts about Asymptotic Fields}
\label{Miscellaneous Facts about Asymptotic Fields}

\noindent
This section indicates some possible and impossible features
of asymptotic fields. Only Lemma~\ref{sine cosine} will be
used later. 
Can a differentially closed field 
have a nontrivial valuation making it an asymptotic field? This was 
answered negatively by Scanlon~\cite{Scanlon08}. The solution given 
here, inspired by
Rosenlicht~\cite{Rosenlicht2}, is a little different. Scanlon used the 
logarithmic derivative map of an elliptic curve. The 
differential-algebraic properties of this map are used here 
for another purpose:
to construct asymptotic fields that cannot be obtained 
by coarsening pre-$\d$-valued fields. Unrelated to this, we also indicate an
algebraically closed $\d$-valued field of $H$-type that is not an
algebraic closure of a real closed $\d$-valued field, and we finish with
the result that~$\T$ as a subset of $\T[\imag]$ is definable in
the differential field~$\T[\imag]$.

\subsection*{Differentially closed fields cannot be asymptotic} 
In this subsection $K$ is a differential field, with constant field $C$, and
$y$ ranges over $K$.  

\index{differential field!differentially closed}
\index{closed!differentially}

\begin{prop}\label{Scanlon-Prop} If $K$ is differentially closed, 
then there is no valuation ring $\mathcal{O}\neq K$ of $K$ with 
$\der \mathcal{O} \subseteq  \mathcal{O}$. 
\end{prop}
\begin{proof} Let $K$ be differentially closed
and suppose for a contradiction that $\mathcal{O}\neq K$ is a
valuation ring of $K$ with $\der\mathcal{O}\subseteq \mathcal{O}$. Since $K$ remains
 differentially closed upon replacing $\der$
by $a\der$ with $a\in  \smallo^{\neq}$, we can assume
$\der \mathcal{O}\subseteq \smallo$. Take $y$ such that $y+(y')^2=y^3$ and $y\ne y^3$.
The argument below uses the second part of Lemma~\ref{Cohn-Lemma} several times.
If $y\succ 1$, then $y+(y')^2 \prec y^3$, a contradiction. Thus $y\in \mathcal{O}$, and so
$y\equiv y^3\bmod \smallo$, hence $y\equiv -1$, $0$, or $1 \bmod \smallo$. The case $y\equiv 0 \bmod \smallo$ is impossible, since $y\ne 0$ gives $(y')^2 \prec y$ and $y^3\prec y$. If $y\equiv \-1 \bmod \smallo$, set $y=z+1$, and
we get a similar contradiction from $-2z + (z')^2 =z^3+3z^2$ and $0\ne z\prec 1$,
and the case $y\equiv -1 \bmod \smallo$ is likewise impossible.    
\end{proof}

\subsection*{The logarithmic derivative map on an elliptic curve}
In this subsection $K$ is a differential field with constant field $C$. 
Let $c\in C$, $c\neq 0,1$, and put
$$P(X) := X(X-1)(X-c) = X^3 - (c+1)X^2 +cX \in C[X].$$
Consider the projective plane curve $E$ defined over $C$ by the equation
$$Y^2Z = X^3 - (c+1)X^2Z + cXZ^2.$$
The affine part of $E$ (in standard coordinates) is given by $Y^2=P(X)$, and
$E$ is an elliptic curve whose group law
has $(0:1:0)$, the unique point at infinity on $E$, as its zero element. 
 
We define the logarithmic derivative map $\ell_E \colon E(K)\to K$ as follows:
$$\ell_E(0:1:0):=0, \quad \ell_E(x:y:1) := x'/y\ \text{ if $y\ne 0$,} \quad
\ell_E(x:0:1):= 0.$$ By Lemma~2 on p.~805 
of \cite{Kolchin-GTDF}, 
the map $\ell_E$ is a group homomorphism from $E(K)$ to the additive group 
of $K$; its kernel is $E(C)$. 

\begin{lemma}\label{ell_E maps into O}
Let $\mathcal{O}$ be a valuation ring of $K$ such that $\der \smallo \subseteq \smallo$ and $c,c-1\in \mathcal{O}\setminus  \smallo$. 
Then $\ell_E\big(E(K)\big)\subseteq \mathcal{O}$. 
\end{lemma}
\begin{proof}
Let $x,y\in K$ with $(x:y:1)\in E(K)$. If $y=0$, then 
$\ell_E(x:y:1)=0\in \mathcal{O}$. Assume $y\neq 0$; 
we need to show $x'/y\preceq 1$. We distinguish two cases. 
Suppose first that $x\preceq 1$; then $y^2= P(x) \preceq 1$, hence 
$y\preceq 1$. Differentiating both sides of the equality $y^2=P(x)$ we 
obtain $2yy'=P'(x)\cdot x'$. Hence if $P'(x)\asymp 1$ then 
$x'/y \asymp 2y'\preceq 1$. Otherwise $P'(x)\prec 1$ and therefore 
$y^2=P(x)\asymp 1$, since the reduced polynomial 
$\bar{P}(X)\in (\res K)[X]$ has no multiple zeros; hence 
$x'/y \asymp x'\preceq 1$. 
Now suppose that $x\succ 1$. Then $y^2=P(x)\asymp x^3$ and $(x')^2 \preceq x^3$ 
by Lemma~\ref{Cohn-Lemma}, therefore $x'/y\preceq 1$.
\end{proof}

\begin{cor} If $\ell_E\colon E(K)\to K$ is surjective, then $K$ has no nontrivial 
valuation making it an asymptotic field.
\end{cor}

\noindent
For use in the next subsection we prove:

\begin{prop}\label{making ell_E surjective}
Suppose $C$ is algebraically closed. There exists a differential field 
extension $L$ of $K$, $\d$-algebraic over $K$ with $C_L=C$, 
such that
\begin{enumerate}
\item[\textup{(i)}] $\ell_E\colon E(L)\to L$ is surjective \textup{(}so $L$ has 
no non\-trivial valuation making it an asymptotic 
field\textup{)};
\item[\textup{(ii)}] for each $y\in L^\times$ with $y^\dagger\in C^\times$ there exists $n\ge 1$
with $y^n\in K$.
\end{enumerate}
\end{prop}

\noindent
The main work goes into proving the next two lemmas, in both of which
$C$ is assumed to be algebraically closed.

\begin{lemma}
There is no group homomorphism $E(C)\to C^\times$ with finite kernel.
\end{lemma}
\begin{proof}
Suppose 
$\varphi\colon E(C)\to C^\times$ is a group homomorphism with 
$k:=\abs{\ker\varphi}<\infty$. For an additive abelian group $G$ 
and $n\ge 1$, put
$G[n] := \{g\in G:ng=0\}$, a subgroup of $G$. Then, for $n\ge 1$, $E(C)[n]$ is 
finite with $\abs{E(C)[n]}=n^2$, see \cite[Theorem~VI.6.1]{Silverman}, 
and 
$\varphi\big(E(C)[n]\big)\subseteq \{\zeta\in C^\times: \zeta^n=1\}$. 
Take a prime number $p$ not dividing $k$. Then $\varphi|E(C)[p]$
%\colon E(C)[p\to \{\zeta\in C^\times: \zeta^p=1\}$ 
is injective, contradicting $\abs{\{\zeta\in C^\times: \zeta^p=1\}}\leq p$.
\end{proof}

\begin{lemma}\label{power of y}
Let $a\in K^\times$, and let $f$ be an element in a differential field 
extension of $K$ with the same constant field $C$ as $K$, satisfying 
$(f')^2=a^2P(f)\neq 0$. Let $y$ be a nonzero element of the differential 
field $K\langle f\rangle=K(f,f')$ with $y^\dagger\in C^\times$. Then 
$y^n\in K$ for some $n\ge 1$.
\end{lemma}
\begin{proof}
From Section~\ref{Differential Fields and Differential Polynomials} recall that given a differential field extension $L$ of $K$, we denote by 
$\Aut_\der(L|K)$ the group of differential automorphisms of $L$ which are 
the identity on $K$. 
Suppose for a contradiction that $y^n\notin K$ for all $n>0$. Then $y$ is 
transcendental over $K$, so $f$ is transcendental over $K$. 
%From $y^\dagger \in C^\times$ we obtain that 
Also $\sigma(y)/y\in C^\times$ for every $\sigma\in\Aut_\der\!\big(K(y)|K\big)$,
and the map 
\begin{equation}\label{isomorphism with multiplicative group}
\sigma\mapsto \sigma(y)/y\colon\Aut_\der\!\big(K(y)|K\big)\to C^\times
\end{equation}
is a group isomorphism \cite[p.~803]{Kolchin-GTDF}.
Writing the group operation on $E(K)$ additively, we have for
$\sigma\in\Aut_\der\!\big(K\langle f\rangle|K\big)$,
$$p_\sigma := \big(\sigma(f):\sigma(f'/a):1\big) - (f:f'/a:1)\in E(C),$$
and the map 
\begin{equation}\label{embedding into E(C)}
\sigma\mapsto p_\sigma\colon \Aut_\der\!\big(K\langle f\rangle|K\big)\to E(C)
\end{equation} 
is an injective group homomorphism \cite[p.~807]{Kolchin-GTDF}. The 
homomorphism 
\begin{equation}\label{restriction homomorphism}
\sigma\mapsto{\sigma\restrict K(y)}\colon
\Aut_\der\!\big(K\langle f\rangle|K\big)\to\Aut_\der\!\big(K(y)|K\big)
\end{equation}
is surjective by \cite[Theorem~3 on p.~797]{Kolchin-GTDF}; in particular, 
$\Aut_\der\!\big(K\langle f\rangle|K\big)$ is infinite. Hence 
the map \eqref{embedding into E(C)} is an 
isomorphism, by \cite[p.~807]{Kolchin-GTDF}. Moreover, the kernel of~\eqref{restriction homomorphism} is finite by \cite[p.~796, Theorem~2]{Kolchin-GTDF}. Composing the inverse of the map~\eqref{embedding into E(C)} 
with~\eqref{restriction homomorphism} and~\eqref{isomorphism with 
multiplicative group} yields a homomorphism $E(C)\to C^\times$ with finite 
kernel; this is impossible by the previous lemma.
\end{proof}

\noindent
To prove Proposition~\ref{making ell_E surjective}, fix a differential 
closure $K^{\operatorname{dc}}$ of $K$; then
 $K^{\operatorname{dc}}$  is $\d$-algebraic over $K$ and  $C_{K^{\operatorname{dc}}}=C$.
(See Section~\ref{sec:DCF}.) By Zorn we can take a differential subfield
$L$ of $K^{\operatorname{dc}}$ 
containing $K$ which is maximal with respect to the
property that for each $y\in L^\times$ 
with $y^\dagger\in C^\times$ there exists $n\ge 1$ with $y^n\in K$.
Then $\ell_E\big(E(L)\big)=L$. To see this, let $a\in L^\times$ and 
$f\in (K^{\operatorname{dc}})^\times$ be such that $(f')^2=a^2P(f)\neq 0$. Suppose that $y$ is a nonzero
element of $L\langle f\rangle$ such that $y^\dagger\in C^\times$. 
Lemma~\ref{power of y} (applied to $L$ in place of $K$) yields an 
$n\ge 1$ with $y^n\in L$; since $(y^n)^\dagger=ny^\dagger\in C^\times$, we
have $y^{nm}\in K$ for some $m\ge 1$. So $L\langle f\rangle=L$ by
maximality of $L$; in particular 
$f\in L$ and thus $(f:f'/a:1)\in E(L)$ with
$\ell_E(f:f'/a:1)=a$. \qed

\subsection*{An asymptotic field that is not a coarsening of a 
pre-$\d$-valued field}

Let~$\k$ be a differential field
such that $a'+na\ne 0$ for all $a\in \k^\times$ and $n\geq 1$. Equip~$\k\(( t\)) $ with the valuation $v$ that has $\k[[t]]$ as its valuation ring,
and with the unique derivation extending that 
of $\k$ such that $t'=t$ and $\k[[t]]'\subseteq \k[[t]]$. 
Then $\k\(( t\)) $ is an asymptotic field.  (To see this, 
note that 
$(at^n)'=(a'+na)t^n$ for $a\in \k^\times$ and $n\geq 1$.)
It has small derivation, with $\Psi=\{0\}$ (so its asymptotic
couple is of $H$-type); its differential residue field is
$\k$, and its constant field is $C_\k$. The valued differential subfield 
$K:=\k(t)$ of~$\k\(( t\)) $ is an asymptotic field with valuation ring 
$\mathcal{O}=\k[t]_{(t)}$. If $\k$ has an element $b$ with $b'\ne 0$,
then $b\asymp b'\asymp 1$, and so $K$ is not pre-$\d$-valued.

\begin{lemma}\label{asd}
Suppose $\k$ has an element $b$ with $b'\neq 0$, and $\k$ has no 
nontrivial valuation making it a pre-$\d$-valued field. 
Then $\k(t)$ has no valuation $v_1$ with~${v_1(t)>0}$ that makes $\k(t)$ a 
pre-$\d$-valued field. 
%such that $t\in\mathcal O_1$, $t^{-1}\notin\mathcal O_1$. 
\end{lemma}
\begin{proof}
Suppose $v_1$ is a valuation on $\k(t)$ with $v_1(t)>0$ that makes $\k(t)$ a 
pre-$\d$-valued field. Let $\mathcal{O}_1$ be the valuation ring of $v_1$.
Then $\mathcal{O}_1\cap\k$ is a valuation ring of~$\k$ making $\k$ a 
pre-$\d$-valued field, so $\mathcal{O}_1\cap \k=\k$. 
Then $\mathcal{O}_1=\k[t]_{(t)}=\mathcal{O}$, which 
contradicts the observation preceding the lemma.
\end{proof}

\noindent
If $\k$ satisfies the hypothesis of the lemma, then
$K$ cannot be obtained by coarsening a pre-differential valuation of $\k(t)$:
there is no valuation $v_1\colon \k(t)^\times\to\Gamma_1$ 
making~$\k(t)$ a pre-$\d$-valued field such that for some
convex subgroup $\Delta$ of $\Gamma_1$ the coarsened valuation
$\dot{v}_1\colon \k(t)^\times \to \Gamma_1/\Delta$ has valuation ring $\mathcal{O}$.

\medskip\noindent
To obtain a differential field $\k$ satisfying the 
requirements above we first take any algebraically closed field $C$ 
of characteristic zero. Let 
$C(x)$ be a field extension of $C$ with~$x$ transcendental over $C$,
equipped with the unique derivation with constant field $C$ and $x'=1$.
Proposition~\ref{making ell_E surjective} yields a differential 
field extension $\k$ of $C(x)$ such that
for every $y\in\k^\times$ with $y^\dagger\in C^\times$ there exists $n>0$ 
with $y^n\in C(x)$, and such that $\k$ has no nontrivial 
valuation making it an asymptotic field. 
Since there is no $y\in C(x)^\times$ with $y^\dagger\in C^\times$ 
(by Corollary~\ref{cor:no new exps under integration}), 
there is no $y\in\k^\times$ with $y^\dagger\in C^\times$. 
Hence, as required, $a'+na\ne 0$ for all $a\in \k^\times$ and $n\geq 1$, 
$C_\k\neq \k$, and~$\k$ has no nontrivial valuation making it 
a pre-$\d$-valued field. Thus the asymptotic field $K$ defined just 
before Lemma~\ref{asd} cannot be obtained by coarsening a 
pre-differential valuation.

\subsection*{An algebraically closed $\d$-valued field 
that is not an algebraic closure of a real closed $\d$-valued field}
Every algebraically closed field $L$ of characteristic zero has a real 
closed subfield $K$ such that $L=K(\imag)$, where $\imag^2=-1$. This suggests the 
following question: 

\medskip\noindent
\textit{Let $L$ be an algebraically closed 
$\d$-valued field of $H$-type; 
is there a real closed $\d$-valued subfield $K$ of $L$ such that 
$L=K(\imag)$?}\/  

\medskip\noindent
We show that the answer is 
negative for $L$ as in Corollary~\ref{newcor}.
This result follows from the next lemma, also used in the next
subsection.
Let $K$ be a real closed $\d$-valued field and $K[\imag]$ its algebraic closure, $\imag^2=-1$. 
The residue field of $K$ is isomorphic to~$C$, hence real closed. Thus the valuation ring $\mathcal O$ of $K$ is convex by Theorem~\ref{thm:KW}. 
For $y=a+b\imag \in K[\imag]$ ($a,b\in K$) we let $\abs{y}:=\sqrt{a^2+b^2}\in K^{\geq}$ be the absolute value of $y$. %For all $x,y\in K[\imag]$,
%$$\abs{x}=0 \Longleftrightarrow x=0, \qquad \abs{xy}=\abs{x}\abs{y}, \qquad \abs{x+y}\leq \abs{x}+\abs{y}.$$
We have the subgroup
$$S\ :=\ \big\{ y\in K[\imag]:\ \abs{y}=1\big\}$$
of the multiplicative group $K[\imag]^\times$, with
$S\subseteq \mathcal O+\mathcal O\imag$, $K[\imag]^\times = K^>\cdot S$, $K^>\cap S=\{1\}$.

\begin{lemma}\label{sine cosine} Let
$a+b\imag\in S$ \textup{(}$a,b\in K$\textup{)}. Then
$(a+b\imag)^\dagger = \wr(a,b)\imag$.
Thus
\begin{align*} \big(K[\imag]^\times\big){}^\dagger\ &=\ (K^>)^\dagger \oplus S^\dagger,\ \text{ an internal direct sum of subgroups of $\big(K[\imag]^\times\big){}^\dagger$,}\\
S^\dagger\ &\subseteq\ \{f\in K:\ f\preceq g' \text{ for some $g\in \smallo$}\}\cdot \imag. 
\end{align*}
\end{lemma}
\begin{proof} Since $(a+b\imag)(a-b\imag)=1$ we have
\begin{align*}
(a+b\imag)^\dagger\	&=\ (a'+b'\imag)(a-b\imag)\ =\ 
			 (aa'+bb')+(ab'-a'b)\imag \\ 
			&=\ \textstyle\frac{1}{2}\big( a^2+b^2 \big)' + (ab'-a'b)\imag\ =\ (ab'-a'b)\imag\ =\ \wr(a,b)\imag.
\end{align*}					
From $a,b\in\mathcal O=C+\smallo$ we get $\wr(a,b)\preceq g'$ for some $g\in \smallo$. It remains to note that $K[\imag]^\times = K^>\cdot S$ gives $\big(K[\imag]^\times\big){}^\dagger = (K^>)^\dagger + S^\dagger$.
\end{proof}

%\begin{lemma}\label{sine cosine}
%Suppose $K$ is a real closed $\d$-valued field with an element $a\in K^\times$ satisfying $a \succ \der \smallo$.  Then the algebraic closure $K(\imag)$ of $K$, where $\imag^2=-1$, has no nonzero element $y$ such that $y^\dagger=a\imag$.
%\end{lemma}
%\begin{proof} The residue field of $K$ is isomorphic to $C$, hence real closed, so $K^{\preceq 1}$ is convex in $K$  by Theorem~\ref{thm:KW}. Suppose $y\in K(\imag)^\times$ satisfies $y^\dagger=a\imag$.  We have $y=f+g\imag$ with $f,g\in K$; think of  $y$ as $\ex^{\imag\int{a}}$ and of $f$ and $g$ as $\cos \int{a}$ and $\sin \int{a}$.  So  $f'=-ag$ and $g'=af$, hence $f,g\neq 0$ and
%$$(f^2+g^2)'=2ff'+2gg'=0.$$
%Therefore $f^2+g^2\in C^\times$, so $f\preceq 1$ and $g\preceq 1$, and either $f\succeq 1$ or $g\succeq 1$. Suppose $f\succeq 1$. (The case  $g\succeq 1$ is similar.)  Since $g'=af\succeq a\succ b'$ for  all $b\in K^{\preceq 1}$, we have $g\succ 1$, a contradiction.
%\end{proof}

\begin{cor}\label{corsinecosine}
If $L$ is a $\d$-valued field with small derivation
and an element $y\in L^\times$ such that  $y^\dagger=\imag$,  $\imag^2=-1$,
then $L$ has no real closed 
$\d$-valued subfield $K$ with $L=K(\imag)$.
\end{cor}

\subsection*{Definability of $\T$ in $\T[\imag]$} It is well-known that
$\R$ as a subset of $\C=\R[\imag]$ is not definable (even allowing parameters) in the field $\C$ of complex numbers; see~\ref{cor:ACF sm}. In stark contrast, the {\em differential\/} field structure of $\T[\imag]$ is rich enough to define~$\T$:

\begin{prop}\label{prop:defining T} The subset $\T$ of $\T[\imag]$ is definable in the
differential field $\T[\imag]$.
\end{prop} 
\begin{proof} Below, ``definable'' means definable without parameters in the differential field $\T[\imag]$. Let $\mathcal{O}$ be the valuation ring of $\T$
with maximal ideal $\smallo$. Then $\mathcal{O}+\mathcal{O}\imag$
is the valuation ring of $\T[\imag]$ with maximal ideal $\smallo + \smallo \imag$. We begin with noting the definability of $\R$, as a consequence of Lemma~\ref{sine cosine}:
$$\R\ =\ \big\{y\in \T[\imag]:\ \text{$y'=0$ and $y=f^\dagger$  for some 
$f\in \T[\imag]^\times$}\big\}.$$
Indeed, this lemma gives the more
precise result $\big(\T[\imag]^\times\big){}^\dagger \subseteq \T+\imag \der\smallo$. We claim that $\big(\T[\imag]^\times\big){}^\dagger = \T+\imag \der\smallo$. To see this, note that for $b\in \smallo$ we have
$$\sin b\ :=\  b-\frac{b^3}{6} + \frac{b^5}{120} -\dots\in \T, \qquad \cos b\ :=\ 1-\frac{b^2}{2} + \frac{b^4}{24} -\dots\in \T^{\ne}$$
with $(\sin b)'=b'\cos b$ and $(\cos b)'=-b'\sin b$, and so
$$f\ :=\ \exp(\imag b)\ :=\ \cos b + \imag \sin b\in \T[\imag]^\times$$ satisfies
$f^\dagger= ib'$. Hence $\imag \der \smallo\subseteq (\T[\imag]^\times)^\dagger$, and thus $\T+ \imag \der \smallo\subseteq (\T[\imag]^\times)^\dagger$, and the claim follows. Think of $\T+\imag \der \smallo$ as a thin strip around the ``real axis'' $\T$ in the ``complex
plane''~$\T[\imag]$. Intersecting $\T+\imag \der \smallo$ with
its multiple $\imag(\T+\imag \der\smallo)$ yields the definable
set~${\der\smallo + \imag\der\smallo}$, and taking integrals gives
the definability of
$\mathcal{O} + \mathcal{O}\imag$, and of its maximal ideal $\smallo + \smallo\imag$.
Hence
$\big(\T[\imag]^\times\big){}^\dagger+ (\smallo + \smallo\imag)=\T + \smallo\imag$ is definable. Thus $\mathcal{O} + \smallo\imag=(\mathcal{O}+\mathcal{O}\imag)\cap (\T+ \smallo\imag)$ is definable. It is easy to check that for $f\in \mathcal{O}+\smallo\imag$, 
$$f\cdot (\T+\smallo \imag)\ \subseteq\ \T+\smallo\imag\ \Longleftrightarrow\ f\in \mathcal{O},$$
so $\mathcal{O}$ is definable, and therefore its
fraction field $\T$ in $\T[\imag]$ is definable.
%$\T= \bigcap_{f\in \mathcal{O}^{\ne}} f\cdot (\T+\smallo \imag)$ is definable.
\end{proof}

%\subsection*{Notes and comments} Except for 
%Corollary~\ref{corsinecosine} and 
%Proposition~\ref{prop:defining T},
%this section is mainly based on 
%Rosenlicht~\cite{Rosenlicht2} and Scanlon~\cite{Scanlon08}.

%% file: mt-11.tex
\chapter{Eventual Quantities, Immediate Extensions, and Special Cuts}\label{evq}

\setcounter{theorem}{0}

\noindent
Our main interest is ultimately in $H$-fields with asymptotic integration
and small derivation such as $\mathbb{T}$. The induced derivation on the residue field of such asymptotic fields is trivial, however, so these asymptotic
fields are not covered by Corollary~\ref{cor2asdiftrzda} on
spherically complete immediate asymptotic extensions. 
One goal in the present chapter is to remedy this
defect by establishing the following result in Section~\ref{sec:construct imm exts}:
 
\begin{theorem}\label{immHasint} Every asymptotically maximal 
$H$-asymptotic field with 
rational asymptotic integration is spherically complete. 
\end{theorem}

\noindent
Proving Theorem~\ref{immHasint}
requires some tools that are also important later. These tools arise from the
fact that for asymptotic fields $K$ with
asymptotic integration and $P\in K\{Y\}^{\ne}$, certain quantities
associated to its compositional conjugates $P^{\phi}$, such as its dominant degree $\ddeg P^{\phi}$, become constant for sufficiently high 
$v\phi\in \Psi^{\downarrow}$. These ``eventual quantities'' will be 
studied in Sections~\ref{evtbeh}, \ref{sec:ndeg and nval}, \ref{sec:nwt}.

In Sections~\ref{sec:special cuts}, \ref{sec:uplfree}, \ref{sec:behupo} we consider special (definable) cuts in 
$H$-asymptotic fields $K$ with
asymptotic integration, and introduce some key elementary properties 
of~$K$, namely $\upl$-freeness and $\upo$-freeness, which express
that these cuts are not realized in~$K$. We show that $\mathbb{T}$ has these properties, but a full exploitation
of these subtle but powerful elementary properties must be left to Chapter~\ref{ch:The Dominant Part and the Newton Polynomial}, where the machinery of
triangular automorphisms of $K\{Y\}$ from Chapter~\ref{ch:triangular automorphisms} is available.

In Section~\ref{sec:special sets} we consider certain special existentially 
definable subsets of Liouville closed $H$-fields $K$, 
and the behavior of the functions $\omega$ and
$\sigma$ (introduced in Section~\ref{sec:secondorder}) on these sets. This will play a
role in our Elimination of Quantifiers for $\mathbb{T}$ in Chapter~\ref{ch:QE}. 

%Terminology: from now on, an 
%{\em asymptotic extension}\/ of an asymptotic field 
%$K$ is a valued differential field extension of $K$ 
%that is asymptotic.

\section{Eventual Behavior}\label{evtbeh}

\noindent
{\em In this section $K$ is an asymptotic field with value group $\Gamma\ne \{0\}$.}\/
We let $\phi$ range over~$K^\times$, and $\bsigma$, $\btau$ over $\N^*$. We also 
fix a differential polynomial $P\in K\{Y\}^{\ne}$. 

The function $y\mapsto P(y)$ defined
by $P$ on $K$ does not give up
its secrets easily, but we do have some tricks up our sleeve. First, things are
more transparent if $\der\smallo\subseteq \smallo$, and the ``best case''
is when $\sup \Psi = 0$. If $\sup \Psi$ exists we can
reduce to that best case by compositional conjugation. In general 
we try to simulate this best case by working in compositional conjugates
$K^{\phi}$ with small derivation, that is, $v\phi < (\Gamma^{>})'$, but 
$v\phi$ as high as possible. 
When $K$ is ungrounded, 
it turns out that certain quantities
associated to~$P^{\phi}$ such as $\dwt(P^{\phi})$ eventually
stabilize for increasing $v\phi\in \Psi^{\downarrow}$. 
These ``eventual'' quantities associated to $P$ turn out to be important 
invariants. This section is devoted to proving their existence.

\subsection*{Behavior of $v F^{n}_{k}(\phi)$} 
In order to better understand $v(P^{\phi})$ as a function of $\phi$
we use from Lemma~\ref{conj-lemma} the identity 
\begin{equation}\label{eq:coeffs of Pphi}
(P^\phi)_{[\bsigma]}\  =\ \sum_{\btau \geq \bsigma} F^{\btau}_{\bsigma}(\phi) P_{[\btau]}. 
\end{equation}
This leads to the study of $vF^{\btau}_{\bsigma}(\phi)$. 
Recall that for $\btau=\tau_1\cdots\tau_d\ge \bsigma=\sigma_1\cdots\sigma_d$,
$$F^{\btau}_{\bsigma}\ :=\ F^{\tau_1}_{\sigma_1}\cdots F^{\tau_d}_{\sigma_d}.$$

\begin{lemma}\label{compconjval, general}
If $\der\mathcal{O}\subseteq\mathcal{O}$ and $\phi\preceq 1$, then $v(P^\phi)\geq v(P)$, with equality if $\phi\asymp 1$.
\end{lemma}
\begin{proof} Assume $\der\mathcal{O}\subseteq\mathcal{O}$ and $\phi\preceq 1$.
For $0\leq k\leq n$ we have $F^n_k\in \Q\{X\}$, so $F^n_k(\phi)\preceq 1$. 
Thus $v(P^\phi)\geq v(P)$ by \eqref{eq:coeffs of Pphi}.
If $\phi\asymp 1$, use
$P=(P^\phi)^{\phi^{-1}}$. 
\end{proof}

\noindent 
In studying $v F^{n}_{k}(\phi)$ we consider the case $\phi^\dagger \preceq \phi$
in the next lemma, and take up the case $\phi^\dagger \succ \phi$ 
in the next subsection. We set $\derdelta=\phi^{-1}\der$.

\begin{lemma}\label{compconjval}
Suppose that $\derdelta\smallo\subseteq \smallo$, and let 
$0\leq k\leq n$.
\begin{enumerate}
\item[\textup{(i)}] If $\phi^\dagger\preceq \phi$, then $v\big(F^n_k(\phi)\big)\geq nv\phi$, with equality if $k=n$.
\item[\textup{(ii)}] If $\phi^\dagger\prec \phi$ and $k<n$, then 
$v\big(F^n_k(\phi)\big)>nv\phi$.
\end{enumerate}
\end{lemma}

\begin{proof}
Note that $\phi^\dagger\preceq \phi$ means 
$v\big(\derdelta(\phi)\big)\geq v\phi$. 
Now use \eqref{recursion for F} and induction on~$n$, 
and Lemma~\ref{closed under der}
applied to $K^\phi$ in place of $K$.
\end{proof}

\begin{cor}\label{compconjval, cor 2}
Suppose that $\derdelta\smallo\subseteq 
\smallo$ and $\phi^\dagger\preceq \phi$, and  $\btau\geq\bsigma$. Then
$v\big(F^\btau_{\bsigma}(\phi)\big)\geq
\|\btau\| v\phi$, with equality if $\btau=\bsigma$.
\end{cor}

\subsection*{Behavior of $v(P^{\phi})$}
Note that for $\derdelta:= \phi^{-1}\der$ we have 
$$ \derdelta\smallo\subseteq \smallo\ \Longleftrightarrow\ 
v\phi < (\Gamma^{>})'.$$
{\em Accordingly we restrict $\phi$ in the rest of
this subsection to satisfy $v\phi < (\Gamma^{>})' $}. 
(Thus $\phi^\dagger \prec 1$ if
$\phi \prec 1$.) The main goal of this subsection is:

\begin{prop}\label{compconjval, prop} Suppose that $\der\smallo \subseteq \smallo$. Then for all $\phi\preceq 1$,
$$v(P) + \Pnu(P^\phi)v\phi\ \leq\ v(P^\phi)\ \leq\ v(P)+\Pmu(P)v\phi.$$
\end{prop}

\noindent
Here $\phi\preceq 1$ is in addition to the standing assumption
that $v\phi < (\Gamma^{>})' $.
We first prove some
lemmas on the behavior of $v F^{n}_{k}(\phi)$ for $\phi^\dagger \succ \phi$. 
Combining this
with the case $\phi^\dagger \preceq \phi$ from the previous subsection,
we shall derive Proposition~\ref{compconjval, prop}.

\medskip\noindent
Assume $P$ has order at most $r$, let 
$\bsigma, \btau\in \{0,\dots,r\}^*$ have equal 
length, and set 
$${\btau\brack \bsigma} :={ \tau_1\brack \sigma_1}\cdots {\tau_d\brack \sigma_d}\ \quad \text{for $\bsigma=\sigma_1\cdots\sigma_d$, $\btau=\tau_1\cdots\tau_d$,}$$
where ${n\brack m}$ is the signless Stirling number of the first kind
from Section~\ref{Compositional Conjugation}. 

\nomenclature[Cc]{${\btau\brack \bsigma}$}{${ \tau_1\brack \sigma_1}\cdots {\tau_d\brack \sigma_d}$,
for $\bsigma=\sigma_1\cdots\sigma_d,\btau=\tau_1\cdots\tau_d\in\N^*$}

\begin{lemma}\label{Riccatipow} Assume $\der\smallo\subseteq \smallo$. Let
$z\in K$, $z\succ 1$. Then 
$$R_n(z)=z^n(1+ \epsilon)\ \text{  with $\epsilon \preceq z^\dagger/z \prec 1$.}$$
\end{lemma}
\begin{proof} This is clear for $n=0$ and $n=1$. Suppose $n\ge 1$ and
$R_n(z)=z^n(1+ \epsilon)$ with $\epsilon$ as in the lemma. Then
\begin{align*} R_{n+1}(z)&=zR_n(z) + R_n(z)'=
z^{n+1}(1+ \epsilon)+nz^{n-1}z'(1+\epsilon) + z^n\epsilon'\\
  &=z^{n+1}\left(1 + \epsilon + n\frac{z^\dagger}{z}(1+\epsilon) +         
                     \frac{\epsilon'}{z}\right).
\end{align*}
Next, use that $\epsilon'\prec z^\dagger$, hence $\epsilon'/z \prec
z^\dagger/z$. 
\end{proof}

\begin{lemma}\label{derivatives of equal weight and degree}
Assume $\der\smallo\subseteq \smallo$. Let $f\in K^\times$ be such that 
$f^\dagger\succ 1$. Then 
$$f^{[\bsigma]}\sim f^{|\bsigma|}(f^\dagger)^{\|\bsigma\|}.$$
Hence if $Q\in \mathcal{O}\{Y\}$ is homogeneous of degree $d$, isobaric of 
weight $w$, and $c\asymp 1$ where $c:=\sum_{\i} Q_{\i}$ is the sum of its 
coefficients, then
$$Q(f)\sim c f^{d}(f^\dagger)^w.$$
\end{lemma}
\begin{proof}
It suffices to do the case $|\bsigma|=1$ of a single factor, that is,
it suffices to show  $f^{(n)}/f\sim
(f^\dagger)^n$. Now $R_n(f^\dagger)=f^{(n)}/f$, so Lemma~\ref{Riccatipow}
applied to $z=f^\dagger$ yields the desired result.  
\end{proof}

\begin{remark}
With the assumptions of this lemma, the proof gives
$$f^{[\bsigma]}=f^{|\bsigma|}(f^\dagger)^{\|\bsigma\|}(1+\varepsilon)\qquad\text{where $\varepsilon\preceq f^{\dagger\dagger}/f^\dagger\prec 1$.}$$
\end{remark}

\begin{lemma}\label{derivatives of equal weight and degree, 2}
Suppose $\phi^\dagger\succ \phi$. Then for $0 < k\leq n$ we have
$$F^n_k(\phi) \sim {n\brack k} \phi^{k}(\phi^\dagger)^{n-k}.$$
\end{lemma}
\begin{proof} Let $0 < k\leq n$.
Recall from Section~\ref{Compositional Conjugation} the definition of the differential polynomial $G^n_k\in K^\phi\{Y\}$ with nonnegative integer coefficients, homogeneous of degree~$n$ and isobaric of weight~$n-k$, satisfying $F^n_k(\phi)=G^n_k(\phi)$. By the previous lemma applied to $Q=G^n_k$ and with $K^\phi$ in place of $K$ we obtain 
$$F^n_k(\phi)=G^n_k(\phi)\sim c \phi^n\left(\frac{\derdelta(\phi)}{\phi}\right)^{n-k}=c\phi^{k}(\phi^\dagger)^{n-k},$$
where $c$ is the sum of the coefficients of $G^n_k$, that is, $c={n\brack k}$
by Lemma~\ref{lem:Stirling}.
\end{proof}

\nomenclature[Cc]{$\supp \bsigma$}{set of $i\in\{1,\dots,d\}$ with $\sigma_i\neq 0$, for $\bsigma=\sigma_1\cdots\sigma_d\in\N^*$}

\noindent
We set $\supp\bsigma:=\{i:\sigma_i\neq 0\}$. 
So if $\btau\geq\bsigma$, then $\supp\bsigma\subseteq\supp\btau$. If $\btau\geq\bsigma$ and $\phi^\dagger\succ \phi$, then we have the equivalences 
$$F^\btau_{\bsigma}(\phi)\neq 0 \quad\Longleftrightarrow\quad 
{\btau \brack \bsigma} \neq 0 \quad\Longleftrightarrow\quad \supp\btau=\supp\bsigma.$$

\begin{lemma}\label{compconjval, 3}
Suppose $\btau\geq\bsigma$ and $\supp\btau=\supp\bsigma$, and $\phi^\dagger\succ \phi$. Then 
$$F^\btau_{\bsigma}(\phi) \sim {\btau \brack\bsigma} \phi^{\|\bsigma\|}
(\phi^\dagger)^{\|\btau\|-\|\bsigma\|}.$$
Also, if $\phi\not\asymp 1$, then 
$v\big(F^\btau_\bsigma(\phi)\big)=
\|\btau\|v\phi+o(v\phi)$. 
\end{lemma}
\begin{proof}
The first statement follows from the previous lemma. Note that
$$\phi^{\|\bsigma\|}(\phi^\dagger)^{\|\btau\|-\|\bsigma\|}\ =\ 
\phi^{\|\btau\|} \big(\derdelta(\phi)/\phi\big)^{\|\btau\|-\|\bsigma\|}.$$
If $\phi\nasymp 1$, then $\psi^\phi(v\phi)=
v\big(\phi^\dagger/\phi\big)<0$, so
$v\big(\derdelta(\phi)/\phi\big)=o(v\phi)$ by Lemma~\ref{PresInf-Lemma2}(iv).
\end{proof}

\begin{lemma}\label{compconjval, 4}
Suppose $\phi \prec 1$ and $\btau\geq\bsigma$.
Then $v\big(F^\btau_{\bsigma}(\phi)\big)\geq\|\bsigma\|v\phi$. Moreover, 
$$v\big(F^\btau_{\bsigma}(\phi)\big)\ =\| \bsigma\|v\phi\ \Longleftrightarrow\ \btau=\bsigma.$$
\end{lemma}
\begin{proof} From  $\phi \prec 1$ we obtain $\phi^\dagger \prec 1$.
If $\phi^\dagger\succ\phi$, then we appeal to 
Lemma~\ref{compconjval, 3}. If $\phi^\dagger\preceq\phi$, then we use Corollary~\ref{compconjval, cor 2}.  
\end{proof}

\begin{lemma}\label{compconjval, 5} Assume $\phi \prec 1$.
Then $v\big((P^\phi)_{[\bsigma]}\big)\geq\|\bsigma\|v\phi+v(P)$, and
$$\begin{cases}
(P^\phi)_{[\bsigma]}\sim\phi^{\|\bsigma\|}P_{[\bsigma]} &\text{if $v(P_{[\bsigma]})=v(P)$,} \\ 
v\big((P^\phi)_{[\bsigma]}\big)>\|\bsigma\|v\phi+v(P) &\text{otherwise.}
\end{cases} $$
\end{lemma}
\begin{proof}
Suppose $\btau\geq\bsigma$. 
Then $v\big(F^\btau_\bsigma(\phi) P_{[\btau]}\big) \geq\|\bsigma\|v\phi+v(P)$
by Lemma~\ref{compconjval, 4}, hence
$$v\big((P^\phi)_{[\bsigma]}\big)\ \geq\ \|\bsigma\|v\phi+v(P)$$ 
by \eqref{eq:coeffs of Pphi},
showing the first statement. 
Lemma~\ref{compconjval, 4} yields that 
if $v(P_{[\bsigma]})=v(P)$ and $\btau\neq\bsigma$, then 
$$v\big(F^\btau_\bsigma(\phi) P_{[\btau]}\big) > 
\|\bsigma\|v\phi+v(P)=v(\phi^{\|\bsigma\|}P_{[\bsigma]}).$$ 
Hence
$(P^\phi)_{[\bsigma]}\sim\phi^{\|\bsigma\|}P_{[\bsigma]}$ if $v(P_{[\bsigma]})=v(P)$.
Suppose that $v(P_{[\bsigma]})>v(P)$. Then
$v(\phi^{\|\bsigma\|}P_{[\bsigma]})>\|\bsigma\|v\phi+v(P)$, and
by the previous lemma again
$$v\big(F^\btau_\bsigma(\phi) P_{[\btau]}\big) > 
\|\bsigma\|v\phi+v(P)
\quad\text{if $\btau\neq\bsigma$,}$$ 
hence $v\big((P^\phi)_{[\bsigma]}\big)>\|\bsigma\|v\phi+v(P)$.
\end{proof}

\begin{proof}[Proof of Proposition~\ref{compconjval, prop}]
Assume $\der\smallo\subseteq \smallo$ and $\phi\preceq 1$. We need to show 
$$v(P) + \Pnu(P^\phi)v\phi\ \leq\  v(P^\phi)\ \leq\ v(P)+\Pmu(P)v\phi.$$
If $\phi \asymp 1$, this holds by Lemma~\ref{compconjval, general}. Let $\phi\prec 1$ and take 
$\bsigma$ with
$\|\bsigma\|=\Pnu(P^\phi)$. Then  by Lemma~\ref{compconjval, 5} we have
$$v(P)+\Pnu(P^\phi)v\phi\ =\ \|\bsigma\|v\phi+v(P)\ \leq\ v\big((P^\phi)_{[\bsigma]}\big)\ =\ v(P^\phi).$$  For
$\bsigma$ such that $\|\bsigma\|=\Pmu(P)$, this same lemma gives
$$v(P^\phi)\ \leq\  v\big((P^\phi)_{[\bsigma]}\big)\ =\ v(P)+\Pmu(P)v\phi,$$
as required.
\end{proof}

\medskip\noindent
We record a few consequences of Proposition~\ref{compconjval, prop} and Lemma~\ref{compconjval, 5}.

\begin{cor} \label{compconjval, cor 3}
 Assume $\der\smallo\subseteq \smallo$ and set $w:=\Pmu(P)$. Then 
\begin{enumerate}
\item[\textup{(i)}] If $\phi\preceq 1$ and $\Pnu(P^\phi)=\Pmu(P)$, then 
$v(P^\phi)=v(P)+w v\phi$.
\item[\textup{(ii)}] If $\phi\prec 1$, then $\Pmu(P^\phi)\le \Pnu(P^\phi) \le \Pmu(P)\le \Pnu(P)$.
\item[\textup{(iii)}] If $\phi\prec 1$ and $\Pmu(P^\phi)=\Pnu(P)$, then for all $\bsigma$,
\begin{align*} 
v(P_{[\bsigma]})=v(P) \qquad &\Longleftrightarrow\qquad v\big( (P^\phi)_{[\bsigma]} \big)=v(P^\phi),\\
 v(P_{[\bsigma]})=v(P) \qquad &\, \Longrightarrow\qquad 
(P^\phi)_{[\bsigma]}\sim\phi^{w}P_{[\bsigma]}.
\end{align*}
\end{enumerate}
\end{cor}
\begin{proof}
Item (i) follows from Proposition~\ref{compconjval, prop}.
Suppose $\phi\prec 1$. Then by the same Proposition we have $\Pnu(P^\phi) \le \Pmu(P)$, which gives (ii). Assume also that 
$\Pmu(P^\phi)=\Pnu(P)$. Then the four numbers in (ii) are all equal to $w$.
Therefore, if $v(P_{[\bsigma]})=v(P)$, then $\|\bsigma\|=\Pmu(P)=\Pnu(P)=w$, so 
$(P^\phi)_{[\bsigma]}\sim\phi^w P_{[\bsigma]}$ by Lemma~\ref{compconjval, 5}, and 
thus $v\big( (P^\phi)_{[\bsigma]} \big)=v(P^\phi)$ by (i).  
Conversely, if $v\big( (P^\phi)_{[\bsigma]} \big)=v(P^\phi)$, then 
$\|\bsigma\|=\Pmu(P^\phi)=\Pnu(P^\phi)=w$, so $v(P^\phi)=v(P)+\|\bsigma\|v\phi$ by (i), and thus $v(P_{[\bsigma]})=v(P)$ by Lemma~\ref{compconjval, 5}.
\end{proof}

\noindent
Now let $\phi_1, \phi_2\in K^\times$ and $v\phi_1, v\phi_2< (\Gamma^{>})'$, 
$\phi_1\preceq \phi_2$. Then $\phi_3:= \phi_1\phi_2^{-1}$ satisfies
$K^{\phi_1}=(K^{\phi_2})^{\phi_3}$ and  $P^{\phi_1}=(P^{\phi_2})^{\phi_3}$, in particular,
the derivation of $(K^{\phi_2})^{\phi_3}$ is small. 
Thus if $\phi_1\prec\phi_2$, then $\phi_3\prec 1$, so we can apply 
Corollary~\ref{compconjval, cor 3} with $K^{\phi_2}$ and $P^{\phi_2}$  instead of 
$K$ and $P$, with $\phi_3$ in the role of $\phi$:

\begin{cor}\label{phi_1 phi_2}
Suppose $\phi_1, \phi_2\in K^\times$ and $v\phi_1, v\phi_2< (\Gamma^{>})'$, $\phi_1\prec\phi_2$. Then
$$\Pmu(P^{\phi_1})\ \leq\ \Pnu(P^{\phi_1})\ \leq\ \Pmu(P^{\phi_2})\ \leq\ \Pnu(P^{\phi_2}),$$ and
if 
$\Pnu(P^{\phi_1})=\Pmu(P^{\phi_2})=w$,
then $$v(P^{\phi_1})-wv(\phi_1)\ =\
v(P^{\phi_2}) - wv(\phi_2).$$
\end{cor}

\begin{cor}\label{cor:Pphi flat equivalence}
Assume $K$ is of $H$-type and $\der\smallo\subseteq\smallo$. Let $\Phi\in K^\times$ be such that $\Phi\nasymp^\flat 1$. Then for all
$\phi\preceq 1$ we have $P^\phi \asymp_\Phi P$. 
\end{cor}
\begin{proof} This is clear for $\phi\asymp 1$. Assume
$\phi\prec 1$. Then $\phi^\dagger\prec 1 \preceq\Phi^\dagger$, so $\phi\flatter\Phi$, and hence the claim follows from Proposition~\ref{compconjval, prop}.
\end{proof}

\subsection*{Newton weight, Newton degree, and Newton multiplicity} We recall that  
$\Psi^{\downarrow}< (\Gamma^{>})'$. Moreover, 
if $K$ has no gap, then 
$$\Psi^{\downarrow}\ =\ \big\{\gamma\in \Gamma:\ \gamma <  (\Gamma^{>})'\big\}.$$
Call $\phi$
{\bf active} (in $K$) if $v\phi \in \Psi^{\downarrow}$. 
If $\phi$ is active in $K$, then $\phi$ remains active in every asymptotic field extension of $K$. 

\index{element!active}
\index{active}

\medskip\noindent
{\em In the rest of this subsection we assume that $K$ is ungrounded,}\/ and
we restrict~$\phi$ and~$\phi_0$ to be active in $K$, in particular, $\phi, \phi_0\in K^\times$. This restriction
implies the restriction on~$\phi$ in the previous subsection, and coincides with it if $K$ has no gap.

\begin{lemma}\label{evpmupnu} There exists $\phi_0$ such that for all $\phi\preceq \phi_0$,
$$ \Pmu(P^{\phi})\ =\ \Pnu(P^{\phi})\ =\ \Pmu(P^{\phi_0})\ =\ \Pnu(P^{\phi_0}).$$
\end{lemma}
\begin{proof} Clear from the first part of Corollary~\ref{phi_1 phi_2}. \end{proof}

\noindent
We say that a property $S(\phi)$ of (active) elements $\phi$ holds 
{\bf eventually} if there exists $\phi_0$ such that $S(\phi)$ holds for
all~$\phi\preceq \phi_0$. Let $\phi_0$ be as in Lemma~\ref{evpmupnu} and note that $\Pnu(P^{\phi_0})$ does not depend on
the choice of such $\phi_0$. We set 
$$\nwt(P)\ :=\ \Pnu(P^{\phi_0}) = \text{eventual value of $\dwt(P^{\phi})$}  = \text{eventual value of $\dwv(P^{\phi})$,}$$ and call it the {\bf Newton weight\/}
of $P$. Thus  
$v(P^{\phi})-\nwt(P)v(\phi)$ is independent of $\phi\preceq \phi_0$, by the second part of Corollary~\ref{phi_1 phi_2}, and \nomenclature[X]{$v^{\ev}(P)$}{eventual value of $v(P^{\phi})-\nwt(P)v(\phi)$} we set $$v^{\ev}(P)\ :=\  v(P^{\phi})-\nwt(P)v(\phi) \qquad(\phi\preceq \phi_0).$$

\index{eventually}
\index{differential polynomial!Newton weight}
\index{Newton!weight}
\index{weight!Newton}
\nomenclature[X]{$\nwt(P)$}{Newton weight of $P$}

\begin{cor}\label{eventualequality} $v(P^{\phi}) = 
v^{\ev}(P)+\nwt(P)v(\phi)$, eventually.
\end{cor}

\noindent
If $P$ is homogeneous, then $P^\phi\in K^\phi\{Y\}$ is homogeneous with $\deg P^\phi = \deg P$, by Corollary~\ref{corconj-lemma}. 
Therefore $v(P^\phi)=\min_{i\in\N} v(P_i^\phi)$ for each $\phi$.
If $P_i\neq 0$, then $v(P_i^\phi)=v^{\ev}(P_i)+\nwt(P_i)v\phi$ eventually, and thus:

\begin{cor}\label{cor:nwt homog}
There is an $i\in\N$ such that 
$v(P^\phi)=v(P_i^\phi)$ eventually, and for any such $i$ we have
$v^{\ev}(P)=v^{\ev}(P_i)$ and $\nwt(P)=\nwt(P_i)$.
\end{cor}

\noindent
It follows from Lemmas~\ref{dwtP} and~\ref{comp-conjugation lemma} that for $g\in K^\times$ the
natural number $\nwt(P_{\times g})$ depends only on $vg$, so we can define a
function $\nwt_P\colon \Gamma \to \N$ by \nomenclature[X]{$\nwt_P(\gamma)$}{Newton weight of $P_{\times g}$, where $vg=\gamma$}
$$\nwt_P(\gamma)\ :=\ \nwt(P_{\times g})\ \quad \text{
for $g\in K^\times$ with $vg=\gamma$}.$$ 
Let $g\in K^\times$. If $g\asymp 1$, then
$v(P_{\times g}^\phi)=v(P^\phi)$ and $\nwt(P)=\nwt(P_{\times g})$, so 
$v^{\ev}(P)=v^{\ev}(P_{\times g})$. Thus 
$v^{\ev}(P_{\times g})$ depends only on $\gamma=vg$, and so we can set 
$v_P^{\ev}(\gamma):= v^{\ev}(P_{\times g})$ when $vg=\gamma$. Note also that 
$0\le \nwt(P_{\times g})\le \wt(P)$ for all $g\in K^\times$.

\nomenclature[X]{$v^{\ev}_P(\gamma)$}{$v^{\ev}(P_{\times g})$ where $vg=\gamma$}

\medskip\noindent
Let $\k$ be the residue field of $K$. Since $\Psi$ has no largest element, the $\Psi$-set $\Psi-v\phi$ of any compositional conjugate $K^{\phi}$ contains positive elements, so by the equivalence~(i)~$\Longleftrightarrow~$(iv) of Proposition~\ref{CharacterizationAsymptoticFields} we have $\derdelta\mathcal{O}\subseteq \smallo$ for the derivation $\derdelta=\phi^{-1}\der$ of $K^{\phi}$.
Thus for all $\phi$ the differential residue field of $K^{\phi}$ is $\k$ with the trivial derivation.

\medskip\noindent
With $\phi_0$ as above, Corollary~\ref{compconjval, cor 3}(iii) yields that the dominant part $D_{P^{\phi_0}}$ of $P^{\phi_0}$ is isobaric of weight~$\nwt(P)$, and that
for all $\phi\prec \phi_0$,
$$D_{P^{\phi}}\ =\ u(\phi)D_{P^{\phi_0}} \text{ for some $u(\phi)\in \k^{\times}$.}$$
%This equality has to be interpreted with care: 
%the residue $\d$-fields of $K^{\phi_0}$ and $K^{\phi}$ may be %different but have the same
%underlying field $\k$. Therefore, the differential polynomial 
%rings in
%the indeterminate $Y$ over these residue $\d$-fields have the same %underlying ring, and the equality refers to equality in this ring. 
Thus the degree of
$D_{P^{\phi_0}}$ is independent of the choice of $\phi_0$; we set
$$\ndeg P\ :=\ \deg D_{P^{\phi_0}}\ =\ \text{eventual value of $\ddeg P^{\phi}$,}$$  and call it the {\bf Newton degree} of $P$. \label{p:ndeg}
Likewise, we set 
$$\nval P\ :=\ \val D_{P^{\phi_0}}\ =\ \text{eventual value of $\dval P^{\phi}$,}$$ 
and call it the {\bf Newton multiplicity} of $P$ at $0$. Thus 
$$\val P\ \le\ \nval P\ \le\ \ndeg P\ \leq\ \deg P.$$
%Note that $\nval P^{\phi}=\nval P$ and $\ndeg P^{\phi}=\ndeg P$ for %all $\phi$. 
Given $P\in K\{Y\}^{\ne}$, the definitions of $v^{\ev}P\in \Gamma$ and $\nwt P, \ndeg P, \nval P\in \N$ involve $K$, but these quantities do not change when $K$ is replaced by an asymptotic extension $L$ such that $\Psi$ is cofinal in $\Psi_L$.
In particular, these quantities remain the same when passing to 
the algebraic closure of $K$.

%In the next section we consider more closely the Newton degree and
%Newton differential valuation. 

\index{differential polynomial!Newton degree}
\index{Newton!degree}
\index{degree!Newton}

\index{differential polynomial!Newton multiplicity}
\index{Newton!multiplicity}
\index{multiplicity!Newton}

\nomenclature[X]{$\ndeg(P)$}{Newton degree of $P$}
\nomenclature[X]{$\nval(P)$}{Newton multiplicity of $P$}

\subsection*{Behavior under compositional conjugation}
We continue to assume that $K$ is ungrounded, and $\phi$
continues to range over active elements in $K$. Let~$f\in K^\times$. The assumptions on $K$ remain valid for $K^f$, and the active elements in $K^f$ are the quotients $\phi/f$, hence
$$\nwt P^f = \nwt P, \quad \nval P^{f}=\nval P,\quad \ndeg P^{f}=\ndeg P. $$

\begin{lemma}\label{vepa} 
$v^{\ev}(P^f)=v^{\ev}(P)+\nwt(P)vf$.
\end{lemma}
\begin{proof} By Corollary~\ref{eventualequality} we have, eventually,
$$v((P^f)^{\phi/f})\ =\ v^{\ev}(P^f)+\nwt(P)(v\phi-vf)\
=\ v^{\ev}(P)+\nwt(P)v\phi,$$ 
which gives the desired result.
\end{proof}

\subsection*{Newton weight of linear differential operators}
\textit{In this subsection we assume that $K$ is ungrounded \textup{(}and so pre-$\d$-valued\textup{)}. We let $\phi$, $\phi_1$, $\phi_2$ range over active elements of $K$, $g$  over~$K^\times$, and $\gamma$ over $\Gamma$. }\/

\index{linear differential operator!Newton weight}
\index{Newton!weight}
\index{weight!Newton}

\medskip\noindent
Let $A=a_0+a_1\der+\cdots+a_r\der^r\in K[\der]^{\neq}$ ($a_0,\dots,a_r\in K$, $a_r\neq 0$).
Recall that in Section~\ref{sec:ldopv} we defined 
$\Pmu(A), \dwt(A)\in \N$ if $K$ has
small derivation, and so $\Pmu(A^\phi)$ and $\dwt(A^\phi)$ are defined. 
Set
$P:=a_0Y+a_1Y'+\cdots+a_rY^{(r)}\in K\{Y\}^{\ne}$. Then
$\Pmu(A^\phi)=\Pmu(P^\phi)$ and $\dwt(A^\phi)=\dwt(P^\phi)$. We call
$\nwt(A) := \nwt(P)$ the {\bf Newton weight\/} of $A$; so
$$\nwt(A)\  =\ \text{eventual value of $\dwt(A^{\phi})$}\  =\ \text{eventual value of $\dwv(A^{\phi})$.}$$  
Thus eventually $D_{P^\phi} = u(\phi)Y^{(w)}$ where $u(\phi)\in\k^\times$ and $w:=\nwt(A)$.
We also define the function  $\nwt_A\colon\Gamma\to\N$ by
$$\nwt_A(\gamma)\ :=\ \nwt(Ag)\ =\ \nwt_P(\gamma)\quad \text{for $\gamma=vg$,}$$
so for $a\in K^\times$ we have $\nwt_{aA}(\gamma)=\nwt_A(\gamma)$ and $\nwt_{Aa}(\gamma)=\nwt_A(va+\gamma)$.
We set 
$v^{\operatorname{e}}(A):=v^{\operatorname{e}}(P)$ and define 
$v_A^{\operatorname{e}}:=v_P^{\operatorname{e}}: \Gamma\to \Gamma$.
So, given any $\gamma$, 
\begin{equation}\label{eq:def of vAe}
v_{A^\phi}(\gamma)\ =\ v_A^{\operatorname{e}}(\gamma)+\nwt_A(\gamma)v\phi, \quad\text{ eventually,}
\end{equation}
by Corollary~\ref{eventualequality}. 
It follows that for $a\in K^\times$,
$$v_{aA}^{\operatorname{e}}(\gamma)\ =\ va+v_A^{\operatorname{e}}(\gamma),\qquad
v_{Aa}^{\operatorname{e}}(\gamma)\ =\ v_A^{\operatorname{e}}(va+\gamma).
$$
\begin{example}
We have $(\der g)^\phi=g'+\phi g\derdelta$ in $K^\phi[\derdelta]$. If $g\asymp 1$, then
$(\der g)^\phi \sim \phi g\derdelta$ for all~$\phi$ (since $K$ is pre-$\d$-valued), and
if $g\nasymp 1$, then eventually  $(\der g)^\phi \sim g'$.
Hence $\nwt_{\der}(0)=1$ and $\nwt_{\der}(\gamma)=0$ for $\gamma\ne 0$. From this and \eqref{eq:def of vAe} we get 
$$v_\der^{\operatorname{e}}(\gamma)\ =\ \begin{cases}
0	& \text{if $\gamma= 0$,} \\
\gamma'		& \text{if $\gamma\neq 0$.}\end{cases}$$
\end{example}

\noindent
Recall from Section~\ref{sec:ldopv} that if $\der\mathcal{O}\subseteq \smallo$, then $\exc(A)=\big\{\gamma:\ \Pmu_A(\gamma) \ge 1\big\}$. Thus for each $\phi$
we have $\exc(A^\phi)=\big\{\gamma:\ \Pmu_{A^\phi}(\gamma) \ge 1\big\}$. Define
$$\exc^{\operatorname{e}}(A)\ :=\ \big\{ \gamma:\ \nwt_A(\gamma) \geq 1 \big\}.$$
If $\phi_1\prec\phi_2$, then $\exc(A^{\phi_1})\subseteq\exc(A^{\phi_2})$ by Corollary~\ref{phi_1 phi_2}, so
$$\exc^{\operatorname{e}}(A)\ =\ \bigcap_\phi \exc(A^\phi)\ =\ \big\{ vg:\ \text{$A(g)\prec A^\phi g$ for all $\phi$} \big\}.$$
Call $\exc^{\operatorname{e}}(A)$ the set of {\bf eventual exceptional values} for $A$. By Lem\-ma~\ref{lem:exc(A) and ker} we have
$v(\ker^{\neq} A)=v(\ker^{\neq} A^\phi)\subseteq \exc(A^\phi)$ for each $\phi$, so
$v(\ker^{\neq} A) \subseteq \exc^{\operatorname{e}}(A)$. For $a\in K^\times$ we have
$\exc^{\operatorname{e}}(aA) =\exc^{\operatorname{e}}(A)$ and
$\exc^{\operatorname{e}}(Aa) = 
\exc^{\operatorname{e}}(A) - va$.

\index{exceptional value!eventual}
\index{eventual!exceptional value}
\index{linear differential operator!eventual exceptional value}
\nomenclature[Q]{$\exc^{\operatorname{e}}(A)$}{set of eventual exceptional values of $A$}

\begin{example}
If $r=0$, then $\exc^{\operatorname{e}}(A)=\emptyset$. Suppose $r=1$, and set $a:=a_0/a_1$. Then 
$$vg\in\exc^{\operatorname{e}}(A)
\quad\Longleftrightarrow\quad  \phi^{-1}a+\phi^{-1} g^\dagger\prec 1
\text{ for all }\phi\quad \Longleftrightarrow\quad v(a+g^\dagger)> \Psi,$$
where the first equivalence follows from a similar equivalence on exceptional values in Section~\ref{sec:ldopv}. Thus 
$\exc^{\operatorname{e}}(A)$ has at most one element.
In particular, $\exc^{\operatorname{e}}(\der)=\{0\}$. 
\end{example}

\noindent
In Chapter~\ref{ch:newtonian fields} we show that for 
certain $K$ the set $\exc^{\operatorname{e}}(A)$ is finite of size at most~$r$.

\section{Newton Degree and Newton Multiplicity}  \label{sec:ndeg and nval}

\noindent
{\em In this section $K$ is an ungrounded asymptotic field with value group 
$\Gamma\ne \{0\}$}.
We also restrict $\phi\in K^\times$ to be active, and we fix a differential polynomial $P\in K\{Y\}^{\ne}$. We let $\fm$ and
$\fn$ range over $K^\times$.

\begin{lemma}\label{newtonzero} Let $\fm$ be given, and suppose $P(f)=0$ for some $f\preceq \fm$
in some asymptotic field extension of $K$. Then $\ndeg P_{\times \fm}\ge 1$.
\end{lemma}
\begin{proof} Let $L$ be an asymptotic field extension of $K$ and suppose $f\in L$, $f\preceq \fm$, and $P(f)=0$. Then $f=a\fm$ with
$a\preceq 1$. We have $Q(a)=0$ for $Q:= P_{\times \fm}$, so $Q^{\phi}(a)=0$, and as $L^{\phi}$ has small derivation, this gives $D_{Q^{\phi}}(\bar{a})=0$,
and thus $\deg D_{Q^{\phi}} \ge 1$. As this holds for all $\phi$, we get
$\ndeg Q\ge 1$.   
\end{proof}

\begin{lemma}\label{ndevcon} There are $\fm_0, \fm_1\in K^\times$ such that  
$\nval P_{\times \fm}=\ndeg P_{\times \mathfrak m}=\val P$ for all $\fm\preceq \fm_0$, and
$\nval P_{\times \fm}=\ndeg P_{\times \mathfrak m}=\deg P$ for all $\fm\succeq \fm_1$. 
\end{lemma}
\begin{proof} Replacing $K$ and $P$ by $K^{\phi_0}$ and $P^{\phi_0}$ for some 
active $\phi_0$ in $K$ we arrange that $\der\smallo \subseteq \smallo$. Then by Proposition~\ref{compconjval, prop} we have for all $\phi \preceq 1$,
$$v_P(\gamma)\ \le\ v_{P^{\phi}}(\gamma)\ \le\ v_P(\gamma) + \wt(P)v\phi.$$ 
This holds in particular also for each nonzero homogeneous part of $P$ 
in place of~$P$. Let $d=\val(P)$, set $F:= P_d$. If $d< e$ with 
$G:=P_e\ne 0$, then by the above inequalities and Corollary~\ref{equalizer} there is $\gamma_0\in \Gamma$ such that
$v_{F^{\phi}}(\gamma)< v_{G^{\phi}}(\gamma)$, for all $\phi\preceq 1$ and all $\gamma\in \Gamma^{\ge \gamma_0}$. This gives the first part of the lemma, and the second part follows in the same way.
\end{proof}

\noindent
Next we state some results on {\em Newton degree\/} that follow easily from corresponding facts on {\em dominant degree} in Section~\ref{sec: dominant}, using
Lemma~\ref{comp-conjugation lemma}.

% or have similar proofs as these facts; 
%see Section~\ref{sec: dominant}. 

\begin{lemma}\label{ndnm} Let $a\in K$, $a\preceq 1$. Then we have
\begin{enumerate}
\item[\textup{(i)}] $\ndeg {P_{+a}}=\ndeg{P}$;
\item[\textup{(ii)}] if $a\prec 1$, then $\ndeg P_{\times a} \le \nval P=\nval P_{+a}$;
\item[\textup{(iii)}] if $a\asymp 1$, then
$\nval P_{\times a} = \nval P$, $\ndeg P_{\times a} = \ndeg P$. 
\end{enumerate}
\end{lemma}

\noindent
By (iii) of this lemma we have maps $\nval_P, \ndeg_P\colon\Gamma\to\N$
given by
$$\nval_P(\gamma)=\nval(P_{\times g}),\quad \ndeg_P(\gamma)=\ndeg(P_{\times g})\quad\text{for $g\in K^\times$ with $vg=\gamma$}.$$

\begin{cor} \label{aftranew} Let $a,b\in K$, $g\in K^\times$ 
be such that $a-b\preceq g$. Then 
$$\ndeg P_{+a,\times g}\ =\ \ndeg P_{+b,\times g}.$$
\end{cor}
%\begin{proof} Just note that $f:= (b-a)/g\in \mathcal{O}$ 
%and
%$P_{+b,\times g}=P_{+a, \times g, +f}$.
%\end{proof}

\noindent
Here is an analogue of Corollary~\ref{dn1cor,new}:

\nomenclature[X]{$\ndeg_P(\gamma)$}{Newton degree of $P_{\times g}$, where $vg=\gamma$}
\nomenclature[X]{$\nval_P(\gamma)$}{Newton multiplicity of $P_{\times g}$, where $vg=\gamma$}

\begin{cor}\label{dn1cor,new,Newton} 
$\val{P} =\val P_{\times \fm} \le \ndeg P_{\times\fm}\le \deg P_{\times \fm}=\deg P$, and 
$$\fm \prec \fn\ \Longrightarrow\ \nval P_{\times\mathfrak m}\ \leq\
\ndeg P_{\times\fm}\ \leq\ \nval P_{\times\mathfrak n}\  
\le\ \ndeg P_{\times\fn}.$$ 
\end{cor}

\noindent
Thus the functions $\nval_P, \ndeg_P\colon \Gamma \to \N$ are decreasing,
taking values in the finite set of $d\in \N$ for which $P_d\ne 0$.

\subsection*{The Newton degree on a set $\E$} 
In this subsection $\E\subseteq K^\times$ is 
$\preceq$-closed, that is, $\E\ne \emptyset$, and $\fm\in \E$ whenever 
$\fm \preceq \fn\in \E$. 
The {\bf Newton degree} of $P$ on $\E$, $\ndeg_{\E}P$, is the natural
number 
$$ \ndeg_{\E}P\ :=\ \max\!\big\{\!\ndeg(P_{\times\fm}):\ \fm\in \E\big\}.$$
By Corollary~\ref{dn1cor,new,Newton} we have $\val P\le \ndeg_{\E}P$.
If $Q\in K\{Y\}$, $Q\ne 0$, then clearly 
$$\ndeg_{\E}PQ\ =\ \ndeg_{\E}P + \ndeg_{\E}Q.$$ 
We have $\ndeg_{\E}P^{\phi}=\ndeg_{\E}P$. For $a\in K^\times$ we have
$\ndeg_{\E} P_{\times a} = \ndeg_{a\E}P$. 

\index{differential polynomial!Newton degree!on $\E$}
\index{Newton!degree!on $\E$}
\index{degree!Newton!on $\E$}
\nomenclature[X]{$\ndeg_{\E}(P)$}{Newton degree of $P$ on $\E$}

\begin{lemma} 
If $v(\mathcal E)$ does not have a smallest element, then 
$$\ndeg_{\mathcal E} P\ =\ \max\!\big\{\!\nval(P_{\times\mathfrak m}):\ \mathfrak m\in\mathcal E\big\}.$$
\end{lemma}
 
\begin{lemma}\label{lem:ndeg add conj} If $f\in \E$, then 
$\ndeg_{\E}P_{+f}=\ndeg_{\E} P$.
\end{lemma}

\noindent
Let $\E'\subseteq \E$ be $\preceq$-closed. Then for $f\in \E$ we have 
$\ndeg_{\E'}P_{+f}\le \ndeg_{\E} P$ by Lemma~\ref{lem:ndeg add conj}.
For $\gamma\in \Gamma$ and $\E=\{\fn: v\fn\ge \gamma\}$, we set 
$\ndeg_{\geq\gamma} P:= \ndeg_{\E} P$, so if $v\fm=\gamma$, then
$\ndeg_{\geq\gamma} P = \ndeg P_{\times \fm}$. 

\nomenclature[X]{$\ndeg_{\geq\gamma}(P)$}{Newton degree of $P$ on $\{\fn: v\fn\ge \gamma\}$}

\begin{cor}\label{cor:ndeg add conj} Let $a, b\in K$ and $\alpha, \beta\in \Gamma$ 
be such that $v(b - a)\ge \alpha$ as well as~$\beta \ge \alpha$. 
Then $\ndeg_{\geq\beta} P_{+b}\ \le\ \ndeg_{\geq\alpha}P_{+a}$. 
\end{cor}

\nomenclature[X]{$\ndeg_{\prec\fm}(P)$}{Newton degree of $P$ on $\{\fn: \fn \prec \fm\}$}

\noindent
For $\E:=\{\fn: \fn \prec \fm\}$ we set 
$\ndeg_{\prec \fm}P:= \ndeg_{\E}P$. If $\Gamma^{>}$ has no least element, then $v\big(\{\fn: \fn\prec\fm\}\big)$ has no least element, and so 
$$ \ndeg_{\prec \fm}P\ =\ \max\!\big\{\!\ndeg P_{\times\fn}:\ \fn\prec \fm\big\}\ =\ \max\!\big\{\!\nval P_{\times\fn}:\ \fn\prec \fm\big\}.$$

\begin{definition}
An algebraic differential equation with asymptotic side condition (for short: an {\bf asymptotic equation}) over $K$
is of the form
\begin{equation}\tag{E} \label{eq:asympt equ 0}
P(Y)=0, \qquad Y\in \E.
\end{equation}
For $g\in K^\times$ and $\E=\{y\in K^\times:\  y \prec g\}$ we also indicate \eqref{eq:asympt equ 0}
by
$$
P(Y)=0, \qquad Y\prec g.
$$
Likewise with $\preceq$ in place of $\prec$.

\medskip
\noindent
A {\bf solution} of \eqref{eq:asympt equ 0} is a $y\in \E$ such that $P(y)=0$.
The {\bf Newton degree} of \eqref{eq:asympt equ 0} is defined to be $\ndeg_{\E} P$.
\end{definition}

\index{equation!asymptotic}
\index{asymptotic equation}
\index{asymptotic equation!solution}
\index{subset!$\preceq$-closed}
\index{closed!$\preceq$-closed}

\noindent
Let $f$ be an element of a valued differential field extension of $K$.
We say that a solution~$y$ of~\eqref{eq:asympt equ 0} 
{\bf best approximates $f$}---tacitly: among solutions of  \eqref{eq:asympt equ 0}---if $y-f\preceq z-f$ for each solution~$z$ of \eqref{eq:asympt equ 0}. 
Of course, if $f\in K^\times$ is a solution of~\eqref{eq:asympt equ 0}, then $f$ is the unique solution of~\eqref{eq:asympt equ 0} that best approximates $f$.
Also, if $f\succ\E$, then $y-f\asymp f$ for all~$y\in\E$, so $f$ is best approximated by each solution of~\eqref{eq:asympt equ 0} in~$K$. 

\index{asymptotic equation!solution!best approximating an element}
\index{best approximation}

\begin{lemma}\label{lem: bestapp10-2}
Let $f$ be an element of a valued differential field extension of $K$ and $\fm\in\E$ with $f\preceq\fm$.
Suppose $y$ is a solution of the  asymptotic equation
\begin{equation}\label{eq:asympt equ, best approx, 1}
P_{\times \fm}(Y)\ =\ 0, \qquad Y\preceq 1
\end{equation}
that best approximates $\fm^{-1}f$. Then the solution $\fm y$ of~\eqref{eq:asympt equ 0} best approximates $f$.
\end{lemma}
\begin{proof}
Let $z$ be a solution of \eqref{eq:asympt equ 0}.
If $z\succ\fm$, then $z-f\sim z\succ\fm\succeq \fm y-f$,
and if $z\preceq\fm$, then $\fm^{-1}z$ is a solution of \eqref{eq:asympt equ, best approx, 1} and so
$y-\fm^{-1}f\preceq \fm^{-1}z-\fm^{-1}f$ and hence $\fm y-f\preceq z-f$. 
\end{proof}

\subsection*{The Newton degree in a cut}
In the next lemma, let $(a_{\rho})$ be a pc-sequence in~$K$, and
put $\gamma_{\rho}=v(a_{s(\rho)}-a_{\rho})$, where $s(\rho)$ is the immediate successor of $\rho$.   \nomenclature[X]{$d(P,(a_\rho))$}{eventual value of $\ndeg_{\geq\gamma_{\rho}}P_{+a_{\rho}}$}

\begin{lemma}\label{ndegeq}  
There is an index $\rho_0$ and $d\in \N$ such that for all $\rho > \rho_0$ 
we have
$\gamma_{\rho}\in \Gamma$ and $\ndeg_{\geq\gamma_{\rho}}P_{+a_{\rho}} = d$.
Denoting this number $d$ by $d\big(P,(a_{\rho})\big)$ to indicate its dependence on $P$
and  $(a_{\rho})$, we have $d\big(P,(a_{\rho})\big)=d\big(P,(b_{\sigma})\big)$
whenever $(b_{\sigma})$ is a pc-sequence in $K$ equivalent to $(a_{\rho})$.
\end{lemma}
\begin{proof} Take $\rho_0$ such that for all $\rho'>\rho\ge \rho_0$ we have $$\gamma_{\rho'}>\gamma_{\rho}, \qquad v(a_{\rho'}-a_{\rho})=\gamma_{\rho}\in \Gamma.$$ Then for such $\rho'$ and $\rho$ we have by Corollary~\ref{cor:ndeg add conj}, 
$$\ndeg_{\geq\gamma_{\rho'}}P_{+a_{\rho'}}\ \le\  \ndeg_{\geq\gamma_{\rho}}P_{+a_{\rho}}.$$
This yields the existence of $d(P,(a_{\rho}))$. For the second part, 
take $\rho_0$ as above such that also $\ndeg_{\geq\gamma_{\rho}}P_{+a_{\rho}} = d$ 
for all $\rho\ge \rho_0$. Let $(b_{\sigma})$ be a pc-sequence in $K$ equivalent 
to~$(a_{\rho})$, and take $\sigma_0$ such that
for all $\sigma'>\sigma \ge \sigma_0$ we have 
$$v(b_{\sigma'}-b_{\sigma})\ =\
\delta_{\sigma}\ :=\  v(b_{\sigma +1}-b_{\sigma})\in \Gamma.$$ 
We can also assume that $e\in \N$ and $\rho_0$,~$\sigma_0$
are such that $\ndeg_{\geq\delta_{\sigma}}P_{+b_{\sigma}} = e$ for all $\sigma\ge \sigma_0$, and $b_{\sigma}-a_{\rho} \prec a_{\rho}-a_{\rho_0}$ for all $\rho > \rho_0$
and $\sigma > \sigma_0$. Finally, we can assume
that $\delta_{\sigma} \ge \gamma_{\rho_0}$ for all 
$\sigma > \sigma_0$. Then for $\sigma > \sigma_0$ we have $v(b_{\sigma}-a_{\rho_0})=\gamma_{\rho_0}$ and
so $$e\ =\ \ndeg_{\geq\delta_{\sigma}}P_{+b_{\sigma}}\ \le\ \ndeg_{\geq\gamma_{\rho_0}}P_{+a_{\rho_0}}\ =\  d$$ 
by Corollary~\ref{cor:ndeg add conj}.
By symmetry we also have $d\le e$, so $d=e$.
\end{proof}

\noindent
We now associate to each pc-sequence $(a_{\rho})$ in $K$ an object $c_K(a_{\rho})$,
the {\bf cut} \index{pc-sequence!cut defined by}\index{cut!defined by a pc-sequence}\nomenclature[K]{$c_K(a_\rho)$}{cut defined by the pc-sequence $(a_\rho)$} defined by~$(a_{\rho})$ in $K$, such that if $(b_{\sigma})$ is 
also a pc-sequence in $K$, then 
$$c_K(a_{\rho})=c_K(b_{\sigma}) \ \Longleftrightarrow\ \text{$(a_{\rho})$ and $(b_{\sigma})$ are equivalent.}$$
We do this in such a way that the cuts $c_K(a_{\rho})$, with $(a_{\rho})$ a
pc-sequence in $K$, are the elements of a set $c(K)$. \nomenclature[K]{$c(K)$}{set of cuts defined by pc-sequences in the valued field $K$} For ${\bf a}\in c(K)$ we define
$$\ndeg_{\bf{a}}P\ :=\ d\big(P,(a_\rho)\big)\ =\ \text{eventual value of $\ndeg_{\geq\gamma_{\rho}}P_{+a_{\rho}}$,}$$
where $(a_{\rho})$ is any pc-sequence in $K$ with ${\bf a}=c_K(a_{\rho})$, using the notations of Lem\-ma~\ref{ndegeq}. 
We call $\ndeg_{\bf{a}}P$ the {\bf Newton degree of $P$ in the cut $\bf a$.} 
Let~$(a_{\rho})$ be a pc-sequence in~$K$ and ${\bf a}=c_K(a_{\rho})$.
For $y\in K$ the cut $c_K(a_{\rho} + y)$ depends only on~$({\bf a}, y)$, and so we can set ${\bf a}+y:=c_K(a_\rho+y)$. Likewise, for $y\in K^\times$ the cut~$c_K(a_{\rho}y)$ depends only on $({\bf a}, y)$, and so we can set 
${\bf a}\cdot y:=c_K(a_\rho y)$. We record some basic facts about $\ndeg_{\bf a} P$: 

\index{differential polynomial!Newton degree!in a cut defined by a pc-sequence}
\index{Newton!degree!in a cut defined by a pc-sequence}
\index{degree!Newton!in a cut defined by a pc-sequence}
\nomenclature[X]{$\ndeg_{\bf a}(P)$}{Newton degree of $P$ in the cut $\bf a$}

\begin{lemma} \label{lem:basic facts on deg_a}
Let $(a_{\rho})$ be a pc-sequence in $K$, ${\bf a}=c_K(a_\rho)$. Then
\begin{enumerate}
\item[\textup{(i)}] $\ndeg_{\bf a} P \leq \deg P$;
\item[\textup{(ii)}] $\ndeg_{\bf a} P^\phi = \ndeg_{\bf a} P$;
\item[\textup{(iii)}] $\ndeg_{\bf a} P_{+y} = \ndeg_{{\bf a}+y} P$  for $y\in K$;
\item[\textup{(iv)}] if $y\in K$ and $vy$ is in the width of $(a_\rho)$, then $\ndeg_{\bf a} P_{+y} = \ndeg_{\bf a} P$;
\item[\textup{(v)}] $\ndeg_{\bf a} P_{\times y} = \ndeg_{{\bf a}\cdot y} P$ for $y\in K^\times$;
\item[\textup{(vi)}] if $Q\in K\{Y\}^{\neq}$, then $\ndeg_{\bf a} PQ = \ndeg_{\bf a} P + \ndeg_{\bf a} Q$;
\item[\textup{(vii)}] if $P(\ell)=0$ for some pseudolimit $\ell$ of $(a_{\rho})$ in an asymptotic field extension of~$K$, then $\ndeg_{\bf a} P \ge 1$;
\item[\textup{(viii)}] if $L$ is an asymptotic field extension of $K$ and $\Psi$ is cofinal in $\Psi_L$, then $\Psi_L$ has no largest element and 
$\ndeg_{\mathbf a} P = \ndeg_{\mathbf a_L} P$, where
$\mathbf a_L=c_L(a_\rho)$.
\end{enumerate}
\end{lemma}

\begin{proof} Most of these items are routine or follow easily from earlier lemmas; in particular, 
item (vii) from  Lemma~\ref{newtonzero}. Item (iv) follows from (iii). 
\end{proof}

%\begin{lemma}
%Suppose $a_\rho\leadsto 0$.  Let $\E=\{y\in K^\times:vy\in W\}$ 
%where $W\subset%eq\Gamma_\infty$ is the width of $(a_\rho)$. 
%If $\E\neq\emptyset$, then $\ndeg_{\bf a} P \geq \ndeg_{\E} P$. 
%\marginpar{Joris claims that we always have equality, but I don't see this.} 
%If $\E=\emptyset$, then $(a_\rho)$ is a c-sequence in $K$, and
%$\ndeg_{\bf a} P = \val P$.
%\end{lemma}
%\begin{proof}
%Set $\gamma_\rho=v(a_{\rho+1}-a_\rho)\in\Gamma_\infty$. Removing some 
%initial terms we arrange that $(\gamma_\rho)$ is strictly increasing, and 
%$\gamma_\rho=v(a_\rho)\in \Gamma$ for each $\rho$.
%Then for each $\rho$ we have 
%$\mathcal E\subseteq \{y\in K^\times:\ vy \geq \gamma_\rho\}$, so 
%$\ndeg_{\geq\gamma_\rho} P_{+a_\rho} = \ndeg_{\geq\gamma_\rho} P$. 
%Therefore, if  $\E\neq\emptyset$, then $\ndeg_{\bf a} P \geq \ndeg_{\E} P$.
%The last part follows from 
%Lemma~\ref{ndevcon}. 
%\end{proof}

\subsection*{The case of order $1$} In this case we derive
properties of the Newton multiplicity and Newton degree that for higher order
require stronger assumptions on $K$ and a lot more work:
see the introduction to Chapter~\ref{ch:The Dominant Part and the Newton Polynomial}. {\em In this subsection we
assume that $K$ is $H$-asymptotic with rational asymptotic integration.}\/
Let $P\in K[Y, Y']^{\ne}$.

\begin{lemma}\label{newtordone} There are $w\in \N$, $\alpha\in \Gamma^{>}$, and
$A\in K[Y]^{\ne}$, such that eventually
$$P^{\phi}\ =\ \phi^w A(Y)(Y')^w + R_{\phi},\quad R_{\phi}\in K^{\phi}[Y, Y'],\quad v(R_{\phi}) > v(P^{\phi}) + \alpha.$$
\end{lemma}
\begin{proof} Let $P=\sum_{j\in J} A_{j}(Y)(Y')^j$ with finite nonempty
$J\subseteq \N$ and $A_j\in K[Y]^{\ne}$ for all $j\in J$. Then
$$P^{\phi}\ =\ \sum_{j\in J} \phi^j A_j(Y)(Y')^j,\qquad v(\phi^jA_j)\ =\ jv(\phi) + v(A_j).$$
Since $\Psi^{\downarrow}$ has no largest element we have $w\in J$
such that eventually 
$$w v(\phi) + v(A_w)\  <\  jv(\phi) + v(A_j)\ \text{ for all $j\in J\setminus \{w\}$.}$$ As $K$ has rational asymptotic integration we also have
$\alpha\in \Gamma^{>}$ such that eventually  $w v(\phi) + v(A_w) +\alpha < jv(\phi) + v(A_j)$ for all $j\in J\setminus \{w\}$. Then the conclusion
of the lemma holds with this $w$, $\alpha$ and $A:= A_w$. 
\end{proof}

\noindent
Note that $\nwt(P)=w$ for $w$ as in the lemma above. 

\begin{prop}\label{propnewtordone} There exists $\gamma\in \Gamma^{>}$ such that
for all $g\in K$,
\begin{align*} 0<vg< \gamma\ &\Longrightarrow\  \nval(P)\ =\ \nval(P_{\times g})\ =\ \ndeg(P_{\times g}),\\
-\gamma < vg < 0\  &\Longrightarrow\  \ndeg(P)\ =\ \ndeg(P_{\times g})\ =\ \nval(P_{\times g}).
\end{align*}
\end{prop}
\begin{proof} Take $w$,~$\alpha$,~$A$ as in Lemma~\ref{newtordone}. Subtracting
from $A$ a polynomial $D\in K[Y]$ with $D\prec  A$ and decreasing $\alpha$ if necessary we arrange that
$\val A = \dval A$ and $\deg A = \ddeg A$, so
$$\nval P = d+w, \quad \ndeg P= e+ w, \text{ where $d:=\val A$, $e:=\deg A$.}$$ Set $B:=A(Y)(Y')^w$ and take active $\phi_0$ in $K$ such that for all $\phi\preceq \phi_0$,
$$P^{\phi}\ =\ \phi^w B + R_{\phi},\quad R_{\phi}\in K^{\phi}[Y, Y'],\quad v(R_{\phi})\ >\ v(P^{\phi}) + \alpha.$$
So for $\phi\preceq \phi_0$ and $g\in K^\times$ with $g\prec 1$ we have, in $K^{\phi}[Y,Y']$,
\begin{align*} 
P_{\times g}^{\phi}\ &=\ \phi^w B_{\times g} + (R_{\phi})_{\times g},\\
v(\phi^w B_{\times g})\ &\le\ w v(\phi) +v(B)+ (d+w+1)v(g)\ =\ v(P^{\phi}) + (d+w+1)v(g),\\
v\big((R_{\phi})_{\times g}\big)\ &\ge\ v(R_{\phi}).
\end{align*} 
Fix $g\in K^\times$ with
$0<(d+w+1)v(g) < \alpha$. Then we have for $\phi\preceq \phi_0$,
$$ v(P^{\phi}_{\times g})\ =\ v(\phi^w B_{\times g})\ <\ v\big((R_{\phi})_{\times g}\big).$$
Now in $K^{\phi}[Y,Y']$ we have $(Y')_{\times g}=\phi^{-1}g'Y+gY'$. 
From $\phi \prec g^{\dagger}$, eventually, we get
$(Y')_{\times g}\sim \phi^{-1}g'Y$, eventually. Now $A=\sum_{i=d}^ea_iY^i$
with all $a_i\in K$, $v(A)=v(a_d)=v(a_e)$. So eventually in $K^{\phi}[Y,Y']$,
$$ B_{\times g}\ =\ \left(\sum_{i=d}^ea_ig^iY^i\right)\cdot \big((Y')_{\times g}\big)^w\ \sim\ \phi^{-w}a_dg^d(g')^wY^{d+w},$$
so $P^{\phi}_{\times g}\sim a_d g^d(g')^w Y^{d+w}$, eventually. Thus
$\nval P_{\times g}=\ndeg P_{\times g}=d+w$.

Next, let $g\in K^\times$, $g\succ 1$, $\phi\preceq \phi_0$. Then 
$v(\phi^w B_{\times g})\le wv(\phi)+v(B)= v(P^{\phi})$ and $v\big((R_{\phi})_{\times g}\big)\ge v(R_{\phi}) + (N+1)vg$ where 
$N:= \deg P$. A similar computation as for $g\prec 1$ also
gives $B_{\times g}\sim \phi^{-w}a_e g^e(g')^w Y^{e+w}$, eventually. It follows
that $\ndeg(P_{\times g})\ =\ \nval(P_{\times g})=e+w$ for 
$-\alpha < (N+1)vg < 0$.  
\end{proof}

\noindent
If the assumption of rational asymptotic integration  is weakened to $K$ having asymptotic integration, then
Lemma~\ref{newtordone}, Proposition~\ref{propnewtordone}, and their proofs
go through for $P\in K[Y,Y']^{\ne}$ of degree $\le 1$ in $Y'$.

\section{Using Newton Multiplicity and Newton Weight}\label{sec:nwt}

\noindent
{\em In this section $K$ is an $H$-asymptotic field with rational 
asymptotic integration.}\/ Also $P\in K\{Y\}^{\ne}$ and $\phi$ ranges over the
active elements of $K$.

\subsection*{Behavior of $vP(y)$}
We establish here a key fact
needed to construct immediate extensions in 
Section~\ref{sec:construct imm exts}, namely
Proposition~\ref{Z5a}. First we observe:
$$\nval P \ge 1\quad \Longleftrightarrow\quad  \text{$D_{P^{\phi}}(0)=0$,  eventually}\quad \Longleftrightarrow\quad
\text{$P(0)\prec P^{\phi}$, eventually.}$$

\begin{prop}\label{Z5a} Suppose $\nval(P)\ge 1$. Then there are $\beta\in \Gamma^{<}$ and a strictly increasing function $i\colon (\beta,0) \to \Gamma$ 
such that $vP(y)=i(vy)$ for all $y\in K$ with $\beta < vy < 0$. If $\Gamma$ is divisible there is such $i$ with the intermediate value property. 
\end{prop}

\noindent
Before establishing this proposition, we prove two lemmas.

\begin{lemma}\label{Z4} Let $Q\in K\{Y\}^{\neq}$ and suppose $P^\phi \succ Q^\phi$, eventually. Then there is $\beta\in \Gamma^{<}$ such that $P(y) \succ Q(y)$ for all $y\in K$ with $\beta < vy <0$. 
\end{lemma}
\begin{proof} For $\mu={\nwt}(P)$ and $\nu={\nwt}(Q)$ we have eventually
$$vP^\phi\ =\ v^{\ev}P + \mu v\phi, \ \qquad vQ^\phi\ =\ v^{\ev}Q + \nu v\phi.$$
The algebraic closure of $K$ has asymptotic integration, and
the above eventual equalities remain true there, 
as does $vP^\phi < vQ^\phi$, eventually. 
It follows that we have $\gamma\in \Gamma^{>}$
such that $vP^\phi + \gamma \le vQ^\phi$, eventually. 
Take $\phi$ such that
$vP^\phi + \gamma \le vQ^\phi$ and $v\phi \ge \psi(\gamma)$. Replace $K$, $P$, $Q$ by
$K^\phi$, $P^\phi$, $Q^\phi$ to arrange that $K$ has small derivation and $v^\flat P < v^\flat Q$. Next, replace
$P$, $Q$ by $P/a$, $Q/a$ for suitable $a\in K^\times$ to get $v^\flat P=0$, in
particular, $P,Q\in \mathcal{O}^\flat\{Y\}$. Let 
$$y\mapsto y^\flat\colon \mathcal{O}^\flat \to K^\flat$$ be the residue map, and
let $P^\flat$, $Q^\flat$ be the images
of $P$, $Q$ under the differential ring morphism 
$\mathcal{O}^\flat\{Y\} \to K^\flat\{Y\}$ that extends this residue map by sending each $Y^{(n)}$ to $Y^{(n)}$. Then $P^\flat\ne 0$ and $K^\flat$ is an $H$-asymptotic field with asymptotic couple $\big(\Gamma^\flat, \psi|(\Gamma^\flat)^{\ne}\big)$.
% so $K^\flat$ is $H$-asymptotic with small derivation, and 
%has asymptotic integration by Lemma~\ref{specasymptint}. 
Applying Corollary~\ref{betterdifpol} to $K^\flat$ and $P^\flat$ instead of $K$ and~$P$, we get 
$\beta\in (\Gamma^\flat)^{<}$
such that $P^\flat(y^\flat)\ne 0$ for all $y\in K$ with $\beta < vy < 0$. 
Also $Q^\flat=0$, so $Q^\flat(y^\flat)=0$ for those~$y$. Thus
$vP(y) \in \Gamma^\flat$ and $vQ(y)>\Gamma^\flat$ for those~$y$. 
\end{proof}

\begin{lemma}\label{Z5} Suppose $\nval(P)\ge 1$ and set $Q:= P-P(0)$.
There is $\beta\in \Gamma^{<}$ such that $Q(y) \succ P(0)$ for all
$y\in K$ with $\beta < vy <0$, and so $P(y)\sim Q(y)$ for those~$y$.
\end{lemma}
\begin{proof} Note that $Q\ne 0$ and $Q(0)=0$. From $\nval(P)\ge 1$ we get
$P(0)\prec P^{\phi}$, eventually, so
$P^{\phi}\sim Q^{\phi}\succ P(0)$, eventually. It remains to apply Corollary~\ref{betterdifpol} if $P(0)=0$, and Lemma~\ref{Z4} if $P(0)\ne 0$. 
\end{proof}

\noindent
Proposition~\ref{Z5a} now follows from Lemma~\ref{Z5} and 
Corollaries~\ref{ps6e} and~\ref{betterdifpol}.

\subsection*{More on $vP(y)$}
In this subsection
we let $y$,~$z$ range over $K$ and $\beta$,~$\gamma$ over $\Gamma$, and
we assume that $P\in K\{Y\}^{\ne}$ has order~$r$. 
Corollary~\ref{betterdifpol} gives us an element $\alpha$ of $\Gamma$ and natural numbers 
$d, d_0, e_1,\dots, e_{r-1}$ such that for 
some $\beta<0$ we have, for all~$y$,
$$ \beta < vy < 0\ \Longrightarrow\  vP(y)\ =\ \alpha + d\psi(vy) + d_0 vy - \sum_{i=1}^{r-1} e_i\chi^{i}(vy).$$ 
Moreover, this property determines $\alpha, d, d_0,e_1,\dots, e_{r-1}$
uniquely. Recall also that if~$P$ is homogeneous, then $d_0=\deg P$,
and that if $P(0)=0$, then $d_0\ge 1$.  
We consider $\alpha + d\psi(vy)$ as the main term in the
sum above as $vy$ tends to $0$, and our goal is to show that
$\alpha$ and $d$ coincide with the ``eventual'' quantities from 
Section~\ref{evtbeh}:

\begin{prop}\label{alpha=vev} We have $\alpha=v^{\ev}(P)$ and $d=\nwt(P)$. 
If $P(0)=0$, then
for some $\beta<0$ and all $y,z\in K$ with
$\beta < vy < vz=0$ we have $P(y) \succ P(z)$. 
\end{prop}
\begin{proof} By the proof of Corollary~\ref{betterdifpol}
the statements in this subsection preceding Proposition~\ref{alpha=vev}
remain true when replacing $K$ by its algebraic closure, with the same
tuple $(\alpha, d, d_0, e_1,\dots, e_{r-1})$. So we assume in the rest of
the proof that $\Gamma$ is divisible.
We now combine tricks
from Section~\ref{Applicationtodifferentialpolynomials}:
passing to an elementary extension, compositional conjugation, and
coarsening.  Set $B:= \Gamma^{<}$. 
Let $L$ be an elementary extension of $K$ and let $\theta\in L^\times$
be such that $\Psi < v\theta\in \Psi_L$. Let
$B_L$ be the convex hull of $B$ in $\Gamma_L$. As in the proof of 
Proposition~\ref{ps6b}, our $\alpha$ and $d$ are related by
$\alpha=\alpha_0-dv\theta$, with $\alpha_0\in \Gamma_L$ obtained by applying
Lemma~\ref{ps6a} to $L^{\theta}$,~$B_L$,~$P^{\theta}$ in the role of
$K$,~$B$,~$P$. We now also apply
Lemma~\ref{alphadelta} to $L^{\theta}$,~$B_L$,~$P^{\theta}$ in the role of
$K$,~$B$,~$P$. In our situation we have
$$ \Delta=\Delta(B_L):=\big\{\delta\in \Gamma_L:\ |\delta|< \epsilon \text{ for all $\epsilon\in \Gamma^{>}$}\big\},$$
a convex subgroup of $\Gamma_L$, and so 
$\alpha_0\equiv v(P^\theta) \bmod \Delta$ by Lemma~\ref{alphadelta}. Therefore, 
$$\alpha \equiv v(P^\theta)-dv\theta \mod \Delta.$$
In terms of the original structure $K$ and using terminology from 
Section~\ref{evtbeh}, it follows that 
for each $\epsilon\in \Gamma^{>}$ we have
$$\left|\alpha - \big(v(P^\phi)-dv\phi\big)\right| < \epsilon, \quad \text{ eventually.}$$ 
Also $v(P^\phi)=v^{\ev}(P) + \nwt(P)v\phi$, eventually. Thus for each $\epsilon\in \Gamma^{>}$,
$$\left|\big(\alpha-v^{\ev}(P)\big) + \big(d-\nwt(P)\big)v\phi\right|< \epsilon, \quad \text{ eventually.}$$
Since $K$ has asymptotic integration, we get $d=\nwt(P)$, 
and then $\alpha=v^{\ev}(P)$. The second part of the proposition follows from the 
second part of 
Lemma~\ref{alphadelta}.
\end{proof} 

\noindent
Now we have the following variant of Corollary~\ref{ps6c}:

\begin{cor}\label{cleanve} There is $\beta<0$ such that 
for all $y\in K$ with $\beta<vy<0$,
\begin{align*} 
vP(y)\ &=\ v^{\ev}(P)+\nwt(P)\psi(vy) + d_0vy -\sum_{i=1}^{r-1}e_i\chi^i(vy),\ \text{ and thus}\\
  vP(y)\ &=\ v^{\ev}(P)+\nwt(P)\psi(vy) + \gamma(y), \quad 
|\gamma(y)| \le (\deg P)\cdot|vy|.
\end{align*}
\end{cor}

\noindent
Here is another interesting consequence of Proposition~\ref{alpha=vev}:

\begin{cor}\label{nwtzero} There exists $\beta<0$ such that $\nwt_{P}(\gamma)= 0$
for all $\gamma\in (\beta,0)$. 
\end{cor}
\begin{proof} With $\beta<0$ and $d, d_0, e_1,\dots, e_{r-1}$ as above, we have
for all $y$,
$$ \beta < vy < 0\ \Longrightarrow\  vP(y)\ =\ v^{\ev}(P) + \nwt(P)\psi(vy) + d_0vy-\sum_{i=1}^{r-1}e_i\chi^i(vy). $$
Now let $\beta< \gamma < 0$. We claim that then $\nwt_{P}(\gamma)=0$.
To see this, take $g\in K$ with $vg=\gamma$, and let $y\in K$ with $vy<0$ 
be so small that 
$\beta < vy + \gamma < \gamma$ and $\psi(\gamma) < \psi(vy)$. 
Then $\psi(v(gy))=\psi(\gamma)$, and thus $\chi(v(gy))=\chi(\gamma)$, so
$$vP_{\times g}(y)\ =\ v^{e}(P) + \nwt(P)\psi(\gamma) + d_0\gamma -\sum_{i=1}^{r-1}e_i\chi^i(\gamma) + d_0 vy.$$
Here only the term $d_0vy$ depends on $vy$, so by Proposition~\ref{alpha=vev} applied to $P_{\times g}$,
$$ v^{\ev}(P_{\times g})\ =\  v^{\ev}(P) + \nwt(P)\psi(\gamma) + d_0\gamma -\sum_{i=1}^{r-1}e_i\chi^i(\gamma), \qquad 
     \nwt(P_{\times g})\ =\ 0.$$
Note also that in replacing $P$ by $P_{\times g}$, the quantity $d_0$ does not change and the quantities $e_1,\dots, e_{r-1}$ become $0$.  
\end{proof} 

\subsection*{Evaluation at pc-sequences} In the discussion following Corollary~\ref{eventualequality} we defined for $P\in K\{Y\}^{\ne}$ a function $v_P^{\ev}\colon \Gamma\to \Gamma$. This function behaves 
much like~$v_P$, at least piecewise:  

\begin{lemma}\label{partitionve} The set
$\Gamma$ is a finite union of subsets $\Gamma(\mu,d)$ with $0\le \mu\ \le \wt(P)$ and $P_d\ne 0$, 
such that for all distinct
$\alpha,\beta\in \Gamma(\mu,d)$ we have
$$v^{\ev}_{P}(\alpha)-v_{P}^{\ev}(\beta)\ =\ d\cdot(\alpha-\beta) +o(\alpha-\beta).$$
Thus if $P(0)=0$, then $v_P^{\ev}$ is strictly
increasing on each set $\Gamma(\mu,d)$.
\end{lemma}
\begin{proof} Since $\Psi$ has no largest element, we can take an elementary 
extension $K^*$ of the asymptotic field $K$ with an element $\theta\in K^*$
such that $v\theta\in \Psi_{K^*}$, and $v\theta>v\phi$ for all $\phi$.
%with an active $\theta\in K^*$ \marginpar{``active'' not yet defined?} such 
%that $\theta \prec \phi$ for all $\phi$. 
For
$\gamma\in \Gamma$ and $g\in K^\times$ such that $vg=\gamma$ we have
$$v_P^{\ev}(\gamma)\ =\ v(P_{\times g}^\theta)-\nwt(P_{\times g})v\theta\ =\ 
v_{P^{\theta}}(\gamma)-\nwt(P_{\times g})v\theta.$$
The nonzero homogeneous parts of $P^{\theta}$ have the same degrees as the
nonzero homogeneous parts of $P$. For
$0\le \mu \le \wt(P)$ and $P_d\ne 0$, define $\Gamma(\mu,d)$ to be the set 
of $\gamma\in \Gamma$ such that $\nwt(P_{\times g})=\mu$ for
$g\in K^\times$ with $vg=\gamma$, and $v_{P^{\theta}}(\gamma)= v_Q(\gamma)$ with
$Q:= P_d^{\theta}$. It remains to use Corollary~\ref{vP-Lemma}. 
\end{proof}

\noindent
In the next lemma $(a_{\rho})$ is a pc-sequence 
in $K$ with a pseudolimit $e$ in an immediate asymptotic extension 
$E$ of $K$. Let $G(Y)\in E\{Y\}\setminus E$ and set
$\gamma_{\rho}:= v\big(a_{s(\rho)}-a_\rho\big)$. 

\begin{lemma}\label{flpc1new} There is a pc-sequence $(b_\lambda)$ in $K$ equivalent to 
$(a_\rho)$ such that $\big(G(b_{\lambda})\big)$ is a pc-sequence in $E$ with $G(b_\lambda) \leadsto G(e)$. 
\end{lemma}
\begin{proof}
By removing some initial $\rho$'s we arrange that 
$\gamma_{\rho}=v(e-a_\rho)\in \Gamma$ for all~$\rho$, and $(\gamma_{\rho})$
is strictly increasing. By removing also some indices $\rho$ corresponding to limit ordinals we arrange in addition that for each $\rho$ we have
$\delta_{\rho}\in \Gamma^{>}$ such that $\gamma_{s(\rho)}-\gamma_\rho > \delta_{\rho}$ and $\gamma_{\rho}- \gamma_{\rho'} > \delta_{\rho}$
whenever $\rho > \rho'$. 
Take $g_\rho\in K$ with
$v(g_\rho)=\gamma_\rho$ and define $u_\rho\in E$ by 
$a_\rho - e = g_\rho u_\rho$, so $u_\rho\asymp 1$. Take $x_{\rho}\in K$, subject 
for now only to $-\delta_{\rho} < vx_{\rho} < 0$, and put 
$b_{\rho}:= a_{\rho} + g_{\rho} x_{\rho}\in K$ and $y_{\rho}: = u_{\rho} + x_{\rho}\in E$. Then
$$ y_{\rho} \sim x_{\rho}, \quad b_\rho - e =g_{\rho}y_{\rho}, \quad
\gamma_{\rho} - \delta_{\rho} < v(b_{\rho}-e) < \gamma_{\rho}.$$ Thus
$(b_{\rho})$ pseudoconverges to $e$ and has the same width as $(a_\rho)$, so
 by Lemma~\ref{pc5} it is a pc-sequence in $K$ equivalent to $(a_{\rho})$. We have 
$$G(b_\rho) - G(e)\ =\ 
\sum_{|\i|\ge 1} G_{(\i)}(e)(g_\rho y_\rho)^{\i}\qquad\text{where $G_{(\i)}=\frac{G^{(\i)}}{\i !}$.}$$
Put $g_{\i}:=G_{(\i)}(e)\in E$ for 
$|\i|\ge 1$. Then
\begin{align*} G(b_\rho)-G(e)\ &=\ 
\sum_{|\i|\ge 1} g_{\i}(g_\rho y_\rho)^{\i}\
=\ P(g_{\rho}y_{\rho})\ =\ P_{\times g_{\rho}}(y_\rho), \text{ where }\\
 P(Y)\ &:=\  
\sum_{|\i|\ge 1} g_{\i}Y^{\i}\in E\{Y\}, \text{ so $\deg P\ge 1$, $P(0)=0$.}
\end{align*}
Corollary~\ref{cleanve} with $E$ instead of $K$ gives for each 
$\rho$ a
$\mu_{\rho}\in \big\{0,\dots, \wt(P)\big\}$
and an $\epsilon_{\rho}\in \Gamma^{>}$ 
such that for all $y\in E$ with
$-\epsilon_\rho < vy < 0$,
$$ v\big(P(g_\rho y)\big)\ =\ v_P^{\ev}(\gamma_{\rho})+
\mu_{\rho}\psi(vy) + \gamma(\rho,y), \quad |\gamma(\rho, y)|\le (\deg P)|vy|.$$
Take these $\epsilon_{\rho}$ so small that $(\deg P)\epsilon_{\rho} < \delta_{\rho}/4$.
Lemma~\ref{partitionve} with $E$ in place of~$K$ gives sets 
$\Gamma(\mu,d)$ with 
$0\le \mu\le \wt(P)$ and $1\le d\le \deg P$. 
Passing to suitable cofinal subsequences we arrange that for a fixed such 
$\mu$ and $d$ we have
$\gamma_{\rho}\in \Gamma(\mu,d)$ for all $\rho$, and $\mu_{\rho}$ is constant 
as a function of $\rho$. Then
$$v_P^{\ev}(\gamma_{\rho'})-v_P^{\ev}(\gamma_{\rho})\ =\ 
d(\gamma_{\rho'}-\gamma_{\rho}) + o(\gamma_{\rho'}-\gamma_{\rho})\ \text{ when $\rho'>\rho$.}$$ 
For $\rho' > \rho$ we have $\gamma_{\rho'}-\gamma_{\rho} > \delta_{\rho'}$ as well as $\gamma_{\rho'}-\gamma_{\rho}> \delta_{\rho}$, so 
$|vy_{\rho'}-vy_{\rho}|< \gamma_{\rho'}-\gamma_{\rho}$, and thus
$\psi(vy_{\rho'})-\psi(vy_{\rho})=o(\gamma_{\rho'}-\gamma_{\rho})$ by 
Lemma~\ref{BasicProperties}(ii).  
We now impose on~$x_{\rho}$ the further restriction 
$-\epsilon_{\rho} < vx_{\rho}<0$. Since $vx_{\rho}=vy_{\rho}$, it follows that for 
$\rho'>\rho$ we have 
$|vy_{\rho}|+ |vy_{\rho'}|< \epsilon_{\rho}+\epsilon_{\rho'}$, so in view of
$(\deg P)\epsilon_{\rho}< \delta_{\rho}/4$, 
$$v\big(P(g_{\rho'}y_{\rho'})\big)-v\big(P(g_{\rho}y_{\rho})\big)\ =\ 
d(\gamma_{\rho'}-\gamma_{\rho}) + o(\gamma_{\rho'}-\gamma_{\rho})+ \epsilon, \quad |\epsilon|< (\gamma_{\rho'}-\gamma_{\rho})/2.$$
Then $v\big(G(b_{\rho})-G(e)\big)=v\big(P(g_{\rho}y_{\rho})\big)$ is strictly 
increasing as a function of $\rho$, so $G(b_{\rho}) \leadsto G(e)$, as promised. 
\end{proof}

\subsection*{Substituting powers of $Y'$} In this subsection we assume $q\in \N^{\ge 1}$. For use in proving Lemma~\ref{lem:eval at uporho} we establish a lower bound on the Newton weight of the differential polynomial $P^{\times q}=P\big((Y')^q\big)\in K\{Y\}$ introduced in Section~\ref{Compositional Conjugation}. 

\begin{lemma}\label{lem:bound on nwt(Ptimesq)}
$\nwt(P^{\times q})\geq dq$ for some $d\in\N$ with $P_d\neq 0$.
\end{lemma}
\begin{proof}
We have $(P^{\times q})_{qi}=(P_i)^{\times q}$ for each $i\in\N$ and $(P^{\times q})_{j}=0$ for $j\in\N\setminus q\N$, so
$\nwt(P^{\times q})=\nwt\!\big( (P_i)^{\times q} \big)$ for some $i\in\N$ with $P_i\neq 0$, by Corollary~\ref{cor:nwt homog}. Hence we may assume that $P$ is homogeneous. Setting $d=\deg P$ we then need to show that $\nwt(P^{\times q})\geq dq$.
By Corollary~\ref{cor:derivatives of powers} we have a homogeneous $E\in K\{Y\}^{\neq}$ of degree $w:=\wt(P)$ such that 
$$P^{\times q}\  =\ (Y')^{dq-w} \cdot E^{\times}.$$
The differential polynomial $E^{\times \phi}=E^\phi_{\times\phi}$ is homogeneous of degree $w$, and so its dominant part $D_{E^{\times\phi}}$ is homogeneous of degree $w$,
with 
$$D_{(E^\times)^\phi}\ =\ D_{E^{\times \phi}(Y')}\ =\ u\,D_{E^{\times \phi}}(Y')\ \text{ for some $u\in \k^\times$.}$$ 
By Corollary~\ref{wtundercomp} we have $$\wt\!\big(D_{E^{\times \phi}}(Y')\big)\ =\ \deg D_{E^{\times\phi}} + 
\wt D_{E^{\times\phi}}\ =\ w+\dwt(E^{\times\phi}).$$
Hence $\nwt(E^\times)\geq w$ and so $$\nwt(P^{\times q})\ =\ \nwt\!\big((Y')^{dq-w}\big)+\nwt(E^\times)\ \geq\ (dq-w)+w\ =\  dq,$$
as desired.
\end{proof}

\begin{cor}\label{cor:bound on nwt(Ptimesq)}
Suppose $P(0)=0$ and $q\geq 2$. Then there is $\beta\in \Gamma^{<}$ 
 such that 
$$\big\{vP\big((y')^q\big):\  y\in K,\ \beta<vy<0\big\}$$ is disjoint
from some interval $(\gamma_0, \delta_0)$ with $\gamma_0\in \Psi$, $\delta_0\in (\Gamma^{>})'$. 
\end{cor}
\begin{proof} Set $d:= \nwt(P^{\times q})$. Then $d\ge 2$ by   
Lemma~\ref{lem:bound on nwt(Ptimesq)}. Next, Corollary~\ref{cleanve} gives
$\alpha,\beta\in \Gamma$ with $\beta<0$, $d_0,e_1,\dots, e_n\in \N$ such that for all $y\in K$,
$$\beta\ <\  vy\  <\  0\ \Longrightarrow\ vP((y')^q)\ =\ \alpha + d_0vy + d\psi(vy)
-\sum_{i=1}^n e_i\chi^i(vy).$$
From the remarks preceding Proposition~\ref{alpha=vev} %at the beginning of the subsection ``More on $vP(y)$''
it follows that $d_0\ge 1$.  
It now remains to apply Corollary~\ref{cor:not cofinal in Psi}. 
\end{proof}

\section{Constructing Immediate Extensions}\label{sec:construct imm exts}

\noindent
{\em In this section $K$ is an $H$-asymptotic field with rational asymptotic 
integration}. Here is the main result of this section:

\begin{theorem}\label{thm:immediate}
Every pc-sequence in $K$ has a pseudolimit in some immediate asymptotic extension of $K$.
\end{theorem}

\noindent
This is of course just an alternative formulation of Theorem~\ref{immHasint}. 
We do not know a direct proof of Theorem~\ref{thm:immediate} 
using evaluation at pc-sequences, along the lines of Sections~\ref{sec:pceq} 
and~\ref{sec:cimex}. Instead we depend heavily on the preceding facts on 
Newton weight, Newton degree, and Newton multiplicity,  
in particular, on Proposition~\ref{Z5a}.

\medskip
\noindent
We let $\phi$ range over the active elements of $K$, and
$a$, $b$, $y$ over $K$. Also,  
$\fm$, $\fn$, $\fd$, $\fv$, and~$\fw$ range over $K^\times$,
and $P$,~$Q$ range over $K\{Y\}^{\neq}$. 

\subsection*{Vanishing}
Recall from Section~\ref{sec:ndeg and nval} that
$$\ndeg_{\prec \fv}P\ =\ \max\!\big\{\!\ndeg P_{\times \fm}:\ \fm \prec \fv\big\}\ =\ \max\!\big\{\!\nval P_{\times \fm}:\ \fm \prec \fv\big\}.$$   
Let $\ell$ be an element in some asymptotic extension of $K$ 
such that $\ell\notin K$ and $v(K-\ell)=\big\{v(a-\ell):\ a\in K\big\}$ has no largest element (so $\ell$ has no best approximation in $K$). 
Note that then $v(K-\ell)\subseteq \Gamma$, and that 
there is a divergent pc-sequence in~$K$ with pseudolimit $\ell$. 

We say that {\bf $P$ vanishes at $(K,\ell)$}
if for all $a$ and $\fv$ with $a-\ell \prec \fv$ we have 
$\ndeg_{\prec \fv}P_{+a} \ge 1$. (Intuitively, ``$P$ vanishes at $(K,\ell)$'' means that $K$ thinks $P$ could have a zero near 
$\ell$.)
Let $Z(K,\ell)$ be the set of all $P$ that vanish
at $(K,\ell)$. 
Here are some frequently used basic facts: \begin{enumerate}
\item  $P\in Z(K,\ell)\ \Longleftrightarrow\ P_{+b}\in Z(K, \ell-b)$;
\item  $P\in Z(K,\ell)\ \Longleftrightarrow\ P_{\times \fm}\in Z(K,\ell/\fm)$;
\item  $P\in Z(K,\ell)\ \Longrightarrow\ PQ\in Z(K,\ell)$;  
\item $P\in K\ \Longrightarrow\ P\notin Z(K,\ell)$.
\end{enumerate}

\noindent
If $P\notin Z(K,\ell)$, we have $a$,~$\fv$ with 
$a-\ell \prec \fv$ and
$\ndeg_{\prec \fv}P_{+a} =0$, and then also $\ndeg_{\prec \fv}P_{+b} =0$ for any
$b$ with $b-\ell\prec \fv$, by Lemma~\ref{lem:ndeg add conj}.  

\index{differential polynomial!vanishing at $(K,\ell)$} 
\nomenclature[Z]{$Z(K,\ell)$}{set of $P\in K\{Y\}^{\neq}$ which vanish at $(K,\ell)$}

\begin{lemma}\label{Z0} $Y-b\notin Z(K,\ell)$.
\end{lemma}
\begin{proof} Take $a$ and $\fv$ such that $a-\ell \prec \fv \asymp b-\ell$. Then
for $P:= Y-b$ and $\fm\prec \fv$ we have $P_{+a,\times \fm}=\fm Y + (a-b)$ and $\fm \prec a-b$,
so $\ndeg_{\prec \fv}P_{+a}=0$. \end{proof}

\begin{lemma}\label{Z1} Suppose $P\notin Z(K,\ell)$, and let $a$,~$\fv$ be such that 
$a-\ell \prec \fv$ and
$\ndeg_{\prec \fv}P_{+a} =0$. Then $P(f)\sim P(a)$ for all 
$f$ in all asymptotic extensions of~$K$ with 
$f-a\asymp \fm \prec \fv$ for some $\fm$. \textup{(}Recall: $\fm\in K^\times$ by convention.\textup{)}
\end{lemma}
\begin{proof} Let $f$ in an asymptotic extension $E$ of $K$
satisfy $f-a\asymp \fm \prec \fv$, so
$f=a+\fm u$ with $u\asymp 1$ in $E$. Now 
$$P_{+a,\times \fm}\ =\ P(a) + R\qquad\text{ with $R\in K\{Y\}$, $R(0)=0$,}$$ 
so
$$ P^\phi_{+a,\times\fm}\ =\ P(a) +R^{\phi}.$$
From $\ndeg P_{+a,\times \fm}=0$ we get $R^{\phi}\prec P(a)$, eventually.
Thus 
$$P(f)\ =\ P_{+a,\times\fm}(u)\ =\ P(a) + R^{\phi}(u)\quad \text{in $E^{\phi}$,}$$ with $R^{\phi}(u)\preceq R^{\phi}\prec P(a)$, eventually, in $E^{\phi}$, so $P(f)\sim P(a)$.  
\end{proof}

\noindent
Note that the conclusion applies to $f=\ell$, and so
for $P$ and $a$,~$\fv$ as in the lemma we have $P(\ell)\sim P(a)$, 
hence $P(\ell)\ne 0$ and $vP(\ell)\in \Gamma$.
In particular, if $K$ is also an $H$-field and $P\notin Z(K,\ell)$, then 
$\sgn P(a)$ is independent of $a\in K$ and $\fv$ subject to $a-\ell \prec \fv$ and
$\ndeg_{\prec \fv} P_{+a}=0$, and so $\sgn(P,\ell):= \sgn P(a)$
does not depend on the choice of such $a$,~$\fv$.  

\begin{lemma}\label{Zv} Suppose that $P,Q\notin Z(K,\ell)$. Then $PQ\notin Z(K,\ell)$.
\end{lemma}
\begin{proof} Take $a$,~$b$,~$\fv$,~$\fw$ such that $a-\ell \prec \fv$, $b-\ell \prec \fw$
and $$\ndeg_{\prec \fv} P_{+a}\ =\ \ndeg_{\prec \fw} Q_{+b}\  =\ 0.$$ 
We can assume $a-\ell \preceq b-\ell$. Take $\fn\asymp a-\ell$ and $d\in K$ with $d-\ell \prec \fn$.
Then $d-\ell \prec \fv$ and $d-\ell \prec \fw$, so 
$\ndeg_{\prec \fv} P_{+a}=\ndeg_{\prec \fv} P_{+d}=0$, and so 
$\ndeg_{\prec \fn} P_{+d}=0$. In the same way we obtain 
$\ndeg_{\prec \fn} Q_{+d}=0$. Hence 
$\ndeg_{\prec \fn} (PQ)_{+d} =0$.
\end{proof}

\begin{lemma}\label{Zu} Suppose $P\in Z(K,\ell)$. Then for each $b$ there is $a$ with
$a-\ell \prec b-\ell$ and $P(a)\not\asymp P(b)$.
\end{lemma}
\begin{proof} Let $\fv\asymp b-\ell$ and take $a_1\in K$ with $a_1-\ell \prec \fv$,
so $\ndeg_{\prec \fv} P_{+a_1}\ge 1$, so we have $\fm\prec \fv$ with
$\nval(P_{+a_{1},\times \fm})\ge 1$. Then by Proposition~\ref{Z5a}, the set 
$$\big\{vP(a_1+\fm y):\ \beta < vy <0\big\}$$ 
is infinite for each $\beta\in \Gamma^{<}$
, so we can take
$y$ with $vy<0$ and $a_1+\fm y -\ell \prec \fv$ and 
$P(a_1 +\fm y)\not\asymp P(b)$. Then $a:= a_1+\fm y$ has the desired property. 
\end{proof}

\begin{lemma}\label{Zw, newton} Suppose $P,Q\notin Z(K,\ell)$ and $P-Q\in Z(K,\ell)$.
Then $P(\ell)\sim Q(\ell)$.
\end{lemma}
\begin{proof} By Lemma~\ref{Zv} we have $b$ and $\fv$ such that
$$\ell - b \prec \fv, \qquad \ndeg_{\prec \fv}P_{+b}\ =\ \ndeg_{\prec \fv}Q_{+b}\ =\ 0.$$
Replacing $\ell$ by $\ell-b$ and $P, Q$ by $P_{+b}, Q_{+b}$ we arrange $b=0$, that is,
$$\ell \prec \fv, \qquad \ndeg_{\prec \fv}P\ =\ 
\ndeg_{\prec \fv}Q\ =\ 0, $$
in particular, $P(0)\ne 0$ and $Q(0)\ne 0$.
By Lemma~\ref{Z1} we have for all $a\prec \fv$,
$$P(a)\sim P(0)\sim P(\ell), \qquad Q(a) \sim Q(0) \sim Q(\ell).$$
If $P(\ell)\not\sim Q(\ell)$, then $P(0)\not\sim Q(0)$, so $(P-Q)(a)\asymp (P-Q)(0)$
for all $a\prec \fv$, contradicting $P-Q\in Z(K, \ell)$ by Lemma~\ref{Zu}.
Thus $P(\ell)\sim Q(\ell)$. 
\end{proof}

\subsection*{Constructing immediate extensions}
As in the previous subsection, $\ell$ is an element in some asymptotic extension of $K$ 
such that $\ell\notin K$ and $v(K-\ell)$ has no largest element.

\begin{lemma}\label{zdf, newton} Suppose $Z(K,\ell)=\emptyset$. Then $P(\ell)\ne 0$ for all $P$, 
and $K\<\ell\>$ is an immediate asymptotic extension of $K$.
Suppose that $g$ in an asymptotic extension~$L$ of~$K$ satisfies
$v(a-g)=v(a-\ell)$ for all $a$. Then there is a unique valued differential field
embedding $K\<\ell \> \to L$ over $K$ that sends $\ell$ to $g$. 
\end{lemma}
\begin{proof} Clearly $P(\ell)\ne 0$ for all $P$. Let any nonzero element $f=P(\ell)/Q(\ell)$ of the asymptotic field extension $K\<\ell\>$ of $K$ be given.
Lemma~\ref{Zv} gives $a$ and $\fv$ such that
$$a-\ell \prec \fv,\qquad \ndeg_{\prec \fv} P_{+a}\ =\ \ndeg_{\prec \fv} Q_{+a}\ =\ 0,$$
and so $P(\ell) \sim P(a)$ and $Q(\ell) \sim Q(a)$ by Lemma~\ref{Z1}, and thus $f\sim P(a)/Q(a)$.
It follows that $K\<\ell\>$ is an immediate extension of $K$.

It is clear that $Z(K,g)=Z(K,\ell)=\emptyset$, so $g$ is 
$\d$-transcendental over~$K$ and~$K\<g\>$ is an immediate extension of $K$, by the first part of the proof.
Given any~$P$ we take $a$ and $\fv$ such that $a-\ell \prec \fv$ and
$\ndeg_{\prec \fv} P_{+a} =0$. Then $P(a)\sim P(g)$ and $P(a) \sim P(\ell)$, and
thus $vP(g)=vP(\ell)$.  Hence the unique differential field embedding $K\<\ell\> \to L$ over $K$ that sends $\ell$ to $g$ is also a valued field embedding.
\end{proof}

\begin{lemma}\label{zda, newton} Suppose that $Z(K,\ell) \ne \emptyset$ and $P$ is an element of
$Z(K,\ell)$ of minimal complexity. Then $K$ has an immediate asymptotic extension $K\<f\>$ with $P(f)=0$ and
$v(a-f)=v(a-\ell)$ for all $a$, and such that if $E$ is an asymptotic extension of~$K$ and $e\in E$ satisfies $P(e)=0$ and
$v(a-e)=v(a-\ell)$ for all $a$, then there is a unique valued differential field
embedding $K\<f \> \to E$ over~$K$ that sends $f$ to $e$.
\end{lemma}
\begin{proof} Let $P$ have order $r$ and take $p\in K[Y_0,\dots,Y_r]$ such that
$$P\ =\ p\big(Y, Y',\dots, Y^{(r)}\big).$$ It is clear from Lemma~\ref{Zv} that $p$ is irreducible. Consider the domain 
$$K[y_0,\dots, y_r]\ =\ K[Y_0,\dots, Y_r]/(p), \quad \text{$y_i=Y_i + (p)$ for $i=0,\dots,r$,}$$
and let $K(y_0,\dots,y_r)$ be its fraction field. We extend $v\colon K^\times \to \Gamma$ to $$v \colon\  K(y_0,\dots, y_r)^\times \to \Gamma$$ as follows. Let $s\in K(y_0,\dots, y_r)^\times$ and take
$g\in K[Y_0,\dots, Y_r]$, $h\in K[Y_0,\dots, Y_{r-1}]$ with 
$g\big(Y, Y',\dots, Y^{(r)}\big)\notin Z(K,\ell)$ (so $g\notin pK[Y_0,\dots, Y_r]$) and $h\ne 0$ such that 
$s=g(y_0,\dots, y_r)/h(y_0,\dots, y_{r-1})$. Then $vg\big(\ell,\dots,\ell^{(r)}\big)$ and $vh\big(\ell, \dots, \ell^{(r-1)}\big)$ lie in~$\Gamma$ by the comments following Lemma~\ref{Z1}.
We claim that 
$$vg\big(\ell, \ell',\dots, \ell^{(r)}\big)-vh\big(\ell,\dots, \ell^{(r-1)}\big)\in \Gamma$$
depends only on $s$ and not on the choice of $g$ and $h$. To see this, let
$g_1\in K[Y_0,\dots, Y_r]$, $h_1\in K[Y_0,\dots, Y_{r-1}]$ be such that 
$g_1\big(Y,\dots, Y^{(r)}\big)\notin Z(K,\ell)$, 
$h_1\ne 0$, and $s=g_1(y_0,\dots, y_r)/h_1(y_0,\dots, y_{r-1})$.
Then 
$$gh_1-g_1h \in pK[Y_0,\dots, Y_r], \quad (gh_1)\big(Y,\dots, Y^{(r)}\big), (g_1h)\big(Y,\dots, Y^{(r)}\big)\notin Z(K,\ell),$$ which
yields the claim by Lemma~\ref{Zw, newton}. We now set, for $g$,~$h$ as above, 
$$vs\ :=\  vg\big(\ell, \ell',\dots, \ell^{(r)}\big)-vh\big(\ell,\dots, \ell^{(r-1)}\big),$$
or more suggestively, 
$$vs\ =\ v\big(G(\ell)/H(\ell)\big)\in \Gamma,\ \text{ with $G=g\big(Y,\dots,Y^{(r)}\big)$, $H=h\big(Y,\dots, Y^{(r-1)}\big)$.}$$
Let $s_1, s_2\in K(y_0,\dots, y_r)^\times$. Then $v(s_1s_2)=vs_1+vs_2$ follows
easily by means of Lemma~\ref{Zv}. Next, assume also
$s_1+s_2\ne 0$; to prove that $v\colon K(y_0,\dots, y_r)^\times \to \Gamma$ is a valuation it remains to show that then $v(s_1+s_2) \ge \min(vs_1, vs_2)$.
For $i=1,2$ we have $s_i=g_i(y_0,\dots, y_r)/h_i(y_0,\dots, y_{r-1})$ where
$$ 0\ne g_i\in K[Y_0,\dots, Y_r],\quad 0\ne h_i\in K[Y_0,\dots, Y_{r-1}],$$
and $g_i$ has lower degree in $Y_r$ than $p$. Then for $s:= s_1+s_2$ we have 
$$s\ =\ g(y_0,\dots, y_r)/h(y_0,\dots, y_{r-1}),\quad  g\ :=\ g_1h_2 + g_2h_1,\quad h\ :=\ h_1h_2,$$
and so $g\ne 0$ (because $s\ne 0$) and $g$ also has lower degree in $Y_r$ than $p$.
In particular, $g\big(Y,\dots, Y^{(r)}\big)\notin Z(K,\ell)$, hence
$$vs\ =\ v\big(g(\ell,\dots, \ell^{(r)})/h(\ell,\dots, \ell^{(r-1)})\big),$$ 
and so by working in the valued field $K\<\ell\>$ we see that $vs\ge \min(vs_1, vs_2)$, as promised.
Thus we now have $K(y_0,\dots, y_r)$ as a valued field extension of $K$.
To show that $K(y_0,\dots, y_r)$ has the same residue field as $K$, consider an
element $s=g(y_0,\dots, y_r)\notin K$ with nonzero $g\in K[Y_0,\dots, Y_r]$ of lower degree in $Y_r$ than $p$;
it suffices to show that $s\sim b$ for some $b$. Set 
$G:= g(Y, \dots, Y^{(r)})$ and take $a$ and $\fv$ with
$a-\ell \prec \fv$ and $\ndeg_{\prec \fv}G_{+a}=0$. Then 
$G(\ell)\sim G(a)$ by Lemma~\ref{Z1}, so for $b:= G(a)$ we have 
$$v(s-b)\ =\ v\big(g(y_0,\dots, y_r)-b\big)\ =\ v(G(\ell)-b)\  >\  vb,$$ 
that is, $s\sim b$. 
This finishes the proof that the valued field $K(y_0,\dots, y_r)$ is an immediate extension of $K$.

\medskip\noindent
Next we equip $K(y_0,\dots, y_r)$ with the derivation extending the derivation of $K$
such that $y_i'=y_{i+1}$ for $0\le i < r$. Setting $f:= y_0$ we have $f^{(i)}=y_i$
for $i=0,\dots,r$, $K\<f\>=K(y_0,\dots, y_r)$, and $P(f)=0$. Note that
$v(f-a)=v(\ell -a)$ by Lemma~\ref{Z0}.

\medskip\noindent
In order for $K\langle f \rangle= K\big(f, \dots,f^{(r-1)}, f^{(r)}\big)$ 
to be an asymptotic field, it suffices by Proposition~\ref{pa5} to show 
that $K\big(f,\dots,f^{(r-1)}\big)$ is asymptotic in $K\langle f \rangle$. 
To do that we are going to apply Lemma~\ref{pa4} with 
$$L=K\big(f,\dots,f^{(r-1)}\big), \quad 
F= K\langle f \rangle,\ \text{ and }\  U=K\big[f,\dots,f^{(r-1)}\big].$$
Consider elements $s\in K\big[f,\dots,f^{(r-1)}\big]$ and 
$b\in K^\times$ such that $b\prec 1$; it suffices to show that then
$s' \prec b^\dagger$ if $s \prec 1$, and $s' \preceq b^\dagger$ 
if $s \asymp 1$. If $s\asymp 1$, take $d\in K$ with $s\sim d$ and use
$d'\preceq  b^\dagger$ to reduce to the case $s\prec 1$. So we assume
$s \prec 1$ and wish to show that $s'\prec b^\dagger$.
This inequality certainly holds when $s'=0$, so assume 
$s'\ne 0$. Take $g\in K[Y_0,\dots,Y_{r-1}]$, $g\ne 0$ with
$s=g\big(f,\dots,f^{(r-1)}\big)$, and 
take  $g_1, g_2\in K[Y_0,\dots,Y_{r-1}]$ such that, in $K\{Y\}$,
$$g\big(Y,\dots,Y^{(r-1)}\big)'\ =\  
g_1\big(Y,\dots,Y^{(r-1)}\big) + g_2\big(Y,\dots,Y^{(r-1)}\big)Y^{(r)}.$$
Then $$s'\ =\ g\big(f,\dots,f^{(r-1)}\big)'\ =\  
g_1\big(f,\dots,f^{(r-1)}\big) + g_2\big(f,\dots,f^{(r-1)}\big)f^{(r)},$$ and
for all $a$,
$$g\big(a,\dots,a^{(r-1)}\big)'\ =\  
g_1\big(a,\dots,a^{(r-1)}\big) + g_2\big(a,\dots,a^{(r-1)}\big)a^{(r)}.$$
There are two cases to consider:

\case[1]{$p$ has degree $>1$ in $Y_r$, or $g_2=0$.}
Then we can take $a$ such that
$$s \sim g\big(a,\dots,a^{(r-1)}\big),\qquad 
s' \sim g\big(a,\dots,a^{(r-1)}\big)'.
$$
Since $s\prec 1$, the left-hand side gives 
$g\big(a,\dots,a^{(r-1)}\big)' \prec b^\dagger$, and then the right-hand side yields
$s' \prec  b^\dagger$.

\case[2]{$p$ has degree $1$ in  $Y_r$, and $g_2 \ne 0$.}
Then 
$$g_1 + g_2Y_r= \frac{h_1p + h_2}{h}, \quad h,h_1,h_2\in K[Y_0,\dots,Y_{r-1}],\ h, h_1\ne 0,$$ so
$0\ne s'=h_2\big(f,\dots,f^{(r-1)}\big)/h(f,\dots, f^{(r-1)})$, so $h_2\ne 0$.
Let 
\begin{align*} G:= g(Y, \dots, Y^{(r-1)}), \quad  &H:= h(Y, \dots, Y^{(r-1)}),\\
\quad H_1:=h_1(Y,\dots, Y^{(r-1)}), \quad &H_2:= h_2(Y, \dots, Y^{(r-1)}).
\end{align*}
By Lemma~\ref{Zv} there is $\fv$ such that for some $a$,
$$ a-\ell \prec \fv,\quad  \ndeg_{\prec \fv}G_{+a}=\ndeg_{\prec \fv}H_{+a}=\ndeg_{\prec \fv}(H_1)_{+a}=\ndeg_{\prec \fv}(H_2)_{+a}=0.$$
Fix such $\fv$, and let $A\subseteq K$ be the set of all $a$ satisfying the above. Then for $a\in A$ we can apply Lemma~\ref{Z1} with $\ell$ in the role of $f$, so
$h(a,\dots, a^{(r-1)})\ne 0$ and
\begin{align*}
s\ =\ g\big(f,\dots,f^{(r-1)}\big)\ &\sim\ g\big(a,\dots,a^{(r-1)}\big),\\
 \frac{h_1\big(f,\dots,f^{(r-1)}\big)}{h\big(f,\dots, f^{(r-1)}\big)}\ 
 &\sim\ \frac{h_1\big(a,\dots,a^{(r-1)}\big)}{h\big(a,\dots, a^{(r-1)}\big)},\\
s'\ =\ \frac{h_2\big(f,\dots,f^{(r-1)}\big)}{h\big(f,\dots, f^{(r-1)}\big)}\ &\sim \frac{h_2\big(a,\dots,a^{(r-1)}\big)}{h\big(a,\dots, a^{(r-1)}\big)},\\
g\big(a,\dots,a^{(r-1)}\big)'\ &=\
\frac{h_1\big(a,\dots,a^{(r-1)}\big)}{h\big(a,\dots, a^{(r-1)}\big)}P(a)+
\frac{h_2\big(a,\dots,a^{(r-1)}\big)}{h\big(a,\dots, a^{(r-1)}\big)}.
\end{align*}
Now $g\big(a,\dots,a^{(r-1)}\big)\asymp s \prec 1$, and so $vg\big(a,\dots,a^{(r-1)}\big)'$  as well as 
$$ v\frac{h_1\big(a,\dots,a^{(r-1)}\big)}{h\big(a,\dots, a^{(r-1)}\big)}\   \text{ and }\ v\frac{h_2\big(a,\dots,a^{(r-1)}\big)}{h\big(a,\dots, a^{(r-1)}\big)}=vs'$$   do not depend on $a\in A$. On the other hand $\ndeg_{\prec \fv}P_{+a} >0$,
so by Lem\-ma~\ref{Zu},
$vP(a)$ is not constant as a function of $a\in A$. 
Hence
$s'\asymp g\big(a,\dots,a^{(r-1)}\big)'$, and thus $s'\prec b^\dagger$. 
This finishes the proof
that $K\langle f \rangle$ is an asymptotic field. 

\medskip\noindent
Suppose now that $e$ in an asymptotic field extension $E$ of $K$ satisfies
$P(e)=0$ and $v(a-e)=v(a-\ell)$ for all $a$.  By
Lemma~\ref{Z1} we have $vQ(e)=vQ(f)$ for all $Q\notin Z(K,\ell)$, in particular,
$Q(e)\ne 0$ for all $Q$ of lower complexity than $P$. Thus we have a
differential field embedding $K\langle f \rangle \to E$ over $K$ sending $f$ to $e$,
and this is also a valued field embedding.
\end{proof}

\noindent
Here are two immediate consequences of Lemmas~\ref{zdf, newton} and~\ref{zda, newton}.

\begin{cor}\label{zdcor-} If $K$ is asymptotically maximal, then $K$ is spherically complete.
\end{cor}
\begin{proof} Suppose $K$ is not spherically complete. Then we have a
divergent pc-se\-quence~$(a_{\rho})$ in $K$, and thus a pseudolimit $\ell$
in an asymptotic extension of $K$ with $\ell\notin K$. 
Now Lemmas~\ref{zdf, newton} (if $Z(K,\ell)=\emptyset$) and~\ref{zda, newton} (if $Z(K,\ell)\ne \emptyset$)
provide a proper immediate asymptotic extension
of~$K$.
\end{proof}

\begin{cor}\label{zdcor} The asymptotic field $K$ has a spherically complete immediate asymptotic
extension. 
\end{cor}
\begin{proof} This follows by Zorn from Corollary~\ref{zdcor-}.
\end{proof}

\subsection*{Relation to the Newton degree in a cut}
As before, $\ell$ is an element in some asymptotic extension of $K$ such that $\ell\notin K$ and $v(K-\ell)$ has no largest element.  Let
$(a_{\rho})$ be a divergent pc-sequence in $K$ with pseudolimit 
$\ell$.

\begin{lemma}\label{ppkell} If $P(a_{\rho}) \leadsto 0$, then $P\in Z(K,\ell)$.
\end{lemma}
\begin{proof} Suppose $P\notin Z(K,\ell)$. Take $a$ and $\fv$ such that
$a-\ell\prec \fv$ and $\ndeg_{\prec \fv}P_{+a}=0$. Now 
$v(a-a_{\rho})=v(a-\ell)$, eventually, so by Lemma~\ref{Z1} we have
$P(a_{\rho}) \sim P(a)$ eventually, so $v\big(P(a_{\rho})\big)= v\big(P(a)\big)\ne \infty$ eventually. 
\end{proof}

\noindent
We now connect the notion of $P$ vanishing at $(K,\ell)$ with the Newton degree $\ndeg_{\bf a} P$ of $P$ in the cut ${\bf a}=c_K(a_\rho)$.

\begin{lemma}\label{dpkell} $\ndeg_{\bf a} P \ge 1\ \Longleftrightarrow\ 
P\in Z(K,\ell)$. More precisely,
$$\ndeg_{\bf a} P\  =\ 
\min\!\big\{\!\ndeg_{\prec \fv}P_{+a}:a-\ell\prec\fv\big\}.$$ 
\end{lemma}
\begin{proof} We may assume $v(\ell- a_{\rho})$
is strictly increasing with $\rho$. 
Given any index~$\rho$, take
$\fv\asymp \ell -a_{\rho}$, take $\rho'>\rho$, and set $a:= a_{\rho'}$. Then
$a-\ell \prec \fv$. Now
$\gamma_{\rho}:= v(\ell - a_{\rho})=v(a-a_{\rho})$, and
thus by Lemma~\ref{lem:ndeg add conj},
$$\ndeg_{\prec \fv}P_{+a}\ \le\ \ndeg_{\preceq \fv}P_{+a}\  =\ \ndeg_{\geq \gamma_{\rho}}P_{+a_{\rho}}.$$
It follows that $\min\!\big\{\!\ndeg_{\prec \fv}P_{+a}:a-\ell\prec\fv\big\}\leq \ndeg_{\bf a} P$.
For the reverse inequality,  let $a$ and $\fv$ be such that
$a-\ell \prec \fv$. Let $\rho$ be such that $\ell - a_{\rho}\preceq \ell-a$. 
Then $a_{\rho}-a \prec \fv$ and 
$\gamma_{\rho}:= v(\ell -a_{\rho})> v(\fv)$, so by 
Lemma~\ref{lem:ndeg add conj}:
$$\ndeg_{\geq\gamma_{\rho}}P_{+a_{\rho}}\ \le\  \ndeg_{\prec \fv}P_{+a_{\rho}}\ =\ \ndeg_{\prec \fv} P_{+a}.$$
Therefore $\ndeg_{\bf a} P\leq \min\!\big\{\!\ndeg_{\prec \fv}P_{+a}:a-\ell\prec\fv\big\}$.
\end{proof}

\noindent
Recall from the end of Section~\ref{Valdifcon} the notion of a minimal
differential polynomial over~$K$ of a pc-sequence in $K$.

\begin{cor}\label{zmindifpol} The following conditions on $P$ are equivalent:
\begin{enumerate}
\item[\textup{(i)}] $P$ is an element of $Z(K,\ell)$ of minimal complexity;
\item[\textup{(ii)}] $P$ is a minimal differential polynomial of $(a_{\rho})$ over $K$. 
\end{enumerate}
\end{cor}
\begin{proof} Assume (i). Then Lemma~\ref{zda, newton} provides 
an immediate extension
$K\<f\>$ of $K$ with $P(f)=0$ and $v(a-\ell)=v(a-f)$ for all $a$.
Hence $a_{\rho} \leadsto f$, and so
Lemma~\ref{flpc1new} gives a pc-sequence $(b_{\lambda})$ in $K$ 
equivalent to $(a_{\rho})$ such that $P(b_{\lambda})\leadsto 0$.
Conversely, if $(b_{\lambda})$ is a pc-sequence in $K$ equivalent to
$(a_{\rho})$ and $Q(b_{\lambda})\leadsto 0$, then $b_{\lambda}\leadsto \ell$, so $Q\in Z(K, \ell)$ by Lemma~\ref{ppkell}. Thus we have established (i)~$\Rightarrow$~(ii). 

Now assume (ii). Again, from Lemma~\ref{ppkell} we get $P\in Z(K,\ell)$.
The direction~(i)~$\Rightarrow$~(ii) shows that $P$ is of minimal complexity
in $Z(K,\ell)$. \end{proof}

\section{Special Cuts in $H$-Asymptotic Fields} \label{sec:special cuts}

\noindent
The cuts referred to in the title of this section
are given by certain jammed pc-sequences in 
$H$-asymptotic fields $K$ such as $\mathbb{T}$. Their divergence in $K$ is equivalent to $K$ having a rather subtle but important elementary property. To explain this for $K=\mathbb{T}$, consider the iterated logarithms 
$\ell_n\in \mathbb{T}$ defined by $\ell_0:=x$ and $\ell_{n+1}:=\log \ell_n$.
Then the corresponding sequence $(\upl_n)$ given by
$$\upl_n\ :=\ -(\ell_n{}^{\dagger}{}^\dagger)\ =\ 
\frac{1}{\ell_0} + \frac{1}{\ell_0\ell_1} + \cdots +  
\frac{1}{\ell_0\ell_1\cdots\ell_n}$$
is a jammed pc-sequence in $\mathbb{T}$ by Example~\ref{ex:iterated log sequence for T} and Lemma~\ref{pcz}
below. It is easy to check that this pc-sequence has no pseudolimit in $\mathbb{T}$.
By Corollary~\ref{psacgap} below, the divergence of this pc-sequence in $\mathbb{T}$ implies that for all 
$s\in \mathbb{T}$ there is~$g\succ 1$ in $\mathbb{T}$ 
such that $s-g^{\dagger\dagger} \succeq g^\dagger$.
A related elementary property of $\mathbb{T}$ plays a key role
in our work. It
corresponds to the fact that the sequence $(\upo_n)$  with
$$ \upo_n\ :=\ -2\upl_n' - \upl_n^2\ =\ \frac{1}{\ell_0^2} + \frac{1}{\ell_0^2\ell_1^2} + \cdots +  
\frac{1}{\ell_0^2\ell_1^2\cdots\ell_n^2}$$
is a divergent (jammed) pc-sequence in $\mathbb{T}$; see 
Corollary~\ref{psacrho}.

\medskip
\noindent
In this section and the next two
we treat this material for ungrounded $H$-asymptotic fields. More precisely, we assume in this section: 
\begin{quote}
\textit{$K$ is an ungrounded $H$-asymptotic field with $\Gamma\neq\{0\}$. \textup{(}Thus $\Psi:=(\Gamma^{\neq})^\dagger$ is nonempty and has no largest element.\textup{)}}\/
\end{quote}
\noindent
Hence $K$ is a pre-differential-valued field by Corollary~\ref{asymp-predif},
with contraction map
$$\chi\colon\ \Gamma^{<}\to \Gamma^{<}.$$ 

\subsection*{Transfinitely iterated logarithms}
For $f\in K^{\succ 1}$, take $\Log f\in K^{\succ 1}$ such that $(\Log f)'\asymp f^\dagger$, that is, $v(\Log f)=\chi(vf)$.  
If~$K$ is closed under logarithms, then we can choose such $\Log f$ with $(\Log f)' = f^\dagger$.\index{logarithm!on an asymptotic field} We now introduce ``iterated logarithms''~$\ell_\rho$, 
for possibly transfinite $\rho$.
More precisely, we construct a sequence
$(\ell_\rho)$ in $K^{\succ 1}$, \nomenclature[Z]{$(\ell_\rho)$}{logarithmic sequence} indexed by the ordinals $\rho$ less than some infinite limit ordinal $\kappa$, such that \begin{enumerate}
\item $\ell_{\rho'} \prec \ell_{\rho}$ whenever $\rho' > \rho$;
\item $(\ell_\rho)$ is coinitial in $K^{\succ 1}$: for each $f\in K^{\succ 1}$ there is an index $\rho$ with $\ell_\rho \preceq f$.
\end{enumerate}
We construct this sequence by transfinite
recursion: take any element $\ell_0\succ 1$ in~$K$, and
take $\ell_{\rho+1}:=\Log \ell_\rho$; if $\lambda$ is an infinite limit ordinal
such that all $\ell_\rho$ with $\rho<\lambda$ have already been chosen,
then we pick $\ell_\lambda$ to be any element in $K^{\succ 1}$ such that $\ell_\lambda\prec \ell_\rho$ for
all $\rho<\lambda$, if there is such an $\ell_\lambda$, while if there is no
such $\ell_\lambda$, we put $\kappa:=\lambda$. We call $(\ell_\rho)$ a 
{\bf logarithmic sequence} for $K$. \index{logarithmic!sequence} Note that if $L$ is an $H$-asymptotic field extension of~$K$ such that $\Gamma^<$ is cofinal in $\Gamma_L^<$, then $\Psi$ is cofinal in $\Psi_L$ (hence $\Psi_L$ also does not have a largest element), and a logarithmic sequence for~$K$ remains a logarithmic sequence for~$L$.

{\sloppy

\medskip
\noindent
Set $e_\rho:=v(\ell_\rho)\in\Gamma^{<}$, so $v(\ell_\rho^\dagger)=
e_\rho^\dagger$, and $e_{\rho+1}=\chi(e_\rho)$, and thus
$$ \qquad e_{\rho+1}^\dagger\ =\ e_{\rho}^\dagger-\chi(e_\rho)\
=\ e_{\rho}^\dagger - e_{\rho +1}\ >\ e_{\rho}^\dagger.$$ 
Therefore the sequence 
$(e_\rho)$ is strictly
increasing and cofinal in $\Gamma^{<}$, and the se\-quence~$(e_\rho^\dagger)$ is strictly increasing and cofinal in $\Psi$. 
%In particular,  $\ell_{\rho'} \flatter \ell_\rho$
%for $\rho<\rho'<\kappa$.
From $(\ell_\rho)$ we define the sequences~$(\upg_\rho)$ in~$K^\times$ and~$(\upl_\rho)$ in $K$ as follows:
$$
\upg_\rho\ :=\ \ell_{\rho}^\dagger, \qquad
\upl_\rho\ :=\  -\upg_\rho^\dagger\ =\ -\ell_\rho^\dagger{}^\dagger\ :=\ -(\ell_\rho^\dagger{}^\dagger).$$
% \quad b_\rho :=(1/l_\rho)'{}^\dagger=a_\rho+y_\rho,$$
Then $v(\upg_\rho)=e_\rho^\dagger$, so $\upg_{\rho'}\prec \upg_{\rho}$ 
for $\rho<\rho'<\kappa$.
}

\nomenclature[Z]{$(\upg_\rho)$}{$(\ell_{\rho}^\dagger)$}
\nomenclature[Z]{$(\upl_\rho)$}{$\big( {-(\ell_\rho^\dagger{}^\dagger)} \big)$}

\begin{exampleNumbered}\label{ex:iterated log sequence for T}
Suppose $K=\mathbb T$. Then the sequence $(\ell_n)$ given by $\ell_0=x$ and $\ell_{n+1}=\log(\ell_n)$ is a logarithmic sequence for $K$, 
and for each $n$ we have
$$\upg_n\ =\ \frac{1}{\ell_0\ell_1\cdots\ell_n}, \qquad \upl_n\ =\ \frac{1}{\ell_0} + \frac{1}{\ell_0\ell_1} + \cdots +  
\frac{1}{\ell_0\ell_1\cdots\ell_n}.$$
\end{exampleNumbered}

\begin{lemma}\label{pcz}
For  $\rho<\rho'<\kappa$ we have
$\upl_{\rho'}-\upl_{\rho} \sim \upg_{\rho+1}$. The sequence $(\upl_\rho)$
is a jammed pc-sequence of width $\big\{\gamma\in \Gamma_{\infty}: \gamma > \Psi\big\}$.
\end{lemma}
\begin{proof}
We have 
$$\upg_{\rho+1}\ =\ \ell_{\rho+1}'/\ell_{\rho+1}\ =\
(\Log \ell_{\rho})'/\ell_{\rho+1}\ \asymp\ 
\ell_{\rho}^\dagger/\ell_{\rho +1}\ =\  \upg_\rho/\ell_{\rho+1},$$
so $\upg_{\rho}/\upg_{\rho+1} \asymp \ell_{\rho+1}$. Hence
$$\upl_{\rho+1}-\upl_{\rho}\ =\ \upg_{\rho}^\dagger-\upg_{\rho+1}^\dagger\ =\ (\upg_{\rho}/\upg_{\rho+1})^\dagger\ \sim\
\ell_{\rho+1}^\dagger=\upg_{\rho+1}.$$
Let $\rho<\rho'<\kappa$; we want
$\upl_{\rho'}-\upl_\rho \sim \upg_{\rho+1}$, and have
just shown this for $\rho'=\rho+1$. When $\rho+1<\rho' < \kappa$, then
by Lemma~\ref{BasicProperties}(i) and
$\upl_{\rho'}-\upl_{\rho+1}=
(\upg_{\rho+1}/\upg_{\rho'})^\dagger$,
$$v(\upl_{\rho'}-\upl_{\rho+1})=\psi\bigl(v(\upg_{\rho+1}/\upg_{\rho'})\bigr)=
\psi\bigl(e_{\rho+1}^\dagger-e_{\rho'}^\dagger\bigr)>e_{\rho+1}^\dagger=v(\upg_{\rho +1}),$$
so $\upl_{\rho'}-\upl_\rho=(\upl_{\rho'}-\upl_{\rho+1}) +(\upl_{\rho+1}-\upl_\rho) \sim \upg_{\rho+1}$, as claimed.
Thus $(\upl_\rho)$ is a pc-sequence with $v(\upl_{\rho'} - \upl_\rho) =e_{\rho+1}^\dagger$ for
$\rho < \rho' < \kappa$. 
By Lemma~\ref{BasicProperties}(ii),
$$ e_{\rho' +1}^\dagger - e_{\rho+1}^\dagger\ =\ 
o(e_{\rho' +1}-e_{\rho +1})\ <\ |e_{\rho}|.$$
Given any convex subgroup $\Delta\ne \{0\}$ of $\Gamma$, we can take an index $\rho(\Delta)$ such that $e_{\rho(\Delta)}\in \Delta$, and then by the
last displayed inequality, 
$$ e_{\rho' +1}^\dagger - e_{\rho+1}^\dagger\in \Delta\ \text{ whenever }\ 
 \kappa > \rho' > \rho > \rho(\Delta).$$
Thus $(\upl_\rho)$ is jammed.
\end{proof}

\begin{remark}
By Lemmas~\ref{pcz} and~\ref{pc5} the sequence $(\upl_\rho+\upg_\rho)$ is also a pc-sequence in~$K$, and is equivalent to $(\upl_\rho)$.
\end{remark}

\noindent
By Proposition~\ref{indcof} below, pc-sequences 
$(\upl_{\rho})$ obtained from different choices of logarithmic sequence~$(\ell_{\rho})$ are
equivalent. The proof uses the properties of $\psi(*-\alpha)$ for $\alpha\in\Psi^\downarrow$; see Corollary~\ref{corhasc}.

\medskip\noindent
By part~(i) of Lemma~\ref{BasicProperties}, if $f, g\in K$ are active 
and $f\succ g$, then $f^\dagger - g^\dagger \prec f$. Below we frequently
use this fact without further mention.
Each $\upg_{\rho}$
is active in $K$, and if $a\in K$ is active, then $\upg_\rho\prec a$ for all big enough~$\rho$. 

\begin{prop}\label{indcof} Let $(u_{\sigma})$ be a well-indexed sequence 
of active elements in~$K$ such that
$\big(v(u_\sigma)\big)$ is strictly increasing and cofinal in 
$\Psi^{\downarrow}$. Then
some cofinal subsequence of $(-u_{\sigma}^\dagger)$ is a pc-sequence
equivalent to $(\upl_\rho)$. 
\end{prop}
\begin{proof} 
Take a pseudolimit $\upl$ of $(\upl_\rho)$ in some 
$H$-asymptotic field extension of $K$. We claim that
$v(\upl+u_{\sigma}^\dagger)\ =\ \psi\big({*-v(u_{\sigma}})\big)$ for all $\sigma$. 
Let any~$\sigma$ be given. Then $\upl+u_{\sigma}^\dagger= (\upl+\upg_{\rho}^{\dagger}) + 
(u_{\sigma}^\dagger-\upg_{\rho}^{\dagger})$ with
$\upl+\upg_{\rho}^{\dagger}=\upl-\upl_{\rho}\sim 
\upg_{\rho+1}$. Also, Lemma~\ref{hasc} shows that for all big enough $\rho$, 
$$v(\upg_{\rho}^{\dagger}-u_{\sigma}^\dagger)\ =\ 
\psi\big(v(\upg_{\rho})-v(u_\sigma)\big)\ =\ \psi\big({*-v(u_{\sigma})}\big)\ <\ v(\upg_{\rho+1}),$$
so $\upl+u_{\sigma}^\dagger\sim u_{\sigma}^\dagger-\upg_{\rho}^{\dagger}$, which proves our claim. 
%Thus by Corollary~\ref{corhasc},
%$$ \sigma\le \sigma'\ \Longrightarrow\ 
%v(\upl+u_{\sigma}^\dagger)\ \le\ 
%v(\upl+u_{\sigma'}^\dagger).$$
Next, given $\sigma$, take $\sigma'> \sigma$ such that 
$u_{\sigma'}\prec \upl+u_{\sigma}^\dagger$. By the claim above (with $\sigma'$ in place of $\sigma$) and Corollary~\ref{corhasc} we have 
$\upl + u_{\sigma'}^\dagger\prec u_{\sigma'}$,
%The argument above 
%(with $\sigma'$ in place of $\sigma$) shows that 
%for big enough $\rho$ we have 
%$\upl+u_{\sigma'}^\dagger \sim 
%u_{\sigma'}^\dagger-\upg_{\rho}^{\dagger}$. 
%Moreover, for
%big enough $\rho$ we have  
%$\upg_{\rho} \prec u_{\sigma'}$, so 
%$u_{\sigma'}^\dagger - \upg_{\rho}^{\dagger} \prec  %u_{\sigma'}$,   
and thus $\upl+u_{\sigma'}^\dagger \prec \upl+u_{\sigma}^\dagger$.

It now follows easily that some cofinal subsequence of 
$(-u_{\sigma}^\dagger)$ 
pseudoconverges to~$\upl$. Replacing $(u_{\sigma})$ by
a suitable cofinal 
subsequence we arrange $-u_{\sigma}^\dagger \leadsto \upl$. We show that then
$(-u_{\sigma}^\dagger)$ is equivalent to $(\upl_{\rho})$.
By the above, $(-u_{\sigma}^\dagger)$ and $(\upl_{\rho})$ have
a common pseudolimit $\upl$. Given $\rho$ we can take $\sigma$ such that
$u_{\sigma}\prec \upg_{\rho}$. Then $u_{\sigma}^\dagger-\upg_{\rho}^{\dagger}\prec \upg_{\rho}$,
which in view of $\upl+\upg_{\rho}^\dagger \prec \upg_{\rho}$ yields
$\upl + u_{\sigma}^\dagger\prec  \upg_{\rho}$. Thus the width of the pc-sequence~$(-u_{\sigma}^\dagger)$ 
equals the width of $(\upl_{\rho})$, and it remains only to apply 
Lemma~\ref{pc5}.
\end{proof}

\begin{cor}\label{cor:uplrho alternative}
Set $\upl_\rho^*:=-(\upl_{\rho+1}-\upl_\rho)^\dagger$ for $\rho<\kappa$. Then
$(\upl_\rho^*)$ is a pc-sequence in $K$ equivalent to $(\upl_\rho)$.
\end{cor}
\begin{proof}
Let $u_\rho:=\upl_{\rho+1}-\upl_\rho$. Then $\upl_\rho^*=-u_\rho^\dagger$, and $u_\rho\sim\upg_{\rho+1}$ by Lemma~\ref{pcz}.
For $\rho<\rho'<\kappa$ we obtain
$$\upl_{\rho'}^*-\upl_\rho^*\ 	=\ (u_\rho/u_{\rho'})^\dagger\
	\sim\ \left(\upg_{\rho+1}/\upg_{\rho'+1}\right)^\dagger\ 
		=\ \upl_{\rho'+1}-\upl_{\rho+1}\ \sim\ \upg_{\rho+2},
$$
hence $(\upl_\rho^*)$ is a pc-sequence in $K$.  Now apply Proposition~\ref{indcof} to $(u_{\rho})$.
\end{proof}

\begin{lemma}\label{lem:eval at uplrho}
Suppose $K$ has asymptotic integration, and $\upl_{\rho}\leadsto \upl$
with $\upl$ in an immediate asymptotic extension~$E$ of~$K$. Let
$G\in E\{Y\}\setminus E$. Then 
$G(\upl_\rho)\leadsto G(\upl)$.
\end{lemma}
\begin{proof} Use Corollary~\ref{cor:value group of dv(K)} to arrange that
$E$ is $\d$-valued of $H$-type, and then Proposition~\ref{asintint}
to get $E$ also closed under integration. 
For each $\rho$, take $y_\rho\in E$ with $y_\rho'=\upl_{\rho}- \upl$.
Then $y_{\rho}'\asymp\upg_{\rho+1}$ by
Lemma~\ref{pcz}. As $\upg_{\rho+1}=\ell_{\rho+1}^\dagger \asymp (\Log\ell_{\rho+1})'=\ell_{\rho+2}'$, and $\ell_{\rho+2}\succ 1$, this gives
$y_\rho \asymp \ell_{\rho+2}$ for all $\rho$.
Put 
$$P\ :=\ G(\upl+Y')-G(\upl)\in E\{Y\}.$$ Then $P\neq 0$, $P(0)=0$, and $P(y_\rho)=G(\upl_\rho)-G(\upl)$ for each $\rho$.  
Corollary~\ref{betterdifpol} gives $\beta\in\Gamma^<$ and a strictly increasing map $i\colon (\beta,0)\to\Gamma$ such that $vP(y)=i(vy)$ for all $y\in E$ with $vy\in (\beta,0)$. Thus $P(y_{\rho}) \leadsto 0$.
\end{proof}

\noindent
We now characterize pseudolimits of~$(\upl_\rho)$ in terms of active elements:

\begin{lemma}\label{lem:gap}
Let $\upl$ be an element of a valued differential field extension of $K$. Then the following conditions on $\upl$ are equivalent: 
\begin{enumerate}
\item[\textup{(i)}] $\upl_\rho \leadsto \upl$;
\item[\textup{(ii)}] for all active $a\in K$ there is an active $b\prec a$ in $K$
with $\upl+b^\dagger\prec a$;
\item[\textup{(iii)}] for all active $a\in K$ we have $\upl+a^\dagger\prec a$;
\item[\textup{(iv)}] for all $g\succ 1$ in $K$ we have $\upl + g^{\dagger\dagger}\prec g^\dagger$.
\end{enumerate}
\end{lemma}
\begin{proof}
Assume (i). Then for all $\rho$, 
$$\upl+\upg_{\rho}^\dagger\ =\ \upl - \upl_{\rho}\ \sim\ \upg_{\rho+1}\ 
\prec\ \upg_{\rho},$$
which gives (ii). To prove (ii)~$\Rightarrow$~(iii), assume (ii), and
suppose $a\in K$ is active. Take active $b\prec a$ in $K$ with 
$\upl+b^\dagger \prec a$. Then
$a^\dagger - b^\dagger \prec a$ by the remark preceding Proposition~\ref{indcof}, 
so
$\upl+a^\dagger\ =\ (\upl+b^\dagger) + (a^\dagger - b^\dagger)\prec a$, 
which gives (iii). The direction (iii)~$\Rightarrow$~(iv) is obvious.
For (iv)~$\Rightarrow$~(i), assume (iv). Let any $\rho$ be given,
and take $g\succ 1$ in $K$ such that for $a=g^\dagger$ we have
$a\prec \upg_{\rho +1}$ in $K$. Then $\upl+a^\dagger \prec a$ and
$\upg_{\rho+1}^\dagger - a^\dagger\prec \upg_{\rho+1}$. In view of
$$\upl+\upg_{\rho}^\dagger\ =\ (\upl+a^\dagger) + 
(\upg_{\rho+1}^\dagger - a^\dagger) + (\upg_{\rho}^\dagger-\upg_{\rho+1}^\dagger),$$
we get $\upl+\upg_{\rho}^\dagger \sim \upg_{\rho+1}$, so (i) holds.
\end{proof}

\begin{cor}\label{cor:gap}
Let $\upg$ be an element in an $H$-asymptotic field extension of $K$ with $\Psi<v\upg<(\Gamma^>)'$, and set $\upl:=-\upg^\dagger$. Then $\upl$ and $\upl+\upg$ are pseudolimits of $(\upl_\rho)$. 
\end{cor}
\begin{proof}
Let $a\in K$ be active. Then $\upl+a^\dagger = (a/\upg)^\dagger \prec a$
by Lemma~\ref{hasc} and the remark following its proof. Hence 
$(\upl+\upg)+a^\dagger = \upg+(a/\upg)^\dagger \prec a$. Thus $\upl_\rho\leadsto \upl$ and $\upl_\rho\leadsto \upl+\upg$ by the equivalence of (i) and (iii) in Lemma~\ref{lem:gap}. 
\end{proof}

\noindent
We now discuss the effect of compositional conjugation on the above. Let
$\phi\in K^\times$, and consider the compositional conjugate
$K^\phi$; its derivation is $\derdelta = \phi^{-1}\der$. 
The logarithmic sequence $(\ell_{\rho})$
for $K$ was obtained via a map 
$\Log \colon K^{\succ 1} \to K^{\succ 1}$. 
The dependence of $\Log$ on the contraction map $\chi$ 
(which is invariant under compositional 
conjugation) shows that we can use the same map $\Log$ to obtain the
same logarithmic sequence $(\ell_{\rho})$
for the compositional conjugate 
$K^\phi$. Let $(\upg_{\rho}^\phi):=   (\derdelta\ell_{\rho}/\ell_{\rho})$ be 
the corresponding sequence of active elements in $K^{\phi}$, 
so $\upg_{\rho}^\phi=\upg_{\rho}/\phi$ and
$$\upl_\rho^\phi\ :=\ -\frac{\derdelta \upg_{\rho}^\phi}{\upg_{\rho}^\phi}\ =\  -\left(\frac{\upg_\rho^\dagger}{\phi} - \frac{\phi^\dagger}{\phi} \right)\  =\ \frac{\upl_\rho}{\phi}+ \frac{\phi^\dagger}{\phi}.$$
Thus, given $\upl\in K$ we have 
$$\upl_\rho \leadsto \upl\ \Longleftrightarrow\ 
  \upl_\rho^\phi \leadsto (\upl/\phi)+(\phi^\dagger/\phi).$$

\begin{cor} \label{cor:lambda under comp conj}
The pc-sequence $(\upl_{\rho})$ has a pseudolimit in $K$ 
if and only if the corresponding pc-sequence
$(\upl_\rho^\phi)$  has a pseudolimit in $K^\phi$.
\end{cor}

\noindent
In the next subsection we relate the above to gap creation.

\subsection*{Gap creation} Recall from Section~\ref{AbstractAsymptoticCouples} that 
a gap in $K$ is an element $\beta$ of its
value group $\Gamma$ such that 
$\Psi < \beta < (\Gamma^{>})'$. Thus by our standing assumption and 
Corollary~\ref{trich}, 
%Since $\Psi$ has no largest element by our standing assumption,
either $K$ has a gap, or $K$ has asymptotic integration.  

\begin{lemma}\label{lem:gap active}
Suppose $K$ has a gap
$vf$ where $f\in K^\times$, and let $s:=f^\dagger$. Then for all active $a\in K$ we have $s-a^\dagger\prec a$. 
\end{lemma}
\begin{proof}
Let $a\in K$ be active. Then $va < vf < (\Gamma^{>})'$, so by Lemma~\ref{BasicProperties}(i),
\equationqed{v(s-a^\dagger)\ =\ v(f^\dagger-a^\dagger)\ =\ v\big((f/a)^\dagger\big)\ 
=\ \psi(vf-va)\ >\ va.}
\end{proof}

\noindent
Thus by Lemmas~\ref{lem:gap} and~\ref{lem:gap active}, if $K$ has a gap $vf$, $f\in K^\times$, then $\upl_{\rho} \leadsto -f^\dagger$.

\medskip
\noindent
Let $s\in K$. We say that $s$ {\bf creates a gap over $K$} if some
exponential integral of~$s$ introduces a gap,
that is, $vf$ is a gap
in $K(f)$, for some element $f\neq 0$ in some
$H$-asymptotic field extension of $K$ with $f^\dagger=s$.  

\index{element!creating a gap}
\index{gap!creating}

\begin{lemma}\label{CofinalLemma}
Suppose that $K$ has asymptotic integration and that $s\in K$ creates a gap over $K$. Then $s-a^\dagger$ is active, for all $a\in K^\times$.
\end{lemma}
\begin{proof} Take $f$ as in the definition above, and let 
$a\in K^{\times}$. Then $f \nasymp a$, so
$v(s-a^\dagger)=v\big((f/a)^\dagger\big)< (\Gamma^{>})'$, and thus $v(s-a^\dagger)\in \Psi^{\downarrow}$. 
\end{proof}

\begin{lemma}\label{lem:gap and compconj}
Suppose $s\in K$ creates a gap over $K$, and let $\phi\in K^\times$. Then 
$(s/\phi)-(\phi^\dagger/\phi)\in K^\phi$ creates a gap over $K^\phi$. 
\end{lemma}
\begin{proof}
Take some  $H$-asymptotic field extension $L=K(f)$ of $K$ with $f\ne 0$,
$f^\dagger = s$, and
$\Psi_L < vf < (\Gamma_L^{>})'$.
Then with $\derdelta=\phi^{-1}\der$
we have, in $L^\phi$, 
%$\derdelta(y/f)/(y/f)\ =\ (s/f)-(f'/f^2)\in K$, and
$$\Psi^\phi_L < v(f/\phi) < (\operatorname{id} +\psi^\phi_L)(\Gamma_L^{>}), \qquad
\derdelta(f/\phi)/(f/\phi)\ =\ (s/\phi)-(\phi^\dagger/\phi),$$
so $(s/\phi)-(\phi^\dagger/\phi)\in K^\phi$ creates a gap over $K^\phi$.    
\end{proof}

\begin{lemma}\label{cg1} Suppose 
$K$ has asymptotic integration and
$s\in K$ creates a gap over~$K$. Then $\upl_\rho\leadsto -s$.
Also, $y'+sy\not\asymp 1$ for all $y\in K$.
\end{lemma}
\begin{proof} Take an $H$-asymptotic field extension $L=K(f)$ of $K$ with gap $vf$,
$f\ne 0$, $f^\dagger=s$.  
Then $s-a^\dagger\prec a$ for all active $a\in L$, by Lemma~\ref{lem:gap active}; so
$s-a^\dagger\prec a$ for all active $a\in K$, and thus $\upl_{\rho}\leadsto -s$ by Lemma~\ref{lem:gap}.

Next, assume towards a contradiction that $y\in K$ satisfies
$y'+sy\asymp 1$. Then $y\ne 0$ and $y'+sy=y(y^\dagger + s)$, so for $a=y^{-1}$ we have $s-a^\dagger \asymp a$. Then $a$ is active by Lemma~\ref{CofinalLemma}, so
$s-a^\dagger \prec a$ by the first part of the proof, and we have a 
contradiction. 
\end{proof}
 
\noindent
Towards a partial converse to the first part of Lemma~\ref{cg1} we have:

\begin{lemma}\label{lem:active S}
Suppose $K$ has asymptotic integration, and $\upl_\rho\leadsto\upl\in K$. 
Set
$$S\ :=\ \big\{v(\upl+y^\dagger):\ y\in K^\times\big\}\ \subseteq\ \Gamma_\infty.$$ 
Then $S$ is a cofinal subset of $\Psi^{\downarrow}$, and thus $\upl\ne 0$. 
Moreover, if $f$ is a nonzero element of some $H$-asymptotic field extension of $K$ and $f^\dagger=-\upl$, then 
$vf\notin\Gamma$.
\end{lemma}
\begin{proof}
We first show that $\upl+y^\dagger$ is active, for each 
$y\in K^\times$. To prove this, assume towards a contradiction
that $y\in K^\times$ and $\upl+y^\dagger$ is not active. Set $\alpha:= vy$ and
take $\beta\in \Gamma^{\ne}$ with $\alpha=\beta'$. Take an active 
$a\not\asymp y$ with $va\ge \beta^\dagger$. Then $\upl+a^\dagger \prec a$, so
$$\psi(va-\alpha)\ =\ v(a^\dagger-y^\dagger)\ =\ v(a^\dagger-y^\dagger+\upl+y^\dagger)\ =\ v(\upl+a^\dagger)\ >\ va,$$
contradicting Lemma~\ref{lem:phipsi}, for $\alpha=vy$ and $\gamma=va$.
Thus $S\subseteq \Psi^{\downarrow}$, and $S$ is cofinal in~$\Psi^{\downarrow}$ by~(iii) of 
Lemma~\ref{lem:gap}. Now let $f\ne 0$ in some $H$-asymptotic field extension
of~$K$ satisfy $f^\dagger=-\upl$, and suppose towards a contradiction that
$vf\in \Gamma$. Then $f=uy$ with $u\in K(f)$, $u\asymp 1$ and $y\in K^\times$, so
$$v(\upl + y^\dagger)\ =\ v(f^\dagger - y^\dagger)\ =\ v(u^\dagger)\ =\ v(u') > \Psi$$
where the last inequality follows from Proposition~\ref{CharacterizationAsymptoticFields}. This is in contradiction to~${S\subseteq \Psi^{\downarrow}}$.
\end{proof} 

\noindent
Here is a partial converse to the first part of Lemma~\ref{cg1}:

\begin{lemma}\label{cg2} Assume that $K$ has asymptotic integration, that $\Gamma$ is divisible, and that $\upl_\rho\leadsto\upl\in K$. Then 
$s=-\upl$ creates a gap over $K$. 
\end{lemma}
\begin{proof} 
Lemmas~\ref{newprop, lemma 1} and \ref{lem:active S} give an $H$-asymptotic field extension
$L=K(f)$ with $vf\notin \Gamma$, $f^\dagger =s$, and $\Psi_L=\Psi$. 
We show that
$\alpha:=vf$ is a gap in $L$. Suppose it is not. Since $\Psi_L=\Psi$ has no maximum, Lemma~\ref{have1} and its Corollary~\ref{trich} then 
give a nonzero $\beta\in \Gamma_L$ such that
$\alpha=\beta'$. Then $\beta^\dagger\in \Psi$. Take active $a$ in $K$
with $a\not\asymp f$ and $va\ge \beta^\dagger$. Then $s-a^\dagger \prec a$, hence
$$\psi_L(\alpha -va)\ =\ v(f^\dagger-a^\dagger)\ =\ v(s-a^\dagger)\ >\ va,$$
which contradicts Lemma~\ref{lem:phipsi}, for $\gamma=va$.  
\end{proof}

\begin{remark}
Suppose that $K$ has asymptotic integration, $\Gamma$ is divisible, and
$s\in K$ creates a gap over $K$.
Then $vf$ is a gap in $K(f)$ for {\em any}\/ nonzero element~$f$ of any $H$-asymptotic field extension of $K$ with $f^\dagger = s$.
(To see this, note that by Lemmas~\ref{cg1} and~\ref{lem:active S}, any such $f$ satisfies $vf\notin\Gamma$ and hence is transcendental over $K$;
it remains to use the uniqueness part of Lemma~\ref{newprop, lemma 1}.)
\end{remark}

%\begin{cor}\label{cor:extend to smallest comp class}
%Every $H$-asymptotic field has a grounded $H$-asymptotic field extension.
%\end{cor}
%\begin{proof}
%Let $E$ be an ungrounded $H$-asymptotic field.
%Then $E$ has a gap or $E$ has asymptotic integration. 
%If $E$ has a gap, then Lemma~\ref{extensionpdvfields1} yields 
%a grounded $H$-asymptotic field extension $E(y)$ of $E$. In general we can 
%pass to an algebraic extension and arrange that $\Gamma_E$ is divisible 
%(Proposition~\ref{Algebraic-Extensions-Asymptotic}).
%For such $E$ it
%remains to consider the case that~$E$ has asymptotic integration. 
%By Corollary~\ref{zdcor}
%we can pass from $E$ to an immediate extension and arrange
%that $E$ is also spherically complete. Then some element of $E$ 
%creates a gap over $E$, by Lemma~\ref{cg2}, and therefore $E$ 
%has an $H$-asymptotic extension $E(f)$ with a gap. The case of a 
%gap has been treated in the beginning of the proof. 
%\end{proof}

\begin{cor}\label{cor:extend to smallest comp class}
Every pre-$\d$-valued field of $H$-type has a grounded $\d$-valued extension
of $H$-type. Every pre-$H$-field has a grounded $H$-field extension.
\end{cor}
\begin{proof}
Let $E$ be a pre-$\d$-valued field of $H$-type.  
%$\Gamma_{\dv(E)} \ne \Gamma_E$, then $\dv(E)$ is grounded by 
%Corollary~\ref{cor:value group of dv(K)}, and $\dv(E)$ is an extension 
%as claimed. Suppose 
%$\Gamma_{\dv(E)} = \Gamma_E$ (so $\dv(E)$ is an immediate extension of $E$).
Replacing $E$ by $\operatorname{dv}(E)$ we assume below that $E$ is $\d$-valued.   
If $E$ is grounded, we are done, so we assume that
$E$ has a gap or has asymptotic integration. 
If $E$ has a gap, then Lemma~\ref{extensionpdvfields1} yields a grounded 
$\d$-valued extension $E(y)$ of $H$-type of $E$. In general we can arrange 
(for example by passing to the algebraic closure using
Corollary~\ref{cor:Algebraic-Extensions-pdv}) that 
$\Gamma_E$ is divisible.
For such $E$ it
remains to consider the case that~$E$ has asymptotic integration. 
In that case Corollary~\ref{zdcor} allows us
to pass from $E$ to an immediate extension and arrange
that $E$ is also spherically complete. Then some element of $E$ creates 
a gap over $E$ by Lemma~\ref{cg2}, and so $E$ has a $\d$-valued 
extension $E(f)$ of $H$-type with a gap, by the remark preceding 
this corollary, and Lemma~\ref{newprop, lemma 1}. The case of a gap 
has been treated earlier in the proof. 

For a pre-$H$-field $E$, replace $\operatorname{dv}(E)$ in the proof 
above by $H(E)$, and appeal to ~\ref{extensionpdvfields1, cor}, ~\ref{cor:H-field rc}, 
~\ref{prehexpint} instead of ~\ref{extensionpdvfields1},
~\ref{cor:Algebraic-Extensions-pdv}, ~\ref{newprop, lemma 1}.
\end{proof}

%\begin{remark} Likewise, every pre-$\d$-valued field of $H$-type 
%has a grounded $\d$-valued field extension of $H$-type, and every 
%pre-$H$-field has a grounded
%$H$-field extension.
%This follows by the arguments in the proof of
%Corollary~\ref{cor:extend to smallest comp class},
%also using the remarks after the proofs of Lemma~\ref{extensionpdvfields1} 
%and Proposition~\ref{Algebraic-Extensions-Asymptotic}, and taking 
%into account much of the above subsection, as well as 
%Lemmas~\ref{newprop, lemma 1} and ~\ref{prehexpint}.
%\end{remark}

\subsection*{Notes and comments}
The $\upg_n$ occur in classical logarithmic
criteria for convergence of infinite series~\cite{AbelOeuvre,Bertrand,deMorgan} and integrals~\cite[p.~229]{Riemann}; see~\cite[\S\S{}4.1,~3.2]{Bourbaki}. Attempts to sharpen these logarithmic criteria provided the initial impetus for du~Bois-Reymond's  ``orders of infinity''~\cite{dBR73}; see~\cite[Section~1]{Fisher}.
In this connection, the role of gaps as an ``ideal boundary between convergence and divergence'' was anticipated in \cite[\S\S{}10--17]{dBR77}.
The~$\upo_n$ appear in \cite{Hartman,Hille48} (investigating the nature of solutions $y$ to ${4y'' + fy =0}$),  and implicitly already in \cite[\S{}25]{Weber}; see
\cite[Section~17]{Boshernitzan} and~\cite{Rosenlicht8} for the setting of Hardy fields.

\section{The Property of $\upl$-Freeness}\label{sec:uplfree}

\noindent
The ultimate goal of our work is
to analyze the elementary (=~first-order) theory of~$\mathbb{T}$. Accordingly we introduce and study
in this section the elementary property of being $\upl$-free and show
that $\mathbb{T}$ is $\upl$-free. {\em We keep the standing assumption from the
previous section that $K$ is an ungrounded $H$-asymptotic field with 
$\Gamma\ne \{0\}$}. We begin with a 
consequence of the material in
the previous section.

\begin{samepage}
\begin{cor}\label{psacgap} The following four conditions on $K$ are equivalent:
\begin{enumerate}
\item[\textup{(i)}] $(\upl_\rho)$ has no pseudolimit in $K$;
\item[\textup{(ii)}] for all $s$ there is $g\succ 1$ such that $s-g^{\dagger\dagger}\succeq g^\dagger$;
\item[\textup{(iii)}] for all $s$ there is an active $a$ such that $s-a^\dagger \succeq a$;
\item[\textup{(iv)}] for all $s$ there is an active $a$ such that $s-b^\dagger \succeq a$ for all active $b\prec a$.
\end{enumerate}
Here $a$,~$b$,~$g$,~$s$ range over $K$. Also, condition \textup{(i)} implies condition \textup{(v)} given by
\begin{enumerate}
\item[\textup{(v)}]
$K$ has asymptotic integration, and no element of~$K$ creates a gap over $K$.
\end{enumerate}
Finally, if $\Gamma$ is divisible, then \textup{(i)} is equivalent to \textup{(v)}.
\end{cor}
\end{samepage}

\begin{proof} The equivalence of (i), (ii), (iii), (iv) follows from Lemma~\ref{lem:gap}. For the im\-pli\-ca\-tion~(i)~$\Rightarrow$~(v), use the remark following Lemma~\ref{lem:gap active}, and Lemma~\ref{cg1}. For the last claim, also use Lemma~\ref{cg2}.
\end{proof}

\noindent
Condition (ii) in Corollary~\ref{psacgap} is a first-order condition on the asymptotic field $K$. 
Conditions (iii) and (iv) are also first-order conditions on $K$, but are logically more complex.
We say that $K$ is {\bf $\upl$-free} if it satisfies condition (ii)
of Corollary~\ref{psacgap}. Thus
if $K$ is $\upl$-free, then $K$ has asymptotic integration, and no element
of $K$ creates a gap over $K$. From 
Lemma~\ref{cg1} we conclude:

\label{p:upl-free}
\index{lambda-free@$\upl$-free}
\index{H-asymptotic@$H$-asymptotic!field!$\upl$-free}

\begin{cor} \label{cor:1-ls and divisible value group implies upl-free}
If $K$ has asymptotic integration and is $1$-linearly surjective, and~$\Gamma$ is divisible, then $K$ is $\upl$-free.
\end{cor}

%\begin{examples}
\noindent
Every Liouville closed $H$-field is $\upl$-free, by Corollary~\ref{cor:1-ls and divisible value group implies upl-free}.
More generally, if $K$ has asymptotic integration and $(K^\times)^\dagger=K$, then $K$ is $\upl$-free, by Lemma~\ref{lem:active S}. By Proposition~\ref{vsur3} below we can drop ``$\Gamma$ is divisible''  in Corollary~\ref{cor:1-ls and divisible value group implies upl-free}.
%\end{examples}
The following also helps in identifying $\upl$-free $K$.

\begin{lemma}\label{conwam} Let $s\in K$ lie in a differential subfield
$F$ of $K$ for which $\max\Psi_F$ exists. Then there is an active $a\in K$ such that $s-a^\dagger \succeq a$.
\end{lemma}
\begin{proof} Take $a\in F$ such that $va=\max \Psi_F$. Then $a\in K$ is 
active, so if
$s-a^\dagger \succeq a$, we are done. Suppose
$s-a^\dagger \prec a$. Then $v(s-a^\dagger)\in (\Gamma_F^{>})'\cup\{\infty\}$, so 
$v(s-a^\dagger) > \Psi$.  
Take $b\succ 1$ in $K$ such that $vb'=va$. Then $vb+ vb^\dagger=va$, 
so $v(a/b)=vb^\dagger$, and thus $a/b\in K$ is active, and 
$s-(a/b)^\dagger\ =\ (s-a^\dagger) + b^\dagger\ \sim\ b^\dagger\ \asymp\ a/b$.
\end{proof}

\begin{cor}\label{recwam} If $K$ is a union of grounded $H$-asymptotic 
subfields, then $K$ is $\upl$-free.
\end{cor} 

\noindent
If $K$ has a $\upl$-free $H$-asymptotic field extension $L$ such that $\Gamma^<$ is cofinal in $\Gamma_L^<$, then~$K$ is $\upl$-free. 
In particular,  if $K$ has an immediate $\upl$-free $H$-asymptotic field extension, then $K$ is $\upl$-free. In the other direction we have:

\begin{lemma}\label{lem:lambda-free completion} If $K$ is $\upl$-free, then so is its completion $K^{\operatorname{c}}$.
\end{lemma}
\begin{proof} If $\upl_{\rho} \leadsto \upl\in K^{\operatorname{c}}$,
then $\upl_{\rho} \leadsto a$ for any $a\in K$ with $v(a-\upl)> \Psi$. 
\end{proof}

\noindent
In connection with the next result, see Proposition~\ref{asfl}.

\begin{lemma}\label{lem:lambda-free fluentcompletion} If $K$ is $\upl$-free, then so is any asymptotic extension which, as a valued field extension of $K$, is
a fluent completion of $K$.
\end{lemma}
\begin{proof} Immediate from 
Lemma~\ref{lem:no jammed pseudolimits}. 
\end{proof}

\begin{lemma}\label{lem:lambda-free henselization}
If $K$ is $\upl$-free, then so is its henselization $K^{\operatorname{h}}$. 
\end{lemma}
\begin{proof}
If $(\upl_\rho)$ pseudoconverges in $K^{\operatorname{h}}$, then
it pseudoconverges in some fluent completion of $K$, and thus in $K$,
by Lemma~\ref{lem:no jammed pseudolimits}.
\end{proof}

\begin{cor}\label{cor:lambda-free alg closure} If $K$ is $\upl$-free, then
$K$ has rational asymptotic integration. Also:
$$  \text{$K$ is $\upl$-free}\ \Longleftrightarrow\  \text{the algebraic closure $K^\alg$ of $K$ is $\upl$-free}.$$
\end{cor}
\begin{proof} The first claim is a consequence of the subsequent equivalence, so we prove the equivalence.
Suppose $(\upl_\rho)$ pseudoconverges in $K^{\alg}$. Then by Lemma~\ref{lem:Kaplansky, 1}, the pc-sequence $(\upl_\rho)$ is of algebraic type over $K$, and
so pseudoconverges in $K^{\operatorname{h}}$, by Corollaries~\ref{cor:algebraically maximal} and \ref{cor:alg max equals henselian}, and thus in $K$ by Lemma~\ref{lem:lambda-free henselization}.
This gives the forward direction of the equivalence, and
the backward direction holds because $\Gamma^<$ is cofinal 
in~$(\Q\Gamma)^<$. 
\end{proof}

\begin{lemma}\label{lem:s+lambda active}
Suppose $K$ is $\upl$-free, $\upl$ is a pseudolimit of $(\upl_\rho)$ in an $H$-asymptotic field extension of $K$, and $s\in K$. Then $s+\upl \sim b$ for some active $b\in K$.
\end{lemma}
\begin{proof}
Take active $a\in K$ as in (iv) of  Corollary~\ref{psacgap}, that is, $s-b^\dagger \succeq a$ for all active $b\prec a$ in $K$.
For such $b$ we have $s-a^\dagger=s-b^\dagger+(b/a)^\dagger\sim s-b^\dagger$, since $(b/a)^\dagger\prec a\preceq s-b^\dagger$. Take $\upg_{\rho_0}\prec a$. Then $s+\upl_\rho \sim s+\upl_{\rho_0}\succ \upg_{\rho_0}$ for $\rho\ge\rho_0$, so $s+\upl_\rho$ is active for $\rho\ge\rho_0$.
Now $\upl-\upl_\rho \sim \upg_{\rho+1}$ for all $\rho$, and thus $\upl-\upl_\rho \prec \upg_{\rho_0} \prec s+\upl_\rho$ for~$\rho\ge\rho_0$.
Hence $s+\upl=(s+\upl_\rho) + (\upl-\upl_\rho) \sim s+\upl_{\rho}$ for $\rho\ge \rho_0$.
\end{proof}

\noindent
The following technical consequence of Corollary~\ref{cor:lambda-free alg closure} and Lemma~\ref{lem:s+lambda active} will be used in Section~\ref{fcuplfr}. 
Recall the derivation $R\mapsto \partial R/\partial Y$ on $K(Y)$ from Section~\ref{Differential Fields and Differential Polynomials}. 

\begin{cor}\label{cor:val of formal log der at upl}
Suppose $K$ is $\upl$-free, and $\upl$ is a pseudolimit of $(\upl_\rho)$ in an $H$-asymptotic field extension of $K$. Then $\upl$ is transcendental over $K$, and for each $R\in K(Y)^\times$ there is an active $a\in K$ such
that
$$\left(\frac{\partial R/\partial Y}{R}\right)(\upl)\ \preceq\ \frac{1}{a}.$$
\end{cor}
\begin{proof} Take an algebraically closed $H$-asymptotic extension
of $K$ containing $\upl$, and replace $K$ by its algebraic closure
in this extension. It remains to appeal to Corollary~\ref{cor:lambda-free alg closure}
and to use Lemmas~\ref{lem:log der of rational fn} and \ref{lem:s+lambda active}.
\end{proof}

\noindent
By Corollary~\ref{cor:lambda under comp conj}, if $K$ is $\upl$-free and $\phi\in K^\times$, then the compositional conjugate~$K^\phi$ is $\upl$-free. 
Being $\upl$-free is also preserved under suitable specializations:

\begin{lemma}\label{specuplfree}
Suppose $K$ is $\upl$-free with small derivation. Let $\Delta\neq\{0\}$ be a convex subgroup of $\Gamma$ with $\psi(\Delta^{\ne})\subseteq\Delta$. Then the asymptotic residue field $\dot K$ of the coarsened valuation $\dot v=v_\Delta$ is of $H$-type, $\psi(\Delta^{\ne})\ne \emptyset$, $\psi(\Delta^{\neq})$ is cofinal in $\Psi$ \textup{(}and so has no largest element\textup{)} and $\dot K$ is $\upl$-free.
\end{lemma}
\begin{proof} It is clear that $\dot{K}$ is $H$-asymptotic and that
the claims about $\psi(\Delta^{\ne})$ hold.
As usual, let $\dot{\mathcal O}=\{a\in K:\text{$va \geq \delta$ for some $\delta\in\Delta$}\}$ denote the valuation ring of $\dot v$, with maximal ideal $\dot{\smallo}=\{a\in K: va>\Delta\}$ and residue map $a\mapsto\dot a\colon\dot{\mathcal O}\to \dot K=\dot{\mathcal O}/\dot{\smallo}$.
Let $s\in \dot{\mathcal O}$; it is enough to get $a\in K^\times$ such that $va\in \Delta^{\ne}$, $a$ is active in~$K$, and $s-a^\dagger\succeq a$. 
By Corollary~\ref{psacgap} and since $\Psi^{>0}\ne \emptyset$, we can take active $b\in K$ such that $vb>0$ and $s-b^\dagger\succeq b$.
Then $vb\in \Delta$ by Lemma~\ref{asymp-lemma1}, and so $a:=b$ works.
\end{proof}

\noindent
Being $\upl$-free is preserved under adjunction of certain exponential integrals:

\begin{lemma}\label{lem:lambda-free, exp int, 1}
Suppose $K$ is $\upl$-free and $\Gamma$ is divisible. Let $s\in K$ be such that $s-a^\dagger$ is active for each $a\in K^\times$, and set $S:=\big\{v(s-a^\dagger):a\in K^\times\big\}$. Let $f^\dagger=s$, where $f\ne 0$ lies in an $H$-asymptotic field extension of~$K$.
Suppose that
\begin{enumerate}
\item[\textup{(i)}] $S$ does not have a largest element,  or
\item[\textup{(ii)}] $S$ has a largest element and $[\gamma + vf]\notin [\Gamma]$ for some $\gamma\in\Gamma$.
\end{enumerate}
Then $K(f)$ is $\upl$-free.
\end{lemma}

\begin{proof}
Let $L$ be the algebraic closure of $K$ in an algebraic closure of $K(f)$.
Then~$L$ is an $H$-asymptotic field extension of $K$ with $\Gamma_L=\Gamma$. We claim that
$$S\ =\ \big\{v(s-b^\dagger):\ b\in L^\times\big\}.$$
To see this, let $b\in L^\times$, and take $a\in K^\times$, $u\asymp 1$ with $b=au$; then 
$$v(u^\dagger)\ =\ v(u')\ >\ v(s-a^\dagger)\ =\ v(s-a^\dagger-u^\dagger)\ =\ v(s-b^\dagger)$$ 
since $s-a^\dagger$ is active in $L$. This proves the claim. In particular, $s-b^\dagger$ is active in~$L$, for each $b\in L^\times$. Since $L$ is $\upl$-free by Corollary~\ref{cor:lambda-free alg closure}, we have a reduction to the case that 
$K$ is algebraically closed and~$f$ is transcendental over~$K$.

 If $S$ satisfies condition (i), then
$\Psi_{K(f)}=\Psi$ by Lemma~\ref{newprop, lemma 1}. 
If $S$ satisfies condition (ii), then
$\Psi_{K(f)}=\Psi\cup \big\{v(s-a^\dagger)\big\}$ for some $a\in K^\times$,
by Lemma~\ref{newprop, lemma 2}(i) and its proof.  
In each case, $vf\notin \Gamma$ and
$\Psi$ is a cofinal subset of $\Psi_{K(f)}$, so our logarithmic sequence $(\ell_\rho)$ for~$K$ remains a logarithmic sequence for~$K(f)$. 
Thus  $(\upl_\rho)$ has no pseudolimit in $K(f)$ by Corollary~\ref{cor:pc-limits in simple trans exts}.
\end{proof}

\begin{lemma}\label{lem:lambda-free, exp int, 2}
Suppose $K$ is $\upl$-free and is equipped with an ordering making it a real closed $H$-field. Let $s\in K^{<}$ be such that
$s-a^\dagger$ is active, for all $a\in K^\times$. Let~$f$ be a nonzero element in a pre-$H$-field extension of $K$ such that $f^\dagger=s$.  Then $K(f)$ is $\upl$-free.
\end{lemma}
\begin{proof}
This is clear if $f\in K$, so suppose $f\notin K$; then $f$ is transcendental over~$K$. Lemma~\ref{lemma62} gives $vf\notin \Gamma$ and
$\Psi$ is cofinal in $\Psi_{K(f)}$, so $\Psi_{K(f)}$ does not have a largest element, and
$\Gamma^<$ is cofinal in $\Gamma_{K(f)}^<$.  Hence the 
logarithmic sequence $(\ell_\rho)$ for~$K$
remains a logarithmic sequence for $K(f)$. 
Now use Corollary~\ref{cor:pc-limits in simple trans exts}.
\end{proof}

\subsection*{Application to eventual equalizing} {\em In this subsection
we assume that $K$ has asymptotic integration}.
For use in Section~\ref{cordone} we prove here Proposition~\ref{vsur3}. 

\index{linear differential operator!$v$-surjective}
\index{v-surjective@$v$-surjective}

We say that $A\in K[\der]^{\ne}$ is {\bf $v$-surjective}
if for each $g\in K^\times$ there is an $f\in K^\times$ such that 
$A(f) \asymp g$. For example, $\der$ is $v$-surjective since
$K$ has asymptotic integration. 
If $A\in K[\der]^{\ne}$ is $v$-surjective and $a\in K^\times$, then so 
are $aA$ and $Aa$ for any
$a\in K^\times$. If $A\in K[\der]^{\ne}$ has order $1$, then
$aA=\der-s$ for some $a\in K^\times$ and $s\in K$. Note that
if $A=\der-s$ with $s\in K$, then $a^{-1}Aa=\der-(s-a^{\dagger})$ for $a\in K^\times$.   

\begin{lemma}\label{vsur1} Let $s\in K$ and $A=\der -s$. Then $A$ is $v$-surjective
if $v(s-a^\dagger)> \Psi$ for some $a\in K^\times$, 
or $\big\{v(s-a^\dagger):\ a\in K^\times\big\}$ is a subset of 
$\Psi^{\downarrow}$ with a largest element.
\end{lemma}
\begin{proof} Consider first the case that $v(s-a^\dagger)> \Psi$ with 
$a\in K^\times$. Renaming $a^{-1}Aa$ and $s-a^\dagger$ as $A$ and $s$
we have 
$vs> \Psi$. Let $g\in K^\times$ and take $f\in K^\times$ with $f\not\asymp 1$ and 
$f'\asymp g$.
Then $A(f)=f'-sf$, and $f^\dagger \succ s$ gives $f'\succ sf$, so 
$A(f)\sim f'\asymp g$. Next, assume $\big\{v(s-a^\dagger):\ a\in K^\times\big\}$ is a subset of 
$\Psi^{\downarrow}$ with a largest element. Take $a\in K^\times$ such that
$v(s-a^\dagger)$ is maximal. Renaming $a^{-1}Aa$ and $s-a^\dagger$ as $A$ and $s$ 
we have $vs\in \Psi^{\downarrow}$ and $s-b^\dagger\succeq s$
for all $b\in K^\times$. Let $g\in K^\times$ and take $\beta\in \Gamma^{\ne}$ 
such that
$\beta'=vg$. If $\beta^\dagger \le v(s)$, then any $f\in K^\times$ with 
$vf=\beta$ satisfies $f^\dagger - s \asymp f^\dagger$, so 
$A(f)=f(f^\dagger -s)\asymp f'\asymp g$. 
Suppose $\beta^\dagger > v(s)$.
Then we take $\alpha\in \Gamma$ with $\alpha+vs=vg$. Since 
$vg=\beta+\beta^\dagger > \beta+ vs$, this gives $\alpha > \beta$, so
$\alpha' > \beta'=vg=\alpha+vs$ (with $0':= \infty$ by convention).
Take $f\in K^\times$ with $vf=\alpha$. Then $A(f)=f'-sf\sim -sf\asymp g$, as desired.   
\end{proof}

\begin{lemma}\label{vsur2} Assume $\Gamma$ is divisible, $s\in K$, and
$\big\{v(s-a^\dagger):\ a\in K^\times\big\}$ is a subset of 
$\Psi^{\downarrow}$ without largest element. Then the following are equivalent for $g\in K^\times$:
\begin{enumerate}
\item[\textup{(i)}] $vg\notin v\big(A(K)\big)$ for $A:=\der -s$; 
\item[\textup{(ii)}] $g^\dagger-s$ creates a gap over $K$.
\end{enumerate}
\end{lemma}
\begin{proof} 
Lemma~\ref{newprop, lemma 1} yields an $H$-asymptotic extension $L=K(b)$ with
$b\ne 0$, $b^\dagger = s$, $\eta:= vb\notin \Gamma$, and
$\Gamma_L=\Gamma \oplus \Z \eta$, $[\Gamma_L]=[\Gamma]$ (and thus
$\Psi_L = \Psi$). Set $A:=\der - s$ and let $f\in K^\times$. Then $A(f)=f\cdot(f/b)^\dagger$, so 
$$v(A(f))-\eta\ =\ (vf-\eta) + \psi_L(vf-\eta)\ =\ (vf-\eta)'.$$ 
Also, $\Psi_L \subseteq \Gamma$ gives $(\Gamma-\eta)\cap (\Gamma_L^{\ne})'\subseteq (\Gamma-\eta)'$, and thus  
for all $\gamma\in \Gamma$,
$$\gamma\in v(A(K))\ \Longleftrightarrow\ \gamma-\eta\in (\Gamma_L^{\ne})'\ \Longleftrightarrow\ \text{$\gamma-\eta$ is not a gap in $L$.}$$
Now assume $\gamma:= vg\notin v\big(A(K)\big)$. Then $\gamma-\eta=v(g/b)$ is a gap
in $L$, and so $(g/b)^\dagger=g^\dagger - s$
creates a gap over $K$. This proves (i)~$\Rightarrow$~(ii). The converse is
part of Lemma~\ref{cg1}. 
\end{proof}

\begin{cor}
Suppose $\Gamma$ is divisible and $A\in K[\der]$ has order $1$. Then 
$$\big|\Gamma\setminus v\big(A(K)\big)\big|\leq 1.$$
\end{cor}
\begin{proof} We first reduce to the case that
$A=\der-s$, $s\in K$, and then set
$$S\ :=\ \big\{v(s-a^\dagger):\ a\in K^\times\big\}.$$
If $S\not\subseteq\Psi^\downarrow$ or $S$ is a subset of 
$\Psi^{\downarrow}$ with a largest element, then $A$ is $v$-surjective by Lemma~\ref{vsur1}.
Suppose $S$ is a subset of 
$\Psi^{\downarrow}$ without a largest element, and $g,h\in K^\times$, $vg,vh\notin v\big(A(K)\big)$.
Then $g^\dagger-s$ and $h^\dagger-s$ create a gap over $K$ by Lemma~\ref{vsur2}, and so $v\big((g/h)^\dagger\big)=v(g^\dagger-h^\dagger)>\Psi$ by
Lemmas~\ref{pcz} and \ref{cg1}, thus $vg=vh$.
\end{proof}

\begin{example}
If $\Gamma$ is divisible and  $\upl_\rho\leadsto\upl\in K$, then $v\big(A(K)\big)=\Gamma^{\neq}$ for $A= \der-\upl$.
\end{example}

\noindent
We let $\phi$ range over the active elements in $K$. 
Let $P=aY+bY'$ with $a,b\in K$ not both zero. We say that {\bf $P$ has eventual
equalizers} if for every $g\in K^\times$ there is $f\in K^\times$ such that
eventually $P_{\times f}^{\phi} \asymp g$. Note that if $P$ has eventual equalizers
and $h\in K^\times$, then so do $hP$ and $P_{\times h}$. The next result improves on Corollary~\ref{cor:1-ls and divisible value group implies upl-free}. 

\index{eventual!equalizer}
\index{equalizer!eventual}

\begin{prop}\label{vsur3} The following conditions on $K$ are equivalent:
\begin{enumerate}
\item[\textup{(i)}] $\der-s$ is $v$-surjective for every $s\in K$;
\item[\textup{(ii)}] every $P=aY + bY'$ with $a,b\in K$ not both zero has eventual equalizers;
\item[\textup{(iii)}] $K$ is $\upl$-free. 
\end{enumerate}
\end{prop}
\begin{proof} Assume (i) and let $P=aY + bY'$ with 
$a,b\in K$ not both zero; our job is to show that $P$ has eventual equalizers.
The case $b=0$ being trivial, assume $b\ne 0$. 
Replacing $P$ by $b^{-1}P$ we get $P=Y'-sY$ with
$s\in K$. Then for $f\in K^\times$,
$$ P_{\times f}^{\phi}\ =\ f\phi Y' + (f'-sf)Y.$$
Consider first the case that $v(s-h^\dagger)> \Psi$ for some $h\in K^\times$.
For such $h$ we have $h^{-1}P_{\times h}=Y'-(s-h^\dagger)Y$, so by suitable renaming
we arrange that $v(s)> \Psi$. Now, let $g\in K^\times$, and take $f\in K^\times$
with $f\not\asymp 1$ such that $f'\asymp g$. Then $f'-sf\sim f'\asymp g$, and
eventually $f\phi\prec f'$, so by the identity above for 
$P_{\times f}^{\phi}$ we get  $P_{\times f}^{\phi} \asymp g$, eventually.

Next, assume that $v(s-h^\dagger)\in \Psi^{\downarrow}$ for all $h\in K^\times$. 
Take $f\in K^\times$ such that
$f'-sf\asymp g$. Then $f'-sf=f(f^\dagger-s)$ and  $f^\dagger-s \succ \phi$, 
eventually, so $f'-sf\succ f\phi$ eventually, and thus  
$P_{\times f}^{\phi} \asymp g$, eventually. This finishes the proof of (i)~$\Rightarrow$~(ii).

\medskip\noindent
For  (ii)~$\Rightarrow$~(i), assume (ii), and let $s\in K$. Then $P:=Y'-sY$ 
has eventual equalizers. Let $g\in K^\times$ and take $f\in K^\times$ such that
$P_{\times f}^{\phi} \asymp g$, eventually. Then by the
identity above we must have $f'-sf\asymp g$.

\medskip\noindent
Next we prove the contrapositive of (i)~$\Rightarrow$~(iii). 
Suppose $K$ is not $\upl$-free.
Take $\upl\in K$ such that $\upl_{\rho}\leadsto \upl$.
If $y\in K^\times$ and $y'-\upl y \asymp 1$, then $a:= y^{-1}$ gives
$a^{\dagger} + \upl \asymp a$, so~$a$ is active by Lemma~\ref{lem:active S}, but
this contradicts the equivalence (i)~$\Longleftrightarrow$~(iii) of 
Lemma~\ref{lem:gap}.
Thus $\der-\upl$ is not $v$-surjective.

\medskip\noindent
Finally, we prove (iii)~$\Rightarrow$~(ii), so assume (iii). Let $P=aY + bY'$ with 
$a,b\in K$ not both zero. Now the algebraic closure
$K^\alg$ of $K$ is also $\upl$-free. Suppose $g\in K^\times$ and 
$f\in (K^\alg)^\times$ are such that eventually $P_{\times f}^{\phi} \asymp g$. 
By taking
$\phi$ (in $K$) with high enough~$v\phi$ and applying the Equalizer Theorem 
to $P^{\phi}\in K^{\phi}\{Y\}$ we see that $vf\in \Gamma$. 
Thus in trying to show that $P$ has eventual
equalizers, we can assume that~$\Gamma$ is divisible. But then the
$\upl$-freeness of $K$ and Lemmas~\ref{vsur1} and~\ref{vsur2} show that~(i) holds, and we already know that  (i)~$\Rightarrow$~(ii).  
\end{proof}

\subsection*{Notes and comments} Section~12 of \cite{AvdD3}
contains some aspects of $\upl$-freeness, but without naming this property. This occurs in connection with gap creation in real closed $H$-fields with asymptotic integration.   
Lemma~\ref{lem:lambda-free henselization} can also be deduced
from Lemma~\ref{specfluent} and Proposition~\ref{prop:approx henselization}, since~$(\upl_\rho)$ is jammed. Allen Gehret has shown us a proof that if $K$ is
$\upl$-free, then so is $\operatorname{dv}(K)$.

\section{Behavior of the Function $\omega$} \label{sec:behupo}

\noindent
{\em In this section we keep the assumption that $K$ is an ungrounded $H$-asymptotic field with $\Gamma\ne \{0\}$}.
Recall from Section~\ref{sec:secondorder} the function
$$\omega\colon K \to K,\qquad \omega(z)\ =\  -(2z'+z^2).$$

\begin{lemma}\label{pcrho}
For $\rho<\rho'<\kappa$ we have $
\omega(\upl_{\rho'})-\omega(\upl_\rho)\ \sim\ \upg_{\rho+1}^2$.
\end{lemma}
\begin{proof} Let $\rho < \rho' < \kappa$. Then 
$\upl_{\rho'}-\upl_\rho=\upg_{\rho+1}+\varepsilon$ with $\varepsilon \prec \upg_{\rho+1}$.
By Lemma~\ref{varrholemma, claim, 1}:
\begin{align*}
\omega(\upl_\rho)-\omega(\upl_{\rho'}) &=\ (\upg_{\rho+1}+\varepsilon)\cdot\bigl(2(\upg_{\rho+1}+\varepsilon)^\dagger+2\upl_\rho+(\upg_{\rho+1}+\varepsilon)\bigr)\\
&=\ (\upg_{\rho+1}+\varepsilon)\cdot \bigl(2\upg_{\rho+1}^\dagger + 2(1+\tau)^\dagger
+2\upl_\rho+(\upg_{\rho+1}+\varepsilon)\bigr),
\end{align*}
where  $\tau:= \varepsilon/\upg_{\rho+1}$. Thus $\tau\prec 1$ and so 
$(1+\tau)^\dagger  \prec \upg_{\rho+1}$, hence 
\begin{align*}
\omega(\upl_\rho)-\omega(\upl_{\rho'})\ &=\ (\upg_{\rho+1}+\varepsilon)\cdot
\bigl(-2\upl_{\rho+1} +2(1+\tau)^\dagger+2\upl_\rho+(\upg_{\rho+1}+\varepsilon)\bigr)\\ &=\ (\upg_{\rho+1}+\varepsilon)\cdot\bigl(2(\upl_{\rho}-\upl_{\rho+1})+2(1+\tau)^\dagger+(\upg_{\rho+1}+\varepsilon)\bigr)\\
&=\ (\upg_{\rho+1}+\varepsilon)\cdot\bigl(-\upg_{\rho+1}+2(1+\tau)^\dagger+\varepsilon_1\bigr) \text{ (with $\varepsilon_1\prec \upg_{\rho+1}$)}\\ 
&\sim\ -\upg_{\rho+1}^2,
\end{align*}
which finishes the proof.
\end{proof}

\noindent
In view of Lemma~\ref{pcz} and its proof, this gives:

\begin{cor} \label{pcrhocor} $(\upo_\rho):=\big(\omega(\upl_\rho)\big)$ 
is a jammed pc-sequence of width 
$$\big\{\gamma\in \Gamma_{\infty}:\ \gamma > 2\Psi\big\}.$$
The sequence $(\upo_\rho+\upg_\rho^2)$ is also a pc-sequence in $K$, equivalent to $(\upo_\rho)$.
\end{cor}

\nomenclature[Z]{$(\upo_\rho)$}{$\big(\omega\big((-\ell_{\rho})^\dagger\big)\big)$}

\begin{example}
Suppose $K=\mathbb T$ and let $(\ell_n)$ be the logarithmic sequence for $K$ from Example~\ref{ex:iterated log sequence for T}. Using Corollary~\ref{cor:formula for omega} we see that for each $n$,
$$\upo_n\ =\ \omega(\upl_n)\ =\ \frac{1}{\ell_0^2} + \frac{1}{(\ell_0\ell_1)^2} + \cdots +  
\frac{1}{(\ell_0\ell_1\cdots\ell_n)^2}.$$ 
\end{example}

\begin{cor}\label{pczrho} 
Suppose $\upl_\rho \leadsto \upl$, with $\upl$ in an asymptotic 
field extension $L$ of~$K$. Then $\omega(\upl_\rho)\leadsto \omega(\upl)$ in $L$.
\end{cor}
\begin{proof} Eventually we have $\upl-\upl_\rho \sim \upg_{\rho+1}$, so a computation as in the proof of Lemma~\ref{pcrho} gives $\omega(\upl)-\omega(\upl_\rho)\ \sim\ \upg_{\rho+1}^2$, eventually.
\end{proof}

\noindent
It follows from Corollaries~\ref{pcrhocor} and \ref{pczrho}, Lemma~\ref{pc5}, and Proposition~\ref{indcof} that the pc-se\-quen\-ces $(\upo_\rho)$ obtained from different choices of 
the logarithmic se\-quence $(\ell_\rho)$
for $K$ are all equivalent.

\begin{lemma}\label{lem:upl upo infinitesimal}
Suppose $K$ has small derivation and asymptotic integration. Then $\upl_\rho,\upo_\rho\prec 1$ eventually.
\end{lemma}
\begin{proof}
Eventually $0<v(\upg_\rho)<(\Gamma^>)'$, so $v(\upl_{\rho})=\psi\big(v(\upg_{\rho})\big)>0$ eventually, by Lemma~\ref{BasicProperties}(i). Hence $\upl_\rho\prec 1$ eventually. 
Now use $\upo_\rho=-(2\upl_\rho'+\upl_\rho^2)$.
\end{proof}

\noindent
The next lemma will be used in proving
Proposition~\ref{cor:adding upo preserves upl-free}.

\begin{lemma}\label{lem:eval at uporho}
Suppose $K$ has rational asymptotic integration, and $\upo_{\rho} \leadsto
\upo$ with~$\upo$ in an immediate asymptotic extension~$F$ of~$K$. Let
$H\in F\{Y\}\setminus F$. Then $H(\upo_\rho)\leadsto H(\upo)$, and
there are $\gamma_0\in \Psi$ and $\delta_0\in (\Gamma^{>})'$ and an index 
$\rho_0$ such that the subset $\big\{v\big(H(\upo_\rho)-H(\upo)\big):\ \rho > \rho_0\big\}$ of $\Gamma_{\infty}$ is disjoint from $(\gamma_0, \delta_0)$. 
\end{lemma}
\begin{proof} Replacing $F$ by an algebraic closure $F^\alg$ and $K$ by the algebraic closure of~$K$ inside $F^\alg$, we arrange that $F$ is algebraically closed. Passing from $F$ to a suitable immediate asymptotic extension, we may further assume that $F$ is closed under integration, using 
Corollaries~\ref{achensacdiv} and~\ref{cor:value group of dv(K)} and Proposition~\ref{asintint}.
For each $\rho$, take $y_\rho\in F$ such that $(y_\rho')^2=\upo_\rho-\upo$.
Then $(y_\rho')^2 \asymp \upg_{\rho+1}^2 \asymp
(\ell_{\rho+2}')^2$, by Lemma~\ref{pcrho},
so
$y_\rho \asymp \ell_{\rho+2}$. Now put $Q:=H(Y+\upo)-H(\upo)$ and $P:=Q\big((Y')^2\big)$.
Then $P\neq 0$, $P(0)=0$, and $P(y_\rho)=H(\upo_\rho)-H(\upo)$ for each $\rho$.  
By Corollary~\ref{betterdifpol} we can take $\beta\in\Gamma^<$ and a strictly increasing map $i\colon (\beta,0)\to\Gamma$ such that $vP(y)=i(vy)$ for all $y\in F$ with $vy\in (\beta,0)$, hence $H(\upo_\rho)\leadsto H(\upo)$.
It now remains to apply Corollary~\ref{cor:bound on nwt(Ptimesq)} to $F$,~$Q$ in the role of $K$,~$P$, using
$H(\upo_\rho)-H(\upo)=Q\big((y_\rho')^2\big)$.
\end{proof}

\noindent
In connection with Proposition~\ref{indcof} we noted that if $a,b\in K$ are active with $a\succ b$, then $a^\dagger - b^\dagger\prec a$. For the function $\omega$ we have likewise:

\begin{lemma}\label{varrholemma, claim, 2}
Let $a,b\in K$ be active with $a\succ b$. Then 
$$\omega(-a^\dagger)-\omega(-b^\dagger)\ \prec\ a^2.$$
\end{lemma} 
\begin{proof}
Apply Lemma~\ref{varrholemma, claim, 1} to $w:= -a^\dagger$, $z:= -b^\dagger$, and $y:=z-w=a^\dagger-b^\dagger$, and note that $y^\dagger + w = y^\dagger - a^\dagger \prec a$ since $a\succ y$ and $y=(a/b)^\dagger$ is active. 
\end{proof}

\noindent
Here is an analogue of Lemma~\ref{lem:gap} for the sequence $(\upo_\rho)$:

\begin{lemma}\label{varrholemma} Let $\upo\in K$. Then the following are 
equivalent: \begin{enumerate}
\item[\textup{(i)}] $\upo_\rho \leadsto \upo$;
\item[\textup{(ii)}] for all active $a\in K$ we have $\upo-\omega(-a^\dagger) \prec a^2$; 
\item[\textup{(iii)}] for all active $a\in K$ there is an active $b\prec a$ in $K$ such that 
$$\upo - \omega(-b^\dagger) \prec a^2;$$
\item[\textup{(iv)}] for all $g\succ 1$ in $K$ we have 
$\upo - \omega(-g^{\dagger\dagger}) \prec (g^\dagger)^2$.
\end{enumerate}
\end{lemma}
\begin{proof} 
Assume (i). To prove (ii), let $a\in K$ be active. Take $\rho$ such that
$$\upg_{\rho}\prec a, \qquad \upo-\omega(\upl_{\rho})\sim \omega(\upl_{\rho'})-\omega(\upl_{\rho})\ \text{ for $\rho < \rho' < \kappa$.}$$
Then $\upo-\omega(\upl_{\rho}) \sim \upg_{\rho +1}^2\prec a^2$. Also $\omega(\upl_{\rho})-\omega(-a^\dagger)\prec a^2$ by Lemma~\ref{varrholemma, claim, 2} above for $b:=\upg_{\rho}$, hence $\upo-\omega(-a^\dagger) \prec a^2$. 

From Lemma~\ref{varrholemma, claim, 2}
we get (ii)~$\Rightarrow$~(iii). To get 
(iii)~$\Rightarrow$~(iv), assume~(iii), and let $g\in K$, $g\succ 1$. Take active $b\prec a:= g^\dagger$ in $K$ with $\upo-\omega(-b^\dagger) \prec a^2$.
By the above we also have $\omega(-b^\dagger)-\omega(-a^\dagger)\ \prec a^2$, so
$\upo-\omega(-a^\dagger) \prec a^2$. 

Finally, to prove (iv)~$\Rightarrow$~(i), assume (iv), and let $\rho < \kappa$. Take $g\succ 1$ in $K$ such that $a:= g^\dagger\prec \upg_{\rho+1}$ in $K$. 
Then $\upo-\omega(-a^{\dagger})\prec a^2$, so 
$$\upo-\omega(\upl_{\rho})\ =\ \big(\upo-\omega(-a^\dagger)\big) + \big(\omega(-a^\dagger)-\omega(\upl_{\rho +1})\big) + \big(\omega(\upl_{\rho+1})-\omega(\upl_{\rho})\big),$$
and so $\upo-\omega(\upl_{\rho})\sim \upg_{\rho +1}^2$, which gives (i).
\end{proof}

\begin{cor}\label{psacrho} The following four conditions on $K$ are equivalent: \begin{enumerate}
\item[\textup{(i)}] $(\upo_\rho)$ has no pseudolimit in $K$;
\item[\textup{(ii)}] for every $f\in K$ there is $g\succ 1$ in $K$  with 
$f-\omega(-g^{\dagger\dagger}) \succeq (g^\dagger)^2$; 
\item[\textup{(iii)}] for every $f\in K$ there is an active $a\in K$ with $f-\omega(-a^\dagger) \succeq a^2$;
\item[\textup{(iv)}] for every $f\in K$ there is an active $a\in K$ such that $f-\omega(-b^\dagger) \succeq a^2$ for all active $b\prec a$ in $K$. 
\end{enumerate}
\end{cor}

\subsection*{Relating $(\upl_{\rho})$ and $(\upo_{\rho})$} \nomenclature[Z]{$\upl(K)$}{cut defined by the pc-sequence $(\upl_\rho)$ in $K$}\nomenclature[Z]{$\upo(K)$}{cut defined by the pc-sequence $(\upo_\rho)$ in $K$} We introduce here
the cuts
$$\upl(K)\ :=\ c_K(\upl_\rho), \qquad \upo(K)\ :=\ c_K(\upo_\rho),$$ 
in $K$ defined by $(\upl_\rho)$ and $(\upo_\rho)$. By Lemmas~\ref{pcz} and \ref{pcrho},
\begin{align*}
\ndeg_{\upl(K)} P\  &=\ \text{eventual value of $\ndeg P_{+\upl_\rho,\times\upg_{\rho+1}}$, and}\\
\ndeg_{\upo(K)} P\ &=\ \text{eventual value of $\ndeg P_{+\upo_\rho,\times\upg^2_{\rho+1}}\ $ for $P\in K\{Y\}^{\ne}$}.
\end{align*}

\begin{lemma} \label{lem:ndeglambda P=2}
Suppose $\upo_{\rho} \leadsto \upo\in K$. Set $P(Z) := 2Z'+Z^2+\upo \in K\{Z\}$.  Let $\phi$ range over the elements of $K^\times$ with 
$v\phi\in\Psi^\downarrow$. Then for all $\rho$,
$$P^\phi_{+\upl_\rho,\times\upg_{\rho+1}}\ \sim\ \upg_{\rho+1}^2(Z-1)^2,\quad\text{ eventually.}$$ 
In particular, $\ndeg_{\upl(K)} P = 2$. 
\end{lemma}
\begin{proof}
Recalling that $\upo_\rho=-(2\upl_\rho'+\upl_\rho^2)$, we get
\begin{align*}
P_{+\upl_\rho,\times\upg_{\rho+1}}\  &=\ 2(\upg_{\rho+1}Z+\upl_\rho)'+
(\upg_{\rho+1}Z+\upl_\rho)^2+\upo \\
	&=\ (2\upg_{\rho+1}'Z+2\upg_{\rho+1}Z'+2\upl_\rho')+(\upg_{\rho+1}^2Z^2+2\upg_{\rho+1}\upl_\rho Z+\upl_\rho^2)+\upo \\
	&=\  \upg_{\rho+1}\big(2Z'+\upg_{\rho+1}Z^2+2(\upl_\rho-\upl_{\rho+1})Z\big)+(\upo-\upo_\rho).
\end{align*}
Hence
$$P_{+\upl_\rho,\times\upg_{\rho+1}}^\phi\ =\ \upg_{\rho+1}\big(2\phi Z'+\upg_{\rho+1}Z^2+2(\upl_\rho-\upl_{\rho+1})Z\big)+(\upo-\upo_\rho),$$
and for fixed $\rho$, we have $\phi \prec \upg_{\rho+1}\sim (\upl_{\rho+1}-\upl_\rho)$, eventually.
Since $\upo-\upo_\rho\sim\upg_{\rho+1}^2$ for all $\rho$, the lemma follows.
\end{proof}

%\noindent
%By Lemmas~\ref{dpkell} and ~\ref{zda, newton}, this has the following 
%consequence needed later:

\begin{cor}\label{uplupo} If $(\upo_{\rho})$ has a pseudolimit in $K$, 
then $K$ has an immediate asymptotic 
extension $L$ such that~$L$ is $\d$-algebraic over $K$ and not $\upl$-free.
\end{cor} 
\begin{proof} If $K$ is not $\upl$-free, take $L=K$. If
$K$ is $\upl$-free, then $K$ has rational asymptotic integration, and so we
are in the setting of Section~\ref{sec:construct imm exts}, 
and the desired result follows
from Lemmas~\ref{lem:ndeglambda P=2}, \ref{dpkell}, and~\ref{zda, newton}.
\end{proof}   

\noindent
To sharpen the above we first prove a lemma: 

\begin{lemma}\label{lem:upl on linear diff poly}
Suppose $K$ is $\upl$-free and $\upl$ is a pseudolimit of $(\upl_\rho)$ in an asymptotic field extension of $K$. Let $a,b\in K$ and set $Q(Z):=Z'+aZ+b\in K\{Z\}$. Then there is an $\alpha\in\Psi^\downarrow$ and $\rho_0<\kappa$ such that
$$
v\big(Q(\upl)-Q(\upl_\rho)\big)\ =\ v(\upg_{\rho+1}) + \alpha \quad\text{for $\rho_0<\rho<\kappa$.}$$
\end{lemma}
\begin{proof}
Let $\upl_\rho^*:=-(\upl_{\rho+1}-\upl_\rho)^\dagger$; then $\upl_\rho^*\not\leadsto a$
by Corollary~\ref{cor:uplrho alternative}, and so
we have $\alpha\in\Psi^\downarrow$ and $\rho_0<\kappa$ such
that $v(a-\upl_\rho^*)=\alpha$ for $\rho_0<\rho<\kappa$.
Also $\upl-\upl_\rho \sim \upl_{\rho+1}-\upl_\rho$ for all $\rho$. With $u_\rho:=(\upl-\upl_\rho)/(\upl_{\rho+1}-\upl_\rho)$ 
we get
$(\upl-\upl_\rho)^\dagger=u_\rho^\dagger+(\upl_{\rho+1}-\upl_\rho)^\dagger$ and $u_\rho\asymp 1$, so
$v(u_\rho^\dagger)>\Psi$, for $\rho_0<\rho<\kappa$.  Therefore
\begin{align*}
Q(\upl)-Q(\upl_\rho)\ 	&=\ (\upl-\upl_\rho) \big((\upl-\upl_\rho)^\dagger+a\big) \\
						&=\ (\upl-\upl_\rho) \big( u_\rho^\dagger + (a-\upl_\rho^*) \big) \\
						&\sim\ \upg_{\rho+1}  (a-\upl_\rho^*)		
\end{align*}
for $\rho_0<\rho<\kappa$. 
\end{proof}

\begin{prop}\label{prop:min diff poly for uplrho}
Suppose $K$ is $\upl$-free and $\upo_\rho\leadsto\upo\in K$. Then
$$P(Z)\ :=\ 2Z'+Z^2+\upo\in K\{Z\}$$
is a minimal differential polynomial of the pc-sequence $(\upl_\rho)$ over $K$.
\end{prop}
\begin{proof}
We have $P(\upl_\rho)=-\upo_\rho+\upo\leadsto 0$.
Towards a contradiction, suppose that $Q\in K\{Z\}\setminus K$ is a differential polynomial of
smaller complexity than $P$ such that ${Q(\tilde{\upl}_\sigma)\leadsto 0}$ for some pc-sequence $(\tilde{\upl}_\sigma)$ in $K$, equivalent to $(\upl_\rho)$.
Then there are two possibilities: either $\order(Q)=0$ (so $Q\in K[Z]$), or ${\order(Q)=\deg(Q)=1}$.
In the first case,~$(\upl_\rho)$ has a pseudolimit in the henselization of $K$, contradicting
Lem\-ma~\ref{lem:lambda-free henselization}. So we are in the second case, and may assume that  $Q=Z'+aZ+b$ where $a,b\in K$.
Now $Q$ is a minimal differential polynomial of $(\upl_\rho)$ over $K$, hence 
Lemma~\ref{zda, newton} and Corollary~\ref{zmindifpol}
give a pseudolimit $\upl$ of
$(\upl_\rho)$ in an immediate asymptotic field extension of $K$ with $Q(\upl)=0$.
Lemma~\ref{lem:upl on linear diff poly} yields 
$\alpha\in\Psi^\downarrow$ and~$\rho_0<\kappa$ such that
$v\big(Q(\upl_\rho)\big)=v(\upg_{\rho+1}) + \alpha$ for $\rho_0<\rho<\kappa$.
Increasing $\rho_0$ if necessary, we may assume that also $\alpha<v(\upg_{\rho+1})$ for $\rho_0<\rho<\kappa$.
Then by Lemma~\ref{pcrho},
$$P(\upl_\rho)\ =\  \upo-\upo_\rho\ \sim\ \upg_{\rho+1}^2\ 
\prec\ Q(\upl_\rho) \qquad (\rho_0 < \rho < \kappa)$$
and thus for $R:=P-2Q=Z^2-2aZ+(\upo-2b)\in K[Z]$ we get 
$R(\upl_\rho)\asymp Q(\upl_\rho)$, for $\rho_0<\rho<\kappa$. So $R(\upl_\rho)\leadsto 0$,
contradicting what we showed above.
\end{proof}

\noindent
Combining Proposition~\ref{prop:min diff poly for uplrho} with  Lemma~\ref{zda, newton} and Corollary~\ref{zmindifpol}
yields:

\begin{cor}\label{cor:adjoin upl}
Suppose $K$ is $\upl$-free and $\upo_\rho\leadsto\upo\in K$. 
Then $K$ has an immediate asymptotic field extension $K(\upl)$ with $\upl_\rho\leadsto\upl$
and $\omega(\upl)=\upo$, such that if
$L$ is any asymptotic field extension of $K$ and 
$\upl^*\in L$ satisfies $\upl_\rho\leadsto\upl^*$ and $\omega(\upl^*)=\upo$, then there is an embedding $K(\upl)\to L$ over $K$ sending $\upl$ to $\upl^*$.
\end{cor} 

\noindent
{\em Embedding\/} here means of course {\em embedding of valued differential fields}.

\subsection*{$\upo$-freeness}
Let us say that our $H$-asymptotic field $K$ is {\bf $\upo$-free} \label{p:upo-free} if $(\upo_{\rho})$ has no
pseudolimit in $K$ (so~$K$ satisfies the conditions of Corollary~\ref{psacrho}). It follows from Corollary~\ref{pczrho} that if
$K$ is $\upo$-free, then $K$ is $\upl$-free. Thus for the $H$-asymptotic fields $K$ considered in this section we have the implications
$$\text{$\upo$-free}\ \Rightarrow\ \text{$\upl$-free}\ \Rightarrow\ \text{rational asymptotic integration}\ \Rightarrow\ \text{asymptotic integration,}$$
where the last two properties are determined by the asymptotic couple of $K$. The pc-sequence $(\upo_n)$ in the Liouville closed $H$-field $\T$ is divergent, so $\T$ is $\upo$-free. 
See~\cite{ADH} for an example of a Liouville closed $H$-field
that is not $\upo$-free.

\medskip
\noindent
The $\upl$-free
$K$ in Corollary~\ref{recwam} are actually $\upo$-free: 

\index{H-asymptotic@$H$-asymptotic!field!$\upo$-free}
\index{omega-free@$\upo$-free}

\begin{lemma}\label{conam} Let $f\in K$ lie in a differential subfield
$F$ of $K$ for which $\max\Psi_F$ exists. Then there is an active $a\in K$ 
such that $f-\omega(-a^\dagger) \succeq a^2$.
\end{lemma}
\begin{proof} Take $a\in F$ such that $va=\max \Psi_F$. Then $a\in K$ is 
active, so in case $f-\omega(-a^\dagger) \succeq a^2$, we are done.
Suppose $f-\omega(-a^\dagger) \prec a^2$. In the algebraic closure~$K^\alg$ of $K$ with its differential subfield $F^\alg$  this gives 
$$v\left(\sqrt{f-\omega(-a^\dagger)}\,\right)\  >\ \Psi_F\ =\ \Psi_{F^\alg},$$
so $v\left(\sqrt{f-\omega(-a^\dagger)}\right) > \Psi$ by Corollary~\ref{asy1}, and thus
$v\big(f-\omega(-a^\dagger)\big)> 2\Psi$. 
Take $b\succ 1$ in $K$ with $vb'=va$. Then $vb+ vb^\dagger=va$, 
so $(a/b) \asymp b^\dagger$. Hence $(a/b)\in K$ is active, and $a\succ (a/b)$.
Then Lemma~\ref{varrholemma, claim, 1} applied to $w=-a^\dagger$, $z=-(a/b)^\dagger$, $y=z-w=b^\dagger$ gives
$$ \omega(-a^\dagger) - \omega(-(a/b)^\dagger) \ =\ 
b^\dagger\cdot \big(2(b^{\dagger\dagger}-a^\dagger) + b^\dagger\big).$$
We have $b^\dagger=(a/b)u$ with $u\asymp 1$, so $b^{\dagger\dagger}=a^\dagger - b^\dagger + u^\dagger$, hence $b^{\dagger\dagger}-a^\dagger=-b^\dagger + u^\dagger \sim -b^\dagger$ (since $K$ is pre-differential-valued), and thus
$$\omega(-a^\dagger) - \omega(-(a/b)^\dagger) \ \sim\ -(b^\dagger)^2\ \asymp\ (a/b)^2\ \prec\ a^2.$$
Together with $f-\omega(-a^\dagger)\prec (a/b)^2$ (a consequence of 
$v(f-\omega(-a^\dagger))> 2\Psi$),  it follows that
$f-\omega(-(a/b)^\dagger)\ \asymp (a/b)^2$. 
\end{proof}

\begin{cor}\label{recam} If $K$ is a union of grounded 
$H$-asymptotic subfields, then $K$ is $\upo$-free. 
\end{cor} 

\noindent
In view of how $\T$ is constructed, this shows again that
$\T$ is $\upo$-free.

\subsection*{Constructing $\upo$-free $H$-asymptotic fields}
Let $F$ be a grounded pre-$\d$-valued field of $H$-type, so we have $f\in F^\times$ with
$$f\succ 1, \qquad  vf^\dagger\ =\ \max \Psi.$$
We shall associate canonically to the pair $F$,~$f$ a pre-$\d$-valued extension $F_{\upo}$ of $F$ with $\res(F_{\upo})=\res(F)$, such that $F_{\upo}$ is of $H$-type, has asymptotic integration, and is $\upo$-free. In fact, by 
Lemma~\ref{pre-extas2} we
can take an increasing sequence 
$$F\ =\ F_0\ \subseteq\ F_1\ \subseteq\ F_2\ \subseteq\ F_3\ \subseteq\ \cdots$$ 
of pre-$\d$-valued extensions $F_n$ of $F$ of $H$-type, with distinguished elements $f_n\in F_n^\times$ such that $f_0=f$, and
for each $n$ we have 
$$f_n\succ 1,\qquad f_n^\dagger\ =\ \max \Psi_{F_n},\qquad
F_{n+1}\ =\ F_n(f_{n+1}),\qquad  f_{n+1}'\ =\ f_n^\dagger.$$
Then $F_{\upo}:= \bigcup_n F_n$ is an extension of $F$ with the properties
we announced, its $\upo$-freeness being a consequence of
each $F_n$ being grounded. From the universal property of Lemma~\ref{pre-extas2} we obtain:

\begin{lemma}\label{logupo} Any embedding of $F$ into a 
pre-$\d$-valued field $L$ of $H$-type that is closed under logarithms extends to an embedding $F_{\upo}\to L$.
\end{lemma}

%If $L$ is a pre-differential-valued extension of $F$ of $H$-type 
%and $L$ is closed under logarithms, then there is
%an $F$-embedding $F_{\upo} \to L$.
%\end{lemma}

\noindent
If $F$ is $\d$-valued, then each $F_n$ in the construction above is $\d$-valued, and so $F_{\upo}$ is $\d$-valued. This applies in particular to $\operatorname{dv}(F)$, since
$\operatorname{dv}(F)$ is a grounded $\d$-valued field of $H$-type. 

\medskip\noindent
Suppose $F$ is also equipped with an ordering making it a pre-$H$-field.
Then by Corollary~\ref{pre-extas2, pre-H}
there is a unique ordering on $F_{\upo}$ making it a pre-$H$-field extension
of $F$. Let $F_{\upo}$ be equipped with this ordering. Note that if $F$ is 
an $H$-field, then so is $F_{\upo}$. Lemma~\ref{logupo} goes through (and ``embedding'' in the setting of pre-$H$-fields means ``embedding of ordered valued differential fields''):

\begin{lemma}\label{Hlogupo} Any embedding of $F$ into a pre-$H$-field $L$ closed under logarithms extends to
an embedding $F_{\upo} \to L$. 
\end{lemma}

\noindent
For any grounded pre-$H$-field $F$, we let $F_{\upo}$ be a pre-$H$-field extension of $F$ as constructed above.
In particular, if $F$ is a grounded $H$-field, then $F_{\upo}$ is an $H$-field extension of $F$ with $C_{F_{\upo}}=C_F$.

\begin{cor}\label{cor:embed into upo-free}
Every pre-$\d$-valued field of $H$-type  has an $\upo$-free $\d$-valued extension of $H$-type, and every pre-$H$-field has an $\upo$-free $H$-field extension. 
\end{cor}
\begin{proof}
Let $E$ be a pre-$\d$-valued field of $H$-type. Then $E$ has a grounded $\d$-valued extension $F$ of $H$-type, by Corollary~\ref{cor:extend to smallest comp class}, and thus $F_{\upo}$ is an $\upo$-free
$\d$-valued extension of $H$-type of $E$. The argument for
pre-$H$-fields is the same.
\end{proof}

\begin{cor}\label{cor:embedintogap} Every pre-$H$-field has an $H$-field extension with a gap.
\end{cor}
\begin{proof} To extend a pre-$H$-field $E$ to an $H$-field with a gap, 
we can arrange by Corollary~\ref{cor:embed into upo-free} that $E$ is an 
$H$-field with rational asymptotic integration. Taking the real closure we 
arrange further that $E$ has divisible value group. The rest is done 
in the proof of Corollary~\ref{cor:extend to smallest comp class}.
%Corollary~\ref{zdcor} then allows us
%to pass from $E$ to an immediate extension and arrange
%that $E$ is in addition spherically complete. Then some element of 
%$E$ creates a gap over $E$, by Lemma~\ref{cg2}, and so~$E$ has an 
%$H$-field extension $E(f)$ with a gap $vf$ by results in 
%Section~\ref{sec:special cuts} 
%%and ~\ref{Solving Homogeneous First-Order Equations} 
%relating to gap creation. 
\end{proof}

\noindent
If $K$ has an $H$-asymptotic field extension $L$ with $\Gamma^<$  cofinal in $\Gamma_L^<$ and $\upo$-free $L$, then~$K$ is $\upo$-free. Analogues of Lemmas~\ref{lem:lambda-free completion}, \ref{lem:lambda-free fluentcompletion}, \ref{lem:lambda-free henselization}  
have similar proofs:

\begin{lemma}\label{upocompletion} If $K$ is $\upo$-free, then so is its completion $K^{\operatorname{c}}$.
\end{lemma}

\begin{lemma}\label{lem:upo-free fluentcompletion} If $K$ is $\upo$-free, then so is any immediate asymptotic extension  of~$K$ that, as a valued field extension of $K$, is a fluent completion of $K$.
\end{lemma}

\begin{lemma}\label{lem:omega-free fluent completion}
If $K$ is $\upo$-free, then so is its henselization $K^{\operatorname{h}}$.
\end{lemma}

\noindent
As in Corollary~\ref{cor:lambda-free alg closure} this yields:

\begin{cor}\label{cor:omega-free alg closure}
$K$ is $\upo$-free if and only if $K^\alg$ is $\upo$-free.
\end{cor}

\noindent
One direction of Corollary~\ref{cor:omega-free alg closure} will be vastly generalized in Section~\ref{consupofree}. 
 
\medskip\noindent
Specialization of $K$ with small derivation preserves $\upo$-freeness:

\begin{lemma}\label{specupofree}
Suppose $K$ is $\upo$-free with small derivation. Let $\Delta\neq\{0\}$ be a convex subgroup of $\Gamma$ with $\psi(\Delta^{\ne})\subseteq\Delta$. Then the asymptotic residue field $\dot K$ of $\dot v=v_\Delta$ is also $\upo$-free.
\end{lemma}
\begin{proof} We mimic the proof of Lemma~\ref{specuplfree} and use its notations.
Let $f\in \dot{\mathcal O}$; it is enough to get $a\in K^\times$ such that $va\in \Delta^{\ne}$, $a$ is active in~$K$, and $f-\omega(-a^\dagger)\succeq a^2$. 
By Corollary~\ref{psacrho} and since $\Psi^{>0}\ne \emptyset$, we can take active $b\in K$ such that $vb>0$ and $f-\omega(-b^\dagger)\succeq b^2$.
Then $vb\in \Delta$ by Lemma~\ref{asymp-lemma1}, and so $a:=b$ works.
\end{proof}

\noindent
Next, we consider the effect of compositional conjugation on the above. Let $\phi\in K^\times$ and let $\derdelta = \phi^{-1}\der$ be the derivation of the compositional conjugate $K^\phi$. As we saw before, the analogue of the sequence $(\upg_{\rho})$ in $K$ is the sequence
$(\upg^\phi_{\rho})=(\upg_\rho/\phi)$ in~$K^\phi$, and the analogue of the pc-sequence $(\upl_{\rho})$ in $K$ is 
the pc-sequence $(\upl^\phi_{\rho})$ in $K^\phi$, where 
$$\upl_\rho^\phi\  =\ (\upl_\rho/\phi) + (\phi^\dagger/\phi). $$
Also,   the role of $\omega=\omega_K \colon K \to K$ is taken over in 
$K^\phi$ by 
$$\omega^\phi:= \omega_{K^\phi} \colon\ K^\phi \to K^\phi, \quad \omega^\phi(z)\ :=\ -\big(2\derdelta(z) + z^2\big)\ =\ -(2/\phi) z' - z^2,$$
so 
\begin{align*} \omega^\phi(\upl_\rho^\phi)\ &=\  -\frac{2}{\phi}\cdot \left(\frac{\upl_\rho}{\phi} + \frac{\phi^\dagger}{\phi} \right)' - \left(\frac{\upl_\rho}{\phi} +\frac{\phi^\dagger}{\phi} \right)^2 \\
&=\ -\frac{2}{\phi}\cdot \left[ \frac{\upl_\rho'}{\phi} - \frac{\upl_\rho\phi^\dagger}{\phi} +  \left(\frac{\phi^\dagger}{\phi}\right)' \right] - 
\left(\frac{\upl_\rho}{\phi}\right)^2- \left(\frac{\phi^\dagger}{\phi}\right)^2-2\,\frac{\upl_\rho\phi^\dagger}{\phi^2}
  \\
              &=\ \frac{\omega(\upl_\rho)}{\phi^2} + \omega^\phi\left(\frac{\phi^\dagger}{\phi}\right).
\end{align*}
Hence, the sequence $(\upo_\rho)=\big(\omega(\upl_\rho)\big)$ has a pseudolimit in $K$ iff the
analogous sequence
$(\upo^\phi_\rho):=\big(\omega^\phi(\upl_\rho^\phi)\big)$  has a pseudolimit in $K^\phi$.
In particular, $K$ is $\upo$-free iff~$K^\phi$ is $\upo$-free.    
For later use we note that
\begin{align*}
 \omega^\phi(\phi^\dagger/\phi)\ 	&=\  -(2/\phi)\cdot (\phi'/\phi^2)' - (\phi'/\phi^2)^2  \\
									&=\   -(2/\phi)\cdot \big((\phi''/\phi^2) - 2(\phi')^2/\phi^3 \big)  - (\phi'/\phi^2)^2 \\
									&=\ -2(\phi''/\phi^3) + 3(\phi'/\phi^2)^2 =\ -\phi^{-2}\omega(-\phi^\dagger),
\end{align*}
so we can also express $\upo^\phi_\rho=\omega^\phi(\upl_\rho^\phi)$ as 
$$\upo^\phi_\rho\ =\ \phi^{-2}\big(\upo_\rho - \omega(-\phi^\dagger)\big),$$
in analogy to the equality
$$\upl_\rho^\phi\ =\ \phi^{-1}\big(\upl_\rho + \phi^\dagger\big).$$
More generally, let $E$ be any differential field, $\phi\in E^\times$, and 
define $\omega\colon E \to E$ by $\omega(z)=-2z'-z^2$ and
likewise $\omega^{\phi}\colon E^{\phi} \to E^{\phi}$ by
$\omega^{\phi}(z)=-2(z'/\phi)-z^2$. With $z^{\phi}:= \phi^{-1}(z+\phi^\dagger)$ for $z\in E$,
computations as before give the identity
$$\omega^{\phi}(z^{\phi})\ =\ \phi^{-2}\big(\omega(z) - \omega(-\phi^\dagger)\big).$$
Related to $\omega$ is the function $\sigma\colon E^\times \to E$
given by $\sigma(y)=\omega(-y^\dagger) + y^2$. The corresponding function
for $E^{\phi}$ is $\sigma^{\phi}\colon (E^{\phi})^\times \to E^{\phi}$ given by
$$\sigma^{\phi}(y)\ =\ \omega^{\phi}(-y^\dagger/\phi) + y^2.$$
It is routine to check that it satisfies a similar transformation formula:
$$ \sigma^{\phi}(y/\phi)\ =\ \phi^{-2}\big(\sigma(y)-\omega(-\phi^\dagger)\big)\qquad (y\in E^{\times}).$$

\section{Some Special Definable Sets}\label{sec:special sets} 

\noindent
{\em Throughout this section $K$ is a pre-$H$-field. If $K$ has asymptotic integration, then we fix a logarithmic sequence 
$(\ell_{\rho})$ for $K$ as in Section~\ref{sec:special cuts}, and corresponding sequences~$(\upg_{\rho})$, $(\upl_{\rho})$, and $(\upo_{\rho})$.}\/  
We single out the $\mathcal O$-submodule
$$ \I(K)\ :=\ 
\{y\in K:\ \text{$y\preceq f'$ for some $f\in \mathcal O$}\}$$
of $K$. We have $\der \mathcal O\subseteq \I(K)$ and $(\mathcal O^\times)^\dagger\subseteq\I(K)$.
If the derivation $\der$ of $K$ is small and~$K$ is an $H$-field, then $\I(K)\subseteq \smallo$. 
If $K$ has no gap or $K$ is an $H$-field, then 
$$ \I(K)\ =\ \{y\in K:\ \text{$y\preceq f'$ for some $f\in \smallo$}\}.$$

\begin{lemma} Suppose $K$ has a gap $\beta$. Then: 
$$\begin{cases}
\ \text{$v(b')\ne \beta$ for all $b\asymp 1$ in $K$} &   \Rightarrow\ \ \I(K)\ =\ \{ y\in K:\  vy>\beta \}; \\
\ \text{$v(b')=\beta$ for some $b\asymp 1$ in $K$} &   \Rightarrow\ \ \I(K)\ =\ \{y\in K:\ vy\ge \beta\}.
\end{cases}$$
\end{lemma}

\noindent
For $\phi\in K^>$ we have $\phi \I(K^\phi)=\I(K)$.
The following facts are also easy to verify:

\nomenclature[Z]{$\I(K)$}{$\{y\in K:\text{$y\preceq f'$ for some $f\in \mathcal O$}\}$}

\begin{lemma}\label{lem:I(K) under as int}
Suppose $K$ has asymptotic integration. Then for $y\in K$,
\begin{align*}
y\in\I(K) 	&\quad\Longleftrightarrow\quad  \text{$y \prec f^\dagger$ for all nonzero $f\in \smallo$} \\
			&\quad\Longleftrightarrow\quad  \text{$y$ is not active in $K$} \\
			&\quad\Longleftrightarrow\quad  \text{$y\prec \upg_\rho$ for all $\rho$}.
\end{align*}
If $\I(K)\subseteq \der K$, then
$\I(K)=\der\mathcal O$. If $\I(K)\subseteq (K^\times)^\dagger$, then
$\I(K)= (\mathcal O^\times)^\dagger$. Moreover, $(1/\ell_\rho)'\in \I(K)$ for all $\rho$, and for each $y\in\I(K)$ there is a $\rho$ such that $y\prec (1/\ell_\rho)'$.
Finally, if $L$ is a pre-$H$-field extension of~$K$, then $\I(K)=\I(L)\cap K$.
\end{lemma}

\begin{samepage}
\begin{cor}\label{cor:criterion for dv}
Suppose $K$ has asymptotic integration and $\I(K)\subseteq \der\smallo$. Then~$K$ is an $H$-field.
\end{cor}
\begin{proof}
Let $u\in \mathcal O$. Then $u'\in \I(K)$, so $u'=y'$ with $y\in \smallo$, and thus $u-y\in  C$.
\end{proof}
\end{samepage}

\begin{example}
$\I(\mathbb T) = \left\{y\in \mathbb T: \text{$y\prec \frac{1}{\ell_0\ell_1\cdots\ell_n}$ for all $n$}\right\}$,
where  $(\ell_n)$ is the logarithmic sequence for $\mathbb T$ from 
Example~\ref{ex:iterated log sequence for T}. 
\end{example}

\noindent
We also define
\begin{align*}
\Upg(K)	&\ :=\  \big\{ a^{\dagger}:\ a\in K^{\succ 1}\big\}, \\
\Upl(K)	&\ :=\  \big\{ {-a^{\dagger\dagger}}:\ a\in K^{\succ 1}\big\},\\
\Upd(K)	&\ :=\  \big\{ {-a^{\prime\dagger}}:\  0\ne a\in K^{\prec 1}\big\}.
\end{align*}
Note  that 
$$\Upl(K)\ =\  -\Upg(K)^\dagger\ =\ \big\{ {-a^{\dagger\dagger}}:\
 a\in K^\times,\ a\nasymp 1\big\}\ =\ 
  \big\{{-a^{\dagger\dagger}}:\ 0\ne a\in K^{\prec 1}\big\}.$$
For $\phi\in K^{>}$ we have 
$$\phi\Upg(K^\phi)\ =\ \Upg(K),\qquad \phi\Upl(K^\phi)\ =\ \phi^\dagger+\Upl(K), \qquad
\phi\Upd(K^\phi)\ =\ \phi^\dagger+\Upd(K).$$
 The set $\Upg(K)$ is closed under addition. The elements of $\Upg(K)$ are active in $K$. 
The sets~$\I(K)$, $\Upg(K)$, and $-\Upg(K)$ are 
pairwise disjoint. If $K$ has asymptotic integration, then there is for
each $y\in \Upg(K)$ an index $\rho$ with $\upg_\rho\prec y$.
If  $K$ has asymptotic integration and $(K^\times)^\dagger=K$, 
then $$K=\I(K)\cup \Upg(K) \cup \big(\!-\Upg(K)\big),$$ 
by Lemma~\ref{lem:I(K) under as int}.
It is easy to see that if $(K^\times)^\dagger=K$ and $L$ is a pre-$H$-field extension of
$K$, then $\Upg(L)\cap K=\Upg(K)$.

\nomenclature[Z]{$\Upg(K)$}{$(K^{\succ 1})^\dagger$}
\nomenclature[Z]{$\Upl(K)$}{$-(K^{\succ 1})^{\dagger\dagger}$}
\nomenclature[Z]{$\Upd(K)$}{$-(K^{\ne,\prec 1})^{\prime\dagger}$}

\medskip
\noindent
Let $a\in K^{\succ 1}$, and set
$y:=a^\dagger\in\Upg(K)$ and $z=-y^\dagger=-a^{\dagger\dagger}\in\Upl(K)$. Then 
$$z+y\ =\ -\big(a^{\dagger\dagger}+(1/a)^\dagger\big)\ =\  - \big( (1/a)^\dagger\, (1/a) \big)^\dagger\ =\ -(1/a)^{\prime\dagger}\in\Upd(K).$$
In particular, if $K$ has asymptotic integration, then 
$\upl_\rho\in\Upl(K)$  and $\upl_\rho + \upg_\rho \in \Upd(K)$ for all $\rho$, and
$\Upl(K)\cup\Upd(K)$ does not contain a pseudolimit of
$(\upl_\rho)$, by Lemma~\ref{lem:active S}. 
The sets $\Upl(K)$ and $\Upd(K)$ are disjoint.

\begin{lemma}\label{lem:Upl}
Suppose  $K$ has asymptotic integration, $\I(K)=\der\smallo$, and $(K^\times)^\dagger=K$.
Then $K=\Upl(K)\cup\Upd(K)$.
\end{lemma}
\begin{proof}
Let $f\in K$, and take $a\in K^\times$ with $a^\dagger\neq 0$ and
$f=-a^{\dagger\dagger}$. If $a\nasymp 1$, then $f\in\Upl(K)$, so suppose $a\asymp 1$.
Then $a^\dagger \in (\mathcal O^\times)^\dagger \subseteq \I(K)=\der\smallo$. Take
$b\in\smallo$ with $a^\dagger=b'$; then $b\neq 0$, and $f=-b^{\prime\dagger}\in\Upd(K)$.
\end{proof}

\begin{lemma}\label{lem:Upl, 1}
Let $f\in K^\times$. Then
$$f\in\I(K)\quad\Longrightarrow\quad -f^\dagger \notin\Upl(K),$$
and if $(K^\times)^\dagger=K$, then
$$f\in\I(K)\quad\Longleftrightarrow\quad -f^\dagger \notin\Upl(K).$$
\end{lemma}
\begin{proof}
Suppose that  $-f^\dagger\in\Upl(K)$.
Take $a\in K^{\times}$, $a\nasymp 1$, such that $-f^\dagger=-a^{\dagger\dagger}$. Then $f\in C^\times a^\dagger$
and hence $vf=va^\dagger\in\Psi$, so $f\notin\I(K)$.
Now suppose $K$ has asymptotic integration, $(K^\times)^\dagger=K$, and
$-f^\dagger\notin\Upl(K)$. Take $a\in K^\times$ such that $f=a^\dagger$;
then $-a^{\dagger\dagger}=-f^\dagger\notin\Upl(K)$ and hence $a\asymp 1$, so
$f=a^\dagger\in (\mathcal O^\times)^\dagger\subseteq \I(K)$.
\end{proof}

\begin{cor}\label{prupl}
Suppose $K$ has asymptotic integration and $(K^\times)^\dagger=K$, and
let~$L$ be a pre-$H$-field extension of $K$. Then $\Upl(K)=\Upl(L)\cap K$.
\end{cor}
\begin{proof}
Clearly
$\Upl(K)\subseteq\Upl(L)\cap K$, so let $g\in \Upl(L)\cap K$, and take $f\in K^{\times}$ with
$g=-f^{\dagger}$. Then $f\notin\I(L)$ by the first part of Lemma~\ref{lem:Upl, 1},
hence $f\notin\I(K)$ and so $g\in\Upl(K)$, by the second part.
\end{proof}

\begin{lemma}\label{lem:Upl, 2}
Let  $f\in K^\times$. Then
$$f\in\der\smallo\quad\Longleftrightarrow\quad -f^\dagger \in\Upd(K).$$
\end{lemma}
\begin{proof}
If $f\in\der\smallo$, say $f=a'$ with $a\in \smallo$, then 
$-f^\dagger=-a^{\prime\dagger}\in\Upd(K)$. Conversely, if
$-f^\dagger\in\Upd(K)$, then taking $a\in \smallo^{\ne}$ with $f^\dagger=a^{\prime\dagger}$ gives
$f\in Ca'\subseteq\der\smallo$.
\end{proof}

\begin{cor}\label{asupik}
Suppose that $K$ has asymptotic integration, $\I(K)=\der\smallo$, and 
$(K^\times)^\dagger=K$.
Then for $f\in K^\times$ we have
$$f\in\I(K)\ \Longleftrightarrow\ -f^\dagger\in\Upd(K),\qquad
  f\notin\I(K) \ \Longleftrightarrow\ -f^\dagger\in\Upl(K).$$
\end{cor}

\begin{lemma}\label{lem:invariance gp}
Suppose $(K^\times)^\dagger=K$. Then
$\Upl(K)+\I(K)\ \subseteq\ \Upl(K)$.
\end{lemma}
\begin{proof}
Let $a\in K$, $a\succ 1$ and $y\in \I(K)$; we want to show $-a^{\dagger\dagger}+y\in\Upl(K)$. Take $b\in K\setminus C$ with $b^{\dagger\dagger}=a^{\dagger\dagger}-y$.
Then $(a^\dagger/b^\dagger)^\dagger=y\in\I(K)$, hence $a^\dagger/b^\dagger\asymp 1$ and therefore $b\nasymp 1$.
\end{proof}

\subsection*{First-order linear differential equations
and the predicate $\operatorname{I}$}{\em Throughout this subsection
$K$ is a Liouville closed $H$-field and $y$ ranges over $K$}.

\begin{lemma}\label{Ifirst} Suppose the derivation of $K$ is small. Let
$A=aY + bY'$ where $0\ne a\succeq b$ in $K$.
Then $A(\mathcal{O}) = a\mathcal{O}$.
\end{lemma}
\begin{proof} Dividing by $a$ we arrange $a=1$. Let $g\in \mathcal{O}$; our job is to find $f\in \mathcal{O}$ such that
$A(f)=g$. Take $f_0\in K$ with $A(f_0)=g$, and $y_0\in K^\times$ with $A(y_0)=0$. Then 
$$A^{-1}(g)\ :=\ \big\{f\in K:\ A(f)=g\big\}\ =\ f_0 + Cy_0,$$
so $\big|v\big(A^{-1}(g)\big)\big|\le 2$. Let $\phi\in K^{>}$ be active, $\phi\preceq 1$. Then the compositional conjugate~$K^{\phi}$ is a Liouville closed $H$-field with small derivation, and $A^{\phi}=Y+b\phi Y'$ with $b\phi\preceq 1$. We equip $K^{\phi}$ with the coarsening
$v^{\flat}_{\phi}$. Then 
$\max \Psi^{\flat}_{\phi}=0$ in $\Gamma/\Gamma^{\flat}_{\phi}$, where~$\Psi^{\flat}_{\phi}$ is the $\Psi$-set of the asymptotic
couple of $(K^{\phi}, v^{\flat}_{\phi})$. By Lemma~\ref{Specialization of asymptotic fields} the residue map
from~$\mathcal{O}^{\flat}_{\phi}$ onto its residue field maps the constant field
of $K$ onto the constant field of that (differential) residue field. Thus we can apply Lemmas~\ref{sup=0neatsur} and~\ref{nslinsur} to~${(K^{\phi}, v^{\flat}_{\phi})}$ to get
$f_{\phi}\in K$ with
$A^{\phi}(f_{\phi})=A(f_{\phi})=g$ (so 
$f_{\phi}\in A^{-1}(g)$) and 
$v^{\flat}_{\phi}(f_{\phi})\ge 0$. In view of 
$\big|v\big(A^{-1}(g)\big)\big|\le 2$ we get $v(f_{\phi})\ge 0$, eventually.
\end{proof}

\begin{lemma}\label{lem:qelinI} 
Let $f,g\in K$, $g\neq 0$. Then 
$$\exists y\,\big[ y'=fy+g\ \&\ y\prec 1\big]\quad \Longleftrightarrow\quad 
f,g\in\operatorname{I}(K) \ \vee\ \big[f\notin \operatorname{I}(K) \wedge g\prec f\big].$$
\end{lemma}
\begin{proof}
Put $P:=Y'-fY-g\in K\{Y\}$.
Take $b\ne 0$ with $b^\dagger = f$ and then 
$a\ne 0$ with $a\nasymp b$ and
$va+\psi(va-vb)=vg$. For active $\phi$ in $K$ we have 
$$ P^\phi_{\times a}\ =\ a\phi Y'+a(a/b)^\dagger Y- g$$
with $a(a/b)^\dagger\asymp g$. Since eventually 
$\phi\preceq (a/b)^\dagger$, Lemma~\ref{Ifirst} yields a zero of $P_{\times a}$ in $\mathcal O$. But $P_{\times a}$ 
has no zero in $\smallo$, so $P$ has a zero $z\in K$ with $z\asymp a$, and for every zero~$y$ of $P$ we have~$y\succeq a$.
Thus $P$ has a zero $y\prec 1$ iff $a\prec 1$.
Now suppose $f\in\operatorname{I}(K)$. Then $a\nasymp b\asymp 1$ and $va+\psi(va)=vg$, that is, $a'\asymp g$; so $a\prec 1$ iff
$g\in\operatorname{I}(K)$. Next assume $f\notin\operatorname{I}(K)$, so $b\nasymp 1$. Then
\begin{align*}
a\prec 1	&\quad\Longleftrightarrow\quad va-vb>-vb \\
			&\quad\Longleftrightarrow\quad \underbrace{va-vb+\psi(va-vb)}_{=\,vg-vb}>\underbrace{-vb+\psi(-vb)}_{=\,vf-vb} \\
			&\quad\Longleftrightarrow\quad g\prec f.\qedhere
\end{align*}
\end{proof}

\begin{cor}\label{qelinI}
Let $f,g,h\in K$, $g,h\neq 0$. Then the following are  equivalent: \begin{enumerate}
\item[\textup{(i)}] there exists $y$ such that $y'=fy+g$ and $y\prec h$;
\item[\textup{(ii)}]
$[f-h^\dagger\in \operatorname{I}(K) \text{ and }g/h\in\operatorname{I}(K)]\ \ \text{ or }\ \big[ f-h^\dagger\notin \operatorname{I}(K) \text{ and }  (g/h)\prec f-h^\dagger\big]$.
\end{enumerate}
\end{cor}
\begin{proof}
Set $f_*:=f-h^\dagger$, $g_*:=g/h$. 
Then $(y/h)'=y'/h-h^\dagger(y/h)$, so
$y'=fy+g$ and $y\prec h$ is equivalent to
$(y/h)'=f_*(y/h)+g_*$ and $y/h\prec 1$. 
It remains to appeal to
Lemma~\ref{lem:qelinI}.
\end{proof}

\subsection*{Interaction with the ordering}
So far the field ordering of our pre-$H$-field $K$ has not played a role.
As a first comment referring to this ordering we mention that $\I(K)$
is convex in the ordered set $K$. 
Below we characterize $\upl$-freeness and $\upo$-freeness in terms of certain cuts in the ordered set $K$, and we explore further the behavior of the differential polynomial function~$\omega$ in this setting.
We often tacitly use the fact that for $a,b\in K^\times$ we have $a\prec b\Rightarrow a^\dagger < b^\dagger$. One consequence of this fact is that $\Upl(K)<\Upd(K)$, which in view of Lemma~\ref{lem:Upl} yields:

\begin{cor}\label{cor:Upl 2}
If $K$ is Liouville closed, then  $\Upl(K)$ is downward closed,  $\Upd(K)$ is 
upward closed, and $K=\Upl(K)\cup\Upd(K)$.
\end{cor}
%\begin{proof}
%Immediate from Lemma~\ref{lem:Upl} and Corollary~\ref{cor:Upl 1}.
%\end{proof}

\noindent
Here is another use of  
$a\prec b\Rightarrow a^\dagger < b^\dagger$, for $a,b\in K^\times$:

\begin{lemma}\label{lem:coinitial Upl Upd}
Assume $K$ is ungrounded and $L\supseteq K$ is a pre-$H$-field extension such that $\Gamma^>$ is coinitial in~$\Gamma_L^>$.
Then $\Upd(K)$ is coinitial in $\Upd(L)$, $\Upg(K)$ is coinitial in $\Upg(L)$, and 
$\Upl(K)$ is cofinal in $\Upl(L)$.
\end{lemma}
\begin{proof}
Let
$a\in L^{\prec 1}$, $a\ne 0$. Take $b\in K$ with $a\prec b\prec 1$. Then $a'\prec b'$ and hence $-a^{\prime\dagger} > -b^{\prime\dagger}$. This shows
$\Upd(K)$ is coinitial in $\Upd(L)$. Next, let
$a\in L^{\succ 1}$. Take $b\in K^{\succ 1}$ with $a^\dagger \succ b^\dagger$.
Then $a^\dagger,b^\dagger>0$, hence $a^\dagger>b^\dagger$
and
$-a^{\dagger\dagger} < -b^{\dagger\dagger}$.
\end{proof}

\begin{lemma}\label{lem:uplrho strictly increasing} Let $K$ have asymptotic 
integration. Then
$(\upl_\rho)$ is strictly increasing and cofinal in $\Upl(K)$,
and $(\upl_\rho+\upg_\rho)$ is strictly decreasing and coinitial in 
$\Upd(K)$. If~$K$ has small derivation, then $\upl_{\rho}>0$, eventually.
\end{lemma}
\begin{proof}
We have $\Upl(K) = -\Upg(K)^\dagger$ and $\upl_\rho=-\upg_\rho^\dagger$, so the 
first claim follows from the remark preceding the previous lemma and the properties of
$(\upg_\rho)$.
Let $a\in K$, $0\neq a\prec 1$. Take $\rho$ with $1/\ell_\rho \succ a$; then $(1/\ell_\rho)' \succ a'$, so
$-(1/\ell_\rho)^{\prime\dagger} < -a^{\prime\dagger}$.
Thus $(\upl_\rho+\upg_\rho)$ is coinitial in $\Upd(K)$. For $\rho<\rho'$ we have $1/\ell_\rho \prec 1/\ell_{\rho'}\prec 1$, hence  $(1/\ell_\rho)' \prec (1/\ell_{\rho'})'$ and thus $-(1/\ell_\rho)^{\prime\dagger} > -(1/\ell_{\rho'})^{\prime\dagger}$; so $(\upl_\rho+\upg_\rho)$ is strictly decreasing. Suppose now that $K$ has small derivation. Then $1$ is active, so $\upg_{\rho}\prec 1$ eventually, and thus 
$\upl_{\rho}>0$ eventually. 
\end{proof}

\begin{cor}\label{cor:Upl 1} Assume $K$ has asymptotic integration. Let $L$ be a pre-$H$-field extension of $K$ and $\upl\in L$. Then:
$\upl_{\rho} \leadsto \upl\ \Longleftrightarrow\ 
\Upl(K)<\upl<\Upd(K)$.
\end{cor}
\begin{proof}
Use the previous lemma in combination with Lemma~\ref{pcz}.
\end{proof}

\begin{cor}\label{cor:invariance gp}
Suppose $K$ has asymptotic integration and $L$ is a pre-$H$-field extension of $K$ with $(L^\times)^\dagger=L$. Then $\Upl(L)^\downarrow\cap K$
equals $\Upl(K)^\downarrow$ or $K\setminus\Upd(K)^\uparrow$.
\end{cor} 
\begin{proof}
Suppose $\Upl(L)^\downarrow\cap K \neq \Upl(K)^\downarrow$.
Take $f\in \big(\Upl(L)^\downarrow\cap K\big) \setminus \Upl(K)^\downarrow$.
Then $f< \Upd(L)$, so $\Upl(K)<f<\Upd(K)$.
To get $\Upl(L)^\downarrow\cap K = K\setminus\Upd(K)^\uparrow$, 
let $g\in K\setminus\Upd(K)^\uparrow$; it is enough to show that then
$g\in \Upl(L)^\downarrow$. As the case
$g\in\Upl(K)^\downarrow$ is obvious, assume $g\notin\Upl(K)^\downarrow$.
Then $\Upl(K)<g<\Upd(K)$.
Then by Corollary~\ref{cor:Upl 1} both $f$ and $g$ are pseudolimits of $(\upl_\rho)$, hence $v(g-f)>\Psi$ by
Lemma~\ref{pcz} and thus $g-f\in \I(K)\subseteq \I(L)$. Therefore 
$g\in f+\I(L)\subseteq \Upl(L)^{\downarrow}$ by Lemma~\ref{lem:invariance gp}.
\end{proof}

%\begin{cor}\label{cor:Upl 2}
%If $K$ is Liouville closed, then  $\Upl(K)$ is closed downward,  
%$\Upd(K)$ is closed upward,
%and $K=\Upl(K)\cup\Upd(K)$.
%\end{cor}
%\begin{proof}
%Immediate from Lemma~\ref{lem:Upl} and Corollary~\ref{cor:Upl 1}.
%\end{proof}

\noindent
Our eventual quantifier elimination for $\T$ requires predicates 
for the sets $\Upl(\T)$ and $\omega(\T)$, and for this reason we
pay attention to properties of these sets expressible by
universal sentences. In this connection it is convenient to
include also $\I(\T)$ and $\Upd(\T)$ in an auxiliary role. 
Such universal properties are 
contained in Lemmas~\ref{lem:Upl, 1},~\ref{lem:invariance gp},
and Co\-rol\-lary~\ref{cor:Upl 2}. Here is another useful one:

\begin{lemma} Assume $K$ has asymptotic integration, and let  $a,\phi\in K$. Then
$$ a, \phi>0,\ a\succeq 1\ \Rightarrow\ \phi a - \phi^\dagger\in \Upd(K)^{\uparrow}.$$
\end{lemma}
\begin{proof} \begin{samepage} We first establish the following: 
\begin{enumerate}
\item if the derivation of $K$ is small, then $\smallo\cap \Upd(K)\ne \emptyset$;
\item if the derivation of $K$ is not small, then $K^{<}\cap \Upd(K)\ne \emptyset$.
\end{enumerate}
\end{samepage}
For (1), assume the derivation of $K$ is small.
Take $1\in\Gamma^{>}$ with $\psi(1)=1$.
Then $\Psi< 1+1$. Pick $a\in \smallo$ with $v(a')=1+1$. Then
$v(-a'^\dagger)=1$ and so $-a'^\dagger\in \smallo \cap \Upd(K)$. 

For (2), assume the derivation of $K$ is not small. Then for each
$\rho$ we have $\upg_{\rho}\succ 1$, so $\upl_{\rho}$ is active
and $\upl_{\rho}< 0$. By Lemma~\ref{lem:active S}, $v(\upl_{\rho})$
is eventually constant, and in $\Psi$. Hence
$\upl_{\rho} \succ \upg_{\rho}$, eventually, and thus  
$\upl_{\rho}+ \upg_{\rho}\in K^{<} \cap \Upd(K)$, eventually. 
     
Now let $a, \phi > 0$ and $a\succeq 1$. Then 
$a\in \Upd(K^{\phi})^{\uparrow}$ by (1) and (2).  Also
$\Upd(K^{\phi})=\phi^{-1}\big(\Upd(K) + \phi^\dagger\big)$, which gives 
$\phi a - \phi^\dagger\in \Upd(K)^{\uparrow}$.
\end{proof}

\begin{lemma}\label{lem:Upg, Liouville closed}
Suppose $K$ is Liouville closed. Then
$$\Upg(K)\ =\  (K^{>C})^\dagger\ =\  \{ y\in K^>: vy\in\Psi \}\ =\  K^> \setminus \I(K),$$
and hence $\Upg(K)$ is upward closed.  
\end{lemma}
\begin{proof} The set equalities are easy consequences of the results above.
To get $\Upg(K)$ upward closed, use that $\I(K)$ is convex.  
\end{proof}

\begin{example}
Let $(\ell_n)$ be the logarithmic sequence for $\mathbb T$ from Example~\ref{ex:iterated log sequence for T}. Then 
for $y\in\mathbb T$ we have
$$y\in\Upg(\mathbb T) \quad\Longleftrightarrow\quad \text{$y\ \geq\  \frac{1}{\ell_0\ell_1\cdots\ell_n}\ $ for some $n$}.$$
Also, for $z\in\mathbb T$ we have
$$z\in\Upl(\mathbb T)	\quad\Longleftrightarrow\quad \text{$z\ \leq\ \frac{1}{\ell_0}  + \frac{1}{\ell_0\ell_1} + \cdots + \frac{1}{\ell_0\ell_1\cdots\ell_n}\ $ for some $n$}$$
and
\begin{align*}
z\in \Upd(\mathbb T)	&\ \Longleftrightarrow\   \text{$z\ \geq\ \left(\frac{1}{\ell_0} +  \frac{1}{\ell_0\ell_1} +\cdots +  \frac{1}{\ell_0\ell_1\cdots\ell_{n}}\right) + \frac{1}{\ell_0\ell_1\cdots\ell_n}\ $ for some $n$} \\
					&\ \Longleftrightarrow\   
\text{$z\ >\ \frac{1}{\ell_0}  + \frac{1}{\ell_0\ell_1} + \cdots + \frac{1}{\ell_0\ell_1\cdots\ell_n}\ $ for all $n$.}
\end{align*}
Figure~\ref{fig:sets} shows the sets $\I(\T)$, $\Upgamma(\T)$, $\Uplambda(\T)$, and $\Updelta(\T)$.
\end{example}

\begin{figure}[h!]
\begin{tikzpicture}[scale=0.70]

\draw (2.5,0) -- (11.5,0);
\draw (13,0) -- (13.5,0);
\draw[dashed] (11.5,0) -- (13,0);

\draw (14.75,0) -- (16,0); 
\draw[dashed] (16,0) -- (17,0);

\draw (0,-0.5) node {$0$};
\fill (0,0) circle (0.2em);

\draw[densely dotted] (1.5,0) -- (1.75,0);
\draw[densely dotted] (2.25,0) -- (2.5,0);
\draw[dotted] (1.75,0) -- (2.25,0);
\draw (2.0,-2.5) node {\mbox{\fontsize{9}{0}\selectfont $\frac{1}{\ell_0\ell_1\cdots}$}};
\draw[-latex, gray, thick] (2,-2) -- (2,-0.1);

\draw[densely dotted] (2,-3) -- (2,-5);

\fill[white] (2,0) circle (0.2em);
\draw (2,0) circle (0.2em);

\draw (-0.5,0) -- (1.5,0);
\draw[dashed] (-1,0) -- (-0.5,0);

\draw[densely dotted] (13.5,0) -- (14.75,0);
\draw[densely dotted] (14.25,0) -- (14.5,0);
\draw[dotted] (13.75,0) -- (14.25,0);
\draw (15.25,-2.5) node {\mbox{\fontsize{9}{0}\selectfont $\frac{1}{\ell_0}+\frac{1}{\ell_0\ell_1}+\cdots$}};
\draw[-latex, gray, thick] (14,-2) -- (14,-0.1);
\fill[white] (14,0) circle (0.2em);
\draw (14,0) circle (0.2em);

\draw[densely dotted] (14,0.1) -- (14,2);

\fill (4,0) circle (0.1em);
\draw (7,-2.5) node {\mbox{\fontsize{9}{0}\selectfont $\upg_0=\upl_0$}};
\draw[-latex, gray, thick] (6.5,-2) -- (4,-0.1);

\fill (3.5,0) circle (0.1em);
\draw (5.75,-2.5) node {\mbox{\fontsize{9}{0}\selectfont $\upg_1$}};
\draw[-latex, gray, thick] (5.5,-2) -- (3.5,-0.1);

\fill (3.25,0) circle (0.1em);
\draw (4.75,-2.5) node {\mbox{\fontsize{9}{0}\selectfont $\upg_2$}};
\draw[-latex, gray, thick] (4.5,-2) -- (3.25,-0.1);

\draw (3.5,-2.5) node {\mbox{\fontsize{9}{0}\selectfont $\cdots$}};

\fill (7.5,0) circle (0.1em);
\draw (8.75,-2.5) node {\mbox{\fontsize{9}{0}\selectfont $\upl_1$}};
\draw[-latex, gray, thick] (8.5,-2) -- (7.5,-0.1);

\fill (10.75,0) circle (0.1em);
\draw (10,-2.5) node {\mbox{\fontsize{9}{0}\selectfont $\upl_2$}};
\draw[-latex, gray, thick] (10,-2) -- (10.75,-0.1);

\draw (11.5,-2.5) node {\mbox{\fontsize{9}{0}\selectfont $\cdots$}};

\draw (1,-4) node {$\I(\T)$};
\draw[-latex, gray, thick] (1.5,-4.5) -- (0.25,-4.5);
\draw (3.125,-4) node {$\Upgamma(\T)$};
\draw[-latex, gray, thick] (2.5,-4.5) -- (3.75,-4.5);

\draw (13,1.5) node {$\Uplambda(\T)$};
\draw[-latex, gray, thick] (13.5,1) -- (12.25,1);
\draw (15.125,1.5) node {$\Updelta(\T)$};
\draw[-latex, gray, thick] (14.5,1) -- (15.75,1);

\end{tikzpicture}
\caption{The sets $\I(\T)$, $\Upgamma(\T)$, $\Uplambda(\T)$, and $\Updelta(\T)$.}\label{fig:sets}
\end{figure}

\noindent
Next we study the behavior of $\omega\colon K \to K$, $\omega(z):=-(2z'+z^2)$.

\begin{prop}\label{prop:omega for H-fields}
Suppose $K$ is Liouville closed. Then 
the restriction of $\omega$ to~$\Upl(K)$ is strictly increasing,
and $\omega\big(\Upl(K)\big)=\omega\big(\Upd(K)\big)$.
\end{prop}

\begin{proof}
Pick for each $f\in K^{\succ 1}$ an $\Log f\in K^{\succ 1}$ with 
$(\Log f)'=f^\dagger$. Thus we have:
$f\in K^{>C} \Rightarrow \Log f>C$. 
We also pick for all $f\in K^{>}$ and
$r\in \Q$ an element
$f^r\in K^{>}$ such that $(f^r)^\dagger=rf^\dagger$.
It follows easily that for $f\in K$, $r\in \Q$, 
$$f> C,\ r>0\ \Rightarrow\ f^r > C, \qquad
f>C,\ r<0\ \Rightarrow\ f^r \prec 1.$$ 
Now let $z\in\Upl(K)$. We shall find an element $f\in \Upd(K)$ 
(depending on~$z$) such that 
$\omega(z)<\omega(z^*)$ for all $z^*\in K$ with $z<z^*\leq f$. As $z$ is arbitrary and $f>\Upl(K)$, this will prove that
$\omega$ is strictly increasing on $\Upl(K)$.
Take $a\in K^{>C}$ such that $z=-a^{\dagger\dagger}$, and let $t:=(\Log a)^{-r}$ where $r$ is any rational
number with $0<r\leq 1$. (The value of $r$ will be further specified later.) Then $0\ne t\prec 1$, so $f:=-t^{\prime\dagger}\in\Upd(K)$. 
From $t'=-r(\Log a)^\dagger t= -r\frac{a^\dagger}{\Log a}t$ we get,
with $y:=a^\dagger=(\Log a)'$:
$$f\  =\ -y^\dagger + (r+1)(\Log a)^\dagger\  =\ z + (r+1)\frac{y}{\Log a}.$$
Put
$$\varepsilon\ :=\ h(r+1)\frac{y}{\Log a}\qquad\text{with $h\in K$, $0<h\leq 1$,}$$
so $\varepsilon>0$ and $z+\varepsilon$ runs through  all
$w\in K$ with $z<w\leq f$, as $h$ varies. Thus it suffices to show that
$\omega(z+\varepsilon)>\omega(z)$. 
Since
\begin{equation}\label{eq:omega1}
\omega(z) - \omega(z+\varepsilon)\  =\ \varepsilon \big( 2(\varepsilon^\dagger+z)+\varepsilon \big)
\end{equation}
by Lemma~\ref{varrholemma, claim, 1},
this reduces to showing that $2(\varepsilon^\dagger+z)+\varepsilon <0$.
Now
$$\varepsilon^\dagger\ =\ h^\dagger + y^\dagger - (\Log a)^\dagger\ =\ h^\dagger - z - \frac{y}{\Log a}$$
and thus
\begin{equation}\label{eq:omega2}
2(\varepsilon^\dagger+z)+\varepsilon\ =\ 2h^\dagger - \big(2-h(r+1)\big) \frac{y}{\Log a},
\end{equation}
so the inequality $2(\varepsilon^\dagger+z)+\varepsilon <0$ becomes
$$2h^\dagger\  <\ \big(2-h(r+1)\big) \frac{y}{\Log a}.$$
This holds if $h\prec 1$, since then $h^\dagger<0$,  while $\frac{y}{\Log a}>0$. Suppose $h\asymp 1$.
Then $v(h^\dagger) = v(h') > \Psi$, in particular $\frac{y}{\Log a} = (-1/r)\, t^\dagger \succ h^\dagger$.
Taking $0<r<1$ we get $0<2-h(r+1)$ and $2-h(r+1) \asymp 1$, so $2h^\dagger \prec \big(2-h(r+1)\big) \frac{y}{\Log a}$. The desired inequality follows, as its right-hand side is positive.
Hence for any rational $r$ with $0<r<1$, we obtain $f>\Upl(K)$ such that $\omega(z)<\omega(z^*)$ for all $z^*\in K$ with $z<z^*\leq f$. Thus $\omega$ is indeed strictly increasing on $\Upl(K)$. 

For $r=1$, we obtain an element $f=z+\frac{2y}{\Log a}>\Upl(K)$, and equations~\eqref{eq:omega1} and~\eqref{eq:omega2} with $h=1$ then yield that $\omega(z)=\omega(f)$. Conversely, let $f\in\Upd(K)$.
We shall construct $a\in K^{>C}$ such that $f=z+\frac{2y}{\Log a}$ for $y=a^\dagger$ and $z=-a^{\dagger\dagger}$. This will finish the proof of the second claim. 
Take $g\in K^<$ with $f=-g^\dagger$. Then
 by Lemma~\ref{lem:Upl, 2} we can take $t\prec 1$ in $K$ with $g=t'$. Then 
$t>0$, hence $(1/t)'=a^\dagger$ where $a\in K^{>C}$. By readjusting our choice of
$\Log a$ we arrange $t= \frac{1}{\Log a}$. Then, with  $y=a^\dagger$ and $z=-a^{\dagger\dagger}$ we have $z\in\Upl(K)$ and
$f=z+\frac{2y}{\Log a}$, as desired.
\end{proof}

\noindent
Recall from Section~\ref{sec:secondorder} that 
$$\Upo(K):=\{f\in K:\ \text{$4y''+fy=0$ for some $y\in K^\times$}\}.$$
In view of that section and Corollary~\ref{cor:Upl 2} this gives for Liouville closed $K$:
$$ \Upo(K)\ =\ \omega(K)\ =\ \omega\big(\Upl(K)\big)\ =\ \omega\big(\Upd(K)\big).$$

\noindent
Since every pre-$H$-field can be embedded into a Liouville closed $H$-field, 
the function~$\omega$ is strictly increasing on $\Upl(K)^\downarrow$ even without assuming that $K$ is Liouville closed.
Together with Lemma~\ref{lem:uplrho strictly increasing}, we therefore obtain:

\begin{cor}\label{cor:omega for H-fields, 1} If $K$ has asymptotic 
integration, then the sequence $(\upo_\rho)=\big(\omega(\upl_\rho)\big)$ is strictly increasing and cofinal in $\omega\big(\Upl(K)\big)$.
\end{cor}

\noindent
The following is a consequence of Lemma~\ref{pcrho}:
% and the previous corollary:

\begin{cor}\label{cor:omega for H-fields, 2} Let $f\in K$ and suppose
$f\ge \omega(z)+cy^2$ for some $y\in\Upg(K)$ and some constant 
$c\in C^{>}$, with $z=-y^\dagger$. Then $f>\omega\big(\Upl(K)\big)$.
\end{cor}
\begin{proof} By passing to an extension if necessary we can arrange that
$K$ is Liouville closed.
Let $y\in\Upg(K)$ and $c\in C^>$ be such that
$f\ge \omega(z)+cy^2$, with $z=-y^\dagger$. Take $a\in K^{>C}$ such that $y=a^\dagger$.
We may choose our logarithmic sequence $(\ell_\rho)$ for~$K$ beginning with 
$\ell_0=a$; then $\upg_0=\ell_0^\dagger=y$ and $\upl_0=-\upg_0^\dagger=z$.
For $\rho>0$ we have $\upg_1^2\sim \omega(\upl_\rho)-\omega(\upl_0)>0$  by 
Lemma~\ref{pcrho}.
% and Corollary~\ref{cor:omega for H-fields, 1}.
Hence $\omega(\upl_\rho)-\omega(\upl_0) < c\upg_0^2$, so
$$f\ \ge\ \omega(z)+cy^2\ =\ \big(\omega(\upl_0)-\omega(\upl_\rho)+c\upg_0^2\big)+\omega(\upl_\rho)\ >\ \omega(\upl_\rho).$$
Since this holds for all $\rho>0$, we conclude that $f>\omega\big(\Upl(K)\big)$.
\end{proof}

\begin{cor} Assume $K$ is Liouville closed, and
 $f,g\in K^>$ are such that $-\frac{1}{2}f^\dagger \in \Upl(K)$ and $f\preceq g$. Then $\omega(-\textstyle\frac{1}{2}f^\dagger)+g
> \omega\big(\Upl(K)\big)$.
\end{cor}
\begin{proof} 
Take $y\in K^>$ with $y^2=f$. Then $z:=-y^\dagger = -\frac{1}{2}f^\dagger \in \Upl(K)$, so $y> \I(K)$ by Lemma~\ref{lem:Upl, 1},
hence $y\in\Upg(K)$ by Lemma~\ref{lem:Upg, Liouville closed}. Thus, with
$c\in C^{>}$ such that $cf\le g$, we have
$\omega(-\frac{1}{2}f^\dagger)+g \ge \omega(z)+cy^2>\omega\big(\Upl(K)\big)$ by Corollary~\ref{cor:omega for H-fields, 2}.
\end{proof}

\begin{cor}\label{cor:omega for H-fields, 3} Suppose $K$ is Liouville closed and has small derivation. Then  $\upo_\rho>0$ eventually.
\end{cor}
\begin{proof} Take
$x\in K$ with $x'=1$. Then $x\succ 1$, so $1/x=-x^{\dagger\dagger}\in\Upl(K)$, and $\omega(1/x)=1/x^2>0$.
Now use Corollary~\ref{cor:omega for H-fields, 1}.
\end{proof}

\begin{cor}\label{cor:linconstcoeff, m(A)=dim ker A}
Suppose $K$ is Liouville closed with small derivation. Then for each $A\in C[\der]^{\neq}$ we have
$\operatorname{m}(A)=\dim_C \ker A$.
\end{cor}
\begin{proof}
Lemma~\ref{lem:upl upo infinitesimal} and Corollary~\ref{cor:omega for H-fields, 3} give $\omega(K)\cap C \subseteq C^{\leq}$. 
It remains to use Lemma~\ref{lem:linconstcoeff, omega}.
\end{proof}

\begin{cor}\label{cor:linconstcoeff, Schwarz closed}
Suppose $K$ is Liouville closed with small derivation, and let $P(Y)\in C[Y]$. Set $A:= P(\der)\in C[\der]$. Then 
\begin{align*}
\exists y\in K^\times \big[ y \prec 1 \ \& \ A(y)=0 \big] & \quad\Longleftrightarrow\quad
\exists c\in C \big[ c < 0 \ \&\ P(c) = 0 \big],\\
\exists y\in K^\times \big[ y \preceq 1 \ \& \ A(y)=0 \big] & \quad\Longleftrightarrow\quad
\exists c\in C \big[ c \le 0 \ \&\ P(c) = 0 \big],\\
\exists y\in K^\times \big[ y \succeq 1 \ \& \ A(y)=0 \big] & \quad\Longleftrightarrow\quad
\exists c\in C \big[ c \ge 0 \ \&\ P(c) = 0 \big].
\end{align*} 
\end{cor}
\begin{proof} Take $x\in K$ with $x'=1$. As in
the subsection on linear differential equations with constant coefficients in Section~\ref{Linear Differential Operators} we pick for each $c\in C$ an element of $K^\times$, denoted by $\ex^{cx}$, such that
$(\ex^{cx})'=c\ex^{cx}$. Thus $\ex^{0x}\in C^\times$, and 
it is routine to show that if $c\in C^{<}$, then 
$\ex^{cx}\prec x^{-n}$ for all $n$, and if $c\in C^{>}$, then 
$\ex^{cx}\succ x^n$ for all~$n$. It follows easily
that $\ex^{c_1x}x^{k_1}\not\asymp\ex^{c_2x}x^{k_2}$
for all distinct $(c_1, k_1), (c_2, k_2)\in C\times \Z$.
It now remains to use Proposition~\ref{prop:linconstcoeff} and Corollary~\ref{cor:linconstcoeff, m(A)=dim ker A}.
\end{proof}

\noindent
Here is a universal property similar to Lemma~\ref{lem:invariance gp}:

\begin{lemma}\label{uprupo} Suppose $K$ is Liouville closed. Then for all $g,h\in \I(K)$,
$$ \omega(K)^{\downarrow} + gh\ \subseteq\ \omega(K)^{\downarrow}.$$
\end{lemma}
\begin{proof} Let $a\in \omega(K)^{\downarrow}$ and $g,h\in \I(K)$. Take $\rho$ so large that $a < \upo_{\rho} =\omega(\upl_{\rho})$. Then 
$g,h\prec \upg_{\rho+2}$, so
$a+gh < \upo_{\rho} + \upg_{\rho+2}^2 < \upo_{\rho+1}$, and thus
$a+gh \in \omega(K)^{\downarrow}$.
\end{proof} 

\begin{cor}\label{cor:uprupo}
Suppose $K$ is Liouville closed, and $f\in K^\times$, $-\frac{1}{2}f^\dagger\notin\Lambda(K)$.
Then $\omega(K)^\downarrow+f\subseteq\omega(K)^\downarrow$.
\end{cor}
\begin{proof}
This is clear if $f<0$, so assume $f>0$, and take $y\in K^>$ with~${y^2=f}$. Then 
$-y^\dagger = -\frac{1}{2}f^\dagger \notin \Upl(K)$, so $y\in \I(K)$
by Lemma~\ref{lem:Upl, 1}. Now use Lemma~\ref{uprupo}. 
\end{proof}

\noindent
We now consider also the function
$$y\mapsto\sigma(y) = \omega(z)+y^2  \colon\  K^\times \to K
\qquad\text{where $z=-y^\dagger$.}$$
Note that $\sigma(y)=\sigma(-y)$ for all $y\in K^\times$. By Corollary~\ref{cor:omega for H-fields, 2},
$\omega\big(\Upl(K)\big) <\sigma\big(\Upg(K)\big)$,
and by embedding $K$ into a Liouville closed $H$-field, this gives
$$  \omega(K)\ <\ \sigma\big(\Upg(K)\big).$$

\begin{lemma}\label{lem:sigma increasing}
The restriction of $\sigma$ to $\Upg(K)$ is strictly increasing.
\end{lemma}
\begin{proof}
Let $y\in \Upg(K)$, $z=-y^\dagger$ and $a>1$ in $K$; we shall derive
$\sigma(ay)>\sigma(y)$. We have $-(ay)^\dagger = z+b$ with $b:=-a^\dagger$,
so
$$\sigma(ay) - \sigma(y)\ =\ \omega(z+b)-\omega(z) + (a^2-1)y^2.$$
We have $a\in C$ iff $b=0$, and in this case
$\sigma(ay)=\sigma(y)+(a^2-1)y^2>\sigma(y)$ as required.
Suppose $a\notin C$. Then by Lemma~\ref{varrholemma, claim, 1},
$$\sigma(ay)-\sigma(y)\ =\ -b\big(2(b/y)^\dagger +b \big)+ (a^2-1)y^2.$$
%We claim that $\sigma(ay)-\sigma(y) \sim (a^2-1)y^2$ and 
%hence $\sigma(ay)>\sigma(y)$ as required.
We now distinguish three cases:

\case[1]{$a \asymp 1$.}
Then $\Psi<va'=vb$, in particular $y\succ b$ and $(b/y)^\dagger\succ b$,
hence $-b\big(2(b/y)^\dagger + b\big) \sim -2b(b/y)^\dagger$.
From $vy \in \Psi < v(b/y)'=v(b/y)+v(b/y)^\dagger$, we obtain
$v(y^2)< vb+v(b/y)^\dagger$. If $a\not\sim 1$, this gives
$-b\big(2(b/y)^\dagger + b\big)\prec y^2\asymp (a^2-1)y^2$, so
$\sigma(ay)-\sigma(y) \sim (a^2-1)y^2>0$.
If $a\sim 1$, then $a=1+\varepsilon$ with $0< \varepsilon\prec 1$, so
$b=-\varepsilon'/(1+\varepsilon) >0$, hence $0 < b/y \prec 1$, which gives
 $-b\big(2(b/y)^\dagger + b\big)\sim -2b(b/y)^\dagger>0$, and thus
$\sigma(ay)-\sigma(y)>0$. 

\case[2]{$a\succ 1$ and $y\prec b$.}
Then $vy\in \Psi <  v(y/b)'=v(y/b)+v(y/b)^\dagger$, so 
$vb < v(y/b)^\dagger=v(b/y)^\dagger$ and thus
$-b\big(2(b/y)^\dagger +b \big) \sim -b^2$.
Also, $v(1/a)'>\Psi$, in particular 
$vb=v(1/a)^\dagger=va+v(1/a)'>v(ay)$ and hence
$b^2\prec (ay)^2\sim (a^2-1)y^2$.
This gives $\sigma(ay)-\sigma(y) \sim (a^2-1)y^2>0$.

\case[3]{$a\succ 1$ and $y\succeq b$.}
Then $(b/y)^\dagger \prec y$ by Lemma~\ref{BasicProperties}(i), so
$$-b\big(2(b/y)^\dagger +b \big)\ \preceq\ by\ \preceq\ y^2\ \prec\ (a^2-1)y^2,$$
and thus once again, $\sigma(ay)-\sigma(y) \sim (a^2-1)y^2>0$.
\end{proof}

\noindent
As any pre-$H$-field
embeds into some Liouville closed $H$-field, 
Lemmas~\ref{lem:Upg, Liouville closed} and \ref{lem:sigma increasing}
give that $\sigma$ is strictly increasing on 
$\Upg(K)^{\uparrow}$,
even without assuming that~$K$ itself is Liouville closed.
By earlier comments, if $K$ has asymptotic integration, then
the sequence~$(\upg_{\rho})$ in $\Upg(K)$ is
strictly decreasing and coinitial in $\Upg(K)$.

{\sloppy
\begin{cor}\label{cor:omega-freeness and cut} Assume that $K$ has asymptotic 
integration. Then
the se\-quence~$\big(\sigma(\upg_\rho)\big)=(\upo_\rho+\upg_\rho^2)$ is strictly decreasing and coinitial in $\sigma\big(\Upg(K)\big)$. 
Thus for any~$\upo$ in any pre-$H$-field extension of $K$ we have the equivalence
$$ \upo_\rho \leadsto \upo\ \Longleftrightarrow\  
\omega\big(\Upl(K)\big)<\upo<\sigma\big(\Upg(K)\big).$$
\end{cor}
}

\noindent
As to the behavior of $\sigma$ on the complement $K^>\setminus\Upg(K)$, we note:

\begin{lemma}\label{lem:sigma, Liouville closed}
Suppose $K$ is Liouville closed. Then 
\begin{align*} &\sigma\big(K^>\setminus\Upg(K)\big)\ \subseteq\ \omega(K)^\downarrow, \text{ so} \\
\sigma(K^\times)\ =\ &\sigma\big(K^>\setminus\Upg(K)\big)\cup\sigma\big(\Upg(K)\big)\quad\text{ with }
\sigma\big(K^>\setminus\Upg(K)\big) < \sigma\big(\Upg(K)\big).
\end{align*}
\end{lemma}

\begin{proof}
Let $s\in K^>\setminus\Upg(K)$. Using Lemma~\ref{lem:Upg, Liouville closed}, take $a\in K^{\succ 1}$
with $s=(1/a)'$ and set $y:=a^\dagger\in\Upg(K)$, $z:=-y^\dagger\in\Upl(K)$. 
Then $s=-a'/a^2=-y/a$, hence $-s^\dagger=-y^\dagger+a^\dagger=z+y$ and thus, using
Corollary~\ref{cor:formula for omega},
\equationqed{
\sigma(s)\ =\ \omega(z+y)+(y/a)^2\ =\ \omega(z)-y^2+(y/a)^2\ <\ \omega(z).}
\end{proof}

\begin{cor}\label{lisplofree}
Suppose $K$ is Liouville closed, and for every $a\in K$ 
the operator $\der^2-a\in K[\der]$ splits over $K[\imag]$.
Then $K\setminus\omega(K)^\downarrow=\sigma\big(\Upg(K)\big)$. In particular, $\sigma\big(\Upg(K)\big)$ is upward closed and $K$ is $\upo$-free.
\end{cor}
\begin{proof} From \eqref{eq:A splits over K[i]} we get
$K=\omega(K)\cup\sigma(K^\times)$, so
$$K\setminus\omega(K)^\downarrow\ =\ \sigma(K^\times)\setminus\omega(K)^\downarrow\ 
=\ \sigma\big(\Upg(K)\big)$$
by Lemma~\ref{lem:sigma, Liouville closed}. The rest now follows, using  Corollary~\ref{cor:omega-freeness and cut}.
\end{proof}

\noindent
In Section~\ref{cordone} we consider additional restrictions on $K$ ensuring that
$\omega(K)$ is downward closed. 
Figure~\ref{fig:omega} is a sketch of the functions $\omega$ on
$\Uplambda(\T)$  and $\sigma$ on~$\Upgamma(\T)$.

\begin{figure}[h!]
\begin{center}
\begin{tikzpicture}[scale=0.74]
\draw (4.5,0) -- (8.5,0);
\draw (9.5,0) -- (12,0); 
\draw[densely dotted] (0,2) -- (12,2);
\draw[densely dotted] (-5,2) -- (-3.75,2);
\draw[densely dotted] (4,7) -- (4,0);
\draw[densely dotted] (4,-1) -- (4,-5);
\draw[densely dotted] (9,7) -- (9,0);
\draw[densely dotted] (9,-1) -- (9,-5);
\draw[densely dotted] (-4,7) -- (-4,0);
\draw[densely dotted] (-4,-1) -- (-4,-5);

\draw (0,1.5) -- (0,-3);
\draw (0,-4.5) -- (0,-5);
\draw[densely dotted] (0,1.5) -- (0,1.75);
\draw[densely dotted] (0,2.25) -- (0,2.5);

\draw[densely dotted] (3.5,0) -- (3.75,0);
\draw[densely dotted] (4.25,0) -- (4.5,0);
\draw[dotted] (3.75,0) -- (4.25,0);
\draw (4.0,-0.5) node {\mbox{\fontsize{9}{0}\selectfont $\frac{1}{\ell_0\ell_1\cdots}$}};
\fill[white] (4,0) circle (0.2em);
\draw (4,0) circle (0.2em);

\draw (-5,0) -- (-4.5,0);
\draw (-3.5,0) -- (3.5,0);
\draw (-4.0,-0.5) node {\mbox{\fontsize{9}{0}\selectfont $-\frac{1}{\ell_0\ell_1\cdots}$}};
\draw[densely dotted] (-3.5,0) -- (-3.75,0);
\draw[densely dotted] (-4.25,0) -- (-4.5,0);
\draw[dotted] (-3.75,0) -- (-4.25,0);
\fill[white] (-4,0) circle (0.2em);
\draw (-4,0) circle (0.2em);

\draw[densely dotted] (8.5,0) -- (8.75,0);
\draw[densely dotted] (9.25,0) -- (9.5,0);
\draw[dotted] (8.75,0) -- (9.25,0);
\draw (9.0,-0.5) node {\mbox{\fontsize{9}{0}\selectfont $\frac{1}{\ell_0}+\frac{1}{\ell_0\ell_1}+\cdots$}};
\fill[white] (9,0) circle (0.2em);
\draw (9,0) circle (0.2em);

\draw[dotted] (0,2.25) -- (0,1.75);
\draw (-2,2) node {\mbox{\fontsize{9}{0}\selectfont $\frac{1}{\ell_0^2}+\frac{1}{(\ell_0\ell_1)^2}+\cdots$}};
\draw (0,2.5) -- (0,7);
\fill[white] (0,2) circle (0.2em);
\draw (0,2) circle (0.2em);

\draw (8,0.5) node {$\Uplambda(\T)$};
\draw[-latex, gray, thick] (7.3,0.55) -- (6,0.55);
  
\draw (10,-1.25) node {$\Updelta(\T)$};
\draw[-latex, gray, thick] (10.75,-1.25) -- (12,-1.25);

\draw [decoration={brace,amplitude=0.7em,mirror},decorate,thick,gray] (-3.75,-3.3) -- (3.75,-3.3);
\draw (0,-4) node {$\I(\T)$};

\draw (5,-1.25) node {$\Upgamma(\T)$};
\draw[-latex, gray, thick] (5.7,-1.25) -- (7,-1.25);

\draw (0.5,6) node {$\T$};
\draw (11,0.5) node {$\T$};
\draw[-latex, gray, thick] (11.5,0.5) -- (12,0.5);
\draw[-latex, gray, thick] (0.5,6.5) -- (0.5,7);

\draw (0.5,-0.5) node {$0$};

\draw[densely dotted, domain=4.2:4.5] plot (\x, { ( 0.1*(\x-4) * (\x-4) + 2 ) });
\draw[densely dotted, domain=8.5:9.5] plot (\x, {  ( 0.1*(\x-4) * (\x-4) + 2 ) });
\draw[line width=0.1em, domain=9.5:11] plot (\x, {  ( 0.1*(\x-4) * (\x-4) + 2 ) });
\draw[line width=0.1em, domain=4.5:8.5] plot (\x, {  ( 0.1*(\x-4) * (\x-4) + 2 ) });
\draw (11,6) node {$\sigma$};
\fill[white] (4,2) circle (0.4em);
\draw (4,2) circle (0.4em);

\draw[densely dotted, domain=3.5:4.5] plot (\x, { ( -0.025*(\x-9) * (\x-9) + 2 ) });
\draw[densely dotted, domain=8.5:9] plot (\x, {  ( -0.025*(\x-9) * (\x-9) + 2 ) });
\draw[line width=0.1em, domain=-3.5:3.5] plot (\x, {  ( -0.025*(\x-9) * (\x-9) + 2 ) });
\draw[densely dotted, domain=-4.5:-3.5] plot (\x, {  ( -0.025*(\x-9) * (\x-9) + 2 ) });
\draw[line width=0.1em, domain=-5:-4.5] plot (\x, {  ( -0.025*(\x-9) * (\x-9) + 2 ) });
\draw[line width=0.1em, domain=4.5:8.5] plot (\x, {  ( -0.025*(\x-9) * (\x-9) + 2 ) });
\draw (-2.5,-2) node {$\ome$};
\fill[white] (9,2) circle (0.4em);
\draw (9,2) circle (0.4em);

\end{tikzpicture}
\end{center}
\caption{The functions $\omega$ and $\sigma$ on $\T$.}\label{fig:omega}
\end{figure}

\begin{cor}\label{eqschwarz}
The following are equivalent for a Liouville closed $H$-field $K$:
\begin{enumerate}
\item[\textup{(i)}] $K=\omega\big(\Upl(K)\big)\cup 
       \sigma\big(\Upg(K)\big)$;
\item[\textup{(ii)}] $K$ is $\upo$-free, $\omega\big(\Upl(K)\big)$ is downward closed, 
and $\sigma\big(\Upg(K)\big)$ is upward closed;
\item[\textup{(iii)}] for every $a\in K$ the operator $\der^2-a\in K[\der]$ splits over $K[\imag]$, and $\omega(K)$ is downward closed.
\end{enumerate}  
\end{cor}
\begin{proof} From  \eqref{eq:A splits over K[i]} and $\omega(K)< \sigma\big(\Upg(K)\big)$ we get (i)~$\Rightarrow$~(iii).
The implication (iii)~$\Rightarrow$~(ii) follows from Corollary~\ref{lisplofree}, 
and (ii)~$\Rightarrow$~(i) from Corollary~\ref{cor:omega-freeness and cut}.
\end{proof}

\noindent
We say that a pre-$H$-field $K$ is {\bf Schwarz closed} if $K$ is Liouville closed and satisfies the equivalent conditions in the previous corollary.  This terminology is motivated by the role of the functions $\omega$ and $\sigma$ in the Schwarzian derivative; see Section~\ref{sec:secondorder}. 

\label{p:Schwarz closed}
\index{Schwarz!closed}
\index{H-field@$H$-field!Schwarz closed}
\index{closed!Schwarz}

\begin{lemma}\label{lem:Schwarzian closed under compconj}
If $K$ is Schwarz closed, then so is $K^\phi$ for all $\phi\in K^{>}$.
\end{lemma}
\begin{proof} Let $K$ be Schwarz closed and $\phi\in K^{>}$. Then 
$K^\phi$ is Liouville closed, and with $\derdelta=\phi^{-1}\der$,  
every operator $\derdelta^2-a\in K^\phi[\derdelta]$ ($a\in K^\phi$) splits over $K^\phi[\imag]=K[\imag]^{\phi}$.
Let $\omega^\phi$ be the map $z\mapsto -\big(2\derdelta(z)+z^2\big)\colon K^\phi\to K^\phi$. We saw in
Section~\ref{sec:behupo} that it
plays the role of $\omega$ in the differential field $K^\phi$. The computations there give
$$\omega^\phi(K^\phi)\ =\ \phi^{-2}\big(\omega(K) - \omega(-\phi^\dagger)\big), $$
hence $\omega^\phi(K^\phi)$ is downward closed.
\end{proof}

\begin{lemma}\label{intersectschwarz} Let $K$ and $L$ be Schwarz closed $H$-subfields
of an $H$-field $M$ such that $C_L=C_M$. Then $K\cap L$ is a Schwarz
closed $H$-subfield of $M$.
\end{lemma}
\begin{proof} First,  $K\cap L$ is a Liouville closed $H$-subfield of $M$ by Lemma~\ref{intersectliou}. It remains to show that
$\omega\big(\Upl(K)\big)\cap L\subseteq \omega\big(\Upl(K\cap L)\big)$ and
$\sigma\big(\Upg(K)\big)\cap L \subseteq 
\sigma\big(\Upg(K\cap L)\big)$.
Let $f\in \omega\big(\Upl(K)\big)\cap L$. Then $f=\omega(z_1)$ with 
$z_1\in \Upl(K)$. Since $L$ is Schwarz closed, 
$$ f\in \omega\big(\Upl(M)\big)\cap L\ \subseteq\ \omega\big(\Upl(L)\big),$$ so $f=\omega(z_2)$ with
$z_2\in \Upl(L)$. Now $\omega(z_1)=\omega(z_2)$ with 
$z_1, z_2\in \Upl(M)$, so 
$$z\ :=\ z_1\ =\ z_2\in K\cap L\cap \Upl(K)\ =\ \Upl(K\cap L),$$ 
where the last step uses Corollary~\ref{prupl}. Hence $f=\omega(z)\in \omega\big(\Upl(K\cap L)\big)$,
as required. The other inclusion follows in the same way.
\end{proof}

\begin{lemma}\label{lem:intersect Schwarz closed H-fields}
Suppose $K$ is Schwarz closed and $(K_i)_{i\in I}$, $I\ne \emptyset$, is a family of Schwarz closed $H$-subfields of $K$, each with the same constant field as $K$. Then $F:=\bigcap_{i\in I} K_i$ is a Schwarz closed $H$-subfield of $K$.
\end{lemma}
\begin{proof}
By Lemma~\ref{intersectlioumany}, $F$ is  a Liouville closed $H$-subfield of $K$.
We need to show that $F=\omega\bigl(\Upl(F)\big)\cup\sigma\big(\Upg(F)\big)$, and for this
it is enough to show that $F\cap\omega\bigl(\Upl(K)\big) \subseteq \omega\bigl(\Upl(F)\big)$ and
$F\cap \sigma\big(\Upg(K)\big)\subseteq \sigma\big(\Upg(F)\big)$.
Let $f\in F\cap\omega\big(\Upl(K)\big)$; then $f\notin\sigma\big(\Upg(K)\big)$. Hence for each $i\in I$
we have $f\in K_i\setminus \sigma\big(\Upg(K_i)\big)=\omega\big(\Upl(K_i)\big)$, so
we may take $z_i\in \Upl(K_i)$ with $\omega(z_i)=f$.
Since the map $z\mapsto\omega(z)\colon\Upl(K)\to K$ is injective,
we have $z_i=z_j$ for all $i,j\in I$. Denoting the common value of $z_i$ ($i\in I$) by~$z$, we obtain $z\in F\cap\Upl(K)=\Upl(F)$ by Corollary~\ref{prupl}, 
hence $f=\omega(z)\in\omega\bigl(\Upl(F)\big)$ as required. Similarly one shows that
$F\cap \sigma\big(\Upg(K)\big)\subseteq \sigma\big(\Upg(F)\big)$.
\end{proof}

%% file: mt-12.tex
\chapter{Triangular Automorphisms}\label{ch:triangular automorphisms}

\noindent
{\em Throughout this chapter $K$ is a commutative ring containing $\Q$ as a subring}.

\medskip\noindent
When $K$ carries also a derivation and $\phi\in K^\times$, then we saw 
in Section~\ref{Compositional Conjugation} that compositional conjugation 
by $\phi$ induces an automorphism of the $K$-algebra
$$K[Y_0,Y_1, Y_2, \dots],$$ where $Y_n:= Y^{(n)}$ and $Y$ is a differential indeterminate. These automorphisms are triangular (as defined below)
and this facilitates their finer analysis by Lie techniques, as we show in
the present chapter. An endomorphism $\sigma$ of this $K$-algebra
is uniquely determined by the sequence of images $\sigma(Y_0),\ \sigma(Y_1),\ \sigma(Y_2),\dots$, and is said to
be {\em triangular\/} if for each $n$ there are $\sigma_{0n}, \dots, 
\sigma_{nn}\in K$
with
$$\sigma(Y_n)\ =\  \sigma_{0n} Y_0 + \sigma_{1n}Y_1 +\cdots + \sigma_{nn}Y_n.$$
In Sections~\ref{triangularlinear}--\ref{sec:Derivations on polynomials rings} of this chapter we introduce a formalism to
analyze triangular automorphisms of such a polynomial algebra by means of 
their \textit{logarithms}\/, the triangular derivations.  
In Section~\ref{appdifpol} we apply this to   
compositional conjugation in differential polynomial rings. From this chapter,
only Corollaries~\ref{cor:companion of Stirling} and~\ref{cor:companion of Stirling, variant} as well as Section~\ref{appdifpol} will 
be needed later.

\section{Filtered Modules and Algebras}\label{sec:Filtered and Graded Algebras}

\noindent
In this preliminary section we collect some definitions and simple facts 
about filtered modules, filtered algebras, and graded algebras.

\subsection*{Filtered modules}
Let $V$ be a $K$-module. A {\bf filtration} of $V$ is a family $(V^i)_{i\in\Z}$ of $K$-submodules of $V$ such that 
\begin{enumerate}
\item $\bigcup_{i} V^i=V$, 
\item $\bigcap_{i} V^i=\{0\}$, and
\item $V^i\supseteq V^{i+1}$ for all $i$.
\end{enumerate}
We call a filtration $(V^i)$ of $V$ {\bf nonnegative} if $V^0=V$; equivalently,
$V^i=V$ for every~$i\leq 0$. In this case we also call 
$(V^n)$ a nonnegative filtration. The {\bf trivial} filtration of $V$ is the nonnegative filtration $(V^n)$ of $V$ with $V^n=\{0\}$ for $n>0$.
If~$(V^i)$ is a filtration of $V$ and $W$ is a submodule of $V$, then $(V^i\cap W)$ is a filtration of $W$,
called the filtration of $W$ {\bf induced} by the filtration $(V^i)$.

\index{filtration!module}
\index{filtration!module!trivial}
\index{filtration!module!nonnegative}
\index{filtration!module!induced on a submodule}

\index{filtered!module}
\index{filtered!submodule}
\index{module!filtered}

Let $(V^i)$ be a filtration of $V$. To keep notations simple we denote 
$V$ with this filtration also by $V$; the combined object is called a  
{\bf filtered $K$-module}, and a 
submodule of $V$ equipped with the filtration induced by $(V^i)$ is called
a {\bf filtered submodule} of~$V$.
The filtration $(V^i)$ is a fundamental system of neighborhoods of~$0$ for a unique
topology on $V$ making the additive group of $V$ a topological group. For $x\in V$, let 
$|x|\in \R^{\ge 0}$ be given by
$$|x|:= 2^{-i} \text{ if $x\in V^i\setminus V^{i+1}$,} \quad |0|:= 0,$$
so $|x|\le 2^{-i}\ \Leftrightarrow\ x\in V^i$.
Thus $|\cdot|$ is an ultrametric norm on $V$ (and a bit more): 
for all $x,y\in V$ and $a\in K$ and nonzero $k\in \Z$, 
$$|x|=0\ \Longleftrightarrow\ x=0,\quad \abs{ax}\le \abs{x}, \quad \abs{kx}= \abs{x}, \quad |x+y|\le \max(|x|,|y|).$$ 
The topology defined above is induced by the metric $(x,y)\mapsto |x-y|$ on $V$.
We call the filtered $K$-module $V$ {\bf complete} if $V$ is complete
with respect to this metric. \index{filtered!module!complete}\index{complete!filtered module} A family 
$(v_{\lambda})_{\lambda\in \Lambda}$ in $V$
is said to be {\bf summable} \index{family!summable}
\index{summable!family of elements} if there is a (necessarily unique) $v\in V$
with the following property: for every $\epsilon>0$ there is a finite $I(\epsilon)\subseteq \Lambda$ such that $\left|v-\sum_{\lambda\in I}v_{\lambda}\right|< \epsilon$ for all finite $I\subseteq \Lambda$ with $I\supseteq I(\epsilon)$; in that case $v_{\lambda}\ne 0$ for only
countably many $\lambda\in \Lambda$, we call this $v$ the {\bf sum} \index{family!sum}\index{sum!family} of the
family $(v_{\lambda})$, and denote it by $\sum_{\lambda\in\Lambda} v_\lambda$, or by
$\sum v_{\lambda}$ if the index set $\Lambda$ is understood from the context.
If~$(u_{\lambda})$ and~$(v_{\lambda})$ are summable 
families in $V$ with the same index set and $a\in K$, then 
$(u_{\lambda}+v_{\lambda})$ and $(av_{\lambda})$ are summable, with
$$\sum u_{\lambda}+v_{\lambda}\ =\ \sum u_{\lambda} + \sum v_{\lambda}, \qquad 
\sum av_{\lambda}\ =\ a\sum v_{\lambda}.$$
Suppose $V$ is complete and $(v_{\lambda})_{\lambda\in \Lambda}$ is a family in $V$. 
Then $(v_{\lambda})$ is summable iff for each  $\epsilon>0$ we have 
$|v_\lambda| < \epsilon$ for all but finitely many $\lambda$. In particular, a sequence~$(v_n)$ in~$V$ is summable iff $v_n \to 0$ as $n\to \infty$. Suppose
also that $\Lambda$ is the disjoint union of $\Lambda(1)$ and $\Lambda(2)$; then
$(v_{\lambda})_{\lambda\in \Lambda}$ is summable iff $(v_{\lambda})_{\lambda\in \Lambda(1)}$
and  $(v_{\lambda})_{\lambda\in \Lambda(2)}$ are summable; moreover, if $(v_{\lambda})_{\lambda\in \Lambda}$ is summable, then 
$$\sum_{\lambda\in \Lambda} v_{\lambda}\ =\ \sum_{\lambda\in \Lambda(1)} v_{\lambda}  +  \sum_{\lambda\in \Lambda(2)} v_{\lambda}.$$
We leave it to the reader to state and prove a similar statement for $\Lambda=\bigcup_{i\in I} \Lambda(i)$ where $I$
is any index set with $\Lambda(i)\cap \Lambda(j)=\emptyset$
for all $i\ne j$.

\medskip\noindent
Let $W$ be a filtered $K$-module, with filtration $(W^j)$, and let $\Phi\colon V\to W$ be $K$-linear. 
Then $\Phi$ is continuous iff for each $j$ there is an $i$ such that $\Phi(V^i)\subseteq W^j$. Thus if~$\Phi$ is continuous and $V$ and $W$ are complete, then $\Phi$ is ``strongly additive'' in the following sense: if $(v_\lambda)$ is a summable family in $V$, then $\big(\Phi(v_\lambda)\big)$ is summable in~$W$, and $\Phi\left(\sum v_\lambda\right)=\sum \Phi(v_\lambda)$.
For $d\in \Z$ we say that $\Phi$ is {\bf of rank $d$} if $\Phi(V^i)\subseteq W^{i+d}$ for each $i\in\Z$, and we call $\Phi$ {\bf ranked} if $\Phi$ is of 
rank $d$ for some $d$. (For example, if~$V$ is a filtered submodule of $W$, then the natural inclusion $V\to W$ is of rank $0$.) If $\Phi$ is ranked, then $\Phi$ is continuous. If $\Phi$ is of rank $d\ge 0$, then $\Phi$ is of rank $0$.
If $\Phi\colon V\to V$ is of rank $d$, then $\Phi^n$  is of rank $nd$; here and below, $\Phi^n$ denotes the $n$-fold composition of $\Phi$ with itself: $\Phi^0=\id$, $\Phi^{n+1}=\Phi\circ\Phi^n$ for each $n$.

\index{summable!family of endomorphisms}
\index{rank!linear map between filtered modules}
\index{ranked}
\index{topologically nilpotent}
\index{nilpotent!topologically}
\index{weakly!nilpotent}
\index{nilpotent!weakly}

\nomenclature[A]{$\End(V)$}{algebra of endomorphisms of a module $V$}

\medskip\noindent
{\em In the rest of this subsection we assume that $V$ is complete.}\/
We have the $K$-al\-ge\-bra~$\End(V)$ of endomorphisms of $V$, whose elements are the $K$-linear maps~${V\to V}$, with multiplication given by composition of maps. The endomorphisms of~$V$ of rank~$0$ form a subalgebra of $\End(V)$. 
A family $(\Phi_\lambda)$ of endomorphisms of $V$ is said to be 
{\bf summable} if for each $v\in V$ the family $\big(\Phi_\lambda(v)\big)$ is summable in $V$; then
the sum of~$(\Phi_\lambda)$ \index{sum!family} is the endomorphism $\sum \Phi_\lambda$ of $V$ given by 
$$\left(\sum  \Phi_{\lambda}\right)(v) := \sum\Phi_{\lambda}(v).$$
Below $\Phi$ ranges over $\End(V)$. We say that $\Phi$ is {\bf topologically nilpotent} if for each~$j$ we have $\Phi^n(V)\subseteq V^j$ for all sufficiently large $n$. Note that if the filtration of~$V$ is nonnegative and $\Phi$ is 
of positive rank, then $\Phi$ is topologically nilpotent. We say that
$\Phi$ is {\bf weakly nilpotent} if for each $v\in V$ we have $\Phi^n(v)\to 0$ as~$n\to \infty$. For example, if $\Phi$ is of positive rank or topologically
nilpotent, then $\Phi$ is weakly nilpotent. If $\Phi$ is weakly nilpotent, then for each formal power series $f=\sum f_n z^n\in K[[z]]$ ($f_n\in K$ for all $n$),
the sequence $(f_n\Phi^n)$ of endomorphisms of~$V$ is summable, and we set
$$f(\Phi)\ :=\ \sum f_n\,\Phi^n \in \End(V). $$
We equip the $K$-module $K[[z]]$ with the nonnegative filtration $\big(z^n K[[z]]\big)$. Then $K[[z]]$ is complete. If $\Phi$ is weakly nilpotent and the family $(f_\lambda)$ in $K[[z]]$ is summable, then the family $\big(f_\lambda(\Phi)\big)$ in $\End(V)$ is summable, and 
$$\left(\sum f_\lambda\right)(\Phi)\ =\ \sum f_\lambda(\Phi).$$ 
In Section~\ref{triangularlinear}  we use the following.

\begin{lemma}\label{lem:linear maps into power series}
Suppose $\Phi$ is continuous and weakly nilpotent. Then:
\begin{enumerate}
\item[\textup{(i)}] the map $f\mapsto f(\Phi)\colon K[[z]]\to\End(V)$ is a morphism of $K$-algebras;  
\item[\textup{(ii)}] if $\Phi$ is of rank $0$, then $f(\Phi)\in \End(V)$ is of rank $0$, 
for all $f\in K[[z]]$;
\item[\textup{(iii)}] for $g\in zK[[z]]$, the endomorphism $g(\Phi)$ of $V$ is weakly nilpotent and
$$(f\circ g)(\Phi)\ =\ f\big(g(\Phi)\big)\ \text{ for all}\ f\in K[[z]].$$
\end{enumerate}
\end{lemma}

\begin{proof}
It is easy to check that $f\mapsto f(\Phi)$ is $K$-linear. Let $f,g\in K[[z]]$, $f=\sum_m f_m z^m$, $g=\sum_n g_n z^n$ with $f_m,g_n\in K$ for all $m$, $n$. Then  for each $v\in V$,
\begin{align*}
\big(f(\Phi)\circ g(\Phi)\big)(v)\	&=\ \sum_{m} f_m\,\Phi^m\left(\sum_{n} g_n\,\Phi^n(v)\right) \\
									&=\ \sum_{m} \,\left(\sum_{n} f_mg_n\,\Phi^{m+n}(v) \right) \\
									&=\ \sum_{k=0}^\infty \left(\sum_{m+n=k} f_mg_n\right) \,\Phi^k(v)\  =\ (f\cdot g)(\Phi)(v),
\end{align*}
where we used continuity of $\Phi^m$ for the second equality. This shows (i). 

To get (ii), just use that each $V^i$ is closed. 
For (iii), let $g=\sum_{i\geq 1} g_iz^i$ where $g_i\in K$ for all $i\geq 1$.
Note that $g(\Phi)^n=g^n(\Phi)$ by (i), and 
$$g^n\ =\ g_{1}^nz^{n}+\text{higher degree terms in $z$}.$$
Let $v\in V$ and $j\ge 1$, and take $n_1\geq 1$ such that $\Phi^n(v)\in V^j$ for all $n\geq n_1$. Then for all $n\geq n_1$ we have
$g^n(\Phi)(v)=g_{1}^n\Phi^{n}(v)+\cdots\in V^{j}$, and so $g(\Phi)^n(v)\in V^{j}$. Thus~$g(\Phi)$ is weakly nilpotent. Let also $f=\sum f_nz^n$ where $f_n\in K$ for all $n$.
Then the family $(f_n g^n)$ in $K[[z]]$ is summable, and 
$f\circ g= \sum f_ng^n$ in $K[[z]]$, so
$$(f\circ g)(\Phi)\ =\  \left(\sum f_n g^n\right)(\Phi)\ =\ \sum f_n g^n(\Phi)\ 
=\ \sum f_n g(\Phi)^n\ =\ f\big(g(\Phi)\big)$$
as claimed.
\end{proof}

\noindent
Thus if $\Phi$ is continuous and weakly nilpotent, and $f,g\in K[[z]]$ are such that $fg=1$, then $f(\Phi)$ is a $K$-module automorphism of $V$ with inverse $g(\Phi)$. 

\subsection*{Filtered algebras} {\em Throughout the rest of this section $A$ is a 
\textup{(}not necessarily commutative\textup{)} 
$K$-algebra with $1\ne 0$.}\/ As usual, this includes the requirement that 
$$\lambda(xy)=(\lambda x)y=x(\lambda y) \qquad (\lambda\in K,\ x,y\in A).$$ So
$\lambda \mapsto \lambda1\colon K \to A$ is ring morphism taking values in
$\text{center}(A)$; via this morphism we identify $\Q$ with 
a subring of $\text{center}(A)$.
A {\bf filtration} of the $K$-algebra $A$ is a fil\-tra\-tion~$(A^i)_{i\in\Z}$ of the $K$-module $A$ with $1\in A^0$ and $A^i A^j\subseteq A^{i+j}$ for all $i$, $j$.
Thus $A^0$ is a $K$-subalgebra of $A$, and the $A^i$ are (left and right) $A^0$-submodules of~$A$.
Given a nonnegative filtration 
$(A^n)$ of $A$, each $A^n$ is a two-sided ideal of $A$. 

\index{filtration!algebra}
\index{algebra!filtered}
 
Let $(A^i)$ be a filtration of the $K$-algebra $A$. To keep notations simple 
we denote~$A$ with this filtration also by $A$; the combined object is called a  
{\bf filtered $K$-algebra}. A subalgebra of $A$ with the induced 
filtration is then also a filtered $K$-algebra. The norm on 
the additive group of $A$ obtained
from the filtration satisfies 
$ |xy| \le |x|\cdot |y|$ for all $x,y\in A$, so $A$ with the topology
given by the filtration is a topological ring. If~$(u_{\lambda})_{\lambda\in \Lambda}$ and $(v_{\omega})_{\omega\in \Omega}$ are summable families in $A$ (with respect to this filtration), then so is the family
$\big(u_{\lambda}v_{\omega}\big)$ indexed by $\Lambda\times \Omega$ and $$\sum u_{\lambda}v_{\omega}\ =\ \sum u_{\lambda} \cdot \sum v_{\omega}.$$

\begin{example}
Let $R$ be a $K$-algebra with $1_R\ne 0$, and let $A:=R[[z]]$ be the $K$-algebra
of power series in one commuting indeterminate $z$ with coefficients in $R$.
(As usual we identify $R$ with a subring of $A$ via $r\mapsto rz^0$.) 
Then $(z^nA)$ is a
nonnegative filtration of~$A$ making $A$ a complete filtered $K$-algebra.
\end{example}

\subsection*{Logarithms}
Let $A$ be a complete filtered $K$-algebra with respect to the fil\-tra\-tion~$(A^i)$.
Then every sequence $(a_n)$ in $A$ with $a_n\in A^n$ for each $n$ is summable in~$A$. In particular,
given a formal power series $f=\sum_{n} f_n z^n\in K[[z]]$ ($f_n\in K$ for each~$n$) and $a\in A^1$, the sequence $\left(f_n a^n\right)$ is summable in $A$, and we set $f(a):=\sum f_n a^n$. Here is an analogue of Lemma~\ref{lem:linear maps into power series}: \index{logarithm!on an algebra}\index{exponential!on an algebra}

\begin{lemma}\label{substitutionlemma} Let $a\in A^1$. Then the map
$f\mapsto f(a)\colon K[[z]] \to A$ is a $K$-algebra morphism of rank $0$. 
If $f\in K[[z]]$, $g\in zK[[z]]$, then the element $f\circ g\in K[[z]]$
satisfies $(f\circ g)(a)=f\big(g(a)\big)$. Also $f(bab^{-1})=bf(a)b^{-1}$ for
$f\in K[[z]]$, $b\in A^{\times}$.
\end{lemma}

\noindent
As a consequence,
each element of $1+A^1$ is a unit in $A$: if $a\in A^1$, then $1-a$
has a multiplicative (two-sided) inverse given by
$(1-a)^{-1} = \sum_{i=0}^\infty a^i$.
Thus the multiplicative group $1+A^1$ is a subgroup of the group 
$A^\times$ of units of $A$. Applying Lemma~\ref{substitutionlemma} to
the series $\exp(z)=\sum z^n/n!$ and $\log(1+z)=\sum_{n=1}^{\infty} (-1)^{n+1}z^n/n$ we  see that 
the maps
\begin{align*} A^1\ \to\ 1+A^1\  &\colon\ a\mapsto \exp(a)\ :=\ \sum \frac{a^n}{n!},\\
 1+A^1\ \to\ A^1\ &\colon\  a\mapsto \log(a)\ :=\ \sum_{n=1}^\infty (-1)^{n+1}\frac{(a-1)^n}{n}
 \end{align*}
are each other's inverse, and thus bijective. Note that $\exp(0)=1$.

\begin{lemma}\label{ab=ba} Let $a,b\in A$ and $ab=ba$. Then: 

\begin{tabular}{lll}
\hskip-2em\textup{(i)} & $\exp(a+b)\ =\ \exp(a)\exp(b)$ & for $a,b\in A^1$; \\
\hskip-2em\textup{(ii)} &  $\log(ab)\ =\ \log(a) + \log b$ & for $a,b\in 1+A^1$; \\
\hskip-2em\textup{(iii)} & $(\log a)(\log b)\ =\ (\log b)(\log a)$ & for $a,b\in 1+A^1$.
\end{tabular}
%\begin{enumerate}
%\item $\ \exp(a+b)\ =\ \exp(a)\exp(b)\ \quad $ for $a,b\in A^1$;
%\item  $\  \log(ab)\ =\ \log(a) + \log b\quad \qquad $ for $a,b\in 1+A^1$;
%\item $\ (\log a)(\log b)\ =\ (\log b)(\log a)\quad $ for $a,b\in 1+A^1$.
%\end{enumerate}
\end{lemma}
\begin{proof} We leave (i) and (ii) as routine verifications. As to (iii), this follows from the (almost obvious) fact that for $f,g\in K[[z]]$ and $s,t\in A^1$ with $st=ts$ we have $f(s)g(t)=g(t)f(s)$:
for $a,b\in 1+A^1$ we have $a=1+s$, $b=1+t$ with $s,t\in A^1$, and
then $ab=ba$ gives $st=ts$, and so
$\log(1+s)\log(1+t)=\log(1+t)\log(1+s)$.
\end{proof}

\noindent
The mutually inverse nature of $\exp$ and $\log$ gives
$$\exp(A^n)\ =\ 1+A^n,\qquad \log(1+A^n)\ =\ A^n\qquad \text{for $n\ge 1$.}$$
For $r\in K$ and $a\in 1+A^1$ we set $a^r:=\exp(r\log a)\in 1+A^1$. 
Thus for $a\in 1+A^1$ and $r,s\in K$ we have  $\log(a^r) = r\log a$ and
$$ a^0\ =\ 1,  \qquad a^{r+s}\ =\ a^r\cdot a^s, \qquad a^{-r}\ =\ (a^r)^{-1}, \qquad a^n\ =\ \underbrace{a\cdots a}_{\text{$n$ times}}.$$
Also, $\exp(ra)= \exp(a)^r$ for $a\in A^1$ and $r\in K$, and 
$$(ab)^r\ =\ a^r b^r\qquad\text{for $a,b\in 1+A^1$ with $ab=ba$ and $r\in K$.}$$
We also write $\ex^a$ instead of $\exp(a)$ when $a\in A^1$. Thus $\ex^z=\exp(z)=\sum z^n/n!$ in the filtered $K$-algebra $K[[z]]$ of the example above.

\subsection*{Lie algebras and filtered Lie algebras}
Recall that a Lie algebra over $K$ is a $K$-module $L$ equipped with a 
$K$-bilinear operation $[\ \, ,\ ]\colon L\times L\to L$ (the Lie bracket of $L$) satisfying 
$[x,x]=0$ for all $x\in L$, as well as the Jacobi Identity:
$$[x,[y,z]] + [y,[z,x]] + [z,[x,y]]\ =\ 0\ \qquad\text{for all $x,y,z\in L$.}$$
(Thus $[x,y]=-[y,x]$ for $x,y\in L$, in view of $[x+y,x+y]=0$.)
A Lie al\-ge\-bra~$L$ over $K$ is said to be {\bf abelian} if $[x,y]=0$ for all~$x,y\in L$. Given a Lie algebra $L$ over $K$, a  {\bf Lie subalgebra} of $L$
is a $K$-submodule $M$ of $L$ such that $[x,y]\in M$ for all~$x,y\in M$, and an {\bf ideal} of $L$ is a $K$-submodule $M$ of $L$ such that $[x,y]\in M$ for all~$x\in L$ and $y\in M$.  
%Given a Lie algebra $L$ over $K$, we put $[x]:=x$ for $x\in L$, and inductively we define the iterated Lie brackets:
%$$[x_1,\dots,x_n] := [x_1,[x_2,\dots,x_n]]\qquad\text{for $x_1,\dots,x_n\in L$, $n\geq 3$.}$$
A {\bf filtration} of a Lie algebra $L$ over $K$ is a filtration $(L^i)_{i\in\Z}$ of~$L$ as a $K$-module such that additionally $[L^i,L^j]\subseteq L^{i+j}$ for all indices $i$, $j$.

\index{Lie algebra}
\index{Lie algebra!abelian}
\index{Lie algebra!subalgebra}
\index{Lie algebra!ideal}
\index{ideal!Lie algebra}
\index{algebra!Lie}

\index{filtration!Lie algebra}
\index{filtered!Lie algebra}

\medskip
\noindent
Let $(L^i)_{i\in\Z}$ be a filtration of the Lie algebra $L$ over $K$. 
If $i\ge 0$, then $[L^i,L^i]\subseteq L^i$, so~$L^i$ is a Lie subalgebra of $L$. If the filtration $(L^i)$ of $L$ is nonnegative, i.e., if $L^0=L$, then each $L^i$ is even an ideal of $L$, since then $[L,L^i]=[L^0,L^i]\subseteq L^i$.
The norm on~$A$ defined by the filtration satisfies
$|[x,y]|\le |x|\cdot |y|$ for all $x,y\in L$, and so the 
Lie bracket operation $[\ \, ,\ ]\colon L\times L\to L$ is continuous 
with respect to the topology on $L$ given by the filtration (with the product topology on $L\times L$).

\begin{example}
The binary operation on $A$ given by
$$[a,b]\ :=\  ab-ba\qquad (a,b\in A)$$
turns $A$ into a Lie algebra $A_{\Lie}$ over $K$. Every filtration of the $K$-algebra $A$ is also a filtration of $A_{\Lie}$. 
\end{example} 

\nomenclature[Y]{$A_{\Lie}$}{Lie algebra associated to the algebra $A$}

\noindent
A {\bf filtered Lie algebra over $K$} is a Lie algebra $L$ over $K$ together with a filtration of~$L$.
Let $(L^i)$ be the filtration of a filtered Lie algebra $L$ over $K$. Given a Lie sub\-al\-ge\-bra~$M$ of $L$, $(L^i\cap M)$ is a filtration of $M$, called the filtration of $M$ {\bf induced} by~$(L^i)$.
A Lie subalgebra of $L$ equipped with the filtration induced by $(L^i)$ is called a {\bf filtered Lie subalgebra} of $L$.

\index{filtration!Lie algebra}
\index{filtration!Lie algebra!induced on a Lie subalgebra}
\index{filtered!Lie algebra}
\index{filtered!Lie subalgebra}
\index{Lie algebra!filtered}

\subsection*{Derivations}  
A {\bf $K$-derivation} on $A$ is a $K$-linear map $\Delta\colon A\to A$ such that
$$\Delta(xy)\ =\ \Delta(x)y+x\Delta(y)\qquad\text{for all $x,y\in A$}$$
(and thus $\Delta(\lambda\cdot1)=0$ for $\lambda\in K$).
 
\index{derivation!on an algebra} 
\index{K-derivation@$K$-derivation}
\index{derivation!$K$-derivation}
\index{K-derivation@$K$-derivation!internal}
\nomenclature[Y]{$\operatorname{ad}_a$}{adjoint derivation}

\begin{example}
Given $a\in A$, the adjoint $x\mapsto \operatorname{ad}_a(x):=[a,x]=ax-xa$ is a $K$-derivation on $A$. 
The $K$-derivations $\operatorname{ad}_a$  on $A$  (where $a\in A$) are called {\bf internal.}
If $A$ is a filtered $K$-algebra and $a\in A^d$, then $\operatorname{ad}_a$ is of rank $d$.
\end{example}

\noindent
We denote the set of all $K$-derivations on $A$ by $\Der_K(A)$. 
If $\Delta_1,\Delta_2\in\Der_K(A)$ and $r\in K$, then 
$\Delta_1+\Delta_2\in\Der_K(A)$ and $r\Delta_1\in\Der_K(A)$, and so 
$\Der_K(A)$ is naturally a left $K$-module.
One also verifies easily that if $\Delta_1, \Delta_2\in \Der_K(A)$, then 
$$[\Delta_1,\Delta_2]\ :=\ \Delta_1\Delta_2 - \Delta_2\Delta_1\in \Der_K(A). $$
With this operation $[\ \, ,\ ]$, the $K$-module $\Der_K(A)$ is  a Lie algebra over $K$, in fact, a Lie subalgebra of the Lie algebra
$\End(A)_{\Lie}$ over $K$, where $\End(A)$ is the $K$-algebra of endomorphisms of $A$ as a $K$-module.
If $\sigma$ is an automorphism of the $K$-algebra~$A$, then for every 
$\Delta\in \Der_K(A)$ we have  $\sigma\Delta\sigma^{-1}\in \Der_K(A)$, and 
$\Delta\mapsto \sigma\Delta\sigma^{-1}$ is an automorphism of the Lie algebra 
$\Der_K(A)$ over $K$, with inverse 
$\Delta\mapsto \sigma^{-1}\Delta\sigma$.

\nomenclature[Y]{$\Der_K(A)$}{Lie algebra of $K$-derivations of the $K$-algebra $A$}

\medskip
\noindent
Given $\Delta\in\Der_K(A)$, we have for all $x,y\in A$,
\begin{equation}\label{eq:Leibniz}
\Delta^n(xy)\ =\ \sum_{i+j=n}\binom{n}{i}\Delta^i(x)\Delta^j(y)  \qquad \text{\textup{(}Leibniz rule\textup{)},}
\end{equation}
and more generally, for $m\ge 1$ and $x_1,\dots,x_m\in  A$,
$$\Delta^n(x_1\cdots x_m)\ =\  \sum_{i_1+\cdots+i_m=n} \frac{n!}{i_1!\cdots i_m!} \Delta^{i_1}(x_1)\cdots\Delta^{i_m}(x_m).
$$

\noindent
We shall often use the following facts, the second of which follows from \eqref{eq:Leibniz} and the remark after Lemma~\ref{lem:linear maps into power series}:

\begin{lemma}\label{lem:exp Delta} Let $A$ be a complete filtered $K$-algebra with respect to the filtration~$(A^i)$. Suppose $\Delta\in \Der_K(A)$ is continuous. Then \begin{enumerate}
\item[\textup{(i)}] if $A$ is commutative and $a\in A^1$, then 
$\Delta(\ex^a)=\ex^a\Delta(a)$; 
\item[\textup{(ii)}] if $\Delta$ is weakly nilpotent, then 
the $K$-module endomorphism $\ex^{\Delta}:= \exp(\Delta)$ of $A$ is a $K$-algebra automorphism of~$A$ with inverse $\ex^{-\Delta}$.
\end{enumerate}
\end{lemma}

\subsection*{Graded algebras}
A {\bf grading} \index{grading!algebra} of $A$ is a family $(A_i)_{i\in\Z}$ of $K$-submodules of $A$ such that the following two conditions are satisfied:
\begin{enumerate}
\item $A=\bigoplus_{i\in\Z} A_i$ (internal direct sum of $K$-submodules of $A$); 
\item $A_i A_j\subseteq A_{i+j}$ for all $i,j\in\Z$.
\end{enumerate}
A little argument shows that then $1\in A_0$, so $A_0$ is a $K$-subalgebra of $A$, and each~$A_i$ is a left-and-right $A_0$-submodule of $A$. 
A {\bf graded $K$-algebra} \index{graded algebra}\index{algebra!graded} is a $K$-algebra $A$ together with a grading 
of $A$. Let $A$ be a graded $K$-algebra and let~$(A_i)$ 
be its grading. 
The elements of $A_i$ are said to be \index{homogeneous!element of a graded algebra} {\bf homogeneous of degree~$i$}. For every $a\in A$ there is a unique family
$(a_i)$ with $a_i\in A_i$ for each $i$ and $a_i=0$ for all but finitely 
many $i$ such that $a=\sum_{i} a_i$. For this family $(a_i)$ we call
 $a_i$ the {\bf homogeneous part of $a$ of degree $i$}. 
For $a\in A$, $a\neq 0$, we define the {\bf degree of~$a$} \index{degree!element of a graded algebra} as the largest $i$ such that $a_i\neq 0$, denoted by $\degree(a)$, or by $\degree_{(A_i)}(a)$ 
if we want to indicate the dependence on $(A_i)$. We also set $\degree(0):=-\infty<\Z$.
The grading $(A_i)$ of $A$ is said to be {\bf nonnegative} \index{grading!algebra!nonnegative} if $A_i=\{0\}$ for all $i<0$. From the grading $(A_i)$ we obtain the filtration $(A^i)$ \index{grading!algebra!associated filtration}\index{filtration!algebra!associated to a grading} of the $K$-algebra $A$ by
$$A^i\ :=\ \bigoplus_{j\geq i} A_j,$$ 
the filtration of $A$ {\bf associated} to $(A_i)$. Clearly $(A_i)$ is nonnegative iff $(A^i)$ is.

%Given graded $K$-algebras $(A,\mathbf A)$ and $(B,\mathbf B)$ and some $d\in\Z$, a $K$-linear map $\Phi\colon A\to B$ is said to be {\bf homogeneous of degree $d$} if $\Phi(A_i)\subseteq B_{i+d}$ for every $i\in\Z$.
%Given a further graded $K$-algebra $(C,\mathbf C)$ and a homogeneous $K$-linear map $\Psi\colon B\to C$ of degree $e\in\Z$, the $K$-linear map $\Psi\circ\Phi\colon A\to C$ is homogeneous of degree $d+e$.

\subsection*{Gradings of polynomial algebras}
{\em In the rest of this section $A=K[Y_0, Y_1,\dots]$ is the \textup{(}commutative\textup{)} $K$-algebra of polynomials in the distinct indeterminates $Y_n$, $n=0,1,2,\dots$, and
$\i$ ranges over the set~$\N^{(\N)}$ of 
sequences $\i=(i_0,i_1,\dots)\in\N^\N$ such that $i_n=0$ for all but finitely many $n$.}\/ For each $\i$ we set
$$Y^{\i}\ :=\  Y_0^{i_0}Y_1^{i_1} \cdots Y_n^{i_n}\cdots\in A.$$
Given $P\in A$ we have a unique family $(P_{\i})$
in $K$ such that $P = \sum_{\i} P_{\i}\, Y^{\i}$
with $P_{\i}=0$ for all but finitely many $\i$.
Below, $Y^\diamond:=\{Y^{\i}:\i\in\N^{(\N)}\}$ is the multiplicative monoid 
of monomials. 

\begin{definition}
A {\bf degree function on $A$} is a function $\degree\colon Y^\diamond\to\Z$ 
such that $\degree(1)=0$ and $\degree(s\cdot t)=\degree(s)+\degree(t)$ for all $s,t\in Y^\diamond$.
\end{definition}

\begin{example}
Given a sequence $(d_n)$ of integers, define
$$\degree(Y^\i)\ :=\ \sum_n d_ni_n\in \Z \qquad\text{for each $\i$.}$$
Then $\degree$ is a degree function on $A$. Any degree function $\degree$
on $A$ arises from a sequence~$(d_n)$ of integers in this manner, by setting $d_n:=\degree(Y_n)$ for each $n$.
\end{example}

\index{degree function!on a polynomial algebra}
\index{homogeneous!with respect to a degree function}

\noindent
Any degree function $\degree$ on $A$ yields a grading 
$(A_d)_{d\in\Z}$ of $A$, 
$$A_d\ :=\ \big\{P\in A:\ \text{$P_{\i}=0$ if $\degree(Y^\i)\neq d$}\big\}\ =\
\sum_{\degree(Y^{\i})=d} KY^{\i} .$$
Suppose $(A_i)$ is a grading of $A$ for which each indeterminate $Y_n$ 
is homogeneous. This grading is induced by a degree function on $A$ as above, namely the restriction of $\degree_{(A_i)}$ to $Y^\diamond$. 
To simplify notation denote
$\degree_{(A_i)}$ by $\degree$. Then
$$\degree(P)\ =\ \max\!\big\{\!\operatorname{d}(Y^\i):\ P_\i\neq 0\big\}\in\Z \text{ if $P\neq 0$,}
\qquad \degree(0)=-\infty<\Z,$$
and the elements of $A_i$ are said to be {\bf $\degree$-homogeneous of degree 
$i$.} Given $P\in A$ and $i\in\Z$, the homogeneous part of $P$ of degree $i$ with respect to $(A_i)$ is called the {\bf $\degree$-homogeneous part of $P$ of degree $i$.} Given $d\in\Z$, a $K$-linear map $\Phi\colon A\to A$ is said to be {\bf $\degree$-homogeneous of degree $d$} if $\Phi(A_i)\subseteq A_{i+d}$ for all $i$; 
given also a $\degree$-ho\-mo\-ge\-ne\-ous $K$-linear map $\Psi\colon A\to A$ of degree $e\in\Z$, the $K$-linear map~${\Psi\circ\Phi\colon A\to A}$ is $\degree$-homogeneous of degree $d+e$.

Clearly the grading induced by a degree function $\degree$ on $A$ is nonnegative if and only if~${\degree(Y_n)\geq 0}$ for each $n$; in this case we say that $\degree$ is {\bf nonnegative}.
We now define two important nonnegative degree functions on $A$:
 
\index{degree!in a polynomial algebra} 
\index{weight!in a polynomial algebra} 
\index{isobaric!element of a polynomial algebra}
\index{homogeneous!element of a polynomial algebra}

\begin{example}
The usual (total) degree 
$$\deg(Y^\i)\ :=\ \abs{\i}\ =\  i_0+ i_1 + i_2 + \cdots + i_n + \cdots\in\N$$ 
yields the degree function $\deg$ on $A$. For $P\in A$ we have
$$\deg(P)\ =\ \max\big\{\abs{\i}:\ P_\i\neq 0\big\}\in\N \text{ if $P\neq 0$,}
\qquad \deg(0)=-\infty<\Z.$$
We denote the grading associated to $\deg$ by $(A_{d})$.
In the rest of this chapter, the term {\bf homogeneous}  (no mention of a degree function or grading) is synonymous with \textit{$\deg$-homogeneous.}\/  For $P\in A$ and $d\in\N$ we let
$$P_{d}\ =\ \sum_{\abs{\i}=d} P_{\i}\, Y^{\i}\in A_d$$
be the homogeneous part of  $P$ of degree $d$. Thus 
$A_d = \{P\in A: P=P_{d}\}$ for each~$d\in\N$. For example,
$A_0=K$, $A_1=\bigoplus_n K\,Y_n$, $A_2=\bigoplus_{m\le n}KY_mY_n$.
\end{example} 

\begin{example}
Setting 
$$\wt(Y^\i)\ :=\ \dabs{\i}\ =\  i_1 + 2i_2 + \cdots + ni_n + \cdots \in\N$$ 
gives a degree function $\wt$ on $A$. Note that $\wt(Y_i)=i$ for each $i\in\N$.
For $P\in A$ we call $\wt(P)$ the {\bf weight} of $P$:
$$\wt(P)\ =\  \max\big\{\dabs{\i}:\ P_\i\neq 0\big\}\in\N \text{ if $P\neq 0$,}
\qquad \wt(0)=-\infty<\Z.$$
We denote the grading associated to $\wt$ by $(A_{[w]})$.
In the rest of this chapter, {\bf isobaric} is synonymous with \textit{$\wt$-homogeneous.}\/ For $P\in A$ and $w\in\N$ we let
$$P_{[w]}\ =\  \sum_{\dabs{\i}=w} P_{\i}\, Y^{\i}\in A_{[w]}$$
be the {\bf isobaric part of $P$ of weight $w$.} Thus 
$$A_{[w]}\ =\  \{P\in A:\  P=P_{[w]}\}\qquad\text{for each $w\in\N$,}$$
so $A_{[0]}=K[Y_0]$, $\ A_{[1]}=K[Y_0]\,Y_1$, $\ A_{[2]}=K[Y_0]Y_1^2 + K[Y_0]Y_2$.
\end{example}

\noindent
If $K$ is equipped with a derivation, $Y$ is a differential indeterminate over $K$, and the $K$-algebra $K\{Y\}=K[Y,Y',\dots]$ is identified with 
$K[Y_0,Y_1,\dots]$ by setting $Y_n=Y^{(n)}$, then these notions of degree and weight agree with the ones for differential polynomials introduced in Section~\ref{Decompositions of Differential Polynomials}.

\section{Triangular Linear Maps}\label{triangularlinear}

\noindent
{\em In this section $V$ is a $K$-module.}\/ We equip $V$ with the trivial filtration. This makes~$V$ a complete filtered $K$-module.
Recall from Section~\ref{sec:Filtered and Graded Algebras} that $\End(V)$ is the $K$-algebra  of
endomorphisms of $V$. 
 Thus a family $(\Phi_i)_{i\in I}$ of endomorphisms of~$V$ is summable (as defined in Section~\ref{sec:Filtered and Graded Algebras}) iff for each $v\in V$ we have $\Phi_i(v)=0$ for all but finitely many $i$; recall that then $\sum\Phi_i$ is the endomorphism of $V$ given~by
$$\left(\sum \Phi_i\right)(v)\  =\ \sum \Phi_i(v).$$
If $(\Phi_i)_{i\in I}$ and $(\Psi_j)_{j\in J}$ are summable families of 
endomorphisms of $V$, then so is $(\Phi_i\Psi_j)_{(i,j)\in I\times J}$ with
$$  \left(\sum_{i} \Phi_i\right)\left(\sum_{j} \Psi_j\right)\ =\ \sum_{i,j}\Phi_i\Psi_j .$$ 
Throughout this section $\Phi,\Psi\in\End(V)$.

\subsection*{Locally nilpotent and locally unipotent endomorphisms}
We call $\Phi$ {\bf locally nilpotent} if for every $v\in V$ there is some $n$ such that $\Phi^n(v)=0$. 
Note that $\Phi$ is locally nilpotent iff $\Phi$ is weakly nilpotent with respect to the trivial filtration on $V$.
If $\Phi$ is locally nilpotent, then so is every power 
$\Phi^n$ with $n\ge 1$ and every scalar multiple $\lambda\Phi$, 
and for each family $(\lambda_n)$ of scalars the family $(\lambda_n\Phi^n)$ 
of endomorphisms is summable, and 
$\sum_{n\ge 1} \lambda_n\Phi^n$ is locally nilpotent. 
We call $\Phi$  {\bf locally unipotent} if $\Phi-1$ is locally nilpotent.
If $\sigma$ is an automorphism of the $K$-module $V$, and $\Phi$ is locally nilpotent, respectively locally unipotent, then so is~$\sigma\Phi\sigma^{-1}$.

\index{nilpotent!locally}
\index{locally!nilpotent}
\index{locally!unipotent}

\medskip\noindent
\textit{In the rest of this subsection we assume that $\Phi$ is locally nilpotent.}\/
Then 
 $(\frac{1}{n!}\Phi^n)$ is summable, so we can define the endomorphism 
$\exp\Phi$ of $V$ by
$$\exp\Phi\ :=\ \sum_{n=0}^\infty \frac{1}{n!}\Phi^n\ =\ 1+\Phi+\frac{1}{2}\Phi^2+\cdots+\frac{1}{n!}\Phi^n+\cdots.$$
Since $\exp \Phi -1=\sum_{n=1}^\infty \frac{1}{n!}\Phi^n$ is locally nilpotent, 
$\exp \Phi$ is locally unipotent.

\medskip\noindent
Conversely, assume $\Psi$ is locally unipotent. Then we define the
locally nilpotent endomorphism $\log\Psi$ of $V$ by
$$
\log\Psi\ :=\ \sum_{n=1}^\infty \frac{(-1)^{n+1}}{n}(\Psi-1)^n\ =\  
(\Psi-1)-\frac{1}{2}(\Psi-1)^2+\frac{1}{3}(\Psi-1)^3 -\cdots.
$$
We have $\exp(\log\Psi)=\Psi$ and $\log(\exp\Phi)=\Phi$. If
$\Phi_1, \Phi_2$ are commuting locally nilpotent endomorphisms of $V$, then
$\Phi_1\Phi_2$ and $\Phi_1+\Phi_2$ are locally nilpotent and  
$$\exp(\Phi_1)\exp(\Phi_2)\ =\ \exp(\Phi_1+\Phi_2).$$ 
Thus $\exp\Phi$ is an automorphism of the $K$-module $V$ with
inverse $\exp(-\Phi)$, and 
$$%\begin{equation}\label{eq:powers of exp Phi}
(\exp\Phi)^k\ =\ \exp(k\Phi) \qquad (k\in \Z).
$$%\end{equation}
If $\sigma$ is an automorphism of the $K$-module $V$, then
$$\exp(\sigma\Phi\sigma^{-1})\ =\ \sigma\exp(\Phi)\sigma^{-1}, \qquad
  \log(\sigma\Psi\sigma^{-1})\ =\ \sigma\log(\Psi)\sigma^{-1}.$$
We also note that for $v\in V$ we have
\begin{equation}\label{eq:invariants}
\Phi(v)=0 \quad\Longleftrightarrow\quad (\exp \Phi)(v)=v.
\end{equation}

\subsection*{Triangular matrices}
We construe $K^{\N\times\N}$ as a $K$-module with the componentwise addition and scalar multiplication. The elements $M=(M_{ij})_{i,j\in\N}$ of $K^{\N\times\N}$ may be visualized as infinite square matrices with entries in $K$:
$$M\ =\  \begin{pmatrix}
M_{00} & M_{01} & M_{02} &  \cdots &\\
M_{10} & M_{11} & M_{12} &  \cdots &\\
M_{20} & M_{21} & M_{22} &  \cdots &\\
\vdots & \vdots & \vdots &  \ddots &
\end{pmatrix}.$$
We say that $M=(M_{ij})\in K^{\N\times\N}$ is {\bf column-finite} if for each $j$ there are only finitely many $i$ with $M_{ij}\neq 0$. Given column-finite matrices
$M=(M_{ij})$ and $N=(N_{ij})$ we can define their matrix product 
$MN\in K^{\N\times\N}$ by 
$$  (MN)_{ij}\  :=\   \sum_{k} M_{ik} N_{kj}.$$
Then $MN$ is again column-finite.
With this product operation, the $K$-submodule of~$K^{\N\times\N}$ consisting of all column-finite matrices is a $K$-algebra with multiplicative identity $1$ given by the identity matrix. 
We identify $K$ with a subring of this $K$-algebra via 
$\lambda\mapsto \lambda\cdot 1$.

\index{matrix!column-finite}
\index{column-finite}
\nomenclature[Y]{$K^{\N\times\N}$}{$K$-module of matrices  $M=(M_{ij})_{i,j\in\N}$  over $K$}

We say that $M=(M_{ij})\in K^{\N\times\N}$ is (upper) {\bf triangular} \index{matrix!triangular}\index{triangular!matrix} if $M_{ij}=0$ for all $i,j$ with $i>j$. %, and {\bf strictly triangular} if $M_{ij}=0$ for all $i,j\in\N$ with $i\geq j$. 
The set $\tr_K$ of triangular matrices is a subalgebra of the $K$-algebra 
of column-finite matrices. \nomenclature[Y]{$\tr_K$}{$K$-algebra of triangular matrices over $K$} For every $n$ we set
$$\tr_K^n\  :=\ \big\{ M\in\tr_K :\ \text{$M_{ij}=0$ for all $i,j$ with $j< i+n$}\big\},$$
so $\tr_K^0=\tr_K$, and $\tr_K^1=\{ M\in\tr_K:\  \text{$M_{ii}=0$ for all $i$}\}$.
It is easily verified that~$(\tr_K^n)$ is a complete nonnegative
filtration of the $K$-algebra $\tr_K$; in particular, each~$\tr_K^n$ is an ideal of the Lie algebra $(\tr_K)_{\Lie}$.
We say that $M\in K^{\N\times \N}$ is {\bf diagonal} if~$M_{ij}= 0$ for all $i\ne j$. Then 
$$D_K\ :=\ \{M\in K^{\N\times \N}:\ \text{$M$ is diagonal}\}$$ is a (commutative)
subalgebra of the $K$-algebra $\tr_K$. For $M\in \tr_K$ we de\-fine the matrix~${M_0\in D_K}$ by
$(M_0)_{ii}=M_{ii}$. Then $M\mapsto M_0 \colon \tr_K\to D_K$ is a $K$-algebra
morphism that is the identity on $D_K$. The multiplicative group of units of 
$D_K$ is
$$D^\times_K\ =\ \big\{ M\in D_K :\ 
\text{$M_{ii}\in K^\times$ for all $i$}\big\}.$$  
The group morphism $M\mapsto M_0 \colon \tr_K^\times \to D_K^\times$
from the group $\tr_K^\times$ of units of $\tr_K$ onto~$D_K^\times$ is the identity on $D_K^\times$ and has kernel $1+\tr_K^1$, so  $1+\tr_K^1$ 
is a normal subgroup of $\tr_K^\times$ with $(1+\tr_K^1)D_K^\times = \tr_K^\times$ and 
$(1+\tr_K^1)\cap D_K^\times=\{1\}$. Thus $\tr_K^\times$ is the internal semidirect
product of $(1+\tr_K^1)$ with $D_K^\times$. 
We set $\mathcal U:= 1+\tr_K^1$ (also denoted by~$\mathcal U_K$ if we need to
indicate the dependence on $K$) and call its elements {\bf unitriangular} matrices. A {\bf unitriangular} group over $K$ is a subgroup of $\mathcal U$. 

\index{matrix!unitriangular}
\index{matrix!diagonal}
\index{diagonal!matrix}
\index{unitriangular!matrix}
\index{unitriangular!group}
\index{group!unitriangular}

\index{triangular!endomorphism}
\index{endomorphism!triangular}

\nomenclature[Y]{$D_K$}{$K$-algebra of diagonal matrices over $K$}
\nomenclature[Y]{$\TrEnd_K(V)$}{$K$-algebra of triangular endomorphisms of $V$}

\subsection*{Triangular linear maps}
We now assume that $V$ is a free $K$-module on the basis~$(Y_n)$. In particular,
$$V\ =\ \bigoplus_{n} \,K\, Y_n \qquad\text{(internal direct sum of $K$-submodules of $V$).}$$
The set of $K$-linear maps $V\to V$, under (pointwise) addition and composition, forms the $K$-algebra $\End(V)$.
We say that $\Phi$ is {\bf triangular} if
$$\Phi(Y_j)\ =\ \Phi_{0j} Y_0 + \Phi_{1j} Y_{1} + \cdots + \Phi_{jj} Y_j\qquad\text{where $\Phi_{ij}\in K$ for $i,j\in\N$, $i\leq j$.}$$ 
Triangular endomorphisms of the $K$-module $V$ may be conveniently represented by triangular bi-infinite  matrices with entries in $K$: for every triangular $\Phi$ define $$M_\Phi\ :=\ (\Phi_{ij})_{i,j\in\N}\ =\  
\begin{pmatrix}
\Phi_{00} 	& \Phi_{01} 	& \Phi_{02} 	& \Phi_{03}  & \cdots  \\
  	& \Phi_{11}				& \Phi_{12}	& \Phi_{13} & \cdots  \\
	&				& \Phi_{22}			& \Phi_{23} & \cdots  \\
	&				&				& \Phi_{33}			  & \cdots  \\
	&				&				&			& \ddots
\end{pmatrix}
\in \tr_K.$$
Here and below, given triangular $\Phi$ we set $\Phi_{ij}:=0$ for all $i,j\in\N$ with~${i>j}$. It is easily verified that if $\Phi$ and $\Psi$ are triangular, then the composition $\Phi\Psi$ and the sum $\Phi+\Psi$ are also triangular, with $M_{\Phi\Psi}=M_\Phi \cdot M_{\Psi}$ and $M_{\Phi+\Psi}=M_\Phi + M_\Psi$.
Hence the triangular endomorphisms of the $K$-module $V$ form a $K$-subalgebra~$\TrEnd_K(V)$ of~$\End(V)$, isomorphic to the $K$-algebra $\tr_K$ via the isomorphism $\Phi\mapsto M_\Phi$. 
If the $K$-module $V$ and the basis $(Y_n)$ are clear from the context, we also abbreviate~$\TrEnd_K(V)$ as~$\TrEnd_K$. 
For each $n$ we set
$$\TrEnd^n_K\ :=\ \{ \Phi\in\TrEnd_K :\  M_\Phi \in \tr_K^n \}.$$
Then $(\TrEnd^n_K)$ is a complete nonnegative filtration of the $K$-algebra $\TrEnd_K$.

\begin{lemma} If the sequence $(\Phi_n)$ in $\TrEnd_K$ is summable
in the sense of the filtration~$(\TrEnd^n_K)$, then
$(\Phi_n)$ is summable as defined in the beginning of this section, and
the sum $\sum \Phi_n\in \TrEnd_K$ in the sense of the filtration $(\TrEnd^n_K)$ equals $\sum \Phi_n$ as defined
in the beginning of this section.
\end{lemma}

\noindent
We shall use this fact tacitly in what follows.

\medskip
\noindent
Suppose $\Phi$ is triangular and $\Phi_{ii} = 0$ for all $i$. Then 
$\Phi^n\in \TrEnd^n_K$ for all $n$, so $\Phi$ is locally nilpotent.
Moreover, $\exp\Phi$ is triangular with  
$(\exp\Phi)_{ii}=1$ for all $i$, and $M_{\exp \Phi}=\exp M_{\Phi}$. 
We can reverse this as follows. Suppose
$\Psi$ is triangular and $\Psi_{ii}=1$ for all $i$. Then $(\Psi-1)^n\in \TrEnd^n_K$ for all $n$, so $\Psi$ is locally unipotent. Moreover, 
$\log\Psi$ is triangular with
$(\log\Psi)_{ii}=0$ for all $i$, and $M_{\log\Psi}=\log M_\Psi$.

\subsection*{Diagonals}
Let $M\in\tr_K$. For each $n$ we call the triangular matrix 
$$M_n\  =\ 
%\begin{tikzpicture}[baseline=(current bounding box.center)]
%\matrix (m) [matrix of math nodes,nodes in empty cells,right delimiter={)},left delimiter={(} ]{
%0\ &  	 	&   	& 0   	& M_{0n}	& 0 \     	&      		&     			& \phantom{0} 	& \phantom{0} \\
%  & 0\     	&   	&   	& 0      	& M_{1,n+1}	& 0 \       &   			& \phantom{0} 	& \phantom{0}\\
%  &        	& 0\   	&   	&   		& \ 0   	& M_{2,n+2} & 0 \  			& \phantom{0} 	& \phantom{0}\\
%  &	 	   	&      	&  \ddots 	&        	&          	&  \ddots    		& \ddots	& \ddots	& \phantom{0} \\  
%} ;
%\draw[line width=0.1em, line cap=round, dash pattern=on 0pt off 0.5em] (m-1-1)-- (m-1-4);
%\draw[line width=0.1em, line cap=round, dash pattern=on 0pt off 0.5em] (m-2-2)-- (m-2-5);
%\draw[line width=0.1em, line cap=round, dash pattern=on 0pt off 0.5em] (m-3-3)-- (m-3-6);
%\draw[line width=0.1em, line cap=round, dash pattern=on 0pt off 0.5em] (m-1-6)-- (m-1-10);
%\draw[line width=0.1em, line cap=round, dash pattern=on 0pt off 0.5em] (m-2-7)-- (m-2-10);
%\draw[line width=0.1em, line cap=round, dash pattern=on 0pt off 0.5em] (m-3-8)-- (m-3-10);
%\end{tikzpicture} 
\begin{pmatrix*}[l]
0 & \cdots & \cdots & \cdots & 0\phantom{M}    & M_{0n}\, & \ \ \ 0          & \cdots    &        & \\
  & 0      & \cdots & \cdots & \cdots 		   & \ \ 0      & \!M_{1,n+1}  & \ \ \ 0         & \cdots & \\
  &        & 0      & \cdots & \cdots 		   & \cdots &  \ \ \   0       & M_{2,n+2} &  \ \ \ \, 0      & \cdots \\
  &		   &        & \hskip-0.5em\ddots &        		   &        &              & \hskip-1.5em\ddots    & \ddots & \ddots  
\end{pmatrix*} 
\in \tr_K^n
$$
the {\bf $n$-diagonal} of $M$;  i.e., 
$$(M_n)_{ij}\ =\ 0\ \text{ if $j\neq i+n$,} \qquad
(M_n)_{i,i+n}\ =\ M_{i,i+n}.
$$
In the metric given by the 
filtration $(\tr_K^n)$ we have $\sum_{i=0}^n M_i \to M$ as $n\to \infty$, so 
$(M_n)$ is summable with
$M=M_0+M_1+\cdots+M_n+\cdots$.
We say that $M$ is {\bf $n$-diagonal} if $M=M_n$.
Thus $M$ is diagonal as defined earlier iff $M$ is $0$-diagonal.

\index{matrix!$n$-diagonal}
\index{diagonal!$n$-diagonal!matrix}
\nomenclature[Y]{$\diag_n a$}{$n$-diagonal matrix with sequence $a$ on its $n$th diagonal}

\begin{notation}
For a sequence $a=(a_i)\in K^\N$, define $\diag_n a\in \tr_K^n$ by
$$(\diag_n a)_{i,i+n}\ =\ a_i, \qquad (\diag_n a)_{ij}\ =\ 0\ \text{ for $j\neq i+n$.}$$  
Thus the $n$-diagonal matrices are precisely the matrices $\diag_n a$ with 
$a\in K^\N$. We also abbreviate $\diag_0 a$ as $\diag a$.
\end{notation}

\noindent
The sum of two $n$-diagonal matrices is $n$-diagonal. As for products, we have:

\begin{lemma}\label{lem:formulas for diagonals}
Let $M=\diag_m a$ be $m$-diagonal and $N=\diag_n b$ be $n$-diagonal, where $a=(a_i), b=(b_i)\in K^\N$. Then $MN$ is $(m+n)$-diagonal, in fact
$$MN\  =\  \diag_{m+n}( a_i\cdot b_{i+m} )_{i}.$$
Therefore $[M,N]$ is $(m+n)$-diagonal, with
$$[M,N]\ =\ \diag_{m+n}( a_i\cdot b_{i+m} - b_i\cdot a_{i+n} )_{i},$$
and for each $k\in\N$, $M^k$ is $km$-diagonal, with
$$M^k\  =\  \diag_{km} (a_i \cdot a_{i+m} \cdots a_{i+(k-1)m} )_{i}.$$
\end{lemma}

\begin{cor}\label{cor:formulas for diagonals}
Let $M=\diag_1 a$ where $a=(a_i)\in K^\N$. Then 
$$(\exp M)_{ij}\ =\ \frac{1}{(j-i)!}\, a_i\cdot a_{i+1} \cdots a_{j-1} \qquad\text{for all $i\leq j$.}$$
\end{cor}

\noindent
Suppose $\Phi\in\TrEnd_K$. For each $n$ we call the triangular endomorphism $\Phi_n\in\TrEnd^n_K$ of $V$ with associated matrix $M_{\Phi_n}=(M_\Phi)_n$ the {\bf $n$-diagonal} of $\Phi$.
We also say that $\Phi$ is {\bf $n$-diagonal} if $\Phi=\Phi_n$ and {\bf diagonal} if $\Phi=\Phi_0$. Thus $(\Phi_n)$ is summable with
$$\Phi\ =\  \Phi_0 + \Phi_1 + \cdots + \Phi_n + \cdots,$$
and $\Phi\in\TrEnd^m_K$ if and only if $\Phi_0=\cdots=\Phi_{m-1}=0$.
For $\Psi\in\TrEnd_K$ we have
\begin{equation}\label{eq:diagonals of products}
(\Phi\Psi)_n\ =\ \sum_{i+j=n} \Phi_i\Psi_j.
\end{equation}
Hence if $\Psi$ is a diagonal automorphism of $V$ then
$$(\Psi\Phi\Psi^{-1})_n\ =\ \Psi\Phi_n\Psi^{-1}.$$
The identity \eqref{eq:diagonals of products} also implies, for $k\in\N$:
$$(\Phi^k)_{n}\ =\ \sum_{i_1+\cdots+i_k=n} \Phi_{i_1}\Phi_{i_2}\cdots\Phi_{i_k};$$
thus if  $\Phi\in\TrEnd^m_K$ and $n<mk$, then $(\Phi^k)_n = 0$. This immediately yields:

\index{endomorphism!$n$-diagonal}
\index{diagonal!$n$-diagonal!endomorphism}

\index{endomorphism!diagonal}
\index{diagonal!endomorphism}

\begin{lemma}\label{lem:formulas for diagonals of exp and log}
Let $\Phi,\Psi\in\TrEnd^1_K$ with $\Phi=\log(1+\Psi)$, and $m\ge 1$. Then $\Phi\in\TrEnd^m_K$ if and only if $\Psi\in\TrEnd^m_K$. Also, for all $n\ge 1$:
$$\Phi_n\ =\ \sum_{i_1,\dots, i_k} \frac{(-1)^{k+1}}{k}\,\Psi_{i_1}\cdots \Psi_{i_k},\ \qquad
\Psi_n\ =\ \sum_{i_1,\dots,i_k} \frac{1}{k!}\,\Phi_{i_1}\cdots \Phi_{i_k},$$
both summed over the $(i_1,\dots, i_k)$ with $k\ge 1$, $i_1,\dots, i_k\ge 1$ and $i_1+ \cdots + i_k=n$.
\end{lemma}

\begin{example} 
Let $\Phi,\Psi\in\TrEnd^1_K$ with $\Phi=\log(1+\Psi)$. Then
$$\Phi_1 = \Psi_1, \quad 
\Phi_2 = \Psi_2 - \textstyle\frac{1}{2} (\Psi_1)^2, \quad 
\Phi_3 = \Psi_3 - \textstyle\frac{1}{2}(\Psi_1\Psi_2+\Psi_2\Psi_1) + \textstyle\frac{1}{3} (\Psi_1)^3,\quad \dots\  .
$$
\end{example}

\section{The Lie Algebra of an Algebraic Unitriangular Group}\label{sec:algebraic groups}

\index{algebra!nontrivial}
\index{group!unitriangular!algebraic}

\noindent
{\em In this section $K$ is an integral domain.}\/ A $K$-algebra $A$ is said to be
{\bf nontrivial} if~${1_A\ne 0}$. Given a morphism $\phi\colon
A \to B$ of nontrivial commutative $K$-algebras, we extend $\phi$ to
a morphism $\phi\colon \tr_A \to \tr_B$ of $K$-algebras by
$\phi(M):= \big(\phi(M_{ij})\big)$. This extended $\phi$ maps 
$\tr_A^n$ to $\tr_B^n$, restricts to a group
morphism $\mathcal U_A \to \mathcal U_B$, and 
$$ \phi\big(f(M)\big)\ =\ f\big(\phi(M)\big) \quad \text{ for 
$f\in K[[z]]$, $M\in \tr_A^1$.}$$
Let $X=(X_{ij})_{i,j\in\N}$ be a family of distinct indeterminates. A unitriangular group~$\mathcal G$ over $K$ is {\bf algebraic} if for some
family $(P_\alpha)$ of polynomials $P_{\alpha}\in K[X]$: \begin{enumerate}
\item[(1)] $ \mathcal G\ =\ \big\{ G=(G_{ij}) \in K^{\N\times\N}:\  \text{$P_\alpha(G)=0$ for all $\alpha$} \big\}$,
\item[(2)] for each nontrivial commutative $K$-algebra $A$, the set
$$\mathcal G_A\ :=\ \big\{ G=(G_{ij}) \in A^{\N\times\N}:\  \text{$P_\alpha(G)=0$ for all $\alpha$} \big\}$$ 
of common zeros of the $P_\alpha$ in $A$ is a subgroup of $\mathcal U_A=1+\tr_A^1$.
\end{enumerate}
In particular, $\mathcal U=\mathcal U_K$ is algebraic.
Below $\mathcal G$ is such an algebraic unitriangular group over $K$ and~$(P_\alpha)$ is a family of polynomials $P_{\alpha}\in K[X]$ as above. Thus a morphism $\phi\colon
A \to B$ of nontrivial commutative $K$-algebras induces a group morphism
$$\mathcal G_A \to \mathcal G_B\colon  \quad G \mapsto \phi(G).$$ 
We also let $t$ be an indeterminate, and let $\varepsilon$ be the image of $t$ under the natural map $K[t]\to K[t]/(t^2)$, so $K[t]/(t^2)=K[\varepsilon]=K\oplus K\varepsilon$ with $\varepsilon^2=0$ is the $K$-algebra of dual numbers over $K$. 
Thus for $M\in\tr_{K[\varepsilon]}^1$ we have $\exp(\varepsilon M)=1+\varepsilon M$ in $\tr_{K[\varepsilon]}$. 

\begin{lemma}
Let $M\in\tr_K^1$. Then $\exp(tM)\in \mathcal G_{K[t]}\Longleftrightarrow 1+\varepsilon M \in \mathcal G_{K[\varepsilon]}$.
\end{lemma}
\begin{proof}
The forward direction is clear by applying the $K$-algebra morphism 
$$K[t]\to K[\varepsilon]\colon \quad t\mapsto\varepsilon. 
$$
For the converse suppose $1+\varepsilon M \in \mathcal G_{K[\varepsilon]}$. Let $n\geq1$ and let $t_1,\dots,t_n$ be distinct indeterminates with respective images $\varepsilon_1,\dots,\varepsilon_n$ under the natural map 
$$K[t_1,\dots,t_n]\to K[t_1,\dots,t_n]/(t_1^2,\dots,t_n^2).$$ 
Then $\varepsilon_i^2=0$ for $i=1,\dots,n$, and with
$$R_n\ :=\ K[t_1,\dots,t_n]/(t_1^2,\dots,t_n^2)\ =\ K[\varepsilon_1,\dots,\varepsilon_n],$$
the kernel of the $K$-algebra morphism $K[t]\to R_n$ with $t\mapsto \tau:=\varepsilon_1+\cdots+\varepsilon_n$ is generated as an ideal by $t^{n+1}$.
Moreover, the image of $\exp(tM)$ under this 
$K$-algebra morphism  is
\begin{align*}
1+\tau M + \frac{\tau^2}{2!} M^2 + \cdots + \frac{\tau^n}{n!}M^n\ &=\ 1+(\varepsilon_1+\cdots+\varepsilon_n)M+\cdots+(\varepsilon_1\cdots\varepsilon_n)M^n\\
&=\ (1+\varepsilon_1 M)\cdots (1+\varepsilon_n M) \in \mathcal G_{K[\tau]}.
\end{align*}
Since this holds for all $n\ge 1$, we obtain 
$$\exp(tM)\ =\ 1+t M + \frac{t^2}{2!} M^2 + \cdots + \frac{t^n}{n!}M^n + \cdots \in \mathcal G_{K[t]}$$ 
as required.
\end{proof}

\begin{lemma}\label{lem:lie, 1}
The set 
$$\mathfrak g\ :=\ \big\{ M\in\tr_K^1 :\  \exp(tM)\in \mathcal G_{K[t]}\big\}\ =\ \big\{  M\in\tr_K^1 :\  1+\varepsilon M \in \mathcal G_{K[\varepsilon]}\big\}$$
is a Lie subalgebra of $\tr_K^1$.
\end{lemma}
\begin{proof}
Let  $M\in\tr_K^1$. Given a polynomial $P\in K[X]$ with $X=(X_{ij})$ as before, we have in $K[\epsilon]$, by Taylor expansion and $\varepsilon^2=0$,
$$P(1+\varepsilon M)\  =\  P(1) + \varepsilon\sum_{i,j} \frac{\partial P}{\partial X_{ij}}(1)\cdot M_{ij}.$$
Note that $\frac{\partial P}{\partial X_{ij}}\ne 0$ for only finitely
many $(i,j)\in \N^2$. 
Now $1\in \mathcal G$ gives $P_\alpha(1)=0$ for each $\alpha$, and so
$$M\in\mathfrak g \quad\Longleftrightarrow\quad 1+\varepsilon M\in\mathcal G_{K[\varepsilon]} \quad\Longleftrightarrow\quad \sum_{i,j} \frac{\partial P_{\alpha}}{\partial X_{ij}}(1)\cdot M_{ij}=0 \text{ for all $\alpha$.}$$
Thus $\mathfrak g$ is a submodule of the $K$-module 
$\tr_K^1$. 
Let $M,N\in\mathfrak g$; to show that  $[M,N]\in\mathfrak g$, we let $R:=R_2=K[\varepsilon_1,\varepsilon_2]$ be as in the proof of the previous lemma, so $G:=1+\varepsilon_1 M \in \mathcal G_{K[\varepsilon_1]}\subseteq\mathcal G_{R}$ and 
$H:=1+\varepsilon_2 N \in \mathcal G_{K[\varepsilon_2]}\subseteq\mathcal G_{R}$. Then in $\tr^1_{R}$ we have
$$GH = 1 + \varepsilon_1 M + \varepsilon_2 N + \varepsilon_1\varepsilon_2 MN, \quad
  HG = 1 + \varepsilon_1 M + \varepsilon_2 N + \varepsilon_1\varepsilon_2 NM,$$
so
$GH = HG\big(1+\varepsilon_1\varepsilon_2 [M,N]\big)$.
Then $GH,HG\in\mathcal G_{R}$ gives $1+\varepsilon_1\varepsilon_2 [M,N]\in \mathcal G_{K[\varepsilon_1\varepsilon_2]}$. Applying the  $K$-algebra isomorphism $K[\varepsilon]\to K[\varepsilon_1\varepsilon_2]$ with $\varepsilon\mapsto\varepsilon_1\varepsilon_2$ then yields $[M,N]\in \mathfrak g$. Hence $\mathfrak g$ is a Lie subalgebra of $\tr^1_K$. 
\end{proof}

\noindent
The next lemma shows that $\mathfrak g$ depends only on $\mathcal G$, not on the particular family $(P_\alpha)$. We call $\mathfrak g$ the {\bf Lie algebra} of $\mathcal G$, and consider it as a Lie subalgebra of $\tr_K^1$. Note that if $M\in\mathfrak g$ then $\exp(tM)\in \mathcal G_{K[t]}$, and substitution of $t=1$ yields $\exp(M)\in\mathcal G$. Thus $\exp(\mathfrak g)\subseteq\mathcal G$. Here is how $\mathcal G$ and $\mathfrak g$ determine each other:

\index{Lie algebra!algebraic unitriangular group}

\begin{lemma}\label{lem:lie, 2}
$\exp(\mathfrak g)=\mathcal G$ and $\log(\mathcal G)=\mathfrak g$.
\end{lemma}
\begin{proof}
Let $G\in\mathcal G$; it suffices to show that then $\log(G)\in\mathfrak g$, i.e., $\exp\!\big(t\log(G)\big)\in\mathcal G_{K[t]}$. Now for each $\alpha$, the polynomial $P_\alpha\big(\!\exp(t\log(G))\big)\in K[t]$ vanishes upon substitution of integers for $t$, since $\exp(k\log(G))=\exp(\log(G^k))=G^k\in\mathcal G$ for each $k\in\Z$. Since $K$ is assumed to be an integral domain, we therefore have $P_\alpha\big(\!\exp(t\log(G))\big)=0$ for each $\alpha$ and so $\exp\!\big(t\log(G)\big)\in\mathcal G_{K[t]}$.
\end{proof}

\noindent
The Lie algebra of the algebraic unitriangular  group $\mathcal U$ over $K$ is $\tr_K^1$; we denote this Lie algebra also by $\mathfrak u$. We equip the $K$-module $\mathfrak u$ with the complete nonnegative fil\-tra\-tion~$(\mathfrak u^n)$ given by $\mathfrak u^n:=\tr_K^{n+1}$ for each $n$. This makes $\mathfrak u$ a {\em filtered\/} Lie algebra over $K$. The Lie subalgebra $\mathfrak g$ of $\mathfrak u$ is made into a filtered Lie algebra over $K$ by giving it the filtration induced by 
$(\mathfrak u^n)$.

\begin{lemma}\label{lem:lie, 3} $\mathcal G$ and $\mathfrak g$ are closed in $\tr_K$. In particular,
the filtered Lie algebra $\mathfrak g$ over $K$ is complete.
\end{lemma}

\begin{samepage}
\begin{proof} For each $P\in K[X]$, the function $G\mapsto P(G)\colon \tr_K \to K$ is clearly locally constant, and so its zero set is closed (and open) in $\tr_K$. Thus $\mathcal G$ is closed in $\tr_K$, and so is 
$\mathfrak g$ in view of the equivalence
$$M\in \mathfrak g\ \Longleftrightarrow\ \sum_{i,j}\frac{\partial P_{\alpha}}{\partial X_{ij}}(1)\cdot M_{ij}=0\ \text{ for all  $\alpha$}$$
from the proof of Lemma~\ref{lem:lie, 1}.
\end{proof}
\end{samepage}

\noindent
Section~\ref{sec:Iteration matrices} below is devoted to the investigation of the Lie algebra of a certain algebraic unitriangular  group over $\Q$ which plays an important role in the study of triangular automorphisms of differential polynomial rings in Section~\ref{appdifpol}.

\subsection*{Notes and comments}
The results in this section are analogues of well-known facts about algebraic groups of (finite-size) matrices. The proof of Lemma~\ref{lem:lie, 1} follows \cite[Theorem~I.5]{Serre}.

\section{Derivations on the Ring of Column-Finite Matrices}\label{dercolfinma}

\noindent 
In this section we investigate two kinds of derivations on the ring of column-finite matrices: those induced by derivations on $K$ and the internal $K$-derivations.

\medskip
\noindent
First, let a derivation $\der$ on $K$ be given. This is a $\Q$-derivation
in the sense of Section~\ref{sec:Filtered and Graded Algebras}. For $M=(M_{ij})\in K^{\N\times\N}$ we set 
$\der(M):=\big(\der(M_{ij})\big)$.  If $M$ is column-finite, then so is $\der(M)$, and $M\mapsto\der(M)$ is a derivation on the $\Q$-algebra of column-finite matrices which restricts to a derivation on its $\Q$-subalgebra $\tr_K$ consisting of the triangular matrices over $K$. 

\medskip\noindent
Now let $t$ be an indeterminate.
Then $\tr_K$ is a $K$-subalgebra of $\tr_{K[t]}$. We equip~$K[t]$ with the derivation $\frac{d}{dt}$, and accordingly we define 
$\frac{dM}{dt}\in\tr_{K[t]}$ for $M\in\tr_{K[t]}$ as just explained. This gives a $K$-derivation $M\mapsto \frac{dM}{dt}$ of rank 
$0$ on the $K$-algebra~$\tr_{K[t]}$ equipped with the filtration $(\tr_{K[t]}^n)$. 
The following two lemmas about this $K$-derivation 
are used in Section~\ref{sec:Iteration matrices}.

\begin{lemma}\label{lem:derivative of exp(tM)}
Let $M\in\tr^1_K$. Then
$$\frac{d}{dt} \ex^{tM}\ =\ \ex^{tM}M.$$
\end{lemma}
\begin{proof}
We have $(tM)^n=t^nM^n$ for every $n$, hence
$$\ex^{tM}\ =\ \sum\frac{(tM)^n}{n!}\ =\ \sum \frac{t^n M^n}{n!}$$
and thus
\equationqed{\frac{d}{dt}\ex^{tM}\ =\ \sum \frac{d}{dt}\left(\frac{t^n M^n}{n!}\right)\ =\ \sum_{n=1}^{\infty} \frac{t^{n-1}M^n}{(n-1)!}\ =\ \ex^{tM}M.}
\end{proof}

%(Similarly, of course, one also sees $\frac{d}{dt} \exp(tM) = M\exp(tM)$, but we won't need this fact.)

\noindent
The following lemma is a familiar fact about systems of linear differential equations with constant coefficients:

\begin{lemma}\label{lem:linear DE}
Let $M\in\tr^1_K$ and $Y\in\tr_{K[t]}$, so $Y(0)\in \tr_K$. Then
$$\frac{dY}{dt}\ =\ YM \qquad\Longleftrightarrow\qquad Y\ =\ Y(0)\ex^{tM}.$$
\end{lemma}
\begin{proof}
Lemma~\ref{lem:derivative of exp(tM)} shows that if $Y=Y(0)\ex^{tM}$, then $\frac{dY}{dt}=YM$. Conversely, suppose $\frac{dY}{dt}=YM$. Replacing
$Y$ by $Y-Y(0)\ex^{tM}$ we arrange $Y(0) = 0$; we need to show that then $Y=0$. Towards a contradiction, suppose $Y\neq 0$. Now $Y(0) = 0$, so for all $i$, $j$ with $Y_{ij}\neq 0$ we have
$Y_{ij}=t^{n_{ij}}Z_{ij}$ with $n_{ij}\in\N$, $n_{ij}\ge 1$, and $Z_{ij}\in K[t]$, $Z_{ij}(0)\neq 0$.
Pick $i$, $j$ with minimal $n_{ij}$. Then $\frac{dY}{dt}=YM$ gives
$$n_{ij}t^{n_{ij}-1}Z_{ij}+t^{n_{ij}}\frac{dZ_{ij}}{dt}=\frac{dY_{ij}}{dt} = \sum_k Y_{ik} M_{kj} = \sum_{Y_{ik}\neq 0} t^{n_{ik}} Z_{ik} M_{kj}. $$
Here the right-hand side is divisible in $K[t]$ by $t^{n_{ij}}$ while the
left-hand side is not, a contradiction. 
\end{proof}

\index{matrix!shift}
\index{shift!matrix}
\index{derivative!matrix}
\index{matrix!derivative}

\noindent
Secondly, given a column-finite matrix $A\in K^{\N\times\N}$, the adjoint $M\mapsto \operatorname{ad}_A M=[A,M]$ of $A$ is a $K$-derivation on the $K$-algebra of column-finite matrices. 
% that is, if  $M,N\in K^{\N\times\N}$ are column-finite 
%and $c\in K$, then
%\begin{align*}
%\operatorname{ad}_A(cM)&=c\operatorname{ad}_A M,  \quad \operatorname{ad}_A(M+N)=\operatorname{ad}_A M+\operatorname{ad}_A N, \\  \operatorname{ad}_A(MN)&=\operatorname{ad}_A(M)N+M\operatorname{ad}_A(N), \quad  \operatorname{ad}_A [M,N]=[\operatorname{ad}_A M,N]+[M,\operatorname{ad}_A N].
%\end{align*}
In the rest of this section we employ the derivation $\operatorname{ad}_A$, for a particular choice of $A$, to establish commutator identities for certain diagonal matrices used in later sections. 

\begin{definition} \label{def:shift matrix}
The (column-finite) {\bf shift matrix} $S\in K^{\N\times\N}$ is
given by 
$$S_{j+1,j}\ =\ 1, \qquad S_{i,j}\ =\ 0\ \text{  for $i\ne j+1$.}$$  
A column-finite matrix $M\in K^{\N\times\N}$ has {\bf derivative} 
$M'\in K^{\N\times\N}$ given by
$$M'\ :=\ \operatorname{ad}_{-S} M\ =\ MS-SM$$
and so $M'$ is also column-finite.
\end{definition}

\noindent
Multiplying a column-finite matrix $M\in K^{\N\times\N}$ by $S$   has the following effect:
\begin{align*}
(MS)_{ij}=M_{i,j+1}: & \quad \text{shifts $M$ one column to the left, cancels the leftmost column,} \\
(SM)_{ij}=M_{i-1,j}: & \quad \text{shifts $M$ one row downwards, adds a top row of zeros.} 
\end{align*}
Here and below $M_{ij}:=0$ if $i<0$ or $j<0$. In particular, 
$$(M')_{ij}\ =\ M_{i,j+1}-M_{i-1,j}.$$

%$$TS=\begin{pmatrix} 
%1 &   &   & \\
%  & 1 &   & \\
%  &   & 1 & \\
%  &   &   & \ddots
% \end{pmatrix}, \qquad ST = \begin{pmatrix} 
%0 &   &   & \\
%  & 1 &   & \\
%  &   & 1 & \\
%  &   &   & \ddots
% \end{pmatrix}.$$

\begin{example}
Let $n\ge 1$, $a=(a_i)\in K^\N$, and $M=\diag_n a$. Then
$$(\diag_n a)'\  =\  \diag_{n-1}(a_0,a_1-a_0,a_2-a_1,\dots).$$
%and if $k\leq 0$ then
%$$(\diag_k a)' = \diag_{k-1}(a_1-a_0,a_2-a_1,a_3-a_2,\dots).$$
\end{example}

\noindent
The next lemma lists some properties of the derivation $M\mapsto M'$.

\begin{lemma}\label{lem:properties of matrix derivative}
Let $M\in\tr_K$. Then 
$$M'\ =\ 0\qquad\Longleftrightarrow\qquad\text{$M\in K$,}$$
and for
$M\in\tr_K^\times$, we have 
$$%\begin{equation}\label{eq:derivative of inverse matrix}
(M^{-1})'\ =\ -M^{-1}M'M^{-1}.
$$%\end{equation}
Suppose now that $M\in \tr^1_K$. Then $M'\in \tr_K$, and for $n\ge 1$,
$$%\begin{equation}\label{eq:derivative of diagonals}
 (M_n)'\ =\ (M')_{n-1}, \qquad   M\in\tr_K^n\ \Longrightarrow\ M'\in\tr_K^{n-1}.
$$%\end{equation}
\end{lemma} 

%\index{finite!column-vector}

\noindent
Below, elements of $K^\N$ are column vectors, and $a=(a_i)\in K^\N$ is called 
{\bf  finite\/} if $a_i=0$ for all but finitely many $i$. So  
$e:=(1,0,0,\dots)^{\operatorname{t}}\in K^\N$ is finite. If  $M\in K^{\N\times\N}$ 
is column-finite and $a\in K^\N$ is finite, then  $Ma\in K^{\N}$ (defined in the obvious way) is finite, and in particular, 
$Me$ is the leftmost column of $M$.

\begin{lemma}\label{lem:recurrence}
Suppose $a\in K^\N$ is finite, $B\in\tr_K$, and $A,C\in K^{\N\times\N}$ are 
column-finite. Then there is a unique column-finite matrix $X\in K^{\N\times\N}$ 
such that
\begin{align}
Xe\ &=\ a, \label{eq:X1} \\ 
X'\  &=\ AXB+C. \label{eq:X2} 
\end{align}
If also $A,C\in\tr_K$ and $a_i=0$ for all $i\ge 1$, then $X\in\tr_K$.
\end{lemma}
\begin{proof} Suppose $X\in K^{\N\times\N}$ is column-finite and 
satisfies \eqref{eq:X1} and \eqref{eq:X2}. Then the leftmost column of 
$X$ is $a$. Let $j$ be given.  Then for each $i$,
$$X_{i,j+1}\ =\ (XS)_{ij}\ =\ (X'+SX)_{ij}\ =\ (AXB)_{ij}+C_{ij}+X_{i-1,j}.$$
Since $B$ is triangular, the sum $(AXB)_{ij} = \sum_{k,l} A_{ik}X_{kl}B_{lj}$  
only involves entries of~$X$ from its columns with indices  $l=0,\dots,j$. 
Thus the column $(X_{i,j+1})_i$ of $X$ with index $j+1$  is determined by 
its columns with lower index, and so there is at most one $X$ as claimed. 
If in addition $A, C\in\tr_K$ and $a_i=0$ for all $i\ge 1$, 
then an induction on $j$ shows that $X_{ij}=0$ for all $i>j$. 

Reversing these considerations we construct an $X$ as claimed.
\end{proof}

\noindent
For the next lemmas and corollaries (used in Sections~\ref{sec:Iteration matrices} and~\ref{sec:riordan}), recall that 
$$\binom{X}{n}\ :=\ \frac{X(X-1)\cdots (X-n+1)}{n!}\in \Q[X]$$
is a polynomial of degree $n$, with $\binom{X}{0}=1$. For $k\in \Z$ we let $\binom{k}{n}$ be
its value at $k$; so $\binom{k}{n}=0$ if $0\le k < n$. It is also
convenient to set $\binom{n}{-1}:= 0$. 
%for any $k\in \Z$ we let  we set, for $k\in \Z$, 
%$$\binom{n}{k}:= 0\ \text{ if $k<0$ or $k>n$.}$$
Consider now
$$A(n)\ :=\ \diag_n {i+n \choose n},$$
so $A(0)=1$ and $A(1)=\diag_1(1,2,3,\dots)$. 
Note that for all $n$,  
\begin{equation}\label{eq:derivative of A(n)}
A(n+1)'\ =\ A(n), 
\end{equation}
by the familiar recurrence relations for binomial coefficients. 
It is easy to verify, using Lemma~\ref{lem:formulas for diagonals}, that
$\big[A(m),A(n)\big]=0$ for all $m$, $n$.
We also set
$$B(n)\ :=\ \diag_n \binom{i+n}{n+1},$$
so $B(0)=\diag(0,1,2,3,\dots)$. For all $n$,
\begin{equation}\label{eq:derivative of B(n)}
B(n+1)'\ =\ B(n),
\end{equation}
again by the recurrence relation for binomial coefficients.

\begin{lemma}\label{lem:binomial diagonal matrices}
For all $m$, $n$,
\begin{equation}
\big[B(m),B(n)\big]\	 =\  \left(\binom{m+n}{m-1}-\binom{m+n}{m+1}\right)\,B(m+n) \label{eq:Lie bracket of B(n)}
\end{equation}
and
\begin{equation} \label{eq:Lie bracket of A(m), B(n)}
\big[A(m),B(n)\big]\ =\ {m+n\choose n+1} A(m+n).
\end{equation}
\end{lemma}

\begin{proof} A simple computation yields $[B(0),B(n)]=-nB(n)$, and 
likewise we have $[B(m),B(0)]=-[B(0),B(m)]=mB(m)$ and
$[A(0),B(n)]=0$, $[A(m),B(0)]=mA(m)$. Thus \eqref{eq:Lie bracket of B(n)} and \eqref{eq:Lie bracket of A(m), B(n)}
hold if $m=0$ or $n=0$. Let $m\ge 1$ and $n\ge 1$. Then by \eqref{eq:derivative of B(n)},
\begin{align*}
\big[B(m),B(n)\big]'\  &=\ \big[B(m)',B(n)\big]+\big[B(m),B(n)'\big] \\
		&=\ \big[B(m-1),B(n)\big] + \big[B(m),B(n-1)\big].
\end{align*}
Inductively we can assume that the last sum equals the sum of
\begin{align*}&\left(\binom{m+n-1}{m-2}-\binom{m+n-1}{m}\right)\,B(m+n-1)\ \text{ and}\\  &\left(\binom{m+n-1}{m-1}-\binom{m+n-1}{m+1}\right)\,B(m+n-1),
\end{align*}
so $\big[B(m),B(n)\big]'=\left(\binom{m+n}{m-1}-\binom{m+n}{m+1}\right)\, B(m+n)'$ by \eqref{eq:derivative of B(n)}. Now \eqref{eq:Lie bracket of B(n)} follows from Lemma~\ref{lem:properties of matrix derivative}. 
Similarly for~\eqref{eq:Lie bracket of A(m), B(n)}, we have by  \eqref{eq:derivative of B(n)} and \eqref{eq:derivative of A(n)}:
\begin{align*}
\big[A(m),B(n)\big]'\  &=\ \big[A(m)',B(n)\big]+\big[A(m),B(n)'\big] \\
		&=\ \big[A(m-1),B(n)\big] + \big[A(m),B(n-1)\big].
\end{align*}
Inductively, we can assume that the last sum equals 
$$ \binom{m+n-1}{n+1}\,A(m+n-1) +   \binom{m+n-1}{n}\,A(m+n-1),$$
so $\big[A(m),B(n)\big]'=\binom{m+n}{n+1}\, A(m+n)'$ by \eqref{eq:derivative of A(n)}.  Now use Lemma~\ref{lem:properties of matrix derivative}.
\end{proof}

\noindent
From \eqref{eq:Lie bracket of B(n)} we obtain:

\begin{cor}\label{cor:binomial diagonal matrices}
Let $c_1, c_2\in K$, and define the sequence $(C(n))_{n\geq 1}$ in $\tr_K$ by 
$$\begin{cases}
C(n)=c_nB(n)    &\text{ for $n=1,2$}, \\ 
C(n+1)=\big[C(1),C(n)\big]&\text{ for $n\ge 2$.}
\end{cases}$$ 
Then
$C(n)=c_nB(n)$ for all $n\ge 1$, where 
$c_{n+1}=\left(1-\binom{n+1}{2}\right)c_1c_n$ for $n\ge 2$.
\end{cor}

\begin{remark}
An easy induction on $n\geq 2$ shows that the terms $c_2,c_3,\dots$ of the se\-quence~$(c_n)_{n\ge 1}$ in the previous lemma are explicitly given by
$$c_n = (-c_1)^{n-2} c_2 \cdot \frac{(n-2)!\,(n+1)!}{3\cdot 2^{n-1}}\qquad\text{for $n\geq 2$.}$$
\end{remark}

\subsection*{Notes and comments}
The derivative of a column-finite matrix  is defined in \cite{Kemeny}, and used there to give simple proofs for combinatorial identities.

\section{Iteration Matrices}\label{sec:Iteration matrices}

\noindent
In this section we introduce special types of triangular matrices, called iteration matrices, and study their matrix logarithms. 

\subsection*{Bell polynomials} Let $x, y_1, y_2, y_3,\dots, z$ be distinct
indeterminates, and set 
$$R\ :=\ \Q[x,y_1,y_2, y_3,\dots], \qquad A\ :=\ R[[z]].$$ We view the power series ring $A$ as a complete filtered $A$-algebra with respect to
the nonnegative filtration $(A^n)$ given by $A^n=z^nA$.
Set
$$y\ :=\ \sum_{n= 1}^{\infty} y_n \frac{z^n}{n!}\in zR[[z]]\ =\ A^1,$$
so $xy\in A^1$, and $\exp(xy)\in 1+A^1$. Here is an explicit
formula for $\exp(xy)$, where for $\k=(k_1,\dots,k_{d})\in\N^{d}$,  $d\geq 1$, we set 
$$\abs{\k}:=k_1+\cdots+k_d, \qquad \dabs{\k}:=k_1+2k_2+\cdots+dk_d.$$

\begin{prop} \label{prop:explicit fm for Bell polys}
In $A$ we have the identity
\begin{align*} \exp(xy)\ &=\ \sum_{j=0}^{\infty} \left(\sum_{i=0}^j B_{ij}x^i\right)\frac{z^j}{j!},\ \text{ where}\\ 
B_{ij}\ =\ \sum_{\substack{\k=(k_1,\dots,k_{d})\in\N^{d} \\ \abs{\k} = i, \dabs{\k}=j}} &\frac{j!}{k_1!k_2!\cdots k_{d}!\cdot (1!)^{k_1}(2!)^{k_2}\cdots (d!)^{k_d}}\, y_1^{k_1}y_2^{k_2}\cdots y^{k_{d}}_{d},
\end{align*} 
for $i\leq j$ and $d=j-i+1$. In particular, for such $i,j,d$ we have $B_{ij}\in\Q[y_1,\dots,y_{d}]$, and the coefficients of $B_{ij}$ are in $\N$. 
\end{prop}
\begin{proof} With $k_1,k_2,\dots$ ranging over $\N$, the multinomial identity gives:
\begin{align*}
\exp(x\cdot y) 	&= \sum_{i=0}^{\infty} \frac{x^i}{i!} \left(\sum_{n=1}^{\infty} y_n\frac{z^n}{n!}\right)^i \\
				&= \sum_{i= 0}^{\infty} \frac{x^i}{i!} \left(\sum_{k_1+k_2+\cdots=i} \frac{i!}{k_1!k_2!\cdots} \left(\frac{y_1z^1}{1!}\right)^{k_1} \left(\frac{y_2z^2}{2!}\right)^{k_2}\cdots \right) \\
				&= \sum_{j=0}^{\infty}\left(\sum_{i=0}^j x^i\!\left(\sum_{\substack{k_1+k_2+\cdots = i \\ k_1+2k_2+\cdots=j}} \frac{j!}{k_1!k_2!\cdots \cdot (1!)^{k_1}(2!)^{k_2}\cdots }\, y_1^{k_1}y_2^{k_2}\cdots\! \right)\!\right)\!\frac{z^j}{j!}.
\end{align*}
The condition on the infinite sequence $k_1, k_2,\dots$ in the innermost
sum forces $k_n=0$ for $n> d:= j-i+1$, which 
yields the displayed formula for $B_{ij}$. 
By the following lemma, the coefficients of $B_{ij}$ are in $\N$.
%The coefficients of $B_{ij}$ are in $\N$, since for 
%$\k=(k_1,\dots,k_d)\in\N^d$, $d\geq 1$, with $\dabs{\k}=j$, 
%the fraction  $\frac{j!}{k_1!\cdots k_{d}!\cdot (1!)^{k_1}\cdots (d!)^{k_d}}$ 
%is the number of partitions of a $j$-element set  into exactly $k_1$ 
%sets of cardinality $1$, $k_2$ sets of cardinality $2$, 
%and so on; see for example \cite{Comtet-Book}.
\end{proof}

\begin{lemma}
Let $\k=(k_1,\dots,k_d)\in\N^d$, $d\geq 1$, and set $n=\dabs{\k}$. The number of partitions of an $n$-element set
into exactly $k_1$ sets of cardinality $1$, $k_2$ sets of car\-di\-na\-li\-ty~$2$, etc., is
$\displaystyle\frac{n!}{k_1!\cdots k_{d}!\cdot (1!)^{k_1}\cdots (d!)^{k_d}}$.
\end{lemma}
\begin{proof} Set $[n]=\{1,\dots,n\}$. Consider first the special case
$n=kd$ with $k,d\ge 1$. 
We claim that the number of partitions of $[n]$ into $k$ sets of size $d$
is $\frac{n!}{k!(d!)^k}$. This claim is obviously true for $k=1$.
Let $k>1$ and assume inductively that the claim holds with $k-1$ instead of $k$.
Then the claim follows by the inductive assumption and the fact that
there are $\binom{n-1}{d-1}$ subsets of $[n]$ of size $d$ containing $1$. 
As to the general case,
the lemma clearly holds for $d=1$, so let $d>1$, and assume the lemma
holds for $d-1$ instead of $d$. This inductive assumption takes care of the case $k_d=0$, so
let $k_d\ge 1$. Then use that there are $\binom{n}{k_d d}$
subsets of $[n]$ of size $k_d d$, apply the claim above and the
inductive assumption, and perform a routine computation. 
\end{proof}

\begin{cor}\label{corprop} Let $b,c_1,c_2, c_3,\ldots\in K$ and set
$$c\ :=\ \sum_{n=1}^{\infty}c_n\frac{z^n}{n!}\in zK[[z]].$$
Then $\ex^{bc}\ =\ \sum_{j=0}^{\infty}\left(\sum_{i=0}^j
B_{ij}(c_1,\dots,c_{j-i+1})b^i\right)\cdot\frac{z^j}{j!}$ in $K[[z]]$.
\end{cor}

\noindent
The $B_{ij}=B_{ij}(y_1,y_2,\dots)\in\Q[y_1,y_2,\dots]$ are the (partial) {\bf Bell polynomials.} \index{polynomial!Bell}\index{Bell polynomials}\nomenclature[Y]{$B_{ij}$}{Bell polynomials} We also set $B_{ij}:= 0$ for $i>j$,
so the family $(B_{ij}x^iz^j/j!)_{i,j\in \N}$ is summable in $A$, with
\begin{equation}\label{eq:Bell polynomials}
\sum_{i\in \N} \frac{x^iy^i}{i!}\  =\ 
\sum_{i,j\in\N} B_{ij}x^i\frac{z^j}{j!}.
\end{equation}
Note that $B_{0j}=0$ and $B_{1j}=y_j$ for $j\geq 1$, and $B_{jj}=y_1^j$ for all $j$. Next we establish the following identity in $\Q[y_1, y_2,\dots][[z]]$, holding for all $i\in \N$:
\begin{equation}\label{eq:power of y}
\frac{y^i}{i!} \ =\ \sum_{j=0}^{\infty} B_{ij} \frac{z^j}{j!}.
\end{equation}
To see why this identity holds we try to view both sides in~\eqref{eq:Bell polynomials} as power series in $x$ and then compare the coefficients of $x^i$. To justify this idea, we note that $R$ is a subring of $S:=\Q[y_1,y_2,\dots][[x]]$, and accordingly, $A=R[[z]]$ is a subring of 
$$B\ :=\  S[[z]]\ =\ \Q[y_1, y_2, \dots][[x,z]]\ =\ \big(\Q[y_1, y_2,\dots][[z]]\big)[[x]],$$
a complete filtered $B$-algebra with respect to the nonnegative
filtration $(B^n)$ given by $B^n=(x,z)^nB$. Then $A^n\subseteq B^n$ for all $n$ and $xy\in B^2$, so
$\sum (x^iy^i)/i!$ takes the value in $B$ that it has in $A$. Likewise,
the right-hand sum in~\eqref{eq:Bell polynomials} is defined in~$B$, and
takes the value in $B$ it has in $A$. Now~\eqref{eq:power of y} follows as
indicated above.

\medskip\noindent
The next recursion formula facilitates the computation of the Bell polynomials: 
\begin{lemma}
Suppose $i_1\leq i\leq j$. Then in $\Q[y_1, y_2,\dots]$ we have
$$B_{ij}\ =\ \frac{1}{{i\choose i_1}} \sum_{j_1=0}^{j} {j\choose j_1}  B_{i_1,j_1} B_{i-i_1,j-j_1}.$$
\end{lemma}
\begin{proof}
By \eqref{eq:power of y},
$$\frac{1}{i_1!}y^{i_1}\ =\ \sum_{j=0}^{\infty} B_{i_1,j} \frac{z^j}{j!}, \qquad 
  \frac{1}{(i-i_1)!}y^{i-i_1}\ =\ \sum_{j=0}^{\infty} B_{i-i_1,j} \frac{z^j}{j!},$$
hence in $\Q[y_1, y_2,\dots][[z]]$ we have
\begin{align*}
{i\choose i_1} \cdot \frac{1}{i!}y^i\ =\ \frac{1}{i_1!}y^{i_1} \cdot \frac{1}{(i-i_1)!}y^{i-i_1} 
									&=\ \left(\sum_{j= 0}^{\infty} B_{i_1,j} \frac{z^j}{j!}\right)\cdot\left(\sum_{j=0}^{\infty} B_{i-i_1,j} \frac{z^j}{j!}\right) \\
									&=\ \sum_{j=0}^{\infty} \left( \sum_{j_1=0}^{j} {j\choose j_1}  B_{i_1,j_1} B_{i-i_1,j-j_1}\right) \frac{z^j}{j!},
\end{align*}
and the lemma follows.
\end{proof}

\noindent
Using this with $i_1=1$ we get easily 
\begin{align*} B_{23}\ &=\ 3y_1y_2,\quad B_{24}\ =\ 4y_1y_3+3y_2^2,\quad \qquad B_{25}\ =\ 5y_1y_4+10y_2y_3,\\ 
B_{34}\ &=\ 6y_1^2y_2,\quad B_{35}\ =\ 10y_1^2y_3+15y_1y_2^2,\quad B_{45}\ =\ 10y_1^3y_2.
\end{align*}

\subsection*{Iteration matrices}
Given a power series $f\in z K[[z]]$,
\begin{equation}\label{eq:f}
f\ =\ \sum_{n\geq 1} f_n\frac{z^n}{n!}\qquad\text{($f_n\in K$ for each $n\geq 1$),}
\end{equation}
we introduce the triangular matrix
\begin{multline*}
\llb f\rrb\ :=\ \big(\llb f \rrb_{ij}\big)_{i,j\in\N}\ =\ \big(B_{ij}(f_1,f_2,\dots,f_{j-i+1})\big)_{i,j\in\N}\ = \\ 
\begin{pmatrix}
1 & 0   & 0     & 0   		& 0   			& 0   					& \cdots \\
  & f_1 & f_2   & f_3 		& f_4 			& f_5					& \cdots \\
  &     & f_1^2 & 3f_1f_2 	& 4f_1f_3+3f_2^2	& 5f_1f_4+10f_2f_3		& \cdots \\
  &     &       & f_1^3      & 6f_1^2f_2		& 10f_1^2f_3+15f_1f_2^2	& \cdots \\
  &     &       &            & f_1^4         & 10f_1^3f_2				& \cdots \\
  &		&		&			&				& f_1^5					& \cdots \\
  &		&		&			&				&						&\ddots
\end{pmatrix}\in\tr_K.
\end{multline*}
Note that \eqref{eq:power of y} gives $\frac{f^i}{i!}=\sum_{j=0}^{\infty}\, \llb f \rrb_{ij}\frac{z^j}{j!}$ in $K[[z]]$. As we see from the second row of the display, the map $f\mapsto \llb f\rrb\colon zK[[z]]\to \tr_K$ is injective. It is also easy to check that $\llb z\rrb = 1$. 
The matrix $\llb f\rrb$ is called the {\bf iteration matrix of $f$}, since $f\mapsto \llb f  \rrb$ converts composition of power series into matrix multiplication: 

\index{matrix!iteration}
\index{iteration matrix}
\nomenclature[Y]{$\llb f\rrb$}{iteration matrix of the power series $f$}

\begin{lemma}\label{lem:iteration matrix}
Let $f,g\in zK[[z]]$. Then
$\ \llb f\circ g \rrb\ =\ \llb f\rrb\cdot \llb g\rrb$.
\end{lemma}

\begin{proof} The above identity for powers of elements in $zK[[z]]$ gives
\begin{align*}
\sum_{k= 0}^{\infty}\ \llb f\circ g\rrb_{ik}\frac{z^k}{k!}\ 	&=\  \frac{1}{i!} (f\circ g)^i\ =\  \frac{1}{i!}\ f^i\circ g\ 
													=\  \sum_{j=0}^{\infty}\ \llb f\rrb_{ij}\, \frac{g^j}{j!} \\
	&=\  \sum_{j= 0}^{\infty} \left(\llb f\rrb_{ij} \sum_{k=0}^{\infty}\ \llb g\rrb_{jk}\ \frac{z^k}{k!}\right)\ 
													=\ \sum_{k= 0}^{\infty} \left(\sum_{j=0}^{\infty} \ \llb f\rrb_{ij}\llb g\rrb_{jk}\right) \frac{z^k}{k!}. 
\end{align*}
Now compare the coefficients of $z^k/k!$ in the first and the last sum.
\end{proof}

\index{unitary!power series}
\index{group!iteration matrices}
\nomenclature[Y]{$\mathcal I$}{group of iteration matrices}

\noindent
The subset $zK^\times + z^2K[[z]]$ of $zK[[z]]$ is a group under formal composition with identity element $z$, and it admits a group embedding 
$$f\mapsto \llb f\rrb\quad \colon\quad zK^\times + z^2K[[z]]\ \to\ \tr_K^\times$$ 
into the group $\tr_K^\times$ of units of $\tr_K$. 
We say that $f$ as in \eqref{eq:f} is {\bf unitary} if $f_1=1$. The set $z+z^2 K[[z]]$ of unitary power series in $K[[z]]$ is a subgroup of $zK^\times + z^2K[[z]]$ under composition, whose image under  $f\mapsto \llb f\rrb$ 
is a subgroup of $\mathcal U=1+\tr_K^1$ which we denote by $\mathcal I$ and call  the {\bf group of iteration matrices over $K$}. 
%Given $f\in K[[z]]$ of the form $f=z+z^{n+1}g$ with $n\ge 1$ and $g\in K[[z]]$, $g\notin zK[[z]]$, we say that the {\bf iterative valuation} of $f$ is $n$; in symbols: $n=\operatorname{itval}(f)$. 
It is easy to see that for $f\in zK[[z]]$ and $n\ge 1$, we have 
$$f\in z+z^{n+1}K[[z]]\ \Longleftrightarrow\ \llb f\rrb\in 1+\tr_K^n.$$

%For each $n>0$ we define the subgroup
%$$\mathcal I_K^{n} := \mathcal I_K\cap (1+\tr_K^n) = \big\{ \llb f\rrb: f\in z+z^{n+1}K[[z]] \big\}$$
%of $\mathcal I_K$. Then 
%$$\mathcal I_K=\mathcal I_K^{1}\supseteq \mathcal I_K^{2}\supseteq\cdots\supseteq \mathcal I_K^{n}\supseteq\cdots\quad\text{and}\quad \bigcap_{n>0} \mathcal I_K^{n}=\{1\},$$ 
%and if $f\in zK[[z]]$ is unitary with $f\neq z$, then $n=\operatorname{itval}(f)$ is the unique $n>0$ such that $\llb f\rrb\in\mathcal I_K^{n}\setminus \mathcal I_K^{n+1}$.

\subsection*{The Lie algebra of the group of iteration matrices} {\em In the
rest of this section~$K$ is an integral domain}.  It is easy to check
that then the  unitriangular group~$\mathcal I$ over $K$ is algebraic:
use the way that the entries of an arbitrary element are given by polynomials in the entries of its second row. Thus $\mathcal I$ 
has an associated Lie algebra by Lemma~\ref{lem:lie, 2} and the remarks preceding it.
Our next goal is to give an explicit description of this Lie algebra.

\index{matrix!iteration!infinitesimal}
\index{iteration matrix!infinitesimal}
\nomenclature[Y]{$\lla h\rra$}{infinitesimal iteration matrix of the power series $h$}

\begin{definition}\label{definfiteration}
Let $h=\sum_{n=1}^{\infty} h_n \frac{z^{n}}{n!}\in zK[[z]]$, $h_n\in K$ for  $n\geq 1$. The {\bf infinitesimal iteration matrix} of $h$ is the triangular matrix
$$\lla h\rra\ =\ \big(\lla h\rra_{ij}\big)\ =\ 
\left( \begin{array}{rrrrrl}
0 & 0  &   0    &  0 \   &   0 \   & \cdots \\
  & h_1&    h_2  &    h_3  &   h_4  & \cdots \\
  &   & 2 h_1 &   3 h_2 &  4 h_3 &\cdots \\[0.25em]
  &   &    &  3 h_1 &  6 h_2 & \cdots \\[0.25em]
  &   &    &      &  4 h_1 & \cdots \\
  &   &    &      &      & \ddots 
\end{array}\right)\in\tr_K 
$$
where $\lla h\rra_{ij}=\textstyle{j\choose j-i+1}\, h_{j-i+1}$ for $i\leq j$.
\end{definition}

\medskip\noindent
Thus $h\mapsto \lla h \rra \colon zK[[z]]\to \tr_K$ is an injective continuous $K$-linear map, and 
$$h\in z^{n+1} K[[z]]\Longleftrightarrow\lla h\rra \in \tr_K^n.$$
We introduce the $K$-submodule
$$\mathfrak i\ :=\  \big\{ \lla h\rra:\  h\in z^{2}K[[z]]\big\}$$
of $\mathfrak u=\tr_K^1$. For each $n$, the matrix 
$$\blla \textstyle\frac{z^{n+1}}{(n+1)!} \brra\ =\ \diag_n {i+n\choose n+1}\in\tr_K^n$$
is $n$-diagonal. In the notation from
Section~\ref{dercolfinma}, we have $\blla \textstyle\frac{z^{n+1}}{(n+1)!} \brra =B(n)$.
Note that for $h\in zK[[z]]$ we have 
$\lla h \rra\ =\ \sum_{n= 1}^{\infty} h_{n} \lla z^{n}/n!\rra$. The map 
$$h\ \mapsto\ \lla h\rra\quad \colon\quad  z^2K[[z]] \to \mathfrak i$$ 
is a $K$-module isomorphism. 
Lemma~\ref{lem:binomial diagonal matrices} implies that $\mathfrak i$ is a Lie subalgebra of the Lie algebra $\mathfrak u$ over $K$. We equip $\mathfrak i$ with the filtration $(\mathfrak i^n)$ induced by the filtration of~$\mathfrak u$: $\mathfrak i^n=\mathfrak i\cap\mathfrak u^n=\mathfrak i\cap\tr_K^{n+1}$. Then the isomorphism above and its inverse
are both of rank $0$ with respect to the filtration $(z^{n+2} K[[z]])$ of the $K$-module $z^2K[[z]]$ and the filtration~$(\mathfrak i^n)$ of $\mathfrak i$. Here is the main result of this subsection:

\nomenclature[Y]{$\mathfrak i$}{Lie algebra of iteration matrices}

\begin{theorem}\label{thm:iteration matrices}
$\exp(\mathfrak i) = \mathcal I$; in other words, $\mathfrak i$ is the Lie algebra of $\mathcal I$.
\end{theorem}

\noindent
We give the proof of this theorem after some preparation. Let $t$ be an
indeterminate and set $K^*=K[t]$. We extend the $K$-derivation $\frac{d}{dt}$
of $K^*$ to the $K$-derivation, also denoted by 
$\frac{d}{dt}$,
 of the power series ring $K^*[[z]]$ by 
$$\frac{d}{dt}\left(\sum f_nz^n\right)\ =\ \sum\frac{d f_n}{dt}z^n  \qquad\text{(all $f_n\in K^*$).}$$
Of course, we also have the usual $K^*$-derivation $\frac{\partial}{\partial z}$ on $K^*[[z]]$. 

\begin{lemma}\label{lem:matrix differential}
Let $f\in zK^*[[z]]$ and $h\in zK[[z]]$. Then 
$$\frac{d f}{dt}\ =\ \frac{\partial f}{\partial z}\, h 
\text{ in }K^*[[z]]\ \Longrightarrow\ 
\frac{d}{dt} \llb f\rrb\ =\ \llb f\rrb\,\lla h\rra \text{ in }\tr_{K^*}.$$
\end{lemma}
\begin{proof} Assume 
$\frac{d f}{dt}\ =\ \frac{\partial f}{\partial z}\, h$;
we need to show that for all $i,j$, 
$$\frac{d}{dt}\llb f\rrb_{ij}\ =\ \big(\llb f\rrb\,\lla h\rra\big)_{ij}.$$
For $i=0$, both sides are $0$, so let $i\ge 1$. By the formula for
$f^i$,
$$\frac{1}{i!}\frac{\partial f^i}{\partial z}\ =\ \sum_{k= 1}^{\infty}\, \llb f\rrb_{ik}\frac{z^{k-1}}{(k-1)!}$$
and hence with $h=\sum_{n=1}^{\infty} h_nz^n/n!$ (all $h_n\in K$),
\begin{align*}
\frac{1}{i!}\frac{\partial f^i}{\partial z}\,h\ &=\
\sum_{j=1}^{\infty} \left(\sum_{k=1}^j\, \frac{\llb f \rrb_{ik}\cdot h_{j-k+1}}{(k-1)!(j-k+1)!}\right) z^j\\
&=\  
\sum_{j=1}^{\infty} \left(\sum_{k=1}^j\, \llb f \rrb_{ik} \lla h\rra_{kj}\right) \frac{z^j}{j!}\
=\ \sum_{j=1}^{\infty} \big(\llb f\rrb\,\lla h\rra\big)_{ij} \frac{z^j}{j!}.
\end{align*}
Moreover
$$\frac{1}{i!}\frac{d f^i}{d t}\ =\ \sum_{j= 1}^{\infty} \frac{d}{dt} \llb f\rrb_{ij}\frac{z^j}{j!}.$$
By the hypothesis of the lemma
$$
\frac{d f^i}{d t}\ =\ if^{i-1}\frac{d f}{d t}\ =\ if^{i-1} \frac{\partial f}{\partial z}\, h\ =\
\frac{\partial f^i}{\partial z}\,h,$$
hence
$\frac{d}{dt}\llb f\rrb_{ij} =\big(\llb f\rrb\,\lla h\rra\big)_{ij}$
for all $j$, as required. 
\end{proof}

\noindent
We can now prove the following important fact:

\begin{prop}\label{prop:iteration matrices}
Let $n\ge 1$, $h\in z^{n+1}K[[z]]$, so 
$$t\lla h\rra \in \tr_{K[t]}^n, \qquad \ex^{t\lla h\rra}\ =\ 1+t\lla h\rra + \cdots\in 1+\tr_{K[t]}^n.$$ Set
$f^{[t]}\ :=\ \sum_{j=1}^{\infty}\, \big(\ex^{t\lla h\rra}\big)_{1j}\, \frac{z^j}{j!} \in K^*[[z]]$. Then
$$f^{[t]}\in z+z^{n+1}tK^*[[z]], \qquad \frac{d f^{[t]}}{d t}\ =\ \frac{\partial f^{[t]}}{\partial z}\,h, \qquad
\big\llb f^{[t]}\big\rrb\ =\ \ex^{t\lla h\rra}.$$
\end{prop}
\begin{proof} That $f^{[t]}\in z+z^{n+1}tK^*[[z]]$ is an easy verification, and gives $\llb f^{[t]}\rrb(0)=1$. Let $h=\sum_{k=1}^{\infty}h_kz^k/k!$ with all $h_k\in K$.
Using Lemma~\ref{lem:derivative of exp(tM)}, we get
\begin{align*}
\frac{df^{[t]}}{d t}\ 	&=\ 
\sum_{j=1}^{\infty} \, \left(\frac{d}{dt}\ex^{t\lla h\rra}\right)_{1j}\cdot\frac{z^j}{j!} \\
&=\ \sum_{j= 1}^{\infty} \, \big(\ex^{t\lla h\rra}\lla h\rra\big)_{1j}\cdot\frac{z^j}{j!} \\
&=\ \sum_{j= 1}^{\infty} \left( \sum_{i=1}^j (\ex^{t\lla h\rra})_{1i}\lla h\rra_{ij}\frac{1}{j!}\right) z^j  \\
&=\ \sum_{j=1}^{\infty} \left( \sum_{i=1}^j \frac{(\ex^{t\lla h\rra})_{1i}}{(i-1)!} \frac{h_{j-i+1}}{(j-i+1)!} \right) z^j\ =\ \frac{\partial f^{[t]}}{\partial z}\,h.
\end{align*}
So $\frac{d}{dt} \llb f^{[t]}\rrb = \llb f^{[t]}\rrb\lla h\rra$ by Lemma~\ref{lem:matrix differential}. Thus both $\llb f^{[t]}\rrb$ and $\ex^{t\lla h\rra}$ satisfy $\frac{dY}{dt}=Y \,\lla h\rra$ and $Y(0) =1$.
Hence $\llb f^{[t]}\rrb = \ex^{t\lla h\rra}$ by Lemma~\ref{lem:linear DE}.
\end{proof}

\noindent
For later use we give another description of the power series $f^{[t]}$ from the previous proposition. 
Let $h\in z^{2}K[[z]]$. Then we have the $K^*$-derivations 
$\Delta\ :=\ h\frac{\partial}{\partial z}$ and $t\Delta$ on $K^*[[z]]$.
These are both of rank $1$, and so give rise to the $K^*$-algebra automorphisms $\ex^{\Delta}$ and $\ex^{t\Delta}$ of $K^*[[z]]$, by Lemma~\ref{lem:exp Delta}. With these notations:

\begin{lemma}\label{lem:iteration matrices}
$\quad  \ex^{t\Delta}(z)\in z+tz^2K^*[[z]], \qquad \llb \ex^{t\Delta}(z)\rrb\ =\ \ex^{t\lla h\rra}$.
\end{lemma}
\begin{proof}
Set $f\ :=\ \ex^{t\Delta}(z)\ =\ \sum_{n= 0}^{\infty} \frac{t^n}{n!} \Delta^n(z)$. Then
$$\frac{df}{dt}\ =\ \sum_{n=1}^{\infty} \frac{t^{n-1}}{(n-1)!} \Delta^n(z)\ 
								=\  \Delta\left(\sum_{n=1}^{\infty} \frac{t^{n-1}}{(n-1)!} \Delta^{n-1}(z)\right)\ 
	=\ \Delta(f)\ =\  h\frac{\partial f}{\partial z}.$$
As in the proof of Proposition~\ref{prop:iteration matrices} it now follows that $\llb f\rrb = \ex^{t\lla h\rra}$. 
\end{proof}

\noindent
The following corollary, with $h\in z^2K[[z]]$ and $\Delta$ as before, is obtained by set\-ting~${t=1}$ in the last identity
of Proposition~\ref{prop:iteration matrices} and in Lemma~\ref{lem:iteration matrices} above. It shows in particular that $\exp(\mathfrak i)\subseteq\mathcal I$:

\begin{cor}\label{cor:iteration matrices}
Set $f:= \sum_{j=1}^{\infty}\, (\ex^{\lla h\rra})_{1j}\, \frac{z^j}{j!}\in z+z^2K[[z]]$. Then
$$\ex^{\lla h\rra}\ =\ \llb f \rrb, \qquad \ex^{\Delta}(z)\ =\  f.$$
\end{cor}

\noindent
Recall that $B(k)=\blla \frac{z^{k+1}}{(k+1)!}\brra$ for $k\in\N$.
Given $k_1,\dots,k_n\in\N$ and $k=k_1+\cdots+k_n$, we have
$$B(k_1) \cdots B(k_n)\ =\ \diag_{k} \left(\textstyle {i+k_1\choose k_1+1}{i+k_1+k_2\choose k_2+1}\cdots {i+k_1+\cdots+k_n\choose k_n+1}
\right)_{i\geq 0}$$
by Lemma~\ref{lem:formulas for diagonals}.
Now let $h=\sum_{k=0}^{\infty} h_{k+1}\frac{z^{k+1}}{(k+1)!}\in zK[[z]]$ (all $h_{k+1}\in K$). Then
$$M\ :=\ \lla h\rra\ =\ h_1 \lla z\rra + h_2 \blla  \textstyle\frac{z^2}{2} \rra + \cdots + h_{k+1} \blla  \textstyle\frac{z^{k+1}}{(k+1)!} \rra + \cdots\in \tr_K,$$
and hence
$$M^n\ =\ \sum_{k_1,\dots,k_n\in\N} h_{k_1+1}\cdots h_{k_n+1}\, \blla \textstyle\frac{z^{k_1+1}}{(k_1+1)!}\brra \cdots \blla \textstyle\frac{z^{k_n+1}}{(k_n+1)!}\brra,$$
and so for  $i,j\in\N$ and $n\ge 1$:
\begin{equation}\label{eq:Mn}
(M^n)_{ij}=\!\!\!\sum_{\substack{k_1,\dots,k_n\in\N \\ k_1+\cdots+k_n=j-i}} \!\!\!h_{k_1+1}\cdots h_{k_n+1}\, \textstyle{i+k_1\choose k_1+1}{i+k_1+k_2\choose k_2+1}\cdots {i+k_1+\cdots+k_n\choose k_n+1}.
\end{equation}
This observation leads to:

\begin{lemma}\label{lem:3.10 corrected}
Let $n\ge 1$. Then
$$(M^n)_{11}\ =\ h_1^n,\qquad (M^n)_{1j}\ =\ \frac{j^n-1}{j-1}\,h_1^{n-1}h_{j} + P_{nj}(h_1,\dots,h_{j-1})\quad\text{for $j\geq 2$,}$$
where $P_{nj}(Y_0,\dots,Y_{j-2})\in\Q[Y_0,\dots,Y_{j-2}]$ has all its coefficients in $\Z$, is homogeneous of degree $n$, isobaric of weight $j-1$, and independent of $h$. 
%\textup{(}Here each $Y_i$ is assigned weight $i$.\textup{)}
\end{lemma}
\begin{proof}
Set $i=1$ in \eqref{eq:Mn}. Then the only terms involving $h_{j}$ in this sum are those of the form
$h_1^{n-1}h_{j} \,j^{n-m}$ where $m\in\{1,\dots,n\}$.
This yields the lemma.
\end{proof}

\noindent
For example, $P_{24}(Y_0, Y_1, Y_2)=10Y_1Y_2$ and $P_{34}(Y_0, Y_1, Y_2)= 76Y_0Y_1Y_2 + 18 Y_1^3$.
With $h$ and $M=\lla h \rra$ as above, we get:

\begin{cor}\label{pjexp}
Suppose $h\in z^2K[[z]]$. Then $M\in \tr_K^1$, and for $j\ge 2$,
$$ (\ex^{M})_{1j}\ =\ h_{j} + P_j(h_2,\dots,h_{j-1})$$
where $P_{j}(Y_1,\dots,Y_{j-2})\in\Q[Y_1,\dots,Y_{j-2}]$ is  independent of $h$. Moreover, $P_2=0$, and for $j>2$, $P_j$ has degree at most 
$j-1$ and is isobaric of weight $j-1$.
\end{cor}
\begin{proof}
Let $j\ge 2$. We have $h_1=0$, so by Lemma~\ref{lem:3.10 corrected},
$$(M^n)_{1j}\ =\ 
\begin{cases} 
h_{j}		&\text{if $n=1$,} \\
P_{nj}(h_1,\dots,h_{j-1})	&\text{if $1<n<j$,} \\
0								&\text{if $n\ge j$.}
\end{cases}$$
Hence
$$(\ex^M)_{1j}\ =\ \sum_{n=1}^{j-1} \frac{1}{n!} (M^n)_{1j}\ =\ h_{j} + \sum_{n=2}^{j-1} \frac{1}{n!}P_{nj}(h_1,\dots,h_{j-1}).$$
Thus in view of $h_1=0$, $$P_j(Y_1,\dots,Y_{j-2})\ :=\ \sum_{n=2}^{j-1} \frac{1}{n!}P_{nj}(0,Y_1,\dots,Y_{j-2})$$ has the right properties.
\end{proof}

\noindent
Theorem~\ref{thm:iteration matrices} now follows immediately from Corollary~\ref{cor:iteration matrices} and the following:

\begin{prop}\label{prop:logit}
Let $f\in zK[[z]]$ be unitary. Then $\log\, \llb f\rrb\in\mathfrak i$. 
\end{prop}
\begin{proof}
We have $f=\sum_{j\geq 1} f_j\frac{z^j}{j!}$ (all $f_j\in K$).
We define recursively the se\-quence~$(h_j)_{j\geq 1}$ in $K$ by $h_1:=0$, and $h_{j+1}:=f_{j+1}-P_{j+1}(h_2,\dots,h_{j})$ for $j\ge 1$. Then $h:=\sum_{j=1}^{\infty} h_j\frac{z^{j}}{j!}\in z^{2}K[[z]]$, and by 
Corollary~\ref{pjexp} 
we have $\big(\!\ex^{\lla h\rra}\big)_{1j}=f_j$ for $j\ge 1$.
Corollary~\ref{cor:iteration matrices} now yields $\ex^{ \lla h\rra} = \llb f\rrb$, hence $\log\,\llb f\rrb=\lla h\rra\in \mathfrak i$.
\end{proof}

\noindent
Note: if $f\in z+z^{n+1}K[[z]]$, $n\ge 1$, then $
\llb f\rrb\in 1+\tr_K^n$, and so $\log\, \llb f\rrb\in\mathfrak i^{n-1}$.

\subsection*{The iterative logarithm}
In this subsection we fix a unitary power series
$$f\ =\ z+\sum_{n=2}^{\infty} f_n \frac{z^n}{n!}\in z+z^2K[[z]] \qquad \text{(all $f_n\in K$).}$$
By Theorem~\ref{thm:iteration matrices} there is a (unique) power series $h\in z^2K[[z]]$ with $\log\,\llb f\rrb=\lla h\rra$; we call this
$h$ the {\bf iterative logarithm of $f$} and  denote it by $\operatorname{itlog}(f)$.
The proof of Proposition~\ref{prop:logit} gives
$\operatorname{itlog}(f)=\sum_{n= 2}^{\infty} h_n \frac{z^n}{n!}$ where
$h_n = H_n(f_2,\dots,f_{n})$
and $H_n\in\Q[Y_1,\dots,Y_{n-1}]$ is isobaric of weight $n-1$, for all $n\ge 2$: the ``isobaric'' statement follows inductively from the recursion 
$$H_{n+1}(Y_1,\dots, Y_n)\ =\ Y_n-P_{n+1}\big(H_2(Y_1), \dots, H_n(Y_1,\dots, Y_{n-1})\big)  \quad(n\ge 2).$$
The $H_n$ for $n=2,3,4$ are easily determined: 
\begin{align} H_2\ &=\ Y_1,\quad H_3\ =\ Y_2-\textstyle\frac{3}{2}Y_1^2,\quad 
H_4\ =\ Y_3-5Y_1Y_2 + \textstyle\frac{9}{2}Y_1^3,\quad \text{ and so} \notag   \\
\label{eq:h2, h3}  h_2\ &=\ f_2,\quad  h_3\ =\ f_3-\textstyle\frac{3}{2}f_2^2, \quad h_4\ =\ f_4-5f_2f_3 + \textstyle\frac{9}{2}f_2^2,\quad \text{ and thus}  
\end{align}
$$\lla h \rra\ =\ \log\, \llb f\rrb\ =\  
\begin{pmatrix}
0 	& 0	& 0 		& 0      					& 0 									& \cdots	 	\\                                                                                                                                                  
  	& 0 	& f_2 	& f_3 - \frac{3}{2}f_2^2 	&    f_4 - 5f_2f_3 + \frac{9}{2}f_2^3	& \cdots 	\\
	& 	&  0   	& 3f_2 						& 4f_3 - 6f_2^2						& \cdots		\\	
	&	&		& 0							& 6f_2   							& \cdots 	\\  
	&	&		&							&									& \ddots                                        
\end{pmatrix}.$$
%in particular,
%\begin{align}
%h_2 &= f_2, 								& H_2 &=% Y_1 \notag \\
%h_3 &= f_3-\textstyle\frac{3}{2}f_2^2, 	& H_3 &= Y_2-\textstyle\frac{3}{%2}Y_1^2.\label{eq:h2, h3}
%\end{align}
Here is another way to obtain $\operatorname{itlog}(f)$:

\index{logarithm!iterative}
\index{iterative logarithm}
\nomenclature[Y]{$\operatorname{itlog}(f)$}{iterative logarithm of the power series $f$}

\begin{lemma}\label{lem:itlog alternative} Define $f[n]\in z^{n+1}K[[z]]$
recursively by
$$  f[0]\ =\ z, \qquad f[n+1]\ =\ f[n]\circ f-f[n].$$
Then
$$\operatorname{itlog}(f)\ =\ \sum_{n= 1}^{\infty} \frac{(-1)^{n-1}}{n}f[n].$$
\end{lemma}
\begin{proof}
Set $h:=\operatorname{itlog}(f)=\sum_{j=1}^{\infty} h_jz^j/j!$ (all $h_j\in K$), and let $j\geq 1$. Then
$$h_j\ =\ \lla h\rra_{1j}\ =\ \sum_{n\geq 1} \frac{(-1)^{n-1}}{n}\left(\big(\llb f\rrb-1\big)^n\right)_{1j}.$$
By Lemma~\ref{lem:iteration matrix},
$$\big(\llb f\rrb-1\big)^n\ =\ \sum_{k=0}^n(-1)^{n-k}{n\choose k} \llb f\rrb^k\ =\ \sum_{k=0}^n(-1)^{n-k}{n\choose k}\big\llb f^{[k]}\big\rrb$$
where $f^{[k]}$ denotes the $k$th compositional iterate of $f$, defined recursively by $f^{[0]}=z$ and $f^{[k+1]}=f^{[k]}\circ f$. Hence
$$\left(\big(\llb f\rrb-1\big)^n\right)_{1j}\ =\  \left(\sum_{k=0}^n(-1)^{n-k}{n\choose k}\big\llb f^{[k]}\big\rrb\right)_{1j}\ =\ 
\left\llb\sum_{k=0}^n(-1)^{n-k}{n\choose k}f^{[k]}\right\rrb_{1j}.$$
An easy induction gives $f[n]=\sum_{k=0}^n(-1)^{n-k}{n\choose k}f^{[k]}$. The lemma follows.
\end{proof}

\begin{cor}\label{cor:subst itlog}
Let $\Delta$ be the $K$-derivation $\operatorname{itlog}(f)\textstyle\frac{d}{d z}$ on $K[[z]]$ of rank $1$. Then the $K$-algebra automorphism
$\ex^{\Delta}$ of $K[[z]]$ satisfies
$$g\circ f\ =\ \ex^{\Delta}(g), \text{ for all $g\in K[[z]]$.}$$ 
%In particular, $f[n]=\big(\ex^z-1\big)^n(\Delta)(z)$ for each $n$.
\end{cor}
\begin{proof}
This follows from Corollary~\ref{cor:iteration matrices} and the continuity of $\ex^{\Delta}$.
% The second statement follows from the first and
%Lemma~\ref{lem:linear maps into power series},~(1).
\end{proof}

\noindent
Set $h:=\operatorname{itlog}(f)$, let $A$ be a commutative ring extension of $K$, and $a\in A$. Then
$\ex^{a\lla h\rra}\ =\ 1+a\lla h\rra + \cdots\in 1+a\tr_{A}^1$ and we set
$$f^{[a]}\ :=\ \sum_{j=1}^{\infty}\, \big(\ex^{a\lla h\rra}\big)_{1j}\, \frac{z^j}{j!} \in z+z^2aA[[z]].$$ (For $A=K^*$ and $a=t$ this is just
the $f^{[t]}$ from Proposition~\ref{prop:iteration matrices}.)  
The $K$-algebra morphism $g=g(t) \mapsto g(a)\colon K^*\to A$ sending $t$ to 
$a$ extends to the $K$-algebra morphism 
$M \mapsto M(a)\colon \tr_{K^*} \to \tr_A$ by substituting $a$ for $t$ in each entry, and it extends also to the $K$-algebra morphism
$$ K^*[[z]]\to A[[z]], \qquad g\ =\ \sum g_nz^n\mapsto g\big|_{t=a}\ :=\ \sum g_n(a)z^n  \quad\text{(all $g_n\in K^*$).}$$
It is easy to check that then
$$\ex^{a\lla h\rra}\ =\ \ex^{t\lla h\rra}(a), \qquad f^{[a]}\ =\ f^{[t]}\big|_{t=a}.$$ 
Then Lemma~\ref{lem:iteration matrix} and Proposition~\ref{prop:iteration matrices} give
\begin{equation}\label{eq:EJ}
f^{[0]}\ =\ z, \qquad f^{[1]}\ =\ f, \qquad f^{[a+b]}\ =\ f^{[a]}\circ f^{[b]}\quad(a,b\in A).
\end{equation}
Thus the power series $f^{[a]}$ with $a\in K$ form a subgroup of $z+z^2K[[z]]$ under composition which contains $f$; they may be thought of as fractional iterates of~$f$. 
In fact,~$f^{[t]}$ is \emph{unique}\/ in the sense that
if $g\in K^*[[z]]$ and $g\big|_{t=n}=f^{[n]}$ for all $n$,  
then $g= f^{[t]}$. (Use that $K$ is an integral domain.) 

\nomenclature[Y]{$f^{[a]}$}{fractional iterates of the power series $f$}
 
\begin{cor}\label{cor:itlog}
\begin{equation}\label{eq:itlog diff}
\operatorname{itlog}(f)\ =\ \left.\frac{d f^{[t]}}{d t}\right\lvert_{t=0}
\end{equation}
and  for $a\in K$,
\begin{equation}\label{eq:itlog it}
\operatorname{itlog}(f^{[a]})\ =\ a\cdot \operatorname{itlog}(f).
\end{equation}
If $g\in z+z^2K[[z]]$ is unitary with $f\circ g=g\circ f$, then 
$$\operatorname{itlog}(f\circ g)\ =\ \operatorname{itlog}(f)+\operatorname{itlog}(g).$$
\end{cor}

\begin{proof} Proposition~\ref{prop:iteration matrices} yields \eqref{eq:itlog diff}. Let
$a\in K$. Then $f^{[at]}\big|_{t=n}=f^{[an]}$ for all~$n$, as is easily checked, and $f^{[an]}=(f^{[a]})^{[n]}$ for all $n$ by \eqref{eq:EJ}.
Thus $f^{[at]}=(f^{[a]})^{[t]}$ by the uniqueness property of $(f^{[a]})^{[t]}$. Together with  \eqref{eq:itlog diff} this gives  \eqref{eq:itlog it}. The last statement follows from Lemma~\ref{ab=ba}(ii). 
\end{proof}

\begin{prop}[Acz\'el and Jabotinsky]\label{prop:Jab} In $K^*[[z]]$ we have
\begin{equation}\label{eq:Jab1}
\operatorname{itlog}(f) \cdot \frac{\partial f^{[t]}}{\partial z}\ =\ \frac{df^{[t]}}{d t}\ =\ \operatorname{itlog}(f)\circ f^{[t]}
\end{equation}
and by evaluating at $t=1$ we get in $K[[z]]$,
\begin{equation}\label{eq:Jab2}
\operatorname{itlog}(f)\cdot f'\ =\ \operatorname{itlog}(f)\circ f.
\end{equation}
\end{prop}

\begin{proof}
The first equality in \eqref{eq:Jab1} is from Proposition~\ref{prop:iteration matrices}. 
To get the second one, let $s$ be a new indeterminate, distinct from $t$ and $z$. In $(K[s,t])[[z]]$ we have
$$f^{[s+t]}\ =\ f^{[s]} \circ f^{[t]}.$$
The $K$-derivations $\partial/ \partial s$ and $\partial/ \partial t$ on $K[s,t]$ extend to continuous
$K$-derivations on $K[s,t][[z]]$ with $\partial z/\partial s=\partial z/\partial t = 0$. Then in $K[s,t][[z]]$, 
$$
\frac{\partial f^{[s+t]}}{\partial s}\ =\ \frac{\partial (f^{[s]}\circ f^{[t]})}{\partial s}\ =\ \frac{\partial f^{[s]}}{\partial s}\circ f^{[t]}.
$$
At $s=0$ the left-hand side becomes $\frac{df^{[t]}}{dt}$, and the 
right-hand side $\operatorname{itlog}(f)\circ f^{[t]}$. This gives the
second equality of \eqref{eq:Jab1}.
\end{proof}

\begin{figure}[h]
$$\xymatrix@L=+1em@C+=15em@R+9em{
\big(z+z^2K[[z]],{\,\circ\,},z\big) \ar@{-}[r]_(.55){\begin{minipage}[r]{7em}\small $f\, \longmapsto\, \operatorname{itlog}(f)$\end{minipage}}^(.55){\begin{minipage}{13em}\small $\ \ex^{h\frac{\partial}{\partial z}}(z)\, \longmapsfrom\, h\ \ \ \  $\end{minipage} } 
\ar@{-}[d]_{\begin{minipage}{3em}\small \begin{center}  \ \\[0.5em] $f$\\[1em] \rotatebox[origin=c]{270}{$\longmapsto$}\\[0.5em] $\llb f\rrb$ \end{center} \end{minipage}}
^{\begin{minipage}{4em}\small \begin{center} $\sum\limits_{j=1}^\infty G_{1j} \frac{z^j}{j!}$\\[0.1em] \rotatebox[origin=c]{90}{$\longmapsto$}\\[0.5em] $G$ \end{center} \end{minipage} } 
& \big(z^2K[[z]],[\ \, ,\ ]\big) 
\ar@{-}[d]_{\begin{minipage}{3em}\small \begin{center}  \ \\[0.5em] $h$\\[1em] \rotatebox[origin=c]{270}{$\longmapsto$}\\[0.5em] $\lla h\rra$ \end{center} \end{minipage}}
^{\begin{minipage}{4.5em}\small \begin{center} $\sum\limits_{j=1}^\infty M_{1j} \frac{z^j}{j!}$\\[0.1em] \rotatebox[origin=c]{90}{$\longmapsto$}\\[0.5em] $M$ \end{center} \end{minipage} } 
 \\
(\mathcal I,{\,\cdot\,},1) \ar@{-}[r]_(.55){\begin{minipage}{6em}\small$G\, \longmapsto\, \log G$\end{minipage}}^(.55){\begin{minipage}{10em}\small$\exp M\, \longmapsfrom\, M\ \ \ \,$\end{minipage} }& \big(\mathfrak i,[\ \, ,\ ]\big)
}$$
\caption{Relationship between power series and (infinitesimal) iteration matrices.}\label{fig:itlog}
\end{figure}

\subsection*{Notes and comments}
A general reference for properties of the Bell polynomials (named after E.~T.~Bell~\cite{Bell}) is~\cite{Comtet-Book}.
(Our notation differs slightly from that in~\cite{Comtet-Book}: our~$B_{ij}$ here is $\text{\bf B}_{ji}$ there.)  Lemma~\ref{lem:iteration matrix} is from Jabotinsky~\cite{J2,J1}.
Theorem~\ref{thm:iteration matrices} for $K=\mathbb C$ is in \cite{Schippers}; the proof given here follows \cite{A-it}.
The Lie subalgebra $\bigoplus_{n\geq 1} K\lla z^{n+1}\rra$ of $\mathfrak i_K$, with $\big[\lla z^{m+1}\rra, \lla z^{n+1}\rra\big]=(m-n)\lla z^{m+n+1}\rra$ for all $m,n\geq 1$ (by \eqref{eq:Lie bracket of B(n)}), is a variant of the  ``Witt algebra'' \cite[pp.~206--212]{AmayoStewart}. 
The terminology   ``iterative logarithm'' was coined in \cite{Ecalle-iterative}.
The construction of a family of power series~$(f^{[a]})_{a\in K}$ sa\-tis\-fy\-ing~\eqref{eq:EJ} goes back to \cite{EJ}.
The functional equation~\eqref{eq:Jab2} satisfied by the iterative logarithm is known as Julia's equation in iteration theory. (See \cite[\S{}8.5A]{KCG}.) It was found by
Acz\'el~\cite{Aczel} and Jabotinsky~\cite{J1}, although~\cite{Gronau} suggests that G.~Frege was already aware of it.

\section{Riordan Matrices}\label{sec:riordan}

\noindent
{\em In this section $K$ is an integral domain.}\/ We enlarge the group of iteration matrices to the group of so-called Riordan matrices. 

\subsection*{The Riordan group}
A {\bf Riordan pair} \index{Riordan!pair} over $K$ is a pair $(f,g)$ where $f\in zK[[z]]$  and $g\in K[[z]]$.
Let $(f,g)$ be a Riordan pair,  $f=\sum f_n\frac{z^n}{n!}$ and $g=\sum g_n\frac{z^n}{n!}$ with $f_n,g_n\in K$ for all $n$, and~${f_0=0}$. Then by \eqref{eq:power of y}, with $i,j$ ranging over $\N$,
$$\frac{1}{i!} f^ig\ =\ \sum_{j\geq i} R_{ij}\frac{z^j}{j!}$$
where for $i\le j$,
$$R_{ij}\ =\ \sum_{k=i}^j  {j\choose k} B_{ik}(f_1,\dots,f_{k-i+1})\cdot g_{j-k}.$$
We also set $R_{ij}=0$ for $i>j$.
We call the triangular matrix
$$\llb f,g \rrb\ :=\ (R_{ij})\in \tr_K$$
the {\bf Riordan matrix} of $(f,g)$. Note that
$$R_{0j}\ =\ g_j\text{ for all $j$,}\qquad R_{1j}\ =\ \sum_{k=1}^j \binom{j}{k} f_kg_{j-k}\text{ for $j\geq 1$.}$$
Clearly $\llb f,1\rrb$ is the iteration matrix $\llb f\rrb$ of $f$ as introduced in Section~\ref{sec:Iteration matrices}. We have $R_{ii}=f_1^i g_0$ for each $i$, hence $\llb f,g \rrb\in 1+\tr_K^1$ iff $f_1=g_0=1$. We say that the Riordan pair $(f,g)$ is {\bf unitary} if $f_1=g_0=1$.

\index{unitary!Riordan pair}
\index{Riordan!matrix}
\nomenclature[Y]{$\llb f,g \rrb$}{Riordan matrix of the Riordan pair $(f,g)$}

\medskip
\noindent
In the next lemma we identify a power series $a=\sum a_n \frac{z^n}{n!}\in K[[z]]$, $a_n\in K$ for each~$n$, with the element $(a_n)$ of $K^\N$, viewed as a row vector. With this convention:

\begin{lemma}\label{lem:Riordan}
Let $(f,g)$ be a Riordan pair over $K$ and $a,b\in K[[z]]$. Then
$$a\cdot \llb f,g\rrb\ =\ b \text{ in $K^{\N}$}\quad\Longleftrightarrow\quad (a\circ f) \cdot g\ =\ b \text{ in $K[[z]]$.}$$
\end{lemma}
\begin{proof}
This is immediate from
\equationqed{(a\circ f) \cdot g\ =\ \left(\sum a_i\frac{f^i}{i!}\right)\cdot g\ =\ \sum_{i} a_i\left(\sum_{j\geq i} R_{ij}\frac{z^j}{j!}\right)\ =\ \sum_{j} \left(\sum_{i\leq j} a_iR_{ij}\right)\frac{z^j}{j!}.}
\end{proof}

\noindent
As a consequence the set of Riordan matrices of unitary Riordan pairs over $K$ is a subgroup of the unitriangular group $\mathcal U=1+\tr^1_K$:

\begin{cor}\label{cor:Riordan}
Let $(f,g)$ and $(f^*,g^*)$ be Riordan pairs over $K$. Then $$\big(f\circ f^*,(g\circ f^*)\cdot g^*\big)$$ is a  Riordan pair over $K$, and
$$\llb f,g\rrb \cdot  \llb f^*,g^* \rrb\ =\ \llb f\circ f^*,(g\circ f^*)\cdot g^*\rrb.$$
If $(f,g)$, $(f^*,g^*)$ are unitary, then so is $\big(f\circ f^*,(g\circ f^*)\cdot g^*\big)$.
Moreover, if $(f,g)$ is unitary, then $\big(f^{[-1]},1/(g\circ f^{[-1]})\big)$ is a unitary Riordan pair over $K$, and for $(f^*,g^*)=\big(f^{[-1]},1/(g\circ f^{[-1]})\big)$ we have 
$$\llb f,g \rrb\cdot \llb f^*,g^*\rrb\ =\ \llb f^*,g^*\rrb \cdot \llb f,g\rrb\ =\ \llb z,1\rrb\ =\ 1.$$
\end{cor}
\begin{proof}
The row with index $i$ of $\llb f,g\rrb$ is $\frac{1}{i!}f^i g$, hence by the previous lemma, the row with index $i$ of
$\llb f,g \rrb \cdot  \llb f^*,g^*\rrb$ is $\big((\frac{1}{i!}f^ig)\circ f^*\big)\cdot g^* = \frac{1}{i!}(f\circ f^*)^i\cdot \big((g\circ f^*)\cdot g^*\big)$.
\end{proof}

\noindent
We call the subgroup $\mathcal R$ of $\mathcal U$ consisting of the Riordan matrices of unitary Riordan pairs over $K$ the {\bf Riordan group} over $K$. Note that the unitriangular group $\mathcal R$ over $K$ is algebraic.
The group $\mathcal I$ of iteration matrices is a subgroup of~$\mathcal R$.

\index{group!Riordan}
\index{Riordan!group}
\nomenclature[Y]{$\mathcal R$}{Riordan group}
\nomenclature[Y]{$[g]$}{Appell matrix $[g]=\llb z,g\rrb$ of the power series $g$}

\subsection*{The Appell group and its Lie algebra}
For $g=\sum g_n\frac{z^n}{n!}$ ($g_n\in K$ for all $n$), we set $[g]:=\llb z,g\rrb$, that is,
$$[g]\ =\  
\left( \begin{array}{rrrrrl}
g_0	& g_1	& g_2	   & g_3	    & g_4	  & \cdots \\
    & g_0	&  2 g_1  & 3g_2	    & 4 g_3	  &    \cdots \\
[0.25em]
  	&   		& g_0 	&  3 g_1 	& 6 g_2 	& \cdots \\
[0.25em]
  	&   		&    	&  g_0 		& 4 g_1 	& \cdots \\
  	&   		&    	&      				& g_0 				& \cdots \\
  	&   		&    	&      				&    	  			& \ddots 
\end{array}\right)\quad
\text{ where $ [g]_{ij}=\textstyle{j\choose i}\, g_{j-i}$ for $i\leq j$.}$$
Note that $[g]\cdot \llb f \rrb=\llb f,g\circ f\rrb$ for $f\in zK[[z]]$, $g\in K[[z]]$. By Corollary~\ref{cor:Riordan}, the map $g\mapsto [g]\colon K[[z]]\to\tr_K$ is an embedding of $K$-algebras. Moreover,
$$g\in z^nK[[z]]\ \Longleftrightarrow\ [g]\in\tr_K^n 
\qquad(g\in K[[z]]).$$
Hence
the image of this embedding, with the filtration induced by $(\tr_K^n)$,
is a complete filtered subalgebra of $\tr_K$, and $g\mapsto [g]\colon K[[z]]\to\tr_K$ is continuous. In particular,
for $g\in zK[[z]]$ and $h\in 1+zK[[z]]$ we have
\begin{equation}\label{eq:explog Appell}
[\exp g]\ =\ \exp\,[g], \qquad [\log h]\ =\ \log\,[h].
\end{equation}
 The embedding $g\mapsto [g]\colon K[[z]]\to\tr_K$ maps the subgroup $1+zK[[z]]$ of $K[[z]]^\times$ onto a commutative normal subgroup of $\mathcal R$, called the {\bf Appell group} over $K$ and denoted here by $\mathcal A$.\index{group!Appell}\index{Appell group}\nomenclature[Y]{$\mathcal A$}{Appell group} 
Now $\mathcal A \cap \mathcal I= \big\{1\big\}$, and from
$$ \llb f,g\rrb\ =\ \llb f\rrb\cdot [g]\ \text{ for $\llb f,g\rrb\in\mathcal R$,}$$ 
we get $\mathcal R=\mathcal I\cdot \mathcal A=\mathcal A \cdot \mathcal I$, so the group $\mathcal R$ is the internal semidirect product 
of its normal subgroup $\mathcal A$ with its subgroup $\mathcal I$. 
Now set 
$$\mathfrak a\ :=\ \big\{ [g] : g\in z K[[z]] \big\},$$ 
an abelian Lie subalgebra of the Lie algebra $\mathfrak u=\tr^1_K$ over $K$;
with the filtration induced by $\mathfrak u$, the filtered Lie algebra $\mathfrak a$ over $K$ is complete.
For $n\geq 1$, the matrix
$$\big[\textstyle\frac{z^n}{n!}\big]\ =\ \diag_n \textstyle{i+n\choose n}\ =\ A(n)\in\mathfrak a^{n-1}$$
is $n$-diagonal, and for $g=\sum_{n= 1}^{\infty}g_n\frac{z^n}{n!}\in zK[[z]]$, where $g_n\in K$ for $n\geq 1$, we have
$$[g]\ =\ \sum_{n=1}^{\infty}g_n \big[\textstyle\frac{z^n}{n!}\big]\ =\ \sum_{n=1}^{\infty} g_nA(n)\ \text{ in $\mathfrak a$.}$$
By \eqref{eq:explog Appell} we have $\exp(\mathfrak a)=\mathcal A$ and $\log(\mathcal A)=\mathfrak a$, so $\mathfrak a$ is the Lie algebra of the algebraic unitriangular group $\mathcal A$ over $K$.

\medskip
\noindent
For each unitary Riordan pair $(f,g)$ and $\phi\in K^\times$, the pair
$\big(\phi^{-1}f(\phi z),g(\phi z)\big)$ is also a unitary Riordan pair, and the diagonal matrix $D=\diag(\phi^i)=\llb \phi z,1\rrb$ satisfies
$$D^{-1} \llb f,g \rrb D\ =\ \big\llb \phi^{-1}f(\phi z),g(\phi z)\big\rrb.$$
Hence for all such $D$,
$$D\mathcal R D^{-1}\ =\ \mathcal R,\quad D\mathcal I D^{-1}\ =\ \mathcal I,\quad D\mathcal A D^{-1}\ =\ \mathcal A.$$

\subsection*{Some identities involving Riordan matrices} Recall 
from Section~\ref{dercolfinma}   the shift matrix $S\in K^{\N\times\N}$ 
and the derivative $M'$ of a column-finite
matrix $M\in K^{\N\times \N}$.

\begin{lemma}\label{lem:shift it matrix}
Let $(f,g)$ be a Riordan pair over $K$. Then
$$\llb f,g \rrb \cdot S\ =\ S\cdot \llb f,f'g \rrb + \llb f,g'\rrb.$$
In particular, $\llb f\rrb \cdot S = S\cdot \llb f,f'\rrb$ and $[g]'=[g']$.
\end{lemma}
\begin{proof}
Let $i\in\N$. Then by definition of $\llb f,g\rrb$ we have
$$\frac{1}{i!}f^i g\ =\ \sum_{j\geq i}\, \llb f,g\rrb_{ij} \frac{z^j}{j!}.$$ 
Differentiating both sides with respect to $z$ yields for $i\ge 1$,
$$\frac{1}{(i-1)!}f^{i-1}f'g + \frac{1}{i!}f^ig'\ =\ \sum_{j \geq i-1 }\, \llb f,g\rrb_{i,j+1} \frac{z^{j}}{j!}.$$
The definition of $\llb f,f'g\rrb$ and $\llb f,g'\rrb$ gives for $i\ge 1$,
$$\frac{1}{(i-1)!}f^{i-1}f'g + \frac{1}{i!}f^ig'\  =\ \left(\sum_{j \geq i-1}\, \llb f,f'g\rrb_{i-1,j}  \frac{z^{j}}{j!}\right) + \left(\sum_{j \geq i}\, \llb f,g'\rrb_{i,j}  \frac{z^{j}}{j!}\right).$$
These last two identities together give for $i\ge 1$ and $j\geq i-1$:
$$\big(\llb f,g\rrb \cdot S\big)_{ij}\ =\ \llb f,g\rrb_{i,j+1}\ =\ \llb f,f'g\rrb_{i-1,j}+\llb f,g'\rrb_{ij}\ =\ \big(S\cdot \llb f,f'g\rrb + \llb f,g'\rrb\big)_{ij}.$$
These equalities actually hold for all $i,j\in \N$, as is easily verified using
$B_{0j}=0$ for $j\ge 1$, and 
thus $\llb f,g \rrb \cdot S = S\cdot \llb f,f'g \rrb + \llb f,g'\rrb$.
Taking $g=1$ yields  $\llb f \rrb\cdot S=S \cdot \llb f,f' \rrb $. Taking $f=z$ gives $[g]'=[g']$.
\end{proof} 

\begin{cor}\label{cor:shift it matrix}
Let $f\in z+z^2K[[z]]$ and $g\in K[[z]]$. Then 
\begin{align*} f'&\in 1+zK[[z]]\subseteq K[[z]]^\times,\\
\llb f,g\rrb'\ &=\ \llb f,g \rrb\cdot S\cdot\big[ 1-(1/f') \big] + \llb f,g'/f'\rrb,\quad \text{in particular,}\\
\llb f\rrb'\ &=\  \llb f\rrb\cdot S\cdot\big[ 1-(1/f') \big].
\end{align*}
\end{cor}
\begin{proof}
By the previous lemma, 
$$\llb f,g\rrb \cdot S\ =\ \big( S\cdot \llb f,g\rrb + \llb f,g'/f'\rrb\big)\cdot [f'].$$
Note that $f'\in 1+zK[[z]]$, so $[f']\in 1+\tr^1_K= \mathcal{U}$,  
and hence 
$$S\cdot\llb f,g\rrb\ =\  \llb f,g\rrb \cdot S\cdot [f']^{-1} - \llb f,g'/f'\rrb,$$
and thus
\begin{align*}
\llb f,g\rrb'\ 	&=\ \llb f,g\rrb \cdot S - S\cdot \llb f,g\rrb \\
		&=\ \llb f,g \rrb\cdot S\cdot\big(1-[f']^{-1}\big) +  \llb f,g'/f'\rrb \\
		&=\ \llb f,g \rrb\cdot S\cdot\big[ 1-(1/f') \big] + \llb f,g'/f'\rrb
\end{align*}
as claimed.
\end{proof}

\subsection*{The Lie algebra of $\mathcal R$}
The material in this subsection is not used later.

\begin{lemma}\label{lem:commutators in r_K}
Let $g\in zK[[z]]$ and $h\in z^2K[[z]]$.
Then
$\ \big[ [g], \lla h\rra \big]\ =\  [g'h]$.
\end{lemma}
\begin{proof}
Let $m,n\ge 1$. With the notation from Section~\ref{dercolfinma},  
$$\big[ \textstyle\frac{z^m}{m!}\big]\ =\ A(m), \qquad \blla\textstyle\frac{z^n}{n!}\brra\ =\ B(n-1),$$ 
and so by Lemma~\ref{eq:Lie bracket of A(m), B(n)},
\begin{align*}
\bigg[\big[ \textstyle\frac{z^m}{m!}\big], \blla\textstyle\frac{z^n}{n!}\brra\bigg]\ &=\ \big[A(m),B(n-1)\big]\ =\ 
{m+n-1 \choose n} A(m+n-1)\\
    &=  {m+n-1 \choose n}\big[ \textstyle\frac{z^{m+n-1}}{(m+n-1)!}\big].
\end{align*}
The general case now follows from this special case and the continuity of the Lie bracket operation.
%Write $g=\sum_{m\geq 1} g_m \frac{z^m}{m!}$ and $h=\sum_{n\geq 2} h_n \frac{z^n}{n!}$ with $g_m,h_n\in K$ for $m\geq 1$, $n\geq 2$.  
%\begin{align*}
%\big[ [g], \lla h\rra \big]	&= \sum_{m\geq 1,\ n\geq 2} g_m h_n \textstyle{m+n-1\choose n} \big[\textstyle\frac{z^{m+n-1}}{(m+n-1)!}\big] \\
%							&= \sum_{k\geq 2} \left(\sum_{\substack{m\geq 1,\ n\geq 2 \\ m-1+n=k}} g_m h_n \textstyle{k\choose m-1}\right) \big[\textstyle\frac{z^{k}}{k!}\big] \\
%							&= \sum_{k\geq 2} \left(\sum_{\substack{m\geq 0,\ n\geq 2 \\ m+n=k}} g_{m+1} h_n \textstyle%{k\choose m}\right) \big[\textstyle\frac{z^{k}}{k!}\big] = [g'h].
%\end{align*}
\end{proof}

\noindent
Clearly $\mathfrak a\cap\mathfrak i=\{0\}$, where  $\mathfrak i$ is the Lie algebra of $\mathcal I$. Let 
$\mathfrak r$ be the Lie algebra of the algebraic unitriangular group $\mathcal R$ over $K$. 
From the previous lemma, it follows that $\mathfrak{a} \oplus \mathfrak{i}$ is a Lie subalgebra
of $\mathfrak{r}$ and $\mathfrak{a}$ is an ideal of $\mathfrak{a} \oplus \mathfrak{i}$. 
In fact, $\mathfrak r=\mathfrak{a} \oplus \mathfrak{i}$, but we will not 
prove this here. 

\subsection*{Notes and comments}
The Riordan group (introduced in \cite{Shapiro} and named in honor of J.~Riordan for his work~\cite{Riordan} on combinatorial identities) is connected to Rota's ``umbral calculus'' dealing with sequences of polynomials. In this calculus,
an Appell sequence (over $K$) is a sequence $(P_j)_{j\geq 0}$ of  polynomials $P_j\in K[z]$ with~$P_0$ of degree $0$ and $P_j' = jP_{j-1}$ for $j>0$ (see \cite[Theorem~2.5.6]{Roman}).
The subgroup~$\mathcal A$ of the Riordan group is called the Appell group because 
for $R=(R_{ij})\in\mathcal A$, the sequence $\big(\sum_{i=0}^j R_{ij}  z^{i}\big)_{j\geq 0}$ is an Appell sequence. The subgroup $\mathcal I$ of $\mathcal R$ is also called the ``associated subgroup'' in the combinatorics literature.

\section{Derivations on Polynomial Rings}\label{sec:Derivations on polynomials rings}

\nomenclature[Y]{$A^\partial$}{algebra of constants of the derivation $\partial$ on $A$}
\nomenclature[Y]{$A^{\Delta,\Lambda}$}{algebra of common constants of $\Delta$, $\Lambda$}

\noindent
{\em In this section $A$ is a commutative ring containing $K$ \textup{(}and thus $\Q$\textup{)} 
as a subring.}\/ Then a $K$-derivation on~$A$ as defined in Section~\ref{sec:Filtered and Graded Algebras} is the same as a derivation on~$A$ whose ring of constants contains $K$. Recall that $\Der_K(A)$ is the Lie algebra over $K$ consisting of the $K$-derivations on $A$, with Lie bracket
$$[\Delta,\Lambda]\ =\ \Delta\Lambda - \Lambda\Delta \qquad (\Delta,\Lambda\in \Der_K(A)).$$
Since $A$ is commutative, $a\Delta\in\Der_K(A)$ for all $a\in A$ and $\Delta\in\Der_K(A)$, and so
$\Der_K(A)$ is naturally a left $A$-module. For $\partial\in\Der_K(A)$, let $A^\partial$ be the ring of constants of $\partial$, so $A^\partial$ is a $K$-subalgebra of $A$. Let $\Delta, \Lambda\in \Der_K(A)$. Then we set $A^{\Delta,\Lambda}:={A^{\Delta}\cap A^{\Lambda}}$ and,
with $\operatorname{Lie}(\Delta,\Lambda)$ the Lie subalgebra of $\Der_K(A)$
over $K$ generated by~$\Delta$,~$\Lambda$, we have
\begin{equation}\label{eq:ring of constants}
A^{\Delta,\Lambda}\ =\ \bigcap_{\partial\in \operatorname{Lie}(\Delta,\Lambda)} A^{\partial}.
\end{equation}

\subsection*{Exponential automorphisms}
If $\Delta\in\Der_K(A)$ is locally nilpotent, then
the (locally unipotent) automorphism $\exp\Delta$ of the $K$-module $A$, with inverse $\exp(-\Delta)$, is a $K$-algebra automorphism (Lem\-ma~\ref{lem:exp Delta}).
The $K$-automorphisms of $A$ of the form $\exp\Delta$ with locally nilpotent $\Delta\in\Der_K(A)$ are said to be {\bf exponential.} The above can be reversed: 

\index{automorphism!exponential}
\index{exponential!automorphism}

\begin{prop}
\label{prop:vdEssen}
Let $\sigma$ be a locally unipotent $K$-algebra endomorphism. 
Then $\log\sigma$ is a locally nilpotent $K$-derivation on $A$, and so 
$\sigma=\exp(\log\sigma)$ is an exponential automorphism. 
\end{prop}

\begin{proof} Set $\Delta:= \log \sigma$, so $\Delta$ is a locally nilpotent
endomorphism of the $K$-module~$A$, and $\sigma=\exp \Delta$. 
Let $a,b\in A$ and take
$N\in \N^{\ge 1}$ with $\Delta^i(a)=\Delta^i(b)=\Delta^i(ab)=0$ for all $i>N/2$. Then $\sigma(ab)=\sigma(a)\sigma(b)$ gives
$$  \sum_{k=0}^N \frac{1}{k!}\Delta^k(ab)\ =\ 
\sum_{k=0}^N \sum_{i+j=k} \frac{1}{i!}\frac{1}{j!} \Delta^i(a)\Delta^j(b).$$
Replacing $\sigma$ by $\sigma^n$ changes $\Delta$ to $n\Delta$, and
the above equality remains valid for $\Delta$ replaced by $n\Delta$, that is, 
$\sum_{k=0}^Nc_k n^k=0$ with 
$$c_k\ :=\ \frac{1}{k!}\Delta^k(ab)-\sum_{i+j=k} \frac{1}{i!}\frac{1}{j!}\Delta^i(a)\Delta^j(b)\in A.$$ 
The $(N+1)\times (N+1)$
Vandermonde matrix with rows $(n^0, n^1,\dots, n^N)$ for $n=0,\dots, N$
is invertible, but annihilates the column vector 
$(c_0,\dots, c_N)^{\operatorname{t}}$, so $c_k=0$ for $k=0,\dots,N$. For $k=1$ this
yields $\Delta(ab)=\Delta(a)b + a\Delta(b)$.
\end{proof}

\noindent
If $\Delta\in\Der_K(A)$ is locally nilpotent and $a\in A^\Delta$ then $a\Delta$ is locally nilpotent.
Using the Leibniz rule one also shows:

\begin{lemma}\label{lem:locally nilpotent on generators}
If $S\subseteq A$ generates the $K$-algebra $A$, and $\Delta\in\Der_K(A)$ is 
such that for every $s\in S$ there is $n$ with $\Delta^n(s)=0$, 
then $\Delta$ is locally nilpotent.
\end{lemma}

\noindent
If $\Delta,\Lambda\in\Der_K(A)$ are locally nilpotent and $[\Delta,\Lambda]=0$, then 
$$\exp(\Delta)\exp(\Lambda)\ =\ \exp(\Delta+\Lambda)$$ is an exponential automorphism of $A$. In particular, for each locally nilpotent $\Delta\in\Der_K(A)$ and $k\in\Z$, $\exp(\Delta)^k$ is an exponential automorphism of $A$ with $\exp(\Delta)^k=\exp(k\Delta)$.

\begin{remark}
The composition of two exponential automorphisms is not in general exponential. For example, take $A=K[Y_0,Y_1]$. Then the $K$-derivations $\Delta=Y_0\frac{\partial}{\partial Y_1}$, $\Lambda=Y_1\frac{\partial}{\partial Y_0}$ on $A$ are locally nilpotent. A computation shows that 
$$\Phi\ :=\ \exp(\Delta)\exp(\Lambda)-1$$ 
satisfies $\Phi(Y_0)=Y_0+Y_1$, $\Phi(Y_1)=Y_0$ and hence $\Phi^n(Y_0)\ne 0$
for all $n$. So $\Phi$ is not locally nilpotent, and the $K$-algebra automorphism $\exp(\Delta)\exp(\Lambda)$ of $A$ is not locally unipotent, and thus not
 exponential.
\end{remark}

\noindent
In the next subsection we study a class of exponential automorphisms of polynomial algebras over $K$ which do form a group under composition.

\subsection*{Triangular derivations}
{\em In the rest of this section $A=K[Y_0, Y_1, Y_2,\dots]$ where $(Y_n)$ is a sequence of distinct indeterminates.}\/ For each $n$ we also put  
\begin{align*} Y_{<n}\ &=\ (Y_0,\dots,Y_{n-1}), \text{ so}\\ 
      K[Y_{<m}]\ \subseteq\ K[Y_{<n}]\ &\subseteq\ K[Y]\ =\ A\qquad\text{for $m\leq n$, and}\\ 
A\ &=\ \bigcup_{n} K[Y_{<n}].
\end{align*}
We abbreviate $\Der_K(A)$ by $\Der_K$, and we denote by $\partial_m$ the $K$-derivation $\frac{\partial}{\partial Y_m}$ of $A$ as well as any restriction of this derivation to a $K$-subalgebra $K[Y_{<n}]$ where $m < n$. The following well-known fact has a routine proof:

\begin{lemma}\label{lem:free basis of derivation module}
$\Der_K\!\big(K[Y_{<n}]\big)$ is a free $K[Y_{<n}]$-module with basis $\partial_0,\dots,\partial_{n-1}$, and $[\partial_i,\partial_j]=0$ for $0\leq i,j< n$. For each $\Delta\in\Der_K\!\big(K[Y_{<n}]\big)$ we have
$$\Delta\ =\ \Delta(Y_0)\partial_0+\cdots+\Delta(Y_{n-1})\partial_{n-1}.$$
\end{lemma}

\noindent
In a similar vein: for each sequence $(f_n)$ in $A$ the sequence of 
en\-do\-mor\-phisms $\big(f_n \frac{\partial}{\partial Y_n}\big)$ of the $K$-module $A$ is
summable (with respect to the trivial filtration on the $K$-mo\-dule~$A$ as in Section~\ref{triangularlinear}), and $\Delta:= \sum_{n=0}^\infty f_n \frac{\partial}{\partial Y_n}$ is
the unique $K$-derivation of $A$ with $\Delta(Y_n)=f_n$ for all $n$. 

\medskip\noindent
If $\Delta\in\Der_K$ and $\Delta(Y_{n-1})\in K[Y_{<n}]$ for all 
$n\ge 1$, then  $\Delta\big(K[Y_{<n}]\big)\subseteq K[Y_{<n}]$ for all $n$. 
%in which case the previous lemma suggests that we write
%$$\Delta = \Delta(Y_0)\partial_0 + \cdots + \Delta(Y_n)\partial_n + \cdots.$$
Note also that each $K$-derivation $\partial_n\in\Der_K$ is locally nilpotent. 
More generally, every $\Delta\in\Der_K$ with $\Delta(Y_n)\in K[Y_{<n}]$ for 
all $n$ is locally nilpotent, by the following:

\begin{lemma}\label{lem:generalized triangular implies locally nilpotent}
Let $\Delta$ be a $K$-derivation of $K[Y_{<n}]$ with $\Delta(Y_j)\in K[Y_{<j}]$ for $j=0,\dots,n-1$. Then $\Delta$ is locally nilpotent.
\end{lemma}
\begin{proof}
We proceed by induction on $n$. The case $n=0$ is clear, so let $n>0$. 
Then the derivation $\Delta$ restricts to a $K$-derivation of $K[Y_{<n-1}]$, 
and so inductively we can assume this restriction to be locally nilpotent. 
Now $\Delta(Y_{n-1})\in K[Y_{<n-1}]$ and hence 
$\Delta^{m+1}(Y_{n-1})=\Delta^m\big(\Delta(Y_{n-1})\big)=0$ for some $m$. 
Thus $\Delta$ is locally nilpotent, 
by Lemma~\ref{lem:locally nilpotent on generators}.
\end{proof}

\begin{definition}\label{def:triangular derivation}
Let $\Delta$ be a $K$-derivation of $A$. Then we call $\Delta$ {\bf triangular} (relative to $Y_0, Y_1,\dots$) if $\Delta$ restricts to a triangular 
endomorphism of the $K$-submodule $A_1=\bigoplus_{j} K\,Y_j$ of $A$ with respect to its basis $Y_0, Y_1,\dots$; that is, 
$$\Delta(Y_j)\ =\ \Delta_{0j} Y_{0} + \Delta_{1j} Y_{1} + \cdots + \Delta_{jj} Y_j$$
where $\Delta_{ij}\in K$ for $i\leq j$. If also $\Delta_{jj}=0$ for all $j$, then $\Delta$ is {\bf strictly triangular}.
Every triangular endomorphism of the $K$-submodule $A_1$ with 
respect to the basis  $Y_0, Y_1,\dots$ extends uniquely to a 
(necessarily triangular) $K$-derivation of $A$.  
\end{definition}

\index{derivation!triangular}
\index{derivation!strictly triangular}
\index{triangular!derivation}
\index{strictly!triangular}
\nomenclature[Y]{$\TrDer_K=\TrDer_K(A)$}{Lie algebra of triangular $K$-derivations of~$A$}

\noindent
Let $\TrDer_K=\TrDer_K(A)$ be the set of triangular $K$-derivations of $A$.
The Lie brack\-et~$[\Delta,\Lambda]$ of triangular $K$-derivations $\Delta$, $\Lambda$ of $A$ is also triangular; hence $\TrDer_K$ is a 
Lie sub\-al\-ge\-bra (over $K$) of $\Der_K$.
For $\Delta\in\TrDer_K$ with the $\Delta_{ij}$ as above, put 
$$M_\Delta\ :=\ M_{\Delta|A_1}\ =\ (\Delta_{ij})_{i,j\in\N}\in\tr_K,$$
with $\Delta_{ij}:=0$ for $i>j$, by convention.
We obtain a complete filtration $(\TrDer^n_K)_{n\geq 0}$ of the Lie algebra $\TrDer_K$ over $K$ by setting
$$\TrDer^n_K\ :=\ \{ \Delta\in \TrDer_K :\  M_\Delta\in\tr_K^n\}.$$
(So  $\TrDer^1_K$ consists of all strictly triangular $K$-derivations.) 
We have a commuting diagram of isomorphisms of Lie algebras over $K$:
$$\xymatrixcolsep{5pc}
\xymatrix{
\tr_{K}   \\
\TrDer_{K} \ar[u]^{\Delta\mapsto M_\Delta}\ar[r]^{\Delta\mapsto \Delta|A_1} & \TrEnd_{K} \ar[ul]_{\Phi\mapsto M_\Phi}
}$$
Every triangular $K$-derivation of $A=K[Y]$ restricts to a derivation of $K[Y_{<n}]$, for each $n$. Hence by Lemma~\ref{lem:generalized triangular implies locally nilpotent}, every strictly triangular $K$-derivation of $A$ is locally nilpotent. 
Note that if $\Delta\in\TrDer^n_K$ and $P\in K[Y_{<n}]$, then $\Delta(P)=0$.

\subsection*{Diagonals}
In the next lemma $\degree$ is a degree function on $A$ and $\Delta\in\Der_K$.

\begin{lemma}\label{lem:isobaric derivation}
Let  $d\in\Z$ be such that for each $n$, $\Delta(Y_n)$ is $\degree$-homogeneous of degree $\degree(Y_n)+d$. Then $\Delta$ is $\degree$-homogeneous of degree $d$.
\end{lemma}
\begin{proof}
Let $\i\in\N^{(\N)}$ with $\degree(Y^\i)=i$. It is enough to show that then $\Delta(Y^{\i})$ is $\degree$-homogeneous of degree $i+d$. For such $\i$ we have
$$\Delta(Y^{\i})\ =\ \sum_{n\geq 0} \Delta(Y_n) \partial_n(Y^{\i})$$
with $\Delta(Y_n)$ $\degree$-homogeneous of degree $d_n+d$ and $\partial_n(Y^{\i})$ $\degree$-homogeneous of degree $i-d_n$, where $d_n:=\degree(Y_n)$.
\end{proof}

\noindent
Thus if $d\in\Z$ and $\Delta(Y_n)\in A_{d+1}$ for every $n$, then $\Delta$ is homogeneous of degree $d$. If
$w\in\Z$ and $\Delta(Y_n)\in A_{[n+w]}$ for every $n$, then $\Delta$ is isobaric of weight $w$.

\index{diagonal!derivation}
\index{diagonal!$n$-diagonal!derivation}
\index{derivation!diagonal}
\index{derivation!$n$-diagonal}

\begin{definition}\label{def:n-diagonal}
Let $\Delta\in\TrDer_K$. We call the triangular $K$-derivation $\Delta_n$ of $A$ with associated matrix $M_{\Delta_n}=(M_\Delta)_n$ the $n$\text{th} {\bf diagonal} or the {\bf $n$-diagonal} of $\Delta$.
We also say that $\Delta$ is {\bf $n$-diagonal} if $\Delta=\Delta_n$ and {\bf diagonal} if $\Delta=\Delta_0$.
\end{definition}

\noindent
In the metric on $\TrDer_K$ given by the filtration $(\TrDer^n_K)$ we have
$\sum_{i=0}^n \Delta_i \to \Delta$ as $n\to \infty$; more suggestively,
$$\Delta\ =\ \Delta_0 + \Delta_1 + \cdots + \Delta_n + \cdots.$$
Also, $\Delta\in\TrDer^m_K$ if and only if $\Delta_0=\cdots=\Delta_{m-1}=0$.
Each $n$-diagonal derivation is homogeneous of degree $0$ and isobaric of weight $-n$, by Lemma~\ref{lem:isobaric derivation}. 
In particular, for all $P\in A$ and
$\Delta\in\TrDer^m_K$ we have 
$$\deg\!\big(\Delta(P)\big)\ \leq\ \deg(P), \qquad \wt\!\big(\Delta(P)\big)\ \leq\ \wt(P)-m.$$
From Lemma~\ref{lem:formulas for diagonals} we obtain:

{\sloppy
\begin{lemma}\label{lem:bracket of diagonals}
Let $\Delta,\Lambda\in\TrDer_K$.
If $\Delta$ is $m$-diagonal and $\Lambda$ is $n$-diagonal, then~$[\Delta,\Lambda]$ is $(m+n)$-diagonal. % with associated matrix
%$$M_{[\Delta,\tilde\Delta]} = \diag_{m+n}(\Delta_{i,i+m}\cdot \tilde\Delta_{i+m,i+m+n}-\tilde\Delta_{i,i+n}\cdot \Delta_{i+n,i+n+m})_{i\geq 0}.$$
For each $k\in \N$ we have 
$$[\Delta,\Lambda]_k\ =\ \sum_{m+n=k} [\Delta_m,\Lambda_n].$$
\end{lemma}
}
\noindent
Also, if $\Delta\in\TrDer_K$ is $n$-diagonal, then $A^\Delta \supseteq K[Y_{<n}]$ and
\begin{equation}\label{eq:expansion of diagonal derivation}
\Delta\ =\ \Delta_{0n}Y_0\partial_n + \Delta_{1,n+1}Y_1\partial_{n+1}+ 
\Delta_{2,n+2}Y_2\partial_{n+2} +\cdots .
\end{equation}

\subsection*{Triangular algebra endomorphisms}
These are defined as follows.

\begin{definition}
An endomorphism $\sigma$ of the $K$-algebra $A$ is called {\bf triangular} (relative to $Y_0, Y_1,\dots$) if $\sigma$  restricts to a triangular 
endomorphism of the $K$-submodule $A_1=\bigoplus_{j} K\,Y_j$ of $A$
with respect to its basis $Y_0, Y_1,\dots$; that is, 
$$\sigma(Y_j)\ =\ \sigma_{0j} Y_0 + \sigma_{1j}Y_{1}+\cdots+\sigma_{jj}Y_j$$
where $\sigma_{ij}\in K$ for $i\leq j$. If in addition $\sigma_{jj}=1$ for all 
$j$, then $\sigma$ is {\bf unitriangular}. If
$\sigma$ is triangular with $\sigma_{ij}=0$ for all $i<j$, then $\sigma$ is called {\bf diagonal.}
\end{definition}

\index{endomorphism!triangular}
\index{endomorphism!unitriangular}
\index{endomorphism!diagonal}
\index{triangular!endomorphism}
\index{unitriangular!endomorphism}
\index{diagonal!endomorphism}

\noindent
Given a triangular $K$-algebra endomorphism $\sigma$ of $A$, set
$$M_\sigma\ :=\ M_{\sigma|A_1}\ =\ (\sigma_{ij})_{i,j\in\N}\in\tr_K,$$
with $\sigma_{ij}:=0$ for $i>j$, by convention.
The following is easy to verify:

\begin{lemma}
Let $\sigma$ be a  triangular $K$-algebra endomorphism of $A$. Then
\begin{enumerate}
\item[\textup{(i)}] $\sigma$ is bijective $\Longleftrightarrow$ $\sigma_{jj}\in K^\times$ for all $j$ $\Longleftrightarrow$ $M_\sigma\in\tr_K^\times$;
\item[\textup{(ii)}] $\sigma$ is unitriangular $\Longleftrightarrow$ $M_\sigma\in 1+\tr_K^1$;
\item[\textup{(iii)}] $\sigma$ is bijective and diagonal $\Longleftrightarrow$ $M_\sigma\in D^\times_K$.
\end{enumerate}
\end{lemma}

\nomenclature[Y]{$\TrAut_K$}{group of triangular $K$-algebra automorphisms}

\noindent
The triangular $K$-algebra automorphisms of $A$ form a subgroup  
$\TrAut_K$ 
of the group of all automorphisms of the $K$-algebra $A$. The map
$$\sigma\mapsto M_\sigma\colon\TrAut_K \to \tr_K^\times$$
is a group isomorphism. Let $m\ge 1 $. Define $\TrAut^m_K$ to be the subgroup of $\TrAut_K$ consisting of the triangular $K$-algebra automorphisms $\sigma$ of $A$ with $M_\sigma \in 1+\tr^m_K$.
For each $\Delta\in\TrDer^m_K$ we have
$\exp\Delta\in\TrAut^m_K$ and 
$\exp M_\Delta = M_{\exp\Delta}$. 
So we have a commuting diagram
$$\xymatrixcolsep{5pc}\xymatrix{
\TrDer^m_K \ar[d]^\exp \ar[r]^{\Delta\mapsto M_\Delta} & \tr^m_K \ar[d]^\exp \\
\TrAut^m_K \ar[r]^{\sigma\mapsto M_\sigma} & 1+\tr^m_K
}$$
where the horizontal arrows are group isomorphisms and the right (and hence also the left) vertical arrow is a bijection. If $\sigma\in\TrAut^m_K$, then 
$\log\sigma\in\TrDer^m_K$, $\log M_\sigma=M_{\log\sigma}$, and the automorphism 
$\sigma=\exp(\log\sigma)$ is exponential. 

\medskip\noindent
Let $\sigma$ be a triangular $K$-automorphism of $A$. Then for every triangular $K$-deri\-vation~$\Delta$ on $A$, the $K$-derivation $\sigma\Delta\sigma^{-1}$ of $A$ is triangular, the map $\Delta\mapsto \sigma\Delta\sigma^{-1}$ is an automorphism of the Lie algebra $\TrDer_K$ over $K$, and
$$M_{\sigma\Delta\sigma^{-1}}\ =\  M_\sigma\,M_\Delta\,(M_{\sigma})^{-1}.$$
If $m\ge 1$, $\Delta\in\TrDer^m_K$, then $\sigma\Delta\sigma^{-1}\in\TrDer^m_K$ 
and
$$\exp(\sigma\Delta\sigma^{-1})\ =\ \sigma\exp(\Delta)\sigma^{-1}.$$

%From Corollary~\ref{cor:BCH} and Lemma~\ref{lem:bracket of diagonals} we obtain:

%\begin{cor}
%Let $\Delta,\tilde\Delta\in\TrDer^1_K$, and set $\Sigma=\exp(\Delta)$, $\tilde\Sigma=\exp(\tilde\Delta)$. Then
%$$\log(\Sigma\,\tilde\Sigma) \equiv \Delta+\tilde\Delta+\textstyle\frac{1}{2}[\Delta,\tilde\Delta]\mod \TrDer_{K,3}.$$
%In particular
%$$\log(\Sigma\,\tilde\Sigma)_1 = \Delta_1+\tilde{\Delta}_1, \qquad  \log(\Sigma\,\tilde\Sigma)_2 = \Delta_2+\tilde{\Delta}_2+\textstyle\frac{1}{2}[\Delta_1,\tilde{\Delta}_1].$$
%\end{cor}

%Every $\Sigma\in\TrAut_K$ can be uniquely written in the form
%$$\Sigma = \Sigma_{\operatorname{e}} \Sigma_{\operatorname{d}}$$
%where $\Sigma_{\operatorname{e}}\in\TrAut_K$ is strictly triangular and $\Sigma_{\operatorname{d}}\in\TrAut_K$ is diagonal. We call $\Sigma_{\operatorname{e}}$ the {\bf exponential part} of $\Sigma$ and $\Sigma_{\operatorname{d}}$ the {\bf diagonal part} of $\Sigma$. If $\Sigma,\tilde\Sigma\in\TrAut_K$, then
%\begin{equation}\label{eq:e and d, 1}
%(\Sigma\tilde\Sigma)_{\operatorname{e}} = \Sigma_{\operatorname{e}}\, \big(\Sigma_{\operatorname{d}}\tilde{\Sigma}_{\operatorname{e}}(\Sigma_{\operatorname{d}})^{-1}\big), \qquad  (\Sigma\tilde\Sigma)_{\operatorname{d}} = \Sigma_{\operatorname{d}} \tilde{\Sigma}_{\operatorname{d}};
%\end{equation}
%in particular,
%\begin{equation}\label{eq:e and d, 2}
%(\Sigma^{-1})_{\operatorname{e}} = (\Sigma_{\operatorname{d}})^{-1} (\Sigma_{\operatorname{e}})^{-1} \Sigma_{\operatorname{d}}, \qquad (\Sigma^{-1})_{\operatorname{d}} = (\Sigma_{\operatorname{d}})^{-1}.
%\end{equation}

\subsection*{Companion derivations}
Let $\Delta\in\TrDer^1_K$ and $\sigma=\exp\Delta$. Then $\sigma$ is homogeneous of degree $0$; in fact $\deg\sigma(P)=\deg P$ for every $P\in A$.
Next, we investigate the isobaric parts of the images of polynomials under $\Delta$. The following lemma for~$P$ of the form $P=Y_n$ is already implicit in Lemma~\ref{lem:formulas for diagonals of exp and log}:

\begin{lemma}\label{lem:isobaric parts of Sigma(P)}
Let $\Delta\in\TrDer^1_K$ and put $\sigma:= \exp(\Delta)$. Suppose also that
$P\in A_{[w]}$, $P\ne 0$, with $w\in\N$. Then $\wt{\sigma(P)}=w$, and  
\begin{align*}
\sigma(P)_{[w]}\quad &=\ P, \\
\sigma(P)_{[w-1]}\ &=\ \Delta_1(P), \\
\sigma(P)_{[w-2]}\ &=\ \Delta_2(P) + \textstyle\frac{1}{2}(\Delta_1)^2(P), \\
\sigma(P)_{[w-3]}\ &=\ \Delta_3(P) + \textstyle\frac{1}{2}(\Delta_1\Delta_2+\Delta_2\Delta_1)(P) + \textstyle\frac{1}{6} (\Delta_1)^3(P) \\
&\ \vdots \\
\sigma(P)_{[w-n]}\ &=\ \sum_{i_1,\dots,i_k} \frac{1}{k!}(\Delta_{i_1}\cdots \Delta_{i_k})(P) \quad\qquad(n\ge 1),
\end{align*}
summed over the $(i_1,\dots, i_k)$ with $k\ge 1$, $i_1,\dots, i_k\ge 1$ and $i_1+ \cdots + i_k=n$.
\end{lemma}
\begin{proof} Note that $\Delta^k(P)=0$ for all big enough $k\in \N$. 
For all $i$,
$$\sigma(P)_{[i]}\ =\ 
P_{[i]} + \Delta(P)_{[i]} + \frac{1}{2}\Delta^2(P)_{[i]} + \cdots + \frac{1}{k!} \Delta^{k}(P)_{[i]}+\cdots.$$
Now $\wt\!\big(\Delta^k(P)\big)\leq\wt(P)-k <\wt P=w$ for $k\ge 1$, hence 
$\wt\!\big(\sigma(P)\big)=w$ and $\sigma(P)_{[w]}=P$. The family 
$(\Delta_i)_{i\ge 1}$ of endomorphisms of the $K$-module $A$ is summable, with
$\Delta=\sum_{i\ge 1}\Delta_i$. Let $k\ge 1$. Then
$$\Delta^k\ =\ \sum_{i_1,\dots,i_k} \Delta_{i_1}\cdots \Delta_{i_k}$$
summed over the $(i_1,\dots, i_k)$ with $i_1,\dots, i_k\ge 1$.
Each operator $\Delta_{i_1}\cdots \Delta_{i_k}$ in this sum is isobaric of weight
$-(i_1+ \cdots +i_k)$, so for all $n$,
$$\Delta^k(P)_{[w-n]}\ =\  \sum_{i_1,\dots,i_k} (\Delta_{i_1}\cdots \Delta_{i_k})(P),$$
summed over the $(i_1,\dots, i_k)$ with $i_1,\dots, i_k\ge 1$ and 
$i_1+ \cdots + i_k=n$.
This yields the lemma. Note also that $\big(\Delta^k(P)\big)_{[w-n]}=0$ if $k> n$. 
\end{proof} 

\index{invariant}

\noindent
As before, $\Delta\in\TrDer^1_K$ and $\sigma= \exp(\Delta)$. 
We say that an element $P$ of $A$ is {\bf $\sigma$-invariant} if~${\sigma(P)=P}$. 
We already observed that $P\in A$ is $\sigma$-invariant iff
$\Delta(P)=0$; see~\eqref{eq:invariants}.
We obtain a refinement of this fact for isobaric polynomials:

\begin{cor}\label{cor:invariants}
Let $P\in A_{[w]}$ where $w\in\N$. Then for $n=1,\dots,w$
$$\sigma(P)_{[w-1]} = \cdots = \sigma(P)_{[w-n]} = 0 \qquad
\Longleftrightarrow\qquad
\Delta_1(P) = \cdots = \Delta_n(P) = 0.$$ 
In particular:\  $\sigma(P)=P\ \Longleftrightarrow\ \text{$\Delta_n(P)=0$ for $n=1,\dots,w$.}$
\end{cor}
\begin{proof}
From the previous lemma we have, for $n=1,\dots,w$:
\begin{align*}\sigma(P)_{[w-n]}\ =\  \Delta_n(P)+\  &\Q \text{-linear combination of terms }(\Delta_{i_1}\cdots\Delta_{i_k})(P)\\
&\quad \text{ with $k>1$ and $1\leq i_1,\dots,i_k<n$.}
\end{align*}
%\text{\parbox{17.5em}{$\Q-$linear combination of the $(\Delta_{i_1}\cdots\Delta%_{i_k})(P)$ with $k>1$ and $1\leq i_1,\dots,i_k<n$.}} \end{cases}
%$$
This yields the corollary.
\end{proof}

\noindent
Thus if $P\in A$ is isobaric and $\sigma$-invariant, then
\begin{equation}\label{eq:Delta1 and Delta2}
\Delta_1(P)\ =\ \Delta_2(P)\ =\ 0.
\end{equation}
The polynomials $P\in A$, not necessarily isobaric, that satisfy \eqref{eq:Delta1 and Delta2}, comprise the $K$-subalgebra $A^{\Delta_1,\Delta_2}$ of $A$. Note that if $P\in A^{\Delta_1,\Delta_2}$ then $\Lambda(P)=0$ for all $\Lambda$ in the Lie subalgebra $\operatorname{Lie}(\Delta_1,\Delta_2)$ (over $K$) of $\TrDer_{K}$ generated by $\Delta_1$, $\Delta_2$; cf.~\eqref{eq:ring of constants}. 
This Lie algebra contains in particular the so-called 
{\bf companion} derivations
\begin{align*}
\Delta_1^\companion\ &=\ \Delta_1 \\
\Delta_2^\companion\ &=\ \Delta_2 \\
\Delta_3^\companion\ &=\ [\Delta_1,\Delta_2]\ =\ \operatorname{ad}_{\Delta_1}(\Delta_2) \\
\Delta_4^\companion\ &=\ [\Delta_1,\Delta_3^\companion]\ =\ [\Delta_1,\Delta_1,\Delta_2]\ =\ \operatorname{ad}^2_{\Delta_1}(\Delta_2)\\
\Delta_5^\companion\ &=\ [\Delta_1,\Delta_4^\companion]\ =\ [\Delta_1,\Delta_1,\Delta_1,\Delta_2]\ =\ \operatorname{ad}^3_{\Delta_1}(\Delta_2) \\
& \ \vdots \\
\Delta_n^\companion\ &=\ [\Delta_1,\Delta_{n-1}^\companion]\ =\ [\underbrace{\Delta_1,\dots,\Delta_1}_{\text{$n-2$ times}},\Delta_2]\ =\ \operatorname{ad}^{n-2}_{\Delta_1}(\Delta_2) \qquad\text{for $n\geq 2$.}
\end{align*}
Each  $\Delta_n^\companion$  is $n$-diagonal. We combine the $\Delta_n^\companion$ into the strictly triangular $K$-derivation
$$\Delta^\companion\ :=\ \Delta_1^\companion + \Delta_2^\companion + \cdots + \Delta_n^\companion + \cdots,$$
which we call the {\bf companion derivation} of $\Delta$. Note that
$\Delta_{0n}^\companion=(\Delta_n^\companion)_{0n}$ for $n\ge 1$, in particular,
$\Delta_{01}^\companion=\Delta_{01}$ and $\Delta_{02}^\companion= \Delta_{02}$. 

%For the next proposition we set $Y_{\leq n}:=(Y_0,\dots,Y_n)$, for each $n$.

\index{derivation!companion}
\index{companion!derivations}
\nomenclature[Y]{$\Delta^\companion$}{companion derivation of the derivation $\Delta$}

\begin{prop}\label{prop:constants of Delta1 and Delta2}
Suppose
$\Delta^\companion_{0n} \neq 0$ for all $n\ge 1$.
Then $A^{\Delta_1,\Delta_2}=K[Y_0]$.
\end{prop}
\begin{proof}
In connection with \eqref{eq:expansion of diagonal derivation} we already observed that $K[Y_0]\subseteq A^{\Delta_1,\Delta_2}$.
For the reverse inclusion, let $P\in A^{\Delta_1,\Delta_2}$. Take $n$ minimal with $P\in K[Y_0,\dots, Y_n]$. Towards a contradiction, 
suppose that $n\ge 1$. By \eqref{eq:expansion of diagonal derivation} and the hypothesis of the proposition, we have
$$\Delta_{n}^\companion\ =\ \Delta^\companion_{0,n} Y_{0}\, \partial_{n} + \Delta^\companion_{1,n+1} Y_{1}\, \partial_{n+1} + \cdots, \qquad \Delta^\companion_{0,n}\ \neq\ 0.$$
Thus $\Delta_{n}^\companion(P)=\Delta^\companion_{0,n} Y_{0}\frac{\partial P}{\partial Y_{n}}\neq 0$, a contradiction. 
\end{proof}

\begin{cor}\label{cor:constants of Delta1 and Delta2}
Suppose
$\Delta_{01} = \Delta_{02} = 0$ and $\Delta^\companion_{1,n} \neq 0$ for all $n\ge 2$.
Then $A^{\Delta_1,\Delta_2}=K[Y_0, Y_1]$, so all isobaric
$\sigma$-invariants belong to $K[Y_0, Y_1]$.
\end{cor}
\begin{proof}
Put $K^*:=K[Y_0]$ and view $A$ as the $K^*$-algebra of polynomials in the  indeterminates $Y_1, Y_2,\dots$ over $K^*$. Since 
$\Delta_{01}=\Delta_{02}=0$, $\Delta_1$ and $\Delta_2$ 
are triangular $K^*$-derivations on $A$
relative to $Y_1,Y_2,\dots$,  and the claim now follows from Proposition~\ref{prop:constants of Delta1 and Delta2} and \eqref{eq:Delta1 and Delta2}.
\end{proof}

\subsection*{The Stirling automorphism}
{\em In the rest of this section $K$ is an integral domain}. The Stirling automorphism
$\Upsilon$ of $A$ is the unitriangular 
automorphism of $A$ whose matrix $M_\Upsilon\in\tr_\Q$ has the signed Stirling numbers of the first kind as its entries:
$$M_\Upsilon = (\Upsilon_{ij}) :=
\left(\begin{array}{rrrrrrr}
1 & 0 & 0  & 0  & 0  & 0   & \cdots \\
  & 1 & -1 & 2  & -6 & 24  & \cdots \\
  &   & 1  & -3 & 11 & -50 & \cdots \\
  &   &    & 1  & -6 & 35  & \cdots \\
  &   &    &    &  1 & -10 & \cdots \\
  &   &    &    &    & 1   & \cdots \\
  &   &    &    &    &     & \ddots  
\end{array}\right)\text{ where $\Upsilon_{ij}=(-1)^{j-i}{j\brack i}$.}$$
By the recurrence relation \eqref{eq:recurrence for s(n,k)} for Stirling numbers of the first kind,
\begin{equation}\label{eq:derivative of M_Upsilon}
M_\Upsilon'\  =\  -M_\Upsilon\,D \qquad\text{where $D=\diag(0,1,2,3,\dots)$.}
\end{equation} 
The matrix $M_\Upsilon$ is an iteration matrix:

\index{automorphism!Stirling}
\index{Stirling!automorphism}
\nomenclature[Y]{$\Upsilon$}{Stirling automorphism}

\begin{lemma}\label{lem:Stirling as it matrix}
$\quad M_\Upsilon\ = \big\llb \log(1+z)\big\rrb$.
\end{lemma}
\begin{proof}
Let $f:=\log(1+z)\in z+z^2\Q[[z]]$. Then 
$$\big[ 1-(1/f') \big]\ =\ [-z]\ =\  -\diag_1(1,2,3,\dots)$$ and thus
$$S\cdot \big[ 1-(1/f') \big]\ =\  -\diag(0,1,2,3,\dots)\ =\ -D,$$
where $D$ is as in \eqref{eq:derivative of M_Upsilon}.
So by Corollary~\ref{cor:shift it matrix} and by \eqref{eq:derivative of M_Upsilon}, respectively, both $X=\llb f\rrb$ and $X=M_\Upsilon$ satisfy the differential equation $X'=-XD$, and have the same leftmost column. Thus
$\llb f\rrb=M_{\Upsilon}$ by Lemma~\ref{lem:recurrence}.
\end{proof}

\noindent
Since $\ex^z-1$ is the formal compositional inverse of $\log(1+z)$ in 
$z+z^2K[[z]]$, we get 
$$M_{\Upsilon^{-1}}\ =\ (M_\Upsilon)^{-1}\ =\ \llb \ex^z-1\rrb.$$
By Lemma~\ref{lem:properties of matrix derivative} and \eqref{eq:derivative of M_Upsilon}, its derivative satisfies the identity
$$M_{\Upsilon^{-1}}'\  =\  DM_{\Upsilon^{-1}}, \qquad D\ :=\  \diag(0,1,2,3,\dots),$$
which gives a recurrence relation for the entries of $M_{\Upsilon^{-1}}=(\Upsilon_{ij}^{-1})$: for all $i,j$
$$\Upsilon_{i,j+1}^{-1}\ =\ \Upsilon^{-1}_{i-1,j}+i\Upsilon^{-1}_{ij} 
\qquad(\Upsilon_{-1,j}:= 0 \text{ by convention})$$
with side conditions $\Upsilon_{00}^{-1}=1$ and $\Upsilon_{i0}^{-1}=\Upsilon_{0j}^{-1}=0$ for $i,j>0$.
Thus 
$$M_{\Upsilon^{-1}}\ =\
\begin{pmatrix}
1 & 0 & 0  & 0  & 0  & 0   & \cdots \\
  & 1 &  1 & 1  & 1  & 1   & \cdots \\
  &   & 1  & 3  & 7  & 15  & \cdots \\
  &   &    & 1  &  6 & 25  & \cdots \\
  &   &    &    &  1 &  10 & \cdots \\
  &   &    &    &    & 1   & \cdots \\
  &   &    &    &    &     & \ddots  
\end{pmatrix}.$$ 
From the recurrence relations we get 
$\Upsilon^{-1}_{ij}={j\brace i}$, where ${j\brace i}$ is by definition
the number of equivalence relations on a $j$-element set with exactly $i$ equivalence classes. These numbers are called Stirling numbers of the 
second kind. (See \cite[\S\S{}5.1, 5.3]{Comtet-Book}, \cite[\S{}6.1]{GKP}, \cite[p.~8]{Stirling}.) 

\index{Stirling!numbers of the second kind}
\nomenclature[Cd]{${j\brace i}$}{Stirling numbers of the second kind}

\medskip\noindent
The {\bf Stirling derivation} $\nabla:=\log\Upsilon\in\TrDer^1_K$  of $A$ has matrix 
$$M_\nabla\ =\ \left(\begin{array}{rrrrrrrr}
0 & 0 & 0  &           0  & 0            & 0   			& 0              &  \cdots \\[0.15em]
  & 0 & -1 & \frac{1}{2}  & -\frac{1}{2} & \frac{2}{3}  	& -\frac{11}{12} & \cdots \\[0.15em]
  &   & 0  & -3           & 2            & -\frac{5}{2} 	& 4              &  \cdots \\[0.15em]
  &   &    &  0           & -6      	    &  5  			& -\frac{15}{2}  &  \cdots \\[0.15em]
  &   &    &              &  0 			& -10 			& 10             &  \cdots \\[0.15em]
  &   &    &              &    			& 0   			& -15            &  \cdots \\[0.15em]
  &   &    &              &    			&     			& 0              &  \cdots \\[0.15em]
  &   &    &              &    			&     			&                &  \ddots 
\end{array}\right).$$\marginpar{only toprow and 0,1,2-diagonals of $M_{\nabla}$ have been checked. But the unchecked entries are not used}Since 
$$M_{\nabla}\ =\ \log M_{\Upsilon}\ =\ \log\, \llb \log(1+z)\rrb\ =\ \blla \!\operatorname{itlog}\!\big(\!\log(1+z)\big) \brra,$$
this gives
$$\operatorname{itlog}\!\big(\!\log(1+z)\big)\ =\ 
 - \frac{z^2}{2!} + \frac{1}{2} \frac{z^3}{3!}  -\frac{1}{2} \frac{z^4}{4!} + \frac{2}{3} \frac{z^5}{5!} -\frac{11}{12} \frac{z^6}{6!} + \cdots.$$
Note that by \eqref{eq:itlog it} in Corollary~\ref{cor:itlog} we have
$$%\begin{equation}\label{eq:itlog of log and exp}
\operatorname{itlog}\!\big(\!\log(1+z)\big)\ =\  -\operatorname{itlog}(\ex^z-1).
$$%\end{equation}
\noindent
In view of Definition~\ref{definfiteration} and setting
$$h\ :=\ \operatorname{itlog}\!\big(\!\log(1+z)\big)\ =\ \sum_{n=1}^{\infty} h_n \frac{z^n}{n!}\qquad\text{ where $h_n\in\Q$, $h_1=0$,}$$
the coefficients $h_n$ determine the diagonals of $\nabla$ as follows:

\index{Stirling!derivation}
\index{derivation!Stirling}
\nomenclature[Y]{$\nabla$}{Stirling derivation}

\begin{lemma}\label{diagnabla} 
$M_{\nabla_n}\ =\  h_{n+1}\diag_n\displaystyle\binom{i+n}{n+1}$ for each $n$. In particular,
$$   M_{\nabla_1}\ =\ -\diag_1\binom{i+1}{2},\qquad
M_{\nabla_2}\ =\ \frac{1}{2}\diag_2\binom{i+2}{3}.
$$
\end{lemma}

\noindent
The companion derivation $\nabla ^\companion$ of $\nabla$ has matrix
$$M_{\nabla^\companion}\ =\ \left(\begin{array}{rrrrrrrr}
0 & 0 & 0  &           0  & 0            & 0   			& 0   & \cdots \\
  & 0 & -1 & \frac{1}{2}  & 1            & 5          	& 45  & \cdots \\
  &   & 0  & -3           & 2            & 5          	& 30  & \cdots \\
  &   &    &  0           & -6      	    & 5  			& 15  & \cdots \\
  &   &    &              &  0 			& -10 			& 10  & \cdots \\
  &   &    &              &    			& 0   			& -15 & \cdots \\
  &   &    &              &    			&     			& 0   & \cdots \\
  &	  &    &			     &  				&				&	  & \ddots 
\end{array}\right).$$  
Here is an explicit formula for the entries of the matrix $M_{\nabla^\companion}$:

\begin{prop}\label{prop:companion of Stirling}
For $n\ge 1$ we have
\begin{equation}\label{eq:nabla'}
M_{\nabla^\companion_n}\ =\ \diag_n  c_n \binom{i+n}{n+1}
\end{equation}
where $c_1=-1$ and $c_n=\frac{(n-2)!(n+1)!}{3\cdot 2^n}$ for $n\geq 2$, so $c_2=\frac{1}{2}$. In particular, 
we have $\nabla^\companion_{0,n}=0$ for all $n$, and 
$\nabla^\companion_{1,n}=c_{n-1}\ne 0$ for all $n\ge 2$.
\end{prop}

\begin{proof}
Lemma~\ref{diagnabla} gives \eqref{eq:nabla'} for $n=1,2$, and then
Corollary~\ref{cor:binomial diagonal matrices} and the subsequent remark gives
\eqref{eq:nabla'} for $n>2$. 
\end{proof}

\noindent
From Proposition~\ref{prop:companion of Stirling} and Corollary~\ref{cor:constants of Delta1 and Delta2} we obtain:

\begin{cor}\label{cor:companion of Stirling}
Suppose $P\in A$ is isobaric. Then 
$$\Upsilon(P)=P\ \Longleftrightarrow\ \nabla_1(P)=\nabla_2(P)=0\ 
\Longleftrightarrow\  P\in K[Y_0,Y_1].$$
\end{cor}

\noindent
Since $(\nabla_2)_{i2}=0$ for all $i$ (see the matrix of $\nabla$), we have
$\nabla_2(Y_2)=0$, and so $K[Y_0, Y_1, Y_2]\subseteq A^{\nabla_2}$.
Thus by Corollary~\ref{cor:companion of Stirling}:

\begin{cor}\label{cor:companion of Stirling, variant}
Suppose $P\in K[Y_0,Y_1,Y_2]\subseteq A$ is isobaric. Then
$$ \nabla_1(P)=0\  \Longleftrightarrow\  P\in K[Y_0,Y_1].$$
\end{cor}

\subsection*{The algebra of partial differential operators}
The $K$-derivations $\partial_n=\frac{\partial}{\partial Y_n}$ of $A$
satisfy $[\partial_m,\partial_n]=0$ for all $m$,~$n$, so they 
generate a commutative subalgebra $K[\partial]:=K[\partial_0, \partial_1,\dots]$ of the $K$-algebra $\End(A)$ of
endomorphisms of the $K$-module~$A$. 
Let $\i$ range over the set $\N^{(\N)}$ of sequences $\i=(i_0,i_1,\dots)\in\N^\N$ such that 
$i_n=0$ for all but finitely many $n$, and likewise with $\j$ and 
$\k$. For such $\i$, set 
$$\partial^{\i}\ :=\ \partial_0^{i_0}\partial_1^{i_1} \cdots
\partial_n^{i_n}\cdots\in K[\partial], \  \abs{\i}:= i_0 + i_1+ \cdots, \ 
\i!\ :=\ i_0!i_1!\cdots i_n!\cdots\in\N^{\geq 1}.$$ 
The $\partial^{\i}$ generate the $K$-module $K[\partial]$. We have a $K$-bilinear map
\begin{equation}\label{eq:diff opp pairing}
(\Delta,P)\mapsto \<\Delta,P\>\ :=\ (\Delta P)\big\lvert_{Y_0=Y_1=\cdots=Y_n=\cdots=0}\ \colon\quad K[\partial]\times A\to K.
\end{equation}
It is easy to check that
$\<\partial^{\i},Y^{\j}\> = \i!$	if $\i=\j$ and $\<\partial^{\i},Y^{\j}\> = 0$ otherwise. Thus the pairing \eqref{eq:diff opp pairing} is non-degenerate, and $(\partial^{\i})$ is a basis for the $K$-module $K[\partial]$.
So for each $\Delta \in K[\partial]$ there is a unique family $(a_\i)$ in $K$ such that $a_{\i}=0$ for all but finitely many~$\i$ and
$\Delta = \sum_{\i} a_{\i} \partial^{\i}$. 
Put
\begin{align*} K[\partial_{<n}]\ &=\ K[\partial_0,\dots,\partial_{n-1}], \text{ so}\\ 
      K[\partial_{<m}]\ \subseteq\ K[\partial_{<n}]\ &\subseteq\ K[\partial]\qquad\text{for $m\leq n$, and}\\ 
K[\partial]\ &=\ \bigcup_{n} K[\partial_{<n}].
\end{align*}
Note that if $\Delta\in K[\partial]$, then 
$\Delta\big(K[Y_{<n}]\big)\subseteq K[Y_{<n}]$, and the $K$-linear map
$$%\begin{equation}\label{eq:restrict partial ders}
\Delta\mapsto \Delta|{K[Y_{<n}]}\ \colon\ K[\partial] \to \End\!\big(K[Y_{<n}]\big)
$$%\end{equation}
is injective on $K[\partial_{<n}]$.

\medskip
\noindent
For each family  $(a_\i)$ in $K$, the family $(a_{\i}\partial^\i)$ of endomorphisms of the $K$-module~$A$
is summable with respect to the trivial filtration on the $K$-module $A$, so we have an endomorphism $\Delta:=\sum_{\i} a_{\i} \partial^{\i}\in\End(A)$. We let $K[[\partial]]$ be the $K$-submodule of~$\End(A)$ consisting of the endomorphisms $\sum_{\i} a_{\i} \partial^{\i}$ where $(a_{\i})$ is a family in $K$.
By the properties of the pairing \eqref{eq:diff opp pairing} there is for each $\Delta\in K[[\partial]]$ a unique 
family~$(a_{\i})$ in $K$ with $\Delta=\sum_\i a_{\i}\partial^\i$.
For any families $(a_{\i})$ and $(b_{\j})$ in $K$,
$$ \left(\sum_{\i} a_{\i}\partial^{\i}\right)\left(\sum_{\j}b_{\j}\partial^{\j}\right)\ =\ \sum_{\k} \left(\sum_{\i+\j=\k}a_{\i}b_{\j} \right)\partial^{\k}.$$
Thus $K[[\partial]]$ is a commutative $K$-subalgebra of $\End(A)$ which contains $K[\partial]$ as a subalgebra.  
We call the elements of $K[[\partial]]$ {\bf partial differential operators} on $A$.
We say that $\Delta=\sum_{\i}a_{\i}\partial^{\i}\in K[[\partial]]$
as above is  {\bf homogeneous  of order $r\in\N$} if $a_\i= 0$ whenever $\abs{\i}\neq r$.
Note that if $P\in A$ has degree $\le d$ and $\Delta\in K[[\partial]]$ is homogeneous of order $r$, 
then $\deg \Delta(P) \leq d-r$; so if $d<r$, then $\Delta(P)=0$.

\index{operator!partial differential}
\index{order!partial differential operator}
\index{homogeneous!partial differential operator}
\nomenclature[M]{$K[[\partial]]$}{algebra of partial differential operators on $K\{Y\}$}

\medskip
\noindent
Let a sequence $\Phi=(\Phi_n)$ in $K[[\partial]]$ be given;  put
$$\Phi^\i\ :=\ \Phi_0^{i_0}\Phi_1^{i_1}\cdots \Phi_n^{i_n}\cdots\in K[[\partial]]. $$
For $\Delta=\sum_{\i}a_{\i}\partial^{\i}\in K[\partial]$ (all
$a_{\i}\in K$, and $a_{\i}=0$ for all but finitely many $\i$), set 
$$\Delta(\Phi)\ :=\ \sum_\i a_\i \Phi^\i\in K[[\partial]].$$
Routine arguments show:

\begin{lemma}
The map $\Delta\mapsto\Delta(\Phi)$ is the unique $K$-algebra morphism from $K[\partial]$ into $K[[\partial]]$ with $\partial_n(\Phi)=\Phi_n$ for each $n$. If $\Delta\in K[\partial]$ is homogeneous of order~$r$ and each $\Phi_n$ is homogeneous of order~$s$, then $\Delta(\Phi)$ is homogeneous of order~$rs$.
\end{lemma}

\subsection*{Notes and comments}
General references on locally nilpotent derivations and exponential automorphisms of \textit{finite-dimensional}\/ polynomial rings are \cite{vdEssen,Freudenburg,Nowicki}. Our usage of ``triangular'' differs from these sources, where a $K$-derivation~$\Delta$ of~$K[Y_{<n}]$ is called triangular if it satisfies the hypothesis of Lemma~\ref{lem:generalized triangular implies locally nilpotent}. 
Proposition~\ref{prop:vdEssen} is \cite[Proposition~2.1.3]{vdEssen}, with a different proof.
In connection with Pro\-po\-si\-tion~\ref{prop:constants of Delta1 and Delta2}, it follows from a  theorem of Maurer~\cite{Maurer} and Weitzen\-b\"ock~\cite{Weitzenboeck} that if $K$ is a field and $\Delta$ is a strictly triangular $K$-derivation of~$A$, then for each~$n$ the $K$-algebra $A^\Delta \cap K[Y_{<n}]$ is finitely generated; see also~\cite[Theorem~6.2.1]{Nowicki}.

 The power series $\operatorname{itlog}(\ex^z-1)\in\Q[[z]]\subseteq \C\(( z\)) $ is 
$\d$-transcendental over the differential subfield of 
$\C\(( z\)) $ consisting  of the 
convergent Laurent series, with $d/dz$ as the derivation on $\C\(( z\)) $; see~\cite{A-it2}.  \marginpar{results stated in these ``Notes and comments'' taken on faith}

\section{Application to Differential Polynomials}\label{appdifpol}

\nomenclature[Y]{$\Upsilon_\phi$}{triangular automorphism with matrix $\big( \phi^{-j} F^j_i(\phi) \big)_{i,j}$}

\noindent
{\em In this section $K$ is a differential field and $Y$ is a differential indeterminate over~$K$}. We apply the material in Section~\ref{sec:Derivations on polynomials rings}  to the $K$-algebra $A=K\{Y\}=K[Y_0,Y_1,\dots]$ where $Y_n=Y^{(n)}$ for each $n$. {\em Throughout this section $P\in A$ and $\phi\in K^\times$}. As in Section~\ref{Compositional Conjugation},
let $\derdelta=\phi^{-1}\der$ be the derivation of the compositional conjugate 
$K^\phi$ of~$K$. Consider the unitriangular $K$-automorphism $\Upsilon_\phi$ of $A$ defined by
\begin{multline}\label{eq:Upsilon_phi}
M_{\Upsilon_\phi}\  =\  \big( \phi^{-j} F^j_i(\phi) \big)_{i,j}\ =\\  
   \left(
\begin{array}{llllll}
1 & 0 & 0            & 0              & 0                                 & \cdots \\
  & 1 & \phi'/\phi^2 & \phi''/\phi^3  & \phi^{(3)}/\phi^4                 & \cdots \\
  &   & 1            & 3\phi'/\phi^2  & 4\phi''/\phi^3+3(\phi')^2/\phi^4  & \cdots \\
  &   &              & 1              & 6\phi'/\phi^2                     & \cdots \\
  &   &              &                &  1                                & \cdots \\
  &   &              &                &                                   & \ddots 
\end{array}\right).
\end{multline}
Expressed in terms of the $\derdelta^n(\phi)$, this is
\begin{multline*}
M_{\Upsilon_\phi}\ =\ \big( \phi^{-j} G^j_i(\phi) \big)_{i,j}\ = \   \left(
\begin{array}{llllll}
1 & 0 & 0                      & 0                                                                                & \cdots \\
  & 1 & \derdelta(\phi)/\phi   & \derdelta^2(\phi)/\phi+\derdelta(\phi)^2/\phi^2  & \cdots \\
  &   & 1                      & 3\derdelta(\phi)/\phi                                          & \cdots \\
  &   &                                                                              & 1       & \cdots \\
    &                        &                                                       &                                   & \ddots 
\end{array}\right).
\end{multline*}
We also define the diagonal $K$-automorphism $\Xi_\phi$ of $A$ by \nomenclature[Y]{$\Xi_\phi$}{diagonal automorphism with matrix $\diag(\phi^i)$}
$$M_{\Xi_\phi}\  =\  \diag (\phi^i)\  =\ 
\begin{pmatrix} 
1 &		&		&		& \\
  & \phi	&		&		& \\
  &		& \phi^2	&		& \\
  &		&		& \phi^3	& \\
  &		&		&		& \ddots
\end{pmatrix}.$$ 
Thus $P^\phi = (\Upsilon_\phi\circ\Xi_\phi)(P)$.
If $P$ is isobaric of weight $w\in\N$, then $\Xi_\phi(P)=\phi^w P$ and hence $P^\phi=\phi^w\,\Upsilon_\phi(P)$. 
We set
$$   \nabla_{\phi}\ :=\ \log\Upsilon_\phi \in \TrDer^1_K.$$
Then by Lemma~\ref{lem:isobaric parts of Sigma(P)} we have: 

\nomenclature[Y]{$\nabla_\phi$}{logarithm of $\Upsilon_\phi$}

\begin{lemma}\label{lem:isobaric components of Pphi}
Suppose $P$ is isobaric of weight $w\in\N$. Then $(P^\phi)_{[w]}=\phi^w P$ and

\smallskip
\noindent
for $n\ge 1$:\hfil  $(P^\phi)_{[w-n]}\  =\  \phi^w \displaystyle\sum_{i_1,\dots,i_k} \frac{1}{k!}(\nabla_{\phi,i_1}\cdots \nabla_{\phi,i_k})(P)$,  \hfil

\smallskip
\noindent
summed over the $(i_1,\dots, i_k)$ with $k\ge 1$, $i_1,\dots, i_k\ge 1$ and $i_1+ \cdots + i_k=n$. So
$$(P^\phi)_{[w-1]}\  =\  \phi^w \nabla_{\phi,1}(P), \qquad
  (P^\phi)_{[w-2]}\ =\  \phi^w \big(\nabla_{\phi,2}+\textstyle\frac{1}{2}(\nabla_{\phi,1})^2\big)(P).$$
\end{lemma}

\noindent
In the next lemma $x\in K$ satisfies $x'=1$, and $t:=\frac{1}{x}$, so
$\Upsilon_{t}=\Upsilon$ by \eqref{Stirling}, and thus $\nabla_t=\nabla$, and
$\nabla_{t,1}=\nabla_1$ as well as $\nabla_{t,2}=\nabla_2$. 

\begin{lemma}\label{cor:P vanish}
Suppose $P$ is isobaric of weight $w\in\N$.
Then the following conditions are equivalent:
\begin{enumerate}
\item[\textup{(i)}] $P\in K[Y,Y']$,
\item[\textup{(ii)}] $P\in K[Y]\cdot (Y')^w$,
\item[\textup{(iii)}] $P^\phi$ is isobaric of weight $w$, for every $\phi$,
\item[\textup{(iv)}] $P^t$ is isobaric of weight $w$,
\item[\textup{(v)}] $P^t=t^wP$.
\end{enumerate}
\end{lemma}

\noindent
In (v), $t^wP$ lies in the differential ring $K\{Y\}$ and $P^t$
in its compositional conjugate~$K^t\{Y\}$, but the equality makes sense as these differential rings have the same underlying ring $A=K[Y,Y',\dots]$.

\begin{proof}
(i)~$\Rightarrow$~(ii)~$\Rightarrow$~(iii)~$\Rightarrow$~(iv)
is clear, and (iv)~$\Rightarrow$~(v) follows from Lemma~\ref{lem:isobaric components of Pphi}.  To show (v)~$\Rightarrow$~(i), suppose $P^t=t^wP$. Then $\Upsilon(P)=\Upsilon_t(P)=t^{-w}P^t=P$ and hence $P\in K[Y,Y']$ by Corollary~\ref{cor:companion of Stirling}.
\end{proof} 

\noindent
By Lemma~\ref{lem:Stirling as it matrix}, the matrix $M_\Upsilon$ is the iteration matrix of the formal power series~$\log(1+z)$. More generally, the matrix $M_{\Upsilon_\phi}$ from \eqref{eq:Upsilon_phi} is an iteration matrix as a consequence of the next proposition:  

\begin{prop}\label{prop:itmatrix F}
Set 
$$f_\phi\ :=\ \sum_{n=1}^{\infty} \phi^{(n-1)}\frac{z^n}{n!}\ =\ \phi z + \phi' \frac{z^2}{2!} + \phi'' \frac{z^3}{3!} + \cdots \in zK[[z]]$$ 
and
$$F_\phi\ :=\ \big(F^j_i(\phi)\big)_{i,j\in\N}\ =\
 \left(
\begin{array}{llllll}
1 & 0 	& 0     & 0			& 0                                 & \cdots \\
  & \phi	& \phi' & \phi''  	& \phi^{(3)}                 & \cdots \\
  &   	& \phi^2& 3\phi\phi'	& 4\phi\phi''+3(\phi')^2  & \cdots \\
  &   	&       & \phi^3     & 6\phi^2\phi'                     & \cdots \\
  &   	&       &            & \phi^4                                & \cdots \\
  &   	&       &            &                                   & \ddots 
\end{array}\right)\in\tr_K.$$
Then $\llb f_\phi\rrb=F_\phi$.
\end{prop}
\begin{proof}
This amounts to showing that for $i\le j$,
$$B_{ij}\big(\phi,\phi',\dots,\phi^{(j-i)}\big)\ =\  F^j_i(\phi)$$
where $B_{ij}(y_1,\dots,y_{j-i+1})\in\Q[y_1,\dots,y_{j-i+1}]$ is the Bell polynomial defined in Section~\ref{sec:Iteration matrices}. Let $x,z$ be distinct indeterminates, and note that $\frac{\partial}{\partial x}$ and
$\partial:=\frac{\partial}{\partial z}$ are commuting $K$-derivations on
$K[[x,z]]$. With $R:=K[x]$, these two derivations map the subring $R[[z]]$
of $K[[x,z]]$ into $R[[z]]$. With the nonnegative filtration $(z^nR[[z]])$ on $R[[z]]$, the restriction of $\partial$ to $R[[z]]$ is a continuous $R$-derivation on $R[[z]]$.  Set $a:= \ex^{xf_{\phi}}\in 1+zR[[z]]$, so
\begin{equation}\label{abell} a\ =\ \sum_{j=0}^{\infty}\left(\sum_{i=0}^j B_{ij}(\phi, \phi',\dots)x^i\right)\frac{z^j}{j!}\ \text{ in $R[[z]]$}\end{equation}
by Corollary~\ref{corprop}. Also 
$\frac{\partial a}{\partial x}=f_{\phi}a$. Below we set
$$ \theta\ :=\ \partial(f_\phi)\ =\  \phi+\phi'z + \phi'' \frac{z^2}{2!}+\cdots \in K[[z]]^\times\subseteq R[[z]]^\times,$$ 
so
$\partial(a)=x\theta a$. Hence with $\delta$ denoting the derivation $\theta^{-1}\partial$ of $R[[z]]$, we have $\delta(a)=\theta^{-1}\partial(a)=xa$, so
$\delta^n(a) = x^n a$ for all $n$ and thus, by the definition of the differential polynomials $F_n^j$:
\begin{equation}\label{eq:powers of x}
a^{-1} \partial^j(a)\ =\ \sum_{n=0}^j F_n^j(\theta) \,x^n.
\end{equation}
Let $i\le j$. Applying $\frac{\partial^i}{\partial x^i}$ to the left-hand side of \eqref{eq:powers of x} and using 
$\frac{\partial a}{\partial x}=f_{\phi}a$ yields
\begin{align*}
\frac{\partial^i}{\partial x^i} \left(a^{-1}\partial^j(a)\right)\ &=\ 
\frac{\partial^i}{\partial x^i} \left(a^{-1}\frac{\partial^j a}{\partial z^j}\right)\ =\ 
 \sum_{k=0}^i {i\choose k} \frac{\partial^k a^{-1} }{\partial x^k} \cdot \frac{\partial^{i+j-k}  a}{\partial x^{i-k}\partial z^j}\\ 
&=\ \sum_{k=0}^i {i\choose k} (-1)^k f_{\phi}^k a^{-1} \cdot \frac{\partial^{i+j-k}  a}{\partial x^{i-k}\partial z^j}
\end{align*}
in $K[[x,z]]$, and so by~\eqref{abell},
\begin{equation}\label{eq:eval, 1}
\frac{1}{i!}\frac{\partial^i}{\partial x^i} \left(a^{-1}\partial^j(a)\right) \Bigg\lvert_{x=z=0}\ =\ 
\frac{1}{i!}\frac{\partial^{i+j}a}{\partial x^i\partial z^j} \Bigg\lvert_{x=z=0}\ =\ B_{ij}(\phi,\phi',\dots).
\end{equation}
Applying $\frac{1}{i!}\frac{\partial^i}{\partial x^i}$ to the right-hand side of \eqref{eq:powers of x} yields
$$ \frac{1}{i!} \frac{\partial^i}{\partial x^i} \left(\sum_{n=0}^j F_n^j(\theta) \,x^n\right)\ =\ F_i^j(\theta)+\text{terms of positive degree in $x$, in $K[[x,z]]$,}$$
hence
$$ \frac{1}{i!} \frac{\partial^i}{\partial x^i} \left(\sum_{n=0}^j F_n^j(\theta) \,x^n\right) \Bigg\lvert_{x=0}\ =\  F_i^j(\theta).$$
Also, for $k\in \N$,
$$\partial^k\theta\ =\ \partial^{k+1}(f_\phi)\ =\ \phi^{(k)} + \phi^{(k+1)} z + \phi^{(k+2)} \frac{z^2}{2!} + \cdots,$$
hence $\partial^k\theta\big\lvert_{z=0}=\phi^{(k)}$, and thus
\begin{equation}\label{eq:eval, 2}
\frac{1}{i!} \frac{\partial^i}{\partial x^i} \left(\sum_{n=0}^j F_n^j(\theta) \,x^n\right) \Bigg\lvert_{x=z=0}\ =\  F_i^j(\theta)\big\lvert_{z=0}\ =\ F_i^j(\phi).
\end{equation}
Equality \eqref{eq:powers of x} together with \eqref{eq:eval, 1} and \eqref{eq:eval, 2} now yields the claim.
\end{proof} 

\noindent
As desired, we can now interpret $M_{\Upsilon_\phi}$ as an iteration matrix. Recall from Section~\ref{The Riccati Transform} the definition of the Riccati polynomials $R_n$.

\begin{cor}\label{cor:itmatrix F}
Let
$$g_\phi:\ =\ \sum_{n=1}^{\infty} \left(\frac{\phi^{(n-1)}}{\phi^n}\right) \frac{z^n}{n!}\ =\ \sum_{n\geq 1} \left(\frac{R_{n-1}(\phi^\dagger)}{\phi^{n-1}}\right) \frac{z^n}{n!}\in z+z^2K[[z]].$$
Then $\llb g_\phi\rrb=M_{\Upsilon_\phi}$.
\end{cor}
\begin{proof}
Let $f_\phi$ and $F_\phi$ be as in Proposition~\ref{prop:itmatrix F}. Then we have $g_\phi=f_\phi\circ (\phi^{-1}z)$, hence
$\llb g_\phi\rrb=\llb f_\phi\rrb\cdot \llb\phi^{-1}z\rrb = F_\phi\cdot \diag(\phi^{-i})=M_{\Upsilon_\phi}$, as claimed.
\end{proof}

\noindent
The next proposition shows that the derivation $\nabla_{\phi,n}$ is
a certain scalar multiple of~$\nabla_n$ for $n=1,2$. This fact will be very useful 
in Chapter~\ref{ch:The Dominant Part and the Newton Polynomial}:

\begin{prop}\label{prop:nabla_phi,i}
For $n=1,2$ there is a polynomial $G_n\in\Q[Y_1,\dots,Y_{n}]$, isobaric of weight $n$ and independent of $\phi$ and $K$,
such that
$$\nabla_{\phi,n}\ =\ \phi^{-n} G_n\big(R_1(\phi^\dagger),\dots,R_n(\phi^\dagger)\big) \nabla_n.$$
In more detail, with $\lambda:=-\phi^\dagger$ and $\omega:=-(2\lambda'+\lambda^2)=2(\phi^{\dagger})'-(\phi^\dagger)^2$, we have
\begin{align*}
\nabla_{\phi,1}\ 	&=\ \phi^{-1}\big(-R_1(\phi^\dagger)\big) \,\nabla_1\ =\ (\lambda/\phi) \,\nabla_1, \\
\nabla_{\phi,2}\		&=\ \phi^{-2}\big(2R_2(\phi^\dagger)-3R_1(\phi^\dagger)^2\big)\,\nabla_2\ =\ (\omega/\phi^2)\,\nabla_2.
\end{align*}
\end{prop}

\noindent 
The proof below shows: if $n\ge 1$ and the coefficient of $z^{n+1}$ in $\operatorname{itlog}\!\big(\!\log(1+z)\big)$ is not zero, then there is a $G_n$ as in the first sentence of the lemma. 

\begin{proof}
Let $g=g_\phi=\sum_{n=1}^{\infty} g_n \frac{z^n}{n!}$, where $g_n=\phi^{(n-1)}/\phi^n$ for $n\geq 1$, be the unitary power series from Corollary~\ref{cor:itmatrix F}, and let $h=\sum_{n=2}^{\infty} h_n\frac{z^n}{n!}$ ($h_n\in K$ for $n\geq 2$) be the iterative logarithm of $g$. Then $M_{\nabla_{\phi}}=\log\, \llb g \rrb = \lla h \rra $. Let $n\ge 1$. Then
$$M_{\nabla_{\phi,n}}\ =\ h_{n+1}\diag_n {i+n\choose n+1}.$$ 
Recall that we have an isobaric polynomial $H_{n+1}\in\Q[Y_1,\dots,Y_n]$ of weight $n$ with
$$h_{n+1}\ =\ H_{n+1}(g_2,\dots,g_{n+1})\ =\ \phi^{-n}H_{n+1}\big(R_1(\phi^\dagger),\dots,R_n(\phi^\dagger)\big).$$
It remains to use Lemma~\ref{diagnabla} and \eqref{eq:h2, h3}. 
\end{proof}

\begin{remark}
It is sometimes convenient to express the transformation factors in the previous proposition
in terms of the  derivation $\derdelta=\phi^{-1}\der$ of $K$:
$$\phi^{-1}\lambda\ =\ -\phi^{-1}\derdelta(\phi)\ \qquad \phi^{-2}\omega\ =\
\phi^{-2}\big(2\phi\derdelta^2(\phi)-\derdelta(\phi)^2\big).$$
\end{remark}

\noindent
From the previous proposition and Lemma~\ref{lem:isobaric components of Pphi} we obtain:

\begin{cor}\label{cor:nabla_phi,i}
Suppose $P$ is nonzero and set $w=\wt(P)$ and  $Q:=P^\phi$. Then with $\lambda=-\phi^\dagger$ and $\omega=-(2\lambda'+\lambda^2)$,
\begin{align*}
Q_{[w]}\ &=\ \phi^w P_{[w]}, \\
Q_{[w-1]}\ &=\ \phi^{w-1}\big[ P_{[w-1]} + \lambda \nabla_1(P_{[w]}) \big], \\
Q_{[w-2]}\ &=\ \phi^{w-2}\big[ P_{[w-2]} + \lambda \nabla_1(P_{[w-1]}) + 
\textstyle \left(\omega\nabla_2+\frac{1}{2}\lambda^2\nabla_1^2\right)(P_{[w]})\big].
\end{align*}
\end{cor}

\subsection*{Additive and multiplicative conjugates of partial differential
operators} This is for use in Section~\ref{sec:unravelers} below 
and concerns partial differential operators on the $K$-algebra $A=K\{Y\}$, using the notions defined at the end of Section~\ref{sec:Derivations on polynomials rings}.  

\begin{lemma}\label{addconjDelta} Let $\Delta\in K[[\partial]]$ and $h\in K$. Then
$$(\Delta  P)_{+ h}\ =\ \Delta (P_{+ h}).$$
\end{lemma}
\begin{proof} For $\Delta=\partial_i$, this is
\eqref{eq:partial ders and operations, 1}. It extends easily, first to
products $\Delta=a\partial_0^{i_0}\cdots \partial_n^{i_n}$ with $a\in K$, and next to $\Delta$ as an infinite sum of such products.  
\end{proof}

\noindent 
As to multiplicative conjugation, let $\Delta$ range over $K[\partial]$ rather than over $K[[\partial]]$. 

\begin{lemma}\label{multconjDelta}
Suppose $h\in K^\times$.
There is a unique $K$-algebra morphism 
$$\Delta\mapsto\Delta_{\times h}\ \colon\ K[\partial] \to K[[\partial]]$$ such that for all $\Delta$,~$P$ we have
\begin{equation}\label{eq:diff op mult conj}
(\Delta_{\times h} P)_{\times h}\ =\ \Delta( P_{\times h}).
\end{equation}
If $\Delta$ is homogeneous of order $r$, then so is $\Delta_{\times h}$.
\end{lemma}
\begin{proof}
If $\Delta\mapsto\Delta_{\times h}\colon K[\partial]\to K[[\partial]]$ is a $K$-algebra morphism satisfying \eqref{eq:diff op mult conj} for
all $\Delta$,~$P$, then
$\Delta_{\times h} P = \big(\Delta (P_{\times h})\big)_{\times h^{-1}}$  for all $\Delta, P$, so there is at most one
such $K$-algebra morphism. For existence,
let $\Delta\mapsto\Delta_{\times h}$  be the unique $K$-algebra morphism from $K[\partial]$ to $K[[\partial]]$ sending each $\partial_i$ to
$$(\partial_i)_{\times h}\ :=\ \sum_{j\geq i} {j\choose i} h^{(j-i)} \partial_j.$$
Then \eqref{eq:diff op mult conj} holds for all $\Delta$,~$P$ by the identity \eqref{eq:partial ders and operations, 2}. 
\end{proof}

\noindent
For $h=1$ we have $\Delta_{\times 1}P=\Delta P$, and for $h\in K^\times$, 
$$\deg \Delta_{\times h}(P)\ =\ \deg \Delta(P_{\times h}), \qquad \wt \Delta_{\times h}(P)\ =\  \wt \Delta(P_{\times h}).$$

\marginpar{chapter checked except as indicated in marginal comments}

\subsection*{Notes and comments} Bank~\cite[Lemma~13]{Bank66}
proves
Lemma~\ref{cor:P vanish} for~$P\in C\{Y\}$, $C=\C$,
by analytic techniques. Babakhanian~\cite[Theorem~8.3]{Babakhanian68} has
an algebraic proof, different from ours, for $P\in C\{Y\}$.
Similar results
play a role in 
the Newton diagram method for differential polynomials developed by Bank~\cite{Bank66} and Strodt~\cite{Strodt54,Strodt57}. It might be interesting
to relate the Bank-Strodt method to our 
Newton diagram method from Chapters~\ref{ch:The Dominant Part and the Newton Polynomial} and~\ref{ch:newtonian fields}.

%% file: mt-13.tex
\chapter{The Newton Polynomial}
\label{ch:The Dominant Part and the Newton Polynomial}
\setcounter{theorem}{0}

\noindent
{\em In this chapter $K$ is a $\d$-valued field of $H$-type with asymptotic
integration and small derivation. We also assume that $K$ is equipped with
a monomial group $\frak{M}$}, and let~$\fm$,~$\fn$ range over $\frak{M}$. As usual, 
$(\Gamma,\psi)$ is the asymptotic couple of $K$. The unique element of
$\Gamma^{>}$ fixed by $\psi$ is denoted by $1$, so $\Psi < 1+1$.  We let 
$\gamma$ range over $\Gamma$, and~$\phi$ over the active elements
of $\frak{M}$ in $K$, so $K^\phi$ inherits the properties we imposed on $K$. Throughout, $P\in K\{Y\}^{\ne}$. We now present
an overview of the main results to be established in this chapter.

\medskip\noindent
Since $K$ is $\d$-valued, we have an isomorphism
$$c\mapsto \bar{c}=c+\smallo\colon\ C \to \k=\mathcal{O}/\smallo$$
from its constant field $C$ (with its trivial derivation) onto the differential residue field~$\k$ of $K$, and below we identify
$C$ with $\k$ via this isomorphism. We also extend
the residue map $a\mapsto\bar{a}\colon \mathcal{O}\to C$ to
the differential ring morphism 
$$Q\mapsto\bar{Q}\colon\  \mathcal{O}\{Y\}\to C\{Y\}$$
that sends $Y^{(n)}$ to $Y^{(n)}$ for each $n$. 
The constant field $C$ of $K$ is also the constant field of
$K^{\phi}$, and so $C\{Y\}$ is a common differential subring of
all $K^{\phi}\{Y\}$. As we did with~$K$, we identify $C$ with the differential residue field of $K^{\phi}$, and extend the residue map $a\mapsto\bar{a}\colon \mathcal{O}^{\phi}\to C$ to
the differential ring morphism 
$$Q\mapsto\bar{Q} \colon\ \mathcal{O}^{\phi}\{Y\}\to C\{Y\},$$ 
where $\mathcal{O}^{\phi}$ is the valuation ring $\mathcal{O}$ of $K^{\phi}$ viewed as a differential subring of $K^{\phi}$.

\medskip\noindent
We now use $\frak{M}$ to associate to 
$P\in K\{Y\}^{\ne}$ its dominant monomial $\frak d_P\in \frak{M}$ and its dominant part $D_P\in C\{Y\}$:
$$\frak d_P\in \frak{M}\ \text{ with $\frak d_P\asymp P$,} \qquad \quad D_P = \bar{\frak d_P^{-1}P}\in C\{Y\}.$$
(This is in agreement with Section~\ref{sec: dominant}. Another choice of $\mathfrak M$ would multiply $D_P$ by a factor
in $C^\times$.)
Likewise, $\frak d_{P^{\phi}}\in \frak{M}$ and 
$D_{P^{\phi}}\in C\{Y\}$ for each $\phi$.
As in Chapter~\ref{evq} a condition $S(\phi)$ on elements $\phi$ is said to hold {\em eventually\/} if there is an active $\phi_0$ in~$K$ such that $S(\phi)$ holds for all $\phi\preceq \phi_0$.

\begin{prop}\label{newexist} Given any $P$ there exists a differential polynomial $N\in C\{Y\}$
such that eventually $D_{P^{\phi}}=N$.
\end{prop}

\index{Newton!polynomial}
\index{polynomial!Newton}
\index{differential polynomial!Newton polynomial}
\nomenclature[X]{$N_P$}{Newton polynomial of $P$}

\noindent
This fact is derived in Section~\ref{sec:The Dominant Part of a Differential Polynomial}. We define the {\bf Newton polynomial\/} of $P$ to be the differential polynomial $N$ in Proposition~\ref{newexist}, and denote it by $N_P$. Section~\ref{basicnewton} contains some elementary results about these Newton polynomials. In Section~\ref{The Shape of the Newton Polynomial} we establish a key consequence of $\upo$-freeness: 

\begin{theorem} Suppose $K$ is $\upo$-free. Then $N_P\in C[Y](Y')^{\N}$
for all $P$. 
\end{theorem} 

\noindent
For $\upo$-free $K$ this allows us to describe explicitly  
the behavior of $P(y)$ near the constant field,   
not just in $K$, but in any $\d$-valued field extension of $K$ of $H$-type.

Another key consequence of $\upo$-freeness proved in this chapter 
concerns \textit{eventual equalizers.}\/ To explain this,
let $\Gamma$ be divisible, and let $P,Q\in K\{Y\}^{\neq}$ be homogeneous of 
degrees $d>e$. 
By the Equalizer Theorem from Chapter~\ref{ch:valueddifferential} there is 
for each $\phi$ an $a\in K^\times$ such that 
$P^\phi_{\times a}\asymp Q^\phi_{\times a}$. We show in Section~\ref{eveqth} 
that for sufficiently high $v\phi$ we can take such $a$ independent of 
$\phi$, provided $K$ is $\upo$-free:

\begin{theorem}\label{evtequalizer} If $K$ is $\upo$-free, $\Gamma$ is 
divisible,
and $P,Q\in K\{Y\}^{\neq}$ are homogeneous of degrees $d>e$, then there exists $a\in K^\times$ such that, eventually, 
$ P^\phi_{\times a}\  \asymp\  Q^\phi_{\times a}$.
\end{theorem}

\index{theorem!Eventual Equalizer Theorem}

\noindent
In Section~\ref{consupofree} we consider more generally any ungrounded
$H$-asymptotic field $E$ with $\Gamma_E\ne \{0\}$ and prove for such $E$:

\begin{theorem}\label{upopreservation} If $E$ is $\upo$-free and $F$ is a 
$\d$-algebraic
$\d$-valued field extension of~$E$ of $H$-type, then $\Gamma_{E}^{<}$ is cofinal in $\Gamma_F^{<}$, and $F$ is $\upo$-free.
\end{theorem}

\noindent
In Section~\ref{consupofree} we also extend the Eventual Equalizer Theorem~\ref{evtequalizer} to $\upo$-free $E$. 
Section~\ref{fcuplfr} contains the construction, for
any $E$ that is $\upl$-free but not $\upo$-free, of a canonical extension $E\<\upg\>$ generated over $E$
by a solution $\upg$ of a certain second-order differential equation. Section~\ref{sec:aseq} on unraveling asymptotic equations is technical, but
crucial in the next chapter. The last section describes some
concrete $H$-fields $\R\<\upo\>\subseteq \R\<\upl\>\subseteq \R\<\upg\>$ with
interesting generic features.

\section{Revisiting the Dominant Part}
\label{sec:The Dominant Part of a Differential Polynomial}

\noindent
We derived some basic facts on the dominant part of $P$ in Section~\ref{sec: dominant},
and here we add to this in the more special setting of this chapter.
Recall that in this setting $D_P$ is defined using our monomial group $\frak{M}$ of $K$. Thus 
$\frak d_{\m P}=\m \frak d_P$ and $D_{\m P}=D_P$. If $a\in\mathcal O$ and $P(a)=0$, then $D_P(\overline{a})=0$.

\subsection*{Elementary facts on the dominant part} 
We represent $P$ as
$$P\ =\ \frak d_PD_P + R_P \quad\text{with $R_P\in K\{Y\}$, $R_P\prec P$.}$$
Note that if $R_P\ne 0$, then 
$$\wv(P)\ \leq\ \wv(R_{P})\ \leq\ \wt(R_{P})\ \leq\  \wt(P).$$
It will also be convenient to define $\frak d_Q$ and $D_Q$ for $Q=0\in K\{Y\}$ by $\frak d_0:=0\in K$ and $D_0:= 0\in C\{Y\}$. If $Q\in K\{Y\}$, then
$\frak d_{PQ}=\frak d_P\frak d_Q$ and $D_{P Q} = D_P D_Q$.

\begin{lemma}\label{dn1} 
Let $Q\in C\{Y\}$, $Q\notin C$. Then $P(Q)\ne 0$, and 
$$\frak d_{P(Q)}\ =\ \frak d_P, \qquad D_{P(Q)}\ =\ D_P(Q),\qquad R_{P(Q)}\ =\ R_P(Q).$$
\textup{(}In particular,   $D_{P_{+c}}=(D_P)_{+c}$ for  $c\in C$ and 
$D_{P_{\times c}}=(D_P)_{\times c}$ for  $c\in C^\times$.\textup{)}
\end{lemma}
\begin{proof}
We have
$$P(Q)\ =\ \frak d_P D_P(Q) + R_P(Q)$$
where $R_P(Q)\prec P$ and $D_P(Q)\in C\{Y\}^{\neq}$ by Lemma~\ref{lem:comp3}.
So $P(Q)\asymp\frak d_P$, hence $\frak d_{P(Q)}=\frak d_P$, and thus
$D_{P(Q)}=D_P(Q)$, and $R_{P(Q)}=R_P(Q)$.
\end{proof}

\begin{lemma}\label{dn1a} Suppose $P$ is homogeneous. Then $D_P$ is 
homogeneous and 
$$       D_{\Ric(P)}\ =\ \Ric(D_P)\ \text{ in $C\{Z\}$.}$$
\end{lemma}
\begin{proof} Clearly $D_P$ is homogeneous. From $P=\frak{d}_P D_P + R$, we get
$$ \Ric(P)\ =\  \frak{d}_P\Ric(D_P)+ \Ric(R),$$
with $\Ric(D_P)\in C\{Z\}$, $\Ric(D_P)\ne 0$, and 
$v(\Ric(R))=v(R) > v(P) = v(\frak{d}_P)$. This gives the desired result.
\end{proof}

\noindent
The case that $D_P\in C[Y](Y')^\N$ is important later in this chapter.
The next lemma concerns the even more special case $D_P\in C[Y]$, except that we use $Z$ rather than~$Y$ as the indeterminate, since the lemma will
be applied to differential polynomials obtained from Riccati 
transforms. 

\begin{lemma}\label{ricnew} Assume
$Q\in  K\{Z\}^{\neq}$ and $D_Q\in C[Z]$.
Then $D_{Q^\phi}=D_Q$ for $\phi\preceq 1$.
\end{lemma}
\begin{proof} We have $Q=\frak{d}_Q D_Q + R_Q$, with $Q \succ R_Q$. Then for
$\phi\preceq 1$ we have $Q^{\phi}=\frak{d}_Q D_Q + (R_Q)^\phi$, with
$(R_Q)^\phi\preceq R_Q\prec \frak{d}_Q$ by 
Lemma~\ref{compconjval, general}, so 
$D_{Q^\phi}=D_Q$.  
\end{proof}

\subsection*{From the dominant part to the Newton polynomial}
If $\phi_0\in \frak{M}$ is active in~$K$ and 
$\phi \preceq\phi_0$ is such that $w=\Pnu(P^{\phi})=\Pmu(P^{\phi_0})$, then
\begin{equation}\label{dominant monomial of P^phi}
\frak d_{P^{\phi}}\ =\ (\phi/\phi_0)^{w}\frak d_{P^{\phi_0}},
\end{equation}
by Corollary~\ref{compconjval, cor 3}(i), with $K^{\phi_0}$ in the role of $K$. Moreover: 

\begin{lemma}\label{dominant part ultimately constant}
Suppose $\phi\preceq 1$ and $\Pmu(P^\phi)=\Pnu(P)=w$. Then $\frak d_{P^\phi}=\phi^{w}\frak d_P$, and $D_P$ is isobaric of weight $w$ with $D_{P^\phi}=D_P$.
\end{lemma}
\begin{proof}
If $\phi=1$, then the lemma holds trivially, so assume $\phi \prec 1$.  
Then $\frak d_{P^\phi}=\phi^{w}\frak d_P$ by
Corollary~\ref{compconjval, cor 3}(ii) and
\eqref{dominant monomial of P^phi}.
Also, by Corollary~\ref{compconjval, cor 3}(iii),
\begin{align*}
D_{P^\phi}\ &=\ \sum_{v((P^\phi)_{[\bsigma]})=v(P^\phi)} \bar{\big((P^\phi)_{[\bsigma]}/\frak d_{P^\phi}\big)}\ Y^{[\bsigma]}\ 
 =\ \sum_{v(P_{[\bsigma]}) =v(P)} \bar{\phi^w P_{[\bsigma]}/\phi^w\frak d_{P}}\ Y^{[\bsigma]} \\
&=\   \sum_{v(P_{[\bsigma]}) =v(P)} \bar{P_{[\bsigma]}/\frak d_{P}}\ Y^{[\bsigma]}\  =\  D_P,
\end{align*}
and $D_P$ is isobaric of weight $w$.  
\end{proof}

\noindent
Recall from Section~\ref{evtbeh} the definition of the Newton weight $\nwt(P)$
of $P$:
$$ \nwt(P)\ =\ \dwt(P^{\phi})\ =\ \dwv(P^\phi),\quad \text{ eventually,}$$ and so Lemma~\ref{dominant part ultimately constant} yields a differential polynomial $N\in C\{Y\}$ such that
$$D_{P^\phi}\ =\ N,\quad \text{ eventually.}$$ 

\begin{definition}
The {\bf Newton polynomial}
of $P$ is the unique $N_P\in C\{Y\}$ such that 
eventually $D_{P^\phi}=N_P$. By convention, $N_Q:=0\in C\{Y\}$ for $Q=0\in K\{Y\}$.
\end{definition}

\noindent
Clearly $N_P$ is isobaric of weight $\nwt(P)$, and if $\dwt(P)=\nwt(P)$, then $D_P=N_P$. In particular, if $P\in K[Y,Y']$, then 
$N_P\in C[Y](Y')^{\N}$. 
Our $N_P$ depends on the monomial group $\frak{M}$ of $K$, but the Newton polynomial of $P$ obtained with another choice of $\frak{M}$ equals $cN_P$ for some $c\in C^\times$. If we want to stress the dependence of~$N_P$ on~$K$ equipped with its monomial group $\frak{M}$, we write $N^K_P$ for $N_P$. 

\begin{example} 
Suppose $P=Q(Y)\cdot (Y')^w$ where $Q\in K[Y]^{\ne}$ and $w\in\N$. Then
$P^\phi=\phi^w P$ for every $\phi$, hence $\nwt(P)=w$
and $N_P=D_Q\cdot (Y')^w$.
\end{example}

\begin{example} 
Suppose $D_P\in C[Y]$. Then $D_{P^\phi}=D_P$ for every $\phi\preceq 1$, by Lem\-ma~\ref{ricnew}, hence
$N_P=D_P$.
\end{example}

\noindent
The Newton polynomials $N_{P_{\times \fm}}$ of the
multiplicative conjugates $P_{\times \fm}$ play a role in detecting
zeros of $P$. To explain this, let $f\in K^\times$, and let 
$(c,\fm)$ be the unique pair with $c\in C^\times$ such that $f\sim c\fm$.
With these notations:

\begin{lemma}\label{newtondetection} Suppose $P(f)=0$.  Then $N_{P_{\times \fm}}(c)=0$, and thus
$\nwt(P_{\times \fm}) \ge 1$ or~$N_{P_{\times \fm}}$ is not homogeneous.
\end{lemma}
\begin{proof} We can reduce to the case that $v(f)=0$, $v(P)=0$
and $N_P=D_P$, so $\fm=1$ and $f=c+\epsilon$, $\epsilon\prec 1$, hence
$$0\ =\ P(f)\ =\ N_P(c+\epsilon) + R_P(f)\ =\ N_P(c)+g,\ \text{ with $g\prec 1$,}$$
so $N_P(c)=0$. If $\wt(N_P)=0$, then $N_P\in C[Y]$ and $N_P$ is not homogeneous.  
\end{proof} 

\noindent
Motivated by this lemma we define a monomial $\fm$ to be a 
{\bf starting monomial for $P$} \index{differential polynomial!starting monomial}\index{starting monomial} if $\nwt(P_{\times \fm}) \ge 1$ or 
$N_{P_{\times \fm}}$ is not homogeneous; equivalently, 
$N_{P_{\times \fm}}\notin CY^{\N}$. We call $\fm$ 
%a {\bf differential starting monomial for $P$} if 
%$\nwt(P_{\times \fm}) \ge 1$, and
an {\bf algebraic starting monomial for $P$} if $N_{P_{\times \fm}}$ is not 
homogeneous. \index{starting monomial!algebraic} Note: if $\fm$ is a starting monomial for $P$, then $\ndeg P_{\times \fm} \ge 1$. Also, 
$\fm$ is an algebraic starting monomial for $P$ iff $\fm/\fn$ is an algebraic starting monomial for $P_{\times \fn}$.
By Corollary~\ref{dn1cor,new,Newton}, $P$ has at most $\deg P-\val P$
{\em algebraic\/} starting monomials. But some $P$  
have infinitely many starting monomials:

\begin{exampleNumbered}\label{ex:infinitely many sm}
Let $K$ be the $H$-subfield $\R(\ex^{\R x}, \ell_0^{\R},\ell_1^{\R},\dots)$ of $\T$, with monomial group
$\fM=\bigcup_n \ex^{\R x}\ell_0^{\R}\cdots \ell_n^{\R}$. Then
$K$ is $\upo$-free, and for  
$P:=Y''Y-(Y')^2$ and every~$r\in \R$ we have $P(\ex^{rx})=0$, so $\ex^{rx}$ is a starting monomial for $P$.
\end{exampleNumbered}

\noindent
Call $f\in K^\times$ an {\bf approximate zero} of $P$ if $N_{P_{\times\fm}}(c)=0$, where
$(c,\fm)$ is the unique pair in $C^\times \times\fM$ with $f\sim c\fm$; the {\bf multiplicity} of~$f$ as an
approximate zero of $P$ is then by definition the multiplicity of $N_{P_{\times\fm}}$ at $c$ as defined just before Lemma~\ref{lem:multiplicity}. If $P(f)=0$, then $f$ is an 
approximate zero of $P$ by Lemma~\ref{newtondetection}.

\index{approximate!zero}
\index{multiplicity}

\medskip\noindent
In the next section we derive various useful properties of these Newton polynomials. % for example Corollary~\ref{npfinite} to the effect that $P$ has only finitely many 
We now continue with technicalities about dominant parts as needed later.

\subsection*{Decomposing $P^{\phi}$} 
The identities~\eqref{eq:Phi for isobaric P} below provide a useful decomposition of $P^{\phi}$ for isobaric $P$.
Accordingly, the asymptotic equivalence from Lemma~\ref{compconjval, 3}
will be improved in Lem\-ma~\ref{lem:properties of P<i>}.   
For $0\leq k\leq n$ we define $\epsilon^n_k(\phi)\in K$
by 
\begin{align*} \epsilon^n_0(\phi)\ &=\ \epsilon^n_n(\phi)\ =\ 0, \quad\epsilon_k^n(\phi)\ =\ 0 \text{ if $\phi'=0$,} \\
F^n_k(\phi)\ &=\ {n\brack k}\phi^k(\phi^\dagger)^{n-k}\big(1+\epsilon^n_k(\phi)\big),
\end{align*}
so $\epsilon^n_k(\phi)\prec 1$ if $\phi^\dagger\succ \phi$. Given $\btau=\tau_1\cdots \tau_d\geq\bsigma=\sigma_1\cdots \sigma_d$ we put 
\begin{align*}\epsilon_{\bsigma}^{\btau}(\phi)\ :&=\
-1+\prod_{i=1}^d \big(1+\epsilon^{\tau_i}_{\sigma_i}(\phi)\big),\ \text{  so}\\ 
 F_{\bsigma}^{\btau}(\phi)\ &=\ {\btau\brack\bsigma} \phi^{\|\bsigma\|} 
(\phi^\dagger)^{\|\btau\|-\|\bsigma\|}\big(1+\epsilon_{\bsigma}^{\btau}(\phi)\big),
\end{align*}
and if $\phi^\dagger\succ \phi$, then $\epsilon^{\btau}_{\bsigma}(\phi)\prec 1$.  By Lemma~\ref{conj-lemma},
$$(P^\phi)_{[\bsigma]}\ =\ \phi^{\|\bsigma\|}\sum_{\btau\geq\bsigma} {\btau\brack\bsigma} (\phi^\dagger)^{\|\btau\|-\|\bsigma\|} 
\big(1+\epsilon^{\btau}_{\bsigma}(\phi)\big)P_{[\btau]}.$$
For $i\in \N$ we define
$$P^{\phi,i}\ :=\ \sum_{\|\bsigma\|=i} \left(\sum_{\btau\geq\bsigma} {\btau\brack\bsigma}  \big(1+\epsilon^\btau_{\bsigma}(\phi)\big)P_{[\btau]} \right)\ Y^{[\bsigma]}\in K\{Y\},
$$
so $P^{\phi,i}$ is isobaric of weight $i$. Note also that 
$P^{\phi,0}=P_{[0]}$ is the isobaric part of $P$ of weight $0$.
If $P$ is isobaric of weight $w$, then $P^{\phi,w}=P$ and
\begin{equation}\label{eq:Phi for isobaric P} 
P^\phi\ =\ \sum_{i=0}^{w} \phi^i(\phi^\dagger)^{w-i} P^{\phi,i}, \quad
(P^{\phi})_{[i]}\ =\ \phi^i(\phi^\dagger)^{w-i}P^{\phi,i}\ \text{ for $i=0,\dots,w$.}
\end{equation}
Suppose $x\in K$ satisfies $x\succ 1$ and $x'=1$, and set $t:=1/x$.
Then $t=x^\dagger$ and
\begin{equation}\label{eq:Pti}
P^{t,i}\ =\  \sum_{\|\bsigma\|=i} \left(\sum_{\btau\geq\bsigma} {\btau\brack\bsigma} P_{[\btau]}  \right)\ Y^{[\bsigma]}.
\end{equation}
(To see this, note that Lemma~\ref{lem:Stirling} gives $\epsilon^n_k(t)=0$ for all $n$ and $k=0,\dots,n$.) The following lemma compares $P^{\phi,i}$ and $P^{t,i}$.

\begin{lemma}\label{lem:properties of P<i>} The $P^{\phi,i}$ have the following properties: \begin{enumerate}
\item[\textup{(i)}] if $\wt(P)=w>0$, then $P^{\phi,w-1}=P^{t,w-1}$;
\item[\textup{(ii)}] if $P$ is isobaric of weight $w>0$, then $P^{\phi,w-1}=-\nabla_1(P)$;
\item[\textup{(iii)}] if $\phi^\dagger\succ\phi$, then for all $i\in \N$,
$$v\big(P^{\phi,i}-P^{t,i}\big)\ \geq\ 
v(P)+\psi\big(\psi(v\phi)-v\phi\big)-\psi(v\phi)\ >\ v(P).$$
\end{enumerate}
\end{lemma}
\begin{proof} For (i), Example~\ref{ex:example for G} shows that for 
$\bsigma\leq\btau$ with $\|\btau\|\leq\|\bsigma\|+1$,
$$F^{\btau}_{\bsigma}(\phi)\ =\ 
{\btau\brack\bsigma}\phi^{\|\bsigma\|}
(\phi^\dagger)^{\|\btau\|-\|\bsigma\|}.$$
For (ii), use (i), the second identity in~\eqref{eq:Phi for isobaric P} for $\phi=t$ and $i=w-1$, and Lemma~\ref{lem:isobaric components of Pphi}
%Corollary~\ref{cor:nabla_phi,i} 
for $\phi=t$, taking into account
that $\nabla_{t,1}=\nabla_1$. For (iii), use the remark following 
Lemma~\ref{derivatives of equal weight 
and degree} and the proof of 
Lemma~\ref{derivatives of equal weight and degree, 2}. 
\end{proof}

\noindent
Even if $K$ does not contain an element $x\succ 1$ with $x'=1$, we
define $P^{t,i}\in K\{Y\}$ by \eqref{eq:Pti}, and then
Lemma~\ref{lem:properties of P<i>}
goes through. To see this, use that by Proposition~\ref{asintint} there is
an $x\succ 1$ with $x'=1$ in an immediate asymptotic extension of~$K$.

\subsection*{More on the case that $D_P\in C[Y](Y')^{\N}$ and $R_P \prec^{\flat} P$} We recall from the subsection on flattening in Section~\ref{Asymptotic-Fields-Basic-Facts} the convention on using~$\prec^{\flat}$, etcetera, to
denote, not only a certain binary relation on $K$, but also its extension to~$K\<Y\>$.

\begin{lemma}\label{cle} Assume $D_P\in C[Y](Y')^{\N}$, 
$R_P \prec^{\flat} P$, and $\phi \preceq 1$. Then 
$$  D_{P^\phi}\ =\ D_P\ =\ N_P, \qquad (R_P)^\phi\ =\ R_{P^\phi}, 
\qquad  R_{P^\phi}\ \prec^{\flat}\ P^{\phi}.$$ 
\end{lemma}
\begin{proof} From $P=\frak{d}_P D_P + R_P$ with $D_P\in C[Y](Y')^n$ we obtain
$$ P^\phi\ =\ \phi^n \frak{d}_P D_P + (R_P)^\phi \quad \text{ in $K^\phi\{Y\}$.}$$
From $0\le v\phi < 1+1$ we get $\phi\asymp^{\flat} 1$. By Lemma~\ref{compconjval, general} with $v^\flat$ in place of $v$, 
$$ (R_P)^\phi \asymp^\flat R_P \prec^\flat P \asymp^\flat \phi^n \frak{d}_P  \asymp^\flat
P^\phi,$$
so $D_{P^\phi}=D_P$. As this holds for all $\phi\preceq 1$, we get 
$D_P=N_P$. 
\end{proof}

\noindent
Here is a slight variant, with almost the same proof:

\begin{lemma}\label{cl} Assume $D_P\in C[Y](Y')^{\N}$,  
$\phi\preceq 1$, and $R_P\prec \phi^n P$ for all $n$. Then 
$$  D_{P^\phi}\ =\ D_P,\ \qquad (R_P)^\phi\ =\ R_{P^\phi}, 
\qquad  R_{P^\phi}\ \prec\ \phi^n P^{\phi} \text{ for all $n$.}$$
\end{lemma}
\begin{proof}  From $P=\frak{d}_P D_P + R_P$ with $D_P\in C[Y](Y')^m$ we obtain
$$ P^\phi\ =\ \phi^m \frak{d}_P D_P + (R_P)^\phi \quad \text{ in $K^\phi\{Y\}$.}$$
By Lemma~\ref{compconjval, general} we have 
$$(R_P)^\phi\ \preceq\ R_P\ \prec\ \phi^nP\ \text{ for all $n$,}$$ 
hence $(R_P)^\phi \prec \phi^m \frak{d}_P  \asymp P^\phi,$ so $D_{P^\phi}=D_P$. 
\end{proof}

\noindent
The next lemma and its corollary will only be needed 
in Section~\ref{sec:unravelers}.  For $i\in\N$ and $Q\in K\{Y\}$, let $\partial_i Q:=\frac{\partial Q}{\partial Y^{(i)}}$ denote the partial derivative of $Q$ with respect to~$Y^{(i)}$. With this notation, we have:

\begin{lemma}
Suppose $D_P\in C[Y](Y')^\N$ and $R_P\prec^\flat P$. Let $i\in\{0,1\}$ be such that 
$\partial_i D_P\neq 0$. Then $\partial_i P \asymp P$ and
$$
D_{\partial_i P}\ =\ \partial_i D_P\in C[Y](Y')^\N,\quad R_{\partial_i P}\ \prec^\flat\  \partial_i P.$$
\end{lemma}
\begin{proof}
From $\partial_i P= \fd_P \, \partial_i D_P  + \partial_i R_P$ and
$\partial_i R_P  \preceq R_P  \prec^{\flat} P$, we obtain 
$\fd_{\partial_i P}=\fd_P$, $D_{\partial_i P}=\partial_i D_P$, and
$R_{\partial_i P} = \partial_i R_P \prec^\flat \partial_i P$.
\end{proof}

\begin{cor}\label{cor:partials and Newton pol} 
Suppose $D_P\in C[Y](Y')^\N$ and $R_P\prec^\flat P$. Let $k,l\in\N$ be such that
$\frac{\partial^{k+l} D_P}{\partial Y^k \partial (Y')^l}\neq 0$. Then for
$Q:=\frac{\partial^{k+l} P}{\partial Y^k \partial (Y')^l}$ we have $Q\asymp P$ and
$$\frac{\partial^{k+l} D_P}{\partial Y^k \partial (Y')^l}\ =\ D_Q\ =\ N_Q, \qquad R_Q\ \prec^{\flat}\ Q.
$$
\end{cor}

\noindent
In the next lemma we use notation introduced in the preceding subsection.

\begin{lemma}\label{R_P^phi}
Suppose $\phi\preceq 1$ and $\Pmu(P^\phi)=\Pnu(P)$. Then 
$$
R_{P^\phi}\ =\  (R_P)^\phi + \sum_{i=0}^{w-1} \frak d_P \cdot \phi^i (\phi^\dagger)^{w-i}\cdot (D_{P})^{\phi,i} \qquad \text{where $w:=\Pmu(P)$.}$$
\end{lemma}
\begin{proof}
By \eqref{eq:Phi for isobaric P} we have 
$$(D_P)^\phi\ =\ \phi^w D_P + \sum_{i=0}^{w-1} \phi^i (\phi^\dagger)^{w-i}\cdot
(D_{P})^{\phi,i}.$$
Using \eqref{dominant monomial of P^phi} we obtain
\begin{align*}
P^\phi\  &=\  \frak d_P (D_P)^\phi + (R_P)^\phi  \\
       &=\ \frak d_{P^\phi} D_P + (R_P)^\phi + \sum_{i=0}^{w-1} \frak d_P\cdot\phi^i (\phi^\dagger)^{w-i}\cdot (D_{P})^{\phi,i}.
\end{align*}
This yields the displayed formula for $R_{P^\phi}$, since $D_{P^\phi}=D_P$ 
by Lemma~\ref{dominant part ultimately constant}.
\end{proof}

\noindent
Recall from Section~\ref{Asymptotic-Fields-Basic-Facts} that $\prec^\flat_\phi$ refers to the flattening $v^\flat_{\phi}$ of the valuation $v$ of $K^\phi$. 

\begin{lemma}\label{cleanliness}
Suppose $\phi\preceq 1$, $R_P \prec^\flat_\phi P$ and $\Pmu(P^\phi)=\Pnu(P)$.
Then 
$$ (R_{P^\phi})_{[w]}\  \prec^\flat_\phi\ P^\phi\qquad\text{for all $w\geq\Pmu(P)$.}$$
\end{lemma}
\begin{proof} Let $w\geq\Pmu(P)$. By the last lemma
$(R_{P^\phi})_{[w]}=\big( (R_P)^\phi \big)_{[w]}$. Hence
$$v\big((R_{P^\phi})_{[w]}\big)\ \geq\  wv\phi+v(R_P)\ \geq\ v(P^\phi)+v(R_P)-v(P)$$
by Lemma~\ref{compconjval, 5} and Corollary~\ref{compconjval, cor 3}. By assumption  $v^\flat_\phi(R_P)>v^\flat_\phi(P)$, so
\equationqed{v^\flat_\phi\big((R_{P^\phi})_{[w]}\big)\ \geq\  v^\flat_\phi(P^\phi)+v^\flat_\phi(R_P)-v^\flat_\phi(P)\ >\ v^\flat_\phi(P^\phi).}
\end{proof}

\subsection*{Behavior of $P(y)$}
If $D_P\in C[Y](Y')^{\N}$ and $R_P\prec^{\flat} P$, then we have a 
good description of $vP(y)$ in the region 
$1\nasymp y\asymp^\flat 1$, even in suitable
extensions of $K$: 

\begin{lemma}\label{lem:sign and valuation}
Suppose $D_P\in C[Y](Y')^{\N}$ and $R_P\prec^{\flat} P$. Then we have for every
$\d$-valued field extension $L$ of $H$-type of $K$ and all $y\in L$,
\begin{align*} 1\prec y\asymp^\flat 1\  &\Longrightarrow\ v\big(P(y)\big)\ =\ v(P) + \dd(P)\,vy + \dwt(P)\,\psi_L(vy),\\
 1\succ y\asymp^\flat 1\  &\Longrightarrow\ v\big(P(y)\big)\ =\ v(P) + \dv(P)\,vy + \dwt(P)\,\psi_L(vy).
\end{align*}
Moreover, if $K$ is equipped with an ordering making $K$ an $H$-field, then there are $\sigma,\tau\in\{-1,+1\}$ such that for all $H$-field extensions $L$ of $K$ and all $y\in L^>$,
\begin{align*}
 1\prec y\asymp^\flat 1\ &\Longrightarrow\ \sgn P(y) = \sigma, \\
 1\succ y\asymp^\flat 1\ &\Longrightarrow\ \sgn P(y) = \tau.
\end{align*}
\end{lemma}
\begin{proof}
After dividing $P$ by $\frak d_P$ we may assume that $v(P)=0$, so
$$P=D_P+R_P\qquad\text{where $v^\flat(R_P)>0$.}$$
We have $D_P=D(Y)\cdot (Y')^{w}$ where $D\in C[Y]$, $w=\wt(D_P)$. Let $L$ be
a $\d$-valued field extension of $H$-type of $K$. Now $1\in \Gamma^\flat\subseteq \Gamma_L^\flat$, so 
$\psi_L\big((\Gamma_L^\flat)^{\ne}\big)\subseteq \Gamma_L^\flat$. Let $y\in L$. 
By these facts about $L$,
if $y\preceq^\flat 1$, then $R_P(y)\prec^\flat 1$. Also 
$\dd P=\deg D + w$ and $\dv P = \val D + w$. If $y\succ 1$, then
$$v\big(D_P(y)\big)\ =\ v(D(y))+w\,v(y')\ =\ (\deg D)\,vy + w\,\big(vy+\psi_L(vy)\big),$$
and for $y\prec 1$ these equalities hold with $\val D$ instead of 
$\deg D$. 
Suppose now that $1\not\asymp y\asymp^\flat 1$. Then $vy\in (\Gamma_L^\flat)^{\ne}$, and 
so $\psi_L(vy)\in\Gamma_L^\flat$, hence 
$$D_P(y)\ \asymp^\flat\ 1, \qquad P(y)\ \simflat\ D_P(y).$$
Suppose now that $L$ is an $H$-field and $y>0$.
If $1\prec y\asymp^\flat 1$, then $y'>0$, so 
$\sgn P(y)=\sgn D(y)$, which equals the sign of the coefficient of the highest degree term of $D$.
If $1\succ y\asymp^\flat 1$, 
then $y'<0$, so $\sgn P(y)=(-1)^w\sgn D(y)$, and $\sgn D(y)$ equals the sign of the coefficient of the lowest degree term of $D$. 
\end{proof}

\begin{remark}
With the same assumptions as in Lemma~\ref{lem:sign and valuation}, its proof also shows that if
$\dd P >0$ and $L$ is a $\d$-valued field extension of $H$-type of 
$K$, then
$$ u,y\in L,\ u\asymp 1\prec y\asymp^\flat 1\ \Longrightarrow\ 
P(u)\prec P(y),$$
whereas if $\dd P=0$ and  $L$ is a $\d$-valued field extension of 
$H$-type of $K$, then 
$$y\in L,\ y\preceq ^\flat 1 \ \Longrightarrow\  P(y)\sim P(0) \sim P.$$
\end{remark}

\subsection*{Notes and comments}
In connection with Example~\ref{ex:infinitely many sm} we mention that for~$P$ of
order~$1$ there are only finitely many starting monomials;\marginpar{taken on faith for now} we omit the proof,  
since we do not use this fact later.
In Section~\ref{cordone} we show that if $K$ is $\upo$-free and~$P$ has degree~$1$, then 
$P$ has only finitely many starting monomials.

\section{Elementary Properties of the Newton Polynomial}
\label{basicnewton}

\noindent
Note that from the definition of $N_P$ we get
$$P^\phi\ =\ \frak d_{P^\phi} N_P + R_{P^\phi}, \qquad \text{ eventually.}$$
It is clear that $N_{P^\phi}=N_P$. We get $N_{PQ}=N_PN_Q$ for $Q\in K\{Y\}^{\ne}$ from the corresponding properties of dominant parts. In particular, $$N_{\frak{m} P}\ =\ N_P, \qquad N_{uP}\ =\ \bar{u}\,N_P\ \text{ for $u\in K$, $u\asymp 1$.}$$ 
Below we prove some other basic facts about Newton polynomials.

\subsection*{Lemmas on Newton polynomials} We begin with an easy consequence of the definitions and Corollary~\ref{compconjval, cor 3}.  

\begin{lemma}\label{DNNWT} Let $w=\nwt(P)$. Then 
$$ D_{P^\phi}=N_P \text{ for all $\phi\preceq 1$}\ \Longleftrightarrow\  D_P=N_P\ \Longleftrightarrow\ \Pnu(P)=w.$$
\end{lemma}

\begin{lemma}\label{newred} Let $P=Q+R$ where $Q,R\in K\{Y\}$ and
$R\prec^{\flat} P$. Then $N_P=N_Q$. If also $D_P=N_P$, then 
$D_{Q^\phi}=N_P$ for all $\phi\preceq 1$.
\end{lemma}
\begin{proof} Note that $Q\ne 0$ and $v(P)=v(Q)$. Let $\phi\preceq 1$. Then 
$0\le v\phi < 1+1$, so
$$ v(Q^\phi)\ \le\ v(Q)+ \wt(Q)v\phi\ <\  v(R)\ \le\ v(R^{\phi}). $$
Since $P^{\phi}=Q^{\phi}+ R^{\phi}$, this gives $v(P^\phi)=v(Q^\phi) < v(R^{\phi})$,
and so $D_{P^\phi}=D_{Q^\phi}$. This holds for all $\phi\preceq 1$, so
$N_P=N_Q$. If in addition $D_P=N_P$, then by Lemma~\ref{DNNWT} we obtain that
$D_{Q^\phi}=N_P$ for all $\phi\preceq 1$. 
\end{proof} 

\begin{lemma}\label{lem:Newton poly, precfn}
Suppose $\fm\flatter \fn\succ^\flat 1$ and $P=Q+R$ where $R\prec_{\fn} P$.
Then 
$$N_{P_{\times\fm}}\ =\ N_{Q_{\times\fm}}.$$
\end{lemma}
\begin{proof}
For $\phi\preceq 1$ we have $R^\phi \asymp_{\fn} R \prec_{\fn} Q \asymp_{\fn} Q^\phi$,
by Corollary~\ref{cor:Pphi flat equivalence}.
Hence replacing $K$, $P$, $Q$, $R$ by $K^\phi$, $P^\phi$, $Q^\phi$, $R^\phi$, respectively, for suitable $\phi\preceq 1$, we
arrange that $D_{P_{\times\fm}}=N_{P_{\times\fm}}$ and $D_{Q_{\times\fm}}=N_{Q_{\times\fm}}$. 
Then by Corollary~\ref{cor:Pg flat equivalence, 3}:
$$R_{\times\fm}\ \asymp_{\fn}\ R
\ \prec_{\fn}\ Q \  \asymp_{\fn}\ 
Q_{\times\fm},$$
so  $R_{\times\fm}\prec Q_{\times\fm}$ and hence $N_{P_{\times\fm}}=D_{P_{\times\fm}}=D_{Q_{\times\fm}}=N_{Q_{\times\fm}}$.
\end{proof}

\begin{cor}\label{cor:Newton poly, precfn}
Suppose $\fn\succ 1$ and $\ndeg P=\ndeg P_{\times\fn}=d$. Set $Q:=P_{\leq d}$.
Then for all $\fm\flatter \fn$ and all $g\preceq 1$ in $K$ we have
$$N_{P_{+g,\times\fm}}\ =\ N_{Q_{+g,\times\fm}}.$$
\end{cor}
\begin{proof}
After replacing $K$, $P$ by $K^\phi$, $P^\phi$, respectively, for suitable $\phi$, we may assume that
$D_P=N_P$, $D_{P_{\times\fn}}=N_{P_{\times\fn}}$, and $\fn\succ^\flat 1$. Let $R:=P-Q=P_{>d}$. Then by Corollary~\ref{cor:Pg flat equivalence, 2}
we have $R\preceq_\fn \fn^{-1} P\prec_\fn P$. Thus given $g\in K^{\preceq 1}$, 
Lemma~\ref{lem:Newton poly, precfn} applies to $P_{+g}$, $Q_{+g}$, $R_{+g}$ in place of $P$, $Q$, $R$, respectively.
\end{proof}

\noindent
Recall that $P_{|i|'}$ is the subhomogeneous part
of $P$ of subdegree $i$ (see Section~\ref{Decompositions of Differential Polynomials}). By Corollary~\ref{corconj-lemma}
we have  $(P^\phi)_{|i|'}=(P_{|i|'})^\phi$, and $P^\phi_{|i|'}$
denotes either of these without ambiguity.
When studying $N_P$, the following lemma sometimes allows us to 
reduce to the case where $P$ is homogeneous or subhomogeneous. 

\begin{lemma}\label{lem:NP components}
$$N_P\ =\ \sum_i N_{P_{|i|'}}\ =\ \sum_i N_{P_i},$$
where the first sum ranges over all $i\in\N$ such that $P^\phi\asymp P^\phi_{|i|'}$ eventually, and the second sum ranges over all $i\in\N$ such that $P^\phi\asymp P^\phi_{i}$ eventually.
\end{lemma}
\begin{proof}
We will only prove $N_P = \sum_i N_{P_{|i|'}}$. (To show $N_P=\sum_i N_{P_i}$ one argues in an analogous way.) Below $i$ ranges over elements of $\N$ with 
$P_{|i|'}\ne 0$, and likewise with $j$.
First, after replacing $P$ by $P^\phi$ for suitable $\phi\preceq 1$, we may assume that for all $\phi\preceq 1$ we have
$v(P_{|i|'}^\phi)=v(P_{|i|'})+\nwt(P_{|i|'})v\phi$ and $N_{P_{|i|'}}=D_{P_{|i|'}^\phi}$.
Therefore, for all $i$,~$j$ with $\nwt(P_{|i|'}) \neq \nwt(P_{|j|'})$, 
either $P_{|i|'}^\phi \prec P_{|j|'}^\phi$ eventually, or $P_{|i|'}^\phi\succ P_{|j|'}^\phi$ eventually. Hence, after replacing $P$ by $P^\phi$ for suitable $\phi\preceq 1$, we may assume that for all $i$,~$j$  and all $\phi\preceq 1$ we have 
$P_{|i|'}\prec P_{|j|'}$ iff
$P_{|i|'}^\phi \prec P_{|j|'}^\phi$. So for all~$i$ and $\phi\preceq 1$, $P\asymp P_{|i|'}$ iff $P^\phi\asymp P_{|i|'}^\phi$, hence for all $i$, $P\asymp P_{|i|'}$ iff $P^\phi\asymp P_{|i|'}^\phi$ eventually.
Thus for each $\phi\preceq 1$,
$$D_{P^\phi}\ =\  \sum_{P_{|i|'}^\phi\asymp P^\phi} D_{P_{|i|'}^\phi}\ 
=\ \sum_{P\asymp P_{|i|'}} N_{P_{|i|'}},$$
and this yields the claim.
\end{proof}

\begin{lemma}\label{npq} Let $P,Q\in K\{Y\}^{\ne}$ be homogeneous of 
different degrees. Then $N_{P+Q}\in \{N_P, N_Q, N_P + N_Q\}$. 
Also, $N_{(P+Q)_{\times \fm}}$ is homogeneous for 
every mo\-no\-mi\-al~$\fm$, with at most one exception.
\end{lemma}
\begin{proof} Take $\phi$ such that 
$$N_P\ =\ D_{P^\phi},\quad   N_Q\ =\ D_{Q^\phi},\quad
N_{P+Q}\ =\ D_{P^\phi + Q^\phi}.$$ If
$v(P^\phi) < v(Q^\phi)$, then $N_{P+Q}= D_{P^\phi}=N_P$, and
if $v(P^\phi) > v(Q^\phi)$, then $N_{P+Q}= D_{Q^\phi}=N_Q$. If
 $v(P^\phi) = v(Q^\phi)$, then  
$$N_{P+Q}\ =\ D_{P^\phi + Q^\phi}\ =\  D_{P^\phi}+ D_{Q^\phi}\ =\ N_P+N_Q,$$
since $P^\phi$ and $Q^\phi$ are homogeneous of different degrees.

For the second claim of the lemma, assume $\deg P= d < \deg Q=e$.
Suppose~$N_{(P+Q)_{\times \fm}}$ 
is not homogeneous. It suffices to show that then
$N_{(P+Q)_{\times g}}= N_{P_{\times g}}$ for all nonzero $g \prec \fm$ in $K$.
Towards proving this, we can arrange $\fm=1$, so by the argument above we
have $v(P^\phi)=v(Q^\phi)$ for $\phi$ as above. Let $g\in K^\times$, $g\prec 1$ 
and set 
$\gamma:= vg$, so $\gamma >0$. Take $\phi$ as above such that in addition
$$ N_{P_{\times g}}\ =\ D_{P_{\times g}^\phi},\qquad  N_{Q_{\times g}}\ =\ 
D_{Q_{\times g}^\phi},\qquad
N_{(P+Q)_{\times g}}\ =\ D_{P_{\times g}^\phi + Q_{\times g}^\phi}.$$ 
By Corollary~\ref{vP-Lemma} we have
$$v(P_{\times g}^\phi)\ =\ v(P^\phi) + d\gamma + o(\gamma)\  <\
v(Q_{\times g}^\phi)\ =\  v(Q^\phi) + e\gamma + o(\gamma),$$
so $N_{(P+Q)_{\times g}}= N_{P_{\times g}}$ by the proof of the first claim
of the lemma.  
\end{proof}

\noindent
Let $J$ be the finite nonempty set of $j\in \N$ such that $P_j\ne 0$; then $P=\sum_{j\in J} P_j$.

\begin{cor}\label{npsum} There is a unique set $I\subseteq J$ such that
$N_P=\sum_{i\in I} N_{P_i}$. This set~$I$ is determined by the condition that
for all $i\in J$,
$$i\in I \ \Longleftrightarrow\  v(P_i^\phi)\le v(P_j^\phi), \text{ eventually, for each $j\in J$.}$$
\end{cor}
\begin{proof} Use that for all $i,j\in J$,
either $v(P_i^\phi)< v(P_j^\phi)$, eventually, or 
$v(P_i^\phi)=v(P_j^\phi)$, eventually, or  $v(P_i^\phi)>v(P_j^\phi)$, eventually.
\end{proof}

\begin{cor}\label{npfinite} For all but finitely many $\fm$ there is 
$i\in J$
such that $$N_{P_{\times \fm}}\ =\ N_{P_{i, \times \fm}}.$$
\end{cor}
\begin{proof} Let $\fm$ be such that for all distinct $i,j\in J$ the
Newton polynomial $$N_{(P_i + P_j)_{\times \fm}}$$ is 
homogeneous. Then there is a (necessarily unique) $i\in J$ such that
$v(P_{i, \times \fm}^\phi)< v(P_{j, \times \fm}^\phi)$, eventually, for every $j\in J\setminus \{i\}$.
For this $i$ we have $N_{P_{\times \fm}}=N_{P_{i, \times \fm}}$.
\end{proof}

\subsection*{Cleanness}
In Section~\ref{sec:The Dominant Part of a Differential Polynomial}, having both $D_P\in C[Y](Y')^\N$ and $R_P\prec^{\flat} P$ turned out to be very strong. The next lemma shows that for divisible $\Gamma$ 
the eventual form of the first condition implies the eventual form of the
second condition.

\begin{lemma}\label{divclean}
Suppose $\Gamma$ is divisible and $N_P\in C[Y](Y')^\N$. Then
$$ R_{P^{\phi}}\ \prec^{\flat}_{\phi}\  P^\phi,\quad \text{ eventually.}$$ 
\end{lemma}
\begin{proof}
Set $w:=\nwt(P)$.
After replacing $P$ by $P^\phi$ for suitable $\phi\preceq 1$, we may assume $w=\dwv(P) =\dwv(P^\phi)=\dwt(P)$ for all $\phi\preceq 1$, and thus for all
such $\phi$,
$$D_P\ =\ D_{P^\phi}\ =\ N_P, \qquad v(P^\phi)\ =\ v(P)+w\,v\phi.$$  
In particular, $(D_P)^{\phi,i}=0$ for $i=0,\dots,w-1$ for $\phi\preceq 1$, so $(R_P)^\phi=R_{P^\phi}$ for $\phi\preceq 1$ by Lemma~\ref{R_P^phi}. If $R_P=0$, this gives $R_{P^\phi}=0$ for $\phi\preceq 1$. So assume $R_P\ne 0$. Replacing $P$ by $P^\phi$ for suitable $\phi\preceq 1$, we arrange in addition:  
$$v(R_{P^\phi})\ =\ v(R_P)+\nwt(R_P)v\phi \quad \text{ for all $\phi\preceq 1$.}$$ 
We need to show that $\psi^\phi\big(v(R_{P^\phi})-v(P^\phi)\big)\leq 0$, eventually. For this, we distinguish three cases. Suppose first that $\nwt(R_P)>w$. Then
we have for $\phi\prec 1$,
$$0\ <\ v\phi\ \leq\ v(R_P)-v(P)+\big(\!\nwt(R_P)-w\big)v\phi\ =\ v(R_{P^\phi})-v(P^\phi).$$
Hence, if $v\phi\geq 1$, then
$$0\ \geq\ 1-v\phi\ =\ \psi^\phi(v\phi)\ \geq\ \psi^\phi\big(v(R_{P^\phi})-v(P^\phi)\big).$$
Now assume $\nwt(R_P)< w$. Let $\alpha:=\frac{v(R_P)-v(P)}{w-\nwt(R_P)}\in
\Gamma^{>}$, and take $\beta\in\Gamma^{\ne}$ such that $\beta+\psi(\beta)=\alpha$. If $v\phi\geq\psi(\beta)$, then by Lemma~\ref{lem:phipsi},
$$\psi^\phi\big(v(R_{P^\phi})-v(P^\phi)\big)\ =\ \psi(\alpha-v\phi)-v\phi\ \leq\  0.$$
Finally, if $\nwt(R_P)=w$ and $v\phi\geq \psi\big(v(R_P)-v(P)\big)$, then clearly 
\equationqed{\psi^\phi\big(v(R_{P^\phi})-v(P^\phi)\big)\ \leq\  0.}
\end{proof}

\begin{cor}\label{ord1clean} If $\Gamma$ is divisible and $P\in K[Y,Y']$,
then 
$$N_P\in C[Y](Y')^{\N}, \text{ and eventually $R_{P^{\phi}}\ 
\prec^{\flat}_{\phi}\  P^{\phi}$.}$$
\end{cor}

\noindent
Let $L$ be a $\d$-valued field extension of $H$-type of $K$ with asymptotic integration, furnished with a monomial group $\mathfrak M_L\supseteq \mathfrak M$. Note that then $L$ has small derivation, and each $\phi$ is active in $L$. If in addition $\Gamma^>$ is coinitial in $\Gamma_L^>$, then
$N_P^L=N_P$. 
This fact yields a variant of Lemma~\ref{divclean}:

\begin{lemma}\label{uplfreeclean} 
If $K$ has rational asymptotic integration and $N_P\in C[Y](Y')^\N$,
then $ R_{P^{\phi}}\ \prec^{\flat}_{\phi}\  P^\phi$, eventually.
\end{lemma}
\begin{proof} The algebraic closure $L$ of $K$ has by  
Lemma~\ref{exmongr} a monomial group $\mathfrak{M}_L\supseteq \mathfrak M$,
and $\Gamma^>$ is coinitial in $\Gamma_{L}^>=(\Q\Gamma)^>$. Thus if $L$ has
asymptotic integration, then
$N_P^{L}=N_P\in C[Y](Y')^\N$, and Lemma~\ref{divclean} 
applies to $L$ instead of $K$.
\end{proof}

\noindent
In Section~\ref{The Shape of the Newton Polynomial} we shall prove: $K$ is $\upo$-free $\ \Longleftrightarrow\ $ $N_P\in C[Y](Y')^\N$ for all $P$.

\begin{cor}\label{cor:sign and valuation, 1}
Assume $N_P\in C[Y](Y')^{\N}$ and eventually $R_{P^{\phi}}\ \prec^{\flat}_{\phi}\  P^\phi$. Then there exists $\phi$ such that for every $\d$-valued field extension $L$ of $H$-type of $K$,
\begin{align*} 
y\in L,\ 1\prec y\asymp^\flat_\phi 1\  &\Longrightarrow\ v\big(P(y)\big)\ =\ v^{\ev}(P) + \ndeg(P)\,vy + \nwt(P)\,\psi_L(vy),\\
y\in L,\ 1\succ y\asymp^\flat_\phi 1\  &\Longrightarrow\ v\big(P(y)\big)\ =\ v^{\ev}(P) + \nval(P)\,vy + \nwt(P)\,\psi_L(vy).
\end{align*}
If $K$ is equipped with an ordering making $K$ an $H$-field, then there are
$\phi$ and
$\sigma,\tau\in\{-1,+1\}$ such that for all $H$-field extensions $L$ of $K$
and all $y\in L^>$,
\begin{align*}
 1\prec y\asymp^\flat_\phi 1\ &\Longrightarrow\ \sgn P(y) = \sigma, \\
 1\succ y\asymp^\flat_\phi 1\ &\Longrightarrow\ \sgn P(y) = \tau.
\end{align*}
\end{cor}
\begin{proof}
This follows from Lemma~\ref{lem:sign and valuation} and Corollary~\ref{eventualequality}, since for $y\nasymp 1$ in any
$\d$-valued field extension $L$ of $H$-type of $K$ we have
\equationqed{v(P^\phi) + \dwt(P^\phi)\psi_L^\phi(vy)\ =\ v^{\ev}(P) + \nwt(P)\,\psi_L(vy),\ \text{ eventually.}}
\end{proof}

\begin{definition}
We say that $K$ is {\bf clean} if for every $P$ we have
$$N_P\in C[Y](Y')^{\N},   \text{ and eventually $R_{P^{\phi}}\ \prec^{\flat}_{\phi}\  P^\phi$.}$$
\end{definition}

\noindent
Note that if $K$ is clean, then so is each compositional conjugate $K^\phi$.

\index{clean}

\subsection*{Behavior under additive and multiplicative conjugation}

\begin{lemma}\label{dn2} Let $c\in C$ and $\epsilon\in K$, $\epsilon\prec 1$. 
Then 
$$N_{P_{+c}}=(N_P)_{+c}, \qquad N_{P_{+\epsilon}}=N_P.$$
\end{lemma}
\begin{proof} Eventually
$N_{P_{+c}}=D_{(P_{+c})^\phi}$. Since  
$(P_{+c})^\phi= (P^\phi)_{+c}$, we have 
$N_{P_{+c}}=D_{(P^\phi)_{+c}}$, eventually. By Lemma~\ref{dn1},
$D_{(P^\phi)_{+c}}=\big(D_{P^\phi}\big)_{+c}$, and eventually we have
$\big(D_{P^\phi}\big)_{+c}=(N_P)_{+c}$. The other displayed item follows likewise from part~(iii) of Lemma~\ref{dn1,new}. 
\end{proof}

\begin{cor}\label{dn2cor}
Let $f,g,h\in K$ and $f-g\prec h$. Then $N_{P_{+f,\times h}}=N_{P_{+g,\times h}}$.
\end{cor}
\begin{proof} Use that $P_{+g,\times h}=P_{+f,\times h,+\epsilon}$ for $\epsilon:= \frac{g-f}{h}\prec 1$.
\end{proof}

\noindent
For use in Sections~\ref{eveqth} and~\ref{sec:unravelers} we need: 

\begin{lemma}\label{dn2a} 
Suppose that
$D_P\in C[Y](Y')^{\N}$, and suppose that $\gamma>0$ is such that $
v(R_P) > v(P)+ m\gamma + n\gamma'$ for all $m,n$.
Then for $g\in K$ with $vg=\gamma$ we have
$$ N_{P_{\times g}} \in C^\times\cdot  Y^\mu, \quad 
\mu:= \val D_P.$$ 
\end{lemma}
\begin{proof} We have $D_P=D(Y)\cdot (Y')^j$ where 
$$D\in C[Y],\qquad  D\ =\ cY^i + \text{terms of higher degree, $c\in C^\times$.}$$ 
We can arrange $v(P)=0$, so $P=D_P+R$, $R= R_P$. Let $g\in K$, $vg=\gamma$. Then 
$$P_{\times g}\ =\ D(g Y)\cdot g^j 
(g^\dagger Y + Y')^j + R_{\times g}.$$
Now $D(g Y)=g^icY^i(1+E)$ with $E\in YK[Y]$, $v(E)\ge \gamma$. 
Hence
$$ P_{\times g}\ =\ 
g^{i+j}cY^i(1+E)(g^\dagger Y + Y')^j +R_{\times g}.$$
Let $\phi$ be such that $v\phi > \gamma^\dagger= v(g^\dagger)$. Then, 
in view of $E^\phi=E$,
\begin{align*} P_{\times g}^\phi\ &=\ g^{i+j}cY^i(1+E)(g^\dagger Y + \phi Y')^j +R^\phi_{\times g}, \\
(g^\dagger Y + \phi Y')^j\ &=\ (g^\dagger Y)^j + F, \quad v(F) > j\gamma^\dagger, \text{ so}\\
P^\phi_{\times g}\ &=\ g^{i+j}(g^\dagger)^jcY^{i+j} + G, \quad v(G) > (i+j)\gamma+j\gamma^\dagger=i\gamma + j\gamma'.
\end{align*} 
Since $i+j=\val D_P$, this gives the desired result.
\end{proof}  

\begin{cor}\label{dn2a, cor}
Suppose that $D_P\in C[Y](Y')^\N$ and $R_P\prec^\flat P$.
Let $f,g\in K$ be such that $f\prec^\flat 1$ and 
$g\prec 1$, $g\asymp^\flat 1$. Then  
$$N_{P_{\times g}}, N_{P_{+f,\times g}}\in C^\times \cdot Y^\mu, \quad\text{where $\mu:=\val D_P$.}$$
\end{cor}
\begin{proof}
From $g\prec 1$ and $g\asymp^\flat 1$ we get $N_{P_{\times g}}\in C^\times \cdot Y^\mu$ by Lemma~\ref{dn2a}.
Lemma~\ref{v-under-conjugation}(i) and $f\prec^\flat 1$ yield
$(D_P)_{+f}\simflat D_P$ and $(R_P)_{+f}\asymp^\flat R_P$, and
so in view of $P_{+f}=\mathfrak d_P (D_P)_{+f}+(R_P)_{+f}$ we have
$\mathfrak d_{P_{+f}}=\mathfrak d_P$, $D_{P_{+f}}=D_P\in C[Y](Y')^\N$, and $R_{P_{+f}}\prec^\flat P_{+f}$.
Hence $N_{P_{+f,\times g}}\in C^\times \cdot Y^\mu$ by Lemma~\ref{dn2a}.
\end{proof}

\section{The Shape of the Newton Polynomial}
\label{The Shape of the Newton Polynomial}

\noindent
In this section we combine the material from the previous two
sections with results from Chapters~\ref{evq} and~\ref{ch:triangular automorphisms}.

\subsection*{Statement of results}
Recall from Section~\ref{sec:Derivations on polynomials rings} that $\nabla_1$ and $\nabla_2$ are the $1,2$-diagonals of the triangular (Stirling) derivation $\nabla$  of the polynomial $K$-algebra $K\{Y\}=K[Y,Y',Y'',\dots]$. 
By Co\-rol\-la\-ry~\ref{cor:companion of Stirling} and Lemma~\ref{cor:P vanish} we have
$$N_P\in C[Y](Y')^\N\ \Longleftrightarrow\ \nabla_1(N_P)=\nabla_2(N_P)=0.$$
The main goal of this section is to prove the following two results:

\begin{theorem}\label{thm:nabla1}
If $K$ is $\upl$-free, then $\nabla_1(N_P)=0$. Conversely, if $\nabla_1(N_Q)=0$ for each 
homogeneous  $Q\in K\{Y\}^{\neq}$ of degree $1$, then $K$ is $\upl$-free.
\end{theorem}

\begin{theorem}\label{thm:nabla2}
If $K$ is $\upo$-free, then $\nabla_1(N_P)=\nabla_2(N_P)=0$. Conversely, if $\nabla_1(N_Q)=\nabla_2(N_Q)=0$ for each 
homogeneous  $Q\in K\{Y\}^{\neq}$ of degree~$2$, then $K$ is $\upo$-free.
\end{theorem}

\noindent
Before we get to the proofs, let us first deduce some consequences.
Theorem~\ref{thm:nabla2} and Lemma~\ref{uplfreeclean} immediately yield a characterization of cleanness showing that~$K$ being clean does not depend on the choice of monomial group $\mathfrak{M}$:

\begin{cor}\label{cor:nabla2}
$K$ is clean iff $K$ is $\upo$-free.
\end{cor}

\noindent
The previous corollary in conjunction with Corollary~\ref{recam} gives:

\begin{cor}\label{cor:nabla3}
If $K$ is a union of asymptotic subfields, each with a smallest comparability class, then $K$ is clean.
\end{cor}

\noindent
Another consequence: if $K$ is spherically complete, then $K$ is not clean.

\subsection*{The $\upl$-free case}
We precede the proof of Theorem~\ref{thm:nabla1} with a lemma.
We identify the ordered group~$\Z$ with the ordered subgroup $\Z\cdot 1$ of $\Gamma$ via 
$k\mapsto k\cdot 1$. 

\begin{lemma}\label{W(P) = w}
Suppose $\dwv(P^\phi)=\dwt(P)$ and $v\phi > 1$. Then 
$D_{P^\phi}=D_{Q^\phi}$, where $Q:=P_{[0]}+\cdots+P_{[w]}$, $w:=\dwt(P)$. 
\end{lemma}
\begin{proof}
Replacing $P$ by $\frak d_P^{-1}P$, we may assume $\frak d_P=1$,
so $\frak d_{P^\phi}=\phi^w$. From the identity $(P-Q)^\phi=\sum_{i>w}(P_{[i]})^\phi$ and
\eqref{eq:Phi for isobaric P} we obtain
$$(P-Q)^\phi\ =\ 
\left(\sum_{i>w\geq j} \phi^j(\phi^\dagger)^{i-j}(P_{[i]})^{\phi,j} \right) + 
\left(\sum_{i\geq j>w} \phi^j(\phi^\dagger)^{i-j}(P_{[i]})^{\phi,j} \right).$$
It is clear that $v(P_{[i]})^{\phi,j}\geq v(P_{[i]})\geq v(P) = 0$ for all $i$,~$j$. Moreover, 
if $i>w\geq j$ then $i-j>(w-j)v\phi$ since $v\phi=1+o(1)$ by Corollary~\ref{1+o(1)}.  
Hence each term $\phi^j(\phi^\dagger)^{i-j}(P_{[i]})^{\phi,j}$ 
($i>w\geq j$) in the first sum has valuation 
$$jv\phi+i-j+v\big(P_{[i]})^{\phi,j}\ \geq\  jv\phi+i-j>wv\phi.$$
Clearly each term $\phi^j(\phi^\dagger)^{i-j}(P_{[i]})^{\phi,j}$ 
($i\geq j>w$) in the
second sum has valuation~${>wv\phi}$.
Hence $v(P-Q)^\phi>wv\phi=vP^\phi$ and thus
$D_{P^\phi}=D_{Q^\phi}$.
\end{proof}

\noindent
Let us also recall the transformation formulas deduced in Corollary~\ref{cor:nabla_phi,i}: with $w=\wt(P)$, $\lambda=-\phi^\dagger$ and $\omega=-(2\lambda'+\lambda^2)$,   
\begin{equation}\label{eq:trans fms}
\left\{\begin{aligned}
(P^\phi)_{[w]} &= \phi^w P_{[w]}, \\
(P^\phi)_{[w-1]} &= \phi^{w-1}\big[ P_{[w-1]} + \lambda \nabla_1(P_{[w]}) \big], \\
(P^\phi)_{[w-2]} &= \phi^{w-2}\big[ P_{[w-2]} + \lambda \nabla_1(P_{[w-1]}) +  \textstyle \left(\omega\nabla_2+\frac{1}{2}\lambda^2\nabla_1^2\right)(P_{[w]})\big].  \hskip-2em
\end{aligned}
\right.\ 
\end{equation}

\noindent
The next proposition and its corollary yield Theorem~\ref{thm:nabla1}:

\begin{prop}\label{prop:nabla1}
Suppose $K$ is $\upl$-free.
Then $\nabla_1(N_P)=0$.
\end{prop}
\begin{proof} 
Let $w:=\nwt(P)$. Then $N_P\in C\{Y\}$ is isobaric of weight $w$; if $w=0$, then $N_P\in C[Y]$ and hence $\nabla_1(N_P)=\nabla_2(N_P)=0$, so we can assume that~${w>0}$.
Replacing~$K$ and $P$ by $K^\phi$ and $P^\phi$ for suitable $\phi\preceq 1$, we arrange~${D_P=N_P}$.
(Here we use the invariance of $\upl$-freeness under compositional conjugation of~$K$.)
Then Lemma~\ref{DNNWT} gives $D_{P^\phi}=N_P$ and $\dwv(P^\phi)=\dwt(P^\phi)=w$ for all~${\phi\preceq 1}$. 
By Lemma~\ref{W(P) = w} we can further reduce to the case that $\wt(P)=w$, and by Lem\-ma~\ref{cleanliness}
we arrange $(R_P)_{[w]} \prec^{\flat} P$. From $P=\frak{d}_PN_P + R_P$
we get $P_{[w]}=\frak{d}_PN_P +(R_P)_{[w]}$. Multiplying 
$P$ by $\frak{d}_P^{-1}$, we obtain in addition
$P\asymp 1$ and $P_{[w]}=N_P + R$ with $R\in K\{Y\}$, $R\prec^\flat 1$.

Suppose towards a contradiction that $\nabla_1 N_P\neq 0$. Then
$$\nabla_1 P_{[w]}\ =\ \nabla_1 N_P + \nabla_1R,  \qquad
\nabla_1 N_P\in C\{Y\}^{\ne},\quad  \nabla_1R\prec^{\flat} 1.$$ 
Recall the pc-sequence $(\upl_\rho)$ of width $\{\gamma\in\Gamma_\infty:\gamma>\Psi\}$ introduced in Section~\ref{sec:special cuts}, and
take a pseudolimit $\upl$ of $(\upl_\rho)$ in some immediate asymptotic field extension of~$K$ (Theorem~\ref{thm:immediate}).
Since $(\upl_\rho)$ has no pseudolimit in~$K$, there is no $a\in K$ with ${v(a+\upl)>\Psi}$, which gives $\phi\preceq 1$ with $P_{[w-1]}+\upl\nabla_1 P_{[w]}\succ\phi$. With $\lambda=-\phi^\dagger$ we have
$\upl-\lambda\prec\phi$ by Lemma~\ref{lem:gap}, so
$$P_{[w-1]} + \lambda \nabla_1 P_{[w]}\ \sim\  P_{[w-1]} + \upl\nabla_1 P_{[w]}\ 
\succ\ \phi.$$
By Proposition~\ref{compconjval, prop} we have $P^{\phi} \asymp \phi^w$, but the second identity
of~\eqref{eq:trans fms} gives
\begin{align*}
(P^{\phi})_{[w-1]}\  &=\	\phi^{w-1}\big(P_{[w-1]}+\lambda\nabla_1 P_{[w]}\big) \\
	&\sim\ 	\phi^{w-1}\big(P_{[w-1]}+\upl\nabla_1 P_{[w]}\big)\ \succ\ \phi^w,
\end{align*}
a contradiction.  
\end{proof}

\begin{remark} 
For later use we record a variant of Proposition~\ref{prop:nabla1}:

\medskip\noindent
\textit{Let $K$ be an immediate extension of its valued differential subfield 
$E$, and assume~$E$ is $\upl$-free, $\frak{M}\subseteq E^\times$, and $P\in E\{Y\}^{\neq}$. Then $\nabla_1 N_P=0$.}

\medskip\noindent
We did not assume here that $K$ is $\upl$-free, and while $C_E\subseteq C$ and 
$N_P\in C\{Y\}$, it does not follow from
$P\in E\{Y\}$ that $N_P\in C_E\{Y\}$. 
Nevertheless, the proof of 
this variant is the same as that of Proposition~\ref{prop:nabla1}, apart 
from minor changes.
\end{remark}

\noindent
The above remark and Corollary~\ref{cor:companion of Stirling, variant}
give a result to be used in Section~\ref{fcuplfr}:
%Section~\ref{consupofree}:

\begin{cor}\label{cor:nabla1, order 2}
Let $K$ be an immediate extension of its valued differential subfield 
$E$, and assume $E$ is $\upl$-free, $\frak{M}\subseteq E^\times$, and $P\in E\{Y\}^{\neq}$ has order at most~$2$. Then $N_P\in C[Y](Y')^\N$.
\end{cor}
 
\noindent
Next we relate $\upl$-freeness to properties of homogeneous 
$P$ of degree $1$: 

\begin{samepage}
\begin{cor}\label{cor:gapp}
The following are equivalent:
\begin{enumerate}
\item[\textup{(i)}] $K$ is $\upl$-free;
\item[\textup{(ii)}] $\nabla_1(N_Q)=0$ for every homogeneous $Q\in K\{Y\}^{\ne}$ of degree $1$;
\item[\textup{(iii)}] for every homogeneous $Q\in K\{Y\}^{\neq}$ of degree $1$ there is $c\in C^\times$ such that 
$$ N_Q\ =\ cY\ \text{ or }\ N_Q\ =\ cY';$$
\item[\textup{(iv)}] for every $a\in K$ and $Q(Y)=aY'+Y''$,
we have $\nwt(Q)\le 1$. 
\end{enumerate}
\end{cor}
\end{samepage}
\begin{proof} Proposition~\ref{prop:nabla1} gives 
(i)~$\Rightarrow$~(ii).
Assume (ii), and let $Q\in K\{Y\}^{\neq}$ be homogeneous of degree $1$. Then $N_Q=cY^{(w)}$, $c\in C^\times$, $w=\nwt(Q)$. If $w>1$, then 
$$\nabla_1(N_Q)\ =\ -c\cdot {w\choose 2}\,Y^{(w-1)}\ \neq\ 0$$ by Lemma~\ref{diagnabla}, contradicting (ii).
Thus $w=0$ or $w=1$. 
This shows (ii)~$\Rightarrow$~(iii), and (iii)~$\Rightarrow$~(iv) is obvious. To show the contrapositive of (iv)~$\Rightarrow$~(i), suppose $\upl\in K$ is a pseudolimit of $(\upl_\rho)$.
Consider $$Q(Y)\ =\ \upl Y'+Y''\in K\{Y\}.$$  Then $Q^\phi=(\phi'+\phi \upl)Y'+\phi^2 Y''$, and 
by Lemma~\ref{lem:gap},
$$ v(\phi'+\phi \upl)\ =\ v\phi+v(\upl+\phi^\dagger)>v(\phi^2),$$ 
hence $\nwt(Q)=2$.
\end{proof}

\noindent
The following variant is also useful:

\begin{cor}\label{Euplone} Suppose $K$ is an immediate extension of its valued differential subfield $E$, and $E$ is $\upl$-free and 
$\frak{M}\subseteq E^\times$. Let $Q\in E\{Y\}^{\ne}$ be homogeneous of degree~$1$. Then
$N_Q=cY$ or $N_Q=cY'$ for some $c\in C^\times$.
\end{cor} 
\begin{proof} Since $N_Q$ is homogeneous of degree $1$ and isobaric of weight $w:= \nwt(Q)$, 
we have $N_Q=cY^{(w)}$ with $c\in C^\times$.  We have $\nabla_1 Q=0$ by the remark after the proof of Proposition~\ref{prop:nabla1}. Thus
$w=0$ or $w=1$ as in the proof of Corollary~\ref{cor:gapp}. 
\end{proof}

% \begin{cor}\label{Euplonemore} With the same assumptions 
%on $K$, $E$, $\frak{M}$ as in Corollary~\ref{Euplone}, 
%let $Q\in E\{Y\}^{\ne}$ have degree $1$. Then
%$N_Q=a+bY$ for some $a,b\in C$, not both zero, or  
%$N_Q=cY'$ for some $c\in C^\times$.
%\end{cor} 
%\begin{proof} Use Corollary~\ref{Euplone},  Lemma~\ref{npq}, %and the fact that $N_Q$ is isobaric. 
%\end{proof}

\noindent
The next lemma shows that $\upl$-freeness imposes restrictions on the shape of $N_P$ for any~$P$. Its corollary
shows this to be decisive for $\nwt(P)\le 3$.

\begin{lemma}\label{prop:nabla1, additional}
Suppose $K$ is $\upl$-free,  and $\nwt(P)>1$. Then 
$(N_P)_{|1|'}=0$. Also $(N_P)^{\phi,1}=0$ for every $\phi$ and $(N_P)^{t,1}=0$.
\end{lemma}
\begin{proof}
Set $w:=\nwt(P)>1$ and $N:=N_P$. Suppose towards a contradiction that $N_{|1|'}\neq 0$. Then $P_{|1|'}\ne 0$ and $N_{|1|'}=N_{P_{|1|'}}$ by Lemma~\ref{lem:NP components}, so after re\-pla\-cing~$P$ by~$P_{|1|'}$ we may assume that $P$ is subhomogeneous of subdegree~$1$. Using Lemma~\ref{lem:NP components} we may further reduce to the case that $P$ is homogeneous. Then ${P=Y^{d}Q}$ where~$d\in\N$ and $Q\in K\{Y'\}\subseteq K\{Y\}$ is homogeneous of degree $1$. Now by Corollary~\ref{cor:gapp}, $N=Y^{d}N_Q=cY^{d}Y'$, $c\in C^\times$; hence $\nwt(P)=\wt(N)=1$, a contradiction.
Thus $N_{|1|'} = 0$. 

Recall from Section~\ref{evtbeh} that ${\btau\brack\bsigma}=0$ for 
$\btau\geq\bsigma$ with $\supp \btau\neq\supp\bsigma$. Also, if 
$\|\bsigma\|=1$, $\btau\geq\bsigma$ and  
$\supp \btau=\supp\bsigma$, then $N_{[\btau]}=0$ (since $N_{|1|'}=0$). Hence
$$N^{\phi,1}\ =\ \sum_{\|\bsigma\|=1}\left(\sum_{\btau\geq\bsigma} {\btau\brack\bsigma} N_{[\btau]}\big(1+\varepsilon^{\btau}_{\bsigma}(\phi)\big)\right)Y^{[\bsigma]}\ =\ 0,$$ 
and
\equationqed{N^{t,1}\ =\ \sum_{\|\bsigma\|=1}\left(\sum_{\btau\geq\bsigma} {\btau\brack\bsigma} N_{[\btau]}\right)Y^{[\bsigma]}\ =\ 0.}
\end{proof}

\begin{cor}\label{cor:nwv leq 3} 
Suppose $K$ is $\upl$-free and $\nwt(P)\leq 3$. Then $N_P\in C[Y](Y')^{\N}$.
\end{cor}
\begin{proof} Assume $\nwt(P)=3$, and set $N:= N_P$. By \eqref{eq:Phi for isobaric P},
$$N^{\phi}\ =\ (\phi^\dagger)^3 N^{\phi,0} + \phi(\phi^\dagger)^2 N^{\phi,1} + \phi^2\phi^\dagger N^{\phi,2} + \phi^3 N^{\phi,3}.$$
Now $N^{\phi,0}$ is the isobaric part of $N$ of weight $0$, so $N^{\phi,0}=0$. Next, $N^{\phi,1}=0$ by Lemma~\ref{prop:nabla1, additional}. Also,
$N^{\phi,2}=0$ by Lemma~\ref{lem:properties of P<i>}(ii) and
Proposition~\ref{prop:nabla1}, and $N^{\phi,3}=N$. This gives $N^{\phi}=\phi^3 N$. Take $x\succ 1$ in an immediate asymptotic extension of $K$ with
$x'=1$ and set $t=1/x$. Replacing $\phi$ by $t$ in the above arguments still gives us $N^t=t^3N$, even if this immediate extension is not $\upl$-free. 
Apply Lemma~\ref{cor:P vanish} to get $N\in C[Y](Y')^3$.
The case $\nwt(P)\le 2$ is similar.
\end{proof}

\subsection*{The $\upo$-free case} Towards proving
Theorem~\ref{thm:nabla2} we need:

\begin{lemma}\label{newtred} Suppose
$\wt(P)=\wt(N_P)=w$, $D_{P}=N_P$, $P\asymp 1$, and 
$$ P_{[w]}\ =\ N_P+R, \quad \nabla_1 N_P=0, \quad v(R) > \Psi.$$ 
%with $\nabla_1P_{[w]} \prec^{\flat} 1$. 
Then there exists $\phi\preceq 1$ such that for
$Q:=P^{\phi}-(P^\phi)_{[w-1]}$ we have: $D_{Q^{\theta}}=N_P$ for all
$\theta\in \frak{M}^{\preceq 1}$ active in $K^{\phi}$; note that $Q_{[w-1]}=0$ and 
$N_P=N_Q$ for such $Q$.
\end{lemma}
\begin{proof} For all $\phi\preceq 1$ we have $P^\phi\sim (P^\phi)_{[w]}$ and by \eqref{eq:trans fms},
\begin{align*}
(P^{\phi})_{[w]}\ &=\ \phi^w P_{[w]}, \\
(P^\phi)_{[w-1]}\ &=\ \phi^{w-1}\cdot\big(P_{[w-1]}-\phi^\dagger\nabla_1P_{[w]}\big)\ =\ \phi^{w-1}\cdot\big( P_{[w-1]} - \phi^\dagger \nabla_1 R\big),
\end{align*}
with $v(\nabla_1R) > \Psi$, so
$v(\phi^\dagger\nabla_1R)> \Psi$. 
We claim that $v(P_{[w-1]})> \Psi$. Assuming towards a contradiction that
this claim is false, we get $\phi$ with $\phi\preceq P_{[w-1]}\preceq 1$, so
$$P^{\phi}\ \sim\ (P^{\phi})_{[w]}\ \asymp\ \phi^w\ \preceq\ (P^\phi)_{[w-1]}$$
by the identities above.
But $N_P=D_{P^{\phi}}$ is isobaric of weight $w$, so 
$(P^{\phi})_{[w-1]} \prec P^{\phi}$. This contradiction proves the claim. 
Take $\gamma > \Psi$ with 
$v(P_{[w-1]})>\gamma$ and $v(R)>\gamma$, 
and then take $\beta$ with $\Psi < \beta < \gamma$. Then $\alpha:= \gamma-\beta>0$, and for all $\phi\preceq 1$, 
$$v\big((P^{\phi})_{[w-1]}\big) - v(P^{\phi})\  >\ \gamma-v\phi\ >\ \alpha$$ 
by the identities above. Thus we can take $\phi\preceq 1$ such that 
$(P^\phi)_{[w-1]}\prec_{\phi}^{\flat} P^{\phi}$.
Then $Q:=P^{\phi}-(P^\phi)_{[w-1]}$ has by Lemma~\ref{newred} applied to $P^{\phi}$ the desired property. 
\end{proof}

\noindent
The next proposition and its corollary imply Theorem~\ref{thm:nabla2}.

\begin{prop}\label{prop:nabla1,nabla2}
Suppose $K$ is $\upo$-free.
Then $\nabla_1(N_P)=\nabla_2(N_P) =0$.
\end{prop}
\begin{proof} As in the proof of Proposition~\ref{prop:nabla1} we 
%arrange that $\Gamma$ is divisible, and 
first reduce to the case that
$w:=\nwt(P)=\wt(P)>0$ and $D_P=N_P$, as well as
$$ P\ \asymp\ 1, \qquad P_{[w]}\ =\ N_P+R, \quad R\prec^{\flat} 1.$$ 
Since $K$ is 
$\upo$-free, $K$ is $\upl$-free, so
by Proposition~\ref{prop:nabla1} we have $\nabla_1 N_P=0$.  If~$w=1$, then $N_P\in K[Y,Y']$, so $\nabla_2 N_P=0$ by
Corollary~\ref{cor:companion of Stirling}. 
Assume $w>1$ in the rest of the proof. Lemma~\ref{newtred} provides
a further reduction to the case $P_{[w-1]}=0$: this involves replacing
$K$ for a suitable $\phi\preceq 1$ by $K^{\phi}$ and 
$P$ by $\phi^{-w}Q$ with $Q:=P^{\phi}- (P^{\phi})_{[w-1]}$. 

Assume towards a contradiction that $\nabla_2 N_P\neq 0$. Then 
$$ \nabla_2 P_{[w]}\ =\ \nabla_2 N_P + \nabla_2 R, \qquad \nabla_2 N_P\in C\{Y\}^{\ne},\quad \nabla_2 R \prec^{\flat} 1.$$ 
Recall the pc-sequence $(\upo_\rho)$ of width $\{\gamma\in\Gamma_\infty:\gamma>2\Psi\}$ introduced in Section~\ref{sec:special cuts}. 
Take a pseudolimit $\upo$ of $(\upo_\rho)$
in an immediate asymptotic field extension of 
$K$. Since~$(\upo_\rho)$
has no pseudolimit in $K$, there is no $b\in K$ with
$b+\upo  \prec a^2$ for all active $a\in K$. This gives
$\phi\preceq 1$ with $P_{[w-2]} + \upo \nabla_2 P_{[w]} \succ \phi^2$.
From 
$P_{[w-1]}=0$, $\nabla_1 P_{[w]}=\nabla_1 R \prec^{\flat} 1$, and \eqref{eq:trans fms} we obtain, with $\lambda=-\phi^\dagger$, $\omega=-(2\lambda'+\lambda^2)$:
\begin{align*}  (P^{\phi})_{[w]}\ &=\ \phi^wP_{[w]}, \qquad (P^\phi)_{[w-1]}\ \prec^{\flat}\ 1, \\
  (P^\phi)_{[w-2]}\ &=\  \phi^{w-2}\big[P_{[w-2]} +
\omega \nabla_2 P_{[w]} + S\big], \quad S\ \prec^{\flat} 1. 
\end{align*} 
By Lemma~\ref{varrholemma} we have $\upo-\omega\prec\phi^2$, so
$$P_{[w-2]}+\omega\nabla_2P_{[w]} + S\ \sim\
P_{[w-2]}+\upo\nabla_2P_{[w]}\ \succ\ \phi^2.$$
By Proposition~\ref{compconjval, prop} we have $P^{\phi} \asymp \phi^w$, but
the above also gives
$$(P^{\phi})_{[w-2]}\ =\  \phi^{w-2}\big[P_{[w-2]} +
\omega \nabla_2 P_{[w]} + S\big]\  \succ\ \phi^w,$$
a contradiction.
\end{proof}

\begin{remark}
For later use we record a variant of Proposition~\ref{prop:nabla1,nabla2}:

\medskip\noindent
\textit{If $K$ is an immediate extension of its valued differential subfield 
$E$, and $E$ is $\upo$-free, $\frak{M}\subseteq E^\times$, and $P\in E\{Y\}^\times$, then $\nabla_1 N_P=\nabla_2N_P = 0$, so $N_P\in C[Y](Y')^{\N}$.}\/

\medskip\noindent
Taking into account the remark following the proof of Proposition~\ref{prop:nabla1}, the proof of 
this variant is the same as that of Proposition~\ref{prop:nabla1,nabla2}, apart 
from routine changes.
\end{remark} 

\begin{cor}
The following are equivalent:
\begin{enumerate}
\item[\textup{(i)}] $K$ is $\upo$-free;
\item[\textup{(ii)}] $\nabla_1(Q)=\nabla_2(Q)=0$ for every homogeneous $Q\in K\{Y\}^{\ne}$
of degree $2$;
\item[\textup{(iii)}] for every homogeneous  $Q\in K\{Y\}^{\neq}$ of degree $2$ we have:  
$$N_Q=cY^i(Y')^j \text{ for some $c\in C^\times$ and some $i,j\in \N$  with  $i+j=2$;}$$
\item[\textup{(iv)}]  for every $a\in K$ and $Q(Y)=a(Y')^2+2Y'Y^{(3)}-3(Y'')^2$,
we have $\nwt(Q)\le 3$. 
\end{enumerate}
\end{cor}
\begin{proof}
Proposition~\ref{prop:nabla1,nabla2} gives (i)~$\Rightarrow$~(ii). Suppose (ii)
holds, and let $Q\in K\{Y\}^{\neq}$ be homogeneous of degree~$2$. Then $N_Q$ is homogeneous of degree~$2$, and isobaric, so $N_Q\in C[Y](Y')^\N$ by 
Corollary~\ref{cor:companion of Stirling} and Lemma~\ref{cor:P vanish}. Thus $N_Q=cY^i(Y')^j$ with $c\in C^\times$, $i+j= 2$. This shows (ii)~$\Rightarrow$~(iii), and (iii)~$\Rightarrow$~(iv) is obvious. To show the contrapositive of (iv)~$\Rightarrow$~(i), suppose $\upo\in K$ is a pseudolimit of $(\upo_\rho)$. Let
$$N(Y)\ =\  2Y'Y^{(3)} - 3(Y'')^2 \in \Q\{Y'\}\subseteq K\{Y'\} \subseteq K\{Y\}.$$
Note that $N$ is homogeneous of degree $2$, isobaric of weight $4$, and
$$N^\phi\ =\ \phi^4 N + (2\phi\phi'' - 3(\phi')^2) (Y')^2 \in  K^\phi\{Y'\}.$$
Consider the differential polynomial 
$$Q\ :=\  -\upo \cdot (Y')^2 + N  \in K\{Y'\}.$$
We have
$$Q^\phi\ =\ \big(2\phi\phi'' - 3(\phi')^2 - \upo\phi^2\big) \cdot (Y')^2 + \phi^4 N.$$
Setting $\lambda:=-\phi^\dagger$ and  $\omega:=-(2\lambda'+\lambda^2)$, we get 
$$(\phi')^2\ =\ \lambda^2 \phi^2, \qquad \phi\phi''\ =\ 
-\phi(\lambda\phi)'\  =\  (-\lambda' + \lambda^2)\phi^2,$$
so by Lemma~\ref{varrholemma}, 
$$2\phi\phi''-3(\phi')^2 - \upo\phi^2\ =\ 
(-2\lambda'+2\lambda^2-3\lambda^2 - \upo)\phi^2\ =\ 
(\omega - \upo)\phi^2\ \prec\ \phi^4.$$
Thus $D_{Q^\phi}=N$ for all $\phi$, and hence $\nwt(Q)=\wt(N)=4$. 
\end{proof}  

\noindent
Suppose $K$ is $\upo$-free. Then $N_{P_{\times\fm}}=A(Y)(Y')^j$ with $A\in C[Y]$ and $j\in \N$. Let~${y\sim c\fm}$ ($c\in C^\times$) be an approximate zero of $P$. Then the multiplicity of $y$ as an approximate zero of $P$ is $i+j$, where $i$ is the
multiplicity of $A$ at $c$; we call~$i$ the {\bf algebraic multiplicity} of the  approximate zero $y$ of $P$. We say that
$y$ is an {\bf algebraic approximate zero of $P$} if $A(c)=0$, 
equivalently, its algebraic multiplicity is $\geq 1$. \index{multiplicity!algebraic}\index{approximate!zero}\index{approximate!zero!algebraic}

\begin{lemma}\label{lem:sum of alg mult}
Suppose $K$ is $\upo$-free. Let
$P$ and $\fm$ be given. Then there
are distinct $c_1,\dots,c_n\in C^\times$ such that $c_1\fm,\dots,c_n\fm$ are algebraic approximate zeros of~$P$ with respective algebraic multiplicities $\mu_1,\dots,\mu_n$,
and such that for any algebraic approximate zero~$y\asymp \fm$ of $P$ we have 
$y\sim c_k\fm$ for some $k$. These properties uniquely determine the set 
$\big\{(c_1, \mu_1),\dots, (c_n,\mu_n)\big\}$. If $C$ is algebraically closed, then 
$$\mu_1+\cdots+\mu_n\ =\ \ndeg P_{\times\fm} - \nval P_{\times\fm}, \phantom{\mod 2}$$
and if $C$ is real closed, then
$$\mu_1+\cdots+\mu_n\ \equiv\ \ndeg P_{\times\fm} - \nval P_{\times\fm}\mod 2.$$
\end{lemma}
\begin{proof} Take $A\in C[Y]$ without zeros in $C$, $i,j\in \N$, distinct $c_1,\dots, c_n\in C^{\times}$, and $\mu_1,\dots, \mu_n\in \N^{\ge 1}$ 
such that
$$N_{P_{\times\fm}}\ =\ A(Y)\cdot Y^i\cdot \left(\,\prod_{k=1}^n (Y-c_k)^{\mu_k}\right) \cdot (Y')^j.$$
Then $c_1,\dots, c_n, \mu_1,\dots, \mu_n$ have the desired
property. If $C$ is algebraically closed, then $\deg A=0$. If $C$ is real closed, then $\deg A$ is even.
\end{proof}

\subsection*{Newton polynomials and upward shift} This subsection assumes familiarity with Appendix~\ref{app:trans}.
Consider $\T$ equipped with its usual ordering, valuation and derivation $\der=\frac{d}{dx}$. Then $\T$ is an $\upo$-free $H$-field with small derivation, and constant field $\R$. The group of transmonomials is a monomial group of~$\T$, and we equip~$\T$ with this monomial group. We have the usual logarithmic sequence~$(\ell_n)$ for $\T$ given by $\ell_0=x$ and $\ell_{n+1}=\log(\ell_n)$, with corresponding
$\upg_n:= \frac{1}{\ell_0\ell_1\cdots \ell_n}$. It is easy to check that the map $f\mapsto f{\uparrow}\colon\T \to \T$ is for each $n$ an
isomorphism $\T^{\upg_{n}}\to \T^{\upg_{n-1}}$ of $H$-fields,
where $\upg_{-1}:=1$. Defining $f{\uparrow^0}:=f$ and $f{\uparrow^{n+1}}:=(f{\uparrow^n}){\uparrow}$ for $f\in \T$, we obtain for each $n$ an isomorphism $f\mapsto f{\uparrow^n}\colon \T^{\upg_{n-1}}\to \T$ of $H$-fields. 

Let $P\in \T\{Y\}$. Recall that $P{\uparrow}$ denotes the differential polynomial in $\T\{Y\}$ obtained by applying $f\mapsto f{\uparrow}$ to the coefficients of $P^{1/x}\in \T^{1/x}\{Y\}$, and that $P(y){\uparrow}=P{\uparrow}(y{\uparrow})$ for $y\in\mathbb T$. 
Defining inductively $P{\uparrow^0}:=P$ and $P{\uparrow^{n+1}}:=(P{\uparrow^n}){\uparrow}$, we have $P(y){\uparrow^n}=P{\uparrow^n}(y{\uparrow^n})$ for all $n$ and $y\in\mathbb T$.

\begin{lemma} Let $P\in \T\{Y\}$. Then $P{\uparrow^n}\in \T\{Y\}$ is obtained by applying 
the operation $f\mapsto f{\uparrow^n}$
to the coefficients of $P^{\upg_{n-1}}\in\T^{\upg_{n-1}}\{Y\}$, and
$D_{P{\uparrow^n}} = D_{P^{\upg_{n-1}}}$.
\end{lemma}
\begin{proof}
Let $Q\in\T\{Y\}$ be  obtained by applying $f\mapsto f{\uparrow^n}$
to the coefficients of $P^{\upg_{n-1}}$. By Lemma~\ref{notthin}, it suffices to show that 
$P(y){\uparrow^n}=Q(y{\uparrow^n})$ for all $y\in\T$. 
For this we may assume $P=Y^{(i)}$, $i\geq 1$, so 
$P^{\upg_{n-1}}=\sum_{j=1}^i F^i_j(\upg_{n-1})Y^{(j)}$. With $\derdelta$ the derivation of~$\T^{\upg_{n-1}}$, we then have $\der^j(y{\uparrow^n})=\derdelta^j(y){\uparrow^n}$ for all $j\in \N$, so
\begin{align*}
Q(y{\uparrow^n}) &= \sum_{j=1}^i F^i_j(\upg_{n-1}){\uparrow^n}\, \der^j(y{\uparrow^n}) 
 = \sum_{j=1}^i F^i_j(\upg_{n-1}){\uparrow^n}\, \derdelta^j(y){\uparrow^n} \\
&= P^{\upg_{n-1}}(y){\uparrow^n} = P(y){\uparrow^n}
\end{align*}
for $y\in \T$, as required. To prove the identity, assume $P\ne 0$, so 
$$ P^{\upg_{n-1}}\ =\ \fd D_{P^{\upg_{n-1}}} + R_n,\quad \fd:= \fd_{P^{\upg_{n-1}}},\ R_n\in \T^{\upg_{n-1}}\{Y\},\ R_n\prec P^{\upg_{n-1}},$$
so $P{\uparrow^n}= (\fd{\uparrow^n})\cdot D_{P^{\upg_{n-1}}} + S_n$
with $S_n\in \T\{Y\}$, $S_n\prec P{\uparrow^n}$. It remains to note that~$\fd{\uparrow^n}$ is a transmonomial. 
\end{proof}

\noindent
Since $\T$ is $\upo$-free, we may conclude:

\begin{cor}\label{cor:newton poly, upward shift}
Let $P\in\T\{Y\}^{\neq}$. Then
there is an $n_0\in\N$ such that 
$$D_{P{\uparrow^n}}=N_P\in \R[Y](Y')^\N\qquad\text{ for all $n\geq n_0$.}$$
\end{cor}

\subsection*{Notes and comments}
Corollary~\ref{cor:newton poly, upward shift} is a translation to our ``compositional conjugation'' setting of \cite[Section~8.3.1]{JvdH}.
Some of the other basic properties of the Newton polynomial established in Section~\ref{basicnewton} 
(for example, Lemma~\ref{dn2}) 
were first shown for $K=\T_{\operatorname{g}}$ in \cite[Section~8.3]{JvdH}.
We do not need this, but it is worth mentioning that in Corollary~\ref{cor:newton poly, upward shift} one can take
$n_0=2\dwv(P)$. \marginpar{taken on faith}

\section{Realizing Cuts in the Value Group}\label{cutsvalgrp}

\noindent
{\em Throughout this section $K$ is $\upo$-free}\/ (and thus clean, by Corollary~\ref{cor:nabla2}). For later use we study realizations
of the following three cuts in the value group $\Gamma$ of $K$:
$$\Gamma^{<},\quad\Gamma^{\leq },\quad (\Gamma^<)'=\Psi^{\downarrow}.$$
Given any element $y$ in any $\d$-valued field extension $L$ of $H$-type of $K$ such that~$vy$ realizes one of these cuts, we shall derive some informative results about 
the  $H$-asymptotic field extension $K\langle y\rangle$ of $K$. 
If $vy$ realizes the cut $\Gamma^<$ or the cut $\Gamma^{\leq}$ (that is, $\Gamma^< < vy<0$ or
$0<vy<\Gamma^{>}$), then this is fairly straightforward, based on earlier results in this chapter.
The case where $vy$ realizes $\Psi^{\downarrow}$ (that is, $\Psi<vy<(\Gamma^>)'$) requires some facts about the operation $Q\mapsto Q^{\times\phi}$ on differential polynomials.

\subsection*{The cut $\Gamma^<$} 
Model-theoretic compactness yields
a $\d$-valued field extension $L$ of $H$-type of $K$ with an element $y\in L$ such that
$\Gamma^< < vy<0$. Let such $L$ and $y\in L$ be given.
We determine the asymptotic couple of $K\<y\>$:

\begin{lemma}  \label{lem:infinitesimal cut}
We have $v\big( P(y) \big) = v^{\ev}(P) + \ndeg(P)\,vy + \nwt(P)\,\psi_L(vy)$.
Moreover, $\Gamma_{K\<y\>} = \Gamma \oplus \Z vy \oplus \Z\psi_L(vy)$ \textup{(}internal direct sum\textup{)} with $\max \Psi_{K\<y\>}=\psi_L(vy)$.
\end{lemma}
\begin{proof}
The first identity follows from Corollary~\ref{cor:sign and valuation, 1}. As this holds for every~$P$, we get
$\Gamma_{K\<y\>} = \Gamma + \Z vy + \Z\psi_L(vy)$. Set $\alpha:= vy$, 
$\beta:= \psi_L(vy)$, and $\Gamma_1=\Gamma+\Z\beta$. Since $K$ has rational asymptotic integration by Corollary~\ref{cor:lambda-free alg closure}, and $\Psi < \beta < (\Gamma^{>})'$, we get from Corollary~\ref{addgapmore} that 
$k\beta\notin \Gamma$ for all nonzero $k\in \Z$, and
$[\Gamma_1]=[\Gamma]$, so $\psi_L\big(\Gamma_1^{\ne}\big)=\Psi$.
With $\psi_1$ the restriction of $\psi_L$ to $\Gamma_1^{\ne}$ we get an
$H$-asymptotic couple~$(\Gamma_1, \psi_1)$ with gap $\beta$.
Since $0< n|\alpha| < \Gamma_1^{>}$ for all $n\ge 1$, and $\psi_L\big(|\alpha|\big)=\beta$, this gives the desired result in view of %Corollary~{5.2.24} and 
$[\Gamma_1+\Z\alpha]=[\Gamma] \cup \big\{[\alpha]\big\}$.
\end{proof}

\noindent
The extension $K\langle y \rangle$ of $K$ is determined up to isomorphism: 

\begin{lemma} \label{lem:infinitesimal cut, uniqueness}
Let $y^*$ also be an element of a $\d$-valued field extension $L^*$ of $H$-type of~$K$
such that  $\Gamma^{<}<vy^*<0$. Then 
there is a unique valued differential field embedding 
$K\langle y \rangle\to L^*$ that is the identity on $K$ and
sends $y$ to $y^*$. 
\end{lemma}
\begin{proof}
Note that $P(y)\ne 0$ for all $P$ by the first part of Lemma~\ref{lem:infinitesimal cut},
and likewise for $y^*$, that is, $y$ and $y^*$ are 
$\d$-transcendental over $K$. The proof of that lemma 
also gives an ordered group isomorphism
$\Gamma_{K\<y\>} \to \Gamma_{K\<y^*\>}$ that is the identity on $\Gamma$, sends~$vy$ to~$vy^*$, and $\psi_L(vy)$ to $\psi_{L^*}(vy^*)$. This yields the desired result.  
\end{proof}

\subsection*{The cut $\Gamma^{\le}$} There also exists
a $\d$-valued field extension $L$ of $H$-type of $K$ with an element $y\in L$ such that
$0 < vy <\Gamma^{>}$. Let such $L$ and $y\in L$ be given. Then the following
analogues of the lemmas above hold, with similar proofs:

\begin{lemma}\label{lem: otherinfcut}  
We have $v\big( P(y) \big) = v^{\ev}(P) + \nval(P)\,vy + \nwt(P)\,\psi_L(vy)$.
Moreover, $\Gamma_{K\<y\>} = \Gamma \oplus \Z vy \oplus \Z\psi_L(vy)$ \textup{(}internal direct sum\textup{)} with $\max \Psi_{K\<y\>}=\psi_L(vy)$.
If $y^*$ is an element of a $\d$-valued field extension $L^*$ of 
$H$-type of $K$
such that  $0 < vy^* < \Gamma^{>}$, then 
there is a unique valued differential field embedding 
$K\langle y \rangle\to L^*$ that is the identity on $K$ and
sends $y$ to $y^*$. 
\end{lemma}

\subsection*{The residue field of $K\langle y \rangle$} Let $L$ be a $\d$-valued field extension of $H$-type of~$K$ and let $y\in L$.

\begin{lemma}  \label{lem:infinitesimal cut, residue field}
Suppose $0<\abs{vy}<\Gamma^{>}$. 
Then $\res K\langle y \rangle=\res K$, and so 
$K\langle y \rangle$ is $\d$-valued with $C_{K\langle y \rangle}=C$.
\end{lemma}
\begin{proof}
Let $f\in K\langle y \rangle^\times$, so $f=P(y)/Q(y)$ where 
$P,Q\in K\{Y\}^{\ne}$. We have $N_P=D(Y)\cdot (Y')^w$ with $D(Y)\in C[Y]^{\neq}$ and $w=\nwt(P)$. The proof of Lemma~\ref{lem:sign and valuation} shows that 
eventually, with $\derdelta=\phi^{-1}\der$,
$$P(y)\ =\ P^\phi(y)\ \sim\ \frak{d}_{P^\phi}D(y)\derdelta(y)^w\ =\
\frak{d}_{P^\phi}D(y)(y'/\phi)^w\ =:\ p(\phi)\in K(y,y')^\times,$$
and likewise, $Q(y)\sim q(\phi)\in K(y,y')^\times$, eventually, 
so $f\sim p(\phi)/q(\phi) \in K(y,y')^\times$.
Thus $\res K\langle y \rangle=\res K(y,y')$. 
By Lemmas~\ref{lem:infinitesimal cut} and \ref{lem: otherinfcut} we have 
$$v K(y,y')^\times\ =\ \Gamma \oplus \Z vy \oplus \Z \psi_L(vy),$$
so $\res K(y,y')=\res K$ by Lemma~\ref{lem:lift value group ext} applied to the valued
field extensions $K(y)|K$ and $K(y,y')|K(y)$. 
\end{proof}

\noindent
Before we study the extension $K\<y\>$ of $K$ in the case where $\Psi<vy<(\Gamma^>)'$, we turn our attention to the operation $Q\mapsto Q^{\times\phi}=Q^\phi_{\times\phi}$.

\index{dominant!part!eventual, of $P^{\times\phi}$}
\nomenclature[X]{${v^{\times\hskip-0.1em\operatorname{e}}(P)}$}{${v^{\operatorname{e}}\big(P(Y')\big)}$}
\nomenclature[X]{$N^\times_P$}{eventual dominant part of $P^{\times\phi}$}
\nomenclature[X]{$\nwt^\times(P)$}{Newton weight of $P(Y')$}
%\nomenclature[X]{$v^{\times\!\operatorname{ev}}(P)$}{$v^{\operatorname{ev}}\big(P(Y')\big)$}

\subsection*{Combining multiplicative and compositional conjugation}
In Section~\ref{Compositional Conjugation} we defined 
$P^{\times \phi}:= P^{\phi}_{\times \phi}$, so
$P^{\times\phi}(Y')=P(Y')^\phi$. Hence
$D_{P^{\times\phi}}(Y') = D_{P(Y')^\phi}$, and so  
$D_{P^{\times\phi}}(Y') = N_{P(Y')} \in C\{Y'\}\subseteq C\{Y\}$, eventually.  Also $N_{P(Y')}\in C[Y](Y')^{\N}$, so
$N_{P(Y')}=c(Y')^w$ with $c\in C^\times$ and $w\in \N$. We set $N^\times_P:= cY^w$, so $N^\times_P(Y')=N_{P(Y')}$,  so
$D_{P^{\times\phi}}=N^\times_P$, eventually. We call $N^\times_P$ the {\bf eventual dominant part of~$P^{\times\phi}$}. For the above $w$ we set 
$$\nwt^\times(P)\ :=\ w\ =\ \nwt\!\big(P(Y')\big).$$
If $P$ is homogeneous of degree $d$, then so is each
$P^{\times\phi}$, hence $\nwt^\times(P)=d$.
Set
$$v^{\times\!\ev}(P)\ :=\ v^{\ev}\big(P(Y')\big),$$
so in view of $v(P^{\times \phi})=v\big(P(Y')^{\phi}\big)$ we get
$$v(P^{\times\phi})\ =\  v^{\times\!\ev}(P) + \nwt^\times(P)v\phi,\quad\text{eventually.}$$
Next an analogue of Corollary~\ref{cor:sign and valuation, 1} for
$v\big(P(z)\big)$ for active $z$. 

\begin{cor}\label{cor:mult-comp, 1}
There is an active $\phi_0=\phi_0(P)$ in $K$ such that for every $\d$-valued field extension $L$ of $H$-type of $K$ and every active 
$z\prec \phi_0$ in $L$,
$$vP(z) \ =\ v^{\times\!\ev}(P)+\nwt^\times(P)\,vz.$$
\end{cor}
\begin{proof}
Take active $\phi_0$ in $K$ such that 
$$D_{P^{\times\phi_0}}=N^\times_P, \quad  R_{P(Y')^{\phi_0}}\prec^{\flat}_{\phi_0} P(Y')^{\phi_0}, \quad v(P^{\times\phi_0})=  v^{\times\!\ev}(P) + \nwt^\times(P)v(\phi_0).$$
Using Lemma~\ref{dn1}, the second relation gives $R_{P^{\times \phi_0}}\prec^{\flat}_{\phi_0} P^{\times \phi_0}$. 
Let $L$ be a $\d$-valued field extension of $H$-type of $K$ and $z$ an
active element of $L$ with $z\prec\phi_0$. Then $y:=z\phi_0^{-1}\prec 1$ is active in $L^{\phi_0}$, so $y \asymp^\flat_{\phi_0} 1$ and
$P(z)=P^{\times\phi_0}(y)$. Hence by Lemma~\ref{lem:sign and valuation}
and using $\wt(N_{P^\times})=0$:
\alignqed{
v\big(P(z)\big)\ 	&=\  v(P^{\times \phi_0})+\dv(P^{\times\phi_0})\,v(y)+\dwt(P^{\times\phi_0})\,\psi_{L^{\phi_0}}(vy) \\
					&=\  v(P^{\times\phi_0})+\nwt^\times(P)\,v(y) \\
					&=\  v^{\times\!\ev}(P)+\nwt^\times(P)\,vz.}
\end{proof}

\noindent
Here is an immediate consequence:

\begin{cor}\label{cor:mult-comp, 2}
There exists $y_0\in K^\times$ with $y_0\prec 1$ such that for every $\d$-valued field extension $L$ of $H$-type of $K$ and any $y\in L$,
$$ 0<\abs{vy}<v(y_0)\ \Longrightarrow\ 
v\big(P(y^\dagger)\big) \ = \ v^{\times\!\ev}(P)+\nwt^\times(P)\psi_L(vy).$$
\end{cor}

\begin{cor}\label{cor:mult-comp, 3}
Suppose $P$ is homogeneous. Then 
$\Ric(N_P)=N_{\Ric(P)}^\times$, and
$$\nwt^\times(P)\ =\ \deg(P),\quad \nwt(P)\ =\ \nwt^\times\!\big(\!\Ric(P)\big),\quad v^{\ev}(P)\ =\ v^{\times\!\ev}\big(\!\Ric(P)\big).$$
\end{cor}

\begin{proof}
By Lemma~\ref{Ri and compositional conjugation} we have $\Ric(P^\phi)=\Ric(P)^{\times\phi}$, and by Lemma~\ref{dn1a} we have
$D_{\Ric(P^\phi)}=\Ric(D_{P^\phi})$, hence, eventually,
$$\Ric(N_P)\ =\ \Ric(D_{P^\phi})\ =\  D_{\Ric(P^\phi)}\ =\  D_{\Ric(P)^{\times\phi}}\ =\ N^\times_{\Ric(P)}.$$
Set $d:=\deg(P)$ and $R:=\Ric(P)\in K\{Z\}$. Since $P$ is homogeneous, we have $\nwt^\times(P)=d$.
Moreover, $N_P\in C^\times\, Y^{d-w}(Y')^w$ where $w=\nwt(P)$,
so $N^\times_{R}=\Ric(N_P)\in C^\times\,Z^w$, and thus $\nwt^\times(R)=w$.
To get $v^{\ev}(P)=v^{\times\!\ev}(R)$ we apply Corollary~\ref{cor:mult-comp, 2} to $R$ in place of $P$. This gives $y_0\in K^\times$ with $y_0\prec 1$ such that for all $y\in K$ with $y_0\prec y\prec 1$ we have
$$v\big(R(y^\dagger)\big)\ =\  v^{\times\!\ev}(R) + \nwt^\times(R)\psi(vy).$$
By Corollary~\ref{cor:sign and valuation, 1} we can arrange that in addition, for all $y\in K$ with $y_0\prec y\prec 1$,
$$v\big(P(y)\big)\ =\ v^{\ev}(P) + d\,vy + \nwt(P)\psi(vy).$$
Since $R(y^\dagger)=P(y)/y^d$ and $\nwt^\times(R)=\nwt(P)$, we get $v^{\times\!\ev}(R)=v^{\ev}(P)$.
\end{proof}

\noindent
We also have a version of Corollary~\ref{cor:mult-comp, 1} for the behavior of $v\big(P(z)\big)$
for $z\in K$ as~$vz$ approaches $\Psi$ from above:  

\begin{cor}\label{cor:mult-comp, 4}
There exists $\delta_0\in (\Gamma^>)'$
such that for every $\d$-valued field extension $L$ of $H$-type
of $K$ and every 
$z\in L$ with 
$\Psi<vz \le \delta_0$,
$$v\big(P(z)\big) \ =\ v^{\times\!\ev}(P)+\nwt^\times(P)\,vz.$$
\end{cor}
\begin{proof}
Take $\phi_0$ as in the proof of Corollary~\ref{cor:mult-comp, 1}.  
By Lemma~\ref{hasc} and the remark following it, we can take
$\gamma_0\in\Psi^{>v\phi_0}$ and $\delta_0\in (\Gamma^>)'$ such that
the map 
$$\gamma\mapsto \psi_L(\gamma-v\phi_0)$$ is constant on
the set $[\gamma_0,\delta_0]_{\Gamma_L}$, for each $\d$-valued field extension
$L$ of $H$-type of~$K$. 
Given such an $L$  and $z\in L$ with
$\Psi<vz\le\delta_0$, set $y:=z\phi_0^{-1}\prec 1$. Then $\psi_L(vy)=\psi_L(vz-v\phi_0)=\psi(\gamma_0-v\phi_0)$,
and the latter is $>v\phi_0$ by Lemma~\ref{BasicProperties}(i).
Thus $y \asymp^\flat_{\phi_0} 1$. Since $P(z)=P^{\times\phi_0}(y)$, the
claims now follow as in the proof of Corollary~\ref{cor:mult-comp, 1}.
\end{proof}

\noindent
Here is an easy consequence of Corollaries~\ref{cor:mult-comp, 1} and~\ref{cor:mult-comp, 4}:

\begin{cor}\label{cor:mult-comp, 5}
There exists $y_0\in K^\times$ with $y_0\prec 1$ such that for every $\d$-valued field extension $L$ of $H$-type of $K$ and all $y\in L$, 
$$0<\abs{vy}<v(y_0)\ \Longrightarrow\ v\big(P(y')\big) \ = \ v^{\times\!\ev}(P)+\nwt^\times(P)v(y').$$
\end{cor}

\subsection*{The cut $\Psi^\downarrow$}
We now return to the setting of the beginning of this section. 
Model-theoretic compactness gives
a $\d$-valued field extension $L$ of $H$-type of $K$ with an element $z\in L$ such that $\Psi<vz<(\Gamma^>)'$. Let such an $L$ and $z\in L$ be given. Then $v\big(P(z)\big)  =  v^{\times\!\ev}(P) + \nwt^\times(P)\,vz$ by 
Corollary~\ref{cor:mult-comp, 4}. 

\begin{lemma}\label{lem:Psi-cut} We have
$\Gamma_{K\<z\>} = \Gamma \oplus \Z vz$ \textup{(}internally\textup{)} with $[\Gamma_{K\<z\>}]=[\Gamma]$.
\end{lemma}
\begin{proof}
Clearly $\Gamma_{K\<z\>} = \Gamma + \Z vz$.
From Corollary~\ref{addgapmore} we get $\Gamma + \Z vz = \Gamma \oplus \Z vz$ and $[\Gamma + \Z vz]=[\Gamma]$. 
\end{proof}

\noindent
Corollary~\ref{addgapmore} and Lemma~\ref{lem:Psi-cut} give an analogue of Lemma~\ref{lem:infinitesimal cut, uniqueness}: 

\begin{lemma}\label{zstargap} 
Let $z^*$ also be an element of a $\d$-valued field extension $L^*$
of $H$-type of $K$
such that  $\Psi<vz^*<(\Gamma^>)'$. Then 
there is a unique valued differential field embedding 
$K\langle z \rangle\to L^*$ that is the identity on $K$ and
sends $z$ to $z^*$. 
\end{lemma}

\begin{cor}\label{corzstargap} The extension $K\<z\>$ is $\d$-valued with $C_{K\<z\>}=C$.
\end{cor}
\begin{proof} Let $\phi_0\in K$ be as in the proof of Corollary~\ref{cor:mult-comp, 1}. So $D_{P^{\times\phi_0}}=N^\times_P=c Y^w$ with $c\in C^\times$ and $w=\nwt^\times(P)$. For $y:= z\phi_0^{-1}$ we have $P(z)=P^{\times \phi_0}(y)$. The proofs of Lemmas~\ref{lem:sign and valuation} and Corollary~\ref{cor:mult-comp, 4} show that then
$$P(z)\ =\ P^{\times\phi_0}(y)\ \sim\ \frak{d}_{P^{\times\phi_0}}cy^w \in K(z)^\times.$$
Thus $\res K\langle z \rangle=\res K(z)$, and as $v K(z)^\times=\Gamma \oplus \Z vz$ by Lemma~\ref{lem:Psi-cut}, we have
$\res K(z)=\res K$ by Lemma~\ref{lem:lift value group ext}. It remains to
appeal to Lemma~\ref{dv}.
\end{proof}

\subsection*{Notes and comments} The following is worth mentioning. \marginpar{taken on faith}
There is a $\delta_0\in (\Gamma^{>})'$ such that
$N_{P_{\times \fm}}=N^{\times}_{P}$ for all $\fm$ with $\Psi < v\fm < \delta_0$. Thus $\ndeg_{\E}(P) = \nwt^\times(P)$ for $\E:=\{y\in K^\times:\ vy > \Psi\}$. (This will not be used
in the present volume.) 

\section{Eventual Equalizers}\label{eveqth}

\noindent
{\em Throughout this section $K$ is $\upo$-free}. The goal of this section is to 
prove the Eventual Equalizer Theorem~\ref{evtequalizer} and to derive from it
a kind of Newton diagram for our differential polynomial 
$P\in K\{Y\}^{\ne}$: Proposition~\ref{newtpolygon}.

\subsection*{More general multiplicative-compositional conjugations}
In the proof of the Eventual Equalizer Theorem we need to consider more general combinations of multiplicative and compositional conjugation than in Section~\ref{cutsvalgrp}. To define these we assume in this subsection that $K$ is algebraically closed, and we also   
fix an algebraic closure $F=K\<Y\>^\alg$ of $K\<Y\>$. We equip $F$ with the unique extension of the derivation $\der$ of $K\<Y\>$ to a derivation, also denoted by $\der$, of $F$. Next, extend the (gaussian) valuation of $K\<Y\>$ to a valuation of $F$.
This makes $F$ a valued differential field, and its derivation is small by  
the remarks at the beginning of Section~\ref{sec:resext} and by
Proposition~\ref{Algebraic-Extensions}. Also by Section~\ref{sec:resext}, the image $y$ of $Y\in \mathcal{O}_F$ in the differential residue field $\k_F$ of $F$ is $\d$-transcendental over the (trivial) differential residue field~$\k$ of $K$, and $\k\<y\>$ is the differential residue field of $K\<Y\>$. Thus
$\k_F=\k\<y\>^\alg$. Since~$\k$ is algebraically closed, Lemma~\ref{lem:composition in K<Y>a} gives:

\begin{lemma}\label{lem:composition in res(F)}
If $Q\in \mathcal O_F$ and $\overline{Q}\notin\k$, then $v\big(P(Q)\big)=v(P)$.
\end{lemma}

\noindent
Next, let any $\phi$ be given. Then $F$ is also an algebraic closure of
$K^{\phi}\<Y\>$, and the unique extension of the derivation $\derdelta=\phi^{-1}\der$ of $K^{\phi}\<Y\>$ to a derivation $\derdelta$ of $F$
is again small with respect to the given valuation of $F$, by the same
arguments as before. Note that~$K^{\phi}$ has the same (trivial) differential residue field $\k$ as $K$, and that the lemma above goes through as follows: 

\begin{lemma}\label{lem:phicomposition in res(F)}
If $Q\in \mathcal O_F$ and $\overline{Q}\notin\k$, then $v\big(P^{\phi}(Q)\big)=v(P^{\phi})$.
\end{lemma}

\noindent
As in Section~\ref{Compositional Conjugation} we extend the ring isomorphism $P\mapsto P^\phi\colon K\{Y\}\to K^\phi\{Y\}$ to an automorphism $R\mapsto R^\phi$ of the field 
$F$, and choose a 
map $$(R,q)\mapsto R^q\colon F^\times\times\Q \to F^\times$$
extending the usual map $(R,k)\mapsto R^k\colon F^\times \times \Z \to F^\times$, subject to $(R^q)^\dagger=qR^\dagger$ for $R\in F^\times$, $q\in \Q$. Thus $a^q\in K^\times$ for $a\in K^\times$. In that same section we defined
$$P^{\times q, \phi}\ :=\  P^{\phi}_{\times \phi^q}\in K^{\phi}\{Y\}\quad (q\in \Q), \qquad P^{\times \phi}\ :=\  P^{\times 1,\phi} = P^\phi_{\times\phi}.$$ 

\begin{prop}\label{prop:valuation of Ptimesqphi}
Suppose $P$ is homogeneous, $d=\deg(P)$, and $q\in\Q^\times$. Then there exists an $\alpha\in\Gamma$ such that
$$v(P^{\times q,\phi})\ =\ \alpha+dq\,v\phi,\qquad\text{eventually.}$$
\end{prop}
\begin{proof}
Lemma~\ref{lem:Ptimesqphi} gives homogeneous $E\in K\{Y\}^{\neq}$ of degree $w:= \wt(P)$, with 
$$P^{\times q,\phi}\big( (Y')^q \big)\ =_{\operatorname{c}}\ (\phi Y')^{dq-w}\cdot E^{\times \phi}(Y')\qquad\text{for each $\phi$.}  $$
Set $\alpha:=v^{\times\!\ev}(E)$. By Corollary~\ref{cor:mult-comp, 3} we have $\nwt^{\times}(E)=\deg(E)=w$. It is obvious that $v\big(E^{\times \phi}(Y')\big)=v(E^{\times \phi})$ for all $\phi$,
hence
$$v\big(E^{\times \phi}(Y')\big)\ =\  v^{\times\!\ev}(E)+\nwt^{\times}(E)v\phi\ =\  \alpha + w\,v\phi,\qquad\text{eventually.}$$
By Lemma~\ref{lem:phicomposition in res(F)} we also have
$v\big(P^{\times q,\phi}\big( (Y')^q \big) \big)=v(P^{\times q,\phi})$. Hence 
\begin{align*}
v\big(P^{\times q,\phi}\big( (Y')^q \big) \big)\   &=\  v\big((\phi Y')^{dq-w}\big) + v\big(E^{\times \phi}(Y')\big) \\
	&=\  (dq-w)\,v\phi + \big( \alpha + w\,v\phi \big)\ 
											=\ \alpha + dq\,v\phi,
\end{align*}
eventually.
\end{proof}

\subsection*{Proof of the Eventual Equalizer Theorem} 
Complementing the ``eventual'' terminology, we say that a property $S(\phi)$ of ele\-ments~$\phi$ holds {\bf cofinally} if for every active $\phi_0$ in $K$ there is a $\phi\preceq \phi_0$ such that $S(\phi)$ holds.

\index{cofinally}

\begin{proof} [Proof of Theorem~\ref{evtequalizer}] We assume that $\Gamma$ is divisible and $P,Q\in K\{Y\}^{\neq}$ are homogeneous of degrees $d>e$.
Our job is to show that there exists~${a\in K^\times}$ such that, eventually, 
$ P^\phi_{\times a}\  \asymp\  Q^\phi_{\times a}$. Note that $v(P^{\phi}_{\times a})$ depends only on $va$ for~$a\in K^\times$, where~$P$ and 
$\phi$ are given. Thus by passing to the algebraic closure of $K$ we can arrange that~$K$ is algebraically closed. We extend the usual map ${(k,a)\mapsto
a^k\colon \Z\times K^\times \to K^\times}$ to a map $(q,a)\mapsto
a^q\colon \Q\times K^\times \to K^\times$ such that $(a^q)^\dagger=qa^\dagger$ for all~$a\in K^\times$ and $q\in \Q$.

Take an elementary extension $K_*$ of $K$ with an active $z$ in $K_*$ such that $z\prec\phi$ for every $\phi$. Then $\Gamma_{K\<z\>}=\Gamma\oplus\Z v(z)$, by Lemma~\ref{lem:Psi-cut}.
Let $L:=K\<z\>^\alg$ be the algebraic closure of $K\<z\>$ in $K_*$. Then $\Gamma_{L}=\Q\Gamma_{K\<z\>}=\Gamma\oplus\Q v(z)$.
Now, working in the compositional conjugate $L^{z}$, Theorem~\ref{theq} gives an element $\alpha+q\,v(z)$ in its value group ($\alpha\in\Gamma$, $q\in\Q$) such that
$$v_{P^{z}}\big(\alpha+q\,v(z)\big)\ =\  
v_{Q^{z}}\big(\alpha+q\,v(z)\big),$$
so $v_{P^{\phi}}\big(\alpha+q\,v(\phi)\big) =  
v_{Q^{\phi}}\big(\alpha+q\,v(\phi)\big)$, cofinally.
Take $a\in K^\times$ with $va=\alpha$. Then 
$$ P^{\phi}_{\times a\phi^q}\ \asymp\ Q^{\phi}_{\times a\phi^q},\qquad\text{cofinally.}$$ 
We claim that $q=0$. Towards a contradiction, assume $q\neq 0$. Renaming $P_{\times a}$ and $Q_{\times a}$ as $P$ and $Q$, respectively, we get
$P^{\times q,\phi}\ \asymp\  Q^{\times q,\phi}$, cofinally.
However, Proposition~\ref{prop:valuation of Ptimesqphi} gives
$\beta,\gamma\in\Gamma$
such that eventually $$v(P^{\times q,\phi})\ =\ \beta+dq\,v\phi, \qquad v(Q^{\times q,\phi})\ =\ \gamma+eq\,v\phi,$$ 
so $\beta-\gamma=(e-d)q\,v\phi$, cofinally, contradicting $q\ne 0$.
Hence $q=0$. Thus $P^\phi_{\times a}\asymp Q^\phi_{\times a}$, cofinally, that is,
$$v^{\ev}(P_{\times a})+\nwt(P_{\times a})v\phi\ =\  v^{\ev}(Q_{\times a})+\nwt(Q_{\times a})v\phi, \qquad\text{cofinally.}$$
Thus $v^{\ev}(P_{\times a})=v^{\ev}(Q_{\times a})$ and $\nwt(P_{\times a})=\nwt(Q_{\times a})$, so the above holds not only cofinally, but even eventually, and we are done.
\end{proof}

\noindent
The Eventual Equalizer Theorem is a key to the Newton diagram of $P$, 
but for this we also need Corollary~\ref{dn5cor} from the next subsection.

\subsection*{Transition from $\ndeg$ to $\nval$} Such a transition
is given by the next lemma. Its proof uses that $K$ is clean. 
Let $f,g$ range over $K$. For $f\preceq \fm$ we let $f_{\fm}$
be the unique element of $C$ such that $f=f_{\fm}\fm + g$ with $g\prec \fm$, so $f_{\fm}=0$ if $f\prec \fm$.

\begin{lemma}\label{transition} Suppose that $f\preceq \fm$. Then with $c:=f_{\m}$,
$$\ndeg_{\prec \m}P_{+f}\ =\ \val{(N_{P_{\times \m}})_{+c}}.$$ 
\end{lemma}
\begin{proof} For $\fn \prec \fm$ and $\fn=\mathfrak{e}\fm$ we have
$$ P_{+f,\times \fn}\ =\ P_{\times \fm,+\fm^{-1}f, \times \mathfrak{e}},$$ so
replacing $P$ by $P_{\times \fm}$ and $f$ by $\fm^{-1}f$ we arrange 
$\fm=1$. Set $Q:=P_{+f}$, so by Lemma~\ref{dn2},
$N_Q=(N_P)_{+c}$. Set $\mu:= 
\val {(N_{P})_{+c}}=\val N_Q \in \N$. 
Then for~${\fn\prec 1}$ and $i>\mu$ we have eventually 
$v(Q^\phi_i)\ge v(Q^\phi_\mu)$, hence eventually  
$$v\big((Q^\phi)_{\times \fn,i}\big)\ =\ v(Q^\phi_i)+iv\fn+o(v\fn)\ >\ 
v\big((Q^\phi)_{\times \fn,\mu}\big)\ =\ v(Q^\phi_\mu)+\mu v\fn+o(v\fn),$$ so
$\deg N_{Q_{\times \fn}}\le \mu$. Take $\phi\preceq 1$ such that 
$N_Q=D_{Q^\phi}$ and
$R_{Q^\phi} \prec^{\flat}_{\phi} Q^{\phi}$. Next, take~${\fn\prec 1}$ with $\psi^\phi(v\fn)=v\fn$,
so $v^\flat_{\phi}(\fn)=0$. Then by Lemma~\ref{dn2a},
$$  N_{Q_{\times \fn}}\ =\ N_{Q_{\times \fn}^\phi}\in C^\times \cdot Y^\mu,$$ 
in particular, $\deg N_{Q_{\times \fn}}= \mu$. 
\end{proof}

\begin{cor}\label{dn5cor} $\ndeg_{\prec \fm} P = \nval P_{\times \fm}$. 
\end{cor}

\subsection*{The Newton diagram of $P$} {\em In this subsection we assume
that $\Gamma$ is divisible.}\/ We begin with a reformulation of the Eventual 
Equalizer Theorem:

\begin{cor}\label{monevq} Assume $P,Q\in K\{Y\}^{\ne}$ are
homogeneous of different degrees. Then there is a unique monomial
$\fm$ such that $N_{(P+Q)_{\times \fm}}$ is not homogeneous.
\end{cor}
\begin{proof} There is at most one such monomial by Lemma~\ref{npq}.
The proof of that lemma shows that if $\fm$ is such that
eventually $P_{\times \fm}^{\phi}\asymp Q_{\times \fm}^{\phi}$, then
 $N_{(P+Q)_{\times \fm}}= N_{P_{\times \fm}}+ N_{Q_{\times \fm}}$ is not homogeneous.
Thus it remains to appeal to the Eventual Equalizer Theorem for the
existence of such an $\fm$. 
\end{proof}

\noindent
For homogeneous
$P,Q\in K\{Y\}^{\ne}$ of different degrees we call the unique $\fm$
as in Corollary~\ref{monevq} the {\bf eventual equalizer} for $P$,~$Q$
and denote it by $\frak{e}(P,Q)$. 

\index{eventual!equalizer}
\index{equalizer!eventual}
\nomenclature[X]{$\frak{e}(P,Q)$}{eventual equalizer for $P$,~$Q$}

{\sloppy
We now focus on $P$, and let $J$ be the finite nonempty 
set of $j\in \N$ with~${P_j\ne 0}$.
Thus $\ndeg P_{\times \fm}\in J$ for all $\fm$. For distinct $i,j\in J$
we let $\frak{e}(P,i,j)$ be the eventual equalizer for $P_i, P_j$.
Therefore any algebraic starting monomial for $P$ is of the form~$\frak{e}(P,i,j)$ with distinct $i,j\in J$. Note that 
$\fM$ is (totally) ordered by $\preceq$. 
Below we fix a $\preceq$-closed set
$\E\subseteq K^\times$.
}

\begin{prop}\label{newtpolygon} There are $i_0, \dots, i_n\in J$ and 
eventual equalizers
$$  \frak{e}(P,i_0,i_1)\ \prec\  \frak{e}(P,i_1,i_2)\ \prec\ \cdots\ \prec\
\frak{e}(P,i_{n-1},i_n)$$
with $\val P= i_0 < \dots < i_n = \ndeg_{\E}P$, such that: \begin{enumerate}
\item[\textup{(i)}] the 
algebraic starting monomials for $P$ in $\E$ are the $\frak{e}(P, i_m, i_{m+1})$, $m<n$;
\item[\textup{(ii)}] for $\fm=\frak{e}(P, i_m, i_{m+1})$, $m<n$, we have
$\val N_{P_{\times \fm}}=i_m$, $\deg N_{P_{\times \fm}}=i_{m+1}$.
\end{enumerate}
\end{prop}    
\begin{proof} Let $i$,~$j$ range over $J$, and set $d:= \ndeg_{\E}P$, 
so  $\val P \le d\le \deg P$. We proceed by induction on
$d-\val P$. If $d=\val P$, then all $N_{P_{\times \fm}}$ with $\fm\in \E$ are 
homogeneous of degree $d$, and so there are no
algebraic starting monomials of~$P$ in~$\E$. Assume that $d> \val P$, and take $i<d$
such that the eventual equalizer $\frak{e}=\frak{e}(P,i,d)$ is maximal 
with respect to $\preceq$.
We claim that then $\frak{e}\in \E$. 
Towards proving this, let~${\fn\prec \frak{e}}$.  
Since $P^{\phi}_{i, \times \frak{e}}\asymp P^{\phi}_{d,\times \frak{e}}$, eventually,
this gives $P^{\phi}_{i, \times \fn}\succ P^{\phi}_{d,\times \fn}$, eventually, and thus $\ndeg P_{\times \fn} \ne d$. Taking $\fn\in \E$
such that $\ndeg P_{\times \fn}=d$ therefore gives $\frak{e}\preceq \fn$, and 
so we get $\frak{e}\in \E$, as claimed. Also $\deg N_{P_{\times \frak{e}}}=d$:
otherwise $\deg N_{P_{\times \frak{e}}}=j < d$, so 
$P^{\phi}_{j,\times \frak{e}}\succ P^{\phi}_{d,\times \frak{e}}$ eventually, hence
$\frak{e}(P,j,d) \succ \frak{e}$, contradicting the maximality of~$\frak{e}$. 
Since $i< d$ and 
$\big(N_{P_{\times \frak{e}}}\big)_i= N_{P_{i,\times \frak{e}}}\ne 0$, it follows
that $\frak{e}$ is an algebraic starting monomial for $P$. 

\claim{$\frak{e}$ is the largest algebraic starting monomial for $P$ in
$\E$.}

\noindent
To see this, suppose towards a contradiction that $\fn\succ \frak{e}$ is an algebraic starting monomial for $P$ in $\E$. Then 
$\ndeg P_{\times \frak{e}}\le \ndeg P_{\times \fn}$, so $\ndeg P_{\times \fn}=d$ by the 
maximality property of $\E$. So $\fn=\frak{e}(P,j,d)$ with $j< d$, but this
contradicts the maximality of $\frak{e}$, and so proves the claim. 

By decreasing $i$ if necessary
we arrange $i=\val N_{P_{\times \frak{e}}}$. Then 
$\ndeg_{\prec\frak{e}}P= i$ by Corollary~\ref{dn5cor}. It remains to apply the inductive assumption with $\E$ replaced by the set 
$\{g\in K^\times:\  g\prec \frak{e}\}$. 
\end{proof}

\noindent
Let $(i_0,\dots, i_n)$ be as in Proposition~\ref{newtpolygon}. This tuple
is uniquely determined by the data $K$,~$P$,~$\E$. If $\val P = \ndeg_{\E}P$, 
then $n=0$ and this tuple is just $(\val P)$. To simplify notation, set 
$\frak{e}_m:= \frak{e}(P, i_{m-1},i_m)$ for $1\le m \le n$. 
We now have a complete
description of the behavior of $\nval P_{\times g} $ and $\ndeg P_{\times g}$ for
$g\in \E$:

\begin{cor}\label{newtpoldes} Assume $\val P \ne \ndeg_{\E}P$, so $n\ge 1$. Let $g$ range over $\E$. Then $\nval P_{\times g}$ and $\ndeg P_{\times g}$ lie
in the set $\{i_0,\dots, i_n\}$, and we have:
\begin{align*} \nval P_{\times g}\ =\ i_0\ &\Longleftrightarrow\ g\preceq 
\frak{e}_1;\\ 
 \ndeg P_{\times g}\ =\ i_0\ &\Longleftrightarrow\ g\prec \frak{e}_1;\\
\nval P_{\times g}\ =\ i_m\ &\Longleftrightarrow\ 
\frak{e}_m \prec g \preceq \frak{e}_{m+1},  \qquad (1\le m < n);\\
      \ndeg P_{\times g}\ =\ i_m\  &\Longleftrightarrow\ 
\frak{e}_m \preceq g \prec \frak{e}_{m+1}, \qquad (1\le m < n);\\
\nval P_{\times g}\ =\ i_n\ &\Longleftrightarrow\ \frak{e}_n\prec g;\\
\ndeg P_{\times g}\ =\ i_n\ &\Longleftrightarrow\  \frak{e}_n \preceq g.
\end{align*}
\end{cor}
\begin{proof} Let $1\le m < n$, and set $\frak{e}:=\frak{e}(P,i_{m-1}, i_m)$
and $\frak{e}^*:= \frak{e}(P,i_{m}, i_{m+1})$. For $\frak{e} \prec g \prec \frak{e}^*$ we obtain from Proposition~\ref{newtpolygon} and Corollary~\ref{dn1cor,new,Newton} that
$$i_m\  =\ \ndeg P_{\times \frak{e}}\ \le\ \nval P_{\times g}\  =\ \ndeg P_{\times g}\ \le\ \nval P_{\times \frak{e}^*}\  =\ i_m,$$
where $\nval P_{\times g}= \ndeg P_{\times g}$ because $P$ has no algebraic
starting monomial $\fm$ with
$\frak{e} \prec \fm \prec \frak{e}^*$. Thus $i_m= \nval P_{\times g}\  =\ \ndeg P_{\times g}$, from which we obtain the ``middle'' equivalences. The ``end''
equivalences are derived in the same way.
\end{proof}

\noindent
Of course, if  $\val P = \ndeg_{\E}P$, then $\nval P_{\times g}=\ndeg P_{\times g}=\val P$ for all $g\in \E$, and $P$ has no algebraic starting monomials in $\E$. 

Suppose now that $\val P\neq\ndeg_{\E} P$, and let $n$ be as above.
So $n$ is by Proposition~\ref{newtpolygon} the number of algebraic starting monomials in $\E$ for $P$, and 
$1\leq n\leq\ndeg_{\E} P - \val P$.  Thus we can take algebraic
approximate zeros $y_1,\dots, y_N$ of~$P$ in~$\E$, $N\in \N$, such that each algebraic approximate zero $y$ of $P$ in $\E$ satisfies
 $y\sim y_i$
for exactly one $i\in \{1,\dots, N\}$.  
Combining Proposition~\ref{newtpolygon} with Lemma~\ref{lem:sum of alg mult} yields:

\begin{cor}\label{cor:sum of alg mult} 
For $i=1,\dots, N$, let $\mu_i$ be the algebraic multiplicity of the approximate zero
$y_i$ of $P$, and set $\mu:=\mu_1+\cdots+\mu_N$.
Then $\mu\leq \ndeg_\E P - \val P$. If $C$ is algebraically closed, then 
$\mu = \ndeg_\E P - \val P$. If $C$ is real closed, then
$\mu \equiv \ndeg_\E P - \val P\bmod 2$.
\end{cor}

\subsection*{Notes and comments}
Corollary~\ref{monevq} and Proposition~\ref{newtpolygon} for $K=\T_{\operatorname{g}}$ are Proposition~8.14(c) 
and Proposition~8.17, respectively,
in \cite{JvdH}. Corollary~\ref{cor:sum of alg mult} for $K=\T_{\operatorname{g}}[\imag]$ is \cite[Exercise~8.14]{JvdH}.

\section{Further Consequences of $\upo$-Freeness}\label{consupofree}

\noindent
{\em In this section $E$ is an ungrounded $H$-asymptotic field with $\Gamma_E\ne \{0\}$}. Thus $E$ is 
pre-$\d$-valued but not necessarily 
$\d$-valued. Also in contrast to $K$, we do not assume the derivation of $E$ small, or that the valued field $E$ has 
a monomial group. 
Let~$(\Gamma_E, \psi_E)$ be the asymptotic couple of $E$, so $\Psi_E\ne \emptyset$ and
$\Psi_E$ has no largest element. We fix a
logarithmic sequence for $E$ and corresponding sequences~$(\upl_{\rho})$ and $(\upo_{\rho})$ in~$E$ as in
Section~\ref{sec:special cuts} (with $E$ instead of~$K$).

 \begin{theorem}\label{upoalgebraic} Suppose that $E$ is $\upo$-free, and that $F$ is a
pre-$\d$-valued field extension of $E$ of $H$-type which is 
$\d$-algebraic over $E$. Then: \begin{enumerate}
\item[\textup{(i)}] there is no $y\in F$ with 
      $0 < |vy| < \Gamma_E^{>}$;  
\item[\textup{(ii)}] there is no $z\in F$ with $\Psi_E < vz < (\Gamma_E^{>})'$; and
\item[\textup{(iii)}] $F$ is $\upo$-free.
\end{enumerate}
\end{theorem}

\noindent
In particular, if $E$ is $\upo$-free, then so is its
differential-valued hull 
$\operatorname{dv}(E)$. (See  Section~\ref{sec:dv(K)} for the basic facts
on differential-valued hulls.)

\begin{proof} If (i) holds, then $\Gamma_E^{<}$ is cofinal in
$\Gamma_F^{<}$, and thus the logarithmic sequence for $E$ can also serve 
as such for
$F$. This remark helps to justify the reduction steps we make 
in the proof, especially in connection with (iii). 

First we take some active 
$a\in E$ and replace $E$ and $F$ 
by their compositional conjugates 
$E^{a}$ and $F^a$. In this way we arrange that the derivations of $E$ and $F$ 
are small.
Next, replacing $F$ by  
$\operatorname{dv}(F)$ we arrange that $F$ is $\d$-valued of $H$-type.
Finally, replacing $E$ and $F$ by their algebraic closures,
we can also assume that $E$ and $F$ are algebraically closed.
This last replacement doesn't change $\Psi_E$, so we can keep the same
logarithmic sequence for $E$.

Lemma~\ref{divmon} gives a monomial group $\frak{M}$ for $E$. Since $E$ is
$\upo$-free, it has asymptotic integration, so
$K:=\operatorname{dv}(E)$ is an immediate extension of $E$. We equip~$K$ with
the monomial group $\frak{M}$, so
$K$ satisfies our standing assumptions at the beginning of this chapter,
and $\Gamma:= \Gamma_K=\Gamma_E$ is divisible. 
Take the unique valued differential field embedding of
$K$ into $F$ that is the identity on $E$, and identify $K$ with a
valued differential subfield of $F$ via this embedding. Then
$F$ is a $\d$-valued field extension of $H$-type of~$K$. Now let
$y\in F$ and take $P\in E\{Y\}^{\ne}$ with $P(y)=0$. By the remark
following the proof of Proposition~\ref{prop:nabla1,nabla2} we have
$N_P\in C[Y](Y')^{\N}$, and thus by Lemma~\ref{divclean} and 
Corollary~\ref{cor:sign and valuation, 1} we cannot have
$0<|vy|< \Gamma^{>}$. This proves (i).

For (ii), suppose towards a contradiction that $z\in F$ and $\Psi < vz < (\Gamma^{>})'$. Corollary~\ref{newcor} provides a $\d$-algebraic 
$\d$-valued field
extension $L$ of $H$-type of $F$ with an element
 $y\in L^\times$ such that  
$y^\dagger=z$. If $y\asymp 1$, then $y=c+y_1$ with~$c\in C_L^\times$ and
$v(y_1)\ge \gamma$ for some $\gamma\in \Gamma^{>}$, by (i) applied to $L$
in the role of~$F$, so
${v(z)=v(y_1') \ge \gamma'}$, contradicting the assumption on~$z$. 
Thus $y\not\asymp 1$, but then the assumption on~$z$ gives $0 < |vy| < \Gamma^{>}$, contradicting (i). This proves (ii). Therefore,~$F$ has asymptotic
integration. For (iii), we first show that
$F$ is $\upl$-free. Suppose it is not. Since we arranged~$\Gamma_F$ to
be divisible, $F$ has a gap creator $s$ by Lemma~\ref{cg2}. Proposition~\ref{newprop}
gives a $\d$-valued field extension $F(f)$ of $H$-type of $F$ 
with $f\ne 0$ and 
$f^\dagger=s$, but then $\Psi < vf < (\Gamma^{>})'$ by the remark preceding Corollary~\ref{psacgap}, and this contradicts (ii). Thus $F$ is indeed $\upl$-free,
and this holds not just for the present~$F$, which is algebraically closed, etcetera, but for any $F$ satisfying the
assumptions in the theorem.
To show that our present~$F$ is even $\upo$-free, assume towards a contradiction that 
$\upo_{\rho} \leadsto \upo\in F$.
Then~$F$ has by Corollary~\ref{uplupo} an immediate asymptotic extension
that is $\d$-algebraic over~$F$ (and thus over $E$) but not
$\upl$-free. This contradicts what we just proved, namely, that any $F$ as in the hypothesis of the 
theorem is $\upl$-free. \end{proof}

\begin{cor}\label{upoliou}
Let $E$ be an $\upo$-free $H$-field. Then $E$ has exactly one Liouville closure, up to isomorphism over $E$.
\end{cor}
\begin{proof} Every Liouville $H$-field extension of $E$ is
$\d$-algebraic over $E$, and thus no such extension has a gap, by 
Theorem~\ref{upoalgebraic}. This gives the desired conclusion
in view of the remarks following Theorem~\ref{thm:Liouville closures}. 
\end{proof}

\begin{cor}\label{differentialtranscendence} Suppose $E$ is $\upo$-free. 
Then the pc-sequences $(\upl_{\rho})$ and $(\upo_{\rho})$ are of 
$\d$-transcendental type over $E$.
\end{cor}
\begin{proof}  
If $(\upl_{\rho})$ is 
of $\d$-algebraic type over $E$, then $(\upl_{\rho})$
pseudoconverges in some immediate asymptotic extension of $E$ 
that is $\d$-algebraic over~$E$, by Corollary~\ref{zmindifpol} and
Lemma~\ref{zda, newton},  but this is impossible by Theorem~\ref{upoalgebraic}. 
Likewise, $(\upo_{\rho})$ is of $\d$-transcendental type over $E$. 
\end{proof}

\noindent
We now combine this with Lemma~\ref{lem:eval at uporho} to get: 

\begin{prop}\label{cor:adding upo preserves upl-free}
Suppose $E$ is $\upo$-free, and $\upo$ is an element in an $H$-asymp\-totic field extension of $E$ such that $\upo_\rho\leadsto\upo$. Then  
$E\<\upo\>$ is $\upl$-free.
\end{prop}
\begin{proof} Since $(\upo_{\rho})$ is of $\d$-transcendental type over $E$, the asymptotic extension $F:=E\<\upo\>$ of $E$ is immediate
and $\upo$ is $\d$-transcendental over $E$. Moreover, $(\upl_{\rho})$ is a divergent pc-sequence in $E$.
Towards a contradiction, suppose $P,Q\in E\{Y\}^{\neq}$ are such that 
$$\upl_{\rho} \leadsto \upl:=\frac{P(\upo)}{Q(\upo)}\in F\setminus E.$$  
Put $R:=P/Q\in E\<Y\>\setminus E$; so $R(\upo)=\upl$.  Set
$$G(Y)\ :=\ P(Y)Q(\upo)-Q(Y)P(\upo)\in F\{Y\}.$$
From $R\notin F$ we get $G\ne 0$. In view of $G(\upo)=0$, this  gives 
$G\notin F$. Lemma~\ref{lem:eval at uporho} then yields $G(\upo_{\rho}) \leadsto 0$. Also $Q(\upo_{\rho}) \not\leadsto 0$,
and thus $Q(\upo_{\rho})\asymp Q(\upo)$, eventually, by Lemma~\ref{lem:eval at uporho}. 
Take $\rho_0$ such that $Q(\upo_\rho)\neq 0$ for $\rho>\rho_0$. Then
$$Q(\upo)Q(\upo_{\rho})\big(R(\upo_{\rho})-\upl\big)\ =\ G(\upo_{\rho})\qquad(\rho>\rho_0),$$
so for $\alpha:= -2v\big(Q(\omega)\big)\in \Gamma_E$ we get 
\begin{equation}\label{eq:approx upl by R(uporho)}
v\big(R(\upo_{\rho})-\upl\big)\ =\ \alpha + v\big(G(\upo_{\rho})\big),\ \text{ eventually.}
\end{equation}
In particular, $R(\upo_{\rho}) \leadsto \upl$ ($\rho> \rho_0$).
If $R(\upo_{\rho}) \leadsto b\in E$ ($\rho>\rho_0$), then
$$P(\upo_{\rho})-bQ(\upo_{\rho})\ =\ Q(\upo_{\rho})\big(R(\upo_{\rho}) - b \big) \leadsto 0,$$
contradicting that $(\upo_{\rho})$ is of $\d$-transcendental
type over $E$. 
%and so $P-bQ=0$, that is, $R=b\in E$, a contradiction.
Thus the pc-se\-quence $\big(R(\upo_\rho)\big)_{\rho > \rho_0}$ in $E$ is divergent. Then by Corollary~\ref{eqpccor}, the pc-sequences $(\upl_\rho)$ and~$\big(R(\upo_\rho)\big)_{\rho > \rho_0}$ in $E$ are equivalent, and so
the latter has width 
$\{\gamma\in\Gamma_\infty:\gamma>\Psi\}$.
Thus for any big enough $\rho_0$ 
the set $\big\{v\big(R(\upo_{\rho})-\upl\big):\ \rho> \rho_0\big\}$
is a cofinal subset of $\Psi_E^{\downarrow}$. 
However, taking $a\in E^\times$ with $va=\alpha$ and setting 
$H:= aG\in F\{Y\}$, it follows from~\eqref{eq:approx upl by R(uporho)} that for big enough 
$\rho_0$ this set equals
$\big\{v\big(H(\upo_{\rho})\big):\ \rho> \rho_0\big\}$. 
As $H(\upo)=0$, this contradicts Lemma~\ref{lem:eval at uporho}. 
\end{proof}

\noindent
Towards extending Lemmas~\ref{lem:infinitesimal cut} and~\ref{lem:infinitesimal cut, uniqueness} to the present setting, we first note:

\begin{lemma}\label{Eval} Suppose $E$ is $\upo$-free. Let $P\in E\{Y\}^{\ne}$. Then there exists an active~$a$ in $E$ such that for every pre-$\d$-valued field extension $F$ of $E$ of $H$-type,
$$y\in F,\ 1\prec y \asymp^{\flat}_{a} 1\ \Longrightarrow\ v\big(P(y)\big)\ =\ v^{\ev}(P) + \ndeg(P)vy + \nwt(P)\psi_F(vy).$$
\end{lemma}
\begin{proof} In view of Lemma~\ref{vepa} we can arrange by compositional conjugation that~$E$ has small derivation. Next, replacing $E$ by $\operatorname{dv}(E)$, and any $F$ by $\operatorname{dv}(F)$, we also
arrange that $E$ is differential-valued, and then passing to algebraic closures we can even assume that
$E$ (and any $F$) is algebraically closed. Then $E$ has a monomial group~$\mathfrak{M}$, so $E$ satisfies the assumptions made on
$K$ in this chapter. As~$E$ is $\upo$-free, we can apply Corollary~\ref{cor:sign and valuation, 1} to get the desired result. 
\end{proof}

\begin{lemma}\label{Evalo} Suppose $E$ is an $\upo$-free $H$-field. Let 
$P\in E\{Y\}^{\ne}$. Then there are active $a$ in $E$ and
$\sigma\in \{-1,+1\}$ such that for every $H$-field extension $F$ of $E$,
$$ y\in F^{>},\ 1\prec y \asymp^{\flat}_{a} 1\ \Longrightarrow\ \operatorname{sign} P(y)\ =\ \sigma.$$
\end{lemma}
\begin{proof} As Lemma~\ref{Eval} follows from the first part of Corollary~\ref{cor:sign and valuation, 1}, this follows 
from its second part. Instead of algebraic closures, take real closures. 
\end{proof}

\begin{prop}\label{Evalcor} Suppose $E$ is $\upo$-free, $F$ is a pre-$\d$-valued field extension of $E$ of $H$-type, and $y\in F$ satisfies
$\Gamma_E^{<} < vy < 0$. Then for all $P\in E\{Y\}^{\ne}$,
$$v\big(P(y)\big)\ =\ v^{\ev}(P) + \ndeg(P)vy + \nwt(P)\psi_F(vy).$$
Moreover, $\Gamma_{E\<y\>}=\Gamma_E \oplus \Z vy \oplus \Z\psi_F(vy)$
\textup{(}internal direct sum\textup{)} and $\max \Psi_{E\<y\>}=\psi_F(vy)$.
For any $y^*$ in any pre-$\d$-valued extension
$F^*$ of $E$ of $H$-type satisfying $\Gamma_E^{<} < vy^* < 0$, there is a
unique valued differential field embedding $E\<y\> \to F^*$ over~$E$ sending $y$ to $y^*$.
\end{prop}
\begin{proof} The first part is immediate from Lemma~\ref{Eval}, and
the proof of the rest is like that of 
Lemmas~\ref{lem:infinitesimal cut} and 
~\ref{lem:infinitesimal cut, uniqueness}.
\end{proof}

\begin{cor}\label{Evalocor} Let $E$ be an $\upo$-free $H$-field.
Let $F$ and $F^*$ be $H$-field extensions of $E$ with elements
$y\in F^{>}$ and $y^*\in F^{* >}$ such that 
$\Gamma_E^{<} < vy < 0$ and $\Gamma_E^{<} < vy^* < 0$. Then there is a
unique pre-$H$-field embedding $E\<y\> \to F^*$ over $E$ sending $y$ to $y^*$.
\end{cor}
\begin{proof} This follows from Proposition~\ref{Evalcor} in view of
Lemma~\ref{Evalo}.
\end{proof}

\noindent
Next we extend the Eventual Equalizer Theorem to the present setting:

\begin{cor}\label{hevq} Assume $E$ is $\upo$-free, $P,Q\in E\{Y\}^{\ne}$ are 
homogeneous of degrees $d> e$, and $(d-e)\Gamma_E=\Gamma_E$.
Then for some $a\in E^{\times}$ and active $f_0$ in $E$,
$$ P_{\times a}^{f}\ \asymp\ Q_{\times a}^{f}\ \text{ for all active $f\preceq f_0$ 
in $E$.}$$
\end{cor} 
\begin{proof} Note: for active $f,g\in E$ with $f\asymp g$ and 
$a,b\in E^\times$ with $a\asymp b$ we have
$$P^f_{\times a}\ \asymp\ P^f_{\times b}\ \asymp\ P^g_{\times b},\ 
\text{ and likewise }Q^f_{\times a}\ \asymp\ Q^g_{\times b}.$$ 
An initial compositional conjugation by an active element of
$E$ arranges that the derivation of $E$ is small.
Since $\operatorname{dv}(E)$ is
still $\upo$-free, by Theorem~\ref{upoalgebraic}, 
and is an immediate extension of $E$, we can also arrange that  
$E$ is $\d$-valued of $H$-type. Then the algebraic closure $K$ of $E$ 
is $\d$-valued of $H$-type and $\upo$-free, and has a monomial group 
$\frak{M}$, and so~$K$ with $\frak{M}$ satisfies the conditions imposed
at the beginning of this chapter, and in addition 
$\Gamma:= \Gamma_K$
is divisible. By the Eventual Equalizer Theorem we can take 
$a\in K^\times$ 
and active $\phi_0$ in $K$ such that for all 
active~$\phi\preceq\phi_0$ we have
$P_{\times a}^{\phi}\ \asymp\ Q_{\times a}^{\phi}$, where, as always in this chapter,
$\phi$ ranges over~$\frak{M}$.
Now let $f\in E$ be active with 
$f\preceq \phi_0$. Taking~$\phi$ such that $\phi\asymp f$ we note that
$P_{\times a}^f\asymp Q_{\times a}^f$. Since $P^f, Q^f\in E^f\{Y\}^{\ne}$ are 
homogeneous of degrees $d>e$, and $(d-e)\Gamma_E=\Gamma_E$,
we have by the Equalizer Theorem of Chapter~\ref{ch:valueddifferential} 
a unique $\alpha\in \Gamma_E$ such that 
$v_{P^f}(\alpha)= v_{Q^f}(\alpha)$. It follows that $\alpha=va$. Thus
with $a$ replaced by any element of $E$ with valuation $\alpha$, and
with $f_0$ any active element of $E$ with $f_0\preceq \phi_0$, we have 
the desired conclusion.
\end{proof}

\noindent
For $d=1$ and $e=0$ this means:

\begin{cor}\label{hevqone} If $E$ is $\upo$-free and 
$P\in E\{Y\}^{\ne}$ is
homogeneous of degree $1$, then for each $a\in E^\times$ there are
elements $g\in E^{\times}$ and active $f_0$ in $E$ such that:
$$ P_{\times g}^{f}\ \asymp\ a\ \text{ for all active $f\preceq f_0$ 
in $E$.}$$
\end{cor}

\noindent
The results in Section~\ref{The Shape of the Newton Polynomial} yield in the present setting initially the following:

\begin{lemma}\label{lem:huposhape}
Assume $E$ is $\upo$-free and $P\in E\{Y\}^{\ne}$ has Newton weight 
$w$. Then there are $A\in E[Y]^{\ne}$ with $A\asymp 1$, $w+\val A=\nval P$, and $w+\deg A=\ndeg P$, an active
$e$ in $E$, an $a\in E^{\times}$, and an $R\in E^{e}\{Y\}$, such that
$$P^{e}\ =\ a\cdot A\cdot (Y')^w + R, \qquad a\asymp P^{e}, \quad R\prec^{\flat}_{e} P^{e}.$$  
\end{lemma}
\begin{proof}
We have $\order(R)\le \order(P)$ and
$\deg R \le \deg P$ for $R$ as above, so we are dealing with
an elementary statement about $E$ and the tuple of
coefficients of $P$. Hence we may pass to an elementary extension of $E$ and
arrange in this way that $E$ is $\aleph_1$-saturated, and 
thus has a monomial group $\frak{M}$, by Lemma~\ref{xs1}. 
Moreover, an initial 
compositional conjugation 
by an active element of $E$ arranges that the derivation of $E$ is small.
Let $K:=\operatorname{dv}(E)$, so $K$ is an immediate extension of $E$
and~$K$ is an $\upo$-free $\d$-valued field of $H$-type with monomial group
$\frak{M}\subseteq E^\times$. This leads to the Newton polynomial 
$N_P\in C[Y](Y')^w$, so
$N_P= B \cdot (Y')^w$, with $B\in C[Y]$, so $w + \deg B = \ndeg P$.
Take $A\in E[Y]$ with $A\sim B$, $\val A = \val B$, $\deg A = \deg B$. Take
$\gamma>0$ in $\Gamma= \Gamma_E$ such that $v(A-B)>\gamma$. 
Take active $e\preceq 1$ in $\frak{M}$ such that  
$$\Gamma^{\flat}_{e}\ <\ \gamma, \quad P^{e}\ =\ \fd_{P^{e}}N_P + R_{P^{e}}, \qquad  R_{P^{e}}\prec^{\flat}_{e} P^{e}.$$
Since $
N_P=B\cdot (Y')^w=A\cdot (Y')^w+G$ with $vG>\gamma$, we obtain 
$$P^e\ =\   
\fd_{P^{e}}\cdot A\cdot (Y')^w + \fd_{P^{e}}G +  R_{P^{e}},$$
so the desired conclusion holds with 
$a:=  \fd_{P^{e}}$ and $R:=  \fd_{P^{e}}G +  R_{P^{e}}$.  
\end{proof}

\noindent
In the next two corollaries of Lemma~\ref{lem:huposhape} we let
$F$ range over $H$-asymptotic field extensions
of $E$. Recall that for such $F$, if $e\in E$ is active in $E$, then 
$e$ remains active in $F$, and if $f\in F$ is active, then $F^f$ has small derivation.

\begin{cor}\label{huposhape}
Assume $E$ is $\upo$-free and $P\in E\{Y\}^{\ne}$ has Newton weight 
$w$. Let $A\in E[Y]^{\ne}$ and $e\in E$ be as in Lemma~\ref{lem:huposhape}. Then for all $F$ and active $f\preceq e$ in $F$
there  are $a_f\in F^\times$ and $R_f\in F^f\{Y\}$ such that
$$P^f\ =\ a_f\cdot A\cdot (Y')^w + R_f, \qquad a_f\asymp P^f, \quad R_f \prec^{\flat}_{e} P^f.$$
\end{cor}
\begin{proof} To simplify notation we
replace $E$, $P$ by $E^{e}$, $P^{e}$ (and each $F$ by $F^e$),
and rename accordingly, so that $e=1$. With 
$a$ and $R$ as in Lemma~\ref{lem:huposhape} we have
$$P\ =\ a\cdot A\cdot (Y')^w+R, \quad a\asymp P,\quad A\asymp 1, \quad R\prec^\flat P  \quad \text{(in $E\{Y\}$).}$$
Note that $R\prec^\flat P$ remains true in $F\{Y\}$.
An easy computation as in the proof of Lemma~\ref{cle} now shows that if $f\preceq 1$ is active in $F$, then
$$P^f\ =\ af^w\cdot A\cdot (Y')^w+R^f,\qquad R^f\prec^\flat af^w\asymp^\flat P^f,$$
so $a_f:=af^w$ and $R_f:=R^f$ works.
\end{proof}
 
\begin{cor}\label{cor:Newton pol under ext}
Suppose $E$ is $\upo$-free. Let $P\in E\{Y\}^{\ne}$. Then there is an active $e\in E$ such that for all $F$ and active $f\preceq e$ in $F$,
$$\ddeg P^f\ =\ \ndeg P, \quad \dval P^f\ =\ \nval P, \quad 
\dwt P^f\ =\ \nwt P.$$
\end{cor}

\noindent
Independent of whether $E$ is $\upo$-free, when we have a decomposition of $P$ as in Lem\-ma~\ref{lem:huposhape}, the following is relevant:

\begin{lemma}\label{decndegnval} Assume the derivation of $E$ is small, and 
$P\in E\{Y\}^{\ne}$, 
$$ P\ =\ a\cdot A \cdot (Y')^w + R, \qquad a\in E^\times,\ A\in E[Y]^{\ne},\ w\in \N,\ R\in E\{Y\},$$
with $a\asymp P$, $A\asymp 1$, $R \prec^{\flat} P$.
Let $\overline{A}$ be the image of $A$ in $\k_E[Y]$. Then $\nval P= \val{\overline {A}} +w$, $\ndeg P= \deg{\overline{A}} + w$,
and there exists $\gamma\in \Gamma_E^{>}$ such that
for all $g\in E$,
\begin{align*} 0<vg< \gamma\ &\Longrightarrow\  \nval(P)\ =\ \nval(P_{\times g})\ =\ \ndeg(P_{\times g}),\\
-\gamma < vg < 0\  &\Longrightarrow\  \ndeg(P)\ =\ \ndeg(P_{\times g})\ =\ \nval(P_{\times g}).
\end{align*}
\end{lemma}
\begin{proof} Let $f\prec 1$ be active in $E$. Then $vf\in \Gamma^{\flat}$
and 
$$   P^f\ =\  af^w\cdot A \cdot (Y')^w + R^f$$
with $v(R^f) \ge v(R) > v(af^w)$, so
$\dval P^f= \val{\overline{A}} + w$ and $\ddeg P^f=\deg{\overline{A}} + w$. This proves the claim about $\nval P$ and $\ndeg P$. Take
$B, G\in E[Y]$ such that $A=B+G$, all nonzero coefficients of $B$ are $\asymp 1$,
and $G\prec 1$, so $B = b_mY^m + \cdots + b_n Y^n$ with $m=\val{\overline{A}}$ and $n=\deg{\overline{A}}$, $b_m,\dots, b_n\in E$,
$b_m, b_n\ne 0$. 
Next, let $g\in E^\times$ be such that $g\not\asymp 1$. 
Then in $E^f\{Y\}$ we have
$$P^f_{\times g}\ =\ af^w \cdot A_{\times g} \cdot (Y'_{\times g})^w + R^f_{\times g}\ =\ af^w\cdot B_{\times g}\cdot (Y'_{\times g})^w + af^w\cdot G_{\times g} \cdot (Y'_{\times g})^w + R^f_{\times g}.$$
If $f \prec g^{\dagger}$, then
$Y'_{\times g}= f^{-1}g'Y + gY' \sim f^{-1}g'Y$ in $E^f\{Y\}$. 
To derive the first displayed implication, assume $0 < v(g) \in \Gamma^{\flat}$ and $m\,v(g) < v(G)$. Then $v(g')\in \Gamma^{\flat}$
and $B_{\times g}\sim b_mg^m Y^m\succ G_{\times g}$. Therefore, eventually with respect to $f$,
$$ P^f_{\times g}\ \sim \ ab_m g^m(g')^w\cdot Y^{m+w}, $$
and so $\nval(P_{\times g}) = \ndeg(P_{\times g})=m+w=\nval(P)$.
For the second implication, let $g$ satisfy $vg< 0$, 
$vg\in \Gamma^{\flat}$,
and $v(g^n)=nv(g) < v(G_{\times g})$. Then $v(g')\in \Gamma^{\flat}$, $B_{\times g}\sim b_ng^nY^n\succ G_{\times g}$, so eventually with
respect to $f$,
$$  P^f_{\times g}\ \sim ab_ng^n(g')^w\cdot Y^{n+w},$$
and thus $\ndeg(P_{\times g})=\nval(P_{\times g})=n+w=\ndeg(P)$. 
\end{proof}

\begin{cor}\label{decupondegnval} Suppose $E$ is $\upo$-free, and $P\in E\{Y\}^{\ne}$. Then there exists $\gamma\in \Gamma_E^{>}$ such that
for all $g\in E$,
\begin{align*} 0<vg< \gamma\ &\Longrightarrow\  \nval(P)\ =\ \nval(P_{\times g})\ =\ \ndeg(P_{\times g}),\\
-\gamma < vg < 0\  &\Longrightarrow\  \ndeg(P)\ =\ \ndeg(P_{\times g})\ =\ \nval(P_{\times g}).
\end{align*}
\end{cor}
\begin{proof} Use Lemma~\ref{huposhape} and compositional conjugation
by some active element of~$E$
to arrange that $E$ has small derivation and $P$ has a decomposition as in
the hypothesis of Lemma~\ref{decndegnval}. Now apply Lemma~\ref{decndegnval}.
\end{proof}

\noindent
Note also that $\deg A = \ddeg A$ for $A$ as in Lemma~\ref{lem:huposhape}.
The device in the proof of passing to an elementary extension with a 
monomial group and then to the differential-valued hull can also
be used in other situations. For example, the results on the 
Newton diagram of $P$ in Section~\ref{eveqth} extend in 
a similar way to the present setting. Here we state what is
needed later:

\begin{cor}\label{newtpolext} Assume $E$ is $\upo$-free, $\Gamma_E$
is divisible, and $P\in E\{Y\}^{\ne}$ is not homogeneous. 
Then there are $i_0 < \cdots < i_n$ in $\{i\in \N:\ P_i\ne 0\}$ with $i_0=\val P$ and $i_n=\deg P$
and elements
$\frak{e}_1  \prec  \cdots \prec \frak{e}_n$ in $E^\times$ such that for all
$g\in E^\times$:
\begin{align*} \nval P_{\times g}\ =\ i_0\ &\Longleftrightarrow\ g\preceq 
\frak{e}_1;\\ 
 \ndeg P_{\times g}\ =\ i_0\ &\Longleftrightarrow\ g\prec \frak{e}_1;\\
\nval P_{\times g}\ =\ i_m\ &\Longleftrightarrow\ 
\frak{e}_m \prec g \preceq \frak{e}_{m+1},  \qquad (1\le m < n);\\
      \ndeg P_{\times g}\ =\ i_m\  &\Longleftrightarrow\ 
\frak{e}_m \preceq g \prec \frak{e}_{m+1}, \qquad (1\le m < n);\\
\nval P_{\times g}\ =\ i_n\ &\Longleftrightarrow\ \frak{e}_n\prec g;\\
\ndeg P_{\times g}\ =\ i_n\ &\Longleftrightarrow\  \frak{e}_n \preceq g.
\end{align*}
\end{cor}

\noindent
This follows from Corollary~\ref{newtpoldes} for $\E=K^\times$.

\begin{lemma}\label{upondegcut1} Assume $E$ is $\upo$-free, $\Gamma_E$ is divisible, and 
$(g_{\rho})$ is a pc-sequence in~$E$
with $g_\rho\leadsto 0$. Let ${\bf g} := c_E(g_{\rho})$ be the corresponding
cut in $E$ and 
$$\E\ :=\ \{g\in E^\times: g\prec g_{\rho},\ \text{eventually}\}.$$
Let $P\in E\{Y\}^{\ne}$.
If $\E\neq\emptyset$, then $\ndeg_{\bf g} P = \ndeg_{\E} P$. 
If $\E=\emptyset$, then $(g_\rho)$ is a c-se\-quence in $E$, and
$\ndeg_{\bf g} P = \val P$.
\end{lemma}
\begin{proof} Set $\gamma_\rho=v(g_{\rho+1}-g_\rho)\in\Gamma_\infty$. Removing some initial 
terms we arrange that~$(\gamma_\rho)$ is strictly increasing, and 
$\gamma_\rho=v(g_\rho)\in \Gamma$ for each $\rho$.
Then for all $\rho$,
$$\ndeg_{\geq\gamma_\rho} P_{+g_\rho}\ =\ \ndeg_{\geq\gamma_\rho} P\ =\ \ndeg P_{\times g_{\rho}}.$$
If $\E\neq\emptyset$, the desired result follows easily from Corollary~\ref{newtpolext}:
the critical case is when $v(\E)$ has a least element of the form 
$v(\frak{e}_m)$ using notation from that corollary.   
If  $\E=\emptyset$, use Lemma~\ref{ndevcon}. 
\end{proof}

\begin{cor}\label{cor:ndegE}
Suppose that $E$ is $\upo$-free and $\Gamma_E$ is divisible. Let $(h_\rho)$ be a pc-sequence in $E$ with
pseudolimit $h\in E$. Let $\mathbf h:=c_E(h_\rho)$  and
$$\E\ :=\ \{g\in E^\times: \text{$g\prec h_{\rho}-h$, eventually}\}.$$
Let $P\in E\{Y\}^{\ne}$.
If $\E\neq\emptyset$, then 
$\ndeg_{\bf h} P = \ndeg_{\E} P_{+h}$,
and 
if $\E=\emptyset$, then 
$\ndeg_{\bf h} P = \val P_{+h}$.
\end{cor}
\begin{proof}
Put $g_\rho := h_\rho-h$. Then $(g_\rho)$ is a pc-sequence in $E$ with $g_\rho\leadsto 0$. By 
Lemma~\ref{lem:basic facts on deg_a}(iii)
we have $\ndeg_{\bf h} P = \ndeg_{{\bf g}+h} P = \ndeg_{\bf g} P_{+h}$, with $\mathbf g = c_E(g_{\rho})$. It
remains to
apply Lemma~\ref{upondegcut1} to $P_{+h}$ in place of $P$.
\end{proof}

\subsection*{Notes and comments}
The conclusion of Corollary~\ref{hevq} 
goes through for $\upl$-free $E$ and homogeneous $P,Q \in E\{Y\}^{\neq}$ of order $\leq 1$
and degrees $d>e$ with $(d-e)\Gamma_E=\Gamma_E$.\marginpar{taken on faith for now} This fact will not be used in this volume, so we omit the proof.  %We do not know if for every $\upl$-free $E$ its differential-valued hull  is $\upl$-free.

\section{Further Consequences of $\upl$-Freeness}\label{fcuplfr}

\noindent
{\em In this section $E$ is a $\upl$-free $H$-asymptotic field}\/ (and thus
pre-$\d$-valued with rational asymptotic integration).  Our main goal here is to establish Proposition~\ref{prop:sigma(upg)=upo}. The basic tool is Lemma~\ref{huplshape}, which rests on 
Corollary~\ref{cor:nabla1, order 2}.

Recall from Section~\ref{sec:secondorder} that 
for elements $y\neq 0$ and $z$ in a differential field, we let $\omega(z)=-(2z'+z^2)$ and $\sigma(y)=\omega(-y^\dagger)+y^2$. If
$F$ is a differential field, $f\in F$, and~$y$ is a nonzero
solution in $F$ of the equation $\sigma(y)=f$, then multiplication by
$y^2$ yields $2y''y= 3(y')^2 - y^4 + fy^2$, so $y$ satisfies the differential equation $A(y)y'' =B(y)$, with $A:=2Y$ and $B:=3(Y')^2-Y^4+fY^2$. This is why we consider such equations, for example in Corollary~\ref{cor:2nd order, uniqueness} below.   

\begin{prop}\label{prop:sigma(upg)=upo}
Suppose $\upo_\rho\leadsto\upo\in E$. Then there is an element $\upg$ of a pre-$\d$-valued field extension of $E$ of $H$-type such that $\Psi_E<v\upg< (\Gamma_E^>)'$ and $\sigma(\upg)=\upo$.
For any such elements $\upg$ and $\upg^*$, in possibly different extensions, there is an isomorphism $E\<\upg\>\to E\<\upg^*\>$ over $E$
sending $\upg$ to~$\upg^*$.
\end{prop}

\noindent
Throughout this section we assume $P\in E\{Y\}^{\ne}$. 

\begin{lemma}\label{huplshape}  Suppose  
$P$ has order at most~$2$ and $\nwt P = w$. There are $a\in E$,
$A\in E[Y]^{\ne}$ with $A\asymp 1$ and $w+\deg A= \ndeg P$, 
and active $e$ in $E$ such that
$$ P^e\ =\ a\cdot A \cdot(Y')^w + R, \qquad a\asymp P^e,\ R\in E^e\{Y\},\ R\prec^{\flat}_e P^e.$$
Given such $a, A, e$ and any active $f\preceq e$ in $E$ there is an $a_f\in E$ such that
$$ P^f\ =\ a_f\cdot A\cdot (Y')^w + R_f,\qquad a_f\asymp P^f, \quad R_f\in E^{f}\{Y\},\ R_f\prec^{\flat}_e P^f.$$ 
\end{lemma}
\begin{proof} For the first part, follow the proof of Lemma~\ref{lem:huposhape}
to reduce to the case that~$E$ has a monomial group $\frak{M}$
and small derivation. Then $K:= \operatorname{dv}(E)$ is $\d$-valued of $H$-type and an immediate extension of $E$.
Whether or not $K$ is $\upl$-free, Corollary~\ref{cor:nabla1, order 2}
yields $N_P\in C[Y](Y')^{\N}$
where $C=C_K$. The rest of the proof of Lemma~\ref{lem:huposhape} can be copied verbatim, except that to get $R_{P^e}\prec^{\flat}_{e} P^{e}$ as in that proof, we appeal to 
Lemma~\ref{uplfreeclean}, using that $K$ has rational
asymptotic integration. For the second part, use the proof of Corollary~\ref{huposhape}. 
\end{proof} 

%\begin{lemma}\label{huplshape}  Suppose  
%$P$ has order at most~$2$ and Newton weight~$w$. Then there is an
%$A\in E[Y]^{\ne}$ with $A\asymp 1$ and $w+\deg A= \ndeg P$, 
%and an active $f_0$ in $E$
%such that for all active $f\preceq f_0$ in $E$ we have $a_f\in E$ and  
%$R_f\in E^f\{Y\}$ with
%$$ P^f\ =\ a_f\cdot A\cdot (Y')^w + R_f,\qquad a_f\asymp P^f, \quad %R_f\prec^{\flat}_f P^f.$$ 
%\end{lemma}
%\begin{proof} Following the proof of Corollary~\ref{huposhape}
%we reduce to the case that $E$ has a monomial group $\frak{M}$
%and small derivation. Then $K:= \operatorname{dv}(E)$ is $\d$-valued %of $H$-type and an immediate extension of $E$.
%Whether or not $K$ is $\upl$-free, 
%Corollary~\ref{cor:nabla1, order 2}
%yields $N_P\in C[Y](Y')^{\N}$
%where $C=C_K$. The rest of the proof of Corollary~\ref{huposhape} can %be copied verbatim. However, to get $R_{P^\phi}\prec^{\flat}_{\phi} %P^{\phi}$ as in that proof, we appeal to 
%Lemma~\ref{uplfreeclean}, using that $K$ has rational
%asymptotic integration.
%\end{proof} 

\noindent
Next we establish variants of some results in Section~\ref{cutsvalgrp}.
Recall from Section~\ref{Compositional Conjugation}
that for $f\in E^\times$ we defined $P^{\times f}:=P^f_{\times f}\in E^f\{Y\}$,
so $P^{\times f}(Y')=P(Y')^f$. Let us generalize notation introduced in Section~\ref{cutsvalgrp}, and put
$$\nwt^\times(P)\ :=\ \nwt\!\big(P(Y')\big),\qquad v^{\times\!\ev}(P)\ :=\ v^{\ev}\big(P(Y')\big).$$
Since $v(P^{\times f})=v\big(P(Y')^f\big)$, there is an active $f_0$ in $E$ such that
$$v(P^{\times f}) \ =\ v^{\times\!\ev}(P)+\nwt^\times(P)vf\quad\text{for active $f\preceq f_0$ in $E$.}$$
Here is a variant of Corollary~\ref{cor:mult-comp, 1}:

\begin{cor}\label{cor:v(P(y)) in gap, 1}
Suppose $\order(P)\leq 1$.
Then there exists an active $f$ in $E$ such that for every $H$-asymptotic
field extension $F$ of $E$ and active $z\prec f$ in $F$,
$$v\big(P(z)\big)\ =\ v^{\times\!\ev}(P) + \nwt^\times(P)vz.$$
\end{cor}
\begin{proof}
Let $w=\nwt^\times(P)$.
Apply Lemma~\ref{huplshape} to $P(Y')$ in place of $P$ to get
active~$f$ in $E$, $a\in E^\times$, and $R\in E^f\{Y\}$ such that 
$$P^{\times f}\ =\  a \cdot (Y^w+R),\quad va\ =\ v(P^{\times f})\ =\  v^{\times\!\ev}(P) + w\,vf,\quad
R\prec^\flat_{f} 1.$$
Let $F$ be an $H$-asymptotic field extension of $E$ and $z$ an active element of $F$ with $z\prec f$.
Then $y:=zf^{-1}\prec 1$ is active in $F^{f}$, so $y\asymp^\flat_{f} 1$. Now $F^f$ has small derivation, so
$R(y)\prec^\flat_{f} 1 \asymp^\flat_{f} y^w$.
Thus $P(z)=P^{\times f}(y)\simflat_f ay^w$,
and hence
$v\big(P(z)\big) = v(a \cdot y^w) = 
v^{\times\!\ev}(P) + w\,vz$ as claimed.
\end{proof}

\noindent
In a similar way we obtain an analogue of Corollary~\ref{cor:mult-comp, 4}:

\begin{cor}\label{cor:v(P(y)) in gap, 2}
Suppose $\order(P)\leq 1$.
Then there exists $\delta_0\in (\Gamma_E^>)'$ such that  for every $H$-asymptotic
field extension~$F$ of~$E$ and $z\in F$ with $\Psi_E<vz\leq\delta_0$,
$$v\big(P(z)\big) = v^{\times\!\ev}(P) + \nwt^\times(P)vz.$$
\end{cor}
\begin{proof}
Let $w$, $f$, $a$ and $R$ be as in the proof of 
Corollary~\ref{cor:v(P(y)) in gap, 1}.
By Lemma~\ref{hasc} and the remark following it, we can take $\gamma_0\in\Gamma_E^{>vf}$ and $\delta_0\in (\Gamma_E^>)'$ such that the map $\gamma\mapsto\psi_F(\gamma-vf)\colon [\gamma_0,\delta_0]_{\Gamma_F}\to\Gamma_F$ is constant, for each $H$-asymptotic field extension~$F$ of~$E$. Given such an $F$ and $z\in F$ with
$\Psi_E<vz\leq\delta_0$, set $y:={zf^{-1}\prec 1}$. Then $y\asymp^\flat_f 1$. (See proof of  Corollary~\ref{cor:mult-comp, 4}.) Hence  we get $P(z)=P^{\times f}(y)\simflat_f ay^w=af^{-w}z^w$ as in the proof of Corollary~\ref{cor:v(P(y)) in gap, 1}.
\end{proof}

\begin{remark}
Suppose $\order(P)\leq 1$, and  $E$ is equipped with an ordering making $E$ a pre-$H$-field. The proofs of Corollaries~\ref{cor:v(P(y)) in gap, 1} and~\ref{cor:v(P(y)) in gap, 2} show that there are $\gamma_0\in(\Gamma_E^<)'$, $\delta_0\in(\Gamma_E^>)'$, and $\sigma\in\{-1,+1\}$, such that  for every pre-$H$-field extension~$F$ of~$E$ and $z\in F^>$ with $\gamma_0\leq vz\leq\delta_0$, we have $\sgn P(z)=\sigma$.
\end{remark}

\noindent
Note that by Corollary~\ref{cor:v(P(y)) in gap, 2}, if  $z$ is an element in an $H$-asymp\-totic field extension of~$E$
with $\Psi_E < vz < (\Gamma_E^>)'$, then $z$,~$z'$ are algebraically independent over $E$.

\begin{cor}\label{cor:2nd order, uniqueness}
Let $A,B\in E\{Y\}$ be of order at most~$1$, $A\ne 0$,
and let $y$,~$z$ be elements of \textup{(}possibly different\textup{)} $H$-asymptotic field extensions of $E$ 
such that 
\begin{align}
\label{eq:y''}&\Psi_E<vy<(\Gamma_E^>)',\quad A(y)y''=B(y),\\
\label{eq:z''}&\Psi_E<vz<(\Gamma_E^>)',\quad A(z)z''=B(z).
\end{align}
%\begin{equation}\label{eq:y''}
%\Psi_E<vy<(\Gamma_E^>)',\quad A(y)y''=B(y),
%\end{equation}
%\begin{equation}\label{eq:z''}
%\Psi_E<vz<(\Gamma_E^>)',\quad A(z)z''=B(z).
%\end{equation}
Then there is an isomorphism $E\<y\>\to E\<z\>$ over $E$ sending $y$ to $z$.
\end{cor}
\begin{proof} 
By the remark preceding the corollary,  $y$,~$y'$ are algebraically
independent over $E$, and so are $z$,~$z'$. Thus $A(y)\ne 0$, $A(z)\ne 0$, and so $E\<y\>=E(y,y')$ and $E\<z\>=E(z,z')$ as fields,
and we have an isomorphism
$h\colon E\<y\>\to E\<z\>$ of differential fields over $E$ sending
$y$ to $z$.  If $P$ has order $\le 1$, then
$$v\big(P(y)\big)\ =\ v^{\times\!\ev}(P)+ \nwt^{\times }(P)vy$$
by Corollary~\ref{cor:v(P(y)) in gap, 2}, and similarly with $y$ replaced by $z$.
In view of Corollary~\ref{addgapmore} it follows that $h$ is also a valued field isomorphism.
\end{proof}

\begin{remark}
Let $A$ and $B$ be as in Corollary~\ref{cor:2nd order, uniqueness}, and suppose additionally that~$E$ is equipped with an ordering making $E$ a pre-$H$-field. Let $y>0$ and $z>0$ be elements of pre-$H$-field extensions of $E$ satisfying \eqref{eq:y''} and \eqref{eq:z''}, respectively.
The isomorphism $E\<y\>\to E\<z\>$ of that corollary is also an isomorphism of ordered fields, in view of the remark following the proof of 
Corollary~\ref{cor:v(P(y)) in gap, 2}.
\end{remark}

\begin{lemma}\label{lem:E(z,z')}
Let $z$ be an element of an $H$-asymptotic field extension of $E$ with
$\Psi_E < vz < (\Gamma_E^>)'$, and consider the valued field extension $F:=E(z,z')$ of $E$.
Then $\Gamma_F = \Gamma_E\oplus \Z vz$ \textup{(}internally\textup{)}, $[\Gamma_F]=[\Gamma_E]$, and
$\res(F)=\res(E)$.
\end{lemma}
\begin{proof}
Clearly $\Gamma_F = \Gamma_E + \Z vz$ by Corollary~\ref{cor:v(P(y)) in gap, 2}.
Since $E$ is $\upl$-free, the algebraic closure $E^\alg$ of $E$ 
(inside the algebraic closure of $E\<z\>$) has asymptotic integration.
Since $\Psi_{E^\alg}=\Psi_E < vz< (\Gamma_{E^\alg}^>)'$, we have $vz\notin \Gamma_{E^\alg}=\Q\Gamma_E$. Thus $\Gamma_F=\Gamma_E\oplus \Z vz$. From Lemma~\ref{gendense} applied to $G:= \Q\Gamma_E$ and $b:= vz$ we obtain $[\Gamma_F]=[\Gamma_E]$. 
To show $\res(F)=\res(E)$, suppose $P$ has order $\le 1$.
Taking $w$, $f$, $a$, and $R$ as in the proof of Corollary~\ref{cor:v(P(y)) in gap, 1}, that proof yields 
$P(z) \sim af^{-w} z^w$. 
Now let $g\in F^\times$. Then the above gives
$e\in E^\times$ and $k\in \Z$ with $g\sim ez^k$. If $g\asymp 1$ this
forces $k=0$, and so $g\sim e\in E^\times$. Thus $\res(F)=\res(E)$. 
\end{proof}

\noindent
We can now prove Proposition~\ref{prop:sigma(upg)=upo}.

\begin{proof}[Proof of Proposition~\ref{prop:sigma(upg)=upo}]
The uniqueness statement is a consequence of Cor\-ol\-lary~\ref{cor:2nd order, uniqueness}, so it suffices to prove the existence of $\upg$ with the desired properties. We arrange that
$\Gamma_E$ is divisible, for example by passing to the algebraic closure
of~$E$. 
Corollary~\ref{cor:adjoin upl} yields an immediate asymptotic
extension $E(\upl)$ of $E$ with $\upl_{\rho}\leadsto \upl$
and $\omega(\upl)=\upo$.
Note that $\upl$ is transcendental over $E$. By  Lemma~\ref{cg2}, $-\upl$ creates a gap over~$E(\upl)$, and the proof of that lemma
yields an element 
$\upg\neq 0$ in an $H$-asymptotic field extension of $E(\upl)$ such that $v\upg\notin \Gamma_E$, $\upg^\dagger=-\upl$, $\Psi_{E(\upl, \upg)}=\Psi_E$
and $v\upg$ is a gap in~${E(\upl, \upg)}$. Note that then
$\upg$ is transcendental over $E(\upl)$, and
$$\Gamma_{E(\upl,\upg)}\ =\ \Gamma_E\oplus \Z v\upg,\qquad 
\res E(\upl,\upg)\ =\ \res E, \qquad \sigma(\upg)\ =\ \upo+\upg^2.$$
In the field extension $E(\upl,\upg)=E(\upg,\upg')$ of $E$ the elements $\upl$, $\upg$ are algebraically independent over $E$. 
Using Corollary~\ref{lem:derivative on trans ext}, let $E_{\upg}$ be the valued differential field extension
of $E$ with $E(\upl, \upg)$ as its underlying valued field, and 
derivation $\der_{\upg}$ given by
$$
\der_{\upg}(\upl)\ =\ \frac{1}{2}\big(\upg^2-\upl^2-\upo\big), \quad
\der_\upg(\upg)\ =\ -\upg\cdot\upl.
$$
Since $-2\upl'-\upl^2=\upo$, the first equality amounts to
$\der_{\upg}(\upl)=\upl'+\frac{1}{2}\upg^2$. The second equality
just says that $\der_{\upg}(\upg)=\upg'$.  
Let $\omega_\upg$, $\sigma_\upg$ be the analogues of  $\omega$,~$\sigma$ in the differential field $E_\upg$. Then 
$\omega_\upg(z)=-2\der_\upg(z)-z^2$ for $z\in E_\upg$ and 
$$\sigma_\upg(y)\ =\ \omega_\upg\big({-\der_\upg(y)/y}\big)+y^2 \qquad(y\in E_\upg^\times),$$
and thus $\sigma_\upg(\upg)=\upo$, as is easily checked. 
To show that $E_{\upg}$ is pre-$\d$-valued, we use Lemma~\ref{Lemma3.3} with
$E$,~$E_\upg$ in place of $K$,~$L$, and with
the multiplicative
subgroup $T:=E(\upl)^\times\cdot\upg^\Z$ of~$E_\upg^\times$. In order to apply that lemma, we need:

\claim{If $s,t\in T$ and $s\preceq 1$, $t\prec 1$, then $\der_\upg(s)\prec \der_\upg(t)/t$ in~$E_{\upg}$.}

\noindent
To prove this claim, let $t\in T$. Take $a\in E(\upl)^\times$ and 
$k\in\Z$ such that $t=a\upg^k$. Take $A\in E(Y)^\times$ with 
$A(\upl)=a$.
Let $R\mapsto R^\der$ be the derivation on the field $E(Y)$ extending that of $E$ with $Y^\der=0$.
Then by Corollary~\ref{cor:derpol},
\begin{align*} a'\, &=\, A^\der(\upl)+(\partial A/\partial Y)(\upl)\cdot\upl',\qquad \der_\upg(a)\, =\, A^\der(\upl)+(\partial A/\partial Y)(\upl)\cdot\der_\upg(\upl), \text{ so}\\
\frac{\der_\upg(t)}{t} &- \frac{t'}{t}= \frac{\der_\upg(a)}{a} - \frac{a'}{a}  =    \left(\frac{\partial A/\partial Y}{A}\right)(\upl)\cdot\big(\der_\upg(\upl)-\upl'\big)  =  \left(\frac{\partial A/\partial Y}{A}\right)(\upl)\cdot \textstyle\frac{1}{2}\upg^2.
\end{align*}
Then Corollary~\ref{cor:val of formal log der at upl} gives
$$v\left(\frac{\der_\upg(t)}{t} - \frac{t'}{t}\right)\ >\  v\upg\ >\  \Psi_E = \Psi_{E(\upl,\upg)}.$$  If $t\not\asymp 1$, this gives 
$\der_\upg(t)/t \sim t'/t$, and so $\der_\upg(t)\sim t'$.
Suppose now that $t\asymp 1$. Then $k=0$ and so $A(\upl)=a=t\asymp 1$, hence
$$\der_\upg(t)-t'\  =\ (\partial A/\partial Y)(\upl)\cdot {\textstyle\frac{1}{2}}\upg^2\ \asymp\ \left(\frac{\partial A/\partial Y}{A}\right)(\upl) \cdot \textstyle\frac{1}{2}\upg^2.$$
As before this yields $v\big(\der_\upg(t)-t'\big)> \Psi_E = \Psi_{E(\upl,\upg)}$, and thus $v\big(\der_\upg(t)\big) >\Psi_{E\<\upl,\upg\>}$. 
The claim now follows. 
\end{proof}

\begin{remarks}
Let $E$ and $E\<\upg\>$ be as in Proposition~\ref{prop:sigma(upg)=upo}.
\begin{enumerate}
\item By Lemma~\ref{lem:E(z,z')} we have $\res E\<\upg\> = \res E$.
In view of Lemma~\ref{dv} it follows that if 
$E$ is $\d$-valued, then so is $E\<\upg\>$ with 
$C_{E\<\upg\>}=C_E$. 
\item Lemma~\ref{lem:E(z,z')} also gives $[\Gamma_E]=[\Gamma_{E\<\upg\>}]$,
and so $v\upg$ is a gap in $E\<\upg\>$. 
\item Suppose
$b\asymp 1$ in $E\<\upg\>$; then $b' \prec \upg$. This is because
$\res E\<\upg\>=\res E$ gives
$u\asymp 1$ in $E$ with $b\sim u$, so $u' \prec \upg$ and
$b'-u'\prec \upg$, and thus $b'\prec \upg$.  
\end{enumerate}
\end{remarks}

\begin{lemma}\label{cor:sigma(upg)=upo}
Let $E$ and $\upg$ be as in Proposition~\ref{prop:sigma(upg)=upo} and let 
$E$ be equipped with an ordering making it a pre-$H$-field.
Then $E\<\upg\>$ has a unique field ordering extending that of $E$ in which $\upg>0$ and $\mathcal{O}_{E\<\upg\>}$ 
is convex. Moreover, $E\<\upg\>$ with this ordering is a
pre-$H$-field in which $\mathcal{O}_{E\<\upg\>}$ is the convex hull of 
$\mathcal{O}_E$. 
\end{lemma}
\begin{proof} Let $f\in E\<\upg\>^\times$. As $E\<\upg\>=E(\upg, \upg')$,  Lemma~\ref{lem:E(z,z')} gives 
$$f\ =\ g\upg^k (1+ a),\ \text{ with $g\in E^\times$, $k\in \Z$,  $a\in E\<\upg\>$, $a\prec 1$.}$$
Thus for any field ordering on $E\<\upg\>$ extending that of $E$ in
which $\upg>0$ and the valuation ring of $E\<\upg\>$
is convex we have: $f>0\Longleftrightarrow g>0$. Thus
at most one such ordering on $E\<\upg\>$ exists.
For any such ordering the valuation ring of $E\<\upg\>$ is
the convex hull of $\mathcal{O}$, since $\res(E)=\res(E\<\upg\>)$ by
Lemma~\ref{lem:E(z,z')}.   

It remains only to show that $E\<\upg\>$ has a field ordering 
making it a pre-$H$-field extension of $E$ in which
$\upg>0$. For this we use the construction of $E\<\upg\>$ in the proof of 
Proposition~\ref{prop:sigma(upg)=upo}. In detail, we first arrange that
$\Gamma_E$ is divisible, for example, by passing to the real closure of
$E$. Then we take an immediate asymptotic
extension~$E(\upl)$ of $E$ with $\upl_{\rho}\leadsto \upl$
and $\omega(\upl)=\upo$, and use Lemma~\ref{prehimm} to make it a pre-$H$-field
extension of $E$. Next we consider a Liouville closed $H$-field extension
of $E(\upl)$, and take an element $\upg>0$ in this extension
such that $\upg^{\dagger}=-\upl$. Then $E\<\upg\>=E(\upl, \upg)$ is a
pre-$H$-field extension of $E(\upl)$, and $v\upg$ is a gap in $E\<\upg\>$
by the remark following the proof of Lemma~\ref{cg2}. It now follows from
Lemma~\ref{lem:E(z,z')} that $E(\upl, \upg)$ as a pre-$\d$-valued
field extension of $E$ is exactly as constructed in the proof of
Proposition~\ref{prop:sigma(upg)=upo}, and in what follows we define
$\der_{\upg}$ and $E_{\upg}$ as in that proof, and
use its notations. We also consider $E_{\upg}$ to be
equipped with the field ordering of~$E(\upl, \upg)$. It is enough to show that~$E_{\upg}$ is a pre-$H$-field. Recall that 
$T:=E(\upl)^\times\cdot\upg^\Z$. 
By Lemma~\ref{prehunr}
it suffices to show: 

\claim{Suppose $t\in T$ and $t\succ 1$; then $\der_\upg(t)/t>0$.}

\noindent
To establish this claim, note that $\der_\upg(t)/t \sim t'/t$ by the proof of Proposition~\ref{prop:sigma(upg)=upo}. Since $E(\upl,\upg)$ is a pre-$H$-field, we have $t'/t>0$, and the claim follows.
\end{proof}

\begin{remark}
If in addition to the hypotheses of Lemma~\ref{cor:sigma(upg)=upo}, $E$ is an $H$-field, then~$E\<\upg\>$ ordered as in that lemma is also an $H$-field, with $C_{E\<\upg\>}=C_E$. 
\end{remark}

\subsection*{Cases of low complexity} For use in the next chapter we derive here results for~$P$ of low complexity, under our standing assumption that $E$ is $\upl$-free. 

\begin{cor}\label{upltwondegnval} If $\order P\le 2$, then there is $\gamma\in \Gamma_E^{>}$ such that
for all $g\in E$,
\begin{align*} 0<vg< \gamma\ &\Longrightarrow\  \nval(P)\ =\ \nval(P_{\times g})\ =\ \ndeg(P_{\times g}),\\
-\gamma < vg < 0\  &\Longrightarrow\  \ndeg(P)\ =\ \ndeg(P_{\times g})\ =\ \nval(P_{\times g}).
\end{align*}
\end{cor}
\begin{proof} Suppose $\order P \le 2$. Use Lemma~\ref{huplshape} and compositional conjugation
by some active element of $E$
to arrange that $E$ has small derivation and $P$ is as in
the hypothesis of Lemma~\ref{decndegnval}. Now apply Lemma~\ref{decndegnval}.
\end{proof}
 
\begin{lemma}\label{ndzeroupl} If $\ndeg P=0$ and $f_0$ is active in $E$, then for some active $f\preceq f_0$ in~$E$
we have $P^f \simflat_{f} a$ in $E^{f}\{Y\}$ with 
$a\in E$.
\end{lemma}

\noindent
We omit the proof: it is like that of the next lemma but shorter.

\begin{lemma}\label{ndoneupl} Suppose $\ndeg P= 1$ and $f_0$ is active in $E$. Then for some
active $f\preceq f_0$ in $E$ we have either $P^{f}\simflat_{f} a+bY$ in $E^{f}\{Y\}$ with $a\preceq b$ in $E$,
or $P^f\simflat_{f} bY'$ in $E^f\{Y\}$ with $b\in E$. In particular,
$\nwt(P)\le 1$.  
\end{lemma}
\begin{proof} By compositional conjugation we arrange that $E$ has small derivation.
By passing to an elementary extension of $E$ we further arrange
that $E$ has a monomial group $\frak{M}$. Then $K:= \operatorname{dv}(E)$
is an immediate extension of $E$. 
With $\frak{M}$ as its monomial group $K$ satisfies the assumptions at the beginning of this chapter. Then $\ndeg P=1$ gives $P_1\ne 0$, and
either $N_P=N_{P_0} + N_{P_1}$ or $N_P=N_{P_1}$. 
Applying Corollary~\ref{Euplone} to
$Q:= P_1$ gives $N_{P_1}=cY$ or $N_{P_1}=cY'$ with $c\in C^\times$. Since $N_P$ is isobaric, this gives
$N_P=c_0+cY$ or $N_P=cY'$, with $c_0, c\in C$, $c\ne 0$.
In particular, $N_P\in C[Y](Y')^{\N}$. We have $\frak{d}_{P^{\phi}}\in \frak{M}\subseteq E^\times$, and 
$P^{\phi}\simflat_{\phi} \frak{d}_{P^{\phi}} N_P$, eventually, by 
Lemma~\ref{uplfreeclean}.
We just do the case $N_P=c_0 + cY$. (The other case is similar.) Take $e_0\preceq e_1\asymp 1$ in $E$ such that
$c_0-e_0\prec 1$ and $c_1-e_1\prec 1$ in $K$. Recall that 
$\phi$ ranges over the elements of $\frak{M}$ active in $K$, but $\frak{M}\subseteq E^\times$, so $\phi\in E$, and eventually
$$P^\phi\ \simflat_{\phi}\  \frak{d}_{P^{\phi}}N_P\ =\ \frak{d}_{P^{\phi}}(e_0+ e_1Y) + \frak{d}_{P^{\phi}}\big((c_0-e_0) + (c_1-e_1)Y\big).$$
Also $(c_0-e_0)+(c_1-e_1)Y\prec^{\flat}_{\phi} 1$, eventually, so $P^\phi\ \simflat_{\phi}\  \frak{d}_{P^{\phi}}(e_0+ e_1Y)$, eventually.  
 \end{proof}

\begin{cor}\label{lem:npol linear}
If $\ndeg P \le 1$, then there is $\gamma\in \Gamma_E^{>}$ such that
for all $g\in E$,
\begin{align*} 0<vg< \gamma\ &\Longrightarrow\  \nval(P)\ =\ \nval(P_{\times g})\ =\ \ndeg(P_{\times g}),\\
-\gamma < vg < 0\  &\Longrightarrow\  \ndeg(P)\ =\ \ndeg(P_{\times g})\ =\ \nval(P_{\times g}).
\end{align*}
\end{cor}
\begin{proof} Assume $\ndeg P \le 1$. By Lemmas~\ref{ndzeroupl} and~\ref{ndoneupl} and compositional conjugation we arrange
that $E$ has small derivation and $P$ is as
in the hypothesis of Lemma~\ref{decndegnval}, with $w=0$ or $w=1$. Now apply
Lemma~\ref{decndegnval}. 
\end{proof}

\begin{lemma}\label{uplnvallemma} If $\nval P = 1$, then there is $\gamma\in \Gamma_E^{>}$ such that for all $g\in E^\times$,
$$ 0 < vg < \gamma\ \Longrightarrow\ \nval P_{\times g}\ =\ \ndeg P_{\times g} =1.$$
\end{lemma}
\begin{proof} Assume $\nval P =1$. 
As in the proof of Lemma~\ref{ndoneupl} we arrange that
$E$ has small derivation, a monomial group, and an immediate ($\d$-valued) extension $K$. With 
$P=P_0 + P_1 + R$, where $R:=\ \sum_{d\ge 2} P_d$, 
we have $P_0\prec P_1^\phi\succeq R^\phi$ in $K^{\phi}\{Y\}$, eventually,
and $v(P_1^\phi)=v^{\ev}(P) + \nwt(P_1)v\phi$, eventually.
By Corollary~\ref{Euplone} we have $\nwt(P_1)\le 1$, and as
$E$ has asymptotic integration, this yields
$\alpha\in \Gamma_E^{>}$ such that $v(P_0) > \alpha + v^{\ev}(P) + \nwt(P_1)v\phi$ for all $\phi$. Then for $\gamma\in \Gamma_E^{>}$ with $2\gamma<\alpha$ we have for all $g\in E^\times$,
$$0 < vg < \gamma\ \Longrightarrow\ v(P_0)\ >\ v(P_{1,\times g}^\phi)\ <\ 
v(R_{\times g}^{\phi}), \text{ eventually,}$$
so any such $\gamma$ has the desired property.   
\end{proof}

\subsection*{Notes and comments}
Corollary~\ref{cor:v(P(y)) in gap, 2} 
improves on~\cite[Proposition~12.12]{AvdD3}.

\section{Asymptotic Equations}\label{sec:aseq}

\noindent
This section is somewhat technical, but indispensable. The main facts
we establish here are Proposition~\ref{prop:unravel} and Lemma~\ref{lem:neglect >d}.
These are needed in Chapter~\ref{ch:newtonian fields} in proving Theorems~\ref{newtmax} and \ref{thm:solve asymptotic equ}, which are of independent interest. Our quantifier elimination for~$\T$
depends essentially on Theorem~\ref{newtmax}.

Recall the assumptions on $K$ imposed in the beginning of this chapter.
{\em In this section we assume in addition that $K$ is $\upo$-free,}\/ so that
we can use results of Section~\ref{eveqth}. Also, as before,  
$P$ ranges over $K\{Y\}^{\ne}$ and $\phi$ over the elements
of~$\frak{M}$ that are active in~$K$. 
Let $\E\subseteq K^\times$ be  $\preceq$-closed,
and consider the asymptotic equation over $K$ given~by
\begin{equation}\tag{E} \label{eq:asympt equ 1}
P(Y)=0, \qquad Y\in \E.
\end{equation}
An {\bf approximate solution} of 
\eqref{eq:asympt equ 1} \index{asymptotic equation!approximate solution}\index{approximate!solution} is an approximate zero~$y$ of $P$ such that $y\in \E$, and
the {\bf multiplicity} \index{multiplicity} of an approximate solution~$y$ of~\eqref{eq:asympt equ 1} is its multiplicity
as an approximate zero of $P$. An {\bf algebraic approximate solution} of 
\eqref{eq:asympt equ 1} \index{approximate!solution!algebraic} is an algebraic approximate zero~$y$ of $P$,
as defined at the end of Section~\ref{The Shape of the Newton Polynomial}, such that~$y\in \E$. 
Let $y\in\E$ with $y\sim c\fm$ ($c\in C^\times$). Then by  Lemma~\ref{transition},
$$\text{$y$  is an approximate solution of 
\eqref{eq:asympt equ 1}} \ \Longleftrightarrow\ 
N_{P_{\times\fm}}(c)=0  \ \Longleftrightarrow\ 
\ndeg_{\prec\fm} P_{+y} \ge 1,$$
and if $y$ is an approximate solution of 
\eqref{eq:asympt equ 1}, then its multiplicity equals $\ndeg_{\prec\fm} P_{+y}$, which is $\le \ndeg_{\E}P$. If $y$ is an approximate solution of  \eqref{eq:asympt equ 1}, then so is every  $z\in K^\times$ with $z\sim y$, with the same multiplicity.
By Lemma~\ref{newtondetection}, if $y$ is a solution  of \eqref{eq:asympt equ 1}, then~$y$ is an approximate solution of  \eqref{eq:asympt equ 1}.  
For each $\phi$, the asymptotic equation
\begin{equation}  \label{eq:asympt equ 1phi} \tag{E$^\phi$}
P^\phi(Y)=0, \qquad Y\in\E
\end{equation}
has the same solutions and the same approximate solutions as  \eqref{eq:asympt equ 1}, with the same multiplicities.

\index{Newton!degree!asymptotic equation}
\index{asymptotic equation!refinement}
\index{refinement}
\index{starting monomial}

\medskip
\noindent
A {\bf starting monomial} for \eqref{eq:asympt equ 1} is a starting monomial $\fm$ for $P$ with $\fm\in\E$;
we define ``algebraic starting monomial  for \eqref{eq:asympt equ 1}'' likewise. 
%and ``differential starting monomial for \eqref{eq:asympt equ 1}.''
If $y\sim c\fm$ ($c\in C^\times$) is an approximate solution of \eqref{eq:asympt equ 1}, then
$\fm$ is a starting monomial for  \eqref{eq:asympt equ 1}. Hence if~\eqref{eq:asympt equ 1} has no starting monomial (in particular, if $\ndeg_{\E}P=0$),
then \eqref{eq:asympt equ 1} has no approximate solution.  By Proposition~\ref{newtpolygon}, if
$\Gamma$ is divisible and $\val P < \ndeg_{\E} P$, then there is an
algebraic starting monomial for
\eqref{eq:asympt equ 1}, and
$\ndeg P_{\times \mathfrak{e}}=\ndeg_{\E} P$, 
$\mathfrak{e}:=\text{largest algebraic starting monomial for  \eqref{eq:asympt equ 1}}$.

\medskip
\noindent
Let $\E'\subseteq \E$ be $\preceq$-closed and let $f\in \E\cup\{0\}$. We call the asymptotic equation
\begin{equation}\tag{E$'$} \label{eq:asympt equ 2}
 P_{+f}(Y)=0, \qquad Y\in \E'
\end{equation}
a {\bf refinement} of  \eqref{eq:asympt equ 1}. 
By Lemma~\ref{lem:ndeg add conj} we have
$\ndeg_{\E'}P_{+f}\le \ndeg_{\E} P$. 
Note that if~$y$ is a solution of 
\eqref{eq:asympt equ 2} and $f+y\ne 0$, then $f+y$ is a solution of \eqref{eq:asympt equ 1}. Moreover:

\begin{lemma}
Let $y \not\sim -f$ be an approximate solution of 
\eqref{eq:asympt equ 2} of multiplicity $\mu$. Then $f+y$ is an approximate solution of \eqref{eq:asympt equ 1}
of multiplicity $\geq\mu$.
\end{lemma}
\begin{proof}
Since $y \not\sim -f$ we have $y\preceq f+y$ and thus
$$1\ \le\ \mu\ =\ \ndeg_{\prec y} (P_{+f})_{+y}\ =\ \ndeg_{\prec y} P_{+(f+y)}\ \leq\ \ndeg_{\prec f+y} P_{+(f+y)},$$
hence $f+y$ is an approximate solution of \eqref{eq:asympt equ 1}
of multiplicity $\geq\mu$.  
\end{proof}

\begin{lemma}\label{lem:equal ndeg}
Suppose $f\ne 0$, $\E'\subseteq K^{\prec f}$ and $\ndeg_{\E} P = \ndeg_{\E'} P_{+f}\geq 1$. Then $f$ is an approximate solution of 
\eqref{eq:asympt equ 1}.
\end{lemma}
\begin{proof} Using Lemma~\ref{lem:ndeg add conj}
we have
$$\ndeg_{\E'} P_{+f}\  \leq\ \ndeg_{\prec f}P_{+f}\ \leq\ \ndeg_{\preceq f} P_{+f}\ =\ \ndeg_{\preceq f} P\ \leq\ \ndeg_{\E} P,$$
and hence $\ndeg_{\prec f} P_{+f}=\ndeg_{\E} P\geq 1$. Thus $f$ is an approximate solution of
\eqref{eq:asympt equ 1} by the equivalence displayed earlier in this section.
% Lemma~\ref{transition}. 
\end{proof}

\noindent
Let an asymptotic equation \eqref{eq:asympt equ 1} be given, with Newton degree
$d=\ndeg_{\E} P$.
Then $\ndeg_{\prec f} P_{+f} \leq  d$ for all $f\in\E$.
Moreover:

\begin{lemma}\label{unraveledconditions}
Suppose $d\ge 1$. Then the following are equivalent:
\begin{enumerate}
\item[\textup{(i)}] $\ndeg_{\prec f} P_{+f} <  d$ for all $f\in\E$;
\item[\textup{(ii)}] $\ndeg_{\prec f} P_{+f} <  d$ for all $f\in\E$ 
with $\ndeg P_{\times f}=d$;
\item[\textup{(iii)}] there is no approximate solution of \eqref{eq:asympt equ 1} of multiplicity $d$.
\end{enumerate}
\end{lemma}
\begin{proof}
Let $f\in\E$ and suppose $\ndeg P_{\times f}<d$. Then
$$\ndeg_{\prec f} P_{+f}\ \leq\ \ndeg_{\preceq f} P_{+f}\ =\ 
\ndeg_{\preceq f} P\ =\ \ndeg P_{\times f}<d.$$
Thus (i)~$\Longleftrightarrow$~(ii). An earlier equivalence in this section gives
(i)~$\Longleftrightarrow$~(iii). 
\end{proof}

\index{asymptotic equation!unraveled}
\index{unraveled}

\noindent
We say that \eqref{eq:asympt equ 1} is {\bf unraveled} if $d\ge 1$ and one of the equivalent conditions in Lem\-ma~\ref{unraveledconditions} holds. So if $d\ge 1$ and  \eqref{eq:asympt equ 1} does not have an approximate solution, then~\eqref{eq:asympt equ 1} is unraveled. If \eqref{eq:asympt equ 1} is unraveled and has an approximate solution, then~${d\ge 2}$ by condition~(iii) of Lemma~\ref{unraveledconditions}.  

\begin{example} \marginpar{skipped for the moment}
Suppose $P\in K[Y](Y')^\N$. Then $d= \ndeg_\E P = \dd_\E P$. Assume that~$d\geq 1$. 
If $P\in K[Y]$,
then~\eqref{eq:asympt equ 1} is not unraveled if and only  if 
there are $\m\in\E$ and $a,b\in C^\times$  such that 
$D_{P_{\times\m}} = a\cdot (Y-b)^d$.
Suppose  $\wt P \geq 1$. Then \eqref{eq:asympt equ 1} is not unraveled if and only if 
$1\in\E$ and $D_P\in C[Y](Y')^\N$ has a zero
in $C^\times$ of multiplicity~$d$.
\end{example}

\noindent
If \eqref{eq:asympt equ 1} is unraveled, then so is
\eqref{eq:asympt equ 1phi} for each $\phi$; moreover:

%\begin{lemma}
%Suppose the asymptotic equation \eqref{eq:asympt equ 1} is 
%unraveled,
%and $f$ is an algebraic approximate solution to 
%\eqref{eq:asympt equ 1}. Then
%$$P_{+f}(Y)=0, \qquad Y\in\E$$
%is also unraveled.
%\end{lemma} \marginpar{lemma and proof checked but replaced by
%next result}
%\begin{proof}
%We have $\ndeg_\E P_{+f} = d > 0$.
%Let $g\in\E$; we need to show that
%$$\ndeg_{\prec g} P_{+(f+g)}<d.$$
%If $g\succ f$, then
%$$\ndeg_{\prec g} P_{+(f+g)} = \ndeg_{\prec g} P_{+g} < d$$
%since   \eqref{eq:asympt equ 1} is unraveled. If
%$g\prec f$, then similarly
%$$\ndeg_{\prec g} P_{+(f+g)} \leq  \ndeg_{\prec f} P_{+(f+g)} 
%= \ndeg_{\prec f} P_{+f} < d.$$
%Next suppose that $g\asymp f$. If $f+g\asymp g$, then
%$$\ndeg_{\prec g} P_{+(f+g)} = \ndeg_{\prec f+g} P_{+(f+g)} < d,$$
%and if $f+g\prec g$, then
%$$\ndeg_{\prec g} P_{+(f+g)} = \ndeg_{\prec g} P = \nval P_{\times f} %< \ndeg P_{\times f}\leq d,$$
%since $f$ is an algebraic approximate zero of $P$.
%\end{proof} 

\begin{lemma}\label{lem:approx unraveler}
Suppose the asymptotic equation \eqref{eq:asympt equ 1} is unraveled. Let  $f\in K$ be such that $f\preceq\frak e$ for some algebraic starting
monomial $\frak e$ for \eqref{eq:asympt equ 1}. Then
$$P_{+f}(Y)\ =\ 0, \qquad Y\in\E$$
is also unraveled.
\end{lemma}
\begin{proof}
We have $\ndeg_\E P_{+f} = d \geq 1$.
Let $g\in\E$; we need to show that
$$\ndeg_{\prec g} P_{+(f+g)}\ <\ d.$$
If $g\succ f$, then
$$\ndeg_{\prec g} P_{+(f+g)}\ =\ \ndeg_{\prec g} P_{+g} < d$$
since   \eqref{eq:asympt equ 1} is unraveled. If
$g\prec f$, then similarly
$$\ndeg_{\prec g} P_{+(f+g)}\ \leq\  \ndeg_{\prec f} P_{+(f+g)}\ =\ \ndeg_{\prec f} P_{+f}\ <\ d.$$
Next, suppose that $g\asymp f$. If $f+g\asymp g$, then
$$\ndeg_{\prec g} P_{+(f+g)}\ =\ \ndeg_{\prec f+g} P_{+(f+g)}\ <\  d,$$
and if $f+g\prec g$, then we take an algebraic starting monomial $\mathfrak e$ for \eqref{eq:asympt equ 1} with $f\preceq\mathfrak e$, and get from Corollary~\ref{dn5cor} that 
 $$\ndeg_{\prec g} P_{+(f+g)}\ =\ \ndeg_{\prec g} P\ =\ \nval P_{\times f}\ \leq\ \nval P_{\times\mathfrak e}\ <\  \ndeg P_{\times\mathfrak e}\ \le\ d,
 $$
so $\ndeg_{\prec g} P_{+(f+g)}< d$.
\end{proof} 

\noindent
Lemma~\ref{lem:approx unraveler} has a converse of sorts:

\begin{lemma}
Suppose $\Gamma$ is divisible and $\ndeg_{\E}P > \val P$. 
Let $\mathfrak e$ be the largest algebraic
starting monomial for \eqref{eq:asympt equ 1}. Suppose $\E_1\subseteq K^\times$ is $\preceq$-closed and 
$f\in\E_1$ is such that $\ndeg_{\E_1}P=\ndeg_{\E} P$ and
$P_{+f}(Y)=0, Y\in \E_1,$
%\eqref{eq:asympt equ 2} 
is unraveled. Then $f\preceq\mathfrak e$.
\end{lemma}
\begin{proof}
We have $d:=\ndeg_{\E} P=\ndeg_{\E_1} P_{+f}$. If $f\succ\mathfrak e$, then
$$d\ =\ \ndeg P_{\times \mathfrak{e}}\ \leq\ \ndeg_{\prec f} P\ =\ \ndeg_{\prec f} (P_{+f})_{+(-f)}\ <\  d,$$
a contradiction.
\end{proof}

\noindent
Assume $\ndeg_{\E} P = d\ge 1$.
Let $f\in \E\cup\{0\}$ and let $\E'\subseteq \E$ be
$\preceq$-closed.
We call the pair $(f,\E')$ a {\bf partial unraveler} for \eqref{eq:asympt equ 1} if  $\ndeg_{\E'} P_{+f} = d$. Thus $(f,\E)$ is a partial unraveler for \eqref{eq:asympt equ 1}. 
If $(f,\E')$ is a partial unraveler for \eqref{eq:asympt equ 1} and $(f_1,\E_1)$ is a partial unraveler for \eqref{eq:asympt equ 2}, then
$(f+f_1,\E_1)$ is a partial unraveler for \eqref{eq:asympt equ 1}.
Moreover, if  $(f,\E')$ is a partial unraveler for \eqref{eq:asympt equ 1},
then $(f,\E')$ is a partial unraveler for \eqref{eq:asympt equ 1phi}.
An {\bf unraveler} for \eqref{eq:asympt equ 1} is a partial unraveler 
$(f, \E')$ for \eqref{eq:asympt equ 1} with unraveled 
 \eqref{eq:asympt equ 2}. In the next easy lemma we continue assuming $\ndeg_{\E} P\ge 1$.

\begin{lemma}\label{lem:unraveler mult conj}
Let $a\in K^\times$ and $a\E :=\{ay\in K^\times:y\in\E\}$, and consider 
\begin{equation}\label{eq:asympt equ 1a} \tag{$a$E}
P_{\times a^{-1}}(Y)\ =\ 0,\qquad Y\in a\E.
\end{equation}
The Newton degree of \eqref{eq:asympt equ 1a} equals that of \eqref{eq:asympt equ 1}.
If 
$(f,\E')$ is a partial unraveler for~\eqref{eq:asympt equ 1}, then
$(af,a\E')$ is a partial unraveler for~\eqref{eq:asympt equ 1a}; 
similarly with {\rm unraveler} instead of {\rm partial unraveler}.
If $a\in \mathfrak{M}$, then the algebraic starting monomials for
~\eqref{eq:asympt equ 1a} are exactly the $a\mathfrak{e}$ with $\mathfrak{e}$ an
algebraic starting monomial for \eqref{eq:asympt equ 1}.
\end{lemma}

\index{asymptotic equation!partial unraveler}
\index{asymptotic equation!unraveler}
\index{unraveler}
\index{unraveler!partial}

\noindent
We also note an easy consequence of Lemma~\ref{lem:approx unraveler}:

\begin{cor}\label{cor:approx unraveler}
Assume $\ndeg_{\E}P\ge 1$. Let $(f,\E')$ be an unraveler for \eqref{eq:asympt equ 1}, and let
$g\in K$ be such that $f-g\preceq\mathfrak e$ for some algebraic starting monomial $\mathfrak e$ of
the refinement  \eqref{eq:asympt equ 2} of  \eqref{eq:asympt equ 1}. Then $(g,\E')$ is also an unraveler for \eqref{eq:asympt equ 1}.
\end{cor}

\noindent
Recall that throughout this section $K$ is $\upo$-free. In the next result we use the notion \textit{asymptotically $\d$-algebraically  maximal}\/ defined in
Section~\ref{As-Fields,As-Couples}.

\begin{prop}\label{prop:unravel}
Suppose $K$ is asymptotically $\d$-algebraically maximal, $\Gamma$ is divisible, $d:=\ndeg_{\E} P\ge 1$, but there is no $f\in\E\cup\{0\}$ with $\val P_{+f}=d$.
Then there exists an unraveler for \eqref{eq:asympt equ 1}.
\end{prop}
\begin{proof}
Let $\big((f_\lambda,\E_\lambda)\big)_{\lambda<\rho}$ be a sequence
of partial unravelers for \eqref{eq:asympt equ 1}, indexed by the ordinals
less than an ordinal $\rho>0$, such that $(f_0, \E_0)=(0,\E)$ and 
\begin{enumerate}
\item $\E_\lambda \supseteq \E_{\mu}$ for all $\lambda<\mu<\rho$,
\item $f_\mu-f_{\lambda} \succ f_\nu-f_\mu$ for all $\lambda<\mu<\nu<\rho$,  
\item  
$f_{\lambda+1}-f_{\lambda}\in\E_\lambda\setminus \E_{\lambda+1}$ for all $\lambda$ with $\lambda+1<\rho$.
\end{enumerate}
Note that for $\rho=1$ we have such a sequence. 
By (2) we have $f_\lambda-f_\mu \asymp f_\lambda-f_{\lambda+1}$ for
$\lambda<\mu<\rho$. Take $\fm_\lambda\in\fM$ with $\fm_\lambda\asymp  f_{\lambda+1}-f_{\lambda}$
for $\lambda+1<\rho$. Then by (3),
\begin{align*}
d\ =\ \ndeg_{\E_{\lambda+1}} P_{+f_{\lambda+1}}\		& \leq\  \ndeg_{\preceq\fm_\lambda} P_{+f_{\lambda+1}} \\
											& =\  \ndeg_{\preceq\fm_\lambda} (P_{+f_\lambda})_{+(f_{\lambda+1}-f_\lambda)} \\
											& =\  \ndeg_{\preceq\fm_\lambda} P_{+f_\lambda} \\
											& \leq\  \ndeg_{\E_{\lambda}} P_{+f_{\lambda}}\ =\  d,
\end{align*}
and thus $\ndeg_{\preceq\fm_\lambda} P_{+f_\lambda}=d$, for all ordinals $\lambda$ with $\lambda+1<\rho$.

Suppose first that $\rho$ is a successor ordinal, $\rho=\sigma+1$.
Consider the refinement
\begin{equation}\label{eq:asympt equ 3}\tag{E$_\sigma$}
P_{+f_\sigma}(Y)\ =\ 0, \qquad Y\in\E_\sigma
\end{equation}
of  \eqref{eq:asympt equ 1}. If \eqref{eq:asympt equ 3} is unraveled, then
$(f_{\sigma}, \E_{\sigma})$ is an unraveler of \eqref{eq:asympt equ 1}, and we are done. Assume \eqref{eq:asympt equ 3} is not unraveled, and take $f\in\E_\sigma$ such that
$\ndeg_{\prec f}\, (P_{+f_\sigma})_{+f} = d$.
The subset 
$$\E_\rho\ :=\ \{y\in K^\times : y\prec f\}$$
of $\E_\sigma$ is $\preceq$-closed,  with 
$\ndeg_{\E_\rho}(P_{+f_\sigma})_{+f} = d$.
Hence $(f,\E_\rho)$ is a partial unraveler for~\eqref{eq:asympt equ 3}, so
$(f_\rho,\E_\rho)$, where $f_\rho:=f_\sigma+f$, is a partial unraveler for~\eqref{eq:asympt equ 1}. Then
$\big((f_\lambda,\E_\lambda)\big)_{\lambda<\rho+1}$
satisfies (1)--(3) with $\rho+1$ instead of $\rho$.

Now suppose $\rho$ is a limit ordinal. Then $(f_\lambda)_{\lambda<\rho}$ is a pc-sequence in $K$. Let ${\bf f}=c_K(f_\lambda)$ be the corresponding cut in $K$.
Then $\ndeg_{\bf f} P = d$. The $\d$-valued field~$K$ of
$H$-type is asymptotically $\d$-algebraically maximal and has rational asymptotic integration, so we can take a pseudolimit $f_\rho$ of $(f_\lambda)_{\lambda<\rho}$ in $K$, by 
Lemmas~\ref{zda, newton} and~\ref{dpkell}. Consider the subset
$$\E_\rho\ :=\ \bigcap_{\lambda<\rho} \E_\lambda\ =\ \{ y\in K^\times:\  \text{$y\prec\fm_\lambda$ for all $\lambda<\rho$} \}$$
of $K^\times$. If $\E_\rho=\emptyset$, then $\val P_{+f_\rho} = \ndeg_{\bf f} P = d$ by Corollary~\ref{cor:ndegE}, contradicting the hypothesis of the proposition.
Thus $\E_\rho\neq\emptyset$, and so $\E_{\rho}$ is $\preceq$-closed. We have 
$\ndeg_{\E_\rho} P_{+f_\rho}=\ndeg_{\bf f} P =d$ by Corollary~\ref{cor:ndegE}, 
hence $(f_\rho,\E_\rho)$ is a partial unraveler  for \eqref{eq:asympt equ 1}, and $\big((f_\lambda,\E_\lambda)\big)_{\lambda<\rho+1}$
satisfies (1)--(3) with $\rho+1$ instead of $\rho$. 

This building process must end in producing an unraveler for \eqref{eq:asympt equ 1}.
\end{proof}

\subsection*{Behavior of unravelers under immediate extensions}
{\em In this subsection we fix an $\upo$-free immediate $H$-asymptotic extension $L$ of $K$.}\/ By Lemma~\ref{dv}, $L$ is $\d$-valued; our monomial group $\fM$ for $K$ continues to serve as a monomial
group for~$L$. Let $\E\subseteq K^\times$ be $\preceq$-closed. Then the subset
$$\E_L\ :=\ \{y\in L^\times:\ vy\in v\E\}$$ 
of $L^\times$ is $\preceq$-closed with $\E_L\cap K=\E$.
We consider the asymptotic equation
\begin{equation}\label{eq:asympt equ, L}\tag{E$_L$}
P(Y)\ =\ 0,\qquad Y\in\E_L
\end{equation}
over $L$, which has the same Newton degree $\ndeg_{\E_L} P=\ndeg_{\E} P$ as \eqref{eq:asympt equ 1}.
An element~$y$ of $K$ is an approximate solution of \eqref{eq:asympt equ 1} if and only if
$y$ is an approximate solution of~\eqref{eq:asympt equ, L}, and in this case the  multiplicity of $y$ 
as an approximate solution of~\eqref{eq:asympt equ 1} agrees with the  multiplicity of $y$ 
as an approximate solution of~\eqref{eq:asympt equ, L}.
Hence~\eqref{eq:asympt equ 1}  is unraveled if and only if \eqref{eq:asympt equ, L} is.
This yields:

\begin{lemma}\label{lem:unravelers under immediate exts}
Assume $\ndeg_{\E} P\ge 1$. Let $f\in \E\cup\{0\}$ and let $\E'\subseteq\E$ be $\preceq$-closed. Then 
$(f,\E')$ is a partial unraveler for  \eqref{eq:asympt equ 1} if and only if $(f,\E'_L)$ is a partial unraveler for 
\eqref{eq:asympt equ, L}; similarly with {\rm unraveler} instead of {\rm partial unraveler}.
\end{lemma}

\noindent
Let $(a_\rho)$ be a divergent pc-sequence in $K$ with
minimal $\d$-polynomial~$P$ over $K$, and suppose
$a_{\rho}\leadsto \ell\in L$.

\begin{lemma}\label{lem:small mult}
 $\val(P_{+\ell})\leq 1$.
\end{lemma}
\begin{proof}
By Lemma~\ref{flpc1new}
%the results in Section~\ref{sec:construct imm exts}, 
we have $Q(\ell)\neq 0$ for all $Q\in K\{Y\}^{\neq}$ of smaller complexity than $P$. Thus $S_P(\ell)\neq 0$, and so $\val(P_{+\ell})\leq 1$.
\end{proof}

\noindent
Let $\mathbf a=c_K(a_\rho)$ be the cut defined by $(a_\rho)$ in $K$,  
and let $a\in K$, $\mathfrak v\in K^\times$ be such that 
$a-\ell\prec\mathfrak v$ and 
$\ndeg_{\prec\mathfrak v} P_{+a}=\ndeg_{\mathbf a} P$. 
Consider the asymptotic equation 
\begin{equation}\label{eq:unravel add conj}
P_{+a}(Y) = 0, \qquad Y\prec\mathfrak v
\end{equation}
over $K$. In the next lemma we also consider it as an asymptotic equation over $L$.
The following consequence of Proposition~\ref{prop:unravel} is needed in Section~\ref{sec:newtonization} below:

\begin{lemma}\label{lem:unravel min pol}
Suppose that $\Gamma$ is divisible, $L$ is asymptotically $\d$-algebraically maximal, and $\ndeg_{\mathbf a} P\geq 2$.
Then there exists an unraveler $(f,\E)$ for~\eqref{eq:unravel add conj} over $L$ 
such that $f\ne 0$, $\ndeg_{\prec f} P_{+(a+f)}=\ndeg_{\mathbf a} P$, and $a_\rho\leadsto a+f+z$ for all $z\in\E\cup\{0\}$.
\end{lemma}
\begin{proof}
Take $g\in K$ such that $a-\ell\sim -g$. Then $0\ne g\prec \fv$, so
$$\ndeg_{\mathbf a} P\ =\ \ndeg_{\prec \fv}P_{+a}\ =\ \ndeg_{\prec \fv} P_{+(a+g)}\ \ge\ \ndeg_{\prec g} P_{+(a+g)}.$$
Also 
$(a+g)-\ell \prec g$, so $\ndeg_{\mathbf a} P\le \ndeg_{\prec g} P_{+(a+g)}$ by  Lemma~\ref{dpkell}, and thus $\ndeg_{\mathbf a} P= \ndeg_{\prec g} P_{+(a+g)}$. Now $P_{+(a+g)}$ is a minimal differential polynomial for $\big(a_\rho-(a+g)\big)$ over $K$
and $\ndeg_{\mathbf a-(a+g)} P_{+(a+g)} = \ndeg_{\mathbf a} P$.
Suppose $\E\subseteq L^{\times,\prec g}$ is $\preceq$-closed in the sense of $L$, and $(h,\E)$ is an unraveler for the asymptotic equation
$$P_{+(a+g)}(Y)\ =\ 0,\qquad Y\prec g$$
over $L$ with $a_\rho-(a+g)\leadsto h+z$ for all $z\in\mathcal E\cup\{0\}$. Then $(f,\E)$, $f:=g+h\ne 0$, is an unraveler for~\eqref{eq:unravel add conj} in $L$ 
with $\ndeg_{\prec f} P_{+(a+f)}=\ndeg_{\mathbf a} P$ and $a_\rho\leadsto a+f+z$ for all $z\in\mathcal E\cup\{0\}$. 
Thus, after replacing $P$,~$(a_\rho)$,~$\ell$,~$\mathfrak v$ by 
$P_{+(a+g)}$, $\big(a_\rho-(a+g)\big)$, $\ell-(a+g)$, $g$, respectively, we may assume $a=0$, and only need to show
the existence of an unraveler $(f,\E)$ for~\eqref{eq:unravel add conj} in $L$ 
with $a_\rho\leadsto f+z$ for all $z\in\E\cup\{0\}$. For this, consider the subset
$$\mathcal Z\ :=\ \{z\in L^\times:\ \text{$z\prec a_\rho-\ell$, eventually}\}$$
of $L^\times$. By 
Lemma~\ref{lem:small mult} and since $\ndeg_{\mathbf a}P\geq 2$, there is no $z\in \mathcal Z\cup\{0\}$
such that $\val(P_{+(\ell+z)})=\ndeg_{\mathbf a} P$. So $\mathcal Z\neq\emptyset$, and thus $\mathcal Z$ is $\preceq$-closed, and
$\ndeg_{\mathcal Z} P_{+\ell} = \ndeg_{\mathbf a} P$, by Corollary~\ref{cor:ndegE}. Proposition~\ref{prop:unravel} now provides us with an unraveler~$(g,\E)$ for the asymptotic equation
$$P_{+\ell}(Y)\ =\ 0,\qquad Y\in\mathcal Z$$
over $L$. Then $(f,\E)$, $f:=\ell+g$, is an unraveler for \eqref{eq:unravel add conj} (with $a=0$), and $a_\rho\leadsto f+z$ for all $z\in \mathcal Z\cup\{0\}$.
\end{proof}

\subsection*{Neglecting terms of high degree} {\em In this subsection we assume $\Gamma$ is divisible and $d:= \ndeg_{\E} P \ge 1$}. 
Let $(f,\E')$ be an unraveler for~\eqref{eq:asympt equ 1}, and set  $\mathfrak f:=\mathfrak d_f$. Suppose also that $d >\ \val(P_{+f})$. 
Then \eqref{eq:asympt equ 2} has an algebraic starting monomial
by Proposition~\ref{newtpolygon}, and 
we let  $\mathfrak e$ be the largest algebraic starting monomial for~\eqref{eq:asympt equ 2}.
Let $g\in K$, $\mathfrak g:=\mathfrak d_g$, and suppose $\mathfrak e\prec \mathfrak g\prec \mathfrak f$.
Put $\tilde f:=f-g$ (so $\tilde f\sim f$), and consider the refinement
\begin{equation}\label{eq:Etilde}\tag{$\tilde{\operatorname{E}}$}
P_{+\tilde f}(Y)\ =\ 0,\qquad Y\preceq  \mathfrak g
\end{equation}
of \eqref{eq:asympt equ 1}.
Set $\tilde\E':=\{y\in\E':\ y \prec\mathfrak g\}$, so $\mathfrak e\in \tilde\E'$. 

\begin{lemma}\label{Etildeunr}
The asymptotic equation \eqref{eq:Etilde} has Newton degree $d$, and
$(g,\tilde\E')$ is an unraveler for  \eqref{eq:Etilde}.
\end{lemma}
\begin{proof}
We have
$$d\, =\, \ndeg_{\preceq\mathfrak e} P_{+f}\, \leq\,  \ndeg_{\preceq\mathfrak g} P_{+f}\, =\, \ndeg_{\preceq\mathfrak g} P_{+\tilde f}\ \leq\ \ndeg_{\E} P_{+\tilde f}\, =\, \ndeg_{\E} P\, =\, d$$
and hence \eqref{eq:Etilde} has Newton degree $d$.
Also
$$d\ =\ \ndeg_{\preceq\mathfrak e} P_{+f}\ \leq\ \ndeg_{\tilde\E'} P_{+f}\ \leq\ \ndeg_{\E} P_{+ f}\ =\ \ndeg_{\E} P\ =\ d.$$
Hence the asymptotic equation
$$P_{+f}(Y)\ =\ 0,\qquad Y \in \tilde\E',$$
which is a refinement of both \eqref{eq:asympt equ 2} and \eqref{eq:Etilde}, has Newton degree $d$, and as
\eqref{eq:asympt equ 2} is unraveled, the pair $(g,\tilde\E')$ is an unraveler for  \eqref{eq:Etilde}.
\end{proof}

\noindent
Recall that for $F\in K\{Y\}$ and $e\in\N$, we defined $F_{\leq e}:=F_0+F_1+\cdots+F_e$. Note that if $e\geq\ndeg F$, then $N_F=N_{F_{\leq e}}$ by Corollary~\ref{npsum} and its proof.
Set $F:= P_{+\tilde f}$. Then $d\ge \ndeg F_{\times \fm}$ for all
$\fm\preceq \mathfrak g$. Consider the ``truncation'' 
\begin{equation}\tag{$\tilde{\operatorname{E}}_{\leq d}$}\label{eq:Etildeleqd}
F_{\leq d}(Y)\ =\ 0,\qquad Y\preceq\mathfrak g
\end{equation}
of \eqref{eq:Etilde} as an asymptotic equation over $K$.
We have 
$$N_{F_{\times\fm}}\ =\ N_{(F_{\times\fm})_{\leq d}}\ =\ N_{(F_{\leq d})_{\times\fm}}
\qquad\text{for $\fm\preceq\mathfrak g$,}
$$ 
so \eqref{eq:Etildeleqd} has the same algebraic starting monomials and the same
Newton degree~$d$ as~\eqref{eq:Etilde}. In the next lemma we show that under suitable conditions the unraveler~$(g,\tilde\E')$ for~\eqref{eq:Etilde} is also an unraveler for~\eqref{eq:Etildeleqd}. This will be crucial in Section~\ref{sec:unravelers}.

\begin{lemma} \label{lem:neglect >d}
Suppose
$(\mathfrak e/\mathfrak g)\flatter (\mathfrak g/\mathfrak f)$.
Then  
$(g,\tilde\E')$ is an unraveler for \eqref{eq:Etildeleqd}, and~$\mathfrak e$ is the largest algebraic starting monomial
of the unraveled asymptotic equation
\begin{equation}\label{eq:neglect >d}\tag{$\tilde{\operatorname{E}}{}'_{\leq d}$}
\big(F_{\leq d}\big)_{+g}(Y)\ =\ 0, \qquad Y\in\tilde\E'.
\end{equation}
\end{lemma}
\begin{proof}
For $(g,\tilde\E')$ to be an unraveler for \eqref{eq:Etildeleqd} it is
enough to show:
\begin{enumerate}
\item $\ndeg_{\tilde\E'} (F_{\leq d})_{+g}=d$;
\item $\ndeg_{\prec h} (F_{\leq d})_{+(g+h)}<d$ for all $h\in\tilde\E'$.
\end{enumerate}
Until further notice we assume $\mathfrak g=1$, so $\mathfrak e\prec 1\prec\mathfrak f$, $\mathfrak e\flatter\mathfrak f$. At the end we reduce to this special case. We have
$\ndeg F=\ndeg_{\preceq \mathfrak g}F=d$ by an equality from the proof
of Lemma~\ref{Etildeunr}, so 
$d\le \ndeg F_{\times \mathfrak f}\le \ndeg_{\E} F =d$, and thus
$\ndeg F=\ndeg F_{\times\mathfrak f}=d$.  
If $\fm\flatter\mathfrak f$,  then by Corollary~\ref{cor:Newton poly, precfn} with $F, \mathfrak f$ in the role of $P, \fn$ we have
$$N_{P_{+f,\times\fm}}\ =\ N_{F_{+g,\times\fm}}\ =\ N_{(F_{\leq d})_{+g,\times\fm}}. $$
In particular, if $\mathfrak e\preceq\fm\prec 1$, then
$N_{P_{+f,\times\fm}} = N_{(F_{\leq d})_{+g,\times\fm}}$, and so
$\mathfrak e$ is the largest algebraic starting monomial
of \eqref{eq:neglect >d}. Also, 
$$\mathfrak e\preceq\fm\in \tilde\E'\ \Longrightarrow\ 
\ndeg (F_{\leq d})_{+g,\times\fm}\ =\ \ndeg P_{+f,\times\fm}\ =\ d,$$ so (1) holds.
For~(2), let $h\in\tilde\E'$, so $h\in \E'$, $h\prec 1$, and 
$h\sim c\,\mathfrak h$ with $c\in C^\times$, $\mathfrak h:=\mathfrak d_{h}$. Applying Lemma~\ref{transition} twice gives
\begin{align*} \ndeg_{\prec h} (F_{\le d})_{+(g+h)}\ &=\ \val\!\big(N_{(F_{\leq d})_{+g,\times\mathfrak h}}\big)_{+c}\\
 \ndeg_{\prec h}P_{+(f+h)}\ &=\ \val\!\big(N_{P_{+f,\times\mathfrak h}}\big)_{+c}.
\end{align*}
Now \eqref{eq:asympt equ 2} is unraveled, so if $\mathfrak e\preceq\mathfrak h$, then 
$\ndeg_{\prec h}P_{+(f+h)}<d$, and thus 
$$\ndeg_{\prec h} (F_{\le d})_{+(g+h)}\ <\ d$$ by
combining various equalities above. 
If $\mathfrak{e}^2\preceq \mathfrak h\prec\mathfrak e$, then
$\mathfrak{h}^\dagger\asymp \mathfrak{e}^\dagger$, so
$\mathfrak{h} \flatter  \mathfrak{f}$, hence 
$$\ndeg\ (F_{\leq d})_{+g,\times\mathfrak h}\ =\ 
\ndeg\ P_{+f,\times\mathfrak h}\ <\ \ndeg\ P_{+f,\times\mathfrak e}\ =\ d,$$
as $\mathfrak e$ is the largest algebraic starting monomial for \eqref{eq:asympt equ 2}. 
Thus if
$\mathfrak h\prec\mathfrak e$, then
$$\ndeg_{\prec h} (F_{\le d})_{+(g+h)}\ =\ \val \big(N_{(F_{\leq d})_{+g,\times\mathfrak h}}\big)_{+c}\ \leq\ \ndeg\ (F_{\leq d})_{+g,\times\mathfrak h}\  <\  d.$$
This gives (2) when $\mathfrak g=1$. To reduce to the case $\mathfrak g=1$, replace $P$, $f$, $g$, $\E$, $\E'$ by $P_{\times\mathfrak g}$, $f/\mathfrak g$, $g/\mathfrak g$, 
$\mathfrak{g}^{-1}\E$, $\mathfrak{g}^{-1}\E'$, respectively, 
and use Lem\-ma~\ref{lem:unraveler mult conj}.
\end{proof}

\section{Some Special $H$-Fields} \label{sec:specialH}

\noindent
This section will not be used later in this volume but is included for its intrinsic interest. We assume familiarity with Appendix~\ref{app:trans}. We construct here $H$-fields,
$$ \R\<\upo\>\ \subseteq\ \R\<\upl\>\ \subseteq\ \R\<\upg\>\  \text{ (with $H$-field inclusions)},$$
each generated as a differential field over their common 
constant field $\R$
by a single element, where $\R\<\upo\>$ is $\upl$-free but not $\upo$-free, $\R\<\upl\>$ is not $\upl$-free but has rational asymptotic integration, and $\R\<\upg\>$
has a gap $v\upg$. These $H$-fields and their asymptotic couples have certain canonical features that are worth documenting. Moreover, they can be realized as Hardy fields as we show at the end of this section.

\subsection*{The ambient $H$-field $\mathbb{L}$}   
Let $\mathfrak{L}_n:= \ell_0^{\R}\cdots \ell_n^{\R}$ be the subgroup of the ordered multiplicative group
$\mathbb{T}^{>}$ generated by the real powers of the iterated logarithms $\ell_i$ for $i=0,\dots,n$. This yields ordered group inclusions 
$$ \mathfrak{L}_0\ \subseteq\ \mathfrak{L}_1\ \subseteq\ \mathfrak{L}_2\ \subseteq\ \cdots \subseteq\ \bigcup_n \mathfrak{L}_n\ =:\ \mathfrak{L}\ \subseteq\  G^{\operatorname{LE}}.$$
We view $\mathbb{L}:=\R[[\mathfrak{L}]]$ in the natural way
as an ordered subfield of the ordered Hahn field $\R[[G^{\operatorname{LE}}]]$.  The latter also contains $\mathbb{T}$ as an ordered subfield, with 
$$ \mathbb{L} \cap \mathbb{T}\ =\ 
\bigcup_n \R[[\mathfrak{L}_n]]\ =\ \mathbb{T}_{\log}.$$ 
We equip $\mathbb{L}$ with the unique strongly additive $\R$-linear derivation such that 
$$ (\ell_0^r)'\ =\ r\,\ell_0^{r-1}, \qquad (\ell_{n+1}^r)'\ =\ r\,\ell_{n+1}^{r-1}(\ell_0\cdots \ell_n)^{-1}\qquad (r\in \R).$$
This makes $\mathbb{L}$ a spherically complete immediate real
closed $H$-field extension of~$\mathbb{T}_{\log}$, with constant field~$\R$. Thus $\mathbb{L}$ and $\mathbb{T}_{\log}$ have the same asymptotic couple~$\big(v(\mathfrak{L}), \psi\big)$. Moreover, $v(\mathfrak{L})$  
is an ordered vector space
over $\R$. We set $e_n:=v(\ell_n)$, so $e_n<0$, $[e_{n}]>[e_{n+1}]$, and $v(\mathfrak{L}) = \bigoplus_{n} \R\, e_n$ (internal direct sum) with
$$ \big[v(\mathfrak{L})^{\ne}\big]\ =\ \big\{[e_n]:\ n=0,1,2,\dots\big\}, \qquad e_n^\dagger\ =\ -(e_0+e_1+\cdots+e_n).$$
Note that $\Psi:=\psi\big(v(\mathfrak{L})^{\ne}\big)$ has no supremum in $v(\mathfrak{L})$, so the $H$-asymptotic couple~$\big(v(\mathfrak{L}), \psi\big)$ has rational asymptotic integration.  

\subsection*{The elements $\upl$ and $\upo$ of
$\mathbb{L}$} Let $\mathfrak M$ be the subgroup of $\mathfrak L$ generated by $\ell_0,\ell_1,\dots$, so
$\mathfrak M= \bigcup_n \ell_0^{\Z}\cdots \ell_n^{\Z}$. The ordered subfield $\R[[\mathfrak M]]$ of $\mathbb L$ is closed under the derivation of $\mathbb L$, which makes it an $H$-subfield of $\mathbb L$. 
Now $\R[[\mathfrak M]]$ has special elements 
$$\sum_{n=1}^{\infty}\ell_n,\quad \upl\, :=\, \left(\sum_{n=1}^{\infty}\ell_n\right)'\, =\,\sum_{n=0}^{\infty}(\ell_0\cdots\ell_n)^{-1},\ \quad \upo\, :=\, \omega(\upl)\, =\, \sum_{n=0}^{\infty}(\ell_0\cdots\ell_n)^{-2}, $$
none lying in $\mathbb{T}_{\log}$. This gives the
$H$-subfields $\R\<\upo\>$ and $\R\<\upl\>$ of 
$\R[[\mathfrak M]]$.

\begin{prop}\label{propimmuplupo} $\R[[\mathfrak M]]$ is an immediate extension 
of $\R\<\upo\>$ and of $\R\<\upl\>$. Al\-so,~$\R\<\upo\>$ is 
$\upl$-free and not $\upo$-free. 
\end{prop} 

\noindent 
Towards the proof, first note that the asymptotic couple $\big(v(\mathfrak{L}),\psi\big)$ of $\mathbb{L}$ extends the
asymptotic couple $\big(v(\mathfrak{M}),\psi\big)$ of 
$\R[[\mathfrak M]]$, with $v(\mathfrak{M})=\bigoplus_{n} \Z\, e_n$. It follows that $\big[v(\mathfrak{M})\big]=\big[v(\mathfrak{L})\big]$, the two
asymptotic couples have the same $\Psi$-set, namely $\Psi$, and the $H$-asymptotic couple $\big(v(\mathfrak{M}), \psi\big)$ has rational asymptotic integration.

\begin{lemma}\label{G=Gamma}
Let $G\neq\{0\}$ be an ordered subgroup of $v(\mathfrak{M})$ such that
$\psi(G^{\ne})\subseteq G$ and
$G^{<}$ is coinitial in $v(\mathfrak{M})^{<}$. Then $G=v(\mathfrak{M})$.
\end{lemma}
\begin{proof} Note: if $[e_n]\in [G]$, then 
$e_n^\dagger\in G$.
Suppose $m<n$ and $[e_m], [e_n]\in [G]$. Then $-e_m^\dagger=e_0+\cdots+e_m\in G$ and $-e_n^\dagger=e_0+\cdots+e_n\in G$, hence 
$$e_{m+1}+\cdots+e_n=e_m^\dagger-e_n^\dagger\in G,$$ so 
$[e_{m+1}]=[e_{m+1}+\cdots+e_n]\in [G]$, and thus $e_{m+1}^\dagger\in G$. Inductively it follows that~$G$ contains $e_m^\dagger,
e_{m+1}^\dagger,\dots,e_n^\dagger$ and hence $G$ contains 
$$e_{m+1}=e_{m}^\dagger-e_{m+1}^\dagger,\ e_{m+2}=e_{m+1}^\dagger-e_{m+2}^\dagger,\ \dots,\ e_n=e_{n-1}^\dagger-e_{n}^\dagger.$$
Take $m$ with $[e_m]\in [G]$. Then
$-\psi(e_m^\dagger)=e_0\in G$.
Thus $e_n\in G$ for all $n$.
\end{proof}

\noindent
Since $\mathbb{T}_{\log}$ is $\upo$-free, it follows from Proposition~\ref{cor:adding upo preserves upl-free} that the $H$-subfield $\mathbb{T}_{\log}\<\upo\>$ of $\mathbb{L}$ is $\upl$-free. In order to conclude that 
$\R\<\upo\>$ is $\upl$-free, it is clearly enough to get
$v\big(\R\<\upo\>^\times\big)=v(\fM)$, and that is part of Corollary~\ref{cor:uplgamma} below.  The $H$-field $\mathbb{L}$ has asymptotic integration, with 
  divisible value group, and $\upl_n \leadsto \upl$, so
  $-\upl$ creates a gap over $\mathbb{L}$, by Lemma~\ref{cg2}.
We take some Liouville closed $H$-field extension~$L$ of $\mathbb{L}$ and take $\upg\in L^{>}$ with 
$\upg^\dagger=-\upl$. Then $v\upg$ is a gap in 
$\mathbb{L}\<\upg\>$ by the remark following the proof of Lemma~\ref{cg2}. 
We use this gap to prove:

\begin{lemma}\label{lem:uplgamma} $v\big(\Q\<\upl\>^\times\big)=v(\mathfrak{M})$, and
$\R[[\mathfrak M]]$ is an immediate extension 
of $\R\<\upl\>$.
\end{lemma}
\begin{proof} Take $z\in L$ with $z'=\upg$.
Subtracting a constant in $L$ from $z$ we arrange
$z\not\asymp 1$, and then 
$v(z^\dagger) > \Psi$. 
Let $\Delta$ be the value group of $\Q\<z\>$. Then $\psi_L(\Delta^{\ne})\subseteq \Delta$. We apply 
Lemma~\ref{Deltacontract} to $\Delta$ as a subgroup of the value group of $L$. From $\upl\in \Q\<z\>$ we get
$v(1/\upl)=v(\ell_0)\in \Delta^{<}$, and as
$v(\ell_0^\dagger) < v(z^\dagger)\in \psi_L(\Delta^{\ne})$, that lemma yields
$$v(\ell_1)\ =\ \chi\big(v(\ell_0)\big)\in \Delta.$$
As $v(\ell_1^\dagger) < v(z^\dagger)$, we likewise get 
$v(\ell_2)\in \Delta$. Continuing this way we get 
$v(\ell_n)\in \Delta$ for all $n\ge 1$. Hence
$[\Delta]$ is infinite, and since
$\operatorname{trdeg}\!\big(\Q\<z\>|\Q\<\upl\>\big)\le 2$, also $\big[v(\Q\<\upl\>^\times)\big]$
is infinite. Thus $v\big(\Q\<\upl\>^\times\big)=v(\mathfrak{M})$ by Lemma~\ref{G=Gamma}.
\end{proof}

\begin{cor}\label{cor:uplgamma} Let $E$ be a differential subfield of 
$\R\<\upl\>$ not contained in $\R$. Then $v(E^\times)=v(\mathfrak{M})$, and $E$ as a pre-$H$-subfield of 
$\R\<\upl\>$ is not $\upo$-free. 
\end{cor}
\begin{proof} Take $a\in E$, $a\notin \R$. Then $\upl$ is
$\d$-algebraic over $\R\<a\>$ by Lemma~\ref{dtrdalg}. By 
Lemma~\ref{lem:uplgamma}, the set
$\big[v(\R\<\upl\>^\times)\big]$ is infinite, so 
$\big[v(\R\<a\>^\times)\big]$ is infinite, hence
$\big[v(\Q\<a\>^\times)\big]$ is infinite by 
Corollary~\ref{cor:dconstantsval}
and Lemma~\ref{archextclass}, and thus $v(E^\times)=v(\mathfrak{M})$ by 
Lemma~\ref{G=Gamma}. Since $\R\<\upl\>$ is not $\upl$-free, it is not $\upo$-free, and as $\R\<\upl\>$ is $\d$-algebraic over $E$, Theorem~\ref{upoalgebraic} yields that
$E$ is not $\upo$-free. 
\end{proof}

\noindent
Applying this to $E=\R\<\upo\>$ yields Proposition~\ref{propimmuplupo}.

\subsection*{Properties of $\R\<\upg\>$} Just before Lemma~\ref{lem:uplgamma} we introduced a pre-$H$-field extension $\mathbb{L}\<\upg\>$ of $\mathbb{L}$. It is generated as a differential field over $\mathbb{L}$ by an element~${\upg>0}$ with $\upg^\dagger=-\upl$. Since 
$\mathbb{L}\<\upg\>=\mathbb{L}(\upg)$, as fields, and $v\upg\notin v(\mathfrak{L})$, we have
$v\big(\mathbb{L}(\upg)^\times\big) = v(\mathfrak{L}) \oplus \Z v\upg$ (internal direct sum) by Lemma~\ref{lem:lift value group ext}. It follows that $\mathbb{L}(\upg)$ has the same residue field as $\mathbb{L}$, and so it is an $H$-field with the same constant field~$\R$ as $\mathbb{L}$. We think of $\upg$ informally as an infinite product,
$$\upg\ =\ \exp\!\left({-\sum_{n=1}^{\infty}\ell_n}\right)\ =\  1/\ell_0\ell_1\ell_2\cdots,$$
which suggests $v\upg= -(e_0 + e_1 + e_2 + \cdots)$. At this point we attach no formal meaning to these identities, but
they suggest other identities that do have meaning and that are easy to prove.
For example, let
$\alpha=\big(\sum_{i=0}^\infty r_i e_i\big) + k v\upg\in v\big(\mathbb{L}(\upg)^\times\big)$, $\alpha\ne 0$,
where $r_i\in\R$ for $i=0,1,\dots$, $r_i=0$ for all but finitely many $i$, and $k\in\Z$. Then in the asymptotic couple of $\mathbb L\<\upg\>$ we have
\begin{equation}\label{eq:alphadagger}
\alpha^\dagger\ =\ e_m^\dagger\ =\ -(e_0+\cdots+e_m)\quad\text{where $m=\min\{i\in\N: r_i\neq k\}$.}
\end{equation}
This is because in $\mathbb L\<\upg\>$ we have:
$$\left( \upg^k\prod_{i=0}^\infty \ell_i^{r_i} \right)^\dagger \ =\ k\upg^\dagger + \sum_{i=0}^\infty r_i\ell_i^\dagger\ =\ \sum_{i=0}^\infty (r_i-k) \frac{1}{\ell_0\ell_1\cdots\ell_i}.$$
We now turn to the $H$-subfield $\R\<\upg\>=\R\<\upl\>(\upg)$ of 
$\mathbb{L}\<\upg\>$. Since $v\upg$ is a gap in~$\mathbb{L}\<\upg\>$, it is a gap in $\R\<\upg\>$. 
The asymptotic couple $\big(v(\mathfrak{M}\big), \psi)$ of $\R\<\upl\>$ has rational asymptotic integration, so $v\big(\R\<\upg\>^\times\big)=v(\fM)\oplus\Z v\upg$ (internal direct sum) by
Corollaries~\ref{cor:ZA inequality} and~\ref{addgapmore}, and
$\alpha^\dagger$ for $0\neq\alpha\in   v\big(\R\<\upg\>^\times\big)$ is given by \eqref{eq:alphadagger}.

\subsection*{Realizing $\R\<\upg\>$ as a Hardy field} The $H$-field $\R\<\upg\>$ is isomorphic over $\R$ to a Hardy field extension of $\R$, and thus the same holds for the
$H$-subfields $\R\<\upo\>$ and~$\R\<\upl\>$ of $\R\<\upg\>$. 
To see this, recall that~$\mathcal G$ is the ring of germs at $+\infty$
of one-variable real-valued functions defined on half-lines~$(a,+\infty)$, $a\in\R$.
Define $l_n, e_n\in \mathcal G$ by recursion on $n$ such that
$l_0(t)=e_0(t)=t$, $l_{n+1}(t)=\log l_n(t)$, and $e_{n+1}(t)=\exp e_n(t)$, eventually. 
Every Hardy field extends to one that contains all $l_n$, $e_n$, and all real numbers.
Boshernitzan~\cite{Boshernitzan86} constructs a Hardy field \marginpar{result from \cite{Boshernitzan86} taken on faith} with an element~$e_\omega$ such that for every $n$,
eventually $e_\omega(t)>e_n(t)$. In particular, $e_{\omega}$ is eventually strictly increasing and $e_{\omega}(t) \to +\infty$ as $t\to +\infty$. Let~$l_\omega$ be the inverse of $e_{\omega}$: the germ in $\mathcal G$ such that $l_\omega\big(e_\omega(t)\big)=t$ eventually. By~\cite{SalvyShackell} 
(see also~\cite[Theo\-rem~1.7]{AvdD4}),~$l_\omega$ lies in a Hardy field extension of
$\R$, and for each $n$ and~$r\in \R$ we have eventually
 $r<l_\omega(t)<l_n(t)$. Then $g:=l_\omega^\dagger$ lies in the same Hardy field extension of~$\R$. 

\begin{lemma}\label{roslem}
There is an isomorphism $\R\<\upg\>\to\R\<g\>$ of ordered differential fields which is the identity on $\R$ and sends $\upg$ to $g$. \textup{(}Here $\R\<g\>$ is a Hardy field.\textup{)}
\end{lemma}
\begin{proof}
Consider the $H$-fields 
$$E\ :=\ \R(\ell_0,\ell_1,\dots)\ \subseteq\ \R[[\fM]], \qquad F\ :=\ \R(l_0,l_1,\dots),$$ the latter a Hardy field extension of $\R$. Lemma~\ref{pre-extas2} and the remarks following it yield an isomorphism $h\colon E\to F$ of ordered differential fields that is the identity on $\R$ and sends~$\ell_n$ to~$l_n$, for each $n$. 
Since $E$ is $\d$-valued of $H$-type with small derivation and $\upo$-free, and has monomial group  $\fM$, both $E$ and $F$ may be used in place of $K$ in Section~\ref{cutsvalgrp}. 
Since $\Psi_E<v\upg<(\Gamma_E^>)'$ and $\Psi_F<vg<(\Gamma_F^>)'$, it follows from Lemma~\ref{zstargap} that~$h$ extends to an
isomorphism $E\<\upg\>\to  F\<g\>$ of valued differential fields sending $\upg$ to $g$. Since
$\upg>0$ and $g>0$, the proof of Corollary~\ref{corzstargap}
shows that this isomorphism is also order preserving. 
\end{proof}

\subsection*{Notes and comments} Rosenlicht~\cite[p.~831]{Rosenlicht6}
states that he has no example of a
Hardy field extension $K\subseteq K(u)$ such that $u\neq 0$, $u^\dagger\in
K$, and $v(u^\dagger)$ is a gap in $K$, except when~$K\subseteq\R$. 
We note here that $K=\R\<g\>$ as in Lemma~\ref{roslem} and $u=\ex^{\int g}$ in a
Hardy field extension of~$K$ furnish such an example.

%\noindent
%Does every Hardy field extend to one into which 
%$\R\<\upg\>$ can be embedded over $\R$?  (See 
%\cite[Question~4, p.~279]{Boshernitzan86}.)

%% file: mt-14.tex
\chapter{Newtonian Differential Fields}
\label{ch:newtonian fields}

\setcounter{theorem}{0}

\noindent
{\em In this chapter $K$ is an ungrounded $H$-asymptotic field with $\Gamma:=v(K^\times)\ne \{0\}$.}\/ So the subset
$\Psi$ of $\Gamma$ is nonempty and has no largest element, and thus
$K$ is pre-$\d$-valued by Corollary~\ref{asymp-predif}. We let $\phi$ range over the active elements of
$K$. Since $\Psi^{\phi}$ contains positive elements, the differential residue field $\k^{\phi}$ of $K^{\phi}$ is
just the residue field~$\k$ of the valued field $K$ with the
trivial derivation. As $\k^{\phi}$ does not depend on~$\phi$, we 
  let~$\k$ stand for $\k^{\phi}$.  
We also fix a ``monomial'' set $\frak{M}\subseteq K^\times$ that is
mapped bijectively by $v$ onto $\Gamma$ and which gives us
the dominant monomial $\frak{d}_P\in \frak{M}\cup\{0\}$ and the dominant part $D_P\in \k\{Y\}$ of any 
$P\in K^{\phi}\{Y\}$. (We do not require $\frak{M}$ to be a monomial group of $K$; such a group might not even exist.)

By an {\em extension\/} of $K$ we mean an $H$-asymptotic
field extension of $K$.

\index{differential polynomial!quasilinear}
\index{quasilinear!differential polynomial}
\index{asymptotic field!newtonian}
\index{newtonian}
\medskip\noindent
Our main interest is in the case that $K$ has asymptotic integration, but then $K$ cannot be $\d$-henselian in the sense of Chapter~\ref{sec:dh1},
by Corollary~\ref{dhns3}. The correct notion in that situation, called {\em newtonian,}\/ is an
eventual variant of {\em $\d$-henselian.}\/ We define
$P\in K\{Y\}$ to be {\bf quasilinear\/} if $\ndeg P =1$, and
we define $K$ to be {\bf newtonian\/} if every quasilinear $P\in K\{Y\}$ has a zero in the valuation ring $\mathcal{O}$ of~$K$.
If $K$ is newtonian, then $K$ is henselian as a valued field, by Lemma~\ref{lem:char henselian}, and $K^{f}$ is newtonian
for every~$f\in K^\times$.  
In Section~\ref{sec:reldifhens} we show that for $\upl$-free $K$,
$$\text{$K$ is newtonian}\ \Longleftrightarrow\ \text{for each $\phi$ the flattening of $K^{\phi}$ is $\d$-henselian,}$$
and derive some consequences of this link between newtonianity and differential-hen\-sel\-ianity. 
In Section~\ref{cordone} we consider weak forms of newtonianity
and apply this to $P\in K\{Y\}$ of low complexity. Among results needed later we show there that if~$K$ is a newtonian Liouville closed $H$-field, then 
the subset $\omega(K)$ of $K$ is downward closed, and the subset $\sigma\big(\Upg(K)\big)$ of $K$ is upward closed. 

In Section~\ref{sqe} we prove newtonian versions of 
$\d$-henselian results in Chapter~\ref{sec:dh1}, leading 
to the following important analogue of Theorem~\ref{thm:damdh}:

\begin{theorem}\label{maxnewt} If $K$ is $\upl$-free and asymptotically $\d$-algebraically maximal, 
then~$K$ is $\upo$-free and newtonian. 
\end{theorem}

\noindent
One (minor) part of this is an immediate consequence of Corollary~\ref{cor:adjoin upl}: if
$K$ is $\upl$-free and asymptotically $\d$-algebraically maximal, then $K$ is $\upo$-free. 
Note also that by Zorn any $\upo$-free
$K$ has an immediate asymptotically $\d$-algebraically maximal $\d$-algebraic extension, and that by Section~\ref{consupofree} any such extension is also 
$\upo$-free, and thus newtonian by Theorem~\ref{maxnewt}.

 The main result of this chapter is almost a converse to Theorem~\ref{maxnewt}:

\begin{theorem}\label{newtmax} If $K$ is $\upo$-free and newtonian with divisible value group, then~$K$ is asymptotically $\d$-algebraically maximal.
\end{theorem}

\noindent
This is key to eliminating quantifiers for
$\T$ in Chapter~\ref{ch:QE}. After a rather technical
Section~\ref{sec:unravelers} on unraveling, we prove Theorem~\ref{newtmax} in 
Section~\ref{sec:newtonization}.

\section{Relation to Differential-Henselianity}\label{sec:reldifhens}

\noindent
We let $a$,~$b$,~$y$ range over $K$, and $P$ over $K\{Y\}^{\ne}$. Recall that throughout this chapter $\phi$ ranges over the active elements of $K$.  

\begin{lemma}\label{newtp1} Suppose $\Psi^{>0}\ne \emptyset$, 
$P$ has order $\le r$, 
$v(P_1)=0$ and $v(P_i)> r\Psi$ for all $i\ne 1$. Then $D_{P^{\phi}}=D_{P_1^{\phi}}$ for all $\phi\preceq 1$, and thus 
$\ndeg P =1$. 
\end{lemma}
\begin{proof} For $\phi\preceq 1$ in $K$ we have
$v(P_1^{\phi}) \le rv\phi < v(P_i^\phi)$ for $i\ne 1$. 
\end{proof} 

\noindent
If $\Psi^{>0}\ne \emptyset$, then we let $1$ denote the unique element of 
$\Gamma^{>}$ with $\psi(1)=1$, and we
identify $\Z$ with a subgroup of $\Gamma$ via $k\mapsto k\cdot 1$.

\begin{lemma}\label{diffnewt} Suppose $K$ is newtonian with $\Psi^{>0}\ne \emptyset$ and $\Delta$ is a convex subgroup of $\Gamma$ with $1\in \Delta$.
Then $(K, v_{\Delta})$ is $\d$-henselian.
\end{lemma}
\begin{proof} Let $\dot{\mathcal{O}}$ be the valuation ring of $\dot{v}=v_{\Delta}$. 
Let $P\in \dot{\mathcal{O}}\{Y\}$ of order $\le r$ be such that
$P_1 \dotasymp 1$ and $P_i \dotprec 1$ for $i\ge 2$; by 
Lemma~\ref{Hdhencor} 
it is enough
to show that $P$ has a zero in $\dot{\mathcal{O}}$. We can arrange $P_1\asymp 1$.
Take $\gamma< 0$ in $\Delta$ such that 
$(\gamma/2) + (r+1) < v(P_0)$.
Take $g\in K^\times$ with $vg=\gamma$.
Then $P_{\times g}= P_0 + L + R$ with
$L=P_{1,\times g}$ and $R \dotprec 1$, so $vL=\gamma+ o(\gamma)$.
Take $a$ with $va=-v(L)$. Then $Q:=aP_{\times g}=aP_0+ aL + aR$ with
$$v(aP_0)=-\gamma + o(\gamma) + v(P_0) > r+1 > r\Psi, \quad v(aL)=0, 
\quad aR \dotprec 1,$$ 
so $\ndeg Q=1$ by Lemma~\ref{newtp1}. Since $K$ is 
newtonian, $Q$ has a zero $y\in \mathcal{O}$, and then~$gy$ is a 
zero of $P$ in $\dot{\mathcal{O}}$.
\end{proof}

\noindent
Lemma~\ref{diffnewt} and its proof go through if the assumption that $K$ is
newtonian is replaced by: every $P\in K\{Y\}^{\ne}$ with $\ndeg P=1$ has a zero in the valuation ring~$\dot{\mathcal{O}}$ of~$v_{\Delta}$. 
Let $(K^{\phi}, v_{\phi}^{\flat})$ be the differential field $K^{\phi}$
with valuation $v^{\flat}_{\phi}$. Then:

\begin{lemma}\label{newt2} Let $\phi\preceq \theta\in K$ $($so $\theta$ is active$)$. If
$(K^{\phi}, v_{\phi}^{\flat})$ is $\d$-henselian, then so is~${(K^{\theta}, v_{\theta}^{\flat})}$.
\end{lemma}
\begin{proof} First, $(K^{\phi}, v_{\theta}^{\flat})$ is isomorphic to 
a coarsening of $(K^{\phi}, v_{\phi}^{\flat})$. Now $K^{\theta}= (K^{\phi})^u$
with $u=\theta/\phi\  \asymp_{\theta}^{\flat}\ 1$, by Lemma~\ref{BasicProperties}(i).
It remains to appeal to Lemma~\ref{inva3} and the subsection on compositional conjugation immediately preceding it.
\end{proof}

\noindent
Thus the following two conditions on $K$ are equivalent: 
\begin{enumerate}
\item for every active $\theta$ in $K$ there is a
$\phi\preceq \theta$ such that $(K^{\phi}, v_{\phi}^{\flat})$ is $\d$-henselian;
\item $(K^{\phi}, v_{\phi}^{\flat})$ is $\d$-henselian, for
every~$\phi$.
\end{enumerate}
We already saw that these conditions are satisfied if $K$ is newtonian.
If $K$ is $\upl$-free we can reverse this implication:

\begin{lemma}\label{difhdifnewt} Suppose $K$ is $\upl$-free and
$(K^{\phi}, v_{\phi}^{\flat})$ is $\d$-henselian, for
every~$\phi$. Then~$K$ is newtonian.
\end{lemma}
\begin{proof} Let $P\in K\{Y\}^{\ne}$ and
$\ndeg P=1$; we need to show that then $P$ has a zero in $\mathcal{O}$. 
By Corollary~\ref{betterdifpol}
we can take $\gamma>0$ such that
$P$ has no zero in the region $-\gamma < vy < 0$.
By Lemma~\ref{ndoneupl} we can take 
$\phi$ such that 
$\Gamma_{\phi}^{\flat}\subseteq (-\gamma, \gamma)$, and either
$P^{\phi}\simflat_{\phi} a+bY$ in $K^{\phi}\{Y\}$, $a\preceq b$ in $K$,
or $P^{\phi}\simflat_{\phi} bY'$ in $K^{\phi}\{Y\}$, $b\in K$. Then $P^{\phi}$ has a zero 
$y\preceq_{\phi}^{\flat} 1$ by Lemma~\ref{difhen1}, 
and the choice of $\gamma$ and $\phi$ gives
$y\in \mathcal{O}$. 
\end{proof}

\begin{cor}\label{newtoniancompletion} Suppose $K$ is $\upl$-free and newtonian, and $\operatorname{cf}(\Gamma)=\omega$. Then the completion $K^{\c}$ of
$K$ is also $\upl$-free and newtonian.
\end{cor}
\begin{proof} Let a flattening $(K^{\phi}, v_{\phi}^{\flat})$ of $K$ be given with $v\phi\in \Psi$.
It is $\d$-henselian. By Lem\-ma~\ref{coarseningcompletion} it has as a completion the flattening
$\big((K^{\c})^\phi, v_{\phi}^{\flat}\big)$, so the
latter is $\d$-henselian by Proposition~\ref{dhenscompletion}.
Now $K^{\c}$ is $\upl$-free by Lemma~\ref{lem:lambda-free completion}, so the desired conclusion
follows from Lemma~\ref{difhdifnewt} applied to $K^{\c}$ instead of $K$. 
\end{proof}

\noindent
In Chapter~\ref{ch:QE} we show that the $\upo$-free newtonian Liouville closed $H$-fields are exactly the 
existentially closed $H$-fields; see Appendix~\ref{app:modeltheory} for the concept {\em existentially closed.}\/ This makes
the next result interesting, in particular for
$K=\T$. 

\begin{sloppy}

\begin{cor}\label{uponewlicomp} If $K$ is an $\upo$-free newtonian Liouville closed $H$-field with $\operatorname{cf}(\Gamma)=\omega$, then  $K^{\c}$  is also an $\upo$-free newtonian Liouville closed $H$-field.
\end{cor}

\end{sloppy}

\noindent
This is immediate from Lemmas~\ref{lioucompletion} and~\ref{upocompletion}, and  Corollary~\ref{newtoniancompletion}.

\subsection*{Preparing for newtonization}
We begin with an analogue of Lemma~\ref{dhli3}: 

\begin{lemma}\label{nlj3} Let $r\ge 1$, suppose $K$ is newtonian, and
let $G\in K\{Y\}\setminus K$ have order $\le r$. Then there do not exist 
$y_0,\dots ,y_{r+1}\in K$ such that \begin{enumerate}
\item[\textup{(i)}] $y_{i-1}-y_{i}\succ y_{i}-y_{i+1}$ for all $i\in\{1,\dots, r\}$, and
$y_r\ne y_{r+1}$;
\item[\textup{(ii)}] $G(y_0)=\cdots = G(y_{r+1})=0$;
\item[\textup{(iii)}] $\ndeg G_{+y_{r+1},\times g}=1$ and $y_0- y_{r+1}\preceq g$ for some $g\in K^\times$. 
\end{enumerate}
\end{lemma}
\begin{proof} Towards a contradiction, suppose $y_0,\dots, y_{r+1}\in K$
satisfy (i), (ii), (iii). By taking $\phi$ with sufficiently high $v\phi$ we arrange that 
$y_{i-1}-y_{i}\succ^{\flat}_\phi y_{i}-y_{i+1}$ for all 
$i\in\{1,\dots, r\}$ and $\ddeg G^{\phi}_{+y_{r+1}, \times g} =1$, where
$g$ witnesses (iii). But $(K^{\phi}, v^{\flat}_{\phi})$ is $\d$-henselian, and so we contradict Lemma~\ref{dhli3}.
\end{proof}

\noindent
This yields also the newtonian version of Proposition~\ref{noextrazeros}:

\begin{lemma}\label{newtonnoextrazeros} Suppose $K$ is newtonian and $G\in K\{Y\}$ satisfies $\ndeg G=1$. Let~$E$ be an immediate extension of $K$.
Then $G$ has the same zeros in $\mathcal{O}$ as in $\mathcal{O}_E$.
\end{lemma}
\begin{proof} Note first that $\ndeg G_{+y}=1$ for all 
$y\in \mathcal{O}_E$. Towards a contradiction, suppose $G(\ell)=0$
with $\ell\in \mathcal{O}_E\setminus\mathcal{O}$. Now 
$\ell\preceq 1$ gives $\ndeg G_{+\ell}=1$, and from $G(\ell)=0$
it follows easily that $\ndeg G_{+\ell, \times g}=1$ for all $g\preceq 1$ in $K^\times$.  

\claim{Let $\gamma\in v(\ell-K)$, $\gamma\ge 0$. Then
$G(y)=0$ for some $y\in \mathcal{O}$ with $v(\ell-y)\ge \gamma$.} 

\noindent
To prove this claim, take $a\in K$ and $g\in K^\times$ such that 
$v(\ell -a)=v(g)=\gamma$. Then by Corollary~\ref{aftranew} and the observation
preceding the claim,
$$\ndeg G_{+a,\times g}\ =\ \ndeg G_{+\ell,\times g}\ =\ 1,$$ so
we get $b\in \mathcal{O}$ such that $G(a+gb)=0$, so $y:=a+gb$ satisfies the claim. Taking~$r\ge 1$ with $G$ of order $\le r$, the claim yields
$y_0,\dots, y_r, y_{r+1}\in \mathcal{O}$ such that 
$$\ell -y_0\succ \ell - y_1 \succ\cdots \succ  \ell-y_{r+1}, \quad 
G(y_0)=G(y_1)=\cdots =G(y_{r+1})=0,$$ 
contradicting Lemma~\ref{nlj3}: take $g=1$ in (iii).    
\end{proof}

\begin{cor}\label{newtonization} If $K$ has a newtonian immediate 
extension, then $K$ has a newtonian immediate extension $L$ such that: \begin{enumerate}
\item[\textup{(i)}] $L$ is $\d$-algebraic over $K$;
\item[\textup{(ii)}] no proper differential subfield of $L$ containing $K$ is newtonian.
\end{enumerate}
\end{cor}
\begin{proof} Suppose $F$ is a newtonian immediate extension of $K$. Let $L$ be the intersection inside
$F$ of the collection of all differential subfields $E$ of $F$ that
contain~$K$ and are newtonian. Applying 
Lemma~\ref{newtonnoextrazeros} to these extensions $E\subseteq F$
shows that~$L$ is newtonian. That (i) holds is because the differential
subfield $E\supseteq K$ of~$F$ consisting of all $y\in F$ that are 
$\d$-algebraic over $F$ is newtonian. It is obvious that~(ii) holds.
\end{proof} 

\noindent
The condition that an extension $L$ of $K$ is newtonian
includes $L$ being ungrounded (which is automatic if $L|K$ is immediate). 
A {\bf newtonization} of $K$ is a newtonian extension of~$K$ that embeds over $K$ into every newtonian extension of $K$.  
At this stage this is just a definition. In Section~\ref{sec:newtonization} we can say more.  

\index{newtonization}
\index{asymptotic field!strongly newtonian}

\subsection*{Strong newtonianity} This notion will be useful
in Section~\ref{sec:newtonization}.  

\begin{lemma}\label{cor:dals}
Let $K$ be newtonian, $(a_\rho)$ a  pc-sequence in $K$, $G\in K\{Y\}^{\neq}$,
$\ndeg_{\mathbf a} G=1$, and $\mathbf a:=c_K(a_\rho)$. Then $G(a)=0$ and $a_\rho\leadsto a$ for some $a\in K$.
\end{lemma}
\begin{proof} We can assume that we have a strictly increasing sequence
$(\gamma_{\rho})$ in $\Gamma$ such that $v(a_\sigma-a_{\rho})=\gamma_{\rho}$ for all $\sigma> \rho$. 
For each $\rho$ take $g_\rho\in K$ with $v(g_\rho)=\gamma_{\rho}$. We can further
assume that $\ndeg G_{+a_\rho,\times g_\rho}=1$ for all $\rho$. Then we get for each $\rho$ an element $z_\rho\in K$ with $G(z_\rho)=0$ and $z_{\rho}-a_\rho\preceq g_\rho$.
Let 
$$B_{\rho}\ :=\ \big\{z\in K:\ v(z-a_\rho)\geq\gamma_\rho\big\}$$ 
be the closed ball in $K$ centered at $a_{\rho}$
with valuation radius $\gamma_{\rho}$, so $z_\rho\in B_\rho$.
If $\rho<\sigma$ and $z\in B_\sigma$, then 
$v(z-a_\sigma)\geq\gamma_\sigma>\gamma_\rho=v(a_\sigma-a_\rho)$ and so
$v(z-a_\rho)=\gamma_\rho$; in particular $B_{\rho}\supseteq B_{\sigma}$ whenever $\rho\le \sigma$.

\claim{There is an index $\rho_0$ such that $z_{\rho_0}\in B_{\rho}$ for all $\rho\ge\rho_0$.}

\noindent
Suppose not. Then we get an infinite sequence $\rho_0<\rho_1<\cdots$ of indices such that $y_m:=z_{\rho_m}\notin B_{\rho_n}$ for all $m<n$.
Taking each $y_m$ as center of $B_{\rho_m}$ we see that
$$v(y_m-y_{m+1})\  <\  v(y_{m+1}-y_{m+2})$$ 
for all $m$.
Take $r\geq 1$ such that $G$ has order $\leq r$.  Then 
conditions~(i) and~(ii) of Lemma~\ref{nlj3} are satisfied.  
Set $g:=g_{\rho_0}$.
Then also 
$y_0-y_{r+1} \preceq g$ and $a_{\rho_0}-y_{r+1}\preceq g$. In view of 
Corollary~\ref{aftranew} the latter gives 
$\ndeg G_{+y_{r+1},\times g}=\ndeg G_{+a_{\rho_0}, \times g}=1$, so
condition~(iii) in Lemma~\ref{nlj3} also holds, which contradicts
that lemma.

\medskip
\noindent
This proves the claim.
Let $\rho_0$ be an index as in the claim; then $v(z_{\rho_0}-a_\rho)=\gamma_\rho$ for all~${\rho\geq\rho_0}$ and thus
$a_{\rho} \leadsto z_{\rho_0}$,
and so  $a:=z_{\rho_0}$ has the required property.
\end{proof}

\noindent
Define
$K$ to be {\bf strongly newtonian\/} if it is newtonian and for every divergent pc-sequence $(a_{\rho})$ in $K$ with minimal $\d$-polynomial~$G(Y)$ over $K$ we have $\ndeg_{\bf a} G=1$, where ${\bf a} = c_K(a_{\rho})$. For this notion 
we have an analogue of Theorem~\ref{thm:fcdifhensalgmax}:

\begin{lemma}\label{lem:strongly newtonian} Suppose $K$ has rational asymptotic integration and is strongly newtonian. Then $K$ is asymptotically $\d$-algebraically maximal.
\end{lemma}
\begin{proof} Towards a contradiction, assume $K$ has a proper immediate
$\d$-algebraic extension.
Then Lemma~\ref{flpc1new} yields a divergent pc-sequence $(a_{\rho})$ in $K$ with a minimal 
$\d$-polynomial
$G(Y)$ over $K$. This contradicts
Lemma~\ref{cor:dals}.
\end{proof} 

\index{asymptotic field!strongly newtonian}
\index{newtonian!strongly}
\index{strongly!newtonian}

\subsection*{Quasilinear asymptotic equations} {\em In this subsection $\E\subseteq K^\times$ is $\preceq$-closed.}\/

\begin{lemma}\label{newt5} If $K$ is newtonian and
$\ndeg_{\E} P =1$, then 
$P$ has a zero in $\E\cup \{0\}$.
\end{lemma}
\begin{proof} Assume $K$ is newtonian and
$\ndeg_{\E} P =1$.
Take $g\in \E$ with
$\ndeg P_{\times g}=1$. Then 
$P_{\times g}$ has a zero in $\mathcal{O}$, and thus~$P$ has a zero in
$g\mathcal{O}\subseteq \E\cup \{0\}$. 
\end{proof}

\noindent
Next we
consider an asymptotic equation
\begin{equation}\tag{E} \label{eq:asympt equ, best approx}
P(Y)\ =\ 0, \qquad Y\in \E
\end{equation}
over $K$. 
We say that \eqref{eq:asympt equ, best approx} is {\bf quasilinear} if $\ndeg_{\E} P = 1$. Recall from Section~\ref{sec:ndeg and nval} the notion of a solution $y$ of
\eqref{eq:asympt equ, best approx} best approximating a 
given element $f$ of a valued differential field extension of $K$.

\index{asymptotic equation!quasilinear}
\index{quasilinear!asymptotic equation}

\begin{lemma}\label{lem:best approx}
Suppose $K$ is newtonian, and \eqref{eq:asympt equ, best approx} is quasilinear and has a solution. Let $f$ be an element of a valued differential field extension of $K$. Then $f$ is best approximated by some solution of~\eqref{eq:asympt equ, best approx}.
\end{lemma}
\begin{proof} We may assume that $f\not\succ\E$. Take $\fm\in\E$ such that $f\preceq\fm$, $\ndeg P_{\times\fm}=\ndeg_{\E} P=1$, and
\eqref{eq:asympt equ, best approx} has a solution $y\preceq \fm$. 
By Lemma~\ref{lem: bestapp10-2} we may replace~$P$ by~$P_{\times\fm}$ and $\E$ by $\mathcal{O}^{\neq}$, and thus assume $\E=\mathcal{O}^{\neq}$. Suppose $f$ is not best approximated by any solution of~\eqref{eq:asympt equ, best approx}. 
Then we get an infinite sequence $y_0,y_1,y_2,\dots$ of solutions of~\eqref{eq:asympt equ, best approx}
with $y_0-f\succ y_1-f \succ y_2-f\succ\cdots$. For each $i$
we have $\operatorname{ndeg} P_{+y_i}=\operatorname{ndeg} P=1$, and this leads to a contradiction with Lemma~\ref{nlj3}. 
\end{proof}

\begin{lemma}\label{cor:best approx is small}
Suppose $K$ is newtonian, \eqref{eq:asympt equ, best approx} is quasilinear, and $f\in \E$ satisfies $\ndeg_{\prec f} P_{+f}\ge 1$. 
%$f\in K$ is an approximate solution of 
%\eqref{eq:asympt equ, best approx}. 
Then \eqref{eq:asympt equ, best approx} has a solution 
$y\sim f$,
and every solution $y$ of~\eqref{eq:asympt equ, best approx} that  best approximates~$f$ satisfies $y\sim f$.
\end{lemma}
\begin{proof} Since \eqref{eq:asympt equ, best approx} is quasilinear
and $f\in \E$,
we have 
$$\ndeg_{\prec f} P_{+f}\ \le\ \ndeg_{\E}P_{+f}\ =\ \ndeg_{\E}P\ =\ 1, $$ 
so $\ndeg_{\prec f}P_{+f}=1$. By Lemma~\ref{newt5} we get
$z\prec f$ in $K$ with $P(f+z)=0$, so
$f+z$ is a solution of \eqref{eq:asympt equ, best approx} with $f+z\sim f$. 
Given any solution $y$ of \eqref{eq:asympt equ, best approx} that best approximates $f$, we get
$y-f\preceq (f+z)-f=z\prec f$, so $y\sim f$. 
\end{proof}

\section{Cases of Low Complexity}\label{cordone}

\noindent
In this section we consider weak forms of newtonianity and differential polynomials of low complexity: degree~$1$ or order at most~$2$.
At the end we apply this to the
differential polynomial function $\omega$ and to the related 
function $\sigma$. Throughout this section, $r$ ranges over $\N$, $a$,~$b$,~$y$ over $K$, and $\gamma$ over $\Gamma$.

Define $K$ to be {\bf $r$-newtonian\/} if
every quasilinear $P\in K\{Y\}$ of order~$\le r$ has a zero in 
$\mathcal{O}$. Define $K$ to be {\bf $(1,1)$-newtonian}
if every quasilinear $P\in K[Y,Y']$ of degree~$\le 1$ in $Y'$
has a zero in $\mathcal{O}$. 
Define $K$ to be {\bf $r$-linearly newtonian} if every
quasilinear $P\in K\{Y\}$ 
with $\deg P = 1$ and $\order(P)\leq r$ has a zero in $\mathcal O$.
Define~$K$ to be {\bf linearly newtonian} if $K$ is $r$-linearly newtonian for every~$r$. Each of these conditions on $K$ is clearly invariant under compositional conjugation by any $f\in K^\times$. These notions are mainly used
for $r=1$ and $r=2$. Thus for $K$ we~have: 
$${\xymatrixcolsep{1em}\xymatrix{
\text{newtonian} \ar@{=>}[r] \ar@{=>}[d]
& \hskip0.5em\text{$2$-newtonian} \ar@{=>}[r] \ar@{=>}[d] & \text{$1$-newtonian} \ar@{=>}[r] & \text{$(1,1)$-newtonian}  \ar@{=>}[r] \ar@{=>}[d] & \text{henselian} \\
\ \parbox{5em}{linearly newtonian} \ar@{=>}[r] & \ \parbox{5em}{$2$-linearly newtonian} \ar@{=>}[rr] & & 
\ \parbox{5em}{$1$-linearly newtonian}
}}$$

\index{asymptotic field!$1$-newtonian}
\index{asymptotic field!$2$-newtonian}
\index{asymptotic field!$(1,1)$-newtonian}
\index{asymptotic field!$r$-linearly newtonian}
\index{asymptotic field!linearly newtonian}
\index{newtonian!linearly}
\index{newtonian!$1$-newtonian}
\index{newtonian!$2$-newtonian}
\index{newtonian!$r$-linearly newtonian}
\index{linearly!newtonian}
\index{newtonian!$(1,1)$-newtonian}

\subsection*{Linear newtonianity} Recall from Section~\ref{sec:maxdh} the notion of a valuation ring of a valued differential field with small derivation being linearly surjective. The following result relates this notion to linear newtonianity:

\begin{lemma}\label{variantdiffnewt} Suppose $K$ is $r$-linearly newtonian with $\Psi^{>0}\ne \emptyset$ and $\Delta$ is a convex subgroup of $\Gamma$ with $1\in \Delta$. Then the valuation ring of $(K, v_{\Delta})$ is $r$-linearly surjective.
\end{lemma}
\begin{proof} Like that of Lemma~\ref{diffnewt}, but simpler since here $R=0$.
\end{proof}

\noindent
The conclusion of Lemma~\ref{variantdiffnewt} implies in particular the $r$-linear surjectivity of $K$.

\begin{cor}\label{11linsur}
%\label{cor:linnewt => linsurj}
If $K$ is $r$-linearly newtonian, then $K$ is $r$-linearly surjective.
\end{cor}
\begin{proof} Apply the above to any compositional conjugate $K^\phi$ of $K$.
\end{proof}

\noindent
In combination with Proposition~\ref{vsur3} this yields:

\begin{cor}\label{aslinnewtuplfree} If $K$ has asymptotic integration and is $1$-linearly newtonian, then~$K$ is $\upl$-free.
\end{cor}

\noindent
Asymptotic integrability plus $1$-linear newtonianity
has further nice consequences:

\begin{lemma}\label{lem:linnewt} If $K$ has asymptotic integration, the following are equivalent:
\begin{enumerate}
\item[\textup{(i)}] $K$ is $1$-linearly newtonian;
\item[\textup{(ii)}] every $P\in K\{Y\}$ with $\nval P=\deg P=1$ and $\order P\leq 1$ has
a zero in $\smallo$.
\end{enumerate}
\end{lemma}
\begin{proof} Assume $K$ has asymptotic integration. Then $\Gamma^{>}$ has no least element. 

To prove (i)~$\Rightarrow$~(ii), let $K$ be $1$-linearly newtonian and let $P\in K\{Y\}$ be such that $\nval P=\deg P =1$, $\order P \le 1$. By the remark
following the proof of Proposition~\ref{propnewtordone} 
we have
$g_0\prec 1$ in $K$ such that $\nval P_{\times g}=\ndeg P_{\times g}=1$
for all~$g\in K$ with $g_0 \prec g \prec 1$. For such $g$ we have a zero
of $P$ in $g\mathcal{O}\subseteq \smallo$.

Next, assume (ii), and let $P\in K\{Y\}$,  $\ndeg P=\deg P =1$, $\order P \le 1$; our job is to show that $P$ has a zero in $\mathcal{O}$.   
By the remark
following the proof of Proposition~\ref{propnewtordone} we have
$f_0\succ 1$ in $K$ such that $\nval P_{\times f}=\ndeg P_{\times f}=1$
for all~$f\in K$ with $f_0 \succ f \succ 1$. By Corollary~\ref{betterdifpol} 
we can take $f_0$ such that also $P(y)\ne 0$ for all $y$ with 
$f_0 \succ y \succ 1$. Take any 
$f\in K$ with $f_0\succ f \succ 1$, and take a zero $y\in \smallo$ of~$P_{\times f}$. Then $P(fy)=0$, and so $fy\in \mathcal{O}$, since otherwise
$f_0 \succ fy \succ 1$.  
\end{proof}

\noindent
In Section~\ref{sec:special sets} we defined the set $\I(E)\subseteq E$ for pre-$H$-fields $E$. We now do this for any asymptotic field $E$ by 
$$ \I(E)\ :=\ 
\{f\in E:\ \text{$f\preceq g'$ for some $g\in \mathcal O_E$}\}.$$
Then $\I(E)$ is an $\mathcal{O}_E$-submodule of $E$ with 
$\der \mathcal O_E\subseteq \I(E)$ and $(\mathcal O_E^\times)^\dagger\subseteq\I(E)$. 

\begin{lemma}\label{nli0} Assume $K$ has asymptotic integration
and is $1$-linearly newtonian. Then $K$ is $\d$-valued and 
$\der\smallo=\operatorname{I}(K)=(1+\smallo)^{\dagger}$.
\end{lemma}
\begin{proof} Let $a\in \operatorname{I}(K)$, 
and $P:= Y'-a$. Then
$P^{\phi}=\phi Y' -a$, so $\nval P=1$, and thus
$P$ has a zero $y\in \smallo$ by Lemma~\ref{lem:linnewt}.
This gives $\der\smallo=\operatorname{I}(K)$, and so $K$ is $\d$-valued.
Next, take $Q:= Y'-(a+aY)$. Then $Q^{\phi}=\phi Y'-(a+aY)$, so
 $\nval Q=1$, and thus $Q$ has a zero $y\in \smallo$ by 
 Lemma~\ref{lem:linnewt}, 
and then
 $a=y'/(1+y)=(1+y)^{\dagger}$. This proves $\operatorname{I}(K)=(1+\smallo)^{\dagger}$.
\end{proof}

\noindent
Here is a generalization of Lemma~\ref{lem:linnewt}:

\begin{lemma}\label{lem:rlinnewt} Suppose $K$ has asymptotic integration and $r\ge 1$.  Then the two conditions below are equivalent:
\begin{enumerate}
\item[\textup{(i)}] $K$ is $r$-linearly newtonian;
\item[\textup{(ii)}] every $P\in K\{Y\}$ with $\nval P=\deg P=1$ and $\order P\leq r$ has
a zero in $\smallo$.
\end{enumerate}
\end{lemma}
\begin{proof} By Corollary~\ref{aslinnewtuplfree}, (i) implies 
$K$ is $\upl$-free, and using also
Lemma~\ref{lem:linnewt}, so does (ii).  Now follow the proof of 
Lemma~\ref{lem:linnewt}, using
Corollary~\ref{lem:npol linear} instead of the remark following the proof of Proposition~\ref{propnewtordone}.   
\end{proof}

\subsection*{Application to linear differential equations}  We
determine here the dimension of the kernel of a linear differential operator over $K$ for linearly newtonian $\d$-valued~$K$. The results in this 
subsection will not be used in this volume.
 
\medskip\noindent
Let $A=a_0+a_1\der+\cdots+a_r\der^r\in K[\der]^{\neq}$ where  $a_0,\dots,a_r\in K$, $a_r\neq 0$. In Section~\ref{evtbeh} we defined the set
$$\exc^{\operatorname{e}}(A)\ =\ \big\{ \gamma :\  \nwt_A(\gamma) \geq 1\big\}\  =\ \bigcap_\phi \exc(A^\phi)$$
of eventual exceptional values of $A$. 
 If $K$ is $\upl$-free, then by Lemma~\ref{ndoneupl} we have 
$\nwt_A(\gamma)=1$ for all $\gamma\in\exc^{\operatorname{e}}(A)$.
In  Section~\ref{evtbeh} we also defined the map $v_A^{\operatorname{e}}\colon\Gamma\to\Gamma$, and mentioned  the fact that
\begin{equation}\label{eq:vAphi}
v_{A^\phi}(\gamma)\ =\  v_A^{\operatorname{e}}(\gamma) + \nwt_A(\gamma)v\phi,\quad\text{eventually.}
\end{equation}
Thus if $\gamma\notin  \exc^{\operatorname{e}}(A)$, then
eventually $v_{A^\phi}(\gamma) = v_A^{\operatorname{e}}(\gamma)$.

\begin{lemma}\label{lem:vAe onto}
The map $\gamma\mapsto v_A^{\operatorname{e}}(\gamma)\colon \Gamma\setminus\exc^{\operatorname{e}}(A)\to\Gamma$  is strictly increasing, and if~$K$ is $\upo$-free, then this map is surjective.
\end{lemma}
\begin{proof}
The first part holds because $v_{A^\phi}$ is strictly increasing for each $\phi$. Suppose~$K$ is $\upo$-free and
let $\alpha\in\Gamma$. Corollary~\ref{hevqone} yields 
$\gamma$ such that $v_{A^\phi}(\gamma)=\alpha$ eventually. 
Then $\nwt_A(\gamma)=0$ by \eqref{eq:vAphi}, so $\gamma\notin\exc^{\operatorname{e}}(A)$ and $v_A^{\operatorname{e}}(\gamma)=\alpha$.
\end{proof}

\noindent
We have $v(\ker^{\neq} A)\subseteq \exc^{\operatorname{e}}(A)$; moreover:

\begin{prop}\label{prop:dim ker A}
Suppose $K$ is $r$-linearly newtonian and $K$ has asymptotic integration. Then  $v(\ker^{\neq} A) = \exc^{\operatorname{e}}(A)$.
\end{prop}

\begin{proof} The case $r=0$ is trivial, so assume $r\ge 1$.
Let $\gamma\in\exc^{\operatorname{e}}(A)$, and take $g\in K^\times$ with $vg=\gamma$; our job is to find
$y\in\ker A$ with $y\asymp g$. Replacing $A$ by $Ag$ we arrange $\gamma=0$ and $g=1$. Since $K$ is $\upl$-free, we have $\nwt(A)=1$. Put $P:=a_0Y+a_1Y'+\cdots+a_rY^{(r)}\in K\{Y\}$. Then 
$D_{P^\phi}\in\k^\times\cdot Y'$, so
by Lemma~\ref{dn1,new}(i),
$$D_{P^\phi_{+1}} \in \k^\times\cdot (D_{P^\phi})_{+1}\ =\ \k^\times \cdot Y',$$
and hence $\nval P_{+1}=1$. By Lemma~\ref{lem:rlinnewt} we get $z\in\smallo$ with $P(1+z)=0$,
so we can take $y:=1+z$.
\end{proof}

\noindent
Corollary~\ref{11linsur}, Proposition~\ref{prop:dim ker A}, and Lemma~\ref{valuation-of-basis} yield:

\begin{cor}\label{cor:dim ker A}
If $K$ is $\d$-valued and $r$-linearly newtonian, then
$$ \dim_C \ker A\ =\ \abs{\exc^{\operatorname{e}}(A)}.$$ 
\end{cor}

\noindent
Suppose $K$ is $\d$-valued. If $K$ has an immediate $r$-linearly newtonian 
extension, then $\abs{\exc^{\operatorname{e}}(A)}\leq r$,
by Corollary~\ref{cor:dim ker A}. Thus 
$\abs{\exc^{\operatorname{e}}(A)}\leq r$ if $K$ is $\upo$-free, 
by the remark following Theorem~\ref{maxnewt} 
(to be proved in Section~\ref{sqe}). 
Also, if $K$ is $r$-linearly newtonian and $L$ is an immediate extension
of $K$, then $\ker_L A=\ker A$.

\begin{prop}
Suppose $K$ is $\upo$-free and $r$-linearly newtonian. Then for every $a\ne 0$ there is $y\neq 0$ such that 
$A(y)=a$, $vy\notin\exc^{\operatorname{e}}(A)$, and $v_A^{\operatorname{e}}(vy)=va$.
\end{prop}
\begin{proof}
It is enough to do the case $a=1$.
Lemma~\ref{lem:vAe onto} gives $g\in K^\times$ with 
$vg\notin \exc^{\operatorname{e}}(A)$ and
$v_A^{\operatorname{e}}(vg)=0$. Replacing $A$ by $Ag$ we arrange $0\notin\exc^{\operatorname{e}}(A)$ and ${v^{\operatorname{e}}(A)=0}$; our job is now
to find $y\asymp 1$ with $A(y)=1$. Eventually, $\dwt(A^\phi)=0$ and $v(A^\phi)=0$. Hence eventually $A^\phi(1)\asymp 1$
and thus $a_0=A(1)=A^\phi(1)\asymp 1$.  Put $P:=a_0Y+a_1Y'+\cdots+a_rY^{(r)}\in K\{Y\}$. 
Eventually $P^\phi\sim a_0 Y$ in~$K^{\phi}\{Y\}$, so eventually
$P^{\phi}_{+(1/a_0)}\sim a_0Y+1$. Then for
$Q:=(-1+P)_{+(1/a_0)}$ we have $Q^{\phi}=-1+P^{\phi}_{+(1/a_0)}\sim a_0 Y$, eventually, so $\nval Q=1$.  Lem\-ma~\ref{lem:rlinnewt} gives $z\in\smallo$ with ${(-1+P)\big((1/a_0)+z\big)}=Q(z)=0$,
and thus $y:=(1/a_0)+z$ works.
\end{proof}

\subsection*{Newtonianity of order $r$} Lemmas~\ref{lem:linnewt}
and~\ref{lem:rlinnewt} extend as follows:

\begin{lemma}\label{ndegnvalnewt} Suppose $K$ has asymptotic integration and $r\ge 1$. Then the two conditions on $K$,~$r$ below are equivalent: 
\begin{enumerate}
\item[\textup{(i)}] $K$ is $r$-newtonian;
\item[\textup{(ii)}] every $P\in K\{Y\}$ with $\nval P =1$ and $\order P \le r$ has a zero in $\smallo$.
\end{enumerate}
\end{lemma}
\begin{proof} Follow the proof of Lemma~\ref{lem:rlinnewt},
but for the direction (i)~$\Rightarrow$~(ii), appeal to Lemma~\ref{uplnvallemma} instead of Corollary~\ref{lem:npol linear}.
\end{proof} 

\begin{cor}\label{cor:ndegnvalnewt}
Let $K$ have asymptotic integration and be $r$-newtonian, $r\ge 1$, and let $P\in K\{Y\}$,  
$u\in \mathcal{O}$, and $A\in \k[Y]$ be such that $\order P \le r$, $A(\bar{u})=0$, $A'(\bar{u})\ne 0$, and $D_{P^{\phi}}\in \k^{\times}\cdot A$, eventually. 
Then $P$ has a zero in $u+\smallo$.
\end{cor}
\begin{proof} By Lemma~\ref{dn1,new}(i) we have
$D_{P^{\phi}_{+u}} \in \k^{\times} \cdot \big(D_{P^{\phi}}\big)_{+\bar{u}}$, hence $\nval P_{+u}=1$. It remains to apply Lemma~\ref{ndegnvalnewt}. 
\end{proof}

\begin{cor}\label{cor:y'=Q(y)}
Suppose $K$ has small derivation and asymptotic integration and is $1$-newtonian. Let $Q\in K[Y, Y']$ satisfy $vQ>\Psi$. Then there is a unique $y\in\smallo$ such that $y'=Q(y)$.
\end{cor}
\begin{proof}
With $P=Y'-Q$ we have $P^\phi\sim \phi Y'$ for $\phi\preceq 1$, so $\nval P=1$. This gives $y\in \smallo$ with $y'=Q(y)$
by Lemma~\ref{ndegnvalnewt}. Suppose $z\in \smallo$, $y\ne z$, and $z'=Q(z)$. With $\varepsilon:= z-y$ we get by Taylor expansion
\begin{align*} \varepsilon'\ &=\ Q(z)-Q(y)\ =\ \sum_{i+j\ge 1} Q_{(i,j)}(y)\varepsilon^i(\varepsilon')^j,\ \text{ so}\\
\varepsilon^\dagger\ &=\ \sum_{(i,j)}Q_{(i+1,j)}(y)\varepsilon^{i}(\varepsilon')^j +\sum_{j\ge 1}Q_{(0,j)}(y)(\varepsilon')^{j-1}\varepsilon^\dagger.
\end{align*}
The valuation of the first sum on the last line is $>\Psi$, and the valuation of the second term is 
$>v(\varepsilon^\dagger)\in \Psi$, and we have a contradiction. 
\end{proof}

\noindent
For our application to $\omega(K)$ it is enough to consider
$(1,1)$-newtonianity, where the proof of Lemma~\ref{lem:linnewt} gives the following:

\begin{lemma}\label{11newtval} If $K$ has asymptotic integration, the following are equivalent: \begin{enumerate}
\item[\textup{(i)}] $K$ is $(1,1)$-newtonian;
\item[\textup{(ii)}] every $P\in K[Y,Y']$ with $\deg_{Y'}P\le 1$ and $\nval P =1$ has a zero in $\smallo$.
\end{enumerate}
\end{lemma}

\begin{cor}\label{11newtsimple} Assume $K$ has asymptotic integration and is $(1,1)$-newtonian. Let $P\in K[Y,Y']$ have degree $\le 1$ in $Y'$,
let $u\in \mathcal{O}$, and let $A\in \k[Y]$ be such that $A(\bar{u})=0$, $A'(\bar{u})\ne 0$, and $D_{P^{\phi}}\in \k^{\times}\cdot A$, eventually. 
Then $P$ has a zero in $u+\smallo$.
\end{cor}
\begin{proof} Like that of Corollary~\ref{cor:ndegnvalnewt}; use
Lemma~\ref{11newtval} instead of~\ref{ndegnvalnewt}.
\end{proof}

%\begin{cor}\label{cor:y'=Q(y)}
%Suppose $K$ has asymptotic integration and is 
%$(1,1)$-newtonian. Let $Q\in K[Y]$ with $vQ>\Psi$. Then 
%there %is a unique $y\in\smallo$ such that $y'=Q(y)$.
%\end{cor}
%\begin{proof}
%With $P=Y'-Q$ we have $P^\phi\sim \phi Y'$, so $\nval P=1$.
%Existence hence follows from Lemma~\ref{11newtval}. For %uniqueness, suppose $y$,~$z$ are distinct elements of 
%$\smallo$ such that $y'=Q(y)$ and $z'=Q(z)$. Then 
%$(y-z)'=Q(y)-Q(z)$, hence by Taylor expansion, for 
%some $n$ we have
%$$(y-z)^\dagger = \sum_{i=1}^n 
%\frac{Q^{(i)}(z)}{i!}(y-z)^{i-1}$$
%where the valuation of the left-hand side is in $\Psi$ and %that of the right-hand side is~$>\Psi$, a contradiction. 
%\end{proof}

%\noindent
%Let $E$ be a valued differential
%field with small derivation. Then $E$ is said to be 
%{\bf $(1,1)$-$\d$-henselian} if every 
%$P\in E[Y,Y']^{\ne}$ with $\deg_{Y'}P\le 1$ and 
%$\operatorname{ddeg} P = 1$ has a zero
%in $\mathcal{O}_E$. Now Lemma~\ref{diffnewt} goes 
%through with basically 
%the same proof if we replace ``newtonian'' and ``$\d$-%henselian''
%by ``$(1,1)$-newtonian'' and 
%``$(1,1)$-$\d$-henselian,'' respectively. 
%Corollary~\ref{newtlinsur} goes through likewise:  

%\index{$(1,1)$-$\d$-henselian}
%\index{asymptotic field!$(1,1)$-$\d$-henselian}

\subsection*{Application to $\omega$}
Let $f\in K$. When is $f\in \omega(K)$, that is, when is there $y$ such that $\omega(y)=f$, equivalently, when does $P(Y):= Y^2+2Y'+f$ have a zero in
$K$? Before we give a partial answer, let $b\ne 0$ and note that
\begin{align*} P_{+a}\ &=\  Y^2 + 2aY + 2Y' + P(a), \text{ so }\\
P_{+a,\times b}\ &=\ b^2Y^2 + (2ab + 2b')Y + 2bY' + P(a), \text{ and thus}\\
b^{-1}P_{+a, \times b}^\phi\ &=\ bY^2 + (2a+ 2b^\dagger)Y + 2\phi Y' + P(a)/b.
\end{align*}
This leads to the following result:

\begin{lemma} \label{lem:omega(K) downward closed}
Suppose $K$ is a $(1,1)$-newtonian real closed $H$-field. Then the
subset $\omega\big(\Upl(K)^\downarrow\big)$ of $K$ is downward closed.
\end{lemma}
\begin{proof} Corollary~\ref{11linsur} for $r=1$ 
shows $K$ to be $1$-linearly surjective. As $K$ is also $\d$-valued, it has asymptotic integration. Thus $K$ is $\upl$-free by Corollary~\ref{aslinnewtuplfree}.

Let $f\in \omega\big(\Upl(K)\big)^\downarrow$; we need to show that
$f\in \omega\big(\Upl(K)^\downarrow\big)$.
By Corollary~\ref{cor:omega for H-fields, 1},
we may take $\rho$ so that $f<\upo_\rho$. Then 
$f-\upo_{\rho}\preceq f-\upo_{\rho'}$ for $\rho < \rho'$, so by increasing~$\rho$ and using Lemma~\ref{pcrho} 
we arrange that $f-\upo_\rho\succ\upg_\rho^2$. 
Set $a:=\upl_{\rho}=-\upg^{\dagger}_{\rho}$. 
Take $b\in K^>$ such that $f-\upo_{\rho}=-b^2$. Then with $P$ as above,
$$ b^{-1}P^{\phi}_{+a, \times b}\ =\ bY^2 + 2(b/\upg_{\rho})^\dagger Y + 2\phi Y' -b.$$
We have $b\succ\upg_\rho$, hence $b$ is active and $b\succ (b/\upg_{\rho})^\dagger$,
and eventually $(b/\upg_{\rho})^\dagger \succ \phi$.
It follows that the
hypothesis of Corollary~\ref{11newtsimple} holds with $P_{+a,\times b}$ in
the role of $P$ and $u=-1$, $A=Y^2-1$. Then that corollary provides a zero
$g\in -1+\smallo$ of $P_{+a,\times b}$, and so $\omega(a+bg)=f$.
It remains to note that $a+bg< a\in \Upl(K)$. 
\end{proof}

\begin{cor}\label{11newtliomegadown} If $K$ is a $(1,1)$-newtonian Liouville closed $H$-field, 
then the subset $\omega(K)=\Upo(K)$ of $K$ is downward closed.
\end{cor}

\subsection*{Application to the function $\sigma$}
In this subsection we assume that $f\in K$. Solving the equation $\sigma(y)=f$ in $K^\times$ means finding a zero of $$S(Y)\ :=\ \sigma(Y)-f\ =\ \omega(-Y^\dagger)+Y^2-f \in K\<Y\>$$ 
in $K^{\times}$. 
Let $b\in K^\times$. Extending the automorphism $P\mapsto P_{\times b}$
of the ring $K\{Y\}$ to an automorphism $R \mapsto R_{\times b}$ of its fraction field $K\<Y\>$, we get
$$S_{\times b}\ =\ \omega( -b^\dagger -Y^\dagger ) + b^2Y^2 - f.$$
Lemma~\ref{varrholemma, claim, 1} applied to $w=-b^\dagger -Y^\dagger$, $z=-b^\dagger$,
gives
\begin{align*}
\omega( -b^\dagger -Y^\dagger )\ &=\ \omega(-b^\dagger)+Y^\dagger\cdot\big( 2(Y^{\dagger\dagger}-b^\dagger-Y^\dagger) + Y^\dagger\big) \\ 
&=\ \omega(-b^\dagger)+Y^\dagger\cdot\big( 2(Y')^\dagger-2b^\dagger-3Y^\dagger\big)
\end{align*}
and so, with $Q:=Y^2S_{\times b}\in K\{Y\}$, we have
\begin{align*}
Q\ &=\ b^2Y^4 + \big(\omega(-b^\dagger)-f\big) Y^2 + Y'Y\cdot\big( 2(Y')^\dagger-2b^\dagger-3Y^\dagger\big) \\
&=\ b^2Y^4 + \big(\omega(-b^\dagger)-f\big) Y^2 + R,\ \text{ where} \\
R\ :&=\ 2Y''Y-2b^\dagger Y'Y-3(Y')^2\in K\{Y\}.
\end{align*}
Now $(Y'')^\phi = \phi^2 Y'' + \phi' Y'$ and thus
\begin{align*}
R^\phi\ &=\ 2(\phi^2 Y'' + \phi'Y')Y - 2b^\dagger\phi Y'Y - 3\phi^2 (Y')^2 \\
&=\ \phi\cdot \big(
2\phi Y''Y +  2(\phi/b)^\dagger Y'Y - 3\phi (Y')^2 \big).
\end{align*}
We use this computation in the proof of the next result:

\begin{prop}\label{prop:sigma(K) upward closed} Let $K$ be a 
$2$-newtonian
real closed $H$-field with asymptotic integration. Then the subset~$\sigma\big(\Upg(K)^\uparrow\big)$ of $K$ is upward closed.
\end{prop}
\begin{proof} 
%Since $K$ is $\upl$-free, $K$ has asymptotic integration.
Assume $f> \sigma(\upg_\rho)$; by Corollary~\ref{cor:omega-freeness and cut}
it suffices to show that then
$f\in \sigma\big(\Upg(K)^\uparrow\big)$. Since $\sigma(\upg_{\rho})>\sigma(\upg_{\rho +1})$ and
$\sigma(\upg_{\rho})-\sigma(\upg_{\rho +1})\sim \upg_{\rho}^2$, we can increase $\rho$ and arrange that $f-\sigma(\upg_\rho) \succ \upg_\rho^2$.
Take $b\in K^>$ with $b^2=f-\sigma(\upg_\rho)$; then $b\succ\upg_\rho$, so~$b$ is active.
Moreover, using Lemma~\ref{varrholemma, claim, 2} in the last step in
the next line,
$$\sigma(b)-f\ =\ \sigma(b)-\sigma(\upg_\rho)-b^2\ =\ \omega(-b^\dagger) - \omega(-\upg_\rho^\dagger) - \upg_\rho^2\ \prec\ b^2,$$
and so $\omega(-b^\dagger) -f \sim -b^2$. 
Eventually $\phi\prec b$, so 
$(\phi/b)^\dagger \prec  b$,
hence
$R^\phi \prec b^2$, and thus
$Q^\phi  \sim b^2Y^2(Y^2-1)$. 
Corollary~\ref{cor:ndegnvalnewt} gives $u\in 1+\smallo$ with $Q(u)=0$.
Then $\sigma(bu)=f$ and $bu\in\Upg(K)^\uparrow$.
\end{proof}

\begin{cor}\label{cor: sigmaupwardclosed} If $K$ is a $2$-newtonian Liouville closed $H$-field, then the subset $\sigma\big(\Upg(K)\big)$ 
of $K$ is upward closed. 
\end{cor}

\noindent
In combination with Corollaries~\ref{eqschwarz} and 
~\ref{11newtliomegadown} this yields:

\begin{cor}\label{2newtschwarz} If $K$ is a $2$-newtonian $\upo$-free Liouville closed $H$-field, then $K$ is Schwarz closed.
\end{cor}

\subsection*{Notes and comments} All newtonian $K$ known to us are $\upo$-free. Note that if $K$ has asymptotic integration and is newtonian, then $K$ is $\upl$-free by 
Corollary~\ref{aslinnewtuplfree}.  This leaves open whether there are newtonian $K$ without asymptotic integration, that is, with a gap, and whether there are $\upl$-free newtonian $K$ that are not $\upo$-free. 
Corollary~\ref{cor:y'=Q(y)} for $K=\T$ and $Q\in K[Y]$ is in \cite[Corollary~63]{FVK-hensel}.

The intersection $E$ of all maximal Hardy fields is clearly a Hardy field that contains~$\R$ and is Liouville closed. 
By Boshernitzan~\cite{Boshernitzan},
$E$ is $\d$-algebraic over $\R$, and by~\cite[Proposition~3.7]{Boshernitzan87} there is no $y\in E$ with $y''+y=\ex^{x^2}$. Hence $E$ is not
$2$-linearly surjective and thus (by Corollary~\ref{11linsur}) $E$ is not newtonian.

\section{Solving Quasilinear Equations} \label{sqe}

\noindent
{\em In this section $K$ is $\upo$-free, $a$,~$b$,~$y$ range over $K$, and~$P$ over $K\{Y\}^{\ne}$.} 

%\begin{lemma}\label{ndegnvalnewt} The following conditions 
%on $K$ are equivalent:
%\begin{enumerate}
%\item $K$ is newtonian;
%\item every $P$ with $\nval P=1$ has a zero in $\smallo$.
%\end{enumerate}
%\end{lemma}
%\begin{proof} Assume (1), and suppose $\nval P=1$. By 
%Corollary~\ref{decupondegnval}
%we have
%$g_0\prec 1$ in $K$ such that $\nval P_{\times g}=\ndeg %P_{\times g}=1$
%for all $g\in K$ with $g_0 \prec g \prec 1$. For such 
%$g$ we have a zero
%of $P$ in $g\mathcal{O}\subseteq \smallo$.
%Next, assume (2), and suppose $\ndeg P=1$; our job is 
%to show that $P$ has a zero in $\mathcal{O}$.   
%By Corollary~\ref{decupondegnval} we have
%$f_0\succ 1$ in $K$ such that $\nval P_{\times f}=\ndeg %P_{\times f}=1$
%for all $f\in K$ with $f_0 \succ f \succ 1$. By 
%Corollary~\ref{betterdifpol} 
%we can take $f_0$ such that also $P(y)\ne 0$ for all 
%$y$ with $f_0 \succ y \succ 1$. Take any 
%$f\in K$ with $f_0\succ f \succ 1$, and take a 
%zero $y\in \mathcal{O}$ of 
%$P_{\times f}$. Then $P(fy)=0$, and so 
%$fy\in \mathcal{O}$, since otherwise
%$f_0 \succ fy \succ 1$.  
%\end{proof}
 
%\begin{cor}\label{cor:ndegnvalnewt}
%Let $K$ be newtonian, and let  
%$u\in \mathcal{O}$ and $A\in \k[Y]$ be such that  
%$A(\bar{u})=0$, $A'(\bar{u})\ne 0$, and 
%$D_{P^{\phi}}\in \k^{\times}\cdot A$, eventually. 
%Then $P$ has a zero in $u+\smallo$.
%\end{cor}
%\begin{proof} By Lemma~\ref{dn1,new}(1) we have
%$D_{P^{\phi}_{+u}} \in 
%\k^{\times} \cdot \big(D_{P^{\phi}}\big)_{+\bar{u}}$.
%It follows that $\nval P_{+u}=1$, so we can 
%apply Lemma~\ref{ndegnvalnewt}. 
%\end{proof}

\subsection*{Newton position, and proof of Theorem~\ref{maxnewt}}
Suppose $\nval P =1$ with $P_0\ne 0$. By Corollary~\ref{hevqone} (which assumes $\upo$-freeness) 
we can take $g\in K^\times$
such that eventually $P_0 \asymp P^{\phi}_{1,\times g}$.
Since $P_0\prec  P^{\phi}_1$, eventually, we have 
$g\prec 1$. Let $i\ge 2$. Since
$P_1^{\phi}\succeq  P_i^{\phi}$, eventually, 
we get $P_{1,\times g}^{\phi}\succ P_{i,\times g}^{\phi}$, eventually.
Thus $\ndeg P_{\times g}=1$, and so if $K$ is newtonian, then $P$ has a zero in $g\mathcal{O}$, 
but $P$ has no zero in $g\smallo$.   

\medskip\noindent
Define $P$ to be {\bf in newton position at $a$\/} \index{Newton!position}
\index{differential polynomial!in newton position at $a$} if $\nval P_{+a}=1$. Suppose
$P$ is in newton position at $a$ and set $Q:= P_{+a}$, so $Q(0)=P(a)$. 
If $P(a)\ne 0$, then by the above there is $g\in K^\times$
such that eventually $P(a)=Q(0) \asymp Q^{\phi}_{1,\times g}$, and as
$vg$ does not depend on the choice of such $g$,
we set $v^{\ev}(P,a):= vg$. If $P(a)=0$ we set $v^{\ev}(P,a)=\infty\in \Gamma_{\infty}$.

\begin{lemma}\label{newt3}
Suppose $K$ is newtonian and $P$ is in newton position at $a$. 
Then $P(b)=0$ 
and $v(a-b)\ge v^{\ev}(P,a)$ for some $b$; any such $b$ satisfies
$v(a-b)=v^{\ev}(P,a)$.
\end{lemma} 
\begin{proof} If $P(a)=0$, then we must take $b=a$. Assume $P(a)\ne 0$, set
$\gamma:= v^{\ev}(P,a)\in \Gamma$, take $g\in K$ with 
$vg=\gamma$. With $Q:= P_{+a}$ we have 
$P(a+gY)\ =\ Q_{\times g}$ and $\ndeg Q_{\times g}=1$, so 
we get $y\asymp 1$ with $Q(gy)=0$, and so for $b:=a+gy$ we have $P(b)=0$ and
$v(a-b)=v^{\ev}(P,a)$. Conversely, if $v(a-b)\ge \gamma$ and $P(b)=0$, then $b=a+gy$ with $y\preceq 1$, so $Q(gy)=0$, hence $y\asymp 1$, and thus
$v(a-b)=\gamma$.  
\end{proof}

\noindent
Without assuming $K$ is newtonian, we get:

\begin{lemma}\label{newt4}
Suppose $P$ is in newton position at $a$ and $P(a)\ne 0$. 
Then there exists~$b$ with the following properties: \begin{enumerate}
\item[\textup{(i)}] $P$ is in newton position at $b$, $v(a-b)= v^{\ev}(P,a)$, and $P(b)\prec P(a)$; 
\item[\textup{(ii)}] for all $b^*\in K$ with $v(a-b^*)\ge v^{\ev}(P,a)$: $\ P(b^*)\prec P(a)\Leftrightarrow a-b\sim a-b^*$;
\item[\textup{(iii)}] for all $b^*\in K$, if $a-b\sim a-b^*$, then $P$ is in newton position at $b^*$ and
$v^{\ev}(P,b^*)> v^{\ev}(P,a)$.
\end{enumerate}
\end{lemma}
\begin{proof} With $Q=P_{+a}$, $\gamma=v^{\ev}(P,a)$, and $g\in K$
with $vg=\gamma$, we have 
$$P(a)\ \asymp\ Q^{\phi}_{1,\times g}\ \succ\  Q^{\phi}_{i, \times g}\ \text{ for
$\ i\ge 2$, eventually}.$$ 
Thus $D_{Q^{\phi}_{\times g}}$ is isobaric of weight $0$, eventually, so we get $d\in K$ with $d\asymp P(a)$
such that eventually we have 
$$Q_{\times g}^{\phi}\ =\ P(a)+ dY + R_{\phi}, \quad R_{\phi}\in K^{\phi}\{Y\},\quad R_{\phi} \prec P(a).$$
Taking $y\sim -P(a)/d$ gives $y\asymp 1$ and
$Q(gy)\prec P(a)$, so with $b:= a+gy$ we obtain $P(b)\prec P(a)$
and $v(a-b)=\gamma$. Now $P_{+b}=P_{+a+gy}$ with $gy\prec 1$, so 
$\nval P_{+b}=\nval P_{+a}=1$ by Lemma~\ref{ndnm}, so $P$ is in newton position at $b$.
Conversely, if $b^*\in K$, $v(a-b^*)\ge \gamma$ and $P(b^*) \prec P(a)$, then
$b^*=a+gy^*$ with $y^*\sim -P(a)/d$.  

With $y$ and $b$ as above it remains to show that $v^{\ev}(P,b) > v^{\ev}(P,a)$. 
Let $P$ have order $\le r$, let
$\i$ and $\j$ range over $\mathbb{N}^{1+r}$, and recall that
$P_{(\i)}= \frac{P^{(\i)}}{\i!}\in K\{Y\}$, so 
$$P_{+a}^{\phi}(Y)\ =\ P^{\phi}(a) + \sum_{|\i|\ge 1} (P^{\phi})_{(\i)}(a)Y^{\i}, \qquad 
\big(P_{+a}^{\phi}\big)_{1}\ =\ \sum_{|\i|=1} (P^{\phi})_{(\i)}(a)Y^{(\i)} ,$$
and likewise with $b$ instead of $a$. 
Taylor expanding $(P^{\phi})_{(\i)}$ at $a$ for $|\i|=1$ gives
$$(P^{\phi})_{(\i)}(b)\ =\ (P^{\phi})_{(\i)}(a)+
\sum\limits_{|\j|\geq 1}(P^{\phi})_{(\i)(\j)}(a)\cdot (gy)^{\j} \quad \text{in $K^{\phi}$.}$$
As $(P^{\phi})_{(\i)(\j)}(a)=
{{\i+\j}\choose {\j}}(P^{\phi})_{(\i+\j)}(a)$, we get for  $A_{\phi}:=\big(P_{+a}^{\phi}\big)_1$ and $B_{\phi}:= \big(P_{+b}^{\phi}\big)_{1}$ that
 $B_{\phi}= A_{\phi}  + E_{\phi}$,
where eventually $v(E_{\phi}) \ge v(A_{\phi})+\gamma+ o(\gamma)$. Together with $v_{A_{\phi}}(\gamma)= v(A_{\phi}) + \gamma + o(\gamma)$, we get $$v_{E_{\phi}}(\gamma)\ =\ v(E_{\phi})+ \gamma+ o(\gamma)\ \ge\ v(A_{\phi})+ 
2\gamma + o(\gamma)\ >\ v_{A_{\phi}}(\gamma),\ \text{ eventually}.$$
Hence $v_{B_{\phi}}(\gamma) = v_{A_{\phi}}(\gamma)$, and so 
$P(b) \prec P(a)$ forces
$v^{\ev}(P,b) > v^{\ev}(P,a)$.  
\end{proof}

\begin{lemma}\label{newt.imm} Suppose $P$ is in newton position at $a$
and there is no $b$ with $P(b)=0$ and 
$v(a-b) = v^{\ev}(P,a)$.
Then there exists a divergent pc-sequence $(a_\rho)$ in~$K$ such that
$P(a_\rho) \leadsto 0$.  
\end{lemma}  
\begin{proof} Let $(a_\rho)_{\rho<\lambda}$ be a sequence 
in $K$ with $\lambda$ an ordinal $>0$, $a_0=a$, and
\begin{enumerate}
  \item $P$ is in newton position at $a_\rho$, for all $\rho < \lambda$,
  \item $v(a_{\rho'} -a_\rho)=v^{\ev}(P,a_\rho)$ 
whenever $\rho<\rho'<\lambda$,
  \item $P(a_{\rho'})\prec P(a_\rho)$ and 
$v^{\ev}(P,a_{\rho'})>v^{\ev}(P,a_\rho)$ whenever $\rho<\rho'<\lambda$.
\end{enumerate}
Note that there is such a sequence if $\lambda=1$. 
Suppose $\lambda= \mu +1$ is a successor ordinal. Then Lemma~\ref{newt4} 
yields $a_\lambda\in K$ such that $v(a_\lambda - a_\mu)=v^{\ev}(P,a_\mu)$,
$P(a_\lambda)\prec P(a_\mu)$ and 
$v^{\ev}(P,a_\lambda)>v^{\ev}(P,a_\mu)$. Then the extended
sequence $(a_\rho)_{\rho<\lambda +1}$ has the above properties 
with $\lambda+1$ instead of $\lambda$. 

Suppose $\lambda$ is a limit ordinal. Then $(a_\rho)$ is
a pc-sequence and 
$P(a_\rho) \leadsto 0$. If~$(a_\rho)$ has no pseudolimit in $K$ 
we are done. Assume otherwise, and take
a pseudolimit $a_\lambda\in K$ of $(a_\rho)$. The extended 
sequence $(a_\rho)_{\rho<\lambda +1}$ clearly satisfies condition 
(2) with $\lambda+1$ instead of $\lambda$. Applying Lemma~\ref{newt4} to
$a_\rho$, $a_{\rho+1}$ and $a_\lambda$ in the place of $a$, $b$ and~$b^*$, where
$\rho < \lambda$, we see that conditions (1) and (3) are also satisfied 
with $\lambda+1$ instead of~$\lambda$.
This building process must come to an end.  
\end{proof}

\begin{proof}[Proof of Theorem~\ref{maxnewt}] Assume that $K$ has no 
proper immediate $\d$-algebraic extension. In order for $K$ to be newtonian, it suffices by Lemma~\ref{ndegnvalnewt}
to show that every $P$ with $\nval P =1$ has a zero in $\smallo$.
So let $\nval P=1$ and suppose towards a contradiction that
$P$ has no zero in $\smallo$. Then $P$ is in newton position at $0$,
and so by Lemma~\ref{newt.imm} there exists a divergent pc-sequence $(a_{\rho})$
in $K$ with $P(a_{\rho}) \leadsto 0$. Then~$K$ has a proper immediate
$\d$-algebraic extension by Section~\ref{sec:construct imm exts}, a contradiction.
This concludes the proof of Theorem~\ref{maxnewt}.
\end{proof}

%\begin{cor}\label{cor:dals} 
%Suppose that $K$ is $\d$-algebraically maximal and $\Gamma$ is %divisible. Let
%$(a_\rho)$ be a pc-sequence in~$K$, $G(Y)\in K\{Y\}^{\ne}$ and 
%$\ndeg_{\bf a} G=1$, ${\bf a}:= c_K(a_{\rho})$.
%Then there exists $b\in K$
%such that $a_\rho \leadsto b$ and $G(b)=0$.
%\end{cor}
%\begin{proof} It follows from Lemma~\ref{dpkell} and other results 
%in Section~\ref{sec:construct imm exts} that $(a_{\rho})$ is  
%of $\d$-algebraic type
%over $K$. So $a_{\rho} \leadsto a\in K$. Set  
%$$\E\ =\ \{f\in K^\times:\ \text{$f\prec a-a_\rho$, eventually}\}.$$
%If $\E\ne \emptyset$, then  Corollary~\ref{cor:ndegE} 
%gives $f\in \E$ with $\ndeg G_{+a, \times f}=1$, and so we get
%$u\in \mathcal{O}$ with $G(a+uf)=0$, so that $b:= a+uf$ has the %desired property. If $\E=\emptyset$, then
%$\val G_{+a}=1$, so $G(a)=0$, and $b:=a$ has the
%desired property.
%\end{proof}

\subsection*{Application to solving asymptotic equations}  
{\em In this subsection
$K$ is $\d$-valued with small derivation, and 
$\fM$ is a monomial group for~$K$.}\/

\begin{lemma}\label{lem:newtonian approx zero}
Suppose $K$ is newtonian.
Let $g\in K^\times$ be an approximate zero of $P$ such that $\ndeg P_{\times g}=1$.
Then there exists $y\sim g$ in $K$ such that $P(y)=0$.
\end{lemma}
\begin{proof}
Take $c\in C^\times$ and $\fm\in\fM$ with $g\sim c\fm$. Then
$N_{P_{\times\fm}}(c)=0$, so 
$$\nval P_{\times\fm, +c}\ =\ \val N_{P_{\times\fm,+c}}\ =\ \val (N_{P_{\times\fm}})_{+c}\ \geq\ 1.$$ Now $\nval P_{\times\fm,+c}\le
\ndeg P_{\times\fm,+c}=\ndeg P_{\times\fm}=1$, so
$\nval P_{\times\fm,+c}=1$, giving $z\in\smallo$ with $P_{\times\fm,+c}(z)=0$ by
Lemma~\ref{ndegnvalnewt}.
Thus $y\sim g$ and $P(y)=0$ for $y:=(c+z)\fm$.
\end{proof}

\noindent
Next we consider an asymptotic equation
\begin{equation}\tag{E} \label{eq:asympt equ solve}
P(Y)\ =\ 0, \qquad Y\in \E
\end{equation}
where $P\in K\{Y\}^{\neq}$ and $\E\subseteq K^\times$ is $\preceq$-closed.

\medskip\noindent
{\em In the next theorem and its corollaries~\ref{cor:solve asymptotic equ} and~\ref{cor:solve asymptotic equ, 2}, $C$ is algebraically closed, $\Gamma$ is divisible, and $K$ is asymptotically $\d$-algebraically maximal.}\/

\begin{theorem}\label{thm:solve asymptotic equ}
If $\ndeg_{\E} P>\val P=0$, then~\eqref{eq:asympt equ solve} has a solution.
\end{theorem}
\begin{proof}
We proceed by induction on $d=\ndeg_{\E} P$. If $d=1$ and $\val P=0$, then we take
$g\in \E$ with $\ndeg P_{\times g}=1$, so Theorem~\ref{maxnewt}
gives $y\preceq 1$ with $P(gy)=0$, and then $gy$ is a solution of \eqref{eq:asympt equ solve}.  Suppose $d>1$, $\val P=0$, and~\eqref{eq:asympt equ solve} does not have a solution; we shall derive a contradiction. Proposition~\ref{prop:unravel} gives an unraveler~$(f,\E')$ for \eqref{eq:asympt equ solve}. The corresponding refinement  \eqref{eq:asympt equ 2} of  \eqref{eq:asympt equ solve} has Newton degree $d$ and is unraveled, with $\val P_{+f}=0$. Replacing $P$,~$\E$ by $P_{+f}$,~$\E'$, we may therefore assume that  \eqref{eq:asympt equ solve} is unraveled. 
As $C$ is algebraically closed, we get from Corollary~\ref{cor:sum of alg mult} an algebraic approximate solution $f$ of  \eqref{eq:asympt equ solve}.
Since \eqref{eq:asympt equ solve} is unraveled, we have $\val P_{+f}=0< \ndeg_{\prec f} P_{+f}<d$,
so by the inductive hypothesis, the refinement
$$P_{+f}(Y)\ =\ 0, \qquad Y\prec f$$
of \eqref{eq:asympt equ solve}, and hence \eqref{eq:asympt equ solve} itself, has a solution. This is the desired contradiction.
\end{proof}

\begin{cor}\label{cor:solve asymptotic equ}
$K$ is weakly differentially closed. 
\end{cor}
\begin{proof}
Let $P\in K\{Y\}\setminus K$. If $\val P>0$, then $P(0)=0$. Otherwise we have  $\deg P> \val P=0$. It remains to apply Theorem~\ref{thm:solve asymptotic equ} with $\E=K^\times$.
\end{proof}

\begin{cor}\label{cor:solve asymptotic equ, 2} Suppose $g\in K^\times$ is an approximate zero of $P$. Then there exists $y\sim g$ such that $P(y)=0$.
\end{cor}
\begin{proof} An
equivalence at the beginning of Section~\ref{sec:aseq} gives $\ndeg_{\prec g} P_{+g}\geq 1$, so
$$P_{+g}(Y)\ =\ 0,\qquad Y\prec g$$
has positive Newton degree. If $\val P_{+g} \geq 1$, then $P(g)=0$, so $y:=g$ works. Suppose $\val P_{+g}=0$.
Then Theorem~\ref{thm:solve asymptotic equ} yields $z\prec g$ in $K^\times$ with $P_{+g}(z)=0$, and so
$P(y)=0$ and $y\sim g$ for $y:=g+z$. 
\end{proof}

\noindent
\textit{In the next theorem and its corollary we assume
that $C$ is real closed, $\Gamma$ is divisible, and $K$ is asymptotically $\d$-algebraically maximal.}\/

\begin{theorem}\label{thm:oddsolution}
If $\val P=0$ and $\ndeg_\E P$ is odd, then~\eqref{eq:asympt equ solve} has a solution.
\end{theorem}
\begin{proof}
We argue by induction on $d=\ndeg_{\E} P$ as in the proof of Theorem~\ref{thm:solve asymptotic equ}.
The case $d=1$ is dealt with as in the proof of that theorem.
Suppose $d>1$ is odd, $\val P=0$; assume towards a contradiction that~\eqref{eq:asympt equ solve} does not have a solution. As in the proof of Theorem~\ref{thm:solve asymptotic equ} we arrange that~\eqref{eq:asympt equ solve} is unraveled.

 Suppose $f=c\fm$ ($c\in C^\times$, $\fm\in\E$) is an approximate solution  of~\eqref{eq:asympt equ solve} of odd multiplicity. This multiplicity equals
 $\mu+w$ where $\mu$ is the algebraic multiplicity of~$f$ as an approximate zero of $P$ and $w:=\nwt P_{\times\fm}$. We have $\val P_{+f}=0$. Also $\ndeg_{\prec f} P_{+f}=\val (N_{P_{\times\fm}})_{+c}=\mu+w<d$ by Lemma~\ref{transition} and~\eqref{eq:asympt equ solve} being unraveled. Then by the inductive hypothesis, the refinement
$$P_{+f}(Y)=0,\qquad Y\prec f$$
of~\eqref{eq:asympt equ solve}, and hence~\eqref{eq:asympt equ solve} itself, has a solution, and we have a contradiction. So it is enough to find 
such an $f$. In order to do so,  we distinguish two cases. 

\case[1]{There exists $\fm\in\E$ such that $w:=\nwt P_{\times\fm}$ is odd.}
Take such $\fm$ and take $A\in C[Y]$ with $N_{P_{\times\fm}}=A(Y)\cdot(Y')^w$; pick
$c\in C^\times$ with $A(c)\neq 0$. Then the approximate solution $f=c\fm$  of~\eqref{eq:asympt equ solve} has odd multiplicity $w$.

\case[2]{$\nwt P_{\times\fm}$ is even for all $\fm\in\E$.}
Corollary~\ref{cor:sum of alg mult} gives an approximate zero $f\in \E$
of $P$ with odd algebraic multiplicity, and thus with odd multiplicity. 
\end{proof}

\noindent
As Theorem~\ref{thm:solve asymptotic equ} gives Corollary~\ref{cor:solve asymptotic equ}, so we get from Theorem~\ref{thm:oddsolution}:

\begin{cor}\label{cor:solve asymptotic equ, real closed}
If $\deg P$ is odd, then $P$ has a zero in $K$.
\end{cor}

\noindent
Recall from Section~\ref{sec:The Dominant Part of a Differential Polynomial} that $\mathfrak m\in \mathfrak M$ is a starting monomial for $P$ iff
$\nwt(P_{\times\fm})\geq 1$ or $N_{P_{\times\fm}}$ is not homogeneous; equivalently, $N_{P_{\times\fm}}\notin CY^\N$.

\begin{cor}\label{cor:solve asymptotic equ, 3}
Let $\fm\in\fM$. Then $\fm$ is a starting monomial for $P$ if
and only if $P(f)=0$ and $f\asymp\fm$ for some $f$ in some $\d$-valued extension of $K$.
\end{cor}
\begin{proof}
Take an algebraic closure $K^\alg$ of the 
$\d$-valued field $K$.
Use Lemma~\ref{exmongr} to equip $K^\alg$ with a monomial group $\fM^\alg$
containing $\fM$. By Zorn we have an immediate
asymptotically $\d$-algebraically maximal $\d$-algebraic extension $L$
of~$K^\alg$. Then $L$ is $\d$-valued, $\upo$-free, $C_L$ is an algebraic closure of $C$, $\Gamma_L=\Q\Gamma$, and $\fM^\alg$ remains a monomial group for $L$.
If $\fm$ is a starting monomial for $P$, then we can take $c\in C_L^\times$ such that $g:=c\fm$ is an approximate zero of~$P$, hence by Corollary~\ref{cor:solve asymptotic equ, 2} applied to~$L$ in place of $K$ there is $f\sim g$ in $L$ such that $P(f)=0$.
Conversely, suppose~$L$ is a $\d$-valued  extension of $K$ and $f\in L$, $P(f)=0$, and $f\asymp\fm$. Take $c\in C_L^\times$ with $f\sim c\fm$. As in the proof of Lemma~\ref{newtondetection} we get $N_{P_{\times\fm}}(c)=0$, and so~$\fm$ is a starting monomial for~$P$.
\end{proof}

\begin{cor}\label{cor:solve asymptotic equ, 4}
Let $\fm\in\fM$. Then there are $\fm_0,\fm_1\in\fM$ with $\fm_0\prec\fm\prec\fm_1$ such that
there is no starting monomial $\fn$ for $P$ with $\fm_0\prec\fn\prec\fm_1$ and $\fn\neq\fm$.
\end{cor}
\begin{proof} Applying
Corollary~\ref{cor:sign and valuation, 1} to $P_{\times\fm}$ we get $\fm_0,\fm_1\in\fM$ with 
 $\fm_0\prec\fm\prec\fm_1$ such that $P(f)\neq 0$ for all $f$ in all $\d$-valued extensions of $K$ with $\fm_0\prec f\prec\fm_1$ and $f\nasymp\fm$. Then $\fm_0$,~$\fm_1$ have the desired property by 
 Corollary~\ref{cor:solve asymptotic equ, 3}.
\end{proof}

\noindent
Recall the concept of {\em newtonization}\/ defined in Section~\ref{sec:reldifhens}.

\begin{cor}\label{newtonization, unique} If $L$ is a newtonization of $K$, then
$L$ is an immediate $\d$-algebraic extension of
$K$, no proper differential subfield 
of~$L$ containing~$K$ is newtonian, and any newtonization of $K$ is $K$-isomorphic to~$L$.
\end{cor}
\begin{proof} By Zorn there exist immediate
asymptotically $\d$-algebraically maximal
$\d$-algebraic extensions of $K$. By
Theorem~\ref{maxnewt} any such
extension is newtonian. It remains to appeal to Corollary~\ref{newtonization}. 
\end{proof}

\noindent
{\em In the rest of this subsection $L$ is an $\upo$-free immediate extension of $K$.}\/ To the asymptotic equation \eqref{eq:asympt equ, best approx} above
corresponds the asymptotic equation
\begin{equation}\tag{E$_L$} \label{eq:asympt equ, best approx, L}
P(Y)\ =\ 0, \qquad Y\in \E_L
\end{equation}
over $L$. If  \eqref{eq:asympt equ, best approx} is quasilinear, then so is \eqref{eq:asympt equ, best approx, L}.
Moreover, if $K$ is newtonian and \eqref{eq:asympt equ, best approx} is quasilinear, then every solution of
\eqref{eq:asympt equ, best approx, L} is a solution of \eqref{eq:asympt equ, best approx}, by Lemma~\ref{newtonnoextrazeros}.
We say that $f\in L$ is {\bf quasilinear} over $K$ if $Q(f)=0$ for some $Q\in K\{Y\}^{\neq}$ with $\ndeg Q_{\times f}=1$.
Every element of $L$ that is linear over $K$ is quasilinear over $K$.
If~$K$ is newtonian, then no element of $L\setminus K$  is  quasilinear over $K$.

\index{element!quasilinear}
\index{quasilinear!element}

\begin{lemma}\label{lem:best approx is small}
Suppose $K$ is newtonian, \eqref{eq:asympt equ, best approx} is quasilinear.
Let $\E'\subseteq\E$ be $\preceq$-closed and 
let $f\in\E_L$ be such that the refinement
\begin{equation}\tag{E$'_L$}\label{eq:asympt equ, best approx, E'L}
P_{+f}(Y)\ =\ 0,\qquad Y\in\E_L'
\end{equation}
of \eqref{eq:asympt equ, best approx, L} is also quasilinear.
Let $y\preceq f$ be a solution of~\eqref{eq:asympt equ, best approx} that best approximates~$f$. Then $f-y\in\E'_L\cup\{0\}$.
\end{lemma}
\begin{proof}
Suppose $f\neq y$, and set $\fm:=\mathfrak d_{f-y}$.
Towards a contradiction, suppose $f-y\notin\E'_L$. Then $\E'_L\prec \fm\in\E$, so by quasilinearity of 
 \eqref{eq:asympt equ, best approx}:
$$1 = \ndeg_{\E_L'} P_{+f} \leq \ndeg_{\preceq \fm} P_{+f} = \ndeg_{\preceq \fm} P_{+y} 
\leq  \ndeg_{\E} P_{+y} =\ndeg_{\E} P = 1.$$
Hence the asymptotic equation
\begin{equation} \label{eq:asympt equ, best approx, proof}
P_{+y}(Y)=0,\qquad Y\preceq \fm
\end{equation}
over $K$ is quasilinear.
By  quasilinearity of \eqref{eq:asympt equ, best approx, E'L} we have
$$\ndeg_{\prec\fm} (P_{+y})_{+(f-y)}\ =\ \ndeg_{\prec\fm} P_{+f}\ \geq\ \ndeg_{\E'_L} P_{+f}\ =\ 1,$$
so $f-y$ is an approximate solution of \eqref{eq:asympt equ, best approx, proof} over $L$, by an equivalence at the beginning of Section~\ref{sec:aseq}. Take $g\in K^\times$ with
$g\sim f-y\asymp \fm$. Then $g$ is an approximate solution of 
\eqref{eq:asympt equ, best approx, proof}, and 
$\ndeg P_{+y, \times g}=\ndeg_{\preceq \fm}P_{+y}=1$. 
Then Lemma~\ref{lem:newtonian approx zero} gives $z\sim f-y$ in~$K$ such that $P(y+z)=P_{+y}(z)=0$. Using $y\preceq f$ here for the first time, we get $y+z\ne 0$, and thus $y+z$ is a better approximation to $f$ by a solution of \eqref{eq:asympt equ, best approx} than $y$, a contradiction.
\end{proof}

\subsection*{Notes and comments} Theorem~\ref{maxnewt}, together with Zorn, is one source of $\upo$-free newtonian 
$H$-asymptotic fields. The next chapter provides another source, and a more constructive procedure to build such objects.

\section{Unravelers}
\label{sec:unravelers}

\noindent
This section is somewhat technical. Our aim is 
Proposition~\ref{prop:unravelers, newtonian} below. In the next section we derive Theorem~\ref{newtmax} from Proposition~\ref{prop:unravelers, newtonian}.

{\em In this section $K$ is $\upo$-free, $\d$-valued, with
divisible value group $\Gamma$ and small derivation, and $\mathfrak M$ is a monomial {\em group}\/ of $K$.}\/ We let $\fm$, $\fn$ range over $\fM$, and~$\phi$ over the elements of 
$\fM$ that are active in $K$.
Also, $\E\subseteq K^\times$ is $\preceq$-closed, and $P\in K\{Y\}^{\neq}$. 
Consider the asymptotic equation
\begin{equation}\tag{E} \label{eq:asympt equ, unravelers}
P(Y)\ =\ 0, \qquad Y\in \E
\end{equation}
over $K$, and assume $d:=\ndeg_{\E} P\geq 2$. 
In addition we fix an $\upo$-free newtonian immediate extension $\hat K$ of~$K$, and use $\mathfrak M$ as a monomial group for $\hat K$. Set 
$$\hat\E\ :=\ \E_{\hat K}\ =\ \big\{y\in \hat{K}^\times\ : vy\in v\E\big\}.$$
Associated to \eqref{eq:asympt equ, unravelers}  we have the asymptotic equation
\begin{equation}\tag{$\hat{\operatorname{E}}$} \label{eq:hat E}
P(Y)\ =\ 0, \qquad Y\in \hat\E
\end{equation}
over $\hat K$, with the same Newton degree~$d$ as  \eqref{eq:asympt equ, unravelers}.
We assume there is given an unraveler~$(\hat f,{\hat\E}')$ for \eqref{eq:hat E} 
with $\hat f \ne 0$, $\ndeg_{\prec\hat f} P_{+\hat f}=d$.
Then $\hat{f}$ is an approximate solution of~\eqref{eq:hat E} 
of multiplicity $d$, by an equivalence at the beginning of Section~\ref{sec:aseq}. Note also that
${\hat\E}'={\E}'_{\hat K}$ for some $\preceq$-closed subset 
${\E}'$ of $\E$, and $\ndeg_{{\hat\E}'} P_{+\hat f}=d$. 
It follows that \eqref{eq:hat E} is not unraveled, and so \eqref{eq:asympt equ, unravelers} is not unraveled.

Finally, we assume $\val P_{+\hat f}<d$. Then by 
Proposition~\ref{newtpolygon}, the refinement
\begin{equation}\label{eq:hat E'}\tag{$\widehat{\operatorname{E}}{}'$}
P_{+\hat f}(Y)\ =\ 0, \qquad Y\in\hat{\E'}
\end{equation}
of \eqref{eq:hat E} has an algebraic starting monomial. 
Let $\mathfrak e$ be the largest
algebraic starting monomial for \eqref{eq:hat E'}.
The goal of this section is to prove the following:

\begin{prop}\label{prop:unravelers, newtonian}
There exists $f\in \hat K$  such that \textup{(i)} or \textup{(ii)} below holds:
\begin{enumerate}
\item[\textup{(i)}] $\hat f-f\preceq\mathfrak e$ and $A(f)=0$ for some $A\in K\{Y\}$ with $\c(A) < \c(P)$ and $\deg A =1$;
\item[\textup{(ii)}] $\hat{f}\sim f$,  $\hat f-a \preceq f-a$ for all $a\in K$, and $A(f)=0$ for some
$A\in K\{Y\}$ with $\c(A)< \c(P)$ and $\ndeg A_{\times f} =1$.
\end{enumerate}
\end{prop}

%\begin{prop}\label{prop:unravelers, newtonian}
%There exists an $f\in \hat K$  such that:
%\begin{enumerate}
%\item $f$ is linear over $K$ with $\hat f-f\preceq\mathfrak e$; or
%\item $f$ is quasilinear over $K$ with minimal annihilator over $K$ %of smaller complexity than $P$ 
%and $\hat f-a \preceq f-a$ for all $a\in K$.
%\end{enumerate}
%\end{prop}

\noindent
Below we first prove Proposition~\ref{prop:unravelers, newtonian} 
in the special case $d=\deg P$, and then show how to reduce the
general case to this special case via Lemma~\ref{lem:neglect >d}.

Note that for each $f\in\hat K$ with $\hat f-f\preceq\mathfrak e$, the pair $(f,{\hat\E}')$ is an unraveler for~\eqref{eq:hat E}, by Corollary~\ref{cor:approx unraveler}. 
So if $f\in K$ and $\hat f-f\preceq\mathfrak e$, then 
$(f,{\E}')$ with ${\E}':= {\hat \E}'\cap K^\times$ is an unraveler for 
\eqref{eq:asympt equ, unravelers}, by Lemma~\ref{lem:unravelers under immediate exts}.

Note also that towards proving Proposition~\ref{prop:unravelers, newtonian}
we may, for any given $\phi$, replace $K$,~$P$,~$\hat{K}$ by $K^{\phi}$,~$P^\phi$,~$\hat{K}^{\phi}$, without changing 
$\E$,~$d$,~$\hat{f}$,~$\hat{\E}'$,~$\mathfrak{e}$.

\subsection*{A special case} 
Set $G:=P_{+\hat f,\times\mathfrak e}\in \hat{K}\{Y\}$. Then
$\ndeg G=d$ and for $w:= \nwt(G)$ we have $N_G\in C[Y](Y')^w$.
As $N_G$ is not homogeneous, we have $w \le d-1$.
Towards proving Proposition~\ref{prop:unravelers, newtonian}
we may compositionally conjugate by an active element of
$\mathfrak{M}$, and thus arrange
$D_G=N_G$ and $R_G\prec^\flat G$. Let $\partial_0$,~$\partial_1$ be the $\hat{K}$-derivations~$\partial/\partial Y$ and $\partial/\partial Y'$ on $\hat{K}\{Y\}$.
Using notations from the end of Section~\ref{appdifpol} we
set
$$\Delta\ :=\ \big(\partial_0^{d-1-w} \partial_1^{w}\big)_{\times\mathfrak e} \in \hat{K}[[\partial]],
\qquad Q\ :=\ \Delta P \in K\{Y\}.
$$
Then by Lemmas~\ref{addconjDelta} and~\ref{multconjDelta} we get
$$Q_{+\hat f,\times\mathfrak e}\ =\ \partial_0^{d-1-w}\partial_1^w G\ \ne\ 0,$$
hence by Corollary~\ref{cor:partials and Newton pol},
$$\partial_0^{d-1-w}\partial_1^{w}N_G\ =\ D_{Q_{+\hat{f},\times \mathfrak{e}}}\ =\ N_{Q_{+\hat{f},\times \mathfrak{e}}}, \qquad
R_{Q_{+\hat{f},\times \mathfrak{e}}}\ \prec^{\flat}\ Q_{+\hat{f},\times \mathfrak{e}},$$ 
so $N_{Q_{+\hat{f},\times \mathfrak{e}}}\in C[Y]$ has degree $1$, and thus $\ndeg_{\preceq\mathfrak e} Q_{+\hat f} = 
\ndeg Q_{+\hat f,\times\mathfrak e}=1$.
Then the asymptotic equation
$$Q_{+\hat f}(Y)=0,\qquad Y\preceq\mathfrak e$$
over $\hat K$ is quasilinear.
Note that 
$$1\ =\ \ndeg_{\preceq\mathfrak e} Q_{+\hat f}\ \leq\ \ndeg_{\hat\E} Q_{+\hat f}\ =\ \ndeg_{\hat \E} Q\ =\ \ndeg_{\E} Q.$$

\begin{lemma}\label{lem:unravelers, newtonian, 1}
Suppose $\mathfrak{e}\prec \hat{f}$ and 
the asymptotic equation
\begin{equation}\label{eq:unravelers, newtonian, 1}
Q(Y)\ =\ 0, \qquad Y\in\hat\E
\end{equation}
over $\hat K$ is quasilinear. Then \eqref{eq:unravelers, newtonian, 1} has a 
solution $y\sim \hat{f}$, and if $f$ is any solution of \eqref{eq:unravelers, newtonian, 1} that best approximates $\hat f$, then $f-\hat f\preceq\mathfrak e$.
\end{lemma} 
\begin{proof} We have $\ndeg_{\prec \hat{f}}Q_{+\hat{f}}\le \ndeg_{\hat{\E}}Q_{+\hat{f}}=\ndeg_{\hat{\E}}Q=1$, and from $\mathfrak{e}\prec \hat{f}$ we get 
$1= \ndeg_{\preceq e}Q_{+\hat{f}}\le \ndeg_{\prec \hat{f}}Q_{+\hat{f}}$, and so $\ndeg_{\prec\hat{f}}Q_{+\hat{f}}=1$.
Thus \eqref{eq:unravelers, newtonian, 1} has a 
solution $y\sim \hat{f}$.  It remains to apply
Lemma~\ref{lem:best approx is small} with $\hat{K}$ in the role of both~$L$ and~$K$, and $Q$,~$\hat{f}$,~$f$ in the role of $P$,~$f$,~$y$
in that lemma.
\end{proof} 

\noindent
If $\deg P=d$, then $\deg Q=1$ and so \eqref{eq:unravelers, newtonian, 1} is automatically quasilinear. It follows that we are in case~(i)
of Proposition~\ref{prop:unravelers, newtonian} when $\deg P=d$:

\begin{cor}\label{cor:unravelers, newtonian, deg P=d}
Suppose $\deg P=d$. Then there exist $f\in \hat K$ and $A\in K\{Y\}$ such that $\hat f-f\preceq\mathfrak e$, $A(f)=0$, $\c(A)< \c(P)$, and $\deg A =1$. 
\end{cor}
\begin{proof} If $\hat{f}\preceq \mathfrak{e}$, then we can take
$f=0$ and $A=Y$. Assume $\mathfrak{e} \prec \hat{f}$. Then 
Lemmas~\ref{lem:unravelers, newtonian, 1} and~\ref{lem:best approx}  give a solution $f$ of  \eqref{eq:unravelers, newtonian, 1} with $f-\hat{f}\preceq\mathfrak e$. Now~${Q(f)=0}$, so $A:= Q$ works. (Recall that $d\ge 2$.) 
\end{proof}

%\begin{cor}\label{cor:unravelers, newtonian, deg P=d}
%Suppose $\deg P=d$. Then there exists $f\in \hat K$ such that
%$f$ is linear over $K$ and 
%$\hat f-f\preceq\mathfrak e$. 
%\end{cor}
%\begin{proof} If $\hat{f}\preceq \mathfrak{e}$, take
%$f=0$. Assume $\mathfrak{e} \prec \hat{f}$. Then 
%Lemmas~\ref{lem:unravelers, newtonian, 1} and~\ref{lem:best approx}  %give a solution $f$ of  \eqref{eq:unravelers, newtonian, 1} with 
%$f-\hat{f}\preceq\mathfrak e$. Now $Q(f)=0$, so 
%$f$ is linear over $K$.  
%\end{proof}

\noindent
Strictly speaking, Corollary~\ref{cor:unravelers, newtonian, deg P=d}
has been derived assuming that $D_G=N_G$ and $R_G \prec^{\flat} G$, but it holds without this assumption, using a
compositional conjugation as indicated in the beginning of this subsection. 

\subsection*{Tschirnhaus refinements} 
Let 
$\mathfrak f:=\mathfrak d_{\hat f}$,
put $H:=P_{\times\mathfrak f}\in K\{Y\}$ and $w:= \nwt(H)$, so
$w\le d$. If $\mathfrak{e} \succeq \mathfrak f$, then case~(i)
of Proposition~\ref{prop:unravelers, newtonian} holds for $f=0$,
with $A:= Y$. For the rest of the proof of this proposition
we assume $\mathfrak e \preceq \mathfrak{f}$. Then
$$ d\ =\ \ndeg_{\preceq \mathfrak e} P_{+\hat{f}}\ \le\ \ndeg_{\preceq \mathfrak f}P_{+\hat{f}}\ =\ \ndeg_{\preceq \mathfrak f} P\ \le\ d,$$
so $\ndeg H=d$. 
Towards proving Proposition~\ref{prop:unravelers, newtonian}
we may compositionally conjugate by an active element
of $\frak{M}$ to arrange that 
$D_H=N_H\in C[Y](Y')^w$ and $R_H\prec^\flat H$.  
Then 
$D_{H^\phi}=N_H$ and $R_{H^\phi} \prec^\flat H^\phi$ for $\phi\preceq 1$,
by Lemma~\ref{cle}. 
Let $\partial_0$ and $\partial_1$ be the $K$-derivations $\partial/\partial Y$ and $\partial/\partial Y'$ on $K\{Y\}$ and define $\Delta\in K[[\partial]]$ by
$$%\begin{equation}\label{eq:def of Delta}
\Delta:=
\begin{cases}
\displaystyle\left(\partial_0^{d-1-w}\partial_1^{w}\right)_{\times\mathfrak f} & \text{if $w\le d-1$,} \\[1em]
\displaystyle\left(\partial_1^{d-1}\right)_{\times\mathfrak f} & \text{if $w=d$,}
\end{cases}
$$%\end{equation}
and set $Q:=\Delta P \in K\{Y\}$. If $w\le d-1$, 
then $Q_{\times \mathfrak f}=\partial_0^{d-1-w}\partial_1^wH$, while if $w=d$, then $Q_{\times \mathfrak f}=\partial_1^{d-1}H$. 
Thus by Corollary~\ref{cor:partials and Newton pol},
$$N_{Q_{\times \mathfrak f}}\ =\ D_{Q_{\times \mathfrak f}}, \qquad
R_{Q_{\times \mathfrak f}}\ \prec^{\flat}\ Q_{\times \mathfrak f},\qquad 
\ndeg Q_{\times \mathfrak f}\ =\ \ddeg Q_{\times\mathfrak f}\ =\ 1,$$  so $\ddeg Q_{\times \mathfrak f}^{\phi}=1$ for all $\phi\preceq 1$, and the asymptotic equation
\begin{equation}\label{eq:Qf}
Q(Y)\ =\ 0,\qquad Y\preceq\mathfrak f
\end{equation}
over $\hat K$ is quasilinear. Thus by Corollary~\ref{cor:NP homog deg 1}:

\begin{lemma}\label{lem:NP homog deg 1, applied} \marginpar{holds also for $f=0$ if $Q(0)=0$}
Suppose $f\in\hat K$ is a solution of  \eqref{eq:Qf}. Then for all
$\mathfrak g\in\fM$ with $\mathfrak g\preceq \mathfrak f$ we have 
$\val Q_{+f,\times\mathfrak g} = 1$ and
$N_{Q_{+f,\times\mathfrak g}}=D_{Q_{+f,\times\mathfrak g}}\in C^\times Y\cup C^\times Y'$. 
Thus $Q_{+f}$ has no algebraic starting monomials $\preceq \mathfrak{f}$, and if $\mathfrak g\in\fM$ with $\mathfrak g\preceq\mathfrak f$ is a starting monomial for $Q_{+f}$, then each $\hat g\asymp\mathfrak g$ in $\hat K$ is an approximate zero of $Q_{+f}$.
\end{lemma}

\begin{lemma}\label{lem:hatf}
The element $\hat f$ of $\hat K$ is an approximate solution of \eqref{eq:Qf}.
\end{lemma}
\begin{proof}
Take $c\in C^\times$ such that $\hat f\sim c\,\mathfrak f$. Since $\hat f$ is an approximate zero of $P$ of multiplicity~$d=\ndeg P_{\times\mathfrak f}$, Corollary~\ref{cor:multiplicity} gives
$N_H = a\, (Y-c)^{d-w}\,(Y')^w$ where $a\in C^\times$. If $w\le d-1$,
then 
$N_{Q_{\times\mathfrak f}} = a\,(d-w)!w!\,(Y-c)$ by Corollary~\ref{cor:partials and Newton pol}, so $N_{Q_{\times\mathfrak f}}(c)=0$. If $w=d$, then $N_{Q_{\times\mathfrak f}}=a\,d!\,Y'$, 
and again $N_{Q_{\times\mathfrak f}}(c)=0$. 
\end{proof} 

\noindent
Let $f \in \hat K$ with $f\sim\hat f$; then $\ndeg_{\prec\mathfrak f} P_{+f}=\ndeg_{\prec\mathfrak f} P_{+\hat f}=d$, that is,
the refinement
\begin{equation}\tag{T}\label{eq:Tsch}
P_{+f}(Y)\ =\ 0,\qquad Y\prec\mathfrak f
\end{equation}
of \eqref{eq:hat E} still has Newton degree~$d$.
Moreover, the refinement
\begin{equation}\tag{$\Delta$T}\label{eq:TschQ}
Q_{+f}(Y)\ =\ 0,\qquad Y\prec\mathfrak f
\end{equation}
of \eqref{eq:Qf} is also quasilinear: use Lemmas~\ref{lem:best approx} and~\ref{lem:hatf}, and
Corollary~\ref{cor:best approx is small}
to get a solution $f_0\in\hat K$ of \eqref{eq:Qf} that best approximates $\hat f$; then $f_0\sim\hat f\sim f$, and thus
$\ndeg_{\prec\mathfrak f} Q_{+f}=\ndeg_{\prec\mathfrak f} Q_{+f_0}=1$ by Lemma~\ref{lem:NP homog deg 1, applied}, as claimed.

\begin{definition} A {\bf Tschirnhaus refinement} of \eqref{eq:hat E} is an asymptotic equation~\eqref{eq:Tsch}
over $\hat{K}$ with $\hat{f}\sim f\in \hat K$ such that some solution $f_0\in\hat K$ of  \eqref{eq:Qf} (taken over~$\hat{K}$) best approximates $\hat f$ and satisfies
$f_0-\hat f\sim f-\hat f$.
Given $f, \hat g\in\hat K$ and $\fm$ with $\fm\prec f-\hat f\preceq \hat g\prec\mathfrak f$ (so $f\sim \hat f$),
the  refinement
\begin{equation}\tag{TC}\label{eq:TschComp}
P_{+(f+\hat g)}(Y)\ =\ 0,\qquad Y\preceq\fm
\end{equation}
of \eqref{eq:Tsch}  is said to be {\bf compatible} 
with \eqref{eq:Tsch} if it has Newton degree~$d$ and $\hat g$ is not an approximate solution of~\eqref{eq:TschQ}. In this compatibility definition we do not require~\eqref{eq:Tsch} to be a Tschirnhaus refinement of \eqref{eq:hat E}.
\end{definition}

\index{refinement!Tschirnhaus}
\index{Tschirnhaus refinement}
\index{compatible refinement}
\index{refinement!compatible}

\begin{lemma}\label{lem:approx hat f, 1}
Let $f,f_0,\hat g\in\hat K$ and $\fm$ be such that
$\fm \prec f_0-\hat f\sim f-\hat f \preceq \hat g\prec\mathfrak f$, and
\eqref{eq:TschComp} has Newton degree $d$. Then $\hat g$ is an approximate solution of \eqref{eq:Tsch} and of
\begin{equation}\tag{T$_0$}\label{eq:Tsch0}
P_{+f_0}(Y)\ =\ 0,\qquad Y\prec\mathfrak f.
\end{equation}
\end{lemma}
\begin{proof}
Since $\ndeg_{\preceq\fm} P_{+(f+\hat g)}=d=\ndeg_{\prec\mathfrak f} P_{+f}$,
Lemma~\ref{lem:equal ndeg} yields that $\hat g$ is an  approximate solution of \eqref{eq:Tsch}.
Since $f_0-f\prec \hat f-f\preceq\hat g$, we have
$\ndeg_{\prec\hat g} P_{+(f_0+\hat g)}=\ndeg_{\prec\hat g} P_{+(f+\hat g)}\geq 1$. Hence
$\hat g$ is also an  approximate solution of \eqref{eq:Tsch0}.
\end{proof}

\begin{lemma}\label{lem:approx hat f, 2}
Let $f,f_0,\hat g\in\hat K$ satisfy $f_0-\hat f\sim f-\hat f\preceq \hat g\prec\mathfrak f$. Then $\hat g$ is an approximate solution of \eqref{eq:TschQ} if and only if $\hat g$ is an approximate solution of
\begin{equation}\tag{$\Delta$T$_0$}\label{eq:TschQ0}
Q_{+f_0}(Y)\ =\ 0,\qquad Y\prec\mathfrak f.
\end{equation}
\end{lemma}
\begin{proof}
We have
$\ndeg_{\prec\hat g} Q_{+(f_0+\hat g)}=\ndeg_{\prec\hat g} Q_{+(f+\hat g)}$,
since ${f_0-f\prec \hat f-f \preceq \hat g}$.
\end{proof}

\noindent
In the rest of this subsection we fix a Tschirnhaus refinement  \eqref{eq:Tsch}  of \eqref{eq:hat E}. If
$\mathfrak e\prec f-\hat f$, then
compatible refinements of \eqref{eq:Tsch} are easy to come by:
 
\begin{lemma}\label{lem:how to get compatible}
Suppose   
$\mathfrak e\prec f-\hat f$. Then  taking $\hat g:=\hat f-f$ and $\fm:=\mathfrak e$, the re\-fine\-ment~\eqref{eq:TschComp} of \eqref{eq:Tsch} is compatible.
\end{lemma}
\begin{proof} \eqref{eq:TschComp} has Newton degree~$d$, since
$\ndeg_{\preceq\mathfrak e} P_{+(f+\hat g)}=\ndeg_{\preceq\mathfrak e} P_{+\hat f}=d$.
Take a solution $f_0\in\hat K$ of  \eqref{eq:Qf} that best approximates $\hat f$ with
$f-\hat f\sim f_0-\hat f$. By Lemma~\ref{lem:approx hat f, 2}, if $\hat g$ is an approximate solution of \eqref{eq:TschQ},
then $\hat g$ is also an approximate solution of \eqref{eq:TschQ0}.
But \eqref{eq:TschQ0} has no approximate solution~$\sim\hat f-f_0$ in~$\hat K$: if it had one, then by Lemma~\ref{cor:best approx is small}
it has a solution $y\sim \hat{f}-f_0$, so $Q(f_0+y)=0$,
and $f_0+y$ would be a solution of \eqref{eq:Qf} with
$f_0+y-\hat{f}=y-(\hat{f}-f_0) \prec \hat{f}-f_0$, contradicting
that $f_0$ is a best approximation. 
Thus~\eqref{eq:TschComp} is compatible.
\end{proof}

\noindent
We now describe the effect of multiplicatively conjugating by $\mathfrak f$ on the above:

\begin{remarkNumbered}\label{rem:reduce to f=1}
Consider the asymptotic equation
\begin{equation}\label{eq:f-1 E}\tag{$\mathfrak f^{-1}$E}
P_{\times \mathfrak f}(Y)\ =\ 0, \qquad Y \in \mathfrak f^{-1}\E
\end{equation}
over $K$.
By Lemma~\ref{lem:unraveler mult conj}, $(\mathfrak f^{-1}\hat f,\mathfrak f^{-1}\hat \E')$ is an unraveler for 
\begin{equation}\label{eq:f-1 hat E}\tag{$\mathfrak f^{-1}\hat{\operatorname{E}}$}
P_{\times \mathfrak f}(Y)\ =\ 0, \qquad Y \in \mathfrak f^{-1}\hat\E
\end{equation}
over $\hat K$, and
$\ndeg_{\prec 1} (P_{\times \mathfrak f})_{+ \mathfrak f^{-1}\hat f}=\ndeg_{\mathfrak f^{-1}\hat\E} P_{\times \mathfrak f}=d$. Moreover,
\begin{equation}\label{eq:f-1 Tsch}\tag{$\mathfrak f^{-1}$T}
(P_{\times\mathfrak f})_{+\mathfrak f^{-1}f}(Y)\ =\ 0,\qquad Y\prec 1
\end{equation}
is a  Tschirnhaus refinement of \eqref{eq:f-1 hat E}, and if \eqref{eq:TschComp} is a compatible refinement of~\eqref{eq:Tsch}, then
\begin{equation}\tag{$\mathfrak f^{-1}$TC}
(P_{\times\mathfrak f})_{+\mathfrak f^{-1}f+\mathfrak f^{-1}\hat g}(Y)\ =\ 0,\qquad Y\preceq \mathfrak f^{-1}\fm
\end{equation}
is a compatible refinement of \eqref{eq:f-1 Tsch}.
\end{remarkNumbered}

\subsection*{The Slowdown Lemma}
%As in the previous subsection we let 
%$\mathfrak f=\mathfrak d_{\hat f}$
%and assume
%$D_H=N_H\in C[Y](Y')^\N$ and $R_H\prec^\flat H$ for 
%$H=P_{\times\mathfrak f}\in K\{Y\}$. We also let 
%$Q=\Delta P$ where $\Delta$ is as in \eqref{eq:def of Delta}.
Let \eqref{eq:Tsch} be a Tschirnhaus refinement of \eqref{eq:hat E} and \eqref{eq:TschComp} a compatible refinement of  \eqref{eq:Tsch}, and set $\mathfrak g=\mathfrak d_{\hat g}$. 
The goal of this subsection is the proof of the following lemma
to the effect that the step from~\eqref{eq:asympt equ, unravelers}  to
 \eqref{eq:Tsch} is much larger than the step from
\eqref{eq:Tsch} to \eqref{eq:TschComp}:

\begin{lemma}[Slowdown]\label{lem:Tsch}
$\displaystyle\frac{\fm}{\mathfrak g}\flatter \frac{\mathfrak g}{\mathfrak f}$, where $\fm$ is the monomial appearing in \eqref{eq:TschComp}.
\end{lemma}

\noindent
We first establish an auxiliary result used in the proof of Lemma~\ref{lem:Tsch}:

\begin{lemma}\label{lemlem:Tsch}
Suppose $\mathfrak f=1$. Then $Q_{+f}(\hat g)\asymp_{\mathfrak g} \mathfrak{g}\, Q_{+f}$.
\end{lemma}
\begin{proof}
Let $f_0\in\hat K$ be a solution of  \eqref{eq:Qf} that best approximates $\hat f$ and satisfies
$f-\hat f\sim f_0-\hat f$. Then $f_0\sim f \sim \hat f \asymp \mathfrak f= 1$. Take
$c\in C^\times$ such that $f_0\sim c$.
%We establish three claims.

\claim[1]{$f_0-c \prec^\flat 1$.} 

\noindent
Suppose otherwise. Then 
$f_0-c \asymp^\flat 1$. By Corollary~\ref{cor:partials and Newton pol}, 
$$Q\ =\ \mathfrak d_Q D_Q + R_Q\qquad\text{where $D_Q\in C^\times (Y-c)\cup C^\times Y'$,  and $R_Q\prec^\flat Q$,}$$
hence $Q_{+c}=\mathfrak{d}_Q(D_Q)_{+c}+(R_Q)_{+c}$ with $(R_Q)_{+c}\asymp R_Q \prec^{\flat}Q\asymp Q_{+c}$, so
$$\mathfrak d_{Q_{+c}}\ =\ \mathfrak d_Q, \quad D_{Q_{+c}}\ =\ (D_Q)_{+c}\in C^\times Y\cup C^\times Y', \quad
R_{Q_{+c}}\ =\ (R_Q)_{+c}\ \prec^\flat\ Q_{+c}.$$
Therefore, with $\mathfrak d=\mathfrak d_{f_0-c}$ we have
$N_{Q_{+c,\times \mathfrak d}}\in C^\times  Y$ by Lemma~\ref{dn2a}, so
$\mathfrak d$ is not a starting monomial for $Q_{+c}$, contradicting $Q_{+c}(f_0-c)=Q(f_0)=0$.

\claim[2]{$\mathfrak g\prec^\flat 1$.} 

\noindent
Suppose not. Then $\mathfrak g \asymp^\flat 1$. 
As in the proof of Claim~1 we get $D_{P_{+c}}\in C[Y](Y')^\N$ and $R_{P_{+c}}\prec^\flat P_{+c}$, 
so by  Claim~1 and Corollary~\ref{dn2a, cor}, 
$$N_{P_{+f_0,\times\mathfrak g}}\ =\ N_{P_{+c,+(f_0-c),\times\mathfrak g}}\ =_{\operatorname{c}}\  N_{P_{+c,\times \mathfrak g}}\in C^\times Y^\N.$$
Therefore $\hat g$ is not an approximate zero of $P_{+f_0}$, contradicting Lemma~\ref{lem:approx hat f, 1}.  

\claim[3]{$Q_{+f}(\hat g)\asymp_{\mathfrak g} Q_{+f,\times\mathfrak g}$.}

\noindent
Since $\hat g$ is not an approximate solution of \eqref{eq:TschQ}, 
it is not an approximate solution of~\eqref{eq:TschQ0} by Lemma~\ref{lem:approx hat f, 2}, and therefore $N_{Q_{+f,\times \mathfrak g}}=N_{Q_{+f_0,\times\mathfrak g}}\in C^\times Y$ by
Co\-rol\-lary~\ref{dn2cor} and Lemma~\ref{lem:NP homog deg 1, applied}.
Take $\phi\preceq 1$ such that $D_{Q^\phi_{+f,\times \mathfrak g}}=N_{Q_{+f,\times \mathfrak g}}$. Then
$$Q^\phi_{+f,\times \mathfrak g}\ =\ aY+R\qquad \text{where $a\in K^\times$, $R\in \hat{K}^{\phi}\{Y\}$, $R\prec a$,}$$
so with $\hat{g}=u\mathfrak{g}$, $u\asymp 1$ in $\hat{K}$,
$$Q_{+f}(\hat g)\ =\  Q_{+f,\times\mathfrak g}(u)\ =\ Q^\phi_{+f,\times\mathfrak g}(u)\ =\ au+R(u)\ \sim\ au$$
and hence $Q_{+f}(\hat g)\asymp a \asymp   Q^\phi_{+f,\times\mathfrak g}$.
We have $Q^\phi_{+f,\times\mathfrak g}\asymp_{\mathfrak g}  Q_{+f,\times\mathfrak g}$ by Claim~2 and
Corollary~\ref{cor:Pphi flat equivalence}.

\medskip
\noindent
By Claim~3 and Lemma~\ref{v-under-conjugation}(i)
we have 
$$Q_{+f}(\hat g)\ \asymp_{\mathfrak g}\ Q_{+f,\times\mathfrak g}\ =\ Q_{\times\mathfrak{g},+f/\mathfrak{g}}\ \sim\ Q_{\times\mathfrak{g},+f_0/\mathfrak{g}}\ =\  Q_{+f_0,\times\mathfrak g}.$$
Corollary~\ref{cor:Pg flat equivalence, 1} and $\val(Q_{+f_0})=\ddeg(Q_{+f_0})=1$ give
$Q_{+f_0,\times\mathfrak g} \asymp_{\mathfrak g} \mathfrak{g}\,Q_{+f_0}$.
Putting everything together, we get $Q_{+f}(\hat g)\asymp_{\mathfrak g} \mathfrak{g}\,Q_{+f_0} \sim \mathfrak{g}\,Q_{+f}$.
\end{proof}

\begin{proof}[Proof of Lemma~\ref{lem:Tsch}] Until further notice
we assume $\mathfrak f=1$. 
Put $F:=P_{+f}$ and $G:=Q_{+f}$. Note: $\ddeg F_{+\hat g}=\ddeg F=\ddeg P=d$, by Lemma~\ref{dn1,new}.

\claim[1]{$\mathfrak{g}\,(F_{+\hat g})_d\preceq_{\mathfrak g} (F_{+\hat g})_{d-1}$.}

\noindent 
We have $G(\hat g)\asymp_{\mathfrak g} \mathfrak{g}\,G$ by Lemma~\ref{lemlem:Tsch}.
If $w\le d-1$, then by \eqref{eq:partial ders and operations, 1},
$$G(\hat g)\ =\ Q(f+\hat{g})\ =\ \big(\partial_0^{d-1-w}\partial_1^wP\big)(f+\hat g)$$ 
is, up to a factor from $\Q^{\times}$, 
the coefficient of $Y^{d-1-w}(Y')^w$ in $P_{+(f+\hat g)}=F_{+\hat g}$, and if $w=d$, then
$G(\hat g)$ is likewise, up to a factor from $\Q^{\times}$,
the coefficient of $(Y')^{d-1}$ in~$F_{+\hat g}$. Thus 
$\mathfrak{g}\,G\preceq_{\mathfrak g} (F_{+\hat g})_{d-1}$. 
 Since $G=Q_{+f}\asymp Q\asymp P$ and
$F_{+\hat g}\sim F\asymp P$ with 
$\ddeg F_{+\hat g}=d$,
we  get
$G\asymp F_{+\hat g}\asymp (F_{+\hat g})_d$, and the claim follows.

\claim[2]{$\fn\prec_{\mathfrak g} \mathfrak g\ \Longrightarrow\ 
(F_{+\hat g,\times\fn})_{d}\prec_{\fn} (F_{+\hat g,\times\fn})_{d-1}$.}

\noindent
Suppose $\fn \prec_{\mathfrak g} \mathfrak g$. 
Then $\fn \prec_{\fn}\mathfrak g$ by Lemma~\ref{symflat}. Also 
$(F_{+\hat g})_d\ne 0$, hence
$$   \fn\,(F_{+\hat g})_d\ \prec_{\fn}\ \mathfrak{g}\,(F_{+\hat g})_{d}.$$ From $\fn\preceq \mathfrak{g}\prec 1$ we get
$\mathfrak{g}\flattereq \mathfrak{n}$, and so by Claim~1 and Lemma~\ref{lem:flat equivalence},
$$ \mathfrak{g}\,(F_{+\hat g})_{d}\ \preceq_{\fn}\ (F_{+\hat g})_{d-1}.$$
This, together with the previous display, yields 
$\fn\,(F_{+\hat g})_d \prec_{\fn} (F_{+\hat g})_{d-1}$.

By Lemma~\ref{lem:Pg flat equivalence} we have for all $i\in \N$,
$$(F_{+\hat g,\times\fn})_{i}\ =\ \big((F_{+\hat g})_{i}\big)_{\times\fn}\ 
\asymp_{\fn}\ \fn^{i}\,(F_{+\hat g})_{i}. $$
Using this for $i=d-1$ and $i=d$ yields the claim.

\claim[3]{$\fn \prec_{\mathfrak g} \mathfrak g\ \Longrightarrow\ \ndeg F_{+\hat g,\times\fn} \leq d-1$.}

\noindent
Assume $\fn\prec_{\mathfrak g} \mathfrak g$. First we note that 
$\ndeg F_{+\hat g,\times\fn}\le d$, since
$$\ndeg F_{+\hat g,\times\fn}\ \le\ \ndeg F_{+\hat g}\ =\ \ndeg F\ =\ \ndeg P_{+f}\ =\ \ndeg P\ =\ d.$$ From $\mathfrak{g}\prec^{\flat} 1$ by Claim 2 in the proof of Lemma~\ref{lemlem:Tsch}, we get $\fn\prec^{\flat} 1$. It remains to use Corollary~\ref{cor:Pphi flat equivalence} and Claim~2 in the present proof.

\medskip\noindent
Now $\ndeg F_{+\hat g,\times \fm}=d$, and so $\mathfrak{g} \preceq_{\mathfrak g} \fm$ by Claim~3, hence 
$\fm/\mathfrak{g}\asymp_{\mathfrak{g}} 1$, and thus $\fm/\mathfrak{g} \flatter \mathfrak{g}$, which is the content
of the Slowdown Lemma for $\mathfrak{f}=1$. The general case is reduced to this special case by multiplicative conjugation as in
Remark~\ref{rem:reduce to f=1}. 
\end{proof}

\subsection*{One more lemma} First an immediate consequence of 
Lemmas~\ref{lem:how to get compatible} and~\ref{lem:Tsch}: 

\begin{cor} \label{corol:Tsch}
Suppose \eqref{eq:Tsch} is a Tschirnhaus refinement of \eqref{eq:hat E}. Then
$$ \mathfrak e\ \prec\ \hat{f} - f\quad \Longrightarrow\quad \displaystyle\frac{\mathfrak e}{\hat{f}-f}\ \flatter\ \frac{\hat{f}-f}{\hat{f}}.$$
\end{cor}

\begin{lemma}\label{cor:Tsch}
Suppose \eqref{eq:Tsch} is a Tschirnhaus refinement of \eqref{eq:hat E} and $\mathfrak e\prec\hat{f}-f$. Then, with $F:= P_{+f}$, 
$\hat g:=\hat f-f$, and $\mathfrak g := \mathfrak d_{\hat g}$, the asymptotic equation
\begin{equation}\label{eq:E''}\tag{$\hat{\operatorname{E}}_{\leq d}$}
F_{\leq d}(Y)\ =\ 0,\qquad Y\preceq\mathfrak g
\end{equation}
has Newton degree~$d$. Moreover, 
$(\hat g,\hat\E'')$, where $\hat\E'' := \{ y\in \hat\E' : y \prec \mathfrak g\}$, is an unraveler for \eqref{eq:E''}, and
$\mathfrak e$ is the largest algebraic starting monomial for the unraveled asymptotic equation
\begin{equation}\label{eq:E'''}\tag{$\hat{\operatorname{E}}{}'_{\leq d}$}
(F_{\leq d})_{+\hat g}(Y)\ =\ 0,\qquad Y\in\hat\E''
\end{equation}
over $\hat K$.
\end{lemma}
\begin{proof}
Using Corollary~\ref{corol:Tsch}, apply Lemma~\ref{lem:neglect >d} to 
$\hat K$, $\hat f$, $f$, $\hat g$, $\hat\E$, $\hat\E'$ in the role of  $K$, $f$, $\tilde f$, $g$, $\E$, $\E'$, respectively.
\end{proof} 
 
\noindent
We are now ready to turn to the proof of  Proposition~\ref{prop:unravelers, newtonian}.
 
\subsection*{Proof of  Proposition~\ref{prop:unravelers, newtonian}} 
By Lemmas~\ref{lem:hatf}, \ref{cor:best approx is small}, and
\ref{lem:best approx} we have
a solution $f_0\sim \hat{f}$ in $\hat K$ of \eqref{eq:Qf} that best approximates~$\hat f$. 
If 
 $\hat f-a\preceq f_0-a$ for all $a\in K$, then (ii) holds for 
 $f:=f_0$, witnessed by $A:=Q$. 
Otherwise, take $f=a\in K$ such that
$\hat{f}-f \succ f_0-f$, that is, $f_0-\hat f\sim f-\hat f$. Then
$f\sim \hat{f}$, since $f_0\sim \hat{f}$, and so \eqref{eq:Tsch}
is a Tschirnhaus refinement of \eqref{eq:hat E}.  
 If $\hat{f}-f\preceq \mathfrak{e}$, then
we are in case (i), as witnessed by $A:= Y-f$. So assume from now on that
$\mathfrak{e}\prec \hat{f}-f$, and
set $\hat g:=\hat f-f$, $\mathfrak g := \mathfrak d_{\hat g}$. 
Consider the asymptotic equation~\eqref{eq:E''}  introduced in
Lemma~\ref{cor:Tsch}. By that lemma,~\eqref{eq:E''} 
has Newton degree~$d$, and $(\hat g,\hat\E'')$,
where $\hat\E'' = \{ y\in \hat\E' : y \prec \mathfrak g\}$, is an unraveler for~\eqref{eq:E''}, and thus $\ndeg_{\prec \hat{g}}(F_{\le d})_{+\hat{g}}=d$.  Moreover, $\mathfrak e$ is the largest algebraic starting monomial for
\eqref{eq:E'''}. 
Since $f\in K$, we can view
\eqref{eq:E''} as an asymptotic equation over~$K$. 
Thus Corollary~\ref{cor:unravelers, newtonian, deg P=d} applies to
\eqref{eq:E''} viewed as an asymptotic equation over $K$ in place of~\eqref{eq:asympt equ, unravelers} and with $(\hat g,\hat\E'')$ in place of
$(\hat f,{\hat\E}')$. So we get $g\in \hat{K}$ and $B\in K\{Y\}$
such that $\hat g-g\preceq\mathfrak e$, $B(g)=0$, 
$\c(B) < \c(F_{\le d})$, and $\deg B=1$. 
Then $f+g$ has the property stated for $f$ in (i) as witnessed by
$A:= B_{+(-f)}$. 
 \qed

\subsection*{Easy consequences} We continue in the setting introduced at the beginning of this section. Here are some easy consequences of Proposition~\ref{prop:unravelers, newtonian}:

\begin{cor}\label{cor:unravelers, newtonian, a}
Suppose $(a_\rho)$ is a divergent pc-sequence in $K$ with pseu\-do\-li\-mit~$\hat f$ and minimal $\d$-polynomial $P$ over $K$.
Then there exist $f\in\hat K$ and $A\in K\{Y\}$ such that 
$\hat f-f\preceq\mathfrak e$, $A(f)=0$, $\c(A)< \c(P)$, and $\deg A =1$.
\end{cor}
\begin{proof} Suppose not. Then Proposition~\ref{prop:unravelers, newtonian} gives
$f\in\hat K$ and $A\in K\{Y\}^{\ne}$ such that $\hat f-a\preceq f-a$ for all $a\in K$, $A(f)=0$, and $\c(A) < \c(P)$. 
Then $f\notin K$, since $\hat{f}\notin K$. Take a divergent
pc-sequence $(b_{\sigma})$ in $K$ such that $b_{\sigma} \leadsto f$.
From $\hat{f}-b_{\sigma}\preceq f-b_{\sigma}$ for all $\sigma$, we get
$b_{\sigma} \leadsto \hat{f}$, hence $(a_{\rho})$ and $(b_{\sigma})$
are equivalent, so $a_\rho\leadsto f$, and thus  
$\hat{f}-a\asymp f-a$ for all $a\in K$. Hence $Z(K,f)=Z(K,\hat{f})$, using a notion from Section~\ref{sec:construct imm exts}. 
Then $A\in Z(K,\hat f)$, by Lemma~\ref{Z1}, hence $\c(A)\ge \c(P)$ by Corollary~\ref{zmindifpol}, a contradiction.
\end{proof}

\begin{cor}
If $K$ is newtonian, then there is an $f\in K$ with $\hat f-f\preceq\mathfrak e$. 
\end{cor}
\begin{proof} Assume $K$ is newtonian. 
Take $f\in\hat K$ as in Proposition~\ref{prop:unravelers, newtonian}.
Then $f\in K$ by Lemma~\ref{newtonnoextrazeros}. In 
case (ii) of Proposition~\ref{prop:unravelers, newtonian}, this gives $\hat{f}=f$. 
\end{proof}

\section{Newtonization} \label{sec:newtonization}

\noindent
{\em In this section we assume $K$ is $\upo$-free
and $\d$-valued with divisible value group $\Gamma$.}\/ 
We do not assume here the very special setting
introduced at the beginning of Section~\ref{sec:unravelers}, but
rather reduce to it in establishing Proposition~\ref{thm:unravelers, ndeg=1} below. This will then be used to derive Theorem~\ref{newtmax}
and related results.

\begin{prop}\label{thm:unravelers, ndeg=1}
Suppose $(a_\rho)$ is a divergent pc-sequence in~$K$, with minimal 
differential polynomial $G$ over $K$, and set $\mathbf a:=c_K(a_\rho)$. Then $\ndeg_{\mathbf a} G=1$.
\end{prop}
\begin{proof} We first arrange by compositional conjugation that
$K$ has small derivation. Set $d:= \ndeg_{\mathbf a} G$. 
Zorn gives an asymptotically
$\d$-algebraically maximal immediate $\d$-algebraic
extension $\hat K$ of~$K$. Then $\hat{K}$ is $\upo$-free by Theorem~\ref{upoalgebraic}. Now forget how we got $\hat{K}$: we only use below that it is an $\upo$-free asymptotically $\d$-algebraically maximal immediate extension of $K$. Note that then $\hat{K}$ is newtonian by Theorem~\ref{maxnewt}. 

Take $\ell\in \hat{K}$ such that $a_{\rho}\leadsto \ell$.
Note that $G$ is an element of $Z(K,\ell)$ of minimal complexity,
by Corollary~\ref{zmindifpol}. Lemma~\ref{dpkell} yields
$d \ge 1$, and provides 
$a\in K$ and $\mathfrak v\in K^\times$ such that
$a-\ell\prec\mathfrak v$ and $\ndeg_{\prec\mathfrak v} G_{+a}=d$.

Suppose towards a contradiction that $d\geq 2$. 
Then Lemma~\ref{lem:unravel min pol} provides an unraveler $(\hat f,\E)$ with $\hat{f}\ne 0$ for the asymptotic equation
$$G_{+a}(Y)\ =\ 0, \qquad Y\prec\mathfrak v$$
over $\hat{K}$ such that $\ndeg_{\prec \hat f} G_{+(a+\hat f)}=d$ and $a_\rho\leadsto a+\hat f+g$ for all $g\in\E\cup\{0\}$. 
Then $\val G_{+(a+\hat f)}< d$
by Lemma~\ref{lem:small mult}. 

Assume temporarily that $K$ has a monomial group.
Then we are in the set-up of Section~\ref{sec:unravelers}
with $P:= G_{+a}$, and can apply
Corollary~\ref{cor:unravelers, newtonian, a} to the
divergent pc-sequence $(a_{\rho}-a)$ and $P:= G_{+a}$.
This yields an $f\in \hat K$ and an $A\in K\{Y\}^{\ne}$
such that $f-\hat f\in\E\cup\{0\}$, $A(f)=0$, and $\c(A) < \c(P)$. Then $a_{\rho}-a \leadsto f$. But $P$ is a minimal
$\d$-polynomial of $(a_{\rho}-a)$ over $K$, and so we contradict Section~\ref{sec:construct imm exts}.

To reduce to the case that $K$ has a monomial group, consider
the valued differential field $\hat{K}$ with $K$ as a distinguished
subset. This structure has the first-order property
that $H\notin Z(K,\ell)$ for all $H\in K\{Y\}$ with $\c(H) < \c(G)$.
Passing to an $\aleph_1$-saturated
elementary extension of this structure and using Lemma~\ref{xs1} we arrange that 
$K$ has a monomial group while preserving
the first-order property just mentioned and other
relevant first-order properties, but some minor issues arise:
\begin{enumerate}
\item the updated $\hat{K}$ might not be asymptotically $\d$-algebraically maximal (though it remains an $\upo$-free immediate extension of $K$);
\item the old pc-sequence $(a_\rho)$ might acquire a pseudolimit in the updated $K$.
\end{enumerate}
To deal with (1), replace $\hat{K}$ by an asymptotically $\d$-algebraically maximal immediate 
$\d$-algebraic extension of~$\hat{K}$.
To deal with (2), replace $(a_\rho)$ by some
divergent pc-sequence in~$K$ with pseudolimit $\ell$. 
This updated pc-sequence still has 
minimal differential polynomial~$G$ over~$K$, and the Newton degree of $G$ in the cut defined by
this pc-sequence remains~$d$, by Corollary~\ref{zmindifpol} and Lemma~\ref{dpkell}.  
\end{proof}

\noindent
This proposition is the analogue of the henselian and $\d$-henselian configuration results~\ref{cor:hensel config} and~\ref{dhconf}.
Below we derive similar consequences from it.

\begin{cor}\label{cor:unravelers, newtonian, 0}
The following are equivalent:
\begin{enumerate}
\item[\textup{(i)}] $K$ is newtonian;
\item[\textup{(ii)}] $K$ is strongly newtonian;
\item[\textup{(iii)}] $K$ is asymptotically $\d$-algebraically maximal.
\end{enumerate}
\end{cor}
\begin{proof}
The implication (i)~$\Rightarrow$~(ii) follows from Proposition~\ref{thm:unravelers, ndeg=1}, (ii)~$\Rightarrow$~(iii) holds by Lemma~\ref{lem:strongly newtonian},  and
(iii)~$\Rightarrow$~(i) by Theorem~\ref{maxnewt}.
\end{proof}

\noindent
This result contains Theorem~\ref{newtmax}: by Lemma~\ref{nli0}
the standing assumption of this section that $K$ is $\d$-valued is satisfied under the hypotheses of Theorem~\ref{newtmax}.

\begin{cor}\label{cor:unravelers, newtonian, 1}
Suppose $K$ is newtonian. If $C$ is algebraically closed, then $K$ is weakly differentially closed. If $C$ is real closed, then every $P\in K\{Y\}$ of odd degree has a zero in $K$.
\end{cor}

\begin{proof}
Combine Corollaries~\ref{cor:solve asymptotic equ} and~\ref{cor:solve asymptotic equ, real closed} with Corollary~\ref{cor:unravelers, newtonian, 0}.  
\end{proof}

\begin{cor} Let $L$ be a newtonian $\d$-algebraic immediate extension of $K$. Then $L$ is a newtonization of $K$, as defined in Section~\ref{sec:reldifhens}. 
\end{cor}
\begin{proof} Let $E$ be a newtonian extension of $K$. To embed
$L$ over $K$ into $E$ we can assume $K\ne L$; it suffices to show that then 
$K\<a\>$ can be embedded
over $K$ into~$E$ for some~${a\in L\setminus K}$. Since $L$ is $\d$-algebraic over $K$ we can take a 
divergent pc-sequence~$(a_{\rho})$
in~$K$ with a minimal differential polynomial~$G$ 
over $K$, by Lemmas~\ref{zdf, newton},~\ref{zda, newton} and Corollary~\ref{zmindifpol}. Then $\ndeg_{\mathbf a} G=1$ by Proposition~\ref{thm:unravelers, ndeg=1}, where 
$\mathbf a:=c_K(a_{\rho})$. Then Lemma~\ref{cor:dals} yields $a\in L$ such that $G(a)=0$ and~${a_{\rho} \leadsto a}$.  
Since $K$ is $\upo$-free,  we also have $\ndeg_{\mathbf a_E} G=1$, where $\mathbf a_E:=c_E(a_{\rho})$, by Corollary~\ref{cor:Newton pol under ext}. Then
Lemma~\ref{cor:dals} gives $b\in E$ with $G(b)=0$ and $a_{\rho} \leadsto b$. The results in 
Section~\ref{sec:construct imm exts} now yield an embedding 
$K\<a\>\to E$
over $K$ sending $a$ to~$b$. 
\end{proof}

\noindent
Thus $K$ has a newtonization, and any two newtonizations of $K$ are isomorphic over~$K$, by Corollary~\ref{newtonization, unique}; this permits us to speak of \textit{the}\/ newtonization of~$K$.  

Suppose $L$ is an immediate newtonian extension of $K$. Then any $K$-embedding of the newtonization of $K$ into $L$
has image
$$\{f\in L:\ f \text{ is $\d$-algebraic over $K$}\},$$
and so we refer to this image as the {\em newtonization of $K$ in $L$.}\/

\begin{cor}\label{embimmnewtonization} Let $L$ be an immediate $\upo$-free newtonian extension of $K$. Then~$L$ embeds over $K$ into any $|\Gamma|^+$-saturated 
newtonian extension of $K$. 
\end{cor} 
\begin{proof} By passing to the newtonization of $K$ in $L$
we arrange that $K$ is newtonian. Let $K^*$ be a $|\Gamma|^+$-saturated 
newtonian extension of $K$, and let $y\in L\setminus K$;
it is enough to show that then $K\<y\>$ can be embedded over
$K$ into $K^*$. Take a divergent pc-sequence~$(a_{\rho})$ in $K$ with $a_{\rho} \leadsto y$.  Since $K$ is asymptotically $\d$-algebraically maximal, it follows from Section~\ref{sec:construct imm exts} that $(a_{\rho})$ is of $\d$-transcendental type over $K$. The saturation assumption on $K^*$ gives $z\in K^*$ such that $a_{\rho}\leadsto z$. Then 
Section~\ref{sec:construct imm exts} gives a valued differential field embedding $K\<y\>\to K^*$ over $K$ sending $y$ to $z$. 
\end{proof}

\subsection*{Preservation of newtonianity} The results of this subsection are not used later, but are included for their intrinsic interest. First a descent property of newtonianity. A special case of it says that $K$ is newtonian if $K[\imag]$ is newtonian. 

\begin{cor}\label{cor:unravelers, newtonian, 2}
Let $L$ be an algebraic extension of $K$ such that $L=K(C_L)$. If~$L$ is newtonian, then so is $K$.
\end{cor}
\begin{proof} Suppose $M$ is an immediate $\d$-algebraic extension 
of $K$. 
Then $M$ is $\d$-valued with constant field $C_M=C$. Using Corollary~\ref{amalgdifferfld} we arrange that $L$ and~$M$ are common differential subfields of some differential field. 
Then $M(C_L)$ is a differential field, and by Corollary~\ref{cor:lin disjoint over constants, 4}, the constant field of $M(C_L)$ is $C_L$.
By Lemma~\ref{lem:lin disjoint over constants}, 
$M$ and $C_L$ are linearly disjoint over $C_M=C$. 
We now use Proposition~\ref{prop:constant field ext}: it gives us a valuation on the field $M(C_L)$ that extends the valuation of $M$ with the same value group $\Gamma$ as $M$ and trivial on $C_L$, and
makes~$M(C_L)$ a $\d$-valued field extension of $H$-type of $M$.  Applying Proposition~\ref{prop:constant field ext} to~$L$, in particular its uniqueness part, we see that~$M(C_L)$ contains $L=K(C_L)$ as a $\d$-valued subfield. Thus $M(C_L)$ is an immediate $\d$-algebraic extension of $L$. 

 Assume $L$ is newtonian. Then $M(C_L)=L=K(C_L)$ by Corollary~\ref{cor:unravelers, newtonian, 0}. Since~$L$ is algebraic over $K$, we have $K(C_L)=K[C_L]$ and
 $M(C_L)=M[C_L]$ and so $K[C_L]=M[C_L]$. Since $M$ is linearly disjoint from $C_L$ over $C$, this gives $K=M$. 
 Therefore $K$ is asymptotically $\d$-algebraically maximal, and thus newtonian.
\end{proof}

\noindent
The next result is obtained by a reduction to Corollary~\ref{dhenspresalg}.

\begin{prop}\label{newtalgpres} If $K$ is newtonian, then so is any algebraic extension of $K$.
\end{prop}
\begin{proof} Assume $K$ is newtonian. Recall that in the present chapter $\phi$ ranges over the active elements of 
$K$. Note that each flattening $(K^\phi, v^{\flat}_{\phi})$ is $\d$-henselian. 
 
\claim{Each
$(K^\phi, v^{\flat}_{\phi})$ is $\d$-algebraically maximal in the sense of Chapter~\ref{sec:dh1}.}

\noindent
To show this, assume towards a contradiction that $(K^\phi, v^{\flat}_{\phi})$ is not $\d$-algebraically maximal. Then~$(K^\phi, v^{\flat}_{\phi})$ has a proper immediate
$\d$-algebraic extension in the sense of Chapter~\ref{sec:dh1}. Such an extension is not required to be asymptotic, but has small derivation, and so is actually $H$-asymptotic in view of
Lemmas~\ref{ndsupPsi=0} and~\ref{impreservesas} applied to~${(K^\phi, v^{\flat}_{\phi})}$. By uncoarsening (Lemma~\ref{Un-coarsen} and the considerations preceding it), this yields a proper immediate
$\d$-algebraic $H$-asymptotic
extension of $K^{\phi}$, and so $K^{\phi}$ would not be
asymptotically $\d$-algebraically maximal, contradicting
Corollary~\ref{cor:unravelers, newtonian, 0}. This proves our claim.

  Let $L$ be an algebraic extension of $K$. Since $\Gamma$ is assumed to be divisible, we have~${\Gamma_L=\Gamma}$. Let $\phi$ be given. The flattening
$(L^{\phi}, v^{\flat}_{\phi})$ of $L^{\phi}$ is an extension of $(K^\phi, v^{\flat}_{\phi})$ in the sense of
Chapter~\ref{sec:dh1}, and is also an algebraic extension. 
Thus by our claim and Corollary~\ref{dhenspresalg}, $(L^{\phi}, v^{\flat}_{\phi})$ is $\d$-henselian. Since $\phi$ is arbitrary,
Lemma~\ref{difhdifnewt} and the equivalence preceding it 
(both applied to $L$ instead of $K$) yield that
$L$ is newtonian.   
\end{proof}

\begin{cor}
Suppose $K$ is newtonian and real closed. Then $K[\imag]$, where $\imag^2=-1$, is linearly closed and so
every monic $A\in K[\der]^{\neq}$ is a product of monic
irreducible operators in $K[\der]$ of order~$1$ and order~$2$.
\end{cor}
\begin{proof}
By Proposition~\ref{newtalgpres} the algebraic closure $K[\imag]$ of $K$ is newtonian.
So by Corollary~\ref{cor:unravelers, newtonian, 1}, $K[\imag]$ is weakly differentially closed and hence linearly closed
by Lemma~\ref{fa0}. The rest now follows from Lemma~\ref{fa6}.
\end{proof}

\subsection*{Newton-Liouville closure} 
{\em Besides assuming that $K$ is $\upo$-free and $\d$-valued with divisible value group, we also assume in this subsection that $K$ comes equipped with a field ordering making 
$K$  an $H$-field.}\/ By Corollary~\ref{newtonization, unique}, the newtonization $K^{\operatorname{nt}}$ of~$K$ is an immediate extension of $K$.
Hence by Lemma~\ref{prehimm},
 $K^{\operatorname{nt}}$ has a unique field ordering extending that of $K$ in which the valuation ring of $K^{\operatorname{nt}}$ is convex; below we consider $K^{\operatorname{nt}}$ equipped with this ordering.
Then $K^{\operatorname{nt}}$ is a newtonian $H$-field extension of $K$, and every embedding of~$K$ into a newtonian pre-$H$-field $L$ extends to an embedding  $K^{\operatorname{nt}}\to L$
of pre-$H$-fields. 
%By Corollary~\ref{cor:Prestel}, if $C$ is real closed, 
%then so is $K^{\operatorname{nt}}$. 

\begin{lemma}\label{lem:NL closure}
There exists a newtonian Liouville closed $H$-field extension 
$ K^{\operatorname{nl}}$ of~$K$ which embeds over $K$ into
every newtonian Liouville closed $H$-field extension of $K$.
Any such $K^{\operatorname{nl}}$ is $\d$-algebraic over $K$, hence $\upo$-free, and its constant field is a real closure of $C$.
\end{lemma}
\begin{proof}
We define inductively an infinite tower $K_0\subseteq K_1\subseteq K_2\subseteq \cdots$ of $\upo$-free  $H$-field extensions of $K$ with divisible value group
as follows. Set $K_0:=K$, and assume inductively that
$K_n$ is an $\upo$-free $H$-field extension of $K$ 
with divisible value group.
For even $n$ we let $K_{n+1}$ be the Liouville closure of $K_n$.
% the definite article ``the'' is harmless here since $K_n$ is
%$\upo$-free. 
For odd $n$ we let 
$K_{n+1}:=(K_n)^{\operatorname{nt}}$ be the newtonization of $K_n$.
In both cases, $K_{n+1}$ has divisible value group, and $K_{n+1}$ is $\d$-algebraic over $K_n$ and
thus remains $\upo$-free, by Theorem~\ref{upoalgebraic}. Thus 
$K^{\operatorname{nl}}:=\bigcup_n K_n$ is a newtonian Liouville closed $H$-field extension of $K$ with the desired semiuniversal  property.
The second part of the lemma follows easily from this semiuniversal property. 
\end{proof}

\noindent
%We say that an $H$-field is {\bf Newton-Liouville closed} if it is both newtonian and Liouville closed.
Let  $E$ be an $\upo$-free $H$-field. We extend 
Lemma~\ref{lem:NL closure} to $E$ by applying it to the $\upo$-free real closed $H$-field extension $K:=E^{\operatorname{rc}}$ of $E$:

\begin{cor}\label{cor:NL closure}
There is a newtonian Liouville closed $H$-field extension~$E^{\operatorname{nl}}$ of~$E$ which embeds over $E$ into
every newtonian Liouville closed $H$-field extension of~$E$.
Any such $E^{\operatorname{nl}}$ is $\d$-algebraic over $E$, thus $\upo$-free, and its constant field is a real closure of~$C_E$.
\end{cor} 

\begin{cor} Every pre-$H$-field extends to an $\upo$-free newtonian
Liouville closed $H$-field.
\end{cor}
\begin{proof} Any pre-$H$-field extends to an $\upo$-free $H$-field
by Corollary~\ref{cor:embed into upo-free}. Now apply Corollary~\ref{cor:NL closure}. 
\end{proof}

\noindent
A {\bf Newton-Liouville closure} of $E$ is by definition a newtonian Liouville closed $H$-field extension $E^{\operatorname{nl}}$ of $E$ with the embedding property stated in Corollary~\ref{cor:NL closure}. Thus
$E$ has a Newton-Liouville closure. In Section~\ref{mctnl} below we show that $E$ has up to isomorphism over~$E$
a unique Newton-Liouville closure. 

\index{closure!Newton-Liouville}
\index{Newton-Liouville closure}
\index{H-field@$H$-field!Newton-Liouville closure}

\subsection*{Schwarz closure}
Let $E$ be an $\upo$-free $H$-field. 

\begin{prop}\label{prop:Schwarz closure}
There exists a Schwarz closed $H$-field extension of $E$ that
 embeds over $E$ into any Schwarz closed $H$-field extension of $E$.
\end{prop}
\begin{proof} Let $E^{\operatorname{nl}}$ be a Newton-Liouville closure
of $E$. Then $E^{\operatorname{nl}}$ is Schwarz closed by Corollary~\ref{2newtschwarz}, and its constant field is a real closure of $C_E$. Let
$E^{\operatorname{s}}$ be the intersection of all Schwarz closed
$H$-subfields of $E^{\operatorname{nl}}$ that contain $E$. Then
$E^{\operatorname{s}}$ is a Schwarz closed $H$-field
extension of $E$ by Lemma~\ref{lem:intersect Schwarz closed H-fields}. We show that it has the desired semiuniversal property. Let $F$ be a Schwarz closed $H$-field extension
of~$E$. Take a Newton-Liouville closure $F^{\operatorname{nl}}$ of $F$.
Then we have an $H$-field embedding 
$i\colon E^{\operatorname{nl}}\to F^{\operatorname{nl}}$ over~$E$, so $i(E^{\operatorname{s}})$ and
$F$ are Schwarz closed $H$-subfields of $F^{\operatorname{nl}}$. 
Now~$F$ and $F^{\operatorname{nl}}$
have the same constant field, so $i(E^{\operatorname{s}})\cap F$ is
a Schwarz closed $H$-subfield of $i(E^{\operatorname{s}})$ containing~$E$, by Lemma~\ref{intersectschwarz}. By the minimality property of $E^{\operatorname{s}}$ this gives $i(E^{\operatorname{s}})\subseteq F$. 
\end{proof}

\noindent
Define a {\bf Schwarz closure} \index{closure!Schwarz}\index{Schwarz!closure}\index{H-field@$H$-field!Schwarz closure} of $E$ to be a Schwarz closed $H$-field extension of~$E$ that embeds over $E$ into
every Schwarz closed $H$-field extension of $E$. So $E$ has a Schwarz closure by Proposition~\ref{prop:Schwarz closure},  and
because of the obvious minimality property of the Schwarz closure
constructed in its proof, any two Schwarz closures of $E$ are isomorphic over $E$. This allows us to speak of \textit{the}\/ Schwarz closure $E^{\operatorname{s}}$ of~$E$. Thus by the proof of~\ref{prop:Schwarz closure}: 

\begin{samepage}
\begin{enumerate}
 \item $E^{\operatorname{s}}$ is $\d$-algebraic over $E$,
\item the constant field of $E^{\operatorname{s}}$ is a real closure of $C_{E}$,
\item $E^{\operatorname{s}}$ has no proper Schwarz closed $H$-subfield containing $E$.
\end{enumerate}
\end{samepage}

%\noindent
%Corollary~\ref{cor:embed into upo-free} and 
%Proposition~\ref{prop:Schwarz closure} together show that 
%each pre-$H$-field has a Schwarz closed $H$-field extension. 
%This will be used in Section~\ref{sec:LO-cuts}. 

%% file: mt-15.tex
\chapter{Newtonianity of Directed Unions}\label{ch:newtdirun}

\setcounter{theorem}{0}

\noindent 
In this brief chapter we prove an analogue of Hensel's Lemma for $\upo$-free $\d$-valued fields of $H$-type: Theorem~\ref{thm:newtonianity of directed union}.
{\em Throughout this chapter~$K$ is an $H$-asymptotic field with asymptotic couple $(\Gamma, \psi)$, and $\gamma$ ranges over $\Gamma$.}\/

\begin{theorem}\label{thm:newtonianity of directed union} 
If $K$ is $\d$-valued with
$\der K = K$, and $K$ is a directed union of spherically complete grounded $\d$-valued subfields, then $K$ is newtonian.
\end{theorem}

\noindent
Note that by Corollary~\ref{recam}, any $K$ as in the hypothesis of this theorem is $\upo$-free. The proof of the theorem is given in Section~\ref{sec:prnewt}. As special cases we shall obtain at the end of that section:

\begin{cor}\label{cor:thmnewtonianity} The $\upo$-free $H$-fields $\T$ and $\T_{\log}$ are newtonian. 
\end{cor}

\section{Finitely Many Exceptional Values}\label{sec:excval}

\noindent
{\em In this section we assume $\sup\Psi=0$, and we let $y$ range over $K^\times$.}\/ Let $A\in K[\der]^{\neq}$.
Recall from Section~\ref{sec:ldopv} that 
$$\exc(A)\ =\ \big\{vy:\ A(y)\prec Ay \big\}\ =\ \big\{vy:\ \Pmu(Ay)>0\big\}$$ 
is the set of exceptional values for $A$.
We have
$v(\ker^{\ne} A)\subseteq \exc(A)$, so knowing $\exc(A)$
helps in locating the solutions in $K$ of the differential equation $A(y)=0$. If $K$ is $\d$-valued, then
$$\exc(A)\ =\ \big\{\gamma:\ \Pmu_A(\gamma)>0\big\},$$
and $\dim_C\ker A \le |\exc(A)|$ by Lemma~\ref{valuation-of-basis}.  

\medskip\noindent
If $K$ extends to a $\d$-henselian asymptotic field, then 
$|\exc(A)|\le \order{A}$ for all $A\in K[\der]^{\ne}$, by Lemma~\ref{dhli2}.
Using this fact we prove: 

\begin{lemma}\label{EAfinite} If $K$ is pre-$\d$-valued, then $|\exc(A)|\le \order{A}$ for all $A\in K[\der]^{\ne}$.
\end{lemma}
\begin{proof} Assume $K$ is pre-$\d$-valued. Consider first the case that $K$ is grounded (so $\max \Psi = 0$). Then $K$ has an $\upo$-free $\d$-valued extension $L$ of $H$-type, by 
Corollary~\ref{cor:embed into upo-free}. In view of Theorem~\ref{maxnewt} 
and the subsequent remarks we can pass to an immediate extension of $L$ and arrange in this way
that $L$ is also newtonian. Then the flattening $(L, v^{\flat})$ of $L$ is $\d$-henselian by Lemma~\ref{diffnewt}, and the inclusion map $K \to L$ is actually an embedding of $K$ into this
flattening. If $K$ is ungrounded, then $0\notin \Psi$, and
this case is taken care of by Corollary~\ref{aspsidh}.
\end{proof}

\section{Integration and the Extension $K(\lowercase{x})$}\label{sec:extkx} 

\noindent
{\em In this section $K$ is $\d$-valued with $\sup \Psi=0$.}\/ 
Then the equation $y'=1$ has no solution in $K$, and so we adjoin a solution:

\begin{lemma} \label{K(x)}
Let $K(x)$ be a field extension of $K$ with $x$
transcendental over $K$. Then there is a unique pair consisting of a
derivation of $K(x)$ and a valuation ring of~$K(x)$
that makes $K(x)$ a $\d$-valued extension of $H$-type of $K$ with
$x'=1$ and $x\succ 1$. This extension $K(x)$ has the same constant field as $K$,
has value group $\Gamma + \Z vx$ with  $\Gamma^{<} < nvx < 0$
for all $n\ge 1$, and has small derivation.
\end{lemma}
\begin{proof} If $\Psi < 0$, use Lemma~\ref{variant41} and subsequent remarks. If $\max \Psi = 0$, use 
Lemma~\ref{pre-extas2} and the remarks following its proof.
\end{proof}

\noindent
Note that the extension $K(x)$ of $K$ described in Lemma~\ref{K(x)} is grounded with $\max\Psi_{K(x)}=-vx>0$, has $K[x]$ as a differential subring, and that
$$\Z vx\ =\ 
\big\{\alpha\in\Gamma_{K(x)}:\ \psi(\alpha)>0\big\}\ =\ \Gamma_{K(x)}^\flat.$$

\begin{lemma}\label{lem:K=der(K+Cx)}
Assume $K$ is spherically complete. Then 
$K=\der K + C=\der(K+Cx)$.
\end{lemma}
\begin{proof}
By Lem\-ma~\ref{var43} and the remarks following it we get
$\smallo=\der\smallo$.
Now let $s\in K$; we need to show $s\in \der K + C$. For this we may assume that
$s\notin\der K$ and hence $s-a'\succeq 1$ for all $a\in K$, since  $\smallo=\der\smallo$. From
Lemma~\ref{var51} it follows that
$$S\ :=\ \big\{ v(s-a'): a\in K\big\}\ \subseteq\ \Gamma^{\leq}$$
has a largest element $\beta$, and so $\beta=0$ by Lemma~\ref{lem:remark5.1(2)}. Take $a\in K$ with $v(s-a')=0$, next take  $c\in C^\times$, $\varepsilon\in\smallo$ such that $s-a'=c(1+\varepsilon)$, and then take  $b\in\smallo$ such that
$\varepsilon=b'$. Then
$s=(a+cb)' + c\in \der K + C$.
\end{proof}

\begin{cor}\label{cor:K[x]=der(K[x])}
Suppose $K$ is spherically complete. Then $K[x]=\der\big(K[x]\big)$.
\end{cor}
\begin{proof} Let $a\in K$; we show that
$ax^n=f'$ for some $f\in K[x]$. For $n=0$, use Lemma~\ref{lem:K=der(K+Cx)}. Let $n\geq 1$ and take $b\in K$, $c\in C$ with $a=b'+c$, and
inductively take $g\in K[x]$ with $g'=bx^{n-1}$. Set $f:=\frac{c}{n+1}x^{n+1}+bx^n-ng$. Then $f'=ax^n$.
\end{proof}

\begin{cor}\label{cor:K=derK+C}
Let $E$ be a spherically complete $\d$-valued field of $H$-type with an element
$\phi$ such that $v\phi=\max \Psi_E$. Let $F$ be a 
differential field extension of $E$ such that $\phi\in \der F$. Then
$E\subseteq \der F$.
\end{cor}
\begin{proof} Take $x\in F$ with $\der x=\phi$. Let 
$\derdelta:= \phi^{-1}\der$ be the derivation of $F^\phi$.
Then the standing assumption of this section is satisfied
for $K:= E^\phi$, since $\max \Psi_K=0$. Also
$\derdelta x=1$, hence by Lemma~\ref{lem:K=der(K+Cx)}
with $\derdelta$ instead of $\der$,
$K\subseteq \derdelta(K)+ C=\derdelta(K+Cx)\subseteq \derdelta(F^\phi)$, that is, $E\subseteq \phi^{-1}\der F$,
and thus $E=\phi E\subseteq \der F$. 
\end{proof}

\subsection*{Notes and comments}
In the next volume we extend Corollary~\ref{cor:K[x]=der(K[x])} as follows: if $K$ is spherically complete and $A\in K[\der]^{\neq}$, then $K[x]=A\big(K[x]\big)$.

\section{Approximating Zeros of Differential Polynomials}

\noindent
{\em In this section $K$ is $\d$-valued with
$\sup \Psi=0$, we fix a  
$P\in K\{Y\}^{\ne}$ of order $\le r$, and we let $a$ be an element of $\mathcal O$ such that $\dval P_{+a}=1$.}\/ 

We also consider the condition 
$\ddeg P_{+a}=1$, equivalent to $P_{+a}\sim (P_{+a})_1$ under our standing assumption that $\dval P_{+a}=1$.

We have $P(a)\prec (P_{+a})_1 \asymp P_{+a}\asymp P$.
Take $\fd\in K^\times$ with $\fd\asymp P$. Then $\fd^{-1}P\in \mathcal{O}\{Y\}$ is in
dh-position at $a$ as defined in Section~\ref{preldh}. 
We have elements $a_0,\dots,a_r\in K$ with
$(P_{+a})_1=a_0Y+a_1Y'+\cdots+a_rY^{(r)}$, so 
$$A\ :=\ L_{P_{+a}}\ = a_0+a_1\der+\cdots+a_r\der^r\in K[\der]^{\neq}.$$
As in Section~\ref{sec:maxdh} where $P\asymp 1$, let 
$v(P,a)$ be the unique $\alpha\in \Gamma_{\infty}$ with
$v_A(\alpha)=v\big(P(a)\big)$; thus $v(P,a)=v(\fd^{-1}P,a)$,
and $v(P,a)=\infty$ iff $P(a)=0$. We use $v(P,a)$ as a measure of how close $a$ is to a potential zero of $P$.

{\em Throughout the rest of this section we assume $P(a)\ne 0$.}\/ Thus $v(P,a)\in \Gamma^{>}$.

\begin{lemma}\label{lem:summary} Let $b\in K$, $v(a-b)\geq v(P,a)$,  
$P(b)\prec P(a)$, and $B:= L_{P_{+b}}$. Then: 
\begin{enumerate}
\item[\textup{(i)}] $\dval P_{+b}=1$;
\item[\textup{(ii)}] $v(a-b)=v(P,a)$ and $v(P,b)>v(P,a)$; 
\item[\textup{(iii)}] for all $y\in K^\times$, if $vy=O\big(v(P,a)\big)$
and $A(y)\prec Ay$, then $B(y)\prec By$;
\item[\textup{(iv)}] $\big\{ \alpha\in\exc(A):\ \alpha=O\big(v(P,a)\big)\big\}\ \subseteq\ \exc(B)$;
\item[\textup{(v)}] if $\ddeg P_{+a}=1$, then $\ddeg P_{+b}=1$.
\end{enumerate}
\end{lemma}
\begin{proof} For (i)--(iv), apply Lemma~\ref{henseld1}
to $\fd^{-1}P$;  for (v), use Lemma~\ref{dn1,new}.
\end{proof}

\noindent
In Propositions~\ref{prop:quasilinear 1} and \ref{prop:quasilinear 2} below we indicate how to improve $\gamma:=v(P,a)$.  

\begin{lemma}\label{lem:not exceptional}
If $\gamma\notin\exc(A)$, then $P(b)\prec P(a)$ for some 
$b\in K$ with $v(a-b)=\gamma$. 
\end{lemma}
\begin{proof} We adapt an argument from the proofs of Lemmas~\ref{newtona} and~\ref{newton}. Take  $g\in K^\times$ with $vg=\gamma$, and set $Q:= P_{+a}$,  
$L:=P(a)^{-1}Q_{1,\times g}\in K\{Y\}$. Then $g\prec 1$ and 
$L\asymp 1$, so $L=\sum_{i=0}^r b_i Y^{(i)}$ with
$b_0,\dots,b_r\in\mathcal O$ and $b_i\asymp 1$ for some $i\in\{0,\dots,r\}$. As in the proof of Lemma~\ref{newtona}, we get
$$P(a+gY)\ =\ P(a)\cdot \big( 1 + L(Y) + R(Y)\big)\quad\text{where 
$\val R \geq 2$ and $R \prec 1$.}$$
Now assume that $\gamma\notin \exc(A)$. 
Then $\Pmu(Q_{1,\times g})=\Pmu(Ag)=0$, so $b_0\asymp 1$.
Take $c\in C^\times$ with $b_0c\sim -1$. Then   
$1+L(c)=1+b_0c\prec 1$ and so
$$P(a+gc)\ =\ P(a)\cdot \big( 1+L(c)+R(c)\big)\ \prec\ P(a).$$
Hence $b:=a+gc$ has the desired properties. 
\end{proof}

\begin{prop}\label{prop:quasilinear 1}
Suppose $K$ is spherically complete, $v(P,a)\notin\exc(L_{P_{+a}})$,  and there is no $b\in K$ with $v(a-b)=v(P,a)$ and $P(b)=0$. Then 
for some $b\in \mathcal O$ we have 
$v(a-b)=v(P,a)$, $P(b)\prec P(a)$, and $v(P,b)\in  \exc(L_{P_{+b}})$.
\end{prop}
\begin{proof}
Let $(a_\rho)_{\rho<\lambda}$ be a sequence 
in $\mathcal{O}$ with $\lambda$ an ordinal $>0$, $a_0=a$, and
\begin{enumerate}
  \item $\dval P_{+a_\rho}=1$, for all $\rho < \lambda$,
  \item $v(a_{\rho'} -a_\rho)=v(P,a_\rho)$ 
whenever $\rho<\rho'<\lambda$,
  \item $P(a_{\rho'})\prec P(a_\rho)$ and 
$v(P,a_{\rho'})>v(P,a_\rho)$ whenever $\rho<\rho'<\lambda$, and
\item $v(P,a_\rho)\notin\exc(L_{P_{+a_\rho}})$ for all $\rho<\lambda$.
\end{enumerate}
Note that there is such a sequence for $\lambda=1$. 
We now keep extending this sequence as in the proof of Lemma~\ref{hensel.imm} while preserving (1)--(4), using 
Lemma~\ref{lem:not exceptional} instead of Lemma~\ref{newton}.
This extension procedure cannot go on indefinitely, and it
must end in a violation of clause (4), which
gives an element $b$ as required. 
\end{proof}

\noindent
In what follows $K(x)$ is the $\d$-valued extension of $H$-type of $K$ from Lemma~\ref{K(x)}, so $x$ is transcendental over~$K$, $x'=1$, and $x\succ 1$. 
We also assume that a $\d$-valued extension $F$  of $K(x)$ of $H$-type is given such that $\Psi_F$ has a supremum in $\Gamma_F$.
We use the flattening
$v^\flat\colon F^\times \to\Gamma_F^\sharp=\Gamma_F/\Gamma_F^\flat$
of the valuation $v\colon F^\times\to\Gamma_F$ of $F$, with the associated dominance relation
$\preceq^\flat$. Here, as usual,
$$\Gamma_F^\flat\ =\ 
\big\{\alpha\in\Gamma_F:\ \psi_F(\alpha)>0\big\},$$ a convex subgroup of $\Gamma_F$.
Note that $x\asymp^\flat 1$. Also 
$\Gamma_F^\flat\cap\Gamma=\{0\}$, so
for $b\in K$ we have $b\prec 1\Longleftrightarrow b\prec^\flat 1$. 
In the next lemma we set 
$n:=\Pmu_A(\gamma)$.  

\begin{lemma} \label{lem:xexceptional}
For some $b\in K[x]$
we have $v(a-b)= \gamma+nvx$ and $P(b)\prec^\flat P(a)$.
\end{lemma}
\begin{proof}
With the notations from the proof of Lemma~\ref{lem:not exceptional}
we have  $\Pmu(Q_{1,\times g})=\Pmu(Ag)=n$, hence $b_0,\dots,b_{n-1}\prec 1\asymp b_n$.
Take $c\in C^\times$ with $b_nc\sim -1$ and set $y:=(c/n!)x^n\in K[x]$.
Then $y\asymp^\flat 1$, so $R(y)\preceq^\flat R\prec^\flat 1$, 
and thus
\begin{align*}
P(a+gy)\	 &=\ P(a)\cdot\big( 1+L(y)+R(y)\big) \\
	&=\ P(a)\cdot\left( 1 + \sum_{i=0}^{n-1} b_i\big(c/(n-i)!\big)x^{n-i} + b_nc + R(y) \right)\ \prec^{\flat}\ P(a),
\end{align*}
so that $b:= a+gy$ has the desired property.
\end{proof}

\noindent
Let $(F, v^\flat)$ be the valued differential field  
whose underlying differential field is that of $F$ and whose valuation is $v^\flat$. Let $\mathcal{O}^{\flat}:= \mathcal{O}_F^\flat$ be the valuation ring of $(F, v^\flat)$.   
Note that~$(F, v^\flat)$ is an $H$-asymptotic field extension
of $K$ with 
$\Gamma_{(F, v^\flat)} = \Gamma_F^\sharp$ and
$\sup\Psi_{(F, v^\flat)}=0$. 
Thus $\fd^{-1}P$ remains in dh-position at $a$ with respect to $(F, v^\flat)$. 

Suppose $\fd^{-1}P$ is in dh-position at
$b\in \mathcal{O}^\flat$, with respect to $(F, v^\flat)$. Then  
$v^\flat(P,b)$ is the unique $\beta\in(\Gamma_F^\sharp)_{\infty}$ 
with $v^\flat_B(\beta)=v^\flat\big(P(b)\big)$, where
$B:= L_{P_{+b}}$. An argument in the proof of Lemma~\ref{inva3}
gives that $P$ is in dh-configuration at $b$ with respect to the valuation $v$ of $F$, and so $v(P,b)\in (\Gamma_F)_{\infty}$ is defined.   
By Lemma~\ref{coarsevp} we have
$v^\flat(P,b)=v(P,b)+\Gamma_F^\flat$ if $P(b)\ne 0$.
Thus by Lemma~\ref{lem:xexceptional}, and by Lemma~\ref{henseld1} applied to $(F, v^\flat)$ in the role of~$K$, we obtain:

\begin{cor}\label{cor:xexceptional} 
There is a $b\in K[x]$
with $v^\flat(a-b)= v^\flat(P,a)$ and $P(b)\prec^\flat P(a)$. For any
$b\in F$ with $v^\flat(a-b)\geq v^\flat(P,a)$ and $P(b)\prec^\flat P(a)$, and any $y\in F^\times$:
\begin{enumerate}
\item[\textup{(i)}] $\fd^{-1}P$ is in dh-position at $b$ with respect to $(F, v^\flat)$;
\item[\textup{(ii)}] $v^\flat(a-b)= v^\flat(P,a)$ and $v^\flat(P,b)>v^\flat(P,a)$;
\item[\textup{(iii)}] if $vy=O\big(v^\flat(P,a)\big)$, and $A(y)\prec^\flat Ay$, then $B(y)\prec^\flat By$ for $B:= L_{P_{+b}}$.
\end{enumerate}
\end{cor}

\noindent
Note also that if $b\in F$ and $v^\flat(a-b)\ge v^\flat(P,a)$, then $b\preceq 1$. 
Next, take $\phi\in F$ with $v\phi = \sup \Psi_F$. We have $\phi\preceq 1/x\prec 1$, $\phi\asymp^\flat 1$,
and
$F^\phi$ is again (like $K$) a $\d$-valued field of $H$-type with $\sup\Psi_{F^\phi}=0$.
We can now state and prove another key fact:

\begin{prop}\label{prop:quasilinear 2} Assume $\ddeg P_{+a}=1$. Then for some $b\in K[x]$ with $b\preceq 1$, 
\begin{enumerate}
\item[\textup{(i)}] $P_{+b}^\phi\sim (P_{+b}^\phi)_1$ in $F^\phi\{Y\}$, and
$v(P^\phi,b)>v(P,a)$;
\item[\textup{(ii)}] $\big\{ \alpha\in\exc_K(A):\alpha=O\big(v(P,a)\big)\big\} \subseteq \exc_{F^\phi}(B^\phi)$ where $B:=L_{P_{+b}}$. 
\end{enumerate}
\end{prop}
\begin{proof} We have $\dval P_{+a}=\ddeg P_{+a}=1$ with respect to $K$ and thus with respect to the valued field extension $(F, v^\flat)$ of $K$.
Take $b\in K[x]\subseteq F$ as in Corollary~\ref{cor:xexceptional}.
Then $b\preceq 1$ and $v^\flat(a-b)>0$, so $\dval P_{+b}= \ddeg P_{+b}=1$ with respect to
$(F, v^\flat)$, by Lemma~\ref{dn1,new}, and so for $H:=(P_{+b})_1\in F\{Y\}$ we have
$$P_{+b}\ =\ H+S\quad
\text{where $S\prec^\flat H$.}$$
By Proposition~\ref{compconjval, prop}, 
$v(H^\phi) \equiv v(H)\bmod\Gamma_F^\flat$ and
$v(S^\phi) \equiv v(S) \bmod\Gamma_F^\flat$, so
$$P_{+b}^\phi\ =\ H^\phi+S^\phi\quad\text{where 
$S^\phi\prec^{\flat} H^\phi=(P_{+b}^\phi)_1$},$$
and thus $P^{\phi}_{+b}\sim (P^{\phi}_{+b})_1$. 
Also $v_{H^\phi}\big(v(P^\phi,b)\big)=v\big(P(b)\big)=v_{H}\big(v(P,b)\big)$, hence 
$$v_{H^\phi}\big(v(P^\phi,b)\big) + \Gamma_F^\flat\ =\  v^\flat\big(P(b)\big)\ =\ v_{H}\big(v(P,b)\big) + \Gamma_F^\flat\ =\ v_{H^\phi}\big(v(P,b)\big) + \Gamma_F^\flat,$$
using Proposition~\ref{compconjval, prop} for the last equality. 
Applying Lemma~\ref{coarsevp} to $F^\phi$ in the role of $K$,
using the injectivity of $v^{\flat}_{H^{\phi}}$, and using Corollary~\ref{cor:xexceptional}, we get
$$v(P^\phi,b)+\Gamma_F^\flat\ =\ v(P,b)+\Gamma_F^\flat\ =\ v^\flat(P,b)\ >\ v^\flat(P,a)\ =\ v(P,a)+\Gamma_F^\flat$$ 
and hence $v(P^\phi,b)>v(P,a)$. We have now established (i).  

For (ii), let $\alpha\in\exc_K(A)$ and
$\alpha=O\big(v(P,a)\big)$. Take $y\in K^\times$ such that
$vy=\alpha$ and $A(y)\prec Ay$. Then $A(y)\prec^{\flat} Ay$ with respect to $(F, v^\flat)$, and so $B(y)\prec^{\flat} By$
by (iii) of Corollary~\ref{cor:xexceptional}. In the common underlying valued subfield of $(F, v^\flat)$ and $(F, v^\flat)^\phi=(F^\phi, v^\flat)$ we have $B(y)=B^{\phi}(y)$, and by Proposition~\ref{compconjval, prop} we have
$By\asymp^\flat (By)^\phi=B^\phi y$ .  Hence $B^\phi(y) \prec^\flat B^\phi y$, so $B^\phi(y) \prec B^\phi y$ with respect to $F^\phi$, and thus
$\alpha=vy\in \exc_{F^\phi}(B^\phi)$. This yields~(ii). 
\end{proof}
 
\noindent
For the proof of Theorem~\ref{cor:quasilinear 2} below it is also relevant that for $B=L_{P_{+b}}$ as in~(ii) of Proposition~\ref{prop:quasilinear 2} we have $B^\phi=L_{P^\phi_{+b}}$, by Corollary~\ref{corconj-lemma}.   
 
\section{Proof of Newtonianity}\label{sec:prnewt}

\noindent
We begin with a result of independent interest and with more constructive content than Theorem~\ref{thm:newtonianity of directed union}.

\begin{theorem}\label{cor:quasilinear 2} Let
$K_0\subseteq K_1\subseteq\cdots\subseteq K_r$ be a tower
of spherically complete $\d$-valued fields of $H$-type. For
$i=0,\dots, r$, let 
$\phi_i\in K_i^\times$ be such that $v\phi_i=\max \Psi_{K_i}$. 
Assume that $\phi_i\in\der K_{i+1}$ for $i=0,\dots,r-1$.
Let $P\in K_0\{Y\}^{\ne}$ have order $\le r$, and let $a\in \mathcal{O}_{K_0}$ be such that $P_{+a}^{\phi_0}\sim (P_{+a}^{\phi_0})_1$. Then $P(b)=0$ for some
$b\in \mathcal{O}_{K_r}$.
\end{theorem}
\begin{proof} 
Towards a contradiction, assume there is no such $b$ (so
$r\ge 1$). We shall build a sequence $a_0, b_0, a_1, b_1,\dots, a_r, b_r$ with $a_i,b_i\in\mathcal O_{K_i}$ and $P^{\phi_i}_{+a_i} \sim (P^{\phi_i}_{+a_i})_1$ and $P^{\phi_i}_{+b_i} \sim (P^{\phi_i}_{+b_i})_1$  for $i=0,\dots,r$.
This sequence will be shown to have certain further properties that lead to a contradiction with Lemma~\ref{EAfinite}. 

We set $a_0:= a$, and $K:= K_0^{\phi_0}$.
If $v(P^{\phi_{0}},a_{0})\in \exc_K\big(L_{P^{\phi_0}_{+a_0}}\big)$, take $b_0=a_0$; otherwise, use Proposition~\ref{prop:quasilinear 1} to get $b_0\in \mathcal{O}_K=\mathcal{O}_{K_0}$ such that
$v(a_0-b_0)=v(P^{\phi_0}, a_0)$, $P(b_0) \prec P(a_0)$, and
$v(P^{\phi_0}, b_0) \in 
\exc_K\big(L_{P^{\phi_0}_{+b_0}}\big)$.
By Lemma~\ref{lem:summary} we have $\dval P^{\phi_0}_{+b_0}=\ddeg P^{\phi_0}_{+b_0}=1$, and so 
$P^{\phi_0}_{+b_0} \sim (P^{\phi_0}_{+b_0})_1$. 

To get $a_1$,
first take
$x\in K_1$ with $x'=\phi_0$, that is, $\der_0 x=1$ for the derivation $\der_0:=\phi_0^{-1}\der$ of $K_r^{\phi_0}$ that extends the derivation of $K$.  
Then $x\succ 1$, and $x$ is transcendental over~$K$ by Lemma~\ref{lem:integration simple}.  By Lemma~\ref{pre-extas2} and a subsequent remark, $K(x)$ is a $\d$-valued field extension of $K$ of $H$-type. Hence $K(x)$ is as described in Lemma~\ref{K(x)}. Applying Proposition~\ref{prop:quasilinear 2} to $F:=K_1^{\phi_0}$
with $P^{\phi_0}$ and $b_0$ in the roles of $P$ and $a$  yields an $a_1\in \mathcal{O}_{K_1}$ such that $P^{\phi_1}_{+a_1} \sim (P^{\phi_1}_{+a_1})_1$, $v(P^{\phi_{1}},a_{1})>v(P^{\phi_0},b_0)$, and every
$\alpha\in \exc_{K_0^{\phi_0}}\big(L_{P^{\phi_0}_{+b_0}}\big)$
with $0\le \alpha\le v(P^{\phi_0}, b_0)$ lies in 
$\exc_{K_1^{\phi_1}}\big(L_{P^{\phi_1}_{+a_1}}\big)$. 

With $K_1$,~$a_1$ as a new starting point instead
of $K_0$,~$a_0$, we repeat the above construction, 
and obtain by iteration elements $a_i,b_i\in\mathcal O_{K_i}$  ($i=0,\dots,r$)  
with the properties listed in the beginning of the proof,
such that moreover for $i=0,\dots,r$:
\begin{enumerate}
\item $v(P^{\phi_i},b_i)\ge v(P^{\phi_{i}},a_{i})$;
\item $v(P^{\phi_i},b_i)\in\exc_{K_i^{\phi_i}}\big(L_{P^{\phi_i}_{+b_i}}\big)$; 
\item $v(P^{\phi_{i+1}},a_{i+1})>v(P^{\phi_i},b_i)$, if $i<r$;
\item all
$\alpha\in \exc_{K_i^{\phi_i}}\big(L_{P^{\phi_i}_{+a_i}}\big)$
with $0\le \alpha\le v(P^{\phi_i}, a_i)$ lie in 
$\exc_{K_{i}^{\phi_{i}}}\big(L_{ P^{\phi_{i}}_{+b_{i}} }\big)$;
\item all
$\alpha\in \exc_{K_i^{\phi_i}}\big(L_{P^{\phi_i}_{+b_i}}\big)$
with $0\le \alpha\le v(P^{\phi_i}, b_i)$ lie in 
$\exc_{K_{i+1}^{\phi_{i+1}}}\big(L_{ P^{\phi_{i+1}}_{+a_{i+1}} }\big)$, if $i<r$.
\end{enumerate}
It follows that 
$v(P^{\phi_0},b_0)<v(P^{\phi_1},b_1)<\cdots<v(P^{\phi_r},b_r)$ are $r+1$ distinct elements of
$\exc_{K_r^{\phi_r}}\big(L_{P^{\phi_{r}}_{+b_{r}}}\big)$, contradicting Lemma~\ref{EAfinite}. 
\end{proof}

\noindent
The proof shows that the conclusion of Theorem~\ref{cor:quasilinear 2} can be strengthened to: 
\begin{quote}
 {\em $P(b)=0$ for some $b\in \mathcal{O}_{K_r}$ with $v(a-b)>0$.}\/
\end{quote}

\begin{example}
Let $K_n$ be the spherically complete $H$-subfield $\R[[\ell_0^{\Z}\cdots \ell_{n}^{\Z}]] $ of~$\T_{\log}$. Then
$\max\Psi_{K_n} = v(\ell_n^\dagger)$, and $\phi_n:=\ell_n^\dagger= \frac{1}{\ell_0\ell_1\cdots\ell_{n}} = \ell_{n+1}'\in\der K_{n+1}$.
Therefore Theorem~\ref{cor:quasilinear 2} applies to each tower
$K_n\subseteq K_{n+1}\subseteq\cdots\subseteq K_{n+r}$: if $P\in K_n\{Y\}^{\ne}$ has order $\leq r$ and $D_{P^{\phi_n}}$ is homogeneous of degree $1$, then $P$ has a zero in $\smallo_{K_{n+r}}$. All this goes through if $\ell_n^{\Z}$ is replaced by $\ell_n^{\R}$ for each $n$.
\end{example}

\begin{proof}[Proof of Theorem~\ref{thm:newtonianity of directed union}] Let $K$ be $\d$-valued such that $\der K = K$, and $K$ is a directed union of spherically complete grounded $\d$-valued subfields.
Let $P\in K\{Y\}^{\neq}$ be quasilinear, $r:=\order(P)$; we need to show that $P$ has a zero in $\mathcal O$.
By Corollary~\ref{recam}, $K$ is $\upo$-free. Hence by
Corollary~\ref{huposhape} and after replacing $K$,~$P$ by
suitable compositional conjugates 
we can assume that $K$ has small derivation, and that one of the following holds:
\begin{enumerate}
\item[(1)] for all active $\phi\preceq 1$ in $K$
there is an $f\in K^\times$ such that
$P^\phi\sim f Y'$,
\item[(2)] for all active $\phi\preceq 1$ in $K$
there is an $f\in K^\times$ such that $P^\phi\sim f Y$,
 \item[(3)] for all active $\phi\preceq 1$ in $K$ there are $f,g\in K^\times$ with
$P^\phi\sim f+gY$ and $f \asymp g$. 
\end{enumerate}
Take spherically complete $\d$-valued subfields $K_0\subseteq K_1\subseteq\cdots\subseteq K_r$ of $K$ with elements
$\phi_i\in K_i^{\times}$
such that $P\in K_0\{Y\}$, $v\phi_i=\max \Psi_{K_i}$
\textup{(}$i=0,\dots,r$\textup{)}, and $\phi_i\in\der K_{i+1}$ for $i=0,\dots,r-1$.
Then $\phi_0\preceq 1$ and $\phi_0$ is active. It follows that
either $P^{\phi_0}\sim f Y'$ for some $f\in K_0^\times$, 
or $P^{\phi_0}\sim fY$ for some $f\in K_0^\times$, or
$P^{\phi_0}\sim f+gY$ for some $f,g\in K_0^\times$ with
$f \asymp g$. In the first two cases the hypothesis of
Theorem~\ref{cor:quasilinear 2} is satisfied for $a=0$, and
in the last case it is satisfied for $a=-f/g$. Thus $P(b)=0$
for some $b\in \mathcal{O}_{K_r}$.   
\end{proof}

\noindent
In concrete cases the hypothesis $K=\der K$ in the theorem can often be verified by means of 
Corollary~\ref{cor:K=derK+C}. Taking in the example above
$\ell_n^{\R}$ in place of $\ell_n^{\Z}$ we get 
$\T_{\log}$ as the directed union of the
spherically complete~$K_n$, with $K_n\subseteq \der K_{n+1}$ by 
Corollary~\ref{cor:K=derK+C}, and
therefore $\T_{\log}=\der \T_{\log}$. Thus $\T_{\log}$ is newtonian. 

For $\T$, use the increasing union representation 
$\T = \bigcup_{n} E_{2n}{\downarrow}{}_{n}$ from 
Appendix~\ref{app:trans}. We show there that each
$E_{2n}{\downarrow}{}_n$ is a spherically complete
grounded $H$-subfield of $\T$ with 
$\max \Psi_{E_{2n}{\downarrow}{}_n}=v(\ell_n^\dagger)$. 
Thus $\T$ is newtonian.

\subsection*{Notes and comments} Berarducci and Mantova~\cite{BM} construct a derivation on Conway's field {\bf No} of surreal numbers with $\omega'=1$ that makes it
a Liouville closed $H$-field with constant field
$\R$. They show moreover that this derivation is simplest possible in a certain sense. They ask
whether {\bf No} with this derivation is elementarily equivalent to the differential field~$\T$. This question has a positive answer, based on Theorem~\ref{thm:newtonianity of directed union} and the next chapter; see~\cite{ADHNo}.

%% file: mt-16.tex
\chapter{Quantifier Elimination}\label{ch:QE}

\setcounter{theorem}{0}

\noindent
We are now close to establishing the main result of this volume:
the theory~$T^{\operatorname{nl}}$ of $\upo$-free newtonian Liouville closed
$H$-fields eliminates quantifiers in a certain natural language. This theory has two completions: $T^{\operatorname{nl}}_{\operatorname{small}}$, whose models are the
models of~$T^{\operatorname{nl}}$ with small derivation, and
$T^{\operatorname{nl}}_{\operatorname{large}}$, in whose models the derivation is not small. One can move from
models of $T^{\operatorname{nl}}_{\operatorname{small}}$ to models of $T^{\operatorname{nl}}_{\operatorname{large}}$ by compositional
conjugation. These two ``sides'' of $T^{\operatorname{nl}}$ reflect in a way the gap phenomenon, and we do not wish to obscure this by restriction to the ``small derivation'' case.

This chapter does not depend on the previous one, where the 
$H$-field $\T$ is shown to be a model of $T^{\operatorname{nl}}_{\operatorname{small}}$, but our quantifier 
elimination should of course be viewed in light of that fact
about $\T$.
%in Chapter~\ref{ch:newtdirun} 
 
\nomenclature[Bet1]{$\mathcal L^\iota_{\Upl\Upo}$}{language of $\Upl\Upo$-fields}
\nomenclature[Bi0]{$T^{\operatorname{nl}}$}{theory of $\upo$-free newtonian Liouville closed $H$-fields in the language of ordered valued differential rings}
\nomenclature[Bi1]{$T^{\operatorname{nl}}_{\operatorname{small}}$, $T^{\operatorname{nl}}_{\operatorname{large}}$}{completions of $T^{\operatorname{nl}}$}
\nomenclature[Bi2]{$T^{\operatorname{nl},\iota}_{\Upl\Upo}$}{theory of $\upo$-free newtonian Liouville closed $H$-fields in the language $\mathcal L^\iota_{\Upl\Upo}$}

To state the results of the present chapter with full precision, we introduce some first-order languages. Recall from Section~\ref{sec:DCF} that $\mathcal L_\der=\{0,1,{-},{+},{\cdot},\der\}$ is
the language of differential rings. We augment it here with binary 
relation symbols~$\leq$ and~$\preceq$ to obtain the language
$$\mathcal L\ :=\ \{0,\ 1,\  +,\  -,\ \cdot\ ,\ \der,\ \le,\ \preceq\}$$
of ordered valued differential rings. Each ordered valued differential field is viewed as an $\mathcal L$-structure
in the natural way, interpreting $\le$ as the ordering and $\preceq$ as the dominance relation as suggested by these symbols. It is clear that the $\upo$-free
newtonian Liouville closed $H$-fields are exactly the models of an
$\mathcal{L}$-theory~$T^{\operatorname{nl}}$, but~$T^{\operatorname{nl}}$ does not eliminate
quantifiers in this language, as we shall see in Section~\ref{sec:lHLO}. To achieve quantifier
elimination
we consider
a certain extension~$T^{\operatorname{nl},\iota}_{\Upl\Upo}$ by definitions of~$T^{\operatorname{nl}}$, in a
language~$\mathcal{L}^{\iota}_{\Upl\Upo}$ that augments 
$\mathcal{L}$ by a new unary function symbol $\iota$ and new 
unary relation symbols~$\I$,~$\Upl$ and~$\Upo$. We obtain
defining axioms in~$T^{\operatorname{nl},\iota}_{\Upl\Upo}$ for these new symbols by
requiring that every model $K$ of $T^{\operatorname{nl}}$ expands uniquely to a model $K^{\iota}_{\Upl\Upo}$ of $T^{\operatorname{nl},\iota}_{\Upl\Upo}$ such that for
all~$a\in K$, 
\begin{align*}
\iota(a)\ &=\ a^{-1} \text{ if $a\ne 0$,}\qquad \iota(0)=0,\\
\I(a)\ &\Longleftrightarrow\ a=y' \text{ for some $y\prec 1$  in $K$,}\\
\Upl(a)\ &\Longleftrightarrow\ a=-y^{\dagger\dagger} \text{ for some $y\succ 1$  in $K$,}\\
\Upo(a)\  &\Longleftrightarrow\ 4y''+ay=0 \text{ for some $y\in K^\times$.}
\end{align*}
Note that $\I$, $\Upl$, and $\Upo$ get interpreted in the expansion 
$K^{\iota}_{\Upl\Upo}$ of a model $K$ of $T^{\operatorname{nl}}$ as the
convex additive subgroup $\I(K)$ of $K$ and as the downward closed subsets~$\Upl(K)$ and~$\Upo(K)$ of~$K$ that were introduced
in Sections~\ref{sec:special sets} and~\ref{sec:secondorder}.  
 We can now state our main result, proved in
Section~\ref{sec:embth}:

\begin{theorem}\label{thm:qe} The theory $T^{\operatorname{nl},\iota}_{\Upl\Upo}$ eliminates
quantifiers.                        
\end{theorem}

\noindent
This fails if either $\Upl$ or $\Upo$ is dropped from the language; see Section~\ref{sec:lHLO}.  On the other hand, the predicate $\I$ is only included for convenience, to simplify some later proofs and formulations: in Theorem~\ref{thm:qe} we can drop $\I$ from the language, since 
for $K$ as above and $a\in K^\times$ we have by Lemmas~\ref{lem:Upl, 1} and~\ref{IOlemma}: 
$$\I(a)\ \Longleftrightarrow\  \neg \Upl(-a^\dagger),\ \qquad\  \I(a)\ \Longleftrightarrow\ \Upo\big(\sigma(a)\big).$$
In Section~\ref{sec:embth} we also derive some consequences of Theorem~\ref{thm:qe}:

\begin{cor}\label{cor:qe} Let $K$ be an $\upo$-free newtonian Liouville closed $H$-field. Then: \begin{enumerate}
\item[\textup{(i)}] $K$ has $\operatorname{NIP}$;
\item[\textup{(ii)}] a subset of $C^n$ is definable in $K$ if and only if it is semialgebraic in the sense of the real closed field $C$;
\item[\textup{(iii)}] $K$ is o-minimal at infinity: if $X\subseteq K$ is definable in $K$,
then there exists $a\in K$ such that $(a,+\infty)\subseteq X$ or $(a,+\infty) \cap X = \emptyset$.
\end{enumerate}
\end{cor}

\noindent
We indicate briefly the ideas behind the proof of Theorem~\ref{thm:qe}, and some steps that still need to be taken.
As is well-known, a theory eliminates quantifiers if and only if its models and their substructures have certain
embedding properties. It is for this reason that we have built an arsenal of embedding and extension results in the previous chapters:
the universal property of the $H$-field hull
of a pre-$H$-field in Section~\ref{sec:H-fields}, the algebraic and Liouville extensions of Sections~\ref{sec:H-fields} and
~\ref{sec:Liouville closed}, the construction of immediate extensions in Section~\ref{sec:construct imm exts} 
and of $F_{\upo}$ in Section~\ref{sec:behupo}, Proposition~\ref{Evalcor}
and Corollary~\ref{Evalocor}, and, crucially, the Newton-Liouville closure of Section~\ref{sec:newtonization} (the latter requiring a lengthy development over several chapters). An important consequence of 
Theorem~\ref{newtmax} is that it 
enables us to deal with immediate extensions. It is worth noting that immediate extensions
typically require the most attention in analogous situations---going back to the work by Ax-Kochen~\cite{AxKochen} and Er{\accentv{s}}ov~\cite{Ersov} in the 1960s---and such extensions indeed preoccupied us in several chapters. In contrast to prior model-theoretic work on (enriched) valued fields, however, the models of our theory $T^{\operatorname{nl}}$ are never spherically complete:
both Liouville closedness and $\upo$-freeness prevent that.

Constant field extensions are taken care of
by Propositions~\ref{prop:constant field ext} and~\ref{prop:Hconstants}. It remains to deal with extensions 
that are
completely controlled by the corresponding extension
of asymptotic couples. We handle such extensions in Section~\ref{sec:valgrp}, 
where it
leads to a result of independent interest:

%The case where the value group gets bigger
%in an extension still requires work. Fortunately, we can %reduce to a situation
%where such an extension is controlled by the corresponding %extension
%of asymptotic couples. The relevant extensions of asymptotic %couples 
%were studied in \cite{AvdD} of which we 
%quote in Section~\ref{sec:valgrp} the key fact we need. 
%This leads in that section to the following result of %independent interest:

\begin{theorem}\label{th:noextL} If $K$ is an $\upo$-free newtonian Liouville closed $H$-field, then $K$ has no proper $\d$-algebraic $H$-field extension with the same constant field. 
\end{theorem}

\noindent
In Section~\ref{mctnl} we use this to get uniqueness-up-to-isomorphism of Newton-Liouville closures of
$\upo$-free $H$-fields, and then employ our arsenal 
to prove an embedding result that has the model completeness of $T^{\operatorname{nl}}$ as a consequence.

\medskip\noindent
Going beyond model completeness to quantifier elimination
requires attention to the substructures of models of~$T^{\operatorname{nl},\iota}_{\Upl\Upo}$ rather than of $T^{\operatorname{nl}}$, and this involves the extra predicates 
$\I$,~$\Upl$,~$\Upo$. 
Therefore we first determine
in Section~\ref{sec:LO-cuts} the substructures of models of~$T^{\operatorname{nl},\iota}_{\Upl\Upo}$:
they are the expanded pre-$H$-fields $(K, I, \Lambda, \Omega)$
where $(I, \Lambda, \Omega)$ is a $\Upl\Upo$-cut
in the pre-$H$-field $K$ as defined in that section. It turns out that  a given pre-$H$-field~$K$ has either exactly one $\Upl\Upo$-cut or exactly two $\Upl\Upo$-cuts. Moreover,
an $\upo$-free pre-$H$-field~$K$ has just one
$\Upl\Upo$-cut, and any $(K, I, \Lambda, \Omega)$ 
has an extension $(K^*, I^*, \Lambda^*, \Omega^*)$ with~$K^*$
an $\upo$-free~$H$-field, such that $(K^*, I^*, \Lambda^*, \Omega^*)$
embeds over~$K$ into any model of~$T^{\operatorname{nl},\iota}_{\Upl\Upo}$ extending $(K, I, \Lambda, \Omega)$: Proposition~\ref{prop:upo ext of HLO}. 
These two facts allow us to focus henceforth for embedding purposes on 
$\upo$-free $H$-fields, and forget about $\Upl\Upo$-cuts,
and this is taken care of by the results in Section~\ref{mctnl}.

\subsection*{Notes and comments} The first suggestion that the primitives $\Upl(\T)$
and $\Upo(\T)$ might be enough to eliminate quantifiers for
$\T$, in addition to the usual primitives for ordered valued differential fields, is in \cite{ADHOleron}. 
The very end of that paper indicates a possible further obstruction to QE that would require also a certain partial inverse of the function $\omega$ as an extra primitive. Fortunately, we were able to get rid of this obstruction: this ultimately rests on observing Corollary~\ref{cor:companion of Stirling, variant}, which gives Corollary~\ref{cor:nabla1, order 2} needed in Section~\ref{fcuplfr}, which in turn is used in the proofs of Lemmas~\ref{ILO6} and~\ref{lem:embed into omega-free, 3} below. 

A minor issue is whether to include field inversion among the primitives to get
QE. We found it convenient to do so. If we drop the symbol $\iota$
for inversion from the language, we would probably loose QE, but we would
certainly regain it upon replacing the unary relation symbols $\Upl$ and $\Upo$ 
by binary relation
symbols $\Upl_2$ and $\Upo_2$, to be interpreted in $\T$ according to
$$\Upl_2(a,b)\Leftrightarrow a\in b\cdot\Upl(\T), \qquad \Upo_2(a,b)\Leftrightarrow a\in b\cdot\Upo(\T).$$

\section{Extensions Controlled by Asymptotic Couples}\label{sec:valgrp}

\noindent
In this section we deal with a kind of $H$-field extension
that is in some sense controlled by the corresponding extension of 
asymptotic couples.  At this point
%The relevant extensions of asymptotic couples were studied already in detail in \cite{AvdD}. We only need here one key fact from that paper, namely
Proposition~\ref{cac} about closed $H$-asymptotic couples becomes relevant, and we use it to obtain 
the following:

\begin{lemma}\label{uponewlidim} Let $K$ be an $\upo$-free newtonian Liouville closed $H$-field and $L$ an $H$-field extension with $C_L=C$, and let $f\in L\setminus K$. Suppose $K$ is maximal in $L$ in the sense that there is no $y\in L\setminus K$ for which $K\<y\>$ is an immediate extension of~$K$. Then the vector space $\Q\Gamma_{K\<f\>}/\Gamma$ over $\Q$ is
infinite-dimensional.
\end{lemma}
\begin{proof} We claim that there is no divergent pc-sequence in $K$
with a pseudolimit in~$L$. To see this, let $(y_{\rho})$ be a divergent pc-sequence in $K$. It cannot be of $\d$-algebraic type, since $K$ is asymptotically $\d$-algebraically maximal. So it is of $\d$-transcendental type, and if it had a pseudolimit $y\in L$, then
$K\<y\>$ would be an immediate extension of $K$. This proves our claim.
Thus for each $y\in L\setminus K$ the set
$v(y-K)\subseteq \Gamma_L$ has a largest element: otherwise there would be a divergent pc-sequence in $K$ with pseudolimit~$y$. Given $y\in L\setminus K$, a {\em best
approximation in $K$ to $y$\/} is by definition an element $y_0\in K$ such that $v(y-y_0)=\max v(y-K)$; note that then 
$v(y-y_0)\notin \Gamma$, since $C_L=C$. For convenience we set
$L=K\<f\>$ below.

Pick a best approximation $b_0$ in $K$ to $f_0:=f$, and set
$f_1:= (f_0-b_0)^\dagger\in L$. Then $f_1\notin K$, since
$K$ is Liouville closed and $C=C_L$. Thus we can take a best approximation $b_1$ in $K$ to $f_1$, and continuing this way, we
obtain a sequence $(f_n)$ in~${L\setminus K}$ and a sequence
$(b_n)$ in $K$, such that $b_n$ is a best approximation in $K$ to~$f_n$ and $f_{n+1}=(f_n-b_n)^\dagger$ for all $n$. Thus
$v(f_n-b_n)\notin \Gamma$ for all $n$.

\claim{$v(f_0-b_0),v(f_1-b_1), v(f_2-b_2),\dots$ are $\Q$-linearly independent over $\Gamma$.}

\noindent
To prove this claim, take $a_n\in K^\times$ with $a_n^\dagger=b_n$ for
$n\ge 1$. Then
$$f_n-b_n\ =\ (f_{n-1}-b_{n-1})^\dagger-a_n^\dagger\ =\ \left(\frac{f_{n-1}-b_{n-1}}{a_n}\right)^\dagger \qquad (n\ge 1).$$
With $\psi:= \psi_L$ and $\alpha_n=va_n\in \Gamma$ for $n\ge 1$, we get 
\begin{align*} v(f_n-b_n)\ &=\ \psi\big(v(f_{n-1}-b_{n-1})-\alpha_n\big),\
 \text{ so by an easy induction on $n$,}\\
 v(f_n-b_n)\ &=\ \psi_{\alpha_1,\dots, \alpha_n}\big(v(f_0-b_0)\big),  \qquad (n\ge 1).
 \end{align*} 
Suppose towards a contradiction that $v(f_0-b_0), \dots, v(f_n-b_n)$ are $\Q$-linearly dependent
over $\Gamma$. Then we have $m<n$ and $q_1,\dots, q_{n-m}\in \Q$ such that
$$v(f_m-b_m) + q_1v(f_{m+1}-b_{m+1}) + \cdots + q_{n-m}v(f_n-b_n)\in \Gamma.$$ 
For $\gamma:= v(f_m-b_m)\in \Gamma_L\setminus \Gamma$ this gives
$$\gamma + q_1\psi_{\alpha_{m+1}}(\gamma) + \cdots + q_{n-m}\psi_{\alpha_{m+1},\dots,\alpha_n}(\gamma)\in \Gamma,$$
but this contradicts Proposition~\ref{cac}. 
\end{proof}

{\sloppy

\begin{proof}[Proof of Theorem~\ref{th:noextL}] Let $K$ be an $\upo$-free newtonian Liouville closed $H$-field, and assume towards a 
contradiction that $L$ is a proper $\d$-algebraic $H$-field extension of~$K$ with $C_L=C$. Since $K$ is asymptotically $\d$-algebraically maximal,
there is no $y\in L\setminus K$ for which $K\<y\>$ is an immediate extension of $K$. Taking any~$f$ in $L\setminus K$, the transcendence degree
of $K\<f\>$ over $K$ is finite, but this contradicts Lemma~\ref{uponewlidim} in view of the Zariski-Abhyankar Inequality (Corollary~\ref{cor:ZA inequality}). 
\end{proof} 

}

\subsection*{Description of $K\<f\>$} Let $K$, $L$, and $f$ be as in 
Lemma~\ref{uponewlidim}. Elaborating on the proof of this lemma we shall
obtain a complete description of $K\<f\>$ as an $H$-field extension
of $K$ generated by $f$. For this we use the notations in that proof, and set $\beta_n:= v(f_{n}-b_{n})-\alpha_{n+1}\in \Gamma_{K\<f\>}$.
Thus $\beta_0, \beta_1, \beta_2,\dots$ are $\Q$-linearly independent over $\Gamma$, by the proof of
Lemma~\ref{uponewlidim}. 
%we have $\beta_n\notin \Gamma$ for all $n$ 

\begin{lemma}\label{type5} The asymptotic couple of $K\<f\>$ has the following properties:  \begin{enumerate}
\item[\textup{(i)}] $\Gamma_{K\<f\>}\ =\ \Gamma \oplus \bigoplus_n \Z \beta_n$ \textup{(}internal direct sum\textup{)};
\item[\textup{(ii)}] $\beta_n^\dagger\notin \Gamma$ for all $n$, and $\beta_m^\dagger\ne \beta_n^\dagger$ for all $m\ne n$;
\item[\textup{(iii)}]  $\psi\big(\Gamma_{K\<f\>}^{\ne}\big)\ =\ \Psi \cup \{\beta_n^\dagger:\ n=0,1,2,\dots\}$;
\item[\textup{(iv)}] $\big[\Gamma_{K\<f\>}\big]\ =\ 
              [\Gamma]\cup \big\{[\beta_n]: n=0,1,2,\dots\big\}$;
\item[\textup{(v)}] $\Gamma^{<}$ is cofinal in $\Gamma_{K\<f\>}^{<}$, and 
 $\beta_0^\dagger\ <\ \beta_1^\dagger\ <\ \beta_2^\dagger\ <\ \cdots.$
\end{enumerate}
\end{lemma}
\begin{proof} Consider the ``monomials''
$\fm_n:= \frac{f_n-b_n}{a_{n+1}}$ with $v(\fm_n)=\beta_n$. Then
\begin{align*}\fm_{n+1}\ &=\ \frac{f_{n+1}-b_{n+1}}{a_{n+2}}\ =\ \frac{(f_n-b_n)^\dagger - b_{n+1}}{a_{n+2}}\\ 
&=\ \frac{(a_{n+1}\fm_n)^\dagger - b_{n+1}}{a_{n+2}}\
=\ \frac{a_{n+1}^\dagger + \fm_n^\dagger - b_{n+1}}{a_{n+2}}\ =\ \frac{\fm_n^\dagger}{a_{n+2}},
\end{align*}
and so $\fm_n'=a_{n+2}\fm_n\fm_{n+1}$. Thus $f=b_0 + a_1\fm_0$ gives $f'=b_0' + a_1'\fm_0 + a_1a_2\fm_0\fm_1$,
and continuing by induction on $n$ gives
$$f^{(n)}\ =\ F_n(\fm_0,\dots, \fm_n), \quad F_n(Y_0,\dots, Y_n)\in K[Y_0,\dots, Y_n],\quad \deg F_n \le n+1.$$
Thus for $P\in K\{Y\}^{\ne}$ of order $\le r\in \N$ we have
$$ P(f)\ =\ \sum_{\i\in I} a_{\i}\, \fm_0^{i_0}\cdots \fm_r^{i_r}$$
where the sum is over a finite nonempty set $I$ of tuples $\i=(i_0,\dots,i_r)\in \N^{1+r}$, and $a_{\i}\in K^\times$ for
all $\i\in I$. Since $v(\fm_0)=\beta_0,v(\fm_1)=\beta_1,\dots$ are 
$\Q$-linearly independent over $\Gamma$, we obtain 
$v\big(P(f)\big)\in \Gamma+\sum_n \N\beta_n$, which proves (i).

\medskip\noindent
We have $\beta_n^\dagger\notin \Gamma$ because by the
proof of Lemma~\ref{uponewlidim},
$$\beta_n^\dagger\ =\ \psi\big(v(f_{n}-b_{n})-\alpha_{n+1}\big)\ =\  v(f_{n+1}-b_{n+1})\ =\ \beta_{n+1}+\alpha_{n+2}\notin \Gamma.$$
Since $\beta_1, \beta_2,\beta_3,\beta_4,\dots$ are $\Q$-linearly independent, so are $\beta_0^\dagger, \beta_1^\dagger, \beta_2^\dagger,\beta_3^\dagger,\dots$ by 
these equalities. This proves (ii),
which in view of (i) yields (iii).
From (ii) we get $[\beta_n]\notin [\Gamma]$ for all $n$, and
$[\beta_m]\ne [\beta_n]$ for all $m\ne n$. Again by (i), this
gives (iv). 

\medskip\noindent
To get (v), assume towards a contradiction that 
$\Gamma^{<}$ is not cofinal in $\Gamma_{K\<f\>}^{<}$. Then
by (iv) we get $n$ with $[\beta_n]< [\alpha]$ for all $\alpha\in \Gamma^{\ne}$, hence $\Psi < \beta_n^\dagger < (\Gamma^{>})'$.
Then $[\beta_n^\dagger-\alpha]\in [\Gamma]$ for all $\alpha\in \Gamma$, by Corollary~\ref{addgapmore}. For $\alpha:= \alpha_{n+2}$ this means
$[\beta_n^\dagger-\alpha_{n+2}]=[\beta_{n+1}]\in [\Gamma]$, contradicting (ii). Thus $\Gamma^{<}$ is indeed
cofinal in $\Gamma_{K\<f\>}^{<}$. For any $n$ we can therefore take
$\alpha\in \Gamma^{\ne}$ with $[\alpha]<[\beta_n]$.
Also $[\beta_{n+1}]\notin [\Gamma]$ and $\beta_n^\dagger-\alpha^\dagger\in (\Gamma+\Z \beta_{n+1})\setminus \Gamma$, and thus
by Lemmas~\ref{archextclass} and~\ref{BasicProperties}, 
$$[\beta_{n+1}]\ \le\ [\beta_n^\dagger-\alpha^\dagger]\ <\ [\beta_n-\alpha]\ =\ [\beta_n].$$
So we have a strictly decreasing sequence $[\beta_0]>[\beta_1]>[\beta_2]>\cdots$ in $[\Gamma_{K\<f\>}]$, and thus a strictly increasing sequence $\beta_0^\dagger < \beta_1^\dagger < \beta_2^\dagger < \cdots$ in view of (ii). 
\end{proof}

\noindent
The following consequence is not needed later but worth pointing out:

\begin{cor} $K\<f\>$ is $\upo$-free.
\end{cor}
\begin{proof} 
Assume towards a contradiction that $K\<f\>$ is not $\upo$-free.
Since $\Gamma^{<}$ is cofinal in $\Gamma_{K\<f\>}^{<}$ this gives
an element $\upo\in K\<f\>$ such that $\upo_{\rho}\leadsto \upo$,
where $(\upo_{\rho})$ is the sequence in $K$ obtained in the usual way
from a logarithmic sequence in $K$. Now $(\upo_{\rho})$ is
of $\d$-transcendental type over $K$, so $K\<\upo\>\subseteq K\<f\>$ 
is an immediate 
extension of $K$. Since $\upo\notin K$, this contradicts Lemma~\ref{uponewlidim}.  
\end{proof}

\begin{lemma}\label{fgcorrespondence} Suppose $g$ in some $H$-field extension $M$ of $K$ realizes the same cut in the ordered set $K$ as $f$ does. Then $v(g-b_0)=\max v(g-K)\notin \Gamma$, and $g_1:=(g-b_0)^\dagger$ realizes the same cut in the ordered set $K$ as $f_1=(f-b_0)^\dagger$.
\end{lemma}
\begin{proof} Let $\alpha\in \Gamma$ and $b\in K$. We claim that 
$$v(f-b) <\alpha\ \Longleftrightarrow\ v(g-b) < \alpha, \quad v(f-b) >\alpha\ \Longleftrightarrow\ v(g-b) > \alpha.$$
To prove this claim, take 
$a\in K^{>}$ with $va=\alpha$.
%assume for simplicity that $f> b$, and 
Consider first the case that $v(f-b) < \alpha$. Then $|f-b| > a$, so 
$|g-b| > a$,
and thus $v(g-b)\le \alpha$. 
Since $\Gamma^{<}$ is cofinal in
$\Gamma_{K\<f\>}^{<}$ we have $\alpha_1\in \Gamma$ such that
$v(f-b) < \alpha_1 < \alpha$. Take 
$a_1\in K^{>}$ such that $va_1=\alpha_1$. Then the argument above
with $a_1$ instead of $a$ gives $v(g-b)\le \alpha_1 < \alpha$.
In a similar way we show: $v(f-b)> \alpha\Rightarrow v(g-b)>\alpha$.
Finally, assume $v(f-b)=\alpha$. Then $C=C_{K\<f\>}$ gives $c\in C^{\times}$ with $f-b\sim ca$, so $\frac{|c|}{2}a < |f-b| < 2|c|a$,
hence $\frac{|c|}{2}a < |g-b| < 2|c|a$, and thus $v(g-b)=va=\alpha$.
This finishes the proof of the claim.

It follows from this claim and $v(f-b_0)\notin \Gamma$ that 
$v(g-b_0)\notin \Gamma$. This gives $v(g-b_0)=\max v(g-K)$: otherwise
we have $b\in K$ with $v(g-b_0)< v(g-b)$, and so
$v(g-b_0)=v(b-b_0)\in \Gamma$, a contradiction. Next we get
$(g-b_0)^\dagger\notin K$: otherwise, $(g-b_0)^\dagger=a^\dagger$
with $a\in K^\times$, so $g-b_0=ca$ for some $c\in C_M^\times$, and thus $v(g-b_0)=va\in \Gamma$, a contradiction. 

Next, we show that $(g-b_0)^\dagger$ realizes the same cut in
$K$ as $(f-b_0)^\dagger$. If $f< b_0$, then we can replace $f$,~$g$,~$b_0$
by $-f$,~$-g$,~$-b_0$ to reduce to the case $f> b_0$. So we assume
$f>b_0$, which gives $g>b_0$. Suppose towards a contradiction that 
$h\in K$ is such that $(f-b_0)^\dagger < h$ in $K\<f\>$ and 
$h< (g-b_0)^\dagger$ in $M$. Take $\phi\in K^{>}$ such that
$h=\phi^\dagger$. Then $s:= (f-b_0)/\phi >0$ and $s^\dagger<0$, so
$s'<0$, and thus $s=c+\epsilon$ with $c\in C^{\ge}$ and $0<\epsilon\prec 1$ in $K\<f\>$. If $c\ne 0$,
then $v(f-b_0)=v\phi\in \Gamma$ contradicts $v(f-b_0)\notin \Gamma$. So $c=0$ and thus $f=b_0 + \phi\epsilon$.  
Likewise, $h< (g-b_0)^\dagger$ gives $t:=(g-b_0)/\phi>0$ and
$t^\dagger > 0$, so $t'>0$, and thus
either $t=c^*-\epsilon^*$ with $c^*\in C_{M}^{>}$ and $0<\epsilon^*\prec 1$ in~$M$, or $t>C_{M}$.
The first case would give $v(g-b_0)=v\phi$, a contradiction, so
we get $t>C_{M}$, but that would give 
$$f\ =\ b_0+\phi\epsilon\ <\  b_0+\phi\ <\  b_0+\phi t\ =\ g,$$ 
with $b_0+\phi\in K$, contradicting that
$f$ and $g$ realize the same cut in $K$. 

The other way that $(g-b_0)^\dagger$ does not realize the same cut in $K$ as $(f-b_0)^\dagger$ is that we have $h\in K$ such that
$(g-b_0)^\dagger < h$ in $M$ and $h < (f-b_0)^\dagger$ in $K\<f\>$. Taking as before $\phi\in K^{>}$ with
$h=\phi^\dagger$, a similar argument gives us $g < b_0 + \phi < f$, and we have again a contradiction.
\end{proof}

\begin{prop}\label{cutcorrespondence} Suppose $g$ in some $H$-field extension 
$M$ of $K$ realizes the same cut in the ordered set $K$ as $f$ does. Then there is
an embedding $K\<f\> \to M$ of $H$-fields over $K$ sending $f$ to $g$.
\end{prop}
\begin{proof} By Lemma~\ref{fgcorrespondence} we can recursively define $g_n\in M\setminus K$ to realize the same cut
in $K$ as $f_n$ by $g_0:= g$, and $g_{n+1}:= (g_n-b_n)^\dagger$. 
Then $v(g_n-b_n)\notin \Gamma$ for all $n$. The same argument
as in the proof of Lemma~\ref{uponewlidim} shows: 
$$v(g_0-b_0),\ v(g_1-b_1),\ v(g_2-b_2),\dots \text{ are $\Q$-linearly independent over $\Gamma$}.$$  
Set $\beta_{n}^*:= v(g_n-b_n)-\alpha_{n+1}$, and
$\fm_n^*:= \frac{g_n-b_n}{a_{n+1}}$. Then $\beta_0^*, \beta_1^*, \beta_2^*,\dots$ are $\Q$-linearly independent over $\Gamma$. 
With $F_n\in K[Y_0,\dots, Y_n]$ as in the proof of Lemma~\ref{type5} we get $g^{(n)}=F_n(\fm_0^*,\dots, \fm_n^*)$, and so for 
$P\in K\{Y\}^{\ne}$ of order $\le r\in \N$ we get
$$ P(g)\ =\ 
\sum_{\i\in I} a_{\i}\, {\fm_0^*}^{i_0}\cdots {\fm_r^*}^{i_r}$$
where the sum is over the same finite nonempty set $I$ of tuples $\i=(i_0,\dots,i_r)\in \N^{1+r}$ as in the proof of Lemma~\ref{type5}, and
with the same coefficients $a_{\i}\in K^\times$ as in that proof. 
Since $v(\fm^*_0)=\beta^*_0,v(\fm^*_1)=\beta^*_1,\dots$ are 
$\Q$-linearly independent over $\Gamma$, we obtain 
$v\big(P(g)\big)\in \Gamma+\sum_n \N\beta^*_n$. The rest of that proof then shows that Lemma~\ref{type5} goes through with $f$ replaced by $g$ and each $\beta_n$ by $\beta^*_n$. In particular, 
$$[\beta^*_0]\ >\ [\beta^*_1]\ >\ [\beta^*_2]\ >\ \cdots.$$
Next, $\fm_n\in K\<f\>$ realizes the same cut in $K$ as $\fm_n^*$,
and so $\beta_n$ realizes the same cut in the ordered set
$\Gamma$ as $\beta^*_n$. It follows that we have an ordered abelian group isomorphism $j\colon \Gamma_{K\<f\>} \to \Gamma_{K\<g\>}$
over $\Gamma$ sending $\beta_n$ to $\beta^*_n$ for each $n$.
Using the expressions above for $P(f)$ and $P(g)$ it follows that
$j\big(v(P(f))\big)\ =\ v(P(g))$ for all $P\in K\{Y\}^{\ne}$, so
we have a valued differential field embedding $K\<f\> \to M$
over~$K$ sending $f$ to $g$. This is even an embedding of ordered valued differential fields, since the facts stated about the
$\fm_n$, $\beta_n$ and $\fm^*_n$, $\beta^*_n$ yield: $P(f) > 0\ \Longleftrightarrow\ P(g)>0$, for all $P\in K\{Y\}^{\ne}$. 
\end{proof}

\subsection*{Notes and comments} Lemma~\ref{type5} is a version for 
$H$-fields of Proposition~5.3 in \cite{AvdD}. That proposition analyzes ``simple extensions of $H$-triples of type (V).'' 

In applying Theorem~\ref{th:noextL} one should watch out: replacing {\em $H$-field extension\/} in its statement by {\em pre-$H$-field extension} results in a false statement.

\section{Model Completeness}\label{mctnl}

\noindent
We first show the uniqueness-up-to-isomorphism of
Newton-Liouville closures of $\upo$-free $H$-fields, although
existence of such closures is enough for model completeness.

\subsection*{Uniqueness of Newton-Liouville closure} Let $E$ be an 
$\upo$-free $H$-field. We saw in Section~\ref{sec:newtonization} that $E$ has a Newton-Liouville closure. 

\begin{lemma}\label{embnl} Let $E^{\operatorname{nl}}$ be any Newton-Liouville closure of $E$ and  $i\colon E^{\operatorname{nl}} \to L$ an embedding into an $H$-field $L$ with $C_L\subseteq i(E^{\operatorname{nl}})$. Then 
$$i(E^{\operatorname{nl}})\ =\ \big\{f\in L:\  \text{$f$ is $\d$-algebraic over $i(E)$}\big\}.$$
\end{lemma}

{\sloppy

\begin{proof} Since $E^{\operatorname{nl}}$ is $\d$-algebraic over 
$E$,
every element of $i(E^{\operatorname{nl}})$ is $\d$-algebraic over~$i(E)$. Also, $i(E^{\operatorname{nl}})$ is a newtonian Liouville
closed $H$-subfield of $L$ with the same constants as~$L$,
so every $f\in L$ that is $\d$-algebraic over 
$i(E^{\operatorname{nl}})$ lies in $i(E^{\operatorname{nl}})$ by
Theorem~\ref{th:noextL}.
\end{proof}

}

\begin{cor}\label{uninl} Any two Newton-Liouville closures of $E$ are isomorphic
over~$E$. If $E^{\operatorname{nl}}$ is a Newton-Liouville closure
of $E$, then $E^{\operatorname{nl}}$ does not have any proper newtonian
Liouville closed $H$-subfield containing $E$.
\end{cor}
\begin{proof} Let $E^{\operatorname{nl}}$ and $L$ be Newton-Liouville closures of $E$. Then there exists an
embedding $E^{\operatorname{nl}}\to L$ over $E$, and any such embedding is necessarily surjective by Lemma~\ref{embnl}. This proves
the first part. The minimality property of $E^{\operatorname{nl}}$
also follows from Lemma~\ref{embnl} by considering embeddings 
$E^{\operatorname{nl}}\to E^{\operatorname{nl}}$ over $E$. 
\end{proof}

\subsection*{Model completeness of $T^{\operatorname{nl}}$}
As usual, $|S|$ denotes the cardinality of a set $S$, and $\kappa^+$ the next cardinal after the cardinal~$\kappa$. Here is the decisive embedding property:

\begin{prop}\label{upoembnl} Let $E$ be an $\upo$-free $H$-subfield of
an $\upo$-free newtonian Liouville closed $H$-field $K$ with $C_E=C$, let
$i\colon E \to L$ be an embedding into an $\upo$-free newtonian Liouville closed $H$-field $L$. Assume $\operatorname{cf}(\Gamma_L^{<})>|\Gamma|$ and the 
underlying ordered set of $L$ is $|K|^{+}$-saturated. Then $i$ extends to an 
embedding $K\to L$.
\end{prop}
\begin{proof} Note that every differential subfield of $K$ containing $E$ is an $H$-subfield of~$K$. Assume $E\ne K$; it is enough to show that then $i$ can be extended to an embedding $F\to L$ for some
$\upo$-free $H$-subfield $F$ of $K$ properly containing $E$.   

Consider first the case that $\Gamma_E^{<}$ is not cofinal in $\Gamma^{<}$. Then we have $y\in K^{>}$ such that 
$\Gamma_E^{<} < vy < 0$. By the cofinality assumption on $\Gamma_L^{<}$ 
we also have $y^*\in L^{>}$ such that 
$\Gamma_{iE}^{<} < vy^* < 0$. 
Then Corollary~\ref{Evalocor} yields an extension of $i$ to
an embedding $E\<y\> \to L $ sending $y$ to $y^*$.  Now
$E\<y\>$ is a grounded $H$-field by
Proposition~\ref{Evalcor}. Then by Lemma~\ref{Hlogupo} we can extend
this embedding $E\<y\>\to L$ to an embedding $F:=E\<y\>_{\upo}\to L$, and for the same reason we can identify the extension $F$ of $E\<y\>$ over $E\<y\>$ with an
$\upo$-free $H$-subfield of $K$. This achieves our goal of extending 
$i$ to an embedding $F\to L$.
      
We are left with the case that $\Gamma_E^{<}$ is cofinal in 
$\Gamma^{<}$. Then every differential subfield of $K$ containing $E$ is an $\upo$-free $H$-subfield of $K$.

\subcase[1]{$E$ is not a newtonian Liouville closed $H$-field.} Then we can extend $i$ to an embedding $F\to L$ where 
$$F\ :=\  \{f\in K:\ f \text{ is $\d$-algebraic over $E$}\}$$
is the Newton-Liouville closure of $E$ inside $K$, by Lemma~\ref{embnl}. 

\subcase[2]{$E$ is newtonian and Liouville closed, and $E\<y\>$ is an immediate extension of $E$ for some $y\in K\setminus E$.} For such $y$ we take a divergent pc-sequence~$(a_{\rho})$ in~$E$
such that $a_{\rho}\leadsto y$. Since $E$ is asymptotically $\d$-algebraically maximal, $(a_{\rho})$ is of $\d$-transcendental type over $E$.
The saturation assumption on $L$ gives $z\in L$ such that
$i(a_{\rho})\leadsto z$, by Lemma~\ref{pcordgp}. Then
Section~\ref{sec:construct imm exts} yields a valued differential field embedding $F:= E\<y\>\to L$ that extends $i$ and sends 
$y$ to $z$. By Lemma~\ref{prehimm} this embedding $F\to L$ is an embedding of ordered valued differential fields.

\subcase[3]{$E$ is newtonian and Liouville closed, and there is no $y\in K\setminus E$ such that $E\<y\>$ is an immediate extension of $E$.} Then we take any $f\in K\setminus E$, and take some
$g\in L$ such that for all $a\in E$ we have: $a<f\ \Rightarrow\ i(a)< g$, and $a>f\ \Rightarrow\ i(a)>g$.  Then Proposition~\ref{cutcorrespondence} yields an extension of $i$ to an
embedding $E\<f\> \to L$ sending~$f$ to~$g$.
\end{proof}

\noindent
When $C_E\ne C$, we require an extra hypothesis:

\begin{cor}\label{cor:upoembnl} Let $E$, $K$, $L$, and $i$ be as in Proposition~\ref{upoembnl}, except that we drop the assumption $C_E=C$. Assume in addition that the underlying ordered set of~$C_L$ is $|C|^+$-saturated. Then $i$ extends to an embedding $K \to L$. 
\end{cor}
\begin{proof} The real closed constant field $C_L$ is 
$|C|^{+}$-saturated, so the ordered field embedding $i|C_E\colon C_E \to C_L$ extends to an ordered field embedding $j\colon C \to C_L$.
Then Propositions~\ref{prop:constant field ext} and~\ref{prop:Hconstants} yield an extension of~$i$ to an embedding $E(C) \to L$ that agrees with $j$ on $C$.
The $H$-subfield $E(C)$ of~$K$ is 
$\d$-algebraic over $E$, so is $\upo$-free. Now
 apply Proposition~\ref{upoembnl} with $E(C)$ in place of $E$.
\end{proof}

\noindent 
The cofinality and saturation hypotheses in Proposition~\ref{upoembnl} and 
Corollary~\ref{cor:upoembnl} are of course satisfied if $L$ 
(as an ordered valued differential field) is 
$|K|^+$-saturated, but the weaker
assumption of Proposition~\ref{upoembnl} is useful in \cite{ADHNo}.  
In view of the next result we recall from
Chapter~\ref{ch:newtdirun} that $\T$ is a model of $T^{\operatorname{nl}}$.  

\begin{cor}\label{Lmodelcomp} The $\mathcal{L}$-theory $T^{\operatorname{nl}}$ is model complete. Thus $T^{\operatorname{nl}}$ is the model companion of the $\mathcal{L}$-theory of $H$-fields and of the $\mathcal{L}$-theory of
pre-$H$-fields. %\marginpar{need to define this in model theory appendix}
\end{cor}
\begin{proof} Apply Corollary~\ref{cor:upoembnl} to the models
$E$ of $T^{\operatorname{nl}}$ and use \ref{cor:modelcomplete test}.
\end{proof}

%\subsection*{Tournant dangereux\marginparsaved{\begin{center}\dbend\end{center}}} 
\subsection*{Tournant dangereux}
 Recall that the valuation ring $\mathcal{O}_{\T}$ of $\T$ is existentially definable without parameters in the differential field $\T$, that is, in $\T$ construed as an $\mathcal L_{\der}$-structure. On the
other hand:

\begin{cor}\label{notmodelcomplete} The valuation ring 
$\mathcal{O}_{\T}$ of $\T$ is not universally definable in the differential field $\T$, 
even if we allow parameters.
\end{cor}
\begin{proof} Take an elementary extension $\T_1$ of the
$H$-field $\T$ with an element $\alpha$ in the
value group $\Gamma_1$ of $\T_1$ such that $0<\alpha < \Gamma_{\T}^{>}$. Let $\T_2$
be the $\Delta$-coarsening of $\T_1$, for
$\Delta := \big\{\gamma\in \Gamma_{1}:\ |\gamma| < \Gamma_{\T}^{>}\big\}$.
It follows easily from Proposition~\ref{Coarsenings-Pdv} that $\T_2$ is a pre-$H$-field extension of
the $H$-field $\T$. Extend the pre-$H$-field
$\T_2$ to an $\upo$-free newtonian Liouville 
closed $H$-field $K$. Thus $\T\preceq K$, as $H$-fields. 

Suppose $\phi(y)$
is a universal formula in the language $\mathcal L_{\der}$ augmented by names for elements of $\T$ such that $\phi(y)$ defines
$\mathcal{O}_{\T}$ in $\T$. Then it defines
$\mathcal{O}_{\T_1}$ in $\T_1$ and~$\mathcal{O}_K$ in~$K$, since $\T_1$ and $K$ are elementary
extensions of $\T$. Take $b\in \T_1$ such that 
$v(b)= -\alpha$. Then $b\notin \mathcal{O}_{\T_1}$, so
$\T_1\models \neg \phi(b)$, hence $K\models \neg \phi(b)$, since $\T\subseteq \T_1\subseteq K$ as differential fields.
Hence $b\notin \mathcal{O}_K$. But $b\in \mathcal{O}_{\T_2}$, and so $b\in \mathcal{O}_K$, a contradiction.
\end{proof}

\noindent
In particular, the theory of $\T$ as a differential field
is not model complete.

\subsection*{Notes and comments} In the next volume we expect to pay more attention to various natural $H$-subfields of $\T$ such as the field $\mathbb T_{\operatorname{g}}$ of grid-based transseries. (The \textit{Notes and comments}\/ to Appendix~\ref{app:trans} sketch how to build~$\mathbb T_{\operatorname{g}}$ inside $\T$.) The construction of $\mathbb T_{\operatorname{g}}$ easily yields that it is $\upo$-free and Liouville closed. It is newtonian by~\cite{JvdH}, especially Chapter~5. Thus
$\mathbb T_{\operatorname{g}}\preceq\mathbb T$ by Corollary~\ref{Lmodelcomp}. Since $\R\subseteq \mathbb T_{\operatorname{g}}$, we
can combine these facts with Theorem~\ref{th:noextL} to conclude:
{\em every $f\in\mathbb T\setminus\mathbb T_{\operatorname{g}}$ is $\d$-transcendental over $\mathbb T_{\operatorname{g}}$}. 
For example, given any nonzero real numbers $c_1,c_2,c_3,\dots$, the transseries
$$c_1 + c_2 2^{-x} + c_3 3^{-x}+\cdots\  =\ \sum_{n=1}^\infty c_n \ex^{-x\log n}$$
is $\d$-transcendental over $\mathbb T_{\operatorname{g}}$. Thus
the series
$\zeta(x)=1+2^{-x}+3^{-x}+\cdots$ for the zeta function is $\d$-trans\-cen\-dental over $\mathbb T_{\operatorname{g}}$; this improves the classical fact that $\zeta(x)$ is $\d$-transcendental, stated by Hilbert in his 1900 ICM address~\cite[p.~287]{Hilbert1900} and proved in~\cite{Ostrowski20,Stadigh}; 
see \cite{Rubel89} for the history, and \cite{Komatsu,Steuding} for other $\d$-transcendence results about Dirichlet series
that can be strenghtened likewise.

\medskip\noindent
Corollary~\ref{notmodelcomplete} corrects an error in~\cite{AvdD3}: Lem\-ma~14.1 there is false. The definition of ``existentially closed $H$-field'' in that paper is therefore not equivalent to the usual definition. 
The mistake is in the equivalence ``$z\succ 1 \Longleftrightarrow \cdots$'' claimed in the alleged proof of that lemma. 
This also led to an incorrect model completeness conjecture on p.~279 of~\cite{ADH}: we should have included the valuation ring of $\T$ among the primitives. (This is the ``minor change'' referred to in our Preface.)

In our treatment the proof of model completeness of~$T^{\operatorname{nl}}$ is a key step towards 
QE (quantifier elimination) in a suitably extended language, rather than
model completeness being obtained as a consequence of QE, as often happens. 
This gives hope that in possible extensions of our work, for example to 
$\T_{\log}$, model completeness, rather than~QE, might be a realistic first aim.

\section{$\HLO$-Cuts and $\HLO$-Fields}\label{sec:LO-cuts}

\noindent
{\em Throughout this section $K$ is a pre-$H$-field}. The reader needs to be familiar with
Sections~\ref{sec:H-fields}, ~\ref{sec:Liouville closed} and~\ref{sec:special sets}.  
We shall determine 
the sets $I, \Lambda, \Omega\subseteq K$ for which $(K, I, \Lambda, \Omega)$ can be embedded into 
$(L, \I(L), \Upl(L), \Upo(L)\big)$ for some $\upo$-free newtonian
Liouville closed $H$-field $L$.  
%It turns out that for any~$K$ there is either exactly one such triple $(I, \Lambda, \Omega)$ or there are exactly two. 
We first consider for any given $K$ just the possibilities for
$I$, next we determine the possibilities for $(I, \Lambda)$, and
finally, the possibilities for~$(I, \Lambda, \Omega)$. This section is independent of the previous two.

\subsection*{$\I$-sets}  
By an {\bf $\I$-set} in $K$ we mean an $\mathcal{O}$-submodule $I$ of~$K$
such that $\der \mathcal{O} \subseteq I$  and $a^\dagger\notin I$ for all
$a\succ 1$ in $K$. An $\I$-set in $K$ is in particular a convex additive subgroup of the ordered additive group of $K$. 
Note that 
$$\I(K)\ =\ \{y\in K:\text{$y\preceq f'$ for some $f\in \mathcal O$}\}$$ 
is the smallest $\I$-set in $K$. If 
$K$ has no gap, then $\I(K)$ is the only $\I$-set in $K$. 
If $K$ has a gap $\beta$ and $v(b')=\beta$ for some $b\asymp 1$ in $K$,
then $\I(K)$ is the only $\I$-set in $K$. 

\index{I-set@$\I$-set}
\index{set!$\I$-set}

\begin{lemma}\label{twoIsets} Suppose $K$ has gap $\beta$ and $v(b')\ne \beta$ for all
$b\asymp 1$ in $K$. Then~$K$ has exactly two $\I$-sets, namely
 $\I(K)=\{a\in K: va>\beta\}$ and $\{a\in K: va\ge \beta\}$. 
\end{lemma}

\noindent
Clearly, if $K$ is an $H$-field with gap~$\beta$, then the hypothesis of Lemma~\ref{twoIsets} holds. 
%since then $v(b')>\beta$ for all $b\preceq 1$ in $K$. 

\begin{lemma}\label{indIset} If $L$ is a pre-$H$-field extension
of $K$ and $J$ is an $\I$-set in $L$, then $J\cap K$ is an $\I$-set
in $K$. As a strong converse, if $I$ is an $\I$-set in $K$, then 
$K$ has a Liouville closed $H$-field extension $L$ such that
$I=\I(L)\cap K$. 
\end{lemma}
\begin{proof} The first part of the lemma is obvious. 
For the second part, first note that~$K$ has a
Liouville closed $H$-field extension $L$, and so the desired conclusion follows
from the first part if $K$ has just one $\I$-set. It remains to consider
the case that $K$ has a gap $\beta$ and $v(b')\ne \beta$ for all $b\asymp 1$ in $K$. Take $s\in K$ with $vs=\beta$.

We first deal with $I=\{a\in K:\ va\ge \beta\}$.
Lemma~\ref{extensionpdvfields1} and subsequent remarks, combined with Corollary~\ref{extensionpdvfields1, cor}, give a pre-$H$-field extension $K(y)$ of $K$ such that $y'=s$, $y\prec 1$, and $K(y)$ has
no gap. Then $s\in \I\!\big(K(y)\big)$,
so $\I\!\big(K(y)\big)\cap K = I$ by Lemma~\ref{twoIsets}. Taking a Liouville closed $H$-field extension $L$ of $K(y)$, we get $\I(L)\cap K(y)=\I\!\big(K(y)\big)$, and thus $\I(L)\cap K=I$. 

Next, let $I=\I(K)=\{a\in K:\ va > \beta\}$. 
The $H$-field hull
$H(K)$ of $K$ is an immediate extension of $K$ by Corollary~\ref{cor:value group of dv(K)}. So $\beta$ is still a gap in $H(K)$ and $v(b')\ne \beta$ for all $b\asymp 1$ in $H(K)$, and thus
$\I\!\big(H(K)\big)\cap K=I$.
Since $H(K)$ is $\d$-valued, we can apply Corollary~\ref{variant41, cor} to give an $H$-field extension $F:=H(K)(y)$ 
such that $y'=s$, $y\succ 1$, and $F$ has no gap. Then $s\notin \I(F)$, and so  $\I(F)\cap K=I$ by Lemma~\ref{twoIsets}. Taking a Liouville closed $H$-field extension $L$ of~$F$, we get $\I(L)\cap F= \I(F)$, and
thus $\I(L)\cap K=I$.  
\end{proof}

\subsection*{$\HL$-cuts} To motivate the notion of a $\HL$-cut, 
suppose $K$ is a Liouville closed $H$-field. 
Various definitions and results in Section~\ref{sec:special sets} show that then 
the subsets $I:= \I(K)$ and 
$\Lambda:= \Upl(K) =  -(K^{\succ 1})^{\dagger\dagger}$
of $K$ have the following universal properties:
\begin{list}{}{\leftmargin=3em }
\item[($\HLO 1$)] $I$ is an $\I$-set in $K$;
\item[($\HLO 2$)] $\Upl(K)\subseteq \Lambda$;
\item[($\HLO 3$)] for all $a\in K^\times$:  
$\ a\in I\Longleftrightarrow -a^\dagger\notin \Lambda$;
\item[($\HLO 4$)] $\Lambda + I \subseteq \Lambda$;
\item[($\HLO 5$)] $\Lambda$ is downward closed;
\item[($\HLO 6$)] for all $a, \phi\in K$, if $a,\phi>0$ and $a\succeq 1$,
then $\phi a - \phi^\dagger\notin \Lambda$.
\end{list}
We define a {\bf $\HL$-cut} in $K$ to be a pair $(I, \Lambda)$ of sets $I, \Lambda\subseteq K$ satisfying conditions ($\HLO 1$)--($\HLO 6$) above. Only ($\HLO 5$) and ($\HLO 6$) involve the ordering of~$K$. 
Although~$I$ is determined by $\Lambda$ via ($\HLO 3$), it is convenient to make it part of a $\HL$-cut. 

\index{Lambda-cut@$\HL$-cut}
\index{cut!$\HL$-cut}

Suppose that $(I, \Lambda)$ is a $\HL$-cut in $K$. It easily follows from ($\HLO 1$)~and~($\HLO 3$) that $\Lambda\cap \Upd(K)=\emptyset$, with $\Upd(K)$ as defined in Section~\ref{sec:special sets}.  
%In the next subsection we use that if
%$f,g\in K^\times$, $f\asymp g$, then $f^\dagger - g^\dagger=
%(f/g)^\dagger\in \I(K)\subseteq I$, so by ($\HLO 4$),
%$$-\frac{1}{2}f^\dagger\in\Lambda\ \Longleftrightarrow\ 
%-\frac{1}{2}g^\dagger\in\Lambda.$$ 
For any differential
subfield $E$ of $K$ we have a $\HL$-cut $(I\cap E, \Lambda\cap E)$
in the pre-$H$-subfield~$E$ of $K$. Thus every pre-$H$-field, having a Liouville closed $H$-field extension, has a
$\HL$-cut. It is easy to check that for $\phi\in K^{>}$ the pair
$$(I, \Lambda)^\phi\ :=\ 
\big(\phi^{-1}I, \phi^{-1}(\Lambda + \phi^\dagger)\big)$$
is a $\HL$-cut in $K^\phi$.

\begin{lemma}\label{IL1} Suppose $K$ is grounded. Then $K$ has a unique $\HL$-cut. 
\end{lemma}
\begin{proof} By compositional conjugation we arrange $\max \Psi=0$. Then $\I(K)=\smallo$.
Assume $(\I(K), \Lambda)$ is a $\HL$-cut in $K$. Since 
$1\notin \I(K)$, we have $-1^\dagger=0\in \Lambda$ by ($\HLO 3$), so
$\smallo^{\downarrow}\subseteq \Lambda$ by ($\HLO 4$) and ($\HLO 5$). 
If $a\in K^{>}$ and $a \succeq 1$, then taking $\phi=1$ in ($\HLO 6$) shows that
$a\notin \Lambda$. Thus $\Lambda=\smallo^{\downarrow}$.
\end{proof}

\begin{lemma}\label{IL2} Suppose $K$ has a gap $\beta$ and $v(b')=\beta$ for some
$b\asymp 1$ in $K$. Then~$K$ has a unique $\HL$-cut.
\end{lemma}
\begin{proof} 
We arrange by compositional conjugation that $\beta=0$. Let $(I, \Lambda)$ be a $\HL$-cut in $K$. 
Then $I=\I(K)=\mathcal{O}$, so $1\in I$, hence $0\notin \Lambda$ by ($\HLO 3$), and thus 
$$\Lambda\ \subseteq\ \{a\in K:\ a< I\}$$ by ($\HLO 4$) and ($\HLO 5$). 
Suppose $a\in K$ and $a< I$. Then $-a^\dagger\in \Lambda$ by ($\HLO 3$),
so $-a^\dagger < I$, hence $a^\dagger \succ 1$. Since the derivation of
$K$ is small, we get $v(a^\dagger)=o(va)$ by Lemma~\ref{PresInf-Lemma2}(iv), and so
$a < -a^\dagger < 0$, hence $a\in \Lambda$ by ($\HLO 5$).  
This yields $\Lambda= \{a\in K:\ a< I\}$.  
\end{proof}

\begin{lemma}\label{IL3} Suppose $K$ has a gap $\beta$ and $v(b')\ne \beta$
for all $b\asymp 1$ in $K$. 
Then there are exactly two $\HL$-cuts $(I, \Lambda)$ in $K$,  one with $I~=~\I(K)$, and the other with $I~=~\{y\in K:\ vy\ge \beta\}$. 
\end{lemma}
\begin{proof} Lemma~\ref{indIset}
provides for each of the two $\I$-sets $I$ in $K$
a Liouville closed $H$-field extension
$L$ of $K$ with $\I(L)\cap K=I$, giving rise to a $\HL$-cut
$\big(I, \Upl(L)\cap K\big)$ in~$K$. To get uniqueness, arrange $\beta=0$ by compositional conjugation. Then $\I(K)=\smallo$ and
$\{y\in K:\ vy\ge \beta\}=\mathcal{O}$. If $(\smallo, \Lambda)$ is a
$\HL$-cut in $K$, then we get $\Lambda=\smallo^{\downarrow}$ as
in the proof of Lemma~\ref{IL1}. If $(\mathcal{O}, \Lambda)$ is
a $\HL$-cut in $K$, then  
we get $\Lambda= \{a\in K:\ a< \mathcal{O}\}$ as in the proof of Lemma~\ref{IL2}.
\end{proof}

\noindent
Thus if $K$ is an $H$-field with a gap $\beta$, then 
there are exactly two $\HL$-cuts $(I_1, \Lambda_1)$ and $(I_2,\Lambda_2)$ in $K$, with $I_1=\I(K)=\{y\in K:\ vy> \beta\}$ and 
$I_2=\{y\in K:\ vy\ge \beta\}$. 

\begin{lemma}\label{IL4} If $K$ is $\upl$-free, then $\big(\!\I(K), \Upl(K)^{\downarrow}\big)$ is the only $\HL$-cut in $K$.
\end{lemma}
\begin{proof} Let $K$ be $\upl$-free. Then $K$ has asymptotic integration, so $\I(K)$ is the only $\I$-set in $K$.
Let $\big(\!\I(K), \Lambda\big)$ be a $\HL$-cut in $K$. Then 
$\Upl(K)^{\downarrow}\subseteq \Lambda\subseteq K\setminus \Upd(K)^{\uparrow}$, and so it remains to note that by $\upl$-freeness we have
$\Upl(K)^{\downarrow}= K\setminus \Upd(K)^{\uparrow}$.
\end{proof}

\begin{lemma}\label{IL5} Suppose $K$ has asymptotic integration and is not $\upl$-free. Then there are exactly two $\HL$-cuts in $K$: 
$ \big(\!\I(K), \Upl(K)^{\downarrow}\big)$, and $\big(\!\I(K), K\setminus \Upd(K)^{\uparrow}\big)$.
\end{lemma}
\begin{proof} Let $(I, \Lambda)$ be a $\HL$-cut in $K$. Then $I=\I(K)$, and we claim that
$\Lambda=\Upl(K)^{\downarrow}$ or $\Lambda= K\setminus \Upd(K)^{\uparrow}$.
To see this, first note that $\Upl(K)^{\downarrow}\subseteq \Lambda\subseteq K\setminus \Upd(K)^{\uparrow}$. Suppose
$\Lambda\ne \Upl(K)^{\downarrow}$, and take $a\in \Lambda\setminus \Upl(K)^{\downarrow}$. Then $\Upl(K) < a < \Upd(K)$, and for every $b\in K$ with $\Upl(K) < b < \Upd(K)$ we have $v(b-a)> \Psi$ by Corollary~\ref{cor:Upl 1}, so $b-a\in \I(K)$, and hence
$b\in \Lambda + \I(K)\subseteq \Lambda$. Thus $\Lambda= K\setminus \Upd(K)^{\uparrow}$, which proves our claim.  
  It remains to show that there is more than one $\HL$-cut in $K$.
Let $E:= H(K)^{\operatorname{rc}}$ be the real closure of the
$H$-field hull of $K$. Now $H(K)$ is an immediate extension
of~$K$ by Corollary~\ref{cor:value group of dv(K)}, so $E$ is an
$H$-field extension of $K$ with divisible value group $\Gamma_E=\Q\Gamma$ and so $\Psi_E=\Psi$. We now distinguish two cases:

\case[1]{$E$ has a gap.}
Then we can take $s\in E^\times$ and $n\geq 1$ such that $vs$ is a gap in $E$ and $s^n\in K$. Note that then $s^\dagger\in K$. 
Now $E$ is an $H$-field, so Lemma~\ref{IL3} yields $\HL$-cuts 
$(I_1, \Lambda_1\big)$ and 
 $(I_2, \Lambda_2)$ in $E$, with $I_1=\I(E)=\{y\in E:\ y\prec s\}$ and 
 $I_2=\{y\in E:\ y\preceq s\}$. Now $s\notin I_1$, so
 $-s^\dagger\in \Lambda_1$ by ($\HLO 3$), and $s\in I_2$, so $-s^\dagger\notin \Lambda_2$, also by ($\HLO 3$).
Hence $-s^\dagger\in \Lambda_1\cap K$ and 
$-s^\dagger\notin \Lambda_2\cap K$. So we have two distinct 
$\HL$-cuts in $K$, namely
$$(I_1\cap K, \Lambda_1\cap K)\ =\ \big(\!\I(K), K\setminus \Upd(K)^{\uparrow}\big), \quad (I_2\cap K, \Lambda_2\cap K)\ =\ \big(\!\I(K), \Upl(K)^{\downarrow}\big).$$   

\case[2]{$E$ has no gap.} Then $E$ has asymptotic integration, and
the sequence~$(\upl_{\rho})$ for~$K$ also serves for $E$.
Take $\upl\in K$ such that $\upl_{\rho} \leadsto \upl$. Then 
$-\upl$ creates a gap over~$E$ by Lemma~\ref{cg2}. Take an element $f\ne 0$ in some Liouville closed $H$-field extension of $E$ such that
$f^\dagger=-\upl$. Then $vf$ is a gap in $E(f)$ by the remark following the proof of Lemma~\ref{cg2}. Using Lemma~\ref{dv} it follows that 
the pre-$H$-field
 $E(f)$ is actually an $H$-field. Hence by Lemma~\ref{IL3}, 
$E(f)$ has $\HL$-cuts $(I_1, \Lambda_1)$ and 
 $(I_2, \Lambda_2)$, with $I_1=\I\!\big(E(f)\big)$ and $I_2=\{a\in E(f):\ a\preceq f\}$.
Now $f\notin I_1$, so $-f^\dagger=\upl\in \Lambda_1$ by ($\HLO 3$), and $f\in I_2$, so $\upl\notin \Lambda_2$, again by ($\HLO 3$). Thus
$$(I_1\cap K, \Lambda_1\cap K)\ =\ \big(\!\I(K), K\setminus \Upd(K)^{\uparrow}\big), \quad (I_2\cap K, \Lambda_2\cap K)\ =\ \big(\!\I(K), \Upl(K)^{\downarrow}\big)$$ 
are two distinct $\HL$-cuts in $K$. 
\end{proof}

\noindent
The proofs above have the following byproduct: 

\begin{cor} \label{cor:IL}
Let $(I,\Lambda)$ be a $\HL$-cut in $K$. 
Then $K$ has a Liouville closed $H$-field extension $L$ such that
$(I, \Lambda)=\big(\!\I(L)\cap K, \Upl(L)\cap K\big)$.
\end{cor}
\begin{proof} If $K$ has a unique $\HL$-cut, then the conclusion of
the corollary holds for any Liouville closed $H$-field extension $L$ of $K$.
The case that $K$ has a gap $\beta$ and $v(b')\ne \beta$ for all $b\asymp 1$ in $K$ is treated in the proof of Lemma~\ref{IL3}. The case that
$K$ has asymptotic integration and is not $\upl$-free reduces to
the previous cases by extending $K$ to $E$ or $E(f)$ as
in the proof of Lemma~\ref{IL5}.
\end{proof}

%\begin{remark} \marginpar{Added this remark.}
%Let $(I,\Lambda)$ be a $\HL$-cut in $K$. 
%The proof of Corollary~\ref{cor:IL} shows that $L$ as in that 
%corollary can additionally be chosen so that $L$ is a Liouville 
%extension of~$K$ and $C_L$ is a real closure of $C$.
%\end{remark}
 
\noindent
This corollary, when combined with Proposition~\ref{prop:omega for H-fields} and Lemma~\ref{lem:sigma increasing}, has the following consequence:

\begin{cor}\label{cor:omega and sigma incr}
Let $(I,\Lambda)$ be a $\HL$-cut in $K$. Then the functions 
$$z\mapsto \omega(z)\colon\Lambda\to K,\qquad y\mapsto\sigma(y)\colon K^>\setminus I\to K$$ 
are strictly increasing.
\end{cor}

\subsection*{$\HLO$-cuts} By various results in Section~\ref{sec:special sets} the sets
$I:= \I(K)$, $\Lambda:=\Upl(K)$, and $\Omega:= \Upo(K)$ of any Schwarz closed $K$ have the following universal properties: 
\begin{list}{}{\leftmargin=4em \labelwidth=4em}
\item[($\HLO 7$)] $\omega(K)\subseteq \Omega$ (and so $0\in \Omega$);
\item[($\HLO 8$)]  for all $f,g \in K^>$, if $-\frac{1}{2}f^\dagger\in\Lambda$ and $f\asymp g$,
then $\omega(-\frac{1}{2}f^\dagger)+g\notin\Omega$;
\item[($\HLO 9$)] for all $f\in K^\times$, if $-\frac{1}{2}f^\dagger\notin\Lambda$,
then $\Omega+f\subseteq\Omega$;
\item[($\HLO 10$)] $\Omega$ is downward closed.
\end{list}
Accordingly we define a {\bf $\HLO$-cut} in $K$ to be a triple $(I, \Lambda, \Omega)$
of sets $I, \Lambda, \Omega\subseteq K$ such that $(I, \Lambda)$ is 
a $\HL$-cut in $K$ and conditions ($\HLO 7$)--($\HLO 10$) above are satisfied.
Let us record some easy consequences of the axioms above: 

\index{LambdaOmega-cut@$\HLO$-cut}
\index{cut!$\HLO$-cut}

\begin{lemma}\label{IOlemma}
Let  $(I, \Lambda, \Omega)$ be a $\HLO$-cut in $K$. Then
\begin{list}{}{\leftmargin=4em \labelwidth=4em}
\item[$(\HLO 11)$] 
for all $g,h\in I$, we have $\Omega + gh\subseteq \Omega$;
\item[$(\HLO 12)$] for all $f\in K^{\times}$: $f\in I\Longleftrightarrow \sigma(f)\in \Omega$.
\end{list}
\end{lemma}
\begin{proof}
To show ($\HLO 11$), let $0\neq g,h\in I$, and set $f:=gh$.
By ($\HLO 9$), it suffices to prove $-\frac{1}{2}f^\dagger\notin\Lambda$.
We may assume $g=ha$, $a\preceq 1$. So $f=h^2a$ and hence $-\frac{1}{2}f^\dagger=-h^\dagger-\frac{1}{2}a^\dagger$, and $-h^\dagger\notin\Lambda$ by ($\HLO 3$). 
If $a\asymp 1$, then ($\HLO 4$) and $a^\dagger\in \I(K)\subseteq I$ give $-\frac{1}{2}f^\dagger\notin\Lambda$. If $a\prec 1$, then $-\frac{1}{2}a^\dagger>0$ and hence
$-\frac{1}{2}f^\dagger>-h^\dagger$ and so $-\frac{1}{2}f^\dagger\notin\Lambda$.

For $(\HLO 12)$, let $f\in K^\times$. If $f\in I$, then $\sigma(f)=\omega(-f^\dagger)+f^2\in\Omega$ by $(\HLO 11)$.
Suppose $f\notin I$. Then $-\frac{1}{2}(f^2)^\dagger=-f^\dagger\in\Lambda$ by $(\HLO 3)$,
so $\sigma(f)=\omega(-f^\dagger)+f^2=\omega(-\frac{1}{2}(f^2)^\dagger)+f^2\notin\Omega$ by $(\HLO 8)$. 
\end{proof}
 
\noindent
Suppose $(I, \Lambda, \Omega)$ is a $\HLO$-cut in $K$. It follows easily 
from ($\HLO 10$) and ($\HLO 12$) that if $a\in K\setminus I$, then 
$\Omega < \sigma(a)$. 
In particular, we have $\Omega<\sigma\big(\Gamma(K)\big)$,  so
$$\omega(K)^\downarrow\ \subseteq\ \Omega\ \subseteq\  K\setminus \sigma\big(\Upg(K)\big){}^\uparrow.$$ 
For any differential
subfield $E$ of $K$ we have a $\HLO$-cut $(I\cap E, \Lambda\cap E, \Omega\cap E)$
in the pre-$H$-subfield $E$ of $K$. Thus every pre-$H$-field, having a 
Schwarz closed $H$-field extension by Corollary~\ref{cor:embed into upo-free} and Proposition~\ref{prop:Schwarz closure},  has a
$\HLO$-cut. Moreover, any $\HL$-cut $(I,\Lambda)$ in $K$ is part of
a $\HLO$-cut $(I, \Lambda, \Omega)$ in $K$. (Proof: let $(I,\Lambda)$ be a $\HL$-cut in~$K$. Take a Liouville closed extension $L$ of $K$
with $\I(L)\cap K=I$ and $\Upl(L)\cap K = \Lambda$, and next a Schwarz closed extension $M$ of $L$. Then $\Omega=\Upo(M)\cap K$ gives a
$\HLO$-cut $(I, \Lambda, \Omega)$ in $K$.)
For $\phi\in K^{>}$ one can use the identities at the end of Section~\ref{sec:behupo} to show that the triple
$(I, \Lambda, \Omega)^\phi=(I^{\phi}, \Lambda^{\phi}, \Omega^{\phi})$ with
\begin{equation}\label{eq:HL-cut comp conj}
I^{\phi}\ :=\ 
\phi^{-1}I,\quad \Lambda^{\phi}\ :=\ \phi^{-1}(\Lambda + \phi^\dagger),\quad \Omega^{\phi}\ :=\  \phi^{-2}\big(\Omega-\omega(-\phi^\dagger)\big)
\end{equation}
is a $\HLO$-cut in $K^\phi$.

\begin{lemma}\label{ILO1} Suppose $K$ is grounded. Then $K$ has a unique $\HLO$-cut. 
\end{lemma}
\begin{proof} By compositional conjugation we arrange $\max \Psi=0$. 
Let $(I, \Lambda, \Omega)$ be a $\HLO$-cut in $K$.
Then $I=\smallo$ and $\Lambda=\smallo^{\downarrow}$ by the proof of Lemma~\ref{IL1}. Let $u\asymp 1$ in $K$. Then $u^{\dagger}\prec 1$ and
$(u^{\dagger})'\prec 1$, so $\sigma(u)=\omega(-u^\dagger)+u^2\sim u^2$.
Also $u\notin I$, so $\Omega <\sigma(u)$. Hence
$\Omega < 2u^2$ for all units $u$ of $\mathcal{O}$, and so 
$\Omega\subseteq \smallo^{\downarrow}$. We claim that
$\Omega=\smallo^{\downarrow}$. By ($\HLO 10$) it is enough to show $\smallo\subseteq \Omega$. We distinguish two cases:

\case[1]{$\Gamma$ has no least positive element.} Then every element
of $\smallo$ is a product $gh$ with $g,h\in \smallo=I$, and as 
$0\in \Omega$, we get $\smallo\subseteq \Omega$ from ($\HLO 11$), and so our claim holds.

\case[2]{$\Gamma$ has a least positive element $\alpha$.}
Take $a\in K$ with $va=\alpha$ and set $b:=-2a$. From $\alpha^\dagger=0$ we get
$v(a')=\alpha$ and so $\omega(a)\sim -2a'=b'$,
hence $v(\omega(a)-b')\ge 2\alpha$, and then $b'\in \Omega$ by ($\HLO 7$) and ($\HLO 11$). 
Thus $b'\in\Omega$ for each $b\in K$ with $vb=\alpha$. 
For such $b$ and  $u\asymp 1$ in $K$ we have 
$v\big((ub)'- ub'\big)=v(u'b)\ge 2\alpha$, and thus $ub'\in \Omega$ by ($\HLO 11$).
Fixing $b$ and varying $u$ we see that all elements of $K$ of valuation~$\alpha$ belong to~$\Omega$, and thus $\smallo\subseteq \Omega$, as desired.  
\end{proof}

\begin{lemma}\label{ILO2} Suppose $K$ has a gap $\beta$ and $v(b')=\beta$ for some
$b\asymp 1$ in $K$. Then~$K$ has a unique $\HLO$-cut.
\end{lemma}
\begin{proof} 
We arrange by compositional conjugation that $\beta=0$. Let $(I, \Lambda, \Omega)$ be
a $\HLO$-cut in $K$. Then $I=\mathcal{O}$ and $\Lambda=\{a\in K:\ a<\mathcal{O}\}$ by the proof of Lemma~\ref{IL2}. From $0\in \Omega$
and ($\HLO 10$) and ($\HLO 11$) we obtain $\mathcal{O}^{\downarrow}\subseteq \Omega$.
We claim that $\Omega= \mathcal{O}^{\downarrow}$. This is 
clear if $\Gamma=\{0\}$, so suppose $\Gamma\neq\{0\}$. Then $\Gamma^<$ does not have a largest element.
For
$a\in K^{>}$ with $a\succ 1$ we have
$\sigma(a)\sim a^2$, as well as $a\notin I$, so $\sigma(a)> \Omega$ by
($\HLO 7$) and~($\HLO 12$),
and thus $\Omega< 2a^2$.  This yields the claim. 
\end{proof}

\begin{lemma}\label{ILO3} Suppose $K$ has a gap $\beta$ and 
$v(b')\ne \beta$ for all $b\asymp 1$ in $K$. 
Then there are exactly two $\HLO$-cuts $(I, \Lambda, \Omega)$ in $K$,  one with $I=\I(K)$, and the other with $I=\{y\in K:\, vy\ge \beta\}$. 
\end{lemma}
\begin{proof} Lemma~\ref{indIset}
provides for each of the two $\I$-sets $I$ in $K$
a Liouville closed $H$-field extension
$L$ of $K$ with $\I(L)\cap K=I$, and by further extending
we can arrange~$L$ to be even Schwarz closed, giving rise to a $\HLO$-cut
$\big(I, \Upl(L)\cap K, \omega(L)\cap K\big)$ in~$K$. To get uniqueness, arrange $\beta=0$ by compositional conjugation. Then by the proof
of Lemma~\ref{IL3} we have
a $\HLO$-cut $(\smallo, \smallo^{\downarrow}, \Omega)$ in $K$. The same
arguments as in the proof of Lemma~\ref{ILO1} show that for each such
$\HLO$-cut we have $\Omega=\smallo^{\downarrow}$.  We also have a $\HLO$-cut $(\mathcal{O}, \Lambda, \Omega)$ with $\Lambda=\{a\in K:\ a < \mathcal{O}\}$. The same argument
as in the proof of Lemma~\ref{ILO2} shows that 
$\Omega=\mathcal{O}^{\downarrow}$ for each such $\HLO$-cut. 
\end{proof}

\noindent
In the next lemmas we treat the case where $K$ has asymptotic integration.
Recall Lemma~\ref{IL4}: if $K$ is $\upl$-free, then $\big(\!\I(K),\Upl(K)^{\downarrow}\big)$ is the only $\HL$-cut in $K$.

\begin{lemma}\label{ILO4} Suppose $K$ is $\upo$-free. Then the only $\HLO$-cut
in $K$ is
$$\big(\!\I(K), \Upl(K)^{\downarrow}, \omega(K)^{\downarrow}\big).$$
\end{lemma} 
\begin{proof} Let $(I, \Lambda, \Omega)$ be a $\HLO$-cut in $K$, so
$I=\I(K)$ and $\Lambda=\Upl(K)^{\downarrow}$.
Also, 
$\omega(K)^{\downarrow}\subseteq \Omega\subseteq K\setminus \sigma\big(\Upg(K)\big){}^\uparrow$, and so it remains to note that by $\upo$-freeness we have
$\omega(K)^{\downarrow}= K\setminus \sigma\big(\Upg(K)\big){}^\uparrow$.
\end{proof}

\begin{lemma}\label{ILO5}
Suppose $K$ has asymptotic integration and the set~$2\Psi$ does not have a supremum in $\Gamma$.
Then for each $\HLO$-cut $(I, \Lambda, \Omega)$ in $K$ we have
$$\Omega\ =\ \omega\big(\Upl(K)\big){}^\downarrow\ =\ \omega(K)^\downarrow \quad\text{ or }\quad \Omega\ =\ K\setminus\sigma\big(\Upg(K)\big){}^\uparrow.$$ 
\end{lemma}
\begin{proof}
Let $(I, \Lambda, \Omega)$ be a $\HLO$-cut in $K$. Recall that
$$\omega\big(\Upl(K)\big){}^\downarrow\ \subseteq\ \Omega\ \subseteq\  K\setminus \sigma\big(\Upg(K)\big){}^\uparrow.$$ 
Note that if $\Omega = \omega\big(\Upl(K)\big){}^{\downarrow}$, then
$\Omega=\omega(K)^{\downarrow}$ by~($\HLO 7$).
Suppose
$\Omega\ne \omega\big(\Upl(K)\big){}^\downarrow$, and take $a\in \Omega\setminus \omega\big(\Upl(K)\big){}^\downarrow$. Then 
$$\omega\big(\Upl(K)\big)\  <\  a\  <\  \sigma\big(\Upg(K)\big),$$ and for every $b\in K$ with $\omega\big(\Upl(K)\big) < b < \sigma\big(\Upg(K)\big)$ we have $v(b-a)> 2\Psi$ by Corollaries~\ref{pcrhocor} and~\ref{cor:omega-freeness and cut}, 
so by Lemma~\ref{lem:split gamma} there are $g,h\in \I(K)$ with $b-a=gh$, and thus 
$b\in \Omega + gh \subseteq \Omega$ by $(\HLO 11)$. Thus $\Omega= K\setminus \sigma\big(\Upg(K)\big){}^\uparrow$.
\end{proof}

\begin{lemma}\label{ILO6} Suppose $K$ is $\upl$-free, but not $\upo$-free.
Then there are exactly two $\HLO$-cuts in $K$:
$ \big(\!\I(K), \Upl(K)^{\downarrow}, \omega(K)^{\downarrow}\big)$, and $\big(\!\I(K), \Upl(K)^{\downarrow}, K\setminus \sigma\big(\Upg(K)\big){}^\uparrow\big)$.
\end{lemma}
\begin{proof} 
Since $K$ is $\upl$-free, $K$ has rational asymptotic integration, and so $2\Psi$ does not have a supremum in $\Gamma$.
Let $(I, \Lambda, \Omega)$ be a $\HLO$-cut in $K$. Then $I=\I(K)$, $\Lambda=\Upl(K)^{\downarrow}$, and either $\Omega=\omega(K)^{\downarrow}$ or $\Omega= K\setminus \sigma\big(\Upg(K)\big){}^\uparrow$ by Lemma~\ref{ILO5}.

It remains to show that there is more than one $\HLO$-cut in $K$.
Take $\upo\in K$ such that $\upo_{\rho} \leadsto \upo$. Then
Corollary~\ref{cor:adjoin upl} and Lemma~\ref{prehimm} yield an immediate pre-$H$-field
extension $K(\upl)$ of $K$ with $\upl_{\rho} \leadsto \upl$ and
$\omega(\upl)=\upo$. Now $K(\upl)$ has a $\HLO$-cut  
$(I_{\upl}, \Lambda_{\upl},\Omega_{\upl})$. Then $\upo\in \Omega_{\upl}$,
and so intersecting with $K$ gives a $\HLO$-cut $(I, \Lambda, \Omega)$ in~$K$
with $\upo\in \Omega$.  
It remains to find such a $\HLO$-cut in $K$
with $\upo\notin \Omega$. For that we take a pre-$H$-field extension 
$K\<\upg\>$ as in Lemma~\ref{cor:sigma(upg)=upo} in which $\sigma(\upg)=\upo$ and $v\upg$ is a gap. By 
remark~(3) following the proof of Proposition~\ref{prop:sigma(upg)=upo}
there is no $b\asymp 1$ in $K\<\upg\>$ with $b'\asymp \upg$.
Then by Lemma~\ref{IL3}  we have
a $\HL$-cut $(I_{\upg}, \Lambda_{\upg})$ in $K\<\upg\>$
with $\upg\notin I_{\upg}$. Take a $\HLO$-cut 
$(I_{\upg}, \Lambda_{\upg}, \Omega_{\upg})$ in $K\<\upg\>$.
Then $\upo=\sigma(\upg)\notin \Omega_{\upg}$ by ($\HLO 12$). Intersecting with $K$
gives a $\HLO$-cut $(I, \Lambda, \Omega)$ in $K$ with $\upo\notin \Omega$.
\end{proof}

\begin{lemma}\label{ILO7} Suppose $K$ has asymptotic integration and is not $\upl$-free,  and the set~$2\Psi$ does not have a supremum in $\Gamma$. Then  there are exactly two $\HLO$-cuts in $K$:
$$\big(\!\I(K), \Upl(K)^{\downarrow}, K\setminus \sigma\big(\Upg(K)\big){}^{\uparrow}\big)\  \text{ and }\ \big(\!\I(K), K\setminus \Upd(K)^{\uparrow}, K\setminus \sigma\big(\Upg(K)\big){}^{\uparrow} \big).$$
\end{lemma} 
\begin{proof} Let $(I, \Lambda, \Omega)$ be a $\HLO$-cut in $K$.
Let $\upl_{\rho} \leadsto \upl\in K$. Then
$\upo_{\rho} \leadsto \upo:= \omega(\upl)$ by Corollary~\ref{pczrho} and $\upo\notin \omega\big(\Upl(K)\big){}^\downarrow$
by Corollary~\ref{cor:omega-freeness and cut}, but
$\upo\in \omega(K)\subseteq \Omega$. Hence 
$\Omega = K\setminus \sigma\big(\Upg(K)\big){}^\uparrow$
by Lemma~\ref{ILO5}.  It remains to use Lemma~\ref{IL5}.
\end{proof} 

\noindent
It follows from Lemma~\ref{lem:split gamma} that in Lemma~\ref{ILO7} we can drop the condition that~$2\Psi$ has no supremum in $\Gamma$ if $\Gamma$ is divisible.

\begin{lemma}\label{ILO8}
Suppose $K$ has asymptotic integration and $vf=\sup 2\Psi$, $f\in K^{>}$. Then 
there are exactly two $\HLO$-cuts in $K$, namely
\begin{align*}\big(\!\I(K), \Upl(K)^{\downarrow}&, K\setminus \sigma\big(\Upg(K)\big){}^{\uparrow}\big)\  \text{ and }\ \big(\!\I(K), K\setminus \Upd(K)^{\uparrow}, \Omega_f \big),\ \text{ where}\\
\Omega_f\ &:=\ \left\{ a\in K:\ \text{$a\leq \omega(-\textstyle\frac{1}{2}f^\dagger)+g$ for some $g\prec f$ in $K$} \right\}.
\end{align*}
Moreover, $\Upl(K)^{\downarrow}\neq K\setminus \Upd(K)^{\uparrow}$ and $K\setminus \sigma\big(\Upg(K)\big){}^{\uparrow}\neq \Omega_f$.
\end{lemma}
\begin{proof}
Let $\sqrt{f}$ be a positive element of the real closure of $K$ with $(\sqrt{f})^2=f$. The
pre-$H$-field extension $L=K(\sqrt{f})$ of $K$ has gap $\beta=v(\sqrt{f})=\frac{1}{2}vf$.
Hence $\Upl(K)<-\frac{1}{2}f^\dagger=-(\sqrt{f})^\dagger<\Upd(K)$ by Corollaries~\ref{cor:gap} and~\ref{cor:Upl 1}.
Thus $K$ is not 
$\upl$-free, and hence by Lemma~\ref{IL5}, $K$ has exactly two $\HL$-cuts,
$ \big(\!\I(K), \Upl(K)^{\downarrow}\big)$ and $\big(\!\I(K), K\setminus \Upd(K)^{\uparrow}\big)$.
By Lemma~\ref{alggapcase} there is
no $b\asymp 1$ in $L$ with $b'\asymp \sqrt{f}$, so by Lemma~\ref{IL3},
$L$ also has exactly two $\Upl$-cuts, $(I_1,\Lambda_1)$ and $(I_2, \Lambda_2)$, with 
$$I_1\ =\ \I(L)\ =\ \{y\in L:\ vy>\beta\}, \qquad I_2\ =\ \{y\in L:\ vy\ge \beta\}.$$ 
We have $\sqrt{f}\notin I_1$ and $\sqrt{f}\in I_2$, hence
$-\frac{1}{2}f^\dagger=-(\sqrt{f})^\dagger\in \Lambda_1\cap K$ and $-\frac{1}{2}f^\dagger\notin \Lambda_2\cap K$
by ($\HLO 3$).
Therefore
$$I_1\cap K\ =\ I_2\cap K\ =\ \I(K),\quad 
\Lambda_1\cap K\ =\ K\setminus \Upd(K)^{\uparrow},\quad
\Lambda_2\cap K\ =\  \Upl(K)^{\downarrow}.$$
Note that $\sigma(\sqrt{f}) = \omega(-\frac{1}{2}f^\dagger) + f$ lies in 
$K$.
Also, by Lemma~\ref{ILO3} we have just one $\Upl\Upo$-cut $(I_1,\Lambda_1, \Omega_1)$ in $L$, and just one $\Upl\Upo$-cut $(I_2, \Lambda_2, \Omega_2)$ in $L$. This yields distinct $\Upl\Upo$-cuts
$(I_1\cap K, \Lambda_1\cap K, \Omega_1\cap K)$ and $(I_2\cap K, \Lambda_2\cap K, \Omega_2\cap K)$ in $K$, since $\Lambda_1\cap K\ne \Lambda_2\cap K$; we claim that these are the only $\Upl\Upo$-cuts in $K$, and that $\Omega_f=\Omega_1\cap K\ne \Omega_2\cap K=K\setminus \sigma\big(\Upg(K)\big){}^{\uparrow}$. 

For this, let $(I, \Lambda, \Omega)$ be a $\Upl\Upo$-cut in $K$.
Note that if $g\in K^\times$ and $g\prec f$, 
then Lemma~\ref{lem:split gamma} gives
$a,b\in \I(K) = I$ with $g=ab$, and so by
($\HLO 11$) we obtain $\Omega+g\subseteq\Omega$. In particular, using  ($\HLO 10$) we get $\Omega_f\subseteq\Omega$.
Also as in the proof of Lemma~\ref{ILO7} we  obtain $\omega(-\frac{1}{2}f^\dagger)\in\Omega\setminus\omega\big(\Upl(K)\big){}^\downarrow$.

\case[1]{$\Lambda=\Upl(K)^{\downarrow}$.}
Then 
$-\frac{1}{2}f^\dagger\notin \Lambda_2\cap K=\Lambda$, and thus $\Omega+f\subseteq \Omega$ by ($\HLO 9$). Replacing the role of $f$ by $uf$ for
any $u\in K^{>}$ with $u\asymp 1$ we get
$\Omega+g\subseteq\Omega$ for all $g\in K^\times$ with $vg>2\Psi$.
Then the proof of Lemma~\ref{ILO5} yields
$\Omega  = K\setminus \sigma\big(\Upg(K)\big){}^{\uparrow}$. In particular, $\Omega=\Omega_2\cap K$ and
$\sigma(\sqrt{f}) = \omega(-\frac{1}{2}f^\dagger) + f\in \Omega$. 

\case[2]{$\Lambda=K\setminus \Upd(K)^{\uparrow}$.} Then $-\frac{1}{2}f^\dagger\in \Lambda_1\cap K=\Lambda$, and  $\sqrt{f}\notin I_1=\I(L)$ gives $\sigma(\sqrt{f})\notin \Omega_1\cap K$ by ($\HLO 12$).
We have $\Omega_f\subseteq\Omega$, and  ($\HLO 8$) and ($\HLO 10$) yield
$\Omega \subseteq \Omega_f$, so $\Omega = \Omega_f = \Omega_1\cap K$.
\end{proof}

\noindent
We have now covered all possibilities, and conclude that $K$ has either exactly one $\HLO$-cut, or exactly two. Moreover:

\begin{cor}\label{cor:ILO1}
The pre-$H$-field $K$ has a unique $\HLO$-cut if and only if
\begin{enumerate}
\item[\textup{(i)}] $K$ is grounded, or
\item[\textup{(ii)}] there exists $b\asymp 1$ in $K$ such that $v(b')$ is a gap in $K$, or
\item[\textup{(iii)}] $K$ is $\upo$-free.
\end{enumerate}
\end{cor}

\noindent
The proofs above have the following byproduct:

\begin{cor} \label{cor:ILO2}
Let $(I,\Lambda, \Omega)$ be a $\HLO$-cut in $K$. 
Then $K$ has a Schwarz closed $H$-field extension $L$ such that
$(I, \Lambda, \Omega)=\big(\!\I(L)\cap K, \Upl(L)\cap K, \Upo(L)\cap K\big)$.
\end{cor}
\begin{proof}
If there is a unique $\HLO$-cut in $K$, then the conclusion of the corollary holds for any
Schwarz closed $H$-field extension $L$ of $K$.
The case that $K$ has a gap $\beta$ with $v(b')\neq\beta$ for all $b\asymp 1$ in $K$ is treated in the proof of Lemma~\ref{ILO3}. 
Suppose $K$ has asymptotic integration but is not $\upl$-free. Then by Lemmas~\ref{ILO7} and \ref{ILO8}, there are
exactly two $\HLO$-cuts
$$(I_1,\Lambda_1,\Omega_1)=\big(\!\I(K),\Upl(K)^\downarrow,\dots\big), \quad
  (I_2,\Lambda_2,\Omega_2)=\big(\!\I(K),K\setminus\Upd(K)^\downarrow,\dots \big)$$
in $K$. By Corollary~\ref{cor:IL} we can take for $i=1,2$ a Liouville closed extension $L_i$ of~$K$ such that $\big(\!\I(L_i)\cap K, \Upl(L_i)\cap K\big) = (I_i,\Lambda_i)$.
Extending $L_i$ if necessary, we can even arrange that each $L_i$ is Schwarz closed. Then 
$$\big(\!\I(L_i)\cap K, \Upl(L_i)\cap K, \Upo(L_i)\cap K\big) = (I_i,\Lambda_i,\Omega_i)\quad\text{for $i=1,2$.}$$
Finally, the case that $K$ is $\upl$-free, but not $\upo$-free, reduces to the case that
$K$ is not $\upl$-free or $K$ has a gap by extending $K$ to $K(\upl)$
and to 
$K\<\upg\>$ as in the proof of Lemma~\ref{ILO6}.
\end{proof}

\noindent
Let $(I,\Lambda, \Omega)$ be a $\HLO$-cut in $K$, let 
$L$ be as in Corollary~\ref{cor:ILO2}, and take 
a Newton-Liouville closure $L^{\operatorname{nl}}$ of $L$. Then the conclusion of that corollary remains valid for~$L^{\operatorname{nl}}$ in place of $L$, since
$\I(L^{\operatorname{nl}})\cap L=\I(L)$, 
$\Upl(L^{\operatorname{nl}})\cap L=\Upl(L)$,  
$\Upo(L^{\operatorname{nl}})\cap L=\Upo(L)$.

\subsection*{$\HLO$-fields}
For model-theoretic use we rephrase some of the results above in the terminology of $\HLO$-fields:

\begin{definition}
A {\bf pre-$\HLO$-field} is a quadruple $(K,I,\Lambda,\Omega)$ where $K$ is a pre-$H$-field  (as throughout this section)
and $(I,\Lambda,\Omega)$ is a $\HLO$-cut in $K$. If in addition~$K$ is an $H$-field, then
we call $(K, I,\Lambda,\Omega)$ a {\bf $\HLO$-field}.
\end{definition} 

\index{pre-$\HLO$-field}
\index{LambdaOmega-field@$\HLO$-field}

\noindent
Since our pre-$H$-field $K$ has a $\HLO$-cut $(I,\Lambda,\Omega)$,
we can turn $K$ into a pre-$\HLO$-field $(K,I,\Lambda,\Omega)$. If $\mathbf K=(K,I,\Lambda,\Omega)$ is a pre-$\HLO$-field and $\phi\in K^{>}$, then $\mathbf K^\phi=(K^\phi,I^\phi,\Lambda^\phi,\Omega^\phi)$,
where $I^\phi$,~$\Lambda^\phi$,~$\Omega^\phi$ are as in \eqref{eq:HL-cut comp conj},  is also
a pre-$\HLO$-field.
Below, any qualifier that applies to pre-$H$-fields, such as \textit{has asymptotic integration}\/, when applied to a pre-$\HLO$-field $(K,\dots)$
means that the underlying pre-$H$-field~$K$ has the property in question.

\medskip
\noindent
Let $\mathbf K=(K,I,\Lambda,\Omega)$ and $\mathbf L=(L,I_L,\Lambda_L,\Omega_L)$ be pre-$\HLO$-fields.  An {\bf embedding} 
$h \colon \mathbf K\to\mathbf L$
of  pre-$\HLO$-fields 
is an embedding $h\colon K\to L$ of pre-$H$-fields such that 
$$h(I)\ =\ h(K)\cap I_L, \qquad
h(\Lambda)\ =\ h(K)\cap\Lambda_L, \qquad h(\Omega)\ =\ h(K)\cap\Omega_L.$$
We also say that $\mathbf L$ is an {\bf extension} of $\mathbf K$ if
$L$ is a pre-$H$-field extension of $K$ with
$I=I_L\cap K$, $\Lambda=\Lambda_L\cap K$, and $\Omega=\Omega_L\cap K$ (so the natural inclusion $K\to L$ is an embedding $\mathbf K\to \mathbf L$ 
of  pre-$\HLO$-fields);
notation: $\mathbf K\subseteq \mathbf L$.
If $(L,I_L,\Lambda_L,\Omega_L)$ is a pre-$\HLO$-field and $K$ is a pre-$H$-subfield of $L$, then 
$$(K,I_L\cap K,\Lambda_L\cap K,\Omega_L\cap K)\ \subseteq\ (L,I_L,\Lambda_L,\Omega_L).$$
From Corollary~\ref{cor:ILO1}
we see that $K$ has a unique
expansion $(K, I,\Lambda,\Omega)$ to a pre-$\HLO$-field if and only if
one of the following conditions is satisfied:
\begin{enumerate}
\item[\textup{(i)}] $K$ is grounded;
\item[\textup{(ii)}] there exists $b\asymp 1$ in $K$ such that $v(b')$ is a gap in $K$;
\item[\textup{(iii)}] $K$ is $\upo$-free.
\end{enumerate}
From this equivalence we obtain: 

\begin{cor}\label{cor:unique HLO, embedding}
Let $\mathbf K=(K,\dots)$ and $\mathbf L=(L,\dots)$ be pre-$\HLO$-fields where~$K$ satisfies one of the conditions~\textup{(i), (ii), (iii)} above.
Then any embedding $K\to L$ of pre-$H$-fields is also an embedding of pre-$\HLO$-fields  $\mathbf K\to\mathbf L$.
\end{cor}

\noindent
By the remark following the proof of Corollary~\ref{cor:ILO2}, every pre-$\HLO$-field extends to some $\upo$-free newtonian Liouville closed $\HLO$-field. 
In the next section we establish a more precise result of this kind.

\subsection*{Notes and comments} For Liouville closed
$K$ the set $\operatorname{I}(K)$ lives in some sense in the asymptotic couple of~$K$, and reflects the
extra predicate needed to get
QE for such asymptotic couples in~\cite{AvdD}. Such an ``equivalence'' to 
a definable set in the asymptotic couple no longer exists
for~$\Upl(K)$ and~$\Upo(K)$.

\section{Embedding Pre-$\HLO$-Fields into $\upo$-Free $\HLO$-Fields}\label{sec:HLOemb}

\noindent
In this section we fix a pre-$\HLO$-field $\mathbf K=(K,I,\Lambda,\Omega)$, and construct an $\upo$-free $\HLO$-field extension of $\mathbf K$ with a useful semiuniversal property:

\begin{prop}\label{prop:upo ext of HLO}
There exists an $\upo$-free $\HLO$-field extension $\mathbf K^*$ of $\mathbf K$ such that $\res \mathbf K^*$ is algebraic over $\res \mathbf K$ and any embedding of $\mathbf K$ into a Schwarz closed $\HLO$-field~$\mathbf L$ extends to an embedding of~$\mathbf K^*$ into $\mathbf L$.
\end{prop}

\noindent
This result is contained in the next lemmas with their corollaries. 

\begin{lemma}\label{lem:HL 1}
Suppose $K$ is grounded, or there exists $b\asymp 1$ in $K$ such that $v(b')$ is a gap in $K$. Then $\mathbf K$ has an $\upo$-free $\HLO$-field extension~$\mathbf K^*$ 
such that any embedding of~$\mathbf K$ into a $\HLO$-field~$\mathbf L$ closed under logarithms extends to an embedding $\mathbf K^*\to \mathbf L$.
\end{lemma}
\begin{proof}
The $H$-field hull $F:=H(K)$ of $K$ is grounded by Corollary~\ref{cor:value group of dv(K)}(i)(iii).
Take $K^*=F_{\upo}$ and apply Lemma~\ref{Hlogupo} and Corollary~\ref{cor:unique HLO, embedding}.
\end{proof}
 
\begin{lemma}\label{lem:embed into omega-free, 1}
Suppose $K$ has gap $\beta$ and $v(b')\neq\beta$ for all $b\asymp 1$ in $K$. 
Then there exists a grounded pre-$\HLO$-field extension $\mathbf K_1$ of~$\mathbf K$ 
such that any embedding of $\mathbf K$ into a $\HLO$-field~$\mathbf L$ closed under integration extends to an embedding $\mathbf K_1\to \mathbf L$.
\end{lemma}
\begin{proof} Take $s\in K$ with 
$vs=\beta$. We distinguish two cases:

\case[1]{$s\notin I$.} Then $I=\I(K)=\{a\in K:va>\beta\}$ by Lemma~\ref{ILO3}. 
By Corollary~\ref{cor:value group of dv(K)}, $H(K)$ is an immediate extension of $K$, so $\beta$ remains a gap in $H(K)$.
Since $H(K)$ is $\d$-valued of $H$-type, 
Corollary~\ref{variant41, cor} yields an $H$-field extension
$K_1:=H(K)(y)$ of $K$ such that $y'=s$ and~$y\succ 1$. Then~$K_1$ is grounded, and thus
admits a unique expansion $\mathbf K_1=(K_1,I_1,\Lambda_1,\Omega_1)$ to a pre-$\HLO$-field. From
$I_1=\I(K_1)$ we get $s\notin I_1$, so $\mathbf K_1\supseteq\mathbf K$
by Lemma~\ref{ILO3}.
Let $\mathbf L\supseteq\mathbf K$ be a $\HLO$-field which is closed under integration.
Take $z\in L$ with $z'=s$. Then $z\succ 1$ since $s\notin I$.
The universal property of $H(K)$ and 
Corollary~\ref{variant41, cor} give a unique embedding $K_1=H(K)(y)\to L$ of pre-$H$-fields over $K$ sending $y$ to $z$, and by Corollary~\ref{cor:unique HLO, embedding} this is an embedding $\mathbf K_1\to\mathbf L$ of pre-$\HLO$-fields. 

\case[2]{$s\in I$.} Then $I=\{a\in K:va\ge\beta\}$ by Lemma~\ref{ILO3}. 
Let $K_1=K(y)$ be a pre-$H$-field extension of $K$ with $y'=s$ and $y\prec 1$, as in Corollary~\ref{extensionpdvfields1, cor} and the subsequent remarks. Then $K_1$ is grounded, and thus
admits a unique expansion $\mathbf K_1=(K_1,I_1,\Lambda_1,\Omega_1)$ to a pre-$\HLO$-field.
From $s\in \I(K_1)= I_1$ and Lemma~\ref{ILO3} we get 
$\mathbf K_1\supseteq\mathbf K$.
Let $\mathbf L\supseteq \mathbf K$ be a $\HLO$-field closed under integration. Take $z\in L$
with $z'=s$. From $s\in I$ we get $z\preceq 1$, and by subtracting a constant from $z$ we arrange~${z\prec 1}$. Then Corollary~\ref{extensionpdvfields1, cor} yields a unique embedding $K_1=K(y)\to L$
of pre-$H$-fields over $K$ sending $y$ to $z$. By 
Corollary~\ref{cor:unique HLO, embedding} this is an embedding $\mathbf K_1\to\mathbf L$ of pre-$\HLO$-fields.
\end{proof}

\noindent
These two lemmas and the constructions in their proofs
yield the following: 

\begin{cor}\label{cor:embed into omega-free, 1}
Suppose $K$ does not have asymptotic integration. 
Then $\mathbf K$ has an $\upo$-free $\HLO$-field extension 
$\mathbf K^*$  
such that $\res \mathbf K^*=\res \mathbf K$ and any embedding of~$\mathbf K$ into a $\HLO$-field~$\mathbf L$ closed under integration 
extends to an embedding $\mathbf K^*\to \mathbf L$.
\end{cor}

\noindent
The next two lemmas deal with the case where $K$ has asymptotic integration.

\begin{lemma}\label{lem:embed into omega-free, 2}
Assume $K$ has asymptotic integration and is not $\upl$-free.
Then $\mathbf K$ extends to an $\upo$-free $\HLO$-field
$\mathbf K^*$  such that 
$\res \mathbf K^*=(\res\mathbf K)^{\operatorname{rc}}$ and any embedding of $\mathbf K$ into a Liouville closed 
$\HLO$-field~$\mathbf L$ extends to an embedding 
$\mathbf K^*\to \mathbf L$.
\end{lemma}
\begin{proof}
By Corollary~\ref{cor:embed into omega-free, 1} it is enough to show: 
$\mathbf K$ has a $\HLO$-field extension~$\mathbf K_1$  with a gap 
such that $\res \mathbf K_1=(\res \mathbf K)^{\operatorname{rc}}$ and any embedding of $\mathbf K$ into a Liouville closed $\HLO$-field~$\mathbf L$ extends to an embedding $\mathbf K_1\to \mathbf L$.
By Lemmas~\ref{ILO7} and \ref{ILO8}, we have 
$\Upl(K)^{\downarrow} \ne K\setminus\Upd(K)^\uparrow$, and
the pre-$H$-field~$K$ has precisely two $\HLO$-cuts, $\big(\!\I(K),\Upl(K)^\downarrow,\dots\big)$ and $\big(\!\I(K),K\setminus\Upd(K)^\uparrow,\dots\big)$.
Let $E:=H(K)^{\operatorname{rc}}$. By Corollary~\ref{cor:value group of dv(K)}(i) we have
$\Gamma_{H(K)}=\Gamma$, so $\Gamma_E=\Q\Gamma$.
We distinguish two cases:

\case[1]{$E$ has a gap.} Take $s\in E^\times$ and $n\geq 1$ such that $vs$ is a gap in $E$ and $s^n\in K$. Then $s^{\dagger}=\frac{1}{n}(s^n)^\dagger\in K$. 
By Lemma~\ref{ILO3}, $E$ has exactly two $\HLO$-cuts $(I_1,\Lambda_1,\Omega_1)$ and $(I_2,\Lambda_2,\Omega_2)$, with 
$I_1=\I(E)=\{y\in E:y\prec s\}$ and $I_2=\{y\in E:y\preceq s\}$. 
We have $s\notin I_1$, so $-s^\dagger\in\Lambda_1\cap K$, and
$s\in I_2$, so $-s^\dagger\notin\Lambda_2\cap K$. If $-s^\dagger\in\Lambda$, then we set $\mathbf K_1:=(E,I_1,\Lambda_1,\Omega_1)$, and if $-s^\dagger\notin\Lambda$, then we set $\mathbf K_1:=(E,I_2,\Lambda_2,\Omega_2)$.
Then~$\mathbf K_1$ is an extension of the pre-$\HLO$-field $\mathbf K$.
Given an embedding $i\colon \mathbf K  \to \mathbf L$ into a
Liouville closed $\HLO$-field $\mathbf L$, there is a unique embedding $j\colon E\to L$
of $H$-fields
such that $j(a)=i(a)$ for all~$a\in K$, and it is easy to check that
$j$ is an embedding $\mathbf K_1\to\mathbf L$ of $\HLO$-fields. 

\case[2]{$E$ has no gap.}
Then $E$ has asymptotic integration, and the sequence~$(\upl_\rho)$ for~$K$ also serves for $E$.
Take $\upl\in K$ such that $\upl_\rho\leadsto\upl$. Then 
$-\upl$ creates a gap over~$E$ by Lemma~\ref{cg2}. Take an element $f\ne 0$ in some Liouville closed $H$-field extension of~$E$ such that
$f^\dagger=-\upl$. Then $vf$ is a gap in $E(f)$ by the remark following the proof of Lemma~\ref{cg2} with 
$\res E(f)=\res E$. Using Lemma~\ref{dv} it follows that 
the pre-$H$-field
 $E(f)$ is actually an $H$-field. 
 Therefore, by Lemma~\ref{ILO3}, $E(f)$ has exactly two $\HLO$-cuts $(I_1,\Lambda_1,\Omega_1)$ and $(I_2,\Lambda_2,\Omega_2)$, with 
$$I_1\ =\ \I\!\big(E(f)\big)\ =\ \big\{y\in E(f)\ :y\prec f\big\},\qquad I_2\ =\ \big\{y\in E(f):\ y\preceq f\big\}.$$
We have $f\notin I_1$, so $\upl=-f^\dagger\in\Lambda_1\cap K$, and $f\in I_2$, so $\upl=-f^\dagger\notin\Lambda_2\cap K$.
If $\upl\in\Lambda$, then we set $\mathbf K_1:=\big(E(f),I_1,\Lambda_1,\Omega_1\big)$, and if $\upl\notin\Lambda$, then we
set $\mathbf K_1:=\big(E(f),I_2,\Lambda_2,\Omega_2\big)$. In any case, we have 
$\mathbf K_1\supseteq \mathbf K$. Let $i\colon \mathbf K  \to \mathbf L$ be an
embedding into a
Liouville closed $\HLO$-field $\mathbf L$.  
By Lemma~\ref{lem:active S} the set 
$$S\ :=\ \big\{v(\upl+a^\dagger):\ a\in E^\times\big\}$$
is a cofinal subset of $\Psi_{E}^\downarrow$ and $f$ is transcendental
over $E$. 
Then Lemmas~\ref{newprop, lemma 1}
and \ref{prehexpint} provide an embedding $j\colon E(f)\to L$ of $H$-fields with $j(a)=i(a)$ for all~$a\in K$, and
any such embedding is an embedding $\mathbf K_1\to\mathbf L$ of $\HLO$-fields.
\end{proof}

\begin{lemma}\label{lem:embed into omega-free, 3}
Suppose $K$ is $\upl$-free but not $\upo$-free. Then $\mathbf K$ has an $\upo$-free $\HLO$-field extension $\mathbf K^*$  
such that 
$\res \mathbf K^*$ is algebraic over $\res \mathbf K$ and any embedding of $\mathbf K$ into a Schwarz closed 
$\HLO$-field~$\mathbf L$ extends to an embedding 
of $\mathbf K^*$ into $\mathbf L$.
\end{lemma}
\begin{proof}
Take $\upo\in K$ such that $\upo_\rho\leadsto\upo$, so $\omega\big(\Upl(K)\big){}^\downarrow<\upo<\sigma\big(\Upg(K)\big){}^\uparrow$.
By Lem\-ma~\ref{ILO6}, there are exactly two $\HLO$-cuts in $K$:
$$\big(\!\I(K), \Upl(K)^{\downarrow}, \omega(K)^{\downarrow}\big),\qquad \big(\!\I(K), \Upl(K)^{\downarrow}, K\setminus \sigma\big(\Upg(K)\big){}^\uparrow\big).$$
Since $\upo\in K\setminus \sigma\big(\Upg(K)\big){}^\uparrow$, it follows from the proof of Lem\-ma~\ref{ILO6} that $\upo\notin \omega(K)^\downarrow$. 
We distinguish two cases:

\case[1]{$\Omega=\omega(K)^{\downarrow}$.}
Lemma~\ref{cor:sigma(upg)=upo} yields a pre-$H$-field extension 
$K_\upg:=K\<\upg\>$ of~$K$ such that $\sigma(\upg)=\upo$, $\upg>0$, $v\upg$ is gap in $K_\upg$, and $\res(K_\upg)=\res(K)$. By 
remark~(3) following the proof of Proposition~\ref{prop:sigma(upg)=upo}
there is no $b\asymp 1$ in $K_{\upg}$ with $b'\asymp \upg$, so
by Lemma~\ref{ILO3}  we have exactly two $\HLO$-cuts $(I_1,\Lambda_1,\Omega_1)$,
$(I_2,\Lambda_2,\Omega_2)$
in $K_\upg$, where 
$$I_1\ =\ \I(K_{\upg})\ =\ \{y\in K_\upg:\  y\prec\upg\}, \qquad I_2\ =\ \{y\in K_\upg:\  y\preceq\upg\}.$$
Put $\mathbf K_\upg := (K_\upg, I_1,\Lambda_1,\Omega_1)$.
We have $\upg\notin I_1$, so $\upo=\sigma(\upg)\notin \Omega_1$ by ($\HLO 12$), and
thus $\mathbf K_\upg\supseteq\mathbf K$.  Let $\mathbf K^*$ be an $\upo$-free $\HLO$-field extension of $\mathbf K_\upg$
obtained by applying Corollary~\ref{cor:embed into omega-free, 1} to $\mathbf K_\upg$ instead of $\mathbf K$. 
Let $\mathbf L$ be a Schwarz closed $\HLO$-field extension of $\mathbf K$;
we claim that there is an embedding $\mathbf K^*\to\mathbf L$ over $K$.
Now $\upo\notin{\Omega=\Upo(L)\cap K}$, hence $\upo\in\sigma\big(\Upg(L)\big)$.
By Lemma~\ref{lem:sigma increasing}, the restriction of~$\sigma$ to~$\Upg(L)$ is strictly increasing. Let $\upg^*$ be the unique element of $\Upg(L)$ such that $\sigma(\upg^*)=\upo$. From $\Upg(L) = L^> \setminus\I(L)$ we get $\upg^*>0$ and $\upg^*\notin\I(L)$, and thus 
$v\upg^*<(\Gamma^>)'$. Also $\sigma(\upg^*)=\upo < \sigma\big(\Upg(K)\big)$ gives
$0 < \upg^*<  \Upg(K)$, and so $\Psi < v\upg^*$.  
Thus Proposition~\ref{prop:sigma(upg)=upo} and Lemma~\ref{cor:sigma(upg)=upo} yield an embedding $h\colon K_\upg\to L$
of pre-$H$-fields over~$K$ with $h(\upg)=\upg^*$. 
Since $\upg\notin I_1$ and $\upg^*\notin \I(L)$, $h$ is an embedding
$\bf K_{\upg}\to \bf L$ of pre-$\HLO$-fields. By Corollary~\ref{cor:embed into omega-free, 1} we can extend $h$
to an embedding~$\mathbf K^*\to\mathbf L$. 

\case[2]{$\Omega=K\setminus \sigma\big(\Upg(K)\big){}^\uparrow$.}
Corollary~\ref{cor:adjoin upl} and Lemma~\ref{prehimm} yield an immediate pre-$H$-field
extension $K_\upl:=K(\upl)$ of $K$ with $\upl_{\rho} \leadsto \upl$ and
$\omega(\upl)=\upo$. Then $K_\upl$ has rational asymptotic integration and is not $\upl$-free, and $\Upl(K_\upl)< \upl < \Upd(K_{\upl})$,  so by
Lemma~\ref{ILO7} there are exactly two $\HLO$-cuts in $K_\upl$: 
$$\big(\!\I(K_\upl), \Upl(K_\upl)^{\downarrow}, K_\upl\setminus \sigma\big(\Upg(K_\upl)\big){}^{\uparrow}\big)\  \text{ and }\ \big(\!\I(K_\upl), K_\upl\setminus \Upd(K_\upl)^{\uparrow}, K_\upl\setminus \sigma\big(\Upg(K_\upl)\big){}^{\uparrow} \big).$$
Note that $\upo\in \Omega$ as well as $\upo\in \omega(K_{\upl})\subseteq K_\upl \setminus \sigma\big(\Upg(K_\upl)\big){}^\uparrow$. Therefore, setting
$$\mathbf K_\upl\ :=\ \big(K_\upl,\I(K_\upl),K_\upl\setminus \Upd(K_\upl)^{\uparrow}, K_\upl\setminus \sigma\big(\Upg(K_\upl)\big){}^{\uparrow}\big),$$
we get $\mathbf K_\upl\supseteq\mathbf K$. Let $\mathbf K^*$ be an $\upo$-free $\HLO$-field extension of $\mathbf K_\upl$
obtained by applying Lemma~\ref{lem:embed into omega-free, 2} to $\mathbf K_\upl$ instead of $\mathbf K$.
Let $\mathbf L$ be a Schwarz closed $\HLO$-field extension of $\mathbf K$;
we claim that there is an embedding $\mathbf K^*\to\mathbf L$ over $K$.
Recall from
Corollary~\ref{cor:Upl 2} and Proposition~\ref{prop:omega for H-fields}
that $\omega$
is strictly increasing on $\Upl(L)$ and $\Upl(L)< \Upd(L)$. Since $\upo\in \Omega=\omega\big(\Upl(L)\big)\cap K$, we get a unique
$\upl^*\in \Upl(L)$ such that $\omega(\upl^*)=\upo$. Then
$\Upl(K)<\upl^*<\Upd(K)$, so $\upl_{\rho}\leadsto \upl^*$ and
Corollary~\ref{cor:adjoin upl} yields an embedding $h\colon K_\upl\to L$
of pre-$H$-fields over $K$ with $h(\upl)=\upl^*$. 
Since $\upl\notin\Upl(K_\upl)^\downarrow$ and $\upl^*\in \Upl(L)$, $h$ is an embedding  $\mathbf K_{\upl} \to \mathbf L$ of pre-$\HLO$-fields, and so by Lemma~\ref{lem:embed into omega-free, 2}, $h$ extends to an embedding $\mathbf K^*\to\mathbf L$.
\end{proof}

\begin{lemma}\label{lem:embed into omega-free, 4} If $K$ is $\upo$-free, then $\mathbf K$ has an $\upo$-free $\HLO$-field extension $\mathbf K^*$  
such that any embedding of $\mathbf K$ into a $\HLO$-field~$\mathbf L$ extends to an embedding 
of $\mathbf K^*$ into $\mathbf L$.
\end{lemma}
\begin{proof} Assume $K$ is $\upo$-free. 
Then $H(K)$ is $\upo$-free by Theorem~\ref{upoalgebraic}. Let
$\mathbf K^*$ be an expansion of $H(K)$ to a $\HLO$-field. Then 
$\mathbf K \subseteq \mathbf K^*$ by Corollary~\ref{cor:ILO1}.
It remains to use the universal property of $H(K)$ and Corollary~\ref{cor:unique HLO, embedding}.
\end{proof}

\noindent
Corollary~\ref{cor:embed into omega-free, 1} and Lemmas~\ref{lem:embed into omega-free, 2},~\ref{lem:embed into omega-free, 3}, and~\ref{lem:embed into omega-free, 4} now have Proposition~\ref{prop:upo ext of HLO} as an immediate consequence. Note: this proof yields 
an extension $\mathbf K^*$ of $\mathbf K$ as in 
Proposition~\ref{prop:upo ext of HLO} that is $\d$-algebraic over $\mathbf K$.

\subsection*{The Newton-Liouville closure of a pre-$\HLO$-field}  
Here we extend the results on Newton-Liouville closures of 
$\upo$-free $H$-fields to pre-$\HLO$-fields. 

\begin{prop}\label{HLOnelicl} Let $\mathbf K=(K,I,\Lambda,\Omega)$ be a 
pre-$\HLO$-field.
Then $\mathbf K$ has an $\upo$-free newtonian Liouville closed $\HLO$-field
extension ${\mathbf K}^{\operatorname{nl}}$ that embeds over $\mathbf K$ into any
$\upo$-free newtonian Liouville closed $\HLO$-field
extension of $\mathbf K$.
\end{prop}
\begin{proof} First we take an $\upo$-free $\HLO$-field extension $\mathbf K^*=(K^*,\dots)$ of $\mathbf K$ as in Proposition~\ref{prop:upo ext of HLO}.
Next we take the Newton-Liouville closure $E^{\operatorname{nl}}$ of the $\upo$-free $H$-field 
$E:= K^*$. Then the unique expansion of 
$E^{\operatorname{nl}}$ to a $\HLO$-field is an extension~${\mathbf K}^{\operatorname{nl}}$ of $\mathbf K$ as claimed.
\end{proof}

\noindent
We define a {\bf Newton-Liouville closure\/} of a pre-$\HLO$-field 
$\mathbf K$ to be an extension~${\mathbf K}^{\operatorname{nl}}$ as in 
Proposition~\ref{HLOnelicl}. Thus every pre-$\HLO$-field has a Newton-Liouville closure.
 
\index{closure!Newton-Liouville}
\index{Newton-Liouville closure}
\index{LambdaOmega-field@$\HLO$-field!Newton-Liouville closure}

\begin{prop} Let $\mathbf K$ be a pre-$\HLO$-field. Any two Newton-Liouville closures of~$\mathbf K$ are isomorphic
over~$\mathbf K$. If ${\mathbf K}^{\operatorname{nl}}$ is a Newton-Liouville closure
of $\mathbf K$, then~${\mathbf K}^{\operatorname{nl}}$ does not have any 
proper newtonian $\upo$-free Liouville closed $\HLO$-subfield 
containing~$\mathbf K$ as a substructure.
\end{prop}
\begin{proof} Let ${\mathbf K}^{\operatorname{nl}}$ be the
Newton-Liouville closure of $\mathbf K$ constructed in the proof of Proposition~\ref{HLOnelicl}. Then 
${\mathbf K}^{\operatorname{nl}}$ is $\d$-algebraic over 
$\mathbf K$ and the residue field
of ${\mathbf K}^{\operatorname{nl}}$ is a real closure of
$\res \mathbf K$. Let $\mathbf L$ be any Newton-Liouville closure of $\mathbf K$. Embedding $\mathbf L$ into ${\mathbf K}^{\operatorname{nl}}$ over $\mathbf K$,
we see that $\mathbf L$ is $\d$-algebraic over $\mathbf K$ and
its residue field is a real closure of $\res \mathbf K$.
Consider any embedding 
  $i\colon {\mathbf K}^{\operatorname{nl}}\to \mathbf L$
  over $\mathbf K$. Then 
  $i\big({\mathbf K}^{\operatorname{nl}}\big)= \mathbf L$ by Theorem~\ref{th:noextL}. This proves the first part, and the minimality property of ${\mathbf K}^{\operatorname{nl}}$ is likewise a consequence of 
Theorem~\ref{th:noextL}. 
\end{proof}

\section{The Language of $\HLO$-Fields}\label{sec:lHLO}

\noindent
In the introduction to this chapter we specified the language
$$\mathcal L\ :=\ \{0,\ 1,\  +,\  -,\ \cdot\ ,\ \der,\ \le,\ \preceq\}$$
of ordered valued differential rings. Each ordered valued differential field is viewed as an $\mathcal L$-structure
in the natural way. In this section we show that Theorem~\ref{thm:qe} 
fails if we drop either the
symbol $\Upl$ or the symbol $\Upo$ from the language $\mathcal{L}^{\iota}_{\Upl\Upo}$ of $\HLO$-fields. (We 
prove somewhat sharper versions of this fact.) 

Throughout this section $K$ is a pre-$H$-field. In Section~\ref{sec:special sets} we introduced the special
subsets
$$%\begin{equation}\label{eq:special sets}
\I(K),\quad \Upg(K),\quad \Upl(K),\quad \Upd(K), \quad \omega(K), \quad \sigma\big(\Upg(K)\big)
$$%\end{equation}
of  $K$. If $K$ is Schwarz closed, then each of these sets is clearly existentially
definable as well as universally definable in the $\mathcal L$-structure 
$K$, both forms of definability holding without parameters from $K$
and witnessed by $\mathcal{L}$-formulas independent of~$K$. 
In this section we successively investigate the quantifier-free definability
of these sets in a Schwarz closed $K$, in the  language~$\mathcal L$ 
and some extensions of this language with predicates for some
of these sets. Our first result shows that no $\upo$-free real closed
$H$-field eliminates quantifiers in $\mathcal L$:

\begin{prop}\label{prop:I(K) not qf definable} Suppose $K$ is an $\upo$-free real closed $H$-field. Then
the subset~$\I(K)$ of $K$ is not quantifier-free definable \textup{(}with
parameters\textup{)} in the $\mathcal L$-structure~$K$. 
\end{prop}
\begin{proof}
Take an element $\ell$ in an elementary extension $K^*$ of $K$ with 
$\ell >0$ and
$1\prec \ell\prec \ell_\rho$ for all $\rho$, and set $\upg:=\ell^\dagger$, $\upl:=-\upg^\dagger$. Then $\upl$ and $\upl+\upg$ are pseudolimits of~$(\upl_\rho)$, by Corollary~\ref{cor:gap}, and  
$(\upl_\rho)$ is of $\d$-transcendental type over $K$, by Corollary~\ref{differentialtranscendence}. Hence by Lemma~\ref{zdf, newton} and Corollary~\ref{zmindifpol},
the pre-$H$-subfields~$K\<\upl\>$ and $K\<\upl+\upg\>$ of $K^*$ are immediate extensions of $K$ (so they are $H$-fields), and we have an isomorphism $K\<\upl\>\to K\<\upl+\upg\>$ of $H$-fields over $K$ sending~$\upl$ to~$\upl+\upg$.
By Lemma~\ref{cg2}, the element $-\upl$ creates a gap over $K\<\upl\>$. Likewise, $-(\upl + \upg)$ creates a gap
over $K\<\upl + \upg\>$.  
Let $f:=(1/\ell)^\dagger=-\upg$ and $g:=(1/\ell)'=-\upg/\ell$, so
$f<0$ and $g<0$ (using $\ell>0$). Then $f^\dagger=-\upl$ and $g^\dagger=-(\upl+\upg)$,  so
the above isomorphism $K\<\upl\>\to K\<\upl+\upg\>$ extends by Lemma~\ref{lem:active S} and the uniqueness parts of Lemmas~\ref{newprop, lemma 1} and~\ref{prehexpint} to an isomorphism $K\<\upl,f\>\to K\<\upl+\upg,g\>$ of $\mathcal L$-structures sending $f$ to $g$. Now, if $\I(K)$ were defined in $K$ by a quantifier-free formula $\varphi(y)$ in the language $\mathcal L$ augmented by names for the elements of $K$, then we would have $K^*\models\neg\varphi(f)$ and $K^*\models\varphi(g)$, and so $K\<\upl,f\>\models\neg\varphi(f)$ and $K\<\upl+\upg,g\>\models\varphi(g)$,  violating the above isomorphism between $K\<\upl,f\>$ and $K\<\upl+\upg,g\>$.
\end{proof}

\noindent
We extend $\mathcal L$ by a single unary function symbol $\iota$ 
to the language $\mathcal L^{\iota}$. Any ordered valued differential field $F$ will be construed as an $\mathcal{L}^{\iota}$-structure by interpreting
this new function symbol as the function $F\to F$ that
agrees with $f\mapsto f^{-1}$ on $F^\times$ and sends~$0$ to~$0$. 
Thus the underlying ring of an $\mathcal L^{\iota}$-substructure of an ordered valued differential field is a field.
Passing from $\mathcal L$ to $\mathcal L^{\iota}$ does not increase what we can express quantifier-free in ordered valued differential fields:

\begin{cor}\label{lem:L vs L'}
Given any quantifier-free $\mathcal L^{\iota}$-formula $\varphi^{\iota}(x_1,\dots,x_n)$, there is a 
quanti\-fier-free $\mathcal L$-formula $\varphi(x_1,\dots,x_n)$ 
which in every
ordered valued differential field $F$ defines the same subset of $F^n$ as $\varphi^{\iota}(x_1,\dots,x_n)$.
\end{cor}
\begin{proof}
Let $\operatorname{OVD}^{\iota}$ be the $\mathcal L^{\iota}$-theory of ordered valued differential fields. By \ref{closureprop} it is enough to show that $\operatorname{OVD}^{\iota}$ has closures of $\mathcal L$-substructures.  
Let $E,F\models \operatorname{OVD}^{\iota}$ have a common 
$\mathcal L$-substructure $D$. Thus $D$ is an ordered subring of both~$E$ and~$F$ with a derivation on it that agrees on $D$ with the derivations of $E$ and $F$ such that
for all $f,g\in D$ we have $f\preceq_E g\Longleftrightarrow f\preceq_D g\Longleftrightarrow f\preceq_F g$, where $\preceq_D$, $\preceq_E$, $\preceq_F$ are the
interpretations of the symbol $\preceq$ of $\mathcal L$ in $D$, $E$, $F$, respectively.
Let~$K_E$ and~$K_F$ be the fraction fields of the integral domain~$D$ in~$E$
and~$F$, respectively. Then $K_E$ is the underlying ring of an $\mathcal L^{\iota}$-substructure of $E$, to be denoted also 
by~$K_E$. Likewise,~$K_F$ denotes the corresponding $\mathcal L^{\iota}$-substructure of $F$.
The unique field isomorphism $K_E \to K_F$ over $D$ is clearly an $\mathcal L^{\iota}$-isomorphism. 
\end{proof}

\noindent
Thus if $K$ is an $\upo$-free real closed $H$-field, then $\I(K)$ is  not quantifier-free definable (with parameters) in the $\mathcal L^{\iota}$-structure $K$. 

Let $\mathcal L^{\iota}_{\Upl}$ be the language $\mathcal L^{\iota}$ augmented by unary predicate symbols $\I$ and $\Upl$.
Given a $\HL$-cut $(I,\Lambda)$ in $K$, we have the
$\mathcal L^{\iota}_{\Upl}$-structure
$(K,I,\Lambda)$: interpret~$\I$ and~$\Upl$ by~$I$ and~$\Lambda$.
Recall that if $K$ is $\upl$-free, then there is only one $\HL$-cut in~$K$.
By Lemma~\ref{lem:Upl, 1}, if~$K$ is Liouville closed, then
 $\I(K)$ is quantifier-free definable in $(K,\Upl(K))$, with $K$ construed as an $\mathcal L^{\iota}$-structure;
nevertheless, we include the symbol $\I$ in $\mathcal L^{\iota}_{\Upl}$.
Note that if $K$ is Liouville closed, then $\Upd(K)=K\setminus\Upl(K)$ is quantifier-free definable
in~the $\mathcal L^{\iota}_{\Upl}$-structure $(K,\I(K),\Upl(K))$ as well.
However:

\begin{prop}\label{prop:omega(K) not qf definable}
Suppose $K$ is an $\upo$-free real closed $H$-field. Then the subsets~$\omega(K)$ and~$\omega(K)^{\downarrow}$ of $K$ are not quantifier-free definable \textup{(}even allowing parameters\textup{)} in the $\mathcal L^{\iota}_{\Upl}$-structure $(K,\I(K),\Upl(K)^\downarrow)$.
\end{prop}
\begin{proof}
Take $K^*$  and $\ell$  as in the proof of Proposition~\ref{prop:I(K) not qf definable}, 
and set $\upg:=\ell^\dagger$, $\upl:=-\upg^\dagger$, and $\upo:=\omega(\upl)$. Then $\upl_\rho\leadsto\upl$  by Corollary~\ref{cor:gap}
and hence $\upo_\rho\leadsto \upo$ by Corollary~\ref{pczrho}.
In view of  Corollary~\ref{pcrhocor} and $v(\upg^2) > 2\Psi$, the pc-sequences~$(\upo_\rho)$ and $(\upo_\rho+\upg_\rho^2)$ in $K$ are equivalent, and $\sigma(\upg_\rho)=\upo_\rho+\upg_\rho^2\leadsto \sigma(\upg)=\upo+\upg^2$. 
By Corollary~\ref{differentialtranscendence}, $(\upo_\rho)$ is of $\d$-transcendental type over $K$. Hence by Lemma~\ref{zdf, newton} and Corollary~\ref{zmindifpol},
the pre-$H$-subfields $K\<\upo\>$ and $K\<\upo+\upg^2\>$ of $K^*$ are immediate extensions of $K$ (so they are $H$-fields) and we have an isomorphism $$K\<\upo\>\ \to\  K\<\upo+\upg^2\>$$ 
of $H$-fields over $K$ sending $\upo$ to $\upo+\upg^2$.
By Proposition~\ref{cor:adding upo preserves upl-free}, $K\<\upo\>$ is $\upl$-free, and hence by Lemma~\ref{IL4}, $\big(\!\I(K\<\upo\>), \Upl(K\<\upo\>)^\downarrow\big)$ is the unique $\HL$-cut in $K\<\upo\>$, and
likewise for $K\<\upo+\upg^2\>$ instead of $K\<\upo\>$.  Thus 
$$\Upl(K^*)^\downarrow\cap K\<\upo\>\ =\ \Upl\big(K\<\upo\>\big)^\downarrow, \qquad 
 \Upl(K^*)^\downarrow\cap K\<\upo+\upg^2\>\ =\  \Upl\big(K\<\upo+\upg^2\>\big)^\downarrow,$$ so 
 our isomorphism $K\langle \upo\rangle\to K\langle \upo+\upg^2\rangle$ is an isomorphism between
$\mathcal L^{\iota}_{\Upl}$-substructures of~$(K^*,\I(K^*), \Upl(K^*)^\downarrow)$.  Now $\upo\in\omega(K^*)$ and $\upo+\upg^2\in \sigma\big(\Upg(K^*)\big)$, so $\upo+\upg^2\notin\omega(K^*)$, and thus
$\omega(K)$ is not quantifier-free definable (with parameters) in the $\mathcal L^{\iota}_{\Upl}$-structure~$(K,\I(K),\Upl(K)^\downarrow)$.
Likewise with $\omega(K)^{\downarrow}$ instead of $\omega(K)$.
\end{proof}

\noindent
Let $\mathcal L_{\Upo}^\iota$ be the language  $\mathcal L^{\iota}$ augmented 
by the unary predicate symbols $\I$ and $\Upo$. 
Then we have the following analogue of Proposition~\ref{prop:omega(K) not qf definable}:

\begin{prop}\label{prop:upoversion} 
Suppose $K$ is an $\upo$-free real closed $H$-field. Then the subsets~$\Upl(K)$ and~$\Upl(K)^{\downarrow}$ of $K$ are not quantifier-free definable \textup{(}even allowing parameters\textup{)} in the 
$\mathcal L^{\iota}_{\Upo}$-structure $(K,\I(K),\omega(K)^\downarrow)$.
\end{prop}

{\sloppy

\begin{proof}
Again, take $K^*$  and $\ell$  as in the proof of Proposition~\ref{prop:I(K) not qf definable}, 
and set $\upg:=\ell^\dagger$ and $\upl:=-\upg^\dagger$. Then  $K\<\upl\>$ and $K\<\upl+\upg\>$ are immediate $H$-field extensions of~$K$, and  we have an
isomorphism $K\<\upl\>\to K\<\upl+\upg\>$ of $H$-fields over $K$ which sends~$\upl$ to~$\upl+\upg$.
Since $K$ has asymptotic integration with divisible value group, so 
does~$K\<\upl\>$. As~$K\<\upl\>$ is not $\upl$-free, 
 there are by Lemma~\ref{ILO7} exactly two $\HLO$-cuts in~$K\<\upl\>$, 
$$\big(I, \Lambda_1, \Omega \big)\quad\text{and}\quad\big(I, \Lambda_2,\Omega  \big),$$
the key point being that these two $\HLO$-cuts have the same first
component $I$ and same third component $\Omega$. In particular,
$$\I(K^*)\cap K\<\upl\>\  =\  I,\qquad 
\omega(K^*)^\downarrow\cap K\<\upl\>\ =\ \Omega.$$
Likewise with $K\<\upl+\upg\>$ in place of $K\<\upl\>$. Hence our iso\-mor\-phism~${K\<\upl\>}\to {K\<\upl+\upg\>}$ is an isomorphism
between $\mathcal{L}^{\iota}_{\Upo}$-substructures of $\big(K^*, \I(K^*), \omega(K^*)^\downarrow\big)$. But $\upl=-\ell^{\dagger\dagger}\in \Upl(K^*)$ and
$\upl+\upg=-(1/\ell)^{\prime\dagger}\in \Upd(K^*)$, so 
$\upl+\upg\notin \Upl(K^*)$, and thus~$\Upl(K)$ is not quantifier-free definable (even allowing parameters) in the 
$\mathcal L^{\iota}_{\Upo}$-struc\-ture $(K,\I(K),\omega(K)^\downarrow)$. 
Similarly with $\Upl(K)^\downarrow$ in place of $\Upl(K)$.
\end{proof}
}

\noindent
Thus $\upo$-free newtonian Liouville closed $H$-fields 
do not eliminate quantifiers when viewed in the usual way
either as $\mathcal{L}^{\iota}_{\Upl}$-structures or as 
$\mathcal{L}^{\iota}_{\Upo}$-structures. This goes to explain
our choice of language $\mathcal{L}^{\iota}_{\Upl\Upo}$ (see
the introduction to this chapter). Could a mild ``algebraic'' extension of $\mathcal{L}^{\iota}$, for example by a square root function, allow us to drop one of $\Upl$, $\Upo$? To eliminate this possibility and
make our choice of language more compelling, we
now indicate some stronger negative results. Towards
this end we specify a language $\mathcal{L}^{\a}$ that serves
as a more robust
version of $\mathcal{L}^{\iota}$.
%that will also be useful in deriving in Section~\ref{sec:con} part %(2) of Corollary~\ref{cor:qe} from the Main Theorem~\ref{thm:qe}.    

To define $\mathcal{L}^{\a}$, note that $\mathcal{L}$ has the language of ordered rings as a sublanguage. We consider $\R$ as a structure for the language of ordered rings in the usual way. A function~${\R^n\to \R}$ is said to be {\bf $\Q$-semialgebraic\/} if its graph is defined in the structure $\R$ by a (quantifier-free) formula in the language of ordered rings; we do not allow names for arbitrary real numbers in such formulas. We extend $\mathcal{L}$ to the language $\mathcal{L}^{\a}$ by adding for each $\Q$-semialgebraic function $f\colon \R^n \to \R$ an $n$-ary function symbol~$f$. We construe any real closed valued differential
field $E$ as an $\mathcal{L}^{\a}$-structure by interpreting such~$f$ as the function $E^n\to E$ whose graph is defined in $E$ by any formula in the language of ordered rings
that defines the graph of $f$ in $\R$. For example, the
function $a\mapsto a^{-1}\colon E \to E$ (with $0^{-1}:= 0$ by convention) is named by a function symbol~$\iota$ of~$\mathcal{L}^{\a}$. For each integer $d\ge 1$, the function
$y\mapsto y^{1/d}\colon E \to E$, taking the value~$0$ for~${y\le 0}$ by convention, is also named by a function symbol. With the richer language~$\mathcal{L}^{\a}$ replacing
$\mathcal{L}^{\iota}$ the above results go through. 
For example, Corollary~\ref{lem:L vs L'} extends as follows:

\begin{cor}\label{lem:L vs L^a}
Given any quantifier-free $\mathcal L^{\a}$-formula $\varphi^{\a}(x_1,\dots,x_n)$, there is a 
quanti\-fier-free $\mathcal L$-formula $\varphi(x_1,\dots,x_n)$ 
which defines in every real closed $K$ the same subset of $K^n$ as
$\varphi^{\a}(x_1,\dots,x_n)$. \textup{(}Recall that $K$ ranges over pre-$H$-fields.\textup{)}
\end{cor}

\noindent
(In the proof of Corollary~\ref{lem:L vs L'}, replace 
$\operatorname{OVD}^{\iota}$ by
the $\mathcal L^{\a}$-theory of real closed
ordered valued differential fields whose valuation ring is convex,
use Corollary~\ref{cvrc}, and
take real closures of fraction fields instead of just fraction fields.)

Thus if $K$ is any $\upo$-free real closed $H$-field, 
then $\I(K)$ is still not quantifier-free definable (with parameters) in the $\mathcal L^{\a}$-structure $K$. 

\medskip\noindent
Let
$\mathcal L^{\a}_{\Upl}$ be the language $\mathcal L^{\a}$ augmented
by unary relation symbols $\I$ and $\Upl$.
Then Proposition~\ref{prop:omega(K) not qf definable} goes through with $\mathcal L^{\a}_{\Upl}$
instead of $\mathcal L^{\iota}_{\Upl}$: replace in the proof of that proposition $K\<\upo\>$ and $K\<\upo+\upg^2\>$
by their real closures in $K^*$, and use that these real closures are 
immediate extensions of $K$ and $\mathcal{L}^{\a}$-substructures
of $K^*$.

Likewise, let
$\mathcal L^{\a}_{\Upo}$ be the language $\mathcal L^{\a}$ augmented
by unary relation symbols~$\I$ and~$\Upo$.
Then Proposition~\ref{prop:upoversion} goes through with 
$\mathcal L^{\a}_{\Upo}$
instead of $\mathcal L^{\iota}_{\Upo}$: replace in the proof of that
proposition $K\<\upl\>$ and $K\<\upl+\upg\>$ by their real closures
in $K^*$.

\subsection*{Notes and comments} For $K=\mathbb T$, Proposition~\ref{prop:I(K) not qf definable} and the $\mathcal{L}^{\a}_{\Upl}$-variant of
Proposition~\ref{prop:omega(K) not qf definable} are in \cite{ADHOleron},
with slightly different notation: $\mathcal{L}'$ instead of our
$\mathcal{L}^{\a}$. Corollary~\ref{lem:L vs L^a} also occurs there as
Proposition~5.5, but with a defective proof.

\section{Elimination of Quantifiers with Applications}\label{sec:embth}

\noindent
In the introduction to this chapter we defined the theory $T^{\operatorname{nl},\iota}_{\Upl\Upo}$. Its models
are exactly the $\upo$-free newtonian Liouville closed $\HLO$-fields;
we defined $\HLO$-fields at the end of Section~\ref{sec:LO-cuts}. As noted at the end of that section, the substructures of models of $T^{\operatorname{nl},\iota}_{\Upl\Upo}$ are exactly the pre-$\HLO$-fields.
Thus by the embedding criterion~\ref{prop:QE test, 10} for~QE, Theorem~\ref{thm:qe} is a consequence of the following: 

\begin{theorem}\label{thm:embqe} Let $\mathbf K$ and 
$\mathbf L$ be
$\upo$-free newtonian Liouville closed $\HLO$-fields such that
$\mathbf L$ is $\kappa^{+}$-saturated,
where $\kappa$ is the cardinality of the underlying set of 
$\mathbf K$. Let $\mathbf E$ be a substructure of $\mathbf K$
and let $i\colon \mathbf E \to \mathbf L$ be an embedding. Then $i$ can be
extended to an embedding $\mathbf K \to \mathbf L$.
$$\xymatrix{ \mathbf K \ar@{-->}[r] & \mathbf L \\
\mathbf E \ar[u]^{\subseteq} \ar[ur]^i }$$
\end{theorem}  
 
\noindent
As to the proof of Theorem~\ref{thm:embqe}, note that by Proposition~\ref{prop:upo ext of HLO} we can reduce to the case that $\mathbf E$ is an $\upo$-free $\HLO$-field. In view of Corollary~\ref{cor:unique HLO, embedding} this case is taken care of by Corollary~\ref{cor:upoembnl}.

\begin{cor}
$T^{\operatorname{nl},\iota}_{\Upl\Upo}$ is the model completion of the $\mathcal L^{\operatorname{nl},\iota}_{\Upl\Upo}$-theory of $\Upl\Upo$-fields. 
%viewed as
%$\mathcal L^{\operatorname{nl},\iota}_{\Upl\Upo}$-structures.
\end{cor}

\noindent
Theorem~\ref{thm:qe} has some immediate logical consequences for $T^{\operatorname{nl}}$:

\begin{cor}\label{cor:completionsdecidable} The completions of $T^{\operatorname{nl}}$ are the  $\mathcal{L}$-theories $T^{\operatorname{nl}}_{\operatorname{small}}$ and $T^{\operatorname{nl}}_{\operatorname{large}}$. These two theories as well as $T^{\operatorname{nl}}$ itself are decidable. 
\end{cor}
\begin{proof} 
Consider the Hardy field $E=\Q(x)$ ($x>\Q$, $x'=1$). Any Liouville closed $H$-field $K$ with small derivation has an element $f>C$ with $f'=1$, and this yields an embedding $E \to K$ sending $x$ to $f$. Note also that $E$ is grounded. In view of Corollary~\ref{cor:unique HLO, embedding}, Theorem~\ref{thm:qe} and a well-known completeness criterion  (Corollary~\ref{cor:QE and sentences}) the completeness of $T^{\operatorname{nl}}_{\operatorname{small}}$ then follows. 

Next, set $a:= x^{-2}$ and consider the compositional conjugate $E^a$ of $E$. Its derivation $\derdelta=x^2\der$ is not small, since $x^{-1}\prec 1$ and $\derdelta(x^{-1})=-1\asymp 1$. Let $K$ be any Liouville closed $H$-field whose derivation is not small. Take  $f\in K$ such that $f'=-1$. By subtracting a constant from $f$ we arrange that $f\not\asymp 1$, and thus $f\prec 1$ and $f>0$. This yields an embedding $E^a\to K$ sending $x^{-1}$ to $f$. As with $T^{\operatorname{nl}}_{\operatorname{small}}$, we derive from this the completeness of $T^{\operatorname{nl}}_{\operatorname{large}}$.

Decidability of these theories then follows since we can effectively enumerate a set of first-order axioms for $T^{\operatorname{nl}}$; see~\ref{sec:compactness}.
\end{proof}

\noindent
Let $E=\Q(x)$ be as in the proof of Corollary~\ref{cor:completionsdecidable}. Then $E$ is a grounded $H$-field with constant field $\Q$, and so $E$ has a unique $\HLO$-cut. Therefore $E$ has by Proposition~\ref{HLOnelicl} a Newton-Liouville closure $E^{\operatorname{nl}}$ in the sense that $E^{\operatorname{nl}}$ is an $\upo$-free newtonian Liouville closed $H$-field extension of $E$ that embeds over 
$E$ into any
$\upo$-free newtonian Liouville closed $H$-field extension of $E$. Thus $E^{\operatorname{nl}}$ is a model of $T^{\operatorname{nl}}_{\operatorname{small}}$ that embeds
into any model of 
$T^{\operatorname{nl}}_{\operatorname{small}}$. As $T^{\operatorname{nl}}_{\operatorname{small}}$ is model complete, this means that $E^{\operatorname{nl}}$ is a so-called {\em prime model}\/ of $T^{\operatorname{nl}}_{\operatorname{small}}$,
as defined in \ref{sec:mc}.

The $E^{\operatorname{nl}}$ obtained above is $\d$-algebraic over $E$, hence over $\Q$, and the constant field of $E^{\operatorname{nl}}$
is a real closure of the constant field $\Q$ of $E$. Thus by
Theorem~\ref{th:noextL} any prime model of $T^{\operatorname{nl}}_{\operatorname{small}}$
is isomorphic to $E^{\operatorname{nl}}$.

\medskip\noindent
We also wish to call attention to the $H$-subfield 
$\mathbb T^{\operatorname{da}}$ of $\T$ given by
$$\mathbb T^{\operatorname{da}}\ :=\ \big\{f\in\mathbb T:\ \text{$f$ is $\d$-algebraic (over $\Q$)}\big\}.$$
Note that $\mathbb T^{\operatorname{da}}$
contains the $\upo$-free $H$-subfield $\R(\ell_0, \ell_1, \ell_2,\dots)$ of $\T$. It follows from Lemma~\ref{embnl} that $\mathbb T^{\operatorname{da}}$ is
actually a Newton-Liouville closure of $\R(\ell_0, \ell_1, \ell_2,\dots)$. In particular, we have
 $\mathbb T^{\operatorname{da}}\preceq\mathbb T$.

\medskip\noindent
Eliminating quantifiers is contingent on a good choice of primitives, but a reasonable~QE should have consequences of a more intrinsic nature that would be hard to obtain otherwise. Below we derive such consequences. 

\subsection*{A further reduction} We can eliminate the primitives $\preceq$,~$\Upl$,~$\Upo$,~$\iota$ by introducing some ``ideal'' elements. In this way we reduce quantifier-free formulas to a very simple form, to be used in proving the results in the present section. 

More precisely, let $y=(y_1,\dots, y_n)$ be a tuple of distinct
syntactic variables, and $Y=(Y_1,\dots,Y_n)$ a
corresponding tuple of distinct differential indeterminates. A routine induction on terms shows that for any $\mathcal{L}^{\iota}$-term $t(y)$, there are quantifier-free formulas $\phi_1(y),\dots, \phi_m(y)$ ($m\ge 1$) in the language of differential rings, and differential
polynomials $F_1(Y), G_1(Y),\dots, F_m(Y), G_m(Y)\in \Q\{Y\}$, such that for all differential fields $K$ and $a\in K^n$,
\begin{align*} K& \models \phi_1(a) \vee \cdots \vee \phi_m(a),\ \text{ and for $i=1,\dots,m$,}\\
  \text{ if }K&\models \phi_i(a), \text{ then }\ G_i(a)\ne 0\ \text{ and }\ t(a)\ =\ \frac{F_i(a)}{G_i(a)}.
\end{align*}
Let $K$ be an $H$-field. Then we give $K$ its order topology and  $K^n$ the corresponding product topology. (Note that by Lemma~\ref{lem:Massaza} the order topology on $K$ equals its valuation topology if the valuation is nontrivial.)  

{\sloppy
\begin{cor}\label{c3plus} Suppose $K$ is an $\upo$-free newtonian Liouville closed $H$-field and~$X\subseteq K^n$ is definable in $K$.
Then $X$ has empty interior in $K^n$ if and only if ${X\subseteq \big\{a\in K^n: P(a)=0\big\}}$ for some $P\in K\{Y\}^{\ne}$.
\end{cor}
}
\begin{proof} Note that Lemma~\ref{c3} goes through for differential polynomials over $K$ in~$n$ indeterminates, by induction on $n$. Next, observe that the sets $\Upl(K)$ and 
$\Upo(K)$ are open and closed in $K$. Now use the remarks above and Theorem~\ref{thm:qe}. \end{proof}

\noindent
Now let $K$ be an $\upo$-free newtonian Liouville closed
$H$-field. Take an immediate $H$-field extension $L$ of $K$ with
an element $\upl$ such that $ \Upl(K) < \upl < \Upd(K)$, and set
$\upo:= \omega(\upl)$. Then $\Upo(K) < \upo < K\setminus \Upo(K)$, 
and for $f,g\in K$ with $g>0$, 
$$\frac{f}{g}\in \Upl(K)\ \Longleftrightarrow\ f < \upl\, g,\qquad \frac{f}{g} \in \Upo(K)\ \Longleftrightarrow\ f< \upo\, g.$$
Thus for any $\mathcal{L}^{\iota}$-term $t(y)$ as before the atomic formula $\Upl\big(t(y)\big)$ is
equivalent, in $K$ and for $y_1,\dots, y_n$ ranging over $K$, to a boolean combination of formulas, each of which has one of the following forms:
$$ F(y) < \upl\, G(y), \quad G(y)>0,\quad G(y)=0 \qquad (F, G\in \Q\{Y\}).$$ 
Likewise for the atomic formula $\Upo\big(t(y)\big)$
with $\upo$ instead of $\upl$. 

We can also eliminate occurrences of $\preceq$, but for this we
take a further $H$-field extension $L^*$ of $L$ with an element $c^*$ such that $C < c^* < a$ for all $a\in K^{>C}$. Then 
for all $f_1,g_1, f_2, g_2\in K$ with $g_1,g_2\ne 0$,
$$ \frac{f_1}{g_1}\ \preceq\ \frac{f_2}{g_2}\ \Longleftrightarrow\ 
|f_1g_2|\ \le\ c^*|f_2g_1|.$$
Thus for any $\mathcal{L}^{\iota}$-terms $t_1(y)$ and $t_2(y)$ the atomic formula $t_1(y)\preceq t_2(y)$ is
equivalent, in $K$ and for $y_1,\dots, y_n$ ranging over $K$, to a boolean combination of formulas, each of which has one of the following forms:
$$ F(y) \le c^*\, G(y), \quad G(y)>0,\quad G(y)=0 \qquad (F, G\in \Q\{Y\}).$$ 
To summarize some of the above in a single lemma, let $z_{ij}$ ($i=1,\dots,n$; $j\in \N$) and $v_\upl$, $v_\upo$, $v_{c^*}$, be distinct syntactic variables, and set $z_j:=(z_{1j},\dots,z_{nj})$. 
For $a=(a_1,\dots,a_n)\in K^n$ we set $a^{(i)}:=\big(a_1^{(i)},\dots,a_n^{(i)}\big)$, and as usual $a'=a^{(1)}$.
Recall from Section~\ref{sec:modth val fields} that $\mathcal L_{\operatorname{OR}}:=\{0,1,{-},{+},{\,\cdot\,},{\leq}\}$ is the language of ordered rings.

\begin{lemma}\label{lem:reduction to inequ}
Let $X\subseteq K^n$ be definable without parameters in $K$. Then
there is a quantifier-free $\mathcal L_{\operatorname{OR}}$-formula $\varphi(z_0,z_1,\dots,z_r,v_\upl, v_\upo, v_{c^*})$, for some $r\in\N$, such that
$$X\ =\ \big\{ a\in K^n:\  L^*\models \varphi\big(a,a',\dots,a^{(r)},\upl,\upo,c^*\big)\big\}.$$
\end{lemma}

\subsection*{NIP} We refer to~\ref{sec:structures having NIP} for a definition and
discussion of the
very robust but highly restrictive property NIP that some model-theoretic structures enjoy. 
In this subsection we establish part~(i) of Corollary~\ref{cor:qe}:

\begin{prop}\label{propNIP}
Every $\upo$-free newtonian Liouville closed $H$-field has $\operatorname{NIP}$.
\end{prop}
\begin{proof}
Let $K$ be an  $\upo$-free newtonian Liouville closed $H$-field. Assume towards a contradiction that the relation $R\subseteq K^m\times K^n$ is definable without parameters in $K$ and independent. 
We just do the case $m=n=1$; the general case only involves more notation.
Thus for every $N\geq 1$ there are $a_1,\dots,a_N\in K$ and $b_I\in K$ ($I\subseteq \{1,\dots,N\}$),
such that for $i=1,\dots,N$ and all $I\subseteq \{1,\dots,N\}$,
$$R(a_i,b_I) \quad\Longleftrightarrow\quad i\in I.$$
Let $L^*$ be an $H$-field extension of $K$ as at the beginning of this section, containing~$\upl$,~$\upo$,~$c^*$. By Lemma~\ref{lem:reduction to inequ} we can take a quantifier-free
$\mathcal L_{\operatorname{OR}}$-formula $$\varphi(x_0,x_1,\dots,x_r,y_0,y_1,\dots,y_r,v_\upl, v_\upo, v_{c^*}),$$ 
such that for all $a,b\in K$:
$$R(a,b) \quad\Longleftrightarrow\quad L^*\models \varphi(a,a',\dots,a^{(r)},b,b',\dots,b^{(r)},\upl,\upo,c^*).$$
Thus the relation $R^*\subseteq (L^*)^{r+1}\times (L^*)^{r+4}$ given by
$$R^*(a_0,\dots,a_r,b_0,\dots,b_{r+3})\quad\Longleftrightarrow\quad L^*\models 
 \varphi(a_0,\dots,a_r,b_0,\dots,b_{r+3})$$
is independent and quantifier-free definable in the 
$\mathcal L_{\operatorname{OR}}$-structure $L^*$, that is, in the ordered field
$L^*$. This
contradicts~\ref{cor:RCF NIP}.
\end{proof}

\subsection*{The induced structure on the constant field} The goal of this subsection is to establish the following:

\begin{prop}\label{indconstants} Let $K$ be an $\upo$-free newtonian Liouville closed $H$-field, and let $X\subseteq K^n$ be definable in $K$. Then
$X\cap C^n$ is semialgebraic in the sense of $C$. 
\end{prop}

\noindent
The proof goes by reduction to Proposition~\ref{indsemi}, using our QE and the fact that a real closed $H$-field $K$ with constant field $C\ne K$ yields
a tame pair $(K,C)$ as defined in Section~\ref{sec:modth val fields}.

\begin{proof}[Proof of Proposition~\ref{indconstants}] Take an immediate real closed $H$-field extension~$L$ of $K$ and $\upl,\upo\in L$ as earlier in this section. 
As in the proof of Corollary~\ref{cor:induced structure on C}, our~QE reduces the problem to showing for
polynomials $p,q\in K[Y_1,\dots, Y_n]$ that the following subsets of~$C^n$ are semialgebraic in the sense of $C$:
\begin{align*}
&\big\{c\in C^n:\ p(c)=0\big\},\qquad \big\{c\in C^n:\ p(c)>0\big\},\qquad \big\{c\in C^n:\ p(c) \preceq q(c)\big\},\\
&\big\{c\in C^n:\  p(c) < \upl\, q(c)\big\}, \quad \big\{c\in C^n:\  p(c) < \upo\, q(c)\big\}.
\end{align*}
This holds for the first three sets by a direct application of Proposition~\ref{indsemi} to the tame pair $(K,C)$.
The two sets involving $\upl$ and~$\upo$ are semialgebraic in the sense of $C$ by applying Proposition~\ref{indsemi} to the 
tame pair $(L,C)$.  
\end{proof}

%\noindent
%As a somewhat random illustration of 
%Proposition~\ref{indconstants}, note that the set of 
%real parameters $(\lambda,\mu)\in \R^2$ for which the system
%$$  y + \lambda(y'')^2 + \mu y'y''' +1\ =\ \log x, 
%\qquad 0\ <\ y\ \preceq\ e^x$$ 
%has a solution in $\T$ is a semialgebraic subset of $\R^2$.   

\subsection*{O-minimality at infinity}  By this we mean the following:

\begin{prop}\label{ominatinf} Let $K$ be an $\upo$-free newtonian Liouville closed $H$-field, and let $X\subseteq K$ be definable in $K$. Then there exists an element $a\in K$ such that
$$(a,+\infty)\ \subseteq\ X \ \text{ or }\ (a,+\infty) \cap X\ =\ \emptyset.$$
\end{prop}

\noindent
This is proved by logarithmic decomposition of differential polynomials; we refer to Section~\ref{Decompositions of Differential Polynomials} for how these decompositions
are defined. First some lemmas.

\begin{lemma}\label{lioudagger} Let $K$ be a Liouville closed $H$-field and $K < a$ 
where $a$ lies in some
$H$-field extension of $K$. Then $K <(a^\dagger)^m < a$ for all $m\ge 1$.
\end{lemma}
\begin{proof} Since $(K^\times)^\dagger=K$, the set $\Psi\subseteq \Gamma$ is downward closed. Set $\alpha:= va$. Then $\alpha< \Gamma$, and so $\alpha^\dagger<
\Gamma$: otherwise we have $\gamma\in \Gamma^{<}$ such that
$\alpha^\dagger > \gamma^\dagger$, and so $\alpha > \gamma$, contradicting $\alpha < \Gamma$. Also, $\alpha^\dagger = o(\alpha)$
by Lemma~\ref{PresInf-Lemma2}; to apply this lemma, first shift by an element of $\Gamma$ to reduce to the small derivation case. It remains to note that $a^\dagger>0$ and $v(a^\dagger)=\alpha^\dagger$.
\end{proof}

\begin{lemma}\label{logdecliou} With $K$ and $a$ as in Lemma~\ref{lioudagger}, $a$ is $\d$-transcendental over $K$,
$C_{K\<a\>}=C$, and $K\<a\>$ is an $H$-field extension of $K$ whose value group
$$v\big(K\<a\>^\times\big)\ =\ \Gamma \oplus \bigoplus_{n} \Z v\big(a^{\<n\>}\big)\qquad\text{\textup{(}internal direct sum\textup{)}}$$
contains $\Gamma$  
as a convex subgroup. If $K$ is $\upo$-free, then so is $K\<a\>$.
\end{lemma}
\begin{proof} By induction and Lemma~\ref{lioudagger} we have
$$ K\ <\ (a^{\<n+1\>})^m\  <\  a^{\<n\>} \qquad  (n=0,1,2,\dots,\ m=1,2,\dots). $$
Let $P\in K\{Y\}^{\ne}$ of order $\le r$ have logarithmic
decomposition
$$ P\ =\ \sum_{\i}P_{\<\i\>}Y^{\<\i\>}$$
with $\i$ ranging over $\N^{1+r}$, all $P_{\<\i\>}\in K$, and
$P_{\<\i\>}\ne 0$ for only finitely many $\i$. Take $\j\in \N^{1+r}$
lexicographically maximal with
$P_{\<\j\>}\ne 0$.  It follows from the above that
\begin{align*} 
P(a)\ &\sim\ P_{\<\j\>}\cdot a^{\<\j\>},\ \text{ and thus}\\
P(a)\ne 0,\qquad \sgn\ P(a)\ &=\ \sgn\ P_{\<\j\>},\qquad v\big(P(a)\big)\ =\ v P_{\<\j\>}+\sum_{n=0}^r j_n v\big(a^{\<n\>}\big).
\end{align*}
Thus $a$ is $\d$-transcendental over $K$, and for any 
$f\in K\<a\>^{\ne}$ there are $g\in K^\times$ and 
$k_0,\dots, k_r\in \Z$ such that $f\sim g\cdot \big(a^{\<0\>}\big)^{k_0}\cdots \big(a^{\<r\>}\big)^{k_r}$. Therefore
$\res K\<a\> =\res K$, and so $K\<a\>$ is an $H$-field extension
of $K$ with $C=C_{K\<a\>}$. The statement about the value group of
$K\<a\>$ also follows easily. Suppose now that $K$ is $\upo$-free, and~$K\<a\>$ is not; it remains to derive a contradiction from this assumption. Since $\Gamma^{<}$ is cofinal in~$\Gamma_{K\<a\>}^{<}$ this gives
an element $\upo\in K\<a\>$ such that $\upo_{\rho}\leadsto \upo$,
where $(\upo_{\rho})$ is the sequence in $K$ obtained in the usual way
from a logarithmic sequence in $K$. Now~$(\upo_{\rho})$ is
of $\d$-transcendental type over $K$, so $K\<\upo\>$ 
is an immediate extension of~$K$, and $\upo$ is $\d$-transcendental
over $K$. Now $\upo=P(a)/Q(a)$ with $P,Q\in K\{Y\}^{\ne}$, so
$a$ is a zero of the differential polynomial $\upo\, Q(Y) -P(Y)\in K\<\upo\>\{Y\}^{\ne}$, and thus~$K\<a\>$ is $\d$-algebraic over 
$K\<\upo\>$. It follows that $K\<a\>$ has finite trancendence degree
over $K\<\upo\>$, and so $\Q v\big(K\<a\>^\times\big)/\Gamma$ has finite dimension as a vector space over~$\Q$, contradicting the above
structure of $v\big(K\<a\>^\times\big)$.    
\end{proof} 

{\sloppy

\begin{proof}[Proof of Proposition~\ref{ominatinf}] By a routine translation into model-theoretic terms it is enough to show the following:

\claim{Let $L$ be an elementary extension of $K$ with elements
$a,b>K$. Then there is a pre-$H$-field isomorphism
$i\colon K\<a\> \to K\<b\>$ over $K$ with $i(a)=b$ such that also
$$i\big(\Upl(L)\cap K\<a\>\big)\ =\ \Upl(L)\cap K\<b\>, \qquad i\big(\Upo(L)\cap K\<a\>\big)\ =\ \Upo(L)\cap K\<b\>.$$
} 

\vskip-2em

\noindent 
A pre-$H$-field isomorphism $i\colon K\<a\>\to K\<b\>$ over $K$ with
$i(a)=b$ is obtained from Lemma~\ref{logdecliou} and its proof, in particular,
the equalities for $\sgn\ P(a)$ and $v\big(P(a)\big)$ in that proof. Since $K$ is $\upo$-free, so are $K\<a\>$ and $K\<b\>$ by
the same lemma, and so the additional property claimed for $i$ is now
a consequence of Corollary~\ref{cor:ILO1}. 
\end{proof}
}

\noindent
Using fractional linear transformations we get analogous behavior
of any definable set $X\subseteq K$ to the left as well as to the
right of any point in $K$. In other words, $K$ is locally o-minimal in the sense of Marker and Steinhorn; see~\cite{TV}. Thus:

\begin{cor}\label{opendiscrete} Let $K$ be an $\upo$-free newtonian Liouville closed $H$-field, and let $X\subseteq K$ be definable in $K$.
Then $X$ is the disjoint union of an open definable subset of $K$ and a discrete definable subset of $K$. Moreover,
$X$ is discrete in $K$ iff $X\subseteq \big\{y\in K:\ P(y)=0\big\}$
for some $P\in K\{Y\}^{\ne}$.
\end{cor}
\begin{proof} The interior of $X$ in $K$ is definable, and
so $X$ with its interior removed is discrete, by local o-minimality. For the second part, use Corollary~\ref{c3plus}.
\end{proof}

\subsection*{O-minimality at $C^{\downarrow}$} 
We also have o-minimality at another important cut:

\begin{prop} Let $K$ be an $\upo$-free newtonian Liouville closed $H$-field, and let $X\subseteq K$ be definable in $K$. Then there exists an element $a>C$ in $K$ such that
$$\{f\in K:\ C < f < a\}\ \subseteq\ X \ \text{ or }\ \{f\in K:\ C < f < a\} \cap X\ =\ \emptyset.$$
\end{prop}
\begin{proof} By a routine translation into model-theoretic terms
it suffices to show:

\claim{Let $L$ be an elementary extension
of $K$ and $f,g\in L$ be such that ${C_L < f < a}$ and $C_L<g<a$
for all $a>C$ in $K$. Then there is a pre-$H$-field isomorphism
$i\colon K\<f\> \to K\<g\>$ over $K$ with $i(f)=g$ such that also
$$i\big(\Upl(L)\cap K\<f\>\big)\ =\ \Upl(L)\cap K\<g\>, \qquad i\big(\Upo(L)\cap K\<f\>\big)\ =\ \Upo(L)\cap K\<g\>.$$
}

\vskip-2em

\noindent
To prove this claim, note first that $\Gamma^{<}< vf < 0$ and
$\Gamma^{<} < vg < 0$. So Corollary~\ref{Evalocor} yields a pre-$H$-field
isomorphism $K\<f\> \to K\<g\>$ over $K$ sending $f$ to~$g$. Also, $K\<f\>$ and $K\<g\>$ have a smallest comparability class by Proposition~\ref{Evalcor}, and thus
$K\<f\>$ and $K\<g\>$ have unique $\Upl\Upo$-cuts by 
Corollary~\ref{cor:ILO1}.  This yields
the claim.
\end{proof} 

{\sloppy

\subsection*{Notes and comments} The uniqueness-up-to-isomorphism of prime models of~$T^{\operatorname{nl}}_{\operatorname{small}}$ holds also on general model-theoretic grounds: \cite[Corollary~4.2.16]{Marker}.

%In \cite{Boshernitzan81,Boshernitzan}, Boshernitzan shows that the intersection $E$ of all maximal Hardy fields is Liouville closed and contains $\R$ (and hence also Hardy's field of logarithmico-exponential functions), and is $\d$-algebraic over~$\R$. In the next volume we show that each maximal Hardy field  is $\upo$-free; this together with Lemma~\ref{embnl} (applied to $\R(\ell_0,\ell_1,\dots)$ and a maximal Hardy field in place of $E$ and $L$) implies that there exists an  embedding $E\to\mathbb T^{\operatorname{da}}$, so in particular, $E$ is also $\upo$-free.

}

Some of the arguments in the subsections on the induced structure on the constant field and o-minimality at infinity were already
used in Section~5 of~\cite{ADHOleron} to prove quantifier-free versions of Propositions~\ref{indconstants} and~\ref{ominatinf} for $K=\T$.

Shelah~\cite{Shelah-dp} considers
a strengthening of  NIP, called \textit{dp-minimality}\/; see \cite{DGL} for basic facts about this notion. 
Algebraically closed valued fields, the field of $p$-adic numbers,
and o-minimal structures
are dp-minimal. Simon~\cite[Theorem~3.6]{Simon}  proved that if an expansion $(G;{\leq}, 0,{+},\dots)$ of a divisible ordered abelian group is dp-minimal, then all infinite definable subsets of $G$ have nonempty interior. Thus if~$K$ is any pre-$H$-field and $K\neq C$, then
$K$ is {\em not\/} dp-minimal, since the definable set~$C\subseteq K$ has empty interior.\marginpar{statements about dp-minimality and distality taken on faith} Another strengthening of NIP is the notion of \textit{distality,}\/ due to Simon~\cite{Simon13}. O-minimal structures and the field of $p$-adic numbers are distal, but algebraically closed valued fields are not. We intend to show elsewhere that $\T$ is distal. 

 In view of \cite[\S{}2.25]{vdD89}, Corollary~\ref{c3plus}
 yields a natural notion of dimension for definable sets $X\subseteq \T^n$; for details, see \cite{ADHdim}. 
In connection with Corollary~\ref{opendiscrete}, the standard example of an infinite discrete definable subset
of~$\T$ is of course the set $\R$ of constants. The question arises if this is the source of all discreteness: is every
discrete definable subset of $\T^n$ the image of some semialgebraic set $S\subseteq \R^m$ under a definable map
$S \to \T^n$? It turns out that the answer is negative; see \cite{ADHdim}.

%% file: mt-trans.tex
\let\thetheoremold\thetheorem
\renewcommand{\thetheorem}{%
    \thechapter.\arabic{theorem}%
}

\appendix

\chapter{Transseries}\label{app:trans}

\setcounter{theorem}{0}

\noindent
We assume here familiarity with well-based series and Hahn fields as exposed in Section~\ref{sec:valued fields}. We begin by adding some items to this material.
Our construction of $\T$ is self-contained as to concepts and definitions, but for proofs of some
key properties we refer to \cite{DMM2}, where $\T$ is denoted by $\R\(( x^{- 1}\)) ^{\operatorname{LE}}$, or $\R\(( t\)) ^{\operatorname{LE}}$ with $t=x^{-1}$, and called the {\em field
of logarithmic-exponential series}\/ (in $x$ over $\R$). 
The construction is also very similar to the treatment in
Schmeling's thesis~\cite{Schm01}.   

The reader should be aware that notations and terminology concerning Hahn fields and transseries vary considerably across the literature, even in our own earlier works. (For example, the $\T$ in \cite{ADH} is not the $\T$
constructed here.) In the present volume we have systematized
things by adopting many notations from~\cite{JvdH}.

\subsection*{Summability in Hahn fields}
In what follows, $\fM$ is a multiplicative (totally)
ordered abelian group, ordered by~$\preceq$. Also $C$ will be a (coefficient) field, so that we have the Hahn field 
$\CM$, with the internal direct sum decomposition
$$\CM\ =\ C[[\fM^{\succ 1}]] \oplus C \oplus C[[\fM^{\prec 1}]]$$ into
$C$-linear subspaces. Note that $C[[\fM^{\prec 1}]]=\CM^{\prec 1}$.  

A family $(f_\lambda)_{\lambda\in \Lambda}$ in $\CM$ is said to be
{\bf summable\/} if $\bigcup_\lambda \supp f_\lambda$ is well-based and
for each $\fm\in \fM$ there are only finitely many $\lambda\in \Lambda$ such that $f_{\lambda,\fm}\ne 0$; in that case we define
its sum $\sum_\lambda f_\lambda$ to be the series $f\in \CM$ such that $f_{\fm}=\sum_\lambda f_{\lambda,\fm}$ for each $\fm\in \fM$. (This agrees with the usual notation\marginpar{has  ``usual notation'' been introduced?} for elements of $\CM$: for a series
$f=\sum_{\fm} f_{\fm}\fm\in C[[\fM]]$ the family 
$f_{\fm}\fm$ is indeed summable with sum $f$.)

\index{summable!family of elements}
\index{sum!family}
\index{family!summable}
\index{family!sum}

\medskip\noindent
Let $t=(t_1,\dots, t_n)$ be a tuple of distinct variables and let
$$F\ =\ F(t)\ =\ \sum_{\nu}c_{\nu}t^{\nu}\in C[[t]]\ :=\ C[[t_1,\dots, t_n]]$$
be a formal power series over $C$; here the sum ranges over all multiindices $\nu=(\nu_1,\dots, \nu_n)\in \N^n$, and
$c_{\nu}\in C$, $t^{\nu}:= t_1^{\nu_1}\cdots t_n^{\nu_n}$. For any tuple $\varepsilon=(\varepsilon_1,\dots, \varepsilon_n)$ of elements of $C[[\fM]]^{\prec 1}$ the family $(c_{\nu}\epsilon^{\nu})$ is summable, where $\varepsilon^{\nu}:= \varepsilon_1^{\nu_1}\cdots \varepsilon_n^{\nu_n}$. Put
$$F(\varepsilon):= \ \sum_{\nu}c_{\nu}\varepsilon^{\nu}\in C[[\fM]]^{\preceq 1}\ =\ C[[\fM^{\preceq 1}]].$$
For example, if $C$ has characteristic zero, then for $n=1$ and $t=t_1$ the formal series 
$\exp(t)=\sum_{\nu=0}^\infty t^\nu/\nu!$ yields a partial
 exponential function 
$$\varepsilon\ \mapsto\ \exp(\varepsilon)\ =\ \sum_{\nu=0}^\infty \varepsilon^\nu/\nu!\ :\ C[[\fM]]^{\prec 1} \to 1+C[[\fM]]^{\prec 1},$$
an isomorphism of the additive subgroup $C[[\fM]]^{\prec 1}$ of $C[[\fM]]$
onto the multiplicative subgroup $1+C[[\fM]]^{\prec 1}$ of $C[[\fM]]^\times$, with inverse
$$1+\delta\ \mapsto\ \log (1+\delta)\ :=\ \sum_{\nu=1}^\infty (-1)^{\nu-1}\delta^\nu/\nu\ :\ 1+C[[\fM]]^{\prec 1}\to C[[\fM]]^{\prec 1}.$$
Let $\fN$ also be a multiplicative ordered abelian group.
Then
a map $$\Phi\ :\ \CM \to \CN$$ is said to be
{\bf strongly additive\/} if for each summable family
$(f_\lambda)$ in $\CM$ the family $\big(\Phi(f_\lambda)\big)$ in $\CN$ is summable with $\Phi(\sum_\lambda f_\lambda)=\sum_\lambda \Phi(f_\lambda)$.
% if in addition $\Phi$ is $C$-linear, we say that $\Phi$ is {\bf strongly $C$-linear}.
Note that if $\Phi$ is strongly additive, then it is additive.  
  
\index{map!strongly additive}
%\index{map!strongly linear}
\index{strongly!additive}
%\index{strongly!linear}
%\index{additive!strongly}
%\index{linear!strongly}
 
\subsection*{The case where $\fM$ is a product with convex factor $\fG$} Suppose now that $\fG$ and $\fR$ are ordered subgroups of $\fM$ such that 
$$ \text{$\fG$ is convex in $\fM$}, \qquad \fG\cap \fR\ =\ \{1\}, \quad \fM\ =\ \fG\fR:=\{\fg\fr:\ \fg\in \fG,\ \fr\in \fR\}.$$
Then we have an isomorphism $C[[\fM]]\to C[[\fG]][[\fR]]$
of $C[[\fG]]$-algebras given by 
$$f\ =\ \sum_{\fm}f_\fm \fm\ \mapsto\ \sum_{\fr\in \fR}\left(\sum_{\fg\in \fG}f_{\fg\fr}\fg\right)\fr.$$
For $f\in C[[\fM]]$ we have in fact $f=\sum_{\fr\in \fR}\left(\sum_{\fg\in \fG}f_{\fg\fr}\fg\right)\fr$ where the indicated sums exist in $C[[\fM]]$ according to the definition of summability. Whenever convenient we 
identify below $C[[\fM]]$ and $C[[\fG]][[\fR]]$ via the above isomorphism.

If in addition $C$ is an ordered field, then $C[[\fM]]$ and $C[[\fG]]$ are ordered Hahn fields, and so is
$C[[\fG]][[\fR]]$, and the above isomorphism is also an isomorphism of ordered fields. (In this remark and in what follows the reader is assumed to be familiar with Section~\ref{sec:valued ordered fields}.)

\subsection*{Directed unions of Hahn fields} A key feature of
$\T$ will be its structure as a directed union of Hahn fields over its constant field $\R$. (A Hahn field over $\R$ with its natural valuation and ordering and any derivation is never a Liouville closed $H$-field, by \cite[Corollary~7.2]{AvdD3}, and so cannot have the properties
we expect of $\T$.) It is therefore useful to extend the notions above to such directed unions, and so we consider here a directed family $(\fM)_{i\in I}$ with $I\ne \emptyset$, of ordered subgroups
of the ordered multiplicative group~$\fM$ such that $\fM=\bigcup_i \fM_i$. Here ``directed'' means that for all
$i,j\in I$ there exists~$k\in I$ with $\fM_i, \fM_j\subseteq \fM_k$. This leads to a directed union of
Hahn fields over $C$, namely the valued subfield  
$$K\ :=\ \bigcup_i C[[\fM_i]]$$ 
of $C[[\fM]]$. Define a {\bf $K$-subgroup\/} of $\fM$ to be an ordered subgroup $\fG$ of $\fM$ such that~${C[[\fG]]\subseteq K}$, inside the ambient $C[[\fM]]$; thus each $\fM_i$ is a $K$-subgroup of $\fM$.  
We say that the family $(\fM_i)$ is {\bf healthy} if every 
$K$-subgroup of $\fM$ is contained in some~$\fM_i$; this might depend on $C$. An easy diagonal argument shows: if $I$ is countable (the relevant case for us), then $(\fM_i)$ is healthy. Also, by an easy cofinality argument, 
if every $\fM_i$ is convex in $\fM$,
then $(\fM_i)$ is healthy.  

\index{family!healthy}
\index{healthy}
\index{K-subgroup@$K$-subgroup}
\index{subgroup!$K$-subgroup}

{\em Assume below that $(\fM_i)$ is healthy}. A family
$(f_\lambda)$ in $K$ is said to be {\bf summable\/} if
there exists a $K$-subgroup $\fG$ of $\fM$ such that all $f_{\lambda}\in C[[\fG]]$
and $(f_{\lambda})$ is summable as a family in $C[[\fG]]$; note that then $\sum_\lambda f_{\lambda}$ is defined as an element of $K$ (lying in~$C[[\fG]]$ for $\fG$ as above). 
Thus for  $F = \sum_{\nu}c_{\nu}t^{\nu}\in C[[t]]$ with $t=(t_1,\dots, t_n)$ and $\epsilon\in \smallo_K^n$, the family
$(c_{\nu}\epsilon^{\nu})_{\nu\in \N^n}$ is summable, with
$F(\epsilon):=\sum_{\nu}c_{\nu}\epsilon^{\nu}\in \mathcal{O}_K$. 

Let~$\fN$ also be an ordered abelian group, with
$\fN=\bigcup_j\fN_j$ and $(\fN_j)$ a directed family of ordered subgroups of $\fN$. Let
$L=\bigcup_j C[[\fN_j]]\subseteq C[[\fN]]$, and assume~$(\fN_j)$ is healthy. Then a map $\Phi\colon K \to L$ is said to be {\bf strongly additive\/} if for every summable family 
$(f_{\lambda})$ in $K$ the family $(\Phi(f_{\lambda}))$ is summable in $L$, and $\Phi(\sum_{\lambda}f_{\lambda})=\sum_{\lambda}\Phi(f_{\lambda})$. 
A map $\Phi\colon K\to L$ is said to be {\bf healthy\/} 
if for each $K$-subgroup $\fG$ of~$\fM$ there exists an
$L$-subgroup $\fH$ of $\fN$ such that $\Phi\big(C[[\fG]]\big)\subseteq C[[\fH]]$.\index{map!strongly additive}\index{strongly!additive}

%\index{map!strongly linear}
%\index{strongly!linear}
% and the restriction of $\Phi$ to a map
%$C[[\fG]]\to C[[\fH]]$ is strongly additive. 
%The term {\bf strongly $C$-linear} is used likewise.    

\subsection*{Exponential ordered fields} An {\bf exponentiation\/} on a field $E$ is a group morphism $\exp\colon E \to E^\times$ from the additive group of $E$ into its multiplicative group.
An {\bf exponential ordered field\/} is an ordered field $E$ equipped with a strictly increasing exponentiation on $E$ (denoted by $\exp$ unless specified otherwise,  
necessarily taking values in the multiplicative subgroup~$E^>$ of~$E^\times$). A {\bf logarithmic-exponential ordered field\/} is an exponential ordered field $E$ with $\exp(E)=E^{>}$;
the inverse of the ordered group isomorphism $\exp\colon E \to E^{>}$
is then an ordered group isomorphism 
$\log\colon E^{>} \to E$.
Below we consider the ordered field $\R$ of real numbers as a logarithmic-exponential ordered field with exponentiation $r\mapsto \ex^r$. For any logarithmic-exponential ordered field $E$ we set 
$a^f:= \exp(f\log a)\in E^{>}$
for $a\in E^{>}$ and $f\in E$, so $a^0=1$ and $a^1=a$, and the usual identitites follow:
$$a^{f+g}\ =\ a^fa^g,\quad (ab)^f\ =\ a^fb^f, \quad a^{fg}\ =\ (a^f)^g \qquad(a,b\in E^{>},\ f,g\in E).$$ 
Initially we shall construct $\T$ as a logarithmic-exponential ordered field extension of~$\R$; the definition of the derivation on $\T$ comes later. This construction involves
the following general procedure. We define a {\bf pre-exponential ordered field}
to be a tuple~$(E, A, B, \exp)$ such that: \index{exponential!ordered field}\index{exponentiation}\index{field!exponential ordered}\index{ordered field!exponential}\index{logarithmic-exponential ordered field}\index{field!logarithmic-exponential ordered}\index{ordered field!logarithmic-exponential}\index{pre-exponential ordered field} \index{field!pre-exponential ordered} \index{ordered field!pre-exponential}
\begin{enumerate}
\item $E$ is an ordered field;
\item $A$ and $B$ are additive subgroups of $E$ with $E=A\oplus B$ and $B$
convex in $E$;
\item $\exp\colon B\to E^{\times}$ is a strictly increasing group morphism (so $\exp(B)\subseteq E^>$).
\end{enumerate}
Let $(E,A,B,\exp)$ be a pre-exponential ordered field. We view $A$ as the part of $E$ where exponentiation is not yet defined, and accordingly we introduce
a ``bigger'' pre-exponential ordered field $(E^*,A^*,B^*,\exp^*)$ as follows:
Take a {\em multiplicative\/} copy $\exp^*(A)$ of the ordered additive
group $A$ with order-preser\-ving isomorphism 
$\exp^*\colon A\to\exp^*(A)$,
and put $E^*\ :=\ E[[\exp^*(A)]]$. Viewing $E^*$ as an ordered Hahn field over the ordered coefficient field $E$, we set
$$ A^*\ :=\ E[[\exp^*(A)^{\succ 1}]],\qquad B^*\ :=\ (E^*)^{\preceq 1}\ =\ E\oplus 
(E^*)^{\prec 1}\ =\ A \oplus B \oplus (E^*)^{\prec 1}.$$ 
Note that $\exp^*(A)^{\succ 1}=\exp^*(A^{>})$. Next we extend
$\exp^*$ to $\exp^*\colon B^*\to (E^*)^{\times}$ by
$$\exp^*(a+b+\varepsilon)\ :=\ \exp^*(a)\cdot \exp(b)\cdot \sum_{n=0}^\infty
\frac{\varepsilon^n}{n!} \qquad(a\in A,\ b\in B,\ \varepsilon\in (E^*)^{\prec 1}).$$
Then $E\subseteq B^*=\operatorname{domain}(\exp^*)$, and $\exp^*$ extends $\exp$. Note that $E<(A^*)^>$ (but $\exp^*(E)$ is cofinal in $E^*$ if $A\ne \{0\}$). In particular, for $a\in A^{>}$, we have 
$$\exp^*(a)\in \exp^*(A^{>})\ \subseteq\ (A^*)^>,\ \text{  so $\exp^*(a)\ >\ E$.}$$ 
Suppose now that $E=\R[[\fN]]$, where $\fN$ is a
multiplicative ordered abelian group. We identify
$\fN$ and $\exp(A)$ with subgroups of the
product group $\fN^*=\fN\times \exp(A)$ via
$\fn\mapsto (\fn,1)$ and $e\mapsto (1,e)$ for $\fn\in \fN$ and 
$e\in \exp(A)$.
Then $\fN^*=\fN\,\exp(A)$ and $\fN\cap \exp(A)=\{1\}$ (in $\fN^*$). We make $\fN^*$ into an ordered abelian group so that
$\fN$ and $\exp(A)$ are ordered subgroups and $\fN$ is convex in $\fN^*$. 
The effect is that 
$$E^*\ =\ \R[[\fN]]\,[[\exp(A)]]\ =\ \R[[\fN^*]],$$ 
after the natural identifications, so the inclusion $E\subseteq E^*$ is now the 
inclusion $\R[[\fN]]\subseteq \R[[\fN^*]]$ induced by $\fN\subseteq \fN^*$. 
 Viewing $E$ and $E^*$ as Hahn fields over~$\R$, we get:
if all infinitesimals of $E$ lie in $B$ and  
$\exp(\epsilon)=\sum_{\nu=0}^\infty \epsilon^\nu/\nu!$ for all 
infinitesimal~$\epsilon\in E$, then all infinitesimals of $E^*$ lie in $B^*$ and $\exp^*(\epsilon)=\sum_{\nu=0}^\infty \epsilon^\nu/\nu!$ for all 
infinitesimal~${\epsilon\in E^*}$.

\subsection*{Construction of $\T_{\exp}$} Starting with 
$E_0 :=\ \R[[\G_0]]$, with  $\G_0=x^\R,$
we construct the field $\T_{\exp}=\bigcup_m E_m$
of \textit{exponential transseries}\/ as the union of an increasing sequence of Hahn fields $E_m=\R[[\G_m]]$. First we make the ordered Hahn field $E_0$ over $\R$ into the pre-exponential ordered field 
\begin{align*} (E_0,A_0,B_0,\exp_0), \quad
A_0\ :=\ \R[[\G_0^{\succ 1}]],\quad B_0\ :=\ E_0^{\preceq 1}\ =\
\R\oplus E_0^{\prec 1}, \text{ with}\\ 
\exp_0\colon B_0\to E_0^{\times}\ \text{ given by }\
\exp_0(r+\varepsilon)\ :=\ \ex^r\sum_{n=0}^\infty\varepsilon^n/n! \qquad(r\in \R,\ \varepsilon\in E_0^{\prec 1}).
\end{align*} 
Inductively, we assume given the pre-exponential ordered field 
$(E_m, A_m, B_m, \exp_m)$ with the ordered Hahn field $E_m=\R[[G_m]]$ over 
$\R$, and set 
\begin{align*} (E_{m+1},A_{m+1},B_{m+1},\exp_{m+1})\ :&=\ 
(E_m^*,A_m^*,B_m^*,\exp_m^*), \text{ so}\\[0.5em] 
E_m\ =\ \R[[\G_m]]\ \subseteq\ E_{m+1}\ &=\ \R[[\G_{m+1}]] \ \text{(inclusions of ordered Hahn fields)}
\end{align*} with 
$\G_m$ a convex ordered subgroup of $\G_{m+1}=G_m\exp(A_m)$. We put
$$\G^{\operatorname{E}}\ :=\ \bigcup_m \G_m, \qquad \T_{\exp}\ =\ \RxE\ :=\ 
\bigcup_m E_m,$$
with $\G^{\operatorname{E}}$ construed as the multiplicative ordered abelian group having the $\G_m$ as ordered subgroups, and 
$\T_{\exp}$ as the
ordered field with the $E_m$ as ordered subfields. The elements of $\G^{\operatorname{E}}$ are called {\bf exponential transmonomials} (or $\operatorname{E}$-monomials), and those of
$\T_{\exp}$ are called {\bf exponential transseries} (or $\operatorname{E}$-series).
The alternative notation~$\RxE$ for $\T_{\exp}$ highlights the role of the formal variable~$x$ and the initial Hahn field
$\R[[x^{\R}]]$ in the construction of $\T_{\exp}$, with the superscript $\operatorname{E}$ indicating closure under exponentiation.
Let $\exp\colon \T_{\exp}\to\T^\times_{\exp}$ be the common extension of the~$\exp_m$. Then~$\T_{\exp}$ with $\exp$ is an exponential ordered field extension of~$\R$.
The ordered Hahn field $\R[[\G^{\operatorname{E}}]]$ gives an ordered field inclusion 
${\T_{\exp} \subseteq\R[[\G^{\operatorname{E}}]]}$. We think of any $f\in \T_{\exp}$ as a series $f(x)\in \R[[G^{\operatorname{E}}]]$ with $\supp f \subseteq \G^{\operatorname{E}}$.  Considering~$\T_{\exp}$ also as a {\em valued\/} subfield of the Hahn field $\R[[\G^{\operatorname{E}}]]$ we have
$$\exp(\epsilon)\ =\ \sum_{\nu=0}^\infty \epsilon^\nu/\nu! \qquad \text{ for infinitesimal }\epsilon\in \T_{\exp}.$$
Note: $\T_{\exp}$ is dense in this valued field $\R[[\G^{\operatorname{E}}]]$, since every $\G_m$ is convex in $\G^{\operatorname{E}}$. In order to indicate an element of $\R[[\G^{\operatorname{E}}]]$ outside $\T_{\exp}$, set $\exp_0(x):=x$ and
$\exp_{n+1}(x):= \exp\!\big(\!\exp_n(x)\big)$; note that we abandon here the earlier
meaning of $\exp_n\colon B_n\to E_n^{\times}$. Induction gives
$\exp_{n+1}(x)\in \exp(A_n^{>})$, so $\exp_{n+1}(x)\in \G_{n+1}$
and $\exp_{n+1}(x) > G_n$.  
Thus the series $\sum_{n=0}^\infty 1/\exp_n(x)$ 
lies in $\R[[\G^{\operatorname{E}}]]$ but not in $\T_{\exp}$.

Straightforward inductions on $m$ yield:

\index{transmonomials} \index{group!exponential transmonomials} \index{E-monomials@$\operatorname{E}$-monomials} \index{transseries!exponential} \index{exponential!transseries} \index{E-series@$\operatorname{E}$-series}

\begin{lemma}\label{appendixT1} The $G_m$ and $A_m$ have the following basic properties:
\begin{enumerate}
\item[\textup{(i)}] $A_m\ =\ \{f\in \T_{\exp}:\  \G_{m-1}\prec \supp f 
\subseteq \G_{m}\}$, with $\G_{-1}:=\{1\}$;
\item[\textup{(ii)}]  $|a| > A_{m-1}$ for all $a\in A_m^{\ne}$, with $A_{-1}:= \{0\}$;
\item[\textup{(iii)}] 
$\big\{f\in \R[[\G_m]]:\ \supp f \succ 1\big\}\ =\  A_0 \oplus \cdots \oplus A_m$;
\item[\textup{(iv)}] $\G_m = x^\R\cdot\exp(A_{0}\oplus\cdots\oplus A_{m-1})$ and
$x^\R\cap \exp(A_{0}\oplus\cdots\oplus A_{m-1})=\{1\}$.
\end{enumerate}
\end{lemma}

\begin{cor}\label{appendixTc} $\{f\in \T_{\exp}:\ \supp f \succ 1\}\ =\ \bigoplus_{m=0}^\infty A_m$, and
$$\G^{\operatorname{E}}\ =\ x^\R\cdot \exp\left(\bigoplus_{m=0}^\infty A_m\right),\  \quad
x^\R\cap \exp\left(\bigoplus_{m=0}^\infty A_{m}\right)=\{1\}.
$$
\end{cor}

\begin{lemma}\label{appendixTT} $ \quad \T_{\exp}^{>}\ =\  x^{\R}\cdot\exp(\T_{\exp})$.
\end{lemma}
\begin{proof} Let $f\in E_m^{>}$. By (iv) of Lemma~\ref{appendixT1}, 
$f=cx^r\exp(a)(1+\delta)$ with $c\in \R^{>}$, $r\in \R$, $a\in A_0 + \cdots + A_{m-1}$, and infinitesimal $\delta$ in the Hahn field
$E_m=\R[[\G_m]]$. Since $1+\delta=
\exp\left(\sum_{\nu=1}^\infty (-1)^{\nu-1}\delta^\nu/\nu\right)\in \exp(E_m)$, we get $f\in x^{\R}\exp(E_m)$.
\end{proof}

\noindent
It is easy to check that $x\notin \exp(\T_{\exp})$, so we are still missing $\log x$. Next we show that copying the above procedure with
$\log x$ instead of $x$, and then with $\log \log x$, and so on, and taking a union, is enough to enlarge the exponential ordered field~$\T_{\exp}$ to a logarithmic-exponential ordered field $\T$.

\subsection*{From $\T_{\exp}$ to $\T$} The idea is to use
distinct symbols $\ell_0,\ \ell_1,\ \ell_2,\ \dots$ in the role of~$x,\ \log x,\ \log \log x,\ \dots$.
Replacing for any given $n$ the formal variable $x$ in $\RxE$ by $\ell_n$ changes informally any $\operatorname{E}$-series $f(x)\in \RxE$ into
a series $f(\ell_n)\in \RlognE$. Formally: take for each~$n$ an isomorphism 
$$\fm \mapsto \fm{\downarrow}{}_n\ :\ \G^{\operatorname{E}} \to \G^{\operatorname{E},n}$$ of (multiplicative) ordered abelian groups,
with $x^r{\downarrow}{}_n=\ell_n^r\in \G^{\operatorname{E},n}$ 
for $r\in \R$. Given~$n$, this isomorphism extends 
uniquely to a strongly additive $\R$-linear map $$f \mapsto f{\downarrow}{}_n\ :\ \R[[\G^{\operatorname{E}}]] \to \R[[\G^{\operatorname{E},n}]].$$
This map is the identity on $\R$ and is an isomorphism of ordered (Hahn) fields; we denote the image of $\RxE$ under this isomorphism by $\RlognE$, and make the ordered subfield
$\RlognE$ of the ordered Hahn field $\R[[\G^{\operatorname{E},n}]]$
into an exponential ordered field, with exponentiation
denoted also by $\exp$, in such a way that 
$$f\mapsto f{\downarrow}{}_{n}\ :\ \RxE \to \RlognE$$
is an isomorphism 
of exponential ordered fields. Also $\ell_n:= \ell_n^1\in \ell_n^{\R}$ by notational convention. For $n=0$ we take
$\G^{\operatorname{E},0}=\G^{\operatorname{E}}$, with $\fm{\downarrow}{}_0=\fm$
for $\fm\in \G^{\operatorname{E}}$. Thus $x^r=\ell_0^r$ for $r\in \R$, and $f{\downarrow}{}_0=f$ for $f\in \R[[\G^{\operatorname{E}}]]$.  
Given $n$, we have the increasing sequence $\left(\G_m{\downarrow}{}_{n}\right)_{m=0}^\infty$ of convex subgroups of $\G^{\operatorname{E},n}$ with 
$\G^{\operatorname{E},n}=\bigcup_m G_m{\downarrow}{}_{n}$, and likewise, $\RlognE=\bigcup_m \R[[G_m{\downarrow}{}_{n}]]$.

\medskip\noindent 
It is straightforward to define inductively a strongly additive
$\R$-linear embedding
$$ \RxE\ \to\ \RxE$$
of exponential ordered fields that sends $x^r$ to $\exp(rx)$ for each $r\in \R$, and show that these properties define the embedding uniquely. (It maps  
$\G_m$ into $\G_{m+1}$ and $E_m$ into~$E_{m+1}$.) By transport to isomorphic copies of $\RxE$ we have
for each $n$ a unique strongly additive
$\R$-linear embedding $ \RlognE\ \to\ \RlognplusE$ of exponential ordered fields that sends~$\ell_n^r$ to $\exp(r\ell_{n+1})$ for each $r\in \R$; it maps~$\G_m{\downarrow}{}_{n}$ into $\G_{m+1}{\downarrow}{}_{n+1}$ and thus 
$\G^{\operatorname{E},n}$ into  $\G^{\operatorname{E},n+1}$. We identify $\RlognE$ with its image in~$\RlognplusE$ under this embedding. So $\ell_n^r=\exp(r\ell_{n+1})$ for $r\in \R$, and $\G_m{\downarrow}{}_{n}\subseteq \G_{m+1}{\downarrow}{}_{n+1}$, and we have inclusions 
\begin{align*} \T_{\exp}\ &=\ \R[[\ell_0^\R]]^{\operatorname{E}}\ \subseteq\ \R[[\ell_1^\R]]^{\operatorname{E}}\ \subseteq\ \R[[\ell_2^\R]]^{\operatorname{E}}\ \subseteq \cdots,\\
\G^{\operatorname{E}}\ &=\ \G^{\operatorname{E},0}\ \subseteq\  \G^{\operatorname{E},1}\ \subseteq\ \G^{\operatorname{E},2}\ \subseteq \cdots
\end{align*}
of exponential ordered fields and ordered abelian groups. We
now set
$$\T\ =\ \RxLE\ :=\ \bigcup_n \RlognE,\qquad 
\G^{\operatorname{LE}}\ :=\ \bigcup_n \G^{\operatorname{E},n}\ \subseteq\ \T,$$
with $\T$ construed as an exponential ordered field having the
$\RlognE$ as exponential ordered subfields, and
$\G^{\operatorname{LE}}$ construed as a multiplicative ordered abelian group with the $\G^{\operatorname{E},n}$ as ordered subgroups. In particular, we have the ordered Hahn field~$\RGLE$ with
the ordered field inclusion $\T \subseteq  \RGLE$. We also consider $\T$ as a {\em valued\/} subfield of the Hahn field $\RGLE$. 
Continuing to denote the exponentiation of $\T$ by $\exp$, we have $\exp(\epsilon)=\sum_{\nu=0}^\infty \epsilon^\nu/\nu!$ for infinitesimal $\epsilon\in \T$. From $x^r=\exp(r\ell_1)$ for $r\in \R$, and Lemma~\ref{appendixTT} we obtain  
$$\bigl(\RxE\bigr)^{>}\ \subseteq\ \exp\!\left(\R[[\ell_1^{\R}]]^{\operatorname{E}}\right),$$
and likewise $\bigl(\RlognE\bigr)^{>} \subseteq \exp\!\left(\RlognplusE\right)$ for all $n$. Thus $\T$ is a logarithmic-exponential ordered field. For the inverse $\log\colon \T^{>}\to \T$ of the exponentiation of~$\T$ it is now literally true that $\ell_1=\log x$, $\ell_2= \log \log x$, and so on. 
As in any logarithmic-exponential field we set $a^f:=\exp(f\log a)$ for $a\in \T^{>}$ and $f\in \T$. Thus for $f\in \T$ we have
$\ex^f=\exp(f)$, and so we use $\ex^f$ as an alternative notation for $\exp(f)$. (The above identification 
$\ell_n^r=\exp(r\ell_{n+1})$ for real $r$ agrees with
this definition of powers.)  The elements of $\G^{\operatorname{LE}}$ are the {\bf transmonomials} (or $\operatorname{LE}$-monomials). For $f\in\T^\times$ we have the dominant monomial $\fd(f)\in \G^{\operatorname{LE}}$.
\index{transmonomials} \index{group!transmonomials} \index{LE-monomials@$\operatorname{LE}$-monomials} \index{monomial!dominant} \index{dominant!monomial}

For any $f\in \T$ and $S\subseteq \G^{\operatorname{LE}}$ the subseries
$f|_{S}:= \sum_{\fm\in S}f_{\fm}\fm\in \R[[\G^{\operatorname{LE}}]]$ also lies in $\T$: to see this, one first observes this holds for each $E_m$ instead of $\T$,  and thus for each $\R[[\ell_n^{\R}]]$.
In particular, $\T$ is a truncation closed subfield of $\R[[\G^{\operatorname{LE}}]]$. 

\subsection*{The valuation of $\T$} Note that the valuation ring of $\T$ is 
$\{f\in \T:\ \supp f \preceq 1\}$, which is also the convex hull of $\R$ in the ordered field $\T$. We make the value group~$\Gamma_{\T}$ of the valuation into an ordered vector space over $\R$ by 
$r\gamma:= v(g^r)$ for $r\in \R$, $\gamma\in \Gamma_{\T}$ and
$g\in \T^{>}$ with $vg=\gamma$. Setting $A:=\{f\in \T: 
\supp f \succ 1\}$ we have the internal direct sum decomposition 
$$\T\ =\ A\oplus \R \oplus \smallo_{\T}$$ of $\T$
into $\R$-linear subspaces. 
The valuation has a concrete realization, based on the interesting fact that $\exp(A)\ =\ \G^{\operatorname{LE}}$.
(Proof of this fact: by Corollary~\ref{appendixTc} we have
$$\exp(A^{\operatorname{E}})\ =\ \G^{\operatorname{E}},\ \text{ where } A^{\operatorname{E}}\ :=\ \big\{f\in \RxE:\ \supp f \succ 1\big\} + \R \ell_1.$$ Use the natural analogues of this for $\RlognE$,~$\G^{\operatorname{E},n}$ in the role of $\RxE$,~$\G^{\operatorname{E}}$.)
Thus the surjective map $f\mapsto -\log \fd(f)\colon \T^\times \to A$ can serve as the valuation of $\T$, with the ordered
subgroup $A$ of $\T$ as value group, its structure as $\R$-linear subspace of $\T$ agreeing with the earlier defined vector space structure on the value group. Note also that $\exp(A) = \G^{\operatorname{LE}}$ gives that if $\fm\in \G^{\operatorname{LE}}$ and $r\in \R$, then $\fm^r\in \G^{\operatorname{LE}}$. 
 
We saw in the construction of $\T_{\exp}$ that the sequence $x, \ex^x, \ex^{\ex^x}, \ex^{\ex^{\ex^x}},\dots$ in $\T_{\exp}$ is strictly increasing and cofinal in it. Hence for each $n$ the analogous sequence
$$\ell_n, \ell_{n-1},\dots, \ell_1, x, \ex^x, \ex^{\ex^x},\dots
$$ in $\RlognE$ is strictly increasing and cofinal in $\RlognE$.
Thus the sequence $$x, \ex^x, \ex^{\ex^x}, \ex^{\ex^{\ex^x}},\dots$$ is even cofinal in $\T$. 
 This argument also shows that the sequence $\ell_0, \ell_1, \ell_2, \ell_3,\dots$ in~$\T^{>\R}$ is strictly decreasing. We claim that it is coinitial in $\T^{>\R}$. To see this, we note first that $[v(x)]$ is the smallest nontrivial
archimedean class of the value group of~$\RxE$, with
$\big[v(x)\big] < \big[v(\ex^x)\big]$. Hence $\big[v(\ell_{n})\big]$ is the smallest nontrivial ar\-chi\-me\-dean class of the value group of $\RlognE$, and $\big[v(\ell_n)\big] < \big[v(\ell_{n-1})\big]$ when $n\ge 1$.
In this way we get a sequence
$$\big[v(\ell_0)\big]\ >\ \big[v(\ell_1)\big]\ >\ \big[v(\ell_2)\big]\ >\ \cdots$$
of archimedean classes of the value group of $\T$, which
is moreover coinitial in the set of nontrivial archimedean classes of this value group. It remains to note that we are dealing with a convex valuation on the ordered field $\T$.

\subsection*{Representing $\T$ as a directed union of Hahn fields} This was not done in~\cite{DMM2} in the strong form we need for the
hypothesis in Theorem~\ref{thm:newtonianity of directed union} (in view
of its Corollary~\ref{cor:thmnewtonianity}), so we give full details here.

 By construction,
$\T$ is an increasing union of increasing unions of Hahn fields over~$\R$. To represent it as as a directed (and also as an increasing) union of Hahn fields over~$\R$,
we recall that $\G_m{\downarrow}{}_{n}\subseteq G_{m+1}{\downarrow}{}_{n+1}$. Also $G_m{\downarrow}{}_{n}\subseteq G_{m+1}{\downarrow}{}_{n}$, and from these two kinds of inclusions
we easily obtain $G_m{\downarrow}{}_n\subseteq G_{2\nu}{\downarrow}{}_{\nu}$ for $\nu=\max(m,n)$; moreover, $G_{2m}{\downarrow}{}_m\subseteq G_{2n}{\downarrow}{}_{n}$ for all $m\le n$. 
Thus the countable family $(G_m{\downarrow}{}_{n})_{m,n}$
of ordered subgroups of $G^{\operatorname{LE}}$ is directed, with
$G^{\operatorname{LE}}=\bigcup_{m,n}G_m{\downarrow}{}_{n}$.
This family is
in particular a healthy family (with respect to the coefficient field $\R$), and 
$\T\ =\ \bigcup_{m,n}\R[[G_m{\downarrow}{}_{n}]]$, and so the notion of a strongly additive map $\T \to \T$ makes sense.
We even have an increasing sequence $(G_{2n}{\downarrow}{}_{n})_n$ of ordered subgroups of $G^{\operatorname{LE}}$ with
$\T\ =\ \bigcup_n \R[[G_{2n}{\downarrow}{}_{n}]]$. Note also that
$$\R[[G_m{\downarrow}{}_{n}]]\ =\ 
\R[[G_m]]{\downarrow}{}_{n}\ =\ E_m{\downarrow}{}_{n}.$$

\subsection*{The upward shift operator}  A very useful automorphism of $\T$ is the {\bf upward shift}, informally to be thought of as~$f(x)\mapsto f(\ex^x)$.  \index{upward!shift} \index{shift!upward} Formally it is the 
 unique strongly additive $\R$-linear automorphism~$f\mapsto f{\uparrow}$ of the exponential ordered field~$\T$ that sends $x$ to $\ex^x$. (To construct it, take
for each $n$ the strongly additive $\R$-linear embedding $\RlognE \to \RlognE$
of exponential ordered fields that sends $\ell_n^r$ to 
$\exp(r\ell_n)$ for each $r\in \R$, and show that these maps have a common extension to a map $\T\to \T$; this common extension is the upward shift.) It is easy to check that $$\G^{\operatorname{E}}{\big\uparrow}\ \subseteq\ \G^{\operatorname{E}}\  \subseteq\ 
\G^{\operatorname{E},1}{\big\uparrow},$$ so
the upward shift maps $\G^{\operatorname{LE}}$ onto itself. The inverse of the upward shift operator is the downward shift operator
$f\mapsto f{\downarrow}$.
\index{downward!shift} \index{shift!downward}
The $n$th iterate of the upward (respectively, downward) shift operator is $f\mapsto f{\uparrow^n}$, respectively, 
$f\mapsto f{\downarrow}{}_{n}$. Thus $x{\downarrow}{}_{n}=\ell_n$. If~$f\in \RxE$, then this $f{\downarrow}{}_{n}$ equals  $f{\downarrow}{}_{n}$ as defined in the subsection ``From $\T_{\exp}$ to~$\T$.'' Thus $\big(\RlognE\big){\big\uparrow^n}=\T_{\exp}$. 

\begin{lemma}\label{appendixTup} If $0<g\in E_m{\downarrow}{}_n$, then $\log g\in E_{m+1}{\downarrow}{}_{n+1}$.
\end{lemma}
\begin{proof} Applying ${\uparrow^n}$, this reduces to the case
$n=0$. So let $0<g\in E_m$; we have to show that $\log g\in E_{m+1}{\downarrow}$. By (iv) of Lemma~\ref{appendixT1} we have 
$g=cx^r\exp(a)(1+\delta)$ with $c\in \R^{>}$, $r\in \R$, $a\in A_0+ \cdots + A_{m-1}$, and $\delta\prec 1$ in $E_m$. Then
$\log g=\log c + r\ell_1 + a + \log(1+\delta)$; now use that $a, \log(1+\delta)\in E_{m}\subseteq E_{m+1}{\downarrow}$.
\end{proof}

\subsection*{Differentiating and integrating transseries} In the rest of this appendix we also assume familiarity with Section~\ref{sec:H-fields}. It would be nice if transseries $f=f(x)$ could be differentiated so that the following rules hold: 
\begin{list}{*}{\setlength\leftmargin{2.5em}}
\item[(D1)] $r'=0$ for all $r\in \R$, and $x'=1$;
\item[(D2)] $(\exp f)'=f'\exp f$ for all $f\in \T$ (and thus $(\log f)'=f^\dagger$ for all $f\in \T^{>}$);
\item[(D3)] $f\mapsto f'\colon \T \to \T$ is strongly additive.
\end{list}
The inductive construction of $\T$ makes it plausible that this can be done uniquely:

\begin{prop} There is a unique
derivation on $\T$ satisfying {\rm(D1)},~{\rm(D2)},~{\rm (D3)}. 
\end{prop}

\noindent
Section~3 of \cite{DMM2} constructs a derivation on $\T$, denoted by $\der$ below, satisfying {\rm(D1)}, {\rm(D2)},
and~{\rm (D3)}. (The paper cited denotes $\der$ by $d/d x$, since $\der$ is thought of as differentiation with respect to $x$.) These properties clearly determine $\der$ uniquely, and
below we consider $\T$ as a differential field whose derivation is $\der$. 

The map
$\der\colon \T \to \T$ is healthy. This is because
$\der(E_{m}{\downarrow}{}_{n})\subseteq E_{m}{\downarrow}{}_{n}$ for all~${m\ge n}$. To prove these inclusions, note first that by the construction of $\der$
each $E_m$ is a differential subfield of $\T$ (so
$\T_{\exp}$ is as well). Next, with $\exp_n=\text{$n$th iterate of $\exp$,}$ 
$$\exp_n(x)'\ =\ \prod_{i=1}^n \exp_i(x), \qquad (f{\uparrow^n})'\ =\  f'{\uparrow^n} \exp_n(x)'  \quad(f\in \T),$$
by \cite{DMM2}. In view of $\exp_n(x)\in E_m$ for 
$m\ge n$, this yields the desired inclusion.

A consequence of the strong additivity of the derivation is that for an ordinary power series $F\in \R[[t_1,\dots, t_n]]$ and $\varepsilon=(\varepsilon_1,\dots, \varepsilon_n)\in \smallo_{\T}^n$ we have
$$F(\varepsilon)'\ =\ \sum_{i=1}^n \frac{\partial F}{\partial t_i}(\varepsilon)\cdot \varepsilon_i'.$$
 Here are some basic properties
of $\T$ as an ordered valued differential field:

\begin{prop} $C_{\T}=\R$, and $\T$ is a Liouville closed $H$-field.
\end{prop}

\noindent
Theorem~3.9(4) in \cite{DMM2} gives $C_{\T}=\R$, which together with 
Proposition~4.3 in that paper makes $\T$ an $H$-field.
The surjectivity of $\der\colon \T \to \T$ is
\cite[Theorem~5.6]{DMM2}. Hence~$\T$ is Liouville closed: given $a\in \T$, take $b\in \T$ with $b'=a$; then $(\ex^b)^\dagger=a$.  

%The upper shift is not an automorphism of
%$\T$ as a differential field: $f'{\uparrow}=\ex^{-x} %f{\uparrow}'$ for $f\in \T$. So $\der$ corresponds under 
%${\uparrow}$ to its compositional conjugate $\ex^{-x}\der$. 

The $E_m$ and thus the $E_{2n}{\downarrow}{}_{n}$ are spherically complete $H$-subfields of $\T$. Moreover, $\big[v(x)\big]$ is the smallest element of $\big[\Gamma_{E_m}^{\ne}\big]$, so $E_m$ is grounded with $\max \Psi_{E_m}=v(x^\dagger)=v(1/x)$.
Likewise, $\big[v(\ell_{n})\big]$ is the smallest elements of~$\big[(E_{m}{\downarrow}{}_{n})^{\ne}\big]$, so 
$E_{2n}{\downarrow}{}_n$ is grounded with
$\max \Psi_{E_{2n}{\downarrow}{}_n}=v(\ell_n^\dagger)$. 
Thus $\T=\bigcup_n E_{2n}{\downarrow}{}_{n}$
represents $\T$ as an increasing union of spherically complete
grounded $H$-subfields.

These facts also lead to an alternative proof that 
$\der\colon \T \to \T$ is surjective: using Corollary~\ref{cor:K=derK+C} it suffices to note that
$\ell_n^\dagger=\ell_{n+1}'\in \der(E_{n+1}{\downarrow}{}_{n+1})$.

\subsection*{Composition} Let $f$ range over $\T$ and $g$,~$h$ over $\T^{>\R}$. With this
convention:

\begin{prop} There is a unique operation
$$ (f,g)\mapsto f\circ g\ :\ \T\times \T^{>\R}\to \T$$ 
such that for each $g$ the following conditions are satisfied: \begin{enumerate}
\item[\textup{(i)}] $x\circ g=g$;
\item[\textup{(ii)}] $f\mapsto f\circ g\colon \T\to \T$ is an $\R$-linear embedding of exponential ordered fields;
\item[\textup{(iii)}] $f\mapsto f\circ g\colon \T \to \T$ is strongly additive.
\end{enumerate}
\end{prop}

\noindent
This is Theorem~6.2 in \cite{DMM2}, except that the
uniqueness is stated there with further conditions than just 
(i), (ii), (iii) above. It is rather easy to check, however, that there can be at most one operation as in the above proposition. Thus
the other conditions of that Theorem 6.2, namely (iv), (v), (vi) below, are consequences: \begin{enumerate}
\item[(iv)] $f\circ x=f$; 
\item[(v)] $f{\uparrow}=f\circ \ex^x$;
\item[(vi)] $(f\circ g)\circ h= f\circ(g\circ h)$;
\item[(vii)] $(f\circ g)'=(f'\circ g)\cdot g'$. 
\end{enumerate}  
The last item here, the Chain Rule, is part of Proposition~6.3 in~\cite{DMM2}. Note also that $\ex^x\circ g= \ex^g$, $\ell_1\circ g= \log g$ and
$f{\downarrow}=f\circ \ell_1$. If $f$ and $g$ are thought of as series~$f(x)$ and~$g(x)$, then $f\circ g$ can be thought of as $f\big(g(x)\big)$. For example, in Section~\ref{Compositional Conjugation} the identities $f{\uparrow}=f(\ex^x)$ and $f{\downarrow}=f(\log x)$  use this suggestive notation. 

\begin{prop} Given $g$ the map $f\mapsto f\circ g$ is healthy. More precisely, assume that
$f\in E_m{\downarrow}{}_n$ and $g\in E_{p}{\downarrow}{}_q$, $p,q\in \N$. Then
$f\circ g\in E_{m+n+p}{\downarrow}{}_{n+q}$.
\end{prop}
\begin{proof} We have $f\circ g = f{\uparrow^n}\circ \log^n g$,
with $f{\uparrow^n}\in E_m$, and $\log^n g\in 
E_{p+n}{\downarrow}{}_{q+n}$ by Lemma~\ref{appendixTup}. It remains to use subsection 6.8 in \cite{DMM2}. 
\end{proof}

\subsection*{The subfield $\T_{\log}$ of logarithmic transseries} This subfield is a particularly transparent
part of $\T$, more so than $\T_{\exp}$. It also has much stronger algebraic
closure properties than $\T_{\exp}$. Let $\mathfrak{L}_n:= \ell_0^{\R}\cdots \ell_n^{\R}$ be the subgroup of $\G^{\operatorname{LE}}$ generated by the real powers of the $\ell_i$ for 
$i=0,\dots,n$. Then 
$\mathfrak{L}_n\subseteq \G_n{\downarrow}{}_{n}$, and $\R[[\mathfrak{L}_n]]$ is an $H$-subfield of 
$\R[[\G_n{\downarrow}{}_{n}]]=E_n{\downarrow}{}_{n}$.  The ordered group inclusions 
$$ \mathfrak{L}_0\ \subseteq\ \mathfrak{L}_1\ \subseteq\ \mathfrak{L}_2\ \subseteq\ \cdots \subseteq\ \mathfrak{L}\ 
:=\ \bigcup_n \mathfrak{L}_n\ \subseteq\ \G^{\operatorname{LE}}$$
induce $H$-field inclusions 
$$\R[[x^\R]]=\R[[\mathfrak{L}_0]]\ \subseteq\ \R[[\mathfrak{L}_1]]\ \subseteq\ \R[[\mathfrak{L}_2]]\ \subseteq\ \cdots\ 
\subseteq\ \T,$$
and we set $$\T_{\log}\ :=\ \bigcup_n \R[[\mathfrak{L}_n]].$$ So $\T_{\log}$ is an $H$-subfield of $\T$ with $\G^{\operatorname{LE}}\cap \T_{\log}=\mathfrak{L}$. It is routine to check that for
$f\in \T_{\log}^{>}$ we have $f^r\in \T_{\log}$ for all $r\in \R$, and $\log f\in \T_{\log}$. 

\index{logarithmic!transseries} \index{transseries!logarithmic}

\subsection*{Notes and comments} 
Analytic functions of infinitesimal arguments in a Hahn field were considered by Neumann~\cite[pp.~206--210]{Neumann}. For more on strongly additive maps between Hahn fields, see~\cite{vdH:noeth}.

The exponential field of transseries was introduced by Dahn and G\"oring~\cite{DG} as a natural candidate for an elementary extension of the exponential ordered field of real numbers, and independently in \'Ecalle's work~\cite{Ecalle1} on the Dulac Conjecture, where the derivation
is also prominent. In connection with the construction of~$\T$ as a directed union of Hahn fields, we note that by \cite{KKS}
no nontrivial Hahn field over $\R$ can be made into a logarithmic-exponential ordered field. 
%we remark that \cite{KKS} shows that, for any ordered abelian group~$\fM\neq\{1\}$, there is no  strictly increasing exponentiation $\exp$ on $K=\R[[ \fM]] $ with $\exp(K)=K^>$.

In this appendix we have followed closely \cite{DMM2}, but the reader who consults that paper will note that it allows any
logarithmic-exponential ordered field $\k$ as coefficient field
in constructing $\k\(( x^{- 1}\)) ^{\operatorname{LE}}$ while
here we stick to $\k=\R$. 
%(The notion of ``logarithmic-exponential ordered field'' is defined below.)
Apart from that, the exponential and logarithmic maps $\operatorname{E}$
and $\operatorname{L}$ there are denoted here by $\exp$ and $\log$, the fields $K_n$ and $L_n$ there
are our~$E_n$ and $\R[[x^{\R}]]^{\operatorname{E},n}$, and
$G_{m,n}$ and $L_{m,n}$ there are our $G_n{\downarrow}{}_m$ 
and $E_n{\downarrow}{}_m$; the map $\Phi$ there is our upward shift operator $\uparrow$.

\medskip\noindent
We wish to say a few words on the grid-based setting, which has
been mentioned several times in earlier \textit{Notes and comments.}\/
This appendix rests on the notion of a 
\textit{well-based}\/ subset of an ordered abelian (multiplicative) group $\fM$. We now replace \textit{well-based}\/ by a much stronger restriction, namely \textit{grid-based.}\/ A {\bf grid-based\/} subset of~$\fM$ is one that is contained in $\fm\, \fn_1^{\N}\cdots \fn_n^{\N}$
for some $\fm\in \fM$ and 
$\fn_1,\dots, \fn_n\in \fM^{\prec 1}$. Let~$C$ be a field. Set 
$$C\lgb\fM\rgb\ :=\ \big\{f\in C[[\fM]]:\ \text{$\supp f$  is grid-based}\big\}.$$
Then $C\lgb\fM\rgb$ is a subfield of $C[[\fM]]$. 
For $\fM=x^\Q$ we obtain the field  $C\lgb x^\Q\rgb=\operatorname{P}(C)$ of Puiseux series over $C$ (Example~\ref{ex:Puiseux}).

{\sloppy

Except for the statements mentioning spherical completeness, the appendix goes through if everywhere \textit{well-based}\/ is replaced by \textit{grid-based}\/, with the Hahn fields~$C[[\fM]]$ accordingly replaced by their grid-based versions $C\lgb\fM\rgb$. In particular, the starting point $\R[[x^{\R}]]$ in the construction of 
$\T$  is replaced by the subfield~$\R\lgb x^{\R}\rgb$ of $\R[[x^{\R}]]$. This leads to a logarithmic-exponential ordered subfield~$\T_{\operatorname{g}}$ of $\T$  which is also an $H$-subfield of~$\T$ and shares key algebraic closure properties with~$\T$ like being Liouville closed. 
Every $f\in \T_{\operatorname{g}}$ has grid-based support. 

}
%Grigor$'$ev-Singer~\cite{GrigorievSinger} studied
%algebraic differential equations over  $\R\lgb x^\R\rgb$.

{\sloppy
Grid-based transseries were introduced in \cite{Ecalle1}. 
There $\mathbb T_{\operatorname{g}}$ is denoted by $\R[[[x]]]$, in~\cite{JvdH} it is written as
$\mathbb T$, and in \cite{DMM2} as 
$\R\(( t\)) ^{\operatorname{LE}, \operatorname{ft}}$; in this last paper
\textit{grid-based}\/ was called \textit{of finite type.}\/ 
}

\index{grid-based!transseries} \index{transseries!grid-based} \index{support!grid-based}
 \index{subset!grid-based}

\medskip\noindent
Logarithmic transseries occur in the work of Loeb and Rota~\cite{Loeb,LoebRota} in connection with difference equations.

\medskip\noindent
As already noted, our field $\T$ of transseries does not contain the well-based series
$x + \log x + \log\log x + \cdots$. Van der Hoeven's thesis~\cite{vdH:PhD} shows how to go beyond~$\T$
to include such series by building a strictly increasing ordinal sequence of fields of transseries, alternately closing off under well-based summation and exponentiation,
and defining differentiation and composition in this setting. In Schmeling's thesis~\cite{Schm01} this is put in an axiomatic framework for 
Hahn fields with a logarithm map, leading to
fields~$T$ of transseries with additional structure, such as an {\em iterator\/} 
$\log_{\omega}\colon T^{>\R}\to T^{>\R}$ of the logarithm satisfying the identity $$\log_{\omega} \log y\ =\ (\log_{\omega} y) -1\qquad(y>\R).$$ 
This axiomatic setting for fields of transseries plays a role in
the recent
work by Berarducci and Mantova~\cite{BM}. We   
plan to consider these extensions in more detail in our second volume in the light of the results of the present volume.

%% file: mt-modth.tex
\let\thetheorem\thetheoremold

\chapter{Basic Model Theory} \label{app:modth}\label{app:modeltheory}

\noindent
This appendix is written for readers unfamiliar with model theory.
It provides a rigorous treatment, with examples, of
the model-theoretic tools used in our work:
back-and-forth, compactness, types, model completeness, quantifier elimination, and (in a different vein) NIP. 
We adopt a many-sorted setting from the outset. In the first four sections we 
specify the notion of a {\em model-theoretic structure\/} and
build up a formalism for handling these structures
efficiently in the rest of the appendix.

\section{Structures and Their Definable Sets} \label{sec:structures and definable sets}

\noindent 
The notion of structure we introduce in this section is very general and includes not only
familiar (one-sorted) algebraic objects such as groups and rings, but also two-sorted structures such as group actions and topological spaces. 

We associate to any structure its category of
\textit{definable sets and definable maps}.
In the case of a field,
the definable sets include the solution sets of finite systems of polynomial
equations, but as we shall see, the notion of definable set goes beyond this by allowing boolean operations and  projections.

\subsection*{Model-theoretic structures}
A (model-theoretic) {\bf structure} \index{structure} consists of: 
\begin{list}{}{\leftmargin=2.5em }
\item[(S1)] a family $(M_s)_{s\in S}$ of nonempty sets $M_s$,
\item[(S2)] a family $(R_i)_{i\in I}$ of relations
$R_i\subseteq M_{s_1}\times\cdots \times M_{s_m}$
on these sets, with the tuple $(s_1,\dots,s_m)\in S^m$ determined by the index $i\in I$, and
\item[(S3)] a family $(f_j)_{j\in J}$ of functions 
$f_j\colon M_{s_1}\times\cdots \times M_{s_n}\to M_{s}$
with the tuple $(s_1,\dots,s_n,s)\in S^{n+1}$ determined by the index~$j\in J$. 
\end{list}
Notation:
$$\mathbf M\ =\ \big( (M_s); (R_i), (f_j) \big).$$
The elements of the index set $S$ are called {\bf sorts}, \index{sorts}\index{structure!sorts} and $M_s$ is the {\bf underlying set}  of sort~$s$ of $\mathbf M$. 
Elements of $M_s$ are also called \textit{elements of $\mathbf M$ of sort~$s$.}\/
The $R_i$ and~$f_j$ are called the {\bf primitives} \index{structure!primitives}   \index{primitives} of $\mathbf M$. For $m=0$ the set $S^m$ has the empty tuple as its only element, and the product $M_{s_1}\times\cdots \times M_{s_m}$ is then a singleton, that is, a one-element set. For $n=0$ in (S3), the function~$f_j$ is then identified with its unique value in $M_{s}$ and is called a {\bf constant.} \index{structure!constant} \index{constant} The family $M=(M_s)$ is said to {\bf underlie} the structure $\mathbf M$, or, abusing language, to be the {\bf underlying set} of $\mathbf M$. \index{structure!underlying set}

We  call $\mathbf M$ as above an {\bf $S$-sorted structure}.
If $S$ is finite of size $n$, then we also speak of an {\bf $n$-sorted structure}.
% and if $S$ has
%more than two elements, of a {\bf many-sorted structure.}
Many texts only consider one-sorted structures, so~$S$ doesn't need to be mentioned. But the extra generality adds useful flexibility and is natural, since mathematical structures are often many-sorted to begin with. In fact, virtually anything that mathematicians consider as a structure can be viewed as a structure in the above sense. We only mention a few examples here; more are discussed in~\ref{ex:L-structures} below. 

In specifying one-sorted structures we often take the liberty of using the same capital letter for the underlying set of the structure as for the structure itself. 

\begin{examplesNumbered}\label{ex:structures}\mbox{}

\begin{enumerate} 
\item Any (additively written) abelian group $M$ is construed naturally as a one-sorted structure 
$(M;0,{-},{+})$ with the constant $0\in M$ and the functions
$$a\mapsto -a\colon M\to M, \quad (a,b)\mapsto a+b\colon M\times M\to M.$$
\item Every ring $R$ is viewed naturally as a one-sorted structure $(R;0,1,{-},{+},{\,\cdot\,})$ with
constants $0$, $1$ and functions 
$$\hskip-3em r\mapsto -r\colon R\to R, \ \ (r,s)\mapsto r+s\colon R\times R\to R,\  \
  (r,s)\mapsto r\cdot s\colon R\times R\to R. $$
\item Let $R$ be a ring  and $M$ an $R$-module. Then $M$ can be viewed as a one-sorted structure, equipped with the functions in (1) 
and for each $r\in R$ the function $a\mapsto ra\colon M\to M$.
\item
Sometimes it is convenient to construe a module~$M$ over a ring $R$ as a {\em two-sorted \/}structure, one sort with underlying set $R$ and the same functions as in~(2),
the other sort with underlying set $M$ and the same functions as in~(1), as well as with a function $(r,a)\mapsto ra\colon R\times M\to M$ for the action of $R$ on~$M$.
\item A nonempty ordered set in the sense of Section~\ref{sec:ordered sets} is the same thing as a one-sorted structure
$(M;{\leq})$ where the binary relation $\leq$ on $M$ is reflexive, antisymmetric, transitive, and total.
\item An incidence geometry is a two-sorted structure consisting of a nonempty set~$P$ whose elements are called points, a nonempty set $L$ whose elements are called lines, and a relation $R\subseteq P\times L$ between points and lines.  
For example, a nonempty topological space can be viewed as an incidence geometry, taking the underlying set of the space for $P$, the collection of open sets for $L$, and the membership relation between points and open sets for $R$.  
\end{enumerate}
\end{examplesNumbered}

\subsection*{Definable sets}
Let $\mathbf M=\big( M; (R_i), (f_j) \big)$  be an $S$-sorted structure as above, where $M=(M_s)$. The main role of the primitives of $\mathbf M$ is to generate the \textit{definable sets}\/ of $\mathbf M$. 
Let $\mathbf s=s_1\dots s_m$ and 
$\mathbf t=t_1\dots t_{n}$ be elements of $S^*$, that is, words on the alphabet $S$. 
Then $\mathbf s\mathbf t=s_1\dots s_m t_1\dots t_{n}\in S^*$ is the concatenation of $\mathbf s$ and $\mathbf t$. We set
$M_{\mathbf s}:=M_{s_1}\times\cdots\times M_{s_m}$, and identify $M_{\mathbf s}\times M_{\mathbf t}$ with $M_{\mathbf s\mathbf t}$ in the obvious way. 

\medskip\noindent
There is no need to require $M_s\cap M_{t}=\emptyset$ for distinct $s, t\in S$. (This would be
unnatural in Example (4) above, since $R$ is naturally an
$R$-module.) Instead we impose the convention
that in referring to an element $a\in M_{\mathbf s}$, this is short for a reference to an ordered pair $(a,\mathbf{s})$ such that $a\in M_{\mathbf s}$.
A similar convention is in force when we refer to a set $X\subseteq M_{\mathbf s}$. Given a map $f\colon X\to Y$, where $X\subseteq M_{\mathbf s}$, $Y\subseteq M_{\mathbf t}$, we let
$\Gamma(f)\subseteq M_{\mathbf s\mathbf t}$ be the graph of 
$f$. 

\begin{definition}\label{def:0-definable} \index{structure!definable set} \index{set!$0$-definable} \index{set!absolutely definable}
The {\bf $0$-definable} (or absolutely definable) sets of $\mathbf M$ are the relations $X\subseteq M_{\mathbf s}$
%, with $\mathbf s$ part of the specification 
obtained recursively as follows:
\begin{list}{}{\leftmargin=1em }
\item[(D1)] the relations $R_i\subseteq M_{s_1\dots s_m}$ and the graphs $\Gamma(f_j)\subseteq M_{s_1\dots s_n s}$ are $0$-definable;
\item[(D2)] if $X,Y\subseteq M_{\mathbf s}$ are $0$-definable, then so are $X\cup Y\subseteq M_{\mathbf s}$ and $(M_{\mathbf s}\setminus X)\subseteq M_{\mathbf s}$;
\item[(D3)] if $X\subseteq M_{\mathbf s}$ and $Y\subseteq M_{\mathbf t}$ are $0$-definable, then so is
$X\times Y\subseteq M_{\mathbf s\mathbf t}$;
\item[(D4)] for all distinct $i,j\in \{1,\dots,m\}$ with $s_i=s_j$,  
the diagonal 
$$\Delta_{ij}:=\big\{ (x_1,\dots,x_m)\in M_{\mathbf s}:x_i=x_j\big\}$$ is $0$-definable;
\item[(D5)] if $X\subseteq M_{\mathbf s\mathbf t}$ is $0$-definable, then so is  $\pi(X)\subseteq M_{\mathbf s}$, where $\pi\colon M_{\mathbf s\mathbf t}\to M_{\mathbf s}$ is the obvious projection map.
\end{list}
\end{definition}

\noindent
We extend this notion to \textit{$A$-definability.}\/ Here $A$ is a so-called {\bf parameter set} in $\mathbf M$, that is,  a family $A=(A_s)$ with $A_s\subseteq M_s$ for each~$s$ (shorthand: $A\subseteq M$).
Then the structure $\mathbf M_A$ is obtained from $\mathbf M$ by adding for each $a\in A_s$ the constant $a\in M_s$ as a primitive. The {\bf $A$-definable sets} of $\mathbf M$ are just the $0$-definable sets of~$\mathbf M_A$. 
The parameter set $A$ with $A_s=\emptyset$ for each $s$ is also denoted by $0$, so $\mathbf M_0=\mathbf M$
and the terminology \textit{$0$-definable}\/ is unambiguous.
Instead of \textit{$M$-definable}\/ we just write \textit{definable}\/; so all finite sets $X\subseteq M_{\mathbf s}$ are definable.

For any family $(A_i)_{i\in I}$ of parameter sets in 
$\mathbf M$ we let $\bigcup_i A_i$ be the parameter set~$A$ in $\mathbf M$ such that $A_s=\bigcup_i A_{i,s}$ for every $s\in S$.

\index{parameter set} \index{structure!parameter set} \index{set!parameter} \index{set!$A$-definable} \index{set!definable}\index{definable!set}

\subsection*{Describing definable sets}
For constructing new definable sets  from old ones,  it is convenient to systematize the correspondence between sets and the conditions used to define them. Let $x$ be a variable ranging over a set
$X$, and let $\phi(x)$, $\psi(x)$ be formulas (``conditions'' on $x$) defining the subsets  
$$\Phi := \big\{x \in X  : \text{$\phi (x)$ holds} \big\}\quad\text{and}\quad
  \Psi := \big\{x \in X  : \text{$\psi (x)$ holds} \big\}$$
of $X$. (We will make this precise in Section~\ref{sec:formulas} below.)
Various logical operations on formulas
correspond to operations on the sets that these formulas define:
\begin{list}{}{\leftmargin=2.5em } 
\item[(C1)] $\neg \phi (x)$ defines the complement $X\setminus\Phi$;
\item[(C2)] $\phi(x)\vee\psi(x)$ defines the union  $\Phi \cup \Psi$;
\item[(C3)] $\phi(x) \wedge \psi(x)$, also written as $\phi(x) \, \& \, \psi(x)$, defines the intersection $\Phi \cap \Psi$; thus $\neg\big(\phi(x)\vee \psi(x)\big)$ and
$\neg\phi(x)\wedge\neg\psi(x)$
define the same subset of $X$.
\end{list}
Moreover, if $y$ is a variable ranging over a set $Y$ and $\theta(x,y)$ is a formula defining the subset $\Theta$ of $X\times Y$, then
\begin{list}{}{\leftmargin=2.5em }
\item[(C4)] $\exists  y\,  \theta(x,y)$  defines   $\pi  (\Theta)
\subseteq X$  where $\pi\colon X\times Y\to X$ is the natural projection;
\item[(C5)]
$\forall y \,\theta(x,y)$ defines the set $\big\{  x \in X : \{ x \}
\times Y \subseteq \Theta \big\}$; hence the formulas $\neg \exists y\,  \theta(x,y)$ and $\forall y\, \neg\theta(x,y)$ define the same subset of $X$.
\end{list}
Thus, with $X=M_{\mathbf s}$ and  $Y=M_{\mathbf t}$, if the sets $\Phi$, $\Psi$ defined by $\phi(x)$, $\psi(x)$ are 
$A$-definable, then so are the sets defined by $\neg \phi(x)$,  $\phi(x)\vee\psi(x)$,
$\phi(x) \wedge \psi(x)$, and if
 $\Theta$ is $A$-definable, then so are the
sets defined by $\exists  y \, \theta(x,y)$  and $\forall y\, \theta(x,y)$.
The advantage of the logical formalism is that it is often more suggestive and transparent than
traditional set-theoretic notation, in particular when quantifiers are involved:

\begin{example} Let $A\subseteq \R$, and let $\mathbf M=(\R;\dots)$ be a one-sorted structure  with the usual ordering $\leq$ 
of the real line among the $A$-definable subsets of $\R^2$. Equip $\R^m$ with its usual product topology.
Suppose $X\subseteq \R^m$ is $A$-definable (in $\mathbf M$). Then the interior~$\operatorname{int}(X)$
of $X$ in the ambient space $\R^m$ is $A$-definable. To see this, note that, with $\phi  \to  \psi$ short for $(\neg \phi) \vee \psi$, 
we have the equivalence
$$ x \in  \operatorname{int}(X)\quad \Longleftrightarrow \quad
\exists  u\, \exists  v\,  \big[ ( u <  x <  v)\ \&\ \forall  y\, ( u <  y <  v \to  y \in X)\big]$$
where $ u$, $ v$, $ x$, $ y$ range over $\R^m$ and $ u< x< v$ is short for
$$u_1<x_1<v_1 \ \&\ \cdots\ \& \ u_m<x_m<v_m.$$
Likewise, if $f\colon X\to \R$ is $A$-definable, that is, its graph $\Gamma(f)$ is $A$-definable, then the set 
$\{x\in X:\text{$f$ is continuous at $x$}\}$ is easily seen to 
be $A$-definable. 
\end{example}

\noindent
In general the $0$-definable relations of a structure cannot be described in a way that is significantly more explicit than the  recursive definition in \ref{def:0-definable}. 
For example, this is the case for the ring $\Z$ viewed as  a one-sorted structure as in Example~\ref{ex:structures}(3);
see \cite[\S{}7.5]{Shoenfield} for details. But in some cases, a more explicit description does exist. One example is the field $\C$ of complex numbers; here we construe~$\C$ as a one-sorted structure as in \ref{ex:structures}(2). (Including also as a primitive, say, $x\mapsto x^{-1}\colon \C^\times\to\C$,
extended to a total function $\C\to\C$ by declaring $0^{-1}:=0$, wouldn't add to the $0$-definable relations of~$\C$.)
Then by the Chevalley-Tarski Constructibility Theorem from Section~\ref{sec:some theories} the $0$-definable subsets of $\C^n$ are just the
finite unions of sets  
$$\big\{ a\in \C^n:\  P_1(a) = \cdots = P_m(a) = 0,\ Q(a)\neq 0\big\}$$
with $P_1,\dots,P_m,Q\in \Q[X_1,\dots,X_n]$. For a subfield $A$ of $\C$ as a parameter set we get the same description of the $A$-definable subsets of $\C^n$ but now the polynomials have their coefficients in $A$. The above goes through for any algebraically closed field~$K$ instead of $\C$, with the prime field of $K$ in place of $\Q$. So in this case the notion of ``$A$-definable'' is akin to Weil's notion of an algebraic variety defined over~$A$. (Model-theoretic notions are often similar to foundational items in Weil's algebraic geometry.)
Another prominent example, more relevant to our work, where a concise description of the definable sets is available is the field $\R$ of real numbers.
Here, by the Tarski-Seidenberg Theorem in Section~\ref{sec:some theories}, the $0$-definable subsets of~$\R^n$ are the
\textit{$\Q$-semialgebraic}\/ subsets of~$\R^n$, that is, finite unions of sets
$$\big\{ a\in \R^n:\  P(a) = 0,\ Q_1(a) > 0, \dots, Q_m(a) > 0 \big\}$$
with $P,Q_1,\dots,Q_m\in \Q[X_1,\dots,X_n]$.
Note that the usual ordering relation on $\R$ is $0$-de\-finable in the field $\R$:
$$x\leq y\quad\Longleftrightarrow\quad \exists z\, [ z^2=y-x ].$$
In this appendix we develop some tools which allow us to 
identify  benign structures, like $\C$ and $\R$, 
whose definable sets allow such explicit descriptions. 
%The theorems of Chevalley-Tarski and Tarski-Seidenberg are proved in Section~\ref{sec:some theories}. 

\index{theorem!Chevalley-Tarski}
\index{theorem!Tarski-Seidenberg}

\subsection*{Conventions and notations} Let $M=(M_s)$ and $M'=(M'_s)$ be families of sets indexed by a set $S$. 
The {\bf size}\/~$|M|$ of~$M$ is defined to be the sum 
$\sum |M_s|$ of the cardinalities of the sets $M_s$, in other words, it is the cardinality of the disjoint union 
$\bigcup_s M_s\times\{s\}$. (We also use {\em size\/} as a synonym for {\em cardinality}.) Define $M\subseteq M':\Leftrightarrow M_s\subseteq M'_s$ for all $s$.
A map $h\colon M\to M'$ is a family $h=(h_s)$ of maps $h_s\colon M_s\to M'_s$.
Let $h=(h_s)\colon M\to M'$ be such a map.
We say that~$h$ is \textit{injective}\/ if each $h_s$ is injective, and similarly with \textit{surjective}\/ or \textit{bijective}\/ in place of \textit{injective.}\/
Given $\mathbf s=s_1\dots s_m\in S^*$ and  $a=(a_1,\dots,a_m)\in M_{\mathbf s}$ we put $ha=h_{\mathbf s}a:=(h_{s_1}a_1, \dots , h_{s_m}a_m)\in M'_{\mathbf s}$. For $A=(A_s)\subseteq M$ we put $h(A):=\big(h(A_s)\big)\subseteq M'$.
Given also $M'':= (M''_s)$ and a map $h'\colon M'\to M''$, we let $h'\circ h$ denote the map $M\to M''$ given
by $(h'\circ h)_s=h'_s\circ h_s$ for all~$s$.

\subsection*{Notes and comments}
The idea of ``structure'' emerged in the 19th century in algebra 
(Dedekind) and the foundations of 
mathematics~\cite{Hilbert99,Schroeder}. It led to the structural point of view in algebra, later extended by Bourbaki to other parts of mathematics; see~\cite{Corry}.
In the early literature in mathematical logic, structures were called  systems until 
the current terminology came into widespread use in the 1950s (see for example \cite{Robinson52}), perhaps under Bourbaki's influence. 

The underlying set $M_s$ of sort $s$ of a structure $\mathbf M$ 
is sometimes called the \textit{universe}\/ of sort $s$ of $\mathbf M$, 
after de Morgan~\cite{deMorgan46}. 
The notion of definable set can be traced back to Weyl~\cite{Weyl}, but became more familiar after Kuratowski and Tarski~\cite{KT31} spelled out the above
correspondence between logical connectives and set-theoretic operations for 
routine use in {\em descriptive set theory}.

\section{Languages}\label{sec:languages}

\noindent
In this section we formalize the idea of {\em structures of the same kind}\/ by defining the category of \textit{$\mathcal L$-structures}\/ for a given (model-theoretic) \textit{language $\mathcal L$.}\/ We also discuss various constructions
with $\mathcal L$-structures: products, direct limits.

\subsection*{Languages} \index{language} \index{sorts} \index{language!sorts} \index{language!relation symbols} \index{language!function symbols} \index{language!constant symbols} \index{symbols!function} \index{symbols!relation} 
A {\bf language} $\mathcal L$ is a triple $(S,\mathcal L^{\operatorname{r}},\mathcal L^{\operatorname{f}})$ consisting of
\begin{list}{}{\leftmargin=2.5em}
\item[(L1)] a set $S$ whose elements are the {\bf sorts} of $\mathcal L$,
\item[(L2)] a set $\mathcal L^{\operatorname{r}}$ whose elements are the {\bf relation symbols} of $\mathcal L$, and
\item[(L3)] a set $\mathcal L^{\operatorname{f}}$ whose elements are the {\bf function symbols} of $\mathcal L$,
\end{list}
where $\mathcal L^{\operatorname{r}}$ and $\mathcal L^{\operatorname{f}}$  are disjoint, each $R \in \mathcal L^{\operatorname{r}}$  is equipped with a word $s_1\dots s_m\in S^*$, called its {\bf sort},
and each $f \in \mathcal L^{\operatorname{f}}$ is equipped with a {\bf sort} $s_1\dots s_n s\in S^*$. 
The elements of $\mathcal L^{\operatorname{r}}\cup \mathcal L^{\operatorname{f}}$ are  called {\bf nonlogical symbols} of~$\mathcal L$.  \index{symbols!nonlogical}
We also call~$\mathcal L$ an {\bf $S$-sorted language}.
% with the same linguistic conventions as for structures if~$S$ is finite.
It is customary to present a  language  
as a disjoint union $\mathcal L=\mathcal L^{\operatorname{r}}\cup\mathcal L^{\operatorname{f}}$ 
of its sets of relation and function symbols,
while separately specifying the set $S$ of sorts as well as the kind (relation symbol or function symbol) and sort of each nonlogical symbol of~$\mathcal L$.
If $R\in \mathcal L^{\operatorname{r}}$ has sort $s_1\dots s_m\in S^*$, then we also call $R$
an {\bf $m$-ary} relation symbol of $\mathcal L$;  similarly  a function symbol of $\mathcal L$ of sort $s_1\dots s_n s$
is called an $n$-ary function symbol of $\mathcal L$.
Instead of \textit{$0$-ary}\/, \textit{$1$-ary}\/, \textit{$2$-ary}\/ we say \textit{nullary}\/, \textit{unary}\/, \textit{binary}\/, respectively. A {\bf constant symbol}
is a nullary function symbol.  An $m$-ary relation symbol of sort $s\dots s$ is also said to be of sort $s$;
an $n$-ary function symbol of sort $s\dots ss$  is said to be of sort~$s$. \index{symbols!constant}

\index{language!extension}\index{extension!language}\index{sublanguage}\index{language!sublanguage}

A language
$\mathcal L' = (S',\mathcal L'^{\operatorname{r}},\mathcal L'^{\operatorname{f}})$
is an {\bf extension} of a language
$\mathcal L = (S,\mathcal L^{\operatorname{r}},\mathcal L^{\operatorname{f}})$ 
(and~$\mathcal L$ is a {\bf sublanguage} of $\mathcal L'$)
if $S \subseteq S'$, $\mathcal L^{\operatorname{r}}\subseteq \mathcal L'^{\operatorname{r}}$,  
 $\mathcal L^{\operatorname{f}}\subseteq \mathcal L'^{\operatorname{f}}$, and each nonlogical symbol of $\mathcal L$  has the same sort in $\mathcal L'$ as it has in $\mathcal L$;
 notation: $\mathcal L\subseteq\mathcal L'$.

\begin{examplesNumbered}\label{ex:languages}
\mbox{}

\begin{enumerate}
\item The one-sorted language $\mathcal L_{\operatorname{G}}=\{ 1, {\ }^{-1}, {\, \cdot\,} \}$ of 
groups has constant
symbol $1$, unary function symbol ${\ }^{-1}$, and binary function symbol $\cdot\,$. %\nomenclature[Bex1]{$\mathcal L_{\operatorname{Gr}}$}{language of groups}

\item The one-sorted language $\mathcal L_{\operatorname{A}}=\{ 0, {-}, {+} \}$ of 
(additive) abelian groups has constant
symbol $0$, unary function symbol $-$, and binary function symbol $+$. %\nomenclature[Bex2]{$\mathcal L_{\operatorname{Ab}}$}{language of abelian groups}

\item The one-sorted language $\mathcal L_{\operatorname{R}}=\{ 0, 1, {-}, {+}, {\,\cdot\,}\}$ of rings
is the extension of $\mathcal L_{\operatorname{A}}$ by a constant symbol $1$ and a binary function symbol $\cdot\,$. \nomenclature[Bex3]{$\mathcal L_{\operatorname{R}}$}{language of rings}

\item The one-sorted language $\mathcal L_{\operatorname{O}}=\{ {\leq} \}$ of ordered sets has just one binary relation symbol $\leq$. \nomenclature[Bex4]{$\mathcal L_{\operatorname{O}}$}{language of ordered sets}

\item Combining (2) and (4) yields the one-sorted language $\mathcal L_{\operatorname{OA}}=\{{\leq},0, {-}, {+} \}$ of ordered abelian groups. Combining (3) and (4) yields
the one-sorted language $\mathcal L_{\operatorname{OR}}=\{ {\leq}, 0, 1, {-}, {+}, {\,\cdot\,}\}$
of ordered rings. \nomenclature[Bex51]{$\mathcal L_{\operatorname{OA}}$}{language of ordered abelian groups} %\nomenclature[Bex52]{$\mathcal L_{\operatorname{OR}}$}{language of ordered rings}
\item Let $R$ be a ring. The language $\mathcal L_{R\!\operatorname{-mod}}=
\mathcal L_{\operatorname{A}}\cup\{ {\lambda_r} \colon r\in R\}$ of $R$-modules is the one-sorted language which
extends $\mathcal L_{\operatorname{A}}$ by unary function symbols $\lambda_r$, one for each $r\in R$.
%\nomenclature[Bex6]{$\mathcal L_{\text{$R$-mod}}$}{language of $R$-modules}

\item The language $$\mathcal L_{\operatorname{Mod}}=\{ 0_{\operatorname{R}},1_{\operatorname{R}},{-_{\operatorname{R}}},{+}_{\operatorname{R}},{\,\cdot\,}_{\operatorname{R}},
0_{\operatorname{M}},{-}_{\operatorname{M}},{+}_{\operatorname{M}},\lambda\}$$ of modules (with unspecified scalar ring)
is a two-sorted language with sorts~$\operatorname{R}$ and~$\operatorname{M}$.
 The symbols
$0_{\operatorname{R}}$,~$1_{\operatorname{R}}$, $0_{\operatorname{M}}$ are constant symbols, 
${-_{\operatorname{R}}}$, ${-_{\operatorname{M}}}$ are unary function symbols, and 
${+}_{\operatorname{R}}$, ${\,\cdot\,}_{\operatorname{R}}$, ${+}_{\operatorname{M}}$  are binary function symbols. Here the symbols indexed with a subscript $\operatorname{R}$ are of sort $\operatorname{R}$,
and similarly for $\operatorname{M}$. The symbol $\lambda$ is a function symbol of sort $\operatorname{R}\operatorname{M}\operatorname{M}$. %\nomenclature[Bex7]{$\mathcal L_{\operatorname{Mod}}$}{language of modules}
\end{enumerate}
\end{examplesNumbered}

\noindent
\textit{In the rest of this appendix $\mathcal L$ is a language with $S$ as its set of sorts, unless specified
otherwise. We let $s$ \textup{(}possibly subscripted\textup{)} range over $S$ and $\mathbf s=s_1\dots s_m$ over~$S^*$.}\/ 
The {\bf size} of $\mathcal L$ is the cardinal
$$\abs{\mathcal L}\ :=\  \max\!\big\{\aleph_0, \abs{S}, \abs{\mathcal L^{\operatorname{r}}\cup\mathcal L^{\operatorname{f}}}\big\},$$
and we say that $\mathcal L$ is {\bf countable} if $\abs{\mathcal L} = \aleph_0$, that is, $S$, $\mathcal  L^{\operatorname{r}}$, and $\mathcal L^{\operatorname{f}}$ are countable.  \index{language!size} \index{language!countable}

\subsection*{$\mathcal L$-structures}
An {\bf $\mathcal L$-structure} is an $S$-sorted structure
$$\mathbf M\ =\ \big( M; (R^{\mathbf M})_{R\in \mathcal L^{\operatorname{r}}}, (f^{\mathbf M})_{f\in \mathcal L^{\operatorname{f}}} \big)\qquad \text{where $M=(M_s)$,}$$
such that for every $R\in \mathcal L^{\operatorname{r}}$
of sort $s_1\dots s_m$ its {\bf interpretation $R^{\mathbf M}$} in $\mathbf M$ is a subset of~$M_{s_1\dots s_m}$, and for $f\in \mathcal L^{\operatorname{f}}$ of sort 
$s_1\dots s_n s$ its {\bf interpretation~$f^{\mathbf M}$} in $\mathbf M$ is a function $M_{s_1\dots s_n}\to M_{s}$. For a constant symbol $c$ of $\mathcal L$ of sort $s$ the corresponding $M_s$-valued function $c^{\mathbf M}$ is identified with its unique value in $M_s$, so $c^{\mathbf M}\in M_s$. If~$\mathbf M$ is understood from the context we often omit the superscript $\mathbf M$ in denoting the interpretation in $\mathbf M$ of a nonlogical symbol of $\mathcal L$.  The reader is supposed to keep in mind the distinction between symbols
of $\mathcal L$ and their interpretation in an $\mathcal L$-structure, even if we use the same notation
for both.   \index{structure!$\mathcal L$-structure} \index{structure@$\mathcal L$-structure} \index{interpretation} \index{structure!interpretation}

\begin{examplesNumbered}\label{ex:L-structures}
\mbox{}

\begin{enumerate}
\item Every group is considered as an $\mathcal L_{\operatorname{G}}$-structure by interpreting the symbols $1$, ${\ }^{-1}$, and $\cdot$ as the identity element of the group, its group inverse, and its group multiplication, respectively.
\item Let $\mathbf A=(A; 0,{-},{+})$ be an abelian group; here $0 \in A$ is the zero element
of the group, and ${-}\colon A \to A$ and ${+} \colon A^2 \to A$ denote the group operations
of $\mathbf A$. We consider $\mathbf A$ as an $\mathcal L_{\operatorname{A}}$-structure by taking as interpretations of
the symbols $0$, $-$ and $+$ of $\mathcal L_{\operatorname{A}}$ the group operations $0$, $-$ and $+$ on $A$.
(We took here the liberty of using the same notation for possibly entirely
different things: $+$ is an element of the set $\mathcal L_{\operatorname{A}}^{\operatorname{f}}$, but also denotes in this
context its interpretation as a binary operation on the set $A$. Similarly
with $0$ and $-$.) In fact, any set~$A$ for which we single out an element of $A$, a
unary operation on $A$, and a binary operation on $A$, is
an $\mathcal L_{\operatorname{A}}$-structure if we choose to construe it that way.
\item Likewise, any ring is construed as an $\mathcal L_{\operatorname{R}}$-structure in the obvious way.
\item
$(\N; {\leq})$ is an $\mathcal L_{\operatorname{O}}$-structure where we interpret $\leq$ as the usual ordering
relation on~$\N$. Similarly for $(\Z; {\leq})$, $(\Q; {\leq})$ and $(\R; {\leq})$. (Here we take
even more notational liberties, by letting $\leq$ denote five different things: a
symbol of $\mathcal L_{\operatorname{O}}$, and the usual orderings of $\N$, $\Z$, $\Q$, and $\R$, respectively.)
Again, any nonempty set  equipped with a binary relation on it can be
viewed as an $\mathcal L_{\operatorname{O}}$-structure.
\item Any ordered abelian group is an $\mathcal L_{\operatorname{OA}}$-structure.
\item The ordered ring of integers and any ordered field are $\mathcal L_{\operatorname{OR}}$-structures.
\item Let $R$ be a ring and $M$ be an $R$-module. Then $M$ becomes an $\mathcal L_{R\!\operatorname{-mod}}$-structure by interpreting the symbols of $\mathcal L_{\operatorname{A}}$ as in (2) and the function symbol~$\lambda_r$~($r\in R$) by the function $x\mapsto rx\colon M\to M$.
We can also construe~$M$ as an $\mathcal L_{\operatorname{Mod}}$-structure 
whose underlying set of sort~$\operatorname{R}$ is $R$ and
whose underlying set of sort~$\operatorname{M}$ is $M$, and where
the symbols $0_{\operatorname{R}}$,~$1_{\operatorname{R}}$,~${-_{\operatorname{R}}}$,~${+}_{\operatorname{R}}$,~${\,\cdot\,}_{\operatorname{R}}$ are interpreted by the ring operations on $R$,
the symbols $0_{\operatorname{M}}$,~${-}_{\operatorname{M}}$,~${+}_{\operatorname{M}}$ by the group
operations on $M$, and~$\lambda$ by $(r,x)\mapsto rx\colon R\times M\to M$.
\end{enumerate}
\end{examplesNumbered}

\noindent
Let
$\mathcal L'$ be an $S'$-sorted extension of $\mathcal L$ (so
$S'\supseteq S$). An
$\mathcal L'$-structure $\mathbf M' = \big(M'; \dots\big)$ is said to be an {\bf $\mathcal L'$-expansion} of an $\mathcal L$-structure
 $\mathbf M=(M;\dots)$ (and
$\mathbf M$ is said to be the {\bf $\mathcal L$-reduct} of $\mathbf M'$) if $M_s = M'_s$ for all $s$, and each
nonlogical symbol of $\mathcal L$ has the same interpretation in~$\mathbf M$ as in $\mathbf M'$. 
Given such an $\mathcal L'$-expansion $\mathbf M'$ of $\mathbf M$ we sometimes abuse notation by letting a parameter set $A = (A_s)$ in $\mathbf M$ denote also the parameter set~$(A_{s'})_{s'\in S'}$ in~$\mathbf M'$, where $A_{s'} = \emptyset$ 
for~$s' \in S' \setminus S$.  \index{structure!expansion}\index{structure!reduct} \index{expansion}\index{reduct}

\begin{example} The ordered field of reals is an $\mathcal L_{\operatorname{OR}}$-expansion of 
the real line, where ``the real line'' stands for the $\mathcal L_{\operatorname{O}}$-structure $(\R; {\leq})$.
\end{example}

\subsection*{Substructures} Let $\mathbf M=(M;\dots)$ and
$\mathbf N=(N;\dots)$ be $\mathcal L$-structures.
Then $\mathbf M$ is called a {\bf substructure} of $\mathbf N$  (notation: $\mathbf M\subseteq\mathbf N$) if $M_s\subseteq N_s$ for all $s$, and 
\begin{align*}
R^{\mathbf N}\cap M_{s_1\dots s_m}\ 	&=\  R^{\mathbf M}	&	& \text{whenever $R\in\mathcal L^{\operatorname{r}}$ has sort $s_1\dots s_m$},\\
f^{\mathbf N}|_{M_{s_1\dots s_n}}\ 	&=\  f^{\mathbf M}\colon M_{s_1\dots s_n}\to M_{s}	&	& \text{whenever $f\in\mathcal L^{\operatorname{f}}$ has sort $s_1\dots s_n s$.}
\end{align*}
We also say in this case that {$\mathbf N$ \bf is an extension of $\mathbf M$} or that {$\mathbf N$ \bf extends $\mathbf M$}. For groups, this is just the notion of subgroup, and for rings the notion of subring.\index{structure!substructure}\index{substructure}\index{extension!$\mathcal L$-structures}\index{structure!extension}

\nomenclature[Bj00]{$\mathbf M\subseteq\mathbf N$}{$\mathbf M$ is a substructure of $\mathbf N$}

\subsection*{Morphisms}
Let $\mathbf M$ and $\mathbf N$ be $\mathcal L$-structures. \index{morphism!structures} \index{structure!morphism}
A {\bf morphism} $h \colon \mathbf M \to \mathbf N$ is a map $h\colon M\to N$ that respects the primitives of $\mathbf M$,  that is:
\begin{list}{}{\leftmargin=2.5em}
\item[(M1)] for every  $R\in\mathcal L^{\operatorname{r}}$ of sort $s_1\dots s_m$ and all 
$a \in M_{s_1\dots s_m}$,
$$a\in R^{\mathbf M} \quad\Rightarrow\quad h_{s_1\dots s_m}a \in R^{\mathbf N};$$
\item[(M2)] for every $f\in\mathcal L^{\operatorname{f}}$ of sort $s_1\dots s_n s$
and all $a \in M_{s_1\dots s_n}$,
$$h_{s}\big(f^{\mathbf M}(a)\big)\ =\ f^{\mathbf N}(h_{s_1\dots s_n}a).$$
\end{list}
Replacing $\Rightarrow$ in (M1) by $\Longleftrightarrow$ 
and also requiring $h$ to be injective
yields the notion of an {\bf embedding}. 
An {\bf isomorphism} is a bijective
embedding, and an {\bf automorphism} of~$\mathbf M$ is an isomorphism $\mathbf M\to\mathbf M$. 
If $\mathbf M \subseteq \mathbf N$, then the inclusions $a\mapsto a \colon M_s \to N_s$ yield an embedding $\mathbf M \to \mathbf N$, called the {\bf natural inclusion} of  $\mathbf M$ into $\mathbf N$.
Conversely, a morphism $h \colon \mathbf M \to \mathbf N$ yields a substructure~$h(\mathbf M)$ of~$\mathbf N$ 
whose 
underlying set of sort~$s$ is  $h_s(M_s)$, and if $h$ is an embedding we have an isomorphism $\mathbf M \to h(\mathbf M)$ given by 
$$a\mapsto h_s(a)\colon M_s\to h(M_s).$$
If $h\colon  \mathbf M\to \mathbf M'$ and $h' \colon \mathbf M' \to\mathbf M''$ are morphisms (embeddings, isomorphisms, respectively), then so is $h'\circ h \colon \mathbf M \to\mathbf M''$.
The automorphisms of $\mathbf M$ form a group~$\Aut(\mathbf M)$ under composition. A parameter set~$A$ in $\mathbf M$
yields the subgroup 
$$\Aut(\mathbf M|A)\ :=\ \big\{f\in\Aut(\mathbf M): \text{$f_s(a) = a$ for every $s$ and $a \in A_s$} \big\}$$
of $\Aut(\mathbf M)$. We write $\mathbf M\cong\mathbf N$ if there is an isomorphism $\mathbf M\to\mathbf N$.
\index{structure!isomorphism}\index{isomorphism!$\mathcal L$-structures}\index{automorphism!$\mathcal L$-structure}\index{structure!automorphism}\index{embedding!$\mathcal L$-structures}\index{structure!embedding}
\index{structure!natural inclusion} \index{natural inclusion}

\nomenclature[Bk1]{$\Aut(\mathbf M)$}{group of automorphisms of $\mathbf M$}
\nomenclature[Bk2]{$\Aut(\mathbf M\lvert A)$}{group of automorphisms of $\mathbf M$ over $A$}
\nomenclature[Bj1]{$\mathbf M\cong\mathbf N$}{$\mathbf M$, $\mathbf N$ are isomorphic}

\begin{examples} 
If $\mathbf M$ and $\mathbf N$ are abelian groups, construed as $\mathcal L_{\operatorname{A}}$-structures according to~\ref{ex:L-structures}(2), then a
morphism $\mathbf M\to\mathbf N$ is exactly what in algebra is called a (homo)morphism
from the group $\mathbf M$ into the group $\mathbf N$. Likewise with rings, and other kinds of
algebraic structures. 
\end{examples}

\subsection*{Products} \index{structure!product}
Now let $(\mathbf N_\lambda)_{\lambda\in\Lambda}$ be a family of $\mathcal L$-structures,
where $\Lambda\neq\emptyset$, and
let $\lambda$ range over $\Lambda$. For each $s$ 
put $N_s:=\prod_{\lambda} (N_\lambda)_s$. (The Axiom of Choice guarantees that $N_s\neq\emptyset$.)
We write the typical element of~$N_s$ as $a=\big(a(\lambda)\big)$.
For  $a=(a_1,\dots,a_m)\in N_{\mathbf s}$,   put 
$$a(\lambda)\ :=\ \big(a_1(\lambda),\dots,a_m(\lambda)\big)\in (N_\lambda)_{\mathbf s}.$$
The {\bf product} $\prod_{\lambda} \mathbf N_\lambda$ of $(\mathbf N_\lambda)$ 
is defined to be the $\mathcal L$-structure $\mathbf N$ whose underlying set of sort $s$ is $N_s$, and where the basic relations and functions are defined coordinatewise:
for $R\in\mathcal L^{\operatorname{r}}$ of sort $s_1\dots s_m$ and $a\in N_{s_1\dots s_m}$,
$$a\in R^{\mathbf N} \quad:\Longleftrightarrow\quad \text{$a(\lambda)\in R^{\mathbf N_\lambda}$ for all $\lambda$},$$
and for $f\in\mathcal L^{\operatorname{f}}$ of sort $s_1\dots s_n s$
and $a\in N_{s_1\dots s_n}$,
$$f^{\mathbf N}(a)\  :=\  \big(f^{\mathbf N_\lambda}(a(\lambda)) \big) \in N_s.$$
For each $\lambda$ the projection to the $\lambda$th factor is the morphism $\pi_\lambda=\pi_\lambda^{\mathbf N}\colon\mathbf N\to\mathbf N_\lambda$
given by $(\pi_{\lambda})_s(a)=a(\lambda)$ for $a\in N_s$.
This product construction makes it possible to combine several morphisms
with a common domain into a single one: if for each~$\lambda$ we have a morphism
$h_\lambda \colon \mathbf M \to \mathbf N_\lambda$, then we obtain a morphism $h\colon \mathbf M\to \mathbf N$
with $\pi_\lambda\circ h = h_\lambda$ for each $\lambda$. If each $h_\lambda$ is an embedding, then so is~$h$.
This yields: \index{embedding!diagonal} \index{diagonal!embedding}

\begin{lemma}\label{lem:morphisms on products}
For each $\lambda$, let $h_\lambda\colon\mathbf M_\lambda\to\mathbf N_\lambda$ be a morphism. Then we have a morphism $h\colon \mathbf M:=\prod_\lambda \mathbf M_\lambda\to \mathbf N=\prod_\lambda \mathbf N_\lambda$ such that $\pi_\lambda^{\mathbf N}\circ h = h_\lambda\circ\pi_\lambda^{\mathbf M}$ for all~$\lambda$.
If every $h_\lambda$ is an embedding, then so is $h$.
\end{lemma}

\noindent
If $\mathbf N_\lambda=\mathbf M$ for all~$\lambda$, then $\prod_\lambda \mathbf N_\lambda$ is denoted by $\mathbf M^\Lambda$;
in this case we have an embedding $\Delta\colon\mathbf M\to \mathbf M^\Lambda$ with 
$\big(\Delta_s(a)\big)(\lambda) = a$ for all $a\in M_s$ and $\lambda$, called the 
{\bf diagonal embedding} of~$\mathbf M$ into $\mathbf M^\Lambda$.

\subsection*{Direct limits}
Let $\Lambda$ be a nonempty partially ordered set, and let $\lambda, \lambda',\lambda_1, \lambda_2,\dots$ range over~$\Lambda$. Suppose that $\Lambda$ is {\bf directed}, that is, for all $\lambda_1$,~$\lambda_2$ there exists~$\lambda$ with $\lambda_1, \lambda_2\leq \lambda$.
A {\bf directed system} of $\mathcal L$-structures indexed by $\Lambda$
consists of a family $(\mathbf M_\lambda)$
of $\mathcal L$-structures together with a family $(h_{\lambda\lambda'})_{\lambda\leq\lambda'}$ of
morphisms $h_{\lambda\lambda'}\colon \mathbf M_\lambda \to \mathbf M_{\lambda'}$ such that
 $h_{\lambda\lambda}=\operatorname{id}_{\mathbf M_{\lambda}}$ for all $\lambda$
and 
 $h_{\lambda\lambda''}=h_{\lambda'\lambda''}\circ h_{\lambda\lambda'}$ whenever $\lambda\leq\lambda'\leq\lambda''$.
Suppose that $\big((\mathbf M_\lambda),(h_{\lambda\lambda'})\big)$ is a directed system 
of $\mathcal L$-structures indexed by $\Lambda$. For this situation we have the following routine lemma. \index{directed system} \index{structure!directed system}

{\sloppy
\begin{lemma}
There
exists an $\mathcal L$-structure $\mathbf M$ and a family $(h_\lambda)$ of morphisms $h_\lambda\colon\mathbf M_\lambda\to\mathbf M$ with the following properties:
\begin{enumerate}
\item[\textup{(i)}] if $\lambda\leq\lambda'$, then $h_\lambda = h_{\lambda'} \circ h_{\lambda\lambda'}$;
\item[\textup{(ii)}] if $\mathbf N$ is an $\mathcal L$-structure and $(f_\lambda)$ is a family of morphisms
$f_\lambda\colon \mathbf M_\lambda\to\mathbf N$ such that $f_\lambda=f_{\lambda'}\circ h_{\lambda\lambda'}$ for all $\lambda\leq\lambda'$, then there is a unique morphism $g\colon \mathbf M\to\mathbf N$ such
that $f_\lambda=g\circ h_\lambda$ for all $\lambda$.
\end{enumerate}
\end{lemma}
}
   
\noindent   
We call an $\mathcal L$-structure $\mathbf M$ together with 
a family $(h_\lambda)$ of morphisms $h_\lambda\colon\mathbf M_\lambda\to\mathbf M$
as in the previous lemma a {\bf direct limit} of the directed system
$\big((\mathbf M_\lambda),(h_{\lambda\lambda'})\big)$.
If $\mathbf M$,~$(h_\lambda)$ and $\tilde{\mathbf M}$,~$(\tilde h_\lambda)$ are two direct limits of
$\big((\mathbf M_\lambda),(h_{\lambda\lambda'})\big)$,  then the unique morphism $h\colon \mathbf M\to \tilde{\mathbf M}$ such that
$\tilde h_\lambda = h\circ h_\lambda$ for all $\lambda$ is actually an isomorphism.
This allows us to speak of \textit{the}\/ direct limit of $\big((\mathbf M_\lambda),(h_{\lambda\lambda'})\big)$. 
One verifies easily that 
if all $h_{\lambda\lambda'}$ are embeddings, then so are all $h_\lambda$.
 \index{direct!limit} \index{structure!direct limit}

\medskip
\noindent
Now let $(\mathbf M_\lambda)$ be a family of $\mathcal L$-structures such that for $\lambda\leq\lambda'$,
$\mathbf M_\lambda$ is a substructure of~$\mathbf M_{\lambda'}$, with
natural inclusion $h_{\lambda\lambda'}\colon\mathbf M_{\lambda}\hookrightarrow \mathbf M_{\lambda'}$.
Then $\big((\mathbf M_\lambda),(h_{\lambda\lambda'})\big)$ is a  directed system of $\mathcal L$-structures, and its direct limit $\bigcup_{\lambda\in\Lambda} \mathbf M_\lambda$,~$(h_\lambda)$ is called  the {\bf direct union} of~$(\mathbf M_\lambda)$.
In the following we always identify each $\mathbf M_\lambda$ with a substructure of~$\bigcup_{\lambda\in\Lambda} \mathbf M_\lambda$
via the embedding $h_\lambda$, and we also simply speak of $\bigcup_{\lambda\in\Lambda} \mathbf M_\lambda$ as the direct union of~$(\mathbf M_\lambda)$. \index{direct!union} \index{structure!direct union}

\begin{exampleNumbered}\label{ex:Fpalg}
View fields as  $\mathcal L_{\operatorname{R}}$-structures in the natural way.
The algebraic closure~$\mathbb F_p^\alg$ of the finite field $\mathbb F_p$ with $p$ elements is the direct union
of the family of finite subfields of $\mathbb F_p^\alg$.
\end{exampleNumbered}

\subsection*{Notes and comments}
What we call here a language is also known as a \textit{signature,}\/ or a \textit{vocabulary.}\/ 
Uncountable languages appear in Mal$'$cev \cite{Malcev36},
and many-sorted languages and structures in Herbrand~\cite{Herbrand} and Schmidt~\cite{Schmidt38,Schmidt51}.  Structures in mathematical practice are usually one-sorted or two-sorted, but the case of infinitely many sorts does naturally arise in model theory. For example, even if the structure $\mathbf M$
is one-sorted, the structure~$\mathbf M^{\operatorname{eq}}$ associated to $\mathbf M$ by Shelah~\cite[III.6]{Shelah78} is most naturally viewed as infinitely-sorted; see \cite[Section~4.3]{Hodges}. 

\section{Variables and Terms}

\noindent
\textit{In the rest of this appendix $\mathbf M=(M;\dots)$ and $\mathbf N=(N;\dots)$ are $\mathcal L$-structures unless noted otherwise.}\/ In this section and the next we introduce the syntax (variables, terms, formulas) that helps to specify the definable sets of $\mathbf M$ along the lines of~\ref{sec:structures and definable sets},
uniformly for all $\mathcal L$-structures $\mathbf M$.  

\subsection*{Variables}
We assume that for each sort $s$ we have available  infinitely  many symbols, called {\bf variables of sort $s$}, chosen so that if $s$ and $s'$ are different sorts, then no variable of sort $s$ is a variable of sort $s'$. We also assume that no variable of any sort is a nonlogical symbol of any language.
A {\bf variable} of $\mathcal L$ is a variable of some sort~$s$, and a {\bf multivariable} of~$\mathcal L$ is a tuple $x=(x_i)_{i\in I}$ of \textit{distinct}\/ variables of $\mathcal L$.
The size of the index set $I$ is called the size of $x$, and the $x_i$ are called the variables in $x$.
Often the index set $I$ is finite, say $I=\{1,\dots,n\}$, so that $x=(x_1,\dots,x_n)$;
in this case, if $x_i$ is of sort $s_i$ ($i=1,\dots,n$), we say that $x$ is of sort $s_1\dots s_n$.
Instead of \textit{$x$ has finite size}\/ we also say \textit{$x$ is finite.}\/ Let $x=(x_i)_{i\in I}$ be a multivariable of~$\mathcal L$, with $x_i$ of sort $s_i$ for $i\in I$. We define the {\bf $x$-set} of~$\mathbf M$ as
$$M_x := \prod_{i\in I} M_{s_i},$$
and we think of $x$ as a variable running over $M_x$.  
When $I=\emptyset$, then
$x$ is said to be {\bf trivial}, and
$M_x$ is a 
singleton. If $h\colon \mathbf M\to\mathbf N$ is a morphism, we obtain a map 
$(a_i) \mapsto \big(h_{s_i}(a_i)\big)\colon M_x\to N_x$ between $x$-sets,
which we denote by $h_x$.
For a parameter set $A=(A_s)$ in $\mathbf M$,  we set
$$A_x := \prod_{i\in I} A_{s_i} \subseteq M_x.$$
Multivariables $x=(x_i)_{i\in I}$ and $y=(y_j)_{j\in J}$ of $\mathcal L$ are said to be {\bf disjoint} if $x_i\neq y_j$ for all $i\in I$, $j\in J$, and in that case we put $M_{x,y}:=M_x\times M_y$.
{\em From now on $x$ and~$y$ denote multivariables of $\mathcal L$, unless specified otherwise}.

\index{variable} \index{multivariable} \index{language!variable} \index{language!multivariable} \index{structure!$x$-set} \index{multivariable!disjoint} \index{multivariable!trivial}

\subsection*{Terms}
We define {\bf $\mathcal L$-terms} to be words on the alphabet consisting of the function
symbols and variables of $\mathcal L$, obtained recursively as follows: 
\begin{list}{}{\leftmargin=2.5em}
\item[(T1)] each variable of
sort~$s$, when viewed as a word of length~$1$, is an $\mathcal L$-term of sort~$s$, and
\item[(T2)]
if $f$ is a function symbol of $\mathcal L$ of sort $s_1\cdots s_n s$ and $t_1, \dots , t_n$ are $\mathcal L$-terms of sort $s_1, \dots, s_n$, respectively, then $ft_1 \cdots t_n$ is an $\mathcal L$-term of sort~$s$.
\end{list}
Thus every constant symbol of $\mathcal L$ of sort $s$ is an $\mathcal L$-term of sort $s$.
Usually we write $f(t_1, \dots , t_n)$ to denote $ft_1 \dots t_n$,   
and shun prefix notation if dictated by tradition. \index{language!term} \index{term}

\begin{example} Let $x$,~$y$,~$z$ be $\mathcal L_{\operatorname{R}}$-variables. Then the word ${\,\cdot\,}{+}x{-}yz$ is an $\mathcal L_{\operatorname{R}}$-term in the official prefix notation. For easier reading we indicate
this term instead by $(x + (-y)) \cdot z$ or even $(x - y)z$. 
\end{example}

\noindent
The following allows us to give definitions and proofs by \textit{induction on terms.}\/

{\sloppy
\begin{lemma}\label{lem:unique readability, terms}
Every $\mathcal L$-term of sort $s$ is either a variable of sort $s$, or equals $ft_1\dots t_n$ for a unique tuple
$(f, t_1,\dots, t_n)$ with $f$ an $n$-ary function symbol of some sort $s_1\dots s_ns$, and $t_i$ an $\mathcal L$-term of sort $s_i$ for $i=1,\dots,n$.
\end{lemma}
}

\noindent
For now we shall assume this lemma without proof. At the
end of this section we establish more general results which are also needed in the next section.

\medskip
\noindent
Let $x$ be a multivariable. An $(\mathcal L, x)$-term is a pair $(t, x)$ where $t$ is an $\mathcal L$-term such that each variable
in~$t$ is a variable in $x$.  Such an $(\mathcal L, x)$-term is also written more suggestively
as $t(x)$, and referred to as \textit{the $\mathcal L$-term $t(x)$.}\/ 
(It is not required that each variable in~$x$ actually occurs in $t$; this is like indicating a polynomial in the indeterminates
$X_1, \dots, X_n$ by $P(X_1, \dots , X_n)$, where some of these indeterminates
might not occur in $P$. Note that only finitely many of the variables in $x$ can occur in any $\mathcal L$-term.)
Given an $\mathcal L$-term $t(x)$ of sort $s$, we define a function $$t^{\mathbf M}\colon M_x\to M_s$$ as follows,
with $x=(x_i)$, and with $a=(a_i)$ ranging over $M_x$:
\begin{enumerate}
\item if $t=x_i$, then $t^{\mathbf M}(a):=a_i$;
\item if $t=ft_1\cdots t_n$ where
$f\in\mathcal L^{\operatorname{f}}$ is of sort $s_1\dots s_n s$ and $t_1, \dots , t_n$ are 
$\mathcal L$-terms of sort $s_1, \dots, s_n$, respectively,
then $$t^{\mathbf M}(a):=f^{\mathbf M}\big(t_1^{\mathbf M}(a),\dots,t_n^{\mathbf M}(a)\big)\in M_s.$$
\end{enumerate}
This inductive definition is justified by Lemma~\ref{lem:unique readability, terms}.

\begin{example}
Let $R$ be a commutative ring viewed as an $\mathcal L_{\operatorname{R}}$-structure in the natural way.
Let $x$,~$y$,~$z$ be distinct $\mathcal L_{\operatorname{R}}$-variables, and let $t(x, y, z)$ be the $\mathcal L_{\operatorname{R}}$-term 
$(x - y)z$. Then the function $t^{R}  \colon R^3 \to R$ is given by $t^R(a, b, c) = (a - b)c$.
In fact,
for each $\mathcal L_{\operatorname{R}}$-term $t(x_1,\dots,x_n)$
there is a (unique) polynomial $P^t(X_1,\dots,X_n)$ with integer coefficients such that for every commutative ring $R$ we have $t^R(a)=P^t(a)$ for all $a\in R^n$. Conversely, there is for each polynomial $P\in \Z[X_1,\dots, X_n]$ an 
$\mathcal L_{\operatorname{R}}$-term $t(x_1,\dots, x_n)$ such that $t^R(a)=P(a)$ for every commutative ring $R$ and all $a\in R^n$.
\end{example}

\noindent
Let $t(x)$ be an $\mathcal L$-term of sort $s$ and $a\in M_x$.
If $h\colon\mathbf M\to\mathbf N$ is a morphism of $\mathcal L$-structures, then
$h_s(t^{\mathbf M}(a))=t^{\mathbf N}(h_x(a))$; so  if
$\mathbf M\subseteq \mathbf N$, then $t^{\mathbf M}(a) = t^{\mathbf N}(a)$.
If~$\mathcal L'$ is an extension of~$\mathcal L$ and $\mathbf M'$ is an $\mathcal L'$-expansion of $\mathbf M$, then each $\mathcal L$-term $t(x)$ is also an $\mathcal L'$-term and
$t^{\mathbf M}=t^{\mathbf M'}\colon M_x \to M_s$.

\subsection*{Variable-free terms}
A term is said to be {\bf variable-free} if no variables occur in it. Let $t$ be a variable-free
$\mathcal L$-term of sort $s$. Then the above gives a nullary function~$t^{\mathbf M}$ with value in $M_s$, identified as usual with its value, so
$t^{\mathbf M} \in M_s$. In particular, if $t$ is a constant symbol $c$, then $t^{\mathbf M} = c^{\mathbf M} \in M_s$, where
$c^{\mathbf M}$ is as in Section~\ref{sec:languages}, and if $t = ft_1\dots t_n$ with $f\in\mathcal L^{\operatorname{f}}$ of sort $s_1\dots s_n s$ and
variable-free $\mathcal L$-terms $t_1,\dots,t_n$ of sorts $s_1,\dots,s_n$, respectively, then $t^{\mathbf M} = f^{\mathbf M}(t^{\mathbf M}_1, \dots, t^{\mathbf M}_n)\in M_s$.\index{term!variable-free}
%A family $(t_i)$ of $\mathcal L$-terms is said to be variable-free if each $t_i$ is variable-free.

\subsection*{Names} \index{language!name} \index{name}\index{symbols!names} \nomenclature[Bj0]{$\mathbf M_A$}{expansion of $\mathbf M$ by names for the elements of $A$}
Let $A\subseteq M$ be a parameter set.
We extend $\mathcal L$ to
a language $\mathcal L_A$ by adding a constant symbol $\underline{a}$ of sort $s$ for each $a \in A_s$, called the {\bf name}
of~$a$. These names are symbols not in $\mathcal L$. 
Note that $\abs{\mathcal L_A}=\max\!\big\{\abs{A},\abs{\mathcal L}\big\}$.
We make~$\mathbf M$ into an $\mathcal L_A$-structure by interpreting each name $\underline{a}$ as the element $a\in A_s$. The $\mathcal L_A$-structure thus
obtained is indicated by $\mathbf M_A$. (This is consistent with the notation introduced in Section~\ref{sec:structures and definable sets}.) Hence for each variable-free $\mathcal L_A$-term $t$ of sort $s$ we have
a corresponding element~$t^{\mathbf M_A}$ of~$M_s$, which for simplicity of notation we denote
instead by $t^{\mathbf M}$. All this applies  to the case $M = A$, where in $\mathcal L_M$ we
have a name~$\underline{a}$ for each element~$a$ of~$M_s$.
If $\mathbf M\subseteq\mathbf N$, then we consider $\mathcal L_M$ to be a sublanguage of~$\mathcal L_N$ in such a way
that each $a \in M$ has the same name in~$\mathcal L_M$ as in $\mathcal L_N$.
Then for each variable-free $\mathcal L_M$-term~$t$ we have $t^{\mathbf M} = t^{\mathbf N}$.
Given a multivariable $x=(x_i)$ of $\mathcal L$ and $a=(a_i)\in M_x$,  we set $\underline{a}:=(\underline{a_i})$. 

\subsection*{Substitution} \index{term!substitution} \index{substitution!into terms} 
Let $\alpha$ be an $\mathcal L_M$-term, $x=(x_i)_{i\in I}$ a multivariable, and
$t=(t_i)_{i\in I}$ a family of $\mathcal L_M$-terms
such that  $t_i$ is of the same sort as $x_i$, for each $i$. Then~$\alpha(t/x)$ denotes the word obtained by
replacing all occurrences of $x_i$ in $\alpha$ by $t_i$, \textit{simultaneously}\/ for all $i$.
If $\alpha$ is given in the form $\alpha(x)$, then $\alpha(t)$ is short
for $\alpha(t/x)$ and for $a\in M_x$ we often write~$\alpha(a)$ instead of $\alpha(\underline{a})$.

\begin{lemma}\label{lem:subst lemma, terms}
Let $\alpha(x)$ be an $\mathcal L_M$-term of sort $s$,
and recall that $\alpha$ defines a
map $\alpha^{\mathbf M}\colon M_{x} \to M_s$. Let $t=(t_i)_{i\in I}$ be a family of  $\mathcal L_M$-terms, with each $t_i$
of the same sort as $x_i$. Then $\alpha(t)$ is an $\mathcal L_M$-term of sort $s$.
Moreover, if all $t_i$ are variable-free, then so is~$\alpha(t)$, and 
with $a_i:= t_i^{\mathbf M}$, $a:=(a_i)\in M_x$ %and $\underline{a}:=(\underline{a}_i)$, 
we have
$$\alpha(t)^{\mathbf M}\ =\ \alpha(a)^{\mathbf M}\ =\  \alpha^{\mathbf M}(a).$$
\end{lemma}

\noindent
This follows by a straightforward induction on the length of $\alpha$.

\subsection*{Generators} \index{substructure!generators} \index{structure!generators} \nomenclature[Bj000]{$\langle A\rangle_{\mathbf M}$}{substructure of $\mathbf M$ generated by $A\subseteq M$}
Assume $\mathcal L$ has for each~$s$ a constant symbol of sort $s$. Let $A\subseteq M$ be given.
For each~$s$ we set
$$B_s\ :=\ \big\{t^{\mathbf M}(a):\ \text{$t(x)$ is an $\mathcal L$-term of sort~$s$, and~$a\in A_x$}\big\}\ \subseteq\ M_s .$$
Then $B=(B_s)$ underlies a substructure $\mathbf B=(B;\dots)$
of $\mathbf M$, and clearly $\mathbf B\subseteq \mathbf C$ for all substructures
$\mathbf C$ of $\mathbf M$ with $A\subseteq C$. We call~$\mathbf B$ the
substructure of~$\mathbf M$ {\bf generated} by~$A$; notation: $\mathbf B=\<A\>_{\mathbf M}$. 
Note that $\abs{A}\leq \abs{\<A\>_{\mathbf M}}\leq\max\!\big\{\abs{A},\abs{\mathcal L}\big\}$.
If $\mathbf M = \<A\>_{\mathbf M}$, then we say that $\mathbf M$ is generated
by $A$.
If $\mathbf N$ is an $\mathcal L$-structure, then each map~${A\to N}$ has clearly at most one extension to a morphism $\<A\>_{\mathbf M}\to\mathbf N$.

\subsection*{Unique readability}
We finish this section with the promised general result on unique readability. 
We let~$F$ be a set of symbols with a function $a \colon F \to \N$ (called the
{\bf arity} function). A symbol~$f \in F$ is said to have arity $n$ if $a(f) = n$. 
A word on~$F$ is said to be {\bf admissible} if it can be obtained by applying the following rules:
\begin{enumerate}
\item If $f \in F$ has arity $0$, then $f$ viewed as a word of length $1$ is admissible.
\item If $f \in F$ has arity $n\geq 1$ and $t_1, \dots, t_n$ are admissible words on $F$, then
the concatenation $ft_1\dots t_n$ is admissible.
\end{enumerate}
Below we just write \textit{admissible word}\/ instead of \textit{admissible word on $F$.}\/ Note
that the empty word is not admissible, and that the last symbol of an admissible
word cannot be of arity~$\geq 1$. \index{admissible!word}

\begin{example}
The $\mathcal L$-terms are admissible words on the alphabet consisting of the function symbols of $\mathcal L$ and the variables of $\mathcal L$, where each $n$-ary $f\in \mathcal L^{\operatorname{f}}$ has arity $n$ and each variable has arity $0$.  
\end{example}

\begin{lemma}
Let $t_1, \dots , t_m$ and $u_1, \dots, u_n$ be admissible words and $w$ any
word on $F$ such that $t_1 \dots t_mw = u_1 \dots u_n$. Then $m \leq n$, $t_i = u_i$ for $i =
1, \dots ,m$, and $w = u_{m+1} \dots u_n$.
\end{lemma}
\begin{proof}
By induction on the length $\ell$ of $u_1 \dots u_n$. If $\ell=0$, then $m = n = 0$ and~$w$ is the empty word. Suppose $\ell> 0$, and assume the lemma
holds for smaller lengths. Note that $n > 0$. If $m = 0$, then the conclusion of the
lemma holds, so suppose~$m > 0$. The first symbol of $t_1$ equals the first symbol
of $u_1$. Say this first symbol is $h \in F$ with arity $k$. Then $t_1 = ha_1\dots a_k$ and
$u_1 = hb_1  \dots b_k$ where $a_1, \dots, a_k$ and $b_1, \dots, b_k$ are admissible words. Canceling
the first symbol $h$ gives
$$a_1 \dots a_k t_2 \dots t_m w = b_1 \dots b_k u_2 \dots u_n.$$
(Caution: any of $k$, $m-1$, $n-1$ could be $0$.) We have $\operatorname{length}(b_1\dots b_ku_2 \dots u_n) = \ell-1$, so the induction hypothesis applies. It yields $k + m  -1 \leq k + n - 1$ (so $m \leq n$), $a_1 = b_1, \dots , a_k = b_k$ (so $t_1 = u_1$), $t_2 = u_2$, \dots , $t_m = u_m$,
and $w = u_{m+1} \cdots u_n$.
\end{proof}

\noindent
In particular, if $t_1, \dots, t_m$ and $u_1, \dots , u_n$ are admissible words such that $t_1 \dots t_m =
u_1 \dots u_n$, then $m = n$ and $t_i = u_i$ for $i = 1,\dots ,m$. Thus:

\begin{cor}[unique readability]
Each admissible word equals $ft_1 \dots t_n$ for a unique tuple $(f,t_1,\dots,t_n)$ with $f \in F$ of arity
$n$ and $t_1, \dots , t_n$ admissible words.
\end{cor}

\section{Formulas}\label{sec:formulas}

\noindent
We now fix once and for all the eight {\bf logical symbols}:
$$\top \qquad \bot \qquad \neg  \qquad \vee \qquad \wedge \qquad = \qquad \exists \qquad \forall$$
to be thought of as \textit{true,}\/ \textit{false,}\/ \textit{not,}\/ \textit{or,}\/ \textit{and,}\/ \textit{equals,}\/ \textit{there exists,}\/ and \textit{for all,}\/ respectively.
These are assumed to be distinct from all nonlogical symbols
and variables of every language. 
The symbols $\neg$, $\vee$, $\wedge$ are called (logical) {\bf connectives}
and $\exists$,~$\forall$ are called {\bf quantifiers}.
It will be convenient to fix
once and for all a sequence $v_0, v_1, v_2,\dots$ of distinct
symbols, called {\em unsorted\/} variables, and to declare that for each sort 
$s$ the symbols $v_n^s$ are the {\bf quantifiable variables}\/ of sort $s$. For any $S$-sorted language $\mathcal{L}$ and~$s\in S$, these are
among its variables of sort $s$, but $\mathcal{L}$ can have other (unquantifiable) variables of sort $s$. 
A multivariable is called quantifiable if each variable in it is quantifiable. \index{symbols!logical}\index{logical!symbols}\index{symbols!connectives}\index{symbols!quantifiers} \index{variable!quantifiable} \index{multivariable!quantifiable}

\subsection*{Formulas} \index{formula!atomic} \index{atomic formula}
The {\bf atomic $\mathcal L$-formulas} are the following words  on the alphabet $$\mathcal{L}^{\operatorname{r}}\cup\mathcal{L}^{\operatorname{f}}\cup\{\text{variables of $\mathcal L$}\}\cup\{\top, \bot,{=}\}:$$  
\begin{list}{}{\leftmargin=2.5em }
\item[(A1)] $\top$ and $\bot$ (as words of length $1$);
\item[(A2)] the words $Rt_1\dots t_m$ where  $R\in\mathcal L^{\operatorname{r}}$ is of
sort $s_1\dots s_m$ and $t_1,\dots,t_m$ are $\mathcal L$-terms of sort $s_1,\dots,s_m$, respectively; 
\item[(A3)] the words ${=}t_1t_2$ where $t_1$, $t_2$ are $\mathcal L$-terms of the same sort. 
\end{list}
The {\bf $\mathcal L$-formulas} are the words on the alphabet 
$$ \mathcal{L}^{\operatorname{r}}\cup\mathcal{L}^{\operatorname{f}}\cup\{\text{variables of $\mathcal L$}\}\cup\{\top, \bot,{\neg}, {\vee},{\wedge},{=}, {\exists}, {\forall}\}$$ obtained as follows:
\begin{list}{}{\leftmargin=2.5em }
\item[(F1)] atomic $\mathcal L$-formulas are $\mathcal L$-formulas;
\item[(F2)] if $\varphi$, $\psi$ are $\mathcal L$-formulas, then so are $\neg\varphi$, $\vee\varphi\psi$, and
$\wedge\varphi\psi$;
\item[(F3)] if $\varphi$ is an $\mathcal L$-formula and $x$ is a quantifiable variable, then $\exists x\varphi$ and
$\forall x\varphi$ are $\mathcal L$-formulas.
\end{list}

\index{formula}
 
\begin{samepage} 
\begin{remarks}\mbox{}

\begin{enumerate}
\item Having the connectives $\vee$ and $\wedge$ in front of the $\mathcal L$-formulas they
``connect'' rather than in between, is called \textit{prefix notation}\/ or \textit{Polish notation.}\/
This is theoretically elegant, but for the sake of readability we usually write $\varphi\vee\psi$
and $\varphi\wedge\psi$ to denote $\vee \varphi\psi$ and $\wedge \varphi\psi$, respectively, and we also use parentheses and
brackets if this helps to clarify the structure of an $\mathcal L$-formula.  
\item All $\mathcal L$-formulas are admissible words on the alphabet consisting of the nonlogical
symbols of $\mathcal L$,  the variables of $\mathcal L$, and the eight logical symbols,
where~$\top$ and~$\bot$ have arity $0$, 
$\neg$ has arity $1$, 
$\vee$, $\wedge$, $=$, $\exists$ and $\forall$ have arity 2, and the other symbols have the arities
assigned to them earlier. Thus the results on unique readability
are applicable to $\mathcal L$-formulas. (However, not all admissible words on this alphabet
are $\mathcal L$-formulas: the word $\exists xx$ is admissible but not an $\mathcal L$-formula.)
\item  The reader should distinguish between different ways of using the symbol~$=$.
Sometimes it denotes one of the eight formal logical symbols, but we also use it
to indicate equality of mathematical objects in the usual way. The context should always make it clear what our intention is in
this respect without having to spell it out. To increase readability we usually
write an atomic formula ${=} t_1t_2$ as $t_1 = t_2$ and its negation $\neg{=}t_1t_2$ as $t_1 \neq t_2$,
where $t_1$, $t_2$ are $\mathcal L$-terms of the same sort. 
\end{enumerate}
\end{remarks}
\end{samepage}

\noindent
We shall use the following notational conventions: $\varphi \to \psi$ denotes $\neg\varphi \vee \psi$, and
$\varphi \leftrightarrow \psi$ denotes $(\varphi \rightarrow \psi) \wedge (\psi \rightarrow \varphi)$.   We sometimes write $\varphi\ \&\ \psi$ instead of~$\varphi \wedge \psi$.

\begin{definition}
Let $\varphi$ be an $\mathcal L$-formula. Written as a word on the alphabet introduced above
we have $\varphi = a_1 \dots a_m$. 
A {\bf subformula} of $\varphi$ is a subword of the 
form $a_i\dots a_k$ (where $1 \leq i \leq k \leq m$) which also happens to be an $\mathcal L$-formula.
An occurrence of a variable $x$ in $\varphi$ at the $j$th place (that is, $a_j = x$) is said
to be {\bf bound} if~$\varphi$ has a subformula $a_ia_{i+1} \dots a_k$ with $i \leq j \leq k$ that
is of the form $\exists x\psi$ or $\forall x\psi$. An occurrence which is not bound is said to be {\bf free}. A {\bf free variable} of~$\varphi$ is a variable that occurs free in $\varphi$.  \index{variable!free} \index{formula!subformula} \index{subformula} \index{variable!free occurrence} \index{variable!bound occurrence}
\end{definition}

\begin{example}
In the $\mathcal L_{\operatorname{A}}$-formula 
$\big( \exists x(x = y) \big) \wedge x = 0$, where $x$ and $y$ are distinct variables, the first two occurrences of $x$ are bound, the third is free, and the only occurrence of $y$ is free.
(Note: this formula is actually the string $\wedge \exists x = xy = x0$, and the occurrences of $x$ and $y$ are the occurrences in this string.)
\end{example}

\noindent
An $(\mathcal L,x)$-formula is a pair $(\varphi,x)$ with $\varphi$ an $\mathcal L$-formula and $x$ a multivariable such that
all free variables of $\varphi$ are in $x$. Such an $(\mathcal L,x)$-formula is written more suggestively as~$\varphi(x)$, and referred to as \textit{the $\mathcal L$-formula $\varphi(x)$.}\/ Likewise, when referring to an $\mathcal L$-formula~$\varphi(x,y)$ we really mean a triple $(\varphi,x,y)$ consisting of an $\mathcal L$-formula~$\varphi$ and disjoint multivariables $x$ and $y$ such that all free variables of $\varphi$ are in $x$ or $y$.
If $\mathcal L$ and $x$ both have size $\le \kappa$,
then the set of $(\mathcal L, x)$-formulas has size $\le \kappa$. 

\medskip
\noindent
An {\bf $\mathcal L$-sentence} is an $\mathcal L$-formula without free 
variables. 
The set of all $\mathcal L$-sentences has size $|\mathcal L|$: that is why in 
(F3) we let $x$ be quantifiable. \index{sentence} \index{formula!sentence}

\subsection*{Substitution} Let $\varphi$ be an $\mathcal L$-formula, let $x=(x_i)_{i\in I}$ be a multivariable, and
let~$t=(t_i)_{i\in I}$ be a family of $\mathcal L$-terms with $t_i$ of the same sort as $x_i$ for all $i$. Then~$\varphi(t/x)$ denotes the word obtained by
replacing all the free occurrences of $x_i$ in $\varphi$ by~$t_i$, \textit{simultaneously}\/ for all $i$.
If $\varphi=\varphi(x)$, then we also write $\varphi(t)$ instead of~$\varphi(t/x)$.
We have the following facts whose routine proofs are left to the reader.

\begin{lemma}
The word $\varphi(t/x)$ is an $\mathcal L$-formula. If
 $t$ is variable-free and $\varphi = \varphi(x)$, then $\varphi(t)$ is an $\mathcal L$-sentence.
\end{lemma}

\noindent
Given an $\mathcal L$-structure $\mathbf M$, an $\mathcal L_M$-formula $\varphi(x)$ and $a\in M_x$, we shall avoid many ugly 
expressions by
writing the $\mathcal L_M$-sentence $\varphi(\underline{a})$ as
just $\varphi(a)$.  \index{substitution!into formulas} \index{formula!substitution}

\subsection*{Truth and definability}
We can now define what it means for an $\mathcal L_M$-sentence~$\sigma$ to be {\bf true} \index{sentence!truth}
in the $\mathcal L$-structure~$\mathbf M$ (notation: $\mathbf M\models \sigma$, also read as   $\mathbf M$ {\bf satisfies} $\sigma$, or   $\sigma$ {\bf holds} in $\mathbf M$). \nomenclature[Bl0]{$\mathbf M\models\sigma$}{$\sigma$ is true in $\mathbf M$} First we consider atomic $\mathcal L_M$-sentences:
\begin{list}{}{\leftmargin=2.5em}
\item[(T1)] $\mathbf M \models \top$ and $\mathbf M \not\models \bot$;
\item[(T2)]  $\mathbf M \models Rt_1 \dots t_m$ if and only if $(t_1^{\mathbf M} , \dots , t_m^{\mathbf M})
\in R^{\mathbf M}$, for $R \in \mathcal L^{\operatorname{r}}$ of sort $s_1\dots s_m$ and
variable-free $\mathcal L_M$-terms $t_1, \dots , t_m$ of sort $s_1,\dots,s_m$, respectively;
\item[(T3)]  $\mathbf M \models {t_1 = t_2}$ if and only if $t^{\mathbf M}_1 = t^{\mathbf M}_2$, for variable-free $\mathcal L_M$-terms $t_1$, $t_2$ of the same sort.
\end{list}
We extend the definition inductively to arbitrary $\mathcal L_M$-sentences as follows:
\begin{list}{}{\leftmargin=2.5em}
\item[(T4)] Suppose $\sigma = \neg\tau$; then $\mathbf M\models \sigma$ if and only if $\mathbf M\not\models\tau$.
\item[(T5)] Suppose $\sigma = \sigma_1 \vee \sigma_2$; then $\mathbf M\models \sigma$ if and only if $\mathbf M\models \sigma_1$ or $\mathbf M\models \sigma_2$.
\item[(T6)] Suppose $\sigma = \sigma_1 \wedge \sigma_2$; then $\mathbf M\models \sigma$ if and only if $\mathbf M\models \sigma_1$ and $\mathbf M\models \sigma_2$.
\item[(T7)] Suppose $\sigma = \exists x\varphi$ where $x$ is a quantifiable variable and $\varphi(x)$ is an $\mathcal L_M$-formula; then $\mathbf M\models \sigma$ if and only if $\mathbf M\models \varphi(a)$ for some $a \in M_x$.
\item[(T8)] Suppose $\sigma = \forall x\varphi$  where $x$ is a quantifiable variable and $\varphi(x)$ is an $\mathcal L_M$-formula; then $\mathbf M\models \sigma$ if and only if $\mathbf M\models \varphi(a)$ for all $a \in M_x$.
\end{list}
Even if we just want to define $\mathbf M\models \sigma$  for $\mathcal L$-sentences $\sigma$, one can see that if
$\sigma$ has the form considered in (T7) or (T8), the inductive definition above forces us to
consider $\mathcal L_M$-sentences $\varphi(a)$. This is why we introduced names. (``Inductive'' refers here to induction with respect to the number of
logical symbols in $\sigma$.)

\nomenclature[Bl9]{$\varphi^{\mathbf M}$}{set defined by $\varphi$ in $\mathbf M$}

\begin{definition} 
Given an $\mathcal L_M$-formula $\varphi(x)$  we let
$$\varphi^{\mathbf M}\ :=\ \big\{ a \in M_x:\  \mathbf M \models  \varphi(a)\big\}\subseteq M_x.$$
The formula $\varphi(x)$ is said to {\bf define} the set $\varphi^{\mathbf M}$ in $\mathbf M$. 
Given a parameter set~$A\subseteq M$, a subset 
of $M_{x}$
is said to be {\bf $A$-definable} in $\mathbf M$ if it is of the form $\varphi^{\mathbf M}$ for some $\mathcal L_A$-formula $\varphi(x)$. 
A map $f\colon X\to M_y$, where $X\subseteq M_x$,  is said to be $A$-definable if its graph $\Gamma(f)\subseteq M_{x,y}$ is.
We use {\bf definable} synonymously with  $M$-de\-finable. \index{set!$A$-definable} \index{set!definable} \index{definable!set}\index{formula!defining a set}\index{map!$A$-definable}\index{definable!map}
\end{definition}

\begin{examples}\mbox{}
Let $\mathbf R=(\R; {\leq}, 0, 1, {+}, {-}, {\cdot})$.
\begin{enumerate}
\item The set $\big\{r \in \R : r \leq \sqrt{2}\big\}$ is $0$-definable in $\mathbf R$: it is defined
by the formula $(x^2 \leq 1+1) \vee (x \le 0)$. (Here $x^2$ abbreviates the term $x\cdot x$.)
\item The set $\{r \in \R : r \leq \pi\}$ is definable in $\mathbf R$: it is defined by
the formula $x \leq \pi$. (It takes more effort to show that it is not $0$-definable; see Example~\ref{ex:RCF QE}.)
\end{enumerate}
\end{examples}

\noindent
It is easy to see that given $A\subseteq M$, a subset of $M_{\mathbf s}$ is $A$-definable in the sense of the previous definition 
if and only if it is $A$-definable in the sense of 
Definition~\ref{def:0-definable} (using the correspondences (C1)--(C5) between formulas and sets defined by them
from Section~\ref{sec:structures and definable sets}).
If $\psi(y)$ is an $\mathcal L_A$-formula, then there is a finite multivariable~$x$ disjoint from $y$, an $\mathcal L$-formula~$\varphi(x,y)$ and an~$a\in A_x$, such that $\varphi(a,y)=\psi(y)$.
Thus a set $Y \subseteq M_y$ is $A$-definable iff for some 
finite $x$ disjoint from $y$, some
$a \in A_x$ and some $0$-definable $Z \subseteq M_{x,y}$ we have $Y = Z(a):=\big\{b\in M_y: (a,b)\in Z\big\}$. In this way {\em $A$-definability} reduces to {\em $0$-definability}.

\medskip
\noindent
We defined what it means for an $\mathcal L$-sentence to hold in~$\mathbf M$. It is convenient to extend this to arbitrary $\mathcal L$-formulas.
First, given $A\subseteq M$ we define an {\bf $A$-instance} of an $\mathcal L_A$-formula $\varphi=\varphi(x)$  to be an 
$\mathcal L_A$-sentence of the form $\varphi(a)$ with $a\in A_x$. Of course~$\varphi$ can also
be written as $\varphi(y)$ for another multivariable $y$. Thus for the above to
count as a definition of \textit{$A$-instance}\/, the reader should check that these different
ways of specifying variables (including at least the free variables of $\varphi$) give the same $A$-instances.  \index{formula!$A$-instance}

\begin{definition} An $\mathcal L$-formula $\varphi$ is said to be {\bf valid} in $\mathbf M$ (notation: $\mathbf M\models\varphi$) if all its
$M$-instances are true in $\mathbf M$.  \index{formula!valid in a structure} 
\end{definition}

\noindent
Suppose the multivariable $x=(x_1,\dots,x_m)$ is (finite and) quantifiable, and $\varphi$ is an $\mathcal L$-formula. Then 
$\forall x\varphi$ denotes $\forall x_1\cdots\forall x_m \varphi$, and likewise with $\forall$ replaced by~$\exists$.
The reader should check that if $\varphi=\varphi(x)$, then~$\forall x\varphi$ is an $\mathcal L$-sentence and
$$\mathbf M \models \varphi\quad\Longleftrightarrow\quad \mathbf M \models \forall x\,\varphi.$$
We define $\models\varphi$ to mean: $\mathbf M\models\varphi$ for all $\mathcal L$-structures $\mathbf M$. We call $\mathcal L$-formulas~$\psi$,~$\theta$ {\bf equivalent} if $\models \psi \leftrightarrow \theta$. \index{formula!equivalence} \index{equivalence!formulas} Thus $\neg\forall x\varphi$ and $\exists x\neg\varphi$ are equivalent.

\medskip\noindent
Let $\varphi_1,\dots,\varphi_n$ be $\mathcal L$-formulas, where $n\geq 1$.
We inductively define 
$$\varphi_1\wedge\cdots\wedge\varphi_n := \begin{cases}
\varphi_1 													& \text{if $n=1$,} \\
\varphi_1\wedge\varphi_2									& \text{if $n=2$,} \\
(\varphi_1\wedge\cdots\wedge\varphi_{n-1})\wedge\varphi_n	& \text{if $n\geq 3$.}
\end{cases}$$
Similarly define $\varphi_1\vee\cdots\vee\varphi_n$.
For each permutation $i$ of $\{1,\dots,n\}$, the $\mathcal L$-formulas
$\varphi_1\wedge\cdots\wedge\varphi_n$ and $\varphi_{i1}\wedge\cdots\wedge\varphi_{in}$
are equivalent; similarly with $\wedge$ replaced by $\vee$.

\subsection*{Formulas of a special form} {\em In this subsection {\rm formula} means {\rm $\mathcal L$-formula}}.
We single out formulas by syntactical conditions with semantic counterparts in terms of behavior under
embeddings, as explained in the next subsection.
A formula is said to be {\bf quantifier-free} if it has no occurrences of $\exists$ and
no occurrences of~$\forall$. A formula is said to be {\bf existential} (or an {\bf $\exists$-formula}) if it has the form~$\exists x\,\rho$ with a finite multivariable $x$ and a quantifier-free formula $\rho$, and  {\bf universal} (or a {\bf $\forall$-formula}) if it has the form $\forall x\,  \rho$ with 
finite  $x$ and quantifier-free~$\rho$.
If~$\varphi$ is a $\forall$-formula, then $\neg\varphi$ is equivalent to an $\exists$-formula, and similarly with ``$\forall$'' and ``$\exists$'' interchanged.
A formula is said to be   {\bf universal-existential} (or a {\bf $\forall\exists$-formula}) if
it has the form $\forall x \exists y \,\rho$ with finite disjoint $x$, $y$ and quantifier-free~$\rho$.\index{formula!quantifier-free}\index{formula!existential}\index{formula!universal}\index{formula!universal-existential}\index{existential formula} \index{universal!formula} \index{formula!$\exists$-formula} \index{formula!$\forall$-formula}\index{formula!$\forall\exists$-formula}
 
\begin{lemma}
If $\varphi$, $\psi$ are $\exists$-formulas, then
$\varphi\wedge\psi$ and $\varphi\vee\psi$ are equivalent to $\exists$-formulas;
likewise with $\forall$ and with $\forall\exists$ in place of $\exists$. 
\end{lemma}
\begin{proof}
The easy proofs of these facts use the device of ``renaming variables'' which we often use tacitly below.
Let $\varphi=\exists x\,\rho$, $\psi=\exists y\,\theta$ where  $\rho$, $\theta$ are quantifier-free; 
it is easy to check that if $x$, $y$ are disjoint, then $\varphi\wedge\psi$ and
$\exists x\exists y(\rho\wedge\theta)$ are
 equivalent. 
In the general case, first
choose
disjoint finite quantifiable multivariables~$x'$,~$y'$ of the same sort as $x$, $y$, respectively,
such that no variable occurring in~$\rho$,~$\theta$ is in $x'$ or $y'$, 
and replace $\varphi$, $\psi$ by the respective equivalent formulas
 $\exists x'\,\rho(x'/x)$, $\exists y'\,\theta(y'/y)$.
\end{proof}

\subsection*{Maps preserving formulas}
\textit{In this subsection $A$ is a parameter set in $\mathbf M$ and 
$h\colon A\to N$ is a map.}\/
For an $\mathcal L_A$-term~$t$, let $h(t)$ be the $\mathcal L_N$-term obtained from~$t$ by replacing every occurrence of a name of an element $a\in A_s$ by the name of $h_s a\in N_s$. For an $\mathcal L_A$-formula~$\varphi$, let~$h(\varphi)$ be the $\mathcal L_N$-formula obtained from $\varphi$  by replacing every occurrence of a name of an element~$a\in A_s$ by the name of $h_sa$. If $t$ is variable-free, then so is~$h(t)$, and if $\varphi$ is a sentence, then so is $h(\varphi)$.
If $\varphi(x)$ is an $\mathcal L$-formula and~$a\in A_x$, then $h\big(\varphi(a)\big)=\varphi(ha)$.
We say that {\bf $h$ preserves atomic formulas} if for all 
atomic $\mathcal{L}$-formulas $\varphi(x)$ and $a\in A_x$ such that
$\mathbf{M}\models \varphi(a)$ we have $\mathbf{N}\models \varphi(ha)$. In the same way we define \textit{$h$~preserves quantifier-free formulas}\/ and \textit{$h$ preserves formulas.}\/\index{map!preserving formulas}

\begin{lemma}\label{lem:persist, 2}
Suppose $A=M$. Then we have the following two equivalences:
\begin{enumerate}
\item[\textup{(i)}] $h$ is a morphism $\mathbf M\to\mathbf N$ iff $h$ preserves atomic formulas;
\item[\textup{(ii)}] $h$ is an embedding $\mathbf M\to\mathbf N$ iff $h$ preserves
quantifier-free formulas.
\end{enumerate}
\end{lemma}

\noindent
We leave the proofs of this routine lemma to the reader.

\begin{cor}\label{cor:persist, 1}
Suppose $\mathbf M\subseteq\mathbf N$, and let $\sigma$ be an $\mathcal L_M$-sentence.
\begin{enumerate}
\item[\textup{(i)}] If $\sigma$ is quantifier-free, then we have: $\mathbf M\models\sigma\Longleftrightarrow\mathbf N\models \sigma$.
\item[\textup{(ii)}] If $\sigma$ is existential, then we have: $\mathbf M\models\sigma \Rightarrow \mathbf N\models \sigma$.
\item[\textup{(iii)}] If $\sigma$ is universal, then we have: $\mathbf N\models\sigma \Rightarrow \mathbf M\models \sigma$.
\end{enumerate}
\end{cor}

{\sloppy
\begin{proof}
Part (i) is immediate from Lemma~\ref{lem:persist, 2}(ii).
Part (ii) follows easily from~(i), and~(iii) follows from (ii).
\end{proof}
}

\noindent
In the next corollary we assume that 
$\mathcal L$ has for every $s$ a constant symbol of sort $s$. This is to guarantee that $\<A\>_{\mathbf M}$ is defined.

\begin{cor}\label{cor:persist, 2}
The map
$h$ extends to
a morphism $\<A\>_{\mathbf M}\to\mathbf N$ iff $h$ preserves
atomic formulas. The map
$h$ extends to
an embedding $\<A\>_{\mathbf M}\to\mathbf N$ iff $h$ preserves 
quantifier-free formulas.
\end{cor}
\begin{proof}
The forward directions in both statements follow from Lem\-ma~\ref{lem:persist, 2} and
Corollary~\ref{cor:persist, 1}(i). Suppose for every
atomic $\mathcal L_A$-sentence $\sigma$ true in $\mathbf M$, its image~$h(\sigma)$ is true in $\mathbf N$. Now every element of sort $s$ of $\<A\>_{\mathbf M}$ has the form $t^{\mathbf M}(a)$ for
some $\mathcal L$-term $t(x)$ of sort $s$ and some $a\in A_x$. Moreover, if $t_1(x)$ and $t_2(y)$ are $\mathcal L$-terms of sort $s$ and $t_1(a)=t_2(b)$, where $a\in A_x$ and $b\in A_y$,
then $\mathbf M\models t_1(a) = t_2(b)$, and hence
$\mathbf N\models t_1\big(h_x(a)\big) = t_2\big(h_y(b)\big)$, and thus
$t_1^{\mathbf N}\big(h_x(a)\big) = t_2^{\mathbf N}\big(h_y(b)\big)$. These two facts easily yield the backward directions of the two equivalences.
\end{proof}

\index{map!elementary}\index{elementary!map}\index{embedding!elementary}\index{elementary!embedding}\index{structure!elementary map}

\noindent
We say that $h\colon A \to N$ is {\bf elementary} if it preserves formulas, that is, for all $\mathcal L$-formulas $\varphi(x)$ and all $a\in A_x$, 
$$\mathbf M\models\varphi(a)\quad\Longleftrightarrow\quad\mathbf N\models\varphi(ha).$$
By Lemmas~\ref{lem:persist, 2}, 
every elementary map $M\to N$ is an embedding $\mathbf M\to\mathbf N$, and
every isomorphism of $\mathcal L$-structures is elementary.
Suppose the $\mathcal L_A$-formula $\varphi(x)$ defines the set
$X\subseteq M_x$. If $h\in\Aut(\mathbf M)$, then the $\mathcal L_{h(A)}$-formula $h(\varphi)(x)$ defines the set
$h(X)\subseteq M_x$, so if $h\in\Aut(\mathbf M|A)$, then  $X=h(X)$.
This observation is often used to show that certain relations are not definable:

\begin{example}
The usual ordering relation $\leq$ on the set of real numbers is not definable in the  
struc\-ture~$(\R; 0, {-},{+})$: for any $r_1,\dots,r_n\in\R$ there is an automorphism $\sigma$ of this structure such that $\sigma(r_i)=r_i$ for $i=1,\dots,n$ and $\sigma(r)<0<r$ for some $r\in\R$.
\end{example}

\subsection*{More on substitution} The next lemma shows that
substitution and evaluation in terms and formulas behave correctly.
%We finish this section with a technical result used in Sections~\ref{sec:ultra products} and \ref{sec:LS} below.
Let $x=(x_i)_{i\in I}$ and $y=(y_j)_{j\in J}$ be multivariables,
and let $t=(t_i)$ be a family of $\mathcal L$-terms with $t_i=t_i(y)$ of the same sort as $x_i$, for all $i\in I$.
Then for $a\in M_y$ we put $t^{\mathbf M}(a):=\big(t_i^\mathbf M(a)\big)\in M_x$.
 
\begin{lemma} \label{lem:subst lemma, formulas} 
Let $a\in M_y$, and let $\alpha(x)$ be an $\mathcal L_M$-term of sort $s$.
Then we have $\alpha(t)^{\mathbf M}(a)=\alpha^{\mathbf M}(t^{\mathbf M}(a))\in M_s$. Let $\varphi(x)$ be a
quantifier-free $\mathcal L_M$-formula. Then
$$  \mathbf M \models \varphi(t)(a)\quad \Longleftrightarrow\quad \mathbf M \models \varphi\big(t^{\mathbf M}(a)\big).$$
\end{lemma}
\begin{proof} The claim about $\alpha$ follows by induction on terms.    
Suppose $\varphi$ is atomic, say $\varphi=
R\alpha_1 \cdots \alpha_m$ with $m$-ary $R\in\mathcal L^{\operatorname{r}}$ and $\mathcal L_M$-terms $\alpha_1(x), \dots, \alpha_m(x)$. Then
\begin{align*} \varphi(t)\ =\ R\alpha_1(t)\cdots\alpha_m(t),&\qquad\varphi\big(t^{\mathbf M}(a)\big)\ =\  
R\alpha_1\big(t^{\mathbf M}(a)\big) \cdots 
\alpha_m\big(t^{\mathbf M}(a)\big), \text{ so}\\ 
\mathbf M\models \varphi(t)(a)\quad &\Longleftrightarrow\quad
\big(\alpha_1(t)^{\mathbf M}(a), \dots, \alpha_m(t)^{\mathbf M}(a)\big) \in R^{\mathbf M}\\
    &\Longleftrightarrow\quad 
\big(\alpha_1^{\mathbf M}(t^{\mathbf M}(a)), \dots, \alpha_m^{\mathbf M}(t^{\mathbf M}(a))\big) \in R^{\mathbf M}\\
 &\Longleftrightarrow\quad \mathbf M\models \varphi\big(t^{\mathbf M}(a)\big).
\end{align*} 
The case that $\varphi$ is $\alpha = \beta$ is handled 
the same way. The desired property is clearly inherited by disjunctions, conjunctions
and negations. 
\end{proof}

\subsection*{Notes and comments}
The definition (T1)--(T8) of the satisfaction relation goes back to  
Tarski's paper~\cite{Tarski35a}, but in the form above is closer to Tarski-Vaught~\cite{TarskiVaught}.

\section{Elementary Equivalence and Elementary Substructures}\label{sec:eleqelsub}

\noindent
The syntax (terms, formulas, sentences) and semantics (truth, definability) from the previous section will be used in this section to compare $\mathcal L$-structures.

\subsection*{Elementary equivalence}
We say that $\mathbf M$ and $\mathbf N$ are {\bf elementarily equivalent} (notation: $\mathbf M\equiv\mathbf N$) if they satisfy the same $\mathcal L$-sentences.
So isomorphic $\mathcal L$-structures 
are elementarily equivalent. By~\ref{ex:DLO} below, however, the non-isomorphic ordered sets $(\Q;{\leq})$ and $(\R;{\leq})$ are elementarily equivalent as well. This uses 
the \textit{back-and-forth method,}\/ a general tool that we also relied on in Chapter~\ref{ch:monotonedifferential}. \index{structure!elementary equivalence} \index{equivalence!elementary} \index{elementary!equivalence}

\nomenclature[Bj2]{$\mathbf M\equiv\mathbf N$}{$\mathbf M$, $\mathbf N$ are elementarily equivalent}

\begin{definition} \index{isomorphism!partial}\index{structure!partial isomorphism}
A {\bf partial isomorphism} from $\mathbf M$ to $\mathbf N$ is a bijection $\gamma\colon A\to B$
with $A\subseteq M$, $B\subseteq N$, such that
\begin{enumerate}
\item for each $R\in\mathcal L^{\operatorname{r}}$ of sort $s_1\dots s_m$ and $a\in A_{s_1\dots s_m}$,
$$a\in R^{\mathbf M}\quad\Longleftrightarrow\quad \gamma a \in R^{\mathbf N};$$
\item for each $f\in\mathcal L^{\operatorname{f}}$ of sort $s_1\dots s_n  s$ and $a\in A_{s_1\dots s_n}$, $b\in A_s$,
$$f^{\mathbf M}(a)=b \quad\Longleftrightarrow\quad f^{\mathbf N}(\gamma a)=\gamma b.$$
\end{enumerate}
\end{definition}

\noindent
Note that in this definition we do not assume that $A$ and $B$ are the underlying sets of substructures of   $\mathbf M$ and
$\mathbf N$, respectively. If
$A$ and $B$ \textit{are}\/ the underlying sets of substructures $\mathbf A$ of $\mathbf M$ and $\mathbf B$ of
$\mathbf N$, respectively, then  a partial isomorphism~$A\to B$ from $\mathbf M$ to $\mathbf N$
is the same thing as an isomorphism $\mathbf A\to\mathbf B$.

Given a partial isomorphism $\gamma\colon A\to B$ from $\mathbf M$ to $\mathbf N$, we set 
$\operatorname{domain}(\gamma_s):=A_s$ and
$\operatorname{domain}(\gamma):=A$, and likewise $\operatorname{codomain}(\gamma_s):=B_s$ and $\operatorname{codomain}(\gamma):=B$.
If~$(\gamma_\lambda)_{\lambda\in \Lambda}$ is a family of partial isomorphisms indexed by a
directed set $(\Lambda,{\leq})$, such that $\gamma_{\lambda'}$ extends $\gamma_\lambda$ for 
$\lambda\leq\lambda'$ in $\Lambda$, then there is a unique partial isomorphism $\gamma:=\bigcup_{\lambda\in\Lambda} \gamma_\lambda$
from $\mathbf M$ to $\mathbf N$ with domain $\bigcup_{\lambda\in\Lambda} \operatorname{domain}(\gamma_\lambda)$ such that $\gamma(a)=\gamma_\lambda(a)$ for 
all $s\in S$, $\lambda\in\Lambda$, and 
$a\in\operatorname{domain}(\gamma_{\lambda,s})$.

\begin{example}
Let $\mathbf M = (M; {\leq})$, $\mathbf N = (N; {\leq})$ be   ordered sets, and let
$a_1, \dots , a_n \in M$ and $b_1, \dots , b_n \in N$ with $a_1 < a_2 < \cdots < a_n$ and
$b_1 < b_2 < \cdots < b_n$; then the map $a_i \mapsto b_i \colon \{a_1, \dots , a_n\} \to \{b_1, \dots , b_n\}$ is a
partial isomorphism from~$\mathbf M$ to~$\mathbf N$.
\end{example}

\begin{definition} \index{back-and-forth system}\index{structure!back-and-forth system}
A {\bf back-and-forth system} from $\mathbf M$ to $\mathbf N$ is a collection $\Gamma$ of partial
isomorphisms from $\mathbf M$ to $\mathbf N$ such that $\Gamma\ne \emptyset$ and:
\begin{list}{}{\leftmargin=1em}
\item[(``Forth'')] for each $\gamma \in  \Gamma$, $s\in S$, and $a \in  M_s$ there is a $\gamma' \in  \Gamma$ extending
$\gamma$ such that $a \in  \operatorname{domain}(\gamma'_s)$;
\item[(``Back'')] for each $\gamma \in  \Gamma$, $s\in S$, and $b \in  N_s$ there is a $\gamma' \in  \Gamma$  extending
$\gamma$ such that $b \in  \operatorname{codomain}(\gamma'_s)$.
\end{list}
We say that $\mathbf M$ and $\mathbf N$ are {\bf back-and-forth equivalent} (notation: $\mathbf M \equiv_{\operatorname{bf}} \mathbf N$) if there
exists a back-and-forth system from $\mathbf M$ to $\mathbf N$. \index{structure!back-and-forth equivalent} \index{equivalence!back-and-forth}\nomenclature[Bj3]{$\mathbf M\equiv_{\operatorname{bf}}\mathbf N$}{$\mathbf M$, $\mathbf N$ are back-and-forth equivalent}
\end{definition}

\begin{prop}\label{prop:bf, 1}
If  $\mathbf M$, $\mathbf N$ are countable and $\mathbf M \equiv_{\operatorname{bf}} \mathbf N$, then~$\mathbf M\cong\mathbf N$.
\end{prop}
\begin{proof} 
Suppose $\Gamma$ is a back-and-forth system from $\mathbf M$ to $\mathbf N$. 
Note first that for $\gamma\in\Gamma$ and $s\in S$ we have:
$\operatorname{domain}(\gamma_s)=M_s$ iff $\operatorname{codomain}(\gamma_s)=N_s$.
In view of this fact, and assuming $M$ and $N$ are countable, we can start with any $\gamma_0\in \Gamma$ and build recursively
a sequence $(\gamma_n)$ in $\Gamma$ such that $\gamma_{n+1}$ extends $\gamma_n$ for all $n$ (going forth if $n$ is even, and going back if $n$ is odd) such that 
$\bigcup_n \operatorname{domain}(\gamma_n) = M$ and 
$\bigcup_n \operatorname{codomain}(\gamma_n) = N$.  
Then we have an isomorphism $\bigcup_n \gamma_n\colon \mathbf M\to\mathbf N$.
\end{proof}

\noindent
In applying this proposition and the next one in a concrete situation, the
key is to guess a back-and-forth system. That is where insight and imagination
(and experience) come in. In the following result we do not assume
countability.

\begin{prop}\label{prop:bf, 2}
If $\mathbf M \equiv_{\operatorname{bf}} \mathbf N$, then $\mathbf M \equiv  \mathbf N$.
\end{prop}

\noindent
Towards the proof, define an $\mathcal L$-formula to be {\bf unnested} \index{formula!unnested} \index{unnested formula}  if every atomic subformula of it is either $\top$, or $\bot$, or has one of the following forms: 
\begin{enumerate}
\item $Rx_1 \dots x_m$; here $R \in \mathcal L^{\operatorname{r}}$  and $x_1,\dots, x_m$ are distinct variables;
\item $x=y$; here $x$ and $y$ are distinct variables; 
\item $fx_1\dots x_n=y$; here $f \in \mathcal L^{\operatorname{f}}$ and $x_1,\dots, x_n,y$ are distinct
variables.
\end{enumerate}

\begin{lemma}\label{unnested}
Each atomic $\mathcal L$-formula $\varphi(x)$ is equivalent to an unnested existential $\mathcal L$-formula~$\varphi_\exists(x)$,
and also to an unnested universal $\mathcal L$-formula $\varphi_\forall(x)$. Each $\mathcal L$-for\-mu\-la~$\varphi(x)$ is equivalent
to an unnested $\mathcal L$-formula $\varphi_{\operatorname{u}}(x)$. 
\end{lemma}

\noindent
We leave the proof of this lemma to the reader.

\begin{proof}[Proof of Proposition~\ref{prop:bf, 2}]
Let $\Gamma$ be a back-and-forth system from $\mathbf M$ to~$\mathbf N$, and let $\varphi(x)$ be an unnested $\mathcal L$-formula. By induction on the
number of logical symbols in~$\varphi$ we show that
for all $\gamma\in\Gamma$ with domain $A$ and $a \in A_x$,
$$\mathbf M\models \varphi(a) \quad\Longleftrightarrow\quad \mathbf N \models\varphi(\gamma a).$$
This yields $\mathbf M \equiv  \mathbf N$, using Lemma~\ref{unnested} for sentences.
The stated equivalence follows from the definitions if $\varphi$ is atomic, and its validity is
preserved under~$\wedge$,~$\vee$,~$\neg$. Suppose $\varphi=\exists y\psi$ where $y$ is a quantifiable
variable of sort $s$ not occurring in $x$ and $\psi(x,y)$ is unnested, and let $\gamma\in\Gamma$ with domain $A$
and $a\in A_x$. If $\mathbf M\models \varphi(a)$, then we
take $b\in M_s$ such that $\mathbf M\models\psi(a,b)$,
and then $\gamma'\in\Gamma$ extending~$\gamma$ with $b\in\operatorname{domain}(\gamma'_s)$;
by inductive hypothesis $\mathbf N\models\psi(\gamma a,\gamma' b)$
and hence $\mathbf N\models\varphi(\gamma a)$.
Similarly one shows that $\mathbf N\models\varphi(\gamma a)\Rightarrow\mathbf M\models \varphi(a)$.
\end{proof}

\begin{exampleNumbered}\label{ex:DLO}
Let $\mathbf M=(M;{\leq})$ and $\mathbf N=(N;{\leq})$ be dense ordered sets without endpoints.
Then the collection of all strictly increasing bijections $A\to B$, where $A\subseteq M$
and $B\subseteq N$ are finite, is a back-and-forth system from $\mathbf M$ to $\mathbf N$.
Hence $\mathbf M\equiv\mathbf N$, and if $M$, $N$ are countable, then $\mathbf M\cong\mathbf N$.
\end{exampleNumbered}

\subsection*{Elementary substructures} \index{substructure!elementary} \index{elementary!substructure} \index{extension!elementary} \index{elementary!extension} Let $\mathbf M\subseteq\mathbf N$.
One says that $\mathbf M$ is an {\bf elementary substructure} of $\mathbf N$ (and that the extension $\mathbf M\subseteq\mathbf N$ is {\bf elementary}) if the natural inclusion
$\mathbf M\hookrightarrow\mathbf N$ is elementary; notation: $\mathbf M\preceq\mathbf N$.
We have  
$$\mathbf M\preceq\mathbf N\quad\Longleftrightarrow\quad \text{$\varphi^{\mathbf M}=\varphi^{\mathbf N}\cap M_x$
for each $\mathcal L_M$-formula $\varphi(x)$}.$$ 
Also, $\mathbf M\preceq\mathbf N \Longleftrightarrow \mathbf M_M\equiv\mathbf N_M$, so  $\mathbf M\preceq\mathbf N\Rightarrow\mathbf M\equiv\mathbf N$.
\nomenclature[Bj01]{$\mathbf M\preceq\mathbf N$}{$\mathbf M$ is an elementary substructure of $\mathbf N$}

\begin{exampleNumbered}\label{ex:simple gps}
View groups as $\mathcal L_{\operatorname{G}}$-structures  as in~\ref{ex:L-structures}(1).
Let  $G$, $H$ be groups with $G\preceq H$, and suppose that $H$ is simple. Then $G$ is also simple: to see this, let 
$g,g'\in G$, $g\neq 1$; it suffices to show that $g'$ is in the normal subgroup of~$G$ generated by~$g$.
Now $g'=h_1 g^{k_1} h_1^{-1}\cdots h_n g^{k_n} h_n^{-1}$
where $n\ge 1$, $h_1,\dots,h_n\in H$, $k_1,\dots,k_n\in\Z$. Therefore
$H$ satisfies the $\mathcal L_G$-sentence 
$$\exists x_1\cdots \exists x_n \big( g' = x_1 g^{k_1} x_1^{-1} \cdots x_n g^{k_n} x_n^{-1}\big),$$
and so does $G$, since $G\preceq H$.
\end{exampleNumbered}

\noindent
We have the following useful criterion for a substructure to be elementary: \index{test!Tarski-Vaught}

\begin{prop}[Tarski-Vaught test]\label{prop:TV}
Let $A\subseteq N$. Suppose that for every $\mathcal L_A$-formula~$\varphi(x)$ with $x$
a quantifiable variable, if $\mathbf N\models\exists x\varphi(x)$, then
$\mathbf N\models\varphi(a)$ for some $a\in A_x$.
Then $A$ underlies an elementary substructure of $\mathbf N$.
\end{prop}

\begin{proof} Note first that $N_s\ne \emptyset$ gives 
$A_s\neq\emptyset$, for all $s$.
If $f$ is a function symbol of~$\mathcal L$ of sort $s_1\dots s_n s$ and $a\in A_{s_1\dots s_n}$,
consider the $\mathcal L_A$-formula $\varphi(x)$ given by~${f(a)=x}$, with $x$ a quantifiable variable of sort~$s$;
then $\mathbf N\models\exists x\varphi$ and hence 
$f^{\mathbf N}(a)\in A_x$. Thus~$A$ is the underlying set of a substructure $\mathbf A$ of $\mathbf N$.
We show that for each $\mathcal L_A$-sentence~$\sigma$ we have
$\mathbf A\models\sigma\Longleftrightarrow\mathbf N\models\sigma$ by induction on the construction of $\sigma$.
If~$\sigma$ is atomic, then this holds by Corollary~\ref{cor:persist, 1}(i), and it is clear that the desired property is preserved under taking negations, conjunctions, and disjunctions.
It remains to treat the case where $\sigma=\exists x\varphi$ with
 an $\mathcal L_A$-formula $\varphi(x)$. Then
\begin{align*}
\mathbf A \models \sigma 	& \quad\Longleftrightarrow\quad \text{$\mathbf A\models\varphi(a)$
for some $a\in A_x$} \\
							& \quad\Longleftrightarrow\quad \text{$\mathbf N\models\varphi(a)$
for some $a\in A_x$}\\
							& \quad\Longleftrightarrow\quad \text{$\mathbf N\models\sigma$,}
\end{align*}
using the hypothesis in the proposition for the third equivalence.
\end{proof}

\begin{cor}
If $\mathbf M\subseteq \mathbf N$ and for all finite $A\subseteq M$ and all $b\in N_s$ there exists $h\in\Aut(\mathbf N|A)$
with $h(b)\in M_s$, then $\mathbf M\preceq\mathbf N$.
\end{cor}

\begin{example}
The previous corollary easily yields $(\Q;{\leq})\preceq (\R;{\leq})$.
\end{example}

\noindent
The Tarski-Vaught test can be used to construct small elementary substructures: \index{theorem!L\"owenheim-Skolem}

\begin{prop}[Downward L\"owenheim-Skolem]\label{prop:LS downwards}
Let $A\subseteq N$ and suppose~$\kappa$ is a cardinal such that
$\max\big\{ |A|,|\mathcal L|\big\} \leq \kappa \leq |N|$.
Then $\mathbf N$ has an elementary substructure~$\mathbf M$ of 
size $\kappa$ with $A\subseteq M$.
\end{prop}
\begin{proof}
After enlarging $A$ to a parameter set in $\mathbf N$ of size $\kappa$,
we may assume that $|A|=\kappa$. For every $B\subseteq N$ and $s$, let $\Phi_{B,s}$ be the set of all $\mathcal L_B$-formulas~$\varphi(x)$ with~$x$ a quantifiable variable of sort $s$ such that $\mathbf N\models \exists x\varphi$.
For every  $\varphi\in\Phi_{B,s}$ choose a $b_\varphi\in N_s$ with $\mathbf N\models\varphi(b_\varphi)$, and let $B'_s:=\{b_\varphi:\varphi\in\Phi_{B,s}\}$. 
For each~$b\in B_s$ and quantifiable variable $x$ of sort $s$ the $\mathcal L_B$-formula $x=b$ is in~$\Phi_{B,s}$,
so $b=b_{x=b}\in B'_s$, hence $B_s\subseteq B_s'$. 
Setting $B':=(B'_s)$, we have $|B|\leq |B'| \leq \max\big\{ |B|, |\mathcal L| \big\}$. We now
inductively define an increasing sequence $A_0\subseteq A_1\subseteq\cdots$ of parameter sets in~$\mathbf M$
by $A_0:=A$ and $A_{n+1}:=A_n'$, and put $M:=\bigcup_n A_n$. Then $|A_n|=\kappa$ for each $n$, hence also 
$|M|=\kappa$. By Proposition~\ref{prop:TV}, $M$ is the underlying set of an elementary substructure of~$\mathbf N$.
\end{proof}

\begin{exampleNumbered}\label{ex:LS downwards}
Any infinite simple group $H$ has a simple subgroup
of any given infinite size $\leq |H|$, by Example~\ref{ex:simple gps} and Proposition~\ref{prop:LS downwards}.
\end{exampleNumbered}

\subsection*{Direct unions}
Suppose $\mathbf M$ is the direct union of the family $(\mathbf M_\lambda)_{\lambda\in\Lambda}$  of substructures of $\mathbf M$. Let~$\lambda$,~$\lambda'$ range over $\Lambda$.
By (i) and (ii) of  Corollary~\ref{cor:persist, 1}, if $\sigma$ is a $\forall\exists$-sentence such that
$\mathbf M_\lambda\models\sigma$ for all~$\lambda$, then $\mathbf M\models\sigma$.
 
\begin{exampleNumbered}\label{ex:forallexists Fpalg}
Let $p$ be a prime number and $\sigma$ a universal-existential $\mathcal L_{\operatorname{R}}$-sentence. If $\sigma$ holds in
all sufficiently large finite fields of characteristic $p$, then $\sigma$ holds in the algebraic closure 
$\mathbb F_p^{\alg}$
of~$\mathbb F_p$; see also Example~\ref{ex:Fpalg}.
\end{exampleNumbered}

\begin{lemma}\label{lem:elementary chains}
Suppose  $\mathbf M_\lambda\preceq \mathbf M_{\lambda'}$ for all $\lambda\leq\lambda'$. Then $\mathbf M_\lambda\preceq\mathbf M$ for all $\lambda$.
\end{lemma} 
\begin{proof}
By induction on $n$ we show that for each $\lambda$ and each $\mathcal L_{M_\lambda}$-sentence $\sigma$ of length~$n$ we have $\mathbf M_\lambda\models\sigma\Longleftrightarrow \mathbf M\models\sigma$.
This is clear if~$\sigma$ is atomic by Corollary~\ref{cor:persist, 1}(i), and  the desired property is preserved under taking negations, conjunctions, and disjunctions.
So let $\varphi(x)$ be an $\mathcal L_{M_\lambda}$-formula, where~$x$ is a single quantifiable variable,
and $\sigma=\exists x\varphi$. 
Suppose first that $\mathbf M\models\sigma$, and take $a\in M_x$ with
$\mathbf M\models\varphi(a)$.
Then we can take some $\lambda'\geq\lambda$ such that  $a\in (M_{\lambda'})_x$.
By inductive hypothesis $\mathbf M_{\lambda'}\models\varphi(a)$ and so
$\mathbf M_{\lambda'}\models\sigma$, hence $\mathbf M_\lambda\models \sigma$ since~${\mathbf M_\lambda\preceq \mathbf M_{\lambda'}}$. Conversely, suppose 
 $\mathbf M_\lambda\models \sigma$, and take $a\in (M_\lambda)_x$ with
 $\mathbf M_\lambda\models\varphi(a)$; then $\mathbf M\models\varphi(a)$ by inductive
 hypothesis, hence $\mathbf M\models\sigma$.
 \end{proof}

\subsection*{Algebraic closure and definable closure} 
Let $A$ be a parameter set in $\mathbf M$ and~$b \in M_{\mathbf s}$. The tuple $b$ is said to be {\bf $A$-definable} in $\mathbf M$ (or {\bf definable over~$A$}
in~$\mathbf M$) if $\{b\}\subseteq M_{\mathbf s}$ is $A$-definable, and is said to be {\bf $A$-algebraic}
in $\mathbf M$ (or {\bf algebraic over $A$} in $\mathbf M$) if $b\in X$ for some finite $A$-definable set
$X\subseteq M_{\mathbf s}$. If $\mathbf M\preceq\mathbf N$, then $b$ is definable over $A$ in $\mathbf M$
iff $b$ is definable over $A$ in $\mathbf N$, and similarly with \textit{algebraic}\/ in place of \textit{definable.}\/
In the above we omit ``in $\mathbf M$'' if $\mathbf M$ is clear from the context.
Clearly \textit{$A$-definable}\/ implies \textit{$A$-algebraic.}\/\index{element!definable}\index{element!algebraic} \index{definable!tuple}\index{algebraic!tuple}

\begin{lemma}\label{lem:dcl and def fns} 
The tuple $b$ is $A$-definable if and only if $f(a)=b$ for some 
finite multivariable $x$, $0$-definable $X\subseteq M_x$, $a\in X\cap A_x$, and $0$-definable $f\colon X\to M_{\mathbf s}$. 
\end{lemma}
\begin{proof}
The ``if'' direction is obvious. 
For the ``only if'' direction, suppose $b$ is $A$-definable. Take a finite multivariable $x$, a $0$-definable set
$Y\subseteq M_x\times M_{\mathbf s}$, and an $a\in A_x$ such that $Y(a)=\{b\}$. Then the set 
$$X\ :=\ \big\{a'\in M_x:\ |Y(a')|=1\big\}\ \subseteq\ M_x$$
is $0$-definable, and
$f\colon X \to M_{\mathbf s}$ given by $Y(a')=\big\{f(a')\big\}$ does the job.
\end{proof}

\noindent
Let 
$\operatorname{dcl}(A)$
be the parameter set in $\mathbf M$ such that for every $s$,
$$\operatorname{dcl}(A)_s\ =\ \{b \in M_s:\  \text{$b$ is definable over $A$}\},$$ 
and define   $\operatorname{acl}(A)$ likewise, with 
\textit{algebraic}\/
instead of \textit{definable.}\/\nomenclature[Bm1]{$\operatorname{acl}(A)$}{algebraic closure of $A$}\nomenclature[Bm2]{$\operatorname{dcl}(A)$}{definable closure of $A$}\index{closure!definable}\index{closure!algebraic}\index{definable!closure}\index{algebraic!closure} Call $\operatorname{dcl}(A)$ (respectively, $\operatorname{acl}(A)$) the 
{\bf definable
closure} of $A$ in $\mathbf M$ (respectively, the {\bf algebraic closure} of~$A$ in $\mathbf M$). We say
that $A$ is {\bf definably closed}\index{closed!definably}  in $\mathbf M$ (respectively, {\bf algebraically closed}\index{algebraically!closed}\index{closed!algebraically} in $\mathbf M$) if
$\operatorname{dcl}(A) = A$ (respectively, $\operatorname{acl}(A) = A$).
It is easy to check that 
$\operatorname{dcl}(A)$ is
definably closed in $\mathbf M$, and that $\operatorname{acl}(A)$ is algebraically closed in $\mathbf M$.
Both $\operatorname{dcl}(A)$ and $\operatorname{acl}(A)$ do not change if we pass from $\mathbf M$
to an elementary extension.
Clearly $\operatorname{dcl}(A)\subseteq\operatorname{acl}(A)$, and
if $A$ is algebraically closed in~$\mathbf M$, then $A$ is definably closed in~$\mathbf M$. Note that if $\mathbf M$ is
one-sorted and there exists an $A$-definable total ordering
on the underlying set $M$ of $\mathbf M$, then $\operatorname{dcl}(A)=\operatorname{acl}(A)$.

\medskip
\noindent
Suppose $\mathcal L$ has for each $s$ a constant symbol of sort~$s$.
If $A$ is definably closed in~$\mathbf M$, then
$A$ underlies a substructure of $\mathbf M$. Hence both  
$\operatorname{dcl}(A)$ and
$\operatorname{acl}(A)$ underlie substructures of~$\mathbf M$, also denoted by
$\operatorname{dcl}(A)$ and $\operatorname{acl}(A)$, respectively, with $\<A\>_{\mathbf M}\subseteq \operatorname{dcl}(A)\subseteq\operatorname{acl}(A)$.

\medskip
\noindent
If $a\in M_{\mathbf s}$ is $A$-definable, then $a$ is $A_0$-definable for some finite subset $A_0$ of $A$, and similarly with \textit{algebraic}\/ in place of \textit{definable.}\/
If $A\subseteq B\subseteq M$, then $\operatorname{dcl}(A)\subseteq\operatorname{dcl}(B)$ and
$\operatorname{acl}(A)\subseteq\operatorname{acl}(B)$.

\begin{lemma}\label{lem:dcl, acl}
Let $b\in M_s$ and set $X:=\big\{\sigma(b):\sigma\in\Aut(\mathbf M|A)\big\}\subseteq M_s$.
If $b\in\operatorname{dcl}(A)_s$, then $X=\{b\}$, and if
$b\in\operatorname{acl}(A)_s$, then $X$ is finite.
\end{lemma}

\begin{example}
Let $K$ be a field,   let
$V$ be an infinite $K$-linear space viewed as an $\mathcal L_{K\!\operatorname{-mod}}$-structure
as in Example~\ref{ex:L-structures}(7), and let $A\subseteq V$. Then  $$\operatorname{dcl}(A)\ =\ \operatorname{acl}(A)\ =\ \text{the subspace of $V$ generated by $A$.}$$ To see this, let $W$ be the subspace of~$V$ generated by $A$.
If $V\neq W$, then $V\setminus W$ is infinite and is an orbit
for the action of $\Aut(V|A)$ on $V$. Now use 
Lemma~\ref {lem:dcl, acl}.
\end{example}

\subsection*{Notes and comments}
The notion of elementary equivalence is from~\cite{Tarski36}.
The fact about countable dense ordered sets without endpoints from~\ref{ex:DLO} is due to Cantor~\cite{Cantor}, with the back-and-forth proof found by Huntington~\cite{Huntington} and Hausdorff~\cite{Hausdorff14}. Proposition~\ref{prop:bf, 2} goes back to Ehrenfeucht~\cite{Ehrenfeucht} and Fra\"\i{}ss\'e~\cite{Fraisse}. (See \cite{Felgner02, Plotkin} for the history of \textit{back-and-forth.}\/)
Elementary extensions as well as Propositions~\ref{prop:TV} and~\ref{prop:LS downwards}  and Lem\-ma~\ref{lem:elementary chains} are from \cite{TarskiVaught}. 
 L\"owenheim~\cite{Loewenheim} and Skolem~\cite{Skolem20} had shown 
Proposition~\ref{prop:LS downwards} in the case where $\mathcal L$ is countable; in the stated form it appears in \cite{TarskiVaught}.
Example~\ref{ex:LS downwards} is from \cite[Section~3.1]{Hodges}.
Model-theoretic algebraic closure was introduced by A.~Robinson~\cite[p.~157]{Robinson65} and gained prominence through the work of
Morley~\cite{Morley} and Baldwin-Lachlan~\cite{BaldwinLachlan}.

\section{Models and the Compactness Theorem}\label{sec:compactness}

\noindent
\textit{In the rest of this appendix, unless indicated otherwise, $t$ is an $\mathcal L$-term, $\varphi$, $\psi$, and~$\theta$ are $\mathcal L$-formulas,
$\sigma$ is an $\mathcal L$-sentence, and $\Sigma$ is a set of $\mathcal L$-sentences.}\/ We drop the prefix
$\mathcal L$ in \textit{$\mathcal L$-term}\/, \textit{$\mathcal L$-formula}\/ and so on, unless this would cause confusion. 

\subsection*{Models}
We say that $\mathbf M$ is a {\bf model of $\Sigma$} or {\bf $\Sigma$ holds in $\mathbf M$} (notation: $\mathbf M\models\Sigma$)
if $\mathbf M\models\sigma$ for all $\sigma\in\Sigma$. \index{model} \index{structure!model}

\nomenclature[Bl1]{$\mathbf M\models\Sigma$}{$\mathbf M$ is a model of $\Sigma$}

To discuss examples it is convenient to introduce some notation. Suppose $\mathcal L$
contains (at least) the constant symbol $0$ and the binary function symbol $+$.
Given any terms~$t_1,\dots,t_n$ we define the term $t_1+\cdots+t_n$ inductively as follows:
it is the term~$0$ if $n = 0$, the term $t_1$ if $n = 1$, and the term $(t_1 + \cdots + t_{n-1})+t_n$
for~$n \geq 2$. We write~$nt$ for the term $t+\cdots+t$ with $n$ summands, in particular, $0t$
and $1t$ denote the terms $0$ and $t$, respectively. Suppose $\mathcal L$ contains the constant
symbol~$1$ and the binary function symbol $\cdot$ (the multiplication sign). Then we
have similar notational conventions for $t_1 \cdot \ldots \cdot t_n$ and $t^n$; in particular, for $n = 0$
both stand for the term~$1$, and $t^1$ is just $t$.

\begin{examplesNumbered}\label{ex:theories}
Fix three distinct quantifiable variables $x$, $y$, $z$.  
\begin{enumerate}
\item \textit{Groups}\/ are the  $\mathcal L_{\operatorname{G}}$-structures that are models of
\begin{flalign*}
\operatorname{Gr}\ &:=\ \big\{\,\forall x(x \cdot 1 = x),\ \forall x\big(x \cdot x^{-1} = 1\big),  \ 
 \forall x \forall y \forall z\big((x \cdot y) \cdot z = x \cdot (y \cdot z)\big)\big\}. &
\end{flalign*}
\item \textit{Abelian groups}\/ are the $\mathcal L_{\operatorname{A}}$-structures that are models of
\begin{flalign*}
\operatorname{Ab}\ &:=\ \big\{\, \forall x(x + 0 = x),\ \forall x\big(x + (-x) = 0\big),\ \forall x \forall y(x + y = y + x), &\\
&\hskip3em  \forall x \forall y \forall z\big((x + y) + z = x + (y + z)\big)\big\}. &
\end{flalign*}
\item \textit{Torsion-free abelian groups}\/ are the $\mathcal L_{\operatorname{A}}$-structures that are models of
\begin{flalign*}
\operatorname{Tf}\ &:=\ \operatorname{Ab}\ \cup\ \big\{ \forall x(nx = 0 \to x = 0) : n = 1, 2, 3, \dots \big\},&
\end{flalign*}
and \textit{divisible abelian groups}\/ are the $\mathcal L_{\operatorname{A}}$-structures that are models of
\begin{flalign*}
\operatorname{Div}\ &:=\ \operatorname{Ab}\ \cup\ \big\{ \forall x\exists y(ny = x) : n = 1, 2, 3, \dots \big\}.&
\end{flalign*}
\item \textit{Ordered sets}\/ are the $\mathcal L_{\operatorname{O}}$-structures that are models of
\begin{flalign*}
\operatorname{Or}\ &:=\ \big\{ \,
\forall x(x\leq x), \ 
\forall x\forall y\forall z \big( (x\leq y \wedge y \leq z) \to x\leq z\big),& \\
&\hskip3em \forall x\forall y \big( (x\leq y \wedge y\leq x) \to x=y\big),\ \forall x\forall y( x\leq y \vee y \leq x) \big\}.&
\end{flalign*}
Abbreviating $x\leq y\ \wedge \neg x = y$ by $x<y$, \textit{dense ordered sets without endpoints}\/ are the $\mathcal L_{\operatorname{O}}$-structures that are models of
\begin{flalign*}
\operatorname{DLO}\ &:=\ \operatorname{Or}\ \cup\ \big\{\,\forall x\forall y\exists z(x<y\to x<z \wedge z<y), \
\forall x\exists y\exists z ( y<x \wedge x<z) \big\}.&
\end{flalign*}
\item \textit{Ordered abelian groups}\/ are the $\mathcal L_{\operatorname{OA}}$-structures that are models of
\begin{flalign*}
\operatorname{OAb}\ &:=\ \operatorname{Or}\ \cup\ \operatorname{Ab}\ \cup\ \big\{\,\forall x \forall y \forall z
(x\leq y \to x+z\leq y+z) \big\}.&
\end{flalign*}
\item \textit{Rings}\/ are the $\mathcal L_{\operatorname{R}}$-structures that are models of
\begin{flalign*}
\operatorname{Ri}\ &:=\ \operatorname{Ab}\, \cup\, \big\{\forall x\forall y\forall z
\big( (x \cdot  y) \cdot  z = x \cdot  (y \cdot  z) \big),\ \forall x ( x \cdot  1 = x \wedge 1 \cdot  x = x),& \\
&\hskip5em \forall x\forall y\forall z \big(
x \cdot  (y + z) = x \cdot  y + x \cdot  z \wedge (x + y) \cdot  z = x \cdot  z + y \cdot  z\big)\big\}.&
\end{flalign*}
\item \textit{Fields}\/ are the $\mathcal L_{\operatorname{R}}$-structures that are models of
\begin{flalign*}
\operatorname{Fl}\ &:=\ \operatorname{Ri}\ \cup\ \big\{\,\forall x\forall y(x \cdot  y = y \cdot  x),\ 1 \neq 0,\ \forall x\big(x \neq 0 \to \exists y (x \cdot  y = 1)\big) \big\}.&
\end{flalign*}
\item \textit{Ordered rings}\/ are the $\mathcal L_{\operatorname{OR}}$-structures that are models of
\begin{flalign*}
\operatorname{ORi}\ &:=\ \operatorname{OAb}\ \cup\ \operatorname{Ri}\ \cup\ \big\{\forall x\forall y( 0 \leq x \wedge 0 \leq y \to 0\leq x\cdot y)\big\},&
\end{flalign*}
and \textit{ordered fields}\/ are the $\mathcal L_{\operatorname{OR}}$-structures that are models of
$\operatorname{OFl}:=\operatorname{ORi}\cup\operatorname{Fl}$.
\item \textit{Fields of characteristic $0$}\/ are the $\mathcal L_{\operatorname{R}}$-structures that are models of
\begin{flalign*}
\operatorname{Fl}(0)\ &:=\  \operatorname{Fl}\ \cup\ \{n1 \neq 0 : n = 2, 3, 5, 7, 11, \dots \},&
\end{flalign*}
and given a prime number $p$, \textit{fields of characteristic $p$}\/ are the $\mathcal L_{\operatorname{R}}$-structures that are models of $\operatorname{Fl}(p) := \operatorname{Fl}\ \cup\ \{p1=0\}$.
\item \textit{Algebraically closed fields}\/ are the $\mathcal L_{\operatorname{R}}$-structures that are models of
\begin{flalign*}
\operatorname{ACF}\ &:=\ \operatorname{Fl}\ \cup\ \big\{\forall u_1 \cdots \forall u_n \exists x(x^n+u_1x^{n-1}+\cdots+u_n = 0) : n \geq 2\big\}.&
\end{flalign*}
Here $u_1, u_2, u_3, \dots$ is some fixed infinite sequence of distinct quantifiable variables, distinct
also from $x$, and $u_ix^{n-i}$ abbreviates $u_i \cdot  x^{n-i}$, for $i = 1, \dots , n$.
\item Given a prime number $p$ or $p=0$, \textit{algebraically closed fields of characteristic~$p$}\/ are the $\mathcal L_{\operatorname{R}}$-structures that are
models of $\operatorname{ACF}(p)\ :=\  \operatorname{ACF} \cup \operatorname{Fl}(p)$.
\end{enumerate}
\end{examplesNumbered}

\nomenclature[Bi11]{$\operatorname{ACF}$}{theory of algebraically closed fields in the language $\mathcal L_{\operatorname{R}}$}
\nomenclature[Bi12]{$\operatorname{ACF(p)}$}{theory of algebraically closed fields of characteristic $p$ in the language $\mathcal L_{\operatorname{R}}$}

\noindent
Here is the important Compactness Theorem: \index{theorem!Compactness Theorem}

\begin{theorem}\label{thm:compactness, 1}
If every finite subset of~$\Sigma$ has a model, then $\Sigma$ has a model.
\end{theorem}

\noindent
We give the proof of this theorem in the next section. The rest of this section contains some
reformulations and 
simple but instructive applications of this theorem.

\subsection*{Logical consequence}
We say that $\sigma$ is a {\bf logical consequence} of $\Sigma$ (writ\-ten~${\Sigma\models\sigma}$) if $\sigma$
is true in every model of $\Sigma$. More generally, we say that a formula~$\varphi$ is 
a {\bf logical consequence} of~$\Sigma$ (notation: $\Sigma\models\varphi$) if $\mathbf M\models\varphi$ for all models $\mathbf M$ of~$\Sigma$.
%Thus $\models\varphi$ iff $\emptyset\models\varphi$.
We also write $\sigma\models\varphi$ instead of $\{\sigma\}\models\varphi$.
Note that $\Sigma\cup\{\sigma_1,\dots,\sigma_n\}\models\varphi$
iff $\Sigma\models (\sigma_1\wedge\cdots\wedge\sigma_n)\to\varphi$. \index{logical!consequence} 

\nomenclature[Bl2]{$\Sigma\models\sigma$}{$\sigma$ is a logical consequence of $\Sigma$}

\begin{example}
It is well-known that in any ring $R$ we have $r \cdot 0 = 0$ for all $r\in R$.
This can now be expressed as $\operatorname{Ri}\models  x \cdot 0 = 0$.
\end{example}

\noindent
Here is a version of the Compactness Theorem in terms of {\em logical consequence}:

{\sloppy

\begin{theorem} \label{thm:compactness, 2}
If $\Sigma\models\sigma$, then $\Sigma_0\models\sigma$
for some finite $\Sigma_0\subseteq \Sigma$.
\end{theorem}
\begin{proof} Suppose $\Sigma_0\not\models \sigma$ for all finite $\Sigma_0\subseteq\Sigma$. Then every finite subset of~${\Sigma\cup\{\neg\sigma\}}$ has a model, so by Theorem~\ref{thm:compactness, 1}, $\Sigma\cup\{\neg\sigma\}$ has a model, and thus~${\Sigma\not\models\sigma}$.
\end{proof}

}
\noindent
A similar argument as in the proof of Theorem~\ref{thm:compactness, 2} shows:

\begin{cor}\label{cor:compactness, 2}
Let $(\sigma_i)_{i\in I}$ and $(\tau_j)_{j\in J}$ be families of sentences such that $$ \bigwedge_{i\in I}\sigma_i \models\bigvee_{j\in J} \tau_j,$$
that is, in each structure where all sentences $\sigma_i$ are true, one of the sentences $\tau_j$ is true. Then there are
$i_1,\dots,i_m\in I$ and $j_1,\dots,j_n\in J$ such that $${\sigma_{i_1}\wedge\cdots\wedge\sigma_{i_m}}\models {\tau_{j_1}\vee\cdots\vee\tau_{j_n}}.$$
\end{cor}
\begin{proof}
The hypothesis expresses that $\{\sigma_i:i\in I\}\cup\{\neg\tau_j:j\in J\}$ has no model. By the Compactness Theorem 
there are $i_1,\dots,i_m\in I$ and $j_1,\dots,j_n\in J$  such that $\{\sigma_{i_1},\cdots,\sigma_{i_m},\neg \tau_{j_1},\dots,\neg \tau_{j_n}\}$ has no model, in other words, ${\sigma_{i_1}\wedge\cdots\wedge\sigma_{i_m}}\models {\tau_{j_1}\vee\cdots\vee\tau_{j_n}}$.
\end{proof}

\noindent
Here is one of many routine applications of the Compactness Theorem:

\begin{cor}\label{cor:Robinson}
If the $\mathcal L_{\operatorname{R}}$-sentence $\sigma$ is true in all fields of characteristic $0$, then~$\sigma$ is true in all
fields of sufficiently high prime characteristic.
\end{cor}
\begin{proof} If $\operatorname{Fl} \cup \{n1 \neq 0 : n \geq 1\} \models \sigma$, then 
compactness yields $N \in \N$ such that $\operatorname{Fl} \cup \{n1 \neq 0 : n = 1, \dots,N\} \models\sigma$, so
$\sigma$ holds in all fields of characteristic~${p > N}$.
\end{proof}

\noindent
Note that $\operatorname{Fl}(0)$ is
infinite; the previous corollary implies that there is no finite set of $\mathcal L_{\operatorname{R}}$-sentences whose models are
exactly the fields of characteristic $0$.
Here is a typical application of compactness (via Corollary~\ref{cor:Robinson}) in algebra:

\begin{exampleNumbered}[Noether-Ostrowski]\label{ex:NO}
Let $T=(T_1,\dots,T_n)$ and $P\in \Z[T]$.
Given a prime number $p$, let $P\bmod p$ denote the image of $P$
under the ring mor\-phism $\Z[T]\to\mathbb F_p[T]$ which extends
$a\mapsto a+p\Z\colon \Z\to\Z/p\Z=\mathbb F_p$ and sends~$T_i$ to~$T_i$~($i=1,\dots,n$).
\textit{Suppose $P$ is irreducible over $\C$. Then there is some $N\in\N$
such that for all primes~${p>N}$ and every field $F$ of characteristic $p$,
the polynomial~${P\bmod p}$ 
is irreducible over $F$.}\/
To see this note that every finitely generated field of characteristic zero can be embedded in $\C$,
so the hypothesis implies that $P$ is irreducible over every field of characteristic zero,
and irreducibility of $P$ over a given field  can be expressed by a sentence: there exists an $\mathcal L_{\operatorname{Ri}}$-sentence $\sigma_P$ such that for each field~$F$, if $F$ has characteristic zero then
$F\models\sigma_P\Longleftrightarrow$ $P$ is irreducible over $F$, and if~$F$ has characteristic~$p>0$ then
$F\models\sigma_P\Longleftrightarrow$ $P\bmod p$ is irreducible over~$F$.
\end{exampleNumbered}

\subsection*{Completeness and compactness}
We call $\Sigma$ {\bf complete} if $\Sigma$ has a model and for all $\sigma$, 
either $\Sigma\models\sigma$ or $\Sigma\models\neg\sigma$;
equivalently, $\Sigma$ has a model and $\mathbf M\equiv\mathbf N$ for all models $\mathbf M$, $\mathbf N$
of $\Sigma$.
Completeness is a strong property and it can be hard to show that a given
set $\Sigma$ is complete. \index{complete!set of sentences}

\begin{examples}
The set of $\mathcal L_{\operatorname{A}}$-sentences $\operatorname{Ab}$ (the set of axioms for abelian groups)
is not complete: consider $\exists x(x\neq 0)$.
The set of $\mathcal L_{\operatorname{O}}$-sentences~$\operatorname{DLO}$ (the set of axioms for dense ordered
sets without endpoints) is complete; see Example~\ref{ex:DLO}.
\end{examples}

\nomenclature[Bm3]{$\operatorname{Th}(\mathbf M)$}{theory of $\mathbf M$}
\nomenclature[Bm31]{$\operatorname{Th}(\mathcal C)$}{theory of a class $\mathcal C$ of $\mathcal L$-structures}
\nomenclature[Bm32]{$\operatorname{Th}(\Sigma)$}{set of logical consequences of $\Sigma$}

\noindent
An {\bf $\mathcal L$-theory} is a set $T$ of 
$\mathcal L$-sentences such that for all
$\sigma$, if $T \models \sigma$, then $\sigma \in T$.\index{theory}
In particular, we have for any $\Sigma$ the $\mathcal{L}$-theory generated by
$\Sigma$:
$$\operatorname{Th}(\Sigma)\ :=\ \{\sigma:\ \Sigma\models\sigma\}.$$
It has the same models as $\Sigma$. An {\bf axiomatization\/} of an 
$\mathcal L$-theory $T$ is a set $\Sigma$ such that 
$\operatorname{Th}(\Sigma)=T$. We use the abbreviation 
$\operatorname{Th}$ also as follows: for any $\mathbf M$, the set 
$$\operatorname{Th}(\mathbf M)\ :=\ \{\sigma:\mathbf M\models\sigma\}$$ is a complete $\mathcal L$-theory,
called the {\bf theory of~$\mathbf M$}.\index{theory!structure}  (Thus if $\Sigma$ has a model, then $\Sigma\subseteq T$ for some complete $\mathcal{L}$-theory $T$.)\index{axiomatization}\index{theory!axiomatization} Given a class $\mathcal C$ of $\mathcal L$-structures, we set
$$\operatorname{Th}(\mathcal C)\ :=\ \bigcap_{\mathbf M\in\mathcal C} \operatorname{Th}(\mathbf M)\ =\ \{\sigma:\ \sigma \text{ is true in every $\mathbf M\in \mathcal C$}\}.$$
This is an $\mathcal L$-theory (not necessarily complete), called the {\bf theory of $\mathcal C$}. Thus $\operatorname{Th}(\Sigma)=\operatorname{Th}(\mathcal C)$,
where $\mathcal C$ is the class of models of $\Sigma$.\index{theory!class of structures} 

\medskip
\noindent
A complete $\mathcal L$-theory that contains $\Sigma$ is called a {\bf completion} of $\Sigma$. 
Let $\operatorname{S}(\Sigma)$ be the set of completions of $\Sigma$; thus
$\abs{\operatorname{S}(\Sigma)}\le 2^{\abs{\mathcal L}}$. We set
$\<\sigma\>:=\{ T \in \operatorname{S}(\Sigma): \sigma\in T\}$.
One verifies easily that for sentences $\sigma$, $\tau$ we have
$\<\sigma\>=\<\tau\>$ iff $\Sigma\models\sigma\leftrightarrow\tau$, and
$$\<\sigma\wedge\tau\>=\<\sigma\>\cap\<\tau\>,\quad \<\sigma\vee\tau\>=\<\sigma\>\cup\<\tau\>,\quad
\<\neg\sigma\>= \operatorname{S}(\Sigma)\setminus\<\sigma\>.$$
The topology on $\operatorname{S}(\Sigma)$ with
the sets $\<\sigma\>$ as a basis is called the {\bf Stone topology}. \index{topology!Stone}
Its open sets are the unions
$\bigcup_{\sigma\in\Delta} \<\sigma\>$ with $\Delta$ a set
of sentences. Note that the basic open sets $\<\sigma\>$ are also
closed, and so the Stone topology on~$\operatorname{S}(\Sigma)$  is hausdorff.  
We have $\operatorname{S}(\Sigma)\neq
\bigcup_{\sigma\in\Delta} \<\sigma\>$ iff $\Sigma\cup\{\neg\sigma:\sigma\in\Delta\}$ has a model; 
hence Theorem~\ref{thm:compactness, 1} also has the following reformulation, explaining the name \textit{Compactness Theorem.}\/\index{completion!set of sentences}

\begin{theorem}\label{thm:compactness, 3}
The hausdorff space $\operatorname{S}(\Sigma)$ is compact.
\end{theorem}

\noindent
Let $\Sigma'$ be a set of $\mathcal{L}$-sentences.
Then $\Sigma$ and $\Sigma'$ are said to be {\bf equivalent} if
they have the same logical consequences. Thus $\Sigma$ and $\Sigma'$ are equivalent iff they  have the same models, iff $\operatorname{S}(\Sigma)=\operatorname{S}(\Sigma')$. For example, $\Sigma$ and
$\operatorname{Th}(\Sigma)$ are equivalent. 

\index{equivalence!sets of sentences}

\subsection*{Completeness and decidability} This subsection concerns the relation between completeness and (algorithmic) decidability, a logical issue that is hardly model-theoretic in nature. We just give
an outline and refer to the literature for details, since decidability 
only makes an appearance in Corollary~\ref{cor:completionsdecidable}. 

First, one should distinguish the notion of \textit{logical consequence of}\/ from that of
\textit{provable from.}\/ To make the latter concept precise requires a \textit{proof system}\/, which specifies logical axioms and
inference rules for generating certain finite sequences of formulas, called (formal) proofs:
a \textit{proof of $\varphi$ from $\Sigma$}\/ is a sequence $\varphi_1,\dots,\varphi_n$
of formulas with~$n \geq 1$ and $\varphi_n = \varphi$, such that for $k = 1, \dots, n$,
either $\varphi_k\in\Sigma$ or $\varphi_k$ is a logical axiom, or $\varphi_k$ is ``inferred'' from some of the earlier formulas
$\varphi_1,\dots,\varphi_{k-1}$ by applying an inference rule. (For example, among the logical axioms might be all formulas of the form $\varphi\vee \neg \varphi$, and among the inference rules is usually {\em Modus Ponens,}\/ which allows
one to infer $\psi$ from $\phi$ and $\phi\to \psi$; 
see for example \cite[Section~2.6]{Shoenfield} for an explicit proof system.)
We call $\varphi$ \textit{provable from $\Sigma$}\/  (in symbols: $\Sigma\vdash \varphi$)
if there exists a proof of $\varphi$
from~$\Sigma$.
The logical axioms and inference rules are chosen so that the logical axioms are valid in all $\mathcal L$-structures, and 
if $\varphi$ is inferred from $\varphi_1,\dots,\varphi_{k}$ by an inference rule, then $\varphi$ is valid in all 
$\mathcal{L}$-structures where
$\varphi_1,\dots,\varphi_{k}$ are valid. Hence $\Sigma\vdash\varphi\Rightarrow \Sigma\models\varphi$. We now fix
some traditional proof system for classical predicate logic.
Then the converse also holds: \index{theorem!Completeness Theorem}

\begin{theorem}[G\"odel's Completeness Theorem]\label{thm:Completeness}
If $\Sigma\models\varphi$, then $\Sigma\vdash\varphi$.
\end{theorem}

\noindent
See \cite[Chapter~4]{Shoenfield} for a proof. 
``Completeness'' here refers to the proof system, not to~$\Sigma$.
Theorem~\ref{thm:Completeness} and its converse immediately yield version \ref{thm:compactness, 2}
of the Compactness Theorem (of which we give
an independent proof in the next section).

\medskip
\noindent
Suppose now that the language $\mathcal L$ has only finitely many nonlogical symbols. Then formulas can be made
into inputs of computer programs, and the logical axioms and inference rules of the proof system hiding behind the notation $\vdash$ can also be effectively given.  
One says that $\Sigma$ is \textit{effectively enumerable}\/ if there is an effective procedure that enumerates all elements of $\Sigma$. 
We say that an $\mathcal{L}$-theory $T$ is \textit{decidable}\/ if there is 
an algorithm (program) that takes any sentence $\sigma$ as input and decides whether or not~$\sigma\in T$.
If $\Sigma$ is effectively enumerable,
then so is the set $\operatorname{Th}(\Sigma)$ of all sentences provable from~$\Sigma$, and we obtain:

\begin{cor}\label{cor:Completeness}
If $\mathcal L$ has only finitely many nonlogical symbols and $\Sigma$ is complete and effectively enumerable, then $\operatorname{Th}(\Sigma)$ is decidable.
\end{cor}
 
\subsection*{Notes and comments}
The notion of ``model'' goes back to the Hilbert school. (But it was Tarski~\cite{Tarski54} who first spoke of ``model theory'' as a subject in its own right.)
Theorem~\ref{thm:Completeness} was  proved by  G\"odel~\cite{Goedel30} (1930)
for countable languages and by Hen\-kin~\cite{Henkin49} (1949) in general.
The Compactness Theorem  was shown independently by Mal$'$cev~\cite{Malcev36}~(1936). The formulation in Theorem~\ref{thm:compactness, 3} is due to Tarski~\cite{Tarski52}.
Corollary~\ref{cor:Robinson} is due to A.~Robinson~\cite{Robinson65}.
Algebraic proofs of the statement in Example~\ref{ex:NO} were given by Noether~\cite{Noether22} and Ostrowski~\cite{Ostrowski19}.

\section{Ultraproducts and Proof of the Compactness Theorem}\label{sec:ultra products}

\noindent
{\em In this section $\Lambda$ is a nonempty index set and
$\lambda$, $\lambda'$ range over $\Lambda$.}

\subsection*{Filters and ultrafilters}
A {\bf proper filter} on $\Lambda$ is a nonempty collection $\mathcal F$ of subsets of $\Lambda$ such that $\emptyset\notin \mathcal F$ and for all $A,B\subseteq \Lambda$:
\begin{list}{}{\leftmargin=2.5em}
\item[(Fi1)] if $A,B\in\mathcal F$, then $A\cap B\in\mathcal F$; 
\item[(Fi2)] if $A\subseteq B\subseteq \Lambda$ and $A\in\mathcal F$, then $B\in\mathcal F$.
\end{list}
Note that then $\Lambda\in\mathcal F$, and
that for all $A\subseteq \Lambda$,  $A\notin\mathcal F$ or $\Lambda\setminus A\notin\mathcal F$.  

\begin{example}
If $\Lambda$ is infinite, then the set of all cofinite subsets of $\Lambda$ is a proper filter on~$\Lambda$, called the {\bf Fr\'echet filter} on $\Lambda$. \index{filter!proper}\index{filter!Fr\'echet}
\end{example}

\noindent
A  proper filter on $\Lambda$ that is maximal with respect to inclusion is called an {\bf ultrafilter} on~$\Lambda$. By Zorn, every proper filter on $\Lambda$ is included in an ultrafilter on $\Lambda$. \index{filter!ultrafilter}\index{ultrafilter}

\begin{lemma}\label{lem:char ultrafilters}
Suppose $\mathcal F$ is a proper filter on $\Lambda$. Then the following are equivalent:
\begin{enumerate}
\item[\textup{(i)}] $\mathcal F$ is an ultrafilter on $\Lambda$;
\item[\textup{(ii)}] for all $A,B\subseteq \Lambda$ with $A\cup B\in\mathcal F$ we have $A\in\mathcal F$ or $B\in\mathcal F$;
\item[\textup{(iii)}] for all $A\subseteq \Lambda$
we have $A\in\mathcal F$ or $\Lambda\setminus A\in\mathcal F$.
\end{enumerate}
\end{lemma}
\begin{proof}  
Suppose $\mathcal F$ is an ultrafilter on $\Lambda$ and $A,B\subseteq \Lambda$, $A\cup B\in\mathcal F$. Then $A\cap F\neq\emptyset$ for all
$F\in\mathcal F$, or $B\cap F\neq\emptyset$ for all
$F\in\mathcal F$; we may assume that the first alternative holds. Then
$\{C\subseteq \Lambda:\   \text{$C\supseteq A\cap F$ for some $F\in \mathcal F$}\}$ is a proper filter on~$\Lambda$ which includes
$\mathcal F\cup \{A\}$, so $A\in\mathcal F$. This shows (i)~$\Rightarrow$~(ii), and
(ii)~$\Rightarrow$~(iii) follows by taking $B=\Lambda\setminus A$. The direction  (iii)~$\Rightarrow$~(i) is obvious.
\end{proof}

\noindent
A nonempty collection $\mathcal F$ of subsets of $\Lambda$ has the {\bf finite intersection property}~(FIP) if for all $A_1,\dots,A_n\in\mathcal F$ with $n\geq 1$ we have $A_1\cap\cdots\cap A_n\neq\emptyset$. In particular, any proper filter
on $\Lambda$ has the FIP. \index{finite intersection property}

\begin{lemma}\label{lem:extend to uf}
Suppose the nonempty collection $\mathcal F$ of subsets of $\Lambda$ has the \textup{FIP}.
Then there exists an ultrafilter $\mathcal U\supseteq\mathcal F$ on~$\Lambda$.
\end{lemma}
\begin{proof}
Let $\mathcal F^*$ be the collection of all subsets of $\Lambda$ that contain a finite in\-ter\-sec\-tion ${A_1\cap\cdots\cap A_n}$
with $A_1,\dots,A_n\in\mathcal F$, $n\ge 1$. Then~$\mathcal F^*$ is a proper filter and $\mathcal F^*\supseteq \mathcal F$. By Zorn there exists an ultrafilter $\mathcal U\supseteq \mathcal F^*$.
\end{proof}

\subsection*{Ultraproducts}
Let  $(\mathbf M_\lambda)$ be a family of $\mathcal L$-structures and $\mathbf M:=\prod_\lambda \mathbf M_\lambda$ be its product.  (See Section~\ref{sec:languages}.)
Let $\sigma$ be an $\mathcal L_M$-sentence. Take
an $\mathcal L$-formula $\varphi(x)$ and $a\in M_{x}$ such that $\sigma=\varphi(a)$, and set
$$\| \sigma \|\ :=\  \big\{ \lambda:\ \mathbf M_\lambda \models \varphi(a(\lambda)) \big\}.$$
This notation is justified since the $\mathcal L_{M_{\lambda}}$-sentence $\varphi(a(\lambda))$ depends only on $\sigma$ and $\lambda$, not on the choice of $\varphi(x)$ and $a$ such that 
$\sigma=\varphi(a)$. 
If  $\sigma$, $\tau$ are $\mathcal L_M$-sentences, then
\begin{equation}\label{eq:boolean extension}
\|\sigma\ \&\ \tau\| = \|\sigma\|\cap\|\tau\|,\quad \|\sigma \vee \tau\| = \|\sigma\|\cup\|\tau\|, \quad
\|\neg\sigma\| = \Lambda \setminus\|\sigma\|,
\end{equation}
\begin{equation}\label{eq:implication-inclusion}
\models \sigma\to\tau\quad\Rightarrow\quad\|\sigma\|\subseteq\|\tau\|.
\end{equation}
Moreover:

\begin{lemma}\label{lem:witnesses}
Let $\psi(y)$ be an $\mathcal L_M$-formula, where $y$ is a quantifiable variable of sort~$s$. Then
for all $b\in M_s$ we have
$\|\psi(b)\| \subseteq \| \exists y\psi\|$,
and there exists $b\in M_s$ such that equality holds.
\end{lemma}
\begin{proof}
For all $b\in M_s$ we have
$\|\psi(b)\| \subseteq \| \exists y\psi\|$ by \eqref{eq:implication-inclusion}.
To obtain equality,  choose $b=(b(\lambda))\in M_s$ as follows:
Take an $\mathcal L$-formula $\varphi(x,y)$
and $a\in M_{x}$
such that $\psi(y)=\varphi(a,y)$. 
Then we have for any $b\in M_s$,
$$
\|\psi(b)\| = \big\{\lambda: 
\mathbf M_\lambda \models \varphi(a(\lambda),b(\lambda)) \big\}, \quad
\|  \exists y\psi\|=\big\{\lambda: 
\mathbf M_\lambda \models (\exists y \varphi)(a(\lambda))\big\}.$$
For $\lambda\in \|  \exists y\psi\|$ we pick $b(\lambda)\in (M_\lambda)_s$ such that
$\mathbf M_\lambda\models \varphi\big(a(\lambda),b(\lambda)\big)$, and for $\lambda\notin \|  \exists y\psi\|$ we let $b(\lambda)\in (M_\lambda)_s$ be arbitrary.
Then 
$\|\psi(b)\|\supseteq \|  \exists y\psi\|$
as required.
\end{proof}

\noindent
Let $\mathcal F$ be a proper filter on $\Lambda$.
For  $s\in S$ we define a binary relation $\sim_s$ on $M_s$ by
$$a \sim_s b \quad :\Longleftrightarrow \quad \|a = b \| \in \mathcal F
\quad \Longleftrightarrow \quad \big\{\lambda:a(\lambda)=b(\lambda)\big\}\in\mathcal F.$$
For $a=(a_1,\dots,a_m),b=(b_1,\dots,b_m)\in M_{s_1\dots s_m}$
we set
$$a \sim_{s_1\dots s_m} b \quad :\Longleftrightarrow \quad a_1\sim_{s_1} b_1\ \&\ \cdots \ \&\ a_m\sim_{s_m} b_m.$$
Using \eqref{eq:implication-inclusion} one easily shows:

\begin{lemma} \mbox{}
\begin{enumerate}
\item[\textup{(i)}] The relation $\sim_{s_1\dots s_m}$ is an equivalence relation on $M_{s_1\dots s_m}$.
\item[\textup{(ii)}] If $R \in \mathcal L^{\operatorname{r}}$ has
sort $s_1\dots s_m$, and $a,b\in M_{s_1\dots s_m}$, 
$a\sim_{s_1\dots s_m} b$, then $$\|R\,a\|\in\mathcal F\quad\Longleftrightarrow\quad\|R\,b\|\in\mathcal F.$$
\item[\textup{(iii)}] If $f \in \mathcal L^{\operatorname{f}}$ has sort $s_1\dots s_n s$ and $a,b\in M_{s_1\dots s_n}$, $a\sim_{s_1\dots s_n} b$, then 
$$f^{\mathbf M}(a)\ \sim_s\  f^{\mathbf M}(b).$$
\end{enumerate}
\end{lemma}

\noindent
For $a\in M_{s_1\dots s_m}$ we let $a^{\sim}$ denote the equivalence class of $a$ with respect to $\sim_{s_1\dots s_m}$, and 
we let~$M^\sim_{s_1\dots s_m}$ be the set of equivalence classes of $\sim_{\mathbf s}$.
We identify $M^\sim_{s_1\dots s_m}$ with $M^\sim_{s_1}\times\cdots\times M^\sim_{s_m}$ in the natural way.
We now define an $\mathcal L$-structure $\mathbf M^{\sim}$ whose underlying set of sort~$s$ is $M^\sim_s$: 
for $R\in\mathcal L^{\operatorname{r}}$ of sort $s_1\dots s_m$ and $a\in M_{s_1\dots s_m}$ we set
$$a^\sim\in R^{\mathbf M^\sim}\quad:\Longleftrightarrow\quad \| R\,a\|\in\mathcal F,$$
and for $f \in \mathcal L^{\operatorname{f}}$ of sort $s_1\dots s_n s$ and $a\in M_{s_1\dots s_n}$ we put
$$f^{\mathbf M^\sim}(a^\sim)\ :=\ f^{\mathbf M}(a)^\sim\in M_s.$$
The $\mathcal L$-structure $\mathbf M^{\sim}$ is called the {\bf reduced product} of
$(\mathbf M_\lambda)$ with respect to $\mathcal F$, and is also denoted by
$\mathbf M/\mathcal F=\left(\prod_{\lambda} \mathbf M_\lambda\right)/\mathcal F$.
If $\mathbf M_\lambda=\mathbf N$ for all $\lambda$, then $\mathbf N^\Lambda/\mathcal F$ is also called the
{\bf reduced power} of $\mathbf N$ with respect to $\mathcal F$. 
If   $\mathcal F$ is an ultrafilter on~$\Lambda$, then we speak of the {\bf ultraproduct~$\mathbf M/\mathcal F$}
of $(\mathbf M_\lambda)$ with respect to $\mathcal F$
and of the {\bf ultrapower~$\mathbf N^\Lambda/\mathcal F$} of $\mathbf N$ with respect to $\mathcal F$.
The maps $a\mapsto a^\sim\colon M_s\to M_s^\sim$ combine to a surjective morphism $\pi\colon\mathbf M\to \mathbf M/\mathcal F$.
\index{structure!reduced product}\index{reduced product}\index{structure!ultraproduct}\index{ultraproduct}\index{structure!ultrapower}\index{ultrapower}

\begin{lemma}\label{lem:Los}
Let $x$ be finite and $a\in M_{x}$. Then for any $\mathcal L$-term  $t(x)$ we have $t(a^\sim)^{\mathbf M/\mathcal F} =t^{\mathbf M/\mathcal F}(a^\sim) =
{t^{\mathbf M}(a)}^\sim$ and for any atomic $\mathcal L$-formula $\varphi(x)$ we have
$$\mathbf M/\mathcal F\models \varphi(a^\sim) \quad\Longleftrightarrow\quad \| \varphi(a) \|\in\mathcal F.$$
\end{lemma}
\begin{proof} The first equality about terms is part of Lemma~\ref{lem:subst lemma, terms}, and the second equality follows by an easy induction on terms. The statement about atomic formulas
is a routine consequence of the equalities about terms. 
%Suppose first that $\varphi=Rt_1\dots t_m$ where  $R\in\mathcal L^{\operatorname{r}}$ has sort $s_1\dots s_m$
%and $t_1,\dots,t_m$ are $\mathcal L$-terms of sort $s_1,\dots,s_m$, respectively. 
%For notational simplicity we assume $m=1$ and set $t_1=t$.
%We have $t(a^\sim)^{\mathbf M/\mathcal F} = t^{\mathbf M/\mathcal F}(a^\sim) =
%{t^{\mathbf M}(a)}^\sim$,
%using Lemma for the first equality, hence
%$$\mathbf M/\mathcal F\models\varphi(a^\sim)	\ \Longleftrightarrow\ 
%t(a^\sim)^{\mathbf M/\mathcal F}={t^{\mathbf M}(a)}^\sim\in
%R^{\mathbf M/\mathcal F} 
%\ \Longleftrightarrow\  \big\| R\, t^{\mathbf M}(a)\big\|\in\mathcal F.$$
%Also, by Lemma~\ref{lem:subst lemma, formulas}:
%$$\big\| R\, t^{\mathbf M}(a)\big\|	=
%\left\{ \lambda: \mathbf M_\lambda \models R\,t^{\mathbf M_\lambda}(a(\lambda)) \right\}
%= \left\{ \lambda: \mathbf M_\lambda \models (R\,t)(a(\lambda)) \right\} =  \| \varphi(a) \|.$$
%If $\varphi$ is of the form $t_1=t_2$ with $\mathcal L$-terms $t_1$, $t_2$,
%then one argues similarly.
\end{proof}
 
 \begin{cor}
Let  $\Delta\colon\mathbf N\to\mathbf N^\Lambda$ be the diagonal embedding.
Then the morphism $\pi\circ\Delta\colon \mathbf N\to\mathbf N^\Lambda/\mathcal F$ is an  embedding, called the
{\bf diagonal embedding} of~$\mathbf N$ into its reduced power $\mathbf N^\Lambda/\mathcal F$. \index{embedding!diagonal}\index{diagonal!embedding}
\end{cor}

\noindent
Here is the main fact about ultraproducts: 

\begin{theorem}[{\L}o\'{s}]\label{thm:Los}
Suppose $\mathcal U$ is an ultrafilter on $\Lambda$.
Let $\varphi(x)$ be a formula with finite $x$, and let 
$a\in M_{x}$. Then
$\mathbf M/\mathcal U\models \varphi(a^\sim) \Longleftrightarrow \| \varphi(a) \|\in\mathcal U$.
\end{theorem}
\begin{proof}
We proceed by induction on the construction of $\varphi$. The case where $\varphi$ is atomic is covered by Lemma~\ref{lem:Los}. The cases $\varphi=\neg\psi$, $\varphi=\psi_1\ \&\ \psi_2$, and $\varphi=\psi_1\vee\psi_2$ 
with formulas $\psi$, $\psi_1$, $\psi_2$
follow by induction, \eqref{eq:boolean extension}, and Lemma~\ref{lem:char ultrafilters}.
Suppose $\varphi=\exists y\psi$ where $\psi(x,y)$ is a formula and $y$ is  of sort $s$; then
\begin{align*}
\mathbf M/\mathcal U\models \varphi(a^\sim)	&\quad\Longleftrightarrow\quad
\text{$\mathbf M/\mathcal U\models \psi(a^\sim,b^\sim)$ for some $b\in M_s$}  \\
&\quad\Longleftrightarrow\quad \text{$\| \psi(a,b) \| \in\mathcal U$ for some $b\in M_s$,}
\end{align*}
by inductive hypothesis. By Lemma~\ref{lem:witnesses}, $\| \psi(a,b) \| \in\mathcal U$ for some $b\in M_s$, if and only if
$\| \varphi(a) \|\in\mathcal U$. The case 
$\varphi=\forall y\psi$ follows from
$\models \forall y\psi \leftrightarrow \neg\exists y\neg\psi$.
\end{proof}

\noindent
In particular, if $\mathcal U$ is an ultrafilter on $\Lambda$ and $\mathbf M_\lambda\models\Sigma$ for all $\lambda$, then $\mathbf M/\mathcal U\models\Sigma$.

\begin{example}
View abelian groups as $\mathcal L_{\operatorname{A}}$-structures in the natural way.
Let $\mathcal U$ be an ultrafilter on $\N^{\geq 1}$ containing the Fr\'echet filter, set
$G:=\prod_{n\geq 1} \left(\Z/n\Z\right)/\mathcal U$, and
$g:=(1+n\Z)^\sim\in G$. Then for every $m\ge 1$ we have $G\not\models m g=0$ by Theorem~\ref{thm:Los}. So $g$ has infinite order in the abelian group $G$.
Thus there is no set of $\mathcal L_{\operatorname{A}}$-sentences whose models are the abelian torsion groups.
\end{example}

\begin{cor} If $\mathcal U$ is an ultrafilter on $\Lambda$, then the diagonal embedding of~$\mathbf N$ into $\mathbf N^\Lambda/\mathcal U$ is elementary.
\end{cor}
\begin{proof}
Let $\varphi(x)$ be a formula with finite $x$ and $a\in N_x$ such that $\mathbf N \models
\varphi(a)$. Then $\| \varphi(\Delta(a)) \| = \Lambda\in\mathcal U$, so
$\mathbf N^\Lambda/\mathcal U\models \varphi(\Delta(a)^\sim)$
by Theorem~\ref{thm:Los}.
\end{proof}

\subsection*{Proof of the Compactness Theorem}
We mean here Theorem~\ref{thm:compactness, 1}.
Thus, suppose every finite subset of $\Sigma$ has a model; we need to show that $\Sigma$ has a model.
We take $\Lambda$ to be the set of all finite subsets of~$\Sigma$;
for each $\lambda$ we take a model $\mathbf M_\lambda$ of~$\lambda$ and set
$F(\lambda):=\{\lambda' : \lambda\subseteq\lambda'\}\subseteq \Lambda$.
Then $F(\lambda_1)\cap F(\lambda_2)=F(\lambda_1\cup\lambda_2)$ for all $\lambda_1,\lambda_2\in\Lambda$,
so $\mathcal F:=\big\{F(\lambda) : \lambda\in\Lambda\big\}$ has the finite intersection property.
By Lemma~\ref{lem:extend to uf} we have an ultrafilter $\mathcal U$ on $\Lambda$ with $\mathcal U \supseteq \mathcal F$. Set $\mathbf M:= \prod_{\lambda} \mathbf M_{\lambda}$. 
We claim that $\mathbf M/\mathcal U\models\Sigma$.
Let $\sigma\in\Sigma$. Then $\{ \sigma\}\in \Lambda$ and
$$F\big(\{\sigma\}\big)\  \subseteq\
\{ \lambda:\ \mathbf M_{\lambda}\models \sigma\}\ =\ \|\sigma\|.$$
Now $F\big(\{\sigma\}\big)\in\mathcal U$, so
$\|\sigma\|\in\mathcal U$, and thus  
$\mathbf M/\mathcal U\models \sigma$ by Theorem~\ref{thm:Los}. \qed

\subsection*{Functoriality of reduced products}
Let $(\mathbf M_\lambda)$ and $(\mathbf N_\lambda)$ be families of $\mathcal L$-structures, and let $h_\lambda\colon \mathbf M_\lambda\to\mathbf N_\lambda$ be a morphism for each $\lambda$. Lemma~\ref{lem:morphisms on products} gives a morphism $h\colon \mathbf M:=\prod_\lambda \mathbf M_\lambda\to \mathbf N:=\prod_\lambda\mathbf N_\lambda$ such that $\pi_\lambda^{\mathbf N}\circ h=h_\lambda\circ \pi_\lambda^{\mathbf M}$ for all $\lambda$.  
Let $\mathcal F$ be a proper filter on $\Lambda$.
Then for all $a,b\in M_s$ we have $a\sim_s b \Rightarrow h_s a \sim_s h_s b$. It easily follows that we have
a morphism $h/{\mathcal F}\colon \mathbf M/{\mathcal F}\to\mathbf N/{\mathcal F}$ making the diagram 
$$\xymatrix{
\mathbf M \ar[d] \ar[r]^h & \mathbf N \ar[d] \\
\mathbf M/{\mathcal F} \ar[r]^{h/{\mathcal F}} & \mathbf N/{\mathcal F}
}$$
commute. Here the vertical arrows are the morphisms $\mathbf M\to\mathbf M/{\mathcal F}$ and  $\mathbf N\to\mathbf N/{\mathcal F}$ defined before Lemma~\ref{lem:Los}.  
If each~$h_\lambda$ is an
embedding, then so is $h/{\mathcal F}$.

We can now prove a statement used in Section~\ref{pafae}. Given $\mathcal L$-structures $\mathbf M$ and~$\mathbf N$, we say that~$\mathbf M$ is
{\bf existentially closed in $\mathbf N$} if $\mathbf M\subseteq\mathbf N$ and every existential $\mathcal L_M$-sentence true in $\mathbf N$ is true in
$\mathbf M$; equivalently, $\mathbf M\subseteq\mathbf N$ and
every universal $\mathcal L_M$-sentence true in $\mathbf M$ is true in $\mathbf N$; notation:
$\mathbf M\preceq_\exists\mathbf N$. Note that if $\mathbf M\preceq_\exists\mathbf N$, then also each universal-existential $\mathcal L_M$-sentence true in $\mathbf N$ is true in
$\mathbf M$.
Clearly $\mathbf M\preceq \mathbf N\Rightarrow \mathbf M\preceq_\exists\mathbf N$.\index{closed!existentially} \index{structure!existentially closed}
\nomenclature[Bj02]{$\mathbf M\preceq_\exists\mathbf N$}{$\mathbf M$ is existentially closed in $\mathbf N$}

\begin{cor}\label{cor:ec in ext} Suppose $\mathbf M$ is the direct union of a directed family $(\mathbf M_\lambda)$ of models $\mathbf M_{\lambda}$ of $\Sigma$. Then $\mathbf M$ is existentially closed in some model of~$\Sigma$.
\end{cor} 
\begin{proof}
Let $\mathcal F$ be the collection of sets $F_{\lambda}:=\{\lambda':\ \lambda'\ge \lambda\}\subseteq \Lambda$. Since $(\Lambda,{\leq})$ is directed, $\mathcal F$ has the FIP. Then Lemma~\ref{lem:extend to uf} yields an ultrafilter 
$\mathcal U\supseteq \mathcal F$ on~$\Lambda$, and so $\mathbf M^*:=(\prod_\lambda \mathbf M_\lambda)/{\mathcal U}\models \Sigma$ by Theorem~\ref{thm:Los}. Let $\iota_\lambda\colon \mathbf M_\lambda\to\mathbf M$ be the natural inclusion, let
$\iota/{\mathcal U}\colon \mathbf M^* \to \mathbf M^\Lambda/{\mathcal U}$ be the embedding
obtained from the family of embeddings~$(\iota_\lambda)$ as described before the corollary, and
let $d\colon \mathbf M\to \mathbf M^\Lambda/{\mathcal U}$ be the diagonal embedding.
Then $d(\mathbf M)\subseteq (\iota/{\mathcal U})(\mathbf M^*)\subseteq \mathbf M^\Lambda/{\mathcal U}$ and $d(\mathbf M)\preceq \mathbf M^\Lambda/{\mathcal U}$, hence
$d(\mathbf M)\preceq_\exists (\iota/{\mathcal U})(\mathbf M^*)$.
\end{proof}

\subsection*{Notes and comments}
Ultrafilters were introduced by H.~Cartan~\cite{Cartan1,Cartan2}, and also appear in Stone~\cite{Stone36}.
The definition of the ultraproduct and Theorem~\ref{thm:Los} are from~\cite{Los55a}, but
versions of ultraproducts had been used already by Skolem~\cite{Skolem31}, Hewitt~\cite{Hewitt},
and Arrow~\cite{Arrow}. Reduced products and the proof of the Compactness Theorem via ultraproducts given above are due to Frayne, Morel and Scott~\cite{FMS}.

\section{Some Uses of Compactness}\label{sec:LS}

\noindent
We first consider  \textit{diagrams,}\/ which provide 
a way to construct embeddings using compactness.
Next we study the relationship between \textit{substructures}\/ and \textit{universal sentences.}\/ We also prove the ``upward'' version of the L\"owenheim-Skolem Theorem.
 
\subsection*{Diagrams}
{\em In this subsection we assume $A\subseteq M$}.
Let
$\operatorname{Diag}(\mathbf M)$ be the set of all quantifier-free sentences $\sigma$ such that  $\mathbf M\models\sigma$, and set $\operatorname{Diag}_A(\mathbf M):=\operatorname{Diag}(\mathbf M_A)$.
We call  $\operatorname{Diag}_M(\mathbf M)$ the (quantifier-free) {\bf diagram} of $\mathbf M$. \index{structure!diagram}\index{diagram}
\nomenclature[Bm4]{$\operatorname{Diag}_M(\mathbf M)$}{diagram of $\mathbf M$}

\begin{lemma}
Let  $h\colon A\to N$, and
$\mathbf N_h$ the expansion of $\mathbf N$ to an $\mathcal L_A$-structure
given by $\underline{a}^{\mathbf N_h}:=h_s(a)$ for $a\in A_s$. 
Then: 
\begin{enumerate}
\item[\textup{(i)}] $h$ preserves quantifier-free formulas iff $\mathbf N_h\models \operatorname{Diag}_A(\mathbf M)$;
\item[\textup{(ii)}] $h$ is elementary iff $\mathbf N_h\models\operatorname{Th}(\mathbf M_A)$.
\end{enumerate}
\end{lemma}

\noindent
This is rather obvious from the definitions, and yields:

\begin{lemma}[Diagram Lemma]\label{lem:diagram}
\mbox{}

\begin{enumerate}
\item[\textup{(i)}] There exists a map $A\to N$ preserving quantifier-free formulas if and only if some $\mathcal L_A$-expansion of
 $\mathbf N$ is a model of $\operatorname{Diag}_A(\mathbf M)$;
\item[\textup{(ii)}] there exists an elementary map $A\to N$ if and only if some $\mathcal L_A$-expansion of
 $\mathbf N$ is elementarily equivalent to $\mathbf M_A$.
\end{enumerate}
\end{lemma}

\noindent
In particular,  there exists an embedding $\mathbf M\to\mathbf N$ iff $\mathbf N$ can be expanded to a model of 
the diagram of $\mathbf M$, and  there exists an elementary embedding $\mathbf M\to\mathbf N$ iff
$\mathbf N$ can be expanded to a model of $\operatorname{Th}(\mathbf M_M)$.
The Diagram Lemma  acquires its power through the Compactness Theorem:

\begin{cor}\label{cor:Henkin}
The following conditions on a structure $\mathbf M$ are equivalent:
\begin{enumerate}
\item[\textup{(i)}] $\mathbf M$ can be embedded into a model of $\Sigma$; 
\item[\textup{(ii)}] every finite subset of $\Sigma\cup\operatorname{Diag}_M(\mathbf M)$ has a model.
\end{enumerate} 
For one-sorted $\mathbf M$, these conditions are also equivalent to: \begin{enumerate}
\item[\textup{(iii)}] every finitely generated substructure of $\mathbf M$ can be embedded into a model of $\Sigma$.
\end{enumerate}
Likewise, $\mathbf M$ has an elementary extension that is a model of $\Sigma$ iff
every finite subset of $\Sigma\cup\operatorname{Th}(\mathbf M_M)$ has a model.
\end{cor}
\begin{proof}
Use the Diagram Lemma and compactness to get (i)~$\Leftrightarrow$~(ii). For one-sor\-ted~$\mathbf M$, use that
every finite subset of $\Sigma\cup\operatorname{Diag}_M(\mathbf M)$ is contained in a set of the form
$\Sigma\cup\operatorname{Diag}_N(\mathbf N)$ where~$\mathbf N$ is a finitely generated
substructure of~$\mathbf M$.
\end{proof}

{\sloppy

\begin{cor}\label{cor:elem equivalence and elem embeddings} Let $\mathbf M$ 
and $\mathbf N$ be given. Then  
$\mathbf M\equiv\mathbf N$ if and only if 
there exists an $\mathcal L$-structure into which both
$\mathbf M$ and~$\mathbf N$ can be ele\-men\-ta\-ri\-ly embedded.
\end{cor}
\begin{proof}
Assume $\mathbf M\equiv\mathbf N$. Extend $\mathcal L$ to
$\mathcal L'$ by adding names for the elements of $\mathbf M$ and $\mathbf N$ such that no name of any $a$ in $\mathbf M$ is the name of any $b$ in $\mathbf N$. We show that the set $\operatorname{Th}(\mathbf M_M)\cup
\operatorname{Th}(\mathbf N_N)$ of $\mathcal L'$-sentences has a model; 
clearly, $\mathbf M$ as well as $\mathbf N$ admits an elementary embedding 
into the $\mathcal L$-reduct of such a model. By compactness
(and replacing a finite subset of $\operatorname{Th}(\mathbf M_M)$ by the 
conjunction of the sentences in it) it suffices that $\big\{\varphi(a)\big\}\cup\operatorname{Th}(\mathbf N_N)$ has a model, for any formula
 $\varphi(x)$ with finite $x$ and $a\in M_x$ such that $\mathbf M \models \varphi(a)$.
For such $\varphi(x)$ and $a$ we have $\mathbf M \models \exists x\varphi(x)$,
so  $\mathbf N \models \exists x\varphi(x)$, and thus $\big\{\varphi(a)\big\}\cup\operatorname{Th}(\mathbf N_N)$ has indeed a model.   

The other direction is obvious.
\end{proof}

}

\noindent
In the same way one shows:

\begin{cor}\label{cor:ex closed and embeddings} Let $\mathbf M \subseteq \mathbf N$. Then $\mathbf M\preceq_\exists\mathbf N$ if and only if $\mathbf N$ embeds over~$\mathbf M$ into some elementary extension of $\mathbf M$.
\end{cor}

\subsection*{Substructures and  universal sentences} 
Call $\Sigma$   {\bf universal} if all sentences in $\Sigma$ are universal. 
\textit{In this subsection we assume that~$\mathcal L$ has for each $s$ a constant symbol of sort~$s$.}\/ Thus
$\<A\>_{\mathbf M}$ is defined for any parameter set~$A$ in~$\mathbf M$.\index{universal!set of sentences}

\begin{prop}\label{prop:Herbrand}
Suppose $\Sigma$ is universal, $x=(x_1,\dots,x_m)$, $y=(y_1,\dots,y_n)$ are disjoint, and $\varphi(x,y)$ is quantifier-free
such that $\Sigma\models\forall x\exists y\, \varphi(x,y)$. Then there are $n$-tuples
$t_1=(t_{11},\dots,t_{1n}),\dots,t_k=(t_{k1},\dots,t_{kn})$  of terms $t_{ij}(x)$, $k\geq 1$,
with
$$\Sigma\models  \varphi\big(x,t_1(x)\big)\vee\cdots\vee\varphi\big(x,t_k(x)\big).$$
\end{prop}

\noindent
The device we use in proving this is often applied:
we extend the language~$\mathcal L$ by new constant symbols and let
them play the role of free parameters. Let $x_1,\dots, x_m$ have sort $s_1,\dots, s_m$, respectively. Let  $\mathcal L_c:=\mathcal L\cup \{c_1,\dots,c_m\}$ where
$c_1,\dots,c_m$ are distinct new constant symbols of sort $s_1,\dots,s_m$, respectively. 
An $\mathcal L_c$-structure $(\mathbf M,a)$ is just an $\mathcal L$-structure $\mathbf M$ together with any $m$-tuple $a=(a_1,\dots,a_m)\in M_{s_1,\dots, s_m}$.
Thus, given any $\mathcal L$-formula
$\psi(x_1,\dots,x_m)$, 
$$\Sigma\models\psi(c_1,\dots,c_m)
\text{ (relative to $\mathcal L_c$)}\ \Longleftrightarrow\ \Sigma\models\psi(x_1,\dots, x_m) \text{ (relative to $\mathcal L$)}.$$ 

\begin{proof}[Proof of Proposition~\ref{prop:Herbrand}]  
Let 
$\mathbf M\models\Sigma$ and $a=(a_1,\dots,a_m)\in M_x$.
Since~$\Sigma$ is universal, $\mathbf N:=\<a_1,\dots,a_m\>_{\mathbf M}$ is also a model
of $\Sigma$. Hence $\mathbf N\models \forall x\exists y\, \varphi$ and so
$\mathbf N\models\exists y\, \varphi(a,y)$. The elements of $\mathbf N$ of sort~$s$
are of the form $t(a)$ where~$t(x)$ is an $\mathcal L$-term of sort~$s$. Hence there
is an $n$-tuple $t=(t_1,\dots,t_n)$ of $\mathcal L$-terms with~$t_i(x)$ of the same sort as~$y_i$, such that $\mathbf N\models\varphi\big(a,t(a)\big)$, and 
thus $\mathbf M\models\varphi\big(a,t(a)\big)$. 
Hence every model of $\Sigma$, viewed as a set of
$\mathcal L_c$-sentences, satisfies some $\mathcal L_c$-sentence
$\varphi\big(c,t(c)\big)$  with $t=(t_1,\dots,t_n)$ an $n$-tuple of $\mathcal L$-terms with $t_i(x)$ of the
same sort as $y_i$ and $c=(c_1,\dots,c_m)$. Now use~\ref{cor:compactness, 2}.
\end{proof}

\noindent
The substructures of fields viewed as $\mathcal L_{\operatorname{R}}$-structures are exactly the integral domains,
that is, the models of the set of universal $\mathcal L_{\operatorname{R}}$-sentences
$$\operatorname{Ri}\cup\big\{ 0\neq 1,\ \forall x\forall y\, (x  y=y  x),\ \forall x\forall y\,
(xy=0 \to x=0 \vee y=0) \big\}.$$
This is an instance of the following general fact:

\begin{prop}\label{prop:LT}
Let 
$$\Sigma_\forall\ :=\ \big\{ \sigma : \text{$\sigma$ is a universal sentence with $\Sigma\models\sigma$} \big\}$$
be the set of universal logical consequences of $\Sigma$. Then
for all $\mathbf M$,
$$\mathbf M\models\Sigma_\forall\ \Longleftrightarrow\ 
\mathbf M \text{ is a substructure of a model of }\Sigma.$$
\end{prop}
\begin{proof} The direction $\Leftarrow$ is clear from Corollary~\ref{cor:persist, 1}(iii). For $\Rightarrow$, 
suppose $\mathbf M\models\Sigma_\forall$.
To show that $\mathbf M$ embeds into a model of $\Sigma$,  let~$\Delta$ be a finite
subset of $\operatorname{Diag}_M(\mathbf M)$; by Corollary~\ref{cor:Henkin} it suffices to show that the set
$\Sigma\cup\Delta$ of $\mathcal L_M$-sentences has a model.
Replacing the sentences in $\Delta$ by their conjunction
we arrange $\Delta=\big\{\varphi(a)\big\}$ where $\varphi(x)$ is a quantifier-free $\mathcal L$-formula,
$x=(x_1,\dots,x_m)$, and $a=(a_1,\dots,a_m)\in M_x$ with distinct $a_1,\dots,a_m$,
such that $\mathbf M\models\varphi(a)$. If $\Sigma\cup\Delta$ has no model, then
$\Sigma\models\neg\varphi(a)$ and hence $\Sigma\models\forall x\,\neg\varphi$, so $\forall x\,\neg\varphi\in \Sigma_\forall$, 
and thus $\mathbf M\models \forall x\,\neg\varphi$, contradicting~${\mathbf M\models\varphi(a)}$. \end{proof}

\begin{cor}\label{cor:LT}
The following conditions on $\Sigma$ are equivalent:
\begin{enumerate}
\item[\textup{(i)}] every substructure of every model of $\Sigma$ is a model of $\Sigma$; 
\item[\textup{(ii)}] $\Sigma$ and $\Sigma_\forall$ are equivalent;
\item[\textup{(iii)}] $\Sigma$ is equivalent to some set of universal sentences.
\end{enumerate}
\end{cor}
%\begin{proof}
%Suppose (i) holds. Clearly every model of $\Sigma$ is also a model of $\Sigma_\forall$;
%conversely, if $\mathbf M\models\Sigma_\forall$, then by the preceding proposition,
%$\mathbf M\subseteq\mathbf N$ for some model~$\mathbf N$ of $\Sigma$, and so $\mathbf M\models\Sigma$ by (i).
%This shows (i)~$\Rightarrow$~(ii);
%the implication (ii)~$\Rightarrow$~(iii) is trivial, and
%(iii)~$\Rightarrow$~(i) is obvious from Corollary~\ref{cor:persist, 1}(iii).
%\end{proof}

\noindent
We say that $\varphi$, $\psi$ are {\bf $\Sigma$-equivalent} if $\Sigma\models\varphi\leftrightarrow\psi$. (So  ``$\emptyset$-equivalent'' is the same as ``equivalent'' in the sense of Section~\ref{sec:formulas}.)
Note that $\Sigma$-equivalence is an equivalence relation on the collection of $\mathcal L$-formulas.
If $\varphi(x)$, $\psi(x)$ are $(\mathcal L,x)$-formulas, then~$\varphi$,~$\psi$ are $\Sigma$-equivalent  iff $\varphi^{\mathbf M}=\psi^{\mathbf M}$ for all $\mathbf M\models\Sigma$. \index{equivalence!formulas}

\begin{cor}\label{cor:equivalent to universal} 
The following are equivalent for $\Sigma,\varphi(x), x=(x_1,\dots, x_m)$:
\begin{enumerate}
\item[\textup{(i)}] for all $\mathbf M\subseteq\mathbf N$ and $a\in M_x$, if $\mathbf M,\mathbf N\models\Sigma$ and
$\mathbf N\models\varphi(a)$, then $\mathbf M\models\varphi(a)$;
\item[\textup{(ii)}] $\varphi(x)$ is $\Sigma$-equivalent to a universal formula $\psi(x)$.
\end{enumerate}
\end{cor}
\begin{proof}
Assume (i), and consider the set $\Sigma':=\Sigma\cup\{\varphi(c)\}$
of $\mathcal L_c$-sentences. By~\ref{prop:LT} and (i)
we have $(\Sigma')_\forall\cup\Sigma\models\varphi(c)$.
By compactness, we can take a universal $\mathcal L$-formula~$\psi(x)$ such that 
$\psi(c)\in (\Sigma')_\forall$ and
$\big\{\psi(c)\big\}\cup\Sigma\models\varphi(c)$. Then $\Sigma\models\varphi\leftrightarrow\psi$.
This shows (i)~$\Rightarrow$~(ii). The reverse implication follows from~\ref{cor:persist, 1}(iii). 
\end{proof}

\subsection*{The L\"owenheim-Skolem Theorem}  {\em In this subsection we
assume that $\mathcal L$ is one-sorted.} (We only use Theorem~\ref{thm:LS}
in this case, and we wish 
to keep formulations simple.) Call an $\mathcal L$-structure {\em infinite\/} 
if its underlying set is infinite.  

First-order logic cannot limit the size of 
an infinite structure: \index{theorem!L\"owenheim-Skolem}

\begin{theorem}[Upward L\"owenheim-Skolem]\label{thm:LS}
Suppose $\mathbf M$ is infinite and
%that one of the underlying sets~$M_s$ of~$\mathbf M$ is infinite.
$\kappa$ is a cardinal $\geq\max\big\{|\mathcal L|,|M|\big\}$. Then~$\mathbf M$ has an elementary extension of size $\kappa$.
\end{theorem}
\begin{proof}
%Take $s$ such that $M_s$ is infinite, 
Let $C$ be a set of new constant symbols 
%of sort $s$
with $|C|=\kappa$, and set $\mathcal L':=\mathcal L\cup C$. Every finite subset of the set
$$\Sigma'\ :=\ \operatorname{Th}(\mathbf M_M)\cup \{c \neq d:\  c,d\in C \text{ are distinct}\}$$
of $\mathcal L'_M$-sentences has a model; in fact, 
%as $M_s$ is infinite, 
as $M$ is infinite, $\mathbf M_M$ itself can be expanded to a model by
suitably interpreting the constants in $C$. By compactness, we obtain a model~$\mathbf N'$ of $\Sigma'$;
by the Diagram Lemma, there exists an elementary embedding of~$\mathbf M$ into the $\mathcal L$-reduct of $\mathbf N'$.
Thus we obtain an elementary extension~$\mathbf N$ of $\mathbf M$ of size~$\geq\kappa$.
Proposition~\ref{prop:LS downwards} gives an elementary substructure $\mathbf N_1$
of~$\mathbf N$ with $M\subseteq N_1$ and $|N_1|=\kappa$. Then $\mathbf M\preceq\mathbf N_1$.
% since both $\mathbf M$ and $\mathbf N_1$ are elementary
%substructures of $\mathbf N$
%and $M\subseteq N_1$.
\end{proof}

\noindent
The downward and upward L\"owenheim-Skolem Theorems (\ref{prop:LS downwards} and \ref{thm:LS}) yield:

\begin{cor}\label{lsupdown}
%If $\Sigma$ has a model, then it has a model of size at most $|\mathcal L|$.
If~$\Sigma$ has an infinite model, then $\Sigma$ has a model of any given size~$\kappa\geq |\mathcal L|$.
\end{cor}

\index{test!Vaught}

\begin{cor}[Vaught's Test]\label{cor:Vaught's Test}
Suppose $|\mathcal L|\leq\kappa$, $\Sigma$ has a model, all models
of~$\Sigma$ are infinite, and all models of $\Sigma$ of size~$\kappa$ are
isomorphic. Then $\Sigma$ is complete.
\end{cor}
\begin{proof} Assume $\sigma$ is true in some model of $\Sigma$.
Then $\sigma$ is true in some model of $\Sigma$ of size $\kappa$,
by Corollary~\ref{lsupdown}. The isomorphism assumption then gives that
$\sigma$ is true in all models of $\Sigma$ of size $\kappa$, and so again by Corollary~\ref{lsupdown}, $\sigma$ is true in all models of~$\Sigma$. It follows that $\Sigma$ is complete.
\end{proof}

\noindent
Vaught's Test is a test for completeness. Here is an application:

\begin{theorem}\label{thm:ACFp}
Let $p$ be a prime number or $p=0$.
The set $\operatorname{ACF}(p)$ of axioms for algebraically closed fields of characteristic $p$ is complete.
\end{theorem}
\begin{proof}
We use basic facts about transcendence bases; see \cite[Chapter~VIII]{Lang}. 
Let~$K$ and~$L$ be algebraically closed fields of characteristic $p$ of the same
uncountable size; let $\k:=\mathbb F_p$ if $p$ is a prime  and $\k:=\Q$ if $p=0$. View $\k$
as a subfield of $K$ and of~$L$ as usual, and note that $|\k|\leq\aleph_0$. Let~$B$ be a transcendence basis of $K$ 
over~$\k$ and~$C$ a transcendence basis of $L$ over $\k$. Then 
$|K|=|\k(B)|=|B|$ and likewise $|L|=|C|$. So there is a bijection $B\to C$, and this bijection extends to a field isomorphism $\k(B)\to \k(C)$ and then further to an isomorphism $K \to L$ between
their algebraic closures.
Thus $\operatorname{ACF}(p)$ is complete by Vaught's Test.
\end{proof}

\noindent
Theorem~\ref{thm:ACFp} and Corollary~\ref{cor:Completeness} imply that  $\operatorname{ACF}(p)$ is decidable. Applications of Theorem~\ref{thm:ACFp} and another proof of this theorem are given in Section~\ref{sec:some theories}.

\subsection*{Notes and comments} Diagrams and their role are explicit in 
A.~Robinson~\cite{Robinson52}. 
Corollary~\ref{cor:Henkin} is in Hen\-kin~\cite{Henkin53}. A sharper version of
Corollary~\ref{cor:elem equivalence and elem embeddings} is due to
Keisler~\cite{Keisler61} and Shelah~\cite{Shelah71a}: if $\mathbf M\equiv\mathbf N$ then there is an ultrafilter $\mathcal U$
on some nonempty set~$\Lambda$ with $\mathbf M^\Lambda/{\mathcal U}\cong \mathbf N^\Lambda/{\mathcal U}$.
Proposition~\ref{prop:Herbrand} is a weak version of a theorem of Herbrand~\cite{Herbrand}.
Corollary~\ref{cor:LT} is due to {\L}o\'{s}~\cite{Los55b} and Tarski~\cite{Tarski54}, and
Corollary~\ref{cor:Vaught's Test} to {\L}o\'{s}~\cite{Los54} and Vaught~\cite{Vaught54}. 

\section{Types and Saturated Structures}\label{sec:sat}

\noindent
Roughly speaking, a \textit{type}\/ is a set of formulas specifying a potential property of a tuple of elements in a structure, similar to a system of 
equations and inequalities that we wish to solve. 
A \textit{saturated structure}\/ realizes many types.

\subsection*{Types}
Let $\Phi=\Phi(x)$ be a set of $(\mathcal L,x)$-formulas. \index{element!realizing a set of formulas}
We say that $a\in M_x$ {\bf realizes~$\Phi$ in~$\mathbf M$} if $\mathbf M\models\varphi(a)$ for all
$\varphi\in\Phi$. Clearly if $\mathbf M\preceq\mathbf N$ and $a\in M_x$, then $a$ realizes~$\Phi$ in~$\mathbf M$ iff~$a$ realizes $\Phi$ in~$\mathbf N$.
We say that {\bf $\Phi$ is realized in $\mathbf M$} if some $a\in M_x$ realizes~$\Phi$ in~$\mathbf M$.
Assume $x=(x_i)_{i\in I}$ and take a tuple $c=(c_i)_{i\in I}$ of distinct new constant symbols
with each $c_i$ of the same sort as $x_i$.
Let $\mathcal L_c$ be $\mathcal L$ augmented by these new constant symbols $c_i$.
Then $\Phi$ is realized in some structure iff the set 
$$\Phi(c)\ :=\ \big\{\varphi(c):\varphi\in\Phi\big\}$$ 
of $\mathcal L_c$-sentences
has a model. Hence by compactness,
$\Phi$ is realized in some structure iff every finite subset of $\Phi$ is realized in some structure.
If this happens we call $\Phi$ an {\bf $x$-type}. \index{type}
An $x$-type $\Phi$ is said to be {\bf complete}\index{type!complete} if for each $\varphi(x)$, either $\varphi \in \Phi$
or $\neg\varphi  \in \Phi$; equivalently, $\Phi(c)$ is a complete $\mathcal{L}_c$-theory. 
For $a \in M_x$ we let 
$\operatorname{tp}^{\mathbf M}_x(a)$ \nomenclature[Bn1]{$\operatorname{tp}^{\mathbf M}_x(b)$}{complete $x$-type in $\mathbf M$ realized by~$b\in M_x$} be  the complete $x$-type in $\mathbf M$ realized by~$a$
 (and we leave out the superscript $\mathbf M$ or subscript $x$ if $\mathbf M$ or~$x$ are clear from the context); that is, for each $\mathcal L$-formula $\varphi(x)$, we have 
$\varphi\in \operatorname{tp}^{\mathbf M}_x(a)$ iff $\mathbf M\models\varphi(a)$.
Every $x$-type is contained in a complete one, namely one of the form $\operatorname{tp}^{\mathbf M}_x(a)$. 
If $\mathbf M\preceq\mathbf N$ and $a\in M_x$, then 
$\operatorname{tp}^{\mathbf M}_x(a)=\operatorname{tp}^{\mathbf N}_x(a)$.

\begin{definition}
We say that $\Phi$ is {\bf $\Sigma$-realizable} if $\Phi$ is realized in some model of~$\Sigma$;
that is, if $\Sigma\cup\Phi(c)$, as a set of $\mathcal L_c$-sentences, has a model. The set of all complete $\Sigma$-realizable $x$-types
is denoted by $\operatorname{S}_x(\Sigma)$. \index{type!$\Sigma$-realizable}
\end{definition}

\noindent
Thus for
a {\em complete\/} $x$-type $\Phi$, we have: 
$\Phi$ is $\Sigma$-realizable iff $\Sigma\subseteq \Phi(c)$.
A complete $x$-type is usually denoted by a letter like $p$ or $q$.

Let $\Sigma_c$ be  $\Sigma$ viewed as set of $\mathcal L_c$-sentences. Then
for $p=p(x)\in \operatorname{S}_x(\Sigma)$ we have $p(c)\in\operatorname{S}(\Sigma_c)$, and 
the map $p\mapsto p(c)\colon \operatorname{S}_x(\Sigma)\to \operatorname{S}(\Sigma_c)$
is a bijection. The {\bf Stone topology} \index{topology!Stone}  on~$\operatorname{S}_x(\Sigma)$ is the topology
making this bijection a homeomorphism. That is, its basic open sets are the sets 
$\<\varphi\>:=\big\{p\in \operatorname{S}_x(\Sigma):\varphi\in p\big\}$, with
$\varphi=\varphi(x)$. Two formulas~$\varphi(x)$ and~$\psi(x)$ are $\Sigma$-equivalent iff $\<\varphi\>=\<\psi\>$.
By Theorem~\ref{thm:compactness, 3}, the Stone topology makes 
$\operatorname{S}_x(\Sigma)$ a compact hausdorff space.

\subsection*{Separating types}
In this subsection we fix a set $\Theta=\Theta(x)$ of $(\mathcal L,x)$-formulas such that $\top, \bot\in \Theta$, and 
for all $\theta_1, \theta_2\in \Theta$, also 
$\theta_1\wedge \theta_2\in \Theta$ and
$\theta_1\vee \theta_2\in \Theta$. For example, $\Theta$ could be the set of
quantifier-free formulas $\theta(x)$.
In Section~\ref{sec:qe} we shall need the next lemma; its corollary was
already used in Section~\ref{cqe1}.

\begin{lemma}\label{lem:separation}
Let a formula $\psi(x)$ be given. Then the following are equivalent:
\begin{enumerate}
\item[\textup{(i)}] $\psi$ is 
$\Sigma$-equivalent to some formula from $\Theta$;
\item[\textup{(ii)}] for all $p,q\in \operatorname{S}_x(\Sigma)$ with
$\psi\in p$ and $\neg\psi\in q$ there exists $\theta\in\Theta$ such that 
$\theta\in p$ and $\neg\theta\in q$.
\end{enumerate}
\end{lemma}
\begin{proof}
The direction (i)~$\Rightarrow$~(ii) is clear. 
Conversely, assume (ii). Consider the open-and-closed subset
$P:=\<\psi\>$ of $\operatorname{S}_x(\Sigma)$ and its complement 
$P^{\c}=\<\neg\psi\>$. 
Let $p\in P$ be given. Then 
$P^{\c}\subseteq\bigcup_{\theta\in\Theta\cap p} \<\neg\theta\>$ by~(ii).
Compactness of $\operatorname{S}_x(\Sigma)$ and $\Theta\cap p$ being closed under conjunction gives
$\theta\in\Theta\cap p$ such that $P^{\c}\subseteq\<\neg\theta\>$ and
hence $p\in\<\theta\>\subseteq P$.
Since $p\in P$ was arbitrary, this yields $P=\bigcup_{\theta\in\Delta} \<\theta\>$ for some $\Delta\subseteq\Theta$.
By compactness of~$\operatorname{S}_x(\Sigma)$ again and $\Theta$ being closed under disjunction,
we obtain $\theta\in\Theta$ with $P=\<\theta\>$, and then $\psi$ is 
$\Sigma$-equivalent to $\theta$. 
\end{proof}

\begin{cor}\label{cor:separation} Suppose $\neg\theta\in \Theta$
for all $\theta\in \Theta$. Then every formula $\psi(x)$ is 
$\Sigma$-equivalent to one in~$\Theta$ iff $p\cap\Theta\neq q\cap\Theta$
for all $p\neq q$ in~$\operatorname{S}_x(\Sigma)$.
\end{cor}

\subsection*{Types over a parameter set}
\textit{In the rest of this section $A$ is a parameter set in~$\mathbf M$.}\/
An {\bf $x$-type over~$A$ in $\mathbf M$}
is a $\operatorname{Th}(\mathbf M_A)$-realizable $x$-type (in the language $\mathcal L_A$).
Equivalently, a set of $(\mathcal L_A,x)$-formulas is an $x$-type over~$A$ in $\mathbf M$
if every finite subset of it is realized in~$\mathbf M_A$.
For~${b \in M_x}$ we let $\operatorname{tp}^{\mathbf M}_x(b|A):=\operatorname{tp}^{\mathbf M_A}_x(b)$ 
be the complete $x$-type over~$A$ in~$\mathbf M$
realized by $b$ (and we leave out the superscript~$\mathbf M$ or subscript $x$ if $\mathbf M$ or~$x$ are clear from the context).
If $\mathbf M\preceq\mathbf N$ and $b\in M_x$, then $\operatorname{tp}^{\mathbf M}_x(b|A)=\operatorname{tp}^{\mathbf N}_x(b|A)$.
Let~$\operatorname{S}^{\mathbf M}_x(A)$ (or~$\operatorname{S}_x(A)$ if $\mathbf M$ is clear from the context)
denote the space $\operatorname{S}_x\big(\!\operatorname{Th}(\mathbf M_A)\big)$ of complete
$x$-types over $A$ in $\mathbf M$. A basis for the Stone topology on $\operatorname{S}_x(A)$ is given by
the sets~$\<\varphi\>=\big\{p\in \operatorname{S}_x(A):\varphi\in p\big\}$ with~$\varphi(x)$ an $\mathcal L_A$-formula. Note that if $\mathcal L$, $x$, $A$
all have size $\le \kappa$, then $\abs{\operatorname{S}_x(A)}\le 2^\kappa$.
 \index{type}\index{type!complete}
\nomenclature[Bn2]{$\operatorname{tp}^{\mathbf M}_x(b\vert A)$}{complete $x$-type in $\mathbf M$ realized by~$b\in M_x$ over $A$}
\nomenclature[Bn3]{$\operatorname{S}^{\mathbf M}_x(A)$}{space of complete $x$-types in $\mathbf M$ over $A$}

\subsection*{Saturated structures}
{\em In this subsection $\kappa$ is a cardinal $>0$}. 
We declare~$\mathbf M$ to be {\bf $\kappa$-saturated} if for all~$A$ of size
$< \kappa$ and every variable $v$ of $\mathcal L$, each complete $v$-type over~$A$ in $\mathbf M$ is realized in
$\mathbf M$; equivalently, for all~$A$ of size
$< \kappa$ and $s\in S$, each collection of $A$-definable subsets of $M_s$ with the finite intersection property has a nonempty intersection.

\begin{example}
Suppose the ordered abelian group $(G;{+},{-},0,{\leq})$ with $G\neq \{0\}$ is $2$-saturated. Then $[G]$ has no largest element:
given $a\in G^{>}$ a realization $b$ of the $v$-type $\{v>na:n=0,1,2,\dots\}$
over $\{a\}$ yields $[b]>[a]$. In particular, the ordered abelian group $(\R;{+},{-},0,{\leq})$ is not $2$-saturated.
\end{example}

{\sloppy
\noindent
If $M_s$ is finite for all $s$, then
$\mathbf M$ is $\kappa$-saturated for all $\kappa$.
If $\mathbf M$ is $\kappa$-saturated and~$M_s$ is infinite, then $|M_s|\ge \kappa$. (Take a variable $v$ of sort $s$ and
consider the $v$-type ${\{v\neq a:\ a\in M_s\}}$.) If $\kappa\le \kappa'$ and
$\mathbf M$ is $\kappa'$-saturated, then $\mathbf M$ is $\kappa$-saturated.
If~$\mathbf{M}$ is $\kappa$-saturated, then so is any reduct of $\mathbf{M}$.
If $\mathbf M$ is $\kappa$-saturated, $\kappa$ is infinite, and 
$\abs{A}<\kappa$, then $\mathbf M_A$ is $\kappa$-saturated.
}

The definition of ``$\kappa$-saturated'' ostensibly only concerns
families of subsets of $M_s$ for $s\in S$. It is a pleasant feature of model theory, however, that a one-variable property often yields a many-variable analogue, with some effort as in this case: \index{structure!$\kappa$-saturated} \index{saturation}

\begin{lemma}\label{lem:sat}
Suppose $\mathbf M$ is $\kappa$-saturated, $\kappa$ is infinite,
$A$ has size $<\kappa$ and $x$ has
size~${\leq \kappa}$. Then every $x$-type over $A$ in $\mathbf M$ is realized in $\mathbf M$.
\end{lemma}
\begin{proof} Let $x=(x_i)_{i\in I}$, and let $i$, $j$ range over $I$. Let 
$p\in \operatorname{S}^{\mathbf M}_x(A)$;
it suffices to show that $p$ is realized in $\mathbf M$.
Take a well-ordering $\leq$ of $I$ of order type $\le |I|$.
Then each proper downward closed subset of $I$ has cardinality~$<|I|$.
For each $j$ let $x_{\leq j}:=(x_i)_{i\leq j}$ and let $p_{\leq j}$
be the set of all formulas in $p$ with free variables in~$x_{\leq j}$;
then $p_{\leq j}$ is a complete $x_{\leq j}$-type over $A$ in $\mathbf M$.
Similarly, for each $j\in I$ we define $x_{<j}$ and the complete $x_{<j}$-type $p_{<j}$ over $A$ in $\mathbf M$;
then for each $\varphi\in p_{\leq j}$ we have $\exists x_j\varphi\in p_{<j}$.
By recursion on $i$ we construct a point $(b_i)\in M_x$ such that for each~$j$, $(b_i)_{i\leq j}$ realizes $p_{\leq j}$.
Suppose that for a certain $j$ we have already a point
$b=(b_i)_{i<j}\in M_{x_{<j}}$ realizing $p_{<j}$.
Then $p_j:=\big\{\varphi(b/x_{<j}):\varphi\in p_{\leq j}\big\}$ is an $x_j$-type
over $A\cup\{b_i:i<j\}$ in $\mathbf M$. Since $A\cup\{b_i:i<j\}$ has size $<\kappa$ and
$\mathbf M$ is $\kappa$-saturated, we have $b_j\in M_{x_j}$ realizing $p_j$.
Then $(b_i)_{i\leq j}$ realizes~$p_{\leq j}$.
\end{proof}

\noindent
In $\kappa$-saturated structures one can do things that otherwise would require passing to an elementary extension.
As an example, here is a variant of Corollary~\ref{cor:elem equivalence and elem embeddings}:

\begin{cor}\label{cor:sat elem embedd}
Suppose that $\mathbf M\equiv\mathbf N$ and $\mathbf N$ is $\kappa$-saturated for some infinite~${\kappa\ge\abs{M}}$. Then there exists an elementary embedding
$\mathbf M\to\mathbf N$.
\end{cor}
\begin{proof} Assume for simplicity that our language $\mathcal{L}$ is one-sorted. Take a multivariable $x=(x_a)_{a\in M}$, and let $\vec a:=(a)_{a\in M}$. Since $\mathbf M\equiv\mathbf N$, 
$\Phi:=\big\{\varphi(x):\mathbf M\models\varphi(\vec a)\big\}$
is an $x$-type over~$\emptyset$ in $\mathbf N$.
By Lemma~\ref{lem:sat}, $\Phi$ is realized in $\mathbf N$.
If $(a')_{a\in M}$ realizes~$\Phi$ in~$\mathbf N$, then $a\mapsto a'$ ($a\in M$) is an elementary
embedding $\mathbf M\to\mathbf N$.
\end{proof}

\noindent
Every structure has a $\kappa$-saturated elementary extension:

\begin{prop}\label{prop:saturated elemext}
Suppose $|\mathcal L|\le \kappa$ and $|M|\le 2^\kappa$. Then $\mathbf M$ has a $\kappa^+$-saturated elementary
extension $\mathbf N$  with $|N| \leq 2^\kappa$.
\end{prop}

\noindent
We first establish the following lemma, with the same assumptions on $\kappa$
as in Proposition~\ref{prop:saturated elemext}.

\begin{lemma}
There exists $\mathbf M'\succeq\mathbf M$ with
$|M'| \leq 2^\kappa$ such that for all $A\subseteq M$
with $|A|\leq\kappa$ and any variable $v$ of $\mathcal L$,  each type in~$\operatorname{S}_v^{\mathbf M}(A)$ is realized in $\mathbf M'$.  
\end{lemma}
\begin{proof} There are at most $2^\kappa$ many pairs $(A,v)$ with $A\subseteq M$,
$\abs{A}\leq\kappa$, and $v$ a quantifiable variable of~$\mathcal L$.
For such $(A,v)$ we have 
$\abs{\operatorname{S}_v^{\mathbf M}(A)}\le 2^\kappa$;
take for each 
$p\in \operatorname{S}_v^{\mathbf M}(A)$ a new constant symbol $c_p$ 
of the same sort as $v$. Let
$\mathcal L'$ be $\mathcal L$ augmented by these $c_p$. By compactness and Proposition~\ref{prop:LS downwards} applied to the
set 
$$\operatorname{Th}(\mathbf M_M)\cup\bigcup\big\{p(c_p):\ p\in \operatorname{S}_v^{\mathbf M}(A),\ A\subseteq M,\ \abs{A}\leq\kappa,\ v \text{ is quantifiable} \big\}$$
of $\mathcal L'_{M}$-sentences, we obtain an elementary extension $\mathbf M'$ of $\mathbf M$ as desired.
\end{proof}

\noindent
Below we use the fact that the least cardinal $\kappa^{+}$ that is greater than $\kappa$ is a regular ordinal as defined in Section~\ref{sec:ordered sets}; see \cite[Theorem 3.1.11]{Harzheim}.

\begin{proof}[Proof of Proposition~\ref{prop:saturated elemext}] 
Let $\alpha$, $\beta$ range over ordinals $<\kappa^+$. Recursion on~$\alpha$ yields a sequence $(\mathbf M_\alpha)$ of
structures of size $\leq 2^\kappa$ such that:
\begin{enumerate}
\item $\mathbf M_0=\mathbf M$ and $\mathbf M_\alpha\preceq\mathbf M_\beta$ if $\alpha\le\beta$;
\item if $\alpha=\beta+1$ is a successor ordinal, then
$\mathbf M_{\beta}$ is obtained from $\mathbf M_\alpha$ as $\mathbf M'$ is from $\mathbf M$
in the previous lemma;
\item if $\alpha>0$ is a limit ordinal, then $\mathbf M_\alpha=\bigcup_{\beta<\alpha}\mathbf M_\beta$. 
\end{enumerate}
Then $\mathbf N:=\bigcup_{\alpha} \mathbf M_\alpha$  is an elementary extension
of $\mathbf M$ with $|N|\leq 2^\kappa$. Suppose $B\subseteq N$ 
and $\abs{B} < \kappa^+$. Since $\kappa^+$ is regular, we have
$\beta$ such that $B\subseteq M_\beta$. Then every  $p\in\operatorname{S}_v^{\mathbf N}(B)=\operatorname{S}_v^{\mathbf M_\beta}(B)$ is realized in
$\mathbf M_{\beta+1}$ and thus in~$\mathbf N$.
\end{proof}

{\sloppy 
\subsection*{Notes and comments}
The notion of $\kappa$-saturated structure goes back to the $\eta_\alpha$-sets of
Hausdorff~\cite[p.~181]{Hausdorff14}, but only appeared in model theory in the late 1950s
in the work of Morley and Vaught~\cite{Morley,MorleyVaught,Vaught61}, who also introduced types.
See~\cite{Felgner02} for the history. A motivation for
Hausdorff was du~Bois-Reymond's ``infinitary pantachy''~\cite{dBR82} (the partially ordered set of germs at $+\infty$ of continuous real-valued
functions): see~\cite{Ehrlich, Steprans}. 

}

\section{Model Completeness}\label{sec:mc}

\noindent 
If the primitives of an
$\mathcal{L}$-structure $\mathbf M$ are computationally or topologically well-be\-haved, then the 
sets defined in~$\mathbf M$ by quantifier-free formulas are
often fairly tame as well.  For example, the
subsets of $\R$ definable by quantifier-free formulas in the
ordered field of real numbers (viewed as an $\mathcal L_{\operatorname{OR}}$-structure) are finite unions of open intervals 
and singletons. In the ordered ring of integers,
the subsets of $\Z$ definable by quantifier-free formulas are just the traces of the preceding sets in $\Z$. But taking arbitrary $\mathcal L_{\operatorname{OR}}$-formulas, one can define much more complicated sets in the ordered ring of integers: for example, the formula $\pi(x)$ given by
$$1+1\leq x \ \& \ \forall u\forall v\big((1\le u\ \&\ 1\le v\ \&\ u\cdot v = x)\ \rightarrow\ u= 1 \ \vee \ v = 1 )$$
defines the set of prime numbers. In general, the more quantifiers that occur in a formula $\varphi(x)$, the more complicated the set $\varphi^{\mathbf M}$ defined by $\varphi(x)$ in $\mathbf M$ can be. When this typical behavior does {\em not\/} occur, it is worth noting!

This is why in this section we consider
{\em model completeness}. In the next section we study a sharper version of this, called
{\em quantifier elimination}. 
A structure $\mathbf M$ is said to be {\bf model complete} if
for each formula~$\varphi(x)$ there exists an existential formula~$\psi(x)$ such that $\varphi^{\mathbf M}=\psi^{\mathbf M}$. For one-sorted model complete $\mathbf M$, every definable subset of~$M^m$ is, for some $n$, the image of a quantifier-free definable
set in $M^{m+n}$ under the natural projection map $M^{m+n}\to M^m$. Such projections often preserve many desirable 
topological properties, for example, having only finitely many connected components (in a suitable topological environment).

\medskip
\noindent
Model completeness also makes sense for a set of sentences:

\begin{definition}\index{theory!model complete} \index{model completeness} $\Sigma$ is said to be {\bf model complete} if every formula  is $\Sigma$-equi\-va\-lent to an existential formula. (Note that then every formula $\varphi(x)$ is actually $\Sigma$-equivalent to an $\exists$-formula $\varphi'(x)$ with no more free variables than those of $\varphi$.) 
\end{definition}

\noindent
Thus $\mathbf M$ is model complete iff  $\operatorname{Th}(\mathbf M)$ is model complete. The following lemma slightly reduces the job of proving model completeness:\index{structure!model complete}

\begin{lemma} Suppose every universal formula is $\Sigma$-equivalent to an existential formula. Then $\Sigma$ is model complete.  \end{lemma}
\begin{proof}
We show by induction that 
every formula $\varphi$ is $\Sigma$-equivalent to
an $\exists$-formula. This is clear if $\varphi$ is quantifier-free,
and the conclusion is preserved under~$\wedge$,~$\vee$, and~$\exists$.
Suppose $\varphi=\neg\psi$. Assuming inductively that 
$\psi$ is $\Sigma$-equivalent to an $\exists$-formula,
 $\varphi$ is $\Sigma$-equivalent to a $\forall$-formula, and thus by hypothesis, $\varphi$ is  $\Sigma$-equivalent to an $\exists$-formula.
The case where $\varphi=\forall y\psi$ with a single quantifiable variable $y$ follows, since then $\varphi$
is equivalent to $\neg\exists y\neg\psi$.
\end{proof}

\noindent
Here is Robinson's model completeness test:\index{model completeness!test}\index{test!model completeness}

\begin{prop} \label{thm:Robinson test} The following are equivalent:
\begin{enumerate}
\item[\textup{(i)}] $\Sigma$ is model complete;
\item[\textup{(ii)}] for all models $\mathbf M$, $\mathbf N$ of $\Sigma$, if $\mathbf M\subseteq\mathbf N$, then 
$\mathbf M\preceq\mathbf N$;
\item[\textup{(iii)}] for all models $\mathbf M$, $\mathbf N$ of $\Sigma$, if $\mathbf M\subseteq\mathbf N$, then 
$\mathbf M\preceq_\exists\mathbf N$.
\end{enumerate}
\end{prop}
\begin{proof}
Clearly (i)~$\Rightarrow$~(ii)~$\Rightarrow$~(iii).
Assume (iii), and let~$\varphi(x)$ be universal. 
Then for $\mathbf M,\mathbf N\models\Sigma$ with $\mathbf M\subseteq\mathbf N$ and $a\in M_x$ we have $\mathbf M\models\varphi(a)\Rightarrow \mathbf N\models\varphi(a)$.
Hence $\varphi$ is $\Sigma$-equivalent to
an $\exists$-formula, by Corollary~\ref{cor:equivalent to universal}. 
Thus $\Sigma$ is model complete by the previous lemma. 
\end{proof}

\noindent
Condition (ii) in Proposition~\ref{thm:Robinson test} means that
$\Sigma\cup\operatorname{Diag}_M(\mathbf M)$ is complete, for all~$\mathbf M\models\Sigma$;
this explains the terminology \textit{model complete.}\/ 
Model completeness often entails completeness:
A {\bf prime model} of $\Sigma$ is a model of $\Sigma$ that embeds elementarily into every model of $\Sigma$. Note that if $\Sigma$ has a prime model, then $\Sigma$ is complete.
If~$\Sigma$ is model complete and $\mathbf M$ is a model of $\Sigma$ that embeds into every model
of $\Sigma$, then $\mathbf M$ is a prime model of $\Sigma$.\index{model!prime}\index{prime model}\index{theory!prime model}

\medskip
\noindent
Here is a variant of the above test for model completeness:  

\begin{cor}\label{cor:modelcomplete test}
The following are equivalent:
\begin{enumerate}
\item[\textup{(i)}] $\Sigma$ is model complete;
\item[\textup{(ii)}] for all models $\mathbf M$, $\mathbf N$ of $\Sigma$ with $\mathbf M\subseteq\mathbf N$ and every elementary extension $\mathbf M^*$ of $\mathbf M$ that is $\kappa$-saturated for some $\kappa>\abs{N}$, there is an embedding $\mathbf N\to\mathbf M^*$
that extends the natural inclusion $\mathbf M\to\mathbf M^*$.
\end{enumerate}
\end{cor}
\begin{proof}
Suppose $\Sigma$ is model complete, and $\mathbf M$, $\mathbf N$, $\mathbf M^*$ are as in the hypothesis of~(ii).
Then $\Sigma\cup\operatorname{Diag}_M(\mathbf M)$ is complete, with models
$\mathbf M^*_M$, $\mathbf N_M$. Since $\mathbf M^*_M$ is $|N|$-saturated,
Corollary~\ref{cor:sat elem embedd} yields an embedding $\mathbf N_M\to\mathbf M^*_M$.
This shows (i)~$\Rightarrow$~(ii).
The implication (ii)~$\Rightarrow$~(i) follows from Corollary~\ref{cor:ex closed and embeddings}, Proposition~\ref{prop:saturated elemext} and the equivalence of (i) and (iii) in Theorem~\ref{thm:Robinson test}.
\end{proof}

\begin{exampleNumbered}\label{ex:TfDiv, 1}
Let $\mathcal L:=\mathcal L_{\operatorname{A}}$ be the language of  additive groups and
$$\Sigma\ :=\ \operatorname{Tf}\,\cup\,\operatorname{Div}\,\cup\,\operatorname{Ab}\,\cup\,\big\{\exists x\,(x\neq 0)\big\}.$$ 
Then the models of $\Sigma$ are the  nontrivial divisible torsion-free abelian groups: see Example~\ref{ex:theories}(2),(3). 
A model of $\Sigma$ may be viewed as a $\Q$-linear space in the natural way,
and if the model is $\kappa$-saturated, then it has dimension
$\geq\kappa$. Together with Corollary~\ref{cor:modelcomplete test} this yields model completeness of $\Sigma$.
\end{exampleNumbered}

\noindent
Before trying to show that a certain theory is model complete, 
one better check that it has an axiomatization by
$\forall\exists$-sentences:

\begin{prop}\label{prop:forallexists}
Suppose $\Sigma$ is model complete. Then $\Sigma$ is equivalent to a set of $\forall\exists$-sentences.
\end{prop}
\begin{proof}
Let $\Sigma_\forall$ be as in Proposition~\ref{prop:LT} (so the models of $\Sigma_\forall$ are the substructures of models of $\Sigma$). 
Every sentence of the form $\forall x(\varphi\to\psi)$, where $\varphi(x)$ is universal and~$\psi(x)$ is existential, is equivalent to a $\forall\exists$-sentence.
Let $\Sigma_{\forall\exists}$ be the set of all $\forall\exists$-sentences  
equivalent to  a sentence $\forall x(\varphi\to\psi)$, where $\varphi(x)$ is universal, $\psi(x)$ is existential, and $\varphi$, $\psi$  are $\Sigma$-equivalent. We claim that $\Sigma$ and $\Sigma':=\Sigma_\forall\cup \Sigma_{\forall\exists}$ are equivalent.
Note first that $\Sigma'$ is also model complete: Let  $\varphi(x)$ be universal, and
take  an $\exists$-formula $\psi(x)$ which is $\Sigma$-equivalent to $\varphi$; then 
$\varphi$ is also $\Sigma'$-equivalent to $\psi$, since  up to equivalence
$\forall x(\psi\to\varphi)$ lies in $\Sigma_\forall$ and
$\forall x(\varphi\to\psi)$  in~$\Sigma_{\forall\exists}$. Let $\mathbf M\models\Sigma'$; we need to show $\mathbf M\models \Sigma$. Now $\mathbf M\models \Sigma_\forall$ yields  some $\mathbf N\models \Sigma$ with $\mathbf M\subseteq\mathbf N$. 
Then $\mathbf N\models\Sigma'$ and so $\mathbf M\preceq\mathbf N$ by model completeness of~$\Sigma'$, thus $\mathbf M\models\Sigma$.
\end{proof}

\noindent
Call $\Sigma$ {\bf inductive} if the
direct union of any directed family
of models of $\Sigma$ is a model of $\Sigma$. If $\Sigma$ is a set of $\forall\exists$-sentences, then $\Sigma$ is inductive. Hence by \ref{prop:forallexists}:\index{inductive}\index{theory!inductive}
 
\begin{cor}\label{cor:model complete=>inductive}
If $\Sigma$ is model complete, then $\Sigma$ is inductive.
\end{cor}

\subsection*{Existentially closed models} 
An {\bf existentially closed model of $\Sigma$} is a
mo\-del~$\mathbf M$ of $\Sigma$ that is existentially closed
in every extension $\mathbf N\models\Sigma$ of $\mathbf M$.
Thus by Robinson's test, $\Sigma$ is model complete iff every model of $\Sigma$ is an existentially closed model of $\Sigma$.
Inductive theories have existentially closed models: \index{model!existentially closed} \index{closed!existentially}

\begin{lemma}\label{lem:ec models}
Suppose $\Sigma$ is inductive and $\mathbf M\models\Sigma$. Then $\mathbf M$ extends to an existentially closed model of $\Sigma$.
\end{lemma}
\begin{proof}
We first show that $\mathbf M$ has an extension $\mathbf M^*\models\Sigma$ which 
satisfies every existential $\mathcal L_M$-sentence that holds in some 
extension of $\mathbf M^*$ to a model of $\Sigma$.

Let $(\sigma_\lambda)_{\lambda<\kappa}$ be an enumeration of all $\mathcal L_M$-sentences,
where $\kappa$ is the cardinality of the set of $\mathcal L_M$-sentences.
By recursion on $\mu<\kappa$ we construct a sequence $(\mathbf M_\mu)_{\mu\leq\kappa}$ of
extensions of $\mathbf M$ to models of $\Sigma$ 
such that $\mathbf M_\lambda\subseteq\mathbf M_{\lambda'}$ for all $\lambda\leq\lambda'\leq\kappa$
as follows:
Set $\mathbf M_0:=\mathbf M$. Suppose $0<\mu\leq\kappa$ and
$(\mathbf M_\lambda)_{\lambda<\mu}$ is a sequence
of extensions of  $\mathbf M$ to models of $\Sigma$ 
such that $\mathbf M_\lambda\subseteq\mathbf M_{\lambda'}$ for all $\lambda\leq\lambda'<\mu$. If $\mu$ is a limit ordinal, then set $\mathbf M_\mu:=\bigcup_{\lambda<\mu}\mathbf M_\lambda$; since $\Sigma$ is inductive, we have $\mathbf M_\mu\models\Sigma$.
Suppose $\mu=\mu'+1$; if some
extension of $\mathbf M_{\mu'}$ is a model of $\Sigma\cup\{\sigma_{\mu'}\}$, then we let $\mathbf M_{\mu}$ be such an extension, and otherwise $\mathbf M_\mu:=\mathbf M_{\mu'}$.
Finally, put $\mathbf M^*:=\mathbf M_\kappa$. It is routine to
check that $\mathbf M^*$ has the desired property.

Now inductively define $\mathbf M^0:=\mathbf M$ and $\mathbf M^{n+1}:=(\mathbf M^n)^*$. Then
$\bigcup_n \mathbf M^n$ is an existentially closed model of $\Sigma$ which extends $\mathbf M$.
\end{proof}

\subsection*{Model companions and model completions}
\textit{In this subsection we assume that $T$ is an inductive 
$\mathcal{L}$-theory $T$.}\/
Even though $T$ may not be model complete, we may hope for the existence of a theory whose
models are exactly the existentially closed models of $T$. This suggests the following notion:

{\sloppy
\begin{definition}\index{model companion}\index{theory!model companion}\index{companion!model} A {\bf model companion} of $T$ is a
model complete $\mathcal{L}$-theo\-ry~${T^*\supseteq T}$ such that every model of $T$ embeds into
a model of $T^*$.
\end{definition}
}

\begin{lemma} Let $T^*$ be an $\mathcal{L}$-theory. Then
$T^*$ is a model companion of $T$ iff the models of $T^*$ are exactly the  existentially closed models of $T$. 
\end{lemma}
\begin{proof}
Suppose $T^*$ is a model companion of $T$. Then every model $\mathbf M^*$ of $T^*$ 
is an existentially closed model of $T$: suppose
$\mathbf N\models T$ extends $\mathbf M^*\models T^*$; extend~$\mathbf N$ to a model~$\mathbf N^*$ of $T^*$; then $\mathbf M^*\preceq\mathbf N^*$ by model completeness of $T^*$ and thus~${\mathbf M^*\preceq_\exists\mathbf N}$ by Corollary~\ref{cor:persist, 1}(ii). Conversely, suppose $\mathbf M$ is an
existentially closed model of~$T$. Take $\mathbf M^*\models T^*$ with $\mathbf M\subseteq\mathbf M^*$;
then $\mathbf M\preceq_\exists \mathbf M^*$ and thus $\mathbf M\models T^*$ by
Proposition~\ref{prop:forallexists}.
This shows the ``only if'' direction of the lemma; the ``if'' direction follows from Theorem~\ref{thm:Robinson test} and Lemma~\ref{lem:ec models}.\end{proof}

\begin{cor}
There is at most one model companion of $T$.
\end{cor}

\noindent
Thus if a model companion of $T$ exists, we may speak of \textit{the}\/ model companion of~$T$.
Clearly $T$ is model complete iff it is its own model companion.

\begin{exampleNumbered}\label{ex:TfDiv, 2}
The $\mathcal L_{\operatorname{A}}$-theory $T$ of torsion-free abelian groups is inductive. We let~$T^*$ be the 
$\mathcal L_{\operatorname{A}}$-theory of nontrivial divisible torsion-free abelian groups. 
Then $T^*$ is model complete (Example~\ref{ex:TfDiv, 1}), and every $A\models T$
extends to some $A^*\models T^*$: if $A=\{0\}$, take $A^*:=\Q$, and if $A\neq\{0\}$, let
$A^*:=A\otimes_\Z\Q$ be the divisible hull of~$A$. Hence $T^*$ is the model companion of $T$.
\end{exampleNumbered}

\begin{definition}
A {\bf model completion} of $T$ is a model companion $T^*$ of $T$ such that in addition
$T^*\cup\operatorname{Diag}_A(\mathbf A)$ is complete, for all  
$\mathbf A\models T$.
%If it exists, we may speak of {\it the}\/ model completion of $T$. 
\index{model completion} \index{theory!model completion}\index{completion!model}
\end{definition}

\noindent
We say that {\bf $\Sigma$ has AP\/} (short for: {\bf $\Sigma$ has
the amalgamation property}) if for all models~$\mathbf A$,~$\mathbf M_1$,~$\mathbf M_2$ of~$\Sigma$ and embeddings
$f_i\colon \mathbf A\to \mathbf M_i$ ($i=1,2$) there exist a model~$\mathbf N$ of~$\Sigma$ and embeddings
$g_i\colon \mathbf M_i\to\mathbf N$ ($i=1,2$) making the diagram
$$\xymatrix{ 
& \mathbf N&  \\
\mathbf M_1 \ar@{-->}[ru]^{g_1} & & \mathbf M_2 \ar@{-->}[lu]_{g_2}   \\
& \mathbf A \ar@{->}[lu]^{f_1} \ar@{->}[ru]_{f_2} & }$$
commute. Using the definitions and Corollary~\ref{cor:elem equivalence and elem embeddings}
one easily shows: \index{amalgamation property} \index{theory!amalgamation property}
 
\begin{lemma}\label{lem:amalg prop}
Let $T^*$ be an $\mathcal{L}$-theory. Then
$T^*$ is a model completion of $T$ iff~$T^*$ is a model companion of~$T$ and $T$ has \textup{AP}.
\end{lemma}

\begin{exampleNumbered}
Let $T$, $T^*$ be as in Example~\ref{ex:TfDiv, 2}.
Then $T$ has AP. To see why, let embeddings $f_i\colon A\to M_i$ ($i=1,2$) of torsion-free abelian groups be given. Set
$$F\ :=\ \big\{(b,c)\in M_1\times M_2:\ nb=f_1(a),\ nc=-f_2(a) \text{ for some $n\ge 1$ and $a\in A$}\big\},$$ 
a subgroup of the direct sum $M_1\oplus M_2$ of
$M_1$ and $M_2$. Put $N:=(M_1\oplus M_2)/F$; it is easy to check that $N$ is torsion-free and that the natural morphisms $g_i\colon M_i\to N$ ($i=1,2$) are embeddings.
Hence $T^*$ is the model completion of $T$.
\end{exampleNumbered}

\subsection*{Notes and comments}
The notion of model completeness and Proposition~\ref{thm:Robinson test} are due to A.~Robinson~
\cite{Robinson53,Robinson55c,Robinson56,Robinson65}. 
See~\cite{Macintyre77} for a survey of
classical model completeness results.
In connection with~\ref{prop:forallexists} and~\ref{cor:model complete=>inductive} we should mention that  
$\Sigma$ is inductive iff $\Sigma$ is equivalent to 
a set of $\forall\exists$-sentences;
this is due to Chang~\cite{Chang59} and {\L}o\'{s}-Suszko~\cite{Los-Suszko}; see \cite[Theorem~6.5.9]{Hodges} for a proof. Existentially closed models were introduced by Rabin~\cite{Rabin62}, 
model companions by Barwise and Robinson~\cite{BarwiseRobinson}, and the connection between the two
notions is from~\cite{EklofSabbagh}.
By Cherlin~\cite{Cherlin73}, the $\mathcal L_{\operatorname{R}}$-theory of
commutative rings has no model companion.

\section{Quantifier Elimination}\label{sec:qe}

\noindent  
To keep formulations simple \textit{we assume in this section that our language~$\mathcal L$ has for each~$s$ a constant symbol of sort $s$.}\/ 
Quantifier elimination is often the first key step to understanding the category of definable sets of a structure.\index{structure!quantifier elimination}\index{quantifier elimination}\index{theory!quantifier elimination}

\begin{definition}
We say that {\bf $\Sigma$ admits quantifier elimination} (QE) if  each formula $\varphi(x)$ with finite $x$ is $\Sigma$-equivalent to a quantifier-free formula $\varphi'(x)$.
We also express this by saying that \textit{$\Sigma$ has \textup{QE}}\/ or \textit{$\Sigma$ eliminates quantifiers.}\/
We say that $\mathbf M$ admits~QE if~$\operatorname{Th}(\mathbf M)$ admits QE.
\end{definition}

\noindent
(Of course, the restriction to finite $x$ in this definition is superfluous.)
Note that if $\Sigma$ has QE, then $\Sigma$ is model complete.  
In this section we give criteria for~$\Sigma$ to have QE. 
For a given structure $\mathbf M$ of interest, one tries to find a~$\Sigma$ that has~$\mathbf M$ as a model and passes such a test.
By Corollary~\ref{cor:separation}, $\Sigma$ admits QE iff for every multivariable $x$, each type
in $\operatorname{S}_x(\Sigma)$ is determined by its subset of quantifier-free formulas; this fact underlies the tests.
To verify that $\Sigma$ admits QE, it suffices to show that  formulas of a   special
form are $\Sigma$-equivalent to quantifier-free formulas:

\begin{lemma}\label{lem:QE simplified}
Suppose each formula $\exists y\, \theta(x,y)$ with
quantifier-free $\theta$, finite $x$, and a single variable 
$y$ is $\Sigma$-equivalent to a quantifier-free formula.
Then $\Sigma$ has \textup{QE}.
\end{lemma}
\begin{proof} We proceed by induction on the
number of logical symbols in a formula.
To show that any formula $\varphi(x)$ is $\Sigma$-equivalent to a quantifier-free formula $\varphi'(x)$, we first note that
this is obvious for quantifier-free $\varphi$, and that this property is
preserved when taking disjunctions, conjunctions and negations.
Suppose $\varphi=\exists y\psi$. We can assume inductively that
$\psi(x,y)$ is $\Sigma$-equivalent to a quantifier-free formula~$\psi'(x,y)$.
The hypothesis of the lemma then yields a quantifier-free formula $\varphi'(x)$ that is
$\Sigma$-equivalent to $\exists y\psi'(x,y)$, and hence to $\varphi(x)$.
The case $\varphi=\forall y\psi$ reduces to the previous cases, since 
$\varphi$ is equivalent to $\neg\exists y\neg\psi$.
\end{proof}

\noindent
Before trying to show that a certain theory has QE, it might be a good idea to
verify the following necessary condition, analogous to 
Proposition~\ref{prop:forallexists}:

\begin{cor} Suppose $\Sigma$ has \textup{QE}. Then $\Sigma$ is equivalent to
a set of sentences of the form $\forall x \exists y \theta(x,y)$ with finite 
$x$, a single variable $y$, and quantifier-free $\theta$.
\end{cor}
\begin{proof} Let $\Sigma'$ be the set of sentences 
$\forall x \exists y \theta(x,y)$ of the indicated form that
are logical consequences of $\Sigma$. An easy exercise using 
Lemma~\ref{lem:QE simplified} shows
that $\Sigma'$ has~QE. Let any sentence $\sigma$ be given. Then we have a 
quantifier-free sentence 
$\sigma'$ with $\Sigma'\models \sigma\leftrightarrow \sigma'$, and thus
$\Sigma\models \sigma\leftrightarrow \sigma'$.  Assume also
that $\sigma\in \Sigma$. Then $\Sigma\models \sigma'$, hence
$\Sigma'\models \sigma'$ (as $\sigma'$ is quantifier-free) and thus $\Sigma'\models \sigma$. 
\end{proof}

\subsection*{A general elimination theorem}
Suppose $\mathcal L$ is a sublanguage of
the language~${\mathcal L}^*$. Let $S^*\supseteq S$ be the set of sorts for $\mathcal{L}^*$; assume that $\mathcal{L}^*$ has for each $s^*\in S^*$ a constant symbol of sort $s^*$. 
Let ${\mathbf M}^*$ and ${\mathbf N}^*$
range over ${\mathcal L}^*$-structures, and let~${\mathbf M}$ and~${\mathbf N}$
be their $\mathcal L$-reducts. Let $\Sigma^*$ be a set of ${\mathcal L}^*$-sentences.   
Here is a criterion for an $\mathcal L^*$-formula to be $\Sigma^*$-equivalent to a quantifier-free $\mathcal L$-formula:

{\sloppy
\begin{lemma}\label{lem:key lemma} Let $x$ be a multivariable of $\mathcal{L}$
and $\varphi^*(x)$ an 
${\mathcal L}^*$-for\-mu\-la. Then~$\varphi^*(x)$ is $\Sigma^*$-equivalent to a
quantifier-free ${\mathcal L}$-formula~$\varphi(x)$ iff for all~${{\mathbf M}^*,{\mathbf N}^*\models \Sigma^*}$, every
common $\mathcal L$-sub\-struc\-ture ${\mathbf A}=(A;\dots)$ of~${\mathbf M}$ and~${\mathbf N}$, and all~$a\in A_x$:
\begin{equation}\label{eq:key lemma}
{\mathbf M}^*\models \varphi^*(a) \quad\Longleftrightarrow\quad {\mathbf N}^*\models \varphi^*(a).
\end{equation}
\end{lemma}

}
\begin{proof}
Forward direction: use Corollary~\ref{cor:persist, 1}(i), and $\mathbf M^*\models\sigma$ iff $\mathbf M\models\sigma$, when~$\sigma$ is an $\mathcal L_M$-sentence.
Conversely, suppose \eqref{eq:key lemma} holds for all $\mathbf M^*$, $\mathbf N^*$, $\mathbf A$ and~$a$ as specified in the lemma. 
Let $p,q\in\operatorname{S}_x(\Sigma^*)$ contain the same quantifier-free $\mathcal L$-formulas.
By Lemma~\ref{lem:separation} it is enough to show that then $\varphi^*\in p$ iff $\varphi^*\in q$.
Take $\mathbf M^*,\mathbf N^*\models\Sigma^*$ and $a\in M_x$ realizing $p$ in $\mathbf M^*$ and $b\in N_x$ realizing $q$ in $\mathbf N^*$. Put $\mathbf A:=\<a\>_{\mathbf M}$
and $\mathbf B:=\<b\>_{\mathbf N}$. Since $a$ and $b$ realize the same quantifier-free $\mathcal L$-formulas
in $\mathbf M$ and $\mathbf N$, respectively, we get an isomorphism $\mathbf A\to\mathbf B$ sending $a$ to $b$. Then ${\mathbf M}^*\models \varphi^*(a)$ iff ${\mathbf N}^*\models \varphi^*(b)$, by
\eqref{eq:key lemma}, so $\varphi^*\in p$ iff $\varphi^*\in q$.
\end{proof}

\noindent
Typically,  for $\mathcal L^*\neq\mathcal L$, the criterion in the lemma above gets used via its corollary below. To state that corollary, we
define $\Sigma^*$ to have
{\bf closures of $\mathcal L$-substructures} if 
for all ${\mathbf M}^*,{\mathbf N}^*\models \Sigma^*$ with a common 
${\mathcal L}$-substructure $\mathbf A=(A;\dots)$ of
${\mathbf M}$ and~${\mathbf N}$, there is a (necessarily unique) 
isomorphism over $A$ from the $\mathcal L^*$-substructure of ${\mathbf M}^*$ generated by~$A$ onto the $\mathcal L^*$-substructure of ${\mathbf N}^*$ generated by $A$.  \index{closure!$\mathcal L$-substructures}

\begin{cor}\label{closureprop}
If $\Sigma^*$ has closures of $\mathcal L$-substructures, then every
quan\-ti\-fier-free ${\mathcal L}^*$-formula is $\Sigma^*$-equivalent to a quantifier-free
$\mathcal L$-formula.
\end{cor}

{\sloppy
\begin{example}
Let $\mathcal L=\mathcal L_{\operatorname{R}}$ be the language of rings
and let $\mathcal L^*$ be $\mathcal L$ augmented by a unary function symbol $\iota$.
Any field is naturally an $\mathcal{L}$-structure and
expands to an $\mathcal L^*$-structure by $\iota(0):=0$ and $\iota(a):= a^{-1}$ for nonzero $a$ in the field. Note that
$$\Sigma^*\ :=\ \operatorname{Fl}\,\cup\, \big\{ \iota(0)=0,\ \forall x \big( x\neq 0 \rightarrow x\cdot \iota(x)=1 \big) \big\}$$
has closures of $\mathcal L$-substructures, by the universal property of the fraction field of an integral domain. Hence every quantifier-free $\mathcal L^*$-formula is
$\Sigma^*$-equivalent to a quantifier-free $\mathcal L$-formula. (This also follows easily
by a simple direct argument; Corol\-lar\-ies~\ref{lem:L vs L'} and \ref{lem:L vs L^a} are more substantial applications.)
\end{example}
}

\subsection*{Application to QE} The expression 
$\mathbf{M}\equiv_A \mathbf{N}$ is short for $\mathbf M_A\equiv\mathbf N_A$; 
here $A$ is a parameter set in both $\mathbf{M}$ and
$\mathbf{N}$, that is, $A\subseteq M$ and $A\subseteq N$. Typically we use
this notation when $A$ underlies a common substructure 
$\mathbf{A}$ of
$\mathbf{M}$ and $\mathbf{N}$.
Taking $\mathcal L^*=\mathcal L$ in Lemma~\ref{lem:key lemma} yields:

\nomenclature[Bj4]{$\mathbf M \equiv_A\mathbf N$}{$\mathbf M_A$, 
$\mathbf N_A$ are elementarily equivalent}

\begin{cor}\label{cor:substructure completeness} Given $\Sigma$, the following conditions are equivalent: \begin{enumerate}
\item[\textup{(i)}] \text{$\Sigma$ has \textup{QE}}; 
\item[\textup{(ii)}] \text{$\mathbf M\equiv_A\mathbf N$,
for all $\mathbf M,\mathbf N\models\Sigma$ and any common substructure $\mathbf A$ of $\mathbf M$, $\mathbf N$};
\item[\textup{(iii)}] \text{$\Sigma\cup\operatorname{Diag}_A(\mathbf A)$ is complete,
for all $\mathbf A\models\Sigma_\forall$.}
\end{enumerate}
\end{cor}

\begin{cor}\label{cor:QE and sentences}
Suppose $\Sigma$ admits \textup{QE} and has a model. Then 
 $\Sigma$ is complete iff some $\mathcal{L}$-structure  \textup{(}not necessarily a model of $\Sigma$\textup{)} embeds into every model of $\Sigma$. 
\end{cor}

\subsection*{A test for quantifier elimination} We now come to our first QE test: \index{test!QE}
\index{quantifier elimination!test}

\begin{prop} \label{prop:QE test} Let $\Sigma$ be given. Then the following are equivalent: \begin{enumerate}
\item[\textup{(i)}] $\Sigma$ has \textup{QE};
\item[\textup{(ii)}] whenever $\mathbf A$ is a substructure of a model $\mathbf M$ of $\Sigma$, any
embedding of $\mathbf A$ into any model $\mathbf N$ of $\Sigma$ extends to an embedding of $\mathbf M$ into some $\mathbf N^*\succeq\mathbf N$.
$$\xymatrix{ 
& \mathbf N^* \\
\mathbf M \ar@{-->}[ur] & \mathbf N  \ar@{-->}[u]_{\preceq} \\
\mathbf A \ar[u]^{\subseteq} \ar[ur] }$$
\end{enumerate}
\end{prop}
\begin{proof}
Suppose $\Sigma$ has QE,  let $\mathbf M,\mathbf N\models\Sigma$,
$\mathbf A\subseteq\mathbf M$, and let $i\colon\mathbf A\to\mathbf N$ be an embedding;
we need to show that $i$ extends to an embedding of $\mathbf M$ into some elementary extension of $\mathbf N$. Replacing $\mathbf N$ by an isomorphic structure we can assume that~${\mathbf A\subseteq\mathbf N}$ and
$i$ is the natural inclusion $\mathbf A\to\mathbf N$. Then by Corollary~\ref{cor:substructure completeness} we have
$\mathbf M\equiv_A\mathbf N$. Now use
Corollary~\ref{cor:elem equivalence and elem embeddings} to obtain an extension of $i$ to an (elementary) embedding of~$\mathbf M$ into some elementary extension of $\mathbf N$.

For the converse, suppose (ii) holds.  
Let $\varphi(x)$ be an existential formula; by Lem\-ma~\ref{lem:QE simplified}
it suffices to show that $\varphi(x)$ is $\Sigma$-equivalent to a quantifier-free for\-mu\-la~$\varphi'(x)$. For this we use Lemma~\ref{lem:key lemma} in the case $\mathcal L^*=\mathcal L$.
So let $\mathbf M,\mathbf N\models\Sigma$, let~$\mathbf A$ be a common substructure of $\mathbf M$ and $\mathbf N$,
and $a\in A_x$; by symmetry it is enough to show: $\mathbf M\models\varphi(a)\Rightarrow \mathbf N\models\varphi(a)$. Extend the natural inclusion $\mathbf A\to\mathbf N$ to an embedding of~$\mathbf M$ into an elementary extension $\mathbf N^*$ of $\mathbf N$. Then
$$\mathbf M\models\varphi(a)\quad\Rightarrow\quad \mathbf N^*\models\varphi(a)
\quad\Rightarrow\quad \mathbf N\models\varphi(a),$$ 
where we used Corollary~\ref{cor:persist, 1}(ii) for the first implication.
\end{proof}

\subsection*{QE tests using saturation}
Using saturated structures yields variants of the above test that are often easier to use:

\begin{cor}[QE test, first variant]\label{prop:QE test, 10} Let $\Sigma$ be 
given. Suppose that for all $\mathbf{M},\mathbf{N}\models \Sigma$ with $\abs{M}^{+}$-saturated $\mathbf{N}$
and all substructures~$\mathbf A$ of $\mathbf M$, every
embedding $\mathbf{A} \to \mathbf{N}$ extends to an embedding $\mathbf{M}\to\mathbf{N}$. Then $\Sigma$ has \textup{QE}.
\end{cor}
\begin{proof} Any structure has a $\kappa$-saturated elementary extension,
for any $\kappa>0$, so condition (ii) of Proposition~\ref{prop:QE test} is satisfied.
\end{proof}

\noindent
The freedom to choose $b$ in the next variant is sometimes convenient.

\begin{cor}[QE test, second variant]\label{prop:QE test, 3}\index{test!QE}\index{quantifier elimination!test}
Let $\Sigma$ be given. Suppose that for every $\mathbf{M}\models \Sigma$ 
and substructure~$\mathbf A$ of $\mathbf M$ with
$A\neq M$ and every embedding $i$ of~$\mathbf A$ into an $\abs{A}^+$-saturated model $\mathbf N$ of $\Sigma$ there exist $s\in S$ and $b\in M_s\setminus A_s$
such that $i$ extends to an embedding $\mathbf A\<b\>\to \mathbf N$. Then $\Sigma$ has \textup{QE}.
\end{cor}
\begin{proof} By Zorn the hypothesis of Corollary~\ref{prop:QE test, 10} is satisfied.   
\end{proof}

\noindent
A {\bf $\Sigma$-closure} of an $\mathcal{L}$-structure $\mathbf A$ is a model $\mathbf A_\Sigma$ of $\Sigma$ with
$\mathbf A\subseteq\mathbf A_\Sigma$ such that every embedding of $\mathbf A$ into a model $\mathbf M$ of $\Sigma$ extends (not necessarily uniquely) to an embedding $\mathbf A_\Sigma\to\mathbf M$. \index{closure!$\Sigma$-closure} 
Existence of $\Sigma$-closures simplifies the above test:

\begin{cor}[QE test, third variant]\label{prop:QE test, 2} \index{test!QE}\index{quantifier elimination!test} Let
$\Sigma$ be given. Assume:
\begin{enumerate}
\item[\textup{(i)}] every substructure of every  model of $\Sigma$ has a $\Sigma$-closure; and
\item[\textup{(ii)}] for all models $\mathbf M$, $\mathbf N$ of $\Sigma$ with $\mathbf M\subseteq\mathbf N$ and
$M\neq N$,  and any $\abs{M}^+$-saturated model $\mathbf M^*$ extending $\mathbf M$
there exist $s\in S$, $b\in N_s\setminus M_s$, and an embedding of $\mathbf M\<b\>$ into $\mathbf M^*$ over $\mathbf M$.
$$\xymatrix{ 
\mathbf N&  \\
\mathbf M\<b\> \ar[u]^{\subseteq} \ar@{-->}[r] & \mathbf M^*   \\
\mathbf M \ar@{-}[u]^{\subseteq} \ar@{-}[ur]_\subseteq }$$
\end{enumerate}
Then $\Sigma$ has \textup{QE}.
\end{cor}
\begin{proof} Check that the hypothesis of \ref{prop:QE test, 3} is satisfied.
\end{proof}

\begin{exampleNumbered}\label{ex:DOAb}
Suppose $\mathcal L=\mathcal L_{\operatorname{OA}}$ is the language of ordered abelian groups and
$\Sigma=\operatorname{Div}\cup\operatorname{OAb}\cup\{\exists x\,x\neq 0\}$, so the models of $\Sigma$ are the  nontrivial divisible ordered abelian groups; see Example~\ref{ex:theories}(3),(5). 
Substructures of models of~$\Sigma$
are ordered abelian groups. Each ordered abelian group $A$ has a $\Sigma$-closure $A_\Sigma$: 
if $A=\{0\}$,  let~$A_\Sigma$ be the ordered additive group~$\Q$; 
if $A\neq\{0\}$, let  $A_\Sigma:=\Q\otimes_\Z A$ be the ordered divisible hull  of~$A$; see Section~\ref{sec:oag}. 

Let $M,N\models\Sigma$ with $M\subseteq N$, let $M^*$ be an $|M|^+$-saturated model of $\Sigma$ extending $M$, and let $a\in N\setminus M$. Take an element $a^*\in M^*$ realizing the same cut in~$M$ as~$a$. We have
$M\<a\>=M+\Z a=M\oplus\Z a$ and $M+\Z a^*=M\oplus\Z a^*$ (internal direct sums), so we get an embedding $M\<a\>\to M^*$ over $M$ sending $a$ to $a^*$. 
Thus~$\Sigma$ has QE by~\ref{prop:QE test, 2}.
\end{exampleNumbered}

\subsection*{QE and definable closure}
\textit{In this subsection we assume that $\Sigma$ is universal, and has \textup{QE}, and $\mathbf M\models\Sigma$.}\/
This has several useful consequences:

\begin{cor} Let $A\subseteq M$. Denoting the underlying set of the substructure~$\<A\>_{\mathbf M}$ of $\mathbf M$ by $\<A\>$,
we have $\<A\>= \operatorname{dcl}(A)$. 
\end{cor}
\begin{proof} This is because $\<A\>\subseteq \operatorname{dcl}(A)$, and 
$\<A\>_{\mathbf M}\preceq \mathbf M$. 
\end{proof}    

\noindent
The next result is known as: $\mathbf M$ has definable Skolem functions.

\begin{cor} Let $x=(x_1,\dots, x_m)$ and $y=(y_1,\dots, y_n)$
be disjoint finite multivariables, and suppose the relation $Z\subseteq M_{x}\times M_{y}$ is $0$-definable in $\mathbf M$.
Let $\pi\colon M_x\times M_y \to M_x$ be the natural projection map.
Then there is a $0$-definable map $f\colon M_x\to M_y$
such that $\big(a,f(a)\big)\in Z$ for all $a\in \pi(Z)$.
\end{cor}
\begin{proof} Replacing $Z$ by $Z\cup \big((M_x\setminus \pi(Z))\times M_y\big)$ we arrange that $\pi(Z)=M_x$. Take a quantifier-free formula $\varphi(x,y)$ that defines $Z$ in~$\mathbf M$. Set 
$$\Sigma(\mathbf M)\ :=\ \Sigma\cup\{\sigma:\   \text{$\sigma$ is quantifier-free and $\mathbf{M}\models \sigma$}\}.$$
Then $\Sigma(\mathbf M)$ is complete, so
$\Sigma(\mathbf M)\models\forall x\exists y\, \varphi(x,y)$. Now use Proposition~\ref{prop:Herbrand}.
\end{proof}

\noindent
The proof above also shows that
definable functions are piecewise given by terms:

\begin{cor} \label{lem:def fns piecewise as terms}
Suppose $f \colon X \to M_s$ with $X\subseteq M_x$ is $0$-definable. Then there are
terms $t_1(x),\dots,t_k(x)$ of sort $s$, $k\in \N^{\ge 1}$, such that for every $a\in X$
we have $f(a)=t_i^{\mathbf M}(a)$ for some $i\in\{1,\dots,k\}$.
\end{cor}

\subsection*{QE and model completeness} Let $T$ be an $\mathcal{L}$-theory. By Corollary~\ref{cor:substructure completeness}, 
$T$ has~\textup{QE} iff
it is the model completion of~$T_\forall$. Thus:

\begin{lemma}\label{lem:model completion with QE} If $T$ has a
universal axiomatization and has a model completion, then this model completion has \textup{QE}.
\end{lemma}

\begin{example}
Let $\mathcal L=\mathcal L_{\operatorname{OA}}$
and $\Sigma$ be as in Example~\ref{ex:DOAb}: the models of $\Sigma$ are the nontrivial divisible 
ordered abelian groups viewed as $\mathcal L$-structures.
Then~$\operatorname{Th}(\Sigma)$ is the model completion of the theory
$\operatorname{Th}(\operatorname{OAb})$ of ordered abelian groups. 
\end{example}

\begin{cor}\label{cor:ES}
$T$ has \textup{QE} $\Longleftrightarrow$
$T$ is model complete and $T_\forall$ has \textup{AP}.
\end{cor}

\subsection*{Notes and comments}
QE goes back to Skolem~\cite{Skolem19,Skolem20}, Langford~\cite{Langford1,Langford2}, Presburger~\cite{Presburger}, and Tars\-ki~\cite{Tarski31,Tarski35b,Ta}. 
Co\-rol\-la\-ry~\ref{cor:substructure completeness} is from Sacks~\cite{Sacks}.
Corollary~\ref{prop:QE test, 10} is close to 
Shoenfield~\cite{Shoenfield71,Shoenfield77} and to~\cite[Theorem 17.2]{Sacks},
the latter attributed there to L.~Blum. Example~\ref{ex:DOAb} is related to ``Fourier-Motzkin elimination''~\cite[III, pp.~209--223]{Schrijver}.
See \cite{DMM0} for an application of Lemma~\ref{lem:def fns piecewise as terms}.
Lemma~\ref{lem:model completion with QE} is from~\cite{Robinson58} and
Corollary~\ref{cor:ES} from~\cite{EklofSabbagh}.
A survey of QE results from before 1984 is \cite{Weispfenning84}.
 
\section{Application to Algebraically Closed and Real Closed Fields}\label{sec:some theories}

\noindent
To demonstrate how the material of the previous sections is used in practice, we now apply it to \textit{algebraically closed fields}\/ and
\textit{real closed fields.}\/  

\subsection*{Algebraically closed fields}
We establish here the basic model-theoretic facts about algebraically closed
fields: eli\-mi\-na\-tion of quantifiers, the Nullstellensatz, strong minimality, definably closed~=~perfect subfield, and definable functions.

In this subsection $K$, $E$, $F$ are fields, and 
\textit{algebraic over}\/, \textit{algebraically closed}\/, and \textit{algebraic closure}\/ are
taken in the sense of field theory; these notions will turn out to agree for
algebraically closed ambient fields with the model-theoretic notions. By
$\operatorname{ACF}$ we mean here the set of axioms for algebraically closed fields in the language
of rings, and for $p$ a prime number or $p = 0$ we let~$\operatorname{ACF}(p)$ be the set of axioms for
algebraically closed fields of characteristic~$p$, as in Examples~\ref{ex:theories}(10),(11). Below~$p$ ranges over
the set $\{0, 2, 3, 5, \dots \}$ of possible characteristics.

\begin{theorem}\label{thm:ACF QE}
$\operatorname{ACF}$ has \textup{QE.}
\end{theorem}
\begin{proof}
Every integral domain has an $\operatorname{ACF}$-closure, namely the algebraic closure of its fraction field. 
Let $E$  be algebraically closed, let $K$ be  a proper 
algebraically closed subfield of $E$, and let $F$ be algebraically closed, $\abs{K}^+$-saturated, also with~$K$ as a subfield.
Take 
any $a \in E \setminus K$. Then $P(a) \neq 0$ for all  $P(T) \in K[T]^{\neq}$,
that is,~$a$ is transcendental over $K$. By saturation we can take $b \in F\setminus K$, so $b$ is transcendental over~$K$. Then the natural inclusion $K\to F$ extends to an embedding
$K[a] \to F$ sending~$a$ to~$b$. Thus $\operatorname{ACF}$ has QE by Corollary~\ref{prop:QE test, 2}.
\end{proof}

\noindent
The substructures of algebraically closed fields are exactly the integral domains, so
by the above, $\operatorname{Th}(\operatorname{ACF})$ is the model completion of the theory of integral domains. 
Here is a well-known manifestation (but not a special case) of Theorem~\ref{thm:ACF QE}:

\begin{example}
Let $X_1,\dots,X_n,Y$ be distinct indeterminates, $X=(X_1,\dots,X_n)$, and
$P,Q\in \Z[X,Y]^{\neq}$ be monic in $Y$. Then there is a polynomial $R\in\Z[X]$ (the \textit{resultant}\/ of $P$ and~$Q$)
such that for every algebraically closed $E$ and $a\in E^n$: $R(a)=0$ iff $P(a,Y),Q(a,Y)\in E[Y]$
have a common zero in $E$; see \cite[IV,~\S{}8]{Lang}.
\end{example}

\noindent
The following consequences of QE make up the Constructibility Theorem
(Che\-val\-ley-Tarski) and the Nullstellensatz (Hilbert).
Let $T=(T_1,\dots,T_n)$ be an $n$-tuple of distinct indeterminates. \index{theorem!Chevalley-Tarski}\index{theorem!Nullstellensatz}

\begin{cor}\label{cor:ACF QE}
Let $E$ be algebraically closed with subfield $K$.
\begin{enumerate}
\item[\textup{(i)}] A subset of $E^n$ is $K$-definable in $E$ if and only if it is a boolean
combination of zero sets
$\big\{a \in  E^n : P(a) = 0\big\}$ of polynomials $P \in  K[T]$.
\item[\textup{(ii)}] If $P_1, \dots , P_k, Q_1, \dots,Q_l \in  K[T]$ and there is an overfield $F$ of
$K$ with a point $a \in  F^n$ such that
$$P_1(a) = \cdots = P_k(a) = 0,\ Q_1(a) \neq 0, \dots , Q_l(a) \neq 0,$$
then there is such a point $a \in  E^n$.
\item[\textup{(iii)}] Let $P_1, \dots, P_m \in  K[T]$. Then $P_1, \dots , P_m$ have no common
zero in $E$ iff there are $Q_1, \dots , Q_m \in  K[T]$ such
that $P_1Q_1 + \cdots + P_mQ_m = 1$   \textup{(}in $K[T]$\textup{)}.
\end{enumerate}
\end{cor}
\begin{proof}
Item (i) is immediate from QE. In (ii), extend $F$ to be algebraically
closed, and use that then $E \equiv_K F$  by QE. In (iii), suppose there
are no $Q_1,\dots, Q_m\in K[T]$ such that $P_1Q_1 + \cdots + P_mQ_m = 1$. Then
the ideal of $K[T]$ generated by $P_1,\dots,P_m$ is a proper ideal, and thus
contained in a maximal ideal $\m$ of $K[T]$. 
Put $t_i := T_i + \m$ for $i = 1,\dots, n$ and $t:=(t_1,\dots,t_n)$.
Then $\m \cap K = \{0\}$, so
$K[T]/\m = K[t]$
is a field extension of $K$, and $P(t) = P(T) + \m$ for each
$P\in K[T]$; in particular, $P_1,\dots, P_m$ have $t$ as a common
zero in an extension field of~$K$, and thus $P_1,\dots, P_m$ have a common zero in
$E$ by (ii).
\end{proof}

\noindent
The following is known as the {\em strong minimality\/} of algebraically closed fields. (This property alone has already many consequences.) 

\begin{cor}\label{cor:ACF sm} Let $E$ be algebraically closed. Then a set $X\subseteq E$ is definable in~$E$ iff $X$ is finite or cofinite.
\end{cor}
\begin{proof} For a single indeterminate $T$ and $P\in E[T]^{\ne}$ the set $\big\{a\in E:\ P(a)=0\big\}$ is finite. Now use Corollary~\ref{cor:ACF QE}(i) for $n=1$.
\end{proof}

\noindent
In Section~\ref{sec:LS} we showed that 
$\operatorname{ACF}(p)$ is complete. This is also a consequence of  Corollary~\ref{cor:QE and sentences} and Theorem~\ref{thm:ACF QE}, since for $p > 0$ the field $\mathbb F_p$ embeds into every
model of $\operatorname{ACF}(p)$, and the ring $\Z$ embeds into every model of $\operatorname{ACF}(0)$.
Thus an $\mathcal L_{\operatorname{R}}$-sentence holds in \textit{some}\/
algebraically closed field of characteristic $p$ iff it holds in \textit{every}\/ 
algebraically closed field of characteristic $p$. Moreover:

\begin{cor}\label{cor:Lefschetz}
Let $\sigma$ be an $\mathcal L_{\operatorname{R}}$-sentence.  Then the following are equivalent:
\begin{enumerate}
\item[\textup{(i)}] $\sigma$ holds in some algebraically closed field of characteristic zero;
\item[\textup{(ii)}] for all but finitely many prime numbers $p$, $\sigma$ holds in all algebraically closed fields of characteristic $p$.
\end{enumerate}
\end{cor}
\begin{proof}
The direction (i)~$\Rightarrow$~(ii) is clear from the remark
preceding the corollary.
% Compactness Theorem~\ref{thm:compactness, 2}.
Suppose conversely that $n_0\in\N$ is such that  $\sigma$ holds in all algebraically closed fields of characteristic $p>n_0$. To get (i), it suffices 
to show that
the set
$\Sigma := \{\sigma\} \cup \operatorname{ACF}(0)$ 
of $\mathcal L_{\operatorname{R}}$-sentences has a model.
Every finite  $\Sigma_0\subseteq\Sigma$  contains only finitely many  sentences $p_1 1\neq 0,\dots,p_n 1\neq 0$
from $\operatorname{Fl}(0)\subseteq \operatorname{ACF}(0)$; since~$\sigma$ holds in every algebraically closed
field of characteristic $p>\max\{n_0,p_1,\dots,p_n\}$, every such algebraically closed field 
is a model of $\Sigma_0$. Hence $\Sigma$ has a model by compactness. This shows (ii)~$\Rightarrow$~(i).
\end{proof}

\noindent
Combining Example~\ref{ex:forallexists Fpalg} and Corollary~\ref{cor:Lefschetz} yields:

\begin{cor}\label{cor:forallexists FpAlg}
Let $\sigma$ be a universal-existential $\mathcal L_{\operatorname{R}}$-sentence, and
suppose that for all but finitely many prime numbers $p$, $\sigma$ holds in all finite fields of characteristic~$p$.
Then $\sigma$ holds in every algebraically closed field of characteristic zero.
\end{cor}

\begin{exampleNumbered}[Serre]\label{ex:Serre}
Suppose $E$ is algebraically closed and $\operatorname{char}(E)\neq 2$.
For $P= (P_1,\dots,P_n )\in E[T]^n$ and $a\in E^n$  set
$P(a):=\big(P_1(a),\dots,P_n(a)\big)\in E^n$.
Then \textit{for any $P\in E[T]^n$ such that $P\big(P(a)\big)=a$ for all $a\in E^n$, there exists an $a\in E^n$ with $P(a)=a$.}\/
To see this note that for fixed $n$ and a bound $d$ on the total degree of the~$P_i$,
the claim can be expressed as a universal-existential $\mathcal L_{\operatorname{R}}$-sentence~$\sigma_{n,d}$. 
Obviously $\sigma_{n,d}$  holds in all finite fields of characteristic $> 2$, and so it 
holds in all algebraically closed
fields of characteristic~$0$ by Corollary~\ref{cor:forallexists FpAlg}, and  in all
algebraically closed
fields of characteristic~$> 2$ by   Example~\ref{ex:forallexists Fpalg} and   completeness of~$\operatorname{ACF}(p)$,~$p>2$.
\end{exampleNumbered}

\noindent
If $E$ is algebraically closed
with subset $A$, then the (model-theoretic) algebraic clo\-sure~$\operatorname{acl}(A)$ of $A$ in $E$
contains obviously the field-theoretic algebraic closure in~$E$ of the subfield
of $E$ generated by $A$, and is in fact equal to (the underlying set of) this
field-theoretic algebraic closure:

\begin{lemma} Let $E$ be algebraically closed with algebraically closed sub\-field~$K$. Then $K$ is algebraically closed in $E$ in the model theory sense.
\end{lemma}
\begin{proof}
Obvious from $K\preceq E$.
\end{proof}

\noindent
If $E$ is algebraically closed and $A\subseteq E$, then the
definable closure $\operatorname{dcl}(A)$ of $A$ in~$E$ contains at least the subfield of $E$ generated by $A$, and equals (the underlying set of) this subfield when $E$ has
characteristic zero:

\begin{prop} \label{prop:ACF(0) dcl} Let $K$ be a subfield of
the algebraically closed field $E$ of characteristic zero. Then $K$ is definably closed in $E$.
\end{prop}
\begin{proof}
Let $a \in E\setminus K$; we claim that then $\sigma a \neq a$ for some $\sigma \in\Aut(E|K)$. (By Lemma~\ref{lem:dcl, acl} the proposition follows from this claim.) If $a$ is transcendental
over~$K$, take a transcendence basis $B$ of $E$ over $K$ with $a \in B$, take the
automorphism of~$K(B)$ over $K$ that sends each $b \in B$ to $b + 1$, and then
extend it to an automorphism of the algebraic closure $E$
of $K(B)$. Suppose $a$ is algebraic over~$K$. Since $a \notin K$, the minimum
polynomial of $a$ over $K$ is of degree $> 1$, so has a zero $b \in E$ with $b \neq a$
(here we use that $E$ has characteristic zero). Take an automorphism $\sigma$ of the
algebraic closure $K^\alg$ of $K$ in $E$ over $K$ that sends $a$ to~$b$, take a transcendence
basis $B$ of $E$ over~$K^\alg$, and extend $\sigma$ to the automorphism of~$K^\alg(B)$ that is
the identity on $B$, and then extend further to an automorphism of $E$.
\end{proof}

\noindent
Characterizations of definable closures of this type lead to corresponding
descriptions of definable functions. In this case definable functions are piecewise rational:

\begin{cor}\label{LR1}
Let $E$ be algebraically closed of characteristic zero, with
subfield $K$, and let $X \subseteq E^n$ and $f \colon X \to E$ be $K$-definable in $E$. Then
there are $P_1, \dots , P_k, Q_1, \dots , Q_k \in K[T]$ such that for each $x \in X$
there is $i \in \{1, \dots , k\}$ with $Q_i(x) \neq 0$ and $f(x) = P_i(x)/Q_i(x)$.
\end{cor}
\begin{proof}
Extending $E$ if necessary we can assume $E$ is $|K|^+$-saturated. Let~${x \in X}$. Then $f(x) \in \operatorname{dcl}(K \cup \{x\}) = K(x)$ by the proposition above, that
is, $f(x) = P(x)/Q(x)$ with polynomials $P, Q \in K[T]$, $Q(x) \neq 0$. Now
use saturation.
\end{proof}

\noindent
Suppose $E$ is algebraically closed of characteristic $p > 0$. Then we have a
$0$-definable automorphism $x \mapsto x^p$ of $E$, the \textit{Frobenius map,}\/ and the inverse
$y \mapsto y^{1/p}$
of this map is not given piecewise by rational functions. The $n$th iterate
$x \mapsto x^{p^n}$
of the Frobenius map has inverse $y \mapsto y^{1/p^n}$, and as we shall
see, these inverse maps are the only obstructions in getting an analogue in
positive characteristic of the above. Recall that a field $K$ of characteristic
$p > 0$ is said to be \textit{perfect}\/ if every element of~$K$ is a $p$th power $x^p$ of some
$x \in K$. So every finite field is perfect. For any subfield $K$ of $E$
there is a smallest perfect subfield of $E$ containing $K$, namely
$$K^{1/p^\infty} := \bigcup_n K^{1/p^n} \quad\text{ where $K^{1/p^n}:= \big\{x^{1/p^n}: x \in K\big\} \subseteq E$,}$$
and by the next result $K^{1/p^\infty}$
is the definable closure of $K$ in $E$.

\begin{prop} Let $E$ be algebraically closed of characteristic $p > 0$, with
perfect subfield $K$. Then $K$ is definably closed in $E$.
\end{prop}

\noindent
The proof is identical to that of Proposition~\ref{prop:ACF(0) dcl}, using the fact that an
irreducible polynomial in one variable over a perfect field is separable.

\begin{cor}\label{cor:def fn ACF(p)}
Let $E$ be algebraically closed of characteristic $p > 0$, with
perfect subfield $K$, and let $X \subseteq E^n$ and $f \colon X \to E$ be $K$-definable in
$E$. Then there are  $P_1, \dots , P_k, Q_1, \dots , Q_k \in K[T]$ and an $m$ with the
following property: 
for each $x \in X$   there is $i \in \{1, \dots , k\}$ such that
$$Q_i\big(x^{1/p^m}\big) \neq 0\ \text{ and }\ f(x) = P_i\big(x^{1/p^m}\big)\big/Q_i\big(x^{1/p^m}\big),$$
where $x^{1/p^m} := \big(x_1^{1/p^m},\dots,x_n^{1/p^m}\big)$ for $x = (x_1, \dots, x_n) \in E^n$.
\end{cor}

\noindent
Here is an application, usually stated only for injective endomorphisms of
algebraic varieties as a theorem of Ax: \index{theorem!Ax}

\begin{cor} \label{cor:Ax}
Let $E$ be algebraically closed, and suppose $X \subseteq E^n$ and
$f \colon X \to X$ are definable in $E$ and $f$ is injective. Then $f$ is surjective.
\end{cor}

{\sloppy

\begin{proof} Consider first the case that $E$ is an algebraic closure of a finite field~$K$
of characteristic $p > 0$. After increasing $K$ we can assume that $X$ and $f$ are
definable over~$K$. Now $E$ is the union of the intermediate finite fields $F$ with
$K \subseteq F \subseteq E$, and all such~$F$ being perfect, it follows from Corollary~\ref{cor:def fn ACF(p)} that~$f$ maps~${X(F) := X \cap F^n}$ into~$X(F)$, so $f\big(X(F)\big) = X(F)$ by injectivity of $f$ and finiteness of $F$.
Taking the union over all these $F$ we get $f(X) = X$, so we are done for
this particular~$E$. The corollary is equivalent to certain $\mathcal L_{\operatorname{Ri}}$-sentences  holding in all algebraically closed fields; we have shown
these sentences hold in all algebraic closures of finite fields. Therefore they
hold in all algebraically closed fields.
\end{proof}

}
 
\subsection*{Real closed fields}
In this subsection we view ordered fields as structures in the language $\mathcal L_{\operatorname{OR}}$, and we let $K$, $E$ and $F$ denote ordered fields.
We assume here that the reader is familiar with the basic 
{\em algebraic\/} facts concerning real closed fields as exposed in Section~\ref{sec:valued ordered fields}, partly
based on \cite[Chapter~XI]{Lang}.
The models of 
\begin{align*}
\operatorname{RCF}  := \operatorname{OFl}\ \cup\ \big\{ &\, \forall x \exists y\, ( x\geq 0 \to x=y^2 ), \\
& \,\forall u_1\cdots \forall u_{2n+1} \exists x \,( x^{2n+1} + u_{1} x^{2n} + \cdots + u_{2n+1} = 0 ) : n\geq 1 \big\},
\end{align*}
are the real closed ordered fields.
We have an analogue of Theorem~\ref{thm:ACF QE}: 

\nomenclature[Bi12]{$\operatorname{RCF}$}{theory of real closed ordered fields in the language $\mathcal L_{\operatorname{OR}}$}

\begin{theorem}\label{thm:RCF QE}
$\operatorname{RCF}$ has \textup{QE.}
\end{theorem}
\begin{proof}
We use Corollary~\ref{prop:QE test, 2}.
Every ordered integral domain has an $\operatorname{RCF}$-closure, namely the
real closure of its ordered fraction field.
Let $E$ and $F$ be real closed and $K$ be a real closed
ordered subfield of both $E$ and $F$. Suppose  $K\neq E$  and~$F$ is $|K|^+$-saturated. It is enough to show that some ordered subfield of $E$ properly containing $K$ embeds
over $K$ into $F$.
Take any $a\in E\setminus K$. Then~$a$ is transcendental over~$K$ by Corollary~\ref{cor:AS}, and realizes a certain cut $A$
in $K$. Saturation gives an element $b\in F$ realizing the same
cut $A$. Then $b$ is transcendental over $K$, so we have a
field embedding $K(a)\to F$  over $K$ sending $a$ to $b$. This embedding is also order preserving by Corollary~\ref{cor:propAS}. 
\end{proof}

\begin{example}
Let $\varphi(x_1,x_2)=\exists y\,( y^2+x_1y+x_2=0)$.
Then $\varphi$ is $\operatorname{RCF}$-equivalent to the quantifier-free
$\mathcal L_{\operatorname{OR}}$-formula $4x_2 \leq x_1^2$.
\end{example}

\noindent
The substructures of real closed ordered fields are exactly the ordered integral domains, so
$\operatorname{Th}(\operatorname{RCF})$ is the model completion of the theory of ordered integral domains.
The ordered integral domain $\Z$ embeds into every model of $\operatorname{RCF}$,
hence $\operatorname{RCF}$ is complete and thus decidable.  

\begin{remark}
If $E$ is real closed, then $E\models x\geq 0 \leftrightarrow \exists y(y^2=x)$, so
the ordering of~$E$ is  definable by an existential $\mathcal L_{\operatorname{R}}$-formula.
Let $\operatorname{RCF}'$ be a set of $\mathcal L_{\operatorname{R}}$-sentences whose models are the real closed fields.
Then $\operatorname{RCF}'$ is model complete, but
does not have~QE: if $\{a\in \R:\ a\ge 0\}$ were definable in the field $\R$ by the quantifier-free $\mathcal{L}_{\operatorname{R}}$-formula~$\varphi(x)$, then $\R\models \varphi(\sqrt{2})$,
so $\R\models \varphi(-\sqrt{2})$, a contradiction.
\end{remark}

\begin{cor}\label{cor:Tarski Principle}
The following are equivalent, for an $\mathcal L_{\operatorname{OR}}$-sentence $\sigma$:
\begin{enumerate}
\item[\textup{(i)}] $\sigma$ holds in the ordered field $\mathbb R$;
\item[\textup{(ii)}] $\sigma$ holds in some real closed ordered field;
\item[\textup{(iii)}] $\sigma$ holds in every real closed ordered field.
\end{enumerate}
\end{cor}

\noindent
In the next three corollaries of Theorem~\ref{thm:RCF QE}, $E$ is real closed with ordered sub\-field~$K$, and
$T=(T_1,\dots,T_n)$ is an $n$-tuple of distinct indeterminates.
First an analogue of part (i) of Corollary~\ref{cor:ACF QE}
for $\operatorname{RCF}$: \index{theorem!Tarski-Seidenberg}

\begin{cor}\label{cor:RCF QE}
A set $X\subseteq E^n$ is $K$-definable in $E$ iff $X$ is a 
boolean
combination  of   sets of the form
$\big\{a \in  E^n : P(a) \geq 0\big\}$ where $P \in  K[T]$.
\textup{(}In particular, if~$K$ is real closed, then the $K$-definable subsets of $E$ are exactly the finite unions
of singletons and intervals $(a,b)$ where $a<b$ are in $K\cup\{-\infty,\infty\}$.\textup{)}
\end{cor}

\begin{exampleNumbered}\label{ex:RCF QE}
The set $X=\{r\in\R:r\leq \pi\}$ is $\R$-definable in $\R$ but not $0$-definable in~$\R$.
(Take $E=\R$ and $K=\Q^{\operatorname{rc}}\subseteq\R$ in the previous corollary.) 
\end{exampleNumbered}

\begin{prop}\label{prop:dcl RCF} The definably closed subsets of $E$ are exactly \textup{(}the underlying sets of\textup{)} the real closed subfields of $E$.
\end{prop}
\begin{proof}
If $K$ is a real closed subfield of $E$, then $K\preceq E$ by QE, hence $K$, as a subset of~$E$, is definably closed in~$E$.
Conversely, suppose~$K$ is a subset of~$E$ and~$K$ is definably closed in $E$. Then 
$K$ is (the underlying set of) a subfield of $E$. The $0$-definable ordering on $E$ guarantees that $K$ is algebraically closed in~$E$ in the field-theoretic sense. Hence~$K$ is a real closed subfield of $E$, by Corollary~\ref{cor:AS}.
\end{proof}

\noindent
In the next corollary $U$ is an indeterminate different from $T_1,\dots,T_n$.

\begin{cor}
Let $X \subseteq E^n$ and $f \colon X \to E$ be $K$-definable in $E$. Then
there exists $P\in K[T,U]^{\neq}$  such that $P\big(x,f(x)\big)=0$ for all $x \in X$.
\end{cor}

{\sloppy

\begin{proof}
We can assume $E$ is $|K|^+$-saturated. Let
$x \in X$. Then $\operatorname{dcl}\!\big(K \cup \{x\}\big)$ equals the real closure of $K(x)$ in $E$, by
Proposition~\ref{prop:dcl RCF}, and contains
$f(x)$. So $Q\big(x,f(x)\big)=0$ for some  $Q\in K[T,U]^{\neq}$.  By saturation we get $Q_1,\dots,Q_m\in K[T,U]^{\neq}$ such that for each $x\in X$ we have $Q_i\big(x,f(x)\big)=0$ for some $i\in\{1,\dots,m\}$; now put $P:=Q_1\cdots Q_m$.
\end{proof}
}

\noindent
In the proof of the corollary above
we tacitly used that for a real closed ordered field extension $F$ of $E$, and $X\subseteq E^n$ and $f\colon X \to E$ that are $K$-definable in $E$, we have a $K$-definable
set $X_F\subseteq F^n$ with $X_F\cap E^n=X$, and an extension
of $f$ to a $K$-definable function $f_F\colon X_F\to F$:
define $X_F$ in $F$ by any formula defining~$X$ in $E$, and likewise with $f_F$, and note that $X_F$, $f_F$ do not
depend on the choice of these formulas.
We finish this section with an application to the 
asymptotics of functions definable in the ordered field~$\R$; in the proof we employ some facts about Hahn fields from Section~\ref{sec:valued ordered fields}. 

\begin{cor}
Let  $f\colon \R^{>}\to \R$ be definable in the ordered field~$\R$.
Then either $f(r)= 0$ for all sufficiently large $r>0$, or there are 
$c\in \R^\times$ and $q\in\Q$ such that $f(r)/cr^q \to 1$ as $r\to\infty$.
\end{cor}
\begin{proof} We can assume $f(r)\ne 0$ for arbitrarily large $r>0$. Then $f(r)\ne 0$ for all sufficiently large $r>0$ by Corollary~\ref{cor:RCF QE}. Consider the real closed ordered field extension $E=\R[[ x^{\Q} ]] $ of $\R$. From $x>\R$ we get
$f_E(x)\ne 0$, so 
$$f_E(x)\ =\ cx^q(1+\delta) \qquad(c\in\R^\times,\ q\in\Q,\ \delta\in E^{\prec 1}).$$ Hence for every
$\varepsilon\in \R^{>}$ we have 
$\left|f_E(x)/cx^q - 1\right|<\varepsilon$, and thus
$\left|f(r)/cr^q - 1\right|<\varepsilon$ for all sufficiently large $r\in \R^{>}$.
\end{proof}

\begin{exampleNumbered}[H\"ormander]\label{ex:Hormander}
Let $P\in \R[T]^{\ne}$ and define $f\colon \R^>\to\R$ by
$$f(r) := \min \big\{ |P(a)|:\  a=(a_1,\dots,a_n)\in\R^n,\ |a_1|+\cdots+|a_n|=r \big\}.$$
Then there are $c\in \R^\times$ and $q\in\Q$ such that $f(r)/cr^q \to 1$ as $r\to\infty$. 
\end{exampleNumbered}

\subsection*{Notes and comments}
Theorem~\ref{thm:ACF QE} was known to Tarski (see \cite[p.~54, Note~16]{Ta}), but
the first published proof was given by Seidenberg \cite[p.~373]{Se}.
Part~(i) of Corollary~\ref{cor:ACF QE} is also due independently to Chevalley~\cite[Th\'eo\-r\`eme~3, Corollaire]{Chevalley55}; see also \cite{Seidenberg69}. A. Robinson~\cite{Robinson57} saw
model-theoretic significance in Hilbert's Nullstellensatz~\cite{Hilbert93}.  
This led him to
the notion of model completeness. 
Corollary~\ref{cor:Lefschetz} may be seen as a model-theoretic formulation of a heuristic principle in algebraic geometry named after Lefschetz~\cite[p.~242f]{Weil}. See~\cite{Seidenberg58} for a discussion of this principle and \cite{BarwiseEklof, Eklof} for other attempts to make it
precise.
Example~\ref{ex:Serre} is from~\cite{Serre09}. Corollaries~\ref{LR1} and ~\ref{cor:def fn ACF(p)} are from~\cite{Lightstone-Robinson} with a
different proof. 
Corollary~\ref{cor:Ax} for endomorphisms of algebraic varieties $X$ is due to Ax~\cite{Ax68,Ax69}, and, independently,
to Grothendieck \cite[\S{}10.4.11]{EGAIV}; 
Borel~\cite{Borel69} gave a different proof of this result for~${K=\C}$. 

Theorem~\ref{thm:RCF QE} is due to Tarski~\cite{Tarski31,Ta}. The proof above was given by A.~Robinson~\cite{Robinson56}, who also used this theorem to derive in a few lines a solution of
Hilbert's 17th Problem, originally solved by Artin~\cite{Artin27}; see \cite{Robinson55a,Robinson55b}. Other proofs of Theorem~\ref{thm:RCF QE}
were given by {\L}ojasiewicz~\cite{Lojasiewicz} and Seidenberg~\cite{Se}. See \cite{vdDries88} for the history of Tarski's theorem and some applications.
A real Nullstellensatz (analogue of Corollary~\ref{cor:ACF QE}(iii)) is in Dubois~\cite{Dubois},
Krivine~\cite{Krivine}, and Risler~\cite{Risler}; see \cite[Section~4.1]{BCR}.
Corollary~\ref{cor:Tarski Principle} is sometimes called the \textit{Tarski Principle.}\/
Example~\ref{ex:Hormander} is from~\cite[Lemma~3.9]{Hormander}.

\section{Structures without the Independence Property}
\label{sec:structures having NIP}

\noindent
In the previous sections we assembled a toolbox that 
 helps in finding an intelligible description of  
 the sets and maps definable in a structure $\mathbf M$. Such a description often leads to useful geometric invariants (dimensions, Euler characteristics, \dots) of these objects. 
 But in this final section of the appendix we take another path by considering a robust dividing line in the realm of all structures discovered by Shelah. The tame side of this dividing line is called:  \textit{not having the independence property}\/~(NIP).
The NIP condition forbids certain combinatorial patterns in the definable binary relations of a structure. (As is often the case in this subject, a ``tame'' property is introduced here as the negation of a ``wild'' one.) 
Motivating NIP is a striking combinatorial dichotomy (Theorem~\ref{thm:SS}), which we prove first.
This dichotomy leads to the definition of a measure of complexity for families of sets called
\textit{VC dimension}\/ and a dual quantity known as \textit{independence dimension.}\/
After discussing these quantities we introduce the NIP property, and show that real closed fields have NIP.
(This fact is used in our proof that $\mathbb T$ has NIP in Section~\ref{sec:embth}.) \index{NIP} \index{structure!independence property}

\subsection*{VC dimension} Given a set $A$ we let $2^A$ be the power set of $A$ and
${A\choose n}$ the set of $n$-element subsets of~$A$.
\textit{Throughout this subsection $X$ is an infinite set and $\mathcal S$ a collection of subsets of~$X$.}\/ 
For $A\subseteq X$ we set $\mathcal S\cap A:=\{S\cap A:S\in\mathcal S\}$, and we say that~$A$ is {\bf shattered} by~$\mathcal S$ if $\mathcal S\cap A=2^A$. \index{shattering}
We have a function $\pi_{\mathcal S}: \N \to \N$ given by
$$\pi_{\mathcal S}(n)\ :=\ \max\left\{\abs{\mathcal S\cap A}:\ A\in{X\choose n}\right\}.$$
Thus $0\leq\pi_{\mathcal S}(n)\leq 2^n$, and
$\pi_{\mathcal S}(n)=2^n$ means: some $A\in {X\choose n}$ is shattered by $\mathcal S$.
If $\pi_{\mathcal S}(n)=2^n$, then $\pi_{\mathcal S}(m)=2^m$ for all $m\leq n$. Here is the promised dichotomy:

\begin{theorem}[Sauer, Shelah]\label{thm:SS} \index{theorem!Sauer-Shelah}
Either $\pi_{\mathcal S}(n)=2^n$ for all $n$, or else there exists~${d\in\N}$ such that $\pi_{\mathcal S}(n)\le n^d$ for all sufficiently large $n$. 
\end{theorem}

\noindent 
Here is the result from finite combinatorics that underlies
this dichotomy:

\begin{lemma}\label{lem:SS} Suppose $|A|=n$, and $\mathcal C\subseteq 2^A$ and $d\in \{0,\dots,n\}$ are such that $|\mathcal C| > \sum_{i<d} {n\choose i}$. Then $A$ has a subset $B$ with $|B|=d$ that is shattered by $\mathcal C$.
\end{lemma}

\noindent
The hypothesis is sharp since the collection of subsets of $A$ of size $<d$ has cardinality equal to the indicated sum of binomial coefficients, and this particular collection violates the conclusion of the lemma.
Let us denote the indicated sum of binomial coefficients by $p_d(n)$. There is clearly a unique polynomial $p_d(X)\in \Q[X]$ (of degree~${d-1}$ if $d\ge 1$, and
$p_0(X)=0$), whose value at $n$ is $p_d(n)$ for $n\ge d$. Note:
$$p_{d-1}(X-1)+p_d(X-1)\ =\ p_d(X) \qquad (d\ge 1).$$

\begin{proof}[Proof of Lemma~\ref{lem:SS}] 
By induction on $n$. The desired result holds trivially for
$d=0$ and $d=n$, so let $0 < d < n$. Pick a point $a\in A$ and set $A':= A\setminus \{a\}$. Also put $C':= C\setminus \{a\}$
for $C\in \mathcal C$, and set $\mathcal{C}':=\{C': C\in \mathcal{C}\}$. Under the map $C\mapsto C'\colon \mathcal{C} \to \mathcal{C}'$ a set $D\in \mathcal{C}'$ has either exactly one preimage or exactly two preimages; in the latter case these two preimages are $D$ and $D\cup\{a\}$. So $\mathcal{C}'=\mathcal{C}_1\cup \mathcal{C}_2$ (disjoint union) where $\mathcal{C}_1$ contains those $D\in \mathcal{C}'$ having one preimage in $\mathcal{C}$, and $\mathcal{C}_2$ those with two preimages. If $|\mathcal{C}'|>p_d(n-1)$, then by the inductive
assumption applied to~$A'$ and~$\mathcal{C}'$ there exists
$B\subseteq A'$ of size $d$ that is shattered by $\mathcal{C}'$
and thus by $\mathcal{C}$. 
So assume $|\mathcal{C}'|\leq p_d(n-1)$. Then
\begin{align*} 
p_d(n-1)+p_{d-1}(n-1)\  =\  p_d(n)\ &<\ |\mathcal{C}|\ =\ |\mathcal{C}_1| + 2|\mathcal{C}_2|\ =\ \big(|\mathcal{C}_1| + |\mathcal{C}_2|\big) + |\mathcal{C}_2|\\ 
&=\ |\mathcal{C}'| + |\mathcal{C}_2|,
\end{align*}
hence $|\mathcal{C}_2|>p_{d-1}(n-1)$, so again by the inductive assumption applied to $A'$ and~$\mathcal{C}_2$ we get 
$B'\subseteq A'$ of size $d-1$ that is shattered by $\mathcal{C}_2$. Since for each $D\in \mathcal{C}_2$ we have
$D\in \mathcal{C}$ and $D\cup\{a\}\in \mathcal{C}$, the set
$B:= B'\cup\{a\}$ is shattered by $\mathcal{C}$.
\end{proof}

\noindent
This yields Theorem~\ref{thm:SS} in a slightly stronger form:

\begin{theorem}\label{thm:SS+} If $d\in \N$ and $\pi_{\mathcal{S}}(d)< 2^d$, then
$\pi_{\mathcal{S}}(n)\le p_d(n)$ for all $n$.
\end{theorem}
\begin{proof} Assume $\pi_{\mathcal{S}}(d)< 2^d$. Then the desired inequality holds for $n<d$, since $p_d(n)=2^n$ for such $n$. Let $n\ge d$ and $A\in {X \choose n}$. If 
$|\mathcal{S}\cap A|>p_d(n)$, then the above lemma would give a set $B\subseteq A$ of size $d$ shattered by $\mathcal{S}$, a contradiction. 
\end{proof}

\noindent 
We define $$\operatorname{VC}(\mathcal S)\ :=\ \sup\big\{n:\pi_{\mathcal S}(n)=2^n\big\}\in\N\cup\{-\infty, \infty\},$$
so $\operatorname{VC}(\mathcal{S}):=-\infty$ iff $\mathcal{S}=\emptyset$, and $\operatorname{VC}(\mathcal S)=\infty$ iff subsets of~$X$ of arbitrarily large finite size are shattered by~$\mathcal S$.
We call $\operatorname{VC}(\mathcal S)$
 the {\bf VC dimension}  of $\mathcal S$ and say that~$\mathcal S$ is a {\bf VC class} if 
$\operatorname{VC}(\mathcal S)<\infty$. The letters $V$ and $C$ here stand for the initials of Vapnik and Chervonenkis, the authors of ~\cite{VC}.
\index{VC!dimension}\index{VC!class}\index{dimension!VC}

\subsection*{Independence dimension}
Let $X$ be a nonempty set.
Given $A_1,\dots,A_n\subseteq X$ we let $S(A_1,\dots,A_n)$ be the set of atoms of the boolean algebra of subsets of $X$ that is generated by the $A_j$: the ``nonempty fields in the Venn diagram of $A_1,\dots,A_n$''; that is, the elements of
$S(A_1,\dots,A_n)\subseteq 2^X$ are the nonempty sets 
$$\bigcap_{i\in I} A_i \cap \bigcap_{i\in [n]\setminus I} X\setminus A_i\qquad\text{where $I\subseteq [n]:=\{1,\dots,n\}$.}$$ 
Thus $S(A_1,\dots, A_n)=\{X\}$ when $n=0$, and always $S(A_1,\dots, A_n)\ne \emptyset$. 

\begin{lemma}\label{lem:atoms BC}
Let $A_1,\dots,A_m,B_1,\dots,B_n\subseteq X$. Then
\begin{enumerate}
\item[\textup{(i)}] $|S(A_1,\dots,A_m,B_1,\dots,B_n)| \leq |S(A_1,\dots,A_m)|\cdot |S(B_1,\dots,B_n)|$, and
\item[\textup{(ii)}] if $A_1,\dots,A_m$ are boolean combinations of $B_1,\dots,B_n$, then
$$|S(A_1,\dots,A_m)|\ \leq\ |S(B_1,\dots,B_n)|.$$
\end{enumerate}
\end{lemma}
\begin{proof}
Part (i) is clear. For (ii) note that under the hypothesis of (ii),
 the boolean algebra of subsets of $X$ generated by the $B_j$ contains that generated by the $A_i$ as a subalgebra, so every atom of the latter is a disjoint union
of atoms of the former. 
\end{proof}

\noindent
We have $1\leq \abs{S(A_1,\dots,A_n)}\leq 2^n$, and we say that the sequence $A_1,\dots,A_n$ is {\bf independent}~(in~$X$) if $\abs{S(A_1,\dots,A_n)}=2^n$, and call $A_1,\dots,A_n$ {\bf dependent}~(in~$X$) otherwise. \index{independent}\index{dependent} Next, let $\mathcal S\subseteq 2^X$, $\mathcal S \ne \emptyset$. Define  
$$\pi^{\mathcal S}(n)\ :=\ \max\big\{\abs{S(A_1,\dots,A_n)}:A_1,\dots,A_n\in\mathcal S \big\}.$$
Note that 
$1\leq \pi^{\mathcal S}(n)\leq 2^n$ for each $n$. 
We say that $\mathcal S$ is {\bf independent} (in $X$) if $\pi^{\mathcal S}(n)=2^n$ for every $n$, that is, if for every $n$ there is an independent sequence of elements of $\mathcal S$ of length $n$. Otherwise, we say that $\mathcal S$ is {\bf dependent} (in $X$). If $\mathcal S$ is dependent, we define the {\bf independence dimension}  of $\mathcal S$ as the largest $n$ such that $\pi^{\mathcal S}(n)=2^n$; notation: $n=\operatorname{IND}(\mathcal S)$. If $\mathcal S$ is independent, we set $\operatorname{IND}(\mathcal S)=\infty$.
If $\mathcal S$ is finite, then clearly $\operatorname{IND}(\mathcal S)\leq \abs{\mathcal S}$. \index{dimension!independence}

\subsection*{VC duality} \index{VC!duality}
Often our collection $\mathcal S$ will be indexed by elements of a parameter
space or index set, and then the set $X$ and this index set play dual roles. To
make this duality explicit, let~$X$ and~$Y$ be infinite sets and $\Phi\subseteq X\times Y$. For $x\in X$,~$y\in Y$,
\begin{align*} \Phi_x\   &:=\ \big\{ y\in Y:\ (x,y) \in\Phi \big\},  \qquad  \Phi_X\ :=\ \{\Phi_x:\ x\in X\}\subseteq 2^Y\\
\Phi^y\   &:=\ \big\{ x\in X:\ (x,y) \in\Phi \big\},  \qquad  \Phi^Y\ :=\ \{\Phi^y:\ y\in Y\}\subseteq 2^X.
\end{align*} 
One verifies easily that, given a finite set $A=\{x_1,\dots, x_n\}\subseteq X$, the map
$$B \mapsto \bigcap_{x\in B} \Phi_x \cap \bigcap_{x\in A\setminus B} Y\setminus \Phi_x\ \colon\ \Phi^Y\cap A \to S(\Phi_{x_1},\dots, \Phi_{x_n})$$
is a bijection. It follows that $\pi_{\Phi^Y} = \pi^{\Phi_X}$.
Thus $\operatorname{VC}(\Phi^Y)=\operatorname{IND}(\Phi_X)$. 
Reversing the role of $X$ and $Y$ also yields $\operatorname{VC}(\Phi_X)=\operatorname{IND}(\Phi^Y)$.

\begin{lemma}\label{lem:dual shatter, 2} If 
$\operatorname{VC}(\Phi^Y)< \infty$, then
$\operatorname{VC}(\Phi_X)<2^{1+\operatorname{VC}(\Phi^Y)}$. 
\end{lemma} 
\begin{proof}
Let $\varepsilon$ range over the  set of functions $[n]\to \{-1,+1\}$.
Suppose  $\operatorname{VC}(\Phi_X)=\operatorname{IND}(\Phi^Y)\geq 2^n$. This yields an independent family 
$\big(\Phi^{y_{\varepsilon}}\big)_\varepsilon$ of $2^n$ subsets of $X$ (with $y_\varepsilon\in Y$ for all $\varepsilon$). So for each $m\in [n]$
 we have 
 $$\bigcap_{\varepsilon(m)=+1} \Phi^{y_\varepsilon} \cap \bigcap_{\varepsilon(m)=-1} (X\setminus\Phi^{y_\varepsilon})\ \neq\ \emptyset,$$ so we can take an element $x(m)$ of the intersection
 on the left. Then
 for each $\varepsilon$ and $m\in [n]$ we have: $x(m)\in\Phi_{y_\varepsilon}\Longleftrightarrow \varepsilon(m)=+1$,
 so $x(1),\dots, x(n)$ are distinct and the set
 $\big\{x(1),\dots,x(n)\big\}$ is shattered by $\Phi^Y$. Hence
 $\operatorname{VC}(\Phi^Y)\geq n$.
 \end{proof}

\noindent
We therefore have the equivalences
$$\text{$\Phi^Y$ is dependent}\ \Leftrightarrow\ 
\text{$\Phi_X$ is a VC class}\ \Leftrightarrow\ 
\text{$\Phi^Y$ is a VC class}\ \Leftrightarrow\ 
\text{$\Phi_X$ is dependent}.$$
We say that $\Phi$ is {\bf dependent} if $\Phi^Y$ is dependent. \index{dependent}
We also set $\pi^\Phi:=\pi^{\Phi^Y}$.

\medskip\noindent
Let $\neg\Phi$ be the relative complement $(X\times Y)\setminus \Phi$ of $\Phi$ in~$X\times Y$. It is clear that $\pi^{\neg\Phi}=\pi^{\Phi}$. 
For
 $\Phi,\Psi\subseteq X\times Y$ we have $\pi^{\Phi\cup\Psi},\pi^{\Phi\cap\Psi} \leq \pi^{\Phi}\cdot\pi^{\Psi}$, by Lemma~\ref{lem:atoms BC}. By Theorem~\ref{thm:SS},
this yields:

\begin{lemma} \label{lem:VC for boolean combinations}
If $\Phi,\Psi\subseteq X\times Y$ are dependent, then so are $\neg\Phi$, $\Phi\cup\Psi$, and~${\Phi\cap\Psi}$.
\end{lemma}

\subsection*{NIP} Let $\mathbf M$ be an $\mathcal L$-structure 
each of whose underlying sets $M_s$ is infinite. Then we define~$\mathbf M$ to have {\bf NIP} (short for:
{\bf the non-independence property}) if every $0$-definable
relation $\Phi\subseteq M_x\times M_y$ with 
finite nonempty multivariables $x$,~$y$ 
is dependent.
Note that if $\mathbf M$ has NIP, then so does $\mathbf M_A$ for each parameter set $A$ in~$\mathbf M$, and so does
every reduct of $\mathbf M$. If
$\mathbf M$ has NIP and  $\mathbf M\equiv\mathbf N$, then
$\mathbf N$ has~NIP. Thus $\mathbf M$ having NIP is really a property of its theory $\operatorname{Th}(\mathbf M)$.

\index{NIP} \index{structure!independence property}\index{structure!NIP} \index{theory!NIP}

Our goal in this subsection is to show that real closed fields have NIP. We obtain this as a consequence of the following:

\begin{prop}\label{prop:Dudley}
Let $\mathcal F$ be an $m$-dimensional real vector space of real-valued functions on an infinite set $X$, and for
each $f\in \mathcal F$, put $$\operatorname{pos}(f):=\big\{x\in X:f(x)>0\big\}.$$ 
Then $\operatorname{pos}(\mathcal F):=\big\{\!\operatorname{pos}(f):f\in\mathcal F\}$ is a VC class 
of VC dimension $m$.
\end{prop}
\begin{proof}
Let $A\in {X\choose m+1}$. The restriction map
$f\mapsto f|A\colon\mathcal F\to\R^A$ is not surjective, 
since
$\dim\R^A=m+1>\dim\mathcal F$. Therefore we can take
a nonzero $w\in\R^A$ that is orthogonal to all restrictions
$f|A$ ($f\in\mathcal F$) with respect to the standard
inner product $\<u,v\>:=\sum_{a\in A} u(a)\cdot v(a)$ on $\R^A$.
Replacing $w$ by $-w$ if necessary we can assume that
$A^+:=\big\{a\in A:w(a)>0\big\}$ is nonempty. 
If there were $f\in\mathcal F$ with ${A^+=A\cap\operatorname{pos}(f)}$,
then we would have $0=\<w,f|A\>=\sum_{a\in A} w(a)\cdot f(a)>0$,
a contradiction. Hence $A$ is not shattered by $\operatorname{pos}(\mathcal F)$.
Thus $\operatorname{VC}\!\big(\!\operatorname{pos}(\mathcal F)\big)\leq m$,
and we leave the proof that equality holds as an exercise.
\end{proof}

{\sloppy

\noindent
The proposition above applies to the vector space of all real polynomial functions on~${X=\R^n}$ (${n\geq 1}$) of degree $\leq d$. 
Together with~\ref{thm:RCF QE} and~\ref{lem:VC for boolean combinations} this yields at once:
}

\begin{cor}\label{cor:RCF NIP}
Real closed fields have \textup{NIP}.
\end{cor}

\subsection*{Notes and comments}
Theorem~\ref{thm:SS} was found independently by Sauer~\cite{Sauer} and Shelah~\cite{Shelah72b}.
The notion of VC class first arose in probability theory in the work of Vapnik and Chervonenkis~\cite{VC}. The independence property was introduced into model theory by Shelah~\cite{Shelah71a}, who also proved in \cite{Shelah71b}
by a curious set-theoretic argument the very useful fact that in the above definition of NIP
it suffices to consider the case that $x$ is a single variable.
Laskowski~\cite{Laskowski} realized the connection between  VC classes  and NIP and 
 gave a more elementary proof of the reduction to a single variable;  see \cite[Chapter~5]{vdDries98} for an exposition.
Proposition~\ref{prop:Dudley} was shown in \cite[Theorem~7.2]{Dudley78} and used in~ 
\cite{StengleYukich} to prove that the ordered field $\R$ has NIP.
More generally, all \textit{o-minimal}\/ structures have NIP~\cite{PillaySteinhorn},
and so do many other structures: for example,    \textit{stable}\/ structures 
(such as algebraically closed fields), ordered abelian groups \cite{GurevichSchmitt},
and the field $\Q_p$ of $p$-adic numbers~\cite{Belair}. 
For more about the significance of NIP in model theory, see~\cite{SimonNIP}.

%% file: mt-bib.tex
\ifbool{PUP}{
\def\bysame{\leavevmode\hbox to3em{\hrulefill}\thinspace}
}{}

\bibliographystyle{amsplain}

%% file: mt-index.tex
%\ifbool{PUP}{\cleartorecto}{\clearpage}

% Redefine \nomgroup
\renewcommand{\nomgroup}[1]{%
 \ifstrequal{#1}{A}{\vspace{1em}\item[]\hspace*{-\leftmargin}%
        \hfill \textbf{\large Algebra} \hfill \hbox{} \vspace{1em}}{} 
 \ifstrequal{#1}{C}{\vspace{1em}\item[]\hspace*{-\leftmargin}%
        \hfill \textbf{\large Combinatorics} \hfill \hbox{} \vspace{1em}}{} 
 \ifstrequal{#1}{E}{\vspace{1em}\item[]\hspace*{-\leftmargin}%
        \hfill \textbf{\large Ordered Sets} \hfill \hbox{} \vspace{1em}}{} 
 \ifstrequal{#1}{G}{\vspace{1em}\item[]\hspace*{-\leftmargin}%
        \hfill \textbf{\large Valued Abelian Groups and Vector Spaces} \hfill \hbox{} \vspace{1em}}{} 
 \ifstrequal{#1}{I}{\vspace{1em}\item[]\hspace*{-\leftmargin}%
        \hfill \textbf{\large Ordered Abelian Groups and Vector Spaces} \hfill \hbox{} \vspace{1em}}{} 
 \ifstrequal{#1}{K}{\vspace{1em}\item[]\hspace*{-\leftmargin}%
        \hfill \textbf{\large Valued Fields} \hfill \hbox{} \vspace{1em}}{} 
 \ifstrequal{#1}{R}{\vspace{1em}\item[]\hspace*{-\leftmargin}%
        \hfill \textbf{\large Asymptotic Relations} \hfill \hbox{} \vspace{1em}}{} 
 \ifstrequal{#1}{M}{\vspace{1em}\item[]\hspace*{-\leftmargin}%
        \hfill \textbf{\large Differential Algebra} \hfill \hbox{} \vspace{1em}}{} 
 \ifstrequal{#1}{O}{\vspace{1em}\item[]\hspace*{-\leftmargin}%
        \hfill \textbf{\large Differential Polynomials} \hfill \hbox{} \vspace{1em}}{} 
 \ifstrequal{#1}{Q}{\vspace{1em}\item[]\hspace*{-\leftmargin}%
        \hfill \textbf{\large Linear Differential Operators} \hfill \hbox{} \vspace{1em}}{} 
 \ifstrequal{#1}{T}{\vspace{1em}\item[]\hspace*{-\leftmargin}%
        \hfill \textbf{\large Asymptotic Couples} \hfill \hbox{} \vspace{1em}}{} 
 \ifstrequal{#1}{V}{\vspace{1em}\item[]\hspace*{-\leftmargin}%
        \hfill \textbf{\large Dominant Quantities} \hfill \hbox{} \vspace{1em}}{} 
 \ifstrequal{#1}{X}{\vspace{1em}\item[]\hspace*{-\leftmargin}%
        \hfill \textbf{\large Eventual Quantities} \hfill \hbox{} \vspace{1em}}{} 
 \ifstrequal{#1}{Y}{\vspace{1em}\item[]\hspace*{-\leftmargin}%
        \hfill \textbf{\large Triangular Automorphisms} \hfill \hbox{} \vspace{1em}}{} 
 \ifstrequal{#1}{Z}{\vspace{1em}\item[]\hspace*{-\leftmargin}%
        \hfill \textbf{\large Asymptotic Fields} \hfill \hbox{} \vspace{1em}}{} 
 \ifstrequal{#1}{B}{\vspace{1em}\item[]\hspace*{-\leftmargin}%
        \hfill \textbf{\large Model Theory} \hfill \hbox{} \vspace{1em}}{} 
}        

%\ifbool{PUP}{
%\markboth{LIST OF SYMBOLS}{LIST OF SYMBOLS}
%}{}

\printnomenclature[10em]

%\ifbool{PUP}{\cleartorecto}{\clearpage}

\printindex